\documentclass[a4paper,openright,openbib,oneside]{report}%
\usepackage{amsfonts}
\usepackage{svsing2e}
\usepackage{amsmath}
\usepackage[a4paper,asymmetric,nomarginpar,verbose,includehead]{geometry}
\usepackage{amssymb}
\usepackage{float}
\usepackage{graphicx}%
\providecommand{\U}[1]{\protect\rule{.1in}{.1in}}
\newtheorem{theorem}{Theorem}

\newtheorem{conclusion}[theorem]{Conclusion}
\newtheorem{condition}[theorem]{Condition}
\newtheorem{conjecture}[theorem]{Conjecture}
\newtheorem{corollary}[theorem]{Corollary}

\newtheorem{definition}[theorem]{Definition}
\newtheorem{example}[theorem]{Example}

\newtheorem{lemma}[theorem]{Lemma}
\newtheorem{notation}[theorem]{Notation}

\newtheorem{proposition}[theorem]{Proposition}
\newtheorem{remark}[theorem]{Remark}

\newtheorem{summary}[theorem]{Summary}
\newenvironment{proof}[1][Proof]{\noindent\textbf{#1.} }{\ \rule{0.5em}{0.5em}}
\makeatletter
\def\@crosshairs{\vbox to0pt{}}
\makeatother
\ifx\pdfoutput\relax\let\pdfoutput=\undefined\fi
\newcount\msipdfoutput
\ifx\pdfoutput\undefined\else
\ifcase\pdfoutput\else
\ifx\paperwidth\undefined\else
\ifdim\paperheight=0pt\relax\else\pdfpageheight\paperheight\fi
\ifdim\paperwidth=0pt\relax\else\pdfpagewidth\paperwidth\fi
\fi\fi\fi
\begin{document}
%
\begin{titlepage}%
%

\title
{A weight function theory of zero order basis function interpolants and smoothers}%
%

\author{\\Phillip Y. Williams\\\\
5/31 Moulden Court, Canberra, ACT 2617, Australia.\\\\
\texttt{phil.y.williams@bigpond.com.au}}%
%

\date{\today}%
%

\maketitle
%

\thanks
{Thanks to my Masters degree supervisors Dr Markus Hegland and Dr Steve Roberts
of the Centre for Mathematics and its Applications (CMA) at the Australian National University.}%
%

\end{titlepage}%
%

\begin{abstract}%

In this document I develop a weight function theory of zero order basis
function interpolants and smoothers.

Chapter 1 Basis functions and data spaces are defined directly using weight
functions. The data spaces are used to formulate the variational problems
which define the interpolants and smoothers discussed in later chapters. The
theory is illustrated using some standard examples of radial basis functions
and a class of weight functions I will call the tensor product extended B-splines.

Chapter 2 The minimal norm interpolant: pointwise convergence of the
interpolant to its data function and orders of convergence. Some classes of
data functions are characterized locally as Sobolev spaces and the results of
several numerical experiments are presented.

In Chapter 3 local interpolation errors are obtained using tempered
distribution Taylor series. Regarding the extended B-spline weight functions,
we obtain orders of convergence that are better than those obtained
previously. We rely on a careful analysis of a remainder of a tempered
distribution Taylor series expansion based on the function \$exp(i(a,x))\$.

In Chapter 4 we derive some local interpolation errors for data functions
which have bounded first derivatives.

In Chapter 5 a new class of weight functions is introduced which I call the
tensor product central difference weight functions. These are based on a 1-dim
central difference operator acting on an \$L\symbol{94}1\$ function and are
closely related to the extended B-splines and have similar properties.
Convolution formulas are derived for the basis function which are based on the
Taylor series and the central difference operator. The theory of this document
is then applied to obtain interpolation convergence results. As with the
extended B-splines the global data functions are characterized locally as
Sobolev spaces.

In Chapter 6 another class central difference weight functions is introduced.
This class uses a full-dimension central difference and includes radial weight
functions. Multiplicative convolution and partial moment formulas are derived
for the basis functions.

In Chapter 7 the Exact smoother: a non-parametric variational smoothing
problem will be studied using the theory of this document with special
interest in its order of pointwise convergence of the smoother to its data
function. This smoothing problem is the minimal norm interpolation problem
stabilized by a smoothing coefficient. The results are obtained using norms
and seminorms based on the smoothing operator.

In Chapter 8 The scalable Approximate smoother: a non-parametric, scalable,
variational smoothing problem will be studied, again with special interest in
its order of pointwise convergence to its data function. We discuss the
\textit{SmoothOperator} software (freeware) package which implements the
Approximate smoother algorithm. It has a full user manual which has several
tutorials and data experiments.

Chapter 9: The goal of this chapter is to use a bilinear form to characterize
the bounded linear functionals on the data spaces generated by the weight
functions used as examples in this document and to explore the properties of
some related operators.

Chapter 10: We derive an upper bound for the derivative of the 1-dimensional
(scaled) hat basis function smoother assuming the data function has a bounded
derivative and sufficiently large support w.r.t. the data region.

In Chapter 11 we work in one dimension only. The local data functions are
assumed to have bounded derivatives on the data region and we consider a
scaled hat basis function. If the basis function has large enough support
w.r.t. the data region then we show the order of convergence of the
interpolant is 1. It seems to me an unpromising approach for higher dimensions
and smaller supports.

Chapter 12: Explicit extension operators based on Wloka but using the
rectangle condition. Not used in previous chapters.%

\end{abstract}%
%

\tableofcontents

\section{Introduction}

\pagenumbering{arabic}In this document I develop a weight function theory of
zero order basis function interpolants and smoothers. Note that the Appendix
of this document contains a list of basic notation, definitions and properties
also used in this document.

A \textbf{region} will be an open connected subset of $\mathbb{R}^{d}$. We use
the `symmetric' Fourier transform i.e. both the transform and it's inverse
have a constant factor $\left(  2\pi\right)  ^{-d/2}$. This will mean that the
basis function formulas may be more complicated. Also, some weight functions
have simple forms at the expense of a more complicated basis function formula.

This document had its genesis in the development of a \textbf{scalable
algorithm for Data Mining} applications. Data Mining is the extraction of
complex information from large databases, often having tens of millions of
records. Scalability means that the time of execution is linearly dependent on
the number of records processed and this is necessary for the algorithms to
have practical execution times. One approach is to develop additive regression
models and these require the approximation of large numbers of data points by
surfaces. Here one is concerned with approximating data by surfaces of the
form $y=f\left(  x\right)  $, where $x\in\mathbb{R}^{d}$, $y\in\mathbb{R}$ and
$d$ is any dimension. Smoothing algorithms are one way of approximating
surfaces and in particular we have decided to use a class of non-parametric
smoothers called basis function smoothers, which solve a variational smoothing
problem over a semi-Hilbert space of continuous \textbf{data functions} and
express the solution in terms of a single \textbf{basis function}.

I started my Masters degree (supervised by Dr. Markus Hegland and Dr. Steve
Roberts at the ANU, Canberra, Australia) searching for a scalable basis
function smoothing algorithm and had the good fortune to devise such an
algorithm (unpublished) by approximating, on a regular grid, the convolution
in Definition 38 of the space $J_{G}$ in Dyn's review article \cite{Dyn89}. In
this document I develop some theoretical tools to construct and analyze this
algorithm for the case of \textbf{zero order} basis functions. Note that in
the document Williams \cite{WilliamsPosOrdSmthV3} the positive order theory is
studied. This theoretical approach applies in any dimension but the smoothing
algorithm is only practical up to about three dimensions. In higher dimensions
the matrices are too large to put into computer memory.

Dr. Hegland was particularly interested in using the tensor product hat
(triangular) function as a basis function. At about this time we had a visit
by the late Professor Will Light who showed me his paper
\cite{LightWayneX98Weight} which defined basis functions in terms of weight
functions using the Fourier transform. Light and Wayne's weight function
properties were designed for positive order basis functions which excluded the
tensor product hat function. They were designed for the well-known `classical'
radial weight functions. I therefore have developed a version of his theory
designed to generate zero order basis functions, including both tensor product
and radial types. This theory is developed in Chapter 1 and requires that the
basis functions have Fourier transforms which can take zero values outside the
origin since this is a property of hat functions.\medskip

\textbf{Chapter by chapter in brief:}

\begin{enumerate}
\item \textbf{Chapter} \ref{Ch_wtfn_basisfn_datasp} \textbf{Weight functions,
data spaces and basis functions} The goal of this chapter is to extend the
theoretical work of Light and Wayne in \cite{LightWayneX98Weight} to allow
classes of \textbf{zero order weight functions} analogous to the positive
order weight functions developed in Chapter 1 of Williams
\cite{WilliamsPosOrdSmthV3}. Both the inner product \textbf{data space} and
the \textbf{basis functions} are defined in terms of the weight function. We
then prove the completeness and smoothness properties of the data space as
well as the continuity and positive definiteness properties of the basis
function. We use as examples the radial weight functions: the thin-plate
splines, the Gaussian and the Sobolev spline. Special attention is paid to a
class of tensor product weight functions we call the \textbf{extended
B-splines} which are the convolutions of hat functions.

\item \textbf{Chapter} \ref{Ch_Interpol} \textbf{The minimum norm
interpolation problem}. Topics include the existence and uniqueness of the
basis function solution, a matrix equation for the solution and the pointwise
convergence to the interpolant to its data function. Pointwise order of
convergence results are derived using two techniques. The first approach is
actually equivalent to using minimal unisolvent data sets of order 1, which
consist of single points and obtains constants which can be easily calculated.
The second approach is much more complex and uses Lagrange interpolation
techniques which correspond to minimal unisolvent data sets of order
$>$
1. Numerical experiments were undertaken using the extended B-splines and
special data functions were constructed for these experiments. This led to the
extended B-splines data functions being characterized locally as Sobolev spaces.

\item \textbf{Chapter} \ref{Ch_interp_err_Taylor_temper_distrib} \textbf{Local
interpolation errors using tempered distribution Taylor series} We will derive
another set of error estimates for the interpolant. Regarding the radial
Sobolev spline weight functions, we obtain orders of convergence that are
better than those obtained previously. We rely on a careful analysis of a
remainder of a tempered distribution Taylor series expansion based on
$exp\left(  i\left(  a,x\right)  \right)  $ and the Fourier transform. This
approach unifies that of Chapter \ref{Ch_Interpol}.

\item \textbf{Chapter} \ref{Ch_interpol_err_H1inf_data} \textbf{Local
interpolation error for 1-dimensional data functions in} $W^{1,\infty}\left(
\Omega\right)  \cap X_{w}^{0}\left(  \Omega\right)  $ Here we work in one
dimension only. The local data functions are assumed to have bounded
derivatives on the data region and we consider a scaled hat basis function
$\Lambda\left(  x./\lambda\right)  $. If the basis function have "large"
support w.r.t. the data region $\Omega$ i.e. if $\operatorname*{diam}%
\Omega\leq\lambda$, then we show the order of convergence of the interpolant
is 1. The bounds of Chapter \ref{Ch_bnd_deriv_hat_smth_large_supp} are used to
obtain realistic estimates for the scaled hat basis function smoother in 1-dimension.

\item \textbf{Chapter} \ref{Ch_cent_diff_wt_fn_ten_prod} \textbf{The}
\textbf{central difference tensor product weight functions} This chapter
introduces a large new class of tensor product weight functions which I call
the \textbf{central difference weight functions} because they are based on
\textbf{1-dimensional central difference operators} and a non-negative
$L^{1}\left(  \mathbb{R}^{1}\right)  $ generating function $q$. They are
closely related to the tensor product extended B-splines. The central
difference basis functions are calculated, first using a
\textbf{multiplicative convolution formula} and then a more convenient
\textbf{partial moment formula}, and smoothness estimates and polynomial upper
and lower bounds are derived for the weight functions and the basis functions.
If $q\in C_{0}^{\infty}$ then the basis function is $C^{\infty}\left(
\mathbb{R}^{d}\setminus0\right)  $ with bounded support. Like the extended
B-splines their data functions are characterized locally as Sobolev spaces.
Pointwise convergence results are derived for the basis function interpolant
using results from Chapter \ref{Ch_Interpol}. These give orders of convergence
identical to those obtained for the extended B-splines. We do not present the
results of any numerical experiments.

\item \textbf{Chapter} \ref{Ch_cent_diff_wt_fn_multivar} \textbf{Central
difference weight functions: }$q$\textbf{\ multivariate} Unlike in the
previous chapter the weight function is calculated using a
\textbf{multivariate central difference operator} and an $L^{1}\left(
\mathbb{R}^{d}\right)  $ generating function. Here $q$ radial implies the
basis function is radial. A formula based on the central difference operator
is derived which expresses the basis function as the sum of a thin-plate
spline $T_{n}$ and several convolutions of this spline with the scaled
generating function $q$, plus a polynomial. Another formula is derived based
on the Taylor series expansion. Formulas are also derived which use the Bessel
and MacDonald's functions.

\item \textbf{Chapter} \ref{Ch_Exact_smth} \textbf{The Exact smoother} (our
terminology) stabilizes the interpolant by adding a smoothing coefficient to
the seminorm functional. Topics include the existence and uniqueness of the
basis function solution, a matrix equation for the solution and the pointwise
convergence of the solution to its data function. Orders of convergence (also
called error orders) are derived which, in special sense, are the same as
those obtained for the interpolant.

\textbf{Note that my smoothing functional is different} from that used in
Narcowich, Ward and Wendland \cite{NarcWardWend2004} - \textbf{I take the
average of the sum of squares}.

\item \textbf{Chapter} \ref{Ch_Approx_smth} \textbf{The scalable\ Approximate
smoother} (our terminology)\textit{\ }turns out to have a discretization
process similar to that described in Garcke and Griebel \cite{GarckGrieb05}.
Topics include the existence and uniqueness of the smoother, matrix equations
for the smoother, and the pointwise convergence of this smoother to the Exact
smoother and to its data function. Orders of pointwise convergence are derived
for the convergence of the Approximate to the Exact smoother which are then
combined with the Exact smoother error formulas of Chapter \ref{Ch_Exact_smth}
to obtain error orders for the Approximate smoother.

\item \textbf{Chapter} \ref{Ch_Xo1/w_and_fnal_on_Xow} \textbf{The spaces
}$X_{1/w}^{0}$\textbf{, }$\widetilde{X}_{1/w}^{0}$\textbf{\ and the bounded
linear functionals on }$X_{w}^{0}$ The goal of this chapter is to use a
bilinear form to characterize the bounded linear functionals on the data space
of Chapter \ref{Ch_wtfn_basisfn_datasp} for the weight functions used as
examples in this document and to explore the properties of some related
operators. A semi-Hilbert space is defined using the reciprocal of the weight
function and isometric isomorphisms are constructed between this space and
$L^{2}$. These isometries and those constructed in Chapter
\ref{Ch_wtfn_basisfn_datasp} between the data space and $L^{2}$ are used to
construct an isometric isomorphism between functional space and the data space.

\item \textbf{Chapter} \ref{Ch_bnd_deriv_hat_smth_large_supp} \textbf{An upper
bound for the derivative of the 1-dimensional hat basis function smoother} We
derive an upper bound for the derivative of the 1-dimensional (scaled) hat
basis function smoother assuming the data function has a bounded derivative
and sufficiently large support w.r.t. the data region.
\end{enumerate}

\textbf{In the APPENDIX:}

\begin{enumerate}
\item[A] \textbf{Chapter} \ref{Ch_Appendx_basic_notation} \textbf{Basic
notation, definitions and symbols} Basic function spaces; Multi-index and
vector notation; Topology; Tempered distributions; Fourier Transforms;
Convolutions; Some simple (combinatorial) binomial sums; Spherical
coordinates; Taylor series expansions with integral remainder; limits and
differentiation under the integral sign etc.

\item[B] \textbf{Chapter} \ref{Ch_quot_sp_reprod_kern} \textbf{Quotient spaces
and reproducing kernels }General quotient space theory for Banach spaces; the
general theory of reproducing kernel (r.k.) Hilbert spaces; the relationship
between restrictions of r.k. Hilbert spaces and quotient spaces.

\item[C] \textbf{Chapter} \ref{Ch_NotesOnSchmeisser} \textbf{Notes on
Schmeisser \cite{Schmeis2006}} This is a very short chapter which just
contains some notes on the 2006 survey concerning Sobolev spaces with
dominating mixed derivatives by Schmeisser \cite{Schmeis2006}.

\item[D] \textbf{Chapter} \ref{Ch_Hm1_proofs_Adams} \textbf{Proofs of the
claims made in Remark \ref{Rem_SobolevSpace2}} This chapter contains the
proofs of the claims made for the Sobolev spaces with dominating mixed
derivatives $W^{m\mathbf{1}}\left(  \Omega\right)  $ etc. in Remark
\ref{Rem_SobolevSpace2}. They are modified results from Adams and Fournier
\cite{AdamFour2003} and Wloka \cite{Wloka87} concerning the Sobolev spaces
$W^{m}\left(  \Omega\right)  $.
\end{enumerate}

\textbf{Chapter by chapter in more detail:}

\subsection{Chapter \ref{Ch_wtfn_basisfn_datasp} \ Weight functions, data
spaces and basis functions}

In his paper \cite{Duchon77} Duchon studied positive order basis function
interpolation. Basically he began by generalizing the functional
$\int\limits_{\mathbb{R}^{d}}%
{\displaystyle\sum\limits_{\left\vert \alpha\right\vert =m}}
\left\vert D^{\alpha}u\right\vert ^{2}$ using a positive weight function $w$
and the Fourier transform to obtain $\int\limits_{\mathbb{R}^{d}}w%
{\displaystyle\sum\limits_{\left\vert \alpha\right\vert =m}}
\left\vert \widehat{D^{m}u}\right\vert ^{2}$, and then devoted the paper to
studying the special case $w=\left\vert \cdot\right\vert ^{2s}$. In
\cite{LightWayneX98Weight} Light and Wayne carried on the study of the
positive order interpolation problem where Duchon left off. Their weight
function theory used a strictly positive weight function to directly define
the basis functions and Hilbert spaces.

In this document I adapt the work of Light and Wayne to the study the much
simpler zero order problem. Here the weight function is only assumed to be
a.e. positive and the conditions placed on the weight functions are expressed
in terms of integrals. The weaker condition on the weight function allows for
basis functions whose Fourier transforms have zeros e.g. the hat or triangular
function. Happily, the zero order theory presented here only requires simple
$L^{1}$ function space theory, the positive order results being considerably
more complex - see the positive order document Williams
\cite{WilliamsPosOrdSmthV3}.

The theory is illustrated using the standard radial basis functions: the
\textbf{shifted thin-plate splines}, the \textbf{Gaussian} and the
\textbf{Sobolev splines} but I am especially interested in basis functions
that are the tensor product of the derivatives of hat (triangle) functions,
denoted $\Lambda$, which I call the \textbf{extended B-splines}. However, as
mentioned previously, the Fourier transform of these basis functions has zero
values on a set of measure zero and the theory has to allow for this. Thus I
begin by using a class of a.e. positive and continuous \textbf{zero order
weight functions} $w$ satisfying \textbf{property W02}: $\int_{\mathbb{R}^{d}%
}\frac{\left\vert x\right\vert ^{2\lambda}}{w\left(  x\right)  }dx<\infty$ for
$0\leq\lambda\leq\kappa\in\mathbb{R}^{1}$, or \textbf{property W03}:
$\int_{\mathbb{R}^{d}}\frac{x^{2\lambda}}{w\left(  x\right)  }dx<\infty$ for
$0\leq\lambda\leq\kappa\in\mathbb{R}^{d}$, to define the \textbf{zero order
data space}
\[
X_{w}^{0}=\left\{  f\in S^{\prime}:\widehat{f}\in L_{loc}^{1}\text{ }and\text{
}\sqrt{w}\widehat{f}\in L^{2}\right\}  .
\]
This turns out to be a \textbf{reproducing kernel Hilbert space} of continuous
functions with norm $\left\Vert f\right\Vert _{w,0}$ and inner product
$\left(  f,g\right)  _{w,0}=\int w\widehat{f}\overline{\widehat{g}}$. It is
shown, for example, that $X_{w}^{0}\hookrightarrow C^{\left(  \left\lfloor
\kappa\right\rfloor \right)  }$ (see Theorem \ref{Thm_X_smooth}) and that
$S\cap X_{w}^{0}$ is dense in $X_{w}^{0}$ .

The weight function is then used to define a zero order basis function $G$ by
the simple Fourier transform formula $\widehat{G}=\frac{1}{w}$. The smoothness
obtained is $G\in C_{B}^{\left(  \left\lfloor 2\kappa\right\rfloor \right)  }$
and several results are proved regarding the tensor products of weight
functions and basis functions, as well as the products and convolutions of
basis functions.

The \textbf{extended B-splines}\textit{\ }are discussed in detail. It turns
out that the corresponding basis functions are the derivatives of convolutions
of hat (triangle) functions and hence their name.

\subsection{Chapter \ref{Ch_Interpol} The minimal norm interpolant}

In this chapter the minimal norm interpolation problem is solved and several
pointwise convergence results are derived. These are illustrated with
numerical results obtained using several extended B-spline basis functions and
special classes of data functions. In more detail, the functions from the
Hilbert data space $X_{w}^{0}$ are used to define the standard minimal norm
interpolation problem with independent data $X=\left\{  x^{(i)}\right\}
_{i=1}^{N}$ and dependent data $y=\left\{  y_{i}\right\}  _{i=1}^{N}$, the
latter being obtained by evaluating a data function at $X$. This problem is
then shown to have a unique basis function solution of the form $\sum
\limits_{i=1}^{N}\alpha_{i}G\left(  \cdot-x^{(i)}\right)  $ with $\alpha
_{i}\in\mathbb{C}$. The standard matrix equation for the $\alpha_{i}$ is then derived.

We will then consider several categories of estimates for the pointwise
convergence of the interpolant to its data function when the data is confined
to a bounded data region. In all these results the order of convergence
appears as the power of the \textbf{maximum cavity radius }i.e. the radius of
the largest ball centered in the data region that can be fitted between the
$X$ data points. The convergence results of this document can be divided into
those which use Taylor series expansions and Lagrange interpolation theory
i.e. use non-trivial minimal unisolvent subsets of $X$ (Definition
\ref{Def_unisolv}), and those which use Riesz representers. The non-unisolvent
proofs are \textit{much} simpler than the unisolvency-based results and if the
data functions are chosen appropriately the constants can be
calculated.\medskip

\fbox{\textbf{Type 1} interpolation error estimates} Neither the unisolvency
property nor a Taylor series expansion is used. Instead we use the Riesz
representer $R_{x}$ of the evaluation functional. Suppose that the weight
function has property W02 for some $\kappa\geq0$, that the data $X$ is
contained in a closed bounded infinite data set $K$, and that the basis
function $G$ satisfies the smoothness estimate
\begin{equation}
\left\vert G\left(  0\right)  -\operatorname{Re}G\left(  x\right)  \right\vert
\leq C_{G}\left\vert x\right\vert ^{2s},\text{\quad}\left\vert x\right\vert
<h_{G},\label{1.9}%
\end{equation}

for some $h_{G}>0$. Then in Theorem \ref{Thm_|f(x)-f(y)|_inequal_2} it is
shown that the interpolant $\mathcal{I}_{X}f$ of a data function $f\in
X_{w}^{0}$ satisfies%
\begin{equation}
\left\vert f\left(  x\right)  -\mathcal{I}_{X}f\left(  x\right)  \right\vert
\leq k_{G}\left\Vert f\right\Vert _{w,0}\left(  h_{X,K}\right)  ^{s}%
,\text{\quad}x\in K,\label{1.51}%
\end{equation}

\begin{enumerate}
\item when $h_{X,K}=\sup\limits_{x\in K}\operatorname*{dist}\left(
x,X\right)  \leq h_{G}$ and $k_{G}=\left(  2\pi\right)  ^{-\frac{d}{4}}%
\sqrt{2C_{G}}$. This implies an order of convergence of at least $s$.
\end{enumerate}

These results are summarized in Table \ref{Tbl_intro_NonUnisolvTyp1Converg} (a
copy of Table \ref{Tbl_NonUnisolvTyp1Converg}).%

\begin{table}[htbp] \centering
$%
\begin{tabular}
[c]{|c|c||c|c|c|}\hline
\multicolumn{5}{|c|}{\textbf{Type 1} Interpolant error estimates.}\\
\multicolumn{5}{|c|}{Smoothness condition on basis function near
origin.}\\\hline
& Parameter & Converg. &  & \\
Weight function & constraints & order $s$ & $C_{G}$ & $h_{G}$\\\hline\hline
\multicolumn{1}{|l|}{Sobolev splines} & $\frac{1}{2}<v-\frac{d}{2}\leq1$ &
$\frac{1}{2}$ & \multicolumn{1}{|l|}{$\frac{\left\Vert \rho^{v-\frac{d}{2}%
}K_{1-\left(  v-\frac{d}{2}\right)  }\right\Vert _{\infty}}{2^{v-1}%
\Gamma\left(  v\right)  }$} & $\infty$\\\cline{2-4}%
$\left(  v>d/2\right)  $ & $v-\frac{d}{2}>1$ & $1$ &
\multicolumn{1}{|l|}{$\frac{\left\Vert D^{2}\widetilde{K}_{v-\frac{d}{2}%
}\right\Vert _{\infty}}{2^{v}\Gamma\left(  v\right)  }$} & $\infty$\\\hline
\multicolumn{1}{|l|}{Shifted thin-plate} & - & $1$ & \multicolumn{1}{|l|}{eq.
(\ref{1.70})} & $\infty$\\
$\left(  -d/2<v<0\right)  $ &  &  &  & \\\hline
\multicolumn{1}{|l|}{Gaussian} & - & $1$ & \multicolumn{1}{|l|}{$2e^{-3/2}$} &
$\infty$\\\hline
\multicolumn{1}{|l|}{Extended B-spline} & - & $\frac{1}{2}$ &
\multicolumn{1}{|l|}{$G_{1}\left(  0\right)  ^{d-1}\left\Vert DG_{1}%
\right\Vert _{\infty}\sqrt{d}$ $^{\left(  1\right)  }$} & $\infty
$\\\hline\hline
\multicolumn{5}{|l|}{$^{\left(  1\right)  }${\small \ }$G_{1}${\small \ is the
univariate basis function used to form the tensor product.}}\\\hline
\end{tabular}
$\caption{}\label{Tbl_intro_NonUnisolvTyp1Converg}%
\end{table}%
\smallskip

\fbox{\textbf{Type 2} interpolation error estimates} Neither the unisolvency
property nor a Taylor series expansion is used. Instead we use the Riesz
representer $R_{x}$ of the evaluation functional. Here we avoid making any
assumptions about $G$ and instead assume $\kappa\geq1$ and use the properties
of the Riesz representer $R_{x}$ of the evaluation functional $f\rightarrow
f\left(  x\right)  $. A consequence of this approach is the estimate (Theorem
\ref{Thm_interpol_error_in_terms_of_wt_fn}): when $h_{X,K}<\infty$
\begin{equation}
\left\vert f\left(  x\right)  -\mathcal{I}_{X}f\left(  x\right)  \right\vert
\leq k_{G}\left\Vert f\right\Vert _{w,0}h_{X,K},\text{\quad}x\in K,\text{
}f\in X_{w}^{0},\label{1.35}%
\end{equation}

where%
\[
k_{G}=\left(  2\pi\right)  ^{-d/4}\sqrt{-\left(  \left\vert D\right\vert
^{2}G\right)  \left(  0\right)  }.
\]

so the order of convergence is \textbf{always} at least $1$. Table
\ref{Tbl_intro_NonUnisolvTyp2Conv} (a copy of Table
\ref{Tbl_InterpNonUnisolvTyp2Conv}) summarizes the results for the weight
function examples:%

\begin{table}[htbp] \centering
$%
\begin{tabular}
[c]{|c|c||c|c|}\hline
\multicolumn{4}{|c|}{\textbf{Type 2} Interpolant error estimates.}\\
\multicolumn{4}{|c|}{We assume W02 for $\kappa\geq1$.}\\\hline
Weight function & Parameter & Converg. & $\left(  2\pi\right)  ^{d/4}%
k_{G}/\sqrt{d}$\\
& constraints & order & \\\hline\hline
\multicolumn{1}{|l|}{Sobolev splines} & $v-\frac{d}{2}\geq2$ & $1$ &
$\sqrt{\frac{\Gamma\left(  v-d/2-1\right)  }{2^{d/2+1}\Gamma\left(  v\right)
}}$\\\cline{2-4}\cline{2-4}%
$\left(  v>\frac{d}{2}\right)  $ & $1<v-\frac{d}{2}<2$ & $1$ & $\sqrt
{\frac{\Gamma\left(  v-d/2-1\right)  }{2^{2v-d/2-3}\Gamma\left(  v\right)  }}%
$\\\hline
\multicolumn{1}{|l|}{Shifted thin-plate} & - & $1$ & $\sqrt{-2v}$\\
$\left(  -d/2<v<0\right)  $ &  &  & \\\hline
\multicolumn{1}{|l|}{Gaussian} & - & $1$ & $\sqrt{2}$\\\hline
Extended B-spline & $n\geq2$ & $1$ & $\sqrt{-G_{1}\left(  0\right)
^{d-1}D^{2}G_{1}\left(  0\right)  }$ $^{\left(  1\right)  }$\\
$\left(  1\leq n\leq l\right)  $ &  &  & \\\hline\hline
\multicolumn{4}{|l|}{$^{\left(  1\right)  }${\small \ }$G_{1}${\small \ is the
univariate basis function used to form the tensor product.}}\\\hline
\end{tabular}
$\caption{}\label{Tbl_intro_NonUnisolvTyp2Conv}%
\end{table}%

Observe that when $1\leq\kappa<2$ the convergence orders obtained are the same
as those obtained for the unisolvency estimates below but that here the
constants are easily calculated.\medskip

\fbox{Interpolation error estimates using \textbf{unisolvency} and
\textbf{multipoint Taylor series}} The data region $\Omega$ satisfies the cone
condition and minimal unisolvent sets of $X$ order $\left\lfloor
\kappa\right\rfloor $ are used. Note that we do not\textit{\ }assume \textit{a
priori} that $X$ is unisolvent. In this case it follows from Theorem
\ref{Thm_converg_interpol_ord_gte_1} that there exist constants $h_{\Omega
,\kappa},k_{G}>0$ such that%
\begin{equation}
\left\vert f\left(  x\right)  -\mathcal{I}_{X}f\left(  x\right)  \right\vert
\leq k_{G}\left\Vert f\right\Vert _{w,0}\left(  h_{X,\Omega}\right)
^{\left\lfloor \kappa\right\rfloor },\quad x\in\overline{\Omega},\text{ }f\in
X_{w}^{0},\label{1.55}%
\end{equation}

when $h_{X,\Omega}=\sup\limits_{\omega\in\Omega}\operatorname*{dist}\left(
\omega,X\right)  <h_{\Omega,\kappa}$. Table \ref{Tbl_intro_UnisolvConverg} (a
copy of Table \ref{Tbl_UnisolvConverg}) summarizes our results:%

\begin{table}[htbp] \centering
$%
\begin{tabular}
[c]{|c|c|c|}\hline
\multicolumn{3}{|c|}{Interpolant error estimates.}\\
\multicolumn{3}{|c|}{Uses Lagrange interpolation/Taylor series; W02 for
$\kappa\geq1$.}\\\hline
& Parameter & Convergence\\
Weight function & constraints & orders $\left(  \left\lfloor \kappa
\right\rfloor \right)  $\\\hline\hline
\multicolumn{1}{|l|}{Sobolev splines} & \multicolumn{1}{|l|}{$v-\frac{d}%
{2}=2,3,4,\ldots$} & \multicolumn{1}{|l|}{$\left\lfloor v-\frac{d}%
{2}\right\rfloor -1$}\\\cline{2-3}\cline{2-3}%
\multicolumn{1}{|l|}{\quad$\left(  v>\frac{d}{2}\right)  $} &
\multicolumn{1}{|l|}{$v-\frac{d}{2}>1,$ $v-\frac{d}{2}\notin\mathbb{Z}_{+}$} &
\multicolumn{1}{|l|}{$\left\lfloor v-\frac{d}{2}\right\rfloor $}\\\hline
\multicolumn{1}{|l|}{Shifted thin-plate spline} & - &
\multicolumn{1}{|l|}{$1,2,3,4,\ldots$}\\\hline
\multicolumn{1}{|l|}{Gaussian} & - & \multicolumn{1}{|l|}{$1,2,3,4,\ldots$%
}\\\hline
\multicolumn{1}{|l|}{Extended B-spline} & $n\geq2$ &
\multicolumn{1}{|l|}{$n-1$}\\
\multicolumn{1}{|l|}{\quad$\left(  1\leq n\leq l\right)  $} &  &
\multicolumn{1}{|l|}{}\\\hline
\end{tabular}
$\caption{}\label{Tbl_intro_UnisolvConverg}%
\end{table}%

Observe that when $1\leq\kappa<2$ the convergence orders obtained are the same
as those obtained for the Type 2 estimates above but that in the Type 2 case
the constants are easily calculated.

\textbf{Numerical results are only presented for the Type 1 and 2 formulas in
1-dimension}. Numerical results are presented which illustrate the convergence
of the interpolant to its data function. We will only be interested in the
convergence of the interpolant to it's data function and not in the
algorithm's performance as an interpolant. Only the extended B-splines will be
considered and we will also restrict ourselves to one dimension so that the
data density parameter $h_{X,\Omega}$ can be easily calculated.

In \textbf{Subsection} \ref{Sect_local_data_space} we \textbf{characterize of
the restriction spaces }$X_{w}^{0}\left(  \Omega\right)  $\textbf{\ for
several classes of weight function, especially tensor products}. If we can
calculate the data function norm we can calculate the error estimate. All the
extended B-spline basis weight functions have a power of $\sin^{2}x$ in the
denominator and so we have derived special classes of data functions for which
the data function norm can be calculated. The derivation of these special
classes of data functions led to the characterization of the restriction
spaces $X_{w}^{0}\left(  \Omega\right)  $ for various classes of weight
function. For example, Theorem \ref{Thm_int_Xow(O)_eq_Hn(O)_dim1} shows that
the restrictions of the B-spline data functions are a member of the class of
\textbf{Sobolev spaces with dominating mixed derivatives}:
\begin{equation}
W^{m\mathbf{1}}\left(  \Omega\right)  =\left\{  u\in L^{2}\left(
\Omega\right)  :D^{\alpha}u\in L^{2}\left(  \Omega\right)  \text{ }for\text{
}\alpha\leq m\mathbf{1}\right\}  ,\text{\quad}m=1,2,3,\ldots.\label{1.062}%
\end{equation}

It is easy to construct functions that are in these spaces. A more general
result is derived which can be applied to \textbf{tensor product weight
functions} and specifically to the \textbf{central difference weight
functions} of Chapter \ref{Ch_cent_diff_wt_fn_ten_prod}.

As expected interpolant instability is evident and because our error estimates
assume an infinite precision we filter the error to remove spikes which are a
manifestation of the instability.

\subsection{Chapter \ref{Ch_interp_err_Taylor_temper_distrib} Local
interpolation errors using tempered distribution Taylor series}

In this chapter we will derive another set of error estimates for the
interpolant which are then applied to the Sobolev splines and the extended
B-splines. Regarding the Sobolev splines, for $\kappa<1$ we obtain orders of
convergence that are better than those obtained previously: compare Table
\ref{Tbl_NonUnisolvInterpolConvergRevisit_intro_copy} with Table
\ref{Tbl_NonUnisolvTyp1Converg} above.

We start by deriving a \textbf{tempered distribution Taylor series} expansion
with a Fourier transform remainder (\ref{a1.55}, \ref{a50.5}) using a
1-dimensional Taylor series expansion applied to $e^{ia\xi}$. Considering the
cases $\kappa<1$ and $\kappa\geq1$ separately the remainder is then estimated
in terms of the Fourier transform norm $\left\Vert f\right\Vert _{w,0}$ and a
factor involving $1/\sqrt{w}$. In Section \ref{Sect_rem_estim_radial_wt_fn}
special remainder estimates are obtained for radial weight functions.

In the next two sections global and local orders of convergence formulas for
the interpolant are derived with separate sections being devoted to the cases
$\kappa<1$ and $\kappa\geq1$ respectively. Unlike the case $\kappa\geq1$, the
case $\kappa<1$ does not make (explicit) use of unisolvency. These results are
summarize in Table \ref{Tbl_NonUnisolvInterpolConvergRevisit_intro_copy} for
$\kappa<1$ and Table \ref{Tbl_UnisolvInterpolConvergRevisit_intro_copy} for
$\kappa\geq1$. Finally, some local interpolation error estimates are derived.

This is a copy of Table \ref{Tbl_NonUnisolvConverg_revisit}:%

\begin{table}[htbp] \centering
$%
\begin{tabular}
[c]{|c|c||c|c|c|}\hline
\multicolumn{5}{|c|}{Interpolant error estimates, W02 or W03 for
$\underline{\kappa}<1$.}\\
\multicolumn{5}{|c|}{Technique: tempered distribution Taylor series.}\\\hline
Weight function & Parameter & Converg. &  & \\
& constraints & order & $c_{w}$ & $h_{w}$\\\hline\hline
\multicolumn{1}{|l|}{Sobolev spline} & $v-\frac{d}{2}\leq1$ & $<v-\frac{d}{2}$
& \multicolumn{1}{|l|}{eq.\thinspace\ref{a9.6}} & $\infty$\\
\quad$\left(  v>d/2\right)  $ &  &  &  & \\\hline
\multicolumn{1}{|l|}{Extended B-spline} & $n=1$ & $1/2$ &
\multicolumn{1}{|l|}{eq.\thinspace\ref{a9.7}} & $\infty$\\\hline
\end{tabular}
$\caption{}\label{Tbl_NonUnisolvInterpolConvergRevisit_intro_copy}%
\end{table}%
\medskip

A copy of Table \ref{Tbl_UnisolvConverg_revisit}:%

\begin{table}[htbp] \centering
$%
\begin{tabular}
[c]{|c|c||c|c|c|}\hline
\multicolumn{5}{|c|}{Interpolant error estimates assuming W02 or W03 for
$\kappa\geq1$.}\\
\multicolumn{5}{|c|}{Techniques: unisolvency, multipoint distribution Taylor
series.}\\\hline
Weight function & Parameter & Converg. &  & \\
& constraints & order & $c_{w}$ & $h_{w}$\\\hline\hline
\multicolumn{1}{|l|}{Sobolev spline (W02)} &  & $<v-d/2$ &
\multicolumn{1}{|l|}{\ref{a110}} & $\infty$\\\cline{2-3}%
\quad$\left(  v>d/2\right)  $ &  &  &  & \\\hline
\multicolumn{1}{|l|}{Extended B-spline (W03)} & $n\geq2$ & $n-1/2$ &  &
$\infty$\\\hline
$\left(  n\leq l\right)  $ &  &  &  & \\\hline
\end{tabular}
$\caption{}\label{Tbl_UnisolvInterpolConvergRevisit_intro_copy}%
\end{table}%
\medskip

Observe that the convergence order estimate for the Sobolev spline in Table
\ref{Tbl_NonUnisolvInterpolConvergRevisit_intro_copy} is better than that
given in Table \ref{Tbl_intro_NonUnisolvTyp1Converg}.

\subsection{Chapter \ref{Ch_interpol_err_H1inf_data} Local interpolation error
for data functions in $W^{1,\infty}\left(  \Omega\right)  \cap X_{w}%
^{0}\left(  \Omega\right)  $}

Here the local data functions are assumed to have bounded derivatives on a
bounded data region $\Omega$ i.e. functions in $W^{1,\infty}\left(
\Omega\right)  \cap X_{w}^{0}\left(  \Omega\right)  $ where $X_{w}^{0}$ is a
global data space.

In Section \ref{Lem_interpol_estim_dim1_W1inf_Xow} we will consider
1-dimensional estimates and the main result is Theorem
\ref{Thm_err_interpol_H1inf_data}: We will choose our data functions $f_{d}\in
W^{1,\infty}\left(  \Omega\right)  \cap X_{w}^{0}\left(  \Omega\right)  $ and
assume that there exist constants $c_{0},c_{1}>0$ such that the basis function
interpolant $\mathcal{I}_{X}f_{d}$ satisfies%
\[
\left\vert D\mathcal{I}_{X}f_{d}\left(  x\right)  \right\vert \leq
c_{0}\left\Vert f_{d}\right\Vert _{\infty,\Omega}+c_{1}\left\Vert
Df_{d}\right\Vert _{\infty,\Omega},\quad x\in\Omega,\text{ }f_{d}\in
W^{1,\infty}\left(  \Omega\right)  \cap X_{w}^{0}\left(  \Omega\right)  .
\]

Then we have the following pointwise \textbf{order 1 spherical cavity} error
estimate%
\[
\left\vert \mathcal{I}_{X}f_{d}\left(  x\right)  -f_{d}\left(  x\right)
\right\vert \leq\left(  c_{0}\left\Vert f_{d}\right\Vert _{\infty,\Omega
}+\left(  1+c_{1}\right)  \left\Vert Df_{d}\right\Vert _{\infty,\Omega
}\right)  h_{\Omega,X},\quad x\in\Omega,
\]

An example is provided by Corollary \ref{vCor_bound_deriv_scal_hat_interpol}
of the Appendix where it is shown that in \textbf{one dimension} the
interpolant generated by the \textbf{scaled hat basis function} $\Lambda
\left(  x/\lambda\right)  $ satisfies the local estimate:
\[
\left\Vert D\mathcal{I}_{X}f_{d}\right\Vert _{\infty;\Omega}\leq\frac
{2}{\lambda}\left\Vert f_{d}\right\Vert _{\infty;\Omega}+\left\Vert
Df_{d}\right\Vert _{\infty;\Omega},\quad f_{d}\in W^{1,\infty}\left(
\Omega\right)  \cap X_{w}^{0}\left(  \Omega\right)  ,
\]

where the scaled hat basis function (extended B-spline $n=1$) is assumed to
have large support w.r.t. the data region i.e. $\operatorname*{diam}\Omega
\geq\lambda$. \textbf{This order 1 convergence is twice the convergence order
of 1/2 obtained above in the previous chapter} but it is only valid for a
subspace of $X_{w}^{0}\left(  \Omega\right)  $. In Section
\ref{Lem_interpol_estim_multivar_W1inf_Xow} an analogous multivariate error
estimate is proven where the interpolant is assumed to satisfy the inequality%
\[
\left\Vert D^{\alpha}\mathcal{I}_{X}f_{d}\right\Vert _{W^{\infty}\left(
\Omega\right)  }\leq c\left\Vert f_{d}\right\Vert _{W^{1,\infty}\left(
\Omega\right)  },\quad x\in\Omega,\text{ }f_{d}\in W^{1,\infty}\left(
\Omega\right)  \cap X_{w}^{0}\left(  \Omega\right)  ,\text{ }\left\vert
\alpha\right\vert =1.
\]

I have not been able to prove that the scaled tensor product hat function with
large support satisfies this inequality.

Note that in 1-dimension data even smoothed functions which have infinite
derivatives show instabilities - subsection
\ref{SbSect_notC1b_data_fns_H1inf_1dim}.

\subsection{Chapter \ref{Ch_cent_diff_wt_fn_ten_prod} The central difference
tensor product weight functions}

To understand the theory of interpolants and smoothers presented in this
document it is not necessary to read this chapter which introduces the
\textit{central difference weight functions}. However, we note that Section
\ref{Sect_CntDifWtFn_DataFuncs} contains local data space results specifically
designed for tensor product weight functions and that central difference
weight functions are used as examples in the chapters dealing with smoothers.

This chapter introduces a large new class of weight functions based on central
difference operators. The basis functions are calculated and smoothness
estimates and upper and lower bounds are derived for the weight functions and
the basis functions.

Like the extended B-splines, their data functions are characterized locally as
Sobolev spaces.

Pointwise convergence results are derived for the basis function interpolant
using results from Chapter \ref{Ch_Interpol}. These give orders of convergence
identical to those obtained for the extended B-splines in. This chapter
presents no numerical interpolation experiments.

The central difference weight functions are defined as follows: Suppose $q\in
L^{1}\left(  \mathbb{R}^{1}\right)  $, $q\neq0$, $q\left(  \xi\right)  \geq0$
and $l\geq n\geq1$ are integers. The univariate \textit{central difference
weight function} with parameters $n,l$ is defined by
\[
w\left(  \xi\right)  =\frac{\xi^{2n}}{\Delta_{2l}\widehat{q}\left(
\xi\right)  },\quad\xi\in\mathbb{R}^{1},
\]

where $\Delta_{2l}$ is the central difference operator%
\[
\Delta_{2l}f\left(  \xi\right)  =\sum_{k=-l}^{l}\left(  -1\right)  ^{k}%
\tbinom{2l}{k+l}f\left(  -k\xi\right)  ,\text{\quad}l=1,2,3,\ldots,\text{ }%
\xi\in\mathbb{R}^{1}.
\]

For example, $\Delta_{2}f\left(  \xi\right)  =-\left(  f\left(  \xi\right)
-2f\left(  0\right)  +f\left(  -\xi\right)  \right)  $. The multivariate
central difference weight function is defined by tensor product.

It is shown in Theorem \ref{Thm_cdiffwt_2} that $w$ belongs to the class of
zero order weight functions introduced in Chapter \ref{Ch_wtfn_basisfn_datasp}
for some $\kappa$ iff $\int\limits_{\left\vert \xi\right\vert \geq
R}\left\vert \xi\right\vert ^{2n-1}q\left(  \xi\right)  d\xi<\infty$ for some
$R\geq0$. Here $\kappa$ satisfies $\kappa+1/2<n$. The central difference
weight functions are closely related to the \textit{extended B-splines}
defined by \ref{1.032} and a discussion of their genesis is given in
Subsection \ref{SbSect_cent_Motivation}.

Various bounds are derived for $w$. For example, in Corollary
\ref{Cor_cdiffwt_bnd_on_wt_fn} it is shown that if $\int\nolimits_{\left\vert
t\right\vert \geq R}t^{2l}q\left(  t\right)  dt<\infty$ then for any $r>0$
there exist constants $c_{r},c_{r}^{\prime},k_{r},k_{r}^{\prime}>0$ such that%
\begin{align*}
k_{r}\xi^{2n}  & \leq w\left(  \xi\right)  \leq k_{r}^{\prime}\xi^{2n}%
,\quad\left\vert \xi\right\vert \geq r,\\
\frac{c_{r}}{\xi^{2\left(  l-n\right)  }}  & \leq w\left(  \xi\right)
\leq\frac{c_{r}^{\prime}}{\xi^{2\left(  l-n\right)  }},\quad\left\vert
\xi\right\vert \leq r.
\end{align*}

By Theorems \ref{Thm_G_basis_def_2} and \ref{Thm2_G_basis_def_2_new} the
univariate \textit{central difference basis function} is given by the
\textit{multiplicative convolution formula}%
\[
G_{c}\left(  s\right)  =\tfrac{2^{2\left(  l-n\right)  +1}}{\sqrt{2\pi}}%
\int_{\mathbb{R}^{1}}G_{s}\left(  \frac{2s}{t}\right)  \left\vert t\right\vert
^{2n-1}q\left(  t\right)  dt,\text{\quad}s\in\mathbb{R}^{1},
\]

where $G_{s}$ is the extended B-spline basis function defined by \ref{1.49}.
We show that if $q$ is bounded then $G_{c}\in C_{B}^{2n-2}\left(
\mathbb{R}^{1}\right)  $ and $D^{2n-1}G_{c}\in C^{\left(  0\right)  }\left(
\mathbb{R}^{1}\setminus0\right)  \cap L^{\infty}$. The multivariate basis
function is defined as a tensor product. For $k\leq2n-2$, $D^{k}G_{c}$ is
uniformly Lipschitz continuous of order $1$ (Theorem
\ref{Thm_cdiffbasis_Lips_dim_d}).

Another \textit{multiplicative convolution formula} is given in Theorem
\ref{Thm_cdiffbasis_mult_convol_ql}. This formula does not involve a formula
for the B-spline.

The Theorem \ref{Thm_cdiffbasis_part_moment_formula} and its corollary give
the more convenient \textit{partial moment formulas} for the basis function
and its derivatives that do not involve calculating the extended B-spline
basis functions $G_{s}$ but instead involves the direct calculation of partial
moments of $q$\textbf{. }In fact%
\begin{align*}
G_{c}\left(  x\right)   & =\left\{
\begin{array}
[c]{ll}%
\frac{\left(  -1\right)  ^{n}}{\left(  2n-1\right)  !}\int_{\left\vert
x\right\vert }^{R_{q}l}\left(  s-\left\vert x\right\vert \right)  ^{2n-1}%
q_{l}\left(  s\right)  ds, & \left\vert x\right\vert \leq R_{q}l,\\
0, & \left\vert x\right\vert \geq R_{q}l,
\end{array}
\right. \\
q_{l}\left(  s\right)   & =\sum\limits_{j=-l,j\neq0}^{l}\frac{\left(
-1\right)  ^{j}}{\left\vert j\right\vert }\tbinom{2l}{j+l}q\left(  \frac{s}%
{j}\right)  ,
\end{align*}

where $\operatorname*{supp}q\subseteq\overline{B}_{R_{q}}$. We also show that
$G_{c}\in C_{B}^{2n-2}\left(  \mathbb{R}^{1}\right)  $ and $D^{2n-1}G_{c}\in
C^{\left(  0\right)  }\left(  \mathbb{R}^{1}\setminus0\right)  \cap L^{\infty
}$ without assuming that $q$ is bounded. It is also shown that if
$\operatorname*{supp}q\subseteq\overline{B}_{R_{q}}$, where possibly
$R_{q}=\infty$, it follows that $\operatorname*{supp}G_{c}\subseteq
\overline{B}_{R_{q}l}$ and that $G_{c}$ satisfies the differential equation
\[
D^{2n}G_{c}\left(  s\right)  =\left(  -1\right)  ^{n}\tbinom{2l}{l}\left(
\int q\right)  \delta+\left(  -1\right)  ^{n}\sum\limits_{j=-l,j\neq0}%
^{l}\frac{\left(  -1\right)  ^{j}}{\left\vert j\right\vert }\tbinom{2l}%
{j+l}q\left(  \frac{s}{j}\right)  ,
\]

implying that $D^{2n}G_{c}$ has the same smoothness as $q$ on $\mathbb{R}%
^{d}\setminus0$. In fact if $q\in C_{0}^{\infty}\left(  \mathbb{R}^{1}\right)
$ then $\ G_{c}\in C^{\infty}\left(  \mathbb{R}^{1}\setminus0\right)  \cap
C_{B}^{\left(  2n-2\right)  }\left(  \mathbb{R}^{1}\right)  $ and $G_{c}$ has
bounded support.

In Section \ref{Sect_CntDifWtFn_DataFuncs} we will characterize the
\textit{data function space} locally as a mixed-derivative Sobolev space
\ref{1.062}. A localization result specifically for tensor product weight
functions is derived in Theorem \ref{Thm_data_fn_tensor_prod}. This supplies
important information about the data functions and makes it easy to choose
data functions for numerical experiments.

Results for the pointwise convergence of the interpolant to its data function
on a bounded data domain are derived using results from Chapter
\ref{Ch_wtfn_basisfn_datasp}. These results are summarized here in Table
\ref{Tbl_intro_ConvergCentral} (a copy of Table \ref{Tbl_ConvergCentral}) and
give pointwise orders of convergence identical to those obtained for the
extended B-splines.%

\begin{table}[htbp] \centering
$%
\begin{tabular}
[c]{|c||c||c|c|c|}\hline
\multicolumn{5}{|c|}{\textbf{Interpolant convergence to data function }-
central difference weight functions.}\\\hline
Estimate & Parameter & Converg. &  & \\
name & constraints & order & $\left(  2\pi\right)  ^{d/4}k_{G}$ & $h_{G}%
$\\\hline\hline
\multicolumn{1}{|l|}{Type 1 (smoothness} & $q$ b'nded if $n=1$ & $1/2$ &
$\sqrt{2G_{1}\left(  0\right)  ^{d-1}\left\Vert DG_{1}\right\Vert _{\infty}%
}\sqrt[4]{d}$ $^{\left(  1\right)  }$ & $\infty$\\
\multicolumn{1}{|l|}{\quad constraint on basis fn)} & else $n\geq2$ &  &  &
\\\hline
\multicolumn{1}{|l|}{Type 2 ($\kappa\geq1$)} & $n\geq2$ & $1$ & $\sqrt
{-G_{1}\left(  0\right)  ^{d-1}D^{2}G_{1}\left(  0\right)  }\sqrt{d}$
$^{\left(  1\right)  }$ & $\infty$\\\hline
Unisolvency/Taylor series & $n\geq2$ & $n-1$ & - & $\infty$\\\hline
\multicolumn{1}{|l|}{Tempered distrib Taylor} & $q$ b'nded a.e. & $n-1/2$ &
- & $\infty$\\\hline\hline
\multicolumn{5}{|l|}{$^{\left(  1\right)  }${\small \ }$G_{1}${\small \ is the
univariate basis function used to form the tensor product.}}\\\hline
\end{tabular}
$\caption{}\label{Tbl_intro_ConvergCentral}%
\end{table}%

This chapter does not present the results of any numerical experiments.

\subsection{Chapter \ref{Ch_cent_diff_wt_fn_multivar} Central difference
weight functions: $q$ multivariate}

In the last chapter we generated multivariate central difference weight
functions using the tensor product. In this (incomplete) chapter we generalize
to higher dimensions the 1-dimensional difference operator $\Delta_{2l}$ used
in the one-dimensional case and thus obtain a new class of weight functions
which seems most suited to radial weight functions.

Using the Fourier transform to directly calculate the basis function from the
weight function formula \ref{a969} may be very awkward.

These central difference weight functions are defined as follows: suppose
$q\in L^{1}\left(  \mathbb{R}^{d}\right)  ,q\neq0,q\left(  \xi\right)  \geq0$
and $l,n\geq0$ are integers. Then the multivariate \textit{central difference
weight function} with parameters $n,l$ is defined by
\[
w\left(  \xi\right)  =\frac{\left\vert \xi\right\vert ^{2n}}{\Delta
_{2l}\widehat{q}\left(  \xi\right)  },\quad\xi\in\mathbb{R}^{d},
\]

where $\Delta_{2l}$ is the central difference operator%
\[
\Delta_{2l}f\left(  \xi\right)  =\sum_{k=-l}^{l}\left(  -1\right)  ^{k}%
\tbinom{2l}{k+l}f\left(  -k\xi\right)  ,\text{\quad}l=1,2,3,\ldots;\text{ }%
\xi\in\mathbb{R}^{d}.
\]

If $q$ is radial then $w$ is radial. It is shown in Theorem
\ref{Thm_wt_func_mdim} that $w$ belongs to the class of zero order weight
functions for some $\kappa$ iff $0\leq\kappa<n-\frac{d}{2}<l$ and
$\int_{\left\vert \tau\right\vert \geq R}\left\vert \tau\right\vert
^{2n-d}q\left(  \tau\right)  d\tau<\infty$.

Using the theory of $S_{\emptyset,m}$ tempered distributions (Definition
\ref{Def_So,n} etc.) Theorem \ref{Thm_centdiffbasis_centdiff_formula} shows us
how to use the central difference operator to express the basis function as
the sum of a thin-plate spline $T_{n}$ and several convolutions of this spline
with the scaled generating function $q$:%
\[
G_{c}\left(  x\right)  =\frac{1}{\left(  2\pi\right)  ^{d}c_{n,d}}\left(
\left(  -1\right)  ^{j}\tbinom{2l}{l}\left(  \int q\right)  T_{n}\left(
x\right)  +\sum\limits_{j=-l,\text{ }j\neq0}^{l}\left(  -1\right)  ^{j}%
\tbinom{2l}{j+l}\int T_{n}\left(  x-jy\right)  q\left(  y\right)  dy\right)
+p_{c}\left(  x\right)  ,
\]

the unique $p_{c}$ is a polynomial which satisfies $\left\vert D\right\vert
^{2n}p_{c}=0$. It remains to determine this polynomial by applying conditions
to the function $q$. For $d$ odd, using the $L_{loc}^{1}$ Taylor series
expansion with integral remainder described in Lemma
\ref{Lem_Taylor_estim_loc_C_L1loc}, I have derived Theorem
\ref{Thm_CentDiff_eqn_thinps} which supplies conditions on $q$ under which the
polynomial $p_{c}$ is zero. Here we assume $q$ has bounded support or that $q
$ is radial with an extra condition. I have not done the case where $d$ is
even but have made partial attempt on a more general result which does not
distinguish between odd and even dimensions.

In \textbf{Section} \ref{Sect_Bessel_centdiff_basis_formula} we derive several
formulas for the central difference basis function which involve the Bessel
functions $J_{\frac{d-2}{2}}$ and the MacDonald's functions $K_{n-\frac{d}{2}%
}$. For example Theorem \ref{Thm_cdiffbasis_Bessel} shows that%
\[
G_{c}\left(  x\right)  =\tfrac{1}{\left(  2\pi\right)  ^{d/2}}\int%
\int\limits_{0}^{\infty}\left(  \sum_{j=-l}^{l}\left(  -1\right)  ^{j}%
\tbinom{2l}{j+l}\frac{J_{\frac{d-2}{2}}\left(  \left\vert x-j\tau\right\vert
r\right)  }{\left(  \left\vert x-j\tau\right\vert r\right)  ^{\frac{d-2}{2}}%
}\right)  \frac{dr}{r^{2n-d+1}}q\left(  \tau\right)  d\tau,
\]

and Theorem \ref{Thm_CentDiffBasisFormula_Bessel} claims that%
\[
G_{c}\left(  x\right)  =\tfrac{1}{\left(  2\pi\right)  ^{\frac{d}{2}}%
2^{n-1}\Gamma\left(  n\right)  }\int\lim_{\varepsilon\rightarrow0^{+}}%
\sum_{j=-l}^{l}\left(  -1\right)  ^{j}\tbinom{2l}{j+l}\left\vert
x-j\tau\right\vert ^{n-1}\frac{K_{n-\frac{d}{2}}\left(  \varepsilon\left\vert
x-j\tau\right\vert \right)  }{\varepsilon^{n-\frac{d}{2}}}q\left(
\tau\right)  d\tau.
\]

In \textbf{Section} \ref{Sect_CentDiffBasis_Taylor} we obtain results which
are an extension of the 1-dimensional of Subsection
\ref{SbSect_tempdistrib_1_dim_centdiff_basis} to an arbitrary dimension.

We use the tempered distribution Taylor series expansion introduced in Section
\ref{Sect_Taylor_series_data_fn} and the theory of the Schwartz subspaces
$S_{\emptyset,k}\subset S$ introduced in Definition \ref{Def_So,n} to expand a
convolution of the thin plate spline and the function $q$ and thus prove the
Taylor series basis function formula given in Theorem
\ref{Thm_centdiff_basis_Tn_q_suppbnd}. The formula is true modulo a polynomial
with degree $\leq d-2$. We will make use of the basic Fourier transform
properties given in Definition \ref{Def_Fourier} of the Appendix. The
$L_{loc}^{1}$ Taylor series expansion of Lemma \ref{Lem_Taylor_estim_L1loc_L1}
is then used to obtain the Taylor series integral remainder formula \ref{a017}
for the basis function which is again true modulo a polynomial with degree
$\leq d-2$.

In \textbf{Section} \ref{Sect_LocDataSpace} I have made a start at
characterizing the data space locally.

I have not yet derived any bounds for $w$ or any orders of convergence for the
smoother. I have begun studying the local structure of the data space.

\subsection{Chapter \ref{Ch_Exact_smth} The Exact smoother}

We call this well-known parameter-stabilized basis function interpolant the
\textit{Exact smoother} because the smoother studied in the next chapter
approximates it. We will assume the basis function is real-valued. The Exact
smoother minimizes the functional%
\begin{equation}
\rho\left\Vert f\right\Vert _{w,0}^{2}+\frac{1}{N}\sum_{i=1}^{N}\left\vert
f(x^{(i)})-y_{i}\right\vert ^{2},\quad f\in X_{w}^{0},\label{1.8}%
\end{equation}

over the data (function) space $X_{w}^{0}$ where $\rho>0$ is termed the
smoothing coefficient. \textbf{Note} that this smoothing functional is
different from that used in Narcowich, Ward and Wendland
\cite{NarcWardWend2004} which is%
\[
\rho\left\Vert f\right\Vert _{w,0}^{2}+\sum_{i=1}^{N}\left\vert f(x^{(i)}%
)-y_{i}\right\vert ^{2}.
\]

To obtain their error estimates you simply replace $\rho N$ by $\rho$ in our
error formulas.

We show that this problem, like the interpolation problem, has a unique basis
function solution $s_{e}=\mathcal{S}_{X}f$ in the space $W_{G,X}$. The finite
dimensionality of the solution allows us to derive matrix equations for the
coefficients $\alpha_{i}$ of the data-translated basis functions.

Following the example of the interpolation semi-inner product $\left(
f-\mathcal{I}_{X}f,g\right)  _{w,0}$ of Definition \ref{Def_interp_seminorm}
of \ref{Ch_Interpol} and it's sequel I will define the semi-inner product
$\left(  \mathcal{S}_{X}f,g\right)  _{w,0}$ and inner product $\left(
f-\mathcal{S}_{X}f,g\right)  _{w,0}$ and use these to study the value and the
error of the Exact smoother respectively. These expressions are useful because
$\left(  \mathcal{S}_{X}f,R_{x}\right)  _{w,0}=\left(  \mathcal{S}%
_{X}f\right)  \left(  x\right)  $ and $\left(  f-\mathcal{S}_{X}%
f,R_{x}\right)  _{w,0}=\left(  f-\mathcal{S}_{X}f\right)  \left(  x\right)  $
and because of the power of Hilbert space theory. An outcome of this approach
are the Exact smoother bounds of Theorem \ref{Thm_bound_Exact_smth}:
\[
\left\vert s_{e}\left(  x\right)  \right\vert \leq\left\Vert f\right\Vert
_{w,0}\sqrt{R_{0}\left(  0\right)  }\min\left\{  1,\frac{R_{0}\left(
0\right)  }{\rho}\right\}  ,
\]

and the convergence estimates of Corollary
\ref{Cor_Thm_Rx(x)minus(SRx)(x)_basis_fn_bound}:%
\begin{equation}
\left\vert f\left(  x\right)  -s_{e}\left(  x\right)  \right\vert
\leq\left\Vert f\right\Vert _{w,0}\left(  \sqrt{\rho N}+\frac{\sqrt{2}%
}{\left(  2\pi\right)  ^{\frac{d}{4}}}\sqrt{G\left(  0\right)  -G\left(
x-x^{\left(  k\right)  }\right)  }\right)  .\label{2.22}%
\end{equation}

When $\rho=0$ these estimates correspond to interpolation results of The Exact
smoother convergence orders and the constants are the same as those for the
interpolation case which are given in the interpolation tables.\medskip

\fbox{Type 1 error estimates} in Subsection \ref{SbSect_ex_smth_typ1_error}.
When the weight function has property W02 for $\kappa\geq0$ it is assumed that
the basis function satisfies an inequality of the form \ref{1.9} and that the
data region is a closed bounded infinite set $K$. In this case it is shown in
Corollary \ref{Cor_OrdConvergExactSmth_k=0_Fd} that the Exact smoother $s_{e}$
of the data function $f$ satisfies the error estimate%
\begin{equation}
\left\vert f\left(  x\right)  -s_{e}\left(  x\right)  \right\vert
\leq\left\Vert f\right\Vert _{w,0}\left(  \sqrt{\rho N}+k_{G}\left(
h_{X,K}\right)  ^{s}\right)  ,\quad x\in K,\label{1.58}%
\end{equation}

when $h_{X,K}=\sup\limits_{x\in K}\operatorname*{dist}\left(  x,X\right)  \leq
h_{G}$ and $k_{G}=\left(  2\pi\right)  ^{-d/4}\sqrt{2C_{G}}$.\medskip

\fbox{Type 2 error estimates} in Subsection \ref{SbSect_ex_smth_typ2_error}.
If it only assumed that $\kappa\geq1$ then it is shown in Theorem
\ref{Thm_G(x)minusG(0)_bound} that%
\begin{equation}
\left\vert f\left(  x\right)  -s_{e}\left(  x\right)  \right\vert
\leq\left\Vert f\right\Vert _{w,0}\left(  \sqrt{\rho N}+k_{G}h_{X,K}\right)
,\quad x\in K,\label{1.86}%
\end{equation}

when $h_{X,K}<\infty$ and $k_{G}=\left(  2\pi\right)  ^{-d/4}\sqrt{-\left(
\left\vert D\right\vert ^{2}G\right)  \left(  0\right)  }\sqrt{d}$.\medskip

\fbox{Unisolvent error estimates} in Section \ref{Sect_ex_unisolv_error}. If
the weight function has property W02 for some parameter $\kappa\geq1$, the
independent data $X$ contains a unisolvent set of order $\left\lfloor
\kappa\right\rfloor $ and $X$ is contained in a bounded data region $\Omega$.
Then using results from the Lagrange theory of interpolation we show in
Theorem \ref{Thm_Exact_smth_ord_gte_1} that there exist constants
$K_{\Omega,m}^{\prime},k_{G}>0$ such that
\begin{equation}
\left\vert f\left(  x\right)  -s_{e}\left(  x\right)  \right\vert
\leq\left\Vert f\right\Vert _{w,0}\left(  K_{\Omega,m}^{\prime}\sqrt{\rho
N}+k_{G}\left(  h_{X,\Omega}\right)  ^{\left\lfloor \kappa\right\rfloor
}\right)  ,\quad x\in\overline{\Omega},\label{1.87}%
\end{equation}

when $h_{X,\Omega}=\sup\limits_{\omega\in\Omega}\operatorname*{dist}\left(
x,X\right)  \leq h_{G}$ is the maximum spherical cavity size (radius).\medskip

\fbox{$W^{1,\infty}\left(  \Omega\right)  \cap X_{w}^{0}\left(  \Omega\right)
$ spaces} in Section \ref{Sect_Ex_smth_err_H1inf_data}. The local data
functions are assumed to have bounded derivatives a.e. on the bounded data
region. The main \textbf{one-dimensional} result is Theorem
\ref{Thm_err_exsmth_H1inf_1dim}: if there exist constants $c_{0},c_{1}\geq0$,
independent of $f_{d}$, such that%
\[
\left\vert Ds_{e}\left(  x\right)  \right\vert \leq c_{0}\left\Vert
f_{d}\right\Vert _{\infty;\Omega}+c_{1}\left\Vert Df_{d}\right\Vert
_{\infty;\Omega},\quad x\in\Omega,
\]

then%
\begin{align*}
\left\vert f_{d}\left(  x\right)  -s_{e}\left(  x\right)  \right\vert
\leq\left\Vert f_{d}\right\Vert _{w,0,\Omega} &  \min\left\{  \sqrt{\rho
N},\frac{G\left(  0\right)  }{\left(  2\pi\right)  ^{d/2}}\right\}  +\\
&  +\left(  c_{0}\left\Vert f_{d}\right\Vert _{\infty;\Omega}+\left(
1+c_{1}\right)  \left\Vert Df_{d}\right\Vert _{\infty;\Omega}\right)
h_{X,\Omega},\quad x\in\Omega.
\end{align*}

and we consider the example of the scaled hat basis function $\Lambda\left(
x./\lambda\right)  $. If the basis function has "large" support w.r.t. the
data region $\Omega$ i.e. if $\operatorname*{diam}\Omega\leq\lambda$, then we
show the order of \textbf{convergence of the interpolant is 1}.

The \textbf{analogous multivariate result} is Theorem
\ref{Thm_err_exsmth_H1inf_multivar} which assumes the cone condition but I
have not yet been able to extend my hat function example to higher dimensions.

The above theoretical results will be illustrated using the weight function
examples from the interpolation chapter, namely the radial \textit{shifted
thin-plate splines}, \textit{Gaussian} and \textit{Sobolev splines} and the
tensor product \textit{extended B-splines}. We will also use the
\textit{central difference} weight functions from Chapter
\ref{Ch_cent_diff_wt_fn_ten_prod}.

\textbf{Numerical results are only presented for the 1-dim non-unisolvent data
cases}. Numeric experiments are carried out using the same 1-dimensional
B-splines and data functions that were used for the interpolants. We restrict
ourselves to one dimension so that the data density parameters $h_{X,\Omega}$
and $h_{X,K}$ can be calculated.

\subsection{Chapter \ref{Ch_Approx_smth} The scalable Approximate smoother}

In \textbf{Section} \ref{Sect_space_Jg} the \textbf{convolution space}
$J_{G}=G\ast S$ is introduced and then some applications of these spaces are
presented. The fact that $J_{G}\subset X_{w}^{0}$ implies some interesting
embedding results which involve supplying necessary and sufficient conditions
for one data space $X_{w}^{0}$ to be \textbf{continuously embedded} in another
e.g. Theorem \ref{Thm_Xw1o_embed_Xw2o_iff}: $X_{w_{1}}^{0}\overset{\iota
}{\hookrightarrow}X_{w_{2}}^{0}$ iff $\frac{w_{2}}{w_{1}}\in L^{\infty}$ and
$\left\Vert \iota\right\Vert =\left\Vert \frac{w_{2}}{w_{1}}\right\Vert
_{\infty}$.

I have also begun a study of the operator $\int_{\Omega}\overline{R}_{x}u$
where $u\in X_{w}^{0}$ and $R_{x}$ is the Riesz representer of the evaluation
functional. This operator might be useful in deriving the Approximate smoother.

We introduce a smoother which we call the \textit{Approximate smoother}
because it approximates the Exact smoother. This is a non-parametric, scalable
smoother. Here \textit{scalable} means the numeric effort to calculate the
Approximate smoother depends linearly on the number of data points.\textit{\ }%
We assume the basis function is real-valued.

Two different approaches will be taken to defining the Approximate smoother,
and both involve formulating the smoother as the solution of a variational
problem. One of these problems will involve minimizing the Exact smoother
functional \ref{1.8} over $W_{G,X^{\prime}}$ where $X^{\prime}=\left\{
x_{i}^{\prime}\right\}  _{i=1}^{N^{\prime}}$ is an arbitrary set of distinct
points in $\mathbb{R}^{d}$. The other, equivalent problem, involves finding
the function in $W_{G,X^{\prime}}$ which is nearest to Exact smoother $s_{e}$
w.r.t. the norm $\left\Vert \cdot\right\Vert _{w,0}$. If
\[
s_{a}\left(  x\right)  =\sum\limits_{i=1}^{N^{\prime}}\alpha_{i}^{\prime
}R_{x_{i}^{\prime}}\left(  x\right)  =\left(  2\pi\right)  ^{-d/2}%
\sum\limits_{i=1}^{N^{\prime}}\alpha_{i}^{\prime}G\left(  x-x_{i}^{\prime
}\right)  ,
\]

denotes the Approximate smoother and $y$ is the dependent data then solving
the second problem yields%
\[
s_{a}=\mathcal{I}_{X^{\prime}}s_{e},
\]

which implies the matrix equation%
\[
\left(  N\rho R_{X^{\prime},X^{\prime}}+R_{X,X^{\prime}}^{T}R_{X,X^{\prime}%
}\right)  \alpha^{\prime}=R_{X,X^{\prime}}^{T}y,
\]

where $R_{X,X^{\prime}}=\left(  R_{x_{j}^{\prime}}\left(  x^{\left(  i\right)
}\right)  \right)  $. The size of the Approximate smoother matrix is
$N^{\prime}\times N^{\prime}$ which is independent of the number of data
points and suggests scalability.

The error estimates for the pointwise convergence of the Approximate smoother
to its data function $f\in X_{w}^{0}$ are based on the simple triangle
inequality%
\[
\left\vert f\left(  x\right)  -s_{a}\left(  x\right)  \right\vert
\leq\left\vert f\left(  x\right)  -s_{e}\left(  x\right)  \right\vert
+\left\vert s_{e}\left(  x\right)  -s_{a}\left(  x\right)  \right\vert ,
\]

and so Section \ref{Sect_ap_Ex_smth_minus_App_smth} will be devoted to
estimating $\left\vert s_{e}\left(  x\right)  -s_{a}\left(  x\right)
\right\vert $.

As with the minimal interpolant and the Exact smoother, we will obtain
estimates that involve unisolvent data sets as well as the Type 1 and Type 2
estimates that do not explicitly involve unisolvency. The Approximate smoother
convergence orders and the constants are the same as those for the
interpolation case which are given in the interpolation tables
\ref{Tbl_intro_NonUnisolvTyp1Converg}, \ref{Tbl_intro_NonUnisolvTyp2Conv},
\ref{Tbl_intro_UnisolvConverg} and \ref{Tbl_intro_ConvergCentral}.\medskip

\fbox{Type 1 error estimates} No \textit{a priori} assumption is made
concerning the weight function parameter $\kappa$ but it will be assumed that
the basis function satisfies an inequality of the form \ref{1.9}. For example,
if it is assumed that the data region $K$ is closed bounded and infinite then
Theorem \ref{Thm_converg_arb_func_K=0} establishes that%
\begin{equation}
\left\vert s_{e}\left(  x\right)  -s_{a}\left(  x\right)  \right\vert
\leq\left\Vert f\right\Vert _{w,0}k_{G}\left(  h_{X^{\prime},K}\right)
^{s},\quad x\in K,\label{1.88}%
\end{equation}

and Theorem \ref{Thm_ap_Appr_smth_err_typ1} shows that%
\begin{equation}
\left\vert f\left(  x\right)  -s_{a}\left(  x\right)  \right\vert
\leq\left\Vert f\right\Vert _{w,0}\left(  \sqrt{\rho N}+k_{G}\left(
h_{X,K}\right)  ^{s}+k_{G}\left(  h_{X^{\prime},K}\right)  ^{s}\right)  ,\quad
x\in K,\label{1.91}%
\end{equation}

when $h_{X,K}=\sup\limits_{x\in K}\operatorname*{dist}\left(  x,X\right)  \leq
h_{G}$ and $h_{X^{\prime},K}=\sup\limits_{x\in K}\operatorname*{dist}\left(
x,X^{\prime}\right)  \leq h_{G}$.\medskip

\fbox{Type 2 error estimates} If it only assumed that $\kappa\geq1$ then by
Theorem \ref{Thm_ap_G(0)minusG(x)_bound}%
\begin{equation}
\left\vert s_{e}\left(  x\right)  -s_{a}\left(  x\right)  \right\vert
\leq\left\Vert f\right\Vert _{w,0}k_{G}\left(  h_{X^{\prime},K}\right)
^{s},\quad x\in\mathbb{R}^{d},\label{1.92}%
\end{equation}

and by Theorem \ref{Thm_ap_Appr_smth_err_typ2}%
\begin{equation}
\left\vert f\left(  x\right)  -s_{a}\left(  x\right)  \right\vert
\leq\left\Vert f\right\Vert _{w,0}\left(  \sqrt{\rho N}+k_{G}h_{X,K}%
+k_{G}h_{X^{\prime},K}\right)  ,\quad x\in\mathbb{R}^{d},\label{1.93}%
\end{equation}

where $k_{G}=\left(  2\pi\right)  ^{-\frac{d}{4}}\sqrt{-\left(  \left\vert
D\right\vert ^{2}G\right)  \left(  0\right)  }\sqrt{d}$.\medskip

\fbox{Unisolvent data error estimates} If $X$ is a unisolvent set of order
$m\geq1$ contained in a bounded data region $\Omega$ then by Theorem
\ref{Thm_ConvergSaToSe_not_unisolv_II},
\begin{equation}
\left\vert s_{e}\left(  x\right)  -s_{a}\left(  x\right)  \right\vert
\leq\left\Vert f\right\Vert _{w,0}k_{G}\left(  h_{X^{\prime},\Omega}\right)
^{m},\quad x\in\overline{\Omega},\label{1.94}%
\end{equation}

and by Theorem \ref{Thm_Err_unisolv_arbitrary_fn},%
\begin{equation}
\left\vert f\left(  x\right)  -s_{a}\left(  x\right)  \right\vert
\leq\left\Vert f\right\Vert _{w,0}\left(  K_{\Omega,m}^{\prime}\sqrt{\rho
N}+k_{G}\left(  h_{X,K}\right)  ^{m}+k_{G}\left(  h_{X^{\prime},K}\right)
^{m}\right)  ,\quad x\in\overline{\Omega},\label{1.98}%
\end{equation}

for some constants $K_{\Omega,m}^{\prime},k_{G}>0$. We say the orders of
convergence are at least $m$.

These theoretical results will be illustrated using the weight function
examples from the interpolation chapter, namely the radial \textit{shifted
thin-plate splines}, \textit{Gaussian} and \textit{Sobolev splines} and the
tensor product \textit{extended B-splines}. We will also use the
\textit{central difference} weight functions from Chapter
\ref{Ch_cent_diff_wt_fn_ten_prod}.

\textbf{Numerical results are only presented for the Type 1 and 2 estimates}.
Numeric experiments are carried out using the same 1-dimensional B-splines and
data functions that were used for the interpolants. We restrict ourselves to
one dimension so that the data density parameters $h_{X,\Omega} $ and
$h_{X,K}$ can be easily calculated.

We discuss the \textit{SmoothOperator} software (freeware) package which
implements the Approximate smoother algorithm. It has a full user manual which
describe several tutorials and data experiments.

\subsection{Chapter \ref{Ch_Xo1/w_and_fnal_on_Xow} The spaces $X_{1/w}^{0}$,
$\protect\widetilde{X}_{1/w}^{0}$ and the bounded linear functionals on
$X_{w}^{0}$}

The goal of this chapter is to use a bilinear form to characterize the bounded
linear functionals $\left(  X_{w}^{0}\right)  ^{\prime}$ on the data space
$X_{w}^{0}$ introduced in Chapter \ref{Ch_wtfn_basisfn_datasp} for the weight
functions used as examples in this document and to explore the properties of
some related operators.

A semi-Hilbert space $X_{1/w}^{0}$ is defined using the reciprocal of the
weight function and isometric isomorphisms are constructed between this space
and $L^{2}$. These isometries and those constructed in Chapter
\ref{Ch_wtfn_basisfn_datasp} between the data space and $L^{2}$ are used to
construct an isometric isomorphism between functional space and the data space.

This chapter involves successively larger classes of weight functions with
increasingly more general definitions of $X_{1/w}^{0}$.

For more details see the introduction to the chapter.

\subsection{Chapter \ref{Ch_bnd_deriv_hat_smth_large_supp} An upper bound for
the derivative of the 1-dimensional hat basis function smoother}

In this chapter we derive an upper bound for the derivative of the
1-dimensional (scaled) hat basis function smoother under the assumption that
the data function has a bounded derivative on the data region $f_{d}\in
X_{w}^{0}\left(  \Omega\right)  \cap W^{1,\infty}\left(  \Omega\right)  $.
This will be used in Example \ref{Ex_ex_smth_hat_larg_supp_1dim} to obtain a
convergence estimate for the Exact smoother when the scaled hat function has
"large" support w.r.t. the data region i.e. when the \textbf{scaled hat basis
function} $\Lambda_{\lambda}\left(  x\right)  =\Lambda\left(  x/\lambda
\right)  $ has been scaled so that $\operatorname*{diam}\Omega\leq\frac{1}%
{2}\operatorname*{diam}\operatorname*{supp}\Lambda_{\lambda}$ i.e.
$\operatorname*{diam}\Omega\leq\lambda$. This convergence estimate is
\ref{a1.15}, namely%
\[
\left\vert s_{e}\left(  x\right)  -f_{d}\left(  x\right)  \right\vert
\leq\left\Vert E\right\Vert \left\Vert f_{d}\right\Vert _{\infty;\Omega}%
\sqrt{\rho N}+\left(  \frac{2}{\lambda}\left\Vert f_{d}\right\Vert
_{\infty;\Omega}+\left(  1+\min\left\{  5,2+\rho N\right\}  \right)
\left\Vert Df_{d}\right\Vert _{\infty;\Omega}\right)  h_{X,\Omega},
\]

where $s_{e}=\mathcal{S}_{X}f_{d}=\mathcal{S}_{X}Ef_{d}$ is the Exact
smoother, $E$ is any continuous linear extension of the data function $f_{d}$
from $\Omega$ to $\mathbb{R}^{d}$ and $h_{X,\Omega}$ is the maximum spherical
cavity size (radius).

Also see introduction to Chapter \ref{Ch_bnd_deriv_hat_smth_large_supp}.

\subsection{Chapter \ref{Ch_exten_rect_condit} Explicit extension operators
based on Wloka but using the rectangle condition}

I needed these explicit extension operators to characterize of the restriction
spaces $X_{w}^{0}\left(  \Omega\right)  $\textbf{\ }for several classes of
weight function, especially the tensor products in Subsection
\ref{Sect_local_data_space}, but $\Omega$ has ??? very awkward constraints.
Subsequently I discovered the very general theoretical extension operator
$r_{\Omega}^{\ast}$ (Theorem \ref{Thm_canon_exten_op}) and so the explicit
extension became unnecessary.

In this chapter, motivated by the work of Wloka \cite{Wloka87} and the Russian
mathematicians described in the Background section, we will derive
\textbf{four extension operators}. Technically:

\begin{enumerate}
\item We will assume that the bounded domain $\Omega$ satisfies the ??? very
restrictive rectangle condition instead of the cone condition - a sphere does
not satisfy the rectangle condition.

\item I will modify the Calderon-Zygmund extension result used in Theorem 5.4,
Section 5.2 of Wloka \cite{Wloka87} by replacing the cone condition by the
rectangle condition but I will still use the integral representation and the
Fourier transform.

\item We adapt the (convolution) integral representation technique. This is in
the form of the integral operators $J_{\delta}\left[  v\right]  \left(
z\right)  =\int_{\mathcal{O}_{\mathbf{1}}}\lambda^{\delta}e^{-\lambda
\mathbf{1}}v\left(  \frac{z}{\lambda}\right)  d\lambda$ where $v$ is weakened
to be an $L^{1}$ function.

\item Develop an extension operator which extends a function from a single
orthant $\mathcal{O}_{\theta}\subset\mathbb{R}^{d}$ to the entire space.

\item Still uses a $C^{\infty}$ partition of unity to define an extension on a
domain $\Omega$.

\item The basic norms are the $L^{2}$ and $L^{\infty}$ norms.
\end{enumerate}

The four extension operators are:

\begin{enumerate}
\item In Section \ref{Sect_Exten_FromOrthant} I construct continuous extension
$\mathcal{E}_{\alpha}^{\mathbf{1}}:C_{0}^{\left(  \alpha\right)  }\left(
\overline{\mathcal{O}_{\mathbf{1}}}\right)  \rightarrow C_{B}^{\left(
\alpha\right)  }\left(  \mathbb{R}^{d}\right)  $ for each $\alpha
\geq\mathbf{1}$. This is then generalized to a continuous extension
$\mathcal{E}_{\alpha}^{\theta}:C_{0}^{\left(  \alpha\right)  }\left(
\overline{\mathcal{O}_{\theta}}\right)  \rightarrow C_{B}^{\left(
\alpha\right)  }\left(  \mathbb{R}^{d}\right)  $. The results of this section
are not used elsewhere and are not used to construct any of the other
extensions in this Chapter.

\item In Section \ref{Sect_ExtenLocWn1_to_Wn1_Fourier} a partition of unity
and the integral representation of Lemma \ref{Lem_SmthFuncIntegRepInOrthant}
to construct continuous convolution extension operators $E_{\Omega
}^{n\mathbf{1}}:W^{n\mathbf{1}}\left(  \Omega\right)  \rightarrow
W^{n\mathbf{1}}\left(  \mathbb{R}^{d}\right)  $ for $n\geq1$. This is done in
Theorem \ref{Thm_ExtenOrthantSobolFourier} and continuity is demonstrated
using a Fourier transform argument. Here $J_{\left(  n-2\right)  \mathbf{1}%
}^{\theta}\left[  v\right]  $ satisfies \ref{a030}.

\item In Section \ref{Sect_ExtenLocCatoBndCa} we generalize the $E_{\Omega
}^{n\mathbf{1}}$ to extension operators $E_{\Omega}^{\alpha}:C^{\left(
\alpha\right)  }\left(  \overline{\Omega}\right)  \rightarrow C_{B}^{\left(
\alpha\right)  }\left(  \mathbb{R}^{d}\right)  $ which are continuous under
the supremum norm. To do this we constrain the function $v$ to have bounded
support in $\mathcal{O}_{\mathbf{1}}$ and to be a tensor product function
which satisfies \ref{X551}. The main extension result is Theorem
\ref{Thm_ExtenContinFuncs_OrthantProp}.

\item In this section we show that the particular extension operators
$\left\{  E_{\Omega}^{n\mathbf{1}}\right\}  _{n\geq1}$ are such that each
$E_{\Omega}^{n\mathbf{1}}:W^{n\mathbf{1}}\left(  \Omega\right)  \rightarrow
W^{n\mathbf{1}}\left(  \mathbb{R}^{d}\right)  $ is continuous if we assume
that in the integral operator $J_{\left(  n-2\right)  \mathbf{1}}^{\theta
}\left[  v\right]  $ the function $v\in L_{0}^{1}\left(  \mathcal{O}%
_{1}\right)  $ is a tensor product with property.
\end{enumerate}

\textbf{In the Appendix\medskip}

\subsection{Chapter \ref{Ch_Appendx_basic_notation} Basic notation,
definitions and symbols}

Basic function spaces, multi-index and vector notation, topology and tempered
distributions, Fourier transforms, convolutions, Taylor series expansions etc.
See the table of contents.

\subsection{Chapter \ref{Ch_quot_sp_reprod_kern} Quotient spaces and
reproducing kernels}

We present the basic theory Banach space quotient spaces and then the general
theory of reproducing kernels. Finally, we characterize restriction spaces in
terms of quotient spaces.

\subsection{Chapter \ref{Ch_NotesOnSchmeisser} Notes on Schmeisser
\cite{Schmeis2006}}

This is a very short chapter which just contains some notes on the 2006 survey
concerning Sobolev spaces with dominating mixed derivatives by Schmeisser
\cite{Schmeis2006}.

\subsection{Chapter \ref{Ch_Hm1_proofs_Adams} Proofs of the claims made in
Remark \ref{Rem_SobolevSpace2}.}

This chapter contains the proofs of the claims made for the Sobolev spaces
with dominating mixed derivatives $W^{m\mathbf{1}}\left(  \Omega\right)  $
etc. in Remark \ref{Rem_SobolevSpace2}. They are modified results from Adams
and Fournier \cite{AdamFour2003} and Wloka \cite{Wloka87} concerning the
Sobolev spaces $W^{m}\left(  \Omega\right)  $.

\chapter{Weight functions, data spaces and basis
functions\label{Ch_wtfn_basisfn_datasp}}

\section{Introduction}

In \cite{LightWayneX98Weight} Light and Wayne developed a weight function
theory of positive order basis function interpolation in which the weight
function directly generated both the basis function and the Hilbert space of
continuous data functions. We intend to do the same thing here for the zero
order case.

We start by introducing the theory of weight functions using as examples the
radial basis functions: the shifted thin-plate splines, the Gaussian, the
Sobolev splines as well as the tensor product extended B-spline weight
functions. Results are also proved regarding the tensor products and products
of weight functions.

The weight functions are then used to define both the Hilbert data spaces and
the basis functions which will be used to formulate and solve the variational
interpolation and smoothing problems discussed in later Chapters. Various
continuity results are derived for the basis functions and data space
functions. The data space is shown to be a reproducing kernel Hilbert space
but instead of using the reproducing kernel we will use the Riesz representers
of the evaluation functionals $f\rightarrow D^{\alpha}f\left(  x\right)  $
because the kernel is only defined for $\alpha=0$.

\section{The weight functions\label{Sect_gen_wt_funcs}}

\subsection{Motivation for the weight function
properties\label{SbSect_motiv_weight_fn}}

We want to modify Light and Wayne's weight function properties and basis
function definition so that the \textbf{tensor product hat functions} are
basis functions of zero order generated by weight functions..

Light and Wayne defined their weight function $w$ to have the three
properties:\medskip

\begin{description}
\item[A3.1] $w$ is continuous and positive outside the origin.

\item[A3.2] $\frac{1}{w}\in L_{loc}^{1}$.

\item[A3.3] For some $R>0$ and $\mu$ real, $\frac{1}{w\left(  x\right)  }\leq
c\left\vert x\right\vert ^{\mu}$ for $\left\vert x\right\vert >R$.\medskip
\end{description}

These properties imply $\frac{1}{w}\in S^{\prime}$ and Light and Wayne's
definition of a \textbf{tempered basis distribution} $G\in S^{\prime}$ of
order $k\geq1$ generated by $w$ (see Section 3 of \cite{LightWayneX98Weight})
imply that as distributions
\[
\widehat{G}=\frac{1}{\left\vert \xi\right\vert ^{2k}w},\text{ }on\text{
}\mathbb{R}^{d}\setminus0.
\]

Part 11 of Remark \ref{Rem_Def_extend_wt_fn} regarding the weight function properties.

Now, as mentioned above, we want to define a basis function so the
multivariate hat function is a basis function of order zero. The one
dimensional hat function will be denoted by $\Lambda$ and is given by
\begin{equation}
\Lambda\left(  x\right)  =\left\{
\begin{array}
[c]{ll}%
0, & \left\vert x\right\vert >1,\text{\quad}x\in\mathbb{R},\\
1-\left\vert x\right\vert , & \left\vert x\right\vert \leq1,\text{\quad}%
x\in\mathbb{R}.
\end{array}
\right. \label{1.05}%
\end{equation}

In higher dimensions the hat function is defined as the tensor product
\begin{equation}
\Lambda\left(  x\right)  =\prod_{i=1}^{d}\Lambda\left(  x_{i}\right)
,\text{\quad}x\in\mathbb{R}^{d}.\label{1.06}%
\end{equation}

It is well known that
\begin{equation}
\widehat{\Lambda}\left(  \xi\right)  =\left(  2\pi\right)  ^{-1/2}\left(
\dfrac{\sin\left(  \xi/2\right)  }{\xi/2}\right)  ^{2},\text{\quad}\xi
\in\mathbb{R}^{1},\label{1.008}%
\end{equation}

so that in higher dimensions
\begin{equation}
\widehat{\Lambda}\left(  \xi\right)  =\prod\limits_{i=1}^{d}\widehat{\Lambda
}\left(  \xi_{i}\right)  =\left(  2\pi\right)  ^{-d/2}\prod\limits_{i=1}%
^{d}\left(  \dfrac{\sin\xi_{i}/2}{\xi_{i}/2}\right)  ^{2},\text{\quad}\xi
\in\mathbb{R}^{d}.\label{1.13}%
\end{equation}

In one dimension, if $\Lambda$ is to be basis function of order zero then we
must have $\widehat{\Lambda}\left(  \xi\right)  =\frac{1}{w\left(  \xi\right)
}$ i.e.
\begin{equation}
w\left(  \xi\right)  =\left(  2\pi\right)  ^{d/2}\prod\limits_{i=1}^{d}\left(
\dfrac{\xi_{i}/2}{\sin\left(  \xi_{i}/2\right)  }\right)  ^{2},\text{\quad}%
\xi\in\mathbb{R}^{d},\label{1.30}%
\end{equation}

and $w$ has discontinuities outside the origin, which violates weight function
property A3.1. Denoting the set of discontinuities by $\mathcal{A}$ we have
that $\mathcal{A}$ is the union of hyperplanes%
\begin{equation}
\mathcal{A}=\bigcup\limits_{i=1}^{d}\left\{  \xi:\xi_{i}\in2\pi\mathbb{Z}%
\setminus0\right\}  .\label{1.46}%
\end{equation}

and this is a closed set of measure zero. This will become weight function
property W01 in Definition \ref{Def_extend_wt_fn} below.

Additionally, from equation \ref{1.13} we have that $\widehat{\Lambda}%
=\frac{1}{w}\in L^{1}$ and we want to only use $L^{1}$ Fourier transform
results. So we will say that $G$ is a basis function of order $0$ if $\frac
{1}{w}\in L^{1}$ and $\widehat{G}=\frac{1}{w}$. A well-known $L^{1}$ result
then implies $G\in C_{B}^{\left(  0\right)  }\subset S^{\prime}$ and so $G$ is
a function. That $\frac{1}{w}\in L^{1}$ is ensured by properties W02 and W03
in Definition \ref{Def_extend_wt_fn}. Note that our definition of a basis
function gives a unique basis function whereas in Light and Wayne a basis
function of higher order $k$ is unique modulo a polynomial of order $\leq k$.

Finally we will incorporate the parameter $\kappa$ into properties W02 and W03
to increase the smoothness of the basis function and to allow a larger choice
of weight functions. The parameter $\kappa$ will be allowed to take real
non-integer values or a vector with non-integer components.

\subsection{The weight function properties}

\begin{definition}
\label{Def_extend_wt_fn}\textbf{Weight functions with parameter} $\kappa$

A weight function's properties are defined with reference to a\textbf{\ }set
$\mathcal{A}\subset\mathbb{R}^{d}$. The introduction of $\mathcal{A}$ is
motivated directly by the properties of the Fourier transform of the hat function.

In this document a weight function $w$ is a mapping $w:\mathbb{R}%
^{d}\rightarrow\mathbb{R}$ which has at least the property W01:

\begin{description}
\item[W01] There exists a closed set $\mathcal{A}$ with measure zero such that
$w$ is continuous, positive and finite outside $\mathcal{A}$ i.e. $w\in
C^{\left(  0\right)  }\left(  \mathbb{R}^{d}\setminus\mathcal{A}\right)  $ and
$w>0$ on $\mathbb{R}^{d}\setminus\mathcal{A}$.

It may also have property W02 or W03:

\item[W02] For some real number \fbox{$\kappa\in\mathbb{R}^{1}$}, $\kappa
\geq0$ we have%
\begin{equation}
\int\dfrac{\left\vert x\right\vert ^{2\lambda}}{w\left(  x\right)  }%
dx<\infty,\quad0\leq\lambda\leq\kappa.\label{1.036}%
\end{equation}

\item[W03] For some \fbox{$\kappa\in\mathbb{R}^{d}$}, $\kappa\geq\mathbf{0}$
we have%
\begin{equation}
\int\dfrac{x^{2\lambda}}{w\left(  x\right)  }dx<\infty,\quad\mathbf{0}%
\leq\lambda\leq\kappa.\label{1.051}%
\end{equation}

\end{description}
\end{definition}

\begin{remark}
\label{Rem_Def_extend_wt_fn}\ 

\begin{enumerate}
\item Regarding property W01: There exists a unique minimal $\mathcal{A}_{c}$
which is the intersection of all closed sets$\mathcal{A}$ with measure zero
such that $w\in C^{\left(  0\right)  }\left(  \mathbb{R}^{d}\setminus
\mathcal{A}_{c}\right)  $ and a unique minimal $\mathcal{A}_{p}$ which is the
intersection of all closed sets $\mathcal{A}$ with measure zero satisfying
$w>0$ on $\mathbb{R}^{d}\setminus\mathcal{A}$. Further, $\mathcal{A}%
=\mathcal{A}_{c}\cap\mathcal{A}_{p}$ is the unique minimal set for which
property W01 is valid.

Also, the countable union of closed sets of measure zero is also a closed set
of measure zero.

If $w$ is a weight function then $1/w$ is a weight function.

\item In Section \ref{Sect_Xow} property W01 is used to define the inner
product space of distributions $X_{w}^{0}$. Property W01 is needed to ensure
that $X_{w}^{0}$ is not empty - see Theorem \ref{Thm_Xwth_non_empty}. It is
not enough just to have a weight function that is positive a.e.

\item Property W01 allows basis functions to have Fourier transforms with
zeros outside the origin e.g. the tensor product hat function.

\item Property W02 is especially suitable for \textbf{radial} weight
functions. I have chosen the integral form \ for W02 because it is more
general than the form A3.3 of Light and Wayne.

Property W03 is especially suitable for \textbf{tensor product} weight
functions. In general it allows more smoothness information to be obtained
about $X_{w}^{0}$ and the basis function.

\item Properties W02 and W03 coincide in one dimension.

\item Properties W02 and W03 allow the definition of a continuous basis
function of order zero.

\item The parameter $\kappa$ of properties W02 and W03 could be called the
\textbf{smoothness parameter} because in Theorem \ref{Thm_X_smooth} it is
shown that $X_{w}^{0}\subset C_{B}^{\left(  \left\lfloor \kappa\right\rfloor
\right)  }$ and in Theorems \ref{Thm_basis_fn_properties_all_m_W2} and
\ref{Thm_basis_fn_properties_all_m_W3} it is shown that the basis functions
are in $C_{B}^{\left(  \left\lfloor 2\kappa\right\rfloor \right)  }$. The
allowance of non-integer values of $\kappa$ sometimes permits an extra degree
of differentiability to be estimated e.g. The smoothness of Sobolev spline
basis functions in Subsubsection \ref{Ex_BasisFuncExSobSpline}.

\item Scaling (or dilation): scaled weight functions i.e. $w\left(  \lambda
x\right)  $ where $\lambda\in\mathbb{R}^{d}$ and $\lambda.>0$, are also weight
functions and properties W02 and W03 are also valid for the same parameter
$\kappa$.

\item Properties such as: $\int\dfrac{\left\vert x\right\vert ^{2\lambda}%
}{w\left(  x\right)  ^{2}}dx<\infty$ for $0\leq\lambda\leq\kappa$, and
$\frac{\left\vert x\right\vert ^{\mu}}{w\left(  x\right)  }$ bounded for
$0\leq\mu\leq v$, might be investigated in order to obtain Sobolev space results.

\item See Section \ref{Sect_space_Jg} for results concerning necessary and
sufficient conditions on weight function/s for the equivalence of $X_{w}^{0}$
spaces and for the embedding of $X_{w}^{0}$ spaces in Sobolev spaces.

\item \textbf{Life can be more complicated}. Instead of the above properties
we could instead use the positive order weight function properties (see
Subsection 1.2.3 of Williams \cite{WilliamsPosOrdSmthV3}) but with $\theta=0
$. Then properties W2.1 and properties W2.2 would stay the same i.e.%
\begin{equation}
\left.
\begin{array}
[c]{l}%
\mathbf{W2.1}:1/w\in L_{loc}^{1},\\
\medskip\\
\mathbf{W2.2}:\int\limits_{\left\vert \cdot\right\vert \geq r_{2}}\frac
{1}{w\left\vert \cdot\right\vert ^{2\sigma}}<\infty\text{ }for\text{
}some\text{ }\sigma>0\text{ }and\text{ }some\text{ }r_{2}>0,
\end{array}
\right\} \label{3.00}%
\end{equation}

and the W3 properties would become W02 or W03.

However W2.1 and W2.2 would still allow us to construct basis distributions
and the data function space would no longer necessarily be a space of
continuous functions but it would still be a Hilbert space with the operators
$\mathcal{I}:X_{w}^{0}\rightarrow L^{2}$ and $\mathcal{J}:L^{2}\rightarrow
X_{w}^{0}$ of Definition \ref{Def_I_J} still being isometric isomorphisms and inverses.
\end{enumerate}
\end{remark}

The following result gives two equivalent criteria for the weight function
property W02.

\begin{theorem}
\label{Thm_equiv_W2}The the following criteria are equivalent to weight
function property W02:

\begin{enumerate}
\item $\int\dfrac{x^{2\beta}\left\vert x\right\vert ^{\kappa-\left\lfloor
\kappa\right\rfloor }}{w\left(  x\right)  }dx<\infty$,$\quad\left\vert
\beta\right\vert \leq\kappa$.

\item $\frac{1}{w}\in L_{loc}^{1}$ and $\int_{\left\vert \cdot\right\vert \geq
R}\dfrac{\left\vert \cdot\right\vert ^{2\kappa}}{w}<\infty$ for some $R\geq0$.
\end{enumerate}
\end{theorem}

\begin{proof}
Part 1 is easily proved using the identity $\left\vert x\right\vert ^{2k}%
=\sum\limits_{\left\vert \beta\right\vert =k}\frac{k!}{\beta!}x^{2\beta}$, and
clearly W02 implies part 2. Finally, if part 2 holds for $0\leq\lambda
\leq\kappa$ and some $R\geq0$
\begin{align*}
\int\dfrac{\left\vert \cdot\right\vert ^{2\lambda}}{w}=\int\limits_{\left\vert
\cdot\right\vert \leq R}\dfrac{\left\vert \cdot\right\vert ^{2\lambda}}%
{w}+\int\limits_{\left\vert \cdot\right\vert \geq R}\dfrac{\left\vert
\cdot\right\vert ^{2\lambda}}{w} &  \leq R^{2\lambda}\int\limits_{\left\vert
\cdot\right\vert \leq R}\dfrac{1}{w}+\int\limits_{\left\vert \cdot\right\vert
\geq R}\frac{1}{\left\vert \cdot\right\vert ^{2\left(  \kappa-\lambda\right)
}}\dfrac{\left\vert \cdot\right\vert ^{2\kappa}}{w}\\
&  \leq R^{2\lambda}\int\limits_{\left\vert \cdot\right\vert \leq R}\dfrac
{1}{w}+\int\limits_{\left\vert \cdot\right\vert \geq R}\frac{1}{R^{2\left(
\kappa-\lambda\right)  }}\dfrac{\left\vert \cdot\right\vert ^{2\kappa}}{w}\\
&  <\infty,
\end{align*}

and so property W02 is satisfied.
\end{proof}

The next theorem relates properties W02 and W03:

\begin{theorem}
\label{Thm_equiv_W3}Weight function properties W02 and W03 are related as
follows: if $w$ is a weight function on $\mathbb{R}^{d}$ then:

\begin{enumerate}
\item If $w$ has property W02 for $\kappa=s$ then $w$ has property W03 for
$\kappa=\frac{s}{d}\mathbf{1}$.

\item If $w$ has property W03 for $\kappa=\mu$ then $w$ has property W02 for
$\kappa=\underline{\mu}=\min\mu=\min_{i}\mu_{i}$.
\end{enumerate}
\end{theorem}

\begin{proof}
\textbf{Part 1} $\int\frac{x^{2\frac{s}{d}\mathbf{1}}}{w}=\int\frac{\left\vert
x^{2\frac{s}{d}\mathbf{1}}\right\vert }{w}\leq\int\frac{\left\vert
x\right\vert ^{\left\vert 2\frac{s}{d}\mathbf{1}\right\vert }}{w}=\int%
\frac{\left\vert x\right\vert ^{2s}}{w}<\infty$.\medskip

\textbf{Part 2} If $w$ has property W03 for $\kappa=\mu$ then $1/w\in L^{1} $
and $\frac{x^{2\lambda}}{w}\in L^{1}$ when $\lambda\leq\mu$. Now%
\[
\left\vert \cdot\right\vert ^{2\underline{\mu}}=\left\vert \cdot\right\vert
^{2\left(  \left\lfloor \underline{\mu}\right\rfloor +\underline{\mu
}-\left\lfloor \underline{\mu}\right\rfloor \right)  }=\left\vert
\cdot\right\vert ^{2\left\lfloor \underline{\mu}\right\rfloor }\left\vert
\cdot\right\vert ^{2\left(  \underline{\mu}-\left\lfloor \underline{\mu
}\right\rfloor \right)  },
\]

and since from \ref{1.57},%
\[
\left\vert x\right\vert ^{2\left\lfloor \underline{\mu}\right\rfloor
}=\left\lfloor \underline{\mu}\right\rfloor !\sum_{\left\vert \alpha
\right\vert =\left\lfloor \underline{\mu}\right\rfloor }\frac{1}{\alpha
!}x^{2\alpha},
\]

and also%
\begin{align*}
\left\vert x\right\vert ^{2\left(  \underline{\mu}-\left\lfloor \underline{\mu
}\right\rfloor \right)  }  & =\left(  x_{1}^{2}+x_{2}^{2}+\ldots+x_{d}%
^{2}\right)  ^{\underline{\mu}-\left\lfloor \underline{\mu}\right\rfloor }\\
& \leq x_{1}^{2\left(  \underline{\mu}-\left\lfloor \underline{\mu
}\right\rfloor \right)  }+x_{2}^{2\left(  \underline{\mu}-\left\lfloor
\underline{\mu}\right\rfloor \right)  }+\ldots+x_{d}^{2\left(  \underline{\mu
}-\left\lfloor \underline{\mu}\right\rfloor \right)  }\\
& =\sum_{k=1}^{d}x_{k}^{2\left(  \underline{\mu}-\left\lfloor \underline{\mu
}\right\rfloor \right)  },
\end{align*}

it follows that%
\begin{align*}
\left\vert x\right\vert ^{2\underline{\mu}}  & \leq\left\lfloor \underline{\mu
}\right\rfloor !\sum_{\left\vert \alpha\right\vert =\left\lfloor
\underline{\mu}\right\rfloor }\frac{x^{2\alpha}}{\alpha!}\sum_{k=1}^{d}%
x_{k}^{2\left(  \underline{\mu}-\left\lfloor \underline{\mu}\right\rfloor
\right)  }\\
& =\left\lfloor \underline{\mu}\right\rfloor !\sum_{k=1}^{d}\sum_{\left\vert
\alpha\right\vert =\left\lfloor \underline{\mu}\right\rfloor }\frac{1}%
{\alpha!}x^{2\alpha}x_{k}^{2\left(  \underline{\mu}-\left\lfloor
\underline{\mu}\right\rfloor \right)  }\\
& =\left\lfloor \underline{\mu}\right\rfloor !\sum_{k=1}^{d}\sum_{\left\vert
\alpha\right\vert =\left\lfloor \underline{\mu}\right\rfloor }\frac{1}%
{\alpha!}x^{2\left(  \alpha+\left(  \underline{\mu}-\left\lfloor
\underline{\mu}\right\rfloor \right)  \mathbf{e}_{k}\right)  }.
\end{align*}

Thus%
\[
\int\frac{\left\vert \cdot\right\vert ^{2\underline{\mu}}}{w}\leq\left\lfloor
\underline{\mu}\right\rfloor !\sum_{k=1}^{d}\sum_{\left\vert \alpha\right\vert
=\left\lfloor \underline{\mu}\right\rfloor }\frac{1}{\alpha!}\int%
\frac{x^{2\left(  \alpha+\left(  \underline{\mu}-\left\lfloor \underline{\mu
}\right\rfloor \right)  \mathbf{e}_{k}\right)  }}{w},
\]

but $\left\vert \alpha+\left(  \underline{\mu}-\left\lfloor \underline{\mu
}\right\rfloor \right)  \mathbf{e}_{k}\right\vert =\left\vert \alpha
\right\vert +\underline{\mu}-\left\lfloor \underline{\mu}\right\rfloor
=\underline{\mu}$ so that $\alpha+\left(  \underline{\mu}-\left\lfloor
\underline{\mu}\right\rfloor \right)  \mathbf{e}_{k}\leq\underline{\mu
}\mathbf{1}\leq\mu$ and we have $\int\frac{\left\vert \cdot\right\vert
^{2\underline{\mu}}}{w}<\infty$ as desired.
\end{proof}

\subsection{Interpolated weight functions}

Here we are concerned with the interpolated weight functions $w=w_{1}%
^{1-t}w_{2}^{t}$ where $t\in\left[  0,1\right]  $.

\begin{theorem}
\label{Thm_interpol_wt_fns}\textbf{Interpolated weight functions} Suppose
$w_{1}$ and $w_{2}$ are weight functions with property W01 for sets
$\mathcal{A}_{1}$ and $\mathcal{A}_{2}$ respectively. Then for all
$t\in\left[  0,1\right]  $:

\begin{enumerate}
\item $w_{1}^{1-t}w_{2}^{t}$ has property W01 on the set $\mathcal{A}%
=\mathcal{A}_{1}\cup\mathcal{A}_{2}$.

\item If $w_{i}\in W02$ for $\kappa_{i}$ then $w_{1}^{1-t}w_{2}^{t}\in W02$
for $\kappa\leq\min\left\{  \kappa_{1},\kappa_{2}\right\}  $.

\item If $w_{i}\in W03$ for $\kappa_{i}$ then $w_{1}^{1-t}w_{2}^{t}\in W03$
for $\kappa\leq\min\left\{  \kappa_{1},\kappa_{2}\right\}  $.
\end{enumerate}
\end{theorem}

\begin{proof}
\textbf{Part 1}. Clearly, outside $\mathcal{A}$, both $w_{1}$ and $w_{2}$ are
continuous, positive and finite.\medskip

\textbf{Part 2}. Since $1-t+t=1$ by H\"{o}lder's theorem ($p=1/\left(
1-t\right)  $, $q=1/t$) and part 2 of Theorem \ref{Thm_equiv_W2}:
\begin{align*}
\int\limits_{\left\vert \cdot\right\vert \geq1}\frac{\left\vert \cdot
\right\vert ^{2\kappa}}{w_{1}^{1-t}w_{2}^{t}}=\int\limits_{\left\vert
\cdot\right\vert \geq1}\frac{\left\vert \cdot\right\vert ^{2\left(
1-t\right)  \kappa}}{w_{1}^{1-t}}\frac{\left\vert \cdot\right\vert ^{2t\kappa
}}{w_{2}^{t}}  & \leq\left(  \int\limits_{\left\vert \cdot\right\vert \geq
1}\frac{\left\vert \cdot\right\vert ^{2\kappa}}{w_{1}}\right)  ^{1-t}\left(
\int\limits_{\left\vert \cdot\right\vert \geq1}\frac{\left\vert \cdot
\right\vert ^{2\kappa}}{w_{2}}\right)  ^{t}\\
& \leq\left(  \int\limits_{\left\vert \cdot\right\vert \geq1}\frac{\left\vert
\cdot\right\vert ^{2\kappa_{1}}}{w_{1}}\right)  ^{1-t}\left(  \int%
\limits_{\left\vert \cdot\right\vert \geq1}\frac{\left\vert \cdot\right\vert
^{2\kappa_{2}}}{w_{2}}\right)  ^{t}\\
& <\infty,
\end{align*}

and if $r>0$,%
\[
\int\limits_{\left\vert \cdot\right\vert \leq r}\frac{1}{w_{1}^{1-t}w_{2}^{t}%
}\leq\left(  \int\limits_{\left\vert \cdot\right\vert \leq r}\frac{1}{w_{1}%
}\right)  ^{1-t}\left(  \int\limits_{\left\vert \cdot\right\vert \leq r}%
\frac{1}{w_{2}}\right)  ^{t}<\infty.
\]
\smallskip

\textbf{Part 3} Again by H\"{o}lder's theorem, when $\mathbf{0}\leq\lambda
\leq\kappa$,
\[
\int\dfrac{x^{2\lambda}}{w_{1}^{1-t}w_{2}^{t}}dx=\int\frac{\left\vert
\cdot\right\vert ^{2\left(  1-t\right)  \lambda}}{w_{1}^{1-t}}\frac{\left\vert
\cdot\right\vert ^{2t\lambda}}{w_{2}^{t}}\leq\left(  \int\frac{\left\vert
\cdot\right\vert ^{2\lambda}}{w_{1}}\right)  ^{1-t}\left(  \int\frac
{\left\vert \cdot\right\vert ^{2\lambda}}{w_{2}}\right)  ^{t}<\infty.
\]

\end{proof}

\subsection{Radial weight functions\label{SbSect_radial_wt_funcs}}

Suppose $w$ is a weight function w.r.t. the closed set $\mathcal{A}$ of
measure zero. Then $1/w$ is positive, finite and continuous outside
$\mathcal{A}$. Motivated by this fact I prove: ?? \textbf{What about
continuity}?

\begin{lemma}
\label{Lem_pos_fin_cont_partial_integ}Suppose $d\geq2$ and $f\in L^{1}\left(
\mathbb{R}^{d}\right)  $ is positive outside a closed set $\mathcal{A}$ of
measure zero. Suppose $\mathbb{R}^{d}=\mathbb{R}^{d_{1}}\times\mathbb{R}%
^{d_{2}}$.

Then the function
\[
g\left(  \xi^{\prime}\right)  :=\int_{\mathbb{R}^{d_{2}}}f\left(  \xi^{\prime
},\xi^{\prime\prime}\right)  d\xi^{\prime\prime},
\]

is in $L^{1}\left(  \mathbb{R}^{d_{1}}\right)  $ with
\[
\int g=\int f.
\]

Further, $g$ is positive outside the closed set of measure zero
\[
\mathcal{A}_{1}:=\left(  \pi_{1}\mathcal{A}^{c}\right)  ^{c}.
\]

If $\mathcal{A}$ is bounded then $\mathcal{A}_{1}$ is empty.
\end{lemma}

\begin{proof}
Since $f$ is positive a.e. we have $\int g=\int f$.

Clearly $\mathcal{A}_{1}$ is closed and if $\mathcal{A}$ is bounded then
$\mathcal{A}_{1}$ is empty.

We must show that $\mathcal{A}_{1}$ has measure zero when $\mathcal{A}$ is
unbounded:%
\begin{align*}
\mathcal{A}_{1}  & =\left(  \pi_{1}\mathcal{A}^{c}\right)  ^{c}=\left\{
x^{\prime}\in\mathbb{R}^{d_{1}}:x^{\prime}\notin\pi_{1}\mathcal{A}%
^{c}\right\}  =\\
& =\left\{  x^{\prime}\in\mathbb{R}^{d_{1}}:\left(  x^{\prime},x^{\prime
\prime}\right)  \notin\mathcal{A}^{c}\text{ }\forall x^{\prime\prime}%
\in\mathbb{R}^{d_{2}}\right\} \\
& =\left\{  x^{\prime}\in\mathbb{R}^{d_{1}}:\left(  x^{\prime},x^{\prime
\prime}\right)  \in\mathcal{A}\text{ }\forall x^{\prime\prime}\in
\mathbb{R}^{d_{2}}\right\} \\
& =\left\{  x^{\prime}\in\mathbb{R}^{d_{1}}:\pi_{1}^{-1}x^{\prime}%
\subset\mathcal{A}\right\}  .
\end{align*}

Thus $\mathcal{A}_{1}$ is not empty iff $\pi_{1}^{-1}\mathcal{A}_{1}%
\subset\mathcal{A}$. This implies that when $\mathcal{A}_{1}$ is not empty we
have
\[
\mathcal{A}_{1}\subset\pi_{1}\mathcal{A}.
\]

We thus have the unique representation%
\[
\mathcal{A}=\left(  \pi_{1}^{-1}\mathcal{A}_{1}\right)  \cup\mathcal{B}%
,\quad\left(  \pi_{1}^{-1}\mathcal{A}_{1}\right)  \cap\mathcal{B}=\emptyset,
\]

with%
\[
\left(  \pi_{1}\mathcal{B}^{c}\right)  ^{c}=\emptyset,\quad\mathcal{A}_{1}%
=\pi_{1}\mathcal{A}\setminus\pi_{1}\mathcal{B}.
\]

Now $0=\operatorname*{meas}\mathcal{A}=\operatorname*{meas}\pi_{1}%
^{-1}\mathcal{A}_{1}=\operatorname*{meas}\mathcal{A}_{1}$.

Choose $\xi^{\prime}\in\mathcal{A}_{1}^{c}=\pi_{1}\mathcal{A}^{c}$. Then there
exists $\xi^{\prime\prime}\in\mathbb{R}^{d_{2}}$ such that $\xi=\left(
\xi^{\prime},\xi^{\prime\prime}\right)  \in\mathcal{A}^{c}$. Thus
$B_{r}\left(  \xi\right)  \subset\mathcal{A}^{c}$ where
$r=\operatorname*{dist}\left(  \xi,\mathcal{A}^{c}\right)  $. Hence if
$\left\vert \xi^{\prime}\right\vert ^{2}+\left\vert \xi^{\prime\prime
}\right\vert ^{2}\leq r^{2}$ then $\left\vert \xi^{\prime\prime}\right\vert
^{2}\leq r^{2}-\left\vert \xi^{\prime}\right\vert ^{2}$ and consequently%
\[
g\left(  \xi^{\prime}\right)  =\int_{\mathbb{R}^{d_{2}}}f\left(  \xi^{\prime
},\xi^{\prime\prime}\right)  d\xi^{\prime\prime}\geq\int_{\left\vert
\xi^{\prime\prime}\right\vert ^{2}\leq r^{2}-\left\vert \xi^{\prime
}\right\vert ^{2}}f\left(  \xi^{\prime},\xi^{\prime\prime}\right)
d\xi^{\prime\prime}>0,
\]

since $f$ is positive outside $\mathcal{A}$.
\end{proof}

We begin with two basic lemmas on radial weight functions:

?? Note the radial result of Theorem \ref{Thm_Integ_u(xy)f(|x|)dx}.

\begin{lemma}
\label{Lem_wt_fn_part_integral}Suppose $w$ has property W01 w.r.t. the closed
set $\mathcal{A}\subset\mathbb{R}^{d}$ and property W02 for parameter $\kappa
$. Write $\xi=\left(  \xi^{\prime},\xi^{\prime\prime}\right)  \in
\mathbb{R}^{d}$ where $\xi^{\prime}\in\mathbb{R}^{m}$ and $\xi^{\prime\prime
}\in\mathbb{R}^{n}$. Then%
\[
\frac{1}{w^{\prime}\left(  \xi^{\prime}\right)  }=\int\frac{d\xi^{\prime
\prime}}{w\left(  \xi^{\prime},\xi^{\prime\prime}\right)  },
\]

defines an a.e. positive function $w^{\prime}$ on $\mathbb{R}^{m}$ and it
follows that:
\end{lemma}

\begin{enumerate}
\item if $w$ is radial then $w^{\prime}$ is radial.

\item If $w$ has property W02 for parameter $\kappa$ then $w^{\prime}$ also
has property W02 for parameter $\kappa$.
\end{enumerate}

\begin{proof}
\textbf{Part 1} Suppose $w\left(  \xi\right)  =w_{\odot}\left(  \left\vert
\xi\right\vert \right)  $. Then%
\[
\frac{1}{w^{\prime}\left(  \xi^{\prime}\right)  }=\int\frac{d\xi^{\prime
\prime}}{w_{\odot}\left(  \sqrt{\left\vert \xi^{\prime}\right\vert
^{2}+\left\vert \xi^{\prime\prime}\right\vert ^{2}}\right)  }.
\]
\smallskip

\textbf{Part 2} If $0\leq\lambda\leq\kappa$ then $\int\frac{\left\vert
\xi^{\prime}\right\vert ^{2\lambda}d\xi^{\prime}}{w^{\prime}\left(
\xi^{\prime}\right)  }=\int\int\frac{\left\vert \xi^{\prime}\right\vert
^{2\lambda}d\xi^{\prime}d\xi^{\prime\prime}}{w\left(  \xi\right)  }\leq
\int\frac{\left\vert \xi\right\vert ^{2\lambda}}{w\left(  \xi\right)  }%
d\xi<\infty$.
\end{proof}

\begin{lemma}
\label{Lem_wt_func_radial}\textbf{Basic radial weight function results}
Suppose $u\geq0$ is a measurable function on $\mathbb{R}^{1}$, and the weight
function $w\left(  \xi\right)  $ is radial on $\mathbb{R}^{d}$ $\left(
d\geq1\right)  $ with property W02 for the parameter $\kappa$. Write
$\xi=\left(  \xi_{1},\xi^{\prime\prime}\right)  $ and define $w_{\odot}\left(
\left\vert x\right\vert \right)  :=w\left(  x\right)  $.

Then if $\int\frac{u\left(  x\xi\right)  }{w\left(  \xi\right)  }d\xi$ exists
it is a radial function of $x$ and%
\begin{equation}
\int_{\mathbb{R}^{d}}\frac{u\left(  x\xi\right)  }{w\left(  \xi\right)  }%
d\xi=\int_{\mathbb{R}^{1}}\frac{u\left(  \left\vert x\right\vert s\right)
}{\overset{\circ}{w}\left(  s\right)  }ds,\label{a118}%
\end{equation}

where%
\begin{equation}
\frac{1}{\overset{\circ}{w}\left(  s\right)  }:=\int_{\mathbb{R}^{d-1}}%
\frac{d\xi^{\prime\prime}}{w\left(  s,\xi^{\prime\prime}\right)  }%
=\int_{\mathbb{R}^{d-1}}\frac{d\xi^{\prime\prime}}{w_{\odot}\left(
\sqrt{s^{2}+\left\vert \xi^{\prime\prime}\right\vert ^{2}}\right)  },\quad
s\in\mathbb{R}^{1}.\label{a119}%
\end{equation}

Here $\overset{\circ}{w}$ is positive a.e. and is an even function with
property W02 for parameter $\kappa$.

?? \textbf{NEED} Lemma \ref{Lem_pos_fin_cont_partial_integ}? ?? Further,
\textbf{if we assume that} $\overset{\circ}{w}$ is continuous and positive
outside a closed set $\mathcal{B}$ of measure zero then $\overset{\circ}{w}$
is a weight function i.e. it has weight function property W01 w.r.t.
$\mathcal{B}$.
\end{lemma}

\begin{proof}
This lemma will employ the technique of Section 4.1 Stein and Weiss
\cite{SteinWeiss71} which defines radial functions in terms of orthogonal
transformations. Stein and Weiss observed that a function $f$ is radial if and
only if $f\left(  \mathcal{O}x\right)  =f\left(  x\right)  $ for any linear,
orthogonal transformation $\mathcal{O}:\mathbb{R}^{d}\rightarrow\mathbb{R}%
^{d}$ and any $x\in\mathbb{R}^{d}$. Now an orthogonal transformation
$\mathcal{O}$ satisfies $\mathcal{O}^{T}=\mathcal{O}^{-1}$ where
$\mathcal{O}xy=x\mathcal{O}^{T}y$ for the Euclidean inner product, and an
orthogonal transformation has a Jacobian of one.

Now $w\left(  \xi\right)  =w_{\odot}\left(  \left\vert \xi\right\vert \right)
$ so that
\begin{align*}
\int\frac{u\left(  \mathcal{O}x\xi\right)  }{w\left(  \xi\right)  }d\xi
=\int\frac{u\left(  x\mathcal{O}^{-1}\xi\right)  }{w_{\odot}\left(  \left\vert
\xi\right\vert \right)  }d\xi=\int\frac{u\left(  x\eta\right)  }{w_{\odot
}\left(  \left\vert \mathcal{O}\eta\right\vert \right)  }\left\vert J%
\genfrac{(}{)}{}{}{\xi}{\eta}%
\right\vert d\eta & =\int\frac{u\left(  x\eta\right)  }{w_{\odot}\left(
\left\vert \eta\right\vert \right)  }d\eta\\
& =\int\frac{u\left(  x\eta\right)  }{w\left(  \eta\right)  }d\eta,
\end{align*}

and so $\int\frac{u\left(  x\xi\right)  }{w\left(  \xi\right)  }d\xi$ is a
radial function of $x$ and we can define the radial function $\mu_{r}\left(
\left\vert x\right\vert \right)  :=\int\frac{u\left(  x\xi\right)  }{w\left(
\xi\right)  }d\xi$.

Now suppose $\rho\geq0$ and set $x=\left(  \rho,0^{\prime\prime}\right)  $.
Then $\left\vert x\right\vert =\rho$ and%
\begin{align*}
\mu_{r}\left(  \rho\right)  =\mu_{r}\left(  \left\vert x\right\vert \right)
=\int\frac{u\left(  \rho\xi_{1}\right)  }{w\left(  \xi\right)  }d\xi & =\int
u\left(  \rho\xi_{1}\right)  \int\frac{d\xi^{\prime\prime}}{w\left(
\xi\right)  }d\xi_{1}\\
& =\int u\left(  \rho\xi_{1}\right)  \int\frac{d\xi^{\prime\prime}}{w\left(
\xi\right)  }d\xi_{1}\\
& =\int_{\mathbb{R}^{1}}\frac{u\left(  \rho\xi_{1}\right)  }{\overset{\circ
}{w}\left(  \xi_{1}\right)  }d\xi_{1},
\end{align*}

by definition \ref{a119} of $\overset{\circ}{w}$. Since $w\left(  x\right)
=w_{\odot}\left(  \left\vert x\right\vert \right)  $ the second equation of
\ref{a119}\ is true and clearly $\overset{\circ}{w}$ is positive a.e. and an
even function. Suppose $0\leq\lambda\leq\kappa$. Then%
\[
\int\frac{\left\vert \xi_{1}\right\vert ^{2\lambda}}{\overset{\circ}{w}\left(
\xi_{1}\right)  }d\xi_{1}=\int\int_{\mathbb{R}^{d-1}}\frac{\left\vert \xi
_{1}\right\vert ^{2\lambda}}{w\left(  \xi_{1},\xi^{\prime\prime}\right)  }%
d\xi^{\prime\prime}d\xi_{1}\leq\int_{\mathbb{R}^{d}}\frac{\left\vert
\xi\right\vert ^{2\lambda}}{w\left(  \xi\right)  }d\xi<\infty,
\]

so that $\overset{\circ}{w}$ has property W02 for $\kappa$.
\end{proof}

\subsection{Examples of radial weight functions\label{SbSect_wt_func_examples}%
}

\begin{example}
\label{Ex_shft_thn_plt_spln_wt}\textbf{The shifted thin-plate splines} From
equations 25, 26 and 27 of Dyn \cite{Dyn89} the shifted thin-plate spline
functions%
\begin{equation}
H\left(  x\right)  =\left\{
\begin{array}
[c]{ll}%
\frac{\left(  -1\right)  ^{\nu+1}}{2}\left(  1+\left\vert x\right\vert
^{2}\right)  ^{\nu}\log\left(  1+\left\vert x\right\vert ^{2}\right)  , &
\nu=1,2,3,\ldots,\\
\left(  -1\right)  ^{\left\lceil \nu\right\rceil }\left(  1+\left\vert
x\right\vert ^{2}\right)  ^{\nu}, & \nu>-d/2,\text{ }\nu\neq0,1,2,\ldots,
\end{array}
\right. \label{1.043}%
\end{equation}

have the distribution Fourier transforms%
\[
\widehat{H}\left(  \xi\right)  =\widetilde{e}\left(  v\right)  \widetilde{K}%
_{v+d/2}\left(  \left\vert \xi\right\vert \right)  \left\vert \xi\right\vert
^{-2\left(  v+d/2\right)  }\;on\text{ }\mathbb{R}^{d}\setminus0,\text{\quad
}v>-d/2,
\]

where%
\[
\widetilde{e}\left(  \nu\right)  =\left\{
\begin{array}
[c]{ll}%
\left(  -1\right)  ^{\nu+1}\widetilde{c}^{\prime}\left(  2\nu\right)  , &
\nu=1,2,3,\ldots,\\
\left(  -1\right)  ^{\left\lceil \nu\right\rceil }\widetilde{c}\left(
2\nu\right)  , & \nu>-d/2,\text{ }\nu\neq1,2,3\ldots,
\end{array}
\right.
\]

and
\[
\widetilde{c}\left(  t\right)  =\left(  2\pi\right)  ^{d/2}2^{\frac{t+2}{2}%
}/\Gamma\left(  -t/2\right)  ,\text{\quad}\widetilde{c}^{\prime}\left(
t\right)  =\dfrac{d\widetilde{c}\left(  t\right)  }{dt},\text{\quad}t>0.
\]

Here $\widetilde{e}\left(  v\right)  >0$, $\widetilde{K}_{\lambda}\left(
t\right)  :=t^{\lambda}K_{\lambda}\left(  t\right)  ,$ $t\geq0,\lambda>0$,
$\widetilde{K}_{\lambda}\in C^{\infty}\left(  \mathbb{R}^{1}\setminus0\right)
$ and $K_{\lambda}$ is called a modified Bessel function or MacDonald's
function. Here $\widetilde{K}_{\lambda}$ has the properties
\begin{equation}
\widetilde{K}_{\lambda}\in C^{\left(  0\right)  }\left(  \mathbb{R}%
^{1}\right)  ,\text{ }\lambda>0;\ \widetilde{K}_{\lambda}\left(  t\right)
>0;\text{ }\lim\limits_{t\rightarrow\infty}\widetilde{K}_{\lambda}\left(
t\right)  =0\text{ }exponentially.\label{1.17}%
\end{equation}

See for example Abramowitz and Stegun \cite{AbramowStegun70}, Watson
\cite{Watson95} and Section 8.1 of Nikol'ski\u{\i} \cite{Nikol75}. See also
Theorem \ref{Thm_bnds_modif_MacDonald} below.

Now if $w=1/\widehat{H}$ then%
\begin{equation}
w\left(  \xi\right)  =\frac{1}{\widetilde{e}\left(  v\right)  }\frac
{\left\vert \xi\right\vert ^{2v+d}}{\widetilde{K}_{v+d/2}\left(  \left\vert
\xi\right\vert \right)  },\label{1.19}%
\end{equation}

and $w$ has property W01 for $\mathcal{A}=\left\{  0\right\}  $. Further,
condition \ref{1.036} holds iff $2\lambda-2\nu-d>-d$ for $0\leq\lambda
\leq\kappa$ i.e. iff
\begin{equation}
-d/2<\nu<0.\label{1.042}%
\end{equation}

and so when $\nu$ satisfies \ref{1.042}, $w$ has property W02 for all
$\kappa\geq0$.
\end{example}

\begin{example}
\label{Ex_Gaussian_wt}\textbf{The Gaussian} The Gaussian function%
\begin{equation}
H\left(  x\right)  =\exp\left(  -\left\vert x\right\vert ^{2}\right)
,\text{\quad}x\in\mathbb{R}^{d},\label{1.044}%
\end{equation}

has Fourier transform%
\[
\widehat{H}\left(  \xi\right)  =\frac{\sqrt{\pi}}{2}\exp\left(  -\frac
{\left\vert \xi\right\vert ^{2}}{4}\right)  ,\text{\quad}\xi\in\mathbb{R}^{d},
\]

and so if $\frac{1}{w}=\widehat{H}$%
\begin{equation}
w\left(  \xi\right)  =\frac{2}{\sqrt{\pi}}\exp\left(  \frac{\left\vert
\xi\right\vert ^{2}}{4}\right)  ,\label{1.029}%
\end{equation}

and $w$ has property W01 for empty $\mathcal{A}$ and property W02 for all
$\kappa\geq0$.
\end{example}

\begin{example}
\label{Ex_Sobolev_splin_wt}\textbf{The Sobolev splines} (Theorem 6.13 of
Wendland \cite{Wendland05}) If $\nu>d/2$ the multivariate modified Bessel
function%
\begin{equation}
H\left(  x\right)  =\frac{1}{2^{\nu-1}\Gamma\left(  \nu\right)  }%
\widetilde{K}_{\nu-\frac{d}{2}}\left(  \left\vert x\right\vert \right)
,\text{\quad}x\in\mathbb{R}^{d},\label{1.045}%
\end{equation}

has Fourier transform%
\[
\widehat{H}\left(  \xi\right)  =\frac{1}{\left(  1+\left\vert \xi\right\vert
^{2}\right)  ^{\nu}},\text{\quad}\xi\in\mathbb{R}^{d},
\]

where we have used the fact that $K_{-\lambda}=K_{\lambda}$ when $\lambda>0$.

Hence if $w:=1/\widehat{H}$,%
\begin{equation}
w\left(  \xi\right)  =\left(  1+\left\vert \xi\right\vert ^{2}\right)  ^{\nu
},\text{\quad}\xi\in\mathbb{R}^{d},\label{1.030}%
\end{equation}

then $w$ has property W01 for empty $\mathcal{A}$ and property W02 for
$0\leq\kappa<\nu-d/2$.
\end{example}

The next two examples of weight functions come from Subsubsection
\ref{SbSbSect_wt_fn_rad_L1loc}.

\begin{example}
\label{Ex_wt_func_L1loc}If $0\leq v<d/2$ and $v^{\prime}>d/2$ then
\[
w\left(  \xi\right)  =\left\{
\begin{array}
[c]{ll}%
\left\vert \xi\right\vert ^{2v}, & \left\vert \xi\right\vert \leq1,\\
\left\vert \xi\right\vert ^{2v^{\prime}}, & \left\vert \xi\right\vert \geq1,
\end{array}
\right.
\]

and this weight function has property W02 for $\kappa<v^{\prime}-d/2$.
\end{example}

\begin{example}
\label{Ex_wt_func_|x|^m_L1loc}Suppose the weight function $w$ has property W01
w.r.t. the set $\mathcal{A}=\left\{  \mathbf{0}\right\}  $ and there exists an
integer $m\geq0$, constants $r,c_{1},c_{2},s>0$ and a bounded function $v$
such that:%
\begin{equation}
w\left(  x\right)  =\frac{v\left(  x\right)  }{\left\vert x\right\vert
^{2\left(  m+s\right)  }},\quad\left\vert x\right\vert \leq r;\quad
2s<d;\quad0<c_{1}\leq v\left(  x\right)  \leq c_{2},\label{a1.09}%
\end{equation}

and for some $\tau\geq0$,%
\begin{equation}
\int\limits_{\left\vert \cdot\right\vert \geq r}\frac{w}{\left\vert
\cdot\right\vert ^{2\tau}}<\infty.\label{a112}%
\end{equation}

Are there such weight functions with property W02? The answer is yes because
Example \ref{Ex_wt_func_L1loc} is a specific case and a more general example
is%
\begin{equation}
w\left(  x\right)  =\left\{
\begin{array}
[c]{ll}%
\frac{1}{\left\vert x\right\vert ^{2\left(  m+s\right)  }}, & \left\vert
x\right\vert \leq1,\\
\left\vert x\right\vert ^{d+2\overline{\kappa}}, & \left\vert x\right\vert
\geq1,
\end{array}
\right. \label{a1.000}%
\end{equation}

where $\mathcal{A}=\left\{  \mathbf{0}\right\}  $, $s,\overline{\kappa}>0$ and
$\tau\geq d+\overline{\kappa}$. Since%
\[
\frac{1}{w}=\left\{
\begin{array}
[c]{ll}%
\left\vert x\right\vert ^{2\left(  m+s\right)  }, & \left\vert x\right\vert
\leq1,\\
\frac{1}{\left\vert x\right\vert ^{d+2\overline{\kappa}}}, & \left\vert
x\right\vert \geq1,
\end{array}
\right.
\]

$w$ has property W02 for $\kappa<\overline{\kappa}$.
\end{example}

\begin{theorem}
\label{Thm_bnds_modif_MacDonald}There exist positive constants $c_{\lambda}$
and $c_{\lambda}^{\prime}$ such that%
\[
c_{\lambda}e^{-t}\leq\widetilde{K}_{\lambda}\left(  t\right)  \leq c_{\lambda
}^{\prime}e^{-t},\text{\quad}t\geq0,\text{ }\lambda>0.
\]

Here $\widetilde{K}_{\lambda}\left(  t\right)  :=t^{\lambda}K_{\lambda}\left(
t\right)  $ and $K_{\lambda}$ is the MacDonald's function or the modified
Bessel function of the third kind.
\end{theorem}

\begin{proof}
This proof is based on estimating the well-known equation%
\[
\widetilde{K}_{\lambda}\left(  t\right)  =\frac{\pi^{1/2}2^{-\lambda}}{\left(
\lambda-\frac{1}{2}\right)  !}e^{-t}\int_{0}^{\infty}e^{-s}\left(  1+\frac
{2s}{t}\right)  ^{2\lambda-1}ds,\quad t,\lambda\geq0.
\]

See, for example, equation 11.137 in Arfken \cite{Arfken70}.\medskip

\textbf{UPPER\ BOUND}\medskip

\fbox{Assume $\lambda\geq1/2$ $and$ $t>0$}%
\begin{align*}
\int_{0}^{\infty}e^{-s}\left(  1+\frac{2s}{t}\right)  ^{2\lambda-1}ds  &
=\int_{0}^{t/2}e^{-s}\left(  1+\frac{2s}{t}\right)  ^{2\lambda-1}ds+\int%
_{t/2}^{\infty}e^{-s}\left(  1+\frac{2s}{t}\right)  ^{2\lambda-1}ds\\
& \leq\int_{0}^{t/2}e^{-s}\left(  1+1\right)  ds+\int_{t/2}^{\infty}%
e^{-s}\left(  \frac{2s}{t}+\frac{2s}{t}\right)  ^{2\lambda-1}ds\\
& =2\int_{0}^{t/2}e^{-s}ds+\int_{t/2}^{\infty}e^{-s}\left(  \frac{4s}%
{t}\right)  ^{2\lambda-1}ds\\
& =2\left(  1-e^{-t/2}\right)  +\left(  \frac{4}{t}\right)  ^{2\lambda-1}%
\int_{t/2}^{\infty}e^{-s}s^{2\lambda-1}ds\\
& <2+\left(  \frac{4}{t}\right)  ^{2\lambda-1}\int_{0}^{\infty}e^{-s}%
s^{2\lambda-1}ds\\
& =2+\left(  \frac{4}{t}\right)  ^{2\lambda-1}\left(  2\lambda-1\right)  !.
\end{align*}

\fbox{Assume $0<\lambda\leq1/2$ $and$ $t>0$}%
\[
\int_{0}^{\infty}e^{-s}\left(  1+\frac{2s}{t}\right)  ^{2\lambda-1}ds=\int%
_{0}^{\infty}\frac{e^{-s}ds}{\left(  1+\frac{2s}{t}\right)  ^{1-2\lambda}%
}<\int_{0}^{\infty}e^{-s}ds=1.
\]

Thus%
\[
\int_{0}^{\infty}e^{-s}\left(  1+\frac{2s}{t}\right)  ^{2\lambda-1}ds<\left\{
\begin{array}
[c]{ll}%
1, & 0<\lambda<1/2\text{ }and\text{ }t>0,\\
2+\left(  \frac{4}{t}\right)  ^{2\lambda-1}\left(  2\lambda-1\right)  !, &
\lambda\geq1/2\text{ }and\text{ }t>0.
\end{array}
\right.
\]

and so for $t>0$%
\[
\widetilde{K}_{\lambda}\left(  t\right)  <\frac{\pi^{1/2}2^{-\lambda}}{\left(
\lambda-\frac{1}{2}\right)  !}e^{-t}\times\left\{
\begin{array}
[c]{ll}%
1, & 0<\lambda\leq1/2,\\
2+\left(  \frac{4}{t}\right)  ^{2\lambda-1}\left(  2\lambda-1\right)  !, &
\lambda\geq1/2,
\end{array}
\right.
\]

and for $t\geq1$,%
\[
\widetilde{K}_{\lambda}\left(  t\right)  <\frac{\pi^{1/2}2^{-\lambda}}{\left(
\lambda-\frac{1}{2}\right)  !}e^{-t}\times\left\{
\begin{array}
[c]{ll}%
1, & 0<\lambda\leq1/2,\\
2+4^{2\lambda-1}\left(  2\lambda-1\right)  !, & \lambda\geq1/2.
\end{array}
\right.
\]

Now from \ref{1.17}, $\widetilde{K}_{\lambda}\in C^{\left(  0\right)  }\left(
\mathbb{R}^{1}\right)  ,$ $\lambda>0\ and$ $\widetilde{K}_{\lambda}\left(
t\right)  >0$ and so $0<d^{\prime}\leq\widetilde{K}_{\lambda}\left(  t\right)
\leq d^{\prime\prime}$ when $0\leq t\leq1$. Thus for some constant
$c_{\lambda}^{\prime}>0$,
\[
\widetilde{K}_{\lambda}\left(  t\right)  \leq c_{\lambda}^{\prime}e^{-t},\quad
t\geq0,\text{ }\lambda>0.
\]

\textbf{LOWER BOUND}\medskip

\fbox{Assume $\lambda\geq1/2$} Then%
\[
\int_{0}^{\infty}e^{-s}\left(  1+\frac{2s}{t}\right)  ^{2\lambda-1}ds\geq
\int_{0}^{\infty}e^{-s}ds=1,
\]

and so%
\[
\widetilde{K}_{\lambda}\left(  t\right)  \geq\pi^{1/2}\frac{2^{-\lambda}%
}{\left(  \lambda-\frac{1}{2}\right)  !}e^{-t},\quad t\geq0,\text{ }%
\lambda\geq1/2.
\]

\fbox{Assume $0<\lambda<1/2$ $and$ $t\geq1$} Then%
\begin{align*}
\int_{0}^{\infty}e^{-s}\left(  1+\frac{2s}{t}\right)  ^{2\lambda-1}ds  &
=\int_{0}^{\infty}\frac{e^{-s}ds}{\left(  1+\frac{2s}{t}\right)  ^{1-2\lambda
}}>\int_{0}^{t}\frac{e^{-s}ds}{\left(  1+\frac{2s}{t}\right)  ^{1-2\lambda}%
}>\\
& >\int_{0}^{t}\frac{e^{-s}ds}{\left(  1+2\right)  ^{1-2\lambda}}%
=\frac{1-e^{-t}}{3^{1-2\lambda}}\geq\frac{1-e^{-1}}{3^{1-2\lambda}}>\frac
{1}{6},
\end{align*}

so that%
\[
\widetilde{K}_{\lambda}\left(  t\right)  \geq\frac{\pi^{1/2}}{6}%
\frac{2^{-\lambda}}{\left(  \lambda-\frac{1}{2}\right)  !}e^{-t},\quad
t\geq1,\text{ }0<\lambda\leq1/2.
\]

Now from \ref{1.17}, $\widetilde{K}_{\lambda}\in C^{\left(  0\right)  }\left(
\mathbb{R}^{1}\right)  ,$ $\lambda>0\ and$ $\widetilde{K}_{\lambda}\left(
t\right)  >0$ and so $0<d^{\prime}\leq\widetilde{K}_{\lambda}\left(  t\right)
\leq d^{\prime\prime}$ when $0\leq t\leq1$. Thus for some constant
$c_{\lambda}>0$,
\[
\widetilde{K}_{\lambda}\left(  t\right)  \geq c_{\lambda}e^{-t},\quad
t\geq0,\text{ }\lambda>0.
\]

\end{proof}

\subsection{Products of weight functions}

See Theorem \ref{Thm_convol_basis_funcs} where the product of weight functions
is discussed in relation to the convolution of basis functions.

\subsection{Tensor products of weight functions\label{SbSect_wt_fn_ten_prod}}

The next theorem gives some general conditions under which the tensor product
of two weight functions is a weight function.

\begin{theorem}
\label{Thm_ten_prod_two_wt_fns}Suppose $w_{1}$ and $w_{2}$ are weight
functions i.e. they have property W01.

\begin{enumerate}
\item Then the tensor product $w_{1}\otimes w_{2}$ also has property W01.

\item The tensor product $w_{1}\otimes w_{2}$ has property W02 for parameter
$\kappa\in\mathbb{R}^{1}$ iff both $w_{1}$ and $w_{2}$ have property W02 for
parameter $\kappa$.

\item The tensor product $w_{1}\otimes w_{2}$ has property W03 for parameter
$\kappa=\left(  \kappa_{i}\right)  $ iff each $w_{i}$ has property W02 for
parameter $\kappa_{i}$.
\end{enumerate}
\end{theorem}

\begin{proof}
\textbf{Part 1} In this proof $\operatorname*{meas}$ denotes the measure of a
set. Now suppose $w_{i}:\mathbb{R}^{d_{i}}\rightarrow\mathbb{R}$. Then each
$w_{i}$ is associated with a weight function set, call it $\mathcal{A}_{i}$,
which is closed and has measure zero. The first step will be to construct the
set $\mathcal{A}$ for $w_{1}\otimes w_{2}$ from the sets $\mathcal{A}_{i} $.
The candidate is
\begin{equation}
\mathcal{A}=\left(  \mathcal{A}_{1}\times\mathbb{R}^{d_{2}}\right)
\cup\left(  \mathbb{R}^{d_{1}}\times\mathcal{A}_{2}\right)  .\label{2.00}%
\end{equation}

We show $\mathcal{A}$ is closed and has measure zero.\textbf{\ }We use the
continuous projection operators $p_{i}:\mathbb{R}^{d}\rightarrow
\mathbb{R}^{d_{i}}$ defined by $p_{1}\left(  x\right)  =x^{\prime}$ and
$p_{2}\left(  x\right)  =x^{\prime\prime}$. Next, since $\mathcal{A}%
=\bigcup\limits_{i=1}^{2}p_{i}^{-1}\left(  \mathcal{A}_{i}\right)  $ and each
$\mathcal{A}_{i}$ is closed $\mathcal{A}$ is also closed since $p_{i}$ is
continuous. Finally it must be shown that $\operatorname*{meas}\mathcal{A}=0$.
But since
\[
\operatorname*{meas}\mathcal{A}=\operatorname*{meas}\left(  \mathcal{A}%
_{1}\times\mathbb{R}^{d_{2}}\cup\mathbb{R}^{d_{1}}\times\mathcal{A}%
_{2}\right)  \leq\operatorname*{meas}\left(  \mathcal{A}_{1}\times
\mathbb{R}^{d_{2}}\right)  +\operatorname*{meas}\left(  \mathbb{R}^{d_{1}%
}\times\mathcal{A}_{2}\right)  ,
\]

it is sufficient to show that $\operatorname*{meas}\left(  \mathcal{A}%
_{1}\times\mathbb{R}^{d_{2}}\right)  =0$. To do this we use the countable
additivity property of the Lebesgue measure on $\mathbb{R}^{d}$:
\begin{align*}
\operatorname*{meas}\left(  \mathcal{A}_{1}\times\mathbb{R}^{d_{2}}\right)   &
=\operatorname*{meas}\bigcup_{n=1}^{\infty}\left(  \mathcal{A}_{1}%
\times\left(  \overline{B\left(  0;n\right)  }-B\left(  0;n-1\right)  \right)
\right) \\
& =\sum_{n=1}^{\infty}\operatorname*{meas}\left(  \mathcal{A}_{1}\times\left(
\overline{B\left(  0;n\right)  }-B\left(  0;n-1\right)  \right)  \right) \\
& =0.
\end{align*}

\textbf{Part 2} Suppose $w_{1}$ and $w_{2}$ have property W02 for parameter
$\kappa$. We show $w_{1}\otimes w_{2}$ also has property W02 for parameter
$\kappa$. Recall that property W02 requires that
\[
\int\dfrac{\left\vert x\right\vert ^{2\lambda}}{w_{i}\left(  x\right)
}dx<\infty,\quad0\leq\lambda\leq\kappa;\text{ }i=1,2.
\]

Now%
\begin{equation}
\int\dfrac{\left\vert x\right\vert ^{2\lambda}}{w\left(  x\right)  }%
dx=\int\dfrac{\left(  \left\vert x^{\prime}\right\vert ^{2}+\left\vert
x^{\prime\prime}\right\vert ^{2}\right)  ^{\lambda}}{w_{1}\left(  x^{\prime
}\right)  w_{2}\left(  x^{\prime\prime}\right)  }dx,\label{1.96}%
\end{equation}

and the inequality
\begin{equation}
\left\vert x\right\vert ^{2\lambda}\leq2^{\lambda}\left(  \left\vert
x^{\prime}\right\vert ^{2\lambda}+\left\vert x^{\prime\prime}\right\vert
^{2\lambda}\right)  ,\label{1.97}%
\end{equation}

allows us to write%
\begin{align*}
\int\dfrac{\left\vert x\right\vert ^{2\lambda}}{w\left(  x\right)  }dx  &
\leq\int\dfrac{2^{\lambda}\left(  \left\vert x^{\prime}\right\vert ^{2\lambda
}+\left\vert x^{\prime\prime}\right\vert ^{2\lambda}\right)  }{w_{1}\left(
x^{\prime}\right)  w_{2}\left(  x^{\prime\prime}\right)  }dx\\
& =2^{\lambda}\left(  \int\dfrac{\left\vert x^{\prime}\right\vert ^{2\lambda}%
}{w_{1}\left(  x^{\prime}\right)  w_{2}\left(  x^{\prime\prime}\right)
}dx\right)  +2^{\lambda}\left(  \int\dfrac{\left\vert x^{\prime\prime
}\right\vert ^{2\lambda}}{w_{1}\left(  x^{\prime}\right)  w_{2}\left(
x^{\prime\prime}\right)  }dx\right) \\
& =2^{\lambda}\left(  \int\dfrac{\left\vert x^{\prime}\right\vert ^{2\lambda
}dx}{w_{1}\left(  x^{\prime}\right)  }\int\dfrac{dx}{w_{2}\left(
x^{\prime\prime}\right)  }\right)  +2^{\lambda}\left(  \int\dfrac{dx^{\prime}%
}{w_{1}\left(  x^{\prime}\right)  }\int\dfrac{\left\vert x^{\prime\prime
}\right\vert ^{2\lambda}dx^{\prime\prime}}{w_{2}\left(  x^{\prime\prime
}\right)  }\right) \\
& <\infty.
\end{align*}
\medskip

Now suppose $w_{1}\otimes w_{2}$ has property W02 for parameter $\kappa$.
Thus, by \ref{1.96}%
\[
\int\dfrac{\left(  \left\vert x^{\prime}\right\vert ^{2}+\left\vert
x^{\prime\prime}\right\vert ^{2}\right)  ^{\lambda}}{w_{1}\left(  x^{\prime
}\right)  w_{2}\left(  x^{\prime\prime}\right)  }dx<\infty,\quad0\leq
\lambda\leq\kappa,
\]

so that $\int\dfrac{\left\vert x^{\prime}\right\vert ^{2\lambda}}{w_{1}\left(
x^{\prime}\right)  w_{2}\left(  x^{\prime\prime}\right)  }dx<\infty$ and
$\int\dfrac{\left\vert x^{\prime\prime}\right\vert ^{2\lambda}}{w_{1}\left(
x^{\prime\prime}\right)  w_{2}\left(  x^{\prime\prime}\right)  }dx<\infty$.

Hence $\int\dfrac{\left\vert x^{\prime}\right\vert ^{2\lambda}}{w_{1}\left(
x^{\prime}\right)  }dx^{\prime}<\infty$ and $\int\dfrac{\left\vert
x^{\prime\prime}\right\vert ^{2\lambda}}{w_{2}\left(  x^{\prime\prime}\right)
}dx^{\prime\prime}<\infty$ i.e. $w_{1}$ and $w_{2}$ both have property W02 for
parameter $\kappa$.

\textbf{Part 3} Follows directly from the definition of property W03.
\end{proof}

The next result generalizes the previous theorem to an arbitrary number of
weight functions.

\begin{corollary}
\label{Cor_Thm_ten_prod_two_wt_fns}Suppose the functions $\left\{
w_{i}\right\}  _{i=1}^{n}$ satisfy property W01 of a weight function.

\begin{enumerate}
\item Then $\bigotimes\limits_{i=1}^{n}w_{i}$ has property W01.

\item $\bigotimes\limits_{i=1}^{n}w_{i}$ has property W02 for real parameter
$\kappa$ iff each $w_{i}$ has property W02 for $\kappa$.

\item $\bigotimes\limits_{i=1}^{n}w_{i}$ has property W03 for parameter
$\kappa=\left(  \kappa_{i}\right)  $ iff each $w_{i}$ has property W02 for
parameter $\kappa_{i}$.
\end{enumerate}
\end{corollary}

\begin{proof}
\textbf{Parts 1 and 2} An easy proof by induction using the two weight
function case given by Theorem \ref{Thm_ten_prod_two_wt_fns}. \textbf{Part 3}
follows directly from the definition of property W03.
\end{proof}

\subsection{Example: the hat weight function\label{SbSect_HatWeightFunc}}

The next theorem shows under what conditions the multivariate hat function is
a basis function of order zero.

\begin{theorem}
\label{Thm_hat_wt_extend_props}Using the multivariate hat function $\Lambda$
define the function $w_{s}$ by
\begin{equation}
w_{s}\left(  \xi\right)  =\dfrac{1}{\widehat{\Lambda}\left(  \xi\right)
}.\label{1.007}%
\end{equation}

Then $w_{s}$ satisfies the weight function properties W01, W02 and W03 in the
following manner:

\begin{enumerate}
\item $w_{s}$ satisfies property W01 for all $\kappa$.

\item $w_{s}$ satisfies property W02 for real parameter $\kappa\geq0$ iff
$\kappa<\frac{1}{2}$\textbf{.}

\item $w_{s}$ satisfies property W03 for parameter $\kappa=\left(  \kappa
_{i}\right)  $ iff $\kappa<\frac{1}{2}\mathbf{1}$, where $\mathbf{1}=\left(
1,\ldots,1\right)  \in\mathbb{R}^{d}$.
\end{enumerate}
\end{theorem}

\begin{proof}
\textbf{Part 1} By \ref{1.46}, $w_{s}$ is the set $\mathcal{A}$ which is
closed with measure zero.\medskip

\textbf{Part 2} Part 2 of Corollary \ref{Cor_Thm_ten_prod_two_wt_fns} tells us
that we need only establish property W02 in one-dimension. Now W03 holds for
$\kappa$ iff
\[
\int_{\mathbb{R}^{1}}\frac{x^{2\lambda}}{w_{s}\left(  x\right)  }=\int
x^{2\lambda}\frac{\sin^{2}\left(  x/2\right)  }{\left(  x/2\right)  ^{2}%
}dx<\infty,\text{\quad}0\leq\lambda\leq\kappa.
\]

But these integrals will exist iff they exist near infinity iff $2-2\lambda>1$
i.e. iff $\lambda<1/2$.\medskip

\textbf{Part 3} Part 3 of Corollary \ref{Cor_Thm_ten_prod_two_wt_fns} tells us
that we need only establish property W03 in one-dimension. But from Part 5 of
Remark \ref{Rem_Def_extend_wt_fn} properties W02 and W03 coincide in one
dimension so Part 2 gives this result.
\end{proof}

\subsection{Example: the extended B-spline weight
functions\label{SbSect_wt_fn_examples_2}}

We now introduce a new class of weight functions which satisfy property W02
and property W03. We call this class the \textbf{extended B-spline weight
functions} because in Subsection \ref{SbSect_basis_fns} it will be shown that
they generate basis functions which are the derivatives of the B-splines.

\begin{lemma}
\label{Lem_exist_Integ_pow(s)_evenpow(sin(s))}Suppose $a>0$ and
$l=0,1,2,\ldots$.

Then $\int_{a}^{\infty}t^{\lambda}\sin^{2l}tdt<\infty$ iff $\lambda<-1$ and
$\int_{0}^{a}t^{\mu}\sin^{2l}tdt<\infty$ iff $2l+\mu>-1$.
\end{lemma}

\begin{theorem}
\label{Thm_ex_splin_wt_fn_properties}\textbf{The (homogeneous) extended
B-spline weight functions} For given integers $l,n\geq1$ define the extended
B-spline weight function $w_{s}$ by
\begin{equation}
w_{s}\left(  x\right)  =\prod_{i=1}^{d}\frac{x_{i}^{2n}}{\sin^{2l}x_{i}%
},\text{\qquad}x=\left(  x_{1},\ldots,x_{d}\right)  \in\mathbb{R}%
^{d}.\label{1.032}%
\end{equation}

Define the closed set of measure zero $\mathcal{A}$ to be the union of
hyperplanes%
\[
\mathcal{A}=\bigcup\limits_{k\in\mathbb{Z}}\bigcup\limits_{i=1}^{d}\left\{
x:x_{i}=\pi k\right\}  .
\]

Then the function $w_{s}$ is a weight function with property W01 w.r.t.
$\mathcal{A}$.

Also, the weight function $w_{s}$ has \textbf{property W02} for real
$\kappa\geq0$ iff $n$ and $l$ satisfy
\begin{equation}
\kappa+1/2<n\leq l,\label{1.028}%
\end{equation}

and $w_{s}$ has \textbf{property W03} for $\kappa\in\mathbb{R}^{d}$,
$\kappa\geq0$ iff%
\begin{equation}
\kappa+1/2<n\mathbf{1};\text{ }n\leq l\mathbf{.}\label{1.054}%
\end{equation}

\end{theorem}

\begin{proof}
Clearly \textbf{property W01} is satisfied for $\mathcal{A}$ so $w$ is a
weight function. Corollary \ref{Cor_Thm_ten_prod_two_wt_fns} tells us that we
need only establish the criterion \ref{1.028} in one dimension.\medskip

Now \textbf{property W02} holds for real $\kappa\geq0$ iff
\begin{equation}
\int_{\mathbb{R}^{1}}\frac{x^{2\lambda}}{w\left(  x\right)  }=\int
x^{2\lambda}\frac{\sin^{2l}x}{x^{2n}}dx<\infty,\text{\quad}0\leq\lambda
\leq\kappa.\label{1.43}%
\end{equation}

But these integrals will exist if and only if they exist near the origin and
near infinity. By Lemma \ref{Lem_exist_Integ_pow(s)_evenpow(sin(s))} they will
exist near the origin iff $2\lambda+2l-2n>-1$ i.e. iff $\lambda>n-l-1/2$, and
they will exist near infinity iff $2\lambda-2n<-1$ i.e. iff $\lambda<n-1/2$.
Thus the integrals all exist iff $n-l-1/2<\lambda<n-1/2$ for $0\leq\lambda
\leq\kappa$ iff $\kappa<n-1/2$ and $n<l+1/2$ iff $\kappa+1/2<n\leq l$.\medskip

Regarding \textbf{property W03}: Part 3 of Corollary
\ref{Cor_Thm_ten_prod_two_wt_fns} tells us that we need only establish
property W03 in one-dimension. But from Part 5 of Remark
\ref{Rem_Def_extend_wt_fn} properties W02 and W03 coincide in one dimension so
by \ref{1.028} we have that $w$ has property W03 for $\kappa\in\mathbb{R}^{d}%
$, $\kappa\geq0$ iff $\kappa_{i}+1/2<n\leq l$ for all $i$ iff $\kappa
+1/2<n\mathbf{1}\leq l\mathbf{1}$ iff $\kappa+1/2<n\mathbf{1}$ and $n\leq
l$\textbf{.}
\end{proof}

\begin{remark}
\label{Rem_Thm_ex_splin_wt_fn_properties}These weight functions can be
extended to%
\[
w_{s}\left(  x\right)  =\prod_{i=1}^{d}\frac{x_{i}^{2v_{i}}}{\sin
^{2\lambda_{i}}x_{i}}=\frac{x^{2v}}{\sin^{2\lambda}x},\quad v\leq\lambda,
\]

where we have employed the notation $\sin^{2\lambda}x=\left(  \sin
^{2\lambda_{i}}x_{i}\right)  $.
\end{remark}

\section{The data spaces $X_{w}^{0}$\label{Sect_Xow}}

We now introduce the zero order data spaces $X_{w}^{0}$ which will later be
used to define the minimum norm interpolant. Here $w$ denotes the weight
function. In Chapters \ref{Ch_Interpol}, \ref{Ch_Exact_smth} and
\ref{Ch_Approx_smth}, $X_{w}^{0}$ will be used to define several variational
interpolation and smoothing problems. In a manner analogous to the Sobolev
spaces, mappings between $X_{w}^{0}$ and $L^{2}$ are used to show that
$X_{w}^{0}$ is a Hilbert space. Various smoothness and $C^{\infty}$ density
results are then established.

\begin{definition}
\label{Def_Xow}\textbf{The zero order inner product space} $X_{w}^{0}$

Suppose $w$ is a weight function i.e. it only has property W01 of Definition
\ref{Def_extend_wt_fn}. Then define
\begin{equation}
X_{w}^{0}=\left\{  f\in S^{\prime}:\widehat{f}\in L_{loc}^{1}\text{ }and\text{
}\sqrt{w}\widehat{f}\in L^{2}\right\}  ,\label{1.033}%
\end{equation}

and endow it with the norm and inner product
\[
\left\Vert f\right\Vert _{w,0}^{2}=\int w\left\vert \widehat{f}\right\vert
^{2},\text{\qquad}\left(  f,g\right)  _{w,0}=\int w\widehat{f}\overline
{\widehat{g}}.
\]

That $\left\Vert \cdot\right\Vert _{w,0}$ is a norm is simple to prove.
Suppose $\left\Vert f\right\Vert _{w,0}=0$. Then, since $w>0$ a.e., $\int
w\left\vert \widehat{f}\right\vert ^{2}=0$ implies that $\widehat{f}\left(
x\right)  =0$ a.e. and thus $f=0$ in the distribution sense.
\end{definition}

The next result reassures us that the space $X_{w}^{0}$ is non-empty.

\begin{theorem}
\label{Thm_Xwth_non_empty}If $w$ is a weight function w.r.t. the closed set
$\mathcal{A}$ of measure zero then
\begin{equation}
\left\{  \overset{\vee}{u}:u\in C_{0}^{\infty}\text{ }and\text{ }%
\operatorname*{supp}u\subset\mathbb{R}^{d}\setminus\mathcal{A}\right\}
\subset X_{w}^{0},\label{1.0571}%
\end{equation}

where the set on the left is not empty.
\end{theorem}

\begin{proof}
Firstly, $\mathbb{R}^{d}\setminus\mathcal{A}$ is not empty since
$\mathcal{A}=\mathbb{R}^{d}$ implies $\operatorname*{meas}\mathcal{A}\neq0$. Clearly

$f\in\left\{  \overset{\vee}{u}:u\in C_{0}^{\infty}\text{ }and\text{
}\operatorname*{supp}u\subset\mathbb{R}^{d}\setminus\mathcal{A}\right\}  $
implies $f\in S\subset S^{\prime}$ and $\widehat{f}\in C_{0}^{\infty}\subset
L_{loc}^{1}$. Also, since $\widehat{f}$ has bounded support in $\mathbb{R}%
^{d}\setminus\mathcal{A}$ and $w$ is positive and continuous on $\mathbb{R}%
^{d}\setminus\mathcal{A}$
\[
\int\left\vert \sqrt{w}\widehat{f}\right\vert ^{2}=\int w\left\vert
\widehat{f}\right\vert ^{2}<\infty.
\]

\end{proof}

\begin{remark}
We mention the interesting result
\[
X_{w}^{0}\hookrightarrow L^{2}\Longleftrightarrow\frac{1}{w}\in L^{\infty},
\]

proven in Theorem \ref{Thm_Xow_embed_L2_iff} of Section \ref{Sect_space_Jg}.
The are also more general results concerning necessary and sufficient
conditions on weight function/s for the equivalence of $X_{w}^{0}$ spaces and
for the embedding of $X_{w}^{0}$ spaces in Sobolev spaces.
\end{remark}

\subsection{The completeness of $X_{w}^{0}$\label{SbSect_Xow_complete}}

In a manner analogous to Sobolev space theory, the completeness of $X_{w}^{0}
$ will be established by constructing an isometric homeomorphism
$\mathcal{I}:X_{w}^{0}\rightarrow L^{2}$ and then making use of the
completeness of $L^{2}$.

\begin{definition}
\label{Def_I_J}\textbf{The linear mappings} $\mathcal{I}$ \textbf{and}
$\mathcal{J}$

Suppose $w$ is a weight function. Using the definition of the space $X_{w}%
^{0}$ we can define the linear mapping $\mathcal{I}:X_{w}^{0}\rightarrow
L^{2}$ by
\[
\mathcal{I}f=\left(  \sqrt{w}\widehat{f}\right)  ^{\vee},\text{\quad}f\in
X_{w}^{0}.
\]

If, in addition, we assume property W02 or W03 then, since $\frac{1}{w}\in
L^{1}$, $\frac{\widehat{g}}{\sqrt{w}}\in L^{1}\subset S^{\prime}$ when $g\in
L^{2}$, and thus we can define the linear mapping $\mathcal{J}:L^{2}%
\rightarrow S^{\prime}$ by
\[
\mathcal{J}g=\left(  \frac{\widehat{g}}{\sqrt{w}}\right)  ^{\vee},\text{\quad
}g\in L^{2}.
\]

\end{definition}

The linear mappings $\mathcal{I}$ and $\mathcal{J}$ have the following properties:

\begin{theorem}
\label{Thm_I_J_property}Suppose $w$ has property W01 of a weight function. Then:

\begin{enumerate}
\item $\mathcal{I}:X_{w}^{0}\rightarrow L^{2}$ is an isometry.

\item $\mathcal{I}$ is one-to-one.

Now suppose that in addition to property W01, $w$ has property W02 or W03 for
parameter $\kappa$. Then:

\item $\mathcal{J}:L^{2}\rightarrow X_{w}^{0}$ and is an isometry.

\item $\mathcal{J}$ is one-to-one.

\item $\mathcal{J}\circ\mathcal{I}=I$ on $X_{w}^{0}$, and $\mathcal{I}%
\circ\mathcal{J}=I$ on $L^{2}$.

\item The mapping $\mathcal{I}:X_{w}^{0}\rightarrow L^{2}$ is an isometric
homeomorphism with inverse $\mathcal{J}$.

\item The mappings $\mathcal{I}$ and $\mathcal{J}$ are adjoints.

\item When $w$ has property W02 and $\left\vert \alpha\right\vert \leq\kappa$:%
\[
D^{a}\mathcal{J}g\left(  x\right)  =\frac{i^{\left\vert \alpha\right\vert }%
}{\left(  2\pi\right)  ^{d/2}}\int\left(  \frac{\xi^{\alpha}}{\sqrt{w}%
}\right)  ^{\vee}\left(  x-y\right)  g\left(  y\right)  dy,\text{\quad}g\in
L^{2}.
\]

\item When $w$ has property W03 and $\alpha\leq\kappa$:%
\[
D^{a}\mathcal{J}g\left(  x\right)  =\frac{i^{\left\vert \alpha\right\vert }%
}{\left(  2\pi\right)  ^{d/2}}\int\left(  \frac{\xi^{\alpha}}{\sqrt{w}%
}\right)  ^{\vee}\left(  x-y\right)  g\left(  y\right)  dy,\text{\quad}g\in
L^{2}.
\]

\end{enumerate}
\end{theorem}

\begin{proof}
The proofs are straight forward and will be omitted.
\end{proof}

The fact that $X_{w}^{0}$ is complete, and hence is a Hilbert space, is a
simple consequence of the last theorem:

\begin{corollary}
\label{Cor_Xow_complete}Suppose $w$ is a weight function with property W02 or
W03. Then $X_{w}^{0}$ is complete and hence a Hilbert space.
\end{corollary}

\begin{proof}
By part 6 of the previous Theorem \ref{Thm_I_J_property}, $\mathcal{I}%
:X_{w}^{0}\rightarrow L^{2}$ is an isometric homeomorphism. Hence $X_{w}^{0}$
is complete since $L^{2}$ is a Hilbert space.
\end{proof}

\begin{remark}
??? Suppose $\frac{1}{\sqrt{w}}L^{2}\subset S^{\prime}$ and $\left(  \frac
{1}{\sqrt{w}}L^{2}\right)  ^{\vee}\in L_{loc}^{1}$. Then it follows that
$\mathcal{J}:L^{2}\rightarrow X_{w}^{0}$ is an isometry and $X_{w}^{0}$ is a
Hilbert space. But is it a reproducing kernel Hilbert space?

\textbf{What if we assume that} $\frac{1}{\sqrt{w}}L^{2}\subset L^{2}$? Then
$f\in L^{2}$ implies $\frac{1}{\sqrt{w}}f\in L^{2}$ and so $f\in\overline
{X}_{1/w}^{0}\subset X_{1/w}^{0}$ i.e. $L^{2}\subset\overline{X}_{1/w}%
^{0}\subset X_{1/w}^{0}$.

Also, $g\in L^{2}$ implies $\left\Vert \mathcal{J}g\right\Vert _{w,0}^{2}=\int
w\left\vert \left(  \left(  \frac{\widehat{g}}{\sqrt{w}}\right)  ^{\vee
}\right)  ^{\wedge}\right\vert ^{2}=\int\left\vert \widehat{g}\right\vert
^{2}=\left\Vert g\right\Vert _{2}^{2}$.

$u\in L_{0}^{2}$ and $f\in X_{w}^{0}$ implies $\int\left\vert \widehat{u\ast
f}\right\vert =\int\left\vert \widehat{u}\right\vert \left\vert \widehat{f}%
\right\vert =\int\frac{1}{\sqrt{w}}\left\vert \widehat{u}\right\vert \sqrt
{w}\left\vert \widehat{f}\right\vert \leq\left\Vert \frac{\widehat{u}}%
{\sqrt{w}}\right\Vert _{2}\left\Vert f\right\Vert _{w,0}<\infty$.

Thus $u\ast f\in C_{B}^{0}$ $\forall u\in L_{0}^{2}$.

Also $\frac{1}{\sqrt{w}}L^{2}\subset L^{2}$ implies $L^{2}\subset\sqrt{w}%
L^{2}$.
\end{remark}

\subsection{The smoothness of functions in $X_{w}^{0}$%
\label{SbSect_Xw,th_smooth}}

We begin with a lemma of basic $L^{1}$ Fourier transform results.

\begin{lemma}
\label{Lem_L1_Fourier_contin}\ 

\begin{enumerate}
\item (Theorem 4.2 of Malliavin \cite{Malliavin95}) If $f\in S^{\prime}$ and
$\widehat{f}\in L^{1}$, then $f\in C_{B}^{\left(  0\right)  }$, $\left\vert
f\left(  x\right)  \right\vert \rightarrow0$ as $\left\vert x\right\vert
\rightarrow\infty$, and
\[
f\left(  x\right)  =\left(  2\pi\right)  ^{-d/2}\int e^{ix\xi}\widehat{f}%
\left(  \xi\right)  d\xi.
\]

\item (Corollary 2.12 of Petersen \cite{Petersen83}) If $f\in L^{1}$ and for
some integer $n>1$, $\left\vert \cdot\right\vert ^{n}f\in L^{1}$, then
$\overset{\vee}{f}\in C_{B}^{\left(  n\right)  }$. Further, $x^{\beta}f\in
L^{1}$ when $\left\vert \beta\right\vert \leq n$ and
\begin{equation}
D^{\beta}\left(  \overset{\vee}{f}\right)  \left(  x\right)  =\left(
i\right)  ^{\left\vert \beta\right\vert }\left(  x^{\beta}f\right)  ^{\vee
}\left(  x\right)  ,\text{\quad}\left\vert \beta\right\vert \leq
n.\label{1.031}%
\end{equation}

\end{enumerate}
\end{lemma}

The next theorem presents our first smoothness result for functions in
$X_{w}^{0}$.

\begin{theorem}
\label{Thm_X_smooth}Suppose that $w$ is a weight function. Then:

\begin{enumerate}
\item If $w$ has property W02 for parameter $\kappa\in\mathbb{R}^{1}$, then
$X_{w}^{0}\subset C_{B}^{\left(  \left\lfloor \kappa\right\rfloor \right)  }$.

\item If $w$ has property W03 for parameter $\kappa=\left(  \kappa_{i}\right)
$ then $X_{w}^{0}\subset C_{B}^{\left(  \left\lfloor \kappa\right\rfloor
\right)  }$.

Here $\left\lfloor \kappa\right\rfloor =\left(  \left\lfloor \kappa
_{i}\right\rfloor \right)  $ and we define for all multi-indexes $\alpha$:%
\begin{equation}
C_{B}^{\left(  \alpha\right)  }:=\left\{  u\in C_{B}^{\left(  0\right)
}:D^{\beta}u\in C_{B}^{\left(  0\right)  }\text{ }for\text{ }\beta\leq
\alpha\right\}  .\label{1.064}%
\end{equation}

\end{enumerate}
\end{theorem}

\begin{proof}
\textbf{Part 1} Since $w$ has property W02, by Theorem \ref{Thm_equiv_W2}, if
$\left\vert \beta\right\vert \leq\kappa$%
\begin{equation}%
{\displaystyle\int}
\dfrac{\left\vert \cdot\right\vert ^{2\left\vert \beta\right\vert }}{w}=%
{\displaystyle\int\limits_{\left\vert \cdot\right\vert \leq1}}
\dfrac{\left\vert \cdot\right\vert ^{2\left\vert \beta\right\vert }}{w}+%
{\displaystyle\int\limits_{\left\vert \cdot\right\vert \geq1}}
\dfrac{\left\vert \cdot\right\vert ^{2\left\vert \beta\right\vert }}{w}\leq%
{\displaystyle\int\limits_{\left\vert \cdot\right\vert \leq1}}
\dfrac{1}{w}+%
{\displaystyle\int\limits_{\left\vert \cdot\right\vert \geq1}}
\dfrac{\left\vert \cdot\right\vert ^{2\kappa}}{w}<\infty.\label{1.50}%
\end{equation}

Hence when $f\in X_{w}^{0}$ the Cauchy-Schwartz inequality yields
\[
\int\left\vert \xi^{\beta}\widehat{f}\right\vert \leq\int\frac{\left\vert
\cdot\right\vert ^{\left\vert \beta\right\vert }}{\sqrt{w}}\sqrt{w}\left\vert
\widehat{f}\right\vert \leq\left(  \int\frac{\left\vert \cdot\right\vert
^{2\left\vert \beta\right\vert }}{w}\right)  ^{1/2}\left\Vert f\right\Vert
_{w,0}<\infty,
\]

i.e. $f\in X_{w}^{0}$ implies $\widehat{D^{\beta}f}\in L^{1}$ when $\left\vert
\beta\right\vert \leq\kappa$ and hence $f\in C_{B}^{\left(  \left\lfloor
\kappa\right\rfloor \right)  }$ by part 1 of Lemma \ref{Lem_L1_Fourier_contin}%
.\medskip

\textbf{Part 2} If $f\in X_{w}^{0}$ and $\beta\leq\kappa$ then the
Cauchy-Schwartz inequality yields
\begin{equation}
\int\left\vert \xi^{\beta}\widehat{f}\right\vert \leq\int\frac{\left\vert
\xi^{\beta}\right\vert }{\sqrt{w}}\sqrt{w}\left\vert \widehat{f}\right\vert
\leq\left(  \int\frac{\xi^{2\beta}}{w}\right)  ^{1/2}\left\Vert f\right\Vert
_{w,0}<\infty,\label{2.0}%
\end{equation}

i.e. $f\in X_{w}^{0}$ implies $\widehat{D^{\beta}f}\in L^{1}$ when $\beta
\leq\kappa$ and hence $f\in C_{B}^{\left(  \left\lfloor \kappa\right\rfloor
\right)  }$ by part 1 of Lemma \ref{Lem_L1_Fourier_contin}.
\end{proof}

The next theorem derives some inverse Fourier transform formulas and
inequalities for derivatives of functions in $X_{w}^{0}$.

\begin{theorem}
\label{Thm_X_invers_Fourier_W2}Suppose that $w$ is a weight function with
property W02 for parameter $\kappa$. Then for $f\in X_{w}^{0}$ we have the
inverse Fourier transform formulas
\begin{equation}
D^{\beta}f\left(  x\right)  =\left(  2\pi\right)  ^{-d/2}\int e^{ix\xi
}\widehat{D^{\beta}f}\left(  \xi\right)  d\xi,\text{\quad}\left\vert
\beta\right\vert \leq\kappa.\label{1.021}%
\end{equation}

and the derivatives satisfy
\begin{equation}
\left\vert D^{\beta}f\left(  x\right)  \right\vert \leq\left(  2\pi\right)
^{-d/2}\left(  \int\frac{\left\vert \cdot\right\vert ^{2\left\vert
\beta\right\vert }}{w}\right)  ^{1/2}\left\Vert f\right\Vert _{w,0}%
,\text{\quad}\left\vert \beta\right\vert \leq\kappa,\label{1.022}%
\end{equation}

so that $X_{w}^{0}\hookrightarrow C_{B}^{\left(  \left\lfloor \kappa
\right\rfloor \right)  }$ when $C_{B}^{\left(  \left\lfloor \kappa
\right\rfloor \right)  }$ is endowed with the supremum norm $\sum
\limits_{\left\vert \beta\right\vert \leq\kappa}\left\vert D^{\beta
}f\right\vert _{\sup}$.
\end{theorem}

\begin{proof}
From the proof of the previous theorem we have: $f\in X_{w}^{0}$ implies
$\widehat{D^{\beta}f}\in L^{1}$ when $\left\vert \beta\right\vert \leq\kappa$.
The inverse Fourier transform formulas \ref{1.021} are now simple consequences
of part 1 of Lemma \ref{Lem_L1_Fourier_contin}.

The bounds on the derivatives are derived from \ref{1.021} using the
Cauchy-Schwartz inequality:
\[
\left\vert D^{\beta}f\left(  x\right)  \right\vert \leq\left(  2\pi\right)
^{-\frac{d}{2}}\int\left\vert \xi^{\beta}\right\vert \left\vert \widehat{f}%
\left(  \xi\right)  \right\vert d\xi\leq\left(  2\pi\right)  ^{-\frac{d}{2}%
}\int\frac{\left\vert \cdot\right\vert ^{\left\vert \beta\right\vert }}%
{\sqrt{w}}\sqrt{w}\left\vert \widehat{f}\right\vert \leq\left(  2\pi\right)
^{-\frac{d}{2}}\left(  \int\frac{\left\vert \cdot\right\vert ^{2\left\vert
\beta\right\vert }}{w}\right)  ^{\frac{1}{2}}\left\Vert f\right\Vert _{w,0}.
\]

\end{proof}

\begin{theorem}
\label{Thm_X_invers_Fourier_W3}Suppose that $w$ is a weight function with
property W03 for parameter $\kappa\in\mathbb{R}^{d}$. Then for $f\in X_{w}%
^{0}$ we have the inverse Fourier transform formulas
\begin{equation}
D^{\beta}f\left(  x\right)  =\left(  2\pi\right)  ^{-\frac{d}{2}}\int
e^{ix\xi}\widehat{D^{\beta}f}\left(  \xi\right)  d\xi,\text{\quad}\beta
\leq\kappa.\label{1.065}%
\end{equation}

and the derivatives satisfy
\begin{equation}
\left\vert D^{\beta}f\left(  x\right)  \right\vert \leq\left(  2\pi\right)
^{-\frac{d}{2}}\left(  \int\frac{\xi^{2\beta}}{w}\right)  ^{1/2}\left\Vert
f\right\Vert _{w,0},\text{\quad}\beta\leq\kappa,\label{1.066}%
\end{equation}

so that $X_{w}^{0}\hookrightarrow C_{B}^{\left(  \left\lfloor \kappa
\right\rfloor \right)  }$ when $C_{B}^{\left(  \left\lfloor \kappa
\right\rfloor \right)  }$ is endowed with the supremum norm $\sum
\limits_{\beta\leq\kappa}\left\vert D^{\beta}f\right\vert _{\sup}$.
\end{theorem}

\begin{proof}
From the proof of the previous theorem we have: $f\in X_{w}^{0}$ implies
$\widehat{D^{\beta}f}\in L^{1}$ when $\beta\leq\kappa$. The inverse Fourier
transform formulas \ref{1.065} are now simple consequences of part 1 of Lemma
\ref{Lem_L1_Fourier_contin}.

The bounds on the derivatives are derived from \ref{1.065} using the
Cauchy-Schwartz inequality:
\[
\left\vert D^{\beta}f\left(  x\right)  \right\vert \leq\left(  2\pi\right)
^{-\frac{d}{2}}\int\left\vert \xi^{\beta}\right\vert \left\vert \widehat{f}%
\left(  \xi\right)  \right\vert d\xi\leq\left(  2\pi\right)  ^{-\frac{d}{2}%
}\int\frac{\left\vert \xi^{\beta}\right\vert }{\sqrt{w}}\sqrt{w}\left\vert
\widehat{f}\right\vert \leq\left(  2\pi\right)  ^{-\frac{d}{2}}\left(
\int\frac{\xi^{2\beta}}{w}\right)  ^{\frac{1}{2}}\left\Vert f\right\Vert
_{w,0}.
\]

\end{proof}

\begin{remark}
\label{Rem_reprod_Hilbert}Inequalities \ref{1.022} and \ref{1.066} for
$\beta=0$ implies that $X_{w}^{0}$ is a reproducing kernel Hilbert space and
in Section \ref{Sect_reps_of_eval_fns_W2} we will calculate the Riesz
representers of the evaluation functionals $f\rightarrow D^{\alpha}f\left(
x\right)  $. We have in fact established \textbf{an important link between
pointwise processes and Hilbert space theory}.
\end{remark}

\subsection{Examples\label{Ex_data_space}}

\begin{example}
\textbf{The shifted thin-plate splines} From Examples
\ref{SbSect_wt_func_examples}, $w$ has property W02 for all $\kappa\geq0$, and
thus Theorem \ref{Thm_X_invers_Fourier_W2} implies $X_{w}^{0}\subset
C_{B}^{\infty}$ as sets. From \ref{1.19} we can write%
\[
w\left(  \xi\right)  =\frac{1}{\widetilde{e}\left(  v\right)  }\frac
{\left\vert \xi\right\vert ^{s}}{\widetilde{K}_{s}\left(  \left\vert
\xi\right\vert \right)  },\quad s=2v+d>0,
\]

where from Theorem \ref{Thm_bnds_modif_MacDonald},%
\[
c_{s}e^{-t}\leq\widetilde{K}_{s}\left(  t\right)  \leq c_{s}^{\prime}%
e^{-t},\quad s>0,\text{ }t\geq0.
\]

so that%
\[
\frac{1}{\widetilde{e}\left(  v\right)  c_{s}^{\prime}}\left\vert
\xi\right\vert ^{s}e^{\left\vert \xi\right\vert }\leq w\left(  \xi\right)
\leq\frac{1}{\widetilde{e}\left(  v\right)  c_{s}}\left\vert \xi\right\vert
^{s}e^{\left\vert \xi\right\vert },\quad\xi\in\mathbb{R}^{d}.
\]

Now expanding the exponentials about the origin using Taylor series we have
for $u\in X_{w}^{0}$,%
\[
\frac{1}{\widetilde{e}\left(  v\right)  c_{s}^{\prime}}\int\left\vert
\cdot\right\vert ^{s}\sum_{k=0}^{\infty}\frac{\left\vert \cdot\right\vert
^{k}}{k!}\left\vert \widehat{u}\right\vert ^{2}\leq\int w\left\vert
\widehat{u}\right\vert ^{2}\leq\frac{1}{\widetilde{e}\left(  v\right)  c_{s}%
}\int\left\vert \cdot\right\vert ^{s}\sum_{k=0}^{\infty}\frac{\left\vert
\cdot\right\vert ^{k}}{k!}\left\vert \widehat{u}\right\vert ^{2},
\]

which means that as sets%
\begin{equation}
X_{w}^{0}=\bigcap\limits_{k=0}^{\infty}\left\{  u\in S^{\prime}:\widehat{u}\in
L_{loc}^{1}\text{ }and\text{ }\int\left\vert \cdot\right\vert ^{s+k}\left\vert
\widehat{u}\right\vert ^{2}<\infty\right\}  .\label{a7.29}%
\end{equation}

\end{example}

\begin{example}
\textbf{The Gaussian }It was shown in Example \ref{Ex_Gaussian_wt} that the
Gaussian weight function has property W02 for all $\kappa\geq0$, and so
Theorem \ref{Thm_X_smooth} implies $X_{w}^{0}\subset C_{B}^{\infty}$ as sets.
Further%
\[
X_{w}^{0}=\left\{  u\in S^{\prime}:\widehat{u}\in L_{loc}^{1}\text{ }and\text{
}\int e^{\left\vert \cdot\right\vert ^{2}}\left\vert \widehat{u}\right\vert
^{2}<\infty\right\}  ,
\]

so that by expanding the Gaussian about the origin using a Taylor series (see
previous example) we can show easily that as sets
\[
X_{w}^{0}\subset W^{\infty}=\bigcap\limits_{n=0}^{\infty}W^{n},
\]

where $W^{n}$ is the Sobolev space of order $n$ given below in Definition
\ref{Def_SobolevSpace}.
\end{example}

\begin{example}
\label{Ex_Xow_example_SobSplines}\textbf{The Sobolev splines} It was shown in
Example \ref{Ex_Sobolev_splin_wt} that property W02 holds for $0\leq
\kappa<v-d/2$. Now by Remark \ref{Rem_reprod_Hilbert}, $X_{w}^{0}%
\hookrightarrow C_{B}^{\left(  \max\left\lfloor \kappa\right\rfloor \right)
}$ and we see that $\max\left\lfloor \kappa\right\rfloor =\left\lfloor
v-d/2\right\rfloor $ when $v-d/2$ is not an integer and $\max\left\lfloor
\kappa\right\rfloor =\left\lfloor v-d/2-1\right\rfloor =\left\lfloor
v-d/2\right\rfloor -1$ when $v-d/2$ is an integer.

Also, from \ref{1.030}, $w\left(  \xi\right)  =\left(  1+\left\vert
\xi\right\vert ^{2}\right)  ^{v}$ and so \ref{1.033} implies
\[
X_{w}^{0}=W^{v},
\]
where $W^{v}$ is the Sobolev space of positive order $v$ defined using the
Fourier transform.
\end{example}

\begin{example}
These splines are defined in Theorem \ref{Thm_ex_splin_wt_fn_properties}.
Since property W03 holds for all $0\leq\kappa<\left(  n-1/2\right)
\mathbf{1}$ it follows that $\max\kappa_{i}=\left(  n-1\right)  \mathbf{1}$,
and thus $X_{w}^{0}\hookrightarrow C_{B}^{\left(  \left(  n-1\right)
\mathbf{1}\right)  } $.

In Lemma \ref{Lem_Thm_Xow(O)_eq_Hn(O)_dim1} we will show that $X_{w}%
^{0}\hookrightarrow W^{n\mathbf{1}}$. Here $W^{n\mathbf{1}}$ is the Sobolev
space introduced in Definition \ref{Def_SobolevSpace2} below and it consists
of all $L^{2}$ functions such that $D^{\alpha}f\in L^{2}$ when each
$\alpha\leq n\mathbf{1}$.
\end{example}

\subsection{Some dense $C^{\infty}$ subspaces of $X_{w}^{0}$}

The results of this subsection are currently not used in this document.

\subsubsection{\protect\underline{Results to handle the weight function
discontinuity set}}

In this section we prove some results needed when we want to use a partition
of unity to handle the discontinuity or zero values of a weight function on a
set of measure zero. Recall that weight function property W01 introduced in
Definition \ref{Def_extend_wt_fn} required that a weight function be
continuous and positive outside a closed set $\mathcal{A}$ of measure zero.
Light and Wayne assumed their basis functions were continuous and positive
outside the origin i.e. $\mathcal{A}=\left\{  0\right\}  $. In order to allow
hat functions to be basis functions we have had to allow the weight set to a
closed unbounded set of measure zero.

The discontinuity of the weight function needs to be taken into account when
we prove the density of some $C^{\infty}$ subspaces in $L^{2}$ and $X_{w}^{0}%
$, and when we prove the smoothness of basis functions generated by weight
functions which have property W02 or W03.

In this document we will use the next lemma whenever we want to construct a
partition of unity using points `near' a weight function set. This lemma will
be used with $\mathcal{F}=\mathcal{A}$, where $\mathcal{A}$ is the weight
function set. Lemma \ref{Lem_Ae_intersect_sphere} will be applied to a
function of the `near' points. Nearness can be measured using the concept of a
\textit{(radial)} \textit{neighborhood of a set}: if $\mathcal{F}%
\subset\mathbb{R}^{d}$ and $\eta\in\mathbb{R}^{1}$, $\eta>0$ then define
$\mathcal{F}_{\eta}=\bigcup\limits_{x\in A}B\left(  x;\eta\right)  $. The set
$\mathcal{F}_{\eta}$ is referred to as the (radial)\textbf{\ }$\eta
-$neighborhood of the set $\mathcal{F}$.

\begin{lemma}
\label{Lem_func_eq_1_nbhd_set}(based on Lemma 1, \S 5.2 of Vladimirov
\cite{Vladimirov}) Let $\mathcal{F}$ be any set of points in $\mathbb{R}^{d}$.
Then for any $\eta\in\mathbb{R}_{+}^{1}$ there exists a real-valued function
$f_{\eta}\in C^{\infty}\left(  \mathbb{R}^{d}\right)  $ such that:

\begin{enumerate}
\item $0\leq f_{\eta}\left(  x\right)  \leq1$,

\item $f_{\eta}\left(  x\right)  =1$ when $x\in\mathcal{F}_{\eta}$,

\item $f_{\eta}\left(  x\right)  =0$ when $x\notin\mathcal{F}_{3\eta}$.
\end{enumerate}

Further, if we emphasize the dependency on $\mathcal{F}$ by writing
$f_{\mathcal{F};\eta}=f_{\eta}$, we have for scalar dilations and
translations:
\begin{equation}
f_{\mathcal{F};\eta}\left(  x/\lambda\right)  =f_{\lambda\mathcal{F}%
;\lambda\eta}\left(  x\right)  ,\quad\lambda\in\mathbb{R}_{+}^{1}%
,\label{a1.62}%
\end{equation}

and%
\begin{equation}
f_{\mathcal{F};\eta}\left(  x-c\right)  =f_{\mathcal{F}+c;\eta}\left(
x\right)  ,\quad c\in\mathbb{R}^{d}.\label{a1.63}%
\end{equation}

\end{lemma}

\begin{proof}
There exists a mollifier $\omega\in C_{0}^{\infty}$ satisfying
\[
\operatorname*{supp}\omega\subset B\left(  0;1\right)  ;\quad\int%
\omega=1;\quad\omega\geq0.
\]

Let $\chi_{\Omega}\left(  x\right)  $ denote the \textbf{characteristic
function} of an arbitrary set $\Omega\subset\mathbb{R}^{d}$. Now define
$\omega_{\eta}\left(  x\right)  =\eta^{-d}\omega\left(  x/\eta\right)  $. Then
it is shown in Lemma 1, \S 5.2 of Vladimirov \cite{Vladimirov} that the
function
\[
f_{\eta}\left(  x\right)  =\int\chi_{\mathcal{F}_{2\eta}}\left(  y\right)
\omega_{\eta}\left(  x-y\right)  dy=\int_{\mathcal{F}_{2\eta}}\omega_{\eta
}\left(  x-y\right)  dy,
\]

has the required properties.

Now regarding dilations: since the mollifier satisfies%
\[
\omega_{\eta}\left(  x/\lambda\right)  =\eta^{-d}\omega\left(  x/\lambda
\eta\right)  =\lambda^{d}\left(  \lambda\eta\right)  ^{-d}\omega\left(
x/\lambda\eta\right)  =\lambda^{d}\omega_{\lambda\eta}\left(  x\right)  ,
\]

we have%
\begin{align*}
f_{\eta}\left(  x/\lambda\right)  =\int\chi_{\mathcal{F}_{2\eta}}\left(
y\right)  \omega_{\eta}\left(  \frac{x}{\lambda}-y\right)  dy  & =\int%
\chi_{\mathcal{F}_{2\eta}}\left(  y\right)  \omega_{\eta}\left(
\frac{x-\lambda y}{\lambda}\right)  dy\\
& =\lambda^{-d}\int\chi_{\mathcal{F}_{2\eta}}\left(  \frac{z}{\lambda}\right)
\omega_{\eta}\left(  \frac{x-z}{\lambda}\right)  dz\\
& =\lambda^{-d}\int\chi_{\mathcal{F}_{2\eta}}\left(  \frac{z}{\lambda}\right)
\lambda^{d}\omega_{\lambda\eta}\left(  x-z\right)  dz\\
& =\int\chi_{\mathcal{F}_{2\eta}}\left(  \frac{z}{\lambda}\right)
\omega_{\lambda\eta}\left(  x-z\right)  dz,
\end{align*}

and because $\left(  \chi_{\mathcal{F}_{2\eta}}\right)  \left(  \frac
{z}{\lambda}\right)  =\left(  \chi_{\lambda\mathcal{F}_{2\eta}}\right)
\left(  z\right)  =\left(  \chi_{\left(  \lambda\mathcal{F}\right)
_{2\lambda\eta}}\right)  \left(  z\right)  $,%
\[
f_{\mathcal{F};\eta}\left(  x/\lambda\right)  =\int\chi_{\left(
\lambda\mathcal{F}\right)  _{2\lambda\eta}}\left(  z\right)  \omega
_{\lambda\eta}\left(  x-z\right)  dz=f_{\lambda\mathcal{F};\lambda\eta}\left(
x\right)  .
\]

Regarding translations:%
\begin{align*}
f_{\mathcal{F};\eta}\left(  x-c\right)  =\int\chi_{\mathcal{F}_{2\eta}}\left(
y\right)  \omega_{\eta}\left(  x-c-y\right)  dy  & =\int\chi_{\mathcal{F}%
_{2\eta}}\left(  z-c\right)  \omega_{\eta}\left(  x-z\right)  dz\\
& =\int\chi_{\left(  \mathcal{F+}c\right)  _{2\eta}}\left(  z\right)
\omega_{\eta}\left(  x-z\right)  dz\\
& =f_{\mathcal{F+}c;\eta}\left(  x\right)  .
\end{align*}

\end{proof}

The next theorem will require the following standard measure theory results
which are stated without proof.

\begin{lemma}
\label{Lem_Ae_intersect_sphere}\ 

\begin{enumerate}
\item Suppose the set $\mathcal{A}$ is closed and has measure zero. If
$\mathcal{A}_{\varepsilon}$ is the $\varepsilon-$neighborhood of the set
$\mathcal{A}$ then for any open ball $B$, $\operatorname*{meas}\left(
B\cap\mathcal{A}_{\varepsilon}\right)  \rightarrow0$ as $\varepsilon
\rightarrow0$.

\item Suppose the set $\mathcal{A}$ is closed and has measure zero.

Then $f\in L^{1}$ implies $\lim\limits_{\varepsilon\rightarrow0}%
\int_{\mathcal{A}_{\varepsilon}}\left\vert f\right\vert =0$.
\end{enumerate}
\end{lemma}

\subsubsection{\protect\underline{Some spaces of continuous functions which
are dense in $L^{2}$}}

The next theorem is an extension of Theorem 2.5 of Light and Wayne
\cite{LightWayneX98Weight}. The difference is that instead of the set
$\left\{  0\right\}  $, we are dealing with a set $\mathcal{A}$ which is
closed and has measure zero. There is no need to assume that $0\in\mathcal{A}%
$. To handle this set we use Lemma \ref{Lem_func_eq_1_nbhd_set} and Lemma
\ref{Lem_Ae_intersect_sphere}.

\begin{theorem}
\label{Thm_sqrtwCoo_dense_L2}Suppose $w$ is a weight function with property
W01 with respect to the set $\mathcal{A}$.

Then the set $\sqrt{w}C_{0}^{\left(  0\right)  }\cap L^{2}=\left\{  \sqrt
{w}f:f\in C_{0}^{\left(  0\right)  }\right\}  \cap L^{2}$ is dense in $L^{2}$.

Here $C_{0}^{\left(  0\right)  }$ denotes the continuous functions with
compact support.
\end{theorem}

\begin{proof}
By Lemma \ref{Lem_func_eq_1_nbhd_set}, given $h>0$ there exists a function
$\psi_{h}\in C^{\infty}$ satisfying $0\leq\psi_{h}\leq1$, $\psi_{h}=1$ on
$\mathcal{A}_{h}$ and $\psi_{h}=0$ outside $\mathcal{A}_{3h}$. Now define the
mapping $\Psi_{h}:$ $C_{0}^{\left(  0\right)  }\rightarrow C_{0}^{\left(
0\right)  }$ by
\[
\Psi_{h}g=\left(  1-\psi_{h}\right)  g,\text{\quad}g\in C_{0}^{\left(
0\right)  },\text{ }h>0.
\]

Since $w$ is positive and continuous on $\mathbb{R}^{d}\setminus\mathcal{A}$,
we have $\frac{1-\psi_{h}}{\sqrt{w}}\in C^{\left(  0\right)  }$ and hence
$\left(  1-\psi_{h}\right)  g\in\sqrt{w}C_{0}^{\left(  0\right)  }\cap L^{2}$.
In other words, $\Psi_{h}:C_{0}^{\left(  0\right)  }\rightarrow\sqrt{w}%
C_{0}^{\left(  0\right)  }\cap L^{2}$.

Now $C_{0}^{\left(  0\right)  }$ is dense in $L^{2}$ e.g. Theorem 2.13 Adams
\cite{Adams75}, and given $g_{0}\in L^{2}$ and $\varepsilon>0$ we can choose
$g_{\varepsilon}\in C_{0}^{\left(  0\right)  }$ such that $\left\Vert
g_{0}-g_{\varepsilon}\right\Vert _{L^{2}}<\varepsilon/2$. Then
\[
\left\Vert g_{0}-\Psi_{h}g_{\varepsilon}\right\Vert _{L^{2}}\leq\left\Vert
g_{0}-g_{\varepsilon}\right\Vert _{L^{2}}+\left\Vert g_{\varepsilon}-\Psi
_{h}g_{\varepsilon}\right\Vert _{L^{2}}\leq\varepsilon/2+\left\Vert \psi
_{h}g_{\varepsilon}\right\Vert _{L^{2}}\leq\varepsilon/2+\left(
\int\limits_{\mathcal{A}_{3h}}\left\vert g_{\varepsilon}\right\vert
^{2}\right)  ^{\frac{1}{2}}.
\]

Finally, since $\left\vert g_{\varepsilon}\right\vert ^{2}\in L^{1}$, it
follows from Lemma \ref{Lem_Ae_intersect_sphere} that $\lim
\limits_{h\rightarrow0}\int_{\mathcal{A}_{3h}}\left\vert g_{\varepsilon
}\right\vert ^{2}=0$. Thus we can choose $h>0$ so that $\left\Vert g_{0}%
-\Psi_{h}g_{\varepsilon}\right\Vert _{L^{2}}\leq\varepsilon$.
\end{proof}

The next theorem improves the previous theorem. This is Light and Wayne's
Proposition 2.7 and our proof is a variation of Light and Wayne's proof.

\begin{theorem}
\label{Thm_sqrtwCoinf_dense_L2}Suppose $w$ is a weight function with respect
to the set $\mathcal{A}$.

Then the set $\sqrt{w}C_{0}^{\infty}\cap L^{2}=\left\{  \sqrt{w}f:f\in
C_{0}^{\infty}\right\}  \cap L^{2}$ is dense in $L^{2}$.
\end{theorem}

\begin{proof}
Now $\sqrt{w}C_{0}^{\infty}\cap L^{2}\subset\sqrt{w}C_{0}^{\left(  0\right)
}\cap L^{2}$. Hence, by Theorem \ref{Thm_sqrtwCoo_dense_L2}, if it can be
shown that $\sqrt{w}C_{0}^{\infty}\cap L^{2}$ is dense in $\sqrt{w}%
C_{0}^{\left(  0\right)  }\cap L^{2}$, then it follows that this theorem is true.

First note that $C_{0}^{\infty}$ is dense in $L^{2}$ and hence dense in
$C_{0}^{\left(  0\right)  }$. Next select $\varepsilon>0$ and define a mapping
$\Theta_{\varepsilon}:C_{0}^{\left(  0\right)  }\rightarrow C_{0}^{\infty}$
where $\Theta_{\varepsilon}f$ is an element of $C_{0}^{\infty}$ such that
$\left\Vert f-\Theta_{\varepsilon}f\right\Vert _{L^{2}}<\varepsilon$.

We want to use $\Theta_{\varepsilon}$ to construct a mapping (see left side of
\ref{a9.1}) from $\sqrt{w}C_{0}^{\left(  0\right)  }\cap L^{2}$ to $\sqrt
{w}C_{0}^{\infty}\cap L^{2}$, which can be used to prove that $\sqrt{w}%
C_{0}^{\infty}\cap L^{2}$ is dense in $\sqrt{w}C_{0}^{\left(  0\right)  }\cap
L^{2}$. If $f\in\sqrt{w}C_{0}^{\left(  0\right)  }\cap L^{2}$, then $\frac
{f}{\sqrt{w}}\in C_{0}^{\left(  0\right)  }$ and $\sqrt{w}\Theta_{\delta
}\left(  \frac{f}{\sqrt{w}}\right)  \in\sqrt{w}C_{0}^{\infty}$, $\delta>0$.
But we still need $\sqrt{w}\Theta_{\delta}\left(  \frac{f}{\sqrt{w}}\right)
\in L^{2}$. To do this we use the function $\psi_{h}$ defined in the previous
Theorem \ref{Thm_sqrtwCoo_dense_L2}. Since $\left(  1-\psi_{h}\right)
\sqrt{w}\in C^{\left(  0\right)  }$ we have
\[
\left(  1-\psi_{h}\right)  \sqrt{w}\Theta_{\delta}\left(  \frac{f}{\sqrt{w}%
}\right)  \in\sqrt{w}C_{0}^{\infty}\cap L^{2},\text{ }when\text{ }f\in\sqrt
{w}C_{0}^{\left(  0\right)  }\cap L^{2}.
\]

Finally, it will turn out that the support of $\left(  1-\psi_{h}\right)
\sqrt{w}\Theta_{\delta}\left(  \frac{f}{\sqrt{w}}\right)  $ needs to be
restricted by multiplying it by a function $\phi_{R}\in C_{0}^{\infty}$,
$R\geq1$, such that $0\leq\phi_{R}\leq1$, $\operatorname*{supp}\phi_{R}%
\subset\left[  -R-1,R+1\right]  $ and $\phi_{R}=1$ on $\left[  -R,R\right]  $.

We now assert that there exist $R$, $h$ and $\delta$ such that

$\left\Vert \left(  1-\psi_{h}\right)  \phi_{R}\sqrt{w}\Theta_{\delta}\left(
\frac{f}{\sqrt{w}}\right)  -f\right\Vert _{L^{2}}<\varepsilon$. To prove this
assertion write
\begin{align*}
\left(  1-\psi_{h}\right)  \phi_{R}\sqrt{w}\Theta_{\delta}\left(  \frac
{f}{\sqrt{w}}\right)  -f  & =\left(  1-\psi_{h}\right)  \phi_{R}\sqrt{w}%
\Theta_{\delta}\left(  \frac{f}{\sqrt{w}}\right)  -\left(  1-\psi_{h}\right)
\phi_{R}\sqrt{w}\left(  \frac{f}{\sqrt{w}}\right)  +\\
& \qquad+\left(  1-\psi_{h}\right)  \phi_{R}\sqrt{w}\left(  \frac{f}{\sqrt{w}%
}\right)  -\phi_{R}f+\phi_{R}f-f\\
& =\left(  1-\psi_{h}\right)  \phi_{R}\sqrt{w}\left(  \Theta_{\delta}\left(
\frac{f}{\sqrt{w}}\right)  -\frac{f}{\sqrt{w}}\right)  +\psi_{h}\phi
_{R}f+\left(  \phi_{R}-1\right)  f.
\end{align*}

Observe that $\left(  1-\psi_{h}\right)  \sqrt{w}\in C^{\left(  0\right)  }$
implies $\left(  1-\psi_{h}\right)  \phi_{R}\sqrt{w}\in C_{0}^{\left(
0\right)  }$. Hence
\begin{align*}
\left\Vert \left(  1-\psi_{h}\right)  \phi_{R}\sqrt{w}\Theta_{\delta}\left(
\frac{f}{\sqrt{w}}\right)  -f\right\Vert _{L^{2}}  & \leq\left\Vert \left(
1-\psi_{h}\right)  \phi_{R}\sqrt{w}\left(  \Theta_{\delta}\left(  \frac
{f}{\sqrt{w}}\right)  -\frac{f}{\sqrt{w}}\right)  \right\Vert _{L^{2}}+\\
& \qquad+\left\Vert \psi_{h}f\right\Vert _{L^{2}}+\left\Vert \left(  \phi
_{R}-1\right)  f\right\Vert _{L^{2}}\\
& \leq\left\Vert \left(  1-\psi_{h}\right)  \phi_{R}\sqrt{w}\right\Vert
_{\infty}\left\Vert \Theta_{\delta}\left(  \frac{f}{\sqrt{w}}\right)
-\frac{f}{\sqrt{w}}\right\Vert _{L^{2}}+\\
& \qquad+\left\Vert \psi_{h}f\right\Vert _{L^{2}}+\left\Vert \left(  \phi
_{R}-1\right)  f\right\Vert _{L^{2}}\\
& \leq\left\Vert \left(  1-\psi_{h}\right)  \phi_{R}\sqrt{w}\right\Vert
_{\infty}\delta+\left\Vert \psi_{h}f\right\Vert _{L^{2}}+\left\Vert \left(
\phi_{R}-1\right)  f\right\Vert _{L^{2}}.
\end{align*}

We will consider each term on the right side of the last inequality
separately. Regarding the last term,

$\left\Vert \left(  \phi_{R}-1\right)  f\right\Vert _{L^{2}}=\left(
\int_{\left\vert \cdot\right\vert \geq R}\left\vert f\right\vert ^{2}\right)
^{1/2}$ and so $R$ can be fixed so that $\left\Vert \left(  \phi_{R}-1\right)
f\right\Vert _{L^{2}}\leq\varepsilon/3$. Next, by definition of $\psi_{h}$,
$\left\Vert \psi_{h}f\right\Vert _{L^{2}}=\left(  \int_{\mathcal{A}_{3h}%
}\left\vert f\right\vert ^{2}\right)  ^{1/2}$ and, since $\left\vert
f\right\vert ^{2}\in L^{1}$, it follows from Lemma
\ref{Lem_Ae_intersect_sphere} that $\lim\limits_{h\rightarrow0}\int%
_{\mathcal{A}_{3h}}\left\vert f\right\vert ^{2}=0$. A value of $h$ can now be
chosen so that $\left\Vert \psi_{h}f\right\Vert _{L^{2}}\leq\varepsilon/3$.
Now we have
\begin{equation}
\left\Vert \left(  1-\psi_{h}\right)  \phi_{R}\sqrt{w}\Theta_{\delta}\left(
\frac{f}{\sqrt{w}}\right)  -f\right\Vert _{L^{2}}\leq\left\Vert \left(
1-\psi_{h}\right)  \phi_{R}\sqrt{w}\right\Vert _{\infty}\delta+2\varepsilon
/3,\label{a9.1}%
\end{equation}

for all $f\in\sqrt{w}C_{0}^{\left(  0\right)  }\cap L^{2}.$

The last step is to choose $\delta$ so that $\left\Vert \left(  1-\psi
_{h}\right)  \phi_{R}\sqrt{w}\right\Vert _{\infty}\delta\leq\varepsilon/3 $,
and the theorem follows.
\end{proof}

The next result corresponds to Light and Wayne's Corollary 2.8. Here
$X_{w}^{0}$ is Light and Wayne's space $Y$ and $\widehat{C_{0}^{\infty}}$
denotes the space of Fourier transforms of all functions in $C_{0}^{\infty}$.

\begin{corollary}
\label{Cor_Xow_S_density}Suppose the weight function $w$ also has property W02
or W03.

Then the spaces $X_{w}^{0}\cap\left(  C_{0}^{\infty}\right)  ^{\vee}$,
$X_{w}^{0}\cap\widehat{C_{0}^{\infty}}$ and $X_{w}^{0}\cap S$ are all dense in
$X_{w}^{0}$.
\end{corollary}

\begin{proof}
Since the Fourier transform is a homeomorphism form $L^{2}$ to $L^{2}$ it
follows from Theorem \ref{Thm_sqrtwCoinf_dense_L2} that $\left(  \sqrt{w}%
C_{0}^{\infty}\cap L^{2}\right)  ^{\vee}$ is dense in $L^{2}$.

Now since $w$ has property W02 or W03 by Theorem \ref{Thm_I_J_property} the
operator $\mathcal{J}:L^{2}\rightarrow X_{w}^{0}$ of Definition \ref{Def_I_J}
is a homeomorphism with inverse denoted by $\mathcal{I}$, and we now show
that
\begin{equation}
\mathcal{J}:\left(  \sqrt{w}C_{0}^{\infty}\cap L^{2}\right)  ^{\vee
}\rightarrow X_{w}^{0}\cap\left(  C_{0}^{\infty}\right)  ^{\vee},\label{a9.0}%
\end{equation}

is onto which means that $X_{w}^{0}\cap\left(  C_{0}^{\infty}\right)  ^{\vee}$
is dense in $X_{w}^{0}$. The key result is%
\[%
\begin{array}
[c]{lll}%
u\in\left(  \sqrt{w}C_{0}^{\infty}\cap L^{2}\right)  ^{\vee} & iff &
\widehat{u}\in\sqrt{w}C_{0}^{\infty}\cap L^{2}\\
& iff & \widehat{u}\in\sqrt{w}C_{0}^{\infty}\text{ }and\text{ }\widehat{u}\in
L^{2}\\
& iff & u\in\left(  \sqrt{w}C_{0}^{\infty}\right)  ^{\vee}\text{ }and\text{
}u\in L^{2}\\
& iff & u\in\left(  \sqrt{w}C_{0}^{\infty}\right)  ^{\vee}\cap L^{2}.
\end{array}
\]

Now $\mathcal{J}u=\left(  \widehat{u}/\sqrt{w}\right)  ^{\vee}\in X_{w}%
^{0}\cap\left(  C_{0}^{\infty}\right)  ^{\vee}$. Further, if $v\in X_{w}%
^{0}\cap\left(  C_{0}^{\infty}\right)  ^{\vee}$ then
\[
\mathcal{J}^{-1}v=\mathcal{I}v=\left(  \sqrt{w}\widehat{u}\right)  ^{\vee}%
\in\left(  \sqrt{w}C_{0}^{\infty}\right)  ^{\vee}\cap L^{2}=\left(  \sqrt
{w}C_{0}^{\infty}\cap L^{2}\right)  ^{\vee},
\]

and so \ref{a9.0} is onto.

Finally, $\widehat{C_{0}^{\infty}}=\left(  C_{0}^{\infty}\right)  ^{\vee
}\subset S$ so $X_{w}^{0}\cap\widehat{C_{0}^{\infty}}$ and $X_{w}^{0}\cap S $
are dense in $X_{w}^{0}$.
\end{proof}

\begin{remark}
\label{Rem_Cor_Xow_S_density}This corollary assumes $w$ has property W02 or
W03. However part 11 of Remark \ref{Rem_Def_extend_wt_fn} implies that if the
weight function conditions W02 or W03 are \textbf{weakened to W2.1 and W2.2
then} \textbf{our corollary is still valid}.
\end{remark}

\subsection{The density of $\left(  C_{0}^{\infty}\left(  \mathbb{R}%
^{d}\setminus\mathcal{A}\right)  \right)  ^{\vee}$ in $X_{w}^{0}%
$\label{SbSect_density_invF[Coinf(R/A)]_in_X}.}

The aim of this subsection is to prove that $\left(  C_{0}^{\infty}\left(
\mathbb{R}^{d}\setminus\mathcal{A}\right)  \right)  ^{\vee}$ is dense in
$X_{w}^{0}$. The next result is the analogue of Theorem
\ref{Thm_sqrtwCoo_dense_L2} of the previous subsection.

\begin{theorem}
\label{Thm_sqrtwCoo(R/A)_dense_L2}Suppose $w$ is a weight function with
properties W01 with respect to the set $\mathcal{A}$.

Then the set $\sqrt{w}C_{0}^{\left(  0\right)  }\left(  \mathbb{R}%
^{d}\setminus\mathcal{A}\right)  =\left\{  \sqrt{w}f:f\in C_{0}^{\left(
0\right)  }\left(  \mathbb{R}^{d}\setminus\mathcal{A}\right)  \right\}  $ is
dense in $L^{2}$.

Here $C_{0}^{\left(  0\right)  }\left(  \mathbb{R}^{d}\setminus\mathcal{A}%
\right)  $ denotes the $C_{0}^{\left(  0\right)  }$ functions with compact
support in $\mathbb{R}^{d}\setminus\mathcal{A}$.
\end{theorem}

\begin{proof}
By Lemma \ref{Lem_func_eq_1_nbhd_set}, given $h>0$ there exists a function
$\psi_{h}\in C^{\infty}$ satisfying $0\leq\psi_{h}\leq1$, $\psi_{h}=1$ on
$\mathcal{A}_{h}$ and $\psi_{h}=0$ outside $\mathcal{A}_{3h}$. Now define the
mapping $\Psi_{h}:C_{0}^{\left(  0\right)  }\rightarrow C_{0}^{\left(
0\right)  }\left(  \mathbb{R}^{d}\setminus\mathcal{A}\right)  $ by
\[
\Psi_{h}g=\left(  1-\psi_{h}\right)  g,\text{\quad}g\in C_{0}^{\left(
0\right)  },\text{ }h>0.
\]

Since $w$ is positive and continuous on $\mathbb{R}^{d}\setminus\mathcal{A}$,
we have $\frac{1-\psi_{h}}{\sqrt{w}}\in C_{0}^{\left(  0\right)  }\left(
\mathbb{R}^{d}\setminus\mathcal{A}\right)  $ and hence

$\left(  1-\psi_{h}\right)  g\in\sqrt{w}C_{0}^{\left(  0\right)  }\left(
\mathbb{R}^{d}\setminus\mathcal{A}\right)  \subset L^{2}$. In other words,
$\Psi_{h}:C_{0}^{\left(  0\right)  }\rightarrow\sqrt{w}C_{0}^{\left(
0\right)  }\left(  \mathbb{R}^{d}\setminus\mathcal{A}\right)  $.

Now $C_{0}^{\left(  0\right)  }$ is dense in $L^{2}$ and given $g_{0}\in
L^{2}$ and $\varepsilon>0$ we can choose $g_{\varepsilon}\in C_{0}^{\left(
0\right)  }$ such that $\left\Vert g_{0}-g_{\varepsilon}\right\Vert _{L^{2}%
}<\varepsilon/2$. Then
\[
\left\Vert g_{0}-\Psi_{h}g_{\varepsilon}\right\Vert _{L^{2}}\leq\left\Vert
g_{0}-g_{\varepsilon}\right\Vert _{L^{2}}+\left\Vert g_{\varepsilon}-\Psi
_{h}g_{\varepsilon}\right\Vert _{L^{2}}\leq\varepsilon/2+\left\Vert \psi
_{h}g_{\varepsilon}\right\Vert _{L^{2}}\leq\varepsilon/2+\left(
\int\limits_{\mathcal{A}_{3h}}\left\vert g_{\varepsilon}\right\vert
^{2}\right)  ^{\frac{1}{2}}.
\]

Finally, since $\left\vert g_{\varepsilon}\right\vert ^{2}\in L^{1}$, it
follows from Lemma \ref{Lem_Ae_intersect_sphere} that $\lim
\limits_{h\rightarrow0}\int_{\mathcal{A}_{3h}}\left\vert g_{\varepsilon
}\right\vert ^{2}=0$. Thus we can choose $h>0$ so that $\left\Vert g_{0}%
-\Psi_{h}g_{\varepsilon}\right\Vert _{L^{2}}\leq\varepsilon$.
\end{proof}

\begin{lemma}
\label{Lem_Co,infin(R/A)_dense_L2}Suppose $\mathcal{A}$ is a closed set of
measure zero. Then $C_{0}^{\infty}\left(  \mathbb{R}^{d}\setminus
\mathcal{A}\right)  $ is dense in $L^{2}\left(  \mathbb{R}^{d}\right)  $.
\end{lemma}

\begin{proof}
We know that $C_{0}^{\infty}$ is dense in $L^{2}$ so given $\phi\in
C_{0}^{\infty}$ we set $\phi_{\eta}=\left(  1-f_{\eta}\right)  \phi\in
C_{0}^{\infty}$, $\eta>0$ where $f_{\eta}$ is constructed in Lemma
\ref{Lem_func_eq_1_nbhd_set}. Now $f_{\eta}$ has the properties
$\operatorname*{supp}f_{\eta}\subset\mathcal{A}_{\eta}$ and $0\leq f_{\eta
}\leq1$ so that $\phi_{\eta}\in C_{0}^{\infty}\left(  \mathbb{R}^{d}%
\setminus\mathcal{A}\right)  $ and by so by part 2 of Lemma
\ref{Lem_Ae_intersect_sphere}
\[
\left\Vert \phi-\phi_{\eta}\right\Vert _{2}^{2}=\int_{\mathcal{A}_{\eta}%
}\left\vert f_{\eta}\phi\right\vert ^{2}\leq\int_{\mathcal{A}_{\eta}%
}\left\vert \phi\right\vert ^{2}\rightarrow0,
\]
as $\eta\rightarrow0^{+}$, which confirms our density claim.
\end{proof}

The next result is the analogue of Theorem \ref{Thm_sqrtwCoinf_dense_L2}.

\begin{theorem}
\label{Thm_sqrtwCoinf(R/A)_dense_L2}Suppose $w$ has weight function property
W1 with respect to weight set $\mathcal{A}$.

Then the set $\sqrt{w}C_{0}^{\infty}\left(  \mathbb{R}^{d}\setminus
\mathcal{A}\right)  =\left\{  \sqrt{w}f:f\in C_{0}^{\infty}\left(
\mathbb{R}^{d}\setminus\mathcal{A}\right)  \right\}  $ is dense in $L^{2}$.
\end{theorem}

\begin{proof}
For compactness set $\mathcal{A}^{c}=\mathbb{R}^{d}\setminus\mathcal{A}$. Now
$\sqrt{w}C_{0}^{\infty}\left(  \mathcal{A}^{c}\right)  \subset\sqrt{w}%
C_{0}^{\left(  0\right)  }\left(  \mathcal{A}^{c}\right)  $. Hence, by Theorem
\ref{Thm_sqrtwCoo(R/A)_dense_L2}, if it can be shown that $\sqrt{w}%
C_{0}^{\infty}\left(  \mathcal{A}^{c}\right)  $ is dense in $\sqrt{w}%
C_{0}^{\left(  0\right)  }\left(  \mathcal{A}^{c}\right)  $, it follows that
this theorem is true.

First note that from Lemma \ref{Lem_Co,infin(R/A)_dense_L2}, $C_{0}^{\infty
}\left(  \mathcal{A}^{c}\right)  $ is dense in $L^{2}$ and hence
$C_{0}^{\left(  0\right)  }\left(  \mathcal{A}^{c}\right)  $ is dense in
$L^{2}$. Thus for each $\varepsilon>0$ there exists a mapping $\Theta
_{\varepsilon}:C_{0}^{\left(  0\right)  }\left(  \mathcal{A}^{c}\right)
\rightarrow C_{0}^{\infty}\left(  \mathcal{A}^{c}\right)  $ such that
$\left\Vert f-\Theta_{\varepsilon}f\right\Vert _{2}<\varepsilon$.

Next we want to use $\Theta_{\varepsilon}$ to construct a mapping (see left
side of \ref{a9.1}) from $\sqrt{w}C_{0}^{\left(  0\right)  }\left(
\mathcal{A}^{c}\right)  $ to $\sqrt{w}C_{0}^{\infty}\left(  \mathcal{A}%
^{c}\right)  $ which can be used to prove that $\sqrt{w}C_{0}^{\infty}\left(
\mathcal{A}^{c}\right)  $ is dense in $\sqrt{w}C_{0}^{\left(  0\right)
}\left(  \mathcal{A}^{c}\right)  $. If $f\in\sqrt{w}C_{0}^{\left(  0\right)
}\left(  \mathcal{A}^{c}\right)  $, then $\frac{f}{\sqrt{w}}\in C_{0}^{\left(
0\right)  }\left(  \mathcal{A}^{c}\right)  $ and $\sqrt{w}\Theta_{\delta
}\left(  \frac{f}{\sqrt{w}}\right)  \in\sqrt{w}C_{0}^{\infty}\left(
\mathcal{A}^{c}\right)  $, $\delta>0$. But we still need $\sqrt{w}%
\Theta_{\delta}\left(  \frac{f}{\sqrt{w}}\right)  \in L^{2}$. To do this we
use the function $\psi_{h}$ defined in the previous Theorem
\ref{Thm_sqrtwCoo(R/A)_dense_L2}. Since $\left(  1-\psi_{h}\right)  \sqrt
{w}\in C^{\left(  0\right)  }\left(  \mathcal{A}^{c}\right)  $ we have
\[
\left(  1-\psi_{h}\right)  \sqrt{w}\Theta_{\delta}\left(  \frac{f}{\sqrt{w}%
}\right)  \in\sqrt{w}C_{0}^{\infty}\left(  \mathcal{A}^{c}\right)  ,\text{
}when\text{ }f\in\sqrt{w}C_{0}^{\left(  0\right)  }\left(  \mathcal{A}%
^{c}\right)  .
\]

Finally, it will turn out that the support of $\left(  1-\psi_{h}\right)
\sqrt{w}\Theta_{\delta}\left(  \frac{f}{\sqrt{w}}\right)  $ needs to be
restricted by multiplying it by a function $\phi_{R}\in C_{0}^{\infty}$,
$R\geq1$, such that $0\leq\phi_{R}\leq1$, $\operatorname*{supp}\phi_{R}%
\subset\left[  -R-1,R+1\right]  $ and $\phi_{R}=1$ on $\left[  -R,R\right]  $.

We now assert that there exist $R$, $h$ and $\delta$ such that

$\left\Vert \left(  1-\psi_{h}\right)  \phi_{R}\sqrt{w}\Theta_{\delta}\left(
\frac{f}{\sqrt{w}}\right)  -f\right\Vert _{L^{2}}<\varepsilon$. To prove this
assertion write
\begin{align*}
\left(  1-\psi_{h}\right)  \phi_{R}\sqrt{w}\Theta_{\delta}\left(  \frac
{f}{\sqrt{w}}\right)  -f  & =\left(  1-\psi_{h}\right)  \phi_{R}\sqrt{w}%
\Theta_{\delta}\left(  \frac{f}{\sqrt{w}}\right)  -\left(  1-\psi_{h}\right)
\phi_{R}\sqrt{w}\left(  \frac{f}{\sqrt{w}}\right)  +\\
& \qquad+\left(  1-\psi_{h}\right)  \phi_{R}\sqrt{w}\left(  \frac{f}{\sqrt{w}%
}\right)  -\phi_{R}f+\phi_{R}f-f\\
& =\left(  1-\psi_{h}\right)  \phi_{R}\sqrt{w}\left(  \Theta_{\delta}\left(
\frac{f}{\sqrt{w}}\right)  -\frac{f}{\sqrt{w}}\right)  +\psi_{h}\phi
_{R}f+\left(  \phi_{R}-1\right)  f.
\end{align*}

Observe that $\left(  1-\psi_{h}\right)  \sqrt{w}\in C^{\left(  0\right)
}\left(  \mathcal{A}^{c}\right)  $ implies $\left(  1-\psi_{h}\right)
\phi_{R}\sqrt{w}\in C_{0}^{\left(  0\right)  }\left(  \mathcal{A}^{c}\right)
$. Hence
\begin{align*}
\left\Vert \left(  1-\psi_{h}\right)  \phi_{R}\sqrt{w}\Theta_{\delta}\left(
\frac{f}{\sqrt{w}}\right)  -f\right\Vert _{L^{2}} &  \leq\left\Vert \left(
1-\psi_{h}\right)  \phi_{R}\sqrt{w}\left(  \Theta_{\delta}\left(  \frac
{f}{\sqrt{w}}\right)  -\frac{f}{\sqrt{w}}\right)  \right\Vert _{L^{2}}+\\
&  \qquad+\left\Vert \psi_{h}f\right\Vert _{L^{2}}+\left\Vert \left(  \phi
_{R}-1\right)  f\right\Vert _{L^{2}}\\
&  \leq\left\Vert \left(  1-\psi_{h}\right)  \phi_{R}\sqrt{w}\right\Vert
_{\infty}\left\Vert \Theta_{\delta}\left(  \frac{f}{\sqrt{w}}\right)
-\frac{f}{\sqrt{w}}\right\Vert _{L^{2}}+\\
&  \qquad+\left\Vert \psi_{h}f\right\Vert _{L^{2}}+\left\Vert \left(  \phi
_{R}-1\right)  f\right\Vert _{L^{2}}\\
&  \leq\left\Vert \left(  1-\psi_{h}\right)  \phi_{R}\sqrt{w}\right\Vert
_{\infty}\delta+\left\Vert \psi_{h}f\right\Vert _{L^{2}}+\left\Vert \left(
\phi_{R}-1\right)  f\right\Vert _{L^{2}}.
\end{align*}

We will now consider each term on the right side of the last inequality
separately. Regarding the last term,

$\left\Vert \left(  \phi_{R}-1\right)  f\right\Vert _{L^{2}}=\left(
\int_{\left\vert \cdot\right\vert \geq R}\left\vert f\right\vert ^{2}\right)
^{1/2}$ and so $R$ can be fixed so that $\left\Vert \left(  \phi_{R}-1\right)
f\right\Vert _{L^{2}}\leq\varepsilon/3$. Next, by definition of $\psi_{h}$,
$\left\Vert \psi_{h}f\right\Vert _{L^{2}}=\left(  \int_{\mathcal{A}_{3h}%
}\left\vert f\right\vert ^{2}\right)  ^{1/2}$ and, since $\left\vert
f\right\vert ^{2}\in L^{1}$, it follows from Lemma
\ref{Lem_Ae_intersect_sphere} that $\lim\limits_{h\rightarrow0}\int%
_{\mathcal{A}_{3h}}\left\vert f\right\vert ^{2}=0$. A value of $h$ can now be
chosen so that $\left\Vert \psi_{h}f\right\Vert _{L^{2}}\leq\varepsilon/3$.
Now we have
\[
\left\Vert \left(  1-\psi_{h}\right)  \phi_{R}\sqrt{w}\Theta_{\delta}\left(
\frac{f}{\sqrt{w}}\right)  -f\right\Vert _{L^{2}}\leq\left\Vert \left(
1-\psi_{h}\right)  \phi_{R}\sqrt{w}\right\Vert _{\infty}\delta+2\varepsilon/3,
\]

for all $f\in\sqrt{w}C_{0}^{\left(  0\right)  }\left(  \mathcal{A}^{c}\right)
. $

The last step is to choose $\delta$ so that $\left\Vert \left(  1-\psi
_{h}\right)  \phi_{R}\sqrt{w}\right\Vert _{\infty}\delta\leq\varepsilon/3 $,
and the theorem follows.
\end{proof}

The next result is the analogue of Corollary \ref{Cor_Xow_S_density}.

\begin{corollary}
\label{Cor_Xow_invFCoinf(R/A)_density}Suppose the weight function $w$ has
property W01 w.r.t. the set $\mathcal{A}$ as well as having property W02 or W03.

Then $\left(  C_{0}^{\infty}\left(  \mathbb{R}^{d}\setminus\mathcal{A}\right)
\right)  ^{\vee}$ is dense in $X_{w}^{0}$.
\end{corollary}

\begin{proof}
Since the Fourier transform is a homeomorphism form $L^{2}$ to $L^{2}$ it
follows from Theorem \ref{Thm_sqrtwCoinf_dense_L2} that $\left(  \sqrt{w}%
C_{0}^{\infty}\left(  \mathbb{R}^{d}\setminus\mathcal{A}\right)  \right)
^{\vee}$ is dense in $L^{2}$.

Now \textbf{since }$w$\textbf{\ has property W02 or W03} by Theorem
\ref{Thm_I_J_property} the operator $\mathcal{J}:L^{2}\rightarrow X_{w}^{0}$
of Definition \ref{Def_I_J} is a homeomorphism with inverse denoted by
$\mathcal{I}$ and we now show that
\[
\mathcal{J}:\left(  \sqrt{w}C_{0}^{\infty}\left(  \mathbb{R}^{d}%
\setminus\mathcal{A}\right)  \right)  ^{\vee}\rightarrow\left(  C_{0}^{\infty
}\left(  \mathbb{R}^{d}\setminus\mathcal{A}\right)  \right)  ^{\vee},
\]

is onto which means that $\left(  C_{0}^{\infty}\left(  \mathbb{R}%
^{d}\setminus\mathcal{A}\right)  \right)  ^{\vee}$ is dense in $X_{w}^{0}$.

Now $\mathcal{J}u=\left(  \widehat{u}/\sqrt{w}\right)  ^{\vee}\in\left(
C_{0}^{\infty}\left(  \mathbb{R}^{d}\setminus\mathcal{A}\right)  \right)
^{\vee}$. Further, if $v\in\left(  C_{0}^{\infty}\left(  \mathbb{R}%
^{d}\setminus\mathcal{A}\right)  \right)  ^{\vee}$ then
\[
\mathcal{J}^{-1}v=\mathcal{I}v=\left(  \sqrt{w}\widehat{u}\right)  ^{\vee}%
\in\left(  \sqrt{w}C_{0}^{\infty}\left(  \mathbb{R}^{d}\setminus
\mathcal{A}\right)  \right)  ^{\vee},
\]

and so \ref{a9.0} is onto.
\end{proof}

\begin{remark}
\label{Rem_Cor_Xow_invFCoinf(R/A)_density}We can replace W02 or W03 by the
weaker properties W2.1 and W2.2 discussed in part 11 of Remark
\ref{Rem_Def_extend_wt_fn}.
\end{remark}

\subsection{The Hilbert spaces $X_{w:K}^{0}$}

\begin{definition}
\label{Def_XowK}If $K\subset\mathbb{R}^{d}$ is a closed set then%
\[
X_{w:K}^{0}:=\left(  X_{w}^{0}\right)  _{K}:=\left\{  f\in X_{w}%
^{0}:\operatorname*{supp}f\subseteq K\right\}  ,
\]

which we endow with the subspace norm.

This notation is analogous to the Sobolev space notation of Petersen
\cite{Petersen83} following Exercise 3:4.10 on p243.
\end{definition}

\begin{theorem}
\label{Thm_XowK_Hilbert}$X_{w:K}^{0}$ is a closed subspace of $X_{w}^{0}$ and
hence a Hilbert space.
\end{theorem}

\begin{proof}
If $X_{w:K}^{0}$ is empty then it is automatically closed. Suppose it is not empty.

Suppose $\left\{  f_{k}\right\}  $ is Cauchy series in $\left(  X_{w}%
^{0}\right)  _{K}$. Then, because this space has been endowed with the
subspace norm, $\left\{  f_{k}\right\}  $ is a Cauchy series in $X_{w}^{0}$
and thus converges to some $f\in X_{w}^{0}$. Suppose $f\left(  x\right)  \neq0
$ for some $x\notin K$. But
\begin{align*}
\left\vert f\left(  x\right)  \right\vert \leq\left\vert f\left(  x\right)
-f_{k}\left(  x\right)  \right\vert +\left\vert f_{k}\left(  x\right)
\right\vert =\left\vert f\left(  x\right)  -f_{k}\left(  x\right)
\right\vert  & =\left\vert \left(  f-f_{k},R_{x}\right)  _{w,0}\right\vert \\
& \leq\left\Vert f-f_{k}\right\Vert _{w,0}\left\Vert R_{x}\right\Vert _{w,0}\\
& \rightarrow0,
\end{align*}

as $k\rightarrow\infty$ which contradicts $f\left(  x\right)  \neq0$. Thus
$\operatorname*{supp}f\subseteq K$ and $\left(  X_{w}^{0}\right)  _{K}$ is
closed and hence a Hilbert space.
\end{proof}

\section{Basis functions\label{Sect_basis_distrib}}

In this section we will define a (unique) basis function of order zero
generated by a weight function with property W02/W03, and note that this
contrasts with the positive order case where the basis function is only unique
modulo a space of polynomials. Basis functions are generally denoted by $G$
and will always be members of $X_{w}^{0}$. If a weight function has property
W02/W03 for parameter $\kappa$ it will be shown that $G\in C_{B}^{\left(
\left\lfloor 2\kappa\right\rfloor \right)  }$, whereas in general we know that
$X_{w}^{0}\subset C_{B}^{\left(  \left\lfloor \kappa\right\rfloor \right)  }$.
Basis functions are used to construct solutions to the variational
interpolation problem of Chapter \ref{Ch_Interpol} of this document and to the
two smoothing problems studied subsequently in Chapters \ref{Ch_Exact_smth}
and \ref{Ch_Approx_smth}.

We will also derive results concerning the basis functions generated by the
tensor (or direct) product of weight functions with reference to the example
of the extended B-spline basis functions introduced in Subsection
\ref{SbSect_wt_fn_examples_2}.

\begin{definition}
\label{Def_basis_func}\textbf{Basis function}

Suppose a weight function $w$ has property W02 or W03. Then $\frac{1}{w}\in
L^{1}$ and hence by Lemma \ref{Lem_L1_Fourier_contin}, $\left(  \frac{1}%
{w}\right)  ^{\vee}\in C_{B}^{\left(  0\right)  }$. We now define the unique
basis function $G$ of order $0$ generated by $w$ to be
\begin{equation}
G=\left(  \frac{1}{w}\right)  ^{\vee}.\label{1.33}%
\end{equation}

\end{definition}

A simple consequence of this definition is that
\begin{equation}
\overline{G\left(  x\right)  }=G\left(  -x\right)  ,\label{1.301}%
\end{equation}

from which follow the results:

\begin{theorem}
\label{Thm_basis_property_minus}Suppose $G$ is a basis function of order $0$.
Then as distributions and for all $\alpha$:

\begin{enumerate}
\item $G\left(  0\right)  $ is real.

\item $\left(  D^{\alpha}G\right)  \left(  -x\right)  =\left(  -1\right)
^{\left\vert \alpha\right\vert }\overline{D^{\alpha}G}\left(  x\right)  $;

\item $\left(  D^{\alpha}\operatorname{Re}G\right)  \left(  -x\right)
=\left(  -1\right)  ^{\left\vert \alpha\right\vert }D^{\alpha}%
\operatorname{Re}G\left(  x\right)  $ and $\left(  D^{\alpha}\operatorname{Im}%
G\right)  \left(  -x\right)  =\left(  -1\right)  ^{\left\vert \alpha
\right\vert +1}D^{\alpha}\operatorname{Im}G\left(  x\right)  $;

\item When $\left\vert \alpha\right\vert $ is even $D^{\alpha}%
\operatorname{Re}G$ is even and when $\left\vert \alpha\right\vert $ is odd
$D^{\alpha}\operatorname{Re}G$ is odd;

\item When $\left\vert \alpha\right\vert $ is odd $D^{\alpha}\operatorname{Im}%
G$ is even and when $\left\vert \alpha\right\vert $ is odd $D^{\alpha
}\operatorname{Im}G$ is even;

\item $\left(  \left(  aD\right)  ^{k}G\right)  \left(  -x\right)  =\left(
-1\right)  ^{k}\left(  \left(  aD\right)  ^{k}\overline{G}\right)  \left(
x\right)  =\left(  \left(  -aD\right)  ^{k}\overline{G}\right)  \left(
x\right)  ,\quad k\geq0;$ $a\in\mathbb{R}^{d}$;
\end{enumerate}
\end{theorem}

The next theorem makes two important points:

(1) unlike in the positive order case basis functions are unique and are data
functions i.e. they are members of $X_{w}^{0}$;

(2) scaling can be used to enlarge the classes of weight functions available.

\begin{theorem}
\label{Thm_G_in_Xow}Suppose $w$ has property W02 or W03 for $\kappa$. Then:

\begin{enumerate}
\item The basis function $G$ generated by the weight function $w$ is a member
of $X_{w}^{0}$.

\item \textbf{Vector scaling/dilation} If $\lambda\in\mathbb{R}^{d}$ and
$\lambda.>0$ then the weight function $\lambda^{\mathbf{1}}w\left(
t./\lambda\right)  $ has property W02/W03 for $\kappa$ iff $w$ has property
W02/W03 for $\kappa$, and $G\left(  \lambda.x\right)  $ is the basis function
generated by $\lambda^{\mathbf{1}}w\left(  t./\lambda\right)  $.

\item \textbf{Scalar scaling/dilation} If $\lambda\in\mathbb{R}^{1}$ and
$\lambda>0$ then the weight function $\lambda^{d}w\left(  t/\lambda\right)  $
has property W02/W03 for $\kappa$ iff $w$ has property W02/W03 for $\kappa$,
and $G\left(  \lambda x\right)  $ is the basis function generated by
$\lambda^{d}w\left(  t/\lambda\right)  $.
\end{enumerate}
\end{theorem}

\begin{proof}
\textbf{Parts 1 and 2} Since $C_{B}^{\left(  0\right)  }\subset S^{\prime}$ it
follows that $G\in S^{\prime}$, $\widehat{G}\in L_{loc}^{1}$ and $\sqrt
{w}\widehat{G}=\frac{1}{\sqrt{w}}\in L^{2}$ i.e. $G\in X_{w}^{0}$.

Further, $\widehat{G\left(  \lambda.x\right)  }=\frac{1}{\lambda^{\mathbf{1}}%
}\widehat{G}\left(  t./\lambda\right)  =\frac{1}{\lambda^{\mathbf{1}}w\left(
t./\lambda\right)  }$ so the corresponding weight function is $\lambda
^{\mathbf{1}}w\left(  \frac{t}{\lambda}\right)  $. But by part 8 of Remark
\ref{Rem_Def_extend_wt_fn} this is a weight function with the same parameter
$\kappa$ as $w$.

\textbf{Part 3} Use $\lambda\mathbf{1}$ in part 2 so that $\left(
\lambda\mathbf{1}\right)  ^{\mathbf{1}}=\lambda^{d}$.
\end{proof}

\subsection{The smoothness of basis functions: $L^{1}$
theory\label{SbSect_smooth_basis}}

Not only is the basis function a data function but all basis function
derivatives up to order $\left\lfloor \kappa\right\rfloor $ are data functions:

\begin{theorem}
\label{Thm_basis_fn_properties_all_m_W2}Suppose the weight function $w$ has
property \textbf{W02} for scalar parameter $\kappa$. Then the basis function
$G$ generated by $w$ is in $C_{B}^{\left(  \max\left\lfloor 2\kappa
\right\rfloor \right)  }\left(  \mathbb{R}^{d}\right)  $. Additionally, the
inverse Fourier transform formulas
\begin{equation}
D^{\beta}G\left(  x\right)  =\tfrac{1}{\left(  2\pi\right)  ^{d/2}}\int
e^{ix\xi}\widehat{D^{\beta}G}\left(  \xi\right)  d\xi,\text{\quad}\left\vert
\beta\right\vert \leq2\kappa,\label{1.024}%
\end{equation}

hold and $D^{\beta}G\in X_{w}^{0}$ when $\left\vert \beta\right\vert
\leq\kappa$.
\end{theorem}

\begin{proof}
Since property W02 implies $\frac{1}{w}\in L^{1}$ and $\frac{\left\vert
\cdot\right\vert ^{2\kappa}}{w}\in L^{1}$, an application of the formula of
part 2 of Lemma \ref{Lem_L1_Fourier_contin} with $f=\frac{1}{w}$ proves $G\in
C_{B}^{\left(  \left\lfloor 2\kappa\right\rfloor \right)  }$. Formula
\ref{1.024} then follows from equation \ref{1.031} of Lemma
\ref{Lem_L1_Fourier_contin}.

Now to prove $D^{\beta}G\in X_{w}^{0}$. Firstly, $G\in S^{\prime}$ implies
$D^{\beta}G\in S^{\prime}$. Next, for compact $K$, $\int_{K}\left\vert
\widehat{D^{\beta}G}\right\vert \leq\int_{K}\frac{\left\vert \cdot\right\vert
^{\left\vert \beta\right\vert }}{w}<\infty$ since $\frac{1}{w}\in L^{1}$.
Finally, for $\left\vert \beta\right\vert \leq\kappa$, \ref{1.50} implies
\begin{equation}
\left\Vert D^{\beta}G\right\Vert _{w,0}^{2}=\int w\left\vert \widehat{D^{\beta
}G}\right\vert ^{2}=\int\frac{\xi^{2\beta}}{w}\leq\int\frac{\left\vert
\cdot\right\vert ^{2\left\vert \beta\right\vert }}{w}<\infty,\label{1.31}%
\end{equation}

since $w$ has property W02.
\end{proof}

\begin{remark}
\label{Rem_basis_fn_properties_all_m_W2}\ 

\begin{enumerate}
\item If the set of values of $2\kappa$ has a finite upper bound then although
$\max2\kappa$ may not exist e.g. the extended B-splines example below,
$\max\left\lfloor 2\kappa\right\rfloor $ always exists because $\left\lfloor
2\kappa\right\rfloor $ only takes integer values.

\item Allowing $\kappa$ to take non-integer values in the definition of weight
property W02 (Definition \ref{Def_extend_wt_fn}) can yield an extra order of
smoothness e.g. the Sobolev spline example below.
\end{enumerate}
\end{remark}

\begin{theorem}
\label{Thm_basis_fn_properties_all_m_W3}Suppose the weight function $w$ has
property \textbf{W03} for parameter $\kappa$. Then the basis function $G$
generated by $w$ is in $C_{B}^{\left(  \max\left\lfloor 2\kappa\right\rfloor
\right)  }\left(  \mathbb{R}^{d}\right)  $. Additionally, the inverse Fourier
transform formulas
\begin{equation}
D^{\beta}G\left(  x\right)  =\tfrac{1}{\left(  2\pi\right)  ^{d/2}}\int
e^{ix\xi}\widehat{D^{\beta}G}\left(  \xi\right)  d\xi,\text{\quad}\beta
\leq2\kappa,\label{1.067}%
\end{equation}

hold and $D^{\beta}G\in X_{w}^{0}$ when $\beta\leq\kappa$.
\end{theorem}

\begin{proof}
Since property W03 implies $\frac{1}{w}\in L^{1}$ and $\frac{\xi^{2\kappa}}%
{w}\in L^{1}$, an application of the formula of part 2 of Lemma
\ref{Lem_L1_Fourier_contin} with $f=\frac{1}{w}$ proves $G\in C_{B}^{\left(
\left\lfloor 2\kappa\right\rfloor \right)  }$ and so $G\in C_{B}^{\left(
\max\left\lfloor 2\kappa\right\rfloor \right)  }$. Formula \ref{1.067} then
follows from equation \ref{1.031} of Lemma \ref{Lem_L1_Fourier_contin}.

Now to prove $D^{\beta}G\in X_{w}^{0}$. Firstly, $G\in S^{\prime}$ implies
$D^{\beta}G\in S^{\prime}$. Next, for compact $K$, $\int_{K}\left\vert
\widehat{D^{\beta}G}\right\vert \leq\int_{K}\frac{\left\vert \xi^{\beta
}\right\vert }{w}\leq\int_{K}\frac{\left\vert \xi\right\vert ^{\left\vert
\beta\right\vert }}{w}<\infty$ since $\frac{1}{w}\in L^{1}$. Finally, for
$\beta\leq\kappa$, \ref{2.0} implies
\begin{equation}
\left\Vert D^{\beta}G\right\Vert _{w,0}^{2}=\int w\left\vert \widehat{D^{\beta
}G}\right\vert ^{2}=\int\frac{\xi^{2\beta}}{w}<\infty,\label{1.068}%
\end{equation}

since $w$ has property W03.
\end{proof}

\begin{corollary}
\label{Cor_basis_smth_W2_W3}Suppose the weight function $w$ has property
\textbf{W02 or} \textbf{W03} for parameter $\kappa$. Then the basis function
satisfies $G\in C_{B}^{\left(  \max\left\lfloor \underline{2\kappa
}\right\rfloor \right)  }$.
\end{corollary}

\begin{proof}
From Theorem \ref{Thm_equiv_W3}, if $w$ has property W03 for (vector)
parameter $\kappa$ then $w$ has property W02 for $\underline{\kappa}$ and
hence Theorem \ref{Thm_basis_fn_properties_all_m_W2} implies $G\in
C_{B}^{\left(  \left\lfloor 2\kappa\right\rfloor \right)  }=C_{B}^{\left(
\left\lfloor \underline{2\kappa}\right\rfloor \right)  }$. If $w$ has property
W02 for parameter $\kappa$ then $C_{B}^{\left(  \left\lfloor 2\kappa
\right\rfloor \right)  }=C_{B}^{\left(  \left\lfloor \underline{2\kappa
}\right\rfloor \right)  }$ since $\kappa$ is now a scalar.
\end{proof}

\subsection{The smoothness of basis functions: $L^{2}$ theory}

\begin{theorem}
\label{Thm_basis_DG_in_L2}Suppose the weight function $w$ has property
\textbf{W02 or} \textbf{W03}. Suppose the basis function generated by $w$
satisfies $D^{\beta}G\in L^{2}$.

Then \ref{1.067} also holds in the $L^{2}$ sense.
\end{theorem}

\subsection{Radial basis functions\label{SbSect_radial_basis_func}}

\begin{theorem}
\label{Thm_basis_radial}If a weight function $w$ is radial with property W02
for parameter $\kappa$. Then the basis function, say $G$, is also radial. In
fact%
\[
G\left(  x\right)  =\left(  2\pi\right)  ^{-\frac{d}{2}}\int_{\mathbb{R}^{1}%
}\frac{e^{i\left\vert x\right\vert s}}{\overset{\circ}{w}\left(  s\right)
}ds,
\]

where $\frac{1}{\overset{\circ}{w}\left(  s\right)  }=\int_{\mathbb{R}^{d-1}%
}\frac{d\xi^{\prime\prime}}{w\left(  s,\xi^{\prime\prime}\right)  }$ and
$\frac{1}{\overset{\circ}{w}}\in L^{1}\left(  \mathbb{R}^{1}\right)  $.

Further, if $\overset{\circ}{w}$ is continuous and positive outside a closed
set of measure zero then $\overset{\circ}{w}$ is a weight function with
property W02 for $\kappa$.
\end{theorem}

\begin{proof}
By definition, $G\left(  x\right)  =\left(  2\pi\right)  ^{-\frac{d}{2}}%
\int\frac{e^{-ix\xi}}{w\left(  \xi\right)  }d\xi$ and so Lemma
\ref{Lem_wt_func_radial} with $u\left(  x\xi\right)  =e^{-ix\xi}$ yields

$G\left(  x\right)  =\left(  2\pi\right)  ^{-\frac{d}{2}}\int_{\mathbb{R}^{1}%
}\frac{e^{-i\left\vert x\right\vert s}}{\overset{\circ}{w}\left(  s\right)
}ds $ and $\int_{\mathbb{R}^{1}}\frac{1}{\overset{\circ}{w}}=\int\frac{1}%
{w}<\infty$. Finally%
\[
\int\frac{s^{2\kappa}ds}{\overset{\circ}{w}\left(  s\right)  }=\int\int%
\frac{s^{2\kappa}dsd\xi^{\prime\prime}}{w\left(  s,\xi^{\prime\prime}\right)
}\leq\int\frac{\left\vert \cdot\right\vert ^{2\kappa}}{w}<\infty,
\]

so that $\overset{\circ}{w}$ is a weight function with property W02 for
$\kappa$.
\end{proof}

The radial functions used in Examples \ref{SbSect_wt_func_examples} to
generate weight functions (generally denoted by $G$), namely the shifted
thin-plate splines, Gaussian and Sobolev splines (univariate Laplacian
kernels) can be used as radial basis functions. The basis functions associated
with the extended B-spline weight functions will be derived below in
Subsection \ref{SbSect_basis_fns}.

\begin{example}
\textbf{The Shifted thin-plate splines} The basis functions are given by
\ref{1.043} where $-d/2<v<0$ i.e.%
\begin{equation}
G\left(  x\right)  =\left(  -1\right)  ^{\left\lceil v\right\rceil }\left(
1+\left\vert x\right\vert ^{2}\right)  ^{v},\quad-d/2<v<0.\label{1.046}%
\end{equation}

The weight function has property W02 for all $\kappa\geq0$. Hence by Theorem
\ref{Thm_basis_fn_properties_all_m_W2}, the basis function lies in
$C_{B}^{\left(  \left\lfloor 2\kappa\right\rfloor \right)  }$ for $\kappa
\geq0$ and so is a $C_{B}^{\infty}$ function.
\end{example}

\begin{example}
\textbf{The Gaussian} The weight function \ref{1.029} has property W02 for
$\kappa\geq0$ and the basis function is given by \ref{1.044} i.e. $G\left(
x\right)  =e^{-\left\vert x\right\vert ^{2}}$. Hence by Theorem
\ref{Thm_basis_fn_properties_all_m_W2}, we have $G\in C_{B}^{\left(
\left\lfloor 2\kappa\right\rfloor \right)  }$ for $\kappa\geq0$ and so $G\in
C_{B}^{\infty}$. In fact $G\in S$.
\end{example}

\begin{example}
\label{Ex_BasisFuncExSobSpline}\textbf{Sobolev splines} The basis function is
\ref{1.045} i.e.%
\[
G\left(  x\right)  =\frac{1}{2^{v-1}\Gamma\left(  v\right)  }\widetilde{K}%
_{v-d/2}\left(  \left\vert x\right\vert \right)  ,\text{\quad}x\in
\mathbb{R}^{d},
\]

where $v>d/2$. It was shown in Examples \ref{SbSect_wt_func_examples} that
property W02 holds for all $0\leq\kappa<v-d/2$. Now by Remark
\ref{Rem_basis_fn_properties_all_m_W2}, $G\in C_{B}^{\left(  \overline
{\left\lfloor 2\kappa\right\rfloor }\right)  }$ and we see that $\max
\left\lfloor 2\kappa\right\rfloor =\left\lfloor 2v-d\right\rfloor $ when $2v $
is not an integer and $\overline{\left\lfloor 2\kappa\right\rfloor
}=\left\lfloor 2v-d\right\rfloor -1$ when $2v$ is an integer. This is an
example where allowing $\kappa$ to take non-integer values in the definition
of weight property W02 (Definition \ref{Def_extend_wt_fn}) yields an
\textbf{extra order of smoothness}. An example is when $v-d/2=1%
\frac34
$.
\end{example}

\begin{example}
The weight functions of Examples \ref{Ex_wt_func_L1loc} and
\ref{Ex_wt_func_|x|^m_L1loc}. Suggestion: consider case where $v=1$. Use the
radial function integration formulas from the appendix of the positive order
basis function document. Then, if necessary, use formulas for $\int x^{s}%
J_{v}\left(  x\right)  dx$ and $\int_{0}^{1}x^{s}J_{v}\left(  x\right)  dx$
from the reference text [??] e.g. equations 13 and 14 in 6.552 where
$S_{\mu,v}$ denotes the Lommel functions of 8.57.
\end{example}

\subsection{Basis functions generated by tensor product weight functions}

Here we prove that the basis function of a tensor product weight function is
the tensor product of basis functions.

\begin{theorem}
\label{Thm_basis_fn_of_tensor_prod}Suppose $w_{1}$ and $w_{2}$ are weight
functions which satisfy property W02 or W03 for parameter $\kappa$. Then part
2 of Theorem \ref{Thm_ten_prod_two_wt_fns} implies that $w=w_{1}\otimes w_{2}$
is a tensor product weight function satisfying property W02 or W03 for
parameter $\kappa$.

We prove here that $G=G_{1}\otimes G_{2}$ where $G_{1},G_{2}$ and $G$ are the
basis functions of order $0$ generated by $w_{1},w_{2}$ and $w$ respectively.
\end{theorem}

\begin{proof}
Here each $G_{i}\in C_{B}^{\left(  0\right)  }$ and so each $G_{i}\in
S^{\prime}$ and we are interested in the tensor product of members of
$S^{\prime}$. In this proof we refer to results from Vladimirov
\cite{Vladimirov} where the term \textit{direct products} is used instead of
\textit{tensor products} and $S^{\prime}$ is referred to as the space of
generalized functions of slow growth. In Subsection 2.8.5 of Vladimirov it is
shown that $G_{1}\otimes G_{2}\in S^{\prime}$ and in part (e) of Subsection
2.9.3 it is shown that $\widehat{G}=\widehat{G_{1}}\otimes\widehat{G_{2}}$.
Thus
\[
\widehat{G}=\widehat{G_{1}}\otimes\widehat{G_{2}}=\frac{1}{w_{1}}\otimes
\frac{1}{w_{2}}=\frac{1}{w_{1}\otimes w_{2}}=\frac{1}{w}.
\]

\end{proof}

\begin{corollary}
\label{Cor_basis_fn_of_tensor_prod}Suppose that $\left\{  w_{i}\right\}  $ is
a set of weight functions which satisfy W02 or W03 for parameter $\kappa$.
Then by part 2 of Theorem \ref{Thm_ten_prod_two_wt_fns} $w=\otimes w_{i}$ is a
tensor product weight function satisfying property W02/W03 for parameter
$\kappa$.

Further suppose that $\left\{  G_{i}\right\}  $ is the set of basis functions
of order $0$ generated by the $w_{i}$.

We claim here that $G=\otimes G_{i}$ is a basis function of order zero
generated by $w$.
\end{corollary}

\begin{proof}
By induction using Theorem \ref{Thm_basis_fn_of_tensor_prod}.
\end{proof}

One of the main motivations of this document is to define the weight functions
so that tensor product hat functions are basis functions. The following
corollary shows that the hat function is a basis function in any dimension.

\begin{corollary}
\label{Cor_hat_basis_fn}Suppose $\Lambda$ is the $d$-dimensional tensor
product hat function. Then $\Lambda$ is the basis function of order zero
generated by the weight function $\frac{1}{\widehat{\Lambda}}$.
\end{corollary}

\subsection{The extended B-spline basis functions\label{SbSect_basis_fns}}

In this subsection we will derive the basis functions of zero order generated
by the class of tensor product weight functions studied in Subsection
\ref{SbSect_wt_fn_examples_2}, namely the extended B-spline weight functions.
The extended B-spline weight functions were so named because their basis
functions will turn out to be the derivatives of the convolutions of hat
functions, denoted $\left(  \ast\Lambda\right)  ^{l},$\quad$l=1,2,3,\ldots$
and these are the B-splines. The next theorem gives the unique tensor product
basis function of zero order generated by these weight functions assuming the
have property W03.

\begin{theorem}
\label{Thm_basis_tensor_hat_W3}\textbf{The extended B-spline basis functions.}
Suppose $w_{s}$ is the extended B-spline weight function considered in Theorem
\ref{Thm_ex_splin_wt_fn_properties} and suppose $w_{s}$ has property
\textbf{W03} for $\kappa$ i.e. $\kappa+1/2<n\mathbf{1}\leq l\mathbf{1}$.

Then the basis function $G_{s}$ of order zero generated by $w$ is the tensor
product $G_{s}\left(  x\right)  =%
{\textstyle\prod\limits_{k=1}^{d}}
G_{1}\left(  x_{k}\right)  $ where
\begin{equation}
G_{1}\left(  t\right)  =\left(  -1\right)  ^{l-n}\tfrac{\left(  2\pi\right)
^{l/2}}{2^{2\left(  l-n\right)  +1}}D^{2\left(  l-n\right)  }\left(  \left(
\ast\Lambda\right)  ^{l}\right)  \left(  \tfrac{t}{2}\right)  ,\text{\quad
}t\in\mathbb{R}^{1},\label{1.49}%
\end{equation}

and $\left(  \ast\Lambda\right)  ^{l}$ denotes the convolution of $l$
1-dimensional hat functions. Further, if $n<l$ we have
\[
D^{2\left(  l-n\right)  }\left(  \ast\Lambda\right)  ^{l}=\frac{\left(
-1\right)  ^{l-n}}{\left(  2\pi\right)  ^{\left(  l-n\right)  /2}}\left(
\ast\Lambda\right)  ^{n-1}\ast\sum_{k=-\left(  l-n\right)  }^{l-n}\left(
-1\right)  ^{k}\tbinom{2\left(  l-n\right)  }{l-n+k}\Lambda\left(
\cdot-k\right)  ,
\]

$G_{1}\in C_{0}^{2n-2}\left(  \mathbb{R}^{1}\right)  $ with $\left\Vert
D^{2n-2}G_{1}\right\Vert _{\infty}=\tfrac{\sqrt{2\pi}}{2^{2l+1}}\tbinom
{2l-2}{l-1}$.

Further, $D^{2n-1}G_{1}$ is a bounded step function with bounded support and
$\left\Vert D^{2n-1}G_{1}\right\Vert _{\infty}=\tfrac{\sqrt{2\pi}}{2^{2l+2}%
}\tbinom{2l-1}{l-1}$ and $D^{2n}G_{1}$ is the finite sum of translated delta functions.

Finally, $G_{s}\in C_{0}^{\left(  \left(  2n-2\right)  \mathbf{1}\right)
}\left(  \mathbb{R}^{d}\right)  $, the functions $\left\{  D^{\alpha}%
G_{s}:\left(  2n-2\right)  \mathbf{1}<\alpha\leq\left(  2n-1\right)
\mathbf{1}\right\}  $ are bounded with bounded support, and $D^{2n\mathbf{1}%
}G_{s}$ is the finite sum of translated delta functions.
\end{theorem}

\begin{proof}
In this proof it will be better to use the operator notation $F\left[
f\right]  $ for the Fourier transform instead of $\widehat{f}$.

Define $w_{1}\left(  t\right)  =\frac{t^{2n}}{\sin^{2l}t}$,\quad
$t\in\mathbb{R}$. Since the 1-dimensional hat function $\Lambda$ satisfies
\begin{equation}
F\left[  \Lambda\right]  \left(  t\right)  =\left(  2\pi\right)
^{-1/2}\left(  \frac{\sin\left(  t/2\right)  }{t/2}\right)  ^{2},\text{\qquad
}t\in\mathbb{R},\label{1.443}%
\end{equation}

we have
\begin{align*}
\tfrac{1}{w_{1}\left(  t\right)  }=\tfrac{\sin^{2l}t}{t^{2n}}=t^{2\left(
l-n\right)  }\left(  \tfrac{\sin t}{t}\right)  ^{2l} &  =\left(  2\pi\right)
^{l/2}t^{2\left(  l-n\right)  }\left(  F\left[  \Lambda\right]  \left(
2t\right)  \right)  ^{l}\\
&  =\left(  2\pi\right)  ^{l/2}t^{2\left(  l-n\right)  }F\left[  \left(
\ast\Lambda\right)  ^{l}\right]  \left(  2t\right) \\
&  =\left(  2\pi\right)  ^{l/2}\left(  \tfrac{1}{2}\right)  ^{2\left(
l-n\right)  }\left(  t^{2\left(  l-n\right)  }F\left[  \left(  \ast
\Lambda\right)  ^{l}\right]  \right)  \left(  2t\right) \\
&  =\left(  2\pi\right)  ^{l/2}\left(  \tfrac{i}{2}\right)  ^{2\left(
l-n\right)  }F\left[  D^{2\left(  l-n\right)  }\left(  \ast\Lambda\right)
^{l}\right]  \left(  2t\right) \\
&  =\left(  -1\right)  ^{l-n}\tfrac{\left(  2\pi\right)  ^{l/2}}{2^{2\left(
l-n\right)  }}F\left[  D^{2\left(  l-n\right)  }\left(  \ast\Lambda\right)
^{l}\right]  \left(  2t\right)  ,
\end{align*}

Now since the parameters $n$ and $l$ used to define the weight function $w$
are independent of the dimension $d$, $w_{1}$ satisfies properties W01 and W03
and so $1/w_{1}\in L^{1}$. Also
\begin{align*}
G_{1}\left(  t\right)  =F^{-1}\left[  \frac{1}{w_{1}}\right]  \left(
t\right)   &  =\left(  -1\right)  ^{l-n}\tfrac{\left(  2\pi\right)  ^{l/2}%
}{2^{2\left(  l-n\right)  }}F^{-1}\left[  F\left[  D^{2\left(  l-n\right)
}\left(  \ast\Lambda\right)  ^{l}\right]  \left(  2t\right)  \right]  \left(
t\right) \\
&  =\left(  -1\right)  ^{l-n}\tfrac{\left(  2\pi\right)  ^{l/2}}{2^{2\left(
l-n\right)  }}\frac{1}{2}D^{2\left(  l-n\right)  }\left(  \ast\Lambda\right)
^{l}\left(  \frac{t}{2}\right) \\
&  =\left(  -1\right)  ^{l-n}\tfrac{\left(  2\pi\right)  ^{l/2}}{2^{2\left(
l-n\right)  +1}}D^{2\left(  l-n\right)  }\left(  \ast\Lambda\right)
^{l}\left(  \frac{t}{2}\right)  .
\end{align*}

We will require the convolution identities%
\[
\delta\left(  \cdot-a\right)  \ast f=\left(  2\pi\right)  ^{-1/2}f\left(
\cdot-a\right)  ,\quad f\in\mathcal{D}^{\prime},
\]

and%
\[
\left(  \ast\left(  \delta\left(  \cdot-1\right)  -2\delta+\delta\left(
\cdot+1\right)  \right)  \right)  ^{m}=\frac{\left(  -1\right)  ^{m}}{\left(
2\pi\right)  ^{\left(  m-1\right)  /2}}\sum_{k=-m}^{m}\left(  -1\right)
^{k}\tbinom{2m}{m+k}\delta\left(  \cdot-k\right)  .
\]

Now suppose that $n<l$. Then in one dimension
\begin{align*}
D^{2\left(  l-n\right)  }\left(  \ast\Lambda\right)  ^{l}=\left(  \ast
\Lambda\right)  ^{n}\ast\left(  \ast D^{2}\Lambda\right)  ^{l-n} &  =\left(
\ast\Lambda\right)  ^{n}\ast\left(  \ast\left(  \delta\left(  \cdot-1\right)
-2\delta+\delta\left(  \cdot+1\right)  \right)  \right)  ^{l-n}\\
&  =\frac{\left(  -1\right)  ^{l-n}}{\left(  2\pi\right)  ^{\left(
l-n-1\right)  /2}}\left(  \ast\Lambda\right)  ^{n}\ast\sum_{k=-\left(
l-n\right)  }^{l-n}\left(  -1\right)  ^{k}\tbinom{2\left(  l-n\right)
}{l-n+k}\delta\left(  \cdot-k\right) \\
&  =\frac{\left(  -1\right)  ^{l-n}}{\left(  2\pi\right)  ^{\left(
l-n\right)  /2}}\left(  \ast\Lambda\right)  ^{n-1}\ast\sum_{k=-\left(
l-n\right)  }^{l-n}\left(  -1\right)  ^{k}\tbinom{2\left(  l-n\right)
}{l-n+k}\Lambda\left(  \cdot-k\right)  .
\end{align*}

Further
\begin{align}
D^{2n-2}G_{1}\left(  t\right)   & =\left(  -1\right)  ^{l-n}\tfrac{\left(
2\pi\right)  ^{l/2}}{2^{2l+1}}\left(  D^{2n-2}D^{2\left(  l-n\right)  }\left(
\ast\Lambda\right)  ^{l}\right)  \left(  \frac{t}{2}\right) \nonumber\\
& =\left(  -1\right)  ^{l-n}\tfrac{\left(  2\pi\right)  ^{l/2}}{2^{2l+1}%
}\left(  D^{2l-2}\left(  \ast\Lambda\right)  ^{l}\right)  \left(  \frac{t}%
{2}\right) \nonumber\\
& =\left(  -1\right)  ^{l-n}\tfrac{\left(  2\pi\right)  ^{l/2}}{2^{2l+1}%
}\left(  \left(  \ast D^{2}\Lambda\right)  ^{l-1}\ast\Lambda\right)  \left(
\frac{t}{2}\right) \nonumber\\
& =\left(  -1\right)  ^{l-n}\tfrac{\left(  2\pi\right)  ^{l/2}}{2^{2l+1}%
}\left(  \left(  \ast\left(  \delta\left(  \cdot-1\right)  -2\delta
+\delta\left(  \cdot+1\right)  \right)  \right)  ^{l-1}\ast\Lambda\right)
\left(  \frac{t}{2}\right) \nonumber\\
& =\left(  -1\right)  ^{n-1}\tfrac{\left(  2\pi\right)  ^{l/2}}{2^{2l+1}}%
\frac{\left(  -1\right)  ^{l-1}}{\left(  2\pi\right)  ^{\frac{l-2}{2}}}%
\sum_{\substack{k= \\-\left(  l-1\right)  }}^{l-1}\left(  -1\right)
^{k}\tbinom{2l-2}{l-1+k}\left(  \delta\left(  \cdot-k\right)  \ast
\Lambda\right)  \left(  \frac{t}{2}\right) \nonumber\\
& =\left(  -1\right)  ^{n-1}\tfrac{\left(  2\pi\right)  ^{l/2}}{2^{2l+1}}%
\frac{\left(  -1\right)  ^{l-1}}{\left(  2\pi\right)  ^{\frac{l-1}{2}}}%
\sum_{k=-\left(  l-1\right)  }^{l-1}\left(  -1\right)  ^{k}\tbinom
{2l-2}{l-1+k}\Lambda\left(  \frac{t}{2}-k\right) \nonumber\\
& =\left(  -1\right)  ^{n+l}\tfrac{\sqrt{2\pi}}{2^{2l+1}}\sum_{k=-\left(
l-1\right)  }^{l-1}\left(  -1\right)  ^{k}\tbinom{2l-2}{l-1+k}\Lambda\left(
\frac{t}{2}-k\right)  ,\label{2.76}%
\end{align}

which is a continuous, piecewise linear function with bounded support. thus%
\[
\left\Vert D^{2n-2}G_{1}\right\Vert _{\infty}=\tfrac{\sqrt{2\pi}}{2^{2l+1}%
}\max_{k=-\left(  l-1\right)  }^{l-1}\tbinom{2l-2}{l-1+k}=\tfrac{\sqrt{2\pi}%
}{2^{2l+1}}\tbinom{2l-2}{l-1},
\]

and so%
\begin{align*}
\left\Vert D^{2n-1}G_{1}\right\Vert _{\infty}  & =\tfrac{\sqrt{2\pi}}%
{2^{2l+1}}\max_{k=-\left(  l-1\right)  }^{l-1}\left\vert \left(  -1\right)
^{k+1}\tbinom{2l-2}{l-1+k+1}-\left(  -1\right)  ^{k}\tbinom{2l-2}%
{l-1+k}\right\vert /2\\
& =\tfrac{\sqrt{2\pi}}{2^{2l+2}}\max\left\{  \max_{k=-\left(  l-1\right)
}^{l-1}\left\vert \left(  -1\right)  ^{k+1}\tbinom{2l-2}{l-1+k+1}-\left(
-1\right)  ^{k}\tbinom{2l-2}{l-1+k}\right\vert ,2l-2\right\} \\
& =\tfrac{\sqrt{2\pi}}{2^{2l+2}}\max\left\{  \max_{k=-\left(  l-1\right)
}^{l-1}\left(  \tbinom{2l-2}{l-1+k+1}+\tbinom{2l-2}{l-1+k}\right)
,2l-2\right\} \\
& =\tfrac{\sqrt{2\pi}}{2^{2l+2}}\max\left\{  \max_{k=0}^{l-2}\left(
\tbinom{2l-2}{k}+\tbinom{2l-2}{k+1}\right)  ,2l-2\right\} \\
& =\tfrac{\sqrt{2\pi}}{2^{2l+2}}\max\left\{  \max_{k=0}^{l-2}\tbinom
{2l-1}{k+1},2l-2\right\} \\
& =\tfrac{\sqrt{2\pi}}{2^{2l+2}}\max\left\{  \tbinom{2l-1}{l-1},2l-2\right\}
\\
& =\tfrac{\sqrt{2\pi}}{2^{2l+2}}\tbinom{2l-1}{l-1}.
\end{align*}

The other stated properties of $G_{1}$ and $G_{s}$ now follow directly.
\end{proof}

\begin{remark}
\label{Rem_Thm_basis_tensor_hat_W3}?? \textbf{TRY USING YOUNG's inequality to
estimate} $\left\Vert D^{j}G_{s}\right\Vert _{\infty}$ in 1 dimension!

Use the following convolution estimate from Section 12B of Jones
\cite{Jones2011}:

Suppose $\frac{1}{p_{1}^{\prime}}+\frac{1}{p_{2}^{\prime}}+\ldots+\frac
{1}{p_{N}^{\prime}}=\frac{1}{r^{\prime}}$ where $\forall i$ $\frac{1}%
{p_{i}^{\prime}}=1-\frac{1}{p_{i}}$, $\frac{1}{r^{\prime}}=1-\frac{1}{r}$ and
$p_{i},r\geq1$.

Then%
\[
\left\Vert f_{1}\ast f_{2}\ast\ldots\ast f_{N}\right\Vert _{r}\leq\left\Vert
f_{1}\right\Vert _{p_{1}}\left\Vert f_{2}\right\Vert _{p_{2}}\ldots\left\Vert
f_{k}\right\Vert _{p_{N}}.
\]

Suppose the $p_{i}\,$s are all equal. In this case $Nr^{\prime}=p_{i}^{\prime
}$, $p_{i}=\frac{Nr^{\prime}}{Nr^{\prime}-1}=\frac{N}{N-1+1/r}$ and so%
\begin{equation}
\left\Vert f_{1}\ast f_{2}\ast\ldots\ast f_{N}\right\Vert _{r}\leq\left\Vert
f_{1}\right\Vert _{\frac{N}{N-1+1/r}}\left\Vert f_{2}\right\Vert _{\frac
{N}{N-1+1/r}}\ldots\left\Vert f_{k}\right\Vert _{\frac{N}{N-1+1/r}%
}.\label{2.70}%
\end{equation}

\textbf{Thought}: choose different norms for the hat and the rectangle function.

\textbf{Assume} $l=n$ so that $G_{s}\left(  t\right)  =\tfrac{\left(
2\pi\right)  ^{n/2}}{2}\left(  \left(  \ast\Lambda\right)  ^{n}\right)
\left(  \frac{t}{2}\right)  $, $D^{j}G_{s}\left(  t\right)  =\tfrac{\left(
2\pi\right)  ^{n/2}}{2^{j}}D^{j}\left(  \left(  \ast\Lambda\right)
^{n}\right)  \left(  \frac{t}{2}\right)  $ and hence%
\begin{equation}
\left\Vert D^{j}G_{s}\right\Vert _{\infty}=\tfrac{\left(  2\pi\right)  ^{n/2}%
}{2^{j}}\left\Vert D^{j}\left(  \left(  \ast\Lambda\right)  ^{n}\right)
\right\Vert _{\infty},\quad j\leq2n-1.\label{2.71}%
\end{equation}

\textbf{Further assume} that $j\leq n$. Then%
\begin{equation}
\left\Vert D^{j}\left(  \left(  \ast\Lambda\right)  ^{n}\right)  \right\Vert
_{\infty}=\left\Vert \left(  \ast D\Lambda\right)  ^{j}\ast\left(  \ast
\Lambda\right)  ^{n-j}\right\Vert _{\infty}\leq\left\Vert \left(  \ast
D\Lambda\right)  ^{j}\right\Vert _{p}\left\Vert \left(  \ast\Lambda\right)
^{n-j}\right\Vert _{q},\label{2.72}%
\end{equation}

where $\frac{1}{p}+\frac{1}{q}=1$. Let $t:=1/np$. Equation \ref{2.70} now
implies%
\begin{align}
\left\Vert D^{j}\left(  \left(  \ast\Lambda\right)  ^{n}\right)  \right\Vert
_{\infty}  & \leq\left\Vert D\Lambda\right\Vert _{\frac{n}{n-1+1/p}}%
^{j}\left\Vert \Lambda\right\Vert _{\frac{n}{n-1+1/q}}^{n-j}=\left\Vert
D\Lambda\right\Vert _{\frac{n}{n-1+1/p}}^{j}\left\Vert \Lambda\right\Vert
_{\frac{n}{n-1/p}}^{n-j}=\nonumber\\
& =\left\Vert D\Lambda\right\Vert _{\frac{1}{1-\frac{1}{n}+t}}^{j}\left\Vert
\Lambda\right\Vert _{\frac{1}{1-t}}^{n-j}\nonumber\\
& :a:=\frac{1}{1-\frac{1}{n}+t},\text{ }b:=\frac{1}{1-t}\Rightarrow\nonumber\\
& =\left\Vert D\Lambda\right\Vert _{a}^{j}\left\Vert \Lambda\right\Vert
_{b}^{n-j}\nonumber\\
& =\left(  \int_{-1}^{1}1^{a}\right)  ^{j/a}\left(  \int_{-1}^{1}\left(
1-\left\vert t\right\vert \right)  ^{b}dt\right)  ^{\left(  n-j\right)
/b}\nonumber\\
& =\left(  2\int_{0}^{1}1\right)  ^{j/a}\left(  2\int_{0}^{1}\left(
1-t\right)  ^{b}dt\right)  ^{\left(  n-j\right)  /b}\nonumber\\
& =2^{j/a}2^{\left(  n-j\right)  /b}\left(  \int_{0}^{1}s^{b}ds\right)
^{\left(  n-j\right)  /b}\nonumber\\
& =2^{j/a}2^{\left(  n-j\right)  /b}\left(  \frac{1}{b+1}\right)  ^{\left(
n-j\right)  /b}\nonumber\\
& =2^{j/a}2^{\left(  n-j\right)  /b}\left(  \frac{1}{2}\frac{2}{b+1}\right)
^{\left(  n-j\right)  /b}\nonumber\\
& =2^{j/a}\left(  \frac{2}{b+1}\right)  ^{\left(  n-j\right)  /b}\nonumber\\
& =2^{j\left(  1-\frac{1}{n}+\frac{1}{np}\right)  }\left(  \frac{2}{\frac
{1}{1-t}+1}\right)  ^{\left(  n-j\right)  \left(  1-t\right)  }\nonumber\\
& =2^{j\left(  1-\frac{1}{n}+\frac{1}{np}\right)  }\left(  \frac{2-2t}%
{2-t}\right)  ^{\left(  n-j\right)  \left(  1-t\right)  }\nonumber\\
& =2^{j\left(  1-\frac{1}{n}+\frac{1}{np}\right)  }\left(  1-\frac{t}%
{2-t}\right)  ^{\left(  n-j\right)  \left(  1-t\right)  }\nonumber\\
& =2^{j\left(  1-\frac{1}{n}+\frac{1}{np}\right)  }\left(  1-\frac{1}%
{2np-1}\right)  ^{\left(  n-j\right)  \left(  1-\frac{1}{np}\right)  },\text{
}p\geq1.\label{2.73}%
\end{align}

Thus%
\begin{align*}
\left\Vert D^{j}G_{s}\right\Vert _{\infty}  & \leq\left(  2\pi\right)
^{n/2}2^{-j}\left\Vert D^{j}\left(  \left(  \ast\Lambda\right)  ^{n}\right)
\right\Vert _{\infty}\\
& =\left(  2\pi\right)  ^{n/2}2^{-j}2^{j\left(  1-\frac{1}{n}+\frac{1}%
{np}\right)  }\left(  1-\frac{1}{2np-1}\right)  ^{\left(  n-j\right)  \left(
1-\frac{1}{np}\right)  }\\
& =\left(  2\pi\right)  ^{n/2}2^{-\frac{j}{n}\left(  1-\frac{1}{p}\right)
}\left(  1-\frac{1}{2np-1}\right)  ^{\left(  n-j\right)  \left(  1-\frac
{1}{np}\right)  }.
\end{align*}

From \ref{1.048} and the definition of the basis function, when $j$ is
even\allowbreak\ $\left\Vert D^{j}G_{s}\right\Vert _{\infty}\leq\left\vert
D^{j}G_{s}\left(  0\right)  \right\vert =\left(  2\pi\right)  ^{-1/2}\int%
\frac{\sin^{2n}s}{s^{2n-j}}$ and so%
\begin{align}
\int\frac{\sin^{2n}s}{s^{2n-j}}  & \leq\left(  \sqrt{2\pi}\right)
^{n+1}2^{-\frac{j}{n}\left(  1-\frac{1}{p}\right)  }\left(  1-\frac{1}%
{2np-1}\right)  ^{\left(  n-j\right)  \left(  1-\frac{1}{np}\right)  },\quad
j\text{ }even,\text{ }0\leq j\leq n.\label{2.74}\\
& \Rightarrow\nonumber\\
\int_{0}^{\infty}\frac{\sin^{2n}s}{s^{m}}  & \leq\left(  \sqrt{2\pi}\right)
^{n+1}2^{-\frac{2m-n}{n}\left(  1-\frac{1}{p}\right)  }\left(  1-\frac
{1}{2np-1}\right)  ^{\left(  m-n\right)  \left(  1-\frac{1}{np}\right)
},\quad m\text{ }even,\text{ }n\leq m\leq2n.\label{2.75}%
\end{align}

??? \textbf{The estimate on the RHS is far, far greater than the RHS}.
\end{remark}

\subsection{Convolutions and products of weight and basis functions: W02 case
\label{SbSect_convol_prod_wt_basis_fns}}

The results of this subsection will not be used later in this document. In
this subsection we give several results concerning the product and convolution
of basis functions. The next result gives conditions under which the product
of two basis functions of order zero is a basis function of order zero.

\begin{theorem}
Suppose the weight functions $w_{1}$ and $w_{2}$ have property \textbf{W02}
for parameter $\kappa_{i}$. Suppose $G_{1}$ and $G_{2}$ are the basis
functions of order zero generated by the weight functions $w_{1}$ and $w_{2}$.
Suppose that $1/w_{1}$ is bounded a.e. or $G_{1}G_{2}\in L^{1}$.

Then the product $G_{1}G_{2}$ is a basis function of order zero generated by a
weight function $w$ which satisfies property W02 for parameter $\kappa
=\min\left\{  \kappa_{1},\kappa_{2}\right\}  $.

Further, $w\in C^{\left(  0\right)  }$ and $w\left(  x\right)  >0$ for all $x$
i.e. $\widehat{G_{1}G_{2}}\left(  \xi\right)  >0$ for all $\xi$.
\end{theorem}

\begin{proof}
We first note that from Theorem 2.2 and Exercise 2.4 of Petersen
\cite{Petersen83}, that if $f,g\in L^{1}$ then $f\ast g\in L^{1}$,
$\widehat{f\ast g}=\widehat{f}\widehat{g}$ and, if $f$ is bounded a.e., $f\ast
g$ is a continuous, bounded function.

Now define the functions $G$ and $w$ by $\tfrac{1}{w}=\tfrac{1}{w_{1}}%
\ast\tfrac{1}{w_{2}}\in L^{1}$ and $\widehat{G}=\frac{1}{w}$. We now have
$\left(  \tfrac{1}{w}\right)  ^{\vee}=\left(  \tfrac{1}{w_{1}}\right)  ^{\vee
}\left(  \tfrac{1}{w_{2}}\right)  ^{\vee}$ i.e. $G=G_{1}G_{2}$.

If $G_{1}G_{2}\in L^{1}$ then $\tfrac{1}{w}\in C_{B}^{\left(  0\right)  }$,
and if $\tfrac{1}{w_{1}}$ is bounded a.e. then $\tfrac{1}{w_{1}}\ast\tfrac
{1}{w_{2}}=\tfrac{1}{w}\in C_{B}^{\left(  0\right)  }$. Thus $\left|  w\left(
x\right)  \right|  \geq c>0$ for some constant $c$, and is continuous whenever
$w$ is finite. But $w_{1}\geq0$ and $w_{2}\geq0$ a.e. so $w\left(  x\right)
\geq c$ for all $x$.

The next step is to show that $w\left(  x\right)  <\infty$ for all $x$. By
definition of $w_{1}$ and $w_{2}$ there exist closed sets of measure zero,
$\mathcal{A}_{1}$ and $\mathcal{A}_{2}$ such that $w_{i}$ is continuous and
positive outside $\mathcal{A}_{i}$. For each $x\in\mathbb{R}^{d}$, let
$\mathcal{B}_{x}=\left(  x-\mathcal{A}_{1}\right)  \cup\mathcal{A}_{2}$ and
note that $\mathcal{B}_{x}$ is closed with measure zero, and that, as a
function of $y$, $w_{1}\left(  x-y\right)  w_{2}\left(  y\right)  $ is
continuous and positive on $\mathbb{R}^{d}\setminus\mathcal{B}_{x}$. Thus, we
have for all $x$%
\begin{equation}
\frac{1}{w\left(  x\right)  }=\tfrac{1}{\left(  2\pi\right)  ^{d/2}}%
\int_{\mathbb{R}^{d}\setminus\mathcal{B}_{x}}\frac{dy}{w_{1}\left(
x-y\right)  w_{2}\left(  y\right)  }>0,\label{1.025}%
\end{equation}

and consequently that $w\left(  x\right)  <\infty$ for all $x$ i.e. $w\in
C^{\left(  0\right)  }\left(  \mathbb{R}\right)  $.

It remains to be shown that $w$ satisfies property W02 for $\kappa
=\min\left\{  \kappa_{1},\kappa_{2}\right\}  $. Suppose $0\leq\lambda
\leq\kappa$. Then, for some constant $c_{\lambda}>0$%
\begin{align*}
\int\frac{\left\vert x\right\vert ^{2\lambda}dx}{w\left(  x\right)  }  & \leq
c_{\lambda}\int\int\frac{\left\vert x-y\right\vert ^{2\lambda}dy\,dx}%
{w_{1}\left(  x-y\right)  w_{2}\left(  y\right)  }+c_{\lambda}\int\int%
\frac{\left\vert y\right\vert ^{2\lambda}dy\,dx}{w_{1}\left(  x-y\right)
w_{2}\left(  y\right)  }\\
& =c_{\lambda}\int\frac{\left\vert x-y\right\vert ^{2\lambda}dx}{w_{1}\left(
x-y\right)  }\int\frac{dy}{w_{2}\left(  y\right)  }+c_{\lambda}\int%
\frac{\left\vert y\right\vert ^{2\lambda}dy}{w_{2}\left(  y\right)  }\int%
\frac{dx}{w_{1}\left(  x-y\right)  }\\
& =c_{\lambda}\int\frac{\left\vert \cdot\right\vert ^{2\lambda}}{w_{1}}%
\int\frac{1}{w_{2}}+c_{\lambda}\int\frac{\left\vert \cdot\right\vert
^{2\lambda}}{w_{2}}\int\frac{1}{w_{1}}\\
& <\infty.
\end{align*}

\end{proof}

The last theorem can be expressed as the following weight function result.

\begin{corollary}
Suppose $w_{1}$ and $w_{2}$ are two weight functions which satisfy
\textbf{W02} for $\kappa=\kappa_{1}$ and $\kappa=\kappa_{2}$ respectively.
Further, suppose that $1/w_{1}$ is bounded. Define the function $w$ by
$\frac{1}{w}=\frac{1}{w_{1}}\ast\frac{1}{w_{2}}$.

Then $w$ is a weight function which has property W02 for $\kappa=\min\left\{
\kappa_{1},\kappa_{2}\right\}  $.

Further, $w\in C^{\left(  0\right)  }\left(  \mathbb{R}\right)  $ and
$w\left(  x\right)  >0$ for all $x$.
\end{corollary}

\begin{remark}
\label{Rem_Thm_prod_basis_funcs} Any basis function can be uniformly pointwise
approximated by a sequence of basis functions which have Fourier transforms
that are always positive.

In fact, suppose $G$ is a basis function for which $\widehat{G}\left(
\xi\right)  \ngtr0$ for all $\xi$. Define the sequence of functions
$G_{n}\left(  \xi\right)  =\Lambda\left(  \xi/n\right)  G\left(  \xi\right)
,$ $n=1,2,3,\ldots$. Clearly, by Theorem \ref{Thm_hat_wt_extend_props} and the
previous theorem, $G_{n}$ is a basis function for $\kappa_{3}<\min\left\{
1/2,\kappa_{2}\right\}  $ and $\widehat{G_{n}}\left(  \xi\right)  >0$ for all
$\xi$.

Further, $G_{n}\rightarrow G$ uniformly pointwise. This is true because the
definition of $\Lambda$ implies
\[
\left\vert G_{n}\left(  \xi\right)  -G\left(  \xi\right)  \right\vert
\leq\left\{
\begin{array}
[c]{ll}%
\dfrac{\left\vert \xi\right\vert }{n}\left\vert G\left(  \xi\right)
\right\vert , & \left\vert \xi\right\vert \leq n,\\
& \\
\left\vert G\left(  \xi\right)  \right\vert , & \left\vert \xi\right\vert >n,
\end{array}
\right.
\]
and the definition of $G$ implies $\lim\limits_{\left\vert \xi\right\vert
\rightarrow\infty}G\left(  \xi\right)  =0$.

A specific example is $G=\Lambda$ for which $\left\|  G_{n}-G\right\|
_{\infty}\leq1/n$.
\end{remark}

The next result gives conditions under which the convolution of two basis
functions is a basis function. This is equivalent to a result showing that the
product of two weight functions is a weight function.

\begin{theorem}
\label{Thm_convol_basis_funcs}Suppose $G_{1}$ and $G_{2}$ are the basis
functions of order zero generated by the weight functions $w_{1}$ and $w_{2}$
respectively, and that each $w_{i}$ has property \textbf{W02} for parameter
$\kappa_{i}$. Then the following two results hold.

\begin{enumerate}
\item If $G_{1}\in L^{1}$, the convolution
\[
G_{1}\ast G_{2}=\tfrac{1}{\left(  2\pi\right)  ^{d/2}}\int G_{1}\left(
x-y\right)  G_{2}\left(  y\right)  dy,
\]
is the basis function generated by the weight function $w=w_{1}w_{2}$, and $w
$ satisfies property W02 for $\kappa=\kappa_{2}$.

\item If we further assume that for an integer $n_{1}\leq\kappa_{1}$,
$D^{2\alpha}G_{1}\in L^{1}$ when $\left\vert \alpha\right\vert =n_{1}$, then
$w$ satisfies property W02 for parameter $\kappa=n_{1}+\kappa_{2}$ and
\begin{equation}
D^{\gamma+\delta}\left(  G_{1}\ast G_{2}\right)  =\tfrac{1}{\left(
2\pi\right)  ^{d/2}}\int D^{\gamma}G_{1}\left(  x-y\right)  D^{\delta}%
G_{2}\left(  y\right)  dy,\text{\quad}\left\vert \gamma\right\vert \leq
n_{1},\text{ }\left\vert \delta\right\vert \leq\kappa_{2}.\label{1.023}%
\end{equation}

\end{enumerate}
\end{theorem}

\begin{proof}
\textbf{Part 1}. Define the function $w=w_{1}w_{2}$. Clearly $w$ is a weight
function since if $\mathcal{A}_{1}$ and $\mathcal{A}_{2}$ are the weight sets,
$w\left(  x\right)  >0$ and $w\in C^{\left(  0\right)  }$ outside the set
$\mathcal{A}_{1}\cup\mathcal{A}_{2}$, which is closed and has measure zero.
Since $\widehat{G_{1}}=\tfrac{1}{w_{1}}\in C_{B}^{\left(  0\right)  }$ we have
$\tfrac{1}{w_{1}w_{2}}\in L^{1}$ and are able to define the function
$\widehat{G}=\tfrac{1}{w_{1}w_{2}}$. We now show that $G=G_{1}\ast G_{2}$.
First observe that because $G_{1}\in L^{1}$ and $G_{2}\in C_{B}^{\left(
0\right)  }$ the convolution integral exists. Now, by definition of $G_{2}$%
\[
\int G_{1}\left(  x-y\right)  G_{2}\left(  y\right)  dy=\tfrac{1}{\left(
2\pi\right)  ^{d/2}}\int\int G_{1}\left(  x-y\right)  \frac{e^{iy\xi}}%
{w_{2}\left(  \xi\right)  }d\xi\,dy,
\]

and, since that this integral is absolutely convergent, we can change the
order of integration so that
\begin{align*}
\int G_{1}\left(  x-y\right)  G_{2}\left(  y\right)  dy  & =\tfrac{1}{\left(
2\pi\right)  ^{d/2}}\int\frac{1}{w_{2}\left(  \xi\right)  }\int e^{iy\xi}%
G_{1}\left(  x-y\right)  dy\,d\xi\\
& =\tfrac{1}{\left(  2\pi\right)  ^{d/2}}\int\frac{1}{w_{2}\left(  \xi\right)
}\int e^{i\left(  x-z\right)  \xi}G_{1}\left(  z\right)  dz\,d\xi\\
& =\int\frac{e^{ix\xi}}{w_{2}\left(  \xi\right)  }\tfrac{1}{\left(
2\pi\right)  ^{d/2}}\int e^{-iz\xi}G_{1}\left(  z\right)  dz\,d\xi.
\end{align*}

The final step uses Corollary 3.7 of Petersen \cite{Petersen83}. This states
that if $f\in L^{1}$ and $\widehat{f}\in L^{1}$ then $f\left(  x\right)
=\tfrac{1}{\left(  2\pi\right)  ^{d/2}}\int e^{ix\xi}\widehat{f}\left(
\xi\right)  d\xi$ a.e. In our case we choose $f\left(  x\right)  =\tfrac
{1}{w_{1}\left(  -x\right)  }$ and obtain
\[
\tfrac{1}{\left(  2\pi\right)  ^{d/2}}\int G_{1}\left(  x-y\right)
G_{2}\left(  y\right)  dy=\tfrac{1}{\left(  2\pi\right)  ^{d/2}}\int%
\frac{e^{ix\xi}d\xi}{w_{1}\left(  \xi\right)  w_{2}\left(  \xi\right)
}=G\left(  x\right)  .
\]

Regarding property W02, if $0\leq\lambda\leq\kappa_{2}$ then, $\int%
\frac{\left\vert \cdot\right\vert ^{2\lambda}}{w_{1}w_{2}}\leq\left\Vert
\frac{1}{w_{1}}\right\Vert _{\infty}\int\frac{\left\vert \cdot\right\vert
^{2\lambda}}{w_{2}}<\infty$.\medskip

\textbf{Part 2}. When $\left\vert \alpha\right\vert =n_{1}$, we have
$\widehat{D^{2\alpha}G_{1}}=\left(  -1\right)  ^{\left\vert \alpha\right\vert
}\xi^{2\alpha}\widehat{G_{1}}=\left(  -1\right)  ^{n_{1}}\frac{\xi^{2\alpha}%
}{w_{1}}$ and $\widehat{D^{2\alpha}G_{1}}\in C_{B}^{\left(  0\right)  }$.
Thus
\[
\int\frac{\left\vert \cdot\right\vert ^{2\kappa}}{w}=\int\frac{\left\vert
\cdot\right\vert ^{2n_{1}}\left\vert \cdot\right\vert ^{2\kappa_{2}}}%
{w_{1}w_{2}}=\sum_{\left\vert \alpha\right\vert =n_{1}}\frac{1}{\alpha!}%
\int\frac{\xi^{2\alpha}}{w_{1}\left(  \xi\right)  }\frac{\left\vert
\xi\right\vert ^{2\kappa_{2}}}{w_{2}\left(  \xi\right)  }d\xi<\infty,
\]

and since $\frac{1}{w}\in L^{1}$, $w$ satisfies W02 for $\kappa=n_{1}%
+\kappa_{2}$.

Suppose that $\left\vert \gamma\right\vert \leq n_{1}$and $\left\vert
\delta\right\vert \leq\kappa_{2}$. Since $G_{1}\in C^{\left(  2\kappa
_{1}\right)  }$, and

$\int D^{\gamma}G_{1}\left(  x-y\right)  G_{2}\left(  y\right)  dy$ is
absolutely convergent, we have
\[
D^{\gamma}\left(  G_{1}\ast G_{2}\right)  =\tfrac{1}{\left(  2\pi\right)
^{d/2}}\int D^{\gamma}G_{1}\left(  x-y\right)  G_{2}\left(  y\right)
dy=\tfrac{1}{\left(  2\pi\right)  ^{d/2}}\int D^{\gamma}G_{1}\left(  z\right)
G_{2}\left(  x-z\right)  dz.
\]
Since $G_{2}\in C^{\left(  2\kappa_{2}\right)  }$ and $\int D^{\gamma}%
G_{1}\left(  z\right)  D^{\delta}G_{2}\left(  x-z\right)  dz$ is absolutely
convergent we have our result.
\end{proof}

\subsection{Convolutions and products of weight and basis functions: W03 case}

?? No theorems yet!

\section{The Riesz representer of the evaluation functionals $f\rightarrow
\left(  D^{\alpha}f\right)  \left(  x\right)  $%
\label{Sect_reps_of_eval_fns_W2}}

It was noted in Remark \ref{Rem_reprod_Hilbert} that $X_{w}^{0}$ is a
reproducing kernel Hilbert space. However, we will use the Riesz representer
of the evaluation functional $f\rightarrow f\left(  x\right)  $ instead of the
reproducing kernel because the former can be easily extended to the evaluation
of derivatives. Thus in this section we calculate the Riesz representers of
the evaluation functionals $f\rightarrow\left(  D^{\alpha}f\right)  \left(
x\right)  $ and derive some of their properties in two theorems; Theorem
\ref{Thm_ord0_Riesz_rep_W2} assumes the weight function has property W02 and
Theorem \ref{Thm_ord0_Riesz_rep_W3}

assumes the weight function has property W03. In particular, in part 6 of
these two results we show that the inclusion $X_{w}^{0}\subset C_{B}^{\left(
\left\lfloor \kappa\right\rfloor \right)  }$ established in Theorem
\ref{Thm_X_smooth} is continuous when $C_{B}^{\left(  \left\lfloor
\kappa\right\rfloor \right)  }$ is endowed with the usual supremum norm. The
results of this section will find use in applications where, for example,
higher degrees of differentiability of a smoother are required. The main tools
used here are the inverse Fourier transform formulas of Theorems
\ref{Thm_basis_fn_properties_all_m_W2} and
\ref{Thm_basis_fn_properties_all_m_W3}.

\begin{theorem}
\label{Thm_ord0_Riesz_rep_W2}Suppose the weight function $w$ has property
\textbf{W02} for some $\kappa$. Then:

\begin{enumerate}
\item The unique Riesz representer $R_{x}\in X_{w}^{0}$ of the evaluation
functional $f\rightarrow f\left(  x\right)  $ is
\begin{equation}
R_{x}(z):=\left(  2\pi\right)  ^{-d/2}G(z-x).\label{1.019}%
\end{equation}

\item If $\left\vert \alpha\right\vert \leq\kappa$ then $D^{\alpha}R_{x}\in
X_{w}^{0}$ and the Riesz representer for the evaluation functional
$f\rightarrow\left(  D^{\alpha}f\right)  \left(  x\right)  $ is $\left(
-D\right)  ^{\alpha}R_{x}$.

\item If $\left\vert \alpha\right\vert \leq2\kappa$ then $x\neq x^{\prime}$
implies $D^{\alpha}R_{x}\neq D^{\alpha}R_{x^{\prime}}$.

\item If $\left\vert \alpha\right\vert \leq\kappa$ and $\left\vert
\beta\right\vert \leq\kappa$ then%
\begin{equation}
\left(  D^{\alpha}R_{x},D^{\beta}R_{y}\right)  _{w,0}=\left(  -1\right)
^{\left\vert \beta\right\vert }\left(  D^{\alpha+\beta}R_{x}\right)  \left(
y\right)  =\left(  -1\right)  ^{\left\vert \alpha+\beta\right\vert }\left(
D_{x}^{\alpha}R_{x},D^{\beta}R_{y}\right)  _{w,0}.\label{1.6}%
\end{equation}

\item If $\left\vert \gamma\right\vert \leq2\kappa$ then
\begin{equation}
\left(  D^{\gamma}R_{x}\right)  \left(  y\right)  =\left(  -1\right)
^{\left\vert \gamma\right\vert }\overline{\left(  D^{\gamma}R_{y}\right)
\left(  x\right)  }.\label{1.013}%
\end{equation}

\item $\max\limits_{\left\vert \beta\right\vert \leq\kappa}\left\Vert
D^{\beta}f\right\Vert _{\infty}\leq\left(  2\pi\right)  ^{-\frac{d}{4}}%
\max\limits_{\left\vert \beta\right\vert \leq\kappa}\sqrt{\left(  -1\right)
^{\left\vert \beta\right\vert }D^{2\beta}G\left(  0\right)  }\left\Vert
f\right\Vert _{w,0}$,\quad$f\in X_{w}^{0}$.
\end{enumerate}
\end{theorem}

\begin{proof}
\textbf{Part 1.} From \ref{1.019},
\begin{align}
\widehat{R_{x}}\left(  \xi\right)   & =\left(  2\pi\right)  ^{-d/2}e^{-ix\xi
}\widehat{G}\left(  \xi\right) \label{1.44}\\
& =\left(  2\pi\right)  ^{-d/2}\frac{e^{-ix\xi}}{w\left(  \xi\right)
}.\label{1.45}%
\end{align}

By Theorem \ref{Thm_basis_fn_properties_all_m_W2}, $G\in X_{w}^{0}$ and so
\ref{1.44} implies $R_{x}\in X_{w}^{0}$. Further, using \ref{1.45} and then
the inverse Fourier transform result \ref{1.021} for functions in $X_{w}^{0}$
we obtain
\[
\left(  f,R_{x}\right)  _{w,0}=\int w\widehat{f}\overline{\widehat{R_{x}}%
}=\left(  2\pi\right)  ^{-d/2}\int e^{ix\xi}\widehat{f}\left(  \xi\right)
d\xi=f\left(  x\right)  .
\]
\medskip

\textbf{Part 2.} If $\left\vert \alpha\right\vert \leq\kappa$ then
$\widehat{D^{\alpha}R_{x}}=\left(  i\xi\right)  ^{\alpha}\widehat{R_{x}%
}=\left(  2\pi\right)  ^{-d/2}\frac{\left(  i\xi\right)  ^{\alpha}e^{-ix\xi}%
}{w\left(  \xi\right)  }$ so that%
\[
\left\Vert D^{\alpha}R_{x}\right\Vert _{w,0}^{2}=\int w\left(  \xi\right)
\left\vert \left(  2\pi\right)  ^{-\frac{d}{2}}\frac{\left(  i\xi\right)
^{\alpha}e^{-ix\xi}}{w\left(  \xi\right)  }\right\vert ^{2}d\xi\leq\left(
2\pi\right)  ^{-d}\int\frac{\left\vert \xi\right\vert ^{2\left\vert
\alpha\right\vert }}{w\left(  \xi\right)  }d\xi<\infty,
\]

by weight function property W02. Further, for $f\in X_{w}^{0}$ and using the
inverse Fourier transform result \ref{1.021} for functions in $X_{w}^{0}$
gives
\begin{align*}
\left(  f,\left(  -D\right)  ^{\alpha}R_{x}\right)  _{w,0}=\int w\left(
\xi\right)  \widehat{f}\left(  \xi\right)  \overline{\widehat{\left(
-D\right)  ^{\alpha}R_{x}}}\left(  \xi\right)  d\xi &  =\int w\left(
\xi\right)  \widehat{f}\left(  \xi\right)  \left(  2\pi\right)  ^{-\frac{d}%
{2}}\frac{\left(  i\xi\right)  ^{\alpha}e^{-ix\xi}}{w\left(  \xi\right)  }%
d\xi\\
&  =\left(  2\pi\right)  ^{-\frac{d}{2}}\int e^{-ix\xi}\left(  i\xi\right)
^{\alpha}\widehat{f}\left(  \xi\right)  d\xi\\
&  =\left(  2\pi\right)  ^{-d/2}\int e^{-ix\xi}\widehat{D^{\alpha}f}\left(
\xi\right)  d\xi\\
&  =\left(  D^{\alpha}f\right)  \left(  x\right)  .
\end{align*}
\medskip

\textbf{Part 3.} Suppose for some $x\neq x^{\prime}$ we have $R_{x}%
=R_{x^{\prime}}$. Taking the Fourier transform and substituting $\widehat{G}%
=\frac{1}{w}$ we must have $e^{-i\xi x}=e^{-i\xi x^{\prime}}$ or
$e^{i\xi\left(  x-x^{\prime}\right)  }=1$ for almost all $\xi$, a
contradiction.\medskip

\textbf{Part 4}. Since we know that $G\in C^{\left(  2\kappa\right)  }$ the
definition of $R_{x}$ implies $R_{x}\in C^{\left(  2\kappa\right)  }$ and
\begin{equation}
\widehat{D^{\alpha}R_{x}}\left(  \xi\right)  =\left(  2\pi\right)  ^{\frac
{d}{2}}\left(  -i\right)  ^{\left\vert \alpha\right\vert }e^{-ix\xi}%
\xi^{\alpha}\widehat{G}\left(  \xi\right)  =\left(  2\pi\right)  ^{\frac{d}%
{2}}\left(  -i\right)  ^{\left\vert \alpha\right\vert }\frac{e^{-ix\xi}%
\xi^{\alpha}}{w\left(  \xi\right)  },\quad\left\vert \alpha\right\vert
\leq\kappa.\label{1.037}%
\end{equation}

Next note that Theorem \ref{Thm_basis_fn_properties_all_m_W2} implies that for
all $x$, and $D^{\gamma}R_{x}\in X_{w}^{0}$ when $\left\vert \gamma\right\vert
\leq\kappa$, so the first and last terms of \ref{1.6} make sense. Now using
equation \ref{1.037}
\begin{align*}
\left(  D^{\alpha}R_{x},D^{\beta}R_{y}\right)  _{w,0}=\int w\widehat{D^{\alpha
}R_{x}}\overline{\widehat{D^{\beta}R_{y}}} &  =\left(  -1\right)  ^{\left\vert
\beta\right\vert }\left(  2\pi\right)  ^{-d}\int\left(  -i\xi\right)
^{\alpha+\beta}\frac{e^{i\left(  y-x\right)  \xi}}{w\left(  \xi\right)  }%
d\xi\\
&  =\left(  -1\right)  ^{\left\vert \beta\right\vert }\left(  2\pi\right)
^{-d}\int\left(  -i\xi\right)  ^{\alpha+\beta}e^{i\left(  y-x\right)  \xi
}\widehat{G}\left(  \xi\right)  d\xi\\
&  =\left(  -1\right)  ^{\left\vert \beta\right\vert }\left(  2\pi\right)
^{-d}\int e^{i\left(  y-x\right)  \xi}\widehat{D^{\alpha+\beta}G}\left(
\xi\right)  d\xi\\
&  =\left(  -1\right)  ^{\left\vert \beta\right\vert }\left(  2\pi\right)
^{-d}\left(  D^{\alpha+\beta}G\right)  \left(  y-x\right)  ,
\end{align*}

where the last step used the inverse Fourier transform rule \ref{1.024} for
basis functions. Finally, substituting \ref{1.019} gives $\left(  D^{\alpha
}R_{x},D^{\beta}R_{y}\right)  _{w,0}=\left(  -1\right)  ^{\left\vert
\beta\right\vert }\left(  D^{\alpha+\beta}R_{x}\right)  \left(  y\right)  $.

Regarding the second equation:%
\begin{align*}
\widehat{D_{x}^{\alpha}R_{x}}\left(  \xi\right)  =\widehat{D_{x}^{\alpha
}G\left(  \cdot-x\right)  }\left(  \xi\right)  =\left(  -1\right)
^{\left\vert \alpha\right\vert }\widehat{\left(  D^{\alpha}G\right)  \left(
\cdot-x\right)  }\left(  \xi\right)   &  =\left(  -1\right)  ^{\left\vert
\alpha\right\vert }e^{-ix\xi}\widehat{D^{\alpha}G}\left(  \xi\right) \\
&  =\left(  -1\right)  ^{\left\vert \alpha\right\vert }\widehat{D^{\alpha
}G\left(  \cdot-x\right)  }\left(  \xi\right) \\
&  =\left(  -1\right)  ^{\left\vert \alpha\right\vert }\widehat{D^{\alpha
}R_{x}}\left(  \xi\right)  ,
\end{align*}

so that%
\[
\left(  D^{\alpha}R_{x},D^{\beta}R_{y}\right)  _{w,0}=\left(  -1\right)
^{\left\vert \alpha+\beta\right\vert }\left(  D^{\alpha}R_{x},D^{\beta}%
R_{y}\right)  _{w,0}.
\]
\medskip

\textbf{Part 5}. If $\left\vert \gamma\right\vert \leq2\kappa$ then
$\gamma=\alpha+\beta$ for some $\alpha$ and $\beta$ satisfying $\left\vert
\alpha\right\vert \leq\kappa$ and $\left\vert \beta\right\vert \leq\kappa$.
Applying part 3 twice we get%
\begin{align*}
D^{\gamma}R_{x}\left(  y\right)  =D^{\alpha+\beta}R_{x}\left(  y\right)
=\left(  -1\right)  ^{\left\vert \beta\right\vert }\left(  D^{\alpha}%
R_{x},D^{\beta}R_{y}\right)  _{w,0} &  =\left(  -1\right)  ^{\left\vert
\beta\right\vert }\overline{\left(  D^{\beta}R_{y},D^{\alpha}R_{x}\right)
}_{w,0}\\
&  =\left(  -1\right)  ^{\left\vert \alpha+\beta\right\vert }\overline
{D^{\alpha+\beta}R_{y}\left(  x\right)  }\\
&  =\left(  -1\right)  ^{\left\vert \gamma\right\vert }\overline{D^{\gamma
}R_{y}\left(  x\right)  }.
\end{align*}
\medskip

\textbf{Part 6} From part 2, $\left\Vert D^{\beta}f\right\Vert _{\infty}%
=\sup\limits_{x\in\mathbb{R}^{d}}\left\vert \left(  f,\left(  -D\right)
^{\beta}R_{x}\right)  _{w,0}\right\vert $ when $\left\vert \beta\right\vert
\leq\kappa$ so that
\begin{align*}
\left\Vert D^{\beta}f\right\Vert _{\infty}\leq\sup_{x\in\mathbb{R}^{d}%
}\left\Vert D^{\beta}R_{x}\right\Vert _{w,0}\left\Vert f\right\Vert _{w,0} &
=\sup_{x\in\mathbb{R}^{d}}\sqrt{\left(  D^{\beta}R_{x},D^{\beta}R_{x}\right)
_{w,0}}\left\Vert f\right\Vert _{w,0}\\
&  =\sup_{x\in\mathbb{R}^{d}}\sqrt{\left(  -1\right)  ^{\left\vert
\beta\right\vert }\left(  D^{2\beta}R_{x}\right)  \left(  x\right)
}\left\Vert f\right\Vert _{w,0}.
\end{align*}

But from part 1, $\left(  D^{2\beta}R_{x}\right)  \left(  x\right)  =\left(
2\pi\right)  ^{-\frac{d}{2}}D^{2\beta}G(0)$ and thus%
\[
\left\Vert D^{\beta}f\right\Vert _{\infty}\leq\left(  2\pi\right)
^{-d/4}\sqrt{\left(  -1\right)  ^{\left\vert \beta\right\vert }D^{2\beta}%
G(0)}\left\Vert f\right\Vert _{w,0},
\]

which implies this part.
\end{proof}

In the next theorem we assume that the weight function has property W03
instead of property W02.

\begin{theorem}
\label{Thm_ord0_Riesz_rep_W3}Suppose the weight function $w$ has property
\textbf{W03} for some $\kappa$. Then:

\begin{enumerate}
\item The unique Riesz representer $R_{x}\in X_{w}^{0}$ of the evaluation
functional $f\rightarrow f\left(  x\right)  $ is
\begin{equation}
R_{x}(z)=\left(  2\pi\right)  ^{-d/2}G(z-x).\label{1.069}%
\end{equation}

\item If $\alpha\leq\kappa$ then $D^{\alpha}R_{x}\in X_{w}^{0}$ and the Riesz
representer for the evaluation functional $f\rightarrow\left(  D^{\alpha
}f\right)  \left(  x\right)  $ is $\left(  -D\right)  ^{\alpha}R_{x}$.

\item If $\alpha\leq2\kappa$ then $x\neq x^{\prime}$ implies $D^{\alpha}%
R_{x}\neq D^{\alpha}R_{x^{\prime}}$.

\item If $\alpha\leq\kappa$ and $\beta\leq\kappa$ then%
\begin{equation}
\left(  D^{\alpha}R_{x},D^{\beta}R_{y}\right)  _{w,0}=\left(  -1\right)
^{\left\vert \beta\right\vert }\left(  D^{\alpha+\beta}R_{x}\right)  \left(
y\right)  =\left(  -1\right)  ^{\left\vert \alpha+\beta\right\vert }\left(
D_{x}^{\alpha}R_{x},D^{\beta}R_{y}\right)  _{w,0}.\label{1.072}%
\end{equation}

\item If $\gamma\leq2\kappa$ then
\begin{equation}
\left(  D^{\gamma}R_{x}\right)  \left(  y\right)  =\left(  -1\right)
^{\left\vert \gamma\right\vert }\overline{\left(  D^{\gamma}R_{y}\right)
\left(  x\right)  }.\label{1.073}%
\end{equation}

\item If $\beta\leq\kappa$ and $x\in\mathbb{R}^{d}$ then%
\[
\left\vert D^{\beta}f\left(  x\right)  \right\vert \leq\left(  2\pi\right)
^{-d/4}\sqrt{\left(  -1\right)  ^{\left\vert \beta\right\vert }D^{2\beta}%
G(0)}\left\Vert f\right\Vert _{w,0},\quad f\in X_{w}^{0}.
\]

\end{enumerate}
\end{theorem}

\begin{proof}
\textbf{Part 1.} From \ref{1.069},
\begin{align}
\widehat{R_{x}}\left(  \xi\right)   & =\left(  2\pi\right)  ^{-d/2}e^{-ix\xi
}\widehat{G}\left(  \xi\right) \label{1.441}\\
& =\left(  2\pi\right)  ^{-d/2}\frac{e^{-ix\xi}}{w\left(  \xi\right)
}.\label{1.451}%
\end{align}

By Theorem \ref{Thm_basis_fn_properties_all_m_W3}, $G\in X_{w}^{0}$ and so
\ref{1.441} implies $R_{x}\in X_{w}^{0}$. Further, using \ref{1.451} and then
the inverse Fourier transform result \ref{1.065} for functions in $X_{w}^{0}$
we obtain
\[
\left(  f,R_{x}\right)  _{w,0}=\int w\widehat{f}\overline{\widehat{R_{x}}%
}=\left(  2\pi\right)  ^{-d/2}\int e^{ix\xi}\widehat{f}\left(  \xi\right)
d\xi=f\left(  x\right)  .
\]
\medskip

\textbf{Part 2.} If $\alpha\leq\kappa$ then $\widehat{D^{\alpha}R_{x}}=\left(
i\xi\right)  ^{\alpha}\widehat{R_{x}}=\left(  2\pi\right)  ^{-d/2}%
\frac{\left(  i\xi\right)  ^{\alpha}e^{-ix\xi}}{w\left(  \xi\right)  }$ so
that%
\[
\left\Vert D^{\alpha}R_{x}\right\Vert _{w,0}^{2}=\int w\left(  \xi\right)
\left\vert \left(  2\pi\right)  ^{-\frac{d}{2}}\frac{\left(  i\xi\right)
^{\alpha}e^{-ix\xi}}{w\left(  \xi\right)  }\right\vert ^{2}d\xi\leq\left(
2\pi\right)  ^{-d}\int\frac{\left\vert \xi\right\vert ^{2\left\vert
\alpha\right\vert }}{w\left(  \xi\right)  }d\xi<\infty,
\]

by weight function property W03. Further, for $f\in X_{w}^{0}$ and using the
inverse Fourier transform result \ref{1.065} for functions in $X_{w}^{0}$
gives
\begin{align*}
\left(  f,\left(  -D\right)  ^{\alpha}R_{x}\right)  _{w,0}=\int w\left(
\xi\right)  \widehat{f}\left(  \xi\right)  \overline{\widehat{\left(
-D\right)  ^{\alpha}R_{x}}}\left(  \xi\right)  d\xi &  =\int w\left(
\xi\right)  \widehat{f}\left(  \xi\right)  \left(  2\pi\right)  ^{-\frac{d}%
{2}}\frac{\left(  i\xi\right)  ^{\alpha}e^{-ix\xi}}{w\left(  \xi\right)  }%
d\xi\\
&  =\left(  2\pi\right)  ^{-\frac{d}{2}}\int e^{-ix\xi}\left(  i\xi\right)
^{\alpha}\widehat{f}\left(  \xi\right)  d\xi\\
&  =\left(  2\pi\right)  ^{-d/2}\int e^{-ix\xi}\widehat{D^{\alpha}f}\left(
\xi\right)  d\xi\\
&  =\left(  D^{\alpha}f\right)  \left(  x\right)  .
\end{align*}
\medskip

\textbf{Part 3.} Suppose for some $x\neq x^{\prime}$ we have $R_{x}%
=R_{x^{\prime}}$. Taking the Fourier transform and substituting $\widehat{G}%
=\frac{1}{w}$ we must have $e^{-i\xi x}=e^{-i\xi x^{\prime}}$ or
$e^{i\xi\left(  x-x^{\prime}\right)  }=1$ for almost all $\xi$, a
contradiction.\medskip

\textbf{Part 4}. Since we know that $G\in C_{B}^{\left(  \left\lfloor
2\kappa\right\rfloor \right)  }$ the definition of $R_{x}$ implies $R_{x}\in
C_{B}^{\left(  \left\lfloor 2\kappa\right\rfloor \right)  }$ and
\begin{equation}
\widehat{D^{\alpha}R_{x}}\left(  \xi\right)  =\left(  2\pi\right)  ^{\frac
{d}{2}}\left(  -i\right)  ^{\left\vert \alpha\right\vert }e^{-ix\xi}%
\xi^{\alpha}\widehat{G}\left(  \xi\right)  =\left(  2\pi\right)  ^{\frac{d}%
{2}}\left(  -i\right)  ^{\left\vert \alpha\right\vert }\frac{e^{-ix\xi}%
\xi^{\alpha}}{w\left(  \xi\right)  },\quad\alpha\leq\kappa.\label{1.0371}%
\end{equation}

Next note that Theorem \ref{Thm_basis_fn_properties_all_m_W3} implies that for
all $x$, and $D^{\gamma}R_{x}\in X_{w}^{0}$ when $\gamma\leq\kappa$, so the
first and last terms of \ref{1.072} make sense. Now using equation
\ref{1.0371}
\begin{align*}
\left(  D^{\alpha}R_{x},D^{\beta}R_{y}\right)  _{w,0}=\int w\widehat{D^{\alpha
}R_{x}}\overline{\widehat{D^{\beta}R_{y}}} &  =\left(  -1\right)  ^{\left\vert
\beta\right\vert }\left(  2\pi\right)  ^{-d}\int\left(  -i\xi\right)
^{\alpha+\beta}\frac{e^{i\left(  y-x\right)  \xi}}{w\left(  \xi\right)  }%
d\xi\\
&  =\left(  -1\right)  ^{\left\vert \beta\right\vert }\left(  2\pi\right)
^{-d}\int\left(  -i\xi\right)  ^{\alpha+\beta}e^{i\left(  y-x\right)  \xi
}\widehat{G}\left(  \xi\right)  d\xi\\
&  =\left(  -1\right)  ^{\left\vert \beta\right\vert }\left(  2\pi\right)
^{-d}\int e^{i\left(  y-x\right)  \xi}\widehat{D^{\alpha+\beta}G}\left(
\xi\right)  d\xi\\
&  =\left(  -1\right)  ^{\left\vert \beta\right\vert }\left(  2\pi\right)
^{-d}\left(  D^{\alpha+\beta}G\right)  \left(  y-x\right)  ,
\end{align*}

where the last step used the inverse Fourier transform rule \ref{1.024} for
basis functions. Finally, substituting \ref{1.069} gives $\left(  D^{\alpha
}R_{x},D^{\beta}R_{y}\right)  _{w,0}=\left(  -1\right)  ^{\left\vert
\beta\right\vert }\left(  D^{\alpha+\beta}R_{x}\right)  \left(  y\right)  $.

Regarding the second equation:%
\begin{align*}
\widehat{D_{x}^{\alpha}R_{x}}\left(  \xi\right)  =\widehat{D_{x}^{\alpha
}G\left(  \cdot-x\right)  }\left(  \xi\right)  =\left(  -1\right)
^{\left\vert \alpha\right\vert }\widehat{\left(  D^{\alpha}G\right)  \left(
\cdot-x\right)  }\left(  \xi\right)   &  =\left(  -1\right)  ^{\left\vert
\alpha\right\vert }e^{-ix\xi}\widehat{D^{\alpha}G}\left(  \xi\right) \\
&  =\left(  -1\right)  ^{\left\vert \alpha\right\vert }\widehat{D^{\alpha
}G\left(  \cdot-x\right)  }\left(  \xi\right) \\
&  =\left(  -1\right)  ^{\left\vert \alpha\right\vert }\widehat{D^{\alpha
}R_{x}}\left(  \xi\right)  ,
\end{align*}

so that%
\[
\left(  D^{\alpha}R_{x},D^{\beta}R_{y}\right)  _{w,0}=\left(  -1\right)
^{\left\vert \alpha+\beta\right\vert }\left(  D^{\alpha}R_{x},D^{\beta}%
R_{y}\right)  _{w,0}.
\]
\medskip

\textbf{Part 5}. If $\gamma\leq2\kappa$ then $\gamma=\alpha+\beta$ for some
$\alpha$ and $\beta$ satisfying $\alpha\leq\kappa$ and $\beta\leq\kappa$.
Applying part 3 twice we get%
\begin{align*}
D^{\gamma}R_{x}\left(  y\right)  =D^{\alpha+\beta}R_{x}\left(  y\right)
=\left(  -1\right)  ^{\left\vert \beta\right\vert }\left(  D^{\alpha}%
R_{x},D^{\beta}R_{y}\right)  _{w,0} &  =\left(  -1\right)  ^{\left\vert
\beta\right\vert }\overline{\left(  D^{\beta}R_{y},D^{\alpha}R_{x}\right)
}_{w,0}\\
&  =\left(  -1\right)  ^{\left\vert \alpha+\beta\right\vert }\overline
{D^{\alpha+\beta}R_{y}\left(  x\right)  }\\
&  =\left(  -1\right)  ^{\left\vert \gamma\right\vert }\overline{D^{\gamma
}R_{y}\left(  x\right)  }.
\end{align*}
\medskip

\textbf{Part 6} From part 2, $\left\Vert D^{\beta}f\right\Vert _{\infty}%
=\sup\limits_{x\in\mathbb{R}^{d}}\left\vert \left(  f,\left(  -D\right)
^{\beta}R_{x}\right)  _{w,0}\right\vert $ when $\beta\leq\kappa$ so that
\begin{align*}
\left\Vert D^{\beta}f\right\Vert _{\infty}\leq\sup_{x\in\mathbb{R}^{d}%
}\left\Vert D^{\beta}R_{x}\right\Vert _{w,0}\left\Vert f\right\Vert _{w,0} &
=\sup_{x\in\mathbb{R}^{d}}\sqrt{\left(  D^{\beta}R_{x},D^{\beta}R_{x}\right)
_{w,0}}\left\Vert f\right\Vert _{w,0}\\
&  =\sup_{x\in\mathbb{R}^{d}}\sqrt{\left(  -1\right)  ^{\left\vert
\beta\right\vert }\left(  D^{2\beta}R_{x}\right)  \left(  x\right)
}\left\Vert f\right\Vert _{w,0}.
\end{align*}

But from part 1, $\left(  D^{2\beta}R_{x}\right)  \left(  x\right)  =\left(
2\pi\right)  ^{-\frac{d}{2}}D^{2\beta}G(0)$ and thus%
\[
\left\Vert D^{\beta}f\right\Vert _{\infty}\leq\left(  2\pi\right)
^{-d/4}\sqrt{\left(  -1\right)  ^{\left\vert \beta\right\vert }D^{2\beta}%
G(0)}\left\Vert f\right\Vert _{w,0},
\]

which implies this part.
\end{proof}

\section{More continuity properties of the data functions: $w\in
W02$\label{Sect_loc_estim_fn_in_Xow_W2}}

In this section we assume that the weight function $w$ has property W02. The
Riesz representers of the evaluation functionals $f\rightarrow D^{\alpha
}f\left(  x\right)  $ discussed in the last section can be used to prove some
local, pointwise smoothness properties of the data functions $X_{w}^{0}$ and
their derivatives e.g. Lipschitz continuity. The basis function $G$ generated
by $w$ also lies in $X_{w}^{0}$ and is considered separately. In Chapter
\ref{Ch_Interpol} these local properties will be used to derive the order of
convergence of the variational interpolant to its data function.

\subsection{General results\label{SbSect_gen_results_W2}}

\begin{theorem}
\label{Thm_||DRx_minus_DRy||_W2}Suppose $w$ is a weight function satisfying
property W02 for parameter $\kappa$, $G$ is the basis function and $R_{x}$ is
the Riesz representer of the evaluation functional $f\rightarrow f\left(
x\right)  $ on $X_{w}^{0}$. Then for $\left\vert \alpha\right\vert \leq\kappa$%
\begin{equation}
\left\Vert D^{\alpha}\left(  R_{x}-R_{y}\right)  \right\Vert _{w,0}%
=\tfrac{\sqrt{2}}{\left(  2\pi\right)  ^{\frac{d}{4}}}\sqrt{\left(  -1\right)
^{\left\vert \alpha\right\vert }\left(  D^{2\alpha}G\left(  0\right)
-\operatorname{Re}\left(  D^{2\alpha}G\right)  \left(  y-x\right)  \right)
}.\label{1.038}%
\end{equation}

Further, if $\kappa\geq1$ and $\left\vert \alpha\right\vert \leq\kappa-1 $
then we have the bound%
\begin{equation}
\left\Vert D^{\alpha}\left(  R_{x}-R_{y}\right)  \right\Vert _{w,0}\leq
\tfrac{1}{\left(  2\pi\right)  ^{\frac{d}{2}}}\left(  \int\frac{\xi^{2\alpha
}\left\vert \xi\right\vert ^{2}}{w\left(  \xi\right)  }d\xi\right)  ^{\frac
{1}{2}}\left\vert x-y\right\vert ,\label{1.039}%
\end{equation}

where%
\begin{equation}
\int\frac{\xi^{2\alpha}\left\vert \xi\right\vert ^{2}}{w\left(  \xi\right)
}d\xi=\left(  -1\right)  ^{1+\left\vert \alpha\right\vert }\left(
2\pi\right)  ^{\frac{d}{2}}\left(  D^{2\alpha}\left\vert D\right\vert
^{2}G\right)  \left(  0\right)  .\label{1.048}%
\end{equation}

\end{theorem}

\begin{proof}
If $\left\vert \alpha\right\vert \leq\kappa$ then from the results of Theorem
\ref{Thm_ord0_Riesz_rep_W2}
\begin{align*}
& \left\Vert D^{\alpha}R_{x}-D^{\alpha}R_{y}\right\Vert _{w,0}^{2}\\
& =\left(  D^{\alpha}R_{x}-D^{\alpha}R_{y},D^{\alpha}R_{x}-D^{\alpha}%
R_{y}\right)  _{w,0}\\
& =\left(  D^{\alpha}R_{x},D^{\alpha}R_{x}\right)  _{w,0}-\left(  D^{\alpha
}R_{x},D^{\alpha}R_{y}\right)  _{w,0}-\left(  D^{\alpha}R_{y},D^{\alpha}%
R_{x}\right)  _{w,0}+\left(  D^{\alpha}R_{y},D^{\alpha}R_{y}\right)  _{w,0}\\
& =\left(  -1\right)  ^{\left\vert \alpha\right\vert }\left(  \left(
D^{2\alpha}R_{x}\right)  \left(  x\right)  -\left(  D^{2\alpha}R_{x}\right)
\left(  y\right)  -\overline{\left(  D^{2\alpha}R_{x}\right)  \left(
y\right)  }+\left(  D^{2\alpha}R_{y}\right)  \left(  y\right)  \right) \\
& =\left(  -1\right)  ^{\left\vert \alpha\right\vert }\left(  2\pi\right)
^{-\frac{d}{2}}\left(  D^{2\alpha}G\left(  0\right)  -\left(  D^{2\alpha
}G\right)  \left(  y-x\right)  -\overline{\left(  D^{2\alpha}G\right)  \left(
y-x\right)  }+D^{2\alpha}G\left(  0\right)  \right) \\
& =\left(  -1\right)  ^{\left\vert \alpha\right\vert }\left(  2\pi\right)
^{-\frac{d}{2}}\left(  D^{2\alpha}G\left(  0\right)  -\left(  D^{2\alpha
}G\right)  \left(  y-x\right)  -\overline{\left(  D^{2\alpha}G\right)  \left(
y-x\right)  }+D^{2\alpha}G\left(  0\right)  \right) \\
& =\left(  -1\right)  ^{\left\vert \alpha\right\vert }2\left(  2\pi\right)
^{-\frac{d}{2}}\left(  D^{2\alpha}G\left(  0\right)  -\operatorname{Re}\left(
D^{2\alpha}G\right)  \left(  y-x\right)  \right)  .
\end{align*}

The proof of our second result uses equation \ref{1.037} i.e.
$\widehat{D^{\alpha}R_{x}}\left(  \xi\right)  =\tfrac{\left(  -i\right)
^{\left\vert \alpha\right\vert }}{\left(  2\pi\right)  ^{\frac{d}{2}}}%
\frac{e^{-ix\xi}\xi^{\alpha}}{w\left(  \xi\right)  }$, so that
\begin{align*}
\left\Vert D^{\alpha}R_{x}-D^{\alpha}R_{y}\right\Vert _{w,0}^{2}=\int
w\left\vert \widehat{D^{\alpha}R_{x}-D^{\alpha}R_{y}}\right\vert ^{2} &
=\tfrac{1}{\left(  2\pi\right)  ^{d}}\int\left\vert e^{i\xi x}-e^{i\xi
y}\right\vert ^{2}\tfrac{\xi^{2\alpha}}{w\left(  \xi\right)  }d\xi\\
&  =\tfrac{1}{\left(  2\pi\right)  ^{d}}\int\left(  2\sin\left(
\tfrac{\left(  x-y\right)  \xi}{2}\right)  \right)  ^{2}\tfrac{\xi^{2\alpha}%
}{w\left(  \xi\right)  }d\xi\\
&  =\tfrac{1}{\left(  2\pi\right)  ^{d}}\int\left\vert \left(  x-y\right)
\xi\right\vert ^{2}\left(  \tfrac{\sin\left(  \left(  x-y\right)
\xi/2\right)  }{\left(  x-y\right)  \xi/2}\right)  ^{2}\tfrac{\xi^{2\alpha}%
}{w\left(  \xi\right)  }d\xi\\
&  \leq\tfrac{1}{\left(  2\pi\right)  ^{d}}\left(  \int\tfrac{\xi^{2\alpha
}\left\vert \xi\right\vert ^{2}}{w\left(  \xi\right)  }d\xi\right)  \left\vert
x-y\right\vert ^{2},
\end{align*}

which is finite since $1+\left\vert \alpha\right\vert \leq\kappa$. Continuing%
\[
\int\tfrac{\xi^{2\alpha}\left\vert \xi\right\vert ^{2}}{w\left(  \xi\right)
}d\xi=\int\xi^{2\alpha}\left\vert \xi\right\vert ^{2}\widehat{G}\left(
\xi\right)  d\xi=\left(  -1\right)  ^{1+\left\vert \alpha\right\vert }%
\int\widehat{D^{2\alpha}\left\vert D\right\vert ^{2}G}=\left(  2\pi\right)
^{\frac{d}{2}}\left(  -1\right)  ^{1+\left\vert \alpha\right\vert }\left(
D^{2\alpha}\left\vert D\right\vert ^{2}G\right)  \left(  0\right)  .
\]

\end{proof}

We now prove some uniform pointwise estimates for functions in $X_{w}^{0}$.

\begin{corollary}
\label{Cor_|f(x)-f(y)|_inequal_1_W2}Suppose the weight function $w$ has
property \textbf{W02} for parameter $\kappa$ and that $G$ is the basis
function of order $0$. Then we have the following local pointwise estimates
for functions $f\in X_{w}^{0}$:

\begin{enumerate}
\item If $\kappa\geq1$ and $\left\vert \alpha\right\vert \leq\kappa-1$ then%
\[
\left\vert D^{\alpha}f\left(  x\right)  -D^{\alpha}f\left(  y\right)
\right\vert \leq\tfrac{1}{\left(  2\pi\right)  ^{\frac{d}{2}}}\left\Vert
f\right\Vert _{w,0}\left(  \int\frac{\left\vert \xi\right\vert ^{2}%
\xi^{2\alpha}}{w\left(  \xi\right)  }d\xi\right)  ^{\frac{1}{2}}\left\vert
x-y\right\vert .
\]

\item If $\left\vert \alpha\right\vert \leq\kappa$ then%
\begin{equation}
\left\vert D^{\alpha}f\left(  x\right)  -D^{\alpha}f\left(  y\right)
\right\vert \leq\tfrac{\sqrt{2}}{\left(  2\pi\right)  ^{\frac{d}{4}}%
}\left\Vert f\right\Vert _{w,0}\sqrt{\left(  -1\right)  ^{\left\vert
\alpha\right\vert }\left(  D^{2\alpha}G\left(  0\right)  -\operatorname{Re}%
\left(  D^{2\alpha}G\right)  \left(  y-x\right)  \right)  }.\label{1.15}%
\end{equation}

\item For the hat weight function in one dimension $\kappa<1/2$ and
\[
\left\vert f\left(  x\right)  -f\left(  y\right)  \right\vert \leq\tfrac
{\sqrt{2}}{\left(  2\pi\right)  ^{\frac{d}{4}}}\left\Vert f\right\Vert
_{w,0}\left\vert x-y\right\vert ^{\frac{1}{2}},\text{\quad}x,y\in
\mathbb{R}^{1},\text{ }\left\vert x-y\right\vert \leq1.
\]

\end{enumerate}
\end{corollary}

\begin{proof}
\textbf{Part 1} This is an application of inequality \ref{1.039} of Theorem
\ref{Thm_||DRx_minus_DRy||_W2}
\begin{align*}
\left\vert D^{\alpha}f\left(  x\right)  -D^{\alpha}f\left(  y\right)
\right\vert =\left\vert \left(  f,D^{\alpha}R_{x}-D^{\alpha}R_{y}\right)
_{w,0}\right\vert  &  \leq\left\Vert f\right\Vert _{w,0}\left\Vert D^{\alpha
}R_{x}-D^{\alpha}R_{y}\right\Vert _{w,0}\\
&  \leq\tfrac{1}{\left(  2\pi\right)  ^{\frac{d}{2}}}\left\Vert f\right\Vert
_{w,0}\left(  \int\frac{\left\vert \xi\right\vert ^{2}\xi^{2\alpha}}{w\left(
\xi\right)  }d\xi\right)  ^{\frac{1}{2}}\left\vert x-y\right\vert .
\end{align*}
\medskip

\textbf{Part 2} By part 2 of Theorem \ref{Thm_ord0_Riesz_rep_W2}%
\[
\left\vert D^{\alpha}f\left(  x\right)  -D^{\alpha}f\left(  y\right)
\right\vert =\left\vert \left(  f,D^{\alpha}R_{x}-D^{\alpha}R_{y}\right)
_{w,0}\right\vert \leq\left\Vert f\right\Vert _{w,0}\left\Vert D^{\alpha}%
R_{x}-D^{\alpha}R_{y}\right\Vert _{w,0},
\]

and inequality \ref{1.038} of Theorem \ref{Thm_||DRx_minus_DRy||_W2} completes
the proof.\medskip

\textbf{Part 3} From Theorem \ref{Thm_hat_wt_extend_props} $\kappa<1/2$ so
$\alpha=0$ and since $G\left(  x\right)  =1-\left\vert x\right\vert $ when
$\left\vert x\right\vert \leq1$ our result follows from the bound proved in
part 2.
\end{proof}

The next corollary provides information concerning the smoothness properties
of basis functions.

\begin{corollary}
\label{Cor_bound_DG(z-x)-DG(z-y)_W2}Suppose $G$ is the basis function
generated by a weight function $w$ satisfying property \textbf{W02} for
parameter $\kappa$. Then for each $z\in\mathbb{R}^{d}$, $G$ has the following properties:

\begin{enumerate}
\item If $\left\vert \alpha\right\vert \leq\kappa$ and $\left\vert
\beta\right\vert \leq\kappa$ then%
\begin{equation}
\left\vert D^{\alpha+\beta}G\left(  x\right)  -D^{\alpha+\beta}G\left(
y\right)  \right\vert \leq k\sqrt{\left(  -1\right)  ^{\left\vert
\alpha\right\vert }\left(  D^{2\alpha}G\left(  0\right)  -\operatorname{Re}%
\left(  D^{2\alpha}G\right)  \left(  x-y\right)  \right)  },\label{1.010}%
\end{equation}

where
\begin{equation}
k=\tfrac{\sqrt{2}}{\left(  2\pi\right)  ^{\frac{d}{4}}}\left(  \int\tfrac
{\xi^{2\alpha}d\xi}{w\left(  \xi\right)  }\right)  ^{\frac{1}{2}}=\sqrt
{2}\sqrt{\left(  -1\right)  ^{\left\vert \alpha\right\vert }D^{2\alpha
}G\left(  0\right)  }.\label{1.049}%
\end{equation}

\item If $\left\vert \alpha\right\vert \leq\kappa$ then%
\[
\left\vert D^{2\alpha}G\left(  x\right)  \right\vert \leq\left(  -1\right)
^{\left\vert \alpha\right\vert }D^{2\alpha}G\left(  0\right)  ,
\]

and%
\[
\left(  -1\right)  ^{\left\vert \alpha\right\vert }\operatorname{Re}%
D^{2\alpha}G\left(  x\right)  <\left(  -1\right)  ^{\left\vert \alpha
\right\vert }D^{2\alpha}G\left(  0\right)  ,\quad x\neq0.
\]

\item Suppose $\kappa\geq1$. Then if $\left\vert \alpha\right\vert \leq
\kappa-1$ and $\left\vert \beta\right\vert \leq\kappa-1$
\[
\left\vert D^{\alpha+\beta}G\left(  x\right)  -D^{\alpha+\beta}G\left(
y\right)  \right\vert \leq\tfrac{1}{\left(  2\pi\right)  ^{\frac{d}{2}}%
}\left(  \int\frac{\xi^{2\beta}d\xi}{w\left(  \xi\right)  }\right)  ^{\frac
{1}{2}}\left(  \int\frac{\left\vert \xi\right\vert ^{2}\xi^{2\alpha}d\xi
}{w\left(  \xi\right)  }\right)  ^{\frac{1}{2}}\left\vert x-y\right\vert .
\]

\end{enumerate}
\end{corollary}

\begin{proof}
\textbf{Part 1.} Our starting point is equation \ref{1.6} i.e. for
$x,y\in\mathbb{R}^{d}$,

$\left(  D^{\alpha}R_{x},D^{\beta}R_{y}\right)  _{w,0}=\left(  -1\right)
^{\left\vert \beta\right\vert }\left(  D^{\alpha+\beta}R_{x}\right)  \left(
y\right)  $ when $\left\vert \alpha\right\vert \leq\kappa$, $\left\vert
\beta\right\vert \leq\kappa$. Hence, since

$R_{x}\left(  y\right)  =\left(  2\pi\right)  ^{-\frac{d}{2}}G\left(
y-x\right)  $%
\begin{align*}
D_{z}^{\alpha+\beta}\left(  G\left(  z-x\right)  -G\left(  z-y\right)
\right)   & =\left(  2\pi\right)  ^{-\frac{d}{2}}D^{\alpha+\beta}\left(
R_{x}-R_{y}\right)  \left(  z\right) \\
& =\left(  2\pi\right)  ^{-\frac{d}{2}}\left(  D^{\alpha}\left(  R_{x}%
-R_{y}\right)  ,D^{\beta}R_{z}\right)  _{w,0},
\end{align*}

so that%
\begin{equation}
\left\vert D_{z}^{\alpha+\beta}\left(  G\left(  z-x\right)  -G\left(
z-y\right)  \right)  \right\vert \leq\left(  2\pi\right)  ^{-\frac{d}{2}%
}\left\Vert D^{\alpha}\left(  R_{x}-R_{y}\right)  \right\Vert _{w,0}\left\Vert
D^{\beta}R_{z}\right\Vert _{w,0}.\label{1.041}%
\end{equation}

From equation \ref{1.0371}: $\widehat{D^{\beta}R_{x}}\left(  \xi\right)
=\left(  2\pi\right)  ^{\frac{d}{2}}\left(  -i\right)  ^{\left\vert
\beta\right\vert }\frac{e^{-ix\xi}\xi^{\beta}}{w\left(  \xi\right)  }$ for
$\beta\leq\kappa$, so we can conclude that
\[
\left\Vert D^{\beta}R_{z}\right\Vert _{w,0}=\left(  \int w\left\vert
\widehat{D^{\beta}R_{z}}\right\vert ^{2}\right)  ^{\frac{1}{2}}=\left(
2\pi\right)  ^{-\frac{d}{2}}\left(  \int\frac{\xi^{2\beta}d\xi}{w\left(
\xi\right)  }\right)  ^{\frac{1}{2}},
\]
and by using equation \ref{1.038} for $\left\Vert D^{\alpha}\left(
R_{x}-R_{y}\right)  \right\Vert _{w,0}$ inequality \ref{1.010}
follows.\medskip

\textbf{Part 2.} We have $\left(  -1\right)  ^{\left\vert \alpha\right\vert
}D^{2\alpha}G\left(  x\right)  =\left(  2\pi\right)  ^{-\frac{d}{2}}\int
e^{i\xi x}\tfrac{\xi^{2\alpha}}{w\left(  \xi\right)  }d\xi$, so that
\[
\left\vert D^{2\alpha}G\left(  x\right)  \right\vert \leq\left(  2\pi\right)
^{-\frac{d}{2}}\int\frac{\xi^{2\alpha}}{w\left(  \xi\right)  }d\xi=\left(
-1\right)  ^{\left\vert \alpha\right\vert }D^{2\alpha}G\left(  0\right)  .
\]

From part 1 we know that $\left(  -1\right)  ^{\left\vert \alpha\right\vert
}\operatorname{Re}\left(  D^{2\alpha}G\right)  \left(  x-y\right)  \leq\left(
-1\right)  ^{\left\vert \alpha\right\vert }D^{2\alpha}G\left(  0\right)  $.
Now suppose that

$\left(  -1\right)  ^{\left\vert \alpha\right\vert }\operatorname{Re}\left(
D^{2\alpha}G\right)  \left(  z\right)  =\left(  -1\right)  ^{\left\vert
\alpha\right\vert }D^{2\alpha}G\left(  0\right)  $ for some $z\neq0$. Equation
\ref{1.038} would then imply that $D^{\alpha}R_{x}=D^{\alpha}R_{z+x}$ for all
$x$, and hence $z=0$ by part 3 of Theorem \ref{Thm_ord0_Riesz_rep_W2}.\medskip

\textbf{Part 3.} Our starting point is inequality \ref{1.041} of part 1.
Substitute the expression for $\left\Vert D^{\beta}R_{z}\right\Vert _{w,0}$
given in the proof of part 1 and then use inequality \ref{1.039} to estimate
$\left\Vert D^{\alpha}\left(  R_{x}-R_{y}\right)  \right\Vert _{w,0}$.
\end{proof}

\begin{remark}
\label{Rem_InTermsOfDeven(0)_W2}Observe that by using \ref{1.048} and
\ref{1.049} the right sides of all the estimates involving weight functions in
this subsection can be written in terms of even order derivatives of the basis
function evaluated at the origin.
\end{remark}

\subsection{Better results obtained using Taylor
series\label{SbSect_better_results_W2}}

We can improve on part 3 of Corollary \ref{Cor_bound_DG(z-x)-DG(z-y)_W2} for
the special case $\kappa\geq1$, $y=0$ and $\alpha=\beta=0$; improved in the
sense that $\left\vert x\right\vert $ is replaced by $\left\vert x\right\vert
^{2}$ and this is done by using the Taylor series expansion with integral
remainder given in Subsection \ref{Sect_apx_TaylorSeries} of the Appendix:

\begin{theorem}
\label{Thm_G(0)minusG(x)_bound_w_W2}Suppose the weight function $w$ satisfies
property \textbf{W02} for $\kappa=1$. Then $\left(  \left\vert D\right\vert
^{2}G\right)  \left(  0\right)  $ is real and negative and we have the bound%
\begin{equation}
G\left(  0\right)  -\operatorname{Re}G\left(  x\right)  \leq-\frac{d}%
{2}\left(  \left\vert D\right\vert ^{2}G\right)  \left(  0\right)  \left\vert
x\right\vert ^{2},\quad x\in\mathbb{R}^{d}.\label{1.060}%
\end{equation}

Now suppose $w$ is radial and let $r=\left\vert x\right\vert $. Then:

\begin{enumerate}
\item If $G\left(  x\right)  =g\left(  r\right)  $ then $g\in C_{B}^{\left(
2\right)  }\left(  \left[  0,\infty\right)  \right)  $ and $r^{-1}g^{\prime
}\in C_{B}^{\left(  0\right)  }\left(  \left[  0,\infty\right)  \right)  $.
Also%
\begin{equation}
\left(  \left\vert D\right\vert ^{2}G\right)  \left(  0\right)  =g^{\prime
\prime}\left(  0\right)  d.\label{1.77}%
\end{equation}

\item If $G\left(  x\right)  =f\left(  r^{2}\right)  $ then $f\in C^{\left(
2\right)  }\left(  \left(  0,\infty\right)  \right)  \cap C_{B}^{\left(
1\right)  }\left(  \left[  0,\infty\right)  \right)  $ and $rf^{\prime\prime
}\left(  r\right)  \in C_{B}^{\left(  0\right)  }\left(  \left[
0,\infty\right)  \right)  $ and

$\lim\limits_{r\rightarrow0^{+}}rf^{\prime\prime}\left(  r\right)  =0$. Also%
\begin{equation}
\left(  \left\vert D\right\vert ^{2}G\right)  \left(  0\right)  =2f^{\prime
}\left(  0\right)  d.\label{1.78}%
\end{equation}

\end{enumerate}
\end{theorem}

\begin{proof}
We start by assuming $w$ is an even function. Because $\kappa=1$ it follows
that $G\in C_{B}^{\left(  \left\lfloor 2\kappa\right\rfloor \right)  }\subset
C_{B}^{\left(  2\right)  }$.

Also, since
\[
D^{\beta}G\left(  x\right)  =\left(  2\pi\right)  ^{-\frac{d}{2}}\int
e^{-ix\xi}\frac{\left(  i\xi\right)  ^{\beta}}{w\left(  \xi\right)  }%
d\xi,\quad\left\vert \beta\right\vert \leq2,
\]
\qquad

and $w$ is even, it follows that $G$ is real, $\left\vert G\left(  x\right)
\right\vert \leq G\left(  0\right)  $ and $D_{k}G\left(  0\right)  =0$ for all
$k $. Thus from Appendix \ref{Sect_apx_TaylorSeries} we now have the second
order Taylor series expansion%
\[
G\left(  x\right)  =G\left(  0\right)  +\left(  \mathcal{R}_{2}G\right)
\left(  0,x\right)  ,
\]
with remainder estimate%
\[
\left\vert \left(  \mathcal{R}_{2}G\right)  \left(  0,x\right)  \right\vert
\leq\frac{d}{2}\left(  \max_{\substack{\left\vert \beta\right\vert =2
\\t\in\left[  0,x\right]  }}\left\vert D^{\beta}G\left(  t\right)  \right\vert
\right)  \left\vert x\right\vert ^{2},\quad x\in\mathbb{R}^{d}.
\]

But when $\left\vert \beta\right\vert =2$, $D^{\beta}G\left(  x\right)
=-\left(  2\pi\right)  ^{-\frac{d}{2}}\int e^{-ix\xi}\frac{\xi^{\beta}%
}{w\left(  \xi\right)  }d\xi$ and hence
\[
\left\vert D^{\beta}G\left(  x\right)  \right\vert \leq\left(  2\pi\right)
^{-\frac{d}{2}}\int\frac{\left\vert \xi\right\vert ^{\left\vert \beta
\right\vert }}{w\left(  \xi\right)  }d\xi=\left(  2\pi\right)  ^{-\frac{d}{2}%
}\int\frac{\left\vert \xi\right\vert ^{2}}{w\left(  \xi\right)  }d\xi=-\left(
\left\vert D\right\vert ^{2}G\right)  \left(  0\right)  ,
\]

so that%
\begin{equation}
G\left(  0\right)  -G\left(  x\right)  \leq-\frac{d}{2}\left(  \left\vert
D\right\vert ^{2}G\right)  \left(  0\right)  \left\vert x\right\vert
^{2},\quad x\in\mathbb{R}^{d}.\label{1.009}%
\end{equation}

Now set aside the assumption that $w$ is even. We define the even function
$w_{e}$ by%
\begin{equation}
\frac{1}{w_{e}\left(  \xi\right)  }=\frac{1}{2}\left(  \frac{1}{w\left(
\xi\right)  }+\frac{1}{w\left(  -\xi\right)  }\right)  ,\label{1.66}%
\end{equation}

and it is easy to show that $w_{e}$ has property W01 and satisfies%
\begin{equation}
\int\frac{\left\vert \xi\right\vert ^{\lambda}}{w_{e}\left(  \xi\right)  }%
d\xi=\int\frac{\left\vert \xi\right\vert ^{\lambda}}{w\left(  \xi\right)
}d\xi,\quad0\leq\lambda\leq\kappa,\label{1.012}%
\end{equation}

so that $w_{e}$ has property W02 for $\kappa$. Further, if $w_{e}$ has basis
function $G_{e}$ then $\widehat{G_{e}}=\frac{1}{2}\left(  \widehat{G}%
+\widehat{\overline{G}}\right)  $ and hence $G_{e}=\operatorname{Re}G$. But
$G_{e}$ satisfies \ref{1.009} and since \ref{1.012} implies $\left(
\left\vert D\right\vert ^{2}G_{e}\right)  \left(  0\right)  =\left(
\left\vert D\right\vert ^{2}G\right)  \left(  0\right)  $ we can conclude that
\ref{1.060} is true.\medskip

\fbox{Part 1} Suppose $G\left(  x\right)  =g\left(  r\right)  $. Then
substituting $x=\left(  r,0^{\prime}\right)  $ we have
\[
g\left(  r\right)  =G\left(  r,0^{\prime}\right)  ,
\]

so that it is clear that $g\in C_{B}^{\left(  2\right)  }\left(  \left[
0,\infty\right)  \right)  $ since $G\in C_{B}^{\left(  2\right)  }$. Regarding
$r^{-1}g^{\prime}$, it is clear that $r^{-1}g^{\prime}\in C_{B}^{\left(
2\right)  }\left(  \left(  0,\infty\right)  \right)  $. Further, since
$g^{\prime}\in C_{B}^{\left(  1\right)  }\left(  \left[  0,\infty\right)
\right)  $, $\lim\limits_{r\rightarrow0^{+}}\frac{g^{\prime}\left(  r\right)
}{r}=g^{\prime\prime}\left(  0\right)  $ and $\lim\limits_{r\rightarrow\infty
}\frac{g^{\prime}\left(  r\right)  }{r}=0$ which means that $r^{-1}g^{\prime
}\in C_{B}^{\left(  2\right)  }\left(  \left[  0,\infty\right)  \right)  $.

If $x\neq0$
\[
D_{k}G\left(  x\right)  =\frac{x_{k}}{r}g^{\prime}\left(  r\right)  ,
\]

and
\[
D_{j}D_{k}G\left(  x\right)  =\frac{x_{j}x_{k}}{r^{2}}g^{\prime\prime}\left(
r\right)  -\frac{x_{j}x_{k}}{r^{3}}g^{\prime}\left(  r\right)  +\frac
{\delta_{j,k}}{r}g^{\prime}\left(  r\right)  ,
\]

so that%
\[
\left\vert D\right\vert ^{2}G\left(  x\right)  =\sum_{k=1}^{d}D_{k}%
^{2}G\left(  x\right)  =g^{\prime\prime}\left(  r\right)  -\frac{g^{\prime
}\left(  r\right)  }{r}+d\frac{g^{\prime}\left(  r\right)  }{r}=g^{\prime
\prime}\left(  r\right)  +\left(  d-1\right)  \frac{g^{\prime}\left(
r\right)  }{r}.
\]

But $g\in C_{B}^{\left(  2\right)  }\left(  \left[  0,\infty\right)  \right)
$ so
\[
\lim\limits_{x\rightarrow0}\left\vert D\right\vert ^{2}G\left(  x\right)
=g^{\prime\prime}\left(  0\right)  +\left(  d-1\right)  g^{\prime\prime
}\left(  0\right)  =g^{\prime\prime}\left(  0\right)  d.
\]
\medskip

\fbox{Part 2} Suppose $G\left(  x\right)  =f\left(  r^{2}\right)  $. Then the
equation $f\left(  r\right)  =g\left(  \sqrt{r}\right)  $ and the properties
of $f$ proved in part 1 can be used to easily derive the results of this part.
\end{proof}

For radial basis functions the estimates the last theorem for $G\left(
0\right)  -\operatorname{Re}G\left(  x\right)  $ give a $d^{2}$ dependency.
However, if the basis function is a radial function this dependency can be avoided:

\begin{theorem}
\label{Thm_G(0)minusG(x)_bnd_w_rad_W2}Suppose the weight function $w$
satisfies property W02 for $\kappa=1$ and denote the basis function by $G$. Then:

\begin{enumerate}
\item If $G\left(  x\right)  =f\left(  r^{2}\right)  $ then%
\begin{equation}
G\left(  0\right)  -\operatorname{Re}G\left(  x\right)  \leq\left\Vert
2rf^{\prime\prime}+f^{\prime}\right\Vert _{\infty}\left\vert x\right\vert
^{2},\quad x\in\mathbb{R}^{d}.\label{1.67}%
\end{equation}

\item If $G\left(  x\right)  =g\left(  r\right)  $ then%
\begin{equation}
G\left(  0\right)  -\operatorname{Re}G\left(  x\right)  \leq\frac{1}%
{2}\left\Vert g^{\prime\prime}\right\Vert _{\infty}\left\vert x\right\vert
^{2},\quad x\in\mathbb{R}^{d},\label{1.68}%
\end{equation}

where $\left\Vert 2rf^{\prime\prime}+f^{\prime}\right\Vert _{\infty}=\frac
{1}{2}\left\Vert g^{\prime\prime}\right\Vert _{\infty}$.
\end{enumerate}
\end{theorem}

\begin{proof}
We will start by assuming $w$ is an even function. Because $\kappa\geq1$ it
follows that $G\in C_{B}^{\left(  \left\lfloor 2\kappa\right\rfloor \right)
}\subset C_{B}^{\left(  2\right)  }$. Also, since
\[
D^{\beta}G\left(  x\right)  =\left(  2\pi\right)  ^{-\frac{d}{2}}\int
e^{-ix\xi}\frac{\left(  i\xi\right)  ^{\beta}}{w\left(  \xi\right)  }%
d\xi,\quad\left\vert \beta\right\vert \leq2,
\]

it follows that if $w$ is even then $G$ is real, $\left\vert G\left(
x\right)  \right\vert \leq G\left(  0\right)  $ and $D_{k}G\left(  0\right)
=0$ for all $k$. Thus from Appendix \ref{Sect_apx_TaylorSeries} we now have
the Taylor series expansion%
\[
G\left(  x\right)  =G\left(  0\right)  +\left(  \mathcal{R}_{2}G\right)
\left(  0,x\right)  ,
\]

where%
\begin{equation}
\left(  \mathcal{R}_{2}G\right)  \left(  0,x\right)  =2\sum_{\left\vert
\beta\right\vert =2}\frac{x^{\beta}}{\beta!}\int_{0}^{1}\left(  1-t\right)
\left(  D^{\beta}G\right)  \left(  tx\right)  dt.\label{1.54}%
\end{equation}
\medskip

\fbox{Part 1} Suppose $G\left(  x\right)  =f\left(  r^{2}\right)  $. Then
\[
D_{k}G\left(  x\right)  =2x_{k}f^{\prime}\left(  r^{2}\right)  ,
\]

and
\[
D_{j}D_{k}G\left(  x\right)  =4x_{j}x_{k}f^{\prime\prime}\left(  r^{2}\right)
+2\delta_{j,k}f^{\prime}\left(  r^{2}\right)  .
\]

Set $\delta_{\beta}=\delta_{2,\max\beta}$. Then in multi-index notation
\[
D^{\beta}G\left(  x\right)  =4x^{\beta}f^{\prime\prime}\left(  r^{2}\right)
+2\delta_{\beta}f^{\prime}\left(  r^{2}\right)  ,\quad\left\vert
\beta\right\vert =2,
\]

and so, using the identity \ref{1.57}, the remainder can be written%
\begin{align}
\left(  \mathcal{R}_{2}G\right)  \left(  0,x\right)   & =2\sum_{\left\vert
\beta\right\vert =2}\frac{x^{\beta}}{\beta!}\int_{0}^{1}\left(  1-t\right)
\left(  4x^{\beta}f^{\prime\prime}\left(  r^{2}\right)  \right)  \left(
tx\right)  dt+\nonumber\\
& \qquad+2\sum_{\left\vert \beta\right\vert =2}\frac{x^{\beta}}{\beta!}%
\int_{0}^{1}\left(  1-t\right)  \left(  2\delta_{\beta}f^{\prime}\left(
r^{2}\right)  \right)  \left(  tx\right)  dt\nonumber\\
& =8\sum_{\left\vert \beta\right\vert =2}\frac{x^{\beta}}{\beta!}\int_{0}%
^{1}\left(  1-t\right)  \left(  t^{2}x^{\beta}\right)  f^{\prime\prime}\left(
t^{2}r^{2}\right)  dt+\nonumber\\
& \qquad+4\sum_{\left\vert \beta\right\vert =2}\frac{x^{\beta}}{\beta!}%
\int_{0}^{1}\left(  1-t\right)  \delta_{\beta}f^{\prime}\left(  t^{2}%
r^{2}\right)  dt\nonumber\\
& =8\left(  \sum_{\left\vert \beta\right\vert =2}\frac{x^{2\beta}}{\beta
!}\right)  \int_{0}^{1}\left(  1-t\right)  t^{2}f^{\prime\prime}\left(
t^{2}r^{2}\right)  dt+\nonumber\\
& \qquad+4\left(  \sum_{\left\vert \beta\right\vert =2}\delta_{\beta}%
\frac{x^{\beta}}{\beta!}\right)  \int_{0}^{1}\left(  1-t\right)  f^{\prime
}\left(  t^{2}r^{2}\right)  dt\nonumber\\
& =4r^{4}\int_{0}^{1}\left(  1-t\right)  t^{2}f^{\prime\prime}\left(
t^{2}r^{2}\right)  dt+2r^{2}\int_{0}^{1}\left(  1-t\right)  f^{\prime}\left(
t^{2}r^{2}\right)  dt\nonumber\\
& =r^{2}\int_{0}^{1}\left(  1-t\right)  \left(  4t^{2}r^{2}f^{\prime\prime
}\left(  t^{2}r^{2}\right)  +2f^{\prime}\left(  t^{2}r^{2}\right)  \right)
dt.\label{1.71}%
\end{align}

The remainder can be estimated as%
\begin{align}
\left\vert \left(  \mathcal{R}_{2}G\right)  \left(  0,x\right)  \right\vert
&  \leq r^{2}\int_{0}^{1}\left(  1-t\right)  \left\vert 4t^{2}r^{2}%
f^{\prime\prime}\left(  t^{2}r^{2}\right)  +2f^{\prime}\left(  t^{2}%
r^{2}\right)  \right\vert dt\nonumber\\
&  =r^{2}\left\Vert 4rf^{\prime\prime}+2f^{\prime}\right\Vert _{\infty}%
\int_{0}^{1}\left(  1-t\right)  dt\nonumber\\
&  =\left\Vert 2rf^{\prime\prime}+f^{\prime}\right\Vert _{\infty}\left\vert
x\right\vert ^{2}.\label{1.73}%
\end{align}

so that we have proved the desired estimate \ref{1.67} when $w$ is an even
function. This estimate can now be extended to an arbitrary weight function by
the technique used in Theorem \ref{Thm_G(0)minusG(x)_bound_w_W2} which
involved defining the even weight function \ref{1.66}.\medskip

\fbox{Part 2} Suppose $G\left(  x\right)  =g\left(  r\right)  $. Then
$g\left(  r\right)  =f\left(  r^{2}\right)  $ and
\begin{equation}
g^{\prime\prime}\left(  r\right)  =2f^{\prime}\left(  r^{2}\right)
+4r^{2}f^{\prime\prime}\left(  r^{2}\right)  ,\label{1.72}%
\end{equation}

so the remainder equation \ref{1.71} can be written%
\[
\left(  \mathcal{R}_{2}G\right)  \left(  0,x\right)  =r^{2}\int_{0}^{1}\left(
1-t\right)  \left(  4t^{2}r^{2}f^{\prime\prime}\left(  t^{2}r^{2}\right)
+2f^{\prime}\left(  t^{2}r^{2}\right)  \right)  dt=\frac{1}{2}r^{2}\int%
_{0}^{1}\left(  1-t\right)  g^{\prime\prime}\left(  tr\right)  dt,
\]

and estimated by%
\begin{equation}
\left\vert \left(  \mathcal{R}_{2}G\right)  \left(  0,x\right)  \right\vert
\leq\frac{1}{2}\left\Vert g^{\prime\prime}\right\Vert _{\infty}\left\vert
x\right\vert ^{2}.\label{1.74}%
\end{equation}

Equation \ref{1.72} ensures that the estimates \ref{1.73} and \ref{1.74} are equal.

We have now proved the desired estimate \ref{1.68} when $w$ is an even
function. This estimate can now be extended to an arbitrary weight function by
the technique of Theorem \ref{Thm_G(0)minusG(x)_bound_w_W2} which involved
defining the even weight function \ref{1.66}.
\end{proof}

\begin{remark}
Calculations using part 1 may be easier when the basis function depends on
$r^{2}$: see the Type 1 convergence estimates for the radial functions of
Subsubsection \ref{SbSect_int_examp_rad_basis_ord_converg}.
\end{remark}

\section{More continuity properties of the data functions: $w\in
W03$\label{Sect_loc_estim_fn_in_Xow_W3}}

In this section we assume that the weight function $w$ has property W03. The
Riesz representers of the evaluation functionals $f\rightarrow D^{\alpha
}f\left(  x\right)  $ discussed in the last section can be used to prove some
local, pointwise smoothness properties of the data functions $X_{w}^{0}$ and
their derivatives e.g. Lipschitz continuity. The basis function $G$ generated
by $w$ also lies in $X_{w}^{0}$ and is considered separately.

\subsection{General results\label{SbSect_gen_results_W3}}

In Remark \ref{Rem_InTermsOfDeven(0)_W3} we observe that the right sides of
all the estimates derived in this subsection can be written in terms of even
order derivatives of the basis function evaluated at the origin.

\begin{theorem}
\label{Thm_||DRx_minus_DRy||_W3}Suppose $w$ is a weight function satisfying
property W03 for parameter $\kappa\in\mathbb{R}^{d}$, $G$ is the basis
function and $R_{x}$ is the Riesz representer of the evaluation functional
$f\rightarrow f\left(  x\right)  $ on $X_{w}^{0}$. Then for $\alpha\leq\kappa
$,%
\begin{equation}
\left\Vert D^{\alpha}\left(  R_{x}-R_{y}\right)  \right\Vert _{w,0}%
=\tfrac{\sqrt{2}}{\left(  2\pi\right)  ^{\frac{d}{4}}}\sqrt{\left(  -1\right)
^{\left\vert \alpha\right\vert }\left(  D^{2\alpha}G\left(  0\right)
-\operatorname{Re}\left(  D^{2\alpha}G\right)  \left(  y-x\right)  \right)
}.\label{1.0381}%
\end{equation}

Further, if $\kappa\geq\mathbf{1}$ and $\alpha\leq\kappa-1$ then we have the
bound%
\begin{equation}
\left\Vert D^{\alpha}\left(  R_{x}-R_{y}\right)  \right\Vert _{w,0}\leq\left(
2\pi\right)  ^{-\frac{d}{2}}\left(  \int\frac{\xi^{2\alpha}\left\vert
\xi\right\vert ^{2}}{w\left(  \xi\right)  }d\xi\right)  ^{\frac{1}{2}%
}\left\vert x-y\right\vert ,\label{1.0391}%
\end{equation}

where%
\begin{equation}
\int\frac{\xi^{2\alpha}\left\vert \xi\right\vert ^{2}}{w\left(  \xi\right)
}d\xi=\left(  -1\right)  ^{1+\left\vert \alpha\right\vert }\left(
2\pi\right)  ^{\frac{d}{2}}\left(  D^{2\alpha}\left\vert D\right\vert
^{2}G\right)  \left(  0\right)  .\label{1.0481}%
\end{equation}

\end{theorem}

\begin{proof}
If $\alpha\leq\kappa$ then from the results of Theorem
\ref{Thm_ord0_Riesz_rep_W3}
\begin{align*}
& \left\Vert D^{\alpha}R_{x}-D^{\alpha}R_{y}\right\Vert _{w,0}^{2}\\
& =\left(  D^{\alpha}R_{x}-D^{\alpha}R_{y},D^{\alpha}R_{x}-D^{\alpha}%
R_{y}\right)  _{w,0}\\
& =\left(  D^{\alpha}R_{x},D^{\alpha}R_{x}\right)  _{w,0}-\left(  D^{\alpha
}R_{x},D^{\alpha}R_{y}\right)  _{w,0}-\left(  D^{\alpha}R_{y},D^{\alpha}%
R_{x}\right)  _{w,0}+\left(  D^{\alpha}R_{y},D^{\alpha}R_{y}\right)  _{w,0}\\
& =\left(  -1\right)  ^{\left\vert \alpha\right\vert }\left(  \left(
D^{2\alpha}R_{x}\right)  \left(  x\right)  -\left(  D^{2\alpha}R_{x}\right)
\left(  y\right)  -\overline{\left(  D^{2\alpha}R_{x}\right)  \left(
y\right)  }+\left(  D^{2\alpha}R_{y}\right)  \left(  y\right)  \right) \\
& =\left(  -1\right)  ^{\left\vert \alpha\right\vert }\left(  2\pi\right)
^{-\frac{d}{2}}\left(  D^{2\alpha}G\left(  0\right)  -\left(  D^{2\alpha
}G\right)  \left(  y-x\right)  -\overline{\left(  D^{2\alpha}G\right)  \left(
y-x\right)  }+D^{2\alpha}G\left(  0\right)  \right) \\
& =\left(  -1\right)  ^{\left\vert \alpha\right\vert }\left(  2\pi\right)
^{-\frac{d}{2}}\left(  D^{2\alpha}G\left(  0\right)  -\left(  D^{2\alpha
}G\right)  \left(  y-x\right)  -\overline{\left(  D^{2\alpha}G\right)  \left(
y-x\right)  }+D^{2\alpha}G\left(  0\right)  \right) \\
& =\left(  -1\right)  ^{\left\vert \alpha\right\vert }2\left(  2\pi\right)
^{-\frac{d}{2}}\left(  D^{2\alpha}G\left(  0\right)  -\operatorname{Re}\left(
D^{2\alpha}G\right)  \left(  y-x\right)  \right)  .
\end{align*}

The proof of our second result uses equation \ref{1.037} i.e.
$\widehat{D^{\alpha}R_{x}}\left(  \xi\right)  =\tfrac{\left(  -i\right)
^{\left\vert \alpha\right\vert }}{\left(  2\pi\right)  ^{\frac{d}{2}}}%
\frac{e^{-ix\xi}\xi^{\alpha}}{w\left(  \xi\right)  }$, $\alpha\leq\kappa$, so
that
\begin{align*}
\left\Vert D^{\alpha}R_{x}-D^{\alpha}R_{y}\right\Vert _{w,0}^{2}=\int
w\left\vert \widehat{D^{\alpha}R_{x}-D^{\alpha}R_{y}}\right\vert ^{2} &
=\tfrac{1}{\left(  2\pi\right)  ^{d}}\int\left\vert e^{i\xi x}-e^{i\xi
y}\right\vert ^{2}\tfrac{\xi^{2\alpha}}{w\left(  \xi\right)  }d\xi\\
&  =\tfrac{1}{\left(  2\pi\right)  ^{d}}\int\left(  2\sin\left(
\tfrac{\left(  x-y\right)  \xi}{2}\right)  \right)  ^{2}\tfrac{\xi^{2\alpha}%
}{w\left(  \xi\right)  }d\xi\\
&  =\tfrac{1}{\left(  2\pi\right)  ^{d}}\int\left\vert \left(  x-y\right)
\xi\right\vert ^{2}\left(  \tfrac{\sin\left(  \left(  x-y\right)
\xi/2\right)  }{\left(  x-y\right)  \xi/2}\right)  ^{2}\tfrac{\xi^{2\alpha}%
}{w\left(  \xi\right)  }d\xi\\
&  \leq\tfrac{1}{\left(  2\pi\right)  ^{d}}\left(  \int\tfrac{\xi^{2\alpha
}\left\vert \xi\right\vert ^{2}}{w\left(  \xi\right)  }d\xi\right)  \left\vert
x-y\right\vert ^{2},
\end{align*}

which is finite since $1+\alpha\leq\kappa$. Continuing%
\[
\int\tfrac{\xi^{2\alpha}\left\vert \xi\right\vert ^{2}}{w\left(  \xi\right)
}d\xi=\int\xi^{2\alpha}\left\vert \xi\right\vert ^{2}\widehat{G}\left(
\xi\right)  d\xi=\left(  -1\right)  ^{1+\left\vert \alpha\right\vert }%
\int\widehat{D^{2\alpha}\left\vert D\right\vert ^{2}G}=\left(  2\pi\right)
^{\frac{d}{2}}\left(  -1\right)  ^{1+\left\vert \alpha\right\vert }\left(
D^{2\alpha}\left\vert D\right\vert ^{2}G\right)  \left(  0\right)  .
\]

\end{proof}

We now prove some uniform pointwise estimates for functions in $X_{w}^{0}$.

\begin{corollary}
\label{Cor_|f(x)-f(y)|_inequal_1_W3}Suppose the weight function $w$ has
property W03 for parameter $\kappa$ and that $G$ is the basis function of
order $0$. Then we have the following local pointwise estimates for functions
$f\in X_{w}^{0}$:

\begin{enumerate}
\item If $\kappa\geq1$ and $\alpha\leq\kappa-1$ then%
\[
\left\vert D^{\alpha}f\left(  x\right)  -D^{\alpha}f\left(  y\right)
\right\vert \leq\tfrac{1}{\left(  2\pi\right)  ^{\frac{d}{2}}}\left\Vert
f\right\Vert _{w,0}\left(  \int\frac{\left\vert \xi\right\vert ^{2}%
\xi^{2\alpha}}{w\left(  \xi\right)  }d\xi\right)  ^{\frac{1}{2}}\left\vert
x-y\right\vert .
\]

\item If $\alpha\leq\kappa$ then%
\begin{equation}
\left\vert D^{\alpha}f\left(  x\right)  -D^{\alpha}f\left(  y\right)
\right\vert \leq\tfrac{\sqrt{2}}{\left(  2\pi\right)  ^{\frac{d}{4}}%
}\left\Vert f\right\Vert _{w,0}\sqrt{\left(  -1\right)  ^{\left\vert
\alpha\right\vert }\left(  D^{2\alpha}G\left(  0\right)  -\operatorname{Re}%
\left(  D^{2\alpha}G\right)  \left(  y-x\right)  \right)  }.\label{1.151}%
\end{equation}

\item For the hat weight function in one dimension $\kappa<1/2$ and
\[
\left\vert f\left(  x\right)  -f\left(  y\right)  \right\vert \leq\tfrac
{\sqrt{2}}{\left(  2\pi\right)  ^{\frac{d}{4}}}\left\Vert f\right\Vert
_{w,0}\left\vert x-y\right\vert ^{\frac{1}{2}},\text{\quad}x,y\in
\mathbb{R}^{1},\text{ }\left\vert x-y\right\vert \leq1.
\]

\end{enumerate}
\end{corollary}

\begin{proof}
\textbf{Part 1} This is an application of inequality \ref{1.0391} of Theorem
\ref{Thm_||DRx_minus_DRy||_W3}
\begin{align*}
\left\vert D^{\alpha}f\left(  x\right)  -D^{\alpha}f\left(  y\right)
\right\vert =\left\vert \left(  f,D^{\alpha}R_{x}-D^{\alpha}R_{y}\right)
_{w,0}\right\vert  &  \leq\left\Vert f\right\Vert _{w,0}\left\Vert D^{\alpha
}R_{x}-D^{\alpha}R_{y}\right\Vert _{w,0}\\
&  \leq\tfrac{1}{\left(  2\pi\right)  ^{\frac{d}{2}}}\left\Vert f\right\Vert
_{w,0}\left(  \int\frac{\left\vert \xi\right\vert ^{2}\xi^{2\alpha}}{w\left(
\xi\right)  }d\xi\right)  ^{\frac{1}{2}}\left\vert x-y\right\vert .
\end{align*}
\medskip

\textbf{Part 2} By part 2 of Theorem \ref{Thm_ord0_Riesz_rep_W3}%
\[
\left\vert D^{\alpha}f\left(  x\right)  -D^{\alpha}f\left(  y\right)
\right\vert =\left\vert \left(  f,D^{\alpha}R_{x}-D^{\alpha}R_{y}\right)
_{w,0}\right\vert \leq\left\Vert f\right\Vert _{w,0}\left\Vert D^{\alpha}%
R_{x}-D^{\alpha}R_{y}\right\Vert _{w,0},
\]

and inequality \ref{1.0381} of Theorem \ref{Thm_||DRx_minus_DRy||_W3}
completes the proof.\medskip

\textbf{Part 3} From Theorem \ref{Thm_hat_wt_extend_props} $\kappa<1/2$ so
$\alpha=0$ and since $G\left(  x\right)  =1-\left\vert x\right\vert $ when
$\left\vert x\right\vert \leq1$ our result follows from the bound proved in
part 3.
\end{proof}

The next corollary provides information concerning the smoothness properties
of basis functions.

\begin{corollary}
\label{Cor_bound_DG(z-x)-DG(z-y)_W3}Suppose $G$ is the basis function
generated by a weight function $w$ satisfying property \textbf{W03} for
parameter $\kappa$. Then for each $z\in\mathbb{R}^{d}$, $G$ has the following properties:

\begin{enumerate}
\item If $\alpha\leq\kappa$ and $\beta\leq\kappa$ then%
\begin{equation}
\left\vert D^{\alpha+\beta}G\left(  x\right)  -D^{\alpha+\beta}G\left(
y\right)  \right\vert \leq k\sqrt{\left(  -1\right)  ^{\left\vert
\alpha\right\vert }\left(  D^{2\alpha}G\left(  0\right)  -\operatorname{Re}%
\left(  D^{2\alpha}G\right)  \left(  x-y\right)  \right)  },\label{1.0101}%
\end{equation}

where
\begin{equation}
k=\tfrac{\sqrt{2}}{\left(  2\pi\right)  ^{\frac{d}{4}}}\left(  \int\tfrac
{\xi^{2\alpha}d\xi}{w\left(  \xi\right)  }\right)  ^{\frac{1}{2}}=\sqrt
{2}\sqrt{\left(  -1\right)  ^{\left\vert \alpha\right\vert }D^{2\alpha
}G\left(  0\right)  }.\label{1.0491}%
\end{equation}

\item If $\alpha\leq\kappa$ then%
\[
\left\vert D^{2\alpha}G\left(  x\right)  \right\vert \leq\left(  -1\right)
^{\left\vert \alpha\right\vert }D^{2\alpha}G\left(  0\right)  ,
\]

and%
\[
\left(  -1\right)  ^{\left\vert \alpha\right\vert }\operatorname{Re}%
D^{2\alpha}G\left(  x\right)  <\left(  -1\right)  ^{\left\vert \alpha
\right\vert }D^{2\alpha}G\left(  0\right)  ,\quad x\neq0.
\]

\item Suppose $\kappa\geq1$. Then if $\alpha\leq\kappa-1$ and $\beta\leq
\kappa$,
\begin{align*}
\left\vert D^{\alpha+\beta}G\left(  x\right)  -D^{\alpha+\beta}G\left(
y\right)  \right\vert  & \leq\left(  2\pi\right)  ^{-\frac{d}{2}}\left(
\int\frac{\xi^{2\beta}d\xi}{w\left(  \xi\right)  }\right)  ^{\frac{1}{2}%
}\left(  \int\frac{\left\vert \xi\right\vert ^{2}\xi^{2\alpha}d\xi}{w\left(
\xi\right)  }\right)  ^{\frac{1}{2}}\left\vert x-y\right\vert \\
& =\left(  \left(  -1\right)  ^{\left\vert \beta\right\vert }D^{2\beta
}G\left(  0\right)  \right)  ^{\frac{1}{2}}\left(  \left(  -1\right)
^{1+\left\vert \alpha\right\vert }\left\vert D\right\vert ^{2}D^{2\alpha
}G\left(  0\right)  \right)  ^{\frac{1}{2}}\left\vert x-y\right\vert .
\end{align*}

\end{enumerate}
\end{corollary}

\begin{proof}
\textbf{Part 1.} Our starting point is equation \ref{1.6} i.e. for
$x,y\in\mathbb{R}^{d}$,

$\left(  D^{\alpha}R_{x},D^{\beta}R_{y}\right)  _{w,0}=\left(  -1\right)
^{\left\vert \beta\right\vert }\left(  D^{\alpha+\beta}R_{x}\right)  \left(
y\right)  $ when $\alpha\leq\kappa$, $\beta\leq\kappa$. Hence, since

$R_{x}\left(  y\right)  =\left(  2\pi\right)  ^{-\frac{d}{2}}G\left(
y-x\right)  $%
\begin{align*}
D_{z}^{\alpha+\beta}\left(  G\left(  z-x\right)  -G\left(  z-y\right)
\right)   & =\left(  2\pi\right)  ^{-\frac{d}{2}}D^{\alpha+\beta}\left(
R_{x}-R_{y}\right)  \left(  z\right) \\
& =\left(  2\pi\right)  ^{-\frac{d}{2}}\left(  D^{\alpha}\left(  R_{x}%
-R_{y}\right)  ,D^{\beta}R_{z}\right)  _{w,0},
\end{align*}

so that%
\begin{equation}
\left\vert D_{z}^{\alpha+\beta}\left(  G\left(  z-x\right)  -G\left(
z-y\right)  \right)  \right\vert \leq\left(  2\pi\right)  ^{-\frac{d}{2}%
}\left\Vert D^{\alpha}\left(  R_{x}-R_{y}\right)  \right\Vert _{w,0}\left\Vert
D^{\beta}R_{z}\right\Vert _{w,0}.\label{1.0411}%
\end{equation}

From equation \ref{1.0371} we conclude that $\left\Vert D^{\beta}%
R_{z}\right\Vert _{w,0}=\left(  2\pi\right)  ^{-\frac{d}{2}}\left(  \int%
\frac{\xi^{2\beta}d\xi}{w\left(  \xi\right)  }\right)  ^{\frac{1}{2}}$ and
then by equation \ref{1.0381} for

$\left\Vert D^{\alpha}\left(  R_{x}-R_{y}\right)  \right\Vert _{w,0}$
inequality \ref{1.0101} follows.\medskip

\textbf{Part 2.} We have $\left(  -1\right)  ^{\left\vert \alpha\right\vert
}D^{2\alpha}G\left(  x\right)  =\left(  2\pi\right)  ^{-\frac{d}{2}}\int
e^{i\xi x}\tfrac{\xi^{2\alpha}}{w\left(  \xi\right)  }d\xi$, so that
\[
\left\vert D^{2\alpha}G\left(  x\right)  \right\vert \leq\left(  2\pi\right)
^{-\frac{d}{2}}\int\frac{\xi^{2\alpha}}{w\left(  \xi\right)  }d\xi=\left(
-1\right)  ^{\left\vert \alpha\right\vert }D^{2\alpha}G\left(  0\right)  .
\]

From part 1 we know that $\left(  -1\right)  ^{\left\vert \alpha\right\vert
}\operatorname{Re}\left(  D^{2\alpha}G\right)  \left(  x-y\right)  \leq\left(
-1\right)  ^{\left\vert \alpha\right\vert }D^{2\alpha}G\left(  0\right)  $.
Now suppose that

$\left(  -1\right)  ^{\left\vert \alpha\right\vert }\operatorname{Re}\left(
D^{2\alpha}G\right)  \left(  z\right)  =\left(  -1\right)  ^{\left\vert
\alpha\right\vert }D^{2\alpha}G\left(  0\right)  $ for some $z\neq0$. Equation
\ref{1.0381} would then imply that $D^{\alpha}R_{x}=D^{\alpha}R_{z+x}$ for all
$x$, and hence $z=0$ by part 3 of Theorem \ref{Thm_ord0_Riesz_rep_W3}.\medskip

\textbf{Part 3.} Our starting point is inequality \ref{1.0411} of part 1.
Substitute the expression for $\left\Vert D^{\beta}R_{z}\right\Vert _{w,0}$
given in the proof of part 1 and then use inequality \ref{1.0391} to estimate
$\left\Vert D^{\alpha}\left(  R_{x}-R_{y}\right)  \right\Vert _{w,0}$.
\end{proof}

\begin{remark}
\label{Rem_InTermsOfDeven(0)_W3}Observe that by using \ref{1.0481} and
\ref{1.0491} the right sides of all the estimates involving weight functions
in this subsection can be written in terms of even order derivatives of the
basis function evaluated at the origin.
\end{remark}

Inequality \ref{1.0101} shows that the smoothness of the basis function near
the origin implies its global smoothness. Specifically, for the
one-dimensional hat function $\Lambda$ we have $\kappa<1/2$ and near the
origin $\left\vert \Lambda\left(  x\right)  -\Lambda\left(  y\right)
\right\vert \leq\sqrt{2}\left\vert x-y\right\vert ^{1/2}$ by part 3 of
Corollary \ref{Cor_|f(x)-f(y)|_inequal_1_W3}. However, $\Lambda$ actually
satisfies the stronger estimate $\left\vert \Lambda\left(  x\right)
-\Lambda\left(  y\right)  \right\vert \leq\left\vert x-y\right\vert $
everywhere i.e. it is uniformly Lipschitz continuous on $\mathbb{R}^{1}$ and
Theorem \ref{Thm_ex_nat_spline_basis_Lipschitz} will generalize this result.
First we will need the following extension of the Taylor series expansion with
integral remainder given in Appendix \ref{Sect_apx_TaylorSeries}.
\textbf{Note} that this result only uses the derivatives $\left\{  D^{\beta
}u\right\}  _{\left\vert \beta\right\vert =1}$ and not $\left\{  D^{\beta
}u\right\}  _{\beta\leq1}$.

\begin{lemma}
\label{Lem_Taylor_extension}Suppose $u\in C_{B}^{\left(  0\right)  }\left(
\mathbb{R}^{d}\right)  $ and as distributional derivatives $\left\{  D^{\beta
}u\right\}  _{\left\vert \beta\right\vert =1}\subset L^{\infty}\left(
\mathbb{R}^{d}\right)  $. Then%
\[
u(z+b)=u(z)+\left(  \mathcal{R}_{1}u\right)  \left(  z,b\right)  ,
\]

where $\mathcal{R}_{1}u$ is the integral remainder term
\[
\left(  \mathcal{R}_{1}u\right)  \left(  z,b\right)  =\sum_{\left\vert
\beta\right\vert =1}\frac{b^{\beta}}{\beta!}\int_{0}^{1}(D^{\beta}u)(z+\left(
1-t\right)  b)dt,
\]

which satisfies
\begin{equation}
\left\vert \left(  \mathcal{R}_{1}u\right)  \left(  z,b\right)  \right\vert
\leq\sqrt{d}\max_{\left\vert \beta\right\vert =1}\left\Vert D^{\beta
}u\right\Vert _{\infty;\left[  z,z+b\right]  }\left\vert b\right\vert
,\label{1.18}%
\end{equation}

and%
\begin{equation}
\left\vert \left(  \mathcal{R}_{1}u\right)  \left(  z,b\right)  \right\vert
\leq\left(  \sum_{\left\vert \beta\right\vert =1}\left\Vert D^{\beta
}u\right\Vert _{\infty;\left[  z,z+b\right]  }\right)  \left\vert b\right\vert
.\label{a1.49}%
\end{equation}

\end{lemma}

\begin{proof}
In order to overcome the fact that $D^{\beta}u$ may not be $C^{\left(
0\right)  }\left(  \mathbb{R}^{d}\right)  $ when $\left\vert \beta\right\vert
=1$, we will use a Taylor series expansion with remainder for distributions.
Suppose $\phi\in C_{0}^{\infty}$. Then the conditions on $u$ render all the
integrals absolutely convergent and allows us to apply Fubini's theorem to
swap the order of integration twice in the calculations of this lemma:%
\begin{align}
\left[  u\left(  z+b\right)  ,\phi\left(  z\right)  \right]   & =\left[
u\left(  z\right)  ,\phi\left(  z-b\right)  \right] \nonumber\\
& =\left[  u\left(  z\right)  ,\phi(z)\right]  +\left[  u\left(  z\right)
,\sum_{\left\vert \beta\right\vert =1}\frac{\left(  -b\right)  ^{\beta}}%
{\beta!}\int_{0}^{1}\left(  D^{\beta}\phi\right)  \left(  z+\left(
1-t\right)  \left(  -b\right)  \right)  dt\right] \nonumber\\
& =\left[  u\left(  z\right)  ,\phi\left(  z\right)  \right]  -\sum
_{\left\vert \beta\right\vert =1}\frac{b^{\beta}}{\beta!}\left[  u\left(
z\right)  ,\int_{0}^{1}\left(  D^{\beta}\phi\right)  \left(  z-\left(
1-t\right)  b\right)  dt\right]  .\label{1.36}%
\end{align}

We now analyze the integral remainder term of \ref{1.36}. When $\left\vert
\beta\right\vert =1$,
\begin{align*}
\left[  u\left(  z\right)  ,\int_{0}^{1}\left(  D^{\beta}\phi\right)  \left(
z-\left(  1-t\right)  b\right)  dt\right]   & =\int u\left(  z\right)
\int_{0}^{1}\left(  D^{\beta}\phi\right)  \left(  z-\left(  1-t\right)
b\right)  dt\text{ }dz\\
& =\int_{0}^{1}\int u\left(  z\right)  \left(  D^{\beta}\phi\right)  \left(
z-\left(  1-t\right)  b\right)  dz\text{ }dt\\
& =\int_{0}^{1}\left[  u\left(  z\right)  ,\left(  D^{\beta}\phi\right)
\left(  z-\left(  1-t\right)  b\right)  \right]  \text{ }dt\\
& =\int_{0}^{1}\left[  u\left(  z+\left(  1-t\right)  b\right)  ,D^{\beta}%
\phi\left(  z\right)  \right]  \text{ }dt\\
& =\left(  -1\right)  ^{\left\vert \beta\right\vert }\int_{0}^{1}\left[
\left(  D^{\beta}u\right)  \left(  z+\left(  1-t\right)  b\right)
,\phi\left(  z\right)  \right]  \text{ }dt\\
& =-\int_{0}^{1}\int\left(  D^{\beta}u\right)  \left(  z+\left(  1-t\right)
b\right)  \phi\left(  z\right)  dz\text{ }dt\\
& =-\int\int_{0}^{1}\left(  D^{\beta}u\right)  \left(  z+\left(  1-t\right)
b\right)  dt\text{ }\phi\left(  z\right)  dz\\
& =-\left[  \int_{0}^{1}\left(  D^{\beta}u\right)  \left(  z+\left(
1-t\right)  b\right)  dt,\phi\left(  z\right)  \right]  ,
\end{align*}

so \ref{1.36} now becomes%
\[
\left[  u\left(  z+b\right)  ,\phi\left(  z\right)  \right]  =\left[  u\left(
z\right)  ,\phi\left(  z\right)  \right]  -\sum_{\left\vert \beta\right\vert
=1}\frac{b^{\beta}}{\beta!}\left[  \int_{0}^{1}\left(  D^{\beta}u\right)
\left(  z+\left(  1-t\right)  b\right)  dt,\phi\left(  z\right)  \right]  ,
\]

for $\phi\in C_{0}^{\infty}$. Thus%
\[
u\left(  z+b\right)  =u\left(  z\right)  +\sum_{\left\vert \beta\right\vert
=1}\frac{b^{\beta}}{\beta!}\int_{0}^{1}\left(  D^{\beta}u\right)  \left(
z+\left(  1-t\right)  b\right)  dt,
\]

as claimed. Finally, the estimate \ref{1.18} is proved in a very similar
manner to the integral remainder estimate \ref{1.34} of Appendix
\ref{Ch_Appendx_basic_notation} and the estimate \ref{a1.49} is proved in a
very similar manner to the remainder estimate \ref{a1.37} of Appendix
\ref{Ch_Appendx_basic_notation}.
\end{proof}

\begin{remark}
\label{Rem_Lem_Taylor_extension}\textbf{Local Taylor expansion} Suppose
$\Omega$ is a region. The Taylor series expansion of the above lemma also
holds if $u\in C_{B}^{\left(  0\right)  }\left(  \mathbb{\Omega}\right)  $,
$\left\{  D^{\beta}u\right\}  _{\left\vert \beta\right\vert =1}\subset
L^{\infty}\left(  \mathbb{\Omega}\right)  $ and $\left[  z,z+b\right]
\subset\Omega$.
\end{remark}

We now prove that all derivatives up to order $\left(  2n-2\right)
\mathbf{1}$ of the extended B-spline tensor product basis function are
uniformly Lipschitz continuous of order $1$.

\begin{theorem}
\label{Thm_ex_nat_spline_basis_Lipschitz}Let $G_{s}\left(  x\right)
=\prod\limits_{k=1}^{d}G_{1}\left(  x_{k}\right)  $ be an extended B-spline
tensor product basis function, as introduced in Subsection
\ref{SbSect_basis_fns}. We then have the estimates
\[
\left\vert G_{s}\left(  x\right)  -G_{s}\left(  y\right)  \right\vert
\leq\sqrt{d}G_{1}\left(  0\right)  ^{d-1}\left\Vert DG_{1}\right\Vert
_{\infty}\left\vert x-y\right\vert ,\text{\quad}x,y\in\mathbb{R}^{d},
\]

and if $\alpha\leq\left(  2n-2\right)  \mathbf{1}$,%
\begin{equation}
\left\vert D^{\alpha}G_{s}\left(  x\right)  -D^{\alpha}G_{s}\left(  y\right)
\right\vert \leq\sqrt{d}\max\limits_{k=1}^{d}\left\Vert D_{k}D^{\alpha}%
G_{s}\right\Vert _{\infty}\left\vert x-y\right\vert ,\text{\quad}%
x,y\in\mathbb{R}^{d}.\label{1.027}%
\end{equation}

\end{theorem}

\begin{proof}
From Theorem \ref{Thm_basis_tensor_hat_W3} we know that $G_{s}\in
C_{0}^{\left(  2n-2\right)  \mathbf{1}}\left(  \mathbb{R}^{d}\right)  $ and
that the derivatives

$\left\{  D^{\beta}G_{s}:\left(  2n-2\right)  \mathbf{1}<\beta\leq\left(
2n-1\right)  \mathbf{1}\right\}  $ are bounded functions. Consequently
$G_{s}\in C_{0}^{\left(  0\right)  }\left(  \mathbb{R}^{d}\right)  $ and all
the first derivatives are bounded functions. This implies $u=G_{s}$ satisfies
the conditions of Lemma \ref{Lem_Taylor_extension} with $z=y$ and $b=x-y$ so
the estimate of that lemma holds i.e.
\begin{equation}
\left\vert G_{s}\left(  x\right)  -G_{s}\left(  y\right)  \right\vert
\leq\sqrt{d}\max\limits_{\left\vert \beta\right\vert =1}\left\Vert D^{\beta
}G_{s}\right\Vert _{\infty}\left\vert x-y\right\vert .\label{1.005}%
\end{equation}

Again by Theorem \ref{Thm_basis_tensor_hat_W3}, $G_{1}\in C_{0}^{\left(
0\right)  }\left(  \mathbb{R}^{1}\right)  $ and $DG_{1}$ is a bounded function
so that%
\[
\left\vert G_{s}\left(  x\right)  -G_{s}\left(  y\right)  \right\vert
\leq\sqrt{d}\left\Vert G_{1}\right\Vert _{\infty}^{d-1}\left\Vert
DG_{1}\right\Vert _{\infty}\left\vert x-y\right\vert .
\]

Finally, part 2 of Corollary \ref{Cor_bound_DG(z-x)-DG(z-y)_W3} tells us that
$\left\Vert G_{1}\right\Vert _{\infty}\leq G_{1}\left(  0\right)  $ which
completes the proof of the first inequality. To prove the second inequality we
simply replace $G_{s}$ by $D^{\alpha}G_{s}$ in \ref{1.005}.
\end{proof}

\begin{remark}
\label{Rem_Thm_ex_nat_spline_basis_Lipschitz}I have realized that the
multivariate theorem can be proven by writing%
\begin{align*}
G_{s}\left(  x\right)   &  -G_{s}\left(  y\right) \\
&  =\left\{  G_{s}\left(  x\right)  -G_{1}\left(  y_{1}\right)  G_{s}\left(
x_{2:d}\right)  \right\}  +\left\{  G_{1}\left(  y_{1}\right)  G_{s}\left(
x_{2:d}\right)  -G_{s}\left(  y_{1:2}\right)  G_{s}\left(  x_{3:d}\right)
\right\}  +\ldots\\
&  =\left\{  G_{1}\left(  x_{1}\right)  -G_{1}\left(  y_{1}\right)  \right\}
G_{s}\left(  x_{2:d}\right)  +G_{1}\left(  y_{1}\right)  \left\{  G_{1}\left(
x_{2}\right)  -G_{1}\left(  y_{2}\right)  \right\}  G_{s}\left(
x_{3:d}\right)  +\ldots,
\end{align*}

so that 1-dimensional estimates can be applied to give%
\begin{align*}
\left\vert G_{s}\left(  x\right)  -G_{s}\left(  y\right)  \right\vert  &
\leq\left\vert G_{1}\left(  x_{1}\right)  -G_{1}\left(  y_{1}\right)
\right\vert \left\Vert G_{1}\right\Vert _{\infty}^{d-1}+\left\vert
G_{1}\left(  x_{2}\right)  -G_{1}\left(  y_{2}\right)  \right\vert \left\Vert
G_{1}\right\Vert _{\infty}^{d-1}+\ldots\\
& \leq\left(  \left\Vert DG_{1}\right\Vert _{\infty}\left\Vert G_{1}%
\right\Vert _{\infty}^{d-1}\right)  \left\vert x-y\right\vert _{1}\\
& \leq\sqrt{d}\left\Vert G_{1}\right\Vert _{\infty}^{d-1}\left\Vert
DG_{1}\right\Vert _{\infty}\left\vert x-y\right\vert _{2}.
\end{align*}

\end{remark}

\chapter{The minimal norm interpolant\label{Ch_Interpol}}

\section{Introduction\label{Sect_InroPart2}}

\textbf{The interpolant and smoother convergence results obtained in this
document all assume the weight function has property W02}.\textbf{\ I have not
been able to use the weight function property W03 to obtain better point-wise
orders of convergence. }However, Theorem \ref{Thm_equiv_W3} tells us that if a
weight function has property W03 for $\kappa=\mu$ then it has property W02 for
$\kappa=\underline{\mu}=\min\mu=\min_{i}\mu_{i}$.

In Chapter \ref{Ch_wtfn_basisfn_datasp} of this document we introduced the
basic mathematical machinery: weight function, data space and basis function.
In this chapter this theory is applied to the well-known minimal norm
interpolation problem. This problem is solved and several pointwise
convergence results are derived. These are illustrated with numerical results
obtained using several extended B-spline basis functions and special classes
of data functions. The study of interpolation convergence is continued in the
next three chapters.

In more detail, the functions from the Hilbert data space $X_{w}^{0}$ are used
to define the standard minimal norm interpolation problem with independent
data $X=\left\{  x^{(i)}\right\}  _{i=1}^{N}$ and dependent data $y=\left\{
y_{i}\right\}  _{i=1}^{N}$, the latter being obtained by evaluating a data
function at $X$. This problem is then shown to have a unique basis function
solution of the form $\sum\limits_{i=1}^{N}\alpha_{i}G\left(  \cdot
-x^{(i)}\right)  $ with $\alpha_{i}\in\mathbb{C}$. The standard matrix
equation for the $\alpha_{i}$ is then derived.

We will then consider several categories of estimates for the pointwise
convergence of the interpolant to its data function when the data is confined
to a bounded data region. In all these results the order of convergence
appears as the power of the radius of the largest ball in the data region that
can be fitted between the $X$ data points. The convergence results of this
document can be divided into those which \textbf{explicitly} use Lagrange
interpolation theory, and thus unisolvent subsets of $X$ (Definition
\ref{Def_unisolv}), and those which do not. The non-unisolvency proofs are
\textit{much} simpler than the unisolvency-based results and if the data
functions are chosen appropriately the constants can be calculated. However,
the maximum order of convergence obtained is $1$, no matter how large the
parameter $\kappa$. On the other hand, the constants for the unisolvency-based
results are more difficult to calculate but the order of convergence is at
least equal to the value of $\overline{\left\lfloor \kappa\right\rfloor }%
$.\medskip

\fbox{Type 1 error estimates: no Taylor series expansion} Suppose that the
weight function has property W02 for some $\kappa\geq0$, that the data $X$ is
contained in a closed bounded infinite data set $K$, and that the basis
function $G$ satisfies the estimate
\[
\left\vert G\left(  0\right)  -\operatorname{Re}G\left(  x\right)  \right\vert
\leq C_{G}\left\vert x\right\vert ^{2s},\text{\quad}\left\vert x\right\vert
<h_{G},
\]

for some $h_{G}>0$. Then in Theorem \ref{Thm_|f(x)-f(y)|_inequal_2} it is
shown that the interpolant $\mathcal{I}_{X}f$ of a data function $f\in
X_{w}^{0}$ satisfies%
\[
\left\vert f\left(  x\right)  -\mathcal{I}_{X}f\left(  x\right)  \right\vert
\leq k_{G}\left\Vert f\right\Vert _{w,0}\left(  h_{X,K}\right)  ^{s}%
,\text{\quad}x\in K,
\]

when $h_{X,K}=\sup\limits_{x\in K}\operatorname*{dist}\left(  x,X\right)  \leq
h_{G}$ and $k_{G}=\left(  2\pi\right)  ^{-\frac{d}{4}}\sqrt{2C_{G}}$. This
implies an order of convergence of at least $s$.\medskip

\fbox{Type 2 error estimates: no Taylor series expansion} Here we avoid making
any assumptions about $G$ and instead assume $\kappa\geq1$ and use the
properties of the Riesz representer $R_{x}$ of the evaluation functional
$f\rightarrow f\left(  x\right)  $. A consequence of this approach is the
estimate (Theorem \ref{Thm_interpol_error_in_terms_of_wt_fn}): when
$h_{X,K}<\infty$
\[
\left\vert f\left(  x\right)  -\mathcal{I}_{X}f\left(  x\right)  \right\vert
\leq k_{G}\left\Vert f\right\Vert _{w,0}h_{X,K},\text{\quad}x\in K,\text{
}f\in X_{w}^{0},
\]

where%
\[
k_{G}=\left(  2\pi\right)  ^{-d/4}\sqrt{-\left(  \left\vert D\right\vert
^{2}G\right)  \left(  0\right)  }.
\]

so the order of convergence is \textbf{always} at least $1$.\medskip

\fbox{Estimates based explicitly on Lagrange interpolation/unisolvency} Here
$X$ is contained in a bounded data region $\Omega$ and $X$ is assumed to have
$m$-unisolvent subsets of order $m=\overline{\left\lfloor \kappa\right\rfloor
}$. In this case it follows from Theorem \ref{Thm_converg_interpol_ord_gte_1}
that there exist constants $h_{\Omega,\kappa},k_{G}>0$ such that%
\[
\left\vert f\left(  x\right)  -\mathcal{I}_{X}f\left(  x\right)  \right\vert
\leq k_{G}\left\Vert f\right\Vert _{w,0}\left(  h_{X,\Omega}\right)
^{m},\quad x\in\overline{\Omega},\text{ }f\in X_{w}^{0},
\]

when $h_{X,\Omega}=\sup\limits_{\omega\in\Omega}\operatorname*{dist}\left(
\omega,X\right)  <h_{\Omega,\kappa}$.

\textbf{Numerical results are only presented for the Type 1 and Type 2 cases
in 1 dimension} and these illustrate the convergence of the interpolant to
it's data function. We will only be interested in the convergence of the
interpolant to it's data function and not in the algorithm's performance as an
interpolant. Only the extended B-splines will be considered and we will also
restrict ourselves to one dimension so that the data density parameter
$h_{X,\Omega}$ can be easily calculated. If we can calculate the data function
norm we can calculate the error estimate. All the extended B-spline basis
weight functions have a power of $\sin^{2}x$ in the denominator and so we have
derived special classes of data functions for which the data function norm can
be calculated. The derivation of these special classes of data functions led
to the characterization of the restriction spaces $X_{w}^{0}\left(
\Omega\right)  $ for various class of weight function. For example, Theorem
\ref{Thm_int_Xow(O)_eq_Hn(O)_dim1} shows that the restrictions of the B-spline
data functions are a member of the class of \textbf{Sobolev spaces}:
\begin{align*}
W^{m\mathbf{1}}\left(  \Omega\right)   & =\left\{  u\in L^{2}\left(
\Omega\right)  :D^{\alpha}u\in L^{2}\left(  \Omega\right)  \text{ }for\text{
}\alpha_{i}\leq m,\text{ }0\leq i\leq m\right\}  ,\text{\quad}m=1,2,3,\ldots\\
& =\left\{  u\in L^{2}\left(  \Omega\right)  :D^{\alpha}u\in L^{2}\left(
\Omega\right)  \text{ }for\text{ }\alpha\leq m\mathbf{1}\right\}
,\text{\quad}m=1,2,3,\ldots
\end{align*}

It is easy to construct functions that are in these spaces.

As expected interpolant instability is evident and because our error estimates
assume an infinite precision we filter the error to remove spikes which are a
manifestation of the instability.

\section{The space $W_{G,X}$}

In this section we will define the space $W_{G,X}$ which contains the solution
to the minimal norm interpolation problem \ref{1.42}. Here $X$ is the
independent data and $G$ is the basis function. The next two results will be
required to show that $W_{G,X}$ is a well defined finite dimensional vector
space. The next theorem requires the following lemma which we state without proof.

\begin{lemma}
\label{Lem_lin_comb_exp(i)_equals_0}Suppose $\left\{  x^{\left(  k\right)
}\right\}  _{k=1}^{N}$ is a set of distinct points in $\mathbb{R}^{d}$ and
$v=\left(  \nu_{k}\right)  _{k=1}^{N}\in\mathbb{C}^{N}$. Define the function
$a_{v}$ by
\[
a_{v}\left(  \xi\right)  =\sum_{k=1}^{N}\nu_{k}e^{-ix^{\left(  k\right)  }\xi
}.
\]

Then:

\begin{enumerate}
\item if $a_{v}\left(  \xi\right)  =0$ for all $\xi$ then $v_{k}=0$ for all
$k$.

\item The null space of $a_{v}$ is a closed set of measure zero.
\end{enumerate}
\end{lemma}

\begin{theorem}
\label{Thm_indep_G(x-xi)}Let $X=\left\{  x^{(k)}\right\}  _{k=1}^{N}$ be $N$
distinct points in $\mathbb{R}^{d}$.

Then the set of translated basis functions\allowbreak\ $\left\{  G\left(
\cdot-x^{(k)}\right)  \right\}  _{k=1}^{N}$ is linearly independent with
respect to the complex scalars.
\end{theorem}

\begin{proof}
Suppose for complex $\alpha_{k}$, $\sum\limits_{k=1}^{N}\alpha_{k}G\left(
x-x^{(k)}\right)  =0$. Taking the Fourier transform we obtain%
\[
\frac{1}{w\left(  \xi\right)  }\sum\limits_{k=1}^{N}\alpha_{k}e^{ix^{(k)}\xi
}=0\quad a.e.
\]

Since $w\left(  \xi\right)  >0$ a.e. we have $\sum\limits_{k=1}^{N}\alpha
_{k}e^{ix^{(k)}\xi}=0$ and part 1 of Lemma \ref{Lem_lin_comb_exp(i)_equals_0}
implies $\alpha_{k}=0$ for all $k$.
\end{proof}

We now introduce the space $W_{G,X}$.

\begin{definition}
\label{Def_Wg,x_m_gt_0}\textbf{The finite dimensional vector space }$W_{G,X}$

Suppose the weight function $w$ has property W02 for parameter $\kappa$, let
$G$ be the (continuous) basis function generated by $w$ and let $X=\left\{
x^{(k)}\right\}  _{k=1}^{N}$ be a set of distinct points in $\mathbb{R}^{d}$.
We have shown in Theorem \ref{Thm_indep_G(x-xi)} that the functions $G\left(
\cdot-x^{(k)}\right)  $ are linearly independent so the span of these
functions makes sense. We define the $N-$ dimensional vector space $W_{G,X}$
by
\begin{equation}
W_{G,X}=\left\{  \sum_{k=1}^{N}\alpha_{k}G\left(  \cdot-x^{(k)}\right)
:\alpha_{k}\in\mathbb{C}\right\}  .\label{1.5}%
\end{equation}

When convenient below, functions in $W_{G,X}$ will be written in the form

$f_{\alpha}\left(  x\right)  =\sum\limits_{k=1}^{N}\alpha_{k}G\left(
x-x^{(k)}\right)  $ where $\alpha=\left(  \alpha_{k}\right)  \in\mathbb{C}%
^{N}$.
\end{definition}

\section{The matrices $G_{X,X}$ and $R_{X,X}$}

We will now define the \textit{basis function matrix} $G_{X,X}$ and the
\textit{reproducing kernel matrix} $R_{X,X}$. In this document we deal with
basis functions of order zero so we will deduce the simple relationship
$R_{X,X}=\left(  2\pi\right)  ^{-d/2}G_{X,X}$. These matrices will be used to
construct the matrix equations for the interpolants studied in this document
and the smoothers studied in Chapters \ref{Ch_Exact_smth} and
\ref{Ch_Approx_smth}.

\begin{definition}
\textbf{The basis function matrix }$G_{X,X}$

Let $X=\left\{  x^{\left(  n\right)  }\right\}  _{n=1}^{N}$ be a set of
distinct points in $\mathbb{R}^{d}$ and suppose $G$ is the basis function
generated by $w$. Then the basis function matrix $G_{X,X}$ is defined by
\[
G_{X,X}=\left(  G\left(  x^{\left(  i\right)  }-x^{\left(  j\right)  }\right)
\right)  .
\]

\end{definition}

\begin{remark}
Since $G=\left(  \frac{1}{w}\right)  ^{\vee}$ we have $\overline{G\left(
-x\right)  }=G\left(  x\right)  $ and so the matrix $G_{X,X}$ is Hermitian.
\end{remark}

\begin{definition}
\label{Def_Matrices_from_R}\textbf{The reproducing kernel matrix }$R_{X,X}$

Suppose $R_{x}$ is the Riesz representer of the evaluation functional
$f\rightarrow f\left(  x\right)  $ introduced in Theorem
\ref{Thm_ord0_Riesz_rep_W2}. Suppose that $X=\left\{  x^{\left(  k\right)
}\right\}  _{k=1}^{N}$ is a set of distinct points in $\mathbb{R}^{d}$. Then
the reproducing kernel matrix is
\[
R_{X,X}=\left(  R_{x^{\left(  j\right)  }}(x^{\left(  i\right)  })\right)  .
\]

I call this the reproducing kernel matrix because the \textbf{reproducing
kernel} $K\left(  x,y\right)  $ satisfies $K\left(  x,y\right)  =R_{x}(y)$.

Using a separate symbol for the matrix $R_{X,X}$ fits in with the case of
positive order, where $R_{X,X}$ is not simply a scalar multiple of $G_{X,X}$
but is related by a more complex formula.
\end{definition}

The reproducing kernel matrix\textbf{\ }$R_{X,X}$ has the following properties:

\begin{theorem}
\label{Thm_Rx,x_properties}The reproducing kernel matrix\textbf{\ }$R_{X,X}$
has the following properties:

\begin{enumerate}
\item $R_{X,X}=\left(  2\pi\right)  ^{-d/2}G_{X,X}$.

\item $R_{X,X}$ is a Gram matrix and hence positive definite, Hermitian and regular.

\item The functions $\left\{  R_{x^{(k)}}\right\}  _{k=1}^{N}$ are independent.
\end{enumerate}
\end{theorem}

\begin{proof}
\textbf{Part 1} True since in Theorem \ref{Thm_ord0_Riesz_rep_W2} it was shown
that $R_{x}(z)=\left(  2\pi\right)  ^{-d/2}G(z-x)$.\medskip

\textbf{Part 2} $R_{X,X}=\left(  R_{x^{\left(  j\right)  }}(x^{\left(
i\right)  })\right)  $. But from the definition of $R_{x}$, $R_{x^{\left(
j\right)  }}(x^{\left(  i\right)  })=\left(  R_{x^{\left(  j\right)  }%
},R_{x^{\left(  i\right)  }}\right)  _{w,0}$ so $R_{X,X}$ is a Gram matrix and
hence is positive definite over $\mathbb{C}$, Hermitian and regular.\medskip

\textbf{Part 3} In Theorem \ref{Thm_indep_G(x-xi)} it was shown that the
functions $G(z-x^{\left(  k\right)  })$ are independent and since
$R_{x}(z)=\left(  2\pi\right)  ^{-d/2}G(z-x)$ the functions $R_{x^{(k)}}$ are independent.
\end{proof}

\section{The vector-valued evaluation operator $\protect\widetilde{\mathcal{E}%
}_{X}$}

The vector-valued evaluation operator $\widetilde{\mathcal{E}}_{X}$ and its
Hilbert space adjoint $\widetilde{\mathcal{E}}_{X}^{\ast}$ are the fundamental
operators used to solve the variational interpolation and smoothing problems.

\begin{definition}
\label{Def_vect_val_eval_op}\textbf{The vector-valued evaluation operator
}$\widetilde{\mathcal{E}}_{X}$

Let $X=\left\{  x^{(i)}\right\}  _{i=1}^{N}$ be a set of $N$ distinct points
in $\mathbb{R}^{d}$. Let $u$ be a complex-valued continuous function. Then the
evaluation operator $\widetilde{\mathcal{E}}_{X}$ is defined by
\[
\widetilde{\mathcal{E}}_{X}u=\left(  u\left(  x^{(i)}\right)  \right)
_{i=1}^{N}.
\]

Sometimes we will use the notation $u_{X}$ for $\widetilde{\mathcal{E}}_{X}u$
and when dealing with matrices $\widetilde{\mathcal{E}}_{X}u$ will be regarded
as a column vector.
\end{definition}

\begin{theorem}
\label{Thm_eval_op_properties}\textbf{Properties of }$\widetilde{\mathcal{E}%
}_{X}$ \ Let $X=\left\{  x^{\left(  i\right)  }\right\}  _{i=1}^{N}$ be a set
of distinct points in $\mathbb{R}^{d}$ and suppose that the weight function
$w$ has property W02. Then we know that $X_{w}^{0}$ is a reproducing kernel
Hilbert space of continuous functions and the evaluation functional
$f\rightarrow f(x)$ has a Riesz representer, say $R_{x}$. The evaluation
operator $\widetilde{\mathcal{E}}_{X}$ will be now shown to have the following properties:

\begin{enumerate}
\item $\widetilde{\mathcal{E}}_{X}:\left(  X_{w}^{0},\left\Vert \cdot
\right\Vert _{w,0}\right)  \rightarrow\left(  \mathbb{C}^{N},\left\vert
\cdot\right\vert \right)  $ is continuous, onto and $\operatorname*{null}%
\widetilde{\mathcal{E}}_{X}=W_{G,X}^{\bot}$.

\item The adjoint operator $\widetilde{\mathcal{E}}_{X}^{\ast}:\mathbb{C}%
^{N}\rightarrow X_{w}^{0}$, defined by $\left(  \widetilde{\mathcal{E}}%
_{X}f,g\right)  _{\mathbb{C}^{N}}=\left(  f,\widetilde{\mathcal{E}}_{X}^{\ast
}g\right)  _{w,0}$, satisfies
\[
\widetilde{\mathcal{E}}_{X}^{\ast}\beta=\sum\limits_{i=1}^{N}\beta
_{i}R_{x^{\left(  i\right)  }},\text{\quad}\beta\in\mathbb{C}^{N},
\]

and is a homeomorphism from $\left(  \mathbb{C}^{N},\left\vert \cdot
\right\vert \right)  $ to $\left(  W_{G,X},\left\Vert \cdot\right\Vert
_{w,0}\right)  $.

\item Suppose $R_{X,X}=\left(  R_{x^{\left(  j\right)  }}\left(  x^{\left(
i\right)  }\right)  \right)  $ is the\textbf{\ }reproducing kernel matrix and
$\left\Vert \cdot\right\Vert $ is the matrix norm corresponding to the
Euclidean vector norm. Then%
\[
\left\Vert \widetilde{\mathcal{E}}_{X}^{\ast}\right\Vert =\left\Vert
\widetilde{\mathcal{E}}_{X}\right\Vert =\left\Vert R_{X,X}\right\Vert
\leq\sqrt{N}\sqrt{R_{0}\left(  0\right)  }.
\]

\item The operator $\widetilde{\mathcal{E}}_{X}^{\ast}\widetilde{\mathcal{E}%
}_{X}$ is self-adjoint and we have $\left\Vert \widetilde{\mathcal{E}}%
_{X}^{\ast}\widetilde{\mathcal{E}}_{X}\right\Vert =\left\Vert R_{X,X}%
\right\Vert ^{2}$, as well as the formulas
\begin{equation}
\mathcal{E}_{X}^{\ast}\widetilde{\mathcal{E}}_{X}f=%
{\textstyle\sum\limits_{i=1}^{N}}
f\left(  x^{\left(  i\right)  }\right)  R_{x^{\left(  i\right)  }},\quad f\in
X_{w}^{0},\label{1.52}%
\end{equation}

and%
\begin{equation}
\left(  \widetilde{\mathcal{E}}_{X}^{\ast}\widetilde{\mathcal{E}}_{X}f\right)
\left(  x\right)  =\left(  \widetilde{\mathcal{E}}_{X}f,\widetilde{\mathcal{E}%
}_{X}R_{x}\right)  _{\mathbb{C}^{N}}=\left(  f,\widetilde{\mathcal{E}}%
_{X}^{\ast}\widetilde{\mathcal{E}}_{X}R_{x}\right)  _{w,0},\quad f\in
X_{w}^{0}.\label{1.53}%
\end{equation}

Also, $\widetilde{\mathcal{E}}_{X}^{\ast}\widetilde{\mathcal{E}}_{X}:X_{w}%
^{0}\mathbb{\rightarrow}W_{G,X}$ is onto and $\operatorname*{null}%
\widetilde{\mathcal{E}}_{X}^{\ast}\widetilde{\mathcal{E}}_{X}=W_{G,X}^{\perp}$.

\item The operator $\widetilde{\mathcal{E}}_{X}\widetilde{\mathcal{E}}%
_{X}^{\ast}$ is self-adjoint and if we regard the range of
$\widetilde{\mathcal{E}}_{X}$ as column vectors%
\[
\widetilde{\mathcal{E}}_{X}\widetilde{\mathcal{E}}_{X}^{\ast}\beta
=R_{X,X}\beta,\quad\beta\in\mathbb{C}^{N},
\]

where $R_{X,X}$ is the\textbf{\ }(regular)\textbf{\ }reproducing kernel matrix
introduced in Definition \ref{Def_Matrices_from_R}.

\item If $X=\left\{  x^{\left(  i\right)  }\right\}  _{i=1}^{N}$ and
$Y=\left\{  y^{\left(  j\right)  }\right\}  _{j=1}^{N^{\prime}}$ define
$R_{X,Y}=\left(  R_{y^{\left(  j\right)  }}\left(  x^{\left(  i\right)
}\right)  \right)  $ and $G_{X,Y}=\left(  G\left(  x^{\left(  i\right)
}-y^{\left(  j\right)  }\right)  \right)  $ so that $R_{X,Y}=\left(
2\pi\right)  ^{-\frac{d}{2}}G_{X,Y}$. Now when $\widetilde{\mathcal{E}}_{Y}$
is assumed to be a column vector we have $\widetilde{\mathcal{E}}%
_{X}\widetilde{\mathcal{E}}_{Y}^{\ast}\beta=R_{X,Y}\beta$ when $\beta
\in\mathbb{C}^{N^{\prime}}$.
\end{enumerate}
\end{theorem}

\begin{proof}
\textbf{Parts 1 and 2} That $\widetilde{\mathcal{E}}_{X}$\textbf{\ }is
continuous follows from the Cauchy-Schwartz inequality and
\[
\left\vert \widetilde{\mathcal{E}}_{X}u\right\vert _{\mathbb{C}^{N}}^{2}%
=\sum_{i=1}^{N}\left\vert u(x^{(i)})\right\vert ^{2}=\sum_{i=1}^{N}\left\vert
\left(  u,R_{x^{(i)}}\right)  _{w,0}\right\vert ^{2}\leq\left(  \sum_{i=1}%
^{N}\left\Vert R_{x^{(i)}}\right\Vert _{w,0}^{2}\right)  \left\Vert
u\right\Vert _{w,0}^{2}.
\]

Clearly $\widetilde{\mathcal{E}}_{X}u=0$ iff $\left(  u,R_{x^{(i)}}\right)
_{w,0}=0$ for each $R_{x^{(i)}}$ so that $\operatorname*{null}%
\widetilde{\mathcal{E}}_{X}=W_{G,X}^{\perp}$.

Next we show that $\widetilde{\mathcal{E}}_{X}$ is onto. The Hilbert space
adjoint of $\widetilde{\mathcal{E}}_{X}$ is denoted $\widetilde{\mathcal{E}%
}_{X}^{\ast}$ and is defined by
\[
\left(  \widetilde{\mathcal{E}}_{X}u,\beta\right)  _{\mathbb{C}^{N}}=\left(
u,\widetilde{\mathcal{E}}_{X}^{\ast}\beta\right)  _{w,0}.
\]

The adjoint is calculated using the representer $R_{x}$:
\[
\left(  \widetilde{\mathcal{E}}_{X}u,\beta\right)  _{\mathbb{C}^{N}}%
=\sum_{i=1}^{N}u\left(  x^{\left(  i\right)  }\right)  \overline{\beta_{i}%
}=\sum_{i=1}^{N}\left(  u,R_{x^{\left(  i\right)  }}\right)  _{w,0}%
\overline{\beta_{i}}=\left(  u,\sum_{i=1}^{N}\beta_{i}R_{x^{\left(  i\right)
}}\right)  _{w,0},
\]

so that
\[
\widetilde{\mathcal{E}}_{X}^{\ast}\beta=\sum_{i=1}^{N}\beta_{i}R_{x^{\left(
i\right)  }},\quad\beta\in\mathbb{C}^{N}.
\]

Finally we show that $\widetilde{\mathcal{E}}_{X}^{\ast}$ is 1-1. But
$\widetilde{\mathcal{E}}_{X}^{\ast}\beta=0$ implies $\sum\limits_{i=1}%
^{N}\beta_{i}R_{x^{\left(  i\right)  }}=0$ and,

since $R_{x^{\left(  i\right)  }}=\left(  2\pi\right)  ^{-d/2}G\left(
\cdot-x^{\left(  i\right)  }\right)  $, Theorem \ref{Thm_indep_G(x-xi)}
implies the $R_{x^{\left(  i\right)  }}$ are linearly independent and so
$\beta=0$.\medskip

\textbf{Part 3} That $\left\Vert \widetilde{\mathcal{E}}_{X}^{\ast}\right\Vert
_{op}=\left\Vert \widetilde{\mathcal{E}}_{X}\right\Vert _{op}$ is an
elementary property of the adjoint. Now%
\begin{align*}
\left\Vert \widetilde{\mathcal{E}}_{X}^{\ast}\beta\right\Vert _{w,0}%
^{2}=\left(  \sum_{i=1}^{N}\beta_{i}R_{x^{\left(  i\right)  }},\sum_{j=1}%
^{N}\beta_{j}R_{x^{\left(  j\right)  }}\right)  _{w,0} &  =\sum_{i,j=1}%
^{N}\beta_{i}\overline{\beta_{j}}\left(  R_{x^{\left(  i\right)  }%
},R_{x^{\left(  j\right)  }}\right)  _{w,0}\\
&  =\sum_{i,j=1}^{N}\beta_{i}\overline{\beta_{j}}R_{x^{\left(  j\right)  }%
}\left(  x^{\left(  i\right)  }\right)  =\beta^{T}R_{X,X}\overline{\beta},
\end{align*}

so that, $\left\Vert \widetilde{\mathcal{E}}_{X}^{\ast}\right\Vert
=\max\limits_{\beta\in\mathbb{C}^{N}}\frac{\left\Vert \widetilde{\mathcal{E}%
}_{X}^{\ast}\beta\right\Vert _{w,0}}{\left\vert \beta\right\vert }%
=\max\limits_{\beta\in\mathbb{C}^{N}}\frac{\sqrt{\beta^{T}R_{X,X}%
\overline{\beta}}}{\left\vert \beta\right\vert }$. But the latter expression
is the largest (positive) eigenvalue of the Hermitian matrix $R_{X,X}$ i.e.
the value of $\left\Vert R_{X,X}\right\Vert $. Hence $\left\Vert
\widetilde{\mathcal{E}}_{X}^{\ast}\right\Vert =\left\Vert R_{X,X}\right\Vert
$.%
\begin{align*}
\left\vert \widetilde{\mathcal{E}}_{X}f\right\vert ^{2}=\sum\limits_{k=1}%
^{N}\left\vert f\left(  x^{\left(  k\right)  }\right)  \right\vert ^{2}%
=\sum\limits_{k=1}^{N}\left\vert \left(  f,R_{x^{\left(  k\right)  }}\right)
_{w,0}\right\vert ^{2} &  \leq\sum\limits_{k=1}^{N}\left\Vert f\right\Vert
_{w,0}^{2}\left\Vert R_{x^{\left(  k\right)  }}\right\Vert _{w,0}^{2}\\
&  =\sum\limits_{k=1}^{N}\left\Vert f\right\Vert _{w,0}^{2}R_{0}\left(
0\right)  =NR_{0}\left(  0\right)  \left\Vert f\right\Vert _{w,0}^{2}.
\end{align*}

\textbf{Part 4} The fact that $\widetilde{\mathcal{E}}_{X}^{\ast
}\widetilde{\mathcal{E}}_{X}$ is self-adjoint and the formulas \ref{1.52} and
\ref{1.53} are simple consequences of the definition of
$\widetilde{\mathcal{E}}_{X}^{\ast}$ and $R_{x}$. That $\left\Vert
\widetilde{\mathcal{E}}_{X}^{\ast}\widetilde{\mathcal{E}}_{X}\right\Vert
=\left\Vert \widetilde{\mathcal{E}}_{X}\right\Vert ^{2}$ is an elementary
Hilbert space result and part 3 now implies $\left\Vert \widetilde{\mathcal{E}%
}_{X}^{\ast}\widetilde{\mathcal{E}}_{X}\right\Vert =\left\Vert R_{X,X}%
\right\Vert ^{2}$.

From Definition \ref{Def_Wg,x_m_gt_0} we know that the functions $\left\{
R_{x^{\left(  i\right)  }}\right\}  _{i=1}^{N}$ form a basis for $W_{G,X}$,
where the $x^{\left(  i\right)  }$ are the unique points of $X$. Hence
$\operatorname*{range}\widetilde{\mathcal{E}}_{X}^{\ast}=W_{G,X}$ and
$\operatorname*{null}\widetilde{\mathcal{E}}_{X}^{\ast}=\left\{  0\right\}  $.
Now recall the closed-range theorem, see for example Yosida \cite{Yosida58},
p205, which states that for a continuous linear operator $\mathcal{V}$, the
range of $\mathcal{V}$ is closed iff the range of $\mathcal{V}^{\ast}$ is
closed. Since $\operatorname*{range}\widetilde{\mathcal{E}}_{X}^{\ast}%
=W_{G,X}$ and $W_{G,X}$ is finite dimensional, $\operatorname*{range}%
\widetilde{\mathcal{E}}_{X}^{\ast}$ is closed and so $\operatorname*{range}%
\widetilde{\mathcal{E}}_{X}$ is closed. Consequently, using the result that
$\overline{\operatorname*{range}\mathcal{V}}=\left(  \operatorname*{null}%
\mathcal{V}^{\ast}\right)  ^{\perp}$ for any continuous linear operator
$\mathcal{V}$, it follows that
\[
\operatorname*{range}\widetilde{\mathcal{E}}_{X}=\overline
{\operatorname*{range}\widetilde{\mathcal{E}}_{X}}=\left(
\operatorname*{null}\widetilde{\mathcal{E}}_{X}^{\ast}\right)  ^{\perp
}=\left\{  0\right\}  ^{\perp}=\mathbb{C}^{N}.
\]

\textbf{Part 5}
\[
\widetilde{\mathcal{E}}_{X}\widetilde{\mathcal{E}}_{X}^{\ast}\beta
=\widetilde{\mathcal{E}}_{X}\sum_{i=1}^{N}\beta_{i}R_{x^{\left(  i\right)  }%
}=\sum_{i=1}^{N}\widetilde{\mathcal{E}}_{X}R_{x^{\left(  i\right)  }}\beta
_{i}=\left(  \widetilde{\mathcal{E}}_{X}R_{x^{\left(  1\right)  }}%
,\ldots,\widetilde{\mathcal{E}}_{X}R_{x^{\left(  N\right)  }}\right)
\beta=R_{X,X}\beta,
\]

since $R_{X,X}=\left(  R_{x^{\left(  j\right)  }}\left(  x^{\left(  i\right)
}\right)  \right)  $.\medskip

\textbf{Part 6} The proof is very similar to that of part 5.
\end{proof}

\section{The minimal norm interpolation problem and its solution}

In this subsection we will formulate and solve the \textit{minimal norm
interpolation problem\ }for scattered data. The solution to this problem will
be shown to be unique and its form derived.

Before defining the interpolation problem we will introduce some notation for
the data. We assume the data is \textit{scattered} i.e. irregularly spaced, as
distinct from on a regular rectangular grid. The data is a set of distinct
points $\left\{  \left(  x^{(i)},y_{i}\right)  \right\}  _{i=1}^{N}$ with
$x^{(i)}\in\mathbb{R}^{d}$ and $y_{i}\in\mathbb{C}$. We assume the $x^{(i)}$
must be distinct and set $X=\left\{  x^{(i)}\right\}  _{i=1}^{N}$ and
$y=\left\{  y_{i}\right\}  _{i=1}^{N}$. In this document $X$ will be called
the \textit{independent data} and $y$ the \textit{dependent data}.\medskip%

\begin{equation}%
\begin{tabular}
[c]{|l|}\hline
\textbf{The minimal norm interpolation problem}\\\hline
Suppose $\left\{  \left(  x^{(i)},y_{i}\right)  \right\}  _{i=1}^{N}$ is
scattered data, where the $x^{(i)}$ are distinct.\\
We say $u_{I}\in X_{w}^{0}$ is a solution of the minimal norm interpolation\\
problem if it interpolates the data and satisfies $\left\Vert u_{I}\right\Vert
_{w,0}\leq\left\Vert u\right\Vert _{w,0}$\\
for any other interpolant $u\in X_{w}^{0}$.\\\hline
\end{tabular}
\label{1.42}%
\end{equation}
\medskip

Using Hilbert space techniques we will now do the following:

\begin{enumerate}
\item Show there exists a unique minimal norm interpolant.

\item Show the interpolant is a basis function interpolant i.e. show that it
lies in $W_{G,X}$.

\item Construct a matrix equation for the coefficients of the
\textit{data-translated }basis functions $G\left(  \cdot-x^{\left(  k\right)
}\right)  $.
\end{enumerate}

The next theorem proves the minimal norm interpolation problem has a unique solution.

\begin{theorem}
\label{Thm_var_interpol_norm}For given data there exists a unique minimal norm
interpolant $u_{I}$. If $v$ is any other interpolant satisfying
$\widetilde{\mathcal{E}}_{X}v=y$ then
\begin{equation}
\left\Vert u_{I}\right\Vert _{w,0}^{2}+\left\Vert v-u_{I}\right\Vert
_{w,0}^{2}=\left\Vert v\right\Vert _{w,0}^{2},\label{1.21}%
\end{equation}

or equivalently
\begin{equation}
\left(  v-u_{I},u_{I}\right)  _{w,0}=0.\label{1.22}%
\end{equation}

\end{theorem}

\begin{proof}
Since the evaluation operator $\widetilde{\mathcal{E}}_{X}:X_{w}%
^{0}\rightarrow\mathbb{C}^{N}$ is continuous and onto, and the singleton set
$\left\{  y\right\}  $ is closed, it follows that the set
\[
\left\{  u:\widetilde{\mathcal{E}}_{X}u=y\right\}  =\widetilde{\mathcal{E}%
}_{X}^{-1}\left(  y\right)  ,
\]

is a non-empty, proper, closed subspace of the Hilbert space $X_{w}^{0}$.
Hence this subspace contains a unique element of smallest norm, say $u_{I}$.
If $\widetilde{\mathcal{E}}_{X}v=y$ then
\[
\left\{  u:u\in X_{w}^{0}\text{ and }\widetilde{\mathcal{E}}_{X}u=y\right\}
=\left\{  v-s:s\in\operatorname*{null}\widetilde{\mathcal{E}}_{X}\right\}  .
\]

Now
\[
\min\left\{  \left\Vert v-s\right\Vert _{w,0}:s\in\operatorname*{null}%
\widetilde{\mathcal{E}}_{X}\right\}  =\operatorname*{dist}\left(
v,\operatorname*{null}\widetilde{\mathcal{E}}_{X}\right)  ,
\]

is the distance between $v$ and the closed subspace $\operatorname*{null}%
\widetilde{\mathcal{E}}_{X}$. Therefore there exists a unique $s_{I}%
\in\operatorname*{null}\widetilde{\mathcal{E}}_{X}$ such that
\begin{equation}
u_{I}=v-s_{I},\label{1.23}%
\end{equation}

\[
\left\Vert v-s_{I}\right\Vert _{w,0}=\min\left\{  \left\Vert v-s\right\Vert
_{w,0}:s\in\operatorname*{null}\widetilde{\mathcal{E}}_{X}\right\}  ,
\]

and
\begin{equation}
\left\Vert s_{I}\right\Vert _{w,0}^{2}+\left\Vert v-s_{I}\right\Vert
_{w,0}^{2}=\left\Vert v\right\Vert _{w,0}^{2}.\label{1.24}%
\end{equation}

Substituting for $s_{I}$ in \ref{1.24} using \ref{1.23} yields \ref{1.21}.
Equation \ref{1.22} follows since it is a necessary and sufficient condition
for \ref{1.21} to be true.
\end{proof}

Next we prove that the minimal norm interpolant lies in the space $W_{G,X}$ of
Definition \ref{Def_Wg,x_m_gt_0}.

\begin{theorem}
\label{Thm_min_norm_interpol_in_Wgx}For each data vector $y\in\mathbb{C}^{N} $
the minimal norm interpolant $u_{I}$ is given by%
\begin{equation}
u_{I}=\widetilde{\mathcal{E}}_{X}^{\ast}R_{X,X}^{-1}y=\left(  2\pi\right)
^{d/2}\widetilde{\mathcal{E}}_{X}^{\ast}\left(  G_{X,X}\right)  ^{-1}%
y,\label{1.216}%
\end{equation}

where $\widetilde{\mathcal{E}}_{X}^{\ast}R_{X,X}^{-1}:\mathbb{C}%
^{N}\rightarrow W_{G,X}$ is an isomorphism. We also have%
\begin{equation}
u_{I}\left(  z\right)  =\left(  \widetilde{\mathcal{E}}_{X}\overline{R_{z}%
}\right)  ^{T}R_{X,X}^{-1}y.\label{1.215}%
\end{equation}

\end{theorem}

\begin{proof}
From part 5 of Theorem \ref{Thm_eval_op_properties}, $\widetilde{\mathcal{E}%
}_{X}\widetilde{\mathcal{E}}_{X}^{\ast}=R_{X,X}$. Now by Theorem
\ref{Thm_Rx,x_properties} $R_{X,X}$ is regular but in general $R_{X,X}\neq I$
so in general $\widetilde{\mathcal{E}}_{X}\widetilde{\mathcal{E}}_{X}^{\ast
}y\neq y$ and $\widetilde{\mathcal{E}}_{X}^{\ast}y$ is not an interpolant.
However, $\widetilde{\mathcal{E}}_{X}\left(  \widetilde{\mathcal{E}}_{X}%
^{\ast}R_{X,X}^{-1}\right)  =I$, so $\widetilde{\mathcal{E}}_{X}^{\ast}%
R_{X,X}^{-1}y$ is an interpolant in $W_{G,X}$. For convenience set
$u=\widetilde{\mathcal{E}}_{X}^{\ast}R_{X,X}^{-1}y$. We want to show that
$u=u_{I}$. But
\[
\left(  u,u_{I}-u\right)  _{w,0}=\left(  \widetilde{\mathcal{E}}_{X}^{\ast
}R_{X,X}^{-1}y,u_{I}-u\right)  _{w,0}=\left(  R_{X,X}^{-1}%
y,\widetilde{\mathcal{E}}_{X}\left(  u_{I}-u\right)  \right)  _{w,0}=0,
\]

since $u_{I}$ and $u$ are both interpolants to $y$. Thus
\[
\left\Vert u_{I}\right\Vert _{w,0}^{2}=\left\Vert u_{I}-u+u\right\Vert
_{w,0}^{2}=\left\Vert u_{I}-u\right\Vert _{w,0}^{2}+\left\Vert u\right\Vert
_{w,0}^{2},
\]

since $\left(  u,u_{I}-u\right)  _{w,0}=0$. If $u\neq u_{I}$ then $\left\Vert
u\right\Vert _{w,0}<\left\Vert u_{I}\right\Vert _{w,0}$, contradicting the
fact that $u_{I}$ is the minimal interpolant. Since $R_{X,X}$ is regular,
$R_{X,X}^{-1}$ is an isomorphism from $\mathbb{C}^{N}$ to $\mathbb{C}^{N}$ and
because $\widetilde{\mathcal{E}}_{X}^{\ast}$ is an isomorphism from
$\mathbb{C}^{N}$ to $W_{G,X}$ we have that $\widetilde{\mathcal{E}}_{X}^{\ast
}R_{X,X}^{-1}:\mathbb{C}^{N}\rightarrow W_{G,X}$ is an isomorphism.

Finally, using the definition of the adjoint $\widetilde{\mathcal{E}}%
_{X}^{\ast}$
\[
u_{I}\left(  z\right)  =\left(  u_{I},R_{z}\right)  _{w,0}=\left(
\widetilde{\mathcal{E}}_{X}^{\ast}R_{X,X}^{-1}y,R_{z}\right)  _{w,0}=\left(
R_{X,X}^{-1}y,\widetilde{\mathcal{E}}_{X}R_{z}\right)  =\left(
\widetilde{\mathcal{E}}_{X}\overline{R_{z}}\right)  ^{T}R_{X,X}^{-1}y.
\]

\end{proof}

The last theorem allows us to define the mapping between a data function and
its corresponding interpolant.

\begin{definition}
\label{Def_data_func_interpol_map}\textbf{Data functions and the interpolant
mapping} $\mathcal{I}_{X}:X_{w}^{0}\rightarrow W_{G,X}$

Given an independent data set $X$, we shall assume that each member of
$X_{w}^{0}$ can act as a legitimate data function $f$ and generate the
dependent data vector $\widetilde{\mathcal{E}}_{X}f$.

Equation \ref{1.216} of Theorem \ref{Thm_min_norm_interpol_in_Wgx} enables us
to define the linear mapping $\mathcal{I}_{X}:X_{w}^{0}\rightarrow W_{G,X} $
from the data functions to the corresponding unique minimal norm interpolant
$\mathcal{I}_{X}f=u_{I}$ given by%
\begin{equation}
\mathcal{I}_{X}f=\widetilde{\mathcal{E}}_{X}^{\ast}R_{X,X}^{-1}%
\widetilde{\mathcal{E}}_{X}f,\quad f\in X_{w}^{0}.\label{1.08}%
\end{equation}

Since $\widetilde{\mathcal{E}}_{X}\mathcal{I}_{X}f=\widetilde{\mathcal{E}}%
_{X}f$, Theorem \ref{Thm_eval_op_properties} can be easily used to show that
$\mathcal{I}_{X}$ is a self-adjoint projection onto $W_{G,X}$ with null space
$W_{G,X}^{\bot}$, and that $\mathcal{I}_{X}f=f$ iff $f\in W_{G,X}$. Since $f$
interpolates the data it follows from Theorem \ref{Thm_var_interpol_norm} that%
\begin{equation}
\left\Vert f-\mathcal{I}_{X}f\right\Vert _{w,0}\leq\left\Vert f\right\Vert
_{w,0},\quad\left\Vert \mathcal{I}_{X}f\right\Vert _{w,0}\leq\left\Vert
f\right\Vert _{w,0},\quad f\in X_{w}^{0},\label{1.09}%
\end{equation}

i.e. $\mathcal{I}_{X}$ and $I-\mathcal{I}_{X}$ are contractions.
\end{definition}

We now know that the interpolant lies in $W_{G,X}$ and from Theorem
\ref{Thm_indep_G(x-xi)} we know that $\left\{  G\left(  \cdot-x^{\left(
k\right)  }\right)  \right\}  $ is a basis for $W_{G,X}$. The next step is to
derive a matrix equation for the coefficients of the $G\left(  \cdot
-x^{\left(  k\right)  }\right)  $:

\begin{theorem}
The space $W_{G,X}$ contains only one interpolant to any given independent
data vector $y\in\mathbb{C}^{N}$. This interpolant is the minimal norm
interpolant and is defined uniquely by
\begin{equation}
u_{I}\left(  x\right)  =\sum_{k=1}^{N}v_{k}G\left(  x-x^{(k)}\right)
,\label{1.271}%
\end{equation}

where the coefficient vector $v=\left(  v_{k}\right)  $ satisfies the regular
basis matrix equation
\begin{equation}
G_{X,X}v=y.\label{1.272}%
\end{equation}

\end{theorem}

\begin{proof}
From Theorem \ref{Thm_min_norm_interpol_in_Wgx} we know that for a given
independent data vector $y$, $W_{G,X}$ contains the minimal norm interpolant.
Suppose $W_{G,X}$ contains another interpolant $v_{I}$. Then
$\widetilde{\mathcal{E}}_{X}\left(  u_{I}-v_{I}\right)  =0$ and so
$u_{I}-v_{I}\in\operatorname*{null}\widetilde{\mathcal{E}}_{X}$. But from part
1 of Theorem \ref{Thm_eval_op_properties} $\operatorname*{null}%
\widetilde{\mathcal{E}}_{X}=W_{G,X}^{\bot}$. Hence $u_{I}-v_{I}\in
W_{G,X}^{\bot}$ and so $u_{I}=v_{I}$ since $W_{G,X}^{\bot}\cap W_{G,X}%
=\left\{  0\right\}  $.

Equation \ref{1.271} and the interpolation requirement $u_{I}\left(
x^{(k)}\right)  =y_{k}$ implies $G_{X,X}v=y$ and by Theorem
\ref{Thm_Rx,x_properties} $G_{X,X}$ is regular.
\end{proof}

\section{Error estimates - no Taylor series
expansions\label{Sect_interp_no_Taylor_converg}}

In this subsection we will prove several pointwise estimates of the rate of
convergence of the minimal norm interpolant $\mathcal{I}_{X}f$ to its data
function $f\in X_{w}^{0}$ as the independent data $X$ `fills' a bounded closed
data region $K$. Unlike the results given by Light and Wayne for positive
order in \cite{LightWayne95ErrEst} and \cite{LightWayne98PowFunc}, the results
of this section do not use "explicit" Lagrange interpolation theory and no
Taylor series expansions so the proofs are much simpler and the constants can
be calculated. Instead, we use the properties of the Riesz representer $R_{x}$
and this allows us to avoid assumptions about the boundary.

The order of convergence appears as the power of the radius of the largest
ball that can be fitted between the data points i.e. the radius of the
\textbf{largest spherical cavity}.

For want of better names we will just call our convergence
estimates\textbf{\ Type 1} and \textbf{Type 2}. In the case of Type 1 the
weight function has property W02 for some $\kappa\geq0$ and the basis function
$G$ satisfies an extra smoothness condition near the origin. For Type
2\textit{\ }estimates\textit{\ }the weight function has property W02 for some
$\kappa\geq1$ and we use pointwise estimates based on the Riesz representer.

Here we could say we are actually using order 1 unisolvent subsets of the data
because by part 3 of Theorem \ref{Thm_unisolv} every data point is a minimal
unisolvent point and the results of the next section are valid for $m=1$.

When deriving error estimates we hope that as the number of evaluation points
$X$ increases i.e. as the `minimum density' of the points increases, the
interpolant will converge uniformly pointwise to the data function $f$ on the
data region $K$. The minimum density of the independent data points is
measured using the expression
\begin{equation}
h_{X,K}=\sup\limits_{x\in K}\operatorname*{dist}\left(  x,X\right)
=\sup\limits_{x\in K}\min\limits_{x^{\left(  k\right)  }\in X}\left\vert
x-x^{\left(  k\right)  }\right\vert .\label{1.12}%
\end{equation}

Here $h_{X,K}$ is the maximum distance between a point in $K$ and a point in
$X$, or equivalently $h_{X,K}$ is the radius of the largest open ball centered
in $K$ which does not contain any of the data points. Its use clearly makes
sense intuitively but from a numerical point of view, in dimensions greater
than 1, the calculation of $h_{X,K}$ presents formidable difficulties.

Theorem \ref{Thm_ord0_Riesz_rep_W2} tells us that the representer of the
evaluation functional $f\rightarrow f\left(  x\right)  $ is

$R_{x}=\left(  2\pi\right)  ^{-d/2}G\left(  \cdot-x\right)  \in X_{w}^{0}$.
Thus for any data function $f\in X_{w}^{0}$
\[
f\left(  x\right)  -\left(  \mathcal{I}_{X}f\right)  \left(  x\right)
=\left(  f-\mathcal{I}_{X}f,R_{x}\right)  _{w,0},
\]

where $\mathcal{I}_{X}f$ is the minimal interpolant of the data function on
the independent data set $X$. Hence
\begin{align}
\left\vert f\left(  x\right)  -\left(  \mathcal{I}_{X}f\right)  \left(
x\right)  \right\vert =\left\vert \left(  f-\mathcal{I}_{X}f,R_{x}\right)
_{w,0}\right\vert  &  \leq\left\Vert f-\mathcal{I}_{X}f\right\Vert
_{w,0}\left\Vert R_{x}\right\Vert _{w,0}\nonumber\\
&  =\left(  2\pi\right)  ^{-d/4}\sqrt{G\left(  0\right)  }\left\Vert
f-\mathcal{I}_{X}f\right\Vert _{w,0},\label{1.2}%
\end{align}

where $G\left(  0\right)  $ is real and positive. The error information is
clearly now restricted to the expression $\left\Vert f-\mathcal{I}%
_{X}f\right\Vert _{w,0}$ and if the inequality
\begin{equation}
\left\Vert f-\mathcal{I}_{X}f\right\Vert _{w,0}\leq\left\Vert f\right\Vert
_{w,0},\label{1.1}%
\end{equation}

implied by equation \ref{1.21} is used the error information is completely
lost. We can retain some of this information as follows. Choose any
independent data point $x^{\left(  j\right)  }$. Then the interpolation
property implies
\begin{equation}
f\left(  x\right)  -\left(  \mathcal{I}_{X}f\right)  \left(  x\right)
=\left(  f-\mathcal{I}_{X}f\right)  \left(  x\right)  -\left(  f-\mathcal{I}%
_{X}f\right)  \left(  x^{\left(  j\right)  }\right)  .\label{1.95}%
\end{equation}

In the results proved below we will use the estimates for $\left\vert f\left(
x\right)  -f\left(  y\right)  \right\vert $ and $\left\vert G\left(  x\right)
-G\left(  y\right)  \right\vert $ derived in Section
\ref{Sect_loc_estim_fn_in_Xow_W2}. Before deriving the first error estimate I
will show that there exists a suitable sequence of independent data sets
$X^{\left(  k\right)  }$ with $h_{X^{\left(  k\right)  },K}$ tending to zero:

\begin{theorem}
\label{Thm_seq_data_regions}Suppose $K$ is a bounded closed infinite set
containing all the independent data sets. Then there exists a sequence of
independent data sets $X^{\left(  k\right)  }\subset K$ such that $X^{\left(
k\right)  }\subset X^{\left(  k+1\right)  }$ and $h_{X^{\left(  k\right)  }%
,K}\rightarrow0$ as $k\rightarrow\infty$.
\end{theorem}

\begin{proof}
For each $k=1,2,3,\ldots$ there exists a finite covering of $K$ by the balls

$\left\{  B\left(  a_{k}^{\left(  j\right)  };\frac{1}{k}\right)  \right\}
_{j=1}^{M_{k}}$. Construct $X^{\left(  1\right)  }$ by choosing points from
$K$ so that one point lies in each ball $B\left(  a_{k}^{\left(  j\right)
};1\right)  $. Construct $X^{\left(  k+1\right)  }$ by first choosing the
points $X^{\left(  k\right)  }$ and then at least one extra point so that
$X^{\left(  k+1\right)  }$ contains points from each ball $B\left(
a_{k+1}^{\left(  j\right)  };\frac{1}{k+1}\right)  $.

Then $x\in K\cap X^{\left(  k\right)  }$ implies $x\in B\left(  a_{k}^{\left(
j\right)  };\frac{1}{k}\right)  $ for some $j$ and hence $\operatorname*{dist}%
\left(  x,X^{\left(  k\right)  }\right)  \leq\frac{1}{k}$. Thus $h_{X^{\left(
k\right)  },K}=\sup\limits_{x\in K}\operatorname*{dist}\left(  x,X^{\left(
k\right)  }\right)  \leq\frac{1}{k}$ and $\lim\limits_{k\rightarrow\infty
}h_{X^{\left(  k\right)  },K}=0$.
\end{proof}

We start with some simple bounds for the pointwise error of the interpolant:

\begin{theorem}
\label{Thm_interp_error_const_bound}Suppose the weight function $w$ has basis
function $G$. Then the interpolant $\mathcal{I}_{X}f$ of any data function $f$
satisfies
\begin{equation}
\left\vert f\left(  x\right)  -\mathcal{I}_{X}f\left(  x\right)  \right\vert
\leq\sqrt{\left(  f-\mathcal{I}_{X}f,f\right)  _{w,0}}\sqrt{R_{0}\left(
0\right)  },\text{\quad}x\in\mathbb{R}^{d},\label{1.006}%
\end{equation}

where $\sqrt{\left(  f-\mathcal{I}_{X}f,f\right)  _{w,0}}\leq\left\Vert
f\right\Vert _{w,0}$ and $R_{0}\left(  0\right)  =\left(  2\pi\right)
^{-\frac{d}{2}}G\left(  0\right)  $.

Also%
\begin{equation}
\left\vert R_{y}\left(  x\right)  -\mathcal{I}_{X}R_{y}\left(  x\right)
\right\vert \leq R_{0}\left(  0\right)  ,\text{\quad}x,y\in\mathbb{R}%
^{d}.\label{1.56}%
\end{equation}

\end{theorem}

\begin{proof}
Since $\mathcal{I}_{X}$ is a self-adjoint projection w.r.t. $\left(
\cdot,\cdot\right)  _{w,0}$
\[
\left\vert f\left(  x\right)  -\left(  \mathcal{I}_{X}f\right)  \left(
x\right)  \right\vert =\left\vert \left(  f-\mathcal{I}_{X}f,R_{x}\right)
_{w,0}\right\vert \leq\left\Vert f-\mathcal{I}_{X}f\right\Vert _{w,0}%
\left\Vert R_{x}\right\Vert _{w,0}=\sqrt{\left(  f-\mathcal{I}_{X}f,f\right)
_{w,0}}\sqrt{R_{x}\left(  x\right)  }.
\]

But $R_{y}\left(  x\right)  =\left(  2\pi\right)  ^{-\frac{d}{2}}G\left(
x-y\right)  $ and inequalities \ref{1.09} imply $\left\Vert f-\mathcal{I}%
_{X}f\right\Vert _{w,0}\leq\left\Vert f\right\Vert _{w,0}$ so that

$\sqrt{\left(  f-\mathcal{I}_{X}f,f\right)  _{w,0}}\leq\left\Vert f\right\Vert
_{w,0}$ and $R_{x}\left(  x\right)  =R_{0}\left(  0\right)  =\left(
2\pi\right)  ^{-\frac{d}{2}}G\left(  0\right)  $. Now letting $f=R_{y}$ in
\ref{1.006} we get%
\[
\left\vert R_{y}\left(  x\right)  -\left(  \mathcal{I}_{X}R_{y}\right)
\left(  x\right)  \right\vert \leq\sqrt{\left(  R_{y}-\mathcal{I}_{X}%
R_{y},R_{y}\right)  _{w,0}}\sqrt{R_{0}\left(  0\right)  }=\sqrt{R_{y}\left(
y\right)  -\left(  \mathcal{I}_{X}R_{y}\right)  \left(  y\right)  }\sqrt
{R_{0}\left(  0\right)  },
\]

which implies that when $x=y$, $R_{y}\left(  y\right)  -\left(  \mathcal{I}%
_{X}R_{y}\right)  \left(  y\right)  \leq R_{0}\left(  0\right)  $ and so%
\[
\left\vert R_{y}\left(  x\right)  -\left(  \mathcal{I}_{X}R_{y}\right)
\left(  x\right)  \right\vert \leq R_{0}\left(  0\right)  ,
\]

as required.
\end{proof}

\subsection{Type 1 pointwise estimates (W02, $\kappa\geq0$%
)\label{SbSect_int_estim_Type1}}

Type 1 estimates place no \textit{a priori} restriction on $\kappa$ which
always satisfies $\kappa\geq0$. In the next theorem a smoothness condition is
applied to the basis function near the origin and this will allow a uniform
order of convergence estimate to be obtained for the interpolant in a closed
bounded infinite data region. The next result is based on part 2 of Corollary
\ref{Cor_|f(x)-f(y)|_inequal_1_W2}.

\begin{theorem}
\label{Thm_|f(x)-f(y)|_inequal_2}\textbf{Interpolant convergence} Suppose the
weight function $w$ has property W02 and that $G$ is the basis function
generated by $w$. Assume that for some $s>0$ and constants $C_{G},h_{G}>0$ the
basis function satisfies
\begin{equation}
G\left(  0\right)  -\operatorname{Re}G\left(  x\right)  \leq C_{G}\left\vert
x\right\vert ^{2s},\text{\quad}\left\vert x\right\vert <h_{G}.\label{1.07}%
\end{equation}

Let $\mathcal{I}_{X}$ be the minimal norm interpolant mapping with the
independent data set $X$ contained in the closed bounded infinite set $K$, and
let $k_{G}=\left(  2\pi\right)  ^{-\frac{d}{4}}\sqrt{2C_{G}}$. Then for any
data function $f\in X_{w}^{0}$ it follows that $\sqrt{\left(  f-\mathcal{I}%
_{X}f,f\right)  _{w,0}}\leq\left\Vert f\right\Vert _{w,0}$ and%
\begin{equation}
\left\vert f\left(  x\right)  -\mathcal{I}_{X}f\left(  x\right)  \right\vert
\leq k_{G}\sqrt{\left(  f-\mathcal{I}_{X}f,f\right)  _{w,0}}\left(
h_{X,K}\right)  ^{s},\text{\quad}x\in K,\label{1.25}%
\end{equation}

when $h_{X,K}=\sup\limits_{x\in K}\operatorname*{dist}\left(  x,X\right)
<h_{G}$ i.e. the order of convergence is at least $s$.
\end{theorem}

\begin{proof}
From Theorem \ref{Thm_interp_error_const_bound}, $\sqrt{\left(  f-\mathcal{I}%
_{X}f,f\right)  _{w,0}}\leq\left\Vert f\right\Vert _{w,0}$. Now fix $x\in K$
and let $X=\left\{  x^{\left(  j\right)  }\right\}  _{j=1}^{N}\subset K$ be an
independent data set. Using the fact that $\mathcal{I}_{X}f$ interpolates $f$
on $X$ we can apply part 2 of Corollary \ref{Cor_|f(x)-f(y)|_inequal_1_W2} to
obtain
\begin{align*}
\left\vert f\left(  x\right)  -\mathcal{I}_{X}f\left(  x\right)  \right\vert
& =\left\vert \left(  f-\mathcal{I}_{X}f\right)  \left(  x\right)  -\left(
f-\mathcal{I}_{X}f\right)  \left(  x^{\left(  j\right)  }\right)  \right\vert
\\
& \leq\left(  2\pi\right)  ^{-\frac{d}{4}}\sqrt{2}\left\Vert f-\mathcal{I}%
_{X}f\right\Vert _{w,0}\sqrt{G\left(  0\right)  -\operatorname{Re}G\left(
x-x^{\left(  j\right)  }\right)  }\\
& =\left(  2\pi\right)  ^{-\frac{d}{4}}\sqrt{2}\sqrt{\left(  f-\mathcal{I}%
_{X}f,f\right)  _{w,0}}\sqrt{G\left(  0\right)  -\operatorname{Re}G\left(
x-x^{\left(  j\right)  }\right)  },
\end{align*}

for all $j$, where the last step used the fact that $\mathcal{I}_{X}$ is a
self-adjoint projection. Then, noting that $\operatorname*{dist}\left(
x,X\right)  <h_{G}$, we can apply the upper bound \ref{1.07} and obtain
\begin{align*}
\left\vert f\left(  x\right)  -\mathcal{I}_{X}f\left(  x\right)  \right\vert
& \leq\left(  2\pi\right)  ^{-\frac{d}{4}}\sqrt{2}\sqrt{\left(  f-\mathcal{I}%
_{X}f,f\right)  _{w,0}}\sqrt{C_{G}\left\vert x-x^{\left(  j\right)
}\right\vert ^{2s}}\\
& =k_{G}\sqrt{\left(  f-\mathcal{I}_{X}f,f\right)  _{w,0}}\left\vert
x-x^{\left(  j\right)  }\right\vert ^{s},
\end{align*}

for all $j$ and so%
\begin{align*}
\left\vert f\left(  x\right)  -\mathcal{I}_{X}f\left(  x\right)  \right\vert
&  \leq k_{G}\sqrt{\left(  f-\mathcal{I}_{X}f,f\right)  _{w,0}}\left(
\operatorname*{dist}\left(  x,X\right)  \right)  ^{s}\\
&  \leq k_{G}\sqrt{\left(  f-\mathcal{I}_{X}f,f\right)  _{w,0}}\left(
\sup\limits_{x\in K}\operatorname*{dist}\left(  x,X\right)  \right)  ^{s}\\
&  =k_{G}\sqrt{\left(  f-\mathcal{I}_{X}f,f\right)  _{w,0}}\left(
h_{X,K}\right)  ^{s},
\end{align*}

where the last step used inequalities \ref{1.09}. Since $\sqrt{\left(
f-\mathcal{I}_{X}f,f\right)  _{w,0}}\leq\left\Vert f\right\Vert _{w,0}$ the
order of convergence is at least $s$.
\end{proof}

We now prove a \textbf{double convergence rate} for the interpolant of data
functions which are Riesz representers $R_{x^{\prime}}$. This predicts a
convergence order which is at least double that for an arbitrary data
function. It also predicts a global bound for the error.

\begin{corollary}
\label{Cor_Thm_|f(x)-f(y)|_inequal_2}Under the notation and assumptions of the
previous Theorem \ref{Thm_|f(x)-f(y)|_inequal_2} the Riesz representer data
functions $R_{x^{\prime}}$ satisfy
\[
\left\vert R_{x^{\prime}}\left(  x\right)  -\left(  \mathcal{I}_{X}%
R_{x^{\prime}}\right)  \left(  x\right)  \right\vert \leq\left(  k_{G}\right)
^{2}\left(  h_{X,K}\right)  ^{2s},\text{\quad}x\in K,\text{ }x^{\prime}%
\in\mathbb{R}^{d},
\]

\end{corollary}

when $h_{X,K}\leq h_{G}$ i.e. the order of convergence is at least $2s$ in
$h_{X,K}$.

\begin{proof}
The key to this result is the term $\sqrt{\left(  f,f-\mathcal{I}_{X}f\right)
_{w,0}}$ in the estimate \ref{1.25} proved in the last theorem. Substituting
$f=R_{x^{\prime}}$ in \ref{1.25} gives%
\begin{align*}
\left\vert R_{x^{\prime}}\left(  x\right)  -\left(  \mathcal{I}_{X}%
R_{x^{\prime}}\right)  \left(  x\right)  \right\vert  & \leq k_{G}%
\sqrt{\left(  R_{x^{\prime}}-\mathcal{I}_{X}R_{x^{\prime}},R_{x^{\prime}%
}\right)  _{w,0}}\left(  h_{X,K}\right)  ^{s}\\
& =k_{G}\sqrt{R_{x^{\prime}}\left(  x^{\prime}\right)  -\left(  \mathcal{I}%
_{X}R_{x^{\prime}}\right)  \left(  x^{\prime}\right)  }\left(  h_{X,K}\right)
^{s},
\end{align*}

so that%
\[
R_{x^{\prime}}\left(  x^{\prime}\right)  -\left(  \mathcal{I}_{X}R_{x^{\prime
}}\right)  \left(  x^{\prime}\right)  \leq k_{G}\sqrt{R_{x^{\prime}}\left(
x^{\prime}\right)  -\left(  \mathcal{I}_{X}R_{x^{\prime}}\right)  \left(
x^{\prime}\right)  }\left(  h_{X,K}\right)  ^{s},
\]

which implies%
\[
\left\vert R_{x^{\prime}}\left(  x\right)  -\left(  \mathcal{I}_{X}%
R_{x^{\prime}}\right)  \left(  x\right)  \right\vert \leq\left(  k_{G}\right)
^{2}\left(  h_{X,K}\right)  ^{2s},
\]

as required.
\end{proof}

\subsection{Examples: radial basis
functions\label{SbSect_int_examp_rad_basis_ord_converg}}

For the shifted thin-plate splines, Gaussian, Sobolev and extended B-splines
we will derive the interpolant error parameters $s$, $h_{G}$ and $C_{G}$
defined above in Theorem \ref{Thm_|f(x)-f(y)|_inequal_2}. To do this I will
use the results of Theorem \ref{Thm_G(0)minusG(x)_bnd_w_rad_W2} when
$\kappa\geq1$: namely if $w$ satisfies property W02 for $\kappa=1$ then%
\begin{equation}
\left\vert G\left(  0\right)  -\operatorname{Re}G\left(  x\right)  \right\vert
\leq C_{G}\left\vert x\right\vert ^{2},\quad x\in\mathbb{R}^{d},\label{1.14}%
\end{equation}

where%
\begin{align}
G\left(  x\right)   & =f\left(  r^{2}\right)  \text{ }implies\text{ }%
C_{G}=\left\Vert 2rf^{\prime\prime}+f^{\prime}\right\Vert _{\infty
},\label{1.331}\\
G\left(  x\right)   & =g\left(  r\right)  \text{\quad}implies\text{ }%
C_{G}=\frac{1}{2}\left\Vert g^{\prime\prime}\right\Vert _{\infty
}.\label{1.332}%
\end{align}

If $\kappa<1$ I will use the mean value theorem.

\begin{example}
\textbf{Shifted thin-plate splines }From \ref{SbSect_wt_func_examples} the
basis functions are given by%
\[
G\left(  x\right)  =\left(  -1\right)  ^{\left\lceil v\right\rceil }\left(
1+\left\vert x\right\vert ^{2}\right)  ^{v},\quad-d/2<v<0.
\]

The weight function has property W02 for all $\kappa\geq0$ so we can use the
estimates \ref{1.331} and \ref{1.332}. But $G$ depends on $\left\vert
x\right\vert ^{2}$ so \ref{1.331} is easier algebraically: $f\left(  r\right)
=\left(  -1\right)  ^{\left\lceil v\right\rceil }\left(  1+r\right)  ^{v}$,
$f^{\prime}\left(  r\right)  =\left(  -1\right)  ^{\left\lceil v\right\rceil
}v\left(  1+r\right)  ^{v-1}$, $f^{\prime\prime}\left(  r\right)  =\left(
-1\right)  ^{\left\lceil v\right\rceil }v\left(  v-1\right)  \left(
1+r\right)  ^{v-2}$ and $f^{\prime\prime\prime}\left(  r\right)  =\left(
-1\right)  ^{\left\lceil v\right\rceil }v\left(  v-1\right)  \left(
v-2\right)  \left(  1+r\right)  ^{v-3}$. For the function $\left\vert
2rf^{\prime\prime}+f^{\prime}\right\vert $ to take a maximum $3f^{\prime
\prime}\left(  r\right)  =-2rf^{\prime\prime\prime}\left(  r\right)  $ must be
satisfied i.e. $3v\left(  v-1\right)  \left(  1+r\right)  ^{v-2}=-2v\left(
v-1\right)  \left(  v-2\right)  \left(  1+r\right)  ^{v-3}$ which reduces to
$3\left(  1+r\right)  =-2\left(  v-2\right)  =2\left(  2-v\right)  $ so that
\begin{equation}
C_{G}=\left\vert \left(  2rf^{\prime\prime}+f^{\prime}\right)  \left(
r_{\max}\right)  \right\vert ,\quad where\text{ }f\left(  r\right)  =\left(
1+r\right)  ^{v}\text{ }and\text{ }r_{\max}=\frac{\left(  1-2v\right)  }%
{3}.\label{1.70}%
\end{equation}

Also $h_{G}=\infty$, $s=1$ and $k_{G}=\left(  2\pi\right)  ^{-\frac{d}{4}%
}\sqrt{2C_{G}}$. The order of convergence for an arbitrary data function is at
least $1$ and that for the Riesz representer is at least $2$.
\end{example}

\begin{example}
\textbf{Gaussian} The weight function has property W02 for all $\kappa\geq0 $.
It is easier algebraically to set $G\left(  x\right)  =f\left(  r^{2}\right)
$ where $f\left(  r\right)  =e^{-r}$ and use \ref{1.331}. For the function
$2rf^{\prime\prime}+f^{\prime}$ to take a maximum $3f^{\prime\prime}\left(
r\right)  =-2rf^{\prime\prime\prime}\left(  r\right)  $ must be satisfied i.e.
$3e^{-r}=2re^{-r}$ so that $r_{\max}=3/2$ and hence $C_{G}=\left(
2rf^{\prime\prime}+f^{\prime}\right)  \left(  r_{\max}\right)  =\left(
2re^{-r}-e^{-r}\right)  \left(  r_{\max}\right)  =2e^{-3/2}$:
\[
C_{G}=2e^{-3/2}.
\]

Also $h_{G}=\infty$, $s=1$ and $k_{G}=\left(  2\pi\right)  ^{-\frac{d}{4}%
}\sqrt{2C_{G}}$.

The order of convergence for an arbitrary data function is at least $1$ and
that for the Riesz representer is at least $2$.
\end{example}

\begin{example}
\textbf{Sobolev splines}\label{Ex_Type1_converg} For this case the basis
function is \ref{1.045} i.e.%
\begin{equation}
G\left(  x\right)  =\frac{1}{2^{v-1}\Gamma\left(  v\right)  }\widetilde{K}%
_{v-d/2}\left(  \left\vert x\right\vert \right)  ,\text{\quad}x\in
\mathbb{R}^{d},\label{1.82}%
\end{equation}

constrained by $v>d/2$. The weight function has property W02 for $0\leq
\kappa<v-d/2$. There are two cases: $v-d/2$ is a positive integer and $v-d/2 $
is a positive non-integer. Refer to Sections 3.1, 3.2 and 3.3 of Magnus et al.
\cite{MagOberSoni66} for the properties of $K_{v}$ and $\widetilde{K}_{v} $
given below: the smoothness results have been obtained from the infinite
series representations and at $r=0$, $\widetilde{K}_{v}\left(  r\right)  $ has
the same smoothness as $r^{2v}\ln r$.\medskip

\fbox{Case 1: $v-d/2=1,2,3,$\ldots} Here we are dealing with the modified
Bessel functions of integer order. Now $rK_{0}$ has properties%
\begin{gather*}
rK_{0}\in C^{\left(  0\right)  }\left(  \left[  0,\infty\right)  \right)  \cap
C^{\infty}\left(  0,\infty\right)  ;\text{\quad}\lim_{r\rightarrow0}%
rK_{0}\left(  r\right)  =0;\text{\quad}K_{0}\left(  r\right)  >0,\text{
}r>0;\\
\lim_{r\rightarrow\infty}rK_{0}\left(  r\right)  =0.
\end{gather*}

and numerical experiments lead to the \textit{hypothesis}: $rK_{0}\left(
r\right)  $ has a single turning (maximum) point:%
\[
0\leq rK_{0}\left(  r\right)  \leq r^{\prime}K_{0}\left(  r^{\prime}\right)
\lesssim0.4665\text{ when }0\leq r\leq r^{\prime}\leq\arg\min\rho K_{0}\left(
\rho\right)  \simeq0.595.
\]

In general, for $m=1,2,3,\ldots$, $\widetilde{K}_{m}=r^{m}K_{m}\left(
r\right)  $ satisfies%
\begin{gather*}
\widetilde{K}_{-m}=\widetilde{K}_{m};\text{\quad}\widetilde{K}_{m}\in
C^{\left(  2m-1\right)  }\left(  \left[  0,\infty\right)  \right)  \cap
C^{\infty}\left(  0,\infty\right)  ;\text{\quad}\widetilde{K}_{m}\left(
0\right)  =2^{m-1}\left(  m-1\right)  !;\\
\widetilde{K}_{m}\left(  r\right)  >0,\text{ }r>0,
\end{gather*}

the derivative of $\widetilde{K}_{m}$ satisfies%
\begin{gather*}
D\widetilde{K}_{m}\left(  r\right)  =-r\widetilde{K}_{m-1}\left(  r\right)
;\text{\quad}D\widetilde{K}_{m}\left(  0\right)  =0;\text{\quad}%
D\widetilde{K}_{m}\left(  r\right)  <0,\text{\quad}x>0;\\
\lim_{r\rightarrow\infty}D\widetilde{K}_{m}\left(  r\right)  =0,
\end{gather*}

and if $m\geq2$, the second derivative of $\widetilde{K}_{m}$ satisfies%
\begin{gather}
D^{2}\widetilde{K}_{m}\left(  r\right)  =-\widetilde{K}_{m-1}\left(  r\right)
+r^{2}\widetilde{K}_{m-2}\left(  r\right)  ;\text{\quad}D^{2}\widetilde{K}%
_{m}\left(  0\right)  =-2^{m-2}\left(  m-2\right)  !\text{ };\label{1.69}\\
\lim_{r\rightarrow\infty}D^{2}\widetilde{K}_{m}\left(  r\right)  =0.\nonumber
\end{gather}

Now we want to estimate $\left\vert G\left(  0\right)  -\operatorname{Re}%
G\left(  x\right)  \right\vert $:\smallskip

\underline{When $v-d/2=1$} we have $\kappa<1$ and so we cannot use \ref{1.14}.
However, by the mean value theorem%
\[
\left\vert \widetilde{K}_{1}\left(  r\right)  -\widetilde{K}_{1}\left(
0\right)  \right\vert \leq\max_{\rho\in\left[  0,r\right]  }\left\vert
D\widetilde{K}_{1}\left(  \rho\right)  \right\vert r=\max_{\rho\in\left[
0,r\right]  }\left(  \rho K_{0}\left(  \rho\right)  \right)  r\leq\left\Vert
\rho K_{0}\left(  \rho\right)  \right\Vert _{\infty}r,
\]

so that%
\begin{align}
G\left(  0\right)  -\operatorname{Re}G\left(  x\right)   & =\frac{1}%
{2^{v-1}\Gamma\left(  v\right)  }\left(  \widetilde{K}_{1}\left(  0\right)
-\widetilde{K}_{1}\left(  \left\vert x\right\vert \right)  \right) \nonumber\\
& \leq\frac{1}{2^{v-1}\Gamma\left(  v\right)  }\left\Vert \rho K_{0}\left(
\rho\right)  \right\Vert _{\infty}\left\vert x\right\vert ,\label{1.64}%
\end{align}

and consequently
\begin{equation}
if\text{ }v-d/2=1\text{ }then\text{ }s=1/2,\text{ }C_{G}=\frac{\left\Vert \rho
K_{0}\left(  \rho\right)  \right\Vert _{\infty}}{2^{v-1}\Gamma\left(
v\right)  },\text{ }h_{G}=\infty.\label{1.59}%
\end{equation}

\underline{On the other hand if $v-d/2=2,3,4,\ldots$} then $\kappa\geq1$ and
we can use the estimate \ref{1.14}. Indeed, since%
\begin{equation}
G\left(  x\right)  =g\left(  r\right)  =\frac{1}{2^{v-1}\Gamma\left(
v\right)  }\widetilde{K}_{n-d/2}\left(  r\right)  ,\label{1.61}%
\end{equation}

we have%
\begin{equation}
C_{G}=\frac{1}{2}\left\Vert g^{\prime\prime}\right\Vert _{\infty}=\frac
{1}{2^{v}\Gamma\left(  v\right)  }\left\Vert D^{2}\widetilde{K}_{v-d/2}%
\right\Vert _{\infty},\text{\quad}v-d/2=2,3,4,\ldots,\label{1.76}%
\end{equation}

and thus%
\begin{equation}
if\text{ }v-d/2=2,3,4,\ldots\text{ }then\text{ }s=1,\text{ }C_{G}\text{
}is\text{ }given\text{ }by\text{ \ref{1.76}},\text{ }h_{G}=\infty.\label{1.60}%
\end{equation}

\fbox{Case 2: $v-d/2>0$ and non-integer} Here we are deal with the modified
Bessel functions $K_{\mu}$ of positive non-integer order. We have%
\begin{align*}
\widetilde{K}_{-\mu}  & =\widetilde{K}_{\mu};\text{\quad}\widetilde{K}_{\mu
}\in C^{\left(  \left\lfloor 2\mu\right\rfloor \right)  }\left(  \left[
0,\infty\right)  \right)  \cap C^{\infty}\left(  0,\infty\right)
;\text{\quad}\widetilde{K}_{\mu}\left(  0\right)  =2^{\left\vert
\mu\right\vert -1}\Gamma\left(  \left\vert \mu\right\vert \right)  ;\\
\widetilde{K}_{\mu}\left(  r\right)   & >0,\text{ }r>0;\text{\quad}%
\lim_{r\rightarrow\infty}\widetilde{K}_{\mu}\left(  r\right)  =0;
\end{align*}

If $\frac{1}{2}\leq\mu<1$ then $\widetilde{K}_{\mu}\in C^{\left(  1\right)
}\left(  \left[  0,\infty\right)  \right)  $ and%
\[
D\widetilde{K}_{\mu}\left(  r\right)  =-r\widetilde{K}_{\mu-1}\left(
r\right)  ,
\]

so that%
\[
D\widetilde{K}_{\mu}\left(  0\right)  =0;\text{\quad}D\widetilde{K}_{\mu
}\left(  r\right)  <0,\text{\quad}r>0;\text{\quad}\lim_{r\rightarrow\infty
}D\widetilde{K}_{\mu}\left(  r\right)  =0.
\]

If $\mu\geq1$ then $\widetilde{K}_{\mu}\in C^{\left(  2\right)  }\left(
\left[  0,\infty\right)  \right)  $ and%
\[
D^{2}\widetilde{K}_{\mu}\left(  r\right)  =-\widetilde{K}_{\mu-1}\left(
r\right)  +r^{2}\widetilde{K}_{\mu-2}\left(  r\right)  ,
\]

so that%
\begin{equation}
D^{2}\widetilde{K}_{\mu}\left(  0\right)  =-\widetilde{K}_{\mu-1}\left(
0\right)  =-2^{\left\vert \mu-2\right\vert }\Gamma\left(  \left\vert
\mu-1\right\vert \right)  ,\label{1.83}%
\end{equation}

and%
\[
\lim_{r\rightarrow\infty}D^{2}\widetilde{K}_{\mu}\left(  r\right)  =0.
\]

Now we want to estimate $G\left(  0\right)  -\operatorname{Re}G\left(
x\right)  $.\smallskip

\underline{When $1/2\leq v-d/2<1$} we have $\kappa<1$ and so we cannot use
\ref{1.14}. However, if $1/2\leq\mu<1$ the mean value theorem implies%
\[
\left\vert \widetilde{K}_{\mu}\left(  r\right)  -\widetilde{K}_{\mu}\left(
0\right)  \right\vert \leq r\max_{\rho\in\left[  0,r\right]  }\left\vert
D\widetilde{K}_{\mu}\left(  \rho\right)  \right\vert \leq\left\Vert
\rho\widetilde{K}_{\mu-1}\left(  \rho\right)  \right\Vert _{\infty
}r=\left\Vert \rho^{\mu}K_{1-\mu}\left(  \rho\right)  \right\Vert _{\infty}r,
\]

so that%
\[
G\left(  0\right)  -\operatorname{Re}G\left(  x\right)  =\frac{1}%
{2^{v-1}\Gamma\left(  v\right)  }\left(  \widetilde{K}_{v-d/2}\left(
0\right)  -\widetilde{K}_{v-d/2}\left(  r\right)  \right)  \leq\frac{1}%
{2^{v}\Gamma\left(  v\right)  }\left\Vert \rho^{v-d/2}K_{1-\left(
v-d/2\right)  }\left(  \rho\right)  \right\Vert _{\infty}r,
\]

and%
\begin{equation}
C_{G}=\frac{1}{2^{v}\Gamma\left(  v\right)  }\left\Vert \rho^{v-d/2}%
K_{1-\left(  v-d/2\right)  }\left(  \rho\right)  \right\Vert _{\infty}%
,\quad1/2\leq v-d/2<1.\label{1.75}%
\end{equation}

\underline{When\ $v-d/2>1,$ $v-d/2\notin\mathbb{Z}_{+}$} we have $\kappa\geq1$
and so we can use \ref{1.14} and \ref{1.332} to get%
\begin{equation}
C_{G}=\frac{1}{2}\left\Vert g^{\prime\prime}\right\Vert _{\infty}=\frac
{1}{2^{v}\Gamma\left(  v\right)  }\left\Vert D^{2}\widetilde{K}_{v-d/2}%
\right\Vert _{\infty},\quad v-d/2>1,\label{1.11}%
\end{equation}

and together with \ref{1.75} we can conclude that:
\begin{equation}%
\begin{array}
[c]{ll}%
s=1/2,\text{ }C_{G}\text{ }is\text{ }given\text{ }by\text{ }\ref{1.75},\text{
}h_{G}=\infty, & 1/2\leq v-d/2<1,\\
s=1,\text{ }C_{G}\text{ }is\text{ }given\text{ }by\text{ }\ref{1.11},\text{
}h_{G}=\infty, & v-d/2>1,\text{ }v-d/2\notin\mathbb{Z}_{+}.
\end{array}
\label{1.62}%
\end{equation}

\end{example}

\subsection{Examples: tensor product extended B-splines}

With reference to the error estimate of Theorem
\ref{Thm_|f(x)-f(y)|_inequal_2}:

\begin{corollary}
\label{Cor_converg_ext_nat_splin_wt_func_n=1}Suppose the weight function is a
tensor product extended B-spline weight function with parameters $n$ and $l $.
Then, if $G_{1}$ is the univariate basis function, the order of convergence
$s$ of the minimal norm interpolant to its data function is at least $1/2$ and
$C_{G}=\sqrt{d}G_{1}\left(  0\right)  ^{d-1}\left\Vert DG_{1}\right\Vert
_{\infty}$ and $h_{G}=\infty$.
\end{corollary}

\begin{proof}
Theorem \ref{Thm_ex_nat_spline_basis_Lipschitz} showed that
\[
\left\vert G_{s}\left(  x\right)  -G_{s}\left(  y\right)  \right\vert
\leq\sqrt{d}G_{1}\left(  0\right)  ^{d-1}\left\Vert DG_{1}\right\Vert
_{\infty}\left\vert x-y\right\vert ,\text{\quad}x,y\in\mathbb{R}^{d},
\]
so Theorem \ref{Thm_|f(x)-f(y)|_inequal_2} gives this result.
\end{proof}

\subsection{Examples: summary table}

Table \ref{Tbl_NonUnisolvTyp1Converg} summarizes the convergence results
\ref{1.59}, \ref{1.60} and \ref{1.62} for the Sobolev splines as predicted by
Theorem \ref{Thm_|f(x)-f(y)|_inequal_2}, as well as those for the other
\textit{Type 1} examples.%

\begin{table}[htbp] \centering
$%
\begin{tabular}
[c]{|c|c||c|c|c|}\hline
\multicolumn{5}{|c|}{\textbf{Type 1} interpolant error estimates}\\
\multicolumn{5}{|c|}{Smoothness condition on basis function near
origin.}\\\hline
& Parameter & Converg. &  & \\
Weight function & constraints & order $s$ & $C_{G}$ & $h_{G}$\\\hline\hline
\multicolumn{1}{|l|}{Sobolev splines} & $\frac{1}{2}<v-\frac{d}{2}\leq1$ &
$\frac{1}{2}$ & \multicolumn{1}{|l|}{$\frac{\left\Vert \rho^{v-\frac{d}{2}%
}K_{1-\left(  v-\frac{d}{2}\right)  }\right\Vert _{\infty}}{2^{v-1}%
\Gamma\left(  v\right)  }$} & $\infty$\\\cline{2-4}%
$\left(  v>d/2\right)  $ & $v-\frac{d}{2}>1$ & $1$ &
\multicolumn{1}{|l|}{$\frac{\left\Vert D^{2}\widetilde{K}_{v-\frac{d}{2}%
}\right\Vert _{\infty}}{2^{v}\Gamma\left(  v\right)  }$} & $\infty$\\\hline
\multicolumn{1}{|l|}{Shifted thin-plate} & - & $1$ & \multicolumn{1}{|l|}{eq.
(\ref{1.70})} & $\infty$\\
$\left(  -d/2<v<0\right)  $ &  &  &  & \\\hline
\multicolumn{1}{|l|}{Gaussian} & - & $1$ & \multicolumn{1}{|l|}{$2e^{-3/2}$} &
$\infty$\\\hline
\multicolumn{1}{|l|}{Extended B-spline} & - & $\frac{1}{2}$ &
\multicolumn{1}{|l|}{$G_{1}\left(  0\right)  ^{d-1}\left\Vert DG_{1}%
\right\Vert _{\infty}\sqrt{d}$ $^{\left(  1\right)  }$} & $\infty
$\\\hline\hline
\multicolumn{5}{|l|}{$^{\left(  1\right)  }${\small \ }$G_{1}${\small \ is the
univariate basis function used to form the tensor product.}}\\\hline
\end{tabular}
$\caption{}\label{Tbl_NonUnisolvTyp1Converg}%
\end{table}%

\subsection{Another approach - the interpolation error seminorm}

Here we will show how the definition of a seminorm using the interpolation
operator can be used to obtain the above results concerning interpolation in a
very concise manner. In Chapter \ref{Ch_Exact_smth} which studies a smoother
(a parameter stabilized interpolant) we will also make use of several
seminorms to study the smoother error.

\begin{definition}
\label{Def_interp_seminorm}\textbf{Interpolation error seminorm
and\ semi-inner product}

Suppose we selected the semi-inner product $\left(  f-\mathcal{I}%
_{X}f,g-\mathcal{I}_{X}g\right)  _{w,0}$ for pointwise estimation of
$f-\mathcal{I}_{X}f$. But $\mathcal{I}_{X}$ is a self-adjoint projection so
$\left(  f-\mathcal{I}_{X}f,g-\mathcal{I}_{X}g\right)  _{w,0}=\left(
f-\mathcal{I}_{X}f,g\right)  _{w,0}$ and this is more suited to pointwise
estimation since

$\left(  f-\mathcal{I}_{X}f,R_{x}\right)  _{w,0}=f\left(  x\right)  -\left(
\mathcal{I}_{X}f\right)  \left(  x\right)  $. So for $f,g\in X_{w}^{0}$ let us
define the interpolation error seminorm and\ semi-inner product by%
\[
\left\vert f\right\vert _{I}=\left(  f-\mathcal{I}_{X}f,f\right)  _{w,0}%
,\quad\left\langle f,g\right\rangle _{I}=\left(  f-\mathcal{I}_{X}f,g\right)
_{w,0}.
\]

To show that $\left\langle \cdot,\cdot\right\rangle _{I}$ is a semi-inner
product generated by $\left\vert \cdot\right\vert _{I}$ we first prove that
the norm satisfies the parallelogram law:%
\[
\left\vert f+g\right\vert _{I}^{2}+\left\vert f-g\right\vert _{I}^{2}=2\left(
\left\vert f\right\vert _{I}^{2}+\left\vert g\right\vert _{I}^{2}\right)  .
\]

It then follows that%
\begin{align}
\left\langle f,g\right\rangle _{I}  & =\frac{1}{4}\left(  \left\vert
f+g\right\vert _{I}^{2}-\left\vert f-g\right\vert _{I}^{2}\right)  +\frac
{i}{4}\left(  \left\vert f+ig\right\vert _{I}^{2}-\left\vert f-ig\right\vert
_{I}^{2}\right) \label{1.000}\\
& =\left(  f-\mathcal{I}_{X}f,g\right)  _{w,0},\nonumber
\end{align}

is a seminorm.
\end{definition}

\begin{remark}
It is interesting to note that \ref{1.000} can also be written%
\[
\left\langle f,g\right\rangle _{I}=\frac{1}{4}\sum\limits_{k=0}^{3}%
i^{k}\left\vert f+i^{k}g\right\vert _{I}^{2}.
\]

\end{remark}

The next theorem shows the close relationship between the pointwise error of
the interpolant and the interpolation seminorm.

\begin{theorem}
\label{Thm_property_interpol_seminorm}If $f,g\in X_{w}^{0}$ then some
properties of the interpolation seminorm are:

\begin{enumerate}
\item $\operatorname*{null}\left\vert \cdot\right\vert _{I}=W_{G,X}$.

\item $\left\vert f\right\vert _{I}^{2}=\left\Vert f\right\Vert _{w,0}%
^{2}-\left\Vert \mathcal{I}_{X}f\right\Vert _{w,0}^{2}$ and $\left\langle
f,g\right\rangle _{I}=\left(  f,g\right)  _{w,0}-\left(  \mathcal{I}%
_{X}f,g\right)  _{w,0}$.

\item $\left\langle f,g\right\rangle _{I}=\left(  f,g\right)  _{w,0}-\left(
\widetilde{\mathcal{E}}_{X}f\right)  ^{T}R_{X,X}^{-1}\widetilde{\mathcal{E}%
}_{X}g$.

\item $\left\langle f,R_{x}\right\rangle _{I}=f\left(  x\right)  -\left(
\mathcal{I}_{X}f\right)  \left(  x\right)  $.

\item $\left\vert R_{x}-R_{x^{\left(  k\right)  }}\right\vert _{I}=\left\vert
R_{x}\right\vert _{I}$\quad when $x^{\left(  k\right)  }\in X$.
\end{enumerate}
\end{theorem}

\begin{proof}
Use the comments in Definitions \ref{Def_data_func_interpol_map} and
\ref{Def_interp_seminorm}.
\end{proof}

The pointwise error of the interpolant has the following properties:

\begin{theorem}
\label{Thm_property_interpol}For all $x,y\in\mathbb{R}^{d}$:

\begin{enumerate}
\item $\left\vert R_{x}\left(  y\right)  -\left(  \mathcal{I}_{X}R_{x}\right)
\left(  y\right)  \right\vert ^{2}\leq\left(  R_{x}\left(  x\right)  -\left(
\mathcal{I}_{X}R_{x}\right)  \left(  x\right)  \right)  \left(  R_{y}\left(
y\right)  -\left(  \mathcal{I}_{X}R_{y}\right)  \left(  y\right)  \right)  .$

\item $R_{x}\left(  x\right)  -\left(  \mathcal{I}_{X}R_{x}\right)  \left(
x\right)  =\left\Vert R_{x}-R_{x^{\left(  k\right)  }}\right\Vert _{w,0}%
^{2}-\left\Vert \mathcal{I}_{X}\left(  R_{x}-R_{x^{\left(  k\right)  }%
}\right)  \right\Vert _{w,0}^{2},\quad x^{\left(  k\right)  }\in X.$
\end{enumerate}
\end{theorem}

\begin{proof}
\textbf{Part 1.} This is just the Cauchy-Schwartz inequality $\left\vert
\left\langle R_{x},R_{y}\right\rangle _{I}\right\vert \leq\left\vert
R_{x}\right\vert _{I}\left\vert R_{y}\right\vert _{I}$.\medskip

\textbf{Part 2.} From part 2 of Theorem \ref{Thm_property_interpol_seminorm}

$\left\vert R_{x}-R_{x^{\left(  k\right)  }}\right\vert _{I}=\left\Vert
R_{x}-R_{x_{k}}\right\Vert _{w,0}^{2}-\left\Vert \mathcal{I}_{X}\left(
R_{x}-R_{x_{k}}\right)  \right\Vert _{w,0}^{2}$ and from part 4 of Theorem
\ref{Thm_property_interpol_seminorm} $\left\vert R_{x}\right\vert _{I}%
=R_{x}\left(  x\right)  -\left(  \mathcal{I}_{X}R_{x}\right)  \left(
x\right)  $.
\end{proof}

From these two theorems and the calculation%
\begin{align}
\left\Vert R_{x}-R_{x^{\left(  k\right)  }}\right\Vert _{w,0}^{2}  & =\left(
R_{x}-R_{x^{\left(  k\right)  }},R_{x}-R_{x^{\left(  k\right)  }}\right)
_{w,0}\nonumber\\
& =R_{x}\left(  x\right)  -R_{x}\left(  x^{\left(  k\right)  }\right)
-\overline{R_{x}\left(  x^{\left(  k\right)  }\right)  }+R_{x^{\left(
k\right)  }}\left(  x^{\left(  k\right)  }\right) \nonumber\\
& =2\left(  2\pi\right)  ^{-d/2}\left(  G\left(  0\right)  -\operatorname{Re}%
G\left(  x-x^{\left(  k\right)  }\right)  \right)  ,\label{1.37}%
\end{align}

the interpolation results of Theorem \ref{Thm_|f(x)-f(y)|_inequal_2} and
Corollary \ref{Cor_Thm_|f(x)-f(y)|_inequal_2} can be obtained without much difficulty.

\subsection{Type 2 pointwise estimates $\left(  \kappa\geq1\right)
$\label{SbSect_estim_Type2}}

Type 2 pointwise estimates only assume $\kappa\geq1$. The Type 1 estimates of
Theorem \ref{Thm_|f(x)-f(y)|_inequal_2} considered the case where the weight
function had property W02 for some $\kappa\geq0$ and an extra condition was
applied to the basis function. In the next theorem we only assume that
$\kappa\geq1$ and derive an order of convergence estimate of $1 $, as well as
a doubled convergence estimate of $2$ for the Riesz data functions $R_{y}$.

Substituting $\alpha=\beta=0$ and $y=0$ in part 3 of Corollary
\ref{Cor_bound_DG(z-x)-DG(z-y)_W2} gives the Riesz representer-type estimate%
\[
G\left(  0\right)  -\operatorname{Re}G\left(  x\right)  \leq\tfrac{1}{\left(
2\pi\right)  ^{\frac{d}{2}}}\left(  \int\frac{d\xi}{w\left(  \xi\right)
}\right)  ^{\frac{1}{2}}\left(  \int\frac{\left\vert \xi\right\vert ^{2}d\xi
}{w\left(  \xi\right)  }\right)  ^{\frac{1}{2}}\left\vert x\right\vert ,
\]

since $G\left(  0\right)  $ is real and $\left\vert G\left(  x\right)
\right\vert \leq G\left(  0\right)  $. Hence $G$ always satisfies \ref{1.07}
with $s=1/2$ and the theory of Type 1 convergence estimates can be always be
applied to obtain convergence rates of at least $1/2$. Also, $\kappa\geq1$
implies $G\in C^{\left(  \left\lfloor 2\kappa\right\rfloor \right)  }\subseteq
C^{\left(  2\right)  }$ and the mean value theorem or the Taylor series
expansion can always be applied with $s=1/2$ or $s=1$ to obtain Type 1
convergence estimates.

However, instead we will start by using part 1 of Corollary
\ref{Cor_|f(x)-f(y)|_inequal_1_W2} - a Riesz representer-type estimate for
data functions - to obtain some error estimates. Note that the constant
$k_{G}$ in the next Theorem is defined to match the constant $k_{G}$ for the
Type 1 convergence estimate \ref{1.25}.

\begin{theorem}
\label{Thm_interpol_error_in_terms_of_wt_fn}\textbf{Interpolant convergence}
Suppose the weight function $w$ has property W02 for a parameter $\kappa\geq
1$. Suppose also that $\mathcal{I}_{X}f$ is the minimal interpolant on $X$ of
the data function $f\in X_{w}^{0}$ and that $K$ is a bounded closed infinite
subset of $\mathbb{R}^{d}$, $X\subset K$ and $h_{X,K}=\sup\limits_{x\in
K}\operatorname*{dist}\left(  x,X\right)  $.

Then%
\begin{equation}
\left\vert f\left(  x\right)  -\left(  \mathcal{I}_{X}f\right)  \left(
x\right)  \right\vert \leq k_{G}\sqrt{\left(  f-\mathcal{I}_{X}f,f\right)
_{w,0}}h_{X,K},\text{\quad}x\in K,\label{1.27}%
\end{equation}

where%
\begin{equation}
k_{G}^{2}=\left(  2\pi\right)  ^{-d}\int\frac{\left\vert \cdot\right\vert
^{2}}{w}=-\left(  2\pi\right)  ^{-\frac{d}{2}}\left(  \left\vert D\right\vert
^{2}G\right)  \left(  0\right)  ,\quad\sqrt{\left(  f-\mathcal{I}%
_{X}f,f\right)  _{w,0}}\leq\left\Vert f\right\Vert _{w,0},\label{1.79}%
\end{equation}

and the order of convergence is at least $1$. Further, for the `Riesz' data
functions $R_{y}$, $y\in K$, we have
\begin{equation}
\left\vert R_{y}\left(  x\right)  -\left(  \mathcal{I}_{X}R_{y}\right)
\left(  x\right)  \right\vert \leq\left(  k_{G}\right)  ^{2}\left(
h_{X,K}\right)  ^{2},\text{\quad}x,y\in K,\label{1.16}%
\end{equation}

i.e. a `doubled' estimate of at least $2$ for the rate of convergence.
\end{theorem}

\begin{proof}
Let $X=\left\{  x^{\left(  j\right)  }\right\}  _{j=1}^{N}$. Then since
$\mathcal{I}_{X}$ is an interpolant on $X$
\[
f\left(  x\right)  -\left(  \mathcal{I}_{X}f\right)  \left(  x\right)
=\left(  f-\mathcal{I}_{X}f\right)  \left(  x\right)  -\left(  f-\mathcal{I}%
_{X}f\right)  \left(  x^{\left(  j\right)  }\right)  ,
\]

for all $j$ and applying the inequality of part 1 of Corollary
\ref{Cor_|f(x)-f(y)|_inequal_1_W2} to $f-\mathcal{I}_{X}f\in X_{w}^{0}$ we
get
\begin{align}
\left\vert f\left(  x\right)  -\left(  \mathcal{I}_{X}f\right)  \left(
x\right)  \right\vert  & =\left\vert \left(  f-\mathcal{I}_{X}f\right)
\left(  x\right)  -\left(  f-\mathcal{I}_{X}f\right)  \left(  x^{\left(
j\right)  }\right)  \right\vert \nonumber\\
& \leq k_{G}\left\Vert f-\mathcal{I}_{X}f\right\Vert _{w,0}\left\vert
x-x^{\left(  j\right)  }\right\vert \nonumber\\
& =k_{G}\sqrt{\left(  f-\mathcal{I}_{X}f,f\right)  _{w,0}}\left\vert
x-x^{\left(  j\right)  }\right\vert ,\quad j=1,\ldots,N,\label{1.89}%
\end{align}

Now suppose that $X\subset K$. Then the inequalities \ref{1.89} imply that
\[
\left\vert f\left(  x\right)  -\left(  \mathcal{I}_{X}f\right)  \left(
x\right)  \right\vert \leq k_{G}\sqrt{\left(  f-\mathcal{I}_{X}f,f\right)
_{w,0}}\operatorname*{dist}\left(  x,X\right)  ,
\]

and so if $x$ is restricted to $K$%
\[
\left\vert f\left(  x\right)  -\left(  \mathcal{I}_{X}f\right)  \left(
x\right)  \right\vert \leq k_{G}\sqrt{\left(  f-\mathcal{I}_{X}f,f\right)
_{w,0}}\sup\limits_{x\in K}\operatorname*{dist}\left(  x,X\right)  =k_{G}%
\sqrt{\left(  f-\mathcal{I}_{X}f,f\right)  _{w,0}}h_{X,K}.
\]

Thus when $f=R_{y}$
\begin{align}
\left\vert R_{y}\left(  x\right)  -\left(  \mathcal{I}_{X}R_{y}\right)
\left(  x\right)  \right\vert  & \leq k_{G}\sqrt{\left(  R_{y}-\mathcal{I}%
_{X}R_{y},R_{y}\right)  _{w,0}}h_{X,K}\nonumber\\
& =k_{G}\sqrt{R_{y}\left(  y\right)  -\left(  \mathcal{I}_{X}R_{y}\right)
\left(  y\right)  }h_{X,K},\label{1.29}%
\end{align}

so that when $x=y$%
\[
R_{y}\left(  y\right)  -\left(  \mathcal{I}_{X}R_{y}\right)  \left(  y\right)
\leq k_{G}\sqrt{R_{y}\left(  y\right)  -\left(  \mathcal{I}_{X}R_{y}\right)
\left(  y\right)  }h_{X,K},
\]

and $R_{y}\left(  y\right)  -\left(  \mathcal{I}_{X}R_{y}\right)  \left(
y\right)  \leq\left(  k_{G}h_{X,K}\right)  ^{2}$. Inequality \ref{1.29} now
implies inequality \ref{1.16} as required.

Finally, that $\int\frac{\left\vert \cdot\right\vert ^{2}}{w}=-\left(
2\pi\right)  ^{d/2}\left(  \left\vert D\right\vert ^{2}G\right)  \left(
0\right)  $, follows directly from \ref{1.048}.
\end{proof}

\subsection{Examples: radial basis functions}

It was shown in Subsection \ref{SbSect_wt_func_examples} that the
\textbf{shifted thin-plate splines} and the \textbf{Gaussian} have weight
functions satisfying property W02\ for all $\kappa\geq1$. Thus we can apply
Theorem \ref{Thm_interpol_error_in_terms_of_wt_fn} to obtain the same orders
of convergence that were obtained for the Type 1 results in Subsection
\ref{SbSect_int_estim_Type1} i.e. a convergence of order $1$.

It was shown in Subsection \ref{SbSect_wt_func_examples} that $0\leq
\kappa<v-d/2$ for the \textbf{Sobolev splines}. Thus we can choose $\kappa
\geq1$ iff $v-d/2>1$.

From \ref{1.79}, $k_{G}^{2}=-\left(  2\pi\right)  ^{-d}\left(  \left\vert
D\right\vert ^{2}G\right)  \left(  0\right)  $. Now if $G\left(  x\right)
=f\left(  r^{2}\right)  $ then by part 2 of Theorem
\ref{Thm_G(0)minusG(x)_bound_w_W2}, $\left(  \left\vert D\right\vert
^{2}G\right)  \left(  0\right)  =2f^{\prime}\left(  0\right)  d$ and so%
\begin{equation}
k_{G}=\left(  2\pi\right)  ^{-d/4}\sqrt{-2f^{\prime}\left(  0\right)  }%
\sqrt{d},\label{1.80}%
\end{equation}

and if $G\left(  x\right)  =g\left(  r\right)  $ then by part 1 of Theorem
\ref{Thm_G(0)minusG(x)_bound_w_W2}, $\left(  \left\vert D\right\vert
^{2}G\right)  \left(  0\right)  =g^{\prime\prime}\left(  0\right)  d$ and so%
\begin{equation}
k_{G}=\left(  2\pi\right)  ^{-d/4}\sqrt{-g^{\prime\prime}\left(  0\right)
}\sqrt{d}.\label{1.81}%
\end{equation}

\begin{example}
\textbf{Shifted thin-plate splines} $f\left(  r\right)  =\left(  1+r\right)
^{v} $, $f^{\prime}\left(  r\right)  =v\left(  1+r\right)  ^{v-1}$ and so
\[
k_{G}=\left(  2\pi\right)  ^{-d/4}\sqrt{-2v}\sqrt{d}.
\]

\end{example}

\begin{example}
\textbf{Gaussian} $f\left(  r\right)  =e^{-r}$, $f\prime\left(  r\right)
=-e^{-r}$ and hence $k_{G}=\left(  2\pi\right)  ^{-d/4}\sqrt{2}\sqrt{d}$.
\end{example}

\begin{example}
\textbf{Sobolev splines} Here we use results of Subsubsection
\ref{Ex_Type1_converg}. From \ref{1.82},

$g\left(  r\right)  =\frac{1}{2^{v-1}\Gamma\left(  v\right)  }\widetilde{K}%
_{v-d/2}\left(  r\right)  $ and so $g^{\prime\prime}\left(  0\right)
=\frac{1}{2^{v-1}\Gamma\left(  v\right)  }D^{2}\widetilde{K}_{v-d/2}\left(
0\right)  $.\medskip

\underline{If $v-d/2$ is an integer} then $v-d/2=2,3,4,\ldots$ and by
\ref{1.69}
\[
g^{\prime\prime}\left(  0\right)  =\frac{D^{2}\widetilde{K}_{v-d/2}\left(
0\right)  }{2^{v-1}\Gamma\left(  v\right)  }=-\frac{2^{v-d/2-2}\left(
v-d/2-2\right)  !}{2^{v-1}\Gamma\left(  v\right)  }=-\frac{\left(
v-d/2-2\right)  !}{2^{d/2+1}\Gamma\left(  v\right)  },
\]

and by \ref{1.81}%
\begin{equation}
k_{G}=\left(  2\pi\right)  ^{-d/4}\sqrt{-g^{\prime\prime}\left(  0\right)
}\sqrt{d}=\left(  2\pi\right)  ^{-d/4}\sqrt{\frac{\left(  v-d/2-2\right)
!}{2^{d/2+1}\Gamma\left(  v\right)  }}\sqrt{d}.\label{1.84}%
\end{equation}

\underline{If $v-d/2$ is \textbf{not} an integer} then $v-d/2>1$ and
$v-d/2\neq2,3,4,\ldots$. By \ref{1.83}%
\begin{align*}
g^{\prime\prime}\left(  0\right)  =\frac{D^{2}\widetilde{K}_{v-d/2}\left(
0\right)  }{2^{v-1}\Gamma\left(  v\right)  }=-\frac{2^{\left\vert
v-d/2-2\right\vert }\Gamma\left(  \left\vert v-d/2-1\right\vert \right)
}{2^{v-1}\Gamma\left(  v\right)  } &  =-\frac{2^{\left\vert v-d/2-2\right\vert
}\Gamma\left(  v-d/2-1\right)  }{2^{v-1}\Gamma\left(  v\right)  }\\
&  =-\frac{\Gamma\left(  v-d/2-1\right)  }{2^{v-1-\left\vert
v-d/2-2\right\vert }\Gamma\left(  v\right)  },
\end{align*}

and%
\begin{align}
k_{G}  & =\left(  2\pi\right)  ^{-d/4}\sqrt{\frac{\Gamma\left(
v-d/2-1\right)  }{2^{v-1-\left\vert v-d/2-2\right\vert }\Gamma\left(
v\right)  }}\sqrt{d}\nonumber\\
& =\left\{
\begin{array}
[c]{ll}%
\left(  2\pi\right)  ^{-d/4}\sqrt{\frac{\Gamma\left(  v-d/2-1\right)
}{2^{2v-d/2-3}\Gamma\left(  v\right)  }}\sqrt{d}, &
\begin{array}
[c]{c}%
1<v-d/2<2,
\end{array}
\\
\left(  2\pi\right)  ^{-d/4}\sqrt{\frac{\Gamma\left(  v-d/2-1\right)
}{2^{d/2+1}\Gamma\left(  v\right)  }}\sqrt{d}, &
\begin{array}
[c]{l}%
v-d/2>2,\\
v-d/2\neq2,3,4,\ldots
\end{array}
\end{array}
\right. \label{1.85}%
\end{align}

Finally, \ref{1.84} and \ref{1.85} can be combined to give%
\[
k_{G}=\left\{
\begin{array}
[c]{ll}%
\left(  2\pi\right)  ^{-d/4}\sqrt{\frac{\Gamma\left(  v-d/2-1\right)
}{2^{2v-d/2-3}\Gamma\left(  v\right)  }}\sqrt{d}, & 1<v-d/2<2,\\
\left(  2\pi\right)  ^{-d/4}\sqrt{\frac{\Gamma\left(  v-d/2-1\right)
}{2^{d/2+1}\Gamma\left(  v\right)  }}\sqrt{d}, & v-d/2>2.
\end{array}
\right.
\]
\medskip
\end{example}

\subsection{Examples: tensor product basis functions\medskip}

\begin{example}
\textbf{Extended B-splines }$\left(  1\leq n\leq l\right)  $ Since the weight
function has property W02 for $\kappa$ iff $\kappa+1/2<n$ we can choose
$\kappa=1$ iff $n\geq2$. Then by Theorem
\ref{Thm_interpol_error_in_terms_of_wt_fn}, $k_{G}=\left(  2\pi\right)
^{-d/4}\sqrt{-\left(  \left\vert D\right\vert ^{2}G_{s}\right)  \left(
0\right)  }$ and since
\begin{align*}
\left(  \left\vert D\right\vert ^{2}G_{s}\right)  \left(  0\right)
=\sum\limits_{k=1}^{d}\left(  D_{k}^{2}G_{s}\right)  \left(  0\right)   &
=\sum\limits_{k=1}^{d}D_{k}^{2}\left(  G_{1}\left(  x_{1}\right)  G_{1}\left(
x_{2}\right)  \ldots G_{1}\left(  x_{d}\right)  \right)  \left(  0\right) \\
&  =\sum\limits_{k=1}^{d}G_{1}\left(  0\right)  ^{d-1}D^{2}G_{1}\left(
0\right) \\
&  =G_{1}\left(  0\right)  ^{d-1}D^{2}G_{1}\left(  0\right)  d,
\end{align*}

it follows that
\[
k_{G}=\left(  2\pi\right)  ^{-d/4}\sqrt{-G_{1}\left(  0\right)  ^{d-1}%
D^{2}G_{1}\left(  0\right)  }\sqrt{d}.
\]
\medskip
\end{example}

These convergence results are summarized in the Table
\ref{Tbl_InterpNonUnisolvTyp2Conv}.

\subsection{Examples: summary table}%

\begin{table}[htbp] \centering
$%
\begin{tabular}
[c]{|c|c||c|r|c|}\hline
\multicolumn{5}{|c|}{\textbf{Type 2} interpolant convergence order
estimates.}\\
\multicolumn{5}{|c|}{Only assume W02 for $\kappa\geq1$.}\\\hline
& Parameter & Converg. &  & \\
Weight function & constraints & order & \multicolumn{1}{|c|}{$\left(
2\pi\right)  ^{d/4}k_{G}$} & $h_{G}$\\\hline\hline
\multicolumn{1}{|l|}{Sobolev splines} & $v-\frac{d}{2}\geq2$ & $1$ &
$\sqrt{\frac{\Gamma\left(  v-d/2-1\right)  }{2^{d/2+1}\Gamma\left(  v\right)
}}\sqrt{d}$ & $\infty$\\\cline{2-5}\cline{2-4}%
$\left(  v>\frac{d}{2}\right)  $ & $1<v-\frac{d}{2}<2$ & $1$ & $\sqrt
{\frac{\Gamma\left(  v-d/2-1\right)  }{2^{2v-d/2-3}\Gamma\left(  v\right)  }%
}\sqrt{d}$ & $\infty$\\\hline
\multicolumn{1}{|l|}{Shifted thin-plate} & - & $1$ & $\sqrt{-2v}\sqrt{d}$ &
$\infty$\\
$\left(  -d/2<v<0\right)  $ &  &  &  & \\\hline
\multicolumn{1}{|l|}{Gaussian} & - & $1$ & $\sqrt{2}\sqrt{d}$ & $\infty
$\\\hline
Extended B-spline & $n\geq2$ & $1$ & $\sqrt{-G_{1}\left(  0\right)
^{d-1}D^{2}G_{1}\left(  0\right)  }\sqrt{d}$ $^{\left(  1\right)  }$ &
$\infty$\\
$\left(  1\leq n\leq l\right)  $ &  &  &  & \\\hline\hline
\multicolumn{5}{|l|}{$^{\left(  1\right)  }${\small \ }$G_{1}${\small \ is the
univariate basis function used to form the tensor product.}}\\\hline
\end{tabular}
$\caption{}\label{Tbl_InterpNonUnisolvTyp2Conv}%
\end{table}%

\section{Error estimates using unisolvent data subsets $\left(  \kappa
\geq1\right)  $\label{Sect_unisolv}}

Multipoint Taylor series expansions and Lagrange interpolation using minimal
unisolvent sets of independent data points is fundamental to the theory of
positive order basis function interpolation \cite{WilliamsPosOrdSmthV3}. In
that case the order of the unisolvent set is the order of the basis function.
We will now show that this technique can be used when the order is zero. Here,
however, the order of unisolvency is related to the value of the weight
function parameter $\kappa$.

For the case $1\leq\kappa\leq2$ it is better to use Type 1 or 2 estimates
because theory is much simpler, the constants can be calculated and the order
of convergence obtained is the same.

In Theorem \ref{Thm_interpol_error_in_terms_of_wt_fn} we assumed the weight
function had property W02 for some $\kappa\geq1$ and so obtained an order $1$
estimate for the convergence of the interpolant. In this subsection we will
show that by assuming the data is a unisolvent set of points of order
$\left\lfloor \kappa\right\rfloor $ an order $\left\lfloor \kappa\right\rfloor
$ convergence estimate can be obtained for and arbitrary data function.

\begin{definition}
\label{Def_unisolv}\textbf{Unisolvent sets and minimal unisolvent sets}

Recall that $P_{m}$ is the set of polynomials of \textbf{order} $m$ i.e. of
degree $m-1$ when $m\geq1$. Note that order=1+degree.

Then a finite set of distinct points $X=\left\{  x_{i}\right\}  $ is said to
be a \textbf{unisolvent set} with respect to $P_{m}$ if: $p\in P_{m}$ and
$p\left(  x_{i}\right)  =0$ for all $x_{i}\in X$ implies $p=0$.

Sometimes we say $X$ is unisolvent of order $m$ or that $X$ is $m$%
-\textbf{unisolvent}.

It is known that any unisolvent set has at least $M=\dim P_{m}$ points, and
that any unisolvent set of more than $M$ points has a unisolvent subset with
$M$ points. Consequently, a unisolvent set with $M$ points is called a
\textbf{minimal unisolvent set}.

The polynomial $p\left(  u\right)  =\sum\limits_{\left\vert \alpha\right\vert
\leq m-1}u^{\alpha}$ has order $m$, and from the identity \ref{a1.11} in the
Appendix:
\begin{equation}
M=\dim P_{m}=p\left(  \mathbf{1}\right)  =\sum\limits_{\left\vert
\alpha\right\vert \leq m-1}1=\binom{d+m-1}{d},\quad m\geq1.\label{a1.52}%
\end{equation}

\end{definition}

\begin{theorem}
\label{Thm_unisolv}Suppose $M=\dim P_{m}$ and $A=\left\{  a_{i}\right\}
_{i=1}^{M}\subset\mathbb{R}^{d}$. Then:

\begin{enumerate}
\item If $\left\{  p_{j}\right\}  _{j=1}^{M}$ is any basis of $P_{m}$ then $A$
is minimally $m$-unisolvent iff $\det\left(  p_{j}\left(  a_{i}\right)
\right)  \neq0$.

\item In one dimension any set of $m$ distinct points is $m$-unisolvent.

\item In any dimension $m=1$ implies $M=1$ and $P_{1}$ consists of the
constant polynomials. The minimal unisolvent sets consist of\textbf{\ }single points.
\end{enumerate}
\end{theorem}

\begin{proof}
\textbf{Part 1} follows directly from the definition of unisolvency.
\textbf{Part 2} Choose the monomial basis $\left\{  x^{k}\right\}
_{k=0}^{m-1}$ for $P_{m}$. Then the determinant of part 1 is Vandermonde's
alternant which is singular iff two points coincide. \textbf{Part 3} From
\ref{a1.52}: $M=\dim P_{m}=\binom{d+1-1}{d}=1$.
\end{proof}

\begin{definition}
\label{Def_cardinal_basis}\textbf{Cardinal basis for polynomials }$P_{m}$.

A basis $\left\{  l_{i}\right\}  _{i=1}^{M}$ for $P_{m}$ is a \textbf{cardinal
basis} for the minimal unisolvent set $A=\left\{  a_{i}\right\}  _{i=1}^{M}$
if $l_{i}\left(  a_{j}\right)  =\delta_{i,j}$ and the $l_{i}$ are polynomials
with \textbf{real valued} coefficients.

It is known that a set is minimally unisolvent iff there exists a (unique)
cardinal basis for the set.
\end{definition}

The next step is to introduce the Lagrange interpolation operator
$\mathcal{P}$ and the operator $\mathcal{Q}=I-\mathcal{P}$:

\begin{definition}
\label{Def_Aux_proj_operator}\textbf{The unisolvency operators }%
$\mathcal{P}:C^{\left(  0\right)  }\rightarrow P_{m}$ \textbf{and}
$\mathcal{Q}:C^{\left(  0\right)  }\rightarrow C^{\left(  0\right)  }$.

These operators are only defined for integers $m\geq1$. Suppose the set
$A=\left\{  a_{i}\right\}  _{i=1}^{M}$ is a minimal unisolvent set with
respect to the polynomials $P_{m}$ and by Definition \ref{Def_cardinal_basis}
there is a unique cardinal basis $\left\{  l_{i}\right\}  _{i=1}^{M}$ for $A$.
Then for any continuous function $f$ the operators $\mathcal{P}$ and
$\mathcal{Q}$ are defined by
\[
\mathcal{P}f=\sum\limits_{i=1}^{M}f\left(  a_{i}\right)  l_{i},\qquad
\mathcal{Q}=I-\mathcal{P}.
\]

\end{definition}

\begin{theorem}
\label{Thm_P_op_properties}\textbf{The operators }$\mathcal{P}$ \textbf{and}
$\mathcal{Q}$ have the following elementary properties:

\begin{enumerate}
\item $p\in P_{m}$ implies $\mathcal{P}p=p$ and hence $\mathcal{Q}p=0$.

\item $\mathcal{P}f$ interpolates the data $\left\{  \text{ }\left(
a_{i},f\left(  a_{i}\right)  \right)  \,\right\}  _{i=1}^{M}$. Indeed,
$\mathcal{P}$ is termed the \textbf{Lagrange polynomial interpolation
function}.

\item $\mathcal{P}^{2}=\mathcal{P}$, $\mathcal{PQ}=\mathcal{QP}=0$ and
$\mathcal{Q}^{2}=\mathcal{Q}$ so the operators are projections.

\item $\operatorname*{null}\mathcal{Q}=\operatorname*{range}\mathcal{P}$.
\end{enumerate}
\end{theorem}

\begin{proof}
\textbf{Part 1} is true since each member of the cardinal basis satisfies
$\mathcal{P}l_{j}=\sum\limits_{i=1}^{M}l_{i}\left(  a_{i}\right)  l_{i}=l_{j}%
$. \textbf{Part 2} is true because $l_{i}\left(  a_{j}\right)  =\delta_{i,j}$.
Regarding \textbf{part 3}, $\mathcal{P}$ is a projection since by part 2
\[
\mathcal{P}^{2}f=\mathcal{P}\left(  \sum\limits_{i=1}^{M}f\left(
a_{i}\right)  l_{i}\right)  =\sum\limits_{\left\vert \alpha\right\vert
<k}f\left(  a_{i}\right)  \mathcal{P}\left(  l_{i}\right)  =\sum
\limits_{i=1}^{M}f\left(  a_{i}\right)  l_{i}=\mathcal{P}f,
\]

and so
\[
\mathcal{PQ}=\mathcal{P}\left(  I-\mathcal{P}\right)  =\mathcal{P}%
-\mathcal{P}^{2}=0=\left(  I-\mathcal{P}\right)  \mathcal{P}=\mathcal{QP},
\]
which implies $\mathcal{Q}^{2}=\mathcal{Q}\left(  I-\mathcal{P}\right)
=\mathcal{Q}$ and we have proved part 3. Finally $\mathcal{P}+\mathcal{Q}=I$
so \textbf{part 4} holds.
\end{proof}

To study the pointwise order of convergence of the minimal interpolant to its
data function we will need the following two lemmas. The first lemma provides
a pointwise upper bound for the operator $\mathcal{Q}$ using a
\textbf{multipoint Taylor series expansion} which will express $\mathcal{Q}%
f\left(  x\right)  =u\left(  x\right)  -\mathcal{P}u\left(  x\right)  $ in
terms of the derivatives of order $m$ evaluated on the intervals $\left[
x,a_{i}\right]  $. The basic tool is the Taylor series expansion with integral remainder.

\begin{lemma}
\label{Lem_Q_estim}Suppose $f\in C^{\left(  m\right)  }\left(  \mathbb{R}%
^{d}\right)  $ and $A=\left\{  a_{i}\right\}  _{i=1}^{M}$ is a minimal
unisolvent set of order $m\geq1$. Then we have the upper bound%
\begin{equation}
\left\vert \mathcal{Q}f\left(  x\right)  \right\vert \leq\frac{d^{\frac{m}{2}%
}}{m!}\left(  \sum_{i=1}^{M}\left\vert l_{i}\left(  x\right)  \right\vert
\right)  \left(  \max_{\substack{\left\vert \beta\right\vert =m \\z\in
S_{A,x}}}\left\vert \left(  D^{\beta}f\right)  \left(  z\right)  \right\vert
\right)  \left(  \max_{i=1}^{M}\left\vert a_{i}-x\right\vert \right)
^{m},\quad x\in\mathbb{R}^{d},\label{1.99}%
\end{equation}

where $S_{A,x}=\bigcup\limits_{i=1}^{M}\left[  x,a_{i}\right]  $ is the union
of closed 1-dimensional intervals.
\end{lemma}

\begin{proof}
From Definition \ref{Def_Aux_proj_operator} of $\mathcal{P}$ and $\mathcal{Q}
$%
\[
\mathcal{Q}f\left(  x\right)  =f\left(  x\right)  -\mathcal{P}f\left(
x\right)  =f\left(  x\right)  -\sum_{i=1}^{M}f\left(  a_{i}\right)
l_{i}\left(  x\right)  =f\left(  x\right)  -\sum_{i=1}^{M}f\left(  x+\left(
a_{i}-x\right)  \right)  l_{i}\left(  x\right)  .
\]

Using the Taylor series expansion formula with integral remainder given in
Section \ref{Sect_apx_TaylorSeries} of the Appendix, we have for each $i$:%
\[
f\left(  x+\left(  a_{i}-x\right)  \right)  =\sum_{\left\vert \beta\right\vert
<m}\frac{D^{\beta}f(x)}{\beta!}\left(  a_{i}-x\right)  ^{\beta}+\left(
\mathcal{R}_{m}f\right)  \left(  x,a_{i}-x\right)  ,
\]

where the remainder term satisfies%
\[
\left\vert \left(  \mathcal{R}_{m}f\right)  \left(  x,a_{i}-x\right)
\right\vert \leq\frac{d^{\frac{m}{2}}}{m!}\max_{\substack{\left\vert
\beta\right\vert =m \\t\in\left[  x,a_{i}\right]  }}\left\vert D^{\beta
}f\left(  t\right)  \right\vert \text{ }\left\vert a_{i}-x\right\vert ^{m}.
\]

Thus%
\begin{align*}
\sum_{i=1}^{M}f\left(  x+\left(  a_{i}-x\right)  \right)  l_{i}\left(
x\right)   & =\sum_{i=1}^{M}\left(  \sum_{\left\vert \beta\right\vert <m}%
\frac{D^{\beta}f(x)}{\beta!}\left(  a_{i}-x\right)  ^{\beta}+\left(
\mathcal{R}_{m}f\right)  \left(  x,a_{i}-x\right)  \right)  l_{i}\left(
x\right) \\
& =\sum_{i=1}^{M}\sum_{\left\vert \beta\right\vert <m}\frac{D^{\beta}%
f(x)}{\beta!}\left(  a_{i}-x\right)  ^{\beta}l_{i}\left(  x\right)
+\sum_{i=1}^{M}\left(  \mathcal{R}_{m}f\right)  \left(  x,a_{i}-x\right)
l_{i}\left(  x\right)  .
\end{align*}

But by part 1 of Theorem \ref{Thm_P_op_properties} the operator $\mathcal{P}$
preserves polynomials of degree $<m$. Hence%
\begin{align*}
\sum_{i=1}^{M}\sum_{\left\vert \beta\right\vert <m}\frac{D^{\beta}f(x)}%
{\beta!}\left(  a_{i}-x\right)  ^{\beta}l_{i}\left(  x\right)   &
=\sum_{\left\vert \beta\right\vert <m}\frac{D^{\beta}f(x)}{\beta!}%
\mathcal{P}_{y}\left(  \left(  y-x\right)  ^{\beta}\right)  \left(  y=x\right)
\\
& =\sum_{\left\vert \beta\right\vert <m}\frac{D^{\beta}f(x)}{\beta!}\left(
\left(  y-x\right)  ^{\beta}\right)  \left(  y=x\right) \\
& =f\left(  x\right)  ,
\end{align*}

leaving us with%
\[
\sum_{i=1}^{M}f\left(  x+\left(  a_{i}-x\right)  \right)  l_{i}\left(
x\right)  =f\left(  x\right)  +\sum_{i=1}^{M}\left(  \mathcal{R}_{m}f\right)
\left(  x,a_{i}-x\right)  l_{i}\left(  x\right)  ,
\]

and%
\begin{equation}
\mathcal{Q}f\left(  x\right)  =-\sum_{i=1}^{M}\left(  \mathcal{R}_{m}f\right)
\left(  x,a_{i}-x\right)  l_{i}\left(  x\right)  .\label{1.063}%
\end{equation}

Finally, by applying the remainder estimate \ref{1.34} given in the Appendix,
we get%
\begin{align}
\left\vert \mathcal{Q}f\left(  x\right)  \right\vert  & \leq\left(  \sum
_{i=1}^{M}\left\vert l_{i}\left(  x\right)  \right\vert \right)  \max
_{i=1}^{M}\left\vert \left(  \mathcal{R}_{m}f\right)  \left(  x,a_{i}%
-x\right)  \right\vert \label{7.59}\\
& \leq\frac{d^{\frac{m}{2}}}{m!}\left(  \sum_{i=1}^{M}\left\vert l_{i}\left(
x\right)  \right\vert \right)  \max_{i=1}^{M}\left(  \max
_{\substack{\left\vert \beta\right\vert =m \\t\in\left[  x,a_{i}\right]
}}\left\vert D^{\beta}f\left(  t\right)  \right\vert \text{ }\left\vert
a_{i}-x\right\vert ^{m}\right) \nonumber\\
& \leq\frac{d^{\frac{m}{2}}}{m!}\left(  \sum_{i=1}^{M}\left\vert l_{i}\left(
x\right)  \right\vert \right)  \left(  \max_{\substack{\left\vert
\beta\right\vert =m \\z\in S_{A,x}}}\left\vert D^{\beta}f\left(  z\right)
\right\vert \right)  \left(  \max_{i=1}^{M}\left\vert a_{i}-x\right\vert
\right)  ^{m},\nonumber
\end{align}

as required.
\end{proof}

To study the convergence of the minimal interpolant we will also need the
following lemma which supplies the required results from the theory of
Lagrange interpolation. These results are stated without proof. This lemma has
been created from Lemma 3.2, Lemma 3.5 and Theorem 3.6 of Light and Wayne
\cite{LightWayne98PowFunc}. The results of this lemma do not involve any
reference to weight or basis functions or to functions in $X_{w}^{0}$, but use
the properties of the data region $\Omega$ which contains the independent data
points $X$ and the order of the unisolvency used for the interpolation. Thus
we have separated the part of the proof that involves basis functions from the
part that uses the detailed theory of Lagrange interpolation operators.

\begin{lemma}
\label{Lem_Lagrange_interpol}Suppose first that:

\begin{enumerate}
\item $\Omega$ is a bounded region of $\mathbb{R}^{d}$ having the cone
property e.g. Section 4.3 of Adams \cite{Adams75}.

\item $X$ is a unisolvent subset of $\Omega$ of order $m$.
\end{enumerate}

Suppose $\left\{  l_{j}\right\}  _{j=1}^{M}$ is the cardinal basis of $P_{m}$
with respect to a minimal unisolvent subset of $\Omega$. Using Lagrange
polynomial interpolation techniques, it can be shown there exists a constant
$K_{\Omega,m}^{\prime}>0$ such that
\begin{equation}
\sum\limits_{j=1}^{M}\left\vert l_{j}\left(  x\right)  \right\vert \leq
K_{\Omega,m}^{\prime},\quad x\in\overline{\Omega},\label{1.7}%
\end{equation}

and all minimal unisolvent subsets of $\Omega$. Now define
\[
h_{X,\Omega}=\sup\limits_{\omega\in\Omega}\operatorname*{dist}\left(
\omega,X\right)  ,
\]

and fix $x\in X$. By using Lagrange interpolation techniques it can be shown
there are constants $c_{\Omega,m},h_{\Omega,m}>0$ such that when $h_{X,\Omega
}<h_{\Omega,m}$ there exists a minimal unisolvent set $A\subset X$ satisfying
\[
\operatorname*{diam}\left(  A\cup\left\{  x\right\}  \right)  \leq
c_{\Omega,m}h_{X,\Omega}.
\]

\end{lemma}

Before deriving the first interpolation error estimate I will show that there
exists a suitable sequence of independent data sets $X^{\left(  k\right)
}\subset\Omega$ with $h_{X^{\left(  k\right)  },\Omega}$ tending to zero. The
proof is very close to that of Theorem \ref{Thm_seq_data_regions} which
considered the case of a closed data region.

\begin{theorem}
\label{Thm_seq_data_regions_2}Suppose $\Omega$ is a bounded region containing
all the independent data sets. Then there exists a sequence of independent
data sets $X^{\left(  k\right)  }\subset\Omega$ such that $X^{\left(
k\right)  }\subset X^{\left(  k+1\right)  }$ and $h_{X^{\left(  k\right)
},\Omega}\rightarrow0$ as $k\rightarrow\infty$.
\end{theorem}

\begin{proof}
For $k=1,2,3,\ldots$ there exists a finite covering of $\Omega$ by the balls

$\left\{  B\left(  a_{k}^{\left(  j\right)  };\frac{1}{k}\right)  \right\}
_{j=1}^{M_{k}}$. Construct $X^{\left(  1\right)  }$ by choosing points from
$\Omega$ so that one point lies in each ball $B\left(  a_{k}^{\left(
j\right)  };1\right)  $. Construct $X^{\left(  k+1\right)  }$ by first
choosing the points $X^{\left(  k\right)  }$ and then at least one extra point
so that $X^{\left(  k+1\right)  }$ contains points from each ball $B\left(
a_{k+1}^{\left(  j\right)  };\frac{1}{k+1}\right)  $.

Then $x\in\Omega\cap X^{\left(  k\right)  }$ implies $x\in B\left(
a_{k}^{\left(  j\right)  };\frac{1}{k}\right)  $ for some $j$ and hence
$\operatorname*{dist}\left(  x,X^{\left(  k\right)  }\right)  <\frac{1}{k}$.
Hence $h_{X^{\left(  k\right)  },\Omega}=\sup\limits_{x\in\Omega
}\operatorname*{dist}\left(  x,X^{\left(  k\right)  }\right)  <\frac{1}{k}$
and $\lim\limits_{k\rightarrow\infty}h_{X^{\left(  k\right)  },\Omega}=0$.
\end{proof}

Now we are ready to state our order of convergence result for the minimal norm
interpolant for the case of a weight function with integer parameter
$\kappa\geq1$.

\begin{theorem}
\label{Thm_converg_interpol_ord_gte_1}Let $w$ be a weight function with
property W02 for some $\kappa\in\mathbb{R}^{1}$ such that $\kappa\geq1$ and
let $G$ be the corresponding basis function. Set $m=\left\lfloor
\kappa\right\rfloor $.

Suppose $\mathcal{I}_{X}f$ is the minimal norm interpolant of the data
function $f\in X_{w}^{0}$ on the independent data set $X$ contained in the
data region $\Omega$.

We use the notation and assumptions of Lemma \ref{Lem_Lagrange_interpol} which
means assuming that $X$ is $m$-unisolvent and $\Omega$ is a bounded region
whose boundary satisfies the cone condition.

Now set $k_{G}=\frac{d^{m/2}}{\left(  2\pi\right)  ^{d/2}m!}\left(
c_{\Omega,m}\right)  ^{m}K_{\Omega,m}^{\prime}\max\limits_{\left\vert
\beta\right\vert =m}\left\vert D^{2\beta}G\left(  0\right)  \right\vert $.

Then there exists $h_{\Omega,m}>0$ such that for $x\in\overline{\Omega}$%
\begin{equation}
\left\vert f\left(  x\right)  -\left(  \mathcal{I}_{X}f\right)  \left(
x\right)  \right\vert \leq k_{G}\sqrt{\left(  f-\mathcal{I}_{X}f,f\right)
_{w,0}}\left(  h_{X,\Omega}\right)  ^{m}\leq k_{G}\left\Vert f\right\Vert
_{w,0}\left(  h_{X,\Omega}\right)  ^{m},\label{1.26}%
\end{equation}

when $h_{X,\Omega}=\sup\limits_{\omega\in\Omega}\operatorname*{dist}\left(
\omega,X\right)  <h_{\Omega,m}$ i.e. the order of convergence is at least $m $.

Further, we have the upper bound%
\begin{equation}
\left\vert f\left(  x\right)  -\left(  \mathcal{I}_{X}f\right)  \left(
x\right)  \right\vert \leq k_{G}^{\prime}\sqrt{\left(  f-\mathcal{I}%
_{X}f,f\right)  _{w,0}}\left(  \operatorname*{diam}\Omega\right)  ^{m},\quad
x\in\overline{\Omega},\label{7.37}%
\end{equation}

where $k_{G}^{\prime}=\frac{d^{m/2}}{\left(  2\pi\right)  ^{d/2}m!}%
K_{\Omega,m}^{\prime}\max\limits_{\left\vert \beta\right\vert =m}\left\vert
D^{2\beta}G\left(  0\right)  \right\vert $. The constants $c_{\Omega
,m},K_{\Omega,m}^{\prime}$ and $h_{\Omega,m}$ only depend on $\Omega,m$ and
$d$.

In terms of the integrals which define weight property W02 we have%
\[
\max_{\left\vert \beta\right\vert =m}\left\vert D^{2\beta}G\left(  0\right)
\right\vert \leq\left(  2\pi\right)  ^{-d/2}\int\tfrac{\left\vert
\xi\right\vert ^{2m}d\xi}{w\left(  \xi\right)  }.
\]

\end{theorem}

\begin{proof}
The set $X$ is unisolvent so from Definition \ref{Def_unisolv} it must have a
minimal unisolvent subset, say $A=\left\{  a_{i}\right\}  _{i=1}^{M}$, which
we use to define the Lagrangian operators $\mathcal{P}$ and $\mathcal{Q}%
=I-\mathcal{P}$ of Definition \ref{Def_Aux_proj_operator}. Since $f$ is a data
function, $f\left(  x\right)  -\left(  \mathcal{I}_{X}f\right)  \left(
x\right)  =\mathcal{Q}\left(  f-\mathcal{I}_{X}f\right)  \left(  x\right)  $
and then by estimating $\left\vert \mathcal{Q}\left(  f-\mathcal{I}%
_{X}f\right)  \left(  x\right)  \right\vert $ using Lemma \ref{Lem_Q_estim} we
get%
\begin{equation}
\left\vert f\left(  x\right)  -\mathcal{I}_{X}f\left(  x\right)  \right\vert
\leq\frac{d^{\frac{m}{2}}}{m!}\left(  \sum_{i=1}^{M}\left\vert l_{i}\left(
x\right)  \right\vert \right)  \left(  \max_{\substack{\left\vert
\beta\right\vert =m \\z\in S_{A,x}}}\left\vert D^{\beta}\left(  f-\mathcal{I}%
_{X}f\right)  \left(  z\right)  \right\vert \right)  \max_{i=1}^{M}\left\vert
a_{i}-x\right\vert ^{m}.\nonumber
\end{equation}

The next step is to consider the third factor on the last line. From part 4 of
Theorem \ref{Thm_ord0_Riesz_rep_W2}%
\begin{align}
\left\vert D^{\beta}\left(  f-\mathcal{I}_{X}f\right)  \left(  z\right)
\right\vert  & =\left\vert \left(  f-\mathcal{I}_{X}f,\left(  -D\right)
^{\beta}R_{z}\right)  _{w,0}\right\vert \nonumber\\
& \leq\left\Vert f-\mathcal{I}_{X}f\right\Vert _{w,0}\left\Vert \left(
-D\right)  ^{\beta}R_{z}\right\Vert _{w,0}\nonumber\\
& \leq\left\Vert f-\mathcal{I}_{X}f\right\Vert _{w,0}\left\vert \left(
D^{2\beta}R_{z}\right)  \left(  z\right)  \right\vert ,\label{1.10}%
\end{align}

by part 4 of Theorem \ref{Thm_ord0_Riesz_rep_W2}.

Finally, by part 1 of Theorem \ref{Thm_ord0_Riesz_rep_W2}, $R_{z}=\left(
2\pi\right)  ^{-d/2}G(\cdot-z)$ so that $\left(  D^{2\beta}R_{z}\right)
\left(  z\right)  =\left(  2\pi\right)  ^{-d/2}D^{2\beta}G\left(  0\right)  $
and%
\begin{equation}
\left\vert D^{\beta}\left(  f-\mathcal{I}_{X}f\right)  \left(  z\right)
\right\vert \leq\left(  2\pi\right)  ^{-d/2}\left\Vert f-\mathcal{I}%
_{X}f\right\Vert _{w,0}\left\vert D^{2\beta}G\left(  0\right)  \right\vert
.\label{1.90}%
\end{equation}

Thus%
\[
\max_{\substack{\left\vert \beta\right\vert =m \\z\in S_{A,x}}}\left\vert
D^{\beta}\left(  f-\mathcal{I}_{X}f\right)  \left(  z\right)  \right\vert
\leq\left(  2\pi\right)  ^{-d/2}\left\Vert f-\mathcal{I}_{X}f\right\Vert
_{w,0}\max_{\left\vert \beta\right\vert =m}\left\vert D^{2\beta}G\left(
0\right)  \right\vert ,
\]

and so%
\begin{equation}
\left\vert f\left(  x\right)  -\mathcal{I}_{X}f\left(  x\right)  \right\vert
\leq\frac{d^{\frac{m}{2}}}{\left(  2\pi\right)  ^{\frac{d}{2}}m!}\left\Vert
f-\mathcal{I}_{X}f\right\Vert _{w,0}\max_{\left\vert \beta\right\vert
=m}\left\vert D^{2\beta}G\left(  0\right)  \right\vert \left(  \sum_{i=1}%
^{M}\left\vert l_{i}\left(  x\right)  \right\vert \right)  \max_{i=1}%
^{M}\left\vert a_{i}-x\right\vert ^{m}.\label{7.29}%
\end{equation}

To estimate the last two factors on the right side of the this equation we
will need the previous Lagrangian Lemma \ref{Lem_Lagrange_interpol}. In the
notation of this lemma, if $h_{X,\Omega}=\sup\limits_{\omega\in\Omega
}\operatorname*{dist}\left(  \omega,X\right)  <h_{\Omega,m}$ then for a given
$x$ there exists a minimal unisolvent set $A=\left\{  a_{i}\right\}
_{i=1}^{M}$ such that $\operatorname*{diam}\left(  A\cup\left\{  x\right\}
\right)  \leq c_{\Omega,m}h_{X,\Omega}$ and $\sum\limits_{j=1}^{M}\left\vert
l_{j}\left(  x\right)  \right\vert \leq K_{\Omega,m}^{\prime}$.

Thus $\left\vert a_{i}-x\right\vert \leq c_{\Omega,m}h_{X,\Omega}$ and
\begin{align*}
\left\vert f\left(  x\right)  -\left(  \mathcal{I}_{X}f\right)  \left(
x\right)  \right\vert  & \leq\frac{d^{\frac{m}{2}}}{\left(  2\pi\right)
^{\frac{d}{2}}m!}\left\Vert f-\mathcal{I}_{X}f\right\Vert _{w,0}%
\max_{\left\vert \beta\right\vert =m}\left\vert D^{2\beta}G\left(  0\right)
\right\vert \text{ }K_{\Omega,m}^{\prime}\left(  c_{\Omega,m}h_{X,\Omega
}\right)  ^{m}\\
& =k_{G}\left\Vert f-\mathcal{I}_{X}f\right\Vert _{w,0}\left(  h_{X,\Omega
}\right)  ^{m},
\end{align*}

and because $\left\Vert f-\mathcal{I}_{X}f\right\Vert _{w,0}^{2}=\left(
f-\mathcal{I}_{X}f,f\right)  _{w,0}$ we obtain%
\begin{equation}
\left\vert f\left(  x\right)  -\left(  \mathcal{I}_{X}f\right)  \left(
x\right)  \right\vert \leq k_{G}\sqrt{\left(  f-\mathcal{I}_{X}f,f\right)
_{w,0}}\left(  h_{X,\Omega}\right)  ^{m},\quad x\in\Omega,\label{1.32}%
\end{equation}

which is almost the required inequality. The extension of the last inequality
to $x\in\overline{\Omega}$ is an easy consequence of the fact that $f$ and
$\mathcal{I}_{X}f$ are continuous on $\mathbb{R}^{d}$. The second inequality
now follows directly from \ref{1.09} and so the order of convergence of the
minimal norm interpolant is at least $m=\left\lfloor \kappa\right\rfloor $.

To prove \ref{7.37} we begin with \ref{7.29} and using the inequalities proved
in this proof we obtain directly that%
\begin{align*}
\left\vert f\left(  x\right)  -\mathcal{I}_{X}f\left(  x\right)  \right\vert
& \leq\frac{d^{\frac{m}{2}}}{\left(  2\pi\right)  ^{\frac{d}{2}}m!}%
\sqrt{\left(  f-\mathcal{I}_{X}f,f\right)  _{w,0}}\max_{\left\vert
\beta\right\vert =m}\left\vert D^{2\beta}G\left(  0\right)  \right\vert
K_{\Omega,m}^{\prime}\left(  \operatorname*{diam}\Omega\right)  ^{m}\\
& =k_{G}^{\prime}\sqrt{\left(  f-\mathcal{I}_{X}f,f\right)  _{w,0}}\left(
\operatorname*{diam}\Omega\right)  ^{m},
\end{align*}

when $x\in\Omega$. Again the extension of the last inequality to
$x\in\overline{\Omega}$ is an easy consequence of the fact that $f$ and
$\mathcal{I}_{X}f$ are continuous on $\mathbb{R}^{d}$.

Finally, from \ref{1.024} we have $D^{2\beta}G\left(  x\right)  =\left(
2\pi\right)  ^{-\frac{d}{2}}\int\frac{\left(  i\xi\right)  ^{2\beta}e^{ix\xi}%
}{w\left(  \xi\right)  }d\xi$ and so $D^{2\beta}G\left(  0\right)  =\left(
2\pi\right)  ^{-\frac{d}{2}}\int\frac{\left(  i\xi\right)  ^{2\beta}}{w\left(
\xi\right)  }d\xi$. When $\left\vert \beta\right\vert =m$
\[
\left(  2\pi\right)  ^{\frac{d}{2}}\left\vert D^{2\beta}G\left(  0\right)
\right\vert \leq\int\frac{\left\vert \left(  i\xi\right)  ^{2\beta}\right\vert
}{w\left(  \xi\right)  }d\xi\leq\int\frac{\left\vert \xi\right\vert
^{2\left\vert \beta\right\vert }}{w\left(  \xi\right)  }d\xi=\int%
\frac{\left\vert \xi\right\vert ^{2m}}{w\left(  \xi\right)  }d\xi.
\]

\end{proof}

Thus for an arbitrary data function the order of convergence of the minimal
norm interpolant is at least $\left\lfloor \kappa\right\rfloor $. Using the
same technique as Corollary \ref{Cor_Thm_|f(x)-f(y)|_inequal_2} the following
\textit{double rate} of convergence estimate can be obtained:

\begin{corollary}
\label{Cor_2_Thm_converg_interpol_ord_gte_1}Under the notation and assumptions
of the previous Theorem \ref{Thm_converg_interpol_ord_gte_1} the data
functions $R_{x^{\prime}}$ satisfy
\begin{equation}
\left\vert R_{x^{\prime}}\left(  x\right)  -\left(  \mathcal{I}_{X}%
R_{x^{\prime}}\right)  \left(  x\right)  \right\vert \leq\left(  k_{G}\right)
^{2}\left(  h_{X,\Omega}\right)  ^{2\left\lfloor \kappa\right\rfloor
},\text{\quad}x,x^{\prime}\in K,\label{1.011}%
\end{equation}

when $h_{X,K}<h_{G}$ i.e. the order of convergence is at least $2\left\lfloor
\kappa\right\rfloor $.
\end{corollary}

\subsection{Examples: radial basis function interpolants}

In Subsection \ref{SbSect_wt_func_examples} it was shown that the weight
functions of the \textbf{shifted thin-plate splines} and the \textbf{Gaussian}
satisfy property W02\ for all $\kappa\geq0$ so $\left\lfloor \kappa
\right\rfloor \in\mathbb{Z}_{+}$.

For the \textbf{Sobolev splines} the weight functions satisfy property W02
when $0\leq\kappa<v-d/2$. Thus $\left\lfloor \kappa\right\rfloor =v-d/2-1$
when $v-d/2\in\mathbb{Z}_{+}$ and $v-d/2\geq2$, and $\left\lfloor
\kappa\right\rfloor =\left\lfloor v-d/2\right\rfloor $ when $v-d/2\notin%
\mathbb{Z}_{+}$ and $v-d/2>1$. Compare these with the estimates of
Subsubsection \ref{SbSect_int_examp_rad_basis_ord_converg} as summarized in
Table \ref{Tbl_NonUnisolvTyp1Converg}.

\subsection{Examples: tensor product extended B-spline interpolant}

\begin{corollary}
\label{Cor_interp_converg_B_splin}Suppose the weight function is an extended
B-spline weight function with parameters $2\leq n\leq l$. Then the order of
convergence of the minimal norm interpolant to an arbitrary data function is
at least $n-1$. Further, the order of convergence to a Riesz data function
$R_{y}$ is at least $2n-2$.
\end{corollary}

\begin{proof}
Theorem \ref{Thm_ex_splin_wt_fn_properties} implies that $\left\lfloor
\kappa\right\rfloor =n-1$ and Theorem \ref{Thm_converg_interpol_ord_gte_1}
implies that the order of convergence is at least $n-1$. Corollary
\ref{Cor_2_Thm_converg_interpol_ord_gte_1} now shows the order of convergence
to a Riesz data function $R_{y}$ is at least $2\left\lfloor \kappa
\right\rfloor =2n-2$.
\end{proof}

\subsection{Examples: summary table}%

\begin{table}[htbp] \centering
$%
\begin{tabular}
[c]{|c|c|c|}\hline
\multicolumn{3}{|c|}{Interpolant error estimates: W02 for some $\kappa\geq1
$.}\\
\multicolumn{3}{|c|}{Uses Lagrange interpolation and Taylor series.}\\\hline
& Parameter & Convergence\\
Weight function & constraints & orders $\left(  \left\lfloor \kappa
\right\rfloor \right)  $\\\hline\hline
\multicolumn{1}{|l|}{Sobolev splines} & \multicolumn{1}{|l|}{$v-\frac{d}%
{2}=2,3,4,\ldots$} & \multicolumn{1}{|l|}{$\left\lfloor v-d/2\right\rfloor
-1$}\\\cline{2-3}\cline{2-3}%
\multicolumn{1}{|l|}{\quad$\left(  v>d/2\right)  $} &
\multicolumn{1}{|l|}{$v-\frac{d}{2}>1,$ $v-\frac{d}{2}\notin\mathbb{Z}_{+}$} &
\multicolumn{1}{|l|}{$\left\lfloor v-d/2\right\rfloor $}\\\hline
\multicolumn{1}{|l|}{Shifted thin-plate spline} & - &
\multicolumn{1}{|l|}{$1,2,3,4,\ldots$}\\\hline
\multicolumn{1}{|l|}{Gaussian} & - & \multicolumn{1}{|l|}{$1,2,3,4,\ldots$%
}\\\hline
\multicolumn{1}{|l|}{Extended B-spline} & $n\geq2$ &
\multicolumn{1}{|l|}{$n-1$}\\
\multicolumn{1}{|l|}{\quad$\left(  1\leq n\leq l\right)  $} &  &
\multicolumn{1}{|l|}{}\\\hline
\end{tabular}
$\caption{}\label{Tbl_UnisolvConverg}%
\end{table}%

\section{Data functions and numerical results for type 1 and type 2 estimates
\label{Sect_int_data_fn_exten_Bsplin}}

In this section we will only be interested in the convergence of the
interpolant to it's data function and not in the algorithm's performance as an
interpolant. We will only consider the numerical experiments regarding the
convergence of the interpolants generated by dilations of the one-dimensional
extended B-spline basis functions \ref{1.49} with parameters $n$ and $l$. We
will also restrict ourselves to one dimension so that the data density
parameter can be easily calculated.

We will divide our numerical experiments into those with $n=1$ and those with
$n=2$. For the case $n=1$ we will consider the two cases $n=l=1$ and $n=1$,
$l=2$ so that derivatives are involved, and for $n=2$ we will use $n=2$,
$l=2$. The data region will be the interval $\left[  -1.5,1.5\right]  $.

Because all the scaled extended B-spline basis weight functions have a power
of $\sin x$ in the denominator we will need to derive special classes of data
functions which are convenient for numerical calculations. This will include
characterizing the data functions locally as Sobolev spaces which makes it
easy to choose data functions for numerical experiments.

See also the discussion in Chapter \ref{Ch_cent_diff_wt_fn_ten_prod} regarding
the data functions generated by the \textbf{central difference} weight
functions. There it is shown that locally the sets of data functions for the
extended B-splines and central difference weight functions are identical.

\subsection{Extended B-splines with $n=1$\label{SbSect_int_data_fn_n_eq_1}}

\subsubsection{\protect\underline{The hat function case $n=1$, $l=1$%
}\label{SbSbSect_int_Bspl_neq1_leq1}}

The hat weight function $w_{s}$ is defined by \ref{1.007} and was discussed in
Subsections \ref{SbSect_motiv_weight_fn} and \ref{SbSect_HatWeightFunc}. This
is a scaled, extended B-splines basis function with parameters $n=l=1$ and
Theorem \ref{Thm_basis_tensor_hat_W3} implies $\max\left\lfloor \kappa
\right\rfloor =0$.

Now suppose $\Pi$ is the multivariate tensor product rectangular function
defined using the 1-dimensional function $\Pi\left(  t\right)  =1$ when
$\left\vert t\right\vert <1/2$ and $\Pi\left(  t\right)  =0$ when $\left\vert
t\right\vert >1/2$. We will now justify using data functions of the form
\begin{equation}
f_{d}=u\ast\Pi,\quad u\in L^{2}\left(  \mathbb{R}^{d}\right)  .\label{a9.3}%
\end{equation}

Because $\widehat{\Pi}\left(  t\right)  =\left(  2\pi\right)  ^{-1/2}%
\frac{\sin\left(  t/2\right)  }{t/2}$ and $\widehat{\Lambda}\left(  t\right)
=\left(  2\pi\right)  ^{-1/2}\left(  \dfrac{\sin\left(  t/2\right)  }%
{t/2}\right)  ^{2}$ it follows that%
\[
\left(  \widehat{\Pi}\left(  \xi\right)  \right)  ^{2}=\left(  2\pi\right)
^{-d/2}\widehat{\Lambda}\left(  \xi\right)  ,\quad\xi\in\mathbb{R}^{d}.
\]

Since $\Pi\in L^{1}$, Young's inequality%
\begin{equation}
\left\Vert f\ast g\right\Vert _{r}\leq\left\Vert f\right\Vert _{p}\left\Vert
g\right\Vert _{q},\quad f\in L^{p},g\in L^{q},\text{ }\frac{1}{p}+\frac{1}%
{q}=1+\frac{1}{r},\text{ }1\leq p,q,r\leq\infty,\label{1.056}%
\end{equation}

implies $u\ast\Pi\in L^{2}\left(  \mathbb{R}^{d}\right)  $ and $u\ast
\Pi=\left(  2\pi\right)  ^{-d/2}\int u\left(  y\right)  \Pi\left(
\cdot-y\right)  dy$. Further $u\ast\Pi\in L^{2}$ implies $u\ast\Pi\in
L_{loc}^{1}$ and if $w_{s}=\frac{1}{\widehat{\Lambda}}$ is the hat weight
function%
\begin{equation}
\left\Vert f_{d}\right\Vert _{w_{s},0}=\left(  \int w_{s}\left\vert
\widehat{u}\widehat{\Pi}\right\vert ^{2}\right)  ^{1/2}=\left(  \int\frac
{1}{\widehat{\Lambda}}\left\vert \widehat{u}\widehat{\Pi}\right\vert
^{2}\right)  ^{1/2}=\left(  2\pi\right)  ^{-\frac{d}{4}}\left\Vert
u\right\Vert _{2},\label{a9.4}%
\end{equation}

so that $u\ast\Pi\in X_{w_{s}}^{0}$. Now $u\in L^{2}$ implies $u\in
L_{loc}^{1}$ and we can define $V\in L_{loc}^{1}\left(  \mathbb{R}^{d}\right)
$ by the integral%
\begin{equation}
V\left(  x\right)  =\int_{0}^{x}u\left(  t\right)  dt,\quad u\in L^{2}\left(
\mathbb{R}^{d}\right)  ,\label{1.4}%
\end{equation}

over the volume of the open rectangle $R\left(  0,x\right)  $ so that%
\begin{align}
f_{d}=\left(  2\pi\right)  ^{-\frac{d}{2}}\int u\left(  y\right)  \Pi\left(
x-y\right)  dy &  =\left(  2\pi\right)  ^{-\frac{d}{2}}\int\limits_{x-1/2}%
^{x+1/2}u\left(  y\right)  dy\nonumber\\
&  =\left(  2\pi\right)  ^{-\frac{d}{2}}\left(  V\left(  x+\frac{1}{2}\right)
-V\left(  x-\frac{1}{2}\right)  \right)  ,\label{a9.5}%
\end{align}

and we note that any translate is also a data function.

Observe that by the Cauchy-Schwartz inequality
\[
\left\vert V\left(  x\right)  -V\left(  y\right)  \right\vert =\left\vert
\int_{y}^{x}u\left(  t\right)  dt\right\vert \leq\left\Vert u\right\Vert
_{2}\left\vert x-y\right\vert ^{d/2},
\]

so $V$ is continuous and since $u\in L^{1}\left(  R\left(  0,x\right)
\right)  $ equation \ref{1.4} implies that $D^{\mathbf{1}}V=u$ where
$\mathbf{1}=\left(  1,1,\ldots,1\right)  $.

Given $u\in L^{2}\left(  \mathbb{R}^{d}\right)  $ we can now avoid the
convolution \ref{a9.3} and calculate the data function using \ref{1.4} and
then \ref{a9.5}. The norm $\left\Vert f_{d}\right\Vert _{w_{s},0}$ is
calculated using \ref{a9.4}.

On the other hand, suppose a distribution $V$ satisfies $D^{\mathbf{1}}V\in
L^{2}$ and set $u=D^{\mathbf{1}}V$. Then $V$ satisfies \ref{1.4}, is
continuous and \ref{a9.5} defines the data function. The norm is again given
by \ref{a9.4}.

To obtain our 1-dimensional\ data function $f_{d}$ we will choose
\[
u=e^{-x^{2}},
\]

for which
\begin{equation}
V=\left(  2\pi\right)  ^{\frac{1}{2}}\operatorname{erf},\text{\qquad
}\left\Vert u\right\Vert _{2}=2\left(  2\pi\right)  ^{\frac{1}{4}%
},\text{\qquad}\left\Vert u\ast\Pi\right\Vert _{w_{s},0}=2,\label{1.014}%
\end{equation}

and%
\begin{equation}
f_{d}=u\ast\Pi=\operatorname{erf}\left(  x+\frac{1}{2}\right)
-\operatorname{erf}\left(  x-\frac{1}{2}\right)  .\label{1.002}%
\end{equation}

For the \textit{double rate} convergence experiment we will use $\Lambda
=\sqrt{2\pi}R_{0}$ as the data function. By part 2 of Theorem
\ref{Thm_hat_wt_extend_props}, $\kappa<1/2$ so we must use the Type 1 error
estimates of Theorems \ref{Thm_|f(x)-f(y)|_inequal_2} and
\ref{Thm_interp_error_const_bound}. These imply that if
\[
G\left(  0\right)  -\operatorname{Re}G\left(  x\right)  \leq C_{G}\left\vert
x\right\vert ^{2s},\text{\quad}\left\vert x\right\vert \leq h_{G},
\]

then there exists some $h_{G}>0$ such that
\begin{equation}
\left\vert f_{d}\left(  x\right)  -\left(  \mathcal{I}_{X}f_{d}\right)
\left(  x\right)  \right\vert \leq\left\Vert f_{d}\right\Vert _{w,0}%
\min\left\{  k_{G}\left(  h_{X,K}\right)  ^{s},\sqrt{R_{0}\left(  0\right)
}\right\}  ,\text{\quad}x\in K,\label{1.3}%
\end{equation}

when $h_{X,K}=\sup\limits_{s\in K}\operatorname*{dist}\left(  s,X\right)  \leq
h_{G}$. Here $k_{G}=\left(  2\pi\right)  ^{-\frac{1}{4}}\sqrt{2C_{G}}$,
$R_{0}\left(  0\right)  =\left(  2\pi\right)  ^{-\frac{1}{2}}G\left(
0\right)  $ and $X$ is an independent data set contained in the closed bounded
infinite data region $K$. From Corollary \ref{Cor_Thm_|f(x)-f(y)|_inequal_2}
we have the corresponding double order convergence estimate%
\begin{equation}
\left\vert R_{0}\left(  x\right)  -\left(  \mathcal{I}_{X}R_{0}\right)
\left(  x\right)  \right\vert \leq\min\left\{  \left(  k_{G}\right)
^{2}\left(  h_{X,K}\right)  ^{2s},R_{0}\left(  0\right)  \right\}
,\text{\quad}x\in\mathbb{R}^{1}.\label{1.41}%
\end{equation}

Since we are using the hat basis function Corollary
\ref{Cor_converg_ext_nat_splin_wt_func_n=1} gives: $G\left(  0\right)  =1$,
$C_{G}=1$, $s=1/2$, $h_{G}=\infty$, and from \ref{1.014} and \ref{1.002},
$\left\Vert f_{d}\right\Vert _{w,0}=2$.

\subsubsection{\protect\underline{Numerical results}}

\begin{figure}[th]%
\centering
\fbox{\includegraphics[
natheight=5.086000in,
natwidth=4.083400in,
height=5.086in,
width=4.0834in
]%
{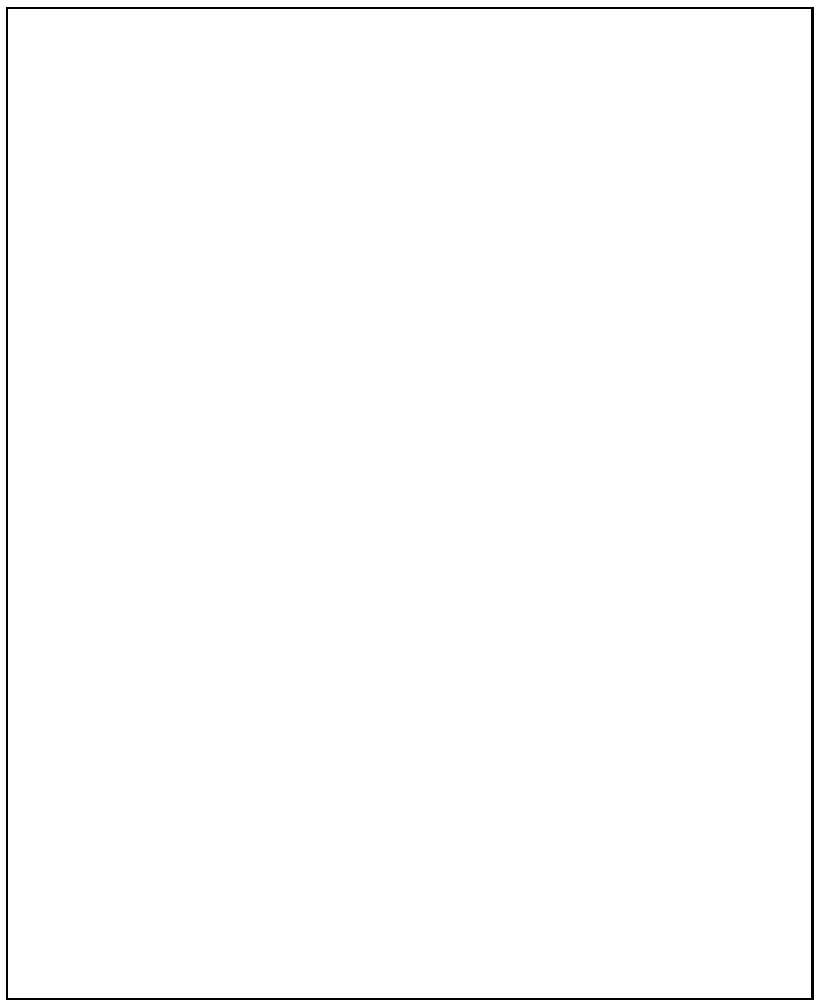}%
}\caption{Interpolant convergence: extended B-spline.}%
\label{Fig_IntepolConverg_N1_L1}%
\end{figure}

Using the functions and parameters discussed in the last subsection the four
subplots displayed in Figure \ref{Fig_IntepolConverg_N1_L1} each display the
superposition of 20 interpolants.%

The two upper subplots relate to the data function \ref{1.002} and the lower
subplots relate to the data function $R_{0}$ i.e. the Riesz representer
\ref{1.019}. The right-hand subplots are filtered versions of the actual
unstable interpolant. The data function is given at the top of the left-hand
plots and the annotation at the bottom of the figure supplies the following
additional information:\medskip

\fbox{Input parameters}\smallskip

$\mathbf{N=L=1}$ - the hat function is a member of the family of scaled
extended B-splines with the indicated parameter values.

\textbf{spl scale 1/2} - changes basis function scale (dilation): $x^{\prime
}=x/spl\_scale$.

\textbf{sm parm 0} - the smoothing parameter is always zero for interpolation.

\textbf{samp 20} - the sample size i.e. the number of test data files
generated. The data function is evaluated on the interval [-1.5,1.5] using a
uniform (statistical) distribution.

\textbf{pts 2:3600} - specifies the smallest number of data points 2 and the
largest number of data points 3600. The other values are given in exponential
steps with a multiplier of approximately 1.2.\medskip

\fbox{Output parameters/messages}\smallskip

\textbf{max ill-cond err -3.5:-4.2} - this relates the ill-conditioning of the
interpolation matrix to the interpolation error. A number or aster * preceding
the colon refers to the data function $f_{d}$ and a number or aster after the
colon refers to the data function $R_{0}$. An aster * will mean there were no
ill-conditioned interpolation matrices generated by the data function. Here
-3.5 means that the largest (unfiltered) interpolation error for which a
matrix was ill-conditioned was $10^{-4.2}$.

\textbf{ill-cond pts 1234:3600 1499:2963} - the first colon-separated group
corresponds to the data function $f_{d}$ and indicates the smallest and
largest numbers of data points for which the interpolant matrix was
ill-conditioned. Asters indicate no ill-conditioning. The second
colon-separated group refers to the data function $R_{0}$.\medskip

Note that all the plots shown in this document have the same format and annotations.

As mentioned above the interpolants are filtered. The filter calculates the
value below which 90\% of the interpolant errors lie. The filter is designed
to remove `large', isolated spikes which dominate the actual errors. The
interpolation error is calculated on a grid with 300 cells applied to the
domain of the data function. No filter is used for the first five interpolants
because there is no instability for small numbers of data points.

As the number of points increases the numerical smoother of $R_{0}$ is
observed to simplify to three very large increasingly narrow spikes at
$\pm1.5$ and $0$, and these dominate by about three orders of magnitude a
residual stable error function of uniform amplitude and zero trend.

The smoother of $f_{d}$ consists of intermingled spikes of various heights
superimposed on a trend curve of amplitude comparable to the average spike
size. The maximum spike height is at most about one order of magnitude of the
average spike height. At the boundary, where the trend is near zero, there are
often two spikes which are narrow w.r.t. the interior spikes.

The (blue) kinked line at the top of each subplot in Figure
\ref{Fig_IntepolConverg_N1_L1} is the theoretical upper bound for the error
given by inequalities \ref{1.3} or \ref{1.41}. Clearly for the data function
$f_{d}$ the theoretical bound of 1/2 underestimates the convergence rate by a
factor of about four - assuming the filtering is valid. The situation for the
data function $R_{0}$ is more complex. The double theoretical convergence rate
of $2s=1$ might suggest that the actual interpolant converges twice as quickly
for the Riesz data function and the filtering reinforces this suspicion.
However, there is no nice linear decrease observed and the error of each
interpolant decreases in large steps to the stable state.

It is also very interesting to note that the theoretical upper error bound for
the data function $R_{0}$ is able to take into account quite closely the
instability of the interpolant.

\subsubsection{\protect\underline{The case $n=1$, $l=2$}}

Amongst other things, the following result will equip us with some tools to
generate data functions for the extended B-splines data spaces, data functions
for which the $X_{w}^{0}$ norm can be calculated \ This result is closely
related to the calculations done above for the hat function.

\begin{definition}
\label{Def_dil_tran_dom_supp}If $x,y\in\mathbb{R}^{d}$ then $x\mathbf{.}%
y=\left(  x_{i}y_{i}\right)  $ denotes the \textbf{component-wise product} of
$x$ and $y$ whereas the scalar or inner product is denoted by $xy$ or $\left(
x,y\right)  $.

\textbf{Component-wise inequality} is denoted $x.\neq y$.

\textbf{Translations of functions} will be defined by $\tau_{c}f\left(
x\right)  =f\left(  x-c\right)  $, $c\in\mathbb{R}^{d}$, and \textbf{dilations
of functions} given by $\sigma_{\lambda}f\left(  x\right)  =f\left(
x./\lambda\right)  $ where $\lambda\in\mathbb{R}^{d}$ and $\lambda.\neq0$.

\textbf{Dilations of sets} will be defined by $\sigma_{\lambda}\left(
\Omega\right)  =\lambda.\Omega$ and \textbf{translations of sets} defined by
$\tau_{c}\left(  \Omega\right)  =\Omega+c$.

Regarding \textbf{function domains}: $\operatorname{dom}\sigma_{\lambda
}f=\lambda.\operatorname{dom}f=\sigma_{\lambda}\operatorname{dom}f$ and
$\operatorname{dom}\tau_{c}f=c+\operatorname{dom}f=\tau_{c}\operatorname{dom}%
f$.

For \textbf{function supports}: $\operatorname*{supp}\sigma_{\lambda}%
f=\lambda.\operatorname*{supp}f=\sigma_{\lambda}\operatorname*{supp}f$ and
$\operatorname*{supp}\tau_{c}f=c+\operatorname*{supp}f=\tau_{c}%
\operatorname*{supp}f$.
\end{definition}

\begin{theorem}
\label{Thm_data_func}\textbf{Central difference operators }and \textbf{global
data functions }If $\alpha$ is a multi-index and $c\in\mathbb{R}^{d} $ we
define the \textbf{central difference operator} $\delta_{c}^{\alpha}$ to be
the composition of the 1-dimensional operators $\delta_{c}^{\alpha}%
=\delta_{c_{1}}^{\alpha_{1}}\delta_{c_{2}}^{\alpha_{2}}\ldots\delta_{c_{d}%
}^{\alpha_{d}}$, where%
\begin{equation}
\delta_{c_{k}}^{\alpha_{k}}=\left(  \delta_{c_{k}}\right)  ^{\alpha_{k}}%
,\quad\delta_{c_{k}}=\tau_{-\frac{c_{k}}{2}\mathbf{e}_{k}}-\tau_{\frac{c_{k}%
}{2}\mathbf{e}_{k}},\quad k=1,\ldots,d,\label{1.001}%
\end{equation}

and if $b\in\mathbb{R}^{1}$ we use the abbreviation $\delta_{b}^{\alpha
}=\delta_{b\mathbf{1}}^{\alpha}$. It then follows (see part 3 of remark below)
that:%
\begin{equation}
\delta_{c}^{\alpha}f=\sum\limits_{\beta\leq\alpha}\left(  -1\right)
^{\left\vert \beta\right\vert }\binom{\alpha}{\beta}\tau_{\left(
2\beta-\alpha\right)  \mathbf{.}\frac{c}{2}}f,\quad f\in\mathcal{D}^{\prime
}\text{ }or\text{ }L^{2},\label{1.015}%
\end{equation}

where $\left(  2\beta-\alpha\right)  \mathbf{.}\frac{c}{2}$ indicates the
component-wise vector product.

Next, suppose $w_{s}$ is the \textbf{extended B-spline} weight function with
parameters $n$ and $l$ given by \ref{1.032}, and that in the sense of
distributions $D^{\alpha}U\in L^{2}\left(  \mathbb{R}^{d}\right)  $ when
$\alpha\leq n\mathbf{1}$.

Then:

\begin{enumerate}
\item
\[
\delta_{c}^{\alpha}U=\left(  \left(  2i\sin\left(  \frac{c}{2}.\xi\right)
\right)  ^{\alpha}\widehat{U}\right)  ^{\vee},
\]

which extends easily to $S^{\prime}$.

\item $\delta_{2}^{l\mathbf{1}}U$ and $\delta_{2}^{2l\mathbf{1}}U$ define
\textbf{data functions} such that%
\begin{equation}
\left\Vert \delta_{2}^{l\mathbf{1}}U\right\Vert _{w_{s},0}=2^{ld}\left\Vert
D^{n\mathbf{1}}U\right\Vert _{2},\quad\left\Vert \delta_{2}^{2l\mathbf{1}%
}U\right\Vert _{w,0}\leq2^{2ld}\left\Vert D^{n\mathbf{1}}U\right\Vert
_{2}.\label{1.003}%
\end{equation}

\item \textbf{Translations} are isometric isomorphisms from $X_{w_{s}}^{0}$ to
$X_{w_{s}}^{0}$.

\item If $m\geq1$ is an integer then the \textbf{dilation} $\sigma_{m}$ is a
continuous mapping from $X_{w}^{0}$ to $X_{w}^{0}$. In fact%
\begin{equation}
\left\Vert \sigma_{m}f\right\Vert _{w_{s},0}\leq m^{\left(  l-n+1/2\right)
d}\left\Vert f\right\Vert _{w_{s},0},\quad f\in X_{w_{s}}^{0}.\label{1.057}%
\end{equation}

\item The \textbf{embedding} $X_{w_{s}}^{0}\hookrightarrow C_{B}^{\left(
n-1\right)  }$ is continuous when $C_{B}^{\left(  n-1\right)  }$ is endowed
with the sup. norm $\max\limits_{\left\vert \alpha\right\vert \leq
n-1}\left\Vert D^{\alpha}u\right\Vert _{\infty}$.

\item Suppose $v\in C_{B}^{\infty}\left(  \mathbb{R}^{d}\right)  $. If $v$ has
periods $\left\{  \mathbf{e}_{k}\right\}  _{k=1}^{d}$ then $v\delta
_{2}^{l\mathbf{1}}U=\delta_{2}^{l\mathbf{1}}\left(  vU\right)  \in X_{w}^{0}$
for all $l$. If $v$ has period $\left\{  2\mathbf{e}_{k}\right\}  _{k=1}^{d}$
then $v\delta_{2}^{l\mathbf{1}}U=\delta_{2}^{l\mathbf{1}}\left(  vU\right)
\in X_{w}^{0}$ when $l$ is even.
\end{enumerate}
\end{theorem}

\begin{proof}
Equation \ref{1.015} can be easily proved in one dimension and then for
arbitrary dimension using the definitions of $\binom{\alpha}{\beta}$ and
$\beta\leq\alpha$.\smallskip

\textbf{Part 1} Since $\widehat{\delta_{c_{k}}U}=e^{-i\frac{c_{k}}{2}\xi_{k}%
}\widehat{U}-e^{+i\frac{c_{k}}{2}\xi_{k}}\widehat{U}=\left(  -2i\sin
\frac{c_{k}}{2}\xi_{k}\right)  \widehat{U}$ i.e. $\delta_{c_{k}}%
^{\mathbf{e}_{k}}U=\left(  \left(  -2i\sin\frac{c_{k}}{2}\xi_{k}\right)
\widehat{U}\right)  ^{\vee}$ we have $\delta_{c_{k}}^{\alpha_{k}}U=\left(
\left(  -2i\sin\frac{c_{k}}{2}\xi_{k}\right)  ^{\alpha_{k}}\widehat{U}\right)
^{\vee}$ and hence
\[
\delta_{c}^{\alpha}U=\left(  \left(  -2i\sin\frac{c_{1}}{2}\xi_{1}\right)
^{\alpha_{1}}\ldots\left(  -2i\sin\frac{c_{d}}{2}\xi_{d}\right)  ^{\alpha_{d}%
}\widehat{U}\right)  ^{\vee}=\left(  \left(  -2i\sin\left(  \frac{c}{2}%
.\xi\right)  \right)  ^{\alpha}\widehat{U}\right)  ^{\vee}.
\]
\smallskip

\textbf{Part 2} The space $X_{w}^{0}$ was introduced in Definition
\ref{Def_Xow}. Now using the Cauchy-Schwartz theorem it is easy to show that
any $L^{2}$ function is $L_{loc}^{1}$. Further, Plancherel's theorem implies
that $\left\Vert \widehat{U}\right\Vert _{2}=\left\Vert U\right\Vert _{2}$ and
so $\widehat{U}\in L_{loc}^{1}$ and $\widehat{U}\in L_{loc}^{1}$. Again by
Plancherel's theorem $D^{\beta}U\in L^{2}$ implies $\left\Vert D^{\beta
}U\right\Vert _{2}^{2}=\left\Vert \widehat{D^{\beta}U}\right\Vert _{2}%
^{2}=\int\xi^{2\beta}\left\vert \widehat{U}\right\vert ^{2}$. Consequently,
from Part 1,%
\begin{align}
\left\Vert \delta_{2}^{l\mathbf{1}}U\right\Vert _{w_{s},0}^{2}=\int
w_{s}\left\vert \widehat{\delta_{2}^{l\mathbf{1}}U}\right\vert ^{2} &
=\int\frac{\xi^{2n\mathbf{1}}}{\left(  \sin\xi_{k}\right)  ^{2l\mathbf{1}}%
}\left\vert \left(  2i\sin\xi_{k}\right)  ^{l\mathbf{1}}\widehat{U}\right\vert
^{2}\label{1.004}\\
&  =2^{2ld}\int\xi^{2n\mathbf{1}}\left\vert \widehat{U}\right\vert
^{2}\nonumber\\
&  =2^{2ld}\left\Vert D^{n\mathbf{1}}U\right\Vert _{2}^{2},\nonumber
\end{align}

and so $\delta_{2}^{l\mathbf{1}}U\in X_{w_{s}}^{0}$. Regarding $\delta
_{2}^{2l\mathbf{1}}U$,
\begin{align*}
\left\Vert \delta_{2}^{2l\mathbf{1}}U\right\Vert _{w_{s},0}^{2}=\int
w_{s}\left\vert \widehat{\delta_{2}^{2l\mathbf{1}}U}\right\vert ^{2}  &
=\int\frac{\xi^{2n\mathbf{1}}}{\left(  \sin\xi_{k}\right)  ^{2l\mathbf{1}}%
}\left\vert \left(  -2i\sin\xi_{k}\right)  ^{2l\mathbf{1}}\widehat{U}%
\right\vert ^{2}\\
& =2^{4ld}\int\frac{\xi^{2n\mathbf{1}}}{\left(  \sin\xi_{k}\right)
^{2l\mathbf{1}}}\left(  \sin\xi_{k}\right)  ^{4l\mathbf{1}}\left\vert
\widehat{U}\right\vert ^{2}\\
& =2^{4ld}\int\left(  \sin\xi_{k}\right)  ^{2l\mathbf{1}}\left\vert
\widehat{D^{n\mathbf{1}}U}\right\vert ^{2}\\
& \leq2^{4ld}\int\left\vert \widehat{D^{n\mathbf{1}}U}\right\vert ^{2}\\
& =2^{4ld}\int\left\vert D^{n\mathbf{1}}U\right\vert ^{2}%
\end{align*}
\smallskip

\textbf{Part 3} The stated properties of translations $\tau_{c}$ follow
directly from the elementary properties: $\tau_{c}\tau_{-c}=\tau_{-c}\tau
_{c}=1$ and $\left\vert \widehat{\tau_{c}u}\right\vert =\left\vert
\widehat{u}\right\vert $.\smallskip

\textbf{Part 4} Regarding the dilations, by using two changes of variable we
obtain%
\begin{align}
\left\Vert \sigma_{m}f\right\Vert _{w_{s},0}^{2}=\int w_{s}\left\vert
\widehat{\sigma_{m\mathbf{1}}f}\right\vert ^{2} &  =\int w_{s}\left\vert
\left(  m\mathbf{1}\right)  ^{\mathbf{1}}\sigma_{\frac{1}{m}\mathbf{1}%
}\widehat{f}\right\vert ^{2}\nonumber\\
&  =m^{2d}\int w_{s}\left(  s\right)  \left\vert \widehat{f}\left(  ms\right)
\right\vert ^{2}ds\nonumber\\
&  =m^{d}\int_{\mathbb{R}^{d}}w_{s}\left(  \frac{t}{m}\right)  \left\vert
\widehat{f}\left(  t\right)  \right\vert ^{2}dt\nonumber\\
&  =m^{d}\int\frac{\left(  \frac{t}{m}\right)  ^{2n\mathbf{1}}}{\left(
\sin\frac{t_{k}}{m}\right)  ^{2l\mathbf{1}}}\left\vert \widehat{f}\left(
t\right)  \right\vert ^{2}dt\nonumber\\
&  =\frac{1}{m^{\left(  2n-1\right)  d}}\int\frac{t^{2n\mathbf{1}}}{\left(
\sin\frac{t_{k}}{m}\right)  ^{2l\mathbf{1}}}\left\vert \widehat{f}\left(
t\right)  \right\vert ^{2}dt\nonumber\\
&  =\frac{1}{m^{\left(  2n-1\right)  d}}\int\frac{\left(  \sin t_{k}\right)
^{2l\mathbf{1}}}{\left(  \sin\frac{t_{k}}{m}\right)  ^{2l\mathbf{1}}}%
\frac{t^{2n\mathbf{1}}}{\left(  \sin t_{k}\right)  ^{2l\mathbf{1}}}\left\vert
\widehat{f}\left(  t\right)  \right\vert ^{2}dt\nonumber\\
&  =\frac{1}{m^{\left(  2n-1\right)  d}}\int\left(  \frac{\sin t_{k}}%
{\sin\frac{t_{k}}{m}}\right)  ^{2l\mathbf{1}}w_{s}\left(  t\right)  \left\vert
\widehat{f}\left(  t\right)  \right\vert ^{2}dt\nonumber\\
&  \leq m^{\left(  2\left(  l-n\right)  +1\right)  d}\left\Vert f\right\Vert
_{w_{s},0}^{2},\label{1.047}%
\end{align}

since $\left\vert \frac{\sin s}{\sin\frac{s}{m}}\right\vert \leq m$ for all
$s\in\mathbb{R}^{1}$ and all integers $m\geq1$. This proves \ref{1.057} and
continuity.\smallskip

\textbf{Part 5} From part 6 Theorem \ref{Thm_ord0_Riesz_rep_W2} the embedding
$X_{w_{s}}^{0}\hookrightarrow C_{B}^{\left(  \kappa\right)  }$ is continuous
and \ref{1.028} allows us to choose $\kappa=n-1$.\smallskip

\textbf{Part 6} Firstly, $U\in L^{2}$ implies $vU\in L^{2}$. Further, since
$L^{2}\subset S^{\prime}$, Leibniz' formula for tempered distributions
$S^{\prime}$ implies
\begin{equation}
D^{\alpha}\left(  vU\right)  =\sum_{\beta\leq\alpha}\binom{\alpha}{\beta
}D^{\alpha-\beta}vD^{\beta}U,\quad U\in L^{2},\text{ }v\in C_{B}^{\infty
}.\label{1.086}%
\end{equation}

Now $D^{\beta}U\in L^{2}$ for $\beta\leq n$, so that $D^{\alpha}\left(
vU\right)  \in L^{2}$ when $\alpha\leq n$. Since $\delta_{2}^{l\mathbf{1}}U\in
X_{w_{s}}^{0}\subset S^{\prime}\subset\mathcal{D}^{\prime}$, if $v$ has period
$\mathbf{e}_{k}$ then in the sense of distributions, in one dimension
\ref{1.015} implies%
\begin{align*}
v\delta_{2}^{l\mathbf{e}_{k}}U  & =v\sum_{j=0}^{l}\left(  -1\right)
^{j}\binom{l}{j}f\left(  \cdot+\left(  l-2j\right)  \mathbf{e}_{k}\right) \\
& =\sum_{j=0}^{l}\left(  -1\right)  ^{j}\binom{l}{j}\left(  vf\right)  \left(
\cdot+\left(  l-2j\right)  \mathbf{e}_{k}\right) \\
& =\delta_{2}^{l\mathbf{e}_{k}}\left(  vU\right)  .
\end{align*}

Consequently $v\delta_{2}^{l\mathbf{1}}U=\delta_{2}^{l\mathbf{1}}\left(
vU\right)  $ and thus \ref{1.003} of this theorem implies $v\delta
_{2}^{l\mathbf{1}}U\in X_{w_{s}}^{0}$. Finally, it is now clear that if $v$
has periods $2\mathbf{e}_{k}$ and $l$ is even then $v\delta_{2}^{l\mathbf{1}%
}U=\delta_{2}^{l\mathbf{1}}\left(  vU\right)  $.
\end{proof}

\begin{remark}
\label{Rem_Thm_data_func}\ 

\begin{enumerate}
\item \textbf{Norms of data functions} Further to the comment preceding the
last theorem: If $D^{\alpha}U\in L^{2}\left(  \mathbb{R}^{d}\right)  $ for
$\alpha\leq n\mathbf{1}$ then $f_{d}=\delta_{2}^{l\mathbf{1}}U$ is a data
function with norm $2^{ld}\left\Vert D^{n\mathbf{1}}U\right\Vert _{2}$ and
translating this function does not change the norm. In addition, if $v\in
C_{B}^{\infty}\left(  \mathbb{R}^{d}\right)  $ has periods $\left\{
\mathbf{e}_{k}\right\}  _{k=1}^{d}$ when $l$ is odd or $\nu$ has period
$\left\{  2\mathbf{e}_{k}\right\}  _{k=1}^{d}$ when $l$ is even, then
$vf_{d}=\delta_{2}^{l\mathbf{1}}\left(  vU\right)  $ is also a data function
with norm $2^{ld}\left\Vert D^{n\mathbf{1}}\left(  vU\right)  \right\Vert
_{2}$.

\item \textbf{Regarding part 6 of the last theorem} The condition $v\in
C_{B}^{\infty}\left(  \mathbb{R}^{d}\right)  $ can be weakened to $v\in
W^{n\mathbf{1},\infty}\left(  \mathbb{R}^{d}\right)  $ where%
\begin{equation}
W^{n\mathbf{1},\infty}\left(  \mathbb{R}^{d}\right)  =\left\{  f\in L^{\infty
}\left(  \mathbb{R}^{d}\right)  :D^{\alpha}f\in L^{\infty},\text{ }\alpha\leq
n\mathbf{1}\right\}  ,\label{a1.25}%
\end{equation}

is the $L^{\infty}$ Sobolev space introduced below in Definition
\ref{Def_SobolevSpace2}. This can be proved using the Leibniz formula
$D\left(  vf\right)  =vDf+\left(  Dv\right)  f$ where $v\in W^{n\mathbf{1}%
,\infty}$, $f\in W^{n\mathbf{1}}$. This formula can in turn be demonstrated by
adapting the proof of Theorem 14.2, Section 1.14 Petersen \cite{Petersen83}.

\item We use the convenient full Fourier transform notation $F\left[
f\right]  $ and partial Fourier transform $F_{x_{k}}\left[  f\left(  x\right)
\right]  $.

Since $F_{x_{k}}\left[  f\left(  x-t\mathbf{e}_{k}\right)  \right]
=e^{-it\xi_{k}}F_{x_{k}}\left[  f\left(  x\right)  \right]  $ and $F\left[
f\left(  x-\tau\right)  \right]  =e^{-i\tau\xi}F\left[  f\left(  x\right)
\right]  $, if $f\in L^{2}$ we have the sequence of equations%
\begin{align*}
F_{x_{k}}\left[  \delta_{c_{k}}f\left(  x\right)  \right]   & =F_{x_{k}%
}\left[  \tau_{-\frac{1}{2}c_{k}\mathbf{e}_{k}}f\right]  -F_{x_{k}}\left[
\tau_{\frac{1}{2}c_{k}\mathbf{e}_{k}}f\right] \\
& =\left(  e^{\frac{i}{2}c_{k}\xi_{k}}-e^{-\frac{i}{2}c_{k}\xi_{k}}\right)
F_{x_{k}}\left[  f\right]  .\\
F_{x_{k}}\left[  \delta_{c_{k}}^{\alpha_{k}}f\right]   & =\left(  e^{\frac
{i}{2}c_{k}\xi_{k}}-e^{-\frac{i}{2}c_{k}\xi_{k}}\right)  ^{\alpha_{k}}%
F_{x_{k}}\left[  f\right]  ,
\end{align*}

so that by means of the binomial theorem%
\begin{align*}
F\left[  \delta_{c}^{\alpha}f\right]   & =\left(  \left(  e^{\frac{i}{2}%
c_{k}\xi_{k}}\right)  +\left(  -e^{-\frac{i}{2}c_{k}\xi_{k}}\right)  \right)
^{\alpha}F\left[  f\right] \\
& =\sum_{\beta\leq\alpha}\tbinom{\alpha}{\beta}\left(  e^{\frac{i}{2}c_{k}%
\xi_{k}}\right)  ^{\alpha-\beta}\left(  -e^{-\frac{i}{2}c_{k}\xi_{k}}\right)
^{\beta}F\left[  f\right] \\
& =\sum_{\beta\leq\alpha}\left(  -1\right)  ^{\left\vert \beta\right\vert
}\tbinom{\alpha}{\beta}\left(  e^{\frac{i}{2}%
{\textstyle\sum}
\left(  \alpha_{k}-\beta_{k}\right)  c_{k}\xi_{k}}\right)  \left(
e^{-\frac{i}{2}%
{\textstyle\sum}
\beta_{k}c_{k}\xi_{k}}\right)  F\left[  f\right] \\
& =\sum_{\beta\leq\alpha}\left(  -1\right)  ^{\left\vert \beta\right\vert
}\tbinom{\alpha}{\beta}e^{\frac{i}{2}%
{\textstyle\sum}
\left(  \alpha_{k}-2\beta_{k}\right)  c_{k}\xi_{k}}F\left[  f\right] \\
& =\sum_{\beta\leq\alpha}\left(  -1\right)  ^{\left\vert \beta\right\vert
}\tbinom{\alpha}{\beta}e^{-i\left(  \left(  2\beta_{k}-\alpha_{k}\right)
\frac{c_{k}}{2}\right)  \xi}F\left[  f\right] \\
& =\sum_{\beta\leq\alpha}\left(  -1\right)  ^{\left\vert \beta\right\vert
}\tbinom{\alpha}{\beta}e^{-i\left(  \left(  2\beta-\alpha\right)  .\frac{c}%
{2}\right)  \xi}F\left[  f\right] \\
& =\sum_{\beta\leq\alpha}\left(  -1\right)  ^{\left\vert \beta\right\vert
}\tbinom{\alpha}{\beta}F\left[  \tau_{\left(  2\beta-\alpha\right)
\mathbf{.}\frac{c}{2}}f\right] \\
& =F\left[  \sum_{\beta\leq\alpha}\left(  -1\right)  ^{\left\vert
\beta\right\vert }\tbinom{\alpha}{\beta}\tau_{\left(  2\beta-\alpha\right)
\mathbf{.}\frac{c}{2}}f\right]  ,
\end{align*}

and applying the inverse Fourier transform gives \ref{1.015} on $L^{2}$. This
immediately means that \ref{1.015} holds for the distribution test functions
$C_{0}^{\infty}$ and hence for all distributions.
\end{enumerate}
\end{remark}

Our basis function is the (unscaled) 1-dimensional extended B-spline $G_{1}$
with parameters $n=1$ and $l=2$ given by \ref{1.49} i.e.%
\[
G_{1}\left(  t\right)  =\left(  -1\right)  ^{l-n}\tfrac{\left(  2\pi\right)
^{l/2}}{2^{2\left(  l-n\right)  +1}}\left(  D^{2\left(  l-n\right)  }\left(
\left(  \ast\Lambda\right)  ^{l}\right)  \right)  \left(  \tfrac{t}{2}\right)
=-\tfrac{\pi}{4}\left(  D^{2}\left(  \Lambda\ast\Lambda\right)  \right)
\left(  \tfrac{t}{2}\right)
\]
for $t\in\mathbb{R}^{1}$. But%
\begin{align}
D^{2}\left(  \Lambda\ast\Lambda\right)  =\Lambda\ast D^{2}\Lambda &
=\Lambda\ast\left(  \delta\left(  \cdot+1\right)  -2\delta+\delta\left(
\cdot-1\right)  \right) \nonumber\\
& =\frac{1}{\sqrt{2\pi}}\left(  \Lambda\left(  \cdot+1\right)  -2\Lambda
+\Lambda\left(  \cdot-1\right)  \right)  ,\label{1.28}%
\end{align}

so that%
\[
G_{1}\left(  t\right)  =-\frac{\sqrt{2\pi}}{8}\left(  \Lambda\left(  \tfrac
{t}{2}+1\right)  -2\Lambda\left(  \tfrac{t}{2}\right)  +\Lambda\left(
\tfrac{t}{2}-1\right)  \right)  ,
\]

and hence%
\[
DG_{1}\left(  t\right)  =-\frac{\sqrt{2\pi}}{16}\left(  \Lambda^{\prime
}\left(  \tfrac{t}{2}+1\right)  -2\Lambda^{\prime}\left(  \tfrac{t}{2}\right)
+\Lambda^{\prime}\left(  \tfrac{t}{2}-1\right)  \right)  ,
\]

i.e. $\left\Vert DG_{1}\right\Vert _{\infty}=\frac{3}{16}\sqrt{2\pi}$. Now by
Theorem \ref{Thm_ex_nat_spline_basis_Lipschitz}%
\[
G_{1}\left(  0\right)  -G_{1}\left(  t\right)  \leq\left\Vert DG_{1}%
\right\Vert _{\infty}\left\vert t\right\vert ,\text{\quad}x\in\mathbb{R}^{1},
\]

which means that
\[
G_{1}\left(  0\right)  =\frac{\sqrt{2\pi}}{4},\text{\quad}C_{G}=\left\Vert
DG_{1}\right\Vert _{\infty}=\frac{3}{16}\sqrt{2\pi},\text{\quad}s=\frac{1}%
{2},\text{\quad}h_{G}=\infty.
\]

With reference to the last theorem we will choose the bell-shaped data
function%
\begin{equation}
f_{d}=\delta_{2}^{2}U\in X_{w}^{0},\label{1.017}%
\end{equation}

where%
\begin{equation}
U\left(  x\right)  =\frac{e^{-k_{1,2}x^{2}}}{\delta_{2}^{2}\left(
e^{-k_{1,2}x^{2}}\right)  \left(  0\right)  }=\frac{e^{-k_{1,2}x^{2}}%
}{2\left(  1-e^{-4k_{1,2}}\right)  },\quad k_{1,2}=0.3,\label{1.035}%
\end{equation}

so that
\[
\left\Vert f_{d}\right\Vert _{w,0}=4\left\Vert DU\right\Vert _{2}%
=\sqrt[4]{2\pi}\frac{\sqrt[4]{4k_{1,2}}}{1-e^{-4k_{1,2}}}.
\]

The interpolation error estimates given by \ref{1.3} and \ref{1.41} are now%
\[
\left\vert f_{d}\left(  x\right)  -\left(  \mathcal{I}_{X}f_{d}\right)
\left(  x\right)  \right\vert \leq\left\Vert f_{d}\right\Vert _{w,0}%
\min\left\{  k_{G}\left(  h_{X,K}\right)  ^{s},\sqrt{R_{0}\left(  0\right)
}\right\}  ,\text{\quad}x\in K,
\]

and%
\[
\left\vert R_{0}\left(  x\right)  -\left(  \mathcal{I}_{X}R_{0}\right)
\left(  x\right)  \right\vert \leq\min\left\{  \left(  k_{G}\right)
^{2}\left(  h_{X,K}\right)  ^{2s},R_{0}\left(  0\right)  \right\}
,\text{\quad}x\in\mathbb{R}^{1},
\]

respectively, where $k_{G}=\left(  2\pi\right)  ^{-\frac{1}{4}}\sqrt{2C_{G}}$
and $R_{0}\left(  0\right)  =\left(  2\pi\right)  ^{-\frac{1}{2}}G_{1}\left(
0\right)  $.

For the theory developed in the previous sections it was convenient to use the
simple, unscaled weight function definition \ref{1.032}. However, I have
observed that scaling can significantly improve the performance of the basis
function interpolant and so I will present the following theorem for the
scaled, extended B-splines.

\begin{theorem}
\label{Cor_Thm_data_func}Suppose $\widetilde{G}_{s}\left(  x\right)  =%
{\textstyle\prod\limits_{k=1}^{d}}
\widetilde{G}_{1}\left(  x_{k}\right)  $ where $\widetilde{G}_{1}=\left(
-1\right)  ^{l-n}D^{2\left(  l-n\right)  }\left(  \left(  \ast\Lambda\right)
^{l}\right)  $ and $\Lambda$ is the univariate hat function and $n,l$ are
integers such that $1\leq n\leq l$.

For given $\lambda>0$, $\widetilde{G}_{s}\left(  \lambda x\right)  $ is called
a scaled extended B-spline basis function. The corresponding weight function
is $\widetilde{w}_{\lambda}\left(  t\right)  =\left(  2\lambda a\right)
^{d}w_{s}\left(  \frac{t}{2\lambda}\right)  $ where $a=\frac{\left(
2\pi\right)  ^{l/2}}{2^{2\left(  l-n\right)  +1}}$ and $w_{s}$ is the extended
B-spline weight function \ref{1.032} with parameters $n,l$. Indeed,
$\widetilde{w}_{\lambda}$ has property W02 for $\kappa$ iff $w_{s}$ has
property W02 for $\kappa$. Further%
\begin{equation}
\widetilde{G}_{s}\left(  0\right)  -\widetilde{G}_{s}\left(  \lambda x\right)
\leq2\lambda a^{-d}\sqrt{d}G_{1}\left(  0\right)  ^{d-1}\left\Vert
DG_{1}\right\Vert _{\infty}\left\vert x\right\vert ,\text{\quad}x\in
\mathbb{R}^{d}.\label{1.016}%
\end{equation}

Finally, if $f_{d}\in X_{w}^{0}$ and $g_{d}\left(  x\right)  =f_{d}\left(
2\lambda x\right)  $, it follows that $g_{d}\in X_{\widetilde{w}_{\lambda}%
}^{0}$ and $\left\Vert g_{d}\right\Vert _{\widetilde{w}_{\lambda},0}%
=a^{d/2}\left\Vert f_{d}\right\Vert _{w,0}$.
\end{theorem}

\begin{proof}
From \ref{1.49}, $\widetilde{G}_{s}\left(  x\right)  =a^{-d}G_{s}\left(
2x\right)  $, where $G_{s}$ is the extended B-spline basis function with
parameters $n$ and $l$. Hence by Theorem \ref{Thm_G_in_Xow} the corresponding
weight function is $w_{\lambda}\left(  t\right)  =\left(  2\lambda a\right)
^{d}w\left(  \frac{t}{2\lambda}\right)  $ with property W02 for $\kappa$. By
Theorem \ref{Thm_ex_nat_spline_basis_Lipschitz}%
\[
G_{s}\left(  0\right)  -G_{s}\left(  x\right)  \leq\sqrt{d}G_{1}\left(
0\right)  ^{d-1}\left\Vert DG_{1}\right\Vert _{\infty}\left\vert x\right\vert
,\text{\quad}x\in\mathbb{R}^{1},
\]

and because $\widetilde{G}_{s}\left(  \lambda x\right)  =a^{-d}G_{s}\left(
2\lambda x\right)  $ we have%
\begin{align*}
\widetilde{G}_{s}\left(  0\right)  -\widetilde{G}_{s}\left(  \lambda x\right)
=a^{-d}\left(  G_{s}\left(  0\right)  -G_{s}\left(  2\lambda x\right)
\right)   &  \leq a^{-d}\sqrt{d}G_{1}\left(  0\right)  ^{d-1}\left\Vert
DG_{1}\right\Vert _{\infty}\left\vert 2\lambda x\right\vert \\
&  =2\lambda a^{-d}\sqrt{d}G_{1}\left(  0\right)  ^{d-1}\left\Vert
DG_{1}\right\Vert _{\infty}\left\vert x\right\vert .
\end{align*}

Finally, $\widehat{g_{d}}\left(  t\right)  =\left(  2\lambda\right)
^{-d}\widehat{f_{d}}\left(  \frac{t}{2\lambda}\right)  $ and so
\[
\left\Vert g_{d}\right\Vert _{w_{\lambda},0}^{2}=\int\left(  2\lambda
a\right)  ^{d}w_{s}\left(  \frac{t}{2\lambda}\right)  \left\vert \left(
2\lambda\right)  ^{-d}\widehat{f_{d}}\left(  \frac{t}{2\lambda}\right)
\right\vert ^{2}dt=a^{d}\int w_{s}\left(  s\right)  \left\vert \widehat{f_{d}%
}\left(  s\right)  \right\vert ^{2}ds=a^{d}\left\Vert f_{d}\right\Vert
_{w,0}^{2}.
\]

\end{proof}

\subsubsection{\protect\underline{Numerical results}}

Using the functions and parameters discussed in the last subsection the four
subplots displayed in Figure \ref{Fig_IntepolConverg_N1_L2} each display the
superposition of 20 interpolants. The two upper subplots relate to the data
function \ref{1.017} and the lower subplot relates to the data function
$R_{0}$ i.e. the Riesz representer \ref{1.019}. The right-hand subplots are
filtered versions of the actual unstable interpolant. The four subplots of
Figure \ref{Fig_IntepolConverg_N1_L2} each display the superposition of 20 interpolants.%

\begin{figure}[th]%
\centering
\includegraphics[
natheight=4.736600in,
natwidth=4.063800in,
height=4.7366in,
width=4.0638in
]%
{C:/Math_SwBasisFunc/InterpolSmthDev/PapersMonog/ZeroOrd/ZeroOrdDev/graphics/figInterpolConverg_N1_L2_samp20_3600pt_scal_1__2.pdf}%
\caption{Interpolant convergence: extended B-spline.}%
\label{Fig_IntepolConverg_N1_L2}%
\end{figure}

As the number of points increases the numerical interpolant of $f_{d}$
gradually simplifies to two spikes at the end points of the data interval with
stable errors in between. The interpolant of $R_{0}$ is observed to simplify
to a large increasingly narrow spike at zero and this dominates a stable error
function with absolute value of the order of $10^{-7}$.

The (blue) kinked line above each interpolant at the top of each subplot in
Figure \ref{Fig_IntepolConverg_N1_L2} is the estimated upper bound for the
error given by the inequalities \ref{1.3} or \ref{1.41}. Clearly for the data
function $f_{d}$ the theoretical bound of 1/2 underestimates the convergence
rate by a factor of approximately four - assuming the filtering is valid. The
situation for the data function $R_{0}$ is more complex. The double
theoretical convergence rate of $2s=1$ might suggest that the actual
interpolant converges twice as quickly for the Riesz data function and the
filtering reinforces this suspicion. However, there is no nice linear decrease
observed and the error of each interpolant decreases in a large step and then
smaller steps to the stable state.

It is also very interesting to note that the theoretical upper error bound for
the data function $R_{0}$ is able to take into account quite closely the
instability of the interpolant.

\subsection{Extended B-splines with $n=2$\label{SbSect_int_data_fn_n_eq_2}}

Since $n\geq2$, Table \ref{Tbl_InterpNonUnisolvTyp2Conv} tells us we can use
the Type 2 error estimates of Theorem
\ref{Thm_interpol_error_in_terms_of_wt_fn} and the estimate of Theorem
\ref{Thm_interp_error_const_bound}.

\subsubsection{\protect\underline{The case: $n=2$, $l=2$}}

From Corollary \ref{Cor_Thm_data_func} the function
\[
G_{2,2}\left(  x\right)  =\frac{\left(  \Lambda\ast\Lambda\right)  \left(
2x\right)  }{\left(  \Lambda\ast\Lambda\right)  \left(  0\right)  }%
=\frac{3\sqrt{2\pi}}{2}\left(  \Lambda\ast\Lambda\right)  \left(  2x\right)  ,
\]

is a scaled extended B-spline basis function with parameters $n=2$, $l=2$ such
that $\operatorname*{supp}G_{2,2}=\left[  -1,1\right]  $ and $G_{2,2}\left(
0\right)  =1$. From \ref{1.49} the scaling factor is $\lambda=4$. To calculate
$G_{2,2}$ we use the convenient formula%
\begin{equation}
G_{2,2}\left(  x\right)  =\left(  1+x\right)  ^{2}\Lambda\left(  2x+1\right)
+\left(  1-2x^{2}\right)  \Lambda\left(  2x\right)  +\left(  1-x\right)
^{2}\Lambda\left(  2x-1\right)  ,\label{1.38}%
\end{equation}

derived by observing that we require the symmetric equivalents of the
B-splines $b^{k}$ studied, for example, in Chapter 3 of H\"{o}llig
\cite{Hollig2003}: Noting that $\Lambda\left(  x\right)  =b^{1}\left(
x+1\right)  $, the symmetric equivalents of the convolution formula of box
3.11 and of the formulas of boxes 3.3, 3.4 can be used to derive \ref{1.38}.
With reference to Theorem \ref{Thm_data_func} choose the data function
\ref{1.017}, \ref{1.035} and for double rate convergence experiments we will
use $R_{0}=\left(  2\pi\right)  ^{-\frac{1}{2}}G_{2,2}$ as the data function.

Since $n\geq2$, Table \ref{Tbl_InterpNonUnisolvTyp2Conv} tells us we can use
the Type 2 error estimates of Theorem
\ref{Thm_interpol_error_in_terms_of_wt_fn} and the estimate of Theorem
\ref{Thm_interp_error_const_bound}. These imply that
\begin{equation}
\left\vert f_{d}\left(  x\right)  -\left(  \mathcal{I}_{X}f_{d}\right)
\left(  x\right)  \right\vert \leq\left\Vert f_{d}\right\Vert _{w,0}%
\min\left\{  k_{G}h_{X,K},\sqrt{R_{0}\left(  0\right)  }\right\}
,\text{\quad}x\in K,\,f_{d}\in X_{w}^{0},\label{1.411}%
\end{equation}

and%
\begin{equation}
\left\vert R_{0}\left(  x\right)  -\left(  \mathcal{I}_{X}R_{0}\right)
\left(  x\right)  \right\vert \leq\min\left\{  \left(  k_{G}\right)
^{2}\left(  h_{X,K}\right)  ^{2},R_{0}\left(  0\right)  \right\}
,\text{\quad}x\in\mathbb{R}^{1},\label{1.111}%
\end{equation}

where $h_{X,K}=\sup\limits_{s\in K}\operatorname*{dist}\left(  s,X\right)  $
and in addition
\begin{align*}
k_{G}  & =\left(  2\pi\right)  ^{-\frac{1}{4}}\sqrt{-G_{2,2}\left(  0\right)
^{d-1}D^{2}G_{2,2}\left(  0\right)  }\sqrt{d}=\left(  2\pi\right)  ^{-\frac
{1}{4}}\sqrt{-D^{2}G_{2,2}\left(  0\right)  },\\
R_{0}\left(  0\right)   & =\left(  2\pi\right)  ^{-\frac{1}{2}}G_{2,2}\left(
0\right)  =\left(  2\pi\right)  ^{-\frac{1}{2}},
\end{align*}
and $K=\left[  -1.5,1.5\right]  $.

Using \ref{1.28} we have $-D^{2}G_{2,2}\left(  0\right)  =-6\sqrt{2\pi}%
D^{2}\left(  \Lambda\ast\Lambda\right)  \left(  0\right)  =6\sqrt{2\pi}%
\frac{2}{\sqrt{2\pi}}=12$ so that $k_{G}=\left(  2\pi\right)  ^{-\frac{1}{4}%
}\sqrt{12}$.

It remains to calculate $\left\Vert f_{d}\right\Vert _{w,0}$. But from
\ref{1.017} and Plancherel's theorem
\[
\left\Vert f_{d}\right\Vert _{w,0}=2^{ld}\left\Vert D^{n\mathbf{1}%
}U\right\Vert _{2}=4\left\Vert D^{2}U\right\Vert _{2}=4\left\Vert \xi
^{2}\widehat{U}\right\Vert _{2}.
\]

By \ref{1.035}, $U\left(  x\right)  =\frac{e^{-k_{1,2}x^{2}}}{2\left(
1-e^{-4k_{1,2}}\right)  }$ so $\widehat{U}\left(  \xi\right)  =\frac
{1}{2\left(  1-e^{-4k_{1,2}}\right)  }\frac{1}{\sqrt{k_{1,2}}}e^{-\frac
{\xi^{2}}{k_{1,2}}}$ and%
\begin{align*}
\left\Vert \xi^{2}\widehat{U}\right\Vert _{2}  & =\frac{4}{2\left(
1-e^{-4k_{1,2}}\right)  \sqrt{k_{1,2}}}\left(  \int\xi^{4}e^{-\frac{2\xi^{2}%
}{k_{1,2}}}d\xi\right)  ^{1/2}\\
& =\frac{2}{\left(  1-e^{-4k_{1,2}}\right)  \sqrt{k_{1,2}}}\frac{k_{1,2}}%
{2}\sqrt[4]{\frac{k_{1,2}}{2}}\left(  \int\xi^{4}e^{-\xi^{2}}d\xi\right)
^{1/2}\\
& =\frac{2}{\left(  1-e^{-4k_{1,2}}\right)  \sqrt{k_{1,2}}}\frac{k_{1,2}}%
{2}\sqrt[4]{\frac{k_{1,2}}{2}}\left(  \frac{3}{4}\sqrt{\pi}\right)  ^{1/2}\\
& =\frac{1}{2\left(  1-e^{-4k_{1,2}}\right)  \sqrt{k_{1,2}}}k_{1,2}%
\sqrt[4]{k_{1,2}}\left(  \frac{9\pi}{2}\right)  ^{1/4}\\
& =\frac{1}{2}\sqrt[4]{\frac{9\pi}{2}}\frac{\left(  k_{1,2}\right)  ^{3/4}%
}{1-e^{-4k_{1,2}}},
\end{align*}

and hence $\left\Vert f_{d}\right\Vert _{w,0}=2\sqrt[4]{\frac{9\pi}{2}}%
\frac{\left(  k_{1,2}\right)  ^{3/4}}{1-e^{-4k_{1,2}}}=\sqrt[4]{72\pi}%
\frac{\left(  k_{1,2}\right)  ^{3/4}}{1-e^{-4k_{1,2}}}$. To summarize:%
\[
k_{G}=\frac{\sqrt{12}}{\left(  2\pi\right)  ^{1/4}};\quad\left\Vert
f_{d}\right\Vert _{w,0}=\sqrt[4]{72\pi}\frac{\left(  k_{1,2}\right)  ^{3/4}%
}{1-e^{-4k_{1,2}}},\quad k_{1,2}=0.3;\quad h_{G}=\infty.
\]

\subsubsection{\protect\underline{Numerical results}}

Using the functions and parameters derived above the 4 subplots of Figure
\ref{Fig_IntepolConverg_N2_L2} were generated. Each plot displays the
superposition of 20 interpolants with unfiltered output on the left and
filtered interpolants on the right. Above the interpolants is a (blue) kinked
line whose slope indicates the theoretically predicted rate of convergence
given by inequalities \ref{1.411} and \ref{1.111}. There is also an adjacent
line which represents an estimate of the slope implied by the filtered
interpolant. For $f_{d}=\delta_{2}^{2}\frac{e^{-k_{1,2}x^{2}}}{2\left(
1-e^{-4k_{1,2}}\right)  }$ the estimated rate of convergence is $n-1=1 $ and
for $R_{0}$ the estimated convergence rate is $2n-2=2$.%

\begin{figure}[th]%
\centering
\includegraphics[
natheight=5.078200in,
natwidth=4.063800in,
height=5.0782in,
width=4.0638in
]%
{C:/Math_SwBasisFunc/InterpolSmthDev/PapersMonog/ZeroOrd/ZeroOrdDev/graphics/figInterpolConverg_N2_L2_samp20_3600pt_scal_inv4__3.pdf}%
\caption{Interpolant convergence: extended B-spline.}%
\label{Fig_IntepolConverg_N2_L2}%
\end{figure}

As the number of points increases both numerical interpolants were observed to
consist of spikes of various sizes but no dominant spikes.

The (blue) kinked line above each interpolant at the top of each subplot in
Figure \ref{Fig_IntepolConverg_N2_L2} is the theoretical upper bound for the
error given by inequalities \ref{1.411} or \ref{1.111}. Clearly for the data
function $f_{d} $ the theoretical bound of $1$ underestimates the convergence
rate by a factor of about 3 - assuming the filtering is valid. For for the
data function $R_{0}$ a theoretical bound of 2 underestimates a convergence
rate of about 5. The double theoretical convergence rate of 2 might suggest
that the actual interpolant converges twice as quickly for the Riesz data
function and the filtering reinforces this suspicion.

\section{Characterizing certain global data spaces
locally\label{Sect_local_data_space}}

The \textbf{global data space} is the Hilbert space $X_{w}^{0}$ of continuous
functions on $\mathbb{R}^{d}$ used to define the minimal norm interpolant, and
later, the Exact and Approximate smoothers. However, these interpolants and
smoothers are independent of the values of the data functions outside the data
region $\Omega$. Therefore we will want pointwise error estimates which
involve the \textbf{local data spaces} i.e. the local restriction space
$X_{w}^{0}\left(  \Omega\right)  $.

In the next subsection we will use the central difference operators
\ref{1.015} to characterize the local data spaces in the case of the scaled
B-spline weight functions: in fact we show that the local data spaces contain
the same set of functions as a special, local, Sobolev space, denoted
$W^{m\mathbf{1}}\left(  \Omega\right)  $. In the subsequent subsection we will
then extend the class of weight functions for which this $W^{m\mathbf{1}%
}\left(  \Omega\right)  $ characterization result is applicable. This is for
the specific purpose of characterizing the local data space of the tensor
product central difference weight functions which will be studied in Chapter
\ref{Ch_cent_diff_wt_fn_ten_prod}. This is useful numerically because it is
easy to choose functions in $W^{m\mathbf{1}}\left(  \Omega\right)  $.

In the last subsection we will discuss the conversion of global pointwise
error estimates to local estimates.

\subsection{The extended class of B-spline weight functions}

Suppose the data region $\Omega\subset\mathbb{R}^{d}$ is bounded and that $w$
is a scaled, extended B-spline weight function. We start with some definitions
and in these definitions we use the notation $\alpha\leq m$ where $m$ is a
non-negative integer, $\alpha$ is a multi-index and
\[
\left\{  \alpha:\alpha\leq m\right\}  =\left\{  \alpha:\alpha_{k}\leq m\text{
}for\text{ }k=1,\ldots,m\right\}  .
\]

See Definition \ref{Def_multi_id} of the Appendix.

Our next important data function result Theorem
\ref{Thm_int_Xow(O)_eq_Hn(O)_dim1} will require some $L^{2}$ Sobolev space
theory based on the comprehensive study \cite{Adams75} by Adams.

\begin{definition}
\label{Def_SobolevSpace}\textbf{Sobolev spaces }$W^{m,\infty}\left(
\Omega\right)  \mathbf{,}W^{m,\infty},W^{m}\left(  \Omega\right)  ,W^{m},$
$H^{m}$, $m=0,1,2,\ldots$

For any open set $\Omega\subset\mathbb{R}^{d}$ and in the sense of
distributions:%
\begin{align*}
W^{m,\infty}\left(  \Omega\right)   & =\left\{  u\in L^{\infty}\left(
\Omega\right)  :D^{\alpha}u\in L^{\infty}\left(  \Omega\right)  \text{
}for\text{ }\left\vert \alpha\right\vert \leq m\right\}  ,\\
W^{m,\infty}  & =W^{m,\infty}\left(  \mathbb{R}^{d}\right)  ,
\end{align*}

with respective norms $\left\Vert u\right\Vert _{m,\infty,\Omega}%
=\sum\limits_{\left\vert \alpha\right\vert \leq m}\binom{m}{\alpha}\left\Vert
D^{\alpha}u\right\Vert _{\infty,\Omega}$ and $\left\Vert u\right\Vert
_{m,\infty}=\left\Vert u\right\Vert _{m,\infty,\mathbb{R}^{d}} $.%
\begin{align*}
W^{m}\left(  \Omega\right)   & =\left\{  u\in L^{2}\left(  \Omega\right)
:D^{\alpha}u\in L^{2}\left(  \Omega\right)  \text{ }for\text{ }\left\vert
\alpha\right\vert \leq m\right\}  ,\\
W^{m}  & =W^{m}\left(  \mathbb{R}^{d}\right)  ,
\end{align*}

with norms $\left\Vert u\right\Vert _{m,2,\Omega}^{2}=\sum\limits_{\left\vert
\alpha\right\vert \leq m}\binom{m}{\alpha}\left\Vert D^{\alpha}u\right\Vert
_{2,\Omega}^{2}$ and $\left\Vert u\right\Vert _{m,2}=\left\Vert u\right\Vert
_{m,2,\mathbb{R}^{d}}$.

Also define%
\[
H^{m}=H^{m}\left(  \mathbb{R}^{d}\right)  =\left\{  u\in L^{2}\left(
\Omega\right)  :\int\left(  1+\left\vert \xi\right\vert ^{2}\right)
^{m}\left\vert \widehat{u}\left(  \xi\right)  \right\vert ^{2}d\xi
<\infty\right\}  .
\]

Using the identity of part 6 of Definition \ref{Def_multi_id} we have: if
$\binom{m}{\alpha}:=\frac{m!}{\alpha!\left(  m-\left\vert \alpha\right\vert
\right)  !}$ then%
\begin{align*}
\int\left(  1+\left\vert \xi\right\vert ^{2}\right)  ^{m}\left\vert
\widehat{u}\left(  \xi\right)  \right\vert ^{2}d\xi & =\int\sum_{\left\vert
\alpha\right\vert \leq m}\binom{m}{\alpha}\xi^{2\alpha}\left\vert
\widehat{u}\left(  \xi\right)  \right\vert ^{2}d\xi\\
& =\sum_{\left\vert \alpha\right\vert \leq m}\binom{m}{\alpha}\int\left\vert
\widehat{D^{\alpha}u}\right\vert ^{2}\\
& =\sum_{\left\vert \alpha\right\vert \leq m}\binom{m}{\alpha}\int\left\vert
D^{\alpha}u\right\vert ^{2}\\
& =\sum_{\left\vert \alpha\right\vert \leq m}\binom{m}{\alpha}\left\Vert
D^{\alpha}u\right\Vert _{L^{2}}^{2}\\
& =\left\Vert u\right\Vert _{m,2}^{2}.
\end{align*}

Here the weight function is $\left(  1+\left\vert \xi\right\vert ^{2}\right)
^{m}$ and we have $W^{m}=H^{m}$ with identical norms. See also Paragraph 7.62
of Adams \& Fournier.
\end{definition}

We will also need the following Sobolev spaces which are described in the
literature as having \textbf{dominating mixed derivatives}:

\begin{definition}
\label{Def_SobolevSpace2}\textbf{Sobolev spaces with dominating mixed
derivatives }$W^{m\mathbf{1}}\left(  \Omega\right)  ,W^{m\mathbf{1}}%
,W_{0}^{m\mathbf{1}}\left(  \Omega\right)  ,$

$W_{0}^{m\mathbf{1}},H^{m\mathbf{1}},$ $m=0,1,2,3,\ldots$

For any open set $\Omega\subset\mathbb{R}^{d}$ and in the sense of
distributions:%
\begin{align*}
W^{m\mathbf{1}}\left(  \Omega\right)   & =\left\{  u\in L^{2}\left(
\Omega\right)  :D^{\alpha}u\in L^{2}\left(  \Omega\right)  \text{ }for\text{
}\alpha\leq m\mathbf{1}\right\}  ,\\
W^{m\mathbf{1}}  & =W^{m\mathbf{1}}\left(  \mathbb{R}^{d}\right)  ,\\
W_{0}^{m\mathbf{1}}\left(  \Omega\right)   & =closure\text{ }of\text{ }%
C_{0}^{\infty}\left(  \Omega\right)  \text{ }in\text{ }W^{m\mathbf{1}}\left(
\Omega\right)  ,\\
W_{0}^{m\mathbf{1}}  & =W_{0}^{m\mathbf{1}}\left(  \mathbb{R}^{d}\right)  ,
\end{align*}

with norms $\left\Vert u\right\Vert _{W^{m\mathbf{1}}\left(  \Omega\right)
}^{2}\equiv\left\Vert u\right\Vert _{m\mathbf{1},2,\Omega}^{2}:=\sum
\limits_{\alpha\leq m\mathbf{1}}\binom{m\mathbf{1}}{\alpha}\left\Vert
D^{\alpha}u\right\Vert _{2,\Omega}^{2}$ and $\left\Vert u\right\Vert
_{W_{0}^{m\mathbf{1}}\left(  \Omega\right)  }=\left\Vert u\right\Vert
_{W^{m\mathbf{1}}\left(  \Omega\right)  }$.

Also define:%
\begin{align*}
W^{m\mathbf{1},\infty}\left(  \Omega\right)   & =\left\{  u\in L^{\infty
}\left(  \Omega\right)  :D^{\alpha}u\in L^{\infty}\left(  \Omega\right)
\text{ }for\text{ }\alpha\leq m\mathbf{1}\right\}  ,\\
W^{m\mathbf{1},\infty}  & =W^{m\mathbf{1},\infty}\left(  \mathbb{R}%
^{d}\right)  ,
\end{align*}

with respective norms $\left\Vert u\right\Vert _{W^{m\mathbf{1},\infty}\left(
\Omega\right)  }^{2}\equiv\left\Vert u\right\Vert _{m\mathbf{1},\infty,\Omega
}^{2}:=\sum\limits_{\alpha\leq m\mathbf{1}}\binom{m\mathbf{1}}{\alpha
}\left\Vert D^{\alpha}u\right\Vert _{\infty,\Omega}^{2}$ and $\left\Vert
u\right\Vert _{m\mathbf{1},\infty}=\left\Vert u\right\Vert _{m\mathbf{1}%
,\infty,\mathbb{R}^{d}}$.

Also define%
\begin{equation}
H^{m\mathbf{1}}\left(  \Omega\right)  :=\left\{  u\in L^{2}:\int\left(
\left(  1+\xi_{1}^{2}\right)  \ldots\left(  1+\xi_{d}^{2}\right)  \right)
^{m}\left\vert \widehat{u}\left(  \xi\right)  \right\vert ^{2}d\xi
<\infty\right\}  ,\label{2.61}%
\end{equation}

with norm $\left\Vert u\right\Vert _{H^{m\mathbf{1}}\left(  \Omega\right)
}^{2}=\left\Vert \left(  \left(  1+\xi_{1}^{2}\right)  \ldots\left(  1+\xi
_{d}^{2}\right)  \right)  ^{m}\widehat{u}\left(  \xi\right)  \right\Vert
_{L^{2}}=\left\Vert \left(  1+\xi.\xi\right)  ^{m\mathbf{1}}\widehat{u}\left(
\xi\right)  \right\Vert _{L^{2}}$, where we have used the

component-wise vector multiplication notation $\xi\centerdot\xi$.

Now using the binomial (multinomial) theorem,%
\begin{align*}
\left\Vert u\right\Vert _{H^{m\mathbf{1}}\left(  \Omega\right)  }^{2}  &
=\int\left(  \left(  1+\xi_{1}^{2}\right)  \ldots\left(  1+\xi_{d}^{2}\right)
\right)  ^{m}\left\vert \widehat{u}\left(  \xi\right)  \right\vert ^{2}d\xi\\
& =\int\left(  1+\xi.\xi\right)  ^{m\mathbf{1}}\left\vert \widehat{u}\left(
\xi\right)  \right\vert ^{2}d\xi\\
& =\int\sum_{\alpha\leq m\mathbf{1}}\binom{m\mathbf{1}}{\alpha}\left(  \xi
.\xi\right)  ^{\alpha}\left\vert \widehat{u}\left(  \xi\right)  \right\vert
^{2}d\xi\\
& =\int\sum_{\alpha\leq m\mathbf{1}}\binom{m\mathbf{1}}{\alpha}\xi^{2\alpha
}\left\vert \widehat{u}\left(  \xi\right)  \right\vert ^{2}d\xi\\
& =\sum_{\alpha\leq m\mathbf{1}}\binom{m\mathbf{1}}{\alpha}\int\left\vert
\widehat{D^{\alpha}u}\right\vert ^{2}\\
& =\sum_{\alpha\leq m\mathbf{1}}\binom{m\mathbf{1}}{\alpha}\int\left\vert
D^{\alpha}u\right\vert ^{2}\\
& =\sum_{\alpha\leq m\mathbf{1}}\binom{m\mathbf{1}}{\alpha}\left\Vert
D^{\alpha}u\right\Vert _{L^{2}}^{2}\\
& =\left\Vert u\right\Vert _{m\mathbf{1},2,\mathbb{R}^{d}}^{2}.
\end{align*}

Thus $H^{m\mathbf{1}}=W^{m\mathbf{1}}$ as sets and the norms are identical.

Note also that when $m\geq1$, Corollary \ref{Cor_Thm_ten_prod_two_wt_fns} and
\ref{2.61} imply $H^{m\mathbf{1}}$ is always a reproducing kernel Hilbert
space with tensor product weight function $\left(  1+\xi_{1}^{2}\right)
^{m}\ldots\left(  1+\xi_{d}^{2}\right)  ^{m}$ and that the weight function
parameter satisfies $\kappa.<\mathbf{1}$.
\end{definition}

\ 

\begin{definition}
\label{Def_SegCondit}\textbf{Segment condition} We say\ a domain $\Omega$
satisfies the segment condition if every $x\in\operatorname*{bdry}\Omega$ has
a neighborhood $U_{x}$ and $y_{x}\in\mathbb{R}^{d}\setminus0$ such that if
$z\in\Omega\cap U_{x}$ then $z+ty_{x}\in\Omega$ for $0<t<1$.

Observe that the segment property is invariant under translations and
dilations of the domain.
\end{definition}

\begin{definition}
\label{Cmb(clos_open_region)}\textbf{Spaces of continuous functions}

In \S 1.28 Adam's and Fournier \cite{AdamFour2003} define $C^{\left(
m\right)  }\left(  \overline{\Omega}\right)  $ by:
\[
C^{\left(  m\right)  }\left(  \overline{\Omega}\right)  =\left\{  u\in
C^{\left(  m\right)  }\left(  \Omega\right)  :D^{\alpha}u\text{ }is\text{
}bounded\text{ }and\text{ }unif.\text{ }contin.\text{ }on\text{ }\Omega\text{
}for\text{ }\left\vert \alpha\right\vert \leq m\right\}  ,
\]

and our analogue is%
\begin{equation}
C^{\left(  m\mathbf{1}\right)  }\left(  \overline{\Omega}\right)  =\left\{
u\in C^{\left(  m\mathbf{1}\right)  }\left(  \Omega\right)  :D^{\alpha}u\text{
}is\text{ }bounded\text{ }and\text{ }unif.\text{ }contin.\text{ }on\text{
}\Omega\text{ }for\text{ }\alpha\leq m\mathbf{1}\right\}  .\label{2.44}%
\end{equation}

We endow them with the respective supremum norms%
\begin{equation}
\left\Vert u\right\Vert _{m,\infty,\Omega}=\max_{\left\vert \alpha\right\vert
\leq m}\sup_{x\in\Omega}\left\vert D^{\alpha}u\left(  x\right)  \right\vert
,\quad\left\Vert u\right\Vert _{m\mathbf{1},\infty,\Omega}=\max_{\alpha\leq
m\mathbf{1}}\sup_{x\in\Omega}\left\vert D^{\alpha}u\left(  x\right)
\right\vert .\label{2.42}%
\end{equation}

These two spaces are Banach spaces. We also define the \textbf{restriction
vector spaces}
\begin{equation}
\left.
\begin{array}
[c]{ll}%
r_{\Omega}C_{B}^{\left(  m\right)  } & :=\left\{  r_{\Omega}u:u\in
C_{B}^{\left(  m\right)  }\right\} \\
r_{\Omega}C_{B}^{\left(  m\mathbf{1}\right)  } & :=\left\{  r_{\Omega}u:u\in
C_{B}^{\left(  m\mathbf{1}\right)  }\right\}
\end{array}
\right\}  .\label{2.43}%
\end{equation}

\end{definition}

\begin{theorem}
If $\Omega$ is a bounded open set then:\smallskip

\begin{enumerate}
\item $r_{\Omega}C^{\left(  m\right)  }\subset C^{\left(  m\right)  }\left(
\overline{\Omega}\right)  $ and $r_{\Omega}C^{\left(  m\mathbf{1}\right)
}\subset C^{\left(  m\mathbf{1}\right)  }\left(  \overline{\Omega}\right)
$.\smallskip

We can now endow $r_{\Omega}C^{\left(  m\right)  }$ and $r_{\Omega}C^{\left(
m\mathbf{1}\right)  }$ with the norms \ref{2.42}. In fact:\smallskip

\item $r_{\Omega}C_{B}^{\left(  m\right)  }\hookrightarrow C^{\left(
m\right)  }\left(  \overline{\Omega}\right)  $ and $r_{\Omega}C_{B}^{\left(
m\mathbf{1}\right)  }\hookrightarrow C^{\left(  m\mathbf{1}\right)  }\left(
\overline{\Omega}\right)  $.\smallskip

\textbf{Using extension operators} we prove:\smallskip

\item If $\Omega$ has the uniform rectangle property of Definition
\ref{Def_UnifRectCondit} then $r_{\Omega}C_{0}^{\left(  m\mathbf{1}\right)
}=C^{\left(  m\mathbf{1}\right)  }\left(  \overline{\Omega}\right)  $ as
Banach spaces.
\end{enumerate}
\end{theorem}

\begin{proof}
\textbf{Part 1} If $u\in C^{\left(  0\right)  }$ then $u$ is continuous on
$\overline{\Omega}$. But $\overline{\Omega}$ is compact so $u$ is uniformly
continuous on $\overline{\Omega}$ and consequently $r_{\Omega}u$ has a unique
extension to $C^{\left(  0\right)  }\left(  \overline{\Omega}\right)  $. Thus
$u\in C^{\left(  m\right)  }$ implies $r_{\Omega}u\in C^{\left(  m\right)
}\left(  \overline{\Omega}\right)  $ and so $r_{\Omega}C^{\left(  m\right)
}\subset C^{\left(  m\right)  }\left(  \overline{\Omega}\right)  $.\medskip

\textbf{Part 2} Clearly $r_{\Omega}C_{B}^{\left(  m\right)  }\hookrightarrow
C^{\left(  m\right)  }\left(  \overline{\Omega}\right)  $ since both spaces
have the same norm. A similar argument applies to $r_{\Omega}C_{B}^{\left(
m\mathbf{1}\right)  }$.\medskip

\textbf{Part 3} Suppose $f\in C^{\left(  m\mathbf{1}\right)  }\left(
\overline{\Omega}\right)  $. From Theorem
\ref{Thm_ExtenContinFuncs_OrthantProp} there exists a linear extension mapping
$\mathcal{E}_{\Omega}:C^{\left(  m\mathbf{1}\right)  }\left(  \Omega\right)
\rightarrow C_{B}^{\left(  m\mathbf{1}\right)  }$ so%
\[
f\in C^{\left(  m\mathbf{1}\right)  }\left(  \overline{\Omega}\right)
\Rightarrow\mathcal{E}_{\Omega}f\in C^{\left(  m\mathbf{1}\right)
}\Rightarrow r_{\Omega}\mathcal{E}_{\Omega}f=f\Rightarrow f\in r_{\Omega
}C^{\left(  m\mathbf{1}\right)  }.
\]

\end{proof}

??? \textbf{SORT\ OUT}! \textbf{QUOTIENT NORM}, $r_{\Omega}C_{B}^{\left(
m\mathbf{1}\right)  }=r_{\Omega}C_{0}^{\left(  m\mathbf{1}\right)  }%
=C_{0}^{\left(  m\mathbf{1}\right)  }\diagup C_{0;\Omega^{c}}^{\left(
m\mathbf{1}\right)  }$?

and%
\[
\left.
\begin{array}
[c]{ll}%
r_{\Omega}C_{0}^{\left(  m\right)  } & =\left\{  r_{\Omega}u:u\in
C_{0}^{\left(  m\right)  }\right\} \\
r_{\Omega}C_{0}^{\left(  m\mathbf{1}\right)  } & =\left\{  r_{\Omega}u:u\in
C_{0}^{\left(  m\mathbf{1}\right)  }\right\}
\end{array}
\right\}  .
\]

We also endow them with the norms \ref{2.42}.

?? Restriction norms:%
\begin{align*}
\left\Vert f\right\Vert _{m,\infty,\Omega}^{B}  & =\min\left\{  \left\Vert
g\right\Vert _{m,\infty}:r_{\Omega}g=r_{\Omega}f,\quad g\in C_{B}^{\left(
m\right)  }\right\}  ,\\
\left\Vert f\right\Vert _{m,\infty,\Omega}^{0}  & =\min\left\{  \left\Vert
g\right\Vert _{m,\infty}:r_{\Omega}g=r_{\Omega}f,\quad g\in C_{0}^{\left(
m\right)  }\right\}  .
\end{align*}

but do these norms exist?

For any closed set $K=\overline{\Omega}\subset\mathbb{R}^{d}$ define%
\begin{align*}
C_{0;K}^{\left(  m\mathbf{1}\right)  }  & =\left\{  u\in C_{0}^{\left(
m\mathbf{1}\right)  }:\operatorname*{supp}u\subseteq K\right\}  ,\\
C_{0;K}^{\left(  m\right)  }  & =\left\{  u\in C_{0}^{\left(  m\right)
}:\operatorname*{supp}u\subseteq K\right\}  ,
\end{align*}

and since the functions in $C_{0}^{\left(  m\mathbf{1}\right)  }$ are ??
uniformly continuous (but not $C_{B}^{\left(  m\mathbf{1}\right)  }$) these
spaces are Banach spaces when endowed with the usual supremum norms. We now
use the quotient space theory of Chapter \ref{Ch_quot_sp_reprod_kern}.

Since $C_{0;\Omega^{c}}^{\left(  m\right)  }$ is a closed subspace of
$C_{0}^{\left(  m\right)  }$ we form the quotient space $C_{0}^{\left(
m\right)  }/C_{0;\Omega^{c}}^{\left(  m\right)  }$ whose elements are the
equivalence classes $\left[  u\right]  $ given by%
\[
v\symbol{126}u\text{ }iff\text{ }u=v\text{ }on\text{ }\Omega,
\]

and endow it with the restriction norm:
\[
\left\Vert \left[  u\right]  \right\Vert _{m,\infty,\Omega}^{\prime}%
=\min_{v\symbol{126}u}\left\Vert v\right\Vert _{m,\infty,\Omega}.
\]

$C_{0}^{\left(  m\right)  }/C_{0;\Omega^{c}}^{\left(  m\right)  }$ is now a
normed Banach space.

?? FIX! Then from Theorem \ref{Thm_KerSeminorm} the quotient mappings
$R_{\Omega}u=\left[  r_{\Omega}u\right]  $ satisfy%
\[
r_{\Omega}C_{0}^{\left(  m\right)  }=C_{0}^{\left(  m\right)  }/C_{0;\Omega
^{c}}^{\left(  m\right)  },\quad r_{\Omega}C_{0}^{\left(  m\mathbf{1}\right)
}=C_{0}^{\left(  m\mathbf{1}\right)  }/C_{0;\Omega^{c}}^{\left(
m\mathbf{1}\right)  },
\]

and are isometries.

if $C_{B;\Omega^{c}}^{\left(  m\right)  }$ and $C_{B;\Omega^{c}}^{\left(
m\mathbf{1}\right)  }$ are closed. ?? If $\Omega$ is bounded then NO!

and by ??, $r_{\Omega}C_{B}^{\left(  m\right)  }$ and $r_{\Omega}%
C_{B}^{\left(  m\mathbf{1}\right)  }$\ are Banach spaces.

We now introduce the local (restriction) data space corresponding to an
arbitrary data space $X_{w}^{0}$:

\begin{definition}
\label{Def_Xow(open_set)}\textbf{The local data space }$X_{w}^{0}\left(
\Omega\right)  $ Suppose that $w$ is a weight function with property W02 or
W03 and suppose $\Omega\subset\mathbb{R}^{d}$ is open. Define the vector space
$X_{w}^{0}\left(  \Omega\right)  =\left\{  r_{\Omega}u:u\in X_{w}^{0}\right\}
$, where $r_{\Omega}u$ denotes the action of the linear operator which
restricts $u$ to $\Omega$. Endow $X_{w}^{0}\left(  \Omega\right)  $ with the
restriction norm of Theorem \ref{Thm_restrict_rkhs} of the Appendix i.e.
\begin{equation}
\left\Vert v\right\Vert _{w,0;\Omega}=\inf\left\{  \left\Vert u\right\Vert
_{w,0}:u\in X_{w}^{0},\text{ }r_{\Omega}u=v\right\}  ,\label{1.026}%
\end{equation}

and by Theorem \ref{Thm_restrict_rkhs}, $r_{\Omega}:X_{w}^{0}\rightarrow
X_{w}^{0}\left(  \Omega\right)  $ is continuous and $X_{w}^{0}\left(
\Omega\right)  $ is a Banach space.

Further, by the same theorem $\left\Vert \cdot\right\Vert _{w,0;\Omega}$
satisfies the parallelogram law and so $X_{w}^{0}\left(  \Omega\right)  $ is a
Hilbert space.

We know that $X_{w}^{0}$ has reproducing kernel $\Phi\left(  x,y\right)
=\left(  2\pi\right)  ^{-d/2}G\left(  x-y\right)  =R_{y}\left(  x\right)  $.
Hence $X_{w}^{0}\left(  \Omega\right)  $ has the unique Riesz representer
$R_{y}^{\Omega}\left(  x\right)  :=r_{\Omega\times\Omega}R_{y}\left(
x\right)  $.
\end{definition}

\begin{remark}
??? Use BLAH from the loc. paper.

To use the notation $X_{w}^{0}\left(  \Omega\right)  $\ or $X_{w}^{0}\left(
\overline{\Omega}\right)  $? Grisvard uses $H^{s}\left(  \overline{\Omega
}\right)  :=r_{\Omega}H^{s}$ in Definition 1.3.2.4 of \cite{Grisvard85}. Wloka
uses $H^{s}\left(  \Omega\right)  :=r_{\Omega}H^{s}$. To be consistent with
the definition $C^{\infty}\left(  \overline{\Omega}\right)  =r_{\Omega
}C^{\infty}$ we should perhaps use $X_{w}^{0}\left(  \overline{\Omega}\right)
=r_{\Omega}X_{w}^{0}$ but if the boundary is such that an extension of
$X_{w}^{0}\left(  \Omega\right)  $ to $X_{w}^{0}$ exists then by Theorem
\ref{Thm_restrict_sp_homeo_quot_space}, $X_{w}^{0}\left(  \Omega\right)  $ is
homeomorphic to $X_{w}^{0}/\left(  X_{w}^{0}\right)  _{\Omega^{c}}$ and in PDE
theory $\Omega$ is used and not $\overline{\Omega}$.
\end{remark}

We now show for $\Omega$ open and bounded there exists an extension operator
based on the adjoint of the restriction operator $r_{\Omega}$. This operator
exists for all local data spaces $X_{w}^{0}\left(  \Omega\right)  $ but no
explicit construction is given.

\begin{theorem}
\label{Thm_canon_exten_op}\textbf{The adjoint extension operator} $r_{\Omega
}^{\ast}$: Suppose $\Omega\subset\mathbb{R}^{d}$ is open and bounded, and the
weight function $w$ has property W02 or W03. Then the Hilbert space adjoint
operator $r_{\Omega}^{\ast}:X_{w}^{0}\left(  \Omega\right)  \rightarrow
X_{w}^{0}$ is a continuous extension with $\left\Vert r_{\Omega}\right\Vert
_{op}=\left\Vert r_{\Omega}^{\ast}\right\Vert _{op}=1$ and this is the minimal
extension operator norm.

Further,%
\begin{equation}
r_{\Omega}r_{\Omega}^{\ast}R_{y}^{\Omega}=R_{y}^{\Omega},\label{1.084}%
\end{equation}

where $R_{y}^{\Omega}$ is the Riesz representer of the evaluation functional
on $X_{w}^{0}\left(  \Omega\right)  $ i.e. $R_{y}^{\Omega}=r_{\Omega}R_{y}$
and $R_{y}\left(  x\right)  =\left(  2\pi\right)  ^{-d/2}G\left(  x-y\right)
$ where $x,y\in\Omega$ and $G$ is the basis function.

Finally, $r_{\Omega}^{\ast}$ commutes with translations $\tau_{c}$ and
dilations $\sigma_{\lambda}$ in the sense that
\begin{equation}
\tau_{c}r_{\Omega}^{\ast}=r_{\tau_{c}\Omega}^{\ast}\tau_{c},\text{\quad}%
\sigma_{\lambda}r_{\Omega}^{\ast}=r_{\lambda\Omega}^{\ast}\sigma_{\lambda
}.\label{1.085}%
\end{equation}

\end{theorem}

\begin{proof}
From the definition of $X_{w}^{0}\left(  \Omega\right)  $ we have $r_{\Omega
}:X_{w}^{0}\overset{contin}{\longrightarrow}X_{w}^{0}\left(  \Omega\right)  $
and so $r_{\Omega}^{\ast}:X_{w}^{0}\left(  \Omega\right)
\overset{contin}{\longrightarrow}X_{w}^{0}$. Now for $y\in\Omega$ and all $z$,%
\begin{align*}
\left(  r_{\Omega}^{\ast}r_{\Omega}R_{y}\right)  \left(  z\right)   & =\left(
r_{\Omega}^{\ast}r_{\Omega}R_{y},R_{z}\right)  _{w,0}=\left(  r_{\Omega}%
R_{y},r_{\Omega}R_{z}\right)  _{w,0;\Omega}=\\
& =\left(  R_{y}^{\Omega},R_{z}^{\Omega}\right)  _{w,0;\Omega}=\overline
{R_{z}^{\Omega}\left(  y\right)  }=\overline{R_{z}\left(  y\right)  }%
=R_{y}\left(  z\right)  ,
\end{align*}

so that
\begin{equation}
r_{\Omega}^{\ast}R_{y}^{\Omega}=R_{y},\label{1.087}%
\end{equation}

and hence $r_{\Omega}r_{\Omega}^{\ast}R_{y}^{\Omega}=r_{\Omega}R_{y}%
=R_{y}^{\Omega}$ as claimed. Define%
\begin{equation}
\left.
\begin{array}
[c]{ll}%
W_{G}^{\Omega} & :=\left\{  \sum\limits_{\substack{finite \\sum}}\alpha
_{k}G\left(  \cdot-x^{\left(  k\right)  }\right)  :x^{\left(  k\right)  }%
\in\Omega\right\}  =\left\{  \sum\limits_{\substack{finite \\sum}}\alpha
_{k}R_{x^{\left(  k\right)  }}:x^{\left(  k\right)  }\in\Omega\right\}  .\\
W_{G}^{\Omega}\left(  \Omega\right)  & :=r_{\Omega}W_{G}^{\Omega}.
\end{array}
\right\} \label{1.088}%
\end{equation}

A standard result is that a subspace $Z$ of a Hilbert space is dense iff its
orthogonal complement $Z^{\bot}$ is $\left\{  0\right\}  $. Hence in general
$W_{G}^{\Omega}$ is not dense in $X_{w}^{0}$ because $\left(  W_{G}^{\Omega
}\right)  ^{\bot}=\left\{  u\in X_{w}^{0}:r_{\Omega}u=0\right\}  $ e.g. when
$G$ has bounded support. But is the restriction space $W_{G}^{\Omega}\left(
\Omega\right)  $ always dense in $X_{w}^{0}\left(  \Omega\right)  $? We aim to
prove that $W_{G}^{\Omega}\left(  \Omega\right)  ^{\bot}=\left\{  0\right\}
$. Indeed, if $F\in X_{w}^{0}\left(  \Omega\right)  $ then $\left(
F,U\right)  _{w,0;\Omega}=0$ $\forall U\in W_{G}^{\Omega}\left(
\Omega\right)  $ implies $\left(  F,R_{y}^{\Omega}\right)  _{w,0;\Omega}=0$
$\forall y\in\Omega$ i.e. $F\left(  y\right)  =0$ $\forall y\in\Omega$ which
means that $F=0$. Thus
\begin{equation}
W_{G}^{\Omega}\left(  \Omega\right)  \text{ }is\text{ }dense\text{ }in\text{
}X_{w}^{0}\left(  \Omega\right)  ,\label{1.089}%
\end{equation}

and we observe that \ref{1.084} implies $r_{\Omega}^{\ast}$ is an extension
operator when applied to $r_{\Omega}W_{G}^{\Omega}$.

Next choose arbitrary $U\in X_{w}^{0}\left(  \Omega\right)  $. Then there must
exist $\left\{  U_{k}\right\}  \subset W_{G}^{\Omega}\left(  \Omega\right)  $
such that $U_{k}\rightarrow U$. Hence $r_{\Omega}^{\ast}U_{k}\rightarrow
r_{\Omega}^{\ast}U$ and so $U_{k}=r_{\Omega}r_{\Omega}^{\ast}U_{k}\rightarrow
r_{\Omega}r_{\Omega}^{\ast}U$. But $U_{k}\rightarrow U$ so $U=r_{\Omega
}r_{\Omega}^{\ast}U$ and therefore $r_{\Omega}^{\ast}$ is an extension operator.

Clearly
\[
\left\Vert r_{\Omega}^{\ast}U\right\Vert _{w,0;\Omega}^{2}=\left(  r_{\Omega
}^{\ast}U,r_{\Omega}^{\ast}U\right)  _{w,0}=\left(  U,r_{\Omega}r_{\Omega
}^{\ast}U\right)  _{w,0;\Omega}=\left\Vert U\right\Vert _{w,0;\Omega}^{2},
\]

so that $\left\Vert r_{\Omega}^{\ast}\right\Vert _{op}=1$. That $\left\Vert
r_{\Omega}^{\ast}\right\Vert _{op}=\left\Vert r_{\Omega}\right\Vert _{op}$ is
an elementary property of adjoints. Further, if $E$ is any continuous
extension then $\left\Vert U\right\Vert _{w,0;\Omega}\leq\left\Vert
EU\right\Vert _{w,0}$ which implies that $\left\Vert E\right\Vert _{op}\geq1$.

If $f\in X_{w}^{0}$ and $U\in X_{w}^{0}\left(  \Omega\right)  $ then%
\begin{align*}
\left(  \tau_{c}r_{\Omega}^{\ast}U,f\right)  _{w,0}=\left(  r_{\Omega}^{\ast
}U,\tau_{-c}f\right)  _{w,0}  & =\left(  U,r_{\Omega}\tau_{-c}f\right)
_{w,0;\Omega}=\left(  U,\tau_{-c}r_{\tau_{c}\Omega}f\right)  _{w,0;\Omega}=\\
& =\left(  \tau_{c}U,r_{\tau_{c}\Omega}f\right)  _{w,0;\tau_{c}\Omega}=\left(
r_{\tau_{c}\Omega}^{\ast}\tau_{c}U,f\right)  _{w,0},
\end{align*}

and%
\begin{align*}
\left(  \sigma_{\lambda}r_{\Omega}^{\ast}U,f\right)  _{w,0}=\left(  r_{\Omega
}^{\ast}U,\sigma_{\lambda}f\right)  _{w,0}  & =\left(  U,r_{\Omega}%
\sigma_{\lambda}f\right)  _{w,0;\Omega}=\left(  U,\sigma_{\lambda}%
r_{\lambda\Omega}f\right)  _{w,0;\Omega}=\\
& =\left(  \sigma_{\lambda}U,r_{\lambda\Omega}f\right)  _{w,0;\lambda\Omega
}=\left(  r_{\lambda\Omega}^{\ast}\sigma_{\lambda}U,f\right)  _{w,0}.
\end{align*}

\end{proof}

Recall that Definition \ref{Def_XowK} introduced the Hilbert spaces
$X_{w:K}^{0}:=\left(  X_{w}^{0}\right)  _{K}:=\left\{  f\in X_{w}%
^{0}:\operatorname*{supp}f\subseteq K\right\}  $.

\begin{corollary}
\label{Cor1_Thm_canon_exten_op}Suppose $\Omega\subset\mathbb{R}^{d}$ is open.
Then $X_{w}^{0}\left(  \Omega^{c}\right)  $ is a Hilbert space and the
equivalence class $u\sim v$ iff $u=v$ on $\Omega$ generates an isometric
isomorphism:%
\[
X_{w}^{0}\left(  \Omega\right)  \approx X_{w}^{0}/\left(  X_{w}^{0}\right)
_{\Omega^{c}}.
\]

\end{corollary}

\begin{proof}
The existence of the continuous extension $r^{\ast}:X_{w}^{0}\left(
\Omega\right)  \rightarrow X_{w}^{0}$ allows us to apply Theorem
\ref{Thm_restrict_sp_homeo_quot_space}.
\end{proof}

\begin{corollary}
\label{Cor2_Thm_canon_exten_op}Under the conditions of Theorem
\ref{Thm_canon_exten_op}, $\operatorname*{supp}r_{\Omega}^{\ast}X_{w}%
^{0}\left(  \Omega\right)  \subset\overline{\Omega}+\operatorname*{supp}G$ and
$r_{\Omega}^{\ast}:X_{w}^{0}\left(  \Omega\right)  \rightarrow\left(
X_{w}^{0}\right)  _{\overline{\Omega}+\operatorname*{supp}G}$ is continuous.
\end{corollary}

\begin{proof}
From \ref{1.087}, $\operatorname*{supp}r_{\Omega}^{\ast}W_{G}^{\Omega}\left(
\Omega\right)  \subset\overline{\Omega}+\operatorname*{supp}G$. Choose $F\in
X_{w}^{0}\left(  \Omega\right)  $. Then by \ref{1.089} there exists a sequence
$\left\{  F_{k}\right\}  $ in $W_{G}^{\Omega}\left(  \Omega\right)  $ which
converges to $F$. We have $\operatorname*{supp}r_{\Omega}^{\ast}F_{k}%
\subset\overline{\Omega}+\operatorname*{supp}G$.

Suppose $\overline{\Omega}+\operatorname*{supp}G\neq\mathbb{R}^{d}$ and choose
$x\notin\overline{\Omega}+\operatorname*{supp}G$. Then%
\[
\left(  r_{\Omega}^{\ast}F\right)  \left(  x\right)  =\left(  r_{\Omega}%
^{\ast}F,R_{x}\right)  _{w,0}=\left(  F,R_{x}^{\Omega}\right)  _{w,0;\Omega
}=\lim\limits_{k\rightarrow\infty}\left(  F_{k},R_{x}^{\Omega}\right)
_{w,0;\Omega}=\lim\limits_{k\rightarrow\infty}F_{k}\left(  x\right)  =0,
\]

and so $\operatorname*{supp}r_{\Omega}^{\ast}F\subset\overline{\Omega
}+\operatorname*{supp}G$.
\end{proof}

\begin{remark}
\label{Rem_SobolevSpace2}\textbf{For certain proofs of the results below} see
Chapter \ref{Ch_Hm1_proofs_Adams} in the Appendix. The work
\cite{AdamFour2003} by Adams and Fournier is devoted to Sobolev space theory
and a review of this study shows that many of it's results can be extended to
the spaces $W^{m\mathbf{1}}\left(  \Omega\right)  $, which are clearly the
"standard" isotropic Sobolev spaces $W^{m}\left(  \Omega\right)  $ in one
dimension. In particular, \textbf{assuming }$\Omega$\textbf{\ is a bounded
region}:

\begin{enumerate}
\item \textbf{Theorem 3.3} can be modified to prove that $W^{m\mathbf{1}%
}\left(  \Omega\right)  $ is a Banach space.

\item \textbf{Theorem 3.6} can be modified to prove that $W^{m\mathbf{1}%
}\left(  \Omega\right)  $ is a separable, reflexive, uniformly convex Hilbert
space with inner products defined in terms of their norms.

\item \textbf{Lemma 3.16} can be easily modified to handle \textbf{mollifiers}
of $W^{m\mathbf{1}}\left(  \Omega\right)  $.

\item Regarding \textbf{approximation by smooth functions}: if $\Omega$ has
the \textbf{segment property} then the restriction space $r_{\Omega}C^{\infty
}$ is dense in $W^{m\mathbf{1}}\left(  \Omega\right)  $ and $r_{\Omega}%
C_{0}^{\infty}=r_{\Omega}C^{\infty}=r_{\Omega}C_{B}^{\infty}$ as sets.
(modified \textbf{Theorem 3.22} - see Theorem \ref{Thm_3.22_AdamsF}).

Also, for \textbf{any open set} $\Omega$, $C^{\left(  m\mathbf{1}\right)
}\left(  \Omega\right)  \cap W^{m\mathbf{1}}\left(  \Omega\right)  $ is dense
in $W^{m\mathbf{1}}\left(  \Omega\right)  $ (modified \textbf{Theorem 3.17} -
see Theorem \ref{Thm_3.17_AdamsF}).

\item Two \textbf{extension operators} are available:

\begin{enumerate}
\item \textbf{The adjoint data space extension }$\mathcal{E}_{\Omega
}:=r_{\Omega}^{\ast}$ \textbf{which\ maps }$X_{w}^{0}\left(  \Omega\right)  $
\textbf{into} $X_{w}^{0}$.

In Theorem \ref{Thm_canon_exten_op} it was shown that for an arbitrary data
space $X_{w}^{0}$ the adjoint $r_{\Omega}^{\ast}$ of the restriction operator
$r_{\Omega}:$ $X_{w}^{0}\rightarrow X_{w}^{0}\left(  \Omega\right)  $ is a
continuous extension operator $r_{\Omega}^{\ast}:X_{w}^{0}\left(
\Omega\right)  \rightarrow X_{w}^{0}$ and some relevant properties proved. It
was noted in Definition \ref{Def_SobolevSpace2} that $W^{m\mathbf{1}}$ is a
data space for all integer $m\geq1$. From Corollary
\ref{Cor2_Thm_canon_exten_op}, $\operatorname*{supp}r_{\Omega}^{\ast}X_{w}%
^{0}\left(  \Omega\right)  \subset\overline{\Omega}+\operatorname*{supp}G$ and
$r_{\Omega}^{\ast}:X_{w}^{0}\left(  \Omega\right)  \rightarrow\left(
X_{w}^{0}\right)  _{\overline{\Omega}+\operatorname*{supp}G}$ is continuous.

\item \textbf{A constructive extension operator} $\mathcal{E}_{\Omega
}^{m\mathbf{1}}$\textbf{\ which maps }$W^{m\mathbf{1}}\left(  \Omega\right)  $
\textbf{into} $W^{m\mathbf{1}}$.

From Theorem 4.26 of Adams \& Fournier \cite{AdamFour2003} there exists a
continuous extension operator $\mathcal{E}_{\Omega}^{\prime}:W^{m}\left(
\Omega\right)  \rightarrow W^{m}$ \textit{which is constructed using the
\textbf{method of reflections}} and assumes $\Omega$ has the \textbf{uniform
}$C^{\left(  m\right)  }$\textbf{-regularity property} of Section 4.6. We have
$\mathcal{E}_{\Omega}^{\prime}:W^{k}\left(  \Omega\right)  \rightarrow W^{k}$
is continuous for $k\leq m$. The proof of Theorem 4.26 involves a
\textbf{two-step} construction. The \textbf{first (reflection) step} involves
an extension from the half-space $\Omega=\mathbb{R}_{+}^{d}$ to $\mathbb{R}%
^{d}$ by reflection and it is easy to show that this extension is also
continuous from $W^{m\mathbf{1}}\left(  \mathbb{R}_{+}^{d}\right)  $ to
$W^{m\mathbf{1}}\left(  \mathbb{R}^{d}\right)  $. The \textbf{second
(localization) step} uses the uniform $C^{\left(  m\right)  }$-regularity
property to locally map neighborhoods of the boundary $\partial\Omega$ into
the unit ball $B_{1}$. It turns out that if we try to modify the localization
step by replacing the $C^{\left(  m\right)  }$-regular map by one having
smoothness $C^{\left(  m\mathbf{1}\right)  }$ then we encounter a basic
problem and our attempt fails, as noted in Section \ref{Ch_exten_rect_condit}.

To construct an \textbf{explicit form} for the extension operator
$\mathcal{E}_{\Omega}^{m\mathbf{1}}:W^{m\mathbf{1}}\left(  \mathbb{\Omega
}\right)  \rightarrow W^{m\mathbf{1}}$ it seems we currently need to severely
restrict the boundary of $\Omega$ and impose the \textbf{rectangle property}
introduced in Definition \ref{Def_UnifRectCondit} of Section
\ref{Ch_exten_rect_condit} below. These extensions are constructed in Theorems
\ref{Thm_ExtenOrthantSobolFourier} and \ref{Thm_exten_Hn1(open)_to_Hn1} and
are denoted there by $E_{\Omega}$ and $E_{\Omega}^{n\mathbf{1}}$ respectively.

If $\Omega$ satisfies the rectangle property then $\sigma_{\lambda}\Omega$ and
$\tau_{c}\Omega$ also satisfy the rectangle property. Here the finite covering
of $\Omega$, the corresponding partition of unity etc. are all transformed appropriately.

Note also that the rectangle property implies the segment property.

\item Finally, both these extension operators commute with translations and
dilations in the sense that
\begin{align*}
\tau_{c}\mathcal{E}_{\Omega}  & =\mathcal{E}_{\Omega+c}\tau_{c},\quad
\sigma_{\lambda}\mathcal{E}_{\Omega}=\mathcal{E}_{\lambda\Omega}%
\sigma_{\lambda}.\\
\tau_{c}\mathcal{E}_{\Omega}^{m\mathbf{1}}  & =\mathcal{E}_{\Omega
+c}^{m\mathbf{1}}\tau_{c},\quad\sigma_{\lambda}\mathcal{E}_{\Omega
}^{m\mathbf{1}}=\mathcal{E}_{\lambda\Omega}^{m\mathbf{1}}\sigma_{\lambda}.
\end{align*}

\end{enumerate}

\item Easy modification of Lemma 3.27 of \cite{AdamFour2003}: let
$\mathcal{E}_{0}u$ denote the \textbf{zero extension} of $u\in W_{0}%
^{m\mathbf{1}}\left(  \Omega\right)  $ outside $\Omega$. Then for $\alpha\leq
m\mathbf{1} $, $D^{\alpha}\mathcal{E}_{0}u=\mathcal{E}_{0}D^{\alpha}u$ in
$\mathcal{D}^{\prime}$, and $\mathcal{E}_{0}u\in W^{m\mathbf{1}}\left(
\mathbb{R}^{d}\right)  $ with $\left\Vert \mathcal{E}_{0}u\right\Vert
_{W^{m\mathbf{1}}}=\left\Vert u\right\Vert _{W_{0}^{m\mathbf{1}}\left(
\Omega\right)  }$.

\item For any closed set $K$ which is the \textbf{closure of an open set}
define%
\begin{align*}
C_{K}^{m\mathbf{1}}  & :=\left\{  \phi\in C^{m\mathbf{1}}:\operatorname*{supp}%
\phi\subseteq K\right\}  .\\
C_{K}^{\infty}  & :=\left\{  \phi\in C^{\infty}:\operatorname*{supp}%
\phi\subseteq K\right\}  .
\end{align*}

\item Special case of Theorem \ref{Thm_norm_Cm1xHm1_loc}: If $\phi\in
r_{\Omega}C^{\infty}$ (definition \ref{2.43}) and $u\in W^{m\mathbf{1}}\left(
\Omega\right)  $ then $\phi u\in W^{m\mathbf{1}}\left(  \Omega\right)  $ and%
\[
\left\Vert \phi u\right\Vert _{m\mathbf{1},2,\Omega}\leq2^{md/2}\left\Vert
\phi\right\Vert _{m\mathbf{1},\infty,\Omega}\left\Vert u\right\Vert
_{m\mathbf{1},2,\Omega}.
\]

where $c_{m,d}$ is independent of both $u$ and $\phi$.

\item For any closed set $K\subset\mathbb{R}^{d}$ define (cf. Definition
\ref{Def_XowK} of $X_{w:K}^{0}$)%
\[
W_{K}^{m\mathbf{1}}:=\left\{  u\in W^{m\mathbf{1}}:\operatorname*{supp}%
u\subseteq K\right\}  ,\text{ }\left\Vert u\right\Vert _{W_{K}^{m\mathbf{1}}%
}=\left\Vert u\right\Vert _{W^{m\mathbf{1}}}.
\]

An easy argument (see Theorem \ref{Thm_XowK_Hilbert}) using the fact that
$W^{m\mathbf{1}}$ is a reproducing kernel Hilbert space, shows that
$W_{K}^{m\mathbf{1}}$ is a closed subspace of $W^{m\mathbf{1}}$ and hence a
Hilbert space.

A modification of Theorem 3.7 of Wloka \cite{Wloka87} is: if $\Omega$ has the
segment property then $r_{\Omega}W_{\overline{\Omega}}^{m\mathbf{1}}%
=W_{0}^{m\mathbf{1}}\left(  \Omega\right)  $ and $C_{0}^{\infty}\left(
\Omega\right)  $ is dense in $W_{\overline{\Omega}}^{m\mathbf{1}}$. This is
Theorem \ref{Thm_Thm3.7_Wloka} in the Appendix. In fact, part 6 implies that
$r_{\Omega}:W_{\overline{\Omega}}^{\mathbf{1}}\rightarrow W_{0}^{\mathbf{1}%
}\left(  \Omega\right)  $ is an \textbf{isometric isomorphism} with inverse
$\mathcal{E}_{0}$.

\item If $\Omega$ is a bounded open set then multiplication by a $r_{\Omega
}C_{\overline{\Omega}}^{\infty}$\ function is a continuous mapping from
$W^{m\mathbf{1}}\left(  \Omega\right)  $ to $W_{0}^{m\mathbf{1}}\left(
\Omega\right)  $. In fact, using part 9 we have $r_{\Omega}C_{\overline
{\Omega}}^{\infty}\subset r_{\Omega}W_{\overline{\Omega}}^{m\mathbf{1}}%
=W_{0}^{m\mathbf{1}}\left(  \Omega\right)  \subset W^{m\mathbf{1}}\left(
\Omega\right)  $ and we can apply Theorem \ref{Thm_norm_Cm1xHm1_loc} to get%
\begin{equation}
\left\Vert \phi u\right\Vert _{W_{0}^{m\mathbf{1}}\left(  \Omega\right)  }%
\leq2^{md/2}\left\Vert \phi\right\Vert _{W^{m\mathbf{1},\infty}\left(
\Omega\right)  }\left\Vert u\right\Vert _{W^{m\mathbf{1}}\left(
\Omega\right)  },\quad\phi\in r_{\Omega}C_{\overline{\Omega}}^{\infty},\text{
}u\in W^{m\mathbf{1}}\left(  \Omega\right)  .\label{2.10}%
\end{equation}

\item \textbf{Grisvard} \cite{Grisvard85} in subsection 1.3.2 gives three
different definitions of Sobolev spaces:

\textbf{Definition 1.3.2.1} is vanilla $W_{p}^{s}\left(  \Omega\right)  $ and
Definition 1.3.2.2. is $\overset{\circ}{W}_{p}^{s}\left(  \Omega\right)
=\operatorname*{cl}C_{0}^{\infty}\left(  \Omega\right)  $.

\textbf{Definition 1.3.2.4} introduces $W_{p}^{s}\left(  \overline{\Omega
}\right)  :=r_{\Omega}W_{p}^{s}$. Grisvard says that H\"{o}rmander
\cite{Hormand63} uses this approach but I checked and he does not - he uses
the \textbf{local }space which exploits the fact that $\phi\in C_{0}^{\infty
}\left(  \Omega\right)  $ and $u\in W_{p}^{s}$ implies $\phi u\in W_{p}%
^{s}\left(  \Omega\right)  $. However, \textbf{Wloka} \cite{Wloka87} does use
this approach in Section 1.5 -\ see part 12.

\textbf{Definition 1.3.2.5} introduces $\widetilde{W}_{p}^{s}\left(
\Omega\right)  :=\left\{  u\in W_{p}^{s}\left(  \Omega\right)  :E_{0}u\in
W_{p}^{s}\left(  \mathbb{R}^{d}\right)  \right\}  $.

\item \textbf{Wloka} \cite{Wloka87}: Avoids interpolation spaces by defining
in Definition 3.1, $W_{2}^{s}\left(  \Omega\right)  $, $s\geq0$, using the
elementary Slobodecki\u{\i} double integral difference method. Definition 3.2
introduces $\overset{\circ}{W}_{2}^{s}\left(  \Omega\right)  :=closure$
$C_{0}^{\infty}\left(  \Omega\right)  $.

Suppose $\Omega$ has the segment property. In Theorem 3.6 he shows that
$r_{\Omega}C_{0}^{\infty}\overset{d}{\hookrightarrow}W_{2}^{s}\left(
\Omega\right)  $ and in Theorem 3.7 he shows that $\overset{\circ}{W}_{2}%
^{s}\left(  \Omega\right)  =\left(  W_{2}^{s}\right)  _{\overline{\Omega}%
}:=r_{\Omega}\left\{  f\in W_{2}^{s}:\operatorname*{supp}f\subseteq
\overline{\Omega}\right\}  $.

Definition 5.1 defines $H^{s}\left(  \mathbb{R}^{d}\right)  $, $s\in
\mathbb{R}^{1}$, using the Fourier transform method. Definition 5.2 defines
$H^{s}\left(  \Omega\right)  :=r_{\Omega}H^{s}\left(  \mathbb{R}^{d}\right)  $
with the restriction norm which is an inner product. Theorem 5.2 shows
$W_{2}^{s}\simeq H^{s}$ i.e. equal as sets and with equivalent norms. Lemma
5.1 says $r_{\Omega}C_{0}^{\infty}\overset{d}{\hookrightarrow}H^{s}\left(
\Omega\right)  $. If we define $\overset{\circ}{H}^{s}\left(  \Omega\right)
:=closure$ $C_{0}^{\infty}\left(  \Omega\right)  $ then $\overset{\circ
}{H}^{s}\left(  \Omega\right)  \simeq\overset{\circ}{W}_{2}^{s}\left(
\Omega\right)  $ when $s\geq0$. Theorem 5.3 states that if $E_{\Omega}%
:H^{s}\left(  \Omega\right)  \rightarrow H^{s}$ is a continuous extension
operator then $H^{s}\left(  \Omega\right)  \simeq W_{2}^{s}\left(
\Omega\right)  $ when $s\geq0$.

\item \textbf{Tr\`{e}ves} \cite{Trev75} defines $H^{s}\left(  K\right)  =??$
after Proposition 25.4.

\item In part 12 Wloka defines $H^{s}\left(  \mathbb{R}^{d}\right)  $ using
the Fourier transform on $L^{2}$ and then let $H^{s}\left(  \Omega\right)
:=r_{\Omega}H^{s}\left(  \mathbb{R}^{d}\right)  $ This will motivate the use
of the restriction space notation $X_{w}^{0}\left(  \Omega\right)  =r_{\Omega
}X_{w}^{0}$ in Definition \ref{Def_Xow(open_set)} below. This is because
$X_{w}^{0}$ is defined globally using the Fourier transform on $\mathbb{R}^{d}
$.
\end{enumerate}
\end{remark}

I next present further results concerning the $W^{m\mathbf{1}}\left(
\Omega\right)  $ spaces. Observe that the proof of part 4 uses both $L^{1}$
and $L^{2}$ Fourier transform theory and introduces an equivalent Fourier
transform norm for $W^{m\mathbf{1}}$.

\begin{lemma}
\label{Lem_SobolevSpProperty}The spaces $W^{m\mathbf{1}}\left(  \Omega\right)
$ and $W_{0}^{m\mathbf{1}}\left(  \Omega\right)  $ have the following properties:

\begin{enumerate}
\item The \textbf{translation} operator $\tau_{c}$ is an isometric isomorphism
from $W^{m\mathbf{1}}\left(  \Omega\right)  $ to $W^{m\mathbf{1}}\left(
\tau_{c}\Omega\right)  $.

The \textbf{dilation} operator $\sigma_{\lambda}$ is a homeomorphism from
$W^{m\mathbf{1}}\left(  \Omega\right)  $ to $W^{m\mathbf{1}}\left(
\sigma_{\lambda}\Omega\right)  $.

\item (\textbf{Copy} of part 10 of Remark \ref{Rem_SobolevSpace2}) If $\Omega$
is a bounded open set then \textbf{multiplication by a }$r_{\Omega
}C_{\overline{\Omega}}^{\infty}$ \textbf{function} is a continuous mapping
from $W^{m\mathbf{1}}\left(  \Omega\right)  $ to $W_{0}^{m\mathbf{1}}\left(
\Omega\right)  $. In fact
\begin{equation}
\left\Vert \phi u\right\Vert _{W_{0}^{m\mathbf{1}}\left(  \Omega\right)  }%
\leq2^{md/2}\left\Vert \phi\right\Vert _{W^{m\mathbf{1},\infty}\left(
\Omega\right)  }\left\Vert u\right\Vert _{W^{m\mathbf{1}}\left(
\Omega\right)  },\quad\phi\in r_{\Omega}C_{\overline{\Omega}}^{\infty},\text{
}u\in W^{m\mathbf{1}}\left(  \Omega\right)  ,\label{1.018}%
\end{equation}

\item \textbf{The extensions} $\mathcal{E}_{\Omega;\varepsilon}:W^{m\mathbf{1}%
}\left(  \Omega\right)  \rightarrow W_{0}^{m\mathbf{1}}\left(  \Omega
_{\varepsilon}\right)  $:

\begin{enumerate}
\item If $\Omega\subset\mathbb{R}^{d}$ has the \textbf{rectangle property}
then given $\varepsilon>0$ we can \textbf{construct} an explicit continuous
linear extension mapping $\mathcal{E}_{\Omega;\varepsilon}:W^{m\mathbf{1}%
}\left(  \Omega\right)  \rightarrow W_{0}^{m\mathbf{1}}\left(  \Omega
_{\varepsilon}\right)  $ using the operator $\mathcal{E}_{\Omega}%
^{m\mathbf{1}}$.

\item Suppose the region $\Omega$ is bounded and satisfies the \textbf{segment
condition}, and $\varepsilon>0$. Then there \textbf{exists} a continuous
linear extension mapping $\mathcal{E}_{\Omega;\varepsilon}:W^{m\mathbf{1}%
}\left(  \Omega\right)  \rightarrow W_{0}^{m\mathbf{1}}\left(  \Omega
_{\varepsilon}\right)  $ constructed using $\mathcal{E}_{\Omega}$.

\item Both these extensions commute with translations and dilations in the
sense that
\[
\tau_{c}\mathcal{E}_{\Omega;\varepsilon}=\mathcal{E}_{\Omega+c;\varepsilon
}\tau_{c},\quad\sigma_{\lambda}\mathcal{E}_{\Omega;\varepsilon}=\mathcal{E}%
_{\lambda\Omega;\lambda\varepsilon}\sigma_{\lambda}.
\]

\end{enumerate}

\item \textbf{Global embedding} If
\begin{equation}
w_{m\mathbf{1}}\left(  \xi\right)  =\left(  1+\xi_{1}^{2}\right)  ^{m}%
\ldots\left(  1+\xi_{d}^{2}\right)  ^{m},\quad m=1,2,3,\ldots,\label{1.052}%
\end{equation}

then $\left\Vert 1/w_{m\mathbf{1}}\right\Vert _{1}=\left(  2^{-2\left(
m-1\right)  }\tbinom{2m-2}{m-1}\pi\right)  ^{d}$ and $w_{m\mathbf{1}}$ is a
weight function with property W03 for $\kappa<\left(  m-\frac{1}{2}\right)
\mathbf{1}$. We have $W^{m\mathbf{1}}=X_{w_{m\mathbf{1}}}^{0}\hookrightarrow
C_{B}^{\left(  m-1\right)  \mathbf{1}}$ where $W^{m\mathbf{1}}$ and
$X_{w_{m\mathbf{1}}}^{0}$ have equivalent norms.

\item \textbf{Local embedding} $X_{w_{m\mathbf{1}}}^{0}\left(  \Omega\right)
\hookrightarrow r_{\Omega}C_{B}^{\left(  m-1\right)  \mathbf{1}}$.

\item \textbf{Poincar\'{e}'s inequality} Suppose $\Omega$ is contained in the
open rectangle $R\left(  a,b\right)  $ and $u\in W_{0}^{n\mathbf{1}}\left(
\Omega\right)  $. Then $\left\Vert D^{m\mathbf{1}}u\right\Vert _{L^{2}\left(
\Omega\right)  }$ is an \textbf{equivalent norm on} $W_{0}^{m\mathbf{1}%
}\left(  \Omega\right)  $. In fact%
\[
\left\Vert D^{\alpha}u\right\Vert _{L^{2}\left(  \Omega\right)  }\leq\left(
\tfrac{b-a}{\sqrt{2}}\right)  ^{m\mathbf{1}-\alpha}\left\Vert D^{m\mathbf{1}%
}u\right\Vert _{L^{2}\left(  \Omega\right)  },\quad\alpha\leq m\mathbf{1}%
,\text{ }u\in W_{0}^{m\mathbf{1}}\left(  \Omega\right)  .
\]

\end{enumerate}
\end{lemma}

\begin{proof}
\textbf{Part 1} $W^{m\mathbf{1}}\left(  \Omega\right)  $ is closed under
translations and dilations, the inverses of translations and dilations are
translations and dilations respectively, and the continuity calculations are
simple changes of variable.\medskip

\textbf{Part 2} Remark \ref{Rem_SobolevSpace2}.\medskip

\textbf{Part 3a} Since $\Omega$ satisfies the uniform rectangle property, part
5a of Remark \ref{Rem_SobolevSpace2} means there exists a continuous extension
operator $\mathcal{E}_{\Omega}^{m\mathbf{1}}:W^{m\mathbf{1}}\left(
\Omega\right)  \rightarrow W^{m\mathbf{1}}\left(  \mathbb{R}^{d}\right)  $.
Using Lemma \ref{Lem_func_eq_1_nbhd_set} we next construct a truncation
function $\phi_{\Omega;\varepsilon}\in C_{\overline{\Omega}_{\varepsilon}%
}^{\infty}$ such that $0\leq\phi_{\Omega;\varepsilon}\leq1$, $\phi
_{\Omega;\varepsilon}=1$ on $\Omega_{\varepsilon/3}$ and $\operatorname*{supp}%
\phi_{\Omega;\varepsilon}\subseteq\overline{\Omega}_{\varepsilon}$. Thus
$r_{\Omega_{\varepsilon}}\phi_{\Omega;\varepsilon}\in r_{\Omega_{\varepsilon}%
}C_{\overline{\Omega}_{\varepsilon}}^{\infty}$ and we define the mapping
$\mathcal{E}_{\Omega;\varepsilon}$ by
\begin{equation}
\mathcal{E}_{\Omega;\varepsilon}u:=r_{\Omega_{\varepsilon}}\left(
\phi_{\Omega;\varepsilon}\mathcal{E}_{\Omega}^{m\mathbf{1}}u\right)  ,\quad
u\in W^{m\mathbf{1}}\left(  \Omega\right)  .\label{1.40}%
\end{equation}

Now by part 5 of Remark \ref{Rem_SobolevSpace2} the rectangle condition
implies the segment property and so by part 2,%
\begin{align*}
\left\Vert \mathcal{E}_{\Omega;\varepsilon}u\right\Vert _{W_{0}^{m\mathbf{1}%
}\left(  \Omega_{\varepsilon}\right)  }  & =\left\Vert r_{\Omega_{\varepsilon
}}\left(  \phi_{\Omega;\varepsilon}\mathcal{E}_{\Omega}^{m\mathbf{1}}u\right)
\right\Vert _{W_{0}^{m\mathbf{1}}\left(  \Omega_{\varepsilon}\right)  }\\
& \leq2^{md/2}\left\Vert r_{\Omega_{\varepsilon}}\phi_{\Omega;\varepsilon
}\right\Vert _{W^{m\mathbf{1},\infty}\left(  \Omega_{\varepsilon}\right)
}\left\Vert \mathcal{E}_{\Omega}^{m\mathbf{1}}u\right\Vert _{W^{m\mathbf{1}%
}\left(  \Omega_{\varepsilon}\right)  }\\
& =2^{md/2}\left\Vert \phi_{\Omega;\varepsilon}\right\Vert _{W^{m\mathbf{1}%
,\infty}\left(  \Omega_{\varepsilon}\right)  }\left\Vert \mathcal{E}_{\Omega
}^{m\mathbf{1}}u\right\Vert _{W^{m\mathbf{1}}}\\
& \leq2^{md/2}\left\Vert \phi_{\Omega;\varepsilon}\right\Vert _{W^{m\mathbf{1}%
,\infty}\left(  \Omega_{\varepsilon}\right)  }\left\Vert \mathcal{E}_{\Omega
}^{m\mathbf{1}}\right\Vert _{op}\left\Vert u\right\Vert _{W^{m\mathbf{1}%
}\left(  \Omega\right)  }.
\end{align*}
\medskip

\textbf{Part 3b} From Definition \ref{Def_SobolevSpace2} we know that
$W^{m\mathbf{1}}$ is a data function space and hence by Theorem
\ref{Thm_canon_exten_op} there exists a continuous extension operator
$\mathcal{E}_{\Omega}:=r_{\Omega}^{\ast}$, $r_{\Omega}^{\ast}:W^{m\mathbf{1}%
}\left(  \Omega\right)  \rightarrow W^{m\mathbf{1}}$ with $\left\Vert
r_{\Omega}^{\ast}\right\Vert _{op}=1$. Following part 3a we define our
extension by \ref{1.40} and by the successive use of Remark
\ref{Rem_SobolevSpace2} (which uses the segment property), part 2 of Lemma
\ref{Lem_SobolevSpProperty} and then Theorem \ref{Thm_canon_exten_op} we
obtain the sequence of continuity estimates%
\begin{align*}
\left\Vert \mathcal{E}_{\Omega;\varepsilon}u\right\Vert _{W_{0}^{m\mathbf{1}%
}\left(  \Omega_{\varepsilon}\right)  }  & =\left\Vert r_{\Omega_{\varepsilon
}}\left(  \phi_{\Omega;\varepsilon}\mathcal{E}_{\Omega}u\right)  \right\Vert
_{W_{0}^{m\mathbf{1}}\left(  \Omega_{\varepsilon}\right)  }=\left\Vert
\phi_{\Omega;\varepsilon}\mathcal{E}_{\Omega}u\right\Vert _{W_{\overline
{\Omega}_{\varepsilon}}^{m\mathbf{1}}}=\\
& \leq2^{md/2}\left\Vert \phi_{\Omega;\varepsilon}\right\Vert _{W^{m\mathbf{1}%
,\infty}\left(  \Omega_{\varepsilon}\right)  }\left\Vert \mathcal{E}_{\Omega
}\right\Vert _{op}\left\Vert u\right\Vert _{W^{m\mathbf{1}}\left(
\Omega\right)  }\\
& =2^{md/2}\left\Vert \phi_{\Omega;\varepsilon}\right\Vert _{W^{m\mathbf{1}%
,\infty}\left(  \Omega_{\varepsilon}\right)  }\left\Vert u\right\Vert
_{W^{m\mathbf{1}}\left(  \Omega\right)  }.
\end{align*}
\medskip

\textbf{Part} \textbf{3c} From \ref{a1.62} and \ref{a1.63} of Lemma
\ref{Lem_func_eq_1_nbhd_set}, $\sigma_{\lambda}\phi_{\Omega;\varepsilon}%
=\phi_{\lambda\Omega;\lambda\varepsilon}$ and $\tau_{c}\phi_{\Omega
;\varepsilon}=\phi_{\Omega+c;\varepsilon}$. By part 5c of Remark
\ref{Rem_SobolevSpace2} $\tau_{c}\mathcal{E}_{\Omega}=\mathcal{E}_{\Omega
+c}\tau_{c}$ and $\sigma_{\lambda}\mathcal{E}_{\Omega}=\mathcal{E}%
_{\lambda\Omega}\sigma_{\lambda}$ so we have the two sequences of equations
\[
\tau_{c}\mathcal{E}_{\Omega;\varepsilon}u=\tau_{c}r_{\Omega_{\varepsilon}%
}\left(  \phi_{\Omega;\varepsilon}\mathcal{E}_{\Omega}u\right)  =r_{\Omega
_{\varepsilon}+c}\left(  \left(  \tau_{c}\phi_{\Omega;\varepsilon}\right)
\tau_{c}\mathcal{E}_{\Omega}u\right)  =r_{\left(  \Omega+c\right)
_{\varepsilon}}\left(  \phi_{\Omega+c;\varepsilon}\mathcal{E}_{\Omega+c}%
\tau_{c}u\right)  =\mathcal{E}_{\Omega+c;\varepsilon}\tau_{c}u,
\]

and%
\begin{equation}
\sigma_{\lambda}\mathcal{E}_{\Omega;\varepsilon}u=\sigma_{\lambda}%
r_{\Omega_{\varepsilon}}\left(  \phi_{\Omega;\varepsilon}\mathcal{E}_{\Omega
}u\right)  =r_{\lambda\Omega_{\varepsilon}}\left(  \phi_{\lambda\Omega
;\lambda\varepsilon}\sigma_{\lambda}\mathcal{E}_{\Omega}u\right)  =r_{\left(
\lambda\Omega\right)  _{\lambda\varepsilon}}\left(  \phi_{\lambda
\Omega;\lambda\varepsilon}\mathcal{E}_{\lambda\Omega}\sigma_{\lambda}u\right)
=\mathcal{E}_{\lambda\Omega;\lambda\varepsilon}\sigma_{\lambda}u.\label{a2.05}%
\end{equation}

An identical argument applies when $\mathcal{E}_{\Omega;\varepsilon}$ is
defined using $\mathcal{E}_{\Omega}^{m\mathbf{1}}$.\medskip

\textbf{Part 4} Noting Theorem \ref{Thm_ten_prod_two_wt_fns} we easily show
that the tensor product function $w_{m\mathbf{1}}$ is a weight function with
property W03 iff $\kappa<\left(  m-\frac{1}{2}\right)  \mathbf{1}$, and hence
by Theorem \ref{Thm_X_invers_Fourier_W3}, $X_{w_{m\mathbf{1}}}^{0}%
\hookrightarrow C_{B}^{\left(  \max\left\lfloor \kappa\right\rfloor \right)
}=C_{B}^{\left(  m-1\right)  \mathbf{1}}$. Now applying the binomial theorem
yields
\begin{align*}
w_{m\mathbf{1}}\left(  \xi\right)  =\left(  1+\left(  \xi_{k}^{2}\right)
\right)  ^{m\mathbf{1}}=\sum_{\alpha\leq m\mathbf{1}}\tbinom{m\mathbf{1}%
}{\alpha}\left(  \xi_{k}^{2}\right)  ^{\alpha}  & =\sum_{\alpha\leq
m\mathbf{1}}\tbinom{m}{\alpha_{1}}\ldots\tbinom{m}{\alpha_{d}}\xi^{2\alpha}\\
& \leq2^{md}\sum_{\alpha\leq m\mathbf{1}}\xi^{2\alpha},
\end{align*}

so that since%
\begin{align*}
\left\Vert u\right\Vert _{W^{m\mathbf{1}}}^{2}=\sum_{\alpha\leq m\mathbf{1}%
}\left\Vert D^{\alpha}u\right\Vert _{2}^{2}  & =\sum_{\alpha\leq m\mathbf{1}%
}\int\left\vert D^{\alpha}u\right\vert ^{2}=\sum_{\alpha\leq m\mathbf{1}}%
\int\left\vert \widehat{D^{\alpha}u}\right\vert ^{2}=\\
& =\sum_{\alpha\leq m\mathbf{1}}\int\xi^{2\alpha}\left\vert \widehat{u}%
\right\vert ^{2}=\int\left(  \sum_{\alpha\leq m\mathbf{1}}\xi^{2\alpha
}\right)  \left\vert \widehat{u}\right\vert ^{2},
\end{align*}

we have the equivalence%
\[
\left\Vert u\right\Vert _{W^{m\mathbf{1}}}^{2}\leq\int w_{m\mathbf{1}%
}\left\vert \widehat{u}\right\vert ^{2}\leq2^{md}\left\Vert u\right\Vert
_{W^{m\mathbf{1}}}^{2}.
\]

That $\left\Vert 1/w_{m\mathbf{1}}\right\Vert _{1}=\left(  2^{-2\left(
m-1\right)  }\tbinom{2m-2}{m-1}\pi\right)  ^{d}$ is just a consequence of a
standard 1-dimensional integral.\medskip

\textbf{Part 5} From part 4, $X_{w_{m\mathbf{1}}}^{0}\hookrightarrow
C_{B}^{\left(  m-1\right)  \mathbf{1}}$ and so from Definition
\ref{Cmb(clos_open_region)}, $X_{w_{m\mathbf{1}}}^{0}\left(  \Omega\right)
\subset r_{\Omega}C_{B}^{\left(  m-1\right)  \mathbf{1}}$. In addition,
$r_{\Omega}:X_{w_{m\mathbf{1}}}^{0}\overset{c}{\rightarrow}X_{w_{m\mathbf{1}}%
}^{0}\left(  \Omega\right)  $ and $r_{\Omega}:C_{B}^{\left(  m-1\right)
\mathbf{1}}\rightarrow r_{\Omega}C_{B}^{\left(  m-1\right)  \mathbf{1}}$ is
continuous under the supremum norm, and so $X_{w_{m\mathbf{1}}}^{0}\left(
\Omega\right)  \hookrightarrow r_{\Omega}C_{B}^{\left(  m-1\right)
\mathbf{1}}$.\medskip

\textbf{Part 6}. Suppose $\phi\in C_{0}^{\infty}\left(  \Omega\right)  $. Then
$\phi\left(  x\right)  =\int_{a_{1}}^{x_{1}}D_{1}\phi\left(  t_{1},x^{\prime
}\right)  dt_{1}$ and so by the Cauchy-Schwartz inequality, when $x_{1}%
\in\left[  a_{1},b_{1}\right]  $,%
\begin{align*}
\left\vert \phi\left(  x\right)  \right\vert \leq\int_{a_{1}}^{x_{1}%
}\left\vert D_{1}\phi\left(  t_{1},x^{\prime}\right)  \right\vert dt_{1} &
\leq\left(  \int_{a_{1}}^{x_{1}}dt_{1}\right)  ^{1/2}\left(  \int_{a_{1}%
}^{x_{1}}\left\vert D_{1}\phi\left(  t_{1},x^{\prime}\right)  \right\vert
^{2}dt_{1}\right)  ^{1/2}\\
&  \leq\left(  x_{1}-a_{1}\right)  ^{1/2}\left(  \int_{a_{1}}^{x_{1}%
}\left\vert D_{1}\phi\left(  t_{1},x^{\prime}\right)  \right\vert ^{2}%
dt_{1}\right)  ^{1/2}\\
&  \leq\left(  x_{1}-a_{1}\right)  ^{1/2}\left(  \int_{a_{1}}^{b_{1}%
}\left\vert D_{1}\phi\left(  t_{1},x^{\prime}\right)  \right\vert ^{2}%
dt_{1}\right)  ^{1/2}.
\end{align*}

Thus%
\begin{align*}
\int_{\Omega}\left\vert \phi\right\vert ^{2}=\int_{\mathcal{R}\left(
a,b\right)  }\left\vert \phi\right\vert ^{2} &  \leq\int_{a_{1}}^{b_{1}%
}\left(  x_{1}-a_{1}\right)  dx_{1}\int_{\mathcal{R}\left(  a^{\prime
},b^{\prime}\right)  }\int_{a_{1}}^{b_{1}}\left\vert D_{1}\phi\left(
t_{1},x^{\prime}\right)  \right\vert ^{2}dt_{1}dx^{\prime}\\
&  =\frac{1}{2}\left(  b_{1}-a_{1}\right)  ^{2}\int_{\mathcal{R}\left(
a,b\right)  }\left\vert D_{1}\phi\right\vert ^{2}\\
&  =\left(  \tfrac{b_{1}-a_{1}}{\sqrt{2}}\right)  ^{2}\int_{\Omega}\left\vert
D_{1}\phi\right\vert ^{2}dx.
\end{align*}

A permutation argument now shows that%
\[
\int_{\Omega}\left\vert \phi\right\vert ^{2}\leq\left(  \tfrac{b_{i}-a_{i}%
}{\sqrt{2}}\right)  ^{2}\int_{\Omega}\left\vert D_{i}\phi\right\vert
^{2},\quad i=1,\ldots,d.
\]

Thus, for all $\alpha$,%
\[
\int_{\Omega}\left\vert D^{\alpha}\phi\right\vert ^{2}\leq\left(  \tfrac
{b-a}{\sqrt{2}}\right)  ^{2\beta}\int_{\Omega}\left\vert D^{\alpha+\beta}%
\phi\right\vert ^{2},\quad\forall\alpha,\beta.
\]

and in particular, when $\beta=m\mathbf{1}-\alpha$,%
\[
\int_{\Omega}\left\vert D^{\alpha}\phi\right\vert ^{2}\leq\left(  \tfrac
{b-a}{\sqrt{2}}\right)  ^{2\left(  m\mathbf{1}-\alpha\right)  }\int_{\Omega
}\left\vert D^{m\mathbf{1}}\phi\right\vert ^{2},\quad\alpha\leq m\mathbf{1}.
\]

Finally, we will extend the result using density of $C_{0}^{\infty}\left(
\Omega\right)  $ in $W_{0}^{m\mathbf{1}}\left(  \Omega\right)  $. Choose $u\in
W_{0}^{m\mathbf{1}}\left(  \Omega\right)  $ and $\phi\in C_{0}^{\infty}\left(
\Omega\right)  $ and $\alpha\leq m\mathbf{1}$. Then%
\begin{align*}
&  \left(  \int_{\Omega}\left\vert D^{\alpha}u\right\vert ^{2}\right)
^{1/2}\\
&  \leq\left(  \int_{\Omega}\left\vert D^{\alpha}\phi\right\vert ^{2}\right)
^{1/2}+\left(  \int_{\Omega}\left\vert D^{\alpha}\left(  u-\phi\right)
\right\vert ^{2}\right)  ^{1/2}\\
&  \leq\left(  \int_{\Omega}\left\vert D^{\alpha}\phi\right\vert ^{2}\right)
^{1/2}+\left\Vert u-\phi\right\Vert _{W^{m\mathbf{1}}\left(  \Omega\right)
}\\
&  \leq\left(  \tfrac{b-a}{\sqrt{2}}\right)  ^{m\mathbf{1}-\alpha}\left(
\int_{\Omega}\left\vert D^{m\mathbf{1}}\phi\right\vert ^{2}\right)
^{1/2}+\left\Vert u-\phi\right\Vert _{W^{m\mathbf{1}}\left(  \Omega\right)
}\\
&  \leq\left(  \tfrac{b-a}{\sqrt{2}}\right)  ^{m\mathbf{1}-\alpha}\left(
\int_{\Omega}\left\vert D^{m\mathbf{1}}u\right\vert ^{2}\right)
^{1/2}+\left(  \tfrac{b-a}{\sqrt{2}}\right)  ^{m\mathbf{1}-\alpha}\left(
\int_{\Omega}\left\vert D^{m\mathbf{1}}\left(  u-\phi\right)  \right\vert
^{2}\right)  ^{1/2}+\left\Vert u-\phi\right\Vert _{W^{m\mathbf{1}}\left(
\Omega\right)  }\\
&  \leq\left(  \tfrac{b-a}{\sqrt{2}}\right)  ^{m\mathbf{1}-\alpha}\left(
\int_{\Omega}\left\vert D^{m\mathbf{1}}u\right\vert ^{2}\right)
^{1/2}+\left(  1+\left(  \tfrac{b-a}{\sqrt{2}}\right)  ^{m\mathbf{1}-\alpha
}\right)  \left\Vert u-\phi\right\Vert _{W^{m\mathbf{1}}\left(  \Omega\right)
},
\end{align*}

and the density of $C_{0}^{\infty}\left(  \Omega\right)  $ in $W_{0}%
^{m\mathbf{1}}\left(  \Omega\right)  $ implies%
\[
\left(  \int_{\Omega}\left\vert D^{\alpha}u\right\vert ^{2}\right)  ^{1/2}%
\leq\left(  \tfrac{b-a}{\sqrt{2}}\right)  ^{m\mathbf{1}-\alpha}\left(
\int_{\Omega}\left\vert D^{m\mathbf{1}}u\right\vert ^{2}\right)  ^{1/2},\quad
u\in W_{0}^{m\mathbf{1}}\left(  \Omega\right)  .
\]

\end{proof}

We will need some more properties of the central difference operator
$\delta_{2}^{2l\mathbf{1}}$ which was introduced in Theorem
\ref{Thm_data_func}. These properties will relate the space $X_{w_{s}}^{0}$ to
the space $W^{n\mathbf{1}}\left(  \Omega\right)  $ locally as illustrated in
Figure \ref{Fig_exten_commut_m1_centdiff} below.

\begin{lemma}
\label{Lem_centdiffop_property_2l}\textbf{Properties of the central difference
operator} $\delta_{2}^{2l\mathbf{1}}$\ Suppose $w_{s}$ is the extended
B-spline weight function with parameters $n$ and $l$. Then:

\begin{enumerate}
\item $\delta_{2}^{2l\mathbf{1}}:W^{n\mathbf{1}}\rightarrow X_{w_{s}}^{0}$ is continuous.

\item If $f\in W^{n\mathbf{1}}$ and $\operatorname*{supp}f\subset\left[
-1,1\right]  ^{d}$ then%
\begin{equation}
\left(  -1\right)  ^{ld}\binom{2l}{l}^{-d}\delta_{2}^{2l\mathbf{1}%
}f=f+\mathcal{A}_{2l}f,\label{a1.60}%
\end{equation}

where%
\begin{equation}
\mathcal{A}_{2l}f=\left(  -1\right)  ^{ld}\tbinom{2l}{l}^{-d}\sum
\limits_{\mathbf{0}\leq\beta\leq2l\mathbf{1};\text{ }\beta\neq l\mathbf{1}%
}\left(  -1\right)  ^{\left\vert \beta\right\vert }\tbinom{2l\mathbf{1}}%
{\beta}\tau_{2\beta-2l}f.\label{a1.16}%
\end{equation}

For $\mathbf{0}\leq\alpha,\beta\leq2l\mathbf{1}$,
\begin{align*}
\operatorname*{supp}\tau_{2\beta-2l}f  & \subset\left[  -1,1\right]
^{d}-\left(  2\beta-2l\right)  \subset\mathbb{R}^{d}\setminus\left(
-1,1\right)  ^{d},\\
\operatorname*{supp}\tau_{2\alpha-2l}f\cap\operatorname*{supp}\tau_{2\beta
-2l}f  & \subset\tau_{2\alpha-2l}\partial\left(  \left[  -1,1\right]
^{d}\right)  \cap\tau_{2\beta-2l}\partial\left(  \left[  -1,1\right]
^{d}\right)  ,\\
& \qquad\qquad when\text{ }\alpha\neq\beta,
\end{align*}

and%
\begin{align}
\operatorname*{supp}\mathcal{A}_{2l}f  & =\bigcup_{\mathbf{0}\leq\beta
\leq2l\mathbf{1};\text{ }\beta\neq l\mathbf{1}}\operatorname*{supp}%
\tau_{2\beta-2l}f,\nonumber\\
\operatorname*{supp}\mathcal{A}_{2l}f  & \subset\left(  2l+1\right)  \left[
-1,1\right]  ^{d}\setminus\left(  -1,1\right)  ^{d}.\label{a1.011}%
\end{align}
\smallskip

If $\Omega\subset\left(  -1,1\right)  ^{d}$ then:\smallskip

\item $\left(  -1\right)  ^{-ld}\tbinom{2l}{l}^{-d}\delta_{2}^{2l\mathbf{1}%
}\mathcal{E}_{0}:W_{0}^{n\mathbf{1}}\left(  \Omega\right)  \rightarrow
X_{w_{s}}^{0}$ is a continuous extension mapping where $\mathcal{E}_{0}$ is
the zero extension operator.

\item If $l$ is even then $\left(  -1\right)  ^{-\frac{l}{2}d}\tbinom{l}%
{l/2}^{-d}\delta_{2}^{l\mathbf{1}}\mathcal{E}_{0}:W_{0}^{n\mathbf{1}}\left(
\Omega\right)  \rightarrow X_{w_{s}}^{0}$ is a continuous extension mapping.
\end{enumerate}
\end{lemma}

\begin{proof}
\textbf{Part 1} If $f\in W^{n\mathbf{1}}$ then by the equality \ref{1.003} of
Theorem \ref{Thm_data_func},%
\begin{equation}
\left\Vert \delta_{2}^{2l\mathbf{1}}f\right\Vert _{w_{s},0}\leq4^{ld}%
\left\Vert D^{n\mathbf{1}}f\right\Vert _{2}\leq4^{ld}\left\Vert f\right\Vert
_{W^{n\mathbf{1}}}.\label{a2.4}%
\end{equation}
\medskip

\textbf{Part 2} From \ref{1.015},%
\begin{align}
\delta_{2}^{2l\mathbf{1}}f &  =\sum\limits_{\beta\leq2l\mathbf{1}}\left(
-1\right)  ^{\left\vert \beta\right\vert }\tbinom{2l\mathbf{1}}{\beta}%
\tau_{2\beta-2l}f\label{a1.00}\\
&  =\left(  -1\right)  ^{\left\vert l\mathbf{1}\right\vert }\tbinom
{2l\mathbf{1}}{l\mathbf{1}}f+\sum\limits_{\beta\leq2l\mathbf{1,}\beta\neq
l\mathbf{1}}\left(  -1\right)  ^{\left\vert \beta\right\vert }\tbinom
{2l\mathbf{1}}{\beta}\tau_{2\beta-2l}f\nonumber\\
&  =\left(  -1\right)  ^{ld}\tbinom{2l}{l}^{d}f+\sum\limits_{\beta
\leq2l\mathbf{1,}\beta\neq l\mathbf{1}}\left(  -1\right)  ^{\left\vert
\beta\right\vert }\tbinom{2l\mathbf{1}}{\beta}\tau_{2\beta-2l}f,\label{a2.00}%
\end{align}

which proves \ref{a1.16}. Further, from \ref{a1.16} and Definition
\ref{Def_dil_tran_dom_supp}, $\beta\neq l\mathbf{1}$ and $\operatorname*{supp}%
f\subset\left[  -1,1\right]  ^{d}$ implies
\begin{align*}
\operatorname*{supp}\tau_{2\beta-2l}f=\tau_{2\beta-2l}\operatorname*{supp}%
f=\operatorname*{supp}f-\left(  2\beta-2l\right)   &  \subset\left[
-1,1\right]  ^{d}-\left(  2\beta-2l\right) \\
&  \subset\mathbb{R}^{d}\setminus\left(  -1,1\right)  ^{d}.
\end{align*}

Since $f=0$ on $\tau_{2\beta-2l}\partial\left(  \left[  -1,1\right]
^{d}\right)  $
\[
\operatorname*{supp}\mathcal{A}_{2l}f\subset\bigcup\limits_{\beta
\leq2l\mathbf{1,}\beta\neq l\mathbf{1}}\operatorname*{supp}\tau_{2\beta
-2l}f\subset\mathbb{R}^{d}\setminus\left(  -1,1\right)  ^{d}.
\]

Also%
\begin{align*}
\operatorname*{supp}\mathcal{A}_{2l}f\subset\bigcup\limits_{\beta
\leq2l\mathbf{1,}\beta\neq l\mathbf{1}}\operatorname*{supp}\tau_{2\beta-2l}f
& \subset\bigcup\limits_{\beta\leq2l\mathbf{1}}\left(  2\beta
-2l+\operatorname*{supp}f\right) \\
& \subset\bigcup\limits_{\beta\leq2l\mathbf{1}}\left(  2\beta-2l+\left[
-1,1\right]  ^{d}\right) \\
& \subset\left[  -\left(  2l+1\right)  ,2l+1\right]  ^{d}\\
& =\left(  2l+1\right)  \left[  -1,1\right]  ^{d}.
\end{align*}
\medskip

\textbf{Part 3} Continuity follows directly from part 6 of Remark
\ref{Rem_SobolevSpace2} and then part 1 above. The extension property follows
from applying \ref{a1.011} to \ref{a1.60}.\medskip

\textbf{Part \ 4} Use \ref{1.003} and the expansion%
\[
\delta_{2}^{l\mathbf{1}}f=\sum\limits_{\beta\leq l\mathbf{1}}\left(
-1\right)  ^{\left\vert \beta\right\vert }\tbinom{l\mathbf{1}}{\beta}%
\tau_{2\beta-l}f=\left(  -1\right)  ^{\frac{l}{2}d}\tbinom{l}{l/2}^{d}%
f+\sum\limits_{\substack{\beta\leq l\mathbf{1} \\\beta\neq\frac{l}%
{2}\mathbf{1}}}\left(  -1\right)  ^{\left\vert \beta\right\vert }%
\tbinom{l\mathbf{1}}{\beta}\tau_{2\beta-l}f.
\]

\end{proof}

We will now use the extension operator $\delta_{2}^{l\mathbf{1}}%
\mathcal{E}_{0}:W_{0}^{n\mathbf{1}}\left(  \left(  -1,1\right)  ^{d}\right)
\rightarrow X_{w_{s}}^{0}$ of part 3 to construct a continuous extension
$E:W^{n\mathbf{1}}\left(  \Omega\right)  \rightarrow X_{w_{s}}^{0}$. The idea
is to map $W^{n\mathbf{1}}\left(  \Omega\right)  $ into $W^{n\mathbf{1}%
}\left(  \Omega^{0}\right)  $ where $\Omega^{0}\subset\left[  -1,1\right]
^{d}$ and then apply Lemma \ref{Lem_centdiffop_property_2l}.

\begin{theorem}
\label{Thm_exten_op_Wloc_to_X0ws}\textbf{The extension operator }%
$E:W^{n\mathbf{1}}\left(  \Omega\right)  \rightarrow X_{w_{s}}^{0}$

We construct our extension operator $E$ in six steps:

\begin{enumerate}
\item Suppose $w_{s}$ is an extended B-spline weight function on
$\mathbb{R}^{d}$ with parameters $n$ and $l$;

\item $\Omega$ is a bounded region satisfying the segment property and so
permits the continuous extension operator $\mathcal{E}_{\Omega;\varepsilon
}:W^{m\mathbf{1}}\left(  \Omega\right)  \rightarrow W_{0}^{m\mathbf{1}}\left(
\Omega_{\varepsilon}\right)  $ introduced in part 3b of Lemma
\ref{Lem_SobolevSpProperty};

\item $C$ is any open cube with $\overline{\Omega}\subset C$ and set
$0<\varepsilon\leq\operatorname*{dist}\left(  \Omega,\mathbb{R}^{d}\setminus
C\right)  $.

\item Then let $c$ be the centre of $C$ and choose an integer $m\geq1$ such
that%
\begin{equation}
\sigma_{1/m}\tau_{-c}C=\frac{1}{m}\left(  C-c\right)  \subseteq\left(
-1,1\right)  ^{d},\label{2.381}%
\end{equation}

and set $\Omega^{0}=\sigma_{1/m}\tau_{-c}\Omega$.

\item By Definitions \ref{Def_UnifRectCondit} and \ref{Def_SegCondit}
$\Omega^{0}$ also has the segment property and so by part 2 of Lemma
\ref{Lem_SobolevSpProperty} there always exists a continuous extension
operator $\mathcal{E}_{\Omega^{0};\varepsilon/m}:W^{n\mathbf{1}}\left(
\Omega^{0}\right)  \rightarrow W_{0}^{n\mathbf{1}}\left(  \Omega
_{\varepsilon/m}\right)  $.

\item We will now define the linear mapping $E$ of Figure
\ref{Fig_exten_commut_m1_centdiff} by:%
\begin{equation}
Ef:=\left(  -1\right)  ^{ld}\tbinom{2l}{l}^{-d}\tau_{c}\sigma_{m}\delta
_{2}^{2l\mathbf{1}}\mathcal{E}_{0}\mathcal{E}_{\Omega^{0};\varepsilon/m}%
\sigma_{1/m}\tau_{-c}f,\quad f\in W^{n\mathbf{1}}\left(  \Omega\right)
,\label{1.0402}%
\end{equation}

where $\mathcal{E}_{0}$ is the extension by zero operator of part 6 of Remark
\ref{Rem_SobolevSpace2}.
\end{enumerate}

\textbf{This theorem will show that}: $E:W^{n\mathbf{1}}\left(  \Omega\right)
\rightarrow X_{w_{s}}^{0}$ satisfies%
\begin{equation}
Ef=\mathcal{E}_{0}\mathcal{E}_{\Omega;\varepsilon}f+\left(  -1\right)
^{ld}\tbinom{2l}{l}^{-d}\sum\limits_{\beta\leq2l\mathbf{1},\beta\neq
l\mathbf{1}}\left(  -1\right)  ^{\left\vert \beta\right\vert }\tbinom
{2l\mathbf{1}}{\beta}\tau_{2m\left(  \beta-l\right)  }\mathcal{E}%
_{0}\mathcal{E}_{\Omega;\varepsilon}f.\label{a1.18}%
\end{equation}

Also we have the following inclusions for the support of $\mathcal{E}%
_{0}\mathcal{E}_{\Omega;\varepsilon}f$:%
\begin{align}
\operatorname*{supp}\mathcal{E}_{0}\mathcal{E}_{\Omega;\varepsilon}f  &
\subset\overline{\Omega_{\varepsilon}}\subset\overline{C}\subset c+m\left[
-1,1\right]  ^{d}.\label{a1.06}\\
\operatorname*{supp}\tau_{2m\left(  \beta-l\right)  }\mathcal{E}%
_{0}\mathcal{E}_{\Omega;\varepsilon}f  & \subset\operatorname*{supp}%
\mathcal{E}_{\Omega;\varepsilon}f-2m\left(  \beta-l\right)  \subset
\mathbb{R}^{d}\setminus\Omega_{\varepsilon},\nonumber\\
when\;\mathbf{0}  & \leq\beta\leq2l\mathbf{1},\text{ }\beta\neq l\mathbf{1}%
.\label{a1.07}\\
\operatorname*{supp}\tau_{2m\left(  \beta-l\right)  }\mathcal{E}%
_{0}\mathcal{E}_{\Omega;\varepsilon}f  & \subset\left(  c+\left(  2l+1\right)
m\left[  -1,1\right]  ^{d}\right)  \setminus\left(  c+m\left(  -1,1\right)
^{d}\right)  ,\nonumber\\
when\;\mathbf{0}  & \leq\beta\leq2l\mathbf{1},\text{ }\beta\neq\mathbf{1}%
,\nonumber\\
\operatorname*{supp}\tau_{2m\left(  \alpha-l\right)  }\mathcal{E}%
_{0}\mathcal{E}_{\Omega;\varepsilon}f\cap\operatorname*{supp}\tau_{2m\left(
\beta-l\right)  }\mathcal{E}_{0}\mathcal{E}_{\Omega;\varepsilon}f  &
\subset\left(  \tau_{2m\left(  \alpha-l\right)  }\partial C\right)
\cap\left(  \tau_{2m\left(  \beta-l\right)  }\partial C\right)  ,\nonumber\\
when\;\mathbf{0}  & \leq\alpha,\beta\leq2l\mathbf{1},\text{ }\alpha\neq
\beta,\label{a1.08}%
\end{align}

and for the support of $Ef$:%
\begin{align}
\operatorname*{supp}Ef  & =\bigcup\limits_{\mathbf{0}\leq\beta\leq
2l\mathbf{1}}\operatorname*{supp}\tau_{2m\left(  \beta-l\right)  }%
\mathcal{E}_{0}\mathcal{E}_{\Omega;\varepsilon}f.\label{a1.56}\\
\operatorname*{supp}Ef  & \subset c+\left(  2l+1\right)  m\left[  -1,1\right]
^{d}.\label{a1.01}%
\end{align}

Finally, $E:W^{n\mathbf{1}}\left(  \Omega\right)  \rightarrow X_{w_{s}}^{0}$
is a continuous linear \textbf{extension }with%
\begin{equation}
\left\Vert Ef\right\Vert _{w_{s},0}\leq\frac{\left(  4m\right)  ^{ld}}%
{\tbinom{2l}{l}^{d}}\left\Vert D^{n\mathbf{1}}\mathcal{E}_{\Omega;\varepsilon
}f\right\Vert _{L^{2}\left(  \Omega_{\varepsilon}\right)  }\leq\frac{\left(
4m\right)  ^{ld}}{\tbinom{2l}{l}^{d}}\left\Vert \mathcal{E}_{\Omega
;\varepsilon}f\right\Vert _{W_{0}^{n\mathbf{1}}\left(  \Omega_{\varepsilon
}\right)  },\label{a1.015}%
\end{equation}

and
\[
\left\Vert Ef\right\Vert _{w_{s},0}\leq\left(  e^{\frac{1}{6l}}\sqrt{\pi
l}m^{l}\right)  ^{d}\left\Vert \mathcal{E}_{\Omega;\varepsilon}f\right\Vert
_{W_{0}^{n\mathbf{1}}\left(  \Omega_{\varepsilon}\right)  }.
\]

\end{theorem}

\begin{proof}
The relevant commutative diagram is Figure \ref{Fig_exten_commut_m1_centdiff}
and the relevant sets and mappings are shown in Figure \ref{Fig_exten_sets_m}.%
\begin{figure}[ptb]%
\centering
\includegraphics[
natheight=1.974300in,
natwidth=5.321000in,
height=1.9743in,
width=5.321in
]%
{C:/Math_SwBasisFunc/InterpolSmthDev/PapersMonog/ZeroOrd/ZeroOrdDev/graphics/figLocData_commut_m_centdif_Wn1__4.pdf}%
\caption{Function spaces and mappings for defining the extension mapping E.}%
\label{Fig_exten_commut_m1_centdiff}%
\end{figure}
%

\begin{figure}[th]%
\centering
\includegraphics[
natheight=2.183600in,
natwidth=5.005300in,
height=2.1836in,
width=5.0053in
]%
{C:/Math_SwBasisFunc/InterpolSmthDev/PapersMonog/ZeroOrd/ZeroOrdDev/graphics/figLocData_SetsMaps_m__5.pdf}%
\caption{Sets and mappings for defining the extension mapping E.}%
\label{Fig_exten_sets_m}%
\end{figure}

We first show that the mapping \ref{1.0402} makes sense:\medskip

\textbf{Result 1} From part 1 of Lemma \ref{Lem_SobolevSpProperty},
$\sigma_{1./m}\tau_{-c}:W^{n\mathbf{1}}\left(  \Omega\right)  \rightarrow
W^{n\mathbf{1}}\left(  \Omega^{0}\right)  $ is a homeomorphism.\smallskip

\textbf{Result 2} From \ref{2.381}, $\Omega^{0}\subset\left[  -1,1\right]
^{d}$ and since $0<\varepsilon\leq\operatorname*{dist}\left(  \Omega
,\mathbb{R}^{d}\setminus C\right)  $, we have

$0<\sigma_{1/m}\varepsilon\leq\operatorname*{dist}\left(  \sigma_{1/m}%
\tau_{-c}\Omega,\mathbb{R}^{d}\setminus\sigma_{1/m}\tau_{-c}C\right)  ,$ i.e.
$0<\varepsilon/m\leq\operatorname*{dist}\left(  \Omega^{0},\mathbb{R}%
^{d}\setminus\left(  -1,1\right)  ^{d}\right)  $ and by part 3 of Lemma
\ref{Lem_SobolevSpProperty}, $\mathcal{E}_{\Omega^{0};\varepsilon
/m}:W^{n\mathbf{1}}\left(  \Omega^{0}\right)  \rightarrow W_{0}^{n\mathbf{1}%
}\left(  \Omega_{\varepsilon/m}^{0}\right)  $.\smallskip

\textbf{Result 3} From part 3 of Lemma \ref{Lem_centdiffop_property_2l},
$\left(  -1\right)  ^{ld}\binom{2l}{l}^{-d}\delta_{2}^{2l\mathbf{1}%
}\mathcal{E}_{0}:W_{0}^{n\mathbf{1}}\left(  \Omega_{\varepsilon/m}^{0}\right)
\rightarrow X_{w_{s}}^{0}$ is a continuous extension.\smallskip

\textbf{Result 4} From part 1 of Lemma \ref{Lem_centdiffop_property_2l},
$\sigma_{m}:X_{w}^{0}\rightarrow X_{w}^{0}$ is continuous.\medskip

These four results imply that $E:W^{n\mathbf{1}}\left(  \Omega\right)
\rightarrow X_{w_{s}}^{0}$ is continuous.

Equation \ref{a1.18} for $E$ follows directly from part 3 of Lemma
\ref{Lem_SobolevSpProperty}.

To prove the extension property of $E$ we start with \ref{a1.18} and
\ref{a1.00} and then compose the translations and dilations by noting that
$\tau_{c}\sigma_{m}\tau_{2\beta-2l}\sigma_{1/m}\tau_{-c}g=\tau_{m\left(
2\beta-2l\right)  }g$:%
\begin{align*}
Ef &  =\left(  -1\right)  ^{ld}\tbinom{2l}{l}^{-d}\tau_{c}\sigma_{m}%
\sum\limits_{\beta\leq2l\mathbf{1}}\left(  -1\right)  ^{\left\vert
\beta\right\vert }\tbinom{2l\mathbf{1}}{\beta}\tau_{2\beta-2l}\sigma_{1/m}%
\tau_{-c}\mathcal{E}_{0}\mathcal{E}_{\Omega;\varepsilon}f\\
&  =\left(  -1\right)  ^{ld}\tbinom{2l}{l}^{-d}\sum\limits_{\beta
\leq2l\mathbf{1}}\left(  -1\right)  ^{\left\vert \beta\right\vert }%
\tbinom{2l\mathbf{1}}{\beta}\tau_{c}\sigma_{m}\tau_{2\beta-2l}\sigma_{1/m}%
\tau_{-c}\mathcal{E}_{0}\mathcal{E}_{\Omega;\varepsilon}f\\
&  =\left(  -1\right)  ^{ld}\tbinom{2l}{l}^{-d}\left(  \sum\limits_{\beta
\leq2l\mathbf{1}}\left(  -1\right)  ^{\left\vert \beta\right\vert }%
\tbinom{2l\mathbf{1}}{\beta}\tau_{m\left(  2\beta-2l\right)  }\right)
\mathcal{E}_{0}\mathcal{E}_{\Omega;\varepsilon}f\\
&  =\mathcal{E}_{0}\mathcal{E}_{\Omega;\varepsilon}f+\left(  -1\right)
^{ld}\tbinom{2l}{l}^{-d}\sum\limits_{\beta\leq2l\mathbf{1},\beta\neq
l\mathbf{1}}\left(  -1\right)  ^{\left\vert \beta\right\vert }\tbinom
{2l\mathbf{1}}{\beta}\tau_{2m\left(  \beta-l\right)  }\mathcal{E}%
_{0}\mathcal{E}_{\Omega;\varepsilon}f,
\end{align*}

which proves \ref{a1.18}.

Obviously $\operatorname*{supp}\tau_{2m\left(  \alpha-l\right)  }%
\mathcal{E}_{0}\mathcal{E}_{\Omega;\varepsilon}f\subset\operatorname*{supp}%
\mathcal{E}_{0}\mathcal{E}_{\Omega;\varepsilon}f-2m\left(  \alpha-l\right)
\subset\Omega_{\varepsilon}-2m\left(  \alpha-l\right)  $. But since
$\Omega_{\varepsilon}\subset C$ and \ref{2.381} implies $C\subset c+m\left[
-1,1\right]  ^{d}$, we have \ref{a1.06}. Because of the scale factor $m$ these
inclusions immediately imply \ref{a1.07} and \ref{a1.08} and subsequently
\ref{a1.56}.

Since $\operatorname*{supp}\mathcal{E}_{0}\mathcal{E}_{\Omega;\varepsilon
}f\subset c+m\left[  -1,1\right]  ^{d}$,%
\begin{align*}
\operatorname*{supp}Ef=\bigcup\limits_{\mathbf{0}\leq\beta\leq2l\mathbf{1}%
}\operatorname*{supp}\tau_{2m\left(  \beta-l\right)  }\mathcal{E}%
_{0}\mathcal{E}_{\Omega;\varepsilon}f &  \subset\bigcup\limits_{\mathbf{0}%
\leq\beta\leq2l\mathbf{1}}\tau_{2m\left(  \beta-l\right)  }C\subset\\
&  \subset\bigcup\limits_{\mathbf{0}\leq\beta\leq2l\mathbf{1}}\left(
C-2m\left(  \beta-l\right)  \right) \\
&  \subset C-\bigcup\limits_{\mathbf{0}\leq\beta\leq2l\mathbf{1}}2m\left(
\beta-l\right)  \left[  -1,1\right]  ^{d}\\
&  =C-2m\bigcup\limits_{\mathbf{0}\leq\beta\leq2l\mathbf{1}}\left(
\beta-l\right) \\
&  =C+2m\left[  -l\mathbf{1},l\mathbf{1}\right] \\
&  =C+2ml\left[  -1,1\right]  ^{d}\\
&  \subset c+m\left[  -1,1\right]  ^{d}+2ml\left[  -1,1\right]  ^{d}\\
&  =c+\left(  2l+1\right)  m\left[  -1,1\right]  ^{d},
\end{align*}

and we have proved the inclusion \ref{a1.56}.\medskip

To prove inequality \ref{a1.015}\textbf{\ }we start with equation \ref{1.0402}
for $E$. When $f\in W^{n\mathbf{1}}\left(  \Omega\right)  $ using the
equations of part 3 of Lemma \ref{Lem_SobolevSpProperty} and then the
estimates \ref{1.057} and \ref{a2.4} yield%
\begin{align*}
\tbinom{2l}{l}^{d}\left\Vert Ef\right\Vert _{w_{s},0} &  =\left\Vert \tau
_{c}\sigma_{m}\delta_{2}^{2l\mathbf{1}}\mathcal{E}_{0}\mathcal{E}_{\Omega
^{0};\varepsilon/m}\sigma_{1/m}\tau_{-c}f\right\Vert _{w_{s},0}\\
&  =\left\Vert \tau_{c}\sigma_{m}\delta_{2}^{2l\mathbf{1}}\sigma_{1/m}%
\tau_{-c}\mathcal{E}_{0}\mathcal{E}_{\Omega;\varepsilon}f\right\Vert
_{w_{s},0}\\
&  =\left\Vert \sigma_{m}\delta_{2}^{2l\mathbf{1}}\sigma_{1/m}\tau
_{-c}\mathcal{E}_{0}\mathcal{E}_{\Omega;\varepsilon}f\right\Vert _{w_{s},0}\\
&  \leq m^{\left(  l-n+1/2\right)  d}\left\Vert \delta_{2}^{2l\mathbf{1}%
}\sigma_{1/m}\tau_{-c}\mathcal{E}_{0}\mathcal{E}_{\Omega;\varepsilon
}f\right\Vert _{w_{s},0}\\
&  =m^{\left(  l-n+1/2\right)  d}4^{ld}\left\Vert D^{n\mathbf{1}}\sigma
_{1/m}\tau_{-c}\mathcal{E}_{0}\mathcal{E}_{\Omega;\varepsilon}f\right\Vert
_{2}\\
&  =m^{\left(  l-n+1/2\right)  d}4^{ld}\left\Vert m^{nd}\sigma_{1/m}%
D^{n\mathbf{1}}\tau_{-c}\mathcal{E}_{0}\mathcal{E}_{\Omega;\varepsilon
}f\right\Vert _{2}\\
&  =m^{\left(  l+1/2\right)  d}4^{ld}\left\Vert \sigma_{1/m}D^{n\mathbf{1}%
}\tau_{-c}\mathcal{E}_{0}\mathcal{E}_{\Omega;\varepsilon}f\right\Vert _{2}\\
&  =m^{\left(  l+1/2\right)  d}4^{ld}m^{-d/2}\left\Vert D^{n\mathbf{1}}%
\tau_{-c}\mathcal{E}_{0}\mathcal{E}_{\Omega;\varepsilon}f\right\Vert _{2}\\
&  =\left(  4m\right)  ^{ld}\left\Vert D^{n\mathbf{1}}\tau_{-c}\mathcal{E}%
_{0}\mathcal{E}_{\Omega;\varepsilon}f\right\Vert _{2}\\
&  =\left(  4m\right)  ^{ld}\left\Vert D^{n\mathbf{1}}\mathcal{E}%
_{0}\mathcal{E}_{\Omega;\varepsilon}f\right\Vert _{2}\\
&  =\left(  4m\right)  ^{ld}\left\Vert D^{n\mathbf{1}}E_{\Omega;\varepsilon
}f\right\Vert _{L^{2}\left(  \Omega_{\varepsilon}\right)  }\\
&  \leq\left(  4m\right)  ^{ld}\left\Vert \mathcal{E}_{\Omega;\varepsilon
}f\right\Vert _{W_{0}^{n\mathbf{1}}\left(  \Omega_{\varepsilon}\right)  }.
\end{align*}

But from \ref{Ap004}, $\frac{^{e^{-\frac{1}{6l}}}}{\sqrt{\pi l}}<2^{-2l}%
\binom{2l}{l}$ so that $\binom{2l}{l}^{-1}<\frac{^{e^{\frac{1}{6l}}\sqrt{\pi
l}}}{2^{2l}}$ and hence%
\begin{align*}
\left\Vert Ef\right\Vert _{w_{s},0}\leq\left(  \binom{2l}{l}^{-1}\left(
4m\right)  ^{l}\right)  ^{d}\left\Vert \mathcal{E}_{\Omega;\varepsilon
}f\right\Vert _{W_{0}^{n\mathbf{1}}\left(  \Omega_{\varepsilon}\right)  } &
<\left(  \frac{^{e^{\frac{1}{6l}}\sqrt{\pi l}}}{2^{2l}}\left(  4m\right)
^{l}\right)  ^{d}\left\Vert \mathcal{E}_{\Omega;\varepsilon}f\right\Vert
_{W_{0}^{n\mathbf{1}}\left(  \Omega_{\varepsilon}\right)  }\\
&  =\left(  e^{\frac{1}{6l}}\sqrt{\pi l}m^{l}\right)  ^{d}\left\Vert
\mathcal{E}_{\Omega;\varepsilon}f\right\Vert _{W_{0}^{n\mathbf{1}}\left(
\Omega_{\varepsilon}\right)  }.
\end{align*}

Finally, an application of \ref{a1.07} and \ref{a1.08} to \ref{a1.18} yields
$\operatorname*{supp}\left(  Ef-\mathcal{E}_{0}\mathcal{E}_{\Omega
;\varepsilon}f\right)  \subset\mathbb{R}^{d}\setminus\Omega_{\varepsilon}$,
which means that $Ef=\mathcal{E}_{\Omega;\varepsilon}f$ on $\Omega
_{\varepsilon}$ and hence that $Ef=\mathcal{E}_{\Omega;\varepsilon}f=f$ on
$\Omega$. Thus $E$ is an extension.
\end{proof}

\begin{remark}
If instead we assume that $\Omega$ satisfies the rectangle condition then we
can either note that $\Omega$ must have the segment property and construct
$\mathcal{E}_{\Omega;\varepsilon}$ using the canonical extension operator
$\mathcal{E}_{\Omega}$ or more directly we can construct $\mathcal{E}%
_{\Omega;\varepsilon}$ using the extension operator $\mathcal{E}_{\Omega
}^{m\mathbf{1}}$.
\end{remark}

\begin{remark}
?? From Theorem \ref{Thm_ord0_Riesz_rep_W2} the Riesz representer is%
\[
R_{x}(y)=\left(  2\pi\right)  ^{-d/2}G(y-x).
\]

From the theorem in Chapter 5 of Part 1 of Aronszajn \cite{Aronszajn50} the
unique local Riesz representer $R_{x}^{\Omega}\left(  y\right)  $ is given by%
\[
R_{x}^{\Omega}\left(  y\right)  =R_{x}\left(  y\right)  ,\quad x,y\in\Omega.
\]

How does the above extension operator $E$ relate to the extension defined by
$R_{x}^{\Omega}\rightarrow R_{x}$ when $x\in\Omega$?

See \ref{1.087} i.e. $r_{\Omega}^{\ast}R_{y}^{\Omega}=R_{y}$.

How does the above extension operator $E$ relate to translations in general?
\end{remark}

If we replace $W^{n\mathbf{1}}\left(  \Omega\right)  $ by $W_{0}^{n\mathbf{1}%
}\left(  \Omega\right)  $ in the last theorem then we do not need to leave
room for the extension $\mathcal{E}_{\Omega;\varepsilon}$ which goes to zero
in a neighborhood of $\Omega$. In this case we can set $\varepsilon=0$ and
adjust the proof to obtain:

\begin{corollary}
\label{Cor_exten_op_Wnoloc_to_X0ws}\textbf{The extension operator }%
$E:W_{0}^{n\mathbf{1}}\left(  \Omega\right)  \rightarrow X_{w_{s}}^{0}$. Suppose:

\begin{enumerate}
\item $w_{s}$ is an extended B-spline weight function on $\mathbb{R}^{d}$ with
parameters $n$ and $l$;

\item $\Omega$ is a bounded region with the segment property;

\item $C$ is any open cube with $\Omega\subseteq C$;

\item Let $c$ be the centre of $C$ and choose an integer $m\geq1$ such that%
\[
\Omega^{0}:=\sigma_{1/m}\tau_{-c}C=\frac{1}{m}\left(  C-c\right)
\subset\left(  -1,1\right)  ^{d}.
\]

Now define the mapping:%
\[
Ef=\left(  -1\right)  ^{ld}\tbinom{2l}{l}^{-d}\tau_{c}\sigma_{m}\delta
_{2}^{2l\mathbf{1}}\mathcal{E}_{0}\sigma_{1/m}\tau_{-c}f,\quad f\in
W_{0}^{n\mathbf{1}}\left(  \Omega\right)  ,
\]

where $\mathcal{E}_{0}$ is the zero extension operator.

Then $E:W_{0}^{n\mathbf{1}}\left(  \Omega\right)  \rightarrow X_{w_{s}}^{0}$
is a continuous linear extension such that%
\[
Ef=\mathcal{E}_{0}f+\left(  -1\right)  ^{ld}\tbinom{2l}{l}^{-d}\sum
\limits_{\beta\leq2l\mathbf{1},\beta\neq l\mathbf{1}}\left(  -1\right)
^{\left\vert \beta\right\vert }\tbinom{2l\mathbf{1}}{\beta}\tau_{2m\left(
\beta-l\right)  }\mathcal{E}_{0}f.
\]

Regarding supports:%
\begin{align*}
\operatorname*{supp}f  & \subset\overline{\Omega}\subset\overline{C}\subset
c+m\left[  -1,1\right]  ^{d},\\
\operatorname*{supp}\tau_{2m\left(  \alpha-l\right)  }f  & \subset
\operatorname*{supp}f-2m\left(  \alpha-l\right)  \subset\mathbb{R}%
^{d}\setminus\Omega,\\
& \qquad when\;\mathbf{0}\leq\alpha\leq2l\mathbf{1},\text{ }\alpha\neq
l\mathbf{1},\\
\operatorname*{supp}\tau_{2m\left(  \alpha-l\right)  }f\cap
\operatorname*{supp}\tau_{2m\left(  \beta-l\right)  }f  & \subset\left(
\tau_{2m\left(  \alpha-l\right)  }\partial C\right)  \cap\left(
\tau_{2m\left(  \beta-l\right)  }\partial C\right)  ,\\
& \qquad when\;\mathbf{0}\leq\alpha,\beta\leq2l\mathbf{1},\text{ }\alpha
\neq\beta,
\end{align*}

and%
\begin{align*}
\operatorname*{supp}Ef  & =\bigcup\limits_{\mathbf{0}\leq\beta\leq
2l\mathbf{1}}\operatorname*{supp}\tau_{2m\left(  \beta-l\right)  }%
\mathcal{E}_{0}f,\\
\operatorname*{supp}Ef  & \subset c+m\left(  2l+1\right)  \left[  -1,1\right]
^{d}.
\end{align*}

Regarding the continuity of $E$:%
\[
\left\Vert Ef\right\Vert _{w_{s},0}\leq\left(  4m\right)  ^{ld}\tbinom{2l}%
{l}^{-d}\left\Vert D^{n\mathbf{1}}f\right\Vert _{L^{2}\left(  \Omega\right)
}\leq\left(  e^{\frac{1}{6l}}\sqrt{\pi l}m^{l}\right)  ^{d}\left\Vert
D^{n\mathbf{1}}f\right\Vert _{L^{2}\left(  \Omega\right)  },
\]

and the right side is an equivalent norm on $W_{0}^{n\mathbf{1}}\left(
\Omega\right)  $.

\item Regarding the Fourier transform of $Ef$,%
\[
\widehat{Ef}=4^{ld}\tbinom{2l}{l}^{-d}\left(  \sin m\xi_{k}\right)
^{2l\mathbf{1}}\widehat{\mathcal{E}_{0}f}.
\]

\item If $f\in W_{0}^{n\mathbf{1}}\left(  \Omega\right)  \simeq\left(
X_{w}^{0}\right)  _{\overline{\Omega}}$ and $g\in X_{w_{s}}^{0}$ then%
\[
\left(  Ef,g\right)  _{w_{s},0}=4^{ld}\tbinom{2l}{l}^{-d}\int_{\Omega
}D^{n\mathbf{1}}f\text{ }\overline{D^{n\mathbf{1}}\mathcal{B}_{m}%
^{2l\mathbf{1}}g},
\]

where for any multi-index $\alpha\geq\mathbf{0}$,%
\begin{align*}
\mathcal{B}_{m}^{\alpha}  & :=\mathcal{B}_{m}^{\alpha_{1}}\mathcal{B}%
_{m}^{\alpha_{2}}\ldots\mathcal{B}_{m}^{\alpha_{d}},\text{\quad}%
\mathcal{B}_{m}^{\alpha_{j}}:=\left(  \mathcal{B}_{m}^{\mathbf{e}_{j}}\right)
^{\alpha_{j}},\text{\quad}\mathcal{B}_{m}^{\mathbf{e}_{j}}:=\sum
\limits_{k=0}^{m}\tau_{\left(  2k-m\right)  \mathbf{e}_{j}}.\\
\widehat{\mathcal{B}_{m}^{\mathbf{e}_{j}}g}\left(  \xi\right)   & =\frac{\sin
m\xi_{j}}{\sin\xi_{j}}\widehat{g},\quad\widehat{\mathcal{B}_{m}^{\alpha}%
g}\left(  \xi\right)  =\left(  \frac{\sin m\xi_{j}}{\sin\xi_{j}}\right)
^{\alpha}\widehat{g}.
\end{align*}

and $\left\{  \mathbf{e}_{1},\mathbf{e}_{2},\ldots,\mathbf{e}_{d}\right\}  $
is the canonical basis for $\mathbb{R}^{d}$.

\item ?? (See part 6 of Lemma \ref{Lem_SobolevSpProperty}) If $f\in
W_{0}^{n\mathbf{1}}\left(  \Omega\right)  \simeq\left(  X_{w_{s}}^{0}\right)
_{\overline{\Omega}}$ then%
\[
f\left(  x\right)  =4^{ld}\tbinom{2l}{l}^{-d}\int_{\Omega}D^{n\mathbf{1}%
}f\text{ }\overline{D^{n\mathbf{1}}\mathcal{B}_{m}^{2l\mathbf{1}}R_{x}^{w_{s}%
}},\quad x\in\Omega,
\]

and%
\[
\left\vert f\left(  x\right)  \right\vert \leq4^{ld}\tbinom{2l}{l}%
^{-d}\left\Vert D^{n\mathbf{1}}\mathcal{B}_{m}^{2l\mathbf{1}}R_{x}^{w_{s}%
}\right\Vert _{L^{2}\left(  \Omega\right)  }\left\Vert D^{n\mathbf{1}%
}f\right\Vert _{L^{2}\left(  \Omega\right)  },\quad x\in\Omega,
\]

and $\left\Vert D^{n\mathbf{1}}f\right\Vert _{L^{2}\left(  \Omega\right)  }$
is a norm on $W_{0}^{n\mathbf{1}}\left(  \Omega\right)  $.
\end{enumerate}
\end{corollary}

\begin{proof}
\textbf{Parts 1 to 4} follow directly from Theorem
\ref{Thm_exten_op_Wloc_to_X0ws} with $\varepsilon=0$ but with part 4 needing
part 6 of Lemma \ref{Lem_SobolevSpProperty}.\medskip

\textbf{Part 5}%
\begin{align*}
\sum\limits_{\beta\leq2l\mathbf{1}}\left(  -1\right)  ^{\left\vert
\beta\right\vert } &  \tbinom{2l\mathbf{1}}{\beta}\left(  \tau_{2m\left(
\beta-l\mathbf{1}\right)  }\mathcal{E}_{0}f\right)  ^{\wedge}\\
&  =\left(  \sum\limits_{\beta\leq2l\mathbf{1}}\tbinom{2l\mathbf{1}}{\beta
}\left(  -1\right)  ^{\left\vert \beta\right\vert }e^{-i2m\left(
\beta-l\mathbf{1}\right)  \xi}\right)  \widehat{\mathcal{E}_{0}f}\\
&  =\left(  \sum\limits_{\beta\leq2l\mathbf{1}}\tbinom{2l\mathbf{1}}{\beta
}\left(  -1\right)  ^{\left\vert \beta\right\vert }e^{i2m\left(
l\mathbf{1}-\beta\right)  \xi}\right)  \widehat{\mathcal{E}_{0}f}\\
&  =e^{-i2ml\mathbf{1}\xi}\left(  \sum\limits_{\beta\leq2l\mathbf{1}}%
\tbinom{2l\mathbf{1}}{\beta}\left(  -1\right)  ^{\left\vert \beta\right\vert
}e^{i2m\left(  2l\mathbf{1}-\beta\right)  \xi}\right)  \widehat{\mathcal{E}%
_{0}f}\\
&  =e^{-i2ml\mathbf{1}\xi}\left(  \sum\limits_{\beta\leq2l\mathbf{1}}%
\tbinom{2l\mathbf{1}}{\beta}\left(  -1\right)  ^{\left\vert \beta\right\vert
}\left(  e^{i2m\xi_{k}}\right)  ^{2l\mathbf{1}-\beta}\right)
\widehat{\mathcal{E}_{0}f}\\
&  =e^{-i2ml\mathbf{1}\xi}\left(  -\mathbf{1}+\left(  e^{i2m\mathbf{\xi}_{k}%
}\right)  \right)  ^{2l\mathbf{1}}\widehat{\mathcal{E}_{0}f}\\
&  =\left(  \left(  -e^{-im\mathbf{\xi}_{k}}\right)  +\left(  e^{im\mathbf{\xi
}_{k}}\right)  \right)  ^{2l\mathbf{1}}\widehat{\mathcal{E}_{0}f}\\
&  =\left(  e^{im\mathbf{\xi}_{k}}-e^{-im\mathbf{\xi}_{k}}\right)
^{2l\mathbf{1}}\widehat{\mathcal{E}_{0}f}\\
&  =\left(  2i\right)  ^{2ld}\left(  \frac{e^{im\mathbf{\xi}_{k}%
}-e^{-im\mathbf{\xi}_{k}}}{2i}\right)  ^{2l\mathbf{1}}\widehat{\mathcal{E}%
_{0}f}\\
&  =\left(  2i\right)  ^{2ld}\left(  \sin m\xi_{k}\right)  ^{2l\mathbf{1}%
}\widehat{\mathcal{E}_{0}f}\\
&  =\left(  2i\right)  ^{2ld}\left(  \sin m\xi_{1}\sin m\xi_{2}\ldots\sin
m\xi_{d}\right)  ^{2l}\widehat{\mathcal{E}_{0}f},
\end{align*}

so that%
\begin{align*}
\widehat{Ef}  & =\left(  -1\right)  ^{ld}\tbinom{2l}{l}^{-d}\left(  2i\right)
^{2ld}\left(  \sin m\xi_{k}\right)  ^{2l\mathbf{1}}\widehat{\mathcal{E}_{0}%
f}\\
& =\left(  i\right)  ^{2ld}\tbinom{2l}{l}^{-d}\left(  2i\right)  ^{2ld}\left(
\sin m\xi_{k}\right)  ^{2l\mathbf{1}}\widehat{\mathcal{E}_{0}f}\\
& =4^{ld}\tbinom{2l}{l}^{-d}\left(  \sin m\xi_{k}\right)  ^{2l\mathbf{1}%
}\widehat{\mathcal{E}_{0}f}.
\end{align*}
\medskip

\textbf{Part 6} ??%
\begin{align*}
\left(  Ef,g\right)  _{w,0}  & =\int w\widehat{Ef}\text{ }\overline
{\widehat{g}}=4^{ld}\tbinom{2l}{l}^{-d}\int w\left(  \sin m\xi_{k}\right)
^{2l\mathbf{1}}\widehat{\mathcal{E}_{0}f}\text{ }\overline{\widehat{g}}=\\
& =4^{ld}\tbinom{2l}{l}^{-d}\int\frac{\xi^{2n\mathbf{1}}}{\left(  \sin\xi
_{k}\right)  ^{2l\mathbf{1}}}\left(  \sin m\xi_{k}\right)  ^{2l\mathbf{1}%
}\widehat{\mathcal{E}_{0}f}\text{ }\overline{\widehat{g}}=\\
& =4^{ld}\tbinom{2l}{l}^{-d}\int\xi^{2n\mathbf{1}}\widehat{\mathcal{E}_{0}%
f}\overline{\left(  \frac{\sin m\xi_{k}}{\sin\xi_{k}}\right)  ^{2l\mathbf{1}%
}\widehat{g}}\\
& =4^{ld}\tbinom{2l}{l}^{-d}\int\xi^{2n\mathbf{1}}\widehat{\mathcal{E}_{0}%
f}\text{ }\overline{\widehat{\mathcal{B}_{m}^{2l\mathbf{1}}g}}\\
& =4^{ld}\tbinom{2l}{l}^{-d}\int\widehat{\mathcal{E}_{0}D^{n\mathbf{1}}%
f}\text{ }\overline{\widehat{D^{n\mathbf{1}}\mathcal{B}_{m}^{2l\mathbf{1}}g}%
}\\
& =4^{ld}\tbinom{2l}{l}^{-d}\int_{\Omega}D^{n\mathbf{1}}f\text{ }%
\overline{D^{n\mathbf{1}}\mathcal{B}_{m}^{2l\mathbf{1}}g}.
\end{align*}

In one dimension%
\begin{align*}
\frac{\sin ms}{\sin s}  & =\frac{\left(  e^{is}\right)  ^{m}-\left(
e^{-is}\right)  ^{m}}{e^{is}-e^{-is}}=\sum\limits_{k=0}^{m}\left(
e^{is}\right)  ^{m-k}\left(  e^{-is}\right)  ^{k}.\\
\frac{\sin ms}{\sin s}\widehat{g}  & =\sum\limits_{k=0}^{m}\left(
e^{is}\right)  ^{m-k}\left(  e^{-is}\right)  ^{k}\widehat{g}.
\end{align*}

But $\left(  \tau_{k}g\right)  ^{\wedge}=e^{-ik\xi}\widehat{g}$ so
\[
\frac{\sin ms}{\sin s}\widehat{g}=\sum\limits_{k=0}^{m}\left(  e^{is}\right)
^{m-k}\widehat{\tau_{k}g}=\left(  \sum\limits_{k=0}^{m}\tau_{2k-m}g\right)
^{\wedge}=\widehat{\mathcal{B}_{m}^{2l}g}.
\]

\textbf{Part 7} Let $g=R_{x}$ in part 6.
\end{proof}

\begin{remark}
\label{Rem_Thm_data_fn_2}\ 

\begin{enumerate}
\item Noting part 4 of Lemma \ref{Lem_centdiffop_property_2l}, when $l$ is
even the operator $\delta_{2}^{l\mathbf{1}}$ can be used instead of
$\delta_{2}^{2l\mathbf{1}}$ to construct an extension $E:W^{n\mathbf{1}%
}\left(  \Omega\right)  \rightarrow X_{w_{s}}^{0}$ which can be used below.

\item The mappings in Figure \ref{Fig_exten_commut_m1_centdiff} which are used
to define the extension $E$ are independent of $n\leq l$.

\item ?? Calculate the Hilbert space adjoint of $E$?
\end{enumerate}
\end{remark}

Noting Definition \ref{Def_dil_tran_dom_supp}:

\begin{theorem}
\label{Thm_transl_locXow_isom_iso}We have $\tau_{a}:X_{w}^{0}\left(
\Omega\right)  \rightarrow X_{w}^{0}\left(  \tau_{a}\Omega\right)  $ is an
isometric isomorphism.
\end{theorem}

\begin{proof}
We use the translation operator facts that for all $b\in\mathbb{R}^{d}$,
$\tau_{b}X_{w}^{0}=X_{w}^{0}$ and%
\begin{equation}
\tau_{b}r_{\Omega}f=r_{\tau_{b}\Omega}\tau_{b}f.\label{1.094}%
\end{equation}

Hence%
\begin{align*}
\left\Vert \tau_{a}f\right\Vert _{w,0;\tau_{a}\Omega}  & =\inf\left\{
\left\Vert u\right\Vert _{w,0}:u\in X_{w}^{0},\text{ }r_{\tau_{a}\Omega}%
u=\tau_{a}f\right\} \\
& =\inf\left\{  \left\Vert u\right\Vert _{w,0}:u\in X_{w}^{0},\text{ }%
\tau_{-a}r_{\tau_{a}\Omega}u=f\right\} \\
& =\inf\left\{  \left\Vert u\right\Vert _{w,0}:u\in X_{w}^{0},\text{
}r_{\Omega}\tau_{-a}u=f\right\} \\
& =\inf\left\{  \left\Vert u\right\Vert _{w,0}:u\in X_{w}^{0},\text{
}r_{\Omega}\tau_{-a}u=f\right\} \\
& =\inf\left\{  \left\Vert v\right\Vert _{w,0}:\tau_{a}v\in X_{w}^{0},\text{
}r_{\Omega}v=f\right\} \\
& =\inf\left\{  \left\Vert v\right\Vert _{w,0}:v\in\tau_{-a}X_{w}^{0},\text{
}r_{\Omega}v=f\right\} \\
& =\inf\left\{  \left\Vert v\right\Vert _{w,0}:v\in X_{w}^{0},\text{
}r_{\Omega}v=f\right\} \\
& =\left\Vert f\right\Vert _{w,0;\Omega}%
\end{align*}

\end{proof}

\begin{theorem}
\label{Thm_dilat_locXow_homeo}Suppose that $w$ is a weight function with
property W02 or W03. Then if $\lambda.>\mathbf{0}$:

\begin{enumerate}
\item The dilation $\sigma_{\lambda}:X_{w}^{0}\rightarrow X_{\sigma
_{\mathbf{1}/\lambda}w}^{0}$ is a homeomorphism with inverse $\sigma
_{\mathbf{1/}\lambda}$ and operator norm $\left\vert \lambda^{\mathbf{1}%
}\right\vert $.

\item Also $\sigma_{\lambda}:X_{w}^{0}\left(  \Omega\right)  \rightarrow
X_{\sigma_{\mathbf{1}/\lambda}w}^{0}\left(  \lambda.\Omega\right)  $ is a
homeomorphism with inverse $\sigma_{\mathbf{1/}\lambda}$ and operator norm
$\left\vert \lambda^{\mathbf{1}}\right\vert $.
\end{enumerate}
\end{theorem}

\begin{proof}
\textbf{Part 1} Note part 8 of Remark \ref{Rem_Def_extend_wt_fn}. From part 12
of Summary \ref{Sum_FourTransf} $\left(  \sigma_{\lambda}f\right)  ^{\wedge
}=\left\vert \lambda^{\mathbf{1}}\right\vert \sigma_{\mathbf{1}./\lambda
}\widehat{f}$, where $\sigma_{\lambda}f\left(  x\right)  :=f\left(
x./\lambda\right)  $. Thus%
\begin{align*}
\left\Vert \sigma_{\lambda}u\right\Vert _{w,0}^{2}  & =\int w\left(
\xi\right)  \left\vert \widehat{\sigma_{\lambda}u}\left(  \xi\right)
\right\vert ^{2}d\xi\\
& =\int w\left(  \xi\right)  \left\vert \left\vert \lambda^{\mathbf{1}%
}\right\vert \widehat{u}\left(  \lambda.\xi\right)  \right\vert ^{2}d\xi\\
& =\int w\left(  \xi\right)  \left(  \lambda^{\mathbf{1}}\right)
^{2}\left\vert \widehat{u}\left(  \lambda.\xi\right)  \right\vert ^{2}d\xi\\
& :\eta=\lambda.\xi,\text{ }d\xi=\frac{1}{\left\vert \lambda^{\mathbf{1}%
}\right\vert }d\eta\Rightarrow\\
& =\int w\left(  \eta./\lambda\right)  \left(  \lambda^{\mathbf{1}}\right)
^{2}\left\vert \widehat{u}\left(  \eta\right)  \right\vert ^{2}\frac
{1}{\left\vert \lambda^{\mathbf{1}}\right\vert }d\eta\\
& =\left\vert \lambda^{\mathbf{1}}\right\vert \int w\left(  \eta
./\lambda\right)  \left\vert \widehat{u}\left(  \eta\right)  \right\vert
^{2}d\eta
\end{align*}

\textbf{Part 2} ??
\end{proof}

\begin{corollary}
\label{Cor_Thm_data_fn_2}Let $w_{s}$ be an extended B-spline weight function
on $\mathbb{R}^{d}$ with parameters $n$ and $l$. Suppose $\Omega$ is a bounded
region with the rectangle or segment property.

Then $W^{n\mathbf{1}}\left(  \Omega\right)  \hookrightarrow X_{\sigma
_{\lambda}w_{s}}^{0}\left(  \Omega\right)  $ where $\sigma_{\lambda}$ is any
scalar dilation operator.
\end{corollary}

\begin{proof}
By part 5b of Remark \ref{Rem_SobolevSpace2} the rectangle condition implies
the segment property. Under this condition Theorem\textbf{\ }%
\ref{Thm_exten_op_Wloc_to_X0ws} showed that a continuous extension
$E:W^{n\mathbf{1}}\left(  \Omega\right)  \hookrightarrow X_{w}^{0}$ exists for
the extended B-splines. Now by Definition \ref{Def_Xow(open_set)}, $r_{\Omega
}E:W^{n\mathbf{1}}\left(  \Omega\right)  \rightarrow X_{w}^{0}\left(
\Omega\right)  $ is a continuous inclusion. Since $\lambda\Omega$ also has the
segment property we also have $W^{n\mathbf{1}}\left(  \lambda\Omega\right)
\hookrightarrow X_{w_{s}}^{0}\left(  \lambda\Omega\right)  $, and an
application of the mapping result of part 2 of Theorem
\ref{Thm_dilat_locXow_homeo} and part 1 of Lemma \ref{Lem_SobolevSpProperty}
enables us to write%
\[
W^{n\mathbf{1}}\left(  \Omega\right)  \overset{\sigma_{\lambda}%
}{\longrightarrow}W^{n\mathbf{1}}\left(  \lambda\Omega\right)  \hookrightarrow
X_{w_{s}}^{0}\left(  \lambda\Omega\right)  \overset{\sigma_{1/\lambda
}}{\longrightarrow}X_{\sigma_{\lambda}w_{s}}^{0}\left(  \Omega\right)  ,
\]

which proves this corollary.
\end{proof}

Observe that the local space $X_{w}^{0}\left(  \Omega\right)  $ can be
regarded as the set of data functions for the minimal norm interpolation
problem where the data is contained in $\Omega$, and Corollary
\ref{Cor_Thm_data_fn_2} proves that $W^{n\mathbf{1}}\left(  \Omega\right)
\subset X_{w}^{0}\left(  \Omega\right)  $ for scaled extended B-splines. The
minimum norm interpolation problem can be defined in terms of extensions of
these spaces to $\mathbb{R}^{d}$ and the matrix equation for the interpolant
only uses values of the data function on $\Omega$. \textbf{Thus }%
$W^{n\mathbf{1}}\left(  \Omega\right)  $\textbf{\ is the source of all data
functions for the interpolation problem} \textbf{when }$w$\textbf{\ is a
scaled extended B-spline weight function}.

The last corollary will now be strengthened to characterize the data functions
for the minimal norm interpolation problem for scaled extended B-splines.

\begin{lemma}
\label{Lem_Thm_Xow(O)_eq_Hn(O)_dim1}Suppose $w_{s}$ is an extended B-spline
weight function with parameters $n$ and $l$, and $\Omega$ is a bounded region
with the rectangle or segment property.

Then $X_{\sigma_{\lambda}w_{s}}^{0}\left(  \Omega\right)  \hookrightarrow
W^{n\mathbf{1}}\left(  \Omega\right)  $ where $\sigma_{\lambda}$ is any scalar
dilation operator.
\end{lemma}

\begin{proof}
Suppose $v\in X_{w_{s}}^{0}\left(  \Omega\right)  $ and let $v^{e}$ be any
extension of $v$ to $X_{w_{s}}^{0}$. Noting that $n\leq l$, if $\alpha\leq n $
we have%
\begin{align*}
\left\Vert v^{e}\right\Vert _{w_{s},0}^{2}=\int w_{s}\left\vert \widehat{v^{e}%
}\right\vert ^{2}=\int\frac{s^{2n\mathbf{1}}}{\left(  \sin s_{i}\right)
^{2l\mathbf{1}}}\left\vert \widehat{v^{e}}\left(  s\right)  \right\vert ^{2}ds
&  =\int s^{2\alpha}\frac{s^{2\left(  n\mathbf{1}-\alpha\right)  }}{\left(
\sin s_{i}\right)  ^{2l\mathbf{1}}}\left\vert \widehat{v^{e}}\left(  s\right)
\right\vert ^{2}ds\\
&  \geq\int s^{2\alpha}\frac{s^{2\left(  n\mathbf{1}-\alpha\right)  }}{\left(
\sin s_{i}\right)  ^{2\left(  n\mathbf{1}-\alpha\right)  }}\left\vert
\widehat{v^{e}}\left(  s\right)  \right\vert ^{2}ds\\
&  \geq\int s^{2\alpha}\left\vert \widehat{v^{e}}\left(  s\right)  \right\vert
^{2}ds\\
&  =\int\left\vert \widehat{D^{\alpha}v^{e}}\left(  s\right)  \right\vert
^{2}ds\\
&  =\int\left\vert D^{\alpha}v^{e}\left(  t\right)  \right\vert ^{2}dt\\
&  \geq\left\Vert D^{\alpha}v\right\Vert _{L^{2}\left(  \Omega\right)  }^{2},
\end{align*}

and so $D^{\alpha}v\in L^{2}\left(  \Omega\right)  $ for $\alpha\leq
n\mathbf{1}$, which means $v\in W^{n\mathbf{1}}\left(  \Omega\right)  $ since
$v\in L^{2}\left(  \Omega\right)  $. Further, the definition \ref{1.026} of
the norm on $X_{w_{s}}^{0}\left(  \Omega\right)  $ implies $\left\Vert
D^{\alpha}v\right\Vert _{L^{2}\left(  \Omega\right)  }\leq\left\Vert
v\right\Vert _{w_{s},0;\Omega}$ when $\alpha\leq n\mathbf{1}$, verifying that
$X_{w_{s}}^{0}\left(  \Omega\right)  \hookrightarrow W^{n\mathbf{1}}\left(
\Omega\right)  $.

By part 5b of Remark \ref{Rem_SobolevSpace2} the rectangle condition implies
the segment property.

Since $\lambda\Omega$ also has the segment property the operator
$\mathcal{E}_{\lambda\Omega}$ is also defined . Hence $X_{w_{s}}^{0}\left(
\lambda\Omega\right)  \hookrightarrow W^{n\mathbf{1}}\left(  \lambda
\Omega\right)  $ and an application of part 2 of Theorem
\ref{Thm_dilat_locXow_homeo} and part 1 of Lemma \ref{Lem_SobolevSpProperty}
enables us to write%
\[
X_{\sigma_{\lambda}w_{s}}^{0}\left(  \Omega\right)  \overset{\sigma_{\lambda
}}{\longrightarrow}X_{w_{s}}^{0}\left(  \lambda\Omega\right)  \hookrightarrow
W^{n\mathbf{1}}\left(  \lambda\Omega\right)  \overset{\sigma_{1/\lambda
}}{\longrightarrow}W^{n\mathbf{1}}\left(  \Omega\right)  ,
\]

which proves this corollary.
\end{proof}

We now have our \textbf{main result}:

\begin{theorem}
\label{Thm_int_Xow(O)_eq_Hn(O)_dim1}Suppose $w_{s}$ is an extended B-spline
weight function with parameters $n,l$ and $\Omega\subset\mathbb{R}^{d}$ is a
bounded region with the rectangle or segment property. Then for any dilation
operator $\sigma_{\lambda}$:

\begin{enumerate}
\item $W^{n\mathbf{1}}\left(  \Omega\right)  =X_{\sigma_{\lambda}w_{s}}%
^{0}\left(  \Omega\right)  $ as sets and their norms are equivalent.

\item $X_{\sigma_{\lambda}w_{s}}^{0}\left(  \Omega\right)  \hookrightarrow
r_{\Omega}C_{B}^{\left(  n-1\right)  \mathbf{1}}$.
\end{enumerate}
\end{theorem}

\begin{proof}
First note that the rectangle property implies the segment property.\medskip

\textbf{Part 1} This follows directly from Lemma
\ref{Lem_Thm_Xow(O)_eq_Hn(O)_dim1} and Corollary \ref{Cor_Thm_data_fn_2}%
.\medskip

\textbf{Part 2} From part 5 of Lemma \ref{Lem_SobolevSpProperty}, $X_{w_{s}%
}^{0}\left(  \Omega\right)  \hookrightarrow r_{\Omega}C_{B}^{\left(
n-1\right)  \mathbf{1}}$.

From part 1, $X_{w_{s}}^{0}\left(  \Omega\right)  =W^{n\mathbf{1}}\left(
\Omega\right)  =X_{\sigma_{\lambda}w_{s}}^{0}\left(  \Omega\right)  $ as sets
so $X_{\sigma_{\lambda}w_{s}}^{0}\left(  \Omega\right)  \subset r_{\Omega
}C_{B}^{\left(  n-1\right)  \mathbf{1}}$.

Also from part 1, $X_{w_{s}}^{0}\left(  \Omega\right)  $ and $X_{\sigma
_{\lambda}w_{s}}^{0}\left(  \Omega\right)  $ must have equivalent norms so
$X_{\sigma_{\lambda}w_{s}}^{0}\left(  \Omega\right)  \hookrightarrow
r_{\Omega}C_{B}^{\left(  n-1\right)  \mathbf{1}}$.
\end{proof}

Recall that for arbitrary $K$ closed the Hilbert space $\left(  X_{w_{s}}%
^{0}\right)  _{K}=\left\{  f\in X_{w}^{0}:\operatorname*{supp}f\subseteq
K\right\}  $ was introduced in Definition \ref{Def_XowK}.

\begin{corollary}
\label{Cor_Thm_int_Xow(O)_eq_Hn(O)_dim1}The mapping $E\iota:X_{w_{s}}%
^{0}\left(  \Omega\right)  \rightarrow\left(  X_{w_{s}}^{0}\right)
_{c+m\left(  2l+1\right)  \left[  -1,1\right]  ^{d}}$ is a continuous extension.

Here $\iota:X_{w_{s}}^{0}\left(  \Omega\right)  \rightarrow W^{n\mathbf{1}%
}\left(  \Omega\right)  $ is the continuous inclusion operator of Lemma
\ref{Lem_Thm_Xow(O)_eq_Hn(O)_dim1} and $E:W^{n\mathbf{1}}\left(
\Omega\right)  \rightarrow X_{w_{s}}^{0}$ is the continuous extension operator
of Theorem \ref{Thm_exten_op_Wloc_to_X0ws}.
\end{corollary}

\begin{proof}
From parts 4 and 5 of Theorem \ref{Thm_exten_op_Wloc_to_X0ws},
$E:W^{n\mathbf{1}}\left(  \Omega\right)  \rightarrow X_{w_{s}}^{0}$ with
$\operatorname*{supp}Ef\subset c+m\left(  2l+1\right)  \left[  -1,1\right]
^{d}$ so that $E:W^{n\mathbf{1}}\left(  \Omega\right)  \rightarrow\left(
X_{w_{s}}^{0}\right)  _{c+m\left(  2l+1\right)  \left[  -1,1\right]  ^{d}}$ is
continuous. Clearly $E\iota f=Ef=f$ on $\Omega$.
\end{proof}

\begin{remark}
\label{Rem_Cor_Thm_int_Xow(O)_eq_Hn(O)_dim1}?? Can the Exact smoother problem
be formulated locally on $\left(  X_{w_{s}}^{0}\right)  _{c+m\left(
2l+1\right)  \left[  -1,1\right]  ^{d}}$? Yes, because the support of the
Exact smoother lies in this space? ?? Minimize the smoother on
$\operatorname*{range}E\iota$?

Since $G_{s}$ has bounded support it follows that for sufficiently large $m$,
$W_{G_{s}}^{\Omega}\subset\left(  X_{w_{s}}^{0}\right)  _{c+m\left(
2l+1\right)  \left[  -1,1\right]  ^{d}}$. Hence, for such an $m$,
\[
\mathcal{S}_{X}f=\operatorname*{argmin}\limits_{f\in\left(  X_{w_{s}}%
^{0}\right)  _{c+m\left(  2l+1\right)  \left[  -1,1\right]  ^{d}}}J_{s}\left[
f\right]  \in W_{G_{s}}^{\Omega}\subset\left(  X_{w_{s}}^{0}\right)
_{c+m\left(  2l+1\right)  \left[  -1,1\right]  ^{d}}.
\]

Define
\[
S_{X}^{\Omega}F:=r_{\Omega}\mathcal{S}_{X}r_{\Omega}^{\ast}F.
\]

Note that \ref{1.084} implies $r_{\Omega}^{\ast}r_{\Omega}g=g$ when $g\in
W_{G_{s}}^{\Omega}$. Hence%
\[
r_{\Omega}^{\ast}\mathcal{S}_{X}^{\Omega}F=\mathcal{S}_{X}r_{\Omega}^{\ast}F.
\]

?? Also, from Corollary \ref{Cor2_Thm_canon_exten_op}, $r_{w_{s};\Omega}%
^{\ast}:X_{w_{s}}^{0}\left(  \Omega\right)  \rightarrow\left(  X_{w_{s}}%
^{0}\right)  _{\overline{\Omega}+\operatorname*{supp}G_{s}}$ is a continuous
extension operator with $\left\Vert r_{w_{s};\Omega}^{\ast}\right\Vert
_{op}=1$. Hence%
\[
\mathcal{S}_{X}f=\operatorname*{argmin}\limits_{f\in\left(  X_{w_{s}}%
^{0}\right)  _{\overline{\Omega}+\operatorname*{supp}G_{s}}}J_{s}\left[
f\right]  \in W_{G_{s}}^{\Omega}\subset\left(  X_{w_{s}}^{0}\right)
_{\overline{\Omega}+\operatorname*{supp}G_{s}}.
\]

Consider $\left(  \mathcal{S}_{X}f-f,g\right)  _{w,0}$ where $f,g\in\left(
X_{w_{s}}^{0}\right)  _{\overline{\Omega}+\operatorname*{supp}G_{s}}$?

Also%
\begin{align*}
\left(  \mathcal{S}_{X}^{\Omega}F-F,F^{\prime}\right)  _{w,0;\Omega}=\left(
r_{\Omega}\mathcal{S}_{X}r_{\Omega}^{\ast}F-F,F^{\prime}\right)
_{w,0;\Omega}  & =\left(  r_{\Omega}^{\ast}r_{\Omega}\mathcal{S}_{X}r_{\Omega
}^{\ast}F-r_{\Omega}^{\ast}F,r_{\Omega}^{\ast}F^{\prime}\right)  _{w,0}\\
& =\left(  \mathcal{S}_{X}r_{\Omega}^{\ast}F-r_{\Omega}^{\ast}F,r_{\Omega
}^{\ast}F^{\prime}\right)  _{w,0}.
\end{align*}

\end{remark}

??

\begin{remark}
??? Characterize $W_{0}^{n\mathbf{1}}\left(  \Omega\right)  ^{\prime}$ (see
characterization of $W_{0}^{n}\left(  \Omega\right)  ^{\prime}$).

p62 of Adams and Fournier \cite{AdamFour2003} - $W^{-m,p^{\prime}}\left(
\Omega\right)  $.

Compare this to $\left(  X_{w_{s}}^{0}\right)  ^{\prime}$ and $X_{1/w_{s}}^{0}
$.
\end{remark}

\begin{theorem}
?? FIX! Suppose that for some integers $n,l>0$,
\[
C_{1}\left(  \sin\xi_{k}\right)  ^{4l\mathbf{1}}\leq\frac{\left(  1+\xi
.\xi\right)  ^{n\mathbf{1}}}{w}\leq C_{2}.
\]

Suppose also that $\Omega$ has the segment property or uniform rectangle property.

Then $X_{w}^{0}\left(  \Omega\right)  \simeq H^{n\mathbf{1}}\left(
\Omega\right)  $.

?? What about the converse?
\end{theorem}

\begin{proof}
Refer to Figure \ref{Fig_commut_embed_Xow_Wn1}.

First note that the uniform property implies the segment property.

?? Theorem \ref{Thm_Xow_embed_Hm1_iff} and the second inequality imply
$X_{w}^{0}\overset{\iota}{\hookrightarrow}H^{n\mathbf{1}}$ with the embedding
$\iota$ having operator norm $\left\Vert \iota\right\Vert =\left\Vert
\frac{\left(  1+\xi.\xi\right)  ^{n\mathbf{1}}}{w\left(  \xi\right)
}\right\Vert _{\infty}^{1/2}$. ??%

\begin{figure}[tbh]%
\centering
\includegraphics[
natheight=1.734200in,
natwidth=5.616000in,
height=1.7342in,
width=5.616in
]%
{C:/Math_SwBasisFunc/InterpolSmthDev/PapersMonog/ZeroOrd/ZeroOrdDev/graphics/figLocData_commut_embed_Xow_in_Wn1__6.pdf}%
\caption{?? ADD BLAH!}%
\label{Fig_commut_embed_Xow_Wn1}%
\end{figure}

By Theorem \ref{Thm_canon_exten_op} there exists a continuous extension
operator $E_{w}:=r_{\Omega}^{\ast}$, $E_{w}:X_{w}^{0}\left(  \Omega\right)
\rightarrow X_{w}^{0}$ so the continuous operator $\iota_{\Omega}=r_{\Omega
}\iota E_{w}$ is an embedding.

From Definition \ref{Def_SobolevSpace2} we know that $W^{n\mathbf{1}}$ is a
data function space and hence by Theorem \ref{Thm_canon_exten_op} there exists
a continuous extension operator $E_{n}:=r_{\Omega}^{\ast}$, $E_{n}%
:W^{n\mathbf{1}}\left(  \Omega\right)  \rightarrow W^{n\mathbf{1}}$ with
$\left\Vert E_{n}\right\Vert _{op}=1$.
\end{proof}

\subsection{A larger class of weight
functions\label{SbSect_loc_data_larger_class}}

In Theorem \ref{Thm_int_Xow(O)_eq_Hn(O)_dim1} it was shown that
$W^{n\mathbf{1}}\left(  \Omega\right)  =X_{w}^{0}\left(  \Omega\right)  $
setwise where $w$ is the scaled extended B-spline weight function with
parameters $n$ and $l$, and $\Omega\subset\mathbb{R}^{d}$ is a bounded region
with the segment (Definition \ref{Def_SegCondit}) or uniform rectangle
property (Definition \ref{Def_UnifRectCondit}). Here we will extend the class
of weight functions to which this $W^{n\mathbf{1}}\left(  \Omega\right)  $
characterization applies, although it still essentially only applies to tensor
product weight functions.

The class of weight functions will be extended for the express purpose of
including the central difference weight functions which will be introduced in
Chapter \ref{Ch_cent_diff_wt_fn_ten_prod}. This will be done by generalizing
the central difference operator $\delta_{2}^{2l\mathbf{1}}$, or $\delta
_{2}^{l\mathbf{1}}$ \textbf{if }$l$\textbf{\ is even}, to a convolution operator.

The central difference operators $\delta_{2}^{l\mathbf{1}}$, $\delta
_{2}^{2l\mathbf{1}}$ have the properties $\widehat{\delta_{2}^{l\mathbf{1}}%
u}\left(  \xi\right)  =\left(  2i\sin\xi_{k}\right)  ^{l\mathbf{1}}%
\widehat{u}\left(  \xi\right)  $ and $\widehat{\delta_{2}^{2l\mathbf{1}}%
u}\left(  \xi\right)  =\left(  2i\sin\xi_{k}\right)  ^{2l\mathbf{1}%
}\widehat{u}\left(  \xi\right)  $, so that when multiplied by the weight
function the zeros of the sines can be used to cancel out the weight function
poles, as in the proof of part 2 of Theorem \ref{Thm_data_func}, leaving a
function of polynomial increase of at most $n$. This ensures the continuity of
$\delta_{2}^{l\mathbf{1}},\delta_{2}^{2l\mathbf{1}}:W^{n\mathbf{1}}\rightarrow
X_{w}^{0}$; so $\left(  2i\sin\xi_{k}\right)  ^{l\mathbf{1}}$ and $\left(
2i\sin\xi_{k}\right)  ^{2l\mathbf{1}}$\ are examples of $g$ below and this
motivates condition \ref{a1.052}. It was also shown in Lemma
\ref{Lem_centdiffop_property_2l} that $f\in W^{n\mathbf{1}}$ and
$\operatorname*{supp}f\subset\left[  -1,1\right]  ^{d}$ implies
\[
\left(  -1\right)  ^{ld}\tbinom{2l\mathbf{1}}{l\mathbf{1}}^{-d}\delta
_{2}^{2l\mathbf{1}}f=f+\mathcal{A}_{2l}f,\text{ }where\text{ }%
\operatorname*{supp}\mathcal{A}_{2l}f\subset\mathbb{R}^{d}\setminus\left[
-1,1\right]  ^{d},
\]

and \textbf{when }$l$\textbf{\ is even}:
\[
\left(  -1\right)  ^{ld}\tbinom{l\mathbf{1}}{\frac{l}{2}\mathbf{1}}^{-d}%
\delta_{2}^{l\mathbf{1}}f=f+\mathcal{A}_{l}f,\text{ }where\text{
}\operatorname*{supp}\mathcal{A}_{l}f\subset\mathbb{R}^{d}\setminus\left[
-1,1\right]  ^{d},
\]

and these examples motivate assumption \ref{1.053} below.

Suppose%
\begin{equation}
\eta\in S^{\prime},\text{\quad}\operatorname*{supp}\eta\subset\mathbb{R}%
^{d}\setminus2\left(  -1,1\right)  ^{d},\text{\quad}\widehat{\eta}\in
L_{loc}^{1}.\label{a1.04}%
\end{equation}

Then we define%
\begin{equation}
\chi=\eta+\left(  2\pi\right)  ^{d/2}\delta,\label{a1.05}%
\end{equation}

which implies $\chi\in S^{\prime}$ and $\widehat{\chi}\in L_{loc}^{1}$.
Suppose also that there is a constant $C_{\chi}\geq0$ and a positive integer
$n$ such that%
\begin{equation}
w\left(  \xi\right)  \left\vert \widehat{\chi}\left(  \xi\right)  \right\vert
^{2}\leq C_{\chi}w_{n}\left(  \xi\right)  \text{ }a.e.,\label{a1.052}%
\end{equation}

where $w_{n}$ is the weight function defined by \ref{1.052}.

Then $\left(  \chi\ast\right)  \mathcal{E}_{0}$ maps $W_{0}^{n\mathbf{1}}$ to
$X_{w}^{0}$ continuously since $u\in W_{0}^{n\mathbf{1}}$ implies $\left(
\chi\ast\right)  \mathcal{E}_{0}u\in S^{\prime}$ and thus $\widehat{\chi
\ast\mathcal{E}_{0}u}=\widehat{\chi}\widehat{\mathcal{E}_{0}u}\in L_{loc}^{1}$
so that
\begin{align}
\left\Vert \left(  \chi\ast\right)  \mathcal{E}_{0}u\right\Vert _{w,0}^{2}  &
=\int w\left\vert \widehat{\left(  \chi\ast\right)  \mathcal{E}_{0}%
u}\right\vert ^{2}=\int w\left\vert \widehat{\chi}\widehat{\mathcal{E}_{0}%
u}\right\vert ^{2}\leq C_{\chi}\int\prod\limits_{k=1}^{d}\left(  1+\xi_{k}%
^{2}\right)  ^{n}\left\vert \widehat{\mathcal{E}_{0}u}\right\vert
^{2}=\nonumber\\
& =C_{\chi}\left\Vert \mathcal{E}_{0}u\right\Vert _{W^{n\mathbf{1}}}%
^{2}=C_{\chi}\left\Vert u\right\Vert _{W_{0}^{n\mathbf{1}}}^{2}.\label{a4.8}%
\end{align}

Also, if $u\in W_{0}^{n\mathbf{1}}\left(  \left(  -1,1\right)  ^{d}\right)  $
then%
\begin{equation}
\operatorname*{supp}\left(  \left(  \chi\ast\right)  \mathcal{E}%
_{0}u-\mathcal{E}_{0}u\right)  =\operatorname*{supp}\left(  \eta
\ast\mathcal{E}_{0}u\right)  \subset\operatorname*{supp}\eta
+\operatorname*{supp}\mathcal{E}_{0}u\subset\mathbb{R}^{d}\setminus\left(
-1,1\right)  ^{d},\label{1.053}%
\end{equation}

which means that $\left(  \chi\ast\right)  \mathcal{E}_{0}:W_{0}^{n\mathbf{1}%
}\left(  \left(  -1,1\right)  ^{d}\right)  \rightarrow X_{w}^{0}$ is a
continuous extension and consequently%
\begin{equation}
\left(  \chi\ast\right)  \mathcal{E}_{0}:W_{0}^{n\mathbf{1}}\left(
\Omega\right)  \rightarrow X_{w}^{0}\text{ }is\text{ }a\text{ }%
continuous\text{ }extension,\label{a1.002}%
\end{equation}

when $\Omega\subset\left[  -1,1\right]  ^{d}$.

\begin{example}
\fbox{The operator $\delta_{2}^{2l\mathbf{1}}$} Set $\left(  -1\right)
^{ld}\binom{2l\mathbf{1}}{l\mathbf{1}}^{-d}\delta_{2}^{2l\mathbf{1}%
}\mathcal{E}_{0}u=\left(  \chi\ast\right)  \mathcal{E}_{0}u$.\ Then
\begin{align*}
\left(  -1\right)  ^{ld}\tbinom{2l\mathbf{1}}{l\mathbf{1}}^{-d}\delta
_{2}^{2l\mathbf{1}}\mathcal{E}_{0}u  & =\left(  2\pi\right)  ^{d/2}\left(
-1\right)  ^{ld}\tbinom{2l\mathbf{1}}{l\mathbf{1}}^{-d}\delta\ast\left(
\delta_{2}^{2l\mathbf{1}}\mathcal{E}_{0}u\right) \\
& =\left(  2\pi\right)  ^{d/2}\left(  \left(  -1\right)  ^{ld}\tbinom
{2l\mathbf{1}}{l\mathbf{1}}^{-d}\delta_{2}^{2l\mathbf{1}}\delta\right)
\ast\mathcal{E}_{0}u,
\end{align*}

so that%
\begin{align*}
\chi & =\left(  2\pi\right)  ^{d/2}\left(  -1\right)  ^{ld}\tbinom
{2l\mathbf{1}}{l\mathbf{1}}^{-d}\delta_{2}^{2l\mathbf{1}}\delta\in
\mathcal{E}^{\prime},\\
\eta & =\chi-\left(  2\pi\right)  ^{d/2}\delta\\
& =\left(  2\pi\right)  ^{d/2}\left(  \left(  -1\right)  ^{ld}\tbinom
{2l\mathbf{1}}{l\mathbf{1}}^{-d}\delta_{2}^{2l\mathbf{1}}-1\right)  \delta
\in\mathcal{E}^{\prime},
\end{align*}

and from \ref{a2.00}%
\[
\eta=\left(  -1\right)  ^{-ld}\tbinom{2l}{l}^{-d}\sum\limits_{\beta
\leq2l\mathbf{1,}\beta\neq l\mathbf{1}}\left(  -1\right)  ^{\left\vert
\beta\right\vert }\tbinom{2l\mathbf{1}}{\beta}\tau_{2\beta-2l}f,
\]

so that $\operatorname*{supp}\eta\subset\mathbb{R}^{d}\setminus2\left(
-1,1\right)  ^{d}$. From \ref{1.004},
\begin{align*}
\widehat{\chi}=\left(  2\pi\right)  ^{d/2}\left(  -1\right)  ^{ld}%
\tbinom{2l\mathbf{1}}{l\mathbf{1}}^{-d}\widehat{\delta_{2}^{2l\mathbf{1}%
}\delta}  & =\left(  2\pi\right)  ^{d/2}\left(  -1\right)  ^{ld}%
\tbinom{2l\mathbf{1}}{l\mathbf{1}}^{-d}\left(  2i\sin\xi_{k}\right)
^{2l\mathbf{1}}\widehat{\delta}\\
& =\left(  -1\right)  ^{ld}\tbinom{2l\mathbf{1}}{l\mathbf{1}}^{-d}\left(
2i\sin\xi_{k}\right)  ^{2l\mathbf{1}}.
\end{align*}

\end{example}

The commutative diagram for the construction of the extension map
$E:W^{n\mathbf{1}}\left(  \Omega\right)  \rightarrow X_{w}^{0}$ is now Figure
\ref{Fig_exten_commut_mu_convol}:%
\begin{figure}[tbh]%
\centering
\includegraphics[
natheight=1.988900in,
natwidth=5.742900in,
height=1.9889in,
width=5.7429in
]%
{C:/Math_SwBasisFunc/InterpolSmthDev/PapersMonog/ZeroOrd/ZeroOrdDev/graphics/figLocData_commut_mu_convol_Wn1__7.pdf}%
\caption{Function spaces and mappings for defining the extension mapping E.}%
\label{Fig_exten_commut_mu_convol}%
\end{figure}

This diagram implies that
\begin{equation}
Ef=\tau_{c}\sigma_{m}\left(  \chi\ast\right)  \mathcal{E}_{0}\mathcal{E}%
_{\Omega^{0};\varepsilon/m}\sigma_{1/m}\tau_{-c}f,\quad f\in W^{n\mathbf{1}%
}\left(  \Omega\right)  .\label{1.071}%
\end{equation}

From the diagram we see that $\sigma_{m}:X_{w}^{0}\rightarrow X_{w}^{0}$ must
be continuous where $m\geq1$ is an integer. This was proven for the spline
weight functions in part 4 of Theorem \ref{Thm_data_func} and inspection of
this proof shows that continuity still holds if there exists an integer
$m\geq1$ and a constant $c_{m}>0$ such that%
\begin{equation}
w\left(  \xi/m\right)  \leq c_{m}w\left(  \xi\right)  ,\quad a.e.\text{
}on\text{ }\mathbb{R}^{d}.\label{a1.055}%
\end{equation}

\begin{lemma}
\label{Lem_Thm_ex_data_fn_3}Suppose the weight function $w$ has property W03
and satisfies \ref{a1.055}. Then:

\begin{enumerate}
\item (cf. \ref{1.057}) The dilation operator $\sigma_{m}:X_{w}^{0}\rightarrow
X_{w}^{0}$ is continuous. In fact,%
\begin{equation}
\left\Vert \sigma_{m}f\right\Vert _{w,0}\leq\sqrt{c_{m}m^{d}}\left\Vert
f\right\Vert _{w,0}.\label{a1.02}%
\end{equation}
\smallskip

If $\chi\in S^{\prime}$ and $u\in\mathcal{E}^{\prime}$ then:\smallskip

\item $\sigma_{m}\left(  \chi\ast u\right)  =m^{-d}\left(  \sigma_{m}%
\chi\right)  \ast\left(  \sigma_{m}u\right)  $ and $\chi\ast\left(  \sigma
_{m}u\right)  =m^{d}\sigma_{m}\left(  \left(  \sigma_{1/m}\chi\right)  \ast
u\right)  $ and,

\item $\tau_{c}\left(  \chi\ast u\right)  =\chi\ast\left(  \tau_{c}u\right)
=\left(  \tau_{c}\chi\right)  \ast u$.

\item The operator $\sigma_{m}\left(  \chi\ast\right)  \mathcal{E}_{0}$ is a
continuous linear mapping from $W_{0}^{n\mathbf{1}}\left(  \left(
-1,1\right)  ^{d}\right)  $ to $X_{w}^{0}$ and%
\[
\left\Vert \sigma_{m}\left(  \chi\ast\right)  \mathcal{E}_{0}u\right\Vert
_{w,0}\leq\sqrt{C_{\chi}c_{m}m^{d}}\left\Vert u\right\Vert _{W_{0}%
^{n\mathbf{1}}\left(  \left(  -1,1\right)  ^{d}\right)  },\quad u\in
W_{0}^{n\mathbf{1}}\left(  \left(  -1,1\right)  ^{d}\right)  .
\]

\item If $\Sigma$ is a region then%
\[
\left\Vert \sigma_{1/m}g\right\Vert _{W^{n\mathbf{1}}\left(  \Sigma\right)
}\leq m^{-\frac{d}{2}}\max\left\{  1,m^{nd}\right\}  \left\Vert g\right\Vert
_{W^{n\mathbf{1}}\left(  m\Sigma\right)  },\quad g\in W^{n\mathbf{1}}\left(
m\Sigma\right)  .
\]

\end{enumerate}
\end{lemma}

\begin{proof}
\textbf{Part 1}%
\begin{align*}
\left\Vert \sigma_{m}f_{d}\right\Vert _{w,0}^{2}=\int w\left(  t\right)
\left\vert \widehat{\sigma_{m}f_{d}}\right\vert ^{2}\left(  t\right)
dt=m^{2d}\int w\left(  t\right)  \left\vert \widehat{f_{d}}\left(  mt\right)
\right\vert ^{2}dt &  =m^{d}\int w\left(  \xi/m\right)  \left\vert
\widehat{f_{d}}\left(  \xi\right)  \right\vert ^{2}d\xi\\
&  \leq c_{m}m^{d}\int w\left(  \xi\right)  \left\vert \widehat{f_{d}}\left(
\xi\right)  \right\vert ^{2}d\xi\\
&  =c_{m}m^{d}\left\Vert f_{d}\right\Vert _{w,0}^{2}.
\end{align*}
\smallskip

\textbf{Parts 2 and 3} Use the fact that $\chi\ast u=\left(  \widehat{u}%
\widehat{\chi}\right)  ^{\vee}$.\smallskip

\textbf{Part 4} From \ref{a4.8} and part 1, continuity is guaranteed by%
\begin{align*}
\left\Vert \sigma_{m}\left(  \chi\ast\right)  \mathcal{E}_{0}u\right\Vert
_{w,0}\leq\sqrt{c_{m}m^{d}}\left\Vert \left(  \chi\ast\right)  \mathcal{E}%
_{0}u\right\Vert _{w,0} &  \leq\sqrt{C_{\chi}c_{m}m^{d}}\left\Vert
\mathcal{E}_{0}u\right\Vert _{W^{n\mathbf{1}}}\\
&  \leq\sqrt{C_{\chi}c_{m}m^{d}}\left\Vert \mathcal{E}_{0}u\right\Vert
_{W^{n\mathbf{1}}\left(  \left(  -1,1\right)  ^{d}\right)  }.
\end{align*}
\smallskip

\textbf{Part 5}%
\begin{align*}
\left\Vert \sigma_{1/m}g\right\Vert _{W^{n\mathbf{1}}\left(  \Sigma\right)
}^{2}=\sum\limits_{\alpha\leq n\mathbf{1}}\left\Vert D^{\alpha}\sigma
_{1/m}g\right\Vert _{L^{2}\left(  \Sigma\right)  }^{2}  & =\sum\limits_{\alpha
\leq n\mathbf{1}}m^{2\left\vert \alpha\right\vert }\left\Vert \sigma
_{1/m}D^{\alpha}g\right\Vert _{L^{2}\left(  \Sigma\right)  }^{2}\\
& =\max\left\{  1,m^{2nd}\right\}  \sum\limits_{\alpha\leq n\mathbf{1}%
}\left\Vert \sigma_{1/m}D^{\alpha}g\right\Vert _{L^{2}\left(  \Sigma\right)
}^{2},
\end{align*}

but%
\[
\left\Vert \sigma_{1/m}D^{\alpha}g\right\Vert _{L^{2}\left(  \Sigma\right)
}^{2}=\int_{\Sigma}\left\vert \left(  D^{\alpha}g\right)  \left(  mx\right)
\right\vert ^{2}dx=m^{-d}\int_{m\Sigma}\left\vert D^{\alpha}g\right\vert
^{2}=m^{-d}\left\Vert D^{\alpha}g\right\Vert _{L^{2}\left(  m\Sigma\right)
}^{2},
\]

so%
\[
\left\Vert \sigma_{1/m}g\right\Vert _{W^{n\mathbf{1}}\left(  \Sigma\right)
}^{2}\leq m^{-d}\max\left\{  1,m^{2nd}\right\}  \sum\limits_{\alpha\leq
n\mathbf{1}}\left\Vert D^{\alpha}g\right\Vert _{L^{2}\left(  m\Sigma\right)
}^{2}.
\]

\end{proof}

Now we are ready to prove:

\begin{theorem}
\label{Thm_ex_data_fn_3}Suppose the weight function $w$ and region $\Omega$ satisfy:

\begin{enumerate}
\item $w\in W03$.

\item There exists $\chi\in S^{\prime}$ with properties \ref{a1.05} and
\ref{a1.04} and also \ref{a1.052} holds for some integer $n\geq1$.

\item $\Omega$ is a bounded region with the segment property.

\item $w$ satisfies \ref{a1.055} for some integer $m\geq1$.
\end{enumerate}

We can choose an open cube $C$ with center $c$ such that $\overline{\Omega
}\subset C$ and $\frac{1}{m}\left(  C-c\right)  \subset\left[  -1,1\right]
^{d}$. Choose $\varepsilon>0$ such that $0<\varepsilon/m\leq
\operatorname*{dist}\left(  \Omega,\mathbb{R}^{d}\setminus C\right)  $ and
let
\begin{equation}
\Omega^{0}=\sigma_{1/m}\tau_{-c}\Omega=\frac{1}{m}\left(  C-c\right)
.\label{1.020}%
\end{equation}

Since $\Omega^{0}$ also has the segment property, by part 3 of Lemma
\ref{Lem_SobolevSpProperty} there exists a continuous extension operator
$\mathcal{E}_{\Omega^{0};\varepsilon/m}:W^{m\mathbf{1}}\left(  \Omega
^{0}\right)  \rightarrow W_{0}^{m\mathbf{1}}\left(  \Omega_{\varepsilon/m}%
^{0}\right)  $. Next define the mapping $E$ of Figure
\ref{Fig_exten_commut_mu_convol} by:%
\begin{equation}
Ef:=\tau_{c}\sigma_{m}\left(  \chi\ast\right)  \mathcal{E}_{0}\mathcal{E}%
_{\Omega^{0};\varepsilon/m}\sigma_{1/m}\tau_{-c}f,\quad f\in W^{n\mathbf{1}%
}\left(  \Omega\right)  .\label{1.040}%
\end{equation}

Then $E:W^{n\mathbf{1}}\left(  \Omega\right)  \rightarrow X_{w}^{0}$ is a
continuous linear extension such that%
\begin{equation}
Ef=\tau_{c}\sigma_{m}\left(  \chi\ast\right)  \sigma_{1/m}\tau_{-c}%
\mathcal{E}_{0}\mathcal{E}_{\Omega;\varepsilon}f,\quad f\in W^{n\mathbf{1}%
}\left(  \Omega\right)  ,\label{a1.03}%
\end{equation}

and%
\begin{equation}
\operatorname*{supp}Ef\subset m\operatorname*{supp}\eta+\Omega_{\varepsilon
}\subset c+m\left(  \left[  -1,1\right]  ^{d}+\operatorname*{supp}\eta\right)
.\label{a1.57}%
\end{equation}

Regarding the continuity of $E$:%
\begin{equation}
\left\Vert Ef\right\Vert _{w,0}\leq\sqrt{C_{\chi}c_{m}}m^{nd}\left\Vert
\mathcal{E}_{\Omega;\varepsilon}f\right\Vert _{W_{0}^{n\mathbf{1}}\left(
\Omega_{\varepsilon}\right)  }\leq\sqrt{C_{\chi}c_{m}}m^{nd}\left\Vert
\mathcal{E}_{\Omega;\varepsilon}\right\Vert \left\Vert f\right\Vert
_{W^{n\mathbf{1}}\left(  \Omega\right)  },\quad f\in W^{n\mathbf{1}}\left(
\Omega\right)  .\label{a1.58}%
\end{equation}

\end{theorem}

\begin{proof}
The relevant commutative diagram is Figure \ref{Fig_exten_commut_mu_convol}
and the relevant sets and mappings are given in Figure \ref{Fig_exten_sets_m}.

We first show that the mapping \ref{1.040} makes sense:\smallskip

\textbf{1}) From part 1 of Lemma \ref{Lem_SobolevSpProperty}, $\sigma
_{1/m}\tau_{-c}:W^{n\mathbf{1}}\left(  \Omega\right)  \rightarrow
W^{n\mathbf{1}}\left(  \Omega^{0}\right)  $ is a homeomorphism.\smallskip

\textbf{2)} The inclusion \ref{1.020} implies $\Omega^{0}\subset\left[
-1,1\right]  ^{d}$ and since $0<\varepsilon<\operatorname*{dist}\left(
\Omega,\mathbb{R}^{d}\setminus C\right)  $, we have

$0<\varepsilon/m<\operatorname*{dist}\left(  \sigma_{1/m}\tau_{-c}%
\Omega,\mathbb{R}^{d}\setminus\sigma_{1/m}\tau_{-c}C\right)  $ i.e.
$0<\varepsilon/m<\operatorname*{dist}\left(  \Omega^{0},\mathbb{R}%
^{d}\setminus\left(  -1,1\right)  \right)  $, and so $\Omega_{\varepsilon
/m}^{0}\subset\left(  -1,1\right)  ^{d}$. Hence by part 3 of Lemma
\ref{Lem_SobolevSpProperty}, $\sigma_{1/m}\tau_{-c}f\in W^{n\mathbf{1}}\left(
\Omega^{0}\right)  $ can be extended by the continuous operator $\mathcal{E}%
_{\Omega^{0};\varepsilon/m}$ to $\mathbb{R}^{d}$ as a function in
$W_{0}^{n\mathbf{1}}\left(  \Omega_{\varepsilon/m}^{0}\right)  $.\smallskip

\textbf{3)} From \ref{a1.002}, $\left(  \chi\ast\right)  \mathcal{E}_{0}%
:W_{0}^{n\mathbf{1}}\left(  \left(  -1,1\right)  ^{d}\right)  \rightarrow
X_{w}^{0}$ is a continuous extension.\smallskip

\textbf{4)} From part 1 of Lemma \ref{Lem_Thm_ex_data_fn_3}, $\tau_{c}%
\sigma_{m}:X_{w}^{0}\rightarrow X_{w}^{0}$ is continuous.\smallskip

These four results imply that $E:W^{n\mathbf{1}}\left(  \Omega\right)
\rightarrow X_{w}^{0}$ is continuous.

Next, part 3 of Lemma \ref{Lem_SobolevSpProperty} enables us to write%
\[
Ef=\tau_{c}\sigma_{m}\left(  \chi\ast\right)  \mathcal{E}_{0}\mathcal{E}%
_{\Omega^{0};\varepsilon/m}\sigma_{1/m}\tau_{-c}f=\tau_{c}\sigma_{m}\left(
\chi\ast\right)  \sigma_{1/m}\tau_{-c}\mathcal{E}_{0}\mathcal{E}%
_{\Omega;\varepsilon}f,
\]

proving \ref{a1.03}. From part 2 of Lemma \ref{Lem_Thm_ex_data_fn_3},
$\chi\ast\left(  \sigma_{1/m}u\right)  =m^{-d}\sigma_{1/m}\left(  \left(
\sigma_{m}\chi\right)  \ast\right)  u$ so%
\begin{align}
Ef=\tau_{c}\sigma_{m}\left(  \chi\ast\right)  \sigma_{1/m}\tau_{-c}%
\mathcal{E}_{0}\mathcal{E}_{\Omega;\varepsilon}f  & =\tau_{c}\sigma_{m}%
\sigma_{1/m}\left(  \left(  m^{-d}\sigma_{m}\chi\right)  \ast\right)
\tau_{-c}\mathcal{E}_{0}\mathcal{E}_{\Omega;\varepsilon}f\nonumber\\
& =\tau_{c}\left(  \left(  m^{-d}\sigma_{m}\chi\right)  \ast\right)  \tau
_{-c}\mathcal{E}_{0}\mathcal{E}_{\Omega;\varepsilon}f\nonumber\\
& =\tau_{c}\tau_{-c}\left(  \left(  m^{-d}\sigma_{m}\chi\right)  \ast\right)
\mathcal{E}_{0}\mathcal{E}_{\Omega;\varepsilon}f\nonumber\\
& =\left(  m^{-d}\sigma_{m}\chi\right)  \ast\left(  \mathcal{E}_{0}%
\mathcal{E}_{\Omega;\varepsilon}f\right)  .\label{a1.010}%
\end{align}

We will now show that the convolution operator $\left(  m^{-d}\sigma_{m}%
\chi\right)  \ast$ is an extension operator by showing that it satisfies the
same conditions as $\chi$ i.e. \ref{a1.04}, \ref{a1.05} and \ref{a1.052}.

But%
\begin{align*}
m^{-d}\sigma_{m}\chi-\left(  2\pi\right)  ^{d/2}\delta & =m^{-d}\sigma
_{m}\left(  \eta+\left(  2\pi\right)  ^{d/2}\delta\right)  -\left(
2\pi\right)  ^{d/2}\delta\\
& =m^{-d}\sigma_{m}\eta+\left(  2\pi\right)  ^{d/2}m^{-d}\sigma_{m}%
\delta-\left(  2\pi\right)  ^{d/2}\delta\\
& =m^{-d}\sigma_{m}\eta+\left(  2\pi\right)  ^{d/2}\delta-\left(  2\pi\right)
^{d/2}\delta\\
& =m^{-d}\sigma_{m}\eta,
\end{align*}

so that%
\begin{align*}
\operatorname*{supp}\left(  m^{-d}\sigma_{m}\chi-\left(  2\pi\right)
^{d/2}\delta\right)  =\operatorname*{supp}\left(  m^{-d}\sigma_{m}\eta\right)
&  =\operatorname*{supp}\left(  \sigma_{m}\eta\right) \\
&  =\sigma_{m}\operatorname*{supp}\eta\\
&  =m\operatorname*{supp}\eta.
\end{align*}

Now \ref{a1.04} and $m\geq1$ imply
\[
\operatorname*{supp}\left(  m^{-d}\sigma_{m}\chi-\left(  2\pi\right)
^{d/2}\delta\right)  \subseteq\mathbb{R}^{d}\setminus2m\left(  -1,1\right)
^{d},
\]

which means that $\left(  m^{-d}\sigma_{m}\chi\right)  \ast$ is an extension
and hence that $Ef=\left(  m^{-d}\sigma_{m}\chi\right)  \ast\left(
\mathcal{E}_{0}\mathcal{E}_{\Omega;\varepsilon}f\right)  $ is an extension of
$f$ beyond $\Omega$.

\textbf{Regarding the support of }$E$:%
\begin{align*}
\operatorname*{supp}Ef  & =\operatorname*{supp}\left(  m^{-d}\sigma_{m}%
\chi\right)  \ast\left(  \mathcal{E}_{0}\mathcal{E}_{\Omega;\varepsilon
}f\right) \\
& \subset\operatorname*{supp}\left(  m^{-d}\sigma_{m}\chi\right)
+\operatorname*{supp}\left(  \mathcal{E}_{0}\mathcal{E}_{\Omega;\varepsilon
}f\right) \\
& =\operatorname*{supp}\left(  m^{-d}\sigma_{m}\chi\right)
+\operatorname*{supp}\mathcal{E}_{\Omega;\varepsilon}f\\
& \subset\operatorname*{supp}\left(  m^{-d}\sigma_{m}\chi\right)
+\Omega_{\varepsilon}\\
& =\operatorname*{supp}\left(  \sigma_{m}\chi\right)  +\Omega_{\varepsilon}\\
& =m\operatorname*{supp}\chi+\Omega_{\varepsilon},
\end{align*}

but from \ref{a1.05} i.e. $\chi=\eta+\left(  2\pi\right)  ^{d/2}\delta$ we
have $\operatorname*{supp}\chi\subset\operatorname*{supp}\eta\cup\left\{
\mathbf{0}\right\}  \subset\operatorname*{supp}\eta$ so that%
\[
\operatorname*{supp}Ef\subset m\operatorname*{supp}\eta+\Omega_{\varepsilon}.
\]

But $\Omega_{\varepsilon}\subset C$ and $\frac{1}{m}\left(  C-c\right)
\subset\left[  -1,1\right]  ^{d}$ implies $\Omega_{\varepsilon}\subset
c+m\left[  -1,1\right]  ^{d}$ and consequently%
\[
\operatorname*{supp}Ef\subset m\operatorname*{supp}\eta+\Omega_{\varepsilon
}\subset c+m\left(  \left[  -1,1\right]  ^{d}+\operatorname*{supp}\eta\right)
.
\]
\smallskip

\textbf{Regarding the continuity of }$E$:%
\[
\left(  m^{-d}\sigma_{m}\chi\right)  ^{\wedge}\left(  \xi\right)
=m^{-d}\left(  \sigma_{m}\chi\right)  ^{\wedge}\left(  \xi\right)
=m^{-d}m^{d}\widehat{\chi}\left(  m\xi\right)  =\widehat{\chi}\left(
m\xi\right)  ,
\]

and inequality \ref{a1.055} yields%
\begin{align*}
w\left(  \xi\right)  \left\vert \left(  m^{-d}\sigma_{m}\chi\right)  ^{\wedge
}\left(  \xi\right)  \right\vert ^{2}=w\left(  \xi\right)  \left\vert
\widehat{\chi}\left(  m\xi\right)  \right\vert ^{2}  & \leq c_{m}w\left(
m\xi\right)  \left\vert \widehat{\chi}\left(  m\xi\right)  \right\vert ^{2}\\
& \leq c_{m}C_{\chi}\prod\limits_{k=1}^{d}\left(  1+m^{2}\xi_{k}^{2}\right)
^{n}\\
& \leq m^{2nd}c_{m}C_{\chi}w_{n}\left(  \xi\right)  .
\end{align*}

Thus%
\[
\left(  m^{-d}\sigma_{m}\chi\right)  \ast:W_{0}^{n\mathbf{1}}\left(  \left(
-1,1\right)  ^{d}\right)  \rightarrow X_{w}^{0}\text{ }is\text{ }a\text{
}continuous\text{ }extension.
\]

Starting from \ref{a1.03}, part 4 of Lemma \ref{Lem_Thm_ex_data_fn_3} implies%
\begin{align*}
\left\Vert Ef\right\Vert _{w,0} &  =\left\Vert \tau_{c}\sigma_{m}\left(
\chi\ast\right)  \sigma_{1/m}\tau_{-c}\mathcal{E}_{0}\mathcal{E}%
_{\Omega;\varepsilon}f\right\Vert _{w,0}\\
&  =\left\Vert \sigma_{m}\left(  \chi\ast\right)  \sigma_{1/m}\tau
_{-c}\mathcal{E}_{0}\mathcal{E}_{\Omega;\varepsilon}f\right\Vert _{w,0}\\
&  \leq\sqrt{C_{\chi}c_{m}m^{d}}\left\Vert \sigma_{1/m}\tau_{-c}%
\mathcal{E}_{0}\mathcal{E}_{\Omega;\varepsilon}f\right\Vert _{W^{n\mathbf{1}}%
}\\
&  =\sqrt{C_{\chi}c_{m}m^{d}}\left\Vert \sigma_{1/m}\tau_{-c}\mathcal{E}%
_{0}\mathcal{E}_{\Omega;\varepsilon}f\right\Vert _{W^{n\mathbf{1}}},
\end{align*}

and, since $m\geq1$, by part 5 of Lemma \ref{Lem_Thm_ex_data_fn_3} and then
part 1 of Lemma \ref{Lem_SobolevSpProperty},%
\begin{align*}
\left\Vert Ef\right\Vert _{w,0}  & \leq\sqrt{C_{\chi}c_{m}m^{d}}m^{-\frac
{d}{2}}\max\left\{  1,m^{nd}\right\}  \left\Vert \tau_{-c}\mathcal{E}%
_{0}\mathcal{E}_{\Omega;\varepsilon}f\right\Vert _{W^{n\mathbf{1}}}\\
& =\sqrt{C_{\chi}c_{m}}m^{nd}\left\Vert \tau_{-c}\mathcal{E}_{0}%
\mathcal{E}_{\Omega;\varepsilon}f\right\Vert _{W^{n\mathbf{1}}}\\
& =\sqrt{C_{\chi}c_{m}}m^{nd}\left\Vert \mathcal{E}_{0}\mathcal{E}%
_{\Omega;\varepsilon}f\right\Vert _{W^{n\mathbf{1}}}\\
& =\sqrt{C_{\chi}c_{m}}m^{nd}\left\Vert \mathcal{E}_{\Omega;\varepsilon
}f\right\Vert _{W_{0}^{n\mathbf{1}}\left(  \Omega_{\varepsilon}\right)  }.
\end{align*}

Finally, from part 3 of Lemma \ref{Lem_SobolevSpProperty},%
\[
\left\Vert Ef\right\Vert _{w,0}\leq\sqrt{C_{\chi}c_{m}}m^{nd}\left\Vert
\mathcal{E}_{\Omega;\varepsilon}f\right\Vert _{W_{0}^{n\mathbf{1}}\left(
\Omega_{\varepsilon}\right)  }\leq\sqrt{C_{\chi}c_{m}}m^{nd}\left\Vert
\mathcal{E}_{\Omega;\varepsilon}\right\Vert _{op}\left\Vert f\right\Vert
_{W^{n\mathbf{1}}\left(  \Omega\right)  }.
\]

\end{proof}

We now have a more general version of Corollary \ref{Cor_Thm_data_fn_2}.

\begin{corollary}
\label{Cor_Thm_data_fn_3}Suppose the weight function $w$ and bounded region
$\Omega$ satisfy assumptions 1 to 4 of Theorem \ref{Thm_ex_data_fn_3}.Then
$W^{n\mathbf{1}}\left(  \Omega\right)  \hookrightarrow X_{\sigma_{\lambda}%
w}^{0}\left(  \Omega\right)  $ where $\sigma_{\lambda}$ is any scalar dilation operator.
\end{corollary}

\begin{proof}
In Theorem\textbf{\ }\ref{Thm_ex_data_fn_3} a continuous extension
$E:W^{n\mathbf{1}}\left(  \Omega\right)  \hookrightarrow X_{w}^{0}$ was
constructed so that by Definition \ref{Def_Xow(open_set)}, $r_{\Omega
}E:W^{n\mathbf{1}}\left(  \Omega\right)  \rightarrow X_{w}^{0}\left(
\Omega\right)  $ is a continuous inclusion. Since $\lambda\Omega$ also permits
one of the extension operators $\mathcal{E}_{\lambda\Omega}$ described in part
5 of Remark \ref{Rem_SobolevSpace2} we also have $W^{n\mathbf{1}}\left(
\lambda\Omega\right)  \hookrightarrow X_{w}^{0}\left(  \lambda\Omega\right)
$, and an application of the dilation results of part 4 of Theorem
\ref{Thm_data_func} and part 1 of Lemma \ref{Lem_SobolevSpProperty} enables us
to write%
\[
W^{n\mathbf{1}}\left(  \Omega\right)  \overset{\sigma_{\lambda}%
}{\longrightarrow}W^{n\mathbf{1}}\left(  \lambda\Omega\right)  \hookrightarrow
X_{w}^{0}\left(  \lambda\Omega\right)  \overset{\sigma_{\lambda}%
^{-1}}{\longrightarrow}X_{\sigma_{\lambda}w}^{0}\left(  \Omega\right)  ,
\]

which proves this corollary.
\end{proof}

The next result characterizes the data functions locally.

\begin{corollary}
\label{Cor_2_Thm_data_fn_3}\textbf{Local data functions} Suppose the
assumptions of Theorem \ref{Thm_ex_data_fn_3} hold and we also assume that
there exist constants $c_{\alpha}>0$ such that%
\begin{equation}
w\left(  \xi\right)  \geq c_{\alpha}\xi^{2\alpha},\quad a.e.\text{ }on\text{
}\mathbb{R}^{d},\text{ }for\text{ }\alpha\leq n\mathbf{1},\label{1.39}%
\end{equation}

or equivalently: there exists a constant $c_{w}>0$ such that%
\begin{equation}
w\left(  \xi\right)  \geq c_{w}w_{n}\left(  \xi\right)  ,\quad a.e.\text{
}on\text{ }\mathbb{R}^{d},\label{1.391}%
\end{equation}

where the weight function $w_{n}$ is defined by \ref{1.052}. Then:

\begin{enumerate}
\item For each $\lambda\in\mathbb{R}_{+}^{1}$, $X_{\sigma_{\lambda}w}%
^{0}\left(  \Omega\right)  =W^{n\mathbf{1}}\left(  \Omega\right)  $ as sets
and the norms are equivalent.

\item $X_{\sigma_{\lambda}w}^{0}\left(  \Omega\right)  \hookrightarrow
r_{\Omega}C_{B}^{\left(  n-1\right)  \mathbf{1}}$.
\end{enumerate}
\end{corollary}

\begin{proof}
\textbf{Part 1} Inspection of the proof of Lemma
\ref{Lem_Thm_Xow(O)_eq_Hn(O)_dim1} reveals that it still holds when the weight
function $w$ satisfies property \ref{1.39}. Hence $X_{\sigma_{\lambda}w}%
^{0}\left(  \Omega\right)  \hookrightarrow W^{n\mathbf{1}}\left(
\Omega\right)  $ where $\sigma_{\lambda}$ is any scalar dilation operator.
Also, by Corollary \ref{Cor_Thm_data_fn_3}, $W^{n\mathbf{1}}\left(
\Omega\right)  \hookrightarrow X_{\sigma_{\lambda}w}^{0}\left(  \Omega\right)
$. These two results imply that $X_{\sigma_{\lambda}w}^{0}\left(
\Omega\right)  =W^{n\mathbf{1}}\left(  \Omega\right)  $ as sets and that the
norms are equivalent.\smallskip

\textbf{Part 2} From the proof of part 2 of Theorem
\ref{Thm_int_Xow(O)_eq_Hn(O)_dim1}, $X_{\sigma_{\lambda}w_{s}}^{0}\left(
\Omega\right)  \hookrightarrow r_{\Omega}C_{B}^{\left(  n-1\right)
\mathbf{1}}$.
\end{proof}

\begin{remark}
Conditions \ref{1.391} and \ref{a1.052} imply that $\left\vert \widehat{\chi
}\left(  \xi\right)  \right\vert \leq\sqrt{C_{\chi}/c_{w}}$ and hence that
$\left\vert \widehat{\eta}\left(  \xi\right)  \right\vert \leq1+\sqrt{C_{\chi
}/c_{w}}$.
\end{remark}

\begin{remark}
\textbf{Interpolated weight functions} The results of this section can be
applied to the interpolated weight function $w_{1}^{1-t}w_{2}^{t}$ introduced
in Theorem \ref{Thm_interpol_wt_fns} i.e. if both $w_{1}$ and $w_{2}$ satisfy
the specified conditions then $w_{1}^{1-t}w_{2}^{t}$ also does for all
$t\in\left[  0,1\right]  $.
\end{remark}

\begin{example}
\textbf{Tensor product weight functions and specifically central difference
weight functions} These are studied in Section \ref{Sect_CntDifWtFn_DataFuncs}.
\end{example}

\subsection{Applications: local pointwise error estimates from global ones
estimates \label{SbSect_global_error_to_local}}

The characterizations of this section can be applied to \textbf{global
pointwise error estimates} of the form%
\[
\left\vert \mathcal{S}f\left(  x\right)  -f\left(  x\right)  \right\vert \leq
c^{\prime}\left\Vert f\right\Vert _{w,0}h^{s},\quad x\in\Omega,\text{ }f\in
X_{w}^{0},
\]

where $\mathcal{S}$ is the minimal norm interpolant operator or the Exact
smoother operator of Chapter \ref{Ch_Exact_smth}, $s$ is the order of
convergence and $h$ is the spherical cavity size of the data. Now if $u\in
W^{\mathbf{1}}\left(  \Omega\right)  $, the extension operator $E$ of Theorem
\ref{Thm_exten_op_Wloc_to_X0ws}, and more generally of Theorem
\ref{Thm_ex_data_fn_3}, can be applied to give the local estimate%
\[
\left\vert \mathcal{S}u\left(  x\right)  -u\left(  x\right)  \right\vert
=\left\vert \mathcal{S}Eu\left(  x\right)  -Eu\left(  x\right)  \right\vert
\leq c^{\prime}\left\Vert Eu\right\Vert _{w,0}h^{s}\leq c^{\prime}\left\Vert
E\right\Vert \left\Vert u\right\Vert _{W^{\mathbf{1}}\left(  \Omega\right)
}h^{s}.
\]

See Remark \ref{Rem_Cor_Thm_int_Xow(O)_eq_Hn(O)_dim1} concerning the
definition of the local smoothing operator $S_{X}^{\Omega}$.

\chapter{Local interpolation errors using tempered distribution Taylor series
\label{Ch_interp_err_Taylor_temper_distrib}}

\section{Introduction}

In this chapter we will derive another set of pointwise error estimates for
the interpolant which are then applied to the Sobolev splines and the extended
B-splines to obtain better convergence estimates.

We start by deriving a \textbf{tempered distribution Taylor series} expansion
with a Fourier transform remainder (\ref{a1.55}, \ref{a50.5}) using a
1-dimensional Taylor series expansion applied to $e^{ia\xi}$. Several
interesting,.weak Taylor series expansions are derived, both global and local
(Lemmas \ref{Lem_Taylor_estim_C_L1} to \ref{Lem_Taylor_estim_L1loc_L1}).
Considering the cases $\underline{\kappa}<1$ and $\underline{\kappa}\geq1$
separately the remainder is then estimated in terms of the Fourier transform
norm $\left\Vert \cdot\right\Vert _{w,0}$ and a factor involving $1/\sqrt{w}$.
In Section \ref{Sect_rem_estim_radial_wt_fn} general remainder estimates are
obtained for the data functions generated by a radial weight function. In
Section \ref{Sect_rem_estim_ten_prod} general remainder estimates are obtained
for data functions generated by tensor product weight functions. In this case
reproducing kernel/Riesz representer techniques are used to bound the
derivatives near the origin e.g. \ref{a3.02}.

In the next two sections these remainder estimates are used to derive global
and local orders of convergence formulas for the interpolant, with separate
sections being devoted to the cases $\underline{\kappa}<1$ and
$\underline{\kappa}\geq1$ respectively. The case $\underline{\kappa}<1$ does
not make (explicit) use of unisolvency. The $\underline{\kappa}\geq1$ case
uses the techniques of multipoint Taylor series expansions combined with
Lagrange interpolation on minimal unisolvent sets of data points. These
results are applied to the radial Sobolev splines and the tensor product
extended B-splines and these are summarized in Table
\ref{Tbl_NonUnisolvConverg_revisit} for $\underline{\kappa}<1$ and Table
\ref{Tbl_UnisolvConverg_revisit} for $\underline{\kappa}\geq1$.

\section{Pointwise Taylor series remainder estimates for data
functions\label{Sect_Taylor_series_data_fn}}

We start by deriving the tempered distribution Taylor series expansion
\ref{a1.55} with Fourier transform remainder \ref{a50.5} by means of the
1-dimensional Taylor series expansion \ref{a1.12} applied to $e^{ia\xi}$. The
Fourier transform is used to form $\left\Vert \cdot\right\Vert _{w,0}$ and
estimates are derived by means of Young's convolution inequality.

Indeed, for $x\in\mathbb{R}^{1}$,%
\begin{align*}
e^{x}=\sum_{k\leq n}\frac{x^{n}}{n!}+\left(  \mathcal{R}_{n+1}e^{x}\right)
\left(  0,x\right)   &  =\sum_{k\leq n}\frac{x^{k}}{k!}+\frac{x^{n+1}}{n!}%
\int_{0}^{1}\left(  1-t\right)  ^{n}\left(  D^{n+1}e^{s}\right)  \left(
tx\right)  dt\\
&  =\sum_{k\leq n}\frac{x^{k}}{k!}+\frac{x^{n+1}}{n!}\int_{0}^{1}\left(
1-t\right)  ^{n}e^{tx}dt,
\end{align*}

so that if we define $g_{n}$ by%
\begin{equation}
g_{n}\left(  t\right)  =\left\{
\begin{array}
[c]{ll}%
0, & t<0,\\
\left(  1-t\right)  ^{n}, & 0\leq t\leq1,\\
0, & t>1,
\end{array}
\right.  \text{\quad}n=0,1,2,\ldots,\label{a117}%
\end{equation}

it follows that%
\begin{align}
e^{ix}=\sum_{k\leq n}\frac{\left(  ix\right)  ^{k}}{k!}+\frac{\left(
ix\right)  ^{n+1}}{n!}\int_{0}^{1}\left(  1-t\right)  ^{n}e^{itx}dt  &
=\sum_{k\leq n}\frac{\left(  ix\right)  ^{k}}{k!}+\frac{\left(  ix\right)
^{n+1}}{n!}\int_{0}^{1}g_{n}\left(  t\right)  e^{itx}dt\nonumber\\
& =\sum_{k\leq n}\frac{\left(  ix\right)  ^{k}}{k!}+\frac{\sqrt{2\pi}}%
{n!}\left(  ix\right)  ^{n+1}\overline{\widehat{g_{n}}}\left(  x\right)
\nonumber\\
& =\sum_{k\leq n}\frac{\left(  ix\right)  ^{k}}{k!}+\frac{\sqrt{2\pi}}%
{n!}\left(  ix\right)  ^{n+1}\overset{\vee}{g_{n}}\left(  x\right)
,\label{a2.35}%
\end{align}

and%
\[
\overset{\vee}{g_{n}}\left(  x\right)  =\frac{n!}{\sqrt{2\pi}}\frac
{e^{ix}-\sum_{k\leq n}\frac{\left(  ix\right)  ^{k}}{k!}}{\left(  ix\right)
^{n+1}}.
\]

Also, for $a,\xi\in\mathbb{R}^{d}$,%
\begin{align}
e^{ia\xi}  & =\sum_{k\leq n}\frac{\left(  ia\xi\right)  ^{k}}{k!}%
+\frac{\left(  ia\xi\right)  ^{n+1}}{n!}\int_{0}^{1}\left(  1-t\right)
^{n}e^{ita\xi}dt\nonumber\\
& =\sum_{k\leq n}\frac{\left(  ia\xi\right)  ^{k}}{k!}+\frac{\sqrt{2\pi}}%
{n!}\left(  ia\xi\right)  ^{n+1}\overline{\widehat{g_{n}}}\left(  a\xi\right)
.\label{a114}%
\end{align}

Clearly%
\[
g_{n}\left(  t\right)  =\left(  1-t\right)  ^{n-k}g_{k}\left(  t\right)  .
\]

More properties of $g_{n}$ and $\widehat{g_{n}}$ are:

\begin{lemma}
\label{Lem_gm_properties}The function $g_{n}$ given by \ref{a117} has the
following properties:

\begin{enumerate}
\item $\left\Vert g_{n}\right\Vert _{1}=\frac{1}{n+1}$ and $\left\Vert
g_{n}\right\Vert _{2}=\frac{1}{\sqrt{2n+1}}$.

\item $\widehat{g_{n}}\in C_{B}^{\infty}$ and%
\begin{equation}
\left\vert \widehat{g_{n}}\left(  t\right)  \right\vert \leq\frac{1}%
{\sqrt{2\pi}}\frac{1}{n+1},\quad n=0,1,2,\ldots\label{a50.3}%
\end{equation}

\item
\[
Dg_{n}=\left\{
\begin{array}
[c]{ll}%
\delta-\delta\left(  \cdot-1\right)  , & n=0,\\
-ng_{n-1}+\delta, & n=1,2,3\ldots,
\end{array}
\right.
\]

and%
\begin{equation}
\widehat{Dg_{n}}=\left\{
\begin{array}
[c]{ll}%
\frac{2i}{\sqrt{2\pi}}e^{-it/2}\sin\frac{t}{2}, & n=0,\\
-n\widehat{g_{n-1}}+\frac{1}{\sqrt{2\pi}}, & n=1,2,3\ldots.
\end{array}
\right. \label{a2.04}%
\end{equation}

\item If we define%
\begin{equation}
c_{n}=2+\frac{1}{n+1},\quad n=0,1,2,\ldots,\label{a2.08}%
\end{equation}

then
\begin{equation}
\left\vert \widehat{g_{n}}\left(  t\right)  \right\vert \leq\frac{c_{n}}%
{\sqrt{2\pi}}\frac{1}{1+\left\vert t\right\vert },\text{\quad}t\in
\mathbb{R}^{1},\text{ }n=0,1,2,\ldots\label{a2.07}%
\end{equation}

\item For each $a\in\mathbb{R}^{d}$ such that $a\centerdot\neq\mathbf{0}$ we
have $\widehat{g_{n}}\left(  a\xi\right)  \in C_{B}^{\infty}\left(
\mathbb{R}^{d}\right)  $.

\item
\[
\left\Vert g_{n}\ast\widetilde{Dg_{n}}\right\Vert _{1}\leq\left\{
\begin{array}
[c]{ll}%
\frac{2}{\sqrt{2\pi}}, & n=0,\\
\left(  1+\frac{1}{\sqrt{2\pi}}\right)  \frac{1}{n+1}, & n\geq1.
\end{array}
\right.
\]

\item
\[
\frac{\left(  -1\right)  ^{k}}{n!}D^{k}g_{n}=\frac{g_{n-k}}{\left(
n-k\right)  !}-\frac{\delta}{\left(  n-k+1\right)  !}+\frac{D\delta}{\left(
n-k+2\right)  !}-\ldots+\left(  -1\right)  ^{k}\frac{D^{k-1}\delta}{n!},\quad
k\leq n.
\]

\[
\frac{\left(  -1\right)  ^{n+1}}{n!}D^{n+1}g_{n}=\delta\left(  \cdot-1\right)
-\delta+\frac{D\delta}{1!}-\frac{D^{2}\delta}{2!}+\ldots+\left(  -1\right)
^{n+1}\frac{D^{n}\delta}{n!}.
\]

\item
\[
\widehat{D^{n+j+1}g_{n}}\left(  t\right)  =\left(  -1\right)  ^{n+j+1}%
\frac{n!}{\sqrt{2\pi}}\left(  -it\right)  ^{j}\left(  e^{-it}-\sum_{k=0}%
^{n}\frac{\left(  -it\right)  ^{k}}{k!}\right)  ,\quad j,n\geq0,
\]

and%
\[
\left\vert \widehat{D^{n+j+1}g_{n}}\left(  t\right)  \right\vert \leq\frac
{n!}{\sqrt{2\pi}}\left\vert t\right\vert ^{j}\left(  1+\sum_{k=0}^{n}%
\frac{\left\vert t\right\vert ^{k}}{k!}\right)  ,\quad j,n\geq0,
\]

and%
\[
\widehat{D^{n+j+1}g_{n}}\in C_{BP}^{\infty}\cap C_{\emptyset,n+j+1}^{\infty}.
\]

\end{enumerate}
\end{lemma}

\begin{proof}
\textbf{Part 1} A simple calculation.\medskip

\textbf{Part 2} Since $g_{n}$ is a distribution with bounded support,
$\widehat{g_{n}}\in C_{BP}^{\infty}$ and the estimates \ref{a50.3} follow
directly from the formula, $\widehat{g_{n}}\left(  t\right)  =\frac{1}%
{\sqrt{2\pi}}\int_{0}^{1}e^{-ist}\left(  1-s\right)  ^{n}ds$.\medskip

\textbf{Part 3} If $n=0$ then $Dg_{0}=\left\{  Dg_{0}\right\}  +\delta
-\delta\left(  \cdot-1\right)  =0+\delta-\delta\left(  \cdot-1\right)  $.

Hence $\widehat{Dg_{0}}=\widehat{\delta}-\widehat{\delta\left(  \cdot
-1\right)  }=\frac{1}{\sqrt{2\pi}}\left(  1-e^{-it}\right)  =\frac{2i}%
{\sqrt{2\pi}}e^{-it/2}\left(  \frac{e^{it/2}-e^{-it/2}}{2i}\right)  =\frac
{2i}{\sqrt{2\pi}}e^{-it/2}\sin\frac{t}{2}$.

If $n\geq1$ then $Dg_{n}=\left\{  Dg_{n}\right\}  +\delta=-n\left(
1-t\right)  ^{n-1}g_{0}+\delta=-ng_{n-1}+\delta$.\medskip

\textbf{Part 4} Suppose $n\geq1$. Then from part 3 and then part 2:%
\begin{align*}
\left\vert t\right\vert \left\vert \widehat{g_{n}\left(  t\right)
}\right\vert =\left\vert \widehat{Dg_{n}}\left(  t\right)  \right\vert  &
\leq\left\{
\begin{array}
[c]{ll}%
\frac{2}{\sqrt{2\pi}}, & n=0,\\
n\left\vert \widehat{g_{n-1}}\left(  t\right)  \right\vert +\frac{1}%
{\sqrt{2\pi}}, & n\geq1,
\end{array}
\right. \\
& \leq\frac{2}{\sqrt{2\pi}},
\end{align*}

so that%
\[
\left(  1+\left\vert t\right\vert \right)  \left\vert \widehat{g_{n}\left(
t\right)  }\right\vert \leq\left\vert \widehat{g_{n}\left(  t\right)
}\right\vert +\frac{2}{\sqrt{2\pi}}\leq\frac{1}{\sqrt{2\pi}}\left(  2+\frac
{1}{n+1}\right)  =\frac{c_{n}}{\sqrt{2\pi}}.
\]
\smallskip

\textbf{Part 5} True since $\widehat{g_{n}}\in C_{B}^{\infty}\left(
\mathbb{R}^{1}\right)  $.\medskip

\textbf{Part 6} From part 3,%
\[
Dg_{n}=\left\{
\begin{array}
[c]{ll}%
\delta-\delta\left(  \cdot-1\right)  , & n=0,\\
-ng_{n-1}+\delta, & n=1,2,3\ldots.
\end{array}
\right.
\]

Hence%
\begin{align*}
\left\Vert g_{0}\ast\widetilde{Dg_{0}}\right\Vert _{1}=\left\Vert g_{0}%
\ast\left(  \delta-\delta\left(  -\cdot-1\right)  \right)  \right\Vert _{1}  &
=\frac{1}{\sqrt{2\pi}}\left\Vert g_{0}-g_{0}\left(  \cdot+1\right)
\right\Vert _{1}\\
& =\frac{1}{\sqrt{2\pi}}\left(  \left\Vert g_{0}\right\Vert _{1}+\left\Vert
g_{0}\left(  \cdot+1\right)  \right\Vert _{1}\right) \\
& =\frac{1}{\sqrt{2\pi}}\left(  \left\Vert g_{0}\right\Vert _{1}+\left\Vert
g_{0}\left(  \cdot+1\right)  \right\Vert _{1}\right) \\
& =\frac{2}{\sqrt{2\pi}}\left\Vert g_{0}\right\Vert _{1}\\
& \leq\frac{2}{\sqrt{2\pi}},
\end{align*}

and if $n\geq1$,%
\begin{align*}
\left\Vert g_{n}\ast\widetilde{Dg_{n}}\right\Vert _{1}  & =\left\Vert
g_{n}\ast\left(  -n\widetilde{g}_{n-1}+\delta\right)  \right\Vert _{1}\\
& =\left\Vert -ng_{n}\ast\widetilde{g}_{n-1}+\frac{1}{\sqrt{2\pi}}%
g_{n}\right\Vert _{1}\\
& \leq n\left\Vert g_{n}\ast\widetilde{g}_{n-1}\right\Vert _{1}+\frac{1}%
{\sqrt{2\pi}}\left\Vert g_{n}\right\Vert _{1}\\
& \leq n\left\Vert g_{n}\right\Vert _{1}\left\Vert \widetilde{g}%
_{n-1}\right\Vert _{1}+\frac{1}{\sqrt{2\pi}}\left\Vert \widetilde{g}%
_{n}\right\Vert _{1}\\
& =n\left\Vert g_{n}\right\Vert _{1}\left\Vert g_{n-1}\right\Vert _{1}%
+\frac{1}{\sqrt{2\pi}}\left\Vert g_{n}\right\Vert _{1}\\
& \leq n\frac{1}{n\left(  n+1\right)  }+\frac{1}{\sqrt{2\pi}}\frac{1}{n+1}\\
& =\left(  1+\frac{1}{\sqrt{2\pi}}\right)  \frac{1}{n+1}.
\end{align*}
\medskip

\textbf{Part 7} From part 3, $Dg_{n}=-ng_{n-1}+\delta$. Hence
\begin{align*}
D^{2}g_{n}  & =-nDg_{n-1}+D\delta=-n\left(  -\left(  n-1\right)
g_{n-2}+\delta\right)  +D\delta=n\left(  n-1\right)  g_{n-2}-n\delta
+D\delta.\\
D^{3}g_{n}  & =n\left(  n-1\right)  Dg_{n-2}-nD\delta+D^{2}\delta=n\left(
n-1\right)  \left(  -\left(  n-2\right)  g_{n-3}+\delta\right)  -nD\delta
+D^{2}\delta=\\
& =-n\left(  n-1\right)  \left(  n-2\right)  g_{n-3}+n\left(  n-1\right)
\delta-nD\delta+D^{2}\delta,
\end{align*}

and in general%
\[
\left(  -D\right)  ^{k}g_{n}=\frac{n!}{\left(  n-k\right)  !}g_{n-k}-\frac
{n!}{\left(  n-k+1\right)  !}\delta+\frac{n!}{\left(  n-k+2\right)  !}%
D\delta-\ldots+\left(  -1\right)  ^{k}D^{k-1}\delta,\quad k\leq n.
\]

\[
\frac{\left(  -1\right)  ^{k}}{n!}D^{k}g_{n}=\frac{g_{n-k}}{\left(
n-k\right)  !}-\frac{\delta}{\left(  n-k+1\right)  !}+\frac{D\delta}{\left(
n-k+2\right)  !}-\ldots+\left(  -1\right)  ^{k}\frac{D^{k-1}\delta}{n!},\quad
k\leq n.
\]

\[
\frac{\left(  -1\right)  ^{n}}{n!}D^{n}g_{n}=g_{0}-\frac{\delta}{1!}%
+\frac{D\delta}{2!}-\ldots+\left(  -1\right)  ^{n}\frac{D^{n-1}\delta}{n!}.
\]

\begin{align*}
\frac{\left(  -1\right)  ^{n}}{n!}D^{n+1}g_{n}  & =Dg_{0}-\frac{D\delta}%
{1!}+\frac{D^{2}\delta}{2!}-\ldots+\left(  -1\right)  ^{n}\frac{D^{n}\delta
}{n!}\\
& =\delta-\delta\left(  \cdot-1\right)  -\frac{D\delta}{1!}+\frac{D^{2}\delta
}{2!}-\ldots+\left(  -1\right)  ^{n}\frac{D^{n}\delta}{n!}.
\end{align*}
\medskip

\textbf{Part 8} From \ref{a2.35},%
\[
e^{ix}=\sum_{k\leq n}\frac{\left(  ix\right)  ^{k}}{k!}+\frac{\sqrt{2\pi}}%
{n!}\left(  ix\right)  ^{n+1}\overset{\vee}{g_{n}}\left(  x\right)  .
\]

Thus%
\begin{align*}
\left(  ix\right)  ^{j}e^{ix}  & =\sum_{k\leq n}\frac{\left(  ix\right)
^{j+k}}{k!}+\frac{\sqrt{2\pi}}{n!}\left(  ix\right)  ^{n+j+1}\overset{\vee
}{g_{n}}\left(  x\right)  ;\\
\left(  ix\right)  ^{j}e^{ix}  & =\sum_{k\leq n}\frac{\left(  ix\right)
^{j+k}}{k!}+\frac{\sqrt{2\pi}}{n!}\left(  ix\right)  ^{n+j+1}\widehat{g_{n}%
}\left(  -x\right)  ;\\
\left(  -it\right)  ^{j}e^{-it}  & =\sum_{k\leq n}\frac{\left(  -it\right)
^{j+k}}{k!}+\frac{\sqrt{2\pi}}{n!}\left(  -it\right)  ^{n+j+1}\widehat{g_{n}%
}\left(  t\right)  ;\\
\left(  -it\right)  ^{j}e^{-it}  & =\sum_{k\leq n}\frac{\left(  -it\right)
^{j+k}}{k!}+\frac{\sqrt{2\pi}}{n!}\left(  -1\right)  ^{n+j+1}%
\widehat{D^{n+j+1}g_{n}}\left(  t\right)  ;\\
\widehat{D^{n+j+1}g_{n}}\left(  t\right)   & =\left(  -1\right)  ^{n+j+1}%
\frac{n!}{\sqrt{2\pi}}\left(  \left(  -it\right)  ^{j}e^{-it}-\sum_{k\leq
n}\frac{\left(  -it\right)  ^{j+k}}{k!}\right)  ;\\
\widehat{D^{n+j+1}g_{n}}\left(  t\right)   & =\left(  -1\right)  ^{n+j+1}%
\frac{n!}{\sqrt{2\pi}}\left(  -it\right)  ^{j}\left(  e^{-it}-\sum_{k\leq
n}\frac{\left(  -it\right)  ^{k}}{k!}\right)  .
\end{align*}

Clearly $\widehat{D^{n+j+1}g_{n}}\left(  t\right)  \in C_{BP}^{\infty}$. Since
$\left(  -it\right)  ^{j}\in C_{\emptyset,j}^{\infty}$ and $e^{-it}%
-\sum_{k\leq n}\frac{\left(  -it\right)  ^{k}}{k!}\in C_{\emptyset
,n+1}^{\infty}$ we have $\widehat{D^{n+j+1}g_{n}}\left(  t\right)  \in
C_{\emptyset,n+j+1}^{\infty}$.

A simple estimate is%
\[
\left\vert \widehat{D^{n+j+1}g_{n}}\left(  t\right)  \right\vert \leq\frac
{n!}{\sqrt{2\pi}}\left\vert t\right\vert ^{j}\left(  1+\sum_{k\leq n}%
\frac{\left\vert t\right\vert ^{k}}{k!}\right)  .
\]

\end{proof}

Thus if $f\in S^{\prime}$ and $\xi$ is the \textbf{action} variable when
$a.\neq\mathbf{0}$ we can write \ref{a114} as%
\begin{equation}
\left(  e^{ia\xi}-\sum_{k\leq n}\frac{\left(  ia\xi\right)  ^{k}}{k!}\right)
\widehat{f}=\frac{\sqrt{2\pi}}{n!}\left(  ia\xi\right)  ^{n+1}\overline
{\widehat{g_{n}}}\left(  a\xi\right)  \widehat{f},\label{a2.11}%
\end{equation}

or on using the first of the multi-index identities \ref{1.57},
\begin{equation}
\frac{\left(  a\xi\right)  ^{k}}{k!}=\sum_{\left\vert \beta\right\vert
=k}\frac{a^{\beta}\xi^{\beta}}{\beta!},\quad\frac{d^{k}}{k!}=\sum_{\left\vert
\beta\right\vert =k}\frac{1}{\beta!},\label{a224}%
\end{equation}

we get%
\[
\left(  f\left(  \cdot+a\right)  -\sum_{k\leq n}\frac{\left(  aD\right)  ^{k}%
}{k!}f\right)  ^{\wedge}=\left(  f\left(  \cdot+a\right)  -\sum_{\left\vert
\beta\right\vert \leq n}\frac{a^{\beta}}{\beta!}D^{\beta}f\right)  ^{\wedge
}=\frac{\sqrt{2\pi}}{n!}\left(  ia\xi\right)  ^{n+1}\overline{\widehat{g_{n}}%
}\left(  a\xi\right)  \widehat{f},
\]

which implies the tempered distribution Taylor series expansion
\begin{equation}
f\left(  \cdot+a\right)  -\sum_{\left\vert \beta\right\vert \leq n}%
\frac{a^{\beta}}{\beta!}D^{\beta}f=f\left(  \cdot+a\right)  -\sum_{k\leq
n}\frac{\left(  aD\right)  ^{k}}{k!}f=\left(  \mathcal{R}_{n+1}f\right)
\left(  \cdot,a\right)  ,\quad f\in S^{\prime},\label{a1.55}%
\end{equation}

where%
\begin{equation}
\left(  \mathcal{R}_{n+1}f\right)  \left(  \cdot,a\right)  :=\frac{\sqrt{2\pi
}}{n!}\left(  \left(  ia\xi\right)  ^{n+1}\overline{\widehat{g_{n}}}\left(
a\xi\right)  \widehat{f}\right)  ^{\vee},\quad f\in S^{\prime},\text{
}n=0,1,2,\ldots\label{a50.5}%
\end{equation}

Suppose $f\in C_{B}^{\left(  n\right)  }$ and $\left(  aD\right)  ^{n+1}f\in
L^{1}$. For clarity set $u=\left(  aD\right)  ^{n+1}f$ where $D=\left(
D_{k}\right)  _{k=1}^{d}$, so that
\begin{align*}
\left(  ia\xi\right)  ^{n+1}\overline{\widehat{g_{n}}}\left(  a\xi\right)
\widehat{f}\left(  \xi\right)  =\overline{\widehat{g_{n}}}\left(  a\xi\right)
\widehat{u}\left(  \xi\right)   &  =\frac{1}{\sqrt{2\pi}}\int_{0}^{1}%
e^{isa\xi}g_{n}\left(  s\right)  \widehat{u}\left(  \xi\right)  ds\\
&  =\frac{1}{\sqrt{2\pi}}\int_{0}^{1}g_{n}\left(  s\right)  F_{x}\left[
u\left(  x+sa\right)  \right]  \left(  \xi\right)  ds\\
&  =\frac{1}{\sqrt{2\pi}}F_{x}\left[  \int_{0}^{1}g_{n}\left(  s\right)
u\left(  x+sa\right)  ds\right]  \left(  \xi\right)  ,
\end{align*}

because $u=\left(  aD\right)  ^{n+1}f\in L^{1}$ implies the Fourier transform
can be transposed allowing an application of Fubini's theorem. Thus%
\begin{equation}
\left(  \mathcal{R}_{n+1}f\right)  \left(  x,a\right)  =\frac{1}{n!}\int%
_{0}^{1}g_{n}\left(  s\right)  \left(  \left(  aD\right)  ^{n+1}f\right)
\left(  x+sa\right)  ds,\quad f\in C_{B}^{\left(  n\right)  },\text{ }\left(
aD\right)  ^{n+1}f\in L^{1}.\label{a51.0}%
\end{equation}

If, in addition, $\left(  \widehat{a}D\right)  ^{n+1}f\in L^{\infty}\left(
\left[  x,x+a\right]  \right)  $ where $\widehat{a}=a/\left\vert a\right\vert
$ then
\begin{align*}
\left\vert \left(  \mathcal{R}_{n+1}f\right)  \left(  x,a\right)  \right\vert
& \leq\frac{1}{n!}\left\Vert \left(  aD\right)  ^{n+1}f\right\Vert _{\infty
}\int_{0}^{1}g_{n}\left(  s\right)  ds\\
& =\left(  \frac{1}{\left(  n+1\right)  !}\left\Vert \left(  \widehat{a}%
D\right)  ^{n+1}f\right\Vert _{\infty,\left[  x,x+a\right]  }\right)
\left\vert a\right\vert ^{n+1}.
\end{align*}

Thus we have proved the global result:

\begin{lemma}
\label{Lem_Taylor_estim_C_L1}Suppose $f\in C_{BP}^{\left(  n\right)  }\left(
\mathbb{R}^{d}\right)  $ and in the distribution sense $\left(  aD\right)
^{n+1}f\in L^{1}\left(  \mathbb{R}^{d}\right)  $ when $\left\vert a\right\vert
=1$.

Then%
\begin{equation}
f\left(  \cdot+a\right)  -\sum_{k\leq n}\frac{\left(  aD\right)  ^{k}}%
{k!}f=\frac{1}{n!}\int_{0}^{1}g_{n}\left(  s\right)  \left(  \left(
aD\right)  ^{n+1}f\right)  \left(  x+sa\right)  ds,\quad a,x\in\mathbb{R}%
^{d}.\label{a111}%
\end{equation}

\end{lemma}

We can easily weaken the condition $\left(  bD\right)  ^{n+1}f\in L^{1}\left(
\mathbb{R}^{d}\right)  $ to $\left(  bD\right)  ^{n+1}f\in L_{loc}^{1}\left(
\mathbb{R}^{d}\right)  $ and obtain:

\begin{lemma}
\label{Lem_Taylor_estim_C_L1loc}Suppose $f\in C_{BP}^{\left(  n\right)
}\left(  \mathbb{R}^{d}\right)  $ and in the distribution sense $\left(
aD\right)  ^{n+1}f\in L_{loc}^{1}\left(  \mathbb{R}^{d}\right)  $ when
$\left\vert a\right\vert =1$. Then for all $a$, $\left(  \mathcal{R}%
_{n+1}f\right)  \left(  \cdot,a\right)  \in L^{1}$ uniformly and%
\begin{equation}
f\left(  \cdot+a\right)  -\sum_{k\leq n}\frac{\left(  aD\right)  ^{k}}%
{k!}f=\frac{1}{n!}\int_{0}^{1}g_{n}\left(  s\right)  \left(  \left(
aD\right)  ^{n+1}f\right)  \left(  \cdot+sa\right)  ds\text{ }a.e.\text{
}on\text{ }\mathbb{R}^{d}.\label{a113}%
\end{equation}

\end{lemma}

\begin{proof}
Choose arbitrary $x$ and $a$. Lemma \ref{Lem_func_eq_1_nbhd_set} then tells us
that we can choose $\phi\in C_{0}^{\infty}$ such $\phi=1$ on a neighborhood of
the closed line segment $\left[  x,x+a\right]  $. Then $\phi f\in
C_{BP}^{\left(  n\right)  }$ and in the distribution sense $\left(  bD\right)
^{n+1}\left(  \phi f\right)  \in L^{1}$ when $\left\vert b\right\vert =1$.
Thus \ref{a111} holds for $\phi f$ i.e.%
\[
\left(  \phi f\right)  \left(  x+a\right)  -\sum_{k\leq n}\frac{\left(
aD\right)  ^{k}}{k!}\left(  \phi f\right)  \left(  x\right)  =\frac{1}{n!}%
\int_{0}^{1}g_{n}\left(  s\right)  \left(  \left(  aD\right)  ^{n+1}\left(
\phi f\right)  \right)  \left(  x+sa\right)  ds\text{\ }a.e.
\]
consequently the remainder formula \ref{a113} holds for each $x$ and $a$.
\end{proof}

A further easy local generalization is:

\begin{lemma}
\label{Lem_Taylor_estim_loc_C_L1loc}Suppose $\Omega\subseteq\mathbb{R}^{d}$ is
an open set, $f\in C_{BP}^{\left(  n\right)  }\left(  \mathbb{\Omega}\right)
$ and in the distribution sense $\left(  bD\right)  ^{n+1}f\in L_{loc}%
^{1}\left(  \mathbb{\Omega}\right)  $ when $\left\vert b\right\vert =1$.
Further, suppose that $x,x+a\in\Omega$ and $\left[  x,x+a\right]
\subset\Omega$. Then%
\begin{equation}
f\left(  x+a\right)  -\sum_{k\leq n}\frac{\left(  aD\right)  ^{k}}{k!}f\left(
x\right)  =\frac{1}{n!}\int_{0}^{1}g_{n}\left(  s\right)  \left(  \left(
aD\right)  ^{n+1}f\right)  \left(  x+sa\right)  ds.\label{a115}%
\end{equation}

\end{lemma}

\begin{proof}
Lemma \ref{Lem_func_eq_1_nbhd_set} tells us that we can choose $\phi\in
C_{0}^{\infty}\left(  \Omega\right)  $ such that $\phi=1$ on a neighborhood of
the closed line segment $\left[  x,x+a\right]  $. Then $\phi f\in
C_{BP}^{\left(  n\right)  }$ and in the distribution sense $\left(  bD\right)
^{n+1}\left(  \phi f\right)  \in L^{1}$ when $\left\vert b\right\vert =1$.
Thus \ref{a111} holds for $\phi f$ and%
\[
\left(  \phi f\right)  \left(  \cdot+a\right)  -\sum_{k\leq n}\frac{\left(
aD\right)  ^{k}}{k!}\left(  \phi f\right)  =\frac{1}{n!}\int_{0}^{1}%
g_{n}\left(  s\right)  \left(  \left(  aD\right)  ^{n+1}\left(  \phi f\right)
\right)  \left(  \cdot+sa\right)  ds,\quad on\text{ }\mathbb{R}^{d}.
\]

consequently \ref{a115} holds for the specific $x$ and $a$.
\end{proof}

Yet another global generalization is:

\begin{lemma}
\label{Lem_Taylor_estim_S_L1}Suppose $f\in S^{\prime}$ and $\left(  aD\right)
^{n+1}f\in L^{1}$ when $\left\vert a\right\vert =1$. Then $\left(
\mathcal{R}_{n+1}f\right)  \left(  \cdot,a\right)  \in L^{1}$ for all $a$, and
in the sense of distributions\
\[
f\left(  \cdot+a\right)  -\sum_{k\leq n}\frac{\left(  aD\right)  ^{k}}%
{k!}f=\frac{1}{n!}\int_{0}^{1}g_{n}\left(  s\right)  \left(  \left(
aD\right)  ^{n+1}f\right)  \left(  \cdot+sa\right)  ds,\quad a\in
\mathbb{R}^{d}.
\]

\end{lemma}

\begin{proof}
The assumption that $f\in S^{\prime}$ yields \ref{a1.55} and \ref{a50.5}. By
inspection we then see that the subsequent calculations still yield equation
\ref{a51.0} for the remainder without making the assumption $f\in
C_{BP}^{\left(  n\right)  }$. Further%
\[
\left\vert \int_{0}^{1}g_{n}\left(  s\right)  \left(  \left(  aD\right)
^{n+1}f\right)  \left(  x+sa\right)  ds\right\vert \leq\left\Vert \left(
aD\right)  ^{n+1}f\right\Vert _{1},
\]

and the remainder is $L^{1}$ for all $a$.
\end{proof}

which in turn leads to the local result:

\begin{lemma}
\label{Lem_Taylor_estim_L1loc_L1}?? \textbf{CHECK AGAIN} - esp. last claim! ??
Suppose $\Omega\subseteq\mathbb{R}^{d}$ is open and that when $\left\vert
b\right\vert =1$, $\left(  bD\right)  ^{k}f\in L_{loc}^{1}\left(
\mathbb{\Omega}\right)  $ when $k\leq n+1$. Further, suppose that
$x,x+a\in\Omega$ and $\left[  x,x+a\right]  \subset\Omega$.

Then in the sense of distributions\
\begin{equation}
f\left(  x+a\right)  -\sum_{k\leq n}\frac{\left(  aD\right)  ^{k}}{k!}f\left(
x\right)  \overset{a.e.}{=}\frac{1}{n!}\int_{0}^{1}g_{n}\left(  s\right)
\left(  \left(  aD\right)  ^{n+1}f\right)  \left(  x+sa\right)
ds.\label{a116}%
\end{equation}

Also, $\left(  \mathcal{R}_{n+1}f\right)  \left(  \cdot,a\right)  \in
L_{loc}^{1}\left(  B_{r}\left(  x\right)  \right)  $ where $r=\min\left\{
\operatorname*{dist}\left(  x,\Omega^{c}\right)  ,\operatorname*{dist}\left(
x+a,\Omega^{c}\right)  \right\}  $.
\end{lemma}

\begin{proof}
Lemma \ref{Lem_func_eq_1_nbhd_set} tells us that we can choose $\phi\in
C_{0}^{\infty}$ such $\phi=1$ on a neighborhood of the closed line segment
$\left[  x,x+a\right]  $. Then $\phi f\in\mathcal{E}^{\prime}\subset
S^{\prime}$ and, since $\left(  bD\right)  ^{k}f\in L_{loc}^{1}\left(
\mathbb{\Omega}\right)  $ when $k\leq n+1$, in the distribution sense $\left(
bD\right)  ^{n+1}\left(  \phi f\right)  \in L^{1}$ when $\left\vert
b\right\vert =1$. Then Lemma \ref{Lem_Taylor_estim_S_L1} implies that in the
distribution sense%
\[
\left(  \phi f\right)  \left(  \cdot+b\right)  -\sum_{k\leq n}\frac{\left(
aD\right)  ^{k}}{k!}\left(  \phi f\right)  =\frac{1}{n!}\int_{0}^{1}%
g_{n}\left(  s\right)  \left(  \left(  bD\right)  ^{n+1}\left(  \phi f\right)
\right)  \left(  \cdot+sb\right)  ds,\quad b\in\mathbb{R}^{d},
\]

which means that%
\[
f\left(  x+a\right)  -\sum_{k\leq n}\frac{\left(  aD\right)  ^{k}}{k!}%
f=\frac{1}{n!}\int_{0}^{1}g_{n}\left(  s\right)  \left(  \left(  aD\right)
^{n+1}f\right)  \left(  x+sa\right)  ds,
\]

as claimed.

The last claim ?? \textbf{TRUE}? ?? follows from the form of the left side of
\ref{a116} and the fact that $\left(  bD\right)  ^{k}f\in L_{loc}^{1}\left(
\mathbb{\Omega}\right)  $ when $k\leq n+1$ and $\left\vert b\right\vert =1$.
\end{proof}

The general remainder result \ref{a50.5} will now be specialized to the
functions in $X_{w}^{0}$ - see the remarks to this theorem:

\begin{theorem}
\label{Thm_Tayor_rem_estim}\textbf{Basic remainder estimates} Suppose the
weight function $w$ has property W02 for some $\kappa$. Set $m=\left\lfloor
\kappa\right\rfloor $ and $\widehat{a}:=a/\left\vert a\right\vert $. Then for
all $f\in X_{w}^{0}$ we have:

\begin{enumerate}
\item the weight function estimates: for all $x,a\in\mathbb{R}^{d}$,%
\begin{align}
\left\vert \left(  \mathcal{R}_{m+1}f\right)  \left(  x,a\right)  \right\vert
& \leq\frac{c_{m}}{m!}\left(  \int\frac{\left\vert \widehat{a}\xi\right\vert
^{2\kappa}}{w\left(  \xi\right)  }d\xi\right)  ^{1/2}\left\Vert f\right\Vert
_{w,0}\left\vert a\right\vert ^{\kappa}\label{a2.06}\\
& \leq\frac{c_{m}}{m!}\left(  \int\frac{\left\vert \cdot\right\vert ^{2\kappa
}}{w}\right)  ^{1/2}\left\Vert f\right\Vert _{w,0}\left\vert a\right\vert
^{\kappa},\nonumber
\end{align}

where $c_{m}$ is given by \ref{a2.08},

\item and the following basis function estimate: for all $x,a\in\mathbb{R}%
^{d}$,
\begin{align}
\left\vert \left(  \mathcal{R}_{m+1}f\right)  \left(  x,a\right)  \right\vert
&  \leq\left(  \int\left\vert e^{ia\xi}-\sum_{k=0}^{m}\frac{\left(
ia\xi\right)  ^{k}}{k!}\right\vert ^{2}\widehat{G}\left(  \xi\right)
d\xi\right)  ^{1/2}\left\Vert f\right\Vert _{w,0}\label{a2.15}\\
&  =\left(  2\pi\right)  ^{\frac{d}{4}}\left(  2G\left(  0\right)  -\sum
_{k=0}^{m}\frac{1}{k!}\left(  \left(  aD\right)  ^{k}G\right)  \left(
-a\right)  -\sum_{k=0}^{m}\frac{1}{k!}\left(  -aD\right)  ^{k}G\left(
a\right)  +\right. \nonumber\\
&  \qquad\qquad\qquad\left.  +%
{\displaystyle\sum\limits_{\substack{j,k\leq m \\j+k>m}}}
\frac{\left(  -1\right)  ^{j}}{j!k!}\left(  \left(  aD\right)  ^{j+k}G\right)
\left(  0\right)  \right)  ^{\frac{1}{2}}\left\Vert f\right\Vert
_{w,0}.\label{a2.31}%
\end{align}

\item When $m=0$ the basis function estimate of part 2 becomes%
\begin{equation}
\left\vert \left(  \mathcal{R}_{1}f\right)  \left(  x,a\right)  \right\vert
\leq\left(  2\pi\right)  ^{\frac{d}{4}}\sqrt{2}\left(  G\left(  0\right)
-\operatorname{Re}G\left(  a\right)  \right)  ^{\frac{1}{2}}\left\Vert
f\right\Vert _{w,0},\quad x,a\in\mathbb{R}^{d}.\label{a54.9}%
\end{equation}

\item For all $m$: suppose $\left(  aD\right)  ^{2m+1}G\in L^{1}$ when
$\left\vert a\right\vert =1$. Then%
\begin{equation}
\left\vert \left(  \mathcal{R}_{m+1}f\right)  \left(  x,a\right)  \right\vert
\leq\left(  2\pi\right)  ^{\frac{d}{4}}\sqrt{\left(  \mathcal{R}%
_{2m+1}G\right)  \left(  a\right)  }\left\Vert f\right\Vert _{w,0},\quad
x,a\in\mathbb{R}^{d},\label{a2.30}%
\end{equation}

where%
\begin{align}
&  \left(  \mathcal{R}_{2m+1}G\right)  \left(  a\right) \nonumber\\
&  =\sum_{k=m+1}^{2m}\frac{1}{k!}\left\{  \left(  \mathcal{R}_{2m+1-k}\left(
\left(  -aD\right)  ^{k}G\right)  \right)  \left(  0,a\right)  +\left(
\mathcal{R}_{2m+1-k}\left(  \left(  aD\right)  ^{k}G\right)  \right)  \left(
0,-a\right)  \right\}  +\nonumber\\
&  \qquad\qquad+\left(  \mathcal{R}_{2m+1}G\right)  \left(  a,-a\right)
+\left(  \mathcal{R}_{2m+1}G\right)  \left(  -a,a\right) \label{a2.32}\\
&  =2\int_{0}^{1}\left(  \sum\limits_{j=0}^{m-1}\frac{\left(  -1\right)
^{j}\left(  1-s\right)  ^{j}}{j!\left(  2m-j\right)  !}-\frac{s^{2m}}{\left(
2m\right)  !}\right)  \left(  \left(  aD\right)  ^{2m+1}\operatorname{Re}%
G\right)  \left(  sa\right)  ds.\label{a2.33}%
\end{align}

\item We have%
\begin{equation}
\int_{0}^{1}\left\vert \sum\limits_{j=0}^{m-1}\frac{\left(  -1\right)
^{j}\left(  1-s\right)  ^{j}}{j!\left(  2m-j\right)  !}-\frac{s^{2m}}{\left(
2m\right)  !}\right\vert ds=\frac{1}{\left(  2m+1\right)  \left(  m!\right)
^{2}}=\frac{\binom{2m}{m}}{\left(  2m+1\right)  !}.\label{a3.01}%
\end{equation}

\end{enumerate}
\end{theorem}

\begin{proof}
\fbox{\textbf{Part 1}} Since $f\in X_{w}^{0}$ we have $\sqrt{w}\widehat{f}\in
L^{2}$ and $\left\Vert f\right\Vert _{w,0}=\left\Vert \sqrt{w}\widehat{f}%
\right\Vert _{2}$ and so from \ref{a50.5} we proceed formally forwards and
then justify each step backwards to get:%
\begin{align}
\left\vert \left(  \mathcal{R}_{m+1}f\right)  \left(  x,a\right)  \right\vert
& \leq\frac{\sqrt{2\pi}}{m!}\left\Vert \left(  \left(  ia\xi\right)
^{m+1}\overline{\widehat{g_{m}}}\left(  a\xi\right)  \widehat{f}\right)
^{\vee}\right\Vert _{\infty}\nonumber\\
& =\frac{\sqrt{2\pi}}{m!}\left\Vert \left(  \frac{\left(  a\xi\right)  ^{m+1}%
}{\sqrt{w}}\overline{\widehat{g_{m}}}\left(  a\xi\right)  \sqrt{w}%
\widehat{f}\right)  ^{\vee}\right\Vert _{\infty}\nonumber\\
& =\frac{\sqrt{2\pi}}{m!}\left\Vert \left(  \frac{\left(  a\xi\right)  ^{m+1}%
}{\sqrt{w}}\overline{\widehat{g_{m}}}\left(  a\xi\right)  \right)  ^{\vee}%
\ast\left(  \sqrt{w}\widehat{f}\right)  ^{\vee}\right\Vert _{\infty
}\label{a2.12}\\
& \leq\frac{\sqrt{2\pi}}{m!}\left\Vert \left(  \frac{\left(  a\xi\right)
^{m+1}}{\sqrt{w}}\overline{\widehat{g_{m}}}\left(  a\xi\right)  \right)
^{\vee}\right\Vert _{2}\left\Vert \left(  \sqrt{w}\widehat{f}\right)  ^{\vee
}\right\Vert _{2}\label{a2.09}\\
& =\frac{\sqrt{2\pi}}{m!}\left\Vert \frac{\left(  a\xi\right)  ^{m+1}}%
{\sqrt{w}}\overline{\widehat{g_{m}}}\left(  a\xi\right)  \right\Vert
_{2}\left\Vert \sqrt{w}\widehat{f}\right\Vert _{2}\label{a2.19}\\
& =\frac{\sqrt{2\pi}}{m!}\left(  \int\frac{\left(  a\xi\right)  ^{2\left(
m+1\right)  }}{w\left(  \xi\right)  }\left\vert \widehat{g_{m}}\left(
a\xi\right)  \right\vert ^{2}d\xi\right)  ^{1/2}\left\Vert f\right\Vert
_{w,0}\label{a2.14}\\
& \leq\frac{\sqrt{2\pi}}{m!}\left(  \int\frac{\left(  a\xi\right)  ^{2\left(
m+1\right)  }}{w\left(  \xi\right)  }\left(  \frac{c_{m}/\sqrt{2\pi}%
}{1+\left\vert a\xi\right\vert }\right)  ^{2}d\xi\right)  ^{1/2}\left\Vert
f\right\Vert _{w,0}\label{a2.10}\\
& =\frac{c_{m}}{m!}\left(  \int\frac{\left(  a\xi\right)  ^{2\left(
m+1\right)  }}{w\left(  \xi\right)  }\left(  \frac{1}{1+\left\vert
a\xi\right\vert }\right)  ^{2}d\xi\right)  ^{1/2}\left\Vert f\right\Vert
_{w,0}\nonumber\\
& =\frac{c_{m}}{m!}\left(  \int\frac{\left\vert a\xi\right\vert ^{2\kappa}%
}{w\left(  \xi\right)  }\left(  \frac{\left\vert a\xi\right\vert ^{m+1-\kappa
}}{1+\left\vert a\xi\right\vert }\right)  ^{2}d\xi\right)  ^{1/2}\left\Vert
f\right\Vert _{w,0}\nonumber\\
& \leq\frac{c_{m}}{m!}\left(  \int\frac{\left\vert a\xi\right\vert ^{2\kappa}%
}{w\left(  \xi\right)  }\left(  \frac{\left(  1+\left\vert a\xi\right\vert
\right)  ^{m+1-\kappa}}{1+\left\vert a\xi\right\vert }\right)  ^{2}%
d\xi\right)  ^{1/2}\left\Vert f\right\Vert _{w,0}\nonumber\\
& \leq\frac{c_{m}}{m!}\left(  \int\frac{\left\vert a\xi\right\vert ^{2\kappa}%
}{w\left(  \xi\right)  }d\xi\right)  ^{1/2}\left\Vert f\right\Vert
_{w,0}\nonumber\\
& =\frac{c_{m}}{m!}\left(  \int\frac{\left\vert \widehat{a}\xi\right\vert
^{2\kappa}}{w\left(  \xi\right)  }d\xi\right)  ^{1/2}\left\Vert f\right\Vert
_{w,0}\left\vert a\right\vert ^{\kappa}\nonumber\\
& \leq\frac{c_{m}}{m!}\left(  \int\frac{\left\vert \cdot\right\vert ^{2\kappa
}}{w}\right)  ^{1/2}\left\Vert f\right\Vert _{w,0}\left\vert a\right\vert
^{\kappa}\nonumber\\
& <\infty.\nonumber
\end{align}

The last integral exists since $w$ has property W02 for $\kappa$. Step
$\ref{a2.12}\Leftrightarrow\ref{a2.09}$ is Young's estimate \ref{1.056} for
convolutions when $p=q=2$ and $r=\infty$. Step $\ref{a2.09}\Leftrightarrow
\ref{a2.19}$ is Plancherel's theorem for $L^{2}$ functions. Step
$\ref{a2.14}\Leftrightarrow\ref{a2.10}$ uses the estimates \ref{a2.07}%
.\medskip

\fbox{\textbf{Part 2}} Our starting point is inequality \ref{a2.14} of part 1
of this theorem:%
\[
\left\vert \left(  \mathcal{R}_{m+1}f\right)  \left(  x,a\right)  \right\vert
\leq\frac{\sqrt{2\pi}}{m!}\left(  \int\frac{\left(  a\xi\right)  ^{2\left(
m+1\right)  }}{w\left(  \xi\right)  }\left\vert \widehat{g_{m}}\left(
a\xi\right)  \right\vert ^{2}d\xi\right)  ^{1/2}\left\Vert f\right\Vert
_{w,0}.
\]

The Taylor series expansion \ref{a114} can be rewritten%
\[
\left(  ia\xi\right)  ^{m+1}\overline{\widehat{g_{m}}}\left(  a\xi\right)
=\frac{m!}{\sqrt{2\pi}}\left(  e^{ia\xi}-\sum_{k=0}^{m}\frac{\left(
ia\xi\right)  ^{k}}{k!}\right)  ,
\]

so that%
\begin{align*}
\left\vert \left(  \mathcal{R}_{m+1}f\right)  \left(  x,a\right)  \right\vert
& \leq\frac{\sqrt{2\pi}}{m!}\left(  \int\frac{\left(  a\xi\right)  ^{2\left(
m+1\right)  }}{w\left(  \xi\right)  }\left\vert \widehat{g_{m}}\left(
a\xi\right)  \right\vert ^{2}d\xi\right)  ^{1/2}\left\Vert f\right\Vert
_{w,0}\\
& =\frac{\sqrt{2\pi}}{m!}\left(  \int\left\vert \left(  a\xi\right)
^{m+1}\widehat{g_{m}}\left(  a\xi\right)  \right\vert ^{2}\widehat{G}\left(
\xi\right)  d\xi\right)  ^{1/2}\left\Vert f\right\Vert _{w,0}\\
& =\frac{\sqrt{2\pi}}{m!}\left(  \int\left\vert \frac{m!}{\sqrt{2\pi}}\left(
e^{ia\xi}-\sum_{k=0}^{m}\frac{\left(  ia\xi\right)  ^{k}}{k!}\right)
\right\vert ^{2}\widehat{G}\left(  \xi\right)  d\xi\right)  ^{1/2}\left\Vert
f\right\Vert _{w,0}\\
& =\left(  \int\left\vert e^{ia\xi}-\sum_{k=0}^{m}\frac{\left(  ia\xi\right)
^{k}}{k!}\right\vert ^{2}\widehat{G}\left(  \xi\right)  d\xi\right)
^{1/2}\left\Vert f\right\Vert _{w,0},
\end{align*}

which proves \ref{a2.15}. But%
\begin{align*}
\left\vert e^{ia\xi}-\sum_{k=0}^{m}\frac{\left(  ia\xi\right)  ^{k}}%
{k!}\right\vert ^{2} &  =\left(  \overline{e^{ia\xi}-\sum_{l=0}^{m}%
\frac{\left(  ia\xi\right)  ^{l}}{l!}}\right)  \left(  e^{ia\xi}-\sum
_{k=0}^{m}\frac{\left(  ia\xi\right)  ^{k}}{k!}\right) \\
&  =\left(  e^{-ia\xi}-\sum_{l=0}^{m}\frac{\left(  -ia\xi\right)  ^{l}}%
{l!}\right)  \left(  e^{ia\xi}-\sum_{k=0}^{m}\frac{\left(  ia\xi\right)  ^{k}%
}{k!}\right) \\
&  =1-e^{-ia\xi}\sum_{k=0}^{m}\frac{\left(  ia\xi\right)  ^{k}}{k!}-e^{ia\xi
}\sum_{l=0}^{m}\frac{\left(  -ia\xi\right)  ^{l}}{l!}+\sum_{k,l=0}^{m}%
\frac{\left(  -ia\xi\right)  ^{l}}{l!}\frac{\left(  ia\xi\right)  ^{k}}{k!}\\
&  =1-e^{-ia\xi}\sum_{k=0}^{m}\frac{\left(  ia\xi\right)  ^{k}}{k!}-e^{ia\xi
}\sum_{l=0}^{m}\frac{\left(  -1\right)  ^{l}\left(  ia\xi\right)  ^{l}}%
{l!}+\sum_{l,k=0}^{m}\frac{\left(  -1\right)  ^{l}\left(  ia\xi\right)
^{l+k}}{k!l!},
\end{align*}

and since we know that $e^{ia\xi}\widehat{G}\left(  \xi\right)  =\left(
G\left(  \cdot+a\right)  \right)  ^{\wedge}\left(  \xi\right)  $ and $\left(
ia\xi\right)  ^{k}\widehat{G}\left(  \xi\right)  =\left(  \left(  aD\right)
^{k}G\right)  ^{\wedge}\left(  \xi\right)  $ it follows that%
\begin{align*}
&  \left\vert e^{ia\xi}-\sum_{k=0}^{m}\frac{\left(  ia\xi\right)  ^{k}}%
{k!}\right\vert ^{2}\widehat{G}\\
&  =\left(  1-\sum_{k=0}^{m}\frac{\left(  ia\xi\right)  ^{k}}{k!}e^{-ia\xi
}-\sum_{l=0}^{m}\frac{\left(  -ia\xi\right)  ^{l}}{l!}e^{ia\xi}+\sum
_{l,k=0}^{m}\frac{\left(  -1\right)  ^{l}}{k!l!}\left(  ia\xi\right)
^{l+k}\right)  \widehat{G}\\
&  =\widehat{G}-\sum_{k=0}^{m}\frac{\left(  ia\xi\right)  ^{k}e^{-ia\xi
}\widehat{G}}{k!}-\sum_{l=0}^{m}\frac{\left(  -ia\xi\right)  ^{l}e^{ia\xi
}\widehat{G}}{l!}+\sum_{l,k=0}^{m}\frac{\left(  -1\right)  ^{l}}{k!l!}\left(
ia\xi\right)  ^{l+k}\widehat{G}\\
&  =\widehat{G}-\sum_{k=0}^{m}\frac{\left(  ia\xi\right)  ^{k}%
\widehat{G\left(  \cdot-a\right)  }}{k!}-\sum_{l=0}^{m}\frac{\left(
-ia\xi\right)  ^{l}\widehat{G\left(  \cdot+a\right)  }}{l!}+\sum_{l,k=0}%
^{m}\frac{\left(  -1\right)  ^{l}}{k!l!}\left(  ia\xi\right)  ^{l+k}%
\widehat{G}\\
&  =\widehat{G}-\sum_{k=0}^{m}\frac{\widehat{\left(  aD\right)  ^{k}G\left(
\cdot-a\right)  }}{k!}-\sum_{l=0}^{m}\frac{\widehat{\left(  -aD\right)
^{l}G\left(  \cdot+a\right)  }}{l!}+\sum_{l,k=0}^{m}\frac{\left(  -1\right)
^{l}}{k!l!}\widehat{\left(  aD\right)  ^{l+k}G}\\
&  =\left(  G-\sum_{k=0}^{m}\frac{\left(  \left(  aD\right)  ^{k}G\right)
\left(  \cdot-a\right)  }{k!}-\sum_{l=0}^{m}\frac{\left(  \left(  -aD\right)
^{l}G\right)  \left(  \cdot+a\right)  }{l!}+\sum_{l,k=0}^{m}\frac{\left(
-1\right)  ^{l}}{k!l!}\left(  aD\right)  ^{l+k}G\right)  ^{\wedge},
\end{align*}

and hence%
\begin{align}
&  \left(  2\pi\right)  ^{-\frac{d}{2}}\int\left\vert e^{ia\xi}-\sum_{k=0}%
^{m}\frac{\left(  ia\xi\right)  ^{k}}{k!}\right\vert ^{2}\widehat{G}\left(
\xi\right)  d\xi\nonumber\\
&  =G\left(  0\right)  -\sum_{k=0}^{m}\frac{\left(  \left(  aD\right)
^{k}G\right)  \left(  -a\right)  }{k!}-\sum_{k=0}^{m}\frac{\left(  \left(
-aD\right)  ^{k}G\right)  \left(  a\right)  }{k!}+\sum_{l,k=0}^{m}%
\frac{\left(  -1\right)  ^{l}}{k!l!}\left(  \left(  aD\right)  ^{l+k}G\right)
\left(  0\right)  ,\label{a2.20}%
\end{align}

which proves \ref{a2.31} and this part.\medskip

\fbox{\textbf{Part 3}} Substitute $m=0$ in the RHS of inequality \ref{a2.31}:%
\begin{align*}
\left\vert \left(  \mathcal{R}_{1}f\right)  \left(  x,a\right)  \right\vert  &
\leq2G\left(  0\right)  -\overline{G\left(  a\right)  }-G\left(  a\right) \\
& =2\left(  G\left(  0\right)  -\operatorname{Re}G\left(  a\right)  \right)  .
\end{align*}
\medskip

\fbox{\textbf{Part 4}} Regarding the last term in \ref{a2.20}:\smallskip%
\begin{align*}
\sum_{l,k=0}^{m}\frac{\left(  -1\right)  ^{l}}{k!l!}\left(  \left(  aD\right)
^{l+k}G\right)  \left(  0\right)   & =%
{\displaystyle\sum\limits_{l+k\leq m}}
\frac{\left(  -1\right)  ^{l}}{k!l!}\left(  \left(  aD\right)  ^{l+k}G\right)
\left(  0\right)  +%
{\displaystyle\sum\limits_{\substack{l+k>m \\l,k\leq m}}}
\frac{\left(  -1\right)  ^{l}}{k!l!}\left(  \left(  aD\right)  ^{l+k}G\right)
\left(  0\right) \\
& =%
{\displaystyle\sum\limits_{q=0}^{m}}
{\displaystyle\sum\limits_{l+k=q}}
\frac{\left(  -1\right)  ^{l}}{k!l!}\left(  \left(  aD\right)  ^{l+k}G\right)
\left(  0\right)  +%
{\displaystyle\sum\limits_{\substack{l+k>m \\l,k\leq m}}}
\frac{\left(  -1\right)  ^{l}}{k!l!}\left(  \left(  aD\right)  ^{l+k}G\right)
\left(  0\right) \\
& =%
{\displaystyle\sum\limits_{q=0}^{m}}
\left(
{\displaystyle\sum\limits_{l+k=q}}
\frac{\left(  -1\right)  ^{l}}{k!l!}\right)  \left(  \left(  aD\right)
^{q}G\right)  \left(  0\right)  +%
{\displaystyle\sum\limits_{\substack{l+k>m \\l,k\leq m}}}
\frac{\left(  -1\right)  ^{l}}{k!l!}\left(  \left(  aD\right)  ^{l+k}G\right)
\left(  0\right)  ,
\end{align*}

and since when $q>0$,%
\begin{align}%
{\displaystyle\sum\limits_{l+k=q}}
\frac{\left(  -1\right)  ^{l}}{k!l!}  & =%
{\displaystyle\sum\limits_{\left\vert \alpha\right\vert =q}}
\frac{\left(  \left(  -1\right)  ^{\alpha_{1}}1^{\alpha_{2}}\right)  \left(
1^{\alpha_{1}}1^{\alpha_{2}}\right)  }{\alpha!}=%
{\displaystyle\sum\limits_{\left\vert \alpha\right\vert =q}}
\frac{\left(  -1,1\right)  ^{\alpha}\left(  1,1\right)  ^{\alpha}}{\alpha
!}=\frac{\left(  \left(  -1,1\right)  \cdot\left(  1,1\right)  \right)  ^{q}%
}{q!}=\nonumber\\
& =0,\label{a2.13}%
\end{align}

we have%
\begin{align*}
\sum_{l,k=0}^{m}\frac{\left(  -1\right)  ^{l}}{k!l!}\left(  \left(  aD\right)
^{l+k}G\right)  \left(  0\right)   & =G\left(  0\right)  +%
{\displaystyle\sum\limits_{\substack{l+k>m \\l,k\leq m}}}
\frac{\left(  -1\right)  ^{l}}{k!l!}\left(  \left(  aD\right)  ^{l+k}G\right)
\left(  0\right) \\
& =G\left(  0\right)  +%
{\displaystyle\sum\limits_{q=m+1}^{2m}}
{\displaystyle\sum\limits_{\substack{l+k=q \\l,k\leq m}}}
\frac{\left(  -1\right)  ^{l}}{k!l!}\left(  \left(  aD\right)  ^{l+k}G\right)
\left(  0\right) \\
& =G\left(  0\right)  +%
{\displaystyle\sum\limits_{q=m+1}^{2m}}
\left(
{\displaystyle\sum\limits_{\substack{l+k=q \\l,k\leq m}}}
\frac{\left(  -1\right)  ^{l}}{k!l!}\right)  \left(  \left(  aD\right)
^{q}G\right)  \left(  0\right)  ,
\end{align*}

so that%
\begin{multline*}
\left(  2\pi\right)  ^{-\frac{d}{2}}\int\left\vert e^{ia\xi}-\sum_{k=0}%
^{m}\frac{\left(  ia\xi\right)  ^{k}}{k!}\right\vert ^{2}\widehat{G}\left(
\xi\right)  d\xi\\
=2G\left(  0\right)  -\sum_{k=0}^{m}\frac{\left(  \left(  aD\right)
^{k}G\right)  \left(  -a\right)  }{k!}-\sum_{k=0}^{m}\frac{\left(  \left(
-aD\right)  ^{k}G\right)  \left(  a\right)  }{k!}+%
{\displaystyle\sum\limits_{\substack{l+k>m \\l,k\leq m}}}
\frac{\left(  -1\right)  ^{l}}{k!l!}\left(  \left(  aD\right)  ^{l+k}G\right)
\left(  0\right)  .
\end{multline*}

But%
\[%
{\displaystyle\sum\limits_{\substack{l+k>m \\l,k\leq m}}}
\frac{\left(  -1\right)  ^{l}}{k!l!}\left(  \left(  aD\right)  ^{l+k}G\right)
=%
{\displaystyle\sum\limits_{q=m+1}^{2m}}
{\displaystyle\sum\limits_{\substack{l+k=q \\l,k\leq m}}}
\frac{\left(  -1\right)  ^{l}}{k!l!}\left(  \left(  aD\right)  ^{l+k}G\right)
=%
{\displaystyle\sum\limits_{q=m+1}^{2m}}
\left(
{\displaystyle\sum\limits_{\substack{l+k=q \\l,k\leq m }}}
\frac{\left(  -1\right)  ^{l}}{k!l!}\right)  \left(  aD\right)  ^{q}G,
\]

so that%
\begin{align}
\left(  2\pi\right)  ^{-\frac{d}{2}} &  \int\left\vert e^{ia\xi}-\sum
_{k=0}^{m}\frac{\left(  ia\xi\right)  ^{k}}{k!}\right\vert ^{2}\widehat{G}%
\left(  \xi\right)  d\xi\nonumber\\
&  =2G\left(  0\right)  -\sum_{k=0}^{m}\frac{\left(  \left(  aD\right)
^{k}G\right)  \left(  -a\right)  }{k!}-\sum_{k=0}^{m}\frac{\left(  \left(
-aD\right)  ^{k}G\right)  \left(  a\right)  }{k!}+\nonumber\\
&  \qquad\qquad+%
{\displaystyle\sum\limits_{q=m+1}^{2m}}
\left(
{\displaystyle\sum\limits_{\substack{l+k=q \\l,k\leq m}}}
\frac{\left(  -1\right)  ^{l}}{k!l!}\right)  \left(  aD\right)  ^{q}G\left(
0\right)  .\label{a337}%
\end{align}

Since $D^{\alpha}G\in L^{1}$ when $\left\vert \alpha\right\vert =2m+1$ we can
apply \ref{a1.55} with remainder \ref{a111} to get%
\begin{align*}
G\left(  0\right)   & =G\left(  a-a\right)  =\sum_{k=0}^{2m}\frac{\left(
\left(  -aD\right)  ^{k}G\right)  \left(  a\right)  }{k!}+\left(
\mathcal{R}_{2m+1}G\right)  \left(  a,-a\right)  ,\\
G\left(  0\right)   & =G\left(  -a+a\right)  =\sum_{k=0}^{2m}\frac{\left(
\left(  aD\right)  ^{k}G\right)  \left(  -a\right)  }{k!}+\left(
\mathcal{R}_{2m+1}G\right)  \left(  -a,a\right)  ,
\end{align*}

so that \ref{a337} becomes%
\begin{align}
&  \left(  2\pi\right)  ^{-\frac{d}{2}}\int\left\vert e^{ia\xi}-\sum_{k\leq
m}\frac{\left(  ia\xi\right)  ^{k}}{k!}\right\vert ^{2}\widehat{G}\left(
\xi\right)  d\xi\nonumber\\
&  =\sum_{k=m+1}^{2m}\frac{\left(  \left(  -aD\right)  ^{k}G\right)  \left(
a\right)  }{k!}+\sum_{k=m+1}^{2m}\frac{\left(  \left(  aD\right)
^{k}G\right)  \left(  -a\right)  }{k!}+%
{\displaystyle\sum\limits_{q=m+1}^{2m}}
\left(
{\displaystyle\sum\limits_{\substack{l+k=q \\l,k\leq m}}}
\frac{\left(  -1\right)  ^{l}}{k!l!}\right)  \left(  aD\right)  ^{q}G\left(
0\right)  +\nonumber\\
&  \qquad\qquad+\left(  \mathcal{R}_{2m+1}G\right)  \left(  a,-a\right)
+\left(  \mathcal{R}_{2m+1}G\right)  \left(  -a,a\right)  .\label{a2.16}%
\end{align}

Since $D^{\alpha}G\in L^{1}$ when $\left\vert \alpha\right\vert =2m+1$ we can
apply \ref{a1.55} with remainder \ref{a111} to $\left(  \left(  -aD\right)
^{k}G\right)  \left(  b\right)  $ about $b=0$ to get%
\[
\left(  \left(  aD\right)  ^{k}G\right)  \left(  b\right)  =\sum_{l=0}%
^{2m-k}\frac{\left(  \left(  bD\right)  ^{l}\left(  aD\right)  ^{k}G\right)
\left(  0\right)  }{l!}+\mathcal{R}_{2m+1-k}\left(  \left(  aD\right)
^{k}G\right)  \left(  0,b\right)  ,\quad0\leq k\leq2m.
\]

When $b=-a$ \ we get%
\begin{equation}
\left(  \left(  aD\right)  ^{k}G\right)  \left(  -a\right)  =\sum_{l=0}%
^{2m-k}\frac{\left(  -1\right)  ^{l}}{l!}\left(  \left(  aD\right)
^{k+l}G\right)  \left(  0\right)  +\mathcal{R}_{2m+1-k}\left(  \left(
aD\right)  ^{k}G\right)  \left(  0,-a\right)  ,\quad0\leq k\leq
2m,\label{a2.17}%
\end{equation}

and when $a\rightarrow-a$ and $b=a$ we get%
\begin{align}
\left(  \left(  -aD\right)  ^{k}G\right)   &  \left(  a\right) \nonumber\\
&  =\sum_{l=0}^{2m-k}\frac{\left(  \left(  aD\right)  ^{l}\left(  -aD\right)
^{k}G\right)  \left(  0\right)  }{l!}+\mathcal{R}_{2m+1-k}\left(  \left(
-aD\right)  ^{k}G\right)  \left(  0,a\right) \nonumber\\
&  =\sum_{l=0}^{2m-k}\frac{\left(  -1\right)  ^{k}}{l!}\left(  \left(
aD\right)  ^{k+l}G\right)  \left(  0\right)  +\mathcal{R}_{2m+1-k}\left(
\left(  -aD\right)  ^{k}G\right)  \left(  0,a\right)  ,\quad0\leq
k\leq2m.\label{a2.18}%
\end{align}

Now substitute the Taylor series expansions \ref{a2.17} and \ref{a2.18}\ into
the terms of the first and second summations of \ref{a2.16} to get%
\begin{align}
\sum_{k=m+1}^{2m} &  \frac{\left(  \left(  -aD\right)  ^{k}G\right)  \left(
a\right)  }{k!}\nonumber\\
&  =\sum_{k=m+1}^{2m}\frac{1}{k!}\left(  \sum_{l=0}^{2m-k}\frac{\left(
-1\right)  ^{k}}{l!}\left(  aD\right)  ^{k+l}G\left(  0\right)  +\mathcal{R}%
_{2m+1-k}\left(  \left(  -aD\right)  ^{k}G\right)  \left(  0,a\right)  \right)
\nonumber\\
&  =\sum_{k=m+1}^{2m}\frac{\left(  -1\right)  ^{k}}{k!}\sum_{l=0}^{2m-k}%
\frac{\left(  aD\right)  ^{k+l}G\left(  0\right)  }{l!}+\sum_{k=m+1}^{2m}%
\frac{1}{k!}\mathcal{R}_{2m+1-k}\left(  \left(  -aD\right)  ^{k}G\right)
\left(  0,a\right) \nonumber\\
&  =\sum_{k=m+1}^{2m}\sum_{l=0}^{2m-k}\frac{\left(  -1\right)  ^{k}}%
{k!l!}\left(  aD\right)  ^{k+l}G\left(  0\right)  +\ldots\nonumber\\
&  =\sum_{q=m+1}^{2m}\sum_{\substack{k+l=q \\k\geq m+1}}\frac{\left(
-1\right)  ^{k}}{k!l!}\left(  aD\right)  ^{k+l}G\left(  0\right)
+\ldots\nonumber\\
&  =\sum_{q=m+1}^{2m}\left(  \sum_{\substack{k+l=q \\k\geq m+1}}\frac{\left(
-1\right)  ^{k}}{k!l!}\right)  \left(  aD\right)  ^{q}G\left(  0\right)
+\ldots\nonumber\\
&  =\sum_{q=m+1}^{2m}\left(  \sum_{\substack{k+l=q \\k\geq m+1}}\frac{\left(
-1\right)  ^{k}}{k!l!}\right)  \left(  aD\right)  ^{q}G\left(  0\right)
+\sum_{k=m+1}^{2m}\frac{1}{k!}\mathcal{R}_{2m+1-k}\left(  \left(  -aD\right)
^{k}G\right)  \left(  0,a\right)  ,\label{a2.21}%
\end{align}

and so%
\begin{align}
&  \frac{\left(  \left(  aD\right)  ^{k}G\right)  \left(  -a\right)  }%
{k!}\nonumber\\
&  =\sum_{q=m+1}^{2m}\left(  \sum_{\substack{k+l=q \\k\geq m+1}}\frac{\left(
-1\right)  ^{k}}{k!l!}\right)  \left(  \left(  -aD\right)  ^{q}G\right)
\left(  0\right)  +\sum_{k=m+1}^{2m}\frac{1}{k!}\mathcal{R}_{2m+1-k}\left(
\left(  aD\right)  ^{k}G\right)  \left(  0,-a\right) \nonumber\\
&  =\sum_{q=m+1}^{2m}\left(  -1\right)  ^{q}\left(  \sum_{\substack{k+l=q
\\k\geq m+1}}\frac{\left(  -1\right)  ^{k}}{k!l!}\right)  \left(  aD\right)
^{q}G\left(  0\right)  +\ldots\nonumber\\
&  =\sum_{q=m+1}^{2m}\left(  -1\right)  ^{q}\left(  \sum_{\substack{k+l=q
\\k\geq m+1}}\frac{\left(  -1\right)  ^{-k}}{k!l!}\right)  \left(  aD\right)
^{q}G\left(  0\right)  +\ldots\nonumber\\
&  =\sum_{q=m+1}^{2m}\left(  \sum_{\substack{k+l=q \\k\geq m+1}}\frac{\left(
-1\right)  ^{q-k}}{k!l!}\right)  \left(  aD\right)  ^{q}G\left(  0\right)
+\ldots\nonumber\\
&  =\sum_{q=m+1}^{2m}\left(  \sum_{\substack{k+l=q \\k\geq m+1}}\frac{\left(
-1\right)  ^{l}}{k!l!}\right)  \left(  aD\right)  ^{q}G\left(  0\right)
+\sum_{k=m+1}^{2m}\frac{1}{k!}\mathcal{R}_{2m+1-k}\left(  \left(  aD\right)
^{k}G\right)  \left(  0,-a\right)  .\label{a2.22}%
\end{align}

Substituting \ref{a2.21} and \ref{a2.22} into \ref{a2.16} gives%
\begin{align*}
\left(  2\pi\right)  ^{-\frac{d}{2}} &  \int\left\vert e^{ia\xi}-\sum_{k\leq
m}\frac{\left(  ia\xi\right)  ^{k}}{k!}\right\vert ^{2}\widehat{G}\left(
\xi\right)  d\xi\\
&  =\sum_{q=m+1}^{2m}\left(  \sum_{\substack{k+l=q \\k\geq m+1}}\frac{\left(
-1\right)  ^{k}}{k!l!}\right)  \left(  aD\right)  ^{q}G\left(  0\right)
+\sum_{k=m+1}^{2m}\frac{1}{k!}\left(  \mathcal{R}_{2m+1}\left(  \left(
-aD\right)  ^{k}G\right)  \right)  \left(  0,a\right)  +\\
&  \quad+\sum_{q=m+1}^{2m}\left(  \sum_{\substack{k+l=q \\k\geq m+1}%
}\frac{\left(  -1\right)  ^{l}}{k!l!}\right)  \left(  aD\right)  ^{q}G\left(
0\right)  +\sum_{k=m+1}^{2m}\frac{1}{k!}\left(  \mathcal{R}_{2m+1}\left(
\left(  aD\right)  ^{k}G\right)  \right)  \left(  0,-a\right)  +\\
&  \quad+%
{\displaystyle\sum\limits_{q=m+1}^{2m}}
\left(
{\displaystyle\sum\limits_{\substack{l+k=q \\l,k\leq m}}}
\frac{\left(  -1\right)  ^{l}}{k!l!}\right)  \left(  aD\right)  ^{q}G\left(
0\right)  +\left(  \mathcal{R}_{2m+1}G\right)  \left(  a,-a\right)  +\left(
\mathcal{R}_{2m+1}G\right)  \left(  -a,a\right) \\
&  =\sum_{q=m+1}^{2m}\left(  \sum_{\substack{k+l=q \\k\geq m+1}}\frac{\left(
-1\right)  ^{k}}{k!l!}\right)  \left(  aD\right)  ^{q}G\left(  0\right)
+\sum_{q=m+1}^{2m}\left(  \sum_{\substack{k+l=q \\k\geq m+1}}\frac{\left(
-1\right)  ^{l}}{k!l!}\right)  \left(  aD\right)  ^{q}G\left(  0\right)  +\\
&  \quad+%
{\displaystyle\sum\limits_{q=m+1}^{2m}}
\left(
{\displaystyle\sum\limits_{\substack{l+k=q \\l,k\leq m}}}
\frac{\left(  -1\right)  ^{l}}{k!l!}\right)  \left(  aD\right)  ^{q}G\left(
0\right)  +\\
&  \quad+\sum_{k=m+1}^{2m}\frac{1}{k!}\left(  \mathcal{R}_{2m+1-k}\left(
\left(  -aD\right)  ^{k}G\right)  \right)  \left(  0,a\right)  +\sum
_{k=m+1}^{2m}\frac{1}{k!}\left(  \mathcal{R}_{2m+1-k}\left(  \left(
aD\right)  ^{k}G\right)  \right)  \left(  0,-a\right)  +\\
&  \quad+\left(  \mathcal{R}_{2m+1}G\right)  \left(  a,-a\right)  +\left(
\mathcal{R}_{2m+1}G\right)  \left(  -a,a\right) \\
&  =\sum_{q=m+1}^{2m}\left(  \sum_{\substack{k+l=q \\k\geq m+1}}\frac{\left(
-1\right)  ^{k}}{k!l!}+\left(  -1\right)  ^{q}\sum_{\substack{k+l=q \\k\geq
m+1}}\frac{\left(  -1\right)  ^{k}}{k!l!}+%
{\displaystyle\sum\limits_{\substack{l+k=q \\l,k\leq m}}}
\frac{\left(  -1\right)  ^{l}}{k!l!}\right)  \left(  aD\right)  ^{q}G\left(
0\right)  +\left(  \mathcal{R}_{2m+1}G\right)  \left(  a\right)  ,
\end{align*}

i.e.%
\begin{align}
&  \left(  2\pi\right)  ^{-\frac{d}{2}}\int\left\vert e^{ia\xi}-\sum_{k\leq
m}\frac{\left(  ia\xi\right)  ^{k}}{k!}\right\vert ^{2}\widehat{G}\left(
\xi\right)  d\xi\nonumber\\
&  =\sum_{q=m+1}^{2m}\left(  \sum_{\substack{k+l=q \\k\geq m+1}}\frac{\left(
-1\right)  ^{k}}{k!l!}+\sum_{\substack{k+l=q \\k\geq m+1}}\frac{\left(
-1\right)  ^{l}}{k!l!}+%
{\displaystyle\sum\limits_{\substack{k+l=q \\k,l\leq m}}}
\frac{\left(  -1\right)  ^{l}}{k!l!}\right)  \left(  aD\right)  ^{q}G\left(
0\right)  +\mathcal{R}_{2m+1}G\left(  a\right)  ,\label{a2.24}%
\end{align}

where $\mathcal{R}_{2m+1}G\left(  a\right)  $ was defined by \ref{a2.32} in
the statement of part 4.

Using \ref{a2.13} we get for $m+1\leq q\leq2m$,%
\begin{align*}
\sum_{\substack{k+l=q \\k\geq m+1}}\frac{\left(  -1\right)  ^{k}}{k!l!}%
+\sum_{\substack{k+l=q \\k\geq m+1}}\frac{\left(  -1\right)  ^{l}}{k!l!}+%
{\displaystyle\sum\limits_{\substack{k+l=q \\k,l\leq m}}}
\frac{\left(  -1\right)  ^{l}}{k!l!}  & =\sum_{\substack{k+l=q \\l\geq
m+1}}\frac{\left(  -1\right)  ^{l}}{k!l!}+\sum_{\substack{k+l=q \\k\geq
m+1}}\frac{\left(  -1\right)  ^{l}}{k!l!}+%
{\displaystyle\sum\limits_{\substack{k+l=q \\k,l\leq m}}}
\frac{\left(  -1\right)  ^{l}}{k!l!}\\
& =\sum_{\substack{k+l=q \\l\geq m+1}}\frac{\left(  -1\right)  ^{l}}%
{k!l!}+\sum_{\substack{k+l=q \\k\geq m+1}}\frac{\left(  -1\right)  ^{l}}%
{k!l!}+%
{\displaystyle\sum\limits_{\substack{k+l=q \\k,l\leq m}}}
\frac{\left(  -1\right)  ^{l}}{k!l!}\\
& =%
{\displaystyle\sum\limits_{k+l=q}}
\frac{\left(  -1\right)  ^{l}}{k!l!}\\
& =0,
\end{align*}

so that \ref{a2.24} becomes%
\begin{equation}
\left(  2\pi\right)  ^{-\frac{d}{2}}\int\left\vert e^{ia\xi}-\sum_{k\leq
m}\frac{\left(  ia\xi\right)  ^{k}}{k!}\right\vert ^{2}\widehat{G}\left(
\xi\right)  d\xi=\mathcal{R}_{2m+1}G\left(  a\right)  ,\label{a1.59}%
\end{equation}

which combined with \ref{a2.15} yields%
\begin{align*}
\left\vert \left(  \mathcal{R}_{m+1}f\right)  \left(  x,a\right)  \right\vert
& \leq\left(  \int\left\vert e^{ia\xi}-\sum_{k\leq m}\frac{\left(
ia\xi\right)  ^{k}}{k!}\right\vert ^{2}\widehat{G}\left(  \xi\right)
d\xi\right)  ^{1/2}\left\Vert f\right\Vert _{w,0}\\
& =\left(  2\pi\right)  ^{\frac{d}{4}}\sqrt{\mathcal{R}_{2m+1}G\left(
a\right)  }\left\Vert f\right\Vert _{w,0},
\end{align*}

which proves \ref{a2.30}.

Since $G\in C_{B}^{\left(  2m\right)  }$ and $\left(  aD\right)  ^{2m+1}G\in
L^{1}$, Lemma \ref{Lem_Taylor_estim_C_L1} implies that%
\begin{align}
&  \sum_{k=m+1}^{2m}\frac{1}{k!}\left(  \mathcal{R}_{2m+1-k}\left(  \left(
-aD\right)  ^{k}G\right)  \right)  \left(  0,a\right) \nonumber\\
&  =\sum_{k=m+1}^{2m}\frac{1}{k!}\left(  \frac{1}{\left(  2m-k\right)  !}%
\int_{0}^{1}g_{2m-k}\left(  s\right)  \left(  \left(  aD\right)
^{2m+1-k}\left(  -aD\right)  ^{k}G\right)  \left(  sa\right)  ds\right)
\nonumber\\
&  =\sum_{k=m+1}^{2m}\frac{1}{k!\left(  2m-k\right)  !}\int_{0}^{1}\left(
-1\right)  ^{k}g_{2m-k}\left(  s\right)  \left(  \left(  aD\right)
^{2m+1}G\right)  \left(  sa\right)  ds\nonumber\\
&  =\int_{0}^{1}\sum_{k=m+1}^{2m}\frac{\left(  -1\right)  ^{k}\left(
1-s\right)  ^{2m-k}}{k!\left(  2m-k\right)  !}\left(  \left(  aD\right)
^{2m+1}G\right)  \left(  sa\right)  ds\nonumber\\
&  =\int_{0}^{1}\sum_{j=0}^{m-1}\frac{\left(  -1\right)  ^{j}\left(
1-s\right)  ^{j}}{j!\left(  2m-j\right)  !}\left(  \left(  aD\right)
^{2m+1}G\right)  \left(  sa\right)  ds\nonumber\\
&  =\int_{0}^{1}\phi_{2m}\left(  s\right)  \left(  \left(  aD\right)
^{2m+1}G\right)  \left(  sa\right)  ds,\label{1.075}%
\end{align}

where%
\begin{equation}
\phi_{2m}\left(  s\right)  :=\left\{
\begin{array}
[c]{ll}%
0, & m=0,\\
\sum\limits_{j=0}^{m-1}\frac{\left(  -1\right)  ^{j}\left(  1-s\right)  ^{j}%
}{j!\left(  2m-j\right)  !}, & m=1,2,3,\ldots,
\end{array}
\right. \label{1.074}%
\end{equation}

and so $a\rightarrow-a$ yields%
\begin{align*}
\sum_{k=m+1}^{2m}\frac{1}{k!}\left(  \mathcal{R}_{2m+1-k}\left(  \left(
aD\right)  ^{k}G\right)  \right)  \left(  0,-a\right)   & =\frac{1}{\left(
2m\right)  !}\int_{0}^{1}\phi_{2m}\left(  s\right)  \left(  \left(
-aD\right)  ^{2m+1}G\right)  \left(  -sa\right)  ds\\
& =\frac{-1}{\left(  2m\right)  !}\int_{0}^{1}\phi_{2m}\left(  s\right)
\left(  \left(  aD\right)  ^{2m+1}G\right)  \left(  -sa\right)  ds\\
& =\frac{1}{\left(  2m\right)  !}\int_{0}^{1}\phi_{2m}\left(  s\right)
\left(  \left(  aD\right)  ^{2m+1}\overline{G}\right)  \left(  sa\right)  ds,
\end{align*}

which means that%
\begin{align}
\sum_{k=m+1}^{2m} &  \frac{1}{k!}\left\{  \left(  \mathcal{R}_{2m+1}\left(
\left(  -aD\right)  ^{k}G\right)  \right)  \left(  0,a\right)  +\mathcal{R}%
_{2m+1}\left(  \left(  aD\right)  ^{k}G\right)  \left(  0,-a\right)  \right\}
\nonumber\\
&  =\frac{1}{\left(  2m\right)  !}\int_{0}^{1}\phi_{2m}\left(  s\right)
\left(  \left(  aD\right)  ^{2m+1}G\right)  \left(  sa\right)  ds+\frac
{1}{\left(  2m\right)  !}\int_{0}^{1}\phi_{2m}\left(  s\right)  \left(
\left(  aD\right)  ^{2m+1}\overline{G}\right)  \left(  sa\right)
ds\nonumber\\
&  =\frac{2}{\left(  2m\right)  !}\int_{0}^{1}\phi_{2m}\left(  s\right)
\left(  \left(  aD\right)  ^{2m+1}\operatorname{Re}G\right)  \left(
sa\right)  ds,\label{a2.26}%
\end{align}

Since $G\in C_{B}^{\left(  2m\right)  }$ and $\left(  aD\right)  ^{2m+1}G\in
L^{1}$, Lemma \ref{Lem_Taylor_estim_C_L1} again implies that%
\begin{align*}
\left(  \mathcal{R}_{2m+1}G\right)  \left(  a,-a\right)   & =\frac{1}{\left(
2m\right)  !}\int_{0}^{1}g_{2m}\left(  s\right)  \left(  \left(  -aD\right)
^{2m+1}G\right)  \left(  a-sa\right)  ds\\
& =\frac{-1}{\left(  2m\right)  !}\int_{0}^{1}g_{2m}\left(  s\right)  \left(
\left(  aD\right)  ^{2m+1}G\right)  \left(  \left(  1-s\right)  a\right)  ds\\
& =\frac{-1}{\left(  2m\right)  !}\int_{0}^{1}\left(  1-s\right)  ^{2m}\left(
\left(  aD\right)  ^{2m+1}G\right)  \left(  \left(  1-s\right)  a\right)  ds\\
& =\frac{1}{\left(  2m\right)  !}\int_{1}^{0}t^{2m}\left(  \left(  aD\right)
^{2m+1}G\right)  \left(  ta\right)  dt\\
& =\frac{-1}{\left(  2m\right)  !}\int_{0}^{1}t^{2m}\left(  \left(  aD\right)
^{2m+1}G\right)  \left(  ta\right)  dt,
\end{align*}

and so $a\rightarrow-a$ yields%
\begin{align}
\left(  \mathcal{R}_{2m+1}G\right)  \left(  -a,a\right)   & =\frac{-1}{\left(
2m\right)  !}\int_{0}^{1}t^{2m}\left(  \left(  -aD\right)  ^{2m+1}G\right)
\left(  -ta\right)  dt\nonumber\\
& =\frac{1}{\left(  2m\right)  !}\int_{0}^{1}t^{2m}\left(  \left(  aD\right)
^{2m+1}G\right)  \left(  -ta\right)  dt\label{a2.25}\\
& =\frac{-1}{\left(  2m\right)  !}\int_{0}^{1}t^{2m}\left(  \left(
aD_{x}\right)  ^{2m+1}\left(  G\left(  -x\right)  \right)  \right)  \left(
ta\right)  dt\nonumber\\
& =\frac{-1}{\left(  2m\right)  !}\int_{0}^{1}t^{2m}\left(  \left(  aD\right)
^{2m+1}\overline{G}\right)  \left(  ta\right)  dt,\nonumber
\end{align}

and hence
\begin{equation}
\left(  \mathcal{R}_{2m+1}G\right)  \left(  a,-a\right)  +\left(
\mathcal{R}_{2m+1}G\right)  \left(  -a,a\right)  =\frac{-2}{\left(  2m\right)
!}\int_{0}^{1}t^{2m}\left(  \left(  aD\right)  ^{2m+1}\operatorname{Re}%
G\right)  \left(  ta\right)  dt.\label{a2.28}%
\end{equation}

Now \ref{a2.26}, \ref{a2.28} substituted into \ref{a2.32}\ gives%
\begin{align*}
\left(  \mathcal{R}_{2m+1}G\right)  \left(  a\right)   & =\sum_{k=m+1}%
^{2m}\frac{1}{k!}\left\{  \left(  \mathcal{R}_{2m+1-k}\left(  \left(
-aD\right)  ^{k}G\right)  \right)  \left(  0,a\right)  +\left(  \mathcal{R}%
_{2m+1-k}\left(  \left(  aD\right)  ^{k}G\right)  \right)  \left(
0,-a\right)  \right\}  +\\
& \qquad\qquad+\left(  \mathcal{R}_{2m+1}G\right)  \left(  a,-a\right)
+\left(  \mathcal{R}_{2m+1}G\right)  \left(  -a,a\right) \\
& =2\int_{0}^{1}\phi_{2m}\left(  s\right)  \left(  \left(  aD\right)
^{2m+1}\operatorname{Re}G\right)  \left(  sa\right)  ds-\frac{2}{\left(
2m\right)  !}\int_{0}^{1}s^{2m}\left(  \left(  aD\right)  ^{2m+1}%
\operatorname{Re}G\right)  \left(  sa\right)  ds\\
& =2\int_{0}^{1}\left(  \sum\limits_{j=0}^{m-1}\frac{\left(  -1\right)
^{j}\left(  1-s\right)  ^{j}}{j!\left(  2m-j\right)  !}-\frac{s^{2m}}{\left(
2m\right)  !}\right)  \left(  \left(  aD\right)  ^{2m+1}\operatorname{Re}%
G\right)  \left(  sa\right)  ds,
\end{align*}

which is \ref{a2.33}.\medskip

\textbf{Part 5} We now calculate%
\[
\int_{0}^{1}\left\vert \sum\limits_{j=0}^{m-1}\frac{\left(  -1\right)
^{j}\left(  1-s\right)  ^{j}}{j!\left(  2m-j\right)  !}-\frac{s^{2m}}{\left(
2m\right)  !}\right\vert ds.
\]

Clearly%
\[
\sum\limits_{j=0}^{m-1}\frac{\left(  -1\right)  ^{j}\left(  1-s\right)  ^{j}%
}{j!\left(  2m-j\right)  !}=\frac{1}{\left(  2m\right)  !}\sum\limits_{j=0}%
^{m-1}\tbinom{2m}{j}\left(  -1\right)  ^{j}\left(  1-s\right)  ^{j},
\]

but%
\begin{align*}
\sum\limits_{j=0}^{m-1}\tbinom{2m}{j}\left(  -1\right)  ^{j}\left(
1-s\right)  ^{j}  & =\sum\limits_{j=0}^{2m}\tbinom{2m}{j}\left(  -1\right)
^{2m-j}\left(  1-s\right)  ^{j}-\sum\limits_{j=m}^{2m}\tbinom{2m}{j}\left(
-1\right)  ^{2m-j}\left(  1-s\right)  ^{j}\\
& =\left(  -1+1-s\right)  ^{2m}-\sum\limits_{j=m}^{2m}\tbinom{2m}{j}\left(
-1\right)  ^{2m-j}\left(  1-s\right)  ^{j}\\
& =s^{2m}-\sum\limits_{j=m}^{2m}\tbinom{2m}{j}\left(  -1\right)
^{2m-j}\left(  1-s\right)  ^{j}\\
& =s^{2m}-\sum\limits_{k=m}^{0}\tbinom{2m}{2m-k}\left(  -1\right)  ^{k}\left(
1-s\right)  ^{2m-k}\\
& =s^{2m}-\sum\limits_{k=0}^{m}\tbinom{2m}{k}\left(  -1\right)  ^{k}\left(
1-s\right)  ^{2m-k}\\
& =s^{2m}-\left(  1-s\right)  ^{m}\sum\limits_{k=0}^{m}\tbinom{2m}{k}\left(
-1\right)  ^{k}\left(  1-s\right)  ^{m-k}\\
& =s^{2m}-\left(  1-s\right)  ^{m}\sum\limits_{k=0}^{m}\tbinom{2m}{m-k}\left(
-1\right)  ^{m-k}\left(  1-s\right)  ^{k}\\
& =s^{2m}+\left(  -1\right)  ^{m+1}\left(  1-s\right)  ^{m}\sum\limits_{k=0}%
^{m}\tbinom{2m}{m-k}\left(  -1\right)  ^{k}\left(  1-s\right)  ^{k},
\end{align*}

so that%
\begin{align*}
\sum\limits_{j=0}^{m-1} &  \frac{\left(  -1\right)  ^{j}\left(  1-s\right)
^{j}}{j!\left(  2m-j\right)  !}-\frac{s^{2m}}{\left(  2m\right)  !}\\
&  =\frac{1}{\left(  2m\right)  !}\sum\limits_{j=0}^{m-1}\tbinom{2m}{j}\left(
-1\right)  ^{j}\left(  1-s\right)  ^{j}-\frac{s^{2m}}{\left(  2m\right)  !}\\
&  =\frac{1}{\left(  2m\right)  !}\left(  s^{2m}+\left(  -1\right)
^{m+1}\left(  1-s\right)  ^{m}\sum\limits_{k=0}^{m}\tbinom{2m}{m-k}\left(
-1\right)  ^{k}\left(  1-s\right)  ^{k}\right)  -\frac{s^{2m}}{\left(
2m\right)  !}\\
&  =\frac{\left(  -1\right)  ^{m+1}}{\left(  2m\right)  !}\left(  1-s\right)
^{m}\sum\limits_{k=0}^{m}\tbinom{2m}{m-k}\left(  -1\right)  ^{k}\left(
1-s\right)  ^{k},
\end{align*}

and hence%
\begin{align*}
\int_{0}^{1}\left\vert \sum\limits_{j=0}^{m-1}\frac{\left(  -1\right)
^{j}\left(  1-s\right)  ^{j}}{j!\left(  2m-j\right)  !}-\frac{s^{2m}}{\left(
2m\right)  !}\right\vert ds  & =\frac{1}{\left(  2m\right)  !}\int_{0}%
^{1}\left(  1-s\right)  ^{m}\sum\limits_{k=0}^{m}\tbinom{2m}{m-k}\left(
-1\right)  ^{k}\left(  1-s\right)  ^{k}ds\\
& =\frac{1}{\left(  2m\right)  !}\int_{0}^{1}t^{m}\sum\limits_{k=0}^{m}%
\tbinom{2m}{m-k}\left(  -1\right)  ^{k}t^{k}dt.
\end{align*}

We consider two cases: $m$ odd and $m$ even.\smallskip

\fbox{\textbf{Case} $m$ is odd} Regarding the integrand:\smallskip%
\begin{align*}
&  t^{m}\sum\limits_{k=0}^{m}\tbinom{2m}{m-k}\left(  -1\right)  ^{k}t^{k}\\
&  =t^{m}\left(  \tbinom{2m}{m}-\tbinom{2m}{m-1}t+\tbinom{2m}{m-2}%
t^{2}-\tbinom{2m}{m-3}t^{3}+\ldots+\tbinom{2m}{1}t^{m-1}-\tbinom{2m}{0}%
t^{m}\right) \\
&  =t^{m}\left(  \left\{  \tbinom{2m}{m}-\tbinom{2m}{m-1}t\right\}  +\left\{
\tbinom{2m}{m-2}t^{2}-\tbinom{2m}{m-3}t^{3}\right\}  +\ldots+\left\{
\tbinom{2m}{1}t^{m-1}-\tbinom{2m}{0}t^{m}\right\}  \right) \\
&  =\left\{  \tbinom{2m}{m}t^{m}-\tbinom{2m}{m-1}t^{m+1}\right\}  +\left\{
\tbinom{2m}{m-2}t^{m+2}-\tbinom{2m}{m-3}t^{m+3}\right\}  +\ldots+\left\{
\tbinom{2m}{1}t^{2m-1}-\tbinom{2m}{0}t^{2m}\right\} \\
&  \geq0,
\end{align*}

since $0\leq t\leq1$ implies this expression is the sum of $m$ non-negative
terms. Further%
\begin{align*}
&  \frac{1}{\left(  2m\right)  !}\int_{0}^{1}t^{m}\sum\limits_{k=0}^{m}%
\tbinom{2m}{m-k}\left(  -1\right)  ^{k}t^{k}ds\\
&  =\frac{1}{\left(  2m\right)  !}\left(
\begin{array}
[c]{c}%
\left\{  \tbinom{2m}{m}\frac{1}{m+1}-\tbinom{2m}{m-1}\frac{1}{m+2}\right\}
+\left\{  \tbinom{2m}{m-2}\frac{1}{m+3}-\tbinom{2m}{m-3}\frac{1}{m+4}\right\}
+\ldots\\
+\left\{  \tbinom{2m}{1}\frac{1}{2m}-\tbinom{2m}{0}\frac{1}{2m+1}\right\}
\end{array}
\right) \\
&  =\left\{  \frac{1}{m!m!}\frac{1}{m+1}-\frac{1}{\left(  m-1\right)  !\left(
m+1\right)  !}\frac{1}{m+2}\right\}  +\\
&  \qquad\qquad+\left\{  \frac{1}{\left(  m-2\right)  !\left(  m+2\right)
!}\frac{1}{m+3}-\frac{1}{\left(  m+1\right)  !\left(  m+3\right)  !}\frac
{1}{m+4}\right\}  +\ldots\\
&  \qquad\qquad+\left\{  \frac{1}{1!\left(  2m-1\right)  !}\frac{1}{2m}%
-\frac{1}{0!\left(  2m\right)  !}\frac{1}{2m+1}\right\} \\
&  =\left\{  \frac{1}{m!\left(  m+1\right)  !}-\frac{1}{\left(  m-1\right)
!\left(  m+2\right)  !}\right\}  +\left\{  \frac{1}{\left(  m-2\right)
!\left(  m+3\right)  !}-\frac{1}{\left(  m-3\right)  !\left(  m+4\right)
!}\right\}  +\ldots\\
&  \qquad\qquad+\left\{  \frac{1}{1!\left(  2m\right)  !}-\frac{1}{0!\left(
2m+1\right)  !}\right\} \\
&  =\frac{1}{\left(  2m+1\right)  !}\left(  \left\{  \tbinom{2m+1}{m}%
-\tbinom{2m+1}{m-1}\right\}  +\left\{  \tbinom{2m+1}{m-2}-\tbinom{2m+1}%
{m-3}\right\}  +\ldots+\left\{  \tbinom{2m+1}{1}-\tbinom{2m+1}{0}\right\}
\right) \\
&  =\frac{1}{\left(  2m+1\right)  !}\left(
\begin{array}
[c]{c}%
\left\{  \tbinom{2m}{m}+\tbinom{2m}{m-1}\right\}  -\left\{  \tbinom{2m}%
{m-1}+\tbinom{2m}{m-2}\right\}  +\left\{  \tbinom{2m}{m-2}+\tbinom{2m}%
{m-3}\right\}  -\ldots\\
+\left\{  \tbinom{2m}{1}+\tbinom{2m}{0}\right\}  -\tbinom{2m}{0}%
\end{array}
\right) \\
&  =\frac{1}{\left(  2m+1\right)  !}\binom{2m}{m},
\end{align*}

which means that when $m$ is odd:%
\begin{align*}
\frac{1}{\left(  2m\right)  !}\int_{0}^{1}\left\vert \sum\limits_{j=0}%
^{m-1}\frac{\left(  -1\right)  ^{j}\left(  1-s\right)  ^{j}}{j!\left(
2m-j\right)  !}-\frac{s^{2m}}{\left(  2m\right)  !}\right\vert ds  & =\frac
{1}{\left(  2m+1\right)  !}\binom{2m}{m}\\
& =\frac{1}{\left(  2m+1\right)  \left(  m!\right)  ^{2}}.
\end{align*}
\medskip

\fbox{\textbf{Case} $m$ is even}%
\begin{align*}
&  t^{m}\sum\limits_{k=0}^{m}\tbinom{2m}{m-k}\left(  -1\right)  ^{k}t^{k}\\
&  =t^{m}\left(  \tbinom{2m}{m}-\tbinom{2m}{m-1}t+\tbinom{2m}{m-2}%
t^{2}-\tbinom{2m}{m-3}t^{3}+\ldots-\tbinom{2m}{1}t^{m-1}+\tbinom{2m}{0}%
t^{m}\right) \\
&  =t^{m}\left(
\begin{array}
[c]{c}%
\left\{  \tbinom{2m}{m}-\tbinom{2m}{m-1}t\right\}  +\left\{  \tbinom{2m}%
{m-2}t^{2}-\tbinom{2m}{m-3}t^{3}\right\}  +\ldots\\
+\left\{  \tbinom{2m}{2}t^{m-2}-\tbinom{2m}{1}t^{m-1}\right\}  +\tbinom{2m}%
{0}t^{m}%
\end{array}
\right) \\
&  =\left\{  \tbinom{2m}{m}t^{m}-\tbinom{2m}{m-1}t^{m+1}\right\}  +\left\{
\tbinom{2m}{m-2}t^{m+2}-\tbinom{2m}{m-3}t^{m+3}\right\}  +\ldots\\
&  \qquad\qquad+\left\{  \tbinom{2m}{2}t^{2m-2}-\tbinom{2m}{1}t^{2m-1}%
\right\}  +\tbinom{2m}{0}t^{2m}\\
&  \geq0,
\end{align*}

since $0\leq t\leq1$ implies this expression is the sum of $m+1$ non-negative
terms. Further%
\begin{align*}
&  \frac{1}{\left(  2m\right)  !}\int_{0}^{1}t^{m}\sum\limits_{k=0}^{m}%
\tbinom{2m}{m-k}\left(  -1\right)  ^{k}t^{k}ds\\
&  =\frac{1}{\left(  2m\right)  !}\left(
\begin{array}
[c]{c}%
\left\{  \tbinom{2m}{m}\frac{1}{m+1}-\tbinom{2m}{m-1}\frac{1}{m+2}\right\}
+\left\{  \tbinom{2m}{m-2}\frac{1}{m+3}-\tbinom{2m}{m-3}\frac{1}{m+4}\right\}
+\ldots\\
+\left\{  \tbinom{2m}{2}\frac{1}{2m-1}-\tbinom{2m}{1}\frac{1}{2m}\right\}
++\tbinom{2m}{0}\frac{1}{2m+1}%
\end{array}
\right) \\
&  =\left\{  \frac{1}{m!m!}\frac{1}{m+1}-\frac{1}{\left(  m-1\right)  !\left(
m+1\right)  !}\frac{1}{m+2}\right\}  +\\
&  \qquad\qquad+\left\{  \frac{1}{\left(  m-2\right)  !\left(  m+2\right)
!}\frac{1}{m+3}-\frac{1}{\left(  m+1\right)  !\left(  m+3\right)  !}\frac
{1}{m+4}\right\}  +\ldots\\
&  \qquad\qquad+\left\{  \frac{1}{2!\left(  2m-2\right)  !}\frac{1}%
{2m-1}-\frac{1}{1!\left(  2m-1\right)  !}\frac{1}{2m}\right\}  +\frac
{1}{0!\left(  2m\right)  !}\frac{1}{2m+1}\\
&  =\left\{  \frac{1}{m!\left(  m+1\right)  !}-\frac{1}{\left(  m-1\right)
!\left(  m+2\right)  !}\right\}  +\left\{  \frac{1}{\left(  m-2\right)
!\left(  m+3\right)  !}-\frac{1}{\left(  m-3\right)  !\left(  m+4\right)
!}\right\}  +\ldots\\
&  \qquad\qquad+\left\{  \frac{1}{2!\left(  2m-1\right)  !}-\frac{1}{1!\left(
2m\right)  !}\right\}  +\frac{1}{0!\left(  2m+1\right)  !}\\
&  =\frac{1}{\left(  2m+1\right)  !}\left(  \tbinom{2m+1}{m}-\tbinom
{2m+1}{m-1}+\tbinom{2m+1}{m-2}-\tbinom{2m+1}{m-3}+\ldots+\tbinom{2m+1}%
{2}-\tbinom{2m+1}{1}+\tbinom{2m+1}{0}\right) \\
&  =\frac{1}{\left(  2m+1\right)  !}\left(
\begin{array}
[c]{c}%
\left\{  \tbinom{2m}{m}+\tbinom{2m}{m-1}\right\}  -\left\{  \tbinom{2m}%
{m-1}+\tbinom{2m}{m-2}\right\}  +\left\{  \tbinom{2m}{m-2}+\tbinom{2m}%
{m-3}\right\}  -\ldots\\
-\left\{  \tbinom{2m}{1}+\tbinom{2m}{0}\right\}  -\tbinom{2m}{0}%
\end{array}
\right) \\
&  =\frac{1}{\left(  2m+1\right)  !}\binom{2m}{m},
\end{align*}

which means that when $m$ is even:%
\begin{align*}
\frac{1}{\left(  2m\right)  !}\int_{0}^{1}\left\vert \sum\limits_{j=0}%
^{m-1}\frac{\left(  -1\right)  ^{j}\left(  1-s\right)  ^{j}}{j!\left(
2m-j\right)  !}-\frac{s^{2m}}{\left(  2m\right)  !}\right\vert ds  & =\frac
{1}{\left(  2m+1\right)  !}\binom{2m}{m}\\
& =\frac{1}{\left(  2m+1\right)  \left(  m!\right)  ^{2}}.
\end{align*}

\end{proof}

\begin{corollary}
If $H=\operatorname{Re}G$ and $f\in X_{w}^{0}$ then%
\begin{align}
\left\vert \left(  \mathcal{R}_{m+1}f\right)  \left(  x,a\right)  \right\vert
&  \leq\left(  \int\left\vert e^{ia\xi}-%
{\displaystyle\sum\limits_{k=0}^{m}}
\frac{\left(  ia\xi\right)  ^{k}}{k!}\right\vert ^{2}\widehat{H}\left(
\xi\right)  d\xi\right)  ^{1/2}\left\Vert f\right\Vert _{w,0}\nonumber\\
&  =\left(  2\pi\right)  ^{\frac{d}{4}}\left(  2H\left(  0\right)  -2%
{\displaystyle\sum\limits_{k=0}^{m}}
\frac{\left(  -1\right)  ^{k}}{k!}\left(  aD\right)  ^{k}H\left(  a\right)
+\right. \nonumber\\
&  \qquad\qquad\left.  +%
{\displaystyle\sum\limits_{j,k\leq m;\text{ }j+k>m}}
\frac{\left(  -1\right)  ^{j}}{j!k!}\left(  \left(  aD\right)  ^{j+k}H\right)
\left(  0\right)  \right)  ^{\frac{1}{2}}\left\Vert f\right\Vert
_{w,0},\label{a1.66}%
\end{align}

for all $x,a\in\mathbb{R}^{d}$.
\end{corollary}

\begin{proof}
Regarding part 2 of Theorem \ref{Thm_Tayor_rem_estim} set%
\[
f\left(  a\right)  =\int\left\vert e^{ia\xi}-%
{\displaystyle\sum\limits_{k=0}^{m}}
\frac{\left(  ia\xi\right)  ^{k}}{k!}\right\vert ^{2}\widehat{G}\left(
\xi\right)  d\xi,
\]

so that%
\begin{align*}
f\left(  -a\right)  =\int\left\vert e^{-ia\xi}-%
{\displaystyle\sum\limits_{k=0}^{m}}
\frac{\left(  -ia\xi\right)  ^{k}}{k!}\right\vert ^{2}\widehat{G}\left(
\xi\right)  d\xi & =\int\left\vert e^{ia\eta}-%
{\displaystyle\sum\limits_{k=0}^{m}}
\frac{\left(  ia\eta\right)  ^{k}}{k!}\right\vert ^{2}\widehat{G}\left(
-\eta\right)  d\eta\\
& =\int\left\vert e^{ia\eta}-%
{\displaystyle\sum\limits_{k=0}^{m}}
\frac{\left(  ia\eta\right)  ^{k}}{k!}\right\vert ^{2}\widehat{G\left(
-x\right)  }\left(  \eta\right)  d\eta\\
& =\int\left\vert e^{ia\eta}-%
{\displaystyle\sum\limits_{k=0}^{m}}
\frac{\left(  ia\eta\right)  ^{k}}{k!}\right\vert ^{2}\widehat{\overline{G}%
}\left(  \eta\right)  d\eta.
\end{align*}

Hence%
\begin{equation}
\int\left\vert e^{ia\xi}-%
{\displaystyle\sum\limits_{k=0}^{m}}
\frac{\left(  ia\xi\right)  ^{k}}{k!}\right\vert ^{2}%
\widehat{\operatorname{Im}G}\left(  \xi\right)  d\xi=0,\label{a1.64}%
\end{equation}

and%
\begin{equation}
\int\left\vert e^{ia\xi}-%
{\displaystyle\sum\limits_{k=0}^{m}}
\frac{\left(  ia\xi\right)  ^{k}}{k!}\right\vert ^{2}\widehat{G}\left(
\xi\right)  d\xi=\int\left\vert e^{ia\xi}-%
{\displaystyle\sum\limits_{k=0}^{m}}
\frac{\left(  ia\xi\right)  ^{k}}{k!}\right\vert ^{2}%
\widehat{\operatorname{Re}G}\left(  \xi\right)  d\xi.\label{a1.65}%
\end{equation}

From Theorem \ref{Thm_basis_property_minus} we know that, as distributions,
\[
\left(  \left(  aD\right)  ^{k}G\right)  \left(  -\cdot\right)  =\left(
-1\right)  ^{k}\left(  aD\right)  ^{k}\overline{G}=\left(  -aD\right)
^{k}\overline{G},\quad k\geq0,\text{ }a\in\mathbb{R}^{d}.
\]

Also $\operatorname{Re}D^{\alpha}G$ is an even distribution when $\left\vert
\alpha\right\vert $ is even and an odd distribution when $\left\vert
\alpha\right\vert $ is odd.

Consequently%
\begin{align*}%
{\displaystyle\sum\limits_{\substack{k\leq m \\\,}}}
\frac{1}{k!} &  \left(  \left(  aD\right)  ^{k}G\right)  \left(  -a\right)  +%
{\displaystyle\sum\limits_{k=0}^{m}}
\frac{1}{k!}\left(  -aD\right)  ^{k}G\left(  a\right) \\
&  =%
{\displaystyle\sum\limits_{k=0}^{m}}
\frac{\left(  -1\right)  ^{k}}{k!}\left(  \left(  aD\right)  ^{k}\overline
{G}\right)  \left(  a\right)  +%
{\displaystyle\sum\limits_{k=0}^{m}}
\frac{\left(  -1\right)  ^{k}}{k!}\left(  aD\right)  ^{k}G\left(  a\right) \\
&  =2%
{\displaystyle\sum\limits_{k=0}^{m}}
\frac{\left(  -1\right)  ^{k}}{k!}\left(  aD\right)  ^{k}\operatorname{Re}%
G\left(  a\right)
\end{align*}

and inequalities \ref{a2.15}\ and \ref{a2.31} become%
\begin{align*}
\left\vert \left(  \mathcal{R}_{m+1}f\right)  \left(  x,a\right)  \right\vert
&  \leq\left(  \int\left\vert e^{ia\xi}-%
{\displaystyle\sum\limits_{k=0}^{m}}
\frac{\left(  ia\xi\right)  ^{k}}{k!}\right\vert ^{2}%
\widehat{\operatorname{Re}G}\left(  \xi\right)  d\xi\right)  ^{1/2}\left\Vert
f\right\Vert _{w,0}\\
&  =\left(  2\pi\right)  ^{\frac{d}{4}}\left(  2G\left(  0\right)
-2\operatorname{Re}%
{\displaystyle\sum\limits_{k=0}^{m}}
\frac{\left(  -1\right)  ^{k}}{k!}\left(  aD\right)  ^{k}G\left(  a\right)
+\right. \\
&  \qquad\qquad\left.  +%
{\displaystyle\sum\limits_{\substack{j,k\leq m \\j+k>m}}}
\frac{\left(  -1\right)  ^{j}}{j!k!}\left(  \left(  aD\right)  ^{j+k}G\right)
\left(  0\right)  \right)  ^{\frac{1}{2}}\left\Vert f\right\Vert _{w,0}\\
&  =\operatorname{Re}\left(  2\pi\right)  ^{\frac{d}{4}}\left(  2G\left(
0\right)  -2\operatorname{Re}%
{\displaystyle\sum\limits_{k=0}^{m}}
\frac{\left(  -1\right)  ^{k}}{k!}\left(  aD\right)  ^{k}G\left(  a\right)
+\right. \\
&  \qquad\qquad\left.  +%
{\displaystyle\sum\limits_{\substack{j,k\leq m \\j+k>m}}}
\frac{\left(  -1\right)  ^{j}}{j!k!}\left(  \left(  aD\right)  ^{j+k}G\right)
\left(  0\right)  \right)  ^{\frac{1}{2}}\left\Vert f\right\Vert _{w,0}\\
&  =\left(  2\pi\right)  ^{\frac{d}{4}}\left(  2\operatorname{Re}G\left(
0\right)  -2%
{\displaystyle\sum\limits_{k=0}^{m}}
\frac{\left(  -1\right)  ^{k}}{k!}\left(  aD\right)  ^{k}\operatorname{Re}%
G\left(  a\right)  +\right. \\
&  \qquad\qquad\left.  +%
{\displaystyle\sum\limits_{\substack{j,k\leq m \\j+k>m}}}
\frac{\left(  -1\right)  ^{j}}{j!k!}\left(  \left(  aD\right)  ^{j+k}%
\operatorname{Re}G\right)  \left(  0\right)  \right)  ^{\frac{1}{2}}\left\Vert
f\right\Vert _{w,0},
\end{align*}

i.e. if $H=\operatorname{Re}G$ then%
\begin{align*}
\left\vert \left(  \mathcal{R}_{m+1}f\right)  \left(  x,a\right)  \right\vert
&  \leq\left(  \int\left\vert e^{ia\xi}-%
{\displaystyle\sum\limits_{k=0}^{m}}
\frac{\left(  ia\xi\right)  ^{k}}{k!}\right\vert ^{2}\widehat{H}\left(
\xi\right)  d\xi\right)  ^{\frac{1}{2}}\left\Vert f\right\Vert _{w,0}\\
&  =\left(  2\pi\right)  ^{\frac{d}{4}}\left(  2H\left(  0\right)  -2%
{\displaystyle\sum\limits_{k=0}^{m}}
\frac{\left(  -1\right)  ^{k}}{k!}\left(  aD\right)  ^{k}H\left(  a\right)
+\right. \\
&  \qquad\qquad\left.  +%
{\displaystyle\sum\limits_{j,k\leq m,\text{ }j+k>m}}
\frac{\left(  -1\right)  ^{j}}{j!k!}\left(  \left(  aD\right)  ^{j+k}H\right)
\left(  0\right)  \right)  ^{\frac{1}{2}}\left\Vert f\right\Vert _{w,0}.
\end{align*}

\end{proof}

We now want to obtain some upper bounds for $\mathcal{R}_{2m+1}\left(
a\right)  $ given by \ref{a2.32}.

\begin{corollary}
\label{Cor_Thm_Tayor_rem_estim_1}Regarding the remainder $\mathcal{R}_{2m+1}$
of part 4 of Theorem \ref{Thm_Tayor_rem_estim}:

\begin{enumerate}
\item If $\left\Vert \left(  \widehat{a}D\right)  ^{2m+1}G\left(  x\right)
\right\Vert _{\infty,B_{r}}<\infty$ then%
\[
\mathcal{R}_{2m+1}G\left(  a\right)  \leq\frac{1}{\left(  2m+1\right)  \left(
m!\right)  ^{2}}\left\Vert \left(  \widehat{\cdot}D\right)  ^{2m+1}%
G\right\Vert _{\infty,B_{r}}\left\vert a\right\vert ^{2m+1},\quad\left\vert
a\right\vert \leq r.
\]

\item If $\max\limits_{\left\vert \alpha\right\vert =2m+1}\left\Vert
D^{\alpha}G\right\Vert _{\infty,B_{r}}<\infty$ then%
\[
\mathcal{R}_{2m+1}G\left(  a\right)  \leq\binom{2m}{m}\left(  \sum_{\left\vert
\beta\right\vert =2m+1}\frac{1}{\beta!}\left\Vert D^{\beta}G\right\Vert
_{\infty,B_{r}}\right)  \left\vert a\right\vert ^{2m+1},\quad\left\vert
a\right\vert \leq r,
\]

and%
\[
\mathcal{R}_{2m+1}G\left(  a\right)  \leq\frac{d^{m+1/2}}{\left(  2m+1\right)
\left(  m!\right)  ^{2}}\left(  \max_{\left\vert \beta\right\vert
=2m+1}\left\Vert D^{\beta}G\right\Vert _{\infty,B_{r}}\right)  \left\vert
a\right\vert ^{2m+1},\quad\left\vert a\right\vert \leq r.
\]

\end{enumerate}
\end{corollary}

\begin{proof}
\textbf{Part 1} From \ref{a2.33} and \ref{a3.01}, if $\widehat{x}%
:=x/\left\vert x\right\vert $ and $\left\vert a\right\vert \leq r$,%
\begin{align*}
\mathcal{R}_{2m+1}G\left(  a\right)   & =2\int_{0}^{1}\left(  \sum
\limits_{j=0}^{m-1}\frac{\left(  -1\right)  ^{j}\left(  1-s\right)  ^{j}%
}{j!\left(  2m-j\right)  !}-\frac{s^{2m}}{\left(  2m\right)  !}\right)
\left(  \left(  aD\right)  ^{2m+1}\operatorname{Re}G\right)  \left(
sa\right)  ds\\
& \leq2\left(  \int_{0}^{1}\left\vert \sum\limits_{j=0}^{m-1}\frac{\left(
-1\right)  ^{j}\left(  1-s\right)  ^{j}}{j!\left(  2m-j\right)  !}%
-\frac{s^{2m}}{\left(  2m\right)  !}\right\vert ds\right)  \max_{s\in\left[
0,1\right]  }\left\vert \left(  aD\right)  ^{2m+1}\operatorname{Re}G\left(
sa\right)  \right\vert \\
& \leq\frac{1}{\left(  2m+1\right)  \left(  m!\right)  ^{2}}\max_{s\in\left[
0,1\right]  }\left\vert \left(  aD\right)  ^{2m+1}\operatorname{Re}G\left(
sa\right)  \right\vert \\
& =\frac{1}{\left(  2m+1\right)  \left(  m!\right)  ^{2}}\max_{s\in\left[
0,1\right]  }\left\vert \left(  \widehat{a}D\right)  ^{2m+1}\operatorname{Re}%
G\left(  sa\right)  \right\vert \left\vert a\right\vert ^{2m+1}\\
& =\frac{1}{\left(  2m+1\right)  \left(  m!\right)  ^{2}}\max_{s\in\left[
0,1\right]  }\left\vert \left(  \widehat{sa}D\right)  ^{2m+1}\operatorname{Re}%
G\left(  sa\right)  \right\vert \left\vert a\right\vert ^{2m+1}\\
& \leq\frac{1}{\left(  2m+1\right)  \left(  m!\right)  ^{2}}\left\Vert \left(
\widehat{\cdot}D\right)  ^{2m+1}\operatorname{Re}G\right\Vert _{\infty,B_{r}%
}\left\vert a\right\vert ^{2m+1}.
\end{align*}
\medskip

\textbf{Part 2} Using the expansion \ref{1.57} we get%
\begin{align*}
\frac{1}{\left(  2m+1\right)  !}\left\vert \left(  \widehat{\cdot}D\right)
^{2m+1}G\left(  x\right)  \right\vert =\left\vert \sum_{\left\vert
\beta\right\vert =2m+1}\frac{\widehat{x}^{\beta}}{\beta!}D^{\beta}G\left(
x\right)  \right\vert  & \leq\sum_{\left\vert \beta\right\vert =2m+1}%
\frac{\left\vert \widehat{x}^{\beta}\right\vert }{\beta!}\left\vert D^{\beta
}G\left(  x\right)  \right\vert \\
& \leq\sum_{\left\vert \beta\right\vert =2m+1}\frac{1}{\beta!}\left\vert
D^{\beta}G\left(  x\right)  \right\vert ,
\end{align*}

so that%
\begin{align*}
\mathcal{R}_{2m+1}G\left(  a\right)   & \leq\frac{\binom{2m}{m}}{\left(
2m+1\right)  !}\left\Vert \left(  \widehat{\cdot}D\right)  ^{2m+1}%
\operatorname{Re}G\right\Vert _{\infty,B_{r}}\left\vert a\right\vert ^{2m+1}\\
& \leq\binom{2m}{m}\left(  \sum_{\left\vert \beta\right\vert =2m+1}\frac
{1}{\beta!}\left\Vert D^{\beta}G\right\Vert _{\infty,B_{r}}\right)  \left\vert
a\right\vert ^{2m+1}.
\end{align*}

Further,%
\begin{align*}
\sum_{\left\vert \beta\right\vert =2m+1}\frac{1}{\beta!}\left\Vert D^{\beta
}G\right\Vert _{\infty,B_{r}}  & \leq\left(  \sum_{\left\vert \beta\right\vert
=2m+1}\frac{1}{\beta!}\right)  \max_{\left\vert \beta\right\vert
=2m+1}\left\Vert D^{\beta}G\right\Vert _{\infty,B_{r}}\\
& =\left(  \sum_{\left\vert \beta\right\vert =2m+1}\frac{\mathbf{1}^{2\beta}%
}{\beta!}\right)  \max_{\left\vert \beta\right\vert =2m+1}\left\Vert D^{\beta
}G\right\Vert _{\infty,B_{r}}\\
& =\frac{\left\vert \mathbf{1}\right\vert ^{2m+1}}{\left(  2m+1\right)  !}%
\max_{\left\vert \beta\right\vert =2m+1}\left\Vert D^{\beta}G\right\Vert
_{\infty,B_{r}}\\
& =\frac{d^{m+1/2}}{\left(  2m+1\right)  !}\max_{\left\vert \beta\right\vert
=2m+1}\left\Vert D^{\beta}G\right\Vert _{\infty,B_{r}},
\end{align*}

and hence%
\begin{align*}
\mathcal{R}_{2m+1}G\left(  a\right)   & \leq\frac{\binom{2m}{m}}{\left(
2m+1\right)  !}d^{m+1/2}\left(  \max_{\left\vert \beta\right\vert
=2m+1}\left\Vert D^{\beta}G\right\Vert _{\infty,B_{r}}\right)  \left\vert
a\right\vert ^{2m+1}\\
& =\frac{d^{m+1/2}}{\left(  2m+1\right)  \left(  m!\right)  ^{2}}\left(
\max_{\left\vert \beta\right\vert =2m+1}\left\Vert D^{\beta}G\right\Vert
_{\infty,B_{r}}\right)  \left\vert a\right\vert ^{2m+1}.
\end{align*}

\end{proof}

The weight function assumptions used to obtain part 4 of Theorem
\ref{Thm_Tayor_rem_estim} imply (Corollary \ref{Cor_basis_smth_W2_W3}) that
the basis function satisfies $G\in C_{B}^{\left(  \left\lfloor 2\kappa
\right\rfloor \right)  }\subset C_{B}^{\left(  2m\right)  }$ if $m\leq
\kappa<m+1/2$ and $G\in C_{B}^{\left(  \left\lfloor 2\kappa\right\rfloor
\right)  }\subset C_{B}^{\left(  2m+1\right)  }$ if $m+1/2\leq\kappa<m+1$. If
$m+1/2\leq\kappa<m+1$ then $G\in C_{B}^{\left(  2m+1\right)  }$ means the
additional assumptions of the next corollary are satisfied. However, part 1 of
Theorem \ref{Thm_Tayor_rem_estim} implies that $\left\vert a\right\vert $ has
the power $\kappa$ which may be greater that $m+\frac{1}{2}$ obtained in the
next corollary.

\begin{corollary}
\label{Cor_Thm_Tayor_rem_estim_2}Let $\left(  \mathcal{R}_{m+1}f\right)
\left(  x,a\right)  $ be the remainder of part 4 of Theorem
\ref{Thm_Tayor_rem_estim}.

\begin{enumerate}
\item If $\left\Vert \left(  \widehat{\cdot}D\right)  ^{2m+1}G\right\Vert
_{\infty,B_{r}}<\infty$ for some $0<r\leq\infty$ then%
\[
\left\vert \left(  \mathcal{R}_{m+1}f\right)  \left(  x,a\right)  \right\vert
\leq\frac{\left(  2\pi\right)  ^{\frac{d}{4}}}{\sqrt{2m+1}m!}\left\Vert
\left(  \widehat{\cdot}D\right)  ^{2m+1}G\right\Vert _{\infty,B_{r}}^{\frac
{1}{2}}\left\Vert f\right\Vert _{w,0}\left\vert a\right\vert ^{m+\frac{1}{2}%
},\quad\left\vert a\right\vert \leq r,\text{ }x\in\mathbb{R}^{d}.
\]

\item If $\max\limits_{\left\vert \beta\right\vert =2m+1}\left\Vert D^{\beta
}G\right\Vert _{\infty,B_{r}}<\infty$ for some $0<r\leq\infty$ then%
\[
\left\vert \left(  \mathcal{R}_{m+1}f\right)  \left(  x,a\right)  \right\vert
\leq\left(  2\pi\right)  ^{\frac{d}{4}}\binom{2m}{m}^{\frac{1}{2}}\left(
\sum_{\left\vert \beta\right\vert =2m+1}\frac{\left\Vert D^{\beta}G\right\Vert
_{\infty,B_{r}}}{\beta!}\right)  ^{\frac{1}{2}}\left\Vert f\right\Vert
_{w,0}\left\vert a\right\vert ^{m+\frac{1}{2}},\quad\left\vert a\right\vert
\leq r,\text{ }x\in\mathbb{R}^{d},
\]

and%
\[
\left\vert \mathcal{R}_{2m+1}\left(  a\right)  \right\vert \leq\left(
2\pi\right)  ^{\frac{d}{4}}\frac{\sqrt{d^{m+\frac{1}{2}}}}{\sqrt{2m+1}%
m!}\left(  \max_{\left\vert \beta\right\vert =2m+1}\left\Vert D^{\beta
}G\right\Vert _{\infty,B_{r}}\right)  ^{\frac{1}{2}}\left\Vert f\right\Vert
_{w,0}\left\vert a\right\vert ^{m+\frac{1}{2}},\quad\left\vert a\right\vert
\leq r,\text{ }x\in\mathbb{R}^{d}.
\]

\item If $m+1/2\leq\kappa<m+1$ then $G\in C_{B}^{\left(  2m+1\right)  }$ and%
\[
\left\vert \left(  \mathcal{R}_{m+1}f\right)  \left(  x,a\right)  \right\vert
\leq\left(  2\pi\right)  ^{\frac{d}{4}}\binom{2m}{m}^{\frac{1}{2}}\left(
\sum_{\left\vert \beta\right\vert =2m+1}\frac{\left\Vert D^{\beta}G\right\Vert
_{\infty,B_{r}}}{\beta!}\right)  ^{\frac{1}{2}}\left\Vert f\right\Vert
_{w,0}\left\vert a\right\vert ^{m+\frac{1}{2}},\quad a,x\in\mathbb{R}^{d},
\]

and%
\[
\left\vert \mathcal{R}_{2m+1}\left(  a\right)  \right\vert \leq\left(
2\pi\right)  ^{\frac{d}{4}}\frac{\sqrt{d^{m+\frac{1}{2}}}}{\sqrt{2m+1}%
m!}\left(  \max_{\left\vert \beta\right\vert =2m+1}\left\Vert D^{\beta
}G\right\Vert _{\infty,B_{r}}\right)  ^{\frac{1}{2}}\left\Vert f\right\Vert
_{w,0}\left\vert a\right\vert ^{m+\frac{1}{2}},\quad a,x\in\mathbb{R}^{d}.
\]

\end{enumerate}
\end{corollary}

\begin{proof}
\textbf{Parts 1 and 2} Substitute the estimates for $\mathcal{R}_{2m+1}\left(
a\right)  $ derived in Corollary \ref{Cor_Thm_Tayor_rem_estim_1} into the
estimate \ref{a2.30}:%
\[
\left\vert \left(  \mathcal{R}_{m+1}f\right)  \left(  x,a\right)  \right\vert
\leq\left(  2\pi\right)  ^{\frac{d}{4}}\sqrt{\mathcal{R}_{2m+1}G\left(
a\right)  }\left\Vert f\right\Vert _{w,0},\quad x,a\in\mathbb{R}^{d}.
\]

\textbf{Part 3} From Corollary \ref{Cor_basis_smth_W2_W3}, $G\in
C_{B}^{\left(  \left\lfloor 2\kappa\right\rfloor \right)  }\subset
C_{B}^{\left(  2m+1\right)  }$ so the additional assumptions of parts 1 and 2
are both satisfied for $r=\infty$.
\end{proof}

\section{Remainder estimates for data functions: radial weight function case
\label{Sect_rem_estim_radial_wt_fn}}

The radial weight functions were introduced in Subsection
\ref{SbSect_wt_func_examples} and the following result is a radial function
application of parts 1 and 4 of Theorem \ref{Thm_Tayor_rem_estim} when the
weight function has property W02.

\begin{theorem}
\label{Thm_radial_estim}\textbf{Basic radial remainder estimates \ }Suppose
the radial weight function $w$ has property W02 for $\kappa$ and that
$\overset{\circ}{w}$ is defined by \ref{a119} when $d\geq2$ and
$\overset{\circ}{w}=w$ when $d=1$. From Lemma \ref{Lem_wt_func_radial},
$\overset{\circ}{w}$ has property W02 for $\kappa$. Set
$\widehat{\overset{\circ}{G}}=1/\overset{\circ}{w}$ and $m=\left\lfloor
\kappa\right\rfloor $.

\begin{enumerate}
\item Then if $f\in X_{w}^{0}$ the estimate of part 1 of Theorem
\ref{Thm_Tayor_rem_estim} becomes:%
\begin{equation}
\left\vert \left(  \mathcal{R}_{m+1}f\right)  \left(  x,a\right)  \right\vert
\leq\frac{c_{m}}{m!}\left(  2\int_{0}^{\infty}\frac{\left\vert s\right\vert
^{2\kappa}ds}{\overset{\circ}{w}\left(  s\right)  }\right)  ^{1/2}\left\Vert
f\right\Vert _{w,0}\left\vert a\right\vert ^{\kappa},\quad x,a\in
\mathbb{R}^{d},\text{ }m=0,1,2,\ldots,\label{a223}%
\end{equation}

where $c_{m}$ is given by \ref{a2.08}.

\item If $\widehat{a}=a/\left\vert a\right\vert $ and $\left(  \widehat{\cdot
}Df\right)  \left(  a\right)  =\widehat{a}_{1}D_{1}+\ldots+\widehat{a}%
_{d}D_{d} $ then%
\begin{align*}
\left(  \widehat{\cdot}D\right)  \left(  \left\vert a\right\vert \right)   &
=1,\\
\left(  \widehat{\cdot}D\right)  g\left(  \left\vert a\right\vert \right)   &
=g^{\prime}\left(  \left\vert a\right\vert \right)  ,\\
G\left(  a\right)   & =\left(  2\pi\right)  ^{-\frac{\left(  d-1\right)  }{2}%
}\overset{\circ}{G}\left(  \left\vert a\right\vert \right)  ,\\
\left(  \widehat{\cdot}D\overset{\circ}{G}\right)  \left(  \left\vert
a\right\vert \right)   & =D\overset{\circ}{G}\left(  \left\vert a\right\vert
\right)  .
\end{align*}

\item Regarding part 4 of Theorem \ref{Thm_Tayor_rem_estim}: if $\left\Vert
D^{2m+1}\overset{\circ}{G}\right\Vert _{\infty,B_{r}}<\infty$ for some $r>0$
then%
\[
\left\vert \left(  \mathcal{R}_{m+1}f\right)  \left(  x,a\right)  \right\vert
\leq\frac{\left(  2\pi\right)  ^{\frac{1}{4}}}{\sqrt{2m+1}m!}\left\Vert
D^{2m+1}\overset{\circ}{G}\right\Vert _{\infty,B_{r}}^{\frac{1}{2}}\left\Vert
f\right\Vert _{w,0}\left\vert a\right\vert ^{m+\frac{1}{2}},\quad\left\vert
a\right\vert \leq r,\text{ }x\in\mathbb{R}^{d}.
\]

\end{enumerate}
\end{theorem}

\begin{proof}
\textbf{Part 1} Since the weight function is radial, when $d>1$ we can apply
Lemma \ref{Lem_wt_func_radial} with $u\left(  s\right)  =\left\vert
s\right\vert ^{2\kappa}$ to the estimate \ref{a2.06} of Theorem
\ref{Thm_Tayor_rem_estim}. In fact%
\[
\int_{\mathbb{R}^{d}}\frac{u\left(  \widehat{x}\xi\right)  }{w\left(
\xi\right)  }d\xi=\int_{\mathbb{R}^{1}}\frac{u\left(  s\right)  }%
{\overset{\circ}{w}\left(  s\right)  }ds=\int_{0}^{\infty}\frac{u\left(
s\right)  +u\left(  -s\right)  }{\overset{\circ}{w}\left(  s\right)  }%
ds=2\int_{0}^{\infty}\frac{\left\vert s\right\vert ^{2\kappa}}{\overset{\circ
}{w}\left(  s\right)  }ds,
\]

so \ref{a2.06} becomes%
\[
\left\vert \left(  \mathcal{R}_{m+1}f\right)  \left(  x,a\right)  \right\vert
\leq\frac{c_{m}}{m!}\left(  2\int_{0}^{\infty}\frac{\left\vert s\right\vert
^{2\kappa}}{\overset{\circ}{w}\left(  s\right)  }ds\right)  ^{1/2}\left\Vert
f\right\Vert _{w,0}\left\vert a\right\vert ^{\kappa},\quad x,a\in
\mathbb{R}^{d}.
\]
\medskip

\textbf{Part 2} We have
\[
\left(  \widehat{\cdot}D\right)  \left(  \left\vert a\right\vert \right)
=\sum\widehat{a}_{i}D_{i}\left\vert a\right\vert =\sum\widehat{a}_{i}%
a_{i}\left\vert a\right\vert ^{-1}=\sum\widehat{a}_{i}\widehat{a}_{i}=1,
\]
and
\[
\left(  \widehat{\cdot}D\right)  g\left(  \left\vert a\right\vert \right)
=\sum\widehat{a}_{i}D_{i}g\left(  \left\vert a\right\vert \right)
=\sum\widehat{a}_{i}g^{\prime}\left(  \left\vert a\right\vert \right)
D_{i}\left\vert a\right\vert =\sum\widehat{a}_{i}g^{\prime}\left(  \left\vert
a\right\vert \right)  a_{i}\left\vert a\right\vert ^{-1}=g^{\prime}\left(
\left\vert a\right\vert \right)  .
\]

From Theorem \ref{Thm_basis_radial}, $G\left(  a\right)  =\left(  2\pi\right)
^{-\frac{d}{2}}\int_{\mathbb{R}^{1}}\frac{e^{i\left\vert a\right\vert s}%
}{\overset{\circ}{w}\left(  s\right)  }ds$ so that%
\[
G\left(  a\right)  =\left(  2\pi\right)  ^{-\frac{d}{2}}\int%
\limits_{\mathbb{R}^{1}}\frac{e^{i\left\vert a\right\vert s}}{\overset{\circ
}{w}\left(  s\right)  }ds=\left(  2\pi\right)  ^{-\frac{d}{2}+\frac{1}{2}%
}\left(  2\pi\right)  ^{-\frac{1}{2}}\int\limits_{\mathbb{R}^{1}%
}e^{i\left\vert a\right\vert t}\widehat{\overset{\circ}{G}}\left(  t\right)
dt=\left(  2\pi\right)  ^{-\frac{\left(  d-1\right)  }{2}}\overset{\circ
}{G}\left(  \left\vert a\right\vert \right)  .
\]

and the other two equations follow from the standard identities%
\[
\left(  aD\overset{\circ}{G}\right)  \left(  \left\vert a\right\vert \right)
=\left\vert a\right\vert D\overset{\circ}{G}\left(  \left\vert a\right\vert
\right)  .
\]
\medskip

\textbf{Part 3} From part 4 of Theorem \ref{Thm_Tayor_rem_estim},
\[
\left\vert \left(  \mathcal{R}_{m+1}f\right)  \left(  x,a\right)  \right\vert
\leq\left(  2\pi\right)  ^{\frac{d}{4}}\sqrt{\mathcal{R}_{2m+1}G\left(
a\right)  }\left\Vert f\right\Vert _{w,0},
\]

and from part 1 Corollary \ref{Cor_Thm_Tayor_rem_estim_1},%
\[
\mathcal{R}_{2m+1}G\left(  a\right)  \leq\frac{1}{\left(  2m+1\right)  \left(
m!\right)  ^{2}}\left\Vert \left(  \widehat{\cdot}D\right)  ^{2m+1}%
G\right\Vert _{\infty,B_{r}}\left\vert a\right\vert ^{2m+1},\quad\left\vert
a\right\vert \leq r,\text{ }x\in\mathbb{R}^{d},
\]

but from part 2%
\[
\left(  \widehat{\cdot}D\right)  ^{2m+1}G\left(  x\right)  =\left(
2\pi\right)  ^{-\frac{\left(  d-1\right)  }{2}}\left(  \widehat{\cdot
}D\right)  ^{2m+1}\left(  \overset{\circ}{G}\left(  \left\vert x\right\vert
\right)  \right)  =\left(  2\pi\right)  ^{-\frac{\left(  d-1\right)  }{2}%
}\left(  D^{2m+1}\overset{\circ}{G}\right)  \left(  \left\vert x\right\vert
\right)  ,
\]

so that%
\begin{align*}
\left\vert \left(  \mathcal{R}_{m+1}f\right)  \left(  x,a\right)  \right\vert
& \leq\left(  2\pi\right)  ^{\frac{d}{4}}\left(  \mathcal{R}_{2m+1}\left(
a\right)  \right)  ^{\frac{1}{2}}\left\Vert f\right\Vert _{w,0}\\
& \leq\frac{\left(  2\pi\right)  ^{\frac{d}{4}}}{\sqrt{2m+1}m!}\left\Vert
\left(  \widehat{\cdot}D\right)  ^{2m+1}G\right\Vert _{\infty,B_{r}}^{\frac
{1}{2}}\left\Vert f\right\Vert _{w,0}\left\vert a\right\vert ^{m+\frac{1}{2}%
}\\
& =\frac{\left(  2\pi\right)  ^{\frac{d}{4}}}{\sqrt{2m+1}m!}\left\Vert \left(
2\pi\right)  ^{-\frac{\left(  d-1\right)  }{2}}\left(  D^{2m+1}\overset{\circ
}{G}\right)  \left(  \left\vert \cdot\right\vert \right)  \right\Vert
_{\infty,B_{r}}^{\frac{1}{2}}\left\Vert f\right\Vert _{w,0}\left\vert
a\right\vert ^{m+\frac{1}{2}}\\
& =\frac{\left(  2\pi\right)  ^{\frac{1}{4}}}{\sqrt{2m+1}m!}\left\Vert
D^{2m+1}\overset{\circ}{G}\right\Vert _{\infty,B_{r}}^{\frac{1}{2}}\left\Vert
f\right\Vert _{w,0}\left\vert a\right\vert ^{m+\frac{1}{2}}.
\end{align*}

\end{proof}

\section{Remainder estimates for data functions: tensor product weight
function case\label{Sect_rem_estim_ten_prod}}

In this section we use part 2 of Corollary \ref{Cor_Thm_Tayor_rem_estim_2} to
estimate the remainder $\left\vert \left(  \mathcal{R}_{m+1}f\right)  \left(
x,a\right)  \right\vert $ in terms of bounds on the derivatives of the
univariant basis function $G_{1}$ To do this we need to impose a bound on
$D^{2m+1}G_{1}$ near zero. It is natural to use the weight function property
W03 and the fact that if a weight function has property W03 for smoothness
parameter $\kappa$ then it has property W02 smoothness parameter
$\underline{\kappa}$.

\begin{theorem}
\label{Thm_rem_estim_ten_prod}\textbf{Data function remainder estimates:
tensor product\ weight function }Suppose the tensor product weight function
$w\left(  x\right)  =w_{1}\left(  x_{1}\right)  w_{1}\left(  x_{2}\right)
\ldots w_{1}\left(  x_{d}\right)  $ has property W03 for $\kappa$. Set
$m=\left\lfloor \underline{\kappa}\right\rfloor $. Suppose $G\left(  x\right)
=G_{1}\left(  x_{1}\right)  \ldots G_{1}\left(  x_{d}\right)  $ is the
corresponding tensor product basis function. Then $G\in C_{B}^{\left(
\left\lfloor 2\kappa\right\rfloor \right)  }\left(  \mathbb{R}^{d}\right)
\subset C_{B}^{\left(  2m\mathbf{1}\right)  }\left(  \mathbb{R}^{d}\right)  $
and $G_{1}\in C_{B}^{\left(  2m\right)  }\left(  \mathbb{R}^{1}\right)  $.

Now suppose \textbf{in addition} we assume that $\left\Vert D^{2m+1}%
G_{1}\right\Vert _{\infty,B_{r}}<\infty$ for some $r>0$.

Then for all $f\in X_{w}^{0}$:%
\begin{equation}
\left\vert \left(  \mathcal{R}_{m+1}f\right)  \left(  x,a\right)  \right\vert
\leq\left(  2\pi\right)  ^{\frac{d}{4}}\binom{2m}{m}^{\frac{1}{2}}\left(
\sum_{j\leq2m+1}\frac{\left\Vert D^{j}G_{1}\right\Vert _{\infty,B_{r}}}%
{j!}\right)  ^{\frac{d}{2}}\left\Vert f\right\Vert _{w,0}\left\vert
a\right\vert ^{m+\frac{1}{2}},\quad\left\vert a\right\vert \leq r,\text{ }%
x\in\mathbb{R}^{d},\label{2.04}%
\end{equation}

and%
\begin{equation}
\left\vert \left(  \mathcal{R}_{m+1}f\right)  \left(  x,a\right)  \right\vert
\leq\frac{\left(  2\pi\right)  ^{\frac{d}{4}}\sqrt{d^{m+\frac{1}{2}}}}%
{\sqrt{2m+1}m!}\left(  \max\limits_{j\leq2m+1}\left\Vert D^{j}G_{1}\right\Vert
_{\infty,B_{r}}\right)  ^{\frac{d}{2}}\left\Vert f\right\Vert _{w,0}\left\vert
a\right\vert ^{m+\frac{1}{2}},\quad\left\vert a\right\vert \leq r,\text{ }%
x\in\mathbb{R}^{d}.\label{2.05}%
\end{equation}

\end{theorem}

\begin{proof}
Since $G\in C_{B}^{\left(  \left\lfloor 2\kappa\right\rfloor \right)  }\left(
\mathbb{R}^{d}\right)  $ we have $G\in C_{B}^{\left(  2m\mathbf{1}\right)
}\left(  \mathbb{R}^{d}\right)  $ and so $\left\Vert D^{\beta}G\right\Vert
_{\infty,B_{r}}<\infty$ when $\beta\leq\left(  2m+1\right)  \mathbf{1}$. Hence
$\left\Vert D^{\beta}G\right\Vert _{\infty,B_{r}}<\infty$ when $\left\vert
\beta\right\vert =2m+1$ and thus from part 2 of Corollary
\ref{Cor_Thm_Tayor_rem_estim_2} we have the two estimates
\begin{equation}
\left\vert \left(  \mathcal{R}_{m+1}f\right)  \left(  x,a\right)  \right\vert
\leq\left(  2\pi\right)  ^{\frac{d}{4}}\binom{2m}{m}^{\frac{1}{2}}\left(
\sum_{\left\vert \beta\right\vert =2m+1}\frac{\left\Vert D^{\beta}G\right\Vert
_{\infty,B_{r}}}{\beta!}\right)  ^{\frac{1}{2}}\left\Vert f\right\Vert
_{w,0}\left\vert a\right\vert ^{m+\frac{1}{2}},\label{2.01}%
\end{equation}

and%
\begin{equation}
\left\vert \left(  \mathcal{R}_{m+1}f\right)  \left(  x,a\right)  \right\vert
\leq\left(  2\pi\right)  ^{\frac{d}{4}}\frac{\sqrt{d^{m+\frac{1}{2}}}}%
{\sqrt{2m+1}m!}\left(  \max_{\left\vert \beta\right\vert =2m+1}\left\Vert
D^{\beta}G\right\Vert _{\infty,B_{r}}\right)  ^{\frac{1}{2}}\left\Vert
f\right\Vert _{w,0}\left\vert a\right\vert ^{m+\frac{1}{2}},\label{2.02}%
\end{equation}

which are valid for $\left\vert a\right\vert \leq r$ and $x\in\mathbb{R}^{d}
$. Now%
\begin{align*}
\left\Vert D^{\beta}G\right\Vert _{\infty,B_{r}}  & =\left\Vert D_{1}%
^{\beta_{1}}G_{1}\text{\thinspace}D_{2}^{\beta_{2}}G_{1}\ldots D_{d}%
^{\beta_{d}}G_{1}\right\Vert _{\infty,B_{r}}\\
& \leq\left\Vert D^{\beta_{1}}G_{1}\right\Vert _{\infty,B_{r}}\left\Vert
D^{\beta_{2}}G_{1}\right\Vert _{\infty,B_{r}}\ldots\left\Vert D^{\beta_{d}%
}G_{1}\right\Vert _{\infty,B_{r}},
\end{align*}

so that\ (note the first inequality)%
\begin{align*}
&  \sum_{\left\vert \beta\right\vert =2m+1}\frac{1}{\beta!}\left\Vert
D^{\beta}G\right\Vert _{\infty,B_{r}}\\
&  \leq\sum_{\beta\leq\left(  2m+1\right)  \mathbf{1}}\frac{1}{\beta
!}\left\Vert D^{\beta}G\right\Vert _{\infty,B_{r}}\\
&  \leq\sum_{\beta\leq\left(  2m+1\right)  \mathbf{1}}\frac{1}{\beta_{1}%
!}\left\Vert D^{\beta_{1}}G_{1}\right\Vert _{\infty,B_{r}}\frac{1}{\beta_{2}%
!}\left\Vert D^{\beta_{2}}G_{1}\right\Vert _{\infty,B_{r}}\ldots\frac{1}%
{\beta_{d}!}\left\Vert D^{\beta_{d}}G_{1}\right\Vert _{\infty,B_{r}}\\
&  =\sum_{\beta_{1}\leq2m+1}\frac{1}{\beta_{1}!}\left\Vert D^{\beta_{1}}%
G_{1}\right\Vert _{\infty,B_{r}}\sum_{\beta_{2}\leq2m+1}\frac{1}{\beta_{2}%
!}\left\Vert D^{\beta_{2}}G_{1}\right\Vert _{\infty,B_{r}}\ldots\sum
_{\beta_{d}\leq2m+1}\frac{1}{\beta_{d}!}\left\Vert D^{\beta_{d}}%
G_{1}\right\Vert _{\infty,B_{r}}\\
&  =\left(  \sum_{j\leq2m+1}\frac{1}{j!}\left\Vert D^{j}G_{1}\right\Vert
_{\infty,B_{r}}\right)  ^{d},
\end{align*}

and substitution into \ref{2.01} gives \ref{2.04}:%
\[
\left\vert \left(  \mathcal{R}_{m+1}f\right)  \left(  x,a\right)  \right\vert
\leq\left(  2\pi\right)  ^{\frac{d}{4}}\binom{2m}{m}^{\frac{1}{2}}\left(
\sum_{j\leq2m+1}\frac{1}{j!}\left\Vert D^{j}G_{1}\right\Vert _{\infty,B_{r}%
}\right)  ^{\frac{d}{2}}\left\Vert f\right\Vert _{w,0}\left\vert a\right\vert
^{m+\frac{1}{2}}.
\]

Regarding \ref{2.02}:%
\begin{align*}
\max\limits_{\left\vert \beta\right\vert =2m+1} &  \left\Vert D^{\beta
}G\right\Vert _{\infty,B_{r}}\\
&  \leq\max\limits_{\beta\leq\left(  2m+1\right)  \mathbf{1}}\left\Vert
D^{\beta}G\right\Vert _{\infty,B_{r}}\\
&  \leq\max\limits_{\beta\leq\left(  2m+1\right)  \mathbf{1}}\left\Vert
D^{\beta_{1}}G_{1}\right\Vert _{\infty,B_{r}}\left\Vert D^{\beta_{2}}%
G_{1}\right\Vert _{\infty,B_{r}}\ldots\left\Vert D^{\beta_{d}}G_{1}\right\Vert
_{\infty,B_{r}}\\
&  =\max\limits_{\beta_{d}\leq2m+1}\ldots\max\limits_{\beta_{2}\leq2m+1}%
\max\limits_{\beta_{1}\leq2m+1}\left\Vert D^{\beta_{1}}G_{1}\right\Vert
_{\infty,B_{r}}\left\Vert D^{\beta_{2}}G_{1}\right\Vert _{\infty,B_{r}}%
\ldots\left\Vert D^{\beta_{d}}G_{1}\right\Vert _{\infty,B_{r}}\\
&  =\left(  \max\limits_{j\leq2m+1}\left\Vert D^{j}G_{1}\right\Vert
_{\infty,B_{r}}\right)  \left(  \max\limits_{j\leq2m+1}\left\Vert D^{j}%
G_{1}\right\Vert _{\infty,B_{r}}\right)  \ldots\left(  \max\limits_{j\leq
2m+1}\left\Vert D^{j}G_{1}\right\Vert _{\infty,B_{r}}\right) \\
&  =\left(  \max\limits_{j\leq2m+1}\left\Vert D^{j}G_{1}\right\Vert
_{\infty,B_{r}}\right)  ^{d},
\end{align*}

and substitution into \ref{2.02} gives \ref{2.05}.
\end{proof}

The next results show that well-known reproducing kernel/Riesz representer
techniques can be used to globally bound the derivatives $\left\{  D^{k}%
G_{1}\right\}  _{k=1}^{2m}$ in terms of $\left\{  D^{k}G_{1}\left(  0\right)
\right\}  _{k=1}^{2m}$.

\begin{lemma}
\label{Lem_Thm_estim_max_DjG1}Suppose $w_{1}$ is a \textbf{1-dimensional}
weight function with smoothness parameter $\kappa$ and the basis function
$G_{1}$ is real-valued. Set $m=\left\lfloor \kappa\right\rfloor $ . Then for
$k\leq2m$ and $x\in\mathbb{R}^{1}$,%
\begin{equation}
\left\vert D^{k}G_{1}\left(  x\right)  \right\vert \leq\left\{
\begin{array}
[c]{ll}%
\left(  -1\right)  ^{\left(  k+1\right)  /2}\left(  D^{k+1}G_{1}\left(
0\right)  \right)  \,\left\vert x\right\vert , & k\text{ }odd,\\
\left(  -1\right)  ^{k/2}D^{k}G_{1}\left(  0\right)  , & k\text{ }even.
\end{array}
\right. \label{a3.02}%
\end{equation}

\end{lemma}

\begin{proof}
From Corollary \ref{Cor_basis_smth_W2_W3}, $G_{1}\in C_{B}^{\left(
\left\lfloor 2\kappa\right\rfloor \right)  }\left(  \mathbb{R}^{1}\right)
\subset C_{B}^{\left(  2m\right)  }\left(  \mathbb{R}^{1}\right)  $ and from
part 2 of Corollary \ref{Cor_bound_DG(z-x)-DG(z-y)_W3},\smallskip

\fbox{so when $k$ is even}%
\begin{equation}
\left\vert D^{k}G_{1}\left(  x\right)  \right\vert \leq\left(  -1\right)
^{k/2}D^{k}G_{1}\left(  0\right)  ,\quad k\text{ }even,\text{ }k\leq
2\kappa,\label{a2.34}%
\end{equation}

and since $\left\lfloor 2\kappa\right\rfloor \leq2m+1$ it follows that
$k\leq2m$.\smallskip

\fbox{When $k$ is odd} $k\leq2m-1$, $k+1\leq2m$ and so $D^{k+1}G_{1}\in
C_{B}^{\left(  0\right)  }\left(  \mathbb{R}^{1}\right)  $. Hence by the mean
value theorem there exists a function $t:\mathbb{R}^{1}\rightarrow\left(
0,1\right)  $ such that%
\[
D^{k}G_{1}\left(  x\right)  =D^{k}G_{1}\left(  0\right)  +x\left(
D^{k+1}G_{1}\right)  \left(  t\left(  x\right)  x\right)  .
\]

But $G$ is even so $D^{k}G_{1}\left(  0\right)  =0$ which means that
$D^{k}G_{1}\left(  x\right)  =\left(  D^{k+1}G_{1}\right)  \left(  t\left(
x\right)  x\right)  \,x$ and by \ref{a2.34},
\[
\left\vert D^{k}G_{1}\left(  x\right)  \right\vert \leq\left(  -1\right)
^{\left(  k+1\right)  /2}\left(  D^{k+1}G_{1}\left(  0\right)  \right)
\,\left\vert x\right\vert ,\quad k\text{ }odd,\text{ }k\leq2\kappa.
\]

\end{proof}

Using Lemma \ref{Lem_Thm_estim_max_DjG1} we now calculate upper bounds for the
terms in Theorem \ref{Thm_rem_estim_ten_prod} which involve the derivatives of
the one dimensional basis function of order at most $2m$.

\begin{theorem}
\label{Thm_estim_max_DjG1}Suppose $w$ is a 1-dimensional weight function with
smoothness parameter $\kappa$ and \textbf{real} basis function $G_{1}$. Set
$m=\left\lfloor \kappa\right\rfloor $. Then regarding \ref{2.05},%
\[
\max\limits_{j=0}^{2m}\left\Vert D^{j}G_{1}\right\Vert _{\infty,B_{r}}\leq
\max\left\{  r\max\limits_{k=1}^{m}\left\vert D^{2k}G_{1}\left(  0\right)
\right\vert ,\max\limits_{k=0}^{m}\left\vert D^{2k}G_{1}\left(  0\right)
\right\vert \right\}  ,
\]

and regarding \ref{2.04},%
\[
\sum_{j=0}^{2m}\frac{1}{j!}\left\Vert D^{j}G_{1}\right\Vert _{\infty,B_{r}%
}\leq G_{1}\left(  0\right)  +\sum_{k=1}^{m}\left(  -1\right)  ^{k}\left(
\frac{1}{\left(  2k\right)  !}+\frac{r}{\left(  2k-1\right)  !}\right)
D^{2k}G_{1}\left(  0\right)  .
\]

\end{theorem}

\begin{proof}
Since $G_{1}$ is real, $G_{1}$ is even and so $D^{j}G_{1}\left(  0\right)  =0$
when $j=1,3,5,\ldots,2m-1$. Thus using Lemma \ref{Lem_Thm_estim_max_DjG1} we
get%
\begin{align*}
\max\limits_{j=0}^{2m}\left\Vert D^{j}G_{1}\right\Vert _{\infty,B_{r}}  &
=\max\left\{  \max\limits_{\substack{j=1 \\j\text{ }odd}}^{2m-1}\left\Vert
D^{j}G_{1}\right\Vert _{\infty,B_{r}},\max\limits_{\substack{j=0 \\j\text{
}even}}^{2m}\left\Vert D^{j}G_{1}\right\Vert _{\infty,B_{r}}\right\} \\
& =\max\left\{  \max\limits_{k=0}^{m-1}\left\Vert D^{2k+1}G_{1}\right\Vert
_{\infty,B_{r}},\max\limits_{k=0}^{m}\left\Vert D^{2k}G_{1}\right\Vert
_{\infty,B_{r}}\right\} \\
& \leq\max\left\{  \max\limits_{k=0}^{m-1}\left\vert D^{2k+2}G_{1}\left(
0\right)  \right\vert r,\max\limits_{k=0}^{m}\left\vert D^{2k}G_{1}\left(
0\right)  \right\vert \right\} \\
& =\max\left\{  r\max\limits_{k=1}^{m}\left\vert D^{2k}G_{1}\left(  0\right)
\right\vert ,\max\limits_{k=0}^{m}\left\vert D^{2k}G_{1}\left(  0\right)
\right\vert \right\}  ,
\end{align*}

and%
\begin{align*}
\sum_{j=0}^{2m} &  \frac{\left\Vert D^{j}G_{1}\right\Vert _{\infty,B_{r}}}%
{j!}\\
&  \leq\sum_{\substack{j=0 \\j\text{ }even}}^{2m}\frac{1}{j!}\left\Vert
D^{j}G_{1}\right\Vert _{\infty,B_{r}}+\sum_{\substack{j=1 \\j\text{ }odd
}}^{2m-1}\frac{1}{j!}\left\Vert D^{j}G_{1}\right\Vert _{\infty,B_{r}}\\
&  =\sum_{k=0}^{m}\frac{1}{\left(  2k\right)  !}\left\Vert D^{2k}%
G_{1}\right\Vert _{\infty,B_{r}}+\sum_{k=0}^{m-1}\frac{1}{\left(  2k+1\right)
!}\left\Vert D^{2k+1}G_{1}\right\Vert _{\infty,B_{r}}\\
&  \leq\sum_{k=0}^{m}\frac{\left(  -1\right)  ^{k}}{\left(  2k\right)
!}D^{2k}G_{1}\left(  0\right)  +r\sum_{k=0}^{m-1}\frac{\left(  -1\right)
^{k+1}}{\left(  2k+1\right)  !}D^{2k+2}G_{1}\left(  0\right) \\
&  =\sum_{k=0}^{m}\frac{\left(  -1\right)  ^{k}}{\left(  2k\right)  !}%
D^{2k}G_{1}\left(  0\right)  +r\sum_{k=1}^{m}\frac{\left(  -1\right)  ^{k}%
}{\left(  2k-1\right)  !}D^{2k}G_{1}\left(  0\right) \\
&  =G_{1}\left(  0\right)  +\sum_{k=1}^{m}\left(  -1\right)  ^{k}\left(
\frac{1}{\left(  2k\right)  !}+\frac{r}{\left(  2k-1\right)  !}\right)
D^{2k}G_{1}\left(  0\right)  .
\end{align*}

\end{proof}

\section{Interpolation error: $w\in W02$ for $\kappa<1$ or $w\in W03 $ for
$\protect\underline{\kappa}<1$}

This section improves on the interpolation error results of Section
\ref{Sect_interp_no_Taylor_converg} of the previous chapter. We want to
estimate the data function remainder $\left(  \mathcal{R}_{m+1}f\right)
\left(  x,a\right)  $ so that the unisolvency data assumption is not required
to remove the polynomial terms of the Taylor series. This implies $m=0$ and we
will use the remainder estimates of parts 1 and 4 of Theorem
\ref{Thm_Tayor_rem_estim}.

The remainder estimates of Theorem \ref{Thm_|f(x)-f(y)|_inequal_2} motivate
the form of the estimate \ref{a335} in the next theorem. This result is our
main interpolation convergence estimate for the case when the weight function
has property W02 for $\kappa<1$ or W03 for $\underline{\kappa}<1$. Whereas in
Section \ref{Sect_interp_no_Taylor_converg} we had two interpolant convergence
estimates i.e. Theorems \ref{Thm_interpol_error_in_terms_of_wt_fn} and
\ref{Thm_|f(x)-f(y)|_inequal_2}, here we need only one.

The results are summarized in Table \ref{Tbl_NonUnisolvConverg_revisit}.

\begin{theorem}
\label{Thm_converg_interpol_not_unisolv2}\textbf{Interpolant convergence
for\ }$\left\{  w\in W02:\kappa<1\right\}  $ Suppose $w\in W02$ for some
$\kappa<1$. Suppose all the independent data sets $X$ are contained in a
closed bounded infinite set $K$.

Assume that for some $s\geq0$ and constants $c_{w},h_{w}>0$ the Taylor series
data function remainder estimate
\begin{equation}
\left\vert f\left(  x+a\right)  -f\left(  x\right)  \right\vert \leq
c_{w}\left\Vert f\right\Vert _{w,0}\left\vert a\right\vert ^{s},\text{\quad
}\left\vert a\right\vert <h_{w},\text{ }x\in K,\label{a335}%
\end{equation}

is valid for all $f\in X_{w}^{0}$.

Let $\mathcal{I}_{X}$ be the minimal norm interpolant mapping for the
independent data set $X$. Then for any data function $f\in X_{w}^{0}$ it
follows that%
\begin{equation}
\left\vert f\left(  x\right)  -\mathcal{I}_{X}f\left(  x\right)  \right\vert
\leq c_{w}\sqrt{\left(  f-\mathcal{I}_{X}f,f\right)  _{w,0}}\left(
h_{X,K}\right)  ^{s}\leq c_{w}\left\Vert f\right\Vert _{w,0}\left(
h_{X,K}\right)  ^{s},\text{\quad}x\in K,\label{a54.3}%
\end{equation}

when $h_{X,K}=\sup\limits_{x\in K}\operatorname*{dist}\left(  x,X\right)
<h_{w}$ i.e. the order of convergence is at least $s$.

Further, for all data and any data function $f\in X_{w}^{0}$ we have%
\begin{equation}
\left\vert f\left(  x\right)  -\mathcal{I}_{X}f\left(  x\right)  \right\vert
\leq c_{w}\sqrt{\left(  f-\mathcal{I}_{X}f,f\right)  _{w,0}}\left(
\operatorname*{diam}K\right)  ^{s},\text{\quad}x\in K.\label{a54.8}%
\end{equation}

\end{theorem}

\begin{proof}
From Theorem \ref{Thm_interp_error_const_bound}, $\sqrt{\left(  f-\mathcal{I}%
_{X}f,f\right)  _{w,0}}\leq\left\Vert f\right\Vert _{w,0}$. Now fix $x\in K$
and let $X=\left\{  x^{\left(  j\right)  }\right\}  _{j=1}^{N}\subset K$ be an
independent data set. Using the fact that $\mathcal{I}_{X}f$ interpolates $f$
on $X$ we can use \ref{a335} to obtain
\begin{align*}
\left\vert f\left(  x\right)  -\mathcal{I}_{X}f\left(  x\right)  \right\vert
=\left\vert \left(  f-\mathcal{I}_{X}f\right)  \left(  x\right)  -\left(
f-\mathcal{I}_{X}f\right)  \left(  x^{\left(  j\right)  }\right)  \right\vert
&  \leq c_{w}\left\Vert f-\mathcal{I}_{X}f\right\Vert _{w,0}\left\vert
x-x^{\left(  j\right)  }\right\vert ^{s}\\
&  =c_{w}\left(  f-\mathcal{I}_{X}f,f\right)  _{w,0}\left\vert x-x^{\left(
j\right)  }\right\vert ^{s},
\end{align*}

because $\max\limits_{j}\left\vert x-x^{\left(  j\right)  }\right\vert
\leq\sup\limits_{x\in K}\operatorname*{dist}\left(  x,X\right)  <h_{w}$. Note
that the last step used the fact that $\mathcal{I}_{X}$ is a self-adjoint
projection. Then, since $\operatorname*{dist}\left(  x,X\right)  <h_{w}$, we
can apply the upper bound \ref{1.07} and obtain
\[
\left\vert f\left(  x\right)  -\mathcal{I}_{X}f\left(  x\right)  \right\vert
\leq\sqrt{c_{w}}\sqrt{\left(  f-\mathcal{I}_{X}f,f\right)  _{w,0}}\left\vert
x-x^{\left(  j\right)  }\right\vert ^{s},
\]

for all $j$ and so%
\begin{align*}
\left\vert f\left(  x\right)  -\mathcal{I}_{X}f\left(  x\right)  \right\vert
&  \leq c_{w}\sqrt{\left(  f-\mathcal{I}_{X}f,f\right)  _{w,0}}\left(
\operatorname*{dist}\left(  x,X\right)  \right)  ^{s}\\
&  \leq c_{w}\sqrt{\left(  f-\mathcal{I}_{X}f,f\right)  _{w,0}}\left(
\sup\limits_{x\in K}\operatorname*{dist}\left(  x,X\right)  \right)  ^{s}\\
&  =c_{w}\sqrt{\left(  f-\mathcal{I}_{X}f,f\right)  _{w,0}}\left(
h_{X,K}\right)  ^{s},
\end{align*}

where the last step used inequalities \ref{1.09}. Since $\sqrt{\left(
f-\mathcal{I}_{X}f,f\right)  _{w,0}}\leq\left\Vert f\right\Vert _{w,0}$ the
order of convergence is at least $s$. Finally, since $h_{X,K}\leq
\operatorname*{diam}K$, \ref{a54.8} follows directly from \ref{a54.3}.
\end{proof}

\begin{remark}
\label{Rem_Thm_converg_interpol_not_unisolv2}The remainder estimate
\ref{a2.06} implies that the order of the interpolant convergence $s$ is at
least $\kappa$.
\end{remark}

\begin{example}
\textbf{Sobolev splines} From Example \ref{Ex_Sobol_splin_k_ge_1}, $w\left(
\xi\right)  =\left(  1+\left\vert \xi\right\vert ^{2}\right)  ^{v}$ has
property W01 for empty $\mathcal{A}$ and property W02 for $\kappa$ iff
$0\leq\kappa<v-d/2$. Here we assume that $\kappa<1$ so $m=0$ and
$0<v-d/2\leq1$. Comparison of \ref{a335} with \ref{a223} yields:
\begin{equation}
c_{w}=2\sqrt{2}\left(  \int\nolimits_{0}^{\infty}\frac{r^{2\kappa}%
}{\overset{\circ}{w}\left(  r\right)  }dr\right)  ^{1/2},\text{\quad}%
h_{w}=\infty,\text{\ }s=\kappa<v-d/2,\label{a9.6}%
\end{equation}

and the integral can be calculated using \ref{a54.2} below.
\end{example}

\begin{example}
\textbf{Extended B-splines with parameter }$\mathbf{n=1}$ From Theorem
\ref{Thm_basis_tensor_hat_W3}, and the weight function has property W03 for
$\kappa<\left(  n-1/2\right)  \mathbf{1}$ and hence property W02 for
$\underline{\kappa}<\left(  n-1/2\right)  $. Hence when $n=1$, $\kappa
<\mathbf{1}/2$ and $m=\left\lfloor \underline{\kappa}\right\rfloor =0$. From
part 1 of Theorem \ref{Thm_Tayor_rem_estim} and \ref{a335} the order of
convergence is at least $\underline{\kappa}$ which is less than the previous
convergence estimates of $1/2$ e.g. Type 1 estimates of Table
\ref{Tbl_NonUnisolvTyp1Converg} and Theorem \ref{Thm_|f(x)-f(y)|_inequal_2}.

However, from part 3 of Theorem \ref{Thm_Tayor_rem_estim},%
\[
\left\vert \left(  \mathcal{R}_{1}f\right)  \left(  x,a\right)  \right\vert
\leq\left(  2\pi\right)  ^{\frac{d}{4}}\sqrt{2}\left(  G\left(  0\right)
-G\left(  a\right)  \right)  ^{\frac{1}{2}}\left\Vert f\right\Vert
_{w,0},\quad x,a\in\mathbb{R}^{d},
\]

and from Theorem \ref{Thm_ex_nat_spline_basis_Lipschitz},%
\[
\left\vert G\left(  x\right)  -G\left(  y\right)  \right\vert \leq\sqrt
{d}G_{1}\left(  0\right)  ^{d-1}\left\Vert DG_{1}\right\Vert _{\infty
}\left\vert x-y\right\vert ,\quad x,y\in\mathbb{R}^{d},
\]

so that%
\begin{align*}
\left\vert \left(  \mathcal{R}_{1}f\right)  \left(  x,a\right)  \right\vert  &
\leq\left(  2\pi\right)  ^{\frac{d}{4}}\sqrt{2}\left(  \sqrt{d}G_{1}\left(
0\right)  ^{d-1}\left\Vert DG_{1}\right\Vert _{\infty}\left\vert a\right\vert
\right)  ^{\frac{1}{2}}\left\Vert f\right\Vert _{w,0}\\
& =\left(  2\left(  2\pi\right)  ^{\frac{d}{2}}\sqrt{d}G_{1}\left(  0\right)
^{d-1}\left\Vert DG_{1}\right\Vert _{\infty}\right)  ^{\frac{1}{2}}\left\Vert
f\right\Vert _{w,0}\left\vert a\right\vert ^{1/2},
\end{align*}

i.e. the order of convergence is $1/2$, and comparison with \ref{a335}
yields:
\begin{equation}
c_{w}=\left(  2\left(  2\pi\right)  ^{\frac{d}{2}}\sqrt{d}G_{1}\left(
0\right)  ^{d-1}\left\Vert DG_{1}\right\Vert _{\infty}\right)  ^{\frac{1}{2}%
},\text{\quad}h_{w}=\infty,\text{\ }s=1/2.\label{a9.7}%
\end{equation}

\end{example}

\subsection{Examples: summary table}%

\begin{table}[htbp] \centering
$%
\begin{tabular}
[c]{|c|c||c|c|c|}\hline
\multicolumn{5}{|c|}{Interpolant error: $w\in W02$ for $\kappa<1$ or $w\in
W03$ for $\underline{\kappa}<1$}\\
\multicolumn{5}{|c|}{Technique: tempered distribution Taylor series.}\\\hline
Weight function & Parameter & Converg. &  & \\
& constraints & order $s$ & $c_{w}$ & $h_{w}$\\\hline\hline
\multicolumn{1}{|l|}{Sobolev spline} & $v-\frac{d}{2}\leq1$ & $s<v-\frac{d}%
{2}$ & \multicolumn{1}{|l|}{eq.\thinspace\ref{a9.6}} & $\infty$\\
\quad$\left(  v>d/2\right)  $ &  &  &  & \\\hline
\multicolumn{1}{|l|}{Extended B-spline} & $n=1$ & $s=1/2$ &
\multicolumn{1}{|l|}{eq.\thinspace\ref{a9.7}} & $\infty$\\\hline
\end{tabular}
$\caption{}\label{Tbl_NonUnisolvConverg_revisit}%
\end{table}%

\subsection{Local estimates}

In Theorem \ref{Thm_converg_interpol_not_unisolv2} we derived error estimates
for all $f\in X_{w}^{0}$. This result can be improved using the localization
results of Section \ref{Sect_local_data_space} which comes, however, at the
price of the hard-to-calculate constant $\left\Vert E\right\Vert $:

\begin{theorem}
\label{Thm_err_interp_Taylor_distrib_k<1}Suppose the weight function $w$ and
the data region $\Omega$ satisfy the localization conditions of Corollary
\ref{Cor_2_Thm_data_fn_3}. Suppose also that $w$ satisfies the convergence
conditions of Theorem \ref{Thm_converg_interpol_not_unisolv2} for
$K=\overline{\Omega}$ and $s$.

Then $X_{w}^{0}\left(  \Omega\right)  =W^{\mathbf{1}}\left(  \Omega\right)  $
as sets and
\begin{equation}
\left\vert u\left(  x\right)  -\mathcal{I}_{X}u\left(  x\right)  \right\vert
\leq c_{w}\left\Vert E\right\Vert \left\Vert u\right\Vert _{W^{\mathbf{1}%
}\left(  \Omega\right)  }\left(  h_{X,\Omega}\right)  ^{s},\quad x\in
\Omega,\text{ }u\in W^{\mathbf{1}}\left(  \Omega\right)  ,\label{a1.14}%
\end{equation}

where $\left\Vert E\right\Vert $ is the norm of the linear extension operator
$E:W^{\mathbf{1}}\left(  \Omega\right)  \rightarrow X_{w}^{0}$ used in Theorem
\ref{Thm_ex_data_fn_3}.
\end{theorem}

\begin{proof}
Theorem \ref{Thm_converg_interpol_not_unisolv2} implies%
\[
\left\vert f\left(  x\right)  -\mathcal{I}_{X}f\left(  x\right)  \right\vert
\leq c_{w}\left\Vert f\right\Vert _{w,0}\left(  h_{X,\overline{\Omega}%
}\right)  ^{s},\text{\quad}x\in\Omega,\text{ }f\in X_{w}^{0}.
\]

But from Theorem \ref{Thm_ex_data_fn_3} the extension operator
$E:W^{\mathbf{1}}\left(  \Omega\right)  \rightarrow X_{w}^{0}$ is continuous
so that if $u\in W^{\mathbf{1}}\left(  \Omega\right)  $,%
\begin{align*}
\left\vert u\left(  x\right)  -\mathcal{I}_{X}u\left(  x\right)  \right\vert
=\left\vert Eu\left(  x\right)  -\mathcal{I}_{X}Eu\left(  x\right)
\right\vert  & \leq c_{w}\left\Vert Eu\right\Vert _{w,0}\left(  h_{X,\Omega
}\right)  ^{s}\\
& \leq c_{w}\left\Vert E\right\Vert \left\Vert u\right\Vert _{W^{\mathbf{1}%
}\left(  \Omega\right)  }\left(  h_{X,\Omega}\right)  ^{s}.
\end{align*}

\end{proof}

\section{Interpolant error: $w\in W02$ for $\kappa\geq1$ or $w\in W03$ for
$\protect\underline{\kappa}\geq1$}

\subsection{Interpolation error}

The pointwise convergence estimates of this section will be obtained by
modifying those of Subsection \ref{Sect_unisolv} which used multipoint Taylor
series expansions and Lagrange interpolation using minimal unisolvent sets of
independent data points. Here we will use the remainder estimates of Theorem
\ref{Thm_radial_estim} and apply them to the Sobolev splines and extended
B-splines to obtain improved convergence rates. Recall that the operator
$\mathcal{Q}$ was introduced in Definition \ref{Def_Aux_proj_operator}.

\begin{lemma}
\label{Lem_Q_estim_2}Suppose $A=\left\{  a_{k}\right\}  _{k=1}^{M}$ is a
minimal unisolvent set of order $m\geq1$ and $f\in X_{w}^{0}$. Then%
\[
\mathcal{Q}f\left(  x\right)  =-\sum_{k=1}^{M}\left(  \mathcal{R}_{m}f\right)
\left(  x,a_{k}-x\right)  l_{k}\left(  x\right)  ,
\]

and we have the upper bound%
\[
\left\vert \mathcal{Q}f\left(  x\right)  \right\vert \leq\left(  \sum
_{k=1}^{M}\left\vert l_{k}\left(  x\right)  \right\vert \right)  \max
_{k=1}^{M}\left\vert \left(  \mathcal{R}_{m}f\right)  \left(  x,a_{k}%
-x\right)  \right\vert ,\quad x\in\mathbb{R}^{d},
\]

where $\mathcal{R}_{m}f$ is the Taylor series remainder.
\end{lemma}

\begin{proof}
From Definition \ref{Def_Aux_proj_operator} of $\mathcal{P}$ and $\mathcal{Q}
$%
\[
\mathcal{Q}f\left(  x\right)  =f\left(  x\right)  -\mathcal{P}f\left(
x\right)  =f\left(  x\right)  -\sum_{k=1}^{M}f\left(  a_{k}\right)
l_{k}\left(  x\right)  =f\left(  x\right)  -\sum_{k=1}^{M}f\left(  x+\left(
a_{k}-x\right)  \right)  l_{k}\left(  x\right)  .
\]

Using the Taylor series expansion formula with remainder we have for each $k$%
\[
f\left(  x+\left(  a_{k}-x\right)  \right)  =\sum_{\left\vert \beta\right\vert
<m}\frac{D^{\beta}f(x)}{\beta!}\left(  a_{k}-x\right)  ^{\beta}+\left(
\mathcal{R}_{m}f\right)  \left(  x,a_{k}-x\right)  ,
\]

so that%
\begin{align*}
\sum_{k=1}^{M}f\left(  x+\left(  a_{k}-x\right)  \right)  l_{k}\left(
x\right)   & =\sum_{k=1}^{M}\left(  \sum_{\left\vert \beta\right\vert <m}%
\frac{D^{\beta}f(x)}{\beta!}\left(  a_{k}-x\right)  ^{\beta}+\left(
\mathcal{R}_{m}f\right)  \left(  x,a_{k}-x\right)  \right)  l_{k}\left(
x\right) \\
& =\sum_{k=1}^{M}\sum_{\left\vert \beta\right\vert <m}\frac{D^{\beta}%
f(x)}{\beta!}\left(  a_{k}-x\right)  ^{\beta}l_{k}\left(  x\right)
+\sum_{k=1}^{M}\left(  \mathcal{R}_{m}f\right)  \left(  x,a_{k}-x\right)
l_{k}\left(  x\right)  .
\end{align*}

But by part 1 of Theorem \ref{Thm_P_op_properties} the operator $\mathcal{P}$
preserves polynomials of degree $<m$. Hence%
\begin{align*}
\sum_{k=1}^{M}\sum_{\left\vert \beta\right\vert <m}\frac{D^{\beta}f(x)}%
{\beta!}\left(  a_{k}-x\right)  ^{\beta}l_{k}\left(  x\right)   &
=\sum_{\left\vert \beta\right\vert <m}\frac{D^{\beta}f(x)}{\beta!}%
\mathcal{P}_{y}\left(  \left(  y-x\right)  ^{\beta}\right)  \left(  y=x\right)
\\
& =\sum_{\left\vert \beta\right\vert <m}\frac{D^{\beta}f(x)}{\beta!}\left(
\left(  y-x\right)  ^{\beta}\right)  \left(  y=x\right) \\
& =f\left(  x\right)  ,
\end{align*}

leaving us with%
\[
\sum_{k=1}^{M}f\left(  x+\left(  a_{k}-x\right)  \right)  l_{k}\left(
x\right)  =f\left(  x\right)  +\sum_{k=1}^{M}\left(  \mathcal{R}_{m}f\right)
\left(  x,a_{k}-x\right)  l_{k}\left(  x\right)  ,
\]

and%
\[
\mathcal{Q}f\left(  x\right)  =-\sum_{k=1}^{M}\left(  \mathcal{R}_{m}f\right)
\left(  x,a_{k}-x\right)  l_{k}\left(  x\right)  ,
\]

so that%
\[
\left\vert \mathcal{Q}f\left(  x\right)  \right\vert \leq\left(  \sum
_{k=1}^{M}\left\vert l_{k}\left(  x\right)  \right\vert \right)  \max
_{k=1}^{M}\left\vert \left(  \mathcal{R}_{m}f\right)  \left(  x,a_{k}%
-x\right)  \right\vert ,
\]

as required.
\end{proof}

The next theorem is our main convergence result. This result corresponds to
Theorem \ref{Thm_converg_interpol_ord_gte_1}.

\begin{theorem}
\label{Thm_converg_interpol_unisolv2}\textbf{Interpolant convergence using
minimal unisolvent data sets}.

Let $w$ be a weight function with property W02 for some $\kappa\geq1$ and let
$G$ be the corresponding basis function. Set $m=\left\lfloor \kappa
\right\rfloor $ so that $m\geq1$.

Suppose $\mathcal{I}_{X}f$ is the minimal norm interpolant of the data
function $f\in X_{w}^{0}$ on the independent data set $X$ contained in the
data region $\Omega$. We use the notation and assumptions of Lemma
\ref{Lem_Lagrange_interpol} which means assuming that $X$ is $m$-unisolvent
and $\Omega$ is a bounded region whose boundary satisfies the cone condition.

Now suppose there exists $\theta\geq0$ and constants $h_{w},c_{w}>0$ such that
the Taylor series remainder of Section \ref{Sect_apx_TaylorSeries} satisfies
the estimate
\begin{equation}
\left\vert \left(  \mathcal{R}_{m}f\right)  \left(  x,a\right)  \right\vert
\leq c_{w}\left\Vert f\right\Vert _{w,0}\left\vert a\right\vert ^{m+\theta
},\quad x\in\Omega,\text{ }\left\vert a\right\vert \leq h_{w},\label{a333}%
\end{equation}

for all $f\in X_{w}^{0}$.

Then there exists $h_{\Omega,m}>0$ such that for $x\in\Omega$:%
\begin{align}
\left\vert f\left(  x\right)  -\mathcal{I}_{X}f\left(  x\right)  \right\vert
& \leq\left(  c_{\Omega,m}\right)  ^{m+\theta}K_{\Omega,m}^{\prime}c_{w}%
\sqrt{\left(  f-\mathcal{I}_{X}f,f\right)  _{w,0}}\left(  h_{X,\Omega}\right)
^{m+\theta}\nonumber\\
& \leq\left(  c_{\Omega,m}\right)  ^{m+\theta}K_{\Omega,m}^{\prime}%
c_{w}\left\Vert f\right\Vert _{w,0}\left(  h_{X,\Omega}\right)  ^{m+\theta
},\label{a338}%
\end{align}

when $h_{X,\Omega}=\sup\limits_{\omega\in\Omega}\operatorname*{dist}\left(
\omega,X\right)  <\min\left\{  h_{\Omega,m},h_{w}\right\}  $ i.e. the order of
convergence is at least $m+\theta$.

Further, for all data sets we have the upper bound%
\begin{equation}
\left\vert f\left(  x\right)  -\mathcal{I}_{X}f\left(  x\right)  \right\vert
\leq K_{\Omega,m}^{\prime}c_{w}\sqrt{\left(  f-\mathcal{I}_{X}f,f\right)
_{w,0}}\left(  \operatorname*{diam}\Omega\right)  ^{m+\theta},\quad
x\in\overline{\Omega}.\label{a54.6}%
\end{equation}

The constants $c_{\Omega,m},K_{\Omega,m}^{\prime}$ and $h_{\Omega,m}$ only
depend on $\Omega,m$ and $d$.
\end{theorem}

\begin{proof}
The set $X$ is unisolvent so from Definition \ref{Def_unisolv} it must have a
minimal unisolvent subset, say $A=\left\{  a_{k}\right\}  _{k=1}^{M}$, which
we use to define the Lagrangian operators $\mathcal{P}$ and $\mathcal{Q}%
=I-\mathcal{P}$ of Definition \ref{Def_Aux_proj_operator}. Now suppose
$\max\limits_{k=1}^{M}\left\vert a_{k}-x\right\vert \leq c_{\Omega
,m}h_{X,\Omega}$. Since $f$ is a data function, $f\left(  x\right)
-\mathcal{I}_{X}f\left(  x\right)  =\mathcal{Q}\left(  f-\mathcal{I}%
_{X}f\right)  \left(  x\right)  $ and by estimating $\left\vert \mathcal{Q}%
\left(  f-\mathcal{I}_{X}f\right)  \left(  x\right)  \right\vert $ using Lemma
\ref{Lem_Q_estim_2} we get the sequence of estimates%
\begin{align}
\left\vert f\left(  x\right)  -\mathcal{I}_{X}f\left(  x\right)  \right\vert
& \leq\left(  \sum\nolimits_{k=1}^{M}\left\vert l_{k}\left(  x\right)
\right\vert \right)  \max_{k=1}^{M}\left\vert \left(  \mathcal{R}_{m}\left(
f-\mathcal{I}_{X}f\right)  \right)  \left(  x,a_{k}-x\right)  \right\vert
\nonumber\\
& \leq c_{w}\left(  \sum\nolimits_{k=1}^{M}\left\vert l_{k}\left(  x\right)
\right\vert \right)  \left\Vert f-\mathcal{I}_{X}f\right\Vert _{w,0}\max
_{k=1}^{M}\left\vert a_{k}-x\right\vert ^{m+\theta}\nonumber\\
& =c_{w}\left(  \sum\nolimits_{k=1}^{M}\left\vert l_{k}\left(  x\right)
\right\vert \right)  \left\Vert f-\mathcal{I}_{X}f\right\Vert _{w,0}\left(
\max_{k=1}^{M}\left\vert a_{k}-x\right\vert \right)  ^{m+\theta}.\label{a54.7}%
\end{align}

To estimate the last two factors on the right side of the last inequality we
will need the previous Lagrangian Lemma \ref{Lem_Lagrange_interpol}. In the
notation of this lemma, if $h_{X,\Omega}=\sup\limits_{\omega\in\Omega
}\operatorname*{dist}\left(  \omega,X\right)  <h_{\Omega,m}$ then for a given
$x$ there exists a minimal unisolvent set $A=\left\{  a_{k}\right\}
_{k=1}^{M}$ such that $\operatorname*{diam}\left(  A\cup\left\{  x\right\}
\right)  \leq c_{\Omega,\kappa}h_{X,\Omega}$ and $\sum\limits_{k=1}%
^{M}\left\vert l_{k}\left(  x\right)  \right\vert \leq K_{\Omega,m}^{\prime}$.

Thus $\left\vert a_{k}-x\right\vert \leq c_{\Omega,m}h_{X,\Omega}$ and
\[
\left\vert f\left(  x\right)  -\mathcal{I}_{X}f\left(  x\right)  \right\vert
\leq c_{w}K_{\Omega,m}^{\prime}\left\Vert f-\mathcal{I}_{X}f\right\Vert
_{w,0}\left(  c_{\Omega,m}h_{X,\Omega}\right)  ^{m+\theta}=k_{G}\left\Vert
f-\mathcal{I}_{X}f\right\Vert _{w,0}\left(  h_{X,\Omega}\right)  ^{m+\theta},
\]

and because $\left\Vert f-\mathcal{I}_{X}f\right\Vert _{w,0}^{2}=\left(
f-\mathcal{I}_{X}f,f\right)  _{w,0}$ we obtain%
\begin{equation}
\left\vert f\left(  x\right)  -\mathcal{I}_{X}f\left(  x\right)  \right\vert
\leq k_{G}\sqrt{\left(  f-\mathcal{I}_{X}f,f\right)  _{w,0}}\left(
h_{X,\Omega}\right)  ^{m+\theta},\quad x\in\Omega,\label{a336}%
\end{equation}

which is almost the required inequality. The extension of the last inequality
to $x\in\overline{\Omega}$ is an easy consequence of the fact that $f$ and
$\mathcal{I}_{X}f$ are continuous on $\mathbb{R}^{d}$. The second inequality
of \ref{a338} now follows directly from \ref{1.09} and so the order of
convergence of the minimal norm interpolant is at least $m=\left\lfloor
\kappa\right\rfloor $.

To prove \ref{a54.6} we start with \ref{a54.7} and write directly%
\begin{align*}
\left\vert f\left(  x\right)  -\mathcal{I}_{X}f\left(  x\right)  \right\vert
& \leq c_{w}K_{\Omega,m}^{\prime}\left\Vert f-\mathcal{I}_{X}f\right\Vert
_{w,0}\left(  \operatorname*{diam}\Omega\right)  ^{m+\theta}\\
& =c_{w}K_{\Omega,m}^{\prime}\sqrt{\left(  f-\mathcal{I}_{X}f,f\right)
_{w,0}}\left(  \operatorname*{diam}\Omega\right)  ^{m+\theta}.
\end{align*}

\end{proof}

\subsection{Examples}

\begin{example}
\textbf{Sobolev splines}\label{Ex_Sobol_splin_k_ge_1} The radial Sobolev
spline weight functions were introduced in Subsection
\ref{SbSect_wt_func_examples}. From \ref{1.030} the Sobolev weight functions
are%
\begin{equation}
w\left(  \xi\right)  =\left(  1+\left\vert \xi\right\vert ^{2}\right)
^{v},\text{\quad}v>d/2,\label{a54.1}%
\end{equation}

and $w$ has property W01 for empty $\mathcal{A}$ and property W02 for $\kappa$
iff $0\leq\kappa<v-d/2$. We have assumed that $\kappa\geq1$ so $v-d/2>1$ and%
\begin{equation}
m=\left\lfloor \kappa\right\rfloor =\left\{
\begin{array}
[c]{ll}%
v-d/2-1, & if\text{ }v-d/2\text{ }is\text{ }an\text{ }integer,\\
\left\lfloor v-d/2\right\rfloor , & otherwise.
\end{array}
\right. \label{a54.5}%
\end{equation}

We will use the radial remainder estimate \ref{a223} i.e.
\[
\left\vert \left(  \mathcal{R}_{m+1}f\right)  \left(  x,a\right)  \right\vert
\leq\frac{2+\frac{1}{m+1}}{m!}\left(  2\int_{0}^{\infty}\frac{s^{2\kappa}%
}{\overset{\circ}{w}\left(  s\right)  }ds\right)  ^{1/2}\left\Vert
f\right\Vert _{w,0}\left\vert a\right\vert ^{\kappa},\quad a,x\in
\mathbb{R}^{d}.
\]

Comparison with \ref{a333} yields:
\begin{equation}
c_{w}=\frac{2+\frac{1}{m+1}}{m!}\left(  2\int_{0}^{\infty}\frac{s^{2\kappa}%
}{\overset{\circ}{w}\left(  s\right)  }ds\right)  ^{1/2},\text{ }h_{w}%
=\infty,\text{ }\theta=\kappa-m,\label{a110}%
\end{equation}

where the constant $\int_{0}^{\infty}\frac{s^{2\kappa}}{\overset{\circ
}{w}\left(  s\right)  }ds$ is given below by \ref{a54.2}. We now show how to
calculate $\int_{0}^{\infty}\frac{s^{2\kappa}}{\overset{\circ}{w}\left(
s\right)  }ds$, starting with the case $d=1$:\medskip

\fbox{\textbf{Case 1} $d=1$} Here $\overset{\circ}{w}=w$. Now the beta
function $B$ and the gamma function $\Gamma$ satisfy%
\begin{align}
\int_{0}^{\pi/2}\sin^{2p}\theta\cos^{2q}\theta\,d\theta & =\frac{1}{2}B\left(
p+\frac{1}{2},q+\frac{1}{2}\right) \nonumber\\
& =\frac{1}{2}\frac{\Gamma\left(  p+\frac{1}{2}\right)  \Gamma\left(
q+\frac{1}{2}\right)  }{\Gamma\left(  p+q+1\right)  },\quad p,q>-\frac{1}%
{2},\label{a50.1}%
\end{align}

and we now show that%
\begin{equation}
\int_{0}^{\infty}\frac{s^{2\kappa}ds}{\left(  1+s^{2}\right)  ^{v}}=\frac
{1}{2}B\left(  \kappa+\frac{1}{2},v-\kappa-\frac{1}{2}\right)  ,\quad
v>\kappa+\frac{1}{2}>-\frac{1}{2}.\label{a229}%
\end{equation}

In fact, the change of variables
\begin{equation}
s=\tan\theta,\text{ }ds=\sec^{2}\theta d\theta,\text{ }1+s^{2}=\sec^{2}%
\theta,\label{a227}%
\end{equation}

yields%
\begin{align*}
\int_{0}^{\infty}\frac{s^{2\kappa}ds}{\left(  1+s^{2}\right)  ^{v}}=\int%
_{0}^{\pi/2}\tan^{2\kappa}\theta\cos^{2v}\theta\sec^{2}\theta d\theta &
=\int_{0}^{\pi/2}\sin^{2\kappa}\theta\cos^{2v-2\kappa-2}\theta d\theta\\
&  =\frac{1}{2}B\left(  \kappa+\frac{1}{2},v-\kappa-\frac{1}{2}\right)  ,
\end{align*}

and so%
\begin{equation}
\int_{0}^{\infty}\frac{s^{2\kappa}}{\overset{\circ}{w}\left(  s\right)
}ds=\frac{1}{2}B\left(  \kappa+\frac{1}{2},v-\kappa-\frac{1}{2}\right)
.\label{a228}%
\end{equation}

\underline{Now suppose $d\geq2$} Here \ref{a223} is the relevant estimate and
we need to calculate $\int_{0}^{\infty}\frac{s^{2\kappa}}{\overset{\circ
}{w}\left(  s\right)  }ds$ where $\overset{\circ}{w}$ is defined by \ref{a119}
i.e.%
\[
\frac{1}{\overset{\circ}{w}\left(  s\right)  }=2^{d-1}\int_{\xi^{\prime}\geq
0}\frac{d\xi^{\prime}}{w\left(  s,\xi^{\prime}\right)  }=2^{d-1}\int%
_{\xi^{\prime}\geq0}\frac{d\xi^{\prime}}{\left(  1+s^{2}+\left\vert
\xi^{\prime}\right\vert ^{2}\right)  ^{v}},
\]

and the change of variables $\xi^{\prime}=\left(  1+s^{2}\right)  ^{\frac
{1}{2}}\eta^{\prime}$ implies $d\xi^{\prime}=\left(  1+s^{2}\right)
^{\frac{d-1}{2}}d\eta^{\prime}$ so that%
\begin{equation}
\frac{1}{\overset{\circ}{w}\left(  s\right)  }=\frac{\lambda_{v}}{\left(
1+s^{2}\right)  ^{v-\frac{d}{2}+\frac{1}{2}}},\qquad\lambda_{v}=2^{d-1}%
\int_{\xi^{\prime}\geq0}\frac{d\eta^{\prime}}{\left(  1+\left\vert
\eta^{\prime}\right\vert ^{2}\right)  ^{v}}.\label{a330}%
\end{equation}

and%
\begin{equation}
\int_{0}^{\infty}\frac{s^{2\kappa}}{\overset{\circ}{w}\left(  s\right)
}ds=\lambda_{v}\int_{0}^{\infty}\frac{s^{2\kappa}}{\left(  1+s^{2}\right)
^{v-\frac{d}{2}+\frac{1}{2}}}ds.\label{a50.6}%
\end{equation}

This means that%
\begin{equation}
\int_{0}^{\infty}\frac{s^{2\kappa}}{\overset{\circ}{w}\left(  s\right)
}ds=\lambda_{v}\int_{0}^{\infty}\frac{s^{2\kappa}}{\left(  1+s^{2}\right)
^{v-\frac{d}{2}+\frac{1}{2}}}ds=\frac{1}{2}\lambda_{v}B\left(  \kappa+\frac
{1}{2},v-\kappa-\frac{d}{2}\right)  .\label{a50.8}%
\end{equation}

The next step is to calculate $\lambda_{v}$. The two cases are $d=2$ and
$d\geq3$:\medskip

\fbox{\textbf{Case 2} $d=2$} Here $v>1$ and from \ref{a330} and \ref{a229}%
\[
\lambda_{v}=2\int_{0}^{\infty}\frac{ds}{\left(  1+s^{2}\right)  ^{v}}=\frac
{1}{2}B\left(  \frac{1}{2},v-\frac{1}{2}\right)  ,
\]

and%
\begin{equation}
\int_{0}^{\infty}\frac{s^{2\kappa}}{\overset{\circ}{w}\left(  s\right)
}ds=\frac{1}{4}B\left(  \frac{1}{2},v-\frac{1}{2}\right)  B\left(
\kappa+\frac{1}{2},v-\kappa-1\right)  ,\quad d=2.\label{a331}%
\end{equation}
\medskip

\fbox{\textbf{Case 3} $d\geq3$} Note first that $v>3/2$. To calculate
$\lambda_{v}$ we use a spherical polar change of variables e.g. Section 5.43
Adams \cite{Adams75}: $\eta^{\prime}=\left(  r,\phi\right)  =\left(
r,\phi_{1},\ldots,\phi_{d-2}\right)  $, $r\geq0$, $0\leq\phi_{j}\leq\frac{\pi
}{2}\mathbf{1}$, which has volume element%
\[
d\eta^{\prime}=r^{d-2}\prod_{j=1}^{d-2}\sin^{j-1}\phi_{j}drd\phi,
\]

where $d\phi=d\phi_{1}\ldots d\phi_{d-2}$. Thus%
\begin{align*}
\lambda_{v}=2^{d-1}\int\limits_{\eta^{\prime}\geq0}\frac{d\eta^{\prime}%
}{\left(  1+\left\vert \eta^{\prime}\right\vert ^{2}\right)  ^{v}} &
=2^{d-1}\int\limits_{0}^{\frac{\pi}{2}\mathbf{1}}\int\limits_{0}^{\infty}%
\frac{r^{d-2}\prod_{j=1}^{d-2}\sin^{j-1}\phi_{j}drd\phi}{\left(
1+r^{2}\right)  ^{v}}\\
&  =2^{d-1}\int\limits_{0}^{\infty}\frac{r^{d-2}}{\left(  1+r^{2}\right)
^{v}}dr\prod_{j=1}^{d-2}\int\limits_{0}^{\frac{\pi}{2}}\sin^{j-1}\phi_{j}%
d\phi_{j},
\end{align*}

and by \ref{a50.1}%
\[
\int_{0}^{\infty}\frac{r^{d-2}}{\left(  1+r^{2}\right)  ^{v}}dr=\frac{1}%
{2}B\left(  \frac{d-1}{2},v-\frac{d-1}{2}\right)  ,
\]

so that%
\begin{align*}
\lambda_{v}  & =2^{d-1}\frac{1}{2}B\left(  \frac{d-1}{2},v-\frac{d-1}%
{2}\right)  \prod_{j=1}^{d-2}\int\limits_{0}^{\frac{\pi}{2}}\sin^{j-1}\phi
_{j}d\phi_{j}\\
& =2^{d-2}B\left(  \frac{d-1}{2},v-\frac{d-1}{2}\right)  \prod_{j=1}%
^{d-2}\frac{1}{2}\frac{\Gamma\left(  \frac{1}{2}\right)  \Gamma\left(
\frac{j}{2}\right)  }{\Gamma\left(  \frac{j+1}{2}\right)  }\\
& =\frac{\Gamma\left(  \frac{1}{2}\right)  ^{d-1}}{\Gamma\left(  \frac{d-1}%
{2}\right)  }B\left(  \frac{d-1}{2},v-\frac{d-1}{2}\right)  ,
\end{align*}

and from \ref{a50.8}%
\begin{align}
\int_{0}^{\infty}\frac{s^{2\kappa}}{\overset{\circ}{w}\left(  s\right)  }ds  &
=\frac{1}{2}\lambda_{v}B\left(  \kappa+\frac{1}{2},v-\kappa-\frac{d}{2}\right)
\nonumber\\
& =\frac{\Gamma\left(  \frac{1}{2}\right)  ^{d-1}}{\Gamma\left(  \frac{d-1}%
{2}\right)  }B\left(  \frac{d-1}{2},v-\frac{d-1}{2}\right)  B\left(
\kappa+\frac{1}{2},v-\kappa-\frac{d}{2}\right)  ,\quad d\geq3.\label{a332}%
\end{align}

We now have our three estimates \ref{a228}, \ref{a331} and \ref{a332}. But
\ref{a332} becomes \ref{a331} when we set $d=2$ so \textbf{in summary:}%
\begin{equation}
\int_{0}^{\infty}\frac{s^{2\kappa}}{\overset{\circ}{w}\left(  s\right)
}ds=\left\{
\begin{array}
[c]{ll}%
\frac{1}{2}B\left(  \kappa+\frac{1}{2},v-\kappa-\frac{1}{2}\right)  , & d=1,\\
\frac{\Gamma\left(  \frac{1}{2}\right)  ^{d-1}}{\Gamma\left(  \frac{d-1}%
{2}\right)  }B\left(  \frac{d-1}{2},v-\frac{d-1}{2}\right)  B\left(
\kappa+\frac{1}{2},v-\kappa-\frac{d}{2}\right)  , & d\geq2,
\end{array}
\right. \label{a54.2}%
\end{equation}

and%
\[
\lambda_{v}=\frac{\Gamma\left(  \frac{1}{2}\right)  ^{d-1}}{\Gamma\left(
\frac{d-1}{2}\right)  }B\left(  \frac{d-1}{2},v-\frac{d-1}{2}\right)  ,\quad
d\geq2.
\]

\end{example}

\begin{example}
\textbf{The extended B-splines }$l,n\geq1$. Introduced in Theorem
\ref{Thm_ex_splin_wt_fn_properties}. They have property W03 for $\kappa$ iff
$\kappa<\left(  n-1/2\right)  \mathbf{1}$. Set $m=\left\lfloor
\underline{\kappa}\right\rfloor $ and so $m\geq1$ when $n\geq2$. They have
property W02 for $\underline{\kappa}$.

Thus, from part 1 of Theorem \ref{Thm_Tayor_rem_estim},%
\[
\left\vert \left(  \mathcal{R}_{m+1}f\right)  \left(  x,a\right)  \right\vert
\leq\frac{c_{m}}{m!}\left(  \int\frac{\left\vert \cdot\right\vert ^{2\kappa}%
}{w_{s}}\right)  ^{1/2}\left\Vert f\right\Vert _{w,0}\left\vert a\right\vert
^{\kappa},\quad x,a\in\mathbb{R}^{d},
\]

and so \ref{a333} is satisfied for
\begin{gather}
n\geq2,\text{ }m=n-1,\text{ }\theta=\kappa-m,\text{\ }h_{w}=\infty,\nonumber\\
c_{w}=\frac{c_{m}}{m!}\left(  \int\frac{\left\vert \cdot\right\vert ^{2\kappa
}}{w_{s}}\right)  ^{1/2}.\label{a9.8}%
\end{gather}

Thus part 1 of Theorem \ref{Thm_Tayor_rem_estim} yields an improvement of
$0<\theta<1/2$ in the estimated order of convergence.

However, we can obtain a further improvement in the order of convergence using
part 4 of Theorem \ref{Thm_Tayor_rem_estim} and the upper bounds of Section
\ref{Sect_rem_estim_ten_prod}. Indeed, by Theorem
\ref{Thm_basis_tensor_hat_W3}, $D^{2m+1}G_{1}\in L^{\infty}$ so by Theorem
\ref{Thm_rem_estim_ten_prod} and Theorem \ref{Thm_estim_max_DjG1} we have
estimates of order $m+1/2=n-1/2$.
\end{example}

\subsection{Examples: summary table}

-%
\begin{table}[htbp] \centering
$%
\begin{tabular}
[c]{|c|c||c|c|c|}\hline
\multicolumn{5}{|c|}{Interpolant error estimates: W02 for $\kappa\geq1$ or W03
for $\underline{\kappa}\geq1$.}\\
\multicolumn{5}{|c|}{Techniques: unisolvency, multipoint tempered distribution
Taylor series.}\\\hline
Weight function & Parameter & Converg. &  & \\
& constraints & order & $c_{w}$ & $h_{w}$\\\hline\hline
\multicolumn{1}{|l|}{Sobolev spline (W02)} &  & $<v-d/2$ &
\multicolumn{1}{|l|}{\ref{a110}} & $\infty$\\\cline{2-3}%
\quad$\left(  v>d/2\right)  $ &  &  &  & \\\hline
\multicolumn{1}{|l|}{Extended B-spline (W03)} & $n\geq2$ & $n-1/2$ &  &
$\infty$\\\hline
$\left(  n\leq l\right)  $ &  &  &  & \\\hline
\end{tabular}
$\caption{}\label{Tbl_UnisolvConverg_revisit}%
\end{table}%

\subsection{Local estimates}

In Theorem \ref{Thm_converg_interpol_unisolv2} we derived error estimates for
all $f\in X_{w}^{0}$. This result can be improved using the localization
results of Section \ref{Sect_local_data_space} which comes, however, at the
price of the possibly hard-to-calculate constant $\left\Vert E\right\Vert $:

\begin{theorem}
\label{Thm_err_interp_Taylor_distrib_k>=1}Suppose the weight function $w$ and
the data region $\Omega$ satisfy the localization conditions of Corollary
\ref{Cor_2_Thm_data_fn_3} for some $n\geq1$. Suppose also that $w$ satisfies
the convergence conditions of Theorem \ref{Thm_converg_interpol_unisolv2} for
$m=\left\lfloor \underline{\kappa}\right\rfloor $ and some $\theta\geq0$.

Then $X_{w}^{0}\left(  \Omega\right)  =W^{n\mathbf{1}}\left(  \Omega\right)  $
as sets and
\[
\left\vert u\left(  x\right)  -\mathcal{I}_{X}u\left(  x\right)  \right\vert
\leq\left(  c_{\Omega,m}\right)  ^{m+\theta}K_{\Omega,m}^{\prime}%
c_{w}\left\Vert E\right\Vert \left\Vert u\right\Vert _{W^{n\mathbf{1}}\left(
\Omega\right)  }\left(  h_{X,\Omega}\right)  ^{m+\theta},\quad x\in
\Omega,\text{ }u\in W^{n\mathbf{1}}\left(  \Omega\right)  ,
\]

where $\left\Vert E\right\Vert $ is the norm of the linear extension operator
$E:W^{n\mathbf{1}}\left(  \Omega\right)  \rightarrow X_{w}^{0}$ introduced in
Theorem \ref{Thm_ex_data_fn_3}.
\end{theorem}

\begin{proof}
Theorem \ref{Thm_converg_interpol_unisolv2} implies%
\[
\left\vert f\left(  x\right)  -\mathcal{I}_{X}f\left(  x\right)  \right\vert
\leq\left(  c_{\Omega,m}\right)  ^{m+\theta}K_{\Omega,m}^{\prime}%
c_{w}\left\Vert f\right\Vert _{w,0}\left(  h_{X,\Omega}\right)  ^{m+\theta
},\text{\quad}x\in\Omega,\text{ }f\in X_{w}^{0}.
\]

But from Theorem \ref{Thm_ex_data_fn_3} the extension operator
$E:W^{n\mathbf{1}}\left(  \Omega\right)  \rightarrow X_{w}^{0}$ is continuous
so that if $u\in W^{n\mathbf{1}}\left(  \Omega\right)  $ and $x\in\Omega$,%
\begin{align*}
\left\vert u\left(  x\right)  -\mathcal{I}_{X}u\left(  x\right)  \right\vert
=\left\vert Eu\left(  x\right)  -\mathcal{I}_{X}Eu\left(  x\right)
\right\vert  & \leq\left(  c_{\Omega,m}\right)  ^{m+\theta}K_{\Omega
,m}^{\prime}c_{w}\left\Vert Eu\right\Vert _{w,0}\left(  h_{X,\Omega}\right)
^{m+\theta}\\
& \leq\left(  c_{\Omega,m}\right)  ^{m+\theta}K_{\Omega,m}^{\prime}%
c_{w}\left\Vert E\right\Vert \left\Vert u\right\Vert _{W^{n\mathbf{1}}\left(
\Omega\right)  }\left(  h_{X,\Omega}\right)  ^{m+\theta}.
\end{align*}

\end{proof}

\begin{remark}
The Sobolev spline weight function of Table \ref{Tbl_UnisolvConverg_revisit}
does not satisfy the localization conditions of Corollary
\ref{Cor_2_Thm_data_fn_3} since it does not satisfy the inequalities
\ref{1.39}.
\end{remark}

\chapter{Local interpolation error for data functions in $W^{1,\infty}\left(
\Omega\right)  \cap X_{w}^{0}\left(  \Omega\right)  $%
\label{Ch_interpol_err_H1inf_data}}

\section{Introduction}

Here the local data functions are assumed to have bounded derivatives on a
bounded data region $\Omega$ i.e. functions in $W^{1,\infty}\left(
\Omega\right)  \cap X_{w}^{0}\left(  \Omega\right)  $ where $X_{w}^{0}$ is a
global data space. See Section \ref{Sect_Ex_smth_err_H1inf_data} for the
corresponding results concerning the Exact smoother.

In Section \ref{Lem_interpol_estim_dim1_W1inf_Xow} we will consider
1-dimensional estimates and the main result is Theorem
\ref{Thm_err_interpol_H1inf_data}. An example is provided by Corollary
\ref{vCor_bound_deriv_scal_hat_interpol} of the Appendix where it is shown
that in \textbf{one dimension} the interpolant generated by the \textbf{scaled
hat basis function} $\Lambda\left(  x/\lambda\right)  $ satisfies the local
estimate:
\[
\left\Vert D\mathcal{I}_{X}f_{d}\right\Vert _{\infty;\Omega}\leq\frac
{2}{\lambda}\left\Vert f_{d}\right\Vert _{\infty;\Omega}+\left\Vert
Df_{d}\right\Vert _{\infty;\Omega},\quad f_{d}\in W^{1,\infty}\left(
\Omega\right)  \cap X_{w}^{0}\left(  \Omega\right)  .
\]

The scaled hat basis function is assumed to have large support w.r.t. the data
region i.e. $\operatorname*{diam}\Omega\geq\lambda$.

In Section \ref{Lem_interpol_estim_multivar_W1inf_Xow} an analogous
multivariate error estimate is proven where the interpolant is assumed to
satisfy the inequality%
\[
\left\Vert D^{\alpha}\mathcal{I}_{X}f_{d}\right\Vert _{W^{\infty}\left(
\Omega\right)  }\leq c\left\Vert f_{d}\right\Vert _{W^{1,\infty}\left(
\Omega\right)  },\quad x\in\Omega,\text{ }f_{d}\in W^{1,\infty}\left(
\Omega\right)  \cap X_{w}^{0}\left(  \Omega\right)  ,\text{ }\left\vert
\alpha\right\vert =1.
\]

I have not been able to prove that the scaled tensor product hat function with
large support satisfies this inequality.

In the light of these observations we will use the following spaces from
Definition \ref{Def_SobolevSpace}:

\begin{definition}
\label{Def_H1inf}\textbf{The Sobolev space} $W^{1,\infty}\left(
\Omega\right)  $: If $\Omega$ is an open set in$\mathbb{\ R}^{d}$ then define
the Sobolev space:%
\[
W^{1,\infty}\left(  \Omega\right)  =\left\{  f\in L^{\infty}\left(
\Omega\right)  :Df\in L^{\infty}\left(  \Omega\right)  \right\}  ,
\]

endowed with the supremum norm: $\left\Vert f\right\Vert _{1,\infty,\Omega
}=\sum\limits_{k=0}^{1}\left\Vert D^{k}f\right\Vert _{\infty,\Omega}$.

It is also well known that if $\Omega$ is bounded then $W^{1,\infty}\left(
\Omega\right)  $ is a Banach space.

Also
\[
W^{1,\infty}\left(  \Omega\right)  :=\left\{  f\in C_{B}^{\left(  0\right)
}\left(  \Omega\right)  :Df\in L^{\infty}\left(  \Omega\right)  \right\}
=\left\{  f\in\mathcal{D}^{\prime}:Df\in L^{\infty}\left(  \Omega\right)
\right\}  .
\]

\end{definition}

Our data function spaces require the derivative to be bounded. This excludes
functions with cusps, vertical inflections/tangents and half-vertical tangents
but allows corners. See Subsection \ref{SbSect_notC1b_data_fns_H1inf_1dim} for
a selection of data functions for numerical experiments.

\section{Local error estimates on $\mathbb{R}^{1}$%
\label{Lem_interpol_estim_dim1_W1inf_Xow}}

We will need the following local form of the Taylor series given in Lemma
\ref{Lem_Taylor_extension}:

\begin{lemma}
\label{Lem_Taylor_rem1_1dim_local}Assume $\Omega\subset\mathbb{R}^{1}$ is
an\textbf{\ open interval}. Suppose $u\in C_{B}^{\left(  0\right)  }\left(
\Omega\right)  $ has distributional derivative $Du\in L^{\infty}\left(
\Omega\right)  $.

Then%
\[
u(z+b)=u(z)+\left(  \mathcal{R}_{1}u\right)  \left(  z,b\right)  ,\quad
z,z+b\in\Omega,
\]

where $\mathcal{R}_{1}u$ is the integral remainder term
\begin{equation}
\left(  \mathcal{R}_{1}u\right)  \left(  z,b\right)  =b\int_{0}^{1}%
(Du)(z+\left(  1-t\right)  b)dt,\label{a1.51}%
\end{equation}

which satisfies
\begin{equation}
\left\vert \left(  \mathcal{R}_{1}u\right)  \left(  z,b\right)  \right\vert
\leq\left\Vert Du\right\Vert _{\infty;\left[  z,z+b\right]  }\left\vert
b\right\vert ,\label{a1.53}%
\end{equation}

and%
\[
\left\vert \left(  \mathcal{R}_{1}u\right)  \left(  z,b\right)  \right\vert
\leq\left(  \left\Vert Du\right\Vert _{\infty;\left[  z,z+b\right]  }\right)
\left\vert b\right\vert .
\]

\end{lemma}

\begin{proof}
Using the test functions $C_{0}^{\infty}\left(  \Omega\right)  $ instead of
$C_{0}^{\infty}\left(  \mathbb{R}^{1}\right)  $ we can easily adapt the proof
of Lemma \ref{Lem_Taylor_extension}.
\end{proof}

The next error estimate is inspired by the (trivial) interpolation Example
\ref{Ex_interp_hat_larg_supp} which considers the scaled hat basis functions.

\begin{theorem}
\label{Thm_err_interpol_H1inf_data}\textbf{One-dimensional result} Suppose the
data region $\Omega\subset\mathbb{R}^{1}$ is an\textbf{\ open interval} and
suppose that the basis function satisfies $G\in W^{1,\infty}\left(
\Omega-\Omega\right)  $. We will choose our data functions $f_{d}\in
W^{1,\infty}\left(  \Omega\right)  \cap X_{w}^{0}\left(  \Omega\right)  $ and
assume that there exist constants $c_{0},c_{1}>0$ such that the basis function
interpolant $\mathcal{I}_{X}f_{d}$ satisfies%
\begin{equation}
\left\vert D\mathcal{I}_{X}f_{d}\left(  x\right)  \right\vert \leq
c_{0}\left\Vert f_{d}\right\Vert _{\infty,\Omega}+c_{1}\left\Vert
Df_{d}\right\Vert _{\infty,\Omega},\quad x\in\Omega,\text{ }f_{d}\in
W^{1,\infty}\left(  \Omega\right)  \cap X_{w}^{0}\left(  \Omega\right)
.\label{a4.6}%
\end{equation}

Then we have the following pointwise error estimate
\[
\left\vert \mathcal{I}_{X}f_{d}\left(  x\right)  -f_{d}\left(  x\right)
\right\vert \leq\left(  c_{0}\left\Vert f_{d}\right\Vert _{\infty,\Omega
}+\left(  1+c_{1}\right)  \left\Vert Df_{d}\right\Vert _{\infty,\Omega
}\right)  \left\vert x-x^{\left(  k\right)  }\right\vert ,\quad x\in
\Omega,\text{ }x^{\left(  k\right)  }\in\Omega.
\]

Further, we have the \textbf{spherical cavity} error estimate%
\begin{equation}
\left\vert \mathcal{I}_{X}f_{d}\left(  x\right)  -f_{d}\left(  x\right)
\right\vert \leq\left(  c_{0}\left\Vert f_{d}\right\Vert _{\infty,\Omega
}+\left(  1+c_{1}\right)  \left\Vert Df_{d}\right\Vert _{\infty,\Omega
}\right)  h_{\Omega,X},\quad x\in\Omega,\label{a6}%
\end{equation}

where
\begin{equation}
h_{\Omega,X}=\sup_{x\in\Omega}\min\limits_{k=1}^{N}\left\vert x-x^{\left(
k\right)  }\right\vert ,\label{a4.7}%
\end{equation}

is the \textbf{maximum spherical cavity size (radius)}.
\end{theorem}

\begin{proof}
Since $G\in W^{1,\infty}\left(  \Omega-\Omega\right)  $ the basis function
formula \ref{1.271} for the interpolant means that $\mathcal{I}_{X}f_{d}\in
C_{B}^{\left(  0\right)  }\cap W^{1,\infty}\left(  \Omega\right)  $. Since
$f_{d}\in X_{w}^{0}\left(  \Omega\right)  $ we have $f_{d}\in C_{B}^{\left(
0\right)  }\left(  \Omega\right)  $ and applying Taylor series expansion
result of Lemma \ref{Lem_Taylor_rem1_1dim_local} with $x\in\Omega$,
$z=x^{\left(  k\right)  }\in X$ and $b=x-x^{\left(  k\right)  }$ yields%
\begin{align*}
\mathcal{I}_{X}f_{d}\left(  x\right)  -f_{d}\left(  x\right)   &  =\left(
\mathcal{I}_{X}f_{d}-f_{d}\right)  \left(  x^{\left(  k\right)  }+\left(
x-x^{\left(  k\right)  }\right)  \right) \\
&  =\left(  \mathcal{I}_{X}f_{d}-f_{d}\right)  \left(  x^{\left(  k\right)
}\right)  +\left(  \mathcal{R}_{1}\left(  \mathcal{I}_{X}f_{d}-f_{d}\right)
\right)  \left(  x^{\left(  k\right)  },x-x^{\left(  k\right)  }\right) \\
&  =\left(  \mathcal{R}_{1}\left(  \mathcal{I}_{X}f_{d}-f_{d}\right)  \right)
\left(  x^{\left(  k\right)  },x-x^{\left(  k\right)  }\right)  ,
\end{align*}

with the remainder estimate \ref{a1.53}, namely
\begin{align*}
\left\vert \left(  \mathcal{R}_{1}\left(  \mathcal{I}_{X}f_{d}-f_{d}\right)
\right)  \left(  x^{\left(  k\right)  },x-x^{\left(  k\right)  }\right)
\right\vert  &  \leq\left(  \max_{y\in\left[  x^{\left(  k\right)  },x\right]
}\left\vert D\left(  \mathcal{I}_{X}f_{d}-f_{d}\right)  (y)\right\vert
\right)  \left\vert x-x^{\left(  k\right)  }\right\vert \\
&  \leq\left(  \left\Vert D\left(  \mathcal{I}_{X}f_{d}-f_{d}\right)
\right\Vert _{\infty,\Omega}\right)  \left\vert x-x^{\left(  k\right)
}\right\vert \\
&  \leq\left(  \left\Vert D\mathcal{I}_{X}f_{d}\right\Vert _{\infty,\Omega
}+\left\Vert Df_{d}\right\Vert _{\infty,\Omega}\right)  \left\vert
x-x^{\left(  k\right)  }\right\vert .
\end{align*}

Hence, using the assumption \ref{a4.6},%
\begin{align}
\left\vert \mathcal{I}_{X}f_{d}\left(  x\right)  -f_{d}\left(  x\right)
\right\vert  & \leq\left(  c_{0}\left\Vert f_{d}\right\Vert _{\infty,\Omega
}+c_{1}\left\Vert Df_{d}\right\Vert _{\infty,\Omega}+\left\Vert Df_{d}%
\right\Vert _{\infty,\Omega}\right)  \left\vert x-x^{\left(  k\right)
}\right\vert \nonumber\\
& \leq\left(  c_{0}\left\Vert f_{d}\right\Vert _{\infty,\Omega}+\left(
1+c_{1}\right)  \left\Vert Df_{d}\right\Vert _{\infty,\Omega}\right)
\left\vert x-x^{\left(  k\right)  }\right\vert ,\label{a5}%
\end{align}

and thus%
\begin{align*}
\left\vert \mathcal{I}_{X}f_{d}\left(  x\right)  -f_{d}\left(  x\right)
\right\vert  & \leq\left(  c_{0}\left\Vert f_{d}\right\Vert _{\infty,\Omega
}+\left(  1+c_{1}\right)  \left\Vert Df_{d}\right\Vert _{\infty,\Omega
}\right)  \min_{k=1}^{N}\left\vert x-x^{\left(  k\right)  }\right\vert \\
& \leq\left(  c_{0}\left\Vert f_{d}\right\Vert _{\infty,\Omega}+\left(
1+c_{1}\right)  \left\Vert Df_{d}\right\Vert _{\infty,\Omega}\right)
\sup_{x\in\overline{\Omega}}\min_{k=1}^{N}\left\vert x-x^{\left(  k\right)
}\right\vert \\
& =\left(  c_{0}\left\Vert f_{d}\right\Vert _{\infty,\Omega}+\left(
1+c_{1}\right)  \left\Vert Df_{d}\right\Vert _{\infty,\Omega}\right)
h_{\Omega,X}.
\end{align*}

\end{proof}

\begin{corollary}
\label{Cor_Thm_err_interpol_H1inf_data}Under the conditions of Theorem
\ref{Thm_err_interpol_H1inf_data} the interpolant operator $\mathcal{I}_{X}$
satisfies
\[
\left\Vert \mathcal{I}_{X}f_{d}\right\Vert _{1,\infty,\Omega}\leq\left(
1+c_{0}+c_{1}\right)  \left(  1+h_{\Omega,X}\right)  \left\Vert f_{d}%
\right\Vert _{1,\infty,\Omega},\quad f_{d}\in W^{1,\infty}\left(
\Omega\right)  \cap X_{w}^{0}\left(  \Omega\right)  .
\]

\end{corollary}

\begin{proof}
From \ref{a6} and then \ref{a4.6} we get:%
\begin{align*}
&  \left\Vert \mathcal{I}_{X}f_{d}\right\Vert _{1,\infty,\Omega}\\
&  \leq\left\Vert \mathcal{I}_{X}f_{d}\right\Vert _{\infty,\Omega}+\left\Vert
D\mathcal{I}_{X}f_{d}\right\Vert _{\infty,\Omega}\\
&  \leq\left\Vert f_{d}\right\Vert _{\infty,\Omega}+\left(  c_{0}\left\Vert
f_{d}\right\Vert _{\infty,\Omega}+\left(  1+c_{1}\right)  \left\Vert
Df_{d}\right\Vert _{\infty,\Omega}\right)  h_{\Omega,X}+\left\Vert
D\mathcal{I}_{X}f_{d}\right\Vert _{\infty,\Omega}\\
&  =\left(  1+c_{0}h_{\Omega,X}\right)  \left\Vert f_{d}\right\Vert
_{\infty,\Omega}+\left(  1+c_{1}\right)  h_{\Omega,X}\left\Vert Df_{d}%
\right\Vert _{\infty,\Omega}+\left\Vert D\mathcal{I}_{X}f_{d}\right\Vert
_{\infty,\Omega}\\
&  \leq\left(  1+c_{0}h_{\Omega,X}\right)  \left\Vert f_{d}\right\Vert
_{\infty,\Omega}+\left(  1+c_{1}\right)  h_{\Omega,X}\left\Vert Df_{d}%
\right\Vert _{\infty,\Omega}+c_{0}\left\Vert f_{d}\right\Vert _{\infty,\Omega
}+c_{1}\left\Vert Df_{d}\right\Vert _{\infty,\Omega}\\
&  \leq\left\{  1+c_{0}h_{\Omega,X}+\left(  1+c_{1}\right)  h_{\Omega,X}%
+c_{0}+c_{1}\right\}  \left\Vert f_{d}\right\Vert _{1,\infty,\Omega}\\
&  =\left\{  1+c_{0}+c_{1}+\left(  1+c_{0}+c_{1}\right)  h_{\Omega,X}\right\}
\left\Vert f_{d}\right\Vert _{1,\infty,\Omega}\\
&  =\left(  1+c_{0}+c_{1}\right)  \left(  1+h_{\Omega,X}\right)  \left\Vert
f_{d}\right\Vert _{1,\infty,\Omega}.
\end{align*}

\end{proof}

\begin{example}
\label{Ex_interp_hat_larg_supp}\textbf{1-dimensional scaled hat basis function
with large support w.r.t. the data region}

Since $r_{\Omega}\Lambda\in$ $X_{w}^{0}\left(  \Omega\right)  $ it is clear
that $X_{w}^{0}\left(  \Omega\right)  \cap W^{1,\infty}\left(  \Omega\right)
\neq\left\{  0\right\}  $.

These basis functions are discussed in Chapter
\ref{Ch_bnd_deriv_hat_smth_large_supp} of the Appendix. Large support w.r.t.
the data region $\Omega$ means that the scaled hat basis function
$\Lambda_{\lambda}\left(  x\right)  =\Lambda\left(  x/\lambda\right)  $ has
been scaled so that
\[
\operatorname*{diam}\Omega\leq\frac{1}{2}\operatorname*{diam}%
\operatorname*{supp}\Lambda_{\lambda}\text{ }i.e.\operatorname*{diam}%
\Omega\leq\lambda.
\]

This assumption has the nice consequence\textbf{\ that }$x,x^{\prime}\in
\Omega$\textbf{\ now implies }$\Lambda_{\lambda}\left(  x-x^{\prime}\right)
=1-\frac{\left\vert x-x^{\prime}\right\vert }{\lambda}$ \textbf{and there are
no zero values.} This will mean that \textbf{in 1-dimension the interpolant is
the (rather trivial) continuous, piecewise linear interpolant} with a finite
number of derivative values which correspond to the linear segments in the
data intervals $\left[  x^{\left(  k\right)  },x^{\left(  k+1\right)
}\right]  $.

In Corollary \ref{vCor_bound_deriv_scal_hat_interpol} it is shown that if
$\operatorname*{diam}\Omega\leq\lambda$ and the data function $f_{d}\in
X_{w}^{0} $ satisfies $Df_{d}\in L^{\infty}\left(  \Omega\right)  $ then the
basis function interpolant $\mathcal{I}_{X}f_{d}$ corresponding to the scaled
hat basis function $\Lambda\left(  \cdot/\lambda\right)  $ satisfies%
\[
\left\Vert D\mathcal{I}_{X}f_{d}\right\Vert _{\infty;\Omega}\leq\frac
{2}{\lambda}\left\Vert f_{d}\right\Vert _{\infty;\Omega}+\left\Vert
Df_{d}\right\Vert _{\infty;\Omega}.
\]

Theorem \ref{Thm_err_interpol_H1inf_data} now implies that%
\begin{align*}
\left\vert \mathcal{I}_{X}f_{d}\left(  x\right)  -f_{d}\left(  x\right)
\right\vert  & \leq\left(  \frac{2}{\lambda}\left\Vert f_{d}\right\Vert
_{\infty,\Omega}+2\left\Vert Df_{d}\right\Vert _{\infty,\Omega}\right)
h_{\Omega,X}\\
& =2\left(  \frac{1}{\lambda}\left\Vert f_{d}\right\Vert _{\infty,\Omega
}+\left\Vert Df_{d}\right\Vert _{\infty,\Omega}\right)  h_{\Omega,X},
\end{align*}

and Corollary \ref{Cor_Thm_err_interpol_H1inf_data} implies that%
\[
\left\Vert \mathcal{I}_{X}f_{d}\right\Vert _{1,\infty,\Omega}\leq C\left\Vert
f_{d}\right\Vert _{1,\infty,\Omega},\quad f_{d}\in W^{1,\infty}\left(
\Omega\right)  ,
\]

where
\[
C=\max\left\{  1+\frac{2}{\lambda}\left(  1+h_{\Omega,X}\right)
,1+2h_{\Omega,X}\right\}  .
\]

\textbf{This interpolant would be more interesting in higher dimensions and
more so for small }$\lambda$.

We are really more interested in the corresponding basis function smoother -
see Section \ref{Sect_Ex_smth_err_H1inf_data} and Example
\ref{Ex_ex_smth_hat_larg_supp_1dim}.
\end{example}

\section{Multivariate local error estimates on $\mathbb{R}^{d}$%
\label{Lem_interpol_estim_multivar_W1inf_Xow}}

We will now derive a \textbf{multivariate analogue} of Theorem
\ref{Thm_err_interpol_H1inf_data} but first we need:

\begin{lemma}
\label{Lem_cone_x_xk}Suppose $\Omega$ is a bounded region satisfying the
(uniform, open) cone condition and the maximum spherical cavity size (radius)
is
\[
h_{\Omega,X}:=\sup_{x\in\Omega}\min\limits_{k=1}^{N}\left\vert x-x^{\left(
k\right)  }\right\vert .
\]

Then there exist constants $c_{\Omega},\varepsilon_{0}>0$ such that: If
$h_{\Omega,X}<\varepsilon_{0}$ then for every $x\in\Omega$ there exists
$x^{\left(  k\right)  }$ such that $\left[  x,x^{\left(  k\right)  }\right]
\subset\Omega$ and $\left\vert x-x^{\left(  k\right)  }\right\vert \leq
c_{\Omega}h_{\Omega,X}$.

More precisely, if $r$ is the radius of the ball of the cone and $l$ is the
distance between the center of the ball and the vertex then we can choose
$\varepsilon_{0}=r$ and $c_{\Omega}=r+l$. Note that $l>r$.
\end{lemma}

\begin{proof}
Choose $\varepsilon_{0}\leq r$. Since $\Omega$ satisfies the cone condition
$x$ lies in some open cone $C\subset\Omega$ and $h_{\Omega,X}\leq
\varepsilon_{0}$ implies there exists an independent data point $x^{\left(
k\right)  }\in C$. This means $\left[  x,x^{\left(  k\right)  }\right]
\subset C\subset\Omega$. Simple geometrical considerations indicate that
$\left\vert x-x^{\left(  k\right)  }\right\vert \leq\left(  \varepsilon
_{0}+l\right)  h_{\Omega,X}$.
\end{proof}

I do not have a multivariate example for the next theorem because I have not
yet been able to prove that the scaled tensor product hat function with large
support satisfies the inequality \ref{a4.60} below. But someone said that hope
springs eternal.

\begin{theorem}
\label{Thm_err_interpol_H1inf_data multivar}\textbf{Multivariate
interpolation} Suppose $\Omega\subset\mathbb{R}^{d}$ is a bounded data region
satisfying the cone condition and suppose that the basis function satisfies
$G\in W^{1,\infty}\left(  \Omega-\Omega\right)  $. We will choose our data
functions $f_{d}\in W^{1,\infty}\left(  \Omega\right)  \cap X_{w}^{0}\left(
\Omega\right)  $ and assume that there exists a constant $c>0 $ such that the
basis function interpolant $\mathcal{I}_{X}f_{d}$ satisfies%
\begin{equation}
\left\Vert D^{\alpha}\mathcal{I}_{X}f_{d}\right\Vert _{\infty,\Omega}\leq
c\left\Vert f_{d}\right\Vert _{1,\infty,\Omega},\quad x\in\Omega,\text{ }%
f_{d}\in W^{1,\infty}\left(  \Omega\right)  \cap X_{w}^{0}\left(
\Omega\right)  ,\text{ }\left\vert \alpha\right\vert =1.\label{a4.60}%
\end{equation}

Then we have the following pointwise error estimate:
\[
\left\vert \mathcal{I}_{X}f_{d}\left(  x\right)  -f_{d}\left(  x\right)
\right\vert \leq\sqrt{d}\left(  1+c\right)  \left\Vert f_{d}\right\Vert
_{1,\infty,\Omega}\left\vert x-x^{\left(  k\right)  }\right\vert ,\quad
x\in\Omega,\text{ }x^{\left(  k\right)  }\in\Omega.
\]

Further, there exist constants $c_{\Omega},\varepsilon_{0}>0$ such that we
have the \textbf{spherical cavity} error estimate:%
\[
\left\vert \mathcal{I}_{X}f_{d}\left(  x\right)  -f_{d}\left(  x\right)
\right\vert \leq\sqrt{d}\left(  1+c\right)  c_{\Omega}\left\Vert
f_{d}\right\Vert _{1,\infty,\Omega}h_{\Omega,X},\quad x\in\Omega,
\]

when $h_{\Omega,X}<\varepsilon_{0}$.
\end{theorem}

\begin{proof}
Since $f_{d}\in X_{w}^{0}\left(  \Omega\right)  $ we have $f_{d}\in
C_{B}^{\left(  0\right)  }\left(  \Omega\right)  $. If $x\in\Omega$ then by
Lemma \ref{Lem_cone_x_xk} there exist constants $r_{\Omega},\varepsilon_{0}>0$
such that when $h_{\Omega,X}<\varepsilon_{0}$ there exists $x^{\left(
k\right)  }\in X$ such that $\left[  x,x^{\left(  k\right)  }\right]
\subset\Omega$ and $\left\vert x-x^{\left(  k\right)  }\right\vert <c_{\Omega
}h_{\Omega,X}$. This permits us to apply the Taylor series expansion result of
Remark \ref{Rem_Lem_Taylor_extension} with $x\in\Omega$, $z=x^{\left(
k\right)  }\in X$ and $b=x-x^{\left(  k\right)  }$ to get%
\begin{align*}
\mathcal{I}_{X}f_{d}\left(  x\right)  -f_{d}\left(  x\right)   &  =\left(
\mathcal{I}_{X}f_{d}-f_{d}\right)  \left(  x^{\left(  k\right)  }+\left(
x-x^{\left(  k\right)  }\right)  \right) \\
&  =\left(  \mathcal{I}_{X}f_{d}-f_{d}\right)  \left(  x^{\left(  k\right)
}\right)  +\left(  \mathcal{R}_{1}\left(  \mathcal{I}_{X}f_{d}-f_{d}\right)
\right)  \left(  x^{\left(  k\right)  },x-x^{\left(  k\right)  }\right) \\
&  =\left(  \mathcal{R}_{1}\left(  \mathcal{I}_{X}f_{d}-f_{d}\right)  \right)
\left(  x^{\left(  k\right)  },x-x^{\left(  k\right)  }\right)  ,
\end{align*}

with the remainder estimate \ref{1.18}, namely
\begin{align*}
\left\vert \left(  \mathcal{R}_{1}\left(  \mathcal{I}_{X}f_{d}-f_{d}\right)
\right)  \left(  x^{\left(  k\right)  },x-x^{\left(  k\right)  }\right)
\right\vert  &  \leq\sqrt{d}\left(  \max_{\left\vert \alpha\right\vert =1}%
\max_{y\in\left[  x^{\left(  k\right)  },x\right]  }\left\vert D^{\alpha
}\left(  \mathcal{I}_{X}f_{d}-f_{d}\right)  (y)\right\vert \right)  \left\vert
x-x^{\left(  k\right)  }\right\vert \\
&  \leq\sqrt{d}\max_{\left\vert \alpha\right\vert =1}\left\Vert D^{\alpha
}\left(  \mathcal{I}_{X}f_{d}-f_{d}\right)  \right\Vert _{\infty,\Omega
}\left\vert x-x^{\left(  k\right)  }\right\vert \\
&  \leq\sqrt{d}\left(  \max_{\left\vert \alpha\right\vert =1}\left\Vert
D^{\alpha}f_{d}\right\Vert _{\infty,\Omega}+\max_{\left\vert \alpha\right\vert
=1}\left\Vert D^{\alpha}\mathcal{I}_{X}f_{d}\right\Vert _{\infty,\Omega
}\right)  \left\vert x-x^{\left(  k\right)  }\right\vert .
\end{align*}

Hence, using the assumption \ref{a4.60},%
\begin{align}
\left\vert \mathcal{I}_{X}f_{d}\left(  x\right)  -f_{d}\left(  x\right)
\right\vert  & \leq\sqrt{d}\left(  \max_{\left\vert \alpha\right\vert
=1}\left\Vert D^{\alpha}f_{d}\right\Vert _{\infty,\Omega}+c\left\Vert
f_{d}\right\Vert _{1,\infty,\Omega}\right)  \left\vert x-x^{\left(  k\right)
}\right\vert \nonumber\\
& \leq\sqrt{d}\left(  1+c\right)  \left\Vert f_{d}\right\Vert _{1,\infty
,\Omega}\left\vert x-x^{\left(  k\right)  }\right\vert ,\label{a551}%
\end{align}

and thus%
\begin{align*}
\left\vert \mathcal{I}_{X}f_{d}\left(  x\right)  -f_{d}\left(  x\right)
\right\vert  & \leq\sqrt{d}\left(  1+c\right)  \left\Vert f_{d}\right\Vert
_{1,\infty,\Omega}\min_{k=1}^{N}\left\vert x-x^{\left(  k\right)  }\right\vert
\\
& \leq\sqrt{d}\left(  1+c\right)  \left\Vert f_{d}\right\Vert _{1,\infty
,\Omega}\sup_{x\in\overline{\Omega}}\min_{k=1}^{N}\left\vert x-x^{\left(
k\right)  }\right\vert \\
& \leq\sqrt{d}\left(  1+c\right)  \left\Vert f_{d}\right\Vert _{1,\infty
,\Omega}c_{\Omega}h_{\Omega,X}.
\end{align*}

\end{proof}

\begin{corollary}
Under the conditions of Theorem \ref{Thm_err_interpol_H1inf_data multivar} the
interpolant operator $\mathcal{I}_{X}$ satisfies
\[
\left\Vert \mathcal{I}_{X}f_{d}\right\Vert _{1,\infty,\Omega}\leq\left\{
1+\sqrt{d}\left(  1+c\right)  c_{\Omega}h_{\Omega,X},c\right\}  \left\Vert
f_{d}\right\Vert _{1,\infty,\Omega},\quad f_{d}\in W^{1,\infty}\left(
\Omega\right)  \cap X_{w}^{0}\left(  \Omega\right)  ,
\]

when $h_{\Omega,X}<\varepsilon_{0}$.
\end{corollary}

\begin{proof}
From \ref{a4.60} and then \ref{a551} we get:%
\begin{align*}
\left\Vert \mathcal{I}_{X}f_{d}\right\Vert _{1,\infty,\Omega} &  =\max\left\{
\left\Vert \mathcal{I}_{X}f_{d}\right\Vert _{\infty,\Omega},\max_{\left\vert
\alpha\right\vert =1}\left\Vert D^{\alpha}\mathcal{I}_{X}f_{d}\right\Vert
_{1,\infty,\Omega}\right\} \\
&  \leq\max\left\{  \left\Vert \mathcal{I}_{X}f_{d}\right\Vert _{\infty
,\Omega},c\left\Vert f_{d}\right\Vert _{1,\infty,\Omega}\right\} \\
&  \leq\max\left\{  \left\Vert f_{d}\right\Vert _{\infty,\Omega}+\sqrt
{d}\left(  1+c\right)  c_{\Omega}\left\Vert f_{d}\right\Vert _{1,\infty
,\Omega}h_{\Omega,X},c\left\Vert f_{d}\right\Vert _{1,\infty,\Omega}\right\}
\\
&  \leq\left\{  1+\sqrt{d}\left(  1+c\right)  c_{\Omega}h_{\Omega
,X},c\right\}  \left\Vert f_{d}\right\Vert _{1,\infty,\Omega}.
\end{align*}

\end{proof}

\chapter{Central difference tensor product weight functions and interpolation
\label{Ch_cent_diff_wt_fn_ten_prod}}

\section{Introduction}

In \textbf{Section} \ref{Sect_wt_fn_central_diff} We define the univariate and
tensor product central difference weight functions using a central difference
operator. Various bounds and smoothness properties are proved for the weight
function. These depend on the assumptions placed upon parameters used to
define the weight function.

In \textbf{Section} \ref{Sect_tenprod_centdiff_basis_fns} we will derive
several formulas for the basis function and its derivatives, as well as
deriving some smoothness properties and bounds. For example, Theorem
\ref{Thm_G_basis_def_2} presents a \textbf{multiplicative convolution formula}
for the basis function in terms of $q$ and the extended B-spline basis
function $G_{s}$. Theorem \ref{Thm_cdiffbasis_part_moment_formula} gives a
basis function formula in terms of the \textbf{partial moments} of $q$.

In \textbf{Section} \ref{Sect_cent_diff_basis_conv_nat_splin_basis} we show
that given an extended B-spline weight function $w_{s}$ and corresponding
basis function $G_{s}$, there exists a sequence of central difference weight
functions $w_{k}$ and corresponding basis functions $G_{k}$ such that
$1/w_{k}\rightarrow1/w_{s}$ in the $L^{1}$ sense and $G_{k}\rightarrow G_{s}$
uniformly pointwise.

In \textbf{Section} \ref{Sect_CntDifWtFn_DataFuncs} we will use the
generalized local data function results of Subsection
\ref{SbSect_loc_data_larger_class} to characterize locally the data functions
for the tensor product central difference basis functions as Sobolev spaces.
These results supply important information about the data functions and makes
it easy to choose data functions for numerical experiments concerning the zero
order basis function interpolation and smoothing problems discussed later in
Chapters \ref{Ch_Interpol}, \ref{Ch_Exact_smth} and \ref{Ch_Approx_smth}. It
also allows us to use the result of Subsection
\ref{SbSect_global_error_to_local} to derive local convergence estimates from
global convergence estimates for interpolants and smoothers.

In \textbf{Section} \ref{Sect_CntlDiffWtFn_InterpolConverg} we will use the
results of \ref{Ch_Interpol} to derive \textbf{orders for the pointwise
converge} of the minimal norm interpolant to its data function in the case of
central difference weight functions.

To understand the theory of interpolants and smoothers presented in this
document it is not necessary to read this chapter which introduces the
\textbf{tensor product central difference weight functions}. However, we note
that Section \ref{Sect_CntDifWtFn_DataFuncs} contains local data space results
specifically designed for tensor product weight functions and that central
difference weight functions are used as examples in the chapters dealing with smoothers.

This chapter introduces a large new class of tensor product functions based on
1-dimensional central difference operators and it is shown that they are
weight functions. The basis functions are calculated and smoothness estimates
and upper and lower bounds are derived for the weight functions and the basis functions.

Like the B-splines, their data functions are characterized locally as Sobolev spaces.

Pointwise convergence results are derived for the basis function interpolant
using results from Chapter \ref{Ch_Interpol}. These give orders of convergence
identical to those obtained for the extended B-splines in. This chapter
presents no numerical interpolation experiments.

This is a large class of weight functions and I have only calculated a couple
of basis functions. The calculations were tedious until I proved Theorem
\ref{Thm_cdiffbasis_part_moment_formula} which gives much easier formulas for
the basis function. I have no results to offer in this regard except Remark
\ref{Rem_Thm_G_basis_def_2} which might mean that it may be possible to
usefully characterize the basis functions or characterize a large subset of
the basis functions.

However, lately I have been able to derive more amenable formulas for the
basis functions - this is Theorem \ref{Thm_cdiffbasis_part_moment_formula}
which calculates the basis function using partial moments of the generating
function $q$.

The central difference weight functions are defined as follows: Suppose $q\in
L^{1}\left(  \mathbb{R}^{1}\right)  ,$ $q\neq0,$ $q\left(  \xi\right)  \geq0$
and $l\geq n\geq1$ are integers. The univariate \textbf{central difference
weight function} with parameters $n,l$ is defined by
\[
w\left(  \xi\right)  =\frac{\xi^{2n}}{\Delta_{2l}\widehat{q}\left(
\xi\right)  },\quad\xi\in\mathbb{R}^{1},
\]

where $\Delta_{2l}$ is central difference operator%
\[
\Delta_{2l}f\left(  \xi\right)  =\sum_{k=-l}^{l}\left(  -1\right)  ^{k}%
\tbinom{2l}{k+l}f\left(  -k\xi\right)  ,\text{\quad}l=1,2,3,\ldots;\text{ }%
\xi\in\mathbb{R}^{1}.
\]

For example, $\Delta_{2}f\left(  \xi\right)  =-\left(  f\left(  \xi\right)
-2f\left(  0\right)  +f\left(  -\xi\right)  \right)  $. The multivariate
central difference weight function is defined by tensor product.

It is shown in Theorem \ref{Thm_cdiffwt_2} that $w$ belongs to the class of
zero order weight functions introduced in Chapter \ref{Ch_wtfn_basisfn_datasp}
for some $\kappa$ iff $\int\limits_{\left\vert \xi\right\vert \geq
R}\left\vert \xi\right\vert ^{2n-1}q\left(  \xi\right)  d\xi<\infty$ for some
$R\geq0$. Here $\kappa$ satisfies $\kappa+1/2<n$. The central difference
weight functions are closely related to the \textbf{extended B-splines}
defined by \ref{1.032} and a discussion of their genesis is given in
Subsection \ref{SbSect_cent_Motivation}.

Various bounds are derived for $w$. For example, in Corollary
\ref{Cor_cdiffwt_bnd_on_wt_fn} it is shown that if $\int\nolimits_{\left\vert
t\right\vert \geq R}t^{2l}q\left(  t\right)  dt<\infty$ then for any $r>0$
there exist constants $c_{r},c_{r}^{\prime},k_{r},k_{r}^{\prime}>0$ such that%
\begin{align*}
k_{r}\xi^{2n}  & \leq w\left(  \xi\right)  \leq k_{r}^{\prime}\xi^{2n}%
,\quad\left\vert \xi\right\vert \geq r,\text{ }\xi\in\mathbb{R}^{1},\\
\frac{c_{r}}{\xi^{2\left(  l-n\right)  }}  & \leq w\left(  \xi\right)
\leq\frac{c_{r}^{\prime}}{\xi^{2\left(  l-n\right)  }},\quad\left\vert
\xi\right\vert \leq r,\text{ }\xi\in\mathbb{R}^{1}.
\end{align*}

By Theorems \ref{Thm_G_basis_def_2} and \ref{Thm2_G_basis_def_2_new} the
univariate central difference basis function is given by the multiplicative
convolution formula%
\[
G_{c}\left(  s\right)  =\left(  -1\right)  ^{\left(  l-n\right)  }%
\int\limits_{\mathbb{R}^{1}}\left(  D^{2\left(  l-n\right)  }\left(
\ast\Lambda\right)  ^{l}\right)  \left(  \frac{s}{t}\right)  \left\vert
t\right\vert ^{2n-1}q\left(  t\right)  dt,\text{\quad}s\in\mathbb{R}^{1},
\]

with $G_{c}\in C_{B}^{\left(  2n-2\right)  }$ and $D^{2n-1}G$ bounded in
$C^{\left(  2n-1\right)  }\left(  \mathbb{R}^{1}\setminus0\right)  $. The
multivariate basis function is defined as a tensor product. For $k\leq2n-2$,
$D^{k}G_{c}$ is uniformly Lipschitz continuous of order $1$ (Theorem
\ref{Thm_cdiffbasis_Lips_dim_d}).

The next theorem Theorem \ref{Thm_cdiffbasis_part_moment_formula} and its
corollary give convenient \textbf{partial moment formulas} for the basis
function and its derivatives that do not involve calculating the extended
B-spline basis functions $G_{s}$ but instead involves the direct calculation
of partial moments of $q$\textbf{. }In fact%
\begin{align*}
G_{c}\left(  x\right)   & =\left\{
\begin{array}
[c]{ll}%
\frac{\left(  -1\right)  ^{n}}{\left(  2n-1\right)  !}\int_{\left\vert
x\right\vert }^{R_{q}l}\left(  s-\left\vert x\right\vert \right)  ^{2n-1}%
q_{l}\left(  s\right)  ds, & \left\vert x\right\vert \leq R_{q}l,\\
0, & \left\vert x\right\vert \geq R_{q}l,
\end{array}
\right. \\
q_{l}\left(  s\right)   & =\sum\limits_{j=-l,j\neq0}^{l}\frac{\left(
-1\right)  ^{j}}{\left\vert j\right\vert }\tbinom{2l}{j+l}q\left(  \frac{s}%
{j}\right)  ,
\end{align*}

where $\operatorname*{supp}q\subseteq\overline{B}_{R_{q}}$.

Also, the improved result that $G_{c}\in C_{B}^{2n-1}\left(  \mathbb{R}%
^{1}\right)  $ is derived. It is also shown that if $\operatorname*{supp}%
q\subseteq\overline{B}_{R_{q}}$, where possibly $R_{q}=\infty$, it follows
that $\operatorname*{supp}G_{c}\subseteq\overline{B}_{R_{q}l}$ and that
$G_{c}$ satisfies the differential equation
\[
D^{2n}G_{c}\left(  s\right)  =\left(  -1\right)  ^{n}\tbinom{2l}{l}\left(
\int q\right)  \delta+\left(  -1\right)  ^{n}\sum\limits_{j=-l,j\neq0}%
^{l}\frac{\left(  -1\right)  ^{j}}{\left\vert j\right\vert }\tbinom{2l}%
{j+l}q\left(  \frac{s}{j}\right)  ,
\]

implying that $D^{2n}G_{c}$ has the same smoothness as $q$ on $\mathbb{R}%
^{d}\setminus0$. In fact if $q\in C_{0}^{\infty}\left(  \mathbb{R}^{1}\right)
$ then $\ G_{c}\in C^{\infty}\left(  \mathbb{R}^{1}\setminus0\right)  \cap
C_{B}^{\left(  2n-1\right)  }\left(  \mathbb{R}^{1}\right)  $ and $G_{c}$ has
bounded support.

In Section \ref{Sect_CntDifWtFn_DataFuncs} we characterize the \textbf{data
function space} locally as the Sobolev space $W^{n\mathbf{1}}\left(
\Omega\right)  $ (defined by \ref{1.062}). A localization result specifically
for tensor product weight functions is derived in Theorem
\ref{Thm_data_fn_tensor_prod}. This supplies important information about the
data functions and makes it easy to choose data functions for numerical experiments.

Results for the pointwise convergence of the interpolant to its data function
on a bounded data domain are derived using results from Chapter
\ref{Ch_wtfn_basisfn_datasp}. These results are summarized in Table
\ref{Tbl_ConvergCentral} and give estimates identical to those obtained for
the extended B-splines.

This chapter does not present the results of any numerical experiments.

\section{Central difference weight functions\label{Sect_wt_fn_central_diff}}

In this section we define the univariate and tensor product central difference
weight functions using a central difference operator. Various bounds and
smoothness properties are proved for the weight function. These depend on the
assumptions placed upon parameters used to define the weight function.

To introduce the central difference weight functions we will need some central
difference operators:

\begin{definition}
\label{Def_op_diff_real}\textbf{The univariate central difference operators
}$\delta_{\xi}$ \textbf{and} $\Delta_{2l}$ \textbf{on} $\mathbb{R}^{1}$.

In Theorem \ref{Thm_data_func}\ we defined the central difference operator%
\[
\delta_{\xi}f\left(  x\right)  =f\left(  x+\frac{\xi}{2}\right)  -f\left(
x-\frac{\xi}{2}\right)  ,\text{\quad}x,\xi\in\mathbb{R}^{1}.
\]

and here we introduce the operator%
\begin{equation}
\Delta_{2l}f\left(  \xi\right)  =\sum_{j=-l}^{l}\left(  -1\right)  ^{j}%
\tbinom{2l}{j+l}f\left(  -j\xi\right)  ,\text{\quad}l=1,2,3,\ldots,\text{ }%
\xi\in\mathbb{R}^{1}.\label{a914}%
\end{equation}

\textbf{This definition of }$\Delta_{2l}$\textbf{\ ensures that}
$\Delta_{2l,\xi}\left(  e^{-i\xi t}\right)  \geq0$, which is part 3 of the
next lemma. For example:%
\begin{align*}
\Delta_{2}f\left(  \xi\right)   & =-\left(  f\left(  \xi\right)  -2f\left(
0\right)  +f\left(  -\xi\right)  \right)  .\\
\Delta_{4}f\left(  \xi\right)   & =f\left(  2\xi\right)  -4f\left(
\xi\right)  +6f\left(  0\right)  -4f\left(  -\xi\right)  +f\left(
-2\xi\right)  .
\end{align*}

Noting that $\Delta_{2l}f\left(  \xi\right)  =\left(  -1\right)  ^{l}\left(
\delta_{\xi}^{2l}f\right)  \left(  0\right)  $ we define%
\[
\Delta_{k}f\left(  \xi\right)  :=i^{k}\left(  \delta_{\xi}^{k}f\right)
\left(  0\right)  ,\text{\quad}k=0,1,2,\ldots
\]

\end{definition}

\begin{lemma}
\label{Lem_central_diff_op}The \textbf{even} central difference operators of
Definition \ref{Def_op_diff_real} have the following properties:

\begin{enumerate}
\item $\Delta_{2k}f\left(  \xi\right)  =\left(  -1\right)  ^{k}\left(
\delta_{\xi}^{2k}f\right)  \left(  0\right)  =\left(  \left(  -\delta_{\xi
}\right)  ^{2k}f\right)  \left(  0\right)  $.

\item Regarding monomials and polynomials:

\begin{enumerate}
\item $\Delta_{2l}\left(  \xi^{k}\right)  =\left(  \sum\limits_{j=-l}%
^{l}\left(  -1\right)  ^{j}\tbinom{2l}{j+l}j^{k}\right)  \left(  -\xi\right)
^{k}$ when $k\geq0$.

\item $\sum\limits_{j=-l}^{l}\left(  -1\right)  ^{j}\tbinom{2l}{j+l}%
j^{k}=\left\{
\begin{array}
[c]{ll}%
0, & if\text{ }k=0,1,2,\ldots,2l-1,\\
\left(  -1\right)  ^{l}\left(  2l\right)  !, & if\text{ }k=2l.
\end{array}
\right.  $

\item $\Delta_{2l}\left(  \xi^{2l+1}\right)  =0$.

\item $\Delta_{2l}\left(  p\left(  \xi\right)  \right)  =0$ when $p$ is a
polynomial with degree less than $2l$.

\item $\Delta_{2l}\left(  \xi^{2l}\right)  =\left(  -1\right)  ^{l}\left(
2l\right)  !\xi^{2l}$.
\end{enumerate}

\item $\Delta_{2l,\xi}\left(  e^{-i\xi t}\right)  =2^{2l}\sin^{2l}\left(  \xi
t/2\right)  $. Here $\Delta_{2l,\xi}$ indicates $\Delta_{2l}$ is acting on the
variable $\xi$ in $e^{-i\xi t}$.

\item If $g\in L^{1}$ then part 12 of Lemma \ref{Sum_FourTransf} implies
\[
\Delta_{2l}\widehat{g}\left(  \xi\right)  =\sum_{j=-l}^{l}\left(  -1\right)
^{j}\tbinom{2l}{j+l}\widehat{g}\left(  -j\xi\right)  =\tbinom{2l}%
{l}\widehat{g}\left(  0\right)  +\left(  \sum_{j=-l,\text{ }j\neq0}^{l}%
\frac{\left(  -1\right)  ^{j}}{\left\vert j\right\vert }\tbinom{2l}%
{j+l}g\left(  -\frac{x}{j}\right)  \right)  ^{\wedge}.
\]

\end{enumerate}
\end{lemma}

\begin{remark}
\label{Rem_fin_diff_Taylor}(to self) In Section 1 of \cite{Bojarski2011}
Bojarski states that it is well known that in $\mathbb{R}^{1}$ the finite
difference expression
\begin{equation}
\sum_{j=0}^{m}\left(  -1\right)  ^{j}\tbinom{m}{j}f\left(  x+jh\right)
,\text{\quad}h=\frac{y-x}{m},\label{a9.9}%
\end{equation}

approximates the Taylor series remainder%
\[
f\left(  y\right)  -\sum_{n=0}^{m-1}D^{n}f\left(  x\right)  \frac{\left(
y-x\right)  ^{n}}{n!},\text{\quad}x,y\in\mathbb{R}^{1}.
\]

The remainder is bounded using the classical Hardy-Littlewood maximal function
$M_{\delta}g$ on $\mathbb{R}^{1}$ e.g. Hardy, Littlewood and P\'{o}lya
\cite{HarLitPol52}. After inequality (2.10) Bojarski considers higher
dimensions with the expression \ref{a9.9} retained and the Taylor series
expansion has the usual higher dimensional form. The estimate of the remainder
is carried out by generalizing the Hardy-Littlewood maximal function to higher
dimensions - see Bojarski and Haj\l asz \cite{BojarHaj93} and Bojarski,
Haj\l asz and Strzelecki \cite{BojHajStr2002}. This involves concepts such as
Bessel capacity and Lebesgue points.

Petersen \cite{Petersen83} mentions Lebesgue points in \S 2.3.

How does this fit in with the derivation of the Approximate smoother?

See Section \ref{SbSect_Comments_multiv_cent_wt}.
\end{remark}

\subsection{Motivation and definition\label{SbSect_cent_Motivation}}

The derivation of the 1-dimensional central difference weight function is
based on the approximation of the reciprocal of the B-spline weight function
$w_{s}$ using a mollifier. From \ref{1.032} we have for parameters $n$ and $l
$,
\[
\frac{1}{w_{s}\left(  \xi\right)  }=\frac{\sin^{2l}\xi}{\xi^{2n}},\quad1\leq
n\leq l.
\]

and so by part 3 Lemma \ref{Lem_central_diff_op} with $t=2$, $\frac{1}%
{w_{s}\left(  \xi\right)  }=\frac{1}{2^{2l}}\frac{\Delta_{2l}\left(
e^{-i2\xi}\right)  }{\xi^{2n}}$ or using the Fourier transform $F\left[
\cdot\right]  \left(  \xi\right)  $
\[
\frac{1}{w_{s}\left(  \xi\right)  }=\frac{\sqrt{2\pi}}{2^{2l}}\frac
{\Delta_{2l}F\left[  \delta\left(  x-2\right)  \right]  }{\xi^{2n}}.
\]

We now want to approximate the delta function $\delta$ in the distribution
sense using the sequence of functions $k\psi\left(  kx\right)  $
$k=1,2,3,\ldots$, where $\psi\in L^{1}\left(  \mathbb{R}^{1}\right)  $,
$\int\psi=1$ and $\int\left\vert x\psi\left(  x\right)  \right\vert dx<\infty
$. Since $k\psi\left(  kx\right)  \rightarrow\delta$ as $k\rightarrow\infty$
in the sense of distributions, $k\psi\left(  k\left(  x-2\right)  \right)
\rightarrow\delta\left(  x-2\right)  $ and so as $k\rightarrow\infty$,
\[
\frac{\sqrt{2\pi}}{2^{2l}}\frac{\Delta_{2l}F\left[  k\psi\left(  k\left(
x-2\right)  \right)  \right]  }{\xi^{2n}}\rightarrow\frac{1}{w_{s}\left(
\xi\right)  }\text{ }on\text{ }\mathcal{D}^{\prime}\left(  \mathbb{R}%
^{1}\setminus\mathbf{0}\right)  .
\]

Here we will not try to extend the convergence to $S^{\prime}$ but observe
that we can write%
\[
\frac{\sqrt{2\pi}}{2^{2l}}\frac{\Delta_{2l}F\left[  k\psi\left(  k\left(
x-2\right)  \right)  \right]  }{\xi^{2n}}=\frac{\Delta_{2l}\widehat{q_{k}}%
}{\xi^{2n}},
\]

where
\begin{equation}
q_{k}\left(  x\right)  =\frac{\sqrt{2\pi}}{2^{2l}}k\psi\left(  k\left(
x-2\right)  \right)  \in L^{1}\left(  \mathbb{R}^{1}\right)  .\label{a978}%
\end{equation}

Thus $\frac{\Delta_{2l}\widehat{q_{k}}}{\xi^{2n}}\rightarrow\frac{1}%
{w_{s}\left(  \xi\right)  }$ on $\mathcal{D}^{\prime}\left(  \mathbb{R}%
^{1}\setminus\mathbf{0}\right)  $ and if we are lucky $\frac{\xi^{2n}}%
{\Delta_{2l}\widehat{q_{k}}\left(  \xi\right)  }$ may be a weight function.
Indeed, we use this idea to define a central difference weight function and
justify calling it a weight function using the next theorem.

Later in Section \ref{Sect_cent_diff_basis_conv_nat_splin_basis} we will give
results which supply conditions on $\psi$ for which $\frac{\Delta
_{2l}\widehat{q_{k}}}{\xi^{2n}}\rightarrow\frac{1}{w_{s}\left(  \xi\right)  }$
in $L^{1}$ and as tempered distributions, and for which the corresponding
sequence of central difference basis functions converges uniformly pointwise
to the extended B-spline basis function.

\begin{definition}
\label{Def_central_diff_wt_func}\textbf{Central difference weight functions}

Suppose that $q\in L^{1}\left(  \mathbb{R}^{1}\right)  $, $q\neq0$ a.e.,
$q\left(  \xi\right)  \geq0$ and $l,n\geq0$ are integers. The
\textbf{univariate} central difference weight function is defined by
\begin{equation}
w\left(  \xi\right)  :=\frac{\xi^{2n}}{\Delta_{2l}\widehat{q}\left(
\xi\right)  },\quad\xi\in\mathbb{R}^{1}.\label{a962}%
\end{equation}

In $d$ dimensions the \textbf{multivariate central difference weight function}
with parameters $n$ and $l$ is the \textbf{tensor product of the univariate
weight function}.
\end{definition}

The next result justifies these definitions.

\begin{theorem}
\label{Thm_cdiffwt_2}Suppose $w$ is the function on $\mathbb{R}^{1}$
introduced in Definition \ref{Def_central_diff_wt_func}. Then $w$ is an even
function satisfying weight function property W01 for $\mathcal{A}=\left\{
0\right\}  $.

Further, $w$ satisfies property W02 for parameter $\kappa$ iff $l$ and $n$
satisfy $\kappa+1/2<n\leq l$ and $q$ satisfies
\begin{equation}
\int\limits_{\left\vert \xi\right\vert \geq R}\left\vert \xi\right\vert
^{2n-1}q\left(  \xi\right)  d\xi<\infty,\quad for\text{ }some\text{ }%
R\geq0.\label{a963}%
\end{equation}

\textbf{Note} that properties W02 and W03 are identical in one dimension and
for tensor products of central differences we can use Theorem
\ref{Thm_ten_prod_two_wt_fns}.
\end{theorem}

\begin{proof}
The weight function $w$ is an even function since it is clear from equation
\ref{a914} that $\Delta_{2l}\widehat{q}$ is even. From part 3 of Lemma
\ref{Lem_central_diff_op} we have $\Delta_{2l,\xi}e^{-i\xi t}=2^{2l}\sin
^{2l}\left(  \xi t/2\right)  $ so that
\begin{align}
\frac{1}{w\left(  \xi\right)  }=\frac{\Delta_{2l}\widehat{q}\left(
\xi\right)  }{\xi^{2n}} &  =\tfrac{1}{\sqrt{2\pi}}\frac{1}{\xi^{2n}}%
\int\left(  \Delta_{2l,\xi}e^{-i\xi t}\right)  q\left(  t\right)
dt\nonumber\\
&  =\tfrac{2^{2l}}{\sqrt{2\pi}}\frac{1}{\xi^{2n}}\int\sin^{2l}\left(  \xi
t/2\right)  q\left(  t\right)  dt,\label{a943}%
\end{align}

and clearly $w\left(  \xi\right)  >0$ when $\xi\neq0$. Since $q\in L^{1}$
implies $\widehat{q}\in C_{B}^{\left(  0\right)  }$, from the definition of
$w$ we have $w\in C^{\left(  0\right)  }\left(  \mathbb{R}^{1}\setminus
\mathbf{0}\right)  $ and $w>0$ on $\mathbb{R}^{1}\setminus\mathbf{0}$.

Regarding property W02, choose $\lambda\geq0$. Then
\begin{align*}
\int\frac{\xi^{2\lambda}d\xi}{w\left(  \xi\right)  }  & =\tfrac{2^{2l}}%
{\sqrt{2\pi}}\int\frac{1}{\xi^{2\left(  n-\lambda\right)  }}\int\sin
^{2l}\left(  \xi t/2\right)  q\left(  t\right)  dtd\xi\\
& =\tfrac{2^{2l}}{\sqrt{2\pi}}\int\int\frac{\sin^{2l}\left(  \xi t/2\right)
}{\xi^{2\left(  n-\lambda\right)  }}d\xi\,q\left(  t\right)  dt.
\end{align*}

where, because the integrand is non-negative, the integrals all exist iff the
last integral exists. Now make the change of variables $\eta=\xi t/2$,
$d\eta=\frac{1}{2}\left\vert t\right\vert d\xi$ so that%
\begin{align}
\int\frac{\xi^{2\lambda}d\xi}{w\left(  \xi\right)  }  & =\tfrac{2^{2l}}%
{\sqrt{2\pi}}\int\int\frac{\sin^{2l}\eta}{\left(  \frac{2\eta}{t}\right)
^{2\left(  n-\lambda\right)  }}\frac{2d\eta}{\left\vert t\right\vert
}\,q\left(  t\right)  dt\nonumber\\
& =\tfrac{2^{2\left(  l-n+\lambda\right)  +1}}{\sqrt{2\pi}}\int\int\frac
{\sin^{2l}\eta}{\eta^{2\left(  n-\lambda\right)  }}d\eta\,\left\vert
t\right\vert ^{2\left(  n-\lambda\right)  -1}q\left(  t\right)  dt\nonumber\\
& =\tfrac{2^{2\left(  l-n+\lambda\right)  +1}}{\sqrt{2\pi}}\int\frac
{\eta^{2\lambda}}{w_{s}\left(  \eta\right)  }d\eta\int\left\vert t\right\vert
^{2\left(  n-\lambda\right)  -1}q\left(  t\right)  dt.\label{a944}%
\end{align}

Since $w_{s}$ is an extended B-spline weight function the first integral of
\ref{a944} exists iff $0\leq\lambda\leq\kappa$ and $\kappa+1/2<n\leq l$.
Further, since $q\in L^{1}$, $\int\left\vert t\right\vert ^{2\left(
n-\lambda\right)  -1}q\left(  t\right)  dt$ exists iff $\int%
\limits_{\left\vert t\right\vert \geq R}\left\vert t\right\vert ^{2\left(
n-\lambda\right)  -1}q\left(  t\right)  dt<\infty$ for some $R\geq0$ iff
$\int\limits_{\left\vert t\right\vert \geq R}\left\vert t\right\vert
^{2\left(  n-\lambda\right)  -1}q\left(  t\right)  dt<\infty$ for $\lambda=0$
and some $R\geq0$. Thus it can be seen that $w$ has property W02 for $\kappa$
iff $\kappa+1/2<n\leq l$ and $\int\limits_{\left\vert t\right\vert \geq
R}\left\vert t\right\vert ^{2n-1}q\left(  t\right)  dt<\infty$.
\end{proof}

\begin{remark}
Subsection 1.2.9 of the positive order document Williams
\cite{WilliamsPosOrdSmthV3} gives necessary and sufficient conditions for a
homogeneous tensor product central difference weight function $w$ to be a
weight function of positive order which satisfies weight function condition
W3.1. In fact, if the weight function is generated by parameters $n,l,q\left(
\cdot\right)  $ then $w$ has property W3.1 for $\theta=\left\vert
\alpha\right\vert \geq1$ and parameter $\kappa=\kappa_{1}\mathbf{1}$ iff%
\[
\kappa_{1}-n+1/2<\underline{\alpha}\leq\overline{\alpha}\leq l-n.
\]

Here $\underline{\alpha}=\min\alpha$ and $\overline{\alpha}=\max\alpha$.

The weight function condition W3.1 is suitable for tensor products and is :
$w$ has property W3.1 for order $\theta$ and some parameter $\kappa
\in\mathbb{R}^{d}$, $\kappa\geq\mathbf{0}$ if there exists a multi-index
$\alpha$ such that $\left\vert \alpha\right\vert =\theta$ and%
\[
\int\dfrac{x^{2\lambda}}{w\left(  x\right)  x^{2\alpha}}dx<\infty,\quad
0\leq\lambda\leq\kappa.
\]

Here $x^{2\lambda}:=\left(  x_{1}^{2}\right)  ^{\lambda_{1}}\times\cdots
\times\left(  x_{d}^{2}\right)  ^{\lambda_{d}}$.

In this case the order $\theta$ basis function $G$ exists and $G\in
C_{B}^{\left(  \left\lfloor 2\kappa_{1}\right\rfloor \mathbf{1}\right)  }$.
For the space of data functions $X_{w}^{\theta}\hookrightarrow C_{BP}^{\infty
}+C_{B}^{\left(  \left\lfloor \kappa_{1}\right\rfloor \mathbf{1}\right)  }$.

The spaces $C_{B}^{\left(  \alpha\right)  }$ are defined by \ref{1.064} and
the spaces $C_{BP}^{\infty}$ are given in Definition
\ref{Def_Some_basic_spaces}.
\end{remark}

The central difference weight function is related to the extended B-spline
weight function as follows:

\begin{theorem}
\label{Thm_cdiffwt_Bspline}Suppose $w_{s}$ is the extended B-spline weight
function with parameters $n,l$. Suppose $w_{c}$ is a central difference weight
function with parameters $n,l,q\left(  \cdot\right)  $. Then these weight
functions are related by the equation%
\begin{equation}
\frac{1}{w_{c}\left(  \xi\right)  }=\tfrac{2^{2\left(  l-n\right)  }}%
{\sqrt{2\pi}}\int_{\mathbb{R}^{1}}\frac{t^{2n}q\left(  t\right)  }%
{w_{s}\left(  \xi t/2\right)  }dt.\label{a968}%
\end{equation}

\end{theorem}

\begin{proof}
From equation \ref{a943} of the proof of Theorem \ref{Thm_cdiffwt_2}%
\begin{align*}
\frac{1}{w_{c}\left(  \xi\right)  }=\tfrac{1}{\sqrt{2\pi}}\frac{2^{2l}}%
{\xi^{2n}}\int\sin^{2l}\left(  \xi t/2\right)  q\left(  t\right)  dt &
=\tfrac{2^{2\left(  l-n\right)  }}{\sqrt{2\pi}}\int\frac{\sin^{2l}\left(  \xi
t/2\right)  }{\left(  \xi t/2\right)  ^{2n}}t^{2n}q\left(  t\right)  dt\\
&  =\tfrac{2^{2\left(  l-n\right)  }}{\sqrt{2\pi}}\int\frac{t^{2n}q\left(
t\right)  }{w_{s}\left(  \xi t/2\right)  }dt.
\end{align*}

\end{proof}

\begin{remark}
Suppose%
\begin{align*}
\frac{1}{w_{c1}\left(  \xi_{1}\right)  }  & :=\tfrac{2^{2\left(  l-n\right)
}}{\sqrt{2\pi}}\int_{\mathbb{R}^{1}}\frac{t_{1}^{2n}q_{1}\left(  t_{1}\right)
}{w_{s}\left(  \xi_{1}t_{1}/2\right)  }dt_{1},\\
\frac{1}{w_{c2}\left(  \xi_{2}\right)  }  & :=\tfrac{2^{2\left(  l-n\right)
}}{\sqrt{2\pi}}\int_{\mathbb{R}^{1}}\frac{t_{2}^{2n}q_{2}\left(  t_{2}\right)
}{w_{s}\left(  \xi_{2}t_{2}/2\right)  }dt_{2}.
\end{align*}

Then the tensor product of these weight functions satisfies%
\begin{align*}
\frac{1}{w_{c1}\left(  \xi_{1}\right)  }\frac{1}{w_{c2}\left(  \xi_{2}\right)
}  & =\tfrac{2^{2\left(  l-n\right)  }}{\sqrt{2\pi}}\int_{\mathbb{R}^{1}}%
\frac{t_{1}^{2n}q_{1}\left(  t_{1}\right)  }{w_{s}\left(  \xi_{1}%
t_{1}/2\right)  }dt_{1}\tfrac{2^{2\left(  l-n\right)  }}{\sqrt{2\pi}}%
\int_{\mathbb{R}^{1}}\frac{t_{2}^{2n}q_{2}\left(  t_{2}\right)  }{w_{s}\left(
\xi_{2}t_{2}/2\right)  }dt_{2}\\
& =\left(  \tfrac{2^{2\left(  l-n\right)  }}{\sqrt{2\pi}}\right)  ^{2}%
\int_{\mathbb{R}^{2}}\frac{t^{2n\mathbf{1}}q_{1}{\small \otimes}q_{2}\left(
\tau\right)  }{w_{s}\left(  \xi\centerdot\tau/2\right)  }d\tau.
\end{align*}

\end{remark}

\subsection{Weight function bounds and smoothness}

We will be mostly concerned with smoothness properties and deriving lower
bounds for the weight function. It turns out that the weight function is
continuous except at the origin where it may have a singularity.

The next theorem will give some properties of the central difference weight
function if condition \ref{a963} is strengthened to \ref{a911}.

\begin{theorem}
\label{Thm_cdiffwt_4}Suppose $w$ is a central difference weight function with
the properties assumed in Theorem \ref{Thm_cdiffwt_2}. Suppose in addition
\begin{equation}
\int\limits_{\left\vert t\right\vert \geq R}t^{2n}q\left(  t\right)
dt<\infty.\label{a911}%
\end{equation}

Then $1/w\in C_{B}^{\left(  0\right)  }\left(  \mathbb{R}^{1}\right)  $ and
\begin{equation}
\lim\limits_{\xi\rightarrow0}\frac{1}{w\left(  \xi\right)  }=\left\{
\begin{array}
[c]{lll}%
0, & \quad if & n<l,\\
\frac{1}{\sqrt{2\pi}}\int t^{2n}q\left(  t\right)  dt, & \quad if & n=l.
\end{array}
\right. \label{a990}%
\end{equation}

Also, there exists a positive constant $C_{1}$ such that%
\begin{equation}
\frac{1}{w\left(  \xi\right)  }\leq\frac{C_{1}}{\xi^{2n}},\quad\xi
\neq0.\label{a950}%
\end{equation}

\end{theorem}

\begin{proof}
From Theorem \ref{Thm_cdiffwt_2} we know that $w$ is continuous and positive
outside the origin. Thus $1/w$ is also continuous and positive outside the
origin. We now show that $1/w$ is also continuous at the origin by proving the
limits \ref{a990}. The existence of the integral $\int\limits_{\left\vert
t\right\vert \geq R}t^{2n}q\left(  t\right)  dt$ implies $t^{2n}q\left(
t\right)  \in L^{1}$ and so $\widehat{q}\in C_{B}^{\left(  2n\right)  }$ by
part 2 of Lemma \ref{Lem_L1_Fourier_contin}. Hence the Taylor series with
integral remainder of Appendix \ref{Sect_apx_TaylorSeries} can be used to
expand each term of \ref{a914} about the origin up to a polynomial of degree
$2n-1$. By part 3 of Lemma \ref{Lem_central_diff_op} $\Delta^{2l}$ annihilates
polynomials of degree less than $2l$ and so:%
\begin{align}
\Delta^{2l}\widehat{q}\left(  \xi\right)   & =\sum_{k=-l}^{l}\left(
-1\right)  ^{\left\vert k\right\vert }\tbinom{2l}{k+l}\widehat{q}\left(
-k\xi\right) \label{a904}\\
& =\sum_{k=-l}^{l}\left(  -1\right)  ^{\left\vert k\right\vert }\tbinom
{2l}{k+l}\sum_{j<2n}\frac{\left(  D^{j}\widehat{q}\right)  (0)}{j!}\left(
-k\xi\right)  ^{j}+\nonumber\\
& \quad+\sum_{k=-l}^{l}\left(  -1\right)  ^{\left\vert k\right\vert }%
\tbinom{2l}{k+l}\frac{\left(  k\xi\right)  ^{2n}}{\left(  2n-1\right)  !}%
\int_{0}^{1}s^{2n-1}\left(  D^{2n}\widehat{q}\right)  \left(  \left(
s-1\right)  k\xi\right)  ds\nonumber\\
& =\sum_{j<2n}\frac{\left(  D^{j}\widehat{q}\right)  (0)}{j!}\left(
-\xi\right)  ^{j}\sum_{k=-l}^{l}\left(  -1\right)  ^{\left\vert k\right\vert
}\tbinom{2l}{k+l}k^{j}+\nonumber\\
& \quad+\sum_{k=-l}^{l}\left(  -1\right)  ^{\left\vert k\right\vert }%
\tbinom{2l}{k+l}\frac{\left(  -k\xi\right)  ^{2n}}{\left(  2n-1\right)  !}%
\int_{0}^{1}s^{2n-1}\left(  D^{2n}\widehat{q}\right)  \left(  \left(
s-1\right)  k\xi\right)  ds\nonumber\\
& =\sum_{k=-l}^{l}\left(  -1\right)  ^{\left\vert k\right\vert }\tbinom
{2l}{k+l}\frac{\left(  k\xi\right)  ^{2n}}{\left(  2n-1\right)  !}\int_{0}%
^{1}s^{2n-1}\left(  D^{2n}\widehat{q}\right)  \left(  \left(  s-1\right)
k\xi\right)  ds\nonumber\\
& =\frac{\xi^{2n}}{\left(  2n-1\right)  !}\sum_{k=-l}^{l}\left(  -1\right)
^{\left\vert k\right\vert }\tbinom{2l}{k+l}k^{2n}\int_{0}^{1}s^{2n-1}\left(
D^{2n}\widehat{q}\right)  \left(  \left(  s-1\right)  k\xi\right)
ds.\label{a960}%
\end{align}

Thus%
\[
\frac{1}{w\left(  \xi\right)  }=\frac{\Delta^{2l}\widehat{q}\left(
\xi\right)  }{\xi^{2n}}=\frac{1}{\left(  2n-1\right)  !}\sum_{k=-l}^{l}\left(
-1\right)  ^{\left\vert k\right\vert }\tbinom{2l}{k+l}k^{2n}\int_{0}%
^{1}s^{2n-1}\left(  D^{2n}\widehat{q}\right)  \left(  \left(  s-1\right)
k\xi\right)  ds.
\]

Because $D^{2n}\widehat{q}\in C_{B}^{\left(  0\right)  }$ we can apply the
Lebesgue-dominated convergence theorem to the sequence

$\left\{  s^{2n-1}\left(  D^{2n}\widehat{q}\right)  \left(  \left(
1-s\right)  \left(  -\frac{k}{j}\right)  \right)  \right\}  _{j=1}^{\infty}$
to obtain%
\begin{align*}
\lim\limits_{\xi\rightarrow0}\frac{1}{w\left(  \xi\right)  }  & =\frac
{1}{\left(  2n-1\right)  !}\sum_{k=-l}^{l}\left(  -1\right)  ^{\left\vert
k\right\vert }\tbinom{2l}{k+l}k^{2n}\int_{0}^{1}s^{2n-1}\lim\limits_{\xi
\rightarrow0}\left(  D^{2n}\widehat{q}\right)  \left(  \left(  s-1\right)
k\xi\right)  ds\\
& =\frac{1}{2n!}\left(  \sum_{k=-l}^{l}\left(  -1\right)  ^{\left\vert
k\right\vert }\tbinom{2l}{k+l}k^{2n}\right)  \left(  D^{2n}\widehat{q}\right)
\left(  0\right) \\
& =\left\{
\begin{array}
[c]{lll}%
0, & \quad if & n<l,\\
\frac{1}{\sqrt{2\pi}}\int t^{2l}q\left(  t\right)  dt, & \quad if & n=l.
\end{array}
\right.
\end{align*}

where the last step used property 4 of Lemma \ref{Lem_central_diff_op}%
.\medskip

The next thing to prove is inequality \ref{a950}. But from equation
\ref{a943},%
\[
\frac{1}{w\left(  \xi\right)  }=\tfrac{1}{\sqrt{2\pi}}\frac{1}{\xi^{2n}}%
\int\sin^{2l}\left(  \xi t/2\right)  q\left(  t\right)  dt\leq\tfrac{1}%
{\sqrt{2\pi}}\left(  \int q\right)  \frac{1}{\xi^{2n}}.
\]

Thus $1/w$ is bounded outside the unit sphere. Since $1/w$ is continuous
everywhere it now follows that $1/w$ is bounded everywhere.
\end{proof}

The next result places the generally stronger condition \ref{a959} on the
function $q$ in exchange for estimates of the weight function behavior near
the origin.

\begin{theorem}
\label{Thm_cdiffwt_3}Suppose $w$ is a \textbf{univariate} central difference
weight function with parameters $n,l$ and generating function $q\left(
\cdot\right)  $ which satisfy the conditions of Theorem \ref{Thm_cdiffwt_2},
and also suppose that%
\begin{equation}
\int\nolimits_{\left\vert t\right\vert \geq R}t^{2l}q\left(  t\right)
dt<\infty.\label{a959}%
\end{equation}

Then%
\begin{equation}
\lim\limits_{\xi\rightarrow0}\frac{1}{w\left(  \xi\right)  \xi^{2\left(
l-n\right)  }}=\frac{1}{2\sqrt{2\pi}}\int\nolimits_{\mathbb{R}^{1}}%
t^{2l}q\left(  t\right)  dt,\quad n\leq l,\label{a961}%
\end{equation}

and given $a>0$ there exist constants $c_{a},c_{a}^{\prime}>0$ such that%
\begin{equation}
c_{a}\leq w\left(  \xi\right)  \xi^{2\left(  l-n\right)  }\leq c_{a}^{\prime
},\quad\left\vert \xi\right\vert \leq a.\label{a964}%
\end{equation}

Regarding smoothness: $w\in C^{\left(  2l\right)  }\left(  \mathbb{R}%
^{1}\setminus0\right)  $ and $w\left(  \xi\right)  \xi^{2\left(  l-n\right)
}\in C^{\left(  0\right)  }$.
\end{theorem}

\begin{proof}
Assumption \ref{a959} implies that $\widehat{q}\in C_{B}^{\left(  2l\right)
}$. Hence we can set $n=l$ in the proof of equation \ref{a960} of Theorem and
obtain%
\[
\Delta^{2l}\widehat{q}\left(  \xi\right)  =\frac{\xi^{2l}}{\left(
2l-1\right)  !}\sum_{k=-l}^{l}\left(  -1\right)  ^{\left\vert k\right\vert
}\tbinom{2l}{k+l}k^{2l}\int_{0}^{1}s^{2l-1}\left(  D^{2l}\widehat{q}\right)
\left(  \left(  s-1\right)  k\xi\right)  ds,
\]

or, since $\frac{1}{w\left(  \xi\right)  }=\xi^{-2n}\Delta^{2l}\widehat{q}%
\left(  \xi\right)  $ for $n\leq l$,
\begin{align*}
\frac{1}{w\left(  \xi\right)  \xi^{2\left(  l-n\right)  }}  & =\xi^{-2l}%
\Delta^{2l}\widehat{q}\left(  \xi\right) \\
& =\frac{1}{\left(  2l-1\right)  !}\sum_{k=-l}^{l}\left(  -1\right)
^{\left\vert k\right\vert }\tbinom{2l}{k+l}k^{2l}\int_{0}^{1}s^{2l-1}\left(
D^{2l}\widehat{q}\right)  \left(  \left(  s-1\right)  k\xi\right)  ds.
\end{align*}

Following the rest of Theorem \ref{Thm_cdiffwt_2} we obtain the limit
\ref{a961}.

To prove \ref{a964} note that since $w$ is continuous, \ref{a961} implies that
there exists $a_{1}>0$ such that%
\[
\left\vert \frac{1}{w\left(  \xi\right)  \xi^{2\left(  l-n\right)  }}-\frac
{1}{\sqrt{2\pi}}\int t^{2l}q\left(  t\right)  dt\right\vert \leq\frac
{1}{2\sqrt{2\pi}}\int t^{2l}q\left(  t\right)  dt,\quad\left\vert
\xi\right\vert \leq a_{1},
\]

and so%
\[
\frac{1}{2\sqrt{2\pi}}\int t^{2l}q\left(  t\right)  dt\leq\frac{1}{w\left(
\xi\right)  \xi^{2\left(  l-n\right)  }}\leq\frac{3}{2\sqrt{2\pi}}\int
t^{2l}q\left(  t\right)  dt,\quad\left\vert \xi\right\vert \leq a_{1},
\]

or%
\[
\frac{1}{\frac{3}{2\sqrt{2\pi}}\int t^{2l}q\left(  t\right)  dt}\leq w\left(
\xi\right)  \xi^{2\left(  l-n\right)  }\leq\frac{1}{\frac{1}{2\sqrt{2\pi}}\int
t^{2l}q\left(  t\right)  dt}.
\]

Now choose
\[
c_{a}=\frac{1}{\frac{3}{2\sqrt{2\pi}}\int t^{2l}q\left(  t\right)  dt},\quad
c_{a}^{\prime}=3c_{a}.
\]

If $a\leq a_{1}$ we are done. If $a>a_{1}$ we use the continuity of $w$ and
\ref{a950} to prove that there exist $C_{1},C_{2}>0$ such that $0<C_{1}\leq
w\left(  \xi\right)  \leq C_{2}$ for $a_{1}\leq\left\vert \xi\right\vert \leq
a$. The inequalities \ref{a964} then follow easily.

Regarding smoothness, we noted at the start of this proof that $\widehat{q}\in
C_{B}^{\left(  2l\right)  }$ and consequently $\Delta^{2l}\widehat{q_{k}}\in
C_{B}^{\left(  2l\right)  }$. Further, equation \ref{a943} implies that
$\Delta_{2l}\widehat{q}\left(  \xi\right)  >0$ when $\xi\neq0$. But $w\left(
\xi\right)  =\frac{\xi^{2n}}{\Delta^{2l}\widehat{q}\left(  \xi\right)  }$ and
so $w\in C^{\left(  2l\right)  }\left(  \mathbb{R}\setminus0\right)  $.
Finally, \ref{a961} of this theorem implies that $\lim\limits_{\xi
\rightarrow0}\frac{1}{w\left(  \xi\right)  \xi^{2\left(  l-n\right)  }}$
exists and is positive. Thus

$\lim\limits_{\xi\rightarrow0}w\left(  \xi\right)  \xi^{2\left(  l-n\right)
}$ exists and so the function $w\left(  \xi\right)  \xi^{2\left(  l-n\right)
} $ is continuous.
\end{proof}

A multi-dimensional analogue of Theorem \ref{Thm_cdiffwt_2} is:

\begin{corollary}
\label{Cor_cdiffwt_W02_W03}\textbf{Tensor product central difference weight
functions} Suppose the functions $\left\{  q_{j}\right\}  _{j=1}^{d}$ satisfy
$q_{j}\in L^{1}\left(  \mathbb{R}^{1}\right)  $, $q_{j}\neq0$, $q_{j}\left(
x\right)  \geq0$. Define the tensor product function $w\left(  \xi\right)
=\prod\limits_{j=1}^{d}\frac{\xi_{j}^{2n}}{\Delta^{2l}\widehat{q_{j}}\left(
\xi_{j}\right)  }$ for integers $l,n\geq0$.

Then:

\begin{enumerate}
\item $w$ is a weight function on $\mathbb{R}^{d}$ with $\mathcal{A}%
\subset\bigcup\limits_{j=1}^{d}\left\{  x:x_{j}=0\right\}  $. We call it a
tensor product central difference weight function.\medskip

Further:

\item $w$ \textbf{also satisfies property W02} for parameter $\kappa
\in\mathbb{R}^{1}$ iff $l$ and $n$ also satisfy the conditions of Theorem
\ref{Thm_cdiffwt_2} i.e.%
\[
\kappa+1/2<n\leq l,
\]

and for some $R_{j}>0$,%
\[
\int_{\left\vert t\right\vert \geq R_{j}}\left\vert t\right\vert ^{2n-1}%
q_{j}\left(  t\right)  dt<\infty,\quad1\leq j\leq d.
\]

\item $w$ \textbf{also satisfies property W03} for parameter $\kappa
\in\mathbb{R}^{d}$ iff $l$ and $n$ also satisfy the conditions%
\[
\kappa+1/2<n\mathbf{1}\leq l\mathbf{1},
\]

and for some $R_{j}>0$,%
\[
\int_{\left\vert t\right\vert \geq R_{j}}\left\vert t\right\vert ^{2n-1}%
q_{j}\left(  t\right)  dt<\infty,\quad1\leq j\leq d.
\]

\end{enumerate}
\end{corollary}

\begin{proof}
\textbf{Part 1} By part 1 of Theorem \ref{Thm_ten_prod_two_wt_fns}, $w$ is a
weight function and by equation \ref{2.00} of Theorem
\ref{Thm_ten_prod_two_wt_fns} and Theorem \ref{Thm_cdiffwt_2} $\mathcal{A}$ is
as given.\medskip

\textbf{Part 2} The weight function $w$ \textbf{satisfies property W02} for
$\kappa\in\mathbb{R}^{1}$ since the tensor product of univariate weight
functions having property W02 for $\kappa$ also satisfies property W02 for
$\kappa$. This is part 2 of Theorem \ref{Thm_ten_prod_two_wt_fns}.\medskip

\textbf{Part 3} The weight function $w$ \textbf{satisfies property W03} for
$\kappa\in\mathbb{R}^{d}$ follows from part 3 of Theorem
\ref{Thm_ten_prod_two_wt_fns}.
\end{proof}

\begin{remark}
\label{RemCor_cdiffwt_W02_W03}The weight functions of Corollary
\ref{Cor_cdiffwt_W02_W03} can be further generalized by choosing different
values of $n$ and $l$ for each dimension i.e. $n$ and $l$ are vectors.
\end{remark}

The next result gives some weight function properties which follow directly
from Theorem \ref{Thm_cdiffwt_3}.

\begin{corollary}
Suppose $w$ is a univariate central difference weight function with parameters
$n,l,q\left(  \cdot\right)  $ which satisfy the conditions of Theorem
\ref{Thm_cdiffwt_3}. Suppose in addition that $l=n$. Then $w$ satisfies:

\begin{enumerate}
\item $w\in C^{\left(  0\right)  }\left(  \mathbb{R}^{1}\right)  \cap
C^{\left(  2n\right)  }\left(  \mathbb{R}^{1}\setminus0\right)  $.\smallskip

There exist constants $c_{1},c_{2}>0$ such that:\smallskip

\item $w\left(  \xi\right)  \geq c_{1}$,\quad for all $\xi$.

\item $w\left(  \xi\right)  \geq c_{2}\xi^{2n}$,\quad for all $\xi$.

\item $\lim\limits_{\xi\rightarrow0}w\left(  \xi\right)  =\left(  \frac
{1}{2^{2n}\sqrt{2\pi}}\int t^{2n}q\left(  t\right)  dt\right)  ^{-1}$.

\item Parts 2 and 3 imply $X_{w}^{0}\left(  \mathbb{R}^{1}\right)
\hookrightarrow W^{n}\left(  \mathbb{R}^{1}\right)  $ where $W^{n}\left(
\mathbb{R}^{1}\right)  $ is the Sobolev space of order $n$. This is because,
$1+\xi^{2n}\leq\left(  1+\xi^{2}\right)  ^{n}\leq2^{n}\left(  1+\xi
^{2n}\right)  $.
\end{enumerate}
\end{corollary}

\subsection{Upper bounds for the weight function}

To derive upper bounds on central difference weight functions we first derive
the weight function for the case when $q$ is a rectangular function.

\begin{theorem}
\label{Thm_cdiffwt_rect_phi}Suppose $b>a$ and $q_{R}$ is the rectangular
function%
\begin{equation}
q_{R}\left(  \xi\right)  =\left\{
\begin{array}
[c]{cc}%
1, & \quad\xi\in\left[  a,b\right]  ,\\
0, & \quad\xi\notin\left[  a,b\right]  .
\end{array}
\right. \label{a930}%
\end{equation}

Then $q=q_{R}$ satisfies the conditions of Theorem \ref{Thm_cdiffwt_2} for
$n=l=1$ and the corresponding weight function is%
\begin{equation}
w_{R}\left(  \xi\right)  =\frac{\sqrt{2\pi}}{8}\frac{\xi^{2}}{\frac{b-a}%
{2}-\frac{\cos\left(  \frac{b+a}{2}\xi\right)  \sin\left(  \frac{b-a}{2}%
\xi\right)  }{\xi}},\label{a954}%
\end{equation}

with bounds%
\begin{equation}
\frac{\sqrt{2\pi}}{6\left(  b-a\right)  }\xi^{2}\leq w_{R}\left(  \xi\right)
\leq\frac{\sqrt{2\pi}}{2\left(  b-a\right)  }\xi^{2},\quad\left\vert
\xi\right\vert \geq\frac{4}{b-a}.\label{a952}%
\end{equation}

\end{theorem}

\begin{proof}
From \ref{a943}%
\begin{align*}
\frac{1}{w_{R}\left(  \xi\right)  }=\tfrac{4}{\sqrt{2\pi}}\frac{1}{\xi^{2}%
}\int_{\mathbb{R}^{1}}\left(  \sin^{2}\frac{\xi t}{2}\right)  q_{R}\left(
t\right)  dt &  =\tfrac{4}{\sqrt{2\pi}}\frac{1}{\xi^{2}}\int_{a}^{b}\sin
^{2}\frac{\xi t}{2}dt\\
&  =\tfrac{8}{\sqrt{2\pi}}\frac{1}{\xi^{2}}\int_{a}^{b}\sin^{2}\frac{\xi t}%
{2}dt.
\end{align*}

The integral is an even function of $\xi$ so now assume $\xi>0$. We now use
the change of variables $s=\xi t/2$, $ds=\left(  \xi/2\right)  dt$ with new
range $\left[  a\xi/2,b\xi/2\right]  $. Hence
\begin{align*}
\int_{a}^{b}\sin^{2}\frac{\xi t}{2}dt=\frac{2}{\xi}\int_{a\xi/2}^{b\xi/2}%
\sin^{2}sds &  =\frac{2}{\xi}\left[  \frac{s}{2}-\frac{\sin2s}{4}\right]
_{a\xi/2}^{b\xi/2}\\
&  =\frac{2}{\xi}\left(  \frac{b-a}{4}\xi-\frac{\sin b\xi-\sin a\xi}{4}\right)
\\
&  =\frac{b-a}{2}-\frac{\sin b\xi-\sin a\xi}{2\xi}\\
&  =\frac{b-a}{2}-\frac{1}{\xi}\cos\left(  \frac{b+a}{2}\xi\right)
\sin\left(  \frac{b-a}{2}\xi\right)  ,
\end{align*}

so that%
\[
\frac{1}{w_{R}\left(  \xi\right)  }=\tfrac{8}{\sqrt{2\pi}}\frac{1}{\xi^{2}%
}\left(  \frac{b-a}{2}-\frac{1}{\xi}\cos\left(  \frac{b+a}{2}\xi\right)
\sin\left(  \frac{b-a}{2}\xi\right)  \right)  ,
\]

as required.

Regarding the inequality \ref{a952}, if $\left\vert \xi\right\vert \geq
\frac{4}{b-a}$ then $\frac{1}{\left\vert \xi\right\vert }\leq\frac{b-a}{4}$
and%
\begin{align*}
w_{R}\left(  \xi\right)  =\frac{\sqrt{2\pi}}{8}\frac{\xi^{2}}{\frac{b-a}%
{2}-\frac{\cos\left(  \frac{b+a}{2}\xi\right)  \sin\left(  \frac{b-a}{2}%
\xi\right)  }{\xi}}\leq\frac{\sqrt{2\pi}}{8}\frac{\xi^{2}}{\frac{b-a}{2}%
-\frac{1}{\left\vert \xi\right\vert }} &  \leq\frac{\sqrt{2\pi}}{8}\frac
{\xi^{2}}{\left(  b-a\right)  /4}\\
&  =\frac{\sqrt{2\pi}}{2\left(  b-a\right)  }\xi^{2}.
\end{align*}

Also, since $\frac{1}{\left\vert \xi\right\vert }\leq\frac{b-a}{4}$%
\[
w_{R}\left(  \xi\right)  =\frac{\sqrt{2\pi}}{8}\frac{\xi^{2}}{\frac{b-a}%
{2}-\frac{\cos\left(  \frac{b+a}{2}\xi\right)  \sin\left(  \frac{b-a}{2}%
\xi\right)  }{\xi}}\geq\frac{\sqrt{2\pi}}{8}\frac{\xi^{2}}{\frac{b-a}{2}%
+\frac{1}{\left\vert \xi\right\vert }}\geq\frac{\sqrt{2\pi}}{6\left(
b-a\right)  }\xi^{2}.
\]

\end{proof}

\begin{remark}
?? \textbf{FINISH}! \textbf{Translation of} $q_{R}$: Define $w_{R}\left(
\xi;a,b\right)  :=w_{R}\left(  \xi\right)  $.%
\begin{align*}
&  \frac{1}{w_{R}\left(  \xi;a+\delta,b+\delta\right)  }\\
&  =\tfrac{8}{\sqrt{2\pi}}\frac{1}{\xi^{2}}\int_{a+\delta}^{b+\delta}\sin
^{2}\frac{\xi t}{2}dt\\
&  =\tfrac{8}{\sqrt{2\pi}}\frac{1}{\xi^{2}}\int_{a}^{b}\sin^{2}\frac
{\xi\left(  t+\delta\right)  }{2}dt=\tfrac{8}{\sqrt{2\pi}}\frac{1}{\xi^{2}%
}\int_{a}^{b}\sin^{2}\frac{\xi t+\xi\delta}{2}dt=\\
&  =\tfrac{8}{\sqrt{2\pi}}\frac{1}{\xi^{2}}\int_{a}^{b}\sin^{2}\frac{\xi
t+\xi\delta}{2}dt=\tfrac{4}{\sqrt{2\pi}}\frac{1}{\xi^{2}}\int_{a}^{b}\left(
1-\cos\left(  \xi t+\xi\delta\right)  \right)  dt=\\
&  =\tfrac{4}{\sqrt{2\pi}}\frac{1}{\xi^{2}}\int_{a}^{b}\left(  1-\cos\xi
t\cos\xi\delta+\sin\xi t\sin\xi\delta\right)  dt\\
&  =\tfrac{4}{\sqrt{2\pi}}\frac{1}{\xi^{2}}\int_{a}^{b}\left(  \cos\xi
\delta-\cos\xi t\cos\xi\delta+1-\cos\xi\delta+\sin\xi t\sin\xi\delta\right)
dt\\
&  =\tfrac{4}{\sqrt{2\pi}}\frac{\cos\xi\delta}{\xi^{2}}\int_{a}^{b}\left(
1-\cos\xi t\right)  dt+\tfrac{4}{\sqrt{2\pi}}\frac{1-\cos\xi\delta}{\xi^{2}%
}\int_{a}^{b}dt+\tfrac{4}{\sqrt{2\pi}}\frac{\sin\xi\delta}{\xi^{2}}\int%
_{a}^{b}\sin\xi tdt\\
&  =\frac{1}{w_{R}\left(  \xi;a,b\right)  }+\tfrac{8}{\sqrt{2\pi}}\left(
b-a\right)  \frac{\sin^{2}\frac{\xi\delta}{2}}{\xi^{2}}+\tfrac{4}{\sqrt{2\pi}%
}\frac{\sin\xi\delta}{\xi^{3}}\left[  -\cos\xi t\right]  _{a}^{b}\\
&  =\frac{1}{w_{R}\left(  \xi;a,b\right)  }+\tfrac{8}{\sqrt{2\pi}}\left(
b-a\right)  \frac{\sin^{2}\frac{\xi\delta}{2}}{\xi^{2}}+\tfrac{4}{\sqrt{2\pi}%
}\frac{\sin\xi\delta}{\xi^{3}}\left(  \cos\xi a-\cos\xi b\right) \\
&  =\frac{1}{w_{R}\left(  \xi;a,b\right)  }+\tfrac{8}{\sqrt{2\pi}}\left(
b-a\right)  \frac{\sin^{2}\frac{\delta\xi}{2}}{\xi^{2}}+\tfrac{8}{\sqrt{2\pi}%
}\frac{\sin\delta\xi}{\xi^{3}}\sin\left(  a+b\right)  \xi\sin\left(
b-a\right)  \xi.
\end{align*}

We can now write%
\begin{align*}
\frac{1}{w_{R}\left(  \xi;a+\delta,b+\delta\right)  } &  =\frac{1}%
{w_{R}\left(  \xi;a,b\right)  }+\tfrac{2}{\sqrt{2\pi}}\left(  b-a\right)
\delta^{2}\frac{\sin^{2}\frac{\delta\xi}{2}}{\left(  \frac{\delta\xi}%
{2}\right)  ^{2}}+\\
&  \qquad+\tfrac{32}{\sqrt{2\pi}}\left(  b^{2}-a^{2}\right)  \delta\frac
{\sin\left(  \frac{1}{2}2\delta\xi\right)  }{\frac{1}{2}2\delta\xi}\frac
{\sin\left(  \frac{1}{2}2\left(  a+b\right)  \xi\right)  }{\frac{1}{2}2\left(
a+b\right)  \xi}\frac{\sin\left(  \frac{1}{2}2\left(  b-a\right)  \xi\right)
}{\frac{1}{2}2\left(  b-a\right)  \xi}\\
&  =??\frac{1}{w_{R}\left(  \xi;a,b\right)  }+\tfrac{2??c_{1}}{\sqrt{2\pi}%
}\left(  b-a\right)  \delta^{2}\widehat{\Lambda}\left(  \delta\xi\right)  +\\
&  \qquad+\tfrac{8??c_{2}}{\sqrt{2\pi}}\delta\widehat{R}\left(  2\delta
\xi\right)  \left(  \frac{a+b}{2}\right)  \widehat{R}\left(  2\left(
a+b\right)  \xi\right)  \left(  \frac{b-a}{2}\right)  \widehat{R}\left(
2\left(  b-a\right)  \xi\right)  ,
\end{align*}

where $R\ast R=\Lambda$.

Now take the inverse Fourier transform to get%
\begin{multline*}
\left(  \frac{1}{w_{R}\left(  \xi;a+\delta,b+\delta\right)  }\right)  ^{\vee
}\left(  x\right)  =\left(  \frac{1}{w_{R}\left(  \xi;a,b\right)  }\right)
^{\vee}\left(  x\right)  +\tfrac{2??c_{1}}{\sqrt{2\pi}}\left(  b-a\right)
\left(  \operatorname*{sgn}\delta\right)  \Lambda\left(  \frac{x}{\delta
}\right)  +\\
+\tfrac{8??c_{2}}{\sqrt{2\pi}}\left(  \operatorname*{sgn}\delta\right)
R\left(  \frac{\cdot}{2\delta}\right)  \ast R\left(  \frac{\cdot}{2\left(
a+b\right)  }\right)  \ast R\left(  \frac{\cdot}{2\left(  b-a\right)
}\right)  .
\end{multline*}

\end{remark}

\begin{corollary}
\label{Cor_cdeffwt_rect_phi}The weight function $w_{R}$ of Theorem
\ref{Thm_cdiffwt_rect_phi} is a $C_{BP}^{\infty}\left(  \mathbb{R}^{1}\right)
$ function.
\end{corollary}

\begin{proof}
We continue $w_{R}$ to $\mathbb{C}^{1}$ as the function of a complex
variable:
\[
w_{R}\left(  z\right)  =\frac{\sqrt{2\pi}}{4}\frac{1}{\frac{\left(
b-a\right)  z-\left(  \sin bz-\sin az\right)  }{z^{3}}},\quad z\in
\mathbb{C}^{1},
\]

and note that%
\[
\frac{\left(  b-a\right)  z-\left(  \sin bz-\sin az\right)  }{z^{3}}%
=\frac{\left(  b^{3}-a^{3}\right)  }{3!}-\frac{\left(  b^{5}-a^{5}\right)
z^{2}}{5!}+\frac{\left(  b^{7}-a^{7}\right)  z^{4}}{7!}-\ldots,
\]

is analytic and non-zero in a neighborhood of the origin. The theory of
\textit{division of a power series} implies that $w_{R}\left(  z\right)  $ is
also analytic in a neighborhood of the origin. Thus $w_{R}$ is $C^{\infty}$ in
a neighborhood of the origin.

To prove that $w_{R}$ is $C^{\infty}$ away from the origin we write%
\[
w_{R}\left(  \xi\right)  =\frac{\sqrt{2\pi}}{8}\frac{1}{\left(  \frac{b-a}%
{2}\right)  ^{3}}f_{R}\left(  \frac{b-a}{2}\xi\right)  ,
\]

where%
\[
f_{R}\left(  \xi\right)  =\frac{\xi^{2}}{1-\cos\left(  \frac{b+a}{b-a}%
\xi\right)  \frac{\sin\xi}{\xi}}.
\]

It is clear that $\cos\left(  \frac{b+a}{b-a}\xi\right)  \frac{\sin\xi}{\xi
}=1$ iff $\xi=0$, so that $f_{R}\in C^{\infty}$ away from the origin and hence
that $w\in C^{\infty}$ away from the origin. Finally, since we can write
$\left\vert f_{R}\left(  \xi\right)  \right\vert \leq\xi^{2}$ we have $w\in
C_{BP}^{\infty}$.
\end{proof}

The previous theorem will now be applied to derive an upper bound for the
central difference weight functions with parameters $n=l=1$.

\begin{corollary}
For the case $n=l=1$ the univariate central difference weight function
$w\left(  \xi\right)  =\frac{\xi^{2}}{\Delta_{2}\widehat{q}\left(  \xi\right)
}$ of Theorem \ref{Thm_cdiffwt_2} satisfies an inequality of the form%
\[
w\left(  \xi\right)  \leq cw_{R}\left(  \xi\right)  ,\quad\xi\in\mathbb{R}%
^{1}.
\]

where $w_{R}$ was defined in Theorem \ref{Thm_cdiffwt_rect_phi}.
\end{corollary}

\begin{proof}
By assumption $q\in L^{1}\left(  \mathbb{R}^{1}\right)  $, $q\left(
\xi\right)  \geq0$ for all $\xi\in\mathbb{R}^{1}$. Therefore, since an
integrable function can be defined using simple functions i.e. the
characteristic functions of intervals, there exist constants $b>a$ and $c>0$
such that $q\left(  \xi\right)  \geq cq_{R}\left(  \xi\right)  $ for $\xi
\in\left[  a,b\right]  $ where $q_{R}$ be the rectangular function defined by
\ref{a930}. Then from \ref{a943} and Theorem \ref{Thm_cdiffwt_rect_phi}%
\[
\frac{1}{w\left(  \xi\right)  }=\tfrac{4}{\sqrt{2\pi}}\frac{1}{\xi^{2}}%
\int\limits_{\mathbb{R}^{1}}\sin^{2}\frac{\xi t}{2}q\left(  t\right)
dt\geq\tfrac{4}{\sqrt{2\pi}}\frac{c}{\xi^{2}}\int_{a}^{b}\sin^{2}\frac{\xi
t}{2}q_{R}\left(  t\right)  dt=\frac{c}{w_{R}\left(  \xi\right)  }.
\]

\end{proof}

For an arbitrary central difference weight function inequality \ref{a950} gave
an order $2n$ lower bound on the growth of a weight function at infinity. We
now derive an upper bound at infinity which is also of order $2n $.

\begin{theorem}
\label{Thm_cdiffwt_bnd_on_wt_fn}There exist constants $0\leq a<b$ such that
the univariate central difference weight function $w\left(  \xi\right)
=\frac{\xi^{2n}}{\Delta_{2l}\widehat{q}\left(  \xi\right)  }$ of Theorem
\ref{Thm_cdiffwt_2} satisfies inequalities of the form%
\begin{equation}
w\left(  \xi\right)  \leq\left\{
\begin{array}
[c]{ll}%
\frac{c_{2}}{\xi^{2\left(  l-n\right)  }},\quad & \left\vert \xi\right\vert
\leq\frac{\pi}{b-a},\\
2c_{1}\xi^{2n},\quad & \left\vert \xi\right\vert >\frac{\pi}{b-a}.
\end{array}
\right. \label{a958}%
\end{equation}

where $c_{1}$ and $c_{2}$ are given by equations \ref{a935} and \ref{a936} respectively.
\end{theorem}

\begin{proof}
The proof will involve calculating a lower bound for $1/w$. Choose $0\leq a<b
$ and $c_{0}$ such that $w\left(  \xi\right)  \geq c_{0}$ when $\xi\in\left[
a,b\right]  $. Set $\delta=b-a$. We start with equation \ref{a943} so that%
\[
\frac{1}{\xi^{2\left(  l-n\right)  }w\left(  \xi\right)  }=\tfrac{2^{2l}%
}{\sqrt{2\pi}}\frac{1}{\xi^{2l}}\int_{\mathbb{R}^{1}}\sin^{2l}\frac{\xi t}%
{2}q\left(  t\right)  dt\geq\tfrac{2^{2l}c_{0}}{\sqrt{2\pi}}\frac{1}{\xi^{2l}%
}\int_{a}^{b}\sin^{2l}\frac{\xi t}{2}dt.
\]

Since $w$ is an even function we can assume $\xi>0$. Applying the change of
variables $s=\xi t/2$ gives%
\begin{equation}
\frac{1}{\xi^{2\left(  l-n\right)  }w\left(  \xi\right)  }\geq\tfrac
{2^{2l+1}c_{0}}{\sqrt{2\pi}}\frac{1}{\xi^{2l+1}}\int_{a\xi/2}^{b\xi/2}%
\sin^{2l}sds.\label{a947}%
\end{equation}

Now $b\xi/2-a\xi/2=\delta\xi/2=k\pi+\sigma\pi$ for some $0\leq\sigma<1$, and
\begin{align*}
\int_{a\xi/2}^{b\xi/2}\sin^{2l}sd  & =\int_{a\xi/2}^{a\xi/2+k\pi}\sin
^{2l}sds+\int_{a\xi/2+k\pi}^{a\xi/2+k\pi+\sigma\pi}\sin^{2l}sds\\
& =\int_{0}^{k\pi}\sin^{2l}sds+\int_{a\xi/2}^{a\xi/2+\sigma\pi}\sin^{2l}sds\\
& =2k\int_{0}^{\pi/2}\sin^{2l}sds+\int_{a\xi/2}^{a\xi/2+\sigma\pi}\sin
^{2l}sds.
\end{align*}

Since $\delta\xi/2=k\pi+\sigma\pi$ implies $\delta\xi/2<\left(  k+1\right)
\pi$ and $\frac{\delta\xi}{2\pi}-1<k$,%
\begin{equation}
\int_{a\xi/2}^{b\xi/2}\sin^{2l}sd\geq2k\int_{0}^{\pi/2}\sin^{2l}sds\geq\left(
\frac{\delta\xi}{\pi}-2\right)  \int_{0}^{\pi/2}\sin^{2l}sds.\label{a948}%
\end{equation}

We now split the range of $\xi$ into the intervals: $0<\xi\leq\frac{\pi
}{\delta}$, $\frac{\pi}{\delta}\leq\xi<\frac{4\pi}{\delta}$ and $\xi\geq
\frac{4\pi}{\delta}$ and consider each as a separate case.\medskip

\fbox{Case 1: $\xi\geq\frac{4\pi}{\delta}$} Here $k\geq2$ and $\frac{2}{\xi
}\leq\frac{\delta}{2\pi}$ so that from \ref{a947} and \ref{a948}
\begin{align*}
\frac{1}{\xi^{2\left(  l-n\right)  }w\left(  \xi\right)  }\geq\tfrac
{2^{2l+1}c_{0}}{\sqrt{2\pi}}\frac{1}{\xi^{2l+1}}\int_{a\xi/2}^{b\xi/2}%
\sin^{2l}sds &  \geq\tfrac{2^{2l+1}c_{0}}{\sqrt{2\pi}}\frac{1}{\xi^{2l+1}%
}\left(  \frac{\delta\xi}{\pi}-2\right)  \int_{0}^{\pi/2}\sin^{2l}sds\\
&  =\tfrac{2^{2l+1}c_{0}}{\sqrt{2\pi}}\frac{1}{\xi^{2l}}\left(  \frac{\delta
}{\pi}-\frac{2}{\xi}\right)  \int_{0}^{\pi/2}\sin^{2l}sds\\
&  \geq\tfrac{2^{2l+1}c_{0}}{\sqrt{2\pi}}\frac{1}{\xi^{2l}}\left(
\frac{\delta}{\pi}-\frac{\delta}{2\pi}\right)  \int_{0}^{\pi/2}\sin^{2l}sds\\
&  =\tfrac{2^{2l+1}c_{0}}{\sqrt{2\pi}}\frac{1}{\xi^{2l}}\frac{\delta}{2\pi
}\int_{0}^{\pi/2}\sin^{2l}sds\\
&  =\tfrac{2^{2l}}{\sqrt{2\pi}}\frac{1}{\pi}\left(  \int_{0}^{\pi/2}\sin
^{2l}sds\right)  \frac{\delta c_{0}}{\xi^{2l}}\\
&  =\tfrac{2^{2l}}{\sqrt{2\pi}}\frac{1}{2}\frac{\left(  2l\right)  !}%
{2^{2l}\left(  l!\right)  ^{2}}\frac{\delta c_{0}}{\xi^{2l}}\\
&  =\tfrac{1}{2\sqrt{2\pi}}\frac{\left(  2l\right)  !}{\left(  l!\right)
^{2}}\frac{\delta c_{0}}{\xi^{2l}},
\end{align*}

and thus%
\begin{equation}
w\left(  \xi\right)  \leq c_{1}\xi^{2n},\quad\xi\geq\frac{4\pi}{\delta
},\label{a955}%
\end{equation}

where%
\begin{equation}
c_{1}=2\sqrt{2\pi}\frac{\left(  l!\right)  ^{2}}{\left(  2l\right)  !}\frac
{1}{\left(  b-a\right)  c_{0}}.\label{a935}%
\end{equation}
\medskip

\fbox{Case 2: $0<\xi<\frac{\pi}{\delta}$} Here $k=0$ and $b\xi/2-a\xi
/2=\delta\xi/2<\pi/2$ so that $\sin s\geq\frac{2}{\pi}s$ and
\begin{align*}
\int_{a\xi/2}^{b\xi/2}\sin^{2l}sds\geq\int_{a\xi/2}^{b\xi/2}\left(  \frac
{2}{\pi}s\right)  ^{2l}ds &  \geq\left(  \frac{2}{\pi}\right)  ^{2l}\left[
\frac{s^{2l+1}}{2l+1}\right]  _{a\xi/2}^{b\xi/2}\\
&  =\frac{1}{2}\left(  \frac{1}{\pi}\right)  ^{2l}\frac{1}{2l+1}\left(
b^{2l+1}-a^{2l+1}\right)  \xi^{2l+1}.
\end{align*}

Thus by \ref{a947}
\begin{align*}
\frac{1}{\xi^{2\left(  l-n\right)  }w\left(  \xi\right)  }\geq\tfrac
{2^{2l+1}c_{0}}{\sqrt{2\pi}}\frac{1}{\xi^{2l+1}}\int_{a\xi/2}^{b\xi/2}%
\sin^{2l}sds &  \geq\tfrac{2^{2l+1}c_{0}}{\sqrt{2\pi}}\frac{1}{2}\left(
\frac{1}{\pi}\right)  ^{2l}\frac{1}{2l+1}\left(  b^{2l+1}-a^{2l+1}\right) \\
&  \geq\frac{1}{\sqrt{2\pi}}\left(  \frac{2}{\pi}\right)  ^{2l}\frac
{b^{2l+1}-a^{2l+1}}{2l+1}c_{0},
\end{align*}

and%
\begin{equation}
w\left(  \xi\right)  \leq\frac{c_{2}}{\xi^{2\left(  l-n\right)  }},\quad
0<\xi<\frac{\pi}{\delta},\label{a956}%
\end{equation}

where%
\begin{equation}
c_{2}=\frac{1}{\sqrt{2\pi}}\left(  \frac{2}{\pi}\right)  ^{2l}\frac{\left(
b^{2l+1}-a^{2l+1}\right)  c_{0}}{2l+1}.\label{a936}%
\end{equation}
\medskip

\fbox{Case 3: $\frac{\pi}{\delta}\leq\xi<\frac{4\pi}{\delta}$} Here%
\[
\int_{a\xi/2}^{b\xi/2}\sin^{2l}sds\geq\int_{0}^{\pi/2}\sin^{2l}sds,
\]

so that by \ref{a947}
\begin{align*}
\frac{1}{\xi^{2\left(  l-n\right)  }w\left(  \xi\right)  }\geq\tfrac
{2^{2l+1}c_{0}}{\sqrt{2\pi}}\frac{1}{\xi^{2l+1}}\int_{a\xi/2}^{b\xi/2}%
\sin^{2l}sds &  \geq\tfrac{2^{2l+1}c_{0}}{\sqrt{2\pi}}\frac{1}{\xi^{2l+1}}%
\int_{0}^{\pi/2}\sin^{2l}sds\\
&  =\tfrac{2^{2l+1}c_{0}}{\sqrt{2\pi}}\frac{1}{\xi^{2l+1}}\frac{\pi}{2}%
\frac{\left(  2l\right)  !}{2^{2l}\left(  l!\right)  ^{2}}\\
&  =\frac{\sqrt{2\pi}}{2}\frac{\left(  2l\right)  !}{\left(  l!\right)  ^{2}%
}\frac{c_{0}}{\xi^{2l+1}}.
\end{align*}

Thus%
\[
w\left(  \xi\right)  \leq\frac{2}{\sqrt{2\pi}}\frac{\left(  l!\right)  ^{2}%
}{\left(  2l\right)  !}\frac{\xi^{2n+1}}{c_{0}}\leq\frac{2}{\sqrt{2\pi}}%
\frac{\left(  l!\right)  ^{2}}{\left(  2l\right)  !}\frac{\xi^{2n}}{c_{0}%
}\frac{4\pi}{\delta}=4\sqrt{2\pi}\frac{\left(  l!\right)  ^{2}}{\left(
2l\right)  !}\frac{1}{c_{0}\delta}\xi^{2n}=2c_{1}\xi^{2n},
\]

and
\begin{equation}
w\left(  \xi\right)  \leq2c_{1}\xi^{2n},\quad\frac{\pi}{\delta}\leq\xi
<\frac{4\pi}{\delta}.\label{a957}%
\end{equation}

Noting that $w$ is an even function the inequalities \ref{a955}, \ref{a956}
and \ref{a957} imply \ref{a958} as desired.
\end{proof}

We have the following results for univariate central difference weight
function near zero and infinity:

\begin{corollary}
\label{Cor_cdiffwt_bnd_on_wt_fn}Suppose $w$ is a central difference weight
function on $\mathbb{R}^{1}$ with parameters $n,l,q\left(  \cdot\right)  $
which satisfies the condition \ref{a959}. Then for any $r>0$ there exist
constants $c_{r},c_{r}^{\prime},k_{r},k_{r}^{\prime}>0$ such that%
\begin{equation}
k_{r}\xi^{2n}\leq w\left(  \xi\right)  \leq k_{r}^{\prime}\xi^{2n}%
,\quad\left\vert \xi\right\vert \geq r,\label{a971}%
\end{equation}

and%
\begin{equation}
\frac{c_{r}}{\xi^{2\left(  l-n\right)  }}\leq w\left(  \xi\right)  \leq
\frac{c_{r}^{\prime}}{\xi^{2\left(  l-n\right)  }},\quad\left\vert
\xi\right\vert \leq r.\label{a972}%
\end{equation}

Also%
\[
X_{w}^{0}\left(  \mathbb{R}^{1}\right)  \hookrightarrow W^{n}\left(
\mathbb{R}^{1}\right)  .
\]

\end{corollary}

\begin{proof}
Inequalities \ref{a971} follow from inequality \ref{a958} of Theorem
\ref{Thm_cdiffwt_bnd_on_wt_fn} and inequality \ref{a950} of Theorem
\ref{Thm_cdiffwt_4}.

Inequalities \ref{a972} are the inequalities \ref{a964} of Theorem
\ref{Thm_cdiffwt_3}.

When $r=1$, \ref{a971} implies $k_{r}\left(  1+\xi^{2n}\right)  \leq w\left(
\xi\right)  $ when $\left\vert \xi\right\vert \geq1$, and \ref{a972} implies
$c_{r}\left(  1+\xi^{2n}\right)  \leq w\left(  \xi\right)  $ when $\left\vert
\xi\right\vert \leq1$. Thus $k_{r}\left(  1+\xi^{2n}\right)  \leq w\left(
\xi\right)  $ for all $\xi$ so that%
\begin{align*}
X_{w}^{0}  & =\left\{  u\in S^{\prime}:\widehat{u}\in L_{loc}^{1},\text{ }\int
w\left\vert \widehat{u}\right\vert ^{2}<\infty\right\} \\
& \subset\left\{  u\in S^{\prime}:\widehat{u}\in L_{loc}^{1},\text{ }%
\int\left(  1+\xi^{2n}\right)  \left\vert \widehat{u}\right\vert ^{2}%
<\infty\right\} \\
& =\left\{  u\in S^{\prime}:\widehat{u}\in L_{loc}^{1},\text{ }\int\left(
1+\xi^{2}\right)  ^{n}\left\vert \widehat{u}\right\vert ^{2}<\infty\right\} \\
& =W^{n}.
\end{align*}

\end{proof}

Later we will be interested in the sequence of central difference weight
functions $\left\{  w_{k}\right\}  $ which corresponds to $b=k^{-1}>0$ and
$a=-k^{-1}$ where $k$ is a positive integer. These weight functions have the
following properties:

\begin{corollary}
\label{Cor_cntl_diff_bnd_on_wt_fn_2}Suppose $k\geq1$ is an integer and $R$ is
the unit rectangular function with support on $\left[  -1,1\right]
\subset\mathbb{R}^{1}$.

Then $q_{k}\left(  \xi\right)  =\frac{k}{2}R\left(  k\xi\right)  $ satisfies
the conditions of Theorem \ref{Thm_cdiffwt_2} for $n=l=1$ and the
corresponding central difference weight function is%
\[
w_{k}\left(  \xi\right)  =\frac{\sqrt{2\pi}}{16}k^{3}f_{R}\left(  \frac{\xi
}{k}\right)  ,
\]

where%
\[
f_{R}\left(  \xi\right)  =\frac{\xi^{2}}{1-\frac{\sin\xi}{\xi}}.
\]

\end{corollary}

\begin{proof}
The function $\frac{k}{2}R\left(  k\xi\right)  $ is a rectangular function
equal to $\frac{k}{2}$ on $\left[  -\frac{1}{k},\frac{1}{k}\right]  $. Thus
$a=-\frac{1}{k}$ and $b=\frac{1}{k}$ so that by Theorem
\ref{Thm_cdiffwt_rect_phi}
\[
w_{k}\left(  \xi\right)  =\frac{\sqrt{2\pi}}{16}\frac{\xi^{2}}{\frac{1}%
{k}-\frac{\sin\left(  \xi/k\right)  }{\xi}}=\frac{\sqrt{2\pi}}{16}\frac
{\xi^{3}}{\frac{\xi}{k}-\sin\frac{\xi}{k}}=\frac{\sqrt{2\pi}}{16}k^{3}%
f_{R}\left(  \xi/k\right)  .
\]

\end{proof}

\section{The central difference basis
functions\label{Sect_tenprod_centdiff_basis_fns}}

In Subsection \ref{SbSect_CentDiffBasis_MultConvolFormula} we derive
multiplicative convolution formulas for the univariate central difference
basis function and its derivatives.

Subsection \ref{SbSect_CentDiffBasis_Lipschitz}: We derive Lipschitz
continuity estimates based on the multiplicative convolution formulas.

In Subsection \ref{SbSect_centdiff_basis_part_mom} we give formulas for the
univariate central difference basis function and its derivatives that does not
involve calculating the extended B-spline basis functions but instead involve
calculating partial moments of $q$. These results also give information about
the support of the basis function and an improved smoothness estimate which
does not involve assuming that $q$ is bounded.

In Subsection \ref{SbSect_tempdistrib_1_dim_centdiff_basis} we will use the
tempered distribution Taylor series expansion introduced in Section
\ref{Sect_Taylor_series_data_fn} and the theory of the spaces $S_{\emptyset
,k}\subset S$ introduced in Definition \ref{Def_So,n} to confirm the partial
moment formulas for the central difference basis function. In the process we
obtain another multiplicative convolution formula which is given in Theorem
\ref{Thm_cdiffbasis_mult_convol_ql}.

Subsection \ref{SbSect_cdiffbasis_Lip_part_mom}: We derive Lipschitz
continuity estimates based on partial moment formulas.

\subsection{A multiplicative convolution formula for the basis
function\label{SbSect_CentDiffBasis_MultConvolFormula}}

We will now derive a multiplicative convolution formula for the basis
functions of zero order generated by the central difference weight functions
of Definition \ref{Def_central_diff_wt_func}. In Theorem
\ref{Thm_basis_tensor_hat_W3} it was shown that a univariate extended B-spline
basis function $G_{s}$ is such that $G_{s}\in C_{0}^{\left(  2n-2\right)
}\left(  \mathbb{R}^{1}\right)  $ and $D^{2n-1}G_{s}$ is a piecewise constant
function. In the next two theorems we will prove analogues of these results
for the central difference basis functions. In fact we show that a central
difference basis function $G_{c}$ is an even function such that $G_{c}\in
C_{B}^{\left(  2n-2\right)  }\left(  \mathbb{R}^{1}\right)  $, $D^{2n-1}%
G_{c}\in C^{\left(  0\right)  }\left(  \mathbb{R}^{1}\setminus0\right)  \cap
L^{\infty}$.

\begin{theorem}
\label{Thm_G_basis_def_2}(\textbf{Multiplicative convolution formula}) Suppose
$w$ is the tensor product central difference weight function with property W02
or W03 defined in Corollary \ref{Cor_cdiffwt_W02_W03}. Set%
\begin{equation}
G_{s}\left(  t\right)  =\left(  -1\right)  ^{l-n}\tfrac{\left(  2\pi\right)
^{l/2}}{2^{2\left(  l-n\right)  +1}}D^{2\left(  l-n\right)  }\left(  \left(
\ast\Lambda\right)  ^{l}\right)  \left(  \tfrac{t}{2}\right)  ,\label{a926}%
\end{equation}

where $G_{s}$ is the extended B-spline basis function defined by \ref{1.49}.

Then this weight function generates the tensor product basis function
$G_{c}\left(  x\right)  =\prod\limits_{i=1}^{d}G_{1}\left(  x_{i}\right)  $
where
\begin{equation}
G_{1}\left(  s\right)  =\tfrac{2^{2\left(  l-n\right)  +1}}{\sqrt{2\pi}}%
\int_{\mathbb{R}^{1}}G_{s}\left(  \frac{2s}{t}\right)  \left\vert t\right\vert
^{2n-1}q\left(  t\right)  dt,\text{\quad}s\in\mathbb{R}^{1},\label{a945}%
\end{equation}

$G_{1}\in C_{B}^{\left(  2n-2\right)  }$ and $G_{1}$ is an \textbf{even
function}. Further, for $1\leq k\leq2n-2$,%
\begin{equation}
D^{k}G_{1}\left(  s\right)  =\tfrac{2^{2\left(  l-n\right)  +k+1}}{\sqrt{2\pi
}}\int\left(  D^{k}G_{s}\right)  \left(  \frac{2s}{\left\vert t\right\vert
}\right)  \left\vert t\right\vert ^{2n-k-1}q\left(  t\right)  dt,\text{\quad
}s\in\mathbb{R}^{1},\label{a940}%
\end{equation}

Finally, for $0\leq\tau\leq k$,%
\[
\left\vert D^{k}G_{1}\left(  s\right)  \right\vert \leq\frac{\left\Vert
\left\vert \cdot\right\vert ^{\tau}D^{k}G_{s}\right\Vert _{\infty}\left\Vert
\left\vert \cdot\right\vert ^{2n-k-1+\tau}q\right\Vert _{1}}{\left\vert
s\right\vert ^{\tau}}.,\text{\quad}s\in\mathbb{R}^{1}\setminus0.
\]

\end{theorem}

\begin{proof}
We first note that from Subsection \ref{SbSect_basis_fns} the multivariate
basis function will just be the product of the univariate basis functions.

Because of the multiplicity of variables involved we will use the notation
$F[\,]$ for the Fourier transform. We want to express $1/w$ as a Fourier
transform. Continuing on from equation \ref{a943},
\begin{align*}
\frac{1}{w\left(  \xi\right)  }  & =\tfrac{1}{\sqrt{2\pi}}\frac{1}{\xi^{2n}%
}\int\limits_{\mathbb{R}^{1}}2^{2l}\sin^{2l}\frac{\xi t}{2}q\left(  t\right)
dt\\
& =\tfrac{1}{\sqrt{2\pi}}\frac{1}{\xi^{2n}}\int\limits_{\mathbb{R}^{1}}%
2^{2l}\left(  \frac{\xi t}{2}\right)  ^{2l}\frac{\sin^{2l}\frac{\xi t}{2}%
}{\left(  \frac{\xi t}{2}\right)  ^{2l}}q\left(  t\right)  dt\\
& =\tfrac{1}{\sqrt{2\pi}}\frac{1}{\xi^{2n}}\int\limits_{\mathbb{R}^{1}}\left(
\xi t\right)  ^{2l}\frac{\sin^{2l}\frac{\xi t}{2}}{\left(  \frac{\xi t}%
{2}\right)  ^{2l}}q\left(  t\right)  dt\\
& =\tfrac{1}{\sqrt{2\pi}}\int\left(  \xi t\right)  ^{2\left(  l-n\right)
}\frac{\sin^{2l}\frac{\xi t}{2}}{\left(  \frac{\xi t}{2}\right)  ^{2l}}%
t^{2n}q\left(  t\right)  dt\\
& =\tfrac{1}{\sqrt{2\pi}}\int\left(  \xi t\right)  ^{2\left(  l-n\right)
}\left(  \frac{\sin^{2}\frac{\xi t}{2}}{\left(  \frac{\xi t}{2}\right)  ^{2}%
}\right)  ^{l}t^{2n}q\left(  t\right)  dt\\
& =\tfrac{1}{\sqrt{2\pi}}\int\left(  \xi t\right)  ^{2\left(  l-n\right)
}\left(  \sqrt{2\pi}\widehat{\Lambda}\left(  \xi t\right)  \right)  ^{l}%
t^{2n}q\left(  t\right)  dt\\
& =\tfrac{1}{\sqrt{2\pi}}\int\left(  2\pi\right)  ^{l/2}\left(  \xi t\right)
^{2\left(  l-n\right)  }F\left[  \left(  \ast\Lambda\right)  ^{l}\right]
\left(  \xi t\right)  t^{2n}q\left(  t\right)  dt\\
& =\tfrac{1}{\sqrt{2\pi}}\int F\left[  \left(  -1\right)  ^{\left(
l-n\right)  }\left(  2\pi\right)  ^{l/2}D^{2\left(  l-n\right)  }\left(
\ast\Lambda\right)  ^{l}\right]  \left(  \xi t\right)  t^{2n}q\left(
t\right)  dt\\
& =\tfrac{2^{2\left(  l-n\right)  +1}}{\sqrt{2\pi}}\int F\left[  G_{s}\left(
2x\right)  \right]  \left(  \xi t\right)  t^{2n}q\left(  t\right)  dt\\
& =\tfrac{2^{2\left(  l-n\right)  +1}}{\sqrt{2\pi}}\int F_{x}\left[
G_{s}\left(  \frac{2x}{t}\right)  \right]  \left(  \xi\right)  \left\vert
t\right\vert ^{2n-1}q\left(  t\right)  dt.
\end{align*}

Now%
\begin{align}
G_{1}\left(  s\right)   & =\tfrac{1}{\sqrt{2\pi}}\int\frac{e^{i\xi s}%
}{w\left(  \xi\right)  }d\xi\nonumber\\
& =\tfrac{2^{2\left(  l-n\right)  +1}}{2\pi}\int\int e^{i\xi s}F_{s}\left[
G_{s}\left(  \frac{2s}{t}\right)  \right]  \left(  \xi\right)  \left\vert
t\right\vert ^{2n-1}q\left(  t\right)  dtd\xi,\label{a951}%
\end{align}

and we want to use Fubini's theorem to change the order of integration. Indeed%
\begin{align*}
\int\int &  \left\vert e^{i\xi s}F_{s}\left[  G_{s}\left(  \frac{2s}%
{t}\right)  \right]  \left(  \xi\right)  \left\vert t\right\vert
^{2n-1}q\left(  t\right)  \right\vert d\xi dt\\
&  \leq\int\int\left\vert F_{s}\left[  G_{s}\left(  \frac{2s}{t}\right)
\right]  \left(  \xi\right)  \right\vert d\xi\text{ }\left\vert t\right\vert
^{2n-1}q\left(  t\right)  dt\\
&  =\int\int\left\vert F\left[  G_{s}\right]  \left(  \frac{1}{2}\xi t\right)
\right\vert \frac{1}{2}\left\vert t\right\vert d\xi\text{ }\left\vert
t\right\vert ^{2n-1}q\left(  t\right)  dt\\
&  =\int\int\frac{\sin^{2l}\left(  \xi t/2\right)  }{\left(  \xi t/2\right)
^{2n}}\frac{1}{2}\left\vert t\right\vert d\xi\text{ }\left\vert t\right\vert
^{2n-1}q\left(  t\right)  dt\\
&  =\int\int\frac{\sin^{2l}\left(  \xi\left\vert t\right\vert /2\right)
}{\left(  \xi\left\vert t\right\vert /2\right)  ^{2n}}\frac{1}{2}\left\vert
t\right\vert d\xi\text{ }\left\vert t\right\vert ^{2n-1}q\left(  t\right)
dt\\
&  :\eta=\xi\left\vert t\right\vert /2,\text{ }d\eta=\frac{1}{2}\left\vert
t\right\vert d\xi\Rightarrow\\
&  =\int\int\frac{\sin^{2l}\eta}{\eta^{2n}}d\eta\text{ }\left\vert
t\right\vert ^{2n-1}q\left(  t\right)  dt\\
&  =\int\int F\left[  G_{s}\right]  \left(  \eta\right)  d\eta\text{
}\left\vert t\right\vert ^{2n-1}q\left(  t\right)  dt\\
&  =\sqrt{2\pi}G_{s}\left(  0\right)  \int\left\vert t\right\vert
^{2n-1}q\left(  t\right)  dt\\
&  <\infty,
\end{align*}

by the assumptions of Corollary \ref{Cor_cdiffwt_W02_W03}. Hence the order of
integration can be reversed to give%
\begin{align*}
G_{1}\left(  s\right)   & =\tfrac{2^{2\left(  l-n\right)  +1}}{2\pi}\int\int
e^{i\xi s}F_{s}\left[  G_{s}\left(  \frac{2s}{t}\right)  \right]  \left(
\xi\right)  d\xi\text{ }\left\vert t\right\vert ^{2n-1}q\left(  t\right)  dt\\
& =\tfrac{2^{2\left(  l-n\right)  +1}}{\sqrt{2\pi}}\int G_{s}\left(  \frac
{2s}{t}\right)  \left\vert t\right\vert ^{2n-1}q\left(  t\right)  dt,
\end{align*}

as required.

By Theorem \ref{Thm_cdiffwt_2} $w$ is even. Hence $G_{1}$ is even since $w$ is
real-valued. Also by Theorem \ref{Thm_cdiffwt_2}, $\kappa+1/2<n$ so that
$\left\lfloor 2\kappa\right\rfloor \leq2n-2$. But by Theorem
\ref{Thm_basis_fn_properties_all_m_W2} $G_{1}\in C_{B}^{\left(  \left\lfloor
2\kappa\right\rfloor \right)  }$ and so $G_{1}\in C_{B}^{\left(  2n-2\right)
} $.

To prove \ref{a940} we start with the equations \ref{a951} for $G_{1}$. If
$k\leq2n-2$
\begin{align*}
D^{k}G_{1}\left(  s\right)   &  =\tfrac{1}{\sqrt{2\pi}}\int\frac{\left(
i\xi\right)  ^{k}e^{i\xi s}}{w\left(  \xi\right)  }d\xi\\
&  =\tfrac{1}{\sqrt{2\pi}}\int e^{i\xi s}\left(  i\xi\right)  ^{k}\left(
\tfrac{2^{2\left(  l-n\right)  +1}}{\sqrt{2\pi}}\int F_{x}\left[  G_{s}\left(
\frac{2x}{t}\right)  \right]  \left(  \xi\right)  \left\vert t\right\vert
^{2n-1}q\left(  t\right)  dt\right)  d\xi\\
&  =\tfrac{2^{2\left(  l-n\right)  +1}}{2\pi}\int\int e^{i\xi s}\left(
i\xi\right)  ^{k}F_{x}\left[  G_{s}\left(  \frac{2x}{t}\right)  \right]
\left(  \xi\right)  \left\vert t\right\vert ^{2n-1}q\left(  t\right)  dtd\xi\\
&  =\tfrac{2^{2\left(  l-n\right)  +1}}{2\pi}\int\int e^{i\xi s}\left(
i\xi\right)  ^{k}F_{x}\left[  G_{s}\left(  \frac{2x}{\left\vert t\right\vert
}\right)  \right]  \left(  \xi\right)  \left\vert t\right\vert ^{2n-1}q\left(
t\right)  dtd\xi\\
&  =\tfrac{2^{2\left(  l-n\right)  +k+1}}{2\pi}\int\int e^{i\xi s}F_{x}\left[
\left(  D^{k}G_{s}\right)  \left(  \frac{2x}{\left\vert t\right\vert }\right)
\right]  \left(  \xi\right)  \left\vert t\right\vert ^{2n-k-1}q\left(
t\right)  dtd\xi\\
&  =\tfrac{2^{2\left(  l-n\right)  +k+1}}{2\pi}\int\int e^{i\xi s}F_{x}\left[
\left(  D^{k}G_{s}\right)  \left(  \frac{2x}{\left\vert t\right\vert }\right)
\right]  \left(  \xi\right)  \left\vert t\right\vert ^{2n-k-1}q\left(
t\right)  dtd\xi.
\end{align*}

An argument similar to that used to change the order of integration in the
case $k=0$ allows us to write%
\begin{align*}
D^{k}G_{1}\left(  s\right)   & =\tfrac{2^{2\left(  l-n\right)  +k+1}}{2\pi
}\int\int e^{i\xi s}F_{x}\left[  \left(  D^{k}G_{s}\right)  \left(  \frac
{2x}{\left\vert t\right\vert }\right)  \right]  \left(  \xi\right)  d\xi\text{
}\left\vert t\right\vert ^{2n-k-1}q\left(  t\right)  dt\\
& =\tfrac{2^{2\left(  l-n\right)  +k+1}}{\sqrt{2\pi}}\int\left(  D^{k}%
G_{s}\right)  \left(  \frac{2s}{\left\vert t\right\vert }\right)  \left\vert
t\right\vert ^{2n-k-1}q\left(  t\right)  dt,
\end{align*}

as required.

Again from \ref{a940}, since $D^{k}G_{s}$ is bounded with bounded support, for
any $0\leq\tau\leq k$,
\begin{align*}
\left\vert s\right\vert ^{\tau}\left\vert D^{k}G_{1}\left(  s\right)
\right\vert  & \leq\int\limits_{\mathbb{R}^{1}}\left\vert \left\vert
s\right\vert ^{\tau}\left(  D^{k}G_{s}\right)  \left(  \frac{2s}{t}\right)
\right\vert \left\vert t\right\vert ^{2n-k-1}q\left(  t\right)  dt\\
& =\int\left\vert \left\vert \frac{s}{t}\right\vert ^{\tau}\left(  D^{k}%
G_{s}\right)  \left(  \frac{2s}{t}\right)  \right\vert \left\vert t\right\vert
^{2n-k-1+\tau}q\left(  t\right)  dt\\
& =2^{-\tau}\int\left\vert \left\vert \frac{2s}{t}\right\vert ^{\tau}\left(
D^{k}G_{s}\right)  \left(  \frac{2s}{t}\right)  \right\vert \left\vert
t\right\vert ^{2n-k-1+\tau}q\left(  t\right)  dt\\
& \leq2^{-\tau}\left\Vert \left\vert \cdot\right\vert ^{\tau}D^{k}%
G_{s}\right\Vert _{\infty}\int\left\vert t\right\vert ^{2n-k-1+\tau}q\left(
t\right)  dt\\
& =2^{-\tau}\left\Vert \left\vert \cdot\right\vert ^{\tau}D^{k}G_{s}%
\right\Vert _{\infty}\left\Vert \left\vert \cdot\right\vert ^{2n-k-1+\tau
}q\right\Vert _{1},
\end{align*}

so that%
\[
\left\vert D^{k}G_{1}\left(  s\right)  \right\vert \leq2^{-\tau}%
\frac{\left\Vert \left\vert \cdot\right\vert ^{\tau}D^{k}G_{s}\right\Vert
_{\infty}\left\Vert \left\vert \cdot\right\vert ^{2n-k-1+\tau}q\right\Vert
_{1}}{\left\vert s\right\vert ^{\tau}}.
\]

\end{proof}

\begin{remark}
\textbf{Check} The formula \ref{a945} can be checked by replacing $q\left(
t\right)  $ by the sequence $q_{k}$ in \ref{a978}.

\textbf{Formally} $\frac{\sqrt{2\pi}}{2^{2l}}k\psi\left(  k\left(  t-2\right)
\right)  \rightarrow\frac{\sqrt{2\pi}}{2^{2l}}\delta\left(  t-2\right)  $ so
\begin{align*}
\tfrac{2^{2\left(  l-n\right)  +1}}{\sqrt{2\pi}}\int G_{s}\left(  \frac{2s}%
{t}\right)  \left\vert t\right\vert ^{2n-1}q_{k}\left(  t\right)  dt  &
\rightarrow\tfrac{2^{2\left(  l-n\right)  +1}}{\sqrt{2\pi}}\int G_{s}\left(
\frac{2s}{t}\right)  \left\vert t\right\vert ^{2n-1}\frac{\sqrt{2\pi}}{2^{2l}%
}\delta\left(  t-2\right)  dt\\
& =2^{-2n+1}\int G_{s}\left(  \frac{2s}{t}\right)  \left\vert t\right\vert
^{2n-1}\delta\left(  t-2\right)  dt\\
& =2^{-2n+1}G_{s}\left(  \frac{2s}{2}\right)  2^{2n-1}\\
& =G_{s}\left(  s\right)  .
\end{align*}

Similarly, the derivative formulas \ref{a940} yield%
\[
\tfrac{2^{2\left(  l-n\right)  +k+1}}{\sqrt{2\pi}}\int\left(  D^{k}%
G_{s}\right)  \left(  \frac{2s}{\left\vert t\right\vert }\right)  \left\vert
t\right\vert ^{2n-k-1}q\left(  t\right)  dt\rightarrow D^{k}G_{s}\left(
s\right)  .
\]

\end{remark}

\begin{remark}
\label{Rem_Thm_G_basis_def_2}\textbf{Calculating }$q_{e}$\textbf{\ from the
central difference basis function }$G_{c}$. Here we will derive the formula
\ref{2.20}.

The expressions \ref{a940} and \ref{a945} are \textbf{multiplicative}
\textbf{convolutions} and exponential substitutions would convert them into
\textbf{additive convolutions}. Then the Fourier transform could be used to
express $q$ in terms of the basis function and perhaps even enable us to
characterize the basis functions.

?? \textbf{UNFINISHED}! \textbf{FIX}! Write%
\begin{align*}
G_{c}\left(  s\right)   & =\tfrac{2^{2\left(  l-n\right)  +1}}{\sqrt{2\pi}%
}\int_{\mathbb{R}^{1}}G_{s}\left(  \frac{2s}{t}\right)  \left\vert
t\right\vert ^{2n-1}q\left(  t\right)  dt=\tfrac{2^{2\left(  l-n\right)  +2}%
}{\sqrt{2\pi}}\int_{0}^{\infty}G_{s}\left(  \frac{2s}{t}\right)  \left\vert
t\right\vert ^{2n-1}q_{e}\left(  t\right)  dt,\\
q_{e}\left(  t\right)   & :=\left(  q\left(  t\right)  +q\left(  -t\right)
\right)  /2.
\end{align*}

Noting that $G_{c}$ is an even function, set $s=e^{u}$ and $t=e^{v}$ so that%
\begin{equation}
G_{c}\left(  e^{u}\right)  =\tfrac{2^{2\left(  l-n\right)  +2}}{\sqrt{2\pi}%
}\int_{-\infty}^{\infty}G_{s}\left(  2e^{u-v}\right)  e^{\left(  2n-1\right)
v}q_{e}\left(  e^{v}\right)  e^{v}dv=\tfrac{2^{2\left(  l-n\right)  +2}}%
{\sqrt{2\pi}}\int_{-\infty}^{\infty}G_{s}\left(  2e^{u-v}\right)  e^{2nv}%
q_{e}\left(  e^{v}\right)  dv.\label{2.21}%
\end{equation}

Define%
\begin{align*}
H_{c}\left(  u\right)   & :=G_{c}\left(  e^{u}\right)  ,\\
H_{\Lambda}\left(  x\right)   & :=2^{2\left(  l-n\right)  +2}G_{s}\left(
2e^{x}\right)  ,\\
p_{e}\left(  v\right)   & :=e^{2nv}q_{e}\left(  e^{v}\right)  ,
\end{align*}

so that formally%
\[
H_{c}=H_{\Lambda}\ast p_{e}.
\]

But%
\[
\int p_{e}\left(  v\right)  dv=\int e^{2nv}q_{e}\left(  e^{v}\right)
dv=\int_{0}^{\infty}t^{2n-1}q_{e}\left(  t\right)  dt<\infty,
\]

so $p_{e}\in L^{1}$ and so $p_{e}\in S^{\prime}$ with $\widehat{p_{e}}\in
C_{B}^{\left(  0\right)  }$. Also, $G_{s}$ has bounded support so $H_{\Lambda
}\in C^{\left(  0\right)  }\cap L^{2}\cap\mathcal{E}^{\prime}$ and
$\widehat{H_{\Lambda}}\in C_{BP}^{\infty}$. This means that $H_{c}\in L^{2}$
and the formula \ref{2.21} is valid and equals the distribution formula in
$S^{\prime}\ast\mathcal{E}^{\prime}$ given in 2.9.6 of Vladimirov
\cite{Vladimirov}. Further, 2.9.5 of Vladimirov \cite{Vladimirov} now allows
us to take the tempered distribution Fourier transform to get%
\begin{align*}
&  \widehat{H_{c}}=\widehat{H_{\Lambda}}\widehat{p_{e}}\Rightarrow
\widehat{p_{e}}=\frac{\widehat{H_{c}}}{\widehat{H_{\Lambda}}}\Rightarrow
\widehat{e^{2nv}q_{e}\left(  e^{v}\right)  }=\frac{\widehat{H_{c}}%
}{\widehat{H_{\Lambda}}}\Rightarrow e^{2nv}q_{e}\left(  e^{v}\right)  =\left(
\frac{\widehat{H}_{c}}{\widehat{H_{\Lambda}}}\right)  ^{\vee}\Rightarrow\\
&  \Rightarrow q_{e}\left(  e^{v}\right)  =e^{-2nv}\left(  \frac
{\widehat{H_{c}}}{\widehat{H_{\Lambda}}}\right)  ^{\vee}\left(  v\right)  ,
\end{align*}

which implies%
\begin{align}
q_{e}\left(  t\right)  =\text{?? }t^{-2n}\left(  \frac{\widehat{H_{c}}%
}{\widehat{H_{\Lambda}}}\right)  ^{\vee}\left(  \ln t\right)   &
=t^{-2n}\left(  \frac{\sqrt{2\pi}\widehat{G_{c}\left(  e^{x}\right)  }%
}{\widehat{G_{s}\left(  e^{x}\right)  }}\right)  ^{\vee}\left(  \ln t\right)
\nonumber\\
&  =\sqrt{2\pi}t^{-2n}\left(  \frac{\widehat{G_{c}\left(  e^{x}\right)  }%
}{\widehat{G_{s}\left(  e^{x}\right)  }}\right)  ^{\vee}\left(  \ln t\right)
.\label{2.20}%
\end{align}

Note that $\frac{\widehat{H_{c}}}{\widehat{H_{\Lambda}}}\in C_{B}^{\left(
0\right)  }$ so the inverse Fourier transform is defined in the tempered
distribution sense.

Now%
\[
\widehat{G_{s}\left(  e^{x}\right)  }\left(  t\right)  =\tfrac{1}{\sqrt{2\pi}%
}\int e^{-its}G_{s}\left(  e^{s}\right)  ds,
\]

and because $G_{s}$ is a piecewise polynomial we are led to calculate
\[
\int_{a}^{b}e^{-its}e^{ks}ds=\int_{a}^{b}e^{\left(  k-it\right)  s}ds=\left[
\frac{e^{\left(  k-it\right)  s}}{k-it}\right]  _{a}^{b}=\frac{e^{\left(
k-it\right)  b}-e^{\left(  k-it\right)  a}}{k-it}.
\]

?? \textbf{NOTE} Could we use the theory of \textbf{Mellin transforms} and the
fact that the Mellin transform of a multiplicative convolution is the product
of the Mellin transforms?
\end{remark}

For the next theorem we will prove the following version of Leibniz's theorem
for differentiating under the integral sign. See the more general Lemma
\ref{Lem_diff_under_integ} in the appendix.

\begin{lemma}
\label{Lem_Leibnitz_diff_under_integ}Suppose we have a function $f:\mathbb{R}%
^{2}\rightarrow\mathbb{C}$ and a point $s_{0}\in\mathbb{R}$. Suppose $f$ has
the following properties :

\begin{enumerate}
\item In some neighborhood of $\mathcal{N}\left(  s_{0}\right)  $ of $s_{0}$
and for almost all $t$, $f\left(  s,t\right)  \in C^{\left(  0\right)  }$ and
$D_{1}f\left(  s,t\right)  $ is piecewise continuous.

\item For each $s\in\mathcal{N}\left(  s_{0}\right)  $, $f\left(
s,\cdot\right)  \in L^{1}$ and $D_{1}f\left(  s,\cdot\right)  \in L^{1}$.

\item There exist constants $c,\varepsilon>0$ such that%
\begin{equation}
\int\left\vert D_{1}f\left(  s,t\right)  -D_{1}f\left(  s_{0},t\right)
\right\vert dt\leq c\left\vert s-s_{0}\right\vert ,\quad\left\vert
s-s_{0}\right\vert \leq\varepsilon.\label{a991}%
\end{equation}

Then we can conclude that in some neighborhood of $s_{0}$, $D_{1}\int f\left(
s,t\right)  dt=\int\left(  D_{1}f\right)  \left(  s,t\right)  dt$ and this is
a continuous function.
\end{enumerate}
\end{lemma}

\begin{proof}
The conditions of parts 1 and 2 imply that for each fixed $t$, $f\left(
s,t\right)  $ is absolutely continuous as a function of $s$ in a neighborhood
of $s_{0}$. This absolute continuity enables us to write
\[
\int\left(  \frac{f\left(  s,t\right)  -f\left(  s_{0},t\right)  }{h}-\left(
D_{1}f\right)  \left(  s_{0},t\right)  \right)  dt=\frac{1}{h}\int\int_{s_{0}%
}^{s}\left(  \left(  D_{1}f\right)  \left(  u,t\right)  -\left(
D_{1}f\right)  \left(  s_{0},t\right)  \right)  du\,dt,
\]

where $h=s-s_{0}$. Condition 3 of this lemma allows us to use Fubini's theorem
to change the order of integration in the following calculation. Assuming that
$s>s_{0}$%
\begin{align*}
\left\vert \int\left(  \frac{f\left(  s,t\right)  -f\left(  s_{0},t\right)
}{h}-\left(  D_{1}f\right)  \left(  s_{0},t\right)  \right)  dt\right\vert  &
\leq\frac{1}{h}\int\int_{s_{0}}^{s}\left\vert \left(  D_{1}f\right)  \left(
u,t\right)  -\left(  D_{1}f\right)  \left(  s_{0},t\right)  \right\vert
du\,dt\\
& =\frac{1}{h}\int_{s_{0}}^{s}\int\left\vert \left(  D_{1}f\right)  \left(
u,t\right)  -\left(  D_{1}f\right)  \left(  s_{0},t\right)  \right\vert
dt\,du\\
& \leq\frac{c}{h}\int_{s_{0}}^{s}\left(  u-s_{0}\right)  du\\
& =\frac{c}{2}\left(  s-s_{0}\right)  .
\end{align*}

The same estimate can be obtained when $s<s_{0}$.
\end{proof}

We have shown that the basis function lies in $C_{B}^{\left(  2n-2\right)  }$.
Now we will use the lemma to consider the derivatives of order $2n-1$.

\begin{theorem}
\label{Thm2_G_basis_def_2_new}Suppose $G_{1}$ is a \textbf{1-dimensional}
central difference basis function generated by a weight function $w\left(
\xi\right)  =\frac{\xi^{2n}}{\Delta^{2l}\widehat{q}\left(  \xi\right)  }$
satisfying property W02. Suppose also that $q$ is bounded a.e. i.e. for some
constant $c_{q}>0$,
\begin{equation}
q\left(  t\right)  \leq c_{q}\text{\quad}a.e.\label{a977}%
\end{equation}

Then, if $G_{s}$ is given by \ref{a926},%
\begin{equation}
D^{2n-1}G_{1}\left(  s\right)  =\tfrac{2^{2l}}{\sqrt{2\pi}}\int\left(
D^{2n-1}G_{s}\right)  \left(  \frac{2s}{\left\vert t\right\vert }\right)
q\left(  t\right)  dt,\text{\quad}s\in\mathbb{R}^{1}\setminus0,\label{a942}%
\end{equation}

and $D^{2n-1}G_{1}\in C^{\left(  0\right)  }\left(  \mathbb{R}^{1}%
\setminus0\right)  $ and $D^{2n-1}G_{1}$ is essentially bounded. Further, for
$0\leq\tau\leq2n-1$,%
\[
\left\vert D^{2n-1}G_{1}\left(  s\right)  \right\vert \leq\frac{2^{-\tau
}\left\Vert x^{\tau}D^{2n-1}G_{s}\right\Vert _{\infty}\left\Vert t^{\tau
}q\right\Vert _{1}}{\left\vert s\right\vert ^{\tau}},\text{\quad}%
s\in\mathbb{R}^{1}\setminus0.
\]

\end{theorem}

\begin{proof}
We will define%
\[
f\left(  s,t\right)  :=\left(  D^{2n-2}G_{s}\right)  \left(  \frac
{2s}{\left\vert t\right\vert }\right)  \left\vert t\right\vert q\left(
t\right)  ,
\]

and verify that $f$ satisfies the three properties of Lemma
\ref{Lem_Leibnitz_diff_under_integ}. First note that from Theorem
\ref{Thm_basis_tensor_hat_W3}, $G_{s}\in C_{0}^{\left(  2n-2\right)  }$ and
$D^{2n-1}G_{s}$ is a piecewise constant function with bounded support. Also,
since $G_{1}$ has weight function property W02, $q\in L^{1}$ and
$\int_{\left\vert \cdot\right\vert \geq R}\left\vert t\right\vert
^{2n-1}q\left(  t\right)  dt<\infty$ so that $\int\left\vert t\right\vert
q\left(  t\right)  dt<\infty$.\medskip

\textbf{Property 1} Thus for each $t\neq0$, $f\left(  s,t\right)  \in
C^{\left(  0\right)  }$ and since
\begin{equation}
D_{1}f\left(  s,t\right)  =2\left(  D^{2n-1}G_{s}\right)  \left(  \frac
{2s}{\left\vert t\right\vert }\right)  q\left(  t\right)  ,\label{a928}%
\end{equation}

it follows that $D_{1}f\left(  \cdot,t\right)  $ is piecewise
continuous.\medskip

\textbf{Property 2} Since $G_{s}\in C_{0}^{\left(  2n-2\right)  }$,
$D^{2n-2}G_{s}$ is bounded and hence%
\[
\int\left\vert f\left(  s,t\right)  \right\vert dt\leq\left\Vert D^{2n-2}%
G_{s}\right\Vert _{\infty}\int\left\vert t\right\vert q\left(  t\right)
dt<\infty,
\]

and%
\[
\int\left\vert D_{1}f\left(  s,t\right)  \right\vert dt\leq2\int\left\vert
\left(  D^{2n-1}G_{s}\right)  \left(  \frac{2s}{\left\vert t\right\vert
}\right)  \right\vert q\left(  t\right)  dt\leq2\left\Vert D^{2n-1}%
G_{s}\right\Vert _{\infty}\int q\left(  t\right)  dt<\infty.
\]
\medskip

\textbf{Property 3} Assume $s_{0}\neq0$. Suppose $\operatorname*{supp}%
D^{2n-1}G_{s}\subset B\left(  0;r\right)  $ and set $t_{\min}=2\max\left\{
\left\vert s_{0}\right\vert ,\left\vert s\right\vert \right\}  /r$. Then%
\begin{align*}
\int\left\vert D_{1}f\left(  s,t\right)  -D_{1}f\left(  s_{0},t\right)
\right\vert dt  & =\int\left\vert \left(  D^{2n-1}G_{s}\right)  \left(
\frac{2s}{\left\vert t\right\vert }\right)  -\left(  D^{2n-1}G_{s}\right)
\left(  \frac{2s_{0}}{\left\vert t\right\vert }\right)  \right\vert q\left(
t\right)  dt\\
& =\int\limits_{\left\vert t\right\vert \geq t_{\min}}\left\vert \left(
D^{2n-1}G_{s}\right)  \left(  \frac{2s}{\left\vert t\right\vert }\right)
-\left(  D^{2n-1}G_{s}\right)  \left(  \frac{2s_{0}}{\left\vert t\right\vert
}\right)  \right\vert q\left(  t\right)  dt.
\end{align*}

Since $D^{2n-1}G_{s}$ is a step function with bounded support it can be
written as the sum of a finite number of rectangular functions, say $\left\{
a_{k}\Pi_{k}\right\}  _{k=1}^{m}$, so that%
\begin{align}
\int\left\vert D_{1}f\left(  s,t\right)  -D_{1}f\left(  s_{0},t\right)
\right\vert dt  & \leq\sum_{k=1}^{m}a_{k}\int\limits_{\left\vert t\right\vert
\geq t_{\min}}\left\vert \Pi_{k}\left(  \frac{2s}{\left\vert t\right\vert
}\right)  -\Pi_{k}\left(  \frac{2s_{0}}{\left\vert t\right\vert }\right)
\right\vert q\left(  t\right)  dt\nonumber\\
& \leq\left\Vert D^{2n-1}G_{s}\right\Vert _{\infty}\sum_{k=1}^{m}%
\int\limits_{\left\vert t\right\vert \geq t_{\min}}\left\vert \Pi_{k}\left(
\frac{2s}{\left\vert t\right\vert }\right)  -\Pi_{k}\left(  \frac{2s_{0}%
}{\left\vert t\right\vert }\right)  \right\vert q\left(  t\right)
dt.\label{a941}%
\end{align}

Suppose $\Pi_{k}$ has support $\left[  x_{1}^{\left(  k\right)  }%
,x_{2}^{\left(  k\right)  }\right]  $. Observing that
\[
\left\{  x:\Pi_{k}\left(  x\right)  -\Pi_{k}\left(  x+h\right)  \neq0\right\}
\subseteq\bigcup\limits_{j=1}^{2}\left\{  x:\left\vert x-x_{j}^{\left(
k\right)  }\right\vert \leq\left\vert h\right\vert \right\}  ,
\]

we have%
\begin{align*}
\left\{  t:\Pi_{k}\left(  \frac{2s_{0}}{\left\vert t\right\vert }\right)
\neq\Pi_{k}\left(  \frac{2s}{\left\vert t\right\vert }\right)  ,\quad
\left\vert t\right\vert \geq t_{\min}\right\}   & \subseteq\bigcup
\limits_{j=1}^{2}\left\{  t:\left\vert \frac{2s_{0}}{\left\vert t\right\vert
}-x_{j}^{\left(  k\right)  }\right\vert \leq\left\vert \frac{2s}{\left\vert
t\right\vert }-\frac{2s_{0}}{\left\vert t\right\vert }\right\vert
,\quad\left\vert t\right\vert \geq t_{\min}\right\} \\
& =\bigcup\limits_{j=1}^{2}\left\{  t:\left\vert \frac{2s_{0}}{\left\vert
t\right\vert }-x_{j}^{\left(  k\right)  }\right\vert \leq\frac{\left\vert
2s-2s_{0}\right\vert }{\left\vert t\right\vert },\quad\left\vert t\right\vert
\geq t_{\min}\right\}  .
\end{align*}

Next apply the further constraint $\left\vert 2s-2s_{0}\right\vert <\left\vert
2s_{0}\right\vert /2$ so that $x_{j}^{\left(  k\right)  }=0$ implies the
corresponding set is empty. If $x_{1},x_{2}\neq0$%
\begin{align*}
\left\{  t:\Pi_{k}\left(  \frac{2s_{0}}{\left\vert t\right\vert }\right)
\neq\Pi_{k}\left(  \frac{2s}{\left\vert t\right\vert }\right)  ,\quad
\left\vert t\right\vert \geq t_{\min}\right\}   & \subseteq\bigcup
\limits_{j=1}^{2}\left\{  t:\left\vert \left\vert t\right\vert -\frac{2s_{0}%
}{x_{j}^{\left(  k\right)  }}\right\vert \leq\frac{\left\vert 2s-2s_{0}%
\right\vert }{\left\vert x_{j}^{\left(  k\right)  }\right\vert }%
,\quad\left\vert t\right\vert \geq t_{\min}\right\} \\
& =I_{1}^{\left(  k\right)  }+I_{2}^{\left(  k\right)  },
\end{align*}

with $0\notin I_{j}^{\left(  k\right)  }$. We can now conclude from equation
\ref{a941} that when $\left\vert s-s_{0}\right\vert <\left\vert s_{0}%
\right\vert /2$%
\begin{align*}
\int\left\vert D_{1}f\left(  s,t\right)  -D_{1}f\left(  s_{0},t\right)
\right\vert dt  & \leq\left\Vert D^{2n-1}G_{s}\right\Vert _{\infty}\sum
_{k=1}^{m}\int_{I_{1}^{\left(  k\right)  }\cup I_{2}^{\left(  k\right)  }}q\\
& =\left\Vert D^{2n-1}G_{s}\right\Vert _{\infty}\int_{J\left(  2s_{0}%
,2s\right)  }q,
\end{align*}

where, $J\left(  2s_{0},2s\right)  =\bigcup\limits_{k=1}^{m}\left(
I_{1}^{\left(  k\right)  }\cup I_{2}^{\left(  k\right)  }\right)
=\bigcup\limits_{k=1}^{m}\left[  a_{k}-b_{k}\left\vert 2s-2s_{0}\right\vert
,a_{k}+b_{k}\left\vert 2s-2s_{0}\right\vert \right]  $ with $b_{k}>0$,
$a_{k}\neq0$, and $0\notin J\left(  2s_{0},2s\right)  $.

Now the additional boundedness assumption \ref{a977} we made about $q$ comes
into play. In fact%
\begin{align*}
\int\limits_{J\left(  s_{0},s\right)  }q\left(  t\right)  dt=\sum_{k=1}%
^{m}\int\limits_{a_{k}-b_{k}\left\vert 2s-2s_{0}\right\vert }^{a_{k}%
+b_{k}\left\vert 2s-2s_{0}\right\vert }q\leq c_{q}\sum_{k=1}^{m}%
\int\limits_{a_{n}-b_{n}\left\vert 2s-2s_{0}\right\vert }^{a_{n}%
+b_{n}\left\vert 2s-2s_{0}\right\vert }dt &  \leq c_{q}\sum_{k=1}^{m}%
4b_{k}\left\vert s-s_{0}\right\vert \\
&  =4c_{q}\left(  \sum_{k=1}^{m}b_{k}\right)  \left\vert s-s_{0}\right\vert ,
\end{align*}

so that%
\[
\int\left\vert D_{1}f\left(  s,t\right)  -D_{1}f\left(  s_{0},t\right)
\right\vert dt\leq4c_{q}\left(  \sum_{k=1}^{m}b_{k}\right)  \left\vert
s-s_{0}\right\vert ,
\]

when $\left\vert s-s_{0}\right\vert <\left\vert s_{0}\right\vert /2$ and
$s_{0}\neq0$. This proves property 3 and hence \ref{a942} and $D^{2n-1}%
G_{1}\in C^{\left(  0\right)  }\left(  \mathbb{R}^{1}\setminus0\right)  $.

From \ref{a942}%
\begin{align*}
\left\vert D^{2n-1}G_{1}\left(  s\right)  \right\vert  & \leq\int%
_{\mathbb{R}^{1}}\left\vert \left(  D^{2n-1}G_{s}\right)  \left(  \frac{2s}%
{t}\right)  \right\vert q\left(  t\right)  dt\\
& \leq\left\Vert D^{2n-1}G_{s}\right\Vert _{\infty}\int q\left(  t\right)
dt<\infty\text{ }a.e.
\end{align*}

Again from \ref{a942}, since $D^{2n-1}G_{s}$ is bounded with bounded support,
for any $0\leq\tau\leq2n-1$
\begin{align*}
\left\vert s\right\vert ^{\tau}\left\vert D^{2n-1}G_{1}\left(  s\right)
\right\vert \leq\int\left\vert \left\vert s\right\vert ^{\tau}\left(
D^{2n-1}G_{s}\right)  \left(  \frac{2s}{t}\right)  \right\vert q\left(
t\right)  dt &  =\int\left\vert \left\vert \frac{s}{t}\right\vert ^{\tau
}\left(  D^{2n-1}G_{s}\right)  \left(  \frac{2s}{t}\right)  \right\vert
\left\vert t\right\vert ^{\tau}q\left(  t\right)  dt\\
&  \leq2^{-\tau}\left\Vert \left\vert x\right\vert ^{\tau}\left(
D^{2n-1}G_{s}\right)  \right\Vert _{\infty}\int\left\vert t\right\vert ^{\tau
}q\left(  t\right)  dt\\
&  =2^{-\tau}\left\Vert \left\vert \cdot\right\vert ^{\tau}D^{2n-1}%
G_{s}\right\Vert _{\infty}\left\Vert \left\vert \cdot\right\vert ^{\tau
}q\right\Vert _{1}.
\end{align*}

\end{proof}

\begin{example}
\label{Ex_central_diff_basis_hat_squ}In this example we show using Theorem
\ref{Thm_G_basis_def_2} that the univariate function $\Lambda\left(  x\right)
^{2}$ is a central difference basis function with weight function $w\left(
\xi\right)  =\frac{\sqrt{2\pi}}{4}\frac{\xi^{2}}{1-\frac{\sin\xi}{\xi}}$.

Suppose $l=n=1$ and $q\left(  x\right)  =R\left(  x\right)  $ where $R$ is the
rectangular function with support $\left[  -1,1\right]  $ and height $1$. From
\ref{a926}, $G_{s}\left(  s\right)  =\tfrac{\sqrt{2\pi}}{2}\Lambda\left(
\tfrac{s}{2}\right)  $ i.e. a scaled hat function. Then using formula
\ref{a945} the central difference basis function is given by%
\[
G\left(  x\right)  =\tfrac{2}{\sqrt{2\pi}}\int_{\mathbb{R}^{1}}\tfrac
{\sqrt{2\pi}}{2}\Lambda\left(  \tfrac{x}{t}\right)  \left\vert t\right\vert
R\left(  t\right)  dt=\int_{\mathbb{-}1}^{1}\Lambda\left(  \tfrac{x}%
{t}\right)  \left\vert t\right\vert dt=2\int_{0}^{1}\Lambda\left(  \frac{x}%
{t}\right)  tdt,
\]

so that $G\left(  x\right)  =0$ when $\left\vert x\right\vert \geq1$, since
$\Lambda\left(  x\right)  =0$ when $\left\vert x\right\vert \geq1$.

When $\left\vert x\right\vert \leq1$
\begin{align*}
G\left(  x\right)  =2\int_{\left\vert x\right\vert }^{1}\left(  1-\frac
{\left\vert x\right\vert }{t}\right)  tdt=2\int_{\left\vert x\right\vert }%
^{1}\left(  t-\left\vert x\right\vert \right)  dt &  =\left[  t^{2}%
-2\left\vert x\right\vert t\right]  _{\left\vert x\right\vert }^{1}\\
&  =\left(  1-2\left\vert x\right\vert \right)  -\left(  \left\vert
x\right\vert ^{2}-2\left\vert x\right\vert ^{2}\right) \\
&  =1-2\left\vert x\right\vert +\left\vert x\right\vert ^{2}\\
&  =\left(  1-\left\vert x\right\vert \right)  ^{2}.
\end{align*}

Thus $G\left(  x\right)  =\Lambda\left(  x\right)  ^{2}$.

Since $\widehat{q}\left(  \xi\right)  =\frac{2}{\sqrt{2\pi}}\frac{\sin\xi}%
{\xi}$, by the definition of $w$,
\begin{align*}
\frac{1}{w\left(  \xi\right)  }=\frac{\Delta_{2l}\widehat{q}\left(
\xi\right)  }{\xi^{2n}}=\frac{\Delta_{2}\widehat{q}\left(  \xi\right)  }%
{\xi^{2}} &  =\frac{-\left(  \widehat{q}\left(  \xi\right)  -2\widehat{q}%
\left(  0\right)  +\widehat{q}\left(  -\xi\right)  \right)  }{\xi^{2}}\\
&  =\frac{2}{\sqrt{2\pi}}\frac{2-2\frac{\sin\xi}{\xi}}{\xi^{2}}\\
&  =\frac{4}{\sqrt{2\pi}}\frac{\xi-\sin\xi}{\xi^{3}},
\end{align*}

and thus $w\left(  \xi\right)  =\frac{\sqrt{2\pi}}{4}\frac{\xi^{2}}%
{1-\frac{\sin\xi}{\xi}}$.
\end{example}

\begin{example}
\label{Ex_ten_prod_basis_central_diff}??? \textbf{DELETE THIS}? Noting the
convergence claims before Definition \ref{Def_central_diff_wt_func} we will
now construct a sequence of central difference basis functions which converges
uniformly pointwise to the B-spline basis function with $n=l=1$. See Theorem
\ref{Thm_centr_diff_basis_conv_to_ext_spl_basis} for a more general result. We
construct the sequence \ref{a978} by choosing $\psi=R$ so that%
\[
q_{k}\left(  s\right)  =\frac{\sqrt{2\pi}}{2}\frac{k}{2}R\left(  \frac{k}%
{2}\left(  s-2\right)  \right)  ,\quad k=1,2,3,\ldots,
\]

where $R$ is the unit rectangular function with support $\left[  -1,1\right]
$. Suppose $l=n=1$ so that the central difference basis function is given by
\ref{a945} and \ref{a926} so we get the sequence%
\begin{align*}
G_{k}\left(  x\right)   & =\tfrac{2}{\sqrt{2\pi}}\int G_{s}\left(  \frac
{2x}{s}\right)  \left\vert s\right\vert q_{k}\left(  s\right)  ds=\tfrac
{2}{\sqrt{2\pi}}\int G_{s}\left(  \frac{2x}{s}\right)  \left\vert s\right\vert
\frac{\sqrt{2\pi}}{2}\frac{k}{2}R\left(  \frac{k}{2}\left(  s-2\right)
\right)  ds=\\
& =\int G_{s}\left(  \frac{2x}{s}\right)  \left\vert s\right\vert \frac{k}%
{2}R\left(  \frac{k}{2}\left(  s-2\right)  \right)  ds\\
& :t=s/2,\text{ }ds=2dt\Rightarrow\\
& =\int G_{s}\left(  \frac{x}{t}\right)  \left\vert t\right\vert kR\left(
k\left(  t-1\right)  \right)  dt=k\int_{1-\frac{1}{k}}^{1+\frac{1}{k}}%
G_{s}\left(  \frac{x}{t}\right)  tdt=\\
& :G_{s}\left(  s\right)  =\tfrac{\sqrt{2\pi}}{2}\Lambda\left(  \tfrac{s}%
{2}\right)  \Rightarrow\\
& =\frac{\sqrt{2\pi}}{2}k\int_{1-\frac{1}{k}}^{1+\frac{1}{k}}\Lambda\left(
\frac{x}{2t}\right)  tdt.
\end{align*}

The corresponding weight functions are derived in Corollary
\ref{Cor_cntl_diff_bnd_on_wt_fn_2}.

We proceed by considering $x$ in three separate sets:\medskip

\fbox{\textbf{Case 1}: $0\leq x<1-\frac{1}{k}$}\smallskip%
\begin{align*}
G_{k}\left(  x\right)   & =\frac{k}{2}\int_{1-\frac{1}{k}}^{1+\frac{1}{k}%
}\Lambda\left(  \frac{x}{t}\right)  tdt=\frac{k}{2}\int_{1-\frac{1}{k}%
}^{1+\frac{1}{k}}\left(  1-\frac{x}{t}\right)  tdt=\frac{k}{2k}\int%
_{1-\frac{1}{k}}^{1+\frac{1}{k}}\left(  t-x\right)  dt=\\
& =\frac{k}{2}\left[  \frac{1}{2}t^{2}-tx\right]  _{1-\frac{1}{k}}^{1+\frac
{1}{k}}\\
& =\frac{k}{2}\left(  \frac{1}{2}\left(  1+\frac{1}{k}\right)  ^{2}-\left(
1+\frac{1}{k}\right)  x\right)  -\frac{1}{k}\left(  -\left(  1-\frac{1}%
{k}\right)  x\right) \\
& =\frac{k}{2}\left(  \frac{1}{2}\left(  1+\frac{1}{k}\right)  ^{2}-\left(
1+\frac{1}{k}\right)  x-\frac{1}{2}\left(  1-\frac{1}{k}\right)  ^{2}+\left(
1-\frac{1}{k}\right)  x\right) \\
& =k\left(  \frac{1}{k}-\frac{x}{k}\right)  =\left(  1-x\right)
=\Lambda\left(  x\right)  .
\end{align*}
\medskip

\fbox{\textbf{Case 2}: $1-\frac{1}{k}\leq x<1+\frac{1}{k}$} Then\smallskip%
\begin{align*}
G_{k}\left(  x\right)  =\frac{k}{2}\int_{x}^{1+\frac{1}{k}}\Lambda\left(
\frac{x}{t}\right)  tdt  & =\frac{k}{2}\int_{x}^{1+\frac{1}{k}}\left(
1-\frac{x}{t}\right)  tdt=\frac{k}{2}\int_{x}^{1+\frac{1}{k}}\left(
t-x\right)  dt\\
& =\frac{k}{2}\left[  \frac{1}{2}t^{2}-tx\right]  _{x}^{1+\frac{1}{k}}\\
& =\frac{k}{2}\left(  \left(  \frac{1}{2}\left(  1+\frac{1}{k}\right)
^{2}-\left(  1+\frac{1}{k}\right)  x\right)  -\left(  \frac{x^{2}}{2}%
-x^{2}\right)  \right) \\
& =\frac{k}{2}\left(  \frac{1}{2}\left(  1+\frac{1}{k}\right)  ^{2}-\left(
1+\frac{1}{k}\right)  x+\frac{1}{2}x^{2}\right) \\
& =\frac{k}{4}\left(  \left(  1+\frac{1}{k}\right)  ^{2}-2\left(  1+\frac
{1}{k}\right)  x+x^{2}\right) \\
& =\frac{k}{4}\left(  1+\frac{1}{k}-x\right)  ^{2}.
\end{align*}

\fbox{\textbf{Case 3}: $1+\frac{1}{k}\leq x$} Clearly $G_{k}\left(  x\right)
=0$. From the three cases we conclude that%
\begin{equation}
G_{k}\left(  x\right)  =\left\{
\begin{array}
[c]{ll}%
\Lambda\left(  x\right)  , & 0\leq\left\vert x\right\vert \leq1-\frac{1}{k},\\
\frac{k}{4}\left(  1+\frac{1}{k}-\left\vert x\right\vert \right)  ^{2}, &
1-\frac{1}{k}\leq\left\vert x\right\vert \leq1+\frac{1}{k},\\
0, & 1+\frac{1}{k}\leq\left\vert x\right\vert ,
\end{array}
\right. \label{a974}%
\end{equation}

and $G_{k}\in C\left(  \mathbb{R}^{1}\right)  \cap C^{\left(  1\right)
}\left(  \mathbb{R}^{1}\setminus0\right)  \cap PWC^{\infty}\left(
\mathbb{R}^{1}\setminus\left\{  -1,0,1\right\}  \right)  $. We can also write%
\[
G_{k}\left(  x\right)  =\Lambda\left(  x\right)  +\Delta_{k}\left(  x\right)
,
\]

where%
\begin{align}
\Delta_{k}\left(  x\right)   & =\left\{
\begin{array}
[c]{ll}%
0, & 0\leq\left\vert x\right\vert \leq1-\frac{1}{k},\\
\frac{k}{4}\left(  1-\left\vert x\right\vert +\frac{1}{k}\right)  ^{2}-\left(
1-\left\vert x\right\vert \right)  , & 1-\frac{1}{k}\leq\left\vert
x\right\vert \leq1,\\
\frac{k}{4}\left(  1-\left\vert x\right\vert +\frac{1}{k}\right)  ^{2}, &
1\leq\left\vert x\right\vert \leq1+\frac{1}{k},\\
0, & 1+\frac{1}{k}\leq\left\vert x\right\vert ,
\end{array}
\right. \nonumber\\
& =\left\{
\begin{array}
[c]{ll}%
0, & 0\leq\left\vert x\right\vert \leq1-\frac{1}{k},\\
\frac{k}{4}\left(  \left\vert x\right\vert -\left(  1-\frac{1}{k}\right)
\right)  ^{2}, & 1-\frac{1}{k}\leq\left\vert x\right\vert \leq1,\\
\frac{k}{4}\left(  1+\frac{1}{k}-\left\vert x\right\vert \right)  ^{2}, &
1\leq\left\vert x\right\vert \leq1+\frac{1}{k},\\
0, & 1+\frac{1}{k}\leq\left\vert x\right\vert .
\end{array}
\right. \label{a979}%
\end{align}

We note that
\[
\operatorname*{supp}\Delta_{k}=\left[  1-\frac{1}{k},1+\frac{1}{k}\right]
;\quad0\leq\left\vert \Delta_{k}\left(  x\right)  \right\vert \leq\frac{1}%
{4k}.
\]

\end{example}

\subsection{Lipschitz continuity\label{SbSect_CentDiffBasis_Lipschitz}}

The next two theorems will show that the continuous derivatives of the basis
function are uniformly Lipschitz continuous of order $1$ on $\mathbb{R}^{1}$.
We start with the univariate case.

\begin{theorem}
\label{Thm_cdiffbasis_Lips_dim_d}\textbf{Univariate Lipschitz continuity} Let
$w\left(  \xi\right)  =\frac{\xi^{2n}}{\Delta^{2l}\widehat{q}\left(
\xi\right)  }$ be the \textbf{univariate} central difference weight function
introduced in Definition \ref{Def_central_diff_wt_func} and which satisfies
property W02/W03. If $G_{1}$ is the basis function of order zero generated by
$w$ the derivatives $D^{k}G_{1}\left(  x\right)  $, $k\leq2n-2$, are uniformly
Lipschitz continuous of order $1$ on $\mathbb{R}^{1}$. In fact, for
$k\leq2n-2$,%
\[
\left\vert D^{k}G_{1}\left(  x\right)  -D^{k}G_{1}\left(  y\right)
\right\vert \leq\tfrac{2^{2\left(  l-n\right)  +k+1}}{\sqrt{2\pi}}\left\Vert
D^{k+1}G_{s}\right\Vert _{\infty}\left\Vert t^{2n-k-2}q\right\Vert
_{1}\left\vert x-y\right\vert ,\text{\quad}x,y\in\mathbb{R}^{1},
\]

where $G_{s}$ is the extended B-spline given by \ref{a926}.
\end{theorem}

\begin{proof}
From \ref{a940} of Theorem \ref{Thm_G_basis_def_2}%
\[
\left\vert D^{k}G_{1}\left(  x\right)  -D^{k}G_{1}\left(  y\right)
\right\vert \leq\tfrac{2^{2\left(  l-n\right)  +k+1}}{\sqrt{2\pi}}%
\int\limits_{\mathbb{R}^{1}}\left\vert \left(  D^{k}G_{s}\right)  \left(
\frac{2x}{\left\vert t\right\vert }\right)  -\left(  D^{k}G_{s}\right)
\left(  \frac{2y}{\left\vert t\right\vert }\right)  \right\vert \left\vert
t\right\vert ^{2n-k-1}q\left(  t\right)  dt.
\]

Since $G_{s}$ is the extended B-spline, from Theorem
\ref{Thm_ex_nat_spline_basis_Lipschitz} we have for dimension $1$ and
$k\leq2n-2 $,
\[
\left\vert D^{k}G_{s}\left(  s\right)  -D^{k}G_{s}\left(  s^{\prime}\right)
\right\vert \leq\left\Vert D^{k+1}G_{s}\right\Vert _{\infty}\left\vert
s-s^{\prime}\right\vert ,\quad s,s^{\prime}\in\mathbb{R}^{1},
\]

and thus%
\begin{align*}
\left\vert D^{k}G_{1}\left(  x\right)  -D^{k}G_{1}\left(  y\right)
\right\vert  & \leq\tfrac{2^{2\left(  l-n\right)  +k+1}}{\sqrt{2\pi}%
}\left\Vert D^{k+1}G_{s}\right\Vert _{\infty}\int\limits_{\mathbb{R}^{1}%
}\left\vert \frac{x}{\left\vert t\right\vert }-\frac{y}{\left\vert
t\right\vert }\right\vert \left\vert t\right\vert ^{2n-k-1}q\left(  t\right)
dt\\
& \leq\tfrac{2^{2\left(  l-n\right)  +k+1}}{\sqrt{2\pi}}\left\Vert
D^{k+1}G_{s}\right\Vert _{\infty}\left\Vert t^{2n-k-2}q\right\Vert
_{1}\left\vert x-y\right\vert .
\end{align*}

\end{proof}

To deal with the multivariate case we will require a lemma concerning a
distribution Taylor series expansion and which follows directly from Lemma
\ref{Lem_Taylor_extension} Chapter \ref{Ch_wtfn_basisfn_datasp} (But see
Remark \ref{Rem_Thm_ex_nat_spline_basis_Lipschitz}).

\begin{lemma}
\label{Lem_estim_G(x)_minus_G(y)}Suppose that $G\in C_{B}^{\left(  0\right)
}\left(  \mathbb{R}^{d}\right)  $ and the (distributional) derivatives
$\left\{  D^{\alpha}G\right\}  _{\left\vert \alpha\right\vert =1}$ are
essentially bounded functions. Then
\[
G\left(  x\right)  -G\left(  y\right)  =\sum_{\left\vert \alpha\right\vert
=1}\left(  x-y\right)  ^{\alpha}\int_{0}^{1}\left(  D^{\alpha}G\right)
\left(  x-t\left(  x-y\right)  \right)  dt,
\]

and
\begin{equation}
\left\vert G\left(  x\right)  -G\left(  y\right)  \right\vert \leq\sqrt{d}%
\max\limits_{\left\vert \alpha\right\vert =1}\left\Vert D^{\alpha}G\right\Vert
_{\infty}\left\vert x-y\right\vert .\label{a901}%
\end{equation}

\end{lemma}

Now we can show that the continuous derivatives of the tensor product central
difference basis functions are also Lipschitz continuous of order $1$.

\begin{theorem}
\label{Thm_cntr_diff_basis_Lips_all_dim}\textbf{Multivariate Lipschitz
continuity} Let $G_{c}\left(  x\right)  =\prod\limits_{k=1}^{d}G_{1}\left(
x_{k}\right)  $ be a central difference tensor product basis function with
parameters $n$, $l$, $q\left(  \cdot\right)  $, and suppose $q$ is bounded if
$n=1$. Then $G_{1}\in C_{B}^{\left(  2n-2\right)  }$, $D^{2n-1}G_{1}$ is
bounded and we have the estimates
\begin{equation}
\left\vert G_{c}\left(  x\right)  -G_{c}\left(  y\right)  \right\vert
\leq\sqrt{d}G_{1}\left(  0\right)  ^{d-1}\left\Vert DG_{1}\right\Vert
_{\infty}\left\vert x-y\right\vert ,\text{\quad}x,y\in\mathbb{R}%
^{d},\label{a946}%
\end{equation}

and if $0<\beta\leq\left(  2n-2\right)  \mathbf{1}$,%
\begin{equation}
\left\vert D^{\beta}G_{c}\left(  x\right)  -D^{\beta}G_{c}\left(  y\right)
\right\vert \leq\sqrt{d}\left(  \max\limits_{k=1}^{d}\left\Vert D_{k}D^{\beta
}G_{c}\right\Vert _{\infty}\right)  \left\vert x-y\right\vert ,\text{\quad
}x,y\in\mathbb{R}^{d}.\label{a949}%
\end{equation}

\end{theorem}

\begin{proof}
By Theorem \ref{Thm_G_basis_def_2}, $G_{1}\in C_{B}^{\left(  2n-2\right)  }$
and by Theorem \ref{Thm2_G_basis_def_2_new}, $D^{2n-1}G_{1}$ is bounded.

Thus $G_{c}\in C_{B}^{\left(  0\right)  }\left(  \mathbb{R}^{d}\right)  $ and
the derivatives $\left\{  D^{\alpha}G_{c}\right\}  _{\left\vert \alpha
\right\vert =1}$ are bounded functions. Thus $G_{c}$ satisfies the conditions
of Lemma \ref{Lem_estim_G(x)_minus_G(y)} and so the estimate of that lemma
holds i.e. for all $x,y\in\mathbb{R}^{d}$,
\[
\left\vert G_{c}\left(  x\right)  -G_{c}\left(  y\right)  \right\vert
\leq\sqrt{d}\left(  \max\limits_{\left\vert \alpha\right\vert =1}\left\Vert
D^{\alpha}G_{c}\right\Vert _{\infty}\right)  \left\vert x-y\right\vert
\leq\sqrt{d}\left\Vert G_{1}\right\Vert _{\infty}^{d-1}\left\Vert
DG_{1}\right\Vert _{\infty}\left\vert x-y\right\vert .
\]

Finally, from \ref{1.33}, $\left\vert G_{1}\left(  x\right)  \right\vert \leq
G_{1}\left(  0\right)  $ which proves \ref{a946}.

To prove \ref{a949} observe that if $0<k\leq2n-2$ then $D^{k}G_{1}\in
C_{B}^{\left(  2n-2-k\right)  }\left(  \mathbb{R}^{d}\right)  $ and
$D^{k+1}G_{1} $ is bounded. Thus if $0<\beta\leq\left(  2n-2\right)
\mathbf{1}$ then $D^{\beta}G_{c}\in C_{B}^{\left(  0\right)  }\left(
\mathbb{R}^{d}\right)  $ and the derivatives $\left\{  D^{\alpha}D^{\beta
}G_{c}\right\}  _{\left\vert \alpha\right\vert =1}$ are bounded functions.
Thus $D^{\beta}G_{c}$ satisfies the conditions of Lemma
\ref{Lem_estim_G(x)_minus_G(y)} and so the estimate of that lemma holds i.e.
for all $x,y\in\mathbb{R}^{d}$,%
\begin{align*}
\left\vert D^{\beta}G_{c}\left(  x\right)  -D^{\beta}G_{c}\left(  y\right)
\right\vert  & \leq\sqrt{d}\left(  \max\limits_{\left\vert \alpha\right\vert
=1}\left\Vert D^{\alpha}D^{\beta}G_{c}\right\Vert _{\infty}\right)  \left\vert
x-y\right\vert \\
& \leq\sqrt{d}\left(  \max\limits_{k=1}^{d}\left\Vert D_{k}D^{\beta}%
G_{c}\right\Vert _{\infty}\right)  \left\vert x-y\right\vert .
\end{align*}

\end{proof}

\subsection{A basis function formula using partial moments of $q$%
\label{SbSect_centdiff_basis_part_mom}}

Theorem \ref{Thm_G_basis_def_2} gave a multiplicative convolution formula
\ref{a945} for the central difference basis function in terms of the B-spline
basis function $G_{s}$. The next theorem gives another formula that does not
involve calculating the extended B-spline basis functions but instead involves
calculating partial moments of $q$. This result also gives information about
the support of the basis function and the corollary gives an improved
smoothness estimate.

\begin{theorem}
\label{Thm_cdiffbasis_part_moment_formula}\textbf{Partial moment formulas}
Suppose $w_{c}$ is a \textbf{1-dimensional} central difference weight function
with parameters $n,l,q$ and $w_{c}\in W02$, and denote the basis function by
$G_{c}$.

Then if $\operatorname*{supp}q\subseteq\overline{B}_{R_{q}}$, where possibly
$R_{q}=\infty$, it follows that $\operatorname*{supp}G_{c}\subseteq
\overline{B}_{R_{q}l}$ and%
\begin{equation}
G_{c}\left(  x\right)  =\left\{
\begin{array}
[c]{ll}%
\frac{\left(  -1\right)  ^{n}}{\left(  2n-1\right)  !}\int_{\left\vert
x\right\vert }^{R_{q}l}\left(  s-\left\vert x\right\vert \right)  ^{2n-1}%
q_{l}\left(  s\right)  ds, & \left\vert x\right\vert \leq R_{q}l,\\
0, & \left\vert x\right\vert \geq R_{q}l,
\end{array}
\right. \label{2.59}%
\end{equation}

where%
\begin{equation}
q_{l}\left(  s\right)  :=\sum\limits_{j=-l,j\neq0}^{l}\frac{\left(  -1\right)
^{j}}{\left\vert j\right\vert }\tbinom{2l}{j+l}q\left(  \frac{s}{j}\right)
.\label{2.60}%
\end{equation}

Further
\begin{equation}
D^{2n}G_{c}=\left(  -1\right)  ^{n}\left(  \Delta_{2l}\widehat{q}\right)
^{\vee}=\left(  -1\right)  ^{n}\tbinom{2l}{l}\left(  \int q\right)
\delta+\left(  -1\right)  ^{n}q_{l},\label{2.62}%
\end{equation}

so that $D^{2n}G_{c}$ has the same smoothness as $q$ on $\mathbb{R}%
^{d}\setminus0$.

Finally%
\begin{equation}
G_{c}\left(  x\right)  =\left\{
\begin{array}
[c]{ll}%
\left(  -1\right)  ^{n}\frac{2}{\left(  2n-1\right)  !}\sum
\limits_{\substack{j=1 \\j\geq\left\vert x\right\vert /R_{q}}}^{l}\left(
-1\right)  ^{j}\tbinom{2l}{j+l}\int_{\left\vert x\right\vert /j}^{R_{q}%
}\left(  jt-\left\vert x\right\vert \right)  ^{2n-1}q_{e}\left(  t\right)
dt, & \left\vert x\right\vert \leq R_{q}l,\\
0, & \left\vert x\right\vert \geq R_{q}l,
\end{array}
\right. \label{2.58}%
\end{equation}

where%
\begin{equation}
q_{e}\left(  t\right)  :=\frac{1}{2}\left(  q\left(  t\right)  +q\left(
-t\right)  \right)  ,\label{2.63}%
\end{equation}

is called the even component of $q$.
\end{theorem}

\begin{proof}
By Definition \ref{Def_central_diff_wt_func}, the weight function is given by%
\begin{equation}
w_{c}\left(  \xi\right)  =\frac{\xi^{2n}}{\Delta_{2l}\widehat{q}\left(
\xi\right)  },\quad\xi\in\mathbb{R}^{1},\label{2.14}%
\end{equation}

and from Theorem \ref{Thm_cdiffwt_2}, $n\leq l$ and the function $q$
satisfies
\[
\int_{\left\vert \cdot\right\vert \geq R}\left\vert \cdot\right\vert
^{2n-1}q<\infty,\quad for\text{ }some\text{ }R\geq0.
\]

This last condition implies%
\begin{equation}
\int\left\vert \cdot\right\vert ^{k}q<\infty,\quad k=0,1,\ldots
,2n-1.\label{2.15}%
\end{equation}

From \ref{2.14}, the basis function $G_{c}$ satisfies%
\begin{align*}
\xi^{2n}\widehat{G_{c}}\left(  \xi\right)   & =\frac{\xi^{2n}}{w_{c}\left(
\xi\right)  }=\Delta_{2l}\widehat{q}\left(  \xi\right)  .\\
\left(  -1\right)  ^{n}\widehat{D^{2n}G_{c}}\left(  \xi\right)   &
=\Delta_{2l}\widehat{q}\left(  \xi\right)  .\\
D^{2n}G_{c}\left(  x\right)   & =\left(  -1\right)  ^{n}\left(  \Delta
_{2l}\widehat{q}\right)  ^{\vee}\left(  x\right)  .
\end{align*}

but from \ref{a914} and part 12 of Definition \ref{Def_Fourier},%
\begin{align*}
\left(  \Delta_{2l}\widehat{q}\right)  ^{\vee}\left(  x\right)   &
=\sum_{j=-l}^{l}\left(  -1\right)  ^{j}\tbinom{2l}{j+l}\left(  \widehat{q}%
\left(  -j\xi\right)  \right)  ^{\vee}\left(  x\right) \\
& =\tbinom{2l}{l}\left(  \widehat{q}\left(  0\right)  \right)  ^{\vee}%
+\sum_{j=-l,j\neq0}^{l}\left(  -1\right)  ^{j}\tbinom{2l}{j+l}\left(
\widehat{q}\left(  -j\xi\right)  \right)  ^{\vee}\left(  x\right) \\
& =\tbinom{2l}{l}\widehat{q}\left(  0\right)  \overset{\vee}{1}+\sum
_{j=-l,j\neq0}^{l}\left(  -1\right)  ^{j}\tbinom{2l}{j+l}\left(
\widehat{q}\left(  -j\xi\right)  \right)  ^{\vee}\left(  x\right) \\
& =\tbinom{2l}{l}\widehat{q}\left(  0\right)  \left(  2\pi\right)  ^{\frac
{d}{2}}\delta+\sum_{j=-l,j\neq0}^{l}\left(  -1\right)  ^{j}\tbinom{2l}%
{j+l}\left(  \widehat{\frac{1}{\left\vert j\right\vert }q\left(  -\frac{\cdot
}{j}\right)  }\right)  ^{\vee}\left(  x\right) \\
& =\tbinom{2l}{l}\left(  \int q\right)  \delta+\sum_{j=-l,j\neq0}^{l}%
\frac{\left(  -1\right)  ^{j}}{\left\vert j\right\vert }\tbinom{2l}%
{j+l}q\left(  -\frac{x}{j}\right)  ,
\end{align*}

so%
\begin{align}
D^{2n}G_{c}\left(  x\right)   & =\left(  -1\right)  ^{n}\left(  \Delta
_{2l}\widehat{q}\right)  ^{\vee}\left(  x\right) \nonumber\\
& =\left(  -1\right)  ^{n}\tbinom{2l}{l}\left(  \int q\right)  \delta+\left(
-1\right)  ^{n}\sum_{j=-l,j\neq0}^{l}\frac{\left(  -1\right)  ^{j}}{\left\vert
j\right\vert }\tbinom{2l}{j+l}q\left(  -\frac{x}{j}\right)  .\label{2.11}%
\end{align}

For compactness, set%
\begin{equation}
\left.
\begin{array}
[c]{l}%
f\left(  x\right)  :=\left(  -1\right)  ^{n}\sum\limits_{j=-l,j\neq0}^{l}%
\frac{\left(  -1\right)  ^{j}}{\left\vert j\right\vert }\tbinom{2l}%
{j+l}q\left(  -\frac{x}{j}\right)  ,\\
a:=\left(  -1\right)  ^{n}\tbinom{2l}{l}\int q,
\end{array}
\right\} \label{2.53}%
\end{equation}

so that \ref{2.11} becomes%
\begin{equation}
D^{2n}G_{c}=a\delta+f.\label{2.54}%
\end{equation}

We want to solve \ref{2.11} for $G_{c}$. Now if $G_{c}\in C^{\left(  0\right)
}$ and $\lim\limits_{\left\vert x\right\vert \rightarrow\infty}G_{c}\left(
x\right)  =0$ satisfies \ref{2.11} then $G_{c}$ is unique. This is because two
different solutions differ by a polynomial.

Our approach will be to solve the two equations
\begin{equation}
D^{2n}G_{\delta}\left(  x\right)  =a\delta,\label{2.13}%
\end{equation}

and
\begin{equation}
D^{2n}G_{q}\left(  x\right)  =f,\label{2.12}%
\end{equation}

in such a way that%
\begin{equation}
G_{c}=G_{\delta}+G_{q},\label{2.56}%
\end{equation}

satisfies%
\begin{equation}
G_{c}\in C^{\left(  0\right)  },\quad G_{c}\left(  \pm\infty\right)
=0.\label{2.57}%
\end{equation}

\fbox{Solution of \ref{2.13}} From Example 2.6(f) of Vladimirov
\cite{Vladimirov},
\begin{equation}
\left.
\begin{array}
[c]{ll}%
G_{\delta}= & HP_{2n-1},\\
P_{k}= & a\frac{x^{k}}{k!},\text{ }k=0,1,2,\ldots,\\
D^{k}G_{\delta}\left(  0\right)  = & 0,\quad k<2n,\\
D^{k}G_{\delta}= & HP_{2n-k-1},\quad k\leq2n-1,
\end{array}
\right\} \label{2.16}%
\end{equation}
\medskip

is a (fundamental) solution of \ref{2.13}. Here $H$ denotes the Heavyside step
function.\medskip

\fbox{Solution of \ref{2.12}} To solve \ref{2.12} (see, for example, Section
II.34 of Zwillinger \cite{Zwill89}) we will integrate it by reversing the
order of integration using \ref{2.46}. Consider the following sequence of
multiple integrals%
\begin{equation}
\left.
\begin{array}
[c]{ll}%
g_{1}\left(  x_{1}\right)  & :=\int_{-\infty}^{x_{1}}f\left(  s\right)  ds,\\
g_{2}\left(  x_{2}\right)  & :=\int_{-\infty}^{x_{2}}\int_{-\infty}^{x_{1}%
}f\left(  s\right)  dsdx_{1},\\
\vdots & \vdots\\
g_{k}\left(  x_{k}\right)  & :=\int_{-\infty}^{x_{k}}\int_{-\infty}^{x_{k-1}%
}\ldots\int_{-\infty}^{x_{1}}f\left(  s\right)  ds\ldots dx_{k-2}dx_{k-1},\\
\vdots & \vdots\\
g_{2n}\left(  x_{2n}\right)  & :=\int_{-\infty}^{x_{2n}}\int_{-\infty
}^{x_{2n-1}}\ldots\int_{-\infty}^{x_{1}}f\left(  s\right)  ds\ldots
dx_{2n-2}dx_{2n-1}.
\end{array}
\right\} \label{2.45}%
\end{equation}

By the theory of differentiation (bounded variation, absolute continuity, the
fundamental theorem of calculus) e.g. \S 5.9 of Kuller \cite{Kuller1969},%
\begin{equation}
g_{1}\in C_{B}^{\left(  0\right)  },\text{ }g_{1}\in BV\left[  a,b\right]
,\text{ }Dg_{1}=f\text{ }a.e.,\text{ }g_{1}\left(  -\infty\right)
=0.\label{2.47}%
\end{equation}

and $g_{1}$ has only step discontinuities - these are countable. Here $\left[
a,b\right]  $ is any finite interval.

We now employ the (reverse the order of integration) identity: if $h,g\geq0$
and the right side of \ref{2.46} exists then%
\begin{equation}
\int_{a}^{x_{2}}h\left(  x_{1}\right)  \int_{a}^{x_{1}}g\left(  s\right)
dsdx_{1}=\int_{a}^{x_{2}}\left(  \int_{s}^{x_{2}}h\left(  t\right)  dt\right)
g\left(  s\right)  ds.\label{2.46}%
\end{equation}

In fact, if we define
\[
\phi\left(  s,x_{1}\right)  :=\left\{
\begin{array}
[c]{ll}%
g\left(  s\right)  h\left(  x_{1}\right)  , & a\leq s\leq x_{1},\\
0, & x_{1}\leq s\leq x,
\end{array}
\right.
\]

then by Fubini's theorem%
\begin{align*}
\int_{a}^{x_{2}}\left(  \int_{s}^{x_{2}}h\left(  t\right)  dt\right)  g\left(
s\right)  ds=\int_{a}^{x_{2}}\int_{s}^{x_{2}}\phi\left(  s,x_{1}\right)
dx_{1}ds &  \overset{Fubini}{=}\int_{a}^{x_{2}}\int_{a}^{x_{2}}\phi\left(
s,x_{1}\right)  dx_{1}ds\\
&  =\int_{a}^{x_{2}}\int_{a}^{x_{2}}\phi\left(  s,x_{1}\right)  dsdx_{1}\\
&  =\int_{a}^{x_{2}}h\left(  x_{1}\right)  \int_{a}^{x_{1}}g\left(  s\right)
dsdx_{1}.
\end{align*}

Regarding the equations \ref{2.45}, the constraints \ref{2.15} imply by the
theory of differentiation that%
\begin{equation}
\left.
\begin{array}
[c]{l}%
g_{2}\left(  x_{2}\right)  =\int_{-\infty}^{x_{2}}\left(  x_{2}-s\right)
f\left(  s\right)  ds,\\
g_{1}\in L^{1},\text{ }Dg_{2}=g_{1}\text{ }a.e.,\\
g_{2}\in C_{B}^{\left(  1\right)  },\text{ }g_{2}\in BV\left[  a,b\right]
,\text{ }g_{2}\left(  -\infty\right)  =0.
\end{array}
\right\} \label{2.48}%
\end{equation}

where $\left[  a,b\right]  $ is any finite interval.

Similarly%
\begin{align*}
g_{3}\left(  x_{3}\right)   & =\int_{-\infty}^{x_{3}}\int_{-\infty}^{x_{2}%
}\int_{-\infty}^{x_{1}}f\left(  s\right)  dsdx_{1}dx_{2}=\int_{-\infty}%
^{x_{3}}g_{2}\left(  x_{2}\right)  dx_{2}=\int_{-\infty}^{x_{3}}\int_{-\infty
}^{x_{2}}\left(  x_{2}-s\right)  f\left(  s\right)  ds=\\
& =\int_{-\infty}^{x_{3}}\left(  \int_{s}^{x_{3}}\left(  t-s\right)
dt\right)  f\left(  s\right)  ds=\frac{1}{2!}\int_{-\infty}^{x_{3}}\left(
x_{3}-s\right)  ^{2}f\left(  s\right)  ds,
\end{align*}

etc. so that%
\begin{equation}
g_{k}\left(  x_{k}\right)  =\frac{1}{\left(  k-1\right)  !}\int_{-\infty
}^{x_{k}}\left(  x_{k}-s\right)  ^{k-1}f\left(  s\right)  ds,\quad
k=1,\ldots,2n,\label{2.52}%
\end{equation}

and%
\begin{equation}
\left.
\begin{array}
[c]{l}%
g_{0}:=f,\\
g_{k-1}\in L^{1},\text{ }Dg_{k}=g_{k-1}\text{ }a.e.,\\
g_{k}\in C_{B}^{\left(  k-1\right)  },\text{ }g_{k}\in BV\left[  a,b\right]
,\text{ }g_{k}\left(  -\infty\right)  =0.
\end{array}
\right\}  ,\quad k=1,\ldots,2n,\label{2.49}%
\end{equation}

where $\left[  a,b\right]  $ is any finite interval. Clearly we can conclude
that%
\begin{equation}
\left.
\begin{array}
[c]{l}%
D^{2n}g_{2n}=f\text{ }a.e.,\\
g_{2n}\in C_{B}^{\left(  2n-1\right)  }\cap BV\left[  a,b\right]  ,\text{
}g_{2n}\left(  -\infty\right)  =0,
\end{array}
\right\} \label{2.51}%
\end{equation}

where $\left[  a,b\right]  $ is any finite interval, and hence $g_{2n}$ is a
solution of \ref{2.12}. Thus we try%
\begin{equation}
G_{q}:=g_{2n}.\label{2.55}%
\end{equation}

From \ref{2.52} and then \ref{2.53},%
\begin{align}
G_{q}\left(  x\right)   &  =\frac{1}{\left(  2n-1\right)  !}\int_{-\infty}%
^{x}\left(  x-s\right)  ^{2n-1}f\left(  s\right)  ds\nonumber\\
&  =\frac{1}{\left(  2n-1\right)  !}\int\limits_{-\infty}^{x}\left(
x-s\right)  ^{2n-1}\left(  \left(  -1\right)  ^{n}\sum_{j=-l,j\neq0}^{l}%
\frac{\left(  -1\right)  ^{j}}{\left\vert j\right\vert }\tbinom{2l}%
{j+l}q\left(  -\frac{s}{j}\right)  \right)  ds\nonumber\\
&  =\frac{\left(  -1\right)  ^{n}}{\left(  2n-1\right)  !}\sum_{j=-l,j\neq
0}^{l}\frac{\left(  -1\right)  ^{j}}{\left\vert j\right\vert }\tbinom{2l}%
{j+l}\int\limits_{-\infty}^{x}\left(  x-s\right)  ^{2n-1}q\left(  -\frac{s}%
{j}\right)  ds.\label{2.50}%
\end{align}

From the statement of this theorem $R_{q}=\max\limits_{x\in
\operatorname*{supp}q}\left\vert x\right\vert $ and thus we have the sequence
of implications
\[
x\leq-R_{q}l\Rightarrow s\leq x\leq-R_{q}l\Rightarrow-s\geq R_{q}%
l\Rightarrow\left\vert -\frac{s}{j}\right\vert \geq R_{q}\frac{l}{\left\vert
j\right\vert }\geq R_{q},
\]
since $\left\vert j\right\vert \leq l$, which means that $G_{q}\left(
x\right)  =0$ when $x\leq-R_{q}l$, and so we can write%
\begin{equation}
G_{q}\left(  x\right)  =\left\{
\begin{array}
[c]{ll}%
0, & x\leq-R_{q}l,\\
\frac{\left(  -1\right)  ^{n}}{\left(  2n-1\right)  !}\sum\limits_{j=-l,j\neq
0}^{l}\frac{\left(  -1\right)  ^{j}}{\left\vert j\right\vert }\tbinom{2l}%
{j+l}\int\limits_{-R_{q}l}^{x}\left(  x-s\right)  ^{2n-1}q\left(  -\frac{s}%
{j}\right)  ds, & x\geq-R_{q}l.
\end{array}
\right. \label{1.083}%
\end{equation}
\medskip

\fbox{Now suppose that $x\geq R_{q}l$} Split the domain of
integration:\smallskip%
\begin{align*}
G_{q}\left(  x\right)  = &  \frac{\left(  -1\right)  ^{n}}{\left(
2n-1\right)  !}\sum_{j=-l,j\neq0}^{l}\frac{\left(  -1\right)  ^{j}}{\left\vert
j\right\vert }\tbinom{2l}{j+l}\int\limits_{-R_{q}l}^{R_{q}l}\left(
x-s\right)  ^{2n-1}q\left(  -\frac{s}{j}\right)  ds+\\
&  +\frac{\left(  -1\right)  ^{n}}{\left(  2n-1\right)  !}\sum_{j=-l,j\neq
0}^{l}\frac{\left(  -1\right)  ^{j}}{\left\vert j\right\vert }\tbinom{2l}%
{j+l}\int\limits_{R_{q}l}^{x}\left(  x-s\right)  ^{2n-1}q\left(  -\frac{s}%
{j}\right)  ds.
\end{align*}

Now $x\geq R_{q}l\Rightarrow R_{q}l\leq s\leq x\Rightarrow\left\vert -\frac
{s}{j}\right\vert =\frac{s}{\left\vert j\right\vert }\geq R_{q}\frac
{l}{\left\vert j\right\vert }\geq R_{q}$ since $\left\vert j\right\vert \leq
l$, which means that%
\[
\frac{\left(  -1\right)  ^{n}}{\left(  2n-1\right)  !}\sum_{j=-l,j\neq0}%
^{l}\frac{\left(  -1\right)  ^{j}}{\left\vert j\right\vert }\tbinom{2l}%
{j+l}\int_{R_{q}l}^{x}\left(  x-s\right)  ^{2n-1}q\left(  -\frac{s}{j}\right)
ds=0,\quad x\geq R_{q}l,
\]

and so%
\[
G_{q}\left(  x\right)  =\frac{\left(  -1\right)  ^{n}}{\left(  2n-1\right)
!}\sum_{j=-l,j\neq0}^{l}\frac{\left(  -1\right)  ^{j}}{\left\vert j\right\vert
}\tbinom{2l}{j+l}\int\limits_{-R_{q}l}^{R_{q}l}\left(  x-s\right)
^{2n-1}q\left(  -\frac{s}{j}\right)  ds,\quad x\geq R_{q}l.
\]

Applying the change of variables $t=-\frac{s}{j}$, $ds=-jdt$ yields%
\begin{align*}
G_{q}\left(  x\right)   &  =\frac{\left(  -1\right)  ^{n}}{\left(
2n-1\right)  !}\sum_{j=-l,j\neq0}^{l}\left(  -1\right)  ^{j}\frac{\left(
-j\right)  }{\left\vert j\right\vert }\tbinom{2l}{j+l}\int_{R_{q}l/j}%
^{-R_{q}l/j}\left(  x+tj\right)  ^{2n-1}q\left(  t\right)  dt\\
&  =\frac{\left(  -1\right)  ^{n}}{\left(  2n-1\right)  !}\sum_{j=-l}%
^{-1}\left(  -1\right)  ^{j}\frac{\left(  -j\right)  }{\left\vert j\right\vert
}\tbinom{2l}{j+l}\int_{R_{q}l/j}^{-R_{q}l/j}\left(  x+tj\right)
^{2n-1}q\left(  t\right)  dt+\\
&  +\frac{\left(  -1\right)  ^{n+1}}{\left(  2n-1\right)  !}\sum_{j=1}%
^{l}\left(  -1\right)  ^{j}\frac{\left(  -j\right)  }{\left\vert j\right\vert
}\tbinom{2l}{j+l}\int_{R_{q}l/j}^{-R_{q}l/j}\left(  x+tj\right)
^{2n-1}q\left(  t\right)  dt\\
&  =\frac{\left(  -1\right)  ^{n}}{\left(  2n-1\right)  !}\sum_{j=-l}%
^{-1}\left(  -1\right)  ^{j}\tbinom{2l}{j+l}\int_{-R_{q}l/\left\vert
j\right\vert }^{R_{q}l/\left\vert j\right\vert }\left(  x+tj\right)
^{2n-1}q\left(  t\right)  dt+\\
&  +\frac{\left(  -1\right)  ^{n}}{\left(  2n-1\right)  !}\sum_{j=1}%
^{l}\left(  -1\right)  ^{j}\tbinom{2l}{j+l}\int_{-R_{q}l/\left\vert
j\right\vert }^{R_{q}l/\left\vert j\right\vert }\left(  x+tj\right)
^{2n-1}q\left(  t\right)  dt\\
&  =\frac{\left(  -1\right)  ^{n}}{\left(  2n-1\right)  !}\sum_{j=-l,j\neq
0}^{l}\left(  -1\right)  ^{j}\tbinom{2l}{j+l}\int_{-R_{q}l/\left\vert
j\right\vert }^{R_{q}l/\left\vert j\right\vert }\left(  x+tj\right)
^{2n-1}q\left(  t\right)  dt\\
&  =\frac{\left(  -1\right)  ^{n}}{\left(  2n-1\right)  !}\sum_{j=-l,j\neq
0}^{l}\left(  -1\right)  ^{j}\tbinom{2l}{j+l}\int_{-R_{q}}^{R_{q}}\left(
x+tj\right)  ^{2n-1}q\left(  t\right)  dt\\
&  =\frac{\left(  -1\right)  ^{n}}{\left(  2n-1\right)  !}\sum_{j=-l,j\neq
0}^{l}\left(  -1\right)  ^{j}\tbinom{2l}{j+l}\int\left(  \sum\limits_{m=0}%
^{2n-1}\tbinom{2n-1}{m}x^{m}\left(  tj\right)  ^{2n-1-m}\right)  q\left(
t\right)  dt\\
&  =\frac{\left(  -1\right)  ^{n}}{\left(  2n-1\right)  !}\sum\limits_{m=0}%
^{2n-1}\tbinom{2n-1}{m}x^{m}\left(  \int t^{2n-1-m}q\left(  t\right)
dt\right)  \sum_{j=-l,j\neq0}^{l}\left(  -1\right)  ^{j}\tbinom{2l}%
{j+l}j^{2n-1-m}.
\end{align*}

But by part 2b of Lemma \ref{Lem_central_diff_op},%
\[
\sum\limits_{j=-l}^{l}\left(  -1\right)  ^{j}\tbinom{2l}{j+l}j^{k}=\left\{
\begin{array}
[c]{ll}%
0, & 0\leq k<2l,\\
\left(  -1\right)  ^{l}\left(  2l\right)  !, & k=2l,
\end{array}
\right.
\]
\medskip

and hence%
\begin{equation}
\sum_{\substack{j=-l \\j\neq0}}^{l}\left(  -1\right)  ^{j}\tbinom{2l}%
{j+l}j^{k}=\left\{
\begin{array}
[c]{ll}%
-\tbinom{2l}{l}, & k=0,\\
0, & 0<k<2l,\\
\left(  -1\right)  ^{l}\left(  2l\right)  !, & k=2l,
\end{array}
\right. \label{a1.202}%
\end{equation}

so that now%
\begin{align*}
G_{q}\left(  x\right)   &  =\left(  -1\right)  ^{n+1}\frac{x^{2n-1}}{\left(
2n-1\right)  !}\left(  \int q\right)  \sum_{j=-l,j\neq0}^{l}\left(  -1\right)
^{j}\tbinom{2l}{j+l}+\\
&  \quad+\frac{\left(  -1\right)  ^{n}}{\left(  2n-1\right)  !}\sum
\limits_{\substack{m=0 \\m\text{ }even}}^{2n-2}\tbinom{2n-1}{m}x^{m}\left(
\int\limits_{-R_{q}}^{R_{q}}t^{2n-1-m}q\left(  t\right)  dt\right)
\sum_{\substack{j=-l \\j\neq0}}^{l}\left(  -1\right)  ^{j}\tbinom{2l}%
{j+l}j^{2n-1-m}\\
&  =\left(  -1\right)  ^{n+1}\tbinom{2l}{l}\left(  \int q\right)
\frac{x^{2n-1}}{\left(  2n-1\right)  !}\\
&  =-a\frac{x^{2n-1}}{\left(  2n-1\right)  !},
\end{align*}

i.e.%
\begin{equation}
G_{q}\left(  x\right)  =-P_{2n-1}\left(  x\right)  ,\quad x\geq R_{q}%
l,\label{a1.025}%
\end{equation}

using the notation of \ref{2.16}. The equations \ref{1.083} now become%
\begin{align*}
G_{q}\left(  x\right)   & =\left\{
\begin{array}
[c]{ll}%
0, & x\leq-R_{q}l,\\
\frac{\left(  -1\right)  ^{n}}{\left(  2n-1\right)  !}\sum\limits_{j=-l,j\neq
0}^{l}\frac{\left(  -1\right)  ^{j}}{\left\vert j\right\vert }\tbinom{2l}%
{j+l}\int\limits_{-R_{q}l}^{x}\left(  x-s\right)  ^{2n-1}q\left(  -\frac{s}%
{j}\right)  ds, & \left\vert x\right\vert \leq R_{q}l,\\
-P_{2n-1}\left(  x\right)  , & x\geq R_{q}l.
\end{array}
\right. \\
& =\left\{
\begin{array}
[c]{ll}%
-\left(  HP_{2n-1}\right)  \left(  x\right)  , & x\leq-R_{q}l,\\
\left.
\begin{array}
[c]{c}%
\frac{\left(  -1\right)  ^{n}}{\left(  2n-1\right)  !}\sum\limits_{j=-l,j\neq
0}^{l}\frac{\left(  -1\right)  ^{j}}{\left\vert j\right\vert }\tbinom{2l}%
{j+l}\int\limits_{-R_{q}l}^{x}\left(  x-s\right)  ^{2n-1}q\left(  -\frac{s}%
{j}\right)  ds+\\
+\left(  HP_{2n-1}\right)  \left(  x\right)  -\left(  HP_{2n-1}\right)
\left(  x\right)
\end{array}
\right\}  , & \left\vert x\right\vert \leq R_{q}l,\\
-\left(  HP_{2n-1}\right)  \left(  x\right)  , & x\geq R_{q}l.
\end{array}
\right. \\
& =-\left(  HP_{2n-1}\right)  \left(  x\right)  +\\
& +\left\{
\begin{array}
[c]{ll}%
0, & x\leq-R_{q}l,\\
\left.
\begin{array}
[c]{c}%
\frac{\left(  -1\right)  ^{n}}{\left(  2n-1\right)  !}\sum\limits_{j=-l,j\neq
0}^{l}\frac{\left(  -1\right)  ^{j}}{\left\vert j\right\vert }\tbinom{2l}%
{j+l}\int\limits_{-R_{q}l}^{x}\left(  x-s\right)  ^{2n-1}q\left(  -\frac{s}%
{j}\right)  ds+\\
+\left(  HP_{2n-1}\right)  \left(  x\right)
\end{array}
\right\}  , & \left\vert x\right\vert \leq R_{q}l,\\
0, & x\geq R_{q}l.
\end{array}
\right. \\
& =-\left(  HP_{2n-1}\right)  \left(  x\right)  +\\
& +\left\{
\begin{array}
[c]{ll}%
\left.
\begin{array}
[c]{c}%
\frac{\left(  -1\right)  ^{n}}{\left(  2n-1\right)  !}\sum\limits_{j=-l,j\neq
0}^{l}\frac{\left(  -1\right)  ^{j}}{\left\vert j\right\vert }\tbinom{2l}%
{j+l}\int\limits_{-R_{q}l}^{x}\left(  x-s\right)  ^{2n-1}q\left(  -\frac{s}%
{j}\right)  ds+\\
+\left(  HP_{2n-1}\right)  \left(  x\right)
\end{array}
\right\}  , & \left\vert x\right\vert \leq R_{q}l,\\
0, & \left\vert x\right\vert \geq R_{q}l.
\end{array}
\right.
\end{align*}

Thus from \ref{2.56} and \ref{2.16},%
\begin{align*}
G_{c}  & =G_{\delta}+G_{q}=HP_{2n-1}+G_{q}=\\
& =\left\{
\begin{array}
[c]{ll}%
\left.
\begin{array}
[c]{c}%
\frac{\left(  -1\right)  ^{n}}{\left(  2n-1\right)  !}\sum\limits_{j=-l,j\neq
0}^{l}\frac{\left(  -1\right)  ^{j}}{\left\vert j\right\vert }\tbinom{2l}%
{j+l}\int\limits_{-R_{q}l}^{x}\left(  x-s\right)  ^{2n-1}q\left(  -\frac{s}%
{j}\right)  ds+\\
+\left(  HP_{2n-1}\right)  \left(  x\right)
\end{array}
\right\}  , & \left\vert x\right\vert \leq R_{q}l,\\
0, & \left\vert x\right\vert \geq R_{q}l.
\end{array}
\right.
\end{align*}

Next we must show that $G_{c}$ has properties \ref{2.57} i.e. $G_{c}$ is
continuous and $G_{c}\left(  \pm\infty\right)  =0$.

From the last display it is clear these conditions hold if $G_{c}$ is
continuous at $\pm R_{q}l$. By inspection $G_{c}\left(  -R_{q}l\right)  =0$ so
we must show that%
\begin{equation}
\sum\limits_{j=-l,j\neq0}^{l}\frac{\left(  -1\right)  ^{j}}{\left\vert
j\right\vert }\tbinom{2l}{j+l}\int_{-R_{q}l}^{R_{q}l}\left(  R_{q}l-s\right)
^{2n-1}q\left(  -\frac{s}{j}\right)  ds=0.\label{a7.31}%
\end{equation}

Indeed,%
\begin{align*}
\sum\limits_{j=-l}^{-1} &  \frac{\left(  -1\right)  ^{j}}{\left\vert
j\right\vert }\tbinom{2l}{j+l}\int_{-R_{q}l}^{R_{q}l}\left(  R_{q}l-s\right)
^{2n-1}q\left(  -\frac{s}{j}\right)  ds\\
&  =\sum\limits_{k=1}^{l}\frac{\left(  -1\right)  ^{-k}}{\left\vert
-k\right\vert }\tbinom{2l}{-k+l}\int_{-R_{q}l}^{R_{q}l}\left(  -R_{q}%
l-s\right)  ^{2n-1}q\left(  \frac{s}{k}\right)  ds\\
&  =\sum\limits_{k=1}^{l}\frac{\left(  -1\right)  ^{k}}{\left\vert
k\right\vert }\tbinom{2l}{k+l}\left(  -1\right)  \int_{R_{q}l}^{-R_{q}%
l}\left(  -R_{q}l+t\right)  ^{2n-1}q\left(  -\frac{t}{k}\right)  dt\\
&  =\sum\limits_{k=1}^{l}\frac{\left(  -1\right)  ^{k}}{\left\vert
k\right\vert }\tbinom{2l}{k+l}\int_{-R_{q}l}^{R_{q}l}\left(  -R_{q}l+t\right)
^{2n-1}q\left(  -\frac{t}{k}\right)  ds\\
&  =-\sum\limits_{k=1}^{l}\frac{\left(  -1\right)  ^{k}}{\left\vert
k\right\vert }\tbinom{2l}{k+l}\int_{-R_{q}l}^{R_{q}l}\left(  R_{q}l-t\right)
^{2n-1}q\left(  -\frac{t}{k}\right)  ds,
\end{align*}

and so \ref{a7.31} holds.

From Theorem \ref{Thm_G_basis_def_2} $G_{c}$ is an even function so when
$-R_{q}l\leq x\leq0$,
\[
G_{c}\left(  x\right)  =\frac{\left(  -1\right)  ^{n}}{\left(  2n-1\right)
!}\sum\limits_{j=-l,j\neq0}^{l}\frac{\left(  -1\right)  ^{j}}{\left\vert
j\right\vert }\tbinom{2l}{j+l}\int_{-R_{q}l}^{x}\left(  x-s\right)
^{2n-1}q\left(  -\frac{s}{j}\right)  ds,
\]

and%
\begin{align}
G_{c}\left(  x\right)   & =\frac{\left(  -1\right)  ^{n}}{\left(  2n-1\right)
!}\sum\limits_{j=-l,j\neq0}^{l}\frac{\left(  -1\right)  ^{j}}{\left\vert
j\right\vert }\tbinom{2l}{j+l}\int_{-R_{q}l}^{-\left\vert x\right\vert
}\left(  -\left\vert x\right\vert -s\right)  ^{2n-1}q\left(  -\frac{s}%
{j}\right)  ds\nonumber\\
& =\frac{\left(  -1\right)  ^{n}}{\left(  2n-1\right)  !}\sum
\limits_{j=-l,j\neq0}^{l}\frac{\left(  -1\right)  ^{j}}{\left\vert
j\right\vert }\tbinom{2l}{j+l}\int_{\left\vert x\right\vert }^{R_{q}l}\left(
-\left\vert x\right\vert +s\right)  ^{2n-1}q\left(  \frac{s}{j}\right)
ds\nonumber\\
& =\frac{\left(  -1\right)  ^{n}}{\left(  2n-1\right)  !}\sum
\limits_{j=-l,j\neq0}^{l}\frac{\left(  -1\right)  ^{j}}{\left\vert
j\right\vert }\tbinom{2l}{j+l}\int_{\left\vert x\right\vert }^{R_{q}l}\left(
s-\left\vert x\right\vert \right)  ^{2n-1}q\left(  \frac{s}{j}\right)
ds\label{2.17}\\
& =\frac{\left(  -1\right)  ^{n}}{\left(  2n-1\right)  !}\int_{\left\vert
x\right\vert }^{R_{q}l}\left(  s-\left\vert x\right\vert \right)  ^{2n-1}%
q_{l}\left(  s\right)  ds,\nonumber
\end{align}

when $\left\vert x\right\vert \leq R_{q}l$, where%
\[
q_{l}\left(  s\right)  :=\sum\limits_{j=-l,j\neq0}^{l}\frac{\left(  -1\right)
^{j}}{\left\vert j\right\vert }\tbinom{2l}{j+l}q\left(  \frac{s}{j}\right)  .
\]

The change of variables $t=s/j$, $ds=jdt$ in \ref{2.17} yields%
\begin{align*}
G_{c}\left(  x\right)   & =\frac{\left(  -1\right)  ^{n}}{\left(  2n-1\right)
!}\sum\limits_{j=-l,j\neq0}^{l}\left(  -1\right)  ^{j}\frac{j}{\left\vert
j\right\vert }\tbinom{2l}{j+l}\int\limits_{\left\vert x\right\vert /j}%
^{R_{q}l/j}\left(  jt-\left\vert x\right\vert \right)  ^{2n-1}q\left(
t\right)  dt\\
& =\frac{\left(  -1\right)  ^{n}}{\left(  2n-1\right)  !}\sum\limits_{j=1}%
^{l}\left(  -1\right)  ^{j}\tbinom{2l}{j+l}\int\limits_{\left\vert
x\right\vert /j}^{R_{q}l/j}\left(  jt-\left\vert x\right\vert \right)
^{2n-1}q\left(  t\right)  dt+\\
& \qquad+\frac{\left(  -1\right)  ^{n}}{\left(  2n-1\right)  !}\sum
\limits_{j=-l}^{-1}\left(  -1\right)  ^{j+1}\tbinom{2l}{j+l}\int%
\limits_{\left\vert x\right\vert /j}^{R_{q}l/j}\left(  jt-\left\vert
x\right\vert \right)  ^{2n-1}q\left(  t\right)  dt.
\end{align*}

The first summation term:%
\begin{align*}
\sum\limits_{j=1}^{l}\left(  -1\right)  ^{j} &  \tbinom{2l}{j+l}%
\int\limits_{\left\vert x\right\vert /j}^{R_{q}l/j}\left(  jt-\left\vert
x\right\vert \right)  ^{2n-1}q\left(  t\right)  dt\\
&  =\sum\limits_{\substack{j=1 \\j\geq\left\vert x\right\vert /R_{q}}%
}^{l}\left(  -1\right)  ^{j}\tbinom{2l}{j+l}\int\limits_{\left\vert
x\right\vert /j}^{R_{q}}\left(  jt-\left\vert x\right\vert \right)
^{2n-1}q\left(  t\right)  dt.
\end{align*}

The second summation term:%
\begin{align*}
\sum\limits_{j=-l}^{-1}\left(  -1\right)  ^{j+1} &  \tbinom{2l}{j+l}%
\int_{\left\vert x\right\vert /j}^{R_{q}l/j}\left(  jt-\left\vert x\right\vert
\right)  ^{2n-1}q\left(  t\right)  dt\\
&  =\sum\limits_{j=1}^{l}\left(  -1\right)  ^{-j+1}\tbinom{2l}{-j+l}%
\int_{-\left\vert x\right\vert /j}^{-R_{q}l/j}\left(  -jt-\left\vert
x\right\vert \right)  ^{2n-1}q\left(  t\right)  dt\\
&  =\sum\limits_{j=1}^{l}\left(  -1\right)  ^{j+1}\tbinom{2l}{j+l}%
\int_{-\left\vert x\right\vert /j}^{-R_{q}l/j}\left(  -jt-\left\vert
x\right\vert \right)  ^{2n-1}q\left(  t\right)  dt\\
&  =\sum\limits_{j=1}^{l}\left(  -1\right)  ^{j}\tbinom{2l}{j+l}%
\int_{\left\vert x\right\vert /j}^{R_{q}l/j}\left(  jt-\left\vert x\right\vert
\right)  ^{2n-1}q\left(  -t\right)  dt\\
&  =\sum\limits_{\substack{j=1 \\j\geq\left\vert x\right\vert /R_{q}}%
}^{l}\left(  -1\right)  ^{j}\tbinom{2l}{j+l}\int_{\left\vert x\right\vert
/j}^{R_{q}}\left(  jt-\left\vert x\right\vert \right)  ^{2n-1}q\left(
-t\right)  dt.
\end{align*}

Combining the terms yields%
\[
G_{c}\left(  x\right)  =\frac{\left(  -1\right)  ^{n}2}{\left(  2n-1\right)
!}\sum\limits_{\substack{j=1 \\j\geq\left\vert x\right\vert /R_{q}}%
}^{l}\left(  -1\right)  ^{j}\tbinom{2l}{j+l}\int_{\left\vert x\right\vert
/j}^{R_{q}}\left(  jt-\left\vert x\right\vert \right)  ^{2n-1}q_{e}\left(
t\right)  dt,
\]

when $-R_{q}l\leq x\leq0$.
\end{proof}

We now need the following lemma.

\begin{lemma}
\label{Lem_DerivIntegAtoX(X-S)f(S)dS}Suppose $x^{j}f\in L^{1}\left(
0,a\right)  $ for $j=0,1,2,\ldots,m$.

Then for $x\in\left[  0,a\right]  $:

\begin{enumerate}
\item
\[
D_{x}^{j}\int_{x}^{a}\frac{\left(  s-x\right)  ^{m}}{m!}f\left(  s\right)
ds=\left\{
\begin{array}
[c]{ll}%
\left(  -1\right)  ^{j}\int_{x}^{a}\frac{\left(  s-x\right)  ^{m-j}}{\left(
m-j\right)  !}f\left(  s\right)  ds, & j\leq m,\\
\left(  -1\right)  ^{m+1}f\left(  x\right)  , & j=m+1,
\end{array}
\right\}  ;
\]

\item $\int_{x}^{a}\left(  s-x\right)  ^{m}f\left(  s\right)  ds\in
C_{B}^{\left(  m\right)  }\left[  0,a\right]  $;

\item $\left\vert \int_{x}^{a}\left(  s-x\right)  ^{k}f\left(  s\right)
ds\right\vert \leq\int_{0}^{a}s^{k}\left\vert f\left(  s\right)  \right\vert
ds$.

\item If $g\in C_{B}^{\left(  2m\right)  }\left[  0,a\right]  $ and
$D^{k}g\left(  0^{+}\right)  =0$ for $k\leq2m-1$ then $g\left(  \left\vert
\cdot\right\vert \right)  \in C_{B}^{\left(  2m-1\right)  }\left[
-a,a\right]  \cap C_{B}^{\left(  2m\right)  }\left(  \left[  -a,a\right]
\setminus0\right)  $ and%
\[
D^{k}\left(  g\left(  \left\vert x\right\vert \right)  \right)  =\left\{
\begin{array}
[c]{ll}%
\left(  \operatorname*{sgn}x\right)  ^{k}\left(  D^{k}g\right)  \left(
\left\vert x\right\vert \right)  , & k\leq2m-1,\\
\left(  \operatorname*{sgn}x\right)  ^{k}\left(  D^{2m}g\right)  \left(
\left\vert x\right\vert \right)  +2\left(  D^{2m}g\right)  \left(
0^{+}\right)  \delta, & k=2m.
\end{array}
\right.
\]

This is a result from the theory of piecewise continuous distributions (e.g.
Vladimirov 2.6.3 \cite{Vladimirov}).

\item By the mean-value theorem:

$\left\vert D^{k}\left(  g\left(  \left\vert x\right\vert \right)  \right)
-D^{k}\left(  g\left(  \left\vert y\right\vert \right)  \right)  \right\vert
\leq\left\vert \left\vert x\right\vert -\left\vert y\right\vert \right\vert
\left\Vert D^{k+1}\left(  g\left(  \left\vert \cdot\right\vert \right)
\right)  \right\Vert _{\left[  \left\vert x\right\vert ,\left\vert
y\right\vert \right]  }$ when $k\leq2m-2$ and $\left\vert x\right\vert
,\left\vert y\right\vert \leq a$.
\end{enumerate}
\end{lemma}

Two formulas for the basis function were given in Theorem
\ref{Thm_cdiffbasis_part_moment_formula}. We give the corresponding formulas
for the derivatives.

\begin{corollary}
\label{Cor_Thm_cdiffbasis_part_moment_formula_2}Suppose $w_{c}$ is a central
difference weight function with parameters $n,l,q$ and $w_{c}\in W02/W03$, and
denote the basis function by $G_{c}$. Suppose $\operatorname*{supp}%
G_{c}\subset\overline{B}_{R_{q}}$. Then for $k\leq2n-1$,%
\begin{equation}
D^{k}G_{c}\left(  x\right)  =\left\{
\begin{array}
[c]{ll}%
\frac{\left(  -1\right)  ^{n+k}}{\left(  2n-k-1\right)  !}\left(
\operatorname*{sgn}x\right)  ^{k}\int_{\left\vert x\right\vert }^{R_{q}%
l}\left(  s-\left\vert x\right\vert \right)  ^{2n-k-1}q_{l}\left(  s\right)
ds, & \left\vert x\right\vert \leq R_{q}l,\\
0, & \left\vert x\right\vert \geq R_{q}l,
\end{array}
\right. \label{2.67}%
\end{equation}

where $q_{l}$ is given by \ref{2.60}, and%
\begin{equation}
D^{k}G_{c}\left(  x\right)  =\left\{
\begin{array}
[c]{ll}%
\frac{\left(  -1\right)  ^{n+k}2}{\left(  2n-k-1\right)  !}\left(
\operatorname*{sgn}x\right)  ^{k}\sum\limits_{\substack{j=1 \\j\geq\left\vert
x\right\vert /R_{q}}}^{l}\left(  -1\right)  ^{j}\tbinom{2l}{j+l}%
\int\limits_{\left\vert x\right\vert /j}^{R_{q}}\left(  jt-\left\vert
x\right\vert \right)  ^{2n-k-1}q_{e}\left(  t\right)  dt, & \left\vert
x\right\vert \leq R_{q}l,\\
0, & \left\vert x\right\vert \geq R_{q}l,
\end{array}
\right. \label{2.69}%
\end{equation}

where $q_{e}\left(  t\right)  :=\frac{1}{2}\left(  q\left(  t\right)
+q\left(  -t\right)  \right)  $ is the even component of $q$.

We also have $G_{c}\in C_{B}^{2n-2}\left(  \mathbb{R}^{1}\right)  \cap
C^{2n-1}\left(  \mathbb{R}^{1}\setminus0\right)  $ with the bounds%
\[
\left\Vert D^{k}G_{c}\right\Vert _{\infty}\leq\frac{1}{\left(  2n-k-1\right)
!}\left(  \sum\limits_{j=1}^{l}\tbinom{2l}{j+l}j^{2n-k-1}\right)  \int
t^{2n-k-1}q\left(  t\right)  dt,\quad k\leq2n-1.
\]

Finally, at zero $D^{2n-1}G_{c}$ has a non-zero jump (or step) of size
$\left(  -1\right)  ^{n}\tbinom{2l}{l}\int q$.
\end{corollary}

\begin{proof}
To obtain these derivative formulas apply Lemma
\ref{Lem_DerivIntegAtoX(X-S)f(S)dS} to equations \ref{2.59} and \ref{2.58}.
Clearly%
\begin{align*}
\left\vert D^{k}G_{c}\left(  x\right)  \right\vert  & \leq\frac{2}{\left(
2n-k-1\right)  !}\sum\limits_{\substack{j=1 \\j\geq\left\vert x\right\vert
/R_{q}}}^{l}\tbinom{2l}{j+l}\int_{\left\vert x\right\vert /j}^{R_{q}}\left(
jt-\left\vert x\right\vert \right)  ^{2n-k-1}q_{e}\left(  t\right)  dt\\
& \leq\frac{2}{\left(  2n-k-1\right)  !}\sum\limits_{\substack{j=1
\\j\geq\left\vert x\right\vert /R_{q}}}^{l}\tbinom{2l}{j+l}\int_{\left\vert
x\right\vert /j}^{R_{q}}\left(  jt\right)  ^{2n-k-1}q_{e}\left(  t\right)
dt\\
& \leq\frac{2}{\left(  2n-k-1\right)  !}\sum\limits_{j=1}^{l}\tbinom{2l}%
{j+l}\int_{0}^{R_{q}}\left(  jt\right)  ^{2n-k-1}q_{e}\left(  t\right)  dt\\
& =\frac{2}{\left(  2n-k-1\right)  !}\left(  \sum\limits_{j=1}^{l}\tbinom
{2l}{j+l}j^{2n-k-1}\right)  \int_{0}^{R_{q}}t^{2n-k-1}q_{e}\left(  t\right)
dt\\
& =\frac{1}{\left(  2n-k-1\right)  !}\left(  \sum\limits_{j=1}^{l}\tbinom
{2l}{j+l}j^{2n-k-1}\right)  \int t^{2n-k-1}q\left(  t\right)  dt.
\end{align*}

Regarding the jump at the origin: when $k=2n-1$ \ref{2.69} becomes%
\[
D^{2n-1}G_{c}\left(  x\right)  =\left(  -1\right)  ^{n+1}\left(
\operatorname*{sgn}x\right)  ^{2n+1}2\sum\limits_{\substack{j=1 \\j\geq
\left\vert x\right\vert /R_{q}}}^{l}\left(  -1\right)  ^{j}\tbinom{2l}%
{j+l}\int\limits_{\left\vert x\right\vert /j}^{R_{q}}q_{e}\left(  t\right)
dt,
\]

so that%
\begin{align*}
D^{2n-1}G_{c}\left(  0^{+}\right)  -D^{2n-1}G_{c}\left(  0^{-}\right)   &
=\left(  -1\right)  ^{n+1}4\left(  \sum\limits_{j=1}^{l}\left(  -1\right)
^{j}\tbinom{2l}{j+l}\right)  \int q_{e}\left(  t\right)  dt\\
& =\left(  -1\right)  ^{n+1}2\left(  \sum\limits_{j=1}^{l}\left(  -1\right)
^{j}\tbinom{2l}{j+l}\right)  \int q\\
& :\ref{a022}\Rightarrow\\
& =\left(  -1\right)  ^{n+1}2\left(  -\frac{1}{2}\tbinom{2l}{l}\right)  \int
q\\
& =\left(  -1\right)  ^{n}\tbinom{2l}{l}\int q.
\end{align*}

\end{proof}

The last corollary can be modified to define the central difference basis
function \textbf{using only the parameters which define the weight function -
and NOT the weight function}:

\begin{remark}
\label{Rem_cdiffbasis_formulas_3}??? \textbf{CHECK}! Suppose $1\leq n\leq l$
are integers, $q\in L^{1}\left(  \mathbb{R}^{1}\right)  $, $q\geq0$ and
$\int_{\left\vert \xi\right\vert \geq R}\left\vert \xi\right\vert
^{2n-1}q\left(  \xi\right)  d\xi<\infty$, for some $R\geq0$. From Theorem
\ref{Thm_cdiffwt_2} $n,l,q$ define a central difference weight function.
Denote the basis function by $G_{c}$.

Then the formulas and upper bounds for the basis function and its derivatives
derived in Corollary \ref{Cor_Thm_cdiffbasis_part_moment_formula_2} hold.
\end{remark}

This remark has the direct consequence:

\begin{corollary}
\label{Cor_Thm_cdiffbasis_part_moment_formula_1}If $q\in C_{0}^{\infty}\left(
\mathbb{R}^{1}\right)  $ and $q\geq0$ then by Theorem \ref{Thm_cdiffwt_2},
$w_{c}=\frac{\xi^{2n}}{\Delta_{2l}\widehat{q}\left(  \xi\right)  }$ is a
central difference weight function for all parameters $1\leq n\leq l$. Denote
the basis function by $G_{c}$.

Then $G_{c}\in C^{\infty}\left(  \mathbb{R}^{1}\setminus0\right)  \cap
C_{B}^{\left(  2n-2\right)  }\left(  \mathbb{R}^{1}\right)  $, $D^{2n-1}%
G_{c}\in L^{\infty}$ and $G_{c}$ has bounded support. Indeed, if
$\operatorname*{supp}q\subseteq\overline{B}_{R_{q}}$ then
$\operatorname*{supp}G_{c}\subseteq\overline{B}_{R_{q}l}$.

The formulas and upper bounds for the basis function and its derivatives
derived in Corollary \ref{Cor_Thm_cdiffbasis_part_moment_formula_2} hold.
\end{corollary}

\begin{remark}
??? Towards creating $C_{0}^{\infty}$ basis functions: We use \ref{2.59}. Set
$n=l$ and choose $q\in C_{0}^{\infty}$ with $\operatorname*{supp}%
q\subseteq\overline{B}_{R_{q}}$. Define $q^{\left(  n\right)  }\left(
s\right)  =q\left(  ns\right)  $. Then $q^{\left(  n\right)  }$ has the same
weight function properties as $q$ and $\operatorname*{supp}q^{\left(
n\right)  }\subseteq\overline{B}_{R_{q}/n}$. Denote the corresponding basis
functions by $G_{c}^{\left(  n\right)  }$. From Corollary
\ref{Cor_Thm_cdiffbasis_part_moment_formula_1} $G_{c}^{\left(  n\right)  }\in
C^{\infty}\left(  \mathbb{R}^{1}\setminus0\right)  \cap C_{B}^{\left(
2n-2\right)  }\left(  \mathbb{R}^{1}\right)  $ and $\operatorname*{supp}%
G_{c}^{\left(  n\right)  }\subseteq\overline{B}_{R_{q}}$.

What happens as $n\rightarrow\infty$? Does the sequence $G_{c}^{\left(
n\right)  }$ converge uniformly to a non-zero $C_{0}^{\infty}$ function with
non-negative Fourier transform? Numerical experiment!

It now occurs to me that the derivatives at the origin may be the key.
Calculate the odd derivatives and choose a sequence of $q$s so that they
converge as $n\rightarrow\infty$. The even derivatives are zero.
\end{remark}

\begin{remark}
??? From \ref{2.59},%
\begin{align*}
G_{c}\left(  x\right)   & =\left\{
\begin{array}
[c]{ll}%
\frac{\left(  -1\right)  ^{n}}{\left(  2n-1\right)  !}\int H\left(
s-\left\vert x\right\vert \right)  \left(  s-\left\vert x\right\vert \right)
^{2n-1}H\left(  s\right)  q_{l}\left(  s\right)  ds, & \left\vert x\right\vert
\leq R_{q}l,\\
0, & \left\vert x\right\vert \geq R_{q}l.
\end{array}
\right. \\
& =\left\{
\begin{array}
[c]{ll}%
\left(  2\pi\right)  ^{d/2}\frac{\left(  -1\right)  ^{n}}{\left(  2n-1\right)
!}\left(  Hs^{2n-1}\ast Hq_{l}\right)  \left(  \left\vert x\right\vert
\right)  , & \left\vert x\right\vert \leq R_{q}l,\\
0, & \left\vert x\right\vert \geq R_{q}l.
\end{array}
\right.
\end{align*}

Since $Hs^{2n-1}$ and $Hq_{l}$ are locally integrable and their supports are
contained in $\left[  0,\infty\right)  $ it follows that the convolution is
locally integrable thus that $G_{c}$ is locally integrable (Vladimirov 2.7.4
\cite{Vladimirov}).
\end{remark}

\begin{example}
\label{Ex_calc_centbasis_from_qrect2}Suppose $q$ is as in Example
\ref{Ex_central_diff_basis_hat_squ} so that $n=l=1$. Here we will use Theorem
\ref{Thm_cdiffbasis_part_moment_formula} to calculate the corresponding basis
function. Here $R_{q}=1$ and so from equation \ref{2.58},
\begin{align*}
G_{c}\left(  x\right)   & =\left\{
\begin{array}
[c]{ll}%
\frac{\left(  -1\right)  ^{n}2}{\left(  2n-1\right)  !}\sum
\limits_{\substack{j=1 \\j\geq\left\vert x\right\vert /R_{q}}}^{l}\left(
-1\right)  ^{j}\tbinom{2l}{j+l}\int_{\left\vert x\right\vert /j}^{R_{q}%
}\left(  jt-\left\vert x\right\vert \right)  ^{2n-1}q_{e}\left(  t\right)
dt, & \left\vert x\right\vert \leq R_{q}l,\\
0, & \left\vert x\right\vert \geq R_{q}l,
\end{array}
\right. \\
& =\left\{
\begin{array}
[c]{ll}%
-2\sum\limits_{\substack{j=1 \\j\geq\left\vert x\right\vert }}^{1}\left(
-1\right)  ^{j}\tbinom{2}{j+1}\int_{\left\vert x\right\vert /j}^{1}\left(
jt-\left\vert x\right\vert \right)  dt, & \left\vert x\right\vert \leq1,\\
0, & \left\vert x\right\vert \geq1,
\end{array}
\right. \\
& =\left\{
\begin{array}
[c]{ll}%
-2\sum\limits_{j=1}\left(  -1\right)  ^{j}\tbinom{2}{j+1}\int_{\left\vert
x\right\vert /j}^{1}\left(  jt-\left\vert x\right\vert \right)  dt, &
\left\vert x\right\vert \leq1,\\
0, & \left\vert x\right\vert \geq1,
\end{array}
\right. \\
& =\left\{
\begin{array}
[c]{ll}%
2\int_{\left\vert x\right\vert }^{1}\left(  t-\left\vert x\right\vert \right)
dt, & \left\vert x\right\vert \leq1,\\
0, & \left\vert x\right\vert \geq1,
\end{array}
\right\}  =\left\{
\begin{array}
[c]{ll}%
2\left[  \frac{1}{2}t^{2}-\left\vert x\right\vert t\right]  _{\left\vert
x\right\vert }^{1}, & \left\vert x\right\vert \leq1,\\
0, & \left\vert x\right\vert \geq1,
\end{array}
\right\}  =\\
& =\left\{
\begin{array}
[c]{ll}%
2\left(  \frac{1}{2}-\left\vert x\right\vert \right)  -2\left(  \frac{1}%
{2}\left\vert x\right\vert ^{2}-\left\vert x\right\vert ^{2}\right)  , &
\left\vert x\right\vert \leq1,\\
0, & \left\vert x\right\vert \geq1,
\end{array}
\right. \\
& =\left\{
\begin{array}
[c]{ll}%
\left(  1-2\left\vert x\right\vert \right)  -\left(  \left\vert x\right\vert
^{2}-2\left\vert x\right\vert ^{2}\right)  , & \left\vert x\right\vert
\leq1,\\
0, & \left\vert x\right\vert \geq1,
\end{array}
\right. \\
& =\left\{
\begin{array}
[c]{ll}%
\left(  1-\left\vert x\right\vert \right)  ^{2}, & \left\vert x\right\vert
\leq1,\\
0, & \left\vert x\right\vert \geq1,
\end{array}
\right.  =\Lambda\left(  x\right)  ^{2},
\end{align*}

which confirms the result of Example \ref{Ex_central_diff_basis_hat_squ}.
\end{example}

\subsection{An alternative proof of the partial moment formula using a
tempered distribution Taylor series expansion and the subspace $S_{\emptyset
,1}\subset S$\label{SbSect_tempdistrib_1_dim_centdiff_basis}}

We use the tempered distribution Taylor series expansion introduced in Section
\ref{Sect_Taylor_series_data_fn} and the theory of the spaces $S_{\emptyset
,k}\subset S$ introduced in Definition \ref{Def_So,n} to confirm the partial
moment formula \ref{2.59} for the central difference basis function. In the
process we obtain another multiplicative convolution formula which is given in
Theorem \ref{Thm_cdiffbasis_mult_convol_ql}.

Recall that this basis function $G_{c}$ is generated by the parameters $n$,
$l$ and a non-negative $L^{1}$ function $q\neq0$. From Theorem
\ref{Thm_cdiffwt_2} we must have $1\leq n\leq l$ and $\int x^{k}q<0$ for
$k=0,1,\ldots,2n-1$.

Suppose $q\in L^{1}\left(  \mathbb{R}^{1}\right)  $ so that $\widehat{q}\in
C_{B}^{\left(  0\right)  }\left(  \mathbb{R}^{1}\right)  $. Then by
definition:%
\[
\Delta_{2l}\widehat{q}\left(  \xi\right)  =\sum_{j=-l}^{l}\left(  -1\right)
^{j}\tbinom{2l}{j+l}\widehat{q}\left(  -j\xi\right)  .
\]

The tempered distribution formula \ref{a1.55} states that if $f\in S^{\prime}$
then%
\[
f\left(  \cdot+\xi\right)  -\sum_{k\leq m}\frac{\xi^{k}}{\beta!}D^{k}%
f=\frac{\sqrt{2\pi}}{m!}\left(  \left(  i\left(  \xi,\cdot\right)  \right)
^{m+1}\overline{\widehat{g_{m}}}\left(  \left(  \xi,\cdot\right)  \right)
\widehat{f}\right)  ^{\vee},\quad m=0,1,2,\ldots,
\]

where the function $g_{n}$ is given by \ref{a117}. When $f$ is replaced by
$\widehat{f}$,%
\[
\widehat{f}\left(  \cdot+\xi\right)  -\sum_{k\leq m}\frac{\xi^{k}}{k!}%
D^{k}\widehat{f}=\frac{\sqrt{2\pi}}{m!}\left(  \left(  i\left(  \xi
,\cdot\right)  \right)  ^{m+1}\overline{\widehat{g_{m}}}\left(  \left(
\xi,\cdot\right)  \right)  f_{\_}\right)  ^{\vee},
\]

where $f_{\_}$ is the distribution extension of $f_{\_}\left(  x\right)
=f\left(  -x\right)  $.

Now, in particular, when $x^{k}q\in L^{1}$ for $k\leq m$ then by Lemma
\ref{Lem_L1_Fourier_contin}, $\widehat{q}\in C_{B}^{\left(  m\right)  }\left(
\mathbb{R}^{1}\right)  $ and
\[
\widehat{q}\left(  \cdot+\xi\right)  -\sum_{k\leq m}\frac{\xi^{k}}{k!}%
D^{k}\widehat{q}=\frac{\sqrt{2\pi}}{m!}F_{\eta}^{-1}\left[  \left(  i\xi
\eta\right)  ^{m+1}\overline{\widehat{g_{m}}}\left(  \xi\eta\right)  q\left(
-\eta\right)  \right]  ,
\]

and consequently%
\[
\widehat{q}\left(  \xi\right)  -\sum_{k\leq m}\frac{\xi^{k}}{k!}%
D^{k}\widehat{q}\left(  0\right)  =\frac{\sqrt{2\pi}}{m!}F_{\eta}^{-1}\left[
\left(  i\xi\eta\right)  ^{m+1}\overline{\widehat{g_{m}}}\left(  \xi
\eta\right)  q\left(  -\eta\right)  \right]  \left(  0\right)  .
\]

From part 4 of Lemma \ref{Lem_gm_properties}, for each $\xi$,
\begin{align*}
\left\vert \left(  i\xi\eta\right)  ^{m+1}\overline{\widehat{g_{m}}}\left(
\xi\eta\right)  q\left(  -\eta\right)  \right\vert  & =\left\vert \xi
\eta\right\vert ^{m+1}\left\vert \widehat{g_{m}}\left(  \xi\eta\right)
\right\vert q\left(  -\eta\right) \\
& \leq\left\vert \xi\eta\right\vert ^{m+1}\frac{c_{m}}{\sqrt{2\pi}}\frac
{1}{1+\left\vert \xi\eta\right\vert }q\left(  -\eta\right) \\
& =\frac{c_{m}}{\sqrt{2\pi}}\left\vert \xi\right\vert ^{m}\left\vert
\eta\right\vert ^{m}q\left(  -\eta\right) \\
& \in L^{1},
\end{align*}

and so%
\[
\widehat{q}\left(  \xi\right)  -\sum_{k\leq m}\frac{\xi^{k}}{k!}%
D^{k}\widehat{q}\left(  0\right)  =\frac{1}{m!}\int\left(  i\xi\eta\right)
^{m+1}\overline{\widehat{g_{m}}}\left(  \xi\eta\right)  q\left(  -\eta\right)
d\eta.
\]

The next step is to use the operator $\Delta_{2l}$. Choose $m=2n-1$ so that
$m<2l$ and%
\begin{align}
\Delta_{2l}\widehat{q}\left(  \xi\right)   & =\Delta_{2l}\left(
\widehat{q}\left(  \xi\right)  -\sum_{k\leq m}\frac{\xi^{k}}{k!}%
D^{k}\widehat{q}\left(  0\right)  \right) \nonumber\\
& =\frac{1}{\left(  2n-1\right)  !}\Delta_{2l,\xi}\int\left(  i\xi\eta\right)
^{2n}\overline{\widehat{g_{2n-1}}}\left(  \xi\eta\right)  q\left(
-\eta\right)  d\eta\nonumber\\
& =\frac{1}{\left(  2n-1\right)  !}\sum_{j=-l,\text{ }j\neq0}^{l}\left(
-1\right)  ^{j}\tbinom{2l}{j+l}\int\left(  -ij\xi\eta\right)  ^{2n}%
\overline{\widehat{g_{2n-1}}}\left(  -j\xi\eta\right)  q\left(  -\eta\right)
d\eta\nonumber\\
& =\frac{\left(  -1\right)  ^{n}\xi^{2n}}{\left(  2n-1\right)  !}%
\sum_{j=-l,\text{ }j\neq0}^{l}\left(  -1\right)  ^{j}\tbinom{2l}{j+l}%
\int\left(  j\eta\right)  ^{2n}\overline{\widehat{g_{2n-1}}}\left(  -j\xi
\eta\right)  q\left(  -\eta\right)  d\eta\nonumber\\
& =\frac{\left(  -1\right)  ^{n}\xi^{2n}}{\left(  2n-1\right)  !}%
\sum_{j=-l,\text{ }j\neq0}^{l}\left(  -1\right)  ^{j}\tbinom{2l}{j+l}%
\int\left(  j\eta\right)  ^{2n}\overline{\widehat{g_{2n-1}}}\left(  j\xi
\eta\right)  q\left(  \eta\right)  d\eta\nonumber\\
& =\frac{\left(  -1\right)  ^{n}\xi^{2n}}{\left(  2n-1\right)  !}%
\sum_{j=-l,\text{ }j\neq0}^{l}\frac{\left(  -1\right)  ^{j}}{\left\vert
j\right\vert }\tbinom{2l}{j+l}\int\tau^{2n}\overline{\widehat{g_{2n-1}}%
}\left(  \xi\tau\right)  q\left(  \tfrac{1}{j}\tau\right)  d\tau\nonumber\\
& =\frac{\left(  -1\right)  ^{n}\xi^{2n}}{\left(  2n-1\right)  !}\int\tau
^{2n}\overline{\widehat{g_{2n-1}}}\left(  \xi\tau\right)  \sum_{j=-l,\text{
}j\neq0}^{l}\frac{\left(  -1\right)  ^{j}}{\left\vert j\right\vert }%
\tbinom{2l}{j+l}q\left(  \tfrac{1}{j}\tau\right)  d\tau\nonumber\\
& =\frac{\left(  -1\right)  ^{n}\xi^{2n}}{\left(  2n-1\right)  !}\int\tau
^{2n}\overline{\widehat{g_{2n-1}}}\left(  \xi\tau\right)  q_{l}\left(
\tau\right)  d\tau,\label{a976}%
\end{align}

where the even function%
\[
q_{l}\left(  \tau\right)  =\sum_{j=-l,\text{ }j\neq0}^{l}\frac{\left(
-1\right)  ^{j}}{\left\vert j\right\vert }\tbinom{2l}{j+l}q\left(  \tfrac
{1}{j}\tau\right)  ,
\]

was introduced in Theorem \ref{Thm_cdiffbasis_part_moment_formula}. Thus%
\[
\widehat{G_{c}}\left(  \xi\right)  =\frac{\Delta_{2l}\widehat{q}\left(
\xi\right)  }{\xi^{2n}}=\frac{\left(  -1\right)  ^{n}}{\left(  2n-1\right)
!}\int\tau^{2n}\overline{\widehat{g_{2n-1}}}\left(  \xi\tau\right)
q_{l}\left(  \tau\right)  d\tau,
\]

and so for $\phi\in S$,%
\[
\left[  \widehat{G_{c}}\left(  \xi\right)  ,\phi\right]  =\frac{\left(
-1\right)  ^{n}}{\left(  2n-1\right)  !}\int\int\tau^{2n}\overline
{\widehat{g_{2n-1}}}\left(  \xi\tau\right)  q_{l}\left(  \tau\right)
d\tau\text{ }d\xi.
\]

I now want to change the order of integration using Fubini's theorem and to do
this I need to show that the integrand is $L^{1}\left(  \mathbb{R}^{2}\right)
$. Write%
\begin{align*}
\int\int &  \left\vert \overline{\widehat{g_{2n-1}}}\left(  \xi\tau\right)
\tau^{2n}q_{l}\left(  \tau\right)  \phi\left(  \xi\right)  \right\vert d\xi
d\tau\\
&  =\int\int\left\vert \tau\right\vert \left\vert \widehat{g_{2n-1}}\left(
\xi\tau\right)  \right\vert \left\vert \phi\left(  \xi\right)  \right\vert
\left\vert \tau^{2n-1}\right\vert q_{l}\left(  \tau\right)  d\xi d\tau\\
&  =\int\int\frac{\left\vert \xi\tau\widehat{g_{2n-1}}\left(  \xi\tau\right)
\right\vert }{\left\vert \xi\right\vert }\left\vert \phi\left(  \xi\right)
\right\vert d\xi\text{ }\left\vert \tau^{2n-1}\right\vert q_{l}\left(
\tau\right)  d\tau\\
&  =\int\int\left\vert \xi\tau\widehat{g_{2n-1}}\left(  \xi\tau\right)
\right\vert \frac{\left\vert \phi\left(  \xi\right)  \right\vert }{\left\vert
\xi\right\vert }d\xi\text{ }\left\vert \tau^{2n-1}\right\vert q_{l}\left(
\tau\right)  d\tau.
\end{align*}

From \ref{a2.07}, $\left\vert \xi\tau\widehat{g_{2n-1}}\left(  \xi\tau\right)
\right\vert \leq\frac{3}{\sqrt{2\pi}}$ so that%
\begin{align*}
\int\int &  \left\vert \overline{\widehat{g_{2n-1}}}\left(  \xi\tau\right)
\tau^{2n}q_{l}\left(  \tau\right)  \phi\left(  \xi\right)  \right\vert d\xi
d\tau\\
&  \leq\frac{3}{\sqrt{2\pi}}\int\int\frac{\left\vert \phi\left(  \xi\right)
\right\vert }{\left\vert \xi\right\vert }d\xi\text{ }\left\vert \tau
^{2n-1}\right\vert q_{l}\left(  \tau\right)  d\tau\\
&  =\frac{3}{\sqrt{2\pi}}\left(  \int\frac{\left\vert \phi\left(  \xi\right)
\right\vert }{\left\vert \xi\right\vert }d\xi\right)  \int\left\vert
\tau^{2n-1}\right\vert q_{l}\left(  \tau\right)  d\tau\\
&  =\frac{3}{\sqrt{2\pi}}\left(  \int_{\left\vert \xi\right\vert \leq1}%
\frac{\left\vert \phi\left(  \xi\right)  \right\vert }{\left\vert
\xi\right\vert }d\xi+\int_{\left\vert \xi\right\vert \geq1}\frac{\left\vert
\phi\left(  \xi\right)  \right\vert }{\left\vert \xi\right\vert }d\xi\right)
\int\left\vert \tau^{2n-1}\right\vert q_{l}\left(  \tau\right)  d\tau.
\end{align*}

From this inequality it is clear that we should restrict $\phi$ to the
subspace $S_{\emptyset,1}=\left\{  \phi\in S:\phi\left(  0\right)  =1\right\}
$ which belongs to the class of subspaces $S_{\emptyset,k}$ of $S $ introduced
in Definition \ref{Def_So,n}. We now have%
\begin{align*}
\int\int\left\vert \overline{\widehat{g_{2n-1}}}\left(  \xi\tau\right)
\tau^{2n}q_{l}\left(  \tau\right)  \phi\left(  \xi\right)  \right\vert d\xi
d\tau & \leq\frac{3}{\sqrt{2\pi}}\left(  \left\Vert D\phi\right\Vert _{\infty
}+\left\Vert \frac{\phi}{\left\vert \cdot\right\vert }\right\Vert _{\infty
}\right)  \int\left\vert \tau^{2n-1}\right\vert q_{l}\left(  \tau\right)
d\tau\\
& <\infty.
\end{align*}

It follows that the integrand is absolutely convergent and Fubini's theorem
allows the order of integration to be reversed. Thus for $\phi\in
S_{\emptyset,1}$,%
\begin{align*}
\left[  \widehat{G_{c}},\phi\right]   & =\frac{\left(  -1\right)  ^{n}%
}{\left(  2n-1\right)  !}\int\left(  \tau\int\overline{\widehat{g_{2n-1}}%
}\left(  \xi\tau\right)  \phi\left(  \xi\right)  d\xi\right)  \tau^{2n-1}%
q_{l}\left(  \tau\right)  d\tau\\
& =\frac{\left(  -1\right)  ^{n}}{\left(  2n-1\right)  !}\int\left(  \tau
\int\overline{\widehat{\frac{1}{\left\vert \tau\right\vert }g_{2n-1}\left(
\frac{x}{\tau}\right)  }}\left(  \xi\right)  \phi\left(  \xi\right)
d\xi\right)  \tau^{2n-1}q_{l}\left(  \tau\right)  d\tau\\
& =\frac{\left(  -1\right)  ^{n}}{\left(  2n-1\right)  !}\int\left(
\int\overline{\widehat{g_{2n-1}\left(  \frac{x}{\tau}\right)  }}\left(
\xi\right)  \phi\left(  \xi\right)  d\xi\right)  \left\vert \tau\right\vert
^{2n-1}q_{l}\left(  \tau\right)  d\tau\\
& :\text{Plancherel's theorem}\Longrightarrow\\
& =\frac{\left(  -1\right)  ^{n}}{\left(  2n-1\right)  !}\int\left(  \int
g_{2n-1}\left(  \frac{\zeta}{\tau}\right)  \widehat{\phi}\left(  \zeta\right)
d\zeta\right)  \left\vert \tau\right\vert ^{2n-1}q_{l}\left(  \tau\right)
d\tau.
\end{align*}

Since $g_{2n-1}$ is bounded the integrand is clearly absolutely convergent and
so%
\[
\left[  \widehat{G_{c}},\phi\right]  =\frac{\left(  -1\right)  ^{n}}{\left(
2n-1\right)  !}\int\left(  \int g_{2n-1}\left(  \frac{\zeta}{\tau}\right)
\left\vert \tau\right\vert ^{2n-1}q_{l}\left(  \tau\right)  d\tau\right)
\widehat{\phi}\left(  \zeta\right)  d\zeta.
\]

Since $\int g_{2n-1}\left(  \frac{\zeta}{\tau}\right)  \left\vert
\tau\right\vert ^{2n-1}q_{l}\left(  \tau\right)  d\tau$ is a bounded function
it lies in $S^{\prime}$ and so

$\widehat{G_{c}}=\frac{\left(  -1\right)  ^{n}}{\left(  2n-1\right)  !}\left(
\int g_{2n-1}\left(  \frac{\zeta}{\tau}\right)  \left\vert \tau\right\vert
^{2n-1}q_{l}\left(  \tau\right)  d\tau\right)  ^{\wedge}$ on $S_{\emptyset,1}%
$. Theorem \ref{Thm_prop_functnl_on_Son} now implies that

$\widehat{G_{c}}-\frac{\left(  -1\right)  ^{n}}{\left(  2n-1\right)  !}\left(
\int g_{2n-1}\left(  \frac{\zeta}{\tau}\right)  \left\vert \tau\right\vert
^{2n-1}q_{l}\left(  \tau\right)  d\tau\right)  ^{\wedge}$ is a delta function
which in turn implies that%
\[
G_{c}\left(  \zeta\right)  =\frac{\left(  -1\right)  ^{n}}{\left(
2n-1\right)  !}\int g_{2n-1}\left(  \frac{\zeta}{\tau}\right)  \left\vert
\tau\right\vert ^{2n-1}q_{l}\left(  \tau\right)  d\tau+const.
\]

If $\zeta>0$ then from \ref{a117},%
\begin{align*}
\int g_{2n-1}\left(  \frac{\zeta}{\tau}\right)  \left\vert \tau\right\vert
^{2n-1}q_{l}\left(  \tau\right)  d\tau & =\int_{\zeta}^{\infty}g_{2n-1}\left(
\frac{\zeta}{\tau}\right)  \left\vert \tau\right\vert ^{2n-1}q_{l}\left(
\tau\right)  d\tau\\
& =\int_{\zeta}^{\infty}\left(  1-\frac{\zeta}{\tau}\right)  ^{2n-1}\left\vert
\tau\right\vert ^{2n-1}q_{l}\left(  \tau\right)  d\tau\\
& =\int_{\zeta}^{\infty}\left(  \tau-\zeta\right)  ^{2n-1}q_{l}\left(
\tau\right)  d\tau.
\end{align*}

Similarly, if $\zeta<0$, then%
\begin{align*}
\int g_{2n-1}\left(  \frac{\zeta}{\tau}\right)  \left\vert \tau\right\vert
^{2n-1}q_{l}\left(  \tau\right)  d\tau & =\int_{-\infty}^{\zeta}%
g_{2n-1}\left(  \frac{\zeta}{\tau}\right)  \left\vert \tau\right\vert
^{2n-1}q_{l}\left(  \tau\right)  d\tau\\
& =\int_{-\infty}^{\zeta}\left(  1-\frac{\zeta}{\tau}\right)  ^{2n-1}%
\left\vert \tau\right\vert ^{2n-1}q_{l}\left(  \tau\right)  d\tau\\
& =\int_{-\zeta}^{\infty}\left(  1+\frac{\zeta}{\tau}\right)  ^{2n-1}%
\left\vert \tau\right\vert ^{2n-1}q_{l}\left(  \tau\right)  d\tau\\
& =\int_{-\zeta}^{\infty}\left(  \tau+\zeta\right)  ^{2n-1}q_{l}\left(
\tau\right)  d\tau,
\end{align*}

and so for all $\zeta$,%
\[
\int g_{2n-1}\left(  \frac{\zeta}{\tau}\right)  \left\vert \tau\right\vert
^{2n-1}q_{l}\left(  \tau\right)  d\tau=\int_{\left\vert \zeta\right\vert
}^{\infty}\left(  \tau-\left\vert \zeta\right\vert \right)  ^{2n-1}%
q_{l}\left(  \tau\right)  d\tau,
\]

i.e.%
\[
G_{c}\left(  \zeta\right)  =\frac{\left(  -1\right)  ^{n}}{\left(
2n-1\right)  !}\int_{\left\vert \zeta\right\vert }^{\infty}\left(
\tau-\left\vert \zeta\right\vert \right)  ^{2n-1}q_{l}\left(  \tau\right)
d\tau+const.
\]

Finally, the observation that%
\[
\left\vert G_{c}\left(  \zeta\right)  \right\vert \leq\frac{1}{\left(
2n-1\right)  !}\int_{\left\vert \zeta\right\vert }^{\infty}\tau^{2n-1}%
q_{l}\left(  \tau\right)  d\tau,
\]

clearly implies that $G_{c}\left(  \infty\right)  =0$ and hence that%
\[
G_{c}\left(  \zeta\right)  =\frac{\left(  -1\right)  ^{n}}{\left(
2n-1\right)  !}\int_{\left\vert \zeta\right\vert }^{\infty}\left(
\tau-\left\vert \zeta\right\vert \right)  ^{2n-1}q_{l}\left(  \tau\right)
d\tau,
\]

and%
\begin{equation}
G_{c}\left(  \zeta\right)  =\frac{\left(  -1\right)  ^{n}}{\left(
2n-1\right)  !}\int g_{2n-1}\left(  \frac{\left\vert \zeta\right\vert }{\tau
}\right)  \left\vert \tau\right\vert ^{2n-1}q_{l}\left(  \tau\right)
d\tau.\label{1.058}%
\end{equation}

The first formula matches the partial moment formula \ref{2.59} of Theorem
\ref{Thm_cdiffbasis_part_moment_formula} and the second formula we present as:

\begin{theorem}
\label{Thm_cdiffbasis_mult_convol_ql}\textbf{Multiplicative convolution
formula} Suppose the conditions of Theorem \ref{Thm_cdiffwt_2} are satisfied.
Then the central difference basis function generated by $n,l,q$ is%
\begin{align*}
G_{c}\left(  \zeta\right)   & =\frac{\left(  -1\right)  ^{n}}{\left(
2n-1\right)  !}\int g_{2n-1}\left(  \frac{\left\vert \zeta\right\vert }{\tau
}\right)  \left\vert \tau\right\vert ^{2n-1}q_{l}\left(  \tau\right)  d\tau\\
& =\frac{\left(  -1\right)  ^{n}}{\left(  2n-1\right)  !}\int_{\left\vert
\zeta\right\vert }^{\infty}\left(  \tau-\left\vert \zeta\right\vert \right)
^{2n-1}q_{l}\left(  \tau\right)  d\tau,
\end{align*}

where $q_{l}$ is given by \ref{2.60} and $g_{2n-1}$ by \ref{a117}.
\end{theorem}

\begin{proof}
Given before this theorem.
\end{proof}

\begin{remark}
Compare this result with the multiplicative convolution formula of Theorem
\ref{Thm_G_basis_def_2}. This formula does not involve a formula for the B-spline.
\end{remark}

\subsection{Lipschitz continuity estimates from partial moment
formulas\label{SbSect_cdiffbasis_Lip_part_mom}}

Some Lipschitz continuity estimates were given in Subsection
\ref{SbSect_CentDiffBasis_Lipschitz} and were derived using a multiplicative
convolution formula for the basis function. Here we derive Lipschitz
continuity estimates base on partial moment formulas.

??? \textbf{WARNING} - the estimate of the next theorem has a factor
$\left\Vert t^{2n-k-1}q\right\Vert _{1}$ but the upper bound of Theorem
\ref{Thm_cdiffbasis_Lips_dim_d} has a factor of $\left\Vert t^{2n-k-2}%
q\right\Vert _{1}$ - CHECK THIS!

\begin{theorem}
\label{Thm_cdiffbasis_Lip_part_moment}\textbf{Lipschitz continuity} Suppose
the conditions of Corollary \ref{Cor_Thm_cdiffbasis_part_moment_formula_2} are
satisfied. Then the central difference basis function $G_{c}$ generated by
parameters $n,l,q$ is Lipschitz continuous and%
\[
\left\vert D^{k}G_{c}\left(  x\right)  -D^{k}G_{c}\left(  y\right)
\right\vert \leq\frac{1}{\left(  2n-k-1\right)  !}\left(  \sum\limits_{j=1}%
^{l}\tbinom{2l}{j+l}j^{2n-k-1}\right)  \left\Vert t^{2n-k-1}q\right\Vert
_{1}\left\vert x-y\right\vert ,
\]

when $k\leq2n-2$.
\end{theorem}

\begin{proof}
By the mean-value theorem and the upper bound of Theorem
\ref{Cor_Thm_cdiffbasis_part_moment_formula_2},%
\begin{align*}
\left\vert D^{k}G_{c}\left(  x\right)  -D^{k}G_{c}\left(  y\right)
\right\vert  & \leq\left\vert x-y\right\vert \left\Vert D^{k+1}G_{c}%
\right\Vert _{\infty}\\
& \leq\frac{\left\vert x-y\right\vert }{\left(  2n-k-1\right)  !}\left(
\sum\limits_{j=1}^{l}\tbinom{2l}{j+l}j^{2n-k-1}\right)  \int t^{2n-k-1}%
q\left(  t\right)  dt\\
& =\frac{1}{\left(  2n-k-1\right)  !}\left(  \sum\limits_{j=1}^{l}\tbinom
{2l}{j+l}j^{2n-k-1}\right)  \left\Vert t^{2n-k-1}q\right\Vert _{1}\left\vert
x-y\right\vert .
\end{align*}

\end{proof}

\section{Weight and basis function convergence
results\label{Sect_cent_diff_basis_conv_nat_splin_basis}}

In the preamble to the definition of a central difference weight function in
Subsection \ref{SbSect_cent_Motivation} several assertions were made regarding
the convergence of sequences of central difference weight and basis functions
to extended B-spline weight and basis functions. This section is devoted to
proving these assertions.

In fact, given an extended B-spline weight function $w_{s}$ and corresponding
basis function $G_{s}$, there exists a sequence of central difference weight
functions $w_{k}$ and corresponding basis functions $G_{k}$ such that
$1/w_{k}\rightarrow1/w_{s}$ in the $L^{1}$ sense and $G_{k}\rightarrow G_{s}$
uniformly pointwise.

\subsection{Uniform pointwise basis function convergence}

The next lemma shows that the tensor product extended B-spline basis functions
are Lipschitz continuous of order 1.

\begin{lemma}
\label{Lem_ex_nat_spline_basis_Lipschitz}(Theorem
\ref{Thm_ex_nat_spline_basis_Lipschitz} Chapter \ref{Ch_wtfn_basisfn_datasp})
Let $G_{s}\left(  x\right)  =\prod\limits_{k=1}^{d}G_{1}\left(  x_{k}\right)
$ be a tensor product \textbf{extended B-spline} basis function, as described
in Theorem \ref{Thm_basis_tensor_hat_W3}. We then have the estimate
\[
\left\vert G_{s}\left(  x\right)  -G_{s}\left(  y\right)  \right\vert
\leq\sqrt{d}G_{1}\left(  0\right)  ^{d-1}\left\Vert DG_{1}\right\Vert
_{\infty}\left\vert x-y\right\vert ,\text{\quad}x,y\in\mathbb{R}^{d}.
\]

\end{lemma}

Using the multiplicative convolution formulas we now show that for any
extended B-spline basis function $G_{s}$ there exists a sequence of central
difference basis functions which converges pointwise to $G_{s}$.

\begin{theorem}
\label{Thm_centr_diff_basis_conv_to_ext_spl_basis}Suppose $l,n\geq1$ are
integers which satisfy $1\leq n\leq l$. Let $\psi\in L^{1}\left(
\mathbb{R}^{1}\right)  $ be a function which satisfies $\psi\left(  x\right)
\geq0$, $\int\psi=1$, $\operatorname*{supp}\psi\subset\left[  -1,1\right]  $.
Thus%
\begin{equation}
\int_{\mathbb{R}^{1}}\left\vert t\right\vert ^{2n}\psi\left(  t\right)
dt<\infty.\label{a910}%
\end{equation}

Then the sequence of functions $q_{k}\left(  t\right)  =\sqrt{2\pi}%
2^{-2l}\frac{k}{2}\psi\left(  \frac{k}{2}\left(  t-2\right)  \right)  ,\quad
k=1,2,3,\ldots$ generates \textbf{a sequence of central difference basis
functions}, say $G_{k}$, which \textbf{converges uniformly pointwise} to the
extended B-spline basis function $G_{s}$ with parameters $l$ and $n$. In fact
\[
\left\vert G_{k}\left(  x\right)  -G_{s}\left(  x\right)  \right\vert
\leq\frac{2^{2n-2}}{k}\left(  \left\Vert DG_{s}\right\Vert _{\infty}\left\vert
x\right\vert +2\left\Vert G_{s}\right\Vert _{\infty}\right)  ,\quad
x\in\left(  1+\frac{1}{k}\right)  \operatorname*{supp}G_{s}.
\]

\textbf{Note}: the use of $k/2$ in the formula for $q_{k}$ instead of $k$,
which was used in the formula \ref{a978}, simplifies the algebra of the proof.
\end{theorem}

\begin{proof}
It is clear that functions $q_{k}$ satisfy the conditions of Theorem
\ref{Thm_cdiffwt_2} so that each function $\frac{\xi^{2n}}{\Delta
^{2l}\widehat{q_{k}}\left(  \xi\right)  }$ is a central difference weight
function with parameters $l$ and $n$. By Theorem \ref{Thm_G_basis_def_2}, the
central difference basis function is%
\[
G_{1}\left(  s\right)  =\tfrac{2^{2\left(  l-n\right)  +1}}{\sqrt{2\pi}}%
\int_{\mathbb{R}^{1}}G_{s}\left(  \frac{2s}{t}\right)  \left\vert t\right\vert
^{2n-1}q\left(  t\right)  dt,
\]

and so the corresponding sequence of basis functions is
\begin{align}
G_{k}\left(  x\right)   & =\frac{1}{2^{2n-1}}\int G_{s}\left(  \frac{2x}%
{t}\right)  \left\vert t\right\vert ^{2n-1}\frac{k}{2}\psi\left(  \frac{k}%
{2}\left(  t-2\right)  \right)  dt\nonumber\\
& :t=2s,\text{ }dt=2ds\Rightarrow\nonumber\\
& =\int G_{s}\left(  \frac{x}{s}\right)  \left\vert s\right\vert ^{2n-1}%
k\psi\left(  k\left(  s-1\right)  \right)  ds,\label{a2.36}%
\end{align}

with $\operatorname*{supp}G_{k}\subset\left(  1+\frac{1}{k}\right)
\operatorname*{supp}G_{s}=\left(  1+\frac{1}{k}\right)  2l$.

Since $\int\psi=1$ it follows that%
\begin{align*}
G_{k}\left(  x\right)  -G_{s}\left(  x\right)   & =\int\left(  G_{s}\left(
\frac{x}{t}\right)  t^{2n-1}-G_{s}\left(  x\right)  \right)  k\psi\left(
k\left(  t-1\right)  \right)  dt\\
& =\int_{1-\frac{1}{k}}^{1+\frac{1}{k}}\left(  G_{s}\left(  \frac{x}%
{t}\right)  t^{2n-1}-G_{s}\left(  x\right)  t^{2n-1}\right)  k\psi\left(
k\left(  t-1\right)  \right)  dt+\\
& \qquad+\int_{1-\frac{1}{k}}^{1+\frac{1}{k}}\left(  G_{s}\left(  x\right)
t^{2n-1}-G_{s}\left(  x\right)  \right)  k\psi\left(  k\left(  t-1\right)
\right)  dt\\
& =\int_{1-\frac{1}{k}}^{1+\frac{1}{k}}\left(  G_{s}\left(  \frac{x}%
{t}\right)  -G_{s}\left(  x\right)  \right)  t^{2n-1}k\psi\left(  k\left(
t-1\right)  \right)  dt+\\
& \qquad+\int_{1-\frac{1}{k}}^{1+\frac{1}{k}}G_{s}\left(  x\right)  \left(
t^{2n-1}-1\right)  k\psi\left(  k\left(  t-1\right)  \right)  dt.
\end{align*}

The Lipschitz continuity estimate on $\mathbb{R}^{1}$ derived for the extended
B-spline in Lemma \ref{Lem_ex_nat_spline_basis_Lipschitz} is%
\[
\left\vert G_{s}\left(  x\right)  -G_{s}\left(  y\right)  \right\vert
\leq\left\Vert DG_{s}\right\Vert _{\infty}\left\vert x-y\right\vert
,\text{\quad}x,y\in\mathbb{R}^{d}.
\]

Hence%
\begin{align*}
&  \left\vert G_{k}\left(  x\right)  -G_{s}\left(  x\right)  \right\vert \\
&  \leq\int_{1-\frac{1}{k}}^{1+\frac{1}{k}}\left\vert G_{s}\left(  \frac{x}%
{t}\right)  -G_{s}\left(  x\right)  \right\vert t^{2n-1}k\psi\left(  k\left(
t-1\right)  \right)  dt+\\
&  \qquad+\int_{1-\frac{1}{k}}^{1+\frac{1}{k}}\left\vert G_{s}\left(
x\right)  \right\vert \left(  t^{2n-1}-1\right)  k\psi\left(  k\left(
t-1\right)  \right)  dt\\
&  \leq\left\Vert DG_{s}\right\Vert _{\infty}\int_{1-\frac{1}{k}}^{1+\frac
{1}{k}}\left\vert \frac{x}{t}-x\right\vert t^{2n-1}k\psi\left(  k\left(
t-1\right)  \right)  dt+\\
&  \qquad+\left\Vert G_{s}\right\Vert _{\infty}\int_{1-\frac{1}{k}}%
^{1+\frac{1}{k}}\left(  t^{2n-1}-1\right)  k\psi\left(  k\left(  t-1\right)
\right)  dt\\
&  =\left\Vert DG_{s}\right\Vert _{\infty}\left\vert x\right\vert
\int_{1-\frac{1}{k}}^{1+\frac{1}{k}}\left(  t-1\right)  t^{2n-2}k\psi\left(
k\left(  t-1\right)  \right)  dt+\\
&  \qquad+\left\Vert G_{s}\right\Vert _{\infty}\int_{1-\frac{1}{k}}%
^{1+\frac{1}{k}}\left(  t-1\right)  \left(  1+t+\ldots+t^{2n-2}\right)
k\psi\left(  k\left(  t-1\right)  \right)  dt\\
&  =\left\Vert DG_{s}\right\Vert _{\infty}\left\vert x\right\vert \int%
_{-\frac{1}{k}}^{\frac{1}{k}}s\left(  1+s\right)  ^{2n-2}k\psi\left(
ks\right)  ds+\\
&  \qquad+\left\Vert G_{s}\right\Vert _{\infty}\int_{-\frac{1}{k}}^{\frac
{1}{k}}s\left(  1+\left(  1+s\right)  +\ldots+\left(  1+s\right)
^{2n-2}\right)  k\psi\left(  ks\right)  ds\\
&  \leq\left\Vert DG_{s}\right\Vert _{\infty}\left\vert x\right\vert \frac
{1}{k}\left(  1+\frac{1}{k}\right)  ^{2n-2}\int_{-\frac{1}{k}}^{\frac{1}{k}%
}k\psi\left(  ks\right)  ds+\\
&  \qquad+\left\Vert G_{s}\right\Vert _{\infty}\frac{1}{k}\left(  1+\left(
1+\frac{1}{k}\right)  +\ldots+\left(  1+\frac{1}{k}\right)  ^{2n-2}\right)
\int_{-\frac{1}{k}}^{\frac{1}{k}}k\psi\left(  ks\right)  ds\\
&  \leq\left\Vert DG_{s}\right\Vert _{\infty}\frac{2^{2n-2}\left\vert
x\right\vert }{k}+\left\Vert G_{s}\right\Vert _{\infty}\frac{1}{k}\left(
1+2+\ldots+2^{2n-2}\right) \\
&  =\frac{2^{2n-2}}{k}\left(  \left\Vert DG_{s}\right\Vert _{\infty}\left\vert
x\right\vert +2\left\Vert G_{s}\right\Vert _{\infty}\right)  ,
\end{align*}

and since $G_{s}$ has compact support we have uniform convergence.
\end{proof}

\begin{theorem}
\label{Thm_cdiffbasis_deriv_convg_to_Bspl_basis}Suppose $l,n\geq1$ are
integers which satisfy $1\leq n\leq l$. Let $\psi\in L^{1}\left(
\mathbb{R}^{1}\right)  $ be a function which satisfies $\psi\left(  x\right)
\geq0$, $\int\psi=1$, $\operatorname*{supp}\psi\subset\left[  -1,1\right]  $.
Thus%
\[
\int_{\mathbb{R}^{1}}\left\vert t\right\vert ^{2n}\psi\left(  t\right)
dt<\infty.
\]

Then the sequence of functions $q_{k}\left(  t\right)  =\sqrt{2\pi}%
2^{-2l}\frac{k}{2}\psi\left(  \frac{k}{2}\left(  t-2\right)  \right)  ,\quad
k=1,2,3,\ldots$ generates \textbf{a sequence of central difference basis
functions}, say $G_{k}$, which \textbf{converges uniformly pointwise} to the
extended B-spline basis function $G_{s}$ with parameters $l$ and $n$. In fact,
for $j\leq2n-2$,
\[
\left\vert D^{j}G_{k}\left(  x\right)  -D^{j}G_{s}\left(  x\right)
\right\vert \leq\frac{2^{2n-2-j}}{k}\left(  \left\Vert D^{j+1}G_{s}\right\Vert
_{\infty}\left\vert x\right\vert +2\left\Vert D^{j}G_{s}\right\Vert _{\infty
}\right)  ,\quad x\in\left(  1+\frac{1}{k}\right)  \operatorname*{supp}G_{s}.
\]

\end{theorem}

\begin{proof}
Start with the derivative formula \ref{a940} of Theorem
\ref{Thm_G_basis_def_2} and adapt the proof of the last theorem.
\end{proof}

\subsection{$L^{1}$ inverse weight function convergence}

In this subsection we prove the inverse weight function convergence theorem
\ref{Thm_CntDiffWtFnConvergToBSplin} i.e. given an extended B-spline weight
function $w_{s}$ there exists a sequence of central difference weight
functions $w_{k}$ such that $1/w_{k}\rightarrow1/w_{s}$ in the $L^{1}$ sense.

Recall that at the start of Subsection \ref{SbSect_cent_Motivation} it was
suggested that if $\psi\in L^{1}\left(  \mathbb{R}^{1}\right)  $ satisfies
$\psi\left(  x\right)  \geq0$ and $\int\psi=1$ and an extra condition then the
sequence of functions $q_{k}\left(  t\right)  =\frac{\sqrt{2\pi}}{2^{2l}}%
\frac{k}{2}\psi\left(  \frac{k}{2}\left(  t-2\right)  \right)  \in L^{1}$ is
such that $\frac{\Delta^{2l}\widehat{q_{k}}\left(  \xi\right)  }{\xi^{2n}}$
converges to $\frac{\sin^{2l}\xi}{\xi^{2n}}$ in the sense of tempered
distributions. Here $\frac{\sin^{2l}\xi}{\xi^{2n}}=\frac{1}{w_{s}\left(
\xi\right)  }\in L^{1}$ where $w_{s}$ is the extended B-spline weight function
with parameters $n,l$.\ The extra condition referred to is \ref{a910} i.e.
$\int_{\left\vert t\right\vert \geq R}\left\vert t\right\vert ^{2n}%
q_{k}\left(  t\right)  dt<\infty$ for some $R\geq0$, which is equivalent to
$\int\left\vert t\right\vert ^{2n}\psi\left(  t\right)  dt<\infty$ since
$\int\psi<\infty$.\ This assumption is stronger than \ref{a963} which was the
used in Theorem \ref{Thm_cdiffwt_2} to ensure that a 1-dimensional central
difference weight function has property W02%
$\backslash$%
W03. We will prove $L^{1}$ convergence and this implies convergence in
$S^{\prime}$. Assume $\xi>0$. Then using the change of variables: $s=\xi t/2$,
$dt=\frac{2}{\xi}ds$%
\begin{align*}
\frac{1}{w_{k}\left(  \xi\right)  } &  =\tfrac{2^{2\left(  l-n\right)  }%
}{\sqrt{2\pi}}\int\frac{t^{2n}q_{k}\left(  t\right)  }{w_{s}\left(  \xi
t/2\right)  }dt\\
&  =\int\frac{\left(  t/2\right)  ^{2n}\frac{k}{2}\psi\left(  \frac{k}%
{2}\left(  t-2\right)  \right)  }{w_{s}\left(  \xi t/2\right)  }dt\\
&  =\int\left(  \frac{s}{\xi}\right)  ^{2n}\frac{\frac{k}{\xi}\psi\left(
\frac{k}{\xi}\left(  s-\xi\right)  \right)  }{w_{s}\left(  s\right)  }ds\\
&  =\frac{1}{\xi^{2n}}\int\frac{k}{\xi}\psi\left(  \frac{k}{\xi}\left(
s-\xi\right)  \right)  \sin^{2l}s\text{ }ds\\
&  =\frac{1}{\xi^{2n}}\int\frac{k}{\xi}\psi\left(  \frac{k}{\xi}\left(
s-\xi\right)  \right)  \left(  \sin^{2l}s-\sin^{2l}\xi\right)  ds+\\
&  \qquad\qquad+\frac{1}{\xi^{2n}}\int\frac{k}{\xi}\psi\left(  \frac{k}{\xi
}\left(  s-\xi\right)  \right)  \sin^{2l}\xi\text{ }ds\\
&  =\frac{1}{\xi^{2n}}\int\frac{k}{\xi}\psi\left(  \frac{k}{\xi}\left(
s-\xi\right)  \right)  \left(  \sin^{2l}s-\sin^{2l}\xi\right)  ds+\frac
{\sin^{2l}\xi}{\xi^{2n}}\int\psi\\
&  =\frac{1}{\xi^{2n}}\int\frac{k}{\xi}\psi\left(  \frac{k}{\xi}\left(
s-\xi\right)  \right)  \left(  \sin^{2l}s-\sin^{2l}\xi\right)  ds+\frac
{\sin^{2l}\xi}{\xi^{2n}},
\end{align*}

so that%
\[
\frac{1}{w_{k}\left(  \xi\right)  }-\frac{1}{w_{s}\left(  \xi\right)  }%
=\frac{1}{\xi^{2n}}\int\frac{k}{\xi}\psi\left(  \frac{k}{\xi}t\right)  \left(
\sin^{2l}\left(  t+\xi\right)  -\sin^{2l}\xi\right)  dt.
\]

Now split the domain of integration about $\left\vert \xi\right\vert =R$:
since $\int\psi=1$,%
\begin{align*}
\int\limits_{\left\vert \cdot\right\vert \geq R}\left\vert \frac{1}{w_{k}%
}-\frac{1}{w_{s}}\right\vert  & \leq\int\limits_{\left\vert \xi\right\vert
\geq R}\frac{1}{\xi^{2n}}\int\frac{k}{\xi}\psi\left(  \frac{k}{\xi}t\right)
\left\vert \sin^{2l}\left(  t+\xi\right)  -\sin^{2l}\xi\right\vert dtd\xi\\
& \leq2\int\limits_{\left\vert \xi\right\vert \geq R}\frac{1}{\xi^{2n}}%
\int\frac{k}{\xi}\psi\left(  \frac{k}{\xi}t\right)  dtd\xi\\
& =2\int\limits_{\left\vert \xi\right\vert \geq R}\frac{d\xi}{\xi^{2n}}\\
& =4\int\limits_{\xi\geq R}\frac{d\xi}{\xi^{2n}}\\
& =\frac{4}{2n-1}\frac{1}{R^{2n-1}}.
\end{align*}

If $G_{s}$ is the B-spline basis function we can write%
\begin{align}
\frac{1}{w_{k}\left(  \xi\right)  }-\frac{1}{w_{s}\left(  \xi\right)  }  &
=\frac{1}{\xi^{2n}}\int\frac{k}{\xi}\psi\left(  \frac{k}{\xi}\left(
s-\xi\right)  \right)  \left(  \sin^{2l}s-\sin^{2l}\xi\right)  ds\nonumber\\
& =\frac{1}{\xi^{2n}}\int\frac{k}{\xi}\psi\left(  \frac{k}{\xi}t\right)
\left(  \sin^{2l}\left(  \xi+t\right)  -\sin^{2l}\xi\right)  dt\nonumber\\
& =\int\frac{k}{\xi}\psi\left(  \frac{k}{\xi}t\right)  \left(  \frac{\sin
^{2l}\left(  \xi+t\right)  }{\xi^{2n}}-\frac{\sin^{2l}\xi}{\xi^{2n}}\right)
dt\nonumber\\
& =\int\frac{k}{\xi}\psi\left(  \frac{k}{\xi}t\right)  \left(  \frac{\left(
t+\xi\right)  ^{2n}}{\xi^{2n}}\widehat{G_{s}}\left(  \xi+t\right)
-\widehat{G_{s}}\left(  \xi\right)  \right)  dt\nonumber\\
& =\int\frac{k}{\xi}\psi\left(  \frac{k}{\xi}t\right)  \left(  \left(
1+\frac{t}{\xi}\right)  ^{2n}\widehat{G_{s}}\left(  \xi+t\right)
-\widehat{G_{s}}\left(  \xi\right)  \right)  dt\nonumber\\
& =\int\frac{k}{\xi}\psi\left(  \frac{k}{\xi}t\right)  \left(  \sum
\limits_{m\leq2n}\binom{2n}{k}\left(  \frac{t}{\xi}\right)  ^{m}%
\widehat{G_{s}}\left(  \xi+t\right)  -\widehat{G_{s}}\left(  \xi\right)
\right)  dt\nonumber\\
& =\sum\limits_{2\leq m\leq2n}\binom{2n}{m}\int\frac{k}{\xi}\left(  \frac
{t}{\xi}\right)  ^{m}\psi\left(  \frac{k}{\xi}t\right)  \widehat{G_{s}}\left(
\xi+t\right)  dt+\nonumber\\
& \qquad+\int\frac{k}{\xi}\psi\left(  \frac{k}{\xi}t\right)  \left(
\widehat{G_{s}}\left(  \xi+t\right)  -\widehat{G_{s}}\left(  \xi\right)
\right)  dt\nonumber\\
& =\sum\limits_{2\leq m\leq2n}\binom{2n}{m}\int\frac{k}{\xi}\left(  \frac
{t}{\xi}\right)  ^{m}\psi\left(  \frac{k}{\xi}t\right)  \widehat{G_{s}}\left(
\xi+t\right)  dt+\nonumber\\
& \qquad+\int\frac{k}{\xi}\psi\left(  \frac{k}{\xi}t\right)  \left(
\widehat{G_{s}}\left(  \xi+t\right)  -\widehat{G_{s}}\left(  \xi\right)
\right)  dt.\label{a913}%
\end{align}

\begin{align*}
\int\left\vert \frac{1}{w_{k}}-\frac{1}{w_{s}}\right\vert \leq\sum
\limits_{2\leq m\leq2n} &  \binom{2n}{m}\int\limits_{\left\vert \xi\right\vert
\leq R}\int\left\vert \frac{k}{\xi}\left(  \frac{t}{\xi}\right)
^{m}\right\vert \psi\left(  \frac{k}{\xi}t\right)  \widehat{G_{s}}\left(
\xi+t\right)  dtd\xi+\\
&  +\int\limits_{\left\vert \xi\right\vert \leq R}\int\frac{k}{\left\vert
\xi\right\vert }\psi\left(  \frac{k}{\xi}t\right)  \left\vert \widehat{G_{s}%
}\left(  \xi+t\right)  -\widehat{G_{s}}\left(  \xi\right)  \right\vert
dtd\xi+\frac{4}{2n-1}\frac{1}{R^{2n-1}}.
\end{align*}

We will consider the last integral in \ref{a913} first. Since%
\begin{align*}
\left\vert \widehat{G_{s}}\left(  \xi+t\right)  -\widehat{G_{s}}\left(
\xi\right)  \right\vert =\left\vert \widehat{e^{i\left(  t,\cdot\right)
}G_{s}}\left(  \xi\right)  -\widehat{G_{s}}\left(  \xi\right)  \right\vert  &
=\left\vert \left(  \left(  e^{i\left(  t,\cdot\right)  }-1\right)
G_{s}\right)  ^{\wedge}\left(  \xi\right)  \right\vert \\
&  =\left\vert t\right\vert \text{ }\left\vert \left(  \left(  \frac
{e^{i\left(  t,\cdot\right)  }-1}{t}\right)  G_{s}\right)  ^{\wedge}\left(
\xi\right)  \right\vert \\
&  \leq\left\vert t\right\vert \text{ }\left\Vert xG_{s}\left(  x\right)
\right\Vert _{1},
\end{align*}

it follows that%
\begin{align*}
-\int\frac{k}{\left\vert \xi\right\vert }\psi\left(  \frac{k}{\xi}t\right)
\left\vert \widehat{G_{s}}\left(  \xi+t\right)  -\widehat{G_{s}}\left(
\xi\right)  \right\vert dt  & \leq\left\Vert xG_{s}\left(  x\right)
\right\Vert _{1}\int\frac{k}{\left\vert \xi\right\vert }\left\vert
t\right\vert \psi\left(  \frac{k}{\xi}t\right)  dt\\
& =\frac{\left\vert \xi\right\vert }{k}\left\Vert xG_{s}\left(  x\right)
\right\Vert _{1}\left\Vert s\psi\left(  s\right)  \right\Vert _{1},
\end{align*}

which exists since $\psi\in L^{1}$ and $\int_{\left\vert s\right\vert \geq
R}\left\vert s\right\vert ^{2n}\psi\left(  s\right)  ds<\infty$. Thus%
\begin{align}
\int\limits_{\left\vert \xi\right\vert \leq R}\int\frac{k}{\left\vert
\xi\right\vert }\psi\left(  \frac{k}{\xi}t\right)  \left\vert \widehat{G_{s}%
}\left(  \xi+t\right)  -\widehat{G_{s}}\left(  \xi\right)  \right\vert dtd\xi
& \leq\int\limits_{\left\vert \xi\right\vert \leq R}\frac{\left\vert
\xi\right\vert }{k}\left\Vert xG_{s}\left(  x\right)  \right\Vert
_{1}\left\Vert s\psi\left(  s\right)  \right\Vert _{1}d\xi\nonumber\\
& =\frac{R^{2}}{k}\left\Vert xG_{s}\left(  x\right)  \right\Vert
_{1}\left\Vert s\psi\left(  s\right)  \right\Vert _{1}.\label{a916}%
\end{align}

Regarding the integrals under the summation sign of \ref{a913}, for
$m=2,3,\ldots,2n$ we apply the change of variable $s=\frac{k}{\xi}t$ to get%
\begin{align*}
\left\vert \int_{0}^{R}\int\frac{k}{\xi}\left(  \frac{t}{\xi}\right)  ^{m}%
\psi\left(  \frac{k}{\xi}t\right)  \widehat{G_{s}}\left(  \xi+t\right)
dtd\xi\right\vert  & <\int_{0}^{R}\int\frac{k}{\xi}\left(  \frac{\left\vert
t\right\vert }{\xi}\right)  ^{m}\psi\left(  \frac{k}{\xi}t\right)  dtd\xi\\
& =\frac{1}{k^{m}}\int_{0}^{R}\int\left\vert s\right\vert ^{m}\psi\left(
s\right)  dsd\xi\\
& =\frac{R}{k^{m}}\left\Vert s^{m}\psi\left(  s\right)  \right\Vert _{1},
\end{align*}

and similarly%
\begin{align*}
\left\vert \int_{-R}^{0}\int\frac{k}{\xi}\left(  \frac{t}{\xi}\right)
^{m}\psi\left(  \frac{k}{\xi}t\right)  \widehat{G_{s}}\left(  \xi+t\right)
dtd\xi\right\vert  & <\frac{R}{k^{m}}\int\left\vert s\right\vert ^{m}%
\psi\left(  s\right)  ds\\
& =\frac{R}{k^{m}}\left\Vert s^{m}\psi\left(  s\right)  \right\Vert _{1},
\end{align*}

so that%
\begin{align}
\left\vert \int_{\left\vert \xi\right\vert <R}\int\frac{k}{\xi}\left(
\frac{t}{\xi}\right)  ^{m}\psi\left(  \frac{k}{\xi}t\right)  \widehat{G_{s}%
}\left(  \xi+t\right)  dtd\xi\right\vert  & <\frac{2R}{k^{m}}\int\left\vert
s\right\vert ^{m}\psi\left(  s\right)  ds\nonumber\\
& =\frac{2R}{k^{m}}\left\Vert s^{m}\psi\left(  s\right)  \right\Vert
_{1}.\label{a917}%
\end{align}

Inequalities \ref{a916} and \ref{a917} can now be used to estimate the right
side of \ref{a913}: if $R>2$ then%
\begin{align}
\int\left\vert \frac{1}{w_{k}}-\frac{1}{w_{s}}\right\vert  & <\sum
\limits_{2\leq m\leq2n}\binom{2n}{m}\frac{2R}{k^{m}}\left\Vert s^{m}%
\psi\left(  s\right)  \right\Vert _{1}+\frac{R^{2}}{k}\left\Vert xG_{s}\left(
x\right)  \right\Vert _{1}\left\Vert s\psi\left(  s\right)  \right\Vert
_{1}+\nonumber\\
& \qquad+\frac{4}{2n-1}\frac{1}{R^{2n-1}}\nonumber\\
& =\frac{R^{2}}{k}\left(  \sum\limits_{2\leq m\leq2n}\binom{2n}{m}+\left\Vert
xG_{s}\left(  x\right)  \right\Vert _{1}\right)  \max_{m=1}^{2n}\left\{
\left\Vert s^{m}\psi\left(  s\right)  \right\Vert _{1}\right\}  +\frac
{4}{R^{2n-1}}\nonumber\\
& =\frac{R^{2}}{k}\left(  2^{2n}-1+\left\Vert xG_{s}\left(  x\right)
\right\Vert _{1}\right)  \max_{m=1}^{2n}\left\{  \left\Vert s^{m}\psi\left(
s\right)  \right\Vert _{1}\right\}  +\frac{4}{R^{2n-1}},\label{a918}%
\end{align}

and it is clear that given $\varepsilon>0$ the right side can be made less
than $\varepsilon$ by first choosing $R$ such that $\frac{4}{R^{2n-1}}%
<\frac{\varepsilon}{2}$ and then choosing $k$ so that the remaining terms also
do not exceed $\varepsilon/2$. Thus we have shown that $\frac{1}{w_{k}%
}\rightarrow\frac{1}{w_{s}}$ in $L^{1}$ and hence in $S^{\prime}$.

Thus we have proved:

\begin{theorem}
\label{Thm_CntDiffWtFnConvergToBSplin}Suppose $w_{s}$ is an extended B-spline
weight function. Suppose $q_{k}$ is the sequence of $q$ functions constructed
in Theorem \ref{Thm_centr_diff_basis_conv_to_ext_spl_basis} and $w_{k}$ is the
corresponding sequence of central difference weight functions. Then $\frac
{1}{w_{k}}\rightarrow\frac{1}{w_{s}}$ in $L^{1}$.
\end{theorem}

\begin{proof}
Precedes this theorem.
\end{proof}

\subsection{Approximating a sequence of B-spline basis functions which
converges to a $C_{0}^{\infty}$ basis function}

A $C_{0}^{\infty}$ basis function can be approximated by an infinite sequence
of B-spline basis functions scaled to have support $\left[  -2n,2n\right]  $
and to take a value of $1$ at the origin - this is the sequence $\mathcal{G}%
_{s}\left(  x;n\right)  =\frac{G_{s}\left(  2nx;n\right)  }{G_{s}\left(
0;n\right)  }$ below. In this section I will construct a sequence of central
difference basis functions which converge to this sequence.

\begin{conjecture}
\label{Conj_CoinfBasis_MatExperim}From \textbf{Matlab experiments} we
\textbf{conjecture} that%
\begin{align}
\int_{0}^{\infty}\frac{\sin^{2n}t}{t^{2n}}dt  & \leq\frac{1}{\sqrt{n}}\int%
_{0}^{\infty}\frac{\sin^{2}t}{t^{2}}dt=\frac{1}{\sqrt{n}}\frac{\pi}%
{2},\text{\quad}n\geq1.\label{2.23}\\
\int_{0}^{\infty}\frac{\sin^{2n}t}{t^{2n}}dt  & \geq\frac{n-1}{n\sqrt{n}}%
\int_{0}^{\infty}\frac{\sin^{2}t}{t^{2}}dt=\frac{n-1}{n}\frac{1}{\sqrt{n}%
}\frac{\pi}{2},\text{\quad}n\geq1.\label{2.24}%
\end{align}

Also, when $n\geq2$,%
\begin{equation}
\max_{m=2:2n-1}\frac{\int_{0}^{\infty}\frac{\sin^{2n}t}{t^{m}}dt}{\int%
_{0}^{\infty}\frac{\sin^{2n}t}{t^{2n}}dt}=\frac{\int_{0}^{\infty}\frac
{\sin^{2n}t}{t^{2}}dt}{\int_{0}^{\infty}\frac{\sin^{2n}t}{t^{2n}}dt}\leq
\frac{\int_{0}^{\infty}\frac{\sin^{4}t}{t^{2}}dt}{\int_{0}^{\infty}\frac
{\sin^{4}t}{t^{4}}dt}\frac{4}{5-\frac{2}{n}}.\label{2.26}%
\end{equation}

\end{conjecture}

\begin{conclusion}
\label{Conc_CoinfBasis_MatExperim}\textbf{Consequences of Conjecture}
\ref{Conj_CoinfBasis_MatExperim}:

From conjectures \ref{2.23} and \ref{2.24},%
\[
0\leq\frac{\pi}{2}-\sqrt{n}\int_{0}^{\infty}\frac{\sin^{2n}t}{t^{2n}}%
dt\leq\frac{1}{n}\frac{\pi}{2},\text{\quad}n\geq1,
\]

and thus%
\begin{equation}
\lim_{n\rightarrow\infty}\sqrt{n}\int_{0}^{\infty}\frac{\sin^{2n}t}{t^{2n}%
}dt=\frac{\pi}{2}.\label{2.29}%
\end{equation}

Easy to verify that%
\begin{equation}
\frac{\left(  2m-1\right)  !!}{\left(  2m\right)  !!}=\frac{1}{2^{2m}}%
\binom{2m}{m},\text{\quad}m\geq1.\label{2.28}%
\end{equation}

This foormula and 3.821.10 of Gradshteyn \cite{GradRyz07} imply%
\begin{equation}
\int_{0}^{\infty}\frac{\sin^{2n}t}{t^{2}}dt=\frac{\left(  2n-3\right)
!!}{\left(  2n-2\right)  !!}\frac{\pi}{2}=\frac{1}{2^{2n-2}}\binom{2n-2}%
{n-1}\frac{\pi}{2},\label{2.27}%
\end{equation}

and so by conjecture \ref{2.23},%
\begin{equation}
\frac{\int_{0}^{\infty}\frac{\sin^{2n}t}{t^{2}}dt}{\int_{0}^{\infty}\frac
{\sin^{2n}t}{t^{2n}}dt}<\frac{\frac{\left(  2n-3\right)  !!}{\left(
2n-2\right)  !!}\frac{\pi}{2}}{\frac{1}{\sqrt{n}}\frac{\pi}{2}}=\sqrt{n}%
\frac{\left(  2n-3\right)  !!}{\left(  2n-2\right)  !!},\text{\quad}%
n\geq2.\label{2.25}%
\end{equation}

From conjecture \ref{2.26} and the equations of 3.827 of Gradshteyn
\cite{GradRyz07},%
\begin{align}
\max_{m=2:2n-1}\frac{\int_{0}^{\infty}\frac{\sin^{2n}t}{t^{m}}dt}{\int%
_{0}^{\infty}\frac{\sin^{2n}t}{t^{2n}}dt}\leq\frac{\int_{0}^{\infty}\frac
{\sin^{4}t}{t^{2}}dt}{\int_{0}^{\infty}\frac{\sin^{4}t}{t^{4}}dt}\frac
{4}{5-\frac{2}{n}}  & =\frac{\pi/4}{\pi/3}\frac{4}{5-\frac{2}{n}}\nonumber\\
& =\frac{3}{5-\frac{2}{n}}\leq\frac{3}{4}.\label{2.31}%
\end{align}

From conjecture \ref{2.26} and \ref{2.27}, when $n\geq2$,%
\[
\int_{0}^{\infty}\frac{\sin^{2n}t}{t^{2n}}dt>\frac{5-\frac{2}{n}}{3}\int%
_{0}^{\infty}\frac{\sin^{2n}t}{t^{2}}dt=\frac{5-\frac{2}{n}}{3}\frac{\left(
2n-3\right)  !!}{\left(  2n-2\right)  !!}\frac{\pi}{2}.
\]

\end{conclusion}

??? It was proven in [???] that:

\begin{theorem}
Suppose $G_{s}\left(  \cdot;n\right)  $ is the extended B-spline basis
function generated by the parameters $n$ and $l=n$. Then the scaled sequence
of basis functions $\mathcal{G}_{s}\left(  x;n\right)  :=\frac{G_{s}\left(
2nx;n\right)  }{G_{s}\left(  0,n\right)  }$ converges uniformly to a
$C_{0}^{\infty}$ function with support $\left[  -1,1\right]  $. Call this
function $\mathcal{G}_{s}$. Further, all the derivatives of $\mathcal{G}%
_{s}\left(  \cdot;n\right)  $ converge uniformly to the corresponding
derivatives of $\mathcal{G}_{s}$.
\end{theorem}

???

\begin{theorem}
The function $G_{s}\left(  2nx;n\right)  $ has Fourier transform%
\[
F\left[  \mathcal{G}_{s}\left(  \cdot;n\right)  \right]  \left(  \xi\right)
=\frac{1}{G_{s}\left(  0;n\right)  }\frac{1}{2n}\frac{\sin^{2n}\frac{\xi}{2n}%
}{\left(  \frac{\xi}{2n}\right)  ^{2n}},\text{\quad}n=1,2,3,\ldots,
\]

where%
\[
G_{s}\left(  0;n\right)  =\frac{2}{\sqrt{2\pi}}\int_{0}^{\infty}\frac
{\sin^{2n}t}{t^{2n}}dt,
\]

and%
\[
\lim_{n\rightarrow\infty}\sqrt{n}G_{s}\left(  0;n\right)  =\frac{\sqrt{2\pi}%
}{2}.
\]

\end{theorem}

If the sequence of central difference basis functions $G_{k}$ and the B-spline
basis function $G_{s}$ are all scaled so they have value 1 at the origin then
we have the following convergence result for the derivatives:

\begin{theorem}
In Theorem \ref{Thm_centr_diff_basis_conv_to_ext_spl_basis} the sequence of
functions
\[
q_{k}\left(  t\right)  =\sqrt{2\pi}2^{-2n}\frac{k}{2}\psi\left(  \frac{k}%
{2}\left(  t-2\right)  \right)  ,\quad k=1,2,3,\ldots
\]

was used to generate the sequence of central difference basis functions
$G_{k}$ given in \ref{a2.36} by%
\[
G_{k}\left(  x\right)  =\int_{\mathbb{R}^{1}}G_{s}\left(  \frac{x}{t}\right)
\left\vert t\right\vert ^{2n-1}k\psi\left(  k\left(  t-1\right)  \right)  dt.
\]

Here I show that for $x\in\mathbb{R}^{1}$,%
\begin{equation}
\left\vert \frac{D^{j}G_{k}\left(  x\right)  }{G_{k}\left(  0\right)  }%
-\frac{D^{j}G_{s}\left(  x\right)  }{G_{s}\left(  0\right)  }\right\vert
\leq\frac{\left\Vert D^{j+1}G_{s}\right\Vert _{\infty}}{G_{s}\left(  0\right)
}\frac{\left\vert x\right\vert }{k}\left(  1-\frac{1}{k}\right)  ^{-j}%
,\quad\left\{
\begin{array}
[c]{l}%
j\leq2n-2,\\
k=2,3,4,\ldots
\end{array}
\right. \label{2.30}%
\end{equation}

\end{theorem}

\begin{proof}
From Theorem \ref{Thm_centr_diff_basis_conv_to_ext_spl_basis},%
\[
\psi\in L^{1}\left(  \mathbb{R}^{1}\right)  ,\text{ }\psi\left(  x\right)
\geq0,\text{ }\int\psi=1,\text{ }\operatorname*{supp}\psi\subset\left[
-1,1\right]  ,
\]

and%
\[
\operatorname*{supp}G_{k}=\left(  1+\frac{1}{k}\right)  2n.
\]

From \ref{a940}: for $1\leq j\leq2n-2$,%
\begin{align}
D^{j}G_{k}\left(  x\right)   & =\tfrac{2^{j+1}}{\sqrt{2\pi}}\int\left(
D^{j}G_{s}\right)  \left(  \frac{2x}{\left\vert t\right\vert }\right)
\left\vert t\right\vert ^{2n-j-1}q_{k}\left(  t\right)  dt\nonumber\\
& :u=t/2,\text{ }dt=2du\Rightarrow\nonumber\\
& =\tfrac{2^{2n+1}}{\sqrt{2\pi}}\int\left(  D^{j}G_{s}\right)  \left(
\frac{x}{\left\vert u\right\vert }\right)  \left\vert u\right\vert
^{2n-j-1}q_{k}\left(  2u\right)  du\nonumber\\
& =\tfrac{2^{2n+1}}{\sqrt{2\pi}}\int\left(  D^{j}G_{s}\right)  \left(
\frac{x}{\left\vert u\right\vert }\right)  \left\vert u\right\vert
^{2n-j-1}\frac{\sqrt{2\pi}}{2^{2n}}\frac{k}{2}\psi\left(  \frac{k}{2}\left(
2u-2\right)  \right)  du\nonumber\\
& =\int\left(  D^{j}G_{s}\right)  \left(  \frac{x}{\left\vert u\right\vert
}\right)  \left\vert u\right\vert ^{2n-j-1}k\psi\left(  k\left(  u-1\right)
\right)  du\nonumber\\
& =\int_{1-\frac{1}{k}}^{1+\frac{1}{k}}\left(  D^{j}G_{s}\right)  \left(
\frac{x}{u}\right)  u^{2n-j-1}k\psi\left(  k\left(  u-1\right)  \right)
du.\label{2.33}%
\end{align}

Thus%
\begin{align*}
&  \frac{D^{j}G_{k}\left(  x\right)  }{G_{k}\left(  0\right)  }-\frac
{D^{j}G_{s}\left(  x\right)  }{G_{s}\left(  0\right)  }\\
&  =\frac{\int_{1-\frac{1}{k}}^{1+\frac{1}{k}}\left(  D^{j}G_{s}\right)
\left(  \frac{x}{t}\right)  t^{2n-j-1}k\psi\left(  k\left(  t-1\right)
\right)  dt}{G_{s}\left(  0\right)  \int_{1-\frac{1}{k}}^{1+\frac{1}{k}%
}t^{2n-1}k\psi\left(  k\left(  t-1\right)  \right)  dt}-\frac{D^{j}%
G_{s}\left(  x\right)  }{G_{s}\left(  0\right)  }\\
&  =\frac{\int\limits_{1-\frac{1}{k}}^{1+\frac{1}{k}}\left(  D^{j}%
G_{s}\right)  \left(  \frac{x}{t}\right)  t^{2n-j-1}k\psi\left(  k\left(
t-1\right)  \right)  dt-\int\limits_{1-\frac{1}{k}}^{1+\frac{1}{k}}\left(
D^{j}G_{s}\right)  \left(  x\right)  t^{2n-j-1}k\psi\left(  k\left(
t-1\right)  \right)  dt}{G_{s}\left(  0\right)  \int_{1-\frac{1}{k}}%
^{1+\frac{1}{k}}t^{2n-1}k\psi\left(  k\left(  t-1\right)  \right)  dt}\\
&  =\frac{\int_{1-\frac{1}{k}}^{1+\frac{1}{k}}\left(  \left(  D^{j}%
G_{s}\right)  \left(  \frac{x}{t}\right)  -\left(  D^{j}G_{s}\right)  \left(
x\right)  \right)  t^{2n-j-1}k\psi\left(  k\left(  t-1\right)  \right)
dt}{G_{s}\left(  0\right)  \int_{1-\frac{1}{k}}^{1+\frac{1}{k}}t^{2n-1}%
k\psi\left(  k\left(  t-1\right)  \right)  dt},
\end{align*}

so that if $k\geq2$, by using Theorem \ref{Thm_ex_nat_spline_basis_Lipschitz}
we get%
\begin{align*}
\left\vert \frac{D^{j}G_{k}\left(  x\right)  }{G_{k}\left(  0\right)  }%
-\frac{D^{j}G_{s}\left(  x\right)  }{G_{s}\left(  0\right)  }\right\vert  &
\leq\frac{\int_{1-\frac{1}{k}}^{1+\frac{1}{k}}\left\vert \left(  D^{j}%
G_{s}\right)  \left(  \frac{x}{t}\right)  -\left(  D^{j}G_{s}\right)  \left(
x\right)  \right\vert t^{2n-j-1}k\psi\left(  k\left(  t-1\right)  \right)
dt}{G_{s}\left(  0\right)  \int_{1-\frac{1}{k}}^{1+\frac{1}{k}}t^{2n-1}%
k\psi\left(  k\left(  t-1\right)  \right)  dt}\\
&  \leq\frac{\left\Vert D^{j+1}G_{s}\right\Vert _{\infty}\int_{1-\frac{1}{k}%
}^{1+\frac{1}{k}}\left\vert \frac{x}{t}-x\right\vert t^{2n-j-1}k\psi\left(
k\left(  t-1\right)  \right)  dt}{G_{s}\left(  0\right)  \int_{1-\frac{1}{k}%
}^{1+\frac{1}{k}}t^{2n-1}k\psi\left(  k\left(  t-1\right)  \right)  dt}\\
&  =\frac{\left\Vert D^{j+1}G_{s}\right\Vert _{\infty}\left\vert x\right\vert
\int_{1-\frac{1}{k}}^{1+\frac{1}{k}}\frac{\left\vert t-1\right\vert }{t^{j}%
}t^{2n-1}k\psi\left(  k\left(  t-1\right)  \right)  dt}{G_{s}\left(  0\right)
\int_{1-\frac{1}{k}}^{1+\frac{1}{k}}t^{2n-1}k\psi\left(  k\left(  t-1\right)
\right)  dt}\\
&  \leq\frac{\left\Vert D^{j+1}G_{s}\right\Vert _{\infty}\frac{\left\vert
x\right\vert }{k}\left(  1-\frac{1}{k}\right)  ^{-j}\int_{1-\frac{1}{k}%
}^{1+\frac{1}{k}}t^{2n-1}k\psi\left(  k\left(  t-1\right)  \right)  dt}%
{G_{s}\left(  0\right)  \int_{1-\frac{1}{k}}^{1+\frac{1}{k}}t^{2n-1}%
k\psi\left(  k\left(  t-1\right)  \right)  dt}\\
&  =\frac{\left\Vert D^{j+1}G_{s}\right\Vert _{\infty}}{G_{s}\left(  0\right)
}\frac{\left\vert x\right\vert }{k}\left(  1-\frac{1}{k}\right)  ^{-j}.
\end{align*}

\end{proof}

\begin{corollary}
(\textbf{Uses several conjectures}) If $k\geq2$ and $n\geq2$ then for all
$j\leq2n-3$,
\begin{align*}
\left\vert \frac{D^{j}G_{k}\left(  x\right)  }{G_{k}\left(  0\right)  }%
-\frac{D^{j}G_{s}\left(  x\right)  }{G_{s}\left(  0\right)  }\right\vert  &
\leq\frac{\int_{0}^{\infty}\frac{\sin^{2n}t}{t^{2n-1-j}}dt}{\int_{0}^{\infty
}\frac{\sin^{2n}t}{t^{2n}}dt}\frac{\left\vert x\right\vert }{k}\left(
1-\frac{1}{k}\right)  ^{-j}\\
& \leq\frac{3}{5-2/n}\frac{\left\vert x\right\vert }{k}\left(  1-\frac{1}%
{k}\right)  ^{-j},
\end{align*}

and for $j=2n-2$,
\begin{align*}
\left\vert \frac{D^{2n-2}G_{k}\left(  x\right)  }{G_{k}\left(  0\right)
}-\frac{D^{2n-2}G_{s}\left(  x\right)  }{G_{s}\left(  0\right)  }\right\vert
& \leq\tfrac{\pi}{2^{2n+1}}\frac{\tbinom{2n-1}{n-1}}{\int_{0}^{\infty}%
\frac{\sin^{2n}t}{t^{2n}}dt}\frac{\left\vert x\right\vert }{k}\left(
1-\frac{1}{k}\right)  ^{-\left(  2n-1\right)  }\\
& \leq\frac{1}{4\sqrt{\pi}}\frac{1}{1-1/n}\frac{\left\vert x\right\vert }%
{k}\left(  1-\frac{1}{k}\right)  ^{-\left(  2n-1\right)  }.
\end{align*}

For all $j\leq2n-2$,
\[
\left\vert \frac{D^{j}G_{k}\left(  x\right)  }{G_{k}\left(  0\right)  }%
-\frac{D^{j}G_{s}\left(  x\right)  }{G_{s}\left(  0\right)  }\right\vert
\leq\frac{3}{4}\frac{\left\vert x\right\vert }{k}\left(  1-\frac{1}{k}\right)
^{-\left(  2n-1\right)  }.
\]

\end{corollary}

\begin{proof}
\fbox{\textbf{Case} $j\leq2n-3$} From Theorem \ref{Thm_basis_tensor_hat_W3}
and Theorem \ref{Thm_basis_fn_properties_all_m_W3},%
\[
D^{p}G_{s}\left(  x\right)  =\frac{1}{\sqrt{2\pi}}\int e^{-ixt}\frac{\left(
-it\right)  ^{p}}{w_{s}\left(  t\right)  }dt,\quad p\leq2n-1,
\]

which implies%
\[
\left\Vert D^{p}G_{s}\right\Vert _{\infty}\leq\frac{1}{\sqrt{2\pi}}\int%
\frac{\left\vert t\right\vert ^{p}}{w_{s}\left(  t\right)  }dt=\frac{2}%
{\sqrt{2\pi}}\int_{0}^{\infty}\frac{\sin^{2n}t}{t^{2n-p}}dt,\quad p\leq2n-1,
\]

and thus from \ref{2.30},%
\begin{align*}
\left\vert \frac{D^{j}G_{k}\left(  x\right)  }{G_{k}\left(  0\right)  }%
-\frac{D^{j}G_{s}\left(  x\right)  }{G_{s}\left(  0\right)  }\right\vert  &
\leq\frac{\left\Vert D^{j+1}G_{s}\right\Vert _{\infty}}{G_{s}\left(  0\right)
}\frac{\left\vert x\right\vert }{k}\left(  1-\frac{1}{k}\right)  ^{-j}\\
& \leq\frac{\frac{2}{\sqrt{2\pi}}\int_{0}^{\infty}\frac{\sin^{2n}t}%
{t^{2n-1-j}}dt}{\frac{2}{\sqrt{2\pi}}\int_{0}^{\infty}\frac{\sin^{2n}t}%
{t^{2n}}dt}\frac{\left\vert x\right\vert }{k}\left(  1-\frac{1}{k}\right)
^{-j}\\
& =\frac{\int_{0}^{\infty}\frac{\sin^{2n}t}{t^{2n-1-j}}dt}{\int_{0}^{\infty
}\frac{\sin^{2n}t}{t^{2n}}dt}\frac{\left\vert x\right\vert }{k}\left(
1-\frac{1}{k}\right)  ^{-j}\\
& \leq\left(  \max_{j\leq2n-3}\frac{\int_{0}^{\infty}\frac{\sin^{2n}%
t}{t^{2n-1-j}}dt}{\int_{0}^{\infty}\frac{\sin^{2n}t}{t^{2n}}dt}\right)
\frac{\left\vert x\right\vert }{k}\left(  1-\frac{1}{k}\right)  ^{-j}\\
& =\left(  \max_{2\leq m\leq2n-1}\frac{\int_{0}^{\infty}\frac{\sin^{2n}%
t}{t^{m}}dt}{\int_{0}^{\infty}\frac{\sin^{2n}t}{t^{2n}}dt}\right)
\frac{\left\vert x\right\vert }{k}\left(  1-\frac{1}{k}\right)  ^{-j}\\
& :conjecture\text{ }\ref{2.31}\Rightarrow\\
& <\frac{3}{5-\frac{2}{n}}\frac{\left\vert x\right\vert }{k}\left(  1-\frac
{1}{k}\right)  ^{-j}.
\end{align*}
\medskip

\fbox{Case $j=2n-2$} Theorem \ref{Thm_basis_tensor_hat_W3} implies%
\[
\left\Vert D^{2n-1}G_{s}\right\Vert _{\infty}=\tfrac{\sqrt{2\pi}}{2^{2n+1}%
}\tbinom{2n-1}{n-1}.
\]

so that \ref{2.30} becomes,%
\begin{align*}
\left\vert \frac{D^{2n-2}G_{k}\left(  x\right)  }{G_{k}\left(  0\right)
}-\frac{D^{2n-2}G_{s}\left(  x\right)  }{G_{s}\left(  0\right)  }\right\vert
& \leq\frac{\left\Vert D^{2n-1}G_{s}\right\Vert _{\infty}}{G_{s}\left(
0\right)  }\frac{\left\vert x\right\vert }{k}\left(  1-\frac{1}{k}\right)
^{-\left(  2n-1\right)  }\\
& =\frac{\tfrac{\sqrt{2\pi}}{2^{2n+2}}\tbinom{2n-1}{n-1}}{\frac{2}{\sqrt{2\pi
}}\int_{0}^{\infty}\frac{\sin^{2n}t}{t^{2n}}dt}\frac{\left\vert x\right\vert
}{k}\left(  1-\frac{1}{k}\right)  ^{-\left(  2n-1\right)  }\\
& =\tfrac{\pi}{2^{2n+2}}\frac{\tbinom{2n-1}{n-1}}{\int_{0}^{\infty}\frac
{\sin^{2n}t}{t^{2n}}dt}\frac{\left\vert x\right\vert }{k}\left(  1-\frac{1}%
{k}\right)  ^{-\left(  2n-1\right)  }.
\end{align*}

But inequality \ref{Ap005} is%
\[
2^{-2n+1}\binom{2n-1}{n-1}<\frac{1}{\sqrt{\pi n}},
\]

and the \textbf{conjectured inequality} \ref{2.24} is%
\[
\int_{0}^{\infty}\frac{\sin^{2n}t}{t^{2n}}dt\geq\frac{n-1}{n}\frac{1}{\sqrt
{n}}\frac{\pi}{2},
\]

so that now%
\begin{align*}
\left\vert \frac{D^{2n-2}G_{k}\left(  x\right)  }{G_{k}\left(  0\right)
}-\frac{D^{2n-2}G_{s}\left(  x\right)  }{G_{s}\left(  0\right)  }\right\vert
& \leq\tfrac{\pi}{2^{2n+2}}\frac{2^{2n-1}\frac{1}{\sqrt{\pi n}}}{\frac{n-1}%
{n}\frac{1}{\sqrt{n}}\frac{\pi}{2}}\frac{\left\vert x\right\vert }{k}\left(
1-\frac{1}{k}\right)  ^{-\left(  2n-1\right)  }\\
& =\tfrac{1}{4}\frac{\frac{1}{\sqrt{\pi}}}{\frac{n-1}{n}}\frac{\left\vert
x\right\vert }{k}\left(  1-\frac{1}{k}\right)  ^{-\left(  2n-1\right)  }\\
& =\frac{1}{4\sqrt{\pi}}\frac{n}{n-1}\frac{\left\vert x\right\vert }{k}\left(
1-\frac{1}{k}\right)  ^{-\left(  2n-1\right)  }.
\end{align*}

\end{proof}

Now we want to make the basis function parameter $n$ explicit and write%
\[
G_{k}\left(  x;n\right)  =G_{k}\left(  x\right)  ,\quad G_{s}\left(
x;n\right)  =G_{s}\left(  x\right)  ,
\]

and define the scaled basis functions%
\[
\mathcal{G}_{n}\left(  x\right)  :=\frac{G_{k_{n}}\left(  2nx;n\right)
}{G_{k_{n}}\left(  0;n\right)  },\quad\mathcal{G}_{s}\left(  x;n\right)
:=\frac{G_{s}\left(  2nx;n\right)  }{G_{s}\left(  0;n\right)  },\quad n\geq2,
\]

where $k_{n}\rightarrow\infty$ is a sequence to be determined. Observe that%
\begin{align*}
\operatorname*{supp}\mathcal{G}_{s}  & \subseteq\left[  -1,1\right]
,\quad\operatorname*{supp}\mathcal{G}_{n}\subseteq\left[  -1,1\right]
+\frac{1}{k_{n}}.\\
\mathcal{G}_{s}\left(  0;n\right)   & =1,\quad\mathcal{G}_{n}\left(  0\right)
=1.
\end{align*}

Now for $j\leq2n-2$,%
\begin{align*}
\left\vert D^{j}\mathcal{G}_{n}\left(  x\right)  -D^{j}\mathcal{G}_{s}\left(
x;n\right)  \right\vert  & =\left\vert D^{j}\frac{G_{k_{n}}\left(
2nx;n\right)  }{G_{k_{n}}\left(  0;n\right)  }-D^{j}\frac{G_{s}\left(
2nx;n\right)  }{G_{s}\left(  0;n\right)  }\right\vert \\
& =\left(  2n\right)  ^{j}\left\vert \frac{\left(  D^{j}G_{k_{n}}\right)
\left(  2nx;n\right)  }{G_{k_{n}}\left(  0;n\right)  }-\frac{\left(
D^{j}G_{s}\right)  \left(  2nx;n\right)  }{G_{s}\left(  0;n\right)
}\right\vert \\
& \leq\left(  2n\right)  ^{j}\frac{3}{4}\frac{\left\vert 2nx\right\vert
}{k_{n}}\left(  1-\frac{1}{k_{n}}\right)  ^{-\left(  2n-1\right)  }\\
& \leq\frac{3}{4}\frac{\left(  2n\right)  ^{j+1}}{k_{n}}\left(  1+\frac
{1}{k_{n}}\right)  \left(  1-\frac{1}{k_{n}}\right)  ^{-\left(  2n-1\right)
}\\
& \leq\frac{3}{4}\frac{\left(  2n\right)  ^{2n-1}}{k_{n}}\left(  1+\frac
{1}{k_{n}}\right)  \left(  1-\frac{1}{k_{n}}\right)  ^{-\left(  2n-1\right)  }%
\end{align*}

We now need the inequality:
\begin{equation}
1+p\varepsilon\leq\left(  1-\varepsilon\right)  ^{-p}\leq1+p\sqrt
{p}\varepsilon\text{ }when\text{ }p\geq1\text{ }and\text{ }0\leq
\varepsilon\leq\frac{1}{2\left(  p+1\right)  },\label{2.32}%
\end{equation}

which is proved by showing that $D_{\varepsilon}\left(  \left(  1-\varepsilon
\right)  ^{p}\left(  1+p\varepsilon\right)  \right)  \leq0$ and
$D_{\varepsilon}\left(  \left(  1-\varepsilon\right)  ^{p}\left(  1+p\sqrt
{p}\varepsilon\right)  \right)  \geq0$.

Next assume that $k_{n}\geq\left(  2n\right)  ^{2n-1}$. Then $\frac{1}{k_{n}%
}\leq\frac{1}{\left(  2n\right)  ^{2n-1}}\leq\frac{1}{4n}$ and \ref{2.32}
implies%
\[
\left(  1-\frac{1}{k_{n}}\right)  ^{-\left(  2n-1\right)  }\leq1+\frac{\left(
2n-1\right)  \sqrt{2n-1}}{k_{n}}\leq1+\frac{2n\sqrt{2n}}{k_{n}},
\]

and hence%
\begin{align*}
\left\vert D^{j}\mathcal{G}_{n}\left(  x\right)  -D^{j}\mathcal{G}_{s}\left(
x;n\right)  \right\vert  & \leq\frac{3}{4}\frac{\left(  2n\right)  ^{2n-1}%
}{k_{n}}\left(  1+\frac{1}{k_{n}}\right)  \left(  1+\frac{2n\sqrt{2n}}{k_{n}%
}\right) \\
& \leq\frac{3}{4}\frac{\left(  2n\right)  ^{2n-1}}{k_{n}}\left(  1+\frac
{1}{4^{3}}\right)  \left(  1+\frac{8}{4^{3}}\right) \\
& <0.86\frac{\left(  2n\right)  ^{2n-1}}{k_{n}}.
\end{align*}

Thus we have shown that:

\begin{theorem}
Suppose $n\geq2$, $j\leq2n-2$ and $k_{n}\geq1$. Then%
\[
\left\vert D^{j}\mathcal{G}_{n}\left(  x\right)  -D^{j}\mathcal{G}_{s}\left(
x;n\right)  \right\vert \leq\frac{3}{4}\frac{\left(  2n\right)  ^{2n-1}}%
{k_{n}}\left(  1+\frac{1}{k_{n}}\right)  \left(  1-\frac{1}{k_{n}}\right)
^{-\left(  2n-1\right)  }.
\]

If, in addition, $k_{n}\geq\left(  2n\right)  ^{2n-1}$ then%
\begin{align*}
\left\vert D^{j}\mathcal{G}_{n}\left(  x\right)  -D^{j}\mathcal{G}_{s}\left(
x;n\right)  \right\vert  & \leq\frac{3}{4}\frac{\left(  2n\right)  ^{2n-1}%
}{k_{n}}\left(  1+\frac{1}{k_{n}}\right)  \left(  1+\frac{2n\sqrt{2n}}{k_{n}%
}\right) \\
& <0.86\frac{\left(  2n\right)  ^{2n-1}}{k_{n}}.
\end{align*}

\end{theorem}

\begin{example}
Set $k_{n}=\left(  2n\right)  ^{2n-1}n^{s}$ for $s\geq1$. Then for $n\geq2$
and $j\leq2n-2$,%
\begin{align*}
\left\vert D^{j}\mathcal{G}_{n}\left(  x\right)  -D^{j}\mathcal{G}_{s}\left(
x;n\right)  \right\vert  & \leq\frac{3}{4}\frac{1}{n^{s}}\left(  1+\frac
{1}{\left(  2n\right)  ^{2n-1}n^{s}}\right)  \left(  1+\frac{2n\sqrt{2n}%
}{\left(  2n\right)  ^{2n-1}n^{s}}\right) \\
& =\frac{3}{4}\frac{1}{n^{s}}\left(  1+\frac{1}{\left(  2n\right)
^{2n-1}n^{s}}\right)  \left(  1+\frac{1}{\left(  2n\right)  ^{2n-5/2}n^{s}%
}\right) \\
& <\frac{0.86}{n^{s}}.
\end{align*}

\end{example}

\begin{remark}
?? An $L^{1}$ convergence result similar to Theorem
\ref{Thm_CntDiffWtFnConvergToBSplin} can be derived for the weight functions
corresponding to the basis functions $\mathcal{G}_{n}$ and $\mathcal{G}%
_{s}\left(  \cdot;n\right)  $.
\end{remark}

\section{Local data function spaces\label{Sect_CntDifWtFn_DataFuncs}}

In this section we will use the generalized local data function results of
Subsection \ref{SbSect_loc_data_larger_class} to characterize locally the data
functions for the tensor product central difference basis functions as Sobolev
spaces. These results supply important information about the data functions
and makes it easy to choose data functions for numerical experiments
concerning the zero order basis function interpolation and smoothing problems
discussed later in Chapters \ref{Ch_Interpol}, \ref{Ch_Exact_smth} and
\ref{Ch_Approx_smth}. It also allows us to use the result of Subsection
\ref{SbSect_global_error_to_local} to derive local convergence estimates from
global convergence estimates for interpolants and smoothers.

\subsection{Tensor product weight functions}

Here we are interested in tensor product weight functions e.g. the
multivariate central difference weight functions and the extended B-splines
\ref{1.032}. We want a result which allows us to prove that there exists
integral $n$ such that the local Sobolev space%
\[
W^{n\mathbf{1}}\left(  \Omega\right)  =\left\{  u\in L^{2}\left(
\Omega\right)  :D^{\alpha}u\in L^{2}\left(  \Omega\right)  \text{ }for\text{
}\alpha\leq n\mathbf{1}\right\}  ,
\]

of Definition \ref{Def_SobolevSpace2} is a subset of the restriction space
$X_{w}^{0}\left(  \Omega\right)  $ (recall Definition \ref{Def_Xow(open_set)}%
). Further we want this result to (essentially) involve \textbf{only proving
estimates about a single univariate weight function}. This result is Theorem
\ref{Thm_data_fn_tensor_prod} below and will allow us to justify using
$W^{n\mathbf{1}}\left(  \Omega\right)  $ functions as data functions in
numerical experiments.

Further, if $w_{1}$ satisfies the assumptions of Corollary
\ref{Cor_Thm_data_fn_tensor_prod} then $W^{n\mathbf{1}}\left(  \Omega\right)
=X_{w}^{0}\left(  \Omega\right)  $ as sets and their norms are equivalent.

\begin{theorem}
\label{Thm_data_fn_tensor_prod}Let $w$ be a \textbf{tensor product} weight
function on $\mathbb{R}^{d}$ generated by the single \textbf{univariate}
weight function $w_{1}\in W03$. Suppose:

\begin{enumerate}
\item $\Omega$ is a bounded region with the segment or uniform rectangle
property described in Definition \ref{Def_UnifRectCondit}.

\item Assume there exists a distribution $\theta_{1}$ such that
\begin{equation}
\theta_{1}\in S^{\prime}\left(  \mathbb{R}^{1}\right)  ,\quad
\operatorname*{supp}\theta_{1}\subset\mathbb{R}^{1}\setminus2\left(
-1,1\right)  ,\quad\widehat{\theta}_{1}\in L_{loc}^{1}\left(  \mathbb{R}%
^{1}\right)  ,\label{a902}%
\end{equation}

and%
\begin{equation}
w_{1}\left(  t\right)  \left\vert \widehat{\theta}_{1}\left(  t\right)
+1\right\vert ^{2}\leq C_{\theta_{1}}\left(  1+t^{2}\right)  ^{n},\quad
a.e.\text{ }on\text{ }\mathbb{R}^{1},\label{a907}%
\end{equation}

for some constant $C_{\theta_{1}}>0$.

\item There exists an integer $m\geq1$ and a constant $c_{1;m}>0$ such that
\begin{equation}
w_{1}\left(  t/m\right)  \leq c_{1;m}w_{1}\left(  t\right)  ,\quad a.e.\text{
}on\text{ }\mathbb{R}^{1}.\label{a906}%
\end{equation}

\end{enumerate}

Then for all scalar $\lambda>0$, $W^{n\mathbf{1}}\left(  \Omega\right)
\hookrightarrow X_{\sigma_{\lambda}w}^{0}\left(  \Omega\right)  $ where
$\sigma_{\lambda}$ is the dilation operator.
\end{theorem}

\begin{proof}
The proof involves showing the assumptions 1 to 4 required by Theorem
\ref{Thm_ex_data_fn_3} are met so that Corollary \ref{Cor_Thm_data_fn_3} can
be applied.\medskip

\textbf{Assumption 1} By Corollary \ref{Cor_Thm_ten_prod_two_wt_fns}, $w$ has
property W03 since $w_{1}$ has property W03.\medskip

\textbf{Assumption 2} Properties\textbf{\ }\ref{a1.05} and \ref{a1.04}: Noting
\ref{a1.05}, we use the tensor (or direct) product of distributions to define
$\chi,\eta\in S^{\prime}\left(  \mathbb{R}^{d}\right)  $ by%
\begin{align}
\chi\left(  x\right)   & =\prod\limits_{k=1}^{d}\left(  \theta_{1}\left(
x_{k}\right)  +\sqrt{2\pi}\delta_{1}\left(  x_{k}\right)  \right)
,\label{a2.01}\\
\eta & =\chi-\left(  2\pi\right)  ^{d/2}\delta,\label{a2.02}%
\end{align}

so that $\operatorname*{supp}\eta\subset\left(  \mathbb{R}^{1}\setminus
2\left(  -1,1\right)  \right)  ^{d}\subset\mathbb{R}^{d}\setminus2\left(
-1,1\right)  ^{d}$, and clearly $\widehat{\eta}\in L_{loc}^{1}\left(
\mathbb{R}^{d}\right)  $.\medskip

Property\textbf{\ }\ref{a1.052}: if $C_{\chi}=C_{\theta_{1}}^{d}$ then
\begin{align*}
w\left(  \xi\right)  \left\vert \widehat{\chi}\left(  \xi\right)  \right\vert
^{2}=\prod\limits_{k=1}^{d}w_{1}\left(  \xi_{k}\right)  \left\vert
\widehat{\chi_{1}}\left(  \xi_{k}\right)  \right\vert ^{2}  & =\prod
\limits_{k=1}^{d}w_{1}\left(  \xi_{k}\right)  \left\vert \widehat{\theta}%
_{1}\left(  \xi_{k}\right)  +1\right\vert ^{2}\\
& \leq C_{\chi}\prod\limits_{k=1}^{d}\left(  1+\xi_{k}^{2}\right)  ^{n},
\end{align*}

which confirms \ref{a1.052}.\medskip

\textbf{Assumption 3} The uniform rectangle property implies the segment
condition which is assumption 3 of Theorem \ref{Thm_ex_data_fn_3}. The uniform
rectangle property was described in Definition \ref{Def_UnifRectCondit} and
the segment condition in Definition \ref{Def_SegCondit}.\medskip

\textbf{Assumption 4} Property\textbf{\ }\ref{a1.055}: If inequality
\ref{a906} is valid then $\eta>\mathbf{1}$ and a.e. on $\mathbb{R}^{d}$,
\[
w\left(  \xi/m\right)  =\prod\limits_{k=1}^{d}w_{1}\left(  \xi_{k}/m\right)
\leq\prod\limits_{k=1}^{d}\left(  c_{1;m}w_{1}\left(  \xi_{k}\right)  \right)
=\left(  c_{1;m}\right)  ^{d}w\left(  \xi\right)  ,
\]

which means that condition \ref{a1.055} holds if we choose $c_{m}=\left(
c_{1;m}\right)  ^{d}$.\medskip

The conclusion of this theorem now follows from Corollary
\ref{Cor_Thm_data_fn_3}.
\end{proof}

This corollary supplies a condition under which $X_{w}^{0}\left(
\Omega\right)  =W^{n\mathbf{1}}\left(  \Omega\right)  $ as sets with
equivalent norms.

\begin{corollary}
\label{Cor_Thm_data_fn_tensor_prod}Suppose the assumptions of Theorem
\ref{Thm_data_fn_tensor_prod} hold. Also suppose there exists a constant
$c_{w_{1}}>0$ such that the univariate weight function $w_{1}$ satisfies%
\begin{equation}
c_{w_{1}}\left(  1+t^{2}\right)  ^{n}\leq w_{1}\left(  t\right)  ,\quad
a.e.\text{ }on\text{ }\mathbb{R}^{1}.\label{a905}%
\end{equation}

Then for each $\lambda>0$, $X_{\sigma_{\lambda}w}^{0}\left(  \Omega\right)
=W^{n\mathbf{1}}\left(  \Omega\right)  $ as sets and their norms are equivalent.
\end{corollary}

\begin{proof}
We must show that condition \ref{1.391} of Corollary \ref{Cor_2_Thm_data_fn_3}
holds. Indeed, if we choose $c_{w}=\left(  c_{w_{1}}\right)  ^{d}$ then a.e.
on $\mathbb{R}^{1}$,%
\[
w\left(  \xi\right)  =\prod\limits_{j=1}^{d}w_{1}\left(  \xi_{j}\right)
\geq\prod\limits_{j=1}^{d}c_{w_{1}}\left(  1+\xi_{j}^{2}\right)  ^{n}%
=c_{w}w_{n}\left(  \xi\right)  ,
\]

as required.
\end{proof}

\subsection{Examples}

We present two applications of the theory of local data functions introduced
in the last subsection to tensor product central difference weight functions.

In the first theorem the strong integral condition \ref{a959} is imposed on
the weight function because this condition is associated with simple upper and
lower bounds on $w$ which imply $X_{w}^{0}\left(  \Omega\right)
=W^{n\mathbf{1}}\left(  \Omega\right)  $ as sets. In the second theorem
(Theorem \ref{Thm_ex2}) condition \ref{a959} is weakened and inequality
\ref{a922} imposed.

\begin{theorem}
\label{Thm_ex1}Suppose $w$ is a tensor product \textbf{central difference}
weight function on $\mathbb{R}^{d}$ with parameters $n,l,q$ which satisfies
the condition \ref{a959} i.e. $\int\nolimits_{\left\vert t\right\vert \geq
R}t^{2l}q\left(  t\right)  dt<\infty$. Further, suppose that $\Omega$ is a
bounded region with the segment property or the uniform rectangle property.

Then for each $\lambda>0$, $X_{\sigma_{\lambda}w}^{0}\left(  \Omega\right)
=W^{n\mathbf{1}}\left(  \Omega\right)  $ as sets and the norms are equivalent.
\end{theorem}

\begin{proof}
\textbf{Assumption 1 of Theorem} \ref{Thm_data_fn_tensor_prod} is clearly
satisfied.\smallskip

We can write $w\left(  \xi\right)  =w_{1}\left(  \xi_{1}\right)  \ldots
w_{1}\left(  \xi_{d}\right)  $ where $w_{1}$ is the univariate central
difference weight function. Starting with definition \ref{a962} we have
$w_{1}\left(  t\right)  \Delta_{2l}\widehat{q}\left(  t\right)  =t^{2n}$ where
$q $ satisfies $q\in L^{1}\left(  \mathbb{R}^{1}\right)  $, $q\neq0$,
$q\left(  \xi\right)  \geq0$ so that $\int q>0$.\smallskip

\textbf{Assumption 2 of Theorem} \ref{Thm_data_fn_tensor_prod}: We will
consider the two cases $n<l$ and $n=l$.\smallskip

\fbox{\textbf{Case} $n<l$} Set $p=l-n$ and choose
\begin{equation}
\widehat{\theta_{1}}\left(  t\right)  =\left(  \left(  -4\right)  ^{p}%
\tbinom{2p}{p}^{-1}\sin^{2p}t\right)  -1,\label{a2.03}%
\end{equation}

so that $\widehat{\theta_{1}}\in C_{B}^{\infty}\subset L_{loc}^{1}$. Hence
$\theta_{1}\left(  t\right)  =\left(  -4\right)  ^{p}\tbinom{2p}{p}%
^{-1}\left(  \sin^{2p}t\right)  ^{\vee}-\sqrt{2\pi}\delta$.

Noting that $\sin^{2p}t=\left(  \frac{e^{it}-e^{-it}}{2i}\right)  ^{2p}$ we
use the identity%
\[
\left(  s+s^{-1}\right)  ^{2p}=\tbinom{2p}{p}+\sum\limits_{k=-p,\text{ }%
k\neq0}^{p}\tbinom{2p}{p+k}s^{2k},
\]

to write%
\begin{align*}
\sin^{2p}t  & =\left(  2i\right)  ^{-2p}\left(  \tbinom{2p}{p}+\sum
\limits_{k=-p,\text{ }k\neq0}^{p}\tbinom{2p}{p+k}e^{2ik}\right) \\
& =\left(  -4\right)  ^{-p}\left(  \tbinom{2p}{p}+\sum\limits_{k=-p,\text{
}k\neq0}^{p}\tbinom{2p}{p+k}e^{2ik}\right)  ,
\end{align*}

and since $\left(  e^{2ik}\right)  ^{\vee}=\sqrt{2\pi}\delta\left(
\xi-2k\right)  $,%
\[
\left(  \sin^{2p}t\right)  ^{\vee}=\sqrt{2\pi}\left(  -4\right)  ^{-p}\left(
\tbinom{2p}{p}\delta+\sum\limits_{k=-p,\text{ }k\neq0}^{p}\tbinom{2p}%
{p+k}\delta\left(  \cdot-2k\right)  \right)  ,
\]

so that%
\begin{align}
\theta_{1}\left(  t\right)   & =\left(  -4\right)  ^{p}\tbinom{2p}{p}%
^{-1}\left(  \sin^{2p}t\right)  ^{\vee}-\sqrt{2\pi}\delta\nonumber\\
& =\left(  -4\right)  ^{p}\tbinom{2p}{p}^{-1}\left(  \sqrt{2\pi}\left(
-4\right)  ^{-p}\left(  \tbinom{2p}{p}\delta+\sum\limits_{k=-p,\text{ }k\neq
0}^{p}\tbinom{2p}{p+k}\delta\left(  t-2k\right)  \right)  \right)  -\sqrt
{2\pi}\delta\nonumber\\
& =\sqrt{2\pi}\left(  -4\right)  ^{p}\tbinom{2p}{p}^{-1}\left(  \left(
-4\right)  ^{-p}\sum\limits_{k=-p,\text{ }k\neq0}^{p}\tbinom{2p}{p+k}%
\delta\left(  t-2k\right)  \right) \nonumber\\
& =\sqrt{2\pi}\tbinom{2p}{p}^{-1}\sum\limits_{k=-p,\text{ }k\neq0}^{p}%
\tbinom{2p}{p+k}\delta\left(  t-2k\right)  ,\label{a204}%
\end{align}

and clearly $\operatorname*{supp}\theta_{1}\subset\mathbb{R}^{1}%
\setminus2\left(  -1,1\right)  $, which means that $\theta_{1}$ satisfies the
conditions \ref{a902}.

We now consider condition \ref{a907}: $w_{1}\left(  t\right)  \left\vert
\widehat{\theta_{1}}\left(  t\right)  +1\right\vert ^{2}\leq C_{\theta_{1}%
}\left(  1+t^{2}\right)  ^{n}$.

From \ref{a2.03},%
\begin{align*}
w_{1}\left(  t\right)  \left\vert \widehat{\theta_{1}}\left(  t\right)
+1\right\vert ^{2}  & \leq w_{1}\left(  t\right)  \left(  \left(  -4\right)
^{p}\tbinom{2p}{p}^{-1}\sin^{2p}t\right)  ^{2}\\
& \leq4^{2p}\tbinom{2p}{p}^{-2}w_{1}\left(  t\right)  \sin^{4p}t.
\end{align*}

In Corollary \ref{Cor_cdiffwt_bnd_on_wt_fn} it was shown that for any $a>0$
there exist constants $c_{a},c_{a}^{\prime},k_{a},k_{a}^{\prime}>0$ such that
for $1\leq n\leq l$,%
\begin{equation}
\frac{c_{a}}{t^{2\left(  l-n\right)  }}\leq w_{1}\left(  t\right)  \leq
\frac{c_{a}^{\prime}}{t^{2\left(  l-n\right)  }},\quad\left\vert t\right\vert
\leq a.\label{a921}%
\end{equation}

and%
\begin{equation}
k_{a}t^{2n}\leq w_{1}\left(  t\right)  \leq k_{a}^{\prime}t^{2n}%
,\quad\left\vert t\right\vert \geq a,\label{a920}%
\end{equation}

In fact, if we now (arbitrarily) choose $a=1$ then $\left\vert t\right\vert
\leq1$ implies
\begin{align*}
w_{1}\left(  t\right)  \left\vert \widehat{\theta_{1}}\left(  t\right)
+1\right\vert ^{2}  & \leq4^{2p}\tbinom{2p}{p}^{-2}w_{1}\left(  t\right)
\sin^{4p}t\\
& \leq4^{2p}\tbinom{2p}{p}^{-2}\frac{c_{1}^{\prime}}{t^{2\left(  l-n\right)
}}\sin^{4p}t\\
& \leq4^{2p}\tbinom{2p}{p}^{-2}c_{1}^{\prime}\frac{\sin^{2\left(  l-n\right)
}t}{t^{2\left(  l-n\right)  }}\\
& \leq4^{2p}\tbinom{2p}{p}^{-2}c_{1}^{\prime}.
\end{align*}
and $\left\vert t\right\vert \geq1$ implies%
\[
w_{1}\left(  t\right)  \left\vert \widehat{\theta_{1}}\left(  t\right)
+1\right\vert ^{2}\leq4^{2p}\tbinom{2p}{p}^{-2}k_{1}^{\prime}t^{2n}\sin
^{4p}t\leq4^{2p}\tbinom{2p}{p}^{-2}k_{1}^{\prime}t^{2n},
\]
so that \ref{a907} holds for $C_{\theta_{1}}=4^{2p}\tbinom{2p}{p}^{-2}%
\max\left\{  c_{1}^{\prime},k_{1}^{\prime}\right\}  $.\smallskip

\fbox{\textbf{Case} $n=l$} Choose $\theta_{1}=0$. Then clearly condition
\ref{a902} holds. From \ref{a921} and \ref{a920} it is clear that \ref{a907}
holds for $C_{\theta_{1}}=\max\left\{  c_{1}^{\prime},k_{1}^{\prime}\right\}
$.\smallskip

\textbf{Assumption 3 of Theorem} \ref{Thm_data_fn_tensor_prod}: From
\ref{a921} and \ref{a920}, if $\eta_{1}\geq1$,%
\[
\frac{1}{w_{1}\left(  t\right)  }\leq\left\{
\begin{array}
[c]{ll}%
\frac{1}{c_{a}}t^{2\left(  l-n\right)  }, & \left\vert t\right\vert \leq a,\\
\frac{1}{k_{a}}\frac{1}{t^{2n}}, & \left\vert t\right\vert \geq a,
\end{array}
\right.
\]

and%
\[
w_{1}\left(  t\right)  \leq\left\{
\begin{array}
[c]{ll}%
\frac{c_{b}^{\prime}}{t^{2\left(  l-n\right)  }}, & \left\vert t\right\vert
\leq b,\\
k_{b}^{\prime}t^{2n}, & \left\vert t\right\vert \geq b,
\end{array}
\right.
\]

and hence%
\[
w_{1}\left(  t/\eta_{1}\right)  \leq\left\{
\begin{array}
[c]{ll}%
c_{b}^{\prime}\eta_{1}^{2\left(  l-n\right)  }\frac{1}{t^{2\left(  l-n\right)
}}, & \left\vert t\right\vert \leq\eta_{1}b,\\
\frac{k_{b}^{\prime}}{\eta_{1}^{2n}}t^{2n}, & \left\vert t\right\vert \geq
\eta_{1}b.
\end{array}
\right.
\]

Now setting $a=\eta_{1}b$ yields%
\begin{align*}
\frac{w_{1}\left(  t/\eta_{1}\right)  }{w_{1}\left(  t\right)  }\leq\left\{
\begin{array}
[c]{ll}%
\frac{c_{b}^{\prime}}{c_{a}}, & \left\vert t\right\vert \leq a,\\
\frac{1}{\eta_{1}^{2n}}\frac{k_{b}^{\prime}}{k_{a}}, & \left\vert t\right\vert
\geq a.
\end{array}
\right.   & =\left\{
\begin{array}
[c]{ll}%
\frac{c_{a/\eta_{1}}^{\prime}}{c_{a}}, & \left\vert t\right\vert \leq a,\\
\frac{1}{\eta_{1}^{2n}}\frac{k_{a/\eta_{1}}^{\prime}}{k_{a}}, & \left\vert
t\right\vert \geq a.
\end{array}
\right. \\
& \leq\min\left\{  \frac{c_{a/\eta_{1}}^{\prime}}{c_{a}},\frac{1}{\eta
_{1}^{2n}}\frac{k_{a/\eta_{1}}^{\prime}}{k_{a}}\right\}  ,
\end{align*}

for all $t$ and $a>0$. Now set $c_{\eta_{1}}=\min\left\{  \frac{c_{a/\eta_{1}%
}^{\prime}}{c_{a}},\frac{1}{\eta_{1}^{2n}}\frac{k_{a/\eta_{1}}^{\prime}}%
{k_{a}}\right\}  $.

\textbf{Assumption} \ref{a905} is $w_{1}\left(  t\right)  \geq c_{w1}\left(
1+t^{2}\right)  ^{n}$ $a.e.$ $on$ $\mathbb{R}^{1}$.

From \ref{a921} and \ref{a920} with $a=1$ we have $w_{1}\left(  t\right)  \geq
k_{1}t^{2n}$ when $\left\vert t\right\vert \geq1$, and $w_{1}\left(  t\right)
\geq c_{1}/t^{2\left(  l-n\right)  }$ when $\left\vert t\right\vert \leq1$ and
$1\leq n\leq l$. The latter condition implies that $w_{1}\left(  t\right)
\geq c_{1}$ when $\left\vert t\right\vert \leq1$. Hence we want to show there
exists a constant $c_{w1}>0$ such that
\begin{align*}
c_{1}  & \geq c_{w1}\left(  1+t^{2}\right)  ^{n},\text{ }\left\vert
t\right\vert \leq1,\\
k_{1}t^{2n}  & \geq c_{w1}\left(  1+t^{2}\right)  ^{n},\text{ }\left\vert
t\right\vert \geq1,
\end{align*}

which is satisfied by%
\begin{align*}
c_{w1}=\min\left\{  \min_{\left\vert t\right\vert \leq1}\frac{c_{1}}{\left(
1+t^{2}\right)  ^{n}},\min_{\left\vert t\right\vert \geq1}\frac{k_{1}t^{2n}%
}{\left(  1+t^{2}\right)  ^{n}}\right\}   & =\min\left\{  \frac{c_{1}}{2^{n}%
},\frac{k_{1}}{2^{n}}\right\} \\
& =\frac{1}{2^{n}}\min\left\{  c_{1},k_{1}\right\}  .
\end{align*}

\end{proof}

The last result assumed that $\int\nolimits_{\left\vert t\right\vert \geq
R}t^{2l}q\left(  t\right)  dt<\infty$. We now weaken this assumption to

$\int\nolimits_{\left\vert t\right\vert \geq R}t^{2n}q\left(  t\right)
dt<\infty$ and see what extra assumption is required to ensure that $X_{w}%
^{0}\left(  \Omega\right)  =W^{n\mathbf{1}}\left(  \Omega\right)  $. In fact:

\begin{theorem}
\label{Thm_ex2}Suppose:

\begin{enumerate}
\item $w$ is a central difference weight function on $\mathbb{R}^{d}$ with
parameters $n,l,q$ which satisfies the condition \ref{a911} i.e.
$\int\nolimits_{\left\vert t\right\vert \geq R}t^{2n}q\left(  t\right)
dt<\infty$.

\item $\Omega$ is a bounded region and satisfies the segment condition or the
uniform rectangle property.

\item There exist constants $r,k_{1}>0$ such that
\begin{equation}
w\left(  \xi\right)  \geq\frac{k_{1}}{\xi^{2\left(  l-n\right)  }}%
,\text{\quad}\left\vert \xi\right\vert \leq r.\label{a922}%
\end{equation}

\end{enumerate}

Then for each $\lambda>0$, $X_{\sigma_{\lambda}w}^{0}\left(  \Omega\right)
=W^{n\mathbf{1}}\left(  \Omega\right)  $ as sets and the norms are equivalent.
\end{theorem}

\begin{proof}
We first verify that \ref{a920} and \ref{a921} imply that the conditions of
Theorem \ref{Thm_ex_data_fn_3} and Corollary \ref{Cor_Thm_data_fn_3} hold for
the case of one dimension. From \ref{a950} there exists a constant $C_{1}>0$
such that%
\begin{equation}
w\left(  \xi\right)  \geq\frac{\xi^{2n}}{C_{1}},\quad\xi\in\mathbb{R}%
^{1},\label{a923}%
\end{equation}

and from \ref{a958} there exist constants $c_{1},c_{2}>0$ and $0\leq a<b$ such
that%
\begin{equation}
w\left(  \xi\right)  \leq\left\{
\begin{array}
[c]{ll}%
\frac{c_{2}}{\xi^{2\left(  l-n\right)  }},\quad & \left\vert \xi\right\vert
\leq\frac{\pi}{b-a},\\
2c_{1}\xi^{2n},\quad & \left\vert \xi\right\vert >\frac{\pi}{b-a}.
\end{array}
\right. \label{a924}%
\end{equation}

Using these two estimates and the additional assumption \ref{a922} we can
derive inequalities of the form \ref{a920} and \ref{a921}, and so prove this
theorem in a similar manner to Theorem \ref{Thm_ex1}.
\end{proof}

\section{Interpolant convergence to its data
function\label{Sect_CntlDiffWtFn_InterpolConverg}}

We will use the results of \ref{Ch_Interpol} to derive orders for the
pointwise converge of the minimal norm interpolant $\mathcal{I}_{X}f$ to its
data function $f$ \textbf{in the case of central difference weight functions}.
These estimates will have the form%
\[
\left\vert f\left(  x\right)  -\left(  \mathcal{I}_{X}f\right)  \left(
x\right)  \right\vert \leq k_{G}\left\Vert f\right\Vert _{w,0}\left(
h_{X,\Omega}\right)  ^{m},\quad x\in\overline{\Omega},
\]

where the spherical cavity size $h_{X,\Omega}=\sup\limits_{\omega\in\Omega
}\operatorname*{dist}\left(  \omega,X\right)  <h_{\Omega,\kappa}$ measures the
minimum density of the independent data $X$ which is contained in the data
region $\Omega$. Here $m$ is the order of convergence.

The results of Chapter \ref{Ch_Interpol} can be classified according to
whether or not they explicitly require the independent data $X$ to be
unisolvent and use Taylor series expansions: a full discussion of unisolvency
is incorporated into Subsection \ref{Sect_unisolv}.

It turns out that the convergence orders obtained for the central difference
weight functions are the same as those obtained for the extended B-splines.

\subsection{Convergence estimates: no Taylor series expansion}

Following the approach of Chapter \ref{Ch_Interpol} we consider two types of
convergence results for the case when the data is not explicitly assumed to be
unisolvent and for which no Taylor series expansion is used. These are
imaginatively called Type 1 and Type 2.

\subsubsection{\protect\underline{Type 1 convergence estimates}}

Recall that Type 1 convergence estimates are those obtained using Theorem
\ref{Thm_|f(x)-f(y)|_inequal_2}, where there are no \textit{a priori}
assumptions made about the weight function parameter $\kappa$ but the
smoothness condition \ref{a912} is applied to the basis function near the
origin, and this will allow a uniform order of convergence estimate to be
obtained for the interpolant in a closed bounded infinite data region.

\begin{theorem}
\label{Thm_|f(x)-f(y)|_inequal_2_2}(Copy of Theorem
\ref{Thm_|f(x)-f(y)|_inequal_2}) Suppose the weight function $w$ has property
W02 and that $G$ is the basis function generated by $w$. Assume that for some
$s>0$ and constants $C_{G},h_{G}>0$ the basis function satisfies
\begin{equation}
G\left(  0\right)  -\operatorname{Re}G\left(  x\right)  \leq C_{G}\left\vert
x\right\vert ^{2s},\text{\quad}\left\vert x\right\vert <h_{G}.\label{a912}%
\end{equation}

Let $\mathcal{I}_{X}$ be the minimal norm interpolant mapping with the
independent data set $X$ contained in the closed bounded infinite set $K$, and
let $k_{G}=\left(  2\pi\right)  ^{-\frac{d}{4}}\sqrt{2C_{G}}$. Then for any
data function $f\in X_{w}^{0}$ it follows that $\sqrt{\left(  f-\mathcal{I}%
_{X}f,f\right)  _{w,0}}\leq\left\Vert f\right\Vert _{w,0}$ and%
\begin{equation}
\left\vert f\left(  x\right)  -\mathcal{I}_{X}f\left(  x\right)  \right\vert
\leq k_{G}\sqrt{\left(  f-\mathcal{I}_{X}f,f\right)  _{w,0}}\left(
h_{X,K}\right)  ^{s},\text{\quad}x\in K,\label{a927}%
\end{equation}

when $h_{X,K}=\sup\limits_{x\in K}\operatorname*{dist}\left(  x,X\right)
<h_{G}$ i.e. the order of convergence is at least $s$.
\end{theorem}

In the next result we derive constants $s$, $C_{G}$ and $h_{G}$ for which the
central difference basis function satisfies equation \ref{a912}.

\begin{theorem}
\label{Thm_converg_ext_nat_splin_wt_func_n=1}Suppose the weight function $w$
is a tensor product central difference weight function with parameters $n$,
$l$, $q\left(  \cdot\right)  $, and $q$ is bounded if $n=1$. Let $G$ be the
basis function generated by $w$

Then the basis function satisfies the estimate \ref{a912} for $s=1/2$,
$h_{G}=\infty$,%
\begin{equation}
C_{G}=G_{1}\left(  0\right)  ^{d-1}\left\Vert DG_{1}\right\Vert _{\infty}%
\sqrt{d},\nonumber
\end{equation}

and%
\[
k_{G}=\left(  2\pi\right)  ^{-\frac{d}{4}}\sqrt{2C_{G}}=\left(  2\pi\right)
^{-\frac{d}{4}}\sqrt{2G_{1}\left(  0\right)  ^{d-1}\left\Vert DG_{1}%
\right\Vert _{\infty}}\sqrt[4]{d}.
\]

\end{theorem}

\begin{proof}
Let $G_{1}$ be the univariate basis function. By Theorem
\ref{Thm_cntr_diff_basis_Lips_all_dim}, $G_{1}\in C_{B}^{\left(  0\right)  }$,
$DG_{1}$ is bounded and
\[
\left\vert G\left(  x\right)  -G\left(  y\right)  \right\vert \leq\sqrt
{d}G_{1}\left(  0\right)  ^{d-1}\left\Vert DG_{1}\right\Vert _{\infty
}\left\vert x-y\right\vert ,\text{\quad}x,y\in\mathbb{R}^{d},
\]

so that we can choose $C_{G}=G_{1}\left(  0\right)  ^{d-1}\left\Vert
DG_{1}\right\Vert _{\infty}\sqrt{d}$ and%
\[
k_{G}=\left(  2\pi\right)  ^{-\frac{d}{4}}\sqrt{2C_{G}}=\left(  2\pi\right)
^{-\frac{d}{4}}\sqrt{2G_{1}\left(  0\right)  ^{d-1}\left\Vert DG_{1}%
\right\Vert _{\infty}}\sqrt[4]{d}.
\]

\end{proof}

\subsubsection{\protect\underline{Type 2 convergence estimates}}

The Type 1 estimates of Theorem \ref{Thm_|f(x)-f(y)|_inequal_2_2} considered
the case where the weight function had property W02 for some $\kappa\geq0$ and
an extra condition was applied to the basis function. Type 2 convergence
estimates only assume $\kappa\geq1$. In the next theorem we only assume that
$\kappa\geq1$ and derive an order of convergence estimate of $1$, as well as a
doubled convergence estimate of $2$ for the Riesz data functions $R_{y}$. Note
that the constant $k_{G}$ in the next Theorem is defined to match the constant
$k_{G}$ for the Type 1 convergence estimate \ref{a927}.

\begin{theorem}
\label{Thm_interpol_error_in_terms_of_wt_fn_2}(Copy of Theorem
\ref{Thm_interpol_error_in_terms_of_wt_fn}) Suppose the weight function $w$
has property W02 for a parameter $\kappa\geq1$. Suppose also that
$\mathcal{I}_{X}f$ is the minimal interpolant on $X$ of the data function
$f\in X_{w}^{0} $ and that $K$ is a bounded closed infinite subset of
$\mathbb{R}^{d}$, $X\subset K$ and $h_{X,K}=\sup\limits_{x\in K}%
\operatorname*{dist}\left(  x,X\right)  $.

Then%
\[
\left\vert f\left(  x\right)  -\left(  \mathcal{I}_{X}f\right)  \left(
x\right)  \right\vert \leq k_{G}\sqrt{\left(  f-\mathcal{I}_{X}f,f\right)
_{w,0}}h_{X,K},\text{\quad}x\in K,
\]

where%
\[
k_{G}^{2}=\left(  2\pi\right)  ^{-d}\int\frac{\left\vert \cdot\right\vert
^{2}}{w}=-\left(  2\pi\right)  ^{-\frac{d}{2}}\left(  \left\vert D\right\vert
^{2}G\right)  \left(  0\right)  ,\quad\sqrt{\left(  f-\mathcal{I}%
_{X}f,f\right)  _{w,0}}\leq\left\Vert f\right\Vert _{w,0},
\]

and the order of convergence is at least $1$. Further, for the `Riesz' data
functions $R_{y}$, $y\in K$, we have
\[
\left\vert R_{y}\left(  x\right)  -\left(  \mathcal{I}_{X}R_{y}\right)
\left(  x\right)  \right\vert \leq k_{G}\left(  h_{X,K}\right)  ^{2}%
,\text{\quad}x,y\in K,
\]

i.e. a `doubled' estimate of at least $2$ for the rate of convergence.
\end{theorem}

Since the weight function has property W02 for $\kappa$ iff $\kappa+1/2<n$ we
can choose $\kappa=1$ iff $n\geq2$. Then by Theorem
\ref{Thm_interpol_error_in_terms_of_wt_fn_2}, $k_{G}=\left(  2\pi\right)
^{-d/4}\sqrt{-\left(  \left\vert D\right\vert ^{2}G\right)  \left(  0\right)
}$ and since
\begin{align*}
\left(  \left\vert D\right\vert ^{2}G\right)  \left(  0\right)  =\sum
\limits_{k=1}^{d}\left(  D_{k}^{2}G\right)  \left(  0\right)   &
=\sum\limits_{k=1}^{d}D_{k}^{2}\left(  G_{1}\left(  x_{1}\right)  G_{1}\left(
x_{2}\right)  \ldots G_{1}\left(  x_{d}\right)  \right)  \left(  0\right) \\
&  =\sum\limits_{k=1}^{d}G_{1}\left(  0\right)  ^{d-1}D^{2}G_{1}\left(
0\right) \\
&  =G_{1}\left(  0\right)  ^{d-1}D^{2}G_{1}\left(  0\right)  d,
\end{align*}

we get the formula%
\begin{equation}
k_{G}=\left(  2\pi\right)  ^{-d/4}\sqrt{-G_{1}\left(  0\right)  ^{d-1}%
D^{2}G_{1}\left(  0\right)  }\sqrt{d}.\label{a973}%
\end{equation}

\subsection{Convergence estimates: unisolvent data sets and $\kappa\geq1$}

In fact comparison of Theorem \ref{Thm_ex_splin_wt_fn_properties}, Chapter
\ref{Ch_wtfn_basisfn_datasp} for the extended B-splines with Theorem
\ref{Thm_cdiffwt_2} of this document shows that for parameters $n$ and $l$ the
smoothness parameters $\kappa$ can take an identical sets of values when
$\int\nolimits_{\mathbb{R}^{1}}\left\vert \cdot\right\vert ^{2n-1}q<\infty$,
since $q\in L^{1}$.

We now require the following convergence estimate from Chapter
\ref{Ch_Interpol}:

\begin{theorem}
\label{Thm_converg_interpol_ord_gte_1_2}(Copy of Theorem
\ref{Thm_converg_interpol_ord_gte_1}) Let $w$ be a weight function with
property W02 for $\kappa\geq1$ and let $G$ be the corresponding basis
function. Set $m=\left\lfloor \kappa\right\rfloor $.

Suppose $\mathcal{I}_{X}f$ is the minimal norm interpolant of the data
function $f\in X_{w}^{0}$ on the independent data set $X$ contained in the
data region $\Omega$. We use the notation and assumptions of Lemma
\ref{Lem_Lagrange_interpol} of Chapter \ref{Ch_Interpol} which means assuming
that $X$ is $m$-unisolvent and $\Omega$ is a bounded region whose boundary
satisfies the cone condition.

Now set $k_{G}=\frac{d^{m/2}}{\left(  2\pi\right)  ^{d/2}m!}\left(
c_{\Omega,\kappa}\right)  ^{m}K_{\Omega,m}^{\prime}\max\limits_{\left\vert
\beta\right\vert =m}\left\vert D^{2\beta}G\left(  0\right)  \right\vert $.
Then there exists $h_{\Omega,\kappa}>0$ such that%
\begin{equation}
\left\vert f\left(  x\right)  -\left(  \mathcal{I}_{X}f\right)  \left(
x\right)  \right\vert \leq k_{G}\left\Vert f\right\Vert _{w,0}\left(
h_{X,\Omega}\right)  ^{m},\quad x\in\overline{\Omega},\label{a934}%
\end{equation}

when $h_{X,\Omega}=\sup\limits_{\omega\in\Omega}\operatorname*{dist}\left(
\omega,X\right)  <h_{\Omega,\kappa}$ i.e. the order of convergence is at least
$m$.

The constants $c_{\Omega,m}$, $K_{\Omega,m}^{\prime}$ and $h_{\Omega,m}$ only
depend on $\Omega,m$ and $d$. In terms of the integrals which define weight
property W02 we have%
\[
\max_{\left\vert \beta\right\vert =m}\left\vert D^{2\beta}G\left(  0\right)
\right\vert \leq\left(  2\pi\right)  ^{-\frac{d}{2}}\int\tfrac{\left\vert
\xi\right\vert ^{2m}d\xi}{w\left(  \xi\right)  }.
\]

\end{theorem}

We now show that the order of convergence of the minimal norm interpolant to
an arbitrary data function is at least $n-1$:

\begin{corollary}
Suppose the weight function is a central difference weight function with
parameters $2\leq n\leq l$. Then the order of convergence of the minimal norm
interpolant to an arbitrary data function is at least $n-1$.
\end{corollary}

\begin{proof}
Theorem \ref{Thm_cdiffwt_2} implies that $m=\left\lfloor \kappa\right\rfloor
=n-1$ and Theorem \ref{Thm_converg_interpol_ord_gte_1_2} implies that the
order of convergence is at least $n-1$.
\end{proof}

\subsection{Convergence estimates using tempered distribution Taylor series}

Suppose $w=%
{\textstyle\bigotimes\limits_{i=1}^{d}}
w_{1}$ is a \textbf{tensor product central difference} weight function defined
in Corollary \ref{Cor_cdiffwt_W02_W03} with parameters $n$, $l$ and (function)
$q$. Then $w_{1} $ satisfies property W02 for parameter $\kappa$ iff $l$ and
$n$ satisfy $\kappa+1/2<n\leq l$. By Corollary
\ref{Cor_Thm_ten_prod_two_wt_fns} $w$ has property W03 for parameter
$\kappa\mathbf{1}$.

Thus the maximum value of $\left\lfloor \kappa\right\rfloor $ is $n-1$. Choose
$\kappa$ so that $\left\lfloor \kappa\right\rfloor =n-1$ and set
$m=\left\lfloor \kappa\right\rfloor $.

By Theorem \ref{Thm_G_basis_def_2} the corresponding 1-dimensional basis
function $G_{1}$ is an even function and satisfies $G_{1}\in C_{B}^{\left(
2n-2\right)  }\left(  \mathbb{R}^{1}\right)  $. If we further assume that
$q$\textbf{\ is bounded} a.e. then by Theorem \ref{Thm2_G_basis_def_2_new},
$D^{2n-1}G_{1}$ is essentially bounded.

We can now easily apply the data function remainder result Theorem
\ref{Thm_rem_estim_ten_prod} of Section \ref{Sect_rem_estim_ten_prod} for
tensor product weight functions which have property W03 for $\kappa\mathbf{1}$
and for which $\left\Vert D^{2m+1}G_{1}\right\Vert _{\infty,B_{r}}<\infty$ for
some $r>0$ where $m=\left\lfloor \kappa\right\rfloor $. The result is a
convergence rate of $m+1/2=n-1/2$ which is an improvement.

\subsection{Summary table}%

\begin{table}[htbp] \centering
$%
\begin{tabular}
[c]{|c||c||c|c|c|}\hline
\multicolumn{5}{|c|}{\textbf{Interpolant convergence to data func:}
\textbf{tensor prod. central diff. weight funcs}.}\\\hline
Estimate & Parameter & Converg. &  & \\
name & constraints & order & $\left(  2\pi\right)  ^{d/4}k_{G}$ & $h_{G}%
$\\\hline\hline
\multicolumn{1}{|l|}{Type 1 (smoothness} & $q$ b'nded if $n=1$ & $1/2$ &
$\sqrt{2G_{1}\left(  0\right)  ^{d-1}\left\Vert DG_{1}\right\Vert _{\infty}%
}\sqrt[4]{d}$ $^{\left(  1\right)  }$ & $\infty$\\
\multicolumn{1}{|l|}{\quad constraint on basis fn)} & else $n\geq2$ &  &  &
\\\hline
\multicolumn{1}{|l|}{Type 2 ($\kappa\geq1$)} & $n\geq2$ & $1$ & $\sqrt
{-G_{1}\left(  0\right)  ^{d-1}D^{2}G_{1}\left(  0\right)  }\sqrt{d}$
$^{\left(  1\right)  }$ & $\infty$\\\hline
Unisolvency/Taylor series & $n\geq2$ & $n-1$ & - & $\infty$\\\hline
\multicolumn{1}{|l|}{Tempered distrib Taylor} & $q$ b'nded a.e. & $n-1/2$ &
- & $\infty$\\\hline\hline
\multicolumn{5}{|l|}{$^{\left(  1\right)  }${\small \ }$G_{1}${\small \ is the
univariate basis function used to form the tensor product.}}\\\hline
\end{tabular}
$\caption{}\label{Tbl_ConvergCentral}%
\end{table}%

\chapter{Central difference weight functions: $q$
multivariate\label{Ch_cent_diff_wt_fn_multivar}}

\section{Introduction}

In the last chapter we generated multivariate central difference weight
functions using the tensor product. In this chapter we generalize to higher
dimensions the 1-dimensional difference operator $\Delta_{2l}$ used in the
one-dimensional case and thus obtain a new class of weight functions which
seems most suited to radial weight functions. This chapter is unfinished.

Using the Fourier transform to directly calculate the basis function from the
weight function formula \ref{a969} may be very awkward.

In \textbf{Section} \ref{Sect_convol_centdiff_basis} we use the theory of
$S_{\emptyset,m}$ tempered distributions (see Definition \ref{Def_So,n} etc.)
and in Theorem \ref{Thm_CentDiff_eqn_thinps} shows us how to use the central
difference operator to express the basis function as the sum of the thin-plate
spline and several scaled convolutions of this spline with the generating
function $q$. The use of a $S_{\emptyset,m}$ technique means that the formula
is true modulo a polynomial. It remains to determine this polynomial by
applying conditions to the function $q$ used to define the weight function.
For $d$ odd, using the $L_{loc}^{1}$ Taylor series expansion with integral
remainder described in Lemma \ref{Lem_Taylor_estim_loc_C_L1loc}, I have
derived Theorem \ref{Thm_CentDiff_eqn_thinps} which supplies conditions on $q$
under which the polynomial $p_{c}$ is zero. Here we assume $q$ has bounded
support or that $q $ is radial with an extra condition. I have not done the
case where $d$ is even but have made partial attempt on a more general result
which does not distinguish between odd and even dimensions.

In \textbf{Section} \ref{Sect_Bessel_centdiff_basis_formula} we will derive
several formulas for the central difference basis function which involve the
Bessel functions $J_{\frac{d-2}{2}}$ and the MacDonald's functions
$K_{n-\frac{d}{2}}$.

In \textbf{Section} \ref{Sect_CentDiffBasis_Taylor} we obtain results which
are an extension of the 1-dimensional of Subsection
\ref{SbSect_tempdistrib_1_dim_centdiff_basis} to an arbitrary dimension.

We use the tempered distribution Taylor series expansion introduced in Section
\ref{Sect_Taylor_series_data_fn} and the theory of the Schwartz subspaces
$S_{\emptyset,k}\subset S$ introduced in Definition \ref{Def_So,n} to expand a
convolution of the thin plate spline and the function $q$ and thus prove the
Taylor series basis function formula given in Theorem
\ref{Thm_centdiff_basis_Tn_q_suppbnd}. The formula is true modulo a polynomial
with degree $\leq d-2$. The $L_{loc}^{1}$ Taylor series expansion of Lemma
\ref{Lem_Taylor_estim_L1loc_L1} is then used to obtain the Taylor series
integral remainder formula \ref{a017} for the basis function which is again
true modulo a polynomial with degree $\leq d-2$.

In Section \ref{Sect_LocDataSpace} I have made a start at characterizing the
data space locally.??

\begin{remark}
?? \textbf{Things to do}:

\begin{enumerate}
\item Does there exist a $C_{0}^{\infty}$ basis function?

\item Derive upper/lower bounds for the weight function.

\item Also consider $d$ even.

\item Calculate interpolant/smoother convergence order.

\item Calculate local data space.
\end{enumerate}
\end{remark}

\section{The weight function\label{Sect_CentDiffWeightFunc_MultiVar}}

In this section I generalize the $\mathbb{R}^{1}$ central difference weight
function of Definition \ref{Def_central_diff_wt_func} to $\mathbb{R}^{d}$
which in turn involves generalizing the central difference operator
$\Delta_{2l}$ of Definition \ref{Def_op_diff_real} to $\mathbb{R}^{d}$.

More precisely, in Definition \ref{Def_op_diff_real} a function $w$ is defined
using the parameters $n,l\geq0$, $q\geq0$ and in Theorem
\ref{Thm_wt_func_mdim} necessary and sufficient conditions are supplied for
$w$ to be a weight function with properties W01 and W02.

The motivation for the 1-dimensional definition of a central difference weight
function, which is given at the start of Subsection
\ref{SbSect_cent_Motivation}, can be easily adapted to the multidimensional
case studied here.

\begin{definition}
\label{Def_op_diff_multivar}\textbf{The multivariate central difference
operators }$\delta_{\xi}$ \textbf{and} $\Delta_{2l}$ \textbf{on}
$\mathbb{R}^{d}$.

We define the central difference operator%
\[
\delta_{\xi}f\left(  x\right)  =f\left(  x+\frac{\xi}{2}\right)  -f\left(
x-\frac{\xi}{2}\right)  ,\text{\quad}x,\xi\in\mathbb{R}^{d}.
\]

and introduce the operator%
\begin{equation}
\Delta_{2l}f\left(  \xi\right)  =\sum_{j=-l}^{l}\left(  -1\right)  ^{j}%
\tbinom{2l}{j+l}f\left(  -j\xi\right)  ,\text{\quad}l=1,2,3,\ldots,\text{ }%
\xi\in\mathbb{R}^{d}.\label{a9.2}%
\end{equation}

and note that $\Delta_{2l}f\left(  -\xi\right)  =\Delta_{2l}f\left(
\xi\right)  $.

This definition of $\Delta_{2l}$ ensures that $\Delta_{2l,\xi}\left(
e^{-i\xi\tau}\right)  \geq0$, which is part 3 of the next lemma. For example:%
\begin{align*}
\Delta_{2}f\left(  \xi\right)   & =-\left(  f\left(  \xi\right)  -2f\left(
0\right)  +f\left(  -\xi\right)  \right)  .\\
\Delta_{4}f\left(  \xi\right)   & =f\left(  2\xi\right)  -4f\left(
\xi\right)  +6f\left(  0\right)  -4f\left(  -\xi\right)  +f\left(
-2\xi\right)  .
\end{align*}

Noting that $\Delta_{2l}f\left(  \xi\right)  =\left(  -1\right)  ^{l}\left(
\delta_{\xi}^{2l}f\right)  \left(  0\right)  $ we define%
\[
\Delta_{k}f\left(  \xi\right)  :=i^{k}\left(  \delta_{\xi}^{k}f\right)
\left(  0\right)  ,\text{\quad}k=0,1,2,\ldots
\]

The operator $\delta_{\xi}f$ can be defined on distributions by%
\[
\delta_{\xi}f=f\left(  \cdot+\frac{\xi}{2}\right)  -f\left(  \cdot-\frac{\xi
}{2}\right)  ,\text{\quad}f\in\mathcal{D}^{\prime},\text{ }\xi\in
\mathbb{R}^{d}.
\]

\end{definition}

\begin{lemma}
\label{Lem_central_diff_op multivar}The \textbf{even} central difference
operators of Definition \ref{Def_op_diff_real} have the following properties:

\begin{enumerate}
\item $\Delta_{2k}f\left(  \xi\right)  =\left(  -1\right)  ^{k}\left(
\delta_{\xi}^{2k}f\right)  \left(  0\right)  =\left(  \left(  -\delta_{\xi
}\right)  ^{2k}f\right)  \left(  0\right)  $.

\item Regarding monomials and polynomials:

\begin{enumerate}
\item $\Delta_{2l}\left(  \xi^{\alpha}\right)  =\left(  \sum\limits_{j=-l}%
^{l}\left(  -1\right)  ^{j}\tbinom{2l}{j+l}j^{\left\vert \alpha\right\vert
}\right)  \left(  -\xi\right)  ^{\alpha}$ when $\alpha\geq\mathbf{0}$.

\item $\sum\limits_{j=-l}^{l}\left(  -1\right)  ^{j}\tbinom{2l}{j+l}%
j^{m}=\left\{
\begin{array}
[c]{ll}%
0, & m=0,1,2,\ldots,2l-1,\\
\left(  -1\right)  ^{l}\left(  2l\right)  !, & m=2l.
\end{array}
\right.  $

\item $\Delta_{2l}\left(  \xi^{2\alpha+1}\right)  =0$ when $\left\vert
\alpha\right\vert \geq l$.

\item $\Delta_{2l}\left(  p\left(  \xi\right)  \right)  =0$ when $p$ is a
polynomial with degree less than $2l$.

\item $\Delta_{2l}\left(  \xi^{2\alpha}\right)  =\left(  -1\right)
^{\left\vert \alpha\right\vert }\left(  2\alpha\right)  !\xi^{2\alpha}$ when
$\left\vert \alpha\right\vert =l$.
\end{enumerate}

\item $\Delta_{2l,\xi}\left(  e^{-i\xi\tau}\right)  =2^{2l}\sin^{2l}\left(
\xi\tau/2\right)  $. Here $\Delta_{2l,\xi}$ indicates $\Delta_{2l}$ is acting
on the variable $\xi$ in $e^{-i\xi\tau}$, and $\xi\tau$ is the inner product
$\sum_{i=1}^{d}\xi_{i}\tau_{i}$.

\item If $f$ is radial then $\Delta_{2l}f$ is also radial. In fact, if
$f\left(  \xi\right)  =f_{\odot}\left(  \left\vert \xi\right\vert \right)  $
then%
\begin{align*}
\Delta_{2l}f\left(  \xi\right)   & =\sum_{j=-l}^{l}\left(  -1\right)
^{j}\tbinom{2l}{j+l}f_{\odot}\left(  \left\vert j\right\vert \left\vert
\xi\right\vert \right)  =\tbinom{2l}{l}f_{\odot}\left(  0\right)  +2\sum
_{j=1}^{l}\left(  -1\right)  ^{j}\tbinom{2l}{j+l}f_{\odot}\left(  j\left\vert
\xi\right\vert \right)  .\\
\left(  \Delta_{2l}f\right)  _{\odot}  & =\Delta_{2l}\left(  f_{\odot}\right)
.
\end{align*}

\item If $\phi\in S$ then $\Delta_{l}\phi=i^{l}\left[  \delta,\delta_{\xi}%
^{l}\phi\right]  $.

\item If $\phi\in S$ and $\psi\in S_{1;2l}$ then%
\[
\Delta_{2l}\phi-\tbinom{2l}{l}\phi\left(  0\right)  \psi\in S_{\emptyset;2l}.
\]

\item If $f\in\mathcal{D}^{\prime}$ then%
\begin{align*}
\left[  \delta_{\xi}f,\phi\right]   & =\left[  f\left(  \cdot+\frac{\xi}%
{2}\right)  -f\left(  \cdot-\frac{\xi}{2}\right)  ,\phi\right]  =\left[
f,\phi\left(  \cdot-\frac{\xi}{2}\right)  -\phi\left(  \cdot+\frac{\xi}%
{2}\right)  \right]  =-\left[  f,\delta_{\xi}\phi\right]  .\\
& \Longrightarrow\\
\left[  \delta_{\xi}^{l}f,\phi\right]   & =\left(  -1\right)  ^{l}\left[
f,\delta_{\xi}^{l}\phi\right]  .
\end{align*}

Thus part 5 implies $\Delta_{l}\phi=i^{l}\left[  \delta,\delta_{\xi}^{l}%
\phi\right]  =\left(  -i\right)  ^{l}\left[  \delta_{\xi}^{l}\delta
,\phi\right]  $ and so%
\begin{equation}
\Delta_{l}=\left(  -i\right)  ^{l}\delta_{\xi}^{l}\delta\text{ }on\text{
}S.\label{a915}%
\end{equation}

\item If $f\in S^{\prime}$, then $\widehat{\delta_{\xi}f}=e^{i\frac{\xi\eta
}{2}}\widehat{f}-e^{-i\frac{\xi\eta}{2}}\widehat{f}=2i\sin\frac{\xi\eta}{2}$
$\widehat{f}$ and $\delta_{\xi}\widehat{f}=-2i\left(  \sin\frac{\xi\eta}%
{2}f\right)  ^{\wedge}$.
\end{enumerate}
\end{lemma}

The motivation for the 1-dimensional definition of a central difference weight
function, which is given at the start of Subsection
\ref{SbSect_cent_Motivation}, can be easily adapted to an arbitrary dimension
as follows. The derivation of the 1-dimensional central difference weight
function was based on the approximation of the reciprocal of the B-spline
weight function $w_{s}$ using a mollifier. From \ref{1.032} we have for
parameters $n$ and $l$,
\[
\frac{1}{w_{s}\left(  \xi\right)  }=\prod_{j=1}^{d}\frac{\sin^{2l}\xi_{j}}%
{\xi_{j}^{2n}},\quad1\leq n\leq l.
\]

Now $\sin^{2l}\xi_{j}=\left(  \frac{e^{i2\xi_{j}}-e^{-i2\xi_{j}}}{2i}\right)
^{2l}$ so that by part 3 Lemma \ref{Lem_central_diff_op multivar} with $t=2$,
$\frac{1}{w_{s}\left(  \xi\right)  }=\frac{1}{2^{2l}}\frac{\Delta_{2l}\left(
e^{-i2\xi\mathbf{1}}\right)  }{\xi^{2n}}$ or using the Fourier transform
$F\left[  \cdot\right]  \left(  \xi\right)  $
\[
\frac{1}{w_{s}\left(  \xi\right)  }=\frac{\left(  2\pi\right)  ^{1/2}}{2^{2l}%
}\frac{\Delta_{2l}F\left[  \delta\left(  x-2\right)  \right]  }{\xi^{2n}}.
\]

We now want to approximate the delta function $\delta$ in the distribution
sense using the mollifier $k\psi\left(  kx\right)  $ where $k=1,2,3,\ldots$,
$\psi\in L^{1}\left(  \mathbb{R}^{d}\right)  $, $\psi\geq0$ and $\int\psi=1$.
Since $k\psi\left(  kx\right)  \rightarrow\delta$ as $k\rightarrow\infty$ in
the sense of distributions, $k\psi\left(  k\left(  x-2\right)  \right)
\rightarrow\delta\left(  x-2\right)  $ and so as $k\rightarrow\infty$,
\[
\frac{\sqrt{2\pi}}{2^{2l}}\frac{\Delta_{2l}F\left[  k\psi\left(  k\left(
x-2\right)  \right)  \right]  }{\xi^{2n}}\rightarrow\frac{1}{w_{s}\left(
\xi\right)  }\text{ }on\text{ }\mathcal{D}^{\prime}\left(  \mathbb{R}%
^{1}\setminus\mathbf{0}\right)  .
\]

Here we will not try to extend the convergence to $S^{\prime}$ but observe
that we can write%
\[
\frac{\sqrt{2\pi}}{2^{2l}}\frac{\Delta_{2l}F\left[  k\psi\left(  k\left(
x-2\right)  \right)  \right]  }{\xi^{2n}}=\frac{\Delta_{2l}\widehat{q_{k}}%
}{\xi^{2n}},
\]

where
\[
q_{k}\left(  x\right)  :=\frac{\sqrt{2\pi}}{2^{2l}}k\psi\left(  k\left(
x-2\right)  \right)  \in L^{1}\left(  \mathbb{R}^{d}\right)  \text{ }and\text{
}q_{k}\geq0.
\]

Thus $\frac{\Delta_{2l}\widehat{q_{k}}}{\xi^{2n}}\rightarrow\frac{1}%
{w_{s}\left(  \xi\right)  }$ on $\mathcal{D}^{\prime}\left(  \mathbb{R}%
^{d}\setminus\mathbf{0}\right)  $ and if we are lucky $\frac{\xi^{2n}}%
{\Delta_{2l}\widehat{q_{k}}\left(  \xi\right)  }$ may be a weight function.
Indeed, we use this idea to define a central difference weight function and
justify calling it a weight function using the next theorem.

?? Later in ?? we will ?? show? ?? that if an extra condition is imposed on
the function $\psi$ then $\frac{\Delta_{2l}\widehat{q_{k}}}{\xi^{2n}%
}\rightarrow\frac{1}{w_{s}\left(  \xi\right)  }$ in $L^{1}$ and as tempered
distributions, and that the corresponding sequence of central difference basis
functions converges uniformly pointwise to the extended B-spline basis function.

\begin{definition}
\label{Def_multi_centdiff_wt_fn}\textbf{Multivariate central difference weight
functions}.

Suppose that $q\in L^{1}\left(  \mathbb{R}^{d}\right)  $, $q\neq0$, $q\left(
\xi\right)  \geq0$ and $l,n\geq0$ are integers. The multivariate central
difference weight function is defined by
\begin{equation}
w\left(  \xi\right)  =\frac{\left\vert \xi\right\vert ^{2n}}{\Delta
_{2l}\widehat{q}\left(  \xi\right)  },\quad\xi\in\mathbb{R}^{d}.\label{a969}%
\end{equation}

Since $q$ radial implies $\widehat{q}$ radial, part 4 of Lemma
\ref{Lem_central_diff_op multivar} implies that if $q$ is radial then $w$ is
radial and consequently the basis function is radial.

The binomial identities 2b of Lemma \ref{Lem_central_diff_op multivar} will
play a central role in the theory of these basis functions.
\end{definition}

\begin{lemma}
\label{Lem_exists_integ_sin(x1)_|x|}Suppose $d$, $n$, $l$ are positive
integers. Then%
\begin{equation}
\int_{\mathbb{R}^{d}}\frac{\sin^{2l}\eta_{1}d\eta}{\left\vert \eta\right\vert
^{2\left(  n-\lambda\right)  }}<\infty\Longleftrightarrow n-l-\frac{d}%
{2}<\lambda<n-\frac{d}{2},\label{a984}%
\end{equation}

and%
\begin{equation}
\int_{\mathbb{R}^{d}}\frac{\sin^{2l}\eta_{1}d\eta}{\left\vert \eta\right\vert
^{2\left(  n-\lambda\right)  }}=\left\{
\begin{array}
[c]{l}%
\omega_{d-1}B\left(  \frac{d-1}{2},n-\lambda-\frac{d-1}{2}\right)  \int%
_{0}^{\infty}\frac{\sin^{2l}t}{t^{2\left(  n-\lambda\right)  -d+1}}dt,\\
2\pi^{\frac{d-1}{2}}\frac{\Gamma\left(  n-\lambda-\frac{d-1}{2}\right)
}{\Gamma\left(  n-\lambda\right)  }\int_{0}^{\infty}\frac{\sin^{2l}%
t}{t^{2\left(  n-\lambda\right)  -d+1}}dt,
\end{array}
\right. \label{a800}%
\end{equation}

where $B$ denotes the beta function.
\end{lemma}

\begin{proof}
We use Theorem \ref{Thm_Integ_u(xy)f(|x|)dx}: first \ref{Ap016} then
\ref{Ap001}. If $e^{\left(  1\right)  }=\left(  1,0,\ldots,0\right)  $ then%
\begin{align}
&  \int_{r\leq\left\vert \eta\right\vert \leq R}\frac{\sin^{2l}\eta_{1}%
}{\left\vert \eta\right\vert ^{2\left(  n-\lambda\right)  }}d\eta\nonumber\\
&  =\int_{r\leq\left\vert \eta\right\vert \leq R}\frac{\sin^{2l}\left(
e^{\left(  1\right)  }\eta\right)  }{\left\vert \eta\right\vert ^{2\left(
n-\lambda\right)  }}d\eta=\int\limits_{r\leq\left\vert x\right\vert \leq
R}\frac{\sin^{2l}\left(  \left\vert e^{\left(  1\right)  }\right\vert \eta
_{1}\right)  }{\left\vert \eta\right\vert ^{2\left(  n-\lambda\right)  }}%
d\eta=\int\limits_{r\leq\left\vert x\right\vert \leq R}\frac{\sin^{2l}\eta
_{1}}{\left\vert \eta\right\vert ^{2\left(  n-\lambda\right)  }}%
d\eta=\nonumber\\
&  =\omega_{d-1}\int_{r}^{R}\left(  \int_{0}^{\pi}\sin^{2l}\left(  \rho
\cos\theta\right)  \sin^{d-2}\theta d\theta\right)  \frac{\rho^{d-1}}%
{\rho^{2\left(  n-\lambda\right)  }}d\rho\nonumber\\
&  =\omega_{d-1}\int_{r}^{R}\left(  \int_{0}^{\pi}\sin^{2l}\left(  \rho
\cos\theta\right)  \sin^{d-2}\theta d\theta\right)  \frac{d\rho}%
{\rho^{2\left(  n-\lambda\right)  -d+1}}\nonumber\\
&  =2\omega_{d-1}\int_{r}^{R}\left(  \int_{0}^{\pi/2}\sin^{2l}\left(  \rho
\cos\theta\right)  \sin^{d-2}\theta d\theta\right)  \frac{d\rho}%
{\rho^{2\left(  n-\lambda\right)  -d+1}}\nonumber\\
&  =2\omega_{d-1}\int_{0}^{\pi/2}\int_{r}^{R}\frac{\sin^{2l}\left(  \rho
\cos\theta\right)  }{\rho^{2\left(  n-\lambda\right)  -d+1}}d\rho\text{ }%
\sin^{d-2}\theta d\theta\nonumber\\
&  =2\omega_{d-1}\lim_{\varepsilon\rightarrow0^{+}}\int_{0}^{\frac{\pi}%
{2}-\varepsilon}\int_{r}^{R}\frac{\sin^{2l}\left(  \rho\cos\theta\right)
}{\rho^{2\left(  n-\lambda\right)  -d+1}}d\rho\text{ }\sin^{d-2}\theta
d\theta.\label{a985}%
\end{align}

The change of variables $t=\rho\cos\theta$, $dt=\cos\theta d\rho$ now gives%
\[
\int_{r}^{R}\frac{\sin^{2l}\left(  \rho\cos\theta\right)  }{\rho^{2\left(
n-\lambda\right)  -d+1}}d\rho=\cos^{2\left(  n-\lambda\right)  -d}\theta
\int_{r\cos\theta}^{R\cos\theta}\frac{\sin^{2l}t}{t^{2\left(  n-\lambda
\right)  -d+1}}dt,
\]

so that%
\[
\int_{0}^{\frac{\pi}{2}-\varepsilon}\int_{r}^{R}\frac{\sin^{2l}\left(
\rho\cos\theta\right)  }{\rho^{2\left(  n-\lambda\right)  -d+1}}d\rho\text{
}\sin^{d-2}\theta d\theta=\int_{0}^{\frac{\pi}{2}-\varepsilon}\cos^{2\left(
n-\lambda\right)  -d}\theta\sin^{d-2}\theta\int_{r\cos\theta}^{R\cos\theta
}\frac{\sin^{2l}t}{t^{2\left(  n-\lambda\right)  -d+1}}dtd\theta,
\]

and hence%
\begin{equation}
\int_{r\leq\left\vert \eta\right\vert \leq R}\frac{\sin^{2l}\eta_{1}%
}{\left\vert \eta\right\vert ^{2\left(  n-\lambda\right)  }}d\eta
=2\omega_{d-1}\lim_{\varepsilon\rightarrow0^{+}}\int_{0}^{\frac{\pi}%
{2}-\varepsilon}\cos^{2\left(  n-\lambda\right)  -d}\theta\sin^{d-2}\theta
\int_{r\cos\theta}^{R\cos\theta}\frac{\sin^{2l}t}{t^{2\left(  n-\lambda
\right)  -d+1}}dtd\theta.\label{a986}%
\end{equation}

Next observe that $\int\frac{\sin^{2l}t}{t^{2\left(  n-\lambda\right)  -d+1}%
}dt<\infty$ iff condition \ref{a984} holds. Thus RHS of \ref{a984} implies%
\begin{align}
\int_{r\leq\left\vert \eta\right\vert \leq R}\frac{\sin^{2l}\eta_{1}%
}{\left\vert \eta\right\vert ^{2\left(  n-\lambda\right)  }}d\eta &
<2\omega_{d-1}\lim_{\varepsilon\rightarrow0^{+}}\int_{0}^{\frac{\pi}%
{2}-\varepsilon}\cos^{2\left(  n-\lambda\right)  -d}\theta\sin^{d-2}\theta
\int\frac{\sin^{2l}t}{t^{2\left(  n-\lambda\right)  -d+1}}dtd\theta\nonumber\\
&  =2\omega_{d-1}\int_{0}^{\frac{\pi}{2}}\cos^{2\left(  n-\lambda\right)
-d}\theta\sin^{d-2}\theta d\theta\int\frac{\sin^{2l}t}{t^{2\left(
n-\lambda\right)  -d+1}}dt\nonumber\\
&  =\omega_{d-1}B\left(  \frac{d-1}{2},n-\lambda-\frac{d-1}{2}\right)
\int\frac{\sin^{2l}t}{t^{2\left(  n-\lambda\right)  -d+1}}dt,\label{a987}%
\end{align}

which implies%
\[
\int_{\mathbb{R}^{d}}\frac{\sin^{2l}\eta_{1}d\eta}{\left\vert \eta\right\vert
^{2\left(  n-\lambda\right)  }}<\omega_{d-1}B\left(  \frac{d-1}{2}%
,n-\lambda-\frac{d-1}{2}\right)  \int\frac{\sin^{2l}t}{t^{2\left(
n-\lambda\right)  -d+1}}dt<\infty.
\]

On the other hand, suppose that $\int_{\mathbb{R}^{d}}\frac{\sin^{2l}\eta
_{1}d\eta}{\left\vert \eta\right\vert ^{2\left(  n-\lambda\right)  }}$ exists.
Fix $0<\varepsilon^{\prime}<\pi/2$ and choose $R>r>0$ such that $R\sin
\varepsilon^{\prime}>r$. Then from \ref{a986},%
\[
\int_{r\leq\left\vert \eta\right\vert \leq R}\frac{\sin^{2l}\eta_{1}d\eta
}{\left\vert \eta\right\vert ^{2\left(  n-\lambda\right)  }}>2\omega_{d-1}%
\int_{0}^{\frac{\pi}{2}-\varepsilon^{\prime}}\cos^{2\left(  n-\lambda\right)
-d}\theta\sin^{d-2}\theta\int_{r}^{R\sin\varepsilon^{\prime}}\frac{\sin^{2l}%
t}{t^{2\left(  n-\lambda\right)  -d+1}}dtd\theta.
\]

Hence $\forall$ $0<\varepsilon^{\prime}<\pi/2$,%
\[
\int\frac{\sin^{2l}\eta_{1}d\eta}{\left\vert \eta\right\vert ^{2\left(
n-\lambda\right)  }}>2\omega_{d-1}\int_{0}^{\frac{\pi}{2}-\varepsilon^{\prime
}}\cos^{2\left(  n-\lambda\right)  -d}\theta\sin^{d-2}\theta d\theta\int%
\frac{\sin^{2l}t}{t^{2\left(  n-\lambda\right)  -d+1}}dt,
\]

and thus%
\begin{align}
\int\frac{\sin^{2l}\eta_{1}d\eta}{\left\vert \eta\right\vert ^{2\left(
n-\lambda\right)  }}  & >2\omega_{d-1}\int_{0}^{\frac{\pi}{2}}\cos^{2\left(
n-\lambda\right)  -d}\theta\sin^{d-2}\theta d\theta\int\frac{\sin^{2l}%
t}{t^{2\left(  n-\lambda\right)  -d+1}}dt\nonumber\\
& =\omega_{d-1}B\left(  \frac{d-1}{2},n-\lambda-\frac{d-1}{2}\right)
\int\frac{\sin^{2l}t}{t^{2\left(  n-\lambda\right)  -d+1}}dt.\label{a988}%
\end{align}

Consequently the RHS of \ref{a984} is true and so \ref{a984} is true. Further
inequalities \ref{a987} and \ref{a988} imply equality and so the proof is complete.
\end{proof}

Part 2 of the next theorem gives an easy necessary and sufficient for the
parameters $n,l\geq0$, $q\geq0$ to define a central difference weight function
\ref{a969}.

\begin{theorem}
\label{Thm_wt_func_mdim}Suppose $w$ is the function on $\mathbb{R}^{d}$
introduced in Definition \ref{Def_multi_centdiff_wt_fn}. Then:

\begin{enumerate}
\item $w$ is an even function satisfying weight function property W01 for
$\mathcal{A}=\left\{  0\right\}  $.

\item $w$ satisfies property W02 for parameter $\kappa>0$ iff
\begin{equation}
0<\kappa<n-\frac{d}{2}<l,\label{a966}%
\end{equation}

and for some $R>0$,
\begin{equation}
\int_{\left\vert \tau\right\vert \geq R}\left\vert \tau\right\vert
^{2n-d}q\left(  \tau\right)  d\tau<\infty.\label{a965}%
\end{equation}

Note \ref{a953}.\medskip

Now suppose $w$ satisfies property W02 for parameter $\kappa>0$. Then:\medskip

\item If $0\leq\lambda\leq\kappa$,%
\begin{align}
\int\frac{\left\vert \xi\right\vert ^{2\lambda}}{w\left(  \xi\right)  }d\xi &
=\tfrac{2^{2\left(  l-n+\lambda\right)  +d}}{\left(  2\pi\right)  ^{d/2}%
}\left(  \int\frac{\sin^{2l}\eta_{1}}{\left\vert \eta\right\vert ^{2\left(
n-\lambda\right)  }}d\eta\right)  \int\left\vert \tau\right\vert ^{2\left(
n-\lambda\right)  -d}q\left(  \tau\right)  d\tau\label{a997}\\
& =\tfrac{2^{2\left(  l-n+\lambda\right)  +d}}{\left(  2\pi\right)  ^{d/2}%
}\omega_{d-1}B\left(  \tfrac{d-1}{2},n-\lambda-\tfrac{d-1}{2}\right)
\int\limits_{0}^{\infty}\frac{\sin^{2l}t}{t^{2\left(  n-\lambda\right)  -d+1}%
}dt\int\left\vert \tau\right\vert ^{2\left(  n-\lambda\right)  -d}q\left(
\tau\right)  d\tau,\label{a801}%
\end{align}

where $2\left(  n-\lambda\right)  -d>0$, $2\left(  l-n+\lambda\right)  +d>0$
and $B$ denotes the beta function.

\item The basis function has at least smoothness $2n-d-1$.

\item ??! $\widehat{q}\in C_{B}^{\left(  2n-d\right)  }$ and $\Delta
_{2l}\widehat{q}\in C_{B}^{\left(  2n-d\right)  }\cap C_{\emptyset
,2n-d+1}^{\left(  2n-d\right)  }$ where in general%
\begin{equation}
C_{\emptyset,k}^{\left(  m\right)  }:=\left\{  u\in C^{\left(  m\right)
}:D^{\alpha}u\left(  0\right)  =0\text{ }for\text{ }\left\vert \alpha
\right\vert <m\right\}  ,\quad k\leq m+1,\text{ }m\geq0.\label{a8.12}%
\end{equation}

(cf. Definition \ref{Def_So,n}).
\end{enumerate}
\end{theorem}

\begin{proof}
\textbf{Part 1} The function $w$ is an even function since it is clear from
equation \ref{a9.2} that $\Delta_{2l}\widehat{q}$ is even.\medskip

\textbf{Part 2} From part 3 of Lemma \ref{Lem_central_diff_op multivar} we
have $\Delta_{2l,\xi}e^{-i\xi\tau}=2^{2l}\sin^{2l}\left(  \xi\tau/2\right)  $
so that
\begin{align}
\frac{1}{w\left(  \xi\right)  }=\frac{\Delta_{2l}\widehat{q}\left(
\xi\right)  }{\left\vert \xi\right\vert ^{2n}} &  =\tfrac{1}{\left(
2\pi\right)  ^{d/2}}\frac{1}{\left\vert \xi\right\vert ^{2n}}\int\left(
\Delta_{2l,\xi}e^{-i\xi\tau}\right)  q\left(  \tau\right)  d\tau\nonumber\\
&  =\tfrac{2^{2l}}{\left(  2\pi\right)  ^{d/2}}\frac{1}{\left\vert
\xi\right\vert ^{2n}}\int\sin^{2l}\left(  \frac{\xi\tau}{2}\right)  q\left(
\tau\right)  d\tau.\label{a953}%
\end{align}

and clearly $w\left(  \xi\right)  >0$ when $\xi\neq0$. Since $q\in L^{1}$
implies $\widehat{q}\in C_{B}^{\left(  0\right)  }$, from the definition of
$w$ we have $w\in C^{\left(  0\right)  }\left(  \mathbb{R}^{1}\setminus
0\right)  $ and $w>0$ on $\mathbb{R}^{1}\setminus0$.

Property W02 was introduced in \ref{1.036}: $w\in W02$ if for some real number
$\kappa\in\mathbb{R}^{1}$, $\kappa>0$ we have%
\begin{equation}
\int\frac{\left\vert \xi\right\vert ^{2\lambda}d\xi}{w\left(  \xi\right)
}<\infty,\quad0\leq\lambda\leq\kappa.\label{a993}%
\end{equation}

Now
\begin{align}
\int\frac{\left\vert \xi\right\vert ^{2\lambda}d\xi}{w\left(  \xi\right)  }  &
=\tfrac{2^{2l}}{\left(  2\pi\right)  ^{d/2}}\int\frac{1}{\left\vert
\xi\right\vert ^{2\left(  n-\lambda\right)  }}\int\sin^{2l}\left(  \xi
\tau/2\right)  q\left(  \tau\right)  d\tau d\xi\nonumber\\
& =\tfrac{2^{2l}}{\left(  2\pi\right)  ^{d/2}}\int\int\frac{\sin^{2l}\left(
\xi\tau/2\right)  }{\left\vert \xi\right\vert ^{2\left(  n-\lambda\right)  }%
}d\xi\,q\left(  \tau\right)  d\tau,\label{a903}%
\end{align}

where, because the integrand is non-negative, the integrals all exist iff the
last integral exists, and the last integral exists iff%
\begin{equation}
\left.
\begin{array}
[c]{l}%
\int\frac{\sin^{2l}\left(  \xi\tau/2\right)  }{\left\vert \xi\right\vert
^{2\left(  n-\lambda\right)  }}d\xi<\infty,\text{ }for\text{ }almost\text{
}all\text{ }\tau,\\
\int\int\frac{\sin^{2l}\left(  \xi\tau/2\right)  }{\left\vert \xi\right\vert
^{2\left(  n-\lambda\right)  }}d\xi\,q\left(  \tau\right)  d\tau<\infty.
\end{array}
\right\} \label{a970}%
\end{equation}

From Theorem \ref{Thm_Integ_u(xy)f(|x|)dx}, in general, for $k=1,2,\ldots,d$,
\begin{align}
\int\limits_{r\leq\left\vert x\right\vert \leq R}u\left(  \xi x\right)
f\left(  \left\vert x\right\vert \right)  dx  & =\int\limits_{r\leq\left\vert
x\right\vert \leq R}u\left(  \left\vert \xi\right\vert x_{k}\right)  f\left(
\left\vert x\right\vert \right)  dx\nonumber\\
& =\omega_{d-1}\int_{r}^{R}\left(  \int_{0}^{\pi}u\left(  \left\vert
\xi\right\vert \rho\cos t\right)  \sin^{d-2}tdt\right)  f\left(  \rho\right)
\rho^{d-1}d\rho\nonumber\\
& =\omega_{d-1}\int_{0}^{\pi}\int_{r}^{R}u\left(  \left\vert \xi\right\vert
\rho\cos t\right)  f\left(  \rho\right)  \rho^{d-1}d\rho\text{ }\sin
^{d-2}tdt\label{a975}%
\end{align}

and so%
\begin{equation}
\int\limits_{r\leq\left\vert \xi\right\vert \leq R}\frac{\sin^{2l}\left(
\xi\tau/2\right)  }{\left\vert \xi\right\vert ^{2\left(  n-\lambda\right)  }%
}d\xi=\int\limits_{r\leq\left\vert \xi\right\vert \leq R}\frac{\sin
^{2l}\left(  \xi_{1}\left\vert \tau\right\vert /2\right)  }{\left\vert
\xi\right\vert ^{2\left(  n-\lambda\right)  }}d\xi=\int\limits_{r\leq
\left\vert \xi\right\vert \leq R}\frac{\sin^{2l}\left(  \xi_{1}\left\vert
\tau\right\vert /2\right)  }{\left\vert \xi\right\vert ^{2\left(
n-\lambda\right)  }}d\xi.\label{a983}%
\end{equation}

The change of variable $\eta=\frac{1}{2}\left\vert \tau\right\vert \xi$,
$d\eta=\frac{1}{2^{d}}\left\vert \tau\right\vert ^{d}d\xi$, $\left\vert
\xi\right\vert =\frac{2}{\left\vert \tau\right\vert }\left\vert \eta
\right\vert $ yields%
\begin{align}
\int\limits_{r\leq\left\vert \xi\right\vert \leq R}\frac{\sin^{2l}\left(
\xi_{1}\left\vert \tau\right\vert /2\right)  }{\left\vert \xi\right\vert
^{2\left(  n-\lambda\right)  }}d\xi & =\int\limits_{\frac{r}{2}\left\vert
\tau\right\vert \leq\left\vert \eta\right\vert \leq\frac{R}{2}\left\vert
\tau\right\vert }\frac{\sin^{2l}\eta_{1}}{\left(  \frac{2}{\left\vert
\tau\right\vert }\left\vert \eta\right\vert \right)  ^{2\left(  n-\lambda
\right)  }}\left(  \frac{2}{\left\vert \tau\right\vert }\right)  ^{d}%
d\eta\nonumber\\
& =\left(  \frac{1}{2}\left\vert \tau\right\vert \right)  ^{2\left(
n-\lambda\right)  -d}\int\limits_{\frac{r}{2}\left\vert \tau\right\vert
\leq\left\vert \eta\right\vert \leq\frac{R}{2}\left\vert \tau\right\vert
}\frac{\sin^{2l}\eta_{1}}{\left\vert \eta\right\vert ^{2\left(  n-\lambda
\right)  }}d\eta.\label{a939}%
\end{align}

Thus%
\begin{equation}
\int\limits_{r\leq\left\vert \xi\right\vert \leq R}\frac{\sin^{2l}\left(
\xi\tau/2\right)  }{\left\vert \xi\right\vert ^{2\left(  n-\lambda\right)  }%
}d\xi=\left(  \frac{1}{2}\left\vert \tau\right\vert \right)  ^{2\left(
n-\lambda\right)  -d}\int\limits_{\frac{r}{2}\left\vert \tau\right\vert
\leq\left\vert \eta\right\vert \leq\frac{R}{2}\left\vert \tau\right\vert
}\frac{\sin^{2l}\eta_{1}}{\left\vert \eta\right\vert ^{2\left(  n-\lambda
\right)  }}d\eta.\label{a999}%
\end{equation}

Thus $\int\frac{\sin^{2l}\left(  \xi\tau/2\right)  }{\left\vert \xi\right\vert
^{2\left(  n-\lambda\right)  }}d\xi$ exists for almost all $\tau$ iff
$\int\frac{\sin^{2l}\eta_{1}}{\left\vert \eta\right\vert ^{2\left(
n-\lambda\right)  }}d\eta$ exists and consequently by Lemma
\ref{Lem_exists_integ_sin(x1)_|x|}, $\int\frac{\sin^{2l}\left(  \xi
\tau/2\right)  }{\left\vert \xi\right\vert ^{2\left(  n-\lambda\right)  }}%
d\xi$ exists for almost all $\tau$ iff
\begin{equation}
n-l-\frac{d}{2}<\lambda<n-\frac{d}{2}.\label{a994}%
\end{equation}

This means that \ref{a970} holds iff \ref{a994} holds and $\int\left\vert
\tau\right\vert ^{2\left(  n-\lambda\right)  -d}q\left(  \tau\right)
d\tau<\infty$ iff \ref{a994} holds and for some $R>0$,
\begin{equation}
\int_{\left\vert \tau\right\vert \geq R}\left\vert \tau\right\vert ^{2\left(
n-\lambda\right)  -d}q\left(  \tau\right)  d\tau<\infty,\label{a996}%
\end{equation}

i.e.
\begin{equation}
\int\frac{\left\vert \xi\right\vert ^{2\lambda}d\xi}{w\left(  \xi\right)
}<\infty\Leftrightarrow\ref{a994}\text{ }and\text{ }\ref{a996}\text{
}hold.\label{a998}%
\end{equation}

Now suppose $\int\frac{\left\vert \xi\right\vert ^{2\lambda}d\xi}{w\left(
\xi\right)  }<\infty$ for $0\leq\lambda\leq\kappa$. Then $\int_{\left\vert
\tau\right\vert \geq R}\left\vert \tau\right\vert ^{2n-d}q\left(  \tau\right)
d\tau<\infty$, $n-l-\frac{d}{2}<0$ and $\kappa<n-\frac{d}{2}$ i.e.
$\int_{\left\vert \tau\right\vert \geq R}\left\vert \tau\right\vert
^{2n-d}q\left(  \tau\right)  d\tau<\infty$ and $0<\kappa<n-\frac{d}{2}<l$ i.e.
the conditions \ref{a966} and \ref{a965} of the lemma hold.

On the other hand, suppose \ref{a966} and \ref{a965} hold i.e. $\int%
_{\left\vert \tau\right\vert \geq R}\left\vert \tau\right\vert ^{2n-d}q\left(
\tau\right)  d\tau<\infty$ and $0<\kappa<n-\frac{d}{2}<l$.

Then $0\leq\lambda\leq\kappa$ implies $\int_{\left\vert \tau\right\vert \geq
R}\left\vert \tau\right\vert ^{2\left(  n-\lambda\right)  -d}q\left(
\tau\right)  d\tau<\infty$ and $\lambda<n-\frac{d}{2}$ and $n-l-\frac{d}%
{2}<0\leq\lambda$ and \ref{a998} implies $\int\frac{\left\vert \xi\right\vert
^{2\lambda}d\xi}{w\left(  \xi\right)  }<\infty$ which means $w$ has property
W02 for $\kappa$.\medskip

\textbf{Part 3} From \ref{a903},%
\[
\int\frac{\left\vert \xi\right\vert ^{2\lambda}d\xi}{w\left(  \xi\right)
}=\tfrac{2^{2l}}{\left(  2\pi\right)  ^{d/2}}\int\int\frac{\sin^{2l}\left(
\xi\tau/2\right)  }{\left\vert \xi\right\vert ^{2\left(  n-\lambda\right)  }%
}d\xi\,q\left(  \tau\right)  d\tau,
\]

and from \ref{a999},%
\begin{equation}
\int\frac{\sin^{2l}\left(  \xi\tau/2\right)  }{\left\vert \xi\right\vert
^{2\left(  n-\lambda\right)  }}d\xi=\left(  \frac{1}{2}\left\vert
\tau\right\vert \right)  ^{2\left(  n-\lambda\right)  -d}\int\frac{\sin
^{2l}\eta_{1}}{\left\vert \eta\right\vert ^{2\left(  n-\lambda\right)  }}%
d\eta,\label{a802}%
\end{equation}

so that%
\begin{align*}
\int\frac{\left\vert \xi\right\vert ^{2\lambda}d\xi}{w\left(  \xi\right)  }  &
=\tfrac{2^{2l}}{\left(  2\pi\right)  ^{d/2}}\int\left(  \frac{1}{2}\left\vert
\tau\right\vert \right)  ^{2\left(  n-\lambda\right)  -d}\int\frac{\sin
^{2l}\eta_{1}}{\left\vert \eta\right\vert ^{2\left(  n-\lambda\right)  }}%
d\eta\,q\left(  \tau\right)  d\tau\\
& =\tfrac{2^{2\left(  l-n+\lambda\right)  +d}}{\left(  2\pi\right)  ^{d/2}%
}\int\frac{\sin^{2l}\eta_{1}}{\left\vert \eta\right\vert ^{2\left(
n-\lambda\right)  }}d\eta\,\int\left\vert \tau\right\vert ^{2\left(
n-\lambda\right)  -d}q\left(  \tau\right)  d\tau,
\end{align*}

and substituting \ref{a800} from the previous lemma completes this
part.\medskip

\textbf{Part 4} From Theorem \ref{Thm_basis_fn_properties_all_m_W2} the
smoothness of the basis function is $\max\limits_{\kappa<n-d/2}\left\lfloor
2\kappa\right\rfloor =2n-d-1$.\medskip

\textbf{Part 5} ??!
\end{proof}

With reference to condition \ref{a965}:

\begin{theorem}
\label{Thm_CentDiffQrad}If $q$ is radial, say $q\left(  \tau\right)
=q_{\odot}\left(  \left\vert \tau\right\vert \right)  $, then%
\begin{equation}
\int\limits_{\left\vert \tau\right\vert \geq R}\left\vert \tau\right\vert
^{2n-d}q\left(  \tau\right)  d\tau=\omega_{d}\int_{R}^{\infty}t^{2n-1}%
q_{\odot}\left(  t\right)  dt.\label{a931}%
\end{equation}

\end{theorem}

\begin{proof}
Use Theorem \ref{Ap148}.
\end{proof}

\section{Convolution formulas for the basis
function\label{Sect_convol_centdiff_basis}}

Using the Fourier transform to directly calculate the basis function using the
weight function formula \ref{a969} may be very awkward. Instead, using the
theory of fundamental solutions in $S^{\prime}$ of the operator $\left\vert
D\right\vert ^{2n}$ and the subspace $S_{\emptyset,m}\subset S$ (see
Definition \ref{Def_So,n} etc.) we derive an expression \ref{p36} for the
basis function which is the sum of the thin-plate spline $T_{n}$ (Lemma
\ref{Lem_thin_plate_splin} below) and several scaled convolutions of this
spline with the generating function $q$. This result is true modulo a unique
polynomial $p_{c}$ which satisfies $\left\vert D\right\vert ^{2n}p_{c}=0$.

Next for $d$ odd, using the $L_{loc}^{1}$ Taylor series expansion with
integral remainder described in Lemma \ref{Lem_Taylor_estim_loc_C_L1loc}, we
derive Theorem \ref{Thm_CentDiff_eqn_thinps} which supplies conditions on $q$
under which the polynomial $p_{c}$ is zero. Here we assume $q$ has bounded
support or that $q$ is radial with an extra condition. I have not done the
case where $d$ is even but have made partial attempt on a more general result
which does not distinguish between odd and even dimensions.

The binomial identities 2b of Lemma \ref{Lem_central_diff_op multivar} will
play a central role here.

Since $\widehat{G}_{c}\in L^{1}$ the central difference basis function
satisfies%
\begin{align}
\widehat{G}_{c}  & =\frac{\Delta_{2l}\widehat{q}\left(  \xi\right)
}{\left\vert \xi\right\vert ^{2n}},\label{p40}\\
& \Rightarrow\nonumber\\
\left\vert \xi\right\vert ^{2n}\widehat{G}_{c}  & =\Delta_{2l}\widehat{q}%
\left(  \xi\right)  .\label{p39}%
\end{align}

Suppose there exists $f\in S^{\prime}$ such that $\left\vert \xi\right\vert
^{2n}\widehat{f}=0$. Then $\operatorname*{supp}\widehat{f}=\left\{  0\right\}
$ which implies $\widehat{f}=\sum_{finite}a_{\alpha}D^{\alpha}\delta$ which in
turn implies $f$ is a polynomial. If we also know that $f\left(
\infty\right)  =0$ then clearly $f=0$. Hence if $G_{c}\in L_{loc}^{1}\cap
C^{\left(  0\right)  }$, $G_{c}\left(  \infty\right)  =0$ also satisfies
\ref{p39} then it is the unique solution with these properties. Further, the
unique solution of \ref{p40} has these properties and it also solves
\ref{p39}. Thus:
\begin{equation}
\text{Solving \ref{p39} for }G_{c}\in L_{loc}^{1}\cap C^{\left(  0\right)
}\text{, }G_{c}\left(  \infty\right)  =0\text{ also solves \ref{p40}%
,}\label{p41}%
\end{equation}

and so modulo a polynomial $p_{c}$ which satisfies $\left\vert D\right\vert
^{2n}p_{c}=0$
\begin{align}
\left\vert D\right\vert ^{2n}G_{c}  & =\left(  -1\right)  ^{n}\left(
\Delta_{2l}\widehat{q}\left(  \xi\right)  \right)  ^{\vee}\label{p00}\\
& =\left(  \left(  -1\right)  ^{n}\sum\limits_{j=-l}^{l}\left(  -1\right)
^{j}\tbinom{2l}{j+l}\widehat{q}\left(  j\xi\right)  \right)  ^{\vee
}\nonumber\\
& =\left(  \left(  -1\right)  ^{n}\tbinom{2l}{l}\widehat{q}\left(  0\right)
\right)  ^{\vee}+\left(  \left(  -1\right)  ^{n}\sum\limits_{j=-l,j\neq0}%
^{l}\left(  -1\right)  ^{j}\tbinom{2l}{j+l}\widehat{q}\left(  j\xi\right)
\right)  ^{\vee}\nonumber\\
& =\left(  -1\right)  ^{n}\tbinom{2l}{l}\widehat{q}\left(  0\right)  \left(
2\pi\right)  ^{d/2}\delta+\left(  -1\right)  ^{n}\sum\limits_{j=-l,j\neq0}%
^{l}\frac{\left(  -1\right)  ^{j}}{\left\vert j\right\vert ^{d}}\tbinom
{2l}{j+l}q\left(  \frac{\cdot}{j}\right) \nonumber\\
& =\left(  -1\right)  ^{n}\tbinom{2l}{l}\left(  \int q\right)  \delta+\left(
-1\right)  ^{n}\sum\limits_{j=-l,j\neq0}^{l}\frac{\left(  -1\right)  ^{j}%
}{\left\vert j\right\vert ^{d}}\tbinom{2l}{j+l}q\left(  \frac{\cdot}%
{j}\right)  .\label{p34}%
\end{align}

Compare with \ref{2.11}.

Now suppose $G_{\delta}$ is a fundamental solution associated with the
operator $\left\vert D\right\vert ^{2n}$ i.e.%
\begin{equation}
\left\vert D\right\vert ^{2n}G_{\delta}=\delta.\label{p01}%
\end{equation}

Now (see, for example, subsection 2.7.4 etc. of Vladimirov \cite{Vladimirov})
if $G_{\delta}\in L_{loc}^{1}$ and $q\in L^{1}$ then $G_{\delta}\ast q\left(
\frac{\cdot}{j}\right)  =\left(  2\pi\right)  ^{-d/2}\int G_{\delta}\left(
x-y\right)  q\left(  \frac{y}{j}\right)  dy\in L_{loc}^{1}$ and%
\begin{align}
G_{c}-p_{c}  & =\left(  2\pi\right)  ^{d/2}\delta\ast G_{c}=\left(
2\pi\right)  ^{d/2}\left\vert D\right\vert ^{2n}G_{\delta}\ast G_{c}=\left(
2\pi\right)  ^{d/2}G_{\delta}\ast\left\vert D\right\vert ^{2n}G_{c}%
=\nonumber\\
& =\left(  2\pi\right)  ^{d/2}G_{\delta}\ast\left(  \left(  -1\right)
^{n}\tbinom{2l}{l}\left(  \int q\right)  \delta+\left(  -1\right)  ^{n}%
\sum\limits_{j=-l,j\neq0}^{l}\frac{\left(  -1\right)  ^{j}}{\left\vert
j\right\vert ^{d}}\tbinom{2l}{j+l}q\left(  \frac{\cdot}{j}\right)  \right)
\nonumber\\
& =\left(  -1\right)  ^{n}\tbinom{2l}{l}\left(  \int q\right)  \left(
2\pi\right)  ^{d/2}G_{\delta}\ast\delta+\left(  -1\right)  ^{n}\sum
\limits_{j=-l,j\neq0}^{l}\frac{\left(  -1\right)  ^{j}}{\left\vert
j\right\vert ^{d}}\tbinom{2l}{j+l}\left(  2\pi\right)  ^{d/2}G_{\delta}%
\ast\left(  q\left(  \frac{\cdot}{j}\right)  \right) \nonumber\\
& =\left(  -1\right)  ^{n}\tbinom{2l}{l}\left(  \int q\right)  G_{\delta
}+\left(  -1\right)  ^{n}\sum\limits_{j=-l,j\neq0}^{l}\frac{\left(  -1\right)
^{j}}{\left\vert j\right\vert ^{d}}\tbinom{2l}{j+l}\left(  2\pi\right)
^{d/2}G_{\delta}\ast\left(  q\left(  \frac{\cdot}{j}\right)  \right)
\label{p37}\\
& =\left(  -1\right)  ^{n}\left(  \tbinom{2l}{l}\left(  \int q\right)
G_{\delta}+\sum\limits_{j=-l,j\neq0}^{l}\frac{\left(  -1\right)  ^{j}%
}{\left\vert j\right\vert ^{d}}\tbinom{2l}{j+l}\int G_{\delta}\left(
\cdot-y\right)  q\left(  \frac{y}{j}\right)  dy\right)  .\label{p21}%
\end{align}

The properties \ref{p38} below show that the convolution in \ref{p37} is a
well-defined distribution of slow growth. Hence we can differentiate the
convolution to obtain \ref{p34} and then \ref{p00} and finally \ref{p39}.

We want%
\begin{equation}
G_{c}\left(  \infty\right)  =0,\label{p22}%
\end{equation}

because then \ref{p41} implies that $G_{c}$ is the desired basis function.

The next lemma comes from Subsection 1.7.1 of my positive order document
Williams \cite{WilliamsPosOrdSmthV3} and uses the subspaces $S_{\emptyset
,2m}\subset S$ introduced above in Definition \ref{Def_So,n}.

\begin{lemma}
\label{Lem_thin_plate_splin}\textbf{Thin-plate spline properties} Suppose
$T_{n}\in S^{\prime}$\ is the thin-plate spline function with integer order
$n>d/2$ i.e.%
\begin{equation}
T_{n}\left(  x\right)  :=\left\{
\begin{array}
[c]{ll}%
\left(  -1\right)  ^{n-\frac{d}{2}+1}\left\vert x\right\vert ^{2n-d}%
\log\left\vert x\right\vert , & d\text{ }is\text{ }even,\\
\left(  -1\right)  ^{n-\frac{1}{2}\left(  d-1\right)  }\left\vert x\right\vert
^{2n-d}, & d\text{ }is\text{ }odd,
\end{array}
\right. \label{p24}%
\end{equation}

for all $x\in\mathbb{R}^{d}$. Then:

\begin{enumerate}
\item If $c_{n,d}:=e\left(  n-d/2\right)  $ where the function $e$ is defined
by%
\begin{equation}
e\left(  s\right)  =\left\{
\begin{array}
[c]{ll}%
\pi^{\frac{d}{2}}2^{2s+d-1}\Gamma\left(  s+\frac{d}{2}\right)  \Gamma\left(
s+1\right)  , & s=1,2,3,\ldots,\\
\left.
\begin{array}
[c]{l}%
\pi^{\frac{d}{2}}2^{2s+d}\Gamma\left(  s+\frac{d}{2}\right)  \Gamma\left(
s+1\right)  \frac{\left(  -1\right)  ^{1+\left\lceil s\right\rceil }\sin\pi
s}{\pi},\\
\left(  -1\right)  ^{\left\lceil s\right\rceil }\pi^{\frac{d}{2}}%
2^{2s+d}\Gamma\left(  s+\frac{d}{2}\right)  /\Gamma\left(  -s\right)
\end{array}
\right\}  , & s>0,\text{ }s\neq1,2,3\ldots,
\end{array}
\right\} \label{p32}%
\end{equation}

we have%
\begin{equation}
\widehat{T_{n}}=\frac{c_{n,d}}{\left\vert \cdot\right\vert ^{2n}}\text{
}on\text{ }S_{\emptyset,2n}.\label{p25}%
\end{equation}

\item $T_{n}\in C_{BP}^{\left(  2n-d-1\right)  }\cap C^{\infty}\left(
\mathbb{R}^{d}\setminus0\right)  $ and $\left(  \widehat{a}D\right)  ^{k}%
T_{n}\in L_{loc}^{1}$ when $k\leq2n-1$.

\item When $d$ \textbf{is odd}
\[
\left\vert \left(  \widehat{a}D_{x}\right)  ^{k}T_{n}\left(  x\right)
\right\vert \leq k_{k,2n-d}\left\vert x\right\vert ^{2n-d-k},
\]

where $k_{k,2n-d}=k_{k,2\left(  n-\frac{d+1}{2}\right)  +1}$ is given by the
\textbf{conjectures} in the positive order document. Also, using the
Pochhammer symbol $\left(  \cdot\right)  _{k}$,%
\[
\left\vert D\right\vert ^{2k}T_{n}\left(  x\right)  =2^{2k}\left(  -\left(
n-\frac{d-1}{2}\right)  \right)  _{k}\left(  -\left(  n-\frac{1}{2}\right)
\right)  _{k}\left\vert x\right\vert ^{2n-d-2k},\quad k\geq0.
\]

and%
\begin{equation}
\left(  \widehat{\cdot}D\right)  ^{k}T_{n}\left(  x\right)  =\left\vert
x\right\vert ^{2n-d-k}\times\left\{
\begin{array}
[c]{ll}%
\left(  -1\right)  ^{k+n-\frac{d-1}{2}}\left(  -\left(  2n-d\right)  \right)
_{k}, & k\leq2n-d,\\
0, & k>2n-d,
\end{array}
\right\}  .\label{p03}%
\end{equation}

\item When $d$ \textbf{is even}, using the Pochhammer symbol,%
\begin{equation}
\frac{1}{k!}\left\vert \left(  \widehat{a}D\right)  ^{k}T_{n}\left(  x\right)
\right\vert \leq\left(  \frac{k_{k,2n-d}}{k!}\left\vert \log\left\vert
x\right\vert \right\vert +\sum\limits_{j=0}^{k-1}\left(  2-\frac{1}{j}\right)
\frac{k_{j,2n-d}}{j!}\right)  \left\vert x\right\vert ^{2n-d-k},\quad
k\geq1,\label{p02}%
\end{equation}

where the $k_{k-j,2n-d}$ are given by the \textbf{conjectures} in the positive
order document. Also%
\[
\left\vert D\right\vert ^{2k}T_{n}\left(  x\right)  =\left(  -1\right)
^{n-\frac{d}{2}+1}2^{2k}\left(  -n+\frac{d}{2}\right)  _{k}\left(
-n+1\right)  _{k}\left\vert x\right\vert ^{2n-d-2k},\text{ }k\geq1,
\]

and for $k\geq0$,%
\begin{align*}
\frac{1}{k!} &  \left(  \widehat{\cdot}D\right)  ^{k}T_{n}\left(  x\right) \\
&  =\left\vert x\right\vert ^{2n-d-k}\times\left\{
\begin{array}
[c]{ll}%
\left(  \tbinom{2n-d}{k}\log\left\vert x\right\vert +\sum\limits_{j=1}%
^{k}\frac{\left(  -1\right)  ^{j+1}}{j}\tbinom{2n-d}{k-j}\right)  , & 1\leq
k\leq2n-d,\\
\left(  -1\right)  ^{k-1}B\left(  2n-d+1,k-2n+d\right)  , & k>2n-d.
\end{array}
\right.  .
\end{align*}

\end{enumerate}
\end{lemma}

\begin{lemma}
\label{Lem_eval_thin_plate_splin_const}When $d$ is odd,%
\begin{align*}
c_{n,d}=e\left(  n-d/2\right)   & =\pi^{\frac{d-1}{2}}2^{d}\left(  n-1\right)
!\left(  2n-d\right)  !!\\
& =\pi^{\frac{d-1}{2}}2^{-\left(  n-\frac{3d-1}{2}\right)  }\frac{\left(
n-1\right)  !\left(  2n-d+1\right)  !}{\left(  n-\frac{d-1}{2}\right)  !}.
\end{align*}

\end{lemma}

\begin{proof}
Define $m$ by
\begin{equation}
2n-d=:2m-1,\label{a2.23}%
\end{equation}

Set $n-\frac{d}{2}=m-\frac{1}{2}$ so that
\begin{align*}
e\left(  n-\frac{d}{2}\right)   & =\pi^{\frac{d}{2}}2^{2\left(  n-\frac{d}%
{2}\right)  +d}\Gamma\left(  \left(  n-\frac{d}{2}\right)  +\frac{d}%
{2}\right)  \Gamma\left(  \left(  n-\frac{d}{2}\right)  +1\right)  \times\\
& \qquad\times\frac{\left(  -1\right)  ^{1+\left(  n-\frac{d}{2}+\frac{1}%
{2}\right)  }\sin\pi\left(  n-\frac{d}{2}\right)  }{\pi}\\
& =\pi^{\frac{d}{2}}2^{2n}\Gamma\left(  n\right)  \Gamma\left(  m+\frac{1}%
{2}\right)  \frac{\left(  -1\right)  ^{1+m}\sin\pi\left(  m-\frac{1}%
{2}\right)  }{\pi}\\
& =\pi^{\frac{d}{2}-1}2^{2n}\left(  n-1\right)  !\Gamma\left(  m+\frac{1}%
{2}\right)  \left(  -1\right)  ^{m-1}\sin\pi\left(  m-\frac{1}{2}\right) \\
& =\pi^{\frac{d}{2}-1}2^{2n}\left(  n-1\right)  !\left(  m-\frac{1}{2}\right)
\Gamma\left(  m-\frac{1}{2}\right)  \left(  -1\right)  ^{m}\cos\pi m\\
& =\pi^{\frac{d}{2}-1}2^{2n}\left(  n-1\right)  !\left(  m-\frac{1}{2}\right)
\Gamma\left(  m-\frac{1}{2}\right) \\
& =\pi^{\frac{d}{2}-1}2^{2n}\left(  n-1\right)  !\left(  m-\frac{1}{2}\right)
\left(  2m-3\right)  !!\pi^{\frac{1}{2}}2^{-\left(  m-1\right)  }\\
& =\pi^{\frac{d-1}{2}}2^{d+1}\left(  n-1\right)  !\left(  m-\frac{1}%
{2}\right)  \left(  2m-3\right)  !!\\
& =\pi^{\frac{d-1}{2}}2^{d}\left(  n-1\right)  !\left(  2m-1\right)  !!\\
& =\pi^{\frac{d-1}{2}}2^{d}\left(  n-1\right)  !\left(  2n-d\right)  !!,\text{
}as\text{ }claimed,\\
& =\pi^{\frac{d-1}{2}}2^{d}\left(  n-1\right)  \ldots\left(  n-\frac{d-3}%
{2}\right)  \left(  n-\frac{d-1}{2}\right)  !\left(  2n-d\right)  !!\\
& =\pi^{\frac{d-1}{2}}2^{d}\left(  n-1\right)  \ldots\left(  n-\frac{d-3}%
{2}\right)  2^{-\left(  n-\frac{d-1}{2}\right)  }\left(  2n-d+1\right)
!!\left(  2n-d\right)  !!\\
& =\pi^{\frac{d-1}{2}}2^{d-\left(  n-\frac{d-1}{2}\right)  }\left(
n-1\right)  \ldots\left(  n-\frac{d-3}{2}\right)  \left(  2n-d+1\right)  !\\
& =\pi^{\frac{d-1}{2}}2^{-\left(  n-\frac{3d-1}{2}\right)  }\frac{\left(
n-1\right)  !\left(  2n-d+1\right)  !}{\left(  n-\frac{d-1}{2}\right)  !}.
\end{align*}

\end{proof}

A consequence of part 1 of Lemma \ref{Lem_thin_plate_splin} is%
\[
\left\vert \cdot\right\vert ^{2n}\widehat{T_{n}}=c_{n,d}\text{ }in\text{
}S^{\prime}.
\]

i.e.%
\[
\left(  -1\right)  ^{n}\left\vert D\right\vert ^{2n}T_{n}=c_{n,d}%
\overset{\vee}{1}=\left(  2\pi\right)  ^{d/2}c_{n,d}\delta,
\]

so that we can choose%
\begin{equation}
G_{\delta}=\frac{\left(  -1\right)  ^{n}}{\left(  2\pi\right)  ^{d/2}c_{n,d}%
}T_{n}.\label{p23}%
\end{equation}

Lemma \ref{Lem_thin_plate_splin} now implies $G_{\delta}\in L_{loc}^{1}$.

In fact, $G_{\delta}\ast q\in L_{loc}^{1}$ and is also a function of
polynomial increase i.e. $G_{\delta}\ast q\in S^{\prime}$, as we now show by
\textbf{considering separately the two cases} $d$ is odd and $d$ is
even.\medskip

\fbox{\textbf{Case} $d$ is odd:} Now $d>1$ means $3\leq d\leq2n-1$ and using
Peetre's inequality:
\begin{equation}
\frac{1}{1+\left\vert b\right\vert }\leq\frac{1+\left\vert a\right\vert
}{1+\left\vert a-b\right\vert },\quad a,b\in\mathbb{R}^{d},\label{p33}%
\end{equation}

we get%
\begin{align*}
\left(  1+\left\vert x\right\vert \right)  ^{-\left(  2n-d\right)  }%
\int\left\vert T_{n}\left(  x-y\right)  \right\vert q\left(  \frac{y}%
{j}\right)  dy  & =\int\frac{\left\vert T_{n}\left(  x-y\right)  \right\vert
}{\left(  1+\left\vert x\right\vert \right)  ^{2n-d}}q\left(  \frac{y}%
{j}\right)  dy\\
& \leq\int\frac{\left\vert T_{n}\left(  x-y\right)  \right\vert }{\left(
1+\left\vert x-y\right\vert \right)  ^{2n-d}}\left(  1+\left\vert y\right\vert
\right)  ^{2n-d}q\left(  \frac{y}{j}\right)  dy\\
& =\int\frac{\left\vert x-y\right\vert ^{2n-d}}{\left(  1+\left\vert
x-y\right\vert \right)  ^{2n-d}}\left(  1+\left\vert y\right\vert \right)
^{2n-d}q\left(  \frac{y}{j}\right)  dy\\
& <\int\left(  1+\left\vert y\right\vert \right)  ^{2n-d}q\left(  \frac{y}%
{j}\right)  dy\\
& \left(  z=y/j\Longrightarrow\right) \\
& \leq\left\vert j\right\vert ^{2n}\int\left(  1+\left\vert \cdot\right\vert
\right)  ^{2n-d}q\\
& <\infty,
\end{align*}

where the last step employs Theorem \ref{Thm_wt_func_mdim}. Thus%
\begin{equation}
\left.
\begin{array}
[c]{c}%
d\text{ }odd,\\
\ \\
\int\left\vert T_{n}\left(  x-y\right)  \right\vert q\left(  \frac{y}%
{j}\right)  dy<\left\vert j\right\vert ^{2n}\left(  \int\left(  1+\left\vert
\cdot\right\vert \right)  ^{2n-d}q\right)  \left(  1+\left\vert x\right\vert
\right)  ^{2n-d},\quad x\in\mathbb{R}^{d},\\
\ \\
T_{n}\ast\left(  q\left(  \frac{\cdot}{j}\right)  \right)  \in S^{\prime}\cap
L_{loc}^{1},\quad\left\vert j\right\vert =1,2,3,\ldots
\end{array}
\right\} \label{p26}%
\end{equation}
\medskip

\fbox{\textbf{Case} $d$ is even:} This means $2\leq d\leq2n-2$ and following
the argument used for $d$ odd,%
\begin{align*}
\left(  1+\left\vert x\right\vert \right)  ^{-\left(  2n-d+1\right)  } &
\int\left\vert T_{n}\left(  x-y\right)  \right\vert q\left(  \frac{y}%
{j}\right)  dy\\
&  =\int\frac{\left\vert T_{n}\left(  x-y\right)  \right\vert }{\left(
1+\left\vert x\right\vert \right)  ^{2n-d+1}}q\left(  \frac{y}{j}\right)  dy\\
&  =\int\frac{\left\vert T_{n}\left(  x-y\right)  \right\vert }{\left(
1+\left\vert x-y\right\vert \right)  ^{2n-d+1}}\frac{\left(  1+\left\vert
x-y\right\vert \right)  ^{2n-d+1}}{\left(  1+\left\vert x\right\vert \right)
^{2n-d+1}}q\left(  \frac{y}{j}\right)  dy\\
&  \leq\int\frac{\left\vert T_{n}\left(  x-y\right)  \right\vert }{\left(
1+\left\vert x-y\right\vert \right)  ^{2n-d+1}}\left(  1+\left\vert
y\right\vert \right)  ^{2n-d+1}q\left(  \frac{y}{j}\right)  dy\\
&  =\int\frac{\left\vert x-y\right\vert ^{2n-d}\log\left\vert x-y\right\vert
}{\left(  1+\left\vert x-y\right\vert \right)  ^{2n-d+1}}\left(  1+\left\vert
y\right\vert \right)  ^{2n-d+1}q\left(  \frac{y}{j}\right)  dy\\
&  <\int\frac{\left\vert x-y\right\vert ^{2n-d+1}}{\left(  1+\left\vert
x-y\right\vert \right)  ^{2n-d+1}}\left(  1+\left\vert y\right\vert \right)
^{2n-d+1}q\left(  \frac{y}{j}\right)  dy\\
&  <\int\left(  1+\left\vert y\right\vert \right)  ^{2n-d+1}q\left(  \frac
{y}{j}\right)  dy\\
&  <\int\left(  1+\left\vert y\right\vert \right)  ^{2n-d+1}q\left(  \frac
{y}{j}\right)  dy\\
&  \left(  z=y/j\Longrightarrow\right) \\
&  \leq\left\vert j\right\vert ^{2n+1}\int\left(  1+\left\vert \cdot
\right\vert \right)  ^{2n-d+1}q\\
&  <\infty,
\end{align*}

where the last step employed Theorem \ref{Thm_wt_func_mdim}. Thus%
\begin{equation}
\left.
\begin{array}
[c]{c}%
When\text{ }d\text{ }is\text{ }even:\\
\ \\
\int\left\vert T_{n}\left(  x-y\right)  \right\vert q\left(  \frac{y}%
{j}\right)  dy<\left\vert j\right\vert ^{2n+1}\left(  \int\left(  1+\left\vert
\cdot\right\vert \right)  ^{2n-d+1}q\right)  \left(  1+\left\vert x\right\vert
\right)  ^{2n-d+1},\quad x\in\mathbb{R}^{d},\\
\ \\
T_{n}\ast\left(  q\left(  \frac{\cdot}{j}\right)  \right)  \in S^{\prime}\cap
L_{loc}^{1},\quad\left\vert j\right\vert =1,2,3,\ldots
\end{array}
\right\} \label{p28}%
\end{equation}

Noting equation (16) of \S \S 2.7.4 of Vladimirov we have shown that
\begin{equation}
\left.
\begin{array}
[c]{l}%
G_{\delta}\ast q\in\mathcal{D}^{\prime}\cap S^{\prime}\cap L_{loc}^{1}.\\
G_{\delta}\ast q\text{ has polynomial growth.}%
\end{array}
\right\} \label{p38}%
\end{equation}

From \ref{p21} and \ref{p23},%
\[
G_{\delta}=\frac{\left(  -1\right)  ^{n}}{\left(  2\pi\right)  ^{d/2}c_{n,d}%
}T_{n},
\]

and so%
\begin{align}
&  G_{c}\left(  x\right)  -p_{c}\left(  x\right) \nonumber\\
&  =\frac{\left(  -1\right)  ^{n}}{\left(  2\pi\right)  ^{d/2}}\left(
\tbinom{2l}{l}\left(  \int q\right)  G_{\delta}\left(  x\right)
+\sum\limits_{j=-l,j\neq0}^{l}\frac{\left(  -1\right)  ^{j}}{\left\vert
j\right\vert ^{d}}\tbinom{2l}{j+l}\int G_{\delta}\left(  x-y\right)  q\left(
\frac{y}{j}\right)  dy\right) \nonumber\\
&  =\frac{1}{\left(  2\pi\right)  ^{d}c_{n,d}}\left(  \tbinom{2l}{l}\left(
\int q\right)  T_{n}\left(  x\right)  +\sum\limits_{j=-l,j\neq0}^{l}%
\frac{\left(  -1\right)  ^{j}}{\left\vert j\right\vert ^{d}}\tbinom{2l}%
{j+l}\int T_{n}\left(  x-y\right)  q\left(  \frac{y}{j}\right)  dy\right)
\label{p30}\\
&  =\frac{1}{\left(  2\pi\right)  ^{d}c_{n,d}}\left(  \tbinom{2l}{l}\left(
\int q\right)  T_{n}\left(  x\right)  +\sum\limits_{j=-l,j\neq0}^{l}\left(
-1\right)  ^{j}\tbinom{2l}{j+l}\int T_{n}\left(  x-jy\right)  q\left(
y\right)  dy\right) \nonumber\\
&  =\frac{1}{\left(  2\pi\right)  ^{d}c_{n,d}}\int\left(  \sum\limits_{j=-l}%
^{l}\left(  -1\right)  ^{j}\tbinom{2l}{j+l}T_{n}\left(  x-jy\right)  \right)
q\left(  y\right)  dy,\label{p35}%
\end{align}

so that we have shown:

\begin{theorem}
\label{Thm_centdiffbasis_centdiff_formula}Suppose the function $q$ and
parameters $n$ and $l$ generate a multivariate central difference weight
function. Set%
\begin{equation}
G_{c}\left(  x\right)  =\frac{1}{\left(  2\pi\right)  ^{d}c_{n,d}}\int\left(
\sum\limits_{j=-l}^{l}\left(  -1\right)  ^{j}\tbinom{2l}{j+l}T_{n}\left(
x-jy\right)  \right)  q\left(  y\right)  dy+p_{c}\left(  x\right)
,\label{p36}%
\end{equation}

where $p_{c}$ is a polynomial which satisfies $\left\vert D\right\vert
^{2n}p_{c}=0$. Then $G_{c}\left(  \infty\right)  =0$ implies $G_{c}$ is the
(unique) basis function.

Also, if $d$ is odd then%
\[
G_{c}\left(  0\right)  =\frac{1}{\left(  2\pi\right)  ^{d}c_{n,d}}\left(
\sum\limits_{j=-l}^{l}\left(  -1\right)  ^{j}\tbinom{2l}{j+l}j^{2n-d}\right)
\int T_{n}\left(  y\right)  q\left(  y\right)  dy+p_{c}\left(  0\right)  .
\]

\end{theorem}

\begin{proof}
Note \ref{p41}.%
\begin{align*}
G_{c}\left(  0\right)   & =\frac{1}{\left(  2\pi\right)  ^{d}c_{n,d}}%
\int\left(  \sum\limits_{j=-l}^{l}\left(  -1\right)  ^{j}\tbinom{2l}{j+l}%
T_{n}\left(  -jy\right)  \right)  q\left(  y\right)  dy+p_{c}\left(  0\right)
\\
& =\frac{1}{\left(  2\pi\right)  ^{d}c_{n,d}}\int\left(  \sum\limits_{j=-l}%
^{l}\left(  -1\right)  ^{j}\tbinom{2l}{j+l}j^{2n-d}T_{n}\left(  y\right)
\right)  q\left(  y\right)  dy+p_{c}\left(  0\right) \\
& =\frac{1}{\left(  2\pi\right)  ^{d}c_{n,d}}\left(  \sum\limits_{j=-l}%
^{l}\left(  -1\right)  ^{j}\tbinom{2l}{j+l}j^{2n-d}\right)  \int T_{n}\left(
y\right)  q\left(  y\right)  dy+p_{c}\left(  0\right)  .
\end{align*}

\end{proof}

\begin{remark}
Note that this theorem does not say that there exists a polynomial $p_{c}$
such that \ref{p36} is the basis function.

Compare this formula with Theorem \ref{Thm_centdiff_basis_Tn_q_suppbnd} below
which IS an existence result.
\end{remark}

\textbf{We consider two cases}: $d$ is odd and $d$ is even. First the odd
case. The even case will be studied later in Subsection
\ref{SbSect_Radial_q_d_even}.

\subsection{Case: $d$ is odd}

Here we will impose conditions on $q$ such that the polynomial $p_{c}$ is
zero. This will involve using the $L_{loc}^{1}$ Taylor series expansion with
integral remainder described in Lemma \ref{Lem_Taylor_estim_loc_C_L1loc}.

From Theorem \ref{Thm_wt_func_mdim},%
\begin{equation}
0\leq\kappa<n-\frac{d}{2}<l,\text{ }i.e.\text{ }0<2n-d<2l,\label{p31}%
\end{equation}

so set (as in \ref{a2.23})
\begin{equation}
2m:=2n-d+1\geq2,\label{a031}%
\end{equation}

and note that%
\begin{equation}
1\leq m\leq l.\label{a900}%
\end{equation}

For fixed $a\in\mathbb{R}^{d}$, Lemma \ref{Lem_thin_plate_splin} implies%
\begin{align*}
\left\vert \left(  \left(  \widehat{a}D\right)  ^{2m}T_{n}\right)  \left(
x\right)  \right\vert  & \leq k_{2m,2n-d}\left\vert x\right\vert
^{2n-d-2m}=\frac{k_{2m,2n-d}}{\left\vert x\right\vert },\quad x\in
\mathbb{R}^{d},\\
\left(  \widehat{a}D\right)  ^{2m}T_{n}  & \in L_{loc}^{1},\\
\left(  \widehat{a}D\right)  ^{2m}T_{n}\left(  x\right)   & \rightarrow0\text{
}as\text{ }\left\vert x\right\vert \rightarrow\infty.
\end{align*}

This means we can employ the $L_{loc}^{1}$ Taylor series expansion about $x$
with integral remainder described in Lemma \ref{Lem_Taylor_estim_loc_C_L1loc}
(here $\Omega=\mathbb{R}^{d}$ and $f=T_{n}$) to write%
\[
T_{n}\left(  x-jy\right)  =\sum_{k\leq2m-1}j^{k}\frac{\left(  -yD\right)
^{k}}{k!}T_{n}\left(  x\right)  +\frac{j^{2m}}{\left(  2m-1\right)  !}\int%
_{0}^{1}g_{2m-1}\left(  s\right)  \left(  \left(  yD\right)  ^{2m}%
T_{n}\right)  \left(  x-sjy\right)  ds,
\]

where the function $g_{n}$ is given by \ref{a117} and we clarify that
\[
\left(  \left(  yD\right)  ^{2m}T_{n}\right)  \left(  x-sjy\right)  =\left(
\left(  yD_{z}\right)  ^{2m}T_{n}\left(  z\right)  \right)  \left(
z=x-sjy\right)  .
\]

The first summation in \ref{p36} now becomes,%
\begin{align*}
&  \sum\limits_{j=-l}^{l}\left(  -1\right)  ^{j}\tbinom{2l}{j+l}T_{n}\left(
x-jy\right) \\
&  =\sum\limits_{j=-l}^{l}\left(  -1\right)  ^{j}\tbinom{2l}{j+l}\left(
\sum_{k\leq2m-1}j^{k}\frac{\left(  -yD\right)  ^{k}}{k!}T_{n}\left(  x\right)
+j^{2m}\int_{0}^{1}\frac{g_{2m-1}\left(  s\right)  }{\left(  2m-1\right)
!}\left(  \left(  yD\right)  ^{2m}T_{n}\right)  \left(  x-sjy\right)
ds\right) \\
&  =\sum_{k\leq2m-1}\left(  \sum\limits_{j=-l}^{l}\left(  -1\right)
^{j}\tbinom{2l}{j+l}j^{k}\right)  \frac{\left(  -yD\right)  ^{k}}{k!}%
T_{n}\left(  x\right)  +\\
&  \qquad\qquad\qquad+\sum\limits_{j=-l}^{l}\left(  -1\right)  ^{j}\tbinom
{2l}{j+l}j^{2m}\int_{0}^{1}\frac{g_{2m-1}\left(  s\right)  }{\left(
2m-1\right)  !}\left(  \left(  yD\right)  ^{2m}T_{n}\right)  \left(
x-sjy\right)  ds.
\end{align*}

But from Lemma \ref{Lem_central_diff_op multivar},%
\begin{equation}
\sum\limits_{j=-l}^{l}\left(  -1\right)  ^{j}\tbinom{2l}{j+l}j^{k}=\left\{
\begin{array}
[c]{ll}%
0, & k=0,1,2,\ldots,2l-1,\\
\left(  -1\right)  ^{l}\left(  2l\right)  !, & k=2l,
\end{array}
\right. \label{a013}%
\end{equation}

and from \ref{a900}, $2m-1\leq2l-1$. Hence%
\begin{align*}
\sum\limits_{j=-l}^{l}\left(  -1\right)  ^{j}\tbinom{2l}{j+l}T_{n}\left(
x-jy\right)   & =\sum\limits_{j=-l}^{l}\left(  -1\right)  ^{j}\tbinom{2l}%
{j+l}j^{2m}\int_{0}^{1}\frac{g_{2m-1}\left(  s\right)  }{\left(  2m-1\right)
!}\left(  \left(  yD\right)  ^{2m}T_{n}\right)  \left(  x-sjy\right)  ds\\
& =\sum\limits_{j=-l,j\neq0}^{l}\left(  -1\right)  ^{j}\tbinom{2l}{j+l}%
j^{2m}\int_{0}^{1}\frac{g_{2m-1}\left(  s\right)  }{\left(  2m-1\right)
!}\left(  \left(  yD\right)  ^{2m}T_{n}\right)  \left(  x-sjy\right)  ds,
\end{align*}

and equation \ref{p36} becomes%
\begin{align}
&  G_{c}\left(  x\right)  -p_{c}\left(  x\right) \nonumber\\
&  =\frac{1}{\left(  2\pi\right)  ^{d}c_{n,d}}\left(  \int\left(
\sum\limits_{j=-l}^{l}\left(  -1\right)  ^{j}\tbinom{2l}{j+l}T_{n}\left(
x-jy\right)  \right)  q\left(  y\right)  dy\right) \nonumber\\
&  =\frac{1}{\left(  2\pi\right)  ^{d}c_{n,d}}\int\left(  \sum
\limits_{j=-l,j\neq0}^{l}\left(  -1\right)  ^{j}\tbinom{2l}{j+l}j^{2m}\int%
_{0}^{1}\frac{g_{2m-1}\left(  s\right)  }{\left(  2m-1\right)  !}\left(
\left(  yD\right)  ^{2m}T_{n}\right)  \left(  x-sjy\right)  ds\right)
q\left(  y\right)  dy.\label{a01}%
\end{align}

Further%
\begin{align*}
&  \left\vert G_{c}\left(  x\right)  -p_{c}\left(  x\right)  \right\vert \\
&  \leq\frac{1}{\left(  2\pi\right)  ^{d}c_{n,d}}\left(  \sum
\limits_{j=-l,j\neq0}^{l}\tbinom{2l}{j+l}j^{2m}\right)  \max_{j=-l,j\neq0}%
^{l}\int\int_{0}^{1}\frac{\left(  1-s\right)  ^{2m-1}}{\left(  2m-1\right)
!}\left\vert \left(  \left(  \widehat{y}D\right)  ^{2m}T_{n}\right)  \left(
x-sjy\right)  \right\vert ds\text{\thinspace}\left\vert y\right\vert
^{2m}q\left(  y\right)  dy,
\end{align*}

but%
\[
\sum\limits_{j=-l,j\neq0}^{l}\tbinom{2l}{j+l}j^{2m}=2\sum\limits_{j=1}%
^{l}\tbinom{2l}{j+l}j^{2m}=2\sum\limits_{j=1}^{l}\tbinom{2l}{l-j}j^{2m}%
=2\sum\limits_{j=0}^{l}\tbinom{2l}{l-j}j^{2m}=2\sum\limits_{k=0}^{l}%
\tbinom{2l}{k}\left(  l-k\right)  ^{2m},
\]

so%
\begin{equation}
\left.
\begin{array}
[c]{l}%
\left\vert G_{c}\left(  x\right)  -p_{c}\left(  x\right)  \right\vert \\
\leq\frac{2\sum\limits_{k=0}^{l}\tbinom{2l}{k}\left(  l-k\right)  ^{2m}%
}{\left(  2\pi\right)  ^{d}c_{n,d}}\max\limits_{j=-l,j\neq0}^{l}\int\int%
_{0}^{1}\frac{\left(  1-s\right)  ^{2m-1}}{\left(  2m-1\right)  !}\left\vert
\left(  \left(  \widehat{y}D\right)  ^{2m}T_{n}\right)  \left(  x-sjy\right)
\right\vert ds\text{ }\left\vert y\right\vert ^{2m}q\left(  y\right)  dy,\\
\ \\
\qquad p_{c}=0,\text{ }2m=2n-d+1,\text{ }1\leq m\leq l,\text{ }m\leq n-1.
\end{array}
\right\} \label{a02}%
\end{equation}

From part 3 of Theorem \ref{Lem_thin_plate_splin},%
\begin{equation}
\left\vert \left(  \widehat{a}D\right)  ^{2m}T_{n}\left(  x\right)
\right\vert \leq k_{2m,2m-1}\left\vert x\right\vert ^{-1},\quad x,a\in
\mathbb{R}^{d}\setminus0.\label{a03}%
\end{equation}

where $k_{k,2n-d}$ is independent of $x$. Hence%
\begin{equation}
\int\int_{0}^{1}\left(  1-s\right)  ^{2m-1}\left\vert \left(  \left(
\widehat{y}D\right)  ^{2m}T_{n}\right)  \left(  x-sjy\right)  \right\vert
ds\text{ }\left\vert y\right\vert ^{2m}q\left(  y\right)  dy\leq
k_{2m,2m-1}\int\left(  \int_{0}^{1}\frac{\left(  1-s\right)  ^{2m-1}%
}{\left\vert x-sjy\right\vert }ds\right)  \left\vert y\right\vert
^{2m}q\left(  y\right)  dy.\label{a04}%
\end{equation}
\medskip

\underline{\textbf{Subcase }$\operatorname*{supp}q$ is bounded,} say
$\operatorname*{supp}q\subset B_{R}$. Then if $\left\vert x\right\vert
\geq2lR$ and $\left\vert j\right\vert \leq l$,%
\[
\left\vert x-sjy\right\vert \geq\left\vert \left\vert x\right\vert
-sl\left\vert y\right\vert \right\vert =\left\vert x\right\vert -sl\left\vert
y\right\vert \geq\left\vert x\right\vert -lR,
\]

and so
\begin{align}
\int\int_{0}^{1}\left(  1-s\right)  ^{2m-1}\left\vert \left(  \left(
\widehat{y}D\right)  ^{2m}T_{n}\right)  \left(  x-sjy\right)  \right\vert
ds\text{ }\left\vert y\right\vert ^{2m}q\left(  y\right)  dy  & \leq
\frac{k_{2m,2m-1}}{\left\vert x\right\vert -lR}\int\left(  \int_{0}^{1}\left(
1-s\right)  ^{2m-1}ds\right)  \left\vert y\right\vert ^{2m}q\left(  y\right)
dy\nonumber\\
& =\frac{k_{2m,2m-1}}{\left\vert x\right\vert -lR}\frac{1}{2m}\int\left\vert
y\right\vert ^{2m}q\left(  y\right)  dy\label{a050}\\
& <\infty,\nonumber
\end{align}

by Theorem \ref{Thm_wt_func_mdim}.

Thus if $\left\vert x\right\vert \rightarrow\infty$ then $G_{c}\left(
x\right)  -p_{c}\left(  x\right)  \rightarrow0$and hence $p_{c}=0$.

From \ref{a02} and \ref{a050}:
\begin{align}
If\text{ }\left\vert x\right\vert \geq2lR, & \nonumber\\
\left\vert G_{c}\left(  x\right)  \right\vert  &  \leq\frac{2\sum
\limits_{k=0}^{l}\tbinom{2l}{k}\left(  l-k\right)  ^{2m}}{\left(  2\pi\right)
^{d}c_{n,d}}\max\limits_{j=-l,j\neq0}^{l}\frac{k_{2m,2m-1}}{\left\vert
x\right\vert -lR}\frac{1}{2m}\int\left\vert y\right\vert ^{2m}q\left(
y\right)  dy\nonumber\\
&  =\frac{\sum\limits_{k=0}^{l}\tbinom{2l}{k}\left(  l-k\right)  ^{2m}%
}{\left(  2\pi\right)  ^{d}c_{n,d}}\frac{k_{2m,2m-1}}{m}\left(  \int\left\vert
y\right\vert ^{2m}q\left(  y\right)  dy\right)  \frac{1}{\left\vert
x\right\vert -lR},\label{a064}%
\end{align}

where $2m=2n-d+1$, $c_{n,d}$ is given by Lemma
\ref{Lem_eval_thin_plate_splin_const} and $k_{2m,2m-1}$ is given by part 3 of
Lemma \ref{Lem_thin_plate_splin}.\medskip

\underline{\textbf{Subcase }$q$ is radial,} say $q\left(  y\right)  =q_{\odot
}\left(  \left\vert y\right\vert \right)  $. Starting with the right side of
\ref{a04}, write%
\begin{align}
\int\left(  \int_{0}^{1}\frac{\left(  1-s\right)  ^{2m-1}}{\left\vert
x-sjy\right\vert }ds\right)  \left\vert y\right\vert ^{2m}q\left(  y\right)
dy  & \leq\int\int_{0}^{1}\frac{ds}{\left\vert x-sjy\right\vert }\left\vert
y\right\vert ^{2m}q_{\odot}\left(  \left\vert y\right\vert \right)
dy\nonumber\\
& =\int_{0}^{1}\int\frac{\left\vert y\right\vert ^{2m}q_{\odot}\left(
\left\vert y\right\vert \right)  dy}{\sqrt{\left\vert x-sjy\right\vert ^{2}}%
}ds\nonumber\\
& =\int_{0}^{1}\int\frac{\left\vert y\right\vert ^{2m}q_{\odot}\left(
\left\vert y\right\vert \right)  dy}{\sqrt{\left\vert x\right\vert
^{2}-2js\left(  xy\right)  +j^{2}s^{2}\left\vert y\right\vert ^{2}}%
}ds,\label{a062}%
\end{align}

and this is an even function of $j$, \textbf{so we need only assume that}
$j>0$.

From Theorem \ref{Thm_Integ_u(xy,|x|)dx},%
\begin{equation}
\int_{\left\vert x\right\vert \leq r}\Phi\left(  \xi x,\left\vert x\right\vert
\right)  dx=\omega_{d-1}\int_{0}^{r}\rho^{d-1}\int_{0}^{\pi}\Phi\left(
\left\vert \xi\right\vert \rho\cos\theta,\rho\right)  \sin^{d-2}\theta d\theta
d\rho,\label{a10}%
\end{equation}

so%
\begin{align*}
\int &  \frac{\left\vert y\right\vert ^{2m}q\left(  y\right)  dy}%
{\sqrt{\left\vert x\right\vert ^{2}-2js\left(  xy\right)  +j^{2}%
s^{2}\left\vert y\right\vert ^{2}}}\\
&  =\omega_{d-1}\int_{0}^{\infty}\rho^{d-1}\int_{0}^{\pi}\frac{\rho
^{2m}q_{\odot}\left(  \rho\right)  \sin^{d-2}\theta}{\sqrt{\left\vert
x\right\vert ^{2}-2js\left\vert x\right\vert \rho\cos\theta+j^{2}s^{2}\rho
^{2}}}d\theta d\rho\\
&  =\omega_{d-1}\int_{0}^{\infty}\rho^{2m+d-1}q_{\odot}\left(  \rho\right)
\int_{0}^{\pi}\frac{\sin^{d-2}\theta d\theta}{\sqrt{\left\vert x\right\vert
^{2}-2js\left\vert x\right\vert \rho\cos\theta+j^{2}s^{2}\rho^{2}}}d\rho\\
&  =\omega_{d-1}\int_{0}^{\infty}\rho^{2m+d-1}q_{\odot}\left(  \rho\right)
\int_{0}^{\pi}\frac{-\sin^{d-3}\theta\text{ }d\cos\theta}{\sqrt{\left\vert
x\right\vert ^{2}-2js\left\vert x\right\vert \rho\cos\theta+j^{2}s^{2}\rho
^{2}}}d\rho\\
&  =\omega_{d-1}\int_{0}^{\infty}\rho^{2m+d-1}q_{\odot}\left(  \rho\right)
\int_{-1}^{1}\frac{\left(  1-t^{2}\right)  ^{d-3}dt}{\sqrt{\left\vert
x\right\vert ^{2}-2js\left\vert x\right\vert \rho t+j^{2}s^{2}\rho^{2}}}%
d\rho\\
&  =\omega_{d-1}\int_{0}^{\infty}\rho^{2m+d-1}q_{\odot}\left(  \rho\right)
\int_{-1}^{1}\frac{\left(  1-t^{2}\right)  ^{d-3}dt}{\sqrt{\left\vert
x\right\vert ^{2}-2js\left\vert x\right\vert \rho t+j^{2}s^{2}\rho^{2}}}%
d\rho\\
&  \leq\omega_{d-1}\int_{0}^{\infty}\rho^{2m+d-1}q_{\odot}\left(  \rho\right)
\int_{-1}^{1}\frac{dt}{\sqrt{\left\vert x\right\vert ^{2}-2js\left\vert
x\right\vert \rho t+j^{2}s^{2}\rho^{2}}}d\rho\\
&  =\omega_{d-1}\int_{0}^{\infty}\rho^{2m+d-1}q_{\odot}\left(  \rho\right)
\int_{-1}^{1}\frac{dt}{\sqrt{B-At}}d\rho,
\end{align*}

where $A=-2js\left\vert x\right\vert \rho$ and $B=\left\vert x\right\vert
^{2}+j^{2}s^{2}\rho^{2}$. Noting that $B+A=\left(  \left\vert x\right\vert
-js\rho\right)  ^{2}$ and $B-A=\left(  \left\vert x\right\vert +js\rho\right)
^{2}$ we continue:%
\begin{align*}
\int_{-1}^{1}\frac{dt}{\sqrt{B-At}}=\left[  -\frac{2}{A}\sqrt{B-At}\right]
_{-1}^{1}  & =-\frac{2}{A}\left(  \sqrt{B-A}-\sqrt{B+A}\right) \\
& =\frac{\left(  \left\vert x\right\vert +js\rho\right)  -\left\vert
\left\vert x\right\vert -js\rho\right\vert }{js\left\vert x\right\vert \rho},
\end{align*}

which is,as noted above, an even function of $j$. Hence%
\[
\int\frac{\left\vert y\right\vert ^{2m}q\left(  y\right)  dy}{\sqrt{\left\vert
x\right\vert ^{2}-2js\left(  xy\right)  +j^{2}s^{2}\left\vert y\right\vert
^{2}}}=\omega_{d-1}\int_{0}^{\infty}\frac{\left(  \left\vert x\right\vert
+js\rho\right)  -\left\vert \left\vert x\right\vert -js\rho\right\vert
}{js\left\vert x\right\vert \rho}\rho^{2m+d-1}q_{\odot}\left(  \rho\right)
d\rho,
\]

and thus
\[
\int_{0}^{1}\int\frac{y^{2m}q_{\odot}\left(  \left\vert y\right\vert \right)
dyds}{\sqrt{\left\vert x\right\vert ^{2}-2js\left(  xy\right)  +j^{2}%
s^{2}\left\vert y\right\vert ^{2}}}=\omega_{d-1}\int_{0}^{1}\int_{0}^{\infty
}\frac{\left(  \left\vert x\right\vert +js\rho\right)  -\left\vert \left\vert
x\right\vert -js\rho\right\vert }{js\left\vert x\right\vert \rho}\rho
^{2m+d-1}q_{\odot}\left(  \rho\right)  d\rho ds.
\]

\begin{align*}
\int_{0}^{1}\int &  \frac{y^{2m}q_{\odot}\left(  \left\vert y\right\vert
\right)  dyds}{\sqrt{\left\vert x\right\vert ^{2}-2js\left(  xy\right)
+j^{2}s^{2}\left\vert y\right\vert ^{2}}}\\
&  =\omega_{d-1}\int_{0}^{1}\int_{0}^{\infty}\frac{\left(  j^{-1}\left\vert
x\right\vert s^{-1}+\rho\right)  -\left\vert j^{-1}\left\vert x\right\vert
s^{-1}-\rho\right\vert }{\left\vert x\right\vert \rho}\rho^{2m+d-1}q_{\odot
}\left(  \rho\right)  d\rho ds\\
&  =\omega_{d-1}\int_{0}^{1}\int_{0}^{\infty}\frac{\left(  j^{-1}\left\vert
x\right\vert s^{-1}+\rho\right)  -\left\vert j^{-1}\left\vert x\right\vert
s^{-1}-\rho\right\vert }{\left\vert x\right\vert \rho}\rho^{2n}q_{\odot
}\left(  \rho\right)  d\rho ds.
\end{align*}

Now apply the change of variables: $t=j^{-1}\left\vert x\right\vert s^{-1}$,
$ds=\left\vert j\right\vert ^{-1}\left\vert x\right\vert t^{-2}dt$ so that%
\begin{align*}
\int_{0}^{1}\int &  \frac{\left\vert y\right\vert ^{2m}q_{\odot}\left(
\left\vert y\right\vert \right)  dyds}{\sqrt{\left\vert x\right\vert
^{2}-2js\left(  xy\right)  +j^{2}s^{2}\left\vert y\right\vert ^{2}}}\\
&  =\omega_{d-1}\int_{j^{-1}\left\vert x\right\vert }^{\infty}\int_{0}%
^{\infty}\frac{t+\rho-\left\vert t-\rho\right\vert }{\left\vert x\right\vert
\rho}\rho^{2n}q_{\odot}\left(  \rho\right)  d\rho\text{ }\left\vert
j\right\vert ^{-1}\left\vert x\right\vert t^{-2}dt\\
&  =\omega_{d-1}j^{-1}\int_{j^{-1}\left\vert x\right\vert }^{\infty}\int%
_{0}^{\infty}\frac{t+\rho-\left\vert t-\rho\right\vert }{t^{2}}\rho
^{2n}q_{\odot}\left(  \rho\right)  d\rho dt.
\end{align*}

Because of the absolute value term we need to split the domain of integration
using the curve $\rho=t$. Thus%
\begin{align*}
&  \int_{0}^{1}\int\frac{y^{2m}q_{\odot}\left(  \left\vert y\right\vert
\right)  \text{ }dyds}{\sqrt{\left\vert x\right\vert ^{2}-2js\left(
xy\right)  +j^{2}s^{2}\left\vert y\right\vert ^{2}}}\\
&  =\omega_{d-1}j^{-1}\int_{j^{-1}\left\vert x\right\vert }^{\infty}\int%
_{t}^{\infty}\frac{t+\rho-\left\vert t-\rho\right\vert }{t^{2}}\rho
^{2n}q_{\odot}\left(  \rho\right)  d\rho dt+\\
&  \qquad+\omega_{d-1}j^{-1}\int_{j^{-1}\left\vert x\right\vert }^{\infty}%
\int_{0}^{t}\frac{t+\rho-\left\vert t-\rho\right\vert }{t^{2}}\rho
^{2n}q_{\odot}\left(  \rho\right)  d\rho dt\\
&  =\omega_{d-1}j^{-1}\int_{j^{-1}\left\vert x\right\vert }^{\infty}\int%
_{t}^{\infty}\frac{t+\rho-\left(  \rho-t\right)  }{t^{2}}\rho^{2n}q_{\odot
}\left(  \rho\right)  d\rho dt+\\
&  \qquad+\omega_{d-1}j^{-1}\int_{j^{-1}\left\vert x\right\vert }^{\infty}%
\int_{0}^{t}\frac{t+\rho-\left(  t-\rho\right)  }{t^{2}}\rho^{2n}q_{\odot
}\left(  \rho\right)  d\rho dt\\
&  =\omega_{d-1}j^{-1}\int_{j^{-1}\left\vert x\right\vert }^{\infty}\int%
_{t}^{\infty}\frac{2}{t}\rho^{2n}q_{\odot}\left(  \rho\right)  d\rho
dt+\omega_{d-1}j^{-1}\int_{j^{-1}\left\vert x\right\vert }^{\infty}\int%
_{0}^{t}\frac{2\rho}{t^{2}}\rho^{2n}q_{\odot}\left(  \rho\right)  d\rho dt\\
&  =\omega_{d-1}j^{-1}\int_{j^{-1}\left\vert x\right\vert }^{\infty}\frac
{2}{t}\int_{t}^{\infty}\rho^{2n}q_{\odot}\left(  \rho\right)  d\rho
dt+\omega_{d-1}j^{-1}\int_{j^{-1}\left\vert x\right\vert }^{\infty}\frac
{2}{t^{2}}\int_{0}^{t}\rho^{2n+1}q_{\odot}\left(  \rho\right)  d\rho dt\\
&  <\ldots+\omega_{d-1}\left\vert j\right\vert ^{-1}\int_{\left\vert
j\right\vert ^{-1}\left\vert x\right\vert }^{\infty}\frac{2}{t^{2}}\int%
_{0}^{\infty}\rho^{2n+1}q_{\odot}\left(  \rho\right)  d\rho dt\\
&  =\ldots+\omega_{d-1}\left\vert j\right\vert ^{-1}\frac{2}{\left\vert
j\right\vert ^{-1}\left\vert x\right\vert }\int_{0}^{\infty}\rho
^{2n+1}q_{\odot}\left(  \rho\right)  d\rho\\
&  =\omega_{d-1}\left\vert j\right\vert ^{-1}\int_{\left\vert j\right\vert
^{-1}\left\vert x\right\vert }^{\infty}\frac{2}{t}\int_{t}^{\infty}\rho
^{2n}q_{\odot}\left(  \rho\right)  d\rho dt+\frac{2\omega_{d-1}}{\left\vert
x\right\vert }\int_{0}^{\infty}\rho^{2n+1}q_{\odot}\left(  \rho\right)  d\rho.
\end{align*}

Regarding the first iterated integral:%
\begin{align*}
\int_{\left\vert j\right\vert ^{-1}\left\vert x\right\vert }^{\infty} &
\frac{2}{t}\int_{t}^{\infty}\rho^{2n}q_{\odot}\left(  \rho\right)  d\rho dt\\
&  =\int_{\left\vert j\right\vert ^{-1}\left\vert x\right\vert }^{\infty}%
\int_{\left\vert j\right\vert ^{-1}\left\vert x\right\vert }^{\rho}\frac{2}%
{t}\rho^{2n}q_{\odot}\left(  \rho\right)  dtd\rho\\
&  =2\int_{\left\vert j\right\vert ^{-1}\left\vert x\right\vert }^{\infty
}\left(  \int_{\left\vert j\right\vert ^{-1}\left\vert x\right\vert }^{\rho
}\frac{dt}{t}\right)  \rho^{2n}q_{\odot}\left(  \rho\right)  d\rho\\
&  =2\int_{\left\vert j\right\vert ^{-1}\left\vert x\right\vert }^{\infty
}\left(  \ln\rho-\ln\left(  \left\vert j\right\vert ^{-1}\left\vert
x\right\vert \right)  \right)  \rho^{2n}q_{\odot}\left(  \rho\right)  d\rho\\
&  =2\int_{\left\vert j\right\vert ^{-1}\left\vert x\right\vert }^{\infty}%
\rho^{2n}\left(  \ln\rho\right)  q_{\odot}\left(  \rho\right)  d\rho
-2\ln\left(  \left\vert j\right\vert ^{-1}\left\vert x\right\vert \right)
\int_{\left\vert j\right\vert ^{-1}\left\vert x\right\vert }^{\infty}\rho
^{2n}q_{\odot}\left(  \rho\right)  d\rho.
\end{align*}

If $\left\vert x\right\vert \geq\left\vert l\right\vert $ then $\left\vert
j\right\vert ^{-1}\left\vert x\right\vert \geq1$ and%
\begin{align*}
2\int_{\left\vert j\right\vert ^{-1}\left\vert x\right\vert }^{\infty}%
\rho^{2n}\left(  \ln\rho\right)  q_{\odot}\left(  \rho\right)  d\rho &
=2\int_{\left\vert j\right\vert ^{-1}\left\vert x\right\vert }^{\infty}%
\rho^{2n+1}\frac{\ln\rho}{\rho}q_{\odot}\left(  \rho\right)  d\rho\\
& <2\int_{\left\vert j\right\vert ^{-1}\left\vert x\right\vert }^{\infty}%
\rho^{2n}q_{\odot}\left(  \rho\right)  d\rho\\
& <2\int_{\left\vert j\right\vert ^{-1}\left\vert x\right\vert }^{\infty}%
\rho^{2n}q_{\odot}\left(  \rho\right)  d\rho\\
& <\infty,
\end{align*}

by Theorem \ref{Thm_CentDiffQrad}.

Thus when $\left\vert x\right\vert \geq\left\vert l\right\vert $,%
\begin{align}
\int_{0}^{1}\int\frac{\left\vert y\right\vert ^{2m}q_{\odot}\left(  \left\vert
y\right\vert \right)  \text{ }dyds}{\sqrt{\left\vert x\right\vert
^{2}-2js\left(  xy\right)  +j^{2}s^{2}\left\vert y\right\vert ^{2}}}%
<2\int_{\left\vert j\right\vert ^{-1}\left\vert x\right\vert }^{\infty}  &
\rho^{2n}q_{\odot}\left(  \rho\right)  d\rho+2\ln\left(  \left\vert
j\right\vert ^{-1}\left\vert x\right\vert \right)  \int_{\left\vert
j\right\vert ^{-1}\left\vert x\right\vert }^{\infty}\rho^{2n}q_{\odot}\left(
\rho\right)  d\rho+\nonumber\\
& +\frac{2\omega_{d-1}}{\left\vert x\right\vert }\int_{0}^{\infty}\rho
^{2n+1}q_{\odot}\left(  \rho\right)  d\rho.\label{a061}%
\end{align}

If we assume that
\begin{equation}
\left.
\begin{array}
[c]{l}%
\int_{0}^{\infty}\rho^{2n+1}q_{\odot}\left(  \rho\right)  d\rho<\infty,\\
i.e.\text{ }\left\vert \cdot\right\vert ^{2n-d+2}q\in L^{1}.
\end{array}
\right\} \label{a065}%
\end{equation}

instead of the inequality \ref{a931} then \textbf{the first and third terms
clearly tend to zero} as $\left\vert x\right\vert \rightarrow0$. The behavior
of \textbf{the second integral} hinges on the behavior of the limit%
\begin{equation}
\lim_{a\rightarrow\infty}\left(  \left(  \ln a\right)  \int_{a}^{\infty}%
\rho^{2n}q_{\odot}\left(  \rho\right)  d\rho\right)  .\label{a05}%
\end{equation}

However%
\begin{align*}
\lim_{a\rightarrow\infty}\left(  \left(  \ln a\right)  \int_{a}^{\infty}%
\rho^{2n}q_{\odot}\left(  \rho\right)  d\rho\right)   & =\lim_{a\rightarrow
\infty}\left(  \ln a\right)  \int_{a}^{\infty}\frac{1}{\rho}\rho
^{2n+1}q_{\odot}\left(  \rho\right)  d\rho\\
& \leq\lim_{a\rightarrow\infty}\left(  \frac{\ln a}{a}\int_{a}^{\infty}%
\rho^{2n+1}q_{\odot}\left(  \rho\right)  d\rho\right) \\
& =0.
\end{align*}

Thus $G_{c}\left(  \infty\right)  =0$ and so $p_{c}=0$.

Combining the successive estimates \ref{a02}, \ref{a04}, \ref{a062} and
\ref{a061} we get:
\begin{align}
&  When\text{ }\left\vert x\right\vert \geq l:\nonumber\\
&  \ \nonumber\\
&  \left\vert G_{c}\left(  x\right)  \right\vert \nonumber\\
&  \leq\frac{2\sum\limits_{k=0}^{l}\tbinom{2l}{k}\left(  l-k\right)  ^{2m}%
}{\left(  2\pi\right)  ^{d}c_{n,d}}\max\limits_{j=-l,j\neq0}^{l}\int\int%
_{0}^{1}\frac{\left(  1-s\right)  ^{2m-1}}{\left(  2m-1\right)  !}\left\vert
\left(  \left(  \widehat{y}D\right)  ^{2m}T_{n}\right)  \left(  x-sjy\right)
\right\vert ds\text{ }\left\vert y\right\vert ^{2m}q\left(  y\right)
dy\nonumber\\
&  \leq\frac{2\sum\limits_{k=0}^{l}\tbinom{2l}{k}\left(  l-k\right)  ^{2m}%
}{\left(  2\pi\right)  ^{d}c_{n,d}}\frac{k_{2m,2m-1}}{\left(  2m-1\right)
!}\max\limits_{j=-l,j\neq0}^{l}\int\left(  \int_{0}^{1}\frac{\left(
1-s\right)  ^{2m-1}}{\left\vert x-sjy\right\vert }ds\right)  \left\vert
y\right\vert ^{2m}q\left(  y\right)  dy\nonumber\\
&  \leq\frac{2\sum\limits_{k=0}^{l}\tbinom{2l}{k}\left(  l-k\right)  ^{2m}%
}{\left(  2\pi\right)  ^{d}c_{n,d}}\frac{k_{2m,2m-1}}{\left(  2m-1\right)
!}\max\limits_{j=-l,j\neq0}^{l}\int_{0}^{1}\int\frac{\left\vert y\right\vert
^{2m}q_{\odot}\left(  \left\vert y\right\vert \right)  dy}{\sqrt{\left\vert
x\right\vert ^{2}-2js\left(  xy\right)  +j^{2}s^{2}\left\vert y\right\vert
^{2}}}ds\nonumber\\
&  \leq\frac{2\sum\limits_{k=0}^{l}\tbinom{2l}{k}\left(  l-k\right)  ^{2m}%
}{\left(  2\pi\right)  ^{d}c_{n,d}}\frac{k_{2m,2m-1}}{\left(  2m-1\right)
!}\max\limits_{j=-l,j\neq0}^{l}\left(
\begin{array}
[c]{c}%
2\int_{\left\vert x\right\vert /j}^{\infty}\rho^{2n}q_{\odot}\left(
\rho\right)  d\rho+2\ln\left(  \left\vert j\right\vert ^{-1}\left\vert
x\right\vert \right)  \int_{\left\vert x\right\vert /j}^{\infty}\rho
^{2n}q_{\odot}\left(  \rho\right)  d\rho+\\
+\frac{2\omega_{d-1}}{\left\vert x\right\vert }\int_{0}^{\infty}\rho
^{2n+1}q_{\odot}\left(  \rho\right)  d\rho
\end{array}
\right) \nonumber\\
&  \leq\frac{4\sum\limits_{k=0}^{l}\tbinom{2l}{k}\left(  l-k\right)  ^{2m}%
}{\left(  2\pi\right)  ^{d}c_{n,d}}\frac{k_{2m,2m-1}}{\left(  2m-1\right)
!}\left(
\begin{array}
[c]{c}%
\int_{\left\vert x\right\vert /l}^{\infty}\rho^{2n}q_{\odot}\left(
\rho\right)  d\rho+\ln\left(  \left\vert x\right\vert /l\right)
\int_{\left\vert x\right\vert /l}^{\infty}\rho^{2n}q_{\odot}\left(
\rho\right)  d\rho+\\
+\frac{\omega_{d-1}}{\left\vert x\right\vert }\int_{0}^{\infty}\rho
^{2n+1}q_{\odot}\left(  \rho\right)  d\rho
\end{array}
\right)  ,\label{a063}%
\end{align}

where $2m=2n-d+1$, $c_{n,d}$ is given by Lemma
\ref{Lem_eval_thin_plate_splin_const} and $k_{2m,2m-1}$ is given by part 3 of
Lemma \ref{Lem_thin_plate_splin}.

We have now almost demonstrated:

\begin{theorem}
\label{Thm_CentDiff_eqn_thinps}\fbox{\textbf{Suppose }$d$ \textbf{is odd}.}
Suppose the central difference basis function $G_{c}$ is generated by the
function $q$ and has parameters $n$ and $l$. If $q$ has bounded support, or
that $q$ is radial and satisfies condition \ref{a065} then
\begin{align}
&  G_{c}\left(  x\right) \nonumber\\
&  :=\frac{1}{\left(  2\pi\right)  ^{d}c_{n,d}}\sum\limits_{j=-l}^{l}\left(
-1\right)  ^{j}\tbinom{2l}{j+l}\int T_{n}\left(  x-jy\right)  q\left(
y\right)  dy\text{ }\left(  a\text{ }copy\text{ }of\text{ \ref{p36}}\right)
\label{a09}\\
&  =\frac{1}{\left(  2\pi\right)  ^{d/2}c_{n,d}}\left(  \tbinom{2l}%
{l}\widehat{q}\left(  0\right)  T_{n}\left(  x\right)  +\sum
\limits_{j=-l,j\neq0}^{l}\left(  -1\right)  ^{j}\tbinom{2l}{j+l}\left\vert
j\right\vert ^{2n-d}\left(  T_{n}\ast q\right)  \left(  \frac{x}{j}\right)
\right) \label{a011}\\
&  =\frac{1}{\left(  2\pi\right)  ^{d/2}c_{n,d}}\sum\limits_{j=-l,j\neq0}%
^{l}\left(  -1\right)  ^{j}\tbinom{2l}{j+l}j^{2m}\int\left(  \int_{0}^{1}%
\frac{\left(  1-s\right)  ^{2m-1}}{\left(  2m-1\right)  !}\left(  \left(
\widehat{y}D\right)  ^{2m}T_{n}\right)  \left(  x-sjy\right)  ds\right)
\left\vert y\right\vert ^{2m}q\left(  y\right)  dy,\label{a016}%
\end{align}

where $2m=2n-d+1\geq2$ and $T_{n}$ is the thin plate spline formula given in
\ref{p24}.

The integral in \ref{a016} is absolutely convergent and also:%
\begin{align*}
&  G_{c}\left(  x\right) \\
&  =\frac{1}{\left(  2\pi\right)  ^{d/2}c_{n,d}}\sum\limits_{\substack{j=-l
\\j\neq0}}^{l}\left(  -1\right)  ^{j}\tbinom{2l}{j+l}\left\vert j\right\vert
^{2m-1}\left(  \int_{0}^{1}\frac{\left(  1-s\right)  ^{2m-1}}{\left(
2m-1\right)  !}\frac{1}{s}\int\left(  \left(  \widehat{y}D\right)  ^{2m}%
T_{n}\right)  \left(  \frac{x}{sj}-y\right)  \left\vert y\right\vert
^{2m}q\left(  y\right)  dyds\right)  ,
\end{align*}

with this convolution integral also being absolutely convergent.

Finally, if $q$ has bounded support in $B_{R_{q}}$ then $G_{c}$ is bounded
near infinity by \ref{a064}. If $q$ is radial and satisfies condition
\ref{a065} then $G_{c}$ is bounded near infinity by \ref{a063}.
\end{theorem}

\begin{proof}
To complete the proof we define%
\[
H_{j}\left(  x\right)  :=\int T_{n}\left(  x-jy\right)  q\left(  y\right)  dy,
\]

then%
\[
H_{j}\left(  jx\right)  =\int T_{n}\left(  jx-jy\right)  q\left(  y\right)
dy=\left\vert j\right\vert ^{2n-d}\int T_{n}\left(  x-y\right)  q\left(
y\right)  dy=\left\vert j\right\vert ^{2n-d}H_{1}\left(  x\right)  ,
\]

and so $H_{j}\left(  x\right)  =\left\vert j\right\vert ^{2n-d}H_{1}\left(
\frac{x}{j}\right)  $. Hence, starting with \ref{p30},%
\begin{align*}
& G_{c}\left(  x\right)  -p_{c}\left(  x\right) \\
& =\frac{1}{\left(  2\pi\right)  ^{d}c_{n,d}}\left(  \tbinom{2l}{l}\left(
\int q\right)  T_{n}\left(  x\right)  +\sum\limits_{j=-l,j\neq0}^{l}%
\frac{\left(  -1\right)  ^{j}}{\left\vert j\right\vert ^{d}}\tbinom{2l}%
{j+l}\int T_{n}\left(  x-y\right)  q\left(  \frac{y}{j}\right)  dy\right) \\
& =\frac{1}{\left(  2\pi\right)  ^{d}c_{n,d}}\left(  \tbinom{2l}{l}\left(
\int q\right)  T_{n}\left(  x\right)  +\sum\limits_{j=-l,j\neq0}^{l}\left(
-1\right)  ^{j}\tbinom{2l}{j+l}\int T_{n}\left(  x-jz\right)  q\left(
z\right)  dz\right) \\
& =\frac{1}{\left(  2\pi\right)  ^{d}c_{n,d}}\left(  \tbinom{2l}{l}\left(
\int q\right)  T_{n}\left(  x\right)  +\sum\limits_{j=-l,j\neq0}^{l}\left(
-1\right)  ^{j}\tbinom{2l}{j+l}H_{j}\left(  x\right)  \right) \\
& =\frac{1}{\left(  2\pi\right)  ^{d}c_{n,d}}\left(  \tbinom{2l}{l}\left(
\int q\right)  T_{n}\left(  x\right)  +\sum\limits_{j=-l,j\neq0}^{l}\left(
-1\right)  ^{j}\tbinom{2l}{j+l}\left\vert j\right\vert ^{2n-d}\int
T_{n}\left(  \frac{x}{j}-y\right)  q\left(  y\right)  dy\right) \\
& =\frac{1}{\left(  2\pi\right)  ^{d/2}c_{n,d}}\left(  \tbinom{2l}%
{l}\widehat{q}\left(  0\right)  T_{n}\left(  x\right)  +\sum
\limits_{j=-l,j\neq0}^{l}\left(  -1\right)  ^{j}\tbinom{2l}{j+l}\left\vert
j\right\vert ^{2n-d}\left(  T_{n}\ast q\right)  \left(  \frac{x}{j}\right)
\right)  ,
\end{align*}

which is \ref{a011}.

The next equation in this theorem, equation \ref{a016}, is merely a copy of
\ref{a01}.

If $q$ has bounded support then inequality \ref{a050} holds. Hence the
integral \ref{a016} is absolutely convergent and $G_{c}\left(  \infty\right)
=0$ i.e. $p_{c}=0$.

If $q$ is radial and satisfies conditions \ref{a065} then the successive
estimates \ref{a02}, \ref{a04}, \ref{a062} and \ref{a061} imply the integral
\ref{a016} is absolutely convergent and that $G_{c}\left(  \infty\right)  =0$
i.e. $p_{c}=0$.

Because the integral \ref{a016} is absolutely convergent Fubini's theorem
allows us to change the order of integration to get%
\begin{align*}
&  G_{c}\left(  x\right) \\
&  =\frac{1}{\left(  2\pi\right)  ^{d/2}c_{n,d}}\sum\limits_{j=-l,j\neq0}%
^{l}\left(  -1\right)  ^{j}\tbinom{2l}{j+l}j^{2m}\left(  \int_{0}^{1}%
\frac{\left(  1-s\right)  ^{2m-1}}{\left(  2m-1\right)  !}\int\left(  \left(
\widehat{y}D\right)  ^{2m}T_{n}\right)  \left(  x-sjy\right)  \left\vert
y\right\vert ^{2m}q\left(  y\right)  dyds\right)  .
\end{align*}

Now define%
\[
f\left(  x\right)  :=\int\left(  \left(  \widehat{y}D\right)  ^{2m}%
T_{n}\right)  \left(  x-sjy\right)  \left\vert y\right\vert ^{2m}q\left(
y\right)  dy=\int\left(  \left(  yD\right)  ^{2m}T_{n}\right)  \left(
x-sjy\right)  q\left(  y\right)  dy.
\]

Then%
\begin{align*}
f\left(  sjx\right)  =\int\left(  \left(  yD\right)  ^{2m}T_{n}\right)
\left(  sjx-sjy\right)  q\left(  y\right)  dy  & =\left(  s\left\vert
j\right\vert \right)  ^{2n-d-2m}\int\left(  \left(  yD\right)  ^{2m}%
T_{n}\right)  \left(  x-y\right)  q\left(  y\right)  dy\\
& =\frac{1}{s\left\vert j\right\vert }\int\left(  \left(  yD\right)
^{2m}T_{n}\right)  \left(  x-y\right)  q\left(  y\right)  dy,
\end{align*}

which implies that%
\[
f\left(  x\right)  =\frac{1}{s\left\vert j\right\vert }\int\left(  \left(
yD\right)  ^{2m}T_{n}\right)  \left(  \frac{x}{sj}-y\right)  q\left(
y\right)  dy=\frac{1}{s\left\vert j\right\vert }\int\left(  \left(
\widehat{y}D\right)  ^{2m}T_{n}\right)  \left(  \frac{x}{sj}-y\right)
\left\vert y\right\vert ^{2m}q\left(  y\right)  dy,
\]

and hence that%
\begin{align}
&  G_{c}\left(  x\right) \nonumber\\
&  =\frac{1}{\left(  2\pi\right)  ^{d/2}c_{n,d}}\sum\limits_{j=-l,j\neq0}%
^{l}\left(  -1\right)  ^{j}\tbinom{2l}{j+l}j^{2m}\left(  \int_{0}^{1}%
\frac{\left(  1-s\right)  ^{2m-1}}{\left(  2m-1\right)  !}\frac{1}{s\left\vert
j\right\vert }\int\left(  \left(  \widehat{y}D\right)  ^{2m}T_{n}\right)
\left(  \frac{x}{sj}-y\right)  \left\vert y\right\vert ^{2m}q\left(  y\right)
dyds\right)  .\label{a051}%
\end{align}

\end{proof}

?? Remember \ref{p41}.

\begin{corollary}
\label{Cor_Thm_CentDiff_eqn_thinps}\fbox{\textbf{Suppose }$d$ \textbf{is
odd}.} Suppose the conditions of Theorem \ref{Thm_CentDiff_eqn_thinps} hold
and that $q$ \textbf{is radial}, say $q\left(  x\right)  =q_{\odot}\left(
\left\vert x\right\vert \right)  $. Then $G_{c}$ is a radial basis function
and \ref{a011} can be written%
\begin{equation}
G_{c\odot}\left(  r\right)  =\frac{1}{\left(  2\pi\right)  ^{\frac{d}{2}%
}c_{n,d}}\left(  \frac{\tbinom{2l}{l}}{\left(  2\pi\right)  ^{d/2}}\left(
\int q\right)  \left(  T_{n}\right)  _{\odot}\left(  r\right)  +2\sum
\limits_{j=1}^{l}\left(  -1\right)  ^{j}\tbinom{2l}{j+l}j^{2n-d}\left(
T_{n}\ast q\right)  _{\odot}\left(  \frac{r}{j}\right)  \right)  ,\quad
r>0,\label{a054}%
\end{equation}

and%
\begin{equation}
\int q=\omega_{d}\int_{0}^{\infty}t^{d-1}q_{\odot}\left(  t\right)
dt.\label{a049}%
\end{equation}

\end{corollary}

\begin{proof}
From Definition \ref{Def_multi_centdiff_wt_fn} $q$ radial implies $G_{c}$ is radial.

From \ref{a011},%
\[
G_{c}\left(  x\right)  =\frac{1}{\left(  2\pi\right)  ^{d/2}c_{n,d}}\left(
\frac{\tbinom{2l}{l}}{\left(  2\pi\right)  ^{d/2}}\left(  \int q\right)
T_{n}\left(  x\right)  +\sum\limits_{j=-l,j\neq0}^{l}\left(  -1\right)
^{j}\tbinom{2l}{j+l}\left\vert j\right\vert ^{2n-d}\left(  T_{n}\ast q\right)
\left(  \frac{x}{j}\right)  \right)  .
\]

Since $T_{n}$ and $q$ are radial $T_{n}\ast q$ is also radial. Thus%
\begin{align*}
G_{c\odot}\left(  r\right)   & =\frac{1}{\left(  2\pi\right)  ^{d/2}c_{n,d}%
}\left(  \frac{\tbinom{2l}{l}}{\left(  2\pi\right)  ^{d/2}}\left(  \int
q\right)  \left(  T_{n}\right)  _{\odot}\left(  r\right)  +\sum
\limits_{j=-l,j\neq0}^{l}\left(  -1\right)  ^{j}\tbinom{2l}{j+l}\left\vert
j\right\vert ^{2n-d}\left(  T_{n}\ast q\right)  _{\odot}\left(  \frac
{r}{\left\vert j\right\vert }\right)  \right) \\
& =\frac{1}{\left(  2\pi\right)  ^{d/2}c_{n,d}}\left(
\begin{array}
[c]{r}%
\frac{\tbinom{2l}{l}}{\left(  2\pi\right)  ^{d/2}}\left(  \int q\right)
\left(  T_{n}\right)  _{\odot}\left(  r\right)  +\sum\limits_{j=1}^{l}\left(
-1\right)  ^{j}\tbinom{2l}{j+l}\left\vert j\right\vert ^{2n-d}\left(
T_{n}\ast q\right)  _{\odot}\left(  \frac{r}{\left\vert j\right\vert }\right)
+\\
+\sum\limits_{j=-l}^{-1}\left(  -1\right)  ^{j}\tbinom{2l}{j+l}\left\vert
j\right\vert ^{2n-d}\left(  T_{n}\ast q\right)  _{\odot}\left(  \frac
{r}{\left\vert j\right\vert }\right)
\end{array}
\right) \\
& =\frac{1}{\left(  2\pi\right)  ^{d/2}c_{n,d}}\left(
\begin{array}
[c]{r}%
\frac{\tbinom{2l}{l}}{\left(  2\pi\right)  ^{d/2}}\left(  \int q\right)
\left(  T_{n}\right)  _{\odot}\left(  r\right)  +\sum\limits_{j=1}^{l}\left(
-1\right)  ^{j}\tbinom{2l}{j+l}\left\vert j\right\vert ^{2n-d}\left(
T_{n}\ast q\right)  _{\odot}\left(  \frac{r}{\left\vert j\right\vert }\right)
+\\
+\sum\limits_{j=1}^{l}\left(  -1\right)  ^{-j}\tbinom{2l}{-j+l}\left\vert
j\right\vert ^{2n-d}\left(  T_{n}\ast q\right)  _{\odot}\left(  \frac
{r}{\left\vert j\right\vert }\right)
\end{array}
\right) \\
& =\frac{1}{\left(  2\pi\right)  ^{d/2}c_{n,d}}\left(  \frac{\tbinom{2l}{l}%
}{\left(  2\pi\right)  ^{d/2}}\left(  \int q\right)  \left(  T_{n}\right)
_{\odot}\left(  r\right)  +2\sum\limits_{j=1}^{l}\left(  -1\right)
^{j}\tbinom{2l}{j+l}j^{2n-d}\left(  T_{n}\ast q\right)  _{\odot}\left(
\frac{r}{j}\right)  \right)  .
\end{align*}

Finally, since $q$ is radial, Corollary \ref{Cor_Thm_Integ_u(xy)f(|x|)dx}
implies $\int q=\omega_{d}\int_{0}^{\infty}t^{d-1}q_{\odot}\left(  t\right)
dt$.
\end{proof}

\begin{remark}
??? Because $T_{n}\notin L^{1}$ we cannot use Corollary
\ref{Cor_Thm_convol_rad_fn} or Theorem \ref{Thm_convol_rad_fn} to calculate
the convolution $T_{n}\ast q$ in \ref{a054}.

?? CHECK! However, by writing $T_{n}=\phi_{0}T_{n}+\phi_{\infty}T_{n}$ where
$\left\{  \phi_{0},\phi_{\infty}\right\}  $ is radial $C^{\infty}$ partition
of unity, and then expanding the terms corresponding to $\phi_{\infty}T_{n}$
using Taylor series so that the remainder is $L^{1}$ and then assuming that
$2n+1<2l$, to enable the polynomial part of the Taylor series to
be\ eliminated, we may be able to use Corollary \ref{Cor_Thm_convol_rad_fn} or
Theorem \ref{Thm_convol_rad_fn} to calculate the resulting convolutions. This
will work if $q$ has compact support or $\left\vert .\right\vert ^{k}q\in
L^{1}$ for sufficiently large $k$.
\end{remark}

\subsection{Case: $d$ is even\label{SbSect_Radial_q_d_even}}

?? To be done. Use the estimate \ref{p02}.

\subsection{?? General case\textbf{\ - unfinished}%
\label{SbSect_CentDiffBasisGeneral}}

We will need to make the assumptions $\int\left\vert \cdot\right\vert
^{2n-d+1}q$ as well as the assumptions \ref{Ap011} or less generally
\ref{Ap037} or \ref{Ap038} when integrating over $U_{2}$. The latter two
assumptions might be related back to the first.

We want to estimate the remainder term%
\[
\int\left(  \int_{0}^{1}\frac{\left(  1-s\right)  ^{m-1}}{\left\vert
x-jsy\right\vert }ds\right)  \left\vert y\right\vert ^{m}q\left(  y\right)
dy.
\]

of the estimate \ref{a04}.

We will use equation 12.41 of Arfken \cite{Arfken70}: if $\left\vert
b\right\vert >\left\vert a\right\vert $ then%
\begin{equation}
\left.
\begin{array}
[c]{l}%
\frac{1}{\left\vert b-a\right\vert }=\frac{1}{\left\vert b\right\vert }%
\sum\limits_{k=0}^{\infty}\left(  \frac{\left\vert a\right\vert }{\left\vert
b\right\vert }\right)  ^{k}P_{k}\left(  \cos\theta\right)  \text{ }is\text{
}uniform\text{ }convergent,\\
\cos\theta=\frac{ab}{\left\vert a\right\vert \left\vert b\right\vert
}=\widehat{a}\widehat{b},\\
The\text{ }P_{n}\text{ }are\text{ }the\text{ }Legendre\text{ }polynomials.\\
P_{0}=1.
\end{array}
\right\} \label{a08}%
\end{equation}

The change of variables $\tau=1/s$ yields:%
\begin{align*}
\int_{0}^{1}\frac{\left(  1-s\right)  ^{m-1}}{\left\vert x-jsy\right\vert }ds
& =\int_{0}^{1}\frac{\left(  1-s\right)  ^{m-1}}{js\left\vert j^{-1}%
s^{-1}x-y\right\vert }ds=\int_{1}^{\infty}\frac{\left(  1-t^{-1}\right)
^{m-1}}{jt^{-1}\left\vert j^{-1}tx-y\right\vert }\frac{dt}{t^{2}}\\
& =\frac{1}{j}\int_{1}^{\infty}\frac{t^{-1}\left(  1-t^{-1}\right)  ^{m-1}%
}{\left\vert j^{-1}tx-y\right\vert }dt=\frac{1}{j}\int_{1}^{\infty}%
\frac{t^{-m}\left(  t-1\right)  ^{m-1}}{\left\vert j^{-1}tx-y\right\vert }dt,
\end{align*}

so that%
\[
\int\left(  \int_{0}^{1}\frac{\left(  1-s\right)  ^{m-1}}{\left\vert
x-jsy\right\vert }ds\right)  \left\vert y\right\vert ^{m}q\left(  y\right)
dy=\frac{1}{j}\int\int_{1}^{\infty}\frac{t^{-m}\left(  t-1\right)  ^{m-1}%
}{\left\vert j^{-1}tx-y\right\vert }dt\text{ }\left\vert y\right\vert
^{m}q\left(  y\right)  dy.
\]

Partition the domain of integration using the conical surfaces%
\begin{align*}
\mu & =j^{-1}\left\vert x\right\vert ,\\
\left\vert y\right\vert  & =\overline{\mu}t,\text{ }\overline{\mu}=\left(
1+\varepsilon\right)  \mu,\\
\left\vert y\right\vert  & =\underline{\mu}t,\text{ }\underline{\mu}=\left(
1-\varepsilon\right)  \mu.
\end{align*}

These two surfaces partition the domain of integration into three volumes;%
\begin{equation}
\left.
\begin{array}
[c]{l}%
U_{1}=\left\{  \left(  y,t\right)  :\left\vert y\right\vert \leq
\underline{\mu}t,\text{ }t\geq1\right\}  ,\\
U_{2}=\left\{  \left(  y,t\right)  :\underline{\mu}t<\left\vert y\right\vert
<\overline{\mu}t,\text{ }t\geq1\right\}  ,\\
U_{3}=\left\{  \left(  y,t\right)  :\left\vert y\right\vert \geq\overline{\mu
}t,\text{ }t\geq1\right\}  .
\end{array}
\right\} \label{a025}%
\end{equation}

Thus%
\begin{align}
\int\left(  \int_{0}^{1}\frac{\left(  1-s\right)  ^{m-1}}{\left\vert
x-jsy\right\vert }ds\right)  \left\vert y\right\vert ^{m}q\left(  y\right)
dy  & =\frac{1}{j}\int\int_{1}^{\infty}\frac{t^{-m}\left(  t-1\right)  ^{m-1}%
}{\left\vert j^{-1}tx-y\right\vert }dt\text{ }\left\vert y\right\vert
^{m}q\left(  y\right)  dy\nonumber\\
& =\frac{1}{j}\sum_{k=1}^{3}\iint_{U_{k}}\frac{t^{-m}\left(  t-1\right)
^{m-1}}{\left\vert j^{-1}tx-y\right\vert }dt\text{ }\left\vert y\right\vert
^{m}q\left(  y\right)  dy,\label{a029}%
\end{align}

The results are the inequalities \ref{a032} and \ref{a033} and \ref{Ap011} or
less generally \ref{Ap037} or \ref{Ap038}.

?? \textbf{Integration over }$U_{2}$\textbf{\ must be checked very, very
carefully}!.

Use \ref{a08} to expand $\frac{1}{\left\vert j^{-1}tx-y\right\vert }$ on
$U_{3}$ and $U_{1}$. To integrate over $U_{2}$ I had to use spherical
coordinates.\medskip

\fbox{\textbf{Integration over} $U_{1}$} On $U_{1}$, $\left\vert y\right\vert
\leq\underline{\mu}t=\left(  1-\varepsilon\right)  \mu t=\left(
1-\varepsilon\right)  j^{-1}t\left\vert x\right\vert $ and so%
\[
\frac{1}{\left\vert j^{-1}tx-y\right\vert }=\frac{1}{j^{-1}t\left\vert
x\right\vert }\sum\limits_{k=0}^{\infty}\left(  \frac{\left\vert y\right\vert
}{j^{-1}t\left\vert x\right\vert }\right)  ^{k}P_{k}\left(  \widehat{x}%
\widehat{y}\right)  ,
\]

with uniform convergence, and so%
\begin{align}
\iint_{U_{1}} &  \frac{t^{-m}\left(  t-1\right)  ^{m-1}}{\left\vert
j^{-1}tx-y\right\vert }dt\text{ }\left\vert y\right\vert ^{m}q\left(
y\right)  dy\nonumber\\
&  =\iint_{U_{1}}t^{-m}\left(  t-1\right)  ^{m-1}\left(  \frac{1}%
{j^{-1}t\left\vert x\right\vert }\sum\limits_{k=0}^{\infty}\left(
\frac{\left\vert y\right\vert }{j^{-1}t\left\vert x\right\vert }\right)
^{k}P_{k}\left(  \widehat{x}\widehat{y}\right)  \right)  dt\text{ }\left\vert
y\right\vert ^{m}q\left(  y\right)  dy\nonumber\\
&  \leq\iint_{U_{1}}t^{-m}\left(  t-1\right)  ^{m-1}\left(  \frac{1}%
{j^{-1}t\left\vert x\right\vert }\sum\limits_{k=0}^{\infty}\left(
\frac{\left\vert y\right\vert }{j^{-1}t\left\vert x\right\vert }\right)
^{k}\left\vert P_{k}\left(  \widehat{x}\widehat{y}\right)  \right\vert
\right)  dt\text{ }\left\vert y\right\vert ^{m}q\left(  y\right)
dy\nonumber\\
&  \leq\iint_{U_{1}}t^{-m}\left(  t-1\right)  ^{m-1}\left(  \frac{1}%
{j^{-1}t\left\vert x\right\vert }\sum\limits_{k=0}^{\infty}\left(
\frac{\left\vert y\right\vert }{j^{-1}t\left\vert x\right\vert }\right)
^{k}\right)  dt\text{ }\left\vert y\right\vert ^{m}q\left(  y\right)
dy\nonumber\\
&  \leq\iint_{U_{1}}t^{-m}\left(  t-1\right)  ^{m-1}\left(  \frac{1}%
{j^{-1}t\left\vert x\right\vert }\sum\limits_{k=0}^{\infty}\left(
1-\varepsilon\right)  ^{k}\right)  dt\text{ }\left\vert y\right\vert
^{m}q\left(  y\right)  dy\nonumber\\
&  =\frac{1}{\varepsilon}\iint_{U_{1}}t^{-m}\left(  t-1\right)  ^{m-1}\frac
{1}{j^{-1}t\left\vert x\right\vert }dt\text{ }\left\vert y\right\vert
^{m}q\left(  y\right)  dy\nonumber\\
&  =\frac{1}{\varepsilon j^{-1}\left\vert x\right\vert }\iint_{U_{1}%
}t^{-\left(  m+1\right)  }\left(  t-1\right)  ^{m-1}dt\text{ }\left\vert
y\right\vert ^{m}q\left(  y\right)  dy\nonumber\\
&  \leq\frac{j}{\varepsilon\left\vert x\right\vert }\left(  \int_{1}^{\infty
}t^{-\left(  m+1\right)  }\left(  t-1\right)  ^{m-1}dt\right)  \left(
\int\left\vert y\right\vert ^{2n-d+1}q\left(  y\right)  dy\right)
\label{a032}\\
&  <\infty.\nonumber
\end{align}

We know the basis function is continuous so we need only obtain an inequality
valid for say $\left\vert x\right\vert \geq1$.\medskip

\fbox{\textbf{Integration over} $U_{3}$} On $U_{3}$, $\left\vert y\right\vert
\geq\overline{\mu}t=\left(  1+\varepsilon\right)  \mu t=\left(  1+\varepsilon
\right)  j^{-1}t\left\vert x\right\vert $ and so%
\[
\frac{1}{\left\vert j^{-1}tx-y\right\vert }=\frac{1}{\left\vert y\right\vert
}\sum\limits_{k=0}^{\infty}\left(  \frac{j^{-1}t\left\vert x\right\vert
}{\left\vert y\right\vert }\right)  ^{k}P_{k}\left(  \widehat{x}%
\widehat{y}\right)  ,
\]

and so%
\begin{align}
\iint_{U_{3}} &  \frac{t^{-m}\left(  t-1\right)  ^{m-1}}{\left\vert
j^{-1}tx-y\right\vert }dt\text{ }\left\vert y\right\vert ^{m}q\left(
y\right)  dy\nonumber\\
&  =\iint_{U_{3}}t^{-m}\left(  t-1\right)  ^{m-1}\left(  \frac{1}{\left\vert
y\right\vert }\sum\limits_{k=0}^{\infty}\left(  \frac{j^{-1}t\left\vert
x\right\vert }{\left\vert y\right\vert }\right)  ^{k}P_{k}\left(
\widehat{x}\widehat{y}\right)  \right)  dt\text{ }\left\vert y\right\vert
^{m}q\left(  y\right)  dy\nonumber\\
&  \leq\iint_{U_{3}}t^{-m}\left(  t-1\right)  ^{m-1}\left(  \frac
{1}{\left\vert y\right\vert }\sum\limits_{k=0}^{\infty}\left(  \frac
{j^{-1}t\left\vert x\right\vert }{\left\vert y\right\vert }\right)
^{k}\left\vert P_{k}\left(  \widehat{x}\widehat{y}\right)  \right\vert
\right)  dt\text{ }\left\vert y\right\vert ^{m}q\left(  y\right)
dy\nonumber\\
&  \leq\iint_{U_{3}}t^{-m}\left(  t-1\right)  ^{m-1}\left(  \frac
{1}{\left\vert y\right\vert }\sum\limits_{k=0}^{\infty}\left(  \frac
{j^{-1}t\left\vert x\right\vert }{\left\vert y\right\vert }\right)
^{k}\right)  dt\text{ }\left\vert y\right\vert ^{m}q\left(  y\right)
dy\nonumber\\
&  \leq\iint_{U_{3}}t^{-m}\left(  t-1\right)  ^{m-1}\left(  \frac
{1}{\left\vert y\right\vert }\sum\limits_{k=0}^{\infty}\left(  \frac
{1}{1+\varepsilon}\right)  ^{k}\right)  dt\text{ }\left\vert y\right\vert
^{m}q\left(  y\right)  dy\nonumber\\
&  =\frac{1}{1-\frac{1}{1+\varepsilon}}\iint_{U_{3}}t^{-m}\left(  t-1\right)
^{m-1}\frac{1}{\left\vert y\right\vert }dt\text{ }\left\vert y\right\vert
^{m}q\left(  y\right)  dy\nonumber\\
&  =\left(  1+\frac{1}{\varepsilon}\right)  \iint_{U_{3}}t^{-m}\left(
t-1\right)  ^{m-1}\frac{1}{\left\vert y\right\vert }dt\text{ }\left\vert
y\right\vert ^{m}q\left(  y\right)  dy\nonumber\\
&  \leq\frac{1+\varepsilon}{\varepsilon}\iint_{U_{3}}t^{-m}\left(  t-1\right)
^{m-1}\frac{1}{\left(  1+\varepsilon\right)  j^{-1}t\left\vert x\right\vert
}dt\text{ }\left\vert y\right\vert ^{m}q\left(  y\right)  dy\nonumber\\
&  =\frac{j}{\varepsilon\left\vert x\right\vert }\iint_{U_{3}}t^{-\left(
m+1\right)  }\left(  t-1\right)  ^{m-1}dt\text{ }\left\vert y\right\vert
^{m}q\left(  y\right)  dy\nonumber\\
&  \leq\frac{j}{\varepsilon\left\vert x\right\vert }\left(  \int_{1}^{\infty
}t^{-\left(  m+1\right)  }\left(  t-1\right)  ^{m-1}dt\right)  \int\left\vert
y\right\vert ^{2n-d+1}q\left(  y\right)  dy\label{a033}\\
&  <\infty.\nonumber
\end{align}

We know the basis function is continuous so we need only obtain an inequality
valid for say $\left\vert x\right\vert \geq1$.\medskip

\fbox{??! \textbf{Integration over} $U_{2}$} \textbf{INVESTIGATION IS
INCOMPLETE}! \textbf{We have assumed that} $d\geq3$. In this case%
\begin{align*}
U_{2}  & =\left\{  \left(  y,t\right)  :\underline{\mu}t<\left\vert
y\right\vert <\overline{\mu}t,\text{ }t\geq1\right\} \\
& =\left\{  \left(  y,t\right)  :\left(  1-\varepsilon\right)  j^{-1}%
\left\vert x\right\vert t<\left\vert y\right\vert <\left(  1+\varepsilon
\right)  j^{-1}\left\vert x\right\vert t,\text{ }t\geq1\right\}  ,
\end{align*}

and so%
\begin{align*}
\iint_{U_{2}}\frac{t^{-m}\left(  t-1\right)  ^{m-1}}{\left\vert j^{-1}%
tx-y\right\vert }dt\text{ }\left\vert y\right\vert ^{m}q\left(  y\right)  dy
& =\iint_{U_{2}}\frac{t^{-m}\left(  t-1\right)  ^{m-1}}{\left\vert
j^{-1}tx-y\right\vert }dt\text{ }\left\vert y\right\vert ^{m}q\left(
y\right)  dy\\
& =\frac{j}{\left\vert x\right\vert }\iint_{U_{2}}\frac{t^{-\left(
m+1\right)  }\left(  t-1\right)  ^{m-1}}{\left\vert \widehat{x}-\frac
{jy}{t\left\vert x\right\vert }\right\vert }dt\text{ }\left\vert y\right\vert
^{m}q\left(  y\right)  dy
\end{align*}

Try the change of variables $z=jy/\left\vert x\right\vert $, $dz=\left(
\frac{j}{\left\vert x\right\vert }\right)  ^{d}dy$ so that%
\begin{align*}
\chi\left(  x\right)   &  :=\iint_{U_{2}}\frac{t^{-m}\left(  t-1\right)
^{m-1}}{\left\vert j^{-1}tx-y\right\vert }dt\text{ }\left\vert y\right\vert
^{m}q\left(  y\right)  dy\\
&  =\frac{j}{\left\vert x\right\vert }\iint_{U_{2}^{\prime}}\frac{t^{-\left(
m+1\right)  }\left(  t-1\right)  ^{m-1}}{\left\vert \widehat{x}-\frac{z}%
{t}\right\vert }dt\text{ }\left(  \frac{\left\vert x\right\vert }{j}\left\vert
z\right\vert \right)  ^{m}q\left(  \frac{\left\vert x\right\vert }{j}z\right)
\left(  \frac{\left\vert x\right\vert }{j}\right)  ^{d}dz,
\end{align*}

where%
\[
U_{2}^{\prime}:=\left\{  \left(  z,t\right)  :\left(  1-\varepsilon\right)
t<\left\vert z\right\vert <\left(  1+\varepsilon\right)  t,\text{ }%
t\geq1\right\}  .
\]

Let $\mathcal{O}$ be a rotation such that $\mathcal{O}x=\left\vert
x\right\vert e^{\left(  d\right)  }=\left\vert x\right\vert \left(
0,0,\ldots,0,1\right)  $. Then%
\begin{align}
\chi\left(  \mathcal{O}x\right)   & =\frac{j}{\left\vert x\right\vert }%
\iint_{U_{2}^{\prime}}\frac{t^{-\left(  m+1\right)  }\left(  t-1\right)
^{m-1}}{\left\vert e^{\left(  d\right)  }-\frac{z}{t}\right\vert }dt\text{
}\left(  \frac{\left\vert x\right\vert }{j}\left\vert z\right\vert \right)
^{m}q\left(  \frac{\left\vert x\right\vert }{j}z\right)  \left(
\frac{\left\vert x\right\vert }{j}\right)  ^{d}dz\nonumber\\
& :s:=\frac{\left\vert x\right\vert }{j}\Rightarrow\label{Ap013}\\
& =\frac{1}{s}\int_{1}^{\infty}\int\limits_{\left\vert z\right\vert =\left(
1-\varepsilon\right)  t}^{\left\vert z\right\vert =\left(  1+\varepsilon
\right)  t}\frac{\left(  s\left\vert z\right\vert \right)  ^{m}q\left(
sz\right)  s^{d}}{\left\vert e^{\left(  d\right)  }-\frac{z}{t}\right\vert
}dz\text{ }t^{-\left(  m+1\right)  }\left(  t-1\right)  ^{m-1}dt.\label{Ap010}%
\end{align}

Note that $m=2n-d+1$.

We could choose the constraint%
\begin{align*}
\int\limits_{\left\vert z\right\vert =\left(  1-\varepsilon\right)
t}^{\left\vert z\right\vert =\left(  1+\varepsilon\right)  t}\frac{\left(
s\left\vert z\right\vert \right)  ^{2n-d+1}q\left(  sz\right)  s^{d}%
}{\left\vert e^{\left(  1\right)  }-\frac{z}{t}\right\vert }dz  & \leq
const,\quad s,t\geq1.\\
s^{2n+1}\int\limits_{\left\vert z\right\vert =\left(  1-\varepsilon\right)
t}^{\left\vert z\right\vert =\left(  1+\varepsilon\right)  t}\frac{\left\vert
z\right\vert ^{2n-d+1}q\left(  sz\right)  }{\left(  \left\vert z^{\prime
}\right\vert ^{2}+\left(  1-z_{d}\right)  ^{2}\right)  ^{1/2}}dz  & \leq
const,\quad s,t\geq1.
\end{align*}

Now try the substitution $\zeta=z/t$. This gives%
\begin{align}
s^{2n+1}\int\limits_{\left\vert \zeta\right\vert =1-\varepsilon}^{\left\vert
\zeta\right\vert =1+\varepsilon}\frac{\left\vert t\zeta\right\vert
^{2n-d+1}q\left(  st\zeta\right)  }{\left(  t^{2}\left\vert \zeta^{\prime
}\right\vert ^{2}+\left(  1-t\zeta_{d}\right)  ^{2}\right)  ^{1/2}}t^{d}d\zeta
& \leq C_{\varepsilon},\quad s,t\geq1.\nonumber\\
\left(  ts\right)  ^{2n+1}\int\limits_{\left\vert \zeta\right\vert
=1-\varepsilon}^{\left\vert \zeta\right\vert =1+\varepsilon}\frac{\left\vert
\zeta\right\vert ^{2n-d+1}q\left(  st\zeta\right)  }{\left(  t^{2}\left\vert
\zeta^{\prime}\right\vert ^{2}+\left(  1-t\zeta_{d}\right)  ^{2}\right)
^{1/2}}d\zeta & \leq C_{\varepsilon},\quad s,t\geq1.\nonumber\\
\int\limits_{\left\vert \zeta\right\vert =1-\varepsilon}^{\left\vert
\zeta\right\vert =1+\varepsilon}\frac{\left\vert \zeta\right\vert
^{2n-d+1}q\left(  u\zeta\right)  }{\left(  t^{2}\left\vert \zeta^{\prime
}\right\vert ^{2}+\left(  1-t\zeta_{d}\right)  ^{2}\right)  ^{1/2}}d\zeta &
\leq\frac{C_{\varepsilon}}{u^{2n+1}},\quad u,t\geq1.\nonumber\\
\int\limits_{\left\vert \zeta\right\vert =1-\varepsilon}^{\left\vert
\zeta\right\vert =1+\varepsilon}\frac{\left\vert \zeta\right\vert
^{2n-d+1}q\left(  u\zeta\right)  }{\left(  t^{2}\left\vert \zeta\right\vert
^{2}-2t\zeta_{d}+1\right)  ^{1/2}}d\zeta & \leq\frac{C_{\varepsilon}}%
{u^{2n+1}},\quad u,t\geq1.\label{Ap011}%
\end{align}

?? Remark ***************************************

Assume $\left\vert q\left(  \zeta\right)  \right\vert \leq c\left\vert
\zeta\right\vert ^{-\left(  2n+1\right)  }$ for $1-\varepsilon\leq\left\vert
\zeta\right\vert $ and calculate%
\[
\int\limits_{\left\vert \zeta\right\vert =1-\varepsilon}^{\left\vert
\zeta\right\vert =1+\varepsilon}\frac{\left\vert \zeta\right\vert ^{2n-d+1}%
}{\left(  t^{2}\left\vert \zeta\right\vert ^{2}-2te^{\left(  d\right)  }%
\zeta+1\right)  ^{1/2}}d\zeta,
\]

using Theorem \ref{Thm_Integ_u(xy,|x|)dx}:%
\begin{align*}
\int_{\left\vert \zeta\right\vert \leq r}\Phi\left(  \left\vert \zeta
\right\vert ^{2},\widehat{a}\zeta\right)  d\zeta & =\omega_{d-1}\int_{0}%
^{r}\rho^{d-1}\int_{0}^{\pi}\Phi\left(  \rho^{2},\rho\cos\theta\right)
\sin^{d-2}\theta d\theta d\rho,\\
\omega_{t}  & :=\frac{2\pi^{t/2}}{\Gamma\left(  t/2\right)  },\text{ }t>0.
\end{align*}

ETC.!

?? ****************************************

The numerator of the integrand can become zero when $1-t\zeta_{d}=0$ and
$\zeta^{\prime}=0$ i.e. $1=t\zeta_{d}$ and $1-\varepsilon\leq\zeta_{d}%
\leq1+\varepsilon$ i.e. $t\in\left[  1,\frac{1}{1-\varepsilon}\right]  $.

Also $\xi=u\zeta$ implies inequality \ref{Ap011} becomes%
\begin{align}
\frac{1}{u^{2n-d+1}}\int\limits_{\left\vert \xi\right\vert =\left(
1-\varepsilon\right)  u}^{\left\vert \xi\right\vert =\left(  1+\varepsilon
\right)  u}\frac{\left\vert \xi\right\vert ^{2n-d+1}q\left(  \xi\right)
}{\left(  \left(  \frac{t}{u}\right)  ^{2}\left\vert \xi\right\vert
^{2}+1-2\frac{t}{u}\xi_{d}\right)  ^{1/2}}\frac{d\xi}{u^{d}}  & \leq
\frac{C_{\varepsilon}}{u^{2n+1}}.\nonumber\\
& :=>\nonumber\\
\int\limits_{\left\vert \xi\right\vert =\left(  1-\varepsilon\right)
u}^{\left\vert \xi\right\vert =\left(  1+\varepsilon\right)  u}\frac
{\left\vert \xi\right\vert ^{2n-d+1}q\left(  \xi\right)  }{\left(  \left(
\frac{t}{u}\right)  ^{2}\left\vert \xi\right\vert ^{2}+1-2\frac{t}{u}\xi
_{d}\right)  ^{1/2}}d\xi & \leq\frac{C_{\varepsilon}u^{2n-d+1}u^{d}}{u^{2n+1}%
}=C_{\varepsilon}.\label{Ap014}%
\end{align}

Suppose that $\xi_{d}<0$. Then%
\begin{equation}
\int\limits_{\left\vert \xi\right\vert =\left(  1-\varepsilon\right)
u}^{\left\vert \xi\right\vert =\left(  1+\varepsilon\right)  u}\frac
{\left\vert \xi\right\vert ^{2n-d+1}q\left(  \xi\right)  }{\left(  \left(
\frac{t}{u}\right)  ^{2}\left\vert \xi\right\vert ^{2}+1-2\frac{t}{u}\xi
_{d}\right)  ^{1/2}}d\xi\leq\int\limits_{\left\vert \xi\right\vert =\left(
1-\varepsilon\right)  u}^{\left\vert \xi\right\vert =\left(  1+\varepsilon
\right)  u}\frac{\left\vert \xi\right\vert ^{2n-d+1}q\left(  \xi\right)  }%
{1}d\xi<\int\left\vert \cdot\right\vert ^{2n-d+1}q<\infty.\label{Ap012}%
\end{equation}

Thus we are only interested in $\xi_{d}>0$.

Use spherical coordinates as discussed in the positive order document:

Suppose $\xi\in\mathbb{R}^{d}$ and $d\geq2$. The \textbf{spherical polar
coordinate representation} is%
\[
\xi=\xi\left(  \rho,\phi\right)  =\xi\left(  \rho,\phi_{1},\phi_{2}%
,\ldots,\phi_{d-1}\right)  ,
\]

where%
\begin{equation}
\left.
\begin{array}
[c]{rr}%
\xi_{1}= & \rho\sin\phi_{1}\sin\phi_{2}\ldots\sin\phi_{d-1},\\
\xi_{2}= & \rho\cos\phi_{1}\sin\phi_{2}\ldots\sin\phi_{d-1},\\
\xi_{3}= & \rho\cos\phi_{2}\ldots\sin\phi_{d-1},\\
\vdots\quad & \vdots\qquad\\
\xi_{d}= & \rho\cos\phi_{d-1},
\end{array}
\right\} \label{Ap051}%
\end{equation}

and on $\mathbb{R}^{d}$,%
\[
\rho\geq0,\text{ }-\pi\leq\phi_{1}\leq\pi,\text{ }0\leq\phi_{2},\ldots
,\phi_{d-1}\leq\pi.
\]

The volume element is%
\[
d\xi=d\xi_{1}d\xi_{2}\ldots d\xi_{d}=\rho^{d-1}%
{\textstyle\prod\limits_{j=2}^{d-1}}
\sin^{j-1}\phi_{j}d\rho d\phi,
\]

where $d\phi=d\phi_{1}\cdots d\phi_{d-1}$. Thus%
\begin{align*}
\int\limits_{\left\vert \xi\right\vert =\left(  1-\varepsilon\right)
u}^{\left\vert \xi\right\vert =\left(  1+\varepsilon\right)  u} &
\frac{\left\vert \xi\right\vert ^{2n-d+1}q\left(  \xi\right)  d\xi}{\left(
\left(  \frac{t}{u}\right)  ^{2}\left\vert \xi\right\vert ^{2}+1-2\frac{t}%
{u}\xi_{d}\right)  ^{1/2}}\\
&  =\int\limits_{\rho=\left(  1-\varepsilon\right)  u}^{\rho=\left(
1+\varepsilon\right)  u}\int_{0}^{\pi/2}\ldots\int_{0}^{\pi}\int_{-\pi}^{\pi
}\frac{\rho^{2n-d+1}q\left(  \xi\right)  \rho^{d-1}d\rho}{\left(  \left(
\frac{t}{u}\right)  ^{2}\rho^{2}-2\frac{t}{u}\rho\cos\phi_{d-1}+1\right)
^{1/2}}%
{\textstyle\prod\limits_{j=1}^{d-1}}
\sin^{j-1}\phi_{j}d\phi\\
&  =\int\limits_{\rho=\left(  1-\varepsilon\right)  u}^{\rho=\left(
1+\varepsilon\right)  u}\int_{0}^{\pi}\ldots\int_{0}^{\pi}\int_{-\pi}^{\pi
}\frac{\rho^{2n}q\left(  \xi\right)  }{\left(  \left(  \frac{t}{u}\rho\right)
^{2}-2\frac{t}{u}\rho\cos\phi_{d-1}+1\right)  ^{1/2}}d\rho%
{\textstyle\prod\limits_{j=1}^{d-1}}
\sin^{j-1}\phi_{j}d\phi\\
&  =\int_{0}^{\pi}\ldots\int_{0}^{\pi}\int_{-\pi}^{\pi}\left(  \int%
\limits_{\left(  1-\varepsilon\right)  u}^{\left(  1+\varepsilon\right)
u}\int_{0}^{\pi/2}\frac{\rho^{2n}q\left(  \xi\right)  d\theta d\rho}{\left(
\left(  \frac{t}{u}\rho\right)  ^{2}-2\frac{t}{u}\rho\cos\theta+1\right)
^{1/2}}\right)
{\textstyle\prod\limits_{j=1}^{d-2}}
\sin^{j-1}\phi_{j}d\phi^{\prime},
\end{align*}

where $d\phi^{\prime}:=d\phi_{1}\cdots d\phi_{d-2}$ and $\theta:=\phi_{d-1}$.
Thus we could apply the constraint%
\begin{equation}
I:=\int\limits_{\left(  1-\varepsilon\right)  u}^{\left(  1+\varepsilon
\right)  u}\int_{0}^{\pi/2}\frac{\rho^{2n}q\left(  \xi\right)  d\theta d\rho
}{\left(  \left(  \frac{t}{u}\rho\right)  ^{2}-2\frac{t}{u}\rho\cos
\theta+1\right)  ^{1/2}}\leq C_{\varepsilon}^{\prime},\quad u,t\geq1,\text{
}\forall\phi^{\prime}.\label{Ap031}%
\end{equation}

Noting \ref{Ap051}, we can write%
\[
I\left(  \eta^{\prime};t,u\right)  =\int\limits_{\left(  1-\varepsilon\right)
u}^{\left(  1+\varepsilon\right)  u}\int_{0}^{\pi/2}\frac{\rho^{2n}q\left(
\eta_{1}\rho\sin\theta,\eta_{2}\rho\sin\theta,\ldots,\eta_{d-2}\rho\sin
\theta,\rho\cos\theta\right)  d\theta d\rho}{\left(  \left(  \frac{t}{u}%
\rho\right)  ^{2}-2\frac{t}{u}\rho\cos\theta+1\right)  ^{1/2}},
\]

where%
\begin{equation}
\left.
\begin{array}
[c]{c}%
\begin{array}
[c]{rr}%
\eta_{1}= & \sin\phi_{1}\sin\phi_{2}\ldots\sin\phi_{d-2},\\
\eta_{2}= & \cos\phi_{1}\sin\phi_{2}\ldots\sin\phi_{d-2},\\
\eta_{3}= & \cos\phi_{2}\ldots\sin\phi_{d-2},\\
\vdots\quad & \vdots\qquad\\
\eta_{d-1}= & \cos\phi_{d-2},
\end{array}
\\
\eta^{\prime}=\left(  \eta_{1},\ldots,\eta_{d-1}\right)  ,\\
\left\vert \eta^{\prime}\right\vert =1.
\end{array}
\right\}  .\label{Ap036}%
\end{equation}

Thus%
\begin{equation}
I\left(  \eta^{\prime};t,u\right)  =\int\limits_{\left(  1-\varepsilon\right)
u}^{\left(  1+\varepsilon\right)  u}\int_{0}^{\pi/2}\frac{\rho^{2n}q\left(
\left(  \rho\sin\theta\right)  \eta^{\prime},\rho\cos\theta\right)  d\theta
d\rho}{\left(  \left(  \frac{t}{u}\rho\right)  ^{2}-2\frac{t}{u}\rho\cos
\theta+1\right)  ^{1/2}}\leq C_{\varepsilon}^{\prime},\quad t,u\geq1,\text{
}\left\vert \eta^{\prime}\right\vert =1.\label{Ap037}%
\end{equation}

We have%
\begin{align*}
\left(  \frac{t}{u}\rho\right)  ^{2}-2\frac{t}{u}\rho\cos\theta+1  & =\left(
\frac{t}{u}\rho\right)  ^{2}-2\frac{t}{u}\rho+1+2\frac{t}{u}\rho-2\frac{t}%
{u}\rho\cos\theta\\
& =\left(  \frac{t}{u}\rho-1\right)  ^{2}+4\frac{t}{u}\rho\sin^{2}\frac
{\theta}{2},
\end{align*}

and%
\begin{align*}
\left(  \frac{t}{u}\rho\right)  ^{2}-2\frac{t}{u}\rho\cos\theta+1  & =\left(
\frac{t}{u}\rho\right)  ^{2}-2\frac{t}{u}\rho\cos\theta+\cos^{2}\theta
+1-\cos^{2}\theta\\
& =\left(  \frac{t}{u}\rho-\cos\theta\right)  ^{2}+\sin^{2}\theta,
\end{align*}

so that
\[
\left(  \frac{t}{u}\rho\right)  ^{2}-2\frac{t}{u}\rho\cos\theta+1=0\text{
}iff\text{ }\rho=u/t\text{ }and\text{ }\theta=0.
\]

?? Determine for what values of $t$ and $u$ the denominator is positive.

?? Thus if $u/t<\left(  1-\varepsilon\right)  u$ i.e. $t>\frac{1}%
{1-\varepsilon}$ then $\rho>u/t$ and so $\left(  \frac{t}{u}\rho\right)
^{2}-2\frac{t}{u}\rho\cos\theta+1>0$. ETC.

Now employ 2-dim polar coordinates: $x_{1}=\rho\cos\theta$, $x_{2}=\rho
\sin\theta$, to get%
\begin{align}
I\left(  \eta^{\prime};t,u\right)   & =\int\limits_{\substack{\left\vert
x\right\vert =\left(  1-\varepsilon\right)  u \\x\geq0}}^{\left\vert
x\right\vert =\left(  1+\varepsilon\right)  u}\frac{\left\vert x\right\vert
^{2n-1}q\left(  \eta_{1}x_{2},\eta_{2}x_{2},\ldots,\eta_{d-2}x_{2}%
,x_{1}\right)  dx_{1}dx_{2}}{\left(  \left(  \frac{t}{u}\right)
^{2}\left\vert x\right\vert ^{2}-2\frac{t}{u}x_{1}+1\right)  ^{1/2}}\leq
C_{\varepsilon}^{\prime},\quad t,u\geq1,\text{ }\left\vert \eta^{\prime
}\right\vert =1.\nonumber\\
I\left(  \eta^{\prime};t,u\right)   & =\int\limits_{\substack{\left\vert
x\right\vert =\left(  1-\varepsilon\right)  u \\x\geq0}}^{\left\vert
x\right\vert =\left(  1+\varepsilon\right)  u}\frac{\left\vert x\right\vert
^{2n-1}q\left(  x_{2}\eta^{\prime},x_{1}\right)  dx_{1}dx_{2}}{\left(  \left(
\frac{t}{u}\right)  ^{2}\left\vert x\right\vert ^{2}-2\frac{t}{u}%
x_{1}+1\right)  ^{1/2}}\leq C_{\varepsilon}^{\prime},\quad t,u\geq1,\text{
}\left\vert \eta^{\prime}\right\vert =1.\label{Ap038}%
\end{align}

\section{Bessel \ and MacDonald's function formulas for the basis function.
\label{Sect_Bessel_centdiff_basis_formula}}

In this section we derive several formulas for the central difference basis
function which involve the Bessel functions $J_{\frac{d-2}{2}}$ and the
MacDonald's functions $K_{n-\frac{d}{2}}$.

\begin{theorem}
\label{Thm_cdiffbasis_Bessel}Suppose $w$ is a central difference weight
function generated by $n$, $l$, $q$. Then when $d\geq2$ the basis function is
given by:%
\[
G_{c}\left(  x\right)  =\tfrac{1}{\left(  2\pi\right)  ^{d/2}}\int%
\int\limits_{0}^{\infty}\left(  \sum_{j=-l}^{l}\left(  -1\right)  ^{j}%
\tbinom{2l}{j+l}\frac{J_{\frac{d-2}{2}}\left(  \left\vert x-j\tau\right\vert
r\right)  }{\left(  \left\vert x-j\tau\right\vert r\right)  ^{\frac{d-2}{2}}%
}\right)  \frac{dr}{r^{2n-d+1}}q\left(  \tau\right)  d\tau.
\]

\end{theorem}

\begin{proof}
From \ref{a953},%
\[
\frac{1}{w\left(  \xi\right)  }=\tfrac{2^{2l}}{\left(  2\pi\right)  ^{d/2}%
}\frac{1}{\left\vert \xi\right\vert ^{2n}}\int\sin^{2l}\left(  \frac{\xi\tau
}{2}\right)  q\left(  \tau\right)  d\tau,
\]

and so%
\[
G_{c}\left(  x\right)  =\tfrac{1}{\left(  2\pi\right)  ^{d/2}}\int%
\frac{e^{ix\xi}}{w\left(  \xi\right)  }d\xi=\tfrac{2^{2l}}{\left(
2\pi\right)  ^{d}}\int\frac{e^{ix\xi}}{\left\vert \xi\right\vert ^{2n}}%
\int\sin^{2l}\left(  \frac{\xi\tau}{2}\right)  q\left(  \tau\right)  d\tau
d\xi.
\]

But $1/w\in L^{1}$ so $\int\frac{d\xi}{w\left(  \xi\right)  }=\tfrac{2^{2l}%
}{\left(  2\pi\right)  ^{d}}\int\frac{1}{\left\vert \xi\right\vert ^{2n}}%
\int\sin^{2l}\left(  \frac{\xi\tau}{2}\right)  q\left(  \tau\right)  d\tau
d\xi<\infty$ and thus the integral for $G_{c}$ is absolutely convergent so by
Fubini's theorem%
\[
G_{c}\left(  x\right)  =\tfrac{2^{2l}}{\left(  2\pi\right)  ^{d/2}}\int\left(
\frac{1}{\left(  2\pi\right)  ^{d/2}}\int e^{ix\xi}\frac{\sin^{2l}\left(
\xi\tau/2\right)  }{\left\vert \xi\right\vert ^{2n}}d\xi\right)  q\left(
\tau\right)  d\tau.
\]

From \ref{a802}, $\frac{\sin^{2l}\left(  \xi\tau/2\right)  }{\left\vert
\xi\right\vert ^{2n}}\in L^{1}$, so%
\[
G_{c}\left(  x\right)  =\tfrac{2^{2l}}{\left(  2\pi\right)  ^{d/2}}\int
F_{\xi}^{-1}\left[  \frac{\sin^{2l}\left(  \xi\tau/2\right)  }{\left\vert
\xi\right\vert ^{2n}}\right]  \left(  x\right)  q\left(  \tau\right)  d\tau.
\]

From part 3 of Lemma \ref{Lem_central_diff_op multivar},%
\[
2^{2l}\sin^{2l}\left(  \xi\tau/2\right)  =\Delta_{2l,\tau}\left(  e^{i\xi\tau
}\right)  =\sum_{j=-l}^{l}\left(  -1\right)  ^{j}\tbinom{2l}{j+l}e^{-ij\xi
\tau},
\]

and hence, since $2n>d$,%
\begin{align*}
2^{2l}\int e^{ix\xi}\frac{\sin^{2l}\left(  \xi\tau/2\right)  }{\left\vert
\xi\right\vert ^{2n}}d\xi & =\lim_{\varepsilon\rightarrow0^{+}}\int%
_{\left\vert \xi\right\vert \geq\varepsilon}e^{ix\xi}\frac{2^{2l}\sin
^{2l}\left(  \xi\tau/2\right)  }{\left\vert \xi\right\vert ^{2n}}d\xi\\
& =\lim_{\varepsilon\rightarrow0^{+}}\int_{\left\vert \xi\right\vert
\geq\varepsilon}e^{ix\xi}\sum_{j=-l}^{l}\left(  -1\right)  ^{j}\tbinom
{2l}{j+l}e^{-ij\xi\tau}\frac{d\xi}{\left\vert \xi\right\vert ^{2n}}\\
& =\lim_{\varepsilon\rightarrow0^{+}}\sum_{j=-l}^{l}\left(  -1\right)
^{j}\tbinom{2l}{j+l}\int_{\left\vert \xi\right\vert \geq\varepsilon}%
\frac{e^{i\left(  x-j\tau\right)  \xi}}{\left\vert \xi\right\vert ^{2n}}d\xi.
\end{align*}

?? \textbf{Expand }$e^{-ij\xi\tau}$\textbf{\ using a Taylor series}?

From Corollary \ref{Cor_Thm_IntegBr_exp(ixy)f(|x|)dx} we have%
\[
\tfrac{1}{\left(  2\pi\right)  ^{d/2}}\int\limits_{\left\vert x\right\vert
\geq\varepsilon}e^{i\xi x}f\left(  \left\vert x\right\vert \right)
dx=\int\limits_{\varepsilon}^{\infty}\frac{J_{\frac{d-2}{2}}\left(  \left\vert
\xi\right\vert r\right)  }{\left(  \left\vert \xi\right\vert r\right)
^{\frac{d-2}{2}}}r^{d-1}f\left(  r\right)  dr,\quad f\in L^{1}\left(
B_{r}^{c}\right)  ,
\]

which implies%
\begin{align}
G_{c}\left(  x\right)   & =\tfrac{2^{2l}}{\left(  2\pi\right)  ^{d/2}}%
\int\left(  \tfrac{1}{\left(  2\pi\right)  ^{d/2}}\int e^{ix\xi}\frac
{\sin^{2l}\left(  \xi\tau/2\right)  }{\left\vert \xi\right\vert ^{2n}}%
d\xi\right)  q\left(  \tau\right)  d\tau\nonumber\\
& =\tfrac{1}{\left(  2\pi\right)  ^{d/2}}\int\left(  \tfrac{2^{2l}}{\left(
2\pi\right)  ^{d/2}}\int e^{ix\xi}\frac{\sin^{2l}\left(  \xi\tau/2\right)
}{\left\vert \xi\right\vert ^{2n}}d\xi\right)  q\left(  \tau\right)
d\tau\nonumber\\
& =\tfrac{1}{\left(  2\pi\right)  ^{d/2}}\int\left(  \tfrac{1}{\left(
2\pi\right)  ^{d/2}}\lim_{\varepsilon\rightarrow0^{+}}\sum_{j=-l}^{l}\left(
-1\right)  ^{j}\tbinom{2l}{j+l}\int_{\left\vert \xi\right\vert \geq
\varepsilon}\frac{e^{i\left(  x-j\tau\right)  \xi}}{\left\vert \xi\right\vert
^{2n}}d\xi\right)  q\left(  \tau\right)  d\tau\label{1.097}\\
& =\tfrac{1}{\left(  2\pi\right)  ^{d/2}}\int\left(  \lim_{\varepsilon
\rightarrow0^{+}}\sum_{j=-l}^{l}\left(  -1\right)  ^{j}\tbinom{2l}{j+l}%
\tfrac{1}{\left(  2\pi\right)  ^{d/2}}\int_{\left\vert \xi\right\vert
\geq\varepsilon}\frac{e^{i\left(  x-j\tau\right)  \xi}}{\left\vert
\xi\right\vert ^{2n}}d\xi\right)  q\left(  \tau\right)  d\tau\nonumber\\
& =\tfrac{1}{\left(  2\pi\right)  ^{d/2}}\int\left(  \lim_{\varepsilon
\rightarrow0^{+}}\sum_{j=-l}^{l}\left(  -1\right)  ^{j}\tbinom{2l}{j+l}%
\int\limits_{\varepsilon}^{\infty}\frac{J_{\frac{d-2}{2}}\left(  \left\vert
x-j\tau\right\vert r\right)  }{\left(  \left\vert x-j\tau\right\vert r\right)
^{\frac{d-2}{2}}}\frac{r^{d-1}}{r^{2n}}dr\right)  q\left(  \tau\right)
d\tau\nonumber\\
& =\tfrac{1}{\left(  2\pi\right)  ^{d/2}}\int\left(  \lim_{\varepsilon
\rightarrow0^{+}}\sum_{j=-l}^{l}\left(  -1\right)  ^{j}\tbinom{2l}{j+l}%
\int\limits_{\varepsilon}^{\infty}\frac{J_{\frac{d-2}{2}}\left(  \left\vert
x-j\tau\right\vert r\right)  }{\left(  \left\vert x-j\tau\right\vert r\right)
^{\frac{d-2}{2}}}\frac{dr}{r^{2n-d+1}}\right)  q\left(  \tau\right)
d\tau\label{1.096}\\
& =\tfrac{1}{\left(  2\pi\right)  ^{d/2}}\int\lim_{\varepsilon\rightarrow
0^{+}}\int\limits_{\varepsilon}^{\infty}\left(  \sum_{j=-l}^{l}\left(
-1\right)  ^{j}\tbinom{2l}{j+l}\frac{J_{\frac{d-2}{2}}\left(  \left\vert
x-j\tau\right\vert r\right)  }{\left(  \left\vert x-j\tau\right\vert r\right)
^{\frac{d-2}{2}}}\right)  \frac{dr}{r^{2n-d+1}}q\left(  \tau\right)
d\tau\nonumber\\
& =\tfrac{1}{\left(  2\pi\right)  ^{d/2}}\int\int\limits_{0}^{\infty}\left(
\sum_{j=-l}^{l}\left(  -1\right)  ^{j}\tbinom{2l}{j+l}\frac{J_{\frac{d-2}{2}%
}\left(  \left\vert x-j\tau\right\vert r\right)  }{\left(  \left\vert
x-j\tau\right\vert r\right)  ^{\frac{d-2}{2}}}\right)  \frac{dr}{r^{2n-d+1}%
}q\left(  \tau\right)  d\tau.\label{1.090}%
\end{align}

\end{proof}

\begin{remark}
In \textbf{ONE dimension}, using the Taylor series expansion \ref{a2.35} ??,%

\[
e^{ix\xi}=\sum_{k\leq2l-1}\frac{\left(  ix\xi\right)  ^{k}}{k!}+\frac
{\sqrt{2\pi}}{\left(  2l-1\right)  !}\left(  ix\xi\right)  ^{2l}%
\overline{\widehat{g_{2l-1}}}\left(  x\xi\right)  .
\]

\[
\lim_{\varepsilon\rightarrow0^{+}}\int_{\left\vert \xi\right\vert
\geq\varepsilon}e^{ix\xi}\sum_{j=-l}^{l}\left(  -1\right)  ^{j}\tbinom
{2l}{j+l}e^{-ij\xi\tau}\frac{d\xi}{\xi^{2n}}=\lim_{\varepsilon\rightarrow
0^{+}}\int_{\left\vert \xi\right\vert \geq\varepsilon}e^{ix\xi}\sum_{j=-l}%
^{l}\left(  -1\right)  ^{j}\tbinom{2l}{j+l}\left(  \frac{\sqrt{2\pi}}{\left(
2l-1\right)  !}\left(  ij\tau\xi\right)  ^{2l}\overline{\widehat{g_{2l-1}}%
}\left(  j\tau\xi\right)  \right)  \frac{d\xi}{\xi^{2n}}.
\]

Use \ref{a2.07}.
\end{remark}

If we assume that $q$ is radial then the formula of Theorem
\ref{Thm_cdiffbasis_Bessel} can be expressed as:

\begin{corollary}
\label{Cor_cdiffbasis_RadialFunc_Bessel}Radial version Theorem
\ref{Thm_cdiffbasis_Bessel}:%
\begin{align*}
&  \left(  G_{c}\right)  _{\odot}\left(  \sigma\right) \\
&  =\tfrac{\omega_{d-1}}{\left(  2\pi\right)  ^{d/2}}\int\limits_{0}^{\infty
}\int\limits_{0}^{\pi}\int\limits_{0}^{\infty}\left(  \sum_{j=-l}^{l}\left(
-1\right)  ^{j}\tbinom{2l}{j+l}\frac{J_{\frac{d-2}{2}}\left(  \left(
\sigma^{2}-2j\sigma\rho\cos\theta+j^{2}\rho^{2}\right)  r\right)  }{\left(
\left(  \sigma^{2}-2j\sigma\rho\cos\theta+j^{2}\rho^{2}\right)  r\right)
^{\frac{d-2}{2}}}\right)  \frac{\rho^{d-1}q_{\odot}\left(  \rho\right)
}{r^{2n-d+1}}\sin^{d-2}\theta drd\theta d\rho,
\end{align*}

\end{corollary}

\begin{proof}
From Theorem \ref{Thm_cdiffbasis_Bessel},
\[
G_{c}\left(  x\right)  =\tfrac{1}{\left(  2\pi\right)  ^{d/2}}\int%
\int\limits_{0}^{\infty}\left(  \sum_{j=-l}^{l}\left(  -1\right)  ^{j}%
\tbinom{2l}{j+l}\frac{J_{\frac{d-2}{2}}\left(  \left\vert x-j\tau\right\vert
r\right)  }{\left(  \left\vert x-j\tau\right\vert r\right)  ^{\frac{d-2}{2}}%
}\right)  \frac{dr}{r^{2n-d+1}}q\left(  \tau\right)  d\tau.
\]

Write $\left\vert x-j\tau\right\vert ^{2}=\left\vert x\right\vert ^{2}%
-2jx\tau+j^{2}\left\vert \tau\right\vert ^{2}$. From Theorem
\ref{Thm_Integ_u(xy,|x|)dx},%
\[
\int\Phi\left(  \left\vert \tau\right\vert ^{2},\xi\tau\right)  d\tau
=\omega_{d-1}\int_{0}^{\infty}\rho^{d-1}\int_{0}^{\pi}\Phi\left(  \rho
^{2},\left\vert \xi\right\vert \rho\cos\theta\right)  \sin^{d-2}\theta d\theta
d\rho,
\]

and here%
\[
\Phi\left(  a,b\right)  =\int\limits_{0}^{\infty}\left(  \sum_{j=-l}%
^{l}\left(  -1\right)  ^{j}\tbinom{2l}{j+l}\frac{J_{\frac{d-2}{2}}\left(
\left(  \left\vert x\right\vert ^{2}-2jb+j^{2}a\right)  ^{\frac{1}{2}%
}r\right)  }{\left(  \left(  \left\vert x\right\vert ^{2}-2jb+j^{2}a\right)
^{\frac{1}{2}}r\right)  ^{\frac{d-2}{2}}}\right)  \frac{dr}{r^{2n-d+1}%
}q_{\odot}\left(  \sqrt{a}\right)  ,
\]

so%
\begin{align*}
&  G_{c}\left(  x\right) \\
&  =\tfrac{\omega_{d-1}}{\left(  2\pi\right)  ^{d/2}}\int\limits_{0}^{\infty
}\int\limits_{0}^{\pi}\int\limits_{0}^{\infty}\left(  \sum_{j=-l}^{l}\left(
-1\right)  ^{j}\tbinom{2l}{j+l}\frac{J_{\frac{d-2}{2}}\left(  \left(
\left\vert x\right\vert ^{2}-2j\left\vert x\right\vert \rho\cos\theta
+j^{2}\rho^{2}\right)  ^{\frac{1}{2}}r\right)  }{\left(  \left(  \left\vert
x\right\vert ^{2}-2j\left\vert x\right\vert \rho\cos\theta+j^{2}\rho
^{2}\right)  ^{\frac{1}{2}}r\right)  ^{\frac{d-2}{2}}}\right)  \frac
{dr}{r^{2n-d+1}}\rho^{d-1}q_{\odot}\left(  \rho\right)  \sin^{d-2}\theta
d\theta d\rho,
\end{align*}

and%
\begin{align*}
&  \left(  G_{c}\right)  _{\odot}\left(  \sigma\right) \\
&  =\tfrac{\omega_{d-1}}{\left(  2\pi\right)  ^{d/2}}\int\limits_{0}^{\infty
}\int\limits_{0}^{\pi}\int\limits_{0}^{\infty}\left(  \sum_{j=-l}^{l}\left(
-1\right)  ^{j}\tbinom{2l}{j+l}\frac{J_{\frac{d-2}{2}}\left(  \left(
\sigma^{2}-2j\sigma\rho\cos\theta+j^{2}\rho^{2}\right)  ^{\frac{1}{2}%
}r\right)  }{\left(  \left(  \sigma^{2}-2j\sigma\rho\cos\theta+j^{2}\rho
^{2}\right)  ^{\frac{1}{2}}r\right)  ^{\frac{d-2}{2}}}\right)  \frac
{\rho^{d-1}q_{\odot}\left(  \rho\right)  }{r^{2n-d+1}}\sin^{d-2}\theta
drd\theta d\rho.
\end{align*}

\end{proof}

The next theorem is a result where the basis function is expressed in terms of
the MacDonald's function $K_{\lambda}$ (also called the modified Bessel
function of the third kind). This formula is a bit different to the previous results.

?? \textbf{This result could be redone using the inverse Fourier transform of}
$\frac{1}{\left(  1+\left\vert \xi\right\vert ^{2}\right)  ^{n}}$. \textbf{See
Subsection} \ref{SbSect_wt_func_examples}.

\begin{theorem}
\label{Thm_CentDiffBasisFormula_Bessel}??%
\[
G_{c}\left(  x\right)  =\tfrac{1}{\left(  2\pi\right)  ^{\frac{d}{2}}%
2^{n-1}\Gamma\left(  n\right)  }\int\lim_{\varepsilon\rightarrow0^{+}}%
\sum_{j=-l}^{l}\left(  -1\right)  ^{j}\tbinom{2l}{j+l}\left\vert
x-j\tau\right\vert ^{n-1}\frac{K_{n-\frac{d}{2}}\left(  \varepsilon\left\vert
x-j\tau\right\vert \right)  }{\varepsilon^{n-\frac{d}{2}}}q\left(
\tau\right)  d\tau.
\]

\end{theorem}

\begin{proof}
?? Repeat the proof of Theorem \ref{Thm_cdiffbasis_Bessel} until equation
\ref{1.097} where $\varepsilon$ was introduced. Instead we now write%
\[
G_{c}\left(  x\right)  =\tfrac{1}{\left(  2\pi\right)  ^{d/2}}\int\left(
\tfrac{1}{\left(  2\pi\right)  ^{d/2}}\lim_{\varepsilon\rightarrow0^{+}}%
\sum_{j=-l}^{l}\left(  -1\right)  ^{j}\tbinom{2l}{j+l}\int\frac{e^{i\left(
x-j\tau\right)  \xi}}{\left(  \varepsilon^{2}+\left\vert \xi\right\vert
^{2}\right)  ^{n}}d\xi\right)  q\left(  \tau\right)  d\tau,
\]

and using the same spherical function formulas obtain%
\begin{align*}
G_{c}\left(  x\right)   & =\tfrac{1}{\left(  2\pi\right)  ^{d/2}}\int%
\lim_{\varepsilon\rightarrow0^{+}}\int\limits_{0}^{\infty}\left(  \sum
_{j=-l}^{l}\left(  -1\right)  ^{j}\tbinom{2l}{j+l}\frac{J_{\frac{d-2}{2}%
}\left(  \left\vert x-j\tau\right\vert r\right)  }{\left(  \left\vert
x-j\tau\right\vert r\right)  ^{\frac{d-2}{2}}}\right)  \frac{r^{d-1}%
dr}{\left(  \varepsilon^{2}+r^{2}\right)  ^{n}}q\left(  \tau\right)  d\tau\\
& =\tfrac{1}{\left(  2\pi\right)  ^{d/2}}\int\lim_{\varepsilon\rightarrow
0^{+}}\sum_{j=-l}^{l}\left(  -1\right)  ^{j}\tbinom{2l}{j+l}\int%
\limits_{0}^{\infty}\frac{J_{\frac{d-2}{2}}\left(  \left\vert x-j\tau
\right\vert r\right)  }{\left(  \left\vert x-j\tau\right\vert r\right)
^{\frac{d-2}{2}}}\frac{r^{d-1}dr}{\left(  \varepsilon^{2}+r^{2}\right)  ^{n}%
}q\left(  \tau\right)  d\tau.
\end{align*}

the inner integral can be eliminated using identity 6.565.4 of Gradshteyn and
Ryzhik \cite{GradRyz07}:%
\[
\int\limits_{0}^{\infty}\frac{J_{\nu}\left(  br\right)  r^{\nu+1}}{\left(
a^{2}+r^{2}\right)  ^{\mu+1}}dr=\frac{a^{\nu-\mu}b^{\mu}}{2^{\mu}\Gamma\left(
\mu+1\right)  }K_{\nu-\mu}\left(  ab\right)  ,\quad-1<\nu<2\mu+\frac{3}%
{2},\text{ }a,b>0,
\]

where $K_{\nu-\mu}=K_{\mu-\nu}$ is the modified Bessel function of the third
kind - see Theorem \ref{Thm_bnds_modif_MacDonald}.

In our case $b=\left\vert x-j\tau\right\vert $, $a=\varepsilon$, $\mu=n-1$,
$\nu=\frac{d-2}{2}$, $\nu-\mu=\frac{d-2}{2}-n+1=\frac{d}{2}-n<0$ so%
\begin{align*}
G_{c}\left(  x\right)   & =\tfrac{1}{\left(  2\pi\right)  ^{\frac{d}{2}}}%
\int\lim_{\varepsilon\rightarrow0^{+}}\sum_{j=-l}^{l}\left(  -1\right)
^{j}\tbinom{2l}{j+l}\frac{\varepsilon^{\frac{d}{2}-n}\left\vert x-j\tau
\right\vert ^{n-1}}{2^{n-1}\Gamma\left(  n\right)  }K_{\frac{d}{2}-n}\left(
\varepsilon\left\vert x-j\tau\right\vert \right)  q\left(  \tau\right)
d\tau\\
& =\tfrac{1}{\left(  2\pi\right)  ^{\frac{d}{2}}2^{n-1}\Gamma\left(  n\right)
}\int\lim_{\varepsilon\rightarrow0^{+}}\sum_{j=-l}^{l}\left(  -1\right)
^{j}\tbinom{2l}{j+l}\left\vert x-j\tau\right\vert ^{n-1}\frac{K_{n-\frac{d}%
{2}}\left(  \varepsilon\left\vert x-j\tau\right\vert \right)  }{\varepsilon
^{n-\frac{d}{2}}}q\left(  \tau\right)  d\tau.
\end{align*}

\end{proof}

\subsection{Comments\label{SbSect_Comments_multiv_cent_wt}}

\begin{enumerate}
\item See Remark \ref{Rem_fin_diff_Taylor}.

\item Can the generalized central difference operator \ref{a9.2} be used to
characterize locally the data spaces associated with multivariate central
difference weight functions introduced in Definition
\ref{Def_multi_centdiff_wt_fn} - especially the radial case?

\item Can a combination of tensor product and radial central difference
operators be used to characterize some of the positive order data function spaces?

\item How are the multivariate central difference weight functions, introduced
in Definition \ref{Def_multi_centdiff_wt_fn}, related to the results which
characterize the Sobolev spaces in a pointwise sense - see. Bojarski and
Haj\l asz \cite{BojarHaj93}, Bojarski, Haj\l asz and Strzelecki
\cite{BojHajStr2002}, Bojarski \cite{Bojarski2011}, \cite{Bojarski2012} and
Haj\l asz \cite{Hajl96}, \cite{Hajl2003}. Here a central difference operator
approximates the Taylor series remainder and the maximal function is used.
This involves concepts such as Bessel capacity and Lebesgue points.
\end{enumerate}

\section{A basis function formula derived using a tempered distribution Taylor
series expansion and the subspaces $S_{\emptyset,k}$ of $S$%
\label{Sect_CentDiffBasis_Taylor}}

These results extend the 1-dimensional results of Subsection
\ref{SbSect_tempdistrib_1_dim_centdiff_basis} to higher dimensions. We use the
tempered distribution Taylor series expansion introduced in Section
\ref{Sect_Taylor_series_data_fn} and the theory of the Schwartz subspaces
$S_{\emptyset,k}\subset S$ introduced in Definition \ref{Def_So,n} to expand a
convolution of the thin plate spline and the function $q$ and thus prove the
Taylor series basis function formula given in Theorem
\ref{Thm_centdiff_basis_Tn_q_suppbnd}. The formula is true modulo a polynomial
with degree $\leq d-2$. We will make use of the basic Fourier transform
properties given in Definition \ref{Def_Fourier} of the Appendix. The
$L_{loc}^{1}$ Taylor series expansion of Lemma \ref{Lem_Taylor_estim_L1loc_L1}
is then used to obtain the Taylor series integral remainder formula \ref{a017}
for the basis function which is again true modulo a polynomial with degree
$\leq d-2$.

Recall that this basis function $G_{c}$ is generated by the parameters $n$,
$l$ and a non-negative $L^{1}$ function $q\neq0$. From Theorem
\ref{Thm_wt_func_mdim} we must have $0<2n-d\leq2l$ and $\int\left\vert
\cdot\right\vert ^{k}q<0$ for $0\leq k\leq2n-d$.

Suppose $q\in L^{1}\left(  \mathbb{R}^{d}\right)  $. Then $\widehat{q}\in
C_{B}^{\left(  0\right)  }\left(  \mathbb{R}^{d}\right)  $ and by definition%
\[
\Delta_{2l}\widehat{q}\left(  \xi\right)  =\sum_{j=-l}^{l}\left(  -1\right)
^{j}\tbinom{2l}{j+l}\widehat{q}\left(  -j\xi\right)  .
\]

The tempered distribution formula \ref{a1.55} states that if $f\in S^{\prime}$
then%
\[
f\left(  \cdot+\xi\right)  -\sum_{\left\vert \alpha\right\vert \leq m}%
\frac{\xi^{\alpha}}{\alpha!}D^{\alpha}f=\frac{\sqrt{2\pi}}{m!}\left(  \left(
i\left(  \xi,\cdot\right)  \right)  ^{m+1}\overline{\widehat{g_{m}}}\left(
\left(  \xi,\cdot\right)  \right)  \widehat{f}\right)  ^{\vee},\quad
m=0,1,2,\ldots,
\]

and when $f$ is replaced by $\widehat{f}$,%
\[
\widehat{f}\left(  \cdot+\xi\right)  -\sum_{\left\vert \alpha\right\vert \leq
m}\frac{\xi^{\alpha}}{\alpha!}D^{\alpha}\widehat{f}=\frac{\sqrt{2\pi}}%
{m!}\left(  \left(  i\left(  \xi,\cdot\right)  \right)  ^{m+1}\overline
{\widehat{g_{m}}}\left(  \left(  \xi,\cdot\right)  \right)  f_{\_}\right)
^{\vee},
\]

where $f_{\_}$ is the distribution extension of $f_{\_}\left(  x\right)
=f\left(  -x\right)  $.

In particular, when $\left\vert \cdot\right\vert ^{k}q\in L^{1}$ for $k\leq
m$, Lemma \ref{Lem_L1_Fourier_contin} implies $\widehat{q}\in C_{B}^{\left(
m\right)  }$ and
\[
\widehat{q}\left(  \cdot+\xi\right)  -\sum_{\left\vert \alpha\right\vert \leq
m}\frac{\xi^{\alpha}}{\alpha!}D^{\alpha}\widehat{q}=\frac{\sqrt{2\pi}}%
{m!}F_{\eta}^{-1}\left[  \left(  i\xi\eta\right)  ^{m+1}\overline
{\widehat{g_{m}}}\left(  \xi\eta\right)  q\left(  -\eta\right)  \right]  ,
\]

and consequently%
\[
\widehat{q}\left(  \xi\right)  -\sum_{\left\vert \alpha\right\vert \leq
m}\frac{\xi^{\alpha}}{\alpha!}D^{\alpha}\widehat{q}\left(  0\right)
=\frac{\sqrt{2\pi}}{m!}F_{\eta}^{-1}\left[  \left(  i\xi\eta\right)
^{m+1}\overline{\widehat{g_{m}}}\left(  \xi\eta\right)  q\left(  -\eta\right)
\right]  \left(  0\right)  .
\]

From part 4 of Lemma \ref{Lem_gm_properties}, for each $\xi$,
\begin{align*}
\left\vert \left(  i\xi\eta\right)  ^{m+1}\overline{\widehat{g_{m}}}\left(
\xi\eta\right)  q\left(  -\eta\right)  \right\vert  & =\left\vert \xi
\eta\right\vert ^{m+1}\left\vert \widehat{g_{m}}\left(  \xi\eta\right)
\right\vert q\left(  -\eta\right) \\
& \leq\left\vert \xi\eta\right\vert ^{m+1}\frac{c_{m}}{\sqrt{2\pi}}\frac
{1}{1+\left\vert \xi\eta\right\vert }q\left(  -\eta\right) \\
& =\frac{c_{m}}{\sqrt{2\pi}}\left\vert \xi\eta\right\vert ^{m}q\left(
-\eta\right) \\
& \in L^{1},
\end{align*}

and so%
\[
\widehat{q}\left(  \xi\right)  -\sum_{\left\vert \alpha\right\vert \leq
m}\frac{\xi^{\alpha}}{\alpha!}D^{\alpha}\widehat{q}\left(  0\right)  =\frac
{1}{m!}\int\left(  i\xi\eta\right)  ^{m+1}\overline{\widehat{g_{m}}}\left(
\xi\eta\right)  q\left(  -\eta\right)  d\eta.
\]

The next step is to use the operator $\Delta_{2l}$. Choose $m=2n-d$ so that
$m<2l$ and%

\begin{align}
\Delta_{2l}\widehat{q}\left(  \xi\right)   & =\Delta_{2l}\left(
\widehat{q}\left(  \xi\right)  -\sum_{\left\vert \alpha\right\vert \leq
m}\frac{\xi^{\alpha}}{\alpha!}D^{\alpha}\widehat{q}\left(  0\right)  \right)
\nonumber\\
& =\frac{1}{\left(  2n-d\right)  !}\Delta_{2l,\xi}\int\left(  i\xi\eta\right)
^{2n-d+1}\overline{\widehat{g_{2n-d}}}\left(  \xi\eta\right)  q\left(
-\eta\right)  d\eta\nonumber\\
& =\frac{1}{\left(  2n-d\right)  !}\sum_{j=-l,\text{ }j\neq0}^{l}\left(
-1\right)  ^{j}\tbinom{2l}{j+l}\int\left(  -ij\xi\eta\right)  ^{2n-d+1}%
\overline{\widehat{g_{2n-d}}}\left(  -j\xi\eta\right)  q\left(  -\eta\right)
d\eta\nonumber\\
& =\frac{\left(  -1\right)  ^{n}}{\left(  2n-d\right)  !}\sum_{j=-l,\text{
}j\neq0}^{l}\left(  -1\right)  ^{j}\tbinom{2l}{j+l}\int\left(  j\xi
\eta\right)  ^{2n-d+1}\overline{\widehat{g_{2n-d}}}\left(  j\xi\eta\right)
q\left(  \eta\right)  d\eta\nonumber\\
& =\frac{\left(  -1\right)  ^{n}}{\left(  2n-d\right)  !}\sum_{j=-l,\text{
}j\neq0}^{l}\frac{\left(  -1\right)  ^{j}}{\left\vert j\right\vert ^{d}%
}\tbinom{2l}{j+l}\int\left(  \xi\tau\right)  ^{2n-d+1}\overline
{\widehat{g_{2n-d}}}\left(  \xi\tau\right)  q\left(  \tfrac{1}{j}\tau\right)
d\tau\nonumber\\
& =\frac{\left(  -1\right)  ^{n}}{\left(  2n-d\right)  !}\int\left(  \xi
\tau\right)  ^{2n-d+1}\overline{\widehat{g_{2n-d}}}\left(  \xi\tau\right)
\sum_{j=-l,\text{ }j\neq0}^{l}\frac{\left(  -1\right)  ^{j}}{\left\vert
j\right\vert ^{d}}\tbinom{2l}{j+l}q\left(  \tfrac{1}{j}\tau\right)
d\tau\nonumber\\
& =\frac{\left(  -1\right)  ^{n}}{\left(  2n-d\right)  !}\int\left(  \xi
\tau\right)  ^{2n-d+1}\overline{\widehat{g_{2n-d}}}\left(  \xi\tau\right)
q_{l}\left(  \tau\right)  d\tau,\label{a023}%
\end{align}

where the even function%
\begin{equation}
q_{l}\left(  \tau\right)  :=\sum_{j=-l,\text{ }j\neq0}^{l}\frac{\left(
-1\right)  ^{j}}{\left\vert j\right\vert ^{d}}\tbinom{2l}{j+l}q\left(
\tfrac{1}{j}\tau\right)  ,\label{a012}%
\end{equation}

is the analogue of the $q_{l}$ introduced in Theorem
\ref{Thm_cdiffbasis_part_moment_formula}.

By definition%
\[
\widehat{G_{c}}\left(  \xi\right)  =\frac{\Delta_{2l}\widehat{q}\left(
\xi\right)  }{\left\vert \xi\right\vert ^{2n}}=\frac{\left(  -1\right)  ^{n}%
}{\left(  2n-d\right)  !}\int\frac{\left(  \xi\tau\right)  ^{2n-d+1}%
}{\left\vert \xi\right\vert ^{2n}}\overline{\widehat{g_{2n-d}}}\left(  \xi
\tau\right)  q_{l}\left(  \tau\right)  d\tau,
\]

and so%
\begin{equation}
\left[  \widehat{G_{c}},\phi\right]  =\frac{\left(  -1\right)  ^{n}}{\left(
2n-d\right)  !}\int\int\frac{\left(  \xi\tau\right)  ^{2n-d+1}}{\left\vert
\xi\right\vert ^{2n}}\overline{\widehat{g_{2n-d}}}\left(  \xi\tau\right)
q_{l}\left(  \tau\right)  d\tau\text{ }\phi\left(  \xi\right)  d\xi,\quad
\phi\in S.\label{a024}%
\end{equation}

I now want to change the order of integration using Fubini's theorem.\medskip

\fbox{\textbf{If }$\phi\in S_{\emptyset,1}$:}\medskip%
\begin{align}
\int\int &  \left\vert \frac{\left(  \xi\tau\right)  ^{2n-d+1}}{\left\vert
\xi\right\vert ^{2n}}\overline{\widehat{g_{2n-d}}}\left(  \xi\tau\right)
q_{l}\left(  \tau\right)  \phi\left(  \xi\right)  \right\vert d\xi
d\tau\nonumber\\
&  =\int\int\left\vert \frac{\left(  \xi\tau\right)  ^{2n-d}}{\left\vert
\xi\right\vert ^{2n-d}\left\vert \xi\right\vert ^{d}}\left(  \xi\tau\right)
\overline{\widehat{g_{2n-d}}}\left(  \xi\tau\right)  q_{l}\left(  \tau\right)
\phi\left(  \xi\right)  \right\vert d\xi d\tau\nonumber\\
&  =\int\int\left\vert \left(  \widehat{\xi}\tau\right)  ^{2n-d}\left(
\xi\tau\right)  \widehat{g_{2n-d}}\left(  \xi\tau\right)  q_{l}\left(
\tau\right)  \frac{\phi\left(  \xi\right)  }{\left\vert \xi\right\vert ^{d}%
}\right\vert d\xi d\tau\nonumber\\
&  =\int\int\left\vert \left(  \xi\tau\right)  \widehat{g_{2n-d}}\left(
\xi\tau\right)  \frac{\phi\left(  \xi\right)  }{\left\vert \xi\right\vert
^{d}}\left(  \widehat{\xi}\tau\right)  ^{2n-d}q_{l}\left(  \tau\right)
\right\vert d\xi d\tau\nonumber\\
&  \leq\int\int\left\vert \left(  \xi\tau\right)  \widehat{g_{2n-d}}\left(
\xi\tau\right)  \frac{\phi\left(  \xi\right)  }{\left\vert \xi\right\vert
^{d}}\right\vert d\xi\text{ }\left\vert \tau\right\vert ^{2n-d}q_{l}\left(
\tau\right)  d\tau\nonumber\\
&  \leq\left\Vert s\widehat{g_{2n-d}}\left(  s\right)  \right\Vert _{\infty
}\left(  \int\frac{\left\vert \phi\left(  \xi\right)  \right\vert }{\left\vert
\xi\right\vert ^{d}}d\xi\right)  \int\left\vert \tau\right\vert ^{2n-d}%
q_{l}\left(  \tau\right)  d\tau\nonumber\\
&  =\left\Vert s\widehat{g_{2n-d}}\left(  s\right)  \right\Vert _{\infty
}\left(  \int_{\left\vert \xi\right\vert \leq1}\frac{\left\vert \phi\left(
\xi\right)  \right\vert }{\left\vert \xi\right\vert }\frac{1}{\left\vert
\xi\right\vert ^{d-1}}d\xi+\int_{\left\vert \xi\right\vert \geq1}%
\frac{\left\vert \phi\left(  \xi\right)  \right\vert }{\left\vert
\xi\right\vert ^{d}}d\xi\right)  \int\left\vert \tau\right\vert ^{2n-d}%
q_{l}\left(  \tau\right)  d\tau\nonumber\\
&  \text{: part 4 of Lemma \ref{Lem_gm_properties}, inequality \ref{a995} =%
$>$%
}\nonumber\\
&  <\infty.\label{a028}%
\end{align}
\medskip

\fbox{\textbf{If }$q$ \textbf{has bounded support and }$\phi\in S$:}\medskip%
\begin{align*}
\int\int &  \left\vert \frac{\left(  \xi\tau\right)  ^{2n-d+1}}{\left\vert
\xi\right\vert ^{2n}}\overline{\widehat{g_{2n-d}}}\left(  \xi\tau\right)
q_{l}\left(  \tau\right)  \phi\left(  \xi\right)  \right\vert d\xi d\tau\\
&  =\int\int\left\vert \frac{\left(  \xi\tau\right)  ^{2n-d+1}}{\left\vert
\xi\right\vert ^{2n-d+1}\left\vert \xi\right\vert ^{d-1}}\left(  \xi
\tau\right)  \overline{\widehat{g_{2n-d}}}\left(  \xi\tau\right)  q_{l}\left(
\tau\right)  \phi\left(  \xi\right)  \right\vert d\xi d\tau\\
&  =\int\int\left\vert \left(  \widehat{\xi}\tau\right)  ^{2n-d+1}\left(
\xi\tau\right)  \widehat{g_{2n-d}}\left(  \xi\tau\right)  q_{l}\left(
\tau\right)  \frac{\phi\left(  \xi\right)  }{\left\vert \xi\right\vert ^{d-1}%
}\right\vert d\xi d\tau\\
&  =\int\int\left\vert \left(  \xi\tau\right)  \widehat{g_{2n-d}}\left(
\xi\tau\right)  \frac{\phi\left(  \xi\right)  }{\left\vert \xi\right\vert
^{d-1}}\left(  \widehat{\xi}\tau\right)  ^{2n-d+1}q_{l}\left(  \tau\right)
\right\vert d\xi d\tau\\
&  \leq\int\int\left\vert \left(  \xi\tau\right)  \widehat{g_{2n-d}}\left(
\xi\tau\right)  \frac{\phi\left(  \xi\right)  }{\left\vert \xi\right\vert
^{d-1}}\right\vert d\xi\text{ }\left\vert \tau\right\vert ^{2n-d+1}%
q_{l}\left(  \tau\right)  d\tau\\
&  \leq\left\Vert s\widehat{g_{2n-d}}\left(  s\right)  \right\Vert _{\infty
}\left(  \int\frac{\left\vert \phi\left(  \xi\right)  \right\vert }{\left\vert
\xi\right\vert ^{d-1}}d\xi\right)  \int\left\vert \tau\right\vert
^{2n-d+1}q_{l}\left(  \tau\right)  d\tau\\
&  =\left\Vert s\widehat{g_{2n-d}}\left(  s\right)  \right\Vert _{\infty
}\left(  \int_{\left\vert \xi\right\vert \leq1}\frac{\left\vert \phi\left(
\xi\right)  \right\vert }{\left\vert \xi\right\vert ^{d-1}}d\xi+\int%
_{\left\vert \xi\right\vert \geq1}\frac{\left\vert \phi\left(  \xi\right)
\right\vert }{\left\vert \xi\right\vert ^{d-1}}d\xi\right)  \int\left\vert
\tau\right\vert ^{2n-d+1}q_{l}\left(  \tau\right)  d\tau\\
&  \leq\left\Vert s\widehat{g_{2n-d}}\left(  s\right)  \right\Vert _{\infty
}\left(  \left\Vert \phi\right\Vert _{\infty}\int_{\left\vert \xi\right\vert
\leq1}\frac{d\xi}{\left\vert \xi\right\vert ^{d-1}}+\int_{\left\vert
\xi\right\vert \geq1}\frac{\left\vert \phi\left(  \xi\right)  \right\vert
}{\left\vert \xi\right\vert ^{d-1}}d\xi\right)  \int\left\vert \tau\right\vert
^{2n-d+1}q_{l}\left(  \tau\right)  d\tau\\
&  \text{: }q\text{ has bounded support, part 4 of Lemma
\ref{Lem_gm_properties}, inequality \ref{a995} =%
$>$%
}\\
&  <\infty.
\end{align*}

Thus%
\begin{align}
\left[  \widehat{G_{c}},\phi\right]   & =\frac{\left(  -1\right)  ^{n}%
}{\left(  2n-d\right)  !}\int\int\frac{\left(  \xi\tau\right)  ^{2n-d+1}%
}{\left\vert \xi\right\vert ^{2n}}\overline{\widehat{g_{2n-d}}}\left(  \xi
\tau\right)  q_{l}\left(  \tau\right)  d\tau\text{ }\phi\left(  \xi\right)
d\xi\nonumber\\
& =\frac{\left(  -1\right)  ^{n}}{\left(  2n-d\right)  !}\int\left(  \int%
\frac{1}{\left\vert \xi\right\vert ^{2n}}\left(  \xi\tau\right)
^{2n-d+1}\overline{\widehat{g_{2n-d}}}\left(  \xi\tau\right)  \phi\left(
\xi\right)  d\xi\right)  q_{l}\left(  \tau\right)  d\tau.\label{a00}%
\end{align}
\medskip

Observe that for all $\tau\centerdot\neq0$, $\phi\in S_{\emptyset,d-1}$
implies $\left(  \xi\tau\right)  ^{2n-d+1}\overline{\widehat{g_{2n-d}}}\left(
\xi\tau\right)  \phi\left(  \xi\right)  \in S_{\emptyset,2n}$ and so from part
1 of Lemma \ref{Lem_thin_plate_splin},%
\begin{align}
\int\frac{1}{\left\vert \xi\right\vert ^{2n}} &  \left(  \xi\tau\right)
^{2n-d+1}\overline{\widehat{g_{2n-d}}}\left(  \xi\tau\right)  \phi\left(
\xi\right)  d\xi\nonumber\\
&  =\frac{1}{c_{n,d}}\left[  \widehat{T_{n}},\left(  \xi\tau\right)
^{2n-d+1}\overline{\widehat{g_{2n-d}}}\left(  \xi\tau\right)  \phi\left(
\xi\right)  \right] \nonumber\\
&  =\frac{1}{c_{n,d}}\left[  T_{n},F_{\xi}\left[  \left(  \xi\tau\right)
^{2n-d+1}\overline{\widehat{g_{2n-d}}}\left(  \xi\tau\right)  \phi\left(
\xi\right)  \right]  \right] \label{a034}\\
&  =\frac{\left(  -1\right)  ^{n-\frac{1}{2}\left(  d-1\right)  }}{c_{n,d}%
}\int T_{n}\left(  \zeta\right)  F_{\xi}\left[  \overline{\widehat{D^{2n-d+1}%
g_{2n-d}}}\left(  \xi\tau\right)  \phi\left(  \xi\right)  \right]  \left(
\zeta\right)  d\zeta\nonumber\\
&  =\frac{\left(  -1\right)  ^{n-\frac{1}{2}\left(  d-1\right)  }}{c_{n,d}%
}\int T_{n}\left(  \zeta\right)  F_{\xi}\left[  \overline{\widehat{D^{2n-d+1}%
g_{2n-d}}}\left(  \xi\tau\right)  \phi\left(  \xi\right)  \right]  \left(
\zeta\right)  d\zeta\nonumber\\
&  =\frac{\left(  -1\right)  ^{n}}{c_{n,d}}\int T_{n}\left(  \zeta\right)
F_{\xi}\left[  \overline{\widehat{D^{2n-d+1}g_{2n-d}}}\left(  \xi\tau\right)
\phi\left(  \xi\right)  \right]  \left(  \zeta\right)  d\zeta\nonumber\\
&  =\frac{\left(  -1\right)  ^{n}}{c_{n,d}}\left[  T_{n},F_{\xi}\left[
\overline{\widehat{D^{2n-d+1}g_{2n-d}}}\left(  \xi\tau\right)  \phi\left(
\xi\right)  \right]  \right] \nonumber\\
&  =\frac{\left(  -1\right)  ^{n}}{c_{n,d}}\left[  T_{n},F_{\xi}\left[
\overline{\widehat{D^{2n-d+1}g_{2n-d}}}\left(  \xi\tau\right)  \phi\left(
\xi\right)  \right]  \right]  .\label{a000}%
\end{align}

Note that since $d\geq2$, $\phi\in S_{\emptyset,d-1}$ implies $\phi\in
S_{\emptyset,1}$ which was assumed above to prove \ref{a028}. \textbf{Note
also that it is not sufficient to assume that }$q$\textbf{\ has bounded
support}.

From part 8 of Lemma \ref{Lem_gm_properties}:%
\begin{align*}
\widehat{D^{m+1}g_{m}}\left(  t\right)   & =\left(  -1\right)  ^{m+1}\frac
{m!}{\sqrt{2\pi}}\left(  e^{-it}-\sum_{k\leq m}\frac{\left(  -it\right)  ^{k}%
}{k!}\right)  ,\quad m\geq0.\\
& \Longrightarrow\\
\widehat{D^{2n-d+1}g_{2n-d}}\left(  \xi\tau\right)   & =\frac{\left(
2n-d\right)  !}{\sqrt{2\pi}}\left(  e^{-i\xi\tau}-\sum_{k\leq2n-d}%
\frac{\left(  -i\xi\tau\right)  ^{k}}{k!}\right)  .\\
\overline{\widehat{D^{2n-d+1}g_{2n-d}}}\left(  \xi\tau\right)   &
=\frac{\left(  2n-d\right)  !}{\sqrt{2\pi}}\left(  e^{i\xi\tau}-\sum
_{k\leq2n-d}\frac{\left(  i\xi\tau\right)  ^{k}}{k!}\right)  .
\end{align*}

Thus%
\begin{align*}
F_{\xi} &  \left[  \overline{\widehat{D^{2n-d+1}g_{2n-d}}}\left(  \xi
\tau\right)  \phi\left(  \xi\right)  \right] \\
&  =\frac{\left(  2n-d\right)  !}{\sqrt{2\pi}}F_{\xi}\left[  \left(
e^{i\xi\tau}-\sum_{k\leq2n-d}\frac{\left(  i\xi\tau\right)  ^{k}}{k!}\right)
\phi\left(  \xi\right)  \right] \\
&  =\frac{\left(  2n-d\right)  !}{\sqrt{2\pi}}F_{\xi}\left[  e^{i\xi\tau}%
\phi\left(  \xi\right)  -\sum_{k\leq2n-d}\frac{1}{k!}\left(  i\xi\tau\right)
^{k}\phi\left(  \xi\right)  \right] \\
&  =\frac{\left(  2n-d\right)  !}{\sqrt{2\pi}}\left(  \widehat{\phi}\left(
\zeta-\tau\right)  -\sum_{k\leq2n-d}\frac{1}{k!}\left(  -i\tau D\right)
^{k}\widehat{\phi}\left(  \zeta\right)  \right) \\
&  =\frac{\left(  2n-d\right)  !}{\sqrt{2\pi}}\left(  \widehat{\phi}\left(
\zeta-\tau\right)  -\sum_{k\leq2n-d}\frac{1}{k!}\left(  \tau D\right)
^{k}\widehat{\phi}\left(  \zeta\right)  \right)  ,
\end{align*}

and as a consequence \ref{a000} becomes%
\begin{align}
\int &  \frac{1}{\left\vert \xi\right\vert ^{2n}}\left(  \xi\tau\right)
^{2n-d+1}\overline{\widehat{g_{2n-d}}}\left(  \xi\tau\right)  \phi\left(
\xi\right)  d\xi\nonumber\\
&  =\frac{\left(  -1\right)  ^{n}}{c_{n,d}}\left[  T_{n},\frac{\left(
2n-d\right)  !}{\sqrt{2\pi}}\left(  \widehat{\phi}\left(  \cdot-\tau\right)
-\sum_{k\leq2n-d}\frac{1}{k!}\left(  \tau D\right)  ^{k}\widehat{\phi}\right)
\right] \nonumber\\
&  =\frac{\left(  -1\right)  ^{n}}{c_{n,d}}\frac{\left(  2n-d\right)  !}%
{\sqrt{2\pi}}\left[  T_{n},\widehat{\phi}\left(  \cdot-\tau\right)
-\sum_{k\leq2n-d}\frac{1}{k!}\left(  \tau D\right)  ^{k}\widehat{\phi}\right]
\label{a035}\\
&  =\frac{\left(  -1\right)  ^{n}}{c_{n,d}}\frac{\left(  2n-d\right)  !}%
{\sqrt{2\pi}}\left(  \left[  T_{n},\widehat{\phi}\left(  \cdot-\tau\right)
\right]  -\left[  T_{n},\sum_{k\leq2n-d}\frac{1}{k!}\left(  \tau D\right)
^{k}\widehat{\phi}\right]  \right) \nonumber\\
&  =\frac{\left(  -1\right)  ^{n}}{c_{n,d}}\frac{\left(  2n-d\right)  !}%
{\sqrt{2\pi}}\left(  \left[  T_{n}\left(  \cdot+\tau\right)  ,\widehat{\phi
}\right]  -\sum_{k\leq2n-d}\frac{1}{k!}\left[  T_{n},\left(  \tau D\right)
^{k}\widehat{\phi}\right]  \right) \nonumber\\
&  =\frac{\left(  -1\right)  ^{n}b_{n,k}}{c_{n,d}}\frac{\left(  2n-d\right)
!}{\sqrt{2\pi}}\left(  \left[  T_{n}\left(  \cdot+\tau\right)  ,\widehat{\phi
}\right]  -\sum_{k\leq2n-d}\frac{1}{k!}\left[  \left(  \tau D\right)
^{k}T_{n},\widehat{\phi}\right]  \right) \nonumber\\
&  =\frac{\left(  -1\right)  ^{n}}{c_{n,d}}\frac{\left(  2n-d\right)  !}%
{\sqrt{2\pi}}\left[  T_{n}\left(  \cdot+\tau\right)  -\sum_{k\leq2n-d}\frac
{1}{k!}\left(  \tau D\right)  ^{k}T_{n},\widehat{\phi}\right] \nonumber\\
&  =\frac{\left(  -1\right)  ^{n}}{c_{n,d}}\frac{\left(  2n-d\right)  !}%
{\sqrt{2\pi}}\int\left(  T_{n}\left(  \zeta+\tau\right)  -\sum_{k\leq
2n-d}\frac{1}{k!}\left(  \tau D\right)  ^{k}T_{n}\left(  \zeta\right)
\right)  \widehat{\phi}\left(  \zeta\right)  d\zeta.\label{a018}%
\end{align}

This means \ref{a00} becomes%
\begin{align}
&  \left[  \widehat{G_{c}},\phi\right] \nonumber\\
&  =\frac{\left(  -1\right)  ^{n}}{\left(  2n-d\right)  !}\int\left(
\int\frac{1}{\left\vert \xi\right\vert ^{2n}}\left(  \xi\tau\right)
^{2n-d+1}\overline{\widehat{g_{2n-d}}}\left(  \xi\tau\right)  \phi\left(
\xi\right)  d\xi\right)  q_{l}\left(  \tau\right)  d\tau\nonumber\\
&  =\frac{\left(  -1\right)  ^{n}}{\left(  2n-d\right)  !}\int\left(
\frac{\left(  -1\right)  ^{n}}{c_{n,d}}\frac{\left(  2n-d\right)  !}%
{\sqrt{2\pi}}\int\left(
\begin{array}
[c]{l}%
T_{n}\left(  \zeta+\tau\right)  -\\
-\sum\limits_{k\leq2n-d}\frac{1}{k!}\left(  \tau D\right)  ^{k}T_{n}\left(
\zeta\right)
\end{array}
\right)  \widehat{\phi}\left(  \zeta\right)  d\zeta\right)  q_{l}\left(
\tau\right)  d\tau\nonumber\\
&  =\frac{1}{\sqrt{2\pi}}\frac{1}{c_{n,d}}\int\left(  \int\left(  T_{n}\left(
\zeta+\tau\right)  -\sum_{k\leq2n-d}\frac{1}{k!}\left(  \tau D\right)
^{k}T_{n}\left(  \zeta\right)  \right)  \widehat{\phi}\left(  \zeta\right)
d\zeta\right)  q_{l}\left(  \tau\right)  d\tau.\label{a014}%
\end{align}

Once again we want to use Fubini's theorem to change the order of integration
by showing the integrand is in $L^{1}\left(  \mathbb{R}^{2d}\right)  $. In
fact%
\begin{align}
\int\int &  \left\vert \left(  T_{n}\left(  \zeta+\tau\right)  -\sum
_{k\leq2n-d}\frac{1}{k!}\left(  \tau D\right)  ^{k}T_{n}\left(  \zeta\right)
\right)  \widehat{\phi}\left(  \zeta\right)  \right\vert d\zeta\text{ }%
q_{l}\left(  \tau\right)  d\tau\nonumber\\
&  \leq\int\int\left(  T_{n}\left(  \zeta+\tau\right)  -\sum_{k\leq2n-d}%
\frac{1}{k!}\left(  \tau D\right)  ^{k}T_{n}\left(  \zeta\right)  \right)
\left\vert \widehat{\phi}\left(  \zeta\right)  \right\vert d\zeta\text{ }%
q_{l}\left(  \tau\right)  d\tau\nonumber\\
&  =\left(  -1\right)  ^{n-\frac{d-1}{2}}\int\int\left(  \left\vert \zeta
+\tau\right\vert ^{2n-d}+\sum_{k\leq2n-d}\frac{1}{k!}\left\vert \left(  \tau
D\right)  ^{k}\left\vert \zeta\right\vert ^{2n-d}\right\vert \right)
\left\vert \widehat{\phi}\left(  \zeta\right)  \right\vert d\zeta\text{ }%
q_{l}\left(  \tau\right)  d\tau\nonumber\\
&  \leq\left(  -1\right)  ^{n-\frac{d-1}{2}}\int\int\left(  \left\vert
\zeta+\tau\right\vert ^{2n-d}+\sum_{k\leq2n-d}\frac{\left\vert \tau\right\vert
^{k}}{k!}\left\vert \left(  \widehat{\tau}D\right)  ^{k}\left\vert
\zeta\right\vert ^{2n-d}\right\vert \right)  \left\vert \widehat{\phi}\left(
\zeta\right)  \right\vert d\zeta\text{ }q_{l}\left(  \tau\right)
d\tau\nonumber\\
&  \text{: part 3 of Lemma \ref{Lem_thin_plate_splin} =%
$>$%
}\nonumber\\
&  \leq\int\int\left(  \left\vert \zeta+\tau\right\vert ^{2n-d}+\sum
\limits_{k\leq2n-d}\frac{\left\vert \tau\right\vert ^{k}}{k!}c_{k}^{\prime
}\left\vert \zeta\right\vert ^{2n-d-k}\right)  \left\vert \widehat{\phi
}\left(  \zeta\right)  \right\vert d\zeta\text{ }q_{l}\left(  \tau\right)
d\tau\nonumber\\
&  \leq\int\left(  \int\left\vert \zeta+\tau\right\vert ^{2n-d}\left\vert
\widehat{\phi}\left(  \zeta\right)  \right\vert d\zeta\right)  q_{l}\left(
\tau\right)  d\tau+\nonumber\\
&  \qquad\qquad+\sum\limits_{k\leq2n-d}\frac{c_{k}^{\prime}}{k!}\left(
\int\left\vert \zeta\right\vert ^{2n-d-k}\left\vert \widehat{\phi}\left(
\zeta\right)  \right\vert d\zeta\right)  \int\left\vert \tau\right\vert
^{k}q_{l}\left(  \tau\right)  d\tau\label{a010}\\
&  <\infty.\nonumber
\end{align}

Thus we can apply Fubini's and reverse the order of integration in \ref{a014}
to get%
\begin{align*}
&  \left[  \widehat{G_{c}},\phi\right] \\
&  =\frac{1}{\sqrt{2\pi}}\frac{1}{c_{n,d}}\int\left(  \int\left(  T_{n}\left(
\zeta+\tau\right)  -\sum_{k\leq2n-d}\frac{1}{k!}\left(  \tau D\right)
^{k}T_{n}\left(  \zeta\right)  \right)  \widehat{\phi}\left(  \zeta\right)
d\zeta\right)  q_{l}\left(  \tau\right)  d\tau\\
&  =\frac{1}{\sqrt{2\pi}}\frac{1}{c_{n,d}}\int\int\left(  T_{n}\left(
\zeta+\tau\right)  -\sum_{k\leq2n-d}\frac{1}{k!}\left(  \tau D\right)
^{k}T_{n}\left(  \zeta\right)  \right)  q_{l}\left(  \tau\right)  d\tau\text{
}\widehat{\phi}\left(  \zeta\right)  d\zeta,
\end{align*}

and so%
\[
\left[  G_{c},\widehat{\phi}\right]  =\frac{1}{\sqrt{2\pi}}\frac{1}{c_{n,d}%
}\int\int\left(  T_{n}\left(  \zeta+\tau\right)  -\sum_{k\leq2n-d}\frac{1}%
{k!}\left(  \tau D\right)  ^{k}T_{n}\left(  \zeta\right)  \right)
q_{l}\left(  \tau\right)  d\tau\text{ }\widehat{\phi}\left(  \zeta\right)
d\zeta
\]

for all $\phi\in S_{\emptyset,d-1}$. Part 3 of Theorem
\ref{Thm_prop_functnl_on_Son} now implies that%
\[
G_{c}\left(  \zeta\right)  =\frac{1}{\sqrt{2\pi}}\frac{1}{c_{n,d}}\int\left(
T_{n}\left(  \zeta+\tau\right)  -\sum_{k\leq2n-d}\frac{1}{k!}\left(  \tau
D\right)  ^{k}T_{n}\left(  \zeta\right)  \right)  q_{l}\left(  \tau\right)
d\tau+p_{c}\left(  \zeta\right)  ,
\]

for some polynomial $p_{c}$ with $\deg p_{c}\leq d-2$.

If we can show that%
\[
\lim_{\left\vert \zeta\right\vert \rightarrow\infty}\int\left(  T_{n}\left(
\zeta+\tau\right)  -\sum_{k\leq2n-d}\frac{1}{k!}\left(  \tau D_{\zeta}\right)
^{k}T_{n}\left(  \zeta\right)  \right)  q_{l}\left(  \tau\right)  d\tau=0,
\]

then $G_{c}\left(  \infty\right)  =0$ implies the polynomial $p_{c}$ is zero.

Thus we have proven:

\begin{theorem}
\label{Thm_centdiff_basis_Tn_q_suppbnd}Suppose $q,n$ and $l$ define a
multivariate central difference weight function introduced in Definition
\ref{Def_multi_centdiff_wt_fn}. Then if\textbf{\ }$q_{l}$ is defined by
\ref{a012}, the basis function is given by%
\begin{equation}
G_{c}\left(  \zeta\right)  =\frac{1}{\sqrt{2\pi}}\frac{1}{c_{n,d}}\int\left(
T_{n}\left(  \zeta+\tau\right)  -\sum_{k\leq2n-d}\frac{1}{k!}\left(  \tau
D_{\zeta}\right)  ^{k}T_{n}\left(  \zeta\right)  \right)  q_{l}\left(
\tau\right)  d\tau+p_{c}\left(  \zeta\right)  ,\label{a015}%
\end{equation}

where $c_{n,d}$ is given by Lemma \ref{Lem_eval_thin_plate_splin_const} and
$p_{c}$ is a unique polynomial with $\deg p_{c}\leq d-2$.

Further%
\begin{align}
G_{c}\left(  \zeta\right)   &  =\frac{1}{\sqrt{2\pi}}\frac{1}{c_{n,d}}\frac
{1}{\left(  2n-d\right)  !}\int\left(  \int_{0}^{1}g_{2n-d}\left(  s\right)
\left(  \left(  \widehat{\tau}D\right)  ^{2n-d+1}T_{n}\right)  \left(
\zeta+s\tau\right)  ds\right)  \left\vert \tau\right\vert ^{2n-d+1}%
q_{l}\left(  \tau\right)  d\tau+\nonumber\\
&  \qquad\qquad+p_{c}\left(  \zeta\right)  .\label{a017}%
\end{align}

\end{theorem}

\begin{proof}
The proof of \ref{a015} precedes the statement of this theorem. It remains to
prove \ref{a017}.

From Lemma \ref{Lem_Taylor_estim_L1loc_L1}: Since $D^{\alpha}T_{n}\in
L_{loc}^{1}$ for $\left\vert \alpha\right\vert \leq2n-d+1$, in the sense of
distributions%
\[
T_{n}\left(  \zeta+\tau\right)  -\sum_{k\leq2n-d}\frac{\left(  \tau D\right)
^{k}}{k!}T_{n}\left(  \zeta\right)  =\frac{1}{\left(  2n-d\right)  !}\int%
_{0}^{1}g_{2n-d}\left(  s\right)  \left(  \left(  \tau D\right)
^{2n-d+1}T_{n}\right)  \left(  \zeta+s\tau\right)  ds.
\]

Thus \ref{a015} becomes%
\begin{align*}
&  G_{c}\left(  \zeta\right) \\
&  =\frac{1}{\sqrt{2\pi}}\frac{1}{c_{n,d}}\frac{1}{\left(  2n-d\right)  !}%
\int\left(  \int_{0}^{1}g_{2n-d}\left(  s\right)  \left(  \left(  \tau
D\right)  ^{2n-d+1}T_{n}\right)  \left(  \zeta+s\tau\right)  ds\right)
q_{l}\left(  \tau\right)  d\tau+p_{c}\left(  \zeta\right) \\
&  =\frac{1}{\sqrt{2\pi}}\frac{1}{c_{n,d}}\frac{1}{\left(  2n-d\right)  !}%
\int\left(  \int_{0}^{1}g_{2n-d}\left(  s\right)  \left(  \left(
\widehat{\tau}D\right)  ^{2n-d+1}T_{n}\right)  \left(  \zeta+s\tau\right)
ds\right)  \left\vert \tau\right\vert ^{2n-d+1}q_{l}\left(  \tau\right)
d\tau+p_{c}\left(  \zeta\right)  ,
\end{align*}

as claimed.
\end{proof}

\begin{remark}
Compare this result with Theorem \ref{Thm_centdiffbasis_centdiff_formula} above.
\end{remark}

\begin{remark}
?? To calculate the constant try $\zeta=t\widehat{\mathbf{1}}$. USE \ref{a08}.
Use techniques used to derive Theorem \ref{Thm_CentDiff_eqn_thinps}.

\textbf{SEE CALCULATIONS IN THE \ldots\_2014\_01\_09 DOCUMENT}.
\end{remark}

\begin{remark}
Perhaps Section \textbf{Sect\_invers\_Fourier} of the positive order document
might have useful results. The equation \ref{a014} can be written in terms of
the translation operator and the operator $\mathcal{Q}_{\emptyset,2n-d}$.

From \ref{a00} and \ref{a035}: if $\phi\in S_{\emptyset,d-1}$ then%
\[
\left[  \widehat{G_{c}},\phi\right]  =\frac{1}{c_{n,d}\sqrt{2\pi}}\int\left[
T_{n},\widehat{\phi}\left(  \cdot-\tau\right)  -\sum_{k\leq2n-d}\frac{1}%
{k!}\left(  \tau D\right)  ^{k}\widehat{\phi}\right]  q_{l}\left(
\tau\right)  d\tau.
\]

\end{remark}

\section{??? Characterizing the data space locally -
unfinished\label{Sect_LocDataSpace}}

In this section I have \textbf{made a start} at characterizing the data space
locally. This is by analogy with the results of Section
\ref{Sect_local_data_space}. Indeed, from \ref{a2.4} the tensor product
central difference operator satisfies%
\[
\left\Vert \delta_{2\mathbf{1}}^{2l\mathbf{1}}f\right\Vert _{w_{s},0}%
\leq4^{ld}\left\Vert D^{n\mathbf{1}}f\right\Vert _{2}\leq4^{ld}\left\Vert
f\right\Vert _{W^{n\mathbf{1}}}.
\]

Here the multivariate central difference operator satisfies%
\begin{align*}
w_{c}\left(  \xi\right)   & =\frac{\left\vert \xi\right\vert ^{2n}}%
{\Delta_{2l}\widehat{q}\left(  \xi\right)  },\\
\left\Vert \delta_{2\mathbf{1}}^{2l}f\right\Vert _{w_{c},0}^{2}  & =\int
w_{c}\left\vert \widehat{\delta_{2\mathbf{1}}^{2l}f}\right\vert ^{2}=\int%
\frac{\left\vert \cdot\right\vert ^{2n}}{\Delta_{2l}\widehat{q}}\left\vert
\widehat{\delta_{2\mathbf{1}}^{2l}f}\right\vert ^{2}.
\end{align*}

But%
\begin{align*}
\delta_{2\mathbf{1}}f\left(  x\right)   & =f\left(  x+\frac{2}{2}%
\mathbf{1}\right)  -f\left(  x-\frac{2}{2}\mathbf{1}\right)  =f\left(
x+1\right)  -f\left(  x-1\right)  \Longrightarrow\\
\widehat{\delta_{2\mathbf{1}}f}\left(  \xi\right)   & =\left(  e^{i\mathbf{1}%
\xi}-e^{-i\mathbf{1}\xi}\right)  \widehat{f}\left(  \xi\right)  =2i\sin
\mathbf{1}\xi\text{\thinspace}\widehat{f}\left(  \xi\right)  \Longrightarrow\\
\widehat{\delta_{2\mathbf{1}}^{2l}f}\left(  \xi\right)   & =\left(  -1\right)
^{l}2^{2l}\sin^{2l}\mathbf{1}\xi\text{\thinspace}\widehat{f}\left(
\xi\right)  ,
\end{align*}

\begin{align*}
\left\Vert \delta_{2\mathbf{1}}^{2l}f\right\Vert _{w_{c},0}^{2}  & =\int
w_{c}\left(  \xi\right)  2^{4l}\sin^{4l}\mathbf{1}\xi\text{\thinspace
}\left\vert \widehat{f}\left(  \xi\right)  \right\vert ^{2}d\xi\\
& =2^{4l}\int w_{c}\left(  \xi\right)  \frac{\sin^{4l}\mathbf{1}\xi
}{\left\vert \xi\right\vert ^{2n}}\left\vert \xi\right\vert ^{2n}\left\vert
\widehat{f}\left(  \xi\right)  \right\vert ^{2}d\xi\\
& \leq2^{4l}\left\Vert w_{c}\left(  \xi\right)  \frac{\sin^{4l}\mathbf{1}\xi
}{\left\vert \xi\right\vert ^{2n}}\right\Vert _{\infty}\int\left\vert
\xi\right\vert ^{2n}\left\vert \widehat{f}\left(  \xi\right)  \right\vert
^{2}d\xi\\
& =2^{4l}\left\Vert w_{c}\left(  \xi\right)  \frac{\sin^{4l}\mathbf{1}\xi
}{\left\vert \xi\right\vert ^{2n}}\right\Vert _{\infty}\sum_{\left\vert
\alpha\right\vert =n}\frac{1}{\alpha!}\int\left\vert \widehat{D^{\alpha}%
f}\right\vert ^{2}.
\end{align*}

Also%
\begin{align*}
&  \left\Vert \delta_{2\mathbf{1}}^{2l}f\right\Vert _{w_{c},0}^{2}\\
&  =\int w_{c}\left(  \xi\right)  2^{4l}\sin^{4l}\mathbf{1}\xi\left\vert
\widehat{f}\left(  \xi\right)  \right\vert ^{2}d\xi\\
&  =2^{4l}\int\limits_{\left\vert \xi\right\vert \leq r}w_{c}\left(
\xi\right)  \sin^{4l}\mathbf{1}\xi\left\vert \widehat{f}\left(  \xi\right)
\right\vert ^{2}d\xi+2^{4l}\int\limits_{\left\vert \xi\right\vert \geq r}%
w_{c}\left(  \xi\right)  \sin^{4l}\mathbf{1}\xi\left\vert \widehat{f}\left(
\xi\right)  \right\vert ^{2}d\xi\\
&  \leq\int\limits_{\left\vert \xi\right\vert \leq r}w_{c}\left(  \xi\right)
\left\vert \mathbf{1}\xi\right\vert ^{4l}\left\vert \widehat{f}\left(
\xi\right)  \right\vert ^{2}d\xi+\int\limits_{\left\vert \xi\right\vert \geq
r}\frac{w_{c}\left(  \xi\right)  }{\left\vert \xi\right\vert ^{2n}}\left\vert
\xi\right\vert ^{2n}\left\vert \widehat{f}\left(  \xi\right)  \right\vert
^{2}d\xi\\
&  \leq\left\vert \mathbf{1}\right\vert ^{4l}\int\limits_{\left\vert
\xi\right\vert \leq r}w_{c}\left(  \xi\right)  \left\vert \xi\right\vert
^{4l}\left\vert \widehat{f}\left(  \xi\right)  \right\vert ^{2}d\xi
+\int\limits_{\left\vert \xi\right\vert \geq r}\frac{w_{c}\left(  \xi\right)
}{\left\vert \xi\right\vert ^{2n}}\left\vert \xi\right\vert ^{2n}\left\vert
\widehat{f}\left(  \xi\right)  \right\vert ^{2}d\xi\\
&  \leq\left(  d^{2l}\max_{B_{r}}w_{c}\left\vert \cdot\right\vert
^{4l}\right)  \int\limits_{\left\vert \cdot\right\vert \leq r}\left\vert
\widehat{f}\right\vert ^{2}+\left(  \max_{B_{r}^{c}}\frac{w_{c}}{\left\vert
\cdot\right\vert ^{2n}}\right)  \int\limits_{\left\vert \cdot\right\vert \geq
r}\left\vert \cdot\right\vert ^{2n}\left\vert \widehat{f}\right\vert ^{2}\\
&  \leq\left(  d^{2l}\max_{B_{r}}w_{c}\left\vert \cdot\right\vert
^{4l}\right)  \int\left\vert \widehat{f}\right\vert ^{2}+\left(  \max
_{B_{r}^{c}}\frac{w_{c}}{\left\vert \cdot\right\vert ^{2n}}\right)
\int\left\vert \cdot\right\vert ^{2n}\left\vert \widehat{f}\right\vert ^{2}\\
&  \leq\max\left\{  d^{2l}\max_{B_{r}}w_{c}\left\vert \cdot\right\vert
^{4l},\max_{B_{r}^{c}}\frac{w_{c}}{\left\vert \cdot\right\vert ^{2n}}\right\}
\int\left(  1+\left\vert \cdot\right\vert ^{2}\right)  ^{n}\left\vert
\widehat{f}\right\vert ^{2}\\
&  =\max\left\{  d^{2l}\max_{B_{r}}w_{c}\left\vert \cdot\right\vert ^{4l}%
,\max_{B_{r}^{c}}\frac{w_{c}}{\left\vert \cdot\right\vert ^{2n}}\right\}
\left\Vert f\right\Vert _{W^{n}}^{2}.
\end{align*}

\textbf{If }$q$\textbf{\ is radial} then $\widehat{q}$ is radial and we can
expand $\widehat{q}$\ about the origin using Theorem \ref{Thm_Tay_rem_zeros}
and obtain a remainder containing the factor $\left\vert \cdot\right\vert
^{2n}$.

\chapter{The Exact smoother\label{Ch_Exact_smth}}

\section{Introduction}

We call this well-known basis function smoother the \textbf{Exact smoother}
because the smoother studied in the next chapter approximates it. We will
assume the \textbf{basis function is real-valued} and that the \textbf{weight
function has property W02}. The Exact smoother minimizes the functional%
\[
\rho\left\Vert f\right\Vert _{w,0}^{2}+\frac{1}{N}\sum_{i=1}^{N}\left\vert
f(x^{(i)})-y_{i}\right\vert ^{2},\quad f\in X_{w}^{0},
\]

over the data (function) space $X_{w}^{0}$ where $\rho>0$ is termed the
smoothing coefficient. It is shown that this problem, like the interpolation
problem, has a unique basis function solution in the space $W_{G,X}$. The
finite dimensionality of the solution allows us to derive matrix equations for
the coefficients $\alpha_{i}$ of the data-translated basis functions.

Following the example of the interpolation semi-inner product $\left(
f-\mathcal{I}_{X}f,g\right)  _{w,0}$ of Definition \ref{Def_interp_seminorm}
of \ref{Ch_Interpol} and it's sequel I will define the semi-inner product
$\left(  \mathcal{S}_{X}f,g\right)  _{w,0}$ and inner product $\left(
f-\mathcal{S}_{X}f,g\right)  _{w,0}$ and use these to study the value and the
error of the Exact smoother respectively. These expressions are useful because
$\left(  \mathcal{S}_{X}f,R_{x}\right)  _{w,0}=\left(  \mathcal{S}%
_{X}f\right)  \left(  x\right)  $ and $\left(  f-\mathcal{S}_{X}%
f,R_{x}\right)  _{w,0}=\left(  f-\mathcal{S}_{X}f\right)  \left(  x\right)  $
and because of the power of Hilbert space theory. An outcome of this approach
are the general Exact smoother bounds \ref{7.68} and the convergence estimates
of \ Corollary \ref{Cor_Thm_Rx(x)minus(SRx)(x)_basis_fn_bound}.

Then using these estimates three estimates are derived for the order of the
pointwise error of the Exact smoother w.r.t. its data function: Types 1 and 2
estimates and the distribution Taylor series estimates\ do not assume
unisolvent $X$ data and one does. When $\rho=0$ these estimates correspond to
interpolation results of Section \ref{Sect_interp_no_Taylor_converg} and the
Exact smoother convergence orders and the constants are the same as those for
the interpolation case which are given in the interpolation tables
\ref{Tbl_intro_NonUnisolvTyp1Converg}, \ref{Tbl_intro_NonUnisolvTyp2Conv},
\ref{Tbl_intro_UnisolvConverg} and \ref{Tbl_intro_ConvergCentral}.\medskip

\fbox{Type 1 error estimates} When the weight function has property W02 for
$\kappa\geq0$ it is assumed that the basis function satisfies an inequality of
the form \ref{1.9} and that the data region is a closed bounded infinite set
$K$. In this case it is shown in Corollary
\ref{Cor_OrdConvergExactSmth_k=0_Fd} that the Exact smoother $s_{e}$ of the
data function $f$ satisfies the error estimate%
\[
\left\vert f\left(  x\right)  -s_{e}\left(  x\right)  \right\vert
\leq\left\Vert f\right\Vert _{w,0}\left(  \sqrt{\rho N}+k_{G}\left(
h_{X,K}\right)  ^{s}\right)  ,\quad x\in K,
\]

when $h_{X,K}=\sup\limits_{x\in K}\operatorname*{dist}\left(  x,X\right)  \leq
h_{G}$ and $k_{G}=\left(  2\pi\right)  ^{-d/4}\sqrt{2C_{G}}$.\medskip

\fbox{Type 2 error estimates} If it is only assumed that $\kappa\geq1$ then it
is shown in Theorem \ref{Thm_G(x)minusG(0)_bound} that%
\[
\left\vert f\left(  x\right)  -s_{e}\left(  x\right)  \right\vert
\leq\left\Vert f\right\Vert _{w,0}\left(  \sqrt{\rho N}+k_{G}h_{X,K}\right)
,\quad x\in K,
\]

when $h_{X,K}<\infty$ and $k_{G}=\left(  2\pi\right)  ^{-d/4}\sqrt{-\left(
\left\vert D\right\vert ^{2}G\right)  \left(  0\right)  }\sqrt{d}$.\medskip

\fbox{Error estimates using Taylor series and unisolvency} Suppose the weight
function has property W02 for some parameter $\kappa\geq1$ and $X$ is
contained in a bounded data region $\Omega$. For sufficiently dense data each
sufficiently small ball contains a minimal unisolvent subset of $X$ of order
$\left\lfloor \kappa\right\rfloor $ Then using results from the Lagrange
theory of Lagrange interpolation we show in Theorem
\ref{Thm_Exact_smth_ord_gte_1} that there exist constants $K_{\Omega
,m}^{\prime},k_{G}>0$ such that
\[
\left\vert f\left(  x\right)  -s_{e}\left(  x\right)  \right\vert
\leq\left\Vert f\right\Vert _{w,0}\left(  K_{\Omega,m}^{\prime}\sqrt{\rho
N}+k_{G}\left(  h_{X,\Omega}\right)  ^{\left\lfloor \kappa\right\rfloor
}\right)  ,\quad x\in\overline{\Omega},
\]

when $h_{X,\Omega}=\sup\limits_{x\in\Omega}\operatorname*{dist}\left(
x,X\right)  \leq h_{G}$.

These theoretical results will be illustrated using the weight function
examples from the interpolation chapter, namely the radial \textbf{shifted
thin-plate splines}, \textbf{Gaussian and Sobolev splines} and the tensor
product \textbf{extended B-splines}. We will also use the \textbf{central
difference} weight functions from Chapter \ref{Ch_cent_diff_wt_fn_ten_prod}.

\textbf{Numerical results are only presented for the Type 1 and Type 2
estimates}. Numeric experiments are carried out using the same 1-dimensional
B-splines and data functions that were used for the interpolants. We restrict
ourselves to one dimension so that the data density parameters $h_{X,\Omega}$
and $h_{X,K}$ can be easily calculated. Also note that in one dimension
equation \ref{a1.52} implies that any set of $m$ distinct points is $m$-unisolvent.

No numerical experiments will be carried out for the unisolvent data case
because I did not want to numerically estimate the parameters related to the
Lagrange polynomial interpolation. However satisfactory easy estimates
\textbf{may} be possible.

\section{The Exact Smoothing problem\label{Sect_Exact_smth_prob}}

The Exact smoothing problem will be defined as a variational minimization
problem within the space of continuous functions $X_{w}^{0}$. The functional
is constructed using a real smoothing parameter $\rho>0$ and a set of distinct
scattered data points $\left\{  \left(  x^{(k)},y_{k}\right)  \right\}
_{k=1}^{N}$, $x^{(k)}\in\mathbb{R}^{d}$ and $y_{k}\in\mathbb{C}$. The
functional to be minimized is%
\begin{equation}
J_{e}[f]=\rho\left\Vert f\right\Vert _{w,0}^{2}+\frac{1}{N}\sum_{k=1}%
^{N}\left\vert f(x^{(k)})-y_{k}\right\vert ^{2},\quad f\in X_{w}%
^{0}.\label{7.63}%
\end{equation}

The independent data will be specified by $X=\left(  x^{(k)}\right)
_{k=1}^{N} $ and the dependent data by $y=\left(  y_{k}\right)  _{k=1}^{N}$.
The Exact smoothing problem is now stated as:\medskip%
\begin{equation}%
\begin{tabular}
[c]{|l|}\hline
\textbf{The Exact smoothing problem}\\\hline
Minimize$\text{ }$the$\text{ }$smoothing$\text{ }$functional\ $J_{e}[f]\text{
for\ }f\in X_{w}^{0}$.\\\hline
\end{tabular}
\label{7.51}%
\end{equation}
\medskip

Here the first component of the functional $\rho\left\Vert f\right\Vert
_{w,0}^{2}$ is the \textbf{global smoothing component} and the second
component is the \textbf{localizing least squares component}. As
$\rho\rightarrow0^{+}$ it will be shown that the matrix equation representing
this problem becomes that of the minimal norm interpolation problem
\ref{1.42}. As the value of $\rho$ increases numerical experiments show the
degree of smoothing increases and the correlation between data and the
smoothing function decreases.

Using Hilbert space techniques I will now do the following:

\begin{enumerate}
\item Show there exists a unique Exact smoother.

\item Show the smoother is a basis function smoother i.e. show that it lies in
$W_{G,X}$.

\item Construct a matrix equation for the coefficients of the data-translated
basis functions $G\left(  \cdot-x^{\left(  k\right)  }\right)  $.
\end{enumerate}

The proofs will be carried out within a Hilbert space framework by formulating
the smoothing functional $J_{e}$ in terms of a norm $\left\Vert \cdot
\right\Vert _{V}$ on the Hilbert product space $V=X_{w}^{0}\otimes
\mathbb{C}^{N}$. To this end I first introduce the following definitions.

\begin{definition}
\label{Def_spaceV_opLx}\textbf{The space} $V$ \textbf{and the operator}
$\mathcal{L}_{X}$

\begin{enumerate}
\item Let $V=X_{w}^{0}\otimes\mathbb{C}^{N}$ be the Hilbert product space with
norm $\left\Vert \cdot\right\Vert _{V}$ and inner product $\left(  \cdot
,\cdot\right)  _{V}$ given by
\[
\left(  \left(  u_{1},u_{2}\right)  ,\left(  v_{1},v_{2}\right)  \right)
_{V}=\rho\,\left(  u_{1},v_{1}\right)  _{w,0}+\frac{1}{N}\left(  u_{2}%
,v_{2}\right)  _{\mathbb{C}^{N}}.
\]

\item Let the operator $\mathcal{L}_{X}:X_{w}^{0}\rightarrow V$ be defined by
$\mathcal{L}_{X}f=\left(  f,\widetilde{\mathcal{E}}_{X}f\right)  $ where
$\widetilde{\mathcal{E}}_{X}$ is the vector-valued evaluation operator
$\widetilde{\mathcal{E}}_{X}f=\left(  f\left(  x^{\left(  k\right)  }\right)
\right)  $ of Definition \ref{Def_vect_val_eval_op}.
\end{enumerate}
\end{definition}

\begin{remark}
\label{Rem_SmoothFunc1}\ 

\begin{enumerate}
\item The definitions above were constructed so I have the following result:
If $y=\left(  y_{i}\right)  $ is the (complex) dependent data given in the
smoothing problem and if I let $\varsigma=(0,y)$ then for $f\in X_{w}^{0}$ the
Exact smoothing functional can be rewritten in terms of the norm $\left\Vert
\cdot\right\Vert _{V}$ as
\[
J_{e}[f]=\rho\left\Vert f\right\Vert _{w,0}^{2}+\frac{1}{N}\left\vert
\widetilde{\mathcal{E}}_{X}f-y\right\vert _{\mathbb{C}^{N}}^{2}=\left\Vert
\left(  f,\widetilde{\mathcal{E}}_{X}f-y\right)  \right\Vert _{V}%
^{2}=\left\Vert \left(  f,\widetilde{\mathcal{E}}_{X}f\right)  -\left(
0,y\right)  \right\Vert _{V}^{2}=\left\Vert \mathcal{L}_{X}f-\varsigma
\right\Vert _{V}^{2}.
\]

The Exact smoother is now the (unique) orthogonal projection of the constant
$\left(  0,y\right)  $ onto the infinite dimensional subspace $\mathcal{L}%
_{X}\left(  X_{w}^{0}\right)  $.

\item ??? \textbf{IDEA} Suppose $\Theta_{X}=\left(  \Theta_{x^{\left(
i\right)  }}\right)  _{i=1}^{N}\subset X_{1/w}^{0}\simeq\left(  X_{w}%
^{0}\right)  ^{\prime}$. Reformulate $\mathcal{L}_{X}$ as%
\[
\mathcal{L}_{X}f=\left(  f,\Theta_{X}f\right)  .
\]

\end{enumerate}
\end{remark}

I will need the following properties of the operators $\mathcal{L}_{X}$,
$\mathcal{L}_{X}^{\ast}$ and $\mathcal{L}_{X}^{\ast}\mathcal{L}_{X}$.

\begin{theorem}
\label{Thm_ex_L_op_properties}The operator $\mathcal{L}_{X}$ has the following properties:

\begin{enumerate}
\item $\left\Vert \mathcal{L}_{X}f\right\Vert _{V}$ and $\left\Vert
f\right\Vert _{w,0}$ are equivalent norms on $X_{w}^{0}$.

\item $\mathcal{L}_{X}:X_{w}^{0}\rightarrow V$ is continuous, 1-1 and has
closed range. ?? FINISH PROOF! In fact,
\[
\operatorname{range}\mathcal{L}_{X}=\left(  W_{G,X}^{\bot}\otimes\left\{
\mathbf{0}\right\}  \right)  +\left\{  \left(  \mathcal{I}_{X}\alpha
,\alpha\right)  :\alpha\in\mathbb{C}^{N}\right\}  .
\]

\item $\mathcal{L}_{X}^{\ast}:V\mathbb{\rightarrow}X_{w}^{0}$ and if
$u=\left(  u_{1},u_{2}\right)  \in V$ then
\[
\mathcal{L}_{X}^{\ast}u=\rho u_{1}+\frac{1}{N}\widetilde{\mathcal{E}}%
_{X}^{\ast}u_{2}.
\]

\item $\mathcal{L}_{X}^{\ast}\mathcal{L}_{X}:X_{w}^{0}\rightarrow X_{w}^{0}$
and $\mathcal{L}_{X}^{\ast}\mathcal{L}_{X}:W_{G,X}\rightarrow W_{G,X}$ and
\[
\mathcal{L}_{X}^{\ast}\mathcal{L}_{X}=\rho I+\frac{1}{N}\widetilde{\mathcal{E}%
}_{X}^{\ast}\widetilde{\mathcal{E}}_{X}.
\]

\item $\mathcal{L}_{X}^{\ast}\mathcal{L}_{X}$ is continuous and 1-1 on
$X_{w}^{0}$.
\end{enumerate}
\end{theorem}

\begin{proof}
\textbf{Part 1} By part 3 of Theorem \ref{Thm_eval_op_properties}, $\left\Vert
\widetilde{\mathcal{E}}_{X}\right\Vert =\left\Vert R_{X,X}\right\Vert $ so
\begin{align*}
\left\Vert \mathcal{L}_{X}f\right\Vert _{V}^{2}=\left\Vert \left(
f,\widetilde{\mathcal{E}}_{X}f\right)  \right\Vert _{V}^{2}=\rho\left\Vert
f\right\Vert _{w,0}^{2}+\frac{1}{N}\left\vert \widetilde{\mathcal{E}}%
_{X}f\right\vert _{\mathbb{C}^{N}}^{2} &  \leq\rho\left\Vert f\right\Vert
_{w,0}^{2}+\frac{1}{N}\left\Vert R_{X,X}\right\Vert ^{2}\left\Vert
f\right\Vert _{w,0}^{2}\\
&  =\left(  \rho+\frac{1}{N}\left\Vert R_{X,X}\right\Vert ^{2}\right)
\left\Vert f\right\Vert _{w,0}^{2}.
\end{align*}

Further, since $\rho>0$
\[
\left\Vert f\right\Vert _{w,0}^{2}\leq\left(  \min\left\{  \rho,\frac{1}%
{N}\right\}  \right)  ^{-1}\left(  \rho\left\Vert f\right\Vert _{w,0}%
^{2}+\frac{1}{N}\left\vert \widetilde{\mathcal{E}}_{X}f\right\vert
_{\mathbb{C}^{N}}^{2}\right)  =\left(  \min\left\{  \rho,\frac{1}{N}\right\}
\right)  ^{-1}\left\Vert \mathcal{L}_{X}f\right\Vert _{V}^{2}.
\]
\medskip

\textbf{Part 2} $\mathcal{L}_{X}:X_{w}^{0}\rightarrow V$ is continuous since
$\widetilde{\mathcal{E}}_{X}$ is continuous by part 1 of Theorem
\ref{Thm_eval_op_properties}. From part 1 there exists a constant $C_{1}$ such
that $\left\Vert f\right\Vert _{w,0}\leq C_{1}\left\Vert \mathcal{L}%
_{X}f\right\Vert _{V}$ and hence the range of $\mathcal{L}_{X}$ is closed and
$\mathcal{L}_{X}$ is 1-1. ?? \textbf{FINISH proof!}\medskip

\textbf{Part 3\ }Since $\mathcal{L}_{X}$ is a continuous operator the adjoint
$\mathcal{L}_{X}^{\ast}:V\rightarrow X_{w}^{0}$ exists and is continuous.
Further, if $f\in X_{w}^{0}$ and $u=\left(  u_{1,}u_{2}\right)  \in V$
\begin{align*}
\left(  \mathcal{L}_{X}f,u\right)  _{V}=\left(  \left(
f,\widetilde{\mathcal{E}}_{X}f\right)  ,\left(  u_{1,}u_{2}\right)  \right)
_{V} &  =\rho\,\left(  f,u_{1}\right)  _{w,0}+\frac{1}{N}\left(
\widetilde{\mathcal{E}}_{X}f,u_{2}\right) \\
&  =\rho\,\left(  f,u_{1}\right)  _{w,0}+\frac{1}{N}\left(
f,\widetilde{\mathcal{E}}_{X}^{\ast}u_{2}\right)  _{w,0}\\
&  =\left(  f,\rho u_{1}+\frac{1}{N}\widetilde{\mathcal{E}}_{X}^{\ast}%
u_{2}\right)  _{w,0}.
\end{align*}

Hence $\mathcal{L}_{X}^{\ast}u=\rho u_{1}+\frac{1}{N}\widetilde{\mathcal{E}%
}_{X}^{\ast}u_{2}$, where $u=\left(  u_{1,}u_{2}\right)  \in V$.\medskip

\textbf{Part 4} If $f\in X_{w}^{0}$ then%
\[
\mathcal{L}_{X}^{\ast}\mathcal{L}_{X}f=\rho\left(  \mathcal{L}_{X}f\right)
_{1}+\frac{1}{N}\widetilde{\mathcal{E}}_{X}^{\ast}\left(  \mathcal{L}%
_{X}f\right)  _{2}=\rho\left(  f,\widetilde{\mathcal{E}}_{X}f\right)
_{1}+\frac{1}{N}\widetilde{\mathcal{E}}_{X}^{\ast}\left(
f,\widetilde{\mathcal{E}}_{X}f\right)  _{2}=\rho f+\frac{1}{N}%
\widetilde{\mathcal{E}}_{X}^{\ast}\widetilde{\mathcal{E}}_{X}f.
\]

By part 2 of Theorem \ref{Thm_eval_op_properties} $\widetilde{\mathcal{E}}%
_{X}^{\ast}:\mathbb{C}^{N}\rightarrow W_{G,X}$ and it follows that
$\mathcal{L}_{X}^{\ast}\mathcal{L}_{X}:X_{w}^{0}\rightarrow X_{w}^{0}$ and
$\mathcal{L}_{X}^{\ast}\mathcal{L}_{X}:W_{G,X}\rightarrow W_{G,X}$.\medskip

\textbf{Part 5} Clearly $\mathcal{L}_{X}^{\ast}\mathcal{L}_{X}$ is continuous.
Now suppose $\mathcal{L}_{X}^{\ast}\mathcal{L}_{X}f=0$ for $f\in X_{w}^{0}$.
Then
\[
0=\left(  \mathcal{L}_{X}^{\ast}\mathcal{L}_{X}f,f\right)  _{w,0}=\left(
\mathcal{L}_{X}f,\mathcal{L}_{X}f\right)  _{w,0}=\left\Vert \mathcal{L}%
_{X}f\right\Vert _{V}^{2},
\]

so that $\mathcal{L}_{X}f=0$ and hence $f=0$ since $\mathcal{L}_{X}$ is 1-1 by
part 2. Therefore $\mathcal{L}_{X}^{\ast}\mathcal{L}_{X}$ is 1-1.
\end{proof}

I next show that the Exact smoother, like the minimum norm interpolant of
Chapter \ref{Ch_Interpol}, has the nice property of being in the N-dimensional
basis function space
\[
W_{G,X}=\left\{  \sum\limits_{i=1}^{N}\alpha_{i}G\left(  \cdot-x^{(i)}\right)
:\alpha_{i}\in\mathbb{C}\right\}  ,
\]
introduced in Definition \ref{Def_Wg,x_m_gt_0}. It was shown in Theorem
\ref{Thm_indep_G(x-xi)} that the data-translated functions $G\left(
\cdot-x^{(i)}\right)  $ are linearly independent so $\dim W_{G,X}=N$.

\begin{theorem}
\label{Thm_ex_smooth_var}Fix $y\in\mathbb{C}^{N}$ and let $\varsigma=(0,y)\in
V$. Then there exists a unique solution to the Exact smoothing problem
\ref{7.51} in $X_{w}^{0}$. Denote this solution by $s_{e}$. This solution has
the following properties:

\begin{enumerate}
\item $\left\Vert \mathcal{L}_{X}s_{e}-\varsigma\right\Vert _{V}<\left\Vert
\mathcal{L}_{X}f-\varsigma\right\Vert _{V}$,\quad for all $f\in X_{w}%
^{0}-\{s_{e}\}$.

\item $\left(  \mathcal{L}_{X}s_{e}-\varsigma,\mathcal{L}_{X}\left(
s_{e}-f\right)  \right)  _{V}=0$,\quad for all $f\in X_{w}^{0}$.

\item $\left\Vert \mathcal{L}_{X}s_{e}-\varsigma\right\Vert _{V}%
^{2}+\left\Vert \mathcal{L}_{X}\left(  s_{e}-f\right)  \right\Vert _{V}%
^{2}=\left\Vert \mathcal{L}_{X}f-\varsigma\right\Vert _{V}^{2}$,\quad for all
$f\in X_{w}^{0}$.\smallskip

The last equality is equivalent to the equality of part 2.

\item $s_{e}=\frac{1}{N}\left(  \mathcal{L}_{X}^{\ast}\mathcal{L}_{X}\right)
^{-1}\widetilde{\mathcal{E}}_{X}^{\ast}y$.

\item $s_{e}\in W_{G,X}$ and%
\begin{equation}
s_{e}=\frac{1}{\rho N}\sum\limits_{k=1}^{N}\left(  y_{k}-s_{e}(x^{(k)}%
)\right)  R_{x^{(k)}}.\label{7.54}%
\end{equation}

\end{enumerate}
\end{theorem}

\begin{proof}
\textbf{Parts 1}, \textbf{2}, \textbf{3\ }Because $\mathcal{L}_{X}(X_{w}^{0})
$ is closed it is well known there exists a unique element $s_{e}$ of the
hyperplane $\mathcal{L}_{X}(X_{w}^{0})$ which minimizes the distance between
the hyperplane $\mathcal{L}_{X}(X_{w}^{0})$ and the point $\varsigma$. From
simple geometric considerations it follows that $s_{e}$ satisfies conditions
1, 2 and 3 of this theorem.\medskip

\textbf{Part 4 }From part 2
\[
0=\left(  \mathcal{L}_{X}s_{e}-\varsigma,\mathcal{L}_{X}s_{e}-\mathcal{L}%
_{X}f\right)  _{V}=\left(  \mathcal{L}_{X}s_{e}-\varsigma,\mathcal{L}%
_{X}\left(  s_{e}-f\right)  \right)  _{V}=\left(  \mathcal{L}_{X}^{\ast
}\left(  \mathcal{L}_{X}s_{e}-\varsigma\right)  ,s_{e}-f\right)  _{w,0},
\]

for all $f\in X_{w}^{0}$. Thus
\begin{equation}
\mathcal{L}_{X}^{\ast}\left(  \mathcal{L}_{X}s_{e}-\varsigma\right)
=0.\label{7.49}%
\end{equation}

Hence
\[
\mathcal{L}_{X}^{\ast}\mathcal{L}_{X}s_{e}=\mathcal{L}_{X}^{\ast}%
\varsigma=\mathcal{L}_{X}^{\ast}\left(  0,y\right)  =\frac{1}{N}%
\widetilde{\mathcal{E}}_{X}^{\ast}y,
\]

and since by part 5 Theorem \ref{Thm_ex_L_op_properties} $\mathcal{L}%
_{X}^{\ast}\mathcal{L}_{X}$ is one-to-one, it follows that $\mathcal{L}%
_{X}^{\ast}\mathcal{L}_{X}$ has an inverse and so $s_{e}=\frac{1}{N}\left(
\mathcal{L}_{X}^{\ast}\mathcal{L}_{X}\right)  ^{-1}\widetilde{\mathcal{E}}%
_{X}^{\ast}y$.\medskip

\textbf{Part 5} Starting with equation \ref{7.49}
\begin{align*}
0=\mathcal{L}_{X}^{\ast}\left(  \mathcal{L}_{X}s_{e}-\varsigma\right)
=\mathcal{L}_{X}^{\ast}\mathcal{L}_{X}s_{e}-\frac{1}{N}\widetilde{\mathcal{E}%
}_{X}^{\ast}y &  =\rho s_{e}+\frac{1}{N}\widetilde{\mathcal{E}}_{X}^{\ast
}\widetilde{\mathcal{E}}_{X}s_{e}-\frac{1}{N}\widetilde{\mathcal{E}}_{X}%
^{\ast}y\\
&  =\rho s_{e}+\frac{1}{N}\widetilde{\mathcal{E}}_{X}^{\ast}\left(
\widetilde{\mathcal{E}}_{X}s_{e}-y\right) \\
&  =\rho s_{e}-\frac{1}{N}\sum_{k=1}^{N}\left(  s_{e}(x^{(k)})-y_{k}\right)
R_{x^{(k)}},
\end{align*}

which proves \ref{7.54}. Finally, since $R_{x^{(k)}}=\left(  2\pi\right)
^{-d/2}G\left(  \cdot-x^{(k)}\right)  $ we have $s_{e}\in W_{G,X}$.
\end{proof}

In a manner analogous to the case of the minimal norm interpolant, part 4 of
the last theorem allows the definition of a mapping between a data function
and it's corresponding Exact smoother.

\begin{definition}
\label{Def_data_func_exact_map}\textbf{Data functions and the Exact smoother
operator} $\mathcal{S}_{X}:X_{w}^{0}\rightarrow W_{G,X}$.

Given an independent data set $X$, I shall assume that each member of
$X_{w}^{0}$ can act as a legitimate data function $g$ and generate the data
vector $\widetilde{\mathcal{E}}_{X}g$ -\textbf{\ }see item 4 of Section
\ref{Sect_num_experim}.

The equation of part 4 of Theorem \ref{Thm_ex_smooth_var} enables us to define
the continuous linear mapping $\mathcal{S}_{X}:X_{w}^{0}\rightarrow W_{G,X}$
from the data functions to the corresponding Exact smoother given by%
\begin{equation}
\mathcal{S}_{X}g=\frac{1}{N}\left(  \mathcal{L}_{X}^{\ast}\mathcal{L}%
_{X}\right)  ^{-1}\widetilde{\mathcal{E}}_{X}^{\ast}\widetilde{\mathcal{E}%
}_{X}g,\quad g\in X_{w}^{0}.\label{7.53}%
\end{equation}

\end{definition}

It was shown in Theorem \ref{Thm_ex_L_op_properties} that $\mathcal{L}%
_{X}^{\ast}\mathcal{L}_{X}:X_{w}^{0}\rightarrow X_{w}^{0}$ is continuous and
1-1. I now prove $\mathcal{L}_{X}^{\ast}\mathcal{L}_{X}$ is onto and hence a
homeomorphism. I also prove some important properties of the Exact smoother
mapping $\mathcal{S}_{X}$.

\begin{corollary}
\label{Cor_ex_Lx*Lx_onto_Xw,th}$\mathcal{L}_{X}^{\ast}\mathcal{L}_{X}%
:X_{w}^{0}\rightarrow X_{w}^{0}$ is a homeomorphism and $\mathcal{L}_{X}%
^{\ast}\mathcal{L}_{X}:W_{G,X}\rightarrow W_{G,X}$ is a homeomorphism.

Regarding the properties of $\mathcal{S}_{X}$:

\begin{enumerate}
\item $\mathcal{S}_{X}g=g-\rho\,\left(  \mathcal{L}_{X}^{\ast}\mathcal{L}%
_{X}\right)  ^{-1}g$,$\quad g\in X_{w}^{0}$,

\item $\mathcal{S}_{X}$ is continuous but not a projection,

\item $\mathcal{S}_{X}:X_{w}^{0}\rightarrow W_{G,X}$ is onto,

\item $\operatorname*{null}\mathcal{S}_{X}=W_{G,X}^{\bot}$,

\item $\mathcal{S}_{X}$ is self-adjoint.
\end{enumerate}
\end{corollary}

\begin{proof}
Suppose $g\in X_{w}^{0}$ and let $s_{e}=\mathcal{S}_{X}g$. From part 5 of
Theorem \ref{Thm_ex_L_op_properties} I know that $\mathcal{L}_{X}^{\ast
}\mathcal{L}_{X}$ is 1-1. From part 4 of Theorem \ref{Thm_ex_smooth_var}
$\mathcal{L}_{X}^{\ast}\mathcal{L}_{X}s_{e}=\frac{1}{N}\widetilde{\mathcal{E}%
}_{X}^{\ast}\widetilde{\mathcal{E}}_{X}g$ and from part 4 of Theorem
\ref{Thm_ex_L_op_properties}, $\mathcal{L}_{X}^{\ast}\mathcal{L}_{X}g=\rho
g+\frac{1}{N}\widetilde{\mathcal{E}}_{X}^{\ast}\widetilde{\mathcal{E}}_{X}g$
so that $\mathcal{L}_{X}^{\ast}\mathcal{L}_{X}g=\rho g+\mathcal{L}_{X}^{\ast
}\mathcal{L}_{X}s_{e}$. Clearly this equation implies part 1 and the equation
\begin{equation}
g=\mathcal{L}_{X}^{\ast}\mathcal{L}_{X}\left(  \frac{g-s_{e}}{\rho}\right)
,\qquad g\in X_{w}^{0}.\label{7.56}%
\end{equation}

This latter equation proves $g\in\operatorname*{range}\mathcal{L}_{X}^{\ast
}\mathcal{L}_{X}$. Thus $\mathcal{L}_{X}^{\ast}\mathcal{L}_{X}:X_{w}%
^{0}\rightarrow X_{w}^{0}$ is onto and by the Open Mapping theorem
$\mathcal{L}_{X}^{\ast}\mathcal{L}_{X}$ is a homeomorphism. If $g\in W_{G,X} $
then \ref{7.56} implies $\mathcal{L}_{X}^{\ast}\mathcal{L}_{X}:W_{G,X}%
\rightarrow W_{G,X}$ is onto and so $\mathcal{L}_{X}^{\ast}\mathcal{L}%
_{X}:W_{G,X}\rightarrow W_{G,X}$ is also a homeomorphism. This fact and part 1
clearly imply that $\mathcal{S}_{X}$ is not a projection.

Further, in Theorem \ref{Thm_eval_op_properties} it was shown that
$\widetilde{\mathcal{E}}_{X}^{\ast}\widetilde{\mathcal{E}}_{X}:X_{w}%
^{0}\rightarrow W_{G,X}$ is onto and $\operatorname*{null}%
\widetilde{\mathcal{E}}_{X}^{\ast}\widetilde{\mathcal{E}}_{X}=W_{G,X}^{\bot}$.
The stated properties of $\mathcal{S}_{X}$ now follow from \ref{7.53} and the
fact that $\widetilde{\mathcal{E}}_{X}^{\ast}\widetilde{\mathcal{E}}_{X}$ and
$\mathcal{L}_{X}^{\ast}\mathcal{L}_{X}$ are self-adjoint.
\end{proof}

\begin{remark}
\label{Rem_AnotherSmthProblem}\textbf{Another smoothing problem} The Exact
smoothing functional used in this document is different from that used in
Narcowich, Ward and Wendland \cite{NarcWardWend2004}. They minimize the
functional%
\[
J_{e}^{\prime}[f]=\rho\left\Vert f\right\Vert _{w,0}^{2}+\sum_{k=1}%
^{N}\left\vert f(x^{(k)})-y_{k}\right\vert ^{2},\quad f\in X_{w}^{0}.
\]

The independent data will be specified by $X=\left(  x^{(k)}\right)
_{k=1}^{N} $ and the dependent data by $y=\left(  y_{k}\right)  _{k=1}^{N}$.
The smoothing problem is now stated as:\medskip%
\[%
\begin{tabular}
[c]{|l|}\hline
\textbf{Another smoothing problem}\\\hline
Minimize$\text{ }$the$\text{ }$smoothing$\text{ }$functional\ $J_{e}^{\prime
}[f]\text{ for\ }f\in X_{w}^{0}$.\\\hline
\end{tabular}
\]

Since%
\[
J_{e}^{\prime}[f]=N\left(  \rho^{\prime}\left\Vert f\right\Vert _{w,0}%
^{2}+\frac{1}{N}\sum_{k=1}^{N}\left\vert f(x^{(k)})-y_{k}\right\vert
^{2}\right)  ,\quad\rho^{\prime}=\frac{\rho}{N},
\]

it is clear that the solution to this problem is the solution to the Exact
smoothing problem with smoothing parameter $\rho^{\prime}=\frac{\rho}{N}$.
This eliminates the very unpleasant factor $\sqrt{N}$ in the smoother errors
derived below in this chapter! Things look much better.
\end{remark}

\section{Matrix equations for the Exact smoother}

Part 5 of Theorem \ref{Thm_ex_smooth_var} means that the Exact smoother lies
in the finite dimensional space $W_{G,X}$ and this enables matrix equations to
be derived for the coefficients of the \textit{data-translated }basis
functions $G\left(  x-x^{\left(  k\right)  }\right)  $. These matrix equations
can be written in terms of the \textit{basis function matrix }$G_{X,X}=\left(
G\left(  x^{\left(  i\right)  }-x^{\left(  j\right)  }\right)  \right)  $ or
the \textit{reproducing kernel matrix} $R_{X,X}=\left(  R_{x^{\left(
j\right)  }}\left(  x^{\left(  i\right)  }\right)  \right)  $ where
$R_{X,X}=\left(  2\pi\right)  ^{-\frac{d}{2}}G_{X,X}$. The relevant properties
of $G_{X,X}$ and $R_{X,X}$ were noted in the Introduction to this paper.
$R_{y}\left(  x\right)  \neq G\left(  x-y\right)  $ because the Fourier
transform I have selected is $\widehat{f}(\xi)=\left(  2\pi\right)
^{-d/2}\int e^{-ix\xi}f(x)dx$. For numeric algorithms the use of $G_{X,X}$
might be preferred because $G$ is specified directly and for analysis the use
of $R_{X,X}$ may be preferred because of the analytic properties of $R_{x}$.

\begin{theorem}
\label{Thm_Exact_smth_mat_eqn}Fix $y\in\mathbb{C}^{N}$ and let $\varsigma
=(0,y)\in V$. Then by Theorem \ref{Thm_ex_smooth_var} there exists a unique
solution to the Exact smoothing problem. I denote this solution by $s$ and
note that here, and in the sequel, I will sometimes use the more compact
notation $s_{X}=\widetilde{\mathcal{E}}_{X}s$. The solution $s$ has the
following properties:

\begin{enumerate}
\item $s_{X}$ satisfies the matrix equation%
\begin{equation}
\left(  N\rho I+R_{X,X}\right)  s_{X}=R_{X,X}y,\label{7.47}%
\end{equation}

where $R_{X,X}=\left(  2\pi\right)  ^{-\frac{d}{2}}G_{X,X}$.

\item The matrix $N\rho I+R_{X,X}$ is real, symmetric, positive definite over
$\mathbb{C}^{N}$ and regular.

\item $s\left(  x\right)  =\sum\limits_{k=1}^{N}\alpha_{k}R_{x^{\left(
k\right)  }}\left(  x\right)  $ where $\alpha=\left(  \alpha_{k}\right)  $
satisfies%
\begin{equation}
\left(  N\rho I+R_{X,X}\right)  \alpha=y,\label{7.48}%
\end{equation}

and $\alpha=\frac{1}{N\rho}\left(  y-s_{X}\right)  $ when $\rho>0$.

The next part is convenient for numerical purposes.

\item $s\left(  x\right)  =\sum\limits_{k=1}^{N}\beta_{k}G\left(  x-x^{\left(
k\right)  }\right)  $ where $\beta=\left(  \beta_{k}\right)  $ satisfies%
\[
\left(  \left(  2\pi\right)  ^{\frac{d}{2}}N\rho I+G_{X,X}\right)  \beta=y.
\]

\end{enumerate}
\end{theorem}

\begin{proof}
\textbf{Part 1} Part 4 of Theorem \ref{Thm_ex_smooth_var} and part 4 of
Theorem \ref{Thm_ex_L_op_properties} imply
\begin{equation}
\mathcal{L}_{X}^{\ast}\mathcal{L}_{X}s=\rho Is+\frac{1}{N}%
\widetilde{\mathcal{E}}_{X}^{\ast}\widetilde{\mathcal{E}}_{X}s=\frac{1}%
{N}\widetilde{\mathcal{E}}_{X}^{\ast}y.\label{7.50}%
\end{equation}

From Theorem \ref{Thm_eval_op_properties} I have $\widetilde{\mathcal{E}}%
_{X}\widetilde{\mathcal{E}}_{X}^{\ast}=R_{X,X}$ so that applying the operator
$\left(  2\pi\right)  ^{d/2}N\widetilde{\mathcal{E}}_{X}$ to \ref{7.50} gives
\ref{7.47}.\medskip

\textbf{Part 2} From Theorem \ref{Thm_Rx,x_properties}, $R_{X,X}=\left(
2\pi\right)  ^{-\frac{d}{2}}G_{X,X}$ is Hermitian and positive definite over
$\mathbb{C}^{N}$. Since this paper assumes the basis function $G$ is real
valued, $R_{X,X}$ must be symmetric and thus $N\rho I+R_{X,X}$ is symmetric.
Now suppose $\alpha\in\mathbb{C}^{N}$. Then
\[
\alpha^{T}\left(  N\rho I+R_{X,X}\right)  \overline{\alpha}=N\rho\left\vert
\alpha\right\vert ^{2}+\alpha^{T}R_{X,X}\overline{\alpha},
\]

and since $R_{X,X}$ is positive definite over $\mathbb{C}^{N}$, $\alpha
^{T}\left(  N\rho I+R_{X,X}\right)  \overline{\alpha}=0$ implies $\alpha=0$.
Thus $N\rho I+R_{X,X}$ is positive definite over $\mathbb{C}^{N}$ and
consequently must be regular.\medskip

\textbf{Part 3} From part 5 of Theorem \ref{Thm_ex_smooth_var}, $s=\sum
\limits_{k=1}^{N}\alpha_{k}R_{x^{\left(  k\right)  }}$ where $\alpha_{k}%
=\frac{y_{k}-s(x^{(k)})}{N\rho}$ i.e. $\alpha=\left(  \alpha_{k}\right)
=\frac{y-s_{X}}{N\rho}$. Thus $s_{X}=R_{X,X}\alpha$ and substituting for
$s_{X}$ in \ref{7.47} we get

$\left(  N\rho I+R_{X,X}\right)  R_{X,X}\alpha=R_{X,X}y$ i.e. matrix equation
\ref{7.48} since $R_{X,X}$ is regular.\medskip

\textbf{Part 4} A straightforward consequence of part 3.
\end{proof}

The next corollary will prove very useful.

\begin{corollary}
\label{Cor_Sexg_mat_eqns}If $g\in X_{w}^{0}$ is a (complex) data function then
the Exact smoother operator $\mathcal{S}_{X}$ satisfies%
\begin{equation}
\mathcal{S}_{X}g=\widetilde{\mathcal{E}}_{X}^{\ast}\left(  N\rho
I+R_{X,X}\right)  ^{-1}\widetilde{\mathcal{E}}_{X}g,\label{7.27}%
\end{equation}

and%
\begin{equation}
\left(  \mathcal{S}_{X}g\right)  \left(  x\right)  =\left(
\widetilde{\mathcal{E}}_{X}R_{x}\right)  ^{T}\left(  N\rho I+R_{X,X}\right)
^{-1}\widetilde{\mathcal{E}}_{X}g.\label{7.28}%
\end{equation}

\end{corollary}

\begin{proof}
The data vector is $\widetilde{\mathcal{E}}_{X}g$ and from part 3 of the last
theorem $\left(  \mathcal{S}_{X}g\right)  \left(  x\right)  =\sum
\limits_{k=1}^{N}\alpha_{k}R_{x^{\left(  k\right)  }}\left(  x\right)  $ where
$\alpha=\left(  \alpha_{k}\right)  $ is given by $\alpha=\left(  N\rho
I+R_{X,X}\right)  ^{-1}\widetilde{\mathcal{E}}_{X}g$. Now making use of the
results concerning $\widetilde{\mathcal{E}}_{X}$ and $\widetilde{\mathcal{E}%
}_{X}^{\ast}$ in parts 1 and 2 of Theorem \ref{Thm_eval_op_properties} we have%
\[
s=\mathcal{S}_{X}g=\sum\limits_{k=1}^{N}\alpha_{k}R_{x^{\left(  k\right)  }%
}=\widetilde{\mathcal{E}}_{X}^{\ast}\alpha=\widetilde{\mathcal{E}}_{X}^{\ast
}\left(  N\rho I+R_{X,X}\right)  ^{-1}\widetilde{\mathcal{E}}_{X}g,
\]

and so, since the \textbf{basis function is assumed to be real valued},
\begin{align*}
\mathcal{S}_{X}g\left(  x\right)  =\left(  \widetilde{\mathcal{E}}_{X}^{\ast
}\alpha,R_{x}\right)  _{w,0}=\left(  \alpha,\widetilde{\mathcal{E}}_{X}%
R_{x}\right)  _{\mathbb{C}^{N}} &  =\left(  \left(  N\rho I+R_{X,X}\right)
^{-1}\widetilde{\mathcal{E}}_{X}g,\widetilde{\mathcal{E}}_{X}R_{x}\right)
_{\mathbb{C}^{N}}\\
&  =\left(  \widetilde{\mathcal{E}}_{X}g\right)  ^{T}\left(  N\rho
I+R_{X,X}\right)  ^{-1}\widetilde{\mathcal{E}}_{X}R_{x}\\
&  =\left(  \widetilde{\mathcal{E}}_{X}R_{x}\right)  ^{T}\left(  N\rho
I+R_{X,X}\right)  ^{-1}\widetilde{\mathcal{E}}_{X}g.
\end{align*}

\end{proof}

\section{More properties of the Exact smoother}

\begin{theorem}
\label{Thm_ex_Min_Smth_in_Wgx}By Theorem \ref{Thm_ex_smooth_var} there exists
a unique solution $s$ to the Exact smoothing problem. This solution has the
following properties:

\begin{enumerate}
\item For all $f\in X_{w}^{0}$%
\begin{align}
\rho\left\Vert s\right\Vert _{w,0}^{2}+\frac{1}{N}\sum_{k=1}^{N}\left\vert
s\left(  x^{(k)}\right)  -y_{k}\right\vert ^{2}+\rho\left\Vert s-f\right\Vert
_{w,0}^{2} &  +\frac{1}{N}\sum_{k=1}^{N}\left\vert s\left(  x^{(k)}\right)
-f\left(  x^{(k)}\right)  \right\vert ^{2}\nonumber\\
&  =\rho\left\Vert f\right\Vert _{w,0}^{2}+\frac{1}{N}\sum_{k=1}^{N}\left\vert
f\left(  x^{(k)}\right)  -y_{k}\right\vert ^{2}.\label{7.58}%
\end{align}

\item If $f$ is the data function then $\left\Vert s\right\Vert _{w,0}%
\leq\left\Vert f\right\Vert _{w,0}$ and $\left\Vert s-f\right\Vert _{w,0}%
\leq\left\Vert f\right\Vert _{w,0}$.

\item $\left\Vert s\right\Vert _{w,0}^{2}=\frac{1}{\rho N}\operatorname{Re}%
\sum\limits_{k=1}^{N}\overline{s\left(  x^{(k)}\right)  }\left(
y_{k}-s\left(  x^{(k)}\right)  \right)  $.

\item $J_{e}\left[  s\right]  =\frac{1}{N}\operatorname{Re}\sum\limits_{k=1}%
^{N}\left(  y_{k}-s\left(  x^{(k)}\right)  \right)  \overline{y_{k}}$.
\end{enumerate}
\end{theorem}

\begin{proof}
\textbf{Part 1} From part 3 of Theorem \ref{Thm_ex_smooth_var}
\[
\left\Vert \mathcal{L}_{X}s-\varsigma\right\Vert _{V}^{2}+\left\Vert
\mathcal{L}_{X}\left(  s-f\right)  \right\Vert _{V}^{2}=\left\Vert
\mathcal{L}_{X}f-\varsigma\right\Vert _{V}^{2},
\]

where $\varsigma=\left(  0,y\right)  $. Part 1 follows from Remark
\ref{Rem_SmoothFunc1} and the definition of the norm on the space $V$.\medskip

\textbf{Part 2} If $f$ is the data function then $f\left(  x^{(k)}\right)
=y_{k}$ and part 1 implies part 2.\medskip

\textbf{Part 3} Substituting $f=0$ in equation \ref{7.58} of part 1 yields%
\begin{align*}
2\rho\left\Vert s\right\Vert _{w,0}^{2}  & =\frac{1}{N}\sum_{k=1}%
^{N}\left\vert y_{k}\right\vert ^{2}-\frac{1}{N}\sum_{k=1}^{N}\left\vert
s\left(  x^{(k)}\right)  -y_{k}\right\vert ^{2}-\frac{1}{N}\sum_{k=1}%
^{N}\left\vert s\left(  x^{(k)}\right)  \right\vert ^{2}\\
& =\frac{1}{N}\sum_{k=1}^{N}\left(  \left\vert y_{k}\right\vert ^{2}%
-\left\vert s\left(  x^{(k)}\right)  -y_{k}\right\vert ^{2}-\left\vert
s\left(  x^{(k)}\right)  \right\vert ^{2}\right) \\
& =\frac{1}{N}\sum_{k=1}^{N}\left(  -2\overline{s\left(  x^{(k)}\right)
}s\left(  x^{(k)}\right)  +2\operatorname{Re}\left(  \overline{s\left(
x^{(k)}\right)  }y_{k}\right)  \right) \\
& =\frac{2}{N}\operatorname{Re}\sum_{k=1}^{N}\overline{s\left(  x^{(k)}%
\right)  }\left(  y_{k}-s\left(  x^{(k)}\right)  \right)  .
\end{align*}
\medskip

\textbf{Part 4} By definition, $J_{e}\left[  s\right]  =\rho\left\Vert
s\right\Vert _{w,0}^{2}+\frac{1}{N}\sum\limits_{k=1}^{N}\left\vert s\left(
x^{(k)}\right)  -y_{k}\right\vert ^{2}$.

Substituting for $\left\Vert s\right\Vert _{w,0}^{2}$ using the result of part
2 I get%
\begin{align*}
J_{e}\left[  s\right]   & =\rho\left\Vert s\right\Vert _{w,0}^{2}+\frac{1}%
{N}\sum\limits_{k=1}^{N}\left\vert s\left(  x^{(k)}\right)  -y_{k}\right\vert
^{2}\\
& =\frac{1}{N}\operatorname{Re}\sum_{k=1}^{N}\overline{s\left(  x^{(k)}%
\right)  }\left(  y_{k}-s\left(  x^{(k)}\right)  \right)  +\frac{1}{N}%
\sum\limits_{k=1}^{N}\left\vert s\left(  x^{(k)}\right)  -y_{k}\right\vert
^{2}\\
& =\frac{1}{N}\sum_{k=1}^{N}\left(
\begin{array}
[c]{c}%
\operatorname{Re}\left(  \overline{s\left(  x^{(k)}\right)  }y_{k}%
-\overline{s\left(  x^{(k)}\right)  }s\left(  x^{(k)}\right)  \right)  +\\
+\overline{s\left(  x^{(k)}\right)  }s\left(  x^{(k)}\right)
-2\operatorname{Re}\overline{s\left(  x^{(k)}\right)  }y_{k}+\overline{y_{k}%
}y_{k}%
\end{array}
\right) \\
& =\frac{1}{N}\operatorname{Re}\sum_{k=1}^{N}\left(  -\overline{s\left(
x^{(k)}\right)  }y_{k}+\overline{y_{k}}y_{k}\right) \\
& =\frac{1}{N}\operatorname{Re}\sum_{k=1}^{N}\left(  \overline{y_{k}%
}-\overline{s\left(  x^{(k)}\right)  }\right)  y_{k}\\
& =\frac{1}{N}\operatorname{Re}\sum_{k=1}^{N}\left(  y_{k}-s\left(
x^{(k)}\right)  \right)  \overline{y_{k}}.
\end{align*}

\end{proof}

\begin{remark}
\label{Rem_Thm_Min_Smth_in_Wgx}Part 1 of the last theorem confirms that if $s
$ is the Exact smoother then $J_{e}\left[  s\right]  \leq J_{e}\left[
f\right]  $ for all $f\in X_{w}^{0}$ i.e. that $s$ minimizes the functional
$J_{e}$.
\end{remark}

\section{Pointwise error estimates using smoother inner products and
semi-inner products\label{Sect_ExactSmthErrEstim_InnerProd}}

Following the example of the interpolation semi-inner product $\left(
f-\mathcal{I}_{X}f,g\right)  _{w,0}$ of Definition \ref{Def_interp_seminorm}
of \ref{Ch_Interpol} and it's sequel I will define the semi-inner product
$\left(  \mathcal{S}_{X}f,g\right)  _{w,0}$ and inner product $\left(
f-\mathcal{S}_{X}f,g\right)  _{w,0}$ and use these to study the value and the
error of the Exact smoother respectively. These expressions are useful because
$\left(  \mathcal{S}_{X}f,R_{x}\right)  _{w,0}=\left(  \mathcal{S}%
_{X}f\right)  \left(  x\right)  $ and $\left(  f-\mathcal{S}_{X}%
f,R_{x}\right)  _{w,0}=\left(  f-\mathcal{S}_{X}f\right)  \left(  x\right)  $
and because of the power of Hilbert space theory. An outcome of this approach
is the Exact smoother bounds \ref{7.68} and the convergence estimates of
\ Corollary \ref{Cor_Thm_Rx(x)minus(SRx)(x)_basis_fn_bound}.

To define these smoother norms I will need the lemma:

\begin{lemma}
\label{Lem_ex_smth_err_1}Suppose $f\in X_{w}^{0}$ and $\mathcal{S}_{X}$ is the
Exact smoother mapping where $X=\left\{  x^{\left(  k\right)  }\right\}
_{k=1}^{N}$. Then:

\begin{enumerate}
\item $\left\Vert \mathcal{S}_{X}f\right\Vert _{w,0}^{2}=\left(
\mathcal{S}_{X}f,f\right)  _{w,0}-\rho N\left\vert \left(  \rho NI+R_{X,X}%
\right)  ^{-1}\widetilde{\mathcal{E}}_{X}f\right\vert ^{2}$.

\item $\left\Vert f-\mathcal{S}_{X}f\right\Vert _{w,0}^{2}=\left(
f-\mathcal{S}_{X}f,f\right)  _{w,0}-\rho N\left\vert \left(  \rho
NI+R_{X,X}\right)  ^{-1}\widetilde{\mathcal{E}}_{X}f\right\vert ^{2}$.
\end{enumerate}
\end{lemma}

\begin{proof}
\textbf{Part 1} I start with formula \ref{7.27} for $\mathcal{S}_{X}f$.%
\begin{align}
\left\Vert \mathcal{S}_{X}f\right\Vert _{w,0}^{2}  & =\left(
\widetilde{\mathcal{E}}_{X}^{\ast}\left(  \rho NI+R_{X,X}\right)
^{-1}\widetilde{\mathcal{E}}_{X}f,\widetilde{\mathcal{E}}_{X}^{\ast}\left(
\rho NI+R_{X,X}\right)  ^{-1}\widetilde{\mathcal{E}}_{X}f\right)
_{w,0}\nonumber\\
& =\left(  \widetilde{\mathcal{E}}_{X}\widetilde{\mathcal{E}}_{X}^{\ast
}\left(  \rho NI+R_{X,X}\right)  ^{-1}\widetilde{\mathcal{E}}_{X}f,\left(
\rho NI+R_{X,X}\right)  ^{-1}\widetilde{\mathcal{E}}_{X}f\right)
_{\mathbb{C}^{N}}\nonumber\\
& =\left(  R_{X,X}\left(  \rho NI+R_{X,X}\right)  ^{-1}\widetilde{\mathcal{E}%
}_{X}f,\left(  \rho NI+R_{X,X}\right)  ^{-1}\widetilde{\mathcal{E}}%
_{X}f\right)  _{\mathbb{C}^{N}},\nonumber
\end{align}

since $\widetilde{\mathcal{E}}_{X}\widetilde{\mathcal{E}}_{X}^{\ast}=R_{X,X} $
but
\begin{align*}
R_{X,X}\left(  \rho NI+R_{X,X}\right)  ^{-1}  & =\left(  \rho NI+R_{X,X}%
\right)  \left(  \rho NI+R_{X,X}\right)  ^{-1}-\rho N\left(  \rho
NI+R_{X,X}\right)  ^{-1}\\
& =I-\rho N\left(  \rho NI+R_{X,X}\right)  ^{-1},
\end{align*}

so that%
\begin{align*}
\left\Vert \mathcal{S}_{X}f\right\Vert _{w,0}^{2}  & =\left(
\widetilde{\mathcal{E}}_{X}f,\left(  \rho NI+R_{X,X}\right)  ^{-1}%
\widetilde{\mathcal{E}}_{X}f\right)  -\rho N\left(  \left(  \rho
NI+R_{X,X}\right)  ^{-1}\widetilde{\mathcal{E}}_{X}f,\left(  \rho
NI+R_{X,X}\right)  ^{-1}\widetilde{\mathcal{E}}_{X}f\right) \\
& =\left(  \widetilde{\mathcal{E}}_{X}^{\ast}\left(  \rho NI+R_{X,X}\right)
^{-1}\widetilde{\mathcal{E}}_{X}f,f\right)  _{w,0}-\rho N\left\vert \left(
\rho NI+R_{X,X}\right)  ^{-1}\widetilde{\mathcal{E}}_{X}f\right\vert ^{2}\\
& =\left(  \mathcal{S}_{X}f,f\right)  _{w,0}-\rho N\left\vert \left(  \rho
NI+R_{X,X}\right)  ^{-1}\widetilde{\mathcal{E}}_{X}f\right\vert ^{2}.
\end{align*}
\medskip

\textbf{Part 2} Substituting the equation for $\left\Vert \mathcal{S}%
_{X}f\right\Vert _{w,0}^{2}$ proved in part 1 yields
\begin{align*}
\left\Vert f-\mathcal{S}_{X}f\right\Vert _{w,0}^{2}  & =\left\Vert
f\right\Vert _{w,0}^{2}-2\left(  \mathcal{S}_{X}f,f\right)  _{w,0}+\left\Vert
\mathcal{S}_{X}f\right\Vert _{w,0}^{2}\\
& =\left\Vert f\right\Vert _{w,0}^{2}-2\left(  \mathcal{S}_{X}f,f\right)
_{w,0}+\left(  \mathcal{S}_{X}f,f\right)  _{w,0}-\rho N\left\vert \left(  \rho
NI+R_{X,X}\right)  ^{-1}\widetilde{\mathcal{E}}_{X}f\right\vert ^{2}\\
& =\left\Vert f\right\Vert _{w,0}^{2}-\left(  \mathcal{S}_{X}f,f\right)
_{w,0}-\rho N\left\vert \left(  \rho NI+R_{X,X}\right)  ^{-1}%
\widetilde{\mathcal{E}}_{X}f\right\vert ^{2}\\
& =\left(  f-\mathcal{S}_{X}f,f\right)  _{w,0}-\rho N\left\vert \left(  \rho
NI+R_{X,X}\right)  ^{-1}\widetilde{\mathcal{E}}_{X}f\right\vert ^{2}.
\end{align*}

\end{proof}

\begin{theorem}
\label{Thm_(Sf,g)_(f-Sf,g)_semi_iiner_prod}The expression $\left(
f-\mathcal{S}_{X}f,g\right)  _{w,0}$ is an inner product on $X_{w}^{0}$ and
$\left(  \mathcal{S}_{X}f,g\right)  _{w,0}$ is a semi-norm with null space
$W_{G,X}^{\perp}$.
\end{theorem}

\begin{proof}
From part 1 of Lemma \ref{Lem_ex_smth_err_1}, $\left(  \mathcal{S}%
_{X}f,f\right)  _{w,0}\geq0$ and $\left(  \mathcal{S}_{X}f,f\right)  _{w,0}=0$

iff $\left(  \rho NI+R_{X,X}\right)  ^{-1}\widetilde{\mathcal{E}}_{X}f=0$ iff
$\widetilde{\mathcal{E}}_{X}f$ iff $f\in W_{G,X}^{\perp}$. Thus $\left(
\mathcal{S}_{X}f,f\right)  _{w,0}$ is a seminorm. Set $\left\vert f\right\vert
^{2}=\left(  \mathcal{S}_{X}f,f\right)  _{w,0}$ so that by part 1 of Lemma
\ref{Lem_ex_smth_err_1}%
\begin{align*}
\left\vert f+g\right\vert ^{2}+\left\vert f-g\right\vert ^{2}  & =\left(
\mathcal{S}_{X}\left(  f+g\right)  ,f+g\right)  _{w,0}+\left(  \mathcal{S}%
_{X}\left(  f-g\right)  ,f-g\right)  _{w,0}\\
& =\left\Vert \mathcal{S}_{X}\left(  f+g\right)  \right\Vert _{w,0}^{2}+\rho
N\left\vert \left(  \rho NI+R_{X,X}\right)  ^{-1}\widetilde{\mathcal{E}}%
_{X}\left(  f+g\right)  \right\vert ^{2}+\\
& +\left\Vert \mathcal{S}_{X}\left(  f-g\right)  \right\Vert _{w,0}^{2}+\rho
N\left\vert \left(  \rho NI+R_{X,X}\right)  ^{-1}\widetilde{\mathcal{E}}%
_{X}\left(  f-g\right)  \right\vert ^{2}\\
& =2\left(  \left\Vert \mathcal{S}_{X}f\right\Vert _{w,0}^{2}+\left\Vert
\mathcal{S}_{X}g\right\Vert _{w,0}^{2}\right)  +\\
& +2\rho N\left(  \left\vert \left(  \rho NI+R_{X,X}\right)  ^{-1}%
\widetilde{\mathcal{E}}_{X}f\right\vert ^{2}+\left\vert \left(  \rho
NI+R_{X,X}\right)  ^{-1}\widetilde{\mathcal{E}}_{X}g\right\vert ^{2}\right) \\
& =2\left(  \mathcal{S}_{X}f,f\right)  _{w,0}+2\left(  \mathcal{S}%
_{X}g,g\right)  _{w,0}\\
& =2\left(  \left\vert f\right\vert ^{2}+\left\vert g\right\vert ^{2}\right)
,
\end{align*}

so that $\left\vert \cdot\right\vert $ satisfies the parallelogram law. Thus a
seminorm can be defined by%
\begin{align*}
\left\langle f,g\right\rangle  & =\frac{1}{4}\left(  \left\vert f+g\right\vert
^{2}-\left\vert f-g\right\vert ^{2}\right)  +\frac{i}{4}\left(  \left\vert
f+ig\right\vert ^{2}-\left\vert f-ig\right\vert ^{2}\right) \\
& =\frac{1}{4}\left(  \left(  \mathcal{S}_{X}\left(  f+g\right)  ,f+g\right)
_{w,0}-\left(  \mathcal{S}_{X}\left(  f-g\right)  ,f-g\right)  _{w,0}\right)
+\\
& \qquad+\frac{i}{4}\left(  \left(  \mathcal{S}_{X}\left(  f+ig\right)
,if+g\right)  _{w,0}-\left(  \mathcal{S}_{X}\left(  f-ig\right)  ,f-ig\right)
_{w,0}\right) \\
& =\frac{1}{2}\left(  \left(  \mathcal{S}_{X}f,g\right)  _{w,0}+\left(
\mathcal{S}_{X}g,f\right)  _{w,0}\right)  +\frac{i}{2}\left(  \left(
\mathcal{S}_{X}f,ig\right)  _{w,0}+\left(  \mathcal{S}_{X}ig,f\right)
_{w,0}\right) \\
& =\frac{1}{2}\left(  \left(  \mathcal{S}_{X}f,g\right)  _{w,0}+\left(
\mathcal{S}_{X}g,f\right)  _{w,0}\right)  +\frac{i}{2}\left(  -i\left(
\mathcal{S}_{X}f,g\right)  _{w,0}+i\left(  \mathcal{S}_{X}g,f\right)
_{w,0}\right) \\
& =\left(  \mathcal{S}_{X}f,g\right)  _{w,0}.
\end{align*}
\smallskip

To prove $\left(  f-\mathcal{S}_{X}f,g\right)  _{w,0}$ is an inner product we
use part 2 of Lemma \ref{Lem_ex_smth_err_1}:%
\[
\left\Vert f-\mathcal{S}_{X}f\right\Vert _{w,0}^{2}=\left(  f-\mathcal{S}%
_{X}f,f\right)  _{w,0}-\rho N\left\vert \left(  \rho NI+R_{X,X}\right)
^{-1}\widetilde{\mathcal{E}}_{X}f\right\vert ^{2},
\]

and \ref{7.53}:%
\[
\mathcal{S}_{X}g=\frac{1}{N}\left(  \mathcal{L}_{X}^{\ast}\mathcal{L}%
_{X}\right)  ^{-1}\widetilde{\mathcal{E}}_{X}^{\ast}\widetilde{\mathcal{E}%
}_{X}g,\quad g\in X_{w}^{0}.
\]

In fact, $\left(  f-\mathcal{S}_{X}f,f\right)  _{w,0}=0$ implies
$\widetilde{\mathcal{E}}_{X}f=0$ and $f-\mathcal{S}_{X}f=0$ which implies
$f=\mathcal{S}_{X}f=0$.
\end{proof}

I can now define our Exact smoother seminorms on\textbf{\ }$X_{w}^{0}$.

\begin{definition}
\label{Def_Exact_smth_seminorms}\textbf{Exact smoother value seminorm} and
\textbf{Exact smoother error norm }on $X_{w}^{0}$.

Define the Exact smoother \textbf{value} seminorm and semi-inner product by
\[
\left\vert f\right\vert _{sv}^{2}=\left(  \mathcal{S}_{X}f,f\right)
_{w,0},\qquad\left\langle f,g\right\rangle _{sv}=\left(  \mathcal{S}%
_{X}f,g\right)  _{w,0}.
\]

The qualifier `value' is used because $\left\langle f,R_{x}\right\rangle
_{sv}=\left(  \mathcal{S}_{X}f\right)  \left(  x\right)  .$

Denote the Exact smoother \textbf{error} norm and inner product by
\[
\left\Vert f\right\Vert _{se}^{2}=\left(  f-\mathcal{S}_{X}f,f\right)
_{w,0},\qquad\left(  f,g\right)  _{se}=\left(  f-\mathcal{S}_{X}f,g\right)
_{w,0}.
\]

The qualifier `error' is used because $\left(  f,R_{x}\right)  _{se}=\left(
f-\mathcal{S}_{X}f\right)  \left(  x\right)  $.
\end{definition}

\begin{remark}
\label{Rem_Def_Exact_smth_seminorms}\ 

\begin{enumerate}
\item Since $\mathcal{I}_{X}\mathcal{S}_{X}=\mathcal{S}_{X}\mathcal{I}%
_{X}=\mathcal{S}_{X}$ and $\left(  \mathcal{I}_{X}f,g\right)  _{w,0}=\left(
f,\mathcal{I}_{X}g\right)  _{w,0}$ we have:%
\begin{align*}
\left(  f,g\right)  _{se}=\left(  f-\mathcal{S}_{X}f,g\right)  _{w,0}=\left(
\left(  I-\mathcal{S}_{X}\right)  f,g\right)  _{w,0}  & =\left(
\mathcal{I}_{X}\left(  I-\mathcal{S}_{X}\right)  \mathcal{I}_{X}f,g\right)
_{w,0}\\
& =\left(  \left(  I-\mathcal{S}_{X}\right)  \mathcal{I}_{X}f,\mathcal{I}%
_{X}g\right)  _{w,0}\\
& =\left(  \mathcal{I}_{X}f,\mathcal{I}_{X}g\right)  _{se},
\end{align*}

and%
\begin{align*}
\left\langle f,g\right\rangle _{sv}=\left(  \mathcal{S}_{X}f,g\right)
_{w,0}=\left(  \mathcal{I}_{X}\mathcal{S}_{X}\mathcal{I}_{X}f,g\right)
_{w,0}  & =\left(  \mathcal{S}_{X}\mathcal{I}_{X}f,\mathcal{I}_{X}g\right)
_{w,0}\\
& =\left\langle \mathcal{I}_{X}f,\mathcal{I}_{X}g\right\rangle _{sv}.
\end{align*}

\item ?? $\left(  f,g\right)  _{sv}:=\left\langle f,g\right\rangle
_{sv}+c\left(  \widetilde{\mathcal{E}}_{X}f,\widetilde{\mathcal{E}}%
_{X}g\right)  $ is a norm when $c=\ldots$ ??.
\end{enumerate}
\end{remark}

\begin{remark}
Since $\mathcal{I}_{X}\mathcal{S}_{X}=\mathcal{S}_{X}\mathcal{I}_{X}$ and
$\mathcal{I}_{X}\mathcal{I}_{X}=\mathcal{I}_{X}$, we have $\left(
\mathcal{I}_{X}f,g\right)  _{se}=\left(  f,\mathcal{I}_{X}g\right)  _{se}$.
Also%
\begin{align*}
\left(  f,g\right)  _{se}  & =\left(  f-\mathcal{S}_{X}f,g\right)  _{w,0}\\
& =\left(  f-\mathcal{I}_{X}f,g\right)  _{w,0}+\left(  \mathcal{I}%
_{X}f-\mathcal{S}_{X}f,g\right)  _{w,0}\\
& =\left(  f-\mathcal{I}_{X}f,g\right)  _{w,0}+\left(  \mathcal{I}%
_{X}f-\mathcal{S}_{X}\mathcal{I}_{X}f,g\right)  _{w,0}\\
& =\left(  f-\mathcal{I}_{X}f,g\right)  _{w,0}+\left(  \mathcal{I}%
_{X}\mathcal{I}_{X}f-\mathcal{I}_{X}\mathcal{S}_{X}\mathcal{I}_{X}f,g\right)
_{w,0}\\
& =\left(  f-\mathcal{I}_{X}f,g\right)  _{w,0}+\left(  \mathcal{I}%
_{X}\mathcal{I}_{X}f-\mathcal{I}_{X}\mathcal{S}_{X}\mathcal{I}_{X}f,g\right)
_{w,0}\\
& =\left(  f-\mathcal{I}_{X}f,g\right)  _{w,0}+\left(  \mathcal{I}%
_{X}f-\mathcal{S}_{X}\mathcal{I}_{X}f,\mathcal{I}_{X}g\right)  _{w,0}\\
& =\left(  f-\mathcal{I}_{X}f,g\right)  _{w,0}+\left(  \mathcal{I}%
_{X}f,\mathcal{I}_{X}g\right)  _{se}\\
& =\left(  f,g\right)  _{se;\rho=0}+\left(  \mathcal{I}_{X}f,\mathcal{I}%
_{X}g\right)  _{se}\\
& change\text{ }f\text{ }away\text{ }from\text{ }data\\
& =\left(  f,g\right)  _{se;\rho=0}+\left(  \mathcal{I}_{X}f_{1}%
,\mathcal{I}_{X}g\right)  _{se}\\
& =\left(  f,g\right)  _{se;\rho=0}+\left(  f_{1},\mathcal{I}_{X}g\right)
_{se}%
\end{align*}

\end{remark}

In the next theorem the pointwise error problem for an arbitrary function is
reduced to that of considering the error of the smoother of $R_{x}$ at $x$,
where $R_{x}$ is the Riesz representer of the evaluation functional
$f\rightarrow f\left(  x\right)  $. This error will latter be estimated using
the values of the basis function near the origin. This theorem also gives an
error bound that is uniform on $\mathbb{R}^{d}$.

\begin{theorem}
\label{Thm_smth_err_bound}Since $\mathcal{S}_{X}f_{d}$ is the Exact smoother
of the data function $f_{d}\in X_{w}^{0}$ and $\mathcal{S}_{X}R_{x}$ is the
Exact smoother of the Riesz representer data function $R_{x}$, we have the
bound%
\begin{equation}
\left\vert f_{d}\left(  x\right)  -\left(  \mathcal{S}_{X}f_{d}\right)
\left(  x\right)  \right\vert \leq\sqrt{\left(  f_{d}-\mathcal{S}_{X}%
f_{d},f_{d}\right)  _{w,0}}\sqrt{R_{x}\left(  x\right)  -\left(
\mathcal{S}_{X}R_{x}\right)  \left(  x\right)  },\quad x\in\mathbb{R}%
^{d},\label{7.19}%
\end{equation}

where $\sqrt{\left(  f_{d}-\mathcal{S}_{X}f_{d},f_{d}\right)  _{w,0}}%
\leq\left\Vert f_{d}\right\Vert _{w,0}$, as well as the bound%
\begin{equation}
\left\vert f_{d}\left(  x\right)  -\left(  \mathcal{S}_{X}f_{d}\right)
\left(  x\right)  \right\vert \leq\sqrt{\left(  f_{d}-\mathcal{S}_{X}%
f_{d},f_{d}\right)  _{w,0}}\sqrt{R_{0}\left(  0\right)  },\quad x\in
\mathbb{R}^{d}.\label{7.26}%
\end{equation}

\end{theorem}

\begin{proof}
Estimate \ref{7.19} is just the Cauchy-Schwartz inequality $\left\vert \left(
f_{d},R_{x}\right)  _{se}\right\vert \leq\left\Vert f_{d}\right\Vert
_{se}\left\Vert R_{x}\right\Vert _{se}$.

The second inequality of this theorem follows directly from the inequality

$\left\Vert f_{d}-\mathcal{S}_{X}f_{d}\right\Vert _{w,0}\leq\left\Vert
f_{d}\right\Vert _{w,0}$, proven in part 2 of Theorem
\ref{Thm_ex_Min_Smth_in_Wgx}. Thus%
\begin{align*}
\left\vert f_{d}\left(  x\right)  -\left(  \mathcal{S}_{X}f_{d}\right)
\left(  x\right)  \right\vert  & \leq\sqrt{\left(  f_{d}-\mathcal{S}_{X}%
f_{d},f_{d}\right)  _{w,0}}\sqrt{R_{x}\left(  x\right)  -\left(
\mathcal{S}_{X}R_{x}\right)  \left(  x\right)  }\\
& \leq\sqrt{\left(  f_{d}-\mathcal{S}_{X}f_{d},f_{d}\right)  _{w,0}}%
\sqrt{R_{x}\left(  x\right)  }\\
& =\sqrt{\left(  f_{d}-\mathcal{S}_{X}f_{d},f_{d}\right)  _{w,0}}\sqrt
{R_{0}\left(  0\right)  }.
\end{align*}

\end{proof}

The last result means that convergence results about arbitrary data functions
in $X_{w}^{0}$ can be proved by deriving error results for the Riesz
representer function $R_{x}$ at $x$.

I will now use the smoothing value seminorm to put bounds on the smoother that
show how the smoother behaves for large values of the smoothing parameter
$\rho$: it decreases to zero at at least the rate $1/\rho$.

\begin{theorem}
\label{Thm_bound_Exact_smth}If $f\in X_{w}^{0}$ and $s_{e}$ is the Exact
smoother of $f$ with smoothing coefficient $\rho$ then%
\begin{align}
\left\vert s_{e}\left(  x\right)  \right\vert  & \leq\left\{
\begin{array}
[c]{ll}%
\left\Vert f\right\Vert _{w,0}\sqrt{R_{0}\left(  0\right)  }, & \rho\leq
R_{0}\left(  0\right)  ,\\
\left\Vert f\right\Vert _{w,0}\left(  \sqrt{R_{0}\left(  0\right)  }\right)
^{3}\rho^{-1},\text{ } & \rho\geq R_{0}\left(  0\right)  ,
\end{array}
\right. \nonumber\\
& =\left\Vert f\right\Vert _{w,0}\sqrt{R_{0}\left(  0\right)  }\min\left\{
1,\frac{R_{0}\left(  0\right)  }{\rho}\right\}  ,\label{7.68}%
\end{align}

where $R_{0}\left(  0\right)  =\left(  2\pi\right)  ^{-\frac{d}{2}}G\left(
0\right)  =\left\Vert R_{0}\right\Vert _{w,0}^{2}$.
\end{theorem}

\begin{proof}
If $\mathcal{S}_{X}$ is the Exact smoother mapping then%
\begin{equation}
\left\vert \mathcal{S}_{X}f\left(  x\right)  \right\vert =\left\vert \left(
\mathcal{S}_{X}f,R_{x}\right)  \right\vert \leq\left\Vert \mathcal{S}%
_{X}f\right\Vert _{w,0}\left\Vert R_{x}\right\Vert _{w,0}\leq\left\Vert
f\right\Vert _{w,0}\left\Vert R_{x}\right\Vert _{w,0}=\left\Vert f\right\Vert
_{w,0}\sqrt{R_{0}\left(  0\right)  }.\label{7.32}%
\end{equation}

From part 1 of Lemma \ref{Lem_ex_smth_err_1}%
\[
\left\vert \left(  \rho NI+R_{X,X}\right)  ^{-1}\widetilde{\mathcal{E}}%
_{X}f\right\vert <\frac{1}{\sqrt{\rho N}}\sqrt{\left(  \mathcal{S}%
_{X}f,f\right)  _{w,0}},\quad f\in X_{w}^{0},
\]

and from part 3 of Theorem \ref{Thm_eval_op_properties}, $\left\vert
\widetilde{\mathcal{E}}_{X}f\right\vert \leq\left\Vert f\right\Vert
_{w,0}\sqrt{N}\sqrt{R_{0}\left(  0\right)  }$. Hence
\begin{align*}
\left(  \mathcal{S}_{X}f,f\right)  _{w,0}=\left(  \widetilde{\mathcal{E}}%
_{X}f\right)  ^{T}\left(  \rho NI+R_{X,X}\right)  ^{-1}\widetilde{\mathcal{E}%
}_{X}f &  \leq\left\vert \widetilde{\mathcal{E}}_{X}f\right\vert \left\vert
\left(  \rho NI+R_{X,X}\right)  ^{-1}\widetilde{\mathcal{E}}_{X}f\right\vert
\\
&  \leq\left\Vert f\right\Vert _{w,0}\left(  \sqrt{N}\sqrt{R_{0}\left(
0\right)  }\right)  \frac{1}{\sqrt{\rho N}}\sqrt{\left(  \mathcal{S}%
_{X}f,f\right)  _{w,0}}\\
&  =\left\Vert f\right\Vert _{w,0}\sqrt{R_{0}\left(  0\right)  }\frac{1}%
{\sqrt{\rho}}\sqrt{\left(  \mathcal{S}_{X}f,f\right)  _{w,0}},
\end{align*}

so that $\left(  \mathcal{S}_{X}f,f\right)  _{w,0}\leq\left\Vert f\right\Vert
_{w,0}^{2}R_{0}\left(  0\right)  /\rho$. But by Theorem
\ref{Thm_smth_err_bound} $\left(  \mathcal{S}_{X}f,g\right)  _{w,0}$ is a
semi-inner product and hence by the Cauchy-Schwartz theorem%
\[
\left\vert \left(  \mathcal{S}_{X}f,g\right)  _{w,0}\right\vert \leq
\sqrt{\left(  \mathcal{S}_{X}f,f\right)  _{w,0}}\sqrt{\left(  \mathcal{S}%
_{X}g,g\right)  _{w,0}},\quad f,g\in X_{w}^{0},
\]

which allows us to conclude that
\[
\left\vert \mathcal{S}_{X}f\left(  x\right)  \right\vert =\left\vert \left(
\mathcal{S}_{X}f,R_{x}\right)  _{w,0}\right\vert \leq\sqrt{\left(
\mathcal{S}_{X}f,f\right)  _{w,0}}\sqrt{\left(  \mathcal{S}_{X}R_{x}%
,R_{x}\right)  _{w,0}}\leq\left\Vert f\right\Vert _{w,0}\left(  \sqrt
{R_{0}\left(  0\right)  }\right)  ^{3}\frac{1}{\rho}.
\]

which when combined with inequality \ref{7.32} proves the first estimates of
this theorem.
\end{proof}

\begin{remark}
These bounds on the smoother make more sense if we note that $\left\vert
f\left(  x\right)  \right\vert \leq\left\Vert f\right\Vert _{w,0}\sqrt
{R_{0}\left(  0\right)  }$.
\end{remark}

The goal of the next theorem is to prove the inequality of part 3 of this
theorem. This inequality is proved using the identity of part 2 which is in
turn derived by expanding the two terms $\left\Vert R_{x}-R_{x^{\left(
k\right)  }}\right\Vert _{w,0}^{2}-\left\Vert \mathcal{S}_{X}\left(
R_{x}-R_{x^{\left(  k\right)  }}\right)  \right\Vert _{w,0}^{2}$. This
approach is motivated by the identity of part 2 of Theorem
\ref{Thm_property_interpol}, namely,%
\[
R_{x}\left(  x\right)  -\left(  \mathcal{I}_{X}R_{x}\right)  \left(  x\right)
=\left\Vert R_{x}-R_{x^{\left(  k\right)  }}\right\Vert _{w,0}^{2}-\left\Vert
\mathcal{I}_{X}\left(  R_{x}-R_{x^{\left(  k\right)  }}\right)  \right\Vert
_{w,0}^{2},
\]

which was proved by expanding its right side. Here $\mathcal{I}_{X}$ is the
interpolant mapping.

\begin{theorem}
\label{Thm_|Rx|se_le_|Ry|se_plus_||Rx_minus_Ry||sq_w0}Recall that $\left\Vert
f\right\Vert _{se}^{2}=\left(  f-\mathcal{S}_{X}f,f\right)  _{w,0}$ is the
Exact smoother error norm introduced in Definition
\ref{Def_Exact_smth_seminorms}. Then for $x,y\in\mathbb{R}^{d}$:

\begin{enumerate}
\item
\begin{align*}
\left\Vert \mathcal{S}_{X}\left(  R_{x}-R_{y}\right)  \right\Vert _{w,0}%
^{2}=\left(  \mathcal{S}_{X}R_{x}\right)  \left(  x\right)   &  -2\left(
\mathcal{S}_{X}R_{x}\right)  \left(  y\right)  +\left(  \mathcal{S}_{X}%
R_{y}\right)  \left(  y\right)  -\\
&  -\rho N\left\vert \left(  \rho NI+R_{X,X}\right)  ^{-1}%
\widetilde{\mathcal{E}}_{X}\left(  R_{x}-R_{y}\right)  \right\vert ^{2}.
\end{align*}

\item
\[
\left\Vert R_{x}-R_{y}\right\Vert _{se}^{2}=\left\Vert R_{x}-R_{y}\right\Vert
_{w,0}^{2}-\left\Vert \mathcal{S}_{X}\left(  R_{x}-R_{y}\right)  \right\Vert
_{w,0}^{2}-\rho N\left\vert \left(  \rho NI+R_{X,X}\right)  ^{-1}%
\widetilde{\mathcal{E}}_{X}\left(  R_{x}-R_{y}\right)  \right\vert ^{2}.
\]

\item
\[
\left\vert \left\Vert R_{x}\right\Vert _{se}-\left\Vert R_{y}\right\Vert
_{se}\right\vert \leq\left\Vert R_{x}-R_{y}\right\Vert _{w,0}.
\]

\end{enumerate}
\end{theorem}

\begin{proof}
\textbf{Part 1} Using the equation of part 1 of Lemma \ref{Lem_ex_smth_err_1}
with $f=R_{x}-R_{y}$
\begin{align*}
\left\Vert \mathcal{S}_{X}\left(  R_{x}-R_{y}\right)  \right\Vert _{w,0}^{2}
& =\left(  \mathcal{S}_{X}\left(  R_{x}-R_{y}\right)  ,R_{x}-R_{y}\right)
_{w,0}-\rho N\left\vert \left(  \rho NI+R_{X,X}\right)  ^{-1}%
\widetilde{\mathcal{E}}_{X}\left(  R_{x}-R_{y}\right)  \right\vert ^{2}\\
& =\left(  \mathcal{S}_{X}R_{x},R_{x}\right)  _{w,0}-2\left(  \mathcal{S}%
_{X}R_{x},R_{y}\right)  _{w,0}+\left(  \mathcal{S}_{X}R_{y},R_{y}\right)
_{w,0}-\\
& \qquad-\rho N\left\vert \left(  \rho NI+R_{X,X}\right)  ^{-1}%
\widetilde{\mathcal{E}}_{X}\left(  R_{x}-R_{y}\right)  \right\vert ^{2}\\
& =\left(  \mathcal{S}_{X}R_{x}\right)  \left(  x\right)  -2\left(
\mathcal{S}_{X}R_{x}\right)  \left(  y\right)  +\left(  \mathcal{S}_{X}%
R_{y}\right)  \left(  y\right)  -\\
& \qquad-\rho N\left\vert \left(  \rho NI+R_{X,X}\right)  ^{-1}%
\widetilde{\mathcal{E}}_{X}\left(  R_{x}-R_{y}\right)  \right\vert ^{2}.
\end{align*}
\medskip

\textbf{Part 2} Since $\left\Vert R_{x}-R_{y}\right\Vert _{w,0}^{2}%
=R_{x}\left(  x\right)  -2R_{x}\left(  y\right)  +R_{y}\left(  y\right)  $, by
part 1%
\begin{align*}
\left\Vert R_{x}-R_{y}\right\Vert _{w,0}^{2}-\left\Vert \mathcal{S}_{X}\left(
R_{x}-R_{y}\right)  \right\Vert _{w,0}^{2}  & =\left(  R_{x}\left(  x\right)
-\left(  \mathcal{S}_{X}R_{x}\right)  \left(  x\right)  \right)  -2\left(
R_{x}\left(  y\right)  -\left(  \mathcal{S}_{X}R_{x}\right)  \left(  y\right)
\right)  +\\
& \qquad+\left(  R_{y}\left(  y\right)  -\left(  \mathcal{S}_{X}R_{y}\right)
\left(  y\right)  \right)  +\\
& \qquad+\rho N\left\vert \left(  \rho NI+R_{X,X}\right)  ^{-1}%
\widetilde{\mathcal{E}}_{X}\left(  R_{x}-R_{y}\right)  \right\vert ^{2}\\
& =\left\Vert R_{x}\right\Vert _{se}^{2}+\left\Vert R_{y}\right\Vert _{se}%
^{2}-2\left(  R_{x},R_{y}\right)  _{se}+\\
& \qquad+\rho N\left\vert \left(  \rho NI+R_{X,X}\right)  ^{-1}%
\widetilde{\mathcal{E}}_{X}\left(  R_{x}-R_{y}\right)  \right\vert ^{2}\\
& =\left\Vert R_{x}-R_{y}\right\Vert _{se}^{2}+\rho N\left\vert \left(  \rho
NI+R_{X,X}\right)  ^{-1}\widetilde{\mathcal{E}}_{X}\left(  R_{x}-R_{y}\right)
\right\vert ^{2},
\end{align*}

and rearranging gives the result.\medskip

\textbf{Part 3} From part 2, $\left\Vert R_{x}-R_{y}\right\Vert _{se}%
\leq\left\Vert R_{x}-R_{y}\right\Vert _{w,0}$ and since $\left\Vert
R_{x}-R_{y}\right\Vert _{se}\geq\left\vert \left\Vert R_{x}\right\Vert
_{se}-\left\Vert R_{y}\right\Vert _{se}\right\vert $ we have part 3.
\end{proof}

I now need the following lemma:

\begin{lemma}
\label{Lem_Rxk(xk).minus.(SRxk)(xk)} Suppose the weight function has property
W02 for some $\kappa\geq0$. Suppose $X=\left\{  x^{\left(  k\right)
}\right\}  _{k=1}^{N}$ is the independent data for the Exact smoothing
operator $\mathcal{S}_{X}$. Then%
\begin{equation}
R_{x^{\left(  j\right)  }}\left(  x^{\left(  k\right)  }\right)  -\left(
\mathcal{S}_{X}R_{x^{\left(  j\right)  }}\right)  \left(  x^{\left(  k\right)
}\right)  =\rho N\delta_{j,k}-\left(  \rho N\right)  ^{2}\left(  \left(  \rho
NI+R_{X,X}\right)  ^{-1}\right)  _{j,k},\label{7.35}%
\end{equation}

where $\left(  \left(  \rho NI+R_{X,X}\right)  ^{-1}\right)  _{j,k}$ is the
$j,k$ th element of $\left(  \rho NI+R_{X,X}\right)  ^{-1}$. Also
\[
0<\left(  \left(  \rho NI+R_{X,X}\right)  ^{-1}\right)  _{k,k}\leq\frac
{1}{\rho N}.
\]

\end{lemma}

\begin{proof}
Using equation \ref{7.28} for $\left(  \mathcal{S}_{X}R_{x^{\left(  j\right)
}}\right)  \left(  x^{\left(  k\right)  }\right)  $ implies%
\begin{align*}
R_{x^{\left(  j\right)  }}\left(  x^{\left(  k\right)  }\right)  -\left(
\mathcal{S}_{X}R_{x^{\left(  j\right)  }}\right)  \left(  x^{\left(  k\right)
}\right)   & =R_{x^{\left(  j\right)  }}\left(  x^{\left(  k\right)  }\right)
-\left(  \widetilde{\mathcal{E}}_{X}R_{x^{\left(  k\right)  }}\right)
^{T}\left(  \rho NI+R_{X,X}\right)  ^{-1}\widetilde{\mathcal{E}}%
_{X}R_{x^{\left(  j\right)  }}\\
& =R_{x^{\left(  j\right)  }}\left(  x^{\left(  k\right)  }\right)  -\left(
\widetilde{\mathcal{E}}_{X}R_{x^{\left(  j\right)  }}\right)  ^{T}\left(  \rho
NI+R_{X,X}\right)  ^{-1}\widetilde{\mathcal{E}}_{X}R_{x^{\left(  k\right)  }%
}\\
& =R_{x^{\left(  j\right)  }}\left(  x^{\left(  k\right)  }\right)  -\left(
\widetilde{\mathcal{E}}_{X}R_{x^{\left(  j\right)  }}\right)  ^{T}\left(  \rho
NI+R_{X,X}\right)  ^{-1}\left(  \rho N\mathbf{i}_{k}+\widetilde{\mathcal{E}%
}_{X}R_{x^{\left(  k\right)  }}\right)  +\\
& \qquad\qquad+\left(  \widetilde{\mathcal{E}}_{X}R_{x^{\left(  j\right)  }%
}\right)  ^{T}\left(  \rho NI+R_{X,X}\right)  ^{-1}\left(  \rho N\mathbf{i}%
_{k}\right) \\
& =R_{x^{\left(  j\right)  }}\left(  x^{\left(  k\right)  }\right)  -\left(
\widetilde{\mathcal{E}}_{X}R_{x^{\left(  j\right)  }}\right)  ^{T}%
\mathbf{i}_{k}+\\
& \qquad\qquad+\left(  \widetilde{\mathcal{E}}_{X}R_{x^{\left(  j\right)  }%
}\right)  ^{T}\left(  \rho NI+R_{X,X}\right)  ^{-1}\left(  \rho N\mathbf{i}%
_{k}\right) \\
& =\left(  \widetilde{\mathcal{E}}_{X}R_{x^{\left(  j\right)  }}\right)
^{T}\left(  \rho NI+R_{X,X}\right)  ^{-1}\left(  \rho N\mathbf{i}_{k}\right)
\\
& =\rho N\left(  \widetilde{\mathcal{E}}_{X}R_{x^{\left(  j\right)  }}\right)
^{T}\left(  \rho NI+R_{X,X}\right)  ^{-1}\mathbf{i}_{k}\\
& =\left(  \rho N\right)  \mathbf{i}_{k}^{T}\left(  \rho NI+R_{X,X}\right)
^{-1}\widetilde{\mathcal{E}}_{X}R_{x^{\left(  j\right)  }}\\
& =\left(  \rho N\right)  \mathbf{i}_{k}^{T}\left(  \rho NI+R_{X,X}\right)
^{-1}\left(  \rho N\mathbf{i}_{j}+\widetilde{\mathcal{E}}_{X}R_{x^{\left(
j\right)  }}\right)  -\\
& \qquad\qquad-\left(  \rho N\right)  \mathbf{i}_{k}^{T}\left(  \rho
NI+R_{X,X}\right)  ^{-1}\left(  \rho N\mathbf{i}_{j}\right) \\
& =\rho N\mathbf{i}_{k}^{T}\mathbf{i}_{j}-\left(  \rho N\right)
^{2}\mathbf{i}_{k}^{T}\left(  \rho NI+R_{X,X}\right)  ^{-1}\mathbf{i}_{j}\\
& =\rho N\delta_{j,k}-\left(  \rho N\right)  ^{2}\left(  \left(  \rho
NI+R_{X,X}\right)  ^{-1}\right)  _{j,k},
\end{align*}

where $\left(  \left(  \rho NI+R_{X,X}\right)  ^{-1}\right)  _{j,k}$ is
element $j,k$ of $\left(  \rho NI+R_{X,X}\right)  ^{-1}$.

Since $\rho NI+R_{X,X}$ is positive definite, $\left(  \rho NI+R_{X,X}\right)
^{-1}$ is positive definite and thus

$\mathbf{i}_{k}^{T}\left(  \rho NI+R_{X,X}\right)  ^{-1}\mathbf{i}_{k}>0$.

\begin{remark}
Equations \ref{7.35} are equivalent to the matrix equation%
\[
R_{X,X}-\left(  \left(  \mathcal{S}_{X}R_{x^{\left(  j\right)  }}\right)
\left(  x^{\left(  i\right)  }\right)  \right)  =\rho NI-\left(  \rho
N\right)  ^{2}\left(  \rho NI+R_{X,X}\right)  ^{-1},
\]

i.e.%
\[
\left(  \left(  R_{x^{\left(  i\right)  }},R_{x^{\left(  j\right)  }}\right)
_{se}\right)  =\rho NI-\left(  \rho N\right)  ^{2}\left(  \rho NI+R_{X,X}%
\right)  ^{-1}.
\]

\end{remark}
\end{proof}

The next theorem shows how the error of the smoother of a function $R_{y}$ is
related to the behavior of the basis function near the origin, as well as to
the number of data points and the smoothing coefficient.

\begin{theorem}
\label{Thm_Rx(x)minus(SRx)(x)_basis_fn_bound}Suppose $\mathcal{S}_{X}$ is the
Exact smoother mapping with smoothing coefficient $\rho$ and independent data
$X=\left\{  x^{\left(  k\right)  }\right\}  _{k=1}^{N}$. Then for each
$x^{\left(  k\right)  }\in X$,%
\[
\sqrt{R_{x}\left(  x\right)  -\left(  \mathcal{S}_{X}R_{x}\right)  \left(
x\right)  }\leq\sqrt{\rho N}\sqrt{1-\rho N\left(  \left(  \rho NI+R_{X,X}%
\right)  ^{-1}\right)  _{k,k}}+\frac{\sqrt{2}}{\left(  2\pi\right)  ^{\frac
{d}{4}}}\sqrt{G\left(  0\right)  -G\left(  x-x^{\left(  k\right)  }\right)  },
\]

where $\left(  \left(  \rho NI+R_{X,X}\right)  ^{-1}\right)  _{k,k}$ is the
$k$th (positive) main diagonal element of the smoother

matrix $\left(  \rho NI+R_{X,X}\right)  ^{-1}$.
\end{theorem}

\begin{proof}
From part 3 of Theorem \ref{Thm_|Rx|se_le_|Ry|se_plus_||Rx_minus_Ry||sq_w0},
$\left\Vert R_{x}\right\Vert _{se}\leq\left\Vert R_{y}\right\Vert
_{se}+\left\Vert R_{x}-R_{y}\right\Vert _{w,0}$ where $\left\Vert
\cdot\right\Vert _{se}$ is the smoother error norm. Setting $y=x^{\left(
k\right)  }\in X$ and using equation \ref{7.35} with $j=k$ this inequality
becomes%
\begin{align*}
\left\Vert R_{x}\right\Vert _{se}=\sqrt{R_{x}\left(  x\right)  -\left(
\mathcal{S}_{X}R_{x}\right)  \left(  x\right)  } &  \leq\sqrt{R_{x^{\left(
k\right)  }}\left(  x^{\left(  k\right)  }\right)  -\left(  \mathcal{S}%
_{X}R_{x^{\left(  k\right)  }}\right)  \left(  x^{\left(  k\right)  }\right)
}+\left\Vert R_{x}-R_{x^{\left(  k\right)  }}\right\Vert _{w,0}\\
&  =\sqrt{\rho N-\left(  \rho N\right)  ^{2}\left(  \left(  \rho
NI+R_{X,X}\right)  ^{-1}\right)  _{k,k}}+\left\Vert R_{x}-R_{x^{\left(
k\right)  }}\right\Vert _{w,0},
\end{align*}

and since $R_{x}\left(  y\right)  =\left(  2\pi\right)  ^{-d/2}G\left(
y-x\right)  $,
\begin{align}
\left\Vert R_{x}-R_{x^{\left(  k\right)  }}\right\Vert _{w,0}^{2}=\left(
R_{x}-R_{x^{\left(  k\right)  }},R_{x}-R_{x^{\left(  k\right)  }}\right)
_{w,0} &  =R_{x}\left(  x\right)  -2R_{x}\left(  x^{\left(  k\right)
}\right)  +R_{x^{\left(  k\right)  }}\left(  x^{\left(  k\right)  }\right)
\nonumber\\
&  =\frac{2}{\left(  2\pi\right)  ^{d/2}}\left(  G\left(  0\right)  -G\left(
x-x^{\left(  k\right)  }\right)  \right)  ,\label{7.36}%
\end{align}

and our result follows.
\end{proof}

Combining the last theorem with the smoother error estimate of Theorem
\ref{Thm_smth_err_bound} we now get:

\begin{corollary}
\label{Cor_Thm_Rx(x)minus(SRx)(x)_basis_fn_bound}If $\mathcal{S}_{X}f_{d}$ is
the Exact smoother (parameter stabilized interpolant) of the data function
$f_{d}\in X_{w}^{0}$ evaluated at the points $\left\{  x^{\left(  k\right)
}\right\}  _{k=1}^{N}$, we have the following pointwise bounds for all
$x\in\mathbb{R}^{d}$:
\[
\left\vert f_{d}\left(  x\right)  -\left(  \mathcal{S}_{X}f_{d}\right)
\left(  x\right)  \right\vert \leq\sqrt{\left(  f_{d}-\mathcal{S}_{X}%
f_{d},f_{d}\right)  _{w,0}}\left(  \sqrt{\rho N}+\frac{\sqrt{2}}{\left(
2\pi\right)  ^{\frac{d}{4}}}\sqrt{G\left(  0\right)  -G\left(  x-x^{\left(
k\right)  }\right)  }\right)  ,
\]

and%
\[
\left\vert f_{d}\left(  x\right)  -\left(  \mathcal{S}_{X}f_{d}\right)
\left(  x\right)  \right\vert \leq\sqrt{\left(  f_{d}-\mathcal{S}_{X}%
f_{d},f_{d}\right)  _{w,0}}\sqrt{R_{0}\left(  0\right)  },
\]

where%
\[
\sqrt{\left(  f_{d}-\mathcal{S}_{X}f_{d},f_{d}\right)  _{w,0}}\leq\left\Vert
f_{d}\right\Vert _{w,0}.
\]

Here $G$ is the basis function and $R_{x}$ is the Riesz representer i.e.
$f\left(  x\right)  =\left(  f,R_{x}\right)  _{w,0}$.
\end{corollary}

\section{Pointwise error estimates - no Taylor series expansions
\label{Sect_exsmth_err_no_Taylor}}

In this section we will prove the Exact smoother analogues of the Type 1 and 2
interpolant error estimates of Section \ref{Sect_interp_no_Taylor_converg}.
Recall that for Type 1 estimates no restriction is placed on the weight
function parameter $\kappa$ but a local smoothness condition is imposed on the
basis function near the origin. On the other hand Type 2 conditions only
assume $\kappa\geq1$ and we use the resultant basis function smoothness
properties. The interpolation weight function examples will be used, augmented
by the central difference weight functions.

\subsection{Type 1 estimates $\left(  \kappa\geq0\right)  $%
\label{SbSect_ex_smth_typ1_error}}

In the next two corollaries a smoothness condition is applied to the basis
function near the origin and this will allow an order estimate to be obtained
for the pointwise smoother of data functions error on a closed bounded
infinite data region. We start with Riesz data functions:

\begin{corollary}
\label{Cor_OrdConvergExactSmth_keq0_Rx}Suppose a weight function has property
W02 for some $\kappa\geq0$ and the (real-valued) basis function $G $ has the
smoothness property assumed by Theorem \ref{Thm_|f(x)-f(y)|_inequal_2}%
\textit{\ i.e. }assume for some $s>0$ and $C_{G},$ $h_{G}>0$ that
\begin{equation}
G\left(  0\right)  -G\left(  x\right)  \leq C_{G}\left\vert x\right\vert
^{2s},\text{\quad}\left\vert x\right\vert \leq h_{G}.\label{7.18}%
\end{equation}

Further, suppose $\mathcal{S}_{X}$ is the Exact smoother mapping with
smoothing coefficient $\rho$ and independent data set $X$ contained in the
bounded closed infinite data set $K$.

Then if $h_{X,K}=\sup\limits_{x\in K}\operatorname*{dist}\left(  x,X\right)
<h_{G}$ it follows that%
\begin{equation}
\left\vert R_{y}\left(  x\right)  -\left(  \mathcal{S}_{X}R_{y}\right)
\left(  x\right)  \right\vert \leq\left(  \sqrt{\rho N}+k_{G}\left(
h_{X,K}\right)  ^{s}\right)  ^{2},\quad x,y\in K,\label{7.13}%
\end{equation}

where
\[
k_{G}=\left(  2\pi\right)  ^{-\frac{d}{4}}\sqrt{2C_{G}}.
\]

Also%
\begin{equation}
\left\vert R_{y}\left(  x\right)  -\left(  \mathcal{S}_{X}R_{y}\right)
\left(  x\right)  \right\vert \leq R_{0}\left(  0\right)  ,\quad
x,y\in\mathbb{R}^{d}.\label{7.17}%
\end{equation}

\end{corollary}

\begin{proof}
Let $X=\left\{  x^{\left(  k\right)  }\right\}  _{k=1}^{N}$ and suppose
$h_{X,K}\leq h_{G}$. Then for each $x\in K$ there exists $x^{\left(  k\right)
}\in X$ such that $G\left(  0\right)  -G\left(  x-x^{\left(  k\right)
}\right)  \leq C_{G}\left\vert x-x^{\left(  k\right)  }\right\vert ^{2s}\leq
C_{G}\left(  h_{X,K}\right)  ^{2s}$. Thus Theorem
\ref{Thm_Rx(x)minus(SRx)(x)_basis_fn_bound} implies that for all $x\in K$%
\begin{align*}
R_{x}\left(  x\right)  -\left(  \mathcal{S}_{X}R_{x}\right)  \left(  x\right)
& \leq\left(  \sqrt{\rho N}+\frac{\sqrt{2}}{\left(  2\pi\right)  ^{\frac{d}%
{4}}}\sqrt{G\left(  0\right)  -G\left(  x-x^{\left(  k\right)  }\right)
}\right)  ^{2}\\
& =\left(  \sqrt{\rho N}+k_{G}\left(  h_{X,K}\right)  ^{s}\right)  ^{2}.
\end{align*}

Suppose $y\in K$. Then setting $f_{d}=R_{y}$ in inequality \ref{7.19} I get%
\begin{align*}
\left\vert R_{y}\left(  x\right)  -\left(  \mathcal{S}_{X}R_{y}\right)
\left(  x\right)  \right\vert  & \leq\sqrt{\left(  R_{y}-\mathcal{S}_{X}%
R_{y},R_{y}\right)  _{w,0}}\sqrt{R_{x}\left(  x\right)  -\left(
\mathcal{S}_{X}R_{x}\right)  \left(  x\right)  }\\
& =\sqrt{R_{y}\left(  y\right)  -\left(  \mathcal{S}_{X}R_{y}\right)  \left(
y\right)  }\sqrt{R_{x}\left(  x\right)  -\left(  \mathcal{S}_{X}R_{x}\right)
\left(  x\right)  }\\
& \leq\left(  \sqrt{\rho N}+k_{G}\left(  h_{X,K}\right)  ^{s}\right)  ^{2}.
\end{align*}

Finally \ref{7.17} follows from \ref{7.26} with $f_{d}=R_{y}$.
\end{proof}

When $\rho=0$ the above results give an interpolant error estimates of order
$2s$. These are the interpolant error estimates obtained in
Chapter\ \ref{Ch_Interpol}. The next corollary is a simple consequence of
inequality \ref{7.19} and Corollary \ref{Cor_OrdConvergExactSmth_keq0_Rx}.
This result is an estimate of the smoother error for an arbitrary data
function in $X_{w}^{0}$.

\begin{corollary}
\label{Cor_OrdConvergExactSmth_k=0_Fd}Suppose:

\begin{enumerate}
\item The basis function $G$ has the properties assumed in Corollary
\ref{Cor_OrdConvergExactSmth_keq0_Rx} for the constants $C_{G},$ $h_{G},$
$s>0$.

\item Suppose $s_{e}$ is the Exact smoother of the arbitrary data function
$f_{d}\in X_{w}^{0}$ on the independent data set $X$ contained in the bounded
closed infinite data set $K$.
\end{enumerate}

Then if $h_{X,K}=\sup\limits_{x\in K}\operatorname*{dist}\left(  x,X\right)
\leq h_{G}$ we have
\begin{equation}
\left\vert f_{d}\left(  x\right)  -s_{e}\left(  x\right)  \right\vert
\leq\sqrt{\left(  f_{d}-s_{e},f_{d}\right)  _{w,0}}\left(  \sqrt{\rho N}%
+k_{G}\left(  h_{X,K}\right)  ^{s}\right)  ,\quad x\in K,\text{ }f_{d}\in
X_{w}^{0},\label{7.21}%
\end{equation}

where $k_{G}$ $=\left(  2\pi\right)  ^{-\frac{d}{4}}\sqrt{2C_{G}}$
and$\sqrt{\left(  f_{d}-s_{e},f_{d}\right)  _{w,0}}<\left\Vert f_{d}%
\right\Vert _{w,0}$.
\end{corollary}

\begin{remark}
\label{Rem_Cor_OrdConvergExactSmth_keq0_Fd}\ 

\begin{enumerate}
\item When $\rho=0$ the last two corollaries yield the Type 1 interpolant
error estimates of Section \ref{Sect_interp_no_Taylor_converg}. Thus we will
label the parameter $s$ the \textit{order of convergence}.

\item The factor $\sqrt{\left(  f_{d}-s_{e},f_{d}\right)  _{w,0}}$ on the
right side of \ref{7.21} enables the error estimate \ref{7.13} for $R_{x}$ to
be recovered immediately. This is because
\end{enumerate}

$\left(  R_{x}-\mathcal{S}_{X}R_{x},R_{x}\right)  _{w,0}=\left(
R_{x}-\mathcal{S}_{X}R_{x}\right)  \left(  x\right)  $ where $\mathcal{S}_{X}$
is the Exact smoother mapping.
\end{remark}

\subsubsection{\protect\underline{Error summary}}

Now we combine the smoother bound of Theorem \ref{Thm_bound_Exact_smth} with
the error estimates derived in this Subsection \ref{SbSect_ex_smth_typ1_error}
to obtain the error estimates that will be tested numerically in Section
\ref{Sect_num_experim}.

\begin{theorem}
\label{Thm_ex_typ1_err_summary}Suppose $K$ is a bounded closed infinite set
and there exist constants $C_{G}$, $s$, $h_{G}>0$ such that
\[
G\left(  0\right)  -G\left(  x\right)  \leq C_{G}\left\vert x\right\vert
^{2s},\text{\quad}\left\vert x\right\vert \leq h_{G},\,x\in K.
\]

Set $k_{G}=\left(  2\pi\right)  ^{-\frac{1}{4}}\sqrt{2C_{G}}$. Then the Exact
smoother $\mathcal{S}_{X}f_{d}$ satisfies the error estimates%
\begin{equation}
\left\vert f_{d}\left(  x\right)  -\left(  \mathcal{S}_{X}f_{d}\right)
\left(  x\right)  \right\vert \leq\min\left\{
\begin{array}
[c]{l}%
\sqrt{\left(  f_{d}-s_{e},f_{d}\right)  _{w,0}}\left(  \sqrt{\rho N}%
+k_{G}\left(  h_{X,K}\right)  ^{s}\right)  ,\\
\\
\sqrt{\left(  f_{d}-s_{e},f_{d}\right)  _{w,0}}\sqrt{R_{0}\left(  0\right)
},\\
\\
\left\Vert f_{d}\right\Vert _{\infty,K}+\left\Vert f_{d}\right\Vert
_{w,0}\sqrt{R_{0}\left(  0\right)  }\min\left\{  1,\frac{R_{0}\left(
0\right)  }{\rho}\right\}  ,
\end{array}
\right. \label{7.09}%
\end{equation}

when $x\in K$ and $h_{X,K}=\sup\limits_{s\in K}\operatorname*{dist}\left(
s,X\right)  \leq h_{G}$. Here $R_{0}\left(  0\right)  =\left(  2\pi\right)
^{-\frac{1}{2}}G\left(  0\right)  $ and $X$ is an independent data set
contained in $K$.

We also have the corresponding double order convergence estimate%
\begin{equation}
\left\vert R_{y}\left(  x\right)  -\left(  \mathcal{S}_{X}R_{y}\right)
\left(  x\right)  \right\vert \leq\min\left\{
\begin{array}
[c]{l}%
\left(  \sqrt{\rho N}+k_{G}\left(  h_{X,K}\right)  ^{s}\right)  ^{2},\\
\\
R_{0}\left(  0\right)  ,\\
\\
\left\Vert R_{y}\right\Vert _{\infty,K}+R_{0}\left(  0\right)  \min\left\{
1,\frac{R_{0}\left(  0\right)  }{\rho}\right\}  ,
\end{array}
\right. \label{7.08}%
\end{equation}

for $x,y\in K$.
\end{theorem}

\begin{proof}
Set $s_{e}=\mathcal{S}_{X}f_{d}$. From Corollary
\ref{Cor_OrdConvergExactSmth_k=0_Fd} and \ref{7.26} of Theorem
\ref{Thm_smth_err_bound} we have%
\[
\left\vert f_{d}\left(  x\right)  -s_{e}\left(  x\right)  \right\vert
\leq\sqrt{\left(  f_{d}-s_{e},f_{d}\right)  _{w,0}}\min\left\{  \left(
\sqrt{\rho N}+k_{G}\left(  h_{X,K}\right)  ^{s}\right)  ,\sqrt{R_{0}\left(
0\right)  }\right\}  .
\]

From Theorem \ref{Thm_bound_Exact_smth}%
\[
\left\vert s_{e}\left(  x\right)  \right\vert \leq\left\Vert f_{d}\right\Vert
_{w,0}\sqrt{R_{0}\left(  0\right)  }\min\left\{  1,\frac{R_{0}\left(
0\right)  }{\rho}\right\}  ,
\]

so that%
\[
\left\vert f_{d}\left(  x\right)  -s_{e}\left(  x\right)  \right\vert
\leq\left\vert f_{d}\left(  x\right)  \right\vert +\left\vert s_{e}\left(
x\right)  \right\vert \leq\left\Vert f_{d}\right\Vert _{\infty,K}+\left\Vert
f_{d}\right\Vert _{w,0}\sqrt{R_{0}\left(  0\right)  }\min\left\{
1,\frac{R_{0}\left(  0\right)  }{\rho}\right\}  ,
\]

which yields \ref{7.09}. Thus when $f_{d}=R_{y}$ and $x\in K$%
\begin{align*}
\left\vert R_{y}\left(  x\right)  -s_{e}\left(  x\right)  \right\vert  &
\leq\left\Vert R_{y}\right\Vert _{\infty,K}+\left\Vert R_{y}\right\Vert
_{w,0}\sqrt{R_{0}\left(  0\right)  }\min\left\{  1,\frac{R_{0}\left(
0\right)  }{\rho}\right\} \\
& =\left\Vert R_{y}\right\Vert _{\infty,K}+R_{0}\left(  0\right)  \min\left\{
1,\frac{R_{0}\left(  0\right)  }{\rho}\right\}  ,
\end{align*}

which combined with the estimates of Corollary
\ref{Cor_OrdConvergExactSmth_keq0_Rx} gives the second set of inequalities
\ref{7.08}.
\end{proof}

\subsubsection{\protect\underline{Smoother convergence}}

I will now introduce the concept of the \textit{convergence of a smoother to
its data function}. Theorem \ref{Thm_seq_data_regions_2} showed that for a
bounded closed infinite data region $K$ there exists a nested sequence of
independent data sets $X^{\left(  k\right)  }\subset X^{\left(  k+1\right)
}\subset K$ such that $h_{X^{\left(  k\right)  },K}\rightarrow0$ as
$k\rightarrow\infty$.

An inspection of the error estimates of this subsection clearly shows that the
smoothing error can tend to zero whilst the smoothing coefficient $\rho$
remains positive. Consider the estimate \ref{7.21} of Corollary
\ref{Cor_OrdConvergExactSmth_k=0_Fd} and suppose $X^{\left(  k\right)  }$ is a
nested sequence of independent data sets such that $h_{X^{\left(  k\right)
},K}\rightarrow0$. Denote the number of points in $X^{\left(  k\right)  }$ by
$N_{k}$. Then given $\varepsilon>0$ there exists $k_{\varepsilon}$ such that

$k_{G}\left\Vert f_{d}\right\Vert _{w,0}\left(  h_{X^{\left(  k_{\varepsilon
}\right)  },K}\right)  ^{s}\leq\frac{\varepsilon}{2}$ when $k\geq
k_{\varepsilon}$. Also for $N\geq N_{k_{\varepsilon}}$ we have $h_{X^{\left(
k\right)  },K}\leq h_{X^{\left(  k_{\varepsilon}\right)  },K}$. Choosing
$\rho_{k}$ so that $\sqrt{\rho_{k}N_{k}}=\frac{\varepsilon}{2\left\Vert
f_{d}\right\Vert _{w,0}}$ an easy calculation shows that
\begin{align*}
\left\vert f_{d}\left(  x\right)  -\mathcal{S}_{X^{\left(  k\right)  }}%
f_{d}\left(  x\right)  \right\vert  &  \leq\left\Vert f_{d}\right\Vert
_{w,0}\left(  \sqrt{\rho N}+k_{G}\left(  h_{X,K}\right)  ^{s}\right) \\
&  \leq\left\Vert f_{d}\right\Vert _{w,0}\sqrt{\rho_{k}N_{k}}+\left\Vert
f_{d}\right\Vert _{w,0}k_{G}\left(  h_{X^{\left(  k\right)  },K}\right)
^{s}\\
&  =\varepsilon.
\end{align*}

We call $s$ the \textit{order of convergence of the smoother}.

\subsubsection{\protect\underline{Type 1 examples}%
\label{SbSbSect_ex_typ1_example}}

The interpolant weight function examples are also used here i.e. the radial
shifted thin-plate splines, Gaussian and Sobolev splines, and the tensor
product extended B-spline weight functions, augmented by the tensor product
central difference weight functions introduced in Chapter
\ref{Ch_cent_diff_wt_fn_ten_prod}.

When $\rho=0$ the error estimates of this subsection are the Type 1
interpolation estimates of Subsection \ref{SbSect_int_estim_Type1} and
consequently the values for $s,C_{G},h_{G},k_{G}$ have already been calculated
and are exhibited in Table \ref{Tbl_NonUnisolvTyp1Converg}. These values are
given below in Table \ref{Tbl_ExNonUnisolvTyp1ConvAugm} augmented by the
values obtained for the central difference tensor product weight functions
from Table \ref{Tbl_ConvergCentral}.%

\begin{table}[htbp] \centering
$%
\begin{tabular}
[c]{|c|c||c|c|c|}\hline
\multicolumn{5}{|c|}{\textbf{Type 1} Exact smoother convergence order
estimates}\\\hline\hline
\multicolumn{5}{|c|}{Smoothness assumption on basis function, $\kappa\geq0$%
.}\\\hline
& Parameter & Converg. &  & \\
Weight function & constraints & order $s$ & $C_{G}$ & $h_{G}$\\\hline\hline
\multicolumn{1}{|l|}{Sobolev splines} & $v-\frac{d}{2}=1$ & $\frac{1}{2}$ &
\multicolumn{1}{|l|}{$\frac{\left\Vert \rho K_{0}\left(  \rho\right)
\right\Vert _{\infty}}{2^{v-1}\Gamma\left(  v\right)  }^{\left(  2\right)  }$}
& $\infty$\\\cline{2-4}%
\multicolumn{1}{|l|}{\quad$\left(  v>d/2\right)  $} & $v-\frac{d}{2}\neq1$ &
$1$ & \multicolumn{1}{|l|}{$\frac{\left\Vert D^{2}\widetilde{K}_{v-d/2}%
\right\Vert _{\infty}}{2^{v}\Gamma\left(  v\right)  }^{\left(  2\right)  }$} &
"\\\hline
\multicolumn{1}{|l|}{Shifted thin-plate} & - & $1$ & \multicolumn{1}{|l|}{eq.
(\ref{7.64})} & $\infty$\\
\multicolumn{1}{|l|}{$\left(  -d/2<v<0\right)  $} &  &  &  & \\\hline
\multicolumn{1}{|l|}{Gaussian} & - & $1$ & \multicolumn{1}{|l|}{$2e^{-3/2}$} &
$\infty$\\\hline
\multicolumn{1}{|l|}{Extended B-spline} & - & $\frac{1}{2}$ &
\multicolumn{1}{|l|}{$G_{1}\left(  0\right)  ^{d-1}\left\Vert DG_{1}%
\right\Vert _{\infty}\sqrt{d}$ $^{\left(  1\right)  }$} & $\infty$\\\hline
\multicolumn{1}{|l|}{$%
\begin{array}
[c]{l}%
\text{Central difference}\\
\text{(tensor product)}%
\end{array}
$} & - & $\frac{1}{2}$ & $\left(  2\pi\right)  ^{-\frac{d}{4}}\sqrt
{2G_{1}\left(  0\right)  ^{d-1}\left\Vert DG_{1}\right\Vert _{\infty}}%
\sqrt[4]{d}$ $^{\left(  1\right)  }$ & $\infty$\\\hline
$%
\begin{array}
[c]{l}%
\text{Central difference}\\
\text{(}q\text{ multivar.)}%
\end{array}
$ & $??$ & $??$ & $??$ & $??$\\\hline\hline
\multicolumn{5}{|l|}{$^{\left(  1\right)  }${\small \ }$G_{1}${\small \ is the
univariate basis function used to form the tensor product.}}\\\hline
\multicolumn{5}{|l|}{$^{\left(  2\right)  }${\small \ }$K_{v}${\small \ is the
modified Bessel function and} $\widetilde{K}_{v}\left(  r\right)  =r^{v}%
K_{v}\left(  r\right)  $.}\\\hline
\end{tabular}
$\caption{}\label{Tbl_ExNonUnisolvTyp1ConvAugm}%
\end{table}%

Equations \ref{7.64} are referenced by Table
\ref{Tbl_ExNonUnisolvTyp1ConvAugm}:%

\begin{equation}
C_{G}=\left\vert \left(  2rf^{\prime\prime}+f^{\prime}\right)  \left(
r_{\max}\right)  \right\vert ,\quad where\text{ }f\left(  r\right)  =\left(
1+r^{2}\right)  ^{v}\text{ }and\text{ }r_{\max}=\frac{\left(  1-2v\right)
}{3}.\label{7.64}%
\end{equation}

\subsection{Type 2 error estimates $\left(  \kappa\geq1\right)  $%
\label{SbSect_ex_smth_typ2_error}}

The Corollaries \ref{Cor_OrdConvergExactSmth_keq0_Rx} and
\ref{Cor_OrdConvergExactSmth_k=0_Fd} which were used to obtain Type 2
interpolant error estimates are again used here by Theorems
\ref{Thm_G(x)minusG(0)_bound}, \ref{Thm_ex_G(x)minusG(0)_bound_w_radial} and
\ref{Thm_Final_Type2_estim} to obtain estimates for $G\left(  0\right)
-G\left(  x\right)  $ when $\kappa\geq1$.

\begin{theorem}
\label{Thm_G(x)minusG(0)_bound}Suppose a weight function satisfies property
W02 for $\kappa=1$ and denote the basis function by $G$. Then the smoother
error estimates of Corollaries \ref{Cor_OrdConvergExactSmth_keq0_Rx} and
\ref{Cor_OrdConvergExactSmth_k=0_Fd} hold where
\begin{equation}
C_{G}=-\frac{1}{2}\left(  \left\vert D\right\vert ^{2}G\right)  \left(
0\right)  d,\text{\quad}s=1,\text{\quad}h_{G}=\infty,\text{\quad}k_{G}=\left(
2\pi\right)  ^{-\frac{d}{4}}\sqrt{-\left(  \left\vert D\right\vert
^{2}G\right)  \left(  0\right)  }\sqrt{d}.\label{7.30}%
\end{equation}

\end{theorem}

\begin{proof}
Theorem \ref{Thm_G(0)minusG(x)_bound_w_W2} showed that if a weight function
has $\kappa\geq1$ then%
\begin{equation}
G\left(  0\right)  -G\left(  x\right)  \leq-\frac{d}{2}\left(  \left\vert
D\right\vert ^{2}G\right)  \left(  0\right)  \left\vert x\right\vert
^{2},\quad x\in\mathbb{R}^{d}.\label{7.01}%
\end{equation}

An application of Corollary \ref{Cor_OrdConvergExactSmth_keq0_Rx} now proves
this theorem directly.
\end{proof}

However, if the weight function is radial we can use the estimates:

\begin{theorem}
\label{Thm_ex_G(x)minusG(0)_bound_w_radial}Suppose a \textbf{radial weight
function} satisfies property W02 for $\kappa=1$ and denote the (radial) basis
function by $G$. Set $r=\left\vert x\right\vert $. Then:

\begin{enumerate}
\item If $G\left(  x\right)  =f\left(  r^{2}\right)  $ then%
\begin{equation}
G\left(  0\right)  -G\left(  x\right)  \leq-f^{\prime}\left(  0\right)
d^{2}\left\vert x\right\vert ^{2},\quad x\in\mathbb{R}^{d}.\label{7.03}%
\end{equation}

\item If $G\left(  x\right)  =g\left(  r\right)  $ then%
\begin{equation}
G\left(  0\right)  -G\left(  x\right)  \leq-\frac{1}{2}g^{\prime\prime}\left(
0\right)  d^{2}\left\vert x\right\vert ^{2},\quad x\in\mathbb{R}%
^{d},\label{7.04}%
\end{equation}

Further, the smoother error estimates of Corollaries
\ref{Cor_OrdConvergExactSmth_keq0_Rx} and \ref{Cor_OrdConvergExactSmth_k=0_Fd}
hold for
\begin{equation}
C_{G}=-f^{\prime}\left(  0\right)  d^{2}\text{ or }-\frac{1}{2}g^{\prime
\prime}\left(  0\right)  d^{2},\text{\quad}s=1,\text{\quad}h_{G}%
=\infty,\text{\quad}k_{G}=\left(  2\pi\right)  ^{-\frac{d}{4}}\sqrt{2C_{G}%
}.\label{7.05}%
\end{equation}

\end{enumerate}
\end{theorem}

\begin{proof}
Follows directly from Theorem \ref{Thm_G(0)minusG(x)_bound_w_W2} and Corollary
\ref{Cor_OrdConvergExactSmth_keq0_Rx}.
\end{proof}

If the weight function is a tensor product the following result will be useful:

\begin{theorem}
\label{Thm_Final_Type2_estim}Suppose a \textbf{tensor product weight function}
satisfies property W02 for $\kappa=1$ and denote the univariate basis function
by $G_{1}$. Then%
\[
G\left(  0\right)  -G\left(  x\right)  \leq-\frac{d}{2}G_{1}\left(  0\right)
^{d-1}D^{2}G_{1}\left(  0\right)  \left\vert x\right\vert ^{2},\quad
x\in\mathbb{R}^{d},
\]

and the smoother error estimates of Corollaries
\ref{Cor_OrdConvergExactSmth_keq0_Rx} and \ref{Cor_OrdConvergExactSmth_k=0_Fd}
hold for
\[
C_{G}=-\frac{d}{2}G_{1}\left(  0\right)  ^{d-1}D^{2}G_{1}\left(  0\right)
,\text{\quad}s=1,\text{\quad}h_{G}=\infty,\text{\quad}k_{G}=\left(
2\pi\right)  ^{-\frac{d}{4}}\sqrt{-G_{1}\left(  0\right)  ^{d-1}D^{2}%
G_{1}\left(  0\right)  }\sqrt{d}.
\]

\end{theorem}

\begin{proof}
Since $G_{1}$ is the univariate basis function
\begin{align*}
\left(  \left\vert D\right\vert ^{2}G\right)  \left(  0\right)  =\sum
\limits_{k=1}^{d}\left(  D_{k}^{2}G\right)  \left(  0\right)  =\sum
\limits_{k=1}^{d}D_{k}^{2}\left(  G_{1}\left(  x_{1}\right)  G_{1}\left(
x_{2}\right)  \ldots G_{1}\left(  x_{d}\right)  \right)  \left(  0\right)   &
=\sum\limits_{k=1}^{d}G_{1}\left(  0\right)  ^{d-1}D^{2}G_{1}\left(  0\right)
\\
&  =G_{1}\left(  0\right)  ^{d-1}D^{2}G_{1}\left(  0\right)  d,
\end{align*}

The estimates of this theorem then follow from Theorem
\ref{Thm_G(x)minusG(0)_bound}, \ref{7.01} and an application of Corollary
\ref{Cor_OrdConvergExactSmth_keq0_Rx} .
\end{proof}

\subsection{Type 2 examples}

When $\rho=0$ the smoother error estimates of Corollaries
\ref{Cor_OrdConvergExactSmth_keq0_Rx} and \ref{Cor_OrdConvergExactSmth_k=0_Fd}
become algebraically identical to the Type 2 interpolant error estimates of
Subsubsection \ref{SbSect_estim_Type2}. Further, the weight function examples
used above were also used for the Type 2 interpolation examples. If we use the
`radial' Theorem \ref{Thm_ex_G(x)minusG(0)_bound_w_radial} to do the estimates
for the radial basis functions and Theorem \ref{Thm_G(x)minusG(0)_bound} to do
the estimates for the tensor product basis functions then the values of the
variables $s,C_{G},h_{G},k_{G}$ will match those obtained for the
interpolants. These are given below in Table
\ref{Tbl_ExNonUnisolvTyp2ConvAugm} which is a copy of Table
\ref{Tbl_InterpNonUnisolvTyp2Conv} augmented by the results for the central
difference tensor product weight functions from Table \ref{Tbl_ConvergCentral}.%

\begin{table}[htbp] \centering
$%
\begin{tabular}
[c]{|c|c||c|c|}\hline
\multicolumn{4}{|c|}{\textbf{Type 2} Exact smoother convergence order
estimates.}\\\hline\hline
\multicolumn{4}{|c|}{We assume $\kappa\geq1$.}\\\hline
& Parameter & Converg. & \\
Weight function & constraints & order & $\left(  2\pi\right)  ^{d/4}%
k_{G}/\sqrt{d}$\\\hline\hline
\multicolumn{1}{|l|}{Sobolev splines} & $v-\frac{d}{2}\geq2$ & $1$ &
$\sqrt{\frac{\Gamma\left(  v-d/2-1\right)  }{2^{d/2+1}\Gamma\left(  v\right)
}}$\\\cline{2-4}\cline{2-4}%
$\left(  v>d/2\right)  $ & $1<v-\frac{d}{2}<2$ & $1$ & $\sqrt{\frac
{\Gamma\left(  v-d/2-1\right)  }{2^{2v-d/2-3}\Gamma\left(  v\right)  }}%
$\\\hline
\multicolumn{1}{|l|}{Shifted thin-plate} & - & $1$ & $\sqrt{-2v}$\\
$\left(  -d/2<v<0\right)  $ &  &  & \\\hline
\multicolumn{1}{|l|}{Gaussian} & - & $1$ & $\sqrt{2}$\\\hline
\multicolumn{1}{|l|}{Extended B-spline} & $n\geq2$ & $1$ & $\sqrt
{-G_{1}\left(  0\right)  ^{d-1}D^{2}G_{1}\left(  0\right)  }^{\text{ }\left(
1\right)  }$\\
$\left(  1\leq n\leq l\right)  $ &  &  & \\\hline
\multicolumn{1}{|l|}{Central difference} & $n\geq2$ & $1$ & $\sqrt
{-G_{1}\left(  0\right)  ^{d-1}D^{2}G_{1}\left(  0\right)  }^{\text{ }\left(
1\right)  }$\\
$\left(  1\leq n\leq l\right)  $ &  &  & \\\hline\hline
\multicolumn{4}{|l|}{$^{\left(  1\right)  }${\small \ }$G_{1}${\small \ is the
univariate basis function used to form the tensor product.}}\\\hline
\end{tabular}
$\caption{}\label{Tbl_ExNonUnisolvTyp2ConvAugm}%
\end{table}%

When $\rho=0$ we see that the order of convergence is at least $1$ for an
arbitrary data function and at least $2$ for a Riesz representer data
function, no matter what value $\kappa$ takes. However, in the next section we
will show that by using the concept of unisolvent data of order $\kappa\geq1$
it follows that an order of convergence of at least $\left\lfloor
\kappa\right\rfloor $ can be attained for an arbitrary data function and an
order of convergence of $\left\lfloor 2\kappa\right\rfloor $ for a Riesz
representer data function.

\section{Pointwise error estimates using unisolvent data subsets
\label{Sect_ex_unisolv_error}}

In this section we will derive smoother error results by applying Lagrange
interpolation theory to the pointwise error estimates obtained using smoother
inner products and semi-inner products of Section
\ref{Sect_ExactSmthErrEstim_InnerProd}. These results are based on showing
that the independent data points must constitute a unisolvent set of specified
order and assumes that the boundary of the data region boundary satisfies a
cone condition.

These results achieve the same order of convergence as the interpolation
results of Subsection \ref{Sect_unisolv} in the sense that when we set
$\rho=0$ we get the same interpolation error estimates.

Unisolvent sets were introduced in Subsection \ref{Sect_unisolv}%
\textit{\ }where they were used to define the Lagrange interpolation operators
$\mathcal{P}$ and $\mathcal{Q}=I-\mathcal{P}$. Lemma
\ref{Lem_Lagrange_interpol} supplied the key interpolation theory and I
reproduce it here:

\begin{lemma}
\label{Lem_Lagrange_interpol_2}(Copy of Lemma \ref{Lem_Lagrange_interpol})
Suppose first that:

\begin{enumerate}
\item $\Omega$ is a bounded region of $\mathbb{R}^{d}$ having the cone property.

\item $X$ is a unisolvent subset of $\Omega$ of order $m$.
\end{enumerate}

Suppose $\left\{  l_{j}\right\}  _{j=1}^{M}$ is the cardinal basis of $P_{m}$
with respect to a minimal unisolvent subset of $\Omega$. Again, using Lagrange
polynomial interpolation techniques, it can be shown there exists a constant
$K_{\Omega,m}^{\prime}>0$ such that
\[
\sum\limits_{j=1}^{M}\left\vert l_{j}\left(  x\right)  \right\vert \leq
K_{\Omega,m}^{\prime}\quad,x\in\Omega,
\]

and all minimal unisolvent subsets of $\Omega$.

Now define the data point density measure
\[
h_{X,\Omega}=\sup\limits_{\omega\in\Omega}\operatorname*{dist}\left(
\omega,X\right)  ,
\]

and fix $x\in X$. By using Lagrange interpolation techniques it can be shown
there are constants $c_{\Omega,m},h_{\Omega,m}>0$ such that when $h_{X,\Omega
}<h_{\Omega,m}$ there exists a minimal unisolvent set $A\subset X$ satisfying
\[
\operatorname*{diam}\left(  A\cup\left\{  x\right\}  \right)  \leq
c_{\Omega,m}h_{X,\Omega}.
\]

\end{lemma}

Our main result for arbitrary data functions is:

\begin{theorem}
\label{Thm_Exact_smth_ord_gte_1}\ 

\begin{enumerate}
\item Let $w$ be a weight function with property W02 for some parameter
$\kappa\geq1$ and let $G$ be the basis function generated by $w$. Set
$m=\left\lfloor \kappa\right\rfloor $.

\item Suppose $s_{e}$ is the Exact smoother of an arbitrary data function
$f_{d}\in X_{w}^{0}$ evaluated on the independent data set $X$ contained in
the data region $\Omega$.

\item Suppose the notation and assumptions of Lemma
\ref{Lem_Lagrange_interpol_2} hold for this theorem i.e. $X$ is $m$-unisolvent
and $\Omega$ is a bounded region whose boundary satisfies the cone condition.
\end{enumerate}

Then there exist positive constants $k_{G},c_{\Omega,m},K_{\Omega,m}^{\prime
},h_{\Omega,m}$ such that the smoothing error satisfies
\begin{equation}
\left\vert f_{d}\left(  x\right)  -s_{e}\left(  x\right)  \right\vert
<\sqrt{\left(  f_{d}-s_{e},f_{d}\right)  _{w,0}}\left(  K_{\Omega,m}^{\prime
}\sqrt{\rho N}+k_{G}\left(  h_{X,\Omega}\right)  ^{m}\right)  ,\quad
x\in\overline{\Omega},\label{7.20}%
\end{equation}

when $h_{X,\Omega}=\sup\limits_{\omega\in\Omega}\operatorname*{dist}\left(
\omega,X\right)  <h_{\Omega,m}$.

Here $k_{G}=\frac{d^{m/2}}{\left(  2\pi\right)  ^{d/2}}\left(  c_{\Omega
,m}\right)  ^{m}K_{\Omega,m}^{\prime}\max\limits_{\left\vert \beta\right\vert
=m}\left\vert D^{2\beta}G\left(  0\right)  \right\vert $ and the constants
$c_{\Omega,m},K_{\Omega,m}^{\prime}$ and $h_{\Omega,m}$ only depend on
$\Omega$, $m$ and $d$. In terms of the integrals which define weight property
W02 I have%
\[
\max_{\left\vert \beta\right\vert =m}\left\vert D^{2\beta}G\left(  0\right)
\right\vert \leq\left(  2\pi\right)  ^{-d/2}\int\tfrac{\left\vert
\xi\right\vert ^{2m}d\xi}{w\left(  \xi\right)  }.
\]

\end{theorem}

\begin{proof}
Fix $x\in\Omega$. By Lemma \ref{Lem_Lagrange_interpol_2} there exists a
minimal unisolvent set $A=\left\{  a_{k}\right\}  _{k=1}^{M}$ of order $m$ and
constants $c_{\Omega,m},h_{X,\Omega}>0$ such that
\[
\operatorname*{diam}\left(  A\cup\left\{  x\right\}  \right)  \leq
c_{\Omega,m}h_{X,\Omega}\text{\quad}when\text{ }h_{X,\Omega}<h_{\Omega,m}.
\]

Further, there exists a positive constant $K_{\Omega,m}^{\prime}$ such that
\[
\sum\limits_{j=1}^{M}\left\vert l_{j}\left(  x\right)  \right\vert \leq
K_{\Omega,m}^{\prime},\quad x\in\Omega,
\]

where $\left\{  l_{j}\right\}  _{j=1}^{M}$ is the cardinal basis of $P_{m}$
with respect to $A$. The next step is to partition the error function using
the Lagrange interpolation projections $\mathcal{P}$ and $\mathcal{Q}%
=I-\mathcal{P}$ so that%
\begin{equation}
f_{d}-s_{e}=\mathcal{P}\left(  f_{d}-s_{e}\right)  +\mathcal{Q}\left(
f_{d}-s_{e}\right)  ,\label{7.2}%
\end{equation}

and note that $\mathcal{P}\left(  f_{d}-s_{e}\right)  $ is not zero. However,
I can still apply Lemma \ref{Lem_Lagrange_interpol_2} and Lemma
\ref{Lem_Q_estim} to the second term $\mathcal{Q}\left(  f_{d}-s_{e}\right)
$, in exactly the the same way as in the proof of Theorem
\ref{Thm_converg_interpol_ord_gte_1}, to show that
\begin{equation}
\left\vert \mathcal{Q}\left(  f_{d}-s_{e}\right)  \left(  x\right)
\right\vert \leq k_{G}\left\Vert f_{d}-s_{e}\right\Vert _{w,0}\left(
h_{X,\Omega}\right)  ^{m},\label{7.22}%
\end{equation}

when $h_{X,\Omega}<h_{\Omega,m}$. Part 2 of Lemma \ref{Lem_ex_smth_err_1} then
implies $\left\Vert s_{e}\right\Vert _{w,0}^{2}\leq\left(  s_{e},f_{d}\right)
_{w,0}$ and so%
\[
\left\vert \mathcal{Q}\left(  f_{d}-s_{e}\right)  \left(  x\right)
\right\vert \leq k_{G}\left(  f_{d}-s_{e},f_{d}\right)  _{w,0}\left(
h_{X,\Omega}\right)  ^{m},
\]

when $h_{X,\Omega}<h_{\Omega,m}$. However, it still remains to estimate
$\mathcal{P}\left(  f_{d}-s_{e}\right)  \left(  x\right)  $. By definition of
$\mathcal{P}$%
\[
\mathcal{P}\left(  f_{d}-s_{e}\right)  \left(  x\right)  =\sum_{k=1}%
^{M}\left(  f_{d}\left(  a_{k}\right)  -s_{e}\left(  a_{k}\right)  \right)
l_{k}\left(  x\right)  ,
\]

and by \ref{7.19},%
\begin{equation}
\left\vert f_{d}\left(  a_{k}\right)  -s_{e}\left(  a_{k}\right)  \right\vert
\leq\sqrt{\left(  f_{d}-s_{e},f_{d}\right)  _{w,0}}\sqrt{R_{a_{k}}\left(
a_{k}\right)  -\left(  \mathcal{S}_{X}R_{a_{k}}\right)  \left(  a_{k}\right)
}.\label{7.43}%
\end{equation}

But by Lemma \ref{Lem_Lagrange_interpol_2}\textit{,} $\sum\limits_{k=1}%
^{M}\left\vert l_{k}\left(  x\right)  \right\vert \leq K_{\Omega,m}^{\prime}$
for $x\in\Omega$, so
\begin{align*}
\left\vert \mathcal{P}\left(  f_{d}-s_{e}\right)  \left(  x\right)
\right\vert  & \leq\sum_{k=1}^{M}\left\vert f_{d}\left(  a_{k}\right)
-s_{e}\left(  a_{k}\right)  \right\vert \left\vert l_{k}\left(  x\right)
\right\vert \\
& \leq\left(  \sum\limits_{k=1}^{M}\left\vert l_{k}\left(  x\right)
\right\vert \right)  \max_{k=1}^{M}\left\vert f_{d}\left(  a_{k}\right)
-s_{e}\left(  a_{k}\right)  \right\vert \\
& \leq K_{\Omega,m}^{\prime}\sqrt{\left(  f_{d}-s_{e},f_{d}\right)  _{w,0}%
}\text{\thinspace}\max_{k=1}^{M}\sqrt{\left(  R_{a_{k}}\left(  a_{k}\right)
-\left(  \mathcal{S}_{X}R_{a_{k}}\right)  \left(  a_{k}\right)  \right)  },
\end{align*}

where $\mathcal{S}_{X}$ is Exact smoother operator. Hence Lemma
\ref{Lem_Rxk(xk).minus.(SRxk)(xk)} implies
\begin{equation}
\left\vert \mathcal{P}\left(  f_{d}-s_{e}\right)  \left(  x\right)
\right\vert \leq K_{\Omega,m}^{\prime}\sqrt{\left(  f_{d}-s_{e},f_{d}\right)
_{w,0}}\sqrt{\rho N}.\label{7.23}%
\end{equation}

Combining inequalities \ref{7.22} and \ref{7.23} gives inequality \ref{7.20}
for all $x\in\Omega$. The extension of the inequality to $\overline{\Omega}$
is an easy consequence of the fact that $f_{d}$ and $s_{e}$ are continuous on
$\mathbb{R}^{d}$.
\end{proof}

Observe that when $\rho=0$ the order of convergence derived above is
$\left\lfloor \kappa\right\rfloor $ and this is the same order derived for the
interpolant of an arbitrary data function in Theorem
\ref{Thm_converg_interpol_ord_gte_1}. I now derive a double order of
convergence result for the Riesz representer data functions $R_{y}$. Note that
the factor $\sqrt{\left(  f_{d}-s_{e},f_{d}\right)  _{w,0}}$ in the right side
of the estimate \ref{7.20} for $\left\vert f_{d}\left(  x\right)
-s_{e}\left(  x\right)  \right\vert $ facilitates the following simple proof:

\begin{corollary}
\label{Cor_Thm_Exact_smth_ord_gte_1}Under the conditions and notation of
Theorem \ref{Thm_Exact_smth_ord_gte_1} the estimate%
\[
R_{y}\left(  x\right)  -\left(  \mathcal{S}_{X}R_{y}\right)  \left(  x\right)
\leq\left(  K_{\Omega,m}^{\prime}\sqrt{\rho N}+k_{G}\left(  h_{X,\Omega
}\right)  ^{m}\right)  ^{2},\quad x,y\in\overline{\Omega},
\]

holds when $h_{X,\Omega}<h_{\Omega,m}$ and $m=\left\lfloor \kappa\right\rfloor
$.
\end{corollary}

\begin{proof}
Suppose $h_{X,\Omega}<h_{\Omega,m}$. When $f_{d}=R_{x}$ the inequality proved
in Theorem \ref{Thm_Exact_smth_ord_gte_1} becomes%
\[
\left(  R_{x}-\mathcal{S}_{X}R_{x}\right)  \left(  x\right)  \leq\sqrt{\left(
R_{x}-\mathcal{S}_{X}R_{x}\right)  \left(  x\right)  }\left(  K_{\Omega
,m}^{\prime}\sqrt{\rho N}+k_{G}\left(  h_{X,\Omega}\right)  ^{m}\right)
,\quad x\in\overline{\Omega},
\]

which implies that%
\[
\left(  R_{x}-\mathcal{S}_{X}R_{x}\right)  \left(  x\right)  \leq\left(
K_{\Omega,m}^{\prime}\sqrt{\rho N}+k_{G}\left(  h_{X,\Omega}\right)
^{m}\right)  ^{2},
\]

and an application of the estimate \ref{7.19} with $f_{d}=R_{y}$ gives the
desired result.
\end{proof}

\subsection{Examples}

The weight function examples used here are those used in for interpolation
i.e. the radial shifted thin-plate splines, Gaussian and Sobolev splines, and
the tensor product extended B-spline weight functions, augmented by the tensor
product central difference weight functions from Chapter
\ref{Ch_cent_diff_wt_fn_ten_prod}.

When the smoothing coefficient $\rho$ is zero the error estimates of this
section are the unisolvent interpolation estimates of Subsection
\ref{Sect_unisolv} and consequently values for $\kappa$ and the Lagrangian
constants $k_{G},c_{\Omega,m},K_{\Omega,m}^{\prime},h_{\Omega,m}$ are the
same. We can't calculate the Lagrangian constants but we estimated the
convergence order for the interpolant as $\left\lfloor \kappa\right\rfloor $
which is given in Table \ref{Tbl_UnisolvConverg}. We use this value here and
they are given below in Table \ref{Tbl_ExUnisolvConverg} augmented by the
order values obtained for the central difference tensor product weight
functions from Table \ref{Tbl_ConvergCentral}.%

\begin{table}[htbp] \centering
$%
\begin{tabular}
[c]{|c|c||c|c|}\hline
\multicolumn{4}{|c|}{Exact smoother convergence order estimates for
$\kappa\geq1$.}\\\hline\hline
\multicolumn{4}{|c|}{Uses Lagrange interpolation and Taylor series.}\\\hline
& Parameter & Convergence & \\
Weight function & constraints & order\thinspace$\left(  \left\lfloor
\kappa\right\rfloor \right)  $ & $h_{\Omega,\left\lfloor \kappa\right\rfloor
}$\\\hline\hline
\multicolumn{1}{|l|}{Sobolev splines} & \multicolumn{1}{|l||}{$v-\frac{d}%
{2}=2,3,4,\ldots$} & \multicolumn{1}{||l|}{$\left\lfloor v-\frac{d}%
{2}\right\rfloor -1$} & $\infty$\\\cline{2-4}\cline{2-4}%
\multicolumn{1}{|l|}{\quad$\left(  v>d/2\right)  $} &
\multicolumn{1}{|l||}{$v-\frac{d}{2}>1$, $v-\frac{d}{2}\notin\mathbb{Z}_{+}$}
& \multicolumn{1}{||l|}{$\left\lfloor v-\frac{d}{2}\right\rfloor $} &
"\\\hline
\multicolumn{1}{|l|}{Shifted thin-plate} & - &
\multicolumn{1}{||l|}{$2,3,4,\ldots$} & $\infty$\\\hline
\multicolumn{1}{|l|}{Gaussian} & - & \multicolumn{1}{||l|}{$2,3,4,\ldots$} &
$\infty$\\\hline
\multicolumn{1}{|l|}{Extended B-spline} & $n\geq2$ &
\multicolumn{1}{||l|}{$n-1$} & $\infty$\\
\multicolumn{1}{|l|}{\quad$\left(  1\leq n\leq l\right)  $} &  &
\multicolumn{1}{||l|}{} & \\\hline
\multicolumn{1}{|l|}{Central difference:} & $n\geq2$ &
\multicolumn{1}{||l|}{$n-1$} & $\infty$\\
\multicolumn{1}{|l|}{ten. prod.$\left(  1\leq n\leq l\right)  $} &  &
\multicolumn{1}{||l|}{} & \\\hline
\end{tabular}
$\caption{}\label{Tbl_ExUnisolvConverg}%
\end{table}%

\section{Numerical experiments with the extended
B-splines\label{Sect_num_experim}}

In this section the convergence of the Exact smoother to its data function is
studied numerically using scaled, extended B-splines, an example of which is
the hat function. We will use the same 1-dimensional data functions and the
same extended B-splines, introduced in Theorem
\ref{Thm_ex_splin_wt_fn_properties}, that were used for the numerical
interpolation experiments described in Section
\ref{Sect_int_data_fn_exten_Bsplin}.

We will only consider the numerical experiments in one dimension so that our
standard data density parameter $h_{X,K}=\sup\limits_{x\in K}%
\operatorname*{dist}\left(  x,X\right)  $ can be easily calculated. It is also
easier to test the Type1 and 2 error estimates i.e. those derived without
explicitly assuming unisolvent data sets. This is because in these cases the
constants are given precisely. In the case of the unisolvent data the theory
is complex and it is unclear what are suitable upper bounds for the constants
e.g. $K_{\Omega,m}^{\prime},c_{\Omega,m},h_{\Omega,m}$ in Theorem
\ref{Thm_Exact_smth_ord_gte_1}.

Recall from Theorem \ref{Thm_ex_splin_wt_fn_properties} that the extended
B-splines are a special class of tensor product weight functions which satisfy
property W02. Indeed, for given integers $l,n$ the multivariate weight
function is defined by
\begin{equation}
w_{s}\left(  x\right)  =\prod_{i=1}^{d}\frac{x_{i}^{2n}}{\sin^{2l}x_{i}%
},\text{\quad}x=\left(  x_{1},\ldots,x_{d}\right)  \in\mathbb{R}%
^{d},\label{7.65}%
\end{equation}

where the weight function $w_{s}$ has property W02 for $\kappa$ iff $n$ and
$l$ satisfy
\begin{equation}
\kappa+1/2<n\leq l.\label{7.66}%
\end{equation}

Theorem \ \ref{Thm_basis_tensor_hat_W3} showed that the basis function $G$
generated by $w$ is the tensor product $G\left(  x\right)  =%
{\textstyle\prod\limits_{k=1}^{d}}
G_{1}\left(  x_{k}\right)  $ where
\begin{equation}
G_{1}\left(  t\right)  =\left(  -1\right)  ^{l-n}\tfrac{\left(  2\pi\right)
^{l/2}}{2^{2\left(  l-n\right)  +1}}D^{2\left(  l-n\right)  }\left(  \left(
\ast\Lambda\right)  ^{l}\right)  \left(  \tfrac{t}{2}\right)  ,\text{\quad
}t\in\mathbb{R}^{1},\label{7.67}%
\end{equation}

and $\left(  \ast\Lambda\right)  ^{l}$ denotes the convolution of $l$
1-dimensional hat functions. Further, $G_{1}\in C_{0}^{\left(  2n-2\right)
}\left(  \mathbb{R}^{1}\right)  $ and $D^{2n-1}G_{1}$ is a piecewise constant
function. Finally, $G\in C_{0}^{\left(  2n-2\right)  \mathbf{1}}\left(
\mathbb{R}^{d}\right)  $ and the derivatives $\left\{  D^{\alpha}G\right\}
_{\alpha\leq\left(  2n-1\right)  \mathbf{1}}$ are bounded functions.

Because all the extended B-spline basis weight functions \ref{7.65} have a
power of $\sin x$ in the denominator we will need the special classes of data
functions developed in Section \ref{Sect_int_data_fn_exten_Bsplin}.

The 1-dimensional independent $X$ data set is constructed using a uniform
distribution on the interval $K=\left[  -1.5,1.5\right]  $. Each of 20 data
files is exponentially sampled using a multiplier of approximately 1.2 and a
maximum of 5000 points, and then we plot $\log_{10}h_{X,K}$ against $\log
_{10}N$ where $N=\left\vert X\right\vert $. It then seems quite reasonable to
use a least-squares linear fit and in this case we obtain%
\begin{equation}
h_{X,K}\simeq3.09N^{-0.81}.\label{7.75}%
\end{equation}

For ease of calculation let
\begin{equation}
h_{X,K}=h_{1}N^{-a},\text{\quad}h_{1}=3.09,\text{\quad}a=0.81.\label{7.6}%
\end{equation}

Noting the error estimates of Theorem \ref{Thm_ex_typ1_err_summary} we use
\ref{7.6} to write $h_{X,K}^{s}=h_{1}^{s}N^{-as}$ and define the Exact
smoother error estimates%
\begin{equation}
\left\vert f_{d}\left(  x\right)  -\left(  \mathcal{S}_{X}f_{d}\right)
\left(  x\right)  \right\vert \lesssim E_{N}\left(  \rho\right)  ,\text{\quad
}x\in K,\label{7.1}%
\end{equation}

and%
\begin{equation}
\left\vert R_{y}\left(  x\right)  -\left(  \mathcal{S}_{X}R_{y}\right)
\left(  x\right)  \right\vert \lesssim\varepsilon_{N}\left(  \rho\right)
,\text{\quad}x,y\in K,\label{7.11}%
\end{equation}

where%
\begin{equation}
E_{N}\left(  \rho\right)  =\min\left\{
\begin{array}
[c]{l}%
\left\Vert f_{d}\right\Vert _{w,0}\left(  \sqrt{\rho N}+k_{G}h_{1}^{s}%
N^{-as}\right)  ,\\
\\
\left\Vert f_{d}\right\Vert _{w,0}\sqrt{R_{0}\left(  0\right)  },\\
\\
\left\Vert f_{d}\right\Vert _{\infty,K}+\left\Vert f_{d}\right\Vert
_{w,0}\sqrt{R_{0}\left(  0\right)  }\min\left\{  1,\frac{R_{0}\left(
0\right)  }{\rho}\right\}  ,
\end{array}
\right. \label{7.10}%
\end{equation}

and%
\begin{equation}
\varepsilon_{N}\left(  \rho\right)  =\min\left\{
\begin{array}
[c]{l}%
\left(  \sqrt{\rho N}+k_{G}h_{1}^{s}N^{-as}\right)  ^{2},\\
\\
R_{0}\left(  0\right)  ,\\
\\
\left\Vert R_{y}\right\Vert _{\infty,K}+R_{0}\left(  0\right)  \min\left\{
1,\frac{R_{0}\left(  0\right)  }{\rho}\right\}  .
\end{array}
\right. \label{7.14}%
\end{equation}

\subsection{Extended B-splines with $n=1$}

\subsubsection{\protect\underline{The case $n=1$, $l=1$}}

The 1-dimensional hat weight function is given by $w_{s}\left(  \xi\right)  =$

$\sqrt{2\pi}\left(  \frac{\xi/2}{\sin\left(  \xi/2\right)  }\right)  ^{2}$
(\ref{1.30}) and examining \ref{7.67} we see that it is a scaled, extended
B-spline weight function corresponding to the parameters $n=l=1$. Also
\ref{7.66} implies $\max\left\lfloor \kappa\right\rfloor =0$.

Now suppose $\Pi$ is the 1-dimensional rectangular function given by
$\Pi\left(  x\right)  =1$ when $\left\vert x\right\vert <1/2$ and $\Pi\left(
x\right)  =0$ when $\left\vert x\right\vert >1/2$. Then it was shown in
Subsubsection \ref{SbSbSect_int_Bspl_neq1_leq1} that:

\begin{theorem}
\textbf{Data function} Suppose $w_{s}$ is the univariate hat weight function.
Then $u\ast\Pi\in X_{w_{s}}^{0}$ if $u\in L^{2}\left(  \mathbb{R}^{1}\right)
$. Further, we can define $V\in L_{loc}^{1}\cap C_{BP}^{\left(  0\right)  }$
by%
\[
V\left(  x\right)  =\int_{0}^{x}u\left(  t\right)  dt,\quad u\in L^{2},
\]

so that $DV=u$ a.e. and%
\[
u\ast\Pi=\left(  2\pi\right)  ^{-\frac{1}{2}}\left(  V\left(  x+\frac{1}%
{2}\right)  -V\left(  x-\frac{1}{2}\right)  \right)  .
\]

\end{theorem}

To obtain our data function $f_{d}=u\ast\Pi$ we will follow the interpolation
approach of Subsubsection \ref{SbSbSect_int_Bspl_neq1_leq1} and choose
\[
u=e^{-x^{2}},
\]
for which
\begin{equation}
V=\left(  2\pi\right)  ^{\frac{1}{2}}\operatorname{erf},\text{\qquad
}\left\Vert u\right\Vert _{2}=2\left(  2\pi\right)  ^{\frac{1}{4}%
},\text{\qquad}\left\Vert u\ast\Pi\right\Vert _{w_{s},0}=2,\label{7.014}%
\end{equation}

and%
\begin{equation}
f_{d}=u\ast\Pi=\operatorname{erf}\left(  x+\frac{1}{2}\right)
-\operatorname{erf}\left(  x-\frac{1}{2}\right)  .\label{7.002}%
\end{equation}

Since we are using the hat basis function, $G\left(  0\right)  =1$, $C_{G}=1$,
$s=1/2$, $h_{G}=\infty$, and from \ref{7.014} and \ref{7.002}, $\left\Vert
f_{d}\right\Vert _{w,0}=2$.

For the \textit{double rate} convergence experiment we will use $R_{0}=\left(
2\pi\right)  ^{-\frac{1}{2}}\Lambda$ as the data function.

\subsubsection{\protect\underline{Numerical results}}

We start by plotting the error bounds given by \ref{7.09} and \ref{7.08},
together with the actual absolute smoother error, against the smoothing
parameter. The results for the data function $f_{d}$ are shown in Figure
\ref{7.12} and the results for the data function $R_{0}$ are shown in Figure
\ref{7.15}:%

\begin{figure}[th]%
\centering
\includegraphics[
natheight=4.467600in,
natwidth=6.059700in,
height=4.4676in,
width=6.0597in
]%
{C:/Math_SwBasisFunc/InterpolSmthDev/PapersMonog/ZeroOrd/ZeroOrdDev/graphics/figExactSmthErrVsSmthParm_HatErf__8.pdf}%
\caption{Exact smoother errors vs smooth parm: data func is
$\operatorname{erf}\left(  x+\frac{1}{2}\right)  -\operatorname{erf}\left(
x-\frac{1}{2}\right)  $.}%
\label{7.12}%
\end{figure}
%

\begin{figure}[th]%
\centering
\includegraphics[
natheight=4.060300in,
natwidth=6.062300in,
height=4.0603in,
width=6.0623in
]%
{C:/Math_SwBasisFunc/InterpolSmthDev/PapersMonog/ZeroOrd/ZeroOrdDev/graphics/figExactSmthErrVsSmthParm_HatRiesz__9.pdf}%
\caption{Exact smoother errors vs smooth param: hat basis func, data func is
$R_{0}$.}%
\label{7.15}%
\end{figure}

The Figure \ref{7.12} seems typical for $C^{\infty}$ data functions and it is
clear that except for smoothing parameters larger than 10 the predicted error
bounds given in \ref{7.10} and \ref{7.14} capture only a small part of the
actual convergence rate. Indeed, the error ratio for the function $R_{0}$ is
$10^{1.5}$ and that for $f_{d}$ ranges from $10^{1.5}$ to $2\times10^{3}$.

We now move on to examine the relationship between the smoother error and the
data density. To do this we need to fix a value for the smoothing parameter
and we choose $\rho=10^{-6}$. We select this value because often the Exact
smoother error curve has a minimum near this value. Using the functions and
parameters discussed at the end of the last subsection we obtain the four
subplots displayed in Figure \ref{Fig_ExactConverg_N1_L1}, each display being
the superposition of 20 smoothers.%

\begin{figure}[th]%
\centering
\includegraphics[
natheight=4.785900in,
natwidth=3.461000in,
height=4.7859in,
width=3.461in
]%
{C:/Math_SwBasisFunc/InterpolSmthDev/PapersMonog/ZeroOrd/ZeroOrdDev/graphics/figExactConverg_N1_L1_samp20_5Kpt_scal_half__10.pdf}%
\caption{Convergence of the Exact smoother.}%
\label{Fig_ExactConverg_N1_L1}%
\end{figure}

The two upper subplots relate to the data function \ref{7.002} and the lower
subplot relates to the Riesz representer data function $R_{0}=\left(
2\pi\right)  ^{-d/2}G\left(  x-\cdot\right)  $ (see \ref{1.019}). The
right-hand subplots are filtered versions of the actual error. The data
function is given at the top of the left-hand plots and the annotation at the
bottom of the figure supplies the following additional information:\medskip

\fbox{Input parameters}\smallskip

$\mathbf{N=L=1}$ - the hat function is a member of the family of scaled,
extended B-splines with the indicated parameter values.

\textbf{spl scale 1/2} - the actual scaling of the spline basis function is
1/2 divided by \textbf{spl scale}.

\textbf{sm parm 1e-6} - the smoothing parameter is $10^{-6}$.

\textbf{samp 20} - the sample size. This is the number of test data files
generated. The data function is evaluated on the interval $[-1.5,1.5]$ at
points selected using a uniform (statistical) distribution.

\textbf{pts 2:5K} - the smallest number of data points is 2 and the largest
number of data points is 5000. The other values are given in exponential steps
with a multiplier of approximately 1.3.\medskip

\fbox{Output parameters/messages}\smallskip

\textbf{No ill-condit} - this indicates all Exact smoother matrices were
always properly conditioned.\medskip

Note that all the plots shown below have the same format and
annotations.\medskip

As mentioned above the smoother is filtered. The filter calculates the value
below which 80\% of the errors lie. The filter is designed to remove `large',
isolated spikes which dominate the actual smoother errors. The smoother error
was calculated on a grid with 300 cells applied to the domain of the data
function. No filter is used for the first five smoothers because there is no
instability for small numbers of data points. I conclude that in this case the
`spike' filter can only meaningfully applied to the smoother of $R_{0}$. This
filter will be used in all the cases below.

As the number of points increases the numerical smoother of $R_{0}$ is
observed to simplify to three very large increasingly narrow spikes at $0$ and
$\pm1$ and these dominate by about four orders of magnitude a stable residual
error function of uniform amplitude and zero trend. The smoother of $f_{d}$
consists of intermingled spikes of various heights superimposed on a trend
curve of amplitude comparable to the average residual spike size. The maximum
spike height is at most one order of magnitude of the average spike height so
the smoother is stable.

The (blue) points above each smoother in Figure \ref{Fig_ExactConverg_N1_L1}
represent the theoretical upper bound for the error given by the estimates
\ref{7.09} or \ref{7.08}, and the adjacent (red) line has a slope which gives
the actual rate of convergence, ignoring instability. Clearly for the data
function $f_{d}$ the theoretical error bound substantially underestimates the
actual error. Regarding the data function $R_{0}$, the estimated `double'
convergence rate is $2s=1$ and the estimated (unfiltered) rate is also 1 so
for $\rho=10^{-6}$ the theoretical upper error bound for the data function
$R_{0}$ is able to take into account quite closely the instability of the
smoother. However, there is a large difference when the filtered error is considered.

\subsubsection{\protect\underline{The case $n=1$, $l=2$}}

Amongst other things the following result, extracted from Theorem
\ref{Thm_data_func}, will allow us to generate data functions for
1-dimensional extended B-spline spaces for which the $X_{w}^{0}$ norm can be
calculated \ This result is closely related to the calculations done above for
the hat function.

\begin{theorem}
\label{Thm_data_func_2}(Portion of Theorem \ref{Thm_data_func}) Suppose $w$ is
the extended B-spline weight function with parameters $n$ and $l$ given by
\ref{1.032}, and that $U\in L^{2}\left(  \mathbb{R}^{1}\right)  $, $D^{n}U\in
L^{2}$ in the sense of distributions. Then if we define the distribution%
\[
f_{d}=\delta_{2}^{l}U,\quad l=1,2,3,\ldots,
\]

where $\delta_{2}$ is the central difference operator%
\begin{equation}
\delta_{2}U=U\left(  \cdot+1\right)  -U\left(  \cdot-1\right)  ,\label{7.001}%
\end{equation}

it follows that $f_{d}\in X_{w}^{0}$ and%
\begin{equation}
\left\Vert f_{d}\right\Vert _{w,0}=2^{l}\left\Vert D^{n}U\right\Vert
_{2}.\label{7.003}%
\end{equation}

\end{theorem}

Our basis function is the extended B-spline $G_{1}$ with parameters $n=1$ and
$l=2$ given by \ref{7.67} i.e.
\begin{align*}
G_{1}\left(  t\right)   & =\left(  -1\right)  ^{l-n}\tfrac{\left(
2\pi\right)  ^{l/2}}{2^{2\left(  l-n\right)  +1}}\left(  D^{2\left(
l-n\right)  }\left(  \left(  \ast\Lambda\right)  ^{l}\right)  \right)  \left(
\tfrac{t}{2}\right) \\
& =-\tfrac{\pi}{4}\left(  D^{2}\left(  \Lambda\ast\Lambda\right)  \right)
\left(  \tfrac{t}{2}\right)  ,\text{\quad}t\in\mathbb{R}^{1}.
\end{align*}
But%
\begin{align*}
D^{2}\left(  \Lambda\ast\Lambda\right)  =\Lambda\ast D^{2}\Lambda &
=\Lambda\ast\left(  \delta\left(  \cdot+1\right)  -2\delta+\delta\left(
\cdot-1\right)  \right) \\
&  =\frac{1}{\sqrt{2\pi}}\left(  \Lambda\left(  \cdot+1\right)  -2\Lambda
+\Lambda\left(  \cdot-1\right)  \right)  ,
\end{align*}

so that%
\[
G_{1}\left(  t\right)  =-\frac{\sqrt{2\pi}}{8}\left(  \Lambda\left(  \tfrac
{t}{2}+1\right)  -2\Lambda\left(  \tfrac{t}{2}\right)  +\Lambda\left(
\tfrac{t}{2}-1\right)  \right)  ,
\]

and hence%
\[
DG_{1}\left(  t\right)  =-\frac{\sqrt{2\pi}}{16}\left(  \Lambda^{\prime
}\left(  \tfrac{t}{2}+1\right)  -2\Lambda^{\prime}\left(  \tfrac{t}{2}\right)
+\Lambda^{\prime}\left(  \tfrac{t}{2}-1\right)  \right)  ,
\]

i.e. $\left\Vert DG_{1}\right\Vert _{\infty}=\frac{3}{16}\sqrt{2\pi}$. Since
the (distributional) derivative is bounded the distributional Taylor series
expansion of Lemma \ref{Lem_Taylor_extension} can be used to write%
\[
G_{1}\left(  0\right)  -G_{1}\left(  t\right)  \leq\left\Vert DG_{1}%
\right\Vert _{\infty}\left\vert t\right\vert ,\text{\quad}x\in\mathbb{R}^{1},
\]

which means that
\[
G_{1}\left(  0\right)  =\frac{\sqrt{2\pi}}{4},\text{ }C_{G}=\left\Vert
DG_{1}\right\Vert _{\infty}=\frac{3}{16}\sqrt{2\pi},\text{ }s=\frac{1}%
{2},\text{ }h_{G}=\infty.
\]

With reference to the last theorem we will choose the bell-shaped data
function%
\begin{equation}
f_{d}=\delta_{2}^{2}U\in X_{w}^{0},\label{7.70}%
\end{equation}

where%
\[
U\left(  x\right)  =\frac{e^{-k_{1,2}x^{2}}}{\delta_{2}^{2}\left(
e^{-k_{1,2}x^{2}}\right)  \left(  0\right)  }=\frac{e^{-k_{1,2}x^{2}}%
}{2\left(  1-e^{-4k_{1,2}}\right)  },\quad k_{1,2}=0.3,
\]

and
\[
\left\Vert f_{d}\right\Vert _{w,0}=4\left\Vert DU\right\Vert _{2}%
=\sqrt[4]{2\pi}\frac{\sqrt[4]{4k_{1,2}}}{1-e^{-4k_{1,2}}}.
\]

The error estimates are given by \ref{7.09} or \ref{7.08} where $k_{G}=\left(
2\pi\right)  ^{-\frac{1}{4}}\sqrt{2C_{G}}$ and $R_{0}\left(  0\right)
=\left(  2\pi\right)  ^{-\frac{1}{2}}G_{1}\left(  0\right)  $.

For the theory developed in the previous sections it was convenient to use the
simple, unscaled weight function definition \ref{7.65} but for computations it
may be easier to use the unscaled version of the extended B-spline basis
function given in the next theorem.

\begin{theorem}
\label{Cor_Thm_data_func_2}Suppose $G=\left(  -1\right)  ^{l-n}D^{2\left(
l-n\right)  }\left(  \left(  \ast\Lambda\right)  ^{l}\right)  $ where
$\Lambda$ is the 1-dimensional hat function and $n,l$ are integers such that
$1\leq n\leq l$.

Then for given $\lambda>0$, $G\left(  \lambda x\right)  $ is called a scaled
extended B-spline basis function. The corresponding weight function is
$w_{\lambda}\left(  t\right)  =2\lambda aw\left(  \frac{t}{2\lambda}\right)  $
where $w$ is the extended B-spline weight function with parameters $n,l$ and
$a=\frac{\left(  2\pi\right)  ^{l/2}}{2^{2\left(  l-n\right)  +1}}$. Indeed,
for $w_{\lambda}$ we can choose $\kappa=n-1$.

Further%
\begin{equation}
G\left(  0\right)  -G\left(  \lambda x\right)  \leq\lambda\left\Vert
DG\right\Vert _{\infty}\left\vert x\right\vert ,\text{\quad}x\in\mathbb{R}%
^{1}.\label{7.016}%
\end{equation}

Finally, if $f_{d}\in X_{w}^{0}$ and $g_{d}\left(  x\right)  =f_{d}\left(
2\lambda x\right)  $, it follows that $g_{d}\in X_{w_{\lambda}}^{0}$ and
$\left\Vert g_{d}\right\Vert _{w_{\lambda},0}=\sqrt{a}\left\Vert
f_{d}\right\Vert _{w,0}$.
\end{theorem}

\subsubsection{\protect\underline{Numerical results}}

As with the previous case we select $\rho=10^{-6}$ and, using the spike filter
of the previous case and the functions and parameters discussed in the last
subsection, we obtain the four subplots displayed in Figure
\ref{Fig_ExactConverg_N1_L2}, each display being the superposition of 20 smoothers.%

\begin{figure}[th]%
\centering
\includegraphics[
natheight=4.792800in,
natwidth=3.410800in,
height=4.7928in,
width=3.4108in
]%
{C:/Math_SwBasisFunc/InterpolSmthDev/PapersMonog/ZeroOrd/ZeroOrdDev/graphics/figExactConverg_N1_L2_samp20_5Kpt_scal_1__11.pdf}%
\caption{Convergence of Exact smoother.}%
\label{Fig_ExactConverg_N1_L2}%
\end{figure}

These results are significantly different to those obtained for the previous
case $n=l=1$. As the number of points increases the numerical smoother of
$R_{0}$ is observed to quickly simplify to a single very large increasingly
narrow spike at the origin, and this dominates a stable error function of
uniform amplitude and zero trend. The smoother of $f_{d}$ typically forms two
major spikes at the boundary $x=\pm3/2$ by 100 data points and then two minor
spikes at $x=\pm1/2$ at about 1000 data points. These spikes are superimposed
on a trend curve of amplitude comparable to the average spike size of a small
stable error function.

The (blue) points above the actual errors in Figure
\ref{Fig_ExactConverg_N1_L2} are the estimated upper bounds for the error
given by the estimates \ref{7.09} or \ref{7.08}, and the adjacent (red) line
gives the average convergence rate. Clearly for the data function $f_{d}$ the
theoretical error bound substantially underestimates the convergence rate.
Regarding the data function $R_{0}$, for $\rho=10^{-6}$ the theoretical upper
error bound for the data function $R_{0}$ is able to take into account quite
closely the instability of the smoother. On the other hand, the filtered
curves for the data function $R_{0}$ are obtained from the corresponding
unfiltered curves by inserting a large step when the data density is about 0.5
and so the theoretical error estimates are at least 4 orders of magnitude
greater than the actual errors.

\subsection{Extended B-splines with $n=2$\label{SbSect_ex_data_fn_n_eq_2}}

Since $n\geq2$, we can use the error estimates of Theorem
\ref{Thm_Final_Type2_estim} as well as the same scaled B-spline and data
function that we used for the interpolant in Subsection
\ref{SbSect_ex_data_fn_n_eq_2}.

\subsubsection{\protect\underline{The case $n=2$, $l=2$}}

The basis function is
\[
G_{2,2}\left(  x\right)  =\frac{\left(  \Lambda\ast\Lambda\right)  \left(
2x\right)  }{\left(  \Lambda\ast\Lambda\right)  \left(  0\right)  }%
=\frac{3\sqrt{2\pi}}{2}\left(  \Lambda\ast\Lambda\right)  \left(  2x\right)  ,
\]
with scaling factor $\lambda=4$ and is such that $\operatorname*{supp}%
G_{2,2}=\left[  -1,1\right]  $ and $G_{2,2}\left(  0\right)  =1$. To calculate
$G_{2,2}$ we use the formula%
\[
G_{2,2}\left(  x\right)  =\left(  1+x\right)  ^{2}\Lambda\left(  2x+1\right)
+\left(  1-2x^{2}\right)  \Lambda\left(  2x\right)  +\left(  1-x\right)
^{2}\Lambda\left(  2x-1\right)  ,
\]

and choose the bell-shaped data function%
\begin{equation}
f_{d}=\delta_{2}^{2}U\in X_{w}^{0},\label{7.71}%
\end{equation}

where%
\begin{equation}
U\left(  x\right)  =\frac{e^{-k_{1,2}x^{2}}}{\delta_{2}^{2}\left(
e^{-k_{1,2}x^{2}}\right)  \left(  0\right)  }=\frac{e^{-k_{1,2}x^{2}}%
}{2\left(  1-e^{-4k_{1,2}}\right)  },\quad k_{1,2}=0.3,\label{7.72}%
\end{equation}

so that%
\begin{equation}
k_{G}=\frac{\sqrt{12}}{\left(  2\pi\right)  ^{1/4}};\quad\left\Vert
f_{d}\right\Vert _{w,0}=\sqrt[4]{72\pi}\frac{\left(  k_{1,2}\right)  ^{3/4}%
}{1-e^{-4k_{1,2}}},\quad k_{1,2}=0.3;\quad h_{G}=\infty.\label{7.73}%
\end{equation}

For double rate convergence experiments $R_{0}=\left(  2\pi\right)
^{-\frac{1}{2}}G_{2,2}$ is the data function and $R_{0}\left(  0\right)
=\left(  2\pi\right)  ^{-\frac{1}{2}}$.

The next step is to calculate $\left\Vert f_{d}\right\Vert _{\infty,K}$ and
$\left\Vert R_{0}\right\Vert _{\infty,K}$ where $K=\left[  -1.5,1.5\right]  $.
Clearly\
\[
\left\Vert R_{0}\right\Vert _{\infty,K}=\left(  2\pi\right)  ^{-\frac{1}{2}%
}G_{2,2}\left(  0\right)  =\left(  2\pi\right)  ^{-\frac{1}{2}},
\]

and from \ref{7.71} and \ref{7.72}, $\left\Vert f_{d}\right\Vert _{\infty
,K}=f_{d}\left(  0\right)  =1$:%
\begin{equation}
\left\Vert R_{0}\right\Vert _{\infty,K}=R_{0}\left(  0\right)  =\left(
2\pi\right)  ^{-\frac{1}{2}},\quad\left\Vert f_{d}\right\Vert _{\infty
,K}=1.\label{7.74}%
\end{equation}

\subsubsection{\protect\underline{Numerical results}}

Using these functions and parameters the four subplots of Figure
\ref{Fig_ExactConverg_N2_L2} each display the superposition of 20 interpolants
or 20 filtered interpolants as well the curves which show the upper bounds for
the error given by the estimates \ref{7.09} and \ref{7.08}.%

\begin{figure}[th]%
\centering
\includegraphics[
natheight=4.054200in,
natwidth=2.922200in,
height=4.0542in,
width=2.9222in
]%
{C:/Math_SwBasisFunc/InterpolSmthDev/PapersMonog/ZeroOrd/ZeroOrdDev/graphics/figExactConverg_N2_L2_samp20_5Kpt_scal_quart__12.pdf}%
\caption{Convergence of the Exact smoother.}%
\label{Fig_ExactConverg_N2_L2}%
\end{figure}

As the number of points increases the numerical smoother of $R_{0}$ is
observed to simplify to a very stable, smooth, symmetrical pattern with a
single spike at the origin and two smaller spikes at $\pm1/2$, each about
$1/3$ the central spike's height. At $\pm1$and $\pm3/2$ there are small
rounded `spikes'. The smoother between the spikes is very smooth and the trend
is zero.

The smoother of $f_{d}$ forms two spikes at the boundary $x=\pm3/2$. These
spikes are superimposed on a trend curve of amplitude which dominates the
amplitude of a stable residual error function. However the trend curve has a
larger amplitude than the boundary spikes.

The (blue) points above each smoother in Figure \ref{Fig_ExactConverg_N2_L2}
are the estimated upper bounds for the error given by \ref{7.09} or
\ref{7.08}. Clearly, in both the filtered and unfiltered cases, the
theoretical bound substantially underestimates the convergence rate.

\section{Local smoother error for data functions in $W^{1,\infty}\left(
\Omega\right)  \cap X_{w}^{0}\left(  \Omega\right)  $%
\label{Sect_Ex_smth_err_H1inf_data}}

\subsection{Introduction}

\textbf{In Subsection} \ref{SbSect_ExactSmthConverg1Dim_locW1inf_Xow} we work
in \textbf{one dimension only} and the data region is an open interval
$\Omega$.\textbf{\ }The local data functions are assumed to have bounded
derivatives i.e. functions in $W^{1,\infty}\left(  \Omega\right)  \cap
X_{w}^{0}\left(  \Omega\right)  $. In Corollary
\ref{vCor_Thm_bound_deriv_scal_hat_smth} of the Appendix, we proved, after a
lot of calculation, the following local bound estimate for the interpolant
derivative:%
\[
\left\Vert Ds_{e}\right\Vert _{\infty;\Omega}\leq\frac{2}{\lambda}\left\Vert
f\right\Vert _{\infty;\Omega}+\min\left\{  5,2+\rho N\right\}  \left\Vert
Df\right\Vert _{\infty;\Omega},\quad N\geq4,\text{ }\rho>0,
\]

where the basis function is the scaled hat function $\Lambda\left(
x/\lambda\right)  $. The hat basis function is assumed to have large support
w.r.t. the data region.

\textbf{Subsection} \ref{SbSect_notC1b_data_fns_H1inf_1dim}: Our data function
spaces will require the derivative to be bounded. This excludes functions with
cusps, vertical inflections/tangents and half-vertical tangents but allows
corners. In this subsection we describe a selection of data functions for
numerical experiments. The smoothers of these functions show instabilities.

\textbf{In Subsection} \ref{SbSect_ExactSmthConvergArbDim_locW1inf_Xow} we
consider the multivariate case but I have no example to present; I have been
unable to prove that the scaled tensor product hat function with large support
satisfies an estimate of the form \ref{a56}:
\[
\left\vert D_{k}s_{e}\left(  x\right)  \right\vert \leq c\left\Vert
f_{d}\right\Vert _{1,\infty;\Omega},\quad x\in\Omega,\text{ }\forall k.
\]

We will need the following spaces from Definition \ref{Def_SobolevSpace}:

\begin{definition}
\label{Def_H1inf_2}\textbf{The Sobolev space} $W^{1,\infty}\left(
\Omega\right)  $: If $\Omega$ is an open interval in$\mathbb{\ R}^{1}$ then
define the Sobolev space:%
\[
W^{1,\infty}\left(  \Omega\right)  =\left\{  f\in L^{\infty}\left(
\Omega\right)  :Df\in L^{\infty}\left(  \Omega\right)  \right\}  ,
\]

endowed with the supremum norm $\left\Vert f\right\Vert _{1,\infty,\Omega
}=\sum\limits_{k=0}^{1}\left\Vert D^{k}f\right\Vert _{\infty,\Omega}$.

It is also well known that if $\Omega$ is bounded then $W^{1,\infty}\left(
\Omega\right)  $ is a Banach space. Also
\[
W^{1,\infty}\left(  \Omega\right)  =\left\{  f\in C_{B}^{\left(  0\right)
}\left(  \Omega\right)  :Df\in L^{\infty}\left(  \Omega\right)  \right\}
=\left\{  f\in\mathcal{D}^{\prime}:Df\in L^{\infty}\left(  \Omega\right)
\right\}  .
\]

\end{definition}

Our data function spaces require the derivative to be bounded. This excludes
functions with cusps, vertical inflections/tangents and half-vertical tangents
but allows corners. See Subsection \ref{SbSect_notC1b_data_fns_H1inf_1dim}
below for a selection of data functions for numerical experiments.

\subsection{Estimates for the convergence of the Exact smoother in one
dimension\label{SbSect_ExactSmthConverg1Dim_locW1inf_Xow}}

Our next result is a smoothing analogue of interpolation result Theorem
\ref{Thm_err_interpol_H1inf_data}.

\begin{theorem}
\label{Thm_err_exsmth_H1inf_1dim}\textbf{1-dim smoother} Suppose:

From Theorem \ref{Thm_canon_exten_op} there exists a continuous, linear
extension operator $r_{\Omega}^{\ast}:W^{1,\infty}\left(  \Omega\right)  \cap
X_{w}^{0}\left(  \Omega\right)  \rightarrow X_{w}^{0}\left(  \mathbb{R}%
^{1}\right)  $ where $W^{1,\infty}\left(  \Omega\right)  \cap X_{w}^{0}\left(
\Omega\right)  \neq\left\{  {}\right\}  $ has norm $\left\Vert \cdot
\right\Vert _{w,0}$. Also $\left\Vert r_{\Omega}^{\ast}\right\Vert _{op}=1$.

\begin{enumerate}
\item $s_{e}:=r_{\Omega}\mathcal{S}_{X}r_{\Omega}^{\ast}f_{d}$ is the Exact
smoother of the arbitrary data function $f_{d}\in W^{1,\infty}\left(
\Omega\right)  \cap X_{w}^{0}\left(  \Omega\right)  $ on the independent data
set $X$ which is contained in the bounded open interval $\Omega\subset
\mathbb{R}^{1}$. Here $\mathcal{S}_{X}$ is the global smoothing operator.

\item The basis function satisfies $G\in W^{1,\infty}\left(  \Omega
-\Omega\right)  $.

\item There exist constants $c_{0},c_{1}\geq0$, independent of $f_{d}$, such
that%
\begin{equation}
\left\vert Ds_{e}\left(  x\right)  \right\vert \leq c_{0}\left\Vert
f_{d}\right\Vert _{\infty;\Omega}+c_{1}\left\Vert Df_{d}\right\Vert
_{\infty;\Omega},\quad x\in\Omega.\label{a1.34}%
\end{equation}

\end{enumerate}

Then if $h_{X,\Omega}=\sup\limits_{x\in\Omega}\operatorname*{dist}\left(
x,X\right)  $ is the \textbf{maximum spherical cavity size (radius) of the
data}:
\begin{align}
\left\vert f_{d}\left(  x\right)  -s_{e}\left(  x\right)  \right\vert
\leq\left\Vert f_{d}\right\Vert _{w,0,\Omega} &  \min\left\{  \sqrt{\rho
N},\frac{G\left(  0\right)  }{\left(  2\pi\right)  ^{1/2}}\right\}
+\nonumber\\
&  +\left(  c_{0}\left\Vert f_{d}\right\Vert _{\infty;\Omega}+\left(
1+c_{1}\right)  \left\Vert Df_{d}\right\Vert _{\infty;\Omega}\right)
h_{X,\Omega},\quad x\in\Omega.\label{a1.54}%
\end{align}

\end{theorem}

\begin{proof}
Let $s_{e}=\mathcal{S}_{X}r_{\Omega}^{\ast}f_{d}$. Since $G\in W^{1,\infty
}\left(  \Omega-\Omega\right)  $ the basis function formula \ref{1.271} for
the smoother means that $s_{e}\in C_{B}^{\left(  0\right)  }\cap W^{1,\infty
}\left(  \Omega\right)  $. Since $f_{d}\in X_{w}^{0}$ we have $f_{d}\in
C_{B}^{\left(  0\right)  }$ and applying Lemma
\ref{Lem_Taylor_rem1_1dim_local} with $x\in\Omega$, $z=x^{\left(  k\right)
}\in X$ and $b=x-x^{\left(  k\right)  }$ yields%
\begin{align*}
s_{e}\left(  x\right)  -f_{d}\left(  x\right)   &  =\left(  s_{e}%
-f_{d}\right)  \left(  x^{\left(  k\right)  }+\left(  x-x^{\left(  k\right)
}\right)  \right) \\
&  =\left(  s_{e}-f_{d}\right)  \left(  x^{\left(  k\right)  }\right)
+\left(  \mathcal{R}_{1}\left(  s_{e}-f_{d}\right)  \right)  \left(
x^{\left(  k\right)  },x-x^{\left(  k\right)  }\right)  .
\end{align*}

We now apply the remainder estimate \ref{a1.53}, namely
\begin{align*}
\left\vert \left(  \mathcal{R}_{1}\left(  s_{e}-f_{d}\right)  \right)  \left(
x^{\left(  k\right)  },x-x^{\left(  k\right)  }\right)  \right\vert  &
\leq\left(  \max_{y\in\left[  x^{\left(  k\right)  },x\right]  }\left\vert
D\left(  s_{e}-f_{d}\right)  (y)\right\vert \right)  \left\vert x-x^{\left(
k\right)  }\right\vert \\
&  \leq\left(  \left\Vert D\left(  s_{e}-f_{d}\right)  \right\Vert
_{\infty,\Omega}\right)  \left\vert x-x^{\left(  k\right)  }\right\vert \\
&  \leq\left(  \left\Vert Ds_{e}\right\Vert _{\infty,\Omega}+\left\Vert
Df_{d}\right\Vert _{\infty,\Omega}\right)  \left\vert x-x^{\left(  k\right)
}\right\vert ,
\end{align*}

and using \ref{a1.34} we get%
\begin{align*}
\left\vert \left(  \mathcal{R}_{1}\left(  s_{e}-f_{d}\right)  \right)  \left(
x^{\left(  k\right)  },x-x^{\left(  k\right)  }\right)  \right\vert  &
\leq\left(  c_{0}\left\Vert f_{d}\right\Vert _{\infty;\Omega}+c_{1}\left\Vert
Df_{d}\right\Vert _{\infty;\Omega}+\left\Vert Df_{d}\right\Vert _{\infty
;\Omega}\right)  \left\vert x-x^{\left(  k\right)  }\right\vert \\
& =\left(  c_{0}\left\Vert f_{d}\right\Vert _{\infty;\Omega}+\left(
1+c_{1}\right)  \left\Vert Df_{d}\right\Vert _{\infty;\Omega}\right)
\left\vert x-x^{\left(  k\right)  }\right\vert
\end{align*}

so that%
\[
\left\vert f_{d}\left(  x\right)  -s_{e}\left(  x\right)  \right\vert
\leq\left\vert f_{d}\left(  x^{\left(  i\right)  }\right)  -s_{e}\left(
x^{\left(  i\right)  }\right)  \right\vert +\left(  c_{0}\left\Vert
f_{d}\right\Vert _{\infty;\Omega}+\left(  1+c_{1}\right)  \left\Vert
Df_{d}\right\Vert _{\infty;\Omega}\right)  \left\vert x-x^{\left(  i\right)
}\right\vert ,
\]

and consequently for all $x\in\Omega$,%
\begin{align*}
\left\vert f_{d}\left(  x\right)  -s_{e}\left(  x\right)  \right\vert  &
\leq\min_{k=1}^{N}\left\vert f_{d}\left(  x^{\left(  k\right)  }\right)
-s_{e}\left(  x^{\left(  k\right)  }\right)  \right\vert +\left(
c_{0}\left\Vert f_{d}\right\Vert _{\infty;\Omega}+\left(  1+c_{1}\right)
\left\Vert Df_{d}\right\Vert _{\infty;\Omega}\right)  \min_{k=1}^{N}\left\vert
x-x^{\left(  k\right)  }\right\vert \\
& \leq\min_{k=1}^{N}\left\vert f_{d}\left(  x^{\left(  k\right)  }\right)
-s_{e}\left(  x^{\left(  k\right)  }\right)  \right\vert +\left(
c_{0}\left\Vert f_{d}\right\Vert _{\infty;\Omega}+\left(  1+c_{1}\right)
\left\Vert Df_{d}\right\Vert _{\infty;\Omega}\right)  \sup_{x\in\Omega
}\operatorname*{dist}\left(  x,X\right) \\
& =\min_{k=1}^{N}\left\vert f_{d}\left(  x^{\left(  k\right)  }\right)
-s_{e}\left(  x^{\left(  k\right)  }\right)  \right\vert +\left(
c_{0}\left\Vert f_{d}\right\Vert _{\infty;\Omega}+\left(  1+c_{1}\right)
\left\Vert Df_{d}\right\Vert _{\infty;\Omega}\right)  h_{X,\Omega}.
\end{align*}

Using \ref{7.19} and then \ref{7.35}:%
\begin{align*}
\left\vert f_{d}\left(  x^{\left(  k\right)  }\right)  -s_{e}\left(
x^{\left(  k\right)  }\right)  \right\vert  & =\left\vert r_{\Omega}^{\ast
}f_{d}\left(  x^{\left(  k\right)  }\right)  -\mathcal{S}_{X}r_{\Omega}^{\ast
}f_{d}\left(  x^{\left(  k\right)  }\right)  \right\vert \\
& \leq\left\Vert r_{\Omega}^{\ast}f_{d}\right\Vert _{w,0}\sqrt{R_{x^{\left(
k\right)  }}\left(  x^{\left(  k\right)  }\right)  -\left(  \mathcal{S}%
_{X}R_{x^{\left(  k\right)  }}\right)  \left(  x^{\left(  k\right)  }\right)
}\\
& \leq\left\Vert f_{d}\right\Vert _{w,0,\Omega}\sqrt{\rho N\delta
_{k,k}-\left(  \rho N\right)  ^{2}\left(  \left(  \rho NI+R_{X,X}\right)
^{-1}\right)  _{k,k}}\\
& \leq\left\Vert f_{d}\right\Vert _{w,0,\Omega}\sqrt{\rho N},
\end{align*}

so that we have the error estimate%
\begin{equation}
\left\vert f_{d}\left(  x\right)  -s_{e}\left(  x\right)  \right\vert
\leq\left\Vert f_{d}\right\Vert _{w,0,\Omega}\sqrt{\rho N}+\left(
c_{0}\left\Vert f_{d}\right\Vert _{\infty;\Omega}+\left(  1+c_{1}\right)
\left\Vert Df_{d}\right\Vert _{\infty;\Omega}\right)  h_{X,\Omega},\quad
x\in\Omega.\label{a1.67}%
\end{equation}

But, from Theorem \ref{Thm_canon_exten_op}, part 1 of Theorem
\ref{Thm_ord0_Riesz_rep_W2} and part 2 of Theorem \ref{Thm_ex_Min_Smth_in_Wgx}%
,%
\begin{align}
\left\vert f_{d}\left(  x^{\left(  k\right)  }\right)  -s_{e}\left(
x^{\left(  k\right)  }\right)  \right\vert =\left\vert \left(  f_{d}%
-s_{e},R_{x^{\left(  k\right)  }}^{\Omega}\right)  _{w,0,\Omega}\right\vert
&  \leq\left\Vert f_{d}-s_{e}\right\Vert _{w,0,\Omega}\left\Vert R_{x^{\left(
k\right)  }}^{\Omega}\right\Vert _{w,0,\Omega}\nonumber\\
&  =\left\Vert f_{d}-s_{e}\right\Vert _{w,0,\Omega}\left\Vert r_{\Omega}%
^{\ast}R_{x^{\left(  k\right)  }}\right\Vert _{w,0,\Omega}\nonumber\\
&  \leq\left\Vert f_{d}-s_{e}\right\Vert _{w,0,\Omega}\left\Vert R_{x^{\left(
k\right)  }}\right\Vert _{w,0}\nonumber\\
&  =\left(  2\pi\right)  ^{-1/2}G\left(  0\right)  \left\Vert f_{d}%
-s_{e}\right\Vert _{w,0,\Omega}\nonumber\\
&  =\left(  2\pi\right)  ^{-1/2}G\left(  0\right)  \left\Vert f_{d}-r_{\Omega
}\mathcal{S}_{X}r_{\Omega}^{\ast}f_{d}\right\Vert _{w,0,\Omega}\nonumber\\
&  =\left(  2\pi\right)  ^{-1/2}G\left(  0\right)  \left\Vert r_{\Omega
}r_{\Omega}^{\ast}f_{d}-r_{\Omega}\mathcal{S}_{X}r_{\Omega}^{\ast}%
f_{d}\right\Vert _{w,0,\Omega}\nonumber\\
&  =\left(  2\pi\right)  ^{-1/2}G\left(  0\right)  \left\Vert r_{\Omega
}\left(  I-\mathcal{S}_{X}\right)  r_{\Omega}^{\ast}f_{d}\right\Vert
_{w,0,\Omega}\nonumber\\
&  \leq\left(  2\pi\right)  ^{-1/2}G\left(  0\right)  \left\Vert \left(
I-\mathcal{S}_{X}\right)  r_{\Omega}^{\ast}f_{d}\right\Vert _{w,0}\nonumber\\
&  \leq\left(  2\pi\right)  ^{-1/2}G\left(  0\right)  \left\Vert r_{\Omega
}^{\ast}f_{d}\right\Vert _{w,0}\nonumber\\
&  \leq\left(  2\pi\right)  ^{-1/2}G\left(  0\right)  \left\Vert
f_{d}\right\Vert _{w,0,\Omega},\label{a1.68}%
\end{align}

and so \ref{a1.67} can be strengthened to%
\begin{align*}
\left\vert f_{d}\left(  x\right)  -s_{e}\left(  x\right)  \right\vert  & \\
& \leq\left\Vert f_{d}\right\Vert _{w,0,\Omega}\min\left\{  \sqrt{\rho
N},\frac{G\left(  0\right)  }{\left(  2\pi\right)  ^{1/2}}\right\}  +\left(
c_{0}\left\Vert f_{d}\right\Vert _{\infty;\Omega}+\left(  1+c_{1}\right)
\left\Vert Df_{d}\right\Vert _{\infty;\Omega}\right)  h_{X,\Omega},\quad
x\in\Omega.
\end{align*}

\end{proof}

\begin{remark}
\label{Rem_Thm_err_exsmth_H1inf_1dim}The estimate \ref{a1.54} is not order 1
convergence because of the factor $\sqrt{\rho N}$ in the first term.
\end{remark}

\begin{example}
\label{Ex_ex_smth_hat_larg_supp_1dim}\textbf{1-dimensional scaled hat basis
function with large support w.r.t. the data region}

From Theorem \ref{Thm_int_Xow(O)_eq_Hn(O)_dim1}, $X_{w}^{0}\left(
\Omega\right)  =W^{\mathbf{1}}\left(  \Omega\right)  $ as sets with equivalent norms.

These basis functions are discussed in Chapter
\ref{Ch_bnd_deriv_hat_smth_large_supp} of the Appendix. Large support w.r.t.
the data region $\Omega$ means that the scaled hat basis function
$\Lambda_{\lambda}\left(  x\right)  =\Lambda\left(  x/\lambda\right)  $ has
been scaled so that
\[
\operatorname*{diam}\Omega\leq\frac{1}{2}\operatorname*{diam}%
\operatorname*{supp}\Lambda_{\lambda}\text{ }i.e.\operatorname*{diam}%
\Omega\leq\lambda.
\]

This assumption has the nice consequence\ that $x,x^{\prime}\in\Omega$\ now
implies $\Lambda_{\lambda}\left(  x-x^{\prime}\right)  =1-\frac{\left\vert
x-x^{\prime}\right\vert }{\lambda}$ and there are no zero values.

In Corollary \ref{vCor_Thm_bound_deriv_scal_hat_smth} of the Appendix it is
shown that if $\operatorname*{diam}\Omega\leq\lambda$ and the data function
$f_{d}\in X_{w}^{0}$ satisfies $Df_{d}\in L^{\infty}\left(  \Omega\right)  $
then the Exact smoother $s_{e}$ corresponding to the basis function
$\Lambda_{\lambda}$ satisfies%
\[
\left\Vert Ds_{e}\right\Vert _{\infty;\Omega}\leq\frac{2}{\lambda}\left\Vert
f\right\Vert _{\infty;\Omega}+\min\left\{  5,2+\rho N\right\}  \left\Vert
Df_{d}\right\Vert _{\infty;\Omega},\quad N\geq4,\text{ }\rho>0,
\]

so that $c_{0}=\frac{2}{\lambda}$ and $c_{1}=\min\left\{  5,2+\rho N\right\}
$. Theorem \ref{Thm_err_exsmth_H1inf_1dim} now implies that when $f_{d}\in
X_{w}^{0}\left(  \Omega\right)  \cap W^{1,\infty}\left(  \Omega\right)  $:%
\begin{align}
&  \left\vert s_{e}\left(  x\right)  -f_{d}\left(  x\right)  \right\vert
\nonumber\\
&  \leq\left\Vert f_{d}\right\Vert _{w,0,\Omega}\min\left\{  \sqrt{\rho
N},\frac{G\left(  0\right)  }{\left(  2\pi\right)  ^{1/2}}\right\}  +\left(
c_{0}\left\Vert f_{d}\right\Vert _{\infty;\Omega}+\left(  1+c_{1}\right)
\left\Vert Df_{d}\right\Vert _{\infty;\Omega}\right)  h_{X,\Omega}\nonumber\\
&  =\left\Vert f_{d}\right\Vert _{w,0,\Omega}\min\left\{  \sqrt{\rho N}%
,\frac{1}{\left(  2\pi\right)  ^{1/2}}\right\}  +\left(  \frac{2}{\lambda
}\left\Vert f_{d}\right\Vert _{\infty;\Omega}+\min\left\{  6,3+\rho N\right\}
\left\Vert Df_{d}\right\Vert _{\infty;\Omega}\right)  h_{X,\Omega
}.\label{a1.15}%
\end{align}

This estimate does not represent order 1 convergence because the ugly factor
$\sqrt{\rho N}$ in the first term depends on $N$.
\end{example}

\subsection{Some non-$PWC_{B}^{\left(  1\right)  }$ data functions for
numerical experiments\label{SbSect_notC1b_data_fns_H1inf_1dim}}

Numerical experiments using the hat basis function with large support relative
to the bounded data region indicate that in 1-dimension the interpolant and
Exact smoother have difficulty handling data functions which are piecewise
$C^{\left(  1\right)  }$ except at a finite number of points where the slope
tends to $\pm\infty$ i.e. \textbf{cusps and vertical tangents/inflections and
half-vertical tangents}.

What about avoiding symmetry about the origin? Use scaling and translation: if
$a\in\left(  -1,1\right)  $ then define%
\begin{equation}
g\left(  x,a\right)  =\left\{
\begin{array}
[c]{ll}%
f\left(  \frac{x-a}{1-a}\right)  , & x\in\left[  a,1\right]  ,\\
f\left(  \frac{x-a}{1+a}\right)  , & x\in\left[  -1,a\right]  .
\end{array}
\right. \label{a8.7}%
\end{equation}

Transform a data function $f$ to a data function $g$ with zero value and
derivative at the boundary by multiplying by $\left(  1-x^{2}\right)  ^{2}$:%
\begin{equation}
g\left(  x\right)  =\left(  1-x^{2}\right)  ^{2}f\left(  x\right)
.\label{a8.8}%
\end{equation}

In the examples below we used 1,000 data points and calculated function values
on the interval $\left[  -1,1\right]  $ using 400 cells. The basis function is
the hat function with support $\left[  -1,1\right]  $.

\begin{example}
\label{Ex_vert_cusp}\textbf{A vertical cusp at the origin} Here we use:%
\[
f_{c}\left(  x;s\right)  =1-\sqrt{\left\vert x\right\vert ^{s}\left(
2-\left\vert x\right\vert ^{s}\right)  },\quad s\geq0,\text{ }x\in
\Omega=\left(  -1,1\right)  .
\]

Hence $f_{c}\left(  x;s\right)  =0$ when $x=\pm1$. But%
\begin{align*}
Df_{c}\left(  x;s\right)   & =\frac{-1}{2\sqrt{\left\vert x\right\vert
^{s}\left(  2-\left\vert x\right\vert ^{s}\right)  }}D\left(  \left\vert
x\right\vert ^{s}\left(  2-\left\vert x\right\vert ^{s}\right)  \right) \\
& =\frac{-1}{2\sqrt{\left\vert x\right\vert ^{s}\left(  2-\left\vert
x\right\vert ^{s}\right)  }}\left(  \left(  D\left\vert x\right\vert
^{s}\right)  \left(  2-\left\vert x\right\vert ^{s}\right)  -\left\vert
x\right\vert ^{s}D\left\vert x\right\vert ^{s}\right) \\
& =-\frac{2-\left\vert x\right\vert ^{s}-\left\vert x\right\vert ^{s}}%
{2\sqrt{\left\vert x\right\vert ^{s}\left(  2-\left\vert x\right\vert
^{s}\right)  }}D\left\vert x\right\vert ^{s}\\
& =\frac{1}{2}\left(  -\frac{2-\left\vert x\right\vert ^{s}}{\sqrt{\left\vert
x\right\vert ^{s}\left(  2-\left\vert x\right\vert ^{s}\right)  }}%
+\frac{\left\vert x\right\vert ^{s}}{\sqrt{\left\vert x\right\vert ^{s}\left(
2-\left\vert x\right\vert ^{s}\right)  }}\right)  D\left\vert x\right\vert
^{s}\\
& =\frac{1}{2}\left(  -\frac{\sqrt{2-\left\vert x\right\vert ^{s}}}%
{\sqrt{\left\vert x\right\vert ^{s}}}+\frac{\sqrt{\left\vert x\right\vert
^{s}}}{\sqrt{2-\left\vert x\right\vert ^{s}}}\right)  D\left\vert x\right\vert
^{s}\\
& =\frac{s}{2}\operatorname*{sgn}\left(  x\right)  \left(  -\frac
{\sqrt{2-\left\vert x\right\vert ^{s}}}{\sqrt{\left\vert x\right\vert ^{s}}%
}+\frac{\sqrt{\left\vert x\right\vert ^{s}}}{\sqrt{2-\left\vert x\right\vert
^{s}}}\right)  \left\vert x\right\vert ^{s-1}\\
& =\frac{s}{2}\operatorname*{sgn}\left(  x\right)  \left(  -\left\vert
x\right\vert ^{\frac{s}{2}-1}\sqrt{2-\left\vert x\right\vert ^{s}}%
+\frac{\left\vert x\right\vert ^{\frac{3s}{2}-1}}{\sqrt{2-\left\vert
x\right\vert ^{s}}}\right) \\
& =\frac{s}{2}\operatorname*{sgn}\left(  x\right)  \frac{-\left\vert
x\right\vert ^{\frac{s}{2}-1}\left(  2-\left\vert x\right\vert ^{s}\right)
+\left\vert x\right\vert ^{\frac{3s}{2}-1}}{\sqrt{2-\left\vert x\right\vert
^{s}}}\\
& =\frac{s}{2}\operatorname*{sgn}\left(  x\right)  \frac{-\left\vert
x\right\vert ^{\frac{s}{2}-1}\left(  2-\left\vert x\right\vert ^{s}\right)
+\left\vert x\right\vert ^{\frac{3s}{2}-1}}{\sqrt{2-\left\vert x\right\vert
^{s}}}\\
& =\frac{s}{2}\operatorname*{sgn}\left(  x\right)  \frac{-2\left\vert
x\right\vert ^{\frac{s}{2}-1}+2\left\vert x\right\vert ^{\frac{3s}{2}-1}%
}{\sqrt{2-\left\vert x\right\vert ^{s}}}\\
& =s\frac{\operatorname*{sgn}\left(  x\right)  }{\sqrt{2-\left\vert
x\right\vert ^{s}}}\left(  -\left\vert x\right\vert ^{\frac{s}{2}%
-1}+\left\vert x\right\vert ^{\frac{3s}{2}-1}\right)  ,
\end{align*}

so that $Df_{c}\left(  x;s\right)  =0$ when $x=\pm1$. Also $Df_{c}\in
L^{2}\left(  \Omega\right)  $ and $\left\vert Df_{c}\left(  0;s\right)
\right\vert =\infty$ iff $1<s<2$.

The properties $f_{c}\left(  \pm1;s\right)  =0$ and $Df_{c}\left(
\pm1;s\right)  =0$ minimize the contribution of the boundary to the
interpolant error since we want to concentrate on the error due to the cusp.

Numerical experiments with the case $s=3/2$ show a large error spike at the origin.

Use asymmetric technique \ref{a8.7} with $a=1/3$?\medskip

\textbf{1}) Asymmetry reduces the spike size?

\textbf{2}) More cusps reduce the spike size?

Measure error spike by max/median? Graph against number of data points?

Several cusps vs. one cusp regarding the error spikes.

Symmetric vs. asymmetric regarding the error spikes.
\end{example}

\begin{example}
\label{Ex_vert_tang}\textbf{A vertical tangent at the origin} Here we use the
functions:%
\[
f_{v}\left(  x;s\right)  =\operatorname*{sgn}\left(  x\right)  \left\vert
x\right\vert ^{s}\sqrt{\left\vert x\right\vert ^{s}\left(  2-\left\vert
x\right\vert ^{s}\right)  }\cos^{2}\frac{\pi}{2}x,\quad s\geq0,\text{ }%
x\in\Omega=\left(  -1,1\right)  .
\]

Here $f_{v}\left(  0;s\right)  =0$. The multiplicative term $\cos^{2}\frac
{\pi}{2}x$ is designed to reduce the contribution of the boundary to the
interpolant error since we want to concentrate on the error due to the
vertical tangent. In fact $f_{v}\left(  \pm1;s\right)  =\left\vert
Df_{v}\left(  \pm1;s\right)  \right\vert =0$. Also $\max\limits_{\overline
{\Omega}}f_{v}\lesssim0.482$.

Similarly to Case 1, $Df_{v}\in L^{2}\left(  \Omega\right)  $ and
$Df_{v}\left(  0;s\right)  =\infty$ iff $1<s<2$.

Numerical experiments with Example \ref{Ex_vert_cusp}: $s=1.1$ show a large
error spike at the origin.

Observations.

Measure error spike by ratio max/median? Graph against number of data points?
\end{example}

\begin{example}
\textbf{A half-vertical tangent at the origin} Here we use the functions:%
\[
f_{h}\left(  x;s\right)  =\left\{
\begin{array}
[c]{ll}%
f_{v}\left(  x;s\right)  , & x\in\left(  -1,0\right]  ,\\
\cos^{2}\frac{\pi}{2}x, & x\in\left[  0,1\right)  .
\end{array}
\right.
\]

Similarly to Case 1, $Df_{v}\in L^{2}\left(  \Omega\right)  $ and
$Df_{v}\left(  0;s\right)  =\infty$ iff $1<s<2$.

Use $s=1.25$.

?? Measure error spike by ratio max/median? Graph against number of data
points?\medskip
\end{example}

Comparing Example \ref{Ex_vert_cusp} with Example \ref{Ex_vert_tang} -
different values of $s$ - no meaningful comment.

\begin{example}
\textbf{A corner at the origin} We use the data function%
\[
f_{d}\left(  x\right)  =(1-\left\vert x\right\vert )^{2}\left(  \left(
m+2\right)  \left\vert x\right\vert +1\right)  ,\quad m=-1.
\]

Here $m$ is the right slope of $f_{d}$ at the origin.

Observations:

1) Spikes in the error at the boundary - width.

2) ??
\end{example}

\begin{example}
\textbf{The functions }$f_{d}$\textbf{\ and }$Df_{d}$\textbf{\ are zero at the
boundary} and $f_{d}\in C^{\left(  2\right)  }\left(  \overline{\Omega
}\right)  $. We use the data functions%
\begin{align*}
f_{d}\left(  x\right)   & =1+\cos\left(  \pi x\right)  ,\\
f_{d}\left(  x\right)   & =\left(  x^{2}-1\right)  ^{2}e^{-x^{2}}.
\end{align*}

Observations: No spikes in the error.
\end{example}

\subsection{Multivariate estimates for convergence of Exact
smoother\label{SbSect_ExactSmthConvergArbDim_locW1inf_Xow}}

A multivariate version of Theorem \ref{Thm_err_exsmth_H1inf_1dim} is:

\begin{theorem}
\label{Thm_err_exsmth_H1inf_multivar}\textbf{Multivariate smoother} From
Theorem \ref{Thm_canon_exten_op} there exists a continuous, linear extension
operator $r_{\Omega}^{\ast}:W^{1,\infty}\left(  \Omega\right)  \cap X_{w}%
^{0}\left(  \Omega\right)  \rightarrow X_{w}^{0}$ where $W^{1,\infty}\left(
\Omega\right)  \cap X_{w}^{0}\left(  \Omega\right)  \neq\left\{  {}\right\}  $
has norm $\left\Vert \cdot\right\Vert _{w,0}$. Also $\left\Vert r_{\Omega
}^{\ast}\right\Vert _{op}=1$.

Now suppose:

\begin{enumerate}
\item $s_{e}:=\mathcal{S}_{X}r_{\Omega}^{\ast}f_{d}$ is the Exact smoother of
the arbitrary data function $f_{d}\in W^{1,\infty}\left(  \Omega\right)  \cap
X_{w}^{0}\left(  \Omega\right)  $ on the independent data set $X$ which is
contained in the bounded region $\Omega$ which satisfies the (uniform) cone condition.

\item The basis function $G$ lies in $W^{1,\infty}\left(  \Omega
-\Omega\right)  $.

\item There exist a constant $c>0$, independent of $f_{d}\in W^{1,\infty
}\left(  \Omega\right)  \cap X_{w}^{0}\left(  \Omega\right)  $, such that%
\begin{equation}
\left\vert D_{k}s_{e}\left(  x\right)  \right\vert \leq c\left\Vert
f_{d}\right\Vert _{1,\infty;\Omega},\quad x\in\Omega,\text{ }\forall
k.\label{a56}%
\end{equation}

\end{enumerate}

Then if $h_{X,\Omega}=\sup\limits_{x\in\Omega}\operatorname*{dist}\left(
x,X\right)  $ is the \textbf{maximum spherical cavity size (radius) of the
data} there exist constants $c_{\Omega}$,$\varepsilon_{0}>0$ such that
\begin{align}
\left\vert f_{d}\left(  x\right)  -s_{e}\left(  x\right)  \right\vert  &
\leq\nonumber\\
&  \left\Vert f_{d}\right\Vert _{w,0,\Omega}\min\left\{  \sqrt{\rho N}%
,\frac{G\left(  0\right)  }{\left(  2\pi\right)  ^{1/2}}\right\}  +\sqrt
{d}\left(  1+c\right)  \left\Vert f_{d}\right\Vert _{1,\infty,\Omega}%
c_{\Omega}h_{\Omega,X},\quad x\in\Omega,\label{a59}%
\end{align}

when $h_{X,\Omega}<\varepsilon_{0}$.
\end{theorem}

\begin{proof}
We will adapt the proof of the multivariate interpolation result Theorem
\ref{Thm_err_interpol_H1inf_data multivar} and the univariate smoothing result
of Theorem \ref{Thm_err_exsmth_H1inf_1dim}.

Since $f_{d}\in X_{w}^{0}\left(  \Omega\right)  $ we have $f_{d}\in
C_{B}^{\left(  0\right)  }\left(  \Omega\right)  $. If $x\in\Omega$ then by
Lemma \ref{Lem_cone_x_xk} there exist constants $r_{\Omega},\varepsilon_{0}>0$
such that when $h_{\Omega,X}<\varepsilon_{0}$ there exists $x^{\left(
k\right)  }\in X$ such that $\left[  x,x^{\left(  k\right)  }\right]
\subset\Omega$ and $\left\vert x-x^{\left(  k\right)  }\right\vert <c_{\Omega
}h_{\Omega,X}$. This permits us to apply the Taylor series expansion result of
Remark \ref{Rem_Lem_Taylor_extension} with $x\in\Omega$, $z=x^{\left(
k\right)  }\in X$ and $b=x-x^{\left(  k\right)  }$ to get%
\begin{align*}
s_{e}\left(  x\right)  -f_{d}\left(  x\right)   &  =\left(  s_{e}%
-f_{d}\right)  \left(  x^{\left(  k\right)  }+\left(  x-x^{\left(  k\right)
}\right)  \right) \\
&  =s_{e}\left(  x^{\left(  k\right)  }\right)  -f_{d}\left(  x^{\left(
k\right)  }\right)  +\left(  \mathcal{R}_{1}\left(  s_{e}-f_{d}\right)
\right)  \left(  x^{\left(  k\right)  },x-x^{\left(  k\right)  }\right)  ,
\end{align*}

and from Remark \ref{Rem_Lem_Taylor_extension} we have the remainder estimate
\ref{1.18}, namely
\begin{align*}
\left\vert \left(  \mathcal{R}_{1}\left(  s_{e}-f_{d}\right)  \right)  \left(
x^{\left(  k\right)  },x-x^{\left(  k\right)  }\right)  \right\vert  &
\leq\sqrt{d}\left(  \max_{\left\vert \alpha\right\vert =1}\max_{y\in\left[
x^{\left(  k\right)  },x\right]  }\left\vert D^{\alpha}\left(  s_{e}%
-f_{d}\right)  (y)\right\vert \right)  \left\vert x-x^{\left(  k\right)
}\right\vert \\
&  \leq\sqrt{d}\max_{\left\vert \alpha\right\vert =1}\left\Vert D^{\alpha
}\left(  s_{e}-f_{d}\right)  \right\Vert _{\infty,\Omega}\left\vert
x-x^{\left(  k\right)  }\right\vert \\
&  \leq\sqrt{d}\left(  \max_{\left\vert \alpha\right\vert =1}\left\Vert
D^{\alpha}f_{d}\right\Vert _{\infty,\Omega}+\max_{\left\vert \alpha\right\vert
=1}\left\Vert D^{\alpha}s_{e}\right\Vert _{\infty,\Omega}\right)  \left\vert
x-x^{\left(  k\right)  }\right\vert .
\end{align*}

Hence, using the assumption \ref{a56},%
\begin{align}
\left\vert s_{e}\left(  x\right)  -f_{d}\left(  x\right)  \right\vert  &
\leq\left\vert s_{e}\left(  x^{\left(  k\right)  }\right)  -f_{d}\left(
x^{\left(  k\right)  }\right)  \right\vert +\sqrt{d}\left(  \max_{\left\vert
\alpha\right\vert =1}\left\Vert D^{\alpha}f_{d}\right\Vert _{\infty,\Omega
}+c\left\Vert f_{d}\right\Vert _{1,\infty,\Omega}\right)  \left\vert
x-x^{\left(  k\right)  }\right\vert \nonumber\\
& \leq\left\vert s_{e}\left(  x^{\left(  k\right)  }\right)  -f_{d}\left(
x^{\left(  k\right)  }\right)  \right\vert +\sqrt{d}\left(  1+c\right)
\left\Vert f_{d}\right\Vert _{1,\infty,\Omega}\left\vert x-x^{\left(
k\right)  }\right\vert ,\label{a55}%
\end{align}

and thus%
\begin{align*}
\left\vert s_{e}\left(  x\right)  -f_{d}\left(  x\right)  \right\vert  &
\leq\min_{k=1}^{N}\left\vert s_{e}\left(  x^{\left(  k\right)  }\right)
-f_{d}\left(  x^{\left(  k\right)  }\right)  \right\vert +\sqrt{d}\left(
1+c\right)  \left\Vert f_{d}\right\Vert _{1,\infty,\Omega}\min_{k=1}%
^{N}\left\vert x-x^{\left(  k\right)  }\right\vert \\
& \leq\min_{k=1}^{N}\left\vert s_{e}\left(  x^{\left(  k\right)  }\right)
-f_{d}\left(  x^{\left(  k\right)  }\right)  \right\vert +\sqrt{d}\left(
1+c\right)  \left\Vert f_{d}\right\Vert _{1,\infty,\Omega}\sup_{x\in\Omega
}\min_{k=1}^{N}\left\vert x-x^{\left(  k\right)  }\right\vert \\
& \leq\min_{k=1}^{N}\left\vert s_{e}\left(  x^{\left(  k\right)  }\right)
-f_{d}\left(  x^{\left(  k\right)  }\right)  \right\vert +\sqrt{d}\left(
1+c\right)  \left\Vert f_{d}\right\Vert _{1,\infty,\Omega}c_{\Omega}%
h_{\Omega,X}.
\end{align*}

Using \ref{7.19} and then \ref{7.35}:%
\begin{align*}
\left\vert s_{e}\left(  x^{\left(  k\right)  }\right)  -f_{d}\left(
x^{\left(  k\right)  }\right)  \right\vert  & =\left\vert \mathcal{S}%
_{X}r_{\Omega}^{\ast}f_{d}\left(  x^{\left(  k\right)  }\right)  -r_{\Omega
}^{\ast}f_{d}\left(  x^{\left(  k\right)  }\right)  \right\vert \\
& =\left\vert \left(  \mathcal{S}_{X}-I\right)  r_{\Omega}^{\ast}f_{d}\left(
x^{\left(  k\right)  }\right)  \right\vert \\
& \leq\left\Vert r_{\Omega}^{\ast}f_{d}\right\Vert _{w,0}\sqrt{R_{x^{\left(
k\right)  }}\left(  x^{\left(  k\right)  }\right)  -\left(  \mathcal{S}%
_{X}R_{x^{\left(  k\right)  }}\right)  \left(  x^{\left(  k\right)  }\right)
}\\
& \leq\left\Vert f_{d}\right\Vert _{w,0,\Omega}\sqrt{\rho N\delta
_{k,k}-\left(  \rho N\right)  ^{2}\left(  \left(  \rho NI+R_{X,X}\right)
^{-1}\right)  _{k,k}}\\
& \leq\left\Vert f_{d}\right\Vert _{w,0,\Omega}\sqrt{\rho N},
\end{align*}

so that we have the error estimate%
\begin{equation}
\left\vert s_{e}\left(  x\right)  -f_{d}\left(  x\right)  \right\vert
\leq\left\Vert f_{d}\right\Vert _{w,0,\Omega}\sqrt{\rho N}+\sqrt{d}\left(
1+c\right)  \left\Vert f_{d}\right\Vert _{1,\infty,\Omega}c_{\Omega}%
h_{\Omega,X},\quad x\in\Omega.\label{a57}%
\end{equation}

But, from Theorem \ref{Thm_canon_exten_op}, part 1 of Theorem
\ref{Thm_ord0_Riesz_rep_W2} and part 2 of Theorem \ref{Thm_ex_Min_Smth_in_Wgx}%
,%
\begin{align*}
\left\vert f_{d}\left(  x^{\left(  k\right)  }\right)  -s_{e}\left(
x^{\left(  k\right)  }\right)  \right\vert =\left\vert \left(  f_{d}%
-s_{e},R_{x^{\left(  k\right)  }}^{\Omega}\right)  _{w,0,\Omega}\right\vert
&  \leq\left\Vert f_{d}-s_{e}\right\Vert _{w,0,\Omega}\left\Vert R_{x^{\left(
k\right)  }}^{\Omega}\right\Vert _{w,0,\Omega}\\
&  =\left\Vert f_{d}-s_{e}\right\Vert _{w,0,\Omega}\left\Vert r_{\Omega}%
^{\ast}R_{x^{\left(  k\right)  }}\right\Vert _{w,0,\Omega}\\
&  \leq\left\Vert f_{d}-s_{e}\right\Vert _{w,0,\Omega}\left\Vert R_{x^{\left(
k\right)  }}\right\Vert _{w,0}\\
&  =\left(  2\pi\right)  ^{-1/2}G\left(  0\right)  \left\Vert f_{d}%
-s_{e}\right\Vert _{w,0,\Omega}\\
&  =\left(  2\pi\right)  ^{-1/2}G\left(  0\right)  \left\Vert f_{d}-r_{\Omega
}\mathcal{S}_{X}r_{\Omega}^{\ast}f_{d}\right\Vert _{w,0,\Omega}\\
&  =\left(  2\pi\right)  ^{-1/2}G\left(  0\right)  \left\Vert r_{\Omega
}r_{\Omega}^{\ast}f_{d}-r_{\Omega}\mathcal{S}_{X}r_{\Omega}^{\ast}%
f_{d}\right\Vert _{w,0,\Omega}\\
&  =\left(  2\pi\right)  ^{-1/2}G\left(  0\right)  \left\Vert r_{\Omega
}\left(  I-\mathcal{S}_{X}\right)  r_{\Omega}^{\ast}f_{d}\right\Vert
_{w,0,\Omega}\\
&  \leq\left(  2\pi\right)  ^{-1/2}G\left(  0\right)  \left\Vert \left(
I-\mathcal{S}_{X}\right)  r_{\Omega}^{\ast}f_{d}\right\Vert _{w,0}\\
&  \leq\left(  2\pi\right)  ^{-1/2}G\left(  0\right)  \left\Vert r_{\Omega
}^{\ast}f_{d}\right\Vert _{w,0}\\
&  \leq\left(  2\pi\right)  ^{-1/2}G\left(  0\right)  \left\Vert
f_{d}\right\Vert _{w,0,\Omega},
\end{align*}

and so \ref{a57} can be strengthened to%
\begin{align*}
\left\vert f_{d}\left(  x\right)  -s_{e}\left(  x\right)  \right\vert  & \\
&  \leq\left\Vert f_{d}\right\Vert _{w,0,\Omega}\min\left\{  \sqrt{\rho
N},\frac{G\left(  0\right)  }{\left(  2\pi\right)  ^{1/2}}\right\}  +\sqrt
{d}\left(  1+c\right)  \left\Vert f_{d}\right\Vert _{1,\infty,\Omega}%
c_{\Omega}h_{\Omega,X},\quad x\in\Omega,
\end{align*}

as claimed.
\end{proof}

\section{Approximation of the Exact smoother\label{Sect_ex_approx_Exact_smth}}

We shall finish this chapter by approximating the Exact smoother by an easily
calculated member of $W_{G,X}$ which does not involve calculating the inverse
of a matrix. This result will be used to motivate the derivation of the
Approximate smoother in Subsection \ref{SbSect_ap_deriv_Approx_smth} of the
next chapter.

The relevant properties of the vector-valued evaluation operator
$\widetilde{\mathcal{E}}_{X}$ and its adjoint $\widetilde{\mathcal{E}}%
_{X}^{\ast}$ w.r.t. $X_{w}^{0}$ were proven in Theorem
\ref{Thm_eval_op_properties}. Now suppose $X=\left\{  x^{\left(  k\right)
}\right\}  _{k=1}^{N}$ is an ordered set of points in $\mathbb{R}^{d}$ and
$f\in X_{w}^{0}$. Define $g\in X_{w}^{0}$ by%
\begin{equation}
f=\frac{1}{NR_{0}\left(  0\right)  }\widetilde{\mathcal{E}}_{X}^{\ast
}\widetilde{\mathcal{E}}_{X}f+g=\frac{1}{NR_{0}\left(  0\right)  }\sum
_{k}f\left(  x^{\left(  k\right)  }\right)  R_{x^{\left(  k\right)  }%
}+g,\label{15.01}%
\end{equation}

so that
\begin{equation}
f-g=\frac{1}{NR_{0}\left(  0\right)  }\widetilde{\mathcal{E}}_{X}^{\ast
}\widetilde{\mathcal{E}}_{X}f,\label{1.65}%
\end{equation}

has the following properties:

\begin{lemma}
\label{Lem_ex_properties_f_minus_g}Suppose $f,g$ satisfy equation \ref{15.01}. Then:

\begin{enumerate}
\item $\left\Vert f-g\right\Vert _{w,0}\leq\left\Vert f\right\Vert _{w,0}$.

\item $\left\vert f\left(  x\right)  -g\left(  x\right)  \right\vert \leq
\sqrt{R_{0}\left(  0\right)  }\left\Vert f\right\Vert _{w,0}$.

\item $\left\vert \widetilde{\mathcal{E}}_{X}f-\widetilde{\mathcal{E}}%
_{X}g\right\vert \leq\frac{1}{\sqrt{NR_{0}\left(  0\right)  }}\left\vert
\widetilde{\mathcal{E}}_{X}f\right\vert \leq\left\Vert f\right\Vert _{w,0}$.

\item The Exact smoother operator $\mathcal{S}_{X}$ satisfies: $\left\vert
\mathcal{S}_{X}f\left(  x\right)  -\mathcal{S}_{X}g\left(  x\right)
\right\vert \leq\sqrt{R_{0}\left(  0\right)  }\left\Vert f\right\Vert _{w,0}$.
\end{enumerate}
\end{lemma}

\begin{proof}
By definition of the Riesz representer
\begin{equation}
\left\Vert R_{x}\right\Vert _{w,0}=\left\Vert R_{0}\right\Vert _{w,0}%
=\sqrt{R_{0}\left(  0\right)  }\label{1.47}%
\end{equation}

so part 3 of Theorem \ref{Thm_eval_op_properties} implies%
\begin{equation}
\left\Vert \widetilde{\mathcal{E}}_{X}^{\ast}\right\Vert =\left\vert
\widetilde{\mathcal{E}}_{X}\right\vert =\left\Vert R_{X,X}\right\Vert
\leq\sqrt{N}\sqrt{R_{0}\left(  0\right)  }.\label{1.48}%
\end{equation}

Thus:\medskip

\textbf{Part 1}%
\[
\left\Vert f-g\right\Vert _{w,0}=\frac{1}{NR_{0}\left(  0\right)  }\left\Vert
\widetilde{\mathcal{E}}_{X}^{\ast}\widetilde{\mathcal{E}}_{X}f\right\Vert
_{w,0}\leq\left\Vert f\right\Vert _{w,0}.
\]
\medskip

\textbf{Part 2}
\begin{align*}
f\left(  x\right)  -g\left(  x\right)  =\frac{1}{NR_{0}\left(  0\right)
}\left\vert \left(  \widetilde{\mathcal{E}}_{X}^{\ast}\widetilde{\mathcal{E}%
}_{X}f,R_{x}\right)  _{w,0}\right\vert  &  =\frac{1}{NR_{0}\left(  0\right)
}\left\Vert \widetilde{\mathcal{E}}_{X}^{\ast}\widetilde{\mathcal{E}}%
_{X}f\right\Vert _{w,0}\left\Vert R_{x}\right\Vert _{w,0}\\
&  \leq\sqrt{R_{0}\left(  0\right)  }\left\Vert f\right\Vert _{w,0}.
\end{align*}
\medskip

\textbf{Part 3} By part 5 of Theorem \ref{Thm_eval_op_properties},
$\widetilde{\mathcal{E}}_{X}\widetilde{\mathcal{E}}_{X}^{\ast}=R_{X,X}$ so
\begin{align*}
\left\vert \widetilde{\mathcal{E}}_{X}f-\widetilde{\mathcal{E}}_{X}%
g\right\vert =\frac{1}{NR_{0}\left(  0\right)  }\left\vert R_{X,X}%
\widetilde{\mathcal{E}}_{X}f\right\vert \leq\frac{1}{NR_{0}\left(  0\right)
}\left\Vert R_{X,X}\right\Vert \left\vert \widetilde{\mathcal{E}}%
_{X}f\right\vert  &  =\frac{1}{\sqrt{NR_{0}\left(  0\right)  }}\left\vert
\widetilde{\mathcal{E}}_{X}f\right\vert \\
&  \leq\left\Vert f\right\Vert _{w,0}.
\end{align*}

\textbf{Part 4} From part 2 of Theorem \ref{Thm_ex_Min_Smth_in_Wgx}, the Exact
smoother operator satisfies, $\left\Vert \mathcal{S}_{X}\right\Vert \leq1$.
Hence
\begin{align*}
\left\vert \mathcal{S}_{X}f\left(  x\right)  -\mathcal{S}_{X}g\left(
x\right)  \right\vert =\frac{1}{NR_{0}\left(  0\right)  }\left\vert \left(
\mathcal{S}_{X}\widetilde{\mathcal{E}}_{X}^{\ast}\widetilde{\mathcal{E}}%
_{X}f,R_{x}\right)  \right\vert  &  \leq\frac{1}{NR_{0}\left(  0\right)
}\left\Vert \mathcal{S}_{X}\widetilde{\mathcal{E}}_{X}^{\ast}%
\widetilde{\mathcal{E}}_{X}f\right\Vert _{w,0}\left\Vert R_{0}\right\Vert
_{w,0}\\
&  \leq\frac{1}{N\sqrt{R_{0}\left(  0\right)  }}\left\Vert
\widetilde{\mathcal{E}}_{X}^{\ast}\widetilde{\mathcal{E}}_{X}f\right\Vert
_{w,0}\\
&  \leq\sqrt{R_{0}\left(  0\right)  }\left\Vert f\right\Vert _{w,0}.
\end{align*}

\end{proof}

We now prove our smoother approximation theorem:

\begin{theorem}
\label{Thm_estim_Exact_smth_by_Wgx}Suppose $f,g\in X_{w}^{0}$ satisfy equation
\ref{15.01}. Then%
\begin{equation}
\mathcal{S}_{X}f\left(  x\right)  -\frac{\left(  \widetilde{\mathcal{E}}%
_{X}R_{x}\right)  ^{T}\widetilde{\mathcal{E}}_{X}f}{\left(  R_{0}\left(
0\right)  +\rho\right)  N}=\frac{R_{0}\left(  0\right)  }{R_{0}\left(
0\right)  +\rho}\mathcal{S}_{X}g\left(  x\right)  ,\quad x\in\mathbb{R}%
^{d},\label{15.19}%
\end{equation}

and%
\begin{equation}
\left\vert \mathcal{S}_{X}f\left(  x\right)  -\frac{\left(
\widetilde{\mathcal{E}}_{X}R_{x}\right)  ^{T}\widetilde{\mathcal{E}}_{X}%
f}{\left(  R_{0}\left(  0\right)  +\rho\right)  N}\right\vert \leq\frac
{2R_{0}\left(  0\right)  ^{3/2}}{R_{0}\left(  0\right)  +\rho}\min\left\{
1,\frac{R_{0}\left(  0\right)  }{\rho}\right\}  \left\Vert f\right\Vert
_{w,0},\label{15.20}%
\end{equation}

and%
\begin{equation}
\left\Vert \mathcal{S}_{X}f-\frac{\left(  \widetilde{\mathcal{E}}_{X}%
R_{x}\right)  ^{T}\widetilde{\mathcal{E}}_{X}f}{\left(  R_{0}\left(  0\right)
+\rho\right)  N}\right\Vert _{w,0}\leq\frac{2R_{0}\left(  0\right)  }%
{R_{0}\left(  0\right)  +\rho}\left\Vert f\right\Vert _{w,0}.\label{15.02}%
\end{equation}

\end{theorem}

\begin{proof}
Starting with \ref{7.28} and using the fact that $\widetilde{\mathcal{E}}%
_{X}\widetilde{\mathcal{E}}_{X}^{\ast}=R_{X,X}$ (part 5 Theorem
\ref{Thm_eval_op_properties}) we have%
\begin{align}
\mathcal{S}_{X}f\left(  x\right)   & =\left(  \widetilde{\mathcal{E}}_{X}%
R_{x}\right)  ^{T}\left(  \rho NI+R_{X,X}\right)  ^{-1}\widetilde{\mathcal{E}%
}_{X}f\nonumber\\
& =\left(  \widetilde{\mathcal{E}}_{X}R_{x}\right)  ^{T}\left(  \rho
NI+R_{X,X}\right)  ^{-1}\left(  \frac{1}{NR_{0}\left(  0\right)  }%
R_{X,X}\widetilde{\mathcal{E}}_{X}f+\widetilde{\mathcal{E}}_{X}g\right)
\nonumber\\
& =\frac{1}{NR_{0}\left(  0\right)  }\left(  \widetilde{\mathcal{E}}_{X}%
R_{x}\right)  ^{T}\left(  \rho NI+R_{X,X}\right)  ^{-1}R_{X,X}%
\widetilde{\mathcal{E}}_{X}f+\left(  \widetilde{\mathcal{E}}_{X}R_{x}\right)
^{T}\left(  \rho NI+R_{X,X}\right)  ^{-1}\widetilde{\mathcal{E}}%
_{X}g\nonumber\\
& =\frac{1}{NR_{0}\left(  0\right)  }\left(  \widetilde{\mathcal{E}}_{X}%
R_{x}\right)  ^{T}\left(  \rho NI+R_{X,X}\right)  ^{-1}R_{X,X}%
\widetilde{\mathcal{E}}_{X}f+\mathcal{S}_{X}g\left(  x\right) \nonumber\\
& =\frac{1}{NR_{0}\left(  0\right)  }\left(  \widetilde{\mathcal{E}}_{X}%
R_{x}\right)  ^{T}\left(  \rho NI+R_{X,X}\right)  ^{-1}\left(  \rho
NI+R_{X,X}\right)  \widetilde{\mathcal{E}}_{X}f-\nonumber\\
& \qquad-\frac{1}{NR_{0}\left(  0\right)  }\left(  \widetilde{\mathcal{E}}%
_{X}R_{x}\right)  ^{T}\left(  \rho NI+R_{X,X}\right)  ^{-1}\left(  \rho
N\widetilde{\mathcal{E}}_{X}f\right)  +\mathcal{S}_{X}g\left(  x\right)
\nonumber\\
& =\frac{\left(  \widetilde{\mathcal{E}}_{X}R_{x}\right)  ^{T}%
\widetilde{\mathcal{E}}_{X}f}{NR_{0}\left(  0\right)  }-\frac{\rho}%
{R_{0}\left(  0\right)  }\mathcal{S}_{X}f\left(  x\right)  +\mathcal{S}%
_{X}g\left(  x\right)  ,\label{15.15}%
\end{align}

and solving for $\mathcal{S}_{X}f\left(  x\right)  $ gives%
\begin{align*}
\mathcal{S}_{X}f\left(  x\right)   & =\frac{1}{1+\frac{\rho}{R_{0}\left(
0\right)  }}\frac{\left(  \widetilde{\mathcal{E}}_{X}R_{x}\right)
^{T}\widetilde{\mathcal{E}}_{X}f}{NR_{0}\left(  0\right)  }+\frac
{\mathcal{S}_{X}g\left(  x\right)  }{1+\frac{\rho}{R_{0}\left(  0\right)  }}\\
& =\frac{\left(  \widetilde{\mathcal{E}}_{X}R_{x}\right)  ^{T}%
\widetilde{\mathcal{E}}_{X}f}{\left(  R_{0}\left(  0\right)  +\rho\right)
N}+\frac{R_{0}\left(  0\right)  }{R_{0}\left(  0\right)  +\rho}\mathcal{S}%
_{X}g\left(  x\right)  ,
\end{align*}

i.e.%
\[
\mathcal{S}_{X}f\left(  x\right)  -\frac{\left(  \widetilde{\mathcal{E}}%
_{X}R_{x}\right)  ^{T}\widetilde{\mathcal{E}}_{X}f}{\left(  R_{0}\left(
0\right)  +\rho\right)  N}=\frac{R_{0}\left(  0\right)  }{R_{0}\left(
0\right)  +\rho}\mathcal{S}_{X}g\left(  x\right)  ,
\]

which is equation \ref{15.19}. To prove the estimate \ref{15.20} we will
estimate $\left\vert \mathcal{S}_{X}g\left(  x\right)  \right\vert $ using the
inequality \ref{7.68}:%
\begin{equation}
\left\vert \mathcal{S}_{X}u\left(  x\right)  \right\vert \leq\sqrt
{R_{0}\left(  0\right)  }\min\left\{  1,\frac{R_{0}\left(  0\right)  }{\rho
}\right\}  \left\Vert u\right\Vert _{w,0},\quad u\in X_{w}^{0}.\label{1.63}%
\end{equation}

Applying the Exact smoother to equation \ref{1.65} gives%
\[
\mathcal{S}_{X}g\left(  x\right)  =\mathcal{S}_{X}f\left(  x\right)  -\frac
{1}{NR_{0}\left(  0\right)  }\mathcal{S}_{X}\widetilde{\mathcal{E}}_{X}^{\ast
}\widetilde{\mathcal{E}}_{X}f\left(  x\right)  ,
\]

and then applying the estimate \ref{1.63} to both terms on the right side
yields%
\begin{align*}
\left\vert \mathcal{S}_{X}g\left(  x\right)  \right\vert  & \leq\left\vert
\mathcal{S}_{X}f\left(  x\right)  \right\vert +\frac{1}{NR_{0}\left(
0\right)  }\left\vert \mathcal{S}_{X}\left(  \widetilde{\mathcal{E}}_{X}%
^{\ast}\widetilde{\mathcal{E}}_{X}f\right)  \left(  x\right)  \right\vert \\
& \leq\sqrt{R_{0}\left(  0\right)  }\min\left\{  1,\frac{R_{0}\left(
0\right)  }{\rho}\right\}  \left\Vert f\right\Vert _{w,0}+\frac{\sqrt
{R_{0}\left(  0\right)  }}{N\left\Vert R_{0}\right\Vert _{w,0}^{2}}%
\min\left\{  1,\frac{R_{0}\left(  0\right)  }{\rho}\right\}  \left\Vert
\widetilde{\mathcal{E}}_{X}^{\ast}\widetilde{\mathcal{E}}_{X}f\right\Vert
_{w,0}\\
& =\sqrt{R_{0}\left(  0\right)  }\min\left\{  1,\frac{R_{0}\left(  0\right)
}{\rho}\right\}  \left(  \left\Vert f\right\Vert _{w,0}+\frac{1}{NR_{0}\left(
0\right)  }\left\Vert \widetilde{\mathcal{E}}_{X}^{\ast}\widetilde{\mathcal{E}%
}_{X}f\right\Vert _{w,0}\right) \\
& \leq2\sqrt{R_{0}\left(  0\right)  }\min\left\{  1,\frac{R_{0}\left(
0\right)  }{\rho}\right\}  \left\Vert f\right\Vert _{w,0},
\end{align*}

where the last step used \ref{1.48}. Finally%
\begin{align*}
\left\vert \mathcal{S}_{X}f\left(  x\right)  -\frac{\left(
\widetilde{\mathcal{E}}_{X}R_{x}\right)  ^{T}\widetilde{\mathcal{E}}_{X}%
f}{\left(  R_{0}\left(  0\right)  +\rho\right)  N}\right\vert  & =\frac
{R_{0}\left(  0\right)  }{R_{0}\left(  0\right)  +\rho}\left\vert
\mathcal{S}_{X}g\left(  x\right)  \right\vert \\
& \leq\frac{2R_{0}\left(  0\right)  ^{3/2}}{R_{0}\left(  0\right)  +\rho}%
\min\left\{  1,\frac{R_{0}\left(  0\right)  }{\rho}\right\}  \left\Vert
f\right\Vert _{w,0},
\end{align*}

which proves \ref{15.20}.
\end{proof}

\begin{remark}
\label{Rem_Thm_estim_Exact_smth_by_Wgx}\ 

\begin{enumerate}
\item If we write
\begin{align*}
\frac{\left(  \widetilde{\mathcal{E}}_{X}R_{x}\right)  ^{T}%
\widetilde{\mathcal{E}}_{X}f}{\left(  R_{0}\left(  0\right)  +\rho\right)  N}
& =\frac{1}{\left(  R_{0}\left(  0\right)  +\rho\right)  N}\sum_{k=1}%
^{N}f\left(  x^{\left(  k\right)  }\right)  R_{x}\left(  x^{\left(  k\right)
}\right) \\
& =\frac{\left(  2\pi\right)  ^{-d/2}}{\left(  R_{0}\left(  0\right)
+\rho\right)  }\frac{1}{N}\sum_{k=1}^{N}f\left(  x^{\left(  k\right)
}\right)  G\left(  x-x^{\left(  k\right)  }\right) \\
& =\frac{1}{\left(  G\left(  0\right)  +\left(  2\pi\right)  ^{d/2}%
\rho\right)  }\frac{1}{N}\sum_{k=1}^{N}f\left(  x^{\left(  k\right)  }\right)
G\left(  x-x^{\left(  k\right)  }\right)  ,
\end{align*}

it is clear that the mapping $f\rightarrow\frac{\left(  \widetilde{\mathcal{E}%
}_{X}R_{x}\right)  ^{T}\widetilde{\mathcal{E}}_{X}f}{\left(  R_{0}\left(
0\right)  +\rho\right)  N}$ is a continuous linear operator from $X_{w}^{0}$
to $W_{G,X}$ with the form of a weighted average, the weights being
data-translated basis functions. Note also that the operator only depends on
the values of the data function in the data region.

\item The approximation \ref{15.20} provides a bound of order $\rho^{-2}$ for
large $\rho$. For all $\rho$ this estimate supplies an upper bound of
$2\sqrt{R_{0}\left(  0\right)  }\left\Vert f\right\Vert _{w,0}$.

\item The normwise approximation \ref{15.02} provides a bound of order
$\rho^{-1}$ as $\rho\rightarrow\infty$. For all $\rho$ this estimate supplies
an upper bound of $2\left\Vert f\right\Vert _{w,0}$.
\end{enumerate}
\end{remark}

\chapter{The scalable Approximate smoother\label{Ch_Approx_smth}}

\section{Introduction\label{Sect_ap_introd}}

In \textbf{Section} \ref{Sect_space_Jg} the \textbf{convolution space}
$J_{G}=G\ast S$ is introduced and then some applications of these spaces are
presented. The fact that $J_{G}\subset X_{w}^{0}$ implies some interesting
embedding results which involve supplying necessary and sufficient conditions
for one data space $X_{w}^{0}$ to be \textbf{continuously embedded} in
another. I have also begun a study of the operator $\int_{\Omega}\overline
{R}_{x}u$ where $u\in X_{w}^{0}$ and $R_{x}$ is the Riesz representer of the
evaluation functional. This operator might be useful in deriving the
Approximate smoother.

In \textbf{Section} \ref{Sect_ApproxSmooth} I \textbf{introduce the
Approximate smoother}. I call this smoother the Approximate smoother because
it approximates the Exact smoother! This is a non-parametric, scalable
smoother. Here \textbf{scalable} means the numeric effort to calculate the
Approximate smoother depends linearly on the number of data points.\textit{\ }%
We assume the basis function is real-valued.

\textbf{Two different approaches} will be taken to defining the Approximate
smoother, and both involve formulating the smoother as the solution of a
variational problem. One of these problems will involve minimizing the Exact
smoother functional \ref{1.8} over $W_{G,X^{\prime}}$ where $X^{\prime
}=\left\{  x_{i}^{\prime}\right\}  _{i=1}^{N^{\prime}}$ is an arbitrary set of
distinct points in $\mathbb{R}^{d}$. The other, equivalent problem, involves
finding the function in $W_{G,X^{\prime}}$ which is nearest to Exact smoother
$s_{e}$ w.r.t. the norm $\left\Vert \cdot\right\Vert _{w,0}$. If
\[
s_{a}\left(  x\right)  =\sum\limits_{i=1}^{N^{\prime}}\alpha_{i}^{\prime
}R_{x_{i}^{\prime}}\left(  x\right)  =\left(  2\pi\right)  ^{-d/2}%
\sum\limits_{i=1}^{N^{\prime}}\alpha_{i}^{\prime}G\left(  x-x_{i}^{\prime
}\right)  ,
\]

denotes the Approximate smoother and $y$ is the dependent data then solving
the second problem yields%
\[
s_{a}=\mathcal{I}_{X^{\prime}}s_{e},
\]

which implies the matrix equation%
\[
\left(  N\rho R_{X^{\prime},X^{\prime}}+R_{X,X^{\prime}}^{T}R_{X,X^{\prime}%
}\right)  \alpha^{\prime}=R_{X,X^{\prime}}^{T}y,
\]

where $R_{X,X^{\prime}}=\left(  R_{x_{j}^{\prime}}\left(  x^{\left(  i\right)
}\right)  \right)  $. The size of the Approximate smoother matrix is
$N^{\prime}\times N^{\prime}$ which is independent of the number of data
points and suggests scalability.

The \textbf{error estimates for the pointwise convergence} of the Approximate
smoother to its data function $f\in X_{w}^{0}$ are based on the simple
triangle inequality%
\[
\left\vert f\left(  x\right)  -s_{a}\left(  x\right)  \right\vert
\leq\left\vert f\left(  x\right)  -s_{e}\left(  x\right)  \right\vert
+\left\vert s_{e}\left(  x\right)  -s_{a}\left(  x\right)  \right\vert ,
\]

and so Section \ref{Sect_ap_Ex_smth_minus_App_smth} will be devoted to
estimating $\left\vert s_{e}\left(  x\right)  -s_{a}\left(  x\right)
\right\vert $.

As with the minimal interpolant and the Exact smoother, we will obtain
estimates that assume unisolvent data sets as well as the Type 1 and Type 2
estimates that do not \textbf{explicitly} involve unisolvency. The Approximate
smoother convergence orders and the constants are the same as those for the
interpolation case which are given in the interpolation tables
\ref{Tbl_intro_NonUnisolvTyp1Converg}, \ref{Tbl_intro_NonUnisolvTyp2Conv},
\ref{Tbl_intro_UnisolvConverg} and \ref{Tbl_intro_ConvergCentral}.\medskip

\fbox{Type 1 error estimates} No \textit{a priori} assumption is made
concerning the weight function parameter $\kappa$ but it will be assumed that
the basis function satisfies an inequality of the form \ref{1.9}. For example,
if it is assumed that the data region $K$ is closed bounded and infinite then
Theorem \ref{Thm_converg_arb_func_K=0} establishes that%
\[
\left\vert s_{e}\left(  x\right)  -s_{a}\left(  x\right)  \right\vert
\leq\left\Vert f\right\Vert _{w,0}k_{G}\left(  h_{X^{\prime},K}\right)
^{s},\quad x\in K,
\]

and Theorem \ref{Thm_ap_Appr_smth_err_typ1} shows that%
\[
\left\vert f\left(  x\right)  -s_{a}\left(  x\right)  \right\vert
\leq\left\Vert f\right\Vert _{w,0}\left(  \sqrt{\rho N}+k_{G}\left(
h_{X,K}\right)  ^{s}+k_{G}\left(  h_{X^{\prime},K}\right)  ^{s}\right)  ,\quad
x\in K,
\]

when $h_{X,K}=\sup\limits_{x\in K}\operatorname*{dist}\left(  x,X\right)  \leq
h_{G}$ and\allowbreak\ $h_{X^{\prime},K}=\sup\limits_{x\in K}%
\operatorname*{dist}\left(  x,X^{\prime}\right)  \leq h_{G}$.\medskip

\fbox{Type 2 error estimates} If it only assumed that $\kappa\geq1$ then by
Theorem \ref{Thm_ap_G(0)minusG(x)_bound}%
\[
\left\vert s_{e}\left(  x\right)  -s_{a}\left(  x\right)  \right\vert
\leq\left\Vert f\right\Vert _{w,0}k_{G}\left(  h_{X^{\prime},K}\right)
^{s},\quad x\in\mathbb{R}^{d},
\]

and by Theorem \ref{Thm_ap_Appr_smth_err_typ2}%
\[
\left\vert f\left(  x\right)  -s_{a}\left(  x\right)  \right\vert
\leq\left\Vert f\right\Vert _{w,0}\left(  \sqrt{\rho N}+k_{G}h_{X,K}%
+k_{G}h_{X^{\prime},K}\right)  ,\quad x\in\mathbb{R}^{d},
\]

where $k_{G}=\left(  2\pi\right)  ^{-\frac{d}{4}}\sqrt{-\left(  \left\vert
D\right\vert ^{2}G\right)  \left(  0\right)  }\sqrt{d}$.\medskip

\fbox{Taylor series and explicit unisolvent data error estimates} If $X$ is a
unisolvent set of order $m\geq1$ contained in a bounded data region $\Omega$
then by Theorem \ref{Thm_ConvergSaToSe_not_unisolv_II}
\[
\left\vert s_{e}\left(  x\right)  -s_{a}\left(  x\right)  \right\vert
\leq\left\Vert f\right\Vert _{w,0}k_{G}\left(  h_{X^{\prime},\Omega}\right)
^{m},\quad x\in\overline{\Omega},
\]

and by Theorem \ref{Thm_Err_unisolv_arbitrary_fn}%
\[
\left\vert f\left(  x\right)  -s_{a}\left(  x\right)  \right\vert
\leq\left\Vert f\right\Vert _{w,0}\left(  K_{\Omega,m}^{\prime}\sqrt{\rho
N}+k_{G}\left(  h_{X,K}\right)  ^{m}+k_{G}\left(  h_{X^{\prime},K}\right)
^{m}\right)  ,\quad x\in\overline{\Omega},
\]

for some constants $K_{\Omega,m}^{\prime},k_{G}>0$. We say the orders of
convergence are at least $m$.

These theoretical error results will be illustrated using the weight function
examples from the interpolation chapter, namely the radial \textbf{shifted
thin-plate splines}, \textbf{Gaussian} and \textbf{Sobolev splines} and the
tensor product \textbf{extended B-splines}. We will also use the
\textbf{central difference} weight functions from Chapter
\ref{Ch_cent_diff_wt_fn_ten_prod}.

\textbf{Numerical results are only presented for the Types 1 and 2 estimates}.
Numeric experiments are carried out using the same 1-dimensional B-splines and
data functions that were used for the interpolants. We restrict ourselves to
one dimension so that the data density parameters $h_{X,\Omega} $ and
$h_{X,K}$ can be easily calculated.

The last section discusses the \textbf{SmoothOperator} software (freeware)
package which implements the Approximate smoother algorithm. It has a full
user manual which describe several tutorials and data experiments.

\section{The convolution space $J_{G}$ and applications\label{Sect_space_Jg}}

Following Dyn \cite{Dyn89} we introduce the convolution space $J_{G}$ and
prove some of its properties, including the fact that $J_{G}$ is dense in
$X_{w}^{0}$. This space will be used in the discretization process which
derives the Approximate smoother problem from the Exact smoother problem in
Subsection \ref{SbSect_ap_deriv_Approx_smth}.

An application of the fact that $J_{G}\subset X_{w}^{0}$ are the embedding
results of Subsection \ref{SbSect_EmbedInXow}, an example of which is Theorem
\ref{Thm_Xw1o_embed_Xw2o_iff}.

\begin{definition}
\label{Def_Jg}\textbf{The convolution space} $J_{G}$

Suppose the weight function $w$ satisfies property W02 for some $\kappa\geq0$
and let $G$ be the zero order basis function generated by $w$. Then define
\[
J_{G}:=G\ast S:=\left\{  G\ast\phi:\phi\in S\right\}  ,
\]

where the convolution $G\ast\phi$ is defined by $G\ast\phi=\left(
\widehat{\phi}\widehat{G}\right)  ^{\vee}$ with $G\in S^{\prime}$ and $\phi\in
S$. We will sometimes write $J_{G}=G\ast S$.
\end{definition}

\begin{theorem}
\label{Thm_Jg_properties}The function space $J_{G}$ has the following properties:

\begin{enumerate}
\item $J_{G}\subset X_{w}^{0}\cap C_{B}^{\infty}$.

\item If $f\in X_{w}^{0}$ then $\left(  f,G\ast\phi\right)  _{w,0}=\left[
f,\overline{\phi}\right]  $ for all $\phi\in S$.

\item $J_{G}$ is dense in $X_{w}^{0}$.

\item $\left(  G\ast\phi\right)  \left(  x\right)  =\int\overline{R}_{x}\phi$
when $\phi\in S$.

\item If $w$ is a tensor product central difference weight function then
$G\ast\left(  S\cap X_{w}^{0}\right)  $ is dense in $X_{w}^{0}$.

\item ??? Is part 5 true in general?
\end{enumerate}
\end{theorem}

\begin{proof}
Suppose $g\in J_{G}$, say $g=G\ast\phi$ where $\phi\in S$.\medskip

\textbf{Part 1} $\widehat{g}=\widehat{\phi}$ $\widehat{G}=\frac{\widehat{\phi
}}{w}\in L^{1}$ since property W02 implies $\frac{1}{w}\in L^{1}$. Thus
\[
\int w\left\vert \widehat{g}\right\vert ^{2}=\int w\frac{\left\vert
\widehat{\phi}\right\vert ^{2}}{w^{2}}=\int\frac{\left\vert \widehat{\phi
}\right\vert ^{2}}{w}\leq\left\Vert \widehat{\phi}\right\Vert _{\infty}%
^{2}\int\frac{1}{w}<\infty,
\]

and so $g\in X_{w}^{0}$. For any multi-index $\alpha\geq0$%
\[
\left(  D^{\alpha}\left(  G\ast\phi\right)  \right)  ^{\wedge}=i^{\left\vert
\alpha\right\vert }\xi^{\alpha}\widehat{G\ast\phi}=i^{\left\vert
\alpha\right\vert }\xi^{\alpha}\widehat{\phi}\widehat{G}=\frac{i^{\left\vert
\alpha\right\vert }\xi^{\alpha}\widehat{\phi}}{w}\in L^{1},
\]

and so by Lemma \ref{Lem_L1_Fourier_contin},$\ D^{\alpha}\left(  G\ast
\phi\right)  \in C_{B}^{\left(  0\right)  }$. Thus $G\ast\phi\in C_{B}%
^{\infty}$.\medskip

\textbf{Part 2} $f\in X_{w}^{0}$ implies $\widehat{f}\in S^{\prime}\cap
L_{loc}^{1}$ and%
\[
\left(  f,G\ast\phi\right)  _{w,0}=\int w\widehat{f}\text{ }\overline
{\widehat{G\ast\phi}}=\int w\widehat{f}\frac{\overline{\widehat{\phi}}}%
{w}=\int\widehat{f}\overline{\widehat{\phi}}.
\]

We now need some results from part (a) of Subsection 2.8.3 of Vladimirov
\cite{Vladimirov}. Here he defines locally integrable functions which have
polynomial growth at infinity i.e. a function $g\in L_{loc}^{1}$ such that
$\int\frac{\left\vert g\right\vert }{\left(  1+\left\vert \cdot\right\vert
\right)  ^{s}}$ for some $s\geq0$. Vladimirov states that these functions have
the property that $g\in S^{\prime}$ and $\left[  g,\phi\right]  =\int g\phi$
for $\phi\in S$.

Now, from the proof of Theorem \ref{Thm_basis_fn_properties_all_m_W2},
$\widehat{f}\in L^{1}$ and so $\widehat{f}$ has polynomial growth at infinity
which implies%
\[
\int\widehat{f}\overline{\widehat{\phi}}=\left[  \widehat{f},\overline
{\widehat{\phi}}\right]  =\left[  f,\overline{\phi}\right]  .
\]
\medskip

\textbf{Part 3} A standard result is that a subspace of a Hilbert space is
dense iff its orthogonal complement is $\left\{  0\right\}  $. In fact, if
$\left(  f,G\ast\phi\right)  _{w,0}=0$ for all $\phi\in S$ then part 2 implies
$f=0$.\medskip

\textbf{Part 4} Let $f=R_{x}$ in part 2.\medskip

\textbf{Part 5} By definition $\overset{\vee}{S}_{w,0}=S\cap X_{w}^{0}$ and
from part 4 of Corollary \ref{Cor_Thm_eta^2/sin(eta)^2.|phi|^2_dim_gt_1},
$G\ast\overset{\vee}{S}_{w,0}$ is dense in $X_{w}^{0}$.\medskip

\textbf{Part 6} ??
\end{proof}

\subsection{Embedding results for $X_{w}^{0}$\label{SbSect_EmbedInXow}}

The fact that $J_{G}\subset X_{w}^{0}$ implies the following interesting
embedding results which involve supplying necessary and sufficient conditions
for one data space to be continuously embedded in another.

\begin{theorem}
\label{Thm_Xow_embed_L2_iff}Suppose the weight function $w$ has property $W01
$ so that $X_{w}^{0}$ is a normed vector space. Then $X_{w}^{0}\hookrightarrow
L^{2}$ iff $\frac{1}{w}\in L^{\infty}$. In fact, if $\frac{1}{w}\in L^{\infty
}$ then%
\[
\left\Vert f\right\Vert _{2}\leq\left\Vert \frac{1}{w}\right\Vert _{\infty
}^{1/2}\left\Vert f\right\Vert _{w,0},\quad f\in X_{w}^{0},
\]

and if the embedding $\iota:X_{w}^{0}\hookrightarrow L^{2}$ is continuous then%
\[
\left\Vert \frac{1}{w}\right\Vert _{\infty}^{1/2}=\left\Vert \iota\right\Vert
.
\]

\end{theorem}

\begin{proof}
\textbf{First suppose} $\frac{1}{w}\in L^{\infty}$. Then, by Plancherel's
theorem for $L^{2}$ functions i.e. $\left\Vert f\right\Vert _{2}=\left\Vert
\widehat{f}\right\Vert _{2}$, we have%
\begin{equation}
\int\left\vert f\right\vert ^{2}=\int\left\vert \widehat{f}\right\vert
^{2}=\int\frac{1}{w}w\left\vert \widehat{f}\right\vert ^{2}\leq\left\Vert
\frac{1}{w}\right\Vert _{\infty}\left\Vert f\right\Vert _{w,0}^{2}%
.\label{1.055}%
\end{equation}

\textbf{Next suppose} $X_{w}^{0}\subset L^{2}$ and the embedding $\iota
:X_{w}^{0}\hookrightarrow L^{2}$ is continuous i.e.%
\[
\left\Vert f\right\Vert _{2}\leq\left\Vert \iota\right\Vert \left\Vert
f\right\Vert _{w,0},\quad f\in X_{w}^{0}.
\]

But from part 1 of Theorem \ref{Thm_Jg_properties} we know that $G\ast\phi\in
X_{w}^{0}$ whenever $\phi\in S$ and thus $G\ast\overset{\vee}{\phi}\in
X_{w}^{0}$ $\forall\phi\in S$. Now by assumption \ref{1.055}, $\left\Vert
G\ast\overset{\vee}{\phi}\right\Vert _{2}^{2}\leq\left\Vert \iota\right\Vert
^{2}\left\Vert G\ast\overset{\vee}{\phi}\right\Vert _{w,0}^{2}$ i.e.
\[
\int\frac{1}{w^{2}}\left\vert \phi\right\vert ^{2}\leq\left\Vert
\iota\right\Vert ^{2}\int\frac{1}{w}\left\vert \phi\right\vert ^{2},\quad
\phi\in S.
\]

With reference to weight function Definition \ref{Def_extend_wt_fn}, choose
arbitrary $x\notin\mathcal{A}$ where $\mathcal{A}$ is a closed set of measure
zero such that $w>0$ and continuous on $\mathbb{R}^{d}\setminus\mathcal{A}$.
Suppose $0<r<R<\operatorname*{dist}\left(  x,\mathcal{A}\right)  $. Then it is
possible to choose $\phi_{r,R}\in S$ such that $0\leq\phi_{r,R}\leq1$,
$\operatorname*{supp}\phi_{r,R}\subseteq\overline{B}_{R}\left(  x\right)  $
and $\phi_{r,R}=1$ on $\overline{B}_{r}\left(  x\right)  $ e.g. the standard
"cap" example of a $C_{0}^{\infty}$ function. Observe now that%
\[
\min_{B_{r}\left(  x\right)  }\frac{1}{w^{2}}\int\limits_{B_{r}\left(
x\right)  }1\leq\int\limits_{B_{R}\left(  x\right)  }\frac{\left\vert
\phi_{r,R}\right\vert ^{2}}{w^{2}}\leq\left\Vert \iota\right\Vert ^{2}%
\int\limits_{B_{R}\left(  x\right)  }\frac{\left\vert \phi_{r,R}\right\vert
^{2}}{w}\leq\left\Vert \iota\right\Vert ^{2}\max_{B_{R}\left(  x\right)
}\frac{1}{w}\int\limits_{B_{R}\left(  x\right)  }1,
\]

i.e. for $0<r<R<\operatorname*{dist}\left(  x,\mathcal{A}\right)  $,%
\[
\frac{\min\limits_{B_{r}\left(  x\right)  }\frac{1}{w^{2}}}{\max
\limits_{B_{R}\left(  x\right)  }\frac{1}{w}}\leq\left\Vert \iota\right\Vert
^{2}\frac{\int_{B_{R}\left(  x\right)  }1}{\int_{B_{r}\left(  x\right)  }%
1}=\left\Vert \iota\right\Vert ^{2}\left(  \frac{R}{r}\right)  ^{d}.
\]

Now since the continuity of $w$ near $x$ implies $\min\limits_{B_{R}\left(
x\right)  }\frac{1}{w^{2}}$ is continuous in $R$ we must have $\frac
{\min\limits_{B_{r}\left(  x\right)  }\frac{1}{w^{2}}}{\max\limits_{B_{r}%
\left(  x\right)  }\frac{1}{w}}\leq\left\Vert \iota\right\Vert ^{2}$ for
$0<r<\operatorname*{dist}\left(  x,\mathcal{A}\right)  $, and hence that
$\frac{1}{w\left(  x\right)  }=\lim\limits_{r\rightarrow0}\frac{\min
\limits_{B_{r}\left(  x\right)  }\frac{1}{w^{2}}}{\max\limits_{B_{r}\left(
x\right)  }\frac{1}{w}}\leq\left\Vert \iota\right\Vert ^{2}$ when
$x\notin\mathcal{A}$ i.e. $1/w\in L^{\infty}$ and $\left\Vert \frac{1}%
{w}\right\Vert _{\infty}^{1/2}\leq\left\Vert \iota\right\Vert $. But from
\ref{1.055}, $\left\Vert \iota\right\Vert \leq\left\Vert \frac{1}%
{w}\right\Vert _{\infty}^{1/2}$ so $\left\Vert \iota\right\Vert =\left\Vert
\frac{1}{w}\right\Vert _{\infty}^{1/2}$.
\end{proof}

More generally, we have:

\begin{theorem}
\label{Thm_Xow_embed_Hm_iff}Suppose the weight function $w$ has property W01
so that $X_{w}^{0}$ is a normed vector space. Then $X_{w}^{0}\overset{\iota
}{\hookrightarrow}H^{m}$ iff $\frac{\xi^{2\alpha}}{w}\in L^{\infty}$ for
$\left\vert \alpha\right\vert \leq m$ iff $\frac{\left(  1+\left\vert
\cdot\right\vert ^{2}\right)  ^{m}}{w}\in L^{\infty}$.

In fact, if $\frac{\xi^{2\alpha}}{w}\in L^{\infty}$ for $\left\vert
\alpha\right\vert \leq m$ then%
\[
\left\Vert f\right\Vert _{H^{m}}\leq\left\Vert \frac{\left(  1+\left\vert
\cdot\right\vert ^{2}\right)  ^{m}}{w}\right\Vert _{\infty}^{1/2}\left\Vert
f\right\Vert _{w,0},\quad f\in X_{w}^{0},
\]

and if the embedding $\iota:X_{w}^{0}\hookrightarrow H^{m}$ is continuous then%
\[
\left\Vert \frac{\left(  1+\left\vert \cdot\right\vert ^{2}\right)  ^{m}}%
{w}\right\Vert _{\infty}^{1/2}=\left\Vert \iota\right\Vert .
\]

\end{theorem}

\begin{proof}
The special case of the identity stated in part 6 in Definition
\ref{Def_multi_id} i.e.%
\[
\left(  1+\left\vert \xi\right\vert ^{2}\right)  ^{m}=\sum_{\left\vert
\alpha\right\vert \leq m}\tbinom{m}{\alpha}\xi^{\alpha},\quad\tbinom{m}%
{\alpha}=\frac{m!}{\alpha!\left(  m-\left\vert \alpha\right\vert \right)  !}.
\]

implies $\frac{\xi^{2\alpha}}{w}\in L^{\infty}$ for $\left\vert \alpha
\right\vert \leq m$ iff $\frac{\left(  1+\left\vert \xi\right\vert
^{2}\right)  ^{m}}{w}\in L^{\infty}$.\smallskip

\textbf{First suppose} $\frac{\xi^{2\alpha}}{w}\in L^{\infty}$ for $\left\vert
\alpha\right\vert \leq m$. Then%
\[
\left\Vert f\right\Vert _{H^{m}}^{2}=\int\left(  1+\left\vert \cdot\right\vert
^{2}\right)  ^{m}\left\vert \widehat{f}\right\vert ^{2}=\int\frac{\left(
1+\left\vert \cdot\right\vert ^{2}\right)  ^{m}}{w}w\left\vert \widehat{f}%
\right\vert ^{2}\leq\left\Vert \frac{\left(  1+\left\vert \cdot\right\vert
^{2}\right)  ^{m}}{w}\right\Vert _{\infty}\left\Vert f\right\Vert _{w,0}%
^{2}<\infty.
\]
\smallskip

\textbf{Next suppose} $X_{w}^{0}\subset H^{m}$ and the embedding $\iota
:X_{w}^{0}\hookrightarrow H^{m}$ is continuous i.e.%
\begin{equation}
\left\Vert f\right\Vert _{H^{m}}\leq\left\Vert \iota\right\Vert \left\Vert
f\right\Vert _{w,0},\quad f\in X_{w}^{0}.\label{1.091}%
\end{equation}

But from part 1 of Theorem \ref{Thm_Jg_properties} we know that if $G$ is the
basis function then for all $\phi\in S$, $G\ast\phi\in X_{w}^{0}$ and thus
$G\ast\overset{\vee}{\phi}\in X_{w}^{0}$. Consequently, by assumption
\ref{1.091},$\allowbreak\left\Vert G\ast\overset{\vee}{\phi}\right\Vert
_{H^{m}}^{2}\leq\left\Vert \iota\right\Vert ^{2}\left\Vert G\ast
\overset{\vee}{\phi}\right\Vert _{w,0}^{2}$ i.e.
\[
\int\left(  1+\left\vert \cdot\right\vert ^{2}\right)  ^{m}\frac{1}{w^{2}%
}\left\vert \phi\right\vert ^{2}\leq\left\Vert \iota\right\Vert ^{2}\int%
\frac{1}{w}\left\vert \phi\right\vert ^{2},\quad\phi\in S.
\]

With reference to weight function Definition \ref{Def_extend_wt_fn}, choose
arbitrary $\eta\notin\mathcal{A}$ where $\mathcal{A}$ is a closed set of
measure zero such that $w>0$ and continuous on $\mathbb{R}^{d}\setminus
\mathcal{A}$. Suppose $0<r<R<\operatorname*{dist}\left(  \eta,\mathcal{A}%
\right)  $. Then it is possible to choose $\phi_{r,R}\in S$ such that
$0\leq\phi_{r,R}\leq1$, $\operatorname*{supp}\phi_{r,R}\subseteq\overline
{B}_{R}\left(  \eta\right)  $ and $\phi_{r,R}=1$ on $\overline{B}_{r}\left(
\eta\right)  $ e.g. the standard "cap" example of a $C_{0}^{\infty}$ function.
Observe now that%
\[
\min_{B_{r}\left(  \eta\right)  }\frac{\left(  1+\left\vert \cdot\right\vert
^{2}\right)  ^{m}}{w^{2}}\int\limits_{B_{r}\left(  \eta\right)  }1\leq
\int\limits_{B_{R}\left(  \eta\right)  }\left(  1+\left\vert \cdot\right\vert
^{2}\right)  ^{m}\frac{\left\vert \phi_{r,R}\right\vert ^{2}}{w^{2}}%
\leq\left\Vert \iota\right\Vert ^{2}\int\limits_{B_{R}\left(  \eta\right)
}\frac{\left\vert \phi_{r,R}\right\vert ^{2}}{w}\leq\left\Vert \iota
\right\Vert ^{2}\max_{B_{R}\left(  \eta\right)  }\frac{1}{w}\int%
\limits_{B_{R}\left(  \eta\right)  }1,
\]

i.e. for $0<r<R<\operatorname*{dist}\left(  \eta,\mathcal{A}\right)  $,%
\[
\frac{\min\limits_{B_{r}\left(  \eta\right)  }\frac{\left(  1+\left\vert
\cdot\right\vert ^{2}\right)  ^{m}}{w^{2}}}{\max\limits_{B_{R}\left(
\eta\right)  }\frac{1}{w}}\leq\left\Vert \iota\right\Vert ^{2}\frac
{\int_{B_{R}\left(  \eta\right)  }1}{\int_{B_{r}\left(  \eta\right)  }%
1}=\left\Vert \iota\right\Vert ^{2}\left(  \frac{R}{r}\right)  ^{d}.
\]

Now since the continuity of $w$ near $\eta$ implies $\min\limits_{B_{R}\left(
\eta\right)  }\frac{1}{w}$ is continuous in $R$ we must have $\frac
{\min\limits_{B_{r}\left(  \eta\right)  }\frac{\left(  1+\left\vert
\cdot\right\vert ^{2}\right)  ^{m}}{w^{2}}}{\max\limits_{B_{r}\left(
\eta\right)  }\frac{1}{w}}\leq\left\Vert \iota\right\Vert ^{2}$ for
$0<r<\operatorname*{dist}\left(  \eta,\mathcal{A}\right)  $, and hence that
\[
\frac{\left(  1+\left\vert \eta\right\vert ^{2}\right)  ^{m}}{w\left(
\eta\right)  }=\frac{\frac{\left(  1+\left\vert \eta\right\vert ^{2}\right)
^{m}}{w\left(  \eta\right)  ^{2}}}{\frac{1}{w\left(  \eta\right)  }}%
=\lim\limits_{r\rightarrow0}\frac{\min\limits_{B_{r}\left(  \eta\right)
}\frac{\left(  1+\left\vert \cdot\right\vert ^{2}\right)  ^{m}}{w^{2}}}%
{\max\limits_{B_{r}\left(  \eta\right)  }\frac{1}{w}}\leq\left\Vert
\iota\right\Vert ^{2}.
\]

when $\eta\notin\mathcal{A}$ i.e. $\frac{\left(  1+\left\vert \eta\right\vert
^{2}\right)  ^{m}}{w\left(  \eta\right)  }\in L^{\infty}$ and $\left\Vert
\frac{\left(  1+\left\vert \cdot\right\vert ^{2}\right)  ^{m}}{w}\right\Vert
_{\infty}^{1/2}\leq\left\Vert \iota\right\Vert $ i.e. $\frac{\xi^{2\alpha}}%
{w}\in L^{\infty}$ for $\left\vert \alpha\right\vert \leq m $.
\end{proof}

\begin{remark}
\label{Rem_Thm_Xow_embed_Hm_iff}Compare the definition of weight function
property W02 for parameter $m$ i.e. $\frac{\xi^{2\alpha}}{w}\in L^{1}$ for
$\left\vert \alpha\right\vert \leq m$, with the condition $\frac{\xi^{2\alpha
}}{w}\in L^{\infty}$ for $\left\vert \alpha\right\vert \leq m$ used in the
last theorem
\end{remark}

\begin{theorem}
\label{Thm_Xow_embed_Hs_iff}Suppose the weight function $w$ has property W01
so that $X_{w}^{0}$ is a normed vector space. Suppose $s\geq0$. Then
$X_{w}^{0}\overset{\iota}{\hookrightarrow}H^{s}$ iff $\frac{\left(
1+\left\vert \cdot\right\vert ^{2}\right)  ^{s}}{w}\in L^{\infty}$.

In fact, if $\frac{\left(  1+\left\vert \cdot\right\vert ^{2}\right)  ^{s}}%
{w}\in L^{\infty}$ then%
\[
\left\Vert f\right\Vert _{H^{s}}\leq\left\Vert \frac{\left(  1+\left\vert
\cdot\right\vert ^{2}\right)  ^{s}}{w}\right\Vert _{\infty}^{1/2}\left\Vert
f\right\Vert _{w,0},\quad f\in X_{w}^{0},
\]

and if the embedding $\iota:X_{w}^{0}\hookrightarrow H^{s}$ is continuous
then$\left\Vert \frac{\left(  1+\left\vert \cdot\right\vert ^{2}\right)  ^{s}%
}{w}\right\Vert _{\infty}^{1/2}=\left\Vert \iota\right\Vert $.
\end{theorem}

Even more generally:

\begin{proof}
\textbf{First suppose} $\frac{\left(  1+\left\vert \cdot\right\vert
^{2}\right)  ^{s}}{w}$. Then%
\[
\left\Vert f\right\Vert _{H^{s}}^{2}=\int\left(  1+\left\vert \cdot\right\vert
^{2}\right)  ^{s}\left\vert \widehat{f}\right\vert ^{2}=\int\frac{\left(
1+\left\vert \cdot\right\vert ^{2}\right)  ^{s}}{w}w\left\vert \widehat{f}%
\right\vert ^{2}\leq\left\Vert \frac{\left(  1+\left\vert \cdot\right\vert
^{2}\right)  ^{s}}{w}\right\Vert _{\infty}\left\Vert f\right\Vert _{w,0}%
^{2}<\infty.
\]

\textbf{Next suppose} $X_{w}^{0}\subset H^{s}$ and the embedding $\iota
:X_{w}^{0}\hookrightarrow H^{s}$ is continuous i.e.%
\begin{equation}
\left\Vert f\right\Vert _{H^{s}}\leq\left\Vert \iota\right\Vert \left\Vert
f\right\Vert _{w,0},\quad f\in X_{w}^{0}.\label{1.093}%
\end{equation}

But from part 1 of Theorem \ref{Thm_Jg_properties} we know that if $G$ is the
basis function then for all $\phi\in S$, $G\ast\phi\in X_{w}^{0}$ and thus
$G\ast\overset{\vee}{\phi}\in X_{w}^{0}$. Consequently, by assumption
\ref{1.093},$\allowbreak\left\Vert G\ast\overset{\vee}{\phi}\right\Vert
_{H^{s}}^{2}\leq\left\Vert \iota\right\Vert ^{2}\left\Vert G\ast
\overset{\vee}{\phi}\right\Vert _{w,0}^{2}$ i.e.
\[
\int\left(  1+\left\vert \cdot\right\vert ^{2}\right)  ^{s}\frac{1}{w^{2}%
}\left\vert \phi\right\vert ^{2}\leq\left\Vert \iota\right\Vert ^{2}\int%
\frac{1}{w}\left\vert \phi\right\vert ^{2},\quad\phi\in S.
\]

With reference to weight function Definition \ref{Def_extend_wt_fn}, choose
arbitrary $\eta\notin\mathcal{A}$ where $\mathcal{A}$ is a closed set of
measure zero such that $w>0$ and continuous on $\mathbb{R}^{d}\setminus
\mathcal{A}$. Suppose $0<r<R<\operatorname*{dist}\left(  \eta,\mathcal{A}%
\right)  $. Then it is possible to choose $\phi_{r,R}\in S$ such that
$0\leq\phi_{r,R}\leq1$, $\operatorname*{supp}\phi_{r,R}\subseteq\overline
{B}_{R}\left(  \eta\right)  $ and $\phi_{r,R}=1$ on $\overline{B}_{r}\left(
\eta\right)  $ e.g. the standard "cap" example of a $C_{0}^{\infty}$ function.
Observe now that%
\[
\min_{B_{r}\left(  \eta\right)  }\frac{\left(  1+\left\vert \cdot\right\vert
^{2}\right)  ^{s}}{w^{2}}\int\limits_{B_{r}\left(  \eta\right)  }1\leq
\int\limits_{B_{R}\left(  \eta\right)  }\left(  1+\left\vert \cdot\right\vert
^{2}\right)  ^{s}\frac{\left\vert \phi_{r,R}\right\vert ^{2}}{w^{2}}%
\leq\left\Vert \iota\right\Vert ^{2}\int\limits_{B_{R}\left(  \eta\right)
}\frac{\left\vert \phi_{r,R}\right\vert ^{2}}{w}\leq\left\Vert \iota
\right\Vert ^{2}\max_{B_{R}\left(  \eta\right)  }\frac{1}{w}\int%
\limits_{B_{R}\left(  \eta\right)  }1,
\]

i.e. for $0<r<R<\operatorname*{dist}\left(  \eta,\mathcal{A}\right)  $,%
\[
\frac{\min\limits_{B_{r}\left(  \eta\right)  }\frac{\left(  1+\left\vert
\cdot\right\vert ^{2}\right)  ^{s}}{w^{2}}}{\max\limits_{B_{R}\left(
\eta\right)  }\frac{1}{w}}\leq\left\Vert \iota\right\Vert ^{2}\frac
{\int_{B_{R}\left(  \eta\right)  }1}{\int_{B_{r}\left(  \eta\right)  }%
1}=\left\Vert \iota\right\Vert ^{2}\left(  \frac{R}{r}\right)  ^{d}.
\]

Now since the continuity of $w$ near $\eta$ implies $\min\limits_{B_{R}\left(
\eta\right)  }\frac{1}{w}$ is continuous in $R$ we must have $\frac
{\min\limits_{B_{r}\left(  \eta\right)  }\frac{\left(  1+\left\vert
\cdot\right\vert ^{2}\right)  ^{s}}{w^{2}}}{\max\limits_{B_{r}\left(
\eta\right)  }\frac{1}{w}}\leq\left\Vert \iota\right\Vert ^{2}$ for
$0<r<\operatorname*{dist}\left(  \eta,\mathcal{A}\right)  $, and hence that
\[
\frac{\left(  1+\left\vert \eta\right\vert ^{2}\right)  ^{s}}{w\left(
\eta\right)  }=\frac{\frac{\left(  1+\left\vert \eta\right\vert ^{2}\right)
^{s}}{w\left(  \eta\right)  ^{2}}}{\frac{1}{w\left(  \eta\right)  }}%
=\lim\limits_{r\rightarrow0}\frac{\min\limits_{B_{r}\left(  \eta\right)
}\frac{\left(  1+\left\vert \cdot\right\vert ^{2}\right)  ^{s}}{w^{2}}}%
{\max\limits_{B_{r}\left(  \eta\right)  }\frac{1}{w}}\leq\left\Vert
\iota\right\Vert ^{2}.
\]

when $\eta\notin\mathcal{A}$ i.e. $\frac{\left(  1+\left\vert \eta\right\vert
^{2}\right)  ^{s}}{w\left(  \eta\right)  }\in L^{\infty}$ and $\left\Vert
\frac{\left(  1+\left\vert \cdot\right\vert ^{2}\right)  ^{s}}{w}\right\Vert
_{\infty}^{1/2}\leq\left\Vert \iota\right\Vert $.
\end{proof}

The analogue of Theorem \ref{Thm_Xow_embed_Hm_iff} for $H^{m\mathbf{1}}$
spaces is:

\begin{theorem}
\label{Thm_Xow_embed_Hm1_iff}Suppose the weight function $w$ has property W01
so that $X_{w}^{0}$ is a normed vector space. Then $X_{w}^{0}\overset{\iota
}{\hookrightarrow}H^{m\mathbf{1}}$ iff $\frac{\xi^{2\alpha}}{w}\in L^{\infty}$
for $\alpha\leq m\mathbf{1}$ iff $\frac{\left(  1+\xi.\xi\right)
^{m\mathbf{1}}}{w\left(  \xi\right)  }\in L^{\infty}$.

In fact, if $\frac{\xi^{2\alpha}}{w}\in L^{\infty}$ for $\alpha\leq
m\mathbf{1}$ then%
\[
\left\Vert f\right\Vert _{H^{m\mathbf{1}}}\leq\left\Vert \frac{\left(
1+\xi.\xi\right)  ^{m\mathbf{1}}}{w\left(  \xi\right)  }\right\Vert _{\infty
}^{1/2}\left\Vert f\right\Vert _{w,0},\quad f\in X_{w}^{0},
\]

where $\xi.\xi$ is the componentwise vector product. Also, if the embedding
$\iota:X_{w}^{0}\hookrightarrow H^{m\mathbf{1}}$ is continuous then%
\[
\left\Vert \frac{\left(  1+\xi.\xi\right)  ^{m\mathbf{1}}}{w\left(
\xi\right)  }\right\Vert _{\infty}^{1/2}=\left\Vert \iota\right\Vert .
\]

\end{theorem}

\begin{proof}
The binomial identity%
\[
\left(  \left(  1+\xi_{1}^{2}\right)  \ldots\left(  1+\xi_{d}^{2}\right)
\right)  ^{m}=\left(  1+\xi.\xi\right)  ^{m\mathbf{1}}=\sum_{\alpha\leq
m\mathbf{1}}\tbinom{m\mathbf{1}}{\alpha}\xi^{2\alpha}.
\]

implies $\frac{\xi^{2\alpha}}{w}\in L^{\infty}$ for $\alpha\leq m\mathbf{1}$
iff $\frac{\left(  1+\xi.\xi\right)  ^{m\mathbf{1}}}{w}\in L^{\infty}%
$.\smallskip

\textbf{Now suppose} $\frac{\xi^{2\alpha}}{w}\in L^{\infty}$ for $\alpha\leq
m\mathbf{1}$. Then%
\[
\left\Vert f\right\Vert _{H^{m\mathbf{1}}}^{2}=\int\left(  1+\xi.\xi\right)
^{m\mathbf{1}}\left\vert \widehat{f}\right\vert ^{2}=\int\frac{\left(
1+\xi.\xi\right)  ^{m\mathbf{1}}}{w}w\left\vert \widehat{f}\right\vert
^{2}\leq\left\Vert \frac{\left(  1+\xi.\xi\right)  ^{m\mathbf{1}}}%
{w}\right\Vert _{\infty}\left\Vert f\right\Vert _{w,0}^{2}<\infty.
\]
\smallskip

\textbf{Next suppose} $X_{w}^{0}\subset H^{m\mathbf{1}}$ and the embedding
$\iota:X_{w}^{0}\hookrightarrow H^{m\mathbf{1}}$ is continuous i.e.%
\begin{equation}
\left\Vert f\right\Vert _{H^{m\mathbf{1}}}\leq\left\Vert \iota\right\Vert
\left\Vert f\right\Vert _{w,0},\quad f\in X_{w}^{0}.\label{1.092}%
\end{equation}

But from part 1 of Theorem \ref{Thm_Jg_properties} we know that for all
$\phi\in S$, $G\ast\phi\in X_{w}^{0}$ and thus $G\ast\overset{\vee}{\phi}\in
X_{w}^{0}$. Consequently, by assumption \ref{1.092},$\allowbreak\left\Vert
G\ast\overset{\vee}{\phi}\right\Vert _{H^{m\mathbf{1}}}^{2}\leq\left\Vert
\iota\right\Vert ^{2}\left\Vert G\ast\overset{\vee}{\phi}\right\Vert
_{w,0}^{2}$ i.e.
\[
\int\left(  1+\xi.\xi\right)  ^{m\mathbf{1}}\frac{1}{w^{2}}\left\vert
\phi\right\vert ^{2}\leq\left\Vert \iota\right\Vert ^{2}\int\frac{1}%
{w}\left\vert \phi\right\vert ^{2},\quad\phi\in S.
\]

With reference to weight function Definition \ref{Def_extend_wt_fn}, choose
arbitrary $\eta\notin\mathcal{A}$ where $\mathcal{A}$ is a closed set of
measure zero such that $w>0$ and continuous on $\mathbb{R}^{d}\setminus
\mathcal{A}$. Suppose $0<r<R<\operatorname*{dist}\left(  \eta,\mathcal{A}%
\right)  $. Then it is possible to choose $\phi_{r,R}\in S$ such that
$0\leq\phi_{r,R}\leq1$, $\operatorname*{supp}\phi_{r,R}\subseteq\overline
{B}_{R}\left(  \eta\right)  $ and $\phi_{r,R}=1$ on $\overline{B}_{r}\left(
\eta\right)  $ e.g. the standard "cap" example of a $C_{0}^{\infty}$ function.
Observe now that \ref{1.092} implies%
\[
\min_{B_{r}\left(  \eta\right)  }\frac{\left(  1+\xi.\xi\right)
^{m\mathbf{1}}}{w^{2}}\int\limits_{B_{r}\left(  \eta\right)  }1\leq
\int\limits_{B_{R}\left(  \eta\right)  }^{m}\left(  1+\xi.\xi\right)
^{m\mathbf{1}}\frac{\left\vert \phi_{r,R}\right\vert ^{2}}{w^{2}}%
\leq\left\Vert \iota\right\Vert ^{2}\int\limits_{B_{R}\left(  \eta\right)
}\frac{\left\vert \phi_{r,R}\right\vert ^{2}}{w}\leq\left\Vert \iota
\right\Vert ^{2}\max_{B_{R}\left(  \eta\right)  }\frac{1}{w}\int%
\limits_{B_{R}\left(  \eta\right)  }1,
\]

i.e. for $0<r<R<\operatorname*{dist}\left(  \eta,\mathcal{A}\right)  $,%
\[
\frac{\min\limits_{B_{r}\left(  \eta\right)  }\frac{\left(  1+\xi.\xi\right)
^{m\mathbf{1}}}{w^{2}}}{\max\limits_{B_{R}\left(  \eta\right)  }\frac{1}{w}%
}\leq\left\Vert \iota\right\Vert ^{2}\frac{\int_{B_{R}\left(  \eta\right)  }%
1}{\int_{B_{r}\left(  \eta\right)  }1}=\left\Vert \iota\right\Vert ^{2}\left(
\frac{R}{r}\right)  ^{d}.
\]

Now, since the continuity of $w$ near $\eta$ implies $\min\limits_{B_{R}%
\left(  \eta\right)  }\frac{1}{w}$ is continuous in $R$, we must have
$\frac{\min\limits_{B_{r}\left(  \eta\right)  }\frac{\left(  1+\xi.\xi\right)
^{m\mathbf{1}}}{w^{2}}}{\max\limits_{B_{r}\left(  \eta\right)  }\frac{1}{w}%
}\leq\left\Vert \iota\right\Vert ^{2}$ for $0<r<\operatorname*{dist}\left(
\eta,\mathcal{A}\right)  $, and hence that
\[
\frac{\left(  1+\eta.\eta\right)  ^{m\mathbf{1}}}{w\left(  \eta\right)
}=\frac{\frac{\left(  1+\eta.\eta\right)  ^{m\mathbf{1}}}{w\left(
\eta\right)  ^{2}}}{\frac{1}{w\left(  \eta\right)  }}=\lim
\limits_{r\rightarrow0}\frac{\min\limits_{B_{r}\left(  \eta\right)  }%
\frac{\left(  1+\eta.\eta\right)  ^{m\mathbf{1}}}{w^{2}}}{\max\limits_{B_{r}%
\left(  \eta\right)  }\frac{1}{w}}\leq\left\Vert \iota\right\Vert ^{2}.
\]

when $\eta\notin\mathcal{A}$ i.e. $\frac{\left(  1+\eta.\eta\right)
^{m\mathbf{1}}}{w\left(  \eta\right)  }\in L^{\infty}$ and $\left\Vert
\frac{\left(  1+\eta.\eta\right)  ^{m\mathbf{1}}}{w}\right\Vert _{\infty
}^{1/2}\leq\left\Vert \iota\right\Vert $ i.e. $\frac{\eta^{2\alpha}}{w}\in
L^{\infty}$ for $\left\vert \alpha\right\vert \leq m$.
\end{proof}

\begin{remark}
\label{Rem_Thm_Xow_embed_Hm1_iff}Compare the definition of weight function
property W03 for parameter $m\mathbf{1}$ i.e. $\frac{\xi^{2\alpha}}{w}\in
L^{1}$ for $\alpha\leq m\mathbf{1}$, with the condition $\frac{\xi^{2\alpha}%
}{w}\in L^{\infty}$ for $\alpha\leq m\mathbf{1}$ used in the last theorem.
\end{remark}

This result is an extension of Theorem \ref{Thm_Xow_embed_Hm1_iff}.

\begin{theorem}
Suppose the weight function $w$ has property W01 so that $X_{w}^{0}$ is a
normed vector space. Then $X_{w}^{0}\overset{\iota}{\hookrightarrow
}H^{s\mathbf{1}}$ iff $\frac{\left(  1+\xi.\xi\right)  ^{s\mathbf{1}}}{w}\in
L^{\infty}$.

In fact, if $\frac{\left(  1+\xi.\xi\right)  ^{s\mathbf{1}}}{w}\in L^{\infty}$
then%
\[
\left\Vert f\right\Vert _{H^{m\mathbf{1}}}\leq\left\Vert \frac{\left(
1+\xi.\xi\right)  ^{s\mathbf{1}}}{w\left(  \xi\right)  }\right\Vert _{\infty
}^{1/2}\left\Vert f\right\Vert _{w,0},\quad f\in X_{w}^{0},
\]

where $\xi.\xi$ is the componentwise vector product. Also, if the embedding
$\iota:X_{w}^{0}\hookrightarrow H^{s\mathbf{1}}$ is continuous then%
\[
\left\Vert \frac{\left(  1+\xi.\xi\right)  ^{s\mathbf{1}}}{w\left(
\xi\right)  }\right\Vert _{\infty}^{1/2}=\left\Vert \iota\right\Vert .
\]

\end{theorem}

For two weight functions:

\begin{theorem}
\label{Thm_Xw1o_embed_Xw2o_iff}Suppose the weight functions $w_{1}$ and
$w_{2}$ have property W01 i.e. $X_{w_{1}}^{0}$ and $X_{w_{2}}^{0}$ are normed
vector spaces. Then $X_{w_{1}}^{0}\overset{\iota}{\hookrightarrow}X_{w_{2}%
}^{0}$ iff $\frac{w_{2}}{w_{1}}\in L^{\infty}$. In fact, $\left\Vert
\iota\right\Vert =\left\Vert \frac{w_{2}}{w_{1}}\right\Vert _{\infty}$.
\end{theorem}

\begin{proof}
\textbf{First suppose} $\frac{w_{2}}{w_{1}}\in L^{\infty}$. Then%
\[
\left\Vert \iota f\right\Vert _{w_{2},0}^{2}=\int w_{2}\left\vert
\widehat{f}\right\vert ^{2}=\int\frac{w_{2}}{w_{1}}w_{1}\left\vert
\widehat{f}\right\vert ^{2}\leq\left\Vert \frac{w_{2}}{w_{1}}\right\Vert
_{\infty}\left\Vert f\right\Vert _{w_{1},0}^{2}<\infty,
\]

so $\left\Vert \iota\right\Vert \leq\left\Vert \frac{w_{2}}{w_{1}}\right\Vert
_{\infty}^{1/2}$.

\textbf{Next suppose} $X_{w_{1}}^{0}\subset X_{w_{2}}^{0}$ and the embedding
$\iota:X_{w_{1}}^{0}\hookrightarrow X_{w_{2}}^{0}$ is continuous i.e.%
\[
\left\Vert f\right\Vert _{w_{2},0}\leq\left\Vert \iota\right\Vert \left\Vert
f\right\Vert _{w_{1},0},\quad f\in X_{w_{1}}^{0}.
\]

But from part 1 of Theorem \ref{Thm_Jg_properties} we know that if $w_{1}$ has
basis function $G_{1}$ then for all $\phi\in S$, $G_{1}\ast\phi\in X_{w_{1}%
}^{0}$ and thus $G_{1}\ast\overset{\vee}{\phi}\in X_{w_{1}}^{0}$.
Consequently, by assumption \ref{1.091},\allowbreak\ $\left\Vert G_{1}%
\ast\overset{\vee}{\phi}\right\Vert _{w_{2},0}^{2}\leq\left\Vert
\iota\right\Vert ^{2}\left\Vert G_{1}\ast\overset{\vee}{\phi}\right\Vert
_{w_{1},0}^{2}$ i.e.
\begin{equation}
\int\frac{w_{2}}{w_{1}^{2}}\left\vert \phi\right\vert ^{2}\leq\left\Vert
\iota\right\Vert ^{2}\int\frac{1}{w_{1}}\left\vert \phi\right\vert ^{2}%
,\quad\phi\in S.\label{X98}%
\end{equation}

With reference to weight function Definition \ref{Def_extend_wt_fn}, suppose
that $\mathcal{A}_{i}$ is a closed set of measure zero such that $w_{i}>0$ and
continuous on $\mathbb{R}^{d}\setminus\mathcal{A}_{i}$. Let $\mathcal{A}%
=\mathcal{A}_{1}\cap$ $\mathcal{A}_{2}$. Choose arbitrary $\eta\notin%
\mathcal{A}$ and suppose $0<r<R<\operatorname*{dist}\left(  \eta
,\mathcal{A}\right)  $. Then it is possible to choose $\phi_{r,R}\in S$ such
that $0\leq\phi_{r,R}\leq1$, $\operatorname*{supp}\phi_{r,R}\subseteq
\overline{B}_{R}\left(  \eta\right)  $ and $\phi_{r,R}=1$ on $\overline{B}%
_{r}\left(  \eta\right)  $ e.g. the standard "cap" example of a $C_{0}%
^{\infty}$ function. Observe now that \ref{X98} implies%
\[
\min_{B_{r}\left(  \eta\right)  }\frac{w_{2}}{w_{1}^{2}}\int\limits_{B_{r}%
\left(  \eta\right)  }1\leq\int\limits_{B_{R}\left(  \eta\right)  }w_{2}%
\frac{\left\vert \phi_{r,R}\right\vert ^{2}}{w_{1}^{2}}\leq\left\Vert
\iota\right\Vert ^{2}\int\limits_{B_{R}\left(  \eta\right)  }\frac{\left\vert
\phi_{r,R}\right\vert ^{2}}{w_{1}}\leq\left\Vert \iota\right\Vert ^{2}%
\max_{B_{R}\left(  \eta\right)  }\frac{1}{w_{1}}\int\limits_{B_{R}\left(
\eta\right)  }1,
\]

i.e. for $0<r<R<\operatorname*{dist}\left(  \eta,\mathcal{A}\right)  $,%
\[
\frac{\min\limits_{B_{r}\left(  \eta\right)  }\frac{w_{2}}{w_{1}^{2}}}%
{\max\limits_{B_{R}\left(  \eta\right)  }\frac{1}{w_{1}}}\leq\left\Vert
\iota\right\Vert ^{2}\frac{\int_{B_{R}\left(  \eta\right)  }1}{\int%
_{B_{r}\left(  \eta\right)  }1}=\left\Vert \iota\right\Vert ^{2}\left(
\frac{R}{r}\right)  ^{d}.
\]

Now, since the continuity of $w_{1}$ near $\eta$ implies $\min\limits_{B_{R}%
\left(  \eta\right)  }\frac{1}{w_{1}}$ is continuous in $R$, we must have
$\frac{\min\limits_{B_{r}\left(  \eta\right)  }\frac{w_{2}}{w_{1}^{2}}}%
{\max\limits_{B_{r}\left(  \eta\right)  }\frac{1}{w_{1}}}\leq\left\Vert
\iota\right\Vert ^{2}$ for $0<r<\operatorname*{dist}\left(  \eta
,\mathcal{A}\right)  $, and hence that
\[
\frac{w_{2}\left(  \eta\right)  }{w_{1}\left(  \eta\right)  }=\frac
{\frac{w_{2}\left(  \eta\right)  }{w_{1}\left(  \eta\right)  ^{2}}}{\frac
{1}{w_{1}\left(  \eta\right)  }}=\lim\limits_{r\rightarrow0}\frac
{\min\limits_{B_{r}\left(  \eta\right)  }\frac{w_{2}}{w_{1}^{2}}}%
{\max\limits_{B_{r}\left(  \eta\right)  }\frac{1}{w_{1}}}\leq\left\Vert
\iota\right\Vert ^{2}.
\]

when $\eta\notin\mathcal{A}$ i.e. $\frac{w_{2}}{w_{1}}\in L^{\infty}$ and
$\left\Vert \frac{w_{2}}{w_{1}}\right\Vert _{\infty}^{1/2}\leq\left\Vert
\iota\right\Vert $.
\end{proof}

\begin{corollary}
Suppose the weight functions $w_{1}$ and $w_{2}$ have property W01 i.e.
$X_{w_{1}}^{0}$ and $X_{w_{2}}^{0}$ are normed vector spaces. Then $X_{w_{1}%
}^{0}\simeq X_{w_{2}}^{0}$ iff $\frac{w_{2}}{w_{1}}\in L^{\infty}$ and
$\frac{w_{1}}{w_{2}}\in L^{\infty}$. The notation $\simeq$ means equal as sets
with equivalent norms.

Further, if $X_{w_{1}}^{0}\simeq X_{w_{2}}^{0}$ then $w_{1}$ has property
W02/W03 for $\kappa$ iff $w_{2}$ has property W02/W03 for $\kappa$.
\end{corollary}

PROVE:

\begin{corollary}
?? \textbf{Show that}: $X_{w_{1}}^{0}\simeq X_{w_{2}}^{0}$ iff there exist
constants $C_{1},C_{2}>0$ such that%
\begin{equation}
C_{1}w_{1}\left(  \xi\right)  \leq w_{2}\left(  \xi\right)  \leq C_{2}%
w_{1}\left(  \xi\right)  ,\quad\xi\notin\mathcal{A}_{1}\cup\mathcal{A}%
_{2}.\label{X100}%
\end{equation}

\end{corollary}

\begin{proof}
From the proof of Theorem \ref{Thm_Xw1o_embed_Xw2o_iff} it is clear that
$X_{w_{1}}^{0}\simeq X_{w_{2}}^{0}$ implies \ref{X100}.

?? FINISH!
\end{proof}

\begin{remark}
??? PROVE!: using the notation of ?? Definition \ref{Def_SobolevSpace2}.

\begin{enumerate}
\item ?? \textbf{Prove analogous results for positive order spaces}!

\item ?? What about the dual spaces $X_{1/w}^{0}$ and $\widetilde{X}_{1/w}%
^{0}$?

\item ?? What about $S_{w,0}$? When does $w_{1}S_{w,0}=S$?

\item ?? What about the quotient space\textbf{\ }$X_{w}^{0}\left(
\Omega\right)  =r_{\Omega}X_{w}^{0}$\textbf{?}

\item What about using the density of $J_{G}$ in $X_{w}^{0}$?
\end{enumerate}
\end{remark}

\subsection{The operator $\protect\int_{\Omega}\overline{R}_{x}u$}

I have begun a study of the operator $\int_{\Omega}\overline{R}_{x}u$ where
$u\in X_{w}^{0}$ and $R_{x}$ is the Riesz representer of the evaluation
functional. These operators might be useful in deriving the Approximate smoother.

\begin{definition}
\label{Def_Op_Rho_u(x)}We define the operator $\rho_{u}\left(  x\right)
:=\int_{\Omega}\overline{R}_{x}u$ where $u\in X_{w}^{0}$ and $\Omega
\subset\mathbb{R}^{d}$ is a bounded open set.

We utilize the space $\overline{X}_{1/w}^{0}$ as introduced in Definition
\ref{Def_barXow} \textbf{below} and used for ??.

For closed $K\subset\mathbb{R}^{d}$, $\left(  \overline{X}_{1/w}^{0}\right)
_{K}:=\left\{  f\in\overline{X}_{1/w}^{0}:\operatorname*{supp}f\subset
K\right\}  $.
\end{definition}

\begin{theorem}
\label{Thm_property_integomegaRxu}\textbf{Properties of the operator }%
$\int_{\Omega}\overline{R}_{x}u$:

\begin{enumerate}
\item $\int_{\Omega}\overline{R}_{x}u$ is a continuous operator from
$X_{w}^{0}$ to $X_{w}^{0}$.

\item ?? \textbf{REDO proof of part 1 using the fact that }$S\cap X_{w}^{0}%
$\textbf{\ is dense in }$X_{w}^{0}$! See Corollary \ref{Cor_Xow_S_density}.

\item $c_{\Omega}u\in\overline{X}_{1/w}^{0}$. In fact, $\left\Vert \rho
_{u}\right\Vert _{w,0}=\left\Vert c_{\Omega}u\right\Vert _{1/w,0}$ and we have
the estimate \ref{Q015}.

\item $\int_{\Omega}\overline{R}_{x}u=0$ iff $\operatorname*{supp}%
u\subset\Omega^{c}$ iff $u\in\left(  X_{w}^{0}\right)  _{\Omega^{c}}$.

\item $\int_{\Omega}\overline{R}_{x}u$ is self-adjoint w.r.t. $X_{w}^{0}$.

\item $\left\vert \int_{\Omega}\overline{R}_{x}u\right\vert \leq\left\Vert
u\right\Vert _{w,0}\left\Vert \rho_{R_{x}}\right\Vert _{w,0}$.

\item For each $x\in\mathbb{R}^{d}$, $\int_{\Omega}\overline{R}_{x}u\in\left(
X_{w}^{0}\right)  ^{\prime}$.

\item ?? $X_{w}^{0}\left(  \Omega\right)  ^{\prime}\simeq\left\{  \rho
_{u}:u\in X_{w}^{0}\right\}  $. Use Theorem \ref{Thm_FnalQuotSp_AnnhilDenom}.
See Subsection \ref{SbSect_fnal_locXow} especially equation \ref{1.101}.
\end{enumerate}
\end{theorem}

\begin{proof}
\textbf{Part 1} If $c_{\Omega}$ denotes the characteristic function of
$\Omega$ then%
\[
\rho_{u}\left(  x\right)  =\int\overline{R}_{x}c_{\Omega}u.
\]

From Theorem \ref{Thm_ord0_Riesz_rep_W2}, $R_{x}(z)=\left(  2\pi\right)
^{-d/2}G(z-x)$ and $\overline{R_{x}(z)}=R_{z}(x)$ so that
\[
\rho_{u}\left(  x\right)  =\int\overline{R}_{x}c_{\Omega}u=\left(
2\pi\right)  ^{-d/2}\int G(x-z)\left(  c_{\Omega}u\right)  \left(  z\right)
dz.
\]

Since $G\in C_{B}^{\left(  0\right)  }\subset L^{\infty}$ and $c_{\Omega}u\in
L^{1}$ we have $\rho_{u}\left(  x\right)  \in L^{\infty}$ and%
\[
\rho_{u}=G\ast\left(  c_{\Omega}u\right)  \in L^{\infty}.
\]

We now need to use the properties of the $S^{\prime}\times O_{c}^{\prime}$
distribution convolution of Section 2.6 of Petersen \cite{Petersen83} which is
summarized below in Section \ref{Sect_temp_convol_multip}. In fact, since
$c_{\Omega}u\in\mathcal{E}^{\prime}\subset O_{c}^{\prime}$ we can use to write%
\begin{align}
\rho_{u}  & =G\ast\left(  c_{\Omega}u\right)  \in S^{\prime},\nonumber\\
\widehat{\rho_{u}}  & =\widehat{G}\left(  c_{\Omega}u\right)  ^{\wedge}%
=\frac{1}{w}\left(  c_{\Omega}u\right)  ^{\wedge}\in L_{loc}^{1},\nonumber\\
\int w\left\vert \widehat{\rho_{u}}\right\vert ^{2}  & =\int\frac{1}%
{w}\left\vert \left(  c_{\Omega}u\right)  ^{\wedge}\right\vert ^{2}%
.\label{Q005}%
\end{align}

Since $X_{w}^{0}\subset C_{B}^{\left(  0\right)  }$ it follows that
$c_{\Omega}u\in L^{1}$. Thus $\left(  c_{\Omega}u\right)  ^{\wedge}\in
C_{B}^{\left(  0\right)  }$ and so
\[
\int\frac{1}{w}\left\vert \left(  c_{\Omega}u\right)  ^{\wedge}\right\vert
^{2}\leq\left\Vert \left(  c_{\Omega}u\right)  ^{\wedge}\right\Vert _{\infty
}^{2}\int\frac{1}{w}<\infty,
\]

which means%
\[
\rho_{u}\in X_{w}^{0},\quad\left\Vert \rho_{u}\right\Vert _{w,0}\leq\left\Vert
\left(  c_{\Omega}u\right)  ^{\wedge}\right\Vert _{\infty}^{2}\int\frac{1}{w}.
\]

Further, since $c_{\Omega}\in\mathcal{E}^{\prime}$ and $u\in S^{\prime}$,
$\left(  c_{\Omega}u\right)  ^{\wedge}=\widehat{c_{\Omega}}\ast\widehat{u}$.
But $c_{\Omega}\in L^{1}$ implies $\widehat{c_{\Omega}}\in C_{B}^{\left(
0\right)  }$, and since the Cauchy-Schwartz inequality and $\widehat{u}%
=\frac{1}{\sqrt{w}}\sqrt{w}\widehat{u}$ implies $\widehat{u}\in L^{1}$, it
follows from Young's convolution inequality that%
\[
\left\Vert \left(  c_{\Omega}u\right)  ^{\wedge}\right\Vert _{\infty
}=\left\Vert \widehat{c_{\Omega}}\ast\widehat{u}\right\Vert _{\infty}%
\leq\left\Vert \widehat{c_{\Omega}}\right\Vert _{\infty}\left\Vert
\widehat{u}\right\Vert _{1}\leq\left(  2\pi\right)  ^{-d/2}\left(
\operatorname*{vol}\Omega\right)  \left(  \int\frac{1}{w}\right)
^{1/2}\left\Vert u\right\Vert _{w,0}.
\]

Hence%
\begin{align*}
\int w\left\vert \widehat{\rho_{u}}\right\vert ^{2}=\int\frac{1}{w}\left\vert
\left(  c_{\Omega}u\right)  ^{\wedge}\right\vert ^{2}  & \leq\int\frac{1}%
{w}\left(  \left(  2\pi\right)  ^{-d/2}\left(  \operatorname*{vol}%
\Omega\right)  \left(  \int\frac{1}{w}\right)  ^{1/2}\left\Vert u\right\Vert
_{w,0}\right)  ^{2}\\
& =\left(  \operatorname*{vol}\Omega\right)  ^{2}\left(  \left(  2\pi\right)
^{-\frac{d}{2}}\int\frac{1}{w}\right)  ^{2}\left\Vert u\right\Vert _{w,0}^{2},
\end{align*}

which means that%
\begin{equation}
\left\Vert \rho_{u}\right\Vert _{w,0}\leq\left(  \operatorname*{vol}%
\Omega\right)  \left(  \left(  2\pi\right)  ^{-\frac{d}{2}}\int\frac{1}%
{w}\right)  \left\Vert u\right\Vert _{w,0}=\left(  \operatorname*{vol}%
\Omega\right)  G\left(  0\right)  \left\Vert u\right\Vert _{w,0},\label{Q015}%
\end{equation}

and continuity is ensured.\medskip

\textbf{Part 2} ?? \textbf{REDO part 1 using the fact that }$S\cap X_{w}^{0}%
$\textbf{\ is dense in }$X_{w}^{0}$, which was proven in Corollary
\ref{Cor_Xow_S_density}.\medskip

\textbf{Part 3}\ Clearly $c_{\Omega}u\in\overline{X}_{1/w}^{0}$.\medskip

\textbf{Part 4 }We have
\[
\rho_{u}=0\Longleftrightarrow\widehat{\rho_{u}}=\frac{1}{w}\left(  c_{\Omega
}u\right)  ^{\wedge}=0\Longleftrightarrow\left(  c_{\Omega}u\right)  ^{\wedge
}=0\Longleftrightarrow c_{\Omega}u=0\Longleftrightarrow\operatorname*{supp}%
u\subset\Omega^{c},
\]

and the second iff follows from the definition of $\left(  X_{w}^{0}\right)
_{\Omega^{c}}$.\medskip

\textbf{Part 5} From Corollary \ref{Cor_Xow_S_density} we have that $X_{w}%
^{0}\cap S$ is dense in $X_{w}^{0}$. Hence it will suffice to prove
self-adjointness on $X_{w}^{0}\cap S$. So choose $u,v\in X_{w}^{0}\cap S$ and
from part 1,%
\[
\left(  \rho_{u},v\right)  _{w,0}=\int w\widehat{\rho_{u}}\overline
{\widehat{v}}=\int w\frac{1}{w}\left(  c_{\Omega}u\right)  ^{\wedge}%
\overline{\widehat{v}}=\int\left(  c_{\Omega}u\right)  ^{\wedge}%
\overline{\widehat{v}}.
\]

Since $c_{\Omega}u\in L^{2}$ and $v\in L^{2}$, by Plancherel's theorem%
\[
\left(  \rho_{u},v\right)  _{w,0}=\int c_{\Omega}u\overline{v}=\int
u\overline{c_{\Omega}v}=\int\widehat{u}\text{\thinspace}\overline
{\widehat{c_{\Omega}v}}=\int w\widehat{u}\overline{\frac{1}{w}%
\widehat{c_{\Omega}v}}=\int w\widehat{u}\overline{\widehat{\rho_{v}}}=\left(
u,\rho_{v}\right)  _{w,0}.
\]

\textbf{Part 6} ?? \textbf{FIX}! From part 1 and then part 5:\allowbreak
\ $\int_{\Omega}\overline{R}_{x}u=\rho_{u}\left(  x\right)  =\left(  \rho
_{u},R_{x}\right)  _{w,0}=??\left(  u,\rho_{R_{x}}\right)  _{w,0}$ and
thus\allowbreak\ $\left\vert \int_{\Omega}\overline{R}_{x}u\right\vert
\leq\left\vert \left(  u,\rho_{R_{x}}\right)  _{w,0}\right\vert \leq\left\Vert
u\right\Vert _{w,0}\left\Vert \rho_{R_{x}}\right\Vert _{w,0}$.\medskip

\textbf{Part 7} Assume $u\in S\cap X_{w}^{0}$ and follow part 2.\medskip

\textbf{Part 8} ??
\end{proof}

??

\begin{theorem}
?? If $1/w\in L^{\infty}$ the results in the last theorem simplify. Note that
we have shown that $1/w\in L^{\infty}$ iff $X_{w}^{0}\hookrightarrow L^{2}$.
What if $\xi^{\alpha}/w\in L^{\infty}$?
\end{theorem}

\begin{proof}
??
\end{proof}

\begin{corollary}
If $v\in W_{G,X}^{\Omega}$ then%
\[
\left\vert \int_{\Omega}\overline{v}u\right\vert \leq??
\]

\end{corollary}

\begin{corollary}
\ \label{Cor_integomegaRxu_locXwo}

\begin{enumerate}
\item $\int_{\Omega}\overline{R}_{x}u$ is a continuous and 1-1 operator from
$X_{w}^{0}\left(  \Omega\right)  $ to $X_{w}^{0}$.

\item $\int_{\Omega}\overline{R}_{x}^{\Omega}u$ is a continuous and 1-1
operator from $X_{w}^{0}\left(  \Omega\right)  $ to $X_{w}^{0}\left(
\Omega\right)  $.

\item The image of $X_{w}^{0}\left(  \Omega\right)  $ under $\int_{\Omega
}\overline{R}_{x}^{\Omega}u$ is a Hilbert space.

\item ?? ADJOINTS?
\end{enumerate}
\end{corollary}

\begin{proof}
\textbf{Part 1} ??

\textbf{Part 2} ??
\end{proof}

??

\begin{definition}
\label{Def_Xow_has_mollif}We say $X_{w}^{0}\left(  \Omega\right)  $
\textbf{has a mollifier} if for some $r>0$ there exists $\phi\in X_{w}%
^{0}\left(  B_{r}\right)  \cap C_{0}^{\infty}\left(  B_{r}\right)  $ such that%
\[
\phi\left(  x\right)  =\phi\left(  -x\right)  ,\text{ }\widehat{\phi}\left(
0\right)  =1,\text{ }\operatorname*{supp}\phi\subseteq\overline{B}_{r/2},
\]

and for which there exists a positive decreasing sequence $\varepsilon
_{j}\rightarrow0$ such that each $\phi_{\varepsilon_{j}}\left(  z\right)
:=\varepsilon_{j}^{-d}\phi\left(  \varepsilon_{j}^{-1}z\right)  $ also lies in
$X_{w}^{0}\left(  B_{r/2}\right)  $.
\end{definition}

\begin{remark}
\textbf{1}. The spaces $X_{w}^{0}$ that are locally Sobolev spaces contain mollifiers.

\textbf{2}. ?? The
\end{remark}

??

\begin{theorem}
\label{Thm_locintegRxS_dense_locXwo}\ 
\end{theorem}

\begin{enumerate}
\item $J_{G}=\left\{  \int\overline{R}_{x}\phi:\phi\in S\right\}  $ is dense
in $X_{w}^{0}$.

\item $J_{G}\left(  \Omega\right)  :=r_{\Omega}J_{G}$ is dense in $X_{w}%
^{0}\left(  \Omega\right)  $.

\item ?? \textbf{FIX PROOF}! If $X_{w}^{0}\left(  \Omega\right)  $ contains a
mollifier then $\left\{  r_{\Omega}\int_{\Omega}\overline{R}_{x}u:u\in
C_{0}^{\infty}\left(  \Omega\right)  \cap X_{w}^{0}\left(  \Omega\right)
\right\}  $ is dense in $X_{w}^{0}\left(  \Omega\right)  $

and $J_{G}^{\Omega}:=\left\{  r_{\Omega}\int_{\Omega}\overline{R}_{x}u:u\in
r_{\Omega}\left(  S\cap X_{w}^{0}\right)  \right\}  $ is dense in $X_{w}%
^{0}\left(  \Omega\right)  $.

\item ?? $\left\{  r_{\Omega}\int_{\Omega}\overline{R}_{x}u:u\in G\ast
S\right\}  $ is dense in $X_{w}^{0}\left(  \Omega\right)  $? Is $C^{\infty
}\left(  \overline{\Omega}\right)  \cap X_{w}^{0}\left(  \Omega\right)  $
dense in $X_{w}^{0}\left(  \Omega\right)  $? Corollary \ref{Cor_Xow_S_density}
states that $C_{0}^{\infty}\cap X_{w}^{0}$ is dense in $X_{w}^{0}$.
\end{enumerate}

\begin{proof}
\textbf{Part 1} By definition $J_{G}=\left\{  G\ast\phi:\phi\in S\right\}
$.$\ $From Theorem \ref{Thm_ord0_Riesz_rep_W2}, $R_{x}(z)=\left(  2\pi\right)
^{-d/2}G(z-x)$ and $\overline{R_{x}(z)}=R_{z}(x)$ so that
\[
\int\overline{R}_{x}\phi=\left(  2\pi\right)  ^{-d/2}\int G(x-z)\phi\left(
z\right)  dz=\left(  G\ast\phi\right)  \left(  x\right)  ,
\]

and $J_{G}=\left\{  \int\overline{R}_{x}\phi:\phi\in S\right\}  $.\medskip

\textbf{Part 2} Choose $U\in X_{w}^{0}\left(  \Omega\right)  $. Then
$r_{\Omega}^{\ast}U\in X_{w}^{0}$ and by part 1 there is a sequence $\phi
_{k}\in J_{G}$ converging to $X_{w}^{0}$. Thus $r_{\Omega}\phi_{k}\in
r_{\Omega}J_{G}$ converges to $r_{\Omega}r_{\Omega}^{\ast}U=U$.\medskip

\textbf{Part 3} Approach: approximate $u$ by $\mathcal{I}_{X}^{\Omega}u\in
W_{G}^{\Omega}\left(  \Omega\right)  $ and then approximate $\mathcal{I}%
_{X}^{\Omega}u$ by $r_{\Omega}\left(  \phi_{\varepsilon}\ast\mathcal{I}%
_{X}^{\Omega}u\right)  \in J_{G}^{\Omega}$.

Choose $u\in X_{w}^{0}\left(  \Omega\right)  $ and $0<\eta<r$. By \ref{1.089}
$W_{G}^{\Omega}\left(  \Omega\right)  $ is dense in $X_{w}^{0}\left(
\Omega\right)  $ so there exists independent data $X_{\eta}=\left(  x^{\left(
k\right)  }\right)  _{k=1}^{N}\subset\Omega$ such that
\begin{equation}
\left\Vert \mathcal{I}_{X_{\eta}}^{\Omega}u-u\right\Vert _{w,0;\Omega}%
<\eta/2.\label{Q008}%
\end{equation}

But
\[
\mathcal{I}_{X_{\eta}}^{\Omega}u=\sum\limits_{k=1}^{N}\alpha_{k}\overline
{R}_{x^{\left(  k\right)  }}^{\Omega}=r_{\Omega}\sum\limits_{k=1}^{N}%
\alpha_{k}\overline{R}_{x^{\left(  k\right)  }}=r_{\Omega}\mathcal{I}%
_{X_{\eta}}u.
\]

?? \textbf{FIX}! fix the $\delta$s and $r$s and $\varepsilon_{j}$s\ etc.

We can pointwise approximate each $\overline{R}_{x^{\left(  k\right)  }}$ by a
mollifier $\phi_{\varepsilon}\ast\overline{R}_{x^{\left(  k\right)  }}\in
X_{w}^{0}$ where $\phi_{\varepsilon}\in C_{0}^{\infty}\cap X_{w}^{0}\left(
\Omega\right)  $ approximates the delta distribution.

Set $\delta=\max\limits_{k=1}^{N}\operatorname*{dist}\left(  x^{\left(
k\right)  },\operatorname*{bdry}\Omega\right)  $ and we can choose\ $\phi
_{\varepsilon}\left(  z\right)  =\varepsilon^{-d}\phi\left(  \varepsilon
^{-1}z\right)  $ where $\phi\left(  x\right)  =\phi\left(  -x\right)  $,
$\widehat{\phi}\left(  0\right)  =1$ and $\operatorname*{supp}\phi
\subseteq\overline{B}_{1}$. Note that $\operatorname*{supp}\phi_{\varepsilon
}\subseteq\overline{B}_{\varepsilon}$. Hence%
\begin{align}
\left\Vert \phi_{\varepsilon}\ast\mathcal{I}_{X_{\eta}}u-\mathcal{I}_{X_{\eta
}}u\right\Vert _{w,0}^{2} &  =\left\Vert \phi_{\varepsilon}\ast\mathcal{I}%
_{X_{\eta}}u-\sum\limits_{k=1}^{N}\alpha_{k}\overline{R}_{x^{\left(  k\right)
}}\right\Vert _{w,0}^{2}\nonumber\\
&  =\left(  \phi_{\varepsilon}\ast\mathcal{I}_{X_{\eta}}u-\sum\limits_{k=1}%
^{N}\alpha_{k}\overline{R}_{x^{\left(  k\right)  }},\phi_{\varepsilon}%
\ast\mathcal{I}_{X_{\eta}}u-\sum\limits_{k=1}^{N}\alpha_{k}\overline
{R}_{x^{\left(  k\right)  }}\right)  _{w,0}\nonumber\\
&  =\left\Vert \phi_{\varepsilon}\ast\mathcal{I}_{X_{\eta}}u\right\Vert
_{w,0}^{2}+\left\Vert \mathcal{I}_{X_{\eta}}u\right\Vert _{w,0}^{2}%
-2\operatorname{Re}\left(  \phi_{\varepsilon}\ast\mathcal{I}_{X_{\eta}}%
u,\sum\limits_{k=1}^{N}\alpha_{k}\overline{R}_{x^{\left(  k\right)  }}\right)
_{w,0}\nonumber\\
&  =\int w\left\vert \widehat{\mathcal{I}_{X_{\eta}}u}\right\vert
^{2}\left\vert \widehat{\phi_{\varepsilon}}\right\vert ^{2}+\left\Vert
\mathcal{I}_{X_{\eta}}u\right\Vert _{w,0}^{2}-2\operatorname{Re}\int
w\widehat{\mathcal{I}_{X_{\eta}}u}\widehat{\phi_{\varepsilon}}\overline
{\left(  \sum\limits_{k=1}^{N}\alpha_{k}\overline{R}_{x^{\left(  k\right)  }%
}\right)  ^{\wedge}}\nonumber\\
&  =\int w\left\vert \widehat{\mathcal{I}_{X_{\eta}}u}\right\vert
^{2}\left\vert \widehat{\phi_{\varepsilon}}\right\vert ^{2}+\left\Vert
\mathcal{I}_{X_{\eta}}u\right\Vert _{w,0}^{2}-2\operatorname{Re}\int
w\left\vert \widehat{\mathcal{I}_{X_{\eta}}u}\right\vert ^{2}\widehat{\phi
_{\varepsilon}}\nonumber\\
&  =\int w\left\vert \widehat{\mathcal{I}_{X_{\eta}}u}\right\vert ^{2}\left(
\left\vert \widehat{\phi_{\varepsilon}}\right\vert ^{2}+1-2\operatorname{Re}%
\widehat{\phi_{\varepsilon}}\right) \nonumber\\
&  =\int w\left\vert \widehat{\mathcal{I}_{X_{\eta}}u}\right\vert
^{2}\left\vert 1-\widehat{\phi_{\varepsilon}}\right\vert ^{2}.\label{Q004}%
\end{align}

Now use the first order Taylor series remainder estimate for functions $g\in
C_{B}^{\left(  1\right)  }$:%
\[
\left\vert g\left(  0\right)  -g\left(  \xi\right)  \right\vert \leq\left\vert
\xi\right\vert \sum_{k=1}^{d}\left\Vert D_{k}g\right\Vert _{\infty},\quad
\xi\in\mathbb{R}^{d},
\]

to obtain%
\begin{align*}
\left\vert 1-\widehat{\phi_{\varepsilon}}\left(  \xi\right)  \right\vert
=\left\vert \widehat{\phi}\left(  0\right)  -\widehat{\phi}\left(
\varepsilon\xi\right)  \right\vert  & \leq\sum_{k=1}^{n}\left\Vert
D_{k}\left(  \widehat{\phi}\left(  \varepsilon\xi\right)  \right)  \right\Vert
_{\infty}\left\vert \xi\right\vert \\
& =\varepsilon\sum_{k=1}^{n}\left\Vert D_{k}\widehat{\phi}\right\Vert
_{\infty}\left\vert \xi\right\vert ,
\end{align*}

so that for all $\xi$,%
\[
\left\vert 1-\widehat{\phi}_{\varepsilon}\left(  \xi\right)  \right\vert
\leq\left\{
\begin{array}
[c]{l}%
\varepsilon\sum\limits_{k=1}^{n}\left\Vert D_{k}\widehat{\phi}\right\Vert
_{\infty}\left\vert \xi\right\vert ,\\
1+\left\Vert \widehat{\phi}\right\Vert _{\infty}.
\end{array}
\right.
\]

Hence the RHS of \ref{Q004} can be estimated as,%
\begin{align*}
&  \int w\left\vert \widehat{\mathcal{I}_{X_{\eta}}u}\right\vert
^{2}\left\vert 1-\widehat{\phi_{\varepsilon}}\right\vert ^{2}\\
&  =\int_{B_{R}}w\left\vert \widehat{\mathcal{I}_{X_{\eta}}u}\right\vert
^{2}\left\vert 1-\widehat{\phi_{\varepsilon}}\right\vert ^{2}+\int_{B_{R}^{c}%
}w\left\vert \widehat{\mathcal{I}_{X_{\eta}}u}\right\vert ^{2}\left\vert
1-\widehat{\phi_{\varepsilon}}\right\vert ^{2}\\
&  \leq\varepsilon R\left(  \sum_{k=1}^{n}\left\Vert D_{k}\widehat{\phi
}\right\Vert _{\infty}\right)  \int_{B_{R}}w\left\vert \widehat{\mathcal{I}%
_{X_{\eta}}u}\right\vert ^{2}+\left(  1+\left\Vert \widehat{\phi}\right\Vert
_{\infty}\right)  \int_{B_{R}^{c}}w\left\vert \widehat{\mathcal{I}_{X_{\eta}%
}u}\right\vert ^{2}\\
&  \leq\varepsilon R\left(  \sum_{k=1}^{n}\left\Vert D_{k}\widehat{\phi
}\right\Vert _{\infty}\right)  \left\Vert \mathcal{I}_{X_{\eta}}u\right\Vert
_{w,0}+\left(  1+\left\Vert \widehat{\phi}\right\Vert _{\infty}\right)
\int_{B_{R}^{c}}w\left\vert \widehat{\mathcal{I}_{X_{\eta}}u}\right\vert ^{2}.
\end{align*}

We can choose $R$ so that $\left(  1+\left\Vert \widehat{\phi}\right\Vert
_{\infty}\right)  \int_{B_{R}^{c}}w\left\vert \widehat{\mathcal{I}_{X_{\eta}%
}u}\right\vert ^{2}<\eta/4$ and then $\varepsilon_{\eta}>0$ so that\
\[
\varepsilon R\left(  \sum_{k=1}^{n}\left\Vert D_{k}\widehat{\phi}\right\Vert
_{\infty}\right)  \left\Vert \mathcal{I}_{X_{\eta}}u\right\Vert _{w,0}%
<\eta/4\text{ }when\text{ }\varepsilon\leq\varepsilon_{\eta}.
\]

If $\varepsilon<r$ then $\operatorname*{supp}\phi_{\varepsilon}\subset\Omega$
and so%
\[
\left\Vert \phi_{\varepsilon}\ast\mathcal{I}_{X_{\eta}}u-\mathcal{I}_{X_{\eta
}}u\right\Vert _{w,0}^{2}\leq\eta/2,\quad\varepsilon<\min\left\{
\varepsilon_{\eta},r\right\}  ,
\]

i.e.%
\begin{equation}
\left\Vert r_{\Omega}\left(  \phi_{\varepsilon}\ast\mathcal{I}_{X_{\eta}%
}u\right)  -\mathcal{I}_{X_{\eta}}^{\Omega}u\right\Vert _{w,0;\Omega}^{2}%
\leq\eta/2,\quad\varepsilon<\min\left\{  \varepsilon_{\eta},r\right\}
.\label{Q009}%
\end{equation}

Combining this estimate with that of \ref{Q008} we get%
\[
\left\Vert r_{\Omega}\left(  \phi_{\varepsilon}\ast\mathcal{I}_{X_{\eta}%
}u\right)  -u\right\Vert _{w,0;\Omega}^{2}\leq\eta,\quad\varepsilon
<\min\left\{  \varepsilon_{\eta},r\right\}  .
\]

Since $\overline{R}_{z}\left(  z^{\prime}\right)  =R_{z^{\prime}}\left(
z\right)  =\left(  2\pi\right)  ^{-\frac{d}{2}}G\left(  z-z^{\prime}\right)
$,%
\begin{align*}
\left(  \phi_{\varepsilon}\ast\mathcal{I}_{X_{\eta}}u\right)  \left(
x\right)   & =\phi_{\varepsilon}\ast\sum\limits_{k=1}^{N}\alpha_{k}%
R_{x^{\left(  k\right)  }}\left(  x\right) \\
& =\sum\limits_{k=1}^{N}\alpha_{k}\left(  \phi_{\varepsilon}\ast R_{x^{\left(
k\right)  }}\right)  \left(  x\right) \\
& =\sum\limits_{k=1}^{N}\alpha_{k}\int\phi_{\varepsilon}\left(  z\right)
R_{x^{\left(  k\right)  }}\left(  x-z\right)  dz\\
& =\left(  2\pi\right)  ^{-\frac{d}{2}}\sum\limits_{k=1}^{N}\alpha_{k}\int%
\phi_{\varepsilon}\left(  z\right)  G\left(  x-z-x^{\left(  k\right)
}\right)  dz.
\end{align*}

Now apply the change of variables $\zeta=z+x^{\left(  k\right)  }$:%
\begin{align*}
\left(  \phi_{\varepsilon}\ast\mathcal{I}_{X_{\eta}}u\right)  \left(
x\right)   & =\left(  2\pi\right)  ^{-\frac{d}{2}}\sum\limits_{k=1}^{N}%
\alpha_{k}\int\phi_{\varepsilon}\left(  \zeta-x^{\left(  k\right)  }\right)
G\left(  x-\zeta\right)  d\zeta\\
& =\sum\limits_{k=1}^{N}\alpha_{k}\int\phi_{\varepsilon}\left(  \zeta
-x^{\left(  k\right)  }\right)  \overline{R}_{x}\left(  \zeta\right)  d\zeta\\
& =\sum\limits_{k=1}^{N}\alpha_{k}\int\phi_{\varepsilon}\left(  x^{\left(
k\right)  }-\zeta\right)  \overline{R}_{x}\left(  \zeta\right)  d\zeta\\
& =\int\left(  \sum\limits_{k=1}^{N}\alpha_{k}\phi_{\varepsilon}\left(
x^{\left(  k\right)  }-\zeta\right)  \right)  \overline{R}_{x}\left(
\zeta\right)  d\zeta.
\end{align*}

Clearly $\psi_{\varepsilon}:=\sum\limits_{k=1}^{N}\alpha_{k}\phi_{\varepsilon
}\left(  x^{\left(  k\right)  }-\cdot\right)  \in C_{0}^{\infty}\left(
\Omega\right)  $ and the assumption that $X_{w}^{0}\left(  \Omega\right)  $
has a mollifier means that $\psi_{\varepsilon}\in X_{w}^{0}\left(
\Omega\right)  $ and
\[
\left(  \phi_{\varepsilon}\ast\mathcal{I}_{X_{\eta}}u\right)  \left(
x\right)  =\int\psi_{\varepsilon}\overline{R}_{x}=\int_{\Omega}\psi
_{\varepsilon}\overline{R}_{x}\in\left\{  r_{\Omega}\int_{\Omega}\overline
{R}_{x}u:u\in C_{0}^{\infty}\left(  \Omega\right)  \cap X_{w}^{0}\left(
\Omega\right)  \right\}  ,
\]

when $\varepsilon<r$.

Hence%
\begin{align*}
\phi_{\varepsilon}\ast\mathcal{I}_{X_{\eta}}u  & \rightarrow\mathcal{I}%
_{X_{\eta}}u,\\
r_{\Omega}\left(  \phi_{\varepsilon}\ast\mathcal{I}_{X_{\eta}}u\right)   &
\rightarrow r_{\Omega}\mathcal{I}_{X_{\eta}}u=\mathcal{I}_{X_{\eta}}^{\Omega
}u\rightarrow u,
\end{align*}

and so%
\[
J_{G}^{\Omega}\overset{d}{\hookrightarrow}X_{w}^{0}\left(  \Omega\right)  .
\]
\medskip

Corollary \ref{Cor_integomegaRxu_locXwo}
\end{proof}

\begin{remark}
If $X=\left\{  x^{\left(  k\right)  }\right\}  _{k=1}^{N}$ then $\sum
\limits_{k=1}^{N}\alpha_{k}\int\limits_{\Omega}\overline{R_{x^{\left(
k\right)  }}}u=\int\limits_{\Omega}\overline{g}u$ with $g\in W_{G,X}^{\Omega
}\subset W_{G}^{\Omega}$. Suppose $g=\mathcal{S}_{X}^{\Omega}f-\mathcal{I}%
_{X}^{\Omega}f\in W_{G,X}^{\Omega}$.

?? Assume we are dealing with B-splines . Thus $1/w\in L^{\infty}$ and
$X_{w}^{0}\hookrightarrow L^{2}$.

?? We use a simple variant of part 3 of Lemma \ref{Lem_centdiffop_property_2l}
i.e. $\left(  -1\right)  ^{-ld}\tbinom{2l}{l}^{-d}\delta_{2}^{2l\mathbf{1}%
}\mathcal{E}_{0}:W_{\overline{\Omega}}^{n\mathbf{1}}\rightarrow X_{w_{s}}^{0}$
is a continuous extension mapping where $\mathcal{E}_{0}$ is the zero
extension operator. Use Corollary \ref{Cor_exten_op_Wnoloc_to_X0ws}?

Choose $u\in W_{\overline{\Omega}}^{n\mathbf{1}}\left(  \Omega\right)  $ and
consider%
\[
\int_{\Omega}\left(  \overline{\mathcal{S}_{X}^{\Omega}f-\mathcal{I}%
_{X}^{\Omega}f}\right)  u,\quad u\in W_{\overline{\Omega}}^{n\mathbf{1}%
}\left(  \Omega\right)  .
\]

?? Choose $u=1$ on $B_{r}$ and $u=0$ on $B_{2r}$.

Use adjoint result?

Generalize proof of part 6 of Theorem \ref{Thm_property_integomegaRxu}?
\end{remark}

??

\begin{remark}
?? From part 1 and then part 5 of Theorem \ref{Thm_property_integomegaRxu}:
$\int_{\Omega}\overline{R}_{x}u=\rho_{u}\left(  x\right)  =\left(  \rho
_{u},R_{x}\right)  _{w,0}=\left(  u,\rho_{R_{x}}\right)  _{w,0}$ and thus
$\left\vert \int_{\Omega}\overline{R}_{x}u\right\vert \leq\left\vert \left(
u,\rho_{R_{x}}\right)  _{w,0}\right\vert \leq\left\Vert u\right\Vert
_{w,0}\left\Vert \rho_{R_{x}}\right\Vert _{w,0}$.

If $u\in X_{w}^{0}\left(  \Omega\right)  $ and $g\in W_{G,X}^{\Omega}\left(
\Omega\right)  $, say $g=\sum\limits_{k=1}^{N}\alpha_{k}R_{x^{\left(
k\right)  }}^{\Omega}$, then%
\[
\int_{\Omega}\overline{g}u=\sum\limits_{k=1}^{N}\alpha_{k}\int\limits_{\Omega
}\overline{R_{x^{\left(  k\right)  }}}u=\sum\limits_{k=1}^{N}\alpha_{k}%
\rho_{u}\left(  x^{\left(  k\right)  }\right)  =\sum\limits_{k=1}^{N}%
\alpha_{k}\left(  \rho_{u},R_{x^{\left(  k\right)  }}\right)  _{w,0;\Omega
}=\left(  u,\sum\limits_{k=1}^{N}\alpha_{k}\rho_{R_{x^{\left(  k\right)  }}%
}\right)  _{w,0;\Omega}=\left(  u,\rho_{g}\right)  _{w,0;\Omega}.
\]

Thus%
\[
\left\vert \int_{\Omega}\overline{g}u\right\vert \leq\left\Vert u\right\Vert
_{w_{s},0;\Omega}\left\Vert \rho_{g}\right\Vert _{w_{s},0;\Omega}.
\]

Use Corollary \ref{Cor_exten_op_Wnoloc_to_X0ws}. Choose $f\in W_{\overline
{\Omega}}^{n\mathbf{1}}$ and set $u=Ef\in X_{w_{s}}^{0}$. Then%
\[
\left\vert \int_{\Omega}\overline{g}f\right\vert \leq\left\Vert Ef\right\Vert
_{w_{s},0}\left\Vert \rho_{g}\right\Vert _{w_{s},0;\Omega},
\]

where%
\begin{align*}
\widehat{Ef}  & =4^{ld}\tbinom{2l}{l}^{-d}\left(  \sin m\xi_{k}\right)
^{2l\mathbf{1}}\widehat{\mathcal{E}_{0}f},\\
\left\Vert Ef\right\Vert _{w_{s},0}  & \leq\left(  4m\right)  ^{ld}\tbinom
{2l}{l}^{-d}\left\Vert D^{n\mathbf{1}}f\right\Vert _{L^{2}\left(
\Omega\right)  }.
\end{align*}

Thus
\[
\left\vert \int_{\Omega}\overline{g}f\right\vert \leq\left(  4m\right)
^{ld}\tbinom{2l}{l}^{-d}\left\Vert D^{n\mathbf{1}}f\right\Vert _{L^{2}\left(
\Omega\right)  }\left\Vert \rho_{g}\right\Vert _{w_{s},0;\Omega},\quad f\in
W_{\overline{\Omega}}^{n\mathbf{1}}.
\]

Choose $f$ to be a tensor product of identical 1-dimensional functions. Scale
using $\sigma_{\varepsilon}f$.

Suppose that $g$ is real-valued.

Choose $x\in\Omega$ and suppose that $g\left(  x\right)  >0$ on
$B_{\varepsilon_{0}}$. Assume $0<\varepsilon\leq\varepsilon_{0}$. Choose
$f\geq0$ so that its support is $B_{2\varepsilon}$ and it is $1$ on
$B_{\varepsilon}$. Then
\[
\int_{B_{\varepsilon}}\overline{g}\leq\left(  4m\right)  ^{ld}\tbinom{2l}%
{l}^{-d}\left\Vert D^{n\mathbf{1}}\left(  f\left(  \frac{x}{\varepsilon
}\right)  \right)  \right\Vert _{B_{2\varepsilon}}\left\Vert \rho
_{g}\right\Vert _{w_{s},0;\Omega}.
\]

From part 3 of ??: $\left\Vert \rho_{u}\right\Vert _{w,0}=\left\Vert
c_{\Omega}u\right\Vert _{1/w,0}$. Thus%
\[
\int_{B_{\varepsilon}}\overline{g}\leq\left(  4m\right)  ^{ld}\tbinom{2l}%
{l}^{-d}\left\Vert D^{n\mathbf{1}}\left(  f\left(  \frac{x}{\varepsilon
}\right)  \right)  \right\Vert _{B_{2\varepsilon}}\left\Vert c_{\Omega
}g\right\Vert _{1/w,0}.
\]

\end{remark}

\section{ The Approximate smoother\label{Sect_ApproxSmooth}}

From part 5 of Summary \ref{Sum_Exact_smth_properties} the matrix equation for
the Exact smoother is $\left(  N\rho I+R_{X,X}\right)  \alpha=y$. The
construction and solution of this system does not generate a scalable matrix
algorithm since the size of the smoothing matrix is $N\times N$ i.e. its size
depends on the number of data points whereas it will turn out that the matrix
size for a scalable algorithm is independent of the number of data points. In
this section we will overcome this limitation and derive the Approximate
smoother problem by discretizing the Exact smoothing problem on a grid. This
will be done by using the space $J_{G}=G\ast S$ which is dense in $X_{w}^{0}$
and was introduced in Definition \ref{Def_Jg}. The space $J_{G}$ is not
necessary for the specification of the Approximate smoother problem but this
is how I derived the Approximate smoother algorithm: I came across the
convolution space $J_{G}$ in Dyn's review paper \cite{Dyn89} and decided to
approximate its functions using a regular rectangular grid $X^{\prime}$ and
the Trapezoidal rule. This led to functions in the finite dimensional basis
function space $W_{G,X^{\prime}}$. Minimizing the Exact smoother functional
$J_{e}\left[  \cdot\right]  $ given by \ref{7.63} over the functions in
$W_{G,X^{\prime}}$ yields the Approximate smoother matrix equation
\ref{a7.094}. The size of this matrix is independent of the number of data
points and hence the construction and solution of the Approximate smoother
matrix equation is a \textbf{scalable algorithm}. In fact, the size of the
matrix is equal to the number of grid points and so increases exponentially
with the number of dimensions. In practice we are limited to two or three dimensions.

In actual fact, the set $X^{\prime}$ will be generalized from a grid to any
set of distinct points. This could be a sparse grid, for example.

\subsection{\textbf{Approach 1} is an analog of the Exact smoothing problem
\label{SbSect_ap_deriv_Approx_smth}}

In this subsection we will provide some justification for approximating the
infinite dimensional Hilbert space $X_{w}^{0}$ by a finite dimensional
subspace $W_{G,X^{\prime}}$, where $X^{\prime}$ is a regular, rectangular grid
of points in $\mathbb{R}^{d}$. The space $W_{G,X^{\prime}}$ will be used to
define the Approximate smoothing problem. The set $X^{\prime}$ will then be
generalized to include any set of distinct points. This discretization process
turns out to be similar to that described in Garcke and Griebel
\cite{GarckGrieb05}.

\begin{definition}
\label{Def_mesh}\textbf{A regular, rectangular grid in }$\mathbb{R}^{d}$

Let the grid occupy a rectangle $R\left(  a,b\right)  $ where $a.<b$. Suppose
the grid has $\mathcal{N}^{\prime}=\left(  N_{1}^{\prime},N_{2}^{\prime
},\ldots,N_{d}^{\prime}\right)  $ points in each dimension and let
$h\in\mathbb{R}^{d}$ denote the grid sizes.

Then $X^{\prime}=\left\{  x_{\alpha}^{\prime}=a+h\alpha\mid\alpha\in
\mathbb{Z}^{d}\text{ and }0\leq\alpha<\mathcal{N}^{\prime}\right\}  $ is the
set of grid points.

Let $N^{^{\prime}}$ be the number of grid points so that $N^{^{\prime}%
}=\left(  \mathcal{N}^{\prime}\right)  ^{\mathbf{1}}=\mathcal{N}_{1}^{\prime
}\mathcal{N}_{2}^{\prime}\ldots\mathcal{N}_{d}^{\prime}$, and of course we
have the constraint $\mathcal{N}^{^{\prime}}h=b-a$.
\end{definition}

By Theorem \ref{Thm_Jg_properties} the space $J_{G}=G\ast S$ is dense in
$X_{w}^{0}$ and so we will approximate the functions in $S$ using the
rectangular grid $X^{\prime}$ defined above. Our analysis will be matrix-based
so we choose an order for the grid points and set $X^{^{\prime}}=\left\{
x_{n}^{\prime}\right\}  _{n=1}^{N^{^{\prime}}}$. Integrals on $\mathbb{R}^{d}$
will be approximated by integrals on the grid region using the trapezoidal
rule i.e.
\begin{equation}
\int\limits_{grid}f\left(  x\right)  dx\simeq h^{\mathbf{1}}\sum
\limits_{n=1}^{N^{^{\prime}}}f\left(  x_{n}^{\prime}\right)  ,\label{a7.096}%
\end{equation}

where $h^{\mathbf{1}}=h_{1}h_{2}\times\ldots\times h_{d}$ is the volume of a
grid element. This will be a two stage approximation: first a restriction to
the grid rectangle and then an application of the trapezoidal rule.\medskip

\textbf{Step1} Approximation of the functions in $J_{G}=G\ast S$ by functions
in $W_{G,X^{\prime}}$.\medskip

Suppose $\phi\in S$. Then the trapezoidal approximation \ref{a7.096} on the
grid $X^{\prime}$ gives%
\begin{equation}
G\ast\phi=\left(  2\pi\right)  ^{-\frac{d}{2}}\int\limits_{\mathbb{R}^{d}%
}G\left(  \cdot-y\right)  \phi\left(  y\right)  dy\simeq\left(  2\pi\right)
^{-\frac{d}{2}}\sum_{n=1}^{N^{\prime}}G\left(  \cdot-x_{n}^{\prime}\right)
\left(  h^{\mathbf{1}}\phi\left(  x_{n}^{\prime}\right)  \right)
.\label{a7.034}%
\end{equation}

The equations and approximations \ref{a7.096} and \ref{a7.034} suggest we
approximate the space $G\ast S$ by functions of the form $\sum\limits_{n=1}%
^{N^{\prime}}G\left(  x-x_{n}^{\prime}\right)  \alpha_{n}$, $\alpha_{n}%
\in\mathbb{C}$ i.e. by functions in $W_{G,X^{\prime}}$.\medskip

\textbf{Step 2} With this motivation we could now specify a smoothing problem
which we will call an Approximate smoothing problem. This would involve
minimizing the Exact smoothing functional \ref{7.63} over $W_{G,X^{\prime}}$
where $X^{\prime}$ is a rectangular grid. However, since the space
$W_{G,X^{\prime}}$ is defined when $X^{\prime}$ is any set of distinct points,
we will define the following more general problem:%

\begin{equation}%
\begin{tabular}
[c]{|c|}\hline
\textbf{The Approximate smoothing problem}\\\hline
Minimize the Exact smoothing functional $J_{e}\left[  f\right]  $\ for $f\in
W_{G,X^{\prime}}$,\\
where $X^{\prime}$ is a set of distinct points in $\mathbb{R}^{d}$.\\\hline
\end{tabular}
\label{a7.5}%
\end{equation}

\begin{remark}
?? Use part 3 of Theorem \ref{Thm_locintegRxS_dense_locXwo}?
\end{remark}

\subsection{\textbf{Approach 2}}

In Theorem \ref{Thm_estim_Exact_smth_by_Wgx} it was shown that for any $f\in
X_{w}^{0}$,%
\[
\left\vert \mathcal{S}_{X}f\left(  x\right)  -\frac{\left(
\widetilde{\mathcal{E}}_{X}R_{x}\right)  ^{T}\widetilde{\mathcal{E}}_{X}%
f}{\left(  R_{0}\left(  0\right)  +\rho\right)  N}\right\vert \leq\frac
{2R_{0}\left(  0\right)  ^{3/2}}{R_{0}\left(  0\right)  +\rho}\min\left\{
1,\frac{R_{0}\left(  0\right)  }{\rho}\right\}  \left\Vert f\right\Vert
_{w,0},\quad x\in\mathbb{R}^{d},
\]

and it was noted in Remark \ref{Rem_Thm_estim_Exact_smth_by_Wgx} that
$\frac{\left(  \widetilde{\mathcal{E}}_{X}R_{x}\right)  ^{T}%
\widetilde{\mathcal{E}}_{X}f}{\left(  R_{0}\left(  0\right)  +\rho\right)
N}\in W_{G,X}$ and this approximation is bounded uniformly on $\mathbb{R}^{d}%
$. Of course $R_{x}\left(  y\right)  =\left(  2\pi\right)  ^{d/2}G\left(
x-y\right)  $. Now suppose our independent data points $X$ are always
contained in a bounded data region $\Omega$. Overlay $\Omega$ with a regular,
rectangular grid and let $X^{\prime}=\left\{  x_{n}^{\prime}\right\}
_{n=1}^{N^{^{\prime}}}\subset\Omega$ denote the grid points in $\Omega$.
Further \textbf{assume that} $G$ and $f$ are Lipschitz continuous on
$\mathbb{R}^{d}$.

Since $f$ is Lipschitz continuous on $\mathbb{R}^{d}$ all the derivatives
exist and are bounded a.e. Noting \ref{av075} we try assuming that $f\in
X_{w}^{0}\cap W^{1,\infty}$ and $G\in W^{1,\infty}$ and for each $x$,
$R_{x}\in W^{1,\infty}$. Indeed
\begin{equation}
R_{x}f\in W^{1,\infty}\cap C_{B}^{\left(  0\right)  }.\label{a7.25}%
\end{equation}

Next we want to express $\left(  \widetilde{\mathcal{E}}_{X}R_{x}\right)
^{T}\widetilde{\mathcal{E}}_{X}f$ in terms of $\left(  \widetilde{\mathcal{E}%
}_{X^{\prime}}R_{x}\right)  ^{T}\widetilde{\mathcal{E}}_{X^{\prime}}f$ and a
remainder term and then estimate the remainder term. But%
\begin{equation}
\frac{1}{N}\left(  \widetilde{\mathcal{E}}_{X}R_{x}\right)  ^{T}%
\widetilde{\mathcal{E}}_{X}f=\frac{1}{N}\sum_{m=1}^{N}R_{x}\left(  x^{\left(
m\right)  }\right)  f\left(  x^{\left(  m\right)  }\right)  =\frac{1}{N}%
\sum_{k=1}^{N^{\prime}}\sum_{x^{\left(  m\right)  }\in X_{k}^{\prime}}%
R_{x}\left(  x^{\left(  m\right)  }\right)  f\left(  x^{\left(  m\right)
}\right)  ,\label{a7.23}%
\end{equation}

where
\[
X_{k}^{\prime}=\left\{  x^{\left(  m\right)  }\in X:x^{\left(  m\right)
}-x_{k}^{\prime}\in R\left[  \mathbf{0};h/2\right)  \right\}  ,
\]

means that $\left\{  X_{k}^{\prime}\right\}  _{k=1}^{N^{\prime}}$ partitions
$X$.

The properties \ref{a7.25} of $R_{x}f$ mean that the distribution Taylor
series expansion of Lemma \ref{Lem_Taylor_extension} can now be brought into
play to yield%
\[
R_{x}\left(  x^{\left(  m\right)  }\right)  f\left(  x^{\left(  m\right)
}\right)  =R_{x}\left(  x_{k}^{\prime}\right)  f\left(  x_{k}^{\prime}\right)
+\mathcal{R}_{1}\left(  R_{x}f\right)  \left(  x_{k}^{\prime},x^{\left(
m\right)  }-x_{k}^{\prime}\right)  ,
\]

and for $x^{\left(  m\right)  }\in X_{k}^{\prime}$,%
\begin{align}
\left\vert \left(  \mathcal{R}_{1}\left(  R_{x}f\right)  \right)  \left(
x_{k}^{\prime},x^{\left(  m\right)  }-x_{k}^{\prime}\right)  \right\vert  &
\leq\sqrt{d}\max_{\left\vert \beta\right\vert =1}\left\Vert D^{\beta}\left(
R_{x}f\right)  \right\Vert _{\infty}\left\vert x^{\left(  m\right)  }%
-x_{k}^{\prime}\right\vert \nonumber\\
&  \leq\frac{\sqrt{d}}{2}\max_{\left\vert \beta\right\vert =1}\left\Vert
D^{\beta}\left(  R_{x}f\right)  \right\Vert _{\infty}\left\vert h\right\vert
\nonumber\\
&  \leq\frac{\sqrt{d}}{2}\max_{\left\vert \beta\right\vert \leq1}\left\Vert
D^{\beta}R_{x}\right\Vert _{\infty}\max_{\left\vert \beta\right\vert \leq
1}\left\Vert D^{\beta}f\right\Vert _{\infty}\left\vert h\right\vert
.\label{a7.59}%
\end{align}

Thus%
\begin{align*}
\sum_{x^{\left(  m\right)  }\in X_{k}^{\prime}}R_{x}\left(  x^{\left(
m\right)  }\right)  f\left(  x^{\left(  m\right)  }\right)   & =\sum
_{x^{\left(  m\right)  }\in X_{k}^{\prime}}R_{x}\left(  x_{k}^{\prime}\right)
f\left(  x_{k}^{\prime}\right)  +\sum_{x^{\left(  m\right)  }\in X_{k}%
^{\prime}}\mathcal{R}_{1}\left(  R_{x}f\right)  \left(  x_{k}^{\prime
},x^{\left(  m\right)  }-x_{k}^{\prime}\right) \\
& =N_{k}^{\prime}R_{x}\left(  x_{k}^{\prime}\right)  f\left(  x_{k}^{\prime
}\right)  +\sum_{x^{\left(  m\right)  }\in X_{k}^{\prime}}\mathcal{R}%
_{1}\left(  R_{x}f\right)  \left(  x_{k}^{\prime},x^{\left(  m\right)  }%
-x_{k}^{\prime}\right)  ,
\end{align*}

and $\frac{1}{N}\left(  \widetilde{\mathcal{E}}_{X}R_{x}\right)
^{T}\widetilde{\mathcal{E}}_{X}f$ can now be written%
\[
\frac{1}{N}\left(  \widetilde{\mathcal{E}}_{X}R_{x}\right)  ^{T}%
\widetilde{\mathcal{E}}_{X}f=\sum_{k=1}^{N^{\prime}}R_{x}\left(  x_{k}%
^{\prime}\right)  \frac{N_{k}^{\prime}}{N}f\left(  x_{k}^{\prime}\right)
+\frac{1}{N}\sum_{k=1}^{N^{\prime}}\sum_{x^{\left(  m\right)  }\in
X_{k}^{\prime}}\mathcal{R}_{1}\left(  R_{x}f\right)  \left(  x_{k}^{\prime
},x^{\left(  m\right)  }-x_{k}^{\prime}\right)  ,
\]

with%
\[
\sum_{k=1}^{N^{\prime}}R_{x}\left(  x_{k}^{\prime}\right)  \frac{N_{k}%
^{\prime}}{N}f\left(  x_{k}^{\prime}\right)  \in W_{G,X^{\prime}}.
\]

Consequently the remainder estimate \ref{a7.59} implies%
\begin{align}
\left\vert \frac{1}{N}\left(  \widetilde{\mathcal{E}}_{X}R_{x}\right)
^{T}\widetilde{\mathcal{E}}_{X}f-\sum_{k=1}^{N^{\prime}}R_{x}\left(
x_{k}^{\prime}\right)  \frac{N_{k}^{\prime}}{N}f\left(  x_{k}^{\prime}\right)
\right\vert  & \leq\frac{1}{N}\sum_{k=1}^{N^{\prime}}\sum_{x^{\left(
m\right)  }\in X_{k}^{\prime}}\left\vert \mathcal{R}_{1}\left(  R_{x}f\right)
\left(  x_{k}^{\prime},x^{\left(  m\right)  }-x_{k}^{\prime}\right)
\right\vert \nonumber\\
& \leq\frac{1}{N}\sum_{k=1}^{N^{\prime}}\sum_{x^{\left(  m\right)  }\in
X_{k}^{\prime}}\frac{\sqrt{d}}{2}\max_{\left\vert \beta\right\vert \leq
1}\left\Vert D^{\beta}R_{x}\right\Vert _{\infty}\max_{\left\vert
\beta\right\vert \leq1}\left\Vert D^{\beta}f\right\Vert _{\infty}\left\vert
h\right\vert \nonumber\\
& =\frac{\sqrt{d}}{2}\max_{\left\vert \beta\right\vert \leq1}\left\Vert
D^{\beta}R_{x}\right\Vert _{\infty}\max_{\left\vert \beta\right\vert \leq
1}\left\Vert D^{\beta}f\right\Vert _{\infty}\left\vert h\right\vert
,\label{a7.19}%
\end{align}

giving a linear dependence on the grid size.

The last estimate can be combined with \ref{15.19} to give%
\begin{align}
\left\vert \mathcal{S}_{X}f\left(  x\right)  -\frac{\sum\limits_{k=1}%
^{N^{\prime}}R_{x}\left(  x_{k}^{\prime}\right)  \frac{N_{k}^{\prime}}%
{N}f\left(  x_{k}^{\prime}\right)  }{R_{0}\left(  0\right)  +\rho}\right\vert
&  \leq\frac{2R_{0}\left(  0\right)  ^{3/2}}{R_{0}\left(  0\right)  +\rho}%
\min\left\{  1,\frac{R_{0}\left(  0\right)  }{\rho}\right\}  \left\Vert
f\right\Vert _{w,0}+\nonumber\\
&  \qquad+\frac{\sqrt{d}}{2}\max_{\left\vert \beta\right\vert \leq1}\left\Vert
D^{\beta}R_{x}\right\Vert _{\infty}\max_{\left\vert \beta\right\vert \leq
1}\left\Vert D^{\beta}f\right\Vert _{\infty}\left\vert h\right\vert ,\quad
x\in\mathbb{R}^{d},\label{a7.28}%
\end{align}

so that as $N\rightarrow\infty$ and $\left\vert h\right\vert \rightarrow0$
this approximation is bounded pointwise on $\mathbb{R}^{d}$.

\textbf{Now if }$w\in W02$\textbf{\ for some }$\kappa\geq1$ then $X_{w}%
^{0}\subset C_{B}^{\left(  1\right)  }$, $G\in C_{B}^{\left(  2\right)  }$ and
so $D^{\beta}f$ and $D^{\beta}G$ are always Lipschitz continuous when
$\left\vert \beta\right\vert \leq1$. This certainly applies to the shifted
thin-plate splines and the Gaussian.

We will now try to improve this approximation by replacing the vector $\left(
\frac{N_{k}^{\prime}}{N}f\left(  x_{k}^{\prime}\right)  \right)  \left(
R_{0}\left(  0\right)  +\rho\right)  ^{-1}$ by an arbitrary complex vector
$\alpha^{\prime}\in\mathbb{C}^{N^{\prime}}$ and then finding a uniform bound
for $\left\vert \mathcal{S}_{X}f\left(  x\right)  -\left(
\widetilde{\mathcal{E}}_{X^{\prime}}^{\ast}\alpha^{\prime}\right)  \left(
x\right)  \right\vert $ when $x\in\mathbb{R}^{d}$ and $f\in X_{w}^{0}$. In
fact we use the standard approach of writing
\[
\left\vert \mathcal{S}_{X}f\left(  x\right)  -\left(  \widetilde{\mathcal{E}%
}_{X^{\prime}}^{\ast}\alpha^{\prime}\right)  \left(  x\right)  \right\vert
=\left\vert \left(  \mathcal{S}_{X}f-\widetilde{\mathcal{E}}_{X^{\prime}%
}^{\ast}\alpha^{\prime},R_{x}\right)  _{w,0}\right\vert \leq\left\Vert
\mathcal{S}_{X}f-\widetilde{\mathcal{E}}_{X^{\prime}}^{\ast}\alpha^{\prime
}\right\Vert _{w,0}\left\Vert R_{0}\right\Vert ,
\]

and then show there is a unique $\alpha^{\prime}$ which minimizes $\left\Vert
\mathcal{S}_{X}f-\widetilde{\mathcal{E}}_{X^{\prime}}^{\ast}\alpha^{\prime
}\right\Vert _{w,0}$. This is just the adjoint formulation of the minimum
seminorm interpolation problem and the calculations are simplified since $G$
is assumed to be real-valued.

\begin{theorem}
If $g\in X_{w}^{0}$ then $R_{X,X}^{-1}\widetilde{\mathcal{E}}_{X^{\prime}%
}g=\operatorname*{argmin}\limits_{\alpha^{\prime}\in\mathbb{C}^{N^{\prime}}%
}\left\Vert g-\widetilde{\mathcal{E}}_{X^{\prime}}^{\ast}\alpha^{\prime
}\right\Vert _{w,0}$ is unique and $\widetilde{\mathcal{E}}_{X^{\prime}}%
^{\ast}R_{X,X}^{-1}\widetilde{\mathcal{E}}_{X^{\prime}}g=\mathcal{I}%
_{X^{\prime}}g$.
\end{theorem}

\begin{proof}
If $\alpha^{\prime}\in\mathbb{C}^{N^{\prime}}$ then%
\begin{align*}
\left\Vert g-\widetilde{\mathcal{E}}_{X^{\prime}}^{\ast}\alpha^{\prime
}\right\Vert _{w,0}^{2}  & =\left(  g-\widetilde{\mathcal{E}}_{X^{\prime}%
}^{\ast}\alpha^{\prime},g-\widetilde{\mathcal{E}}_{X^{\prime}}^{\ast}%
\alpha^{\prime}\right)  _{w,0}\\
& =\left\Vert g\right\Vert _{w,0}^{2}-\left(  g,\widetilde{\mathcal{E}%
}_{X^{\prime}}^{\ast}\alpha^{\prime}\right)  _{w,0}-\left(
\widetilde{\mathcal{E}}_{X^{\prime}}^{\ast}\alpha^{\prime},g\right)
_{w,0}+\left(  \widetilde{\mathcal{E}}_{X^{\prime}}^{\ast}\alpha^{\prime
},\widetilde{\mathcal{E}}_{X^{\prime}}^{\ast}\alpha^{\prime}\right)  _{w,0}\\
& =\left\Vert g\right\Vert _{w,0}^{2}-2\operatorname{Re}\left(
g,\widetilde{\mathcal{E}}_{X^{\prime}}^{\ast}\alpha^{\prime}\right)
_{w,0}+\left(  \widetilde{\mathcal{E}}_{X^{\prime}}^{\ast}\alpha^{\prime
},\widetilde{\mathcal{E}}_{X^{\prime}}^{\ast}\alpha^{\prime}\right)  _{w,0}\\
& =\left\Vert g\right\Vert _{w,0}^{2}-2\operatorname{Re}\left(
\widetilde{\mathcal{E}}_{X^{\prime}}g,\alpha^{\prime}\right)  +\left(
\alpha^{\prime},\widetilde{\mathcal{E}}_{X^{\prime}}\widetilde{\mathcal{E}%
}_{X^{\prime}}^{\ast}\alpha^{\prime}\right) \\
& =\left\Vert g\right\Vert _{w,0}^{2}-2\operatorname{Re}\left(
\widetilde{\mathcal{E}}_{X^{\prime}}g,\alpha^{\prime}\right)  +\left(
\alpha^{\prime},R_{X,X}\alpha^{\prime}\right)  ,
\end{align*}

and if we substitute $\alpha^{\prime}=\beta^{\prime}+i\gamma^{\prime}$ and
differentiate w.r.t. $\beta^{\prime}$ and then $\gamma^{\prime}$ we obtain a
unique minimum at $\alpha^{\prime}=R_{X,X}^{-1}\widetilde{\mathcal{E}%
}_{X^{\prime}}g$.
\end{proof}

This theorem implies that the closest point in $W_{G,X^{\prime}}$ to
$\mathcal{S}_{X}f$ under $\left\Vert \cdot\right\Vert _{w,0}$ is
$\mathcal{I}_{X^{\prime}}\mathcal{S}_{X}f$, and so we define the Approximate
smoother operator $\mathcal{S}_{X}^{a}$ by%
\begin{equation}
\mathcal{S}_{X}^{a}=\mathcal{I}_{X^{\prime}}\mathcal{S}_{X}\label{a7.65}%
\end{equation}

It was observed that the Approximate smoother problem \ref{a7.5} still makes
sense if $X^{\prime}$ is any set of distinct points in $\mathbb{R}^{d}$ and
clearly Definition \ref{a7.65} is still meaningful if $X^{\prime}$ is any set
of distinct points. In part 2 of Corollary \ref{Cor_Apprx_smth_prop} we will
show that the solution of the Approximate smoother problem is the Approximate
smoother \ref{a7.65}.

\section{Preparation: the Exact smoother mapping and data functions}

This section presents a convenient summary of the Exact smoother mapping
properties required for our study of the Approximate smoother problem
\ref{a7.5}.

The Exact smoother problem was described in the Introduction and was solved in
Chapter \ref{Ch_Exact_smth} using a geometric Hilbert space framework which
involved the introduction of the Hilbert product space $V$ and the operator
$\mathcal{L}_{X}:X_{w}^{0}\rightarrow V$:

\begin{definition}
\textbf{The space} $V$ \textbf{and the operator} $\mathcal{L}_{X}$

\begin{enumerate}
\item Let $V=X_{w}^{0}\otimes\mathbb{C}^{N}$ be the Hilbert product space with
norm $\left\Vert \cdot\right\Vert _{V}$ and inner product $\left(  \cdot
,\cdot\right)  _{V}$ given by
\[
\left(  \left(  u_{1},u_{2}\right)  ,\left(  v_{1},v_{2}\right)  \right)
_{V}=\rho\,\left(  u_{1},v_{1}\right)  _{w,0}+\frac{1}{N}\left(  u_{2}%
,v_{2}\right)  _{\mathbb{C}^{N}}.
\]

\item Let the operator $\mathcal{L}_{X}:X_{w}^{0}\rightarrow V$ be defined by
$\mathcal{L}_{X}f=\left(  f,\widetilde{\mathcal{E}}_{X}f\right)  $ where
$\widetilde{\mathcal{E}}_{X}$ is the vector-valued evaluation function
$\widetilde{\mathcal{E}}_{X}f=\left(  f\left(  x^{\left(  k\right)  }\right)
\right)  $ of Definition \ref{Def_vect_val_eval_op}.
\end{enumerate}
\end{definition}

The Exact smoothing functional $J_{e}$ can now be written%
\begin{equation}
J_{e}[f]=\left\Vert \mathcal{L}_{X}f-\left(  0,y\right)  \right\Vert _{V}%
^{2},\quad f\in X_{w}^{0},\label{a7.15}%
\end{equation}

and the Exact smoother $s_{e}$ is the unique orthogonal projection of $\left(
0,y\right)  $ onto the infinite dimensional subspace $\mathcal{L}_{X}\left(
X_{w}^{0}\right)  $. Using this approach it was shown in part 4 of Theorem
\ref{Thm_ex_smooth_var} that%
\begin{equation}
s_{e}=\frac{1}{N}\left(  \mathcal{L}_{X}^{\ast}\mathcal{L}_{X}\right)
^{-1}\widetilde{\mathcal{E}}_{X}^{\ast}y,\quad y\in\mathbb{C}^{N}%
.\label{a7.14}%
\end{equation}

\begin{definition}
\label{Def_data_func_Exact_smth_map}\textbf{Data functions and the Exact
smoother mapping} $\mathcal{S}_{X}:X_{w}^{0}\rightarrow W_{G,X}$.

Given an independent data set $X$ each member of $X_{w}^{0}$ is assumed to act
as a legitimate data function $f$ and generate a dependent data vector
$\widetilde{\mathcal{E}}_{X}f$. Further, equation \ref{a7.14} enables us to
define a mapping $\mathcal{S}_{X}:X_{w}^{0}\rightarrow W_{G,X}$ from the data
functions to their corresponding Exact smoothers defined by%
\begin{equation}
\mathcal{S}_{X}f=\frac{1}{N}\left(  \mathcal{L}_{X}^{\ast}\mathcal{L}%
_{X}\right)  ^{-1}\widetilde{\mathcal{E}}_{X}^{\ast}\widetilde{\mathcal{E}%
}_{X}f,\quad f\in X_{w}^{0}.\label{a7.16}%
\end{equation}

By Corollary \ref{Cor_ex_Lx*Lx_onto_Xw,th}, $\mathcal{L}_{X}^{\ast}%
\mathcal{L}_{X}$ is a homeomorphism from $X_{w}^{0}$ to $X_{w}^{0}$ so
$\mathcal{S}_{X}$ is a continuous linear mapping.
\end{definition}

\begin{summary}
\label{Sum_Exact_smth_properties}This is a list of some of the properties of
the Exact smoother mapping $\mathcal{S}_{X}$ which will be used in this document:

\begin{enumerate}
\item $\mathcal{S}_{X}$ maps $X_{w}^{0}$ onto $W_{G,X}$ and
$\operatorname*{null}\mathcal{S}_{X}=W_{G,X}^{\perp}$. Also, $\mathcal{S}_{X}$
is self-adjoint but not a projection.

\item $\left\Vert \mathcal{S}_{X}f\right\Vert _{w,0}\leq\left\Vert
f\right\Vert _{w,0}$ and $\left\Vert \left(  I-\mathcal{S}_{X}\right)
f\right\Vert _{w,0}\leq\left\Vert f\right\Vert _{w,0}$,$\quad f\in X_{w}^{0} $.

\item $\mathcal{S}_{X}f=f-\rho\left(  \mathcal{L}_{X}^{\ast}\mathcal{L}%
_{X}\right)  ^{-1}f$,$\quad f\in X_{w}^{0}$.

\item $\mathcal{S}_{X}f=\widetilde{\mathcal{E}}_{X}^{\ast}\left(  N\rho
I+R_{X,X}\right)  ^{-1}\widetilde{\mathcal{E}}_{X}f$,$\quad f\in X_{w}^{0}%
$.\medskip

\item $\mathcal{S}_{X}f=\sum\limits_{k=1}^{N}\alpha_{k}R_{x^{\left(  k\right)
}}$, where $\alpha=\left(  \alpha_{k}\right)  $ satisfies the matrix equation

$\left(  N\rho I+R_{X,X}\right)  \alpha=y$, where $y=\widetilde{\mathcal{E}%
}_{X}f$ is the dependent data.
\end{enumerate}
\end{summary}

\begin{proof}
We just list the references from Chapter \ref{Ch_Exact_smth}:\medskip

\textbf{Part 1} Corollary \ref{Cor_ex_Lx*Lx_onto_Xw,th}; \textbf{Part 2} Part
2 of Theorem \ref{Thm_ex_Min_Smth_in_Wgx}; \textbf{Part 3} Corollary
\ref{Cor_ex_Lx*Lx_onto_Xw,th}; \textbf{Parts 4 and 5} Part 3 of Theorem
\ref{Thm_Exact_smth_mat_eqn}.
\end{proof}

\section{Properties of the Approximate smoother problem}

It will be shown in part 1 of the next theorem that the Approximate smoothing
problem has a unique solution, which we will call the Approximate smoother.

\begin{theorem}
\label{Thm_Appr_smth_prop1}Denote the data for the Approximate smoothing
problem by $\left[  X;y\right]  $ and $X^{\prime}$. Set $\varsigma=\left(
0,y\right)  $ so that by \ref{a7.15}, $J_{e}\left[  f\right]  =\left\Vert
\mathcal{L}_{X}f-\varsigma\right\Vert _{V}^{2}$ is the Exact smoothing functional.

Then the Approximate smoothing problem has a unique solution $s_{a}\in
W_{G,X^{\prime}}$ which satisfies:

\begin{enumerate}
\item $J_{e}\left[  s_{a}\right]  <J_{e}\left[  f\right]  $,\quad for all
$f\in W_{G,X^{\prime}}$ and $f\neq s_{a}$.

\item $\left(  \mathcal{L}_{X}s_{a}-\varsigma,\mathcal{L}_{X}s_{a}%
-\mathcal{L}_{X}f\right)  _{V}=0$,\quad for all $f\in W_{G,X^{\prime}}$.

\item $\left\Vert \mathcal{L}_{X}s_{a}-\varsigma\right\Vert _{V}%
^{2}+\left\Vert \mathcal{L}_{X}s_{a}-\mathcal{L}_{X}f\right\Vert _{V}%
^{2}=\left\Vert \mathcal{L}_{X}f-\varsigma\right\Vert _{V}^{2}$,\quad for all
$f\in W_{G,X^{\prime}}$.\medskip

Note that the equations of parts 2 and 3 are directly equivalent.\medskip

\item $\left(  \mathcal{L}_{X}^{\ast}\mathcal{L}_{X}s_{a}-\frac{1}%
{N}\widetilde{\mathcal{E}}_{X}^{\ast}y,\,f\right)  _{w,0}=0$,\quad for all
$f\in W_{G,X^{\prime}}$.
\end{enumerate}
\end{theorem}

\begin{proof}
\textbf{Part 1} We want to show that there is a unique function in
$W_{G,X^{\prime}}$ which minimizes the Exact smoothing functional
$J_{e}\left[  \cdot\right]  $ over $W_{G,X^{\prime}}$. Since $W_{G,X^{\prime}%
}$ is a finite dimensional subspace of $X_{w}^{0}$, its image under
$\mathcal{L}_{X}$ must be a finite dimensional subspace of $V$ and hence a
closed subspace of $V$.

Consequently, there exists a unique element of $\mathcal{L}_{X}\left(
W_{G,X^{\prime}}\right)  $, say $v$, which is the projection of $\zeta$ onto
$\mathcal{L}_{X}\left(  W_{G,X^{\prime}}\right)  $, such that $\left\Vert
v-\varsigma\right\Vert _{V}<\left\Vert \mathcal{L}_{X}f-\varsigma\right\Vert
_{V}$ for all $f\in W_{G,X^{\prime}}$ such that $\mathcal{L}_{X}\left(
f\right)  \neq v$.

Since $\mathcal{L}_{X}$ is 1-1 on $X_{w}^{0}$ there exists a unique element of
$W_{G,X^{\prime}}$, call it $s_{a}$, such that $v=\mathcal{L}_{X}\left(
s_{a}\right)  $.

In terms of $J_{e}$ we have $J_{e}\left[  s_{a}\right]  <J_{e}\left[
f\right]  $\quad for all $f\in W_{G,X^{\prime}}$ and $f\neq s_{a}$.\medskip

\textbf{Parts 2 and 3} Since $v$ is the projection of $\zeta$ onto
$\mathcal{L}_{X}\left(  W_{G,X^{\prime}}\right)  $ simple Hilbert space
geometry yields the equivalent equations of parts 2 and 3.\medskip

\textbf{Part 4} We start with the equation of part 2 and use the fact from
part 3 of Theorem \ref{Thm_ex_L_op_properties} that $\mathcal{L}_{X}^{\ast
}u=\rho u_{1}+\frac{1}{N}\widetilde{\mathcal{E}}_{X}^{\ast}u_{2}$, $\left(
u_{1},u_{2}\right)  \in V$. Noting that $X$ has $N$ data points we have for
all $f\in W_{G,X^{\prime}}$
\begin{align*}
0=\left(  \mathcal{L}_{X}s_{a}-\varsigma,\mathcal{L}_{X}s_{a}-\mathcal{L}%
_{X}f\right)  _{V} &  =\left(  \mathcal{L}_{X}s_{a}-\varsigma,\mathcal{L}%
_{X}\left(  s_{a}-f\right)  \right)  _{V}\\
&  =\left(  \mathcal{L}_{X}^{\ast}\mathcal{L}_{X}s_{a}-\mathcal{L}_{X}^{\ast
}\varsigma,s_{a}-f\right)  _{w,0}\\
&  =\left(  \mathcal{L}_{X}^{\ast}\mathcal{L}_{X}s_{a}-\frac{1}{N}%
\widetilde{\mathcal{E}}_{X}^{\ast}y,s_{a}-f\right)  _{w,0}.
\end{align*}

But $s_{a}\in W_{G,X^{\prime}}$ so%
\[
\left(  \mathcal{L}_{X}^{\ast}\mathcal{L}_{X}s_{a}-\frac{1}{N}%
\widetilde{\mathcal{E}}_{X}^{\ast}y,f\right)  _{w,0}=0,\quad f\in
W_{G,X^{\prime}},
\]

and thus $\mathcal{L}_{X}^{\ast}\mathcal{L}_{X}s_{a}-\frac{1}{N}%
\widetilde{\mathcal{E}}_{X}^{\ast}y\in W_{G,X^{\prime}}^{\bot}$ by definition
of the orthogonal complement $W_{G,X^{\prime}}^{\bot}$.
\end{proof}

Part 4 is used to prove the next two corollaries.

\begin{corollary}
\label{Cor_Apprx_smth_prop}Suppose $s_{e}$ is the Exact smoother of the data
$\left[  X,y\right]  $. Suppose $s_{a}$ is the Approximate smoother of the
data $\left[  X,y\right]  $ generated by the points $X^{\prime}$. Then:

\begin{enumerate}
\item $s_{e}-s_{a}\in W_{G,X^{\prime}}^{\bot}$ i.e. $\widetilde{\mathcal{E}%
}_{X^{\prime}}s_{e}=\widetilde{\mathcal{E}}_{X^{\prime}}s_{a}$.

\item $s_{a}=\widetilde{\mathcal{E}}_{X^{\prime}}^{\ast}\left(
\widetilde{\mathcal{E}}_{X^{\prime}}\widetilde{\mathcal{E}}_{X^{\prime}}%
^{\ast}\right)  ^{-1}\widetilde{\mathcal{E}}_{X^{\prime}}s_{e}%
=\widetilde{\mathcal{E}}_{X^{\prime}}^{\ast}R_{X^{\prime},X^{\prime}}%
^{-1}\widetilde{\mathcal{E}}_{X^{\prime}}s_{e}=\mathcal{I}_{X^{\prime}}s_{e}$.

\item $s_{a}=\frac{1}{N}\widetilde{\mathcal{E}}_{X^{\prime}}^{\ast
}R_{X^{\prime},X^{\prime}}^{-1}\widetilde{\mathcal{E}}_{X^{\prime}}\left(
\mathcal{L}_{X}^{\ast}\mathcal{L}_{X}\right)  ^{-1}\widetilde{\mathcal{E}}%
_{X}^{\ast}y$.
\end{enumerate}
\end{corollary}

\begin{proof}
\textbf{Part 1} From the proof of part 4 of the previous theorem
$\mathcal{L}_{X}^{\ast}\mathcal{L}_{X}s_{a}-\frac{1}{N}\widetilde{\mathcal{E}%
}_{X}^{\ast}y\in W_{G,X^{\prime}}^{\bot}$ so $\mathcal{L}_{X}^{\ast
}\mathcal{L}_{X}s_{a}-\frac{1}{N}\widetilde{\mathcal{E}}_{X}^{\ast}y=g$ for
some $g\in W_{G,X^{\prime}}^{\bot}$ and $s_{a}-\frac{1}{N}\left(
\mathcal{L}_{X}^{\ast}\mathcal{L}_{X}\right)  ^{-1}\widetilde{\mathcal{E}}%
_{X}^{\ast}y=\left(  \mathcal{L}_{X}^{\ast}\mathcal{L}_{X}\right)  ^{-1}g$.
But by equation \ref{a7.16}, $s_{e}=\frac{1}{N}\left(  \mathcal{L}_{X}^{\ast
}\mathcal{L}_{X}\right)  ^{-1}\widetilde{\mathcal{E}}_{X}^{\ast}y$ and so
$s_{a}-s_{e}=\left(  \mathcal{L}_{X}^{\ast}\mathcal{L}_{X}\right)  ^{-1}g$.

Further, from part 3 of Summary \ref{Sum_Exact_smth_properties},
$\mathcal{S}_{X}g=g-\rho\,\left(  \mathcal{L}_{X}^{\ast}\mathcal{L}%
_{X}\right)  ^{-1}g$ and thus

$\left(  \mathcal{L}_{X}^{\ast}\mathcal{L}_{X}\right)  ^{-1}g=\frac{1}{\rho}g
$ since $\operatorname*{null}\mathcal{S}_{X}=W_{G,X^{\prime}}^{\bot}$. Hence
$s_{a}-s_{e}=\frac{1}{\rho}g\in W_{G,X^{\prime}}^{\bot}$ which is
characterized by $\widetilde{\mathcal{E}}_{X^{\prime}}\left(  s_{e}%
-s_{a}\right)  =0$ i.e. $\widetilde{\mathcal{E}}_{X^{\prime}}s_{e}%
=\widetilde{\mathcal{E}}_{X^{\prime}}s_{a}$.\medskip

\textbf{Part 2} From Theorem \ref{Thm_eval_op_properties} we know
$\widetilde{\mathcal{E}}_{X}^{\ast}:\mathbb{C}^{N}\rightarrow W_{G,X}$ and
$\widetilde{\mathcal{E}}_{X^{\prime}}^{\ast}:\mathbb{C}^{N^{\prime}%
}\rightarrow W_{G,X^{\prime}}$ are onto. Hence, since $s_{a}\in W_{G,X^{\prime
}}$ and $s_{e}\in W_{G,X}$ we have $s_{e}=\widetilde{\mathcal{E}}_{X}^{\ast
}\alpha$ and $s_{a}=\widetilde{\mathcal{E}}_{X^{\prime}}^{\ast}\beta$ for some
$\alpha\in\mathbb{C}^{N}$ and $\beta\in\mathbb{C}^{N^{\prime}}$. Parts 1 and 5
of Theorem \ref{Thm_eval_op_properties} now imply
\[
0=\widetilde{\mathcal{E}}_{X^{\prime}}\left(  s_{e}-s_{a}\right)
=\widetilde{\mathcal{E}}_{X^{\prime}}\left(  s_{e}-\widetilde{\mathcal{E}%
}_{X^{\prime}}^{\ast}\beta\right)  =\widetilde{\mathcal{E}}_{X^{\prime}}%
s_{e}-\widetilde{\mathcal{E}}_{X^{\prime}}\widetilde{\mathcal{E}}_{X^{\prime}%
}^{\ast}\beta=\widetilde{\mathcal{E}}_{X^{\prime}}s_{e}-R_{X^{\prime
},X^{\prime}}\beta,
\]

so that $\beta=R_{X^{\prime},X^{\prime}}^{-1}\widetilde{\mathcal{E}%
}_{X^{\prime}}s_{e}$ and $s_{a}=\widetilde{\mathcal{E}}_{X^{\prime}}^{\ast
}\beta=\widetilde{\mathcal{E}}_{X^{\prime}}^{\ast}R_{X^{\prime},X^{\prime}%
}^{-1}\widetilde{\mathcal{E}}_{X^{\prime}}s_{e}$, as required.

The last equality follows from the formula for the minimal norm interpolant
operator $\mathcal{I}_{X^{\prime}}$ given by equation \ref{1.08}, namely
$\mathcal{I}_{X^{\prime}}f=\widetilde{\mathcal{E}}_{X^{\prime}}^{\ast}\left(
R_{X^{\prime},X^{\prime}}\right)  ^{-1}\widetilde{\mathcal{E}}_{X^{\prime}}%
f$.\medskip

\textbf{Part 3} From the proof of part 1, $s_{e}=\frac{1}{N}\left(
\mathcal{L}_{X}^{\ast}\mathcal{L}_{X}\right)  ^{-1}\widetilde{\mathcal{E}}%
_{X}^{\ast}y $. Hence by part 2
\begin{align*}
s_{a}=\widetilde{\mathcal{E}}_{X^{\prime}}^{\ast}R_{X^{\prime},X^{\prime}%
}^{-1}\widetilde{\mathcal{E}}_{X^{\prime}}s_{e} &  =\widetilde{\mathcal{E}%
}_{X^{\prime}}^{\ast}R_{X^{\prime},X^{\prime}}^{-1}\widetilde{\mathcal{E}%
}_{X^{\prime}}\frac{1}{N}\left(  \mathcal{L}_{X}^{\ast}\mathcal{L}_{X}\right)
^{-1}\widetilde{\mathcal{E}}_{X}^{\ast}y\\
&  =\frac{1}{N}\widetilde{\mathcal{E}}_{X^{\prime}}^{\ast}R_{X^{\prime
},X^{\prime}}^{-1}\widetilde{\mathcal{E}}_{X^{\prime}}\left(  \mathcal{L}%
_{X}^{\ast}\mathcal{L}_{X}\right)  ^{-1}\widetilde{\mathcal{E}}_{X}^{\ast}y.
\end{align*}

\end{proof}

The following corollary contains results which are analogous to some of those
derived in Theorem \ref{Thm_ex_Min_Smth_in_Wgx} for the Exact smoother.

\begin{corollary}
\label{Cor_Appr_smth_prop2}Suppose $s_{a}\in W_{G,X^{\prime}}$ is the (unique)
Approximate smoother of the data $\left[  X,y\right]  $ generated by the
points $X^{\prime}$. Then:

\begin{enumerate}
\item $\left(  s_{a},f\right)  _{w,0}=\frac{1}{N\rho}\left(  y-\left(
s_{a}\right)  _{X}\right)  ^{T}\overline{f_{X}}$, when $f\in W_{G,X^{\prime}}$.

\item $\left\Vert s_{a}\right\Vert _{w,0}^{2}=\frac{1}{N\rho}\left(  y-\left(
s_{a}\right)  _{X}\right)  ^{T}\overline{\left(  s_{a}\right)  _{X}}$.

\item $J_{e}\left(  s_{a}\right)  =\frac{1}{N}\left(  y-\left(  s_{a}\right)
_{X}\right)  ^{T}\overline{y}$.
\end{enumerate}
\end{corollary}

\begin{proof}
\textbf{Part 1}. From the definition of $\mathcal{L}_{X}^{\ast}\mathcal{L}%
_{X}$
\[
\mathcal{L}_{X}^{\ast}\mathcal{L}_{X}s_{a}-\frac{1}{N}\widetilde{\mathcal{E}%
}_{X}^{\ast}y=\rho s_{a}+\frac{1}{N}\widetilde{\mathcal{E}}_{X}^{\ast
}\widetilde{\mathcal{E}}_{X}s_{a}-\frac{1}{N}\widetilde{\mathcal{E}}_{X}%
^{\ast}y=\rho s_{a}+\frac{1}{N}\widetilde{\mathcal{E}}_{X}^{\ast}\left(
\left(  s_{a}\right)  _{X}-y\right)  .
\]

Thus, when $f\in W_{G,X^{\prime}}$%
\begin{align*}
0=\left(  \mathcal{L}_{X}^{\ast}\mathcal{L}_{X}s_{a}-\frac{1}{N}%
\widetilde{\mathcal{E}}_{X}^{\ast}y,f\right)  _{w,0} &  =\left(  \rho
s_{a},f\right)  _{w,0}+\frac{1}{N}\left(  \widetilde{\mathcal{E}}_{X}^{\ast
}\left(  \left(  s_{a}\right)  _{X}-y\right)  ,f\right)  _{w,0}\\
&  =\rho\left(  s_{a},f\right)  _{w,0}+\frac{1}{N}\left(  \left(
s_{a}\right)  _{X}-y,\widetilde{\mathcal{E}}_{X}f\right)  _{\mathbb{C}^{N}}\\
&  =\rho\left(  s_{a},f\right)  _{w,0}+\frac{1}{N}\left(  \left(
s_{a}\right)  _{X}-y,f_{X}\right)  _{\mathbb{C}^{N}}\\
&  =\rho\left(  s_{a},f\right)  _{w,0}+\frac{1}{N}\left(  \left(
s_{a}\right)  _{X}-y\right)  ^{T}\overline{f_{X}}.
\end{align*}
\medskip

\textbf{Part 2}. If $f=s_{a}$ in the equation proved in part 1, then
\[
0=\rho\left\Vert s_{a}\right\Vert _{w,0}^{2}+\frac{1}{N}\left(  \left(
s_{a}\right)  _{X}-y,\left(  s_{a}\right)  _{X}\right)  _{\mathbb{C}^{N}}%
=\rho\left\Vert s_{a}\right\Vert _{w,0}^{2}+\frac{1}{N}\left(  \left(
s_{a}\right)  _{X}-y\right)  ^{T}\overline{\left(  s_{a}\right)  _{X}},
\]

so that $\left\Vert s_{a}\right\Vert _{w,0}^{2}=\frac{1}{\rho N}\left(
y-\left(  s_{a}\right)  _{X}\right)  ^{T}\overline{\left(  s_{a}\right)  _{X}%
}$.\medskip

\textbf{Part 3}.
\begin{align*}
J_{e}\left(  s\right)   & =\rho\left\Vert s_{a}\right\Vert _{w,0}^{2}+\frac
{1}{N}\left\vert \left(  s_{a}\right)  _{X}-y\right\vert ^{2}\\
& =\rho\left\Vert s_{a}\right\Vert _{w,0}^{2}+\frac{1}{N}\left(  \left(
s_{a}\right)  _{X}^{T}-y^{T}\right)  \left(  \overline{\left(  s_{a}\right)
_{X}}-\overline{y}\right) \\
& =\rho\left\Vert s_{a}\right\Vert _{w,0}^{2}+\frac{1}{N}\left(  \left(
s_{a}\right)  _{X}^{T}\overline{\left(  s_{a}\right)  _{X}}-\left(
s_{a}\right)  _{X}^{T}\overline{y}-y^{T}\overline{\left(  s_{a}\right)  _{X}%
}+y^{T}\overline{y}\right) \\
& =\frac{1}{N}\left(  y-\left(  s_{a}\right)  _{X}\right)  ^{T}\overline
{\left(  s_{a}\right)  _{X}}+\frac{1}{N}\left(  \left(  s_{a}\right)  _{X}%
^{T}\overline{\left(  s_{a}\right)  _{X}}-\left(  s_{a}\right)  _{X}%
^{T}\overline{y}-y^{T}\overline{\left(  s_{a}\right)  _{X}}+y^{T}\overline
{y}\right) \\
& =\frac{1}{N}\left(  y^{T}\overline{\left(  s_{a}\right)  _{X}}-\left(
s_{a}\right)  _{X}^{T}\overline{\left(  s_{a}\right)  _{X}}+\left(
s_{a}\right)  _{X}^{T}\overline{\left(  s_{a}\right)  _{X}}-\left(
s_{a}\right)  _{X}^{T}\overline{y}-y^{T}\overline{\left(  s_{a}\right)  _{X}%
}+y^{T}\overline{y}\right) \\
& =\frac{1}{N}\left(  -\left(  s_{a}\right)  _{X}^{T}\overline{y}%
+y^{T}\overline{y}\right) \\
& =\frac{1}{N}\left(  y-\left(  s_{a}\right)  _{X}\right)  ^{T}\overline{y}.
\end{align*}

\end{proof}

In a manner analogous to the minimal norm interpolant mapping $\mathcal{I}%
_{X}$ of Definition \ref{Def_data_func_interpol_map} and the Exact smoother
mapping $\mathcal{S}_{X}$ of Definition \ref{Def_data_func_exact_map} which
map data functions to their interpolant and smoother respectively, we will now
define the Approximate smoother mapping:

\begin{definition}
\label{Def_Apprx_smth_map}\textbf{The Approximate smoother mapping}
$\mathcal{S}_{X,X^{\prime}}^{a}:X_{w}^{0}\rightarrow W_{G,X^{\prime}}$

Given an independent data set $X$, we shall assume that each member of
$X_{w}^{0}$ can act as a legitimate data function $f_{d}$ and generate the
data vector $\widetilde{\mathcal{E}}_{X}f_{d}$.

The equation of part 2 of Corollary \ref{Cor_Apprx_smth_prop} enables us to
define a continuous linear mapping $\mathcal{S}_{X,X^{\prime}}^{a}:X_{w}%
^{0}\rightarrow W_{G,X^{\prime}}$ from the data functions to the corresponding
Approximate smoother of the data $\left[  X,\widetilde{\mathcal{E}}_{X}%
f_{d}\right]  $ generated by the points $X^{\prime}$. This mapping is given by%
\[
\mathcal{S}_{X,X^{\prime}}^{a}f_{d}=\frac{1}{N}\widetilde{\mathcal{E}%
}_{X^{\prime}}^{\ast}R_{X^{\prime},X^{\prime}}^{-1}\widetilde{\mathcal{E}%
}_{X^{\prime}}\left(  \mathcal{L}_{X}^{\ast}\mathcal{L}_{X}\right)
^{-1}\widetilde{\mathcal{E}}_{X}^{\ast}\widetilde{\mathcal{E}}_{X}f_{d},\quad
f_{d}\in X_{w}^{0}.
\]

\end{definition}

We now prove some properties of the Approximate smoother mapping.

\begin{corollary}
\label{Cor_Appr_smth_map_prop}Suppose $\mathcal{S}_{X,X^{\prime}}^{a}$ is the
Approximate smoother mapping. Then $\mathcal{S}_{X,X^{\prime}}^{a}%
=\mathcal{I}_{X^{\prime}}\mathcal{S}_{X}$ where $\mathcal{I}_{X^{\prime}}$ is
the minimal norm interpolant mapping with independent data $X^{\prime}$ and
$\mathcal{S}_{X}$ is the Exact smoother mapping with independent data $X$.

Further, $\mathcal{S}_{X,X^{\prime}}^{a}$ is self-adjoint and $\left\Vert
\mathcal{S}_{X,X^{\prime}}^{a}f\right\Vert _{w,0}\leq\left\Vert f\right\Vert
_{w,0}$ for all $f\in X_{w}^{0}$.
\end{corollary}

\begin{proof}
That $\mathcal{S}_{X,X^{\prime}}^{a}=\mathcal{I}_{X^{\prime}}\mathcal{S}_{X}$
follows immediately from \ref{a7.16} for $\mathcal{S}_{X}$ and equation
\ref{1.08} for $\mathcal{I}_{X^{\prime}}$. $\mathcal{S}_{X,X^{\prime}}^{a}$ is
self-adjoint since $\mathcal{I}_{X^{\prime}}$ and $\mathcal{S}_{X}$ are both
self-adjoint. Finally, from part 2 Summary \ref{Sum_Exact_smth_properties}%
,$\left\Vert \mathcal{S}_{X,X^{\prime}}^{a}f\right\Vert _{w,0}\leq\left\Vert
f\right\Vert _{w,0}$ since $\left\Vert \mathcal{S}_{X}f\right\Vert _{w,0}%
\leq\left\Vert f\right\Vert _{w,0}$ and $\left\Vert \mathcal{I}_{X^{\prime}%
}f\right\Vert _{w,0}\leq\left\Vert f\right\Vert _{w,0}$ for all $f\in
X_{w}^{0}$.
\end{proof}

In the last corollary it was shown that the Approximate smoother is an
interpolant of the Exact smoother. Hence the Approximate smoother uses the
same information about the data function as the Exact smoother does i.e. it
involves no additional data function evaluations.

\subsection{The matrix equation for the Approximate
smoother\label{SbSect_Approx_smth_mat_eqn}}

We now know the Approximate smoother exists and is unique. The next step is to
derive its matrix equation.

\begin{theorem}
\label{Thm_Approx_smth_mat_eqn}Suppose $s\in W_{G,X^{\prime}}$ is the (unique)
Approximate smoother of the data $\left[  X,y\right]  $ generated by the
distinct points $X^{\prime}=\left\{  x_{i}^{\prime}\right\}  _{i=1}%
^{N^{\prime}}$.

Then $s\left(  x\right)  =\sum\limits_{i=1}^{N^{\prime}}\alpha_{i}^{\prime
}R_{x_{i}^{\prime}}\left(  x\right)  $, where $\alpha^{\prime}=\left(
\alpha_{i}^{\prime}\right)  \in\mathbb{C}^{N^{\prime}}$ and $\alpha^{\prime}$
satisfies
\begin{equation}
\left(  N\rho R_{X^{\prime},X^{\prime}}+R_{X,X^{\prime}}^{T}R_{X,X^{\prime}%
}\right)  \alpha^{\prime}=R_{X,X^{\prime}}^{T}y,\label{a7.094}%
\end{equation}

where $R_{X,X}=\left(  2\pi\right)  ^{-\frac{d}{2}}G_{X,X}$ and
$R_{X,X^{\prime}}=\left(  2\pi\right)  ^{-\frac{d}{2}}G_{X,X^{\prime}}$.

The matrix on the left of equation \ref{a7.094} will be called the Approximate
smoother matrix and it will usually be denoted by the symbol $\Psi$.
\end{theorem}

\begin{proof}
This proof uses the results $\widetilde{\mathcal{E}}_{X^{\prime}%
}\widetilde{\mathcal{E}}_{X}^{\ast}=R_{X^{\prime},X}$ and
$\widetilde{\mathcal{E}}_{X}\widetilde{\mathcal{E}}_{X^{\prime}}^{\ast
}=R_{X,X^{\prime}}$ which follow from part 6 of Theorem
\ref{Thm_eval_op_properties}. From Theorem \ref{Thm_Appr_smth_prop1}, $\left(
\mathcal{L}_{X}^{\ast}\mathcal{L}_{X}s-\frac{1}{N}\widetilde{\mathcal{E}}%
_{X}^{\ast}y,\,f\right)  _{w,0}=0$ for all $f\in W_{G,X^{\prime}}$.
Substituting $f=R_{x_{i}^{\prime}}$ yields $\widetilde{\mathcal{E}}%
_{X^{\prime}}\left(  \mathcal{L}_{X}^{\ast}\mathcal{L}_{X}s-\frac{1}%
{N}\widetilde{\mathcal{E}}_{X}^{\ast}y\right)  =0$, and since by part 4
Theorem \ref{Thm_ex_L_op_properties}, $\mathcal{L}_{X}^{\ast}\mathcal{L}%
_{X}s=\rho s+\frac{1}{N}\widetilde{\mathcal{E}}_{X}^{\ast}%
\widetilde{\mathcal{E}}_{X}s$, we have
\begin{align}
0=\widetilde{\mathcal{E}}_{X^{\prime}}\left(  \mathcal{L}_{X}^{\ast
}\mathcal{L}_{X}s-\frac{1}{N}\widetilde{\mathcal{E}}_{X}^{\ast}y\right)   &
=\widetilde{\mathcal{E}}_{X^{\prime}}\left(  \rho s+\frac{1}{N}%
\widetilde{\mathcal{E}}_{X}^{\ast}\widetilde{\mathcal{E}}_{X}s-\frac{1}%
{N}\widetilde{\mathcal{E}}_{X}^{\ast}y\right) \nonumber\\
&  =\rho\widetilde{\mathcal{E}}_{X^{\prime}}s+\frac{1}{N}R_{X^{\prime}%
,X}\widetilde{\mathcal{E}}_{X}s-\frac{1}{N}R_{X^{\prime},X}y\nonumber\\
&  =N\rho s_{X^{\prime}}+R_{X^{\prime},X}s_{X}-R_{X^{\prime},X}%
y.\label{a7.093}%
\end{align}

The next step is to write $s_{X^{\prime}}$ and $s_{X}$ in terms of basis
functions. But $s\left(  x\right)  =\sum\limits_{i=1}^{N^{\prime}}\alpha
_{i}^{\prime}R_{x_{i}^{\prime}}\left(  x\right)  =\widetilde{\mathcal{E}%
}_{X^{\prime}}^{\ast}\alpha^{\prime}$ implies%
\begin{equation}
s_{X}=R_{X,X^{\prime}}\alpha^{\prime},\qquad s_{X^{\prime}}=R_{X^{\prime
},X^{\prime}}\alpha^{\prime},\label{a7.084}%
\end{equation}

so that \ref{a7.093} becomes%
\[
0=N\rho R_{X^{\prime},X^{\prime}}\alpha^{\prime}+R_{X^{\prime},X}%
R_{X,X^{\prime}}\alpha^{\prime}-R_{X^{\prime},X}y,
\]

which can be rearranged to give the desired matrix equation \ref{a7.094}.
\end{proof}

\begin{remark}
\ 

\begin{enumerate}
\item Note that the reason the matrix equation \ref{a7.094} can be more
concisely written in matrices based on $R_{x}$ than the basis function $G$ is
that I have chosen the `symmetric' Fourier transform pair $\widehat{f}%
(\xi)=\left(  2\pi\right)  ^{-\frac{d}{2}}\int e^{-ix\xi}f(x)dx$ and
$\overset{\vee}{f}(\xi)=\left(  2\pi\right)  ^{-\frac{d}{2}}\int e^{ix\xi
}f(x)dx$.

\item The `basis function form' of the Approximate smoother matrix equation
\ref{a7.094} is%
\[
\left(  \left(  2\pi\right)  ^{\frac{d}{2}}N\rho G_{X^{\prime},X^{\prime}%
}+G_{X,X^{\prime}}^{T}G_{X,X^{\prime}}\right)  \alpha^{\prime}=G_{X,X^{\prime
}}^{T}y,
\]

with the smoother given by%
\[
s\left(  x\right)  =\sum\limits_{i=1}^{N^{\prime}}\alpha_{i}^{\prime}G\left(
x-x_{i}^{\prime}\right)  .
\]

\end{enumerate}
\end{remark}

In the next theorem we prove some properties of the Approximate smoother matrix.

\begin{theorem}
\label{Thm_ap_mat_properties}The Approximate smoother matrix $\Psi$ specified
in Theorem \ref{Thm_Approx_smth_mat_eqn} has the following properties:

\begin{enumerate}
\item $\Psi$ is real valued, symmetric, positive definite and regular.

\item $\Psi$ has size $N^{\prime}\times N^{\prime}$. Hence the size of $\Psi$
is independent of the number of (scattered) data points.
\end{enumerate}
\end{theorem}

\begin{proof}
\textbf{Part 1}. Since $G$ is a Hermitian function and is assumed to be real
valued in this document, it follows that $R_{X^{\prime},X^{\prime}}$ and
$R_{X^{\prime},X}$ have real elements, $R_{X^{\prime},X^{\prime}}$ is
symmetric and so $\Psi$ is real and symmetric. Further, if $\alpha
\in\mathbb{C}^{N}$
\begin{align*}
\alpha^{T}\Psi\overline{\alpha}=\alpha^{T}\left(  N\rho R_{X^{\prime
},X^{\prime}}+R_{X,X^{\prime}}^{T}R_{X,X^{\prime}}\right)  \overline{\alpha}
&  =N\rho\alpha^{T}R_{X^{\prime},X^{\prime}}\overline{\alpha}+\alpha
^{T}R_{X,X^{\prime}}^{T}R_{X,X^{\prime}}\overline{\alpha}\\
&  =N\rho\alpha^{T}R_{X^{\prime},X^{\prime}}\overline{\alpha}+\left\vert
R_{X,X^{\prime}}\overline{\alpha}\right\vert ^{2}.
\end{align*}

But from the Introduction \ref{Sect_ap_introd} $R_{X^{\prime},X^{\prime}}$ is
positive definite and so $\alpha^{T}\Psi\overline{\alpha}>0$ iff $\alpha=0$
which implies that $\Psi$ is positive definite over $\mathbb{C}$, and since
$G$ is real we have that $\Psi$ is regular.\medskip

\textbf{Part 2.}\ From the block sizes it is clear that $\Psi$ is square with
$N^{\prime}+2M$ rows. Hence the size of $\Psi$ is independent of the number of
(scattered) data points $N$.
\end{proof}

\subsection{The scalability of the Approximate
smoother\label{SbSect_Approx_smth_scalable}}

From the basis function $G$, the $N$ independent data points $X$ and the
$N^{\prime}$ points $X^{\prime}$, we construct the $N^{\prime}\times
N^{\prime}$ matrix $R_{X^{\prime},X^{\prime}}=\left(  2\pi\right)  ^{-\frac
{d}{2}}G_{X^{\prime},X^{\prime}}$ and the $N\times N^{\prime}$ matrix
$R_{X,X^{\prime}}=\left(  2\pi\right)  ^{-\frac{d}{2}}G_{X,X^{\prime}}$. For a
specified smoothing parameter $\rho$ and independent data $y$, we construct
the matrix equation \ref{a7.094} i.e.%
\[
\left(  N\rho R_{X^{\prime},X^{\prime}}+R_{X,X^{\prime}}^{T}R_{X,X^{\prime}%
}\right)  \alpha^{\prime}=R_{X,X^{\prime}}^{T}y.
\]

The system is solved for $\alpha^{\prime}$ and the Approximate smoother is
evaluated at various points $Z$ using the formula $s_{Z}=R_{Z,X^{\prime}%
}\alpha^{\prime}$.

Our next result shows the algorithm is \textit{scalable} i.e. the time of
construction and execution of the solution is proportional to the number of
data points. This is in contrast with the Exact smoother which is not scalable
but which has quadratic dependency on the number of data points.

\begin{corollary}
\label{Cor_scalable}The Approximate smoother algorithm is scalable.
\end{corollary}

\begin{proof}
Suppose the \textit{evaluation cost} for $G\left(  x\right)  $ is $m_{G}$
multiplications and that $N\gg N^{\prime}$. The \textit{construction cost} for
$\Psi$ is
\[
\left(  N^{\prime}\right)  ^{2}m_{G}+N^{\prime}Nm_{G}+N^{\prime}N+N^{\prime
}Nm_{G}\simeq2N^{\prime}Nm_{G}.
\]

The \textit{solution cost} of an $N^{\prime}\times N^{\prime}$ matrix equation
is $\frac{1}{3}\left(  N^{\prime}\right)  ^{3}$ multiplications for a dense
matrix. Thus the \textit{total cost} is $2N^{\prime}Nm_{G}+\frac{1}{3}\left(
N^{\prime}\right)  ^{3}$ which is linearly dependent on the number of data points.

However, by the use of a basis function with support containing only several
points in $X^{\prime}$ e.g. the hat function, the construction and solution
costs can be reduced significantly. However we still have linear dependency on
$N$.
\end{proof}

\section{Convergence of the Approximate smoother to the Exact
smoother\label{Sect_ap_Ex_smth_minus_App_smth}}

In this Section we will prove several results concerning the convergence of
the Exact smoother to the Approximate smoother. We will start by proving some
general results which only assume that the weight function has property W02.
These results yield no convergence orders (hence called order-less) but show
that the Exact smoother always converges to the Approximate smoother uniformly
pointwise and normwise as the density of the set $X^{\prime}$ increases. They
will be applicable to Approximate smoothers generated by regular rectangular
grids and sparse grids.

We then derive pointwise order of convergence results which are analogous to
the results derived for the minimal norm interpolant in Section
\ref{Sect_interp_no_Taylor_converg} and subsection \ref{Sect_unisolv}, and for
the Exact smoother in Sections \ref{Sect_exsmth_err_no_Taylor} and
\ref{Sect_ex_unisolv_error}. Type 1 and Type 2 results are first derived for
the case where no explicit assumption is made regarding the unisolvency
(Definition \ref{Def_unisolv}) of the independent data. Then the concept of
unisolvent data is introduced to obtain higher orders of convergence.

\subsection{Order-less convergence
results\label{SbSect_ap_gen_converg_App_to_Ex}}

The results of this subsection only assume that the weight function has
property W02 so that the data functions are continuous. Uniform pointwise
convergence results are derived but \textbf{no order of convergence} results
are obtained. To start with we will need a definition of convergence for
\textit{sequences of independent data sets}.

\begin{definition}
\label{Def_Xk-->X}\textbf{Convergence of independent data sets}

\begin{enumerate}
\item A sequence of independent data $X_{n}=\left(  x_{n}^{\left(  i\right)
}\right)  _{i=1}^{N_{n}}$ is said to converge to the independent data
$X=\left(  x^{\left(  i\right)  }\right)  _{i=1}^{N}$, denoted $X_{n}%
\rightarrow X $, if there exists an integer $K$ such that $N_{n}=N$ when
$n\geq K$, and for each $i$, $\left\vert x_{n}^{\left(  i\right)  }-x^{\left(
i\right)  }\right\vert \rightarrow0$ as $n\rightarrow\infty$.

\item When $k\geq K$, $X_{n}$ and $X$ can be regarded as members of
$\mathbb{R}^{Nd}$ and convergence as convergence in $\mathbb{R}^{Nd}$ under
the Euclidean norm.
\end{enumerate}
\end{definition}

\begin{theorem}
\label{Thm_Jsd[sig(Zk)]->Jsd(sig(Z))}Suppose $s_{e}$ is the Exact smoother
generated by the data $\left[  X;y\right]  $. Suppose $X_{n}^{\prime}$ is a
sequence of independent data sets which converge to $X$ in the sense of
Definition \ref{Def_Xk-->X} and that $s_{a}^{\left(  n\right)  }$ is the
Approximate smoother generated by $X_{n}^{\prime}$ and $\left[  X;y\right]  $.

Then the Approximate smoothers satisfy $J_{e}\left[  s_{a}^{\left(  n\right)
}\right]  \rightarrow J_{e}\left[  s_{e}\right]  $ as $n\rightarrow\infty$,
where%
\[
J_{e}\left[  f\right]  =\rho\left\Vert f\right\Vert _{w,0}^{2}+\frac{1}%
{N}\left\vert f_{X}-y\right\vert ^{2},
\]

is the Exact smoother functional \ref{7.63}.
\end{theorem}

\begin{proof}
We first note that the definition of the convergence of independent data sets
allows us to assume that the $X_{n}^{\prime}$ have the same number of points
as $X$.

Now suppose $X^{\prime}$ is an arbitrary independent data set with the same
number of points as $X$, and $s_{a}=s_{a}\left(  X^{\prime}\right)  $ is the
corresponding Approximate smoother. If it can be shown that as a function of
$X^{\prime}$, $J_{e}\left[  s_{a}\left(  X^{\prime}\right)  \right]  $ is
continuous everywhere the the theorem holds since $J_{e}\left[  s_{e}\right]
=J_{e}\left[  s_{a}\left(  X\right)  \right]  $.

In fact, by Theorem \ref{Thm_Approx_smth_mat_eqn}
\[
s_{a}\left(  X^{\prime}\right)  \left(  x\right)  =\sum\limits_{i=1}%
^{N}R\left(  x-x_{i}^{\prime}\right)  \alpha_{i}^{\prime},
\]

where $\alpha^{\prime}=\left(  \alpha_{i}^{\prime}\right)  $ satisfies the
matrix equation
\[
\Psi\alpha^{\prime}=R_{X^{\prime},X}y,
\]

and
\[
\Psi=N\rho R_{X^{\prime},X^{\prime}}+R_{X,X^{\prime}}^{T}R_{X,X^{\prime}}.
\]

Starting with part 3 of Corollary \ref{Cor_Appr_smth_prop2} and noting that
$\left(  \alpha^{\prime}\right)  ^{T}=\left(  \Psi^{-1}R_{X^{\prime}%
,X}y\right)  ^{T}=y^{T}\overline{R}_{X,X^{\prime}}\overline{\Psi}^{-1}$ we
have%
\begin{align*}
J_{e}\left[  s_{a}\left(  X^{\prime}\right)  \right]  =\frac{1}{N}\left(
y-s_{a}\left(  X^{\prime}\right)  _{X}\right)  ^{T}\overline{y} &  =\frac
{1}{N}\left\vert y\right\vert ^{2}-\frac{1}{N}s_{a}\left(  X^{\prime}\right)
_{X}^{T}\overline{y}\\
&  =\frac{1}{N}\left\vert y\right\vert ^{2}-\frac{1}{N}\left(  R_{X,X^{\prime
}}\alpha^{\prime}\right)  ^{T}\overline{y}\\
&  =\frac{1}{N}\left\vert y\right\vert ^{2}-\frac{1}{N}\left(  \alpha^{\prime
}\right)  ^{T}\left(  R_{X,X^{\prime}}\right)  ^{T}\overline{y}\\
&  =\frac{1}{N}\left\vert y\right\vert ^{2}-\frac{1}{N}y^{T}R_{X,X^{\prime}%
}\Psi^{-1}R_{X^{\prime},X}\overline{y}.
\end{align*}

If we can show that $R_{X,X^{\prime}}$ and $\Psi_{X^{\prime}}^{-1}$ are
continuous functions of $X^{\prime}$ in a neighborhood of $X$, then we have
$J_{e}\left[  s_{a}\left(  X^{\prime}\right)  \right]  $ is a continuous
function of $X^{\prime}$. But $R_{X^{\prime},X^{\prime}}$ and $R_{X,X^{\prime
}}$ are clearly continuous for all $X^{\prime}$ so $\Psi_{X^{\prime}}$ and
$\det\Psi_{X^{\prime}}$ are continuous for all $X^{\prime}$. Further, since
$\Psi_{X^{\prime}}$ is positive definite and regular for all $X^{\prime}$,
$\det\Psi_{X^{\prime}}>0$ for all $X^{\prime}$ and it is clear from Crammer's
rule that $\Psi_{X^{\prime}}^{-1}$ is continuous everywhere. Thus
$J_{e}\left[  s_{a}\left(  X^{\prime}\right)  \right]  $ is continuous
everywhere and the proof is complete.
\end{proof}

The next corollary shows that the Approximate smoother converges to the Exact
smoother, both pointwise and norm-wise, as the grid size goes to zero.

\begin{corollary}
\label{Cor_Jsd[sig(Zk)]->Jsd(sig(Z))_1}Suppose $s_{e}$ is the Exact smoother
generated by the data $\left[  X;y\right]  $ and that we have a sequence of
regular grids $X_{n}^{\prime}$ with a common rectangular grid boundary and
grid sizes $h_{n}$. Suppose $X$ lies in the interior of the common grid
boundary. Also, let $s_{a}\left(  X_{n}^{\prime}\right)  $ denote the
Approximate smoother generated by $X_{n}^{\prime}$ and $\left[  X;y\right]  $.

Then $\left\vert h_{n}\right\vert \rightarrow0$ implies $\left\Vert
s_{a}\left(  X_{n}^{\prime}\right)  -s_{e}\right\Vert _{w,0}\rightarrow0$ and
$\left\Vert s_{a}\left(  X_{n}^{\prime}\right)  -s_{e}\right\Vert _{\infty
}\rightarrow0$, where the supremum norm $\left\Vert \cdot\right\Vert _{\infty
}$ is defined on $\mathbb{R}^{d}$.
\end{corollary}

\begin{proof}
For each data point $x^{\left(  k\right)  }\in X$ there exists a sequence of
distinct points $\left(  z_{n}^{\left(  k\right)  }\right)  _{n=1}^{\infty}$
such that $z_{n}^{\left(  k\right)  }\in X_{n}^{\prime}$ and $z_{n}^{\left(
k\right)  }\rightarrow x^{\left(  k\right)  }$ in $\mathbb{R}^{d}$ as
$n\rightarrow\infty$.

Set $Z_{n}=\left\{  z_{n}^{\left(  k\right)  }\right\}  _{k=1}^{N}$ and let
$s_{a}\left(  Z_{n}\right)  $ be the Approximate smoother generated by $Z_{n}%
$. Then $Z_{n}\rightarrow X$ as independent data and by Theorem
\ref{Thm_Jsd[sig(Zk)]->Jsd(sig(Z))}, $J_{e}\left[  s_{a}\left(  Z_{n}\right)
\right]  \rightarrow J_{e}\left[  s_{e}\right]  $.

But since $Z_{n}\subset X_{n}^{\prime}$, we have $J_{e}\left[  s_{a}\left(
X_{n}^{\prime}\right)  \right]  \leq J_{e}\left[  s_{a}\left(  Z_{n}\right)
\right]  $. Also from the definition of $s_{e}$ we have $J_{e}\left[
s_{e}\right]  \leq J_{e}\left[  s_{a}\left(  X_{n}^{\prime}\right)  \right]
$. Thus
\[
J_{e}\left[  s_{e}\right]  \leq J_{e}\left[  s_{a}\left(  X_{n}^{\prime
}\right)  \right]  \leq J_{e}\left[  s_{a}\left(  Z_{n}\right)  \right]  ,
\]

and so $J_{e}\left[  s_{a}\left(  X_{n}^{\prime}\right)  \right]  \rightarrow
J_{e}\left[  s_{e}\right]  $. To establish the convergence of $s_{a}\left(
X_{n}^{\prime}\right)  $ we use part 1 of Theorem \ref{7.58} which can be
written%
\[
J_{e}\left[  s_{e}\right]  +\rho\left\Vert s_{e}-f\right\Vert _{w,0}^{2}%
+\frac{1}{N}\sum_{k=1}^{N}\left\vert s_{e}\left(  x^{(k)}\right)  -f\left(
x^{(k)}\right)  \right\vert ^{2}=J_{e}\left[  f\right]  ,
\]

for all $f\in X_{w}^{0}$. Choosing $f=s_{a}\left(  X_{n}^{\prime}\right)  $ we
see that $\left\Vert s_{e}-s_{a}\left(  X_{n}^{\prime}\right)  \right\Vert
_{w,0}\rightarrow0$.

Finally, if $R_{x}$ is the Riesz representer of the functional $f\rightarrow
f\left(  x\right)  $
\begin{align*}
\left\vert s_{e}\left(  x\right)  -s_{a}\left(  X_{n}^{\prime}\right)  \left(
x\right)  \right\vert =\left\vert \left(  s_{e}-s_{a}\left(  X_{n}^{\prime
}\right)  ,R_{x}\right)  _{w,0}\right\vert  &  \leq\left\Vert s_{e}%
-s_{a}\left(  X_{n}^{\prime}\right)  \right\Vert _{w,0}\left\Vert
R_{x}\right\Vert _{w,0}\\
&  =\left\Vert s_{e}-s_{a}\left(  X_{n}^{\prime}\right)  \right\Vert
_{w,0}\sqrt{R_{0}\left(  0\right)  },
\end{align*}

with the right side independent of $x$.
\end{proof}

We can remove the constraint that the sets $X_{n}^{\prime}$ are regular,
regular grids:

\begin{corollary}
\label{Cor_Jsd[sig(Zk)]->Jsd(sig(Z))_2}Suppose $\Omega$ is a bounded region.
Suppose $s_{e}$ is the Exact smoother generated by the data $\left[
X;y\right]  $ and that $X\subset\Omega$. Suppose also that we have a sequence
of sets $X_{n}^{\prime}\subset\Omega$ such that $\sup\limits_{x\in\Omega
}\operatorname*{dist}\left(  x,X_{n}^{\prime}\right)  \rightarrow0$. Finally
let $s_{a}\left(  X_{n}^{\prime}\right)  $ denote the Approximate smoother
generated by $X_{n}^{\prime}$ and $\left[  X;y\right]  $.

Then $\left\Vert s_{a}\left(  X_{n}^{\prime}\right)  -s_{e}\right\Vert
_{w,0}\rightarrow0$ and $\left\Vert s_{a}\left(  X_{n}^{\prime}\right)
-s_{e}\right\Vert _{\infty}\rightarrow0$, where the supremum norm $\left\Vert
\cdot\right\Vert _{\infty}$ is defined on $\mathbb{R}^{d}$.
\end{corollary}

\begin{proof}
Let $X=\left\{  x^{\left(  k\right)  }\right\}  _{k=1}^{N}$. Since
$X\subset\Omega$ the assumption that $\sup_{x\in\Omega}\operatorname*{dist}%
\left(  x,X_{n}^{\prime}\right)  \rightarrow0$ implies

$\operatorname*{dist}\left(  x^{\left(  k\right)  },X_{n}^{\prime}\right)
\rightarrow0$ and we can now use the arguments of the proof of the previous corollary.
\end{proof}

\subsection{General error results\label{SbSect_ap_gen_App_minus_Ex}}

The results of this subsection only assume that the weight function has
property W02. The next result establishes some upper bounds for the pointwise
difference between the Exact and Approximate smoothers which are uniform on
$\mathbb{R}^{d}$. Here $\rho$ is called the smoothing coefficient and is used
to define the Exact smoother functional. No data densities are involved.

\begin{theorem}
If $s_{e}$ is the Exact smoother of the data $\left[  X;\widetilde{\mathcal{E}%
}_{X}f_{d}\right]  $ and $s_{a}$ is the Approximate smoother generated by
$\left[  X;\widetilde{\mathcal{E}}_{X}f_{d}\right]  $ and $X^{\prime}$ then%
\begin{equation}
\left\vert s_{e}\left(  x\right)  -s_{a}\left(  x\right)  \right\vert
\leq\left(  s_{e}-s_{a},s_{a}\right)  _{w,0}\sqrt{R_{0}\left(  0\right)
},\text{\quad}x\in\mathbb{R}^{d},\label{a7.20}%
\end{equation}

and%
\begin{equation}
\left\vert \mathcal{S}_{X}R_{y}\left(  x\right)  -\mathcal{S}_{X,X^{\prime}%
}^{a}R_{y}\left(  x\right)  \right\vert \leq R_{0}\left(  0\right)
\min\left\{  1,\frac{R_{0}\left(  0\right)  }{\rho}\right\}  ,\text{\quad
}x,y\in\mathbb{R}^{d},\label{a7.24}%
\end{equation}

and $\left\Vert s_{e}\right\Vert _{w,0}\leq\left\Vert f_{d}\right\Vert _{w,0}$.
\end{theorem}

\begin{proof}
From Corollary \ref{Cor_Apprx_smth_prop} we have $s_{a}=\mathcal{I}%
_{X^{\prime}}s_{e}$ so that
\begin{align*}
\left\vert s_{e}\left(  x\right)  -s_{a}\left(  x\right)  \right\vert
\leq\left\vert s_{e}\left(  x\right)  -\left(  \mathcal{I}_{X^{\prime}}%
s_{e}\right)  \left(  x\right)  \right\vert =\left\vert \left(  s_{e}%
-\mathcal{I}_{X^{\prime}}s_{e},R_{x}\right)  _{w,0}\right\vert  &
\leq\left\Vert s_{e}-\mathcal{I}_{X^{\prime}}s_{e}\right\Vert _{w,0}\left\Vert
R_{x}\right\Vert _{w,0}\\
&  =\left\Vert s_{e}-s_{a}\right\Vert _{w,0}\sqrt{R_{0}\left(  0\right)  },
\end{align*}

where the last equality follows directly from the definition of $R_{x}$.
Substituting $f_{d}=R_{y}$ into \ref{a7.20} we get
\[
\left\vert \mathcal{S}_{X}R_{y}\left(  x\right)  -\mathcal{S}_{X,X^{\prime}%
}^{a}R_{y}\left(  x\right)  \right\vert \leq\left\Vert \mathcal{S}_{X}%
R_{y}\right\Vert _{w,0}\sqrt{R_{0}\left(  0\right)  },\text{\quad}%
x,y\in\mathbb{R}^{d}.
\]

From equation \ref{7.68}, $\left\Vert \mathcal{S}_{X}f\right\Vert _{w,0}%
\leq\left\Vert f\right\Vert _{w,0}\sqrt{R_{0}\left(  0\right)  }\min\left\{
1,\frac{R_{0}\left(  0\right)  }{\rho}\right\}  $, so that
\[
\left\Vert \mathcal{S}_{X}R_{y}\right\Vert _{w,0}\leq\left\Vert R_{y}%
\right\Vert _{w,0}\left\Vert R_{0}\right\Vert _{w,0}\min\left\{  1,\frac
{R_{0}\left(  0\right)  }{\rho}\right\}  =R_{0}\left(  0\right)  \min\left\{
1,\frac{R_{0}\left(  0\right)  }{\rho}\right\}  .
\]

and hence
\[
\left\vert \mathcal{S}_{X}R_{y}\left(  x\right)  -\mathcal{S}_{X,X^{\prime}%
}^{a}R_{y}\left(  x\right)  \right\vert \leq R_{0}\left(  0\right)
\min\left\{  1,\frac{R_{0}\left(  0\right)  }{\rho}\right\}  .
\]

That $\left\Vert s_{e}\right\Vert _{w,0}\leq\left\Vert f_{d}\right\Vert
_{w,0}$ follows from part 2 of Summary \ref{Sum_Exact_smth_properties}.
\end{proof}

\subsection{Error estimates derived without explicitly assuming unisolvent
data sets}

We will now derive order estimates for the pointwise difference between the
Approximate smoother and the Exact smoother. The derivation of these estimates
will rely strongly on the fact that \textbf{the Approximate smoother is the
interpolant of the Exact smoother on the set }$X^{\prime}$ (see \ref{a7.5})
and so we will use the convergence results for the minimal norm interpolant
derived in Chapter \ref{Ch_Exact_smth}. The convergence results for the
interpolant and Exact smoother all involved one of three assumptions
concerning the unisolvency of the independent data set $X$ and the weight
function $w$ and its parameter $\kappa$:\smallskip

\begin{enumerate}
\item $w$ has property W02 for some $\kappa\geq0$ with no explicit assumption
that $X$ is unisolvent;

\item $\kappa\geq1$ with no explicit assumption that $X$ is unisolvent;

\item $\kappa\geq1$ and $X$ is unisolvent of order $\geq0$.

We will consider the first two cases in this subsection.
\end{enumerate}

\subsubsection{\protect\underline{Type 1 estimates $\left(  \kappa
\geq0\right)  $}}

The next result is an estimate of the difference between the Exact and
Approximate smoothers for an arbitrary data function in $X_{w}^{0}$.

\begin{theorem}
\label{Thm_converg_arb_func_K=0}Suppose:

\begin{enumerate}
\item The weight function $w$ has property W02 and that $G$ is the basis
function generated by $w$. Assume that for some $s$, $C_{G}$, $h_{G}>0$ the
basis function satisfies%
\begin{equation}
\left\vert G\left(  0\right)  -G\left(  x\right)  \right\vert \leq
C_{G}\left\vert x\right\vert ^{2s},\text{\quad}\left\vert x\right\vert \leq
h_{G}.\label{a7.07}%
\end{equation}

\item Let $s_{e}$ be the Exact smoother generated by the data $\left[
X;\widetilde{\mathcal{E}}_{X}f_{d}\right]  $.

\item Suppose that $s_{a}$ is the Approximate smoother generated by the data
$\left[  X;\widetilde{\mathcal{E}}_{X}f_{d}\right]  $ and the points
$X^{\prime}$ contained in $K$, where $K$ is a closed bounded infinite set.

\item Let $\mathcal{I}_{X^{\prime}}$ be the minimal norm interpolant operator
on the set $X^{\prime}$.
\end{enumerate}

Then if $k_{G}=\left(  2\pi\right)  ^{-\frac{d}{4}}\sqrt{2C_{G}}$ we have the
error bound
\begin{equation}
\left\vert s_{e}\left(  x\right)  -s_{a}\left(  x\right)  \right\vert
\leq\sqrt{\left(  s_{e}-s_{a},s_{e}\right)  _{w,0}}k_{G}\left(  h_{X^{\prime
},K}\right)  ^{s},\text{\quad}x\in K,\label{a7.025}%
\end{equation}

when $h_{X^{\prime},K}=\sup\limits_{x\in K}\operatorname*{dist}\left(
x,X^{\prime}\right)  \leq h_{G}$.

Further, $\sqrt{\left(  s_{e}-s_{a},s_{e}\right)  _{w,0}}\leq\left\Vert
s_{e}\right\Vert _{w,0}\leq\left\Vert f_{d}\right\Vert _{w,0}$ and the order
of convergence is at least $s$ in $h_{X^{\prime},K}$.
\end{theorem}

\begin{proof}
From Corollary \ref{Cor_Appr_smth_map_prop} we have $s_{a}=\mathcal{I}%
_{X^{\prime}}s_{e}$. Now we can apply Theorem \ref{Thm_|f(x)-f(y)|_inequal_2}
which estimates the error of the minimal norm interpolant of an arbitrary data
function. In this case the data function is $s_{e}$, the data points are
$X^{\prime}$, the data region is $K$, and so for $x\in K$
\begin{align*}
\left\vert s_{e}\left(  x\right)  -s_{a}\left(  x\right)  \right\vert
=\left\vert s_{e}\left(  x\right)  -\left(  \mathcal{I}_{X^{\prime}}%
s_{e}\right)  \left(  x\right)  \right\vert  &  \leq\sqrt{\left(
s_{e}-\mathcal{I}_{X^{\prime}}s_{e},s_{e}\right)  _{w,0}}k_{G}\left(
h_{X^{\prime},K}\right)  ^{s}\\
&  =\sqrt{\left(  s_{e}-s_{a},s_{e}\right)  _{w,0}}k_{G}\left(  h_{X^{\prime
},K}\right)  ^{s},
\end{align*}

when $h_{X^{\prime},K}<h_{G}$. Finally, from part 2 of Summary
\ref{Sum_Exact_smth_properties}, $I-\mathcal{I}_{X^{\prime}}$ and
$\mathcal{S}_{X}$ are contractions and so%
\[
\sqrt{\left(  s_{e}-\mathcal{I}_{X^{\prime}}s_{e},s_{e}\right)  _{w,0}}%
\leq\sqrt{\left\Vert s_{e}-\mathcal{I}_{X^{\prime}}s_{e}\right\Vert
_{w,0}\left\Vert s_{e}\right\Vert _{w,0}}\leq\left\Vert s_{e}\right\Vert
_{w,0}\leq\left\Vert f_{d}\right\Vert _{w,0}.
\]

\end{proof}

\begin{remark}
\label{Rem_Thm_converg_arb_func_K=0}\ 

\begin{enumerate}
\item Unlike the error estimates for the Approximate smoother and the minimal
norm interpolant, this estimate for $\left\vert s_{e}\left(  x\right)
-s_{a}\left(  x\right)  \right\vert $ does not require the independent data
set $X$ to be in a closed bounded infinite set. Instead, it is $X^{\prime}$
which is required to be in such a set. The density of $X$ is not explicitly
involved in the estimate i.e. there is no term involving $h_{X,K}$ and the
dependency on $X$ is in the formula for $s_{e}$ and can thus be eliminated by
the approximation $\sqrt{\left(  s_{e}-\mathcal{I}_{X^{\prime}}s_{e}%
,s_{e}\right)  _{w,0}}\leq\left\Vert f_{d}\right\Vert _{w,0}$.

\item This result confirms the convergence result of Corollary
\ref{Cor_Jsd[sig(Zk)]->Jsd(sig(Z))_1}. Unlike the Approximate smoother error
estimates, this estimate does not explicitly involve the smoothing coefficient
$\rho$ and has the same form as the interpolation error estimates. However,
the smoothing coefficient is part of the formula for $s_{e}$ and can be
eliminated by the approximation $\sqrt{\left(  s_{e}-\mathcal{I}_{X^{\prime}%
}s_{e},s_{e}\right)  _{w,0}}\leq\left\Vert f_{d}\right\Vert _{w,0}$.

\item The following argument shows the factor $\sqrt{\left(  s_{e}-s_{a}%
,s_{e}\right)  _{w,0}}$ in the estimate \ref{a7.025} can be calculated
numerically. Indeed, since $s_{a}=\mathcal{I}_{X^{\prime}}s_{e}$ we have
$\left(  s_{e}-s_{a},s_{a}\right)  _{w,0}=\left(  s_{e}-s_{a},\mathcal{I}%
_{X^{\prime}}s_{e}\right)  _{w,0}=\left(  \mathcal{I}_{X^{\prime}}%
s_{e}-\mathcal{I}_{X^{\prime}}s_{a},s_{e}\right)  _{w,0}=0$ so that $\left(
s_{e},s_{a}\right)  _{w,0}=\left\Vert s_{a}\right\Vert _{w,0}^{2}$ and

$\left(  s_{e}-s_{a},s_{e}\right)  _{w,0}=\sqrt{\left\Vert s_{e}\right\Vert
_{w,0}^{2}-\left\Vert s_{a}\right\Vert _{w,0}^{2}}$. Part 2 of Corollary
\ref{Cor_Appr_smth_prop2} and part 3 of Theorem \ref{Thm_smth_err_bound} can
then be used to calculate $\left\Vert s_{a}\right\Vert _{w,0}$ and $\left\Vert
s_{e}\right\Vert _{w,0}$ respectively.
\end{enumerate}
\end{remark}

As with the Exact smoother error estimates which were valid for arbitrary data
functions in $X_{w}^{0}$ we can try to improve the last convergence result for
data functions which have the form $R_{y}$, $y\in K$.

\begin{theorem}
\label{Thm_converg_Rx_func_K=0}Suppose we have the same assumptions and
notation as Theorem \ref{Thm_converg_arb_func_K=0} except that now the data
function is specialized to $R_{y}$. Hence $\mathcal{S}_{X}R_{y}$ is the Exact
smoother and $\mathcal{S}_{X,X^{\prime}}^{a}R_{y}$ is the Approximate smoother.

Further, suppose that the data region $K$ that contains $X$ is a bounded
closed and infinite set.

Then if $h_{X,K}\leq h_{G}$ and $h_{X^{\prime},K}\leq h_{G}$ we have the
estimate
\begin{equation}
\left\vert \mathcal{S}_{X}R_{y}\left(  x\right)  -\mathcal{S}_{X,X^{\prime}%
}^{a}R_{y}\left(  x\right)  \right\vert \leq k_{G}\left(  h_{X^{\prime}%
,K}\right)  ^{s}\left(  \sqrt{\rho N}+k_{G}\left(  h_{X,K}\right)  ^{s}%
+k_{G}\left(  h_{X^{\prime},K}\right)  ^{s}\right)  ,\text{\quad}x,y\in
K,\label{a7.17}%
\end{equation}

where $k_{G}=\left(  2\pi\right)  ^{-\frac{d}{4}}\sqrt{2C_{G}}$ and $C_{G}$
satisfies \ref{a7.07}.
\end{theorem}

\begin{proof}
Theorem \ref{Thm_converg_arb_func_K=0} with $f_{d}=R_{y}$ yields the estimate%
\begin{equation}
\left\vert \mathcal{S}_{X}R_{y}\left(  x\right)  -\mathcal{S}_{X,X^{\prime}%
}^{a}R_{y}\left(  x\right)  \right\vert \leq k_{G}\sqrt{\left(  \mathcal{S}%
_{X}R_{y}-\mathcal{I}_{X^{\prime}}\mathcal{S}_{X}R_{y},\mathcal{S}_{X}%
R_{y}\right)  _{w,0}}\left(  h_{X^{\prime},K}\right)  ^{s},\label{a7.04}%
\end{equation}

when $h_{X^{\prime},K}\leq h_{G}$ and $x\in K$. The term $\sqrt{\left(
\mathcal{S}_{X}R_{y}-\mathcal{I}_{X^{\prime}}\mathcal{S}_{X}R_{y}%
,\mathcal{S}_{X}R_{y}\right)  _{w,0}}$ will now be manipulated in order to try
to improve the convergence. Since $\mathcal{I}_{X^{\prime}}$ is a self-adjoint
projection we have
\begin{align}
\left(  \mathcal{S}_{X}R_{y}-\mathcal{I}_{X^{\prime}}\mathcal{S}_{X}%
R_{y},\mathcal{S}_{X}R_{y}\right)  _{w,0}  & =\left(  \left(  I-\mathcal{I}%
_{X^{\prime}}\right)  \mathcal{S}_{X}R_{y},\mathcal{S}_{X}R_{y}\right)
_{w,0}\label{a7.051}\\
& =\left(  \left(  I-\mathcal{I}_{X^{\prime}}\right)  \mathcal{S}_{X}%
R_{y},R_{y}-\left(  I-\mathcal{S}_{X}\right)  R_{y}\right)  _{w,0}\nonumber\\
& =\left(  \left(  I-\mathcal{I}_{X^{\prime}}\right)  \mathcal{S}_{X}%
R_{y},R_{y}\right)  _{w,0}-\left(  \left(  I-\mathcal{I}_{X^{\prime}}\right)
\mathcal{S}_{X}R_{y},\left(  I-\mathcal{S}_{X}\right)  R_{y}\right)
_{w,0}\nonumber\\
& =\left(  \left(  I-\mathcal{I}_{X^{\prime}}\right)  \mathcal{S}_{X}%
R_{y}\right)  \left(  y\right)  -\left(  \left(  I-\mathcal{S}_{X}\right)
\left(  I-\mathcal{I}_{X^{\prime}}\right)  \mathcal{S}_{X}R_{y},R_{y}\right)
_{w,0}\nonumber\\
& =\left(  \left(  I-\mathcal{I}_{X^{\prime}}\right)  \mathcal{S}_{X}%
R_{y}\right)  \left(  y\right)  -\left(  \left(  I-\mathcal{S}_{X}\right)
\left(  I-\mathcal{I}_{X^{\prime}}\right)  \mathcal{S}_{X}R_{y}\right)
\left(  y\right) \nonumber\\
& \leq\left\vert \left(  \left(  I-\mathcal{I}_{X^{\prime}}\right)
\mathcal{S}_{X}R_{y}\right)  \left(  y\right)  \right\vert +\left\vert \left(
\left(  I-\mathcal{S}_{X}\right)  \left(  I-\mathcal{I}_{X^{\prime}}\right)
\mathcal{S}_{X}R_{y}\right)  \left(  y\right)  \right\vert .\label{a7.01}%
\end{align}

The application of the Type 1 interpolation convergence estimate \ref{a927} to
the data function $\mathcal{S}_{X}R_{y}$ gives%
\begin{equation}
\left\vert \left(  \left(  I-\mathcal{I}_{X^{\prime}}\right)  \mathcal{S}%
_{X}R_{y}\right)  \left(  y\right)  \right\vert \leq\left(  \left(
I-\mathcal{I}_{X^{\prime}}\right)  \mathcal{S}_{X}R_{y},\mathcal{S}_{X}%
R_{y}\right)  _{w,0}k_{G}\left(  h_{X^{\prime},K}\right)  ^{s},\label{a7.02}%
\end{equation}

when $h_{X^{\prime},K}<h_{G}$ and $y\in K$. Further, for an arbitrary data
function $f_{d}\in X_{w}^{0}$, the Exact smoother estimate \ref{7.21} implies
\[
\left\vert \left(  \left(  I-\mathcal{S}_{X}\right)  f_{d}\right)  \left(
x\right)  \right\vert \leq\sqrt{\left(  f_{d}-\mathcal{S}_{X}f_{d}%
,f_{d}\right)  _{w,0}}\left(  \sqrt{\rho N}+k_{G}\left(  h_{X,K}\right)
^{s}\right)  ,\quad x\in K,
\]

when $h_{X,K}\leq h_{G}$. But in this case $f_{d}=\left(  I-\mathcal{I}%
_{X^{\prime}}\right)  \mathcal{S}_{X}R_{y}$ where $x,y\in K$, and so%
\[
\left\vert \left(  \left(  I-\mathcal{S}_{X}\right)  \left(  I-\mathcal{I}%
_{X^{\prime}}\right)  \mathcal{S}_{X}R_{y}\right)  \left(  x\right)
\right\vert \leq\sqrt{\left(  f_{d}-\mathcal{S}_{X}f_{d},f_{d}\right)  _{w,0}%
}\left(  \sqrt{\rho N}+k_{G}\left(  h_{X,K}\right)  ^{s}\right)  ,
\]

and
\[
\sqrt{\left(  f_{d}-\mathcal{S}_{X}f_{d},f_{d}\right)  _{w,0}}\leq
\sqrt{\left(  f_{d},f_{d}\right)  _{w,0}}=\sqrt{\left(  \left(  I-\mathcal{I}%
_{X^{\prime}}\right)  \mathcal{S}_{X}R_{y},\mathcal{S}_{X}R_{y}\right)
_{w,0}},
\]

so that%
\begin{align}
\left\vert \left(  \left(  I-\mathcal{S}_{X}\right)  \left(  I-\mathcal{I}%
_{X^{\prime}}\right)  \mathcal{S}_{X}R_{y}\right)  \left(  x\right)
\right\vert  & \leq\sqrt{\left(  \left(  I-\mathcal{I}_{X^{\prime}}\right)
\mathcal{S}_{X}R_{y},\mathcal{S}_{X}R_{y}\right)  _{w,0}}\times\nonumber\\
& \qquad\times\left(  \sqrt{\rho N}+k_{G}\left(  h_{X,K}\right)  ^{s}\right)
.\label{a7.09}%
\end{align}

Next using \ref{a7.02} and \ref{a7.09} to estimate the right side of
\ref{a7.01} we obtain%
\begin{align*}
\left(  \left(  I-\mathcal{I}_{X^{\prime}}\right)  \mathcal{S}_{X}%
R_{y},\mathcal{S}_{X}R_{y}\right)  _{w,0}  & \leq\sqrt{\left(  \left(
I-\mathcal{I}_{X^{\prime}}\right)  \mathcal{S}_{X}R_{y},\mathcal{S}_{X}%
R_{y}\right)  _{w,0}}k_{G}\left(  h_{X^{\prime},K}\right)  ^{s}+\\
& \quad+\sqrt{\left(  \left(  I-\mathcal{I}_{X^{\prime}}\right)
\mathcal{S}_{X}R_{y},\mathcal{S}_{X}R_{y}\right)  _{w,0}}\left(  \sqrt{\rho
N}+k_{G}\left(  h_{X,K}\right)  ^{s}\right)  ,
\end{align*}

so that%
\[
\sqrt{\left(  \left(  I-\mathcal{I}_{X^{\prime}}\right)  \mathcal{S}_{X}%
R_{y},\mathcal{S}_{X}R_{y}\right)  _{w,0}}\leq\sqrt{\rho N}+k_{G}\left(
h_{X,K}\right)  ^{s}+k_{G}\left(  h_{X^{\prime},K}\right)  ^{s},
\]

and as a consequence of equation \ref{a7.051}, inequality \ref{a7.04} implies
\begin{align*}
\left\vert \mathcal{S}_{X}R_{y}\left(  x\right)  -\mathcal{S}_{X,X^{\prime}%
}^{a}R_{y}\left(  x\right)  \right\vert  & \leq k_{G}\left(  h_{X^{\prime}%
,K}\right)  ^{s}\sqrt{\left(  \left(  I-\mathcal{I}_{X^{\prime}}\right)
\mathcal{S}_{X}R_{y},\mathcal{S}_{X}R_{y}\right)  _{w,0}}\\
& \leq k_{G}\left(  h_{X^{\prime},K}\right)  ^{s}\left(  \sqrt{\rho N}%
+k_{G}\left(  h_{X,K}\right)  ^{s}+k_{G}\left(  h_{X^{\prime},K}\right)
^{s}\right)  ,
\end{align*}

as claimed.
\end{proof}

\begin{remark}
\label{Rem_Thm_converg_Rx_func_K=0}Our attempt at improving convergence has
failed and the order of convergence in $h_{X^{\prime},K}$ is still $s$. Unlike
the cases of the interpolation error in Subsection
\ref{Sect_interp_no_Taylor_converg} and the Exact smoother error of Section
\ref{Sect_exsmth_err_no_Taylor} this attempt to improve the estimate for the
order of convergence of the Approximate smoother to the Exact smoother has
resulted in new terms appearing. Essentially the original term $k_{G}\left(
h_{X,K}\right)  ^{s}$ has been squared but the term $k_{G}\left(
h_{X^{\prime},K}\right)  ^{s}\left(  \sqrt{\rho N}+k_{G}\left(  h_{X,K}%
\right)  ^{s}\right)  $ has been added. Hence, to get the benefit of squaring,
the term $\sqrt{\rho N}+k_{G}\left(  h_{X,K}\right)  ^{s}$ should be the much
the same size as $k_{G}\left(  h_{X^{\prime},K}\right)  ^{s}$.
\end{remark}

The next result gives some idea of how the Exact and Approximate smoothers
compare when the smoothing parameter $\rho$ is large.

\begin{theorem}
\label{Thm_Se-Sa_large_smth_parm}Suppose the assumptions and notation of
Theorem \ref{Thm_converg_arb_func_K=0} hold so that $s_{e}$ is the Exact
smoother and $s_{a}$ is the Approximate smoother of the data function $f_{d}$.
Then for $\rho>0$%
\begin{equation}
\left\vert s_{e}\left(  x\right)  -s_{a}\left(  x\right)  \right\vert
\leq\sqrt{\left(  s_{e}-f_{d},f_{d}\right)  _{w,0}}\left(  k_{G}\right)
^{2}\frac{\left(  h_{X^{\prime},K}\right)  ^{2s}}{\sqrt{\rho}},\text{\quad
}x\in K,\text{ }f_{d}\in X_{w}^{0},\label{a7.05}%
\end{equation}

with $\sqrt{\left(  s_{e}-f_{d},f_{d}\right)  _{w,0}}\leq\left\Vert
f_{d}\right\Vert _{w,0}$. For given $\rho$ the order of convergence is $2s$.
\end{theorem}

\begin{proof}
From part 4 of Summary \ref{Sum_Exact_smth_properties}%
\[
s_{e}=\widetilde{\mathcal{E}}_{X}^{\ast}\left(  N\rho I+R_{X,X}\right)
^{-1}\widetilde{\mathcal{E}}_{X}f_{d},
\]

so that%
\begin{align*}
\left(  s_{e}-s_{a},s_{e}\right)  _{w,0}  & =\left(  s_{e}-s_{a}%
,\widetilde{\mathcal{E}}_{X}^{\ast}\left(  N\rho I+R_{X,X}\right)
^{-1}\widetilde{\mathcal{E}}_{X}f_{d}\right)  _{w,0}\\
& =\left(  \widetilde{\mathcal{E}}_{X}\left(  s_{e}-s_{a}\right)  ,\left(
N\rho I+R_{X,X}\right)  ^{-1}\widetilde{\mathcal{E}}_{X}f_{d}\right) \\
& \leq\left\vert \widetilde{\mathcal{E}}_{X}\left(  s_{e}-s_{a}\right)
\right\vert \left\vert \left(  N\rho I+R_{X,X}\right)  ^{-1}%
\widetilde{\mathcal{E}}_{X}f_{d}\right\vert
\end{align*}

But from part 2 Lemma \ref{Lem_ex_smth_err_1}%
\[
\left\Vert s_{e}-f_{d}\right\Vert _{w,0}^{2}=\left(  s_{e}-f_{d},f_{d}\right)
_{w,0}-\rho N\left\vert \left(  \rho NI+R_{X,X}\right)  ^{-1}%
\widetilde{\mathcal{E}}_{X}f_{d}\right\vert ^{2},
\]

so that $\left\vert \left(  \rho NI+R_{X,X}\right)  ^{-1}%
\widetilde{\mathcal{E}}_{X}f_{d}\right\vert \leq\left(  \rho N\right)
^{-1/2}\sqrt{\left(  s_{e}-f_{d},f_{d}\right)  _{w,0}}$.

Hence, if $\left\vert s_{e}-s_{a}\right\vert _{\infty,K}=\sup\limits_{x\in
K}\left\vert s_{e}\left(  x\right)  -s_{a}\left(  x\right)  \right\vert $,%
\begin{align}
\left(  s_{e}-s_{a},s_{e}\right)  _{w,0}  & \leq\left\vert
\widetilde{\mathcal{E}}_{X}\left(  s_{e}-s_{a}\right)  \right\vert
\sqrt{\left(  s_{e}-f_{d},f_{d}\right)  _{w,0}}\left(  \rho N\right)
^{-1/2}\nonumber\\
& \leq\sqrt{N}\left\vert s_{e}-s_{a}\right\vert _{\infty,K}\sqrt{\left(
s_{e}-f_{d},f_{d}\right)  _{w,0}}\left(  \rho N\right)  ^{-1/2}\nonumber\\
& =\left\vert s_{e}-s_{a}\right\vert _{\infty,K}\sqrt{\left(  s_{e}%
-f_{d},f_{d}\right)  _{w,0}}\rho^{-1/2},\label{a7.03}%
\end{align}

and by \ref{a7.025}
\[
\left\vert s_{e}-s_{a}\right\vert _{\infty,K}\leq\sqrt{\left(  s_{e}%
-s_{a},s_{e}\right)  _{w,0}}k_{G}\left(  h_{X^{\prime},K}\right)  ^{s},
\]

which implies%
\begin{align*}
\left\vert s_{e}-s_{a}\right\vert _{\infty,K}  & \leq\sqrt{\left\vert
s_{e}-s_{a}\right\vert _{\infty,K}\sqrt{\left(  s_{e}-f_{d},f_{d}\right)
_{w,0}}\rho^{-1/2}}\left(  k_{G}\left(  h_{X^{\prime},K}\right)  ^{s}\right)
\\
& =\sqrt{\left\vert s_{e}-s_{a}\right\vert _{\infty,K}}\sqrt[4]{\left(
s_{e}-f_{d},f_{d}\right)  _{w,0}}\rho^{-1/4}k_{G}\left(  h_{X^{\prime}%
,K}\right)  ^{s},
\end{align*}

and hence%
\[
\left\vert s_{e}-s_{a}\right\vert _{\infty,K}\leq\sqrt{\left(  s_{e}%
-f_{d},f_{d}\right)  _{w,0}}\left(  k_{G}\right)  ^{2}\frac{\left(
h_{X^{\prime},K}\right)  ^{2s}}{\sqrt{\rho}},
\]

which implies \ref{a7.05}.
\end{proof}

\begin{remark}
\label{Rem_Thm_Se-Sa_large_smth_parm}Theorem \ref{Thm_Se-Sa_large_smth_parm}
suggests that as the smoothing parameter $\rho$ increases the Approximate
smoother converges to the Exact smoother. Further, for any `large' value of
$\rho$ the rate of convergence in $h_{X^{\prime},K}$ is at least $\left(
h_{X^{\prime},K}\right)  ^{2s}$. Of course, as the smoothing parameter
increases the Exact smoother gets further from the data function.
\end{remark}

The next result summarizes the estimates of $\left\vert s_{e}\left(  x\right)
-s_{a}\left(  x\right)  \right\vert $ and $\left\vert \mathcal{S}_{X}%
R_{y}\left(  x\right)  -\mathcal{S}_{X,X^{\prime}}^{a}R_{y}\left(  x\right)
\right\vert $ we have derived above.

\begin{corollary}
\label{Cor_Thm_Se-Sa_large_smth_parm}For an arbitrary data function $f_{d}$ we
have the error estimates
\begin{equation}
\left\vert s_{e}\left(  x\right)  -s_{a}\left(  x\right)  \right\vert \leq
\min\left\{
\begin{array}
[c]{l}%
\sqrt{\left(  s_{e}-s_{a},s_{e}\right)  _{w,0}}k_{G}\left(  h_{X^{\prime}%
,K}\right)  ^{s},\\
\sqrt{\left(  s_{e}-s_{a},s_{e}\right)  _{w,0}}\sqrt{R_{0}\left(  0\right)
},\\
\sqrt{\left(  s_{e}-f_{d},f_{d}\right)  _{w,0}}\left(  k_{G}\right)  ^{2}%
\frac{\left(  h_{X^{\prime},K}\right)  ^{2s}}{\sqrt{\rho}}.
\end{array}
\right.  ,\quad x\in K,\label{a7.18}%
\end{equation}

and for the Riesz data functions $R_{y}$ we have the error estimates:%
\begin{equation}
\left\vert \mathcal{S}_{X}R_{y}\left(  x\right)  -\mathcal{S}_{X,X^{\prime}%
}^{a}R_{y}\left(  x\right)  \right\vert \leq\min\left\{
\begin{array}
[c]{l}%
k_{G}\left(  h_{X^{\prime},K}\right)  ^{s}\left(  \sqrt{\rho N}+k_{G}\left(
h_{X,K}\right)  ^{s}+k_{G}\left(  h_{X^{\prime},K}\right)  ^{s}\right)  ,\\
R_{0}\left(  0\right)  \min\left\{  1,\frac{R_{0}\left(  0\right)  }{\rho
}\right\}  ,\\
\left(  k_{G}\right)  ^{2}\frac{\left(  h_{X^{\prime},K}\right)  ^{2s}}%
{\sqrt{\rho}}\left(  \sqrt{\rho N}+k_{G}\left(  h_{X,K}\right)  ^{s}\right)
,\\
\left(  k_{G}\right)  ^{2}\frac{\left(  h_{X^{\prime},K}\right)  ^{2s}}%
{\sqrt{\rho}}\sqrt{R_{0}\left(  0\right)  },\\
\left(  k_{G}\right)  ^{2}\frac{\left(  h_{X^{\prime},K}\right)  ^{2s}}%
{\sqrt{\rho}}\sqrt{\left\Vert R_{y}\right\Vert _{\infty,K}+R_{0}\left(
0\right)  \min\left\{  1,\frac{R_{0}\left(  0\right)  }{\rho}\right\}  }.
\end{array}
\right. \label{a7.191}%
\end{equation}

for $x,y\in K.$
\end{corollary}

\begin{proof}
The inequalities \ref{a7.18} are the estimates \ref{a7.20}, \ref{a7.025} and
\ref{a7.05}.\medskip

Regarding the inequalities \ref{a7.191}: The first inequality is \ref{a7.17}
and the second inequality is \ref{a7.24}. When $f_{d}=R_{y}$ is substituted
into \ref{a7.05} we get for $x\in K$
\[
\left\vert \mathcal{S}_{X}R_{y}\left(  x\right)  -\mathcal{S}_{X,X^{\prime}%
}^{a}R_{y}\left(  x\right)  \right\vert \leq\left(  k_{G}\right)  ^{2}%
\frac{\left(  h_{X^{\prime},K}\right)  ^{2s}}{\sqrt{\rho}}\sqrt{\mathcal{S}%
_{X}R_{y}\left(  x\right)  -R_{y}\left(  x\right)  },
\]

and then substituting the estimates \ref{a7.22} for $\mathcal{S}_{X}%
R_{y}\left(  x\right)  -R_{y}\left(  x\right)  $ immediately gives the
remaining three inequalities.
\end{proof}

\subsubsection{\protect\underline{Type 1 examples}%
\label{SbSbSect_ap_typ1_Ap_to_Ex_example}}

The weight function examples used here are those used in the interpolation and
Exact smoother documents i.e. the radial shifted thin-plate splines, Gaussian
and Sobolev splines, and the tensor product extended B-spline weight
functions, augmented by the tensor product central difference weight functions
from Chapter \ref{Ch_cent_diff_wt_fn_ten_prod}.

When the smoothing coefficient $\rho$ is zero the error estimates of this
subsection are the Type 1 interpolation estimates of Subsubsection
\ref{Sect_interp_no_Taylor_converg} and consequently the formulas or values
for $s,C_{G},h_{G},k_{G}$ have already been calculated and were given in Table
\ref{Tbl_NonUnisolvTyp1Converg}. We give these values here in Table
\ref{Tbl_NonUnisolvTyp1ConvAugm} together with the Type 1 values from Table
\ref{Tbl_ConvergCentral} for the central difference tensor product weight functions.%

\begin{table}[htbp] \centering
$%
\begin{tabular}
[c]{|c|c||c|c|c|}\hline
\multicolumn{5}{|c|}{\textbf{Type 1} smoother convergence order estimates,
$\kappa\geq0$.}\\\hline\hline
\multicolumn{5}{|c|}{Smoothness assumption on $G$; $k_{G}=\left(  2\pi\right)
^{-d/4}\sqrt{2C_{G}}$.}\\\hline
& Parameter & Converg. &  & \\
Weight function & constraints & order $s$ & $C_{G}$ & $h_{G}$\\\hline\hline
\multicolumn{1}{|l|}{Sobolev splines} & $v-\frac{d}{2}=1$ & $\frac{1}{2}$ &
\multicolumn{1}{|l|}{$\frac{\left\Vert \rho K_{0}\left(  \rho\right)
\right\Vert _{\infty}}{2^{v-1}\Gamma\left(  v\right)  }^{\left(  2\right)  }$}
& $\infty$\\\cline{2-4}%
\multicolumn{1}{|l|}{\quad($v>d/2$)} & $v-\frac{d}{2}\neq1$ & $1$ &
\multicolumn{1}{|l|}{$\frac{\left\Vert D^{2}\widetilde{K}_{v-d/2}\right\Vert
_{\infty}}{2^{v}\Gamma\left(  v\right)  }^{\left(  2\right)  }$} & "\\\hline
\multicolumn{1}{|l|}{Shifted thin-plate} & - & $1$ & \multicolumn{1}{|l|}{eq.
(\ref{a7.641})} & $\infty$\\
\multicolumn{1}{|l|}{($-d/2<v<0$)} &  &  &  & \\\hline
\multicolumn{1}{|l|}{Gaussian} & - & $1$ & \multicolumn{1}{|l|}{$2e^{-3/2}$} &
$\infty$\\\hline
\multicolumn{1}{|l|}{Extended B-spline} & - & $\frac{1}{2}$ &
\multicolumn{1}{|l|}{$G_{1}\left(  0\right)  ^{d-1}\left\Vert DG_{1}%
\right\Vert _{\infty}\sqrt{d}$ $^{\left(  1\right)  }$} & $\infty$\\\hline
\multicolumn{1}{|l|}{Central difference} & - & $\frac{1}{2}$ & $G_{1}\left(
0\right)  ^{d-1}\left\Vert DG_{1}\right\Vert _{\infty}\sqrt{d}$ $^{\left(
1\right)  }$ & $\infty$\\\hline\hline
\multicolumn{5}{|l|}{$^{\left(  1\right)  }${\small \ }$G_{1}${\small \ is the
univariate basis function used to form the tensor product}}\\\hline
\multicolumn{5}{|l|}{$^{\left(  2\right)  }${\small \ }$K_{v}${\small \ is the
modified Bessel function and} $\widetilde{K}_{v}\left(  r\right)  =r^{v}%
K_{v}\left(  r\right)  $.}\\\hline
\end{tabular}
$\caption{}\label{Tbl_NonUnisolvTyp1ConvAugm}%
\end{table}%
%

\begin{equation}
C_{G}=\left\vert \left(  2rf^{\prime\prime}+f^{\prime}\right)  \left(
r_{\max}\right)  \right\vert ,\quad where\text{ }f\left(  r\right)  =\left(
1+r\right)  ^{v}\text{ }and\text{ }r_{\max}=\frac{1-2v}{3}.\label{a7.641}%
\end{equation}

\subsubsection{\protect\underline{Type 2 estimates $\left(  \kappa
\geq1\right)  $}\label{SbSbSect_ap_Typ2_ApprMinusExact}}

These results can be proved by first using the theorems of Subsection
\ref{SbSect_ex_smth_typ2_error} to obtain estimates of the form $\left\vert
G\left(  0\right)  -G\left(  x\right)  \right\vert \leq C_{G}\left\vert
x\right\vert ^{2s}$ and then employing Corollary
\ref{Cor_Thm_Se-Sa_large_smth_parm}.

\begin{theorem}
\label{Thm_ap_G(0)minusG(x)_bound}Suppose a weight function satisfies property
W02 for some $\kappa\geq1$ and denote the basis function by $G$. Then the
smoother error estimates of Corollary \ref{Cor_Thm_Se-Sa_large_smth_parm} hold
for
\begin{equation}
C_{G}=-\frac{1}{2}\left(  \left\vert D\right\vert ^{2}G\right)  \left(
0\right)  d,\text{\quad}s=1,\text{\quad}h_{G}=\infty,\text{\quad}k_{G}=\left(
2\pi\right)  ^{-\frac{d}{4}}\sqrt{-\left(  \left\vert D\right\vert
^{2}G\right)  \left(  0\right)  }\sqrt{d}.\label{a7.301}%
\end{equation}

\end{theorem}

However, if the weight function is radial we can use the estimates of:

\begin{theorem}
\label{Thm_ap_G(0)minusG(x)_bound_w_radial}Suppose a radial weight function
satisfies property W02 for $\kappa=1$ and denote the (radial) basis function
by $G$. Set $r=\left\vert x\right\vert $. Then:

\begin{enumerate}
\item If $G\left(  x\right)  =f\left(  r^{2}\right)  $ the smoother error
estimates of Corollary \ref{Cor_Thm_Se-Sa_large_smth_parm} hold for%
\[
C_{G}=-f^{\prime}\left(  0\right)  d^{2},\text{\quad}s=1,\text{\quad}%
h_{G}=\infty,\text{\quad}k_{G}=\left(  2\pi\right)  ^{-\frac{d}{4}}%
\sqrt{-2f^{\prime}\left(  0\right)  }d.
\]

\item If $G\left(  x\right)  =g\left(  r\right)  $ then the smoother error
estimates of Corollary \ref{Cor_Thm_Se-Sa_large_smth_parm} hold for
\begin{equation}
C_{G}=-\frac{1}{2}g^{\prime\prime}\left(  0\right)  d^{2},\text{\quad
}s=1,\text{\quad}h_{G}=\infty,\text{\quad}k_{G}=\left(  2\pi\right)
^{-\frac{d}{4}}\sqrt{-g^{\prime\prime}\left(  0\right)  }d.\label{a7.35}%
\end{equation}

\end{enumerate}
\end{theorem}

If the weight function is a tensor product the following result will be useful:

\begin{theorem}
\label{Thm_Final_Type2_estim_2}Suppose a tensor product weight function
satisfies property W02 for $\kappa=1$ and denote the univariate basis function
by $G_{1}$. Then the smoother error estimates of Corollary
\ref{Cor_Thm_Se-Sa_large_smth_parm} hold when
\[
C_{G}=-\frac{d}{2}G_{1}\left(  0\right)  ^{d-1}D^{2}G_{1}\left(  0\right)
,\text{\quad}s=1,\text{\quad}h_{G}=\infty,\text{\quad}k_{G}=\left(
2\pi\right)  ^{-\frac{d}{4}}\sqrt{-G_{1}\left(  0\right)  ^{d-1}D^{2}%
G_{1}\left(  0\right)  }\sqrt{d}.
\]

\end{theorem}

\subsubsection{\protect\underline{Type 2 examples}}

When the smoothing coefficient is zero the smoother error estimates of the
Theorems of the previous subsubsection become algebraically identical to the
Type 2 interpolant error estimates of Chapters \ref{Ch_Interpol} and
\ref{Ch_cent_diff_wt_fn_ten_prod}. Further, the weight function examples used
above were also used for the Type 2 interpolation examples of Chapter
\ref{Ch_Interpol}. If we use the `radial' Theorem
\ref{Thm_ap_G(0)minusG(x)_bound_w_radial} to do the estimates for the radial
basis functions and Theorem \ref{Thm_ap_G(0)minusG(x)_bound} to do the
estimates for the tensor product basis functions then the values for $s$,
$C_{G}$, $h_{G}$, $k_{G}$ will match those obtained for the interpolants.
These are given below in Table \ref{Tbl_NonUnisolvTyp2ConvAugm} which is Table
\ref{Tbl_InterpNonUnisolvTyp2Conv} augmented by the results for the central
difference tensor product weight functions from Table \ref{Tbl_ConvergCentral}.%

\begin{table}[htbp] \centering
$%
\begin{tabular}
[c]{|c|c||c|c|}\hline
\multicolumn{4}{|c|}{Type 2 smoother convergence order estimates.}%
\\\hline\hline
\multicolumn{4}{|c|}{Assume $\kappa\geq1$.}\\\hline
& Parameter & Converg. & \\
Weight function & constraints & order & $\left(  2\pi\right)  ^{d/4}%
k_{G}/\sqrt{d}$\\\hline\hline
\multicolumn{1}{|l|}{Sobolev splines} & $v-\frac{d}{2}\geq2$ & $1$ &
$\sqrt{\frac{\Gamma\left(  v-d/2-1\right)  }{2^{d/2+1}\Gamma\left(  v\right)
}}$\\\cline{2-4}\cline{2-4}%
($v>d/2$) & $1<v-\frac{d}{2}<2$ & $1$ & $\sqrt{\frac{\Gamma\left(
v-d/2-1\right)  }{2^{2v-d/2-3}\Gamma\left(  v\right)  }}$\\\hline
\multicolumn{1}{|l|}{Shifted thin-plate} & - & $1$ & $\sqrt{-2v}$\\
($-d/2<v<0$) &  &  & \\\hline
\multicolumn{1}{|l|}{Gaussian} & - & $1$ & $\sqrt{2}$\\\hline
\multicolumn{1}{|l|}{Extended B-spline} & $n\geq2$ & $1$ & $\sqrt
{-G_{1}\left(  0\right)  ^{d-1}D^{2}G_{1}\left(  0\right)  }^{\text{ }\left(
1\right)  }$\\
($1\leq n\leq l$) &  &  & \\\hline
\multicolumn{1}{|l|}{Central difference} & $n\geq2$ & $1$ & $\sqrt
{-G_{1}\left(  0\right)  ^{d-1}D^{2}G_{1}\left(  0\right)  }^{\text{ }\left(
1\right)  }$\\
($1\leq n\leq l$) &  &  & \\\hline\hline
\multicolumn{4}{|l|}{$^{\left(  1\right)  }${\small \ }$G_{1}${\small \ is the
univariate basis function used to form the tensor product.}}\\\hline
\end{tabular}
$\caption{}\label{Tbl_NonUnisolvTyp2ConvAugm}%
\end{table}%

When $\rho=0$ we see that the order of convergence is at least $1$ for an
arbitrary data function and at least $2$ for a Riesz representer data
function, no matter what value $\kappa$ takes. However, in the next section we
will show that by explicitly assuming the independent data is unisolvent of
order $\kappa\geq1$ it follows that an order of convergence of at least
$\left\lfloor \kappa\right\rfloor $ can be attained for an arbitrary data
function and an order of convergence of $\left\lfloor 2\kappa\right\rfloor $
for a Riesz representer data function.

\subsection{Error estimates explicitly using unisolvent data sets}

Unisolvent sets were introduced in Subsection \ref{Sect_unisolv}%
\textit{\ }where they were used to define the Lagrange interpolation operators
$\mathcal{P}$ and $\mathcal{Q}=I-\mathcal{P}$. Lemma
\ref{Lem_Lagrange_interpol} supplied the key interpolation theory and I
reproduce it here:

\begin{lemma}
\label{Lem_Lagrange_interpol_3}(Copy of Lemma \ref{Lem_Lagrange_interpol}%
)Suppose first that:

\begin{enumerate}
\item $\Omega$ is a bounded region of $\mathbb{R}^{d}$ having the cone property.

\item $X$ is a unisolvent subset of $\Omega$ of order $\kappa$.

Now define
\[
h_{X,\Omega}=\sup\limits_{\omega\in\Omega}\operatorname*{dist}\left(
\omega,X\right)  ,
\]

\end{enumerate}

and fix $x\in X$. By using Lagrange polynomial interpolation techniques it can
be shown there are constants $c_{\Omega,\kappa},h_{\Omega,\kappa}>0$ such that
when $h_{X,\Omega}<h_{\Omega,\kappa}$ there exists a minimal unisolvent set
$A\subset X$ satisfying
\[
\operatorname*{diam}\left(  A\cup\left\{  x\right\}  \right)  \leq
c_{\Omega,\kappa}h_{X,\Omega}.
\]

Further, suppose $\left\{  l_{j}\right\}  _{j=1}^{M}$ is the cardinal basis of
$P_{\kappa}$ with respect to a minimal unisolvent subset of $\Omega$. Again,
using Lagrange interpolation techniques, it can be shown there exists a
constant $K_{\Omega,m}^{\prime}>0$ such that
\[
\sum\limits_{j=1}^{M}\left\vert l_{j}\left(  x\right)  \right\vert \leq
K_{\Omega,m}^{\prime},
\]

for all $x\in\Omega$ and all minimal unisolvent subsets of $\Omega$.
\end{lemma}

To derive our estimates for $\left\vert s_{e}\left(  x\right)  -s_{a}\left(
x\right)  \right\vert $ we will rely on the formula $s_{a}=\mathcal{I}%
_{X^{\prime}}s_{e}$ of part 2 Corollary \ref{Cor_Apprx_smth_prop}. Hence we
will also need the minimal norm interpolation convergence results of
Subsection \ref{Sect_unisolv}. We are now ready to state our order of
convergence result for an arbitrary data function in $X_{w}^{0}$:

\begin{theorem}
\label{Thm_ConvergSaToSe_not_unisolv_II}Suppose:

\begin{enumerate}
\item $w$ is a weight function with property W02 for parameter $\kappa\geq1 $
and let $G$ be the basis function generated by $w$. Set $m=\left\lfloor
\kappa\right\rfloor $.

\item $s_{e}$ is the Exact smoother generated by the data $\left[
X;\widetilde{\mathcal{E}}_{X}f_{d}\right]  $ for some data function $f_{d}\in
X_{w}^{0}$.

\item $s_{a}$ is the Approximate smoother generated by $\left[
X;\widetilde{\mathcal{E}}_{X}f_{d}\right]  $ and the set $X^{\prime}%
\subset\Omega$.

\item $\mathcal{I}_{X^{\prime}}$ is the minimal norm interpolant operator on
the set $X^{\prime}$.

\item We use the notation and assumptions of Lemma
\ref{Lem_Lagrange_interpol_3} which means that $X^{\prime}$ is $\kappa
$-unisolvent and $\Omega$ is a bounded region whose boundary satisfies the
cone condition.
\end{enumerate}

Then there exist positive constants $k_{G}=\frac{d^{m/2}}{\left(  2\pi\right)
^{d/2}}\left(  c_{\Omega,\kappa}\right)  ^{m}K_{\Omega,m}^{\prime}%
\max\limits_{\left\vert \beta\right\vert =m}\left\vert D^{2\beta}G\left(
0\right)  \right\vert $ and $h_{\Omega,m}$ such that%
\begin{equation}
\left\vert s_{e}\left(  x\right)  -s_{a}\left(  x\right)  \right\vert
\leq\sqrt{\left(  s_{e}-s_{a},s_{e}\right)  _{w,0}}\text{ }k_{G}\left(
h_{X^{\prime},\Omega}\right)  ^{m},\quad x\in\overline{\Omega},\label{a7.12}%
\end{equation}

when $h_{X^{\prime},\Omega}=\sup\limits_{\omega^{\prime}\in\Omega
}\operatorname*{dist}\left(  \omega^{\prime},X^{\prime}\right)  <h_{G}$.

The constants $c_{\Omega,m}$, $K_{\Omega,m}^{\prime}$ and $h_{\Omega,m}$ only
depend on $\Omega,m$ and $d$. Further,

$\sqrt{\left(  s_{e}-s_{a},s_{e}\right)  _{w,0}}\leq\left\Vert s_{e}%
\right\Vert _{w,0}\leq\left\Vert f_{d}\right\Vert _{w,0}$ and the order of
convergence is at least $\left\lfloor \kappa\right\rfloor $ in $h_{X^{\prime
},\Omega}$.
\end{theorem}

\begin{proof}
From part 2 Corollary \ref{Cor_Apprx_smth_prop}, $s_{a}=\mathcal{I}%
_{X^{\prime}}s_{e}$ and by the interpolation estimate \ref{1.32} there exists
a constant $k_{G}$ such that%
\begin{equation}
\left\vert s_{e}\left(  x\right)  -\left(  \mathcal{I}_{X^{\prime}}%
s_{e}\right)  \left(  x\right)  \right\vert \leq\sqrt{\left(  s_{e}%
-\mathcal{I}_{X^{\prime}}s_{e},s_{e}\right)  _{w,0}}k_{G}\left(  h_{X^{\prime
},\Omega}\right)  ^{m},\quad x\in\overline{\Omega},\label{a7.32}%
\end{equation}

when $h_{X^{\prime},\Omega}<h_{\Omega,m}$. From Theorem
\ref{Thm_converg_arb_func_K=0}, $\sqrt{\left(  s_{e}-\mathcal{I}_{X^{\prime}%
}s_{e},s_{e}\right)  _{w,0}}\leq\left\Vert s_{e}\right\Vert _{w,0}%
\leq\left\Vert f_{d}\right\Vert _{w,0}$ and so the order of convergence is at
least $m$.
\end{proof}

\begin{remark}
\label{Rem_Thm_ConvergSaToSe_not_unisolv_II}The estimate derived in the last
theorem has its dependency on the data points $X$ contained in $s_{e}$ in the
expression $\sqrt{\left(  s_{e}-s_{a},s_{e}\right)  _{w,0}}$.
\end{remark}

We now try to derive an improved convergence result from the last theorem for
a data function which has the form $R_{y}$.

\begin{theorem}
\label{Thm_converg_Rx_func_Kgt0}Suppose we have the same assumptions and
notation as Theorem \ref{Thm_ConvergSaToSe_not_unisolv_II} except that now the
data function is specialized to $R_{y}$, so that $\mathcal{S}_{X}R_{y}$ is the
Exact smoother of $R_{y}$ and $\mathcal{S}_{X,X^{\prime}}^{a}R_{y}$ is the
Approximate smoother of $R_{y}$.

Further, suppose $X$ is a unisolvent subset of $\Omega$ of order $\kappa$, and
suppose that the data region $\Omega$ that contains $X$ is a bounded region
which satisfies the cone condition. Set $m=\left\lfloor \kappa\right\rfloor $.

Then if $h_{X,K}\leq h_{\Omega,m}$ and $h_{X^{\prime},K}\leq h_{\Omega,m}$ we
have the estimate
\begin{equation}
\left\vert \mathcal{S}_{X}R_{y}\left(  x\right)  -\mathcal{S}_{X,X^{\prime}%
}^{a}R_{y}\left(  x\right)  \right\vert \leq k_{G}\left(  h_{X^{\prime}%
,\Omega}\right)  ^{m}\left(  K_{\Omega,m}^{\prime}\sqrt{\rho N}+k_{G}\left(
h_{X,\Omega}\right)  ^{m}+k_{G}\left(  h_{X^{\prime},\Omega}\right)
^{m}\right)  ,\label{a7.06}%
\end{equation}

for all $x,y\in\overline{\Omega}$, where $k_{G}=\frac{d^{m/2}}{\left(
2\pi\right)  ^{d/2}}\left(  c_{\Omega,\kappa}\right)  ^{m}K_{\Omega,m}%
^{\prime}\max\limits_{\left\vert \beta\right\vert =m}\left\vert D^{2\beta
}G\left(  0\right)  \right\vert $.

The order of convergence is at least $m$ in $h_{X^{\prime},\Omega}$.
\end{theorem}

\begin{proof}
Theorem \ref{Thm_ConvergSaToSe_not_unisolv_II} with data function $R_{y}$
implies
\begin{equation}
\left\vert \mathcal{S}_{X}R_{y}\left(  x\right)  -\mathcal{S}_{X,X^{\prime}%
}^{a}R_{y}\left(  x\right)  \right\vert \leq k_{G}\sqrt{\left(  \mathcal{S}%
_{X}R_{y}-\mathcal{S}_{X,X^{\prime}}^{a}R_{y},\mathcal{S}_{X}R_{y}\right)
_{w,0}}\left(  h_{X^{\prime},\Omega}\right)  ^{m},\label{a7.54}%
\end{equation}

for $x\in\overline{\Omega}$, when $h_{X^{\prime},\Omega}\leq h_{\Omega,m}$.

The term $\sqrt{\left(  \mathcal{S}_{X}R_{y}-\mathcal{S}_{X,X^{\prime}}%
^{a}R_{y},\mathcal{S}_{X}R_{y}\right)  _{w,0}}$ will now be manipulated to try
to improve the convergence. Indeed, from inequality \ref{a7.01}
\begin{align}
\left(  \mathcal{S}_{X}R_{y}-\mathcal{S}_{X,X^{\prime}}^{a}R_{y}%
,\mathcal{S}_{X}R_{y}\right)  _{w,0}  & =\left(  \left(  I-\mathcal{I}%
_{X^{\prime}}\right)  \mathcal{S}_{X}R_{y},\mathcal{S}_{X}R_{y}\right)
_{w,0}\label{a7.55}\\
& \leq\left\vert \left(  \left(  I-\mathcal{I}_{X^{\prime}}\right)
\mathcal{S}_{X}R_{y}\right)  \left(  y\right)  \right\vert +\nonumber\\
& \qquad+\left\vert \left(  \left(  I-\mathcal{S}_{X}\right)  \left(
I-\mathcal{I}_{X^{\prime}}\right)  \mathcal{S}_{X}R_{y}\right)  \left(
y\right)  \right\vert ,\label{a7.51}%
\end{align}

and using \ref{1.32} we estimate the interpolation error to be%
\begin{equation}
\left\vert \left(  \left(  I-\mathcal{I}_{X^{\prime}}\right)  \mathcal{S}%
_{X}R_{y}\right)  \left(  y\right)  \right\vert \leq\sqrt{\left(  \left(
I-\mathcal{I}_{X^{\prime}}\right)  \mathcal{S}_{X}R_{y},\mathcal{S}_{X}%
R_{y}\right)  _{w,0}}\text{ }k_{G}\left(  h_{X^{\prime},\Omega}\right)
^{m},\label{a7.52}%
\end{equation}

when $y\in\Omega$ and $h_{X^{\prime},\Omega}\leq h_{\Omega,m}$, and using
\ref{7.20} we estimate the Exact smoothing error for an arbitrary data
function $f_{d}\in X_{w}^{0}$ to be
\[
\left\vert f_{d}\left(  y\right)  -\left(  \mathcal{S}_{X}f_{d}\right)
\left(  y\right)  \right\vert \leq\sqrt{\left(  f_{d}-\mathcal{S}_{X}%
f_{d},f_{d}\right)  _{w,0}}\left(  K_{\Omega,m}^{\prime}\sqrt{\rho N}%
+k_{G}\left(  h_{X,\Omega}\right)  ^{m}\right)  ,\quad y\in\overline{\Omega}.
\]

Consequently, when $f_{d}=\left(  I-\mathcal{I}_{X^{\prime}}\right)
\mathcal{S}_{X}R_{y}$%
\[
\left\vert \left(  \left(  I-\mathcal{S}_{X}\right)  \left(  I-\mathcal{I}%
_{X^{\prime}}\right)  \mathcal{S}_{X}R_{y}\right)  \left(  y\right)
\right\vert \leq\sqrt{\left(  f_{d}-\mathcal{S}_{X}f_{d},f_{d}\right)  _{w,0}%
}\left(  K_{\Omega,m}^{\prime}\sqrt{\rho N}+k_{G}\left(  h_{X,\Omega}\right)
^{m}\right)  ,
\]

and
\begin{align*}
\sqrt{\left(  f_{d}-\mathcal{S}_{X}f_{d},f_{d}\right)  _{w,0}}\leq
\sqrt{\left(  f_{d},f_{d}\right)  _{w,0}}=\left\Vert f_{d}\right\Vert _{w,0}
&  =\left\Vert \left(  I-\mathcal{I}_{X^{\prime}}\right)  \mathcal{S}_{X}%
R_{y}\right\Vert _{w,0}\\
&  =\sqrt{\left(  \left(  I-\mathcal{I}_{X^{\prime}}\right)  \mathcal{S}%
_{X}R_{y},\mathcal{S}_{X}R_{y}\right)  _{w,0}},
\end{align*}

so that%
\begin{equation}
\left\vert \left(  \left(  I-\mathcal{S}_{X}\right)  \left(  I-\mathcal{I}%
_{X^{\prime}}\right)  \mathcal{S}_{X}R_{y}\right)  \left(  y\right)
\right\vert \leq\sqrt{\left(  \left(  I-\mathcal{I}_{X^{\prime}}\right)
\mathcal{S}_{X}R_{y},\mathcal{S}_{X}R_{y}\right)  _{w,0}}\left(  K_{\Omega
,m}^{\prime}\sqrt{\rho N}+k_{G}\left(  h_{X,\Omega}\right)  ^{m}\right)
.\label{a7.53}%
\end{equation}

Next, using \ref{a7.52} and \ref{a7.53} to estimate the right side of
\ref{a7.51} we obtain%
\begin{align*}
\left(  \left(  I-\mathcal{I}_{X^{\prime}}\right)  \mathcal{S}_{X}%
R_{y},\mathcal{S}_{X}R_{y}\right)  _{w,0}  & \leq\sqrt{\left(  \left(
I-\mathcal{I}_{X^{\prime}}\right)  \mathcal{S}_{X}R_{y},\mathcal{S}_{X}%
R_{y}\right)  _{w,0}}\left(  K_{\Omega,m}^{\prime}\sqrt{\rho N}+k_{G}\left(
h_{X,\Omega}\right)  ^{m}\right)  +\\
& \qquad+\sqrt{\left(  \left(  I-\mathcal{I}_{X^{\prime}}\right)
\mathcal{S}_{X}R_{y},\mathcal{S}_{X}R_{y}\right)  _{w,0}}k_{G}\left(
h_{X^{\prime},\Omega}\right)  ^{m},
\end{align*}

for $y\in\overline{\Omega}$ so that%
\[
\sqrt{\left(  \left(  I-\mathcal{I}_{X^{\prime}}\right)  \mathcal{S}_{X}%
R_{y},\mathcal{S}_{X}R_{y}\right)  _{w,0}}\leq K_{\Omega,m}^{\prime}\sqrt{\rho
N}+k_{G}\left(  h_{X,\Omega}\right)  ^{m}+k_{G}\left(  h_{X^{\prime},\Omega
}\right)  ^{m},
\]

and as a consequence of equation \ref{a7.55}, inequality \ref{a7.54} implies
\begin{align*}
\left\vert s_{e}\left(  y\right)  -s_{a}\left(  y\right)  \right\vert  & \leq
k_{G}\left(  h_{X^{\prime},\Omega}\right)  ^{m}\sqrt{\left(  s_{e}-s_{a}%
,s_{e}\right)  _{w,0}}\\
& =k_{G}\left(  h_{X^{\prime},\Omega}\right)  ^{m}\sqrt{\left(  \left(
I-\mathcal{I}_{X^{\prime}}\right)  \mathcal{S}_{X}R_{y},\mathcal{S}_{X}%
R_{y}\right)  _{w,0}}\\
& \leq k_{G}\left(  h_{X^{\prime},\Omega}\right)  ^{m}\left(  K_{\Omega
,m}^{\prime}\sqrt{\rho N}+k_{G}\left(  h_{X,\Omega}\right)  ^{m}+k_{G}\left(
h_{X^{\prime},\Omega}\right)  ^{m}\right)  ,
\end{align*}

as claimed.
\end{proof}

\begin{remark}
\label{Rem_Thm_converg_Rx_func_K>0}This is the same order of convergence in
$h_{X^{\prime},\Omega}$ that was obtained in Theorem
\ref{Thm_ConvergSaToSe_not_unisolv_II} for an arbitrary data function, namely
$m=\left\lfloor \kappa\right\rfloor $.
\end{remark}

\section{Type 1 Exact smoother error estimates $\left(  \kappa\geq0\right)  $}

In order to estimate the pointwise error of the Approximate smoother we will
need the following Type 1 error estimates of the Exact smoother.

\begin{theorem}
\label{Thm_ap_typ1_Exact_smth_err_summary}(Copy of Theorem
\ref{Thm_ex_typ1_err_summary}) Suppose $K$ is a bounded closed infinite set
and there exist constants $C_{G}$, $s$, $h_{G}>0$ such that
\[
G\left(  0\right)  -G\left(  x\right)  \leq C_{G}\left\vert x\right\vert
^{2s},\text{\quad}\left\vert x\right\vert \leq h_{G},\,x\in K.
\]

Set $k_{G}=\left(  2\pi\right)  ^{-\frac{1}{4}}\sqrt{2C_{G}}$. Then the Exact
smoother $\mathcal{S}_{X}f_{d}$ satisfies the error estimates%
\begin{equation}
\left\vert f_{d}\left(  x\right)  -\left(  \mathcal{S}_{X}f_{d}\right)
\left(  x\right)  \right\vert \leq\min\left\{
\begin{array}
[c]{l}%
\sqrt{\left(  f_{d}-s_{e},f_{d}\right)  _{w,0}}\left(  \sqrt{\rho N}%
+k_{G}\left(  h_{X,K}\right)  ^{s}\right)  ,\\
\sqrt{\left(  f_{d}-s_{e},f_{d}\right)  _{w,0}}\sqrt{R_{0}\left(  0\right)
},\\
\left\Vert f_{d}\right\Vert _{\infty,K}+\left\Vert f_{d}\right\Vert
_{w,0}\sqrt{R_{0}\left(  0\right)  }\min\left\{  1,\frac{R_{0}\left(
0\right)  }{\rho}\right\}  ,
\end{array}
\right. \label{a7.11}%
\end{equation}

and $\mathcal{S}_{X}R_{y}$ satisfies the double order convergence estimate%
\begin{equation}
\left\vert R_{y}\left(  x\right)  -\left(  \mathcal{S}_{X}R_{y}\right)
\left(  x\right)  \right\vert \leq\min\left\{
\begin{array}
[c]{l}%
\left(  \sqrt{\rho N}+k_{G}\left(  h_{X,K}\right)  ^{s}\right)  ^{2},\\
R_{0}\left(  0\right)  ,\\
\left\Vert R_{y}\right\Vert _{\infty,K}+R_{0}\left(  0\right)  \min\left\{
1,\frac{R_{0}\left(  0\right)  }{\rho}\right\}  ,
\end{array}
\right. \label{a7.22}%
\end{equation}

for $x,y\in K$, when $h_{X,K}=\sup\limits_{s\in K}\operatorname*{dist}\left(
s,X\right)  \leq h_{G}$.

Here $R_{0}\left(  0\right)  =\left(  2\pi\right)  ^{-\frac{1}{2}}G\left(
0\right)  $ and $X$ is an independent data set contained in $K$.
\end{theorem}

\section{Type 2 Exact smoother error estimates $\left(  \kappa\geq1\right)  $}

These will just be the Type 1 estimates of the previous section but with the
constants from the most appropriate result given in Subsection
\ref{SbSect_ex_smth_typ2_error}.

\section{Convergence of the Approximate smoother to the data function}

In Section \ref{Sect_ap_Ex_smth_minus_App_smth} several estimates were derived
for the difference between the Approximate smoother $s_{a}$ and the Exact
smoother $s_{e}$ and a copy of the estimates for the Exact smoother error
obtained in Chapter \ref{Ch_Exact_smth} was presented. In this subsection we
will combine these estimates by means of the simple triangle inequality:
$\left\vert f_{d}\left(  x\right)  -s_{a}\left(  x\right)  \right\vert
\leq\left\vert f_{d}\left(  x\right)  -s_{e}\left(  x\right)  \right\vert
+\left\vert s_{e}\left(  x\right)  -s_{a}\left(  x\right)  \right\vert $,
where $f_{d}$ is a data function. The convergence results of Subsection
\ref{SbSect_ap_gen_converg_App_to_Ex} do not involve orders of convergence so
they will not be used.

\subsection{Order-less convergence estimates}

The following result is a simple consequence of Corollary
\ref{Cor_Jsd[sig(Zk)]->Jsd(sig(Z))_2}:

\begin{corollary}
\label{Thm_ap_orderless_Appr_smth_err}Suppose:

\begin{enumerate}
\item $\Omega$ is a bounded region and $f_{d}\in X_{w}^{0}$ is a data function.

\item $s_{e}^{\left(  k\right)  }$ is a sequence of Exact smoothers generated
by the data $\left[  X_{k};\widetilde{\mathcal{E}}_{X_{k}}f_{d}\right]  $ with
$X_{k}\subset\Omega$.

\item $s_{e}^{\left(  k\right)  }$ converges uniformly pointwise to $f_{d}$ on
$\overline{\Omega}$.

\item $s_{a}^{\left(  k\right)  }$ is a sequence of Approximate smoothers
generated by $X_{k}^{\prime}\subset\Omega$ and $\left[  X_{k};f_{d}\left(
X_{k}\right)  \right]  $ with the $X_{k}^{\prime}$ satisfying $\sup
\limits_{x\in\Omega}\operatorname*{dist}\left(  x,X_{k}^{\prime}\right)
\rightarrow0$.
\end{enumerate}

Then $s_{a}^{\left(  k\right)  }$ converges uniformly pointwise to $f_{d}$ on
$\overline{\Omega}$.
\end{corollary}

\subsection{Error estimates without explicitly assuming unisolvent data sets}

\subsubsection{\protect\underline{Type 1 estimates $\left(  \kappa
\geq0\right)  $}}

The next result gives estimates for an arbitrary data function and a Riesz
data function $R_{y}$.

\begin{theorem}
\label{Thm_ap_Appr_smth_err_typ1}Suppose:

\begin{enumerate}
\item the weight function $w$ has property W02 for $\kappa=0$ and that $G$ is
the basis function generated by $w$. Assume that for some $0<s<1$ and
$C_{G},h_{G}>0$ the basis function satisfies%
\begin{equation}
G\left(  0\right)  -G\left(  x\right)  \leq C_{G}\left\vert x\right\vert
^{2s},\text{\quad}\left\vert x\right\vert <h_{G}.\label{a7.13}%
\end{equation}

\item Let $s_{e}$ be the Exact smoother generated by the data $\left[
X;\widetilde{\mathcal{E}}_{X}f_{d}\right]  $ for some data function $f_{d}\in
X_{w}^{0}$. The independent data $X$ is contained in a data region $K$ which
is a closed bounded infinite set.

\item Suppose that $s_{a}$ is the Approximate smoother generated by the data
$\left[  X;\widetilde{\mathcal{E}}_{X}f_{d}\right]  $ and the points
$X^{\prime}\subset K$.
\end{enumerate}

Then if $k_{G}=\left(  2\pi\right)  ^{-\frac{d}{4}}\sqrt{2C_{G}}$ we have the
estimates%
\begin{equation}
\left\vert f_{d}\left(  x\right)  -s_{a}\left(  x\right)  \right\vert \leq
\min\left\{
\begin{array}
[c]{l}%
\max\left\{  \sqrt{\left(  f_{d}-s_{e},f_{d}\right)  _{w,0}},\sqrt{\left(
s_{e}-s_{a},s_{e}\right)  _{w,0}}\right\}  \left(  \sqrt{\rho N}+k_{G}\left(
h_{X,K}\right)  ^{s}+k_{G}\left(  h_{X^{\prime},K}\right)  ^{s}\right)  ,\\
\max\left\{  \sqrt{\left(  f_{d}-s_{e},f_{d}\right)  _{w,0}},\sqrt{\left(
s_{e}-s_{a},s_{e}\right)  _{w,0}}\right\}  \sqrt{R_{0}\left(  0\right)  },\\
\left\Vert f_{d}\right\Vert _{\infty,K}+\sqrt{\left(  s_{e}-f_{d}%
,f_{d}\right)  _{w,0}}\left(  k_{G}\right)  ^{2}\frac{\left(  h_{X^{\prime}%
,K}\right)  ^{2s}}{\sqrt{\rho}}+\\
\qquad\qquad+\left\Vert f_{d}\right\Vert _{w,0}\sqrt{R_{0}\left(  0\right)
}\min\left\{  1,\frac{R_{0}\left(  0\right)  }{\rho}\right\}  .
\end{array}
\right. \label{a7.56}%
\end{equation}

and%
\begin{equation}
\left\vert R_{y}\left(  x\right)  -\mathcal{S}_{X,X^{\prime}}^{a}R_{y}\left(
x\right)  \right\vert \leq\eta_{e}\left(  \rho\right)  +\min\left\{
\begin{array}
[c]{l}%
k_{G}\left(  h_{X^{\prime},K}\right)  ^{s}\left(  \sqrt{\rho N}+k_{G}\left(
h_{X,K}\right)  ^{s}+k_{G}\left(  h_{X^{\prime},K}\right)  ^{s}\right)  ,\\
R_{0}\left(  0\right)  \min\left\{  1,\frac{R_{0}\left(  0\right)  }{\rho
}\right\}  ,\\
\left(  k_{G}\right)  ^{2}\frac{\left(  h_{X^{\prime},K}\right)  ^{2s}}%
{\sqrt{\rho}}\sqrt{\eta_{e}\left(  \rho\right)  }.
\end{array}
\right. \label{a7.57}%
\end{equation}

where%
\begin{equation}
\eta_{e}\left(  \rho\right)  =\min\left\{
\begin{array}
[c]{l}%
\left(  \sqrt{\rho N}+k_{G}\left(  h_{X,K}\right)  ^{s}\right)  ^{2},\\
R_{0}\left(  0\right)  ,\\
\left\Vert R_{y}\right\Vert _{\infty,K}+R_{0}\left(  0\right)  \min\left\{
1,\frac{R_{0}\left(  0\right)  }{\rho}\right\}  .
\end{array}
\right. \label{a7.10}%
\end{equation}

when $h_{X,K}=\sup\limits_{x\in K}\operatorname*{dist}\left(  x,X\right)  \leq
h_{G}$ and $h_{X^{\prime},K}=\sup\limits_{x\in K}\operatorname*{dist}\left(
x,X^{\prime}\right)  \leq h_{G}$.
\end{theorem}

\begin{proof}
From Theorem \ref{Thm_ap_typ1_Exact_smth_err_summary} and Corollary
\ref{Cor_Thm_Se-Sa_large_smth_parm},
\[
\left\vert f_{d}\left(  x\right)  -s_{e}\left(  x\right)  \right\vert \leq
\min\left\{
\begin{array}
[c]{l}%
\sqrt{\left(  f_{d}-s_{e},f_{d}\right)  _{w,0}}\left(  \sqrt{\rho N}%
+k_{G}\left(  h_{X,K}\right)  ^{s}\right)  ,\\
\sqrt{\left(  f_{d}-s_{e},f_{d}\right)  _{w,0}}\sqrt{R_{0}\left(  0\right)
},\\
\left\Vert f_{d}\right\Vert _{\infty,K}+\left\Vert f_{d}\right\Vert
_{w,0}\sqrt{R_{0}\left(  0\right)  }\min\left\{  1,\frac{R_{0}\left(
0\right)  }{\rho}\right\}  .
\end{array}
\right.
\]

and%
\[
\left\vert s_{e}\left(  x\right)  -s_{a}\left(  x\right)  \right\vert \leq
\min\left\{
\begin{array}
[c]{l}%
\sqrt{\left(  s_{e}-s_{a},s_{e}\right)  _{w,0}}k_{G}\left(  h_{X^{\prime}%
,K}\right)  ^{s},\\
\sqrt{\left(  s_{e}-s_{a},s_{e}\right)  _{w,0}}\sqrt{R_{0}\left(  0\right)
},\\
\sqrt{\left(  s_{e}-f_{d},f_{d}\right)  _{w,0}}\left(  k_{G}\right)  ^{2}%
\frac{\left(  h_{X^{\prime},K}\right)  ^{2s}}{\sqrt{\rho}}.
\end{array}
\right.
\]

so that%
\begin{align*}
\left\vert f_{d}\left(  x\right)  -s_{a}\left(  x\right)  \right\vert  &
\leq\left\vert f_{d}\left(  x\right)  -s_{e}\left(  x\right)  \right\vert
+\left\vert s_{e}\left(  x\right)  -s_{a}\left(  x\right)  \right\vert \\
& \leq\min\left\{
\begin{array}
[c]{l}%
\max\left\{  \sqrt{\left(  f_{d}-s_{e},f_{d}\right)  _{w,0}},\sqrt{\left(
s_{e}-s_{a},s_{e}\right)  _{w,0}}\right\}  \left(  \sqrt{\rho N}+k_{G}\left(
h_{X,K}\right)  ^{s}+k_{G}\left(  h_{X^{\prime},K}\right)  ^{s}\right)  ,\\
\max\left\{  ,\right\}  \sqrt{R_{0}\left(  0\right)  },\\
\left\Vert f_{d}\right\Vert _{\infty,K}+\left\Vert f_{d}\right\Vert
_{w,0}\left(  k_{G}\right)  ^{2}\frac{\left(  h_{X^{\prime},K}\right)  ^{2s}%
}{\sqrt{\rho}}+\\
\qquad\qquad+\left\Vert f_{d}\right\Vert _{w,0}\sqrt{R_{0}\left(  0\right)
}\min\left\{  1,\frac{R_{0}\left(  0\right)  }{\rho}\right\}  .
\end{array}
\right.
\end{align*}

Using $\eta_{e}\left(  \rho\right)  $ the estimate \ref{a7.191} for
$\left\vert \mathcal{S}_{X}R_{y}\left(  x\right)  -\mathcal{S}_{X,X^{\prime}%
}^{a}R_{y}\left(  x\right)  \right\vert $ can be written
\begin{align*}
\left\vert \mathcal{S}_{X}R_{y}\left(  x\right)  -\mathcal{S}_{X,X^{\prime}%
}^{a}R_{y}\left(  x\right)  \right\vert  & \leq\min\left\{
\begin{array}
[c]{l}%
k_{G}\left(  h_{X^{\prime},K}\right)  ^{s}\left(  \sqrt{\rho N}+k_{G}\left(
h_{X,K}\right)  ^{s}+k_{G}\left(  h_{X^{\prime},K}\right)  ^{s}\right)  ,\\
R_{0}\left(  0\right)  \min\left\{  1,\frac{R_{0}\left(  0\right)  }{\rho
}\right\}  ,\\
\left(  k_{G}\right)  ^{2}\frac{\left(  h_{X^{\prime},K}\right)  ^{2s}}%
{\sqrt{\rho}}\left(  \sqrt{\rho N}+k_{G}\left(  h_{X,K}\right)  ^{s}\right)
,\\
\left(  k_{G}\right)  ^{2}\frac{\left(  h_{X^{\prime},K}\right)  ^{2s}}%
{\sqrt{\rho}}\sqrt{R_{0}\left(  0\right)  },\\
\left(  k_{G}\right)  ^{2}\frac{\left(  h_{X^{\prime},K}\right)  ^{2s}}%
{\sqrt{\rho}}\sqrt{\left\Vert R_{y}\right\Vert _{\infty,K}+R_{0}\left(
0\right)  \min\left\{  1,\frac{R_{0}\left(  0\right)  }{\rho}\right\}  }.
\end{array}
\right. \\
& \leq\min\left\{
\begin{array}
[c]{l}%
k_{G}\left(  h_{X^{\prime},K}\right)  ^{s}\left(  \sqrt{\rho N}+k_{G}\left(
h_{X,K}\right)  ^{s}+k_{G}\left(  h_{X^{\prime},K}\right)  ^{s}\right)  ,\\
R_{0}\left(  0\right)  \min\left\{  1,\frac{R_{0}\left(  0\right)  }{\rho
}\right\}  ,\\
\left(  k_{G}\right)  ^{2}\frac{\left(  h_{X^{\prime},K}\right)  ^{2s}}%
{\sqrt{\rho}}\sqrt{\eta\left(  \rho\right)  }.
\end{array}
\right.
\end{align*}

Now%
\[
\left\vert R_{y}\left(  x\right)  -\mathcal{S}_{X,X^{\prime}}^{a}R_{y}\left(
x\right)  \right\vert \leq\left\vert R_{y}\left(  x\right)  -\left(
\mathcal{S}_{X}R_{y}\right)  \left(  x\right)  \right\vert +\left\vert
\mathcal{S}_{X}R_{y}\left(  x\right)  -\mathcal{S}_{X,X^{\prime}}^{a}%
R_{y}\left(  x\right)  \right\vert ,
\]

and since the estimate for $\left\vert R_{y}\left(  x\right)  -\left(
\mathcal{S}_{X}R_{y}\right)  \left(  x\right)  \right\vert $ from \ref{a7.22}
can be written $\left\vert R_{y}\left(  x\right)  -\left(  \mathcal{S}%
_{X}R_{y}\right)  \left(  x\right)  \right\vert \leq\eta\left(  \rho\right)  $
we have \ref{a7.57}.
\end{proof}

\begin{remark}
\label{Rem_Thm_Err_not_unisolv_k=0_arb}Define $h_{a}=\left(  \left(
h_{X,K}\right)  ^{s}+\left(  h_{X^{\prime},K}\right)  ^{s}\right)  ^{1/s}$. Then:

\begin{enumerate}
\item From \ref{a7.56}%
\[
\left\vert f_{d}\left(  x\right)  -s_{a}\left(  x\right)  \right\vert \leq
\max\left\{  \sqrt{\left(  f_{d}-s_{e},f_{d}\right)  _{w,0}},\sqrt{\left(
s_{e}-s_{a},s_{e}\right)  _{w,0}}\right\}  \left(  \sqrt{\rho N}+k_{G}\left(
h_{a}\right)  ^{s}\right)  ,\quad x\in K,
\]

and the order of convergence is $s$ in $h_{a}$.

\item From \ref{a7.57}: for $x,y\in K$%
\begin{align*}
\left\vert R_{y}\left(  x\right)  -\mathcal{S}_{X,X^{\prime}}^{a}R_{y}\left(
x\right)  \right\vert  & \leq\left(  \sqrt{\rho N}+k_{G}\left(  h_{X,K}%
\right)  ^{s}\right)  ^{2}+\sqrt{\rho N}k_{G}\left(  h_{X^{\prime},K}\right)
^{s}+\\
& \qquad\qquad+\left(  k_{G}\right)  ^{2}\left(  h_{X^{\prime},K}\right)
^{s}\left(  h_{X^{\prime},K}\right)  ^{s}+\left(  k_{G}\right)  ^{2}\left(
h_{X,K}\right)  ^{2s}\\
& \leq\left(  \sqrt{\rho N}+k_{G}\left(  h_{X,K}\right)  ^{s}+k_{G}\left(
h_{X^{\prime},K}\right)  ^{s}\right)  ^{2}\\
& =\left(  \sqrt{\rho N}+k_{G}\left(  h_{a}\right)  ^{s}\right)  ^{2},
\end{align*}

and the order of convergence is $2s$ in $h_{a}$.

\item The last two Approximate smoother error estimates are analogous to the
Exact smoother error estimates%
\[
\left\vert f_{d}\left(  x\right)  -\left(  \mathcal{S}_{X}f_{d}\right)
\left(  x\right)  \right\vert \leq\sqrt{\left(  f_{d}-s_{e},f_{d}\right)
_{w,0}}\left(  \sqrt{\rho N}+k_{G}\left(  h_{X,K}\right)  ^{s}\right)  ,
\]

and%
\[
\left\vert R_{y}\left(  x\right)  -\left(  \mathcal{S}_{X}R_{y}\right)
\left(  x\right)  \right\vert \leq\left(  \sqrt{\rho N}+k_{G}\left(
h_{X,K}\right)  ^{s}\right)  ^{2},
\]

from \ref{a7.11} and \ref{a7.22}.
\end{enumerate}
\end{remark}

\subsubsection{\protect\underline{Type 2 error estimates $\left(  \kappa
\geq1\right)  $}}

We want to derive Type 2 error estimates for $\left\vert f_{d}\left(
x\right)  -s_{a}\left(  x\right)  \right\vert $. In fact, these results can be
proved by combining the Type 2 Exact smoother error estimates from Subsection
\ref{SbSect_ex_smth_typ2_error} with the Type 2 estimates for $\left\vert
s_{e}\left(  x\right)  -s_{a}\left(  x\right)  \right\vert $ from
Subsubsection \ref{SbSbSect_ap_Typ2_ApprMinusExact}.

\begin{theorem}
\label{Thm_ap_Appr_smth_err_typ2}Suppose a weight function satisfies property
W02 for some $\kappa\geq1$ and denote the basis function by $G$. Then the
Approximate smoother error estimates \ref{a7.18} and \ref{a7.191} hold for
\[
C_{G}=-\frac{1}{2}\left(  \left\vert D\right\vert ^{2}G\right)  \left(
0\right)  d,\text{\quad}s=1,\text{\quad}h_{G}=\infty,\text{\quad}k_{G}=\left(
2\pi\right)  ^{-\frac{d}{4}}\sqrt{-\left(  \left\vert D\right\vert
^{2}G\right)  \left(  0\right)  }\sqrt{d}.
\]

\end{theorem}

However, if the weight function is radial we can use the estimates of:

\begin{theorem}
Suppose a radial weight function satisfies property W02 for $\kappa=1$ and
denote the (radial) basis function by $G$. Set $r=\left\vert x\right\vert $. Then:

\begin{enumerate}
\item If $G\left(  x\right)  =f\left(  r^{2}\right)  $ the Approximate
smoother error estimates of Theorem \ref{Thm_ap_Appr_smth_err_typ1} hold for%
\[
C_{G}=-f^{\prime}\left(  0\right)  d^{2},\text{\quad}s=1,\text{\quad}%
h_{G}=\infty,\text{\quad}k_{G}=\left(  2\pi\right)  ^{-\frac{d}{4}}%
\sqrt{-2f^{\prime}\left(  0\right)  }d.
\]

\item If $G\left(  x\right)  =g\left(  r\right)  $ then the Approximate
smoother error estimates of Theorem \ref{Thm_ap_Appr_smth_err_typ1} hold for
\[
C_{G}=-\frac{1}{2}g^{\prime\prime}\left(  0\right)  d^{2},\text{\quad
}s=1,\text{\quad}h_{G}=\infty,\text{\quad}k_{G}=\left(  2\pi\right)
^{-\frac{d}{4}}\sqrt{-g^{\prime\prime}\left(  0\right)  }d.
\]

\end{enumerate}
\end{theorem}

If the weight function is a tensor product we can use the estimates of:

\begin{theorem}
\label{Thm_Final_Type2_estim_3}Suppose a tensor product weight function
satisfies property W02 for $\kappa=1$ and denote the univariate basis function
by $G_{1}$. Then the Approximate smoother error estimates of Theorem
\ref{Thm_ap_Appr_smth_err_typ1} hold when
\[
C_{G}=-\frac{d}{2}G_{1}\left(  0\right)  ^{d-1}D^{2}G_{1}\left(  0\right)
,\text{\quad}s=1,\text{\quad}h_{G}=\infty,\text{\quad}k_{G}=\left(
2\pi\right)  ^{-\frac{d}{4}}\sqrt{-G_{1}\left(  0\right)  ^{d-1}D^{2}%
G_{1}\left(  0\right)  }\sqrt{d}.
\]

\end{theorem}

\subsection{Estimates using minimally unisolvent data sets}

The next result is an error estimate for an arbitrary data function in
$X_{w}^{0}$.

\begin{theorem}
\label{Thm_Err_unisolv_arbitrary_fn}Suppose:

\begin{enumerate}
\item $w$ is a weight function with property W02 for some $\kappa\geq1$ and
let $G$ be the basis function generated by $w$. Set $m=\left\lfloor
\kappa\right\rfloor $.

\item $s_{e}$ is the Exact smoother generated by the data $\left[
X;\widetilde{\mathcal{E}}_{X}f_{d}\right]  $ for some data function $f_{d}\in
X_{w}^{0}$, and $X$ is contained in the data region $\Omega$.

\item $s_{a}$ is the Approximate smoother generated by $\left[
X;\widetilde{\mathcal{E}}_{X}f_{d}\right]  $ and the set $X^{\prime}%
\subset\Omega$.

\item We use the notation and assumptions of Lemma \ref{Lem_Lagrange_interpol}
which means here that $X^{\prime}$ and $X$ are $m$-unisolvent and that
$\Omega$ is a bounded region whose boundary satisfies the cone condition.
\end{enumerate}

Then there exist positive $h_{\Omega,m},c_{\Omega,m},K_{\Omega,m}^{\prime}$
and $k_{G}=\frac{d^{m/2}}{\left(  2\pi\right)  ^{d/2}}\left(  c_{\Omega
,m}\right)  ^{m}K_{\Omega,m}^{\prime}\max\limits_{\left\vert \beta\right\vert
=m}\left\vert D^{2\beta}G\left(  0\right)  \right\vert $ such that%
\begin{multline*}
\left\vert f_{d}\left(  x\right)  -s_{a}\left(  x\right)  \right\vert
\leq\sqrt{\left(  f_{d}-s_{e},f_{d}\right)  _{w,0}}\left(  K_{\Omega
,m}^{\prime}\sqrt{\rho N}+k_{G}\left(  h_{X,\Omega}\right)  ^{m}\right)  +\\
+\sqrt{\left(  s_{e}-s_{a},s_{e}\right)  _{w,0}}k_{G}\left(  h_{X^{\prime
},\Omega}\right)  ^{m},\text{\quad}x\in\overline{\Omega},
\end{multline*}

\end{theorem}

when $h_{X,\Omega}=\sup\limits_{x\in\Omega}\operatorname*{dist}\left(
x,X\right)  \leq h_{G}$ and $h_{X^{\prime},\Omega}=\sup\limits_{x\in\Omega
}\operatorname*{dist}\left(  x,X^{\prime}\right)  \leq h_{G}$.

Further,$\sqrt{\left(  s_{e}-s_{a},s_{e}\right)  _{w,0}}\leq\left\Vert
s_{e}\right\Vert _{w,0}\leq\left\Vert f_{d}\right\Vert _{w,0}$ and
$\sqrt{\left(  f_{d}-s_{e},f_{d}\right)  _{w,0}}\leq\left\Vert f_{d}%
\right\Vert _{w,0}$.

\begin{proof}
From Theorem \ref{Thm_ConvergSaToSe_not_unisolv_II}%
\[
\left\vert s_{e}\left(  x\right)  -s_{a}\left(  x\right)  \right\vert
\leq\sqrt{\left(  s_{e}-s_{a},s_{e}\right)  _{w,0}}k_{G}\left(  h_{X^{\prime
},\Omega}\right)  ^{m},\quad x\in\overline{\Omega},
\]

when $h_{X^{\prime},\Omega}<h_{\Omega,m}$. From inequality \ref{7.20}
\[
\left\vert f_{d}\left(  x\right)  -s_{e}\left(  x\right)  \right\vert
<\sqrt{\left(  f_{d}-s_{e},f_{d}\right)  _{w,0}}\left(  K_{\Omega,m}^{\prime
}\sqrt{\rho N}+k_{G}\left(  h_{X,\Omega}\right)  ^{m}\right)  ,\quad
x\in\overline{\Omega},
\]

when $h_{X,\Omega}\leq h_{\Omega,m}$. Hence when $x\in\overline{\Omega}$%
\begin{align*}
\left\vert f_{d}\left(  x\right)  -s_{a}\left(  x\right)  \right\vert  &
\leq\left\vert f_{d}\left(  x\right)  -s_{e}\left(  x\right)  \right\vert
+\left\vert s_{e}\left(  x\right)  -s_{a}\left(  x\right)  \right\vert \\
& \leq\sqrt{\left(  f_{d}-s_{e},f_{d}\right)  _{w,0}}\left(  K_{\Omega
,m}^{\prime}\sqrt{\rho N}+k_{G}\left(  h_{X,\Omega}\right)  ^{m}\right)  +\\
& \qquad+\sqrt{\left(  s_{e}-s_{a},s_{e}\right)  _{w,0}}k_{G}\left(
h_{X^{\prime},\Omega}\right)  ^{m}.
\end{align*}

The final inequalities follow from Theorem
\ref{Thm_ConvergSaToSe_not_unisolv_II} and Corollary
\ref{Cor_OrdConvergExactSmth_k=0_Fd}.
\end{proof}

\begin{remark}
\label{Rem_Thm_Err_unisolv_arbitrary_fn}We can combine the two density
measures $h_{X,\Omega}$ and $h_{X^{\prime},\Omega}$ into one measure
\[
h_{a}=\left(  \left(  h_{X,\Omega}\right)  ^{m}+\left(  h_{X^{\prime},\Omega
}\right)  ^{m}\right)  ^{1/m}.
\]
From Theorem \ref{Thm_ap_Appr_smth_err_typ1} we have $\sqrt{\left(
f_{d}-s_{e},f_{d}\right)  _{w,0}}\leq\left\Vert f_{d}\right\Vert _{w,0}$ and
$\sqrt{\left(  s_{e}-s_{a},s_{e}\right)  _{w,0}}\leq\left\Vert f_{d}%
\right\Vert _{w,0}$ so that
\begin{align*}
\left\vert f_{d}\left(  x\right)  -s_{a}\left(  x\right)  \right\vert  &
\leq\sqrt{\left(  f_{d}-s_{e},f_{d}\right)  _{w,0}}\left(  K_{\Omega
,m}^{\prime}\sqrt{\rho N}+k_{G}\left(  h_{X,\Omega}\right)  ^{m}\right)  +\\
& \qquad+\sqrt{\left(  s_{e}-s_{a},s_{e}\right)  _{w,0}}k_{G}\left(
h_{X^{\prime},\Omega}\right)  ^{m}\\
& \leq\left\Vert f_{d}\right\Vert _{w,0}\left(  K_{\Omega,m}^{\prime}%
\sqrt{\rho N}+k_{G}\left(  h_{X,\Omega}\right)  ^{m}+k_{G}\left(
h_{X^{\prime},\Omega}\right)  ^{m}\right) \\
& \leq\left\Vert f_{d}\right\Vert _{w,0}\left(  K_{\Omega,m}^{\prime}%
\sqrt{\rho N}+k_{G}\left(  h_{a}\right)  ^{m}\right)  ,
\end{align*}

so that when $\rho=0$ the order of convergence in $h_{a}$ is $m=\left\lfloor
\kappa\right\rfloor $.
\end{remark}

An important application of the Approximate smoother is the case where
$X^{\prime}$ is a regular rectangular grid. To apply the last theorem we must
show that if the grid size decreases then $X^{\prime}$ eventually contains a
$m$-unisolvent subset. The next theorem confirms this. The result of the next
theorem is obvious in one-dimension but results for polynomials in higher
dimensions can be complex and/or unexpected. However, here things are OK. We
will need the following lemma:

\begin{lemma}
We have the following unisolvency results:

\begin{enumerate}
\item The set $\left\{  \gamma\in\mathbb{Z}^{d}:0\leq\gamma<n\right\}  $ is
unisolvent w.r.t. $P_{n}$.

\item Translations of minimal unisolvent sets are minimal unisolvent sets.

\item Dilations of minimal unisolvent sets are minimal unisolvent sets.
\end{enumerate}
\end{lemma}

\begin{proof}
\textbf{Part 1}. From the definition of unisolvency, Definition
\ref{Def_unisolv}, we must show that for each $p\in P_{n}$, $p\left(
\gamma\right)  =0$ for $0\leq\gamma<n$ implies $p=0$. The proof will be by
induction on the order of the polynomial.

Clearly the lemma is true for $n=1$ since $P_{1}$ is the constant polynomials.

Now assume that $n\geq2$ and that if $p\in P_{n}$ and $p\left(  \gamma\right)
=0$ for $0\leq\gamma<n$ then $p=0$. Set $p\left(  x\right)  =\sum
\limits_{\left\vert \beta\right\vert <n}c_{\beta}x^{\beta}$. Then if
$\gamma_{k}=0$ and $\gamma_{i}=1$ when $i\neq k$ then $0=\sum
\limits_{\left\vert \beta\right\vert <n}c_{\beta}\gamma^{\beta}$ implies
$c_{\beta}=0$ when $\beta_{k}=0$. Thus $p\left(  \gamma\right)  =0$ for
$0\leq\gamma<n$ implies $c_{\beta}=0$ when $\beta_{i}=0$ for some $i$.
Consequently, $p=0$ if $n\leq d$ else $p\left(  x\right)  =\sum
\limits_{\substack{\left\vert \beta\right\vert <n \\\beta>0}}c_{\beta}%
x^{\beta}$. If $n>d$ we can write $p\left(  x\right)  =x^{\mathbf{1}}q\left(
x\right)  $ where $q\in P_{n-1}$, so that $q\left(  \gamma\right)  =0$ for
$1\leq\gamma<n$ must imply $q=0$. Finally, if we define $r\in P_{n-1}$ by
$r\left(  x\right)  =q\left(  x+1\right)  $ then $r\left(  \gamma\right)  =0$
for $0\leq\gamma<n-1$ must imply $r=0$, and so the truth of our lemma for
$n-1$ implies the truth of the lemma for $n$, and the lemma is proved.\medskip

\textbf{Parts 2 and 3}. From the definition of a cardinal basis, Definition
\ref{Def_cardinal_basis}, there is a unique cardinal basis $\left\{
l_{i}\right\}  $ of $P_{n}$ associated with a set $A=\left\{  a_{i}\right\}  $
iff $A$ is minimally unisolvent of order $n$. Now by definition $l_{i}\left(
a_{j}\right)  =\delta_{i,j}$. Hence, if $\tau,\delta\in\mathbb{R}^{d}$ and
$\delta$ has positive components, then the cardinal basis associated with the
translation $A+\tau$ is $\left\{  l_{i}\left(  \cdot-\tau\right)  \right\}  $
and the cardinal basis associated with the dilation $\delta A$ is $\left\{
l_{i}\left(  \cdot/\delta\right)  \right\}  $.
\end{proof}

\begin{theorem}
\label{Thm_X_unisolv}Suppose $X^{\prime}=\left\{  x_{\alpha}^{\prime
}=a+h\alpha\mid\alpha\in\mathbb{Z}^{d}\text{ and }0\leq\alpha<\mathcal{N}%
^{\prime}\right\}  $ is the regular, rectangular grid introduced in Definition
\ref{Def_mesh}.

Then $X^{\prime}$ is $m$-unisolvent if $\mathcal{N}^{\prime}\geq m$.
\end{theorem}

\begin{proof}
Since $\left(  X^{\prime}-a\right)  /h=\left\{  \alpha\mid\alpha\in
\mathbb{Z}^{d}\text{ and }0\leq\alpha<\mathcal{N}^{\prime}\right\}  $ and
$\mathcal{N}^{\prime}\geq m$, part 1 of the lemma implies $\left(  X^{\prime
}-a\right)  /h$ is $m$-unisolvent and thus from the definition of unisolvency,
Definition \ref{Def_unisolv}, $\left(  X^{\prime}-a\right)  /h$ must contain a
minimal unisolvent subset. Parts 2 and 3 of the lemma imply that $X^{\prime}$
contains a minimal unisolvent subset and so $X^{\prime}$ is unisolvent.
\end{proof}

We can improve the last convergence result when the data functions have the
form $R_{y}$.

\begin{theorem}
\label{Thm_Err_unisolv_Ry}Suppose:

\begin{enumerate}
\item $w$ is a weight function with property W02 for some $\kappa\geq1$ and
let $G$ be the basis function generated by $w$. Set $m=\left\lfloor
\kappa\right\rfloor $.

\item $\mathcal{S}_{X}R_{y}$ is the Exact smoother generated by the data
$\left[  X;\widetilde{\mathcal{E}}_{X}R_{y}\right]  $ for some data function
$R_{y}$, and $X$ is contained in the data region $\Omega$.

\item $\mathcal{S}_{X,X^{\prime}}^{a}R_{y}$ is the Approximate smoother
generated by $\left[  X;\widetilde{\mathcal{E}}_{X}R_{y}\right]  $ and the set
$X^{\prime}$ is contained in the data region $\Omega$.

\item We use the notation and assumptions of Lemma \ref{Lem_Lagrange_interpol}
which means here that $X^{\prime}$ and $X$ are $m$-unisolvent and that
$\Omega$ is a bounded region whose boundary satisfies the cone condition.
\end{enumerate}

Then there exist positive $h_{\Omega,m},c_{\Omega,m},K_{\Omega,m}^{\prime}$
and $k_{G}=\frac{d^{m/2}}{\left(  2\pi\right)  ^{d/2}}\left(  c_{\Omega
,m}\right)  ^{m}K_{\Omega,m}^{\prime}\max\limits_{\left\vert \beta\right\vert
=m}\left\vert D^{2\beta}G\left(  0\right)  \right\vert $ such that%
\[
\left\vert R_{y}\left(  x\right)  -\mathcal{S}_{X,X^{\prime}}^{a}R_{y}\left(
x\right)  \right\vert \leq\left(  K_{\Omega,m}^{\prime}\sqrt{\rho N}%
+k_{G}\left(  h_{X,\Omega}\right)  ^{m}+k_{G}\left(  h_{X^{\prime},\Omega
}\right)  ^{m}\right)  ^{2},
\]

for $x,y\in\overline{\Omega}$, when $h_{X,\Omega}=\sup\limits_{x\in\Omega
}\operatorname*{dist}\left(  x,X\right)  \leq h_{G}$ and $h_{X^{\prime}%
,\Omega}=\sup\limits_{x\in\Omega}\operatorname*{dist}\left(  x,X^{\prime
}\right)  \leq h_{G}$.
\end{theorem}

\begin{proof}
When $h_{X,\Omega}\leq h_{\Omega,m}$ and $h_{X^{\prime},\Omega}\leq
h_{\Omega,m}$ Theorem \ref{Thm_converg_Rx_func_Kgt0} implies
\[
\left\vert \mathcal{S}_{X}R_{y}\left(  x\right)  -\mathcal{S}_{X,X^{\prime}%
}^{a}R_{y}\left(  x\right)  \right\vert \leq k_{G}^{\prime}\left(
h_{X^{\prime},\Omega}\right)  ^{m}\left(  K_{\Omega,m}^{\prime}\sqrt{\rho
N}+k_{G}\left(  h_{X,\Omega}\right)  ^{m}+k_{G}\left(  h_{X^{\prime},\Omega
}\right)  ^{m}\right)  ,
\]

for $x\in\overline{\Omega}$, and \ref{a7.22} implies that
\[
\left\vert R_{y}\left(  x\right)  -\mathcal{S}_{X}R_{y}\left(  x\right)
\right\vert \leq\left(  K_{\Omega,m}^{\prime}\sqrt{\rho N}+k_{G}\left(
h_{X,\Omega}\right)  ^{m}\right)  ^{2},\quad x,y\in\overline{\Omega}.
\]

Hence when $x\in\overline{\Omega}$%
\begin{align*}
\left\vert R_{y}\left(  x\right)  -\mathcal{S}_{X,X^{\prime}}^{a}R_{y}\left(
x\right)  \right\vert  & \leq\left\vert R_{y}\left(  x\right)  -\mathcal{S}%
_{X}R_{y}\left(  x\right)  \right\vert +\left\vert \mathcal{S}_{X}R_{y}\left(
x\right)  -\mathcal{S}_{X,X^{\prime}}^{a}R_{y}\left(  x\right)  \right\vert \\
& \leq\left(  K_{\Omega,m}^{\prime}\sqrt{\rho N}+k_{G}\left(  h_{X,\Omega
}\right)  ^{m}\right)  ^{2}+\\
& \qquad k_{G}^{\prime}\left(  h_{X^{\prime},\Omega}\right)  ^{m}\left(
K_{\Omega,m}^{\prime}\sqrt{\rho N}+k_{G}\left(  h_{X,\Omega}\right)
^{m}+k_{G}\left(  h_{X^{\prime},\Omega}\right)  ^{m}\right) \\
& <\left(  K_{\Omega,m}^{\prime}\sqrt{\rho N}+k_{G}\left(  h_{X,\Omega
}\right)  ^{m}+k_{G}\left(  h_{X^{\prime},\Omega}\right)  ^{m}\right)  ^{2}.
\end{align*}

\end{proof}

\begin{remark}
\label{Rem_Thm_Err_unisolv_Ry}We have%
\[
\left\vert R_{y}\left(  x\right)  -\mathcal{S}_{X,X^{\prime}}^{a}R_{y}\left(
x\right)  \right\vert \leq\left(  K_{\Omega,m}^{\prime}\sqrt{\rho N}%
+k_{G}\left(  h_{a}\right)  ^{m}\right)  ^{2},\quad x,y\in\overline{\Omega},
\]

where $h_{a}=\left(  \left(  h_{X,\Omega}\right)  ^{m}+\left(  h_{X^{\prime
},\Omega}\right)  ^{m}\right)  ^{1/m}$ so that when $\rho=0$ the order of
convergence is $2m=2\left\lfloor \kappa\right\rfloor $. This compares with
order $m=\left\lfloor \kappa\right\rfloor $ obtained in Remark
\ref{Rem_Thm_Err_unisolv_arbitrary_fn} in the case of an arbitrary data function.
\end{remark}

\subsection{??}

??

\section{Numerical experiments with scaled, extended B-splines}

In this section the convergence of the Approximate smoother to its data
function is studied numerically using scaled, extended B-splines, the first of
which will be the hat function. We will use the same 1-dimensional data
functions and B-spline basis functions \ref{7.67} that were used for the
numerical interpolation experiments described in Section
\ref{Sect_num_experim}. We will only consider the numerical experiments in one
dimension to allow the data density parameter $h_{X,K}=\sup\limits_{x\in
K}\operatorname*{dist}\left(  x,X\right)  $ to be easily calculated.

It is easier to validate the non-unisolvent error estimates i.e. those derived
without explicitly assuming unisolvent data sets. This is because in the
non-unisolvent case the constants can be known precisely. In the case of
unisolvent data the theory is complex and it is unclear what are suitable
upper bounds for the constants $K_{\Omega,m}^{\prime}$, $c_{\Omega,m}$,
$h_{\Omega,m}$ in Theorem \ref{Thm_Err_unisolv_arbitrary_fn} - constants are
"buried" in the analysis.

Because all the extended B-spline basis weight functions have a power of $\sin
x$ in the denominator we will use the special classes of data functions
developed in Section \ref{Sect_int_data_fn_exten_Bsplin}. The multivariate
theory of local data spaces $X_{w}^{0}\left(  \Omega\right)  $ where $\Omega$
is a bounded region was developed in Sections \ref{Sect_local_data_space} and
\ref{Sect_CntDifWtFn_DataFuncs}.

The 1-dimensional data set is constructed using a uniform distribution on the
data region $K=\left[  -1.5,1.5\right]  $.

We now assume that $X^{\prime}$ is a regular grid so that
\begin{equation}
h_{X^{\prime},K}=\frac{3}{2N^{\prime}},\label{a7.21}%
\end{equation}

where $N^{\prime}\geq2$ is the number of grid cells.

Each of 20 data files is exponentially sampled using a multiplier of
approximately 1.3 and then we plot $\log_{10}h_{X,K}$ against $\log_{10}N$
where $N=\left\vert X\right\vert $. It then seems quite reasonable to use a
least-squares linear fit and in this case we obtain%
\begin{equation}
h_{X,K}\simeq3.09N^{-0.81}.\label{a8.75}%
\end{equation}

For ease of calculation let
\begin{equation}
h_{X,K}=h_{1}N^{-a},\text{\quad}h_{1}=3.09,\text{\quad}a=0.81.\label{a8.6}%
\end{equation}

Noting the error estimates of Theorem \ref{Thm_ap_Appr_smth_err_typ1} we use
\ref{a8.6} to write $h_{X,K}^{s}=h_{1}^{s}N^{-as}$ and define the Exact
smoother error estimates%
\begin{equation}
\left\vert f_{d}\left(  x\right)  -\left(  \mathcal{S}_{X}f_{d}\right)
\left(  x\right)  \right\vert \lesssim E_{N}\left(  \rho\right)  ,\text{\quad
}x\in K,\label{a8.1}%
\end{equation}

and%
\begin{equation}
\left\vert R_{y}\left(  x\right)  -\left(  \mathcal{S}_{X}R_{y}\right)
\left(  x\right)  \right\vert \lesssim\varepsilon_{N}\left(  \rho\right)
,\text{\quad}x,y\in K,\label{a8.11}%
\end{equation}

by%
\begin{equation}
E_{N}\left(  \rho\right)  =\min\left\{
\begin{array}
[c]{l}%
\max\left\{  \sqrt{\left(  f_{d}-s_{e},f_{d}\right)  _{w,0}},\sqrt{\left(
s_{e}-s_{a},s_{e}\right)  _{w,0}}\right\}  \left(  \sqrt{\rho N}+k_{G}%
\frac{h_{1}^{s}}{N^{as}}+k_{G}\left(  \frac{3}{2N^{\prime}}\right)
^{s}\right)  ,\\
\max\left\{  \sqrt{\left(  f_{d}-s_{e},f_{d}\right)  _{w,0}},\sqrt{\left(
s_{e}-s_{a},s_{e}\right)  _{w,0}}\right\}  \sqrt{R_{0}\left(  0\right)  },\\
\left\Vert f_{d}\right\Vert _{\infty,K}+\left\Vert f_{d}\right\Vert
_{w,0}\left(  k_{G}\right)  ^{2}\left(  \frac{3}{2N^{\prime}}\right)
^{2s}\frac{1}{\sqrt{\rho}}+\\
\qquad\qquad+\left\Vert f_{d}\right\Vert _{w,0}\sqrt{R_{0}\left(  0\right)
}\min\left\{  1,\frac{R_{0}\left(  0\right)  }{\rho}\right\}  .
\end{array}
\right. \label{a8.10}%
\end{equation}

and%
\begin{equation}
\varepsilon_{N}\left(  \rho\right)  =\eta_{e}\left(  \rho\right)
+\min\left\{
\begin{array}
[c]{l}%
k_{G}\left(  \frac{3}{2N^{\prime}}\right)  ^{s}\left(  \sqrt{\rho N}%
+k_{G}\frac{h_{1}^{s}}{N^{as}}+k_{G}\left(  \frac{3}{2N^{\prime}}\right)
^{s}\right)  ,\\
R_{0}\left(  0\right)  \min\left\{  1,\frac{R_{0}\left(  0\right)  }{\rho
}\right\}  ,\\
\left(  k_{G}\right)  ^{2}\left(  \frac{3}{2N^{\prime}}\right)  ^{2s}%
\sqrt{\frac{\eta_{e}\left(  \rho\right)  }{\rho}}.
\end{array}
\right. \label{a7.26}%
\end{equation}

where%
\begin{equation}
\eta_{e}\left(  \rho\right)  =\min\left\{
\begin{array}
[c]{l}%
\left(  \sqrt{\rho N}+k_{G}\frac{h_{1}^{s}}{N^{as}}\right)  ^{2},\\
R_{0}\left(  0\right)  ,\\
\left\Vert R_{y}\right\Vert _{\infty,K}+R_{0}\left(  0\right)  \min\left\{
1,\frac{R_{0}\left(  0\right)  }{\rho}\right\}  .
\end{array}
\right. \label{a7.27}%
\end{equation}

\subsection{Extended B-splines with $n=1$}

\subsubsection{\protect\underline{The case $n=1$, $l=1$}}

The 1-dimensional hat weight function $w_{s}$ is given by $w_{s}\left(
\xi\right)  =\frac{\xi^{2}}{\sin^{2}\xi}$ and was discussed in Subsections
\ref{SbSect_motiv_weight_fn} and \ref{SbSect_HatWeightFunc}. This is a scaled,
extended B-spline basis function corresponding to $n=l=1$ and thus Theorem
\ref{Thm_ex_splin_wt_fn_properties} and part 7 Remark
\ref{Rem_Def_extend_wt_fn} imply $\kappa<1/2$.

Now suppose $\Pi$ is the 1-dimensional unit-valued rectangular function with
support $\left[  -.5,.5\right]  $. Then it was shown in Subsubsection
\ref{SbSbSect_int_Bspl_neq1_leq1} that:

\begin{theorem}
\label{Thm_data_fns_for_Xhat}If $u\in L^{2}\left(  \mathbb{R}^{1}\right)  $
then $u\ast\Pi$ is a data function i.e. $u\ast\Pi\in X_{w_{s}}^{0}$ where
$w_{s}$ is the hat weight function. Further, if we define $V\in L_{loc}%
^{1}\cap C_{BP}^{\left(  0\right)  }$ by%
\begin{equation}
V\left(  x\right)  =\int_{0}^{x}u\left(  t\right)  dt,\label{a7.74}%
\end{equation}

then $DV=u$ a.e.,%
\[
u\ast\Pi=\left(  2\pi\right)  ^{-\frac{1}{2}}\left(  V\left(  x+\frac{1}%
{2}\right)  -V\left(  x-\frac{1}{2}\right)  \right)  ,
\]

and%
\[
\left\Vert u\ast\Pi\right\Vert _{w_{s},0}=\left(  2\pi\right)  ^{-\frac{d}{4}%
}\left\Vert u\right\Vert _{2}.
\]

\end{theorem}

To obtain our data function $f_{d}=u\ast\Pi$ we will choose
\[
u\left(  x\right)  =e^{-x^{2}},
\]
for which
\begin{equation}
V=\left(  2\pi\right)  ^{\frac{1}{2}}\operatorname{erf},\text{\qquad
}\left\Vert u\right\Vert _{2}=2\left(  2\pi\right)  ^{\frac{1}{4}%
},\text{\qquad}\left\Vert u\ast\Pi\right\Vert _{w_{s},0}=2,\label{a7.914}%
\end{equation}

and%
\begin{equation}
f_{d}\left(  x\right)  =u\ast\Pi\left(  x\right)  =\operatorname{erf}\left(
x+\frac{1}{2}\right)  -\operatorname{erf}\left(  x-\frac{1}{2}\right)
.\label{a7.902}%
\end{equation}

Since we are using the hat basis function, from Subsubsection
\ref{SbSbSect_int_Bspl_neq1_leq1}
\[
G\left(  0\right)  =1,\quad C_{G}=1,\quad s=1/2,\quad h_{G}=\infty,\quad
k_{G}=\left(  2\pi\right)  ^{-1/4}\sqrt{2},
\]

and from \ref{a7.914} and \ref{a7.902},
\[
\left\Vert f_{d}\right\Vert _{w,0}=2.
\]

Finally, $Df_{d}\left(  x\right)  =\left(  2\pi\right)  ^{\frac{1}{2}%
}e^{-\left(  x+\frac{1}{2}\right)  ^{2}}-\left(  2\pi\right)  ^{\frac{1}{2}%
}e^{-\left(  x-\frac{1}{2}\right)  ^{2}}=0$ iff $x=0$ so the maximum value of
the data function on the data region $K=\left[  -1.5,1.5\right]  $ is
\[
\max\limits_{K}f_{d}=\operatorname{erf}\left(  \frac{1}{2}\right)
-\operatorname{erf}\left(  -\frac{1}{2}\right)  \simeq1.041.
\]

For the \textit{double rate} convergence experiment we will use the Riesz data
function $R_{0}=\left(  2\pi\right)  ^{-\frac{1}{2}}\Lambda$.

\subsubsection{\protect\underline{Numerical results}}

\fbox{Absolute smoother error vs. Smoothing parameter} We start by plotting
the absolute error bounds given by \ref{a8.10}, \ref{a7.26} and \ref{a7.27},
together with the actual absolute smoother errors, against the smoothing
parameter. To be precise the actual errors are averaged over a regular
\textit{error grid} covering\textbf{\ }the data region. The results for the
data function $f_{d}$ are shown in Figure
\ref{Fig_ApprSmthErrVsSmthParm_HatErf} and the results for the Riesz data
function $R_{0}$ are shown in Figure \ref{Fig_ApprSmthErrVsSmthParm_HatRiesz}:%

\begin{figure}[th]%
\centering
\includegraphics[
natheight=4.324100in,
natwidth=5.864300in,
height=4.3241in,
width=5.8643in
]%
{C:/Math_SwBasisFunc/InterpolSmthDev/PapersMonog/ZeroOrd/ZeroOrdDev/graphics/figApprSmthErrVsSmthParm_HatErf__13.pdf}%
\caption{Data function: $\operatorname{erf}\left(  x+\frac{1}{2}\right)
-\operatorname{erf}\left(  x-\frac{1}{2}\right)  $. Num smoother grid cells is
16.}%
\label{Fig_ApprSmthErrVsSmthParm_HatErf}%
\end{figure}
%

\begin{figure}[th]%
\centering
\includegraphics[
natheight=4.324100in,
natwidth=5.864300in,
height=4.3241in,
width=5.8643in
]%
{C:/Math_SwBasisFunc/InterpolSmthDev/PapersMonog/ZeroOrd/ZeroOrdDev/graphics/figApprSmthErrVsSmthParm_HatRiesz__14.pdf}%
\caption{Riesz data function: $R_{0}=\Lambda$ (Hat func). Num smoother grid
cells is 16.}%
\label{Fig_ApprSmthErrVsSmthParm_HatRiesz}%
\end{figure}

\fbox{Smoother error on Data region} As the number of data points increases
the smoother error of $f_{d}$ is observed to take the form of stable
alternating positive/negative spikes as shown in Figure
\ref{Fig_ApprSmth_ErrOnDataRegion_HatErf} which corresponds to 30K data
points. The positive errors are modulated by a hat function and the zero error
points have a period of 3/16; note there are 16 cells in the smoother grid.

It is seen in Figure \ref{Fig_ApprSmth_ErrOnDataRegion_HatRiesz} that as the
number of points increases the smoother error of the Riesz data function
$R_{0}$ is observed to stabilize to a wide negative central spike surrounded
by smaller spikes of about 1/5th the amplitude.

It is clear that, except for smoothing parameters larger than 10, the
theoretical error bounds given in \ref{a8.10} and \ref{a7.26} capture only a
small part of the actual error. The convergence of the error bounds to the
average spatial errors for large $\rho$ is typical. The averaging process also
disguises the instability indicated by Figure \ref{Fig_ApprConverg_N1_L1}.%

\begin{figure}[th]%
\centering
\includegraphics[
natheight=4.625000in,
natwidth=6.439400in,
height=4.625in,
width=6.4394in
]%
{C:/Math_SwBasisFunc/InterpolSmthDev/PapersMonog/ZeroOrd/ZeroOrdDev/graphics/figApprSmth_ErrOnDataRegion_HatErf__15.pdf}%
\caption{Data function: $\operatorname{erf}\left(  x+\frac{1}{2}\right)
-\operatorname{erf}\left(  x-\frac{1}{2}\right)  $. Number of smoother grid
cells is 16.}%
\label{Fig_ApprSmth_ErrOnDataRegion_HatErf}%
\end{figure}
%

\begin{figure}[th]%
\centering
\includegraphics[
natheight=4.605100in,
natwidth=6.412600in,
height=4.6051in,
width=6.4126in
]%
{C:/Math_SwBasisFunc/InterpolSmthDev/PapersMonog/ZeroOrd/ZeroOrdDev/graphics/figApprSmth_ErrOnDataRegion_HatRiesz__16.pdf}%
\caption{Riesz data function: $R_{0}=\Lambda$ (hat func). Number of smoother
grid cells is 16.}%
\label{Fig_ApprSmth_ErrOnDataRegion_HatRiesz}%
\end{figure}

\fbox{Absolute smoother error vs. Data density} Following the numerical work
of Chapter \ref{Ch_Exact_smth} which studied the Exact smoother we will also
apply the `spike' filter to the Approximate smoother errors. This filter
calculates the value below which 80\% of the errors lie. The filter is
designed to remove `large', isolated spikes which dominate the actual smoother
errors. The smoother error was calculated on a grid with 200 cells applied to
the domain of the data function. No filter is used for the first five
interpolants because there is no instability for small numbers of data points.
This filter will be used in all the cases below.

We now move on to examine the relationship between the smoother error and the
data density. To do this we need to fix a value for the smoothing coefficient
and here we choose $\rho=10^{-6}$. We select this value because often the
Exact smoother average error curve has a minimum near this value. Using the
functions and parameters discussed at the end of the last subsection we obtain
the four subplots displayed in Figure \ref{Fig_ApprConverg_N1_L1}, each
subplot being the superposition of 20 interpolants.%

\begin{figure}[th]%
\centering
\includegraphics[
natheight=5.840100in,
natwidth=3.646900in,
height=5.8401in,
width=3.6469in
]%
{C:/Math_SwBasisFunc/InterpolSmthDev/PapersMonog/ZeroOrd/ZeroOrdDev/graphics/figApprConverg_N1_L1_samp20_50Kpt_scal_half_16smcell__17.pdf}%
\caption{Error of Approximate smoother vs data density.}%
\label{Fig_ApprConverg_N1_L1}%
\end{figure}

The two upper subplots relate to the data function \ref{a7.902} and the lower
subplots relate to the data function $R_{0}$. The right-hand subplots are
filtered versions of the actual errors on the right. The (blue) points above
the actual errors in Figure \ref{Fig_ApprConverg_N1_L1} are the estimated
upper bounds for the error given by the theoretical estimates \ref{a7.56} and
\ref{a7.57}. The annotation at the bottom of the figure supplies the following
additional information:\medskip

\underline{\textbf{Input parameters}}\smallskip

$\mathbf{N=L=1}$ - the hat function is a member of the family of scaled,
extended splines with the indicated parameter values.

\textbf{spl scale 1/2} - the actual scaling of the spline basis function is
1/2 divided by \textbf{spl scale}.

\textbf{sm parm 1e-6} - the smoothing parameter is zero for interpolation.

\textbf{samp 20} - the sample size which is the number of test data files
generated. The data function is evaluated on the data interval $K=[-1.5,1.5]$
at points selected using a uniform (statistical) distribution.

\textbf{pts 2:50K} - specifies the smallest number of data points to be 2 and
the largest number of data points to be 30,000. The other values are given in
exponential steps with a multiplier of approximately 1.2.

\textbf{sm grid cell} \textbf{16} - number of cells in the regular smoother
grid superimposed on the data region $K$.\medskip

\underline{\textbf{Output parameters/messages}}\smallskip

\textbf{No ill-condit} - this indicates all Exact smoother matrices were
always properly conditioned.\textbf{\medskip}

Note that all the plots of error vs. data density given below have the same
format and annotations.

Clearly both theoretical error bounds substantially underestimates the
convergence rate. The instabilities associated with small numbers of data
points disappear as the number of data points increases. The smoother of
$R_{0}$ is less stable than that of $f_{d}$.\smallskip

\subsubsection{\protect\underline{The case $n=1$,$\ l=2$}}

Here we are interested in the 1-dimensional extended B-spline basis function
with $n=1$ and $l=2$. The following result, which is Theorem
\ref{Thm_data_func} in one dimension, will allow us to generate data functions
for which the $X_{w}^{0}$ norm can be calculated \ This result is closely
related to the calculations done above for the hat function.

\begin{theorem}
\label{Thm_data_func_3}\textbf{Data functions} Suppose $w$ is the B-spline
tensor product weight function with parameters $n$ and $l$ and univariate
\ref{7.65}, and that $U\in L^{2}\left(  \mathbb{R}^{1}\right)  $ with
$D^{n}U\in L^{2}\left(  \mathbb{R}^{1}\right)  $ in the sense of
distributions. Then if we define the distribution%
\[
f_{d}=\delta_{2}^{l}U,\quad l=1,2,3,\ldots,
\]

where $\delta_{2}$ is the central difference operator%
\begin{equation}
\delta_{2}U=U\left(  \cdot+1\right)  -U\left(  \cdot-1\right)  ,\label{a7.901}%
\end{equation}

it follows that $f_{d}\in X_{w}^{0}$ and%
\begin{equation}
\left\Vert f_{d}\right\Vert _{w,0}=2^{l}\left\Vert D^{n}U\right\Vert
_{2}.\label{a7.903}%
\end{equation}

\end{theorem}

Our tensor product basis function $G$ has univariate \ref{7.67} i.e.
\begin{align*}
G\left(  t\right)   & =\left(  -1\right)  ^{l-n}\tfrac{\left(  2\pi\right)
^{l/2}}{2^{2\left(  l-n\right)  +1}}\left(  D^{2\left(  l-n\right)  }\left(
\left(  \ast\Lambda\right)  ^{l}\right)  \right)  \left(  \tfrac{t}{2}\right)
\\
& =-\tfrac{\pi}{4}\left(  D^{2}\left(  \Lambda\ast\Lambda\right)  \right)
\left(  \tfrac{t}{2}\right)  ,\text{\quad}t\in\mathbb{R}^{1}.
\end{align*}
But%
\begin{align*}
D^{2}\left(  \Lambda\ast\Lambda\right)  =\Lambda\ast D^{2}\Lambda &
=\Lambda\ast\left(  \delta\left(  \cdot+1\right)  -2\delta+\delta\left(
\cdot-1\right)  \right) \\
&  =\frac{1}{\sqrt{2\pi}}\left(  \Lambda\left(  \cdot+1\right)  -2\Lambda
+\Lambda\left(  \cdot-1\right)  \right)  ,
\end{align*}

so that%
\[
G\left(  t\right)  =-\frac{\sqrt{2\pi}}{8}\left(  \Lambda\left(  \tfrac{t}%
{2}+1\right)  -2\Lambda\left(  \tfrac{t}{2}\right)  +\Lambda\left(  \tfrac
{t}{2}-1\right)  \right)  ,
\]

and hence%
\[
DG\left(  t\right)  =-\frac{\sqrt{2\pi}}{16}\left(  \Lambda^{\prime}\left(
\tfrac{t}{2}+1\right)  -2\Lambda^{\prime}\left(  \tfrac{t}{2}\right)
+\Lambda^{\prime}\left(  \tfrac{t}{2}-1\right)  \right)  ,
\]

i.e. $\left\Vert DG\right\Vert _{\infty}=\frac{3}{16}\sqrt{2\pi}$. Since the
(distributional) derivative is bounded the distributional Taylor series
expansion of Lemma \ref{Lem_Taylor_extension} can be used to write%
\[
G\left(  0\right)  -G\left(  t\right)  \leq\left\Vert DG\right\Vert _{\infty
}\left\vert t\right\vert ,\text{\quad}x\in\mathbb{R}^{1},
\]

which means that
\[
G\left(  0\right)  =\frac{\sqrt{2\pi}}{4},\text{ }C_{G}=\left\Vert
DG_{1,2}\right\Vert _{\infty}=\frac{3}{16}\sqrt{2\pi},\text{ }s=\frac{1}%
{2},\text{ }h_{G}=\infty.
\]

With reference to the last theorem we will choose the bell-shaped data
function%
\begin{equation}
f_{d}=\delta_{2}^{2}U\in X_{w}^{0},\label{a7.75}%
\end{equation}

where%
\[
U\left(  x\right)  =\frac{e^{-k_{1,2}x^{2}}}{\delta_{2}^{2}\left(
e^{-k_{1,2}x^{2}}\right)  \left(  0\right)  }=\frac{e^{-k_{1,2}x^{2}}%
}{2\left(  1-e^{-4k_{1,2}}\right)  },\quad k_{1,2}=0.3,
\]

so that $\left\Vert f_{d}\right\Vert _{\infty,K}=f_{d}\left(  0\right)  =1$
and
\[
\left\Vert f_{d}\right\Vert _{w,0}=4\left\Vert DU\right\Vert _{2}%
=\sqrt[4]{2\pi}\frac{\sqrt[4]{4k_{1,2}}}{1-e^{-4k_{1,2}}}.
\]

The error estimates are given by \ref{a7.56} and \ref{a7.57} where
$k_{G}=\left(  2\pi\right)  ^{-\frac{1}{4}}\sqrt{2C_{G}}$ and
\[
\left\Vert R_{0}\right\Vert _{\infty,K}=R_{0}\left(  0\right)  =\left(
2\pi\right)  ^{-\frac{1}{2}}G\left(  0\right)  =1/4.
\]

For the theory developed in the previous sections it was convenient to use the
unscaled B-spline weight function definition \ref{7.65} but for computations
it may be easier to use the unscaled version of the extended B-spline basis
function given in the next theorem.

\begin{theorem}
\label{Cor_Thm_data_func_3}Suppose $G_{s}=\left(  -1\right)  ^{l-n}D^{2\left(
l-n\right)  }\left(  \left(  \ast\Lambda\right)  ^{l}\right)  $ where
$\Lambda$ is the 1-dimensional hat function and $n,l$ are integers such that
$1\leq n\leq l$.

Then for given $\lambda>0$, $G_{s}\left(  \lambda x\right)  $ is called a
scaled extended B-spline basis function. The corresponding weight function is
$w_{\lambda}\left(  t\right)  =2\lambda aw\left(  \frac{t}{2\lambda}\right)  $
where $w$ is the extended B-spline weight function with parameters $n,l$ and
$a=\frac{\left(  2\pi\right)  ^{l/2}}{2^{2\left(  l-n\right)  +1}}$. Indeed,
for $w_{\lambda}$ we can choose $\kappa=n-1$. Further%
\begin{equation}
G_{s}\left(  0\right)  -G_{s}\left(  \lambda x\right)  \leq\lambda\left\Vert
DG_{s}\right\Vert _{\infty}\left\vert x\right\vert ,\text{\quad}x\in
\mathbb{R}^{1}.\label{a7.916}%
\end{equation}

Finally, if $f_{d}\in X_{w}^{0}$ and $g_{d}\left(  x\right)  =f_{d}\left(
2\lambda x\right)  $, it follows that $g_{d}\in X_{w_{\lambda}}^{0}$ and
$\left\Vert g_{d}\right\Vert _{w_{\lambda},0}=\sqrt{a}\left\Vert
f_{d}\right\Vert _{w,0}$.
\end{theorem}

\subsubsection{\protect\underline{Numerical results}}

\fbox{Smoother error on Data region} As the number of data points increases
the error function of the smoother of $f_{d}$ stabilizes: it has alternating
positive and negative spikes with an overall mean of zero. The zeros of the
error function have a period of 3/16: note there are 16 cells in the smoother grid.

As the number of data points increases the error of the smoother of the Riesz
data function $R_{0}$ is observed to stabilize: it has a wide negative central
spike at $x=0$ surrounded by much smaller values.\medskip

\fbox{Absolute smoother error vs. Data density} As with the previous case we
select $\rho=10^{-6}$ and, using the `spike filter' of the previous case, we
obtain the four subplots displayed in Figure \ref{Fig_ApprConverg_N1_L2}, each
subplot being the superposition of 20 smoothers. The two upper subplots relate
to the data function \ref{a7.75} and the lower subplots relate to the Riesz
data function $R_{0}$.%

\begin{figure}[th]%
\centering
\includegraphics[
natheight=6.431500in,
natwidth=3.957100in,
height=6.4315in,
width=3.9571in
]%
{C:/Math_SwBasisFunc/InterpolSmthDev/PapersMonog/ZeroOrd/ZeroOrdDev/graphics/figApprConverg_N1_L2_samp20_50Kpt_scal_1_16smcell__18.pdf}%
\caption{Error of Approximate smoother vs data density.}%
\label{Fig_ApprConverg_N1_L2}%
\end{figure}

The (blue) points above the actual errors in Figure
\ref{Fig_ApprConverg_N1_L2} are the estimated upper bounds for the error given
by the estimates \ref{a7.56} and \ref{a7.57}. Clearly both theoretical bounds
substantially underestimates the convergence rate. The instabilities
associated with small numbers of data points eventually disappear as the
number of data points increases. The smoother of $R_{0}$ is less stable than
that of $f_{d}$.

\subsection{Extended B-splines with $n=2$}

Since $n\geq2$ we can use the error estimates of Theorem
\ref{Thm_Final_Type2_estim_3}. We will use the same scaled (dilated) B-spline
and data function that we used for interpolation in
\ref{SbSect_int_data_fn_n_eq_2} and for Exact smoothing in Subsection
\ref{SbSect_ex_data_fn_n_eq_2}.

\subsubsection{\protect\underline{The case $n=2$, $l=2$}}

The basis function is
\[
G_{2,2}\left(  x\right)  =\frac{\left(  \Lambda\ast\Lambda\right)  \left(
2x\right)  }{\left(  \Lambda\ast\Lambda\right)  \left(  0\right)  }%
=\frac{3\sqrt{2\pi}}{2}\left(  \Lambda\ast\Lambda\right)  \left(  2x\right)  ,
\]
with scaling factor $\lambda=4$ and is such that $\operatorname*{supp}%
G_{2,2}=\left[  -1,1\right]  $ and $G_{2,2}\left(  0\right)  =1$. To calculate
$G_{2,2}$ we use the convenient formula%
\[
G_{2,2}\left(  x\right)  =\left(  1+x\right)  ^{2}\Lambda\left(  2x+1\right)
+\left(  1-2x^{2}\right)  \Lambda\left(  2x\right)  +\left(  1-x\right)
^{2}\Lambda\left(  2x-1\right)  ,
\]

and choose the bell-shaped data function%
\begin{equation}
f_{d}=\delta_{2}^{2}U\in X_{w}^{0},\label{a7.71}%
\end{equation}

where%
\begin{equation}
U\left(  x\right)  =\frac{e^{-k_{1,2}x^{2}}}{\delta_{2}^{2}\left(
e^{-k_{1,2}x^{2}}\right)  \left(  0\right)  }=\frac{e^{-k_{1,2}x^{2}}%
}{2\left(  1-e^{-4k_{1,2}}\right)  },\quad k_{1,2}=0.3,\label{a7.72}%
\end{equation}

so that%
\begin{equation}
k_{G}=\frac{\sqrt{12}}{\left(  2\pi\right)  ^{1/4}},\quad\left\Vert
f_{d}\right\Vert _{w,0}=\sqrt[4]{72\pi}\frac{\left(  k_{1,2}\right)  ^{3/4}%
}{1-e^{-4k_{1,2}}},\quad k_{1,2}=0.3,\quad s=1,\quad h_{G}=\infty
.\label{a7.73}%
\end{equation}

For the double rate convergence experiments the data function is
$R_{0}=\left(  2\pi\right)  ^{-\frac{1}{2}}G_{2,2}$ where $R_{0}\left(
0\right)  =\left(  2\pi\right)  ^{-\frac{1}{2}}$.

The next step is to calculate $\left\Vert f_{d}\right\Vert _{\infty,K}$ and
$\left\Vert R_{0}\right\Vert _{\infty,K}$ where $K=\left[  -1.5,1.5\right]  $.
Clearly $\left\Vert R_{0}\right\Vert _{\infty,K}=R_{0}\left(  0\right)  $ and
from \ref{a7.71} and \ref{a7.72}, $\left\Vert f_{d}\right\Vert _{\infty
,K}=f_{d}\left(  0\right)  =1$. Thus%
\begin{equation}
\left\Vert R_{0}\right\Vert _{\infty,K}=R_{0}\left(  0\right)  =\left(
2\pi\right)  ^{-\frac{1}{2}},\quad\left\Vert f_{d}\right\Vert _{\infty
,K}=1.\label{a7.741}%
\end{equation}

\subsubsection{\protect\underline{Numerical results}}

\fbox{Smoother error on Data region} As the number of data points increases
the smoother error of $f_{d}$ is observed to consist of stable alternating
positive and negative spikes with mean zero. The zeros of the error function
have a period of 3/16: note there are 16 cells in the smoother grid.

As the number of data points increases the smoother error of the Riesz data
function $R_{0}$ is observed to stabilize: it has a wide negative central
spike at $x=0$ surrounded by much smaller values.\medskip

\fbox{Absolute smoother error vs. Data density} As in the previous cases we
select $\rho=10^{-6}$ and use the spike filter. The output consists of the
four subplots displayed in Figure \ref{Fig_ApprConverg_N2_L2_16smcell}, each
subplot being the superposition of 20 smoothers. The bottom subplots
correspond to the Riesz data function.%

\begin{figure}[th]%
\centering
\includegraphics[
natheight=4.664800in,
natwidth=3.845000in,
height=4.6648in,
width=3.845in
]%
{C:/Math_SwBasisFunc/InterpolSmthDev/PapersMonog/ZeroOrd/ZeroOrdDev/graphics/figApprConverg_N2_L2_samp20_50Kpt_scal_quart_16smcell__19.pdf}%
\caption{Error of Approx smoother vs data density: 16 smth grid cells.}%
\label{Fig_ApprConverg_N2_L2_16smcell}%
\end{figure}

The (blue) points above the actual errors in Figure
\ref{Fig_ApprConverg_N2_L2_16smcell} are the theoretical upper bounds for the
smoother error given by the estimates \ref{a7.56} and \ref{a7.57}.

Clearly for the data function $f_{d}$ the theoretical error bound
substantially underestimates the convergence rate. Regarding the data function
$R_{0}$, the theoretical upper error bound for the data function $R_{0}$ is a
rather better fit to the unfiltered errors. The instabilities associated with
small numbers of data points eventually disappear as the number of data points increases.

\begin{figure}[th]%
\centering
\fbox{\includegraphics[
natheight=4.747600in,
natwidth=3.041200in,
height=4.7476in,
width=3.0412in
]%
{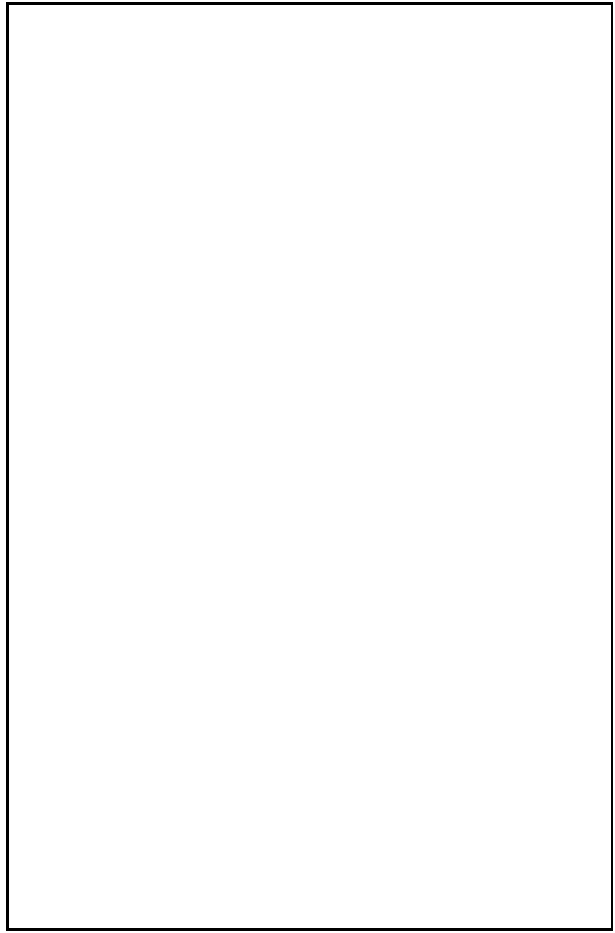}%
}\caption{Error of Approx smoother vs data density: 8 smth grid cells.}%
\label{Fig_ApprConverg_N2_L2_8smcell}%
\end{figure}

It is interesting to note that the case of 8 smoother grid cells shown in
Figure \ref{Fig_ApprConverg_N2_L2_8smcell} is more stable and converges to a
much smaller error value than the case of 16 smooth grid cells in Figure
\ref{Fig_ApprConverg_N2_L2_16smcell}.%

\section{A numerical implementation of the Approximate
smoother\label{Sect_numer_implem_Approx_smth}}

\subsection{The \textit{SmoothOperator} software (freeware)}

In this section we discuss a numerical implementation for the construction and
solution of the zero order Approximate smoother matrix equation \ref{a7.094}
i.e.%
\[
\left(  N\rho R_{X^{\prime},X^{\prime}}+R_{X,X^{\prime}}^{T}R_{X,X^{\prime}%
}\right)  \alpha^{\prime}=R_{X,X^{\prime}}^{T}y.
\]

I called this software \textit{SmoothOperator}. In Subsection
\ref{SbSect_Approx_smth_scalable} it was shown that the construction and
solution of this matrix equation is a scalable process and thus worthy of
numerical implementation. This software also implements the positive order
version of this equation which is derived in Chapter 6 of Williams
\cite{WilliamsPosOrdSmthV3} but the tutorials concentrate on the zero order
basis functions.

The algorithm has been implemented in \textbf{Matlab 6.0} with a GUI interface
but has only been tested on Windows. \textit{SmoothOperator} can be obtained
by emailing the author. However there is a \textbf{short user document }(4
pages) and the potential user can read this document first to decide whether
they want the software. The top-level directory of the software contains the
file \textit{read\_me.txt} which can also be downloaded separately. A
\textbf{full user manual} (82 pages) comes with the software. The main
features of the full user manual are:\medskip

\textbf{1} \textbf{Tutorials and data experiments }To learn about the system
and the behavior of the algorithm, I have prepared three tutorials and five
data experiments.\medskip

\textbf{2} \textbf{Context-sensitive help }Each dialog box incorporates
context-sensitive help which is invoked using the right mouse button. An F1
key facility can be implemented. The actual help text is contained in the text
file \newline\textbf{%
$\backslash$%
Help%
$\backslash$%
ContextHelpText.m} and this can be easily edited.\medskip

\textbf{3} \textbf{Matlab diary facility} When the system is started the
Matlab diary facility is invoked. This means that most user information
generated by \textit{SmoothOperator} and written to the Matlab command line is
also written to the text diary file. Each dialog box has a drop-down diary
menu which allows the user to view the diary file using Notepad. There is also
a facility for you to choose another editor/browser. The file can also be
emptied, if, for example, it gets too big. It can also be disabled.\medskip

\textbf{4} \textbf{Tools for viewing data files} Before generating a smoother
you can view the contents of the data file.

For ASCII delimited text data files, use the (slow) \textbf{View records}
facility to display records and the file header, and then with this
information you can use the high speed \textbf{Study records} facility to
check the records and then obtain detailed information about single fields
e.g. a histogram, and multiple fields e.g. correlation coefficients and
scatter plots.

For binary \textbf{test data files} only the \textbf{View records} facility is
needed. The parameters which generated the file can also be viewed using the
\textbf{Make or View data} option.\medskip

\textbf{5} \textbf{The output data} The output from the experiments and
tutorials mentioned in point 1 above consists of well-documented Matlab one
and two-dimensional plots, command line output and diary output. There is
currently no file output, but this can be implemented on request.\medskip

\textbf{6} \textbf{Reading delimited text files} A MEX C file allows ASCII
delimited text files to be read very quickly. You can specify:

\textbf{6.1} that the file be read in chunks of records, and not in one go.

\textbf{6.2} The ids of the fields to be read. This means that the file can
contain non-numeric fields.

\textbf{6.3} Fields can be checked to ensure they are numeric.\medskip

The tutorials which create smoothers allow the data, smoothed data, and
related functions to be viewed using scatter plots and plots along lines and
planes.\medskip

\textbf{7} Please note that there is no explicit suite of functions -
application programmer interface or API - supplied by \textit{SmoothOperator}
for immediate use by the user. After reflecting on how I would create such an
API, I decided that I lacked the experience to produce it, and that, anyhow,
there were too many possibilities to anticipate. I would expect that possible
users of this software, designed to be applied to perhaps millions of records,
would need to familiarize themselves thoroughly with how this system works. I
urge them to contact me and discuss their application.

\subsection{Algorithms}

The following three algorithms are used in the \textit{SmoothOperator}
software package to calculate the Approximate smoother.

\textbf{Algorithm 1} uses data generated internally according to
specifications supplied by the user using the interface. The smoothing
parameter can be either specified or calculated using an error grid. This
algorithm is \textbf{scalable}.

\textbf{Algorithm 2} uses a `small' subset of actual data to get an idea of a
suitable smoothing parameter to use for the full data set. This algorithm is
\textbf{not scalable}.

\textbf{Algorithm 3} uses all the actual data. The smoothing parameter can be
either specified or calculated using an \textit{error grid}. This algorithm is
\textbf{scalable}.

We will now explain these algorithms in more detail.

\subsubsection{\protect\underline{Algorithm 1: using experimental data}}

\begin{enumerate}
\item Generate the experimental data $\left[  X,y\right]  $. The independent
data\textbf{\ }$X$ is generated by uniformly distributed random numbers on $X
$. The dependent data $y$ is generated by uniformly perturbing an analytic
data function $g_{dat}$.

\item Choose a smoothing grid $X^{\prime}$ whose boundary contains $X$. Choose
an error grid $X_{err}^{\prime}$ which will be used to estimate the optimal
smoothing parameter $\rho$.

\item Read the data $\left[  X,y\right]  $ and construct the matrix equation.
If the matrix $G_{X,X^{\prime}}$ is dense, the matrix $G_{X,X^{\prime}}%
^{T}G_{X,X^{\prime}}$ is constructed using a Matlab MEX file (a compiled C
file). If $G_{X,X^{\prime}}$ is sparse the usual matrix multiplication is used.

\item Given a value for the smoothing parameter $\rho$ we can solve the matrix
equation and evaluate the smoother.

We want to estimate the value of $\rho$ which minimizes the `sum of squares'
error between the smoother and the data function
\[
\delta_{1}\left(  \rho\right)  =\sum\limits_{x^{\prime\prime}\in
X_{err}^{\prime}}\left(  \sigma_{\rho}\left(  x^{\prime\prime}\right)
-g_{dat}\left(  x^{\prime\prime}\right)  \right)  ^{2},\text{\quad}\rho>0.
\]

Empirical work indicates a standard shape for $\delta_{1}\left(  \rho\right)
$, namely decreasing from right to left, reaching a minimum and then
increasing at a decreasing rate. To find the minimum we basically use the
standard iterative algorithm of dividing and multiplying by a factor e.g. 10
and choosing the smallest value. The process is stopped when the percentage
change of one or both of $\delta_{1}\left(  \rho\right)  $ and $\rho$ are less
than prescribed values.
\end{enumerate}

\subsubsection{\protect\underline{Algorithm 2: using a `test' subset of actual
data}}

\begin{enumerate}
\item Perhaps based on results using Algorithm 1, choose an initial value for
smoothing\textbf{\ }parameter $\rho$ and choose a smoothing grid $X^{\prime}$.

\item Step 2 of Algorithm 1. Denote the data by $\left[  X,y\right]  $ where
$X=\left(  x^{\left(  i\right)  }\right)  $ and $y=\left(  y_{i}\right)  $.

\item This is the same as step 4 of Algorithm 2 except we now minimize
\[
\delta_{2}\left(  \rho\right)  =\sum\limits_{i=1}^{N}\left(  \sigma_{\rho
}\left(  x^{\left(  i\right)  }\right)  -y_{i}\right)  ^{2},\text{\quad}%
\rho>0,
\]

because we do not have a data function.
\end{enumerate}

\subsubsection{\protect\underline{Algorithm 3: using all the actual data}}

\begin{enumerate}
\item Choose a value for smoothing parameter $\rho$, based on experiments
using Algorithm 2. Choose a smoothing grid $X^{\prime}$.

\item Read all the data $\left[  X,y\right]  $ and construct the matrix
equation. If the matrix $G_{X,X^{\prime}}$ is dense, the matrix
$G_{X,X^{\prime}}^{T}G_{X,X^{\prime}}$ is constructed using a Matlab MEX file
(a compiled C file). If $G_{X,X^{\prime}}$ is sparse the usual multiplication
is used.

\item Solve the matrix equation to obtain the basis function coefficients and
evaluate the smoother at the desired points.
\end{enumerate}

\subsection{Features of the smoothing algorithm and its
implementation\label{SbSect_AlgorImplem}}

\textbf{1} The Short user manual and the User manual contains a lot of detail
regarding the \textit{SmoothOperator} system and algorithms. So we will
content ourselves here with just some key points.\medskip

\textbf{2} Although the algorithm is scalable there can still be a problem
with rapidly increasing memory usage as the grid size decreases and the
dimension increases. The classical radial basis functions, such as the thin
plate spline functions, have support everywhere. Hence the smoothing matrix is
completely full. To significantly reduce this problem we do the following :

\textbf{a)} We use basis functions with bounded support.

\textbf{b)} We shrink the basis function support to the magnitude of the grid
cells. This makes $G_{X,X^{\prime}}^{T}G_{X,X^{\prime}}$ a \textbf{very
sparse} banded diagonal matrix, and $G_{X^{\prime},X^{\prime}}$ has a small
number of non-zero diagonals. We say this basis function has \textit{small
support}.\medskip

\textbf{3} Instead of using the memory devouring Matlab \texttt{repmat}
function to calculate $G_{X,X^{\prime}}$ and $G_{X^{\prime},X^{\prime}}$, we
directly calculate the arguments of the Matlab function \texttt{sparse}. This
function takes three arrays, namely the row ids, the column ids and the
corresponding matrix elements, as well as the matrix dimensions, and converts
them to the Matlab sparse internal representation. I must mention that
although this is quick and space efficient, the algorithm is much more
complicated than using \texttt{repmat}.

Matlab's sparse multiplication facility is then used to quickly and
efficiently calculate $G_{X,X^{\prime}}^{T}G_{X,X^{\prime}}$.\medskip

\textbf{4} This software has implemented the above techniques for the
smoothing matrix, and a selection of basis functions is supplied.\medskip

\textbf{5} The tutorials and exercises are based around the zero order
\textbf{tensor product hat (triangle) function, }denoted by $\Lambda$. The hat
function, also known as the triangle function, is used because of its
simplicity and its analytic properties. The higher dimensional hat functions
are defined as the tensor product of the one dimensional hat function :
$\Lambda\left(  x_{1}\right)  =1-\left\vert x_{1}\right\vert $, when
$\left\vert x_{1}\right\vert \leq1$, and zero otherwise. This function has
zero order so that the smoothing matrix simplifies to%
\begin{equation}
\Psi_{G}=\rho NG_{X^{\prime},X^{\prime}}+G_{X,X^{\prime}}^{T}G_{X,X^{\prime}%
},\text{\quad}G_{X^{\prime},X^{\prime}}=I,\label{2.39}%
\end{equation}

with matrix equation%
\[
\Psi_{G}\alpha=G_{X^{\prime},X}y.
\]

Smoothers constructed from the hat function are continuous but not smooth. I
have also included the \textbf{B-spline of order 3}, namely $\Lambda
\ast\Lambda$, which is a basis function of order zero or one. This can be used
when a smoother is required to be smooth i.e. continuously
differentiable.\medskip

\textbf{6} \textbf{Note that this software was written before I embarked on
the error analysis contained in this document} so \textit{SmoothOperator}
concentrates on the Approximate smoother and the user cannot study the error
of the Exact smoother or the Approximate smoother or compare the Approximate
smoother with the Exact smoother.

\subsection{An application - predictive modelling of forest cover type}

In this section we demonstrate how the method developed in the previous
sections can be used in data mining for predictive modelling. This application
smooths \textbf{binary-valued} data - I was unable to obtain a large
`continuous-valued' data set for distribution on the web. The source of our
data is the web site file:%
\[%
\begin{tabular}
[c]{l}%
\textit{http://kdd.ics.uci.edu/databases/covertype/covertype.html,}%
\end{tabular}
\]

in the UCI KDD Archive, Information and Computer Science, University of
California, Urvine.

The data gives the forest cover type in $30\times30$ meter cells as a function
of the following cartographic parameters:%

\begin{table}[htbp] \centering
$%
\begin{tabular}
[c]{|l|l|l|}\hline
Id & Independent variable & Description\\\hline\hline
1 & ELEVATION & Altitude above sea level\\
2 & ASPECT & Azimuth\\
3 & SLOPE & Inclination\\
4 & HORIZ\_HYDRO & Horizontal distance to water\\
5 & VERT\_HYDRO & Vertical distance to water\\
6 & HORIZ\_ROAD & Horizontal distance to roadways\\
7 & HILL\_SHADE\_9 & Hill shade at 9am\\
8 & HILL\_SHADE\_12 & Hill shade at noon\\
9 & HILL\_SHADE\_15 & Hill shade at 3pm\\
10 & HORIZ\_FIRE & Horizontal distance to fire points\\\hline
\end{tabular}
$\caption{}\label{Tbl_cart_parms}%
\end{table}%

Forest cover type is the dependent variable\ $y$ and it takes on one of seven values:%

\begin{table}[htbp] \centering
$%
\begin{tabular}
[c]{|l|l|}\hline
Forest cover type & Id\\\hline\hline
Spruce fir & 1\\
Lodge-pole pine & 2\\
Ponderosa pine & 3\\
Cottonwood/Willow & 4\\
Aspen & 5\\
Douglas fir & 6\\
Krummholtz & 7\\\hline
\end{tabular}
$\caption{}\label{Tbl_forest_cover}%
\end{table}%

\subsubsection{\protect\underline{The data file}}

In this study we will use the data to train a model predicting on the
\textbf{presence or absence of the Ponderosa pine forest cover (}id = 3), but
the results will be similar for the other forest types. To this end we have
created from the full web site file a file called \texttt{%
$\backslash$%
UserData%
$\backslash$%
forest\_1\_to\_10\_pondpin.dat} which contains the ten independent variables
of Table \ref{Tbl_cart_parms} and then a binary dependent variable derived
from the variable \textbf{Ponderosa pine} of Table \ref{Tbl_forest_cover}.
This variable is 1 if the cover is Ponderosa pine and zero otherwise.

\subsubsection{\protect\underline{Methodology}}

We recommend you first use the user interface of the \textit{SmoothOperator
}software\textit{\ }to construct artificial data sets to understand the
behavior of the smoother and run the experiments to get a feel for the
influence of the parameters.\medskip

\textbf{1} We chose a small\ subset of the forest-cover data and selected
various smoothing parameters and smoothing grid sizes to study the performance
of the smoother using plots and the value of the error.

Note that \textit{two subsets of data} could be used here, often called the
\textit{training set} and the \textit{test set}. The training set would be
used to calculate the smoother for the initial value of the smoothing
coefficient, and then the test set could be used to determine the smoothing
coefficient which minimizes the least squares\ error.\medskip

\textbf{2} Having chosen our parameters we run the smoother program on the
full data set.

\chapter{The spaces $X_{1/w}^{0}$, $\protect\widetilde{X}_{1/w}^{0}$ and the
bounded linear functionals on $X_{w}^{0}$\label{Ch_Xo1/w_and_fnal_on_Xow}}

\section{Introduction}

??? \textbf{UPDATE INTRO}! ??

The goal of this chapter is to use a bilinear form to characterize the bounded
linear functionals $\left(  X_{w}^{0}\right)  ^{\prime}$ on the data set
$X_{w}^{0}$ for the weight functions used as examples in this document and to
explore the properties of all related operators. Now section by
section:\smallskip

\textbf{Section \ref{Sect_Xo1/w}} If $w$ is a weight function w.r.t. the set
$\mathcal{A}$ this characterization is first done in Section \ref{Sect_Xo1/w}
by applying certain conditions to $w$ and defining the semi-Hilbert space%
\[
X_{1/w}^{0}=\left\{  v\in S^{\prime}:v_{F}\in L_{loc}^{1}\left(
\mathbb{R}^{d}\setminus\mathcal{A}\right)  \text{ }and\text{ }\int%
\frac{\left\vert v_{F}\right\vert ^{2}}{w}<\infty\right\}  ,
\]

where $v_{F}$ is the restriction of $\widehat{v}$ to $\mathbb{R}^{d}%
\setminus\mathcal{A}$ and then constructing an isometric isomorphism from
$L^{2}$ to $X_{1/w}^{0}$, as indicated by:%
\begin{gather}%
\begin{array}
[c]{ccccc}%
X_{1/w}^{0} &
\begin{array}
[c]{c}%
\overset{\mathcal{L}}{\longleftarrow}\\
\underset{\mathcal{M}}{\longrightarrow}%
\end{array}
& L^{2} &
\begin{array}
[c]{c}%
\overset{\mathcal{I}}{\longleftarrow}\\
\underset{\mathcal{J}}{\longrightarrow}%
\end{array}
& X_{w}^{0}%
\end{array}
\label{Fig_X_maps_intro}\\
\text{Isometric isomorphisms between }X_{1/w}^{0}\text{, }L^{2}\text{, }%
X_{w}^{0}\text{.}\nonumber\\
\text{Also, }\mathcal{V}=\mathcal{JM}\text{ and }\mathcal{W}=\mathcal{LI}%
\text{.}\nonumber
\end{gather}

This case\ will include all the weight functions introduced in Section
\ref{Sect_gen_wt_funcs} except the Gaussian and the thin-plate splines.

Now the mapping $\mathcal{M}$ is easy to construct and the work lies in the
mapping $\mathcal{L}$. We first define $\mathcal{L}$ under the conditions
\ref{a1.046}, namely%
\begin{equation}
w\left\vert \cdot\right\vert ^{2m}\in L_{loc}^{1}\left(  \mathbb{R}%
^{d}\right)  ,\quad\int\limits_{\left\vert \cdot\right\vert \geq r}\frac
{w}{\left\vert \cdot\right\vert ^{2\tau}}<\infty,\label{a6.3}%
\end{equation}

which (Theorem \ref{Thm_wt_fn_satisfy_property_1}) includes the
\textit{Sobolev splines} and the univariate \textit{central difference} weight
functions.\textit{\ }We then introduce the subspace
\[
S_{w,0}=\left\{  \phi\in S:\int w\left\vert \phi\right\vert ^{2}%
<\infty\right\}  ,
\]

and we say that $w\in W_{S;0}$ if%
\begin{equation}
\left(  \int w\left\vert \phi\right\vert ^{2}\right)  ^{1/2}\leq\left\vert
\phi\right\vert _{\sigma},\quad\phi\in S_{w,0},\label{a6.2}%
\end{equation}

where $\left\vert \cdot\right\vert _{\sigma}$ is a positive, linear
combination of the seminorms which define the topology of $S$. We define
$\mathcal{L}$ for $w\in W_{S;0}$, which includes the \textit{extended
B-splines} and the \textit{central difference weight} functions.

In Corollary \ref{Cor_X1/w,o_semiHilb_opL2} we show that if $w\in W01\cap
W_{S;0}$ then in general $X_{1/w}^{0}$ is a semi-Hilbert space.

The bounded linear functionals $\left(  X_{w}^{0}\right)  ^{\prime}$ will be
characterized by the isometric, isomorphism $\Phi:X_{1/w}^{0}\rightarrow
\left(  X_{w}^{0}\right)  ^{\prime}$ defined by
\[
\left(  \Phi v\right)  \left(  u\right)  =\left(  u,\mathcal{V}v\right)
_{w,0},\text{\quad}u\in X_{w}^{0},\text{ }v\in X_{1/w}^{0},
\]

where $\mathcal{V}=\mathcal{JM}$. In an analogous fashion to the negative
order Sobolev spaces, this can be written in terms of a bilinear form as%
\[
\left(  \Phi v\right)  \left(  u\right)  =\int\widehat{u}\overline{v_{F}%
},\text{\quad}u\in X_{w}^{0},\text{ }v\in X_{1/w}^{0}.
\]

Later on, in Section \ref{Sect_generalize2_tildXo1/w_T}, the space
$X_{1/w}^{0}$ will be generalized to include the Gaussian and the thin-plate
splines.\smallskip

\textbf{Section \ref{Sect_tildXo1/w_S2dag}} Under the assumption that $w\in
W_{S;0}$, a larger space $\widetilde{X}_{1/w}^{0}$ is defined which allows
similar, but simpler, maps to be defined to and from $L^{2}$. We then set
about defining the isometric maps analogous to those defined for $X_{1/w}^{0}
$. Some of these maps are indicated by:%

\begin{equation}%
\begin{array}
[c]{ccccccc}%
X_{1/w}^{0} &
\begin{array}
[c]{c}%
\overset{E}{\longleftarrow}\\
\underset{R}{\longrightarrow}%
\end{array}
& \widetilde{X}_{1/w}^{0} &
\begin{array}
[c]{c}%
\overset{\widetilde{\mathcal{L}}}{\longleftarrow}\\
\underset{\widetilde{\mathcal{M}}_{1}}{\longrightarrow}%
\end{array}
& L^{2} &
\begin{array}
[c]{c}%
\overset{\mathcal{I}}{\longleftarrow}\\
\underset{\mathcal{J}}{\longrightarrow}%
\end{array}
& X_{w}^{0}%
\end{array}
\label{Fig_tildX_maps_intro}%
\end{equation}

The space $X_{w}^{0}$ will be defined as%
\begin{equation}
\widetilde{X}_{1/w}^{0}=\left\{  u\in\left(  \widehat{S}_{w,0}\right)
^{\prime}:u_{F}\in L_{loc}^{1}\left(  \mathbb{R}^{d}\setminus\mathcal{A}%
\right)  \text{ }and\text{ }\int\frac{\left\vert u_{F}\right\vert ^{2}}%
{w}<\infty\right\}  ,\label{a6.9}%
\end{equation}

where $u_{F}$ is $\widehat{u}$ restricted to $\mathbb{R}^{d}\setminus
\mathcal{A}$. $\left(  \widehat{S}_{w,0}\right)  ^{\prime}$ denotes the
continuous linear functionals on the space $\widehat{S}_{w,0}$ which consists
of the Fourier transforms of functions in $S_{w,0}$, and $\widehat{u}\in
S_{w,0}^{\prime}$.

The relationship between the spaces $\widetilde{X}_{1/w}^{0}$ and $X_{w}^{0}$
is studied and we construct an isometric, isomorphism between $\widetilde{X}%
_{1/w}^{0}$ and $X_{w}^{0}$ using the restriction mapping $R$ and extension
mapping $E$. Next we use Fourier transforms to map $\widetilde{X}_{1/w}^{0}$
and $X_{w}^{0}$ onto Fourier-independent spaces so that the equivalents of the
mappings $\widetilde{\mathcal{L}}$, $\widetilde{\mathcal{M}}_{1}$,
$\mathcal{I}$ and $\mathcal{J}$ become Fourier-independent (see Figure
\ref{Fig_ops_Y_tildY_Fourier_indep}).

Finally, it is shown that the Gaussian and the shifted thin-plate spline
weight functions do not lie in $W_{S;0}$.\smallskip

\textbf{Section \ref{Sect_generalize1_tildXo1/w_T}} We define $X_{1/w}^{0}$
when $w$ satisfies \ref{a6.3} and we define $\widetilde{X}_{1/w}^{0}$ for
$w\in S_{w,0}^{\dag}$. In Condition \ref{Cnd_wt_fn_T} we show how to
generalize our definition of $\widetilde{X}_{1/w}^{0}$ in order to include the
case when $w$ satisfies \ref{a6.3} whilst retaining the form of the integral
in \ref{a6.2} used to define $S_{w,0}^{\dag}$. In fact, we will postulate the
existence of a subspace $T\subset S_{w,0}$ such that \ref{a6.2} holds for
$\phi\in T$, and we also will assume that $C_{0}^{\infty}\left(
\mathbb{R}^{d}\setminus\mathcal{A}\right)  \subset T$.

Given these considerations our approach will be to replace the space $S_{w,0}
$\ by the space $T$ in the results and definitions of Section
\ref{Sect_tildXo1/w_S2dag} and show that they all still hold.

Finally, we will show that the Gaussian and the shifted thin-plate spline
weight functions do not satisfy Condition \ref{Cnd_wt_fn_T}.\smallskip

\textbf{Section \ref{Sect_generalize2_tildXo1/w_T}} We now consider the
problem of characterizing the bounded linear functionals on $X_{w}^{0}$ for
the weight functions which do not have property \ref{a1.046}, do not belong to
$W_{S;0}$ and do not satisfy Condition \ref{Cnd_wt_fn_T}; which means the
results of Sections \ref{Sect_Xo1/w}, \ref{Sect_tildXo1/w_S2dag} and
\ref{Sect_generalize1_tildXo1/w_T} are not available.

To handle these weight functions we will adapt definition \ref{a6.9} of
$\widetilde{X}_{1/w}^{0}$. However, \textbf{in this section }$S_{w,0}%
$\textbf{\ will be considered not as a subspace of }$S$\textbf{\ but as a
topological vector space in itself} and will be endowed with the topology
using the $S$ seminorms \textbf{and the norm }$\int w\left\vert \phi
\right\vert ^{2}$. We will still call the data space $\widetilde{X}_{1/w}^{0}%
$. We will prove analogues of the results of Section
\ref{Sect_tildXo1/w_S2dag}.

The space $\widehat{S}_{w,0}$ will be endowed with the topology that makes the
inverse-Fourier transform to $S_{w,0}$ a homeomorphism. Following the approach
used in Section \ref{Sect_tildXo1/w_S2dag} the mappings
$\widetilde{\mathcal{L}}_{4}:L^{2}\rightarrow\left(  \widehat{S}_{w,0}\right)
^{\prime}$ and $\widetilde{\mathcal{M}}_{4}:\widetilde{X}_{1/w}^{0}\rightarrow
L^{2}$ will be introduced, followed by the operators $\widetilde{\mathcal{V}%
}_{3}:\widetilde{X}_{1/w}^{0}\rightarrow X_{w}^{0}$, $\widetilde{\mathcal{B}%
}_{3}:X_{w}^{0}\otimes\widetilde{X}_{1/w}^{0}\rightarrow C_{B}^{\left(
0\right)  }$ and $\widetilde{\Phi}_{3}:\widetilde{X}_{1/w}^{0}\rightarrow
\left(  X_{w}^{0}\right)  ^{\prime}$, the last operator being used to
characterize the bounded linear functionals on $X_{w}^{0}$ as members of
$\widetilde{X}_{1/w}^{0}$ w.r.t. to the bilinear form $\int\widehat{u}%
\overline{v_{F}}$ where $u\in X_{w}^{0}$ and $v\in\widetilde{X}_{1/w}^{0}$.

Although we do not do so here, analogues of the extension and restriction
mappings between $\widetilde{X}_{1/w}^{0}$ and $X_{1/w}^{0}$ of Subsection
\ref{SbSect_E_R_tildXo1/w_Xo1/w} can be easily proved as can analogues of the
Fourier-independent spaces $\widetilde{Y}_{1/w}^{0}$ and $Y_{w}^{0}$ of
Subsection \ref{SbSect_tildYo1/w_Yow_FTindep}.

\section{The spaces $X_{1/w}^{0}$ and $\overline{X}_{1/w}^{0}$%
\label{Sect_Xo1/w}}

\begin{definition}
\label{Def_Xo1/w}\textbf{The semi-inner product space }$X_{1/w}^{0}$

Suppose that the weight function $w$ has property W01 for the set
$\mathcal{A}$, where $\mathcal{A}$ is a closed set of measure zero. Then we
define the semi-inner product space:%
\begin{equation}
X_{1/w}^{0}=\left\{  u\in S^{\prime}:u_{F}\in L_{loc}^{1}\left(
\mathbb{R}^{d}\setminus\mathcal{A}\right)  \text{ }and\text{ }\int%
\frac{\left\vert u_{F}\right\vert ^{2}}{w}<\infty\right\}  ,\label{a3.3}%
\end{equation}

where $u_{F}$ is the restriction of $\widehat{u}$, regarded as a member of
$\mathcal{D}^{\prime}$, to $\mathbb{R}^{d}\setminus\mathcal{A}$. We endow
$X_{1/w}^{0}$ with seminorm and semi-inner product%
\[
\left\vert u\right\vert _{1/w,0}=\int\frac{\left\vert u_{F}\right\vert ^{2}%
}{w},\quad\left\langle u,v\right\rangle _{1/w,0}=\left(  \int\frac
{u_{F}\overline{v_{F}}}{w}\right)  ^{1/2}.
\]

The notation $u_{F}$ is used extensively for the positive order case in
Williams \cite{WilliamsPosOrdSmthV3}.

In Lemma \ref{Lem_Xoinvw} below we show that $X_{1/w}^{0}$ is not empty, that
the seminorm has null space $\left(  S_{\mathcal{A}}^{\prime}\right)  ^{\vee}
$ where
\begin{equation}
S_{\mathcal{A}}^{\prime}=\left\{  u\in S^{\prime}:\operatorname*{supp}%
u\subseteq\mathcal{A}\text{ }as\text{ }a\text{ }distribution\right\}
,\label{a1.4}%
\end{equation}

and that we can choose $\mathcal{A}$ to be empty iff $X_{1/w}^{0}$ is an inner
product space. Note that in \ref{a1.4} we regard $u$ as a distribution because
tempered distributions are not localizable on $S$.
\end{definition}

\begin{remark}
\begin{enumerate}
\item \textbf{1}. Why not choose the space%
\[
X_{1/w}^{0}=\left\{  u\in S^{\prime}:\widehat{u}\in L_{loc}^{1}\text{
}and\text{ }\int\frac{\left\vert \widehat{u}\right\vert ^{2}}{w}%
<\infty\right\}  ?
\]

See Definition \ref{Def_barXow} below.

??? I looked into this but concluded that I needed \ref{a3.3} but at the
moment I can't remember where.

\item Perhaps we could restrict $\mathcal{A}$ to having the form
$\bigcup_{j=0}^{d-1}\bigcup_{i=1}^{n_{j}}\mathcal{A}_{i,j}$ where $n_{j}%
\leq\infty$ and $\mathcal{A}_{i,j}$ is locally homeomorphic to the ball
$B_{1}\left(  \mathbb{R}^{j}\right)  $ when $j\geq1$ and a single point when
$j=0$.
\end{enumerate}
\end{remark}

\begin{lemma}
\label{Lem_Xoinvw}Suppose the weight function $w$ has property W01 for the set
$\mathcal{A}$. Then:

\begin{enumerate}
\item $X_{1/w}^{0}$ is not empty.

\item The seminorm $\left\vert \cdot\right\vert _{1/w,0}$ has null space
$\left(  S_{\mathcal{A}}^{\prime}\right)  ^{\vee}$, defined by \ref{a1.4}.
Further, noting part 1 of Remark \ref{Rem_Def_extend_wt_fn}, we see that
$X_{1/w}^{0}$ is an inner product space iff we can choose $\mathcal{A}$ to be
the (minimal) empty set.\medskip

Suppose $w$ also has property W02 for parameter $\kappa$. Then:\medskip

\item If $\left\vert v\left(  \xi\right)  \right\vert \leq c\left(
1+\left\vert \xi\right\vert \right)  ^{\kappa}$ then $\overset{\vee}{v}\in
X_{1/w}^{0}$. If $w\in W02$ then $D^{\alpha}\delta\in X_{1/w}^{0}$ when
$\left\vert \alpha\right\vert \leq\kappa$.

\item $S\subset X_{1/w}^{0}$.

\item If $u\in X_{1/w}^{0}$ then $\frac{\left(  D^{\alpha}u\right)  _{F}}%
{w}\in L^{1}$ when $\left\vert \alpha\right\vert \leq\kappa$.

Also, $\int\frac{\left(  D^{\alpha}u\right)  _{F}}{w}\in\left(  X_{1/w}%
^{0}\right)  ^{\prime}$ with $\left\Vert \int\frac{\left(  D^{\alpha}u\right)
_{F}}{w}\right\Vert _{op}\leq\left(  \int\frac{\xi^{2\alpha}}{w}\right)
^{1/2}$.

\item However, suppose instead that $w$ has property W03 for parameter
$\kappa$. Then parts 3 to 5 hold when $\alpha\leq\kappa$.
\end{enumerate}
\end{lemma}

\begin{proof}
\textbf{Part 1} Choose $\phi\in C_{0}^{\left(  0\right)  }$ such that
$\operatorname*{supp}\phi\subset\mathbb{R}^{d}\setminus\mathcal{A}$. Then
$\sqrt{w}\phi\in C_{0}^{\left(  0\right)  }$ with $\operatorname*{supp}%
\sqrt{w}\phi\subset\mathbb{R}^{d}\setminus\mathcal{A}$. Set $u=\left(
\sqrt{w}\phi\right)  ^{\vee}$ so that $u\in C_{BP}^{\infty}\subset S^{\prime}%
$. Then $u_{F}=\sqrt{w}\phi\in L_{loc}^{1}\left(  \mathbb{R}^{d}%
\setminus\mathcal{A}\right)  $ and $\int\frac{\left\vert u_{F}\right\vert
^{2}}{w}=\int\left\vert \phi\right\vert ^{2}<\infty$.\medskip

\textbf{Part 2}
\begin{align}
\operatorname{null}\left\vert \cdot\right\vert _{1/w,0} &  =\left\{  u\in
X_{1/w}^{0}:\left\vert u\right\vert _{1/w,0}=0\right\} \nonumber\\
&  =\left\{  u\in S^{\prime}:\widehat{u}=0\text{ }on\text{ }\mathbb{R}%
^{d}\setminus\mathcal{A}\text{ }as\text{ }a\text{ }distribution\right\}
\nonumber\\
&  =\left\{  \widehat{u}:u\in S^{\prime},\text{ }\widehat{u}=0\text{ }on\text{
}\mathbb{R}^{d}\setminus\mathcal{A}\right\}  ^{\vee}\nonumber\\
&  =\left\{  v\in S^{\prime}:v=0\text{ }on\text{ }\mathbb{R}^{d}%
\setminus\mathcal{A}\text{ }as\text{ }a\text{ }distribution\right\}  ^{\vee
}\label{a5.9}\\
&  =\left\{  v\in S^{\prime}:\operatorname*{supp}u\subset\mathcal{A}\right\}
^{\vee}\nonumber\\
&  =\left(  S_{\mathcal{A}}^{\prime}\right)  ^{\vee}.\nonumber
\end{align}

Now if $\mathcal{A}$ is empty then \ref{a5.9} implies that
$\operatorname{null}\left\vert \cdot\right\vert _{1/w,0}=\left\{  0\right\}
$. On the other hand if $\operatorname{null}\left\vert \cdot\right\vert
_{1/w,0}=\left\{  0\right\}  $ then from \ref{a5.9} it follows that $u\in
S^{\prime}$ and $u=0$ on $\mathbb{R}^{d}\setminus\mathcal{A}$ implies $u=0$ on
$\mathbb{R}^{d}$. Now if $\mathcal{A}$ is not empty and $a\in\mathcal{A}$ then
$\delta\left(  \cdot-a\right)  \in S^{\prime}$ so $\delta\left(
\cdot-a\right)  \in S^{\prime}$ and $\delta\left(  \cdot-a\right)  =0$ on
$\mathbb{R}^{d}\setminus\mathcal{A}$. But $\delta\left(  \cdot-a\right)
\neq0$ on $\mathbb{R}^{d}$ which is a contradiction. Thus $\mathcal{A}$ must
be empty.\medskip

\textbf{Part 3} Let $f=\overset{\vee}{v}$. Then $f\in S^{\prime}$,
$\widehat{f}\in L_{loc}^{1}$ so that $f_{F}\in L_{loc}^{1}\left(
\mathbb{R}^{d}\setminus\mathcal{A}\right)  $, and $\int\frac{\left\vert
f_{F}\right\vert ^{2}}{w}\leq\int\frac{c\left(  1+\left\vert \cdot\right\vert
\right)  ^{\kappa}}{w}<\infty$.\medskip

\textbf{Part 4} If $u\in S$ then $v=\widehat{u}\in S$ is bounded and by part
4, $\overset{\vee}{v}=u\in X_{1/w}^{0}$.\medskip

\textbf{Part 5} $\int\frac{\left\vert \left(  D^{\alpha}u\right)
_{F}\right\vert }{w}=\int\frac{\left\vert \xi^{\alpha}u_{F}\right\vert }%
{w}=\int\frac{\left\vert \xi^{\alpha}\right\vert }{\sqrt{w}}\frac{\left\vert
u_{F}\right\vert }{\sqrt{w}}\leq\left(  \int\frac{\xi^{2\alpha}}{w}\right)
^{1/2}\left(  \int\frac{\left\vert u_{F}\right\vert ^{2}}{w}\right)
^{1/2}=\left(  \int\frac{\left\vert \cdot\right\vert ^{2\left\vert
\alpha\right\vert }}{w}\right)  ^{1/2}\left\vert u\right\vert _{1/w,0}<\infty
$.\medskip

\textbf{Part 6} ?? \textbf{FINISH}!
\end{proof}

\subsection{The operator $\mathcal{M}:X_{1/w}^{0}\rightarrow L^{2}%
$\label{SbSect_op_M}}

In this section we define the mapping $\mathcal{M}:X_{1/w}^{0}\rightarrow
L^{2}$ displayed in Figure \ref{Fig_X_maps_intro} which we hope will have nice
properties in the seminorm sense.

\begin{definition}
\label{Def_op_M}\textbf{The linear operator} $\mathcal{M}:X_{1/w}%
^{0}\rightarrow L^{2}$\textbf{:} Suppose $w$ is a weight function w.r.t. the
set $\mathcal{A}$. From the definition of $X_{1/w}^{0}$, $u\in X_{1/w}^{0}$
implies $\frac{u_{F}}{\sqrt{w}}\in L^{2}$ where $u_{F}$ is the restriction of
$\widehat{u}$ to $\mathbb{R}^{d}\setminus\mathcal{A}$. We can now define the
linear mapping $\mathcal{M}:X_{1/w}^{0}\rightarrow L^{2}$ by%
\begin{equation}
\mathcal{M}u=\left(  \frac{u_{F}}{\sqrt{w}}\right)  ^{\vee},\text{\quad}u\in
X_{1/w}^{0}.\label{a1.067}%
\end{equation}

\end{definition}

The operator $\mathcal{M}$ has the following properties:

\begin{theorem}
\label{Thm_M_property} The operator $\mathcal{M}:X_{1/w}^{0}\rightarrow L^{2}
$ is linear and in the seminorm sense it is isometric, 1-1 with null space
$\left(  S_{\mathcal{A}}^{\prime}\right)  ^{\vee}$.

Also, $\tau_{a}\mathcal{M=M}\tau_{a}$ where $\tau_{a}$ is the translation
operator $\tau_{a}u=u\left(  \cdot-a\right)  $.
\end{theorem}

\begin{proof}
That $\mathcal{M}$ is an isometry is clear from the Definition \ref{Def_Xo1/w}
of $X_{1/w}^{0}$. Since $\mathcal{M}$ is an isometry the null space of
$\mathcal{M}$ is the null space of the seminorm $\left\vert \cdot\right\vert
_{1/w,0}$, namely $\left(  S_{\mathcal{A}}^{\prime}\right)  ^{\vee}$. Finally,%
\[
\tau_{a}\mathcal{M}u=\tau_{a}\left(  \frac{u_{F}}{\sqrt{w}}\right)  ^{\vee
}=\left(  e^{-ia\xi}\frac{u_{F}}{\sqrt{w}}\right)  ^{\vee}=\left(
\frac{\left(  \tau_{a}u\right)  _{F}}{\sqrt{w}}\right)  ^{\vee}=\mathcal{M}%
\tau_{a}u.
\]

\end{proof}

\begin{definition}
\label{Def_barXow}\textbf{The inner product space }$\overline{X}_{1/w}^{0}$
Suppose that the weight function $w$ has property W01 for the set
$\mathcal{A}$, where $\mathcal{A}$ is a closed set of measure zero. Then we
define the inner product space:%
\[
\overline{X}_{1/w}^{0}:=\left\{  u\in S^{\prime}:\widehat{u}\in L_{loc}%
^{1}\text{ }and\text{ }\int\frac{\left\vert \widehat{u}\right\vert ^{2}}%
{w}<\infty\right\}  .
\]

Endow $\overline{X}_{1/w}^{0}$ with the norm and inner product $\left\Vert
u\right\Vert _{w,0}^{2}=\int\frac{\left\vert \widehat{u}\right\vert ^{2}}{w}
$, $\left(  u,v\right)  _{w,0}=\int\frac{\widehat{u}\overline{\widehat{v}}}{w}
$.

Clearly we have $\overline{X}_{1/w}^{0}\hookrightarrow X_{1/w}^{0}$ and the
embedding is isometric.
\end{definition}

\begin{theorem}
\label{Thm_barXow_tenprod_quotient}\ ???

\begin{enumerate}
\item ?? We have the direct sum formula $X_{1/w}^{0}=\overline{X}_{1/w}%
^{0}+\left(  S_{\mathcal{A}}^{\prime}\right)  ^{\vee}$?

\item ?? $\overline{X}_{1/w}^{0}$ is a Hilbert space? \textbf{Construct
mapping to} $L^{2}$!

\item ?? MORE!
\end{enumerate}
\end{theorem}

\begin{proof}
\textbf{Parts} 1 \textbf{and} 2: ??.\medskip

\textbf{Part 3} ??
\end{proof}

\begin{theorem}
\label{Thm_L2_embed_Xoinvw_Xow}Suppose $w\in W01$ so that $X_{w}^{0}$ is an
inner product space and $\overline{X}_{1/w}^{0}$ is an product space
(Definition \ref{Def_barXow}).

Then $X_{w}^{0}\hookrightarrow L^{2}\hookrightarrow\overline{X}_{1/w}^{0}$ iff
$1/w\in L^{\infty}$.
\end{theorem}

\begin{proof}
Theorem \ref{Thm_Xow_embed_L2_iff} implies that $X_{w}^{0}\hookrightarrow
L^{2}$ iff $1/w\in L^{\infty}$. Thus $X_{w}^{0}\hookrightarrow L^{2}%
\hookrightarrow\overline{X}_{1/w}^{0}$ implies $1/w\in L^{\infty}$.

Conversely, if $1/w\in L^{\infty}$ then $X_{w}^{0}\hookrightarrow L^{2}$.
Further, if $1/w\in L^{\infty}$ and $f\in L^{2}$ then $f\in S^{\prime}$,
$\widehat{f}\in L^{2}\subset L_{loc}^{1}\left(  \mathbb{R}^{d}\right)  \subset
L_{loc}^{1}\left(  \mathbb{R}^{d}\setminus\mathcal{A}\right)  $, and
$\int\frac{1}{w}\left\vert \widehat{f}\right\vert ^{2}\leq\left\Vert \frac
{1}{w}\right\Vert _{\infty}\int\left\vert \widehat{f}\right\vert
^{2}=\left\Vert \frac{1}{w}\right\Vert _{\infty}\int\left\vert f\right\vert
^{2}<\infty$. Thus by definition $f\in\overline{X}_{1/w}^{0}$ and
$L^{2}\hookrightarrow\overline{X}_{1/w}^{0}$.
\end{proof}

\begin{theorem}
\label{Thm_L2_embed_Xoinvw_iff_w_in_Linf}Suppose $w\in W01$ so that
$\overline{X}_{1/w}^{0}$ is an inner product space.

Then $L^{2}\hookrightarrow\overline{X}_{1/w}^{0}$ iff $1/w\in L^{\infty}$.

More precisely, if $\frac{1}{w}\in L^{\infty}$ then%
\[
\left\Vert f\right\Vert _{1/w,0}\leq\left\Vert \frac{1}{w}\right\Vert
_{\infty}^{1/2}\left\Vert f\right\Vert _{2},\quad f\in\overline{X}_{1/w}^{0},
\]

and if the embedding $\iota:L^{2}\hookrightarrow\overline{X}_{1/w}^{0}$ is
continuous then%
\[
\left\Vert \frac{1}{w}\right\Vert _{\infty}^{1/2}=\left\Vert \iota\right\Vert
.
\]

\end{theorem}

\begin{proof}
From the theorem above we know that $1/w\in L^{\infty}$ implies $L^{2}%
\hookrightarrow\overline{X}_{1/w}^{0}$.

Now suppose that we have the continuous embedding $L^{2}\overset{\iota
}{\hookrightarrow}\overline{X}_{1/w}^{0}$. Since $S\subset L^{2}$,%
\begin{equation}
\left\vert \phi\right\vert _{1/w,0}^{2}=\int\frac{1}{w}\left\vert
\widehat{\phi}\right\vert ^{2}\leq\left\Vert \iota\right\Vert _{op}%
\int\left\vert \widehat{\phi}\right\vert ^{2},\quad\phi\in S,\label{a50.2}%
\end{equation}

i.e.%
\[
\int\frac{1}{w}\left\vert \phi\right\vert ^{2}\leq\left\Vert \iota\right\Vert
_{op}\int\left\vert \phi\right\vert ^{2},\quad\phi\in S.
\]

With reference to Definition \ref{Def_extend_wt_fn} of the weight function,
choose arbitrary $x\notin\mathcal{A}$ where $\mathcal{A}$ is a closed set of
measure zero such that $w>0$ and continuous on $\mathbb{R}^{d}\setminus
\mathcal{A}$. Suppose $0<r<R<\operatorname*{dist}\left(  x,\mathcal{A}\right)
$. Then it is possible to choose $\phi_{r,R}\in S$ such that $0\leq\phi
_{r,R}\leq1$, $\operatorname*{supp}\phi_{r,R}\subseteq\overline{B}_{R}\left(
x\right)  $ and $\phi_{r,R}=1$ on $\overline{B}_{r}\left(  x\right)  $ e.g.
the standard "cap" example of a $C_{0}^{\infty}$ function. Observe now that%
\[
\min_{B_{r}\left(  x\right)  }\frac{1}{w}\int\limits_{B_{r}\left(  x\right)
}1\leq\int\limits_{B_{R}\left(  x\right)  }\frac{\left\vert \phi
_{r,R}\right\vert ^{2}}{w}\leq\left\Vert \iota\right\Vert ^{2}\int%
\limits_{B_{R}\left(  x\right)  }\left\vert \phi_{r,R}\right\vert ^{2}%
\leq\left\Vert \iota\right\Vert ^{2}\int\limits_{B_{R}\left(  x\right)  }1,
\]

i.e. for $0<r<R<\operatorname*{dist}\left(  x,\mathcal{A}\right)  $,%
\[
\min\limits_{B_{r}\left(  x\right)  }\frac{1}{w}\leq\left\Vert \iota
\right\Vert ^{2}\frac{\int_{B_{R}\left(  x\right)  }1}{\int_{B_{r}\left(
x\right)  }1}=\left\Vert \iota\right\Vert ^{2}\left(  \frac{R}{r}\right)
^{d}.
\]

Now since the continuity of $w$ near $x$ implies $\min\limits_{B_{r}\left(
x\right)  }\frac{1}{w}$ is continuous in $R$ we must have $\frac{1}%
{\max\limits_{B_{r}\left(  x\right)  }\frac{1}{w}}\leq\left\Vert
\iota\right\Vert ^{2}$ for $0<r<\operatorname*{dist}\left(  x,\mathcal{A}%
\right)  $, and hence that $\frac{1}{w\left(  x\right)  }\leq\left\Vert
\iota\right\Vert ^{2} $ when $x\notin\mathcal{A}$, i.e. $1/w\in L^{\infty}$
and $\left\Vert \frac{1}{w}\right\Vert _{\infty}^{1/2}\leq\left\Vert
\iota\right\Vert $. But from \ref{a50.2}, $\left\Vert \iota\right\Vert
\leq\left\Vert \frac{1}{w}\right\Vert _{\infty}^{1/2}$ so $\left\Vert
\iota\right\Vert =\left\Vert \frac{1}{w}\right\Vert _{\infty}^{1/2}$.
\end{proof}

\subsection{The spaces $S_{\emptyset,k}$ and $S_{\emptyset,k}^{\prime}%
$\label{SbSect_Sok_fnal_Sok}}

We begin by preparing the way for the application of the Hahn-Banach extension
theorem to subspaces of the Schwartz functions of rapidly decrease $S$ and
then consider the properties of some special subspaces of functions which have
several zero derivatives at the origin.

\begin{lemma}
\label{Lem_convex_tls_extend_2}\textbf{Schwartz functions and the Hahn-Banach
extension theorem}

\begin{enumerate}
\item (Hahn-Banach) Any continuous linear functional $f$ on a subspace $M$ of
a locally convex topological vector space $T$ can be extended non-uniquely to
a continuous linear functional $f^{e}$ on $T$ ($f^{e}\in T^{\prime}$).

The set of extensions is $f^{e}+M^{\bot}$ where $M^{\bot}$ is the set of
annihilators of $M$ i.e. the members of $T^{\prime}$ which are zero on each
member of $M$.

\item The space of rapidly decreasing Schwartz functions $S$ is a locally
convex topological space when endowed with the topology induced by any of the
(equivalent) countable seminorm sets given in part 1 of Definition
\ref{Def_Distributions}.

\item A linear functional $f$ defined on a subspace $M$ of the space $S$ is
continuous iff%
\[
\left\vert \left[  f,\psi\right]  \right\vert \leq\left\vert \psi\right\vert
_{\sigma},\quad\psi\in M.
\]

where $\left\vert \cdot\right\vert _{\sigma}$ denotes a positive, linear
combination of the seminorms in the (equivalent) sets given in part 1 of
Definition \ref{Def_Distributions}. We then write $f\in M^{\prime}$.
\end{enumerate}
\end{lemma}

\begin{proof}
This lemma can be proved, for example, by using the results and definitions of
Chapter V., Volume I of Reed and Simon \cite{ReedSimon72}.
\end{proof}

We will now introduce a special class of subspaces of the Schwartz spaces $S$
that are widely used in the positive order theory expounded in Williams
\cite{WilliamsPosOrdSmthV3}. These spaces are:

\begin{definition}
\label{Def_So,n}\textbf{The space} $S_{\emptyset,n}$
\begin{equation}
S_{\emptyset,0}=S,\qquad S_{\emptyset,n}=\left\{  \phi\in S:D^{\alpha}%
\phi\left(  0\right)  =0,\text{\ }\left\vert \alpha\right\vert <n\right\}
,\;n=1,2,3,\ldots,\label{a1.50}%
\end{equation}

and we endow $S_{\emptyset,n}$ with the subspace topology induced by the space
$S$. $S$ is the space of $C^{\infty}$ functions of rapid decrease used as test
functions for the tempered distributions. $S$ is endowed with any of the
countable seminorm topology described in Definition \ref{Def_Distributions} of
the Appendix.
\end{definition}

The next inequality gives a simple upper bound for functions in $S_{\emptyset
,n}$ near the origin. This follows directly from the estimate (\ref{a1.37}\ in
the Appendix) of the integral remainder term of the Taylor series expansion.%
\begin{equation}
\left\vert u(x)\right\vert \leq\left(  \sum\limits_{\left\vert \alpha
\right\vert =n}\left\Vert D^{\alpha}u\right\Vert _{\infty,B_{1}}\right)
\left\vert x\right\vert ^{n},\quad u\in S_{\emptyset,n},\text{ }x\in
\overline{B}_{1}.\label{a992}%
\end{equation}

Consequently%
\begin{equation}
\int\frac{\left\vert u(x)\right\vert }{\left\vert x\right\vert ^{n}}%
dx<\max\left\{  \sum\limits_{\left\vert \alpha\right\vert =n}\left\Vert
D^{\alpha}u\right\Vert _{\infty,B_{1}},\left\Vert \frac{u}{\left\vert
\cdot\right\vert ^{n}}\right\Vert _{\infty,\mathbb{R}^{d}\setminus B_{1}%
}\right\}  ,\quad u\in S_{\emptyset,n}.\label{a995}%
\end{equation}

From Appendix \ref{Def_Some_basic_spaces} we know that $P_{n}$ is the space of
polynomials of order at most $n$ ($\deg P_{n}\leq n-1$) with complex
coefficients, and $P$ is the space of all polynomials with complex
coefficients. We will let $\widehat{P}_{n}$ denote the Fourier transforms of
the functions in $P_{n}$ and $\widehat{P}$ denote the Fourier transforms of
the functions in $P$.

\begin{theorem}
\label{Thm_prop_functnl_on_Son}\ 

\begin{enumerate}
\item $\widehat{P}_{n}=\left\{  g\in S^{\prime}:\left[  g,\phi\right]
=0\text{ }for\text{ }all\text{ }\phi\in S_{\emptyset,n}\right\}  $ i.e.
$S_{\emptyset,n}^{\bot}=\widehat{P}_{n}$.

\item If $f\in S_{\emptyset,n}^{\prime}$ then there exists $f^{e}\in
S^{\prime}$ such that $f^{e}=f$ on $S_{\emptyset,n}$. The set of extensions is
$f^{e}+\widehat{P}_{n}$. We sometimes say that $f$ can be extended to $S$ as a
member $f^{e}$ of $S^{\prime}$.

\item $P_{n}=\left\{  g\in S^{\prime}:\left[  g,\widehat{\phi}\right]
=0\text{ }for\text{ }all\text{ }\phi\in S_{\emptyset,n}\right\}  $.
\end{enumerate}
\end{theorem}

\begin{proof}
\textbf{Part 1} Suppose $\left[  u,\phi\right]  =0$ for all $\phi\in
S_{\emptyset,n}$. This implies that the $\operatorname*{supp}u\subset\left\{
0\right\}  $ and by a well known theorem in distribution theory,
$u\in\widehat{P}$. Thus $u=\widehat{p}$ for some polynomial $p$. Suppose $\deg
p>n $. For each $\left\vert \beta\right\vert >n$ choose $\phi_{\beta}\in
S_{\emptyset,n}$ such that $\phi_{\beta}\left(  x\right)  =x^{\beta}/\beta!$
in a neighborhood of zero. Then we have $\left(  D^{\alpha}\phi_{\beta
}\right)  \left(  0\right)  =\delta_{\alpha,\beta}$ and if the coefficients of
$p$ are $b_{\alpha}$,

$0=\left[  \widehat{p},\phi_{\beta}\right]  =\left[  p(-D)\left(  \phi_{\beta
}\right)  \right]  (0)=\left(  -1\right)  ^{\left\vert \beta\right\vert
}b_{\beta}$, and thus $\deg p<n$ and $u\in\widehat{P}_{n}$.\medskip

\textit{Conversely}, suppose $u\in\widehat{P}_{n}$. Then there exists $p\in
P_{n}$ such that $u=\widehat{p}=\left(  2\pi\right)  ^{d/2}p(iD)\delta$.
Hence, if $\phi\in S_{\emptyset,n}$
\begin{align*}
\left[  u,\phi\right]   &  =\left[  \left(  2\pi\right)  ^{d/2}p(iD)\delta
,\phi\right]  =\left(  2\pi\right)  ^{d/2}\left[  p(iD)\delta,\phi\right]
=\left(  2\pi\right)  ^{d/2}\left[  p(iD)\delta,\phi\right] \\
&  =\left(  2\pi\right)  ^{d/2}\left[  \delta,p(-iD)\phi\right]  =\left(
2\pi\right)  ^{d/2}\left[  p(-iD)\phi\right]  \left(  0\right)  =0.
\end{align*}
\medskip

\textbf{Part 2} We use Lemma \ref{Lem_convex_tls_extend_2} with
$M=S_{\emptyset,n}$ and $T=S$. In this case the the set of extensions is
$f^{e}+S_{\emptyset,n}^{\bot}$ and $g\in S_{\emptyset,n}^{\bot}$ iff $g\in
S^{\prime}$ and $\left[  g,\phi\right]  =0$ for all $\phi\in S_{\emptyset,n}$
i.e. $g\in\widehat{P}_{n}$ by part 1.\medskip

\textbf{Part 3 }Use part 1.
\end{proof}

We give the following results without proof:

\begin{theorem}
\label{Thm_product_of_Co,k_funcs}If $\phi\in S_{\emptyset,m}$ and $\psi\in
S_{\emptyset,n}$ then:

\begin{enumerate}
\item $\phi\psi\in S_{\emptyset,m+n}$.

\item If $\left\vert \alpha\right\vert =k$ then $x^{\alpha}\psi\in
S_{\emptyset,n+k}$.

\item $\left\vert \cdot\right\vert ^{2k}\psi\in S_{\emptyset,n+2k}$.
\end{enumerate}
\end{theorem}

\begin{theorem}
\label{Thm_abs_squar_and_Som}If $\phi\in S$ then for all $n\geq1$: $\left\vert
\phi\right\vert ^{2}\in S_{\emptyset,2n}$ iff $\phi\in S_{\emptyset,n}$.
\end{theorem}

\begin{proof}
From part 1 Theorem \ref{Thm_product_of_Co,k_funcs} we know that $\phi\in
S_{\emptyset,n}$ implies $\left\vert \phi\right\vert ^{2}\in S_{\emptyset,2n}$.

The converse will now be proved. Firstly, since $\left\vert \phi\right\vert
^{2}\in S_{\emptyset,2n}$ we have $\phi\left(  0\right)  =0$. Now suppose
there exists a multi-index $\beta$ such that: $0<\left\vert \beta\right\vert
<m$, $D^{\beta}\phi\left(  0\right)  \neq0$ and $D^{\gamma}\phi\left(
0\right)  =0$ when $\left\vert \gamma\right\vert <\left\vert \beta\right\vert
$. We now write%
\begin{align*}
0=D^{2\beta}\left(  \left\vert \phi\right\vert ^{2}\right)  \left(  0\right)
&  =\sum_{\alpha\leq2\beta}\binom{2\beta}{\alpha}D^{\alpha}\phi\left(
0\right)  \overline{D^{2\beta-\alpha}\phi\left(  0\right)  }\\
&  =\left\vert \binom{2\beta}{\beta}D^{\beta}\phi\left(  0\right)  \right\vert
^{2}+\sum_{\substack{\alpha\leq2\beta\\\alpha\neq\beta}}\binom{2\beta}{\alpha
}D^{\alpha}\phi\left(  0\right)  \overline{D^{2\beta-\alpha}\phi\left(
0\right)  },
\end{align*}

and observe that for each term in the last summation either $\left\vert
\alpha\right\vert <\left\vert \beta\right\vert $ or $\left\vert 2\beta
-\alpha\right\vert <\left\vert \beta\right\vert $. Thus our assumptions mean
that the summation term is actually zero and so $D^{\beta}\phi\left(
0\right)  =0$, contradicting our assumptions. We have thus shown that $\phi\in
S_{\emptyset,n}$.
\end{proof}

\subsection{The operator $\mathcal{L}_{1}:L^{2}\rightarrow X_{1/w}^{0}%
$\label{SbSect_op_L1}}

Making certain assumptions about the weight function we introduce the operator
$\mathcal{L}_{1}:L^{2}\rightarrow S^{\prime}$. This is the analogue of the
mapping $\mathcal{J}:L^{2}\rightarrow X_{w}^{0}$ from $X_{1/w}^{0}$ to $L^{2}$
introduced in Definition \ref{Def_I_J}. We show that $\mathcal{L}_{1}%
:L^{2}\rightarrow X_{1/w}^{0}$ and that it has nice properties in the seminorm sense.

We begin by discussing the following weight function property: For some
integer $m\geq0$ and constants $\tau,r\geq0$ the weight function $w$ satisfies
the two conditions
\begin{equation}
w\left\vert \cdot\right\vert ^{2m}\in L_{loc}^{1}\left(  \mathbb{R}%
^{d}\right)  ,\quad\int\limits_{\left\vert \cdot\right\vert \geq r}\frac
{w}{\left\vert \cdot\right\vert ^{2\tau}}<\infty.\label{a1.046}%
\end{equation}

\begin{lemma}
\label{Lem_wt_fn_extra_property_2}Suppose the weight function $w$ satisfies
\ref{a1.046}. Then $w\in L_{loc}^{1}\left(  \mathbb{R}^{d}\setminus0\right)
$. Further:

\begin{enumerate}
\item If $w\in L_{loc}^{1}$ and $g\in L^{2}$ then $\sqrt{w}g\in S^{\prime}\cap
L_{loc}^{1}$.

\item If $g\in L^{2}$ then $\sqrt{w}g\in S_{\emptyset,m}^{\prime}\cap
L_{loc}^{1}\left(  \mathbb{R}^{d}\setminus0\right)  \cap L_{loc}^{1}\left(
\mathbb{R}^{d}\setminus\mathcal{A}\right)  $.

\item If $\psi\in S_{\emptyset,2m}$ the linear functional $\int w\psi$ is a
member of $S_{\emptyset,2m}^{\prime}$ i.e.%
\[
\left\vert \int w\psi\right\vert \leq\left\vert \psi\right\vert _{\sigma
},\quad\psi\in S_{\emptyset,2m},
\]

where $\left\vert \psi\right\vert _{\sigma}$ is a positive, linear combination
of the seminorms \ref{a1.2} which define the topology of the tempered
distributions $S$.

\item If $\psi\in S_{\emptyset,m}$ then the non-linear functional $\left(
\int w\left\vert \psi\right\vert ^{2}\right)  ^{1/2}$ satisfies%
\[
\left(  \int w\left\vert \psi\right\vert ^{2}\right)  ^{1/2}\leq\left\vert
\psi\right\vert _{\sigma},\quad\psi\in S_{\emptyset,m},
\]

where $\left\vert \cdot\right\vert _{\sigma}$ is a positive, linear
combination of the seminorms \ref{a1.2} which define the topology of the
tempered distributions $S$.
\end{enumerate}
\end{lemma}

\begin{proof}
Since $w\left\vert \cdot\right\vert ^{2m}\in L_{loc}^{1}\left(  \mathbb{R}%
^{d}\right)  $ it is clear that $w\in L_{loc}^{1}\left(  \mathbb{R}%
^{d}\setminus0\right)  $.\medskip

\textbf{Part 1} If $K$ is compact then $\int\limits_{K}\sqrt{w}g\leq\left(
\int\limits_{K}w\right)  ^{\frac{1}{2}}\left(  \int\limits_{K}\left\vert
g\right\vert ^{2}\right)  ^{\frac{1}{2}}<\infty$ so $\sqrt{w}g\in L_{loc}^{1}
$. If $\phi\in S$ then%
\begin{align*}
\int\sqrt{w}g\phi & =\int\limits_{\left\vert \cdot\right\vert \leq r}\sqrt
{w}g\phi+\int\limits_{\left\vert \cdot\right\vert \geq r}\sqrt{w}g\phi
=\int\limits_{\left\vert \cdot\right\vert \leq r}\sqrt{w}g\phi+\int%
\limits_{\left\vert \cdot\right\vert \geq r}\frac{\sqrt{w}}{\left\vert
\cdot\right\vert ^{\tau}}g\left\vert \cdot\right\vert ^{\tau}\phi\\
& \leq\left(  \int\limits_{\left\vert \cdot\right\vert \leq r}w\right)
^{1/2}\left(  \int\limits_{\left\vert \cdot\right\vert \leq r}\left\vert
g\phi\right\vert ^{2}\right)  ^{\frac{1}{2}}+\left(  \int\limits_{\left\vert
\cdot\right\vert \geq r}\frac{w}{\left\vert \cdot\right\vert ^{2\tau}}\right)
^{1/2}\left(  \int\limits_{\left\vert \cdot\right\vert \geq r}\left\vert
g\right\vert ^{2}\left\vert \cdot\right\vert ^{2\tau}\left\vert \phi
\right\vert ^{2}\right)  ^{1/2}\\
& \leq\left(  \int\limits_{\left\vert \cdot\right\vert \leq r}w\right)
^{1/2}\left(  \int\limits_{\left\vert \cdot\right\vert \leq r}\left\vert
g\right\vert ^{2}\right)  ^{\frac{1}{2}}\left\Vert \phi\right\Vert _{\infty
}+\left(  \int\limits_{\left\vert \cdot\right\vert \geq r}\frac{w}{\left\vert
\cdot\right\vert ^{2\tau}}\right)  ^{1/2}\left(  \int\limits_{\left\vert
\cdot\right\vert \geq r}\left\vert g\right\vert ^{2}\right)  ^{1/2}\left\Vert
\left\vert \cdot\right\vert ^{\tau}\phi\right\Vert _{\infty}\\
& \leq\left(  \left(  \int\limits_{\left\vert \cdot\right\vert \leq
r}w\right)  ^{1/2}+\left(  \int\limits_{\left\vert \cdot\right\vert \geq
r}\frac{w}{\left\vert \cdot\right\vert ^{2\tau}}\right)  ^{1/2}\right)
\left\Vert g\right\Vert _{2}\left(  \left\Vert \phi\right\Vert _{\infty
}+\left\Vert \left(  1+\left\vert \cdot\right\vert \right)  ^{\left\lceil
\tau\right\rceil }\phi\right\Vert _{\infty}\right)  ,
\end{align*}

which is a positive, linear combination of the seminorms \ref{a1.2}.\medskip

\textbf{Part 2} Since $w\in L_{loc}^{1}\left(  \mathbb{R}^{d}\setminus
0\right)  $ for compact $K\subset\mathbb{R}^{d}\setminus0$ we have by the
Cauchy-Schwartz inequality%
\[
\int\limits_{K}\sqrt{w}\left\vert g\right\vert \leq\left(  \int\limits_{K}%
w\right)  ^{1/2}\left\Vert g\right\Vert _{2}<\infty.
\]

Also, if $K\subset\mathbb{R}^{d}\setminus\mathcal{A}$ is compact then $w$ is
continuous on $K$ and again%
\[
\int\limits_{K}\sqrt{w}\left\vert g\right\vert \leq\left(  \int\limits_{K}%
w\right)  ^{1/2}\left\Vert g\right\Vert _{2}<\infty.
\]

If $\psi\in S$ then
\begin{equation}
\left\vert \left[  \sqrt{w}\left\vert g\right\vert ,\psi\right]  \right\vert
\leq\int\sqrt{w}\left\vert g\right\vert \left\vert \psi\right\vert
\leq\left\Vert g\right\Vert _{2}\left(  \int w\left\vert \psi\right\vert
^{2}\right)  ^{1/2},\label{a1.004}%
\end{equation}

and if $\psi\in S_{\emptyset,m}$
\begin{align}
\left(  \int w\left\vert \psi\right\vert ^{2}\right)  ^{1/2}  & \leq\left(
\int\limits_{\left\vert \cdot\right\vert \leq r}w\left\vert \psi\right\vert
^{2}\right)  ^{1/2}+\left(  \int\limits_{\left\vert \cdot\right\vert \geq
r}w\left\vert \psi\right\vert ^{2}\right)  ^{1/2}\nonumber\\
& =\left(  \int\limits_{\left\vert \cdot\right\vert \leq r}w\left\vert
\cdot\right\vert ^{2m}\frac{\left\vert \psi\right\vert ^{2}}{\left\vert
\cdot\right\vert ^{2m}}\right)  ^{1/2}+\left(  \int\limits_{\left\vert
\cdot\right\vert \geq r}w\left\vert \psi\right\vert ^{2}\right)
^{1/2}\nonumber\\
& \leq\left(  \int\limits_{\left\vert \cdot\right\vert \leq r}w\left\vert
\cdot\right\vert ^{2m}\right)  ^{1/2}\left\Vert \frac{\psi}{\left\vert
\cdot\right\vert ^{m}}\right\Vert _{\infty}+\left(  \int\limits_{\left\vert
\cdot\right\vert \geq r}w\left\vert \psi\right\vert ^{2}\right)
^{1/2}.\label{a1.003}%
\end{align}

Next we estimate the second term of \ref{a1.003}:%
\begin{align*}
\left(  \int\limits_{\left\vert \cdot\right\vert \geq r}w\left\vert
\psi\right\vert ^{2}\right)  ^{1/2}=\left(  \int\limits_{\left\vert
\cdot\right\vert \geq r}\frac{w}{\left\vert \cdot\right\vert ^{2\tau}%
}\left\vert \cdot\right\vert ^{2\tau}\left\vert \psi\right\vert ^{2}\right)
^{1/2} &  =r^{\tau}\left(  \int\limits_{\left\vert \cdot\right\vert \geq
r}\frac{w}{\left\vert \cdot\right\vert ^{2\tau}}\left(  \frac{\left\vert
\cdot\right\vert }{r}\right)  ^{2\tau}\left\vert \psi\right\vert ^{2}\right)
^{1/2}\\
&  \leq r^{\tau}\left(  \int\limits_{\left\vert \cdot\right\vert \geq r}%
\frac{w}{\left\vert \cdot\right\vert ^{2\tau}}\left(  \frac{\left\vert
\cdot\right\vert }{r}\right)  ^{2\left\lceil \tau\right\rceil }\left\vert
\psi\right\vert ^{2}\right)  ^{1/2}\\
&  =r^{\tau-\left\lceil \tau\right\rceil }\left(  \int\limits_{\left\vert
\cdot\right\vert \geq r}\frac{w}{\left\vert \cdot\right\vert ^{2\tau}%
}\left\vert \cdot\right\vert ^{2\left\lceil \tau\right\rceil }\left\vert
\psi\right\vert ^{2}\right)  ^{1/2}\\
&  \leq r^{\tau-\left\lceil \tau\right\rceil }\left(  \int\limits_{\left\vert
\cdot\right\vert \geq r}\frac{w}{\left\vert \cdot\right\vert ^{2\tau}}\right)
^{1/2}\left\Vert \left\vert \cdot\right\vert ^{\left\lceil \tau\right\rceil
}\psi\right\Vert _{\infty},
\end{align*}

so that
\[
\left(  \int w\left\vert \psi\right\vert ^{2}\right)  ^{1/2}\leq\left(
\int\limits_{\left\vert \cdot\right\vert \leq r}w\left\vert \cdot\right\vert
^{2m}\right)  ^{1/2}\left\Vert \frac{\psi}{\left\vert \cdot\right\vert ^{m}%
}\right\Vert _{\infty}+r^{\tau-\left\lceil \tau\right\rceil }\left(
\int\limits_{\left\vert \cdot\right\vert \geq r}\frac{w}{\left\vert
\cdot\right\vert ^{2\tau}}\right)  ^{1/2}\left\Vert \left\vert \cdot
\right\vert ^{\left\lceil \tau\right\rceil }\psi\right\Vert _{\infty}.
\]

But by \ref{a992} we have%
\begin{equation}
\left(  \int w\left\vert \psi\right\vert ^{2}\right)  ^{1/2}\leq\left(
\int\limits_{\left\vert \cdot\right\vert \leq r}w\left\vert \cdot\right\vert
^{2m}\right)  ^{1/2}2^{md/2}\sum\limits_{\left\vert \alpha\right\vert
=m}\left\Vert D^{\alpha}\psi\right\Vert _{\infty}+r^{\tau-\left\lceil
\tau\right\rceil }\left(  \int\limits_{\left\vert \cdot\right\vert \geq
r}\frac{w}{\left\vert \cdot\right\vert ^{2\tau}}\right)  ^{1/2}\left\Vert
\left\vert \cdot\right\vert ^{\left\lceil \tau\right\rceil }\psi\right\Vert
_{\infty},\label{a1.007}%
\end{equation}

which is a positive linear combination of the $L^{\infty}$ seminorms on $S$.
Thus inequality \ref{a1.004} implies that $\sqrt{w}\left\vert g\right\vert \in
S_{\emptyset,m}^{\prime}$.\medskip

\textbf{Part 3}. If $\psi\in S_{\emptyset,2m}$ then the assumptions
\ref{a1.046} imply
\begin{align*}
\left\vert \int w\psi\right\vert  & \leq\int\limits_{\left\vert \cdot
\right\vert \leq r}w\left\vert \cdot\right\vert ^{2m}\frac{\left\vert
\psi\right\vert }{\left\vert \cdot\right\vert ^{2m}}+\int\limits_{\left\vert
\cdot\right\vert \geq r}\frac{w}{\left\vert \cdot\right\vert ^{2\tau}%
}\left\vert \cdot\right\vert ^{2\tau}\left\vert \psi\right\vert \\
& \leq\left(  \int\limits_{\left\vert \cdot\right\vert \leq r}w\left\vert
\cdot\right\vert ^{2m}\right)  \left\Vert \frac{\psi}{\left\vert
\cdot\right\vert ^{2m}}\right\Vert _{\infty}+\left(  \int\limits_{\left\vert
\cdot\right\vert \geq r}\frac{w}{\left\vert \cdot\right\vert ^{2\tau}}\right)
\max_{\left\vert x\right\vert \geq r}\left\{  \left\vert x\right\vert ^{2\tau
}\left\vert \psi\left(  x\right)  \right\vert \right\} \\
& \leq\left(  \int\limits_{\left\vert \cdot\right\vert \leq r}w\left\vert
\cdot\right\vert ^{2m}\right)  \left\Vert \frac{\psi}{\left\vert
\cdot\right\vert ^{2m}}\right\Vert _{\infty}+r^{2\tau-\left\lceil
2\tau\right\rceil }\left(  \int\limits_{\left\vert \cdot\right\vert \geq
r}\frac{w}{\left\vert \cdot\right\vert ^{2\tau}}\right)  \left\Vert \left\vert
\cdot\right\vert ^{\left\lceil 2\tau\right\rceil }\psi\right\Vert _{\infty},
\end{align*}

But by \ref{a992} we have%
\[
\left\vert \int w\psi\right\vert \leq\left(  \int\limits_{\left\vert
\cdot\right\vert \leq r}w\left\vert \cdot\right\vert ^{2m}\right)  2^{md}%
\sum\limits_{\left\vert \alpha\right\vert =2m}\left\Vert D^{\alpha}%
\psi\right\Vert _{\infty}+r^{2\tau-\left\lceil 2\tau\right\rceil }\left(
\int\limits_{\left\vert \cdot\right\vert \geq r}\frac{w}{\left\vert
\cdot\right\vert ^{2\tau}}\right)  \left\Vert \left\vert \cdot\right\vert
^{\left\lceil 2\tau\right\rceil }\psi\right\Vert _{\infty},
\]

which is a positive linear combination of the $L^{\infty}$ seminorms on $S$.
Thus $\int w\psi\in S_{\emptyset,2m}^{\prime}$.\medskip

\textbf{Part 4}. This follows directly from inequality \ref{a1.007} of part 2.
\end{proof}

We now define the operator $\mathcal{L}_{1}:L^{2}\rightarrow S^{\prime}$. It
will be shown that $\mathcal{L}_{1}:L^{2}\rightarrow X_{1/w}^{0}$ and that
$\mathcal{L}_{1}$ is an inverse of $\mathcal{M}:X_{1/w}^{0}\rightarrow L^{2}$.
This will allow us to prove that $\mathcal{M}$ is onto and hence that
$X_{1/w}^{0}$ is complete.

\begin{definition}
\label{Def_op_L1}\textbf{The (class of) operators} $\mathcal{L}_{1}%
:L^{2}\rightarrow S^{\prime}$

Suppose the weight function $w$ has property \ref{a1.046}. If $g\in L^{2}$
then part 2 of Lemma \ref{Lem_wt_fn_extra_property_2} shows that $\sqrt
{w}\widehat{g}\in S_{\emptyset,m}^{\prime}$. By part 3 of Lemma
\ref{Lem_convex_tls_extend_2} $\sqrt{w}\widehat{g}$ can be extended
(non-uniquely) to act on $S$ as a member of $S^{\prime}$. Denote such an
extension by $\left(  \sqrt{w}\widehat{g}\right)  ^{e}$ and define the
operator $\mathcal{L}_{1}:L^{2}\rightarrow S^{\prime}$ by
\[
\mathcal{L}_{1}g=\left(  \left(  \sqrt{w}\widehat{g}\right)  ^{e}\right)
^{\vee},\quad g\in L^{2}.
\]

In general $\mathcal{L}_{1}$ is not unique and\ not linear.
\end{definition}

The next theorem gives some properties of $\mathcal{L}_{1}$.

\begin{theorem}
\label{Thm_L1_propert}Suppose the weight function $w$ has property
\ref{a1.046}. Then:

\begin{enumerate}
\item If $g\in L^{2}$ then $\left(  \mathcal{L}_{1}g\right)  _{F}=\sqrt
{w}\widehat{g}\in L_{loc}^{1}\left(  \mathbb{R}^{d}\setminus0\right)  \cap
L_{loc}^{1}\left(  \mathbb{R}^{d}\setminus\mathcal{A}\right)  $.

\item In the seminorm sense $\mathcal{L}_{1}:L^{2}\rightarrow X_{1/w}^{0}$,
$\mathcal{L}_{1}$ is an isometry and is 1-1.

\item $\mathcal{L}_{1}$ is linear in the seminorm sense i.e.%
\[
\left\vert \mathcal{L}_{1}\left(  \lambda_{1}g_{1}+\lambda_{2}g_{2}\right)
-\mathcal{L}_{1}\left(  \lambda_{1}g_{1}\right)  -\mathcal{L}_{1}\left(
\lambda_{2}g_{2}\right)  \right\vert _{1/w,0}=0.
\]

\item In the seminorm sense, $\tau_{a}\mathcal{L}_{1}=\mathcal{L}_{1}\tau_{a}$
on $L^{2}$ where $\tau_{a}$ is the translation operator $\tau_{a}u=u\left(
\cdot-a\right)  $.
\end{enumerate}
\end{theorem}

\begin{proof}
\textbf{Part 1}. From Lemma \ref{Lem_wt_fn_extra_property_2} we know that if
$g\in L^{2}$ then $\mathcal{L}_{1}g\in S^{\prime}$ and $\left(  \mathcal{L}%
_{1}g\right)  _{F}=\sqrt{w}\widehat{g}\in L_{loc}^{1}\left(  \mathbb{R}%
^{d}\setminus0\right)  \cap L_{loc}^{1}\left(  \mathbb{R}^{d}\setminus
\mathcal{A}\right)  $.\medskip

\textbf{Part 2}. From the definition of $\mathcal{L}_{1}$%
\[
\left\vert \mathcal{L}_{1}g\right\vert _{1/w,0}^{2}=\int\frac{\left\vert
\left(  \mathcal{L}_{1}g\right)  _{F}\right\vert ^{2}}{w}=\int\frac{\left\vert
\sqrt{w}\widehat{g}\right\vert ^{2}}{w}=\left\Vert g\right\Vert _{2}%
^{2}<\infty.
\]

and so $\mathcal{L}_{1}g\in X_{1/w}^{0}$ and $\mathcal{L}_{1}$ is an isometry.
Clearly $\mathcal{L}_{1}g\in\left(  S_{\mathcal{A}}^{\prime}\right)  ^{\vee}$
iff $g=0$ so $\mathcal{L}_{1}$ is 1-1.\medskip

\textbf{Part 3}. Follows from%
\begin{align*}
& \left\vert \mathcal{L}_{1}\left(  \lambda_{1}g_{1}+\lambda_{2}g_{2}\right)
-\mathcal{L}_{1}\left(  \lambda_{1}g_{1}\right)  -\mathcal{L}_{1}\left(
\lambda_{2}g_{2}\right)  \right\vert _{1/w,0}\\
& =\int\frac{\left\vert \left(  \mathcal{L}_{1}\left(  \lambda_{1}%
g_{1}+\lambda_{2}g_{2}\right)  -\mathcal{L}_{1}\left(  \lambda_{1}%
g_{1}\right)  -\mathcal{L}_{1}\left(  \lambda_{2}g_{2}\right)  \right)
_{F}\right\vert ^{2}}{w}\\
& =\int\frac{\left\vert \sqrt{w}\left(  \lambda_{1}g_{1}+\lambda_{2}%
g_{2}\right)  ^{\wedge}-\sqrt{w}\left(  \lambda_{1}g_{1}\right)  ^{\wedge
}-\sqrt{w}\left(  \lambda_{2}g_{2}\right)  ^{\wedge}\right\vert ^{2}}{w}\\
& =0.
\end{align*}
\medskip

\textbf{Part 4} If $g\in L^{2}$ then
\begin{align*}
\left(  \tau_{a}\mathcal{L}_{1}g-\mathcal{L}_{1}\tau_{a}g\right)  _{F}  &
=\left(  \tau_{a}\mathcal{L}_{1}g\right)  _{F}-\left(  \mathcal{L}_{1}\tau
_{a}g\right)  _{F}\\
& =\widehat{\left(  \tau_{a}\mathcal{L}_{1}g\right)  }\mid_{\mathbb{R}%
^{d}\setminus\mathcal{A}}-\sqrt{w}\widehat{\tau_{a}g}\mid_{\mathbb{R}%
^{d}\setminus\mathcal{A}}\\
& =e^{ia\xi}\sqrt{w}\widehat{g}\mid_{\mathbb{R}^{d}\setminus\mathcal{A}%
}-e^{ia\xi}\sqrt{w}\widehat{g}\mid_{\mathbb{R}^{d}\setminus\mathcal{A}}\\
& =0,
\end{align*}

which ensures that $\left\vert \tau_{a}\mathcal{L}_{1}-\mathcal{L}_{1}\tau
_{a}\right\vert _{1/w,0}=0$.
\end{proof}

The following theorem indicates how the operators $\mathcal{L}_{1}%
:L^{2}\rightarrow X_{1/w}^{0}$ and $\mathcal{M}:X_{1/w}^{0}\rightarrow L^{2}$ interact.

\begin{theorem}
\label{Thm_M_L1_propert}Suppose the weight function $w$ has property
\ref{a1.046} which was assumed in Lemma \ref{Lem_wt_fn_extra_property_2}. Then:

\begin{enumerate}
\item $\mathcal{ML}_{1}=I$.

\item $\left\vert \mathcal{L}_{1}\mathcal{M}u-u\right\vert _{1/w,0}=0$ when
$u\in X_{1/w}^{0}$.

\item $\mathcal{M}$ and $\mathcal{L}_{1}$ are onto in the seminorm sense.

\item $\mathcal{M}$ and $\mathcal{L}_{1}$ are adjoints.
\end{enumerate}
\end{theorem}

\begin{proof}
\textbf{Part 1} Now $\mathcal{M}u=\left(  \frac{u_{F}}{\sqrt{w}}\right)
^{\vee}$ for $u\in X_{1/w}^{0}$ and $\mathcal{L}_{1}g=\left(  \sqrt
{w}\widehat{g}\right)  ^{\vee}$ for $g\in L^{2}$. it follows immediately that
$\mathcal{ML}_{1}=I$.\medskip

\textbf{Part 2} That $\mathcal{L}_{1}\mathcal{M}u-u\in\left(  S_{\mathcal{A}%
}^{\prime}\right)  ^{\vee}$ when $u\in X_{1/w}^{0}$ is proved by showing
$\left\vert \mathcal{L}_{1}\mathcal{M}u-u\right\vert _{1/w,0}=0$. In fact,
$\widehat{\mathcal{L}_{1}\mathcal{M}u}=\sqrt{w}\widehat{\mathcal{M}u}=u_{F}$
so that $\left\vert \mathcal{L}_{1}\mathcal{M}u-u\right\vert _{1/w,0}%
=0$.\medskip

\textbf{Part 3} $\mathcal{ML}_{1}=I$ implies $\mathcal{M}$ is onto. That
$\left\vert \mathcal{L}_{1}\mathcal{M}u-u\right\vert _{1/w,0}=0$ when $u\in
X_{1/w}^{0}$ implies $\mathcal{L}_{1}$ is onto.\medskip

\textbf{Part 4} Regarding adjointness, suppose $u\in X_{1/w}^{0}$ and $g\in
L^{2}$. Then from Theorem \ref{Thm_L1_propert} $\left(  \mathcal{L}%
_{1}g\right)  _{F}=\sqrt{w}\widehat{g}$ and from Definition \ref{Def_op_M},
$\widehat{\mathcal{M}u}=\frac{u_{F}}{\sqrt{w}}$. Thus
\[
\left\langle \mathcal{L}_{1}g,u\right\rangle _{1/w,0}=\int\frac{1}{w}\left(
\mathcal{L}_{1}g\right)  _{F}\overline{u_{F}}=\int\frac{1}{w}\sqrt
{w}\widehat{g}\overline{u_{F}}=\int\widehat{g}\frac{\overline{u_{F}}}{\sqrt
{w}}=\int\widehat{g}\overline{\widehat{\mathcal{M}u}}=\left(  \widehat{g}%
,\widehat{\mathcal{M}u}\right)  _{2}=\left(  g,\mathcal{M}u\right)  _{2}.
\]

\end{proof}

Since $L^{2}$ is complete, the mappings of the previous theorem easily yield
the following important result:

\begin{corollary}
\label{Cor_X1/w,o_semiHilb_opL1}Suppose the weight function $w$ has property
\ref{a1.046}. Then, in general, $X_{1/w}^{0}$ is a semi-Hilbert space.
However, $X_{1/w}^{0}$ is a Hilbert space iff $\mathcal{A}$ is empty.
\end{corollary}

\begin{proof}
By part 2 of Theorem \ref{Thm_L1_propert} $\mathcal{M}$ is isometric. Hence if
$\left\{  u_{k}\right\}  $ is Cauchy in $X_{1/w}^{0}$ then $\left\{
\mathcal{M}u_{k}\right\}  $ is Cauchy in $L^{2}$ and so $\mathcal{M}%
u_{k}\rightarrow u$ for some $u\in L^{2}$ since $L^{2}$ is complete. From
Definition \ref{Def_Xo1/w} the seminorm for $X_{1/w}^{0}$ has null space
$\left(  S_{\mathcal{A}}^{\prime}\right)  ^{\vee}$. Hence by Theorem
\ref{Thm_M_L1_propert} $\mathcal{L}_{2}\mathcal{M}u_{k}=u_{k}\rightarrow
\mathcal{L}_{2}u\in X_{1/w}^{0}$, as required.

That $X_{1/w}^{0}$ is a Hilbert space iff $\mathcal{A}$ is empty was stated in
Definition \ref{Def_Xo1/w}.
\end{proof}

Not all the weight functions discussed in Subsections
\ref{SbSect_wt_func_examples}, \ref{SbSect_wt_fn_examples_2} and Section
\ref{Sect_wt_fn_central_diff} satisfy property \ref{a1.046}. In fact we have
the following results:

\begin{theorem}
\label{Thm_wt_fn_satisfy_property_1}Regarding the property \ref{a1.046}:

\begin{enumerate}
\item The radial \textbf{Gaussian} weight function \ref{1.044} does
\textbf{not satisfy} property \ref{a1.046}.

\item The radial \textbf{Thin-plate spline} weight functions \ref{1.19} do
\textbf{not} \textbf{satisfy} property \ref{a1.046}.

\item The radial \textbf{Sobolev spline} weight functions \ref{1.030}
\textbf{satisfy} property \ref{a1.046}.

\item The \textbf{univariate} \textbf{central difference} weight functions
\ref{a962} of Subsection \ref{SbSect_cent_Motivation} \textbf{satisfy} the
condition \ref{a1.046}.

\item When $n=l$ the tensor product \textbf{multivariate} \textbf{central
difference} weight functions of Corollary \ref{Cor_cdiffwt_W02_W03}
\textbf{satisfy} the condition \ref{a1.046}.

When $n<l$ and \ref{a959} holds the tensor product \textbf{multivariate}
\textbf{central difference} weight functions of Corollary
\ref{Cor_cdiffwt_W02_W03} do \textbf{not satisfy} the condition \ref{a1.046}.

\item The tensor product \textbf{extended B-spline} weight functions of
Theorem \ref{Thm_ex_splin_wt_fn_properties} do \textbf{not satisfy} condition
\ref{a1.046}.
\end{enumerate}
\end{theorem}

\begin{proof}
Recall that property \ref{a1.046} is%
\[
w\left\vert \cdot\right\vert ^{2m}\in L_{loc}^{1}\text{ }and\text{ }%
\int\limits_{\left\vert \cdot\right\vert \geq r}\frac{w}{\left\vert
\cdot\right\vert ^{2\tau}}<\infty\text{ }for\text{ }some\text{ }\tau,r\geq0.
\]
\medskip

\textbf{Part 1} Easy.\medskip

\textbf{Part 2} Use the upper bound of Theorem \ref{Thm_bnds_modif_MacDonald}%
.\medskip

\textbf{Part 3} Since $w\left(  \xi\right)  =\left(  1+\left\vert
\xi\right\vert ^{2}\right)  ^{\nu}$ for $\nu>d/2$, choose $\tau>\nu
+d/2$.\medskip

\textbf{Part 4} In the univariate case just choose $m=l-n$.

Now suppose $n=l$ and $d>1$. Then \ref{a958} imply that for $t\in
\mathbb{R}^{1}$ and some $r_{0},c_{1},c_{2}>0$
\[
w\left(  t\right)  \leq\left\{
\begin{array}
[c]{ll}%
c_{2}, & \left\vert t\right\vert \leq r_{0},\\
2c_{1}t^{2n}, & \left\vert t\right\vert >r_{0}.
\end{array}
\right.
\]

Clearly $w$ is bounded on any compact set so $w\left\vert \cdot\right\vert
^{2m}\in L_{loc}^{1}$ when $m=0$. Also, if $r>0$
\[
\int\limits_{\left\vert \xi\right\vert \geq r}\frac{w\left(  \xi\right)
}{\left\vert \xi\right\vert ^{2\tau}}d\xi\leq\int\limits_{\left\vert
\xi\right\vert \geq r}\frac{%
{\textstyle\prod\limits_{k=1}^{d}}
\left(  c_{2}+2c_{1}\xi_{k}^{2n}\right)  }{\left\vert \xi\right\vert ^{2\tau}%
}d\xi\leq\int\limits_{\left\vert \xi\right\vert \geq r}\frac{%
{\textstyle\prod\limits_{k=1}^{d}}
\left(  c_{2}+2c_{1}\left\vert \xi\right\vert ^{2n}\right)  }{\left\vert
\xi\right\vert ^{2\tau}}d\xi\leq\int\limits_{\left\vert \xi\right\vert \geq
r}\frac{\left(  c_{2}+2c_{1}\left\vert \xi\right\vert ^{2n}\right)  ^{d}%
}{\left\vert \xi\right\vert ^{2\tau}}d\xi,
\]

which exists iff $2\tau>\left(  2n+1\right)  d$.

Finally, suppose $n<l$ and \ref{a959} holds. Then by Corollary
\ref{Cor_cdiffwt_bnd_on_wt_fn} for any $r>0$%
\[
\frac{c_{r}}{t^{2\left(  l-n\right)  }}\leq w_{1}\left(  t\right)
,\quad\left\vert t\right\vert \leq r,
\]

so that for all $k$%
\[
\frac{\left(  c_{r}\right)  ^{d}}{\xi^{2\left(  l-n\right)  \mathbf{1}}}\leq
w\left(  \xi\right)  ,\quad\left\vert \xi_{k}\right\vert \leq r.
\]

Thus if $0<\rho<\sqrt{d}r$%
\[
\int\limits_{\left\vert \xi\right\vert \leq\rho}w\left(  \xi\right)
\left\vert \xi\right\vert ^{2m}d\xi\geq r^{2r}\int\limits_{\left\vert
\xi\right\vert \leq\rho}w\left(  \xi\right)  d\xi\geq r^{2r}\left(
c_{r}\right)  ^{d}\int\limits_{\left\vert \xi\right\vert \leq\rho}\frac{d\xi
}{\xi^{2\left(  l-n\right)  \mathbf{1}}},
\]

and the last integral diverges.\medskip

\textbf{Part 5} An easy proof.
\end{proof}

\begin{remark}
The case where $n<l$ and \ref{a959} does not hold is not considered in the
previous theorem because I have not derived any relevant results in Section
\ref{Sect_wt_fn_central_diff}.
\end{remark}

Part 2 of Lemma \ref{Lem_wt_fn_extra_property_2} allowed the operator
$\mathcal{L}_{1}$ of Definition \ref{Def_op_L1}\ to be defined. Part 2 of
Lemma \ref{Lem_wt_fn_extra_property_2} was proved by assuming inequality
\ref{a1.046}. We will finish this section by giving alternative but closely
related conditions which imply part 2 of of Lemma
\ref{Lem_wt_fn_extra_property_2}.

\begin{theorem}
\label{Thm_Som_and_Sw2}Suppose the weight function $w$ satisfies the condition%
\begin{equation}
\left(  \int w\left\vert \psi\right\vert ^{2}\right)  ^{1/2}\leq\left\vert
\psi\right\vert _{\sigma},\quad\psi\in S_{\emptyset,m},\label{a1.31}%
\end{equation}

where $\left\vert \cdot\right\vert _{\sigma}$ is a positive, linear
combination of the seminorms in Definition \ref{Def_Distributions} which
define the topology of the tempered distributions $S$.

Then $w\in L_{loc}^{1}\left(  \mathbb{R}^{d}\setminus0\right)  \cap
L_{loc}^{1}\left(  \mathbb{R}^{d}\setminus\mathcal{A}\right)  $ and part 2 of
Lemma \ref{Lem_wt_fn_extra_property_2} holds i.e. $g\in L^{2}$ implies
$\sqrt{w}g\in S_{\emptyset,m}^{\prime}\cap L_{loc}^{1}\left(  \mathbb{R}%
^{d}\setminus0\right)  \cap L_{loc}^{1}\left(  \mathbb{R}^{d}\setminus
\mathcal{A}\right)  $.
\end{theorem}

\begin{proof}
Suppose $K$ is compact and $0\notin K$. We can choose $\psi\in S_{\emptyset
,m}$ such that $\psi=1$ on $K$. Then $\int_{K}w=\int_{K}w\left\vert
\psi\right\vert ^{2}\frac{1}{\left\vert \psi\right\vert ^{2}}=\int%
_{K}w\left\vert \psi\right\vert ^{2}\leq\left\vert \psi\right\vert _{\sigma
}^{2}$ and $w\in L_{loc}^{1}\left(  \mathbb{R}^{d}\setminus0\right)  $. Since
$w$ is continuous on $\mathbb{R}^{d}\setminus\mathcal{A}$ it follows directly
that $w\in L_{loc}^{1}\left(  \mathbb{R}^{d}\setminus\mathcal{A}\right)  $.

Examination of the proof of part 2 of Lemma \ref{Lem_wt_fn_extra_property_2}
now shows that $\sqrt{w}g\in S_{\emptyset,m}^{\prime}$.

Further, $\int_{K}\sqrt{w}g\leq\int_{K}\sqrt{w}\psi g\leq\left(  \int
w\left\vert \psi\right\vert ^{2}\right)  ^{1/2}\left\Vert g\right\Vert
_{2}<\infty$ and if $K^{\prime}\subset\mathbb{R}^{d}\setminus\mathcal{A}$ is
compact, $w$ is continuous on $K^{\prime}$ and $\int_{K^{\prime}}\sqrt{w}%
g\leq\left(  \int_{K^{\prime}}w\right)  ^{1/2}\left\Vert g\right\Vert
_{2}<\infty$.
\end{proof}

\subsection{The operator $\mathcal{L}_{2}:L^{2}\rightarrow X_{1/w}^{0}$}

In Theorem \ref{Thm_wt_fn_satisfy_property_1} of the previous section it was
shown that several of the weight functions do not satisfy property
\ref{a1.046}. We now tackle these deficiencies by replacing $S^{\prime}$ by a
larger space $S_{w,0}^{\dag}$, where $S_{w,0}$ is a special subspace of $S$
generated by $w$ which will be endowed with a topology generated by $w$. If
the topology on $S_{w,0}$ is the subspace topology induced by $S$ then a
generalization of the mapping $\mathcal{L}_{1}:L^{2}\rightarrow X_{1/w}^{0}$,
denoted by $\mathcal{L}_{2}$, will be constructed and used to prove that
$X_{1/w}^{0}$ is semi-Hilbert space. In Subsection
\ref{SbSect_ext_Nat_splin_wt_fn_Swo} it will be shown that the topologies are
equivalent for weight functions that are radial-like near the origin, the
extended B-spline weight functions and the central difference weight
functions. For several classes of weight function the topologies are not
equivalent and the space $X_{1/w}^{0}$ must be extended. This will be done in
Section \ref{Sect_generalize2_tildXo1/w_T}.

\begin{definition}
\label{Def_Sw2_and_fnal(Sw2)}\textbf{The spaces} $S_{w,0}$\textbf{, }%
$S_{w,0}^{\prime}$\textbf{,} $W_{S;0}$\textbf{\ and }$S_{w,0}^{\dag}$

\begin{enumerate}
\item \fbox{Property $S_{w,0}$} Suppose $w\in W01$ and define the linear
subspace%
\[
S_{w,0}=\left\{  \phi\in S:\int w\left\vert \phi\right\vert ^{2}%
<\infty\right\}  .
\]

We endow $S_{w,0}$ with the subspace topology induced by $S$ (Lemma
\ref{Lem_convex_tls_extend_2}) and observe that $S_{w,0}$ \textbf{is not
empty} because $w\in C^{\left(  0\right)  }\left(  \mathbb{R}^{d}%
\setminus\mathcal{A}\right)  $ implies $C_{0}^{\infty}\left(  \mathbb{R}%
^{d}\setminus\mathcal{A}\right)  \subset S_{w,0}$.

\item \fbox{Property $S_{w,0}^{\prime}$} denotes the bounded linear
functionals on $S_{w,0}$.

\item \fbox{Property $W_{S;0}$} We say $w\in W_{S;0}$ if there exists some
positive, linear combination $\left\vert \cdot\right\vert _{\sigma}$ of
seminorms used to define the topology on $S$ such that the \textbf{non-linear}
functional $\left(  \int w\left\vert \phi\right\vert ^{2}\right)  ^{1/2}$
satisfies%
\begin{equation}
\left(  \int w\left\vert \phi\right\vert ^{2}\right)  ^{1/2}\leq\left\vert
\phi\right\vert _{\sigma},\quad\phi\in S_{w,0}.\label{a1.8}%
\end{equation}

\item \fbox{Property $S_{w,0}^{\dag}$} We say the bilinear functional
$\Psi\left(  \phi,\psi\right)  \rightarrow\mathbb{C}$ belongs to
$S_{w,0}^{\dag} $ w.r.t. $S$ if
\[
\left\vert \Psi\left(  \phi,\psi\right)  \right\vert \leq\left\vert
\phi\right\vert _{\sigma}\left\vert \psi\right\vert _{\sigma},\quad\phi
,\psi\in S_{w,0},
\]

where $\left\vert \cdot\right\vert _{\sigma}$ is some positive linear
combination of the seminorms used to define the topology on $S$.

Observe that the bilinear functional $\int w\phi\overline{\psi}$ is a member
of $S_{w,0}^{\dag}$.
\end{enumerate}
\end{definition}

\begin{remark}
Regarding our examples of $S_{w,0}$:

\begin{enumerate}
\item For the \textbf{Sobolev splines} of Example \ref{Ex_Sobolev_splin_wt}
$S_{w,0}=S$ (Theorem \ref{Thm_Sw2_eq_Son}).

\item For the weight functions of Example \ref{Ex_wt_func_|x|^m_L1loc},
$S_{w,0}=S_{\emptyset,n}$ where $n$ is the smallest integer such that
$n>m-\left(  \frac{d}{2}-s\right)  $ and $n\geq0$ (Theorem
\ref{Thm_Sw2_eq_Son}).

\item For the \textbf{extended B-spline weight} functions of Theorem
\ref{Thm_ex_splin_wt_fn_properties}, $S_{w,0}=\left\{  \phi\in
S:x^{n\mathbf{1}}\phi\in S_{\otimes,l\mathbf{1}}\left(  \pi\mathbb{Z}%
^{d}\right)  \right\}  $ (Theorem \ref{Thm_Sw2_extsplin_wt}).

\item For the \textbf{tensor product central difference weight functions} of
Definition \ref{Def_central_diff_wt_func} satisfying \ref{a959}:
$S_{w,0}=S_{\otimes,l-n}\left(  \mathbf{0}\right)  $ (Theorem
\ref{Thm_Sw2.eq.Sx,l_minus_n}).

\item ??

\item ??
\end{enumerate}
\end{remark}

\begin{remark}
\label{Rem_Def_Sw2_and_fnal(Sw2)}That
\begin{equation}
S_{w,0}^{\bot}\subset S_{\mathcal{A}}^{\prime},\label{a2.7}%
\end{equation}

where $S_{w,0}^{\bot}$ is the subspace of annihilators of $S_{w,0}$,\ follows
immediately from the fact that $C_{0}^{\infty}\left(  \mathbb{R}^{d}%
\setminus\mathcal{A}\right)  \subset S_{w,0}$.
\end{remark}

\begin{definition}
\label{Def_op_L2}\textbf{The operator} $\mathcal{L}_{2}:L^{2}\rightarrow
S^{\prime}$

Suppose $w\in W01\cap W_{S;0}$. Then for $\phi\in S_{w,0}$ and $g\in L^{2}$%
\[
\left\vert \int\sqrt{w}\widehat{g}\phi\right\vert =\left\vert \int\left(
\sqrt{w}\phi\right)  \widehat{g}\right\vert \leq\left\Vert g\right\Vert
_{2}\left(  \int w\left\vert \phi\right\vert ^{2}\right)  ^{1/2}\leq\left\Vert
g\right\Vert _{2}\left\vert \phi\right\vert _{S},
\]

so that $\sqrt{w}\widehat{g}\in S_{w,0}^{\prime}$. Using Lemma
\ref{Lem_convex_tls_extend_2} we can then extend $\sqrt{w}\widehat{g}$ to $S$
as a member of $S^{\prime}$. Such an extension is unique up to a member of the
annihilator $S_{w,0}^{\bot}$ of $S_{w,0}$.

Denote such an extension by $\left(  \sqrt{w}\widehat{g}\right)  ^{e}\in
S^{\prime}$.

The class of operators $\mathcal{L}_{2}:L^{2}\rightarrow S^{\prime}$ can now
be defined by%
\begin{equation}
\mathcal{L}_{2}g=\left(  \left(  \sqrt{w}\widehat{g}\right)  ^{e}\right)
^{\vee},\quad g\in L^{2}.\label{a1.020}%
\end{equation}

In general $\mathcal{L}_{2}$ is not unique.
\end{definition}

The operator $\mathcal{L}_{2}$ has the following properties:

\begin{theorem}
\label{Thm_L2_propert}\textbf{Properties of }$\mathcal{L}_{2}$:

\begin{enumerate}
\item If $g\in L^{2}$ then $\left(  \mathcal{L}_{2}g\right)  _{F}=\sqrt
{w}\widehat{g}\in L_{loc}^{1}\left(  \mathbb{R}^{d}\setminus\mathcal{A}%
\right)  $.

\item $\mathcal{L}_{2}:L^{2}\rightarrow X_{1/w}^{0}$.

\item $\mathcal{L}_{2}$ is an isometry and 1-1.

\item $\mathcal{L}_{2}$ is linear in the seminorm sense i.e.%
\[
\left\vert \mathcal{L}_{2}\left(  \lambda_{1}g_{1}+\lambda_{2}g_{2}\right)
-\mathcal{L}_{2}\left(  \lambda_{1}g_{1}\right)  -\mathcal{L}_{2}\left(
\lambda_{2}g_{2}\right)  \right\vert _{1/w,0}=0.
\]

\item If $\tau_{a}$ is the translation operator $\tau_{a}u=u\left(
\cdot-a\right)  $ then $\tau_{a}\mathcal{L}_{2}=\mathcal{L}_{2}\tau_{a}$ in
the seminorm sense.
\end{enumerate}
\end{theorem}

\begin{proof}
\textbf{Part 1} Suppose $g\in L^{2}$. From the definition of $\mathcal{L}_{2}
$, $\widehat{\mathcal{L}_{2}g}=\sqrt{w}\widehat{g}$ on $S_{w,0}$. Next observe
that because $w\in C^{\left(  0\right)  }\left(  \mathbb{R}^{d}\setminus
\mathcal{A}\right)  $ and $C_{0}^{\infty}\left(  \mathbb{R}^{d}\setminus
\mathcal{A}\right)  \subset S_{w,0}$ it follows that
\[
\left[  \widehat{\mathcal{L}_{2}g},\phi\right]  =\left[  \sqrt{w}%
\widehat{g},\phi\right]  ,\quad\phi\in C_{0}^{\infty}\left(  \mathbb{R}%
^{d}\setminus\mathcal{A}\right)  .
\]

If we can show $\sqrt{w}\widehat{g}\in L_{loc}^{1}\left(  \mathbb{R}%
^{d}\setminus\mathcal{A}\right)  $ it follows that $\left(  \mathcal{L}%
_{2}g\right)  _{F}=\sqrt{w}\widehat{g}\in L_{loc}^{1}\left(  \mathbb{R}%
^{d}\setminus\mathcal{A}\right)  $. But if $K\subset\mathbb{R}^{d}%
\setminus\mathcal{A}$ is compact, $w\in C_{0}^{\infty}\left(  \mathbb{R}%
^{d}\setminus\mathcal{A}\right)  $ implies
\[
\int_{K}\sqrt{w}\left\vert \widehat{g}\right\vert \leq\left(  \int%
_{K}w\right)  ^{1/2}\left\Vert g\right\Vert _{2}\leq\max_{K}\left(  w\right)
\left(  \int_{K}1\right)  \left\Vert g\right\Vert _{2}.
\]

Thus $\widehat{\mathcal{L}_{2}g}\in L_{loc}^{1}\left(  \mathbb{R}^{d}%
\setminus\mathcal{A}\right)  $ and $\left(  \mathcal{L}_{2}g\right)
_{F}=\sqrt{w}\widehat{g}$.\medskip

\textbf{Part 2} $\int\frac{\left\vert \left(  \mathcal{L}_{2}g\right)
_{F}\right\vert ^{2}}{w}=\left\Vert \widehat{g}\right\Vert _{2}=\left\Vert
g\right\Vert _{2}$ so $\mathcal{L}_{2}:L^{2}\rightarrow X_{1/w}^{0}$.\medskip

\textbf{Part 3} $\mathcal{L}_{2}$ is an isometry since $\left\vert
\mathcal{L}_{2}g\right\vert _{1/w,0}=\left(  \int\frac{\left\vert \left(
\mathcal{L}_{2}g\right)  _{F}\right\vert ^{2}}{w}\right)  ^{1/2}=\left\Vert
g\right\Vert _{2}$. Thus $\mathcal{L}_{2}:L^{2}\rightarrow X_{1/w}^{0}$ is an
isometry and clearly $\mathcal{L}_{2}g=0$ implies $g=0$.\medskip

\textbf{Part 4} Linearity in the seminorm sense is true since%
\begin{align*}
& \left\vert \mathcal{L}_{2}\left(  \lambda_{1}g_{1}+\lambda_{2}g_{2}\right)
-\mathcal{L}_{2}\left(  \lambda_{1}g_{1}\right)  -\mathcal{L}_{2}\left(
\lambda_{2}g_{2}\right)  \right\vert _{1/w,0}\\
& =\int\frac{\left\vert \left(  \mathcal{L}_{2}\left(  \lambda_{1}%
g_{1}+\lambda_{2}g_{2}\right)  -\mathcal{L}_{2}\left(  \lambda_{1}%
g_{1}\right)  -\mathcal{L}_{2}\left(  \lambda_{2}g_{2}\right)  \right)
_{F}\right\vert ^{2}}{w}\\
& =\int\frac{\left\vert \sqrt{w}\left(  \lambda_{1}g_{1}+\lambda_{2}%
g_{2}\right)  ^{\wedge}-\sqrt{w}\left(  \lambda_{1}g_{1}\right)  ^{\wedge
}-\sqrt{w}\left(  \lambda_{2}g_{2}\right)  ^{\wedge}\right\vert ^{2}}{w}\\
& =0.
\end{align*}
\medskip

\textbf{Part 5}
\[
\left\vert \tau_{a}\mathcal{L}_{2}g-\mathcal{L}_{2}\tau_{a}g\right\vert
_{1/w,0}^{2}=\int\frac{\left\vert \left(  \tau_{a}\mathcal{L}_{2}%
g-\mathcal{L}_{2}\tau_{a}g\right)  _{F}\right\vert ^{2}}{w}=\int%
\frac{\left\vert \left(  \tau_{a}\mathcal{L}_{2}g\right)  _{F}-\left(
\mathcal{L}_{2}\tau_{a}g\right)  _{F}\right\vert ^{2}}{w}.
\]

But from part 1, $\left(  \mathcal{L}_{2}g\right)  _{F}=\sqrt{w}\widehat{g}$
on $\mathbb{R}^{d}\setminus\mathcal{A}$. Hence on $\mathbb{R}^{d}%
\setminus\mathcal{A}$%
\[
\left(  \mathcal{L}_{2}\tau_{a}g\right)  _{F}=\sqrt{w}\widehat{\tau_{a}%
g}=e^{-ia\xi}\sqrt{w}\widehat{g}=e^{-ia\xi}\left(  \mathcal{L}_{2}g\right)
_{F},
\]

and since $\left(  \tau_{a}\mathcal{L}_{2}g\right)  ^{\wedge}=e^{-ia\xi
}\left(  \mathcal{L}_{2}g\right)  ^{\wedge}$, it follows that%
\[
\left(  \tau_{a}\mathcal{L}_{2}g\right)  _{F}=e^{-ia\xi}\left(  \mathcal{L}%
_{2}g\right)  _{F},
\]

and therefore $\left(  \tau_{a}\mathcal{L}_{2}g\right)  _{F}=\left(
\mathcal{L}_{2}\tau_{a}g\right)  _{F}$.
\end{proof}

The following theorem indicates how the operators $\mathcal{L}_{2}%
:L^{2}\rightarrow X_{1/w}^{0}$ and $\mathcal{M}:X_{1/w}^{0}\rightarrow L^{2}$ interact.

\begin{theorem}
\label{Thm_M_L2_propert}Suppose the weight function $w$ is a member of
$W_{S;0}$. Then:

\begin{enumerate}
\item $\mathcal{ML}_{2}=I$ on $L^{2}$.

\item $\mathcal{L}_{2}\mathcal{M}u-u\in\left(  S_{\mathcal{A}}^{\prime
}\right)  ^{\vee}$ when $u\in X_{1/w}^{0}$ i.e. $\left\vert \mathcal{L}%
_{2}\mathcal{M}u-u\right\vert _{1/w,0}=0$.

\item $\mathcal{M}:X_{1/w}^{0}\rightarrow L^{2}$ is onto.

\item $\mathcal{L}_{2}:L^{2}\rightarrow X_{1/w}^{0}$ is onto in the seminorm sense.

\item $\mathcal{M}$ and $\mathcal{L}_{2}$ are adjoints.
\end{enumerate}
\end{theorem}

\begin{proof}
\textbf{Part 1} By definition $\widehat{\mathcal{M}u}=\frac{u_{F}}{\sqrt{w}%
}\in L^{2}$ when $u\in X_{1/w}^{0}$. From the proof of Theorem
\ref{Thm_L2_propert} $\left(  \mathcal{L}_{2}f\right)  _{F}=\sqrt{w}f$. Thus
from the definition of $\mathcal{M}$, for $f\in L^{2}$%
\[
\left(  \mathcal{ML}_{2}f\right)  ^{\wedge}=\frac{\left(  \mathcal{L}%
_{2}f\right)  _{F}}{\sqrt{w}}=\frac{\sqrt{w}\widehat{f}}{\sqrt{w}}%
=\widehat{f},
\]

and so $\mathcal{ML}_{2}=I$ on $L^{2}$.\medskip

\textbf{Part 2} If $u\in X_{1/w}^{0}$ then $\mathcal{L}_{2}\mathcal{M}u\in
X_{1/w}^{0}$ and $\left(  \mathcal{L}_{2}\mathcal{M}u\right)  _{F}=\sqrt
{w}\widehat{\mathcal{M}u}=u_{F}$. Thus%
\[
\left\vert \mathcal{L}_{2}\mathcal{M}u-u\right\vert _{1/w,0}^{2}=\int%
\frac{\left\vert \left(  \mathcal{L}_{2}\mathcal{M}u-u\right)  _{F}\right\vert
^{2}}{w}=0,
\]

and $\mathcal{L}_{2}\mathcal{M}u-u\in\left(  S_{\mathcal{A}}^{\prime}\right)
^{\vee}$.\medskip

\textbf{Parts 3 and 4} The equation $\mathcal{ML}_{2}=I$ implies $\mathcal{M}
$ is onto and $\mathcal{L}_{2}\mathcal{M}u-u\in\left(  S_{\mathcal{A}}%
^{\prime}\right)  ^{\vee}$ when $u\in X_{1/w}^{0}$ implies that $\mathcal{L}%
_{2}$ are onto.\medskip

\textbf{Part 5} Regarding adjointness, suppose $u\in X_{1/w}^{0}$ and $g\in
L^{2}$. Then from Theorem \ref{Thm_L2_propert} $\left(  \mathcal{L}%
_{1}g\right)  _{F}=\sqrt{w}\widehat{g}$ and from Definition \ref{Def_op_M},
$\widehat{\mathcal{M}u}=\frac{u_{F}}{\sqrt{w}}$. Thus
\[
\left(  \mathcal{L}_{2}g,u\right)  _{1/w,0}=\int\frac{1}{w}\left(
\mathcal{L}_{2}g\right)  _{F}\overline{u_{F}}=\int\frac{1}{w}\sqrt
{w}\widehat{g}\overline{u_{F}}=\int\widehat{g}\frac{\overline{u_{F}}}{\sqrt
{w}}=\int\widehat{g}\overline{\widehat{\mathcal{M}u}}=\left(  \widehat{g}%
,\widehat{\mathcal{M}u}\right)  _{2}=\left(  g,\mathcal{M}u\right)  _{2}.
\]

\end{proof}

Since $L^{2}$ is complete, the mappings of the previous theorem will yield the
following important result:

\begin{corollary}
\label{Cor_X1/w,o_semiHilb_opL2}If the weight function $w\in W_{S;0}$ then in
general $X_{1/w}^{0}$ is a semi-Hilbert space. Indeed, $X_{1/w}^{0}$ is a
Hilbert space iff $\mathcal{A}$ is empty.
\end{corollary}

\begin{proof}
By Theorem \ref{Thm_L2_propert} $\mathcal{M}$ is isometric. Hence if $\left\{
u_{k}\right\}  $ is Cauchy in $X_{1/w}^{0}$ then $\left\{  \mathcal{M}%
u_{k}\right\}  $ is Cauchy in $L^{2}$ and so $\mathcal{M}u_{k}\rightarrow u$
for some $u\in L^{2}$ since $L^{2}$ is complete. From part 2 Lemma
\ref{Lem_Xoinvw} the seminorm for $X_{1/w}^{0}$ has null space $\left(
S_{\mathcal{A}}^{\prime}\right)  ^{\vee}$. Hence by Theorem
\ref{Thm_M_L2_propert} $\mathcal{L}_{2}\mathcal{M}u_{k}=u_{k}\rightarrow
\mathcal{L}_{2}u\in X_{1/w}^{0}$, as required.

Clearly $\mathcal{A}$ is empty iff $S_{\mathcal{A}}^{\prime}=\left\{
0\right\}  $ iff $\left(  S_{\mathcal{A}}^{\prime}\right)  ^{\vee}=\left\{
0\right\}  $.
\end{proof}

Combining Corollary \ref{Cor_X1/w,o_semiHilb_opL1} and Corollary
\ref{Cor_X1/w,o_semiHilb_opL2} we have:

\begin{theorem}
\label{Def_Xo1/w_Hilbert}Suppose the weight function $w$ has property
\ref{a1.046} or is a member of $W_{S;0}$.

Then $X_{1/w}^{0}$ is a semi-Hilbert space or a Hilbert space. Further,
$X_{1/w}^{0}$ is a Hilbert space iff $\mathcal{A}$ is empty.
\end{theorem}

\subsection{Examples of weight functions in $W_{S;0}$ and the corresponding
spaces $S_{w,0}$\label{SbSect_ext_Nat_splin_wt_fn_Swo}}

We start by considering a class of radial basis functions which satisfy
\ref{a1.046} but which also lie in $W_{S;0}$. Then it is shown that the
\textbf{tensor product extended B-spline and central difference weight
functions are of type} $W_{S;0}$ and we characterize the spaces $S_{w,0}$.

\subsubsection{\protect\underline{Weight functions such that $w\left\vert
\cdot\right\vert ^{2m}\in L_{loc}^{1}$} \label{SbSbSect_wt_fn_rad_L1loc}}

\begin{definition}
\label{Def_radial_like_wt_fn}\textbf{Radial-like weight functions} We say that
a weight function is radial-like if $w\in W01$, $w\left\vert \cdot\right\vert
^{2m}\in L_{loc}^{1}$ for some integer $m\geq0$, and for some $\tau\geq0$,%
\begin{equation}
\int\limits_{\left\vert \cdot\right\vert \geq r}\frac{w}{\left\vert
\cdot\right\vert ^{2\tau}}<\infty.\label{a1.005}%
\end{equation}

This class of weight functions includes the Sobolev splines $w=\left(
1+\left\vert \cdot\right\vert ^{2}\right)  ^{\nu}$ of Example
\ref{Ex_Sobolev_splin_wt} and the weight functions introduced in Examples
\ref{Ex_wt_func_L1loc} and \ref{Ex_wt_func_|x|^m_L1loc}.
\end{definition}

We first consider separately the simplest case $m=0$. Here it turns out that
$S_{w,0}=S$. Indeed:

\begin{theorem}
\textbf{Radial-like weight functions with }$\mathbf{m=0}$. Suppose the weight
function $w$ is radial-like i.e. it satisfies Definition
\ref{Def_radial_like_wt_fn} for $m=0$. Then $S_{w,0}=S$ and $w\in W_{S;0}$.
Further%
\[
\int w\left\vert \phi\right\vert ^{2}\leq\left(  \int_{\left\vert
\cdot\right\vert \leq r}w\right)  \left\Vert \phi\right\Vert _{\infty}%
^{2}+\left\lceil \tau\right\rceil !\left(  \int_{\left\vert \cdot\right\vert
\geq r}\frac{w}{\left\vert \cdot\right\vert ^{2\tau}}\right)  \sum
\limits_{\left\vert \alpha\right\vert =\left\lceil \tau\right\rceil }\frac
{1}{\alpha!}\left\Vert \xi^{\alpha}\phi\right\Vert _{\infty}^{2}.
\]

\begin{proof}
If $\phi\in S$ then%
\begin{align*}
\int w\left\vert \phi\right\vert ^{2}  & =\int_{\left\vert \cdot\right\vert
\leq r}w\left\vert \phi\right\vert ^{2}+\int_{\left\vert \cdot\right\vert \geq
r}w\left\vert \phi\right\vert ^{2}\\
& \leq\left(  \int_{\left\vert \cdot\right\vert \leq r}w\right)  \left\Vert
\phi\right\Vert _{\infty}^{2}+\int_{\left\vert \cdot\right\vert \geq
r}w\left\vert \phi\right\vert ^{2}\\
& =\left(  \int_{\left\vert \cdot\right\vert \leq r}w\right)  \left\Vert
\phi\right\Vert _{\infty}^{2}+\int_{\left\vert \cdot\right\vert \geq r}%
\frac{w}{\left\vert \cdot\right\vert ^{2\tau}}\left\vert \cdot\right\vert
^{2\tau}\left\vert \phi\right\vert ^{2}\\
& \leq\left(  \int_{\left\vert \cdot\right\vert \leq r}w\right)  \left\Vert
\phi\right\Vert _{\infty}^{2}+\int_{\left\vert \cdot\right\vert \geq r}%
\frac{w}{\left\vert \cdot\right\vert ^{2\tau}}\left\vert \cdot\right\vert
^{2\left\lceil \tau\right\rceil }\left\vert \phi\right\vert ^{2}\\
& \leq\left(  \int_{\left\vert \cdot\right\vert \leq r}w\right)  \left\Vert
\phi\right\Vert _{\infty}^{2}+\left(  \int_{\left\vert \cdot\right\vert \geq
r}\frac{w}{\left\vert \cdot\right\vert ^{2\tau}}\right)  \left\Vert \left\vert
\cdot\right\vert ^{2\left\lceil \tau\right\rceil }\left\vert \phi\right\vert
^{2}\right\Vert _{\infty}\\
& =\left(  \int_{\left\vert \cdot\right\vert \leq r}w\right)  \left\Vert
\phi\right\Vert _{\infty}^{2}+\left\lceil \tau\right\rceil !\left(
\int_{\left\vert \cdot\right\vert \geq r}\frac{w}{\left\vert \cdot\right\vert
^{2\tau}}\right)  \left\Vert \frac{\left\vert \cdot\right\vert ^{2\left\lceil
\tau\right\rceil }}{\left\lceil \tau\right\rceil !}\left\vert \phi\right\vert
^{2}\right\Vert _{\infty}\\
& =\left(  \int_{\left\vert \cdot\right\vert \leq r}w\right)  \left\Vert
\phi\right\Vert _{\infty}^{2}+\left\lceil \tau\right\rceil !\left(
\int_{\left\vert \cdot\right\vert \geq r}\frac{w}{\left\vert \cdot\right\vert
^{2\tau}}\right)  \left\Vert \sum\limits_{\left\vert \alpha\right\vert
=\left\lceil \tau\right\rceil }\frac{\xi^{2\alpha}}{\alpha!}\left\vert
\phi\right\vert ^{2}\right\Vert _{\infty}\\
& \leq\left(  \int_{\left\vert \cdot\right\vert \leq r}w\right)  \left\Vert
\phi\right\Vert _{\infty}^{2}+\left\lceil \tau\right\rceil !\left(
\int_{\left\vert \cdot\right\vert \geq r}\frac{w}{\left\vert \cdot\right\vert
^{2\tau}}\right)  \sum\limits_{\left\vert \alpha\right\vert =\left\lceil
\tau\right\rceil }\frac{1}{\alpha!}\left\Vert \xi^{2\alpha}\left\vert
\phi\right\vert ^{2}\right\Vert _{\infty}\\
& =\left(  \int_{\left\vert \cdot\right\vert \leq r}w\right)  \left\Vert
\phi\right\Vert _{\infty}^{2}+\left\lceil \tau\right\rceil !\left(
\int_{\left\vert \cdot\right\vert \geq r}\frac{w}{\left\vert \cdot\right\vert
^{2\tau}}\right)  \sum\limits_{\left\vert \alpha\right\vert =\left\lceil
\tau\right\rceil }\frac{1}{\alpha!}\left\Vert \xi^{\alpha}\phi\right\Vert
_{\infty}^{2}.
\end{align*}

Hence $S_{w,0}=S$ and $w\in W_{S;0}$.
\end{proof}
\end{theorem}

We now consider the weight functions of Definition \ref{Def_radial_like_wt_fn}
when $m>0$. The next theorem illustrates the interplay between the pointwise
properties used to define property W01 and the integration conditions
\ref{a1.046} required by our Hilbert space techniques.

\begin{theorem}
\textbf{Radial-like weight functions with }$\mathbf{m>0}$. Suppose the weight
function $w$ satisfies Definition \ref{Def_radial_like_wt_fn} with $m\geq1$. Then:

\begin{enumerate}
\item $w\in L_{loc}^{1}\left(  \mathbb{R}^{d}\setminus0\right)  $.

\item $w\left\vert \cdot\right\vert ^{2m}\in L_{loc}^{1}$.

\item $w$ is a function in $S_{\emptyset,2m}^{\prime}$. Indeed $w\in
L_{loc}^{1}\left(  \mathbb{R}^{d}\setminus0\right)  $ and has polynomial
(slow)increase at infinity i.e. \ref{a1.005} holds for some $\tau\geq0$.
Further, $w$ can be extended (non-uniquely) to $S^{\prime}$.\medskip

\quad Regarding the set $\mathcal{A}$ used to define the weight property
W01:\medskip

\item If $w\in C^{\left(  0\right)  }\left(  B\left(  0;r\right)  \right)  $
and $0<2m+2s<d$ then $0$ is an isolated point of $\mathcal{A}$.

\item If $w\in C^{\left(  0\right)  }\left(  \overline{B\left(  0;r\right)
}\right)  $ and ??$2m+2s\leq0$?? then $\mathcal{A}\cap\overline{B\left(
0;r\right)  }$ is empty.
\end{enumerate}

\begin{proof}
?? FINISH!.
\end{proof}
\end{theorem}

\begin{remark}
Part 3 is analogous to the definition of a function in $S^{\prime}$ e.g.
Definition 2.8.3(a)\thinspace Vladimirov \cite{Vladimirov} which requires that
$w\in L_{loc}^{1}$ and has polynomial (slow)increase at infinity.
\end{remark}

Part 2 of the last theorem shows that the weight function $w$\ satisfies
Definition \ref{Def_radial_like_wt_fn} and so we can define the operator
$\mathcal{L}_{1}$\ without using the lemma below!

The next theorem will characterize the space $S_{w,0}$ in terms of the spaces
$S_{\emptyset,n}$. However, all the work will be done by the next lemma.

\begin{lemma}
\label{Lem_wt_radial_near_0}Suppose $\sigma\geq0$ and $\phi\in S\left(
\mathbb{R}^{d}\right)  $. Let $n$ be the smallest integer such that
$2n>2\sigma-d$ and $n\geq0$.

Then: $\int_{\left\vert \cdot\right\vert \leq1}\frac{\left\vert \phi
\right\vert ^{2}}{\left\vert \cdot\right\vert ^{2\sigma}}<\infty$ iff $\phi\in
S_{\emptyset,n}$.

Further, if $\int_{\left\vert \cdot\right\vert \leq1}\frac{\left\vert
\phi\right\vert ^{2}}{\left\vert \cdot\right\vert ^{2\sigma}}<\infty$ then%
\begin{equation}
\int_{\left\vert \cdot\right\vert \leq1}\frac{\left\vert \phi\right\vert ^{2}%
}{\left\vert \cdot\right\vert ^{2\sigma}}\leq k_{\sigma}\left(  \sum
_{\left\vert \beta\right\vert =n}\left\Vert D^{\beta}\phi\right\Vert _{\infty
}\right)  ^{2},\label{a1.41}%
\end{equation}

where
\begin{equation}
k_{\sigma}=\int_{\left\vert \cdot\right\vert \leq1}\frac{1}{\left\vert
\cdot\right\vert ^{2\sigma-2n}}=\frac{2}{2n-2\sigma+d+1}\frac{\left(
\sqrt{\pi}/2\right)  ^{d/2}}{\Gamma\left(  d/2\right)  }.\label{a1.46}%
\end{equation}

\end{lemma}

\begin{proof}
First suppose $\phi\in S_{\emptyset,n}$. Then by part 1 Theorem
\ref{Thm_product_of_Co,k_funcs}, $\left\vert \phi\right\vert ^{2}\in
S_{\emptyset,2n}$ and by \ref{a992} there exists a constant $C_{\phi,2n}\geq0$
such that $\left\vert \phi\left(  x\right)  \right\vert ^{2}\leq C_{\phi
,2n}\left\vert x\right\vert ^{2n}$ when $\left\vert x\right\vert \leq1$. Hence
$\int_{\left\vert \cdot\right\vert \leq1}\frac{\left\vert \phi\right\vert
^{2}}{\left\vert \cdot\right\vert ^{2\sigma}}\leq C_{\phi,2n}\int_{\left\vert
\cdot\right\vert \leq1}\frac{1}{\left\vert \cdot\right\vert ^{2\sigma-2n}%
}<\infty$ since we have assumed that $2\sigma-2n<d$.\medskip

To prove that $\int_{\left\vert \cdot\right\vert \leq1}\frac{\left\vert
\phi\right\vert ^{2}}{\left\vert \cdot\right\vert ^{2\sigma}}<\infty$ implies
$\phi\in S_{\emptyset,n}$ we will split the domain of $\sigma$ and consider
three cases:\medskip

\fbox{\textbf{Case 1} $0\leq2\sigma<d$} Here $n=0$. This condition means that
$\int_{\left\vert \cdot\right\vert \leq1}\frac{\left\vert \phi\right\vert
^{2}}{\left\vert \cdot\right\vert ^{2\sigma}}$ exists for all $\phi\in
S_{\emptyset,0}=S$. Further
\[
\int_{\left\vert \cdot\right\vert \leq1}\frac{\left\vert \phi\right\vert ^{2}%
}{\left\vert \cdot\right\vert ^{2\sigma}}\leq\left(  \int_{\left\vert
\cdot\right\vert \leq1}\frac{1}{\left\vert \cdot\right\vert ^{2\sigma}%
}\right)  \left\Vert \phi\right\Vert _{\infty}^{2},\qquad\phi\in
S_{\emptyset,0}.
\]
\medskip

\fbox{\textbf{Case 2} $d\leq2\sigma<d+2$} Here $0\leq2\sigma-d<2$ and $n=1$.
If $\phi\left(  0\right)  \neq0$ then $\int_{\left\vert \cdot\right\vert
\leq1}\frac{\left\vert \phi\right\vert ^{2}}{\left\vert \cdot\right\vert
^{2\sigma}}\geq\int_{\left\vert \cdot\right\vert \leq1}\frac{\left\vert
\phi\right\vert ^{2}}{\left\vert \cdot\right\vert ^{d}}$ which diverges. Thus
$\phi\left(  0\right)  =0$ and $\phi\in S_{\emptyset,1}$. Since we have proved
that $\phi\in S_{\emptyset,1}$ implies $\int_{\left\vert \cdot\right\vert
\leq1}\frac{\left\vert \phi\right\vert ^{2}}{\left\vert \cdot\right\vert
^{2\sigma}}<\infty$ and this case is true.\medskip

\fbox{\textbf{Case 3} $d+2\leq2\sigma$} Here $n\geq2$ and
\begin{equation}
2n>2\sigma-d\geq2n-2\geq2.\label{a1.43}%
\end{equation}

Set $\psi=\left\vert \phi\right\vert ^{2}$. Using a similar argument to Case 2
we see immediately that%
\begin{equation}
\phi\left(  0\right)  =0\label{a1.39}%
\end{equation}

Expanding $\psi$ using the Taylor series about the origin we have%
\[
\psi\left(  x\right)  =\sum_{\left\vert \beta\right\vert <2n}\frac{x^{\beta}%
}{\beta!}D^{\beta}\psi\left(  0\right)  +\left(  \mathcal{R}_{2n}\psi\right)
\left(  0,x\right)  ,
\]

where the remainder satisfies the estimate \ref{a1.37}:%
\[
\left\vert \left(  \mathcal{R}_{2n}\psi\right)  \left(  0,x\right)
\right\vert \leq\left\vert x\right\vert ^{2n}\sum_{\left\vert \beta\right\vert
=2n}\max_{y\in\left[  0,x\right]  }\left\vert \left(  D^{\beta}\psi\right)
(y)\right\vert \leq\left\vert x\right\vert ^{2n}\sum_{\left\vert
\beta\right\vert =2n}\left\Vert D^{\beta}\psi\right\Vert _{\infty}.
\]

Thus the definition of $n$ implies that%
\[
\int_{\left\vert \cdot\right\vert \leq1}\frac{\left\vert \left(
\mathcal{R}_{2n}\psi\right)  \left(  0,x\right)  \right\vert }{\left\vert
\cdot\right\vert ^{2\sigma}}\leq\sum_{\left\vert \beta\right\vert
=2n}\left\Vert D^{\beta}\psi\right\Vert _{\infty}\int_{\left\vert
\cdot\right\vert \leq1}\frac{\left\vert \cdot\right\vert ^{2n}}{\left\vert
\cdot\right\vert ^{2\sigma}}<\infty,
\]

and so the Taylor polynomial satisfies%
\[
\int_{\left\vert \cdot\right\vert \leq1}\frac{1}{\left\vert \cdot\right\vert
^{2\sigma}}\sum_{\left\vert \beta\right\vert <2n}\frac{x^{\beta}}{\beta
!}D^{\beta}\psi\left(  0\right)  dx<\infty.
\]

By considering the odd/even symmetries of the integrand we obtain%
\begin{align}
\int\limits_{\left\vert \cdot\right\vert \leq1}\frac{1}{\left\vert
\cdot\right\vert ^{2\sigma}}\sum_{\left\vert \beta\right\vert <2n}%
\frac{x^{\beta}}{\beta!}D^{\beta}\psi\left(  0\right)  dx  & =\lim
_{\varepsilon\rightarrow0^{+}}\sum_{\left\vert \beta\right\vert <2n}%
\frac{D^{\beta}\psi\left(  0\right)  }{\beta!}\int\limits_{\varepsilon
\leq\left\vert x\right\vert \leq1}\frac{x^{\beta}}{\left\vert x\right\vert
^{2\sigma}}dx\nonumber\\
& =2^{d}\lim_{\varepsilon\rightarrow0^{+}}\sum_{\left\vert \alpha\right\vert
<n}\frac{D^{2\alpha}\psi\left(  0\right)  }{\left(  2\alpha\right)  !}%
\int\limits_{\substack{\varepsilon\leq\left\vert x\right\vert \leq1 \\x\geq
0}}\frac{x^{2\alpha}}{\left\vert x\right\vert ^{2\sigma}}dx.\label{a1.32}%
\end{align}

The next step is to study these integrals using multivariate spherical polar
coordinates e.g.\thinspace Section 5.43 Adams \cite{Adams75}:%
\[
x=\left(  \rho,\phi_{1},\ldots,\phi_{d-1}\right)  ,
\]

where $\rho\geq0,$ $0\leq\phi_{i}\leq\pi/2$ and%
\begin{equation}
\left.
\begin{array}
[c]{rr}%
x_{1}= & \rho\sin\phi_{1}\sin\phi_{2}\ldots\sin\phi_{d-1},\\
x_{2}= & \rho\cos\phi_{1}\sin\phi_{2}\ldots\sin\phi_{d-1},\\
x_{3}= & \rho\cos\phi_{2}\ldots\sin\phi_{d-1},\\
\vdots & \\
x_{d}= & \rho\cos\phi_{d-1}.
\end{array}
\right\} \label{a1.42}%
\end{equation}

The volume element is%
\[
dx=\rho^{d}\prod\limits_{j=1}^{d-1}\left(  \sin\phi_{j}\right)  ^{j-1}d\rho
d\phi,
\]

where $d\phi=d\phi_{1}\ldots d\phi_{d-1}$. We now have%
\begin{align*}
x^{2\alpha}=\rho^{2\left\vert \alpha\right\vert }\left(  \sin\phi_{1}\right)
^{2\alpha_{1}}\left(  \cos\phi_{1}\right)  ^{2\alpha_{2}} &  \times\left(
\sin\phi_{2}\right)  ^{2\left(  \alpha_{1}+\alpha_{2}\right)  }\left(
\cos\phi_{2}\right)  ^{2\alpha_{3}}\times\\
&  \times\left(  \sin\phi_{3}\right)  ^{2\left(  \alpha_{1}+\alpha_{2}%
+\alpha_{3}\right)  }\left(  \cos\phi_{3}\right)  ^{2\alpha_{4}}\times\\
&  \vdots\\
&  \times\left(  \sin\phi_{d-1}\right)  ^{2\left(  \alpha_{1}+\ldots
+\alpha_{d-1}\right)  }\left(  \cos\phi_{d-1}\right)  ^{2\alpha_{d}},
\end{align*}

so that%
\begin{align*}
x^{2\alpha}dx=\rho^{2\left\vert \alpha\right\vert +d-1}\left(  \sin\phi
_{1}\right)  ^{2\alpha_{1}}\left(  \cos\phi_{1}\right)  ^{2\alpha_{2}} &
\times\left(  \sin\phi_{2}\right)  ^{1+2\left(  \alpha_{1}+\alpha_{2}\right)
}\left(  \cos\phi_{2}\right)  ^{2\alpha_{3}}\times\\
&  \times\left(  \sin\phi_{3}\right)  ^{2+2\left(  \alpha_{1}+\alpha
_{2}+\alpha_{3}\right)  }\left(  \cos\phi_{3}\right)  ^{2\alpha_{4}}\times\\
&  \vdots\\
&  \times\left(  \sin\phi_{d-1}\right)  ^{d-2+2\left(  \alpha_{1}%
+\ldots+\alpha_{d-1}\right)  }\left(  \cos\phi_{d-1}\right)  ^{2\alpha_{d}},
\end{align*}

and using the gamma function formula
\begin{equation}
\int_{0}^{\pi/2}\left(  \sin\phi\right)  ^{m}\left(  \cos\phi\right)
^{n}d\phi=\frac{\Gamma\left(  \frac{m+1}{2}\right)  \Gamma\left(  \frac
{n+1}{2}\right)  }{2\Gamma\left(  \frac{m+n}{2}+1\right)  },\quad
m,n=0,1,2,\ldots\label{a1.40}%
\end{equation}

we obtain separated variables%
\begin{align*}
\int\limits_{\substack{\varepsilon\leq\left\vert x\right\vert \leq1 \\x\geq
0}}\frac{x^{2\alpha}}{\left\vert x\right\vert ^{2\sigma}}dx  & =\int%
\limits_{\varepsilon}^{1}\frac{d\rho}{\rho^{2\sigma-2\left\vert \alpha
\right\vert -d+1}}\text{ }\int\limits_{0}^{\pi/2}\left(  \sin\phi_{1}\right)
^{2\alpha_{1}}\left(  \cos\phi_{1}\right)  ^{2\alpha_{2}}d\phi_{1}\times\\
& \qquad\times\int\limits_{0}^{\pi/2}\left(  \sin\phi_{2}\right)  ^{1+2\left(
\alpha_{1}+\alpha_{2}\right)  }\left(  \cos\phi_{2}\right)  ^{2\alpha_{3}%
}d\phi_{2}\times\ldots\\
& \qquad\times\int\limits_{0}^{\pi/2}\left(  \sin\phi_{d-1}\right)
^{d-2+2\left(  \alpha_{1}+\ldots+\alpha_{d-1}\right)  }\left(  \cos\phi
_{d-1}\right)  ^{2\alpha_{d}}d\phi_{d-1}\\
& =\int\limits_{\varepsilon}^{1}\frac{d\rho}{\rho^{2\sigma-2\left\vert
\alpha\right\vert -d+1}}\text{ }\frac{\Gamma\left(  \alpha_{1}+\frac{1}%
{2}\right)  \Gamma\left(  \alpha_{2}+\frac{1}{2}\right)  }{2\text{ }%
\Gamma\left(  \alpha_{1}+\alpha_{2}+1\right)  }\text{ }\frac{\Gamma\left(
\alpha_{1}+\alpha_{2}+1\right)  \text{ }\Gamma\left(  \alpha_{3}+\frac{1}%
{2}\right)  }{2\text{ }\Gamma\left(  \alpha_{1}+\alpha_{2}+\alpha_{3}+\frac
{3}{2}\right)  }\times\ldots\\
& \qquad\qquad\times\frac{\Gamma\left(  \alpha_{1}+\alpha_{2}+\ldots
+\alpha_{d-1}+\frac{d-1}{2}\right)  \Gamma\left(  \alpha_{d}+\frac{1}%
{2}\right)  }{2\text{ }\Gamma\left(  \alpha_{1}+\alpha_{2}+\ldots+\alpha
_{d}+\frac{d}{2}\right)  }\\
& =\int\limits_{\varepsilon}^{1}\frac{d\rho}{\rho^{2\sigma-2\left\vert
\alpha\right\vert -d+1}}\text{ }\frac{\Gamma\left(  \alpha_{1}+\frac{1}%
{2}\right)  \text{ }\Gamma\left(  \alpha_{2}+\frac{1}{2}\right)  \ldots
\Gamma\left(  \alpha_{d}+\frac{1}{2}\right)  }{2^{d-1}\Gamma\left(  \left\vert
\alpha\right\vert +\frac{d}{2}\right)  }.
\end{align*}

Now $\Gamma\left(  \frac{1}{2}\right)  =\sqrt{\pi}$ and for $k=1,2,3,\ldots$%
\[
\Gamma\left(  k+\frac{1}{2}\right)  =\frac{1.3.5\ldots2k-1}{2^{k}}\sqrt{\pi
}=\frac{\left(  2k\right)  !}{2.4.6\ldots2k}\frac{\sqrt{\pi}}{2^{k}}%
=\frac{\left(  2k\right)  !}{k!}\frac{\sqrt{\pi}}{2^{2k}},
\]

so that%
\[
\int\limits_{\substack{\varepsilon\leq\left\vert x\right\vert \leq1 \\x\geq
0}}\frac{x^{2\alpha}}{\left\vert x\right\vert ^{2\sigma}}dx=\int%
\limits_{\varepsilon}^{1}\frac{d\rho}{\rho^{2\sigma-2\left\vert \alpha
\right\vert -d+1}}\;\frac{\left(  2\alpha\right)  !}{\alpha!\text{
}2^{2\left\vert \alpha\right\vert }}\pi^{d/2}\text{ }\frac{1}{2^{d-1}%
\Gamma\left(  \left\vert \alpha\right\vert +\frac{d}{2}\right)  },
\]

and \ref{a1.32} becomes%
\begin{align*}
\int\limits_{\left\vert \cdot\right\vert \leq1}\frac{1}{\left\vert
\cdot\right\vert ^{2\sigma}}\sum_{\left\vert \beta\right\vert <2n}%
\frac{x^{\beta}}{\beta!}D^{\beta}\psi\left(  0\right)  dx  & =2^{d}%
\lim_{\varepsilon\rightarrow0^{+}}\sum_{\left\vert \alpha\right\vert <n}%
\frac{D^{2\alpha}\psi\left(  0\right)  }{\left(  2\alpha\right)  !}%
\int\limits_{\substack{\varepsilon\leq\left\vert x\right\vert \leq1 \\x\geq
0}}\frac{x^{2\alpha}}{\left\vert x\right\vert ^{2\sigma}}dx\\
& =\lim_{\varepsilon\rightarrow0^{+}}\sum_{\left\vert \alpha\right\vert
<n}\frac{D^{2\alpha}\psi\left(  0\right)  }{\left(  2\alpha\right)  !}\text{
}\int\limits_{\varepsilon}^{1}\frac{d\rho}{\rho^{2\sigma-2\left\vert
\alpha\right\vert -d+1}}\text{ }\frac{\left(  2\alpha\right)  !}{\alpha!\text{
}2^{2\left\vert \alpha\right\vert }}\pi^{\frac{d}{2}}\frac{2^{d}}%
{2^{d-1}\Gamma\left(  \left\vert \alpha\right\vert +\frac{d}{2}\right)  }\\
& =\lim_{\varepsilon\rightarrow0^{+}}\sum_{\left\vert \alpha\right\vert
<n}\frac{2\pi^{\frac{d}{2}}}{2^{2k}\Gamma\left(  k+\frac{d}{2}\right)  }\text{
}\int\limits_{\varepsilon}^{1}\frac{d\rho}{\rho^{2\sigma-2k-d+1}}\text{ }%
\frac{D^{2\alpha}\psi\left(  0\right)  }{\alpha!}\\
& =\lim_{\varepsilon\rightarrow0^{+}}\sum_{k=1}^{n-1}\sum_{\left\vert
\alpha\right\vert =k}\frac{2\pi^{\frac{d}{2}}}{2^{2k}\Gamma\left(  k+\frac
{d}{2}\right)  }\text{ }\int\limits_{\varepsilon}^{1}\frac{d\rho}%
{\rho^{2\sigma-2k-d+1}}\text{ }\frac{D^{2\alpha}\psi\left(  0\right)  }%
{\alpha!}\\
& =\lim_{\varepsilon\rightarrow0^{+}}\sum_{k=1}^{n-1}\frac{2\pi^{\frac{d}{2}}%
}{2^{2k}\Gamma\left(  k+\frac{d}{2}\right)  }\text{ }\int\limits_{\varepsilon
}^{1}\frac{d\rho}{\rho^{2\sigma-2k-d+1}}\text{ }\sum_{\left\vert
\alpha\right\vert =k}\frac{D^{2\alpha}\psi\left(  0\right)  }{\alpha!}\\
& =\lim_{\varepsilon\rightarrow0^{+}}\sum_{k=1}^{n-1}\frac{2\pi^{\frac{d}{2}}%
}{2^{2k}k!\text{ }\Gamma\left(  k+\frac{d}{2}\right)  }\text{ }\int%
\limits_{\varepsilon}^{1}\frac{d\rho}{\rho^{2\sigma-2k-d+1}}\text{ }\left(
\left\vert D\right\vert ^{2k}\psi\right)  \left(  0\right)  ,
\end{align*}

where the last step used the second identity of \ref{1.57} in the Appendix.

We now consider two subcases: $2\sigma-d\neq2,4,6,\ldots$and $2\sigma
-d=2,4,6,\ldots$.\medskip

\fbox{\textbf{Subcase} $2\sigma-d\neq2,4,6,\ldots$} Here $2\sigma-2k-d>0$ for
$1\leq k\leq n-1$ and%
\[
\int\limits_{\varepsilon}^{1}\frac{d\rho}{\rho^{2\sigma-2k-d+1}}=\frac
{1}{2\sigma-2k-d}\left(  \frac{1}{\varepsilon^{2\sigma-2k-d}}-1\right)  ,
\]

so that if we set%
\begin{equation}
a_{k}=\frac{\pi^{\frac{d}{2}}}{2^{2k-1}k!\text{ }\Gamma\left(  k+\frac{d}%
{2}\right)  \left(  2\sigma-2k-d\right)  }\left(  \left\vert D\right\vert
^{2k}\psi\right)  \left(  0\right)  ,\label{a1.30}%
\end{equation}

then%
\begin{align}
\int\limits_{\left\vert \cdot\right\vert \leq1}\frac{1}{\left\vert
\cdot\right\vert ^{2\sigma}}\sum_{\left\vert \beta\right\vert <2n}%
\frac{x^{\beta}}{\beta!}D^{\beta}\psi\left(  0\right)  dx  & =\lim
_{\varepsilon\rightarrow0^{+}}\sum_{k=1}^{n-1}a_{k}\left(  \frac
{1}{\varepsilon^{2\sigma-2k-d}}-1\right) \nonumber\\
& =\lim_{\varepsilon\rightarrow0^{+}}\sum_{k=1}^{n-1}\frac{a_{k}}%
{\varepsilon^{2\sigma-2k-d}}-\sum_{k=1}^{n-1}a_{k}.\label{a1.33}%
\end{align}

It is now not difficult to show that the limit \ref{a1.33} exists iff $a_{k}=0
$ for $1\leq k\leq n-1$ i.e. iff%
\[
\left(  \left\vert D\right\vert ^{2k}\psi\right)  \left(  0\right)
=0,\quad1\leq k\leq n-1.
\]
\medskip

\fbox{\textbf{Subcase} $2\sigma-d=2,4,6,\ldots$} When $k=n-1$, $2\sigma
-2k-d=0$ and so $2\sigma-d=2\left(  n-1\right)  $. When $k<n-1$,
$2\sigma-2k-d>0$. Thus%
\begin{align*}
& \int\limits_{\left\vert \cdot\right\vert \leq1}\frac{1}{\left\vert
\cdot\right\vert ^{2\sigma}}\sum_{\left\vert \beta\right\vert <2n}%
\frac{x^{\beta}}{\beta!}D^{\beta}\psi\left(  0\right)  dx\\
& =2\text{ }\pi^{\frac{d}{2}}\lim_{\varepsilon\rightarrow0^{+}}\sum
_{k=1}^{n-1}\frac{1}{2^{2k}k!\text{ }\Gamma\left(  k+\frac{d}{2}\right)
}\text{ }\int\limits_{\varepsilon}^{1}\frac{d\rho}{\rho^{2\sigma-2k-d+1}%
}\text{ }\left(  \left\vert D\right\vert ^{2k}\psi\right)  \left(  0\right) \\
& =\lim_{\varepsilon\rightarrow0^{+}}\left(  \sum_{k=1}^{n-2}\frac{\pi
^{\frac{d}{2}}}{2^{2k-1}k!\text{ }\Gamma\left(  k+\frac{d}{2}\right)  }\text{
}\int\limits_{\varepsilon}^{1}\frac{d\rho}{\rho^{2\sigma-2k-d+1}}\text{
}\left(  \left\vert D\right\vert ^{2k}\psi\right)  \left(  0\right)  +\right.
\\
& \qquad\qquad\qquad+\left.  \frac{\pi^{\frac{d}{2}}}{2^{2n-3}\left(
n-1\right)  !\text{ }\Gamma\left(  n-1+\frac{d}{2}\right)  }\text{ }%
\int\limits_{\varepsilon}^{1}\frac{d\rho}{\rho^{2\sigma-2\left(  n-1\right)
-d+1}}\text{ }\left(  \left\vert D\right\vert ^{2\left(  n-1\right)  }%
\psi\right)  \left(  0\right)  \right) \\
& =\lim_{\varepsilon\rightarrow0^{+}}\left(  \sum_{k=1}^{n-2}\frac{\pi
^{\frac{d}{2}}}{2^{2k-1}k!\text{ }\Gamma\left(  k+\frac{d}{2}\right)  }\text{
}\int\limits_{\varepsilon}^{1}\frac{d\rho}{\rho^{2\left(  n-k\right)  -1}%
}\text{ }\left(  \left\vert D\right\vert ^{2k}\psi\right)  \left(  0\right)
+\right. \\
& \qquad\qquad\qquad+\left.  \frac{\pi^{\frac{d}{2}}}{2^{2n-3}\left(
n-1\right)  !\text{ }\Gamma\left(  n-1+\frac{d}{2}\right)  }\text{ }%
\int\limits_{\varepsilon}^{1}\frac{d\rho}{\rho}\text{ }\left(  \left\vert
D\right\vert ^{2\left(  n-1\right)  }\psi\right)  \left(  0\right)  \right) \\
& =\lim_{\varepsilon\rightarrow0^{+}}\left(  \sum_{k=1}^{n-2}\frac{\pi
^{\frac{d}{2}}\left(  \frac{1}{\varepsilon^{2\left(  n-k-1\right)  }%
}-1\right)  }{2^{2k}k!\text{ }\Gamma\left(  k+\frac{d}{2}\right)  \left(
n-k-1\right)  }\text{ }\left(  \left\vert D\right\vert ^{2k}\psi\right)
\left(  0\right)  +\right. \\
& \qquad\qquad\qquad+\left.  \frac{\pi^{\frac{d}{2}}\ln\left(  1/\varepsilon
\right)  }{2^{2n-3}\left(  n-1\right)  !\text{ }\Gamma\left(  n-1+\frac{d}%
{2}\right)  }\left(  \left\vert D\right\vert ^{2\left(  n-1\right)  }%
\psi\right)  \left(  0\right)  \right)  .
\end{align*}

Now if we set%
\[
b_{k}=\left\{
\begin{array}
[c]{ll}%
\frac{\pi^{\frac{d}{2}}}{2^{2k}k!\text{ }\Gamma\left(  k+\frac{d}{2}\right)
\text{ }\left(  n-k-1\right)  }\left(  \left\vert D\right\vert ^{2k}%
\psi\right)  \left(  0\right)  , & k=1,\ldots,n-2,\\
\frac{\pi^{\frac{d}{2}}}{2^{2n-3}\left(  n-1\right)  !\text{ }\Gamma\left(
n-1+\frac{d}{2}\right)  }\left(  \left\vert D\right\vert ^{2\left(
n-1\right)  }\psi\right)  \left(  0\right)  , & k=n-1,
\end{array}
\right.
\]

it follows that
\begin{align*}
& \int\limits_{\left\vert \cdot\right\vert \leq1}\frac{1}{\left\vert
\cdot\right\vert ^{2\sigma}}\sum_{\left\vert \beta\right\vert <2n}%
\frac{x^{\beta}}{\beta!}D^{2\beta}\psi\left(  0\right)  dx\\
& =\left\{
\begin{array}
[c]{ll}%
\lim\limits_{\varepsilon\rightarrow0^{+}}b_{n-1}\ln\left(  1/\varepsilon
\right)  , & n=2,\\
\lim\limits_{\varepsilon\rightarrow0^{+}}\left(  \sum\limits_{k=1}^{n-2}%
b_{k}\left(  \frac{1}{\varepsilon^{2\left(  n-k-1\right)  }}-1\right)
+b_{n-1}\ln\left(  1/\varepsilon\right)  \right)  , & n\geq3,
\end{array}
\right. \\
& =\left\{
\begin{array}
[c]{ll}%
\lim\limits_{\varepsilon\rightarrow0^{+}}b_{n-1}\ln\left(  1/\varepsilon
\right)  , & n=2,\\
\lim\limits_{\varepsilon\rightarrow0^{+}}\left(  \frac{b_{1}}{\varepsilon
^{2\left(  n-2\right)  }}+\frac{b_{2}}{\varepsilon^{2\left(  n-3\right)  }%
}+\ldots+\frac{b_{n-2}}{\varepsilon^{2}}+b_{n-1}\ln\left(  1/\varepsilon
\right)  \right)  -\sum\limits_{k=1}^{n-2}b_{k}, & n\geq3,
\end{array}
\right.
\end{align*}

and these limits exist iff $b_{k}=0$ when $1\leq k\leq n-1$ i.e.%
\begin{equation}
\left(  \left\vert D\right\vert ^{2k}\psi\right)  \left(  0\right)
=0,\quad1\leq k\leq n-1.\label{a1.38}%
\end{equation}

Thus the conclusion of these two subcases is the same, namely \ref{a1.38}.
Combining this result with \ref{a1.39} we deduce that%
\begin{equation}
\left(  \left\vert D\right\vert ^{2k}\psi\right)  \left(  0\right)
=0,\quad0\leq k\leq n-1.\label{a1.44}%
\end{equation}

The final step in this proof is to show by induction on $n\geq2$ that
\ref{a1.44} implies $\phi\in S_{\emptyset,n}$. When $n=2$ we use \ref{a1.44}
to deduce that
\begin{align*}
0  & =\left(  \left\vert D\right\vert ^{2}\psi\right)  \left(  0\right)
=\sum_{i=1}^{d}D_{i}^{2}\psi\left(  0\right)  =\sum_{i=1}^{d}D_{i}^{2}\left(
\phi\overline{\phi}\right)  \left(  0\right) \\
& =\sum_{i=1}^{d}\left(  2\left\vert D_{i}\phi\left(  0\right)  \right\vert
^{2}+\phi\left(  0\right)  D_{i}^{2}\overline{\phi\left(  0\right)
}+\overline{\phi\left(  0\right)  }D_{i}^{2}\phi\left(  0\right)  \right) \\
& =\sum_{i=1}^{d}2\left\vert D_{i}\phi\left(  0\right)  \right\vert ^{2},
\end{align*}

and so $D_{i}\phi\left(  0\right)  =0$ for $i=1,\ldots,d$ i.e. $\phi\in
S_{\emptyset,n}$ when $n=2$.

Now assume the result is true for arbitrary $n$. Then if $\left\vert
\alpha\right\vert =n+1$, applying Leibniz' theorem gives
\begin{align*}
D^{2\alpha}\psi\left(  0\right)   & =\sum_{\beta\leq2\alpha}\binom{2\alpha
}{\beta}\left(  D^{\beta}\phi\right)  \left(  0\right)  \left(  D^{2\alpha
-\beta}\overline{\phi}\right)  \left(  0\right) \\
& =\sum_{\beta<\alpha}\binom{2\alpha}{\beta}\left(  D^{\beta}\phi\right)
\left(  0\right)  \left(  D^{2\alpha-\beta}\overline{\phi}\right)  \left(
0\right)  +\sum_{\beta=\alpha}\binom{2\alpha}{\beta}\left(  D^{\beta}%
\phi\right)  \left(  0\right)  \left(  D^{2\alpha-\beta}\overline{\phi
}\right)  \left(  0\right)  +\\
& \qquad\qquad+\sum_{\beta>\alpha}\binom{2\alpha}{\beta}\left(  D^{\beta}%
\phi\right)  \left(  0\right)  \left(  D^{2\alpha-\beta}\overline{\phi
}\right)  \left(  0\right) \\
& =\sum_{\beta<\alpha}\binom{2\alpha}{\beta}\left(  D^{\beta}\phi\right)
\left(  0\right)  \left(  D^{2\alpha-\beta}\overline{\phi}\right)  \left(
0\right)  +\binom{2\alpha}{\alpha}\left(  D^{\alpha}\phi\right)  \left(
0\right)  \left(  D^{\alpha}\overline{\phi}\right)  \left(  0\right)  +\\
& \qquad\qquad+\sum_{\beta<\alpha}\binom{2\alpha}{\beta}\left(  D^{2\alpha
-\beta}\phi\right)  \left(  0\right)  \left(  D^{\beta}\overline{\phi}\right)
\left(  0\right) \\
& =\binom{2\alpha}{\alpha}\left\vert D^{\alpha}\phi\left(  0\right)
\right\vert ^{2},
\end{align*}

where the last step is valid since $\beta<\alpha$ implies $\left\vert
\beta\right\vert <\left\vert \alpha\right\vert =n+1$ and hence $D^{\beta}%
\phi\left(  0\right)  =0$. We now have%
\[
0=\left(  \left\vert D\right\vert ^{2\left(  n+1\right)  }\psi\right)
=\sum_{\left\vert \alpha\right\vert =n+1}\frac{1}{\alpha!}\binom{2\alpha
}{\alpha}\left\vert D^{\alpha}\phi\left(  0\right)  \right\vert ^{2},
\]

and so $D^{\alpha}\phi\left(  0\right)  =0$ whenever $\left\vert
\alpha\right\vert =n+1$, which completes the inductive argument.\medskip

Finally, we prove inequality \ref{a1.41}. We now know that if $\int%
_{\left\vert \cdot\right\vert \leq1}\frac{\left\vert \phi\right\vert ^{2}%
}{\left\vert \cdot\right\vert ^{2\sigma}}<\infty$ then $\phi\in S_{\emptyset
,n}$. Thus the Taylor series expansion of Section \ref{Sect_apx_TaylorSeries}
can be used to write%
\[
\phi\left(  x\right)  =\left(  \mathcal{R}_{n}\phi\right)  \left(  0,x\right)
,\quad n=1,2,3,\ldots,
\]

with remainder estimate \ref{a1.37}:%
\[
\left\vert \left(  \mathcal{R}_{n}\phi\right)  \left(  0,x\right)  \right\vert
\leq\left\vert x\right\vert ^{n}\sum_{\left\vert \beta\right\vert =n}%
\max_{y\in\left[  0,x\right]  }\left\vert \left(  D^{\beta}\phi\right)
(y)\right\vert \leq\left\vert x\right\vert ^{n}\sum_{\left\vert \beta
\right\vert =n}\left\Vert D^{\beta}\phi\right\Vert _{\infty}.
\]

Hence%
\begin{align*}
\int_{\left\vert \cdot\right\vert \leq1}\frac{\left\vert \phi\right\vert ^{2}%
}{\left\vert \cdot\right\vert ^{2\sigma}}\leq\int_{\left\vert \cdot\right\vert
\leq1}\frac{\left\vert \left(  \mathcal{R}_{n}\phi\right)  \left(
0,\cdot\right)  \right\vert ^{2}}{\left\vert \cdot\right\vert ^{2\sigma}} &
\leq\left(  \int_{\left\vert \cdot\right\vert \leq1}\frac{\left\vert
\cdot\right\vert ^{2n}}{\left\vert \cdot\right\vert ^{2\sigma}}\right)
\left(  \sum_{\left\vert \beta\right\vert =n}\left\Vert D^{\beta}%
\phi\right\Vert _{\infty}\right)  ^{2}\\
&  =\left(  \int_{\left\vert \cdot\right\vert \leq1}\frac{1}{\left\vert
\cdot\right\vert ^{2\sigma-2n}}\right)  \left(  \sum_{\left\vert
\beta\right\vert =n}\left\Vert D^{\beta}\phi\right\Vert _{\infty}\right)
^{2},
\end{align*}

and we set $k_{\sigma}=\int_{\left\vert \cdot\right\vert \leq1}\frac
{1}{\left\vert \cdot\right\vert ^{2\sigma-2n}}$.\medskip

Case 1 implies that this inequality also holds for $n=0$. This proves the
estimate \ref{a1.41}.

Using the spherical polar change of variables \ref{a1.42} and the formula
\ref{a1.40} we obtain

$2n>2\sigma-d$
\begin{align*}
\int_{\left\vert x\right\vert \leq1}\frac{dx}{\left\vert x\right\vert
^{2\sigma-2n}}  & =\int\limits_{0}^{1}\int\limits_{\left[  0,\frac{\pi}%
{2}\right]  ^{d-1}}\frac{1}{\rho^{2\sigma-2n}}\rho^{d}\prod\limits_{j=1}%
^{d-1}\left(  \sin\phi_{j}\right)  ^{j-1}d\phi d\rho\\
& =\left(  \int\limits_{0}^{1}\rho^{2n-2\sigma+d}d\rho\right)  \left(
\int\limits_{0}^{\pi/2}d\phi_{1}\right)  \left(  \int\limits_{0}^{\pi/2}%
\sin\phi_{2}d\phi_{2}\right)  \ldots\left(  \int\limits_{0}^{\pi/2}\left(
\sin\phi_{d-1}\right)  ^{d-2}d\phi_{2}\right) \\
& =\frac{1}{2n-2\sigma+d+1}\left(  \int\limits_{0}^{\pi/2}d\phi_{1}\right)
\left(  \int\limits_{0}^{\pi/2}\sin\phi_{2}d\phi_{2}\right)  \ldots\left(
\int\limits_{0}^{\pi/2}\left(  \sin\phi_{d-1}\right)  ^{d-2}d\phi_{2}\right)
\end{align*}

From \ref{a1.40}
\[
\int_{0}^{\pi/2}\left(  \sin\phi\right)  ^{m}d\phi=\frac{\Gamma\left(
\frac{m+1}{2}\right)  \Gamma\left(  \frac{1}{2}\right)  }{2\text{ }%
\Gamma\left(  \frac{m}{2}+1\right)  },
\]

so%
\begin{align*}
\int_{\left\vert x\right\vert \leq1}\frac{dx}{\left\vert x\right\vert
^{2\sigma-2n}}  & =\frac{1}{2n-2\sigma+d+1}\text{ }\frac{\Gamma\left(
\frac{1}{2}\right)  \Gamma\left(  \frac{1}{2}\right)  }{2\text{ }\Gamma\left(
1\right)  }\text{ }\frac{\Gamma\left(  1\right)  \Gamma\left(  \frac{1}%
{2}\right)  }{2\text{ }\Gamma\left(  \frac{3}{2}\right)  }\text{ }\frac
{\Gamma\left(  \frac{3}{2}\right)  \Gamma\left(  \frac{1}{2}\right)  }{2\text{
}\Gamma\left(  2\right)  }\times\ldots\\
& \qquad\qquad\times\frac{\Gamma\left(  \frac{d-1}{2}\right)  \Gamma\left(
\frac{1}{2}\right)  }{2\text{ }\Gamma\left(  \frac{d}{2}\right)  }\\
& =\frac{1}{2n-2\sigma+d+1}\frac{1}{2^{d-1}}\frac{\Gamma\left(  \frac{1}%
{2}\right)  ^{d}}{\Gamma\left(  \frac{d}{2}\right)  }\\
& =\frac{1}{2n-2\sigma+d+1}\frac{1}{2^{d-1}}\frac{\pi^{\frac{d}{2}}}%
{\Gamma\left(  \frac{d}{2}\right)  }\\
& =\frac{2}{2n-2\sigma+d+1}\left(  \frac{\sqrt{\pi}}{2}\right)  ^{\frac{d}{2}%
}\frac{1}{\Gamma\left(  \frac{d}{2}\right)  }.
\end{align*}

\end{proof}

\begin{remark}
\label{Rem_Lem_wt_radial_near_0}\ 

\begin{enumerate}
\item Does anybody have an easier proof?

\item Since $\sqrt{\pi}/2<1$ the constant $\int_{\left\vert \cdot\right\vert
\leq1}\frac{1}{\left\vert \cdot\right\vert ^{2\sigma-2n}}=\frac{2}%
{2n-2\sigma+d+1}\frac{\left(  \sqrt{\pi}/2\right)  ^{d/2}}{\Gamma\left(
d/2\right)  }$ decreases very rapidly but note that if we wanted to average
the right side we would have to divide the equation by $\left(  \sum
_{\left\vert \beta\right\vert =n}1\right)  ^{2}=\binom{d+n-1}{n}^{2}$ (part 8
of Definition \ref{Def_multi_id} in the Appendix).
\end{enumerate}
\end{remark}

The previous lemma considered some consequences of the existence of
$\int_{\left\vert \cdot\right\vert \leq1}\frac{\left\vert \phi\right\vert
^{2}}{\left\vert \cdot\right\vert ^{2\sigma}}$. We now use simple scaling
arguments to consider $\int_{\left\vert \cdot\right\vert \leq r}%
\frac{\left\vert \phi\right\vert ^{2}}{\left\vert \cdot\right\vert ^{2\sigma}%
}<\infty$ for any $r>0$.

\begin{lemma}
\label{Lem_2_wt_radial_near_0}Suppose $r>0$, $\sigma\geq0$ and $\phi\in
S\left(  \mathbb{R}^{d}\right)  $. Let $n$ be the smallest non-negative
integer such that $2n>2\sigma-d$.

Then $\int_{\left\vert \cdot\right\vert \leq r}\frac{\left\vert \phi
\right\vert ^{2}}{\left\vert \cdot\right\vert ^{2\sigma}}<\infty$ iff $\phi\in
S_{\emptyset,n}$. Further, if $\int_{\left\vert \cdot\right\vert \leq r}%
\frac{\left\vert \phi\right\vert ^{2}}{\left\vert \cdot\right\vert ^{2\sigma}%
}<\infty$ then%
\begin{equation}
\int_{\left\vert \cdot\right\vert \leq r}\frac{\left\vert \phi\right\vert
^{2}}{\left\vert \cdot\right\vert ^{2\sigma}}\leq k_{\sigma}r^{2n-2\sigma
+d}\left(  \sum_{\left\vert \beta\right\vert =n}\left\Vert D^{\beta}%
\phi\right\Vert _{\infty}\right)  ^{2},\label{a1.45}%
\end{equation}

where $k_{\sigma}$ is given by \ref{a1.46}.
\end{lemma}

\begin{proof}
The previous lemma considered the case $r=1$ and so $\int_{\left\vert
\cdot\right\vert \leq r}\frac{\left\vert \phi\right\vert ^{2}}{\left\vert
\cdot\right\vert ^{2\sigma}}<\infty$ iff $\int_{\left\vert \cdot\right\vert
\leq1}\frac{\left\vert \phi\right\vert ^{2}}{\left\vert \cdot\right\vert
^{2\sigma}}<\infty$ iff $\phi\in S_{\emptyset,n}$.

Since $\int_{\left\vert \cdot\right\vert \leq r}\frac{\left\vert
\phi\right\vert ^{2}}{\left\vert \cdot\right\vert ^{2\sigma}}=r^{d-2\sigma
}\int_{\left\vert y\right\vert \leq1}\frac{\left\vert \phi\left(  ry\right)
\right\vert ^{2}}{\left\vert y\right\vert ^{2\sigma}}dy$ we have from the
previous lemma%
\begin{align*}
\int_{\left\vert \cdot\right\vert \leq r}\frac{\left\vert \phi\right\vert
^{2}}{\left\vert \cdot\right\vert ^{2\sigma}}=r^{d-2\sigma}\int_{\left\vert
y\right\vert \leq1}\frac{\left\vert \phi\left(  ry\right)  \right\vert ^{2}%
}{\left\vert y\right\vert ^{2\sigma}}dy &  =r^{d-2\sigma}\int_{\left\vert
\cdot\right\vert \leq1}\frac{1}{\left\vert \cdot\right\vert ^{2\sigma-2n}%
}\left(  \sum_{\left\vert \beta\right\vert =n}\left\Vert D^{\beta}\left(
\phi\left(  ry\right)  \right)  \right\Vert _{\infty}\right)  ^{2}\\
&  =r^{2n-2\sigma+d}\int_{\left\vert \cdot\right\vert \leq1}\frac
{1}{\left\vert \cdot\right\vert ^{2\sigma-2n}}\left(  \sum_{\left\vert
\beta\right\vert =n}\left\Vert D^{\beta}\phi\right\Vert _{\infty}\right)
^{2}\\
&  =k_{\sigma}r^{2n-2\sigma+d}\left(  \sum_{\left\vert \beta\right\vert
=n}\left\Vert D^{\beta}\phi\right\Vert _{\infty}\right)  ^{2}.
\end{align*}

\end{proof}

We now characterize the space $S_{w,0}$ for the weight functions of Example
\ref{Ex_wt_func_|x|^m_L1loc} as well as showing that $w\in W_{S;0}$
(Definition \ref{Def_Sw2_and_fnal(Sw2)}).

\begin{theorem}
\label{Thm_Sw2_eq_Son}Suppose the weight function $w$ is from Example
\ref{Ex_wt_func_|x|^m_L1loc}. Then:
\begin{equation}
S_{w,0}=S_{\emptyset,n},\label{a1.28}%
\end{equation}

where $n$ be the smallest integer such that $n>m-\left(  \frac{d}{2}-s\right)
$ and $n\geq0$. Set $\sigma=m+s$. Further, $w\in W_{S;0}$ and%
\[
\int w\left\vert \phi\right\vert ^{2}\leq\max\left\{  c_{2}k_{\sigma
}r^{2n-2\sigma+d},\int\limits_{\left\vert \cdot\right\vert \geq r}\frac
{w}{\left\vert \cdot\right\vert ^{2\tau}}\right\}  \left(  \left\Vert \left(
1+\left\vert \cdot\right\vert \right)  ^{\left\lceil \tau\right\rceil }%
\phi\right\Vert _{\infty}+\sum_{\left\vert \beta\right\vert =n}\left\Vert
D^{\beta}\phi\right\Vert _{\infty}\right)  ^{2},\quad\phi\in S_{w,0},
\]

where $k_{\sigma}$ is given by \ref{a1.46}.
\end{theorem}

\begin{proof}
The weight function properties are:%
\[
w\left(  x\right)  =\frac{v\left(  x\right)  }{\left\vert x\right\vert
^{2\left(  m+s\right)  }},\quad\left\vert x\right\vert \leq r;\quad
2s<d;\quad0<c_{1}\leq v\left(  x\right)  \leq c_{2},
\]

which imply that for $\phi\in S$,%
\begin{align}
c_{1}\int_{\left\vert \cdot\right\vert \leq r}\frac{1}{\left\vert
\cdot\right\vert ^{2m+2s}}\left\vert \phi\right\vert ^{2} &  \leq
\int_{\left\vert \cdot\right\vert \leq r}w\left\vert \phi\right\vert ^{2}%
=\int_{\left\vert \cdot\right\vert \leq r}w\left\vert \phi\right\vert
^{2}+\int_{\left\vert \cdot\right\vert \geq r}w\left\vert \phi\right\vert
^{2}\leq\nonumber\\
&  \leq\int\limits_{\left\vert \cdot\right\vert \leq r}w\left\vert
\phi\right\vert ^{2}+\int\limits_{\left\vert \cdot\right\vert \geq r}\frac
{w}{\left\vert \cdot\right\vert ^{2\tau}}\left\vert \cdot\right\vert ^{2\tau
}\left\vert \phi\right\vert ^{2}\nonumber\\
&  \leq c_{2}\int\limits_{\left\vert \cdot\right\vert \leq r}\frac
{1}{\left\vert \cdot\right\vert ^{2m+2s}}\left\vert \phi\right\vert
^{2}+\left(  \int\limits_{\left\vert \cdot\right\vert \geq r}\frac
{w}{\left\vert \cdot\right\vert ^{2\tau}}\right)  \left\Vert \left(
1+\left\vert \cdot\right\vert \right)  ^{\left\lceil \tau\right\rceil }%
\phi\right\Vert _{\infty}^{2},\label{a1.36}%
\end{align}

and so $\phi\in S_{w,0}$ iff $\int_{\left\vert \cdot\right\vert \leq r}%
\frac{1}{\left\vert \cdot\right\vert ^{2m+2s}}\left\vert \phi\right\vert
^{2}<\infty$. But by Lemma \ref{Lem_2_wt_radial_near_0}, $\int_{\left\vert
\cdot\right\vert \leq r}\frac{1}{\left\vert \cdot\right\vert ^{2m+2s}%
}\left\vert \phi\right\vert ^{2}<\infty$ iff $\phi\in S_{\emptyset,n}$ and we
can conclude that $S_{\emptyset,n}=S_{w,0}$.

Using \ref{a1.45}, \ref{a1.36} becomes%
\begin{align*}
\int w\left\vert \phi\right\vert ^{2}  & \leq c_{2}k_{\sigma}r^{2n-2\sigma
+d}\left(  \sum_{\left\vert \beta\right\vert =n}\left\Vert D^{\beta}%
\phi\right\Vert _{\infty}\right)  ^{2}+\left(  \int\limits_{\left\vert
\cdot\right\vert \geq r}\frac{w}{\left\vert \cdot\right\vert ^{2\tau}}\right)
\left\Vert \left(  1+\left\vert \cdot\right\vert \right)  ^{\left\lceil
\tau\right\rceil }\phi\right\Vert _{\infty}^{2}\\
& \leq\max\left\{  c_{2}k_{\sigma}r^{2n-2\sigma+d},\int\limits_{\left\vert
\cdot\right\vert \geq r}\frac{w}{\left\vert \cdot\right\vert ^{2\tau}%
}\right\}  \left(  \left\Vert \left(  1+\left\vert \cdot\right\vert \right)
^{\left\lceil \tau\right\rceil }\phi\right\Vert _{\infty}+\sum_{\left\vert
\beta\right\vert =n}\left\Vert D^{\beta}\phi\right\Vert _{\infty}\right)
^{2}.
\end{align*}

\end{proof}

\begin{remark}
?? When is $w_{1}S_{w,0}\subset S$ where $w_{1}^{2}=w$ and $w_{1}\in??$ etc.?
\end{remark}

\begin{lemma}
\label{Lem_Son_Taylor_represent}Suppose $\phi\in S$ and $n\geq1$ is an
integer. Then $\phi\in S_{\emptyset,n}$ iff%
\begin{equation}
\phi\left(  x\right)  =n\sum_{\left\vert \beta\right\vert =n}\frac{x^{\beta}%
}{\beta!}\int_{0}^{1}\left(  1-t\right)  ^{n-1}\left(  D^{\beta}\phi\right)
\left(  tx\right)  dt.\label{a1.205}%
\end{equation}

\end{lemma}

\begin{proof}
If $\phi\in S_{\emptyset,n}$ then the Taylor series expansion of $\phi$ about
the origin%
\[
\phi\left(  x\right)  =\sum_{\left\vert \beta\right\vert <n}\frac{x^{\beta}%
}{\beta!}D^{\beta}\phi\left(  0\right)  +n\sum_{\left\vert \beta\right\vert
=n}\frac{x^{\beta}}{\beta!}\int_{0}^{1}\left(  1-t\right)  ^{n-1}\left(
D^{\beta}\phi\right)  \left(  tx\right)  dt,
\]

yields \ref{a1.205}. On the other hand, if \ref{a1.205} holds then
$\sum_{\left\vert \beta\right\vert <n}\frac{x^{\beta}}{\beta!}D^{\beta}%
\phi\left(  0\right)  =0$ is true for all $x$ and so $D^{\beta}\phi\left(
0\right)  =0$ when $\left\vert \beta\right\vert <n$ i.e. $\phi\in
S_{\emptyset,n}$.
\end{proof}

??? The next theorem is a \textbf{FAILED ATTEMPT} at an analogue of Theorem
\ref{Thm_eta^2/sin(eta)^2.|phi|^2_dim_gt_1}.

\begin{theorem}
Suppose the weight function $w$ is from Example \ref{Ex_wt_func_|x|^m_L1loc}
and $n\geq1$.

\begin{enumerate}
\item Suppose also that $\left(  w_{1}\left(  x\right)  \right)
^{2}:=w\left(  x\right)  $ and $1/w_{1}\in C_{BP}^{\infty}$. Then $\frac
{1}{w_{1}}S\subset S_{w,0}$.

\item \textbf{FAILED}! If, in addition, $w_{1}\left(  x\right)  x^{\beta}\in S
$ when $\left\vert \beta\right\vert =n$ then $S_{w,0}=\frac{1}{w_{1}}S$.
\end{enumerate}
\end{theorem}

\begin{proof}
\textbf{Part 1} Clearly $\frac{1}{w_{1}}S\subset S$ so that if $\phi=\frac
{1}{w_{1}}\psi$ for some $\psi\in S$ then $\int w\left\vert \phi\right\vert
^{2}=\int\left\vert \psi\right\vert ^{2}<\infty$ and $\frac{1}{w_{1}}S\subset
S_{w,0}$.

\textbf{Part 2} On the other hand, assume $\phi\in S_{w,0}$. Now Theorem
\ref{Thm_Sw2_eq_Son} implies $S_{w,0}=S_{\emptyset,n}$ and so $\phi\in
S_{\emptyset,n}$.%
\[
\left(  w_{1}\phi\right)  \left(  x\right)  =n\sum_{\left\vert \beta
\right\vert =n}\frac{w_{1}\left(  x\right)  x^{\beta}}{\beta!}\int_{0}%
^{1}\left(  1-t\right)  ^{n-1}\left(  D^{\beta}\phi\right)  \left(  tx\right)
dt\in S.
\]

?? \textbf{HOWEVER} the assumptions $1/w_{1}\in C_{BP}^{\infty}$ and
$w_{1}\left(  x\right)  x^{\beta}\in S$ when $\left\vert \beta\right\vert =n$
are incompatible. Indeed $w_{1}\left(  x\right)  x^{\beta}\in S$ when
$\left\vert \beta\right\vert =n$ implies $w\left\vert \cdot\right\vert
^{2n}\in S$ and $1/w_{1}\in C_{BP}^{\infty}$ implies $1/w\in C_{BP}^{\infty}$.
Thus $\left\vert \cdot\right\vert ^{2n}\in\frac{1}{w}S\subset C_{BP}^{\infty
}.S\subset S$ which is a contradiction.
\end{proof}

\subsubsection{\protect\underline{The extended B-splines}}

Here we are compelled to use more general extended B-splines here, namely
those with weight functions $w\left(  \xi\right)  =\frac{\xi^{2\nu}}%
{\sin^{2\lambda}\xi}$ with multi-index parameters $\lambda$ and $\nu$
satisfying $\lambda\geq\nu\geq\mathbf{1}$. Here $\sin^{2\lambda}\xi:=\left(
\sin^{2\lambda_{k}}\xi_{k}\right)  =\left(  \sin\xi_{k}\right)  ^{2\lambda}$.

\begin{definition}
\label{Def_space_Sxm}\textbf{The spaces }$S_{\otimes,\alpha}\left(  \left\{
\mathbf{0}\right\}  \right)  $, $S_{\otimes,\alpha}\left(  \left\{  x\right\}
\right)  $ and $S_{\otimes,\alpha}\left(  \mathbb{\Sigma}\right)  $.

For any multi-index $\alpha$,
\begin{align*}
S_{\otimes,\alpha}\left(  \left\{  \mathbf{0}\right\}  \right)   &
=S_{\otimes,\alpha}\left(  \mathbf{0}\right) \\
& =\left\{
\begin{array}
[c]{ll}%
S, & \alpha=\mathbf{0},\\
\left\{  \phi\in S:\left(  D_{k}^{n}\phi\right)  \left(  x^{\prime}%
,0_{k},x^{\prime\prime}\right)  =0,\text{\quad}\left(  x^{\prime}%
,x^{\prime\prime}\right)  \in\mathbb{R}^{d-1},\text{ }n<\alpha_{k}%
\mathbf{,}\text{ }\forall k\right\}  , & \alpha\neq\mathbf{0}.
\end{array}
\right.
\end{align*}

$S_{\otimes,\alpha}$ is endowed with the subspace topology induced by $S$, and
$S_{\otimes,\alpha}^{\prime}$ denotes the corresponding bounded linear
functionals on $S_{\otimes,\alpha}$.

For an arbitrary point%
\[
\phi\in S_{\otimes,\alpha}\left(  x\right)  :=S_{\otimes,\alpha}\left(
\left\{  x\right\}  \right)  \text{ }iff\text{ }\phi\left(  \cdot+x\right)
\in S_{\otimes,\alpha}\left(  \mathbf{0}\right)  .
\]

More generally, for any set $\Sigma\subset\mathbb{R}^{d}$ without a finite
cluster point,%
\[
S_{\otimes,\alpha}\left(  \mathbb{\Sigma}\right)  =\bigcap_{x\in\Sigma
}S_{\otimes,\alpha}\left(  \left\{  x\right\}  \right)  .
\]

\end{definition}

\begin{lemma}
\label{Lem_Sxl_and_Sx,n,l}\ 

\begin{enumerate}
\item If $\phi\in S_{\otimes,\alpha}\left(  \mathbb{\Sigma}\right)  $ then
$D^{\beta}\phi\in S_{\otimes,\alpha-\beta}\left(  \mathbb{\Sigma}\right)  $
when $\beta\leq\alpha$.

\item If $\phi\in S_{\otimes,\alpha}\left(  \mathbf{0}\right)  $ then
$x^{\beta}\phi\in S_{\otimes,\alpha+\beta}\left(  \mathbf{0}\right)  $.

\item If $\phi\in S$ and $\beta\leq\alpha$ then $\phi\in S_{\otimes
,\alpha-\beta}\left(  \mathbf{0}\right)  $ iff $x^{\beta}\phi\in
S_{\otimes,\alpha}\left(  \mathbf{0}\right)  $.

This can be expressed as $S_{\otimes,\alpha-\beta}\left(  \mathbf{0}\right)
=S\cap x^{-\beta}S_{\otimes,\alpha}\left(  \mathbf{0}\right)  $.

\item If $z_{+}.>\mathbf{0}$ and $\phi\in S$ then $\phi\in S_{\otimes
,\alpha-\beta}\left(  \left\{  z\right\}  \right)  $ iff $x^{\beta}\phi\in
S_{\otimes,\alpha}\left(  \left\{  z\right\}  \right)  $.

This can be expressed as $S_{\otimes,\alpha}\left(  \left\{  z\right\}
\right)  =S\cap x^{-\beta}S_{\otimes,\alpha}\left(  \left\{  z\right\}
\right)  $.\smallskip

Parts 3 and 4 can be generalized to:\smallskip

\item Suppose $\phi\in S$, $z\in\mathbb{R}^{d}$ and $\alpha\geq\left(
z\And\mathbf{0}\right)  .\beta$ where%
\[
\left(  z\And\mathbf{0}\right)  _{i}:=\left\{
\begin{array}
[c]{ll}%
1, & if\text{\ }z_{i}=0,\\
0, & if\text{\ }z_{i}\neq0.
\end{array}
\right.
\]

Then $x^{\beta}\phi\in S_{\otimes,\alpha}\left(  \left\{  z\right\}  \right)
$ iff $\phi\in S_{\otimes,\alpha-\left(  z\And\mathbf{0}\right)  .\beta
}\left(  \left\{  z\right\}  \right)  $.

This can be expressed as $S_{\otimes,\alpha-\left(  z\And\mathbf{0}\right)
.\beta}\left(  \left\{  z\right\}  \right)  =S\cap x^{-\beta}S_{\otimes
,\alpha}\left(  \left\{  z\right\}  \right)  $.

\item If $f\in C_{BP}^{\infty}$ and $\phi\in S_{\otimes,\alpha}\left(
\mathbb{\Sigma}\right)  $ then $f\phi\in S_{\otimes,\alpha}\left(
\mathbb{\Sigma}\right)  $. Also $f,\phi\in S_{\otimes,\alpha}\left(
\mathbb{\Sigma}\right)  $ implies $f\phi\in S_{\otimes,\alpha}\left(
\mathbb{\Sigma}\right)  $.
\end{enumerate}
\end{lemma}

\begin{proof}
\textbf{Part 1} Follows directly from the definition of $S_{\otimes,\alpha
}\left(  \mathbb{\Sigma}\right)  $.\smallskip

\textbf{Part 2} Use 1-dimensional binomial theorem.\smallskip

\textbf{Part 3} From part 2, $\phi\in S_{\otimes,\alpha-\beta}\left(
\mathbf{0}\right)  $ implies $x^{\beta}\phi\in S_{\otimes,\alpha}\left(
\mathbf{0}\right)  $.

Conversely, suppose $x^{\beta}\phi\in S_{\otimes,\alpha}\left(  \mathbf{0}%
\right)  $. Now%
\begin{align}
D_{k}^{m}\left(  x^{\beta}\phi\right)  =\sum\limits_{n=0}^{m}\tbinom{m}%
{n}\left(  D_{k}^{n}x^{\beta}\right)  D_{k}^{m-n}\phi & =\sum\limits_{n=0}%
^{\min\left\{  \beta_{k},m\right\}  }\tbinom{m}{n}\left(  D_{k}^{n}x^{\beta
}\right)  D_{k}^{m-n}\phi\nonumber\\
& =\sum\limits_{n=0}^{\min\left\{  \beta_{k},m\right\}  }\tbinom{m}{n}%
\tfrac{\beta_{k}!}{\left(  \beta_{k}-n\right)  !}x^{\beta-n\mathbf{e}_{k}%
}D_{k}^{m-n}\phi,\label{a1.014}%
\end{align}

and we consider three cases:\medskip

\fbox{\textbf{Case} $m<\beta_{k}<\alpha_{k}$} Here%
\[
D_{k}^{m}\left(  x^{\beta}\phi\right)  =\sum\limits_{n=0}^{m}\tbinom{m}%
{n}\tfrac{\alpha_{k}!}{\left(  \alpha_{k}-n\right)  !}x^{\beta-n\mathbf{e}%
_{k}}D_{k}^{m-n}\phi,
\]

so%
\[
\sum\limits_{n=0}^{m}\tbinom{m}{n}\tfrac{\alpha_{k}!}{\left(  \alpha
_{k}-n\right)  !}x^{\beta-n\mathbf{e}_{k}}\left(  D_{k}^{m-n}\phi\right)
\left(  x^{\prime},0_{k},x^{\prime\prime}\right)  =0,
\]

and hence $x^{\beta}\phi\in S_{\otimes,\alpha}\left(  \mathbf{0}\right)  $
supplies no new information.\smallskip

\fbox{\textbf{Case} $m=\beta_{k}<\alpha_{k}$} In this case \ref{a1.014}
becomes%
\[
D_{k}^{m}\left(  x^{\beta}\phi\right)  =\sum\limits_{n=0}^{\beta_{k}}%
\tbinom{\beta_{k}}{n}\tfrac{\beta_{k}!}{\left(  \beta_{k}-\beta_{k}\right)
!}x^{\beta-n\mathbf{e}_{k}}D_{k}^{\beta_{k}-n}\phi,
\]

and hence $x^{\beta}\phi\in S_{\otimes,\alpha}\left(  \mathbf{0}\right)  $
implies%
\[
0=D_{k}^{m}\left(  x^{\beta}\phi\right)  \left(  x^{\prime},0_{k}%
,x^{\prime\prime}\right)  =\beta_{k}!\left(  x^{\prime}\right)  ^{\beta
^{\prime}}\left(  x^{\prime\prime}\right)  ^{\beta^{\prime\prime}}\phi\left(
x^{\prime},0_{k},x^{\prime\prime}\right)  ,
\]

i.e.%
\[
\phi\left(  x^{\prime},0_{k},x^{\prime\prime}\right)  =0,\quad\forall x.
\]
\smallskip

\fbox{\textbf{Case} $\beta_{k}<m<\alpha_{k}$} Observe that \ref{a1.014} is now%
\[
D_{k}^{m}\left(  x^{\beta}\phi\right)  =\sum\limits_{n=0}^{\beta_{k}}%
\tbinom{m}{n}\tfrac{\beta_{k}!}{\left(  \beta_{k}-n\right)  !}x^{\beta
-n\mathbf{e}_{k}}D_{k}^{m-n}\phi,
\]

and so $x^{\beta}\phi\in S_{\otimes,\alpha}\left(  \mathbf{0}\right)  $
implies%
\begin{align*}
0  & =D_{k}^{m}\left(  x^{\beta}\phi\right)  \left(  x^{\prime},0_{k}%
,x^{\prime\prime}\right) \\
& =\sum\limits_{n=0}^{\beta_{k}}\tbinom{m}{n}\tfrac{\beta_{k}!}{\left(
\beta_{k}-n\right)  !}\left(  x^{\prime},0_{k},x^{\prime\prime}\right)
^{\beta-n\mathbf{e}_{k}}\left(  D_{k}^{m-n}\phi\right)  \left(  x^{\prime
},0_{k},x^{\prime\prime}\right) \\
& =\tbinom{m}{\beta_{k}}\tfrac{\beta_{k}!}{\left(  \beta_{k}-\beta_{k}\right)
!}\left(  x^{\prime},0_{k},x^{\prime\prime}\right)  ^{\beta-\beta
_{k}\mathbf{e}_{k}}\left(  D_{k}^{m-\beta_{k}}\phi\right)  \left(  x^{\prime
},0_{k},x^{\prime\prime}\right) \\
& =\tbinom{m}{\beta_{k}}\beta_{k}!\left(  x^{\prime}\right)  ^{\beta^{\prime}%
}\left(  x^{\prime\prime}\right)  ^{\beta^{\prime\prime}}\left(
D_{k}^{m-\beta_{k}}\phi\right)  \left(  x^{\prime},0_{k},x^{\prime\prime
}\right)  ,
\end{align*}

which means that $\left(  D_{k}^{m-\beta_{k}}\phi\right)  \left(  x^{\prime
},0_{k},x^{\prime\prime}\right)  =0$ when $\beta_{k}<m<\alpha_{k}$ i.e.
\[
\left(  D_{k}^{m}\phi\right)  \left(  x^{\prime},0_{k},x^{\prime\prime
}\right)  =0,\quad0<m<\alpha_{k}.
\]

The the second and third cases now yield%
\[
\left(  D_{k}^{m}\phi\right)  \left(  x^{\prime},0_{k},x^{\prime\prime
}\right)  =0,\quad0\leq m<\alpha_{k},
\]

so that $\phi\in S_{\otimes,\alpha-\beta}\left(  \mathbf{0}\right)  $, as
claimed.\smallskip

\textbf{Part 4} ?? \textbf{FINISH}!\smallskip

\textbf{Part 5} ?? \textbf{FINISH}! Approach: since $z\neq\mathbf{0}$ and
$z_{+}>\mathbf{0}$ we can permute $z$ so that it has the form $\left(
\mathbf{0}^{\prime},z^{\prime\prime}\right)  $ where $z_{+}^{\prime\prime
}>\mathbf{0}$. Now apply parts 3 and 4.\smallskip

\textbf{Part 6} Note that $f\phi\in S$.
\end{proof}

Lemma \ref{Lem_deriv_origin} is based on the following result which we supply
without proof.

\begin{lemma}
\label{Lem_space_Sxm}Suppose $p$ is a polynomial of degree at most $n-1$ on
$\mathbb{R}^{1}$. Then the integral $\lim\limits_{\varepsilon\rightarrow0^{+}%
}\int_{\varepsilon}^{1}\frac{p\left(  x\right)  }{x^{n}}dx$ exists iff $p=0$.
\end{lemma}

\begin{lemma}
\label{Lem_deriv_origin}Suppose $m\geq1$ is an integer and $\phi\in S\left(
\mathbb{R}^{d}\right)  $. Then:

\begin{enumerate}
\item If
\begin{equation}
\int_{0}^{1}\frac{\phi\left(  b^{\prime},x_{k},b^{\prime\prime}\right)
}{x_{k}^{2m}}dx_{k},\label{a1.21}%
\end{equation}

exists a.e. for $b=\left(  b^{\prime},b^{\prime\prime}\right)  \in
\mathbb{R}^{d-1}$ it follows that \ref{a1.21} exists for all $b\in
\mathbb{R}^{d-1}$ and that

$\left(  D_{k}^{n}\phi\right)  \left(  b^{\prime},0,b^{\prime\prime}\right)
=0$ for $n<2m$ and all $b\in\mathbb{R}^{d-1}$.

\item The integral \ref{a1.21} exists for some $b=\left(  b^{\prime}%
,b^{\prime\prime}\right)  \in\mathbb{R}^{d-1}$ iff $\left(  D_{k}^{n}%
\phi\right)  \left(  b^{\prime},0,b^{\prime\prime}\right)  =0$ for $n<2m$.
\end{enumerate}
\end{lemma}

\begin{proof}
\textbf{Part 1} Set $\psi_{b}\left(  s\right)  =\phi\left(  b^{\prime
},s,b^{\prime\prime}\right)  $. Hence $\psi_{b}\in S\left(  \mathbb{R}%
^{1}\right)  $ for any $b$. Assume $\lim\limits_{\varepsilon\rightarrow0^{+}%
}\int_{\varepsilon}^{1}\frac{\psi_{b}\left(  s\right)  }{s^{2m}}ds$ exists for
some $b$. Now using the Taylor series expansion with remainder \ref{a1.35} we
have for $\varepsilon>0$
\begin{align}
\int_{\varepsilon}^{1}\frac{\psi_{b}\left(  s\right)  }{s^{2m}}ds  &
=\int_{\varepsilon}^{1}\frac{\sum\limits_{k=0}^{2m-1}\frac{s^{k}}{k!}D^{k}%
\psi_{b}\left(  0\right)  +\frac{s^{2m}}{2m!}\left(  D^{2m}\psi_{b}\right)
\left(  t\left(  s\right)  s\right)  }{s^{2m}}ds\nonumber\\
& =\int_{\varepsilon}^{1}\frac{\sum\limits_{k=0}^{2m-1}\frac{s^{k}}{k!}%
D^{k}\psi_{b}\left(  0\right)  }{s^{2m}}ds+\frac{1}{\left(  2m\right)  !}%
\int_{\varepsilon}^{1}\left(  D^{2m}\psi_{b}\right)  \left(  t\left(
s\right)  s\right)  ds,\label{a1.001}%
\end{align}

where $0<t\left(  s\right)  <1$. Since $D^{2m}\psi_{b}$ is bounded the last
integral exists and Lemma \ref{Lem_space_Sxm} implies $D^{k}\psi_{b}\left(
0\right)  =0$ for $k<2m$ and%
\begin{align*}
\int_{0}^{1}\frac{\psi_{b}\left(  s\right)  }{s^{2m}}ds  & =\frac{1}{\left(
2m\right)  !}\int_{0}^{1}\left(  D^{2m}\psi_{b}\right)  \left(  t\left(
s\right)  s\right)  ds\\
& =\frac{1}{\left(  2m\right)  !}\int_{0}^{1}\left(  D_{k}^{2m}\psi\right)
\left(  b^{\prime},t\left(  s\right)  s,b^{\prime\prime}\right)  ds.
\end{align*}

But the last term is a continuous function of $b$ so \ref{a1.21} exists for
all $b$ and $\left(  D_{k}^{n}\phi\right)  \left(  b^{\prime},0,b^{\prime
\prime}\right)  =0$ for $n<2m$ and all $b\in\mathbb{R}^{d-1}$.\medskip

\textbf{Part 2} An easy consequence of part 1.
\end{proof}

\begin{lemma}
\label{Lem_Sxm_expansion}Suppose $\phi\in S_{\otimes,\alpha}\left(
\mathbf{0}\right)  $. Then:

\begin{enumerate}
\item If $\alpha\geq1$ then there exist functions $t_{k}:\mathbb{R}%
^{d}\rightarrow\left(  0,1\right)  $ such that%
\begin{equation}
\phi\left(  \xi\right)  =\left\{
\begin{array}
[c]{l}%
\frac{\xi^{\alpha}}{\alpha!}\left(  D^{\alpha}\phi\right)  \left(
t_{1}\left(  \xi\right)  \xi_{1},\ldots,t_{d}\left(  \xi\right)  \xi
_{d}\right)  ,\\
\frac{\xi^{\alpha}}{\alpha!}\left(  D^{\alpha}\phi\right)  \left(  t\left(
\xi\right)  .\xi\right)  ,
\end{array}
\right\}  ,\quad\xi\in\mathbb{R}^{d}.\label{a1.231}%
\end{equation}

\item If $\mathbf{0}\leq\delta\leq\mathbf{1}$ and $\alpha\geq1$ there exist
functions $t_{k}:\mathbb{R}^{d}\rightarrow\left(  0,1\right]  $ such that
$t_{k}=1$ when $\delta_{k}=0$ and $t_{k}:\mathbb{R}^{d}\rightarrow\left(
0,1\right)  $ when $\delta_{k}=1$, and such that
\begin{equation}
\phi\left(  \xi\right)  =\frac{\xi^{\alpha.\delta}}{\left(  \alpha
.\delta\right)  !}\left(  D^{\alpha.\delta}\phi\right)  \left(  t_{1}\left(
\xi\right)  \xi_{1},\ldots,t_{d}\left(  \xi\right)  \xi_{d}\right)  ,\quad
\xi\in\mathbb{R}^{d},\label{a1.24}%
\end{equation}

where $.$ is the component-wise product.

\item If $\alpha\geq0$ then there exist functions $t_{k}:\mathbb{R}%
^{d}\rightarrow\left(  0,1\right]  $ such that $t_{k}=1$ when $\alpha_{k}=0$
and $t_{k}:\mathbb{R}^{d}\rightarrow\left(  0,1\right)  $ when $\alpha_{k}>0$,
and such that%
\begin{equation}
\phi\left(  \xi\right)  =\frac{\xi^{\alpha}}{\alpha!}\left(  D^{\alpha}%
\phi\right)  \left(  t_{1}\left(  \xi\right)  \xi_{1},\ldots,t_{d}\left(
\xi\right)  \xi_{d}\right)  ,\quad\xi\in\mathbb{R}^{d}.\label{a1.23}%
\end{equation}

\end{enumerate}
\end{lemma}

\begin{proof}
\textbf{Part 1} For successive coordinates $\xi_{1},\xi_{2},\xi_{3},\ldots$ we
use the Taylor series expansion with remainder \ref{a1.35} about zero. Since
$\phi\in S_{\otimes,\alpha}\left(  0\right)  $, for the first coordinate there
exists a function $0<t_{1}\left(  \xi\right)  <1$ such that
\[
\phi\left(  \xi_{1},\xi^{\prime}\right)  =\sum\limits_{k<\alpha_{1}}\frac
{\xi_{1}^{k}}{k!}D_{1}^{k}\phi\left(  0,\xi^{\prime}\right)  +\frac{\xi
_{1}^{\alpha_{1}}}{\alpha_{1}!}D_{1}^{\alpha_{1}}\phi\left(  t_{1}\left(
\xi\right)  \xi_{1},\xi^{\prime}\right)  =\frac{\xi_{1}^{\alpha_{1}}}%
{\alpha_{1}!}D_{1}^{\alpha_{1}}\phi\left(  t_{1}\left(  \xi\right)  \xi
_{1},\xi^{\prime}\right)  .
\]

For the second coordinate%
\begin{align*}
D_{1}^{\alpha_{1}}\phi\left(  t_{1}\left(  \xi\right)  \xi_{1},\xi_{2}%
,\xi^{\prime\prime}\right)   & =\sum\limits_{k<\alpha_{2}}\frac{\xi_{2}^{k}%
}{k!}D_{1}^{\alpha_{1}}D_{2}^{k}\phi\left(  t_{1}\left(  \xi\right)  \xi
_{1},0,\xi^{\prime\prime}\right)  +\frac{\xi_{2}^{\alpha_{2}}}{\alpha_{2}%
!}D_{1}^{\alpha_{1}}D_{2}^{\alpha_{2}}\phi\left(  t_{1}\left(  \xi\right)
\xi_{1},t_{2}\left(  \xi\right)  \xi_{2},\xi^{\prime\prime}\right) \\
& =\frac{\xi_{2}^{\alpha_{2}}}{\alpha_{2}!}D_{1}^{\alpha_{1}}D_{2}^{\alpha
_{2}}\phi\left(  t_{1}\left(  \xi\right)  \xi_{1},t_{2}\left(  \xi\right)
\xi_{2},\xi^{\prime\prime}\right)  ,
\end{align*}

so that%
\[
\phi\left(  \xi\right)  =\frac{\xi_{1}^{\alpha_{1}}}{\alpha_{1}!}\frac{\xi
_{2}^{\alpha_{2}}}{\alpha_{2}!}D_{1}^{\alpha_{1}}D_{2}^{\alpha_{2}}\phi\left(
t_{1}\left(  \xi\right)  \xi_{1},t_{2}\left(  \xi\right)  \xi_{2},\xi
^{\prime\prime}\right)  ,
\]

and exhausting the coordinates yields \ref{a1.23}.\medskip

\textbf{Part 2} Equation \ref{a1.24} is just the result of only applying the
preceding Taylor series expansion to the coordinates specified by the binary
multi-index $\delta$.\medskip

\textbf{Part 3} We apply part 2. So let $\delta_{k}=\left\{
\begin{array}
[c]{ll}%
0, & \alpha_{k}=0,\\
1, & \alpha_{k}\geq1,
\end{array}
\right\}  $ and define $\beta_{k}=\left\{
\begin{array}
[c]{ll}%
1, & \alpha_{k}=0,\\
\alpha_{k}, & \alpha_{k}\geq1,
\end{array}
\right\}  $. Then $\mathbf{0}\leq\delta\leq\mathbf{1}$, $\beta\geq\mathbf{1}$
and $\beta.\delta=\alpha$ and we can apply part 2 with $\alpha$ replaced by
$\beta$.
\end{proof}

\begin{lemma}
\label{Lem_1}Suppose $m\geq1$ is an integer and $\psi\in S\left(
\mathbb{R}^{d}\right)  $. Then:

\begin{enumerate}
\item If the integral%
\begin{equation}
\int_{0}^{1}\frac{\left\vert \psi\left(  b^{\prime},x_{k},b^{\prime\prime
}\right)  \right\vert ^{2}}{x_{k}^{2m}}dx_{k},\label{a1.47}%
\end{equation}

exists for almost all $b=\left(  b^{\prime},b^{\prime\prime}\right)
\in\mathbb{R}^{d-1}$ it follows that \ref{a1.47} exists for all $b\in
\mathbb{R}^{d-1}$ and that $\left(  D_{k}^{n}\psi\right)  \left(  b^{\prime
},0,b^{\prime\prime}\right)  =0$ for $n<m$ and all $b\in\mathbb{R}^{d-1}$.

\item The integral \ref{a1.47} exists for some $b=\left(  b^{\prime}%
,b^{\prime\prime}\right)  \in\mathbb{R}^{d-1}$ iff $\left(  D_{k}^{n}%
\psi\right)  \left(  b^{\prime},0,b^{\prime\prime}\right)  =0$ for $n<m$.
\end{enumerate}
\end{lemma}

\begin{proof}
\textbf{Part 1} From Lemma \ref{Lem_deriv_origin}, if \ref{a1.47} exists for
almost all $b$ it exists for all $b\in\mathbb{R}^{d-1}$ and

$\left(  D_{k}^{n}\left(  \left\vert \psi\right\vert ^{2}\right)  \right)
\left(  b^{\prime},0,b^{\prime\prime}\right)  =0$ for $n<m$ and all
$b\in\mathbb{R}^{d-1}$. Lemma \ref{Thm_abs_squar_and_Som} now implies that
$\left(  D_{k}^{n}\psi\right)  =0$ for $n<m$ and all $b\in\mathbb{R}^{d-1}%
$.\medskip$\left(  b^{\prime},0,b^{\prime\prime}\right)  $

\textbf{Part 2} If \ref{a1.47} exists then part 1 implies $\left(  D_{k}%
^{n}\psi\right)  \left(  b^{\prime},0,b^{\prime\prime}\right)  =0$ for $n<m$.
The converse is easily proved using a Taylor series expansion.
\end{proof}

The next result characterizes the spaces $S_{\otimes,m}$.

\begin{lemma}
\label{Lem_2}Suppose $\psi\in S\left(  \mathbb{R}^{d}\right)  $ and
$\alpha\geq0$. Then the integral%
\begin{equation}
\int_{\mathbf{0}}^{\mathbf{1}}\frac{\left\vert \psi\left(  x\right)
\right\vert ^{2}}{x^{2\alpha}}dx:=\int_{0}^{1}\ldots\int_{0}^{1}%
\frac{\left\vert \psi\left(  x\right)  \right\vert ^{2}}{x^{2\alpha}%
}dx,\label{a1.48}%
\end{equation}

exists iff $\psi\in S_{\otimes,\alpha}\left(  \mathbf{0}\right)  $.
\end{lemma}

\begin{proof}
Suppose that $\psi\in S_{\otimes,\alpha}\left(  \mathbf{0}\right)  $. Then
equation \ref{a1.23} implies that the integral \ref{a1.48} exists. Conversely,
if \ref{a1.48} exists we have%
\[
\int_{\mathbf{0}}^{\mathbf{1}}\frac{\left\vert \psi\left(  x\right)
\right\vert ^{2}}{x^{2\alpha}}dx\geq\int_{\mathbf{0}^{\prime\prime}%
}^{\mathbf{1}^{\prime\prime}}\int_{0}^{1}\int_{\mathbf{0}^{\prime}%
}^{\mathbf{1}^{\prime}}\frac{\left\vert \psi\left(  x^{\prime},x_{k}%
,x^{\prime\prime}\right)  \right\vert ^{2}}{x_{k}^{2\alpha_{k}}}dx^{\prime
}dx_{k}dx^{\prime\prime},
\]
which shows that $\int_{0}^{1}\frac{\left\vert \psi\left(  x^{\prime}%
,x_{k},x^{\prime\prime}\right)  \right\vert ^{2}}{x_{k}^{2\alpha_{k}}}dx_{k}$
exists for almost all $\left(  x^{\prime},x^{\prime\prime}\right)
\in\mathbb{R}^{d-1}$. Lemma \ref{Lem_1} now implies that

$\int_{0}^{1}\frac{\left\vert \psi\left(  x^{\prime},x_{k},x^{\prime\prime
}\right)  \right\vert ^{2}}{x_{k}^{2\alpha_{k}}}dx_{k}$ exists for all
$\left(  x^{\prime},x^{\prime\prime}\right)  \in\mathbb{R}^{d-1}$ and that
$\left(  D_{k}^{n}\psi\right)  \left(  x^{\prime},0,x^{\prime\prime}\right)
=0$ for $n<\alpha_{k}$ and all $\left(  x^{\prime},x^{\prime\prime}\right)
\in\mathbb{R}^{d-1}$. Since we can choose $k=1,\ldots,d$ it follows that
$\psi\in S_{\otimes,\alpha}\left(  \mathbf{0}\right)  $.
\end{proof}

The next result characterizes the spaces $S_{\otimes,\alpha}$ in terms an
integral-form Taylor series remainder.

\begin{lemma}
\label{Lem_Sxm_integ_expansion}\ Suppose $\phi\in S$. Then:

\begin{enumerate}
\item $\phi\in S_{\otimes,\alpha}\left(  \mathbf{0}\right)  $ for some
$\alpha\geq\mathbf{1}$ iff%
\begin{equation}
\left.
\begin{array}
[c]{ll}%
\phi\left(  \xi\right)  & =\frac{\xi^{\alpha}}{\left(  \alpha-1\right)  !}%
\int_{\mathbf{0}}^{\mathbf{1}}\left(  1-\tau_{1}\right)  ^{\alpha_{1}-1}%
\ldots\left(  1-\tau_{d}\right)  ^{\alpha_{d}-1}\left(  D^{\alpha}\phi\right)
\left(  \tau_{1}\xi_{1},\ldots,\tau_{d}\xi_{d}\right)  d\tau_{1}\ldots
d\tau_{d}\\
& =\frac{\xi^{\alpha}}{\left(  \alpha-1\right)  !}\int_{\mathbf{0}%
}^{\mathbf{1}}\left(  1-\tau\right)  ^{\alpha-1}\left(  D^{\alpha}\phi\right)
\left(  \tau.\xi\right)  d\tau
\end{array}
\right\}  ,\label{1.076}%
\end{equation}

where the second equation uses component-wise notation. See part 7.

\item Suppose $\alpha=\left(  \alpha^{\prime},\mathbf{0}^{\prime\prime
}\right)  $ where $\alpha^{\prime}\geq\mathbf{1}^{\prime}$. Then $\phi\in
S_{\otimes,\alpha}\left(  \mathbf{0}\right)  $ iff%
\[
\phi\left(  \xi\right)  =\frac{\left(  \xi^{\prime}\right)  ^{\alpha^{\prime}%
}}{\left(  \alpha^{\prime}-1\right)  !}\int_{\mathbf{0}^{\prime}}%
^{\mathbf{1}^{\prime}}\left(  \mathbf{1}^{\prime}-\tau^{\prime}\right)
^{\alpha^{\prime}-1}\left(  D^{\alpha^{\prime}}\phi\right)  \left(
\tau^{\prime}.\xi^{\prime},\xi^{\prime\prime}\right)  d\tau^{\prime}.
\]

\item If $\alpha\geq\mathbf{0}$ then $\phi\in S_{\otimes,\alpha}\left(
\mathbf{0}\right)  $ iff%
\begin{equation}
\phi\left(  \xi\right)  =\frac{\xi^{\alpha}}{\left(  \alpha-1\right)  !}%
\int_{\mathbf{0}}^{\mathbf{1}}\left(  \mathbf{1}-\tau\right)  ^{\alpha-\left(
\alpha.>\mathbf{0}\right)  }\left(  D^{\alpha}\phi\right)  \left(
\tau^{\left(  \alpha.>\mathbf{0}\right)  }.\xi\right)  d\tau,\label{1.080}%
\end{equation}

where $\left(  -1\right)  !:=1$; the \textbf{logical operator} $\left(
\alpha.>\mathbf{0}\right)  $ is defined by: $\left(  \alpha.>\mathbf{0}%
\right)  _{k}=1$ if $\alpha_{k}>0$ and $\left(  \alpha>\mathbf{0}\right)
_{k}=0$ if $\alpha_{k}=0$.

\item If $\alpha\geq\mathbf{0}$ and $z\in\mathbb{R}^{d}$ then $\phi\in
S_{\otimes,\alpha}\left(  z\right)  $ iff%
\[
\phi\left(  \xi\right)  =\frac{\left(  \xi-z\right)  ^{\alpha}}{\left(
\alpha-1\right)  !}\int\limits_{\mathbf{0}}^{\mathbf{1}}\left(  \mathbf{1}%
-\tau\right)  ^{\alpha-\left(  \alpha.>\mathbf{0}\right)  }\left(  D^{\alpha
}\phi\right)  \left(  z+\tau^{\left(  \alpha.>\mathbf{0}\right)  }.\left(
\xi-z\right)  \right)  d\tau.
\]

Further, when $\beta\leq\alpha$,%
\[
D^{\beta}\phi\left(  \xi\right)  =\frac{\left(  \xi-z\right)  ^{\alpha-\beta}%
}{\left(  \alpha-\beta-1\right)  !}\int\limits_{\mathbf{0}}^{\mathbf{1}%
}\left(  \mathbf{1}-\tau\right)  ^{\alpha-\beta-\left(  \alpha.>\beta\right)
}\left(  D^{\alpha}\phi\right)  \left(  z+\tau^{\left(  \alpha.>\beta\right)
}.\left(  \xi-z\right)  \right)  d\tau.
\]

\item If $\phi\in S_{\otimes,\alpha}\left(  \mathbf{0}\right)  $ for some
$\alpha\geq\mathbf{1}$ then%
\[
\left\vert \phi\left(  \xi\right)  \right\vert \leq\frac{\xi_{+}^{\alpha}%
}{\alpha!}\max_{R\left[  \mathbf{0},\xi\right]  }\left\vert D^{\alpha}%
\phi\left(  \cdot\right)  \right\vert ,\quad\xi\in\mathbb{R}^{d},
\]

where $R\left[  \mathbf{0},\xi\right]  =R\left[  \min\left\{  \mathbf{0}%
,\xi\right\}  ,\max\left\{  \mathbf{0},\xi\right\}  \right]  =\left\{
\tau.\xi:\mathbf{0}\leq\tau\leq\mathbf{1}\right\}  $ - see Definition
\ref{Def_topol_on_Rd} of the Appendix.

\item If $\phi\in S_{\otimes,\alpha}\left(  z\right)  $ for some $\alpha
\geq\mathbf{1}$ then%
\[
\left\vert \phi\left(  \xi\right)  \right\vert \leq\frac{\left\vert \left(
\xi-z\right)  ^{\alpha}\right\vert }{\alpha!}\max_{R\left[  z,\xi\right]
}\left\vert D^{\alpha}\phi\left(  \cdot\right)  \right\vert ,\quad\xi
\in\mathbb{R}^{d}.
\]

\item If $\phi\in S_{\otimes,\alpha}\left(  \mathbf{0}\right)  $ for some
$\alpha\geq\mathbf{1}$ then%
\[
\phi\left(  \xi\right)  =\frac{\left(  \operatorname*{sgn}\xi\right)
^{\mathbf{1}}}{\left(  \alpha-1\right)  !}\int_{\min\left\{  \mathbf{0}%
,\xi\right\}  }^{\max\left\{  \mathbf{0},\xi\right\}  }\left(  \xi
-\eta\right)  ^{\alpha-1}D^{\alpha}\phi\left(  \eta\right)  d\eta,
\]

where, $R\left[  \mathbf{0},\xi\right]  =R\left(  \min\left\{  \mathbf{0}%
,\xi\right\}  ,\max\left\{  \mathbf{0},\xi\right\}  \right)  $. Further, for
all $\xi\in\mathbb{R}^{d}$ we have the $L^{1}$ estimate%
\[
\left\vert \frac{\phi\left(  \xi\right)  }{\xi^{\alpha}}\right\vert \leq
\frac{1}{\left(  \alpha-1\right)  !}\frac{1}{\left\vert \xi^{\mathbf{1}%
}\right\vert }\int_{\min\left\{  \mathbf{0},\xi\right\}  }^{\max\left\{
\mathbf{0},\xi\right\}  }\left\vert D^{\alpha}\phi\right\vert ,
\]

and the $L^{2}$ estimate%
\[
\left\vert \frac{\phi\left(  \xi\right)  }{\xi^{\alpha}}\right\vert \leq
\frac{1}{\sqrt{2\alpha-1}\left(  \alpha-1\right)  !}\left(  \frac
{1}{\left\vert \xi^{\mathbf{1}}\right\vert }\int_{\min\left\{  \mathbf{0}%
,\xi\right\}  }^{\max\left\{  \mathbf{0},\xi\right\}  }\left\vert D^{\alpha
}\phi\right\vert ^{2}\right)  ^{1/2}.
\]

Here $\int_{\min\left\{  \mathbf{0},\xi\right\}  }^{\max\left\{
\mathbf{0},\xi\right\}  }=\int_{R\left[  \mathbf{0},\xi\right]  }$. See part 5.

\item If $\phi\in S_{\otimes,\alpha}\left(  z\right)  $ for some $\alpha
\geq\mathbf{1}$ then for all $\xi\in\mathbb{R}^{d}$,%
\begin{align*}
\phi\left(  \xi\right)   & =\frac{\left(  \zeta-z\right)  ^{\alpha}}{\left(
\alpha-1\right)  !}\int_{\mathbf{0}}^{\mathbf{1}}\left(  1-\tau\right)
^{\alpha-1}\left(  D^{\alpha}\phi\right)  \left(  \tau.\zeta+\left(
1-\tau\right)  .z\right)  d\tau\\
& =\frac{\left(  \operatorname*{sgn}\left(  \xi-z\right)  \right)
^{\mathbf{1}}}{\left(  \alpha-1\right)  !}\int_{\min\left\{  z,\xi\right\}
}^{\max\left\{  z,\xi\right\}  }\left(  \xi-\sigma\right)  ^{\alpha
-1}D^{\alpha}\phi\left(  \sigma\right)  d\sigma,
\end{align*}

and%
\[
\left\vert \frac{\phi\left(  \xi\right)  }{\left(  \xi-z\right)  ^{\alpha}%
}\right\vert \leq\frac{1}{\left(  \alpha-1\right)  !}\frac{1}{\left\vert
\left(  \xi-z\right)  ^{\mathbf{1}}\right\vert }\int_{\min\left\{
z,\xi\right\}  }^{\max\left\{  z,\xi\right\}  }\left\vert D^{\alpha}%
\phi\right\vert ,
\]

and%
\[
\left\vert \frac{\phi\left(  \xi\right)  }{\left(  \xi-z\right)  ^{\alpha}%
}\right\vert \leq\frac{1}{\sqrt{2\alpha-1}\left(  \alpha-1\right)  !}\left(
\frac{1}{\left\vert \left(  \xi-z\right)  ^{\mathbf{1}}\right\vert }%
\int\limits_{\min\left\{  z,\xi\right\}  }^{\max\left\{  z,\xi\right\}
}\left\vert D^{\alpha}\phi\right\vert ^{2}\right)  ^{1/2}.
\]

Here $\int_{\min\left\{  z,\xi\right\}  }^{\max\left\{  z,\xi\right\}  }%
=\int_{R\left[  z,\xi\right]  }$. See part 5.
\end{enumerate}
\end{lemma}

\begin{proof}
\textbf{Part 1} Since $\phi\in S_{\otimes,\alpha}\left(  \mathbf{0}\right)  $
implies $D_{1}^{k}\phi\left(  0_{1},\cdot\right)  =0$ for $k<\alpha_{1}$, we
have%
\begin{align*}
\phi\left(  \xi_{1},\xi^{\prime}\right)   & =\sum\limits_{k<\alpha_{1}}%
\frac{\xi_{1}^{k}}{k!}D_{1}^{k}\phi\left(  0,\xi^{\prime}\right)  +\frac
{\xi_{1}^{\alpha_{1}}}{\left(  \alpha_{1}-1\right)  !}\int_{0}^{1}\left(
1-\tau_{1}\right)  ^{\alpha_{1}-1}\left(  D_{1}^{\alpha_{1}}\phi\right)
\left(  \tau_{1}\xi_{1},\xi^{\prime}\right)  d\tau_{1}\\
& =\frac{\xi_{1}^{\alpha_{1}}}{\left(  \alpha_{1}-1\right)  !}\int_{0}%
^{1}\left(  1-\tau_{1}\right)  ^{\alpha_{1}-1}\left(  D_{1}^{\alpha_{1}}%
\phi\right)  \left(  \tau_{1}\xi_{1},\xi^{\prime}\right)  d\tau_{1}.
\end{align*}

Since $\psi\in S_{\otimes,\alpha}\left(  \mathbf{0}\right)  $, we have
$\left(  D_{1}^{\alpha_{1}}\phi\right)  \left(  \tau_{1}\xi_{1},\xi^{\prime
}\right)  =D_{2}^{k}\phi\left(  \tau_{1}\xi_{1},0_{1},\xi^{\prime\prime
}\right)  =0$ for $k<\alpha_{2}$, which implies%
\[
\left(  D_{1}^{\alpha_{1}}\phi\right)  \left(  \tau_{1}\xi_{1},\xi_{2}%
,\xi^{\prime\prime}\right)  =\frac{\xi_{2}^{\alpha_{1}}}{\left(  \alpha
_{1}-1\right)  !}\int_{0}^{1}\left(  1-\tau_{2}\right)  ^{\alpha_{1}-1}\left(
D_{1}^{\alpha_{1}}D_{2}^{\alpha_{2}}\phi\right)  \left(  \tau_{1}\xi_{1}%
,\tau_{2}\xi_{2},\xi^{\prime\prime}\right)  d\tau_{2},
\]

and so%
\begin{align*}
&  \phi\left(  \xi_{1},\xi_{2},\xi^{\prime\prime}\right) \\
&  =\frac{\xi_{1}^{\alpha_{1}}}{\left(  m-1\right)  !}\int_{0}^{1}\left(
1-\tau_{1}\right)  ^{\alpha_{1}-1}\frac{\xi_{2}^{\alpha_{2}}}{\left(
m-1\right)  !}\int_{0}^{1}\left(  1-\tau_{2}\right)  ^{\alpha_{2}-1}\left(
D_{1}^{\alpha_{1}}D_{2}^{\alpha_{2}}\phi\right)  \left(  \tau_{1}\xi_{1}%
,\tau_{2}\xi_{2},\xi^{\prime\prime}\right)  d\tau_{2}d\tau_{1}\\
&  =\frac{\xi_{1}^{\alpha_{1}}\xi_{2}^{\alpha_{2}}}{\left(  \alpha
_{1}-1\right)  !\left(  \alpha_{2}-1\right)  !}\int_{0}^{1}\int_{0}^{1}\left(
1-\tau_{1}\right)  ^{\alpha_{1}-1}\left(  1-\tau_{2}\right)  ^{\alpha_{2}%
-1}\left(  D_{1}^{\alpha_{1}}D_{2}^{\alpha_{2}}\phi\right)  \left(  \tau
_{1}\xi_{1},\tau_{2}\xi_{2},\xi^{\prime\prime}\right)  d\tau_{1}d\tau_{2},
\end{align*}

and%
\begin{equation}
\phi\left(  \xi\right)  =\frac{\left(  \xi^{\prime}\right)  ^{\alpha^{\prime}%
}}{\left(  \alpha^{\prime}-\mathbf{1}^{\prime}\right)  !}\int_{\mathbf{0}%
^{\prime}}^{\mathbf{1}^{\prime}}\left(  \mathbf{1}^{\prime}-\tau^{\prime
}\right)  ^{\alpha^{\prime}-1}\left(  D^{\alpha^{\prime}}\phi\right)  \left(
\tau^{\prime}.\xi^{\prime},\xi^{\prime\prime}\right)  d\tau^{\prime
},\label{1.034}%
\end{equation}

and continuing in this manner until exhaustion proves \ref{1.076}.

Conversely, if \ref{1.076} is true then $\left\vert \phi\left(  \xi\right)
\right\vert \leq\frac{\left\vert \xi^{\alpha}\right\vert }{\alpha!}\left\Vert
D^{\alpha}\phi\right\Vert _{\infty}$ and so $\int_{\mathbf{0}}^{\mathbf{1}%
}\frac{\left\vert \phi\left(  \xi\right)  \right\vert ^{2}}{\xi^{2\alpha}}%
d\xi<\infty$ and Lemma \ref{Lem_2} implies $\phi\in S_{\otimes,\alpha}\left(
\mathbf{0}\right)  $.\medskip

\textbf{Part 2} Examining part 1 it is clear that in this case \ref{1.034}
holds.\medskip

\textbf{Part 3} Choose a permutation $\pi$ such that $\pi\alpha=\beta=\left(
\beta^{\prime},\mathbf{0}^{\prime\prime}\right)  $ with $\beta^{\prime}%
\geq\mathbf{1}^{\prime}$. Now set $\pi\xi=\eta$ and $\psi\left(  \eta\right)
=\phi\left(  \xi\right)  $ which means that $\psi\in S_{\otimes,\beta}\left(
\mathbf{0}\right)  $. From part 2, $\psi\in S_{\otimes,\beta}\left(
\mathbf{0}\right)  $ iff%
\begin{align*}
\psi\left(  \eta\right)   & =\frac{\left(  \eta^{\prime}\right)
^{\beta^{\prime}}}{\left(  \beta^{\prime}-1\right)  !}\int_{\mathbf{0}%
^{\prime}}^{\mathbf{1}^{\prime}}\left(  \mathbf{1}-\sigma^{\prime}\right)
^{\beta^{\prime}-1}\left(  D^{\beta^{\prime}}\psi\right)  \left(
\sigma^{\prime}.\eta^{\prime},\eta^{\prime\prime}\right)  d\sigma^{\prime}\\
& =\frac{\eta^{\beta}}{\left(  \beta-1\right)  !}\int_{\mathbf{0}}%
^{\mathbf{1}}\left(  \mathbf{1}-\sigma\right)  ^{\beta-\left(  \beta
.>\mathbf{0}\right)  }\left(  D^{\beta}\psi\right)  \left(  \sigma^{\left(
\beta>\mathbf{0}\right)  }.\eta\right)  d\sigma.
\end{align*}

The change of variables $\tau=\pi\sigma$, $d\tau=d\sigma$ yields%
\begin{align*}
\psi\left(  \eta\right)   & =\frac{\eta^{\beta}}{\left(  \beta-1\right)
!}\int_{\mathbf{0}}^{\mathbf{1}}\left(  \mathbf{1}-\pi\tau\right)
^{\beta-\left(  \beta.>\mathbf{0}\right)  }\left(  D^{\beta}\psi\right)
\left(  \left(  \pi\tau\right)  ^{\left(  \beta.>\mathbf{0}\right)  }%
.\eta\right)  d\tau\\
& =\frac{\xi^{\alpha}}{\left(  \alpha-1\right)  !}\int_{\mathbf{0}%
}^{\mathbf{1}}\left(  \mathbf{1}-\tau\right)  ^{\alpha-\left(  \alpha
.>\mathbf{0}\right)  }\left(  D^{\alpha}\psi\right)  \left(  \tau^{\left(
\alpha.>\mathbf{0}\right)  }.\xi\right)  d\tau\\
& =\phi\left(  \xi\right)  .
\end{align*}

Conversely, suppose \ref{1.080} holds. Then%
\begin{align*}
\left\vert \phi\left(  \xi\right)  \right\vert  & \leq\frac{\left\vert
\xi^{\alpha}\right\vert }{\left(  \alpha-1\right)  !}\int_{\mathbf{0}%
}^{\mathbf{1}}\left(  1-\tau\right)  ^{\alpha-\left(  \alpha.>\mathbf{0}%
\right)  }\left\vert \left(  D^{\alpha}\phi\right)  \left(  \tau^{\left(
\alpha.>\mathbf{0}\right)  }.\xi\right)  \right\vert d\tau\\
& \leq\frac{\left\vert \xi^{\alpha}\right\vert }{\left(  \alpha-1\right)
!}\left\Vert D^{\alpha}\phi\right\Vert _{\infty},
\end{align*}

and so $\int_{\mathbf{0}}^{\mathbf{1}}\frac{\left\vert \phi\left(  \xi\right)
\right\vert ^{2}}{\xi^{\alpha}}d\xi<\infty$. An application of Lemma
\ref{Lem_2} now implies $\phi\in S_{\otimes,\alpha}\left(  \mathbf{0}\right)
$.\medskip

\textbf{Part 4} From Definition \ref{Def_space_Sxm}, $\phi\in S_{\otimes
,\alpha}\left(  z\right)  $ iff $\phi\left(  \cdot+z\right)  \in
S_{\otimes,\alpha}\left(  \mathbf{0}\right)  $, and $\phi\in S_{\otimes
,\alpha}\left(  z\right)  $ implies $D^{\beta}\phi\in S_{\otimes,\alpha-\beta
}\left(  z\right)  $.\medskip

\textbf{Part 5} From part 1,
\begin{align*}
\left\vert \phi\left(  \xi\right)  \right\vert  & \leq\frac{\xi_{+}^{\alpha}%
}{\left(  \alpha-1\right)  !}\int_{\mathbf{0}}^{\mathbf{1}}\left(
1-\tau\right)  ^{\alpha-1}\left\vert \left(  D^{\alpha}\phi\right)  \left(
\tau.\xi\right)  \right\vert d\tau\\
& \leq\frac{\xi_{+}^{\alpha}}{\left(  \alpha-1\right)  !}\max_{R\left[
\mathbf{0},\xi\right]  }\left\vert D^{\alpha}\phi\right\vert \int_{\mathbf{0}%
}^{\mathbf{1}}\left(  1-\tau\right)  ^{\alpha-1}d\tau\\
& =\frac{\xi_{+}^{\alpha}}{\left(  \alpha-1\right)  !}\max_{R\left[
\mathbf{0},\xi\right]  }\left\vert D^{\alpha}\phi\right\vert \int_{\mathbf{0}%
}^{\mathbf{1}}\tau^{\alpha-1}d\tau\\
& =\frac{\xi_{+}^{\alpha}}{\alpha!}\max_{R\left[  \mathbf{0},\xi\right]
}\left\vert D^{\alpha}\phi\right\vert .
\end{align*}
\medskip

\textbf{Part 6} From part 5, if $\phi\in S_{\otimes,\alpha}\left(
\mathbf{0}\right)  $ for some $\alpha\geq\mathbf{1}$ then%
\[
\left\vert \phi\left(  \xi\right)  \right\vert \leq\frac{\xi_{+}^{\alpha}%
}{\alpha!}\max_{\eta\in R\left[  \mathbf{0},\xi\right]  }\left\vert D^{\alpha
}\phi\left(  \eta\right)  \right\vert ,\quad\xi\in\mathbb{R}^{d}.
\]

From Definition \ref{Def_space_Sxm}, $\psi\in S_{\otimes,\alpha}\left(
z\right)  $ iff $\psi\left(  z+\cdot\right)  \in S_{\otimes,\alpha}\left(
\mathbf{0}\right)  $ and so%
\begin{align*}
\left\vert \psi\left(  \xi+z\right)  \right\vert  & \leq\frac{\xi_{+}^{\alpha
}}{\alpha!}\max_{\eta\in R\left[  \mathbf{0},\xi\right]  }\left\vert
D^{\alpha}\psi\left(  \eta+z\right)  \right\vert ,\quad\xi\in\mathbb{R}^{d}.\\
\left\vert \psi\left(  \zeta\right)  \right\vert  & \leq\frac{\left(
\zeta-z\right)  _{+}^{\alpha}}{\alpha!}\max_{\eta\in R\left[  \mathbf{0}%
,\zeta-z\right]  }\left\vert D^{\alpha}\psi\left(  \eta+z\right)  \right\vert
,\quad\zeta\in\mathbb{R}^{d}.\\
& =\frac{\left(  \zeta-z\right)  _{+}^{\alpha}}{\alpha!}\max_{\eta\in R\left[
z,\zeta\right]  }\left\vert D^{\alpha}\psi\left(  \eta\right)  \right\vert
,\quad\zeta\in\mathbb{R}^{d}.
\end{align*}
\medskip

\textbf{Part 7} From part 1,%
\[
\phi\left(  \xi\right)  =\frac{\xi^{\alpha}}{\left(  \alpha-1\right)  !}%
\int_{\mathbf{0}}^{\mathbf{1}}\left(  1-\tau\right)  ^{\alpha-1}\left(
D^{\alpha}\phi\right)  \left(  \tau.\xi\right)  d\tau.
\]

Applying the change of variables: $\eta=\tau.\xi$, $d\eta=\left\vert
\xi^{\mathbf{1}}\right\vert d\tau$ gives for $\xi.\neq\mathbf{0}$,%
\begin{align*}
\phi\left(  \xi\right)   & =\frac{\xi^{\alpha}}{\left(  \alpha-1\right)
!}\int_{\mathbf{0}}^{\mathbf{1}}\left(  1-\frac{\eta}{\xi}\right)  ^{\alpha
-1}D^{\alpha}\phi\left(  \eta\right)  \frac{d\eta}{\left\vert \xi^{\mathbf{1}%
}\right\vert }\\
& =\frac{1}{\left(  \alpha-1\right)  !}\frac{\xi^{\mathbf{1}}}{\left\vert
\xi^{\mathbf{1}}\right\vert }\int_{R\left[  \mathbf{0},\xi\right]  }\left(
\xi-\eta\right)  ^{\alpha-1}D^{\alpha}\phi\left(  \eta\right)  d\eta\\
& =\frac{\left(  \operatorname*{sgn}\xi\right)  ^{\mathbf{1}}}{\left(
\alpha-1\right)  !}\int_{R\left[  \mathbf{0},\xi\right]  }\left(  \xi
-\eta\right)  ^{\alpha-1}D^{\alpha}\phi\left(  \eta\right)  d\eta\\
& =\frac{\left(  \operatorname*{sgn}\xi\right)  ^{\mathbf{1}}}{\left(
\alpha-1\right)  !}\int_{\min\left\{  \mathbf{0},\xi\right\}  }^{\max\left\{
\mathbf{0},\xi\right\}  }\left(  \xi-\eta\right)  ^{\alpha-1}D^{\alpha}%
\phi\left(  \eta\right)  d\eta,
\end{align*}

where $R\left[  \mathbf{0},\xi\right]  =R\left(  \min\left\{  \mathbf{0}%
,\xi\right\}  ,\max\left\{  \mathbf{0},\xi\right\}  \right)  $. Further%
\begin{align}
\left\vert \phi\left(  \xi\right)  \right\vert  & \leq\frac{1}{\left(
\alpha-1\right)  !}\int\limits_{\min\left\{  \mathbf{0},\xi\right\}  }%
^{\max\left\{  \mathbf{0},\xi\right\}  }\left\vert \left(  \xi-\eta\right)
^{\alpha-1}\right\vert \left\vert D^{\alpha}\phi\left(  \eta\right)
\right\vert d\eta\label{1.081}\\
& \leq\frac{\left\vert \xi^{\alpha-1}\right\vert }{\left(  \alpha-1\right)
!}\int\limits_{\min\left\{  \mathbf{0},\xi\right\}  }^{\max\left\{
\mathbf{0},\xi\right\}  }\left\vert D^{\alpha}\phi\right\vert ,\nonumber
\end{align}

and so%
\[
\left\vert \frac{\phi\left(  \xi\right)  }{\xi^{\alpha}}\right\vert \leq
\frac{1}{\left(  \alpha-1\right)  !}\frac{1}{\left\vert \xi^{\mathbf{1}%
}\right\vert }\int\limits_{\min\left\{  \mathbf{0},\xi\right\}  }%
^{\max\left\{  \mathbf{0},\xi\right\}  }\left\vert D^{\alpha}\phi\right\vert .
\]

Also, from \ref{1.081},%
\[
\left\vert \phi\left(  \xi\right)  \right\vert \leq\frac{1}{\left(
\alpha-1\right)  !}\left(  \int\limits_{\min\left\{  \mathbf{0},\xi\right\}
}^{\max\left\{  \mathbf{0},\xi\right\}  }\left(  \xi-\eta\right)  ^{2\left(
\alpha-1\right)  }d\eta\right)  ^{1/2}\left(  \int\limits_{\min\left\{
\mathbf{0},\xi\right\}  }^{\max\left\{  \mathbf{0},\xi\right\}  }\left\vert
D^{\alpha}\phi\right\vert ^{2}\right)  ^{1/2},
\]

where%
\begin{align*}
\int\limits_{\min\left\{  \mathbf{0},\xi\right\}  }^{\max\left\{
\mathbf{0},\xi\right\}  }\left(  \xi-\eta\right)  ^{2\left(  \alpha-1\right)
}d\eta & =\int\limits_{\min\left\{  \mathbf{0},\xi\right\}  }^{\max\left\{
\mathbf{0},\xi\right\}  }\left(  \eta-\xi\right)  ^{2\left(  \alpha-1\right)
}d\eta=\int\limits_{\min\left\{  -\xi,\mathbf{0}\right\}  }^{\max\left\{
-\xi,\mathbf{0}\right\}  }\zeta^{2\left(  \alpha-1\right)  }d\zeta=\\
& =\int\limits_{\min\left\{  \xi,\mathbf{0}\right\}  }^{\max\left\{
\xi,\mathbf{0}\right\}  }\zeta^{2\left(  \alpha-1\right)  }d\zeta
=\frac{\left\vert \xi^{2\alpha-1}\right\vert }{2\alpha-1},
\end{align*}

so that%
\begin{align*}
\left\vert \phi\left(  \xi\right)  \right\vert  & \leq\frac{1}{\left(
\alpha-1\right)  !}\frac{\left\vert \xi^{2\alpha-1}\right\vert ^{1/2}}{\left(
2\alpha-1\right)  ^{\mathbf{1}/2}}\left(  \int\limits_{\min\left\{
\mathbf{0},\xi\right\}  }^{\max\left\{  \mathbf{0},\xi\right\}  }\left\vert
D^{\alpha}\phi\right\vert ^{2}\right)  ^{1/2}\\
& =\frac{1}{\left(  \alpha-1\right)  !}\frac{\left\vert \xi^{\alpha
}\right\vert }{\left(  2\alpha-1\right)  ^{\mathbf{1}/2}}\left(  \frac
{1}{\left\vert \xi^{\mathbf{1}}\right\vert }\int\limits_{\min\left\{
\mathbf{0},\xi\right\}  }^{\max\left\{  \mathbf{0},\xi\right\}  }\left\vert
D^{\alpha}\phi\right\vert ^{2}\right)  ^{1/2},
\end{align*}

i.e.%
\[
\left\vert \frac{\phi\left(  \xi\right)  }{\xi^{\alpha}}\right\vert \leq
\frac{1}{\left(  2\alpha-1\right)  ^{\mathbf{1}/2}\left(  \alpha-1\right)
!}\left(  \frac{1}{\left\vert \xi^{\mathbf{1}}\right\vert }\int\limits_{\min
\left\{  \mathbf{0},\xi\right\}  }^{\max\left\{  \mathbf{0},\xi\right\}
}\left\vert D^{\alpha}\phi\right\vert ^{2}\right)  ^{1/2},
\]

as claimed.\medskip

\textbf{Part 8} From Definition \ref{Def_space_Sxm}, $\phi\in S_{\otimes
,\alpha}\left(  z\right)  $ iff $\phi\left(  \cdot+z\right)  \in
S_{\otimes,\alpha}\left(  \mathbf{0}\right)  $. So from part 1,%
\begin{align*}
\phi\left(  \xi+z\right)   & =\frac{\xi^{\alpha}}{\left(  \alpha-1\right)
!}\int_{\mathbf{0}}^{\mathbf{1}}\left(  1-\tau\right)  ^{\alpha-1}\left(
D^{\alpha}\phi\right)  \left(  \tau.\xi+z\right)  d\tau.\\
\phi\left(  \zeta\right)   & =\frac{\left(  \zeta-z\right)  ^{\alpha}}{\left(
\alpha-1\right)  !}\int_{\mathbf{0}}^{\mathbf{1}}\left(  1-\tau\right)
^{\alpha-1}\left(  D^{\alpha}\phi\right)  \left(  \tau.\left(  \zeta-z\right)
+z\right)  d\tau\\
& =\frac{\left(  \zeta-z\right)  ^{\alpha}}{\left(  \alpha-1\right)  !}%
\int_{\mathbf{0}}^{\mathbf{1}}\left(  1-\tau\right)  ^{\alpha-1}\left(
D^{\alpha}\phi\right)  \left(  \tau.\zeta+\left(  1-\tau\right)  .z\right)
d\tau,
\end{align*}

and from part 7%
\begin{align*}
\phi\left(  \xi+z\right)   & =\frac{\left(  \operatorname*{sgn}\xi\right)
^{\mathbf{1}}}{\left(  \alpha-1\right)  !}\int_{\min\left\{  \mathbf{0}%
,\xi\right\}  }^{\max\left\{  \mathbf{0},\xi\right\}  }\left(  \xi
-\eta\right)  ^{\alpha-1}D^{\alpha}\phi\left(  \eta+z\right)  d\eta.\\
\phi\left(  \zeta\right)   & =\frac{\left(  \operatorname*{sgn}\left(
\zeta-z\right)  \right)  ^{\mathbf{1}}}{\left(  \alpha-1\right)  !}\int%
_{\min\left\{  \mathbf{0},\eta-z\right\}  }^{\max\left\{  \mathbf{0}%
,\eta-z\right\}  }\left(  \zeta-z-\eta\right)  ^{\alpha-1}D^{\alpha}%
\phi\left(  \eta+z\right)  d\eta\\
& :\sigma=\eta+z,\text{ }d\sigma=d\eta\Rightarrow\\
& =\frac{\left(  \operatorname*{sgn}\left(  \zeta-z\right)  \right)
^{\mathbf{1}}}{\left(  \alpha-1\right)  !}\int_{\min\left\{  z,\zeta\right\}
}^{\max\left\{  z,\zeta\right\}  }\left(  \zeta-\sigma\right)  ^{\alpha
-1}D^{\alpha}\phi\left(  \sigma\right)  d\sigma.
\end{align*}

Also from part 7,%
\[
\left\vert \frac{\phi\left(  \xi+z\right)  }{\xi^{\alpha}}\right\vert
\leq\frac{1}{\left(  \alpha-1\right)  !}\frac{1}{\left\vert \xi^{\mathbf{1}%
}\right\vert }\int_{\min\left\{  \mathbf{0},\xi\right\}  }^{\max\left\{
\mathbf{0},\xi\right\}  }\left\vert D^{\alpha}\phi\left(  \cdot+z\right)
\right\vert ,\text{ }\forall\xi,
\]

so that: $\zeta=\xi+z\Rightarrow$%
\begin{align*}
\left\vert \frac{\phi\left(  \zeta\right)  }{\left(  \zeta-z\right)  ^{\alpha
}}\right\vert  & \leq\frac{1}{\left(  \alpha-1\right)  !}\frac{1}{\left\vert
\left(  \zeta-z\right)  ^{\mathbf{1}}\right\vert }\int_{\min\left\{
\mathbf{0},\zeta-z\right\}  }^{\max\left\{  \mathbf{0},\zeta-z\right\}
}\left\vert D^{\alpha}\phi\left(  \eta+z\right)  \right\vert d\eta\\
& :\sigma=\eta+z\Rightarrow\\
& =\frac{1}{\left(  \alpha-1\right)  !}\frac{1}{\left\vert \left(
\zeta-z\right)  ^{\mathbf{1}}\right\vert }\int_{\min\left\{  z,\zeta\right\}
}^{\max\left\{  z,\zeta\right\}  }\left\vert D^{\alpha}\phi\left(
\sigma\right)  \right\vert d\sigma.
\end{align*}

\end{proof}

\begin{remark}
What about replacing $D^{\alpha}$ by a difference?
\end{remark}

?? ADD BLURB!

\begin{lemma}
\label{Lem_Sxm_integ_expansion_2}Suppose $\psi\in S_{\otimes,\alpha}\left(
\pi\mathbb{Z}^{d}\right)  $ for some $\alpha\geq\mathbf{1}$. Then:

\begin{enumerate}
\item For all $\beta\in\mathbb{Z}^{d}$,%
\begin{align*}
\left\vert \psi\left(  \xi\right)  \right\vert  & \leq\frac{\left\vert \left(
\xi-\pi\beta\right)  ^{\alpha}\right\vert }{\alpha!}\max_{\pi\beta+\frac{\pi
}{2}I_{0}}\left\vert D^{\alpha}\psi\left(  \cdot\right)  \right\vert ,\quad
\xi\in\pi\beta+\pi I_{0},\\
\left\vert \psi\left(  \xi\right)  \right\vert  & \leq\frac{\left(
\pi/2\right)  ^{\left\vert \alpha\right\vert }}{\alpha!}\left\Vert D^{\alpha
}\psi\right\Vert _{\infty},\quad\xi\in\mathbb{R}^{d},
\end{align*}

where $I_{0}=R\left[  -\mathbf{1},\mathbf{1}\right]  $.

\item For all $\beta\in\mathbb{Z}^{d}$,%
\[
\left\vert \frac{\psi\left(  \xi\right)  }{\sin^{\alpha}\xi}\right\vert
\leq\frac{\left(  \pi/2\right)  ^{\left\vert \alpha\right\vert }}{\alpha!}%
\max_{\pi\beta+\frac{\pi}{2}I_{0}}\left\vert D^{\alpha}\psi\left(
\cdot\right)  \right\vert ,\quad\xi\in\pi\beta+\frac{\pi}{2}I_{0}.
\]

Further%
\[
\left\vert \frac{\psi\left(  \xi\right)  }{\sin^{\alpha}\xi}\right\vert
\leq\frac{\left(  \pi/2\right)  ^{\left\vert \alpha\right\vert }}{\alpha
!}\left\Vert D^{\alpha}\psi\right\Vert _{\infty},\quad\xi\in\mathbb{R}^{d}.
\]

\item If $f\in C_{BP}^{\infty}$ then%
\[
\left\vert \psi\left(  \xi\right)  \right\vert \leq\frac{\left(  \pi/2\right)
^{\left\vert \alpha\right\vert }}{\alpha!}\left\Vert D^{\alpha}\left(
f\psi\right)  \right\Vert _{\infty}\frac{\sin^{\alpha}\xi}{f\left(
\xi\right)  },\quad\xi\in\mathbb{R}^{d}.
\]

\item Suppose $f\in C_{BP}^{\infty}$ and $\left\vert \frac{f\left(  x\right)
}{f\left(  y\right)  }\right\vert \leq C_{f}$ whenever $\left\vert
x-y\right\vert \leq\pi$.

Then%
\[
\left\vert \psi\left(  \xi\right)  \right\vert \leq C_{f}\frac{\left(
\pi/2\right)  ^{\left\vert \alpha\right\vert }}{\alpha!}\left\Vert fD^{\alpha
}\psi\right\Vert _{\infty}\frac{\sin^{\alpha}\xi}{f\left(  \xi\right)  }%
,\quad\xi\in\mathbb{R}^{d}.
\]

The functions $f\left(  x\right)  =\left(  \left\vert a\right\vert
^{2}+\left\vert x\right\vert ^{2}\right)  ^{\mu}$ where $\mu>0$, satisfy the
specified condition when

\ $C_{f}=\left(  1+\max\left\{  1,\frac{2\pi^{2}}{\left\vert a\right\vert
^{2}}\right\}  \right)  ^{\mu}$.

\item If $\mathbf{0}\leq\nu\leq\mathbf{1}$ then%
\[
\sum_{\beta\in\mathbb{Z}^{d}}\left\vert \psi\left(  \pi\left(  \beta
+\nu\right)  \right)  \right\vert \leq\frac{\pi^{\left\vert \alpha
-1\right\vert }}{\left(  \alpha-1\right)  !}\int\left\vert D^{\alpha}%
\psi\right\vert ,
\]

and%
\[
\sum_{\beta\in\mathbb{Z}^{d}}\left\vert \psi\left(  \pi\left(  \beta
+\nu\right)  \right)  \right\vert ^{2}\leq\frac{\pi^{\left\vert 2\alpha
-1\right\vert }}{\left(  2\alpha-1\right)  ^{\mathbf{1}}\left(  \left(
\alpha-1\right)  !\right)  ^{2}}\int\left\vert D^{\alpha}\psi\right\vert ^{2}.
\]

Observe that the left hand sides of these two inequalities are actually
multi-periodic functions of $\nu$ with periods $\left\{  0,1\right\}
^{d}\setminus\mathbf{0}$.
\end{enumerate}
\end{lemma}

\begin{proof}
\textbf{Part 1} Use part 6 of Lemma \ref{Lem_Sxm_integ_expansion}.\medskip

\textbf{Part 2} From part 1,%
\[
\left\vert \psi\left(  \xi\right)  \right\vert \leq\frac{\left\vert \left(
\xi-\pi\beta\right)  ^{\alpha}\right\vert }{\alpha!}\max_{\pi\beta+\frac{\pi
}{2}I_{0}}\left\vert D^{\alpha}\psi\left(  \cdot\right)  \right\vert ,\quad
\xi\in\pi\beta+\frac{\pi}{2}I_{0}.
\]

Hence, when $\xi\in\pi\beta+\frac{\pi}{2}I_{0}$,%
\begin{align*}
\left\vert \frac{\psi\left(  \xi\right)  }{\sin^{\alpha}\xi}\right\vert
\leq\frac{1}{\alpha!}\left\vert \frac{\left(  \xi-\pi\beta\right)  ^{\alpha}%
}{\sin^{\alpha}\xi}\right\vert \max_{\pi\beta+\frac{\pi}{2}I_{0}}\left\vert
D^{\alpha}\psi\left(  \cdot\right)  \right\vert  &  \leq\frac{1}{\alpha!}%
\max_{\frac{\pi}{2}I_{0}}\left\vert \frac{\zeta^{\alpha}}{\sin^{\alpha}\zeta
}\right\vert \max_{\pi\beta+\frac{\pi}{2}I_{0}}\left\vert D^{\alpha}%
\psi\left(  \cdot\right)  \right\vert \\
&  =\frac{\left(  \pi/2\right)  ^{\left\vert \alpha\right\vert }}{\alpha!}%
\max_{\pi\beta+\frac{\pi}{2}I_{0}}\left\vert D^{\alpha}\psi\left(
\cdot\right)  \right\vert .
\end{align*}
\medskip

\textbf{Part 3} From part 6 of Lemma \ref{Lem_Sxl_and_Sx,n,l}, $f\psi\in
S_{\otimes,\alpha}\left(  \pi\mathbb{Z}^{d}\right)  $ and hence by part 8,%
\[
\left\vert \frac{\left(  f\psi\right)  \left(  \xi\right)  }{\sin^{\alpha}\xi
}\right\vert \leq\frac{\left(  \pi/2\right)  ^{\left\vert \alpha\right\vert }%
}{\alpha!}\left\Vert D^{\alpha}\left(  f\psi\right)  \right\Vert _{\infty
},\quad\xi\in\mathbb{R}^{d}.
\]

i.e.%
\[
\left\vert \psi\left(  \xi\right)  \right\vert \leq\frac{\left(  \pi/2\right)
^{\left\vert \alpha\right\vert }}{\alpha!}\left\Vert D^{\alpha}\left(
f\psi\right)  \right\Vert _{\infty}\frac{\sin^{\alpha}\xi}{f\left(
\xi\right)  },\quad\xi\in\mathbb{R}^{d}.
\]
\medskip

\textbf{Part 4} From part 2, for all $\beta\in\mathbb{Z}^{d}$,%
\[
\left\vert \frac{\psi\left(  \xi\right)  }{\sin^{\alpha}\xi}\right\vert
\leq\frac{\left(  \pi/2\right)  ^{\left\vert \alpha\right\vert }}{\alpha!}%
\max_{\eta\in\pi\beta+\frac{\pi}{2}I_{0}}\left\vert D^{\alpha}\psi\left(
\eta\right)  \right\vert ,\quad\xi\in\pi\beta+\frac{\pi}{2}I_{0}.
\]

Now%
\begin{align*}
\left\vert \frac{\psi\left(  \xi\right)  }{\sin^{\alpha}\xi}\right\vert  &
\leq\frac{\left(  \pi/2\right)  ^{\left\vert \alpha\right\vert }}{\alpha!}%
\max_{\eta\in\pi\beta+\frac{\pi}{2}I_{0}}\frac{\left\vert f\left(  \xi\right)
D^{\alpha}\psi\left(  \eta\right)  \right\vert }{f\left(  \xi\right)  }\\
& =\frac{\left(  \pi/2\right)  ^{\left\vert \alpha\right\vert }}{\alpha!}%
\max_{\eta\in\pi\beta+\frac{\pi}{2}I_{0}}\frac{\left\vert \frac{f\left(
\xi\right)  }{f\left(  \eta\right)  }f\left(  \eta\right)  D^{\alpha}%
\psi\left(  \eta\right)  \right\vert }{f\left(  \xi\right)  }\\
& \leq\frac{\left(  \pi/2\right)  ^{\left\vert \alpha\right\vert }}{\alpha
!}C_{f}\max_{\eta\in\pi\beta+\frac{\pi}{2}I_{0}}\frac{\left\vert f\left(
\eta\right)  D^{\alpha}\psi\left(  \eta\right)  \right\vert }{f\left(
\xi\right)  }\\
& \leq\frac{\left(  \pi/2\right)  ^{\left\vert \alpha\right\vert }}{\alpha
!}C_{f}\frac{\left\Vert fD^{\alpha}\psi\right\Vert _{\infty}}{f\left(
\xi\right)  }.
\end{align*}

Also%
\begin{align*}
\left\vert \frac{f\left(  x\right)  }{f\left(  y\right)  }\right\vert  &
=\left(  \frac{\left\vert a\right\vert ^{2}+\left\vert x\right\vert ^{2}%
}{\left\vert a\right\vert ^{2}+\left\vert y\right\vert ^{2}}\right)  ^{\mu
}\leq\left(  \frac{\left\vert a\right\vert ^{2}+2\left\vert x-y\right\vert
^{2}+2\left\vert y\right\vert ^{2}}{\left(  \left\vert a\right\vert
^{2}+\left\vert y\right\vert ^{2}\right)  ^{\mu}}\right)  ^{\mu}\leq\left(
\frac{\left\vert a\right\vert ^{2}+2\pi^{2}+2\left\vert y\right\vert ^{2}%
}{\left\vert a\right\vert ^{2}+\left\vert y\right\vert ^{2}}\right)  ^{\mu}=\\
& =\left(  1+\frac{2\pi^{2}+\left\vert y\right\vert ^{2}}{\left\vert
a\right\vert ^{2}+\left\vert y\right\vert ^{2}}\right)  ^{\mu}\leq\left(
1+\max\left\{  1,\frac{2\pi^{2}}{\left\vert a\right\vert ^{2}}\right\}
\right)  ^{\mu},
\end{align*}

and so we can choose $C_{f}=\left(  1+\max\left\{  1,\frac{2\pi^{2}%
}{\left\vert a\right\vert ^{2}}\right\}  \right)  ^{\mu}$.\medskip

\textbf{Part 5} From part 8 of Lemma \ref{Lem_Sxm_integ_expansion},%
\[
\left\vert \frac{\psi\left(  \xi\right)  }{\left(  \xi-z\right)  ^{\alpha}%
}\right\vert \leq\frac{1}{\left(  \alpha-1\right)  !}\frac{1}{\left\vert
\left(  \xi-z\right)  ^{\mathbf{1}}\right\vert }\int\limits_{\min\left\{
z,\xi\right\}  }^{\max\left\{  z,\xi\right\}  }\left\vert D^{\alpha}%
\psi\right\vert ,\quad\beta\in\mathbb{Z}^{d},
\]

and when $z=\pi\left(  \beta-\nu^{\prime}\right)  $, $\xi=\pi\left(  \beta
+\nu\right)  $, $\nu+\nu^{\prime}=\mathbf{1}$, $\nu,\nu^{\prime}\geq
\mathbf{0}$ this inequality becomes,%
\begin{align*}
\left\vert \frac{\psi\left(  \pi\left(  \beta+\nu\right)  \right)  }{\left(
\pi\mathbf{1}\right)  ^{\alpha}}\right\vert  & \leq\frac{1}{\left(
\alpha-1\right)  !}\frac{1}{\left\vert \left(  \pi\mathbf{1}\right)
^{\mathbf{1}}\right\vert }\int\limits_{\pi\left(  \beta-\nu^{\prime}\right)
}^{\pi\left(  \beta+\nu\right)  }\left\vert D^{\alpha}\psi\right\vert
,\quad\beta\in\mathbb{Z}^{d}.\\
\left\vert \psi\left(  \pi\left(  \beta+\nu\right)  \right)  \right\vert  &
\leq\frac{\pi^{\left\vert \alpha-1\right\vert }}{\left(  \alpha-1\right)
!}\int\limits_{\pi\left(  \beta-\nu^{\prime}\right)  }^{\pi\left(  \beta
+\nu\right)  }\left\vert D^{\alpha}\psi\right\vert ,\quad\beta\in
\mathbb{Z}^{d}.
\end{align*}

which implies%
\[
\sum_{\beta\in\mathbb{Z}^{d}}\left\vert \psi\left(  \pi\left(  \beta
+\nu\right)  \right)  \right\vert \leq\frac{\pi^{\left\vert \alpha
-1\right\vert }}{\left(  \alpha-1\right)  !}\int\left\vert D^{\alpha}%
\psi\right\vert .
\]

Again from part 8 of Lemma \ref{Lem_Sxm_integ_expansion},%
\[
\left\vert \frac{\psi\left(  \xi\right)  }{\left(  \xi-z\right)  ^{\alpha}%
}\right\vert \leq\frac{1}{\left(  2\alpha-1\right)  ^{\mathbf{1}/2}\left(
\alpha-1\right)  !}\left(  \frac{1}{\left\vert \left(  \xi-z\right)
^{\mathbf{1}}\right\vert }\int\limits_{\min\left\{  z,\xi\right\}  }%
^{\max\left\{  z,\xi\right\}  }\left\vert D^{\alpha}\psi\right\vert
^{2}\right)  ^{1/2}.
\]

When $z=\pi\left(  \beta-\nu^{\prime}\right)  $ and $\xi=\pi\left(  \beta
+\nu\right)  $, we get the sequence of inequalities%
\begin{align*}
\left\vert \frac{\psi\left(  \pi\left(  \beta+\nu\right)  \right)  }{\left(
\pi\mathbf{1}\right)  ^{\alpha}}\right\vert  & \leq\frac{1}{\left(
2\alpha-1\right)  ^{\mathbf{1}/2}\left(  \alpha-1\right)  !}\left(  \frac
{1}{\left(  \pi\mathbf{1}\right)  ^{\mathbf{1}}}\int\limits_{\pi\left(
\beta-\nu^{\prime}\right)  }^{\pi\left(  \beta+\nu\right)  }\left\vert
D^{\alpha}\psi\right\vert ^{2}\right)  ^{1/2}.\\
\left\vert \frac{\psi\left(  \pi\left(  \beta+\nu\right)  \right)  }{\left(
\pi\mathbf{1}\right)  ^{\alpha}}\right\vert ^{2}  & \leq\frac{1}{\left(
2\alpha-1\right)  ^{\mathbf{1}}\left(  \left(  \alpha-1\right)  !\right)
^{2}}\frac{1}{\left(  \pi\mathbf{1}\right)  ^{\mathbf{1}}}\int\limits_{\pi
\left(  \beta-\nu^{\prime}\right)  }^{\pi\left(  \beta+\nu\right)  }\left\vert
D^{\alpha}\psi\right\vert ^{2}.\\
\left\vert \psi\left(  \pi\left(  \beta+\nu\right)  \right)  \right\vert ^{2}
& \leq\frac{\pi^{\left\vert 2\alpha-1\right\vert }}{\left(  2\alpha-1\right)
^{\mathbf{1}}\left(  \left(  \alpha-1\right)  !\right)  ^{2}}\int%
\limits_{\pi\left(  \beta-\nu^{\prime}\right)  }^{\pi\left(  \beta+\nu\right)
}\left\vert D^{\alpha}\psi\right\vert ^{2}.\\
\left\vert \psi\left(  \pi\left(  \beta+\nu\right)  \right)  \right\vert ^{2}
& \leq\frac{\pi^{\left\vert 2\alpha-1\right\vert }}{\left(  2\alpha-1\right)
^{\mathbf{1}}\left(  \left(  \alpha-1\right)  !\right)  ^{2}}\int%
\limits_{\pi\left(  \beta-\nu^{\prime}\right)  }^{\pi\left(  \beta+\nu\right)
}\left\vert D^{\alpha}\psi\right\vert ^{2}.
\end{align*}

Hence%
\[
\sum_{\beta\in\mathbb{Z}^{d}}\left\vert \psi\left(  \pi\left(  \beta
+\nu\right)  \right)  \right\vert ^{2}\leq\frac{\pi^{\left\vert 2\alpha
-1\right\vert }}{\left(  2\alpha-1\right)  ^{\mathbf{1}}\left(  \left(
\alpha-1\right)  !\right)  ^{2}}\int\left\vert D^{\alpha}\psi\right\vert ^{2}.
\]

\end{proof}

We can now characterize the spaces $S_{\otimes,2\alpha}\left(  \mathbf{0}%
\right)  $.

\begin{lemma}
\label{Lem_Sx2m_expansion}Suppose $\phi\in S\left(  \mathbb{R}^{d}\right)  $
and $\alpha\geq0$. Then $\phi\in S_{\otimes,2\alpha}\left(  \mathbf{0}\right)
$ iff there exist functions $t_{k}:\mathbb{R}^{d}\rightarrow\left(
0,1\right)  $ such that%
\begin{equation}
\phi\left(  \xi\right)  =\frac{\xi^{2\alpha}}{\left(  2\alpha\right)
!}\left(  D^{2\alpha}\phi\right)  \left(  t_{1}\left(  \xi\right)  \xi
_{1},\ldots,t_{d}\left(  \xi\right)  \xi_{d}\right)  .\label{a1.29}%
\end{equation}

\end{lemma}

\begin{proof}
If $\phi\in S_{\otimes,2\alpha}\left(  \mathbf{0}\right)  $ then part 3 of
Lemma \ref{Lem_Sxm_expansion} implies \ref{a1.29}.

On the other hand, suppose \ref{a1.29} holds for some $\phi\in S$. Then
$\int_{0}^{1}\frac{\phi\left(  b^{\prime},x_{k},b^{\prime\prime}\right)
}{x_{k}^{2\alpha_{k}}}dx_{k}<\infty$ for all $b=\left(  b^{\prime}%
,b^{\prime\prime}\right)  \in\mathbb{R}^{d-1}$ and thus Lemma
\ref{Lem_deriv_origin} implies that $\left(  D_{k}^{n}\phi\right)  \left(
b^{\prime},0,b^{\prime\prime}\right)  =0$ for all $b\in\mathbb{R}^{d-1}$,
$\forall k$ and $n<2\alpha_{k}$ i.e. $\phi\in S_{\otimes,2\alpha}\left(
\mathbf{0}\right)  $.
\end{proof}

\begin{theorem}
\label{Thm_integ_condits_for_func_in_Sw2_splin_wt_fn}Suppose $\lambda\geq
\nu\geq\mathbf{1}$ are multi-indexes.

Then $\int\frac{\eta^{2\nu}}{\sin^{2\lambda}\eta}\left\vert \phi\left(
\eta\right)  \right\vert ^{2}d\eta<\infty$ iff $\int\nolimits_{-\frac{\pi}%
{2}\mathbf{1}}^{\frac{\pi}{2}\mathbf{1}}\frac{\left\vert \phi\left(
\eta\right)  \right\vert ^{2}}{\eta^{2\left(  \lambda-\nu\right)  }}%
d\eta<\infty$ and $\int\nolimits_{-\frac{\pi}{2}\mathbf{1}}^{\frac{\pi}%
{2}\mathbf{1}}\frac{\left\vert \left(  x^{\nu}\phi\right)  \left(  \pi
\alpha+\eta\right)  \right\vert ^{2}}{\eta^{2\lambda}}d\eta<\infty$ for
$\alpha\in\mathbb{Z}^{d}\setminus\mathbf{0}$.
\end{theorem}

\begin{proof}
We begin by changing the domain of integration to the closed rectangle
$R\left[  -\frac{\pi}{2}\mathbf{1,}\frac{\pi}{2}\mathbf{1}\right]  $. Then,
using the notation $\alpha_{+}:=\left(  \left\vert \alpha_{i}\right\vert
\right)  $ and $\phi_{\nu}\left(  \eta\right)  :=\eta^{\nu}\phi\left(
\eta\right)  $,%
\begin{align}
&  \int\frac{\eta^{2\nu}}{\sin^{2\lambda}\eta}\left\vert \phi\left(
\eta\right)  \right\vert ^{2}d\eta\nonumber\\
&  =\sum\limits_{\alpha\in\mathbb{Z}^{d}}\int_{\pi\alpha-\pi/2}^{\pi\alpha
+\pi/2}\frac{\eta^{2\nu}}{\sin^{2\lambda}\eta}\left\vert \phi\left(
\eta\right)  \right\vert ^{2}d\eta\nonumber\\
&  =\sum\limits_{\alpha\in\mathbb{Z}^{d}}\int_{-\frac{\pi}{2}\mathbf{1}%
}^{\frac{\pi}{2}\mathbf{1}}\frac{\left\vert \phi_{\nu}\left(  \eta-\pi
\alpha\right)  \right\vert ^{2}}{\sin^{2\lambda}\eta}d\eta\nonumber\\
&  =\sum\limits_{\alpha\in\mathbb{Z}^{d}}\int_{-\frac{\pi}{2}\mathbf{1}%
}^{\frac{\pi}{2}\mathbf{1}}\frac{\left\vert \phi_{\nu}\left(  \pi\alpha
+\eta\right)  \right\vert ^{2}}{\sin^{2\lambda}\eta}d\eta\nonumber\\
&  =\int\limits_{-\frac{\pi}{2}\mathbf{1}}^{\frac{\pi}{2}\mathbf{1}}\frac
{\eta^{2\lambda}}{\sin^{2\lambda}\eta}\frac{\left\vert \phi\left(
\eta\right)  \right\vert ^{2}}{\eta^{2\left(  \lambda-\nu\right)  }}d\eta
+\sum\limits_{\alpha\in\mathbb{Z}^{d}\setminus\mathbf{0}}\int\limits_{-\frac
{\pi}{2}\mathbf{1}}^{\frac{\pi}{2}\mathbf{1}}\frac{\eta^{2\lambda}}%
{\sin^{2\lambda}\eta}\frac{\left\vert \phi_{\nu}\left(  \pi\alpha+\eta\right)
\right\vert ^{2}}{\eta^{2\lambda}}d\eta,\label{a1.008}%
\end{align}

which means $\int\frac{\eta^{2\nu}}{\sin^{2\lambda}\eta}\left\vert \phi\left(
\eta\right)  \right\vert ^{2}d\eta<\infty$ implies%
\begin{equation}
\int\limits_{-\frac{\pi}{2}\mathbf{1}}^{\frac{\pi}{2}\mathbf{1}}%
\frac{\left\vert \phi\left(  \eta\right)  \right\vert ^{2}}{\eta^{2\left(
\lambda-\nu\right)  }}d\eta<\infty,\quad\int\limits_{-\frac{\pi}{2}\mathbf{1}%
}^{\frac{\pi}{2}\mathbf{1}}\frac{\left\vert \phi_{\nu}\left(  \pi\alpha
+\eta\right)  \right\vert ^{2}}{\eta^{2\lambda}}d\eta<\infty,\quad\alpha
\in\mathbb{Z}^{d}\setminus\mathbf{0}.\label{a1.013}%
\end{equation}

Conversely, suppose \ref{a1.013} is true. Then Lemma \ref{Lem_2} implies
$\phi\in S_{\otimes,\lambda-\nu}\left(  \mathbf{0}\right)  $ and $\phi_{\nu
}\left(  \pi\alpha+\cdot\right)  \in S_{\otimes,\lambda}\left(  \mathbf{0}%
\right)  $ for $\alpha\in\mathbb{Z}^{d}\setminus\mathbf{0}$. Further, by part
3 of Lemma \ref{Lem_Sxm_expansion}, we have the easy estimate
\[
\left\vert \frac{\phi_{\nu}\left(  \pi\alpha+\eta\right)  }{\eta^{\lambda}%
}\right\vert \leq\frac{1}{\left(  \lambda!\right)  ^{d}}\max_{\xi\in
I_{\alpha}}\left\vert D^{\lambda}\phi_{\nu}\left(  \xi\right)  \right\vert
,\quad\left\{
\begin{array}
[c]{l}%
\alpha\in\mathbb{Z}^{d}\setminus\mathbf{0},\\
\eta\in I_{\alpha},
\end{array}
\right.
\]

where $I_{\alpha}$ is the closed rectangle%
\[
I_{\alpha}=\pi\alpha+R\left[  -\frac{\pi}{2}\mathbf{1,}\frac{\pi}{2}%
\mathbf{1}\right]  =\pi\alpha+I_{\mathbf{0}}.
\]

Thus%
\begin{align*}
\int\limits_{-\frac{\pi}{2}\mathbf{1}}^{\frac{\pi}{2}\mathbf{1}}%
\frac{\left\vert \phi_{\nu}\left(  \pi\alpha+\eta\right)  \right\vert ^{2}%
}{\eta^{2\lambda}}d\eta & \leq\frac{vol\left(  \overline{I}_{\alpha}\right)
}{\left(  \lambda!\right)  ^{2d}}\max_{x\in I_{\alpha}}\left\vert D^{\lambda
}\phi_{\nu}\left(  x\right)  \right\vert ^{2}\\
& =\frac{\pi^{d}}{\left(  \lambda!\right)  ^{2d}}\max_{x\in I_{\alpha}%
}\left\vert D^{\lambda}\phi_{\nu}\left(  x\right)  \right\vert ^{2},
\end{align*}

and consequently%
\begin{align*}
\sum\limits_{\alpha\in\mathbb{Z}^{d}\setminus\mathbf{0}}\int_{-\frac{\pi}%
{2}\mathbf{1}}^{\frac{\pi}{2}\mathbf{1}}\frac{\eta^{2\lambda}}{\sin^{2\lambda
}\eta}\frac{\left\vert \phi_{\nu}\left(  \pi\alpha+\eta\right)  \right\vert
^{2}}{\eta^{2\lambda}}d\eta & \leq\sum\limits_{\alpha\in\mathbb{Z}%
^{d}\setminus\mathbf{0}}\frac{\pi^{d}}{\left(  \lambda!\right)  ^{2d}}%
\max_{x\in I_{\alpha}}\left\vert D^{\lambda}\phi_{\nu}\left(  x\right)
\right\vert ^{2}\\
& \leq\frac{\pi^{d}}{\left(  \lambda!\right)  ^{2d}}\left\Vert D^{\lambda}%
\phi_{\nu}\right\Vert _{\infty}^{2}\\
& <\infty.
\end{align*}

\end{proof}

\begin{theorem}
\label{Thm_Sw2_extsplin_wt}Suppose $w_{\nu,\lambda}\left(  \eta\right)
=\frac{\eta^{2\nu}}{\sin^{2\lambda}\eta}$ where $\lambda\geq\nu\geq\mathbf{1}%
$. Then:

\begin{enumerate}
\item
\begin{align}
S_{w_{\nu,\lambda;2}}  & =\left\{  \phi\in S_{\otimes,\lambda-\nu}\left(
\mathbf{0}\right)  :x^{\nu}\phi\in S_{\otimes,\lambda}\left(  \pi
\mathbb{Z}^{d}\setminus\mathbf{0}\right)  \right\} \label{a1.17}\\
& =S_{\otimes,\lambda-\nu}\left(  \mathbf{0}\right)  \cap x^{-\nu}%
S_{\otimes,\lambda}\left(  \pi\mathbb{Z}^{d}\setminus\mathbf{0}\right)
.\nonumber
\end{align}

\item
\begin{align*}
S_{w_{\nu,\lambda;2}}  & =\left\{  \phi\in S:x^{\nu}\phi\in S_{\otimes
,\lambda}\left(  \pi\mathbb{Z}^{d}\right)  \right\} \\
& =S\cap x^{-\nu}S_{\otimes,\lambda}\left(  \pi\mathbb{Z}^{d}\right)  .
\end{align*}

\item $D^{\beta}S_{w_{\nu,\lambda;2}}\subset S_{w_{\nu,\lambda-\beta;2}}$ when
$\beta\leq\lambda-\nu$.

\item $S_{w_{\nu,\lambda;2}}\subset S_{\otimes,\lambda}\left(  \pi
+2\pi\mathbb{Z}^{d}\right)  $. The set $\pi+2\pi\mathbb{Z}^{d}$ avoids any
hyperplane parallel to the axes which passes through the origin i.e. no point
has a zero component.
\end{enumerate}
\end{theorem}

\begin{proof}
\textbf{Part 1} We apply Lemma \ref{Lem_2} to the result of Theorem
\ref{Thm_integ_condits_for_func_in_Sw2_splin_wt_fn} and conclude that $\phi\in
S_{w_{\nu,\lambda};2}$ iff $\phi\in S_{\otimes,\lambda-\nu}\left(
\mathbf{0}\right)  $ and $\left(  x^{\nu}\phi\right)  \left(  \pi\alpha
+\cdot\right)  \in S_{\otimes,\lambda}\left(  \mathbf{0}\right)  $ for
$\alpha\in\mathbb{Z}^{d}\setminus\mathbf{0}$ i.e.
\begin{align*}
S_{w_{\nu,\lambda};2}  & =\left\{  \phi\in S_{\otimes,\lambda-\nu}\left(
\mathbf{0}\right)  :x^{\nu}\phi\in S_{\otimes,\lambda}\left(  \mathbb{Z}%
^{d}\setminus\mathbf{0}\right)  \right\} \\
& =S_{\otimes,\lambda-\nu}\left(  \mathbf{0}\right)  \cap x^{-\nu}%
S_{\otimes,\lambda}\left(  \mathbb{Z}^{d}\setminus\mathbf{0}\right)  .
\end{align*}
\medskip

\textbf{Part 2} Noting part 2 of Lemma \ref{Lem_Sxl_and_Sx,n,l}, from part 1,
\begin{align*}
\phi\in S_{w_{\nu,\lambda};2} &  \Rightarrow\phi\in S_{\otimes,\lambda-\nu
}\left(  \mathbf{0}\right)  \text{ }and\text{ }x^{\nu}\phi\in S_{\otimes
,\lambda}\left(  \pi\mathbb{Z}^{d}\setminus\mathbf{0}\right) \\
&  \Rightarrow x^{\nu}\phi\in S_{\otimes,\lambda}\left(  \mathbf{0}\right)
\text{ }and\text{ }x^{\nu}\phi\in S_{\otimes,\lambda}\left(  \pi\mathbb{Z}%
^{d}\setminus\mathbf{0}\right) \\
&  \Rightarrow x^{\nu}\phi\in S_{\otimes,\lambda}\left(  \pi\mathbb{Z}%
^{d}\right)  ,
\end{align*}

and by part 3 of Lemma \ref{Lem_Sxl_and_Sx,n,l},%
\begin{align*}
\phi\in S\text{ }and\text{ }x^{\nu}\phi\in S_{\otimes,\lambda}\left(
\pi\mathbb{Z}^{d}\right)   &  \Rightarrow x^{\nu}\phi\in S_{\otimes,\lambda
}\left(  \mathbf{0}\right)  \text{ }and\text{ }x^{\nu}\phi\in S_{\otimes
,\lambda}\left(  \pi\mathbb{Z}^{d}\setminus\mathbf{0}\right) \\
&  \Rightarrow\phi\in S_{\otimes,\lambda-\nu}\left(  \mathbf{0}\right)  \text{
}and\text{ }x^{\nu}\phi\in S_{\otimes,\lambda}\left(  \pi\mathbb{Z}%
^{d}\setminus\mathbf{0}\right) \\
&  \Rightarrow\phi\in S_{w_{\nu,\lambda};2}.
\end{align*}
\medskip

\textbf{Part 3} We use the characterization \ref{a1.17} of $S_{w_{\nu,\lambda
};2}$. Thus $\phi\in S_{w_{\nu,\lambda};2}$ implies $\phi\in S_{\otimes
,\lambda-\nu}\left(  \mathbf{0}\right)  $ and from part 1 of Lemma
\ref{Lem_Sxl_and_Sx,n,l}, $D^{\beta}\phi\in S_{\otimes,\lambda-\nu-\beta
}\left(  \mathbf{0}\right)  $ when $\beta\leq\lambda-\nu$.

Noting \ref{a1.17} we must show $x^{\nu}D^{\beta}\phi\in S_{\otimes
,\lambda-\beta}\left(  \pi\mathbb{Z}^{d}\setminus\mathbf{0}\right)  $. But%
\[
\phi\in S_{w_{\nu,\lambda};2}%
\begin{array}
[t]{ll}%
\Rightarrow x^{\nu}\phi\in S_{\otimes,\lambda}\left(  \pi\mathbb{Z}%
^{d}\setminus\mathbf{0}\right)  =%
{\textstyle\bigcup\limits_{\gamma\neq\mathbf{0}}}
S_{\otimes,\lambda}\left(  \gamma\right)  , & part\text{ }1\text{
}Theorem\text{ }\ref{Thm_Sw2_extsplin_wt},\\
\Rightarrow\phi\in%
{\textstyle\bigcup\limits_{\gamma\neq\mathbf{0}}}
S_{\otimes,\lambda-\left(  \gamma\And\mathbf{0}\right)  .\nu}\left(
\gamma\right)  , & part\text{ }5\text{ }Lemma\text{ }\ref{Lem_Sxl_and_Sx,n,l}%
,\\
\Rightarrow D^{\beta}\phi\in%
{\textstyle\bigcup\limits_{\gamma\neq\mathbf{0}}}
S_{\otimes,\lambda-\beta-\left(  \gamma\And\mathbf{0}\right)  .\nu}\left(
\gamma\right)  , & part\text{ }1\text{ }Lemma\text{ }\ref{Lem_Sxl_and_Sx,n,l}%
,\\
\Rightarrow x^{\nu}D^{\beta}\phi\in%
{\textstyle\bigcup\limits_{\gamma\neq\mathbf{0}}}
x^{\nu}S_{\otimes,\lambda-\beta-\left(  \gamma\And\mathbf{0}\right)  .\nu
}\left(  \gamma\right)  , & \\
\Rightarrow x^{\nu}D^{\beta}\phi\in%
{\textstyle\bigcup\limits_{\gamma\neq\mathbf{0}}}
S_{\otimes,\lambda-\beta}\left(  \gamma\right)  , & part\text{ }5\text{
}Lemma\text{ }\ref{Lem_Sxl_and_Sx,n,l},
\end{array}
\]

so that $x^{\nu}D^{\beta}\phi\in S_{\otimes,\lambda-\beta}\left(
\pi\mathbb{Z}^{d}\setminus\mathbf{0}\right)  $ as claimed.

\textbf{Part 4} ??
\end{proof}

\begin{remark}
\label{Rem_Lem_bnd_Sx,l_cube}What about when $\lambda\ngeq\mathbf{1}$?
\end{remark}

?? We now show that $w\in W_{S;0}$ - see Definition
\ref{Def_Sw2_and_fnal(Sw2)}.

\begin{theorem}
\label{Thm_eta^2/sin(eta)^2.|phi|^2_dim_gt_1}Suppose $\phi\in S\left(
\mathbb{R}^{d}\right)  $, $\nu$ and $\lambda$ are multi-indexes such that
$\lambda\geq\nu\geq1$ and

$\int\frac{\eta^{2\nu}}{\sin^{2\lambda}\eta}\left\vert \phi\left(
\eta\right)  \right\vert ^{2}d\eta<\infty$ i.e. $\phi\in S_{w,0}$ - see
Definition \ref{Def_Sw2_and_fnal(Sw2)}. Note part 3 of
\ref{Thm_Sw2_extsplin_wt}.

Then:

\begin{enumerate}
\item
\begin{align*}
\int &  \frac{\eta^{2\nu}}{\sin^{2\lambda}\eta}\left\vert \phi\left(
\eta\right)  \right\vert ^{2}d\eta\\
&  \leq\left(  \frac{\pi}{2}\right)  ^{2\left\vert \lambda\right\vert }\left(
\left(  \frac{2^{\left\vert \lambda-\nu\right\vert }}{\left(  2\left(
\lambda-\nu\right)  -1\right)  !!}\right)  ^{2}\left\Vert D^{\lambda-\nu}%
\phi\right\Vert _{2}^{2}+\left(  \frac{2^{\left\vert \lambda\right\vert }%
}{\left(  2\lambda-1\right)  !!}\right)  ^{2}\left\Vert D^{\lambda}\left(
x^{\nu}\phi\right)  \right\Vert _{2}^{2}\right)  ,
\end{align*}

where we define $\left(  -1\right)  !!=1$.\medskip

For parts 2 and 3 we need the following \textbf{definition}: given $\xi
\in\mathbb{R}^{d}$ denote by $\mathcal{R}_{\pi}\left(  \xi\right)  $ the
unique \textbf{half-open rectangle} $\left\{  \eta:\pi\beta\leq\eta
\centerdot<\pi\left(  \beta+1\right)  \right\}  $ which contains $\xi$ for
some $\beta\in\mathbb{Z}^{d}$.

\item ??

\begin{enumerate}
\item Derive $L^{\infty}$, $L^{1}$, $L^{2}$ upper bounds for $\phi$. See
attempts in earlier documents.

\item From part 2 of Theorem \ref{Thm_Sw2_extsplin_wt},\allowbreak
\ $S_{w_{\nu,\lambda};2}=\left\{  \phi\in S:\xi^{\nu}\phi\in S_{\otimes
,\lambda}\left(  \pi\mathbb{Z}^{d}\right)  \right\}  $ so the results of Lemma
\ref{Lem_Sxm_integ_expansion_2} apply for $\psi=\xi^{\nu}\phi$ and
$\alpha=\lambda$.
\end{enumerate}

\item $\frac{\xi^{\nu}}{\sin^{\lambda}\xi}\phi\left(  \xi\right)  \in S$ and
for all $\xi\in\mathbb{R}^{d}$,%
\[
\left\vert \frac{\xi^{\nu}}{\sin^{\lambda}\xi}\phi\left(  \xi\right)
\right\vert \leq\frac{\left(  \pi/2\right)  ^{\left\vert \lambda\right\vert }%
}{\lambda!}\left\Vert D^{\lambda}\left(  x^{\nu}\phi\right)  \right\Vert
_{\infty}.
\]

\item (See better result in Corollary
\ref{Cor_Thm_eta^2/sin(eta)^2.|phi|^2_dim_gt_2}) For all $n=0,1,2,\ldots$ and
$\beta\geq\mathbf{0}$:%
\begin{multline*}
\left\vert \xi\right\vert ^{n}\left\vert D^{\beta}\left(  \frac{\xi^{\nu}%
}{\sin^{\lambda}\xi}\phi\left(  \xi\right)  \right)  \right\vert \\
\leq\frac{\left(  1+2\sqrt{d}\right)  ^{n}}{\lambda!}\left(  \sum_{\gamma
\leq\beta}\tbinom{\beta}{\gamma}\left\Vert D^{\gamma}\left(  \frac
{\eta^{\lambda}}{\sin^{\lambda}\eta}\right)  \right\Vert _{\infty,R\left[
-\frac{\pi}{2}\mathbf{1,}\frac{\pi}{2}\mathbf{1}\right]  }\right)  \times\\
\times\max_{\gamma\leq\beta}\left\Vert \left(  1+\left\vert \xi\right\vert
^{n}\right)  D^{\gamma+\lambda}\left(  \xi^{\nu}\phi\right)  \right\Vert
_{\infty},
\end{multline*}

where%
\begin{align*}
\sum_{\gamma\leq\beta}\tbinom{\beta}{\gamma} &  \left\Vert D^{\gamma}\left(
\frac{\eta^{\lambda}}{\sin^{\lambda}\eta}\right)  \right\Vert _{\infty
,R\left[  -\frac{\pi}{2}\mathbf{1,}\frac{\pi}{2}\mathbf{1}\right]  }%
=a_{\beta_{1}}a_{\beta_{2}}\ldots a_{\beta_{d}},\\
a_{k} &  :=\sum_{j=0}^{k}\tbinom{k}{j}\max_{\left\vert s\right\vert \leq
\frac{\pi}{2}}\left\vert D^{j}\text{ }\left(  \frac{\sin s}{s}\right)
^{-\lambda_{k}}\right\vert .
\end{align*}

\item For all $\xi\in\mathbb{R}^{d}$,%
\begin{equation}
\left\vert \phi\left(  \xi\right)  \right\vert \leq\frac{\left(  \pi/2\right)
^{\left\vert \lambda\right\vert }}{\lambda!}\left\Vert D^{\lambda}\left(
x^{\nu}\phi\right)  \right\Vert _{\infty}\left\vert \frac{\sin^{\lambda}\xi
}{\xi^{\nu}}\right\vert .\label{a3.11}%
\end{equation}

\item
\[
\int\left\vert \phi\right\vert \leq\left(  \int\frac{\sin^{2\lambda}\eta}%
{\eta^{2\nu}}\right)  ^{1/2}\left(  \int\frac{\eta^{2\nu}}{\sin^{2\lambda}%
\eta}\left\vert \phi\left(  \eta\right)  \right\vert ^{2}d\eta\right)  ^{1/2}.
\]

\item
\[
\left(  \int\left\vert \phi\right\vert ^{2}\right)  ^{1/2}\leq\left\Vert
\frac{\sin^{\lambda}\eta}{\eta^{\nu}}\right\Vert _{\infty}\left(  \int%
\frac{\eta^{2\nu}}{\sin^{2\lambda}\eta}\left\vert \phi\left(  \eta\right)
\right\vert ^{2}d\eta\right)  ^{1/2}.
\]

\end{enumerate}
\end{theorem}

\begin{proof}
\fbox{\textbf{Part 1}} From the proof of Theorem \ref{Thm_Sw2_extsplin_wt},
$\int\frac{\eta^{2\nu}}{\sin^{2\lambda}\eta}\left\vert \phi\left(
\eta\right)  \right\vert ^{2}d\eta<\infty$ implies \ref{a1.008} i.e.%
\begin{equation}
\int\frac{\eta^{2\nu}}{\sin^{2\lambda}\eta}\left\vert \phi\left(  \eta\right)
\right\vert ^{2}d\eta=\int\limits_{-\frac{\pi}{2}\mathbf{1}}^{\frac{\pi}%
{2}\mathbf{1}}\frac{\eta^{2\nu}}{\sin^{2\lambda}\eta}\left\vert \phi\left(
\eta\right)  \right\vert ^{2}d\eta+\sum\limits_{\alpha\in\mathbb{Z}%
^{d}\setminus0}\int\limits_{-\frac{\pi}{2}\mathbf{1}}^{\frac{\pi}{2}%
\mathbf{1}}\frac{\eta^{2\lambda}}{\sin^{2\lambda}\eta}\frac{\left\vert
\phi_{\nu}\left(  \pi\alpha+\eta\right)  \right\vert ^{2}}{\eta^{2\lambda}%
}d\eta,\label{a1.70}%
\end{equation}

where
\begin{equation}
\phi_{\nu}\left(  \eta\right)  :=\eta^{\nu}\phi\left(  \eta\right)
.\label{a1.012}%
\end{equation}

From part 1 of Theorem \ref{Thm_Sw2_extsplin_wt} we have $\phi\in
S_{\otimes,\lambda-\nu}\left(  \mathbf{0}\right)  $ and $\phi_{\nu}\in
S_{\otimes,\lambda}\left(  \pi\mathbb{Z}^{d}\right)  $. Since $\lambda\geq1$
by \ref{1.076},%
\begin{equation}
\phi_{\nu}\left(  \pi\alpha+\eta\right)  =\frac{\eta^{\lambda}}{\left(
\lambda-1\right)  !}\int_{\mathbf{0}}^{\mathbf{1}}\left(  1-\tau\right)
^{\lambda-1}\left(  D^{\lambda}\phi_{\nu}\right)  \left(  \pi\alpha+\tau
.\eta\right)  d\tau,\quad\alpha\in\pi\mathbb{Z}^{d},\label{a1.13}%
\end{equation}

and%
\begin{equation}
\phi\left(  \eta\right)  =\frac{\eta^{\lambda-\nu}}{\left(  \lambda
-\nu-1\right)  !}\int_{\mathbf{0}}^{\mathbf{1}}\left(  1-\tau\right)
^{\lambda-\nu-1}\left(  D^{\lambda-\nu}\phi\right)  \left(  \tau.\eta\right)
d\tau.\label{a1.131}%
\end{equation}
\medskip

\underline{\textbf{The summation term of} \ref{a1.70}} The next step is to use
equation \ref{a1.13} to estimate the summation term on the right hand side of
\ref{a1.70}. In fact, if $I_{\alpha}$ denotes the closed rectangle%
\[
I_{\alpha}:=R\left[  \pi\alpha-\frac{\pi}{2}\mathbf{1,}\pi\alpha+\frac{\pi}%
{2}\mathbf{1}\right]  =\pi\alpha+I_{\mathbf{0}},\quad\alpha\in\mathbb{Z}^{d},
\]

substituting the estimate \ref{a1.13} into a summation term of \ref{a1.70}
implies%
\begin{align*}
\int_{-\frac{\pi}{2}\mathbf{1}}^{\frac{\pi}{2}\mathbf{1}} &  \frac
{\eta^{2\lambda}}{\sin^{2\lambda}\eta}\frac{\left\vert \phi_{\nu}\left(
\pi\alpha+\eta\right)  \right\vert ^{2}}{\eta^{2\lambda}}d\eta\\
&  =\frac{1}{\left(  \left(  \lambda-1\right)  !\right)  ^{2}}\int_{-\frac
{\pi}{2}\mathbf{1}}^{\frac{\pi}{2}\mathbf{1}}\frac{\eta^{2\lambda}}%
{\sin^{2\lambda}\eta}\left\vert \int_{\mathbf{0}}^{\mathbf{1}}\left(
1-\tau\right)  ^{\lambda-1}\left(  D^{\lambda}\phi_{\nu}\right)  \left(
\pi\alpha+\tau.\eta\right)  d\tau\right\vert ^{2}d\eta\\
&  \leq\frac{1}{\left(  \left(  \lambda-1\right)  !\right)  ^{2}}\left(
\frac{\pi}{2}\right)  ^{2\left\vert \lambda\right\vert }\int_{-\frac{\pi}%
{2}\mathbf{1}}^{\frac{\pi}{2}\mathbf{1}}\left\vert \int_{\mathbf{0}%
}^{\mathbf{1}}\left(  1-\tau\right)  ^{\lambda-1}\left(  D^{\lambda}\phi_{\nu
}\right)  \left(  \pi\alpha+\tau.\eta\right)  d\tau\right\vert ^{2}d\eta\\
&  \leq\left(  \frac{\left(  \pi/2\right)  ^{\left\vert \lambda\right\vert }%
}{\left(  \lambda-1\right)  !}\right)  ^{2}\int_{-\frac{\pi}{2}\mathbf{1}%
}^{\frac{\pi}{2}\mathbf{1}}\left(  \int_{\mathbf{0}}^{\mathbf{1}}\left(
1-\tau\right)  ^{\lambda-1}\left\vert \left(  D^{\lambda}\phi_{\nu}\right)
\left(  \pi\alpha+\tau.\eta\right)  \right\vert d\tau\right)  ^{2}d\eta\\
&  =\left(  \frac{\left(  \pi/2\right)  ^{\left\vert \lambda\right\vert }%
}{\left(  \lambda-1\right)  !}\right)  ^{2}\int_{-\frac{\pi}{2}\mathbf{1}%
}^{\frac{\pi}{2}\mathbf{1}}\left(
\begin{array}
[c]{l}%
\int_{\mathbf{0}}^{\mathbf{1}}\left(  1-\sigma\right)  ^{\lambda-1}\left\vert
\left(  D^{\lambda}\phi_{\nu}\right)  \left(  \pi\alpha+\sigma.\eta\right)
\right\vert d\sigma\times\\
\quad\times\int_{\mathbf{0}}^{\mathbf{1}}\left(  1-\tau\right)  ^{\lambda
-1}\left\vert \left(  D^{\lambda}\overline{\phi_{\nu}}\right)  \left(
\pi\alpha+\tau.\eta\right)  \right\vert d\tau
\end{array}
\right)  d\eta\\
&  =\left(  \frac{\left(  \pi/2\right)  ^{\left\vert \lambda\right\vert }%
}{\left(  \lambda-1\right)  !}\right)  ^{2}\int_{\mathbf{0}}^{\mathbf{1}}%
\int_{\mathbf{0}}^{\mathbf{1}}\left(  1-\sigma\right)  ^{\lambda-1}\left(
1-\tau\right)  ^{\lambda-1}f_{\alpha}^{\left(  \lambda\right)  }\left(
\phi_{\nu};\sigma,\tau\right)  d\sigma d\tau,
\end{align*}

where%
\begin{equation}
f_{\alpha}^{\left(  \lambda\right)  }\left(  \phi_{\nu};\sigma,\tau\right)
:=\int\limits_{-\frac{\pi}{2}\mathbf{1}}^{\frac{\pi}{2}\mathbf{1}}\left\vert
\left(  D^{\lambda}\phi_{\nu}\right)  \left(  \pi\alpha+\sigma.\eta\right)
\right\vert \left\vert \left(  D^{\lambda}\overline{\phi_{\nu}}\right)
\left(  \pi\alpha+\tau.\eta\right)  \right\vert d\eta,\quad\left\{
\begin{array}
[c]{l}%
\lambda\geq\mathbf{1},\\
\alpha\in\mathbb{Z}^{d}\setminus\mathbf{0.}%
\end{array}
\right. \label{1.078}%
\end{equation}

so that%
\begin{align}
\sum\limits_{\alpha\in\mathbb{Z}^{d}\setminus0} &  \int_{-\frac{\pi}%
{2}\mathbf{1}}^{\frac{\pi}{2}\mathbf{1}}\frac{\left\vert \phi_{\nu}\left(
\pi\alpha+\eta\right)  \right\vert ^{2}}{\sin^{2\lambda}\eta}d\eta\nonumber\\
&  \leq\left(  \frac{\left(  \pi/2\right)  ^{\left\vert \lambda\right\vert }%
}{\left(  \lambda-1\right)  !}\right)  ^{2}\int_{\mathbf{0}}^{\mathbf{1}}%
\int_{\mathbf{0}}^{\mathbf{1}}\left(  1-\sigma\right)  ^{\lambda-1}\left(
1-\tau\right)  ^{\lambda-1}\sum\limits_{\alpha\in\mathbb{Z}^{d}\setminus
0}f_{\alpha}^{\left(  \lambda\right)  }\left(  \phi_{\nu};\sigma,\tau\right)
d\sigma d\tau,\label{1.077}%
\end{align}

where%
\[
\sum\limits_{\alpha\in\mathbb{Z}^{d}\setminus0}f_{\alpha}^{\left(
\lambda\right)  }\left(  \phi_{\nu};\sigma,\tau\right)  \leq\sum
\limits_{\alpha\in\mathbb{Z}^{d}\setminus0}\int\limits_{-\frac{\pi}%
{2}\mathbf{1}}^{\frac{\pi}{2}\mathbf{1}}\left\vert \left(  D^{\lambda}%
\phi_{\nu}\right)  \left(  \pi\alpha+\sigma.\eta\right)  \right\vert
\left\vert \left(  D^{\lambda}\overline{\phi_{\nu}}\right)  \left(  \pi
\alpha+\tau.\eta\right)  \right\vert d\eta.
\]

Now to estimate the right side of the last inequality using the
Cauchy-Schwartz inequality on $L^{2}$ and then on infinite sequences of real
numbers i.e. $l^{2}$:%
\begin{align*}
&  \sum\limits_{\alpha\in\mathbb{Z}^{d}\setminus0}f_{\alpha}^{\left(
\lambda\right)  }\left(  \phi_{\nu};\sigma,\tau\right) \\
&  \leq\sum\limits_{\alpha\in\mathbb{Z}^{d}\setminus0}\left(  \int%
\limits_{-\frac{\pi}{2}\mathbf{1}}^{\frac{\pi}{2}\mathbf{1}}\left\vert \left(
D^{\lambda}\phi_{\nu}\right)  \left(  \pi\alpha+\sigma.\eta\right)
\right\vert ^{2}d\eta\right)  ^{\frac{1}{2}}\left(  \int\limits_{-\frac{\pi
}{2}\mathbf{1}}^{\frac{\pi}{2}\mathbf{1}}\left\vert \left(  D^{\lambda}%
\phi_{\nu}\right)  \left(  \pi\alpha+\tau.\eta\right)  \right\vert ^{2}%
d\eta\right)  ^{\frac{1}{2}}\\
&  \leq\left(  \sum\limits_{\alpha\in\mathbb{Z}^{d}\setminus0}\int%
\limits_{-\frac{\pi}{2}\mathbf{1}}^{\frac{\pi}{2}\mathbf{1}}\left\vert \left(
D^{\lambda}\phi_{\nu}\right)  \left(  \pi\alpha+\sigma.\eta\right)
\right\vert ^{2}d\eta\right)  ^{\frac{1}{2}}\left(  \sum\limits_{\alpha
\in\mathbb{Z}^{d}\setminus0}\int\limits_{-\frac{\pi}{2}\mathbf{1}}^{\frac{\pi
}{2}\mathbf{1}}\left\vert \left(  D^{\mathbf{\lambda}}\phi_{\nu}\right)
\left(  \pi\alpha+\tau.\eta\right)  \right\vert ^{2}d\eta\right)  ^{\frac
{1}{2}}%
\end{align*}

With no loss of generality we can assume $\mathbf{0.}<\sigma\leq\mathbf{1}$
and the change of variables $\xi=\sigma.\eta$, $d\xi=\sigma^{\mathbf{1}}d\eta$
yields%
\begin{align*}
\sum\limits_{\alpha\in\mathbb{Z}^{d}\setminus0}\int_{-\frac{\pi}{2}\mathbf{1}%
}^{\frac{\pi}{2}\mathbf{1}}\left\vert \left(  D^{\lambda}\phi_{\nu}\right)
\left(  \pi\alpha+\sigma.\eta\right)  \right\vert ^{2}d\eta & =\sum
\limits_{\alpha\in\mathbb{Z}^{d}\setminus0}\frac{1}{\sigma^{\mathbf{1}}}%
\int_{-\frac{\pi}{2}\sigma}^{\frac{\pi}{2}\sigma}\left\vert \left(
D^{\lambda}\phi_{\nu}\right)  \left(  \pi\alpha+\xi\right)  \right\vert
^{2}d\xi\\
& \leq\frac{1}{\sigma^{\mathbf{1}}}\sum\limits_{\alpha\in\mathbb{Z}%
^{d}\setminus0}\int_{-\frac{\pi}{2}\mathbf{1}}^{\frac{\pi}{2}\mathbf{1}%
}\left\vert \left(  D^{\lambda}\phi_{\nu}\right)  \left(  \pi\alpha
+\xi\right)  \right\vert ^{2}d\xi\\
& \leq\frac{1}{\sigma^{\mathbf{1}}}\int\left\vert D^{\lambda}\phi_{\nu
}\right\vert ^{2}=:\frac{1}{\sigma^{\mathbf{1}}}\left\Vert D^{\lambda}%
\phi_{\nu}\right\Vert _{2}^{2},
\end{align*}

so that%
\[
\sum\limits_{\alpha\in\mathbb{Z}^{d}\setminus0}f_{\alpha}^{\left(
\lambda\right)  }\left(  \eta;\sigma,\tau\right)  \leq\frac{1}{\sigma
^{\mathbf{1}/2}\tau^{\mathbf{1}/2}}\left\Vert D^{\lambda}\phi_{\nu}\right\Vert
_{2}^{2},
\]

and consequently \ref{1.077} becomes%
\begin{align*}
\sum\limits_{\alpha\in\mathbb{Z}^{d}\setminus0} &  \int_{-\frac{\pi}%
{2}\mathbf{1}}^{\frac{\pi}{2}\mathbf{1}}\frac{\left\vert \phi_{\nu}\left(
\pi\alpha+\eta\right)  \right\vert ^{2}}{\sin^{2\lambda}\eta}d\eta\\
&  \leq\frac{\left(  \pi/2\right)  ^{2\left\vert \lambda\right\vert }}{\left(
\left(  \lambda-1\right)  !\right)  ^{2}}\int_{\mathbf{0}}^{\mathbf{1}}%
\int_{\mathbf{0}}^{\mathbf{1}}\left(  1-\sigma\right)  ^{\lambda-1}\left(
1-\tau\right)  ^{\lambda-1}\sum\limits_{\alpha\in\mathbb{Z}^{d}\setminus
0}f_{\alpha}^{\left(  \lambda\right)  }\left(  \phi_{\nu};\sigma,\tau\right)
d\sigma d\tau\\
&  \leq\frac{\left(  \pi/2\right)  ^{2\left\vert \lambda\right\vert }}{\left(
\left(  \lambda-1\right)  !\right)  ^{2}}\left(  \int_{\mathbf{0}}%
^{\mathbf{1}}\int_{\mathbf{0}}^{\mathbf{1}}\left(  1-\sigma\right)
^{\lambda-1}\left(  1-\tau\right)  ^{\lambda-1}\frac{d\sigma d\tau}%
{\sigma^{\mathbf{1}/2}\tau^{\mathbf{1}/2}}\right)  \left\Vert D^{\lambda}%
\phi_{\nu}\right\Vert _{2}^{2}\\
&  =\frac{\left(  \pi/2\right)  ^{2\left\vert \lambda\right\vert }}{\left(
\left(  \lambda-1\right)  !\right)  ^{2}}\left(  \int_{\mathbf{0}}%
^{\mathbf{1}}\frac{\left(  1-\sigma\right)  ^{\lambda-1}}{\sigma
^{\mathbf{1}/2}}d\sigma\right)  ^{2}\left\Vert D^{\lambda}\phi_{\nu
}\right\Vert _{2}^{2}\\
&  =\frac{\left(  \pi/2\right)  ^{2\left\vert \lambda\right\vert }}{\left(
\left(  \lambda-1\right)  !\right)  ^{2}}\left(  \int_{0}^{1}\frac{\left(
1-t\right)  ^{\lambda_{i}-1}}{\sqrt{t}}dt\right)  ^{2\mathbf{1}}\left\Vert
D^{\lambda}\phi_{\nu}\right\Vert _{2}^{2}\\
&  =\frac{\left(  \pi/2\right)  ^{2\left\vert \lambda\right\vert }}{\left(
\left(  \lambda-1\right)  !\right)  ^{2}}\left(  \int_{0}^{1}\frac
{t^{\lambda_{i}-1}}{\sqrt{1-t}}dt\right)  ^{2\mathbf{1}}\left\Vert D^{\lambda
}\phi_{\nu}\right\Vert _{2}^{2}\\
&  =\left(  \frac{\left(  \pi/2\right)  ^{\lambda_{i}}}{\left(  \lambda
_{i}-1\right)  !}\int_{0}^{1}\frac{t^{\lambda_{i}-1}}{\sqrt{1-t}}dt\right)
^{2\mathbf{1}}\left\Vert D^{\lambda}\phi_{\nu}\right\Vert _{2}^{2}.
\end{align*}

From the definite integral 855.34 of Dwight \cite{Dwight61},%
\begin{align}
\int_{0}^{1}\frac{t^{m}dt}{\sqrt{1-t^{n}}}  & =\frac{\sqrt{\pi}}{n}%
\frac{\Gamma\left(  \frac{m+1}{n}\right)  }{\Gamma\left(  \frac{m+1}{n}%
+\frac{1}{2}\right)  }=\frac{1}{n}B\left(  \frac{m+1}{n},\frac{1}{2}\right)
,\label{2.03}\\
\int_{0}^{1}\frac{t^{\lambda_{i}-1}dt}{\sqrt{1-t}}  & =\frac{\sqrt{\pi}}%
{1}\frac{\Gamma\left(  \frac{\lambda_{i}}{1}\right)  }{\Gamma\left(
\frac{\lambda_{i}}{1}+\frac{1}{2}\right)  }=\sqrt{\pi}\frac{\Gamma\left(
\lambda_{i}\right)  }{\Gamma\left(  \lambda_{i}+\frac{1}{2}\right)  }=B\left(
\lambda_{i},\frac{1}{2}\right)  ,\nonumber\\
\frac{\left(  \pi/2\right)  ^{\lambda_{i}}}{\left(  \lambda_{i}-1\right)
!}\int_{0}^{1}\frac{t^{\lambda_{i}-1}dt}{\sqrt{1-t}}  & =\frac{\left(
\pi/2\right)  ^{\lambda_{i}}}{\left(  \lambda_{i}-1\right)  !}\sqrt{\pi}%
\frac{\Gamma\left(  \lambda_{i}\right)  }{\Gamma\left(  \lambda_{i}+\frac
{1}{2}\right)  }=\frac{\left(  \pi/2\right)  ^{\lambda_{i}}\sqrt{\pi}}%
{\Gamma\left(  \lambda_{i}+\frac{1}{2}\right)  }=\nonumber\\
& =\frac{\left(  \pi/2\right)  ^{\lambda_{i}}}{\left(  2\lambda_{i}-1\right)
!!2^{-\lambda_{i}}}=\frac{\pi^{\lambda_{i}}}{\left(  2\lambda_{i}-1\right)
!!},\nonumber
\end{align}

where we have used the standard notation%
\begin{equation}
m!!:=\left\{
\begin{array}
[c]{ll}%
1, & m=-1,\\
1\cdot3\cdot5\cdots\cdot m, & m=1,3,5,\ldots,
\end{array}
\right. \label{2.09}%
\end{equation}

and shown that for $\lambda\geq\nu\geq\mathbf{1}$,%
\begin{align}
\sum\limits_{\alpha\in\mathbb{Z}^{d}\setminus0}\int_{-\frac{\pi}{2}\mathbf{1}%
}^{\frac{\pi}{2}\mathbf{1}}\frac{\left\vert \phi_{\nu}\left(  \pi\alpha
+\eta\right)  \right\vert ^{2}}{\sin^{2\lambda}\eta}d\eta & \leq\left(
\frac{\pi^{\lambda_{i}}}{\left(  2\lambda_{i}-1\right)  !!}\right)
^{2\mathbf{1}}\left\Vert D^{\lambda}\phi_{\nu}\right\Vert _{2}^{2}\nonumber\\
& =\left(  \frac{\pi^{\left\vert \lambda\right\vert }}{\left(  2\lambda
-1\right)  !!}\right)  ^{2}\left\Vert D^{\lambda}\phi_{\nu}\right\Vert
_{2}^{2}.\label{2.07}%
\end{align}

Using \ref{a1.021},%
\[
D^{\beta}x^{\nu}=\left\{
\begin{array}
[c]{ll}%
\beta!\binom{\nu}{\beta}x^{\nu-\beta}, & \beta\leq\nu,\\
0, & otherwise,
\end{array}
\right.
\]

we get%
\begin{align*}
D^{\lambda}\phi_{\nu}  & =D^{\lambda}\left(  \eta^{\nu}\phi\right)
=\sum\limits_{\beta\leq\lambda}\tbinom{\lambda}{\beta}\left(  D^{\lambda
-\beta}\eta^{\nu}\right)  D^{\beta}\phi=\\
& =\sum\limits_{\beta=\lambda-\nu}^{\lambda}\tbinom{\lambda}{\beta}\left(
\lambda-\beta\right)  !\binom{\nu}{\lambda-\beta}\eta^{\beta-\left(
\lambda-\nu\right)  }D^{\beta}\phi\\
& =\nu!\sum\limits_{\beta=\lambda-\nu}^{\lambda}\tbinom{\lambda}{\beta}%
\frac{\eta^{\beta-\left(  \lambda-\nu\right)  }}{\left(  \beta-\left(
\lambda-\nu\right)  \right)  !}D^{\beta}\phi,
\end{align*}

so that%
\[
\left\Vert D^{\lambda}\phi_{\nu}\right\Vert _{2}^{2}\leq\left\vert
\nu\right\vert \left(  \nu!\right)  ^{2}\sum\limits_{\beta=\lambda-\nu
}^{\lambda}\left(  \frac{\tbinom{\lambda}{\beta}}{\left(  \beta-\left(
\lambda-\nu\right)  \right)  !}\right)  ^{2}\left\Vert \eta^{\beta-\left(
\lambda-\nu\right)  }D^{\beta}\phi\right\Vert _{2}^{2},
\]

and hence%
\[
\sum\limits_{\alpha\in\mathbb{Z}^{d}\setminus0}\int_{-\frac{\pi}{2}\mathbf{1}%
}^{\frac{\pi}{2}\mathbf{1}}\frac{\left\vert \phi_{\nu}\left(  \pi\alpha
+\eta\right)  \right\vert ^{2}}{\sin^{2\lambda}\eta}d\eta\leq\left\vert
\nu\right\vert \left(  \frac{\pi^{\left\vert \lambda\right\vert }\nu!}{\left(
2\lambda-1\right)  !!}\right)  ^{2}\sum\limits_{\beta=\lambda-\nu}^{\lambda
}\left(  \frac{\tbinom{\lambda}{\beta}}{\left(  \beta-\left(  \lambda
-\nu\right)  \right)  !}\right)  ^{2}\left\Vert \eta^{\beta-\left(
\lambda-\nu\right)  }D^{\beta}\phi\right\Vert _{2}^{2}.
\]
\smallskip

\underline{\textbf{The first term of} \ref{a1.70}} If $\lambda=\nu$ then%
\begin{align}
\int_{-\frac{\pi}{2}\mathbf{1}}^{\frac{\pi}{2}\mathbf{1}}\frac{\eta^{2\lambda
}}{\sin^{2\lambda}\eta}\frac{\left\vert \phi\left(  \eta\right)  \right\vert
^{2}}{\eta^{2\left(  \lambda-\nu\right)  }}d\eta & \leq\left(  \frac{\pi}%
{2}\right)  ^{2\left\vert \lambda\right\vert }\int_{-\frac{\pi}{2}\mathbf{1}%
}^{\frac{\pi}{2}\mathbf{1}}\left\vert \phi\left(  \eta\right)  \right\vert
^{2}d\eta\nonumber\\
& \leq\left(  \frac{\pi}{2}\right)  ^{2\left\vert \lambda\right\vert
}\left\Vert \phi\right\Vert _{2}^{2},\quad\lambda=\nu.\label{2.08}%
\end{align}

There remains the case $\lambda>\nu$. We use \ref{a1.131} to obtain
\begin{align}
&  \int_{-\frac{\pi}{2}\mathbf{1}}^{\frac{\pi}{2}\mathbf{1}}\frac
{\eta^{2\lambda}}{\sin^{2\nu}\eta}\frac{\left\vert \phi\left(  \eta\right)
\right\vert ^{2}}{\eta^{2\left(  \lambda-\nu\right)  }}d\eta\nonumber\\
&  \leq\left(  \frac{\pi}{2}\right)  ^{2\left\vert \lambda\right\vert }%
\int_{-\frac{\pi}{2}\mathbf{1}}^{\frac{\pi}{2}\mathbf{1}}\frac{\left\vert
\phi\left(  \eta\right)  \right\vert ^{2}}{\eta^{2\left(  \lambda-\nu\right)
}}d\eta\nonumber\\
&  =\left(  \frac{\pi}{2}\right)  ^{2\left\vert \lambda\right\vert }%
\int_{-\frac{\pi}{2}\mathbf{1}}^{\frac{\pi}{2}\mathbf{1}}\frac{1}%
{\eta^{2\left(  \lambda-\nu\right)  }}\left\vert \frac{\eta^{\lambda-\nu}%
}{\left(  \lambda-\nu-1\right)  !}\int_{\mathbf{0}}^{\mathbf{1}}\left(
1-\tau\right)  ^{\lambda-\nu-1}\left(  D^{\lambda-\nu}\phi\right)  \left(
\tau.\eta\right)  d\tau\right\vert ^{2}d\eta\nonumber\\
&  =\tfrac{\left(  \pi/2\right)  ^{2\left\vert \lambda\right\vert }}{\left(
\left(  \lambda-\nu-1\right)  !\right)  ^{2}}\int_{-\frac{\pi}{2}\mathbf{1}%
}^{\frac{\pi}{2}\mathbf{1}}\left\vert \int_{\mathbf{0}}^{\mathbf{1}}\left(
1-\tau\right)  ^{\lambda-\nu-1}\left(  D^{\lambda-\nu}\phi\right)  \left(
\tau.\eta\right)  d\tau\right\vert ^{2}d\eta\nonumber\\
&  =\tfrac{\left(  \pi/2\right)  ^{2\left\vert \lambda\right\vert }}{\left(
\left(  \lambda-\nu-1\right)  !\right)  ^{2}}\int\limits_{-\frac{\pi}%
{2}\mathbf{1}}^{\frac{\pi}{2}\mathbf{1}}\left(  \int\limits_{\mathbf{0}%
}^{\mathbf{1}}\left(  1-\tau\right)  ^{\lambda-\nu-1}\left(  D^{\lambda-\nu
}\phi\right)  \left(  \tau.\eta\right)  d\tau\right)  \left(  \int%
\limits_{\mathbf{0}}^{\mathbf{1}}\left(  1-\sigma\right)  ^{\lambda-\nu
-1}\left(  D^{\lambda-\nu}\overline{\phi}\right)  \left(  \sigma.\eta\right)
d\sigma\right)  d\eta\nonumber\\
&  =\tfrac{\left(  \pi/2\right)  ^{2\left\vert \lambda\right\vert }}{\left(
\left(  \lambda-\nu-1\right)  !\right)  ^{2}}\int\limits_{-\frac{\pi}%
{2}\mathbf{1}}^{\frac{\pi}{2}\mathbf{1}}\int\limits_{\mathbf{0}}^{\mathbf{1}%
}\int\limits_{\mathbf{0}}^{\mathbf{1}}\left(  1-\sigma\right)  ^{\lambda
-\nu-1}\left(  1-\tau\right)  ^{\lambda-\nu-1}\left(  D^{\lambda-\nu}%
\phi\right)  \left(  \tau.\eta\right)  \left(  D^{\lambda-\nu}\overline{\phi
}\right)  \left(  \sigma.\eta\right)  d\sigma d\tau d\eta\nonumber\\
&  =\tfrac{\left(  \pi/2\right)  ^{2\left\vert \lambda\right\vert }}{\left(
\left(  \lambda-\nu-1\right)  !\right)  ^{2}}\int\limits_{\mathbf{0}%
}^{\mathbf{1}}\int\limits_{\mathbf{0}}^{\mathbf{1}}\left(  1-\sigma\right)
^{\lambda-\nu-1}\left(  1-\tau\right)  ^{\lambda-\nu-1}\int\limits_{-\frac
{\pi}{2}\mathbf{1}}^{\frac{\pi}{2}\mathbf{1}}\left(  D^{\lambda-\nu}%
\phi\right)  \left(  \tau.\eta\right)  \left(  D^{\lambda-\nu}\overline{\phi
}\right)  \left(  \sigma.\eta\right)  d\eta\text{ }d\sigma d\tau\nonumber\\
&  \leq\tfrac{\left(  \pi/2\right)  ^{2\left\vert \lambda\right\vert }%
}{\left(  \left(  \lambda-\nu-1\right)  !\right)  ^{2}}\int_{\mathbf{0}%
}^{\mathbf{1}}\int_{\mathbf{0}}^{\mathbf{1}}\left(  1-\sigma\right)
^{\lambda-\nu-1}\left(  1-\tau\right)  ^{\lambda-\nu-1}f_{0}^{\left(
\lambda-\nu\right)  }\left(  \phi;\sigma,\tau\right)  d\sigma d\tau
,\label{1.079}%
\end{align}

where the function $f_{0}^{\left(  \lambda-\nu\right)  }$ is given by
\ref{1.078}. Again applying the Cauchy-Schwartz inequality we get%
\begin{align*}
f_{0}^{\left(  \lambda-\nu\right)  }\left(  \phi;\sigma,\tau\right)   &
\leq\int_{-\frac{\pi}{2}\mathbf{1}}^{\frac{\pi}{2}\mathbf{1}}\left\vert
\left(  D^{\lambda-\nu}\phi\right)  \left(  \sigma.\eta\right)  \right\vert
\left\vert \left(  D^{\lambda-\nu}\overline{\phi}\right)  \left(  \tau
.\eta\right)  \right\vert d\eta\\
& \leq\left(  \int_{-\frac{\pi}{2}\mathbf{1}}^{\frac{\pi}{2}\mathbf{1}%
}\left\vert \left(  D^{\lambda-\nu}\phi\right)  \left(  \sigma.\eta\right)
\right\vert ^{2}d\eta\right)  ^{1/2}\left(  \int_{-\frac{\pi}{2}\mathbf{1}%
}^{\frac{\pi}{2}\mathbf{1}}\left\vert \left(  D^{\lambda-\nu}\phi\right)
\left(  \tau.\eta\right)  \right\vert ^{2}d\eta\right)  ^{1/2}.
\end{align*}

and again assuming $\mathbf{0}<\sigma,\tau\leq\mathbf{1}$ and applying the
change of variables $\xi=\sigma.\eta$, $d\xi=\sigma^{\mathbf{1}}d\eta$ yields%
\begin{align*}
f_{0}^{\left(  \lambda-\nu\right)  }\left(  \eta;\sigma,\tau\right)   &
\leq\frac{1}{\sigma^{\mathbf{1}/2}\tau^{\mathbf{1}/2}}\left(  \int_{-\frac
{\pi}{2}\sigma}^{\frac{\pi}{2}\sigma}\left\vert D^{\lambda-\nu}\phi\right\vert
^{2}\right)  ^{1/2}\left(  \int_{-\frac{\pi}{2}\tau}^{\frac{\pi}{2}\tau
}\left\vert D^{\lambda-\nu}\phi\right\vert ^{2}\right)  ^{1/2}\\
& \leq\frac{1}{\sigma^{\mathbf{1}/2}\tau^{\mathbf{1}/2}}\left(  \int%
_{-\frac{\pi}{2}\mathbf{1}}^{\frac{\pi}{2}\mathbf{1}}\left\vert D^{\lambda
-\nu}\phi\right\vert ^{2}\right) \\
& \leq\frac{1}{\sigma^{\mathbf{1}/2}\tau^{\mathbf{1}/2}}\left\Vert
D^{\lambda-\nu}\phi\right\Vert _{2}^{2}.
\end{align*}

This inequality applied to \ref{1.079} allows the estimate%
\begin{align*}
\int_{-\frac{\pi}{2}\mathbf{1}}^{\frac{\pi}{2}\mathbf{1}} &  \frac
{\eta^{2\lambda}}{\sin^{2\lambda}\eta}\frac{\left\vert \phi\left(
\eta\right)  \right\vert ^{2}}{\eta^{2\left(  \lambda-\nu\right)  }}d\eta\\
&  \leq\frac{\left(  \pi/2\right)  ^{2\left\vert \lambda\right\vert }}{\left(
\left(  \lambda-\nu-1\right)  !\right)  ^{2}}\left(  \int_{\mathbf{0}%
}^{\mathbf{1}}\int_{\mathbf{0}}^{\mathbf{1}}\frac{\left(  1-\sigma\right)
^{\lambda-\nu-1}\left(  1-\tau\right)  ^{\lambda-\nu-1}}{\sigma^{\mathbf{1}%
/2}\tau^{\mathbf{1}/2}}d\sigma d\tau\right)  \left\Vert D^{\lambda-\nu}%
\phi\right\Vert _{2}^{2}\\
&  =\frac{\left(  \pi/2\right)  ^{2\left\vert \lambda\right\vert }}{\left(
\left(  \lambda-\nu-1\right)  !\right)  ^{2}}\left(  \int_{\mathbf{0}%
}^{\mathbf{1}}\frac{\left(  1-\sigma_{i}\right)  ^{\lambda_{i}-\nu_{i}-1}%
}{\sigma_{i}^{1/2}}d\sigma_{i}\right)  ^{2\mathbf{1}}\left\Vert D^{\lambda
-\nu}\phi\right\Vert _{2}^{2}\\
&  =\frac{\left(  \pi/2\right)  ^{2\left\vert \lambda\right\vert }}{\left(
\lambda-\nu-1\right)  !}\left(  \int_{0}^{1}\frac{t^{\lambda_{i}-\nu_{i}-1}%
}{\left(  1-t\right)  ^{1/2}}dt\right)  ^{2\mathbf{1}}\left\Vert D^{\left(
\lambda-\nu\right)  \mathbf{1}}\phi\right\Vert _{2}^{2}.
\end{align*}

But from \ref{2.03},%
\[
\int_{0}^{1}\frac{t^{\lambda_{i}-\nu_{i}-1}dt}{\sqrt{1-t}}=\frac{\sqrt{\pi}%
}{1}\frac{\Gamma\left(  \frac{\lambda_{i}-\nu_{i}}{1}\right)  }{\Gamma\left(
\frac{\lambda_{i}-\nu_{i}}{1}+\frac{1}{2}\right)  }=\sqrt{\pi}\frac
{\Gamma\left(  \lambda_{i}-\nu_{i}\right)  }{\Gamma\left(  \lambda_{i}-\nu
_{i}+\frac{1}{2}\right)  },
\]

so when $\lambda>\nu$, by 850.1 of Dwight \cite{Dwight61}, for example,%
\begin{align}
\int_{-\frac{\pi}{2}\mathbf{1}}^{\frac{\pi}{2}\mathbf{1}}\frac{\eta^{2\nu}%
}{\sin^{2\lambda}\eta}\left\vert \phi\left(  \eta\right)  \right\vert
^{2}d\eta & \leq\left(  \frac{\left(  \pi/2\right)  ^{2\left\vert
\lambda\right\vert }}{\left(  \lambda-\nu-1\right)  !}\sqrt{\pi}\frac
{\Gamma\left(  \lambda-\nu\right)  }{\Gamma\left(  \lambda-\nu+\frac{1}%
{2}\right)  }\right)  ^{2d}\left\Vert D^{\lambda-\nu}\phi\right\Vert _{2}%
^{2}\nonumber\\
& =\left(  \left(  \frac{\pi}{2}\right)  ^{\left\vert \lambda\right\vert
}\frac{\sqrt{\pi}}{\Gamma\left(  \lambda-\nu+\frac{1}{2}\right)  }\right)
^{2}\left\Vert D^{\lambda-\nu}\phi\right\Vert _{2}^{2}\nonumber\\
& =\left(  \left(  \frac{\pi}{2}\right)  ^{\left\vert \lambda\right\vert
}\frac{2^{\left\vert \lambda\right\vert -\left\vert \nu\right\vert }}{\left(
2\left(  \lambda-\nu\right)  -1\right)  !!}\right)  ^{2}\left\Vert
D^{\lambda-\nu}\phi\right\Vert _{2}^{2},\quad\lambda>\nu.\label{2.06}%
\end{align}

The inequalities \ref{2.06} and \ref{2.08} can now be combined as%
\begin{equation}
\int_{-\frac{\pi}{2}\mathbf{1}}^{\frac{\pi}{2}\mathbf{1}}\frac{\eta^{2\nu}%
}{\sin^{2\lambda}\eta}\left\vert \phi\left(  \eta\right)  \right\vert
^{2}d\eta\leq\left(  \left(  \frac{\pi}{2}\right)  ^{\left\vert \lambda
\right\vert }\frac{2^{\left\vert \lambda\right\vert -\left\vert \nu\right\vert
}}{\left(  2\left(  \lambda-\nu\right)  -1\right)  !!}\right)  ^{2}\left\Vert
D^{\lambda-\nu}\phi\right\Vert _{2}^{2},\quad\lambda\geq\nu.\label{2.41}%
\end{equation}

Noting \ref{a1.012}, the inequalities \ref{2.41} and \ref{2.07} can be
combined into \ref{a1.70} to obtain part 1 of this theorem:%
\begin{align*}
\int &  \frac{\eta^{2\nu}}{\sin^{2\lambda}\eta}\left\vert \phi\left(
\eta\right)  \right\vert ^{2}d\eta\\
&  \leq\left(  \left(  \frac{\pi}{2}\right)  ^{\left\vert \lambda\right\vert
}\frac{2^{\left\vert \lambda-\nu\right\vert }}{\left(  2\left(  \lambda
-\nu\right)  -1\right)  !!}\right)  ^{2}\left\Vert D^{\lambda-\nu}%
\phi\right\Vert _{2}^{2}+\left(  \frac{\pi^{\left\vert \lambda\right\vert }%
}{\left(  2\lambda-1\right)  !!}\right)  ^{2}\left\Vert D^{\lambda}\left(
x^{\nu}\phi\right)  \right\Vert _{2}^{2}\\
&  =\left(  \frac{\pi}{2}\right)  ^{2\left\vert \lambda\right\vert }\left(
\left(  \frac{2^{\left\vert \lambda-\nu\right\vert }}{\left(  2\left(
\lambda-\nu\right)  -1\right)  !!}\right)  ^{2}\left\Vert D^{\lambda-\nu}%
\phi\right\Vert _{2}^{2}+\left(  \frac{2^{\left\vert \lambda\right\vert }%
}{\left(  2\lambda-1\right)  !!}\right)  ^{2}\left\Vert D^{\lambda}\left(
x^{\nu}\phi\right)  \right\Vert _{2}^{2}\right)  .
\end{align*}
\bigskip

\fbox{\textbf{Part 2a}} See attempts in previous files.\medskip

\fbox{\textbf{Part 2b}} See statement of this part.\medskip

\fbox{\textbf{Part 3}} From Theorem \ref{Thm_Sw2_extsplin_wt}, $\phi\in
S_{\otimes,\lambda-\nu}\left(  \mathbf{0}\right)  $ and $\phi_{\nu}:=\eta
^{\nu}\phi\in S_{\otimes,\lambda}\left(  \pi\mathbb{Z}^{d}\right)  $, so that
by part 3 of Lemma \ref{Lem_Sxm_expansion}:%
\begin{align}
\phi\left(  \eta\right)   & =\frac{\eta^{\lambda-\nu}}{\left(  \lambda
-\nu\right)  !}\left(  D^{\lambda-\nu}\phi\right)  \left(  t^{\prime}\left(
\eta\right)  .\eta\right)  ,\text{ }\eta\in\mathbb{R}^{d}.\label{a3.06}\\
\phi_{\nu}\left(  \pi\alpha+\eta\right)   & =\frac{\eta^{\lambda}}{\lambda
!}\left(  D^{\lambda}\phi_{\nu}\right)  \left(  \pi\alpha+t\left(
\eta\right)  .\eta\right)  ,\text{\quad}\alpha\in\mathbb{Z}^{d},\text{ }%
\eta\in R\left[  -\frac{\pi}{2}\mathbf{1,}\frac{\pi}{2}\mathbf{1}\right]
.\label{a3.07}%
\end{align}

Thus, using \ref{a3.07},%
\begin{align*}
\left\Vert \frac{\eta^{\nu}}{\sin^{\lambda}\eta}\phi\left(  \eta\right)
\right\Vert _{\infty}  & \leq\max_{\alpha\in\mathbb{Z}^{d}}\max_{\eta\in
R\left[  -\frac{\pi}{2}\mathbf{1,}\frac{\pi}{2}\mathbf{1}\right]  }\left\vert
\frac{\eta^{\mathbf{\lambda}}}{\sin^{\mathbf{\lambda}}\eta}\frac{\phi_{\nu
}\left(  \pi\alpha+\eta\right)  }{\eta^{\lambda}}\right\vert \\
& \leq\left(  \frac{\pi}{2}\right)  ^{\left\vert \lambda\right\vert }%
\max_{\alpha\in\mathbb{Z}^{d}}\max_{\eta\in R\left[  -\frac{\pi}{2}%
\mathbf{1,}\frac{\pi}{2}\mathbf{1}\right]  }\left\vert \frac{\phi_{\nu}\left(
\pi\alpha+\eta\right)  }{\eta^{\lambda}}\right\vert \\
& =\left(  \frac{\pi}{2}\right)  ^{\left\vert \lambda\right\vert }\max
_{\alpha\in\mathbb{Z}^{d}}\max_{\eta\in R\left[  -\frac{\pi}{2}\mathbf{1,}%
\frac{\pi}{2}\mathbf{1}\right]  }\left\vert \frac{\left(  D^{\lambda}\phi
_{\nu}\right)  \left(  \pi\alpha+t\left(  \eta\right)  .\eta\right)  }%
{\lambda!}\right\vert \\
& \leq\frac{\left(  \pi/2\right)  ^{\left\vert \lambda\right\vert }}{\lambda
!}\left\Vert D^{\lambda}\phi_{\nu}\right\Vert _{\infty}.
\end{align*}

See also part 8 of Lemma \ref{Lem_Sxm_integ_expansion}.

Clearly $\frac{\eta^{\nu}}{\sin^{\lambda}\eta}\phi\in C^{\infty}\left(
\mathbb{R}^{d}\setminus\left(  \pi\mathbb{Z}^{d}\setminus0\right)  \right)  $.
We now show that $\frac{\eta^{\nu}}{\sin^{\lambda}\eta}\phi\in C_{B}^{\infty
}\left(  \mathbb{R}^{d}\right)  $.

To this end suppose $\xi\in R\left(  \pi\alpha-\frac{3}{4}\pi\mathbf{1,}%
\pi\alpha+\frac{3}{4}\pi\mathbf{1}\right)  $ and set
\begin{equation}
\xi=\pi\alpha+\eta,\label{a3.09}%
\end{equation}

so that $\eta\in R\left(  -\frac{3}{4}\pi\mathbf{1,}\frac{3}{4}\pi
\mathbf{1}\right)  $. Then from \ref{a1.13},
\begin{align}
\frac{\xi^{\nu}}{\sin^{\lambda}\xi}\phi\left(  \xi\right)   & =\frac{\phi
_{\nu}\left(  \xi\right)  }{\sin^{\lambda}\xi}=\frac{\phi_{\nu}\left(
\pi\alpha+\eta\right)  }{\sin^{\lambda}\eta}=\nonumber\\
& =\frac{\left(  -1\right)  ^{\left\vert \lambda\right\vert }}{\left(
\lambda-1\right)  !}\frac{\eta^{\lambda}}{\sin^{\lambda}\eta}\int_{\mathbf{0}%
}^{\mathbf{1}}\left(  1-\tau\right)  ^{\lambda-1}\left(  D^{\lambda}\phi_{\nu
}\right)  \left(  \pi\alpha+\tau.\eta\right)  d\tau.\label{a3.04}%
\end{align}

But on $R\left[  -\frac{3}{4}\pi\mathbf{1,}\frac{3}{4}\pi\mathbf{1}\right]  $,
$\frac{\eta^{\lambda}}{\sin^{\lambda}\eta}\in C^{\infty}$ and $\int%
_{\mathbf{0}}^{\mathbf{1}}\left(  1-\tau\right)  ^{\lambda-1}\left(
D^{\lambda}\phi_{\nu}\right)  \left(  \pi\alpha+\tau.\eta\right)  d\tau\in
C^{\infty}$ which is easily proved by applying Lemma
\ref{Lem_diff_under_integral_2} to show that differentiation under the
integral sign is valid. Thus, for each $\beta$,%
\begin{align*}
D_{\eta}^{\beta}\int_{\mathbf{0}}^{\mathbf{1}}\left(  1-\tau\right)
^{\lambda-1}\left(  D^{\lambda}\phi_{\nu}\right)  \left(  \pi\alpha+\tau
.\eta\right)  d\tau & =\int_{\mathbf{0}}^{\mathbf{1}}\left(  1-\tau\right)
^{\lambda-1}\tau^{\beta}\left(  D^{\beta+\lambda}\phi_{\nu}\right)  \left(
\pi\alpha+\tau.\eta\right)  d\tau\\
& \in C_{B}^{\left(  0\right)  }\left(  \mathbb{R}^{d}\right)  ,
\end{align*}

and so%
\begin{align}
D^{\beta} &  \left(  \frac{\xi^{\nu}}{\sin^{\lambda}\xi}\phi\left(
\xi\right)  \right) \nonumber\\
&  =\frac{\left(  -1\right)  ^{\left\vert \lambda\right\vert }}{\left(
\lambda-1\right)  !}\sum_{\gamma\leq\beta}\tbinom{\beta}{\gamma}%
D^{\beta-\gamma}\left(  \frac{\eta^{\lambda}}{\sin^{\lambda}\eta}\right)
\int_{\mathbf{0}}^{\mathbf{1}}\left(  1-\tau\right)  ^{\lambda-1}\tau^{\gamma
}\left(  D^{\gamma+\lambda}\phi_{\nu}\right)  \left(  \pi\alpha+\tau
.\eta\right)  d\tau,\label{a2.29}%
\end{align}

and we can conclude that $\frac{\xi^{\nu}}{\sin^{\lambda}\xi}\phi\left(
\xi\right)  \in C_{B}^{\infty}\left(  \mathbb{R}^{d}\setminus R\left(
-\frac{\pi}{2}\mathbf{1,}\frac{\pi}{2}\mathbf{1}\right)  \right)  $.

Thus $\frac{\xi^{\nu}}{\sin^{\lambda}\xi}\phi\left(  \xi\right)  \in
C_{B}^{\infty}\left(  \mathbb{R}^{d}\right)  $ and it remains to be proven
that $\frac{\xi^{\nu}}{\sin^{\lambda}\xi}\phi\left(  \xi\right)  \in S$.
Noting that $\xi=\pi\alpha+\eta$ and $I_{\alpha}=R\left[  \pi\alpha-\frac{\pi
}{2}\mathbf{1,}\pi\alpha+\frac{\pi}{2}\mathbf{1}\right]  $, equation
\ref{a2.29} yields the straight-forward estimate
\[
\left\Vert \left\vert \xi\right\vert ^{n}D^{\beta}\left(  \frac{\xi^{\nu}%
}{\sin^{\lambda}\xi}\phi\left(  \xi\right)  \right)  \right\Vert
_{\infty,I_{\alpha}}\leq\frac{1}{\lambda!}\sum_{\gamma\leq\beta}\tbinom{\beta
}{\gamma}\left\Vert D^{\beta-\gamma}\left(  \frac{\eta^{\lambda}}%
{\sin^{\lambda}\eta}\right)  \right\Vert _{\infty,I_{0}}\left(  \max
_{I_{\alpha}}\left\vert \xi\right\vert _{2}^{n}\right)  \left\Vert
D^{\gamma+\lambda}\phi_{\nu}\right\Vert _{\infty,I_{\alpha}}.
\]
\medskip

\fbox{If $\alpha=\mathbf{0}$,}\medskip%
\begin{align}
&  \left\Vert \xi^{\sigma}D^{\beta}\left(  \frac{\xi^{\nu}}{\sin^{\lambda}\xi
}\phi\left(  \xi\right)  \right)  \right\Vert _{\infty,I_{\alpha}}\nonumber\\
&  \leq\frac{1}{\lambda!}\sum_{\gamma\leq\beta}\tbinom{\beta}{\gamma
}\left\Vert D^{\beta-\gamma}\left(  \frac{\eta^{\lambda}}{\sin^{\lambda}\eta
}\right)  \right\Vert _{\infty,I_{0}}\left(  \max_{I_{0}}\left\vert
\xi\right\vert _{2}\right)  ^{n}\left\Vert D^{\gamma+\lambda}\phi_{\nu
}\right\Vert _{\infty,I_{0}}\nonumber\\
&  \leq\frac{1}{\lambda!}\sum_{\gamma\leq\beta}\tbinom{\beta}{\gamma
}\left\Vert D^{\beta-\gamma}\left(  \frac{\eta^{\lambda}}{\sin^{\lambda}\eta
}\right)  \right\Vert _{\infty,I_{0}}\left(  \frac{\pi}{2}\sqrt{d}\right)
^{n}\left\Vert D^{\gamma+\lambda}\phi_{\nu}\right\Vert _{\infty,I_{0}%
}\nonumber\\
&  \leq\frac{\left(  2\sqrt{d}\right)  ^{n}}{\lambda!}\left(  \sum_{\gamma
\leq\beta}\tbinom{\beta}{\gamma}\left\Vert D^{\gamma}\left(  \frac
{\eta^{\lambda}}{\sin^{\lambda}\eta}\right)  \right\Vert _{\infty,R\left[
-\frac{\pi}{2}\mathbf{1,}\frac{\pi}{2}\mathbf{1}\right]  }\right)  \left(
\frac{\pi}{4}\right)  ^{n}\max_{\gamma\leq\beta}\left\Vert D^{\gamma+\lambda
}\phi_{\nu}\right\Vert _{\infty}.\label{a3.03}%
\end{align}
\medskip

\fbox{If $\alpha\neq\mathbf{0}$,}\medskip%
\[
\max_{I_{\alpha}}\left\vert \cdot\right\vert _{e}=\max_{\eta\in I_{0}%
}\left\vert \pi\alpha+\eta\right\vert _{e}\leq\max_{\eta\in I_{0}}\left(
\pi\left\vert \alpha\right\vert _{e}+\left\vert \eta\right\vert _{e}\right)
\leq\pi\left\vert \alpha\right\vert _{e}+\max_{\eta\in I_{0}}\left\vert
\eta\right\vert _{e}=\pi\left\vert \alpha\right\vert _{e}+\frac{\pi}{2}%
\sqrt{d},
\]

and if%
\[
\left(  \alpha_{+}>\mathbf{0}\right)  _{i}:=\left\{
\begin{array}
[c]{ll}%
0, & \left(  \alpha_{+}\right)  _{i}=0,\\
1, & \left(  \alpha_{+}\right)  _{i}=1.
\end{array}
\right.
\]

then%
\begin{align*}
\min_{I_{\alpha}}\left\vert \cdot\right\vert _{e}=\min_{\eta\in I_{0}%
}\left\vert \pi\alpha+\eta\right\vert =\min_{\eta\in I_{0}}\left\vert
\pi\alpha_{+}+\eta\right\vert = &  \text{??}\left\vert \pi\alpha_{+}-\frac
{\pi}{2}\left(  \alpha_{+}>\mathbf{0}\right)  \right\vert \text{ }\left(
\text{\textbf{justify}!}\right) \\
&  \geq\left\vert \pi\alpha_{+}\right\vert -\frac{\pi}{2}\left\vert \left(
\alpha_{+}>\mathbf{0}\right)  \right\vert \\
&  =\pi\left\vert \alpha\right\vert -\frac{\pi}{2}\left\vert \left(
\alpha_{+}>\mathbf{0}\right)  \right\vert .
\end{align*}

Hence%
\begin{align*}
\frac{\max_{I_{\alpha}}\left\vert \cdot\right\vert }{\min_{I_{\alpha}%
}\left\vert \cdot\right\vert }  & \leq\frac{\pi\left\vert \alpha\right\vert
+\frac{\pi}{2}\sqrt{d}}{\pi\left\vert \alpha\right\vert -\frac{\pi}%
{2}\left\vert \left(  \alpha_{+}>\mathbf{0}\right)  \right\vert }\\
& =\frac{\pi\left\vert \alpha\right\vert -\frac{\pi}{2}\left\vert \left(
\alpha_{+}>\mathbf{0}\right)  \right\vert +\frac{\pi}{2}\left\vert \left(
\alpha_{+}>\mathbf{0}\right)  \right\vert +\frac{\pi}{2}\sqrt{d}}%
{\pi\left\vert \alpha\right\vert -\frac{\pi}{2}\left\vert \left(  \alpha
_{+}>\mathbf{0}\right)  \right\vert }=\\
& =1+\frac{\frac{1}{2}\left\vert \left(  \alpha_{+}>\mathbf{0}\right)
\right\vert +\frac{1}{2}\sqrt{d}}{\left\vert \alpha\right\vert -\frac{1}%
{2}\left\vert \left(  \alpha_{+}>\mathbf{0}\right)  \right\vert }\leq
1+\frac{\frac{1}{2}\sqrt{d}+\frac{1}{2}\sqrt{d}}{\left\vert \alpha\right\vert
-\frac{1}{2}\left\vert \left(  \alpha_{+}>\mathbf{0}\right)  \right\vert }=\\
& =1+\frac{\sqrt{d}}{\left\vert \alpha\right\vert -\frac{1}{2}\left\vert
\left(  \alpha_{+}>\mathbf{0}\right)  \right\vert }\leq1+\frac{\sqrt{d}%
}{\left\vert \left(  \alpha_{+}>\mathbf{0}\right)  \right\vert -\frac{1}%
{2}\left\vert \left(  \alpha_{+}>\mathbf{0}\right)  \right\vert }=\\
& =1+\frac{\sqrt{d}}{\frac{1}{2}\left\vert \left(  \alpha_{+}>\mathbf{0}%
\right)  \right\vert }\\
& \leq1+2\sqrt{d}=:c_{d},
\end{align*}

so when $\alpha\neq\mathbf{0}$,%
\begin{align*}
&  \left\Vert \left\vert \xi\right\vert ^{n}D^{\beta}\left(  \frac{\xi^{\nu}%
}{\sin^{\lambda}\xi}\phi\left(  \xi\right)  \right)  \right\Vert
_{\infty,I_{\alpha}}\\
&  \leq\left(  \frac{1}{\lambda!}\sum_{\gamma\leq\beta}\tbinom{\beta}{\gamma
}\left\Vert D^{\beta-\gamma}\left(  \frac{\eta^{\lambda}}{\sin^{\lambda}\eta
}\right)  \right\Vert _{\infty,I_{0}}\right)  \left(  c_{d}\right)
^{n}\left(  \min_{I_{\alpha}}\left\vert \cdot\right\vert ^{n}\right)
\max_{\gamma\leq\beta}\left\Vert D^{\gamma+\lambda}\phi_{\nu}\right\Vert
_{\infty,I_{\alpha}}\\
&  \leq\frac{\left(  c_{d}\right)  ^{n}}{\lambda!}\left(  \sum_{\gamma
\leq\beta}\tbinom{\beta}{\gamma}\left\Vert D^{\gamma}\left(  \frac
{\eta^{\lambda}}{\sin^{\lambda}\eta}\right)  \right\Vert _{\infty,I_{0}%
}\right)  \max_{\gamma\leq\beta}\left\Vert \left\vert \cdot\right\vert
^{n}D^{\gamma+\lambda}\phi_{\nu}\right\Vert _{\infty,I_{\alpha}}.
\end{align*}

But $\phi_{\nu}\in S$ so $\max\limits_{\gamma\leq\beta}\left\Vert \left\vert
\cdot\right\vert ^{n}D^{\gamma+\lambda}\phi_{\nu}\right\Vert _{\infty
,I_{\alpha}}\leq\max\limits_{\gamma\leq\beta}\left\Vert \left\vert
\cdot\right\vert ^{n}D^{\gamma+\lambda}\phi_{\nu}\right\Vert _{\infty}<\infty$
which implies that $\frac{\xi^{\nu}}{\sin^{\lambda}\xi}\phi\left(  \xi\right)
\in S$.\medskip

\fbox{Part 4} From part 3,%
\begin{align*}
&  \left\Vert \left\vert \xi\right\vert ^{n}D^{\beta}\left(  \frac{\xi^{\nu}%
}{\sin^{\lambda}\xi}\phi\left(  \xi\right)  \right)  \right\Vert
_{\infty,I_{\alpha}}\\
&  \leq\frac{\left(  1+2\sqrt{d}\right)  ^{n}}{\lambda!}\left(  \sum
_{\gamma\leq\beta}\tbinom{\beta}{\gamma}\left\Vert D^{\gamma}\left(
\frac{\eta^{\lambda}}{\sin^{\lambda}\eta}\right)  \right\Vert _{\infty
,R\left[  -\frac{\pi}{2}\mathbf{1,}\frac{\pi}{2}\mathbf{1}\right]  }\right)
\max_{\gamma\leq\beta}\left\Vert \left\vert \cdot\right\vert ^{n}%
D^{\gamma+\lambda}\phi_{\nu}\right\Vert _{\infty},
\end{align*}

and combining this estimate with that of \ref{a3.03} we get%
\begin{align*}
\left\vert \xi\right\vert ^{n}\left\vert D^{\beta}\left(  \frac{\xi^{\nu}%
}{\sin^{\lambda}\xi}\phi\left(  \xi\right)  \right)  \right\vert  &  \leq
\frac{\left(  1+2\sqrt{d}\right)  ^{n}}{\lambda!}\left(  \sum_{\gamma\leq
\beta}\tbinom{\beta}{\gamma}\left\Vert D^{\gamma}\left(  \frac{\eta^{\lambda}%
}{\sin^{\lambda}\eta}\right)  \right\Vert _{\infty,R\left[  -\frac{\pi}%
{2}\mathbf{1,}\frac{\pi}{2}\mathbf{1}\right]  }\right)  \times\\
&  \qquad\qquad\times\left\{  \max_{\gamma\leq\beta}\left(  \frac{\pi}%
{4}\right)  ^{n}\left\Vert D^{\gamma+\lambda}\phi_{\nu}\right\Vert _{\infty
},\max_{\gamma\leq\beta}\left\Vert \left\vert \cdot\right\vert ^{n}%
D^{\gamma+\lambda}\phi_{\nu}\right\Vert _{\infty}\right\} \\
&  \leq\frac{\left(  1+2\sqrt{d}\right)  ^{n}}{\lambda!}\left(  \sum
_{\gamma\leq\beta}\tbinom{\beta}{\gamma}\left\Vert D^{\gamma}\left(
\frac{\eta^{\lambda}}{\sin^{\lambda}\eta}\right)  \right\Vert _{\infty
,R\left[  -\frac{\pi}{2}\mathbf{1,}\frac{\pi}{2}\mathbf{1}\right]  }\right)
\times\\
&  \qquad\qquad\times\max_{\gamma\leq\beta}\left\Vert \left(  1+\left\vert
\cdot\right\vert ^{n}\right)  D^{\gamma+\lambda}\left(  x^{\nu}\phi\right)
\right\Vert _{\infty}.
\end{align*}
\medskip

\fbox{\textbf{Part 5}} From part 3,%
\[
\left\vert \frac{\eta^{\nu}}{\sin^{\lambda}\eta}\phi\left(  \eta\right)
\right\vert \leq\frac{\left(  \pi/2\right)  ^{\left\vert \lambda\right\vert }%
}{\lambda!}\left\Vert D^{\lambda}\phi_{\nu}\right\Vert _{\infty},
\]

so that%
\begin{equation}
\left\vert \phi\left(  \eta\right)  \right\vert \leq\frac{\left(
\pi/2\right)  ^{\left\vert \lambda\right\vert }}{\lambda!}\left\Vert
D^{\lambda}\phi_{\nu}\right\Vert _{\infty}\left\vert \frac{\sin^{\lambda}\eta
}{\eta^{\nu}}\right\vert .\label{a3.05}%
\end{equation}
\medskip

\fbox{\textbf{Part 6}} To prove the inequality we use part 1 to obtain
\[
\int\left\vert \phi\right\vert =\int\frac{\sin^{\lambda}\eta}{\eta^{\nu}}%
\frac{\eta^{\nu}}{\sin^{\lambda}\eta}\left\vert \phi\left(  \eta\right)
\right\vert d\eta\leq\left(  \int\frac{\sin^{2\lambda}\eta}{\eta^{2\nu}%
}\right)  ^{\frac{1}{2}}\left(  \int\frac{\eta^{2\nu}}{\sin^{2\lambda}\eta
}\left\vert \phi\left(  \eta\right)  \right\vert ^{2}d\eta\right)  ^{\frac
{1}{2}}.
\]
\medskip

\fbox{\textbf{Part 7}} If $w=\frac{\sin^{2\lambda}\eta}{\eta^{2\nu}}$ then%
\[
\int\left\vert \phi\right\vert ^{2}=\int\frac{1}{w}w\left\vert \phi\right\vert
^{2}\leq\left\Vert \frac{1}{w}\right\Vert \int w\left\vert \phi\right\vert
^{2}.
\]

\end{proof}

\begin{corollary}
\label{Cor_Thm_eta^2/sin(eta)^2.|phi|^2_dim_gt_3}When $\beta\leq\lambda-\nu$
the inequalities of Theorem \ref{Thm_Sw2_extsplin_wt} also hold with $\phi$
replaced by $D^{\beta}\phi$ and $\lambda$ replaced by $\lambda-\beta$.
\end{corollary}

\begin{proof}
This corollary follows directly from part 3 of Theorem
\ref{Thm_Sw2_extsplin_wt} i.e. $D^{\beta}S_{w_{\nu,\lambda};2}\subset
S_{w_{\nu,\lambda-\beta};2}$ when $\beta\leq\lambda-\nu$.
\end{proof}

\begin{remark}
?? From \ref{a3.04},%
\[
\frac{\xi^{\nu}}{\sin^{\lambda}\xi}\phi\left(  \xi\right)  =\frac{\left(
-1\right)  ^{\left\vert \lambda\right\vert }}{\left(  \lambda-1\right)
!}\frac{\eta^{\lambda}}{\sin^{\lambda}\eta}\int_{\mathbf{0}}^{\mathbf{1}%
}\left(  1-\tau\right)  ^{\lambda-1}\left(  D^{\lambda}\phi_{\nu}\right)
\left(  \pi\alpha+\tau.\eta\right)  d\tau,\text{ }\phi\in S_{w,0}.
\]

Write%
\[
\frac{\xi^{\nu}}{\sin^{\lambda}\xi}\overset{\wedge}{\overset{\vee}{\phi}%
}\left(  \xi\right)  =\frac{\left(  -1\right)  ^{\left\vert \lambda\right\vert
}}{\left(  \lambda-1\right)  !}\frac{\eta^{\lambda}}{\sin^{\lambda}\eta}%
\int_{\mathbf{0}}^{\mathbf{1}}\left(  1-\tau\right)  ^{\lambda-1}\left(
D^{\lambda}\left(  x^{\nu}\overset{\wedge}{\overset{\vee}{\phi}}\right)
\right)  \left(  \pi\alpha+\tau.\eta\right)  d\tau.
\]

Replace $\overset{\vee}{\phi}$ by $u\in X_{w}^{0}$:%
\[
\frac{\xi^{\nu}}{\sin^{\lambda}\xi}\overset{\wedge}{u}\left(  \xi\right)
=\frac{\left(  -1\right)  ^{\left\vert \lambda\right\vert }}{\left(
\lambda-1\right)  !}\frac{\eta^{\lambda}}{\sin^{\lambda}\eta}\int_{\mathbf{0}%
}^{\mathbf{1}}\left(  1-\tau\right)  ^{\lambda-1}\left(  D^{\lambda}\left(
x^{\nu}\overset{\wedge}{u}\right)  \right)  \left(  \pi\alpha+\tau
.\eta\right)  d\tau.
\]

Does this mean anything?
\end{remark}

We now slightly generalize part 4 of Theorem
\ref{Thm_eta^2/sin(eta)^2.|phi|^2_dim_gt_1}.\medskip

\begin{corollary}
\label{Cor_Thm_eta^2/sin(eta)^2.|phi|^2_dim_gt_2}\ 

\begin{enumerate}
\item Suppose $w$ is a B-spline weight function with parameters $\nu$ and
$\lambda$. Suppose $\phi\in S_{w,0}$. By Theorem
\ref{Thm_eta^2/sin(eta)^2.|phi|^2_dim_gt_1} $\frac{\xi^{\nu}}{\sin^{\lambda
}\xi}\phi\in S$, and we show below that for all $\alpha,\beta\geq\mathbf{0}$:%
\begin{multline*}
\left\vert \xi^{\alpha}D^{\beta}\left(  \frac{\xi^{\nu}}{\sin^{\lambda}\xi
}\phi\left(  \xi\right)  \right)  \right\vert \\
\leq\frac{1}{\left(  \lambda-1\right)  !}\sum_{\substack{\gamma\leq
\beta\\\sigma\leq\alpha}}\tfrac{\tbinom{\alpha}{\sigma}\tbinom{\beta}{\gamma}%
}{\tbinom{\alpha-\sigma+\lambda-1+\gamma}{\gamma}\left(  \alpha-\sigma
+\lambda+\gamma\right)  }\left\Vert \eta^{\alpha-\sigma}D^{\beta-\gamma
}\left(  \frac{\eta^{\lambda}}{\sin^{\lambda}\eta}\right)  \right\Vert
_{\infty,I_{0}}\times\\
\times\max_{\substack{\gamma\leq\beta\\\sigma\leq\alpha}}\left\Vert x^{\sigma
}D^{\gamma+\lambda}\left(  x^{\nu}\phi\right)  \right\Vert _{\infty},
\end{multline*}

where $I_{0}=\left[  -\frac{\pi}{2}\mathbf{1},\frac{\pi}{2}\mathbf{1}\right]
$.

\item $w_{1}\phi:S_{w,0}\rightarrow S$ is continuous where $S_{w,0}$ is
endowed with the subspace topology.

\item ?? Derive an $L^{2}$ upper bound for $\int\left\vert \xi^{\alpha
}D^{\beta}\left(  \frac{\xi^{\nu}}{\sin^{\lambda}\xi}\phi\left(  \xi\right)
\right)  \right\vert ^{2}d\xi$.
\end{enumerate}
\end{corollary}

\begin{proof}
\textbf{Part 1} We start with equation \ref{a2.29} i.e.
\begin{align*}
D^{\beta} &  \left(  \frac{\xi^{\nu}}{\sin^{\lambda}\xi}\phi\left(
\xi\right)  \right) \\
&  =\frac{\left(  -1\right)  ^{\left\vert \lambda\right\vert }}{\left(
\lambda-1\right)  !}\sum_{\gamma\leq\beta}\tbinom{\beta}{\gamma}%
D^{\beta-\gamma}\left(  \frac{\eta^{\lambda}}{\sin^{\lambda}\eta}\right)
\int_{\mathbf{0}}^{\mathbf{1}}\left(  1-\tau\right)  ^{\lambda-1}\tau^{\gamma
}\left(  D^{\gamma+\lambda}\phi_{\nu}\right)  \left(  \pi\alpha+\tau
.\eta\right)  d\tau,
\end{align*}

where $\xi=\pi\alpha+\eta$, $\eta\in I_{0}$. Define $I_{\alpha}=\pi
\alpha+I_{0}$ for all $\alpha\in\mathbb{Z}^{d}$. The binomial theorem implies%
\begin{align*}
\xi^{\alpha}=\left(  \pi\alpha+\tau.\eta+\eta-\tau.\eta\right)  ^{\alpha}  &
=\sum_{\sigma\leq\alpha}\tbinom{\alpha}{\sigma}\left(  \pi\alpha+\tau
.\eta\right)  ^{\sigma}\left(  \eta-\tau.\eta\right)  ^{\alpha-\sigma}\\
& =\sum_{\sigma\leq\alpha}\tbinom{\alpha}{\sigma}\left(  \pi\alpha+\tau
.\eta\right)  ^{\sigma}\eta^{\alpha-\sigma}\left(  1-\tau\right)
^{\alpha-\sigma},
\end{align*}

so that%
\begin{align*}
&  \xi^{\alpha}D^{\beta}\left(  \frac{\xi^{\nu}}{\sin^{\lambda}\xi}\phi\left(
\xi\right)  \right) \\
&  =\frac{\left(  -1\right)  ^{\left\vert \lambda\right\vert }}{\left(
\lambda-1\right)  !}\sum_{\gamma\leq\beta}\tbinom{\beta}{\gamma}%
D^{\beta-\gamma}\left(  \frac{\eta^{\lambda}}{\sin^{\lambda}\eta}\right)
\int_{\mathbf{0}}^{\mathbf{1}}\left(  1-\tau\right)  ^{\lambda-1}\tau^{\gamma
}\xi^{\alpha}\left(  D^{\gamma+\lambda}\phi_{\nu}\right)  \left(  \pi
\alpha+\tau.\eta\right)  d\tau\\
&  =\frac{\left(  -1\right)  ^{\left\vert \lambda\right\vert }}{\left(
\lambda-1\right)  !}\sum_{\gamma\leq\beta}\tbinom{\beta}{\gamma}%
D^{\beta-\gamma}\left(  \frac{\eta^{\lambda}}{\sin^{\lambda}\eta}\right)
\times\\
&  \qquad\qquad\times\int_{\mathbf{0}}^{\mathbf{1}}\left(  1-\tau\right)
^{\lambda-1}\tau^{\gamma}\left(  \sum_{\sigma\leq\alpha}\tbinom{\alpha}%
{\sigma}\left(  \pi\alpha+\tau.\eta\right)  ^{\sigma}\eta^{\alpha-\sigma
}\left(  1-\tau\right)  ^{\alpha-\sigma}\right)  \left(  D^{\gamma+\lambda
}\phi_{\nu}\right)  \left(  \pi\alpha+\tau.\eta\right)  d\tau\\
&  =\frac{\left(  -1\right)  ^{\left\vert \lambda\right\vert }}{\left(
\lambda-1\right)  !}\sum_{\substack{\gamma\leq\beta\\\sigma\leq\alpha}%
}\tbinom{\alpha}{\sigma}\tbinom{\beta}{\gamma}\eta^{\alpha-\sigma}%
D^{\beta-\gamma}\left(  \frac{\eta^{\lambda}}{\sin^{\lambda}\eta}\right)
\times\\
&  \qquad\qquad\times\int_{\mathbf{0}}^{\mathbf{1}}\left(  1-\tau\right)
^{\alpha-\sigma+\lambda-1}\tau^{\gamma}\left(  \pi\alpha+\tau.\eta\right)
^{\sigma}\left(  D^{\gamma+\lambda}\phi_{\nu}\right)  \left(  \pi\alpha
+\tau.\eta\right)  d\tau.
\end{align*}

This yields the sequence of estimates:%
\begin{align}
&  \left\vert \xi^{\alpha}D^{\beta}\left(  \frac{\xi^{\nu}}{\sin^{\lambda}\xi
}\phi\left(  \xi\right)  \right)  \right\vert \nonumber\\
&  =\frac{1}{\left(  \lambda-1\right)  !}\sum_{\gamma\leq\beta,\sigma
\leq\alpha}\tbinom{\alpha}{\sigma}\tbinom{\beta}{\gamma}\left\vert
\eta^{\alpha-\sigma}D^{\beta-\gamma}\left(  \frac{\eta^{\lambda}}%
{\sin^{\lambda}\eta}\right)  \right\vert \times\nonumber\\
&  \qquad\qquad\qquad\qquad\times\int_{\mathbf{0}}^{\mathbf{1}}\left(
1-\tau\right)  ^{\alpha-\sigma+\lambda-1}\tau^{\gamma}\left\vert \left(
\pi\alpha+\tau.\eta\right)  ^{\sigma}\left(  D^{\gamma+\lambda}\phi_{\nu
}\right)  \left(  \pi\alpha+\tau.\eta\right)  \right\vert d\tau\nonumber\\
&  \leq\frac{1}{\left(  \lambda-1\right)  !}\sum_{\substack{\gamma\leq
\beta\\\sigma\leq\alpha}}\tbinom{\alpha}{\sigma}\tbinom{\beta}{\gamma
}\left\Vert \eta^{\alpha-\sigma}D^{\beta-\gamma}\left(  \frac{\eta^{\lambda}%
}{\sin^{\lambda}\eta}\right)  \right\Vert _{\infty,I_{0}}\times\nonumber\\
&  \qquad\qquad\qquad\qquad\times\int_{\mathbf{0}}^{\mathbf{1}}\left(
1-\tau\right)  ^{\alpha-\sigma+\lambda-1}\tau^{\gamma}\left\vert \left(
\pi\alpha+\tau.\eta\right)  ^{\sigma}\left(  D^{\gamma+\lambda}\phi_{\nu
}\right)  \left(  \pi\alpha+\tau.\eta\right)  \right\vert d\tau\label{1.082}\\
&  \leq\frac{1}{\left(  \lambda-1\right)  !}\sum_{\substack{\gamma\leq
\beta\\\sigma\leq\alpha}}\tbinom{\alpha}{\sigma}\tbinom{\beta}{\gamma
}\left\Vert \eta^{\alpha-\sigma}D^{\beta-\gamma}\left(  \frac{\eta^{\lambda}%
}{\sin^{\lambda}\eta}\right)  \right\Vert _{\infty,I_{0}}\left(
\int_{\mathbf{0}}^{\mathbf{1}}\left(  1-\tau\right)  ^{\alpha-\sigma
+\lambda-1}\tau^{\gamma}d\tau\right)  \left\Vert x^{\sigma}D^{\gamma+\lambda
}\phi_{\nu}\right\Vert _{\infty,I_{\alpha}}\nonumber\\
&  \Rightarrow beta\text{ }function\Rightarrow\nonumber\\
&  =\frac{1}{\left(  \lambda-1\right)  !}\sum_{\substack{\gamma\leq
\beta\\\sigma\leq\alpha}}\tbinom{\alpha}{\sigma}\tbinom{\beta}{\gamma
}\left\Vert \eta^{\alpha-\sigma}D^{\beta-\gamma}\left(  \frac{\eta^{\lambda}%
}{\sin^{\lambda}\eta}\right)  \right\Vert _{\infty,I_{0}}\tfrac{1}%
{\tbinom{\alpha-\sigma+\lambda-1+\gamma}{\gamma}\left(  \alpha-\sigma
+\lambda+\gamma\right)  }\left\Vert x^{\sigma}D^{\gamma+\lambda}\phi_{\nu
}\right\Vert _{\infty,I_{\alpha}}\nonumber\\
&  =\frac{1}{\left(  \lambda-1\right)  !}\sum_{\substack{\gamma\leq
\beta\\\sigma\leq\alpha}}\tfrac{\tbinom{\alpha}{\sigma}\tbinom{\beta}{\gamma}%
}{\tbinom{\alpha-\sigma+\lambda-1+\gamma}{\gamma}\left(  \alpha-\sigma
+\lambda+\gamma\right)  }\left\Vert \eta^{\alpha-\sigma}D^{\beta-\gamma
}\left(  \frac{\eta^{\lambda}}{\sin^{\lambda}\eta}\right)  \right\Vert
_{\infty,I_{0}}\left\Vert x^{\sigma}D^{\gamma+\lambda}\phi_{\nu}\right\Vert
_{\infty,I_{\alpha}}\nonumber\\
&  \leq\frac{1}{\left(  \lambda-1\right)  !}\left(  \sum_{\substack{\gamma
\leq\beta\\\sigma\leq\alpha}}\tfrac{\tbinom{\alpha}{\sigma}\tbinom{\beta
}{\gamma}}{\tbinom{\alpha-\sigma+\lambda-1+\gamma}{\gamma}\left(
\alpha-\sigma+\lambda+\gamma\right)  }\left\Vert \eta^{\alpha-\sigma}%
D^{\beta-\gamma}\left(  \frac{\eta^{\lambda}}{\sin^{\lambda}\eta}\right)
\right\Vert _{\infty,I_{0}}\right)  \max_{\substack{\gamma\leq\beta
\\\sigma\leq\alpha}}\left\Vert x^{\sigma}D^{\gamma+\lambda}\phi_{\nu
}\right\Vert _{\infty,I_{\alpha}}\nonumber\\
&  \leq\frac{1}{\left(  \lambda-1\right)  !}\left(  \sum_{\substack{\gamma
\leq\beta\\\sigma\leq\alpha}}\tfrac{\tbinom{\alpha}{\sigma}\tbinom{\beta
}{\gamma}}{\tbinom{\alpha-\sigma+\lambda-1+\gamma}{\gamma}\left(
\alpha-\sigma+\lambda+\gamma\right)  }\left\Vert \eta^{\alpha-\sigma}%
D^{\beta-\gamma}\left(  \frac{\eta^{\lambda}}{\sin^{\lambda}\eta}\right)
\right\Vert _{\infty,I_{0}}\right)  \max_{\substack{\gamma\leq\beta
\\\sigma\leq\alpha}}\left\Vert x^{\sigma}D^{\gamma+\lambda}\left(  x^{\nu}%
\phi\right)  \right\Vert _{\infty}.\nonumber
\end{align}
\medskip

\textbf{Part 2} Follows directly from the inequality proved in part 1.\medskip

\textbf{Part 3} ?? Copy proof of part 1 of Theorem
\ref{Thm_eta^2/sin(eta)^2.|phi|^2_dim_gt_1}.
\end{proof}

Recall that $\tau_{a}$ denotes the translation operator $\tau_{a}f\left(
x\right)  =f\left(  x-a\right)  $.

\begin{corollary}
\label{Cor_Thm_eta^2/sin(eta)^2.|phi|^2_dim_gt_1}If $w$ is a B-spline weight
function with parameters $\nu$ and $\lambda$ then:

\begin{enumerate}
\item $w\in W_{S;0}$ and%
\[
\int w\left\vert \phi\right\vert ^{2}\leq\left(  \frac{\pi}{2}\right)
^{2\left\vert \lambda\right\vert }\left(  \left(  \frac{2^{\left\vert
\lambda-\nu\right\vert }}{\left(  2\left(  \lambda-\nu\right)  -1\right)
!!}\right)  ^{2}\left\Vert D^{\lambda-\nu}\phi\right\Vert _{2}^{2}+\left(
\frac{2^{\left\vert \lambda\right\vert }}{\left(  2\lambda-1\right)
!!}\right)  ^{2}\left\Vert D^{\lambda}\left(  x^{\nu}\phi\right)  \right\Vert
_{2}^{2}\right)  ,
\]

when $\phi\in S_{w,0}$. Here we have used the notation \ref{2.09}.

\item $S_{w,0}=\frac{\sin^{\lambda}\xi}{\xi^{\nu}}S$, $\frac{\sin^{\lambda}%
\xi}{\xi^{\nu}}\in L^{2}$ and $\overset{\vee}{S}_{w,0}=\left(  \frac
{\sin^{\lambda}\xi}{\xi^{\nu}}\right)  ^{\vee}\ast S\in L^{r}$ for any $2\leq
r\leq\infty$. Here where we are using the "$L^{p}$" convolution described by
Young's inequality in part 6 of Definition \ref{Def_convol}.

\item $\overset{\vee}{S}_{w,0}\ast\overset{\vee}{S}_{w,0}=G\ast S\ast S\subset
G\ast S$.

\item $G\ast S\subset\overset{\vee}{S}_{w,0}$.

\item $\overset{\vee}{S}_{w,0}\ast\overset{\vee}{S}_{w,0}\subset
\overset{\vee}{S}_{w,0}$.

\item If $\psi\in\overset{\vee}{S}_{w,0}$ then%
\[
\left\Vert \psi\right\Vert _{w,0}^{2}\leq\left(  \frac{\pi}{2}\right)
^{2\left\vert \lambda\right\vert }\left(  \left(  \frac{2^{\left\vert
\lambda-\nu\right\vert }}{\left(  2\left(  \lambda-\nu\right)  -1\right)
!!}\right)  ^{2}\left\Vert x^{\lambda-\nu}\psi\right\Vert _{2}^{2}+\left(
\frac{2^{\left\vert \lambda\right\vert }}{\left(  2\lambda-1\right)
!!}\right)  ^{2}\left\Vert x^{\lambda}D^{\nu}\psi\right\Vert _{2}^{2}\right)
.
\]

\end{enumerate}
\end{corollary}

\begin{proof}
\textbf{Part 1} Part 1 of Theorem \ref{Thm_eta^2/sin(eta)^2.|phi|^2_dim_gt_1}%
.\medskip

\textbf{Part 2} If $\phi\in S_{w,0}$ then by part 3 of Theorem
\ref{Thm_eta^2/sin(eta)^2.|phi|^2_dim_gt_1}, $\frac{\xi^{\nu}}{\sin^{\lambda
}\xi}\phi\in S$. If $\phi\in\frac{\sin^{\lambda}\xi}{\xi^{\nu}}S$ then
$\frac{\xi^{\nu}}{\sin^{\lambda}\xi}\phi\in S$ and so $\left\vert \frac
{\xi^{\nu}}{\sin^{\lambda}\xi}\phi\right\vert ^{2}=w\left\vert \phi\right\vert
^{2}\in S$ and consequently $\int w\left\vert \phi\right\vert ^{2}<\infty$
i.e. $\phi\in S_{w,0}$. Thus $S_{w,0}=\frac{\sin^{\lambda}\xi}{\xi^{\nu}}S$.

By definition $\frac{\sin^{2\lambda}\xi}{\xi^{2\nu}}\in L^{1}$ and so
$\frac{\sin^{\lambda}\xi}{\xi^{\nu}}\in L^{2}$. Also $S\subset L^{q}$ when
$1\leq q\leq\infty$ and so by Young's convolution estimate:
\[
\left(  2\pi\right)  ^{d/2}\left\Vert u\ast v\right\Vert _{r}\leq\left\Vert
u\right\Vert _{p}\left\Vert v\right\Vert _{q},\text{ }where\text{ }\left\{
\begin{array}
[c]{l}%
\frac{1}{p}+\frac{1}{q}=1+\frac{1}{r},\\
1\leq p,q,r\leq\infty,
\end{array}
\right.
\]

the convolution is an $L^{r}$ function for any $2\leq r\leq\infty$.\medskip

\textbf{Part 3} Extend Young's inequality to two convolutions and then write
\begin{align*}
\overset{\vee}{S}_{w,0}\ast\overset{\vee}{S}_{w,0}  & =\left(  \left(
\frac{\sin^{\lambda}\xi}{\xi^{\nu}}\right)  ^{\vee}\ast S\right)  \ast\left(
\left(  \frac{\sin^{\lambda}\xi}{\xi^{\nu}}\right)  ^{\vee}\ast S\right) \\
& =\left(  \left(  \frac{\sin^{\lambda}\xi}{\xi^{\nu}}\right)  ^{\vee}%
\ast\left(  \frac{\sin^{\lambda}\xi}{\xi^{\nu}}\right)  ^{\vee}\right)  \ast
S\ast S\\
& :part\text{ }9\text{ }of\text{ }Definition\text{ }\ref{Def_convol}%
\Rightarrow\\
& =\left(  \frac{\sin^{2\lambda}\xi}{\xi^{2\nu}}\right)  ^{\vee}\ast S\ast S\\
& =G\ast S\ast S\\
& :part\text{ }8\text{ }of\text{ }Definition\text{ }\ref{Def_convol}%
\Rightarrow\\
& \subset G\ast S.
\end{align*}
\medskip

\textbf{Part 4} $S_{w,0}=\frac{\sin^{\lambda}\xi}{\xi^{\nu}}S\Leftrightarrow
\frac{\sin^{\lambda}\xi}{\xi^{\nu}}S_{w,0}=\frac{1}{w}S\Rightarrow\frac{1}%
{w}S\subset S_{w,0}\Leftrightarrow\left(  \frac{1}{w}S\right)  ^{\vee}%
\subset\overset{\vee}{S}_{w,0}\Leftrightarrow G\ast S\subset\overset{\vee
}{S}_{w,0}$.\medskip

\textbf{Part 5} Apply parts 3 and 4.\medskip

\textbf{Part 6} If $\phi\in S_{w,0}$ then from part 1,%
\begin{align*}
&  \left\Vert \overset{\vee}{\phi}\right\Vert _{w,0}^{2}=\int\frac{\eta^{2\nu
}}{\sin^{2\lambda}\eta}\left\vert \widehat{\overset{\vee}{\phi}}\left(
\eta\right)  \right\vert ^{2}d\eta\leq\\
&  \leq\left(  \frac{\pi}{2}\right)  ^{2\left\vert \lambda\right\vert }\left(
\left(  \frac{2^{\left\vert \lambda-\nu\right\vert }}{\left(  2\left(
\lambda-\nu\right)  -1\right)  !!}\right)  ^{2}\left\Vert D^{\lambda-\nu
}\widehat{\overset{\vee}{\phi}}\right\Vert _{2}^{2}+\left(  \frac
{2^{\left\vert \lambda\right\vert }}{\left(  2\lambda-1\right)  !!}\right)
^{2}\left\Vert D^{\lambda}\left(  x^{\nu}\widehat{\overset{\vee}{\phi}%
}\right)  \right\Vert _{2}^{2}\right)  .\\
&  \leq\left(  \frac{\pi}{2}\right)  ^{2\left\vert \lambda\right\vert }\left(
\left(  \frac{2^{\left\vert \lambda-\nu\right\vert }}{\left(  2\left(
\lambda-\nu\right)  -1\right)  !!}\right)  ^{2}\left\Vert \widehat{x^{\lambda
-\nu}\overset{\vee}{\phi}}\right\Vert _{2}^{2}+\left(  \frac{2^{\left\vert
\lambda\right\vert }}{\left(  2\lambda-1\right)  !!}\right)  ^{2}\left\Vert
\widehat{x^{\lambda}D^{\nu}\overset{\vee}{\phi}}\right\Vert _{2}^{2}\right) \\
&  =\left(  \frac{\pi}{2}\right)  ^{2\left\vert \lambda\right\vert }\left(
\left(  \frac{2^{\left\vert \lambda-\nu\right\vert }}{\left(  2\left(
\lambda-\nu\right)  -1\right)  !!}\right)  ^{2}\left\Vert x^{\lambda-\nu
}\overset{\vee}{\phi}\right\Vert _{2}^{2}+\left(  \frac{2^{\left\vert
\lambda\right\vert }}{\left(  2\lambda-1\right)  !!}\right)  ^{2}\left\Vert
x^{\lambda}D^{\nu}\overset{\vee}{\phi}\right\Vert _{2}^{2}\right)  .
\end{align*}

\end{proof}

??? \textbf{ALSO SEE }Definition \ref{Def_tildL} and Corollary
\ref{Cor_invF[wSw2]_dense_tildXo1/w} in the section about the space
$\widetilde{X}_{1/w}^{0}$.

\begin{definition}
\label{Def_I1_J1}\textbf{\ The operators} $\mathcal{I}_{1}:X_{w}%
^{0}\rightarrow L^{2}$ and $\mathcal{J}_{1}:L^{2}\rightarrow X_{w}^{0}$. The
operators $\mathcal{I}:X_{w}^{0}\rightarrow L^{2}$ and $\mathcal{J}%
:L^{2}\rightarrow X_{w}^{0}$ were introduced in Definition \ref{Def_I_J} using
$\sqrt{w}$ and their properties described in Theorem \ref{Thm_I_J_property}.

Suppose instead we use a function $w_{1}$ such that $\left(  w_{1}\right)
^{2}=w$ and define the operators $\mathcal{I}_{1}$ and $\mathcal{J}_{1}$ by%
\begin{align*}
\mathcal{I}_{1}u  & :=\left(  w_{1}\widehat{u}\right)  ^{\vee},\quad u\in
X_{w}^{0}.\\
\mathcal{J}_{1}f  & :=\left(  \frac{1}{w_{1}}\widehat{f}\right)  ^{\vee},\quad
f\in L^{2}.
\end{align*}

Then it is easy to show that $\mathcal{I}_{1}:X_{w}^{0}\rightarrow L^{2}$ is
an isometric isomorphism (as was $\mathcal{I}$) with left/right inverse
$\mathcal{J}_{1}$ - $\mathcal{J}_{1}$ corresponding to the operator
$\mathcal{J}$.

Regarding Theorem \ref{Thm_I_J_property}, $\mathcal{I}_{1}$ and $\mathcal{J}%
_{1}$ have the same properties as $\mathcal{I}$ and $\mathcal{J}$ except with
$\sqrt{w}$ replaced by $w_{1}$.

Since $1/w\in L^{1}$ we have $1/w_{1}\in L^{2}$ and we define%
\[
G_{1}:=\left(  \frac{1}{w_{1}}\right)  ^{\vee}\in L^{2}.
\]

We can now express the basis function $G$ as
\[
G=\left(  \frac{1}{w}\right)  ^{\vee}=\left(  \frac{1}{w_{1}}\frac{1}{w_{1}%
}\right)  ^{\vee}=\left(  \frac{1}{w_{1}}\right)  ^{\vee}\ast\left(  \frac
{1}{w_{1}}\right)  ^{\vee}=G_{1}\ast G_{1},
\]

and clearly%
\[
\mathcal{J}_{1}f:=G_{1}\ast f,\quad f\in L^{2}.
\]

\end{definition}

\begin{lemma}
\label{Lem_S2_eq_S}$S^{2}:=S.S=S$ where $S.S$ indicates the pointwise product
of functions in $S$.
\end{lemma}

\begin{proof}
?? This result is proved in the Theorem in Paul Garrett's PDF document
"Weil-Schwartz envelopes for rapidly decreasing functions." dated December 21,
2004. HOW to reference the document? Garrett suggests cite URL, adding "[after Weil]"?

The document "faculty.gvsu.edu/alayontf/notes/dixmier\_malliavin.pdfShare"
proves the 1-dim case.\medskip
\end{proof}

??? \textbf{The following theorem is UNFINISHED}.

See also Remark ?? \ref{Rem_FT_wt_func} regarding the \textbf{Fourier
transform of a weight function}.

\begin{theorem}
\label{Thm_ops_IJL_invFSw2}Suppose $w$ is any B-spline weight function, say
$w=\frac{\xi^{2\nu}}{\sin^{2\lambda}\xi}$. Now choose $w_{1}=\frac{\xi^{\nu}%
}{\sin^{\lambda}\xi}$. Then:

\begin{enumerate}
\item $G_{1}\ast S=\overset{\vee}{S}_{w,0}$ and $\overset{\vee}{S}_{w,0}%
\ast\overset{\vee}{S}_{w,0}=G\ast S$ which is a vector space.

\item $\mathcal{I}_{1}:\overset{\vee}{S}_{w,0}\rightarrow S$ and is onto and
$\mathcal{J}_{1}:S\rightarrow\overset{\vee}{S}_{w,0}$ and is onto.

\item $\mathcal{I}_{1}:\overset{\vee}{S}_{w,0}\ast\overset{\vee}{S}%
_{w,0}\rightarrow\overset{\vee}{S}_{w,0}$ is onto and $\mathcal{J}%
_{1}:\overset{\vee}{S}_{w,0}\rightarrow\overset{\vee}{S}_{w,0}\ast
\overset{\vee}{S}_{w,0}$ is onto.

\item ?? $\mathcal{L}\widetilde{\mathcal{L}}$ restricted to $\overset{\vee
}{S}_{w,0}$ is $\mathcal{I}_{1}$.

\item ?? The operator $\mathcal{L}\widetilde{\mathcal{L}}:\overset{\vee
}{S}_{w,0}\rightarrow S$ and is 1-1 and onto?

\item ?? $\overset{\vee}{S}_{w,0}$ is dense in $L^{1}$ and $\overset{\vee
}{S}_{w,0}\ast\overset{\vee}{S}_{w,0}$ is dense in $X_{w}^{0}$.

\item \textbf{See summary Figure} \ref{Fig_ops_I1_J1_L_Sw2} below.
\end{enumerate}
\end{theorem}

\begin{proof}
\textbf{Part 1} From part 2 of Corollary
\ref{Cor_Thm_eta^2/sin(eta)^2.|phi|^2_dim_gt_1}, $S_{w,0}=\frac{1}{w_{1}}S$
i.e. $\overset{\vee}{S}_{w,0}=G_{1}\ast S$. Thus\allowbreak\ $\overset{\vee
}{S}_{w,0}\ast\overset{\vee}{S}_{w,0}=\left(  G_{1}\ast S\right)  \ast\left(
G_{1}\ast S\right)  =G\ast S$ which is a vector space.\medskip

\textbf{Part 2} If $\psi\in\overset{\vee}{S}_{w,0}$then $\mathcal{I}_{1}%
\psi=\left(  w_{1}\widehat{\psi}\right)  ^{\vee}\in S$ since $w_{1}S_{w,0}=S$.

If $\phi\in S$ then $\mathcal{J}_{1}\phi=\left(  \frac{1}{w_{1}}\widehat{\phi
}\right)  ^{\vee}\in\overset{\vee}{S}_{w,0}$. The onto-ness of the two
operators now follows from the fact that $\mathcal{J}_{1}\mathcal{I}%
_{1}=\mathcal{I}_{1}\mathcal{J}_{1}=I$.\medskip

\textbf{Part 3} Suppose $\psi_{1},\psi_{2}\in S_{w,0}$. We first show that
$\mathcal{I}_{1}\left(  \overset{\vee}{\psi_{1}}\ast\overset{\vee}{\psi_{2}%
}\right)  \in\overset{\vee}{S}_{w,0}$.

Now $\mathcal{I}_{1}\left(  \overset{\vee}{\psi_{1}}\ast\overset{\vee
}{\psi_{2}}\right)  \in\overset{\vee}{S}_{w,0}$ iff $w_{1}\left(
\overset{\vee}{\psi_{1}}\ast\overset{\vee}{\psi_{2}}\right)  ^{\wedge}\in
S_{w,0} $ iff $w_{1}\psi_{1}\psi_{2}\in S_{w,0}$ iff $w_{1}\psi_{1}w_{1}%
\psi_{2}\in w_{1}S_{w,0}=S$ which is clearly true.

Suppose $\psi_{3}\in S_{w,0}$. Then $\mathcal{J}_{1}\overset{\vee}{\psi_{3}%
}=\left(  \frac{1}{w_{1}}\psi_{3}\right)  ^{\vee}=G_{1}\ast\overset{\vee
}{\psi_{3}}$ and we require the existence of $\psi_{1},\psi_{2}\in S_{w,0}$
such that $G_{1}\ast\overset{\vee}{\psi_{3}}=\overset{\vee}{\psi_{1}}%
\ast\overset{\vee}{\psi_{2}}$ \allowbreak i.e. $\frac{1}{w_{1}}\psi_{3}%
=\psi_{1}\psi_{2}\Leftrightarrow\psi_{3}=w_{1}\psi_{1}\psi_{2}\Leftrightarrow
w_{1}\psi_{3}=w_{1}\psi_{1}w_{1}\psi_{2}\Leftrightarrow S.S=S^{2}=S$ and this
is true by Lemma \ref{Lem_S2_eq_S}.\medskip

\textbf{Part 4} ?? \medskip

\textbf{Part 5} ??\medskip

\textbf{Part 6} ??\medskip

\textbf{Part 7} Summary figure:%
\begin{equation}%
\begin{array}
[c]{ccccc}%
\widetilde{X}_{1/w}^{0},X_{1/w}^{0} & \overset{\widetilde{\mathcal{L}%
},\mathcal{L}}{\longleftarrow} & L^{2} & \overset{\mathcal{I}_{1}%
:i.i.}{\underset{\mathcal{J}_{1}}{{\leftrightarrows}}} & X_{w}^{0}\\
\uparrow &  & \uparrow\; &  & \uparrow\;\\
\uparrow\; &  & S & \overset{\mathcal{I}_{1}}{\longleftarrow} & \overset{\vee
}{S}_{w,0}=G_{1}\ast S\\
\uparrow &  & \uparrow &  & \uparrow\;\\
S & \overset{\widetilde{\mathcal{L}},\mathcal{L}}{\longleftarrow} &
\overset{\vee}{S}_{w,0} & \overset{\mathcal{I}_{1}}{\longleftarrow} &
\overset{\vee}{S}_{w,0}\ast\overset{\vee}{S}_{w,0}=G\ast S
\end{array}
\label{Fig_ops_I1_J1_L_Sw2}%
\end{equation}

\end{proof}

\begin{remark}
\label{Rem_Thm_ops_IJL_invFSw2}\textbf{1}) ?? Define $\Psi:S\times
S\rightarrow S$ by $\Psi\left(  \phi,\psi\right)  =\phi\psi$. Endow $S\times
S$ with a suitable sequence of seminorms. Show $\Psi$ is continuous. Lemma
\ref{Lem_S2_eq_S} implies $\Psi$ is onto.\medskip

\textbf{2}) Given $u\in S^{\prime}$ define $v:S\times S\rightarrow\mathbb{C}$
by $v\left(  \phi,\psi\right)  =\left[  u,\Psi\left(  \phi,\psi\right)
\right]  $. Write $v=\Psi^{\ast}$?\medskip

\textbf{3}) Part 2 of Corollary
\ref{Cor_Thm_eta^2/sin(eta)^2.|phi|^2_dim_gt_2} implies $w_{1}\phi
:S_{w,0}\rightarrow S$ is continuous and so $u\in S^{\prime}$ implies $\left[
u,w_{1}\phi\right]  \in S_{w,0}^{\prime}$.\medskip

\textbf{4}) $\left[  u,w_{1}\phi\right]  :S^{\prime}\rightarrow S_{w,0}%
^{\prime}$ is continuous.%
\[%
\begin{array}
[c]{ccc}%
S_{w,0} & \overset{w_{1}.}{\longrightarrow} & S\\
& \searrow & u\downarrow\in S^{\prime}\\
&  & \mathbb{C}%
\end{array}
\]
\smallskip

\textbf{5}) ?? $\frac{1}{w_{1}}\psi:S\rightarrow S_{w,0}$ is
continuous.\medskip

6) ??
\end{remark}

?? Is this theorem, in particular part 3, in the right place? See Definitions
\ref{Def_Sw2_lin_fnal}, \ref{Def_Sw2_lin_fnal_0} and \ref{Def_T_lin_fnal_0}
and other stuff.

\begin{theorem}
\label{Thm_FSw2_and_invFSw2}?? FIX! Suppose $w$ is a (homogeneous) extended
B-spline weight function i.e. $w=\frac{x^{2n\mathbf{1}}}{\sin^{2l\mathbf{1}}%
x}$. Then:

\begin{enumerate}
\item $\phi\in\widehat{S}_{w,0}$ iff $\overline{\phi}\in\overset{\vee
}{S}_{w,0}$.

\item If $w$ is even then $\widehat{S}_{w,0}=\overset{\vee}{S}_{w,0}$.

\item ?? $\left(  \widehat{S}_{w,0}\right)  ^{\prime}=\left(  S_{w,0}^{\prime
}\right)  ^{\vee}$ and $\left(  \overset{\vee}{S}_{w,0}\right)  ^{\prime
}=\left(  S_{w,0}^{\prime}\right)  ^{\wedge}$.

\item $\int w\phi\overline{\psi}$ is a bilinear form on $S_{w,0}\otimes
S_{w,0}$ and%
\[
\int w\phi\overline{\psi}\leq\left(  \int w\left\vert \phi\right\vert
^{2}\right)  ^{1/2}\left(  \int w\left\vert \psi\right\vert ^{2}\right)
^{1/2}.
\]

\end{enumerate}
\end{theorem}

\begin{proof}
The identities $\widehat{\widehat{\phi}}\left(  \xi\right)  =\phi\left(
-\xi\right)  $ and then $\overset{\vee}{\phi}\left(  -\xi\right)
=\widehat{\phi}\left(  \xi\right)  $ imply:\medskip

\textbf{Part 1} $\overline{\phi}\in\overset{\vee}{S}_{w,0}\Leftrightarrow
\phi\in\overset{\overline{\vee}}{S}_{w,0}=\overset{\wedge}{\overset{\_}{S}%
}_{w,0}=S_{w,0}$.\medskip

\textbf{Part 2} If $\phi\in\widehat{S}_{w,0}$ then $\int w\left\vert
\widehat{\phi}\right\vert ^{2}$, and since $w$ is even
\[
\int w\left\vert \overset{\vee}{\phi}\right\vert ^{2}=\int w\left(
\xi\right)  \left\vert \overset{\vee}{\phi}\left(  -\xi\right)  \right\vert
^{2}d\xi=\int w\left\vert \widehat{\phi}\right\vert ^{2},
\]

so that $\phi\in\overset{\vee}{S}_{w,0}$ and the argument is
reversible.\medskip

\textbf{Part 3} ?? FINISH! If $g\in\left(  \widehat{S}_{w,0}\right)  ^{\prime
}$ and $\psi\in S_{w,0}$ then $\left[  \widehat{g},\psi\right]  =\left[
g,\widehat{\psi}\right]  $.

?? If $g\in\left(  \overset{\vee}{S}_{w,0}\right)  ^{\prime}$ and $\psi\in
S_{w,0}$ then $\left[  \overset{\vee}{g},\psi\right]  =\left[  g,\overset{\vee
}{\psi}\right]  $.\medskip

\textbf{Part 4} Use the Cauchy-Schwartz inequality.
\end{proof}

\begin{remark}
\label{Rem_FT_wt_func}?? \textbf{Define the Fourier transform of a weight
function\medskip}

\textbf{1}) We use part 4 of Theorem \ref{Thm_FSw2_and_invFSw2} to define the
Fourier transform $\widehat{w}$ of a weight function $w$ on the subspace
$\overset{\vee}{S}_{w,0}\ast\overset{\vee}{S}_{w,0}=G\ast S$ of $\overset{\vee
}{S}_{w,0}$. In fact, for $\phi,\psi\in S_{w,0}$,%
\[
\left[  \widehat{w},\overset{\vee}{\phi}\ast\overset{\vee}{\overline{\psi}%
}\right]  :=\left[  w,\phi\overline{\psi}\right]  =\int w\phi\overline{\psi
}=\int w_{1}\phi w_{1}\overline{\psi},
\]

and the integral exists since by the Cauchy-Schwartz theorem%
\[
\left\vert \int w\phi\overline{\psi}\right\vert \leq\int w\left\vert
\phi\right\vert \left\vert \psi\right\vert =\int\sqrt{w}\left\vert
\phi\right\vert \sqrt{w}\left\vert \psi\right\vert \leq\left(  \int
w\left\vert \phi\right\vert ^{2}\right)  ^{1/2}\left(  \int w\left\vert
\psi\right\vert ^{2}\right)  ^{1/2}<\infty.
\]
\medskip

\textbf{2}) Further%
\[
\int w\left\vert \phi\right\vert ^{2}=\left[  w,\phi\overline{\phi}\right]
=\left[  \widehat{w},\overset{\vee}{\phi}\ast\overset{\vee}{\overline{\phi}%
}\right]  ,
\]

so that%
\[
\left\vert \left[  \widehat{w},\overset{\vee}{\phi}\ast\overset{\vee
}{\overline{\psi}}\right]  \right\vert \leq\left[  \widehat{w},\overset{\vee
}{\phi}\ast\overset{\vee}{\overline{\phi}}\right]  \left[  \widehat{w}%
,\overset{\vee}{\psi}\ast\overset{\vee}{\overline{\psi}}\right]  .
\]
\medskip

\textbf{3})%
\begin{align*}
\left\{  \overset{\vee}{\phi}\ast\overset{\vee}{\overline{\phi}}:\phi
\in\overset{\vee}{S}_{w,0}\right\}   & \subset\overset{\vee}{S}_{w,0}%
\ast\overset{\vee}{S}_{w,0}.\\
\left\{  \overset{\vee}{\phi}\ast\overset{\vee}{\overline{\phi}}:\phi
\in\overset{\vee}{S}_{w,0}\right\}   & =\left\{  \left(  \left\vert
\phi\right\vert ^{2}\right)  ^{\vee}:\phi\in\overset{\vee}{S}_{w,0}\right\}
=\left\{  \left\vert \phi\right\vert ^{2}:\phi\in\overset{\vee}{S}%
_{w,0}\right\}  ^{\vee}.
\end{align*}

Hence%
\[
\left[  \widehat{w},\overset{\vee}{\phi}\ast\overset{\vee}{\overline{\phi}%
}\right]  =\left[  \widehat{w},\left(  \left\vert \phi\right\vert ^{2}\right)
^{\vee}\right]  =\int w\left\vert \phi\right\vert ^{2}.
\]
\medskip

\textbf{4}) Noting part 1 we denote the mapping $\phi\rightarrow\left[
\widehat{w},\overset{\vee}{\phi}\ast\overset{\vee}{\overline{\psi}}\right]  $
by $\widehat{w}$ and it's action by $\widehat{w}\left(  \phi\right)  $. Part 1
implies it satisfies:
\[
\widehat{w}:S_{w,0}\rightarrow S_{w,0}^{\prime},
\]

and is continuous in the sense that if $\phi_{k}\rightarrow\phi$ in $S_{w,0}$
then $\widehat{w}\left(  \phi_{k}\right)  \rightarrow\widehat{w}\left(
\phi\right)  $ in $S_{w,0}^{\prime}$ i.e. $\left[  \widehat{w}\left(  \phi
_{k}\right)  ,\psi\right]  \rightarrow\left[  \widehat{w}\left(  \phi\right)
,\psi\right]  $ for all $\psi\in S_{w,0}$.

\textbf{So the Fourier transform of a weight function can be regarded as a
continuous mapping from} $S_{w,0}$ $to$ $S_{w,0}^{\prime}$.
\end{remark}

???

\begin{remark}
\textbf{Characterize} $S_{\mathcal{A}}^{\prime}$ HARD! Necessary? ??

$S_{\mathcal{A}}^{\prime}\left(  \mathbb{R}^{d}\right)  $ contains all
distributions the form $\mathcal{W}g$ where $g\left(  x\right)  =f\left(
x^{\prime}\right)  G\left(  x^{\prime\prime}\right)  $, $f\in\overset{\vee
}{S}_{w,0}\cup G\ast S$, $G$ is the basis function in $\mathbb{R}%
^{d^{\prime\prime}}$ and $d^{\prime\prime}\geq1$.
\end{remark}

??

\begin{remark}
Also contains all functions $g\left(  \pi x\right)  $ where $\pi$ is any permutation.

Can replace $G$ by a function in $W_{G}$.

To characterize tempered distributions with support at isolated points a
Taylor series expansion was used. What about in two dimensions for the case
$\mathcal{A}=\left\{  x:x_{1}=0\right\}  $?

Perhaps $S_{\mathcal{A}}^{\prime}=\sum\limits_{k=1}^{M}\left(  D_{1}^{n_{k}%
}\delta_{1}\right)  \left(  x_{1}\right)  u_{k}\left(  x_{2}\right)  $ with
all $u_{k}\in S^{\prime}$?

Use the translation operator.

What about $\widehat{S_{\mathcal{A}}^{\prime}}$? $\mathcal{A}$ is the union of
translated hyperplanes through the origin.

$S_{\tau_{a}\mathcal{A}}^{\prime}=\tau_{-a}S_{\mathcal{A}}^{\prime}$?

But the infinite sum is a problem? Sequence characterization in Arfken
\cite{Arfken70}?
\end{remark}

\subsubsection{\protect\underline{The tensor product central difference weight
functions}}

We will now characterize the space of functions $S_{w,0}$ for the tensor
product central difference weight functions (Section
\ref{Sect_wt_fn_central_diff}) subject to the constraint \ref{a959}, as well
as showing that $w\in W_{S;0}$ (Definition \ref{Def_Sw2_and_fnal(Sw2)}).

\begin{theorem}
\label{Thm_Sw2.eq.Sx,l_minus_n}Suppose $w$ is a tensor product central
difference weight function on $\mathbb{R}^{d}$ satisfying \ref{a959} for
parameters $n$ and $l$.

Then $S_{w,0}=S_{\otimes,l-n}\left(  \mathbf{0}\right)  $. Further, $w\in
W_{S;0}$ and%
\begin{equation}
\int w\left\vert \phi\right\vert ^{2}\leq2^{d}\left(  k_{1}^{\prime}%
+\frac{c_{1}^{\prime}}{\left(  \left(  l-n\right)  !\right)  ^{2}}\right)
^{d}\max_{0\leq\alpha\leq1}\left\Vert x^{\left(  n+1\right)  \left(
\mathbf{1}-\alpha\right)  }D^{\left(  l-n\right)  \alpha}\phi\left(  x\right)
\right\Vert _{\infty}^{2},\quad\phi\in S_{w,0},\label{a1.27}%
\end{equation}

where $c_{1}^{\prime},k_{1}^{\prime}$ are defined by \ref{a971} and \ref{a972}.
\end{theorem}

\begin{proof}
For $a.<b\in\mathbb{R}^{d}$ let $R\left[  a,b\right]  =\left\{  x:a\leq x\leq
b\right\}  $ denote the closed rectangle with extreme corners $a$ and $b$.

Denote the univariate weight function by $w_{1}$ and set $m=l-n$. Since $w$
satisfies \ref{a959} we can use Corollary \ref{Cor_cdiffwt_bnd_on_wt_fn}.
Specifically, when $r=1$ inequality \ref{a972} becomes%
\begin{equation}
\frac{c_{1}}{s^{2m}}\leq w_{1}\left(  s\right)  \leq\frac{c_{1}^{\prime}%
}{s^{2m}},\quad\left\vert s\right\vert \leq1,\mathcal{6}\label{a1.10}%
\end{equation}

for some constants $c_{1}^{\prime}>c_{1}>0$, and inequality \ref{a971} becomes%
\begin{equation}
k_{1}s^{2n}\leq w_{1}\left(  s\right)  \leq k_{1}^{\prime}s^{2n}%
,\quad\left\vert s\right\vert \geq1,\label{a1.22}%
\end{equation}

for some constants $k_{1}^{\prime}>k_{1}>0$. The lower bound of \ref{a1.10}
implies%
\[
\left(  c_{1}\right)  ^{d}\int_{\mathbf{0}}^{\mathbf{1}}\frac{\left\vert
\phi\left(  \xi\right)  \right\vert ^{2}}{\xi^{2m\mathbf{1}}}d\xi
<\int_{\mathbb{R}^{d}}w\left(  \xi\right)  \left\vert \phi\left(  \xi\right)
\right\vert ^{2}d\xi<\infty,
\]

and Lemma \ref{Lem_2} implies $\phi\in S_{\otimes,m}\left(  \mathbf{0}\right)
=S_{\otimes,l-n}\left(  \mathbf{0}\right)  $ and we have demonstrated that
$S_{w,0}\subset S_{\otimes,l-n}\left(  \mathbf{0}\right)  $.

Now suppose that $\phi\in S_{\otimes,m}\left(  \mathbf{0}\right)  $. We will
require the following covering of $\mathbb{R}^{d}$:%
\[
\mathbb{R}^{d}=\bigcup\limits_{\mathbf{0}\leq\alpha\leq\mathbf{1}}A_{\alpha
}=\bigcup\limits_{\mathbf{0}\leq\alpha\leq\mathbf{1}}\left(  I\left(
\alpha_{1}\right)  \times I\left(  \alpha_{2}\right)  \times\ldots\times
I\left(  \alpha_{d}\right)  \right)  ,
\]

where%
\[
I\left(  0\right)  =\left[  -1,1\right]  ,\text{\quad}I\left(  1\right)
=\left[  -\infty,-1\right]  \cup\left[  1,\infty\right]  .
\]

Hence%
\begin{align*}
\int\limits_{\mathbb{R}^{d}}w\left\vert \phi\right\vert ^{2}  & =\sum
_{0\leq\alpha\leq1}\int\limits_{A_{\alpha}}w\left\vert \phi\right\vert ^{2}\\
& =\sum_{0\leq\alpha\leq1}\int\limits_{I\left(  \alpha_{d}\right)  }\ldots
\int\limits_{I\left(  \alpha_{2}\right)  }\int\limits_{I\left(  \alpha
_{1}\right)  }w\left(  x\right)  \left\vert \phi\left(  x\right)  \right\vert
^{2}dx\\
& =\sum_{0\leq\alpha\leq1}\int\limits_{I\left(  \alpha_{d}\right)  }%
w_{1}\left(  x_{d}\right)  \ldots\int\limits_{I\left(  \alpha_{2}\right)
}w_{1}\left(  x_{2}\right)  \int\limits_{I\left(  \alpha_{1}\right)  }%
w_{1}\left(  x_{1}\right)  \left\vert \phi\left(  x\right)  \right\vert
^{2}dx,
\end{align*}

For each $\alpha$ let $\pi_{\alpha}$ be a permutation such that $\pi_{\alpha
}\alpha=\left(  \mathbf{1}^{\prime},\mathbf{0}^{\prime\prime}\right)  $ and
set $k_{\beta}=\left\vert \mathbf{1}^{\prime}\right\vert $. Then $\beta
=\pi_{\alpha}\alpha$ satisfies $A_{\beta}=I\left(  0\right)  ^{k_{\beta}%
}\times I\left(  1\right)  ^{d-k_{\beta}}$ and for each $\alpha$ we apply
$\pi_{\alpha}$ to the variables of integration so that all the partial
integrations over $I\left(  0\right)  $ are carried out first i.e. let
$y=\pi_{\alpha}x$ so that $x=\pi_{\alpha}y$ and the Jacobian is $\pm1$. In the
sequel we will make use of the notation $\phi_{\alpha}=\phi\circ\pi_{\alpha}$,
$y^{\prime}=\left(  y_{1},\ldots,y_{k_{\alpha}}\right)  $, $y^{\prime\prime
}=\left(  y_{k_{\alpha}+1},\ldots,y_{d}\right)  $, $w^{\prime}\left(
y^{\prime}\right)  =w_{1}\left(  x_{1}\right)  \ldots w_{1}\left(
x_{k_{\alpha}}\right)  $, $w^{\prime\prime}\left(  y^{\prime\prime}\right)
=w_{1}\left(  x_{k_{\alpha}+1}\right)  \ldots w_{1}\left(  x_{d}\right)  $.

Now%
\begin{align*}
\int\limits_{\mathbb{R}^{d}}w\left\vert \phi\right\vert ^{2}  & =\sum
_{0\leq\alpha\leq1}\int\limits_{I\left(  \alpha_{d}\right)  }w_{1}\left(
x_{d}\right)  \ldots\int\limits_{I\left(  \alpha_{2}\right)  }w_{1}\left(
x_{2}\right)  \int\limits_{I\left(  \alpha_{1}\right)  }w_{1}\left(
x_{1}\right)  \left\vert \phi\left(  x\right)  \right\vert ^{2}dx\\
& =\sum_{0\leq\alpha\leq1}\int\limits_{I\left(  1\right)  ^{d-k_{\alpha}}%
}w^{\prime\prime}\left(  y^{\prime\prime}\right)  \int\limits_{I\left(
0\right)  ^{k_{\alpha}}}w^{\prime}\left(  y^{\prime}\right)  \left\vert
\phi_{\alpha}\left(  y\right)  \right\vert ^{2}dy^{\prime}dy^{\prime\prime}\\
& \leq\sum_{0\leq\alpha\leq1}\left(  k_{1}^{\prime}\right)  ^{d-k_{\alpha}%
}\left(  c_{1}^{\prime}\right)  ^{k}\int\limits_{I\left(  1\right)
^{d-k_{\alpha}}}\int\limits_{I\left(  0\right)  ^{k_{\alpha}}}y^{2n\left(
\mathbf{0}^{\prime},\mathbf{1}^{\prime\prime}\right)  }\frac{\left\vert
\phi_{\alpha}\left(  y\right)  \right\vert ^{2}}{y^{2m\left(  \mathbf{1}%
^{\prime},\mathbf{0}^{\prime\prime}\right)  }}dy^{\prime}dy^{\prime\prime}.
\end{align*}

But $\phi\in S_{\otimes,m}\left(  \mathbf{0}\right)  $ implies $\phi_{\alpha
}\in S_{\otimes,m}\left(  \mathbf{0}\right)  $ and by \ref{a1.24} when
$\delta=\left(  \mathbf{1}^{\prime},\mathbf{0}^{\prime\prime}\right)  $ there
exists a vector of functions $t^{\prime}=\left(  t_{i}\right)  _{i=1}%
^{k_{\alpha}}$ such that $\mathbf{0}^{\prime}<t^{\prime}<\mathbf{1}^{\prime}$
and%
\[
\phi_{\alpha}\left(  y\right)  =\frac{1}{\left(  m!\right)  ^{k_{\alpha}}%
}y^{m\left(  \mathbf{1}^{\prime},\mathbf{0}^{\prime\prime}\right)  }\left(
D^{m\left(  \mathbf{1}^{\prime},\mathbf{0}^{\prime\prime}\right)  }%
\phi_{\alpha}\right)  \left(  t^{\prime}\left(  y\right)  .y^{\prime
},y^{\prime\prime}\right)  ,\quad y\in\mathbb{R}^{d},
\]

where $.$ denotes the component-wise product of vectors. Thus%
\begin{align*}
&  \int_{\mathbb{R}^{d}}w\left\vert \phi\right\vert ^{2}\\
&  \leq\sum_{0\leq\alpha\leq1}\frac{\left(  k_{1}^{\prime}\right)
^{d-k_{\alpha}}\left(  c_{1}^{\prime}\right)  ^{k_{\alpha}}}{\left(
m!\right)  ^{2k_{\alpha}}}\int\limits_{I\left(  1\right)  ^{d-k_{\alpha}}}%
\int\limits_{I\left(  0\right)  ^{k_{\alpha}}}y^{2n\left(  \mathbf{0}^{\prime
},\mathbf{1}^{\prime\prime}\right)  }\left\vert \left(  D^{m\left(
\mathbf{1}^{\prime},\mathbf{0}^{\prime\prime}\right)  }\phi_{\alpha}\right)
\left(  t^{\prime}\left(  y\right)  .y^{\prime},y^{\prime\prime}\right)
\right\vert ^{2}dy^{\prime}dy^{\prime\prime}\\
&  \leq\left(  \sum_{0\leq\alpha\leq1}\frac{\left(  k_{1}^{\prime}\right)
^{d-k_{\alpha}}\left(  c_{1}^{\prime}\right)  ^{k_{\alpha}}}{\left(
m!\right)  ^{2k_{\alpha}}}\right)  \max_{0\leq\alpha\leq1}\int%
\limits_{I\left(  1\right)  ^{d-k_{\alpha}}}\int\limits_{I\left(  0\right)
^{k_{\alpha}}}y^{2n\left(  \mathbf{0}^{\prime},\mathbf{1}^{\prime\prime
}\right)  }\left\vert \left(  D^{m\left(  \mathbf{1}^{\prime},\mathbf{0}%
^{\prime\prime}\right)  }\phi_{\alpha}\right)  \left(  t^{\prime}\left(
y\right)  .y^{\prime},y^{\prime\prime}\right)  \right\vert ^{2}dy^{\prime
}dy^{\prime\prime}.
\end{align*}

and regarding the leading constant:%
\begin{align*}
\sum_{0\leq\alpha\leq1}\frac{\left(  k_{1}^{\prime}\right)  ^{d-k_{\alpha}%
}\left(  c_{1}^{\prime}\right)  ^{k_{\alpha}}}{\left(  m!\right)
^{2k_{\alpha}}}=\sum_{0\leq\alpha\leq1}\left(  k_{1}^{\prime}\right)
^{d-k_{\alpha}}\left(  \frac{c_{1}^{\prime}}{\left(  m!\right)  ^{2}}\right)
^{k} &  =\sum_{k=1}^{d}\left(  \sum_{k_{\alpha}=k}1\right)  \left(
k_{1}^{\prime}\right)  ^{d-k_{\alpha}}\left(  \frac{c_{1}^{\prime}}{\left(
m!\right)  ^{2}}\right)  ^{k_{\alpha}}\\
&  =\sum_{k_{\alpha}=1}^{d}\binom{d}{k}\left(  k_{1}^{\prime}\right)
^{d-k_{\alpha}}\left(  \frac{c_{1}^{\prime}}{\left(  m!\right)  ^{2}}\right)
^{k_{\alpha}}\\
&  =\left(  k_{1}^{\prime}+\frac{c_{1}^{\prime}}{\left(  m!\right)  ^{2}%
}\right)  ^{d},
\end{align*}

so that%
\[
\int\limits_{\mathbb{R}^{d}}w\left\vert \phi\right\vert ^{2}\leq\left(
k_{1}^{\prime}+\frac{c_{1}^{\prime}}{\left(  m!\right)  ^{2}}\right)  ^{d}%
\max_{0\leq\alpha\leq1}\int\limits_{I\left(  1\right)  ^{d-k_{\alpha}}}%
\int\limits_{I\left(  0\right)  ^{k_{\alpha}}}y^{2n\left(  \mathbf{0}^{\prime
},\mathbf{1}^{\prime\prime}\right)  }\left\vert \left(  D^{m\left(
\mathbf{1}^{\prime},\mathbf{0}^{\prime\prime}\right)  }\phi_{\alpha}\right)
\left(  t^{\prime}\left(  y\right)  .y^{\prime},y^{\prime\prime}\right)
\right\vert ^{2}dy^{\prime}dy^{\prime\prime},
\]

so that%
\begin{align*}
&  \int_{\mathbb{R}^{d}}w\left\vert \phi\right\vert ^{2}\\
&  \leq\left(  k_{1}^{\prime}+\frac{c_{1}^{\prime}}{\left(  m!\right)  ^{2}%
}\right)  ^{d}\max_{0\leq\alpha\leq1}\int\limits_{I\left(  1\right)
^{d-k_{\alpha}}}\int\limits_{I\left(  0\right)  ^{k_{\alpha}}}y^{2n\left(
\mathbf{0}^{\prime},\mathbf{1}^{\prime\prime}\right)  }\left\vert \left(
D^{m\left(  \mathbf{1}^{\prime},\mathbf{0}^{\prime\prime}\right)  }%
\phi_{\alpha}\right)  \left(  t^{\prime}\left(  y\right)  .y^{\prime
},y^{\prime\prime}\right)  \right\vert ^{2}dy^{\prime}dy^{\prime\prime}\\
&  =\left(  k_{1}^{\prime}+\frac{c_{1}^{\prime}}{\left(  m!\right)  ^{2}%
}\right)  ^{d}\max_{0\leq\alpha\leq1}\int\limits_{I\left(  1\right)
^{d-k_{\alpha}}}\int\limits_{I\left(  0\right)  ^{k_{\alpha}}}\frac
{1}{y^{2\mathbf{1}^{\prime\prime}}}y^{\left(  2n+2\right)  \left(
\mathbf{0}^{\prime},\mathbf{1}^{\prime\prime}\right)  }\left\vert \left(
D^{m\left(  \mathbf{1}^{\prime},\mathbf{0}^{\prime\prime}\right)  }%
\phi_{\alpha}\right)  \left(  t^{\prime}\left(  y\right)  .y^{\prime
},y^{\prime\prime}\right)  \right\vert ^{2}dy^{\prime}dy^{\prime\prime}\\
&  \leq\left(  k_{1}^{\prime}+\frac{c_{1}^{\prime}}{\left(  m!\right)  ^{2}%
}\right)  ^{d}\max_{0\leq\alpha\leq1}\left\Vert y^{\left(  2n+2\right)
\left(  \mathbf{0}^{\prime},\mathbf{1}^{\prime\prime}\right)  }\left\vert
\left(  D^{m\left(  \mathbf{1}^{\prime},\mathbf{0}^{\prime\prime}\right)
}\phi_{\alpha}\right)  \left(  t^{\prime}\left(  y\right)  .y^{\prime
},y^{\prime\prime}\right)  \right\vert ^{2}\right\Vert _{\infty}%
\int\limits_{I\left(  1\right)  ^{d-k_{\alpha}}}\int\limits_{I\left(
0\right)  ^{k_{\alpha}}}\frac{dy^{\prime}dy^{\prime\prime}}{y^{2\mathbf{1}%
^{\prime\prime}}}\\
&  \leq\left(  k_{1}^{\prime}+\frac{c_{1}^{\prime}}{\left(  m!\right)  ^{2}%
}\right)  ^{d}\max_{0\leq\alpha\leq1}\left\Vert y^{\left(  n+1\right)  \left(
\mathbf{0}^{\prime},\mathbf{1}^{\prime\prime}\right)  }D^{m\left(
\mathbf{1}^{\prime},\mathbf{0}^{\prime\prime}\right)  }\phi_{\alpha}\left(
y\right)  \right\Vert _{\infty}^{2}\int\limits_{I\left(  1\right)
^{d-k_{\alpha}}}\int\limits_{I\left(  0\right)  ^{k_{\alpha}}}\frac
{dy^{\prime}dy^{\prime\prime}}{y^{2\mathbf{1}^{\prime\prime}}}\\
&  =\left(  k_{1}^{\prime}+\frac{c_{1}^{\prime}}{\left(  m!\right)  ^{2}%
}\right)  ^{d}\max_{0\leq\alpha\leq1}2^{k_{\alpha}}\left\Vert y^{\left(
n+1\right)  \left(  \mathbf{0}^{\prime},\mathbf{1}^{\prime\prime}\right)
}D^{m\left(  \mathbf{1}^{\prime},\mathbf{0}^{\prime\prime}\right)  }%
\phi_{\alpha}\left(  y\right)  \right\Vert _{\infty}^{2}\\
&  \leq2^{d}\left(  k_{1}^{\prime}+\frac{c_{1}^{\prime}}{\left(  m!\right)
^{2}}\right)  ^{d}\max_{0\leq\alpha\leq1}\left\Vert y^{\left(  n+1\right)
\left(  \mathbf{0}^{\prime},\mathbf{1}^{\prime\prime}\right)  }D^{m\left(
\mathbf{1}^{\prime},\mathbf{0}^{\prime\prime}\right)  }\phi_{\alpha}\left(
y\right)  \right\Vert _{\infty}^{2}.
\end{align*}

But if $y=\pi_{\alpha}x$ then%
\[
y^{\left(  \mathbf{0}^{\prime},\mathbf{1}^{\prime\prime}\right)  }=\left(
\pi_{\alpha}x\right)  ^{\left(  \mathbf{0}^{\prime},\mathbf{1}^{\prime\prime
}\right)  }=x^{\pi_{\alpha}\left(  \mathbf{0}^{\prime},\mathbf{1}%
^{\prime\prime}\right)  }=x^{\pi_{\alpha}\left(  \mathbf{1}-\left(
\mathbf{1}^{\prime},\mathbf{0}^{\prime\prime}\right)  \right)  }%
=x^{\mathbf{1}-\alpha},
\]

so that%
\[
y^{\left(  n+1\right)  \left(  \mathbf{0}^{\prime},\mathbf{1}^{\prime\prime
}\right)  }=x^{\left(  n+1\right)  \left(  \mathbf{1}-\alpha\right)  }.
\]

and also%
\[
D^{m\left(  \mathbf{1}^{\prime},\mathbf{0}^{\prime\prime}\right)  }%
\phi_{\alpha}\left(  y\right)  =\left(  D^{m\left(  \mathbf{1}^{\prime
},\mathbf{0}^{\prime\prime}\right)  }\left(  \phi\circ\pi_{\alpha}\right)
\right)  \left(  \pi_{\alpha}x\right)  =D^{m\alpha}\phi\left(  x\right)  .
\]

Thus%
\begin{align*}
\int_{\mathbb{R}^{d}}w\left\vert \phi\right\vert ^{2}  & \leq2^{d}\left(
k_{1}^{\prime}+\frac{c_{1}^{\prime}}{\left(  m!\right)  ^{2}}\right)  ^{d}%
\max_{0\leq\alpha\leq1}\left\Vert y^{\left(  n+1\right)  \left(
\mathbf{0}^{\prime},\mathbf{1}^{\prime\prime}\right)  }D^{m\left(
\mathbf{1}^{\prime},\mathbf{0}^{\prime\prime}\right)  }\phi_{\alpha}\left(
y\right)  \right\Vert _{\infty}^{2}\\
& =2^{d}\left(  k_{1}^{\prime}+\frac{c_{1}^{\prime}}{\left(  m!\right)  ^{2}%
}\right)  ^{d}\max_{0\leq\alpha\leq1}\left\Vert x^{\left(  n+1\right)  \left(
\mathbf{1}-\alpha\right)  }D^{m\alpha}\phi\left(  x\right)  \right\Vert
_{\infty}^{2}\\
& =2^{d}\left(  k_{1}^{\prime}+\frac{c_{1}^{\prime}}{\left(  \left(
l-n\right)  !\right)  ^{2}}\right)  ^{d}\max_{0\leq\alpha\leq1}\left\Vert
x^{\left(  n+1\right)  \left(  \mathbf{1}-\alpha\right)  }D^{\left(
l-n\right)  \alpha}\phi\left(  x\right)  \right\Vert _{\infty}^{2},
\end{align*}

which proves \ref{a1.27}. We can also conclude that $\phi\in S_{w,0}$ implies
$\phi\in S_{\otimes,l-n}\left(  \mathbf{0}\right)  $ and thus $S_{w,0}%
=S_{\otimes,\left(  l-n\right)  \mathbf{1}}\left(  \mathbf{0}\right)  $.

\textbf{Part 6} From part 1,%
\begin{align*}
\int &  \leq\left\Vert \overset{\vee}{\phi}\right\Vert _{w,0}^{2}\\
&  \leq\left(  \frac{\pi}{2}\right)  ^{2\left\vert \lambda\right\vert }\left(
\left(  \frac{2^{\left\vert \lambda-\nu\right\vert }}{\left(  2\left(
\lambda-\nu\right)  -1\right)  !!}\right)  ^{2}\left\Vert \widehat{x^{\lambda
-\nu}\overset{\vee}{\phi}}\right\Vert _{2}^{2}+\left(  \frac{2^{\left\vert
\lambda\right\vert }}{\left(  2\lambda-1\right)  !!}\right)  ^{2}\left\Vert
\widehat{x^{\lambda}D^{\nu}\overset{\vee}{\phi}}\right\Vert _{2}^{2}\right) \\
&  =\left(  \frac{\pi}{2}\right)  ^{2\left\vert \lambda\right\vert }\left(
\left(  \frac{2^{\left\vert \lambda-\nu\right\vert }}{\left(  2\left(
\lambda-\nu\right)  -1\right)  !!}\right)  ^{2}\left\Vert x^{\lambda-\nu
}\overset{\vee}{\phi}\right\Vert _{2}^{2}+\left(  \frac{2^{\left\vert
\lambda\right\vert }}{\left(  2\lambda-1\right)  !!}\right)  ^{2}\left\Vert
x^{\lambda}D^{\nu}\overset{\vee}{\phi}\right\Vert _{2}^{2}\right)
\end{align*}

\end{proof}

\begin{corollary}
?? $S_{w,0}=\frac{1}{w}S$.
\end{corollary}

\begin{proof}
??
\end{proof}

\begin{remark}
??? \textbf{These calculations relate to trying to prove:} $\int w\left\vert
\chi\eta\right\vert ^{2}\leq C\left(  \int w\left\vert \chi\right\vert
^{2}\right)  \left(  \int w\left\vert \eta\right\vert ^{2}\right)  $ when
$\chi,\eta\in S_{w,0}=S_{\otimes,\left(  l-n\right)  \mathbf{1}}\left(
\mathbf{0}\right)  $.

Suppose the conditions of Theorem \ref{Thm_cdiffwt_3} are satisfied. Then by
Theorem \ref{Thm_cdiffwt_3} if $m=l-n$ then%
\[
\frac{c_{1}}{s^{2m}}\leq w_{1}\left(  s\right)  \leq\frac{c_{1}^{\prime}%
}{s^{2m}},\quad\left\vert s\right\vert \leq1,\mathcal{6}%
\]

for some constants $c_{1}^{\prime}>c_{1}>0$, and%
\[
k_{1}s^{2n}\leq w_{1}\left(  s\right)  \leq k_{1}^{\prime}s^{2n}%
,\quad\left\vert s\right\vert \geq1,
\]

for some constants $k_{1}^{\prime}>k_{1}>0$.

From \ref{a1.23},%
\[
\left\vert \phi\left(  \xi\right)  \right\vert \leq\frac{\left\vert
\xi^{\alpha}\right\vert }{\alpha!}\max_{t\in\left[  0,1\right]  }\left\vert
\left(  D^{\left(  l-n\right)  \mathbf{1}}\phi\right)  \left(  t\xi\right)
\right\vert .
\]

\[
\sqrt{w\left(  \xi\right)  }\leq\frac{c_{1}^{\prime}}{\xi_{+}^{m\mathbf{1}}%
},\quad\xi_{+}\leq\mathbf{1}.
\]

?? Apply above to $\sqrt{w_{k_{1}}\left(  \xi_{k_{1}}\right)  \ldots w_{k_{j}%
}\left(  \xi_{k_{j}}\right)  }\left\vert \eta\left(  \xi\right)  \right\vert $
where $\left\vert \xi_{k_{i}}\right\vert \leq1$. Then use%
\[
k_{1}\xi_{j}^{2n}\leq w_{j}\left(  \xi_{j}\right)  \leq k_{1}^{\prime}\xi
_{j}^{2n},\quad\left\vert \xi_{j}\right\vert \geq1.
\]

USE $\left(  \beta\left(  \xi\right)  \right)  _{k}=\left\{
\begin{array}
[c]{ll}%
l, & \left\vert \xi_{k}\right\vert <1,\\
n, & \left\vert \xi_{k}\right\vert \geq1,
\end{array}
\right\}  $ and $\mathcal{R}\left(  \xi\right)  =$??.%
\[
\sum??+\int w\left\vert \eta\chi\right\vert ^{2}%
\]

Perhaps if $m=0$ i.e. if $l=n$, the inequality holds?
\end{remark}

??

\begin{remark}
??? Try proving that if $w\left(  x\right)  \leq c\left(  1+\left\vert
x\right\vert \right)  ^{m}$ for some $c>0$ and integer $m>d$, and if $w$ has
sufficient smoothness (to allow a suitable local Taylor series expansion),
then $X_{w}^{0}$ is a ?? \textbf{NOT LIKELY}! \textbf{convolution algebra} ??.
Exponential weight functions are suitable candidates.

Try proving this by partitioning $\mathbb{R}^{d}$ into unit cubes and then
proceeding in a manner analogous to the expansion \ref{a1.70} in the proof of
part 1 of Theorem \ref{Thm_eta^2/sin(eta)^2.|phi|^2_dim_gt_1} i.e. we have the
full volume integrals
\begin{align*}
\int w\left\vert \phi\right\vert ^{2}  & =\sum\limits_{\alpha\in\mathbb{Z}%
^{d}}\int_{\alpha-\frac{1}{2}}^{\alpha+\frac{1}{2}}w\left(  \eta\right)
\left\vert \phi\left(  \eta\right)  \right\vert ^{2}d\eta\\
& =\sum\limits_{\alpha\in\mathbb{Z}^{d}}\int_{-\frac{\mathbf{1}}{2}}%
^{\frac{\mathbf{1}}{2}}w\left(  \alpha+\xi\right)  \left\vert \phi\left(
\alpha+\xi\right)  \right\vert ^{2}d\xi\\
& =\int_{-\frac{\mathbf{1}}{2},}^{\frac{\mathbf{1}}{2}}w\left\vert
\phi\right\vert ^{2}+\sum\limits_{\alpha\in\mathbb{Z}^{d}\setminus\mathbf{0}%
}\int_{-\frac{\mathbf{1}}{2}}^{\frac{\mathbf{1}}{2}}w\left(  \alpha
+\xi\right)  \left\vert \phi\left(  \alpha+\xi\right)  \right\vert ^{2}d\xi.
\end{align*}

Assume that $\forall\alpha$, $0<C_{1}\leq\left\vert \xi^{2\lambda_{\alpha}%
}\frac{w\left(  \alpha+\xi\right)  }{\left(  a+\xi\right)  ^{2\nu_{\alpha}}%
}\right\vert \leq C_{2}$ on $\overline{I}_{0}$, $\max\limits_{\alpha}%
\lambda_{\alpha}<\infty$ and $\max\limits_{\alpha}\nu_{\alpha}<\infty$. Now
write%
\[
\int w\left\vert \phi\right\vert ^{2}=\int_{-\frac{\mathbf{1}}{2}}%
^{\frac{\mathbf{1}}{2}}w\left\vert \phi\right\vert ^{2}+\sum\limits_{\alpha
\in\mathbb{Z}^{d}\setminus\mathbf{0}}\int_{-\frac{\mathbf{1}}{2}}%
^{\frac{\mathbf{1}}{2}}\xi^{2\lambda_{\alpha}}\left(  \frac{w}{x^{\nu_{\alpha
}}}\right)  \left(  \alpha+\xi\right)  \left\vert \frac{\left(  x^{\nu
_{\alpha}}\phi\right)  \left(  \alpha+\xi\right)  }{\xi^{\lambda_{\alpha}}%
}\right\vert ^{2}d\xi.
\]

Thus $\int w\left\vert \phi\right\vert ^{2}<\infty$ implies%
\[
\int_{-\frac{\mathbf{1}}{2}}^{\frac{\mathbf{1}}{2}}w\left\vert \phi\right\vert
^{2}<\infty,\quad\int_{-\frac{\mathbf{1}}{2}}^{\frac{\mathbf{1}}{2}}\left\vert
\frac{\left(  x^{\nu_{\alpha}}\phi\right)  \left(  \alpha+\xi\right)  }%
{\xi^{\lambda_{\alpha}}}\right\vert ^{2}d\xi<\infty\quad\forall\alpha
\neq\mathbf{0}.
\]

?? ETC.
\end{remark}

\subsubsection{\protect\underline{?? What about the central difference weight
functions of Definition \ref{Def_multi_centdiff_wt_fn}?}}

???

\subsection{The operator $\mathcal{B}:X_{w}^{0}\otimes X_{1/w}^{0}\rightarrow
C_{B}^{\left(  0\right)  }$\label{SbSect_ops_B_and_V}}

In this section we will assume that $w\in S_{w,0}^{\dag}$ and in Section
\ref{Sect_generalize1_tildXo1/w_T} we will construct a generalization which
includes property \ref{a1.046}.

In this section we will introduce the operators $\mathcal{B}:X_{w}^{0}\otimes
X_{1/w}^{0}\rightarrow C_{B}^{\left(  0\right)  }$, $\mathcal{V}:X_{1/w}%
^{0}\rightarrow X_{w}^{0}$, $\mathcal{W}:X_{w}^{0}\rightarrow X_{1/w}^{0}$ and
$\Phi:X_{1/w}^{0}\rightarrow\left(  X_{w}^{0}\right)  ^{\prime}$.

where $\left(  X_{w}^{0}\right)  ^{\prime}$ denotes the bounded linear
functionals on $X_{w}^{0}$. In summary, if $G$ is the basis function and $u\in
X_{w}^{0}$, $v\in X_{1/w}^{0}$ then:%
\begin{align*}
\mathcal{B}\left(  u,v\right)   & =\left(  \widehat{u}v_{F}\right)  ^{\vee
}=\left(  \mathcal{I}u\right)  \ast\mathcal{M}v,\\
\mathcal{V}v  & =\mathcal{B}\left(  G,v\right)  =\mathcal{JM}v,\\
\mathcal{W}u  & =\mathcal{LI}u,\\
\left(  \Phi v\right)  \left(  u\right)   & =\left(  u,\mathcal{V}v\right)
_{w,0},
\end{align*}

and $\Phi$ characterizes the bounded linear functionals on $X_{w}^{0}$. The
operator $\mathcal{B}$ is actually not needed to define $\Phi$. The bilinear
mapping $\mathcal{B}$ can be regarded as a convolution in the sense of
\ref{1.20} (in Appendix) because%
\[
\mathcal{B}\left(  u,\phi\right)  =u\ast\phi=\left(  2\pi\right)  ^{-\frac
{d}{2}}\left[  u_{y},\phi\left(  \cdot-y\right)  \right]  ,\quad u\in
X_{w}^{0},\text{ }\phi\in S.
\]

We start by defining the mapping $\mathcal{B}$:

\begin{definition}
\label{Def_map_B}\textbf{The mapping }$\mathcal{B}$

If $\left(  u,v\right)  \in X_{w}^{0}\otimes X_{1/w}^{0}$ define the bilinear
mapping $\mathcal{B}:X_{w}^{0}\otimes X_{1/w}^{0}\rightarrow C_{B}^{\left(
0\right)  }$ by:%
\begin{equation}
\mathcal{B}\left(  u,v\right)  =\left(  \widehat{u}v_{F}\right)  ^{\vee},\quad
u\in X_{w}^{0},\text{ }v\in X_{1/w}^{0}.\label{a1.20}%
\end{equation}

Noting that the product of two measurable functions is a measurable function
definition \ref{a1.20} makes sense since from the definitions of $X_{1/w}^{0}
$ and $X_{w}^{0}$, $v_{F}\in L_{loc}^{1}\left(  \mathbb{R}^{d}\setminus
\mathcal{A}\right)  $, $\frac{v_{F}}{\sqrt{w}}\in L^{2}$, $\widehat{u}\in
L_{loc}^{1}$ and $\sqrt{w}\widehat{u}\in L^{2}$ so that $\int\left\vert
\widehat{u}v_{F}\right\vert \leq\left\Vert u\right\Vert _{w,0}\left\vert
v\right\vert _{1/w,0}<\infty$. Thus $\widehat{u}v_{F}\in L^{1}$ and hence
$\left(  \widehat{u}v_{F}\right)  ^{\vee}\in C_{B}^{\left(  0\right)  }$.
\end{definition}

\begin{theorem}
\label{Thm_properties_B}\textbf{Properties of the operator }$\mathcal{B}$

\begin{enumerate}
\item $\mathcal{B}:X_{w}^{0}\otimes X_{1/w}^{0}\rightarrow C_{B}^{\left(
0\right)  }$ is a continuous bilinear mapping when $C_{B}^{\left(  0\right)
}$ is endowed with the supremum norm $\left\Vert \cdot\right\Vert _{\infty}$.
In fact%
\[
\left\Vert \mathcal{B}\left(  u,v\right)  \right\Vert _{\infty}\leq\left(
2\pi\right)  ^{-\frac{d}{2}}\left\Vert u\right\Vert _{w,0}\left\vert
v\right\vert _{1/w,0},\quad u\in X_{w}^{0},\text{ }v\in X_{1/w}^{0}.
\]

\item In the sense of distributions: if $\alpha+\beta=\gamma\geq0$ then
\[
D^{\gamma}\mathcal{B}\left(  u,v\right)  =\left(  \widehat{D^{\alpha}u}\left(
D^{\beta}v\right)  _{F}\right)  ^{\vee},\quad u\in X_{w}^{0},\text{ }v\in
X_{1/w}^{0}.
\]

\item The operator $\mathcal{B}$ commutes with the (distribution) translation
operator $\tau_{a}f=f\left(  \cdot-a\right)  $, $a\in\mathbb{R}^{d}$ in the
sense that%
\[
\tau_{a}\mathcal{B}\left(  u,v\right)  =\mathcal{B}\left(  \tau_{a}u,v\right)
=\mathcal{B}\left(  u,\tau_{a}v\right)  .
\]

\item We have%
\[
\mathcal{B}\left(  u,v\right)  =\left(  \mathcal{I}u\right)  \ast
\mathcal{M}v,\quad u\in X_{w}^{0},\text{ }v\in X_{1/w}^{0},
\]

where $\mathcal{I}$ was introduced in Definition \ref{Def_I_J} and
$\mathcal{M}$ in Definition \ref{Def_op_M}.\medskip

Now suppose $w$ also satisfies property W02 or W03 for parameter $\kappa$.
Then:\medskip

\item $\mathcal{B}$ is a convolution in the sense of \ref{1.20}: $X_{w}%
^{0}\subset S^{\prime}$, $S\subset X_{1/w}^{0}$ and%
\[
\mathcal{B}\left(  u,\phi\right)  =u\ast\phi,\quad u\in X_{w}^{0},\text{ }%
\phi\in S.
\]

Further, $\mathcal{B}\left(  u,\phi\right)  \in C_{B}^{\infty}$ and
$D^{\gamma}\mathcal{B}\left(  u,\phi\right)  =D^{\alpha}u\ast D^{\beta}\phi$
for all $\gamma=\alpha+\beta$.
\end{enumerate}
\end{theorem}

\begin{proof}
\textbf{Part 1} Noting the calculations done in the definition of
$\mathcal{B}$ all that remains to be shown is the continuity of $\mathcal{B}$.
But from \ref{a1.20}
\begin{align*}
\left\vert \mathcal{B}\left(  u,v\right)  \right\vert \leq\left\vert \left(
\widehat{u}v_{F}\right)  ^{\vee}\right\vert =\left(  2\pi\right)  ^{-\frac
{d}{2}}\left\vert \int e^{-ix\xi}\widehat{u}\left(  \xi\right)  v_{F}\left(
\xi\right)  d\xi\right\vert  &  \leq\left(  2\pi\right)  ^{-\frac{d}{2}}%
\int\left\vert \widehat{u}\right\vert \left\vert v_{F}\right\vert \\
&  \leq\left(  2\pi\right)  ^{-\frac{d}{2}}\left\Vert u\right\Vert
_{w,0}\left\vert v\right\vert _{1/w,0}.
\end{align*}
\medskip

\textbf{Part 2}%
\[
D^{\gamma}\mathcal{B}\left(  u,v\right)  =D^{\gamma}\left(  \widehat{u}%
v_{F}\right)  ^{\vee}=\left(  \left(  i\xi\right)  ^{\gamma}\left(
\widehat{u}v_{F}\right)  \right)  ^{\vee}=\left(  \left(  i\xi\right)
^{\alpha}\widehat{u}\left(  i\xi\right)  ^{\beta}v_{F}\right)  ^{\vee}=\left(
\widehat{D^{\alpha}u}\left(  D^{\beta}v\right)  _{F}\right)  ^{\vee}.
\]
\medskip

\textbf{Part 3}%
\[
\tau_{a}\mathcal{B}\left(  u,v\right)  =\tau_{a}\left(  \widehat{u}%
v_{F}\right)  ^{\vee}=\left(  e^{-ia\xi}\widehat{u}v_{F}\right)  ^{\vee
}=\left(  \widehat{\tau_{a}u}v_{F}\right)  ^{\vee}=\mathcal{B}\left(  \tau
_{a}u,v\right)  ,
\]

and%
\[
\tau_{a}\mathcal{B}\left(  u,v\right)  =\tau_{a}\left(  \widehat{u}%
v_{F}\right)  ^{\vee}=\left(  \widehat{u}e^{-ia\xi}v_{F}\right)  ^{\vee
}=\left(  \widehat{u}\text{ }\left(  \tau_{a}v\right)  _{F}\right)  ^{\vee
}=\mathcal{B}\left(  u,\tau_{a}v\right)  .
\]
\medskip

\textbf{Part 4} From the definition of $\mathcal{B}$ and the definition of
$\mathcal{M}$ (Definition \ref{Def_op_M})%
\begin{align*}
\mathcal{B}\left(  u,v\right)  =\left(  \widehat{u}v_{F}\right)  ^{\vee
}=\left(  \sqrt{w}\widehat{u}\frac{v_{F}}{\sqrt{w}}\right)  ^{\vee} &
=\left(  \sqrt{w}\widehat{u}\right)  ^{\vee}\ast\left(  \frac{v_{F}}{\sqrt{w}%
}\right)  ^{\vee}\\
&  =\left(  \mathcal{I}u\right)  \ast\mathcal{M}v.
\end{align*}
\medskip

\textbf{Part 5} From part 4 of Lemma \ref{Lem_Xoinvw}, $S\subset X_{1/w}^{0}$.
So by definition \ref{a1.20} of $\mathcal{B}\left(  u,\phi\right)  $, and the
definition \ref{1.20} of the convolution of a member of $S^{\prime}$ and a
member of $S$, $\mathcal{B}\left(  u,\phi\right)  =\left(  \widehat{u}\phi
_{F}\right)  ^{\vee}=\left(  \widehat{u}\widehat{\phi}\right)  ^{\vee}%
=u\ast\phi\in C_{B}^{\left(  0\right)  }$ so that by \ref{2.40}, $D^{\gamma
}\mathcal{B}\left(  u,\phi\right)  =\left(  D^{\alpha}u\right)  \ast D^{\beta
}\phi$ for all $\gamma$ and $\alpha+\beta=\gamma$. Thus $\mathcal{B}\left(
u,\phi\right)  \in C_{B}^{\infty}$.
\end{proof}

Part 5 of the last theorem shows how $\mathcal{B}$ can be regarded as a
convolution by restricting the domain of the second variable to $S$.

\begin{definition}
\label{Def_Sw2_lin_fnal_0}\textbf{The spaces }$\widehat{S}_{w,0}$,
$\overset{\vee}{S}_{w,0}$ \textbf{and the spaces of continuous functionals}
$S_{w,0}^{\prime}$, $\left(  \widehat{S}_{w,0}\right)  ^{\prime}$, $\left(
\overset{\vee}{S}_{w,0}\right)  ^{\prime}$.

The space $S_{w,0}$ was introduced in Definition \ref{Def_Sw2_and_fnal(Sw2)}.

\begin{enumerate}
\item The spaces $\widehat{S}_{w,0}$ and $\overset{\vee}{S}_{w,0}$ are defined
by
\[
\widehat{S}_{w,0}=\left\{  \widehat{\phi}:\phi\in S_{w,0}\right\}
,\quad\overset{\vee}{S}_{w,0}=\left\{  \overset{\vee}{\phi}:\phi\in
S_{w,0}\right\}  ,
\]

and we observe that from the definition of $S_{w,0}$: $\overset{\vee}{S}%
_{w,0}=S\cap X_{w}^{0}$.\smallskip

$\widehat{S}_{w,0}$ is topologized using the seminorms which define the
topology of $S$.

$\overset{\vee}{S}_{w,0}$ is topologized using the seminorms which define the
topology of $S$.

\item With these topologies the Fourier transform is now a homeomorphism from
$S_{w,0}$ to $\overset{\vee}{S}_{w,0}$ and the inverse-Fourier transform is
now a homeomorphism from $S_{w,0}$ to $\widehat{S}_{w,0}$.

\item Now define $S_{w,0}^{\prime}$, $\left(  \widehat{S}_{w,0}\right)
^{\prime}$ and $\left(  \overset{\vee}{S}_{w,0}\right)  ^{\prime}$ to be the
spaces of continuous, linear functionals on $S_{w,0}$, $\widehat{S}_{w,0} $
and $\overset{\vee}{S}_{w,0}$ respectively.

\item We define the Fourier transform on $S_{w,0}^{\prime}$ and
inverse-Fourier transform on $S_{w,0}^{\prime}$ by%
\begin{align*}
\left[  \widehat{f},\phi\right]   & =\left[  f,\widehat{\phi}\right]  ,\quad
f\in S_{w,0}^{\prime},\text{ }\phi\in\overset{\vee}{S}_{w,0},\\
\left[  \overset{\vee}{f},\phi\right]   & =\left[  f,\overset{\vee}{\phi
}\right]  ,\quad f\in S_{w,0}^{\prime},\text{ }\phi\in\widehat{S}_{w,0},
\end{align*}

so that
\[
\left(  S_{w,0}^{\prime}\right)  ^{\vee}=\left(  \widehat{S}_{w,0}\right)
^{\prime},\quad\left(  S_{w,0}^{\prime}\right)  ^{\wedge}=\left(
\overset{\vee}{S}_{w,0}\right)  ^{\prime}.
\]

These equations imply that the Fourier transform is a homeomorphism from
$S_{w,0}^{\prime}$ to $\left(  \overset{\vee}{S}_{w,0}\right)  ^{\prime}$, and
that the inverse-Fourier transform is a homeomorphism from $S_{w,0}^{\prime}$
to $\left(  \widehat{S}_{w,0}\right)  ^{\prime}$.\medskip

\textbf{Localization}:\medskip

\item Suppose $\Omega\subseteq\mathbb{R}^{d}\setminus\mathcal{A}$ is open and
let $C_{0}^{\infty}\left(  \Omega\right)  =\left\{  \phi\in C_{0}^{\infty
}:\operatorname*{supp}\phi\subset\Omega\right\}  $. Then $C_{0}^{\infty
}\left(  \Omega\right)  \subset S_{w,0}$ since $w$ is continuous on
$\mathbb{R}^{d}\setminus\mathcal{A}$. We now say $u\in S_{w,0}^{\prime}$ is a
member of $L_{loc}^{1}\left(  \Omega\right)  $ if there exists $f\in
L_{loc}^{1}\left(  \Omega\right)  $ such that $\left[  u,\phi\right]  =\int
f\phi$ for all $\phi\in C_{0}^{\infty}\left(  \Omega\right)  $.

\item Suppose $\Omega\subset\mathbb{R}^{d}\setminus\mathcal{A}$ is open and
$A$ is one of the spaces $S_{w,0}^{\prime}$, $\left(  \widehat{S}%
_{w,0}\right)  ^{\prime}$ or $\left(  \overset{\vee}{S}_{w,0}\right)
^{\prime}$.

Then if $u,v\in A$ we say that $u=v$ on $\Omega$ if $\left[  u-v,\phi\right]
=0$ for all $\phi\in C_{0}^{\infty}\left(  \Omega\right)  $.

\item A member $u$ of $S_{w,0}^{\prime}$ has support on $\mathcal{A}$ iff
$\left[  u,\phi\right]  =0$ when $\phi\in C_{0}^{\infty}$ and
$\operatorname*{supp}\phi\subset\mathbb{R}^{d}\setminus\mathcal{A}$.
\end{enumerate}
\end{definition}

We now want to define the convolution of a member of $\overset{\vee}{S}_{w,0}$
and a member of $\left(  \widehat{S}_{w,0}\right)  ^{\prime}$ and we use as a
guide definition \ref{1.20} of the Appendix:
\begin{equation}
\phi\ast f=\left(  \widehat{\phi}\widehat{f}\right)  ^{\vee}=\left(
2\pi\right)  ^{-d/2}\left[  f_{y},\phi\left(  \cdot-y\right)  \right]
,\quad\phi\in S,\text{ }f\in S^{\prime},\label{a2.0}%
\end{equation}

Here $\phi\ast f\in C_{BP}^{\infty}$.

\begin{definition}
\label{Def_Sw2_convol_fnalSw2}\textbf{The convolution }$\overset{\vee
}{S}_{w,0}\ast\left(  \widehat{S}_{w,0}\right)  ^{\prime}\rightarrow\left(
\widehat{S}_{w,0}\right)  ^{\prime}$

If $f\in\left(  \widehat{S}_{w,0}\right)  ^{\prime}$ and $\phi\in
\overset{\vee}{S}_{w,0}$ then $\widehat{f}=S_{w,0}^{\prime}$ and
$\widehat{\phi}\in S$ so that $\widehat{\phi}\widehat{f}\in S_{w,0}^{\prime}$
and $\left(  \widehat{\phi}\widehat{f}\right)  ^{\vee}\in\left(
\widehat{S}_{w,0}\right)  ^{\prime}$. Now by analogy with \ref{a2.0} we define%
\begin{equation}
\phi\ast f=\left(  \widehat{\phi}\widehat{f}\right)  ^{\vee},\quad\phi
\in\overset{\vee}{S}_{w,0},\text{ }f\in\left(  \widehat{S}_{w,0}\right)
^{\prime}.\label{a2.1}%
\end{equation}

\end{definition}

Inspired by \ref{a2.0} we now prove the convolution representation:

\begin{lemma}
\label{Lem_convol_Sw2}\textbf{The convolution of Definition
\ref{Def_Sw2_convol_fnalSw2}} can be written:%
\[
\phi\ast f=\left(  \widehat{\phi}\widehat{f}\right)  ^{\vee}=\left(
2\pi\right)  ^{-\frac{d}{2}}\left[  f_{y},\phi\left(  \cdot-y\right)  \right]
,\quad\phi\in\overset{\vee}{S}_{w,0},\text{ }f\in\left(  \widehat{S}%
_{w,0}\right)  ^{\prime}.
\]

\end{lemma}

\begin{proof}
Now $\phi\in S$ and $f$ can be extended (non-uniquely) to $S$ as, say,
$f^{e}\in S^{\prime}$.

Hence\ $\left(  \widehat{\phi}\widehat{f^{e}}\right)  ^{\vee}=\left(
2\pi\right)  ^{-\frac{d}{2}}\left[  f_{y}^{e},\phi\left(  \cdot-y\right)
\right]  =\left(  2\pi\right)  ^{-\frac{d}{2}}\left[  f_{y},\phi\left(
\cdot-y\right)  \right]  $, since $\phi\left(  \cdot-y\right)  \in
\overset{\vee}{S}_{w,0}$ for all $y$.

Also, if $\psi\in\widehat{S}_{w,0}$ then, $\left[  \left(  \widehat{\phi
}\widehat{f^{e}}\right)  ^{\vee},\psi\right]  =\left[  \widehat{\phi
}\widehat{f^{e}},\overset{\vee}{\psi}\right]  =\left[  \widehat{f^{e}%
},\widehat{\phi}\overset{\vee}{\psi}\right]  =\left[  f^{e},\left(
\widehat{\phi}\overset{\vee}{\psi}\right)  ^{\wedge}\right]  $. But
$\overset{\vee}{\psi}\in S_{w,0}$ and $\widehat{\phi}\in S$ implies
$\widehat{\phi}\overset{\vee}{\psi}\in S_{w,0}$ so that $\left(
\widehat{\phi}\overset{\vee}{\psi}\right)  ^{\wedge}\in\widehat{S}_{w,0}$ and
consequently\allowbreak\ $\left[  \left(  \widehat{\phi}\widehat{f^{e}%
}\right)  ^{\vee},\psi\right]  =\left[  \left(  \widehat{\phi}\widehat{f}%
\right)  ^{\vee},\psi\right]  $ which now allows us to write,\allowbreak
\ $\phi\ast f=\left(  \widehat{\phi}\widehat{f}\right)  ^{\vee}=\left(
2\pi\right)  ^{-\frac{d}{2}}\left[  f_{y},\phi\left(  \cdot-y\right)  \right]
$.
\end{proof}

This lemma can be applied to obtain the following representation of the
operator $\mathcal{B}$:

\begin{theorem}
\label{Thm_properties2_B}Suppose the weight function $w$ has property W01
w.r.t. the set $\mathcal{A}$ and $w\in W_{S;0}$ (Definition
\ref{Def_Sw2_and_fnal(Sw2)}). Then:

\begin{enumerate}
\item If $v\in X_{1/w}^{0}$ then $v\in\left(  S_{w,0}^{\prime}\right)  ^{\vee
}$, $v_{F}\in S_{w,0}^{\prime}\cap L_{loc}^{1}\left(  \mathbb{R}^{d}%
\setminus\mathcal{A}\right)  $ and $\widehat{v}-v_{F}\in S_{w,0;\mathcal{A}%
}^{\prime}$ where%
\begin{equation}
S_{w,0;\mathcal{A}}^{\prime}:=\left\{  g\in S_{w,0}^{\prime}%
:\operatorname*{supp}g\subset\mathcal{A}\right\}  ,\label{a2.2}%
\end{equation}

with the support given by part 7 of Definition \ref{Def_Sw2_lin_fnal_0}.

\item $\overset{\vee}{S}_{w,0}=S\cap X_{w}^{0}$.

\item If $\phi\in\overset{\vee}{S}_{w,0}$ and $v\in X_{1/w}^{0}$ then
\begin{equation}
\mathcal{B}\left(  \phi,v\right)  =\phi\ast\overset{\vee}{v_{F}}=\left(
2\pi\right)  ^{-\frac{d}{2}}\left[  \left(  \overset{\vee}{v_{F}}\right)
_{y},\phi\left(  \cdot-y\right)  \right]  .\label{a2.3}%
\end{equation}

\item If $w$ also has property W02 or W03 then $\overset{\vee}{S}_{w,0}$ is
dense in $X_{w}^{0}$.
\end{enumerate}
\end{theorem}

\begin{proof}
\textbf{Part 1} $S^{\prime}\subset\left(  \widehat{S}_{w,0}\right)  ^{\prime}$
so $v\in\left(  \widehat{S}_{w,0}\right)  ^{\prime}=\left(  S_{w,0}^{\prime
}\right)  ^{\vee}$. Suppose $\phi\in S_{w,0}$.

Then
\[
\left\vert \int v_{F}\phi\right\vert =\left\vert \int\frac{v_{F}}{\sqrt{w}%
}\sqrt{w}\phi\right\vert \leq\left\vert v\right\vert _{1/w,0}\left(  \int
w\left\vert \phi\right\vert ^{2}\right)  ^{1/2}%
\]

and the assumptions of this theorem imply that $\left(  \int w\left\vert
\phi\right\vert ^{2}\right)  ^{1/2}$ is bounded by a positive, linear
combination $\left\vert \cdot\right\vert _{\sigma}$ of the seminorms used to
define the topology on $S$ i.e. $v_{F}\in S_{w,0}^{\prime}$.

If $\phi\in C_{0}^{\infty}\left(  \mathbb{R}^{d}\setminus\mathcal{A}\right)  $
then $\left[  \widehat{v}-v_{F},\phi\right]  =\left[  \widehat{v},\phi\right]
-\left[  v_{F},\phi\right]  =0$ so that $\operatorname*{supp}\left(
\widehat{v}-v_{F}\right)  \subset\mathcal{A}$ and $\widehat{v}-v_{F}\in
S_{w,0;\mathcal{A}}^{\prime}$.\medskip

\textbf{Part 2} From part 1 Definition \ref{Def_Sw2_lin_fnal_0}.\medskip

\textbf{Part 3} From \ref{a1.20}, $\mathcal{B}\left(  \phi,v\right)  =\left(
\widehat{\phi}v_{F}\right)  ^{\vee}=\left(  \widehat{\phi}%
\widehat{\overset{\vee}{v_{F}}}\right)  ^{\vee}$ and by part 1, $\overset{\vee
}{v_{F}}\in\left(  \widehat{S}_{w,0}\right)  ^{\prime}$ so that Lemma
\ref{Lem_convol_Sw2} implies, $\allowbreak\mathcal{B}\left(  \phi,v\right)
=\phi\ast\overset{\vee}{v_{F}}=\left(  2\pi\right)  ^{-d/2}\left[  \left(
\overset{\vee}{v_{F}}\right)  _{y},\phi\left(  \cdot-y\right)  \right]
$.\medskip

\textbf{Part 4} This was shown in Corollary \ref{Cor_Xow_S_density}.
\end{proof}

??

\begin{corollary}
\label{Cor_properties_B}\textbf{Properties of }$\mathcal{B}\left(
G\mathbf{,\cdot}\right)  $. Suppose $w$ has property W02 for $\kappa$ and $G$
is the basis function. Then if $v,v^{\prime}\in X_{1/w}^{0}$:

\begin{enumerate}
\item $D^{\gamma}G\in X_{w}^{0}$ when $\left\vert \gamma\right\vert \leq
\kappa$ and%
\[
D^{\gamma}\mathcal{B}\left(  G,v\right)  =\mathcal{B}\left(  D^{\gamma
}G,v\right)  ,\quad\left\vert \gamma\right\vert \leq\kappa.
\]

\item $\mathcal{B}\left(  G,v\right)  \in X_{w}^{0}$ and $\left\Vert
\mathcal{B}\left(  G,v\right)  \right\Vert _{w,0}=\left\vert v\right\vert
_{1/w,0}$.

\item $\mathcal{B}\left(  G,\phi\right)  =G\ast\phi$ when $\phi\in S$.

\item $\int F\mathcal{B}\left(  \mathcal{B}\left(  G,v\right)  ,\overline
{v^{\prime}}\right)  =\left\langle v,v_{\ast}^{\prime}\right\rangle _{1/w,0}
$, where $F$ is the Fourier transform and $v_{\ast}^{\prime}\left(  x\right)
=v^{\prime}\left(  -x\right)  $.
\end{enumerate}
\end{corollary}

\begin{proof}
\textbf{Part 1} By Theorem \ref{Thm_basis_fn_properties_all_m_W2}, $D^{\gamma
}G\in X_{w}^{0}$ when $\left\vert \gamma\right\vert \leq\kappa$ and so by part
2 of Theorem \ref{Thm_properties_B}, $D^{\gamma}\mathcal{B}\left(  G,v\right)
\in\left(  \widehat{D^{\gamma}G}v_{F}\right)  ^{\vee}=\mathcal{B}\left(
D^{\gamma}G,v\right)  $.\medskip

\textbf{Part 2} Observe that $G\in X_{w}^{0}$ and $\mathcal{B}\left(
G,v\right)  \in S^{\prime}$. Further $F\mathcal{B}\left(  G,v\right)
=\widehat{G}v_{F}=\frac{v_{F}}{w}\in L^{1}$ by Lemma \ref{Lem_Xoinvw}, and
$\int w\left\vert F\mathcal{B}\left(  G,v\right)  \right\vert ^{2}=\int%
\frac{\left\vert v_{F}\right\vert ^{2}}{w}=\left\vert v\right\vert _{1/w,0}$.
Thus $\mathcal{B}\left(  G,v\right)  \in X_{w}^{0}$.\medskip

\textbf{Part 3} From Lemma \ref{Lem_Xoinvw}, $S\subset X_{1/w}^{0}$ and
$\mathcal{B}\left(  G,\phi\right)  =\left(  \widehat{G}\widehat{\phi}\right)
^{\vee}=G\ast\phi$ by the convolution definition \ref{1.20}.\medskip

\textbf{Part 4} $\mathcal{B}\left(  G,v\right)  =\left(  \frac{v_{F}}%
{w}\right)  ^{\vee}$ so
\begin{align*}
\int F\mathcal{B}\left(  \mathcal{B}\left(  G,v\right)  ,\overline{v^{\prime}%
}\right)  =\int F\mathcal{B}\left(  \left(  \frac{v_{F}}{w}\right)  ^{\vee
},\overline{v^{\prime}}\right)  =\int F\left(  \left(  \frac{v_{F}}{w},\left(
\overline{v^{\prime}}\right)  _{F}\right)  ^{\vee}\right)   &  =\int\frac
{1}{w}v_{F}\left(  \overline{v^{\prime}}\right)  _{F}\\
&  =\int\frac{1}{w}v_{F}\overline{\left(  v_{\ast}^{\prime}\right)  _{F}}\\
&  =\left\langle v,v_{\ast}^{\prime}\right\rangle _{1/w,0}.
\end{align*}

\end{proof}

\begin{remark}
??? WHAT IF $w\in W03$?
\end{remark}

\subsection{The operators $\mathcal{V}:X_{1/w}^{0}\rightarrow X_{w}^{0}$,
$\mathcal{W}:X_{w}^{0}\rightarrow X_{1/w}^{0}$ and $\Phi:X_{1/w}%
^{0}\rightarrow\left(  X_{w}^{0}\right)  ^{\prime}$}

The operators $\mathcal{W}=\mathcal{LI}$ and $\mathcal{V}=\mathcal{JM}$ will
be introduced:%

\begin{gather}%
\begin{array}
[c]{ccccc}%
X_{1/w}^{0} &
\begin{array}
[c]{c}%
\overset{\mathcal{L}}{\longleftarrow}\\
\underset{\mathcal{M}}{\longrightarrow}%
\end{array}
& L^{2} &
\begin{array}
[c]{c}%
\overset{\mathcal{I}}{\longleftarrow}\\
\underset{\mathcal{J}}{\longrightarrow}%
\end{array}
& X_{w}^{0}%
\end{array}
\label{Fig_ops_Xo1/w_L2_Xow}\\
\text{Isometric isomorphisms between }X_{1/w}^{0}\text{, }L^{2}\text{, }%
X_{w}^{0}\text{.}\nonumber\\
\text{Also, }\mathcal{V}=\mathcal{JM}\text{ and }\mathcal{W}=\mathcal{LI}%
\text{.}\nonumber
\end{gather}

and used to define the operator $\Phi:X_{1/w}^{0}\rightarrow\left(  X_{w}%
^{0}\right)  ^{\prime}$ which characterizes the bounded linear functionals on
$X_{w}^{0}$ as the space $X_{1/w}^{0}$. In fact, in an analogous fashion to
the negative order Sobolev spaces, Theorem \ref{Thm_op_Phi_bilinear} shows
that $\Phi$ can be written as a bilinear form.

\begin{definition}
\label{Def_op_V}\textbf{The operator} $\mathcal{V}:X_{1/w}^{0}\rightarrow
X_{w}^{0}$

The operator $\mathcal{J}:L^{2}\rightarrow X_{w}^{0}$ was introduced in
Definition \ref{Def_I_J} and studied in Theorem \ref{Thm_I_J_property}.

Suppose the weight function $w$ has property W02 or W03. These conditions
ensure that both $\mathcal{M}$ and $\mathcal{J}$ are defined.

\textbf{Define the operator} $\mathcal{V}$ \textbf{by} $\mathcal{V}%
:=\mathcal{JM}$ where $\mathcal{M}:X_{1/w}^{0}\rightarrow L^{2}$.
\end{definition}

The operator $\mathcal{V}$ can be expressed in terms of the operator
$\mathcal{B}$ and the basis function as follows:

\begin{theorem}
\label{Thm_Vu_=_G*u}Suppose the weight function $w$ has property W02 or W03
and that $G$ is the corresponding basis function. Then
\begin{equation}
\mathcal{V}v=\left(  \frac{v_{F}}{w}\right)  ^{\vee}=\mathcal{B}\left(
G,v\right)  ,\quad v\in X_{1/w}^{0},\label{a6.0}%
\end{equation}

and%
\[
\mathcal{V}v=G\ast v,\quad v\in S,
\]

and $\mathcal{V}:S\rightarrow G\ast S$ is 1-1 and onto.
\end{theorem}

\begin{proof}
By Definition \ref{Def_I_J}, $\widehat{\mathcal{J}g}=\frac{\widehat{g}}%
{\sqrt{w}}\in L_{loc}^{1}$ for $g\in L^{2}$.

By Definition \ref{Def_op_M}, $\widehat{\mathcal{M}v}=\frac{v_{F}}{\sqrt{w}%
}\in L^{2}$ for $v\in X_{1/w}^{0}$. Thus
\[
\widehat{\mathcal{V}v}=\widehat{\mathcal{JM}v}=\frac{\widehat{\mathcal{M}v}%
}{\sqrt{w}}=\frac{v_{F}}{w}=\left(  \mathcal{B}\left(  G,v\right)  \right)
^{\wedge},
\]

and so $\mathcal{V}v=\mathcal{B}\left(  G,v\right)  $. The second equation of
this theorem follows from part 3 of Corollary \ref{Cor_properties_B}.

Clearly $\mathcal{V}$ is onto. Since $\mathcal{V}$ is isometric $\mathcal{V}%
\phi=0$ implies $\left\Vert \mathcal{V}\phi\right\Vert _{w,0}=\left\vert
\phi\right\vert _{1/w,0}=\left(  \int\frac{\left\vert _{\phi_{F}}\right\vert
^{2}}{w}\right)  ^{1/2}=0$, where $\phi_{F}$ is the restriction of
$\widehat{\phi}$ to $\mathbb{R}^{d}\setminus\mathcal{A}$ where $\mathcal{A}$
is a closed set of measure zero. But $\widehat{\phi}\in S$ and so $\phi=0$.
\end{proof}

\begin{definition}
\label{Def_ops_L_W}\textbf{The operators }$\mathcal{L}:L^{2}\rightarrow
X_{1/w}^{0}$ \textbf{and }$\mathcal{W}:X_{w}^{0}\rightarrow X_{1/w}^{0}$.

Suppose $w$ has property \ref{a1.046} or $w\in W_{S;0}$. Then we define:

\textbf{The operator }$\mathcal{L}:L^{2}\rightarrow X_{1/w}^{0}$\textbf{:} if
$w$ has property \ref{a1.046} let $\mathcal{L}=\mathcal{L}_{1}$ otherwise if
$w\in W_{S;0}$ let $\mathcal{L=L}_{2}$.

\textbf{The operator} $\mathcal{W}:=\mathcal{LI}:X_{w}^{0}\rightarrow
X_{1/w}^{0}$.
\end{definition}

The next theorem relates the properties of the operators $\mathcal{V}$ and
$\mathcal{W}$, and for which we will require the properties of the operators
$\mathcal{I}$ and $\mathcal{J}$ given in Theorem \ref{Thm_M_L1_propert}.

\begin{theorem}
\label{Thm_V_W_properties}Suppose the weight function $w$ has property W02 or
W03 as well as property \ref{a1.046} or property \ref{a1.8}. Then in the
seminorm sense the operators $\mathcal{V}=\mathcal{JM}$ and $\mathcal{W}%
=\mathcal{LI}$ have the following properties:

\begin{enumerate}
\item $\mathcal{V}:X_{1/w}^{0}\rightarrow X_{w}^{0}$ is a unique, linear isometry.

\item $\mathcal{W}:X_{w}^{0}\rightarrow X_{1/w}^{0}$ is a class of linear isometries.

\item $\mathcal{V}$ and $\mathcal{W}$ are inverses.

\item $\mathcal{V}$ and $\mathcal{W}$ are onto.

\item $\mathcal{V}$ and $\mathcal{W}$ are 1-1.

\item $\mathcal{V}$ and $\mathcal{W}$ are adjoints.

\item If $f\in X_{w}^{0}$ and $w\widehat{f}\in S^{\prime}$ then $\mathcal{W}%
f=\left(  w\widehat{f}\right)  ^{\vee}$.

\item $\mathcal{V}$ and $\mathcal{W}$ are isometric isomorphisms, inverses and adjoints.
\end{enumerate}
\end{theorem}

\begin{proof}
\textbf{Part 1}\ From Theorem \ref{Thm_I_J_property} $\mathcal{J}$ is a
unique, linear isometry, and from Theorem \ref{Thm_M_property} $\mathcal{M}$
is also a unique, linear isometry.\medskip

\textbf{Part 2} From Theorem \ref{Thm_I_J_property} $\mathcal{I}$ is a unique,
linear isometry. From Theorem \ref{Thm_L1_propert} or Theorem
\ref{Thm_L2_propert} $\mathcal{L}$ is a class of linear isometries.\medskip

\textbf{Part 3} From Theorem \ref{Thm_M_L1_propert} or Theorem
\ref{Thm_M_L2_propert} $\mathcal{ML}=I$\ and Theorem \ref{Thm_I_J_property}
$\mathcal{JI}f=f$. Hence $\mathcal{VW}=\mathcal{JMLI}=\mathcal{JI}=I$.

From Theorem \ref{Thm_I_J_property} $\mathcal{I}\mathcal{J}=I$ and from
Theorem \ref{Thm_M_L1_propert} or Theorem \ref{Thm_M_L2_propert}
$\mathcal{LM}g-g\in\left(  S_{\mathcal{A}}^{\prime}\right)  ^{\vee}$. Hence
$\mathcal{WV}=\mathcal{LIJM}=\mathcal{LM}$ and $\mathcal{WV}g-g\in\left(
S_{\mathcal{A}}^{\prime}\right)  ^{\vee}$.\medskip

\textbf{Parts 4 and 5} follow directly from part 3.\medskip

\textbf{Part 6} Theorem \ref{Thm_M_L1_propert} or Theorem
\ref{Thm_M_L2_propert}\ implies $\mathcal{L}$ and $\mathcal{M}$ are adjoints
and Theorem \ref{Thm_I_J_property} implies $\mathcal{I}$ and $\mathcal{J}$ are
adjoints. Thus
\[
\left\langle g,\mathcal{W}f\right\rangle _{1/w,0}=\left\langle g,\mathcal{LI}%
f\right\rangle _{1/w,0}=\left(  \mathcal{M}g,\mathcal{I}f\right)  _{2}=\left(
\mathcal{JM}g,f\right)  _{w,0}=\left(  \mathcal{V}g,f\right)  _{w,0}.
\]
\medskip

\textbf{Part 7} If $f\in X_{w}^{0}$ then from the definitions of
$\mathcal{L}_{1}$ and $\mathcal{L}_{2}$
\[
\widehat{\mathcal{W}f}=\widehat{\mathcal{LI}f}=\sqrt{w}\widehat{\mathcal{I}%
f}=\left(  \sqrt{w}\left(  \sqrt{w}\widehat{f}\right)  \right)  ^{e}=\left(
w\widehat{f}\right)  ^{e}=w\widehat{f}.
\]
\medskip

\textbf{Part 8} Summary of previous parts.
\end{proof}

\begin{remark}
??? SEE also Definition \ref{Def_tildL}. Regarding part 7, what about the
space:%
\begin{align*}
\left\{  f\in X_{w}^{0}:w\widehat{f}\in S^{\prime}\right\}   & =\left\{  f\in
X_{w}^{0}:\left(  w\widehat{f}\right)  ^{\vee}\in S^{\prime}\right\} \\
& =\left\{  f\in X_{w}^{0}:\mathcal{W}f\in S^{\prime}\right\} \\
& =\mathcal{J}\left\{  u\in L^{2}:\sqrt{w}u\in S^{\prime}\right\}  .
\end{align*}

If $w_{1}$ exists and $\frac{1}{w_{1}}\in C_{BP}^{\infty}$ e.g. extended
B-splines, then%
\begin{align*}
\left\{  f\in X_{w}^{0}:w\widehat{f}\in S^{\prime}\right\}   & =\mathcal{J}%
_{1}\left\{  u\in L^{2}:w_{1}u\in S^{\prime}\right\} \\
& =??\mathcal{J}_{1}\left\{  u\in L^{2}:u\in\frac{1}{w_{1}}S^{\prime}\right\}
\\
& =\mathcal{J}_{1}\left(  L^{2}\cap\frac{1}{w_{1}}S^{\prime}\right)  .
\end{align*}

\end{remark}

Next we prove that $S$ is dense in $X_{1/w}^{0}$.

\begin{corollary}
\label{Cor_S_dense_Xo1/w}Suppose the weight function $w$ has property W02 or
W03 and either property \ref{a1.046} holds or else $w\in W_{S;0}$.

Then $S$ is dense in $X_{1/w}^{0}$ under $\left\vert \mathcal{\cdot
}\right\vert _{1/w,0}$.
\end{corollary}

\begin{proof}
Choose $u\in X_{1/w}^{0}$ so that $\mathcal{V}u\in X_{w}^{0}$. From part 3
Theorem \ref{Thm_Jg_properties}, $G\ast S$ is dense in $X_{w}^{0}$ and hence
given $\varepsilon>0$ there exists $\phi_{\varepsilon}\in G\ast S$ such that
$\left\Vert \mathcal{V}u-\phi_{\varepsilon}\right\Vert _{w,0}<\varepsilon$.
But by Theorem \ref{Thm_Vu_=_G*u}, $\mathcal{V}:S\rightarrow G\ast S$ is an
isomorphism and so $\mathcal{V}^{-1}\phi_{\varepsilon}\in S$. Since
$\mathcal{V}:X_{1/w}^{0}\rightarrow X_{w}^{0}$ is an isometry%
\[
\left\vert u-\mathcal{V}^{-1}\phi_{\varepsilon}\right\vert _{1/w,0}=\left\Vert
\mathcal{V}\left(  u-\mathcal{V}^{-1}\phi_{\varepsilon}\right)  \right\Vert
_{w,0}=\left\Vert \mathcal{V}u-\phi_{\varepsilon}\right\Vert _{w,0}%
<\varepsilon,
\]

which proves this density result.

Since $X_{1/w}^{0}$ is complete and $S$ is dense in $X_{1/w}^{0}$ it follows
that $X_{1/w}^{0}$ is the completion of $S$.
\end{proof}

The operator $\mathcal{V}$ will now be used characterize the bounded linear
functionals on $X_{w}^{0}$, denoted $\left(  X_{w}^{0}\right)  ^{\prime}$, in
the following sense:

\begin{theorem}
\label{Thm_op_Phi}\textbf{The operator }$\Phi$. Suppose the weight function
$w$ has property W02 or W03 and either property \ref{a1.046} holds or $w\in
W_{S;0}$. Denote by $\left(  X_{w}^{0}\right)  ^{\prime}$ the space of bounded
linear functionals on $X_{w}^{0}$.

Then the equation%
\begin{equation}
\left(  \Phi v\right)  \left(  u\right)  :=\left(  u,\mathcal{V}v\right)
_{w,0},\quad u\in X_{w}^{0},\text{ }v\in X_{1/w}^{0},\label{a1.6}%
\end{equation}

defines a linear operator $\Phi:X_{1/w}^{0}\rightarrow\left(  X_{w}%
^{0}\right)  ^{\prime}$ which is an isometric isomorphism in the seminorm sense.
\end{theorem}

\begin{proof}
Since $\mathcal{V}:X_{1/w}^{0}\rightarrow X_{w}^{0}$ is an isometry, given
$v\in X_{1/w}^{0}$ and $u\in X_{w}^{0}$ we have
\[
\left\vert \left(  u,\mathcal{V}v\right)  _{w,0}\right\vert \leq\left\Vert
u\right\Vert _{w,0}\left\Vert \mathcal{V}v\right\Vert _{w,0}=\left\Vert
u\right\Vert _{w,0}\left\vert v\right\vert _{1/w,0},
\]

and for each $v\in X_{1/w}^{0}$, the expression $\left(  u,\mathcal{V}%
v\right)  _{w,0}$, $u\in X_{w}^{0}$ defines a bounded linear functional on
$X_{w}^{0}$. Denote this functional by $\Phi v$ so that \ref{a1.6} holds and
the operator norm is
\[
\left\Vert \Phi v\right\Vert _{op}=\sup\limits_{u\in X_{w}^{0}}\frac
{\left\vert \left(  u,\mathcal{V}v\right)  _{w,0}\right\vert }{\left\Vert
u\right\Vert _{w,0}}\leq\left\vert v\right\vert _{1/w,0}.
\]

Thus $\Phi:X_{1/w}^{0}\rightarrow\left(  X_{w}^{0}\right)  ^{\prime}$. In
fact, we can easily prove that $\left\Vert \Phi v\right\Vert _{op}=\left\vert
v\right\vert _{1/w,0}$ by noting that when $u=\widetilde{\mathcal{V}}v$,
$\frac{\left\vert \left(  u,\mathcal{V}v\right)  _{w,0}\right\vert
}{\left\Vert u\right\Vert _{w,0}}=\left\Vert \mathcal{V}v\right\Vert
_{w,0}=\left\vert v\right\vert _{1/w,0}$. Clearly $\Phi v=0$ implies
$\left\vert v\right\vert _{1/w,0}=0$ so $\Phi$ is 1-1.

To prove $\Phi$ is onto choose $\mathcal{Y}\in\left(  X_{w}^{0}\right)
^{\prime}$. Since $X_{w}^{0}$ is an inner product space there exists
$\lambda\in X_{w}^{0}$ such that $\mathcal{Y}u=\left(  u,\lambda\right)
_{w,0}$ when $u\in X_{w}^{0}$. But by parts 3 and 6 of Theorem
\ref{Thm_V_W_properties}, $\mathcal{VV}^{\ast}=I$ so $\mathcal{Y}u=\left(
u,\lambda\right)  _{w,0}=\left(  u,\mathcal{VV}^{\ast}\lambda\right)  _{w,0}$
and comparison with \ref{a1.6} yields $\left(  \Phi\mathcal{V}^{\ast}%
\lambda\right)  \left(  u\right)  =\mathcal{Y}u$ and $\Phi\mathcal{V}^{\ast
}\lambda=\mathcal{Y}$.
\end{proof}

In an analogous fashion to the negative order Sobolev spaces a bilinear form
can be used to characterize the bounded linear functionals on $X_{w}^{0}$:

\begin{theorem}
\label{Thm_op_Phi_bilinear}If $u\in X_{w}^{0}$ and $v\in X_{1/w}^{0}$ then the
operator $\Phi:X_{1/w}^{0}\rightarrow\left(  X_{w}^{0}\right)  ^{\prime} $ can
be expressed directly in terms of the bilinear form $\int\widehat{u}%
\overline{v_{F}}$ as%
\[
\left(  \Phi v\right)  \left(  u\right)  =\int\widehat{u}\overline{v_{F}%
},\quad u\in X_{w}^{0},\text{ }v\in X_{1/w}^{0}.
\]

\end{theorem}

\begin{proof}
A direct consequence of \ref{a1.6} and \ref{a6.1}.
\end{proof}

If $w$ has property W02 or W03 for parameter $\kappa$ then from part 2 Theorem
\ref{Thm_ord0_Riesz_rep_W2} there exists $R_{x}\in X_{w}^{0}$ such that for
$\left\vert \gamma\right\vert \leq\kappa$, $D^{\gamma}R_{x}\in X_{w}^{0}$ and
$\left(  f,\left(  -D\right)  ^{\gamma}R_{x}\right)  _{w,0}=D^{\gamma}f\left(
x\right)  $. On the other hand, the space $X_{1/w}^{0}$is not a reproducing
kernel Hilbert space because it contains tempered distributions which are not
continuous functions e.g. $\delta\left(  \cdot-x\right)  $. However it does
have the following evaluation property \ref{a1.006}:

\begin{theorem}
\label{Thm_V1*(DRx)=delta}Suppose the weight function $w$ has property W02 for
parameter $\kappa$ as well as property \ref{a1.046} or $w\in W_{S;0}$. Then%
\[
\mathcal{V}^{\ast}D^{\gamma}\left(  R_{x}\right)  =D^{\gamma}\delta\left(
\cdot-x\right)  ,\quad\left\vert \gamma\right\vert \leq\kappa,
\]

where $R_{x}$ is the Riesz representer of the evaluation functional
$f\rightarrow f\left(  x\right)  $. Further%
\begin{equation}
\left\langle u,D^{\gamma}\delta\left(  \cdot-x\right)  \right\rangle
_{1/w,0}=\left(  -1\right)  ^{\left\vert \gamma\right\vert }\left(  D^{\gamma
}\mathcal{V}u\right)  \left(  x\right)  ,\quad u\in X_{1/w}^{0},\text{
}\left\vert \gamma\right\vert \leq\kappa.\label{a1.006}%
\end{equation}

\end{theorem}

\begin{proof}
If $f\in X_{1/w}^{0}$ and $w\widehat{f}\in S^{\prime}$ then by parts 6 and 7
of Theorem \ref{Thm_V_W_properties},
\[
\mathcal{V}^{\ast}f=\mathcal{W}f=\left(  w\widehat{f}\right)  ^{\vee}.
\]

Thus when $f=D^{\gamma}R_{x}$, by using the Fourier transform properties given
in Appendix \ref{Ch_Appendx_basic_notation}, we obtain%
\[
w\widehat{D^{\gamma}R_{x}}=w\left(  \xi\right)  \left(  i\xi\right)  ^{\gamma
}\widehat{R_{x}}\left(  \xi\right)  =\left(  2\pi\right)  ^{-d/2}w\left(
\xi\right)  \left(  i\xi\right)  ^{\gamma}\widehat{G\left(  \cdot-x\right)
}\left(  \xi\right)  =\left(  2\pi\right)  ^{-d/2}\left(  i\xi\right)
^{\gamma}e^{-ix\xi}\in S^{\prime},
\]

and as a consequence
\[
\mathcal{V}^{\ast}D^{\gamma}\left(  R_{x}\right)  =\left(  w\widehat{D^{\gamma
}R_{x}}\right)  ^{\vee}=\left(  2\pi\right)  ^{-d/2}\left(  \left(
i\xi\right)  ^{\gamma}e^{-ix\xi}\right)  ^{\vee}=\left(  2\pi\right)
^{-d/2}D^{\gamma}\left(  \left(  e^{-ix\xi}\right)  ^{\vee}\right)  .
\]

But $\left(  e^{-ia\xi}\right)  ^{\vee}=\left(  2\pi\right)  ^{d/2}%
\delta\left(  \cdot-a\right)  $ so%
\[
\mathcal{V}^{\ast}D^{\gamma}\left(  R_{x}\right)  =D^{\gamma}\delta\left(
\cdot-x\right)  ,
\]

and thus%
\begin{align*}
\left\langle u,D^{\gamma}\delta\left(  \cdot-x\right)  \right\rangle
_{1/w,0}=\left\langle u,\mathcal{V}^{\ast}D^{\gamma}\left(  R_{x}\right)
\right\rangle _{1/w,0}=\left(  \mathcal{V}u,D^{\gamma}R_{x}\right)  _{w,0} &
=\left(  \left(  -D\right)  ^{\gamma}\mathcal{V}u,R_{x}\right)  _{w,0}\\
&  =\left(  -1\right)  ^{\left\vert \gamma\right\vert }\left(  D^{\gamma
}\mathcal{V}u\right)  \left(  x\right)  .
\end{align*}

\end{proof}

\begin{remark}
??? If $w$ has property W03 then this theorem is also true if we replace
$\left\vert \gamma\right\vert \leq\kappa$ by $\gamma\leq\kappa$.
\end{remark}

\subsection{The space $X_{w}^{0}\left(  \Omega\right)  ^{\prime}%
$\label{SbSect_fnal_locXow}}

?? How does the operator $\int_{\Omega}\overline{R}_{x}u$ in Theorem
\ref{Thm_property_integomegaRxu} above affect this?

From Theorem \ref{Thm_FnalQuotSp_AnnhilDenom},%
\begin{align*}
V  & :\left(  X/A\right)  ^{\prime}\rightarrow A^{0}\text{ }is\text{ }an\text{
}isom.\text{ }iso.\\
A^{0}  & =\left\{  u\in X^{\prime}:u\left(  A\right)  =\left\{  0\right\}
\right\}  .\\
\left(  Vu^{\prime}\right)  \left(  x\right)   & :=u^{\prime}\left(  x\right)
.
\end{align*}

Now consider the case: $X=X_{w}^{0}$ and $A=\left(  X_{w}^{0}\right)
_{\Omega^{c}}=\left\{  f\in X_{w}^{0}:\operatorname*{supp}f\subset\Omega
^{c}\right\}  $. $A $ is closed by Theorem \ref{Thm_XowK_Hilbert}. From
Corollary \ref{Cor1_Thm_canon_exten_op} we have%
\[
X_{w}^{0}\left(  \Omega\right)  \overset{\hom}{=}X_{w}^{0}/\left(  X_{w}%
^{0}\right)  _{\Omega^{c}}.
\]

\begin{equation}
X_{w}^{0}\left(  \Omega\right)  ^{\prime}=\operatorname*{anhil}\left(  \left(
X_{w}^{0}\right)  _{\Omega^{c}}\right)  =\left\{  \sigma\in\left(  X_{w}%
^{0}\right)  ^{\prime}:\sigma\left(  \left(  X_{w}^{0}\right)  _{\Omega^{c}%
}\right)  =\left\{  0\right\}  \right\}  .\label{1.101}%
\end{equation}

See the (several) definitions of the operators $\mathcal{B}$, $\mathcal{B}%
\left(  G,v\right)  $, $\mathcal{V}$ and $\Phi$ used in relation to the spaces
$X_{1/w}^{0}$ and $\left(  X_{w}^{0}\right)  ^{\prime}$.

When does%
\[
X_{w}^{0}\left(  \Omega\right)  ^{\prime}=\left\{  \int_{\Omega}\overline
{R}_{x}u:u\in X_{w}^{0}\right\}  ?
\]

From Theorem \ref{Thm_op_Phi_bilinear}: if $u\in X_{w}^{0}$ and $v\in
X_{1/w}^{0}$ then the isometric isomorphism $\Phi:X_{1/w}^{0}\rightarrow
\left(  X_{w}^{0}\right)  ^{\prime}$ introduced in Theorem \ref{Thm_op_Phi}
can be expressed directly in terms of the bilinear form $\int\widehat{u}%
\overline{v_{F}}$ as%
\[
\left(  \Phi v\right)  \left(  u\right)  =\int\widehat{u}\overline{v_{F}%
},\quad u\in X_{w}^{0},\text{ }v\in X_{1/w}^{0}.
\]

Thus%
\begin{align*}
X_{w}^{0}\left(  \Omega\right)  ^{\prime}  & =\left\{  \Phi v:v\in X_{1/w}%
^{0}\text{ }and\text{ }\left(  \Phi v\right)  \left(  \left(  X_{w}%
^{0}\right)  _{\Omega^{c}}\right)  =\left\{  0\right\}  \right\} \\
& =\left\{  \Phi v:v\in X_{1/w}^{0}\text{ }and\text{ }\left(  \Phi v\right)
\left(  u\right)  =0\text{ }\forall u\in\left(  X_{w}^{0}\right)  _{\Omega
^{c}}\right\} \\
& =\left\{  \Phi v:v\in X_{1/w}^{0}\text{ }and\text{ }\int\widehat{u}%
\overline{v_{F}}=0\text{ }\forall u\in\left(  X_{w}^{0}\right)  _{\Omega^{c}%
}\right\} \\
& =\Phi\left\{  v\in X_{1/w}^{0}:\int\widehat{u}\overline{v_{F}}=0\text{
}\forall u\in\left(  X_{w}^{0}\right)  _{\Omega^{c}}\right\}
\end{align*}

What is $C_{\Omega^{c}}^{\infty}$? Choose $C_{\Omega^{c}}^{\infty}=\left\{
\phi\in C_{0}^{\infty}:\operatorname*{supp}\phi\subset\Omega^{c}\right\}  $?
PERHAPS\ $S_{\Omega^{c}}:=\left\{  \phi\in S:\operatorname*{supp}\phi
\subset\Omega^{c}\right\}  $?

??? IS $C_{\Omega^{c}}^{\infty}\cap\left(  X_{w}^{0}\right)  _{\Omega^{c}}$
dense in $\left(  X_{w}^{0}\right)  _{\Omega^{c}}$? Petersen? Is
$\overset{\vee}{S}_{w,0}\cap C_{\Omega^{c}}^{\infty}\cap\left(  X_{w}%
^{0}\right)  _{\Omega^{c}}$ dense in $\left(  X_{w}^{0}\right)  _{\Omega^{c}}%
$? Is $\overset{\vee}{S}_{w,0}\cap C_{\Omega^{c}}^{\infty}\cap\left(
X_{w}^{0}\right)  _{\Omega^{c}}$ non-empty? NOTE that $\overset{\vee}{S}%
_{w,0}=X_{w}^{0}\cap S$ so $\overset{\vee}{S}_{w,0}\cap C_{\Omega^{c}}%
^{\infty}\cap\left(  X_{w}^{0}\right)  _{\Omega^{c}}=X_{w}^{0}\cap S\cap
C_{\Omega^{c}}^{\infty}\cap\left(  X_{w}^{0}\right)  _{\Omega^{c}}=S\cap
C_{\Omega^{c}}^{\infty}\cap\left(  X_{w}^{0}\right)  _{\Omega^{c}}%
=C_{\Omega^{c}}^{\infty}\cap\left(  X_{w}^{0}\right)  _{\Omega^{c}}$. SEE
remark about Theorem \ref{Thm_XowK_Hilbert}.

If true then%
\[
X_{w}^{0}\left(  \Omega\right)  ^{\prime}=\Phi\left\{  v\in X_{1/w}^{0}%
:\int\overline{v_{F}}\widehat{\phi}=0\text{ }\forall\phi\in\overset{\vee
}{S}_{w,0}\cap C_{\Omega^{c}}^{\infty}\right\}  .
\]

From part 1 of Theorem \ref{Thm_properties2_B}: if $v\in X_{1/w}^{0}$ then
$v\in\left(  S_{w,0}^{\prime}\right)  ^{\vee}$, $v_{F}\in S_{w,0}^{\prime}$
and $\widehat{v}-v_{F}\in S_{w,0;\mathcal{A}}^{\prime}$ where%
\[
S_{w,0;\mathcal{A}}^{\prime}:=\left\{  g\in S_{w,0}^{\prime}%
:\operatorname*{supp}g\subset\mathcal{A}\right\}  ,
\]

with the support given by part 7 of Definition \ref{Def_Sw2_lin_fnal_0}.
Consequently
\[
0=\int\overline{v_{F}}\widehat{\phi}=\left[  \overline{v_{F}},\widehat{\phi
}\right]  =\left[  \overline{v_{F}},\widehat{\phi}\right]  =\left[
\overline{\widehat{v}},\widehat{\phi}\right]  =\left[  \overline{v}%
,\phi\right]  ,
\]

and%
\begin{align*}
X_{w}^{0}\left(  \Omega\right)  ^{\prime}  & =\Phi\left\{  v\in X_{1/w}^{0}:
\left[  \overline{v},\phi\right]  =0\text{ }\forall\phi\in\overset{\vee
}{S}_{w,0}\cap C_{\Omega^{c}}^{\infty}\right\} \\
& =\Phi\left\{  v\in X_{1/w}^{0}:\left[  v,\phi\right]  =0\text{ }\forall
\phi\in\overset{\vee}{S}_{w,0}\cap C_{\Omega^{c}}^{\infty}\right\}  .
\end{align*}

?? SEE pages 306 - 309 of Petersen \cite{Petersen83}.

Perhaps define%
\[
X_{1/w}^{0}\left(  \Omega\right)  :=r_{\Omega}X_{1/w}^{0}.
\]

WHEN does%
\[
X_{1/w}^{0}\left(  \Omega\right)  =\left\{  v\in X_{1/w}^{0}:\left[
v,\phi\right]  =0\text{ }\forall\phi\in\overset{\vee}{S}_{w,0}\cap
C_{\Omega^{c}}^{\infty}\right\}  ?
\]

\section{The space $\protect\widetilde{X}_{1/w}^{0}$%
\label{Sect_tildXo1/w_S2dag}}

In this section we assume that the weight function belongs to $W_{S;0}$
(\ref{a1.8}). Examination of the results concerning the space $X_{1/w}^{0}$ of
the last section shows that a larger space $\widetilde{X}_{1/w}^{0}$ can be
defined which allows similar, but simpler, maps to be defined to and from
$L^{2}$.

We then set about defining the isometric maps analogous to $\mathcal{M}$,
$\mathcal{L}$ $\left(  \mathcal{L}_{1}\text{ or }\mathcal{L}_{2}\right)  $,
$\mathcal{B}$, $\mathcal{V}$, $\mathcal{W}$ and $\Phi$ and proving their
properties; the goal being to characterize the bounded linear functionals on
$X_{w}^{0}$. Some of these maps are indicated in Figure
\ref{Fig_ops_Xo1/w_tildXo1/w_Xow_Sw2}:%
\begin{gather}%
\begin{array}
[c]{ccccccc}%
X_{1/w}^{0} &
\begin{array}
[c]{c}%
\overset{E}{\longleftarrow}\\
\underset{R}{\longrightarrow}%
\end{array}
& \widetilde{X}_{1/w}^{0} &
\begin{array}
[c]{c}%
\overset{\widetilde{\mathcal{L}}}{\longleftarrow}\\
\underset{\widetilde{\mathcal{M}}_{1}}{\longrightarrow}%
\end{array}
& L^{2} &
\begin{array}
[c]{c}%
\overset{\mathcal{I}}{\longleftarrow}\\
\underset{\mathcal{J}}{\longrightarrow}%
\end{array}
& X_{w}^{0}%
\end{array}
\label{Fig_ops_Xo1/w_tildXo1/w_Xow_Sw2}\\
\widetilde{\mathcal{W}}_{1}=\widetilde{\mathcal{L}}\mathcal{I}\text{ and
}\widetilde{\mathcal{V}}_{1}=\mathcal{J}\widetilde{\mathcal{M}}_{1}\nonumber
\end{gather}

There is a great deal of commonality in the proofs but putting the proofs in
an appendix makes editing difficult so the proofs are abbreviated where possible.

The relationship between the spaces $\widetilde{X}_{1/w}^{0}$ and $X_{w}^{0}$
is studied in Subsection \ref{SbSect_E_R_tildXo1/w_Xo1/w} where we construct
an isometric, isomorphism between $\widetilde{X}_{1/w}^{0}$ and $X_{w}^{0}$.

Next, as illustrated in Figure \ref{Fig_ops_Y_tildY_Fourier_indep} of
Subsection \ref{SbSect_tildYo1/w_Yow_FTindep} we use Fourier transforms to map
$\widetilde{X}_{1/w}^{0}$ and $X_{w}^{0}$ onto Fourier-independent spaces so
that the equivalents of the mappings $\widetilde{\mathcal{L}}$,
$\widetilde{\mathcal{M}}_{1}$, $\mathcal{I}$ and $\mathcal{J}$ become Fourier-independent.

Finally, it is shown that the Gaussian and the shifted thin-plate splines do
not lie in $W_{S;0}$.

\begin{theorem}
\label{Thm_w_phi_in_fnalSw2}Suppose that $w$ has property W01 w.r.t. the set
$\mathcal{A}$ and $w\in W_{S;0}$. Then:
\end{theorem}

\begin{enumerate}
\item Suppose $\phi\in S$. Then $w\phi\in L_{loc}^{1}\left(  \mathbb{R}%
^{d}\setminus\mathcal{A}\right)  $. Further, $\phi\in S_{w,0}$ iff $w\phi\in
S_{w,0}^{\prime}$ in the sense of part 5 of Definition
\ref{Def_Sw2_lin_fnal_0}.

\item Suppose $\phi\in S$. Then $w\widehat{\phi}\in L_{loc}^{1}\left(
\mathbb{R}^{d}\setminus\mathcal{A}\right)  $. Further, $\phi\in\overset{\vee
}{S}_{w,0}$ iff $\left(  w\widehat{\phi}\right)  ^{\vee}\in\left(
\widehat{S}_{w,0}\right)  ^{\prime}$ and $u=\left(  w\widehat{\phi}\right)
^{\vee}$ satisfies
\begin{equation}
\int w\left\vert \widehat{\phi}\right\vert ^{2}=\int\frac{\left\vert
\widehat{u}\right\vert ^{2}}{w}.\label{1.059}%
\end{equation}

\end{enumerate}

\begin{proof}
\textbf{Part 1} Clearly $w\phi\in L_{loc}^{1}\left(  \mathbb{R}^{d}%
\setminus\mathcal{A}\right)  $ if $\phi\in S$. Now suppose $\phi\in S_{w,0}$.
Then $w\phi\in L_{loc}^{1}\left(  \mathbb{R}^{d}\setminus\mathcal{A}\right)  $
since $w\in L_{loc}^{1}\left(  \mathbb{R}^{d}\setminus\mathcal{A}\right)  $.
Further, if $\psi\in S_{w,0}$ then
\[
\left[  w\phi,\psi\right]  =\int w\phi\psi\leq\int\sqrt{w}\phi\sqrt{w}\psi
\leq\left(  \int w\left\vert \phi\right\vert ^{2}\right)  ^{1/2}\left(  \int
w\left\vert \psi\right\vert ^{2}\right)  ^{1/2}\leq\left\vert \phi\right\vert
_{\sigma}\left\vert \psi\right\vert _{\sigma},
\]

so that $w\phi\in S_{w,0}^{\prime}$.\medskip

On the other hand if $w\phi\in S_{w,0}^{\prime}$ there exists a positive,
linear combination $\left\vert \cdot\right\vert _{\sigma\left(  \phi\right)
}$ of the seminorms which define $S$ such that $\int w\phi\overline{\psi}=\int
w\phi\overline{\psi}\leq\left\vert \overline{\psi}\right\vert _{\sigma\left(
\phi\right)  }$ when $\psi\in S_{w,0}$. Setting $\psi=\phi$ we now obtain
$\int w\left\vert \phi\right\vert ^{2}<\infty$ which means that $\phi\in
S_{w,0}$.\medskip

\textbf{Part 2} A simple consequence of part 1.
\end{proof}

\begin{definition}
\label{Def_tildXo1/w_1}\textbf{The semi-inner product space }$\widetilde{X}%
_{1/w}^{0}$ From part 4 of Definition \ref{Def_Sw2_lin_fnal_0}, $u\in\left(
\widehat{S}_{w,0}\right)  ^{\prime}$ implies $\widehat{u}\in S_{w,0}^{\prime}%
$, and from part 5 it is meaningful to require $u\in S_{w,0}^{\prime}$ be a
member of $L_{loc}^{1}\left(  \Omega\right)  $ for any open $\Omega
\subseteq\mathbb{R}^{d}\setminus\mathcal{A}$.

Now suppose $w$ has property W01 w.r.t. the set $\mathcal{A}$ and that $w\in
W_{S;0}$. Then we define the semi-inner product space:
\begin{equation}
\widetilde{X}_{1/w}^{0}=\left\{  u\in\left(  \widehat{S}_{w,0}\right)
^{\prime}:u_{F}\in L_{loc}^{1}\left(  \mathbb{R}^{d}\setminus\mathcal{A}%
\right)  \text{ }and\text{ }\int\frac{\left\vert u_{F}\right\vert ^{2}}%
{w}<\infty\right\}  ,\label{a2.5}%
\end{equation}

where $u_{F}$ is the restriction of $\widehat{u}\in S_{w,0}^{\prime}$ to
$\mathbb{R}^{d}\setminus\mathcal{A}$. Endow $\widetilde{X}_{1/w}^{0}$ with the
semi-norm and semi-inner product%
\[
\left\vert u\right\vert _{w,0}=\left(  \int\frac{\left\vert u_{F}\right\vert
^{2}}{w}\right)  ^{1/2},\quad\left\langle u,v\right\rangle _{w,0}=\int%
\frac{u_{F}\overline{v_{F}}}{w}.
\]

Part 2 of Theorem \ref{Thm_w_phi_in_fnalSw2} now implies that $\widetilde{X}%
_{1/w}^{0}$ is not empty: in fact%
\begin{equation}
\left\{  \left(  w\widehat{\phi}\right)  ^{\vee}:\phi\in\overset{\vee
}{S}_{w,0}\right\}  =\left(  wS_{w,0}\right)  ^{\vee}\subset\widetilde{X}%
_{1/w}^{0}.\label{a3.7}%
\end{equation}

\end{definition}

\begin{remark}
\label{Rem_Def_tildXo1/w_1}Clearly $S^{\prime}\subset\left(  \widehat{S}%
_{w,0}\right)  ^{\prime}$ and so $X_{1/w}^{0}\subset\widetilde{X}_{1/w}^{0}$
as sets.
\end{remark}

\begin{theorem}
\label{Thm_property1_tildXo1/w}\textbf{Some properties of }$\widetilde{X}%
_{1/w}^{0}$

\begin{enumerate}
\item $\operatorname{null}\left\vert \cdot\right\vert _{1/w,0}=\left(
S_{w,0;\mathcal{A}}^{\prime}\right)  ^{\vee}$ where%
\begin{equation}
S_{w,0;\mathcal{A}}^{\prime}:=\left\{  u\in S_{w,0}^{\prime}%
:\operatorname*{supp}u\subset\mathcal{A}\text{ }as\text{ }a\text{
}distribution\right\}  .\label{a2.8}%
\end{equation}

\item Noting part 1 of Remark \ref{Rem_Def_extend_wt_fn}, $\left\vert
\cdot\right\vert _{1/w,0}$ is a norm iff we can choose $\mathcal{A}$ to be the
(minimal) empty set.

If, in addition, $w$ has property W02 for $\kappa$ then:

\item If $\left\vert v\left(  \xi\right)  \right\vert \leq c\left(
1+\left\vert \xi\right\vert \right)  ^{\kappa}$ then $\overset{\vee}{v}%
\in\widetilde{X}_{1/w}^{0}$.

\item $\overset{\vee}{S}_{w,0}\subset X_{w}^{0}$, $S\subset\widetilde{X}%
_{1/w}^{0}$ and $\left(  wS_{w,0}\right)  ^{\vee}=S\cap\widetilde{X}_{1/w}%
^{0}$.

\item $u\in\widetilde{X}_{1/w}^{0}$ implies $u_{F}\in S_{w,0}^{\prime}\cap
L_{loc}^{1}\left(  \mathbb{R}^{d}\setminus\mathcal{A}\right)  $, $\left(
u_{F}\right)  ^{\vee}\in\widetilde{X}_{1/w}^{0}$ and $u-\left(  u_{F}\right)
^{\vee}\in\left(  S_{w,0;\mathcal{A}}^{\prime}\right)  ^{\vee}$.

\item If $u\in\widetilde{X}_{1/w}^{0}$ then $\frac{\left(  D^{\alpha}u\right)
_{F}}{w}\in L^{1}$ when $\left\vert \alpha\right\vert \leq\kappa$.

Also, $\int\frac{\left(  D^{\alpha}u\right)  _{F}}{w}\in\left(  \widetilde{X}%
_{1/w}^{0}\right)  ^{\prime}$ and $\left\Vert \int\frac{\left(  D^{\alpha
}u\right)  _{F}}{w}\right\Vert _{op}\leq\left(  \int\frac{\xi^{2\alpha}}%
{w}\right)  ^{1/2}$.
\end{enumerate}
\end{theorem}

\begin{proof}
\textbf{Part 1}
\begin{align}
&  \operatorname{null}\left\vert \cdot\right\vert _{1/w,0}\nonumber\\
&  =\left\{  v\in\widetilde{X}_{1/w}^{0}:\left\vert v\right\vert
_{1/w,0}=0\right\} \nonumber\\
&  =\left\{  v\in\left(  \widehat{S}_{w,0}\right)  ^{\prime}:\widehat{v}\in
S_{w,0}^{\prime}\text{ }and\text{ }\widehat{v}=0\text{ }on\text{ }%
\mathbb{R}^{d}\setminus\mathcal{A}\text{ }as\text{ }a\text{ }%
distribution\right\} \nonumber\\
&  =\left\{  v\in\left(  S_{w,0}^{\prime}\right)  ^{\vee}:\widehat{v}\in
S_{w,0}^{\prime}\text{ }and\text{ }\widehat{v}=0\text{ }on\text{ }%
\mathbb{R}^{d}\setminus\mathcal{A}\text{ }as\text{ }a\text{ }%
distribution\right\} \nonumber\\
&  =\left\{  \overset{\vee}{u}:u\in S_{w,0}^{\prime}\text{ }and\text{
}u=0\text{ }on\text{ }\mathbb{R}^{d}\setminus\mathcal{A}\text{ }as\text{
}a\text{ }distribution\right\} \label{a3.2}\\
&  =\left\{  \overset{\vee}{u}:u\in S_{w,0}^{\prime}\text{ }and\text{
}\operatorname*{supp}u\subset\mathcal{A}\text{ }as\text{ }a\text{
}distribution\right\} \nonumber\\
&  =\left(  S_{w,0;\mathcal{A}}^{\prime}\right)  ^{\vee}.\nonumber
\end{align}
\medskip

\textbf{Part 2} If $\mathcal{A}$ is empty then \ref{a3.2} implies that
$\operatorname{null}\left\vert \cdot\right\vert _{1/w,0}=\left\{  0\right\}
$. If $\operatorname{null}\left\vert \cdot\right\vert _{1/w,0}=\left\{
0\right\}  $ then it follows that $u\in S_{w,0}^{\prime}$ and $u=0$ on
$\mathbb{R}^{d}\setminus\mathcal{A}$ implies $u=0$ on $\mathbb{R}^{d}$. Now if
$\mathcal{A}$ is not empty and $a\in\mathcal{A}$ then $\delta\left(
\cdot-a\right)  \in S^{\prime}$ so $\delta\left(  \cdot-a\right)  \in
S_{w,0}^{\prime}$ and $\delta\left(  \cdot-a\right)  =0$ on $\mathbb{R}%
^{d}\setminus\mathcal{A}$. But $\delta\left(  \cdot-a\right)  \neq0$ on
$\mathbb{R}^{d}$ which is a contradiction. Thus $\mathcal{A}$ must be
empty.\medskip

\textbf{Part 3} Clearly $v\in S^{\prime}$ so $\overset{\vee}{v}\in S^{\prime}%
$. Also, $\left(  \overset{\vee}{v}\right)  _{F}=v\in L_{loc}^{1}$ and
property W02 or W03 imply $\int\frac{\left\vert \left(  \overset{\vee
}{v}\right)  _{F}\right\vert ^{2}}{w}<\infty$.\medskip

\textbf{Part 4} That $\overset{\vee}{S}_{w,0}\subset X_{w}^{0}$ was shown in
part 2 of Theorem \ref{Thm_properties2_B}. That $S\subset\widetilde{X}%
_{1/w}^{0}$ follows directly from part 2. That $\left(  wS_{w,0}\right)
^{\vee}=S\cap\widetilde{X}_{1/w}^{0}$ is just \ref{a3.7}.\medskip

\textbf{Part 5} By definition, $u\in\widetilde{X}_{1/w}^{0}$ implies
$u\in\left(  \widehat{S}_{w,0}\right)  ^{\prime}$, $u_{F}\in L_{loc}%
^{1}\left(  \mathbb{R}^{d}\setminus\mathcal{A}\right)  $ and $\int%
\frac{\left\vert u_{F}\right\vert ^{2}}{w}<\infty$.

Thus $\widehat{u}\in S_{w,0}^{\prime}$. If $\phi\in\widehat{S}_{w,0}$ then by
the Cauchy-Schwartz inequality%
\[
\left\vert \int u_{F}\phi\right\vert \leq\int\left\vert u_{F}\right\vert
\left\vert \phi\right\vert =\int\frac{\left\vert u_{F}\right\vert }{\sqrt{w}%
}\sqrt{w}\left\vert \phi\right\vert \leq\left\vert u\right\vert _{1/w,0}%
\left(  \int w\left\vert \phi\right\vert ^{2}\right)  ^{1/2},
\]

and since $w\in W_{S;0}$ it follows that $u_{F}\in S_{w,0}^{\prime}$. We can
now conclude that $\left(  u_{F}\right)  ^{\vee}\in\left(  \widehat{S}%
_{w,0}\right)  ^{\prime}$, $\left(  \left(  u_{F}\right)  ^{\vee}\right)
_{F}=u_{F}\in L_{loc}^{1}\left(  \mathbb{R}^{d}\setminus\mathcal{A}\right)  $
and $\int\frac{\left\vert \left(  \left(  u_{F}\right)  ^{\vee}\right)
_{F}\right\vert ^{2}}{w}<\infty$ so that $\left(  u_{F}\right)  ^{\vee}%
\in\widetilde{X}_{1/w}^{0}$ and
\[
\left\vert u-\left(  u_{F}\right)  ^{\vee}\right\vert _{1/w,0}^{2}=\int%
\frac{\left\vert \left(  u-\left(  u_{F}\right)  ^{\vee}\right)
_{F}\right\vert ^{2}}{w}=\int\frac{\left\vert u_{F}-\left(  \left(
u_{F}\right)  ^{\vee}\right)  _{F}\right\vert ^{2}}{w}=0,
\]

i.e. $u-\left(  u_{F}\right)  ^{\vee}\in\left(  S_{w,0;\mathcal{A}}^{\prime
}\right)  ^{\vee}$.\medskip

\textbf{Part 6} $\int\frac{\left\vert \left(  D^{\alpha}u\right)
_{F}\right\vert }{w}=\int\frac{\left\vert \xi^{\alpha}u_{F}\right\vert }%
{w}=\int\frac{\left\vert \xi^{\alpha}\right\vert }{\sqrt{w}}\frac{\left\vert
u_{F}\right\vert }{\sqrt{w}}\leq\left(  \int\frac{\xi^{2\alpha}}{w}\right)
^{1/2}\left(  \int\frac{\left\vert u_{F}\right\vert ^{2}}{w}\right)
^{1/2}=\left(  \int\frac{\left\vert \cdot\right\vert ^{2\left\vert
\alpha\right\vert }}{w}\right)  ^{1/2}\left\vert u\right\vert _{1/w,0}<\infty$.
\end{proof}

\begin{remark}
\label{Rem_Thm_property1_tildXo1/w}

\textbf{1}. Part 5 forms the basis for Subsection
\ref{SbSect_E_R_tildXo1/w_Xo1/w}.

\textbf{2}. What if $w\in W03$?
\end{remark}

??

\begin{remark}
ADD MARKER ?? Prove analogue of Theorem \ref{Thm_barXow_tenprod_quotient}!
\end{remark}

\subsection{The maps $\protect\widetilde{\mathcal{M}}_{1}:X_{1/w}%
^{0}\rightarrow L^{2}$ and $\protect\widetilde{\mathcal{L}}:L^{2}\rightarrow
X_{1/w}^{0}$}

We will now define the equivalent of the operators $\mathcal{M}:X_{1/w}%
^{0}\rightarrow L^{2}$ (Subsection \ref{SbSect_op_M}) and $\mathcal{L}%
:L^{2}\rightarrow X_{1/w}^{0}$ (Definition \ref{Def_ops_L_W}) and prove their
properties. This will be done in a very similar manner.

\begin{definition}
\label{Def_tildM1}\textbf{The operator} $\widetilde{\mathcal{M}}%
_{1}:\widetilde{X}_{1/w}^{0}\rightarrow L^{2}$

From the definition of $\widetilde{X}_{1/w}^{0}$, $u\in\widetilde{X}_{1/w}%
^{0}$ implies $\frac{u_{F}}{\sqrt{w}}\in L^{2}$ where $u_{F}$ is the
restriction of $\widehat{u}$ to $\mathbb{R}^{d}\setminus\mathcal{A}$. We can
now define the linear mapping $\widetilde{\mathcal{M}}_{1}:\widetilde{X}%
_{1/w}^{0}\rightarrow L^{2}$ by%
\[
\widetilde{\mathcal{M}}_{1}u=\left(  \frac{u_{F}}{\sqrt{w}}\right)  ^{\vee
},\text{\quad}u\in\widetilde{X}_{1/w}^{0}.
\]

\end{definition}

The operator $\widetilde{\mathcal{M}}_{1}$ has the following properties:

\begin{theorem}
\label{Thm_tildM1_property}The operator $\widetilde{\mathcal{M}}%
_{1}:\widetilde{X}_{1/w}^{0}\rightarrow L^{2}$ is linear and in the seminorm
sense it is isometric, 1-1 with null space $\left(  S_{w,0;\mathcal{A}%
}^{\prime}\right)  ^{\vee}$.

Also, $\tau_{a}\widetilde{\mathcal{M}}_{1}=\widetilde{\mathcal{M}}_{1}\tau
_{a}$ where $\tau_{a}$ is the translation operator $\tau_{a}u=u\left(
\cdot-a\right)  $.
\end{theorem}

\begin{proof}
That $\widetilde{\mathcal{M}}_{1}$ is an isometry is clear from the definition
of $\widetilde{X}_{1/w}^{0}$. Since $\widetilde{\mathcal{M}}_{1}$ is an
isometry the null space of $\widetilde{\mathcal{M}}_{1}$ is the null space of
the seminorm $\left\vert \cdot\right\vert _{1/w,0}$, namely $\left(
S_{w,0;\mathcal{A}}^{\prime}\right)  ^{\vee}$. Finally%
\[
\tau_{a}\widetilde{\mathcal{M}}_{1}u=\tau_{a}\left(  \frac{u_{F}}{\sqrt{w}%
}\right)  ^{\vee}=\left(  e^{-ia\xi}\frac{u_{F}}{\sqrt{w}}\right)  ^{\vee
}=\left(  \frac{\left(  \tau_{a}u\right)  _{F}}{\sqrt{w}}\right)  ^{\vee
}=\widetilde{\mathcal{M}}_{1}\tau_{a}u.
\]

\end{proof}

\begin{definition}
\label{Def_tildL}\textbf{The operator} $\widetilde{\mathcal{L}}:L^{2}%
\rightarrow\left(  \widehat{S}_{w,0}\right)  ^{\prime}$

If $w\in W_{S;0}$ then Definition \ref{Def_op_L2} of $\mathcal{L}_{2}$
involves first showing that $\sqrt{w}\widehat{g}\in S_{w,0}^{\prime}$ when
$g\in L^{2}$ and then using Lemma \ref{Lem_convex_tls_extend_2} to extend
$\sqrt{w}\widehat{g}$ to $S$ as a member of $S^{\prime}$. If $\left(  \sqrt
{w}\widehat{g}\right)  ^{e}$ is such an extension then $\mathcal{L}_{2}$ was
defined (non-uniquely) by $\mathcal{L}_{2}g=\left(  \left(  \sqrt
{w}\widehat{g}\right)  ^{e}\right)  ^{\vee}$ when $g\in L^{2}$.

We will want $\widetilde{\mathcal{L}}\left(  L^{2}\right)  \subset
\widetilde{X}_{1/w}^{0}\subset\left(  \widehat{S}_{w,0}\right)  ^{\prime}$ so
we can skip the extension step and \textbf{define the linear operator}
$\widetilde{\mathcal{L}}$ \textbf{uniquely} by%
\begin{equation}
\widetilde{\mathcal{L}}g=\left(  \sqrt{w}\widehat{g}\right)  ^{\vee},\quad
g\in L^{2},\label{a2.9}%
\end{equation}

so that $\left(  \sqrt{w}\widehat{g}\right)  ^{\vee}\in\left(  S_{w,0}%
^{\prime}\right)  ^{\vee}=\left(  \widehat{S}_{w,0}\right)  ^{\prime}$.
\end{definition}

The operator $\widetilde{\mathcal{L}}$ has the following properties:

\begin{theorem}
\label{Thm_tildL_property}\textbf{Properties of }$\widetilde{\mathcal{L}}$

\begin{enumerate}
\item If $g\in L^{2}$ then $\left(  \widetilde{\mathcal{L}}g\right)
_{F}=\sqrt{w}\widehat{g}\in L_{loc}^{1}\left(  \mathbb{R}^{d}\setminus
\mathcal{A}\right)  $.

\item $\widetilde{\mathcal{L}}:L^{2}\rightarrow\widetilde{X}_{1/w}^{0}$.

\item $\widetilde{\mathcal{L}}$ is an isometry and 1-1.

\item If $\tau_{a}$ is the translation operator $\tau_{a}u=u\left(
\cdot-a\right)  $ then $\tau_{a}\widetilde{\mathcal{L}}=\widetilde{\mathcal{L}%
}\tau_{a}$.
\end{enumerate}
\end{theorem}

\begin{proof}
\textbf{Part 1} Suppose $g\in L^{2}$. From the definition of
$\widetilde{\mathcal{L}}$, $\widehat{\widetilde{\mathcal{L}}g}=\sqrt
{w}\widehat{g}$. Next observe that because $w\in C^{\left(  0\right)  }\left(
\mathbb{R}^{d}\setminus\mathcal{A}\right)  $ and $C_{0}^{\infty}\left(
\mathbb{R}^{d}\setminus\mathcal{A}\right)  \subset S_{w,0}$ it follows that
\[
\left[  \widehat{\widetilde{\mathcal{L}}g},\phi\right]  =\left[  \sqrt
{w}\widehat{g},\phi\right]  ,\quad\phi\in C_{0}^{\infty}\left(  \mathbb{R}%
^{d}\setminus\mathcal{A}\right)  .
\]

If we can show $\sqrt{w}\widehat{g}\in L_{loc}^{1}\left(  \mathbb{R}%
^{d}\setminus\mathcal{A}\right)  $ it follows that $\left(
\widetilde{\mathcal{L}}g\right)  _{F}=\sqrt{w}\widehat{g}\in L_{loc}%
^{1}\left(  \mathbb{R}^{d}\setminus\mathcal{A}\right)  $. But if
$K\subset\mathbb{R}^{d}\setminus\mathcal{A}$ is compact, $w\in C^{\left(
0\right)  }\left(  \mathbb{R}^{d}\setminus\mathcal{A}\right)  $ implies
\[
\int_{K}\sqrt{w}\left\vert \widehat{g}\right\vert \leq\left(  \int%
_{K}w\right)  ^{1/2}\left\Vert g\right\Vert _{2}\leq\max_{K}\left(  w\right)
\left(  \int_{K}1\right)  \left\Vert g\right\Vert _{2}.
\]

Thus $\widehat{\widetilde{\mathcal{L}}g}\in L_{loc}^{1}\left(  \mathbb{R}%
^{d}\setminus\mathcal{A}\right)  $ and $\left(  \widetilde{\mathcal{L}%
}g\right)  _{F}=\sqrt{w}\widehat{g}$.\medskip

\textbf{Part 2} $\int\frac{\left\vert \left(  \widetilde{\mathcal{L}}g\right)
_{F}\right\vert ^{2}}{w}=\left\Vert \widehat{g}\right\Vert _{2}=\left\Vert
g\right\Vert _{2}$ so $\widetilde{\mathcal{L}}:L^{2}\rightarrow X_{1/w}^{0}%
$.\medskip

\textbf{Part 3} $\widetilde{\mathcal{L}}$ is an isometry since $\left\vert
\widetilde{\mathcal{L}}g\right\vert _{1/w,0}=\left(  \int\frac{\left\vert
\left(  \widetilde{\mathcal{L}}g\right)  _{F}\right\vert ^{2}}{w}\right)
^{1/2}=\left\Vert g\right\Vert _{2}$. Thus $\widetilde{\mathcal{L}}%
:L^{2}\rightarrow X_{1/w}^{0}$ is an isometry and clearly
$\widetilde{\mathcal{L}}g=0$ implies $g=0$.\medskip

\textbf{Part 4}
\[
\left\vert \tau_{a}\widetilde{\mathcal{L}}g-\widetilde{\mathcal{L}}\tau
_{a}g\right\vert _{1/w,0}^{2}=\int\frac{\left\vert \left(  \tau_{a}%
\widetilde{\mathcal{L}}g-\widetilde{\mathcal{L}}\tau_{a}g\right)
_{F}\right\vert ^{2}}{w}=\int\frac{\left\vert \left(  \tau_{a}%
\widetilde{\mathcal{L}}g\right)  _{F}-\left(  \widetilde{\mathcal{L}}\tau
_{a}g\right)  _{F}\right\vert ^{2}}{w}.
\]

But from part 1, $\left(  \widetilde{\mathcal{L}}g\right)  _{F}=\sqrt
{w}\widehat{g}$ on $\mathbb{R}^{d}\setminus\mathcal{A}$. Hence on
$\mathbb{R}^{d}\setminus\mathcal{A}$%
\[
\left(  \widetilde{\mathcal{L}}\tau_{a}g\right)  _{F}=\sqrt{w}\widehat{\tau
_{a}g}=e^{-ia\xi}\sqrt{w}\widehat{g}=e^{-ia\xi}\left(  \widetilde{\mathcal{L}%
}g\right)  _{F},
\]

and since $\left(  \tau_{a}\widetilde{\mathcal{L}}g\right)  ^{\wedge
}=e^{-ia\xi}\left(  \widetilde{\mathcal{L}}g\right)  ^{\wedge}$, it follows
that%
\[
\left(  \tau_{a}\widetilde{\mathcal{L}}g\right)  _{F}=e^{-ia\xi}\left(
\widetilde{\mathcal{L}}g\right)  _{F},
\]

and therefore $\left(  \tau_{a}\widetilde{\mathcal{L}}g\right)  _{F}=\left(
\widetilde{\mathcal{L}}\tau_{a}g\right)  _{F}$.
\end{proof}

The following theorem indicates how the operators $\widetilde{\mathcal{L}%
}:L^{2}\rightarrow\widetilde{X}_{1/w}^{0}$ and $\widetilde{\mathcal{M}}%
_{1}:\widetilde{X}_{1/w}^{0}\rightarrow L^{2}$ interact.

\begin{theorem}
\label{Thm_tildM1_tildL_property}Suppose the weight function $w$ is a member
of $W_{S;0}$. Then:

\begin{enumerate}
\item $\widetilde{\mathcal{M}}_{1}\widetilde{\mathcal{L}}=I$ on $L^{2}$.

\item $\widetilde{\mathcal{L}}\widetilde{\mathcal{M}}_{1}u-u\in\left(
S_{w,0;\mathcal{A}}^{\prime}\right)  ^{\vee}$ when $u\in X_{1/w}^{0}$ i.e.
$\left\vert \widetilde{\mathcal{L}}\widetilde{\mathcal{M}}_{1}u-u\right\vert
_{1/w,0}=0$.

\item $\widetilde{\mathcal{M}}_{1}:\widetilde{X}_{1/w}^{0}\rightarrow L^{2}$
is onto.

\item $\widetilde{\mathcal{L}}:L^{2}\rightarrow\widetilde{X}_{1/w}^{0}$ is
onto in the seminorm sense.

\item $\widetilde{\mathcal{M}}_{1}$ and $\widetilde{\mathcal{L}}$ are adjoints.
\end{enumerate}
\end{theorem}

\begin{proof}
\textbf{Part 1} By definition $\widehat{\widetilde{\mathcal{M}}_{1}u}%
=\frac{u_{F}}{\sqrt{w}}\in L^{2}$ when $u\in\widetilde{X}_{1/w}^{0}$. From
part 1 of Theorem \ref{Thm_tildL_property} $\left(  \widetilde{\mathcal{L}%
}f\right)  _{F}=\sqrt{w}f$. Thus from the definition of
$\widetilde{\mathcal{M}}_{1}$, for $f\in L^{2}$%
\[
\left(  \widetilde{\mathcal{M}}_{1}\widetilde{\mathcal{L}}f\right)  ^{\wedge
}=\frac{\left(  \widetilde{\mathcal{L}}f\right)  _{F}}{\sqrt{w}}=\frac
{\sqrt{w}\widehat{f}}{\sqrt{w}}=\widehat{f},
\]

and so $\widetilde{\mathcal{M}}_{1}\widetilde{\mathcal{L}}=I$ on $L^{2}%
$.\medskip

\textbf{Part 2} If $u\in\widetilde{X}_{1/w}^{0}$ then $\widetilde{\mathcal{L}%
}\widetilde{\mathcal{M}}_{1}u\in\widetilde{X}_{1/w}^{0}$ and
$\widetilde{\mathcal{L}}\widetilde{\mathcal{M}}_{1}u=\sqrt{w}%
\widehat{\widetilde{\mathcal{M}}_{1}u}=u_{F}$. Thus%
\[
\left\vert \widetilde{\mathcal{L}}\widetilde{\mathcal{M}}_{1}u-u\right\vert
_{1/w,0}^{2}=\int\frac{\left\vert \left(  \widetilde{\mathcal{L}%
}\widetilde{\mathcal{M}}_{1}u-u\right)  _{F}\right\vert ^{2}}{w}=0,
\]

and $\widetilde{\mathcal{L}}\widetilde{\mathcal{M}}_{1}u-u\in\left(
S_{w,0;\mathcal{A}}^{\prime}\right)  ^{\vee}$.\medskip

\textbf{Parts 3 and 4} The equation $\widetilde{\mathcal{M}}_{1}%
\widetilde{\mathcal{L}}=I$ implies $\widetilde{\mathcal{M}}_{1}$ is onto and
$\widetilde{\mathcal{L}}\widetilde{\mathcal{M}}_{1}u-u\in\left(
S_{w,0;\mathcal{A}}^{\prime}\right)  ^{\vee}$ when $u\in\widetilde{X}%
_{1/w}^{0}$ implies that $\widetilde{\mathcal{L}}$ are onto.\medskip

\textbf{Part 5} Regarding adjointness, suppose $u\in\widetilde{X}_{1/w}^{0}$
and $g\in L^{2}$. Then from part 1 Theorem \ref{Thm_tildL_property}, $\left(
\widetilde{\mathcal{L}}g\right)  _{F}=\sqrt{w}\widehat{g}$ and from Definition
\ref{Def_tildM1}, $\widehat{\widetilde{\mathcal{M}}_{1}u}=\frac{u_{F}}%
{\sqrt{w}}$. Thus
\[
\left(  \widetilde{\mathcal{L}}g,u\right)  _{1/w,0}=\int\frac{1}{w}\left(
\widetilde{\mathcal{L}}g\right)  _{F}\overline{u_{F}}=\int\frac{1}{w}\sqrt
{w}\widehat{g}\overline{u_{F}}=\int\widehat{g}\frac{\overline{u_{F}}}{\sqrt
{w}}=\int\widehat{g}\overline{\widehat{\widetilde{\mathcal{M}}_{1}u}}=\left(
\widehat{g},\widehat{\widetilde{\mathcal{M}}_{1}u}\right)  _{2}=\left(
g,\widetilde{\mathcal{M}}_{1}u\right)  _{2}.
\]

\end{proof}

Since $L^{2}$ is complete, the mappings of the previous theorem will yield the
following important result:

\begin{corollary}
\label{Cor_tildXo1/w_semiHilb_Sw2}If the weight function $w\in W_{S;0}$ then
in general $\widetilde{X}_{1/w}^{0}$ is a semi-Hilbert space. Indeed,
$\widetilde{X}_{1/w}^{0}$ is a Hilbert space iff $\mathcal{A}$ is empty.
\end{corollary}

\begin{proof}
By Theorem \ref{Thm_tildM1_property}, $\widetilde{\mathcal{M}}_{1}$ is
isometric. Hence if $\left\{  u_{k}\right\}  $ is Cauchy in $\widetilde{X}%
_{1/w}^{0}$ then $\left\{  \widetilde{\mathcal{M}}_{1}u_{k}\right\}  $ is
Cauchy in $L^{2}$ and so $\widetilde{\mathcal{M}}_{1}u_{k}\rightarrow u$ for
some $u\in L^{2}$ since $L^{2}$ is complete. From part 1 Lemma
\ref{Thm_property1_tildXo1/w} the seminorm for $\widetilde{X}_{1/w}^{0}$ has
null space $\left(  S_{w,0;\mathcal{A}}^{\prime}\right)  ^{\vee}$. Hence by
Theorem \ref{Thm_tildM1_tildL_property}, $\widetilde{\mathcal{L}%
}\widetilde{\mathcal{M}}_{1}u_{k}=u_{k}\rightarrow\widetilde{\mathcal{L}}%
u\in\widetilde{X}_{1/w}^{0}$ and $\widetilde{X}_{1/w}^{0}$ is complete.

Finally, from part 2 of Theorem \ref{Thm_property1_tildXo1/w}, $\widetilde{X}%
_{1/w}^{0}$ is an inner product space iff $\mathcal{A}$ is empty.
\end{proof}

\subsection{The map $\protect\widetilde{\mathcal{B}}_{1}:X_{w}^{0}%
\otimes\protect\widetilde{X}_{1/w}^{0}\rightarrow C_{B}^{\left(  0\right)  }$}

We now define the analogue $\widetilde{\mathcal{B}}_{1}$ of the bilinear
operator $\mathcal{B}$ which was introduced in Definition \ref{Def_map_B}:

\begin{definition}
\label{Def_tildB1}\textbf{The mapping }$\widetilde{\mathcal{B}}_{1}$

If $\left(  u,v\right)  \in X_{w}^{0}\otimes\widetilde{X}_{1/w}^{0}$ define
the bilinear mapping $\widetilde{\mathcal{B}}_{1}:X_{w}^{0}\otimes
\widetilde{X}_{1/w}^{0}\rightarrow C_{B}^{\left(  0\right)  }$ by:%
\begin{equation}
\widetilde{\mathcal{B}}_{1}\left(  u,v\right)  =\left(  \widehat{u}%
v_{F}\right)  ^{\vee},\quad u\in X_{w}^{0},\text{ }v\in\widetilde{X}_{1/w}%
^{0}.\label{a1.201}%
\end{equation}

Noting that the product of two measurable functions is a measurable function
definition \ref{a1.201} makes sense since from the definitions of
$\widetilde{X}_{1/w}^{0}$ and $X_{w}^{0}$, $v_{F}\in L_{loc}^{1}\left(
\mathbb{R}^{d}\setminus\mathcal{A}\right)  $, $\frac{v_{F}}{\sqrt{w}}\in
L^{2}$, $\widehat{u}\in L_{loc}^{1}$ and $\sqrt{w}\widehat{u}\in L^{2}$ so
that $\int\left\vert \widehat{u}v_{F}\right\vert \leq\left\Vert u\right\Vert
_{w,0}\left\vert v\right\vert _{1/w,0}<\infty$. Thus $\widehat{u}v_{F}\in
L^{1}$ and hence $\left(  \widehat{u}v_{F}\right)  ^{\vee}\in C_{B}^{\left(
0\right)  }$.
\end{definition}

\begin{theorem}
\label{Thm_property_tildB1}\textbf{Properties of the operator }%
$\widetilde{\mathcal{B}}_{1}$

\begin{enumerate}
\item $\widetilde{\mathcal{B}}_{1}:X_{w}^{0}\otimes\widetilde{X}_{1/w}%
^{0}\rightarrow C_{B}^{\left(  0\right)  }$ is a continuous bilinear mapping
when $C_{B}^{\left(  0\right)  }$ is endowed with the supremum norm
$\left\Vert \cdot\right\Vert _{\infty}$. In fact%
\begin{equation}
\left\Vert \widetilde{\mathcal{B}}_{1}\left(  u,v\right)  \right\Vert
_{\infty}\leq\left(  2\pi\right)  ^{-\frac{d}{2}}\left\Vert u\right\Vert
_{w,0}\left\vert v\right\vert _{1/w,0},\quad u\in X_{w}^{0},\text{ }%
v\in\widetilde{X}_{1/w}^{0}.\label{a3.0}%
\end{equation}

\item In the sense of distributions: if $\alpha+\beta=\gamma\geq0$ then
\[
D^{\gamma}\widetilde{\mathcal{B}}_{1}\left(  u,v\right)  =\left(
\widehat{D^{\alpha}u}\left(  D^{\beta}v\right)  _{F}\right)  ^{\vee},\quad
u\in X_{w}^{0},\text{ }v\in\widetilde{X}_{1/w}^{0}.
\]

\item The operator $\widetilde{\mathcal{B}}_{1}$ commutes with the
(distribution) translation operator $\tau_{a}f=f\left(  \cdot-a\right)  $,
$a\in\mathbb{R}^{d}$ in the sense that%
\[
\tau_{a}\widetilde{\mathcal{B}}_{1}\left(  u,v\right)  =\widetilde{\mathcal{B}%
}_{1}\left(  \tau_{a}u,v\right)  =\widetilde{\mathcal{B}}_{1}\left(
u,\tau_{a}v\right)  .
\]

\item We have%
\[
\widetilde{\mathcal{B}}_{1}\left(  u,v\right)  =\left(  \mathcal{I}u\right)
\ast\widetilde{\mathcal{M}}_{1}v,\quad u\in X_{w}^{0},\text{ }v\in
\widetilde{X}_{1/w}^{0},
\]

where $\mathcal{I}$ was introduced in Definition \ref{Def_I_J} and
$\widetilde{\mathcal{M}}_{1}$ in Definition \ref{Def_tildM1}.\medskip

Now suppose $w$ also satisfies property W02 or W03 for parameter $\kappa$.
Then:\medskip

\item $\widetilde{\mathcal{B}}_{1}$ is a convolution in the sense of
\ref{1.20}: $X_{w}^{0}\subset S^{\prime}$, $S\subset\widetilde{X}_{1/w}^{0}$
and%
\[
\widetilde{\mathcal{B}}_{1}\left(  u,v\right)  =u\ast v,\quad u\in X_{w}%
^{0},\text{ }v\in S.
\]

Further, $\widetilde{\mathcal{B}}_{1}\left(  u,v\right)  \in C_{B}^{\infty}$
and $D^{\gamma}\widetilde{\mathcal{B}}_{1}\left(  u,v\right)  =D^{\alpha}u\ast
D^{\beta}v$ for all $\gamma=\alpha+\beta$.
\end{enumerate}
\end{theorem}

\begin{proof}
\textbf{Part 1} Noting the calculations done in the definition of
$\widetilde{\mathcal{B}}_{1}$ all that remains to be shown is the continuity
of $\widetilde{\mathcal{B}}_{1}$. But from \ref{a1.201}
\begin{align*}
\left\vert \widetilde{\mathcal{B}}_{1}\left(  u,v\right)  \right\vert
\leq\left\vert \left(  \widehat{u}v_{F}\right)  ^{\vee}\right\vert \leq\left(
2\pi\right)  ^{-\frac{d}{2}}\left\vert \int e^{-ix\xi}\widehat{u}\left(
\xi\right)  v_{F}\left(  \xi\right)  d\xi\right\vert  &  \leq\left(
2\pi\right)  ^{-\frac{d}{2}}\int\left\vert \widehat{u}\right\vert \left\vert
v_{F}\right\vert \\
&  \leq\left(  2\pi\right)  ^{-\frac{d}{2}}\left\Vert u\right\Vert
_{w,0}\left\vert v\right\vert _{1/w,0}.
\end{align*}
\medskip

\textbf{Part 2}%
\[
D^{\gamma}\widetilde{\mathcal{B}}_{1}\left(  u,v\right)  =D^{\gamma}\left(
\widehat{u}v_{F}\right)  ^{\vee}=\left(  \left(  i\xi\right)  ^{\gamma}\left(
\widehat{u}v_{F}\right)  \right)  ^{\vee}=\left(  \left(  i\xi\right)
^{\alpha}\widehat{u}\text{ }\left(  i\xi\right)  ^{\beta}v_{F}\right)  ^{\vee
}=\left(  \widehat{D^{\alpha}u}\text{ }\left(  D^{\beta}v\right)  _{F}\right)
^{\vee}.
\]
\medskip

\textbf{Part 3}%
\[
\tau_{a}\widetilde{\mathcal{B}}_{1}\left(  u,v\right)  =\tau_{a}\left(
\widehat{u}v_{F}\right)  ^{\vee}=\left(  e^{-ia\xi}\widehat{u}v_{F}\right)
^{\vee}=\left(  \widehat{\tau_{a}u}v_{F}\right)  ^{\vee}%
=\widetilde{\mathcal{B}}_{1}\left(  \tau_{a}u,v\right)  ,
\]

and%
\[
\tau_{a}\widetilde{\mathcal{B}}_{1}\left(  u,v\right)  =\tau_{a}\left(
\widehat{u}v_{F}\right)  ^{\vee}=\left(  \widehat{u}e^{-ia\xi}v_{F}\right)
^{\vee}=\left(  \widehat{u}\text{ }\left(  \tau_{a}v\right)  _{F}\right)
^{\vee}=\widetilde{\mathcal{B}}_{1}\left(  u,\tau_{a}v\right)  .
\]
\medskip

\textbf{Part 4} From the definition of $\widetilde{\mathcal{B}}_{1}$ and the
definition of $\widetilde{\mathcal{M}}_{1}$ (Definition \ref{Def_tildM1})%
\begin{align*}
\widetilde{\mathcal{B}}_{1}\left(  u,v\right)  =\left(  \widehat{u}%
v_{F}\right)  ^{\vee}=\left(  \sqrt{w}\widehat{u}\frac{v_{F}}{\sqrt{w}%
}\right)  ^{\vee} &  =\left(  \sqrt{w}\widehat{u}\right)  ^{\vee}\ast\left(
\frac{v_{F}}{\sqrt{w}}\right)  ^{\vee}\\
&  =\left(  \mathcal{I}u\right)  \ast\widetilde{\mathcal{M}}_{1}v.
\end{align*}
\medskip

\textbf{Part 5} From part 4 of Lemma \ref{Lem_Xoinvw}, $S\subset
\widetilde{X}_{1/w}^{0}$. So by definition \ref{a1.201} of
$\widetilde{\mathcal{B}}_{1}\left(  u,v\right)  $, and the definition
\ref{1.20} of the convolution of a member of $S^{\prime}$ and a member of $S$,
$\widetilde{\mathcal{B}}_{1}\left(  u,v\right)  =\left(  \widehat{u}%
v_{F}\right)  ^{\vee}=\left(  \widehat{u}\widehat{v}\right)  ^{\vee}=u\ast
v\in C_{B}^{\left(  0\right)  }$ so that by \ref{2.40}, $D^{\gamma
}\widetilde{\mathcal{B}}_{1}\left(  u,v\right)  =\left(  D^{\alpha}u\right)
\ast D^{\beta}v$ for all $\gamma$ and $\alpha+\beta=\gamma$. Thus
$\widetilde{\mathcal{B}}_{1}\left(  u,v\right)  \in C_{B}^{\infty}$.
\end{proof}

Part 5 of the last theorem shows how $\widetilde{\mathcal{B}}_{1}$ can be
regarded as a convolution by restricting the domain of the second variable to
$S$.

The next result considers the important case where the first argument of
$\widetilde{\mathcal{B}}_{1}$ is the basis function.

\begin{corollary}
\label{Cor_property_tildB1_basis}\textbf{Properties of }%
$\widetilde{\mathcal{B}}_{1}\left(  G\mathbf{,\cdot}\right)  $. Suppose $w$
has property W02 or W03 for $\kappa$ and $G$ is the basis function. Then if
$v,v^{\prime}\in X_{1/w}^{0}$:

\begin{enumerate}
\item $D^{\gamma}G\in X_{w}^{0}$ when $\left\vert \gamma\right\vert \leq
\kappa$ and%
\[
D^{\gamma}\widetilde{\mathcal{B}}_{1}\left(  G,v\right)
=\widetilde{\mathcal{B}}_{1}\left(  D^{\gamma}G,v\right)  ,\quad\left\vert
\gamma\right\vert \leq\kappa.
\]

\item $\widetilde{\mathcal{B}}_{1}\left(  G,v\right)  \in X_{w}^{0}$ and
$\left\Vert \widetilde{\mathcal{B}}_{1}\left(  G,v\right)  \right\Vert
_{w,0}=\left\vert v\right\vert _{1/w,0}$.

\item $\widetilde{\mathcal{B}}_{1}\left(  G,\phi\right)  =G\ast\phi$ when
$\phi\in S$.

\item $\int F\widetilde{\mathcal{B}}_{1}\left(  \widetilde{\mathcal{B}}%
_{1}\left(  G,v\right)  ,\overline{v^{\prime}}\right)  =\left\langle
v,v_{\ast}^{\prime}\right\rangle _{1/w,0}$, where $F$ is the Fourier transform
and $v_{\ast}^{\prime}\left(  x\right)  =v^{\prime}\left(  -x\right)  $.
\end{enumerate}
\end{corollary}

\begin{proof}
\textbf{Part 1} By Theorem \ref{Thm_basis_fn_properties_all_m_W2}, $D^{\gamma
}G\in X_{w}^{0}$ when $\left\vert \gamma\right\vert \leq\kappa$ and so by part
2 of Theorem \ref{Thm_property_tildB1}, $D^{\gamma}\widetilde{\mathcal{B}}%
_{1}\left(  G,v\right)  \in\left(  \widehat{D^{\gamma}G}v_{F}\right)  ^{\vee
}=\widetilde{\mathcal{B}}_{1}\left(  D^{\gamma}G,v\right)  $.\medskip

\textbf{Part 2} We have $G\in X_{w}^{0}$ and $\widetilde{\mathcal{B}}%
_{1}\left(  G,v\right)  \in S^{\prime}$. Further, $F\widetilde{\mathcal{B}%
}_{1}\left(  G,v\right)  =\widehat{G}v_{F}=\frac{v_{F}}{w}\in L^{1}$, by
Theorem \ref{Thm_property1_tildXo1/w} and

$\int w\left\vert F\widetilde{\mathcal{B}}_{1}\left(  G,v\right)  \right\vert
^{2}=\int\frac{\left\vert v_{F}\right\vert ^{2}}{w}=\left\vert v\right\vert
_{1/w,0}$. Thus $\widetilde{\mathcal{B}}_{1}\left(  G,v\right)  \in X_{w}^{0}%
$.\medskip

\textbf{Part 3} From part 4 Theorem \ref{Thm_property1_tildXo1/w}, $S\subset
X_{1/w}^{0}$ and $\widetilde{\mathcal{B}}_{1}\left(  G,\phi\right)  =\left(
\widehat{G}\widehat{\phi}\right)  ^{\vee}=G\ast\phi$ by the convolution
definition \ref{1.20}.\medskip

\textbf{Part 4} $\widetilde{\mathcal{B}}_{1}\left(  G,v\right)  =\left(
\frac{v_{F}}{w}\right)  ^{\vee}$ so
\begin{align*}
\int F\widetilde{\mathcal{B}}_{1}\left(  \widetilde{\mathcal{B}}_{1}\left(
G,v\right)  ,\overline{v^{\prime}}\right)  =\int F\widetilde{\mathcal{B}}%
_{1}\left(  \left(  \frac{v_{F}}{w}\right)  ^{\vee},\overline{v^{\prime}%
}\right)  =\int F\left(  \left(  \frac{v_{F}}{w},\left(  \overline{v^{\prime}%
}\right)  _{F}\right)  ^{\vee}\right)   &  =\int\frac{1}{w}v_{F}\left(
\overline{v^{\prime}}\right)  _{F}\\
&  =\int\frac{1}{w}v_{F}\overline{\left(  v_{\ast}^{\prime}\right)  _{F}}\\
&  =\left\langle v,v_{\ast}^{\prime}\right\rangle _{1/w,0}.
\end{align*}

\end{proof}

The next theorem is the analogue of part 3 of Theorem \ref{Thm_properties2_B}.

\begin{theorem}
\label{Thm_tildB1_invFSw2_convol}Suppose the weight function $w$ has property
W01 w.r.t. the set $\mathcal{A}$ and $w\in W_{S;0}$ (Definition
\ref{Def_Sw2_and_fnal(Sw2)}). Then if $\phi\in\overset{\vee}{S}_{w,0}\subset
X_{w}^{0}$ and $v\in\widetilde{X}_{1/w}^{0}\subset\left(  \widehat{S}%
_{w,0}\right)  ^{\prime}$:%
\[
\widetilde{\mathcal{B}}_{1}\left(  \phi,v\right)  =\phi\ast\overset{\vee
}{v_{F}}=\left(  2\pi\right)  ^{-\frac{d}{2}}\left[  \left(  \overset{\vee
}{v_{F}}\right)  _{y},\phi\left(  \cdot-y\right)  \right]  ,
\]

where the convolution is that of Definition \ref{Def_Sw2_convol_fnalSw2}.
\end{theorem}

\begin{proof}
From \ref{a1.201}, $\widetilde{\mathcal{B}}_{1}\left(  \phi,v\right)  =\left(
\widehat{\phi}v_{F}\right)  ^{\vee}=\left(  \widehat{\phi}%
\widehat{\overset{\vee}{v_{F}}}\right)  ^{\vee}$ and by part 5 Theorem
\ref{Thm_property1_tildXo1/w}, $\overset{\vee}{v_{F}}\in\left(  \widehat{S}%
_{w,0}\right)  ^{\prime}$ so that Lemma \ref{Lem_convol_Sw2} implies,
$\widetilde{\mathcal{B}}_{1}\left(  \phi,v\right)  =\phi\ast\overset{\vee
}{v_{F}}=\left(  2\pi\right)  ^{-d/2}\left[  \left(  \overset{\vee}{v_{F}%
}\right)  _{y},\phi\left(  \cdot-y\right)  \right]  $.
\end{proof}

\subsection{The maps $\protect\widetilde{\mathcal{W}}_{1}:X_{w}^{0}%
\rightarrow\protect\widetilde{X}_{1/w}^{0}$, $\protect\widetilde{\mathcal{V}%
}_{1}:\protect\widetilde{X}_{1/w}^{0}\rightarrow X_{w}^{0}$ and
$\protect\widetilde{\Phi}_{1}:\protect\widetilde{X}_{1/w}^{0}\rightarrow
\left(  X_{w}^{0}\right)  ^{\prime}$}

The operators $\widetilde{\mathcal{W}}_{1}=\widetilde{\mathcal{L}}\mathcal{I}
$ and $\widetilde{\mathcal{V}}_{1}=\mathcal{J}\widetilde{\mathcal{M}}_{1}$
will be introduced - see Figure \ref{Fig_ops_Xo1/w_tildXo1/w_Xow_Sw2} - and
used to define the operator $\widetilde{\Phi}_{1}:\widetilde{X}_{1/w}%
^{0}\rightarrow\left(  X_{w}^{0}\right)  ^{\prime}$ which characterizes the
bounded linear functionals on $X_{w}^{0}$ as the space $\widetilde{X}%
_{1/w}^{0}$. In fact, in an analogous fashion to the negative order Sobolev
spaces, Theorem \ref{Thm_op_tildPhi1_bilinear} shows that $\widetilde{\Phi
}_{1}$ can be written as a bilinear form.

\begin{definition}
\label{Def_op_tildV1}\textbf{The operator} $\widetilde{\mathcal{V}}%
_{1}:\widetilde{X}_{1/w}^{0}\rightarrow X_{w}^{0}$

The operator $\mathcal{J}:L^{2}\rightarrow X_{w}^{0}$ was introduced in
Definition \ref{Def_I_J} and studied in Theorem \ref{Thm_I_J_property}.

Suppose the weight function $w$ has property W02 or W03. These conditions
ensure $\centerdot$that both $\widetilde{\mathcal{M}}_{1}$ and $\mathcal{J}$
are defined.

\textbf{Define the operator} $\widetilde{\mathcal{V}}_{1}$ by
$\widetilde{\mathcal{V}}_{1}=\mathcal{J}\widetilde{\mathcal{M}}_{1}$ where
$\widetilde{\mathcal{M}}_{1}:\widetilde{X}_{1/w}^{0}\rightarrow L^{2}$.
\end{definition}

The operator $\widetilde{\mathcal{V}}_{1}$ can be expressed in terms of the
operator $\widetilde{\mathcal{B}}_{1}$ and the basis function as follows:

\begin{theorem}
\label{Thm_tildV1u_=_G*u}Suppose the weight function $w$ has property W02 or
W03 and that $G$ is the corresponding basis function. Then
\begin{equation}
\widetilde{\mathcal{V}}_{1}v=\left(  \frac{v_{F}}{w}\right)  ^{\vee
}=\widetilde{\mathcal{B}}_{1}\left(  G,v\right)  ,\quad v\in\widetilde{X}%
_{1/w}^{0},\label{a6.1}%
\end{equation}

and%
\[
\widetilde{\mathcal{V}}_{1}\phi=G\ast\phi,\quad\phi\in S,
\]

and $\widetilde{\mathcal{V}}_{1}:S\rightarrow G\ast S$ is 1-1 and onto.
\end{theorem}

\begin{proof}
By Definition \ref{Def_I_J}, $\widehat{\mathcal{J}g}=\frac{\widehat{g}}%
{\sqrt{w}}\in L_{loc}^{1}$ for $g\in L^{2}$.

By Definition \ref{Def_tildM1}, $\widehat{\widetilde{\mathcal{M}}_{1}v}%
=\frac{v_{F}}{\sqrt{w}}\in L^{2}$ for $v\in\widetilde{X}_{1/w}^{0}$. Thus
\[
\widehat{\widetilde{\mathcal{V}}_{1}v}=\widehat{\mathcal{J}%
\widetilde{\mathcal{M}}_{1}v}=\frac{\widehat{\widetilde{\mathcal{M}}_{1}v}%
}{\sqrt{w}}=\frac{v_{F}}{w}=\left(  \widetilde{\mathcal{B}}_{1}\left(
G,v\right)  \right)  ^{\wedge},
\]

and so $\widetilde{\mathcal{V}}_{1}v=\widetilde{\mathcal{B}}_{1}\left(
G,v\right)  $. The second equation of this theorem follows from part 3 of
Corollary \ref{Cor_property_tildB1_basis}.

Clearly $\widetilde{\mathcal{V}}_{1}$ is onto. Since $\widetilde{\mathcal{V}%
}_{1}$ is isometric $\widetilde{\mathcal{V}}_{1}\phi=0$ implies $\left\Vert
\widetilde{\mathcal{V}}_{1}\phi\right\Vert _{w,0}=\left\vert \phi\right\vert
_{1/w,0}=\left(  \int\frac{\left\vert _{\phi_{F}}\right\vert ^{2}}{w}\right)
^{1/2}=0$, where $\phi_{F}$ is the restriction of $\widehat{\phi}$ to
$\mathbb{R}^{d}\setminus\mathcal{A}$ where $\mathcal{A}$ is a closed set of
measure zero. But $\widehat{\phi}\in S$ and so $\phi=0$.
\end{proof}

\begin{definition}
\label{Def_op_tildW1}\textbf{The operator }$\widetilde{\mathcal{W}}%
_{1}=\widetilde{\mathcal{L}}\mathcal{I}$
\end{definition}

The next theorem relates the properties of the operators
$\widetilde{\mathcal{V}}_{1}$ and $\widetilde{\mathcal{W}}_{1}$, and for which
we will require the properties of the operators $\mathcal{I}$ and
$\mathcal{J}$ given in Theorem \ref{Thm_I_J_property}. The next result is an
easy analogue of Theorem \ref{Thm_V_W_properties}.

\begin{theorem}
\label{Thm_tildV1_tildW1_property}Suppose the weight function $w$ has property
W02 or W03 and $w\in W_{S;0}$. Then in the seminorm sense the operators
$\widetilde{\mathcal{V}}_{1}=\mathcal{J}\widetilde{\mathcal{M}}_{1} $ and
$\widetilde{\mathcal{L}}\mathcal{I}$ have the following properties:

\begin{enumerate}
\item $\widetilde{\mathcal{V}}_{1}:\widetilde{X}_{1/w}^{0}\rightarrow
X_{w}^{0}$ is a unique, linear isometry.

\item $\widetilde{\mathcal{W}}_{1}:X_{w}^{0}\rightarrow\widetilde{X}_{1/w}%
^{0}$ is a class of linear isometries.

\item $\widetilde{\mathcal{V}}_{1}$ and $\widetilde{\mathcal{W}}_{1}$ are inverses.

\item $\widetilde{\mathcal{V}}_{1}$ and $\widetilde{\mathcal{W}}_{1}$ are onto.

\item $\widetilde{\mathcal{V}}_{1}$ and $\widetilde{\mathcal{W}}_{1}$ are 1-1.

\item $\widetilde{\mathcal{V}}_{1}$ and $\widetilde{\mathcal{W}}_{1}$ are adjoints.

\item If $f\in X_{w}^{0}$ then $\widetilde{\mathcal{W}}_{1}f=\left(
w\widehat{f}\right)  ^{\vee}$.

\item $\widetilde{\mathcal{V}}_{1}$ and $\widetilde{\mathcal{W}}_{1}$ are
isometric isomorphisms, inverses and adjoints.
\end{enumerate}
\end{theorem}

\begin{proof}
The proof is almost identical to that of Theorem ??.\medskip

\textbf{Part 1}\ From Theorem \ref{Thm_I_J_property} $\mathcal{J}$ is a
unique, linear isometry, and from Theorem \ref{Thm_tildM1_property},
$\widetilde{\mathcal{M}}_{1}$ is also a unique, linear isometry.\medskip

\textbf{Part 2} From Theorem \ref{Thm_I_J_property} $\mathcal{I}$ is a unique,
linear isometry and from Theorem \ref{Thm_tildL_property}
$\widetilde{\mathcal{L}}$ is a linear isometry.\medskip

\textbf{Part 3} From Theorem \ref{Thm_tildM1_tildL_property}
$\widetilde{\mathcal{M}}_{1}\widetilde{\mathcal{L}}=I$\ and Theorem
\ref{Thm_I_J_property} $\mathcal{JI}=I$. Hence $\widetilde{\mathcal{V}}%
_{1}\widetilde{\mathcal{W}}_{1}=\mathcal{J}\widetilde{\mathcal{M}}%
_{1}\widetilde{\mathcal{L}}\mathcal{I}=\mathcal{JI}$ and so
$\widetilde{\mathcal{V}}_{1}\widetilde{\mathcal{L}}\mathcal{I}=I$.

From Theorem \ref{Thm_I_J_property} $\mathcal{I}\mathcal{J}=I$ and from
Theorem \ref{Thm_tildM1_tildL_property} $\widetilde{\mathcal{L}}%
\widetilde{\mathcal{M}}_{1}g-g\in\left(  S_{w,0;\mathcal{A}}^{\prime}\right)
^{\vee}$. Hence $\widetilde{\mathcal{W}}_{1}\widetilde{\mathcal{V}}%
_{1}=\widetilde{\mathcal{L}}\mathcal{IJ}\widetilde{\mathcal{M}}_{1}%
=\widetilde{\mathcal{L}}\widetilde{\mathcal{M}}_{1}$ and so
$\widetilde{\mathcal{W}}_{1}\widetilde{\mathcal{V}}_{1}g-g\in\left(
S_{w,0;\mathcal{A}}^{\prime}\right)  ^{\vee}$.\medskip

\textbf{Parts 4 and 5} follow directly from part 3.\medskip

\textbf{Part 6} Theorem \ref{Thm_tildM1_tildL_property}\ implies
$\widetilde{\mathcal{L}}$ and $\widetilde{\mathcal{M}}_{1}$ are adjoints and
Theorem \ref{Thm_I_J_property} implies $\mathcal{I}$ and $\mathcal{J}$ are
adjoints. Thus
\[
\left\langle g,\widetilde{\mathcal{W}}_{1}f\right\rangle _{1/w,0}=\left\langle
g,\widetilde{\mathcal{L}}\mathcal{I}f\right\rangle _{1/w,0}=\left(
\widetilde{\mathcal{M}}_{1}g,\mathcal{I}f\right)  _{2}=\left(  \mathcal{J}%
\widetilde{\mathcal{M}}_{1}g,f\right)  _{w,0}=\left(  \widetilde{\mathcal{V}%
}_{1}g,f\right)  _{w,0}.
\]
\medskip

\textbf{Part 7} If $f\in X_{w}^{0}$ then from the definition of
$\widetilde{\mathcal{L}}$
\[
\widehat{\widetilde{\mathcal{W}}_{1}f}=\widehat{\widetilde{\mathcal{L}%
}\mathcal{I}f}=\sqrt{w}\widehat{\mathcal{I}f}=\sqrt{w}\left(  \sqrt
{w}\widehat{f}\right)  =w\widehat{f}.
\]

\end{proof}

\begin{corollary}
\label{Cor_invF[wSw2]_dense_tildXo1/w}\textbf{The mapping }$\left(
w\widehat{\phi}\right)  ^{\vee}$ \textbf{where} $\phi\in\overset{\vee
}{S}_{w,0}$\textbf{.}

Suppose the weight function $w$ has property W01 w.r.t. the set $\mathcal{A}$
and that $w\in W_{S;0}$. Then:

\begin{enumerate}
\item The mapping $\left(  w\widehat{\phi}\right)  ^{\vee}$ is an isometric,
isomorphism from $\overset{\vee}{S}_{w,0}=S\cap X_{w}^{0}$ to $\left(
wS_{w,0}\right)  ^{\vee}\subset\widetilde{X}_{1/w}^{0}$.\medskip

If, in addition, $w$ has property W02 or W03 then:\medskip

\item $\left(  wS_{w,0}\right)  ^{\vee}$ is dense in $\widetilde{X}_{1/w}^{0}$.

\item The mapping $\left(  w\widehat{\phi}\right)  ^{\vee}$ can be extended to
an isometric, isomorphism (in the seminorm sense).

\item The isometric isomorphism of part 3 is actually $\widetilde{\mathcal{W}%
}_{1}$: $\widetilde{\mathcal{W}}_{1}=\left(  w\widehat{\phi}\right)  ^{\vee}$.
\end{enumerate}
\end{corollary}

\begin{proof}
\textbf{Part 1} This is a direct consequence of part 2 of Theorem
\ref{Thm_w_phi_in_fnalSw2} and the fact that \ref{Thm_properties2_B} shows
$\overset{\vee}{S}_{w,0}=S\cap X_{w}^{0}$ is dense in $X_{w}^{0}%
$.\textbf{\medskip}

\textbf{Part 2} Follows from part 1 and the fact, proved in Theorem
\ref{Thm_properties2_B}, that $\overset{\vee}{S}_{w,0}$ is dense in $X_{w}%
^{0}$\textbf{.\medskip}

\textbf{Part 3} Follows from parts 1 and 2 using density
arguments.\textbf{\medskip}

\textbf{Part 4} Follows directly from part 7 of Theorem
\ref{Thm_tildV1_tildW1_property}.
\end{proof}

Next we prove that $S$ is dense in $\widetilde{X}_{1/w}^{0}$. This result is
the analogue of Corollary \ref{Cor_S_dense_Xo1/w}.

\begin{corollary}
\label{Cor_S_dense_tildXo1/w}Suppose the weight function $w$ has property W02
or W03 and that $w\in W_{S;0}$.

Then $S$ is dense in $\widetilde{X}_{1/w}^{0}$ under $\left\vert
\mathcal{\cdot}\right\vert _{1/w,0}$.
\end{corollary}

\begin{proof}
Choose $u\in\widetilde{X}_{1/w}^{0}$ so that $\widetilde{\mathcal{V}}_{1}u\in
X_{w}^{0}$. From part 3 Theorem \ref{Thm_Jg_properties}, $G\ast S$ is dense in
$X_{w}^{0}$ and hence given $\varepsilon>0$ there exists $\phi_{\varepsilon
}\in G\ast S$ such that $\left\Vert \widetilde{\mathcal{V}}_{1}u-\phi
_{\varepsilon}\right\Vert _{w,0}<\varepsilon$. But by Theorem
\ref{Thm_tildV1u_=_G*u}, $\widetilde{\mathcal{V}}_{1}:S\rightarrow G\ast S$ is
an isomorphism and so $\widetilde{\mathcal{V}}_{1}^{-1}\phi_{\varepsilon}\in
S$. Since $\widetilde{\mathcal{V}}_{1}:\widetilde{X}_{1/w}^{0}\rightarrow
X_{w}^{0}$ is an isometry%
\[
\left\vert u-\widetilde{\mathcal{V}}_{1}^{-1}\phi_{\varepsilon}\right\vert
_{1/w,0}=\left\Vert \widetilde{\mathcal{V}}_{1}\left(
u-\widetilde{\mathcal{V}}_{1}^{-1}\phi_{\varepsilon}\right)  \right\Vert
_{w,0}=\left\Vert \widetilde{\mathcal{V}}_{1}u-\phi_{\varepsilon}\right\Vert
_{w,0}<\varepsilon,
\]

which proves this density result.

Since $\widetilde{X}_{1/w}^{0}$ is complete and $S$ is dense in $\widetilde{X}%
_{1/w}^{0}$ it follows that $\widetilde{X}_{1/w}^{0}$ is the completion of $S$.
\end{proof}

In the following analogue of Theorem \ref{Thm_op_Phi} the operator
$\widetilde{\mathcal{V}}_{1}$ will now be used characterize the bounded linear
functionals on $X_{w}^{0}$, denoted $\left(  X_{w}^{0}\right)  ^{\prime}$, as
the members of $\widetilde{X}_{1/w}^{0}$:

\begin{theorem}
\label{Thm_op_tildPhi1}\textbf{The operator }$\widetilde{\Phi}_{1}$. Suppose
the weight function $w$ has property W02 or W03 and that $w\in W_{S;0}$.
Denote by $\left(  X_{w}^{0}\right)  ^{\prime}$ the space of bounded linear
functionals on $X_{w}^{0}$.

Then the equation%
\begin{equation}
\left(  \widetilde{\Phi}_{1}v\right)  \left(  u\right)  =\left(
u,\widetilde{\mathcal{V}}_{1}v\right)  _{w,0},\quad u\in X_{w}^{0},\text{
}v\in\widetilde{X}_{1/w}^{0},\label{a1.61}%
\end{equation}

defines a linear operator $\widetilde{\Phi}_{1}:\widetilde{X}_{1/w}%
^{0}\rightarrow\left(  X_{w}^{0}\right)  ^{\prime}$ which is an isometric
isomorphism in the seminorm sense.
\end{theorem}

\begin{proof}
This proof is "identical" to that of Theorem \ref{Thm_op_Phi} and so will be omitted.
\end{proof}

In an analogous fashion to the negative order Sobolev spaces a bilinear form
can be used to characterize the bounded linear functionals on $X_{w}^{0}$:

\begin{theorem}
\label{Thm_op_tildPhi1_bilinear}If $u\in X_{w}^{0}$ and $v\in\widetilde{X}%
_{1/w}^{0}$ then $\widetilde{\Phi}_{1}$ can be expressed directly in terms of
the bilinear form $\int\widehat{u}\overline{v_{F}}$ as%
\[
\left(  \widetilde{\Phi}_{1}v\right)  \left(  u\right)  =\int\widehat{u}%
\overline{v_{F}},\quad u\in X_{w}^{0},\text{ }v\in\widetilde{X}_{1/w}^{0}.
\]

\end{theorem}

\begin{proof}
A direct consequence of \ref{a1.61} and \ref{a6.1}.
\end{proof}

\subsection{Extension and restriction mappings between $\protect\widetilde{X}%
_{1/w}^{0}$ and $X_{1/w}^{0}$\label{SbSect_E_R_tildXo1/w_Xo1/w}}

In this subsection we will still assume that $w\in S_{w,0}^{\dag}$ and we will
relate the space $\widetilde{X}_{1/w}^{0}$ to the space $X_{1/w}^{0}$ by
constructing an isometric isomorphism (in the seminorm sense) between them, as
indicated in Figure \ref{Fig_ops_E_and_R}.%

\begin{equation}%
\begin{array}
[c]{ccccccc}%
X_{1/w}^{0} &
\begin{array}
[c]{c}%
\overset{E}{\longleftarrow}\\
\underset{R}{\longrightarrow}%
\end{array}
& \widetilde{X}_{1/w}^{0} &
\begin{array}
[c]{c}%
\overset{\widetilde{\mathcal{L}}}{\longleftarrow}\\
\underset{\widetilde{\mathcal{M}}_{1}}{\longrightarrow}%
\end{array}
& L^{2} &
\begin{array}
[c]{c}%
\overset{\mathcal{I}}{\longleftarrow}\\
\underset{\mathcal{J}}{\longrightarrow}%
\end{array}
& X_{w}^{0}%
\end{array}
\label{Fig_ops_E_and_R}%
\end{equation}

This isometry involves the null space $\left(  \widehat{S}_{w,0}\right)
^{\bot}$, whereas the isometries $\mathcal{L}$ and $\mathcal{M}$ between
$\widetilde{X}_{1/w}^{0}$ and $L^{2}$ involve the null space $\left(
S_{w,0;\mathcal{A}}^{\prime}\right)  ^{\vee}$. This separation simplifies the calculations.

Recall the definitions \ref{a2.5}:
\[
\widetilde{X}_{1/w}^{0}=\left\{  u\in\left(  \widehat{S}_{w,0}\right)
^{\prime}:u_{F}\in L_{loc}^{1}\left(  \mathbb{R}^{d}\setminus\mathcal{A}%
\right)  \text{ }and\text{ }\int\frac{\left\vert u_{F}\right\vert ^{2}}%
{w}<\infty\right\}  ,
\]

and \ref{a3.3}:%
\[
X_{1/w}^{0}=\left\{  v\in S^{\prime}:v_{F}\in L_{loc}^{1}\left(
\mathbb{R}^{d}\setminus\mathcal{A}\right)  \text{ }and\text{ }\int%
\frac{\left\vert v_{F}\right\vert ^{2}}{w}<\infty\right\}  ,
\]

of the spaces $\widetilde{X}_{1/w}^{0}$ and $X_{1/w}^{0}$ and it is clear that
$X_{1/w}^{0}\subset\widetilde{X}_{1/w}^{0}$ as sets.

Since $\widehat{S}_{w,0}$ is a subspace of $S$, the obvious mappings between
$\widetilde{X}_{1/w}^{0}$ and $X_{1/w}^{0}$ are the \textit{extension mapping}
$E:\left(  \widehat{S}_{w,0}\right)  ^{\prime}\rightarrow S^{\prime}$ (defined
using the Hahn-Banach Lemma \ref{Lem_convex_tls_extend_2}) and the
\textit{restriction mapping} $R:S^{\prime}\rightarrow\left(  \widehat{S}%
_{w,0}\right)  ^{\prime}$. Lemma \ref{Lem_convex_tls_extend_2} tells us that
$E:\left(  \widehat{S}_{w,0}\right)  ^{\prime}\rightarrow S^{\prime}$ is
unique up to an annihilator of $\widehat{S}_{w,0}$ so it is unique and linear
in the seminorm sense. Clearly%
\begin{equation}
ER=RE=I.\label{a3.6}%
\end{equation}

Now if $\phi\in C_{0}^{\infty}\left(  \mathbb{R}^{d}\setminus\mathcal{A}%
\right)  $ and $u\in\widetilde{X}_{1/w}^{0}$ then $u\in\left(  \widehat{S}%
_{w,0}\right)  ^{\prime}$, $\phi\in S_{w,0}$, $\widehat{\phi}\in
\widehat{S}_{w,0}$ and so%
\[
\left[  \left(  Eu\right)  _{F},\phi\right]  =\left[  \widehat{Eu}%
|_{\mathbb{R}^{d}\setminus\mathcal{A}},\phi\right]  =\left[  \widehat{Eu}%
,\phi\right]  =\left[  Eu,\widehat{\phi}\right]  =\left[  u,\widehat{\phi
}\right]  =\left[  \widehat{u},\phi\right]  =\left[  u_{F},\phi\right]  .
\]

This means that $\left(  Eu\right)  _{F}=u_{F}$ on $\mathbb{R}^{d}%
\setminus\mathcal{A}$ and hence that
\begin{equation}
E:\widetilde{X}_{1/w}^{0}\rightarrow X_{1/w}^{0},\quad\left\vert Eu\right\vert
_{1/w,0}=\left\vert u\right\vert _{1/w,0}.\label{a3.4}%
\end{equation}

On the other hand, suppose $\phi\in C_{0}^{\infty}\left(  \mathbb{R}%
^{d}\setminus\mathcal{A}\right)  $, $v\in X_{1/w}^{0}$ and set $u=Rv$. Then
$u\in\left(  \widehat{S}_{w,0}\right)  ^{\prime}$, $\widehat{u}\in
S_{w,0}^{\prime}$, $C_{0}^{\infty}\left(  \mathbb{R}^{d}\setminus
\mathcal{A}\right)  \subset S_{w,0}$ and so%
\[
\left[  u_{F},\phi\right]  =\left[  \widehat{u},\phi\right]  =\left[
u,\widehat{\phi}\right]  =\left[  v,\widehat{\phi}\right]  =\left[
\widehat{v},\phi\right]  =\left[  v_{F},\phi\right]  ,
\]

so that $u_{F}=v_{F}$ on $\mathbb{R}^{d}\setminus\mathcal{A}$ and consequently
$Rv\in\widetilde{X}_{1/w}^{0}$ and%
\begin{equation}
R:X_{1/w}^{0}\rightarrow\widetilde{X}_{1/w}^{0},\quad\left\vert Rv\right\vert
_{1/w,0}=\left\vert v\right\vert _{1/w,0}.\label{a3.5}%
\end{equation}

To summarize:

\begin{theorem}
\label{Thm_exten_restrict_X}The natural extension operator $E:\left(
\widehat{S}_{w,0}\right)  ^{\prime}\rightarrow S^{\prime}$ and the restriction
operator $R:S^{\prime}\rightarrow\left(  \widehat{S}_{w,0}\right)  ^{\prime}$
discussed above are inverses and isometric isomorphisms between $\widetilde{X}%
_{1/w}^{0}$ and $X_{1/w}^{0}$.
\end{theorem}

\begin{proof}
Re-statement of the equations and mappings \ref{a3.6}, \ref{a3.4} and
\ref{a3.5}.
\end{proof}

\subsection{The Fourier-independent spaces $\protect\widetilde{Y}_{1/w}^{0}$
and $Y_{w}^{0}$\label{SbSect_tildYo1/w_Yow_FTindep}}

In this subsection we separate the Fourier transform from defining the
mappings between $\widetilde{X}_{1/w}^{0}$, $L^{2}$ and $X_{w}^{0}$. The
operators between $\widetilde{X}_{1/w}^{0}$, $L^{2}$ and $X_{w}^{0}$ induce
corresponding operators between $\widetilde{Y}_{1/w}^{0}$, $L^{2}$ and
$Y_{w}^{0}$ which are defined \textbf{without using the Fourier transform}.
The maps induced in the bottom row of Figure
\ref{Fig_ops_Y_tildY_Fourier_indep} are Fourier-independent.%

\begin{gather}%
\begin{array}
[c]{ccccc}%
\widetilde{X}_{1/w}^{0} &
\begin{array}
[c]{c}%
\overset{\mathcal{L}}{\longleftarrow}\\
\underset{\mathcal{M}}{\longrightarrow}%
\end{array}
& L^{2} &
\begin{array}
[c]{c}%
\overset{\mathcal{I}}{\longleftarrow}\\
\underset{\mathcal{J}}{\longrightarrow}%
\end{array}
& X_{w}^{0}\\
F^{-1}\uparrow\downarrow\mathcal{F\quad} &  & id\updownarrow\quad &  &
\mathcal{\quad}F\uparrow\downarrow F^{-1}\\
\widetilde{Y}_{1/w}^{0} & \leftrightarrows & L^{2} & \leftrightarrows &
Y_{w}^{0}%
\end{array}
\label{Fig_ops_Y_tildY_Fourier_indep}\\
\text{Fourier-independent maps and spaces }Y_{w}^{0}\text{ and }%
\widetilde{Y}_{1/w}^{0}\text{.}\nonumber
\end{gather}

\begin{definition}
\label{Def_tildYo1/w_Yow}\textbf{The spaces }$\widetilde{Y}_{1/w}^{0}$
\textbf{and }$Y_{w}^{0}$
\begin{equation}
\widetilde{Y}_{1/w}^{0}=\left\{  u\in S_{w,0}^{\prime}:u_{E}\in L_{loc}%
^{1}\left(  \mathbb{R}^{d}\setminus\mathcal{A}\right)  \text{ }and\text{ }%
\int\frac{\left\vert u_{E}\right\vert ^{2}}{w}<\infty\right\}  ,\label{a3.1}%
\end{equation}

where $u_{E}$ is the restriction of $u$ to $\mathbb{R}^{d}\setminus
\mathcal{A}$. We endow $\widetilde{Y}_{1/w}^{0}$ with the seminorm and
semi-inner product%
\[
\left\vert u\right\vert _{1/w,0}^{\prime}=\left(  \int\frac{\left\vert
u_{E}\right\vert ^{2}}{w}\right)  ^{1/2},\quad\left\langle u,v\right\rangle
_{1/w,0}^{\prime}=\int\frac{u_{E}\overline{v_{E}}}{w}.
\]

Also%
\[
Y_{w}^{0}=\left\{  u\in L_{loc}^{1}\cap S^{\prime}:\int w\left\vert
u\right\vert ^{2}<\infty\right\}  ,
\]

endowed with the norm and inner product%
\[
\left\vert u\right\vert _{w,0}^{\prime}=\left(  \int w\left\vert u\right\vert
^{2}\right)  ^{1/2},\quad\left(  u,v\right)  _{w,0}^{\prime}=\int
wu\overline{v}.
\]

\end{definition}

\begin{theorem}
\label{Thm_tildYo1/w_Yow_property}\textbf{Properties of the spaces
}$\widetilde{Y}_{1/w}^{0}$ \textbf{and }$Y_{w}^{0}$ (see Theorem
\ref{Thm_property1_tildXo1/w}).

\begin{enumerate}
\item $\operatorname{null}\left\vert \cdot\right\vert _{1/w,0}^{\prime
}=S_{w,0;\mathcal{A}}^{\prime}$.

\item $\left\vert \cdot\right\vert _{1/w,0}^{\prime}$ is a norm iff we can
choose $\mathcal{A}$ to be the (minimal) empty set.

\item If $\left\vert u\left(  \xi\right)  \right\vert \leq c\left(
1+\left\vert \xi\right\vert \right)  ^{\kappa}$ then $u\in\widetilde{Y}%
_{1/w}^{0}$.

\item $S_{w,0}\subset Y_{w}^{0}$; $S\subset\widetilde{Y}_{1/w}^{0}$ and
$wS_{w,0}\subset\widetilde{Y}_{1/w}^{0}$.

\item $u_{E}\in\widetilde{Y}_{1/w}^{0}$ and $u-u_{E}\in S_{w,0;\mathcal{A}%
}^{\prime}$ i.e. $\left\vert u-u_{E}\right\vert _{1/w,0}^{\prime}=0$.
\end{enumerate}
\end{theorem}

\begin{proof}
Adapt the proofs of Theorem \ref{Thm_property1_tildXo1/w}.
\end{proof}

\begin{definition}
\label{Def_map_F}\textbf{The mapping }$\mathcal{F}:\widetilde{X}_{1/w}%
^{0}\rightarrow S_{w,0}^{\prime}$ is defined by%
\[
\mathcal{F}u=u_{F},
\]

which makes sense since $u_{F}\in S_{w,0}^{\prime}$ by part 5 of Theorem
\ref{Thm_property1_tildXo1/w}.
\end{definition}

We now show that $\mathcal{F}$ is an isometric isomorphism from $\widetilde{X}%
_{1/w}^{0}$ to $\widetilde{Y}_{1/w}^{0}$ in the seminorm sense.

\begin{theorem}
\label{Thm_property_map_F}\textbf{Some properties of }$\mathcal{F}$. In the
seminorm sense:

\begin{enumerate}
\item $\mathcal{F}$ is an isometry from $\widetilde{X}_{1/w}^{0}$ to
$\widetilde{Y}_{1/w}^{0}$.

\item The inverse Fourier transform $F^{-1}$ is an isometry from
$\widetilde{Y}_{1/w}^{0}$ to $\widetilde{X}_{1/w}^{0}$.

\item $\mathcal{F}F^{-1}=F^{-1}\mathcal{F}=I$.

\item $\mathcal{F}$ is an isometric isomorphism from $\widetilde{X}_{1/w}^{0}
$ to $\widetilde{Y}_{1/w}^{0}$.
\end{enumerate}
\end{theorem}

\begin{proof}
From part 5 of Theorem \ref{Thm_property1_tildXo1/w}: $u\in\widetilde{X}%
_{1/w}^{0}$ implies $u_{F}\in S_{w,0}^{\prime}\cap L_{loc}^{1}\left(
\mathbb{R}^{d}\setminus\mathcal{A}\right)  $, $\overset{\vee}{u_{F}}%
\in\widetilde{X}_{1/w}^{0}$ and $u-\overset{\vee}{u_{F}}\in\left(
S_{w,0;\mathcal{A}}^{\prime}\right)  ^{\vee}$.\medskip

\textbf{Part 1} Thus $u\in\widetilde{X}_{1/w}^{0}$ implies $u_{F}\in
S_{w,0}^{\prime}\cap L_{loc}^{1}\left(  \mathbb{R}^{d}\setminus\mathcal{A}%
\right)  $ and $\int\frac{\left\vert u_{F}\right\vert ^{2}}{w}<\infty$ i.e.
$u_{F}\in\widetilde{Y}_{1/w}^{0}$, $\left\vert u\right\vert _{1/w,0}%
=\left\vert u_{F}\right\vert _{1/w,0}^{\prime}$ and $\mathcal{F}$ is an
isometry.\medskip

\textbf{Part 2} If $v\in\widetilde{Y}_{1/w}^{0}$ then $v\in S_{w,0}^{\prime}$,
$v_{E}\in L_{loc}^{1}\left(  \mathbb{R}^{d}\setminus\mathcal{A}\right)  $ and
$\int\frac{\left\vert v_{E}\right\vert ^{2}}{w}<\infty$. Set $u=F^{-1}%
v=\overset{\vee}{v}$. Then $u\in\left(  \widehat{S}_{w,0}\right)  ^{\prime}$,
$u_{F}=v_{E}\in L_{loc}^{1}\left(  \mathbb{R}^{d}\setminus\mathcal{A}\right)
$ and $\int\frac{\left\vert u_{F}\right\vert ^{2}}{w}=\int\frac{\left\vert
v_{E}\right\vert ^{2}}{w}<\infty$ i.e. $u\in\widetilde{X}_{1/w}^{0}$ and
$F^{-1}$ is an isometry.\medskip

\textbf{Part 3} If $v\in\widetilde{Y}_{1/w}^{0}$ then $\mathcal{F}%
F^{-1}v=\left(  F^{-1}v\right)  _{F}=v_{E}$ and by part 5 of Theorem
\ref{Thm_tildYo1/w_Yow_property}, $\left\vert v-v_{E}\right\vert _{1/w,0}=0$.

If $u\in\widetilde{X}_{1/w}^{0}$ then from the preface to the proof,
$u-F^{-1}\mathcal{F}u=u-\overset{\vee}{u_{F}}\in\left(  S_{w,0;\mathcal{A}%
}^{\prime}\right)  ^{\vee}$ i.e. $\left\vert u-\overset{\vee}{u_{F}%
}\right\vert _{1/w,0}=0$.\medskip

\textbf{Part 4} This part summarizes the preceding parts.
\end{proof}

The following theorem shows that the definition of $\widetilde{X}_{1/w}^{0}$
allows the complete separation of the Fourier transform and the definition of
the space:

\begin{theorem}
\textbf{Some properties of }$\widetilde{Y}_{1/w}^{0}$\textbf{\ and }$Y_{w}%
^{0}$\textbf{\ in relation to }$\widetilde{X}_{1/w}^{0}$\textbf{\ and
}$\widetilde{Y}_{1/w}^{0}$:

\begin{enumerate}
\item $\mathcal{F}$ is an isometric isomorphism from $\widetilde{X}_{1/w}^{0}
$ to $\widetilde{Y}_{1/w}^{0}$.

\item $\widetilde{Y}_{1/w}^{0}$ is a semi-Hilbert space.

\item $F^{-1}:X_{w}^{0}\rightarrow Y_{w}^{0}$ is an isometric isomorphism in
the norm sense.

\item $Y_{w}^{0}$ is a Hilbert space.
\end{enumerate}
\end{theorem}

\begin{proof}
\textbf{Part 1} Re-statement of part 4 of Theorem \ref{Thm_property_map_F}%
.\medskip

\textbf{Part 2 }The space $\widetilde{X}_{1/w}^{0}$ was shown to be a
semi-Hilbert space in Corollary \ref{Cor_tildXo1/w_semiHilb_Sw2}.\medskip

\textbf{Part 3} Easy.\medskip

\textbf{Part 4} The space $X_{w}^{0}$ was shown to be a Hilbert space in
Corollary \ref{Cor_Xow_complete} so this result follows from part 3.
\end{proof}

Other Fourier transform-independent results can be easily derived e.g. from
Theorem \ref{Cor_invF[wSw2]_dense_tildXo1/w}.

\subsection{Some weight functions not in $W_{S;0}$%
\label{SbSect_wt_fn_not_in_S2}}

In this section we have so far assumed that $w\in W_{S;0}$. Now we show that
the Gaussian weight function and the shifted thin-plate spline weight function
do not belong to $W_{S;0}$.

The next lemma will be used to show that the Gaussian weight function is not a
member of $W_{S;0}$ (\ref{a1.8}).

\begin{lemma}
\label{Lem_w_not_in_S2dagger_1}Suppose the weight function $w$ satisfies $w\in
C^{\infty}$, $w>0$ and $1/w\in S$. Further, suppose there exists $b>0$ such
that $\frac{1}{w\left(  \lambda x\right)  }\in S_{w,0}$ when $\lambda>b$ and
$\int\frac{w\left(  x\right)  }{w\left(  \lambda x\right)  ^{2}}%
dx\rightarrow\infty$ as $\lambda\rightarrow b^{+}$.

Then $w\notin W_{S;0}$.
\end{lemma}

\begin{proof}
Suppose that $w\in W_{S;0}$. Then there exists a positive, linear combination
of the seminorms \ref{a1.5} such that
\[
\left(  \int w\left\vert \phi\right\vert ^{2}\right)  ^{1/2}\leq\sum
_{\alpha,\beta}c_{\alpha,\beta}\left\Vert x^{\alpha}D^{\beta}\phi\right\Vert
_{\infty},\quad\phi\in S_{w,0}.
\]

Choose $\phi_{\lambda}\left(  x\right)  =\frac{1}{w\left(  \lambda x\right)
}$ where $\lambda>0$ and $\lambda$ is real. Then%
\[
\int w\left\vert \phi_{\lambda}\right\vert ^{2}dx=\int\frac{w\left(  x\right)
}{w\left(  \lambda x\right)  ^{2}}dx<\infty.
\]

On the other hand%
\begin{align*}
\left\Vert x^{\alpha}D_{x}^{\beta}\phi_{\lambda}\right\Vert _{\infty
}=\left\Vert x^{\alpha}D_{x}^{\beta}\frac{1}{w\left(  \lambda x\right)
}\right\Vert _{\infty}=\lambda^{\left\vert \beta\right\vert }\left\Vert
x^{\alpha}\left(  D^{\beta}\frac{1}{w}\right)  \left(  \lambda x\right)
\right\Vert _{\infty} &  \leq\lambda^{\left\vert \beta\right\vert -\left\vert
\alpha\right\vert }\left\Vert \left(  \lambda x\right)  ^{\alpha}\left(
D^{\beta}\frac{1}{w}\right)  \left(  \lambda x\right)  \right\Vert _{\infty}\\
&  =\lambda^{\left\vert \beta\right\vert -\left\vert \alpha\right\vert
}\left\Vert x^{\alpha}D^{\beta}\frac{1}{w}\right\Vert _{\infty},
\end{align*}

so that%
\[
\int\frac{w\left(  x\right)  }{w\left(  \lambda x\right)  ^{2}}dx\leq
\sum_{\alpha,\beta}c_{\alpha,\beta}\lambda^{\left\vert \beta\right\vert
-\left\vert \alpha\right\vert }\left\Vert x^{\alpha}D^{\beta}\frac{1}%
{w}\right\Vert _{\infty}.
\]

But as $\lambda\rightarrow b^{+}$ our assumptions imply that the right side
diverges and the left side converges. Hence $w\notin W_{S;0}$ which is a contradiction.
\end{proof}

It is now easy to apply our lemma to the Gaussian.

\begin{theorem}
The Gaussian weight function $w=e^{\left\vert \cdot\right\vert ^{2}}$ is not a
member of $W_{S;0}$.
\end{theorem}

\begin{proof}
Clearly $w\in C^{\infty}$, $w>0$ and $1/w\in S$. Further,
\[
\int\frac{w\left(  x\right)  }{w\left(  \lambda x\right)  ^{2}}dx=\int
e^{\left\vert x\right\vert ^{2}}e^{-2\lambda^{2}\left\vert x\right\vert ^{2}%
}=\int e^{-\left(  2\lambda^{2}-1\right)  \left\vert x\right\vert ^{2}}%
=\sqrt{\frac{\pi}{2\lambda^{2}-1}},
\]

iff $\lambda>\frac{1}{\sqrt{2}}$ so that $\frac{1}{w\left(  \lambda x\right)
}\in S_{w,0}$ iff $\lambda>b=\frac{1}{\sqrt{2}}$. When $\lambda\rightarrow
b^{+}$ the integrals diverge.
\end{proof}

This lemma will be used to show that the shifted thin-plate splines weight
functions $w$ are not members of $W_{S;0}$. Indeed, these weight functions
lack $C^{\infty}$ smoothness at the origin and $1/w\notin S$. But it is enough
to assume that $1/w$ is like $S$ on $\mathbb{R}^{d}\setminus0$.

\begin{lemma}
\label{Lem_w_not_in_S2dag_2}Suppose that $w\in C^{\infty}\left(
\mathbb{R}^{d}\setminus0\right)  $, $w>0$ on $\mathbb{R}^{d}\setminus0$ and
$\phi_{w}=\frac{\psi}{w}\in S$ for some $\psi\in C^{\infty}$ such that
$0\leq\psi\leq1$, $\psi=0$ near the origin and $\psi\left(  x\right)  =1$ when
$\left\vert x\right\vert \geq1$.

Further, suppose there exists $b>0$ such that $\int\limits_{\left\vert
x\right\vert \geq1}\frac{w\left(  x\right)  }{w\left(  \lambda x\right)  ^{2}%
}dx<\infty$ when $\lambda>b$ and $\int\limits_{\left\vert x\right\vert \geq
1}\frac{w\left(  x\right)  }{w\left(  \lambda x\right)  ^{2}}dx\rightarrow
\infty$ as $\lambda\rightarrow b^{+}$.

Then $w$ does not satisfy Condition \ref{Cnd_wt_fn_T}.
\end{lemma}

\begin{proof}
Suppose that $w$ satisfies Condition \ref{Cnd_wt_fn_T}. Then there exists a
positive, linear combination of the seminorms \ref{a1.5} such that
\[
\left(  \int w\left\vert \phi\right\vert ^{2}\right)  ^{1/2}\leq\sum
_{\alpha,\beta}c_{\alpha,\beta}\left\Vert x^{\alpha}D^{\beta}\phi\right\Vert
_{\infty},\quad\phi\in S_{w,0}.
\]

We now show that $\phi_{w}\left(  \lambda x\right)  \in S_{w,0}$ when
$\lambda>b$ i.e. that
\[
\int w\left(  x\right)  \phi_{w}\left(  \lambda x\right)  ^{2}dx<\infty,
\]

when $\lambda>b$. Indeed%
\begin{align*}
\int w\left(  x\right)  \phi_{w}\left(  \lambda x\right)  ^{2}dx  &
=\int\limits_{\left\vert x\right\vert \leq1}w\left(  x\right)  \phi_{w}\left(
\lambda x\right)  ^{2}dx+\int\limits_{\left\vert x\right\vert \geq1}w\left(
x\right)  \phi_{w}\left(  \lambda x\right)  ^{2}dx\\
& =\int\limits_{\left\vert x\right\vert \leq1}w\left(  x\right)  \phi
_{w}\left(  \lambda x\right)  ^{2}dx+\int\limits_{\left\vert x\right\vert
\geq1}\frac{w\left(  x\right)  }{w\left(  \lambda x\right)  ^{2}}\psi\left(
\lambda x\right)  ^{2}dx\\
& \leq\int\limits_{\left\vert x\right\vert \leq1}w\left(  x\right)  \phi
_{w}\left(  \lambda x\right)  ^{2}dx+\int\limits_{\left\vert x\right\vert
\geq1}\frac{w\left(  x\right)  }{w\left(  \lambda x\right)  ^{2}}dx,
\end{align*}

and since the second integral exists when $\lambda>b$ and the first integrand
is continuous, we can conclude that $\phi_{w}\left(  \lambda x\right)  \in
S_{w,0}$ when $\lambda>b$. On the other hand%
\begin{align*}
\left\Vert x^{\alpha}D_{x}^{\beta}\phi_{w}\left(  \lambda x\right)
\right\Vert _{\infty}=\lambda^{\left\vert \beta\right\vert }\left\Vert
x^{\alpha}\left(  D^{\beta}\phi_{w}\right)  \left(  \lambda x\right)
\right\Vert _{\infty} &  \leq\lambda^{\left\vert \beta\right\vert -\left\vert
\alpha\right\vert }\left\Vert \left(  \lambda x\right)  ^{\alpha}\left(
D^{\beta}\phi_{w}\right)  \left(  \lambda x\right)  \right\Vert _{\infty}\\
&  =\lambda^{\left\vert \beta\right\vert -\left\vert \alpha\right\vert
}\left\Vert x^{\alpha}D^{\beta}\phi_{w}\right\Vert _{\infty},
\end{align*}

so that when $\lambda>b$%
\[
\int\limits_{\left\vert x\right\vert \geq1}\frac{w\left(  x\right)  }{w\left(
\lambda x\right)  ^{2}}dx\leq c\int w\left(  x\right)  \phi_{w}\left(  \lambda
x\right)  ^{2}dx\leq c\sum_{\alpha,\beta}c_{\alpha,\beta}\lambda^{\left\vert
\beta\right\vert -\left\vert \alpha\right\vert }\left\Vert x^{\alpha}D^{\beta
}\phi_{w}\right\Vert _{\infty}.
\]

But as $\lambda\rightarrow b^{+}$ our assumptions imply that the right side
diverges to infinity and the left side converges. Hence $w\notin W_{S;0}$.
\end{proof}

\begin{theorem}
No shifted thin-plate spline weight function is a member of $W_{S;0}$.
\end{theorem}

\begin{proof}
From Subsubsection \ref{Ex_shft_thn_plt_spln_wt}, if we set $\mu=v+d/2$ then
the shifted thin-plate spline weight functions have the form
\[
w\left(  x\right)  =\frac{1}{\widetilde{e}\left(  v\right)  }\frac{\left\vert
x\right\vert ^{2\mu}}{\widetilde{K}_{\mu}\left(  \left\vert x\right\vert
\right)  },\quad0<\mu<d/2,
\]

where $\widetilde{K}_{\mu}\left(  t\right)  =t^{\mu}K_{\mu}\left(  t\right)  ,
$ $t\geq0,\mu>0$ and $K_{\lambda}$ is called a modified Bessel function or
MacDonald's function. From Theorem \ref{Thm_bnds_modif_MacDonald} this
function has the key property that for some constants $c_{\mu},c_{\mu}%
^{\prime}>0$,
\[
c_{\mu}e^{-t}\leq\widetilde{K}_{\mu}\left(  t\right)  \leq c_{\mu}^{\prime
}e^{-t},\text{\quad}\mu>0,\text{ }t\geq0,
\]

so that%
\begin{align*}
\int\limits_{\left\vert x\right\vert \geq1}\frac{w\left(  x\right)  }{w\left(
\lambda x\right)  ^{2}}dx  & =\int\limits_{\left\vert x\right\vert \geq1}%
\frac{1}{\widetilde{e}\left(  v\right)  }\frac{\left\vert x\right\vert ^{2\mu
}}{\widetilde{K}_{\mu}\left(  \left\vert x\right\vert \right)  }\frac
{1}{\left(  \frac{1}{\widetilde{e}\left(  v\right)  }\frac{\left\vert \lambda
x\right\vert ^{2\mu}}{\widetilde{K}_{\mu}\left(  \left\vert \lambda
x\right\vert \right)  }\right)  ^{2}}dx\\
& =\widetilde{e}\left(  v\right)  \int\limits_{\left\vert x\right\vert \geq
1}\frac{\left(  \widetilde{K}_{\mu}\left(  \left\vert \lambda x\right\vert
\right)  \right)  ^{2}}{\widetilde{K}_{\mu}\left(  \left\vert x\right\vert
\right)  }\frac{dx}{\left\vert \lambda x\right\vert ^{2\mu}}\\
& =\frac{\widetilde{e}\left(  v\right)  }{\lambda^{2\mu}}\int%
\limits_{\left\vert x\right\vert \geq1}\frac{\left(  \widetilde{K}_{\mu
}\left(  \left\vert \lambda x\right\vert \right)  \right)  ^{2}}%
{\widetilde{K}_{\mu}\left(  \left\vert x\right\vert \right)  }\frac
{dx}{\left\vert x\right\vert ^{2\mu}}\\
& \geq\frac{\widetilde{e}\left(  v\right)  }{\lambda^{2\mu}}\int%
\limits_{\left\vert x\right\vert \geq1}\frac{\left(  c_{\mu}e^{-\left\vert
\lambda x\right\vert }\right)  ^{2}}{c_{\mu}^{\prime}e^{-\left\vert
x\right\vert }}\frac{dx}{\left\vert x\right\vert ^{2\mu}}\\
& =\frac{\left(  c_{\mu}\right)  ^{2}}{c_{\mu}^{\prime}}\frac{\widetilde{e}%
\left(  v\right)  }{\lambda^{2\mu}}\int\limits_{\left\vert x\right\vert \geq
1}\frac{e^{\left(  1-2\lambda\right)  \left\vert x\right\vert }}{\left\vert
x\right\vert ^{2\mu}}dx,
\end{align*}

and since $2\mu<d$ the last integral exists iff $1-2\lambda<0$ i.e. iff
$\lambda>\frac{1}{2}$. Now set $b=\frac{1}{2}$ and assume $\lambda>b$. The
change of variables $y=\left(  \frac{\lambda}{b}-1\right)  x$ implies
$dy=\left(  \frac{\lambda}{b}-1\right)  ^{d}dx$ and $\left\vert y\right\vert
=\left(  \frac{\lambda}{b}-1\right)  \left\vert x\right\vert $ and hence%
\begin{align*}
\int\limits_{\left\vert x\right\vert \geq1}\frac{e^{\left(  1-2\lambda\right)
\left\vert x\right\vert }}{\left\vert x\right\vert ^{2\mu}}dx=\int%
\limits_{\left\vert x\right\vert \geq1}\frac{e^{-\left(  \frac{\lambda}%
{b}-1\right)  \left\vert x\right\vert }}{\left\vert x\right\vert ^{2\mu}}dx &
=\frac{1}{\left(  \frac{\lambda}{b}-1\right)  ^{d-2\mu}}\int%
\limits_{\left\vert y\right\vert \geq\frac{\lambda}{b}-1}\frac{e^{-\left\vert
y\right\vert }}{\left\vert y\right\vert ^{2\mu}}dy\\
&  \geq\frac{1}{\left(  \frac{\lambda}{b}-1\right)  ^{d-2\mu}}\int%
\limits_{\left\vert y\right\vert \geq1}\frac{e^{-\left\vert y\right\vert }%
}{\left\vert y\right\vert ^{2\mu}}dy,
\end{align*}

when $\lambda\leq2b$. Clearly $\int\limits_{\left\vert x\right\vert \geq
1}\frac{w\left(  x\right)  }{w\left(  \lambda x\right)  ^{2}}dx\rightarrow
\infty$ as $\lambda\rightarrow b^{+}$ and so $w\notin W_{S;0}$ by Lemma
\ref{Lem_w_not_in_S2dag_2}.
\end{proof}

We end this subsection with a (unused) generalization of the above results.

\begin{lemma}
\label{Lem_w_not_in_S2dag_3}Suppose that $w\in C^{\infty}\left(
\mathbb{R}^{d}\setminus\mathcal{A}\right)  $ and $\phi_{w}=\frac{\psi}{w}\in
S$ for some $\psi\in C_{B}^{\infty}$ such that $0\leq\psi\leq1$.

Further, suppose there exists $b>0$ such that $\int w\left(  x\right)
\left\vert \phi_{w}\left(  \lambda x\right)  \right\vert ^{2}dx<\infty$ when
$\lambda>b$ and that

$\int w\left(  x\right)  \left\vert \phi_{w}\left(  \lambda x\right)
\right\vert ^{2}dx\rightarrow\infty$ as $\lambda\rightarrow b^{+}$.

Then $w\notin W_{S;0}$.
\end{lemma}

\begin{proof}
Suppose that $w\in W_{S;0}$. Now $\phi_{w}\left(  \lambda x\right)  \in
S_{w,0}$ when $\lambda>b$ and so there exists a positive, linear combination
of the seminorms \ref{a1.5} such that
\[
\left(  \int w\left(  x\right)  \left\vert \phi_{w}\left(  \lambda x\right)
\right\vert ^{2}dx\right)  ^{1/2}\leq\sum_{\alpha,\beta}c_{\alpha,\beta
}\left\Vert x^{\alpha}D_{x}^{\beta}\phi_{w}\left(  \lambda x\right)
\right\Vert _{\infty},\quad\phi\in S_{w,0}.
\]

But%
\begin{align*}
\left\Vert x^{\alpha}D_{x}^{\beta}\phi_{w}\left(  \lambda x\right)
\right\Vert _{\infty}=\lambda^{\left\vert \beta\right\vert }\left\Vert
x^{\alpha}\left(  D^{\beta}\phi_{w}\right)  \left(  \lambda x\right)
\right\Vert _{\infty} &  \leq\lambda^{\left\vert \beta\right\vert -\left\vert
\alpha\right\vert }\left\Vert \left(  \lambda x\right)  ^{\alpha}\left(
D^{\beta}\phi_{w}\right)  \left(  \lambda x\right)  \right\Vert _{\infty}\\
&  =\lambda^{\left\vert \beta\right\vert -\left\vert \alpha\right\vert
}\left\Vert x^{\alpha}D^{\beta}\phi_{w}\right\Vert _{\infty},
\end{align*}

so that%
\[
\left(  \int w\left(  x\right)  \left\vert \phi_{w}\left(  \lambda x\right)
\right\vert ^{2}dx\right)  ^{1/2}\leq\lambda^{\left\vert \beta\right\vert
-\left\vert \alpha\right\vert }\left\Vert x^{\alpha}D^{\beta}\phi
_{w}\right\Vert _{\infty}.
\]

But as $\lambda\rightarrow b^{+}$ our assumptions imply that the right side
diverges to infinity and the left side converges. Hence $w\notin W_{S;0}$.
\end{proof}

\section{A generalization of $\protect\widetilde{X}_{1/w}^{0}$%
\label{Sect_generalize1_tildXo1/w_T}}

\subsection{Motivation}

In Section \ref{Sect_tildXo1/w_S2dag} we defined $\widetilde{X}_{1/w}^{0}$ for
the case when $w\in S_{w,0}^{\dag}$ and in Section \ref{Sect_Xo1/w} we defined
the space $X_{1/w}^{0}$ for the case when $w$ satisfies \ref{a1.046}. In this
section we show how to generalize our definition of $\widetilde{X}_{1/w}^{0}$
in order to include the case when $w$ satisfies \ref{a1.046}.

We would like to retain a form of the integral condition \ref{a1.8} used to
define $S_{w,0}^{\dag}$ and this is where Lemma
\ref{Lem_wt_fn_extra_property_2} and Theorem \ref{Thm_Som_and_Sw2} are helpful
as they establish results which suggest a generalization. In Subsubsection
\ref{SbSbSect_wt_fn_rad_L1loc} it was assumed that a weight function
satisfying \ref{a1.046} had the form \ref{a1.09} i.e. for some integer
$m\geq0$,%
\[
w\left(  x\right)  =\frac{v\left(  x\right)  }{\left\vert x\right\vert
^{2\left(  m+s\right)  }},\quad\left\vert x\right\vert \leq r;\quad
2s<d;\quad0<c_{1}\leq v\left(  x\right)  \leq c_{2},
\]

and in Theorem \ref{Thm_Sw2_eq_Son} the space $S_{w,0}$ was shown to be of the
form $S_{\emptyset,n}$ and an integral condition \ref{a1.8} derived with $w\in
W_{S;0}$. However, if we just assume the weaker \ref{a1.046} is satisfied then
from part 4 of Lemma \ref{Lem_wt_fn_extra_property_2}
\begin{equation}
\left(  \int w\left\vert \psi\right\vert ^{2}\right)  ^{1/2}\leq\left\vert
\psi\right\vert _{\sigma^{\prime}},\quad\psi\in S_{\emptyset,m},\label{a4.5}%
\end{equation}

for some $m$ and a positive, linear combination $\left\vert \cdot\right\vert
_{\sigma^{\prime}}$ of the seminorms \ref{a1.2} which define the topology of
the tempered distributions $S$. Hence $S_{\emptyset,m}\subset S_{w,0}$. And if
\ref{a4.5} holds then we would expect Theorem \ref{Thm_Som_and_Sw2} to allow
us to define the isometric operator $\mathcal{L}:L^{2}\rightarrow
\widetilde{X}_{1/w}^{0}$ in the manner of Definition \ref{Def_tildL}. So in
the next subsection Condition \ref{Cnd_wt_fn_T} will postulate the existence
of a subspace $T\subset S_{w,0}$ such that \ref{a4.5} holds for $\psi\in T$.

Part 5 of Definition \ref{Def_Sw2_lin_fnal_0} now suggests we assume $T$ has
the property that if $\Omega\subset\mathbb{R}^{d}\setminus\mathcal{A}$ then
$C_{0}^{\infty}\left(  \Omega\right)  \subset T$ i.e. that $C_{0}^{\infty
}\left(  \mathbb{R}^{d}\setminus\mathcal{A}\right)  \subset T$. However, if
for the case \ref{a1.046} we choose $T=S_{\emptyset,m}$, $w=1$ near the origin
and $w>0$ everywhere then the minimal $\mathcal{A}$ is empty. But noting part
1 of Remark \ref{Rem_Def_extend_wt_fn} we can "adjust" $\mathcal{A}$ by
enlarging it so that $\mathcal{A}=\left\{  0\right\}  $ for $T=S_{\emptyset
,m}$.

Given the above considerations our approach will now be to replace the space
$S_{w,0}$\ by the space $T$ in the results and definitions of Section
\ref{Sect_tildXo1/w_S2dag} and show that they all still hold.

However, in the next subsection we will show that the Gaussian and the shifted
thin-plate spline weight functions still do not satisfy Condition
\ref{Cnd_wt_fn_T}.

\subsection{Condition \ref{Cnd_wt_fn_T}}

We encapsulate the properties discussed in the motivation as:

\begin{condition}
\label{Cnd_wt_fn_T}Suppose $w$ is a weight function w.r.t. a closed set
$\mathcal{A}$ which has measure zero. Suppose $T$ is a subspace of $S$ such that:

\begin{enumerate}
\item There exists a positive, linear combination $\left\vert \cdot\right\vert
_{\sigma}$ of the seminorms which define $S$ such that
\begin{equation}
\left(  \int w\left\vert \phi\right\vert ^{2}\right)  ^{1/2}\leq\left\vert
\phi\right\vert _{\sigma},\quad\phi\in T.\label{a4.1}%
\end{equation}

\item There exists an $\mathcal{A}$ such that $C_{0}^{\infty}\left(
\mathbb{R}^{d}\setminus\mathcal{A}\right)  \subset T$ - see part 1 of Remark
\ref{Rem_Def_extend_wt_fn}.

\item There exists a sequence of open balls $B\left(  x^{\left(  k\right)
};r_{k}\right)  \subset C_{0}^{\infty}\left(  \mathbb{R}^{d}\setminus
\mathcal{A}\right)  $ such that $\left\vert x^{\left(  k\right)  }\right\vert
\rightarrow\infty$ and $r_{k}\geq c\left\vert x^{\left(  k\right)
}\right\vert ^{-m}$ for some constant $c>0$ and integer $m\geq0$.
\end{enumerate}
\end{condition}

\begin{remark}
\label{Rem_Cnd_wt_fn_T}\ 

\begin{enumerate}
\item With reference to part 1 of Remark \ref{Rem_Def_extend_wt_fn}, part 2 of
this condition caters for the case \ref{a1.046} where $w=1$ near the origin
and $w>0$ everywhere. Here the minimal $\mathcal{A}$ is empty and noting Lemma
\ref{Lem_wt_fn_extra_property_2} we can "adjust" $\mathcal{A}$ by enlarging it
so that $\mathcal{A}=\left\{  0\right\}  $ and $T=S_{\emptyset,m}$. Note
however that $T_{\mathcal{A}}^{\prime}$ now changes.

\item There is no loss of generality if we assume that the open balls satisfy
$\left\vert x^{\left(  k\right)  }\right\vert \rightarrow\infty$ and $1\geq
r_{k}\geq\left\vert x^{\left(  k\right)  }\right\vert ^{-m}$.

\item An example of $\mathcal{A}$ can be many (tensor product) $\mathbb{R}%
^{d-1}$ hyperplanes perpendicular to each axis which get closer together
towards infinity. Part 3 restricts how quickly these hyperplanes get closer together.

\item If Condition \ref{Cnd_wt_fn_T} holds for $T_{1}$ and $T_{2}$ then it
holds for the subspace $T_{1}+T_{2}$.
\end{enumerate}
\end{remark}

\subsection{Some weight functions which do not satisfy Condition
\ref{Cnd_wt_fn_T}\label{SbSect_condit_T}}

The following results are analogous to those proved in Subsection
\ref{SbSect_wt_fn_not_in_S2} concerning $W_{S;0}$ where it was shown that the
Gaussian and the shifted thin-plate spline weight functions do not belong to
$W_{S;0}$. Here we show the same weight functions do not satisfy Condition
\ref{Cnd_wt_fn_T}.

The next lemma will be used to show that the Gaussian weight function does not
satisfy Condition \ref{Cnd_wt_fn_T}.

\begin{lemma}
\label{Lem_condit_A_empty}Suppose the weight function $w$ satisfies $w\in
C^{\infty}$, $w>0$ and $1/w\in S$. Then $w$ does not satisfy Condition
\ref{Cnd_wt_fn_T}.
\end{lemma}

\begin{proof}
Suppose that $w$ satisfies Condition \ref{Cnd_wt_fn_T} for some $T$ and
$\mathcal{A}$. Then there exists a positive, linear combination of the
seminorms \ref{a1.5} such that
\[
\left(  \int w\left\vert \phi\right\vert ^{2}\right)  ^{1/2}\leq\sum
_{\alpha,\beta}c_{\alpha,\beta}\left\Vert x^{\alpha}D^{\beta}\phi\right\Vert
_{\infty},\quad\phi\in T,
\]

and%
\[
C_{0}^{\infty}\left(  \mathbb{R}^{d}\setminus\mathcal{A}\right)  \subset T.
\]

Further, by part 2 of Remark \ref{Rem_Cnd_wt_fn_T} there exists a sequence a
functions $\phi_{k}\in C_{0}^{\infty}\left(  \mathbb{R}^{d}\setminus
\mathcal{A}\right)  $ defined by
\begin{equation}
\phi_{k}\left(  x\right)  =\phi_{0}\left(  \frac{x-x^{\left(  k\right)  }%
}{r_{k}}\right)  ,\quad\phi_{0}\in C_{0}^{\infty}\left(  B\left(  0;1\right)
\right)  ,\label{a5.0}%
\end{equation}

and such that%
\begin{equation}
\left\vert x^{\left(  k\right)  }\right\vert \rightarrow\infty,\quad1\geq
r_{k}\geq\left\vert x^{\left(  k\right)  }\right\vert ^{-m}.\label{a5.1}%
\end{equation}

Hence $\operatorname*{supp}\phi_{k}\subset B\left(  x^{\left(  k\right)
};r_{k}\right)  \subset\mathbb{R}^{d}\setminus\mathcal{A}$. Thus we must have%
\begin{equation}
\left(  \int w\left\vert \phi_{k}\right\vert ^{2}\right)  ^{1/2}\leq
\sum_{\alpha,\beta}c_{\alpha,\beta}\left\Vert x^{\alpha}D^{\beta}\phi
_{k}\right\Vert _{\infty},\label{a5.2}%
\end{equation}

and the form \ref{a5.0} of $\phi_{k}$ implies%
\begin{align*}
\left\Vert x^{\alpha}D^{\beta}\phi_{k}\right\Vert _{\infty}=\left\Vert
x^{\alpha}D_{x}^{\beta}\phi_{0}\left(  \frac{x-x^{\left(  k\right)  }}{r_{k}%
}\right)  \right\Vert _{\infty} &  =r_{k}^{-\left\vert \beta\right\vert
}\left\Vert x^{\alpha}\left(  D_{x}^{\beta}\phi_{0}\right)  \left(
\frac{x-x^{\left(  k\right)  }}{r_{k}}\right)  \right\Vert _{\infty}\\
&  =r_{k}^{-\left\vert \beta\right\vert }\left\Vert \left(  r_{k}y+x^{\left(
k\right)  }\right)  ^{\alpha}\left(  D^{\beta}\phi_{0}\right)  \left(
y\right)  \right\Vert _{\infty}\\
&  \leq r_{k}^{-\left\vert \beta\right\vert }\left\Vert \left\vert
r_{k}y+x^{\left(  k\right)  }\right\vert ^{\left\vert \alpha\right\vert
}\left(  D^{\beta}\phi_{0}\right)  \left(  y\right)  \right\Vert _{\infty}\\
&  \leq r_{k}^{-\left\vert \beta\right\vert }\left\Vert \left(  r_{k}%
\left\vert y\right\vert +\left\vert x^{\left(  k\right)  }\right\vert \right)
^{\left\vert \alpha\right\vert }\left(  D^{\beta}\phi_{0}\right)  \left(
y\right)  \right\Vert _{\infty}\\
&  \leq r_{k}^{-\left\vert \beta\right\vert }\sum\limits_{j\leq\left\vert
\alpha\right\vert }\left\Vert \binom{\left\vert \alpha\right\vert }{j}%
r_{k}^{j}\left\vert y\right\vert ^{j}\left\vert x^{\left(  k\right)
}\right\vert ^{\left\vert \alpha\right\vert -j}D^{\beta}\phi_{0}\left(
y\right)  \right\Vert _{\infty}\\
&  \leq r_{k}^{-\left\vert \beta\right\vert }\sum\limits_{j\leq\left\vert
\alpha\right\vert }\binom{\left\vert \alpha\right\vert }{j}r_{k}^{j}\left\Vert
\left\vert \cdot\right\vert ^{j}D^{\beta}\phi_{0}\right\Vert _{\infty
}\left\vert x^{\left(  k\right)  }\right\vert ^{\left\vert \alpha\right\vert
-j}\\
&  \leq r_{k}^{-\left\vert \beta\right\vert }\sum\limits_{j\leq\left\vert
\alpha\right\vert }\binom{\left\vert \alpha\right\vert }{j}\left\Vert
\left\vert \cdot\right\vert ^{j}D^{\beta}\phi_{0}\right\Vert _{\infty
}\left\vert x^{\left(  k\right)  }\right\vert ^{\left\vert \alpha\right\vert
-j}\\
&  \leq\left\vert x^{\left(  k\right)  }\right\vert ^{m\left\vert
\beta\right\vert }\sum\limits_{j\leq\left\vert \alpha\right\vert }%
\binom{\left\vert \alpha\right\vert }{j}\left\Vert y^{j}D^{\beta}\phi
_{0}\right\Vert _{\infty}\left\vert x^{\left(  k\right)  }\right\vert
^{\left\vert \alpha\right\vert -j}\\
&  \leq\left(  \sum\limits_{j\leq\left\vert \alpha\right\vert }\binom
{\left\vert \alpha\right\vert }{j}\left\Vert y^{j}D^{\beta}\phi_{0}\right\Vert
_{\infty}\right)  \left\vert x^{\left(  k\right)  }\right\vert ^{m\left\vert
\beta\right\vert +\left\vert \alpha\right\vert },
\end{align*}

where the last two inequalities make use of \ref{a5.1}. Now \ref{a5.1} implies
$\left\vert x^{\left(  k\right)  }\right\vert \geq1$ so by virtue of
\ref{a5.2} we can now establish the upper bound%
\begin{align}
\left(  \int w\left\vert \phi_{k}\right\vert ^{2}\right)  ^{1/2}  & \leq
\sum_{\alpha,\beta}\sum\limits_{\gamma\leq\alpha}c_{\alpha,\beta}\binom
{\alpha}{\gamma}\left\Vert y^{\gamma}D^{\beta}\phi_{0}\right\Vert _{\infty
}\left\vert x^{\left(  k\right)  }\right\vert ^{m\left\vert \beta\right\vert
+\left\vert \alpha\right\vert -\left\vert \gamma\right\vert }\nonumber\\
& \leq\left(  \sum_{\alpha,\beta}\sum\limits_{\gamma\leq\alpha}c_{\alpha
,\beta}\binom{\alpha}{\gamma}\left\Vert y^{\gamma}D^{\beta}\phi_{0}\right\Vert
_{\infty}\right)  \left\vert x^{\left(  k\right)  }\right\vert ^{mb+a}%
,\label{a5.3}%
\end{align}

where $a=\max\left\vert \alpha\right\vert $ and $b=\max\left\vert
\beta\right\vert $.

Next we estimate a lower bound for $\left(  \int w\left\vert \phi
_{k}\right\vert ^{2}\right)  ^{1/2}$. Now $1/w\in S$ implies $\left\vert
\cdot\right\vert ^{2n}/w\in S$ so that $\left\vert \cdot\right\vert ^{2n}/w $
is bounded and there exists a constant $c_{n,d}>0$ such that%
\[
w\left(  x\right)  \geq c_{n}\left\vert x\right\vert ^{2n},\quad
x\in\mathbb{R}^{d}.
\]

Hence, since \ref{a5.1} implies $1\geq r_{k}\geq\left\vert x^{\left(
k\right)  }\right\vert ^{-m}$, we again use \ref{a5.0} to obtain the sequence
of lower bounds
\begin{align*}
\int w\left\vert \phi_{k}\right\vert ^{2}\geq c_{n}\int\left\vert
\cdot\right\vert ^{2n}\left\vert \phi_{k}\right\vert ^{2} &  =c_{n}%
\int\left\vert x\right\vert ^{2n}\left\vert \phi_{0}\left(  \frac{x-x^{\left(
k\right)  }}{r_{k}}\right)  \right\vert ^{2}dx\\
&  =c_{n}r_{k}^{d}\int_{\left\vert y\right\vert \leq1}\left\vert
r_{k}y+x^{\left(  k\right)  }\right\vert ^{2n}\left\vert \phi_{0}\left(
y\right)  \right\vert ^{2}dy\\
&  \geq c_{n}r_{k}^{d}\int_{\left\vert y\right\vert \leq1}\left(  \left\vert
x^{\left(  k\right)  }\right\vert -r_{k}\left\vert y\right\vert \right)
^{2n}\left\vert \phi_{0}\left(  y\right)  \right\vert ^{2}dy\\
&  \geq c_{n}r_{k}^{d}\int_{\left\vert y\right\vert \leq1}\left(  \left\vert
x^{\left(  k\right)  }\right\vert -1\right)  ^{2n}\left\vert \phi_{0}\left(
y\right)  \right\vert ^{2}dy\\
&  =c_{n}r_{k}^{d}\left\Vert \phi_{0}\right\Vert _{2}^{2}\left(  \left\vert
x^{\left(  k\right)  }\right\vert -1\right)  ^{2n}\\
&  \geq c_{n}\left\Vert \phi_{0}\right\Vert _{2}^{2}\left\vert x^{\left(
k\right)  }\right\vert ^{-md}\left(  \left\vert x^{\left(  k\right)
}\right\vert -1\right)  ^{2n}.
\end{align*}

Now $\left\vert x^{\left(  k\right)  }\right\vert \rightarrow\infty$ so we can
assume $\left\vert x^{\left(  k\right)  }\right\vert \geq2$ so that
$\left\vert x^{\left(  k\right)  }\right\vert -1\geq\left\vert x^{\left(
k\right)  }\right\vert /2$ and%
\begin{equation}
\int w\left\vert \phi_{k}\right\vert ^{2}\geq\frac{c_{n}}{4^{n}}\left\Vert
\phi_{0}\right\Vert _{2}^{2}\left\vert x^{\left(  k\right)  }\right\vert
^{2n-md},\quad\left\vert x^{\left(  k\right)  }\right\vert \geq2.\label{a5.6}%
\end{equation}

Hence when $\left\vert x^{\left(  k\right)  }\right\vert \geq2$ the bounds
\ref{a5.3} and \ref{a5.6} imply%
\[
\frac{c_{n}}{4^{n}}\left\Vert \phi_{0}\right\Vert _{2}^{2}\left\vert
x^{\left(  k\right)  }\right\vert ^{2n-md}\leq\left(  \sum_{\alpha,\beta}%
\sum\limits_{\gamma\leq\alpha}c_{\alpha,\beta}\binom{\alpha}{\gamma}\left\Vert
y^{\gamma}D^{\beta}\phi_{0}\right\Vert _{\infty}\right)  ^{2}\left\vert
x^{\left(  k\right)  }\right\vert ^{2mb+2a},
\]

or%
\[
\frac{c_{n}}{4^{n}}\left\Vert \phi_{0}\right\Vert _{2}^{2}\left\vert
x^{\left(  k\right)  }\right\vert ^{2n}\leq\left(  \sum_{\alpha,\beta}%
\sum\limits_{\gamma\leq\alpha}c_{\alpha,\beta}\binom{\alpha}{\gamma}\left\Vert
y^{\gamma}D^{\beta}\phi_{0}\right\Vert _{\infty}\right)  ^{2}\left\vert
x^{\left(  k\right)  }\right\vert ^{md+2mb+2a}.
\]

Since $\left\vert x^{\left(  k\right)  }\right\vert \rightarrow\infty$ if we
choose $2n>md+2mb+2a$ we easily obtain a contradiction by dividing both sides
by $\left\vert x^{\left(  k\right)  }\right\vert ^{2n}$ and letting
$k\rightarrow\infty$.
\end{proof}

We now have the very simple consequence:

\begin{theorem}
The Gaussian weight function $e^{\left\vert x\right\vert ^{2}}$ does not
satisfy Condition \ref{Cnd_wt_fn_T}.
\end{theorem}

This lemma will be used to show that the shifted thin-plate splines weight
functions $w$ do not satisfy Condition \ref{Cnd_wt_fn_T}.

\begin{lemma}
\label{Lem_condit_A_eq_0}Suppose $w$ is a weight function w.r.t. the minimal
set $\mathcal{A}=\left\{  0\right\}  $ and is such that: $w\in C^{\infty
}\left(  \mathbb{R}^{d}\setminus0\right)  $, $w>0$ on $\mathbb{R}^{d}%
\setminus0$ and $\phi_{w}:=\frac{\psi}{w}\in S$ for some $\psi\in C^{\infty}$
such that $0\leq\psi\leq1$, $\psi=0$ near the origin and $\psi\left(
x\right)  =1$ when $\left\vert x\right\vert \geq1$.

Then $w$ does not satisfy Condition \ref{Cnd_wt_fn_T}.
\end{lemma}

\begin{proof}
Suppose that $w$ does satisfy Condition \ref{Cnd_wt_fn_T} for some $T$ and a
set $\mathcal{A}\supseteq\left\{  0\right\}  $. Then there exists a positive,
linear combination of the seminorms \ref{a1.5} such that
\[
\left(  \int w\left\vert \phi\right\vert ^{2}\right)  ^{1/2}\leq\sum
_{\alpha,\beta}c_{\alpha,\beta}\left\Vert x^{\alpha}D^{\beta}\phi\right\Vert
_{\infty},\quad\phi\in T,
\]

and%
\[
C_{0}^{\infty}\left(  \mathbb{R}^{d}\setminus\mathcal{A}\right)  \subset T.
\]

Further, by part 2 of Remark \ref{Rem_Cnd_wt_fn_T} there exists a sequence a
functions $\phi_{k}\in C_{0}^{\infty}\left(  \mathbb{R}^{d}\setminus
\mathcal{A}\right)  $ defined by
\begin{equation}
\phi_{k}\left(  x\right)  =\phi_{0}\left(  \frac{x-x^{\left(  k\right)  }%
}{r_{k}}\right)  ,\quad\phi_{0}\in C_{0}^{\infty}\left(  B\left(  0;1\right)
\right)  ,\label{a5.7}%
\end{equation}

and such that%
\begin{equation}
\left\vert x^{\left(  k\right)  }\right\vert \rightarrow\infty,\quad1\geq
r_{k}\geq\left\vert x^{\left(  k\right)  }\right\vert ^{-m},\quad\left\vert
x^{\left(  k\right)  }\right\vert \geq1.\label{a5.8}%
\end{equation}

Hence $\operatorname*{supp}\phi_{k}\subset B\left(  x^{\left(  k\right)
};r_{k}\right)  \subset\mathbb{R}^{d}\setminus\mathcal{A}$. Thus we must have%
\[
\left(  \int w\left\vert \phi_{k}\right\vert ^{2}\right)  ^{1/2}\leq
\sum_{\alpha,\beta}c_{\alpha,\beta}\left\Vert x^{\alpha}D^{\beta}\phi
_{k}\right\Vert _{\infty},
\]

and now the form \ref{a5.7} of $\phi_{k}$ together with the inequalities
\ref{a5.8} imply%
\begin{align*}
\left\Vert x^{\alpha}D^{\beta}\phi_{k}\right\Vert _{\infty}=\left\Vert
x^{\alpha}D_{x}^{\beta}\phi_{0}\left(  \frac{x-x^{\left(  k\right)  }}{r_{k}%
}\right)  \right\Vert _{\infty} &  =r_{k}^{-\left\vert \beta\right\vert
}\left\Vert x^{\alpha}\left(  D_{x}^{\beta}\phi_{0}\right)  \left(
\frac{x-x^{\left(  k\right)  }}{r_{k}}\right)  \right\Vert _{\infty}\\
&  =r_{k}^{-\left\vert \beta\right\vert }\left\Vert \left(  r_{k}y+x^{\left(
k\right)  }\right)  ^{\alpha}\left(  D^{\beta}\phi_{0}\right)  \left(
y\right)  \right\Vert _{\infty}\\
&  \leq r_{k}^{-\left\vert \beta\right\vert }\left\Vert \left\vert
r_{k}y+x^{\left(  k\right)  }\right\vert ^{\left\vert \alpha\right\vert
}\left(  D^{\beta}\phi_{0}\right)  \left(  y\right)  \right\Vert _{\infty}\\
&  \leq r_{k}^{-\left\vert \beta\right\vert }\left\Vert \left(  r_{k}%
\left\vert y\right\vert +\left\vert x^{\left(  k\right)  }\right\vert \right)
^{\left\vert \alpha\right\vert }\left(  D^{\beta}\phi_{0}\right)  \left(
y\right)  \right\Vert _{\infty}\\
&  \leq r_{k}^{-\left\vert \beta\right\vert }\sum\limits_{j\leq\left\vert
\alpha\right\vert }\left\Vert \binom{\left\vert \alpha\right\vert }{j}%
r_{k}^{j}\left\vert y\right\vert ^{j}\left\vert x^{\left(  k\right)
}\right\vert ^{\left\vert \alpha\right\vert -j}D^{\beta}\phi_{0}\left(
y\right)  \right\Vert _{\infty}\\
&  \leq r_{k}^{-\left\vert \beta\right\vert }\sum\limits_{j\leq\left\vert
\alpha\right\vert }\binom{\left\vert \alpha\right\vert }{j}r_{k}^{j}\left\Vert
\left\vert \cdot\right\vert ^{j}D^{\beta}\phi_{0}\right\Vert _{\infty
}\left\vert x^{\left(  k\right)  }\right\vert ^{\left\vert \alpha\right\vert
-j}\\
&  \leq r_{k}^{-\left\vert \beta\right\vert }\sum\limits_{j\leq\left\vert
\alpha\right\vert }\binom{\left\vert \alpha\right\vert }{j}\left\Vert
\left\vert \cdot\right\vert ^{j}D^{\beta}\phi_{0}\right\Vert _{\infty
}\left\vert x^{\left(  k\right)  }\right\vert ^{\left\vert \alpha\right\vert
-j}\\
&  \leq\left\vert x^{\left(  k\right)  }\right\vert ^{m\left\vert
\beta\right\vert }\sum\limits_{j\leq\left\vert \alpha\right\vert }%
\binom{\left\vert \alpha\right\vert }{j}\left\Vert y^{j}D^{\beta}\phi
_{0}\right\Vert _{\infty}\left\vert x^{\left(  k\right)  }\right\vert
^{\left\vert \alpha\right\vert -j}\\
&  \leq\left(  \sum\limits_{j\leq\left\vert \alpha\right\vert }\binom
{\left\vert \alpha\right\vert }{j}\left\Vert y^{j}D^{\beta}\phi_{0}\right\Vert
_{\infty}\right)  \left\vert x^{\left(  k\right)  }\right\vert ^{m\left\vert
\beta\right\vert +\left\vert \alpha\right\vert },
\end{align*}

By virtue of \ref{a5.2} we can now establish the upper bound%
\begin{align}
\left(  \int w\left\vert \phi_{k}\right\vert ^{2}\right)  ^{1/2}  &
\leq\left(  \sum_{\alpha,\beta}\sum\limits_{j\leq\left\vert \alpha\right\vert
}\binom{\left\vert \alpha\right\vert }{j}\left\Vert y^{j}D^{\beta}\phi
_{0}\right\Vert _{\infty}\right)  \left\vert x^{\left(  k\right)  }\right\vert
^{m\left\vert \beta\right\vert +\left\vert \alpha\right\vert }\nonumber\\
& \leq\left(  \sum_{\alpha,\beta}\sum\limits_{j\leq\left\vert \alpha
\right\vert }\binom{\left\vert \alpha\right\vert }{j}\left\Vert y^{j}D^{\beta
}\phi_{0}\right\Vert _{\infty}\right)  \left\vert x^{\left(  k\right)
}\right\vert ^{mb+a},\label{a5.4}%
\end{align}

where $a=\max\left\vert \alpha\right\vert $ and $b=\max\left\vert
\beta\right\vert $.

Next we estimate a lower bound for $\left(  \int w\left\vert \phi
_{k}\right\vert ^{2}\right)  ^{1/2}$. Now $\psi/w\in S$ implies $\left\vert
\cdot\right\vert ^{2n}\psi/w\in S$ so that $\left\vert \cdot\right\vert
^{2n}\psi/w$ is bounded and there exists a constant $c_{n}>0$ such that%
\[
w\left(  x\right)  \geq c_{n}\left\vert x\right\vert ^{2n}\psi\left(
x\right)  ,\quad x\in\mathbb{R}^{d}.
\]

We again use \ref{a5.7} to obtain the lower bound%
\begin{align*}
\int w\left\vert \phi_{k}\right\vert ^{2}  & \geq c_{n}\int\left\vert
\cdot\right\vert ^{2n}\psi\left\vert \phi_{k}\right\vert ^{2}=c_{n}%
\int\left\vert x\right\vert ^{2n}\psi\left(  x\right)  \left\vert \phi
_{0}\left(  \frac{x-x^{\left(  k\right)  }}{r_{k}}\right)  \right\vert
^{2}dx\\
& =c_{n}r_{k}^{d}\int\limits_{\left\vert y\right\vert \leq1}\left\vert
r_{k}y+x^{\left(  k\right)  }\right\vert ^{2n}\psi\left(  r_{k}y+x^{\left(
k\right)  }\right)  \left\vert \phi_{0}\left(  y\right)  \right\vert ^{2}dy.
\end{align*}

But if $\left\vert y\right\vert \leq1$ and $\left\vert x^{\left(  k\right)
}\right\vert \geq2$ then%
\[
\left\vert r_{k}y+x^{\left(  k\right)  }\right\vert \geq\left\vert x^{\left(
k\right)  }\right\vert -r_{k}\left\vert y\right\vert \geq\left\vert x^{\left(
k\right)  }\right\vert -r_{k}\geq\left\vert x^{\left(  k\right)  }\right\vert
-1\geq\left\vert x^{\left(  k\right)  }\right\vert /2,
\]

and so if $\left\vert x^{\left(  k\right)  }\right\vert \geq2$%
\begin{align}
\int w\left\vert \phi_{k}\right\vert ^{2}  & \geq c_{n}r_{k}^{d}%
\int\limits_{\left\vert y\right\vert \leq1}\left\vert r_{k}y+x^{\left(
k\right)  }\right\vert ^{2n}\psi\left(  r_{k}y+x^{\left(  k\right)  }\right)
\left\vert \phi_{0}\left(  y\right)  \right\vert ^{2}dy\nonumber\\
& =\frac{c_{n}}{4^{n}}r_{k}^{d}\left\Vert \phi_{0}\right\Vert _{2}^{2}\left(
\left\vert x^{\left(  k\right)  }\right\vert /2\right)  ^{2n}\nonumber\\
& \geq\frac{c_{n}}{4^{n}}\left\Vert \phi_{0}\right\Vert _{2}^{2}\left\vert
x^{\left(  k\right)  }\right\vert ^{2n-md},\label{a5.5}%
\end{align}

where last step used \ref{a5.8}. Thus, if $\left\vert x^{\left(  k\right)
}\right\vert \geq2$
\[
\frac{c_{n}}{4^{n}}\left\Vert \phi_{0}\right\Vert _{2}^{2}\left\vert
x^{\left(  k\right)  }\right\vert ^{2n-md}\leq\left(  \sum_{\alpha,\beta}%
\sum\limits_{\gamma\leq\alpha}c_{\alpha,\beta}\binom{\alpha}{\gamma}\left\Vert
y^{\gamma}D^{\beta}\phi_{0}\right\Vert _{\infty}\right)  ^{2}\left\vert
x^{\left(  k\right)  }\right\vert ^{mb+a},
\]

or%
\[
\frac{c_{n}}{4^{n}}\left\Vert \phi_{0}\right\Vert _{2}^{2}\left\vert
x^{\left(  k\right)  }\right\vert ^{2n}\leq\left(  \sum_{\alpha,\beta}%
\sum\limits_{\gamma\leq\alpha}c_{\alpha,\beta}\binom{\alpha}{\gamma}\left\Vert
y^{\gamma}D^{\beta}\phi_{0}\right\Vert _{\infty}\right)  ^{2}\left\vert
x^{\left(  k\right)  }\right\vert ^{md+2mb+2a}%
\]

Since $\left\vert x^{\left(  k\right)  }\right\vert \rightarrow\infty$ if we
choose $2n>md+2mb+2a$ we easily obtain a contradiction by dividing both sides
by $\left\vert x^{\left(  k\right)  }\right\vert ^{2n}$ and letting
$k\rightarrow\infty$.
\end{proof}

Again we have a simple consequence:

\begin{theorem}
No shifted thin-plate spline weight function $w$ (Subsubsection
\ref{Ex_shft_thn_plt_spln_wt}) satisfies Condition \ref{Cnd_wt_fn_T}.
\end{theorem}

\begin{proof}
From Subsubsection \ref{Ex_shft_thn_plt_spln_wt}, $w\in C^{\infty}\left(
\mathbb{R}^{d}\setminus0\right)  $ has property W01 w.r.t. the set $\left\{
0\right\}  $ and $1/w$ decreases exponentially at infinity. Thus $\psi/w\in S
$ where $\psi$ was defined in Lemma \ref{Lem_condit_A_eq_0} i.e. $\psi\in
C^{\infty}$, $0\leq\psi\leq1$, $\psi=0$ near the origin and $\psi\left(
x\right)  =1$ when $\left\vert x\right\vert \geq1$.

It now follows from Lemma \ref{Lem_condit_A_eq_0} that $w$ does not satisfy
Condition \ref{Cnd_wt_fn_T}.
\end{proof}

\subsection{The space $\protect\widetilde{X}_{1/w}^{0}$}

Given a subspace $T\subset S$ we will define the spaces $\widehat{T}$,
$\overset{\vee}{T}$, $T^{\prime}$, $\left(  \widehat{T}\right)  ^{\prime}$,
$\left(  \overset{\vee}{T}\right)  ^{\prime}$ in a similar fashion to the
subspaces $\widehat{S}_{w,0}$, $\overset{\vee}{S}_{w,0}$ $S_{w,0}^{\prime}$,
$\left(  \widehat{S}_{w,0}\right)  ^{\prime}$, $\left(  \overset{\vee
}{S}_{w,0}\right)  ^{\prime}$ which were based on the subspace $S_{w,0}$ and
were introduced in Definition \ref{Def_Sw2_lin_fnal_0}:

\begin{definition}
\label{Def_T_lin_fnal_0}\textbf{The spaces }$\widehat{T}$, $\overset{\vee}{T}$
\textbf{and the spaces of continuous functionals} $T^{\prime}$, $\left(
\widehat{T}\right)  ^{\prime}$, $\left(  \overset{\vee}{T}\right)  ^{\prime}$.

The space $T$ was introduced in Definition \ref{Cnd_wt_fn_T}.

\begin{enumerate}
\item The spaces $\widehat{T}$ and $\overset{\vee}{T}$ are defined by
\[
\widehat{T}=\left\{  \widehat{\phi}:\phi\in T\right\}  ,\quad\overset{\vee
}{T}=\left\{  \overset{\vee}{\phi}:\phi\in T\right\}  .
\]

$\widehat{T}$ is topologized using the seminorms which define the topology of
$S$.

$\overset{\vee}{T}$ is topologized using the seminorms which define the
topology of $S$.

\item With these topologies the Fourier transform is now a homeomorphism from
$T$ to $\overset{\vee}{T}$ and the inverse-Fourier transform is now a
homeomorphism from $T$ to $\widehat{T}$.

\item Now define $T^{\prime}$, $\left(  \widehat{T}\right)  ^{\prime}$ and
$\left(  \overset{\vee}{T}\right)  ^{\prime}$ to be the spaces of continuous,
linear functionals on $T$, $\widehat{T}$ and $\overset{\vee}{T}$ respectively.

\item We define the Fourier transform on $T^{\prime}$ and inverse-Fourier
transform on $T^{\prime}$ by%
\begin{align*}
\left[  \widehat{f},\phi\right]   & =\left[  f,\widehat{\phi}\right]  ,\quad
f\in T^{\prime},\text{ }\phi\in\overset{\vee}{T},\\
\left[  \overset{\vee}{f},\phi\right]   & =\left[  f,\overset{\vee}{\phi
}\right]  ,\quad f\in T^{\prime},\text{ }\phi\in\widehat{T},
\end{align*}

so that
\[
\left(  T^{\prime}\right)  ^{\vee}=\left(  \widehat{T}\right)  ^{\prime}%
,\quad\left(  T^{\prime}\right)  ^{\wedge}=\left(  \overset{\vee}{T}\right)
^{\prime}.
\]

These equations imply that the Fourier transform is a homeomorphism from
$T^{\prime}$ to $\left(  \overset{\vee}{T}\right)  ^{\prime}$, and that the
inverse-Fourier transform is a homeomorphism from $T^{\prime}$ to $\left(
\widehat{T}\right)  ^{\prime}$.\medskip

\textbf{Localization}:\medskip

\item Suppose $\Omega\subset\mathbb{R}^{d}\setminus\mathcal{A}$ is open and
let $C_{0}^{\infty}\left(  \Omega\right)  =\left\{  \phi\in C_{0}^{\infty
}:\operatorname*{supp}\phi\subset\Omega\right\}  $. Then by part 2 of
Condition \ref{Cnd_wt_fn_T}, $C_{0}^{\infty}\left(  \Omega\right)  \subset T$.
We now say $u\in T^{\prime}$ is a member of $L_{loc}^{1}\left(  \Omega\right)
$ if there exists $f\in L_{loc}^{1}\left(  \Omega\right)  $ such that $\left[
u,\phi\right]  =\int f\phi$ for all $\phi\in C_{0}^{\infty}\left(
\Omega\right)  $.

\item Suppose $\Omega\subset\mathbb{R}^{d}\setminus\mathcal{A}$ is open and
$A$ is one of the spaces $T^{\prime}$, $\left(  \widehat{T}\right)  ^{\prime}$
or $\left(  \overset{\vee}{T}\right)  ^{\prime}$.

Then if $u,v\in A$ we say that $u=v$ on $\Omega$ if $\left[  u-v,\phi\right]
=0$ for all $\phi\in C_{0}^{\infty}\left(  \Omega\right)  $.

\item A member $u$ of $T^{\prime}$ has support on $\mathcal{A}$ iff $\left[
u,\phi\right]  =0$ when $\phi\in C_{0}^{\infty}$ and $\operatorname*{supp}%
\phi\subset\mathbb{R}^{d}\setminus\mathcal{A}$.
\end{enumerate}
\end{definition}

Condition \ref{Cnd_wt_fn_T} has the following implications cf. Theorem
\ref{Thm_Som_and_Sw2}:

\begin{theorem}
\label{Thm_T_and_sqrt(wt)L2}Suppose the weight function $w$ satisfies
Condition \ref{Cnd_wt_fn_T} for some $T\subset S$. Then:

\begin{enumerate}
\item $w\in L_{loc}^{1}\left(  \mathbb{R}^{d}\setminus\mathcal{A}\right)  $.

\item $g\in L^{2}$ implies $\sqrt{w}g\in T^{\prime}\cap L_{loc}^{1}\left(
\mathbb{R}^{d}\setminus\mathcal{A}\right)  $.

\item $\left(  C_{0}^{\infty}\left(  \mathbb{R}^{d}\setminus\mathcal{A}%
\right)  \right)  ^{\vee}\subset\overset{\vee}{T}\cap X_{w}^{0}\subset S\cap
X_{w}^{0}=\overset{\vee}{S}_{w,0}$.

\item If $w$ has property W02 or W03 then $\overset{\vee}{T}\cap X_{w}^{0}$ is
dense in $X_{w}^{0}$.
\end{enumerate}
\end{theorem}

\begin{proof}
\textbf{Part 1} Since $w$ is continuous on $\mathbb{R}^{d}\setminus
\mathcal{A}$ it follows directly that $w\in L_{loc}^{1}\left(  \mathbb{R}%
^{d}\setminus\mathcal{A}\right)  $.\medskip

\textbf{Part 2} Examination of the proof of part 2 of Lemma
\ref{Lem_wt_fn_extra_property_2} shows that $\sqrt{w}g\in T^{\prime}$.

Further, if $K\subset\mathbb{R}^{d}\setminus\mathcal{A}$ is compact, $w$ is
continuous on $K$ and $\int_{K}\sqrt{w}g\leq\left(  \int_{K}w\right)
^{1/2}\left\Vert g\right\Vert _{2}<\infty$.\medskip

\textbf{Part 3} $C_{0}^{\infty}\left(  \mathbb{R}^{d}\setminus\mathcal{A}%
\right)  \subset X_{w}^{0}$ because $w\in C^{\left(  0\right)  }\left(
\mathbb{R}^{d}\setminus\mathcal{A}\right)  $ and $C_{0}^{\infty}\left(
\mathbb{R}^{d}\setminus\mathcal{A}\right)  \subset T$ by part 2 of Condition
\ref{Cnd_wt_fn_T}. Thus $\left(  C_{0}^{\infty}\left(  \mathbb{R}^{d}%
\setminus\mathcal{A}\right)  \right)  ^{\vee}\subset\overset{\vee}{T}\cap
X_{w}^{0}$. Also, $T\subset S$ implies $\overset{\vee}{T}\subset S$ so
$\overset{\vee}{T}\cap X_{w}^{0}\subset S\cap X_{w}^{0}$. Finally $S\cap
X_{w}^{0}=\overset{\vee}{S}_{w,0}$ by part 1 of Definition
\ref{Def_Sw2_lin_fnal_0}.\medskip

\textbf{Part 4} By part 3, $\left(  C_{0}^{\infty}\left(  \mathbb{R}%
^{d}\setminus\mathcal{A}\right)  \right)  ^{\vee}\subset\overset{\vee}{T}\cap
X_{w}^{0}\subset X_{w}^{0}$ and in Corollary
\ref{Cor_Xow_invFCoinf(R/A)_density} we showed that $\left(  C_{0}^{\infty
}\left(  \mathbb{R}^{d}\setminus\mathcal{A}\right)  \right)  ^{\vee}$ is dense
in $X_{w}^{0}$.
\end{proof}

The next result is the analogue of Theorem \ref{Thm_w_phi_in_fnalSw2}.

\begin{theorem}
\label{Thm_w_phi_in_fnalT}Suppose that $w$ has property W01 w.r.t. the set
$\mathcal{A}$ and that $w$ and $T\subset S$ satisfy Condition
\ref{Cnd_wt_fn_T}. Then:
\end{theorem}

\begin{enumerate}
\item Suppose $\phi\in S_{w,0}$. Then $w\phi\in L_{loc}^{1}\left(
\mathbb{R}^{d}\setminus\mathcal{A}\right)  $. Further, $\phi\in T$ implies
$w\phi\in T^{\prime}$.

\item Suppose $\phi\in S_{w,0}$. Then $w\widehat{\phi}\in L_{loc}^{1}\left(
\mathbb{R}^{d}\setminus\mathcal{A}\right)  $. Further, $\phi\in\overset{\vee
}{T}$ implies $\left(  w\widehat{\phi}\right)  ^{\vee}\in\left(
\widehat{T}\right)  ^{\prime}$.

In addition, if $\phi\in\overset{\vee}{T}$ then $u=\left(  w\widehat{\phi
}\right)  ^{\vee}$ satisfies
\[
\int w\left\vert \widehat{\phi}\right\vert ^{2}=\int\frac{\left\vert
\widehat{u}\right\vert ^{2}}{w}.
\]

\end{enumerate}

\begin{proof}
\textbf{Part 1} Clearly $w\phi\in L_{loc}^{1}\left(  \mathbb{R}^{d}%
\setminus\mathcal{A}\right)  $ if $\phi\in S_{w,0}$. Now suppose $\phi\in T$.
Clearly $w\phi\in L_{loc}^{1}\left(  \mathbb{R}^{d}\setminus\mathcal{A}%
\right)  $ since $w\in L_{loc}^{1}\left(  \mathbb{R}^{d}\setminus
\mathcal{A}\right)  $. Further, if $\psi\in T$ then
\[
\left[  w\phi,\psi\right]  =\int w\phi\psi\leq\int\sqrt{w}\phi\sqrt{w}\psi
\leq\left(  \int w\left\vert \phi\right\vert ^{2}\right)  ^{1/2}\left(  \int
w\left\vert \psi\right\vert ^{2}\right)  ^{1/2}\leq\left(  \int w\left\vert
\phi\right\vert ^{2}\right)  ^{1/2}\left\vert \psi\right\vert _{\sigma},
\]

so that $w\phi\in T^{\prime}$.\medskip

\textbf{Part 2} A simple consequence of part 1.
\end{proof}

\begin{definition}
\label{Def_tildXo1/w_with_T}\textbf{The semi-inner product space
}$\widetilde{X}_{1/w}^{0}$

Suppose $w$ has property W01 w.r.t. the set $\mathcal{A}$ and satisfies
Condition \ref{Cnd_wt_fn_T}. Then we define the semi-inner product space:
\[
\widetilde{X}_{1/w}^{0}=\left\{  u\in\left(  \widehat{T}\right)  ^{\prime
}:u_{F}\in L_{loc}^{1}\left(  \mathbb{R}^{d}\setminus\mathcal{A}\right)
\text{ }and\text{ }\int\frac{\left\vert u_{F}\right\vert ^{2}}{w}%
<\infty\right\}  ,
\]

where $u_{F}$ is the restriction of $\widehat{u}\in T^{\prime}$ to
$\mathbb{R}^{d}\setminus\mathcal{A}$. Endow $\widetilde{X}_{1/w}^{0}$ with the
seminorm and semi-inner product%
\[
\left\vert u\right\vert _{w,0}=\left(  \int\frac{\left\vert u_{F}\right\vert
^{2}}{w}\right)  ^{1/2},\quad\left\langle u,v\right\rangle _{w,0}=\int%
\frac{u_{F}\overline{v_{F}}}{w}.
\]

Part 2 of Theorem \ref{Thm_w_phi_in_fnalT} below will imply that
$\widetilde{X}_{1/w}^{0}$ is not empty: in fact%
\begin{equation}
\left\{  \left(  w\widehat{\phi}\right)  ^{\vee}:\phi\in\overset{\vee
}{T}\right\}  \equiv\left(  wT\right)  ^{\vee}\subset\widetilde{X}_{1/w}%
^{0}.\label{a3.9}%
\end{equation}

\end{definition}

The next result is the analogue of Theorem \ref{Thm_property1_tildXo1/w}.

\begin{theorem}
\label{Thm_property_tildXo1/w_T}\textbf{Some properties of }$\widetilde{X}%
_{1/w}^{0}$

\begin{enumerate}
\item $\operatorname{null}\left\vert \cdot\right\vert _{1/w,0}=\left(
T_{\mathcal{A}}^{\prime}\right)  ^{\vee}$ where%
\begin{equation}
T_{\mathcal{A}}^{\prime}=\left\{  u\in T^{\prime}:\operatorname*{supp}%
u\subset\mathcal{A}\text{ }as\text{ }a\text{ }distribution\right\}
.\label{a4.0}%
\end{equation}

\item Noting part 1 of Remark \ref{Rem_Def_extend_wt_fn}, $\left\vert
\cdot\right\vert _{1/w,0}$ is a norm iff we can choose $\mathcal{A}$ to be the
(minimal) empty set.

If, in addition, $w$ has property W02 or W03 for $\kappa$ then:

\item If $\left\vert v\left(  \xi\right)  \right\vert \leq c\left(
1+\left\vert \xi\right\vert \right)  ^{\kappa}$ then $\overset{\vee}{v}%
\in\widetilde{X}_{1/w}^{0}$.

\item $\overset{\vee}{T}\subset X_{w}^{0}$, $S\subset\widetilde{X}_{1/w}^{0}$
and $\left(  wT\right)  ^{\vee}\subset\widetilde{X}_{1/w}^{0}$.

\item $u\in\widetilde{X}_{1/w}^{0}$ implies $u_{F}\in T^{\prime}\cap
L_{loc}^{1}\left(  \mathbb{R}^{d}\setminus\mathcal{A}\right)  $, $\left(
u_{F}\right)  ^{\vee}\in\widetilde{X}_{1/w}^{0}$ and $u-\left(  u_{F}\right)
^{\vee}\in\left(  T_{\mathcal{A}}^{\prime}\right)  ^{\vee}$.

\item If $u\in\widetilde{X}_{1/w}^{0}$ then $\frac{\left(  D^{\alpha}u\right)
_{F}}{w}\in L^{1}$ when $\left\vert \alpha\right\vert \leq\kappa$.

Also, $\int\frac{\left(  D^{\alpha}u\right)  _{F}}{w}\in\left(  \widetilde{X}%
_{1/w}^{0}\right)  ^{\prime}$ and $\left\Vert \int\frac{\left(  D^{\alpha
}u\right)  _{F}}{w}\right\Vert _{op}\leq\left(  \int\frac{\xi^{2\alpha}}%
{w}\right)  ^{1/2}$.
\end{enumerate}
\end{theorem}

\begin{proof}
\textbf{Part 1}
\begin{align}
&  \operatorname{null}\left\vert \cdot\right\vert _{1/w,0}\nonumber\\
&  =\left\{  v\in\widetilde{X}_{1/w}^{0}:\left\vert v\right\vert
_{1/w,0}=0\right\} \nonumber\\
&  =\left\{  v\in\left(  \widehat{T}\right)  ^{\prime}:\widehat{v}\in
T^{\prime}\text{ }and\text{ }\widehat{v}=0\text{ }on\text{ }\mathbb{R}%
^{d}\setminus\mathcal{A}\text{ }as\text{ }a\text{ }distribution\right\}
\nonumber\\
&  =\left\{  \overset{\vee}{u}:u\in T^{\prime}\text{ }and\text{ }u=0\text{
}on\text{ }\mathbb{R}^{d}\setminus\mathcal{A}\text{ }as\text{ }a\text{
}distribution\right\} \label{a3.8}\\
&  =\left\{  \overset{\vee}{u}:u\in T^{\prime}\text{ }and\text{ }%
\operatorname*{supp}u\subset\mathcal{A}\text{ }as\text{ }a\text{
}distribution\right\} \nonumber\\
&  =\left(  T_{\mathcal{A}}^{\prime}\right)  ^{\vee}.\nonumber
\end{align}
\medskip

\textbf{Part 2} If $\mathcal{A}$ is empty then \ref{a3.8} implies that
$\operatorname{null}\left\vert \cdot\right\vert _{1/w,0}=\left\{  0\right\}
$. If $\operatorname{null}\left\vert \cdot\right\vert _{1/w,0}=\left\{
0\right\}  $ then it follows that $u\in T^{\prime}$ and $u=0$ on
$\mathbb{R}^{d}\setminus\mathcal{A}$ implies $u=0$ on $\mathbb{R}^{d}$. Now if
$\mathcal{A}$ is not empty and $a\in\mathcal{A}$ then $\delta\left(
\cdot-a\right)  \in S^{\prime}$ so $\delta\left(  \cdot-a\right)  \in
T^{\prime}$ and $\delta\left(  \cdot-a\right)  =0$ on $\mathbb{R}^{d}%
\setminus\mathcal{A}$. But $\delta\left(  \cdot-a\right)  \neq0$ on
$\mathbb{R}^{d}$ which is a contradiction. Thus $\mathcal{A}$ must be
empty.\medskip

\textbf{Part 3} Clearly $v\in S^{\prime}$ so $\overset{\vee}{v}\in S^{\prime}%
$. Also, $\left(  \overset{\vee}{v}\right)  _{F}=v\in L_{loc}^{1}$ and
property W02 or W03 implies $\int\frac{\left\vert \left(  \overset{\vee
}{v}\right)  _{F}\right\vert ^{2}}{w}<\infty$.\medskip

\textbf{Part 4} It was shown in part 2 of Theorem \ref{Thm_properties2_B} that
$\overset{\vee}{S}_{w,0}\subset X_{w}^{0}$ so $\overset{\vee}{T}%
\subset\overset{\vee}{S}_{w,0}\subset X_{w}^{0}$. That $S\subset
\widetilde{X}_{1/w}^{0}$ follows directly from part 2. That $\left(
wT\right)  ^{\vee}\subset\widetilde{X}_{1/w}^{0}$ is just \ref{a3.9}.\medskip

\textbf{Part 5} By definition, $u\in\widetilde{X}_{1/w}^{0}$ implies
$u\in\left(  \widehat{T}\right)  ^{\prime}$, $u_{F}\in L_{loc}^{1}\left(
\mathbb{R}^{d}\setminus\mathcal{A}\right)  $ and $\int\frac{\left\vert
u_{F}\right\vert ^{2}}{w}<\infty$.

Thus $\widehat{u}\in T^{\prime}$. If $\phi\in\widehat{T}$ then by the
Cauchy-Schwartz inequality%
\[
\left\vert \int u_{F}\phi\right\vert \leq\int\left\vert u_{F}\right\vert
\left\vert \phi\right\vert =\int\frac{\left\vert u_{F}\right\vert }{\sqrt{w}%
}\sqrt{w}\left\vert \phi\right\vert \leq\left\vert u\right\vert _{1/w,0}%
\left(  \int w\left\vert \phi\right\vert ^{2}\right)  ^{1/2},
\]

and by \ref{a4.1} $u_{F}\in T^{\prime}$. We can now conclude that $\left(
u_{F}\right)  ^{\vee}\in\widehat{T}$, $\left(  \left(  u_{F}\right)  ^{\vee
}\right)  _{F}=u_{F}\in L_{loc}^{1}\left(  \mathbb{R}^{d}\setminus
\mathcal{A}\right)  $ and $\int\frac{\left\vert \left(  \left(  u_{F}\right)
^{\vee}\right)  _{F}\right\vert ^{2}}{w}<\infty$ so that $\left(
u_{F}\right)  ^{\vee}\in\widetilde{X}_{1/w}^{0}$ and
\[
\left\vert u-\left(  u_{F}\right)  ^{\vee}\right\vert _{1/w,0}^{2}=\int%
\frac{\left\vert \left(  u-\left(  u_{F}\right)  ^{\vee}\right)
_{F}\right\vert ^{2}}{w}=\int\frac{\left\vert u_{F}-\left(  \left(
u_{F}\right)  ^{\vee}\right)  _{F}\right\vert ^{2}}{w}=0,
\]

i.e. $u-\left(  u_{F}\right)  ^{\vee}\in\left(  T_{\mathcal{A}}^{\prime
}\right)  ^{\vee}$.\medskip

\textbf{Part 6} $\int\frac{\left\vert \left(  D^{\alpha}u\right)
_{F}\right\vert }{w}=\int\frac{\left\vert \xi^{\alpha}u_{F}\right\vert }%
{w}=\int\frac{\left\vert \xi^{\alpha}\right\vert }{\sqrt{w}}\frac{\left\vert
u_{F}\right\vert }{\sqrt{w}}\leq\left(  \int\frac{\xi^{2\alpha}}{w}\right)
^{1/2}\left(  \int\frac{\left\vert u_{F}\right\vert ^{2}}{w}\right)
^{1/2}=\left(  \int\frac{\left\vert \cdot\right\vert ^{2\left\vert
\alpha\right\vert }}{w}\right)  ^{1/2}\left\vert u\right\vert _{1/w,0}<\infty$.
\end{proof}

\subsection{The operators $\protect\widetilde{\mathcal{M}}_{2}%
:\protect\widetilde{X}_{1/w}^{0}\rightarrow L^{2}$ and
$\protect\widetilde{\mathcal{L}}_{3}:L^{2}\rightarrow\protect\widetilde{X}%
_{1/w}^{0}$}

The next definition is the analogue of Definition \ref{Def_tildM1}.

\begin{definition}
\label{Def_tildM2}\textbf{The operator} $\widetilde{\mathcal{M}}%
_{2}:\widetilde{X}_{1/w}^{0}\rightarrow L^{2}$

From the definition of $\widetilde{X}_{1/w}^{0}$, $u\in\widetilde{X}_{1/w}%
^{0}$ implies $\frac{u_{F}}{\sqrt{w}}\in L^{2}$ where $u_{F}$ is the
restriction of $\widehat{u}$ to $\mathbb{R}^{d}\setminus\mathcal{A}$. We can
now define the linear mapping $\widetilde{\mathcal{M}}_{2}:\widetilde{X}%
_{1/w}^{0}\rightarrow L^{2}$ by%
\[
\widetilde{\mathcal{M}}_{2}u=\left(  \frac{u_{F}}{\sqrt{w}}\right)  ^{\vee
},\text{\quad}u\in\widetilde{X}_{1/w}^{0}.
\]

\end{definition}

The operator $\widetilde{\mathcal{M}}_{2}$ has the following properties:

\begin{theorem}
\label{Thm_tildM2_property}\textbf{The operator }$\widetilde{\mathcal{M}}%
_{2}:\widetilde{X}_{1/w}^{0}\rightarrow L^{2}$ is linear and in the seminorm
sense it is isometric, 1-1 with null space $\left(  T_{\mathcal{A}}^{\prime
}\right)  ^{\vee}$.

Also, $\tau_{a}\widetilde{\mathcal{M}}_{2}=\widetilde{\mathcal{M}}_{2}\tau
_{a}$ where $\tau_{a}$ is the translation operator $\tau_{a}u=u\left(
\cdot-a\right)  $.
\end{theorem}

\begin{proof}
That $\widetilde{\mathcal{M}}_{2}$ is an isometry is clear from the definition
of $\widetilde{X}_{1/w}^{0}$. Since $\widetilde{\mathcal{M}}_{2}$ is an
isometry the null space of $\widetilde{\mathcal{M}}_{2}$ is the null space of
the seminorm $\left\vert \cdot\right\vert _{1/w,0}$, namely $\left(
T_{\mathcal{A}}^{\prime}\right)  ^{\vee}$. Finally%
\[
\tau_{a}\widetilde{\mathcal{M}}_{2}u=\tau_{a}\left(  \frac{u_{F}}{\sqrt{w}%
}\right)  ^{\vee}=\left(  e^{-ia\xi}\frac{u_{F}}{\sqrt{w}}\right)  ^{\vee
}=\left(  \frac{\left(  \tau_{a}u\right)  _{F}}{\sqrt{w}}\right)  ^{\vee
}=\widetilde{\mathcal{M}}_{2}\tau_{a}u.
\]

\end{proof}

The next definition is the analogue of Definition \ref{Def_tildL}.

\begin{definition}
\label{Def_tildL3}\textbf{The operator} $\widetilde{\mathcal{L}}_{3}%
:L^{2}\rightarrow\left(  \widehat{T}\right)  ^{\prime}$

Suppose the weight function $w$ satisfies Condition \ref{Cnd_wt_fn_T}. Then
noting part 2 of Theorem \ref{Thm_T_and_sqrt(wt)L2} we define the linear
operator $\widetilde{\mathcal{L}}_{3}$ uniquely by%
\begin{equation}
\widetilde{\mathcal{L}}_{3}g=\left(  \sqrt{w}\widehat{g}\right)  ^{\vee},\quad
g\in L^{2},\label{a4.2}%
\end{equation}

so that $\left(  \sqrt{w}\widehat{g}\right)  ^{\vee}\in\left(  T^{\prime
}\right)  ^{\vee}=\left(  \widehat{T}\right)  ^{\prime}$.
\end{definition}

The next result is the analogue of Theorem \ref{Thm_tildL_property}:

\begin{theorem}
\label{Thm_tildL3_property}\textbf{Properties of }$\widetilde{\mathcal{L}}%
_{3}$

\begin{enumerate}
\item If $g\in L^{2}$ then $\left(  \widetilde{\mathcal{L}}_{3}g\right)
_{F}=\sqrt{w}\widehat{g}\in L_{loc}^{1}\left(  \mathbb{R}^{d}\setminus
\mathcal{A}\right)  $.

\item $\widetilde{\mathcal{L}}_{3}:L^{2}\rightarrow\widetilde{X}_{1/w}^{0}$.

\item $\widetilde{\mathcal{L}}_{3}$ is an isometry and 1-1.

\item If $\tau_{a}$ is the translation operator $\tau_{a}u=u\left(
\cdot-a\right)  $ then $\tau_{a}\widetilde{\mathcal{L}}_{3}%
=\widetilde{\mathcal{L}}_{3}\tau_{a}$.
\end{enumerate}
\end{theorem}

\begin{proof}
\textbf{Part 1} Suppose $g\in L^{2}$. From the definition of
$\widetilde{\mathcal{L}}_{3}$, $\widehat{\widetilde{\mathcal{L}}_{3}g}%
=\sqrt{w}\widehat{g}$. Next observe that because $w\in C^{\left(  0\right)
}\left(  \mathbb{R}^{d}\setminus\mathcal{A}\right)  $ and Condition
\ref{Cnd_wt_fn_T} implies $C_{0}^{\infty}\left(  \mathbb{R}^{d}\setminus
\mathcal{A}\right)  \subset T$, it follows that
\[
\left[  \widehat{\widetilde{\mathcal{L}}_{3}g},\phi\right]  =\left[  \sqrt
{w}\widehat{g},\phi\right]  ,\quad\phi\in C_{0}^{\infty}\left(  \mathbb{R}%
^{d}\setminus\mathcal{A}\right)  .
\]

But by part 2 of Theorem \ref{Thm_T_and_sqrt(wt)L2}, $\sqrt{w}\widehat{g}\in
L_{loc}^{1}\left(  \mathbb{R}^{d}\setminus\mathcal{A}\right)  $ and so
$\left(  \widetilde{\mathcal{L}}_{3}g\right)  _{F}=\sqrt{w}\widehat{g}\in
L_{loc}^{1}\left(  \mathbb{R}^{d}\setminus\mathcal{A}\right)  $. But if
$K\subset\mathbb{R}^{d}\setminus\mathcal{A}$ is compact, $w\in C^{\left(
0\right)  }\left(  \mathbb{R}^{d}\setminus\mathcal{A}\right)  $ implies
\[
\int_{K}\sqrt{w}\left\vert \widehat{g}\right\vert \leq\left(  \int%
_{K}w\right)  ^{1/2}\left\Vert g\right\Vert _{2}\leq\max_{K}\left(  w\right)
\left(  \int_{K}1\right)  \left\Vert g\right\Vert _{2}.
\]

Thus $\widehat{\widetilde{\mathcal{L}}_{3}g}\in L_{loc}^{1}\left(
\mathbb{R}^{d}\setminus\mathcal{A}\right)  $ and $\left(
\widetilde{\mathcal{L}}_{3}g\right)  _{F}=\sqrt{w}\widehat{g}$.\medskip

The rest of the proof is identical to that of Theorem \ref{Thm_tildL_property}.
\end{proof}

The following theorem indicates how the operators $\widetilde{\mathcal{L}}%
_{3}:L^{2}\rightarrow\widetilde{X}_{1/w}^{0}$ and $\widetilde{\mathcal{M}}%
_{2}:\widetilde{X}_{1/w}^{0}\rightarrow L^{2}$ interact.

\begin{theorem}
\label{Thm_tildM2_tildL3_property}Suppose the weight function $w$ satisfies
Condition \ref{Cnd_wt_fn_T}. Then:

\begin{enumerate}
\item $\widetilde{\mathcal{M}}_{2}\widetilde{\mathcal{L}}_{3}=I$ on $L^{2}$.

\item $\widetilde{\mathcal{L}}_{3}\widetilde{\mathcal{M}}_{2}u-u\in\left(
T_{\mathcal{A}}^{\prime}\right)  ^{\vee}$ when $u\in X_{1/w}^{0}$ i.e.
$\left\vert \widetilde{\mathcal{L}}_{3}\widetilde{\mathcal{M}}_{2}%
u-u\right\vert _{1/w,0}=0$.

\item $\widetilde{\mathcal{M}}_{2}:\widetilde{X}_{1/w}^{0}\rightarrow L^{2}$
is onto.

\item $\widetilde{\mathcal{L}}_{3}:L^{2}\rightarrow\widetilde{X}_{1/w}^{0}$ is
onto in the seminorm sense.

\item $\widetilde{\mathcal{M}}_{2}$ and $\widetilde{\mathcal{L}}_{3}$ are adjoints.
\end{enumerate}
\end{theorem}

\begin{proof}
Easy adaptation of Theorem \ref{Thm_tildM1_tildL_property}.
\end{proof}

Since $L^{2}$ is complete, the mappings of the previous theorem will yield the
following important result:

\begin{corollary}
\label{Cor_tildXo1/w_semiHilb_L3_M2}Suppose the weight function $w$ satisfies
Condition \ref{Cnd_wt_fn_T}. Then in general $\widetilde{X}_{1/w}^{0}%
$\textbf{\ is a semi-Hilbert space}. Indeed, $\widetilde{X}_{1/w}^{0}$ is a
Hilbert space iff $\mathcal{A}$ is empty.
\end{corollary}

\begin{proof}
Proof of Corollary \ref{Cor_tildXo1/w_semiHilb_Sw2} with $T_{\mathcal{A}%
}^{\prime}$ replaced by $S_{w,0;\mathcal{A}}^{\prime}$.
\end{proof}

\subsection{The operator $\protect\widetilde{\mathcal{B}}_{2}:X_{w}^{0}%
\otimes\protect\widetilde{X}_{1/w}^{0}\rightarrow C_{B}^{\left(  0\right)  }$}

We now define the analogue $\widetilde{\mathcal{B}}_{2}$ of the bilinear
operator $\widetilde{\mathcal{B}}_{1}$ which was introduced in Definition
\ref{Def_tildB1}:

\begin{definition}
\label{Def_tildB2}\textbf{The mapping }$\widetilde{\mathcal{B}}_{2}$

If $\left(  u,v\right)  \in X_{w}^{0}\otimes\widetilde{X}_{1/w}^{0}$ define
the bilinear mapping $\widetilde{\mathcal{B}}_{2}:X_{w}^{0}\otimes
\widetilde{X}_{1/w}^{0}\rightarrow C_{B}^{\left(  0\right)  }$ by:%
\begin{equation}
\widetilde{\mathcal{B}}_{2}\left(  u,v\right)  =\left(  \widehat{u}%
v_{F}\right)  ^{\vee},\quad u\in X_{w}^{0},\text{ }v\in\widetilde{X}_{1/w}%
^{0}.\label{a4.3}%
\end{equation}

Noting that the product of two measurable functions is a measurable function
definition \ref{a4.3} makes sense since from the definitions of $\widetilde{X}%
_{1/w}^{0}$ and $X_{w}^{0}$, $v_{F}\in L_{loc}^{1}\left(  \mathbb{R}%
^{d}\setminus\mathcal{A}\right)  $, $\frac{v_{F}}{\sqrt{w}}\in L^{2}$,
$\widehat{u}\in L_{loc}^{1}$ and $\sqrt{w}\widehat{u}\in L^{2}$ so that
$\int\left\vert \widehat{u}v_{F}\right\vert \leq\left\Vert u\right\Vert
_{w,0}\left\vert v\right\vert _{1/w,0}<\infty$. Thus $\widehat{u}v_{F}\in
L^{1}$ and hence $\left(  \widehat{u}v_{F}\right)  ^{\vee}\in C_{B}^{\left(
0\right)  }$.
\end{definition}

\begin{theorem}
\label{Thm_property_tildB2}\textbf{Properties of the operator }%
$\widetilde{\mathcal{B}}_{2}$

\begin{enumerate}
\item $\widetilde{\mathcal{B}}_{2}:X_{w}^{0}\otimes\widetilde{X}_{1/w}%
^{0}\rightarrow C_{B}^{\left(  0\right)  }$ is a continuous bilinear mapping
when $C_{B}^{\left(  0\right)  }$ is endowed with the supremum norm
$\left\Vert \cdot\right\Vert _{\infty}$. In fact%
\begin{equation}
\left\Vert \widetilde{\mathcal{B}}_{2}\left(  u,v\right)  \right\Vert
_{\infty}\leq\left(  2\pi\right)  ^{-\frac{d}{2}}\left\Vert u\right\Vert
_{w,0}\left\vert v\right\vert _{1/w,0},\quad u\in X_{w}^{0},\text{ }%
v\in\widetilde{X}_{1/w}^{0}.
\end{equation}

\item In the sense of distributions: if $\alpha+\beta=\gamma\geq0$ then
\[
D^{\gamma}\widetilde{\mathcal{B}}_{2}\left(  u,v\right)  =\left(
\widehat{D^{\alpha}u}\left(  D^{\beta}v\right)  _{F}\right)  ^{\vee},\quad
u\in X_{w}^{0},\text{ }v\in\widetilde{X}_{1/w}^{0}.
\]

\item The operator $\widetilde{\mathcal{B}}_{2}$ commutes with the
(distribution) translation operator $\tau_{a}f=f\left(  \cdot-a\right)  $,
$a\in\mathbb{R}^{d}$ in the sense that%
\[
\tau_{a}\widetilde{\mathcal{B}}_{2}\left(  u,v\right)  =\widetilde{\mathcal{B}%
}_{2}\left(  \tau_{a}u,v\right)  =\widetilde{\mathcal{B}}_{2}\left(
u,\tau_{a}v\right)  .
\]

\item We have%
\[
\widetilde{\mathcal{B}}_{2}\left(  u,v\right)  =\left(  \mathcal{I}u\right)
\ast\widetilde{\mathcal{M}}_{2}v,\quad u\in X_{w}^{0},\text{ }v\in
\widetilde{X}_{1/w}^{0},
\]

where $\mathcal{I}$ was introduced in Definition \ref{Def_I_J} and
$\widetilde{\mathcal{M}}_{2}$ in Definition \ref{Def_tildM2}.\medskip

Now suppose $w$ also satisfies property W02 or W03 for parameter $\kappa$.
Then:\medskip

\item $\widetilde{\mathcal{B}}_{2}$ is a convolution in the sense of
\ref{1.20}: $X_{w}^{0}\subset S^{\prime}$, $S\subset\widetilde{X}_{1/w}^{0}$
and%
\[
\widetilde{\mathcal{B}}_{2}\left(  u,v\right)  =u\ast v,\quad u\in X_{w}%
^{0},\text{ }v\in S.
\]

Further, $\widetilde{\mathcal{B}}_{2}\left(  u,v\right)  \in C_{B}^{\infty}$
and $D^{\gamma}\widetilde{\mathcal{B}}_{2}\left(  u,v\right)  =D^{\alpha}u\ast
D^{\beta}v$ for all $\gamma=\alpha+\beta$.
\end{enumerate}
\end{theorem}

\begin{proof}
Identical to the proof of Theorem \ref{Thm_property_tildB1}.
\end{proof}

Part 5 of the last theorem shows how $\widetilde{\mathcal{B}}_{2}$ can be
regarded as a convolution by restricting the domain of the second variable to
$S$.

The next result considers the important case where the first argument of
$\widetilde{\mathcal{B}}_{2}$ is the basis function.

\begin{corollary}
\label{Cor_property_tildB2_basis}\textbf{Properties of }%
$\widetilde{\mathcal{B}}_{2}\left(  G\mathbf{,\cdot}\right)  $. Suppose $w$
has property W02 or W03 for $\kappa$ and $G$ is the basis function. Then if
$v,v^{\prime}\in X_{1/w}^{0}$:

\begin{enumerate}
\item $D^{\gamma}G\in X_{w}^{0}$ when $\left\vert \gamma\right\vert \leq
\kappa$ and%
\[
D^{\gamma}\widetilde{\mathcal{B}}_{2}\left(  G,v\right)
=\widetilde{\mathcal{B}}_{2}\left(  D^{\gamma}G,v\right)  ,\quad\left\vert
\gamma\right\vert \leq\kappa.
\]

\item $\widetilde{\mathcal{B}}_{2}\left(  G,v\right)  \in X_{w}^{0}$ and
$\left\Vert \widetilde{\mathcal{B}}_{2}\left(  G,v\right)  \right\Vert
_{w,0}=\left\vert v\right\vert _{1/w,0}$.

\item $\widetilde{\mathcal{B}}_{2}\left(  G,\phi\right)  =G\ast\phi$ when
$\phi\in S$.

\item $\int F\widetilde{\mathcal{B}}_{2}\left(  \widetilde{\mathcal{B}}%
_{2}\left(  G,v\right)  ,\overline{v^{\prime}}\right)  =\left\langle
v,v_{\ast}^{\prime}\right\rangle _{1/w,0}$, where $F$ is the Fourier transform
and $v_{\ast}^{\prime}\left(  x\right)  =v^{\prime}\left(  -x\right)  $.
\end{enumerate}
\end{corollary}

\begin{proof}
Proof identical to the proof of Corollary \ref{Cor_property_tildB1_basis}.
\end{proof}

Our extension of the convolution Definition \ref{Def_Sw2_convol_fnalSw2} is:

\begin{definition}
\label{Def_invFSb_convol_fnalFSb}\textbf{The convolution }$\overset{\vee
}{T}\ast\left(  \widehat{T}\right)  ^{\prime}\rightarrow\left(  \widehat{T}%
\right)  ^{\prime}$

Suppose $T\subset S$. If $f\in\left(  \widehat{T}\right)  ^{\prime}$ and
$\phi\in\overset{\vee}{T}$ then $\widehat{f}=T^{\prime}$ and $\widehat{\phi
}\in S$ so that $\widehat{\phi}\widehat{f}\in T^{\prime}$ and $\left(
\widehat{\phi}\widehat{f}\right)  ^{\vee}\in\left(  \widehat{T}\right)
^{\prime}$. Now by analogy with \ref{a2.0} we define%
\begin{equation}
\phi\ast f=\left(  \widehat{\phi}\widehat{f}\right)  ^{\vee},\quad\phi
\in\overset{\vee}{T},\text{ }f\in\left(  \widehat{T}\right)  ^{\prime
}.\label{a4.4}%
\end{equation}

\end{definition}

The next convolution result is the analogue of Lemma \ref{Lem_convol_Sw2}.

\begin{lemma}
\label{Lem_convol_T}The convolution of Definition
\ref{Def_invFSb_convol_fnalFSb} can be written:%
\[
\phi\ast f=\left(  \widehat{\phi}\widehat{f}\right)  ^{\vee}=\left(
2\pi\right)  ^{-\frac{d}{2}}\left[  f_{y},\phi\left(  \cdot-y\right)  \right]
,\quad\phi\in\overset{\vee}{T},\text{ }f\in\left(  \widehat{T}\right)
^{\prime}.
\]

\end{lemma}

\begin{proof}
Replace $S_{w,0}$ by $T$ in the proof of Lemma \ref{Lem_convol_Sw2}.
\end{proof}

This convolution result for $\widetilde{\mathcal{B}}_{2}$\ is the analogue of
part 3 of Theorem \ref{Thm_properties2_B}.

\begin{theorem}
\label{Thm_tildB2_invFSb_convol}Suppose the weight function $w$ has property
W01 w.r.t. the set $\mathcal{A}$ and $w$ satisfies Condition \ref{Cnd_wt_fn_T}%
. Then if $\phi\in\overset{\vee}{T}\subset X_{w}^{0}$ and $v\in\widetilde{X}%
_{1/w}^{0}\subset\left(  \widehat{T}\right)  ^{\prime}$:%
\[
\widetilde{\mathcal{B}}_{2}\left(  \phi,v\right)  =\phi\ast\overset{\vee
}{v_{F}}=\left(  2\pi\right)  ^{-\frac{d}{2}}\left[  \left(  \overset{\vee
}{v_{F}}\right)  _{y},\phi\left(  \cdot-y\right)  \right]  ,
\]

where the convolution is that of Definition \ref{Def_invFSb_convol_fnalFSb}.
\end{theorem}

\begin{proof}
From \ref{a4.3}, $\widetilde{\mathcal{B}}_{2}\left(  \phi,v\right)  =\left(
\widehat{\phi}v_{F}\right)  ^{\vee}=\left(  \widehat{\phi}%
\widehat{\overset{\vee}{v_{F}}}\right)  ^{\vee}$ and by part 5 of Theorem
\ref{Thm_property_tildXo1/w_T}, $\overset{\vee}{v_{F}}\in\left(
\widehat{T}\right)  ^{\prime}$ so that Lemma \ref{Lem_convol_T} implies,
$\widetilde{\mathcal{B}}_{2}\left(  \phi,v\right)  =\phi\ast\overset{\vee
}{v_{F}}=\left(  2\pi\right)  ^{-d/2}\left[  \left(  \overset{\vee}{v_{F}%
}\right)  _{y},\phi\left(  \cdot-y\right)  \right]  $.
\end{proof}

\subsection{The operators $\protect\widetilde{\mathcal{V}}_{2}%
:\protect\widetilde{X}_{1/w}^{0}\rightarrow X_{w}^{0}$,
$\protect\widetilde{\mathcal{W}}_{2}:X_{w}^{0}\rightarrow\protect\widetilde{X}%
_{1/w}^{0}$ and $\protect\widetilde{\Phi}_{2}:\protect\widetilde{X}_{1/w}%
^{0}\rightarrow\left(  X_{w}^{0}\right)  ^{\prime}$}

In an analogous fashion to the operator $\widetilde{\mathcal{V}}_{1}$ of
Definition \ref{Def_op_tildV1} we define the operator $\widetilde{\mathcal{V}%
}_{2}$ in terms of the operators $\mathcal{J}$ and $\widetilde{\mathcal{M}%
}_{2}$ by:

\begin{definition}
\label{Def_op_tildV2}\textbf{The operator} $\widetilde{\mathcal{V}}%
_{2}:\widetilde{X}_{1/w}^{0}\rightarrow X_{w}^{0}$

The operator $\mathcal{J}:L^{2}\rightarrow X_{w}^{0}$ was introduced in
Definition \ref{Def_I_J} and studied in Theorem \ref{Thm_I_J_property}.

Suppose the weight function $w$ has property W02 or W03. These conditions
ensure that both $\widetilde{\mathcal{M}}_{2}$ and $\mathcal{J}$ are defined.

\textbf{Define the operator} $\widetilde{\mathcal{V}}_{2}$ by
$\widetilde{\mathcal{V}}_{2}=\mathcal{J}\widetilde{\mathcal{M}}_{2}$ where
$\widetilde{\mathcal{M}}_{2}:\widetilde{X}_{1/w}^{0}\rightarrow L^{2}$.
\end{definition}

The operator $\widetilde{\mathcal{V}}_{2}$ can be expressed in terms of the
operator $\widetilde{\mathcal{B}}_{2}$ and the basis function as follows:

\begin{theorem}
\label{Thm_tildV2u_=_G*u}Suppose the weight function $w$ has property W02 or
W03 and that $G$ is the corresponding basis function. Then
\begin{equation}
\widetilde{\mathcal{V}}_{2}v=\left(  \frac{v_{F}}{w}\right)  ^{\vee
}=\widetilde{\mathcal{B}}_{2}\left(  G,v\right)  ,\quad v\in\widetilde{X}%
_{1/w}^{0},\label{a6.5}%
\end{equation}

and%
\[
\widetilde{\mathcal{V}}_{2}\phi=G\ast\phi,\quad\phi\in S,
\]

and $\widetilde{\mathcal{V}}_{2}:S\rightarrow G\ast S$ is one-to-one and onto.
\end{theorem}

\begin{proof}
By Definition \ref{Def_I_J}, $\widehat{\mathcal{J}g}=\frac{\widehat{g}}%
{\sqrt{w}}\in L_{loc}^{1}$ for $g\in L^{2}$.

By Definition \ref{Def_tildM2}, $\widehat{\widetilde{\mathcal{M}}_{2}v}%
=\frac{v_{F}}{\sqrt{w}}\in L^{2}$ for $v\in\widetilde{X}_{1/w}^{0}$. Thus
\[
\widehat{\widetilde{\mathcal{V}}_{2}v}=\widehat{\mathcal{J}%
\widetilde{\mathcal{M}}_{2}v}=\frac{\widehat{\widetilde{\mathcal{M}}_{2}v}%
}{\sqrt{w}}=\frac{v_{F}}{w}=\left(  \widetilde{\mathcal{B}}_{2}\left(
G,v\right)  \right)  ^{\wedge},
\]

and so $\widetilde{\mathcal{V}}_{2}v=\widetilde{\mathcal{B}}_{2}\left(
G,v\right)  $. The second equation of this theorem follows from part 3 of
Corollary \ref{Cor_property_tildB2_basis}.

Clearly $\widetilde{\mathcal{V}}_{2}$ is onto. Since $\widetilde{\mathcal{V}%
}_{2}$ is isometric $\widetilde{\mathcal{V}}_{2}\phi=0$ implies $\left\Vert
\widetilde{\mathcal{V}}_{2}\phi\right\Vert _{w,0}=\left\vert \phi\right\vert
_{1/w,0}=\left(  \int\frac{\left\vert _{\phi_{F}}\right\vert ^{2}}{w}\right)
^{1/2}=0$, where $\phi_{F}$ is the restriction of $\widehat{\phi}$ to
$\mathbb{R}^{d}\setminus\mathcal{A}$ where $\mathcal{A}$ is a closed set of
measure zero. But $\widehat{\phi}\in S$ and so $\phi=0$.
\end{proof}

\begin{definition}
\label{Def_op_tildW2}\textbf{The operator }$\widetilde{\mathcal{W}}%
_{2}=\widetilde{\mathcal{L}}_{3}\mathcal{I}$
\end{definition}

The next theorem relates the properties of the operators
$\widetilde{\mathcal{V}}_{2}$ and $\widetilde{\mathcal{W}}_{2}$, and for which
we will require the properties of the operators $\mathcal{I}$ and
$\mathcal{J}$ given in Theorem \ref{Thm_I_J_property}. The next result is the
analogue of Theorem \ref{Thm_tildV1_tildW1_property}.

\begin{theorem}
\label{Thm_tildV2_tildW2_properties}Suppose the weight function $w$ has
property W02 or W03 and satisfies Condition \ref{Cnd_wt_fn_T}. Then in the
seminorm sense the operators $\widetilde{\mathcal{V}}_{2}$ and
$\widetilde{\mathcal{W}}_{2}$ have the following properties:

\begin{enumerate}
\item $\widetilde{\mathcal{V}}_{2}:\widetilde{X}_{1/w}^{0}\rightarrow
X_{w}^{0}$ is a unique, linear isometry.

\item $\widetilde{\mathcal{W}}_{2}:X_{w}^{0}\rightarrow\widetilde{X}_{1/w}%
^{0}$ is a class of linear isometries.

\item $\widetilde{\mathcal{V}}_{2}$ and $\widetilde{\mathcal{W}}_{2}$ are inverses.

\item $\widetilde{\mathcal{V}}_{2}$ and $\widetilde{\mathcal{W}}_{2}$ are onto.

\item $\widetilde{\mathcal{V}}_{2}$ and $\widetilde{\mathcal{W}}_{2}$ are 1-1.

\item $\widetilde{\mathcal{V}}_{2}$ and $\widetilde{\mathcal{W}}_{2}$ are adjoints.

\item If $f\in X_{w}^{0}$ then $\widetilde{\mathcal{W}}_{2}f=\left(
w\widehat{f}\right)  ^{\vee}$.

\item $\widetilde{\mathcal{V}}_{2}$ and $\widetilde{\mathcal{W}}_{2}$ are
isometric isomorphisms, inverses and adjoints.
\end{enumerate}
\end{theorem}

\begin{proof}
The proof is identical to that of Theorem \ref{Thm_tildV1_tildW1_property}
except that $S_{w,0;\mathcal{A}}$ is replaced by $T$.
\end{proof}

Next we prove that $S$ is dense in $\widetilde{X}_{1/w}^{0}$. This result is
the analogue of Corollary \ref{Cor_S_dense_tildXo1/w}.

\begin{corollary}
\label{Cor_S_dense_tildXo1/w_T}Suppose the weight function $w$ has property
W02 or W03 and satisfies Condition \ref{Cnd_wt_fn_T}.

Then $S$ is dense in $\widetilde{X}_{1/w}^{0}$ under $\left\vert
\mathcal{\cdot}\right\vert _{1/w,0}$.
\end{corollary}

\begin{proof}
Proof identical to Corollary \ref{Cor_S_dense_tildXo1/w}.
\end{proof}

The next theorem is an analogue of Corollary
\ref{Cor_invF[wSw2]_dense_tildXo1/w}.

\begin{corollary}
\label{Cor_invF[wT]_dense_tildXo1/w}Suppose the weight function $w$ has
property W02 or W03 w.r.t. the set $\mathcal{A}$, and that $w$ satisfies
Condition \ref{Cnd_wt_fn_T} w.r.t. $\mathcal{A}$ and the subspace $T\subset S
$.

Then the mapping $\left(  w\widehat{\phi}\right)  ^{\vee}$ is an isometry from
$X_{w}^{0}\cap\overset{\vee}{T}\subset X_{w}^{0}$ to $\left(  wT\right)
^{\vee}\subset\widetilde{X}_{1/w}^{0}$ which can be extended uniquely to the
isometric isomorphism $\widetilde{\mathcal{W}}_{2}$.

Finally, $\left(  wT\right)  ^{\vee}$ is dense in $\widetilde{X}_{1/w}^{0}$.
\end{corollary}

\begin{proof}
The isometric property is a consequence of Theorem \ref{Thm_w_phi_in_fnalT}.

Now by part 4 of Theorem \ref{Thm_T_and_sqrt(wt)L2}, $X_{w}^{0}\cap
\overset{\vee}{T}$ is dense in $X_{w}^{0}$ so $\left(  w\widehat{\phi}\right)
^{\vee}$ can be extended uniquely to $X_{w}^{0}$.

From part 7 of Theorem \ref{Thm_tildV2_tildW2_properties} we have
$\widetilde{\mathcal{W}}_{2}f=\left(  w\widehat{f}\right)  ^{\vee}$ when $f\in
X_{w}^{0}$ so the extension of $\left(  w\widehat{\phi}\right)  ^{\vee}$ must
coincide with $\widetilde{\mathcal{W}}_{2}$ which was shown to be an isometric
isomorphism in Theorem \ref{Thm_tildV2_tildW2_properties}.

That $\left(  wT\right)  ^{\vee}$ is dense in $\widetilde{X}_{1/w}^{0}$ is a
straightforward consequence of the fact that $\widetilde{\mathcal{W}}_{2}$
maps $X_{w}^{0}$ onto $\widetilde{X}_{1/w}^{0}$.
\end{proof}

In the following analogue of Theorem \ref{Thm_op_tildPhi1} the operator
$\widetilde{\mathcal{V}}_{2}$ will now be used characterize the bounded linear
functionals on $X_{w}^{0}$, denoted $\left(  X_{w}^{0}\right)  ^{\prime}$, as
the members of $\widetilde{X}_{1/w}^{0}$:

\begin{theorem}
\label{Thm_op_tildPhi2}\textbf{The operator }$\widetilde{\Phi}_{2}$. Suppose
the weight function $w$ has property W02 or W03 and satisfies Condition
\ref{Cnd_wt_fn_T}. Denote by $\left(  \widetilde{X}_{1/w}^{0}\right)
^{\prime}$ the space of bounded linear functionals on $\widetilde{X}_{1/w}%
^{0}$. Then the equation%
\begin{equation}
\left(  \widetilde{\Phi}_{2}v\right)  \left(  u\right)  =\left(
u,\widetilde{\mathcal{V}}_{2}v\right)  _{w,0},\quad u\in X_{w}^{0},\text{
}v\in\widetilde{X}_{1/w}^{0},\label{a6.4}%
\end{equation}

defines a linear operator $\widetilde{\Phi}_{2}:\widetilde{X}_{1/w}%
^{0}\rightarrow\left(  X_{w}^{0}\right)  ^{\prime}$ which is an isometric
isomorphism in the seminorm sense.
\end{theorem}

\begin{proof}
This proof is "identical" to that of Theorem \ref{Thm_op_tildPhi1} and so will
be omitted.
\end{proof}

In an analogous fashion to the negative order Sobolev spaces a bilinear form
can be used to characterize the bounded linear functionals on $X_{w}^{0}$:

\begin{theorem}
\label{Thm_op_tildPhi2_bilinear}If $u\in X_{w}^{0}$ and $v\in\widetilde{X}%
_{1/w}^{0}$ then $\widetilde{\Phi}_{2}$ can be expressed directly in terms of
the bilinear form $\int\widehat{u}\overline{v_{F}}$ as%
\[
\left(  \widetilde{\Phi}_{2}v\right)  \left(  u\right)  =\int\widehat{u}%
\overline{v_{F}},\quad u\in X_{w}^{0},\text{ }v\in\widetilde{X}_{1/w}^{0}.
\]

\end{theorem}

\begin{proof}
A direct consequence of \ref{a6.4} and \ref{a6.5}.
\end{proof}

\section{A further generalization of $\protect\widetilde{X}_{1/w}^{0}%
$\label{Sect_generalize2_tildXo1/w_T}}

\subsection{Introduction}

We now consider the problem of characterizing the bounded linear functionals
on $X_{w}^{0}$ for the weight functions which do not have property
\ref{a1.046}, do not belong to $W_{S;0}$ and do not satisfy Condition
\ref{Cnd_wt_fn_T}.

It was shown in Theorem \ref{Thm_wt_fn_satisfy_property_1} that the Gaussian
and the shifted thin-plate spline weight functions do not have property
\ref{a1.046}, in Subsection \ref{SbSect_wt_fn_not_in_S2} it was shown that
they do not belong to $W_{S;0}$ (\ref{a1.8}) and in Subsection
\ref{SbSect_condit_T} it was shown that they do not satisfy Condition
\ref{Cnd_wt_fn_T}. Thus the results of Sections \ref{Sect_Xo1/w},
\ref{Sect_tildXo1/w_S2dag} and \ref{Sect_generalize1_tildXo1/w_T} are not applicable.

To handle these cases we will adapt definition \ref{a2.5} of $\widetilde{X}%
_{1/w}^{0}$ i.e. that
\[
\widetilde{X}_{1/w}^{0}=\left\{  u\in\left(  \widehat{S}_{w,0}\right)
^{\prime}:u_{F}\in L_{loc}^{1}\left(  \mathbb{R}^{d}\setminus\mathcal{A}%
\right)  \text{ }and\text{ }\int\frac{\left\vert u_{F}\right\vert ^{2}}%
{w}<\infty\right\}  ,
\]

where $u_{F}=\widehat{u}$ restricted to $\mathbb{R}^{d}\setminus\mathcal{A}$.
Here $\left(  \widehat{S}_{w,0}\right)  ^{\prime}$ denotes the continuous
linear functionals on the space $\widehat{S}_{w,0}$ which consists of the
Fourier transforms of functions in $S_{w,0}$. The spaces $S_{w,0}$ and
$\widehat{S}_{w,0}$ were introduced in Definition \ref{Def_Sw2_and_fnal(Sw2)}.
However, \textbf{in this section }$S_{w,0}$\textbf{\ will be considered not as
a subspace of }$S$\textbf{\ but as a topological vector space in itself} and
will be endowed with the topology using the $S$ seminorms \textbf{and} the
norm $\int w\left\vert \phi\right\vert ^{2}$. \textbf{We will still call the
modified data space} $\widetilde{X}_{1/w}^{0}$ \textbf{and define it using
same definition as} \ref{a2.5}. We will prove analogues of the results of
Section \ref{Sect_tildXo1/w_S2dag}.

The space $\widehat{S}_{w,0}$ will be endowed with the topology that makes the
inverse-Fourier transform to $S_{w,0}$ a homeomorphism. Following the approach
used in Section \ref{Sect_tildXo1/w_S2dag} the mappings
$\widetilde{\mathcal{L}}_{4}:L^{2}\rightarrow\left(  \widehat{S}_{w,0}\right)
^{\prime}$ and $\widetilde{\mathcal{M}}_{4}:\widetilde{X}_{1/w}^{0}\rightarrow
L^{2}$ will be introduced, followed by the operators $\widetilde{\mathcal{B}%
}_{3}:X_{w}^{0}\otimes\widetilde{X}_{1/w}^{0}\rightarrow C_{B}^{\left(
0\right)  }$, $\widetilde{\mathcal{V}}_{3}:\widetilde{X}_{1/w}^{0}\rightarrow
X_{w}^{0}$, $\widetilde{\mathcal{W}}_{3}:X_{w}^{0}\rightarrow\widetilde{X}%
_{1/w}^{0}$ and $\widetilde{\Phi}_{3}:\widetilde{X}_{1/w}^{0}\rightarrow
\left(  X_{w}^{0}\right)  ^{\prime}$ the last operator being used to
characterize the bounded linear functionals on $X_{w}^{0}$ as members of
$\widetilde{X}_{1/w}^{0}$ w.r.t. to the bilinear form $\int\widehat{u}%
\overline{v_{F}}$ where $u\in X_{w}^{0}$ and $v\in\widetilde{X}_{1/w}^{0}$.

Although we \textbf{do not} do so here, analogues of the extension and
restriction mappings between $\widetilde{X}_{1/w}^{0}$ and $X_{1/w}^{0}$ of
Subsection \ref{SbSect_E_R_tildXo1/w_Xo1/w} can be easily proved as can
analogues of the Fourier-independent spaces $\widetilde{Y}_{1/w}^{0}$ and
$Y_{w}^{0}$ of Subsection \ref{SbSect_tildYo1/w_Yow_FTindep}.

There is also the possibility of generalizing the space $\widetilde{X}%
_{1/w}^{0}$ of Section \ref{Sect_generalize1_tildXo1/w_T} but we \textbf{do
not} do so here.

\subsection{The space $\protect\widetilde{X}_{1/w}^{0}$}

\begin{definition}
\label{Def_Sw2_lin_fnal}\textbf{The spaces of functionals} $S_{w,0}^{\prime}
$, $\left(  \widehat{S}_{w,0}\right)  ^{\prime}$, $\left(  \overset{\vee
}{S}_{w,0}\right)  ^{\prime}$.

\begin{enumerate}
\item The linear subspace $S_{w,0}\subset S$ was introduced in Definition
\ref{Def_Sw2_and_fnal(Sw2)} where it had the subspace topology induced by $S$.
\textbf{Now we will topologize} $S_{w,0}$\ \textbf{by means of the norm}
$\left(  \int w\left\vert \phi\right\vert ^{2}\right)  ^{1/2}$ \textbf{used to
define} $S_{w,0}$ \textbf{as well as the countable set of seminorms used to
topologize} $S$ (Definition \ref{Def_Distributions})\textbf{.} Thus $f\in
S_{w,0}^{\prime}$ iff $\left[  f,\phi\right]  $ is bounded by a positive,
linear combination of the norm $\left(  \int w\left\vert \phi\right\vert
^{2}\right)  ^{1/2}$ and the seminorms which define the topology on $S$.

\item $\widehat{S}_{w,0}$ is topologized using the seminorms which define the
topology of $S$ and the norm $\left(  \int w\left\vert \overset{\vee}{\phi
}\right\vert ^{2}\right)  ^{1/2}=\left\Vert \overline{\phi}\right\Vert _{w,0}$.

$\overset{\vee}{S}_{w,0}$ is topologized using the seminorms which define the
topology of $S$ and the norm $\left(  \int w\left\vert \widehat{\phi
}\right\vert ^{2}\right)  ^{1/2}=\left\Vert \phi\right\Vert _{w,0}$.

\item With these topologies the Fourier transform is now a homeomorphism from
$S_{w,0}$ to $\overset{\vee}{S}_{w,0}$ and the inverse-Fourier transform is
now a homeomorphism from $S_{w,0}$ to $\widehat{S}_{w,0}$. ?? \textbf{MORE
DETAIL}? See Definitions \ref{Def_Sw2_lin_fnal_0} and \ref{Def_T_lin_fnal_0} ??

\item Now define $S_{w,0}^{\prime}$, $\left(  \widehat{S}_{w,0}\right)
^{\prime}$ and $\left(  \overset{\vee}{S}_{w,0}\right)  ^{\prime}$ to be the
spaces of continuous, linear functionals on $S_{w,0}$, $\widehat{S}_{w,0} $
and $\overset{\vee}{S}_{w,0}$ respectively.
\end{enumerate}
\end{definition}

We now \textbf{define }$\widetilde{X}_{1/w}^{0}$\textbf{\ using the same
definition as} \ref{a2.5}:

\begin{definition}
\label{Def_tildXo1/w_Sw2_2}\textbf{The semi-inner product space}
$\widetilde{X}_{1/w}^{0}$

Suppose the weight function $w$ satisfies property W01 and that $S_{w,0}$ is
endowed with the topology of part 1 of Definition \ref{Def_Sw2_lin_fnal}. Then%
\begin{equation}
\widetilde{X}_{1/w}^{0}=\left\{  u\in\left(  \widehat{S}_{w,0}\right)
^{\prime}:u_{F}\in L_{loc}^{1}\left(  \mathbb{R}^{d}\setminus\mathcal{A}%
\right)  \text{ }and\text{ }\int\frac{\left\vert u_{F}\right\vert ^{2}}%
{w}<\infty\right\}  ,\label{a1.26}%
\end{equation}

where $u_{F}=\widehat{u}$ on $\mathbb{R}^{d}\setminus\mathcal{A}$. Note that
from part 3 of Theorem \ref{Thm_FSw2_and_invFSw2}, $u\in\left(  \widehat{S}%
_{w,0}\right)  ^{\prime}=\left(  S_{w,0}^{\prime}\right)  ^{\vee}$ implies
$\widehat{u}\in S_{w,0}^{\prime}$ and so by part 6 of Definition
\ref{Def_Sw2_lin_fnal}, $\widehat{u}\in L_{loc}^{1}\left(  \mathbb{R}%
^{d}\setminus\mathcal{A}\right)  $ means there exists $f\in L_{loc}^{1}\left(
\mathbb{R}^{d}\setminus\mathcal{A}\right)  $ such that $\left[  \widehat{u}%
,\phi\right]  =\int f\phi$ for all $\phi\in C_{0}^{\infty}\left(
\mathbb{R}^{d}\setminus\mathcal{A}\right)  $.

Lemma \ref{Lem_tildXo1/w_Sw2_2_property} below shows that $\widetilde{X}%
_{1/w}^{0}$ is not empty.

We endow $\widetilde{X}_{1/w}^{0}$ with the semi-inner product and seminorm
\[
\left\langle u,v\right\rangle _{1/w,0}=\int\frac{u_{F}\overline{v_{F}}}%
{w},\qquad\left\vert u\right\vert _{1/w,0}=\left(  \int\frac{\left\vert
u_{F}\right\vert ^{2}}{w}\right)  ^{1/2}.
\]

Regarding part 8 of Definition \ref{Def_Sw2_lin_fnal}, Lemma
\ref{Lem_tildXo1/w_Sw2_2_property} below shows that%
\[
\operatorname*{null}\left\vert \cdot\right\vert _{1/w,0}=\left(
S_{w,0;\mathcal{A}}^{\prime}\right)  ^{\vee},
\]

where we have defined $S_{w,0;\mathcal{A}}^{\prime}$ by%
\begin{equation}
S_{w,0;\mathcal{A}}^{\prime}=\left\{  v\in S_{w,0}^{\prime}%
:\operatorname*{supp}v\subset\mathcal{A}\text{ }as\text{ }a\text{
}distribution\right\}  .\label{a1.042}%
\end{equation}

\end{definition}

\begin{lemma}
\label{Lem_Co,inf(Rd/C)}Let $\mathcal{A}$ be the weight function set used in
the Definition \ref{Def_Xo1/w} of the space $X_{1/w}^{0}$. Suppose
\[
C_{0}^{\left(  0\right)  }\left(  \mathbb{R}^{d}\setminus\mathcal{A}\right)
=\left\{  \psi\in C_{0}^{\left(  0\right)  }:\operatorname*{supp}\psi
\subset\mathbb{R}^{d}\setminus\mathcal{A}\right\}  .
\]

Then $C_{0}^{\left(  0\right)  }\left(  \mathbb{R}^{d}\setminus\mathcal{A}%
\right)  \subset S_{w,0}^{\prime}\cap L_{loc}^{1}\left(  \mathbb{R}%
^{d}\setminus\mathcal{A}\right)  $ and $\int\frac{\left\vert u\right\vert
^{2}}{w}<\infty$ when $u\in C_{0}^{\left(  0\right)  }\left(  \mathbb{R}%
^{d}\setminus\mathcal{A}\right)  $.
\end{lemma}

\begin{proof}
Since $\mathcal{A}$ is a closed set, $C_{0}^{\left(  0\right)  }\left(
\mathbb{R}^{d}\setminus\mathcal{A}\right)  \subset L_{loc}^{1}\left(
\mathbb{R}^{d}\setminus\mathcal{A}\right)  $. If $u\in C_{0}^{\left(
0\right)  }\left(  \mathbb{R}^{d}\setminus\mathcal{A}\right)  $ and $\psi\in
S_{w,0}$ then%
\begin{align*}
\left\vert \int u\psi\right\vert \leq\int\frac{\left\vert u\right\vert }%
{\sqrt{w}}\sqrt{w}\left\vert \psi\right\vert  &  \leq\left(  \int%
\frac{\left\vert u\right\vert ^{2}}{w}\right)  ^{1/2}\left(  \int w\left\vert
\psi\right\vert ^{2}\right)  ^{1/2}\\
&  \leq\frac{1}{\min\limits_{\operatorname*{supp}u}\sqrt{w}}\left(
\int\left\vert u\right\vert ^{2}\right)  ^{1/2}\left(  \int w\left\vert
\psi\right\vert ^{2}\right)  ^{1/2}\\
&  <\infty,
\end{align*}

since $w$ is continuous and positive outside $\mathcal{A}$. Thus $u\in
S_{w,0}^{\prime}$ by virtue of the topology we endowed $S_{w,0}$ with in part
1 of Definition \ref{Def_Sw2_lin_fnal}.

It is also clear that $C_{0}^{\left(  0\right)  }\left(  \mathbb{R}%
^{d}\setminus\mathcal{A}\right)  \subset L_{loc}^{1}\left(  \mathbb{R}%
^{d}\setminus\mathcal{A}\right)  $ in the sense of part 6 of Definition
\ref{Def_Sw2_lin_fnal}.
\end{proof}

\begin{lemma}
\label{Lem_tildXo1/w_Sw2_2_property}Suppose the weight function $w$ has
property W01 for the set $\mathcal{A}$ and that $S_{w,0}$ is endowed with the
topology of part 1 of Definition \ref{Def_Sw2_lin_fnal}. Then:

\begin{enumerate}
\item $\widetilde{X}_{1/w}^{0}$ is not empty. Indeed, $\left(  C_{0}^{\infty
}\left(  \mathbb{R}^{d}\setminus\mathcal{A}\right)  \right)  ^{\vee}%
\subset\widetilde{X}_{1/w}^{0}$.

\item The seminorm $\left\vert \cdot\right\vert _{1/w,0}$ has null space
$\left(  S_{w,0;\mathcal{A}}^{\prime}\right)  ^{\vee}$ (given by
\ref{a1.042}), and $\widetilde{X}_{1/w}^{0}$ is an inner product space iff
$\mathcal{A}$ is empty.

\item If $\tau_{a}$ is the translation operator $u\rightarrow u\left(
\cdot-a\right)  $ then $\tau_{a}:\widetilde{X}_{1/w}^{0}\rightarrow
\widetilde{X}_{1/w}^{0}$ is an isometry and $u\in\widetilde{X}_{1/w}^{0}$
implies $\left(  \tau_{a}u\right)  _{F}=e^{-ia\xi}u_{F}\in L_{loc}^{1}\left(
\mathbb{R}^{d}\setminus\mathcal{A}\right)  $.\medskip

Now suppose $w$ also has property W02 or W03 for parameter $\kappa$.
Then:\medskip

\item If $\left\vert v\left(  \xi\right)  \right\vert \leq c\left(
1+\left\vert \xi\right\vert \right)  ^{\kappa}$ for some constant $c$, then
$\overset{\vee}{v}\in\widetilde{X}_{1/w}^{0}\cap S^{\prime}$.

\item $S\subset\widetilde{X}_{1/w}^{0}$.

\item If $u\in\widetilde{X}_{1/w}^{0}$ then $\frac{\xi^{\alpha}u_{F}}{w}\in
L^{1}$ when $\left\vert \alpha\right\vert \leq\kappa$.
\end{enumerate}
\end{lemma}

\begin{proof}
\textbf{Part 1} Suppose $u\in\left(  C_{0}^{\infty}\left(  \mathbb{R}%
^{d}\setminus\mathcal{A}\right)  \right)  ^{\vee}$. Then by Lemma
\ref{Lem_Co,inf(Rd/C)}, $\widehat{u}\in S_{w,0}^{\prime}$ so that $u\in\left(
\widehat{S}_{w,0}\right)  ^{\prime}$ and $u_{F}=\widehat{u}|_{\mathbb{R}%
^{d}\setminus\mathcal{A}}\in L_{loc}^{1}\left(  \mathbb{R}^{d}\setminus
\mathcal{A}\right)  $. Finally, $\int\frac{\left\vert u_{F}\right\vert ^{2}%
}{w}<\infty$ by Lemma \ref{Lem_Co,inf(Rd/C)}.\medskip

\textbf{Part 2} Suppose $\left\vert v\right\vert _{1/w,0}^{2}=\int%
\frac{\left\vert v_{F}\right\vert ^{2}}{w}=0$. Then $v_{F}=0$ a.e. on
$\mathbb{R}^{d}\setminus\mathcal{A}$ and so $\widehat{v}\in S_{w,0}^{\prime}$
and $\operatorname*{supp}\widehat{v}\subset\mathcal{A}$ i.e. $\widehat{v}%
\in\widetilde{S}_{w,0;\mathcal{A}}^{\prime}$. The argument is easily
reversible.\medskip

\textbf{Part 3} Suppose $u\in\widetilde{X}_{1/w}^{0}$ and $\phi\in
C_{0}^{\infty}\left(  \mathbb{R}^{d}\setminus\mathcal{A}\right)  $. Then
$u\in\left(  \widehat{S}_{w,0}\right)  ^{\prime}=\left(  S_{w,0}^{\prime
}\right)  ^{\vee}$ and $u_{F}\in L_{loc}^{1}\left(  \mathbb{R}^{d}%
\setminus\mathcal{A}\right)  $ means there exists $f\in L_{loc}^{1}\left(
\mathbb{R}^{d}\setminus\mathcal{A}\right)  $ such that $\left[  u_{F}%
,\psi\right]  =\int f\psi$ for all $\psi\in C_{0}^{\infty}\left(
\mathbb{R}^{d}\setminus\mathcal{A}\right)  $.

Thus
\begin{align*}
\left[  \left(  \tau_{a}u\right)  _{F},\phi\right]  =\left[  \widehat{\tau
_{a}u},\phi\right]   &  =\left[  \tau_{a}u,\widehat{\phi}\right]  =\left[
u,\tau_{-a}\widehat{\phi}\right]  =\left[  u,\widehat{e^{ia\xi}\phi}\right] \\
&  =\left[  \widehat{u},e^{ia\xi}\phi\right]  =\int e^{ia\xi}f\phi=\left[
e^{ia\xi}u_{F},\phi\right]  ,
\end{align*}

so that $\left(  \tau_{a}u\right)  _{F}=e^{ia\xi}u_{F}$ on $\mathbb{R}%
^{d}\setminus\mathcal{A}$. Thus%
\[
\left\vert \tau_{a}u\right\vert _{1/w,0}=\int\frac{\left\vert \left(  \tau
_{a}u\right)  _{F}\right\vert ^{2}}{w}=\int\frac{\left\vert u_{F}\right\vert
^{2}}{w}=\left\vert u\right\vert _{1/w,0}.
\]
\medskip

\textbf{Part 4} Now $\widetilde{X}_{1/w}^{0}=\left\{  u\in\left(
\widehat{S}_{w,0}\right)  ^{\prime}:u_{F}\in L_{loc}^{1}\left(  \mathbb{R}%
^{d}\setminus\mathcal{A}\right)  \text{ }and\text{ }\int\frac{\left\vert
u_{F}\right\vert ^{2}}{w}<\infty\right\}  $.

Set $u=\overset{\vee}{v}$. Then $u\in S^{\prime}$ and $\widehat{u}=v\in
L_{loc}^{1}$ so that $u_{F}\in L_{loc}^{1}\left(  \mathbb{R}^{d}%
\setminus\mathcal{A}\right)  $. If $\phi\in\widehat{S}_{w,0}$ then%
\[
\left[  u,\phi\right]  =\left[  \widehat{u},\overset{\vee}{\phi}\right]  =\int
v\overset{\vee}{\phi}=\int\frac{v}{\sqrt{w}}\sqrt{w}\overset{\vee}{\phi},
\]

so that%
\[
\left\vert \left[  u,\phi\right]  \right\vert \leq\left(  \int\frac{\left\vert
v\right\vert ^{2}}{w}\right)  ^{1/2}\left(  \int w\left\vert \overset{\vee
}{\phi}\right\vert ^{2}\right)  ^{1/2}\leq c\left(  \int\frac{\left(
1+\left\vert \cdot\right\vert \right)  ^{2\kappa}}{w}\right)  ^{1/2}\left(
\int w\left\vert \overset{\vee}{\phi}\right\vert ^{2}\right)  ^{1/2}.
\]

Thus $\left\vert \left[  \widehat{u},\overset{\vee}{\phi}\right]  \right\vert
\leq c^{\prime}\left(  \int w\left\vert \overset{\vee}{\phi}\right\vert
^{2}\right)  ^{1/2}$ when $\overset{\vee}{\phi}\in S_{w,0}$ and part 5 of
Definition \ref{Def_Sw2_lin_fnal} means that $\widehat{u}\in\left(
S_{w,0}\right)  ^{\prime}$ i.e. $u\in\left(  \widehat{S}_{w,0}\right)
^{\prime}$.\medskip

\textbf{Part 5} Follows directly from part 4.\medskip

\textbf{Part 6} If $u\in\widetilde{X}_{1/w}^{0}$ then from part 1, $\int%
\frac{\left\vert u_{F}\right\vert ^{2}}{w}<\infty$ and thus%
\[
\int\frac{\xi^{\alpha}u_{F}}{w}=\int\frac{\xi^{\alpha}}{\sqrt{w}}\frac{u_{F}%
}{\sqrt{w}}\leq\left(  \int\frac{\xi^{2\alpha}}{w}\right)  ^{1/2}\left(
\int\frac{\left\vert u_{F}\right\vert ^{2}}{w}\right)  ^{1/2}<\infty.
\]

\end{proof}

\subsection{The operators $\protect\widetilde{\mathcal{L}}_{4}:L^{2}%
\rightarrow\left(  \protect\widehat{S}_{w,0}\right)  ^{\prime}$ and
$\protect\widetilde{\mathcal{M}}_{4}:\protect\widetilde{X}_{1/w}%
^{0}\rightarrow L^{2}$}

The operators $\mathcal{L}_{1}$ (Definition \ref{Def_op_L1} assuming property
\ref{a1.046}) and $\mathcal{L}_{2}$ (Definition \ref{Def_op_L2} assuming $w\in
S_{2}^{\dagger}$) were constructed so that they mapped $L^{2}\rightarrow
S^{\prime}$ and then were shown to map $L^{2}\rightarrow X_{1/w}^{0}$. So
noting the definition of $\widetilde{X}_{1/w}^{0}$ (equation \ref{a1.26}) we
will define the operator $\widetilde{\mathcal{L}}_{4}:L^{2}\rightarrow\left(
\widehat{S}_{w,0}\right)  ^{\prime}$ and then show that
$\widetilde{\mathcal{L}}_{4}:L^{2}\rightarrow\widetilde{X}_{1/w}^{0}$ is an isomorphism.

\begin{definition}
\label{Def_op_tildL4}\textbf{The linear operator} $\widetilde{\mathcal{L}}%
_{4}:L^{2}\rightarrow\left(  \widehat{S}_{w,0}\right)  ^{\prime}$\textbf{.}

Suppose that $S_{w,0}$ is endowed with the topology of part of Definition
\ref{Def_Sw2_lin_fnal}. To be consistent with previous definitions of the
mappings $\mathcal{L}_{1}$ (Definition \ref{Def_op_L1}) and $\mathcal{L}_{2}$
(Definition \ref{Def_op_L2}) for $g\in L^{2}$ we try $\widetilde{\mathcal{L}%
}_{4}g=\left(  \sqrt{w}\widehat{g}\right)  ^{\vee}$ and we want
$\widetilde{\mathcal{L}}_{4}g\in\left(  \widehat{S}_{w,0}\right)  ^{\prime}$.
So suppose $\phi\in S_{w,0}$. Then%
\[
\left\vert \int\sqrt{w}\widehat{g}\phi\right\vert \leq\int\left\vert
\widehat{g}\right\vert \left(  \sqrt{w}\left\vert \phi\right\vert \right)
\leq\left\Vert g\right\Vert _{2}\left(  \int w\left\vert \phi\right\vert
^{2}\right)  ^{1/2},
\]

and the topology of $S_{w,0}$ ensures that $\sqrt{w}\widehat{g}\in
S_{w,0}^{\prime}$ and hence that $\left(  \sqrt{w}\widehat{g}\right)  ^{\vee
}\in\left(  \widehat{S}_{w,0}\right)  ^{\prime}$.

So we define $\widetilde{\mathcal{L}}_{4}$ by%
\begin{equation}
\widetilde{\mathcal{L}}_{4}g=\left(  \sqrt{w}\widehat{g}\right)  ^{\vee},\quad
g\in L^{2}.\label{a1.7}%
\end{equation}

\end{definition}

We now prove some properties of $\widetilde{\mathcal{L}}_{4}$ which relate to
the space $\widetilde{X}_{1/w}^{0}$.

\begin{theorem}
\label{Thm_tildL4_property}\textbf{Properties of }$\widetilde{\mathcal{L}}%
_{4}$

\begin{enumerate}
\item If $g\in L^{2}$ then $\left(  \widetilde{\mathcal{L}}_{4}g\right)
_{F}=\sqrt{w}\widehat{g}$ on $\mathbb{R}^{d}\setminus\mathcal{A}$.

\item $\widetilde{\mathcal{L}}_{4}:L^{2}\rightarrow\widetilde{X}_{1/w}^{0}$ is
a linear isometry.

\item $\tau_{a}\widetilde{\mathcal{L}}_{4}=\widetilde{\mathcal{L}}_{4}\tau
_{a}$, where $\tau_{a}$ is the translation operator $\tau_{a}u=u\left(
\cdot-a\right)  $.
\end{enumerate}
\end{theorem}

\begin{proof}
\textbf{Part 1} From the definition of $\widetilde{\mathcal{L}}_{4}$, $\left(
\widetilde{\mathcal{L}}_{4}g\right)  ^{\wedge}=\sqrt{w}\widehat{g}\in
S_{w,0}^{\prime}$. Next observe that because $w\in C^{\left(  0\right)
}\left(  \mathbb{R}^{d}\setminus\mathcal{A}\right)  $ we have $C_{0}^{\infty
}\left(  \mathbb{R}^{d}\setminus\mathcal{A}\right)  \subset S_{w,0}$ and hence
that
\[
\left[  \left(  \widetilde{\mathcal{L}}_{4}g\right)  ^{\wedge},\phi\right]
=\left[  \sqrt{w}\widehat{g},\phi\right]  ,\quad\phi\in C_{0}^{\infty}\left(
\mathbb{R}^{d}\setminus\mathcal{A}\right)  .
\]

If we can show $\sqrt{w}g\in L_{loc}^{1}\left(  \mathbb{R}^{d}\setminus
\mathcal{A}\right)  $ it follows that $\left(  \widetilde{\mathcal{L}}%
_{4}g\right)  ^{\wedge}\in L_{loc}^{1}\left(  \mathbb{R}^{d}\setminus
\mathcal{A}\right)  $. Indeed, if $K\subset\mathbb{R}^{d}\setminus\mathcal{A}$
is compact, $w\in C^{\left(  0\right)  }\left(  \mathbb{R}^{d}\setminus
\mathcal{A}\right)  $ implies
\[
\int_{K}\sqrt{w}\left\vert g\right\vert \leq\left(  \int_{K}w\right)
^{1/2}\left\Vert g\right\Vert _{2}\leq\max_{K}\left(  w\right)  \left(
\int_{K}1\right)  \left\Vert g\right\Vert _{2}.
\]

Thus $\left(  \widetilde{\mathcal{L}}_{4}g\right)  ^{\wedge}\in L_{loc}%
^{1}\left(  \mathbb{R}^{d}\setminus\mathcal{A}\right)  $ and $\left(
\widetilde{\mathcal{L}}_{4}g\right)  _{F}=\sqrt{w}g$ a.e.\medskip

\textbf{Part 2} $\left\vert \widetilde{\mathcal{L}}_{4}g\right\vert
_{1/w,0}=\left(  \int\frac{\left\vert \left(  \widetilde{\mathcal{L}}%
_{4}g\right)  _{F}\right\vert ^{2}}{w}\right)  ^{1/2}=\left\Vert g\right\Vert
_{2}$.\medskip

\textbf{Part 3} $\tau_{a}\widetilde{\mathcal{L}}_{4}u=\tau_{a}\left(  \sqrt
{w}\widehat{g}\right)  ^{\vee}=\left(  e^{-ia\xi}\sqrt{w}\widehat{g}\right)
^{\vee}=\left(  \sqrt{w}\widehat{\tau_{a}g}\right)  ^{\vee}%
=\widetilde{\mathcal{L}}_{4}\tau_{a}u$.
\end{proof}

Our definition of the operator $\widetilde{\mathcal{M}}_{4}:\widetilde{X}%
_{1/w}^{0}\rightarrow L^{2}$ will formally be the same as that of
$\widetilde{\mathcal{M}}_{1}$ in Definition \ref{Def_tildM1}.

\begin{definition}
\label{Def_op_tildM4}\textbf{The operator} $\widetilde{\mathcal{M}}%
_{4}:\widetilde{X}_{1/w}^{0}\rightarrow L^{2}$

Suppose the weight function $w$ has property W01 and that $S_{w,0}$ is endowed
with the topology of part of Definition \ref{Def_Sw2_lin_fnal}.

From the definition of $\widetilde{X}_{1/w}^{0}$, $u\in\widetilde{X}_{1/w}%
^{0}$, $u_{F}\in L_{loc}^{1}\left(  \mathbb{R}^{d}\setminus\mathcal{A}\right)
$ and $\frac{u_{F}}{\sqrt{w}}\in L^{2}$. We can now define the linear mapping
$\widetilde{\mathcal{M}}_{4}:\widetilde{X}_{1/w}^{0}\rightarrow L^{2}$ by%
\[
\widetilde{\mathcal{M}}_{4}u=\left(  \frac{u_{F}}{\sqrt{w}}\right)  ^{\vee
},\text{\quad}u\in\widetilde{X}_{1/w}^{0}.
\]

\end{definition}

$\widetilde{\mathcal{M}}_{4}$ has the following properties:

\begin{theorem}
\label{Thm_tildM4_property}\textbf{Properties of }$\widetilde{\mathcal{M}}%
_{4}$. The operator $\widetilde{\mathcal{M}}_{4}:\widetilde{X}_{1/w}%
^{0}\rightarrow L^{2}$ is linear and isometric with null space $\left(
S_{w,0;\mathcal{A}}^{\prime}\right)  ^{\vee}$.

Further, the translation operator $\tau_{a}$ commutes with
$\widetilde{\mathcal{M}}_{4}$ i.e. $\tau_{a}\widetilde{\mathcal{M}}%
_{4}=\widetilde{\mathcal{M}}_{4}\tau_{a}$ on $\widetilde{X}_{1/w}^{0}$.
\end{theorem}

\begin{proof}
That $\widetilde{\mathcal{M}}_{4}$ is an isometry is clear from the definition
of $\widetilde{X}_{1/w}^{0}$. Since $\widetilde{\mathcal{M}}_{4}$ is an
isometry the null space of $\widetilde{\mathcal{M}}_{4}$ is the null space of
the seminorm $\left\vert \cdot\right\vert _{1/w,0}$, namely $\left(
S_{w,0;\mathcal{A}}^{\prime}\right)  ^{\vee}$.

Now suppose $u\in\widetilde{X}_{1/w}^{0}$. Then from part 3 of Lemma
\ref{Lem_tildXo1/w_Sw2_2_property}%
\[
\tau_{a}\widetilde{\mathcal{M}}_{4}u=\tau_{a}\left(  \frac{u_{F}}{\sqrt{w}%
}\right)  ^{\vee}=\left(  e^{-ia\xi}\frac{u_{F}}{\sqrt{w}}\right)  ^{\vee
}=\left(  \frac{\left(  \tau_{a}u\right)  _{F}}{\sqrt{w}}\right)  ^{\vee
}=\widetilde{\mathcal{M}}_{4}\tau_{a}u.
\]

\end{proof}

The next theorem indicates how the operators $\widetilde{\mathcal{L}}_{4}$ and
$\widetilde{\mathcal{M}}_{4}$ interact.

\begin{theorem}
\label{Thm_tildL4_tildM4}\textbf{Interaction of }$\widetilde{\mathcal{L}}_{4}
$\textbf{\ and }$\widetilde{\mathcal{M}}_{4}$\textbf{.} Suppose the weight
function $w$ has the properties assumed in Definition \ref{Def_op_L2}. Then
the operators $\widetilde{\mathcal{L}}_{4}:L^{2}\rightarrow\widetilde{X}%
_{1/w}^{0}$ and $\widetilde{\mathcal{M}}_{4}:\widetilde{X}_{1/w}%
^{0}\rightarrow L^{2}$ satisfy:

\begin{enumerate}
\item $\widetilde{\mathcal{M}}_{4}\widetilde{\mathcal{L}}_{4}=I$ on $L^{2}$
and $\widetilde{\mathcal{L}}_{4}\widetilde{\mathcal{M}}_{4}u-u\in\left(
S_{w,0;\mathcal{A}}^{\prime}\right)  ^{\vee}$ when $u\in\widetilde{X}%
_{1/w}^{0}$.\medskip

Thus $\widetilde{\mathcal{M}}_{4}$ and $\widetilde{\mathcal{L}}_{4}$ are
inverses in the seminorm sense.\medskip

\item $\widetilde{\mathcal{L}}_{4}$ is 1-1. $\widetilde{\mathcal{L}}_{4}$ is
also onto in the seminorm sense.

\item $\widetilde{\mathcal{M}}_{4}$ is onto. It is also 1-1 in the seminorm sense.

\item $\widetilde{\mathcal{M}}_{4}$ and $\widetilde{\mathcal{L}}_{4}$ are
adjoints in the sense that $\left\langle \widetilde{\mathcal{L}}%
_{4}g,f\right\rangle _{1/w,0}=\left(  g,\widetilde{\mathcal{M}}_{4}f\right)
_{2}$.
\end{enumerate}
\end{theorem}

\begin{proof}
\textbf{Property 1}. From the definition of $\widetilde{\mathcal{L}}_{4}$,
$\left(  \widetilde{\mathcal{L}}_{4}g\right)  _{F}=\sqrt{w}\widehat{g}$. Thus,
from the definition of $\widetilde{\mathcal{M}}_{4}$ and Theorem
\ref{Thm_tildL4_property}, for $g\in L^{2}$
\[
\left(  \widetilde{\mathcal{M}}_{4}\widetilde{\mathcal{L}}_{4}g\right)
^{\wedge}=\frac{\left(  \widetilde{\mathcal{L}}_{4}g\right)  _{F}}{\sqrt{w}%
}=\frac{\sqrt{w}\widehat{g}}{\sqrt{w}}=\widehat{g},
\]

and thus $\widetilde{\mathcal{M}}_{4}\widetilde{\mathcal{L}}_{4}=I$.

Next we show that $\widetilde{\mathcal{L}}_{4}\widetilde{\mathcal{M}}%
_{4}u-u\in\left(  S_{w,0;\mathcal{A}}^{\prime}\right)  ^{\vee}$ on
$\widetilde{X}_{1/w}^{0}$. In fact by Theorem \ref{Thm_tildL4_property},
$u\in\widetilde{X}_{1/w}^{0}$ implies $\left(  \widetilde{\mathcal{L}}%
_{4}\widetilde{\mathcal{M}}_{4}u\right)  _{F}=\sqrt{w}\left(
\widetilde{\mathcal{M}}_{4}u\right)  ^{\wedge}=u_{F}$ so that
\[
\left\vert \widetilde{\mathcal{L}}_{4}\widetilde{\mathcal{M}}_{4}%
u-u\right\vert _{1/w,0}^{2}=\int w\left\vert \left(  \widetilde{\mathcal{L}%
}_{4}\widetilde{\mathcal{M}}_{4}u-u\right)  _{F}\right\vert ^{2}=0,
\]

so that $\widetilde{\mathcal{L}}_{4}\widetilde{\mathcal{M}}_{4}u-u\in\left(
S_{w,0;\mathcal{A}}^{\prime}\right)  ^{\vee}$ i.e. $\widetilde{\mathcal{L}%
}_{4}\widetilde{\mathcal{M}}_{4}u-u\in\left(  S_{w,0;\mathcal{A}}^{\prime
}\right)  ^{\vee}$, where $S_{w,0;\mathcal{A}}^{\prime}$ was defined by
\ref{a2.8}. Since $\left(  S_{w,0;\mathcal{A}}^{\prime}\right)  ^{\vee}$ is
the null space of the seminorm on $\widetilde{X}_{1/w}^{0}$,
$\widetilde{\mathcal{M}}_{4}$ and $\widetilde{\mathcal{L}}_{4}$ are inverses
in the seminorm sense.\medskip

\textbf{Property 2}. That $\widetilde{\mathcal{L}}_{4}$ is 1-1 follows from
$\widetilde{\mathcal{M}}_{4}\widetilde{\mathcal{L}}_{4}=I$. That
$\widetilde{\mathcal{L}}_{4}$ is onto follows from $\widetilde{\mathcal{L}%
}_{4}\widetilde{\mathcal{M}}_{4}u-u\in\left(  S_{w,0;\mathcal{A}}^{\prime
}\right)  ^{\vee}$ when $u\in\widetilde{X}_{1/w}^{0}$.\medskip

\textbf{Property 3}. That $\widetilde{\mathcal{M}}_{4}$ is onto follows from
$\widetilde{\mathcal{M}}_{4}\widetilde{\mathcal{L}}_{4}=I$. That
$\widetilde{\mathcal{L}}_{4}$ is onto follows from $\widetilde{\mathcal{L}%
}_{4}\widetilde{\mathcal{M}}_{4}u-u\in\left(  S_{w,0;\mathcal{A}}^{\prime
}\right)  ^{\vee}$ when $u\in\widetilde{X}_{1/w}^{0}$.\medskip

\textbf{Property 4}. If $g\in L^{2}$ and $f\in\widetilde{X}_{1/w}^{0}$,
\[
\left\langle \widetilde{\mathcal{L}}_{4}g,f\right\rangle _{1/w,0}=\int%
\frac{\left(  \widetilde{\mathcal{L}}_{4}g\right)  _{F}\overline{f_{F}}}%
{w}=\int\frac{\sqrt{w}\widehat{g}\overline{f_{F}}}{w}=\int\widehat{g}%
\frac{\overline{f_{F}}}{\sqrt{w}}=\left(  \widehat{g},\overline{\left(
\widetilde{\mathcal{M}}_{4}f\right)  ^{\wedge}}\right)  _{2}=\left(
g,\overline{\widetilde{\mathcal{M}}_{4}f}\right)  _{2}.
\]

\end{proof}

Since $L^{2}$ is complete, the mappings of the previous theorem will yield the
following important result.

\begin{corollary}
\label{Cor_Xo1/w_Sw2_2_semiHilb}In general $\widetilde{X}_{1/w}^{0}$ is a
semi-Hilbert space and $\widetilde{X}_{1/w}^{0}$ is a Hilbert space iff
$\mathcal{A}$ can be chosen to be empty.
\end{corollary}

\begin{proof}
If $\left\{  u_{k}\right\}  $ is a Cauchy sequence in $\widetilde{X}_{1/w}%
^{0}$ then $\left\{  \widetilde{\mathcal{M}}_{4}u_{k}\right\}  $ is Cauchy in
$L^{2}$ since $\widetilde{\mathcal{M}}_{4}$ is isometric. Thus
$\widetilde{\mathcal{M}}_{4}u_{k}\rightarrow u$ for some $u\in L^{2}$ and
since $\widetilde{\mathcal{M}}_{4}\widetilde{\mathcal{L}}_{4}=I$ by Theorem
\ref{Thm_tildL4_tildM4}%
\[
\left\vert u_{k}-\widetilde{\mathcal{L}}_{4}u\right\vert _{1/w,0}=\left\Vert
\widetilde{\mathcal{M}}_{4}\left(  u_{k}-\widetilde{\mathcal{L}}_{4}u\right)
\right\Vert _{2}=\left\Vert \widetilde{\mathcal{M}}_{4}u_{k}-u\right\Vert
_{2}\rightarrow0.
\]

Thus $\widetilde{X}_{1/w}^{0}$ is complete.

Finally, from part 2 of Lemma \ref{Lem_tildXo1/w_Sw2_2_property},
$\widetilde{X}_{1/w}^{0}$is an inner product space iff $\mathcal{A}$ can be
chosen to be empty.
\end{proof}

\subsection{The operator $\protect\widetilde{\mathcal{B}}_{3}:X_{w}^{0}%
\otimes\protect\widetilde{X}_{1/w}^{0}\rightarrow C_{B}^{\left(  0\right)  }$}

We will now construct the analogues of the operators $\widetilde{\mathcal{V}%
}_{3}$, $\widetilde{\mathcal{B}}_{3}$ and $\widetilde{\Phi}_{3}$ of Subsection
\ref{SbSect_ops_B_and_V} and these will allow us to characterize the bounded
linear functionals on $\widetilde{X}_{1/w}^{0}$ using the bilinear form
$\int\widehat{u}\overline{v_{F}}$ where $u\in X_{w}^{0}$ and $v\in
\widetilde{X}_{1/w}^{0}$.

\begin{definition}
\label{Def_map_B general}\textbf{The mapping }$\widetilde{\mathcal{B}}%
_{3}:X_{w}^{0}\otimes\widetilde{X}_{1/w}^{0}\rightarrow C_{B}^{\left(
0\right)  }$

If $\left(  u,v\right)  \in X_{w}^{0}\otimes\widetilde{X}_{1/w}^{0}$ define
the bilinear mapping $\widetilde{\mathcal{B}}_{3}:X_{w}^{0}\otimes
\widetilde{X}_{1/w}^{0}\rightarrow C_{B}^{\left(  0\right)  }$ by:%
\begin{equation}
\widetilde{\mathcal{B}}_{3}\left(  u,v\right)  =\left(  \widehat{u}%
v_{F}\right)  ^{\vee},\quad u\in X_{w}^{0},\text{ }v\in\widetilde{X}_{1/w}%
^{0}.\label{a2.6}%
\end{equation}

Noting that the product of two measurable functions is a measurable function
this definition makes sense since from the definitions of $\widetilde{X}%
_{1/w}^{0}$ and $X_{w}^{0}$, $v_{F}\in L_{loc}^{1}\left(  \mathbb{R}%
^{d}\setminus\mathcal{A}\right)  $, $\frac{v_{F}}{\sqrt{w}}\in L^{2}$,
$\widehat{u}\in L_{loc}^{1}$ and $\sqrt{w}\widehat{u}\in L^{2}$ so that,
$\int\left\vert \widehat{u}v_{F}\right\vert =\int\left\vert \sqrt
{w}\widehat{u}\frac{v_{F}}{\sqrt{w}}\right\vert \leq\left\Vert u\right\Vert
_{w,0}\left\vert v\right\vert _{1/w,0}<\infty$. Thus $\widehat{u}v_{F}\in
L^{1}$ and hence $\left(  \widehat{u}v_{F}\right)  ^{\vee}\in C_{B}^{\left(
0\right)  }$.
\end{definition}

\begin{theorem}
\label{Thm_op_tildB3_property}\textbf{Properties of the operator
}$\widetilde{\mathcal{B}}_{3}$

\begin{enumerate}
\item $\widetilde{\mathcal{B}}_{3}:X_{w}^{0}\otimes\widetilde{X}_{1/w}%
^{0}\rightarrow C_{B}^{\left(  0\right)  }$ is a continuous bilinear mapping
when $C_{B}^{\left(  0\right)  }$ is endowed with the supremum norm
$\left\Vert \cdot\right\Vert _{\infty}$. In fact%
\[
\left\Vert \widetilde{\mathcal{B}}_{3}\left(  u,v\right)  \right\Vert
_{\infty}\leq\left(  2\pi\right)  ^{-\frac{d}{2}}\left\Vert u\right\Vert
_{w,0}\left\vert v\right\vert _{1/w,0},\quad u\in X_{w}^{0},\text{ }%
v\in\widetilde{X}_{1/w}^{0}.
\]

\item In the sense of distributions, if $\alpha+\beta=\gamma\geq0$ then
\[
D^{\gamma}\widetilde{\mathcal{B}}_{3}\left(  u,v\right)  =\left(
\widehat{D^{\alpha}u}\left(  D^{\beta}v\right)  _{F}\right)  ^{\vee},\quad
u\in X_{w}^{0},\text{ }v\in\widetilde{X}_{1/w}^{0}.
\]

\item Suppose $\tau_{a}$ denotes the (distribution) translation operator
$\tau_{a}f=f\left(  \cdot-a\right)  $, $a\in\mathbb{R}^{d}$. Then
$\widetilde{\mathcal{B}}_{3}$ commutes with $\tau_{a}$ in the sense that%
\[
\tau_{a}\widetilde{\mathcal{B}}_{3}\left(  u,v\right)  =\widetilde{\mathcal{B}%
}_{3}\left(  \tau_{a}u,v\right)  =\widetilde{\mathcal{B}}_{3}\left(
u,\tau_{a}v\right)  ,\quad u\in X_{w}^{0},\text{ }v\in\widetilde{X}_{1/w}^{0}.
\]

\item We have%
\[
\widetilde{\mathcal{B}}_{3}\left(  u,v\right)  =\left(  \mathcal{I}u\right)
\ast\widetilde{\mathcal{M}}_{4}v,\quad u\in X_{w}^{0},\text{ }v\in
\widetilde{X}_{1/w}^{0}.
\]

\item $\widetilde{\mathcal{B}}_{3}$ is a convolution in the sense that:
$X_{w}^{0}\subset S^{\prime}$, $S\subset\widetilde{X}_{1/w}^{0}$ and%
\[
\widetilde{\mathcal{B}}_{3}\left(  u,\phi\right)  =u\ast\phi,\quad u\in
X_{w}^{0},\text{ }\phi\in S.
\]

Further, $\widetilde{\mathcal{B}}_{3}\left(  u,\phi\right)  \in C_{B}^{\infty
}$ and $D^{\gamma}\widetilde{\mathcal{B}}_{3}\left(  u,\phi\right)
=D^{\alpha}u\ast D^{\beta}\phi$ for all $\gamma=\alpha+\beta$.
\end{enumerate}
\end{theorem}

\begin{proof}
\textbf{Part 1} Noting the calculations done in the definition of
$\widetilde{\mathcal{B}}_{3}$ all that remains to be shown is the continuity
of $\widetilde{\mathcal{B}}_{3}$. But from \ref{a1.20}
\begin{align*}
\left\vert \mathcal{B}\left(  u,v\right)  \right\vert \leq\left\vert \left(
\widehat{u}v_{F}\right)  ^{\vee}\right\vert \leq\left(  2\pi\right)
^{-\frac{d}{2}}\left\vert \int e^{-ix\xi}\widehat{u}\left(  \xi\right)
v_{F}\left(  \xi\right)  d\xi\right\vert  &  \leq\left(  2\pi\right)
^{-\frac{d}{2}}\int\left\vert \widehat{u}\right\vert \left\vert v_{F}%
\right\vert \\
&  \leq\left(  2\pi\right)  ^{-\frac{d}{2}}\left\Vert u\right\Vert
_{w,0}\left\vert v\right\vert _{1/w,0}.
\end{align*}
\medskip

\textbf{Part 2}%
\begin{align*}
D^{\gamma}\widetilde{\mathcal{B}}_{3}\left(  u,v\right)  =D^{\gamma}\left(
\widehat{u}v_{F}\right)  ^{\vee}=\left(  \left(  i\xi\right)  ^{\gamma}\left(
\widehat{u}v_{F}\right)  \right)  ^{\vee} &  =\left(  \left(  i\xi\right)
^{\alpha}\widehat{u}\text{ }\left(  i\xi\right)  ^{\beta}v_{F}\right)  ^{\vee
}\\
&  =\left(  \widehat{D^{\alpha}u}\text{ }\left(  D^{\beta}v\right)
_{F}\right)  ^{\vee}.
\end{align*}
\medskip

\textbf{Part 3}%
\[
\tau_{a}\widetilde{\mathcal{B}}_{3}\left(  u,v\right)  =\tau_{a}\left(
\widehat{u}v_{F}\right)  ^{\vee}=\left(  e^{-ia\xi}\widehat{u}v_{F}\right)
^{\vee}=\left(  \widehat{\tau_{a}u}v_{F}\right)  ^{\vee}%
=\widetilde{\mathcal{B}}_{3}\left(  \tau_{a}u,v\right)  ,
\]

and%
\[
\tau_{a}\widetilde{\mathcal{B}}_{3}\left(  u,v\right)  =\tau_{a}\left(
\widehat{u}v_{F}\right)  ^{\vee}=\left(  \widehat{u}e^{-ia\xi}v_{F}\right)
^{\vee}=\left(  \widehat{u}\text{ }\left(  \tau_{a}v\right)  _{F}\right)
^{\vee}.
\]

From the definition of $\widetilde{X}_{1/w}^{0}$ there exists $f\in
L_{loc}^{1}\left(  \mathbb{R}^{d}\setminus\mathcal{A}\right)  $ such that
$\left[  \widehat{v},\phi\right]  =\int f\phi$ for all $\phi\in S_{w,0}$. Thus
for $\phi\in S_{w,0}$,
\[
\left[  \widehat{\tau_{a}v},\phi\right]  =\left[  e^{-ia\xi}\widehat{v}%
,\phi\right]  =\left[  \widehat{v},e^{-ia\xi}\phi\right]  =\int e^{-ia\xi
}f\phi,
\]

and consequently $\left(  \tau_{a}v\right)  _{F}\in L_{loc}^{1}\left(
\mathbb{R}^{d}\setminus\mathcal{A}\right)  $.\medskip

\textbf{Part 4 }From the definition of $\widetilde{\mathcal{B}}_{3}$ and the
definition of $\widetilde{\mathcal{M}}_{4}$ (Definition \ref{Def_op_tildM4})%
\begin{align*}
\widetilde{\mathcal{B}}_{3}\left(  u,v\right)  =\left(  \widehat{u}%
v_{F}\right)  ^{\vee}=\left(  \sqrt{w}\widehat{u}\frac{v_{F}}{\sqrt{w}%
}\right)  ^{\vee} &  =\left(  \sqrt{w}\widehat{u}\right)  ^{\vee}\ast\left(
\frac{v_{F}}{\sqrt{w}}\right)  ^{\vee}\\
&  =\left(  \mathcal{I}u\right)  \ast\widetilde{\mathcal{M}}_{4}v.
\end{align*}

\textbf{Part 5} Part 5 of Lemma \ref{Lem_tildXo1/w_Sw2_2_property} showed that
$S\subset\widetilde{X}_{1/w}^{0}$. So by definition \ref{a2.6} of
$\widetilde{\mathcal{B}}_{3}\left(  u,\phi\right)  $ and the definition
\ref{1.20} of the convolution of a member of $S^{\prime}$ and a member of $S$,
$\widetilde{\mathcal{B}}_{3}\left(  u,\phi\right)  =\left(  \widehat{u}%
\phi_{F}\right)  ^{\vee}=\left(  \widehat{u}\widehat{\phi}\right)  ^{\vee
}=u\ast\phi\in C_{B}^{\left(  0\right)  }$ so that by \ref{2.40}, $D^{\gamma
}\widetilde{\mathcal{B}}_{3}\left(  u,\phi\right)  =\left(  D^{\alpha
}u\right)  \ast D^{\beta}\phi$ for all $\gamma$ and $\alpha+\beta=\gamma$.
Thus $\widetilde{\mathcal{B}}_{3}\left(  u,\phi\right)  \in C_{B}^{\infty}$.
\end{proof}

\begin{remark}
\label{Rem_Thm_op_tildB3_property}\ 

\begin{enumerate}
\item Part 3 also follows directly from part 4 and the translation operator
results of Theorem \ref{Thm_tildM4_property} and Theorem
\ref{Thm_tildL4_property}.

\item Part 5 shows how $\widetilde{\mathcal{B}}_{3}$ can be regarded as a
convolution by restricting the domain of the second variable to $S$.
\end{enumerate}
\end{remark}

\begin{corollary}
\label{Cor_property_tildB3_basis}\textbf{Properties of }%
$\widetilde{\mathcal{B}}_{3}\left(  G\mathbf{,\cdot}\right)  $\textbf{:}
Suppose $w$ has property W02 or W03 for $\kappa$ and $G$ is the basis
function. Then if $v,v^{\prime}\in\widetilde{X}_{1/w}^{0}$:

\begin{enumerate}
\item $D^{\gamma}G\in X_{w}^{0}$ when $\left\vert \gamma\right\vert \leq
\kappa$ and%
\[
D^{\gamma}\widetilde{\mathcal{B}}_{3}\left(  G,v\right)
=\widetilde{\mathcal{B}}_{3}\left(  D^{\gamma}G,v\right)  ,\quad\left\vert
\gamma\right\vert \leq\kappa.
\]

\item $\widetilde{\mathcal{B}}_{3}\left(  G,v\right)  \in X_{w}^{0}$ and
$\left\Vert \widetilde{\mathcal{B}}_{3}\left(  G,v\right)  \right\Vert
_{w,0}=\left\vert v\right\vert _{1/w,0}$.

\item $\widetilde{\mathcal{B}}_{3}\left(  G,\phi\right)  =G\ast\phi$ when
$\phi\in S$.

\item $\int F\widetilde{\mathcal{B}}_{3}\left(  \widetilde{\mathcal{B}}%
_{3}\left(  G,v\right)  ,\overline{v^{\prime}}\right)  =\left\langle
v,v_{\ast}^{\prime}\right\rangle _{1/w,0}$, where $F$ is the Fourier transform
and $v_{\ast}^{\prime}\left(  x\right)  =v^{\prime}\left(  -x\right)  $.
\end{enumerate}
\end{corollary}

\begin{proof}
\textbf{Part 1} By Theorem \ref{Thm_basis_fn_properties_all_m_W2}, $D^{\gamma
}G\in X_{w}^{0}$ when $\left\vert \gamma\right\vert \leq\kappa$ and so by part
2 of Theorem \ref{Thm_op_tildB3_property}, $D^{\gamma}\widetilde{\mathcal{B}%
}_{3}\left(  G,v\right)  \in\left(  \widehat{D^{\gamma}G}v_{F}\right)  ^{\vee
}=\widetilde{\mathcal{B}}_{3}\left(  D^{\gamma}G,v\right)  $.\medskip

\textbf{Part 2} $G\in X_{w}^{0}$ and $\widetilde{\mathcal{B}}_{3}\left(
G,v\right)  \in S^{\prime}$. Further, $F\widetilde{\mathcal{B}}_{3}\left(
G,v\right)  =\widehat{G}v_{F}=\frac{v_{F}}{w}\in L^{1}$ and $\int w\left\vert
F\widetilde{\mathcal{B}}_{3}\left(  G,v\right)  \right\vert ^{2}=\int%
\frac{\left\vert v_{F}\right\vert ^{2}}{w}=\left\vert v\right\vert _{1/w,0}$.
Thus $\widetilde{\mathcal{B}}_{3}\left(  G,v\right)  \in X_{w}^{0}$.\medskip

\textbf{Part 3} From part 5 Lemma \ref{Lem_tildXo1/w_Sw2_2_property},
$S\subset\widetilde{X}_{1/w}^{0}$ and $\widetilde{\mathcal{B}}_{3}\left(
G,\phi\right)  =\left(  \widehat{G}\widehat{\phi}\right)  ^{\vee}=G\ast\phi$
by the convolution definition \ref{1.20}.\medskip

\textbf{Part 4} $\widetilde{\mathcal{B}}_{3}\left(  G,v\right)  =\left(
\frac{v_{F}}{w}\right)  ^{\vee}$ so
\begin{align*}
\int F\widetilde{\mathcal{B}}_{3}\left(  \widetilde{\mathcal{B}}_{3}\left(
G,v\right)  ,\overline{v^{\prime}}\right)  =\int F\widetilde{\mathcal{B}}%
_{3}\left(  \left(  \frac{v_{F}}{w}\right)  ^{\vee},\overline{v^{\prime}%
}\right)  =\int F\left(  \left(  \frac{v_{F}}{w},\left(  \overline{v^{\prime}%
}\right)  _{F}\right)  ^{\vee}\right)   &  =\int\frac{1}{w}v_{F}\left(
\overline{v^{\prime}}\right)  _{F}\\
&  =\int\frac{1}{w}v_{F}\overline{\left(  v_{\ast}^{\prime}\right)  _{F}}\\
&  =\left\langle v,v_{\ast}^{\prime}\right\rangle _{1/w,0}.
\end{align*}

\end{proof}

We will now define the equivalent of the operator $\widetilde{\mathcal{V}}%
_{1}$ of Definition \ref{Def_op_tildV1}.

\subsection{The operators $\protect\widetilde{\mathcal{V}}_{3}%
:\protect\widetilde{X}_{1/w}^{0}\rightarrow X_{w}^{0}$,
$\protect\widetilde{\mathcal{W}}_{3}:X_{w}^{0}\rightarrow\protect\widetilde{X}%
_{1/w}^{0}$ and $\protect\widetilde{\Phi}_{3}:\protect\widetilde{X}_{1/w}%
^{0}\rightarrow\left(  X_{w}^{0}\right)  ^{\prime}$}

\begin{definition}
\label{Def_op_tildV3}\textbf{The operator} $\widetilde{\mathcal{V}}%
_{3}:\widetilde{X}_{1/w}^{0}\rightarrow X_{w}^{0}$

Suppose the weight function $w$ has property W02 or W03. These conditions
ensure that both $\widetilde{\mathcal{M}}_{4}$ (Definition \ref{Def_op_tildM4}%
) and $\mathcal{J}$ (Definition \ref{Def_I_J}) are defined.

The operator $\widetilde{\mathcal{V}}_{3}$ is now defined by
$\widetilde{\mathcal{V}}_{3}=\mathcal{J}\widetilde{\mathcal{M}}_{4}$.
\end{definition}

The operator $\widetilde{\mathcal{V}}_{3}$ can be expressed in terms of the
operator $\widetilde{\mathcal{B}}_{3}$ and the basis function as follows:

\begin{theorem}
\label{Thm_tildV3g_=_G*g}Suppose the weight function $w$ has property W02 or
W03 and that $G$ is the corresponding basis function. Then
\begin{equation}
\widetilde{\mathcal{V}}_{3}v=\left(  \frac{v_{F}}{w}\right)  ^{\vee
}=\widetilde{\mathcal{B}}_{3}\left(  G,v\right)  ,\quad v\in\widetilde{X}%
_{1/w}^{0},\label{a6.7}%
\end{equation}

and%
\[
\widetilde{\mathcal{V}}_{3}v=G\ast v,\quad v\in S,
\]

and $\widetilde{\mathcal{V}}_{3}:S\rightarrow G\ast S$ is 1-1 and onto.
\end{theorem}

\begin{proof}
By Definition \ref{Def_I_J}, $\widehat{\mathcal{J}g}=\frac{\widehat{g}}%
{\sqrt{w}}\in L_{loc}^{1}$ for $g\in L^{2}$.

By Definition \ref{Def_op_tildM4}, $\widehat{\widetilde{\mathcal{M}}_{4}%
v}=\frac{v_{F}}{\sqrt{w}}\in L^{2}$ for $v\in\widetilde{X}_{1/w}^{0}$. Thus
\[
\widehat{\widetilde{\mathcal{V}}_{3}v}=\widehat{\mathcal{J}%
\widetilde{\mathcal{M}}_{4}v}=\frac{\widehat{\widetilde{\mathcal{M}}_{4}v}%
}{\sqrt{w}}=\frac{v_{F}}{w}=\left(  \widetilde{\mathcal{B}}_{3}\left(
G,v\right)  \right)  ^{\wedge},
\]

and so $\widetilde{\mathcal{V}}_{3}v=\widetilde{\mathcal{B}}_{3}\left(
G,v\right)  $. The second equation of this theorem follows from part 3 of
Corollary \ref{Cor_property_tildB3_basis}.

Clearly $\widetilde{\mathcal{V}}_{3}$ is onto. Since $\widetilde{\mathcal{V}%
}_{3}$ is isometric $\widetilde{\mathcal{V}}_{3}\phi=0$ implies $\left\Vert
\widetilde{\mathcal{V}}_{3}\phi\right\Vert _{w,0}=\left\vert \phi\right\vert
_{1/w,0}=\left(  \int\frac{\left\vert _{\phi_{F}}\right\vert ^{2}}{w}\right)
^{1/2}=0$, where $\phi_{F}$ is the restriction of $\widehat{\phi}$ to
$\mathbb{R}^{d}\setminus\mathcal{A}$ where $\mathcal{A}$ is a closed set of
measure zero. But $\widehat{\phi}\in S$ and so $\phi=0$.
\end{proof}

\begin{definition}
\label{Def_op_tildW3}\textbf{The operator }$\widetilde{\mathcal{W}}%
_{3}=\widetilde{\mathcal{L}}_{4}\mathcal{I}$
\end{definition}

The operator $\widetilde{\mathcal{V}}_{3}$ has the following properties:

\begin{theorem}
\label{Thm_tildV3_tildL4}Suppose the weight function $w$ has property W02 or
W03. Then in the seminorm sense:

\begin{enumerate}
\item $\widetilde{\mathcal{V}}_{3}:\widetilde{X}_{1/w}^{0}\rightarrow
X_{w}^{0}$ is a linear isometry.

\item $\widetilde{\mathcal{W}}_{3}:X_{w}^{0}\rightarrow\widetilde{X}_{1/w}%
^{0}$ is an isometry.

\item $\widetilde{\mathcal{V}}_{3}$ and $\widetilde{\mathcal{W}}_{3}$ are inverses.

\item $\widetilde{\mathcal{V}}_{3}$ and $\widetilde{\mathcal{W}}_{3}$ are onto.

\item $\widetilde{\mathcal{V}}_{3}$ and $\widetilde{\mathcal{W}}_{3}$ are 1-1.

\item $\widetilde{\mathcal{V}}_{3}$ and $\widetilde{\mathcal{W}}_{3}$ are adjoints.

\item If $f\in X_{w}^{0}$ then $\widetilde{\mathcal{W}}_{3}f=\left(
w\widehat{f}\right)  ^{\vee}$.

\item $\widetilde{\mathcal{V}}_{3}$ and $\widetilde{\mathcal{W}}_{3}$ are
isometric isomorphisms, inverses and adjoints.
\end{enumerate}
\end{theorem}

\begin{proof}
\textbf{Parts 1 and 2}. True since $\mathcal{I}$, $\mathcal{J}$,
$\widetilde{\mathcal{L}}_{4}$ and $\widetilde{\mathcal{M}}_{4}$ are
isometric.\medskip

\textbf{Part 3}. From Theorem \ref{Thm_tildL4_tildM4} $\widetilde{\mathcal{M}%
}_{4}\widetilde{\mathcal{L}}_{4}=I$\ and Theorem \ref{Thm_I_J_property} showed
that $\mathcal{JI}f=f$. Hence $\widetilde{\mathcal{V}}_{3}%
\widetilde{\mathcal{W}}_{3}=\widetilde{\mathcal{V}}_{3}\widetilde{\mathcal{L}%
}_{4}\mathcal{I}=\mathcal{JM}_{3}\widetilde{\mathcal{L}}_{4}\mathcal{I}%
=\mathcal{JI}$ and so $\widetilde{\mathcal{V}}_{3}\widetilde{\mathcal{W}}%
_{3}=I$.

From Theorem \ref{Thm_I_J_property} $\mathcal{IJ}=I$ and from Theorem
\ref{Thm_tildL4_tildM4} $\widetilde{\mathcal{L}}_{4}\widetilde{\mathcal{M}%
}_{4}g-g\in\left(  \widetilde{S}_{\mathcal{A}}^{\prime}\right)  ^{\vee}$.
Hence $\widetilde{\mathcal{W}}_{3}\widetilde{\mathcal{V}}_{3}%
=\widetilde{\mathcal{L}}_{4}\mathcal{IJM}_{3}=\widetilde{\mathcal{L}}%
_{4}\widetilde{\mathcal{M}}_{4}$ and $\widetilde{\mathcal{W}}_{3}%
\widetilde{\mathcal{V}}_{3}g-g\in\left(  \widetilde{S}_{\mathcal{A}}^{\prime
}\right)  ^{\vee}$.\medskip

\textbf{Parts 4 and 5}. follow directly from part 3.\medskip

\textbf{Part 6}. Theorem \ref{Thm_tildL4_tildM4}\ and Theorem
\ref{Thm_I_J_property} imply $\widetilde{\mathcal{L}}_{4}$,
$\widetilde{\mathcal{M}}_{4}$ are adjoints and $\mathcal{I}$, $\mathcal{J}$
are adjoints. Thus
\[
\left\langle g,\widetilde{\mathcal{W}}_{3}f\right\rangle _{1/w,0}=\left\langle
g,\widetilde{\mathcal{L}}_{4}\mathcal{I}f\right\rangle _{1/w,0}=\left(
\widetilde{\mathcal{M}}_{4}g,\mathcal{I}f\right)  _{2}=\left(  \mathcal{JM}%
_{3}g,f\right)  _{w,0}=\left(  \widetilde{\mathcal{V}}_{3}g,f\right)  _{w,0}.
\]
\medskip

\textbf{Part 7}. If $f\in X_{w}^{0}$ then
\[
\widehat{\widetilde{\mathcal{W}}_{3}f}=\widehat{\widetilde{\mathcal{L}}%
_{4}\mathcal{I}f}=\sqrt{w}\widehat{\mathcal{I}f}=\sqrt{w}\left(  \sqrt
{w}\widehat{f}\right)  =w\widehat{f}.
\]

\end{proof}

\begin{theorem}
\label{Thm_tildV3*(DRx)=delta}Suppose the weight function $w$ has property W02
or W03 for parameter $\kappa$. Then
\[
\left(  2\pi\right)  ^{d/2}\widetilde{\mathcal{V}}_{3}^{\ast}D^{\gamma}%
R_{x}=D^{\gamma}\delta\left(  \cdot-x\right)  ,\quad\left\vert \gamma
\right\vert \leq\kappa,
\]

where $R_{x}$ is the Riesz representer of the evaluation functional
$f\rightarrow f\left(  x\right)  $. Further%
\[
\left\langle u,D^{\gamma}\delta\left(  \cdot-x\right)  \right\rangle
_{1/w,0}=\left(  -1\right)  ^{\left\vert \gamma\right\vert }\left(  D^{\gamma
}\widetilde{\mathcal{V}}_{3}^{\ast}u\right)  \left(  x\right)  ,\quad u\in
X_{1/w}^{0},\text{ }\left\vert \gamma\right\vert \leq\kappa.
\]

\end{theorem}

\begin{proof}
If $f\in\widetilde{X}_{1/w}^{0}$ then by parts 6 and 7 of Theorem
\ref{Thm_tildV3_tildL4}
\[
\widetilde{\mathcal{V}}_{3}^{\ast}f=\widetilde{\mathcal{W}}_{3}f=\left(
w\widehat{f}\right)  ^{\vee}.
\]

By using the Fourier transform properties given in the Appendix
\ref{Ch_Appendx_basic_notation} we obtain
\begin{align*}
\widetilde{\mathcal{V}}_{3}^{\ast}D^{\gamma}R_{x}=\left(  w\widehat{D^{\gamma
}R_{x}}\right)  ^{\vee}=\left(  w\left(  \xi\right)  \left(  i\xi\right)
^{\gamma}\widehat{R_{x}}\left(  \xi\right)  \right)  ^{\vee} &  =\left(
2\pi\right)  ^{-d/2}\left(  w\left(  \xi\right)  \left(  i\xi\right)
^{\gamma}\widehat{G\left(  \cdot-x\right)  }\left(  \xi\right)  \right)
^{\vee}\\
&  =\left(  2\pi\right)  ^{-d/2}\left(  \left(  i\xi\right)  ^{\gamma
}e^{-ix\xi}\right)  ^{\vee}\\
&  =\left(  2\pi\right)  ^{-d/2}D^{\gamma}\left(  \left(  e^{-ix\xi}\right)
^{\vee}\right)  .
\end{align*}

But $\left(  e^{-ia\xi}\right)  ^{\vee}=\left(  2\pi\right)  ^{d/2}%
\delta\left(  \cdot-a\right)  $ so that $\widetilde{\mathcal{V}}_{3}^{\ast
}D^{\gamma}R_{x}=D^{\gamma}\delta\left(  \cdot-x\right)  $ and thus%
\begin{align*}
\left\langle u,D^{\gamma}\delta\left(  \cdot-x\right)  \right\rangle
_{1/w,0}=\left\langle u,\widetilde{\mathcal{V}}_{3}^{\ast}D^{\gamma}%
R_{x}\right\rangle _{1/w,0}=\left(  \widetilde{\mathcal{V}}_{3}^{\ast
}u,D^{\gamma}R_{x}\right)  _{w,0} &  =\left(  \left(  -D\right)  ^{\gamma
}\widetilde{\mathcal{V}}_{3}^{\ast}u,R_{x}\right)  _{w,0}\\
&  =\left(  -1\right)  ^{\left\vert \gamma\right\vert }\left(  D^{\gamma
}\widetilde{\mathcal{V}}_{3}^{\ast}u\right)  \left(  x\right)  .
\end{align*}

\end{proof}

We now prove a density result for $\widetilde{X}_{1/w}^{0}$.

\begin{corollary}
\label{Cor_S_dense_Xo1/w_Sw2_2}The space $S$ is dense in $\widetilde{X}%
_{1/w}^{0}$.
\end{corollary}

\begin{proof}
Choose $u\in\widetilde{X}_{1/w}^{0}$ so that $\widetilde{\mathcal{V}}_{3}u\in
X_{w}^{0}$. From part 3 Theorem \ref{Thm_Jg_properties}, $G\ast S$ is dense in
$X_{w}^{0}$ so given $\varepsilon>0$ there exists $\phi_{\varepsilon}\in G\ast
S$ such that $\left\Vert \widetilde{\mathcal{V}}_{3}u-\phi_{\varepsilon
}\right\Vert _{w,0}<\varepsilon$. But $\widetilde{\mathcal{V}}_{3}%
:S\rightarrow G\ast S$ is an isomorphism so $\widetilde{\mathcal{V}}_{3}%
^{-1}\phi_{\varepsilon}\in S$. Since $\widetilde{\mathcal{V}}_{3}:X_{1/w}%
^{0}\rightarrow X_{w}^{0}$ is an isometry%
\[
\left\vert u-\widetilde{\mathcal{V}}_{3}^{-1}\phi_{\varepsilon}\right\vert
_{1/w,0}=\left\Vert \widetilde{\mathcal{V}}_{3}\left(
u-\widetilde{\mathcal{V}}_{3}^{-1}\phi_{\varepsilon}\right)  \right\Vert
_{w,0}=\left\Vert \widetilde{\mathcal{V}}_{3}u-\phi_{\varepsilon}\right\Vert
_{w,0}<\varepsilon,
\]

which proves this density result.

Since $\widetilde{X}_{1/w}^{0}$ is complete and $S$ is dense in $\widetilde{X}%
_{1/w}^{0}$.
\end{proof}

\begin{corollary}
\label{Cor_invF[wS^]_dense_tildXo1/w}Suppose the weight function $w$ has
property W02 or W03 w.r.t. the set $\mathcal{A}$, and that $S_{w,0}$ is
endowed with the countable seminorm topology of part 1 of Definition
\ref{Def_Sw2_lin_fnal}. Then:

\begin{enumerate}
\item The mapping $\left(  w\widehat{\phi}\right)  ^{\vee}$ is an isometry
from $\overset{\vee}{S}_{w,0}=X_{w}^{0}\cap S$ onto $\left(  wS_{w,0}\right)
^{\vee}\subset\widetilde{X}_{1/w}^{0}$ and is the \textbf{restriction of the
isometric isomorphism} $\widetilde{\mathcal{W}}_{3}$.

\item The space $\left(  wS_{w,0}\right)  ^{\vee}$ is dense in $\widetilde{X}%
_{1/w}^{0}$.

\item $\psi\in S\cap\left(  wS_{w,0}\right)  ^{\vee}\Longleftrightarrow\psi\in
S$ and $G\ast\psi\in\overset{\vee}{S}_{w,0}$.

\item $w\phi$ is a continuous map from $S_{w,0}$ to $S_{w,0}^{\prime}$.

\item $\left(  w\widehat{\phi}\right)  ^{\vee}$ is a continuous map from
$\overset{\vee}{S}_{w,0}$ to $\left(  \widehat{S}_{w,0}\right)  ^{\prime}$.
\end{enumerate}
\end{corollary}

\begin{proof}
\textbf{Part 1} From Theorem \ref{Thm_tildV3_tildL4} $\left(  w\widehat{\phi
}\right)  ^{\vee}$ is an isometry from $\overset{\vee}{S}_{w,0}=X_{w}^{0}\cap
S$ \textbf{into} $\left(  wS_{w,0}\right)  ^{\vee}\subset\widetilde{X}%
_{1/w}^{0}$ which is the restriction of the isometric isomorphism
$\widetilde{\mathcal{W}}_{3}:X_{w}^{0}\rightarrow\widetilde{X}_{1/w}^{0}$.
Clearly $\left(  w\widehat{\phi}\right)  ^{\vee}:\overset{\vee}{S}%
_{w,0}\rightarrow\left(  wS_{w,0}\right)  ^{\vee}$ is onto.\medskip

\textbf{Part 2} Since $\overset{\vee}{S}_{w,0}$ is dense in $X_{w}^{0}$ the
ontoness of $\left(  w\widehat{\phi}\right)  ^{\vee}$ implies the density of
$\left(  wS_{w,0}\right)  ^{\vee}$ in $\widetilde{X}_{1/w}^{0}$.\medskip

\textbf{Part 3} If
\begin{align*}
\psi\in S\cap\left(  wS_{w,0}\right)  ^{\vee}  & \Longleftrightarrow\psi\in
S\text{ }and\text{ }\widehat{\psi}\in wS_{w,0}\\
& \Longleftrightarrow\psi\in S\text{ }and\text{ }\frac{1}{w}\widehat{\psi
}=\widehat{G}\widehat{\psi}=\left(  G\ast\psi\right)  ^{\wedge}\in S_{w,0}\\
& \Longleftrightarrow\psi\in S\text{ }and\text{ }G\ast\psi\in\overset{\vee
}{S}_{w,0}.
\end{align*}

\textbf{Part 4} If $\phi,\psi\in S_{w,0}$ then%
\[
\left\vert \left[  w\phi,\psi\right]  \right\vert =\left\vert \int w\phi
\psi\right\vert \leq\int w\left\vert \phi\right\vert \left\vert \psi
\right\vert =\int\sqrt{w}\left\vert \phi\right\vert \sqrt{w}\left\vert
\psi\right\vert \leq\left(  \int w\left\vert \phi\right\vert ^{2}\right)
^{1/2}\left(  \int w\left\vert \psi\right\vert ^{2}\right)  ^{1/2},
\]

and since $S_{w,0}$ is endowed with the countable seminorm topology of part 1
of Definition \ref{Def_Sw2_lin_fnal}, it follows that $w\phi\in S_{w,0}%
^{\prime}$ and $w\phi$ is a continuous map from $S_{w,0}$ to $S_{w,0}^{\prime
}$.\medskip

\textbf{Part 5} From part 4, $w\widehat{\phi}$ is a continuous map from
$\overset{\vee}{S}_{w,0}$ to $S_{w,0}^{\prime}$, and so by Definition
\ref{Def_Sw2_lin_fnal}, $\left(  w\widehat{\phi}\right)  ^{\vee}$ is a
continuous map from $\overset{\vee}{S}_{w,0}$ to $\left(  S_{w,0}^{\prime
}\right)  ^{\vee}=\left(  \widehat{S}_{w,0}\right)  ^{\prime}$.
\end{proof}

\begin{theorem}
Suppose the weight function $w$ has the properties used to define
$\widetilde{X}_{1/w}^{0}$.

Denote by $\left(  \widetilde{X}_{1/w}^{0}\right)  ^{\prime}$ the space of
bounded linear functionals on $\widetilde{X}_{1/w}^{0}$. Then the equation%
\begin{equation}
\left(  \widetilde{\Phi}_{3}u\right)  g=\left(  g,\widetilde{\mathcal{V}}%
_{3}u\right)  _{w,0},\quad g\in X_{w}^{0},\text{ }u\in\widetilde{X}_{1/w}%
^{0},\label{a6.6}%
\end{equation}

defines a linear operator $\widetilde{\Phi}_{3}:\widetilde{X}_{1/w}%
^{0}\rightarrow\left(  X_{w}^{0}\right)  ^{\prime}$ which is an isometric
isomorphism in the seminorm sense.
\end{theorem}

\begin{proof}
Since $\widetilde{\mathcal{V}}_{3}:\widetilde{X}_{1/w}^{0}\rightarrow
X_{w}^{0}$ is an isometry, given $u\in\widetilde{X}_{1/w}^{0}$ and $g\in
X_{w}^{0}$ we have
\[
\left\vert \left(  g,\widetilde{\mathcal{V}}_{3}u\right)  _{w,0}\right\vert
\leq\left\Vert g\right\Vert _{w,0}\left\Vert \widetilde{\mathcal{V}}%
_{3}u\right\Vert _{w,0}=\left\Vert g\right\Vert _{w,0}\left\vert u\right\vert
_{1/w,0},
\]

and for each $u\in\widetilde{X}_{1/w}^{0}$, the expression $\left(
g,\widetilde{\mathcal{V}}_{3}u\right)  _{w,0}$, $g\in X_{w}^{0}$ defines a
bounded linear functional on $X_{w}^{0}$. Denote this functional by
$\widetilde{\Phi}_{3}u$ so that \ref{a6.6} holds and the operator norm is
\[
\left\Vert \widetilde{\Phi}_{3}u\right\Vert _{op}=\sup\limits_{g\in X_{w}^{0}%
}\frac{\left\vert \left(  g,\widetilde{\mathcal{V}}_{2}u\right)
_{w,0}\right\vert }{\left\Vert g\right\Vert _{w,0}}\leq\left\vert u\right\vert
_{1/w,0}.
\]

Thus $\widetilde{\Phi}_{3}:\widetilde{X}_{1/w}^{0}\rightarrow\left(  X_{w}%
^{0}\right)  ^{\prime}$. In fact, we can easily prove that $\left\Vert
\widetilde{\Phi}_{3}u\right\Vert _{op}=\left\vert u\right\vert _{1/w,0}$ by
noting that when $g=\widetilde{\mathcal{V}}_{3}u$, $\frac{\left\vert \left(
g,\widetilde{\mathcal{V}}_{3}u\right)  _{w,0}\right\vert }{\left\Vert
g\right\Vert _{w,0}}=\left\Vert \widetilde{\mathcal{V}}_{3}u\right\Vert
_{w,0}=\left\vert u\right\vert _{1/w,0}$. Clearly $\widetilde{\Phi}_{3}u=0$
implies $\left\vert u\right\vert _{1/w,0}=0$ so $\widetilde{\Phi}_{3}$ is 1-1.

To prove $\widetilde{\Phi}_{3}$ is onto choose $\mathcal{Y}\in\left(
X_{w}^{0}\right)  ^{\prime}$. Since $X_{w}^{0}$ is an inner product space
there exists $v\in X_{w}^{0}$ such that $\mathcal{Y}g=\left(  g,v\right)
_{w,0}$ when $g\in X_{w}^{0}$. But by Theorem \ref{Thm_tildV3_tildL4}
$\mathcal{Y}g=\left(  g,v\right)  _{w,0}=\left(  g,\widetilde{\mathcal{V}}%
_{3}\widetilde{\mathcal{V}}_{3}^{\ast}v\right)  _{w,0}$ and comparison with
\ref{a6.6} yields $\left(  \widetilde{\Phi}_{3}\widetilde{\mathcal{V}}%
_{3}^{\ast}v\right)  g=\mathcal{Y}g$ and $\widetilde{\Phi}_{3}%
\widetilde{\mathcal{V}}_{3}^{\ast}v=\mathcal{Y}$.
\end{proof}

In an analogous fashion to the negative order Sobolev spaces a bilinear form
can be used to characterize the bounded linear functionals on $X_{w}^{0}$:

\begin{theorem}
\label{Thm_op_tildPhi3}If $u\in X_{w}^{0}$ and $v\in\widetilde{X}_{1/w}^{0}$
then $\widetilde{\Phi}_{3}$ can be expressed directly in terms of the bilinear
form $\int\widehat{u}\overline{v_{F}}$ as%
\begin{equation}
\left(  \widetilde{\Phi}_{3}v\right)  \left(  u\right)  =\int\widehat{u}%
\overline{v_{F}},\quad u\in X_{w}^{0},\text{ }v\in\widetilde{X}_{1/w}%
^{0}.\label{a6.8}%
\end{equation}

\end{theorem}

\begin{proof}
A direct consequence of \ref{a6.6} and \ref{a6.7}.
\end{proof}

\subsection{Examples}

\begin{example}
\textbf{The Gaussian} Here $w\left(  x\right)  =e^{\left\vert x\right\vert
^{2}}$ so $\mathcal{A}$ is empty and $\widetilde{X}_{1/w}^{0}$ is a Hilbert
space. Also, $\widetilde{\Phi}_{3}:\widetilde{X}_{1/w}^{0}\rightarrow\left(
X_{w}^{0}\right)  ^{\prime}$ is an isometric isomorphism in the norm sense.
\end{example}

\begin{example}
\textbf{Shifted thin-plate splines\ }Here $w=\frac{1}{\widetilde{e}\left(
v\right)  }\frac{\left\vert \cdot\right\vert ^{2v+d}}{\widetilde{K}%
_{v+d/2}\left(  \left\vert \cdot\right\vert \right)  }$ Set $s=2v+d$. Then
$s>0$ and by Theorem \ref{Thm_bnds_modif_MacDonald},%
\[
\frac{1}{\widetilde{e}\left(  v\right)  c_{s}^{\prime}}\left\vert
\xi\right\vert ^{s}e^{\left\vert \xi\right\vert }\leq w\left(  \xi\right)
\leq\frac{1}{\widetilde{e}\left(  v\right)  c_{s}}\left\vert \xi\right\vert
^{s}e^{\left\vert \xi\right\vert },\quad\xi\in\mathbb{R}^{d}.
\]

so $\mathcal{A}$ is empty and $\widetilde{X}_{1/w}^{0}$ is a Hilbert space.
Also, $\widetilde{\Phi}_{3}:\widetilde{X}_{1/w}^{0}\rightarrow\left(
X_{w}^{0}\right)  ^{\prime}$ is an isometric isomorphism in the norm sense.
\end{example}

\chapter{An upper bound for the derivative of the 1-dim. hat basis function
smoother smoother\label{Ch_bnd_deriv_hat_smth_large_supp}}

\section{Introduction\label{Sect_chapter_intro}}

In this chapter we derive an upper bound for the derivative of the
1-dimensional (scaled) hat basis function smoother under the assumption that
the data function has a bounded derivative and sufficiently large support
w.r.t. the data region. This will be used above in Example
\ref{Ex_ex_smth_hat_larg_supp_1dim} to prove order 1 convergence of the Exact smoother.

Suppose the data region is $\Omega=\left(  \underline{a},\overline{a}\right)
$ and that there are $N$ \textbf{dependent data} points $X=\left\{  x^{\left(
k\right)  }\right\}  _{k=1}^{N}$ satisfying $\underline{a}<x^{\left(
n\right)  }<x^{\left(  n+1\right)  }<\overline{a}$. Further suppose we are
using a \textbf{scaled hat basis function} $\Lambda_{\lambda}\left(  x\right)
=\Lambda\left(  x/\lambda\right)  $ which has been scaled so that
$\operatorname*{diam}\Omega\leq\frac{1}{2}\operatorname*{diam}%
\operatorname*{supp}\Lambda_{\lambda}$ i.e. $\operatorname*{diam}\Omega
\leq\lambda$. In this case we say that the hat basis function has
\textbf{large support w.r.t. the data region} and \textbf{this chapter will be
based on the fact that }$x,x^{\prime}\in\Omega$\textbf{\ now implies }%
$\Lambda_{\lambda}\left(  x-x^{\prime}\right)  =1-\frac{\left\vert
x-x^{\prime}\right\vert }{\lambda}$ \textbf{and there are no zero values.}
This will mean that the smoother is a continuous, piecewise linear function
with a finite number of derivative values which correspond to the linear
segments of the data intervals $\left[  x^{\left(  k\right)  },x^{\left(
k+1\right)  }\right]  $.

Suppose the \textbf{data functions} are drawn from the functions:%
\begin{equation}
\left\{  f\in C_{B}^{\left(  0\right)  }\left(  \Omega\right)  :Df\in
L^{\infty}\left(  \Omega\right)  \right\}  =C_{B}^{\left(  0\right)  }\left(
\Omega\right)  \cap W^{1,\infty}\left(  \Omega\right)  .\label{av075}%
\end{equation}

Then the \textbf{1-dimensional Exact smoother }(or\textbf{\ parameter
stabilized interpolant}) is given by%
\begin{equation}
s\left(  x\right)  =\sum_{k=1}^{N}\alpha_{k}\Lambda_{\lambda}\left(
x-x^{\left(  k\right)  }\right)  \text{,\quad}x\in\mathbb{R}^{1},\label{av073}%
\end{equation}

where $\alpha$ satisfies the \textbf{basis function matrix equation}%
\begin{equation}
\left.
\begin{array}
[c]{ll}
& \left(  N\rho I_{N}+\left(  \Lambda_{\lambda}\right)  _{X,X}\right)
\alpha=\widetilde{\mathcal{E}}_{X}f,\\
where, & \\
& \left(  \Lambda_{\lambda}\right)  _{X,X}=\left(  \Lambda\left(
\frac{x^{\left(  i\right)  }-x^{\left(  j\right)  }}{\lambda}\right)  \right)
,\text{ }\rho\geq0.
\end{array}
\right\} \label{av074}%
\end{equation}

Here $\rho$ denotes the \textbf{smoothing parameter} and
$\widetilde{\mathcal{E}}_{X}f=\left(  f\left(  x^{\left(  k\right)  }\right)
\right)  $ is the \textbf{dependent data}.

We now define the scaled independent data $X_{\lambda}=X/\lambda=\left\{
x^{\left(  k\right)  }/\lambda\right\}  _{k=1}^{N}$ which is contained in the
data region $\Omega_{\lambda}=\Omega/\lambda$ and satisfies
$\operatorname*{diam}\Omega_{\lambda}\leq1$. Also define the data function
$f_{1/\lambda}\left(  x\right)  =f\left(  \lambda x\right)  $ which is clearly
a member of $C_{B}^{\left(  0\right)  }\left(  \Omega_{\lambda}\right)  \cap
W^{1,\infty}\left(  \Omega_{\lambda}\right)  $. Then \ref{av073} and
\ref{av074} become the equivalent matrix equation%
\begin{equation}
\left.
\begin{array}
[c]{ll}
& \left(  N\rho I_{N}+\Lambda_{X_{\lambda},X_{\lambda}}\right)  \alpha
=\widetilde{\mathcal{E}}_{X_{\lambda}}f_{1/\lambda},\\
where, & \\
& \Lambda_{X_{\lambda},X_{\lambda}}=\left(  \Lambda\left(  \frac{x^{\left(
i\right)  }}{\lambda}-\frac{x^{\left(  j\right)  }}{\lambda}\right)  \right)
,\text{ }\rho\geq0,
\end{array}
\right\} \label{av084}%
\end{equation}

with smoother%
\begin{equation}
s_{1/\lambda}\left(  x\right)  =s\left(  \lambda x\right)  =\sum_{k=1}%
^{N}\alpha_{k}\Lambda\left(  x-\frac{x^{\left(  k\right)  }}{\lambda}\right)
,\text{\quad}x\in\mathbb{R}^{1}.\label{av083}%
\end{equation}

Therefore, \textbf{the great bulk of this chapter} will be devoted to the case
where the basis function is the \textbf{unscaled hat function} $\Lambda$ and
the \textbf{data region satisfies} $\operatorname*{diam}\Omega\leq1$. For this
case we will derive the pointwise upper bounds for the derivative given in
Theorems \ref{vThm_bound_deriv_hat_interpol} and
\ref{vThm_bound_deriv_hat_smth} which consider the interpolant $\iota$ and
smoother $s$ respectively. These bounds will turn out to be
\[
\left\Vert D\mathcal{\iota}\right\Vert _{\infty;\Omega}\leq2\left\Vert
f\right\Vert _{\infty;\Omega}+\left\Vert Df\right\Vert _{\infty;\Omega},\quad
N\geq??,\text{ }\rho=0,
\]

and%
\[
\left\Vert Ds\right\Vert _{\infty;\Omega}\leq2\left\Vert f\right\Vert
_{\infty;\Omega}+\min\left\{  5,2+\rho N\right\}  \left\Vert Df\right\Vert
_{\infty;\Omega},\quad N\geq4,\text{ }\rho>0.
\]

Note that when $\rho=0$ we obtain an interpolant derivative bound that is
weaker than that given above. ?? Have I made a mistake or can a stronger
estimate be obtained fairly easily?

Using the preceding remarks concerning scaling, for the case of the scaled hat
basis function $\Lambda\left(  x/\lambda\right)  $ and a data region
satisfying $\operatorname*{diam}\Omega\leq\lambda$, it will be a simple matter
to obtain the interpolant derivative estimate of Corollary
\ref{vCor_bound_deriv_scal_hat_interpol}:
\[
\left\Vert D\iota\right\Vert _{\infty;\Omega}\leq\frac{2}{\lambda}\left\Vert
f\right\Vert _{\infty;\Omega}+\left\Vert Df\right\Vert _{\infty;\Omega},\quad
N\geq??,
\]
and the smoother derivative estimate of Corollary
\ref{vCor_Thm_bound_deriv_scal_hat_smth}:%
\[
\left\Vert Ds\right\Vert _{\infty;\Omega}\leq\frac{2}{\lambda}\left\Vert
f\right\Vert _{\infty;\Omega}+\min\left\{  5,2+\rho N\right\}  \left\Vert
Df\right\Vert _{\infty;\Omega},\quad N\geq4,\text{ }\rho>0.
\]

\textbf{For higher dimensions} it will suffice to bound $D^{\mathbf{1}}\iota$
and $D^{\mathbf{1}}s$.

\subsection{Steps in proof:}

See introduction to each section.

\section{Calculating the smoother derivatives using Cramer's
rule\label{Sect_Cramer_1dim_deriv_smther}}

In this section, using the fact that $\Lambda$\textbf{\ has large support
w.r.t. the data region}, we start by deriving the matrix equation \ref{av076}:
$\mu\left(  f\right)  =A\alpha,$ which expresses the $N+1$ derivative values
$\mu\left(  f\right)  $ of the smoother in terms of the vector $\alpha$ of
basis function coefficients (see \ref{av073}). The matrix $A$ has a simple
inverse with integer elements. The basis function matrix equation \ref{av074}
is then used to express $\mu\left(  f\right)  $ in the form \ref{av078} which
involves the product of the inverse of the symmetric matrix $v\left(
A^{T}\right)  ^{-1}A^{-1}+\left(  A^{T}\right)  ^{-1}\Lambda_{X,X}A^{-1}$ and
the vector $\left(  A^{T}\right)  ^{-1}\widetilde{\mathcal{E}}_{X}f$, the
latter involving function value differences $\Delta_{1}f\left(  x^{\left(
i\right)  }\right)  $ to which we can apply the Taylor expansion result of
Lemma \ref{vLem_Taylor_extension} if we assume that the data functions are in
\ref{av075}. The actual expression is \ref{av002} which involves the inverse
of a bordered matrix which is the bordering of the tridiagonal symmetric
matrix $T_{v}$ given by \ref{av511}. This in turn leads to the formulas
\ref{av006} and \ref{av110} which express $\mu\left(  f\right)  $ in terms of
the inverse of the matrix $T_{v}$ which plays a central role in the sequel.

In Subsection \ref{SbSect_semihomog_Tv} we will simplify the algebra by taking
advantage of the partial homogeneity of $T\left(  \cdot;v\right)  $ in $v$ to
replace $T_{v}=T\left(  \Delta_{1}X;v\right)  $ by $T=T\left(  \frac{1}%
{v}\Delta_{1}X;1\right)  $. This yields the equations \ref{av008} and
\ref{av007} for $\mu\left(  f\right)  $. Next we substitute the value of the
sparse vector $\beta$ and this enables us to express the components of
\ref{av008} and \ref{av007} in terms of the elements of $T^{-1}$.

For arbitrary $a\in\mathbb{R}^{n}$ studying determinants of the form
$\left\vert T\left(  a;1\right)  \right\vert $ is simpler and more compact
than studying determinants of the form $\left\vert T\left(  a;v\right)
\right\vert $.

\subsection{The derivative values and the basis function matrix equation}

The form of the smoother \ref{av073} and the fact that $\Lambda$ has large
support w.r.t. the data region means that the smoother is a continuous,
piecewise linear function with a finite number of derivative values which
correspond to the linear segments of the data intervals $\left[  x^{\left(
k\right)  },x^{\left(  k+1\right)  }\right]  $. We have three cases:\medskip

\fbox{Suppose $x^{\left(  n\right)  }<x<x^{\left(  n+1\right)  }$}
\begin{align*}
s\left(  x\right)  =\sum_{k=1}^{N}\alpha_{k}\Lambda\left(  x-x^{\left(
k\right)  }\right)   & =\sum_{k=1}^{n}\alpha_{k}\Lambda\left(  x-x^{\left(
k\right)  }\right)  +\sum_{k=n+1}^{N}\alpha_{k}\Lambda\left(  x-x^{\left(
k\right)  }\right) \\
& =\sum_{k=1}^{n}\alpha_{k}\left(  1-\left(  x-x^{\left(  k\right)  }\right)
\right)  +\sum_{k=n+1}^{N}\alpha_{k}\left(  1-\left(  x^{\left(  k\right)
}-x\right)  \right)  ,
\end{align*}

so that%
\[
Ds\left(  x\right)  =-\sum_{k=1}^{n}\alpha_{k}+\sum_{k=n+1}^{N}\alpha
_{k},\quad x^{\left(  n\right)  }\leq x\leq x^{\left(  n+1\right)  }.
\]

\fbox{If $\underline{a}\leq x<x^{\left(  1\right)  }$} then%
\[
s\left(  x\right)  =\sum_{k=1}^{N}\alpha_{k}\left(  1-\left(  x^{\left(
k\right)  }-x\right)  \right)  ,
\]

so that%
\[
Ds\left(  x\right)  =\sum_{k=1}^{N}\alpha_{k}.
\]

\fbox{If $x^{\left(  N\right)  }<x\leq\overline{a}$} then%
\[
s\left(  x\right)  =\sum_{k=1}^{N}\alpha_{k}\left(  1-\left(  x-x^{\left(
k\right)  }\right)  \right)  ,
\]

and hence%
\[
Ds\left(  x\right)  =-\sum_{k=1}^{N}\alpha_{k}.
\]

To summarize:%
\begin{equation}
Ds\left(  x\right)  =\left\{
\begin{array}
[c]{ll}%
\sum\limits_{k=1}^{N}\alpha_{k}, & \underline{a}\leq x<x^{\left(  1\right)
},\\
-\sum\limits_{k=1}^{n}\alpha_{k}+\sum\limits_{k=n+1}^{N}\alpha_{k},\text{ } &
x^{\left(  n\right)  }<x<x^{\left(  n+1\right)  },\\
-\sum\limits_{k=1}^{N}\alpha_{k}, & x^{\left(  N\right)  }<x\leq\overline{a}.
\end{array}
\right. \label{av526}%
\end{equation}

The derivative is thus a step function with $N+1$ values given by \ref{av526}
and so the interpolant is the "simplest" piecewise linear interpolant. In this
case the error is easily seen to be order $1$. From \ref{av526},
$\sup_{\left[  a,x^{\left(  1\right)  }\right)  }\left\vert Ds\left(
x\right)  \right\vert =\sup_{\left(  x^{\left(  N\right)  },\overline
{a}\right]  }\left\vert Ds\left(  x\right)  \right\vert $ so%
\begin{equation}
\left\Vert Ds\right\Vert _{\infty;\Omega}=\left\vert \mu\left(  f\right)
\right\vert _{\max},\label{av089}%
\end{equation}

and therefore we need only consider the first $N$ values of the derivative
which are given by
\begin{equation}
\mu\left(  f\right)  :=A\alpha,\label{av076}%
\end{equation}

where%
\begin{equation}
A:=%
\begin{pmatrix}
1 & 1 & 1 &  & 1 & 1\\
-1 & 1 & 1 &  & 1 & 1\\
-1 & -1 & 1 &  & 1 & 1\\
&  &  & \ddots &  & \\
-1 & -1 & -1 &  & 1 & 1\\
-1 & -1 & -1 &  & -1 & 1
\end{pmatrix}
,\quad\operatorname*{size}A=N.\label{av077}%
\end{equation}

We note that%
\begin{equation}
A^{-1}=\frac{1}{2}%
\begin{pmatrix}
1 & -1 & 0 &  & 0 & 0\\
0 & 1 & -1 &  & 0 & 0\\
0 & 0 & 1 &  & 0 & 0\\
&  &  & \ddots &  & \\
0 & 0 & 0 &  & 1 & -1\\
1 & 0 & 0 &  & 0 & 1
\end{pmatrix}
.\label{av841}%
\end{equation}

Set
\begin{equation}
v:=\rho N.\label{av839}%
\end{equation}

Now we write%
\begin{align}
\mu\left(  f\right)   &  =A\alpha\nonumber\\
&  =A\left(  N\rho I_{N}+\Lambda_{X,X}\right)  ^{-1}\widetilde{\mathcal{E}%
}_{X}f\nonumber\\
&  =\left(  \left(  vI_{N}+\Lambda_{X,X}\right)  A^{-1}\right)  ^{-1}%
\widetilde{\mathcal{E}}_{X}f\nonumber\\
&  =\left(  vA^{-1}+\Lambda_{X,X}A^{-1}\right)  ^{-1}\widetilde{\mathcal{E}%
}_{X}f\nonumber\\
&  =\left(  vA^{-1}+\Lambda_{X,X}A^{-1}\right)  ^{-1}A^{T}\left(
A^{T}\right)  ^{-1}\widetilde{\mathcal{E}}_{X}f\nonumber\\
&  =\left\{  \left(  A^{T}\right)  ^{-1}\left(  vA^{-1}+\Lambda_{X,X}%
A^{-1}\right)  \right\}  ^{-1}\left(  A^{T}\right)  ^{-1}%
\widetilde{\mathcal{E}}_{X}f\nonumber\\
&  =\left\{  v\left(  A^{-1}\right)  ^{T}A^{-1}+\left(  A^{-1}\right)
^{T}\Lambda_{X,X}A^{-1}\right\}  ^{-1}\left(  A^{-1}\right)  ^{T}%
\widetilde{\mathcal{E}}_{X}f,\label{av078}%
\end{align}

which is the product of the inverse of the symmetric matrix $v\left(
A^{-1}\right)  ^{T}A^{-1}+\left(  A^{-1}\right)  ^{T}\Lambda_{X,X}A^{-1}$ and
the vector $\left(  A^{-1}\right)  ^{T}\widetilde{\mathcal{E}}_{X}f$, the
latter involving function value differences $\Delta_{1}f\left(  x^{\left(
i\right)  }\right)  $ to which we can apply the Taylor expansion result of
Lemma \ref{vLem_Taylor_extension}, if we assume that the data functions are in
\ref{av075}.

Now using the compact notation%
\begin{equation}
w_{v}=2+v-\left(  x^{\left(  N\right)  }-x^{\left(  1\right)  }\right)
,\label{av090}%
\end{equation}

\begin{equation}
\eta_{i,j}:=\Lambda\left(  x^{\left(  i\right)  }-x^{\left(  j\right)
}\right)  =1-\left\vert x^{\left(  i\right)  }-x^{\left(  j\right)
}\right\vert ,\label{av001}%
\end{equation}

we have%
\[
\Delta_{0,1}\eta_{i,j}:=\eta_{i,j+1}-\eta_{i,j}=\left\{
\begin{array}
[c]{rr}%
-\Delta_{1}x^{\left(  j\right)  }, & i\leq j,\\
\Delta_{1}x^{\left(  j\right)  }, & i>j,
\end{array}
\right.
\]

so that%
\begin{align*}
\Lambda_{X,X}A^{-1} &  =\frac{1}{2}\left(
\begin{array}
[c]{cccccc}%
1 & \eta_{1,2} & \eta_{1,3} &  & \eta_{1,N-1} & \eta_{1,N}\\
\eta_{2,1} & 1 & \eta_{2,3} &  & \eta_{2,N-1} & \eta_{2,N}\\
\eta_{3,1} & \eta_{3,2} & 1 &  & \eta_{3,N-1} & \eta_{3,N}\\
&  &  & \ddots &  & \\
\eta_{N-1,1} & \eta_{N-1,2} & \eta_{N-1,3} &  & 1 & \eta_{N-1,N}\\
\eta_{N,1} & \eta_{N,2} & \eta_{N,3} &  & \eta_{N,N-1} & 1
\end{array}
\right)  \left(
\begin{array}
[c]{cccccc}%
1 & -1 & 0 &  & 0 & 0\\
0 & 1 & -1 &  & 0 & 0\\
0 & 0 & 1 &  & 0 & 0\\
&  &  & \ddots &  & \\
0 & 0 & 0 &  & 1 & -1\\
1 & 0 & 0 &  & 0 & 1
\end{array}
\right) \\
& \\
&  =\frac{1}{2}%
\begin{pmatrix}
w_{0} & -\Delta_{1}x^{\left(  1\right)  } & -\Delta_{1}x^{\left(  2\right)  }
&  & -\Delta_{1}x^{\left(  N-2\right)  } & -\Delta_{1}x^{\left(  N-1\right)
}\\
w_{0} & \Delta_{1}x^{\left(  1\right)  } & -\Delta_{1}x^{\left(  2\right)  } &
& -\Delta_{1}x^{\left(  N-2\right)  } & -\Delta_{1}x^{\left(  N-1\right)  }\\
w_{0} & \Delta_{1}x^{\left(  1\right)  } & \Delta_{1}x^{\left(  2\right)  } &
& -\Delta_{1}x^{\left(  N-2\right)  } & -\Delta_{1}x^{\left(  N-1\right)  }\\
&  &  & \ddots &  & \\
w_{0} & \Delta_{1}x^{\left(  1\right)  } & \Delta_{1}x^{\left(  2\right)  } &
& \Delta_{1}x^{\left(  N-2\right)  } & -\Delta_{1}x^{\left(  N-1\right)  }\\
w_{0} & \Delta_{1}x^{\left(  1\right)  } & \Delta_{1}x^{\left(  2\right)  } &
& \Delta_{1}x^{\left(  N-2\right)  } & \Delta_{1}x^{\left(  N-1\right)  }%
\end{pmatrix}
,
\end{align*}

and so%
\begin{align}
&  \left(  A^{-1}\right)  ^{T}\Lambda_{X,X}A^{-1}\nonumber\\
& \nonumber\\
&  =\frac{1}{4}%
\begin{pmatrix}
1 & 0 & 0 &  & 0 & 1\\
-1 & 1 & 0 &  & 0 & 0\\
0 & -1 & 1 &  & 0 & 0\\
&  &  & \ddots &  & \\
0 & 0 & 0 &  & 1 & 0\\
0 & 0 & 0 &  & -1 & 1
\end{pmatrix}%
\begin{pmatrix}
w_{0} & -\Delta_{1}x^{\left(  1\right)  } & -\Delta_{1}x^{\left(  2\right)  }
&  & -\Delta_{1}x^{\left(  N-2\right)  } & -\Delta_{1}x^{\left(  N-1\right)
}\\
w_{0} & \Delta_{1}x^{\left(  1\right)  } & -\Delta_{1}x^{\left(  2\right)  } &
& -\Delta_{1}x^{\left(  N-2\right)  } & -\Delta_{1}x^{\left(  N-1\right)  }\\
w_{0} & \Delta_{1}x^{\left(  1\right)  } & \Delta_{1}x^{\left(  2\right)  } &
& -\Delta_{1}x^{\left(  N-2\right)  } & -\Delta_{1}x^{\left(  N-1\right)  }\\
&  &  & \ddots &  & \\
w_{0} & \Delta_{1}x^{\left(  1\right)  } & \Delta_{1}x^{\left(  2\right)  } &
& \Delta_{1}x^{\left(  N-2\right)  } & -\Delta_{1}x^{\left(  N-1\right)  }\\
w_{0} & \Delta_{1}x^{\left(  1\right)  } & \Delta_{1}x^{\left(  2\right)  } &
& \Delta_{1}x^{\left(  N-2\right)  } & \Delta_{1}x^{\left(  N-1\right)  }%
\end{pmatrix}
\nonumber\\
& \nonumber\\
&  =\frac{1}{4}%
\begin{pmatrix}
2w_{0} & 0 & 0 &  & 0 & 0\\
0 & 2\Delta_{1}x^{\left(  1\right)  } & 0 &  & 0 & 0\\
0 & 0 & 2\Delta_{1}x^{\left(  2\right)  } &  & 0 & 0\\
&  &  & \ddots &  & \\
0 &  &  &  & 2\Delta_{1}x^{\left(  N-2\right)  } & 0\\
0 &  &  &  & 0 & 2\Delta_{1}x^{\left(  N-1\right)  }%
\end{pmatrix}
\nonumber\\
& \nonumber\\
&  =\frac{1}{2}%
\begin{pmatrix}
w_{0} & 0 & 0 &  & 0 & 0\\
0 & \Delta_{1}x^{\left(  1\right)  } & 0 &  & 0 & 0\\
0 & 0 & \Delta_{1}x^{\left(  2\right)  } &  & 0 & 0\\
&  &  & \ddots &  & \\
0 &  &  &  & \Delta_{1}x^{\left(  N-2\right)  } & 0\\
0 &  &  &  & 0 & \Delta_{1}x^{\left(  N-1\right)  }%
\end{pmatrix}
,\label{av808}%
\end{align}

and%
\begin{align*}
\left(  A^{-1}\right)  ^{T}A^{-1}  & =\frac{1}{4}%
\begin{pmatrix}
1 & 0 & 0 &  & 0 & 1\\
-1 & 1 & 0 &  & 0 & 0\\
0 & -1 & 1 &  & 0 & 0\\
&  &  & \ddots &  & \\
0 & 0 & 0 &  & 1 & 0\\
0 & 0 & 0 &  & -1 & 1
\end{pmatrix}%
\begin{pmatrix}
1 & -1 & 0 &  & 0 & 0\\
0 & 1 & -1 &  & 0 & 0\\
0 & 0 & 1 &  & 0 & 0\\
&  &  & \ddots &  & \\
0 & 0 & 0 &  & 1 & -1\\
1 & 0 & 0 &  & 0 & 1
\end{pmatrix}
\\
& =\frac{1}{4}%
\begin{pmatrix}
2 & -1 & 0 &  & 0 & 1\\
-1 & 2 & -1 &  & 0 & 0\\
0 & -1 & 2 &  & 0 & 0\\
&  &  & \ddots &  & \\
0 & 0 & 0 &  & 2 & -1\\
1 & 0 & 0 &  & -1 & 2
\end{pmatrix}
.
\end{align*}

Now noting \ref{av078},%
\begin{align*}
&  v\left(  A^{-1}\right)  ^{T}A^{-1}+\left(  A^{-1}\right)  ^{T}\Lambda
_{X,X}A^{-1}\\
& \\
&  =\frac{1}{2}%
\begin{pmatrix}
v & -\frac{v}{2} & 0 &  & 0 & \frac{v}{2}\\
-\frac{v}{2} & v & -\frac{v}{2} &  & 0 & 0\\
0 & -\frac{v}{2} & v &  & 0 & 0\\
&  &  & \ddots &  & \\
0 & 0 & 0 &  & v & -\frac{v}{2}\\
\frac{v}{2} & 0 & 0 &  & -\frac{v}{2} & v
\end{pmatrix}
+\frac{1}{2}%
\begin{pmatrix}
w_{0} & 0 & 0 &  & 0 & 0\\
0 & \Delta_{1}x^{\left(  1\right)  } & 0 &  & 0 & 0\\
0 & 0 & \Delta_{1}x^{\left(  2\right)  } &  & 0 & 0\\
&  &  & \ddots &  & \\
0 &  &  &  & \Delta_{1}x^{\left(  N-2\right)  } & 0\\
0 &  &  &  & 0 & \Delta_{1}x^{\left(  N-1\right)  }%
\end{pmatrix}
\\
& \\
&  =\frac{1}{2}%
\begin{pmatrix}
w_{v} & -\frac{v}{2} & 0 &  & 0 & \frac{v}{2}\\
-\frac{v}{2} & v+\Delta_{1}x^{\left(  1\right)  } & -\frac{v}{2} &  & 0 & 0\\
0 & -\frac{v}{2} & v+\Delta_{1}x^{\left(  2\right)  } &  & -\frac{v}{2} & 0\\
&  &  & \ddots &  & \\
0 & 0 & 0 &  & v+\Delta_{1}x^{\left(  N-2\right)  } & -\frac{v}{2}\\
\frac{v}{2} & 0 & 0 &  & -\frac{v}{2} & v+\Delta_{1}x^{\left(  N-1\right)  }%
\end{pmatrix}
,
\end{align*}

so that%
\begin{align}
\mu\left(  f\right)   &  =2%
\begin{pmatrix}
w_{v} & -\frac{v}{2} & 0 &  & 0 & \frac{v}{2}\\
-\frac{v}{2} & v+\Delta_{1}x^{\left(  1\right)  } & -\frac{v}{2} &  & 0 & 0\\
0 & -\frac{v}{2} & v+\Delta_{1}x^{\left(  2\right)  } &  & -\frac{v}{2} & 0\\
&  &  & \ddots &  & \\
0 & 0 & 0 &  & v+\Delta_{1}x^{\left(  N-2\right)  } & -\frac{v}{2}\\
\frac{v}{2} & 0 & 0 &  & -\frac{v}{2} & v+\Delta_{1}x^{\left(  N-1\right)  }%
\end{pmatrix}
^{-1}\left(  A^{-1}\right)  ^{T}\widetilde{\mathcal{E}}_{X}f\nonumber\\
&  =2%
\begin{pmatrix}
w_{v} & b^{T}\\
b & T_{v}%
\end{pmatrix}
^{-1}\left(  A^{-1}\right)  ^{T}\widetilde{\mathcal{E}}_{X}f,\label{av272}%
\end{align}

where, using the \textbf{forward difference} operator $\Delta_{1}$,%
\begin{align}
T_{v}  & :=%
\begin{pmatrix}
v+\Delta_{1}x^{\left(  1\right)  } & -\frac{v}{2} & 0 &  & 0 & 0\\
-\frac{v}{2} & v+\Delta_{1}x^{\left(  2\right)  } & -\frac{v}{2} &  & 0 & 0\\
0 & -\frac{v}{2} & v+\Delta_{1}x^{\left(  3\right)  } &  & 0 & 0\\
&  &  & \ddots &  & \\
0 & 0 & 0 &  & v+\Delta_{1}x^{\left(  N-2\right)  } & -\frac{v}{2}\\
0 & 0 & 0 &  & -\frac{v}{2} & v+\Delta_{1}x^{\left(  N-1\right)  }%
\end{pmatrix}
,\label{av511}\\
b  & :=\frac{v}{2}\beta,\label{av531}\\
\beta & :=%
\begin{pmatrix}
-1 & 0 &  & 0 & 1
\end{pmatrix}
^{T},\quad\operatorname*{size}\beta=N-1\times1.\label{av530}%
\end{align}

Further, from \ref{av841},%
\begin{align*}
\left(  A^{-1}\right)  ^{T}\widetilde{\mathcal{E}}_{X}f=\frac{1}{2}%
\begin{pmatrix}
1 & 0 & 0 &  & 0 & 1\\
-1 & 1 & 0 &  & 0 & 0\\
0 & -1 & 1 &  & 0 & 0\\
&  &  & \ddots &  & \\
0 & 0 & 0 &  & 1 & 0\\
0 & 0 & 0 &  & -1 & 1
\end{pmatrix}%
\begin{pmatrix}
f\left(  x^{\left(  1\right)  }\right) \\
f\left(  x^{\left(  2\right)  }\right) \\
f\left(  x^{\left(  3\right)  }\right) \\
\vdots\\
f\left(  x^{\left(  N-1\right)  }\right) \\
f\left(  x^{\left(  N\right)  }\right)
\end{pmatrix}
&  =\frac{1}{2}%
\begin{pmatrix}
f\left(  x^{\left(  1\right)  }\right)  +f\left(  x^{\left(  N\right)
}\right) \\
\Delta_{1}f\left(  x^{\left(  1\right)  }\right) \\
\Delta_{1}f\left(  x^{\left(  2\right)  }\right) \\
\vdots\\
\Delta_{1}f\left(  x^{\left(  N-2\right)  }\right) \\
\Delta_{1}f\left(  x^{\left(  N-1\right)  }\right)
\end{pmatrix}
\\
& \\
&  =\frac{1}{2}%
\begin{pmatrix}
f\left(  x^{\left(  1\right)  }\right)  +f\left(  x^{\left(  N\right)
}\right) \\
\Delta_{X}f
\end{pmatrix}
,
\end{align*}

where%
\begin{equation}
\Delta_{X}f:=\Delta_{1}\widetilde{\mathcal{E}}_{X}f:=%
\begin{pmatrix}
\Delta_{1}f\left(  x^{\left(  1\right)  }\right) \\
\Delta_{1}f\left(  x^{\left(  2\right)  }\right) \\
\vdots\\
\Delta_{1}f\left(  x^{\left(  N-2\right)  }\right) \\
\Delta_{1}f\left(  x^{\left(  N-1\right)  }\right)
\end{pmatrix}
,\label{av528}%
\end{equation}

so that%
\begin{equation}
\mu\left(  f\right)  =%
\begin{pmatrix}
w_{v} & b^{T}\\
b & T_{v}%
\end{pmatrix}
^{-1}%
\begin{pmatrix}
f\left(  x^{\left(  1\right)  }\right)  +f\left(  x^{\left(  N\right)
}\right) \\
\Delta_{X}f
\end{pmatrix}
\label{av002}%
\end{equation}

\subsection{Basis function interpolation}

The case $\rho=0$ corresponds to the basis function interpolant
$\mathcal{\iota}$. Considering $\rho=0$ allows us to introduce some notation
and results in a simplified context as well as acting as a check. In fact,
using equations \ref{av090} and \ref{av511} to \ref{av530},%
\begin{align*}
\mu\left(  f\right)   & =%
\begin{pmatrix}
w_{0} & 0^{T}\\
0 & T_{0}%
\end{pmatrix}
^{-1}%
\begin{pmatrix}
f\left(  x^{\left(  1\right)  }\right)  +f\left(  x^{\left(  N\right)
}\right) \\
\Delta_{X}f
\end{pmatrix}
\\
& =%
\begin{pmatrix}
w_{0} & 0^{T}\\
0 & \operatorname*{diag}\left(  \Delta_{1}X\right)
\end{pmatrix}
^{-1}%
\begin{pmatrix}
f\left(  x^{\left(  1\right)  }\right)  +f\left(  x^{\left(  N\right)
}\right) \\
\Delta_{X}f
\end{pmatrix}
\\
& =%
\begin{pmatrix}
\frac{1}{w_{0}} & 0^{T}\\
0 & \operatorname*{diag}\left(  1./\Delta_{1}X\right)
\end{pmatrix}%
\begin{pmatrix}
f\left(  x^{\left(  1\right)  }\right)  +f\left(  x^{\left(  N\right)
}\right) \\
\Delta_{X}f
\end{pmatrix}
\\
& =%
\begin{pmatrix}
\frac{f\left(  x^{\left(  1\right)  }\right)  +f\left(  x^{\left(  N\right)
}\right)  }{2-\left(  x^{\left(  N\right)  }-x^{\left(  1\right)  }\right)
}\\
\widetilde{D}_{X}f
\end{pmatrix}
,
\end{align*}

where we have used the component-wise notation%
\begin{equation}
\widetilde{D}_{X}f:=\frac{\Delta_{X}f}{\Delta_{1}X}=\left(  \frac{f\left(
x^{\left(  k+1\right)  }\right)  -f\left(  x^{\left(  k\right)  }\right)
}{x^{\left(  k+1\right)  }-x^{\left(  k\right)  }}\right)  _{k=1}%
^{N-1}.\label{av082}%
\end{equation}

We will now need the following result to bound the components of
$\widetilde{D}_{X}f$:

\begin{lemma}
\label{vLem_Taylor_extension}\textbf{Taylor series expansion} Suppose
$\Omega=\left(  a,b\right)  \subset\mathbb{R}^{1}$ and $D$ is the
distributional derivative. Then:
\end{lemma}

\begin{lemma}
$f\in C_{B}^{\left(  0\right)  }\left(  \Omega\right)  $ and $Df\in L^{\infty
}\left(  \Omega\right)  $ implies:%
\begin{equation}
f(y)=f(x)+\left(  y-x\right)  \int_{0}^{1}(Df)(\left(  1-t\right)
x+ty)dt,\quad x,y\in\Omega,\label{av091}%
\end{equation}

\begin{enumerate}
\item and%
\begin{equation}
\left\vert \frac{f\left(  y\right)  -f\left(  x\right)  }{y-x}\right\vert
\leq\left\Vert Df\right\Vert _{\infty},\quad x,y\in\Omega.\label{av080}%
\end{equation}

\end{enumerate}
\end{lemma}

\begin{proof}
Since we are working in one dimension $f\in C_{B}^{\left(  0\right)  }\left(
\Omega\right)  $ can be easily extended to $f\in C_{B}^{\left(  0\right)
}\left(  \mathbb{R}^{1}\right)  $ and then an application of Lemma
\ref{Lem_Taylor_extension} yields \ref{av091}. Now \ref{av080} follows
directly from \ref{av091}.
\end{proof}

Since $0<x^{\left(  N\right)  }-x^{\left(  1\right)  }<1$, if $f$ lies in the
space \ref{av075} then
\[
\left\vert \mu_{1}\left(  f\right)  \right\vert \leq\frac{\left\vert f\left(
x^{\left(  1\right)  }\right)  \right\vert +\left\vert f\left(  x^{\left(
N\right)  }\right)  \right\vert }{2-\left(  x^{\left(  N\right)  }-x^{\left(
1\right)  }\right)  }\leq2\left\Vert f\right\Vert _{\infty;\Omega},
\]

and applying Lemma \ref{vLem_Taylor_extension} gives%
\[
\left\vert \mu^{\prime}\left(  f\right)  \right\vert _{\max}\leq\left\Vert
Df\right\Vert _{\infty;\Omega},
\]

so that%
\[
\left\vert \mu\left(  f\right)  \right\vert _{\max}\leq2\left\Vert
f\right\Vert _{\infty;\Omega}+\left\Vert Df\right\Vert _{\infty;\Omega}.
\]

But the case $x^{\left(  N\right)  }<x\leq\overline{a}$ from \ref{av526}
implies that $\sup\limits_{\left[  a,x^{\left(  1\right)  }\right)
}\left\vert D\mathcal{\iota}\left(  x\right)  \right\vert =\sup
\limits_{\left(  x^{\left(  N\right)  },\overline{a}\right]  }\left\vert
D\mathcal{\iota}\left(  x\right)  \right\vert $ and so%
\[
\left\Vert D\mathcal{\iota}\right\Vert _{\infty;\Omega}\leq2\left\Vert
f\right\Vert _{\infty;\Omega}+\left\Vert Df\right\Vert _{\infty;\Omega}.
\]

Thus we have proved:

\begin{theorem}
\label{vThm_bound_deriv_hat_interpol}\textbf{Hat basis function interpolation}
Suppose the data region satisfies $\operatorname*{diam}\Omega\leq1$ and the
data function $f$ lies in the space $\left\{  f\in C_{B}^{\left(  0\right)
}\left(  \Omega\right)  :Df\in L^{\infty}\left(  \Omega\right)  \right\}  $.

Then the basis function interpolant $\mathcal{\iota}$ corresponding to the hat
basis function $\Lambda$ satisfies%
\begin{equation}
\left\Vert D\mathcal{\iota}\right\Vert _{\infty;\Omega}\leq2\left\Vert
f\right\Vert _{\infty;\Omega}+\left\Vert Df\right\Vert _{\infty;\Omega
}.\label{av088}%
\end{equation}

\end{theorem}

Regarding the \textbf{scaled} hat basis function:

\begin{corollary}
\label{vCor_bound_deriv_scal_hat_interpol}\textbf{Scaled hat basis function
interpolation} Suppose the data region satisfies $\operatorname*{diam}%
\Omega\leq\lambda$ and the data function $f$ lies in the space $\left\{  f\in
C_{B}^{\left(  0\right)  }\left(  \Omega\right)  :Df\in L^{\infty}\left(
\Omega\right)  \right\}  $.

Then the basis function interpolant $\iota$ corresponding to the scaled hat
basis function $\Lambda\left(  \cdot/\lambda\right)  $ satisfies%
\begin{equation}
\left\Vert D\mathcal{\iota}\right\Vert _{\infty;\Omega}\leq\frac{2}{\lambda
}\left\Vert f\right\Vert _{\infty;\Omega}+\left\Vert Df\right\Vert
_{\infty;\Omega}.\label{av085}%
\end{equation}

\end{corollary}

\begin{proof}
From \ref{av084} and \ref{av083} of the introduction to this chapter it
follows that $\mathcal{\iota}_{1/\lambda}\left(  x\right)  =\mathcal{\iota
}\left(  \lambda x\right)  $ is the hat basis function interpolant associated
with data function $f_{1/\lambda}\left(  x\right)  =f\left(  \lambda x\right)
$ and data $X_{\lambda}\subset\Omega_{\lambda}$ satisfying
$\operatorname*{diam}X_{\lambda}\leq1$. By Theorem
\ref{vThm_bound_deriv_hat_interpol},%
\begin{align*}
\left\Vert D\mathcal{\iota}_{1/\lambda}\right\Vert _{\infty;\Omega_{\lambda}}
& \leq2\left\Vert f_{1/\lambda}\right\Vert _{\infty;\Omega_{\lambda}%
}+\left\Vert Df_{1/\lambda}\right\Vert _{\infty;\Omega_{\lambda}}\\
& =2\left\Vert f\right\Vert _{\infty;\Omega}+\left\Vert Df_{1/\lambda
}\right\Vert _{\infty;\Omega_{\lambda}}\\
& =2\left\Vert f\right\Vert _{\infty;\Omega}+\lambda\left\Vert \left(
Df\right)  _{1/\lambda}\right\Vert _{\infty;\Omega_{\lambda}}\\
& =2\left\Vert f\right\Vert _{\infty;\Omega}+\lambda\left\Vert Df\right\Vert
_{\infty;\Omega},
\end{align*}

and since%
\[
\left\Vert D\mathcal{\iota}_{1/\lambda}\right\Vert _{\infty;\Omega_{\lambda}%
}=\lambda\left\Vert \left(  D\mathcal{\iota}\right)  _{1/\lambda}\right\Vert
_{\infty;\Omega_{\lambda}}=\lambda\left\Vert D\mathcal{\iota}\right\Vert
_{\infty;\Omega},
\]

it follows that%
\[
\left\Vert D\mathcal{\iota}\right\Vert _{\infty;\Omega}\leq\frac{2}{\lambda
}\left\Vert f\right\Vert _{\infty;\Omega}+\left\Vert Df\right\Vert
_{\infty;\Omega}.
\]

\end{proof}

\begin{remark}
\label{vRem_Cor_bound_deriv_scal_hat_interpol}When $\operatorname*{diam}%
\Omega\leq\lambda$ increasing $\lambda$ only changes the interpolant when
$x<\underline{a}$ and $x>\overline{a}$.
\end{remark}

\subsection{The inverse of a bordered matrix}

To calculate $\mu$ using \ref{av002} we will derive the formula \ref{av513}
which expresses the inverse of the bordered matrix $%
\begin{pmatrix}
w_{v} & b^{T}\\
b & T_{v}%
\end{pmatrix}
$ in terms of $T_{v}^{-1}$. This in turn leads to the formulas \ref{av006} and
\ref{av110} which express $\mu\left(  f\right)  $ in terms of the inverse of
the matrix $T_{v}$ which plays a central role in the sequel. Suppose%
\begin{equation}%
\begin{pmatrix}
w_{v} & b^{T}\\
b & T_{v}%
\end{pmatrix}%
\begin{pmatrix}
z & c^{T}\\
c & M
\end{pmatrix}
=%
\begin{pmatrix}
1 & O\\
O & I
\end{pmatrix}
.\label{av295}%
\end{equation}

Then we have the equations%
\begin{align*}
w_{v}z+b\mathbf{\cdot}c  & =1,\\
zb+T_{v}c  & =O,\\
w_{v}c^{T}+b^{T}M  & =O,\\
bc^{T}+T_{v}M  & =I.
\end{align*}

Thus $c=-zT_{v}^{-1}b$ and so $w_{v}z+b\mathbf{\cdot}c=w_{v}z-zb^{T}T_{v}%
^{-1}b=1$ which implies that%
\begin{equation}%
\begin{array}
[c]{ll}%
z= & \frac{1}{w_{v}-b^{T}T_{v}^{-1}b},\\
c= & -\frac{1}{w_{v}-b^{T}T_{v}^{-1}b}T_{v}^{-1}b,
\end{array}
\label{av509}%
\end{equation}

and%
\begin{align}
M=T_{v}^{-1}\left(  I-bc^{T}\right)   &  =T_{v}^{-1}\left(  I+\frac{1}%
{w_{v}-b^{T}T_{v}^{-1}b}b\left(  T_{v}^{-1}b\right)  ^{T}\right) \nonumber\\
&  =T_{v}^{-1}\left(  I+\frac{1}{w_{v}-b^{T}T_{v}^{-1}b}bb^{T}T_{v}%
^{-1}\right)  .\label{av510}%
\end{align}

Thus%
\begin{align}%
\begin{pmatrix}
w_{v} & b^{T}\\
b & T_{v}%
\end{pmatrix}
^{-1} &  =%
\begin{pmatrix}
z & c^{T}\\
c & M
\end{pmatrix}
\nonumber\\
&  =\left(
\begin{array}
[c]{ccc}%
\frac{1}{w_{v}-b^{T}T_{v}^{-1}b} & \vdots & -\left(  \frac{1}{w_{v}-b^{T}%
T_{v}^{-1}b}T_{v}^{-1}b\right)  ^{T}\\
-\frac{1}{w_{v}-b^{T}T_{v}^{-1}b}T_{v}^{-1}b & \vdots & T_{v}^{-1}\left(
I+\frac{1}{w_{v}-b^{T}T_{v}^{-1}b}bb^{T}T_{v}^{-1}\right)
\end{array}
\right) \nonumber\\
&  =\frac{1}{w_{v}-b^{T}T_{v}^{-1}b}\left(
\begin{array}
[c]{ccc}%
1 & \vdots & -\left(  T_{v}^{-1}b\right)  ^{T}\\
-T_{v}^{-1}b & \vdots & T_{v}^{-1}\left(  \left(  w_{v}-b^{T}T_{v}%
^{-1}b\right)  I+bb^{T}T_{v}^{-1}\right)
\end{array}
\right) \nonumber\\
&  =\frac{1}{w_{v}-b^{T}T_{v}^{-1}b}\left(
\begin{array}
[c]{ccc}%
1 & \vdots & -b^{T}T_{v}^{-1}\\
-T_{v}^{-1}b & \vdots & \left(  w_{v}-b^{T}T_{v}^{-1}b\right)  T_{v}%
^{-1}+T_{v}^{-1}bb^{T}T_{v}^{-1}%
\end{array}
\right)  .\label{av513}%
\end{align}

Check:%
\begin{align*}
&
\begin{pmatrix}
w_{v} & b^{T}\\
b & T_{v}%
\end{pmatrix}
\frac{1}{w_{v}-b^{T}T_{v}^{-1}b}\left(
\begin{array}
[c]{ccc}%
1 & \vdots & -b^{T}T_{v}^{-1}\\
-T_{v}^{-1}b & \vdots & \left(  w_{v}-b^{T}T_{v}^{-1}b\right)  T_{v}%
^{-1}+T_{v}^{-1}bb^{T}T_{v}^{-1}%
\end{array}
\right) \\
&  =\frac{1}{w_{v}-b^{T}T_{v}^{-1}b}%
\begin{pmatrix}
w_{v} & b^{T}\\
b & T_{v}%
\end{pmatrix}
\left(
\begin{array}
[c]{ccc}%
1 & \vdots & -b^{T}T_{v}^{-1}\\
-T_{v}^{-1}b & \vdots & \left(  w_{v}-b^{T}T_{v}^{-1}b\right)  T_{v}%
^{-1}+T_{v}^{-1}bb^{T}T_{v}^{-1}%
\end{array}
\right) \\
&  =\frac{1}{w_{v}-b^{T}T_{v}^{-1}b}\left(
\begin{array}
[c]{ccc}%
w_{v}-b^{T}T_{v}^{-1}b & \vdots & -w_{v}b^{T}T_{v}^{-1}+b^{T}\left(  \left(
w_{v}-b^{T}T_{v}^{-1}b\right)  T_{v}^{-1}+T_{v}^{-1}bb^{T}T_{v}^{-1}\right) \\
O & \vdots & -bb^{T}T_{v}^{-1}+T\left(  w_{v}-b^{T}T_{v}^{-1}b\right)
T_{v}^{-1}+bb^{T}T_{v}^{-1}%
\end{array}
\right) \\
&  =\frac{1}{w_{v}-b^{T}T_{v}^{-1}b}\left(
\begin{array}
[c]{ccc}%
w_{v}-b^{T}T_{v}^{-1}b & \vdots & -w_{v}b^{T}T_{v}^{-1}+\left(  w_{v}%
-b^{T}T_{v}^{-1}b\right)  b^{T}T_{v}^{-1}+b^{T}T_{v}^{-1}bb^{T}T_{v}^{-1}\\
O & \vdots & w_{v}-b^{T}T_{v}^{-1}b
\end{array}
\right) \\
&  =I_{N}.
\end{align*}

Substituting \ref{av513} into \ref{av002}:%
\begin{align*}
\mu\left(  f\right)   & =%
\begin{pmatrix}
w_{v} & b^{T}\\
b & T_{v}%
\end{pmatrix}
^{-1}%
\begin{pmatrix}
f\left(  x^{\left(  1\right)  }\right)  +f\left(  x^{\left(  N\right)
}\right) \\
\Delta_{X}f
\end{pmatrix}
\\
& =\frac{1}{w_{v}-b^{T}T_{v}^{-1}b}\left(
\begin{array}
[c]{ccc}%
1 & \vdots & -b^{T}T_{v}^{-1}\\
-T_{v}^{-1}b & \vdots & \left(  w_{v}-b^{T}T_{v}^{-1}b\right)  T_{v}%
^{-1}+T_{v}^{-1}bb^{T}T_{v}^{-1}%
\end{array}
\right)
\begin{pmatrix}
f\left(  x^{\left(  1\right)  }\right)  +f\left(  x^{\left(  N\right)
}\right) \\
\Delta_{X}f
\end{pmatrix}
\\
& =\frac{1}{w_{v}-b^{T}T_{v}^{-1}b}\left(
\begin{array}
[c]{c}%
f\left(  x^{\left(  1\right)  }\right)  +f\left(  x^{\left(  N\right)
}\right)  -b^{T}T_{v}^{-1}\Delta_{X}f\\
u
\end{array}
\right)  ,
\end{align*}

where%
\begin{align*}
u  & =-\left(  f\left(  x^{\left(  1\right)  }\right)  +f\left(  x^{\left(
N\right)  }\right)  \right)  T_{v}^{-1}b+\left(  \left(  w_{v}-b^{T}T_{v}%
^{-1}b\right)  T_{v}^{-1}+T_{v}^{-1}bb^{T}T_{v}^{-1}\right)  \Delta_{X}f\\
& =-\left(  f\left(  x^{\left(  1\right)  }\right)  +f\left(  x^{\left(
N\right)  }\right)  \right)  T_{v}^{-1}b+\left(  w_{v}-b^{T}T_{v}%
^{-1}b\right)  T_{v}^{-1}\Delta_{X}f+T_{v}^{-1}b\left(  b^{T}T_{v}^{-1}%
\Delta_{X}f\right) \\
& =-\left(  f\left(  x^{\left(  1\right)  }\right)  +f\left(  x^{\left(
N\right)  }\right)  \right)  T_{v}^{-1}b+\left(  w_{v}-b^{T}T_{v}%
^{-1}b\right)  T_{v}^{-1}\Delta_{X}f+\left(  b^{T}T_{v}^{-1}\Delta
_{X}f\right)  T_{v}^{-1}b\\
& =\left(  w_{v}-b^{T}T_{v}^{-1}b\right)  T_{v}^{-1}\Delta_{X}f-\left(
f\left(  x^{\left(  1\right)  }\right)  +f\left(  x^{\left(  N\right)
}\right)  -b^{T}T_{v}^{-1}\Delta_{X}f\right)  T_{v}^{-1}b,
\end{align*}

so that%
\begin{align}
\mu\left(  f\right)   & =%
\begin{pmatrix}
\frac{f\left(  x^{\left(  1\right)  }\right)  +f\left(  x^{\left(  N\right)
}\right)  -b^{T}T_{v}^{-1}\Delta_{X}f}{w_{v}-b^{T}T_{v}^{-1}b}\\
T_{v}^{-1}\Delta_{X}f-\frac{f\left(  x^{\left(  1\right)  }\right)  +f\left(
x^{\left(  N\right)  }\right)  -b^{T}T_{v}^{-1}\Delta_{X}f}{w_{v}-b^{T}%
T_{v}^{-1}b}T_{v}^{-1}b
\end{pmatrix}
\label{av503}\\
& =%
\begin{pmatrix}
\mu_{1}\left(  f\right) \\
T_{v}^{-1}\Delta_{X}f-\mu_{1}\left(  f\right)  T_{v}^{-1}b
\end{pmatrix}
\nonumber\\
& =%
\begin{pmatrix}
\mu_{1}\left(  f\right) \\
T_{v}^{-1}\left(  \Delta_{X}f-\mu_{1}\left(  f\right)  b\right)
\end{pmatrix}
.\nonumber
\end{align}

where%
\begin{align}
\mu_{1}\left(  f\right)   & =\frac{f\left(  x^{\left(  1\right)  }\right)
+f\left(  x^{\left(  N\right)  }\right)  -b^{T}T_{v}^{-1}\Delta_{X}f}%
{w_{v}-b^{T}T_{v}^{-1}b},\label{av006}\\
\mu^{\prime}\left(  f\right)   & :=\left(  \mu_{k}\right)  _{k=2}^{N}%
=T_{v}^{-1}\Delta_{X}f-\mu_{1}\left(  f\right)  T_{v}^{-1}b.\label{av110}%
\end{align}

\subsection{Replacing $T_{v}:=T\left(  \Delta_{1}X;v\right)  $ by $T=T\left(
\frac{1}{v}\Delta_{1}X;1\right)  $\label{SbSect_semihomog_Tv}}

We will simplify the algebra by taking advantage of the partial homogeneity of
$T\left(  \cdot;v\right)  $ in $v$ to replace $T_{v}=T\left(  \Delta
_{1}X;v\right)  $ by $T=T\left(  \frac{1}{v}\Delta_{1}X;1\right)  $. This
yields the equations \ref{av008} and \ref{av007} for $\mu\left(  f\right)  $.
Next we substitute the value of the sparse vector $\beta$ and this enables us
to express the components of \ref{av008} and \ref{av007} in terms of the
elements of $T^{-1}$. Studying determinants of the form $\left\vert T\left(
a;1\right)  \right\vert $ is simpler and more compact than studying
determinants of the form $\left\vert T\left(  a;v\right)  \right\vert $ where
$a\in\mathbb{R}^{n}$.

In the last subsection we use Cramer's rule for matrix inversion to express
the components of the formulas for $\mu\left(  f\right)  $ (\ref{av008} and
\ref{av007}) in terms of the minors of $T$. These are equations \ref{av725} to
\ref{av757}.

Noting the formula \ref{av511} for $T_{v}$ we introduce the following definitions:

\begin{definition}
\label{vDef_T_minors}\textbf{Matrix notation }Suppose\textbf{\ }$a=\left(
a_{1},\ldots,a_{n}\right)  \in\mathbb{R}^{n}$.%
\begin{align}
\Delta_{1}X  & :=\left(  \Delta_{1}x^{\left(  k\right)  }\right)  _{k=1}%
^{N-1},\label{av055}\\
T\left(  a;v\right)   & :=%
\begin{pmatrix}
v+a_{1} & -\frac{v}{2} & 0 &  & 0 & 0\\
-\frac{v}{2} & v+a_{2} & -\frac{v}{2} &  & 0 & 0\\
0 & -\frac{v}{2} & v+a_{3} &  & 0 & 0\\
&  &  & \ddots &  & \\
0 & 0 & 0 &  & v+a_{n-1} & -\frac{v}{2}\\
0 & 0 & 0 &  & -\frac{v}{2} & v+a_{n}%
\end{pmatrix}
,\label{av054}\\
T  & :=T\left(  \frac{1}{v}\Delta_{1}X;1\right)  :=T\left(  \frac{1}{v}%
\Delta_{1}X\right)  ,\nonumber\\
T_{m:n}  & :=T\left(  \frac{1}{v}\left(  \Delta_{1}x^{\left(  k\right)
}\right)  _{k=m}^{n-1};1\right)  :=T\left(  \frac{1}{v}\left(  \Delta
_{1}x^{\left(  k\right)  }\right)  _{k=m}^{n-1}\right)  ,\nonumber\\
T\left(  a\right)   & :=T\left(  a;1\right)  ,\nonumber\\
a_{m:n}  & :=\left(  a_{m},a_{m+1},\ldots,a_{n}\right)  \text{ }if\text{
}1\leq m\leq n,\nonumber\\
a_{m:n}  & :=\left(  a_{m},a_{m-1},\ldots,a_{n}\right)  \text{ }if\text{
}1\leq n\leq m,\nonumber\\
T_{m:n}\left(  a\right)   & :=T\left(  a_{m:n}\right)  ,\nonumber\\
T\left(  a_{1:0}\right)   & :=1,\label{av057}\\
T\left(  a_{n+1:n}\right)   & :=1,\label{av050}\\
0_{0}  & :=\left\{  {}\right\}  ,\text{ the empty sequence},\nonumber\\
T\left(  \left\{  {}\right\}  \right)   & :=0.\nonumber
\end{align}

\textbf{Minor notation}%
\[%
\begin{array}
[c]{ll}%
\underline{Matrix} & \underline{Minors}\\
T\left(  a;v\right)  & M_{i,j}\left(  a;v\right) \\
T & M_{i,j}\\
T\left(  a\right)  & M_{i,j}\left(  a\right)
\end{array}
\]

\end{definition}

\begin{remark}
\label{vRem_Def_T_minors}\textbf{Homogeneity} We have the \textbf{partial
homogeneity} result%
\begin{equation}
\left\vert T\left(  a;v\right)  \right\vert =v^{m}\left\vert T\left(  \frac
{1}{v}a;1\right)  \right\vert ,\quad a\in\mathbb{R}^{m},\label{av826}%
\end{equation}

and the \textbf{homogeneity} result%
\begin{equation}
\left\vert T\left(  ta;tv\right)  \right\vert =t^{m}\left\vert T\left(
a;1\right)  \right\vert ,\quad a\in\mathbb{R}^{m},\label{av825}%
\end{equation}

and so from \ref{av511},%
\begin{equation}%
\begin{array}
[c]{ll}%
T_{v} & =T\left(  \Delta_{1}X;v\right)  =vT\left(  \frac{1}{v}\Delta
_{1}X;1\right)  =vT,\\
T_{v}^{-1} & =\frac{1}{v}T^{-1}.
\end{array}
\label{av011}%
\end{equation}

\end{remark}

We now apply equations \ref{av011} to \ref{av006} and \ref{av110}. Regarding
$\mu^{\prime}$:%
\begin{equation}
T_{v}^{-1}b=\frac{1}{v}T^{-1}b=\frac{1}{v}T^{-1}\left(  \frac{v}{2}%
\beta\right)  =\frac{1}{2}T^{-1}\beta,\label{av009}%
\end{equation}

and%
\begin{equation}
T_{v}^{-1}\Delta_{X}f=\frac{1}{v}T^{-1}\Delta_{X}f,\label{av010}%
\end{equation}

so that%
\[
\mu^{\prime}\left(  f\right)  =T_{v}^{-1}\Delta_{X}f-\mu_{1}\left(  f\right)
T_{v}^{-1}b=\frac{1}{v}T^{-1}\Delta_{X}f-\frac{1}{2}\mu_{1}\left(  f\right)
T^{-1}\beta.
\]

Regarding $\mu_{1}$:%
\[
b^{T}T_{v}^{-1}\Delta_{X}f=\left(  \frac{v}{2}\beta\right)  ^{T}\frac{1}%
{v}T^{-1}\Delta_{X}f=\frac{1}{2}\beta^{T}T^{-1}\Delta_{X}f,
\]

and%
\[
b^{T}T_{v}^{-1}b=\left(  \frac{v}{2}\beta\right)  ^{T}\frac{1}{2}T^{-1}%
\beta=\frac{v}{4}\beta^{T}T^{-1}\beta.
\]

so equation \ref{av006} now becomes:%
\begin{align*}
\mu_{1}\left(  f\right)   & =\frac{f\left(  x^{\left(  1\right)  }\right)
+f\left(  x^{\left(  N\right)  }\right)  -b^{T}T_{v}^{-1}\Delta_{X}f}%
{w_{v}-b^{T}T_{v}^{-1}b}\\
& =\frac{f\left(  x^{\left(  1\right)  }\right)  +f\left(  x^{\left(
N\right)  }\right)  -\frac{1}{2}\beta^{T}T^{-1}\Delta_{X}f}{w_{v}-\frac{v}%
{2}\beta^{T}T^{-1}\beta}\\
& =\frac{f\left(  x^{\left(  1\right)  }\right)  +f\left(  x^{\left(
N\right)  }\right)  -\frac{1}{2}\beta^{T}T^{-1}\Delta_{X}f}{2-\left(
x^{\left(  N\right)  }-x^{\left(  1\right)  }\right)  +v-\frac{v}{4}\beta
^{T}T^{-1}\beta}\\
& =\frac{f\left(  x^{\left(  1\right)  }\right)  +f\left(  x^{\left(
N\right)  }\right)  -\frac{1}{2}\beta^{T}T^{-1}\Delta_{X}f}{2-\mathbf{1}%
\cdot\Delta_{1}X+v-\frac{v}{4}\beta^{T}T^{-1}\beta}\\
& =\frac{1}{v}\frac{f\left(  x^{\left(  1\right)  }\right)  +f\left(
x^{\left(  N\right)  }\right)  -\frac{1}{2}\beta^{T}T^{-1}\Delta_{X}f}{\left(
\frac{2}{v}-\mathbf{1}\cdot\frac{1}{v}\Delta_{1}X+1\right)  -\frac{1}{4}%
\beta^{T}T^{-1}\beta}.
\end{align*}

To summarize:%
\begin{align}
\mu_{1}\left(  f\right)   & =\frac{1}{v}\frac{f\left(  x^{\left(  1\right)
}\right)  +f\left(  x^{\left(  N\right)  }\right)  -\frac{1}{2}\beta^{T}%
T^{-1}\Delta_{X}f}{\left(  \frac{2}{v}-\mathbf{1}\cdot\frac{1}{v}\Delta
_{1}X+1\right)  -\frac{1}{4}\beta^{T}T^{-1}\beta},\label{av008}\\
\mu^{\prime}\left(  f\right)   & =\left(  \mu_{k}\left(  f\right)  \right)
_{k=2}^{N}=\frac{1}{v}T^{-1}\Delta_{X}f-\frac{1}{2}\mu_{1}\left(  f\right)
T^{-1}\beta.\label{av007}%
\end{align}

\textbf{The goal of this chapter} is to uniformly bound $\left\Vert
\mu\right\Vert _{\max}$ and we will bound $T^{-1}\beta$ in Section
\ref{Sect_bound_invTb}, bound $\frac{1}{v}T^{-1}\Delta_{X}f$ in Section
\ref{Sect_bound_invTdiff_data_fn} and bound $\mu_{1}\left(  f\right)  $ in
Subsection \ref{SbSect_bound_mu1}.

Now we substitute the value of $\beta$ into the expressions involving $T^{-1}$
in \ref{av008} and \ref{av007}:
\begin{equation}
T^{-1}\beta=T^{-1}%
\begin{pmatrix}
-1\\
0\\
\vdots\\
0\\
1
\end{pmatrix}
=%
\begin{pmatrix}
\left(  T^{-1}\right)  _{1,N-1}-\left(  T^{-1}\right)  _{1,1}\\
\vdots\\
\left(  T^{-1}\right)  _{k,N-1}-\left(  T^{-1}\right)  _{k,1}\\
\vdots\\
\left(  T^{-1}\right)  _{N-1,N-1}-\left(  T^{-1}\right)  _{N-1,1}%
\end{pmatrix}
,\label{av527}%
\end{equation}

so%
\[
\left(  T^{-1}\beta\right)  _{k}=\left(  \left(  T^{-1}\right)  _{k,N-1}%
-\left(  T^{-1}\right)  _{k,1}\right)  ,\qquad1\leq k\leq N-1.
\]

Also%
\[
\left(  T^{-1}\Delta_{X}f\right)  _{k}=\sum\limits_{n=1}^{N-1}\left(
T^{-1}\right)  _{k,n}\Delta_{1}f\left(  x^{\left(  n\right)  }\right)  .
\]

Further%
\begin{align*}
\beta^{T}T^{-1}\beta &  =%
\begin{pmatrix}
-1 & 0 & \ldots & 0 & 1
\end{pmatrix}
T^{-1}%
\begin{pmatrix}
-1 & 0 & \ldots & 0 & 1
\end{pmatrix}
^{T}\\
&  =\left(  T^{-1}\right)  _{1,1}-\left(  T^{-1}\right)  _{1,N-1}-\left(
T^{-1}\right)  _{N-1,1}+\left(  T^{-1}\right)  _{N-1,N-1}\\
&  =\left(  T^{-1}\right)  _{1,1}-2\left(  T^{-1}\right)  _{1,N-1}+\left(
T^{-1}\right)  _{N-1,N-1},
\end{align*}

and%
\[
\beta^{T}T^{-1}\Delta_{X}f=\left(  \Delta_{X}f\right)  ^{T}T^{-1}\beta
=\sum\limits_{n=1}^{N-1}\left(  \left(  T^{-1}\right)  _{n,N-1}-\left(
T^{-1}\right)  _{n,1}\right)  \Delta_{1}f\left(  x^{\left(  n\right)
}\right)  .
\]

To summarize: for $1\leq k\leq N-1$,%
\begin{align}
\left(  T^{-1}\beta\right)  _{k}  & =\left(  \left(  T^{-1}\right)
_{k,N-1}-\left(  T^{-1}\right)  _{k,1}\right)  ,\label{av534}\\
\left(  T^{-1}\Delta_{X}f\right)  _{k}  & =\sum\limits_{n=1}^{N-1}\left(
T^{-1}\right)  _{k,n}\Delta_{1}f\left(  x^{\left(  n\right)  }\right)
,\label{av539}\\
\beta^{T}T^{-1}\beta & =\left(  T^{-1}\right)  _{1,1}-2\left(  T^{-1}\right)
_{1,N-1}+\left(  T^{-1}\right)  _{N-1,N-1},\label{av514}\\
\beta^{T}T^{-1}\Delta_{X}f  & =\sum\limits_{n=1}^{N-1}\left(  \left(
T^{-1}\right)  _{n,N-1}-\left(  T^{-1}\right)  _{n,1}\right)  \Delta
_{1}f\left(  x^{\left(  n\right)  }\right)  .\label{av529}%
\end{align}

\begin{remark}
\label{vRem_invTdf_invTb_bounded_implications}Suppose $\mu_{1}\left(
f\right)  =\frac{f\left(  x^{\left(  1\right)  }\right)  +f\left(  x^{\left(
N\right)  }\right)  -b^{T}T^{-1}\Delta_{X}f}{w_{v}-b^{T}T^{-1}b}$ is uniformly
bounded for all $N\geq4$.

Then \ref{av503} implies that $\mu$ is uniformly bounded for all $N\geq4$ iff
$T^{-1}\Delta_{X}f$ and $T^{-1}b$ are uniformly bounded for all $N\geq4 $.

From MATLAB experiments it is reasonable to believe that $\frac{v}{N}%
T^{-1}\Delta_{X}f=\rho T^{-1}\Delta_{X}f$ is uniformly bounded for all
$N\geq4$, $\rho\geq0$ and all data in $\Omega$.
\end{remark}

\subsection{Applying Cramer's rule}

In this subsection we use Cramer's rule for matrix inversion to express the
right sides of equations \ref{av534} to \ref{av529} in terms of the minors of
$T$. Further, we will also express several of these minors in terms of
determinants of principal sub-matrices of $T$ e.g. \ref{av533}.

Cramer's rule states that%
\begin{equation}
\left(  T^{-1}\right)  _{m,n}=\frac{A_{m,n}}{\left\vert T\right\vert }%
=\frac{\left(  -1\right)  ^{m+n}M_{m,n}}{\left\vert T\right\vert
},\label{av003}%
\end{equation}

where $A=\left(  A_{m,n}\right)  $ is the adjoint, $A_{m,n}$ is the cofactor
of the element $T_{m,n}$, $M_{m,n}=M_{n,m}$ denotes the minor. Further, from
\ref{av011} and Definition \ref{vDef_T_minors},%
\begin{equation}
M_{1,N-1}=%
\begin{vmatrix}
-\frac{1}{2} & 1+\frac{1}{v}\Delta_{1}x^{\left(  2\right)  } & -\frac{1}{2} &
& 0 & 0 & 0\\
0 & -\frac{1}{2} & 1+\frac{1}{v}\Delta_{1}x^{\left(  3\right)  } &  & 0 & 0 &
0\\
0 & 0 & -\frac{1}{2} & \ddots & 0 & 0 & 0\\
&  &  & \ddots &  &  & \\
0 & 0 & 0 &  & -\frac{1}{2} & 1+\frac{1}{v}\Delta_{1}x^{\left(  N-3\right)  }
& -\frac{1}{2}\\
0 & 0 & 0 &  & 0 & -\frac{1}{2} & 1+\frac{1}{v}\Delta_{1}x^{\left(
N-2\right)  }\\
0 & 0 & 0 &  & 0 & 0 & -\frac{1}{2}%
\end{vmatrix}
.\label{av005}%
\end{equation}

so%
\begin{equation}
M_{1,N-1}=M_{N-1,1}=\left(  -\frac{1}{2}\right)  ^{N-2},\quad N\geq
3,\label{av552}%
\end{equation}

and using the notation of Definition \ref{vDef_T_minors},%
\begin{equation}
M_{N-1,N-1}=\left\vert T_{1:N-2}\right\vert ,\quad M_{1,1}=\left\vert
T_{2:N-1}\right\vert .\label{av533}%
\end{equation}

Equations \ref{av534} to \ref{av529} now become for $1\leq k\leq N-1$:%
\begin{align*}
\left(  T^{-1}\beta\right)  _{k}=\left(  T^{-1}\right)  _{k,N-1}-\left(
T^{-1}\right)  _{k,1} &  =\frac{1}{\left\vert T\right\vert }\left(  \left(
-1\right)  ^{k-1}\left(  \left(  -1\right)  ^{N}M_{k,N-1}-M_{k,1}\right)
\right) \\
&  =\left(  -1\right)  ^{k-1}\frac{1}{\left\vert T\right\vert }\left(  \left(
-1\right)  ^{N}M_{k,N-1}-M_{k,1}\right)  ,
\end{align*}

\[
\left(  T^{-1}\Delta_{X}f\right)  _{k}=\sum\limits_{n=1}^{N-1}\frac{\left(
-1\right)  ^{k+n}M_{k,n}}{\left\vert T\right\vert }\Delta_{1}f\left(
x^{\left(  n\right)  }\right)  =\frac{\left(  -1\right)  ^{k}}{\left\vert
T\right\vert }\sum\limits_{n=1}^{N-1}\left(  -1\right)  ^{n}M_{k,n}\Delta
_{1}f\left(  x^{\left(  n\right)  }\right)  ,
\]

using equations \ref{av552} and \ref{av533},%
\begin{align*}
\beta^{T}T^{-1}\beta & =\left(  T^{-1}\right)  _{1,1}-2\left(  T^{-1}\right)
_{1,N-1}+\left(  T^{-1}\right)  _{N-1,N-1}\\
& =\frac{1}{\left\vert T\right\vert }\left(  M_{1,1}-\left(  -1\right)
^{N}2M_{1,N-1}+M_{N-1,N-1}\right) \\
& =\frac{1}{\left\vert T\right\vert }\left(  \left\vert T_{2:N-1}\right\vert
-\left(  -1\right)  ^{N}2\left(  -\frac{1}{2}\right)  ^{N-2}+\left\vert
T_{1:N-2}\right\vert \right) \\
& =\frac{1}{\left\vert T\right\vert }\left(  \left\vert T_{2:N-1}\right\vert
-\left(  \frac{1}{2}\right)  ^{N-3}+\left\vert T_{1:N-2}\right\vert \right)  ,
\end{align*}

and%
\begin{align*}
\beta^{T}T^{-1}\Delta_{X}f  & =\sum\limits_{n=1}^{N-1}\left(  \left(
T^{-1}\right)  _{n,N-1}-\left(  T^{-1}\right)  _{n,1}\right)  \Delta
_{1}f\left(  x^{\left(  n\right)  }\right) \\
& =\frac{1}{\left\vert T\right\vert }\sum\limits_{n=1}^{N-1}\left(  \left(
-1\right)  ^{n+N-1}M_{n,N-1}-\left(  -1\right)  ^{n+1}M_{n,1}\right)
\Delta_{1}f\left(  x^{\left(  n\right)  }\right) \\
& =\frac{1}{\left\vert T\right\vert }\sum\limits_{n=1}^{N-1}\left(  -1\right)
^{n-1}\left(  \left(  -1\right)  ^{N}M_{n,N-1}-M_{n,1}\right)  \Delta
_{1}f\left(  x^{\left(  n\right)  }\right)  .
\end{align*}

which we summarize as: for $1\leq k\leq N-1$,%
\begin{align}
\left(  T^{-1}\beta\right)  _{k}  & =\frac{\left(  -1\right)  ^{k-1}%
}{\left\vert T\right\vert }\left(  \left(  -1\right)  ^{N}M_{k,N-1}%
-M_{k,1}\right)  ,\label{av725}\\
\left(  T^{-1}\Delta_{X}f\right)  _{k}  & =\frac{\left(  -1\right)  ^{k}%
}{\left\vert T\right\vert }\sum\limits_{n=1}^{N-1}\left(  -1\right)
^{n}M_{k,n}\Delta_{1}f\left(  x^{\left(  n\right)  }\right)  ,\label{av758}\\
\beta^{T}T^{-1}\beta & =\frac{1}{\left\vert T\right\vert }\left(  \left\vert
T_{2:N-1}\right\vert -\left(  \frac{1}{2}\right)  ^{N-3}+\left\vert
T_{1:N-2}\right\vert \right)  ,\label{av756}\\
\beta^{T}T^{-1}\Delta_{X}f  & =\frac{1}{\left\vert T\right\vert }%
\sum\limits_{k=1}^{N-1}\left(  -1\right)  ^{k-1}\left(  \left(  -1\right)
^{N}M_{k,N-1}-M_{k,1}\right)  \Delta_{1}f\left(  x^{\left(  k\right)
}\right)  .\label{av757}%
\end{align}

These equations will later be substituted in equations \ref{av008} and
\ref{av007} for $\mu\left(  f\right)  $.

\section{Calculating the minors of $T$\label{Sect_minors}}

Regarding the equations \ref{av725} to \ref{av757} derived at the end of the
last section,\ in this section we will derive the expression \ref{av787} which
expresses all the minors of $T\left(  a\right)  =T\left(  a;1\right)  $ in
terms of the products of one or two \textbf{principal sub-determinants}
$\left\vert T\left(  a_{i:j};1\right)  \right\vert $. These results are
summarized in Subsection \ref{SbSect_minor_summary}.

Note that not all the determinant proofs of this section are complete e.g.
\ref{av782}.

\subsection{Calculating the minors $M_{n,1}\left(  a_{1:m}\right)  $ and
$M_{n,m}\left(  a_{1:m}\right)  $\label{SbSect_minors_Mm1_MmN-1}}

We use the notation of Definition \ref{vDef_T_minors}.

\begin{example}
If $a\in\mathbb{R}^{10}$ then%
\begin{equation}
\left\vert T\left(  a\right)  \right\vert =%
\begin{vmatrix}
1+a_{1} & -\frac{1}{2} & 0 & 0 & 0 & 0 & 0 & 0 & 0 & 0\\
-\frac{1}{2} & 1+a_{2} & -\frac{1}{2} & 0 & 0 & 0 & 0 & 0 & 0 & 0\\
0 & -\frac{1}{2} & 1+a_{3} & -\frac{1}{2} & 0 & 0 & 0 & 0 & 0 & 0\\
0 & 0 & -\frac{1}{2} & 1+a_{4} & -\frac{1}{2} & 0 & 0 & 0 & 0 & 0\\
0 & 0 & 0 & -\frac{1}{2} & 1+a_{5} & -\frac{1}{2} & 0 & 0 & 0 & 0\\
0 & 0 & 0 & 0 & -\frac{1}{2} & 1+a_{6} & -\frac{1}{2} & 0 & 0 & 0\\
0 & 0 & 0 & 0 & 0 & -\frac{1}{2} & 1+a_{7} & -\frac{1}{2} & 0 & 0\\
0 & 0 & 0 & 0 & 0 & 0 & -\frac{1}{2} & 1+a_{8} & -\frac{1}{2} & 0\\
0 & 0 & 0 & 0 & 0 & 0 & 0 & -\frac{1}{2} & 1+a_{9} & -\frac{1}{2}\\
0 & 0 & 0 & 0 & 0 & 0 & 0 & 0 & -\frac{1}{2} & 1+a_{10}%
\end{vmatrix}
,\label{av778}%
\end{equation}

and the minor corresponding to element $\left(  1,7\right)  $ is
\begin{align*}
M_{1,7}\left(  a\right)   & =%
\begin{vmatrix}
-\frac{1}{2} & 1+a_{2} & -\frac{1}{2} & 0 & 0 & 0 & 0 & 0 & 0\\
0 & -\frac{1}{2} & 1+a_{3} & -\frac{1}{2} & 0 & 0 & 0 & 0 & 0\\
0 & 0 & -\frac{1}{2} & 1+a_{4} & -\frac{1}{2} & 0 & 0 & 0 & 0\\
0 & 0 & 0 & -\frac{1}{2} & 1+a_{5} & -\frac{1}{2} & 0 & 0 & 0\\
0 & 0 & 0 & 0 & -\frac{1}{2} & 1+a_{6} & 0 & 0 & 0\\
0 & 0 & 0 & 0 & 0 & -\frac{1}{2} & -\frac{1}{2} & 0 & 0\\
0 & 0 & 0 & 0 & 0 & 0 & 1+a_{8} & -\frac{1}{2} & 0\\
0 & 0 & 0 & 0 & 0 & 0 & -\frac{1}{2} & 1+a_{9} & -\frac{1}{2}\\
0 & 0 & 0 & 0 & 0 & 0 & 0 & -\frac{1}{2} & 1+a_{10}%
\end{vmatrix}
\\
& \\
& =-\frac{1}{2}%
\begin{vmatrix}
-\frac{1}{2} & 1+a_{3} & -\frac{1}{2} & 0 & 0 & 0 & 0 & 0\\
0 & -\frac{1}{2} & 1+a_{4} & -\frac{1}{2} & 0 & 0 & 0 & 0\\
0 & 0 & -\frac{1}{2} & 1+a_{5} & -\frac{1}{2} & 0 & 0 & 0\\
0 & 0 & 0 & -\frac{1}{2} & 1+a_{6} & 0 & 0 & 0\\
0 & 0 & 0 & 0 & -\frac{1}{2} & -\frac{1}{2} & 0 & 0\\
0 & 0 & 0 & 0 & 0 & 1+a_{8} & -\frac{1}{2} & 0\\
0 & 0 & 0 & 0 & 0 & -\frac{1}{2} & 1+a_{9} & -\frac{1}{2}\\
0 & 0 & 0 & 0 & 0 & 0 & -\frac{1}{2} & 1+a_{10}%
\end{vmatrix}
\\
& \\
& =\left(  -\frac{1}{2}\right)  ^{2}%
\begin{vmatrix}
-\frac{1}{2} & 1+a_{4} & -\frac{1}{2} & 0 & 0 & 0 & 0\\
0 & -\frac{1}{2} & 1+a_{5} & -\frac{1}{2} & 0 & 0 & 0\\
0 & 0 & -\frac{1}{2} & 1+a_{6} & 0 & 0 & 0\\
0 & 0 & 0 & -\frac{1}{2} & -\frac{1}{2} & 0 & 0\\
0 & 0 & 0 & 0 & 1+a_{8} & -\frac{1}{2} & 0\\
0 & 0 & 0 & 0 & -\frac{1}{2} & 1+a_{9} & -\frac{1}{2}\\
0 & 0 & 0 & 0 & 0 & -\frac{1}{2} & 1+a_{10}%
\end{vmatrix}
\\
& \\
& =\left(  -\frac{1}{2}\right)  ^{3}%
\begin{vmatrix}
-\frac{1}{2} & 1+a_{5} & -\frac{1}{2} & 0 & 0 & 0\\
0 & -\frac{1}{2} & 1+a_{6} & 0 & 0 & 0\\
0 & 0 & -\frac{1}{2} & -\frac{1}{2} & 0 & 0\\
0 & 0 & 0 & 1+a_{8} & -\frac{1}{2} & 0\\
0 & 0 & 0 & -\frac{1}{2} & 1+a_{9} & -\frac{1}{2}\\
0 & 0 & 0 & 0 & -\frac{1}{2} & 1+a_{10}%
\end{vmatrix}
\\
& \\
& =\left(  -\frac{1}{2}\right)  ^{4}%
\begin{vmatrix}
-\frac{1}{2} & 1+a_{6} & 0 & 0 & 0\\
0 & -\frac{1}{2} & -\frac{1}{2} & 0 & 0\\
0 & 0 & 1+a_{8} & -\frac{1}{2} & 0\\
0 & 0 & -\frac{1}{2} & 1+a_{9} & -\frac{1}{2}\\
0 & 0 & 0 & -\frac{1}{2} & 1+a_{10}%
\end{vmatrix}
\\
& \\
& =\left(  -\frac{1}{2}\right)  ^{5}%
\begin{vmatrix}
-\frac{1}{2} & -\frac{1}{2} & 0 & 0\\
0 & 1+a_{8} & -\frac{1}{2} & 0\\
0 & -\frac{1}{2} & 1+a_{9} & -\frac{1}{2}\\
0 & 0 & -\frac{1}{2} & 1+a_{10}%
\end{vmatrix}
=\left(  -\frac{1}{2}\right)  ^{6}%
\begin{vmatrix}
1+a_{8} & -\frac{1}{2} & 0\\
-\frac{1}{2} & 1+a_{9} & -\frac{1}{2}\\
0 & -\frac{1}{2} & 1+a_{10}%
\end{vmatrix}
=\\
& \\
& =\left(  -\frac{1}{2}\right)  ^{6}\left\vert T\left(  a_{8:10}\right)
\right\vert ,
\end{align*}

i.e.%
\[
M_{1,7}\left(  a_{1:10}\right)  =\left(  -\frac{1}{2}\right)  ^{6}\left\vert
T\left(  a_{8:10}\right)  \right\vert ,
\]

\end{example}

and in general we assert \textbf{without proof} that%
\begin{equation}
M_{n,1}\left(  a_{1:m}\right)  =M_{1,n}\left(  a_{1:m}\right)  =\left(
-\frac{1}{2}\right)  ^{n-1}\left\vert T\left(  a_{n+1:m}\right)  \right\vert
,\quad2\leq n\leq m-1.\label{av746}%
\end{equation}

When $a=\frac{1}{v}\Delta_{1}X$, in the notation of Definition
\ref{vDef_T_minors}:%
\begin{equation}
M_{n,1}=M_{1,n}=\left(  -\frac{1}{2}\right)  ^{n-1}\left\vert T\left(
a_{n+1:N-1}\right)  \right\vert ,\quad2\leq n\leq N-2.\label{av015}%
\end{equation}

Also%
\begin{align*}
M_{10,7}\left(  a_{1:10}\right)   & =%
\begin{vmatrix}
1+a_{1} & -\frac{1}{2} & 0 & 0 & 0 & 0 & 0 & 0 & 0\\
-\frac{1}{2} & 1+a_{2} & -\frac{1}{2} & 0 & 0 & 0 & 0 & 0 & 0\\
0 & -\frac{1}{2} & 1+a_{3} & -\frac{1}{2} & 0 & 0 & 0 & 0 & 0\\
0 & 0 & -\frac{1}{2} & 1+a_{4} & -\frac{1}{2} & 0 & 0 & 0 & 0\\
0 & 0 & 0 & -\frac{1}{2} & 1+a_{5} & -\frac{1}{2} & 0 & 0 & 0\\
0 & 0 & 0 & 0 & -\frac{1}{2} & 1+a_{6} & 0 & 0 & 0\\
0 & 0 & 0 & 0 & 0 & -\frac{1}{2} & -\frac{1}{2} & 0 & 0\\
0 & 0 & 0 & 0 & 0 & 0 & 1+a_{8} & -\frac{1}{2} & 0\\
0 & 0 & 0 & 0 & 0 & 0 & -\frac{1}{2} & 1+a_{9} & -\frac{1}{2}%
\end{vmatrix}
\\
& \\
& =-\frac{1}{2}%
\begin{vmatrix}
1+a_{1} & -\frac{1}{2} & 0 & 0 & 0 & 0 & 0 & 0\\
-\frac{1}{2} & 1+a_{2} & -\frac{1}{2} & 0 & 0 & 0 & 0 & 0\\
0 & -\frac{1}{2} & 1+a_{3} & -\frac{1}{2} & 0 & 0 & 0 & 0\\
0 & 0 & -\frac{1}{2} & 1+a_{4} & -\frac{1}{2} & 0 & 0 & 0\\
0 & 0 & 0 & -\frac{1}{2} & 1+a_{5} & -\frac{1}{2} & 0 & 0\\
0 & 0 & 0 & 0 & -\frac{1}{2} & 1+a_{6} & 0 & 0\\
0 & 0 & 0 & 0 & 0 & -\frac{1}{2} & -\frac{1}{2} & 0\\
0 & 0 & 0 & 0 & 0 & 0 & 1+a_{8} & -\frac{1}{2}%
\end{vmatrix}
\\
& \\
& =\left(  -\frac{1}{2}\right)  ^{2}%
\begin{vmatrix}
1+a_{1} & -\frac{1}{2} & 0 & 0 & 0 & 0 & 0\\
-\frac{1}{2} & 1+a_{2} & -\frac{1}{2} & 0 & 0 & 0 & 0\\
0 & -\frac{1}{2} & 1+a_{3} & -\frac{1}{2} & 0 & 0 & 0\\
0 & 0 & -\frac{1}{2} & 1+a_{4} & -\frac{1}{2} & 0 & 0\\
0 & 0 & 0 & -\frac{1}{2} & 1+a_{5} & -\frac{1}{2} & 0\\
0 & 0 & 0 & 0 & -\frac{1}{2} & 1+a_{6} & 0\\
0 & 0 & 0 & 0 & 0 & -\frac{1}{2} & -\frac{1}{2}%
\end{vmatrix}
\\
& \\
& =\left(  -\frac{1}{2}\right)  ^{3}%
\begin{vmatrix}
1+a_{1} & -\frac{1}{2} & 0 & 0 & 0 & 0\\
-\frac{1}{2} & 1+a_{2} & -\frac{1}{2} & 0 & 0 & 0\\
0 & -\frac{1}{2} & 1+a_{3} & -\frac{1}{2} & 0 & 0\\
0 & 0 & -\frac{1}{2} & 1+a_{4} & -\frac{1}{2} & 0\\
0 & 0 & 0 & -\frac{1}{2} & 1+a_{5} & -\frac{1}{2}\\
0 & 0 & 0 & 0 & -\frac{1}{2} & 1+a_{6}%
\end{vmatrix}
\\
& =\left(  -\frac{1}{2}\right)  ^{3}\left\vert T_{1,6}\right\vert ,
\end{align*}

so that%
\[
M_{11,7}\left(  a_{1:10}\right)  =\left(  -\frac{1}{2}\right)  ^{3}\left\vert
T\left(  a_{1:5}\right)  \right\vert ,
\]

and in general we assert without proof that%
\begin{equation}
M_{n,m}\left(  a_{1:m}\right)  =M_{m,n}\left(  a_{1:m}\right)  =\left(
-\frac{1}{2}\right)  ^{m-n}\left\vert T\left(  a_{1:n-1}\right)  \right\vert
,\quad2\leq n\leq m.\label{av747}%
\end{equation}

When $a=\frac{1}{v}\Delta_{1}X$, in the notation of Definition
\ref{vDef_T_minors}:%
\begin{equation}
M_{n,N-1}=M_{N-1,n}=\left(  -\frac{1}{2}\right)  ^{N-1-n}\left\vert T\left(
a_{1:n-1}\right)  \right\vert ,\quad2\leq n\leq N-1.\label{av017}%
\end{equation}

\subsection{Calculation of the minors $M_{n,n}$}

\begin{example}
From \ref{av778},%
\begin{align*}
\left\vert M_{4,4}\left(  a_{1:9}\right)  \right\vert  & =%
\begin{vmatrix}
1+a_{1} & -\frac{1}{2} & 0 & 0 & 0 & 0 & 0 & 0\\
-\frac{1}{2} & 1+a_{2} & -\frac{1}{2} & 0 & 0 & 0 & 0 & 0\\
0 & -\frac{1}{2} & 1+a_{3} & 0 & 0 & 0 & 0 & 0\\
0 & 0 & 0 & 1+a_{5} & -\frac{1}{2} & 0 & 0 & 0\\
0 & 0 & 0 & -\frac{1}{2} & 1+a_{6} & -\frac{1}{2} & 0 & 0\\
0 & 0 & 0 & 0 & -\frac{1}{2} & 1+a_{7} & -\frac{1}{2} & 0\\
0 & 0 & 0 & 0 & 0 & -\frac{1}{2} & 1+a_{8} & -\frac{1}{2}\\
0 & 0 & 0 & 0 & 0 & 0 & -\frac{1}{2} & 1+a_{9}%
\end{vmatrix}
\\
& \\
& =%
\begin{vmatrix}
1+a_{1} & -\frac{1}{2} & 0\\
-\frac{1}{2} & 1+a_{2} & -\frac{1}{2}\\
0 & -\frac{1}{2} & 1+a_{3}%
\end{vmatrix}%
\begin{vmatrix}
1+a_{5} & -\frac{1}{2} & 0 & 0 & 0\\
-\frac{1}{2} & 1+a_{6} & -\frac{1}{2} & 0 & 0\\
0 & -\frac{1}{2} & 1+a_{7} & -\frac{1}{2} & 0\\
0 & 0 & -\frac{1}{2} & 1+a_{8} & -\frac{1}{2}\\
0 & 0 & 0 & -\frac{1}{2} & 1+a_{9}%
\end{vmatrix}
\\
& =\left\vert T\left(  a_{1:3}\right)  \right\vert \text{ }\left\vert T\left(
a_{5:9}\right)  \right\vert ,
\end{align*}

\end{example}

so that we hypothesize, \textbf{but do not prove}, that if $m\geq3$ then%
\begin{equation}
M_{n,n}\left(  a_{1:m}\right)  =%
\begin{array}
[c]{ll}%
\left\vert T_{2:m}\left(  a_{1:m}\right)  \right\vert , & n=1,\\
\left\vert T_{1:n-1}\left(  a_{1:m}\right)  \right\vert \left\vert
T_{n+1:m}\left(  a_{1:m}\right)  \right\vert , & 2\leq n\leq m-1,\\
\left\vert T_{1:m-1}\left(  a_{1:m}\right)  \right\vert , & n=m.
\end{array}
\label{av782}%
\end{equation}

\subsection{Calculation of the minors $M_{n,n+1}$}

Also%
\begin{align*}
M_{5,4}\left(  a_{1:9}\right)   & =%
\begin{vmatrix}
1+a_{1} & -\frac{1}{2} & 0 & 0 & 0 & 0 & 0 & 0\\
-\frac{1}{2} & 1+a_{2} & -\frac{1}{2} & 0 & 0 & 0 & 0 & 0\\
0 & -\frac{1}{2} & 1+a_{3} & 0 & 0 & 0 & 0 & 0\\
0 & 0 & -\frac{1}{2} & \frame{$-\frac{1}{2}$} & 0 & 0 & 0 & 0\\
0 & 0 & 0 & \frame{$-\frac{1}{2}$} & 1+a_{6} & -\frac{1}{2} & 0 & 0\\
0 & 0 & 0 & 0 & -\frac{1}{2} & 1+a_{7} & -\frac{1}{2} & 0\\
0 & 0 & 0 & 0 & 0 & -\frac{1}{2} & 1+a_{8} & -\frac{1}{2}\\
0 & 0 & 0 & 0 & 0 & 0 & -\frac{1}{2} & 1+a_{9}%
\end{vmatrix}
\\
& \\
& =-\frac{1}{2}%
\begin{vmatrix}
1+a_{1} & -\frac{1}{2} & 0 & 0 & 0 & 0 & 0 & 0\\
-\frac{1}{2} & 1+a_{2} & -\frac{1}{2} & 0 & 0 & 0 & 0 & 0\\
0 & -\frac{1}{2} & 1+a_{3} & 0 & 0 & 0 & 0 & 0\\
0 & 0 & -\frac{1}{2} & \frame{$1$} & 0 & 0 & 0 & 0\\
0 & 0 & 0 & \frame{$1$} & 1+a_{6} & -\frac{1}{2} & 0 & 0\\
0 & 0 & 0 & 0 & -\frac{1}{2} & 1+a_{7} & -\frac{1}{2} & 0\\
0 & 0 & 0 & 0 & 0 & -\frac{1}{2} & 1+a_{8} & -\frac{1}{2}\\
0 & 0 & 0 & 0 & 0 & 0 & -\frac{1}{2} & 1+a_{9}%
\end{vmatrix}
\\
& \\
& =-\frac{1}{2}%
\begin{vmatrix}
1+a_{1} & -\frac{1}{2} & 0 & 0 & 0 & 0 & 0\\
-\frac{1}{2} & 1+a_{2} & -\frac{1}{2} & 0 & 0 & 0 & 0\\
0 & -\frac{1}{2} & 1+a_{3} & 0 & 0 & 0 & 0\\
0 & 0 & 0 & 1+a_{6} & -\frac{1}{2} & 0 & 0\\
0 & 0 & 0 & -\frac{1}{2} & 1+a_{7} & -\frac{1}{2} & 0\\
0 & 0 & 0 & 0 & -\frac{1}{2} & 1+a_{8} & -\frac{1}{2}\\
0 & 0 & 0 & 0 & 0 & -\frac{1}{2} & 1+a_{9}%
\end{vmatrix}
+\\
& \\
& \qquad\qquad+\frac{1}{2}%
\begin{vmatrix}
1+a_{1} & -\frac{1}{2} & 0 & 0 & 0 & 0 & 0\\
-\frac{1}{2} & 1+a_{2} & -\frac{1}{2} & 0 & 0 & 0 & 0\\
0 & -\frac{1}{2} & 1+a_{3} & 0 & 0 & 0 & 0\\
0 & 0 & -\frac{1}{2} & 0 & 0 & 0 & 0\\
0 & 0 & 0 & \frame{$-\frac{1}{2}$} & 1+a_{7} & -\frac{1}{2} & 0\\
0 & 0 & 0 & 0 & -\frac{1}{2} & 1+a_{8} & -\frac{1}{2}\\
0 & 0 & 0 & 0 & 0 & -\frac{1}{2} & 1+a_{9}%
\end{vmatrix}
\\
& \\
& =-\frac{1}{2}%
\begin{vmatrix}
1+a_{1} & -\frac{1}{2} & 0 & 0 & 0 & 0 & 0\\
-\frac{1}{2} & 1+a_{2} & -\frac{1}{2} & 0 & 0 & 0 & 0\\
0 & -\frac{1}{2} & 1+a_{3} & 0 & 0 & 0 & 0\\
0 & 0 & 0 & 1+a_{6} & -\frac{1}{2} & 0 & 0\\
0 & 0 & 0 & -\frac{1}{2} & 1+a_{7} & -\frac{1}{2} & 0\\
0 & 0 & 0 & 0 & -\frac{1}{2} & 1+a_{8} & -\frac{1}{2}\\
0 & 0 & 0 & 0 & 0 & -\frac{1}{2} & 1+a_{9}%
\end{vmatrix}
-\\
& \\
& \qquad\qquad-\frac{1}{2^{2}}%
\begin{vmatrix}
1+a_{1} & -\frac{1}{2} & 0 & 0 & 0 & 0\\
-\frac{1}{2} & 1+a_{2} & -\frac{1}{2} & 0 & 0 & 0\\
0 & -\frac{1}{2} & 1+a_{3} & 0 & 0 & 0\\
0 & 0 & -\frac{1}{2} & 0 & 0 & 0\\
0 & 0 & 0 & -\frac{1}{2} & 1+a_{8} & -\frac{1}{2}\\
0 & 0 & 0 & 0 & -\frac{1}{2} & 1+a_{9}%
\end{vmatrix}
.
\end{align*}

Regarding the second determinant:%
\begin{align*}%
\begin{vmatrix}
1+a_{1} & -\frac{1}{2} & 0 & 0 & 0 & 0\\
-\frac{1}{2} & 1+a_{2} & -\frac{1}{2} & 0 & 0 & 0\\
0 & -\frac{1}{2} & 1+a_{3} & 0 & 0 & 0\\
0 & 0 & -\frac{1}{2} & 0 & 0 & 0\\
0 & 0 & 0 & \frame{$-\frac{1}{2}$} & 1+a_{8} & -\frac{1}{2}\\
0 & 0 & 0 & 0 & -\frac{1}{2} & 1+a_{9}%
\end{vmatrix}
& =\frac{1}{2}%
\begin{vmatrix}
1+a_{1} & -\frac{1}{2} & 0 & 0 & 0\\
-\frac{1}{2} & 1+a_{2} & -\frac{1}{2} & 0 & 0\\
0 & -\frac{1}{2} & 1+a_{3} & 0 & 0\\
0 & 0 & -\frac{1}{2} & 0 & 0\\
0 & 0 & 0 & \frame{$-\frac{1}{2}$} & 1+a_{9}%
\end{vmatrix}
\\
& \\
& =-\frac{1}{2^{2}}%
\begin{vmatrix}
1+a_{1} & -\frac{1}{2} & 0 & 0\\
-\frac{1}{2} & 1+a_{2} & -\frac{1}{2} & 0\\
0 & -\frac{1}{2} & 1+a_{3} & 0\\
0 & 0 & -\frac{1}{2} & 0
\end{vmatrix}
\\
& =0,
\end{align*}

so that%
\[
M_{5,4}\left(  a_{1:9}\right)  =-\frac{1}{2}%
\begin{vmatrix}
T\left(  a_{1:3}\right)  & O_{3,4}\\
O_{3,4}^{T} & T\left(  a_{6:9}\right)
\end{vmatrix}
=-\frac{1}{2}\left\vert T\left(  a_{1:3}\right)  \right\vert \left\vert
T\left(  a_{6:9}\right)  \right\vert =M_{4,5}\left(  a_{1:9}\right)  .
\]

We hypothesize, \textbf{but do not prove}, that%
\begin{equation}
M_{n,n+1}\left(  a_{1:m}\right)  =M_{n+1,n}\left(  a\right)  =-\frac{1}%
{2}\left\vert T\left(  a_{1:n-1}\right)  \right\vert \left\vert T\left(
a_{n+1:m}\right)  \right\vert ,\quad2\leq n\leq m-2.\label{av776}%
\end{equation}

\begin{align*}
M_{1,2}\left(  a\right)   & =%
\begin{vmatrix}
-\frac{1}{2} & -\frac{1}{2} & 0 & 0 & 0 & 0 & 0 & 0\\
0 & 1+a_{3} & -\frac{1}{2} & 0 & 0 & 0 & 0 & 0\\
0 & -\frac{1}{2} & 1+a_{4} & -\frac{1}{2} & 0 & 0 & 0 & 0\\
0 & 0 & -\frac{1}{2} & 1+a_{5} & -\frac{1}{2} & 0 & 0 & 0\\
0 & 0 & 0 & -\frac{1}{2} & 1+a_{6} & -\frac{1}{2} & 0 & 0\\
0 & 0 & 0 & 0 & -\frac{1}{2} & 1+a_{7} & -\frac{1}{2} & 0\\
0 & 0 & 0 & 0 & 0 & -\frac{1}{2} & 1+a_{8} & -\frac{1}{2}\\
0 & 0 & 0 & 0 & 0 & 0 & -\frac{1}{2} & 1+a_{9}%
\end{vmatrix}
\\
& \\
& \\
& =-\frac{1}{2}%
\begin{vmatrix}
1+a_{3} & -\frac{1}{2} & 0 & 0 & 0 & 0 & 0\\
-\frac{1}{2} & 1+a_{4} & -\frac{1}{2} & 0 & 0 & 0 & 0\\
0 & -\frac{1}{2} & 1+a_{5} & -\frac{1}{2} & 0 & 0 & 0\\
0 & 0 & -\frac{1}{2} & 1+a_{6} & -\frac{1}{2} & 0 & 0\\
0 & 0 & 0 & -\frac{1}{2} & 1+a_{7} & -\frac{1}{2} & 0\\
0 & 0 & 0 & 0 & -\frac{1}{2} & 1+a_{8} & -\frac{1}{2}\\
0 & 0 & 0 & 0 & 0 & -\frac{1}{2} & 1+a_{9}%
\end{vmatrix}
\\
& \\
& =-\frac{1}{2}T\left(  a_{3:9}\right)  ,
\end{align*}

so%
\begin{equation}
M_{2,1}\left(  a_{1:m}\right)  =M_{1,2}\left(  a_{1:m}\right)  =-\frac{1}%
{2}\left\vert T\left(  a_{3:m}\right)  \right\vert ,\quad m\geq3.\label{av777}%
\end{equation}

Further%
\begin{align*}
M_{8,9}\left(  a_{1:9}\right)   & =%
\begin{vmatrix}
1+a_{1} & -\frac{1}{2} & 0 & 0 & 0 & 0 & 0 & 0\\
-\frac{1}{2} & 1+a_{2} & -\frac{1}{2} & 0 & 0 & 0 & 0 & 0\\
0 & -\frac{1}{2} & 1+a_{3} & -\frac{1}{2} & 0 & 0 & 0 & 0\\
0 & 0 & -\frac{1}{2} & 1+a_{4} & -\frac{1}{2} & 0 & 0 & 0\\
0 & 0 & 0 & -\frac{1}{2} & 1+a_{5} & -\frac{1}{2} & 0 & 0\\
0 & 0 & 0 & 0 & -\frac{1}{2} & 1+a_{6} & -\frac{1}{2} & 0\\
0 & 0 & 0 & 0 & 0 & -\frac{1}{2} & 1+a_{7} & 0\\
0 & 0 & 0 & 0 & 0 & 0 & -\frac{1}{2} & -\frac{1}{2}%
\end{vmatrix}
\\
& \\
& \\
& =-\frac{1}{2}%
\begin{vmatrix}
1+a_{1} & -\frac{1}{2} & 0 & 0 & 0 & 0 & 0\\
-\frac{1}{2} & 1+a_{2} & -\frac{1}{2} & 0 & 0 & 0 & 0\\
0 & -\frac{1}{2} & 1+a_{3} & -\frac{1}{2} & 0 & 0 & 0\\
0 & 0 & -\frac{1}{2} & 1+a_{4} & -\frac{1}{2} & 0 & 0\\
0 & 0 & 0 & -\frac{1}{2} & 1+a_{5} & -\frac{1}{2} & 0\\
0 & 0 & 0 & 0 & -\frac{1}{2} & 1+a_{6} & -\frac{1}{2}\\
0 & 0 & 0 & 0 & 0 & -\frac{1}{2} & 1+a_{7}%
\end{vmatrix}
\\
& \\
& \\
& =-\frac{1}{2}T\left(  a_{1:7}\right)  ,
\end{align*}

and we hypothesize, \textbf{but do not prove}, that%
\begin{equation}
M_{m-1,m}\left(  a_{1:m}\right)  =M_{m,m-1}\left(  a\right)  =-\frac{1}%
{2}\left\vert T\left(  a_{1:m-2}\right)  \right\vert ,\quad m\geq
3.\label{av779}%
\end{equation}

Thus the equations \ref{av777}, \ref{av776} and \ref{av779} combine to give
for $m\geq4$,%
\begin{equation}
M_{n,n+1}\left(  a_{1:m}\right)  =\left\{
\begin{array}
[c]{ll}%
-\frac{1}{2}\left\vert T\left(  a_{3:m}\right)  \right\vert , & n=1,\\
-\frac{1}{2}\left\vert T\left(  a_{1:n-1}\right)  \right\vert \left\vert
T\left(  a_{n+2:m}\right)  \right\vert , & 2\leq n\leq m-2,\\
-\frac{1}{2}\left\vert T\left(  a_{1:m-2}\right)  \right\vert , & n=m-1,
\end{array}
\right. \label{av783}%
\end{equation}

and for the special case $m=3$,%
\begin{equation}
M_{n,n+1}\left(  a_{1:3}\right)  =\left\{
\begin{array}
[c]{ll}%
-\frac{1}{2}\left\vert T\left(  a_{3}\right)  \right\vert , & n=1,\\
-\frac{1}{2}\left\vert T\left(  a_{1}\right)  \right\vert , & n=2,
\end{array}
\right.  =\left\{
\begin{array}
[c]{ll}%
-\frac{1}{2}\left(  1+a_{3}\right)  , & n=1,\\
-\frac{1}{2}\left(  1+a_{1}\right)  , & n=2.
\end{array}
\right. \label{av056}%
\end{equation}

\subsection{Calculation of $M_{j,n}$, $j<n$.}

\begin{example}
From \ref{av778},%
\begin{align*}
M_{3,6}\left(  a_{1:9}\right)   & =%
\begin{vmatrix}
1+a_{1} & -\frac{1}{2} & 0 & 0 & 0 & 0 & 0 & 0\\
-\frac{1}{2} & 1+a_{2} & \frame{$-\frac{1}{2}$} & 0 & 0 & 0 & 0 & 0\\
0 & 0 & \frame{$-\frac{1}{2}$} & 1+a_{4} & -\frac{1}{2} & 0 & 0 & 0\\
0 & 0 & 0 & -\frac{1}{2} & 1+a_{5} & 0 & 0 & 0\\
0 & 0 & 0 & 0 & -\frac{1}{2} & -\frac{1}{2} & 0 & 0\\
0 & 0 & 0 & 0 & 0 & 1+a_{7} & -\frac{1}{2} & 0\\
0 & 0 & 0 & 0 & 0 & -\frac{1}{2} & 1+a_{8} & -\frac{1}{2}\\
0 & 0 & 0 & 0 & 0 & 0 & -\frac{1}{2} & 1+a_{9}%
\end{vmatrix}
\\
& \\
& =-\frac{1}{2}%
\begin{vmatrix}
1+a_{1} & -\frac{1}{2} & 0 & 0 & 0 & 0 & 0 & 0\\
-\frac{1}{2} & 1+a_{2} & \frame{$1$} & 0 & 0 & 0 & 0 & 0\\
0 & 0 & \frame{$1$} & 1+a_{4} & -\frac{1}{2} & 0 & 0 & 0\\
0 & 0 & 0 & -\frac{1}{2} & 1+a_{5} & 0 & 0 & 0\\
0 & 0 & 0 & 0 & -\frac{1}{2} & -\frac{1}{2} & 0 & 0\\
0 & 0 & 0 & 0 & 0 & 1+a_{7} & -\frac{1}{2} & 0\\
0 & 0 & 0 & 0 & 0 & -\frac{1}{2} & 1+a_{8} & -\frac{1}{2}\\
0 & 0 & 0 & 0 & 0 & 0 & -\frac{1}{2} & 1+a_{9}%
\end{vmatrix}
\\
& \\
& =\frac{1}{2}%
\begin{vmatrix}
\frame{$1+a_1$} & -\frac{1}{2} & 0 & 0 & 0 & 0 & 0\\
0 & 0 & 1+a_{4} & -\frac{1}{2} & 0 & 0 & 0\\
0 & 0 & -\frac{1}{2} & 1+a_{5} & 0 & 0 & 0\\
0 & 0 & 0 & -\frac{1}{2} & -\frac{1}{2} & 0 & 0\\
0 & 0 & 0 & 0 & 1+a_{7} & -\frac{1}{2} & 0\\
0 & 0 & 0 & 0 & -\frac{1}{2} & 1+a_{8} & -\frac{1}{2}\\
0 & 0 & 0 & 0 & 0 & -\frac{1}{2} & 1+a_{9}%
\end{vmatrix}
-\\
& \\
& \qquad-\frac{1}{2}%
\begin{vmatrix}
1+a_{1} & -\frac{1}{2} & 0 & 0 & 0 & 0 & 0\\
-\frac{1}{2} & 1+a_{2} & 0 & 0 & 0 & 0 & 0\\
0 & 0 & -\frac{1}{2} & 1+a_{5} & 0 & 0 & 0\\
0 & 0 & 0 & -\frac{1}{2} & -\frac{1}{2} & 0 & 0\\
0 & 0 & 0 & 0 & 1+a_{7} & -\frac{1}{2} & 0\\
0 & 0 & 0 & 0 & -\frac{1}{2} & 1+a_{8} & -\frac{1}{2}\\
0 & 0 & 0 & 0 & 0 & -\frac{1}{2} & 1+a_{9}%
\end{vmatrix}
\\
& \\
& =0-\frac{1}{2}%
\begin{vmatrix}
1+a_{1} & -\frac{1}{2} & 0 & 0 & 0 & 0 & 0\\
-\frac{1}{2} & 1+a_{2} & 0 & 0 & 0 & 0 & 0\\
0 & 0 & \fbox{$-\frac{1}{2}$} & 1+a_{5} & 0 & 0 & 0\\
0 & 0 & 0 & -\frac{1}{2} & -\frac{1}{2} & 0 & 0\\
0 & 0 & 0 & 0 & 1+a_{7} & -\frac{1}{2} & 0\\
0 & 0 & 0 & 0 & -\frac{1}{2} & 1+a_{8} & -\frac{1}{2}\\
0 & 0 & 0 & 0 & 0 & -\frac{1}{2} & 1+a_{9}%
\end{vmatrix}
\\
& \\
& =\frac{1}{2^{2}}%
\begin{vmatrix}
1+a_{1} & -\frac{1}{2} & 0 & 0 & 0 & 0 & 0\\
-\frac{1}{2} & 1+a_{2} & 0 & 0 & 0 & 0 & 0\\
0 & 0 & 1 & 1+a_{5} & 0 & 0 & 0\\
0 & 0 & 0 & -\frac{1}{2} & -\frac{1}{2} & 0 & 0\\
0 & 0 & 0 & 0 & 1+a_{7} & -\frac{1}{2} & 0\\
0 & 0 & 0 & 0 & -\frac{1}{2} & 1+a_{8} & -\frac{1}{2}\\
0 & 0 & 0 & 0 & 0 & -\frac{1}{2} & 1+a_{9}%
\end{vmatrix}
,
\end{align*}

but%
\begin{align*}
&
\begin{vmatrix}
1+a_{1} & -\frac{1}{2} & 0 & 0 & 0 & 0 & 0\\
-\frac{1}{2} & 1+a_{2} & 0 & 0 & 0 & 0 & 0\\
0 & 0 & \frame{$1$} & 1+a_{5} & 0 & 0 & 0\\
0 & 0 & 0 & -\frac{1}{2} & -\frac{1}{2} & 0 & 0\\
0 & 0 & 0 & 0 & 1+a_{7} & -\frac{1}{2} & 0\\
0 & 0 & 0 & 0 & -\frac{1}{2} & 1+a_{8} & -\frac{1}{2}\\
0 & 0 & 0 & 0 & 0 & -\frac{1}{2} & 1+a_{9}%
\end{vmatrix}
\\
&  =%
\begin{vmatrix}
1+a_{1} & -\frac{1}{2} & 0 & 0 & 0 & 0\\
-\frac{1}{2} & 1+a_{2} & 0 & 0 & 0 & 0\\
0 & 0 & -\frac{1}{2} & -\frac{1}{2} & 0 & 0\\
0 & 0 & 0 & 1+a_{7} & -\frac{1}{2} & 0\\
0 & 0 & 0 & -\frac{1}{2} & 1+a_{8} & -\frac{1}{2}\\
0 & 0 & 0 & 0 & -\frac{1}{2} & 1+a_{9}%
\end{vmatrix}
\\
\operatorname{col}4 &  =\operatorname{col}4-\operatorname{col}3\Longrightarrow
\\
&  =%
\begin{vmatrix}
1+a_{1} & -\frac{1}{2} & 0 & 0 & 0 & 0\\
-\frac{1}{2} & 1+a_{2} & 0 & 0 & 0 & 0\\
0 & 0 & \fbox{$-\frac{1}{2}$} & 0 & 0 & 0\\
0 & 0 & 0 & 1+a_{7} & -\frac{1}{2} & 0\\
0 & 0 & 0 & -\frac{1}{2} & 1+a_{8} & -\frac{1}{2}\\
0 & 0 & 0 & 0 & -\frac{1}{2} & 1+a_{9}%
\end{vmatrix}
\\
& \\
&  =-\frac{1}{2}%
\begin{vmatrix}
1+a_{1} & -\frac{1}{2} & 0 & 0 & 0 & 0\\
-\frac{1}{2} & 1+a_{2} & 0 & 0 & 0 & 0\\
0 & 0 & \frame{$1$} & 0 & 0 & 0\\
0 & 0 & 0 & 1+a_{7} & -\frac{1}{2} & 0\\
0 & 0 & 0 & -\frac{1}{2} & 1+a_{8} & -\frac{1}{2}\\
0 & 0 & 0 & 0 & -\frac{1}{2} & 1+a_{9}%
\end{vmatrix}
\\
& \\
&  =-\frac{1}{2}%
\begin{vmatrix}
1+a_{1} & -\frac{1}{2} & 0 & 0 & 0\\
-\frac{1}{2} & 1+a_{2} & 0 & 0 & 0\\
0 & 0 & 1+a_{7} & -\frac{1}{2} & 0\\
0 & 0 & -\frac{1}{2} & 1+a_{8} & -\frac{1}{2}\\
0 & 0 & 0 & -\frac{1}{2} & 1+a_{9}%
\end{vmatrix}
\\
&  =-\frac{1}{2}%
\begin{vmatrix}
T\left(  a_{1:2}\right)  & O\\
O & T\left(  a_{7:9}\right)
\end{vmatrix}
\\
&  =-\frac{1}{2}\left\vert T\left(  a_{1:2}\right)  \right\vert \left\vert
T\left(  a_{7:9}\right)  \right\vert ,
\end{align*}

so that%
\[
M_{3,6}\left(  a_{1:9}\right)  =-\frac{1}{2^{3}}\left\vert T\left(
a_{1:2}\right)  \right\vert \left\vert T\left(  a_{7:9}\right)  \right\vert .
\]

\end{example}

This suggests the hypothesis%
\[
M_{j,n}\left(  a_{1:m}\right)  =\left(  -\frac{1}{2}\right)  ^{n-j}\left\vert
T\left(  a_{1:j-1}\right)  \right\vert \left\vert T\left(  a_{n+1:m}\right)
\right\vert ,\quad\left\{
\begin{array}
[c]{l}%
2\leq j<n-1,\\
n\leq m-1,\text{ }m\geq5,
\end{array}
\right.
\]

which we combine with \ref{av783} to give%
\begin{equation}
M_{j,n}\left(  a_{1:m}\right)  =\left(  -\frac{1}{2}\right)  ^{n-j}\left\vert
T\left(  a_{1:j-1}\right)  \right\vert \left\vert T\left(  a_{n+1:m}\right)
\right\vert ,\quad\left\{
\begin{array}
[c]{l}%
2\leq j<n\leq m-1,\\
m\geq3.
\end{array}
\right. \label{av781}%
\end{equation}

\begin{example}
Another example:
\end{example}

\begin{example}
If $a=a_{1:11}$ then (boxes indicate elements involved)%
\begin{align*}
M_{3,9}\left(  a\right)   &  =%
\begin{vmatrix}
1+a_{1} & -\frac{1}{2} & 0 & 0 & 0 & 0 & 0 & 0 & 0 & 0\\
-\frac{1}{2} & 1+a_{2} & \fbox{$-\frac{1}{2}$} & 0 & 0 & 0 & 0 & 0 & 0 & 0\\
0 & 0 & \fbox{$-\frac{1}{2}$} & 1+a_{4} & -\frac{1}{2} & 0 & 0 & 0 & 0 & 0\\
0 & 0 & 0 & -\frac{1}{2} & 1+a_{5} & -\frac{1}{2} & 0 & 0 & 0 & 0\\
0 & 0 & 0 & 0 & -\frac{1}{2} & 1+a_{6} & -\frac{1}{2} & 0 & 0 & 0\\
0 & 0 & 0 & 0 & 0 & -\frac{1}{2} & 1+a_{7} & -\frac{1}{2} & 0 & 0\\
0 & 0 & 0 & 0 & 0 & 0 & -\frac{1}{2} & 1+a_{8} & 0 & 0\\
0 & 0 & 0 & 0 & 0 & 0 & 0 & \fbox{$-\frac{1}{2}$} & \fbox{$-\frac{1}{2}$} &
0\\
0 & 0 & 0 & 0 & 0 & 0 & 0 & 0 & 1+a_{10} & -\frac{1}{2}\\
0 & 0 & 0 & 0 & 0 & 0 & 0 & 0 & -\frac{1}{2} & 1+a_{11}%
\end{vmatrix}
\\
&  =\frac{1}{2^{2}}%
\begin{vmatrix}
1+a_{1} & -\frac{1}{2} & 0 & 0 & 0 & 0 & 0 & 0 & 0 & 0\\
-\frac{1}{2} & 1+a_{2} & \fbox{$1$} & 0 & 0 & 0 & 0 & 0 & 0 & 0\\
0 & 0 & \fbox{$1$} & 1+a_{4} & -\frac{1}{2} & 0 & 0 & 0 & 0 & 0\\
0 & 0 & 0 & -\frac{1}{2} & 1+a_{5} & -\frac{1}{2} & 0 & 0 & 0 & 0\\
0 & 0 & 0 & 0 & -\frac{1}{2} & 1+a_{6} & -\frac{1}{2} & 0 & 0 & 0\\
0 & 0 & 0 & 0 & 0 & -\frac{1}{2} & 1+a_{7} & -\frac{1}{2} & 0 & 0\\
0 & 0 & 0 & 0 & 0 & 0 & -\frac{1}{2} & 1+a_{8} & 0 & 0\\
0 & 0 & 0 & 0 & 0 & 0 & 0 & 1 & 1 & 0\\
0 & 0 & 0 & 0 & 0 & 0 & 0 & 0 & 1+a_{10} & -\frac{1}{2}\\
0 & 0 & 0 & 0 & 0 & 0 & 0 & 0 & -\frac{1}{2} & 1+a_{11}%
\end{vmatrix}
\\
&  =\frac{1}{2^{2}}%
\begin{vmatrix}
1+a_{1} & -\frac{1}{2} & 0 & 0 & 0 & 0 & 0 & 0 & 0\\
0 & 0 & 1+a_{4} & -\frac{1}{2} & 0 & 0 & 0 & 0 & 0\\
0 & 0 & -\frac{1}{2} & 1+a_{5} & -\frac{1}{2} & 0 & 0 & 0 & 0\\
0 & 0 & 0 & -\frac{1}{2} & 1+a_{6} & -\frac{1}{2} & 0 & 0 & 0\\
0 & 0 & 0 & 0 & -\frac{1}{2} & 1+a_{7} & -\frac{1}{2} & 0 & 0\\
0 & 0 & 0 & 0 & 0 & -\frac{1}{2} & 1+a_{8} & 0 & 0\\
0 & 0 & 0 & 0 & 0 & 0 & 1 & 1 & 0\\
0 & 0 & 0 & 0 & 0 & 0 & 0 & 1+a_{10} & -\frac{1}{2}\\
0 & 0 & 0 & 0 & 0 & 0 & 0 & -\frac{1}{2} & 1+a_{11}%
\end{vmatrix}
-\\
&  \qquad\qquad-\frac{1}{2^{2}}%
\begin{vmatrix}
1+a_{1} & -\frac{1}{2} & 0 & 0 & 0 & 0 & 0 & 0 & 0\\
-\frac{1}{2} & 1+a_{2} & 0 & 0 & 0 & 0 & 0 & 0 & 0\\
0 & 0 & -\frac{1}{2} & 1+a_{5} & -\frac{1}{2} & 0 & 0 & 0 & 0\\
0 & 0 & 0 & -\frac{1}{2} & 1+a_{6} & -\frac{1}{2} & 0 & 0 & 0\\
0 & 0 & 0 & 0 & -\frac{1}{2} & 1+a_{7} & -\frac{1}{2} & 0 & 0\\
0 & 0 & 0 & 0 & 0 & -\frac{1}{2} & 1+a_{8} & 0 & 0\\
0 & 0 & 0 & 0 & 0 & 0 & 1 & 1 & 0\\
0 & 0 & 0 & 0 & 0 & 0 & 0 & 1+a_{10} & -\frac{1}{2}\\
0 & 0 & 0 & 0 & 0 & 0 & 0 & -\frac{1}{2} & 1+a_{11}%
\end{vmatrix}
\\
&  =0-\frac{1}{2^{2}}%
\begin{vmatrix}
1+a_{1} & -\frac{1}{2} & 0 & 0 & 0 & 0 & 0 & 0 & 0\\
-\frac{1}{2} & 1+a_{2} & 0 & 0 & 0 & 0 & 0 & 0 & 0\\
0 & 0 & -\frac{1}{2} & 1+a_{5} & -\frac{1}{2} & 0 & 0 & 0 & 0\\
0 & 0 & 0 & -\frac{1}{2} & 1+a_{6} & -\frac{1}{2} & 0 & 0 & 0\\
0 & 0 & 0 & 0 & -\frac{1}{2} & 1+a_{7} & -\frac{1}{2} & 0 & 0\\
0 & 0 & 0 & 0 & 0 & -\frac{1}{2} & 1+a_{8} & 0 & 0\\
0 & 0 & 0 & 0 & 0 & 0 & 1 & 1 & 0\\
0 & 0 & 0 & 0 & 0 & 0 & 0 & 1+a_{10} & -\frac{1}{2}\\
0 & 0 & 0 & 0 & 0 & 0 & 0 & -\frac{1}{2} & 1+a_{11}%
\end{vmatrix}
.
\end{align*}

Further%
\begin{align*}
&
\begin{vmatrix}
1+a_{1} & -\frac{1}{2} & 0 & 0 & 0 & 0 & 0 & 0 & 0\\
-\frac{1}{2} & 1+a_{2} & 0 & 0 & 0 & 0 & 0 & 0 & 0\\
0 & 0 & \fbox{$-\frac{1}{2}$} & 1+a_{5} & -\frac{1}{2} & 0 & 0 & 0 & 0\\
0 & 0 & 0 & -\frac{1}{2} & 1+a_{6} & -\frac{1}{2} & 0 & 0 & 0\\
0 & 0 & 0 & 0 & -\frac{1}{2} & 1+a_{7} & -\frac{1}{2} & 0 & 0\\
0 & 0 & 0 & 0 & 0 & -\frac{1}{2} & 1+a_{8} & 0 & 0\\
0 & 0 & 0 & 0 & 0 & 0 & 1 & 1 & 0\\
0 & 0 & 0 & 0 & 0 & 0 & 0 & 1+a_{10} & -\frac{1}{2}\\
0 & 0 & 0 & 0 & 0 & 0 & 0 & -\frac{1}{2} & 1+a_{11}%
\end{vmatrix}
\\
& \\
&  =%
\begin{vmatrix}
1+a_{1} & -\frac{1}{2}\\
-\frac{1}{2} & 1+a_{2}%
\end{vmatrix}%
\begin{vmatrix}
-\frac{1}{2} & 1+a_{5} & -\frac{1}{2} & 0 & 0 & 0 & 0\\
0 & -\frac{1}{2} & 1+a_{6} & -\frac{1}{2} & 0 & 0 & 0\\
0 & 0 & -\frac{1}{2} & 1+a_{7} & -\frac{1}{2} & 0 & 0\\
0 & 0 & 0 & -\frac{1}{2} & 1+a_{8} & 0 & 0\\
0 & 0 & 0 & 0 & \fbox{$1$} & \fbox{$1$} & 0\\
0 & 0 & 0 & 0 & 0 & 1+a_{10} & -\frac{1}{2}\\
0 & 0 & 0 & 0 & 0 & -\frac{1}{2} & 1+a_{11}%
\end{vmatrix}
\\
& \\
&  =%
\begin{vmatrix}
1+a_{1} & -\frac{1}{2}\\
-\frac{1}{2} & 1+a_{2}%
\end{vmatrix}
\left(
\begin{array}
[c]{c}%
\begin{vmatrix}
-\frac{1}{2} & 1+a_{5} & -\frac{1}{2} & 0 & 0 & 0\\
0 & -\frac{1}{2} & 1+a_{6} & -\frac{1}{2} & 0 & 0\\
0 & 0 & -\frac{1}{2} & 1+a_{7} & 0 & 0\\
0 & 0 & 0 & -\frac{1}{2} & 0 & 0\\
0 & 0 & 0 & 0 & 1+a_{10} & -\frac{1}{2}\\
0 & 0 & 0 & 0 & -\frac{1}{2} & 1+a_{11}%
\end{vmatrix}
-\\
-%
\begin{vmatrix}
-\frac{1}{2} & 1+a_{5} & -\frac{1}{2} & 0 & 0 & 0\\
0 & -\frac{1}{2} & 1+a_{6} & -\frac{1}{2} & 0 & 0\\
0 & 0 & -\frac{1}{2} & 1+a_{7} & -\frac{1}{2} & 0\\
0 & 0 & 0 & -\frac{1}{2} & 1+a_{8} & 0\\
0 & 0 & 0 & 0 & 0 & \fbox{$-\frac{1}{2}$}\\
0 & 0 & 0 & 0 & 0 & 1+a_{11}%
\end{vmatrix}
\end{array}
\right) \\
& \\
&  =%
\begin{vmatrix}
1+a_{1} & -\frac{1}{2}\\
-\frac{1}{2} & 1+a_{2}%
\end{vmatrix}
\left(
\begin{vmatrix}
-\frac{1}{2} & 1+a_{5} & -\frac{1}{2} & 0 & 0 & 0\\
0 & -\frac{1}{2} & 1+a_{6} & -\frac{1}{2} & 0 & 0\\
0 & 0 & -\frac{1}{2} & 1+a_{7} & 0 & 0\\
0 & 0 & 0 & -\frac{1}{2} & 0 & 0\\
0 & 0 & 0 & 0 & 1+a_{10} & -\frac{1}{2}\\
0 & 0 & 0 & 0 & -\frac{1}{2} & 1+a_{11}%
\end{vmatrix}
+0\right) \\
& \\
&  =%
\begin{vmatrix}
1+a_{1} & -\frac{1}{2}\\
-\frac{1}{2} & 1+a_{2}%
\end{vmatrix}%
\begin{vmatrix}
-\frac{1}{2} & 1+a_{5} & -\frac{1}{2} & 0\\
0 & -\frac{1}{2} & 1+a_{6} & -\frac{1}{2}\\
0 & 0 & -\frac{1}{2} & 1+a_{7}\\
0 & 0 & 0 & -\frac{1}{2}%
\end{vmatrix}%
\begin{vmatrix}
1+a_{10} & -\frac{1}{2}\\
-\frac{1}{2} & 1+a_{11}%
\end{vmatrix}
\\
& \\
&  =-\frac{1}{2^{4}}%
\begin{vmatrix}
1+a_{1} & -\frac{1}{2}\\
-\frac{1}{2} & 1+a_{2}%
\end{vmatrix}%
\begin{vmatrix}
1+a_{10} & -\frac{1}{2}\\
-\frac{1}{2} & 1+a_{11}%
\end{vmatrix}
,
\end{align*}

so that%
\[
M_{3,9}\left(  a_{1:11}\right)  =\frac{1}{2^{6}}%
\begin{vmatrix}
1+a_{1} & -\frac{1}{2}\\
-\frac{1}{2} & 1+a_{2}%
\end{vmatrix}%
\begin{vmatrix}
1+a_{10} & -\frac{1}{2}\\
-\frac{1}{2} & 1+a_{11}%
\end{vmatrix}
=\frac{1}{2^{6}}T\left(  a_{1:2}\right)  T\left(  a_{10:11}\right)  ,
\]

which supports \ref{av781}.
\end{example}

\subsection{Minor calculations summary\label{SbSect_minor_summary}}

We summarize the calculations from previous subsections here.

Suppose $a\in\mathbb{R}^{m}$. Then from \ref{av746}, \ref{av552}, when
$m\geq2$,%
\begin{equation}
M_{1,j}\left(  a\right)  =M_{j,1}\left(  a\right)  =\left\{
\begin{array}
[c]{ll}%
\left(  -\frac{1}{2}\right)  ^{j-1}\left\vert T\left(  a_{j+1:m}\right)
\right\vert , & 1\leq j\leq m-1,\\
\left(  -\frac{1}{2}\right)  ^{m-1}, & j=m.
\end{array}
\right. \label{av784}%
\end{equation}

From \ref{av552} and \ref{av747}, when $m\geq2$,%
\[
M_{n,m}\left(  a\right)  =M_{m,n}\left(  a\right)  =\left\{
\begin{array}
[c]{ll}%
\left(  -\frac{1}{2}\right)  ^{m-1}, & n=1,\\
\left(  -\frac{1}{2}\right)  ^{m-n}\left\vert T\left(  a_{1:n-1}\right)
\right\vert , & 2\leq n\leq m.
\end{array}
\right.
\]

From \ref{av781}: if $m\geq4$ then%
\[
M_{j,n}\left(  a_{1:m}\right)  =M_{j,n}\left(  a_{1:m}\right)  =\left(
-\frac{1}{2}\right)  ^{n-j}\left\vert T\left(  a_{1:j-1}\right)  \right\vert
\left\vert T\left(  a_{n+1:m}\right)  \right\vert ,\quad2\leq j<n\leq m-1,
\]

From \ref{av782}, if $m\geq3$ then%
\[
M_{n,n}\left(  a_{1:m}\right)  =%
\begin{array}
[c]{ll}%
\left\vert T_{2:m}\left(  a_{1:m}\right)  \right\vert , & n=1,\\
\left\vert T_{1:n-1}\left(  a_{1:m}\right)  \right\vert \left\vert
T_{n+1:m}\left(  a_{1:m}\right)  \right\vert , & 2\leq n\leq m-1,\\
\left\vert T_{1:m-1}\left(  a_{1:m}\right)  \right\vert , & n=m.
\end{array}
\]

and using the notation $T\left(  a_{1:0}\right)  :=1$ and $T\left(
a_{m+1:m}\right)  :=1$\ of Definition \ref{vDef_T_minors}, these formulas can
be combined compactly to give%
\begin{equation}
\left.
\begin{array}
[c]{l}%
M_{j,n}\left(  a\right)  =\left(  -\frac{1}{2}\right)  ^{n-j}\left\vert
T\left(  a_{1:j-1}\right)  \right\vert \left\vert T\left(  a_{n+1:m}\right)
\right\vert ,\\
M_{n,j}\left(  a\right)  =M_{j:n}\left(  a\right)  .
\end{array}
\right\}  ,\quad1\leq j\leq n\leq m.\label{av787}%
\end{equation}

Thus we can conclude that all the minors of $T\left(  a\right)  $ can all be
written as the product of one or two principal sub-determinants of $T\left(
a\right)  $ and that the variables are always separated.

When $a=\frac{1}{v}\Delta_{1}X:=\left(  \Delta_{1}x^{\left(  k\right)
}\right)  _{k=1}^{N-1}$ the minor formulas \ref{av787} become%
\begin{equation}%
\begin{array}
[c]{ll}%
M_{j,n}=\left(  -\frac{1}{2}\right)  ^{n-j}\left\vert T_{1:j-1}\right\vert
\left\vert T_{n+1:N-1}\right\vert , & 1\leq j\leq n\leq N-1,\\
M_{n,j}=M_{j,n}, &
\end{array}
\label{av023}%
\end{equation}

where $T_{1:0}:=1$ and $T_{N:N-1}:=1$. These expressions for the minors will
later be substituted into equations \ref{av725} to \ref{av757}.

\section{A formula for $\left\vert T\left(  a\right)  \right\vert $}

The determinant $\left\vert T\right\vert $ occurs in the denominator of the
formulas \ref{av725} to \ref{av757}. Now from Definition \ref{vDef_T_minors},
$T:=T\left(  \frac{1}{v}\Delta_{1}X;1\right)  $ and since $\frac{1}{v}%
\Delta_{1}X\in\mathbb{R}^{N-1}$ we will calculate the determinant of $T\left(
a\right)  =T\left(  a;1\right)  $ for $a\in\mathbb{R}^{m}$ and $m\geq1 $. More
precisely, we will derive a formula for the positive coefficients $c_{\beta
}^{\left(  m\right)  }$ of the monomials $a^{\beta}$ in the equation
\ref{av599} for $\left\vert T\left(  a\right)  \right\vert $. This formula is
given in Theorem \ref{vThm_formula_for_C^(m)_idx}.

Observe that when $m\geq3$,%
\begin{align}
\left\vert T\left(  a_{1:m}\right)  \right\vert  & =\left(  1+a_{m}\right)
\left\vert T\left(  a_{1:m-1}\right)  \right\vert +\frac{1}{2}%
\begin{vmatrix}
1+a_{1} & -\frac{1}{2} & 0 &  & 0 & 0 & 0\\
-\frac{1}{2} & 1+a_{2} & -\frac{1}{2} &  & 0 & 0 & 0\\
0 & -\frac{1}{2} & 1+a_{3} &  & 0 & 0 & 0\\
&  &  & \ddots &  &  & \\
0 & 0 & 0 &  & 1+a_{m-3} & -\frac{1}{2} & 0\\
0 & 0 & 0 &  & -\frac{1}{2} & 1+a_{m-2} & -\frac{1}{2}\\
0 & 0 & 0 &  & 0 & 0 & -\frac{1}{2}%
\end{vmatrix}
\nonumber\\
& =\left(  1+a_{m}\right)  \left\vert T\left(  a_{1:m-1}\right)  \right\vert
-\frac{1}{4}\left\vert T\left(  a_{1:m-2}\right)  \right\vert ,\label{av585}%
\end{align}

and also%
\begin{equation}
\left\vert T\left(  a_{1}\right)  \right\vert =1+a_{1},\quad\left\vert
T\left(  a_{1:2}\right)  \right\vert =\left(  1+a_{1}\right)  \left(
1+a_{2}\right)  -\frac{1}{4},\label{av586}%
\end{equation}

so that in the notation of Definition \ref{vDef_T_minors},
\begin{equation}
\left\vert T\left(  a_{1:m}\right)  \right\vert =\left(  1+a_{m}\right)
\left\vert T\left(  a_{1:m-1}\right)  \right\vert -\frac{1}{4}\left\vert
T\left(  a_{1:m-2}\right)  \right\vert ,\quad m\geq1.\label{av540}%
\end{equation}

A well known tridiagonal matrix result is:

\begin{lemma}
\label{vLem_T(0m)}\
\begin{equation}
\left\vert T\left(  0_{m}\right)  \right\vert =\frac{m+1}{2^{m}},\quad
m\geq1.\label{av004}%
\end{equation}

\end{lemma}

\begin{proof}
From \ref{av586} and \ref{av585} we have the linear difference equation,
$\left\vert T\left(  0_{m}\right)  \right\vert =\left\vert T\left(
0_{m-1}\right)  \right\vert -\frac{1}{4}\left\vert T\left(  0_{m-2}\right)
\right\vert $ with initial conditions $\left\vert T\left(  0_{1}\right)
\right\vert =1$ and $\left\vert T\left(  0_{2}\right)  \right\vert =\frac
{3}{4} $. The auxiliary equation has repeated root $\alpha=1/2$ and since the
solution has the form $\left(  A+Bm\right)  \alpha^{m}$ the initial conditions
imply the solution \ref{av004}.
\end{proof}

Compare this with Lemma \ref{vLem_T(s1m)} below which derives a formula for
$T\left(  s1_{m}\right)  $.

The elementary symmetric functions will prove useful in deriving formulas and
bounds for $\left\vert T\left(  a\right)  \right\vert $:

\begin{definition}
\label{vDef_elem_sym_funct}\textbf{Elementary symmetric function }$\sigma
_{k}\left(  x\right)  =\sum\limits_{\substack{\mathbf{0}\leq\alpha
\leq\mathbf{1} \\\left\vert \alpha\right\vert =k}}x^{\alpha}$ \textbf{of
order} $k$.

Some properties are:

\begin{enumerate}
\item $\sigma_{k}\left(  tx\right)  =t^{k}\sigma_{k}\left(  x\right)  $;

\item $\sum\limits_{k=0}^{m}\sigma_{k}\left(  x\right)  =\left(  1+x\right)
^{\mathbf{1}}$;

\item $\sigma_{k}\left(  1_{m}\right)  =\tbinom{m}{k}$.
\end{enumerate}
\end{definition}

Now%
\begin{equation}
\left\vert T\left(  a_{1:2}\right)  \right\vert =\frac{3}{4}+\left(
a_{1}+a_{2}\right)  +a_{1}a_{2}=\frac{3}{4}+\sigma_{1}\left(  a_{1:2}\right)
+\sigma_{2}\left(  a_{1:2}\right)  ,\label{av587}%
\end{equation}

so consequently%
\begin{align}
\left\vert T\left(  a_{1:3}\right)  \right\vert  &  =\left(  1+a_{3}\right)
\left\vert T\left(  a_{1:2}\right)  \right\vert -\frac{1}{4}\left\vert
T\left(  a_{1}\right)  \right\vert \nonumber\\
&  =\left(  1+a_{3}\right)  \left(  \left(  1+a_{1}\right)  \left(
1+a_{2}\right)  -\frac{1}{4}1\right)  -\frac{1}{4}\left(  1+a_{1}\right)
\nonumber\\
&  =\left(  1+a_{1}\right)  \left(  1+a_{2}\right)  \left(  1+a_{3}\right)
-\frac{1}{4}\left(  1+a_{3}\right)  -\frac{1}{4}\left(  1+a_{1}\right)
\nonumber\\
&  =\frac{1}{2}+\left(  a_{1}+a_{2}+a_{3}-\frac{1}{4}a_{1}-\frac{1}{4}%
a_{3}\right)  +\sigma_{2}\left(  a_{1:3}\right)  +\sigma_{3}\left(
a_{1:3}\right) \nonumber\\
&  =\frac{4}{8}+\left(  \frac{3}{4}a_{1}+a_{2}+\frac{3}{4}a_{3}\right)
+\sigma_{2}\left(  a_{1:3}\right)  +\sigma_{3}\left(  a_{1:3}\right)
\label{av589}\\
&  =\frac{4}{8}+\left(  \frac{3}{4}a_{1}+a_{2}+\frac{3}{4}a_{3}\right)
+\sigma_{2}\left(  a_{1:3}\right)  +\sigma_{3}\left(  a_{1:3}\right)
.\nonumber
\end{align}

Further%
\begin{align}
&  \left\vert T\left(  a_{1:4}\right)  \right\vert \nonumber\\
&  =\left(  1+a_{4}\right)  \left\vert T\left(  a_{1:3}\right)  \right\vert
-\frac{1}{4}\left\vert T\left(  a_{1:2}\right)  \right\vert \nonumber\\
&  =\left(  1+a_{4}\right)  \left(  \frac{1}{2}+\left(  \frac{3}{4}a_{1}%
+a_{2}+\frac{3}{4}a_{3}\right)  +\sigma_{2}\left(  a_{1:3}\right)  +\sigma
_{3}\left(  a_{1:3}\right)  \right)  -\nonumber\\
&  \qquad\qquad-\frac{1}{4}\left(  \frac{3}{4}+\left(  a_{1}+a_{2}\right)
+a_{1}a_{2}\right) \nonumber\\
&  =\frac{1}{2}+\left(  \frac{3}{4}a_{1}+a_{2}+\frac{3}{4}a_{3}\right)
+\sigma_{2}\left(  a_{1:3}\right)  +\sigma_{3}\left(  a_{1:3}\right)
+\frac{1}{2}a_{4}+\nonumber\\
&  \qquad+\left(  \frac{3}{4}a_{1}+a_{2}+\frac{3}{4}a_{3}\right)  a_{4}%
+\sigma_{2}\left(  a_{1:3}\right)  a_{4}+\sigma_{3}\left(  a_{1:3}\right)
a_{4}-\frac{3}{16}-\frac{1}{4}\left(  a_{1}+a_{2}\right)  -\frac{1}{4}%
a_{1}a_{2}\nonumber\\
&  =\frac{5}{16}+\left(  \frac{3}{4}a_{1}+a_{2}+\frac{3}{4}a_{3}+\frac{1}%
{2}a_{4}-\frac{1}{4}\left(  a_{1}+a_{2}\right)  \right)  +\nonumber\\
&  \qquad+\left(  \sigma_{2}\left(  a_{1:3}\right)  +\left(  \frac{3}{4}%
a_{1}+a_{2}+\frac{3}{4}a_{3}\right)  a_{4}-\frac{1}{4}a_{1}a_{2}\right)
+\left(  \sigma_{3}\left(  a_{1:3}\right)  +\sigma_{2}\left(  a_{1:3}\right)
a_{4}\right)  +\nonumber\\
&  \qquad+\sigma_{3}\left(  a_{1:3}\right)  a_{4}\nonumber\\
&  =\frac{5}{16}+\left(  \frac{1}{2}a_{1}+\frac{3}{4}a_{2}+\frac{3}{4}%
a_{3}+\frac{1}{2}a_{4}\right)  +\nonumber\\
&  \qquad+\left(  a_{1}a_{2}+a_{1}a_{3}+a_{2}a_{3}+\frac{3}{4}a_{1}a_{4}%
+a_{2}a_{4}+\frac{3}{4}a_{3}a_{4}-\frac{1}{4}a_{1}a_{2}\right)  +\nonumber\\
&  \qquad+\sigma_{3}\left(  a_{1:4}\right)  +\sigma_{4}\left(  a_{1:4}\right)
\nonumber\\
&  =\frac{5}{16}+\left(  \frac{1}{2}a_{1}+\frac{3}{4}a_{2}+\frac{3}{4}%
a_{3}+\frac{1}{2}a_{4}\right)  +\label{av588}\\
&  \qquad+\left(  \frac{3}{4}a_{1}a_{2}+a_{1}a_{3}+\frac{3}{4}a_{1}a_{4}%
+a_{2}a_{3}+a_{2}a_{4}+\frac{3}{4}a_{3}a_{4}\right)  +\nonumber\\
&  \qquad+\sigma_{3}\left(  a_{1:4}\right)  +\sigma_{4}\left(  a_{1:4}\right)
,\nonumber
\end{align}

and%
\begin{align*}
\left\vert T\left(  a_{1:5}\right)  \right\vert  & =\left(  1+a_{5}\right)
\left\vert T\left(  a_{1:4}\right)  \right\vert -\frac{1}{4}\left\vert
T\left(  a_{1:3}\right)  \right\vert \\
& =-\frac{1}{4}\left\vert T\left(  a_{1:3}\right)  \right\vert +\left\vert
T\left(  a_{1:4}\right)  \right\vert +\left\vert T\left(  a_{1:4}\right)
\right\vert a_{5}.
\end{align*}

Now%
\begin{align*}
\left\vert T\left(  a_{1:4}\right)  \right\vert  &  =\frac{5}{16}+\left(
\frac{1}{2}a_{1}+\frac{3}{4}a_{2}+\frac{3}{4}a_{3}+\frac{1}{2}a_{4}\right)
+\\
&  +\left(  \frac{3}{4}a_{1}a_{2}+a_{1}a_{3}+\frac{3}{4}a_{1}a_{4}+a_{2}%
a_{3}+a_{2}a_{4}+\frac{3}{4}a_{3}a_{4}\right)  +\\
&  +\sigma_{3}\left(  a_{1:4}\right)  +\sigma_{4}\left(  a_{1:4}\right)  ,
\end{align*}

and%
\begin{align*}
\left\vert T\left(  a_{1:4}\right)  \right\vert a_{5}  & =0+\frac{5}{16}%
a_{5}+\left(  \frac{1}{2}a_{1}a_{5}+\frac{3}{4}a_{2}a_{5}+\frac{3}{4}%
a_{3}a_{5}+\frac{1}{2}a_{4}a_{5}\right)  +\\
& +\left(  \frac{3}{4}a_{1}a_{2}a_{5}+a_{1}a_{3}a_{5}+\frac{3}{4}a_{1}%
a_{4}a_{5}+a_{2}a_{3}a_{5}+a_{2}a_{4}a_{5}+\frac{3}{4}a_{3}a_{4}a_{5}\right)
+\\
& +\sigma_{3}\left(  a_{1:4}\right)  a_{5}+\sigma_{4}\left(  a_{1:4}\right)
a_{5},
\end{align*}

and%
\begin{align*}
-\frac{1}{4}\left\vert T\left(  a_{1:3}\right)  \right\vert  & =-\frac{1}%
{4}\frac{4}{8}-\frac{1}{4}\left(  \frac{3}{4}a_{1}+a_{2}+\frac{3}{4}%
a_{3}\right)  -\frac{1}{4}\sigma_{2}\left(  a_{1:3}\right)  -\frac{1}{4}%
\sigma_{3}\left(  a_{1:3}\right) \\
& =-\frac{1}{8}+\left(  -\frac{3}{16}a_{1}-\frac{1}{4}a_{2}-\frac{3}{16}%
a_{3}\right)  -\frac{1}{4}\sigma_{2}\left(  a_{1:3}\right)  -\frac{1}{4}%
\sigma_{3}\left(  a_{1:3}\right)
\end{align*}

so that%
\begin{align}
\left\vert T\left(  a_{1:5}\right)  \right\vert  & =\left(  \frac{5}%
{16}+0-\frac{1}{8}\right)  +\nonumber\\
& +\left(  \frac{1}{2}a_{1}+\frac{3}{4}a_{2}+\frac{3}{4}a_{3}+\frac{1}{2}%
a_{4}\right)  +\left(  \frac{5}{16}a_{5}\right)  +\left(  -\frac{3}{16}%
a_{1}-\frac{1}{4}a_{2}-\frac{3}{16}a_{3}\right)  +\nonumber\\
& +\left(
\begin{array}
[c]{l}%
\frac{3}{4}a_{1}a_{2}+a_{1}a_{3}+\frac{3}{4}a_{1}a_{4}+a_{2}a_{3}+a_{2}%
a_{4}+\frac{3}{4}a_{3}a_{4}+\\
+\frac{1}{2}a_{1}a_{5}+\frac{3}{4}a_{2}a_{5}+\frac{3}{4}a_{3}a_{5}+\frac{1}%
{2}a_{4}a_{5}-\frac{1}{4}\sigma_{2}\left(  a_{1:3}\right)
\end{array}
\right)  +\nonumber\\
& +\left(
\begin{array}
[c]{l}%
\sigma_{3}\left(  a_{1:4}\right)  +\\
+\left(  \frac{3}{4}a_{1}a_{2}a_{5}+a_{1}a_{3}a_{5}+\frac{3}{4}a_{1}a_{4}%
a_{5}+a_{2}a_{3}a_{5}+a_{2}a_{4}a_{5}+\frac{3}{4}a_{3}a_{4}a_{5}\right)
-\frac{1}{4}a_{1}a_{2}a_{3}%
\end{array}
\right)  +\nonumber\\
& +\sigma_{4}\left(  a_{1:4}\right)  +\sigma_{3}\left(  a_{1:4}\right)
a_{5}+\nonumber\\
& +\sigma_{4}\left(  a_{1:4}\right)  a_{5}\nonumber\\
& =\frac{3}{16}+\left(  \frac{5}{16}a_{1}+\frac{1}{2}a_{2}+\frac{9}{16}%
a_{3}+\frac{1}{2}a_{4}+\frac{5}{16}a_{5}\right)  +\nonumber\\
& +\left(
\begin{array}
[c]{l}%
\frac{1}{2}a_{1}a_{2}+\frac{3}{4}a_{1}a_{3}+\frac{3}{4}a_{1}a_{4}+\frac{3}%
{4}a_{2}a_{3}+a_{2}a_{4}+\frac{3}{4}a_{3}a_{4}+\\
+\frac{1}{2}a_{1}a_{5}+\frac{3}{4}a_{2}a_{5}+\frac{3}{4}a_{3}a_{5}+\frac{1}%
{2}a_{4}a_{5}%
\end{array}
\right)  +\nonumber\\
& +\left(
\begin{array}
[c]{l}%
a_{1}a_{2}a_{3}+a_{1}a_{2}a_{4}+a_{1}a_{3}a_{4}+a_{2}a_{3}a_{4}+\\
+\left(  \frac{3}{4}a_{1}a_{2}a_{5}+a_{1}a_{3}a_{5}+\frac{3}{4}a_{1}a_{4}%
a_{5}+a_{2}a_{3}a_{5}+a_{2}a_{4}a_{5}+\frac{3}{4}a_{3}a_{4}a_{5}\right)
-\frac{1}{4}a_{1}a_{2}a_{3}%
\end{array}
\right) \nonumber\\
& +\sigma_{4}\left(  a_{1:5}\right)  +\sigma_{5}\left(  a_{1:5}\right)
\nonumber\\
& =\frac{3}{16}+\left(  \frac{5}{16}a_{1}+\frac{1}{2}a_{2}+\frac{9}{16}%
a_{3}+\frac{1}{2}a_{4}+\frac{5}{16}a_{5}\right)  +\nonumber\\
& +\left(
\begin{array}
[c]{l}%
\frac{1}{2}a_{1}a_{2}+\frac{3}{4}a_{1}a_{3}+\frac{3}{4}a_{1}a_{4}+\frac{3}%
{4}a_{2}a_{3}+a_{2}a_{4}+\frac{3}{4}a_{3}a_{4}+\\
+\frac{1}{2}a_{1}a_{5}+\frac{3}{4}a_{2}a_{5}+\frac{3}{4}a_{3}a_{5}+\frac{1}%
{2}a_{4}a_{5}%
\end{array}
\right)  +\nonumber\\
& +\left(
\begin{array}
[c]{l}%
\frac{3}{4}a_{1}a_{2}a_{3}+a_{1}a_{2}a_{4}+a_{1}a_{3}a_{4}+a_{2}a_{3}a_{4}+\\
+\frac{3}{4}a_{1}a_{2}a_{5}+a_{1}a_{3}a_{5}+\frac{3}{4}a_{1}a_{4}a_{5}%
+a_{2}a_{3}a_{5}+a_{2}a_{4}a_{5}+\frac{3}{4}a_{3}a_{4}a_{5}%
\end{array}
\right)  +\nonumber\\
& +\sigma_{4}\left(  a_{1:5}\right)  +\sigma_{5}\left(  a_{1:5}\right)
,\label{av563}%
\end{align}

The expressions for the determinants $\left\{  \left\vert T\left(
a_{1:k}\right)  \right\vert \right\}  _{m=1}^{5}$ obtained above yield the
following bounds in terms of the symmetric functions $\sigma_{k}$ described in
Definition \ref{vDef_elem_sym_funct}:%
\begin{align*}
\left\vert T\left(  a_{1}\right)  \right\vert  & =1+a_{1},\\
\left\vert T\left(  a_{1:2}\right)  \right\vert  & =\frac{3}{4}+\sigma
_{1}\left(  a_{1:2}\right)  +\sigma_{2}\left(  a_{1:2}\right)  ,
\end{align*}%
\[
\frac{4}{8}+\frac{3}{4}\sigma_{1}\left(  a_{1:3}\right)  +\sigma_{2}\left(
a_{1:3}\right)  +\sigma_{3}\left(  a_{1:3}\right)  \leq\left\vert T\left(
a_{1:3}\right)  \right\vert \leq\frac{4}{8}+\sigma_{1}\left(  a_{1:3}\right)
+\sigma_{2}\left(  a_{1:3}\right)  +\sigma_{3}\left(  a_{1:3}\right)  ,
\]

\begin{multline*}
\frac{5}{16}+\frac{1}{2}\sigma_{1}\left(  a_{1:4}\right)  +\frac{3}{4}%
\sigma_{2}\left(  a_{1:4}\right)  +\sigma_{3}\left(  a_{1:4}\right)
+\sigma_{4}\left(  a_{1:4}\right) \\
\leq\left\vert T\left(  a_{1:4}\right)  \right\vert \leq\\
\frac{5}{16}+\frac{1}{2}\sigma_{1}\left(  a_{1:4}\right)  +\sigma_{2}\left(
a_{1:4}\right)  +\sigma_{3}\left(  a_{1:4}\right)  +\sigma_{4}\left(
a_{1:4}\right)
\end{multline*}%
\begin{multline*}
\frac{6}{32}+\frac{1}{2}\sigma_{1}\left(  a_{1:5}\right)  +\frac{1}{2}%
\sigma_{2}\left(  a_{1:5}\right)  +\frac{1}{2}\sigma_{3}\left(  a_{1:5}%
\right)  +\frac{3}{4}\sigma_{4}\left(  a_{1:5}\right)  +\sigma_{5}\left(
a_{1:5}\right) \\
\leq\left\vert T\left(  a_{1:5}\right)  \right\vert \leq\\
\frac{6}{32}+\frac{9}{16}\sigma_{1}\left(  a_{1:5}\right)  +\sigma_{2}\left(
a_{1:5}\right)  +\sigma_{3}\left(  a_{1:5}\right)  +\sigma_{4}\left(
a_{1:5}\right)  +\sigma_{5}\left(  a_{1:5}\right)  .
\end{multline*}

\subsection{A general formula for $\left\vert T\left(  a\right)  \right\vert
$}

From the above examples and \ref{av004} we hypothesize the following
representation for $m\geq1$: if $a\in\mathbb{R}^{m}$ then
\begin{align}
\left\vert T\left(  a\right)  \right\vert  &  =\frac{m+1}{2^{m}}%
+\sum\limits_{\substack{\beta\leq1_{m} \\\left\vert \beta\right\vert
=1}}c_{\beta}^{\left(  m\right)  }a^{\beta}+\sum\limits_{\substack{\beta
\leq1_{m} \\\left\vert \beta\right\vert =2}}c_{\beta}^{\left(  m\right)
}a^{\beta}+\ldots+\sum\limits_{\substack{\beta\leq1_{m} \\\left\vert
\beta\right\vert =m-1}}c_{\beta}^{\left(  m\right)  }a^{\beta}+\sum
\limits_{\substack{\beta\leq1_{m} \\\left\vert \beta\right\vert =m}}c_{\beta
}^{\left(  m\right)  }a^{\beta}\label{av596}\\
&  =\sum_{k=0}^{m}\sum\limits_{\beta\leq1_{m},\left\vert \beta\right\vert
=k}c_{\beta}^{\left(  m\right)  }a^{\beta},\quad where\text{ }c_{\mathbf{0}%
}^{\left(  m\right)  }=\frac{m+1}{2^{m}},\label{av599}%
\end{align}

and it will be shown in Theorem \ref{vThm_formula_for_C^(m)_idx} below that
$0<c_{\beta}^{\left(  m\right)  }\leq1$ for all $\beta$.

Regarding notation we could write in full%
\[
\left\vert T\left(  a\right)  \right\vert =\sum_{k=0}^{m}\sum
\limits_{\substack{\beta\leq\mathbf{1},\left\vert \beta\right\vert =k
\\\beta\in\mathbb{Z}_{\oplus}^{m}}}c_{\beta}^{\left(  m\right)  }a^{\beta
}=\sum_{k=0}^{m}\sum\limits_{\substack{0\mathbf{\leq}\beta\leq1_{m}
\\\left\vert \beta\right\vert =k}}c_{\beta}^{\left(  m\right)  }a^{\beta},
\]

but since the constants $c_{\beta}^{\left(  m\right)  }$ always have $\beta
\in\mathbb{Z}_{\oplus}^{m}$ we will omit this restraint and write%
\[
\left\vert T\left(  a\right)  \right\vert =\sum_{k=0}^{m}\sum
\limits_{\substack{\beta\leq1 \\\left\vert \beta\right\vert =k}}c_{\beta
}^{\left(  m\right)  }a^{\beta},
\]

except when we have something like $c_{\left(  \beta,1\right)  }^{\left(
m\right)  }$ and then we will use a summation such as $\sum
\limits_{\substack{\beta\leq\mathbf{1},\left\vert \beta\right\vert =k
\\\beta\in\mathbb{Z}_{\oplus}^{m-1}}}$ or $\sum\limits_{\substack{\beta
\leq1_{m-1} \\\left\vert \beta\right\vert =k}}$. Since all $c_{\beta}^{\left(
m\right)  }$ also satisfy $\beta\leq\mathbf{1}$ we could omit $\beta
\leq\mathbf{1}$ but I choose not to.

Rearranging the summations in \ref{av599}:%
\begin{align}
\left\vert T\left(  a\right)  \right\vert =\sum_{k=0}^{m}\sum\limits_{\beta
\leq1_{m},\left\vert \beta\right\vert =k}c_{\beta}^{\left(  m\right)
}a^{\beta}=\sum_{k=0}^{m}\sum\limits_{\beta\leq\mathbf{1},\left\vert
\beta\right\vert =m-k}c_{\beta}^{\left(  m\right)  }a^{\beta}  & =\sum
_{k=0}^{m}\sum\limits_{\beta\leq\mathbf{1},\left\vert \mathbf{1-}%
\beta\right\vert =k}c_{\beta}^{\left(  m\right)  }a^{\beta}\nonumber\\
& =\sum_{k=0}^{m}\sum\limits_{\beta\leq\mathbf{1},\left\vert \beta\right\vert
=k}c_{\mathbf{1}-\beta}^{\left(  m\right)  }a^{\mathbf{1}-\beta}.\label{av735}%
\end{align}

From \ref{av586}, \ref{av587} and \ref{av589},%
\begin{equation}%
\begin{array}
[c]{l}%
c_{0}^{\left(  1\right)  }=c_{1}^{\left(  1\right)  }=1,\\
c_{\mathbf{0}}^{\left(  2\right)  }=\frac{3}{4},\quad c_{0,1}^{\left(
2\right)  }=c_{1,0}^{\left(  2\right)  }=1,\quad c_{\mathbf{1}}^{\left(
2\right)  }=1,\\
c_{\mathbf{0}}^{\left(  3\right)  }=\frac{1}{2},\quad c_{0,0,1}^{\left(
3\right)  }=c_{1,0,0}^{\left(  3\right)  }=\frac{3}{4},\text{ }c_{0,1,0}%
^{\left(  3\right)  }=1,\quad c_{0,1,1}^{\left(  3\right)  }=c_{1,1,0}%
^{\left(  3\right)  }=c_{1,0,1}^{\left(  3\right)  }=1,\quad c_{\mathbf{1}%
}^{\left(  3\right)  }=1.
\end{array}
\label{av591}%
\end{equation}

The calculations below will prove Theorem \ref{vThm_formula_for_C^(m)_idx} and
along the way formulas for special $c_{\beta}^{\left(  m\right)  }$ and bounds
for the $c_{\beta}^{\left(  m\right)  }$ e.g. \ref{av820} and \ref{av618} as
well as the bounds in the next subsection.

Using \ref{av599} and the iterative formula \ref{av540} i.e.%
\[
\left\vert T\left(  a_{1:m}\right)  \right\vert =\left(  1+a_{m}\right)
\left\vert T\left(  a_{1:m-1}\right)  \right\vert -\frac{1}{4}\left\vert
T\left(  a_{1:m-2}\right)  \right\vert ,\quad m\geq3,
\]

we have for $m\geq4$:%
\begin{align}
&  \frac{m+1}{2^{m}}+\sum_{k=1}^{m}\sum\limits_{\beta\leq\mathbf{1},\left\vert
\beta\right\vert =k}c_{\beta}^{\left(  m\right)  }a^{\beta}\label{av98}\\
&  =\left(  1+a_{m}\right)  \left(  \frac{m}{2^{m-1}}+\sum_{k=1}^{m-1}%
\sum\limits_{\substack{\gamma\leq\mathbf{1} \\\left\vert \gamma\right\vert
=k}}c_{\gamma}^{\left(  m-1\right)  }a_{1:m-1}^{\gamma}\right)  -\frac{1}%
{4}\left(  \frac{m-1}{2^{m-2}}+\sum_{k=1}^{m-2}\sum\limits_{\substack{\delta
\leq\mathbf{1} \\\left\vert \delta\right\vert =k}}c_{\delta}^{\left(
m-2\right)  }a_{1:m-2}^{\delta}\right) \nonumber\\
&  =\frac{m}{2^{m-1}}+\sum_{k=1}^{m-1}\sum\limits_{\gamma\leq\mathbf{1}%
,\left\vert \gamma\right\vert =k}c_{\gamma}^{\left(  m-1\right)  }%
a_{1:m-1}^{\gamma}+\frac{m}{2^{m-1}}a_{m}+\sum_{k=1}^{m-1}a_{m}\sum
\limits_{\gamma\leq\mathbf{1},\left\vert \gamma\right\vert =k}c_{\gamma
}^{\left(  m-1\right)  }a_{1:m-1}^{\gamma}-\nonumber\\
&  \qquad\qquad-\frac{m-1}{2^{m}}-\frac{1}{4}\sum_{k=1}^{m-2}\sum
\limits_{\delta\leq\mathbf{1},\left\vert \delta\right\vert =k}c_{\delta
}^{\left(  m-2\right)  }a_{1:m-2}^{\delta}\nonumber\\
&  =\left(  \frac{m+1}{2^{m}}+\frac{m}{2^{m-1}}a_{m}\right)  +\sum_{k=1}%
^{m-1}\sum\limits_{\gamma\leq\mathbf{1},\left\vert \gamma\right\vert
=k}c_{\gamma}^{\left(  m-1\right)  }a_{1:m-1}^{\gamma}+\sum_{k=2}^{m}a_{m}%
\sum\limits_{\gamma\leq\mathbf{1,}\left\vert \gamma\right\vert =k-1}c_{\gamma
}^{\left(  m-1\right)  }a_{1:m-1}^{\gamma}-\nonumber\\
&  \qquad\qquad-\frac{1}{4}\sum_{k=1}^{m-2}\sum\limits_{\delta\leq
\mathbf{1},\left\vert \delta\right\vert =k}c_{\delta}^{\left(  m-2\right)
}a_{1:m-2}^{\delta}\nonumber\\
&  =\left(  \frac{m+1}{2^{m}}+\frac{m}{2^{m-1}}a_{m}\right)  +\sum_{k=1}%
\sum\limits_{\gamma\leq\mathbf{1,}\left\vert \gamma\right\vert =k}c_{\gamma
}^{\left(  m-1\right)  }a_{1:m-1}^{\gamma}+\sum_{k=2}^{m-2}\sum\limits_{\gamma
\leq\mathbf{1,}\left\vert \gamma\right\vert =k}c_{\gamma}^{\left(  m-1\right)
}a_{1:m-1}^{\gamma}+\nonumber\\
&  \qquad+\sum_{k=m-1}\sum\limits_{\gamma\leq\mathbf{1,}\left\vert
\gamma\right\vert =k}c_{\gamma}^{\left(  m-1\right)  }a_{1:m-1}^{\gamma
}+\nonumber\\
&  \qquad+\sum_{k=2}^{m-2}a_{m}\sum\limits_{\gamma\leq\mathbf{1,}\left\vert
\gamma\right\vert =k-1}c_{\gamma}^{\left(  m-1\right)  }a_{1:m-1}^{\gamma
}+\sum_{k=m-1}^{m}a_{m}\sum\limits_{\gamma\leq\mathbf{1,}\left\vert
\gamma\right\vert =k-1}c_{\gamma}^{\left(  m-1\right)  }a_{1:m-1}^{\gamma
}-\nonumber\\
&  \qquad-\frac{1}{4}\sum_{k=1}\sum\limits_{\delta\leq\mathbf{1},\left\vert
\delta\right\vert =k}c_{\delta}^{\left(  m-2\right)  }a_{1:m-2}^{\delta}%
-\frac{1}{4}\sum_{k=2}^{m-2}\sum\limits_{\delta\leq\mathbf{1},\left\vert
\delta\right\vert =k}c_{\delta}^{\left(  m-2\right)  }a_{1:m-2}^{\delta
}\nonumber\\
&  =\frac{m+1}{2^{m}}+\frac{m}{2^{m-1}}a_{m}+\sum\limits_{\substack{\gamma
\leq\mathbf{1} \\\left\vert \gamma\right\vert =1}}c_{\gamma}^{\left(
m-1\right)  }a_{1:m-1}^{\gamma}+\sum_{k=2}^{m-2}\sum\limits_{\substack{\gamma
\leq\mathbf{1} \\\left\vert \gamma\right\vert =k}}c_{\gamma}^{\left(
m-1\right)  }a_{1:m-1}^{\gamma}+\sum\limits_{\substack{\gamma\leq\mathbf{1}
\\\left\vert \gamma\right\vert =m-1}}c_{\gamma}^{\left(  m-1\right)
}a_{1:m-1}^{\gamma}+\nonumber\\
&  \quad+\sum_{k=2}^{m-2}a_{m}\sum\limits_{\gamma\leq\mathbf{1,}\left\vert
\gamma\right\vert =k-1}c_{\gamma}^{\left(  m-1\right)  }a_{1:m-1}^{\gamma
}+a_{m}\sum\limits_{\gamma\leq\mathbf{1},\left\vert \gamma\right\vert
=m-2}c_{\gamma}^{\left(  m-1\right)  }a_{1:m-2}^{\gamma}+\nonumber\\
&  \quad+a_{m}\sum\limits_{\gamma\leq\mathbf{1},\left\vert \gamma\right\vert
=m-1}c_{\gamma}^{\left(  m-1\right)  }a_{m-1}^{\gamma}-\frac{1}{4}%
\sum\limits_{\delta\leq\mathbf{1,}\left\vert \delta\right\vert =1}c_{\delta
}^{\left(  m-2\right)  }a_{1:m-2}^{\delta}-\frac{1}{4}\sum_{k=2}^{m-2}%
\sum\limits_{\delta\leq\mathbf{1,}\left\vert \delta\right\vert =k}c_{\delta
}^{\left(  m-2\right)  }a_{1:m-2}^{\delta}\nonumber\\
&  =\frac{m+1}{2^{m}}+\frac{m}{2^{m-1}}a_{m}+\sum\limits_{\gamma
\leq\mathbf{1,}\left\vert \gamma\right\vert =1}c_{\gamma}^{\left(  m-1\right)
}a_{1:m-1}^{\gamma}-\frac{1}{4}\sum\limits_{\delta\leq\mathbf{1,}\left\vert
\delta\right\vert =1}c_{\delta}^{\left(  m-2\right)  }a_{1:m-2}^{\delta
}+\nonumber\\
&  \quad+\sum_{k=2}^{m-2}\left(  \sum\limits_{\gamma\leq\mathbf{1,}\left\vert
\gamma\right\vert =k}c_{\gamma}^{\left(  m-1\right)  }a_{1:m-1}^{\gamma}%
-\frac{1}{4}\sum\limits_{\delta\leq\mathbf{1,}\left\vert \delta\right\vert
=k}c_{\delta}^{\left(  m-2\right)  }a_{1:m-2}^{\delta}+a_{m}\sum
\limits_{\gamma\leq\mathbf{1,}\left\vert \gamma\right\vert =k-1}c_{\gamma
}^{\left(  m-1\right)  }a_{1:m-1}^{\gamma}\right)  +\nonumber\\
&  \quad+\sum\limits_{\gamma\leq\mathbf{1,}\left\vert \gamma\right\vert
=m-1}c_{\gamma}^{\left(  m-1\right)  }a_{1:m-1}^{\gamma}+a_{m}\sum
\limits_{\gamma\leq\mathbf{1,}\left\vert \gamma\right\vert =m-2}c_{\gamma
}^{\left(  m-1\right)  }a_{1:m-1}^{\gamma}+a_{m}\sum\limits_{\gamma
\leq\mathbf{1,}\left\vert \gamma\right\vert =m-1}c_{\gamma}^{\left(
m-1\right)  }a_{1:m-1}^{\gamma},\nonumber
\end{align}

i.e.
\begin{subequations}
\label{av580}%
\begin{multline*}
\frac{m+1}{2^{m}}+\sum_{k=1}^{m}\sum\limits_{\substack{\beta\leq\mathbf{1}
\\\left\vert \beta\right\vert =k}}c_{\beta}^{\left(  m\right)  }a^{\beta}\\
=\frac{m+1}{2^{m}}+\left(  \sum\limits_{\substack{\gamma\leq\mathbf{1}
\\\left\vert \gamma\right\vert =1}}c_{\gamma}^{\left(  m-1\right)  }%
a_{1:m-1}^{\gamma}-\frac{1}{4}\sum\limits_{\substack{\delta\leq\mathbf{1}
\\\left\vert \delta\right\vert =1}}c_{\delta}^{\left(  m-2\right)  }%
a_{1:m-2}^{\delta}+\frac{m}{2^{m-1}}a_{m}\right)  +\\
+\sum_{k=2}^{m-2}\left(  \sum\limits_{\substack{\gamma\leq\mathbf{1}
\\\left\vert \gamma\right\vert =k}}c_{\gamma}^{\left(  m-1\right)  }%
a_{1:m-1}^{\gamma}-\frac{1}{4}\sum\limits_{\substack{\delta\leq\mathbf{1}
\\\left\vert \delta\right\vert =k}}c_{\delta}^{\left(  m-2\right)  }%
a_{1:m-2}^{\delta}+a_{m}\sum\limits_{\substack{\gamma\leq\mathbf{1}
\\\left\vert \gamma\right\vert =k-1}}c_{\gamma}^{\left(  m-1\right)
}a_{1:m-1}^{\gamma}\right)  +\\
+\sum\limits_{\substack{\gamma\leq\mathbf{1} \\\left\vert \gamma\right\vert
=m-1}}c_{\gamma}^{\left(  m-1\right)  }a_{1:m-1}^{\gamma}+a_{m}\sum
\limits_{\substack{\gamma\leq\mathbf{1} \\\left\vert \gamma\right\vert
=m-2}}c_{\gamma}^{\left(  m-1\right)  }a_{1:m-1}^{\gamma}+\\
+a_{m}\sum\limits_{\substack{\gamma\leq\mathbf{1} \\\left\vert \gamma
\right\vert =m-1}}c_{\gamma}^{\left(  m-1\right)  }a_{1:m-1}^{\gamma},\quad
m\geq4,
\end{multline*}

and by replacing $a$ by $a/s$ in the last equation, where $s\in\mathbb{R}$,
and matching powers of $s$ on both sides of this equation we find that we must
have\ the set of equations
\end{subequations}
\begin{equation}
\sum\limits_{\beta\leq\mathbf{1},\left\vert \beta\right\vert =1}c_{\beta
}^{\left(  m\right)  }a^{\beta}=\sum\limits_{\gamma\leq\mathbf{1},\left\vert
\gamma\right\vert =1}c_{\gamma}^{\left(  m-1\right)  }a^{\gamma}-\frac{1}%
{4}\sum\limits_{\delta\leq\mathbf{1},\left\vert \delta\right\vert =1}%
c_{\delta}^{\left(  m-2\right)  }a^{\delta}+\frac{m}{2^{m-1}}a_{m},\quad
m\geq4,\label{av592}%
\end{equation}

and when $2\leq k\leq m-2$ and $m\geq4$,%
\begin{equation}
\sum\limits_{\beta\leq\mathbf{1},\left\vert \beta\right\vert =k}c_{\beta
}^{\left(  m\right)  }a^{\beta}=\sum\limits_{\gamma\leq\mathbf{1},\left\vert
\gamma\right\vert =k}c_{\gamma}^{\left(  m-1\right)  }a^{\gamma}-\frac{1}%
{4}\sum\limits_{\delta\leq\mathbf{1},\left\vert \delta\right\vert =k}%
c_{\delta}^{\left(  m-2\right)  }a^{\delta}+a_{m}\sum\limits_{\gamma
\leq\mathbf{1},\left\vert \gamma\right\vert =k-1}c_{\gamma}^{\left(
m-1\right)  }a^{\gamma},\label{av640}%
\end{equation}

and when $k=m-1$,%
\begin{equation}
\sum\limits_{\beta\leq\mathbf{1},\left\vert \beta\right\vert =m-1}c_{\beta
}^{\left(  m\right)  }a^{\beta}=\sum\limits_{\gamma\leq\mathbf{1},\left\vert
\gamma\right\vert =m-1}c_{\gamma}^{\left(  m-1\right)  }a^{\gamma}+a_{m}%
\sum\limits_{\gamma\leq\mathbf{1},\left\vert \gamma\right\vert =m-2}c_{\gamma
}^{\left(  m-1\right)  }a^{\gamma},\quad m\geq4,\label{av593}%
\end{equation}

and when $k=m$,%
\begin{equation}
\sum\limits_{\beta\leq\mathbf{1},\left\vert \beta\right\vert =m}c_{\beta
}^{\left(  m\right)  }a^{\beta}=a_{m}\sum\limits_{\gamma\leq\mathbf{1}%
,\left\vert \gamma\right\vert =m-1}c_{\gamma}^{\left(  m-1\right)  }a^{\gamma
},\quad m\geq4.\label{av594}%
\end{equation}

\fbox{We start with the simplest case \ref{av594}} which can be written%
\[
c_{\mathbf{1}}^{\left(  m\right)  }a_{1}a_{2}\ldots a_{m}=c_{\mathbf{1}%
}^{\left(  m-1\right)  }a_{1}a_{2}\ldots a_{m-1}a_{m},\quad m\geq4,
\]

so that from \ref{av591}, $c_{\mathbf{1}}^{\left(  m\right)  }=c_{\mathbf{1}%
}^{\left(  m-1\right)  }=\ldots=c_{\mathbf{1}}^{\left(  3\right)
}=c_{\mathbf{1}}^{\left(  2\right)  }=c_{\mathbf{1}}^{\left(  1\right)  }=1$
i.e.
\begin{equation}
c_{\mathbf{1}}^{\left(  m\right)  }=1,\quad m\geq1.\label{av595}%
\end{equation}

\fbox{We next consider the equations \ref{av593}} i.e. for $m\geq4$,%
\begin{align*}
\sum\limits_{\beta\leq\mathbf{1},\left\vert \beta\right\vert =m-1}c_{\beta
}^{\left(  m\right)  }a^{\beta}  & =\sum\limits_{\gamma\leq\mathbf{1}%
,\left\vert \gamma\right\vert =m-1}c_{\gamma}^{\left(  m-1\right)  }a^{\gamma
}+a_{m}\sum\limits_{\gamma\leq\mathbf{1},\left\vert \gamma\right\vert
=m-2}c_{\gamma}^{\left(  m-1\right)  }a_{1:m-1}^{\gamma}\\
& =c_{\mathbf{1}}^{\left(  m-1\right)  }a_{1}a_{2}\ldots a_{m-1}+a_{m}%
\sum\limits_{\gamma\leq\mathbf{1},\left\vert \gamma\right\vert =m-2}c_{\gamma
}^{\left(  m-1\right)  }a_{1:m-1}^{\gamma}\\
& =a_{1}a_{2}\ldots a_{m-1}+a_{m}\sum\limits_{\gamma\leq\mathbf{1},\left\vert
\gamma\right\vert =m-2}c_{\gamma}^{\left(  m-1\right)  }a_{1:m-1}^{\gamma}.
\end{align*}

But from \ref{av589}, \ref{av588} and \ref{av563},%
\[
c_{\beta}^{\left(  m\right)  }=1\text{ }when\text{ }\mathbf{0}\leq\beta
\leq\mathbf{1},\text{ }\left\vert \beta\right\vert =m-1\text{ }and\text{
}4\leq m\leq5,
\]

so the multi-index identity%
\[
\sum\limits_{\beta\leq\mathbf{1},\left\vert \beta\right\vert =m-1}a^{\beta
}=a_{1}a_{2}\ldots a_{m-1}+a_{m}\sum\limits_{\gamma\leq\mathbf{1},\left\vert
\gamma\right\vert =m-2}a_{1:m-1}^{\gamma},\quad m\geq4,
\]

implies%
\begin{equation}
c_{\beta}^{\left(  m\right)  }=1\text{ }when\text{ }\mathbf{0}\leq\beta
\leq\mathbf{1},\text{ }\left\vert \beta\right\vert =m-1\text{ }and\text{
}m\geq4.\label{av613}%
\end{equation}

But from \ref{av591} $c_{0,1,1}^{\left(  3\right)  }=c_{1,1,0}^{\left(
3\right)  }=c_{1,0,1}^{\left(  3\right)  }=1$ and $c_{0,1}^{\left(  2\right)
}=c_{1,0}^{\left(  2\right)  }=1$ and $c_{0}^{\left(  1\right)  }=1$ so
\ref{av613} can be extended to%
\begin{equation}
c_{\beta}^{\left(  m\right)  }=1\text{ }when\text{ }\mathbf{0}\leq\beta
\leq\mathbf{1},\text{ }\left\vert \beta\right\vert =m-1\text{ }and\text{
}m\geq1.\label{av619}%
\end{equation}

\fbox{Next consider \ref{av592}} These equations can be written%
\begin{align*}
c_{\mathbf{i}_{1}}^{\left(  m\right)  } &  a_{1}+c_{\mathbf{i}_{2}}^{\left(
m\right)  }a_{2}+\ldots+c_{\mathbf{i}_{m}}^{\left(  m\right)  }a_{m}\\
&  =\left(  c_{\mathbf{i}_{1}}^{\left(  m-1\right)  }a_{1}+c_{\mathbf{i}_{2}%
}^{\left(  m-1\right)  }a_{2}+\ldots+c_{\mathbf{i}_{m-1}}^{\left(  m-1\right)
}a_{m-1}\right)  -\\
&  \qquad\qquad-\frac{1}{4}\left(  c_{\mathbf{i}_{1}}^{\left(  m-2\right)
}a_{1}+c_{\mathbf{i}_{2}}^{\left(  m-2\right)  }a_{2}+\ldots+c_{\mathbf{i}%
_{m-2}}^{\left(  m-2\right)  }a_{m-2}\right)  +\frac{m}{2^{m-1}}a_{m}\\
&  =\left(  c_{\mathbf{i}_{1}}^{\left(  m-1\right)  }-\frac{1}{4}%
c_{\mathbf{i}_{1}}^{\left(  m-2\right)  }\right)  a_{1}+\left(  c_{\mathbf{i}%
_{2}}^{\left(  m-1\right)  }-\frac{1}{4}c_{\mathbf{i}_{2}}^{\left(
m-2\right)  }\right)  a_{2}+\ldots+\left(  c_{\mathbf{i}_{m-2}}^{\left(
m-1\right)  }-\frac{1}{4}c_{\mathbf{i}_{m-2}}^{\left(  m-2\right)  }\right)
a_{m-2}+\\
&  \qquad\qquad+c_{\mathbf{i}_{m-1}}^{\left(  m-1\right)  }a_{m-1}+\frac
{m}{2^{m-1}}a_{m},
\end{align*}

which implies that for $m\geq4$,%
\begin{align*}
c_{\mathbf{i}_{1}}^{\left(  m\right)  } &  =c_{\mathbf{i}_{1}}^{\left(
m-1\right)  }-\frac{1}{4}c_{\mathbf{i}_{1}}^{\left(  m-2\right)  },\\
c_{\mathbf{i}_{2}}^{\left(  m\right)  } &  =c_{\mathbf{i}_{2}}^{\left(
m-1\right)  }-\frac{1}{4}c_{\mathbf{i}_{2}}^{\left(  m-2\right)  },\\
&  \vdots\\
c_{\mathbf{i}_{m-2}}^{\left(  m\right)  } &  =c_{\mathbf{i}_{m-2}}^{\left(
m-1\right)  }-\frac{1}{4}c_{\mathbf{i}_{m-2}}^{\left(  m-2\right)  },\\
c_{\mathbf{i}_{m-1}}^{\left(  m\right)  } &  =c_{\mathbf{i}_{m-1}}^{\left(
m-1\right)  },\\
c_{\mathbf{i}_{m}}^{\left(  m\right)  } &  =\frac{m}{2^{m-1}}.
\end{align*}

or more compactly%
\begin{equation}
c_{\mathbf{i}_{k}}^{\left(  m\right)  }=\left\{
\begin{array}
[c]{ll}%
c_{\mathbf{i}_{k}}^{\left(  m-1\right)  }-\frac{1}{4}c_{\mathbf{i}_{k}%
}^{\left(  m-2\right)  }, & k\leq m-2,\\
c_{\mathbf{i}_{m-1}}^{\left(  m-1\right)  }, & k=m-1,\\
\frac{m}{2^{m-1}}, & k=m,
\end{array}
\right.  \qquad m\geq4.\label{av597}%
\end{equation}

But from \ref{av591},%
\begin{equation}%
\begin{array}
[c]{l}%
c_{0}^{\left(  1\right)  }=c_{\mathbf{i}_{1}}^{\left(  1\right)  }=1,\\
c_{\mathbf{0}}^{\left(  2\right)  }=\frac{3}{4},\quad c_{\mathbf{i}_{1}%
}^{\left(  2\right)  }=c_{\mathbf{i}_{2}}^{\left(  2\right)  }=1,\quad
c_{\mathbf{1}}^{\left(  2\right)  }=1,\\
c_{\mathbf{0}}^{\left(  3\right)  }=\frac{1}{2},\quad c_{\mathbf{i}_{1}%
}^{\left(  3\right)  }=c_{\mathbf{i}_{3}}^{\left(  3\right)  }=\frac{3}%
{4},\text{ }c_{\mathbf{i}_{2}}^{\left(  3\right)  }=1,\quad c_{1,1,0}^{\left(
3\right)  }=c_{1,0,1}^{\left(  3\right)  }=c_{0,1,1}^{\left(  3\right)
}=1,\quad c_{\mathbf{1}}^{\left(  3\right)  }=1,
\end{array}
\label{av932}%
\end{equation}

so by inspection \ref{av597} can be extended to:\medskip

\fbox{$m=1$} Only $c_{\mathbf{i}_{m}}^{\left(  m\right)  }=\frac{m}{2^{m-1}}$
holds since $c_{\mathbf{i}_{1}}^{\left(  1\right)  }=\frac{1}{2^{1-1}}=1$.

\fbox{$m=2$} $c_{\mathbf{i}_{m-1}}^{\left(  m\right)  }=c_{\mathbf{i}_{m-1}%
}^{\left(  m-1\right)  }$ holds since $c_{\mathbf{i}_{1}}^{\left(  2\right)
}=c_{\mathbf{i}_{1}}^{\left(  1\right)  }=1$ and $c_{\mathbf{i}_{m}}^{\left(
m\right)  }=\frac{m}{2^{m-1}}$ holds since $c_{\mathbf{i}_{2}}^{\left(
2\right)  }=\frac{2}{2}=1$.

\fbox{$m=3$} All equations hold: we have $k=1$ and $c_{\mathbf{i}_{1}%
}^{\left(  3\right)  }=c_{\mathbf{i}_{1}}^{\left(  2\right)  }-\frac{1}%
{4}c_{\mathbf{i}_{1}}^{\left(  1\right)  }=\frac{3}{4}$, $c_{\mathbf{i}_{2}%
}^{\left(  3\right)  }=c_{\mathbf{i}_{2}}^{\left(  2\right)  }=1$ and
$c_{\mathbf{i}_{3}}^{\left(  3\right)  }=\frac{3}{4}$.\medskip

Thus \ref{av597} can be extended to:%
\[
c_{\mathbf{i}_{k}}^{\left(  m\right)  }=\left\{
\begin{array}
[c]{lll}%
c_{\mathbf{i}_{k}}^{\left(  m-1\right)  }-\frac{1}{4}c_{\mathbf{i}_{k}%
}^{\left(  m-2\right)  }, & k\leq m-2, & m\geq3,\\
c_{\mathbf{i}_{m-1}}^{\left(  m-1\right)  }, & k=m-1, & m\geq2,\\
\frac{m}{2^{m-1}}, & k=m, & m\geq1,
\end{array}
\right.
\]

which further simplifies to%
\begin{equation}
c_{\mathbf{i}_{k}}^{\left(  m\right)  }=\left\{
\begin{array}
[c]{lll}%
c_{\mathbf{i}_{k}}^{\left(  m-1\right)  }-\frac{1}{4}c_{\mathbf{i}_{k}%
}^{\left(  m-2\right)  }, & k\leq m-2, & m\geq3,\\
\frac{m-1}{2^{m-2}}, & k=m-1, & m\geq2,\\
\frac{m}{2^{m-1}}, & k=m, & m\geq1.
\end{array}
\right. \label{av644}%
\end{equation}

For each $k$, \ref{av644} can be interpreted as a second order difference
equation for $c_{\mathbf{i}_{k}}^{\left(  m\right)  }$ in terms of $m$. In
fact, for given $k\geq1$,%
\begin{equation}%
\begin{array}
[c]{l}%
c_{\mathbf{i}_{k}}^{\left(  m\right)  }=c_{\mathbf{i}_{k}}^{\left(
m-1\right)  }-\frac{1}{4}c_{\mathbf{i}_{k}}^{\left(  m-2\right)  },\quad m\geq
k+2,\\
c_{\mathbf{i}_{k}}^{\left(  k+1\right)  }=c_{\mathbf{i}_{k}}^{\left(
k\right)  }=\frac{k}{2^{k-1}}.
\end{array}
\label{av612}%
\end{equation}

The auxiliary equation is $x^{2}-x+1/4$ which has the single zero $1/2$ so the
solution has the form%
\[
c_{\mathbf{i}_{k}}^{\left(  m\right)  }=\frac{1}{2^{m}}\left(  A_{k}%
+B_{k}m\right)  .
\]

The initial conditions imply%
\[
\frac{k}{2^{k-1}}=\frac{1}{2^{k}}\left(  A_{k}+B_{k}k\right)  =\frac
{1}{2^{k+1}}\left(  A_{k}+B_{k}\left(  k+1\right)  \right)  ,
\]

so that $2\left(  A_{k}+B_{k}k\right)  =A_{k}+B_{k}\left(  k+1\right)  $ and
$2k=A_{k}+B_{k}k$ i.e. $B_{k}=2k$ and $A_{k}=2k\left(  1-k\right)  $ and hence%
\begin{equation}
c_{\mathbf{i}_{k}}^{\left(  m\right)  }=\frac{1}{2^{m}}\left(  2k\left(
1-k\right)  +2km\right)  =\frac{k}{2^{m-1}}\left(  m+1-k\right)
.\label{av820}%
\end{equation}

Regarding upper and lower\ bounds for $c_{\mathbf{i}_{k}}^{\left(  m\right)
}$ for each $m$. Since $x\left(  m+1-x\right)  $ has a maximum at
$x=\frac{m+1}{2} $ it follows that when $m$ is odd,$\ \operatorname*{argmax}%
\limits_{k}c_{\mathbf{i}_{k}}^{\left(  m\right)  }=\frac{m+1}{2}$ and
$\max\limits_{k}c_{\mathbf{i}_{k}}^{\left(  m\right)  }=\frac{\frac{m+1}{2}%
}{2^{m-1}}\left(  m+1-\frac{m+1}{2}\right)  =\frac{\left(  m+1\right)  ^{2}%
}{2^{m+1}}$. Also, if $m$ is even, $\operatorname*{argmax}\limits_{k}%
c_{\mathbf{i}_{k}}^{\left(  m\right)  }=\left\{  \frac{m}{2},\frac{m}%
{2}+1\right\}  $ and $\max\limits_{k}c_{\mathbf{i}_{k}}^{\left(  m\right)
}=\frac{m/2}{2^{m-1}}\left(  m+1-\frac{m}{2}\right)  =\frac{m}{2^{m}}\left(
\frac{m}{2}+1\right)  $.

Further, $\operatorname*{argmin}\limits_{k}c_{\mathbf{i}_{k}}^{\left(
m\right)  }=\left\{  1,m\right\}  $ and $\min\limits_{k}c_{\mathbf{i}_{k}%
}^{\left(  m\right)  }=\frac{m}{2^{m-1}}$. In summary%
\begin{equation}%
\begin{array}
[c]{ll}%
\max\limits_{k}c_{\mathbf{i}_{k}}^{\left(  m\right)  }=\frac{m}{2^{m}}\left(
\frac{m}{2}+1\right)  , & \operatorname*{argmax}\limits_{k}c_{\mathbf{i}_{k}%
}^{\left(  m\right)  }=\left\{
\begin{array}
[c]{cc}%
\frac{m+1}{2}, & m\text{ }odd,\\
\frac{\left(  m+1\right)  ^{2}}{2^{m+1}}, & m\text{ }even,
\end{array}
\right. \\
\min\limits_{k}c_{\mathbf{i}_{k}}^{\left(  m\right)  }=\frac{m}{2^{m-1}}, &
\operatorname*{argmin}\limits_{k}c_{\mathbf{i}_{k}}^{\left(  m\right)
}=\left\{  1,m\right\}  ,
\end{array}
\label{av615}%
\end{equation}

so that%
\begin{equation}
\frac{m}{2^{m-1}}\leq c_{\mathbf{i}_{k}}^{\left(  m\right)  }\leq\frac
{m}{2^{m}}\left(  \frac{m}{2}+1\right)  ,\quad1\leq k\leq m,\text{ }%
m\geq1.\label{av618}%
\end{equation}

So we have considered the cases $c_{\mathbf{\beta}}^{\left(  m\right)  }$
where $\left\vert \beta\right\vert =1$(\ref{av598}), $\left\vert
\beta\right\vert =m-1$(\ref{av619}) and $\left\vert \beta\right\vert
=m$(\ref{av595}).\medskip

\fbox{Now to consider the most complex equation \ref{av640}} i.e. when $2\leq
k\leq m-2$ and $m\geq4$,%
\[
\sum\limits_{\beta\leq\mathbf{1},\left\vert \beta\right\vert =k}c_{\beta
}^{\left(  m\right)  }a^{\beta}=\sum\limits_{\gamma\leq\mathbf{1},\left\vert
\gamma\right\vert =k}c_{\gamma}^{\left(  m-1\right)  }a^{\gamma}-\frac{1}%
{4}\sum\limits_{\delta\leq\mathbf{1},\left\vert \delta\right\vert =k}%
c_{\delta}^{\left(  m-2\right)  }a^{\delta}+a_{m}\sum\limits_{\gamma
\leq\mathbf{1},\left\vert \gamma\right\vert =k-1}c_{\gamma}^{\left(
m-1\right)  }a^{\gamma}.
\]

\begin{definition}
\label{vDef_alt_notation_e}\textbf{Alternative notation }$e_{\gamma}^{\left(
m\right)  }$: We will use the following alternative notation which is very
useful for examples: since in $c_{\beta}^{\left(  m\right)  }$ we always have
$0\leq\beta\leq1_{m}$ an alternative and more compact notation $e_{\gamma
}^{\left(  m\right)  }=c_{\beta}^{\left(  m\right)  }$ uses the subscripts of
$e_{\gamma}^{\left(  m\right)  }$ to give the locations of the $1$s in $\beta$
e.g. in \ref{av588},%
\begin{equation}
e_{1,2}^{\left(  4\right)  }=\frac{1}{2},\text{ }e_{1,3}^{\left(  4\right)
}=1,\text{ }e_{1,4}^{\left(  4\right)  }=\frac{3}{4},\text{ }e_{2,3}^{\left(
4\right)  }=1,\text{ }e_{2,4}^{\left(  4\right)  }=1,\text{ }e_{3,4}^{\left(
4\right)  }=\frac{3}{4},\label{av623}%
\end{equation}

where $e_{1,2}^{\left(  4\right)  }=c_{\mathbf{i}_{1}+\mathbf{i}_{2}}^{\left(
4\right)  }$, $e_{1,3}^{\left(  4\right)  }=c_{\mathbf{i}_{1}+\mathbf{i}_{3}%
}^{\left(  4\right)  }$ etc. and $e_{m,n}^{\left(  4\right)  }=c_{\mathbf{i}%
_{m}+\mathbf{i}_{n}}^{\left(  4\right)  }$ and in general
\begin{equation}
e_{k_{1},k_{2},\ldots,k_{n}}^{\left(  m\right)  }=c_{\mathbf{i}_{k_{1}%
}+\mathbf{i}_{k_{2}}+\ldots+\mathbf{i}_{k_{n}}}^{\left(  m\right)
}.\label{av625}%
\end{equation}

\end{definition}

From \ref{av563},%
\begin{equation}%
\begin{array}
[c]{l}%
e_{1,2}^{\left(  5\right)  }=\frac{1}{2},\text{\quad}e_{1,3}^{\left(
5\right)  }=\frac{3}{4},\text{\quad}e_{1,4}^{\left(  5\right)  }=\frac{3}%
{4},\text{\quad}e_{1,5}^{\left(  5\right)  }=\frac{1}{2},\\
e_{2,3}^{\left(  5\right)  }=\frac{3}{4},\text{\quad}e_{2,4}^{\left(
5\right)  }=1,\text{\quad}e_{2,5}^{\left(  5\right)  }=\frac{3}{4},\\
e_{3,4}^{\left(  5\right)  }=\frac{3}{4},\text{\quad}e_{3,5}^{\left(
5\right)  }=\frac{3}{4},\\
e_{4,5}^{\left(  5\right)  }=\frac{1}{2}.
\end{array}
\label{av624}%
\end{equation}

\begin{example}
From equation\textbf{\ }\ref{av640}: when $m=4$ and $k=2$ equation \ref{av640}
becomes%
\[
\sum\limits_{\beta\leq\mathbf{1},\left\vert \beta\right\vert =2}c_{\beta
}^{\left(  4\right)  }a^{\beta}=\sum\limits_{\gamma\leq\mathbf{1},\left\vert
\gamma\right\vert =2}c_{\gamma}^{\left(  4-1\right)  }a^{\gamma}+a_{m}%
\sum\limits_{\gamma\leq\mathbf{1},\left\vert \gamma\right\vert =2-1}c_{\gamma
}^{\left(  4-1\right)  }a^{\gamma}-\frac{1}{4}\sum\limits_{\delta
\leq\mathbf{1},\left\vert \delta\right\vert =2}c_{\delta}^{\left(  4-2\right)
}a^{\delta},
\]

i.e.%
\[
\sum\limits_{\beta\leq\mathbf{1},\left\vert \beta\right\vert =2}c_{\beta
}^{\left(  4\right)  }a^{\beta}=\sum\limits_{\gamma\leq\mathbf{1},\left\vert
\gamma\right\vert =2}c_{\gamma}^{\left(  3\right)  }a^{\gamma}+a_{m}%
\sum\limits_{\gamma\leq\mathbf{1},\left\vert \gamma\right\vert =1}c_{\gamma
}^{\left(  3\right)  }a^{\gamma}-\frac{1}{4}\sum\limits_{\delta\leq
\mathbf{1},\left\vert \delta\right\vert =2}c_{\delta}^{\left(  2\right)
}a^{\delta},
\]

or in notation \ref{av625},%
\begin{align*}
e_{1,2}^{\left(  4\right)  }a_{1}a_{2} &  +e_{1,3}^{\left(  4\right)  }%
a_{1}a_{3}+e_{1,4}^{\left(  4\right)  }a_{1}a_{4}+e_{2,3}^{\left(  4\right)
}a_{2}a_{3}+e_{2,4}^{\left(  4\right)  }a_{2}a_{4}+e_{3,4}^{\left(  4\right)
}a_{3}a_{4}\\
&  =e_{1,2}^{\left(  3\right)  }a_{1}a_{2}+e_{1,3}^{\left(  3\right)  }%
a_{1}a_{3}+e_{2,3}^{\left(  3\right)  }a_{2}a_{3}+a_{4}\left(  e_{1}^{\left(
3\right)  }a_{1}+e_{2}^{\left(  3\right)  }a_{2}+e_{3}^{\left(  3\right)
}a_{3}\right)  -\frac{1}{4}e_{1,2}^{\left(  2\right)  }a_{1}a_{2}\\
&  =e_{1,2}^{\left(  3\right)  }a_{1}a_{2}+e_{1,3}^{\left(  3\right)  }%
a_{1}a_{3}+e_{2,3}^{\left(  3\right)  }a_{2}a_{3}+e_{1}^{\left(  3\right)
}a_{1}a_{4}+e_{2}^{\left(  3\right)  }a_{2}a_{4}+e_{3}^{\left(  3\right)
}a_{3}a_{4}-\frac{1}{4}e_{1,2}^{\left(  2\right)  }a_{1}a_{2}\\
&  =\left(  e_{1,2}^{\left(  3\right)  }-\frac{1}{4}e_{1,2}^{\left(  2\right)
}\right)  a_{1}a_{2}+e_{1,3}^{\left(  3\right)  }a_{1}a_{3}+e_{1}^{\left(
3\right)  }a_{1}a_{4}+e_{2,3}^{\left(  3\right)  }a_{2}a_{3}+e_{2}^{\left(
3\right)  }a_{2}a_{4}+e_{3}^{\left(  3\right)  }a_{3}a_{4},
\end{align*}

Noting the definition \ref{av625}, by equating monomials of $a$ and then using
\ref{av589} we get%
\begin{equation}%
\begin{array}
[c]{l}%
e_{1,2}^{\left(  4\right)  }=e_{1,2}^{\left(  3\right)  }-\frac{1}{4}%
e_{1,2}^{\left(  2\right)  }=3/4,\\
e_{1,3}^{\left(  4\right)  }=e_{1,3}^{\left(  3\right)  }=1,\\
e_{1,4}^{\left(  4\right)  }=e_{1}^{\left(  3\right)  }=3/4,\\
e_{2,3}^{\left(  4\right)  }=e_{2,3}^{\left(  3\right)  }=1,\\
e_{2,4}^{\left(  4\right)  }=e_{2}^{\left(  3\right)  }=1,\\
e_{3,4}^{\left(  4\right)  }=e_{3}^{\left(  3\right)  }=3/4,
\end{array}
\label{av626}%
\end{equation}

which matches \ref{av588}.
\end{example}

\begin{example}
\textbf{Regarding equation \ref{av640}}: when $m=5$ and $k=3$ equation
\ref{av640} becomes%
\[
\sum\limits_{\beta\leq\mathbf{1},\left\vert \beta\right\vert =3}c_{\beta
}^{\left(  5\right)  }a^{\beta}=\sum\limits_{\gamma\leq\mathbf{1},\left\vert
\gamma\right\vert =3}c_{\gamma}^{\left(  4\right)  }a_{1:4}^{\gamma}+a_{5}%
\sum\limits_{\gamma\leq\mathbf{1},\left\vert \gamma\right\vert =2}c_{\gamma
}^{\left(  4\right)  }a_{1:4}^{\gamma}-\frac{1}{4}\sum\limits_{\delta
\leq\mathbf{1},\left\vert \delta\right\vert =3}c_{\delta}^{\left(  3\right)
}a_{1:3}^{\delta},
\]

i.e.%
\begin{align*}
&  c_{1,2,3}^{\left(  5\right)  }a_{1}a_{2}a_{3}+c_{1,2,4}^{\left(  5\right)
}a_{1}a_{2}a_{4}+c_{1,2,5}^{\left(  5\right)  }a_{1}a_{2}a_{5}+c_{1,3,4}%
^{\left(  5\right)  }a_{1}a_{3}a_{4}+c_{1,3,5}^{\left(  5\right)  }a_{1}%
a_{3}a_{5}+c_{2,3,4}^{\left(  5\right)  }a_{2}a_{3}a_{4}+\\
&  \qquad\qquad\qquad+c_{2,3,5}^{\left(  5\right)  }a_{2}a_{3}a_{5}%
+c_{2,4,5}^{\left(  5\right)  }a_{2}a_{4}a_{5}+c_{3,4,5}^{\left(  5\right)
}a_{3}a_{4}a_{5}\\
&  =\left(  c_{1,2,3}^{\left(  4\right)  }a_{1}a_{2}a_{3}+c_{1,2,4}^{\left(
4\right)  }a_{1}a_{2}a_{4}+c_{1,3,4}^{\left(  4\right)  }a_{1}a_{3}%
a_{4}+c_{2,3,4}^{\left(  4\right)  }a_{2}a_{3}a_{4}\right)  +\\
&  \qquad+a_{5}\left(  c_{1,2}^{\left(  4\right)  }a_{1}a_{2}+c_{1,3}^{\left(
4\right)  }a_{1}a_{3}+c_{1,4}^{\left(  4\right)  }a_{1}a_{4}+c_{2,3}^{\left(
4\right)  }a_{2}a_{3}+c_{2,4}^{\left(  4\right)  }a_{2}a_{4}+c_{3,4}^{\left(
4\right)  }a_{3}a_{4}\right)  -\frac{1}{4}c_{1,2,3}^{\left(  3\right)  }%
a_{1}a_{2}a_{3}\\
&  =c_{1,2,3}^{\left(  4\right)  }a_{1}a_{2}a_{3}+c_{1,2,4}^{\left(  4\right)
}a_{1}a_{2}a_{4}+c_{1,3,4}^{\left(  4\right)  }a_{1}a_{3}a_{4}+c_{2,3,4}%
^{\left(  4\right)  }a_{2}a_{3}a_{4}+\\
&  \qquad+a_{5}\left(  c_{1,2}^{\left(  4\right)  }a_{1}a_{2}+c_{1,3}^{\left(
4\right)  }a_{1}a_{3}+c_{1,4}^{\left(  4\right)  }a_{1}a_{4}+c_{2,3}^{\left(
4\right)  }a_{2}a_{3}+c_{2,4}^{\left(  4\right)  }a_{2}a_{4}+c_{3,4}^{\left(
4\right)  }a_{3}a_{4}\right)  -\frac{1}{4}c_{1,2,3}^{\left(  3\right)  }%
a_{1}a_{2}a_{3}\\
&  =c_{1,2,3}^{\left(  4\right)  }a_{1}a_{2}a_{3}+c_{1,2,4}^{\left(  4\right)
}a_{1}a_{2}a_{4}+c_{1,3,4}^{\left(  4\right)  }a_{1}a_{3}a_{4}+c_{2,3,4}%
^{\left(  4\right)  }a_{2}a_{3}a_{4}+\\
&  \qquad+\left(  c_{1,2}^{\left(  4\right)  }a_{1}a_{2}a_{5}+c_{1,3}^{\left(
4\right)  }a_{1}a_{3}a_{5}+c_{1,4}^{\left(  4\right)  }a_{1}a_{4}a_{5}%
+c_{2,3}^{\left(  4\right)  }a_{2}a_{3}a_{5}+c_{2,4}^{\left(  4\right)  }%
a_{2}a_{4}a_{5}+c_{3,4}^{\left(  4\right)  }a_{3}a_{4}a_{5}\right)  -\\
&  \qquad-\frac{1}{4}c_{1,2,3}^{\left(  3\right)  }a_{1}a_{2}a_{3}\\
&  =\left(  c_{1,2,3}^{\left(  4\right)  }-\frac{1}{4}c_{1,2,3}^{\left(
3\right)  }\right)  a_{1}a_{2}a_{3}+c_{1,2,4}^{\left(  4\right)  }a_{1}%
a_{2}a_{4}+c_{1,3,4}^{\left(  4\right)  }a_{1}a_{3}a_{4}+c_{2,3,4}^{\left(
4\right)  }a_{2}a_{3}a_{4}+\\
&  \qquad+\left(  c_{1,2}^{\left(  4\right)  }a_{1}a_{2}a_{5}+c_{1,3}^{\left(
4\right)  }a_{1}a_{3}a_{5}+c_{1,4}^{\left(  4\right)  }a_{1}a_{4}a_{5}%
+c_{2,3}^{\left(  4\right)  }a_{2}a_{3}a_{5}+c_{2,4}^{\left(  4\right)  }%
a_{2}a_{4}a_{5}+c_{3,4}^{\left(  4\right)  }a_{3}a_{4}a_{5}\right) \\
&  =\left(  c_{1,2,3}^{\left(  4\right)  }-\frac{1}{4}c_{1,2,3}^{\left(
3\right)  }\right)  a_{1}a_{2}a_{3}+c_{1,2,4}^{\left(  4\right)  }a_{1}%
a_{2}a_{4}+c_{1,2}^{\left(  4\right)  }a_{1}a_{2}a_{5}+c_{1,3,4}^{\left(
4\right)  }a_{1}a_{3}a_{4}+c_{1,3}^{\left(  4\right)  }a_{1}a_{3}a_{5}+\\
&  \qquad+c_{2,3,4}^{\left(  4\right)  }a_{2}a_{3}a_{4}+c_{2,3}^{\left(
4\right)  }a_{2}a_{3}a_{5}+c_{2,4}^{\left(  4\right)  }a_{2}a_{4}%
a_{5}++c_{3,4}^{\left(  4\right)  }a_{3}a_{4}a_{5}.
\end{align*}

Noting \ref{av625} and equating monomials of $a$ we get%
\begin{align*}
e_{1,2,3}^{\left(  5\right)  }  & =e_{1,2,3}^{\left(  4\right)  }-\frac{1}%
{4}e_{1,2,3}^{\left(  3\right)  },\\
e_{1,2,4}^{\left(  5\right)  }  & =e_{1,2,4}^{\left(  4\right)  },\text{\quad
}e_{1,2,5}^{\left(  5\right)  }=e_{1,2}^{\left(  4\right)  },\text{\quad
}e_{1,3,4}^{\left(  5\right)  }=e_{1,3,4}^{\left(  4\right)  },\\
e_{2,3,4}^{\left(  5\right)  }  & =e_{2,3,4}^{\left(  4\right)  },\text{\quad
}e_{2,3,5}^{\left(  5\right)  }=e_{2,3}^{\left(  4\right)  },\text{\quad
}e_{2,4,5}^{\left(  5\right)  }=e_{2,4}^{\left(  4\right)  },\\
e_{3,4,5}^{\left(  5\right)  }  & =e_{3,4}^{\left(  4\right)  }.
\end{align*}

\end{example}

In general we need an equation which expresses $\sum\limits_{\beta
\leq\mathbf{1},\left\vert \beta\right\vert =k}c_{\beta}^{\left(  m\right)
}a^{\beta}$ in terms of $\sum\limits_{\substack{\gamma\leq\mathbf{1}%
,\left\vert \gamma\right\vert =k-1 \\\gamma\in\mathbb{Z}_{\oplus}^{m-1}
}}c_{\gamma,1}^{\left(  m\right)  }a_{1:m-1}^{\gamma}$ etc. However, using the
simple identity%
\[
\sum\limits_{\substack{\beta\leq\mathbf{1},\left\vert \beta\right\vert =k
\\\beta\in\mathbb{Z}_{\oplus}^{m}}}a^{\beta}=\sum\limits_{\substack{\gamma
\leq\mathbf{1},\left\vert \gamma\right\vert =k \\\gamma\in\mathbb{Z}_{\oplus
}^{m-1}}}a_{1:m-1}^{\gamma}+a_{m}\sum\limits_{\substack{\delta\leq
\mathbf{1},\left\vert \delta\right\vert =k-1 \\\delta\in\mathbb{Z}_{\oplus
}^{m-1}}}a_{1:m-1}^{\delta},
\]

for guidance, we write using \ref{av625}:%
\begin{align*}
&  e_{1,2,3}^{\left(  5\right)  }a_{1}a_{2}a_{3}+e_{1,2,4}^{\left(  5\right)
}a_{1}a_{2}a_{4}+e_{1,2,5}^{\left(  5\right)  }a_{1}a_{2}a_{5}+e_{1,3,4}%
^{\left(  5\right)  }a_{1}a_{3}a_{4}+e_{1,3,5}^{\left(  5\right)  }a_{1}%
a_{3}a_{5}+e_{1,4,5}^{\left(  5\right)  }a_{1}a_{4}a_{5}+\\
&  \qquad\qquad+e_{2,3,4}^{\left(  5\right)  }a_{2}a_{3}a_{4}+e_{2,3,5}%
^{\left(  5\right)  }a_{2}a_{3}a_{5}+e_{2,4,5}^{\left(  5\right)  }a_{2}%
a_{4}a_{5}+e_{3,4,5}^{\left(  5\right)  }a_{3}a_{4}a_{5}\\
&  =e_{1,2,3}^{\left(  5\right)  }a_{1}a_{2}a_{3}+e_{1,2,4}^{\left(  5\right)
}a_{1}a_{2}a_{4}+e_{1,3,4}^{\left(  5\right)  }a_{1}a_{3}a_{4}+e_{2,3,4}%
^{\left(  5\right)  }a_{2}a_{3}a_{4}+e_{1,2,5}^{\left(  5\right)  }a_{1}%
a_{2}a_{5}+e_{1,3,5}^{\left(  5\right)  }a_{1}a_{3}a_{5}+\\
&  \qquad\qquad e_{1,4,5}^{\left(  5\right)  }a_{1}a_{4}a_{5}+e_{2,3,5}%
^{\left(  5\right)  }a_{2}a_{3}a_{5}+e_{2,4,5}^{\left(  5\right)  }a_{2}%
a_{4}a_{5}+e_{3,4,5}^{\left(  5\right)  }a_{3}a_{4}a_{5}\\
&  =e_{1,2,3}^{\left(  5\right)  }a_{1}a_{2}a_{3}+e_{1,2,4}^{\left(  5\right)
}a_{1}a_{2}a_{4}+e_{1,3,4}^{\left(  5\right)  }a_{1}a_{3}a_{4}+e_{2,3,4}%
^{\left(  5\right)  }a_{2}a_{3}a_{4}+\\
&  \qquad\qquad+a_{5}\left(  e_{1,2,5}^{\left(  5\right)  }a_{1}%
a_{2}+e_{1,3,5}^{\left(  5\right)  }a_{1}a_{3}+e_{1,4,5}^{\left(  5\right)
}a_{1}a_{4}+e_{2,3,5}^{\left(  5\right)  }a_{2}a_{3}+e_{2,4,5}^{\left(
5\right)  }a_{2}a_{4}+e_{3,4,5}^{\left(  5\right)  }a_{3}a_{4}\right)  ,
\end{align*}

or in $c_{\beta}^{\left(  m\right)  }$ notation:%
\[
\sum\limits_{\beta\leq\mathbf{1},\left\vert \beta\right\vert =3}c_{\beta
}^{\left(  5\right)  }a^{\beta}=\sum\limits_{\substack{\beta\leq
\mathbf{1},\left\vert \beta\right\vert =3 \\\beta\in\mathbb{Z}_{\oplus}^{4}%
}}c_{\beta,0}^{\left(  5\right)  }a_{1:4}^{\beta}+a_{5}\sum
\limits_{\substack{\beta\leq\mathbf{1},\left\vert \beta\right\vert =2
\\\beta\in\mathbb{Z}_{\oplus}^{4}}}c_{\beta,1}^{\left(  5\right)  }%
a_{1:4}^{\beta},
\]

and in general: for $2\leq k\leq m-1$ and $m\geq4$:%
\begin{equation}
\sum\limits_{\beta\leq\mathbf{1},\left\vert \beta\right\vert =k}c_{\beta
}^{\left(  m\right)  }a^{\beta}=\sum\limits_{\substack{\gamma\leq
\mathbf{1},\left\vert \gamma\right\vert =k \\\gamma\in\mathbb{Z}_{\oplus
}^{m-1}}}c_{\gamma,0}^{\left(  m\right)  }a_{1:m-1}^{\gamma}+a_{m}%
\sum\limits_{\substack{\delta\leq\mathbf{1},\left\vert \delta\right\vert =k-1
\\\delta\in\mathbb{Z}_{\oplus}^{m-1}}}c_{\delta,1}^{\left(  m\right)
}a_{1:m-1}^{\delta},\label{av616}%
\end{equation}

With $m$ replaced by $m-1$ \ref{av616} becomes: for $2\leq k\leq m-2,$
$m\geq5$:
\begin{equation}
\sum\limits_{\gamma\leq\mathbf{1},\left\vert \gamma\right\vert =k}c_{\gamma
}^{\left(  m-1\right)  }a_{1:m-1}^{\gamma}=\sum\limits_{\substack{\beta
\leq\mathbf{1},\left\vert \beta\right\vert =k \\\beta\in\mathbb{Z}_{\oplus
}^{m-2}}}c_{\beta,0}^{\left(  m-1\right)  }a_{1:m-2}^{\beta}+a_{m-1}%
\sum\limits_{\substack{\beta\leq\mathbf{1},\left\vert \beta\right\vert =k-1
\\\beta\in\mathbb{Z}_{\oplus}^{m-2}}}c_{\beta,1}^{\left(  m-1\right)
}a_{1:m-2}^{\beta}.\label{av617}%
\end{equation}

Recall that equation \ref{av640} is: for $2\leq k\leq m-2,$ $m\geq4$:%
\begin{multline*}
\sum\limits_{\beta\leq\mathbf{1},\left\vert \beta\right\vert =k}c_{\beta
}^{\left(  m\right)  }a^{\beta}\\
=\sum\limits_{\gamma\leq\mathbf{1},\left\vert \gamma\right\vert =k}c_{\gamma
}^{\left(  m-1\right)  }a_{1:m-1}^{\gamma}-\frac{1}{4}\sum\limits_{\delta
\leq\mathbf{1},\left\vert \delta\right\vert =k}c_{\delta}^{\left(  m-2\right)
}a_{1:m-2}^{\delta}+a_{m}\sum\limits_{\gamma\leq\mathbf{1},\left\vert
\gamma\right\vert =k-1}c_{\gamma}^{\left(  m-1\right)  }a_{1:m-1}^{\gamma}.
\end{multline*}

Substituting \ref{av616} into the LHS of \ref{av640} yields: for $2\leq k\leq
m-2,$ $m\geq4$:
\begin{multline*}
\sum\limits_{\substack{\beta\leq\mathbf{1},\left\vert \beta\right\vert =k
\\\beta\in\mathbb{Z}_{\oplus}^{m-1}}}c_{\beta,0}^{\left(  m\right)  }%
a_{1:m-1}^{\beta}+a_{m}\sum\limits_{\substack{\beta\leq\mathbf{1},\left\vert
\beta\right\vert =k-1 \\\beta\in\mathbb{Z}_{\oplus}^{m-1}}}c_{\beta
,1}^{\left(  m\right)  }a_{m-1}^{\beta}\\
=\sum\limits_{\substack{\gamma\leq\mathbf{1} \\\left\vert \gamma\right\vert
=k}}c_{\gamma}^{\left(  m-1\right)  }a_{1:m-1}^{\gamma}-\frac{1}{4}%
\sum\limits_{\substack{\delta\leq\mathbf{1} \\\left\vert \delta\right\vert
=k}}c_{\delta}^{\left(  m-2\right)  }a_{1:m-2}^{\delta}+a_{m}\sum
\limits_{\substack{\gamma\leq\mathbf{1} \\\left\vert \gamma\right\vert
=k-1}}c_{\gamma}^{\left(  m-1\right)  }a_{1:m-1}^{\gamma},
\end{multline*}

and matching powers of $a$ implies: for $2\leq k\leq m-2,$ $m\geq4$:
\begin{align}
\sum\limits_{\substack{\beta\leq\mathbf{1},\left\vert \beta\right\vert =k-1
\\\beta\in\mathbb{Z}_{\oplus}^{m-1}}}c_{\beta,1}^{\left(  m\right)  }%
a_{1:m-1}^{\beta}  & =\sum\limits_{\beta\leq\mathbf{1},\left\vert
\beta\right\vert =k-1}c_{\beta}^{\left(  m-1\right)  }a_{1:m-1}^{\beta
},\label{av620}\\
\sum\limits_{\substack{\beta\leq\mathbf{1},\left\vert \beta\right\vert =k
\\\beta\in\mathbb{Z}_{\oplus}^{m-1}}}c_{\beta,0}^{\left(  m\right)  }%
a_{1:m-1}^{\beta}  & =\sum\limits_{\gamma\leq\mathbf{1},\left\vert
\gamma\right\vert =k}c_{\gamma}^{\left(  m-1\right)  }a_{1:m-1}^{\gamma}%
-\frac{1}{4}\sum\limits_{\delta\leq\mathbf{1},\left\vert \delta\right\vert
=k}c_{\delta}^{\left(  m-2\right)  }a_{1:m-2}^{\delta},\label{av600}%
\end{align}

The first equations \ref{av620} imply that%
\[
c_{\beta,1}^{\left(  m\right)  }=c_{\beta}^{\left(  m-1\right)  }%
,\quad\left\{
\begin{array}
[c]{l}%
\beta\in\mathbb{Z}_{\oplus}^{m-2},\\
\beta\leq\mathbf{1},\\
\left\vert \beta\right\vert =k-1\\
2\leq k\leq m-2,\\
m\geq4,
\end{array}
\right.
\]

which can be written more simply as%
\begin{equation}
c_{\beta,1}^{\left(  m\right)  }=c_{\beta}^{\left(  m-1\right)  }%
,\quad\left\{
\begin{array}
[c]{l}%
\beta\in\mathbb{Z}_{\oplus}^{m-2},\\
\beta\leq\mathbf{1},\\
1\leq\left\vert \beta\right\vert \leq m-3,\\
m\geq4,
\end{array}
\right. \label{av621}%
\end{equation}

Now to consider \ref{av600}: this can be rewritten%
\[
\sum\limits_{\substack{\beta\leq\mathbf{1},\left\vert \beta\right\vert =k
\\\beta\in\mathbb{Z}_{\oplus}^{m-1}}}\left(  c_{\beta,0}^{\left(  m\right)
}-c_{\beta}^{\left(  m-1\right)  }\right)  a_{m-1}^{\beta}=-\frac{1}{4}%
\sum\limits_{\delta\leq\mathbf{1},\left\vert \delta\right\vert =k}c_{\delta
}^{\left(  m-2\right)  }a^{\delta},\quad2\leq k\leq m-2,\text{ }m\geq4,
\]

and an application to the left hand side of this equation of the formula
analogous to \ref{av617} yields: for $2\leq k\leq m-2$ and $m\geq4$,%
\begin{multline*}
\sum\limits_{\substack{\beta\leq\mathbf{1},\left\vert \beta\right\vert =k
\\\beta\in\mathbb{Z}_{\oplus}^{m-2}}}\left(  c_{\beta,0,0}^{\left(  m\right)
}-c_{\beta,0}^{\left(  m-1\right)  }\right)  a_{1:m-2}^{\beta}+a_{m-1}%
\sum\limits_{\substack{\gamma\leq\mathbf{1},\left\vert \gamma\right\vert =k-1
\\\gamma\in\mathbb{Z}_{\oplus}^{m-2}}}\left(  c_{\gamma,1,0}^{\left(
m\right)  }-c_{\gamma,1}^{\left(  m-1\right)  }\right)  a_{1:m-2}^{\gamma}\\
=-\frac{1}{4}\sum\limits_{\delta\leq\mathbf{1},\left\vert \delta\right\vert
=k}c_{\delta}^{\left(  m-2\right)  }a_{1:m-2}^{\delta}.
\end{multline*}

Again comparing powers of $a$ we get,%
\[
\left.
\begin{array}
[c]{ll}%
c_{\beta,0,0}^{\left(  m\right)  }=c_{\beta,0}^{\left(  m-1\right)  }-\frac
{1}{4}c_{\beta}^{\left(  m-2\right)  }, & \left\vert \beta\right\vert =k,\\
c_{\beta,1,0}^{\left(  m\right)  }=c_{\beta,1}^{\left(  m-1\right)  }, &
\left\vert \beta\right\vert =k-1,
\end{array}
\right\}  \quad\text{when }\left\{
\begin{array}
[c]{l}%
\beta\in\mathbb{Z}_{\oplus}^{m-2},\\
\beta\leq\mathbf{1},\\
2\leq k\leq m-2,\\
m\geq4,
\end{array}
\right.
\]

or equivalently%
\begin{equation}
\left.
\begin{array}
[c]{ll}%
c_{\beta,0,0}^{\left(  m\right)  }=c_{\beta,0}^{\left(  m-1\right)  }-\frac
{1}{4}c_{\beta}^{\left(  m-2\right)  }, & 2\leq\left\vert \beta\right\vert
\leq m-2,\\
c_{\beta,1,0}^{\left(  m\right)  }=c_{\beta,1}^{\left(  m-1\right)  }, &
1\leq\left\vert \beta\right\vert \leq m-3,
\end{array}
\right\}  \quad\text{when }\left\{
\begin{array}
[c]{l}%
\beta\in\mathbb{Z}_{\oplus}^{m-2},\\
\beta\leq\mathbf{1},\\
m\geq4.
\end{array}
\right. \label{av622}%
\end{equation}

Combining \ref{av621} and \ref{av622} gives:%
\begin{equation}
\left.
\begin{array}
[c]{ll}%
c_{\beta,0,0}^{\left(  m\right)  }=c_{\beta,0}^{\left(  m-1\right)  }-\frac
{1}{4}c_{\beta}^{\left(  m-2\right)  }, & 2\leq\left\vert \beta\right\vert
\leq m-2,\\
c_{\gamma,1}^{\left(  m\right)  }=c_{\gamma}^{\left(  m-1\right)  }, &
1\leq\left\vert \gamma\right\vert \leq m-3,\\
c_{\beta,1,0}^{\left(  m\right)  }=c_{\beta,1}^{\left(  m-1\right)  }, &
1\leq\left\vert \beta\right\vert \leq m-3,
\end{array}
\right\}  \quad\text{when }\left\{
\begin{array}
[c]{l}%
\gamma\in\mathbb{Z}_{\oplus}^{m-1},\\
\beta\in\mathbb{Z}_{\oplus}^{m-2},\\
\gamma,\beta\leq\mathbf{1},\\
m\geq4,
\end{array}
\right. \label{av627}%
\end{equation}

Indeed, setting $\gamma=\left(  \beta,0\right)  $ and $\gamma=\left(
\beta,0\right)  $ we get%
\begin{equation}
\left.
\begin{array}
[c]{ll}%
c_{\beta,0,0}^{\left(  m\right)  }=c_{\beta,0}^{\left(  m-1\right)  }-\frac
{1}{4}c_{\beta}^{\left(  m-2\right)  }, & 2\leq\left\vert \beta\right\vert
\leq m-2,\\
c_{\beta,0,1}^{\left(  m\right)  }=c_{\beta,0}^{\left(  m-1\right)  }, &
1\leq\left\vert \beta\right\vert \leq m-3,\\
c_{\beta,1,1}^{\left(  m\right)  }=c_{\beta,1}^{\left(  m-1\right)  }, &
0\leq\left\vert \beta\right\vert \leq m-4.\\
c_{\beta,1,0}^{\left(  m\right)  }=c_{\beta,1}^{\left(  m-1\right)  }, &
1\leq\left\vert \beta\right\vert \leq m-3,
\end{array}
\right\}  \quad\text{when }\left\{
\begin{array}
[c]{l}%
\beta\in\mathbb{Z}_{\oplus}^{m-2},\\
\beta\leq\mathbf{1},\\
m\geq4.
\end{array}
\right. \label{av717}%
\end{equation}

We want to extend \ref{av717} so that the equations are true when
$0\leq\left\vert \beta\right\vert \leq m-2$. There are several cases:\medskip

\textbf{1}) \fbox{$c_{\beta,0,0}^{\left(  m\right)  }=c_{\beta,0}^{\left(
m-1\right)  }-\frac{1}{4}c_{\beta}^{\left(  m-2\right)  }$} If $\left\vert
\beta\right\vert =1$ then $LHS=c_{\beta,0,0}^{\left(  m\right)  }%
=c_{\mathbf{i}_{k}}^{\left(  m\right)  }$ for some $1\leq k\leq m-2$ and by
\ref{av598}, $c_{\beta,0,0}^{\left(  m\right)  }=\frac{k}{2^{m-1}}\left(
m+1-k\right)  $. But also by \ref{av598},
\begin{align*}
RHS  & =c_{\mathbf{i}_{k},0}^{\left(  m-1\right)  }-\frac{1}{4}c_{\mathbf{i}%
_{k}}^{\left(  m-2\right)  }\\
& =c_{\mathbf{i}_{k},0}^{\left(  m-1\right)  }-\frac{1}{4}c_{\mathbf{i}_{k}%
}^{\left(  m-2\right)  }\\
& =\frac{k}{2^{m-2}}\left(  m-k\right)  -\frac{1}{4}\frac{k}{2^{m-3}}\left(
m-1-k\right) \\
& =\frac{k}{2^{m-2}}\left(  m-k\right)  -\frac{k}{2^{m-1}}\left(  m-1-k\right)
\\
& =\frac{k}{2^{m-1}}\left(  2\left(  m-k\right)  -\left(  m-1-k\right)
\right) \\
& =\frac{k}{2^{m-1}}\left(  m+1-k\right)  =LHS.
\end{align*}

If $\left\vert \beta\right\vert =0$ then by \ref{av599} $c_{\mathbf{0}%
}^{\left(  m\right)  }=\frac{m+1}{2^{m}}$ so $LHS=c_{\mathbf{0}}^{\left(
m\right)  }=\frac{m+1}{2^{m}}$ and
\begin{align*}
RHS  & =c_{\mathbf{0}}^{\left(  m-1\right)  }-\frac{1}{4}c_{\mathbf{0}%
}^{\left(  m-2\right)  }=\frac{m}{2^{m-1}}-\frac{1}{4}\frac{m-1}{2^{m-2}%
}=\frac{m}{2^{m-1}}-\frac{m-1}{2^{m}}\\
& =\frac{1}{2^{m}}\left(  2m-\left(  m-1\right)  \right)  =\frac{m+1}{2^{m}%
}=LHS.
\end{align*}
\medskip\smallskip

\textbf{2}) \fbox{$c_{\beta,0,1}^{\left(  m\right)  }=c_{\beta,0}^{\left(
m-1\right)  }$} If $\left\vert \beta\right\vert =m-2$ then $\left\vert \left(
\beta,0,1\right)  \right\vert =m-1$ so we can conclude from \ref{av619} that
$LHS=c_{\beta,0,1}^{\left(  m\right)  }=1$. Also, from \ref{av619},
$RHS=c_{\beta,0}^{\left(  m-1\right)  }=1=LHS$.

If $\left\vert \beta\right\vert =\mathbf{0}$ then by \ref{av598},
$LHS=c_{\mathbf{i}_{m}}^{\left(  m\right)  }=\frac{m}{2^{m-1}}$ and by
\ref{av599}, $RHS=c_{\mathbf{0}}^{\left(  m-1\right)  }=\frac{m}{2^{m-1}}%
=LHS$.\smallskip

\textbf{3}) \fbox{$c_{\beta,1,1}^{\left(  m\right)  }=c_{\beta,1}^{\left(
m-1\right)  }$} If $\left\vert \beta\right\vert =m-3$ then $\left\vert \left(
\beta,1,1\right)  \right\vert =m-1$ so we can conclude from \ref{av619} that
$LHS=c_{\beta,1,1}^{\left(  m\right)  }=1$. Also, from \ref{av619},
$RHS=c_{\beta,1}^{\left(  m-1\right)  }=1=LHS$.

If $\left\vert \beta\right\vert =m-2$ then $\left\vert \left(  \beta
,1,1\right)  \right\vert =m$ which means that $\beta=\mathbf{1}$ and so we can
conclude from \ref{av595} that $LHS=c_{\beta,1,1}^{\left(  m\right)  }=1$.
Also, from \ref{av595}, $RHS=c_{\beta,1}^{\left(  m-1\right)  }=1=LHS$%
.\smallskip

\textbf{4}) \fbox{$c_{\beta,1,0}^{\left(  m\right)  }=c_{\beta,1}^{\left(
m-1\right)  }$} If $\left\vert \beta\right\vert =m-2$ then $\left\vert \left(
\beta,1,0\right)  \right\vert =m-1$ so we can conclude from \ref{av619} that
$LHS=c_{\beta,1,0}^{\left(  m\right)  }=1$. Also, from \ref{av619},
$RHS=c_{\beta,1}^{\left(  m-1\right)  }=1=LHS$.

$c_{\mathbf{i}_{k}}^{\left(  m\right)  }=\frac{k}{2^{m-1}}\left(
m+1-k\right)  $

If $\beta=\mathbf{0}$ then by \ref{av598}, $LHS=c_{\beta,1,0}^{\left(
m\right)  }=c_{\mathbf{i}_{m-1}}^{\left(  m\right)  }=\frac{m-1}{2^{m-1}%
}2=\frac{m-1}{2^{m-2}}$ and by \ref{av598}, $RHS=c_{\mathbf{i}_{m-1}}^{\left(
m-1\right)  }=\frac{m-1}{2^{m-2}}\left(  m-\left(  m-1\right)  \right)
=\frac{m-1}{2^{m-2}}=LHS$.\medskip

These cases can be represented by:%
\begin{equation}
\left.
\begin{array}
[c]{l}%
c_{\beta,0,0}^{\left(  m\right)  }=c_{\beta,0}^{\left(  m-1\right)  }-\frac
{1}{4}c_{\beta}^{\left(  m-2\right)  },\\
c_{\beta,0,1}^{\left(  m\right)  }=c_{\beta,0}^{\left(  m-1\right)  },\\
c_{\beta,1,1}^{\left(  m\right)  }=c_{\beta,1}^{\left(  m-1\right)  },\\
c_{\beta,1,0}^{\left(  m\right)  }=c_{\beta,1}^{\left(  m-1\right)  },
\end{array}
\right\}  \quad m\geq4.\label{av301}%
\end{equation}

Next we extend \ref{av301} to $m=3$ using equations \ref{av591} i.e. the
equations%
\[%
\begin{array}
[c]{l}%
c_{0}^{\left(  1\right)  }=c_{1}^{\left(  1\right)  }=1,\\
c_{\mathbf{0}}^{\left(  2\right)  }=\frac{3}{4},\quad c_{1,0}^{\left(
2\right)  }=c_{0,1}^{\left(  2\right)  }=1,\quad c_{\mathbf{1}}^{\left(
2\right)  }=1,\\
c_{\mathbf{0}}^{\left(  3\right)  }=\frac{1}{2},\quad c_{1,0,0}^{\left(
3\right)  }=c_{0,0,1}^{\left(  3\right)  }=\frac{3}{4},\text{ }c_{0,1,0}%
^{\left(  3\right)  }=1,\quad c_{1,1,0}^{\left(  3\right)  }=c_{1,0,1}%
^{\left(  3\right)  }=c_{0,1,1}^{\left(  3\right)  }=1,\quad c_{\mathbf{1}%
}^{\left(  3\right)  }=1.
\end{array}
\]

Indeed, when $m=3$, \ref{av301} becomes%
\[%
\begin{array}
[c]{l}%
c_{0,0,0}^{\left(  3\right)  }=c_{0,0}^{\left(  2\right)  }-\frac{1}{4}%
c_{0}^{\left(  1\right)  },\\
c_{1,0,0}^{\left(  3\right)  }=c_{1,0}^{\left(  2\right)  }-\frac{1}{4}%
c_{1}^{\left(  1\right)  },\\
c_{0,0,1}^{\left(  3\right)  }=c_{0,0}^{\left(  2\right)  },\\
c_{1,0,1}^{\left(  3\right)  }=c_{1,0}^{\left(  2\right)  },\\
c_{0,1,1}^{\left(  3\right)  }=c_{0,1}^{\left(  2\right)  },\\
c_{1,1,1}^{\left(  3\right)  }=c_{1,1}^{\left(  2\right)  },\\
c_{0,1,0}^{\left(  3\right)  }=c_{0,1}^{\left(  2\right)  },\\
c_{1,1,0}^{\left(  3\right)  }=c_{1,1}^{\left(  2\right)  }.
\end{array}
.
\]

and inspection of \ref{av591} shows that these equations all hold and we can
generalize \ref{av301} to%
\begin{equation}
\left.
\begin{array}
[c]{l}%
c_{\beta,0,0}^{\left(  m\right)  }=c_{\beta,0}^{\left(  m-1\right)  }-\frac
{1}{4}c_{\beta}^{\left(  m-2\right)  },\\
c_{\beta,0,1}^{\left(  m\right)  }=c_{\beta,0}^{\left(  m-1\right)  },\\
c_{\beta,1,1}^{\left(  m\right)  }=c_{\beta,1}^{\left(  m-1\right)  },\\
c_{\beta,1,0}^{\left(  m\right)  }=c_{\beta,1}^{\left(  m-1\right)  },
\end{array}
\right\}  \quad m\geq3.\label{av370}%
\end{equation}

\[
c_{\gamma,0,0}^{\left(  m+1\right)  }=c_{\gamma,0}^{\left(  m\right)  }%
-\frac{1}{4}c_{\gamma}^{\left(  m-1\right)  },
\]

and there are two cases: $\gamma=\left(  \delta,0\right)  $ and $\gamma
=\left(  \delta,1\right)  $. When $\gamma=\left(  \delta,1\right)  $,%
\[
c_{\delta,1,0,0}^{\left(  m+1\right)  }=c_{\delta,1,0}^{\left(  m\right)
}-\frac{1}{4}c_{\delta,1}^{\left(  m-1\right)  }=c_{\delta,1}^{\left(
m-1\right)  }-\frac{1}{4}c_{\delta,1}^{\left(  m-1\right)  }=\frac{3}%
{4}c_{\delta,1}^{\left(  m-1\right)  }>0.
\]

When $\gamma=\left(  \delta,0\right)  $,%
\[
c_{\delta,0,0,0}^{\left(  m+1\right)  }=c_{\delta,0,0}^{\left(  m\right)
}-\frac{1}{4}c_{\delta,0}^{\left(  m-1\right)  }=c_{\delta,0}^{\left(
m-1\right)  }-\frac{1}{4}c_{\delta}^{\left(  m-2\right)  }-\frac{1}%
{4}c_{\delta,0}^{\left(  m-1\right)  }=\frac{3}{4}c_{\delta,0}^{\left(
m-1\right)  }-\frac{1}{4}c_{\delta}^{\left(  m-2\right)  }.
\]

Since%
\begin{align}
c_{\gamma,0,0}^{\left(  m\right)  }  & =c_{\gamma,0}^{\left(  m-1\right)
}-\frac{1}{4}c_{\gamma}^{\left(  m-2\right)  },\label{av630}\\
c_{\gamma,0,0,0}^{\left(  m\right)  }  & =c_{\gamma,0,0}^{\left(  m-1\right)
}-\frac{1}{4}c_{\gamma,0}^{\left(  m-2\right)  }=c_{\gamma,0}^{\left(
m-2\right)  }-\frac{1}{4}c_{\gamma}^{\left(  m-3\right)  }-\frac{1}%
{4}c_{\gamma,0}^{\left(  m-2\right)  }=\frac{3}{4}c_{\gamma,0}^{\left(
m-2\right)  }-\frac{1}{4}c_{\gamma}^{\left(  m-3\right)  },\nonumber
\end{align}

we hypothesize that%
\begin{equation}
c_{\gamma,0_{n}}^{\left(  m\right)  }=a_{n}c_{\gamma,0}^{\left(  m-n+1\right)
}-b_{n}c_{\gamma}^{\left(  m-n\right)  },\quad1\leq n\leq m-1.\label{av648}%
\end{equation}

When $n=1$ we see that $a_{1}=1$ and $b_{1}=0$. When $n=2$ we see that
$a_{2}=1$ and $b_{2}=1/4$. Thus%
\begin{equation}
Initial\text{ }values:\left\{
\begin{array}
[c]{l}%
a_{1}=1,\text{ }a_{2}=1,\\
b_{1}=0,\text{ }b_{2}=\frac{1}{4}.
\end{array}
\right. \label{av659}%
\end{equation}

If $n\geq3$ and $m\geq4$, from \ref{av630} and then \ref{av648},%
\begin{align}
c_{\gamma,0_{n}}^{\left(  m\right)  }  & =c_{\gamma,0_{n-2},0,0}^{\left(
m\right)  }\nonumber\\
& =c_{\gamma,0_{n-2},0}^{\left(  m-1\right)  }-\frac{1}{4}c_{\gamma,0_{n-2}%
}^{\left(  m-2\right)  }\nonumber\\
& =c_{\gamma,0_{n-1}}^{\left(  m-1\right)  }-\frac{1}{4}c_{\gamma,0_{n-2}%
}^{\left(  m-2\right)  }\nonumber\\
& =\left(  a_{n-1}c_{\gamma,0}^{\left(  m-n+1\right)  }-b_{n-1}c_{\gamma
}^{\left(  m-n\right)  }\right)  -\frac{1}{4}\left(  a_{n-2}c_{\gamma
,0}^{\left(  m-n+1\right)  }-b_{n-2}c_{\gamma}^{\left(  m-n\right)  }\right)
\nonumber\\
& =a_{n-1}c_{\gamma,0}^{\left(  m-n+1\right)  }-b_{n-1}c_{\gamma}^{\left(
m-n\right)  }-\frac{1}{4}a_{n-2}c_{\gamma,0}^{\left(  m-n+1\right)  }+\frac
{1}{4}b_{n-2}c_{\gamma}^{\left(  m-n\right)  }\nonumber\\
& =\left(  a_{n-1}-\frac{1}{4}a_{n-2}\right)  c_{\gamma,0}^{\left(
m-n+1\right)  }-\left(  b_{n-1}-\frac{1}{4}b_{n-2}\right)  c_{\gamma}^{\left(
m-n\right)  }\nonumber\\
& =a_{n}c_{\gamma,0}^{\left(  m-n+1\right)  }-b_{n}c_{\gamma}^{\left(
m-n\right)  }\label{av610}\\
& =a_{n}c_{\gamma,0}^{\left(  \left\vert \gamma\right\vert +1\right)  }%
-b_{n}c_{\gamma}^{\left(  \left\vert \gamma\right\vert \right)  },\nonumber
\end{align}

so that we obtain the two second-order difference equations:%
\begin{equation}%
\begin{array}
[c]{c}%
a_{n}=a_{n-1}-\frac{1}{4}a_{n-2},\quad a_{1}=1,\text{ }a_{2}=1,\\
b_{n}=b_{n-1}-\frac{1}{4}b_{n-2},\quad b_{1}=0,\text{ }b_{2}=\frac{1}{4}.
\end{array}
\label{av608}%
\end{equation}

Recall \ref{av612}: the auxiliary equation is $x^{2}-x+1/4$ which has the
single zero $1/2$ so the solutions have the form%
\[
a_{n}=\frac{1}{2^{n}}\left(  A+Bn\right)  ,\text{ }b_{n}=\frac{1}{2^{n}%
}\left(  A^{\prime}+B^{\prime}n\right)  .
\]

The initial conditions imply $1=\frac{1}{2}\left(  A+B\right)  =\frac{1}%
{4}\left(  A+2B\right)  $ so that $A=0$ and $B=2$. Also, $0=\frac{1}{2}\left(
A^{\prime}+B^{\prime}\right)  $ and $\frac{1}{4}=\frac{1}{4}\left(  A^{\prime
}+2B^{\prime}\right)  $ so that $A^{\prime}=-1$ and $B^{\prime}=1 $. Hence%
\begin{equation}
a_{n}=\frac{n}{2^{n-1}},\text{ }b_{n}=\frac{n-1}{2^{n}},\label{av609}%
\end{equation}

and \ref{av610} becomes%
\begin{equation}
c_{\gamma,0_{n}}^{\left(  m\right)  }=\frac{n}{2^{n-1}}c_{\gamma,0}^{\left(
m-n+1\right)  }-\frac{n-1}{2^{n}}c_{\gamma}^{\left(  m-n\right)  },\quad
n\geq1.\label{av611}%
\end{equation}

Thus by \ref{av591}, if $m\geq2$,%
\begin{align}
c_{\mathbf{0}}^{\left(  m\right)  }=c_{0,0_{m=1}}^{\left(  m\right)  }%
=\frac{m-1}{2^{m-2}}c_{0,0}^{\left(  2\right)  }-\frac{m-2}{2^{m-1}}%
c_{0}^{\left(  1\right)  }  & =\frac{m-1}{2^{m-2}}\frac{3}{4}-\frac
{m-2}{2^{m-1}}\nonumber\\
& =\frac{1}{2^{m}}\left\{  3\left(  m-1\right)  -2\left(  m-2\right)  \right\}
\nonumber\\
& =\frac{m+1}{2^{m}}.\label{av598}%
\end{align}

Consider: if $m\geq3$ then $\gamma$ must have one of the unique partitions:
\begin{equation}
\gamma=\left\{
\begin{array}
[c]{l}%
z^{\left(  1\right)  },\\
\left(  z^{\left(  1\right)  }u^{\left(  1\right)  }\right)  \left(
z^{\left(  2\right)  }u^{\left(  2\right)  }\right)  \ldots\left(  z^{\left(
k\right)  }u^{\left(  k\right)  }\right)  ,\\
\left(  z^{\left(  1\right)  }u^{\left(  1\right)  }\right)  \left(
z^{\left(  2\right)  }u^{\left(  2\right)  }\right)  \ldots\left(  z^{\left(
k\right)  }u^{\left(  k\right)  }\right)  z^{\left(  k+1\right)  },\\
u^{\left(  1\right)  },\\
\left(  u^{\left(  1\right)  }z^{\left(  1\right)  }\right)  \ldots\left(
u^{\left(  k\right)  }z^{\left(  k\right)  }\right)  ,\\
\left(  u^{\left(  1\right)  }z^{\left(  1\right)  }\right)  \ldots\left(
u^{\left(  k\right)  }z^{\left(  k\right)  }\right)  u^{\left(  k+1\right)
},\\
0\left(  u^{\left(  1\right)  }z^{\left(  1\right)  }\right)  \ldots\left(
u^{\left(  k\right)  }z^{\left(  k\right)  }\right)  ,\\
0\left(  u^{\left(  1\right)  }z^{\left(  1\right)  }\right)  \ldots\left(
u^{\left(  k\right)  }z^{\left(  k\right)  }\right)  u^{\left(  k+1\right)  },
\end{array}
\right.  ,\quad k\geq1.\label{av646}%
\end{equation}

where the $u^{\left(  k\right)  }$s have only isolated zeros and the
$z^{\left(  k\right)  }$s are maximal sequences of zeros (of at least length 2).

Observe that if $z=0_{k}$ is a sequence of (at least two) zeros, $u=1$ or the
sequence $u=1\ldots1$ has only isolated zeros then applying \ref{av370} to
$c_{z,u}^{\left(  m\right)  }$ to eliminate $u$ and then applying \ref{av598}
we get:%
\begin{equation}
c_{z,u}^{\left(  m\right)  }=c_{0_{k}1\ldots1}^{\left(  m\right)  }=c_{0_{k}%
}^{\left(  m-\left\Vert u\right\Vert \right)  }=\frac{m-\left\Vert
u\right\Vert +1}{2^{m-\left\Vert u\right\Vert }}=\frac{1+\left\Vert
z\right\Vert }{2^{\left\vert z\right\vert }},\label{av303}%
\end{equation}

where $\left\Vert u\right\Vert $ is the length of the sequence $u$. If
$\left\vert \gamma\right\vert \geq0$,%
\begin{align}
c_{\gamma,1,z,u}^{\left(  m\right)  }=c_{\gamma,1,z}^{\left(  m-\left\Vert
u\right\Vert \right)  }  & =\frac{\left\Vert z\right\Vert }{2^{\left\Vert
z\right\Vert -1}}c_{\gamma,1,0}^{\left(  m-\left\Vert u,z\right\Vert
+1\right)  }-\frac{\left\Vert z\right\Vert -1}{2^{\left\Vert z\right\Vert }%
}c_{\gamma,1}^{\left(  m-\left\Vert u,z\right\Vert \right)  }\nonumber\\
& =\frac{\left\Vert z\right\Vert }{2^{\left\Vert z\right\Vert -1}}c_{\gamma
,1}^{\left(  m-\left\Vert u,z\right\Vert \right)  }-\frac{\left\Vert
z\right\Vert -1}{2^{\left\Vert z\right\Vert }}c_{\gamma,1}^{\left(
m-\left\Vert u,z\right\Vert \right)  }\nonumber\\
& =\left(  \frac{\left\Vert z\right\Vert }{2^{\left\Vert z\right\Vert -1}%
}-\frac{\left\Vert z\right\Vert -1}{2^{\left\Vert z\right\Vert }}\right)
c_{\gamma,1}^{\left(  m-\left\Vert u,z\right\Vert \right)  }\nonumber\\
& =\frac{1+\left\Vert z\right\Vert }{2^{\left\Vert z\right\Vert }}c_{\gamma
,1}^{\left(  m-\left\Vert u,z\right\Vert \right)  }.\label{av614}%
\end{align}

Also, from \ref{av611}, $c_{\gamma,0_{n}}^{\left(  m\right)  }=\frac
{n}{2^{n-1}}c_{\gamma,0}^{\left(  m-n+1\right)  }-\frac{n-1}{2^{n}}c_{\gamma
}^{\left(  m-n\right)  }$ so%
\begin{align}
c_{u,z}^{\left(  m\right)  }  & =\frac{\left\Vert z\right\Vert }{2^{\left\Vert
z\right\Vert -1}}c_{u,0}^{\left(  \left\Vert u\right\Vert +1\right)  }%
-\frac{\left\Vert z\right\Vert -1}{2^{\left\Vert z\right\Vert }}c_{u}^{\left(
\left\Vert u\right\Vert \right)  }\nonumber\\
& =\frac{\left\Vert z\right\Vert }{2^{\left\Vert z\right\Vert -1}}%
c_{u}^{\left(  \left\Vert u\right\Vert \right)  }-\frac{\left\Vert
z\right\Vert -1}{2^{\left\Vert z\right\Vert }}c_{u}^{\left(  \left\Vert
u\right\Vert \right)  }\nonumber\\
& =\frac{1+\left\Vert z\right\Vert }{2^{\left\Vert z\right\Vert }}%
c_{u}^{\left(  \left\Vert u\right\Vert \right)  }=\left\{
\begin{array}
[c]{c}%
\frac{1+\left\Vert z\right\Vert }{2^{\left\Vert z\right\Vert }}c_{1}^{\left(
1\right)  }\\
or\\
\frac{1+\left\Vert z\right\Vert }{2^{\left\Vert z\right\Vert }}c_{1?}^{\left(
2\right)  }%
\end{array}
\right\}  =\frac{1+\left\Vert z\right\Vert }{2^{\left\Vert z\right\Vert }%
},\label{av13}%
\end{align}

and when $\left\vert \delta\right\vert \geq1$%
\begin{align}
c_{\delta,u,z}^{\left(  m\right)  }  & =\frac{\left\Vert z\right\Vert
}{2^{\left\Vert z\right\Vert -1}}c_{\delta,u,0}^{\left(  m-\left\Vert
z\right\Vert +1\right)  }-\frac{\left\Vert z\right\Vert -1}{2^{\left\Vert
z\right\Vert }}c_{\delta,u}^{\left(  m-\left\Vert z\right\Vert \right)
}\nonumber\\
& =\frac{\left\Vert z\right\Vert }{2^{\left\Vert z\right\Vert -1}}c_{\delta
,u}^{\left(  m-\left\Vert z\right\Vert \right)  }-\frac{\left\Vert
z\right\Vert -1}{2^{\left\Vert z\right\Vert }}c_{\delta,u}^{\left(
m-\left\Vert z\right\Vert \right)  }\nonumber\\
& =\left(  \frac{\left\Vert z\right\Vert }{2^{\left\Vert z\right\Vert -1}%
}-\frac{\left\Vert z\right\Vert -1}{2^{\left\Vert z\right\Vert }}\right)
c_{\delta,u}^{\left(  m-\left\Vert z\right\Vert \right)  }\nonumber\\
& =\frac{1+\left\Vert z\right\Vert }{2^{\left\Vert z\right\Vert }}c_{\delta
,u}^{\left(  m-\left\Vert z\right\Vert \right)  }\nonumber\\
& =\frac{1+\left\Vert z\right\Vert }{2^{\left\Vert z\right\Vert }}c_{\delta
}^{\left(  m-\left\Vert u\right\Vert -\left\Vert z\right\Vert \right)
},\label{av645}%
\end{align}

so that%
\begin{equation}
c_{\delta,u,z}^{\left(  m\right)  }=\left\{
\begin{array}
[c]{ll}%
\frac{1+\left\Vert z\right\Vert }{2^{\left\Vert z\right\Vert }}c_{\delta
}^{\left(  m-\left\Vert u\right\Vert -\left\Vert z\right\Vert \right)  }, &
\left\vert \delta\right\vert \geq1,\\
\frac{1+\left\Vert z\right\Vert }{2^{\left\Vert z\right\Vert }}, & \delta=0.
\end{array}
\right. \label{av17}%
\end{equation}

The cases in \ref{av646} are:\bigskip

\fbox{Case $\gamma=z^{\left(  1\right)  }$} From \ref{av598},%
\begin{equation}
c_{z^{\left(  1\right)  }}^{\left(  m\right)  }=\frac{1+\left\Vert z^{\left(
1\right)  }\right\Vert }{2^{\left\Vert z^{\left(  1\right)  }\right\Vert }%
}.\label{av477}%
\end{equation}

\fbox{Case $\gamma=\left(  z^{\left(  1\right)  }u^{\left(  1\right)
}\right)  \left(  z^{\left(  2\right)  }u^{\left(  2\right)  }\right)
\ldots\left(  z^{\left(  k\right)  }u^{\left(  k\right)  }\right)  $} From
\ref{av614} and then \ref{av303},%
\begin{align}
c_{\left(  z^{\left(  1\right)  }u^{\left(  1\right)  }\right)  \left(
z^{\left(  2\right)  }u^{\left(  2\right)  }\right)  \ldots\left(  z^{\left(
k\right)  }u^{\left(  k\right)  }\right)  }^{\left(  m\right)  }  &
=\frac{1+\left\Vert z^{\left(  k\right)  }\right\Vert }{2^{\left\Vert
z^{\left(  k\right)  }\right\Vert }}c_{\left(  z^{\left(  1\right)
}u^{\left(  1\right)  }\right)  \left(  z^{\left(  2\right)  }u^{\left(
2\right)  }\right)  \ldots\left(  z^{\left(  k-1\right)  }u^{\left(
k-1\right)  }\right)  }^{\left(  m-\left\Vert z^{\left(  k\right)  }u^{\left(
k\right)  }\right\Vert \right)  }\nonumber\\
& =\left(  \frac{1+\left\Vert z^{\left(  2\right)  }\right\Vert }%
{2^{\left\Vert z^{\left(  2\right)  }\right\Vert }}\times\ldots\times
\frac{1+\left\Vert z^{\left(  k\right)  }\right\Vert }{2^{\left\Vert
z^{\left(  k\right)  }\right\Vert }}\right)  c_{z^{\left(  1\right)
}u^{\left(  1\right)  }}^{\left(  \left\Vert z^{\left(  1\right)  }u^{\left(
1\right)  }\right\Vert \right)  }\nonumber\\
& =\frac{1+\left\Vert z^{\left(  1\right)  }\right\Vert }{2^{\left\Vert
z^{\left(  1\right)  }\right\Vert }}\times\ldots\times\frac{1+\left\Vert
z^{\left(  k\right)  }\right\Vert }{2^{\left\Vert z^{\left(  k\right)
}\right\Vert }}.\label{av305}%
\end{align}

\fbox{Case $\gamma=\left(  z^{\left(  1\right)  }u^{\left(  1\right)
}\right)  \left(  z^{\left(  2\right)  }u^{\left(  2\right)  }\right)
\ldots\left(  z^{\left(  k\right)  }u^{\left(  k\right)  }\right)  z^{\left(
k+1\right)  }$} From \ref{av17} and \ref{av598}:%
\begin{align}
c_{\left(  z^{\left(  1\right)  }u^{\left(  1\right)  }\right)  \left(
z^{\left(  2\right)  }u^{\left(  2\right)  }\right)  \ldots\left(  z^{\left(
k\right)  }u^{\left(  k\right)  }\right)  z^{\left(  k+1\right)  }}^{\left(
m\right)  } &  =c_{z^{\left(  1\right)  }\left(  u^{\left(  1\right)
}z^{\left(  2\right)  }\right)  \ldots\left(  u^{\left(  k\right)  }z^{\left(
k+1\right)  }\right)  }^{\left(  m\right)  }\nonumber\\
&  =\frac{1+\left\Vert z^{\left(  k+1\right)  }\right\Vert }{2^{\left\Vert
z^{\left(  k+1\right)  }\right\Vert }}c_{z^{\left(  1\right)  }\left(
u^{\left(  1\right)  }z^{\left(  2\right)  }\right)  \ldots\left(  u^{\left(
k-1\right)  }z^{\left(  k\right)  }\right)  }^{\left(  m-\left\vert u^{\left(
k\right)  }z^{\left(  k+1\right)  }\right\vert \right)  }\\
&  =\frac{1+\left\Vert z^{\left(  k\right)  }\right\Vert }{2^{\left\Vert
z^{\left(  k\right)  }\right\Vert }}\frac{1+\left\Vert z^{\left(  k+1\right)
}\right\Vert }{2^{\left\Vert z^{\left(  k+1\right)  }\right\Vert }%
}c_{z^{\left(  1\right)  }\left(  u^{\left(  1\right)  }z^{\left(  2\right)
}\right)  \ldots\left(  u^{\left(  k-1\right)  }z^{\left(  k\right)  }\right)
}^{\left(  m\right)  }\nonumber\\
&  =\frac{1+\left\Vert z^{\left(  2\right)  }\right\Vert }{2^{\left\Vert
z^{\left(  2\right)  }\right\Vert }}\times\ldots\times\frac{1+\left\Vert
z^{\left(  k+1\right)  }\right\Vert }{2^{\left\Vert z^{\left(  k+1\right)
}\right\Vert }}c_{z^{\left(  1\right)  }}^{\left(  \left\Vert z^{\left(
1\right)  }\right\Vert \right)  }\nonumber\\
&  =\frac{1+\left\Vert z^{\left(  1\right)  }\right\Vert }{2^{\left\Vert
z^{\left(  1\right)  }\right\Vert }}\times\ldots\times\frac{1+\left\Vert
z^{\left(  k+1\right)  }\right\Vert }{2^{\left\Vert z^{\left(  k+1\right)
}\right\Vert }}.\label{av306}%
\end{align}

\fbox{Case $\gamma=u^{\left(  1\right)  }$} Here $c_{u^{\left(  1\right)  }%
}^{\left(  m\right)  }=c_{1}^{\left(  1\right)  }=1$ or from \ref{av370} and
\ref{av591}, $c_{u^{\left(  1\right)  }}^{\left(  m\right)  }=c_{1,1}^{\left(
2\right)  }=1$.
\begin{equation}
c_{u^{\left(  1\right)  }}^{\left(  m\right)  }=1.\label{av723}%
\end{equation}

This confirms \ref{av595} i.e. $c_{\mathbf{1}}^{\left(  m\right)  }=1,\quad
m\geq1$.\medskip

\fbox{Case $\gamma=\left(  u^{\left(  1\right)  }z^{\left(  1\right)
}\right)  \ldots\left(  u^{\left(  k\right)  }z^{\left(  k\right)  }\right)
$} Next, applying \ref{av17} we get%
\begin{align}
c_{\left(  u^{\left(  1\right)  }z^{\left(  1\right)  }\right)  \ldots\left(
u^{\left(  k\right)  }z^{\left(  k\right)  }\right)  }^{\left(  m\right)  } &
=\frac{1+\left\Vert z^{\left(  k\right)  }\right\Vert }{2^{\left\Vert
z^{\left(  k\right)  }\right\Vert }}c_{\left(  u^{\left(  1\right)
}z^{\left(  1\right)  }\right)  \ldots\left(  u^{\left(  k-1\right)
}z^{\left(  k-1\right)  }\right)  }^{\left(  m-\left\Vert u^{\left(  k\right)
}z^{\left(  k\right)  }\right\Vert \right)  }\nonumber\\
&  =\frac{1+\left\Vert z^{\left(  2\right)  }\right\Vert }{2^{\left\Vert
z^{\left(  2\right)  }\right\Vert }}\times\ldots\times\frac{1+\left\Vert
z^{\left(  k\right)  }\right\Vert }{2^{\left\Vert z^{\left(  k\right)
}\right\Vert }}c_{\left(  u^{\left(  1\right)  }z^{\left(  1\right)  }\right)
}^{\left(  \left\Vert u^{\left(  1\right)  }z^{\left(  1\right)  }\right\Vert
\right)  }\nonumber\\
&  =\frac{1+\left\Vert z^{\left(  1\right)  }\right\Vert }{2^{\left\Vert
z^{\left(  1\right)  }\right\Vert }}\times\ldots\times\frac{1+\left\Vert
z^{\left(  k\right)  }\right\Vert }{2^{\left\Vert z^{\left(  k\right)
}\right\Vert }}.\label{av15}%
\end{align}

\fbox{$\left(  u^{\left(  1\right)  }z^{\left(  1\right)  }\right)
\ldots\left(  u^{\left(  k\right)  }z^{\left(  k\right)  }\right)  u^{\left(
k+1\right)  }$} Next, from \ref{av370} and then \ref{av17},%
\begin{align}
c_{\left(  u^{\left(  1\right)  }z^{\left(  1\right)  }\right)  \ldots\left(
u^{\left(  k\right)  }z^{\left(  k\right)  }\right)  u^{\left(  k+1\right)  }%
}^{\left(  m\right)  } &  =c_{\left(  u^{\left(  1\right)  }z^{\left(
1\right)  }\right)  \ldots\left(  u^{\left(  k\right)  }z^{\left(  k\right)
}\right)  }^{\left(  m-\left\Vert u^{\left(  k+1\right)  }\right\Vert \right)
}\nonumber\\
&  =\frac{1+\left\Vert z^{\left(  k\right)  }\right\Vert }{2^{\left\Vert
z^{\left(  k\right)  }\right\Vert }}c_{\left(  u^{\left(  1\right)
}z^{\left(  1\right)  }\right)  \ldots\left(  u^{\left(  k-1\right)
}z^{\left(  k-1\right)  }\right)  }^{\left(  m-\left\Vert u^{\left(
k+1\right)  }\right\Vert \right)  }\nonumber\\
&  =\frac{1+\left\Vert z^{\left(  2\right)  }\right\Vert }{2^{\left\Vert
z^{\left(  2\right)  }\right\Vert }}\times\ldots\times\frac{1+\left\Vert
z^{\left(  k\right)  }\right\Vert }{2^{\left\Vert z^{\left(  k\right)
}\right\Vert }}c_{\left(  u^{\left(  1\right)  }z^{\left(  1\right)  }\right)
}^{\left(  \left\Vert u^{\left(  1\right)  }z^{\left(  1\right)  }\right\Vert
\right)  }\nonumber\\
&  =\frac{1+\left\Vert z^{\left(  1\right)  }\right\Vert }{2^{\left\Vert
z^{\left(  1\right)  }\right\Vert }}\times\ldots\times\frac{1+\left\Vert
z^{\left(  k\right)  }\right\Vert }{2^{\left\Vert z^{\left(  k\right)
}\right\Vert }}.\label{av669}%
\end{align}

\fbox{$0\left(  u^{\left(  1\right)  }z^{\left(  1\right)  }\right)
\ldots\left(  u^{\left(  k\right)  }z^{\left(  k\right)  }\right)  $} Next,
from \ref{av17},%
\begin{align}
c_{0\left(  u^{\left(  1\right)  }z^{\left(  1\right)  }\right)  \ldots\left(
u^{\left(  k\right)  }z^{\left(  k\right)  }\right)  }^{\left(  m\right)  } &
=\frac{1+\left\Vert z^{\left(  k\right)  }\right\Vert }{2^{\left\Vert
z^{\left(  k\right)  }\right\Vert }}c_{0\left(  u^{\left(  1\right)
}z^{\left(  1\right)  }\right)  \ldots\left(  u^{\left(  k-1\right)
}z^{\left(  k-1\right)  }\right)  }^{\left(  m-\left\Vert u^{\left(  k\right)
}z^{\left(  k\right)  }\right\Vert \right)  }\nonumber\\
&  =\frac{1+\left\Vert z^{\left(  1\right)  }\right\Vert }{2^{\left\Vert
z^{\left(  1\right)  }\right\Vert }}\times\ldots\times\frac{1+\left\Vert
z^{\left(  k\right)  }\right\Vert }{2^{\left\Vert z^{\left(  k\right)
}\right\Vert }}c_{0}^{\left(  1\right)  }\nonumber\\
&  =\frac{1+\left\Vert z^{\left(  1\right)  }\right\Vert }{2^{\left\Vert
z^{\left(  1\right)  }\right\Vert }}\times\ldots\times\frac{1+\left\Vert
z^{\left(  k\right)  }\right\Vert }{2^{\left\Vert z^{\left(  k\right)
}\right\Vert }}.\label{av674}%
\end{align}

\fbox{$0\left(  u^{\left(  1\right)  }z^{\left(  1\right)  }\right)
\ldots\left(  u^{\left(  k\right)  }z^{\left(  k\right)  }\right)  u^{\left(
k+1\right)  }$} From \ref{av370},%
\[
c_{0\left(  u^{\left(  1\right)  }z^{\left(  1\right)  }\right)  \ldots\left(
u^{\left(  k\right)  }z^{\left(  k\right)  }\right)  u^{\left(  k+1\right)  }%
}^{\left(  m\right)  }=c_{0\left(  u^{\left(  1\right)  }z^{\left(  1\right)
}\right)  \ldots\left(  u^{\left(  k\right)  }z^{\left(  k\right)  }\right)
}^{\left(  m-\left\Vert u^{\left(  k+1\right)  }\right\Vert \right)  },
\]

and from \ref{av674},%
\begin{equation}
c_{0\left(  u^{\left(  1\right)  }z^{\left(  1\right)  }\right)  \ldots\left(
u^{\left(  k\right)  }z^{\left(  k\right)  }\right)  u^{\left(  k+1\right)  }%
}^{\left(  m\right)  }=\frac{1+\left\Vert z^{\left(  1\right)  }\right\Vert
}{2^{\left\Vert z^{\left(  1\right)  }\right\Vert }}\times\ldots\times
\frac{1+\left\Vert z^{\left(  k\right)  }\right\Vert }{2^{\left\Vert
z^{\left(  k\right)  }\right\Vert }}.\label{av729}%
\end{equation}

Equations \ref{av477}, \ref{av305}, \ref{av306}, \ref{av723}, \ref{av15},
\ref{av669}, \ref{av674} and \ref{av729} now yield the following
multiplicative representation for the $c_{\gamma}^{\left(  m\right)  }$ in the
expansion \ref{av599}:

\begin{theorem}
\label{vThm_formula_for_C^(m)_idx}(Multiplicative representation) For $m\geq1
$,%
\[
\left\vert T\left(  a;1\right)  \right\vert =\sum_{k=0}^{m}\sum
\limits_{\substack{\beta\leq\mathbf{1} \\\left\vert \beta\right\vert
=k}}c_{\beta}^{\left(  m\right)  }a^{\beta},\quad a\in\mathbb{R}^{m},
\]

where all $0<c_{\beta}^{\left(  m\right)  }\leq1$.

Now suppose the multi-index $0\leq\beta\leq1_{m}$ contains non-trivial
sequences of zeros i.e. sequences with length $\geq2$. If the sequences are
$z^{\left(  1\right)  },\ldots,z^{\left(  k\right)  }$ then%
\begin{equation}
c_{\beta}^{\left(  m\right)  }=\frac{1+\left\Vert z^{\left(  1\right)
}\right\Vert }{2^{\left\Vert z^{\left(  1\right)  }\right\Vert }}\times
\ldots\times\frac{1+\left\Vert z^{\left(  k\right)  }\right\Vert
}{2^{\left\Vert z^{\left(  k\right)  }\right\Vert }},\label{av052}%
\end{equation}

and there are no such zero sequences iff $c_{\beta}^{\left(  m\right)  }=1$.
Here $\left\Vert z\right\Vert $ denotes the length of the sequence $z$.
\end{theorem}

\begin{remark}
\label{vRem_Thm_formula_for_C^(m)_idx} Suppose the multi-index $0_{m}\leq
\beta\leq1_{m}$ contains the (possible trivial) sequences of zeros
$\zeta^{\left(  1\right)  },\ldots,\zeta^{\left(  n\right)  }$ i.e. sequences
with length $\geq1$. Then%
\begin{align*}
c_{\beta}^{\left(  m\right)  }  & =\frac{1+\left\Vert \zeta^{\left(  1\right)
}\right\Vert }{2^{\left\Vert \zeta^{\left(  1\right)  }\right\Vert }}%
\times\ldots\times\frac{1+\left\Vert \zeta^{\left(  n\right)  }\right\Vert
}{2^{\left\Vert \zeta^{\left(  n\right)  }\right\Vert }}\\
& =\frac{1}{2^{\left\Vert \zeta^{\left(  1\right)  }\right\Vert +\ldots
+\left\Vert \zeta^{\left(  n\right)  }\right\Vert }}\left(  1+\left\Vert
\zeta^{\left(  1\right)  }\right\Vert \right)  \times\ldots\times\left(
1+\left\Vert \zeta^{\left(  n\right)  }\right\Vert \right) \\
& =\frac{1}{2^{m-\left\vert \beta\right\vert }}\left(  1+\left\Vert
\zeta^{\left(  1\right)  }\right\Vert \right)  \times\ldots\times\left(
1+\left\Vert \zeta^{\left(  n\right)  }\right\Vert \right)  .
\end{align*}

Can this formula be related to the entropy $H\left(  X\right)  =-\sum
\limits_{i}p_{i}\log_{2}p_{i}$ of some random variable $X$ on a discrete
multivariate distribution?
\end{remark}

\subsection{Some consequences of Theorem \ref{vThm_formula_for_C^(m)_idx}}

\begin{corollary}
\label{vCor_1_formula_for_C^(m)_gam}If $\mathcal{R}$ is the \textbf{reverse
order permutation} then:

\begin{enumerate}
\item $c_{\mathcal{R}\gamma}^{\left(  m\right)  }=c_{\gamma}^{\left(
m\right)  }$.

\item $\left\vert T\left(  \mathcal{R}x\right)  \right\vert =\left\vert
T\left(  x\right)  \right\vert $.
\end{enumerate}
\end{corollary}

\begin{proof}
The reverse order permutation is defined by: $\mathcal{R}\left(  x_{1}%
,x_{2},\ldots,x_{m}\right)  =\left(  x_{m},x_{m-1},\ldots,x_{1}\right)  $ i.e.
$\left(  \mathcal{R}x\right)  _{k}=x_{m+1-k}$.\medskip

\textbf{Part 1} Regarding the formula \ref{av052} for $c_{\gamma}^{\left(
m\right)  }$, sequences of zeros are mapped to sequences of zeros of the same
length.\medskip

\textbf{Part 2} From Part 2,%
\begin{align*}
\left\vert T\left(  \mathcal{R}x;1\right)  \right\vert  & =\sum_{k=0}^{m}%
\sum\limits_{\beta\leq\mathbf{1},\left\vert \beta\right\vert =k}c_{\beta
}^{\left(  m\right)  }\left(  \mathcal{R}x\right)  ^{\beta}=\sum_{k=0}^{m}%
\sum\limits_{\beta\leq\mathbf{1},\left\vert \beta\right\vert =k}c_{\beta
}^{\left(  m\right)  }x^{\mathcal{R}\beta}=\\
& =\sum_{k=0}^{m}\sum\limits_{\beta\leq\mathbf{1},\left\vert \beta\right\vert
=k}c_{\mathcal{R}\beta}^{\left(  m\right)  }x^{\mathcal{R}\beta}=\sum
_{k=0}^{m}\sum\limits_{\beta\leq\mathbf{1},\left\vert \beta\right\vert
=k}c_{\beta}^{\left(  m\right)  }x^{\beta}=\\
& =\left\vert T\left(  x;1\right)  \right\vert .
\end{align*}

Alternatively, regarding formula \ref{av054} for $\left\vert T\left(
x\right)  \right\vert $, first swapping rows as $r_{k}\leftrightarrow
r_{m+1-k}$ and then columns as $c_{k}\leftrightarrow c_{m+1-k}$ yields this part.
\end{proof}

The next corollary requires the following elementary properties of the
function $\frac{1+x}{2^{x}}$.

\begin{remark}
\label{vRem_1_plus_X_div_2^X}\textbf{Elementary properties of} $\frac
{1+x}{2^{x}}$:

\begin{enumerate}
\item \textbf{Monotonic decreasing} When $x\geq1$ the function $\frac
{1+x}{2^{x}}$ decreases monotonically from a maximum of one.

\item \textbf{Subadditivity} $\frac{1+x+y}{2^{x+y}}<\frac{1+x}{2^{x}}%
\frac{1+y}{2^{y}}$ when $x,y>0$.
\end{enumerate}
\end{remark}

\begin{corollary}
\label{vCor_2_formula_for_C^(m)_gam}\ 

\begin{enumerate}
\item $c_{\mathbf{0}}^{\left(  m\right)  }=\frac{m+1}{2^{m}}\leq
\frac{1+m-\left\vert \gamma\right\vert }{2^{m-\left\vert \gamma\right\vert }%
}\leq c_{\gamma}^{\left(  m\right)  }\leq1$ for all $\gamma$.

\item $c_{\mathbf{0}}^{\left(  m\right)  }\left(  1+a\right)  ^{\mathbf{1}%
}\leq\left\vert T\left(  a\right)  \right\vert \leq\left(  1+a\right)
^{\mathbf{1}}$ when $a\in\mathbb{R}_{\oplus}^{m}$. This inequality is sharp.

\item When $a\in\mathbb{R}_{+}^{m}$ we have the polynomial upper bound
\[
0<\left\vert T\left(  a\right)  \right\vert -c_{\mathbf{0}}^{\left(  m\right)
}<\left(  1+a\right)  ^{\mathbf{1}}\left(  \frac{a_{1}}{1+a_{1}}+\frac{a_{2}%
}{1+a_{2}}+\ldots+\frac{a_{m}}{1+a_{m}}\right)  .
\]

Compare this bound with the lower bound \ref{av819}.

\item If $\gamma$ contains $n\geq0$ sequences of zeros of length $\geq2$ then
$c_{\gamma}^{\left(  m\right)  }\leq\left(  \frac{3}{4}\right)  ^{n}$.

\item If $\gamma$ contains $n\geq0$ sequences of zeros of length $\geq1$ then
$c_{\gamma}^{\left(  m\right)  }\geq\frac{1}{2^{m-\left\vert \gamma\right\vert
-n}}$.

\item $\min\limits_{\left\vert \gamma\right\vert =n}c_{\gamma}^{\left(
m\right)  }=\frac{1+m-n}{2^{m-n}}=c_{\beta}^{\left(  m\right)  }$ when
$\beta=0_{m-n}1_{n}$.
\end{enumerate}
\end{corollary}

\begin{proof}
First note the properties of the function $\frac{1+x}{2^{x}}$ described in
Remark \ref{vRem_1_plus_X_div_2^X}.\medskip

\textbf{Part 1 }Equation \ref{av052} of the last theorem and the monotonicity
of $\frac{1+x}{2^{x}}$ implies $c_{\gamma}^{\left(  m\right)  }\leq1$.

Also, \ref{av052} and the subadditivity of $\frac{1+x}{2^{x}}$ implies
\[
c_{\gamma}^{\left(  m\right)  }=\frac{1+\left\Vert z^{\left(  1\right)
}\right\Vert }{2^{\left\Vert z^{\left(  1\right)  }\right\Vert }}\times
\ldots\times\frac{1+\left\Vert z^{\left(  k\right)  }\right\Vert
}{2^{\left\Vert z^{\left(  k\right)  }\right\Vert }}>\frac{1+\left\Vert
z^{\left(  1\right)  }\right\Vert +\ldots+\left\Vert z^{\left(  k\right)
}\right\Vert }{2^{\left\Vert z^{\left(  1\right)  }\ldots z^{\left(  k\right)
}\right\Vert }}=\frac{1+\left\Vert z^{\left(  1\right)  }\ldots z^{\left(
k\right)  }\right\Vert }{2^{\left\Vert z^{\left(  1\right)  }\ldots z^{\left(
k\right)  }\right\Vert }},
\]

where $z^{\left(  1\right)  }\ldots z^{\left(  k\right)  }$ indicates the
concatenation of the sequences.

But $\left\Vert z^{\left(  1\right)  }\ldots z^{\left(  k\right)  }\right\Vert
+\left\vert \gamma\right\vert \leq m$ so if $\gamma$ contains sequences of
zeros having length $\geq2$, \ref{av598} implies that,
\begin{equation}
c_{\gamma}^{\left(  m\right)  }\geq\frac{1+m-\left\vert \gamma\right\vert
}{2^{m-\left\vert \gamma\right\vert }}\geq\frac{1+m}{2^{m}}=c_{\mathbf{0}%
}^{\left(  m\right)  }.\label{av731}%
\end{equation}
\medskip

\textbf{Part 2} From \ref{av599}, Property \ref{vDef_elem_sym_funct} and Part
1 of this corollary:%
\[
c_{\mathbf{0}}^{\left(  m\right)  }\left(  1+a\right)  ^{\mathbf{1}%
}=c_{\mathbf{0}}^{\left(  m\right)  }\sum_{k=0}^{m}\sum
\limits_{\substack{\beta\leq\mathbf{1} \\\left\vert \beta\right\vert
=k}}a^{\beta}<T\left(  a\right)  =\sum_{k=0}^{m}\sum\limits_{\substack{\beta
\leq\mathbf{1} \\\left\vert \beta\right\vert =k}}c_{\beta}^{\left(  m\right)
}a^{\beta}<\sum_{k=0}^{m}\sum\limits_{\substack{\beta\leq\mathbf{1}
\\\left\vert \beta\right\vert =k}}a^{\beta}=\left(  1+a\right)  ^{\mathbf{1}}.
\]
\medskip

\textbf{Part 3} If $a>0$ then%
\[
\left\vert T\left(  a\right)  \right\vert -c_{\mathbf{0}}^{\left(  m\right)
}=\sum_{k=1}^{m}\sum\limits_{\substack{\beta\leq\mathbf{1} \\\left\vert
\beta\right\vert =k}}c_{\beta}^{\left(  m\right)  }a^{\beta}<\sum_{k=1}%
^{m}\sum\limits_{\substack{\beta\leq\mathbf{1} \\\left\vert \beta\right\vert
=k}}a^{\beta}=\sum_{k=0}^{m}\sum\limits_{\substack{\beta\leq\mathbf{1}
\\\left\vert \beta\right\vert =k}}a^{\beta}-1=\left(  1+a\right)
^{\mathbf{1}}-1,
\]

But
\begin{align*}
\left(  1+a\right)  ^{\mathbf{1}}-1  & =\left(  1+a_{1}\right)  \left(
1+a_{2}\right)  \ldots\left(  1+a_{m}\right)  -1\\
& =\left(  1+a_{1}\right)  \left(  1+a_{2}\right)  \ldots\left(
1+a_{m}\right)  -\left(  1+a_{1}\right)  \left(  1+a_{2}\right)  \ldots\left(
1+a_{m-1}\right)  +\\
& +\left(  1+a_{1}\right)  \left(  1+a_{2}\right)  \ldots\left(
1+a_{m-1}\right)  -\left(  1+a_{1}\right)  \left(  1+a_{2}\right)
\ldots\left(  1+a_{m-2}\right)  +\ldots\\
& +\left(  1+a_{1}\right)  -1\\
& =\left\{  \left(  1+a_{1}\right)  \left(  1+a_{2}\right)  \ldots\left(
1+a_{m-1}\right)  \right\}  a_{m}+\\
& +\left\{  \left(  1+a_{1}\right)  \left(  1+a_{2}\right)  \ldots\left(
1+a_{m-2}\right)  \right\}  a_{m-1}+\ldots+\left(  1+a_{1}\right)  a_{2}%
+a_{1},
\end{align*}

so that%
\[
\frac{\left(  1+a\right)  ^{\mathbf{1}}-1}{\left(  1+a\right)  ^{\mathbf{1}}%
}<\frac{a_{m}}{1+a_{m}}+\frac{a_{m-1}}{1+a_{m-1}}+\ldots+\frac{a_{1}}{1+a_{1}%
},
\]

and hence%
\[
\left\vert T\left(  a\right)  \right\vert -c_{\mathbf{0}}^{\left(  m\right)
}<\left(  1+a\right)  ^{\mathbf{1}}\frac{\left(  1+a\right)  ^{\mathbf{1}}%
-1}{\left(  1+a\right)  ^{\mathbf{1}}}=\left(  1+a\right)  ^{\mathbf{1}%
}\left(  \frac{a_{1}}{1+a_{1}}+\frac{a_{2}}{1+a_{2}}+\ldots+\frac{a_{m}%
}{1+a_{m}}\right)  .
\]
\medskip\medskip

\textbf{Part 4} Suppose the sequences of zeros of length $\geq2$ are $\left\{
z^{\left(  j\right)  }\right\}  _{j=1}^{n}$. Then \ref{av052} and the
monotonicity of $\frac{1+x}{2^{x}}$ imply the bound:
\[
c_{\gamma}^{\left(  m\right)  }=\frac{1+\left\Vert z^{\left(  1\right)
}\right\Vert }{2^{\left\Vert z^{\left(  1\right)  }\right\Vert }}\times
\ldots\times\frac{1+\left\Vert z^{\left(  n\right)  }\right\Vert
}{2^{\left\Vert z^{\left(  n\right)  }\right\Vert }}\leq\left(  \frac{3}%
{4}\right)  ^{n}.
\]

\textbf{Part 5} Follows directly from Remark
\ref{vRem_Thm_formula_for_C^(m)_idx} to Theorem
\ref{vThm_formula_for_C^(m)_idx}.\medskip

\textbf{Part 6} If there are no sequences of zeros of length $\geq2$ in
$\gamma$ then by Theorem \ref{vThm_formula_for_C^(m)_idx}, $c_{\gamma
}^{\left(  m\right)  }=1$.

Let $q=\left\Vert z^{\left(  1\right)  }\right\Vert +\left\Vert z^{\left(
2\right)  }\right\Vert +\ldots+\left\Vert z^{\left(  k\right)  }\right\Vert
=\left\Vert z^{\left(  1\right)  }z^{\left(  2\right)  \ldots}z^{\left(
k\right)  }\right\Vert $. Then the subadditivity of $\frac{1+x}{2^{x}}$
implies
\begin{align*}
c_{\gamma}^{\left(  m\right)  }  & =\frac{1+\left\Vert z^{\left(  1\right)
}\right\Vert }{2^{\left\Vert z^{\left(  1\right)  }\right\Vert }}%
\frac{1+\left\Vert z^{\left(  2\right)  }\right\Vert }{2^{\left\Vert
z^{\left(  2\right)  }\right\Vert }}\times\ldots\times\frac{1+\left\Vert
z^{\left(  k\right)  }\right\Vert }{2^{\left\Vert z^{\left(  k\right)
}\right\Vert }}\\
& \geq\frac{1+\left\Vert z^{\left(  1\right)  }\right\Vert +\left\Vert
z^{\left(  2\right)  }\right\Vert +\ldots+\left\Vert z^{\left(  k\right)
}\right\Vert }{2^{\left\Vert z^{\left(  1\right)  }\right\Vert +\left\Vert
z^{\left(  2\right)  }\right\Vert +\ldots+\left\Vert z^{\left(  k\right)
}\right\Vert }}\\
& =\frac{1+\left\Vert z^{\left(  1\right)  }z^{\left(  2\right)  \ldots
}z^{\left(  k\right)  }\right\Vert }{2^{\left\Vert z^{\left(  1\right)
}z^{\left(  2\right)  \ldots}z^{\left(  k\right)  }\right\Vert }}\\
& \geq\frac{1+q+\left(  m-q-\left\vert \gamma\right\vert \right)
}{2^{q+\left(  m-q-\left\vert \gamma\right\vert \right)  }}\\
& =\frac{1+m-\left\vert \gamma\right\vert }{2^{m-\left\vert \gamma\right\vert
}}.
\end{align*}

The sequence $\beta=0_{m-\left\vert \gamma\right\vert }1_{\left\vert
\gamma\right\vert }$ satisfies%
\[
c_{\gamma}^{\left(  m\right)  }\geq c_{\beta}^{\left(  m\right)  }%
=\frac{1+m-\left\vert \gamma\right\vert }{2^{m-\left\vert \gamma\right\vert }%
},\quad\left\vert \beta\right\vert =\left\vert \gamma\right\vert .
\]

Thus $\min\limits_{\left\vert \gamma\right\vert =n}c_{\gamma}^{\left(
m\right)  }=\frac{1+m-n}{2^{m-n}}=c_{\beta}^{\left(  m\right)  }$.
\end{proof}

A result regarding upper bounds for $c_{\gamma}^{\left(  m\right)  }$ was
proven as part 4 of Corollary \ref{vCor_2_formula_for_C^(m)_gam}. We now prove
a more precise result which is the upper bound equivalent of :the problem
$\min\limits_{\left\vert \gamma\right\vert =n}c_{\gamma}^{\left(  m\right)  }$
solved in part 6 of Corollary \ref{vCor_2_formula_for_C^(m)_gam}.

\begin{corollary}
\label{vCor_3_formula_for_C^(m)_gam}\textbf{Solution of the problem}
$\max\limits_{\left\vert \gamma\right\vert =k}c_{\gamma}^{\left(  m\right)  }%
$. There are two cases:

\begin{description}
\item[Case $m/2\leq k\leq m$]
\begin{equation}
\max\limits_{\left\vert \gamma\right\vert =k}c_{\gamma}^{\left(  m\right)
}=1=c_{\beta}^{\left(  m\right)  }\text{ }when\text{ }\beta=\left(  01\right)
_{m-k}1_{2k-m}.\label{av07}%
\end{equation}

\item[\textbf{Case} $0\leq k<m/2$] Define $\underline{\mu}=\left\lfloor
\frac{m+1}{k+1}\right\rfloor $, $\overline{\mu}=\left\lceil \frac{m+1}%
{k+1}\right\rceil $, $k^{\prime\prime}=\operatorname*{rem}\frac{m+1}{k+1}$ and
$k^{\prime}=k+1-k^{\prime\prime}$. Then
\begin{equation}
\max\limits_{\left\vert \gamma\right\vert =k}c_{\gamma}^{\left(  m\right)
}=\frac{1}{2^{m-k}}\underline{\mu}^{k^{\prime}}\overline{\mu}^{k^{\prime
\prime}}=c_{\beta}^{\left(  m\right)  }\text{ }when\text{ }\beta=\left(
0_{\underline{\mu}-1}1\right)  _{k^{\prime}}\left(  0_{\overline{\mu}%
}1\right)  _{k^{\prime\prime}-1}0_{\overline{\mu}},\label{av071}%
\end{equation}

and we have the upper bound%
\begin{equation}
\max\limits_{\left\vert \gamma\right\vert =k}c_{\gamma}^{\left(  m\right)
}\leq\frac{1}{2^{m-k}}\left(  \frac{m+1}{k+1}\right)  ^{k+1}.\label{av072}%
\end{equation}

\end{description}
\end{corollary}

\begin{proof}
\fbox{Case $m/2\leq k\leq m$} When $\beta=\left(  01\right)  _{m-k}1_{2k-m}$
formula \ref{av052} implies $c_{\beta}^{\left(  m\right)  }=1$ and by part 1
of Corollary \ref{vCor_2_formula_for_C^(m)_gam}, $\max\limits_{\beta}c_{\beta
}^{\left(  m\right)  }=1$.\smallskip

\fbox{Case $0\leq k<m/2$} Here. $k<m-k$ so there are more zeros than ones. Now
$0_{m}\leq\gamma\leq1_{m}$ and $\left\vert \gamma\right\vert =k$. Suppose
$\gamma$ has ones in positions $n_{1}<n_{2}<\ldots<n_{k}$. Then from Remark
\ref{vRem_Thm_formula_for_C^(m)_idx} to Theorem
\ref{vThm_formula_for_C^(m)_idx},%
\begin{align*}
c_{\gamma}^{\left(  m\right)  }  & =\frac{1+\left\Vert \zeta^{\left(
1\right)  }\right\Vert }{2^{\left\Vert \zeta^{\left(  1\right)  }\right\Vert
}}\frac{1+\left\Vert \zeta^{\left(  2\right)  }\right\Vert }{2^{\left\Vert
\zeta^{\left(  2\right)  }\right\Vert }}\times\ldots\times\frac{1+\left\Vert
\zeta^{\left(  j\right)  }\right\Vert }{2^{\left\Vert \zeta^{\left(  j\right)
}\right\Vert }}\times\ldots\times\frac{1+\left\Vert \zeta^{\left(  k\right)
}\right\Vert }{2^{\left\Vert \zeta^{\left(  k\right)  }\right\Vert }}%
\frac{1+\left\Vert \zeta^{\left(  k+1\right)  }\right\Vert }{2^{\left\Vert
\zeta^{\left(  k=1\right)  }\right\Vert }}\\
& =\frac{1+\left(  n_{1}-1\right)  }{2^{n_{1}-1}}\frac{1+\left(  n_{2}%
-n_{1}-1\right)  }{2^{n_{2}-n_{1}-1}}\times\ldots\times\frac{1+\left(
n_{j}-n_{j-1}-1\right)  }{2^{n_{j}-n_{j-1}-1}}\times\ldots\\
& \qquad\qquad\qquad\times\frac{1+\left(  n_{k}-n_{k-1}-1\right)  }%
{2^{n_{k}-n_{k-1}-1}}\frac{1+\left(  m-n_{k}\right)  }{2^{m-n_{k}}}\\
& =\frac{n_{1}\left(  n_{2}-n_{1}\right)  \ldots\left(  n_{j}-n_{j-1}\right)
\ldots\left(  n_{k}-n_{k-1}\right)  \left(  1+m-n_{k}\right)  }{2^{\left(
n_{1}-1\right)  +\left(  n_{2}-n_{1}-1\right)  +\ldots+\left(  n_{j}%
-n_{j-1}-1\right)  +\ldots+\left(  n_{k}-n_{k-1}-1\right)  +\left(
m-n_{k}\right)  }}\\
& =\frac{1}{2^{m-k}}n_{1}\left(  n_{2}-n_{1}\right)  \ldots\left(
n_{j}-n_{j-1}\right)  \ldots\left(  n_{k}-n_{k-1}\right)  \left(
1+m-n_{k}\right)  .
\end{align*}

Define $p=\left(  p\right)  _{j=1}^{k+1}$, $p_{1}=n_{1}$, $p_{j}=n_{j}%
-n_{j-1}$ for $2\leq j\leq k$, and $p_{k+1}=1+m-n_{k}$ so that
\begin{equation}
c_{\gamma}^{\left(  m\right)  }=\frac{1}{2^{m-k}}p_{1}p_{2}\ldots
p_{k+1},\label{av060}%
\end{equation}
$p_{j}\geq1$ and $\left\vert p\right\vert =1+m$. Thus the problem
$\max\limits_{\left\vert \gamma\right\vert =k}c_{\gamma}^{\left(  m\right)  }$
is equivalent to the \textbf{nonlinear integer programming problem}:
\begin{equation}
\max_{\substack{p\geq\mathbf{1} \\\left\vert p\right\vert =1+m}}p_{1}%
p_{2}\ldots p_{k+1}.\label{av058}%
\end{equation}

Now if there exists $p_{i}\leq p_{j}-2$ then the replacement: $p_{i}%
\rightarrow p_{i}+1$ and $p_{j}\rightarrow p_{j}-1$, does not change
$\left\vert p\right\vert $ but increases $p_{1}p_{2}\ldots p_{k+1}$. Repeating
this process until exhaustion implies that there exist integers
$\underline{\mu}$ and $\overline{\mu}$ such that $\underline{\mu}%
=\overline{\mu}-1$ and $1\leq k^{\prime},k^{\prime\prime}<k+1$, or
$\underline{\mu}=\overline{\mu}$ and which satisfy%
\begin{align}
p_{j}  & =\left\{
\begin{array}
[c]{ll}%
\underline{\mu}, & 1\leq j\leq k^{\prime},\\
\overline{\mu}, & k^{\prime}+1\leq j\leq k+1,
\end{array}
\right. \label{av066}\\
\operatorname*{argmax}\limits_{_{\substack{p\geq\mathbf{1} \\\left\vert
p\right\vert =1+m}}}p_{1}p_{2}\ldots p_{k+1}  & =\left(  \underset{k^{\prime
}\text{ }repeats}{\underbrace{\underline{\mu},\ldots,\underline{\mu}}%
},\underset{k^{\prime\prime}\text{ }repeats}{\underbrace{\overline{\mu}%
,\ldots,\overline{\mu}}}\right)  =\left(  \underline{\mu}\right)  _{k^{\prime
}}\left(  \overline{\mu}\right)  _{k^{\prime\prime}},\nonumber\\
\max_{\substack{p\geq\mathbf{1} \\\left\vert p\right\vert =1+m}}p_{1}%
p_{2}\ldots p_{k+1}  & =\underline{\mu}^{k^{\prime}}\overline{\mu}%
^{k^{\prime\prime}},\nonumber\\
k^{\prime}+k^{\prime\prime}  & =k+1,\label{av059}\\
\underline{\mu}k^{\prime}+\overline{\mu}k^{\prime\prime}  & =m+1.\label{av06}%
\end{align}

From the definition of the $p_{j}$, $\left\Vert \zeta^{\left(  1\right)
}\right\Vert =n_{1}-1=p_{1}-1,$ $\left\Vert \zeta^{\left(  j\right)
}\right\Vert =n_{j}-n_{j-1}-1=p_{j}-1$ and $\left\Vert \zeta^{\left(
k+1\right)  }\right\Vert =m-n_{k}=p_{k+1}-1$ i.e.%
\begin{equation}
n_{j}=p_{1}+p_{2}+\ldots+p_{j},\quad1\leq j\leq k,\label{av063}%
\end{equation}

and%
\begin{equation}
\left\Vert \zeta^{\left(  j\right)  }\right\Vert =p_{j}-1,\quad1\leq j\leq
k+1,\label{av064}%
\end{equation}

so that from \ref{av060},%
\begin{equation}
\max\limits_{\left\vert \gamma\right\vert =k}c_{\gamma}^{\left(  m\right)
}=\frac{1}{2^{m-k}}\underline{\mu}^{k^{\prime}}\overline{\mu}^{k^{\prime
\prime}},\label{av067}%
\end{equation}

and%
\[
\left(  0_{p_{1}-1}1\right)  \left(  0_{p_{2}-1}1\right)  \ldots\left(
0_{p_{k}-1}1\right)  0_{p_{k+1}}\in\operatorname*{argmax}\limits_{_{\left\vert
\gamma\right\vert =k}}c_{\gamma}^{\left(  m\right)  }.
\]

Substituting \ref{av066} in the last expression now gives%
\[
\underset{k^{\prime}\text{ }repetitions}{\underbrace{\left(  0_{\underline{\mu
}-1}1\right)  \left(  0_{\underline{\mu}-1}1\right)  \ldots\left(
0_{\underline{\mu}-1}1\right)  }}\underset{k^{\prime\prime}-1\text{
}repetitions}{\text{ }\underbrace{\left(  0_{\overline{\mu}}1\right)  \left(
0_{\overline{\mu}}1\right)  \ldots\left(  0_{\overline{\mu}}1\right)  }\text{
}}0_{\overline{\mu}}\in\operatorname*{argmax}\limits_{_{\left\vert
\gamma\right\vert =k}}c_{\gamma}^{\left(  m\right)  },
\]

or more compactly%
\begin{equation}
\left(  0_{\underline{\mu}-1}1\right)  _{k^{\prime}}\left(  0_{\overline{\mu}%
}1\right)  _{k^{\prime\prime}-1}0_{\overline{\mu}}\in\operatorname*{argmax}%
\limits_{_{\left\vert \gamma\right\vert =k}}c_{\gamma}^{\left(  m\right)
},\label{av068}%
\end{equation}

where, since we can have $k^{\prime\prime}-1$, we define $\delta_{0}=1$ for
any sequence $\delta$. Expressions \ref{av067} and \ref{av068} can now be
written%
\begin{equation}
\max\limits_{\left\vert \gamma\right\vert =k}c_{\gamma}^{\left(  m\right)
}=\frac{1}{2^{m-k}}\underline{\mu}^{k^{\prime}}\overline{\mu}^{k^{\prime
\prime}}=c_{\beta}^{\left(  m\right)  },\quad\beta=\left(  0_{\underline{\mu
}-1}1\right)  _{k^{\prime}}\left(  0_{\overline{\mu}}1\right)  _{k^{\prime
\prime}-1}0_{\overline{\mu}}.\label{av069}%
\end{equation}

Suppose $\underline{\mu}=\overline{\mu}$. Then \ref{av059} and \ref{av06}
imply $\underline{\mu}=\frac{m+1}{k+1}$ and so $\frac{m+1}{k+1}$ is an
integer. In fact, we now show that if $\frac{m+1}{k+1}$ is an integer then
$\underline{\mu}=\overline{\mu}=\frac{m+1}{k+1}$. If $\underline{\mu
}=\overline{\mu}-1$ and $k^{\prime},k^{\prime\prime}\geq1$ then from
\ref{av059} and \ref{av06}, $m+1=\underline{\mu}k^{\prime}+\overline{\mu
}k^{\prime\prime}=\underline{\mu}k^{\prime}+\left(  1+\underline{\mu}\right)
k^{\prime\prime}=\underline{\mu}\left(  k+1\right)  +k^{\prime\prime}$ i.e.
$\frac{m+1}{k+1}=\underline{\mu}+\frac{k^{\prime\prime}}{k+1}$. But $1\leq
k^{\prime\prime}<k+1$ so $\frac{m+1}{k+1}=\underline{\mu}$ as claimed.
Observing that the mean of $p$ is $\overline{p}=\frac{m+1}{k+1}$, this
motivates us to consider two subcases:\smallskip

\quad\fbox{\textbf{Subcase} $\overline{p}=\frac{m+1}{k+1}$ is an integer} We
have just shown that:%
\begin{equation}
p_{j}=\underline{\mu}=\overline{\mu}=\frac{m+1}{k+1},\quad1\leq j\leq
k+1,\label{av061}%
\end{equation}

Since $2k<m$: $\underline{\mu}=\frac{1+m}{1+k}>\frac{1+m}{1+m/2}%
=\frac{1+m/2+m/2}{1+m/2}=1+\frac{m/2}{1+m/2}>1$ and so $\underline{\mu}\geq2$
which from \ref{av064} implies that $\left\Vert \zeta^{\left(  j\right)
}\right\Vert =\underline{\mu}-1\geq1$. Expression \ref{av069} becomes
\begin{equation}
\max\limits_{\left\vert \gamma\right\vert =k}c_{\gamma}^{\left(  m\right)
}=\frac{1}{2^{m-k}}\underline{\mu}^{k+1}=c_{\beta}^{\left(  m\right)  }%
,\quad\beta=\left(  0_{\underline{\mu}-1}1\right)  _{k}0_{\underline{\mu}%
-1}.\label{av065}%
\end{equation}

\quad\fbox{\textbf{Subcase} $\overline{p}=\frac{m+1}{k+1}$is not an integer}
In this case we know from above that $\underline{\mu}=\overline{\mu}-1$ and
$1\leq k^{\prime},k^{\prime\prime}<k+1$. Further, from \ref{av059} and
\ref{av06},
\[
m+1=\underline{\mu}k^{\prime}+\overline{\mu}k^{\prime\prime}=\underline{\mu
}k^{\prime}+\left(  1+\underline{\mu}\right)  k^{\prime\prime}=\underline{\mu
}\left(  k+1\right)  +k^{\prime\prime},
\]

so that%
\begin{align}
\underline{\mu}  & =\left\lfloor \frac{m+1}{k+1}\right\rfloor \text{
}and\text{ }\overline{\mu}=\underline{\mu}+1,\label{av062}\\
k^{\prime\prime}  & =\operatorname*{rem}\frac{m+1}{k+1}\text{ }and\text{
}k^{\prime}=k+1-k^{\prime\prime},\label{av070}%
\end{align}

Since $2k<m$ and $k^{\prime\prime}\leq k$: $\underline{\mu}=\frac
{1+m-k^{\prime\prime}}{1+k}>\frac{1+m-k^{\prime\prime}}{1+m/2}=\frac
{1+m/2+m/2-k^{\prime\prime}}{1+m/2}=1+\frac{m/2-k^{\prime\prime}}{1+m/2}>1 $
and so $\underline{\mu}\geq2$. The equivalent of \ref{av065} here is
\ref{av069}. Finally we observe that for both subcases $\underline{\mu
}=\left\lfloor \frac{m+1}{k+1}\right\rfloor $ and $\overline{\mu}=\left\lceil
\frac{m+1}{k+1}\right\rceil $ so that \ref{av071} holds.

It remains to prove the upper bound \ref{av072}. We do this using the
\textbf{theorem of weighted means}\textit{\ }which compares the weighted
geometric and arithmetic means. In fact, regarding \ref{av071} and using
equations \ref{av059} and \ref{av06} we get
\[
\left(  \underline{\mu}^{k^{\prime}}\overline{\mu}^{k^{\prime\prime}}\right)
^{\frac{1}{k^{\prime}+k^{\prime\prime}}}\leq\frac{k^{\prime}\underline{\mu
}+k^{\prime\prime}\overline{\mu}}{k^{\prime}+k^{\prime\prime}}=\frac{m+1}%
{k+1},
\]
so that%
\[
\max\limits_{\left\vert \gamma\right\vert =k}c_{\gamma}^{\left(  m\right)
}=\frac{1}{2^{m-k}}\underline{\mu}^{k^{\prime}}\overline{\mu}^{k^{\prime
\prime}}\leq\frac{1}{2^{m-k}}\left(  \frac{m+1}{k+1}\right)  ^{k+1},
\]
as required.
\end{proof}

\section{Uniform upper bounds for $T^{-1}\beta$ and $\left\vert \mu\left(
f\right)  \right\vert _{\max}$\label{Sect_bound_invTb}}

In this section we will proof Theorem \ref{vThm_inequal_Sigma_C^m_1minusBeta}
which is applied in the succeeding subsections to prove the estimate
\ref{av035} which is a uniform upper bound for $\left\vert T^{-1}%
\beta\right\vert _{\max}$. This bound is applied in turn to \ref{av007} to
obtain the upper bound \ref{av087} for $\left\vert \mu\left(  f\right)
\right\vert _{\max}$.

\subsection{Bounds for $\left(  T^{-1}\beta\right)  _{1}$}

From \ref{av725}, \ref{av552} and \ref{av533},%
\begin{align}
\left(  T^{-1}\beta\right)  _{1}  & =\left(  -1\right)  ^{1-1}\frac
{1}{\left\vert T\right\vert }\left(  \left(  -1\right)  ^{N}M_{1,N-1}%
-M_{1,1}\right) \nonumber\\
& =\frac{1}{\left\vert T\right\vert }\left(  \left(  -1\right)  ^{N}\left(
-\frac{1}{2}\right)  ^{N-2}-\left\vert T_{2:N-1}\right\vert \right)
\nonumber\\
& =\frac{1}{\left\vert T\right\vert }\left(  \left(  \frac{1}{2}\right)
^{N-2}-\left\vert T_{2:N-1}\right\vert \right)  .\label{av032}%
\end{align}

Set
\begin{equation}
x=\frac{1}{v}\Delta_{1}X=\frac{1}{v}\left(  x^{\left(  k+1\right)
}-x^{\left(  k\right)  }\right)  _{k=1}^{N-1}.\label{av012}%
\end{equation}

From \ref{av735} and then using \ref{av598} and \ref{av595} i.e.
$c_{\mathbf{0}}^{\left(  m\right)  }=\frac{m+1}{2^{m}}$, $c_{\mathbf{1}%
}^{\left(  m\right)  }=1 $, we calculate for $N\geq3$:%
\begin{align}
\left(  T^{-1}\beta\right)  _{1}  & =\frac{1}{\left\vert T\right\vert }\left(
\frac{1}{2^{N-2}}-\left\vert T_{2:N-1}\right\vert \right) \nonumber\\
& =\frac{\frac{1}{2^{N-2}}-\sum\limits_{k=0}^{N-2}\sum\limits_{\gamma
\leq\mathbf{1,}\left\vert \gamma\right\vert =k}c_{\gamma}^{\left(  N-2\right)
}x_{2:N-1}^{\gamma}}{\sum\limits_{k=0}^{N-1}\sum\limits_{\beta\leq
\mathbf{1,}\left\vert \beta\right\vert =k}c_{\beta}^{\left(  N-1\right)
}x^{\beta}}\nonumber\\
& =\frac{\frac{1}{2^{N-2}}-\sum\limits_{k=0}\sum\limits_{\gamma\leq
\mathbf{1,}\left\vert \gamma\right\vert =k}c_{\gamma}^{\left(  N-2\right)
}x_{2:N-1}^{\gamma}-\sum\limits_{k=1}^{N-2}\sum\limits_{\gamma\leq
\mathbf{1},\left\vert \gamma\right\vert =k}c_{\gamma}^{\left(  N-2\right)
}x_{2:N-1}^{\gamma}}{\sum\limits_{k=0}\sum\limits_{\beta\leq\mathbf{1}%
,\left\vert \beta\right\vert =k}c_{\beta}^{\left(  N-1\right)  }x^{\beta}%
+\sum\limits_{k=1}^{N-1}\sum\limits_{\beta\leq\mathbf{1},\left\vert
\beta\right\vert =k}c_{\beta}^{\left(  N-1\right)  }x^{\beta}}\nonumber\\
& =\frac{\frac{1}{2^{N-2}}-c_{\mathbf{0}}^{\left(  N-2\right)  }%
-\sum\limits_{k=1}^{N-2}\sum\limits_{\gamma\leq\mathbf{1,}\left\vert
\gamma\right\vert =k}c_{\gamma}^{\left(  N-2\right)  }x_{2:N-1}^{\gamma}%
}{c_{\mathbf{0}}^{\left(  N-1\right)  }+\sum\limits_{k=1}^{N-1}\sum
\limits_{\beta\leq\mathbf{1,}\left\vert \beta\right\vert =k}c_{\beta}^{\left(
N-1\right)  }x^{\beta}}\nonumber\\
& =\frac{\frac{1}{2^{N-2}}-\frac{N-1}{2^{N-2}}-\sum\limits_{k=1}^{N-2}%
\sum\limits_{\gamma\leq\mathbf{1,}\left\vert \gamma\right\vert =k}c_{\gamma
}^{\left(  N-2\right)  }x_{2:N-1}^{\gamma}}{\frac{N}{2^{N-1}}+\sum
\limits_{k=1}^{N-1}\sum\limits_{\beta\leq\mathbf{1,}\left\vert \beta
\right\vert =k}c_{\beta}^{\left(  N-1\right)  }x^{\beta}}\nonumber\\
& =\frac{-\frac{N-2}{2^{N-2}}-\sum\limits_{k=1}^{N-2}\sum\limits_{\gamma
\leq\mathbf{1,}\left\vert \gamma\right\vert =k}c_{\gamma}^{\left(  N-2\right)
}x_{2:N-1}^{\gamma}}{\frac{N}{2^{N-1}}+\sum\limits_{k=1}^{N-1}\sum
\limits_{\beta\leq\mathbf{1,}\left\vert \beta\right\vert =k}c_{\beta}^{\left(
N-1\right)  }x^{\beta}}\nonumber\\
& =\frac{-\frac{N-1}{2^{N-2}}-\sum\limits_{k=0}^{N-2}\sum\limits_{\gamma
\leq\mathbf{1,}\left\vert \gamma\right\vert =k}c_{\gamma}^{\left(  N-2\right)
}x_{2:N-1}^{\gamma}}{\frac{N}{2^{N-1}}+\sum\limits_{k=1}^{N-2}\sum
\limits_{\beta\leq\mathbf{1,}\left\vert \beta\right\vert =k}c_{\beta}^{\left(
N-1\right)  }x^{\beta}+x^{\mathbf{1}}}\nonumber\\
& =-\frac{\frac{N-1}{2^{N-2}}+\sum\limits_{k=0}^{N-2}\sum\limits_{\gamma
\leq\mathbf{1,}\left\vert \gamma\right\vert =k}c_{\gamma}^{\left(  N-2\right)
}x_{2:N-1}^{\gamma}}{\frac{N}{2^{N-1}}+\sum\limits_{k=1}^{N-2}\sum
\limits_{\beta\leq\mathbf{1,}\left\vert \beta\right\vert =k}c_{\beta}^{\left(
N-1\right)  }x^{\beta}+x^{\mathbf{1}}}\label{av753}\\
& =-\frac{\frac{N-1}{2^{N-2}}+\sum\limits_{k=1}^{N-2}\sum\limits_{\gamma
\leq\mathbf{1,}\left\vert \gamma\right\vert =k-1}c_{\gamma}^{\left(
N-2\right)  }x_{2:N-1}^{\gamma}}{\frac{N}{2^{N-1}}+\sum\limits_{k=1}^{N-2}%
\sum\limits_{\beta\leq\mathbf{1,}\left\vert \beta\right\vert =k}c_{\beta
}^{\left(  N-1\right)  }x^{\beta}+x^{\mathbf{1}}},\nonumber
\end{align}

so that%
\begin{equation}
\left(  T^{-1}\beta\right)  _{1}=-\frac{\sum\limits_{k=1}^{N-2}\sum
\limits_{\gamma\leq\mathbf{1,}\left\vert \gamma\right\vert =k-1}c_{\gamma
}^{\left(  N-2\right)  }x_{2:N-1}^{\gamma}+\frac{N-1}{2^{N-2}}}{\frac
{N}{2^{N-1}}+\sum\limits_{k=1}^{N-2}\sum\limits_{\beta\leq\mathbf{1,}%
\left\vert \beta\right\vert =k}c_{\beta}^{\left(  N-1\right)  }x^{\beta
}+x^{\mathbf{1}}}<0,\label{av014}%
\end{equation}

where $x$ is given by \ref{av012}. But from inequality \ref{av020} of
Corollary \ref{vThm_inequal_Sigma_C^m_1minusBeta} below:
\[
\sum\limits_{\beta\leq\mathbf{1,}\left\vert \beta\right\vert =k}c_{\beta
}^{\left(  N-1\right)  }x^{\beta}>\frac{1}{2}\sum\limits_{\gamma
\leq\mathbf{1,}\left\vert \gamma\right\vert =k-1}c_{\gamma}^{\left(
N-2\right)  }x_{2:N-1}^{\gamma},\quad x>0_{N-1},\text{ }1\leq k\leq N-1,
\]

which implies%
\begin{align*}
\frac{N}{2^{N-1}}+\sum\limits_{k=1}^{N-2}\sum\limits_{\beta\leq\mathbf{1,}%
\left\vert \beta\right\vert =k}c_{\beta}^{\left(  N-1\right)  }x^{\beta
}+x^{\mathbf{1}}  & >\frac{N}{2^{N-1}}+\frac{1}{2}\sum\limits_{k=1}^{N-2}%
\sum\limits_{\gamma\leq\mathbf{1,}\left\vert \gamma\right\vert =k-1}c_{\gamma
}^{\left(  N-2\right)  }x_{2:N-1}^{\gamma}+x^{\mathbf{1}}\\
& >\frac{1}{2}\left(  \frac{N-1}{2^{N-2}}+\sum\limits_{k=1}^{N-2}%
\sum\limits_{\gamma\leq\mathbf{1,}\left\vert \gamma\right\vert =k-1}c_{\gamma
}^{\left(  N-2\right)  }x_{2:N-1}^{\gamma}\right)  ,
\end{align*}

and so%
\begin{equation}
0<-\left(  T^{-1}\beta\right)  _{1}<2,\quad N\geq3,\text{ }\rho
>0,\label{av754}%
\end{equation}

where $T=T\left(  x;1\right)  $, $x=\frac{1}{v}\Delta_{1}X$ and from
\ref{av839}, $v=\rho N$.

The first step toward proving Theorem \ref{vThm_inequal_Sigma_C^m_1minusBeta}
are the following estimates:

\begin{theorem}
\label{vThm_C^m_0n_gt_.5C^(m-1)_n}If $0\leq\mu\leq1_{m-1}$ then:

\begin{enumerate}
\item $\frac{1}{2}c_{\mu}^{\left(  m-1\right)  }<c_{0,\mu}^{\left(  m\right)
}\leq c_{\mu}^{\left(  m-1\right)  }$,\quad$m\geq2$.

\item $\frac{1}{2}c_{\lambda,\mu}^{\left(  m-1\right)  }<c_{\lambda,0,\mu
}^{\left(  m\right)  }\leq c_{\lambda,\mu}^{\left(  m-1\right)  }$,\quad
$m\geq3$.

\item $\frac{1}{2}c_{\lambda}^{\left(  m-1\right)  }<c_{\lambda,0}^{\left(
m\right)  }\leq c_{\lambda}^{\left(  m-1\right)  }$,\quad$m\geq2$.

\item $\frac{1}{2^{m-\left\vert \sigma\right\vert }}\leq c_{\sigma}^{\left(
m\right)  }\leq1$,\quad$m\geq1$. (cf. \ref{av731})

The estimate of part1 can be sharpened slightly to:

\item $\frac{1}{2}\frac{m+1}{m}c_{\mu}^{\left(  m-1\right)  }\leq c_{0,\mu
}^{\left(  m\right)  }$,$\quad m\geq2$.

Finally, if $0<\mu\leq1_{m-1}$ then

\item $\frac{1}{2}\frac{m}{m-1}c_{\mu}^{\left(  m-1\right)  }\leq c_{0,\mu
}^{\left(  m\right)  }$,$\quad$when $m\geq2$.
\end{enumerate}
\end{theorem}

\begin{proof}
\textbf{Part 1\smallskip}

\fbox{\textbf{Case 1} $\left(  0,\mu\right)  $ contains no non-trivial
sequence of zeros} This means $\mu$ contains no non-trivial sequence of zeros
and so Theorem \ref{vThm_formula_for_C^(m)_idx} implies $c_{0,\mu}^{\left(
m\right)  }=c_{\mu}^{\left(  m-1\right)  }=1$.\medskip

There remain four cases:\medskip

\fbox{\textbf{Case 2} $\mu=0$} From \ref{av591}, $c_{0,0}^{\left(  2\right)
}=3/4$ and $c_{0}^{\left(  1\right)  }=1$ so that
\[
c_{0,\mu}^{\left(  m\right)  }=\frac{3}{4}c_{\mu}^{\left(  m-1\right)  }.
\]
\medskip

Again with reference to Theorem \ref{vThm_formula_for_C^(m)_idx}: if $\left(
0,\mu\right)  $ contains the non-trivial sequences $z^{\left(  1\right)  }%
$,\ldots,$z^{\left(  k\right)  }$ of zeros then%
\begin{equation}
c_{0,\mu}^{\left(  m\right)  }=\frac{1+\left\Vert z^{\left(  1\right)
}\right\Vert }{2^{\left\Vert z^{\left(  1\right)  }\right\Vert }}\times
\ldots\times\frac{1+\left\Vert z^{\left(  k\right)  }\right\Vert
}{2^{\left\Vert z^{\left(  k\right)  }\right\Vert }}.\label{av751}%
\end{equation}

\fbox{\textbf{Case 3} $\mu=\left(  1,\ldots\right)  $} Here $\mu$ has the same
non-trivial sequences of zeros as $\left(  0,\mu\right)  $ and so
\[
c_{0,\mu}^{\left(  m\right)  }=c_{\mu}^{\left(  m-1\right)  }.
\]
\medskip

\fbox{\textbf{Case 4} $\mu=\left(  0,1,\ldots\right)  $} In this case
$z^{\left(  1\right)  }=\left(  0,0\right)  $. Thus $c_{0,\mu}^{\left(
m\right)  }=\frac{3}{4}\frac{1+\left\Vert z^{\left(  2\right)  }\right\Vert
}{2^{\left\Vert z^{\left(  2\right)  }\right\Vert }}\times\ldots\times
\frac{1+\left\Vert z^{\left(  k\right)  }\right\Vert }{2^{\left\Vert
z^{\left(  k\right)  }\right\Vert }}$ and $c_{\mu}^{\left(  m-1\right)
}=\frac{1+\left\Vert z^{\left(  2\right)  }\right\Vert }{2^{\left\Vert
z^{\left(  2\right)  }\right\Vert }}\times\ldots\times\frac{1+\left\Vert
z^{\left(  k\right)  }\right\Vert }{2^{\left\Vert z^{\left(  k\right)
}\right\Vert }}$ and hence
\[
c_{0,\mu}^{\left(  m\right)  }=\frac{3}{4}c_{\mu}^{\left(  m-1\right)  }.
\]
\medskip

\fbox{\textbf{Case 5} $\mu=\left(  0,0,\ldots\right)  $} In this case
$z^{\left(  1\right)  }=\left(  0,0,0,\ldots\right)  $. Thus $c_{\mu}^{\left(
m-1\right)  }=\frac{\left\Vert z^{\left(  1\right)  }\right\Vert
}{2^{\left\Vert z^{\left(  1\right)  }\right\Vert -1}}\frac{1+\left\Vert
z^{\left(  2\right)  }\right\Vert }{2^{\left\Vert z^{\left(  2\right)
}\right\Vert }}\times\ldots\times\frac{1+\left\Vert z^{\left(  k\right)
}\right\Vert }{2^{\left\Vert z^{\left(  k\right)  }\right\Vert }}$ and so
$c_{0,\mu}^{\left(  m\right)  }=\frac{1}{2}\frac{1+\left\Vert z^{\left(
1\right)  }\right\Vert }{\left\Vert z^{\left(  1\right)  }\right\Vert }c_{\mu
}^{\left(  m-1\right)  }=\frac{1}{2}\left(  1+\frac{1}{\left\Vert z^{\left(
1\right)  }\right\Vert }\right)  c_{\mu}^{\left(  m-1\right)  }$.

But $3\leq\left\Vert z^{\left(  1\right)  }\right\Vert \leq m$ so
\begin{equation}
\frac{1}{2}c_{\mu}^{\left(  m-1\right)  }<\frac{m+1}{2m}c_{\mu}^{\left(
m-1\right)  }\leq c_{0,\mu}^{\left(  m\right)  }\leq\frac{2}{3}c_{\mu
}^{\left(  m-1\right)  }.\label{av931}%
\end{equation}

Noting \ref{av932} we can verify that part 1 also holds for $m=2$.\medskip

\textbf{Part 2} \fbox{\textbf{Case 1} $\left(  \lambda,0,\mu\right)  $
contains no non-trivial sequence of zeros} Then
\[
c_{\lambda,0,\mu}^{\left(  m\right)  }=c_{\lambda,\mu}^{\left(  m-1\right)
}=1.
\]
\medskip

Now suppose has a non-trivial sequence of zeros. By Theorem
\ref{vThm_formula_for_C^(m)_idx},%
\begin{equation}
c_{\lambda,0,\mu}^{\left(  m\right)  }=\frac{1+\left\Vert z^{\left(  1\right)
}\right\Vert }{2^{\left\Vert z^{\left(  1\right)  }\right\Vert }}\times
\ldots\times\frac{1+\left\Vert z^{\left(  k\right)  }\right\Vert
}{2^{\left\Vert z^{\left(  k\right)  }\right\Vert }}.\label{av762}%
\end{equation}

We consider six cases:\medskip

\fbox{\textbf{Case 1} $\lambda=\left(  \ldots,1\right)  $, $\mu=\left(
1,\ldots\right)  $} Then $c_{\lambda,0,\mu}^{\left(  m\right)  }%
=c_{\ldots1,0,1\ldots}^{\left(  m\right)  }$ and $c_{\lambda,\mu}^{\left(
m-1\right)  }=c_{\ldots1,1\ldots}^{\left(  m\right)  }$ and $c_{\lambda,0,\mu
}^{\left(  m\right)  }$ and $c_{\lambda,\mu}^{\left(  m-1\right)  }$ both have
the same non-trivial sequence of zeros. Thus by \ref{av762},
\[
c_{\lambda,0,\mu}^{\left(  m\right)  }=c_{\lambda,\mu}^{\left(  m-1\right)  }.
\]
\medskip

\fbox{\textbf{Case 2} $\lambda=\left(  \ldots,1,0\right)  $, $\mu=\left(
1,\ldots\right)  $} Write $\left(  \lambda,0,\mu\right)  =\left(
\ldots1,0,0,1\ldots\right)  =\left(  \ldots1,z^{\left(  n\right)  }%
,1\ldots\right)  $ where $z^{\left(  n\right)  }=\left(  0,0\right)  $, so
that $\left(  \lambda,\mu\right)  =\left(  \ldots1,0,1\ldots\right)  $ and%
\begin{align*}
c_{\lambda,0,\mu}^{\left(  m\right)  } &  =\frac{1+\left\Vert z^{\left(
1\right)  }\right\Vert }{2^{\left\Vert z^{\left(  1\right)  }\right\Vert }%
}\times\ldots\times\frac{1+\left\Vert z^{\left(  n\right)  }\right\Vert
}{2^{\left\Vert z^{\left(  n\right)  }\right\Vert }}\times\ldots\times
\frac{1+\left\Vert z^{\left(  k\right)  }\right\Vert }{2^{\left\Vert
z^{\left(  k\right)  }\right\Vert }}\\
&  =\frac{3}{4}\frac{1+\left\Vert z^{\left(  1\right)  }\right\Vert
}{2^{\left\Vert z^{\left(  1\right)  }\right\Vert }}\times\ldots\times
\frac{1+\left\Vert z^{\left(  n-1\right)  }\right\Vert }{2^{\left\Vert
z^{\left(  n-1\right)  }\right\Vert }}\frac{1+\left\Vert z^{\left(
n+1\right)  }\right\Vert }{2^{\left\Vert z^{\left(  n+1\right)  }\right\Vert
}}\times\ldots\times\frac{1+\left\Vert z^{\left(  k\right)  }\right\Vert
}{2^{\left\Vert z^{\left(  k\right)  }\right\Vert }}\\
&  =\frac{3}{4}c_{\lambda,\mu}^{\left(  m-1\right)  }.
\end{align*}
\medskip

\fbox{\textbf{Case 3} $\lambda=\left(  \ldots,0,0\right)  $, $\mu=\left(
1,\ldots\right)  $} Write $\left(  \lambda,0,\mu\right)  =\left(
\ldots,0,0,0,1\ldots\right)  =\left(  \ldots,z^{\left(  n\right)  }%
,1\ldots\right)  $ where

$z^{\left(  n\right)  }=\left(  0,\ldots,0,0\right)  $, so that $\left(
\lambda,\mu\right)  =\left(  \ldots,z^{\left(  n\right)  }\setminus
0,1\ldots\right)  $ and hence
\begin{align*}
c_{\lambda,\mu}^{\left(  m-1\right)  }  & =\frac{1+\left\Vert z^{\left(
1\right)  }\right\Vert }{2^{\left\Vert z^{\left(  1\right)  }\right\Vert }%
}\times\ldots\times\frac{1+\left\Vert z^{\left(  n\right)  }\setminus
0\right\Vert }{2^{\left\Vert z^{\left(  n\right)  }\setminus0\right\Vert }%
}\times\ldots\times\frac{1+\left\Vert z^{\left(  k\right)  }\right\Vert
}{2^{\left\Vert z^{\left(  k\right)  }\right\Vert }}\\
& =\frac{\frac{1+\left\Vert z^{\left(  n\right)  }\setminus0\right\Vert
}{2^{\left\Vert z^{\left(  n\right)  }\setminus0\right\Vert }}}{\frac
{1+\left\Vert z^{\left(  n\right)  }\right\Vert }{2^{\left\Vert z^{\left(
n\right)  }\right\Vert }}}c_{\lambda,0,\mu}^{\left(  m\right)  }%
=2\frac{\left\Vert z^{\left(  n\right)  }\right\Vert }{1+\left\Vert z^{\left(
n\right)  }\right\Vert }c_{\lambda,0,\mu}^{\left(  m\right)  }.
\end{align*}

But $3\leq\left\Vert z^{\left(  n\right)  }\right\Vert \leq m-1$ and since
$c_{\lambda,0,\mu}^{\left(  m\right)  }=\frac{1}{2}\frac{1+\left\Vert
z^{\left(  n\right)  }\right\Vert }{\left\Vert z^{\left(  n\right)
}\right\Vert }c_{\lambda,\mu}^{\left(  m-1\right)  }=\frac{1}{2}\left(
1+\frac{1}{\left\Vert z^{\left(  n\right)  }\right\Vert }\right)
c_{\lambda,\mu}^{\left(  m-1\right)  }$ we obtain
\begin{equation}
\frac{1}{2}c_{\lambda,\mu}^{\left(  m-1\right)  }<\frac{1}{2}\frac{m}%
{m-1}c_{\lambda,\mu}^{\left(  m-1\right)  }\leq c_{\lambda,0,\mu}^{\left(
m\right)  }\leq\frac{2}{3}c_{\lambda,\mu}^{\left(  m-1\right)  }.\label{av761}%
\end{equation}
\medskip

\fbox{\textbf{Case 4} $\lambda=\left(  \ldots,1\right)  $, $\mu=\left(
0,1,\ldots\right)  $} Write $\left(  \lambda,0,\mu\right)  =\left(
\ldots1,0,0,1\ldots\right)  =\left(  \ldots1,z^{\left(  n\right)  }%
,1\ldots\right)  $ where $z^{\left(  n\right)  }=\left(  0,0\right)  $, so
that $\left(  \lambda,\mu\right)  =\left(  \ldots1,0,1\ldots\right)  $ and%
\begin{align*}
c_{\lambda,0,\mu}^{\left(  m\right)  } &  =\frac{1+\left\Vert z^{\left(
1\right)  }\right\Vert }{2^{\left\Vert z^{\left(  1\right)  }\right\Vert }%
}\times\ldots\times\frac{1+\left\Vert z^{\left(  n\right)  }\right\Vert
}{2^{\left\Vert z^{\left(  n\right)  }\right\Vert }}\times\ldots\times
\frac{1+\left\Vert z^{\left(  k\right)  }\right\Vert }{2^{\left\Vert
z^{\left(  k\right)  }\right\Vert }}\\
&  =\frac{3}{4}\frac{1+\left\Vert z^{\left(  1\right)  }\right\Vert
}{2^{\left\Vert z^{\left(  1\right)  }\right\Vert }}\times\ldots\times
\frac{1+\left\Vert z^{\left(  n-1\right)  }\right\Vert }{2^{\left\Vert
z^{\left(  n-1\right)  }\right\Vert }}\frac{1+\left\Vert z^{\left(
n+1\right)  }\right\Vert }{2^{\left\Vert z^{\left(  n+1\right)  }\right\Vert
}}\times\ldots\times\frac{1+\left\Vert z^{\left(  k\right)  }\right\Vert
}{2^{\left\Vert z^{\left(  k\right)  }\right\Vert }}\\
&  =\frac{3}{4}c_{\lambda,\mu}^{\left(  m-1\right)  }.
\end{align*}
\medskip

\fbox{\textbf{Case 5} $\lambda=\left(  \ldots,1\right)  $, $\mu=\left(
0,0,\ldots\right)  $} In a manner similar to Case 3 we get%
\[
\frac{1}{2}c_{\lambda,\mu}^{\left(  m-1\right)  }<\frac{1}{2}\frac{m}%
{m-1}c_{\lambda,\mu}^{\left(  m-1\right)  }\leq c_{\lambda,0,\mu}^{\left(
m\right)  }\leq\frac{2}{3}c_{\lambda,\mu}^{\left(  m-1\right)  }.
\]
\medskip\medskip

\fbox{\textbf{Case 6} $\lambda=\left(  \ldots,0\right)  $, $\mu=\left(
0,\ldots\right)  $} Write $\left(  \lambda,0,\mu\right)  =\left(
\ldots,0,0,0,\ldots\right)  =\left(  \ldots,z^{\left(  n\right)  }%
,\ldots\right)  $

where $z^{\left(  n\right)  }=\left(  0,\ldots,0,\ldots,0\right)  $, so that
$\left(  \lambda,\mu\right)  =\left(  \ldots,z^{\left(  n\right)  }%
\setminus0,\ldots\right)  $ and hence if%
\[
c_{\lambda,0,\mu}^{\left(  m\right)  }=\frac{1+\left\Vert z^{\left(  1\right)
}\right\Vert }{2^{\left\Vert z^{\left(  1\right)  }\right\Vert }}\times
\ldots\times\frac{1+\left\Vert z^{\left(  n\right)  }\right\Vert
}{2^{\left\Vert z^{\left(  n\right)  }\right\Vert }}\times\ldots\times
\frac{1+\left\Vert z^{\left(  k\right)  }\right\Vert }{2^{\left\Vert
z^{\left(  k\right)  }\right\Vert }},
\]

then%
\begin{align*}
c_{\lambda,\mu}^{\left(  m-1\right)  }  & =\frac{1+\left\Vert z^{\left(
1\right)  }\right\Vert }{2^{\left\Vert z^{\left(  1\right)  }\right\Vert }%
}\times\ldots\times\frac{1+\left\Vert z^{\left(  n\right)  }\setminus
0\right\Vert }{2^{\left\Vert z^{\left(  n\right)  }\setminus0\right\Vert }%
}\times\ldots\times\frac{1+\left\Vert z^{\left(  k\right)  }\right\Vert
}{2^{\left\Vert z^{\left(  k\right)  }\right\Vert }}\\
& =\frac{\frac{1+\left\Vert z^{\left(  n\right)  }\setminus0\right\Vert
}{2^{\left\Vert z^{\left(  n\right)  }\setminus0\right\Vert }}}{\frac
{1+\left\Vert z^{\left(  n\right)  }\right\Vert }{2^{\left\Vert z^{\left(
n\right)  }\right\Vert }}}c_{\lambda,0,\mu}^{\left(  m\right)  }%
=2\frac{\left\Vert z^{\left(  n\right)  }\right\Vert }{1+\left\Vert z^{\left(
n\right)  }\right\Vert }c_{\lambda,0,\mu}^{\left(  m\right)  }.
\end{align*}

But $3\leq\left\Vert z^{\left(  n\right)  }\right\Vert \leq m$ and since
$c_{\lambda,0,\mu}^{\left(  m\right)  }=\frac{1}{2}\frac{1+\left\Vert
z^{\left(  n\right)  }\right\Vert }{\left\Vert z^{\left(  n\right)
}\right\Vert }c_{\lambda,\mu}^{\left(  m-1\right)  }=\frac{1}{2}\left(
\frac{1}{\left\Vert z^{\left(  n\right)  }\right\Vert }+1\right)
c_{\lambda,\mu}^{\left(  m-1\right)  }$ we obtain%
\[
\frac{1}{2}c_{\lambda,\mu}^{\left(  m-1\right)  }<\frac{1}{2}\frac{m+1}%
{m}c_{\lambda,\mu}^{\left(  m-1\right)  }\leq c_{\lambda,0,\mu}^{\left(
m\right)  }\leq\frac{2}{3}c_{\lambda,\mu}^{\left(  m-1\right)  }.
\]
\medskip

\textbf{Part 3} \smallskip

\fbox{\textbf{Case 1} $\left(  \lambda,0\right)  $ contains no non-trivial
sequence of zeros} This means $\lambda$ contains no non-trivial sequence of
zeros and so Theorem \ref{vThm_formula_for_C^(m)_idx} implies $c_{\lambda
,0}^{\left(  m\right)  }=c_{\lambda}^{\left(  m-1\right)  }=1$.\medskip

There remain four cases:\medskip

\fbox{\textbf{Case 2} $\lambda=0$} From \ref{av591}, $c_{0,0}^{\left(
2\right)  }=3/4$ and $c_{0}^{\left(  1\right)  }=1$ so that
\[
c_{0,\lambda}^{\left(  m\right)  }=\frac{3}{4}c_{\lambda}^{\left(  m-1\right)
}.
\]

Again with reference to Theorem \ref{vThm_formula_for_C^(m)_idx}: if $\left(
\lambda,0\right)  $ contains the non-trivial sequences $z^{\left(  1\right)
}$,\ldots,$z^{\left(  k\right)  }$ of zeros then%
\[
c_{\lambda,0}^{\left(  m\right)  }=\frac{1+\left\Vert z^{\left(  1\right)
}\right\Vert }{2^{\left\Vert z^{\left(  1\right)  }\right\Vert }}\times
\ldots\times\frac{1+\left\Vert z^{\left(  k\right)  }\right\Vert
}{2^{\left\Vert z^{\left(  k\right)  }\right\Vert }}.
\]
\medskip

\fbox{\textbf{Case 3} $\lambda=\left(  \ldots,1\right)  $} Here $\lambda$ has
the same non-trivial sequences of zeros as $\left(  \lambda,0\right)  $ and
so
\[
c_{\lambda,0}^{\left(  m\right)  }=c_{\lambda}^{\left(  m-1\right)  }.
\]
\medskip\medskip

\fbox{\textbf{Case 4} $\lambda=\left(  \ldots1,0\right)  $} In this case
$z^{\left(  k\right)  }=\left(  0,0\right)  $. Thus $c_{\lambda,0}^{\left(
m\right)  }=\frac{3}{4}\frac{1+\left\Vert z^{\left(  1\right)  }\right\Vert
}{2^{\left\Vert z^{\left(  1\right)  }\right\Vert }}\times\ldots\times
\frac{1+\left\Vert z^{\left(  k-1\right)  }\right\Vert }{2^{\left\Vert
z^{\left(  k-1\right)  }\right\Vert }}$ and $c_{\lambda}^{\left(  m-1\right)
}=\frac{1+\left\Vert z^{\left(  1\right)  }\right\Vert }{2^{\left\Vert
z^{\left(  1\right)  }\right\Vert }}\times\ldots\times\frac{1+\left\Vert
z^{\left(  k-1\right)  }\right\Vert }{2^{\left\Vert z^{\left(  k-1\right)
}\right\Vert }}$ and hence
\[
c_{\lambda,0}^{\left(  m\right)  }=\frac{3}{4}c_{\lambda}^{\left(  m-1\right)
}.
\]
\medskip

\fbox{\textbf{Case 5} $\lambda=\left(  \ldots,0,0\right)  $} In this case
$z^{\left(  k\right)  }=\left(  \ldots,0,0,0\right)  $. Thus $c_{\lambda
}^{\left(  m-1\right)  }=\frac{1+\left\Vert z^{\left(  1\right)  }\right\Vert
}{2^{\left\Vert z^{\left(  1\right)  }\right\Vert }}\frac{1+\left\Vert
z^{\left(  2\right)  }\right\Vert }{2^{\left\Vert z^{\left(  2\right)
}\right\Vert }}\times\ldots\times\frac{\left\Vert z^{\left(  k\right)
}\right\Vert }{2^{\left\Vert z^{\left(  k\right)  }\right\Vert -1}}$ and so
$c_{\lambda,0}^{\left(  m\right)  }=\frac{1}{2}\frac{1+\left\Vert z^{\left(
k\right)  }\right\Vert }{\left\Vert z^{\left(  k\right)  }\right\Vert
}c_{\lambda}^{\left(  m-1\right)  }=\frac{1}{2}\left(  1+\frac{1}{\left\Vert
z^{\left(  k\right)  }\right\Vert }\right)  c_{\lambda}^{\left(  m-1\right)
}$.

But $3\leq\left\Vert z^{\left(  k\right)  }\right\Vert \leq m$ so%
\[
\frac{1}{2}c_{\lambda}^{\left(  m-1\right)  }<\frac{1+m}{2m}c_{\lambda
}^{\left(  m-1\right)  }\leq c_{\lambda,0}^{\left(  m\right)  }\leq\frac{2}%
{3}c_{\lambda}^{\left(  m-1\right)  }.
\]

The case $m=2$ follows from \ref{av932}.\medskip

\textbf{Part 4} Since $\sigma$ has $\left\vert \sigma\right\vert $ ones and
$m-\left\vert \sigma\right\vert $ zeros it follows directly from the previous
parts and \ref{av595} that $\frac{1}{2^{m-\left\vert \sigma\right\vert }%
}c_{\mathbf{1}}^{\left(  \left\vert \sigma\right\vert \right)  }<c_{\sigma
}^{\left(  m\right)  }\leq1$ and $c_{\mathbf{1}}^{\left(  \left\vert
\sigma\right\vert \right)  }=1$.\medskip

\textbf{Part 5} This follows directly from \ref{av931} in case 5 of part
1.\medskip

\textbf{Part 6} This also follows from \ref{av931} in case 5 of part 1 except
that now $0<\mu$ implies $m\geq4$ and $3\leq\left\Vert z^{\left(  1\right)
}\right\Vert \leq m-1$. The cases $m=2$ and $m=3$ can be verified with
reference to \ref{av932}.
\end{proof}

We can now prove:

\begin{theorem}
\label{vThm_inequal_Sigma_C^m_1minusBeta}Suppose $x=\left(  x^{\prime
},x_{l:l+k-1},x^{\prime\prime}\right)  \in\mathbb{R}_{+}^{m}$ where $m\geq2$
and $1\leq l\leq m-k+1$: note that either $x^{\prime}$ or $x^{\prime\prime}$
may be empty but not both. Then%
\begin{equation}
\sum\limits_{\beta\leq\mathbf{1},\left\vert \beta\right\vert =n}c_{\beta
}^{\left(  m\right)  }x^{\beta}>\frac{1}{2^{k}}\sum\limits_{\gamma
\leq\mathbf{1},\left\vert \gamma\right\vert =n-k}c_{\gamma}^{\left(
m-k\right)  }\left(  x^{\prime},x^{\prime\prime}\right)  ^{\gamma},\quad1\leq
k\leq n.\label{av767}%
\end{equation}

Suppose $1\leq n\leq m$. Then%
\begin{equation}
\sum\limits_{\beta\leq\mathbf{1},\left\vert \beta\right\vert =n}c_{\beta
}^{\left(  m\right)  }x^{\beta}>\frac{1}{2^{k}}\sum\limits_{\gamma
\leq\mathbf{1},\left\vert \gamma\right\vert =n-k}c_{\gamma}^{\left(
m-k\right)  }x_{1:m-k}^{\gamma},\quad1\leq k\leq n,\label{av019}%
\end{equation}

and%
\begin{equation}
\sum\limits_{\beta\leq\mathbf{1},\left\vert \beta\right\vert =n}c_{\beta
}^{\left(  m\right)  }x^{\beta}>\frac{1}{2^{j}}\sum\limits_{\gamma
\leq\mathbf{1},\left\vert \gamma\right\vert =n-j}c_{\gamma}^{\left(
m-j\right)  }x_{j+1:m}^{\gamma},\quad1\leq j\leq n.\label{av020}%
\end{equation}

\end{theorem}

\begin{proof}
Suppose $x=\left(  x^{\prime},x_{l},x^{\prime\prime}\right)  \in\mathbb{R}%
_{+}^{m}$. Note that either $x^{\prime}$ or $x^{\prime\prime}$ may be empty
but not both. Then for $1\leq n\leq m$:%
\begin{align*}
\sum\limits_{\substack{\beta\leq\mathbf{1} \\\left\vert \beta\right\vert
=n}}c_{\beta}^{\left(  m\right)  }x^{\beta}  & =\sum\limits_{\substack{\beta
\leq\mathbf{1} \\\left\vert \beta\right\vert =n}}c_{\beta}^{\left(  m\right)
}\left(  x^{\prime},x_{l},x^{\prime\prime}\right)  ^{\beta}\\
& =\sum\limits_{\substack{\left(  \beta^{\prime},\beta_{l},\beta^{\prime
\prime}\right)  \leq\mathbf{1} \\\left\vert \beta^{\prime}\right\vert
+\beta_{l}+\left\vert \beta^{\prime\prime}\right\vert =n}}c_{\beta^{\prime
},\beta_{l},\beta^{\prime\prime}}^{\left(  m\right)  }\left(  x^{\prime
}\right)  ^{\beta^{\prime}}x_{l}^{\beta_{l}}\left(  x^{\prime\prime}\right)
^{\beta^{\prime\prime}}\\
& =\sum\limits_{\substack{\beta_{l}=1 \\\left(  \beta^{\prime},\beta
^{\prime\prime}\right)  \leq\mathbf{1}_{m-1} \\\left\vert \beta^{\prime
}\right\vert +\beta_{l}+\left\vert \beta^{\prime\prime}\right\vert
=n}}c_{\beta^{\prime},\beta_{l},\beta^{\prime\prime}}^{\left(  m\right)
}\left(  x^{\prime}\right)  ^{\beta^{\prime}}x_{l}^{\beta_{l}}\left(
x^{\prime\prime}\right)  ^{\beta^{\prime\prime}}+\\
& \qquad+\sum\limits_{\substack{\beta_{l}=0 \\\left(  \beta^{\prime}%
,\beta^{\prime\prime}\right)  \leq\mathbf{1}_{m-1} \\\left\vert \beta^{\prime
}\right\vert +\beta_{l}+\left\vert \beta^{\prime\prime}\right\vert
=n}}c_{\beta^{\prime},\beta_{l},\beta^{\prime\prime}}^{\left(  m\right)
}\left(  x^{\prime}\right)  ^{\beta^{\prime}}x_{l}^{\beta_{l}}\left(
x^{\prime\prime}\right)  ^{\beta^{\prime\prime}}\\
& >\sum\limits_{\substack{\beta_{l}=0 \\\left(  \beta^{\prime},\beta
^{\prime\prime}\right)  \leq\mathbf{1}_{m-1} \\\left\vert \beta^{\prime
}\right\vert +\beta_{l}+\left\vert \beta^{\prime\prime}\right\vert
=n}}c_{\beta^{\prime},\beta_{l},\beta^{\prime\prime}}^{\left(  m\right)
}\left(  x^{\prime}\right)  ^{\beta^{\prime}}x_{l}^{\beta_{l}}\left(
x^{\prime\prime}\right)  ^{\beta^{\prime\prime}}\\
& =\sum\limits_{\substack{\left(  \beta^{\prime},\beta^{\prime\prime}\right)
\leq\mathbf{1}_{m-1} \\\left\vert \beta^{\prime}\right\vert +\left\vert
\beta^{\prime\prime}\right\vert =n-1}}c_{\beta^{\prime},0,\beta^{\prime\prime
}}^{\left(  m\right)  }\left(  x^{\prime}\right)  ^{\beta^{\prime}}\left(
x^{\prime\prime}\right)  ^{\beta^{\prime\prime}}.
\end{align*}

Now from parts 1,2 and 3 of Theorem \ref{vThm_C^m_0n_gt_.5C^(m-1)_n},
$c_{\beta^{\prime},0,\beta^{\prime\prime}}^{\left(  m\right)  }>\frac{1}%
{2}c_{\beta^{\prime},\beta^{\prime\prime}}^{\left(  m-1\right)  }$, so that
now%
\begin{align}
\sum\limits_{\beta\leq\mathbf{1},\left\vert \beta\right\vert =n}c_{\beta
}^{\left(  m\right)  }x^{\beta}  & >\frac{1}{2}\sum\limits_{\substack{\left(
\beta^{\prime},\beta^{\prime\prime}\right)  \leq\mathbf{1} \\\left\vert
\beta^{\prime}\right\vert +\left\vert \beta^{\prime\prime}\right\vert
=n-1}}c_{\beta^{\prime},\beta^{\prime\prime}}^{\left(  m-1\right)  }\left(
x^{\prime}\right)  ^{\beta^{\prime}}\left(  x^{\prime\prime}\right)
^{\beta^{\prime\prime}}\nonumber\\
& =\frac{1}{2}\sum\limits_{\gamma\leq\mathbf{1},\left\vert \gamma\right\vert
=n-1}c_{\gamma}^{\left(  m-1\right)  }\left(  x^{\prime},x^{\prime\prime
}\right)  ^{\gamma}.\label{av766}%
\end{align}

The equations of this theorem now follow easily by successive applications of
\ref{av766}.
\end{proof}

\subsection{Bounds for $\left(  T^{-1}\beta\right)  _{N-1}$}

From \ref{av725}, $\left(  T^{-1}\beta\right)  _{k}=\left(  -1\right)
^{k-1}\frac{1}{\left\vert T\right\vert }\left(  \left(  -1\right)
^{N}M_{k,N-1}-M_{k,1}\right)  $ for $k\leq N-1$. From \ref{av533},
$M_{N-1,N-1}=\left\vert T_{1:N-2}\right\vert $ and from \ref{av552},
$M_{N-1,1}=-\frac{1}{2^{N-2}}$. Then%
\begin{align}
\left(  T^{-1}\beta\right)  _{N-1}  & =\left(  -1\right)  ^{N-2}\frac
{1}{\left\vert T\right\vert }\left(  \left(  -1\right)  ^{N}M_{N-1,N-1}%
-M_{N-1,1}\right) \nonumber\\
& =\frac{1}{\left\vert T\right\vert }\left(  M_{N-1,N-1}-\left(  -1\right)
^{N}M_{N-1,1}\right) \nonumber\\
& =\frac{1}{\left\vert T\right\vert }\left(  \left\vert T_{1:N-2}\right\vert
-\left(  -1\right)  ^{N}-\frac{1}{2^{N-2}}\right) \nonumber\\
& =\frac{1}{\left\vert T\right\vert }\left(  \left\vert T_{1:N-2}\right\vert
-\frac{1}{2^{N-2}}\right)  .\label{av031}%
\end{align}

Set $x=\frac{1}{v}\Delta_{1}X$. Then $x\in\mathbb{R}_{+}^{N-1}$ from
\ref{av599} and then \ref{av598},%
\begin{align*}
\left(  T^{-1}\beta\right)  _{N-1}  & =\frac{\sum\limits_{k=0}^{N-2}%
\sum\limits_{\gamma\leq\mathbf{1},\left\vert \gamma\right\vert =k}c_{\gamma
}^{\left(  N-2\right)  }x_{1:N-2}^{\gamma}-\frac{1}{2^{N-2}}}{\sum
\limits_{k=0}^{N-1}\sum\limits_{\beta\leq\mathbf{1},\left\vert \beta
\right\vert =k}c_{\beta}^{\left(  N-1\right)  }x^{\beta}}\\
& =\frac{c_{0}^{\left(  N-2\right)  }+\sum\limits_{k=1}^{N-2}\sum
\limits_{\gamma\leq\mathbf{1},\left\vert \gamma\right\vert =k}c_{\gamma
}^{\left(  N-2\right)  }x_{1:N-2}^{\gamma}-\frac{1}{2^{N-2}}}{\sum
\limits_{k=0}^{N-1}\sum\limits_{\beta\leq\mathbf{1},\left\vert \beta
\right\vert =k}c_{\beta}^{\left(  N-1\right)  }x^{\beta}}\\
& =\frac{\frac{N-1}{2^{N-2}}+\sum\limits_{k=1}^{N-2}\sum\limits_{\gamma
\leq\mathbf{1},\left\vert \gamma\right\vert =k}c_{\gamma}^{\left(  N-2\right)
}x_{1:N-2}^{\gamma}-\frac{1}{2^{N-2}}}{\sum\limits_{k=0}^{N-1}\sum
\limits_{\beta\leq\mathbf{1},\left\vert \beta\right\vert =k}c_{\beta}^{\left(
N-1\right)  }x^{\beta}}\\
& =\frac{\frac{N-2}{2^{N-2}}+\sum\limits_{k=1}^{N-2}\sum\limits_{\gamma
\leq\mathbf{1},\left\vert \gamma\right\vert =k}c_{\gamma}^{\left(  N-2\right)
}x_{1:N-2}^{\gamma}}{\sum\limits_{k=0}^{N-1}\sum\limits_{\beta\leq
\mathbf{1},\left\vert \beta\right\vert =k}c_{\beta}^{\left(  N-1\right)
}x^{\beta}}\\
& >0.
\end{align*}

We will again use Theorem \ref{vThm_inequal_Sigma_C^m_1minusBeta} but this
time with $n=k$ and $l=m=N-1$ so that the estimate \ref{av019} becomes:%
\[
\sum\limits_{\beta\leq\mathbf{1},\left\vert \beta\right\vert =k}c_{\beta
}^{\left(  N-1\right)  }x^{\beta}>\frac{1}{2}\sum\limits_{\gamma\leq
\mathbf{1},\left\vert \gamma\right\vert =k-1}c_{\gamma}^{\left(  N-2\right)
}x_{1:N-2}^{\gamma},\quad1\leq k\leq N-1,\text{ }N\geq2,
\]

and hence by also using \ref{av598} and \ref{av595} i.e. $c_{\mathbf{1}%
}^{\left(  N-1\right)  }=1$ and $c_{\mathbf{0}}^{\left(  N-1\right)  }%
=\frac{N}{2^{N-1}}$:%
\begin{align*}
\left(  T^{-1}\beta\right)  _{N-1}  & =\frac{\frac{N-2}{2^{N-2}}%
+\sum\limits_{k=1}^{N-2}\sum\limits_{\gamma\leq\mathbf{1},\left\vert
\gamma\right\vert =k}c_{\gamma}^{\left(  N-2\right)  }x_{1:N-2}^{\gamma}}%
{\sum\limits_{k=0}^{N-1}\sum\limits_{\beta\leq\mathbf{1},\left\vert
\beta\right\vert =k}c_{\beta}^{\left(  N-1\right)  }x^{\beta}}\\
& =\frac{\frac{N-2}{2^{N-2}}+\sum\limits_{k=1}^{N-2}\sum\limits_{\gamma
\leq\mathbf{1},\left\vert \gamma\right\vert =k}c_{\gamma}^{\left(  N-2\right)
}x_{1:N-2}^{\gamma}}{c_{\mathbf{0}}^{\left(  N-1\right)  }+\sum\limits_{k=1}%
^{N-2}\sum\limits_{\beta\leq\mathbf{1},\left\vert \beta\right\vert =k}%
c_{\beta}^{\left(  N-1\right)  }x^{\beta}+c_{\mathbf{1}}^{\left(  N-1\right)
}x^{\mathbf{1}}}\\
& =\frac{\frac{N-2}{2^{N-2}}+\sum\limits_{k=1}^{N-2}\sum\limits_{\gamma
\leq\mathbf{1},\left\vert \gamma\right\vert =k}c_{\gamma}^{\left(  N-2\right)
}x_{1:N-2}^{\gamma}}{\frac{N}{2^{N-1}}+\sum\limits_{k=1}^{N-2}\sum
\limits_{\beta\leq\mathbf{1},\left\vert \beta\right\vert =k}c_{\beta}^{\left(
N-1\right)  }x^{\beta}+x^{\mathbf{1}}}\\
& <\frac{\frac{N-2}{2^{N-2}}+\sum\limits_{k=1}^{N-2}\sum\limits_{\gamma
\leq\mathbf{1},\left\vert \gamma\right\vert =k}c_{\gamma}^{\left(  N-2\right)
}x_{1:N-2}^{\gamma}}{\frac{N}{2^{N-1}}+\frac{1}{2}\sum\limits_{k=1}^{N-2}%
\sum\limits_{\gamma\leq\mathbf{1},\left\vert \gamma\right\vert =k-1}c_{\gamma
}^{\left(  N-2\right)  }x_{1:N-2}^{\gamma}+x^{\mathbf{1}}}\\
& <\frac{\frac{N-2}{2^{N-2}}+\sum\limits_{k=1}^{N-2}\sum\limits_{\gamma
\leq\mathbf{1},\left\vert \gamma\right\vert =k}c_{\gamma}^{\left(  N-2\right)
}x_{1:N-2}^{\gamma}}{\frac{N}{2^{N-1}}+\frac{1}{2}\sum\limits_{k=1}^{N-2}%
\sum\limits_{\gamma\leq\mathbf{1},\left\vert \gamma\right\vert =k-1}c_{\gamma
}^{\left(  N-2\right)  }x_{1:N-2}^{\gamma}}\\
& =\frac{\frac{N-2}{2^{N-2}}+\sum\limits_{k=1}^{N-2}\sum\limits_{\gamma
\leq\mathbf{1},\left\vert \gamma\right\vert =k}c_{\gamma}^{\left(  N-2\right)
}x_{1:N-2}^{\gamma}}{\frac{1}{2}\left(  \frac{N}{2^{N-2}}+\sum\limits_{k=1}%
^{N-2}\sum\limits_{\gamma\leq\mathbf{1},\left\vert \gamma\right\vert
=k-1}c_{\gamma}^{\left(  N-2\right)  }x_{1:N-2}^{\gamma}\right)  }\\
& <2,
\end{align*}

and in summary%
\begin{equation}
0<\left(  T^{-1}\beta\right)  _{N-1}<2,\quad N\geq3,\text{ }\rho
>0,\label{av764}%
\end{equation}

where $T=T\left(  x;1\right)  $, $x=\frac{1}{v}\Delta_{1}X$ and from
\ref{av839}, $v=\rho N$.

\subsection{Bounds for $\left(  T^{-1}\beta\right)  _{k}$, $2\leq k\leq N-2$}

The result of this subsection is the uniform bound \ref{av022} calculated
using minor formulas from Subsection \ref{SbSect_minor_summary}.

Set $x=x_{1:N-1}=\frac{1}{v}\Delta_{1}X\in\mathbb{R}_{+}^{N-1}$. Now when
$2\leq k\leq N-2$, from \ref{av015} and \ref{av017},%
\begin{align}
\left(  T^{-1}\beta\right)  _{k} &  =\frac{1}{\left\vert T\right\vert }\left(
-1\right)  ^{k-1}\left(  \left(  -1\right)  ^{N}M_{k,N-1}-M_{k,1}\right)
\nonumber\\
&  =\frac{1}{\left\vert T\right\vert }\left(  -1\right)  ^{k-1}\left(  \left(
-1\right)  ^{N}\left(  -\frac{1}{2}\right)  ^{N-1-k}\left\vert T_{1:k-1}%
\right\vert -\left(  -\frac{1}{2}\right)  ^{k-1}\left\vert T_{k+1:N-1}%
\right\vert \right) \nonumber\\
&  =\frac{1}{\left\vert T\right\vert }\left(  \left(  -1\right)  ^{N}\left(
-1\right)  ^{N}\frac{1}{2^{N-1-k}}\left\vert T_{1:k-1}\right\vert -\frac
{1}{2^{k-1}}\left\vert T_{k+1:N-1}\right\vert \right) \nonumber\\
&  =\frac{1}{2^{N-1-k}}\frac{1}{\left\vert T\right\vert }\left\vert
T_{1:k-1}\right\vert -\frac{1}{2^{k-1}}\frac{1}{\left\vert T\right\vert
}\left\vert T_{k+1:N-1}\right\vert ,\label{av748}%
\end{align}

so that%
\[
-\frac{1}{2^{k-1}}\frac{1}{\left\vert T\right\vert }\left\vert T_{k+1:N-1}%
\right\vert <\left(  T^{-1}\beta\right)  _{k}<\frac{1}{2^{N-1-k}}\frac
{1}{\left\vert T\right\vert }\left\vert T_{1:k-1}\right\vert ,
\]

which implies that%
\begin{equation}
\left\vert \left(  T^{-1}\beta\right)  _{k}\right\vert <\max\left\{  \frac
{1}{2^{N-1-k}}\frac{1}{\left\vert T\right\vert }\left\vert T_{1:k-1}%
\right\vert ,\frac{1}{2^{k-1}}\frac{1}{\left\vert T\right\vert }\left\vert
T_{k+1:N-1}\right\vert \right\}  .\label{av018}%
\end{equation}

Applying Theorem \ref{vThm_formula_for_C^(m)_idx} to \ref{av748} gives
\[
\frac{1}{2^{N-1-k}}\frac{1}{\left\vert T\right\vert }\left\vert T_{1:k-1}%
\right\vert =\frac{1}{2^{N-1-k}}\frac{1}{\left\vert T\right\vert }%
\sum\limits_{n=0}^{k-1}\sum\limits_{\gamma\leq\mathbf{1},\left\vert
\gamma\right\vert =n}c_{\gamma}^{\left(  k-1\right)  }x_{1:k-1}^{\gamma},
\]

and by the estimate \ref{av019} of Theorem
\ref{vThm_inequal_Sigma_C^m_1minusBeta},%
\begin{align*}
\left\vert T\right\vert  & =\sum\limits_{m=0}^{N-1}\sum\limits_{\beta
\leq\mathbf{1},\left\vert \beta\right\vert =m}c_{\beta}^{\left(  N-1\right)
}x^{\beta}\\
& >\sum\limits_{m=N-k}^{N-1}\sum\limits_{\beta\leq\mathbf{1},\left\vert
\beta\right\vert =m}c_{\beta}^{\left(  N-1\right)  }x^{\beta}\\
& >\frac{1}{2^{N-k}}\sum\limits_{m=N-k}^{N-1}\sum\limits_{\substack{\beta
\leq\mathbf{1} \\\left\vert \beta\right\vert =m-\left(  N-k\right)  }%
}c_{\beta}^{\left(  N-1-\left(  N-k\right)  \right)  }x_{1:N-1-\left(
N-k\right)  }^{\beta}\\
& =\frac{1}{2^{N-k}}\sum\limits_{m=N-k}^{N-1}\sum\limits_{\substack{\beta
\leq\mathbf{1} \\\left\vert \beta\right\vert =m-\left(  N-k\right)  }%
}c_{\beta}^{\left(  k-1\right)  }x_{1:k-1}^{\beta}\\
& =\frac{1}{2^{N-k}}\sum\limits_{n=0}^{k-1}\sum\limits_{\beta\leq
\mathbf{1},\left\vert \beta\right\vert =n}c_{\beta}^{\left(  k-1\right)
}x_{1:k-1}^{\beta},
\end{align*}

so that%
\begin{equation}
\frac{1}{2^{N-k}}\frac{1}{\left\vert T\right\vert }\left\vert T_{1:k-1}%
\right\vert <2.\label{av016}%
\end{equation}

Further%
\[
\frac{1}{2^{k-1}}\frac{1}{\left\vert T\right\vert }\left\vert T_{k+1:N-1}%
\right\vert =\frac{1}{2^{k-1}}\frac{1}{\left\vert T\right\vert }%
\sum\limits_{n=0}^{N-k-1}\sum\limits_{\gamma\leq\mathbf{1},\left\vert
\gamma\right\vert =n}c_{\gamma}^{\left(  N-k-1\right)  }x_{k+1:N-1}^{\gamma},
\]

and by the estimate \ref{av020} of Theorem
\ref{vThm_inequal_Sigma_C^m_1minusBeta},%
\begin{align*}
\left\vert T\right\vert  & =\sum\limits_{m=0}^{N-1}\sum\limits_{\beta
\leq\mathbf{1},\left\vert \beta\right\vert =m}c_{\beta}^{\left(  N-1\right)
}x^{\beta}\\
& >\sum\limits_{m=k}^{N-1}\sum\limits_{\beta\leq\mathbf{1},\left\vert
\beta\right\vert =m}c_{\beta}^{\left(  N-1\right)  }x^{\beta}\\
& >\frac{1}{2^{k}}\sum\limits_{m=k}^{N-1}\sum\limits_{\substack{\beta
\leq\mathbf{1} \\\left\vert \beta\right\vert =m-k}}c_{\beta}^{\left(
N-1-k\right)  }x_{k+1:N-1}^{\beta}\\
& =\frac{1}{2^{k}}\sum\limits_{m=0}^{N-k-1}\sum\limits_{\substack{\beta
\leq\mathbf{1} \\\left\vert \beta\right\vert =m}}c_{\beta}^{\left(
N-1-k\right)  }x_{k+1:N-1}^{\beta},
\end{align*}

so that%
\begin{equation}
\frac{1}{2^{k-1}}\frac{1}{\left\vert T\right\vert }\left\vert T_{k+1:N-1}%
\right\vert <2.\label{av021}%
\end{equation}

The inequality \ref{av018} now becomes the uniform bound%
\begin{equation}
\left\vert \left(  T^{-1}\beta\right)  _{k}\right\vert <2,\quad2\leq k\leq
N-2,\text{ }N\geq4,\text{ }\rho>0,\label{av022}%
\end{equation}

where $T=T\left(  x;1\right)  $, $x=\frac{1}{v}\Delta_{1}X$ and from
\ref{av839}, $v=\rho N$.

\subsection{Upper bounds for $\left\vert T^{-1}\beta\right\vert _{\max}$ and
$\left\vert \mu\left(  f\right)  \right\vert _{\max}$}

The estimates \ref{av022}, \ref{av764} and \ref{av754} can be combined
directly to give%
\begin{equation}
\left\vert T^{-1}\beta\right\vert _{\max}<2,\quad N\geq4,\text{ }%
\rho>0,\label{av035}%
\end{equation}

where $T=T\left(  \frac{1}{v}\Delta_{1}X;1\right)  $ and from \ref{av839},
$v=\rho N$.

Applying the bound \ref{av035} to the formula \ref{av007} for $\mu^{\prime
}\left(  f\right)  $ yields%
\begin{align*}
\left\vert \mu^{\prime}\left(  f\right)  \right\vert _{\max}  & =\left\vert
\frac{1}{v}T^{-1}\Delta_{X}f-\frac{1}{2}\mu_{1}\left(  f\right)  T^{-1}%
\beta\right\vert _{\max}\\
& \leq\left\vert \frac{1}{v}T^{-1}\Delta_{X}f\right\vert _{\max}+\left\vert
\mu_{1}\left(  f\right)  \right\vert \left\vert \frac{1}{2}T^{-1}%
\beta\right\vert _{\max}\\
& \leq\left\vert \frac{1}{v}T^{-1}\Delta_{X}f\right\vert _{\max}+\left\vert
\mu_{1}\left(  f\right)  \right\vert ,
\end{align*}

when $N\geq4$ and $\rho>0$. This implies that%
\begin{align}
\left\vert \mu\left(  f\right)  \right\vert _{\max}  & =\max\left\{
\left\vert \mu_{1}\left(  f\right)  \right\vert ,\left\vert \mu^{\prime
}\left(  f\right)  \right\vert _{\max}\right\} \nonumber\\
& \leq\max\left\{  \left\vert \mu_{1}\left(  f\right)  \right\vert ,\left\vert
\frac{1}{v}T^{-1}\Delta_{X}f\right\vert _{\max}+\left\vert \mu_{1}\left(
f\right)  \right\vert \right\} \nonumber\\
& =\left\vert \frac{1}{v}T^{-1}\Delta_{X}f\right\vert _{\max}+\left\vert
\mu_{1}\left(  f\right)  \right\vert .\label{av087}%
\end{align}

\section{Uniform upper bounds of $\left\vert \frac{1}{v}T^{-1}\Delta
_{X}f\right\vert _{\max}$ using $\left\Vert Df\right\Vert _{\infty}%
$\label{Sect_bound_invTdiff_data_fn}}

On the assumption that $f\in C_{B}^{\left(  0\right)  }\left(  \Omega\right)
$ and $Df\in L^{\infty}\left(  \Omega\right)  $, the goal of this section is
Corollary \ref{vCor_2_Thm_low_bnd_detT} which shows that $\frac{1}%
{v}\left\vert T^{-1}\Delta_{X}f\right\vert _{\max}\leq2\left\Vert
Df\right\Vert _{\infty;\Omega}$ and so from \ref{av087}, $\left\vert
\mu\left(  f\right)  \right\vert _{\max}\leq2\left\Vert Df\right\Vert
_{\infty;\Omega}+\mu_{1}\left(  f\right)  $. The steps are:

\textbf{Subsection} \ref{SbSect_bound_invTdiff_subdetT}: We first derive the
upper bound \ref{av029} for the term $\left\vert \frac{1}{v}T^{-1}\Delta
_{X}f\right\vert _{\max}$ which occurs in the estimate \ref{av087} for
$\left\vert \mu\left(  f\right)  \right\vert _{\max}$.

\textbf{Subsection} \ref{SbSect_detT_in_terms_of_detTzeros}: We will proceed
below by deriving the expansions \ref{av039} and \ref{av037} of $\left\vert
T\left(  x_{1:m}\right)  \right\vert $ which involve sums of terms of the form
$\left\vert T\left(  x^{\prime},0_{n}\right)  \right\vert $ and $\left\vert
T\left(  0_{n},x^{\prime\prime}\right)  \right\vert $ respectively.

\textbf{Subsection} \ref{SbSect_bnds_detT_with_zeros}: The next step is to
obtain bounds for determinants with zeros in their arguments e.g. $\left\vert
T\left(  0_{n},x\right)  \right\vert \geq c_{n}\left\vert T\left(  x\right)
\right\vert $ and $\left\vert T\left(  x,0_{n}\right)  \right\vert \geq
c_{n}^{\prime}\left\vert T\left(  x\right)  \right\vert $.

\textbf{Subsection} \ref{SbSecT_max_invT_deltaXFn}:These bounds are then
applied to eliminate the zeros from the expansions \ref{av039} and \ref{av037}
of $\left\vert T\left(  x_{1:m}\right)  \right\vert $ to get Theorem
\ref{vThm_low_bnd_detT}. It is then easy to obtain from \ref{av029} the
estimate \ref{av081}:
\[
\left\vert \mu\left(  f\right)  \right\vert _{\max}\leq2\left\Vert
Df\right\Vert _{\infty;\Omega}+\left\vert \mu_{1}\left(  f\right)  \right\vert
,
\]

\subsection{Bounds in terms of $\left\Vert Df\right\Vert _{\infty}$ and
sub-determinants $\left\vert T_{i:j}\right\vert $%
\label{SbSect_bound_invTdiff_subdetT}}

In this subsection we will derive the upper bounds \ref{av029} for the
components of $\frac{1}{v}T^{-1}\Delta_{X}f$.

We require the minor formulas \ref{av023} and assume that $f\in C_{B}^{\left(
0\right)  }\left(  \Omega\right)  $ and $Df\in L^{\infty}\left(
\Omega\right)  $.\medskip

\fbox{Case $k=1$} From \ref{av023}, $M_{1,n}=\left(  -\frac{1}{2}\right)
^{n-1}\left\vert T_{1:0}\right\vert $ $\left\vert T_{n+1:N-1}\right\vert
=\left(  -\frac{1}{2}\right)  ^{n-1}\left\vert T_{n+1:N-1}\right\vert $, so
that \ref{av758} becomes,%
\begin{align*}
\left\vert T\right\vert \left(  T^{-1}\Delta_{X}f\right)  _{1}  &
=-\sum\limits_{n=1}^{N-1}\left(  -1\right)  ^{n}M_{1,n}\Delta_{1}f\left(
x^{\left(  n\right)  }\right) \\
& =-\sum\limits_{n=1}^{N-1}\left(  -1\right)  ^{n}\left(  \left(  -\frac{1}%
{2}\right)  ^{n-1}\left\vert T_{n+1:N-1}\right\vert \right)  \Delta
_{1}f\left(  x^{\left(  n\right)  }\right) \\
& =\sum\limits_{n=1}^{N-1}\frac{1}{2^{n-1}}\left\vert T_{n+1:N-1}\right\vert
\Delta_{1}f\left(  x^{\left(  n\right)  }\right) \\
& =\sum\limits_{n=1}^{N-2}\frac{1}{2^{n-1}}\left\vert T_{n+1:N-1}\right\vert
\Delta_{1}f\left(  x^{\left(  n\right)  }\right)  +\frac{1}{2^{N-2}}\left\vert
T_{N:N-1}\right\vert \Delta_{1}f\left(  x^{\left(  N-1\right)  }\right) \\
& =\sum\limits_{n=1}^{N-2}\frac{1}{2^{n-1}}\left\vert T_{n+1:N-1}\right\vert
\Delta_{1}f\left(  x^{\left(  n\right)  }\right)  +\frac{1}{2^{N-2}}\Delta
_{1}f\left(  x^{\left(  N-1\right)  }\right) \\
& =\sum\limits_{n=1}^{N-2}\frac{\Delta_{1}x^{\left(  n\right)  }}{2^{n-1}%
}\left\vert T_{n+1:N-1}\right\vert \frac{\Delta_{1}f\left(  x^{\left(
n\right)  }\right)  }{\Delta_{1}x^{\left(  n\right)  }}+\frac{\Delta
_{1}x^{\left(  N-1\right)  }}{2^{N-2}}\frac{\Delta_{1}f\left(  x^{\left(
N-1\right)  }\right)  }{\Delta_{1}x^{\left(  N-1\right)  }},
\end{align*}

so that%
\[
\left\vert T\right\vert \left\vert \left(  T^{-1}\Delta_{X}f\right)
_{1}\right\vert \leq\max_{n=1}^{N-1}\left\vert \frac{\Delta_{1}f\left(
x^{\left(  n\right)  }\right)  }{\Delta_{1}x^{\left(  n\right)  }}\right\vert
\left(  \sum\limits_{n=1}^{N-2}\frac{\Delta_{1}x^{\left(  n\right)  }}%
{2^{n-1}}\left\vert T_{n+1:N-1}\right\vert +\frac{\Delta_{1}x^{\left(
N-1\right)  }}{2^{N-2}}\right)  .
\]

From Lemma \ref{vLem_Taylor_extension},%
\[
\left\vert \left(  T^{-1}\Delta_{X}f\right)  _{1}\right\vert \leq
\frac{\left\Vert Df\right\Vert _{\infty;\Omega}}{\left\vert T\right\vert
}\left(  \sum\limits_{n=1}^{N-2}\frac{\Delta_{1}x^{\left(  n\right)  }%
}{2^{n-1}}\left\vert T_{n+1:N-1}\right\vert +\frac{\Delta_{1}x^{\left(
N-1\right)  }}{2^{N-2}}\right)  .
\]

Now set%
\begin{equation}
x:=\frac{1}{v}\Delta_{1}X=\frac{1}{v}\left(  \Delta_{1}x^{\left(  n\right)
}\right)  _{n=1}^{N-1},\label{av025}%
\end{equation}

so that%
\[
\left\vert \left(  T^{-1}\Delta_{X}f\right)  _{1}\right\vert \leq
\frac{\left\Vert Df\right\Vert _{\infty;\Omega}}{\left\vert T\right\vert
}\left(  \sum\limits_{n=1}^{N-2}\frac{vx_{n}}{2^{n-1}}\left\vert
T_{n+1:N-1}\right\vert +\frac{vx_{N-1}}{2^{N-2}}\right)  ,
\]

and hence%
\begin{equation}
\frac{1}{v}\left\vert \left(  T^{-1}\Delta_{X}f\right)  _{1}\right\vert
\leq\frac{\left\Vert Df\right\Vert _{\infty;\Omega}}{\left\vert T\right\vert
}\left(  \sum\limits_{n=1}^{N-2}\frac{x_{n}}{2^{n-1}}\left\vert T_{n+1:N-1}%
\right\vert +\frac{x_{N-1}}{2^{N-2}}\right)  .\label{av027}%
\end{equation}
\medskip

\fbox{Case $2\leq k\leq N-2$} From \ref{av758} and then \ref{av023},%
\begin{align*}
&  \left\vert T\right\vert \left(  T^{-1}\Delta_{X}f\right)  _{k}\\
&  =\left(  -1\right)  ^{k}\left(  \sum\limits_{n=1}^{k}\left(  -1\right)
^{n}M_{k,n}\Delta_{1}f\left(  x^{\left(  n\right)  }\right)  +\sum
\limits_{m=k+1}^{N-1}\left(  -1\right)  ^{m}M_{k,m}\Delta_{1}f\left(
x^{\left(  m\right)  }\right)  \right) \\
&  =\left(  -1\right)  ^{k}\left(  \sum\limits_{n=1}^{k}\left(  -1\right)
^{n}M_{n,k}\Delta_{1}f\left(  x^{\left(  n\right)  }\right)  +\sum
\limits_{m=k+1}^{N-1}\left(  -1\right)  ^{m}M_{k,m}\Delta_{1}f\left(
x^{\left(  m\right)  }\right)  \right) \\
&  =\left(  -1\right)  ^{k}\sum\limits_{n=1}^{k}\left(  -1\right)  ^{n}\left(
-\frac{1}{2}\right)  ^{k-n}\left\vert T_{1:n-1}\right\vert \text{ }\left\vert
T_{k+1:N-1}\right\vert \Delta_{1}f\left(  x^{\left(  n\right)  }\right)  +\\
&  \qquad\qquad+\left(  -1\right)  ^{k}\sum\limits_{m=k+1}^{N-1}\left(
-1\right)  ^{m}\left(  -\frac{1}{2}\right)  ^{m-k}\left\vert T_{1:k-1}%
\right\vert \text{ }\left\vert T_{m+1:N-1}\right\vert \Delta_{1}f\left(
x^{\left(  m\right)  }\right) \\
&  =\sum\limits_{n=1}^{k}\frac{1}{2^{k-n}}\left\vert T_{1:n-1}\right\vert
\left\vert T_{k+1:N-1}\right\vert \Delta_{1}f\left(  x^{\left(  n\right)
}\right)  +\sum\limits_{m=k+1}^{N-1}\frac{1}{2^{m-k}}\left\vert T_{1:k-1}%
\right\vert \left\vert T_{m+1:N-1}\right\vert \Delta_{1}f\left(  x^{\left(
m\right)  }\right) \\
&  =\left\vert T_{k+1:N-1}\right\vert \sum\limits_{n=1}^{k}\frac{1}{2^{k-n}%
}\left\vert T_{1:n-1}\right\vert \Delta_{1}f\left(  x^{\left(  n\right)
}\right)  +\left\vert T_{1:k-1}\right\vert \sum\limits_{m=k+1}^{N-1}\frac
{1}{2^{m-k}}\left\vert T_{m+1:N-1}\right\vert \Delta_{1}f\left(  x^{\left(
m\right)  }\right) \\
&  =\left\vert T_{k+1:N-1}\right\vert \sum\limits_{n=1}^{k}\frac{vx_{n}%
}{2^{k-n}}\left\vert T_{1:n-1}\right\vert \frac{\Delta_{1}f\left(  x^{\left(
n\right)  }\right)  }{\Delta_{1}x^{\left(  n\right)  }}+\left\vert
T_{1:k-1}\right\vert \sum\limits_{m=k+1}^{N-1}\frac{vx_{m}}{2^{m-k}}\left\vert
T_{m+1:N-1}\right\vert \frac{\Delta_{1}f\left(  x^{\left(  m\right)  }\right)
}{\Delta_{1}x^{\left(  m\right)  }}\\
&  =v\left(  \left\vert T_{k+1:N-1}\right\vert \sum\limits_{n=1}^{k}%
\frac{x_{n}}{2^{k-n}}\left\vert T_{1:n-1}\right\vert \frac{\Delta_{1}f\left(
x^{\left(  n\right)  }\right)  }{\Delta_{1}x^{\left(  n\right)  }}+\left\vert
T_{1:k-1}\right\vert \sum\limits_{m=k+1}^{N-1}\frac{x_{m}}{2^{m-k}}\left\vert
T_{m+1:N-1}\right\vert \frac{\Delta_{1}f\left(  x^{\left(  m\right)  }\right)
}{\Delta_{1}x^{\left(  m\right)  }}\right)  ,
\end{align*}

and applying Lemma \ref{vLem_Taylor_extension} yields%
\begin{align*}
&  \frac{1}{v}\left\vert T\right\vert \left\vert \left(  T^{-1}\Delta
_{X}f\right)  _{k}\right\vert \\
&  \leq\left\vert T_{k+1:N-1}\right\vert \sum\limits_{n=1}^{k}\frac{x_{n}%
}{2^{k-n}}\left\vert T_{1:n-1}\right\vert \left\vert \frac{\Delta_{1}f\left(
x^{\left(  n\right)  }\right)  }{\Delta_{1}x^{\left(  n\right)  }}\right\vert
+\left\vert T_{1:k-1}\right\vert \sum\limits_{m=k+1}^{N-1}\frac{x_{m}}%
{2^{m-k}}\left\vert T_{m+1:N-1}\right\vert \left\vert \frac{\Delta_{1}f\left(
x^{\left(  m\right)  }\right)  }{\Delta_{1}x^{\left(  m\right)  }}\right\vert
\\
&  \leq\max_{n=1}^{N-1}\left\vert \frac{\Delta_{1}f\left(  x^{\left(
n\right)  }\right)  }{\Delta_{1}x^{\left(  n\right)  }}\right\vert \left(
\left\vert T_{k+1:N-1}\right\vert \sum\limits_{n=1}^{k}\frac{x_{n}}{2^{k-n}%
}\left\vert T_{1:n-1}\right\vert +\left\vert T_{1:k-1}\right\vert
\sum\limits_{m=k+1}^{N-1}\frac{x_{m}}{2^{m-k}}\left\vert T_{m+1:N-1}%
\right\vert \right) \\
&  \leq\left\Vert Df\right\Vert _{\infty;\Omega}\left(  \left\vert
T_{k+1:N-1}\right\vert \sum\limits_{n=1}^{k}\frac{x_{n}}{2^{k-n}}\left\vert
T_{1:n-1}\right\vert +\left\vert T_{1:k-1}\right\vert \sum\limits_{m=k+1}%
^{N-1}\frac{x_{m}}{2^{m-k}}\left\vert T_{m+1:N-1}\right\vert \right) \\
&  \leq\left\Vert Df\right\Vert _{\infty;\Omega}\left(
\begin{array}
[c]{l}%
\left\vert T_{k+1:N-1}\right\vert \left(  \frac{x_{1}}{2^{k-1}}+\sum
\limits_{n=2}^{k}\frac{x_{n}}{2^{k-n}}\left\vert T_{1:n-1}\right\vert \right)
+\\
\qquad\qquad\qquad+\left\vert T_{1:k-1}\right\vert \left(  \sum\limits_{m=k+1}%
^{N-2}\frac{x_{m}}{2^{m-k}}\left\vert T_{m+1:N-1}\right\vert +\frac{x_{N-1}%
}{2^{N-1-k}}\right)
\end{array}
\right)  ,
\end{align*}

so that when $2\leq k\leq N-2$,%
\begin{align}
&  \frac{1}{v}\left\vert \left(  T^{-1}\Delta_{X}f\right)  _{k}\right\vert
\nonumber\\
&  \leq\frac{\left\Vert Df\right\Vert _{\infty;\Omega}}{\left\vert
T\right\vert }\left(
\begin{array}
[c]{l}%
\left\vert T_{k+1:N-1}\right\vert \left(  \frac{x_{1}}{2^{k-1}}+\sum
\limits_{n=2}^{k}\frac{x_{n}}{2^{k-n}}\left\vert T_{1:n-1}\right\vert \right)
+\\
\qquad\qquad+\left\vert T\left(  x_{1:k-1}\right)  \right\vert \left(
\sum\limits_{m=k+1}^{N-2}\frac{x_{m}}{2^{m-k}}\left\vert T_{m+1:N-1}%
\right\vert +\frac{x_{N-1}}{2^{N-1-k}}\right)
\end{array}
\right) \label{av026}%
\end{align}
\medskip

\fbox{When $k=N-1$} From \ref{av023},
\begin{align*}
M_{n,N-1}=\left(  -\frac{1}{2}\right)  ^{N-1-n}\left\vert T_{1:n-1}\right\vert
\left\vert T_{N-1+1:N-1}\right\vert  & =\left(  -\frac{1}{2}\right)
^{N-1-n}\left\vert T_{1:n-1}\right\vert \left\vert T_{N:N-1}\right\vert \\
& =\left(  -\frac{1}{2}\right)  ^{N-1-n}\left\vert T_{1:n-1}\right\vert ,
\end{align*}

and so \ref{av758} becomes,
\begin{align*}
\left\vert T\right\vert \left(  T^{-1}\Delta_{X}f\right)  _{N-1}  & =\left(
-1\right)  ^{N-1}\sum\limits_{n=1}^{N-1}\left(  -1\right)  ^{n}M_{N-1,n}%
\Delta_{1}f\left(  x^{\left(  n\right)  }\right) \\
& =\left(  -1\right)  ^{N-1}\sum\limits_{n=1}^{N-1}\left(  -1\right)
^{n}M_{n,N-1}\Delta_{1}f\left(  x^{\left(  n\right)  }\right) \\
& =\left(  -1\right)  ^{N-1}\sum\limits_{n=1}^{N-1}\left(  -1\right)
^{n}\left(  -\frac{1}{2}\right)  ^{N-1-n}\left\vert T_{1:n-1}\right\vert
\Delta_{1}f\left(  x^{\left(  n\right)  }\right) \\
& =\sum\limits_{n=1}^{N-1}\frac{1}{2^{N-1-n}}\left\vert T_{1:n-1}\right\vert
\Delta_{1}f\left(  x^{\left(  n\right)  }\right) \\
& =\frac{1}{2^{N-1-1}}\left\vert T_{1:1-1}\right\vert \Delta_{1}f\left(
x^{\left(  1\right)  }\right)  +\sum\limits_{n=2}^{N-1}\frac{1}{2^{N-1-n}%
}\left\vert T_{1:n-1}\right\vert \Delta_{1}f\left(  x^{\left(  n\right)
}\right) \\
& =\frac{1}{2^{N-2}}\Delta_{1}f\left(  x^{\left(  1\right)  }\right)
+\sum\limits_{n=2}^{N-1}\frac{1}{2^{N-1-n}}\left\vert T_{1:n-1}\right\vert
\Delta_{1}f\left(  x^{\left(  n\right)  }\right) \\
& =\frac{vx_{1}}{2^{N-2}}\frac{\Delta_{1}f\left(  x^{\left(  1\right)
}\right)  }{\Delta_{1}x^{\left(  1\right)  }}+\sum\limits_{n=2}^{N-1}%
\frac{vx_{n}}{2^{N-1-n}}\left\vert T_{1:n-1}\right\vert \frac{\Delta
_{1}f\left(  x^{\left(  n\right)  }\right)  }{\Delta_{1}x^{\left(  n\right)
}}\\
& =v\left(  \frac{x_{1}}{2^{N-2}}\frac{\Delta_{1}f\left(  x^{\left(  1\right)
}\right)  }{\Delta_{1}x^{\left(  1\right)  }}+\sum\limits_{n=2}^{N-1}%
\frac{x_{n}}{2^{N-1-n}}\left\vert T_{1:n-1}\right\vert \frac{\Delta
_{1}f\left(  x^{\left(  n\right)  }\right)  }{\Delta_{1}x^{\left(  n\right)
}}\right)  ,
\end{align*}

Applying Lemma \ref{vLem_Taylor_extension},%
\begin{align*}
\frac{\left\vert T\right\vert }{v}\left\vert \left(  T^{-1}\Delta_{X}f\right)
_{N-1}\right\vert  & \leq\frac{x_{1}}{2^{N-2}}\left\vert \frac{\Delta
_{1}f\left(  x^{\left(  1\right)  }\right)  }{\Delta_{1}x^{\left(  1\right)
}}\right\vert +\sum\limits_{n=2}^{N-1}\frac{x_{n}}{2^{N-1-n}}\left\vert
T_{1:n-1}\right\vert \left\vert \frac{\Delta_{1}f\left(  x^{\left(  n\right)
}\right)  }{\Delta_{1}x^{\left(  n\right)  }}\right\vert \\
& \leq\max_{n=1}^{N-1}\left\vert \frac{\Delta_{1}f\left(  x^{\left(  n\right)
}\right)  }{\Delta_{1}x^{\left(  n\right)  }}\right\vert \left(  \frac{x_{1}%
}{2^{N-2}}+\sum\limits_{n=2}^{N-1}\frac{x_{n}}{2^{N-1-n}}\left\vert
T_{1:n-1}\right\vert \right) \\
& \leq\left\Vert Df\right\Vert _{\infty;\Omega}\left(  \frac{x_{1}}{2^{N-2}%
}+\sum\limits_{n=2}^{N-1}\frac{x_{n}}{2^{N-1-n}}\left\vert T_{1:n-1}%
\right\vert \right)  ,
\end{align*}

or on rearranging%
\begin{equation}
\frac{1}{v}\left\vert \left(  T^{-1}\Delta_{X}f\right)  _{N-1}\right\vert
\leq\frac{\left\Vert Df\right\Vert _{\infty;\Omega}}{\left\vert T\right\vert
}\left(  \frac{x_{1}}{2^{N-2}}+\sum\limits_{n=2}^{N-1}\frac{x_{n}}{2^{N-1-n}%
}\left\vert T_{1:n-1}\right\vert \right)  .\label{av028}%
\end{equation}

In summary, with reference to \ref{av725}, we have the estimates,%
\begin{align}
&  \frac{1}{v}\left\vert \left(  T^{-1}\Delta_{X}f\right)  _{k}\right\vert
\nonumber\\
&  \leq\left\Vert Df\right\Vert _{\infty;\Omega}\times\left\{
\begin{array}
[c]{ll}%
\frac{1}{\left\vert T\right\vert }\left(  \sum\limits_{n=1}^{N-2}\frac{x_{n}%
}{2^{n-1}}\left\vert T_{n+1:N-1}\right\vert +\frac{x_{N-1}}{2^{N-2}}\right)
, & k=1,\\
& \\
\left.
\begin{array}
[c]{l}%
\frac{\left\vert T_{k+1:N-1}\right\vert }{\left\vert T\left(  x\right)
\right\vert }\left(  \frac{x_{1}}{2^{k-1}}+\sum\limits_{n=2}^{k}\left\vert
T_{1:n-1}\right\vert \frac{x_{n}}{2^{k-n}}\right)  +\\
+\frac{\left\vert T_{1:k-1}\right\vert }{\left\vert T\right\vert }\left(
\sum\limits_{n=k+1}^{N-2}\frac{x_{n}}{2^{n-k}}\left\vert T_{n+1:N-1}%
\right\vert +\frac{x_{N-1}}{2^{N-1-k}}\right)
\end{array}
\right\}  , &
\begin{array}
[c]{l}%
k\geq2,\\
k\leq N-2,
\end{array}
\\
& \\
\frac{1}{\left\vert T\right\vert }\left(  \frac{x_{1}}{2^{N-2}}+\sum
\limits_{n=2}^{N-1}\left\vert T_{1:n-1}\right\vert \frac{x_{n}}{2^{N-1-n}%
}\right)  , & k=N-1,
\end{array}
\right. \label{av029}%
\end{align}

where (Definition \ref{vDef_T_minors}) $T=T\left(  x\right)  =T\left(
x;1\right)  $, $x=\frac{1}{v}\Delta_{1}X$, $v=\rho N$ and $T_{i:j}=T\left(
x_{i:j}\right)  =T\left(  x_{i:j};1\right)  $.

\subsection{Expansions of $\left\vert T\left(  x\right)  \right\vert $ in
terms of $\left\vert T\left(  x^{\prime},0_{n}\right)  \right\vert $,
$\left\vert T\left(  0_{n},x^{\prime\prime}\right)  \right\vert $%
\label{SbSect_detT_in_terms_of_detTzeros}}

We will proceed below by deriving the expansions \ref{av039} and \ref{av037}
of $\left\vert T\left(  x_{1:m}\right)  \right\vert $ which involve sums of
terms of the form $\left\vert T\left(  x^{\prime},0_{n}\right)  \right\vert $
and $\left\vert T\left(  0_{n},x^{\prime\prime}\right)  \right\vert $
respectively.
\begin{align}
\left\vert T\left(  x_{1:m}\right)  \right\vert  &  =%
\begin{vmatrix}
1+x_{1} & -\frac{1}{2} &  & 0 & 0 & 0 &  & 0 & 0\\
-\frac{1}{2} & 1+x_{2} &  & 0 & 0 & 0 &  & 0 & 0\\
&  & \ddots &  &  &  &  &  & \\
0 & 0 &  & 1+x_{k-1} & -\frac{1}{2} & 0 &  & 0 & 0\\
0 & 0 &  & -\frac{1}{2} & 1+x_{k} & -\frac{1}{2} &  & 0 & 0\\
0 & 0 &  & 0 & -\frac{1}{2} & 1+x_{k+1} &  & 0 & 0\\
&  &  &  &  &  & \ddots &  & \\
0 & 0 &  & 0 & 0 & 0 &  & 1+x_{m-1} & -\frac{1}{2}\\
0 & 0 &  & 0 & 0 & 0 &  & -\frac{1}{2} & 1+x_{m}%
\end{vmatrix}
\nonumber\\
& \nonumber\\
&  =%
\begin{vmatrix}
1+x_{1} & -\frac{1}{2} &  & 0 & 0 & 0 &  & 0 & 0\\
-\frac{1}{2} & 1+x_{2} &  & 0 & 0 & 0 &  & 0 & 0\\
&  & \ddots &  &  &  &  &  & \\
0 & 0 &  & 1+x_{k-1} & 0 & 0 &  & 0 & 0\\
0 & 0 &  & -\frac{1}{2} & x_{k} & -\frac{1}{2} &  & 0 & 0\\
0 & 0 &  & 0 & 0 & 1+x_{k+1} &  & 0 & 0\\
&  &  &  &  &  & \ddots &  & \\
0 & 0 &  & 0 & 0 & 0 &  & 1+x_{m-1} & -\frac{1}{2}\\
0 & 0 &  & 0 & 0 & 0 &  & -\frac{1}{2} & 1+x_{m}%
\end{vmatrix}
+\nonumber\\
& \nonumber\\
&  \qquad+%
\begin{vmatrix}
1+x_{1} & -\frac{1}{2} &  & 0 & 0 & 0 &  & 0 & 0\\
-\frac{1}{2} & 1+x_{2} &  & 0 & 0 & 0 &  & 0 & 0\\
&  & \ddots &  &  &  &  &  & \\
0 & 0 &  & 1+x_{k-1} & -\frac{1}{2} & 0 &  & 0 & 0\\
0 & 0 &  & -\frac{1}{2} & 1 & -\frac{1}{2} &  & 0 & 0\\
0 & 0 &  & 0 & -\frac{1}{2} & 1+x_{k+1} &  & 0 & 0\\
&  &  &  &  &  & \ddots &  & \\
0 & 0 &  & 0 & 0 & 0 &  & 1+x_{m-1} & -\frac{1}{2}\\
0 & 0 &  & 0 & 0 & 0 &  & -\frac{1}{2} & 1+x_{m}%
\end{vmatrix}
\nonumber\\
& \nonumber\\
&  =%
\begin{vmatrix}
1+x_{1} & -\frac{1}{2} &  & 0\\
-\frac{1}{2} & 1+x_{2} &  & 0\\
&  & \ddots & \\
0 & 0 &  & 1+x_{k-1}%
\end{vmatrix}
x_{k}%
\begin{vmatrix}
1+x_{k+1} &  & 0 & 0\\
& \ddots &  & \\
0 &  & 1+x_{m-1} & -\frac{1}{2}\\
0 &  & -\frac{1}{2} & 1+x_{m}%
\end{vmatrix}
+\nonumber\\
& \nonumber\\
&  \qquad\qquad+\left\vert T\left(  x_{1:k-1},0,x_{k+1:m}\right)  \right\vert
\nonumber\\
& \nonumber\\
&  =\left\vert T\left(  x_{1:k-1}\right)  \right\vert x_{k}\left\vert T\left(
x_{k+1:m}\right)  \right\vert +\left\vert T\left(  x_{1:k-1},0,x_{k+1:m}%
\right)  \right\vert ,\label{av873}%
\end{align}

\fbox{Case $k=1$}
\begin{align}
&  \left\vert T\left(  x_{1:m}\right)  \right\vert \nonumber\\
& \nonumber\\
&  =%
\begin{vmatrix}
1+x_{1} & -\frac{1}{2} & 0 &  & 0 & 0\\
-\frac{1}{2} & 1+x_{2} & -\frac{1}{2} &  & 0 & 0\\
0 & -\frac{1}{2} & 1+x_{3} &  & 0 & 0\\
&  &  & \ddots &  & \\
0 & 0 & 0 &  & 1+x_{m-1} & -\frac{1}{2}\\
0 & 0 & 0 &  & -\frac{1}{2} & 1+x_{m}%
\end{vmatrix}
\nonumber\\
& \nonumber\\
&  =%
\begin{vmatrix}
x_{1} & -\frac{1}{2} & 0 &  & 0 & 0\\
0 & 1+x_{2} & -\frac{1}{2} &  & 0 & 0\\
0 & -\frac{1}{2} & 1+x_{3} &  & 0 & 0\\
&  &  & \ddots &  & \\
0 & 0 & 0 &  & 1+x_{m-1} & -\frac{1}{2}\\
0 & 0 & 0 &  & -\frac{1}{2} & 1+x_{m}%
\end{vmatrix}
+%
\begin{vmatrix}
1 & -\frac{1}{2} & 0 &  & 0 & 0\\
-\frac{1}{2} & 1+x_{2} & -\frac{1}{2} &  & 0 & 0\\
0 & -\frac{1}{2} & 1+x_{3} &  & 0 & 0\\
&  &  & \ddots &  & \\
0 & 0 & 0 &  & 1+x_{m-1} & -\frac{1}{2}\\
0 & 0 & 0 &  & -\frac{1}{2} & 1+x_{m}%
\end{vmatrix}
\nonumber\\
& \nonumber\\
&  =x_{1}\left\vert T\left(  x_{2:m}\right)  \right\vert +\left\vert T\left(
0,x_{2:m}\right)  \right\vert .\label{av936}%
\end{align}

This result will be used to estimate $\mu_{1}\left(  f\right)  $.\medskip

\fbox{Case $k=m$}
\begin{align}
&  \left\vert T\left(  x_{1:m}\right)  \right\vert \nonumber\\
& \nonumber\\
&  =%
\begin{vmatrix}
1+x_{1} & -\frac{1}{2} &  & 0 & 0 & 0\\
-\frac{1}{2} & 1+x_{2} &  & 0 & 0 & 0\\
&  & \ddots &  &  & \\
0 & 0 &  & 1+x_{m-2} & -\frac{1}{2} & 0\\
0 & 0 &  & -\frac{1}{2} & 1+x_{m-1} & -\frac{1}{2}\\
0 & 0 &  & 0 & -\frac{1}{2} & 1+x_{m}%
\end{vmatrix}
\nonumber\\
& \nonumber\\
&  =%
\begin{vmatrix}
1+x_{1} & -\frac{1}{2} &  & 0 & 0 & 0\\
-\frac{1}{2} & 1+x_{2} &  & 0 & 0 & 0\\
&  & \ddots &  &  & \\
0 & 0 &  & 1+x_{m-2} & -\frac{1}{2} & 0\\
0 & 0 &  & -\frac{1}{2} & 1+x_{m-1} & 0\\
0 & 0 &  & 0 & -\frac{1}{2} & x_{m}%
\end{vmatrix}
+%
\begin{vmatrix}
1+x_{1} & -\frac{1}{2} &  & 0 & 0 & 0\\
-\frac{1}{2} & 1+x_{2} &  & 0 & 0 & 0\\
&  & \ddots &  &  & \\
0 & 0 &  & 1+x_{m-2} & -\frac{1}{2} & 0\\
0 & 0 &  & -\frac{1}{2} & 1+x_{m-1} & -\frac{1}{2}\\
0 & 0 &  & 0 & -\frac{1}{2} & 1
\end{vmatrix}
\nonumber\\
& \nonumber\\
&  =\left\vert T\left(  x_{1:m-1}\right)  \right\vert x_{m}+\left\vert
T\left(  x_{1:m-1},0\right)  \right\vert .\label{av840}%
\end{align}

This equation will be used to estimate $\mu_{1}\left(  f\right)  $.\medskip

\fbox{\textbf{Case} $k=1$} Applying \ref{av873} with $k=1$:%
\[
\left\vert T\left(  x\right)  \right\vert =x_{1}\left\vert T\left(
x_{2:m}\right)  \right\vert +\left\vert T\left(  0,x_{2:m}\right)  \right\vert
.
\]

\fbox{\textbf{Case} $k=2$} Applying \ref{av873} with $k=2$:%
\[
\left\vert T\left(  x\right)  \right\vert =x_{1}\left\vert T\left(
0,x_{3:m}\right)  \right\vert +\left\vert T\left(  x_{1}\right)  \right\vert
x_{2}\left\vert T\left(  x_{3:m}\right)  \right\vert +\left\vert T\left(
0_{2},x_{3:m}\right)  \right\vert .
\]

\fbox{\textbf{Case} $3\leq k\leq m-1$} Applying \ref{av873} to exhaustion to
it's right-most term (in order to substitute zeros to the left) we get:
\begin{align}
\left\vert T\left(  x\right)  \right\vert  & =\left\vert T\left(
x_{1:k-1}\right)  \right\vert x_{k}\left\vert T\left(  x_{k+1:m}\right)
\right\vert +\left\vert T\left(  x_{1:k-1},0,x_{k+1:m}\right)  \right\vert
\nonumber\\
& =\left\vert T\left(  x_{1:k-1}\right)  \right\vert x_{k}\left\vert T\left(
x_{k+1:m}\right)  \right\vert +\left\vert T\left(  x_{1:k-2}\right)
\right\vert x_{k-1}\left\vert T\left(  0,x_{k+1:m}\right)  \right\vert
+\nonumber\\
& \qquad\qquad+\left\vert T\left(  x_{1:k-2},0_{2},x_{k+1:m}\right)
\right\vert \nonumber\\
& =\left\vert T\left(  x_{1:k-1}\right)  \right\vert x_{k}\left\vert T\left(
x_{k+1:m}\right)  \right\vert +\left\vert T\left(  x_{1:k-2}\right)
\right\vert x_{k-1}\left\vert T\left(  0,x_{k+1:m}\right)  \right\vert
+\nonumber\\
& \qquad+\left\vert T\left(  x_{1:k-3}\right)  \right\vert x_{k-2}\left\vert
T\left(  0_{2},x_{k+1:m}\right)  \right\vert +\left\vert T\left(
x_{1:k-3},0_{3},x_{k+1:m}\right)  \right\vert \nonumber\\
& =\left\vert T\left(  x_{1:k-1}\right)  \right\vert x_{k}\left\vert T\left(
x_{k+1:m}\right)  \right\vert +\left\vert T\left(  x_{1:k-2}\right)
\right\vert x_{k-1}\left\vert T\left(  0,x_{k+1:m}\right)  \right\vert
+\nonumber\\
& \qquad+\left\vert T\left(  x_{1:k-3}\right)  \right\vert x_{k-2}\left\vert
T\left(  0_{2},x_{k+1:m}\right)  \right\vert +\ldots\nonumber\\
& \qquad\ldots+\left\vert T\left(  x_{1:j-1}\right)  \right\vert
x_{j}\left\vert T\left(  0_{k-j},x_{k+1:m}\right)  \right\vert +\left\vert
T\left(  x_{1:j-1},0_{k-j+1},x_{k+1:m}\right)  \right\vert \nonumber\\
& =\left\vert T\left(  x_{1:k-1}\right)  \right\vert x_{k}\left\vert T\left(
x_{k+1:m}\right)  \right\vert +\left\vert T\left(  x_{1:k-2}\right)
\right\vert x_{k-1}\left\vert T\left(  0,x_{k+1:m}\right)  \right\vert
+\ldots\nonumber\\
& \qquad\ldots+\left\vert T\left(  x_{1:j-1}\right)  \right\vert
x_{j}\left\vert T\left(  0_{k-j},x_{k+1:m}\right)  \right\vert +\ldots
\nonumber\\
& \qquad\ldots+\left\vert T\left(  x_{1}\right)  \right\vert x_{2}\left\vert
T\left(  0_{k-2},x_{k+1:m}\right)  \right\vert +\left\vert T\left(
x_{1},0_{k-1},x_{k+1:m}\right)  \right\vert .\label{av875}%
\end{align}

Finally, applying \ref{av936}: $\left\vert T\left(  u_{1:m}\right)
\right\vert =u_{1}\left\vert T\left(  u_{2:m}\right)  \right\vert +\left\vert
T\left(  0,u_{2:m}\right)  \right\vert $, to the last term gives%
\begin{align*}
\left\vert T\left(  x\right)  \right\vert  & =\left\vert T\left(
x_{1:k-1}\right)  \right\vert x_{k}\left\vert T\left(  x_{k+1:m}\right)
\right\vert +\left\vert T\left(  x_{1:k-2}\right)  \right\vert x_{k-1}%
\left\vert T\left(  0,x_{k+1:m}\right)  \right\vert +\\
& \qquad\ldots+\left\vert T\left(  x_{1:j-1}\right)  \right\vert
x_{j}\left\vert T\left(  0_{k-j},x_{k+1:m}\right)  \right\vert +\\
& \qquad\ldots+\left\vert T\left(  x_{1}\right)  \right\vert x_{2}\left\vert
T\left(  0_{k-2},x_{k+1:m}\right)  \right\vert +x_{1}\left\vert T\left(
0_{k-1},x_{k+1:m}\right)  \right\vert +\left\vert T\left(  0_{k}%
,x_{k+1:m}\right)  \right\vert \\
& =x_{1}\left\vert T\left(  0_{k-1},x_{k+1:m}\right)  \right\vert
+\sum\limits_{j=2}^{k-1}\left\vert T\left(  x_{1:j-1}\right)  \right\vert
x_{j}\left\vert T\left(  0_{k-j},x_{k+1:m}\right)  \right\vert +\\
& \qquad+\left\vert T\left(  x_{1:k-1}\right)  \right\vert x_{k}\left\vert
T\left(  x_{k+1:m}\right)  \right\vert +\left\vert T\left(  0_{k}%
,x_{k+1:m}\right)  \right\vert .
\end{align*}

\fbox{\textbf{Case} $k=m$} Apply \ref{av840} and then, to substitute zeros to
the left, use the consequence

$\left\vert T\left(  x_{1:p},0_{n}\right)  \right\vert =\left\vert T\left(
x_{1:p-1}\right)  \right\vert x_{p}\left\vert T\left(  0_{n}\right)
\right\vert +\left\vert T\left(  x_{1:p-1},0_{n+1}\right)  \right\vert $ of
\ref{av873} to exhaustion on the right-most terms:
\begin{align*}
\left\vert T\left(  x\right)  \right\vert  & =\left\vert T\left(
x_{1:m-1}\right)  \right\vert x_{m}+\left\vert T\left(  x_{1:m-1},0\right)
\right\vert \\
& =\left\vert T\left(  x_{1:m-1}\right)  \right\vert x_{m}+\left\vert T\left(
x_{1:m-2}\right)  \right\vert x_{m-1}\left\vert T\left(  0\right)  \right\vert
+\left\vert T\left(  x_{1:m-2},0_{2}\right)  \right\vert \\
& =\left\vert T\left(  x_{1:m-1}\right)  \right\vert x_{m}+\left\vert T\left(
x_{1:m-2}\right)  \right\vert x_{m-1}\left\vert T\left(  0\right)  \right\vert
+\left\vert T\left(  x_{1:m-3}\right)  \right\vert x_{m-2}\left\vert T\left(
0_{2}\right)  \right\vert +\\
& \qquad+\left\vert T\left(  x_{1:m-3},0_{3}\right)  \right\vert \\
& =\left\vert T\left(  x_{1:m-1}\right)  \right\vert x_{m}+\left\vert T\left(
x_{1:m-2}\right)  \right\vert x_{m-1}\left\vert T\left(  0\right)  \right\vert
+\left\vert T\left(  x_{1:m-3}\right)  \right\vert x_{m-2}\left\vert T\left(
0_{2}\right)  \right\vert +\\
& \qquad\ldots+\left\vert T\left(  x_{1}\right)  \right\vert x_{2}\left\vert
T\left(  0_{m-2}\right)  \right\vert +\left\vert T\left(  x_{1},0_{m-1}%
\right)  \right\vert \\
& =\left\vert T\left(  x_{1:m-1}\right)  \right\vert x_{m}+\left\vert T\left(
x_{1:m-2}\right)  \right\vert x_{m-1}\left\vert T\left(  0\right)  \right\vert
+\left\vert T\left(  x_{1:m-3}\right)  \right\vert x_{m-2}\left\vert T\left(
0_{2}\right)  \right\vert +\\
& \qquad\ldots+\left\vert T\left(  x_{1}\right)  \right\vert x_{2}\left\vert
T\left(  0_{m-2}\right)  \right\vert +\left\vert T\left(  x_{1},0_{m-1}%
\right)  \right\vert .
\end{align*}

Finally, applying \ref{av936} to the last term yields%
\begin{align*}
\left\vert T\left(  x\right)  \right\vert  & =\left\vert T\left(
x_{1:m-1}\right)  \right\vert x_{m}+\left\vert T\left(  x_{1:m-2}\right)
\right\vert x_{m-1}\left\vert T\left(  0\right)  \right\vert +\left\vert
T\left(  x_{1:m-3}\right)  \right\vert x_{m-2}\left\vert T\left(
0_{2}\right)  \right\vert +\\
& \qquad\ldots+\left\vert T\left(  x_{1}\right)  \right\vert x_{2}\left\vert
T\left(  0_{m-2}\right)  \right\vert +x_{1}\left\vert T\left(  0_{m-1}\right)
\right\vert +\left\vert T\left(  0_{m}\right)  \right\vert \\
& =x_{1}\left\vert T\left(  0_{m-1}\right)  \right\vert +\sum\limits_{j=2}%
^{m-1}\left\vert T\left(  x_{1:j-1}\right)  \right\vert x_{j}\left\vert
T\left(  0_{m-j}\right)  \right\vert +\left\vert T\left(  x_{1:m-1}\right)
\right\vert x_{m}+\left\vert T\left(  0_{m}\right)  \right\vert .
\end{align*}

To summarize: if $m\geq4$ then%
\begin{equation}
\left\vert T\left(  x\right)  \right\vert =\left\{
\begin{array}
[c]{ll}%
x_{1}\left\vert T\left(  x_{2:m}\right)  \right\vert +\left\vert T\left(
0,x_{2:m}\right)  \right\vert , & k=1,\\
& \\
\left.
\begin{array}
[c]{l}%
x_{1}\left\vert T\left(  0,x_{3:m}\right)  \right\vert +\left\vert T\left(
x_{1}\right)  \right\vert x_{2}\left\vert T\left(  x_{3:m}\right)  \right\vert
+\\
\qquad\qquad\qquad+\left\vert T\left(  0_{2},x_{3:m}\right)  \right\vert
\end{array}
\right\}  , & k=2,\\
& \\
\left.
\begin{array}
[c]{l}%
x_{1}\left\vert T\left(  0_{k-1},x_{k+1:m}\right)  \right\vert +\\
\quad+\sum\limits_{j=2}^{k-1}\left\vert T\left(  x_{1:j-1}\right)  \right\vert
x_{j}\left\vert T\left(  0_{k-j},x_{k+1:m}\right)  \right\vert +\\
\quad+\left\vert T\left(  x_{1:k-1}\right)  \right\vert x_{k}\left\vert
T\left(  x_{k+1:m}\right)  \right\vert +\left\vert T\left(  0_{k}%
,x_{k+1:m}\right)  \right\vert
\end{array}
\right\}  , & 3\leq k\leq m-1,\\
& \\
\left.
\begin{array}
[c]{l}%
x_{1}\left\vert T\left(  0_{m-1}\right)  \right\vert +\sum\limits_{j=2}%
^{m-1}\left\vert T\left(  x_{1:j-1}\right)  \right\vert x_{j}\left\vert
T\left(  0_{m-j}\right)  \right\vert +\\
\qquad\qquad\qquad+\left\vert T\left(  x_{1:m-1}\right)  \right\vert
x_{m}+\left\vert T\left(  0_{m}\right)  \right\vert
\end{array}
\right\}  , & k=m.
\end{array}
\right. \label{av039}%
\end{equation}

Observe that if we define $\sum\limits_{j=2}^{1}\ldots=\sum\limits_{j=2}%
^{0}\ldots=0$ then substituting $k=1$, $k=2$ and $k=m$ in the case $3\leq
k\leq m-1$ and using the notation \ref{av057} and \ref{av050} of Definition
\ref{vDef_T_minors} we get the results for $k=1$, $k=2$ and $k=m $ given above
in \ref{av039}.\medskip

The next result \ref{av037} can be derived in a similar manner to \ref{av039}
by substituting zeros to the right, or using the reverse order (a reflection)
permutation and part 2 of Corollary \ref{vCor_1_formula_for_C^(m)_gam}.

\fbox{\textbf{Case} $k=1$} From part 2 of Corollary
\ref{vCor_1_formula_for_C^(m)_gam} and case $k=1$ of \ref{av039}:%
\begin{align*}
\left\vert T\left(  x\right)  \right\vert =\left\vert T\left(  \mathcal{R}%
x\right)  \right\vert  & =\left(  \mathcal{R}x\right)  _{1}\left\vert T\left(
\left(  \mathcal{R}x\right)  _{2:m}\right)  \right\vert +\left\vert T\left(
0,\left(  \mathcal{R}x\right)  _{2:m}\right)  \right\vert \\
& =x_{m}\left\vert T\left(  x_{m-1:1}\right)  \right\vert +\left\vert T\left(
0,x_{m-1:1}\right)  \right\vert \\
& =x_{m}\left\vert T\left(  \mathcal{R}\left(  x_{m-1:1}\right)  \right)
\right\vert +\left\vert T\left(  \mathcal{R}\left(  0,x_{m-1:1}\right)
\right)  \right\vert \\
& =\left\vert T\left(  x_{1:m-1}\right)  \right\vert x_{m}+\left\vert T\left(
x_{1:m-1},0\right)  \right\vert \\
& =\left\vert T\left(  x_{1:m-1},0\right)  \right\vert +\left\vert T\left(
x_{1:m-1}\right)  \right\vert x_{m}.
\end{align*}

\fbox{\textbf{Case} $k=2$} From part 2 of Corollary
\ref{vCor_1_formula_for_C^(m)_gam} and case $k=2$ of \ref{av039}:%
\begin{align*}
\left\vert T\left(  x_{1:m}\right)  \right\vert  & =\left\vert T\left(
\mathcal{R}x\right)  \right\vert \\
& =\left(  \mathcal{R}x\right)  _{1}\left\vert T\left(  0,\left(
\mathcal{R}x\right)  _{3:m}\right)  \right\vert +\left\vert T\left(  \left(
\mathcal{R}x\right)  _{1}\right)  \right\vert \left(  \mathcal{R}x\right)
_{2}\left\vert T\left(  \left(  \mathcal{R}x\right)  _{3:m}\right)
\right\vert +\left\vert T\left(  0_{2},\left(  \mathcal{R}x\right)
_{3:m}\right)  \right\vert \\
& =x_{m}\left\vert T\left(  0,x_{m-2:1}\right)  \right\vert +\left\vert
T\left(  x_{m}\right)  \right\vert x_{m-1}\left\vert T\left(  x_{m-2:1}%
\right)  \right\vert +\left\vert T\left(  0_{2},x_{m-2:1}\right)  \right\vert
\\
& =x_{m}\left\vert T\left(  \mathcal{R}\left(  0,x_{m-2:1}\right)  \right)
\right\vert +\left\vert T\left(  x_{m}\right)  \right\vert x_{m-1}\left\vert
T\left(  \mathcal{R}\left(  x_{m-2:1}\right)  \right)  \right\vert +\left\vert
T\left(  \mathcal{R}\left(  0_{2},x_{m-2:1}\right)  \right)  \right\vert \\
& =x_{m}\left\vert T\left(  x_{1:m-2},0\right)  \right\vert +\left\vert
T\left(  x_{m}\right)  \right\vert x_{m-1}\left\vert T\left(  x_{1:m-2}%
\right)  \right\vert +\left\vert T\left(  x_{1:m-2},0_{2}\right)  \right\vert
\\
& =\left\vert T\left(  x_{1:m-2},0\right)  \right\vert x_{m}+\left\vert
T\left(  x_{1:m-2}\right)  \right\vert x_{m-1}\left\vert T\left(
x_{m}\right)  \right\vert +\left\vert T\left(  x_{1:m-2},0_{2}\right)
\right\vert \\
& =\left\vert T\left(  x_{1:m-2}\right)  \right\vert x_{m-1}\left\vert
T\left(  x_{m}\right)  \right\vert +\left\vert T\left(  x_{1:m-2},0\right)
\right\vert x_{m}+\left\vert T\left(  x_{1:m-2},0_{2}\right)  \right\vert .
\end{align*}

\fbox{\textbf{Case} $3\leq k\leq m-1$} Applying part 2 of Corollary
\ref{vCor_1_formula_for_C^(m)_gam} to case $3\leq k\leq m-1$ of \ref{av039}:%
\begin{align*}
\left\vert T\left(  x\right)  \right\vert  &  =\left\vert T\left(
\mathcal{R}\left(  x\right)  \right)  \right\vert \\
&  =\left(  \mathcal{R}x\right)  _{1}\left\vert T\left(  0_{k-1},\left(
\mathcal{R}x\right)  _{k+1:m}\right)  \right\vert +\sum\limits_{j=2}%
^{k-1}\left\vert T\left(  \left(  \mathcal{R}x\right)  _{1:j-1}\right)
\right\vert \left(  \mathcal{R}x\right)  _{j}\left\vert T\left(
0_{k-j},\left(  \mathcal{R}x\right)  _{k+1:m}\right)  \right\vert +\\
&  \qquad\qquad+\left\vert T\left(  \left(  \mathcal{R}x\right)
_{1:k-1}\right)  \right\vert \left(  \mathcal{R}x\right)  _{k}\left\vert
T\left(  \left(  \mathcal{R}x\right)  _{k+1:m}\right)  \right\vert +\left\vert
T\left(  0_{k},\left(  \mathcal{R}x\right)  _{k+1:m}\right)  \right\vert \\
&  =x_{m}\left\vert T\left(  0_{k-1},x_{m-k:1}\right)  \right\vert
+\sum\limits_{j=2}^{k-1}\left\vert T\left(  x_{m:m-j}\right)  \right\vert
x_{m+1-j}\left\vert T\left(  0_{k-j},x_{m-k:1}\right)  \right\vert +\\
&  \qquad\qquad+\left\vert T\left(  x_{m:m-k+2}\right)  \right\vert
x_{m-k+1}\left\vert T\left(  x_{m-k:1}\right)  \right\vert +\left\vert
T\left(  0_{k},x_{m-k:1}\right)  \right\vert \\
&  =x_{m}\left\vert T\left(  \mathcal{R}\left(  0_{k-1},x_{m-k:1}\right)
\right)  \right\vert +\sum\limits_{j=2}^{k-1}\left\vert T\left(
\mathcal{R}\left(  x_{m:m-j}\right)  \right)  \right\vert x_{m+1-j}\left\vert
T\left(  \mathcal{R}\left(  0_{k-j},x_{m-k:1}\right)  \right)  \right\vert +\\
&  \qquad\qquad+\left\vert T\left(  \mathcal{R}\left(  x_{m:m-k+2}\right)
\right)  \right\vert x_{m-k+1}\left\vert T\left(  \mathcal{R}\left(
x_{m-k:1}\right)  \right)  \right\vert +\left\vert T\left(  \mathcal{R}\left(
0_{k},x_{m-k:1}\right)  \right)  \right\vert \\
&  =x_{m}\left\vert T\left(  x_{1:m-k},0_{k-1}\right)  \right\vert
+\sum\limits_{j=2}^{k-1}\left\vert T\left(  x_{m-j:m}\right)  \right\vert
x_{m+1-j}\left\vert T\left(  x_{1:m-k},0_{k-j}\right)  \right\vert +\\
&  \qquad\qquad+\left\vert T\left(  x_{m-k+2:m}\right)  \right\vert
x_{m-k+1}\left\vert T\left(  x_{1:m-k}\right)  \right\vert +\left\vert
T\left(  x_{1:m-k},0_{k}\right)  \right\vert \\
&  =\left\vert T\left(  x_{1:m-k},0_{k-1}\right)  \right\vert x_{m}%
+\sum\limits_{j=2}^{k-1}\left\vert T\left(  x_{1:m-k},0_{k-j}\right)
\right\vert x_{m+1-j}\left\vert T\left(  x_{m-j:m}\right)  \right\vert +\\
&  \qquad\qquad+\left\vert T\left(  x_{1:m-k}\right)  \right\vert
x_{m-k+1}\left\vert T\left(  x_{m-k+2:m}\right)  \right\vert +\left\vert
T\left(  x_{1:m-k},0_{k}\right)  \right\vert \\
&  =\left\vert T\left(  x_{1:m-k}\right)  \right\vert x_{m-k+1}\left\vert
T\left(  x_{m-k+2:m}\right)  \right\vert +\sum\limits_{l=m-k+2}^{m-1}%
\left\vert T\left(  x_{1:m-k},0_{l+k-m-1}\right)  \right\vert x_{l}\left\vert
T\left(  x_{l+1:m}\right)  \right\vert +\\
&  \qquad\qquad+\left\vert T\left(  x_{1:m-k},0_{k-1}\right)  \right\vert
x_{m}+\left\vert T\left(  x_{1:m-k},0_{k}\right)  \right\vert .
\end{align*}

The change of summation variable $k\rightarrow m-k+2$ now yields%
\begin{align*}
\left\vert T\left(  x\right)  \right\vert  & =\left\vert T\left(
x_{1:k-2}\right)  \right\vert x_{k-1}\left\vert T\left(  x_{k:m}\right)
\right\vert +\sum\limits_{l=k}^{m-1}\left\vert T\left(  x_{1:k-2}%
,0_{l-k+1}\right)  \right\vert x_{l}\left\vert T\left(  x_{l+1:m}\right)
\right\vert +\\
& \qquad\qquad+\left\vert T\left(  x_{1:k-2},0_{m-k+1}\right)  \right\vert
x_{m}+\left\vert T\left(  x_{1:k-2},0_{m-k+2}\right)  \right\vert .
\end{align*}
\medskip

\fbox{\textbf{Case} $k=m$} From part 2 of Corollary
\ref{vCor_1_formula_for_C^(m)_gam}: $\left\vert T\left(  \mathcal{R}x\right)
\right\vert =\left\vert T\left(  x\right)  \right\vert $, and so starting with
case $k=m$ of \ref{av039}:%
\begin{align*}
\left\vert T\left(  x\right)  \right\vert  & =\left\vert T\left(
\mathcal{R}x\right)  \right\vert =\left(  \mathcal{R}x\right)  _{1}\left\vert
T\left(  0_{m-1}\right)  \right\vert +\sum\limits_{j=2}^{m-1}\left\vert
T\left(  \left(  \mathcal{R}x\right)  _{1:j-1}\right)  \right\vert \left(
\mathcal{R}x\right)  _{j}\left\vert T\left(  0_{m-j}\right)  \right\vert +\\
& \qquad\qquad+\left\vert T\left(  \left(  \mathcal{R}x\right)  _{1:m-1}%
\right)  \right\vert \left(  \mathcal{R}x\right)  _{m}+\left\vert T\left(
0_{m}\right)  \right\vert \\
& =x_{m}\left\vert T\left(  0_{m-1}\right)  \right\vert +\sum\limits_{j=2}%
^{m-1}\left\vert T\left(  x_{m:m-j+2}\right)  \right\vert x_{m-j+1}\left\vert
T\left(  0_{m-j}\right)  \right\vert +\left\vert T\left(  x_{m:2}\right)
\right\vert x_{1}+\left\vert T\left(  0_{m}\right)  \right\vert \\
& =x_{1}\left\vert T\left(  x_{m:2}\right)  \right\vert +\sum\limits_{j=2}%
^{m-1}\left\vert T\left(  x_{m:m-j+2}\right)  \right\vert x_{m-j+1}\left\vert
T\left(  0_{m-j}\right)  \right\vert +\left\vert T\left(  0_{m-1}\right)
\right\vert x_{m}+\left\vert T\left(  0_{m}\right)  \right\vert \\
& =x_{1}\left\vert T\left(  \mathcal{R}\left(  x_{m:2}\right)  \right)
\right\vert +\sum\limits_{j=2}^{m-1}\left\vert T\left(  \mathcal{R}\left(
x_{m:m-j+2}\right)  \right)  \right\vert x_{m-j+1}\left\vert T\left(
0_{m-j}\right)  \right\vert +\\
& \qquad+\left\vert T\left(  0_{m-1}\right)  \right\vert x_{m}+\left\vert
T\left(  0_{m}\right)  \right\vert \\
& =x_{1}\left\vert T\left(  x_{2:m}\right)  \right\vert +\sum\limits_{j=2}%
^{m-1}\left\vert T\left(  x_{m-j+2:m}\right)  \right\vert x_{m-j+1}\left\vert
T\left(  0_{m-j}\right)  \right\vert +\left\vert T\left(  0_{m-1}\right)
\right\vert x_{m}+\left\vert T\left(  0_{m}\right)  \right\vert \\
& =x_{1}\left\vert T\left(  x_{2:m}\right)  \right\vert +\sum\limits_{l=m-1}%
^{2}\left\vert T\left(  x_{l+1:m}\right)  \right\vert x_{l}\left\vert T\left(
0_{l-1}\right)  \right\vert +\left\vert T\left(  0_{m-1}\right)  \right\vert
x_{m}+\left\vert T\left(  0_{m}\right)  \right\vert \\
& =x_{1}\left\vert T\left(  x_{2:m}\right)  \right\vert +\sum\limits_{l=2}%
^{m-1}\left\vert T\left(  0_{l-1}\right)  \right\vert x_{l}\left\vert T\left(
x_{l+1:m}\right)  \right\vert +\left\vert T\left(  0_{m-1}\right)  \right\vert
x_{m}+\left\vert T\left(  0_{m}\right)  \right\vert \\
& =\sum\limits_{l=1}^{m-1}\frac{x_{l}}{2^{l-1}}\left\vert T\left(
x_{l+1:m}\right)  \right\vert +\frac{m}{2^{m-1}}x_{m}+\frac{m+1}{2^{m}}.
\end{align*}

To summarize: if $m\geq4$ then:%
\begin{equation}
\left\vert T\left(  x\right)  \right\vert =\left\{
\begin{array}
[c]{ll}%
\left\vert T\left(  x_{1:m-1}\right)  \right\vert x_{m}+\left\vert T\left(
x_{1:m-1},0\right)  \right\vert , & k=1,\\
& \\
\left.
\begin{array}
[c]{l}%
\left\vert T\left(  x_{1:m-2}\right)  \right\vert x_{m-1}\left\vert T\left(
x_{m}\right)  \right\vert +\left\vert T\left(  x_{1:m-2},0\right)  \right\vert
x_{m}+\\
\qquad\qquad\qquad+\left\vert T\left(  x_{1:m-2},0_{2}\right)  \right\vert ,
\end{array}
\right\}  , & k=2,\\
& \\
\left.
\begin{array}
[c]{l}%
\left\vert T\left(  x_{1:k-2}\right)  \right\vert x_{k-1}\left\vert T\left(
x_{k:m}\right)  \right\vert +\\
\quad+\sum\limits_{l=k}^{m-1}\left\vert T\left(  x_{1:k-2},0_{l-k+1}\right)
\right\vert x_{l}\left\vert T\left(  x_{l+1:m}\right)  \right\vert +\\
\quad+\left\vert T\left(  x_{1:k-2},0_{m-k+1}\right)  \right\vert
x_{m}+\left\vert T\left(  x_{1:k-2},0_{m-k+2}\right)  \right\vert
\end{array}
\right\}  , & 3\leq k\leq m-1,\\
& \\
\sum\limits_{l=1}^{m-1}\frac{x_{l}}{2^{l-1}}\left\vert T\left(  x_{l+1:m}%
\right)  \right\vert +\frac{m}{2^{m-1}}x_{m}+\frac{m+1}{2^{m}}, & k=m.
\end{array}
\right. \label{av037}%
\end{equation}

\subsection{Upper and lower bounds for $\left\vert T\left(  0_{n},x\right)
\right\vert $ and $\left\vert T\left(  x,0_{n}\right)  \right\vert
$\label{SbSect_bnds_detT_with_zeros}}

The next step is to obtain bounds for determinants with zeros in their
arguments e.g. $\left\vert T\left(  0_{n},x\right)  \right\vert \geq
c_{n}\left\vert T\left(  x\right)  \right\vert $ and $\left\vert T\left(
x,0_{n}\right)  \right\vert \geq c_{n}^{\prime}\left\vert T\left(  x\right)
\right\vert $. In fact we will eventually prove Theorem
\ref{vThm_low_bnd_remove_zeros}.

\begin{theorem}
\label{vThm_estim_prod_C^m_and_b_C^n_c}For $m,n\geq1$:\
\begin{align}
\frac{1+m-\left\vert \beta\right\vert +n-\left\vert \gamma\right\vert
}{\left(  1+m-\left\vert \beta\right\vert \right)  \left(  1+n-\left\vert
\gamma\right\vert \right)  }c_{\beta}^{\left(  m\right)  }c_{\gamma}^{\left(
n\right)  }  & =\frac{1+\left\vert \mathbf{1}-\left(  \beta,\gamma\right)
\right\vert }{\left(  1+\left\vert \mathbf{1}-\beta\right\vert \right)
\left(  1+\left\vert \mathbf{1}-\gamma\right\vert \right)  }c_{\beta}^{\left(
m\right)  }c_{\gamma}^{\left(  n\right)  }\nonumber\\
& \leq c_{\left(  \beta,\gamma\right)  }^{\left(  m+n\right)  }\label{av842}\\
& \leq c_{\beta}^{\left(  m\right)  }c_{\gamma}^{\left(  n\right)
},\label{av848}%
\end{align}

and the inequalities are sharp. Indeed $\frac{1+m+n}{\left(  1+m\right)
\left(  1+n\right)  }c_{\mathbf{0}}^{\left(  m\right)  }c_{\mathbf{0}%
}^{\left(  n\right)  }=c_{\mathbf{0}}^{\left(  m+n\right)  }$ and
$c_{\mathbf{1}}^{\left(  m+n\right)  }=c_{\mathbf{1}}^{\left(  m\right)
}c_{\mathbf{1}}^{\left(  n\right)  }=1$.

Also,%
\begin{equation}
\frac{1}{2}c_{\beta}^{\left(  m\right)  }c_{\gamma}^{\left(  n\right)  }\leq
c_{\left(  \beta,\gamma\right)  }^{\left(  m+n\right)  }\leq c_{\beta
}^{\left(  m\right)  }c_{\gamma}^{\left(  n\right)  },\label{av849}%
\end{equation}

except when $\beta=\left(  \ldots,0,0\right)  $ and $\gamma=\left(
0,0,\ldots\right)  $.
\end{theorem}

\begin{proof}
From Theorem \ref{vThm_formula_for_C^(m)_idx}:
\begin{align}
c_{\beta}^{\left(  m\right)  }  & =\left\{
\begin{array}
[c]{ll}%
1, & only\text{ }isolated\text{ }zeros,\\
\frac{1+\left\Vert z^{\left(  1\right)  }\right\Vert }{2^{\left\Vert
z^{\left(  1\right)  }\right\Vert }}\times\ldots\times\frac{1+\left\Vert
z^{\left(  k\right)  }\right\Vert }{2^{\left\Vert z^{\left(  k\right)
}\right\Vert }}, & otherw_{v}ise,
\end{array}
\right. \label{av794}\\
c_{\gamma}^{\left(  n\right)  }  & =\left\{
\begin{array}
[c]{ll}%
1, & only\text{ }isolated\text{ }zeros,\\
\frac{1+\left\Vert \zeta^{\left(  1\right)  }\right\Vert }{2^{\left\Vert
\zeta^{\left(  1\right)  }\right\Vert }}\times\ldots\times\frac{1+\left\Vert
\zeta^{\left(  l\right)  }\right\Vert }{2^{\left\Vert \zeta^{\left(  l\right)
}\right\Vert }}, & otherw_{v}ise.
\end{array}
\right. \label{av804}%
\end{align}

Suppose $\beta\in\mathbb{R}^{m}$ and $\gamma\in\mathbb{R}^{n}$.\medskip

\fbox{\textbf{Case 1} $\beta$,$\gamma$ both have only isolated zeros.} Here
$c_{\beta}^{\left(  m\right)  }=c_{\gamma}^{\left(  n\right)  }=1$.\medskip

\qquad\fbox{\textbf{Subcase 1.1} $\beta=\left(  \ldots,1\right)  $ and
$\gamma=\left(  1,\ldots\right)  $} $c_{\left(  \beta,\gamma\right)
}^{\left(  m+n\right)  }=1=c_{\beta}^{\left(  m\right)  }c_{\gamma}^{\left(
n\right)  }$.\medskip

\qquad\fbox{\textbf{Subcase 1.2} $\beta=\left(  \ldots,1\right)  $ and
$\gamma=\left(  0,1,\ldots\right)  $} $c_{\left(  \beta,\gamma\right)
}^{\left(  m+n\right)  }=1=c_{\beta}^{\left(  m\right)  }c_{\gamma}^{\left(
n\right)  }$.\medskip

\qquad\fbox{\textbf{Subcase 1.3} $\beta=\left(  \ldots,1,0\right)  $ and
$\gamma=\left(  1,\ldots\right)  $} $c_{\left(  \beta,\gamma\right)
}^{\left(  m+n\right)  }=1=c_{\beta}^{\left(  m\right)  }c_{\gamma}^{\left(
n\right)  }$.\medskip

\qquad\fbox{\textbf{Subcase 1.4} $\beta=\left(  \ldots,1,0\right)  $ and
$\gamma=\left(  0,1,\ldots\right)  $} Here $\left(  \beta,\gamma\right)  $ has
the single non-trivial zero sequence $\left(  0,0\right)  $ so

$c_{\left(  \beta,\gamma\right)  }^{\left(  m+n\right)  }=\frac{3}{4}=\frac
{3}{4}c_{\beta}^{\left(  m\right)  }c_{\gamma}^{\left(  n\right)  }$.\medskip

\fbox{\textbf{Case 2} $\beta$ has only isolated zeros but $\gamma$ has a
non-trivial zero sequence.} Here $c_{\beta}^{\left(  m\right)  }=1$.\medskip

\qquad\fbox{\textbf{Subcase 2.1} $\beta=\left(  \ldots,1\right)  $}
$c_{\left(  \beta,\gamma\right)  }^{\left(  m+n\right)  }=c_{\gamma}^{\left(
n\right)  }=c_{\beta}^{\left(  m\right)  }c_{\gamma}^{\left(  n\right)  }%
$.\medskip

\qquad\fbox{\textbf{Subcase 2.2} $\beta=\left(  \ldots,1,0\right)  $ and
$\gamma=\left(  1,\ldots\right)  $} $c_{\left(  \beta,\gamma\right)
}^{\left(  m+n\right)  }=c_{\gamma}^{\left(  n\right)  }=c_{\beta}^{\left(
m\right)  }c_{\gamma}^{\left(  n\right)  }$.\medskip

\qquad\fbox{\textbf{Subcase 2.3} $\beta=\left(  \ldots,1,0\right)  $ and
$\gamma=\left(  0,1,\ldots\right)  $} $c_{\left(  \beta,\gamma\right)
}^{\left(  m+n\right)  }=\frac{3}{4}c_{\gamma}^{\left(  n\right)  }=\frac
{3}{4}c_{\beta}^{\left(  m\right)  }c_{\gamma}^{\left(  n\right)  }$.\medskip

\qquad\fbox{\textbf{Subcase 2.4} $\beta=\left(  \ldots,1,0\right)  $ and
$\gamma=\left(  0,0,\ldots\right)  $} Here $\left(  \beta,\gamma\right)  $ has
the single non-trivial zero sequence $\left(  0,0\right)  $ so
\begin{align*}
c_{\left(  \beta,\gamma\right)  }^{\left(  m+n\right)  }  & =\frac
{1+1+\left\Vert \zeta^{\left(  1\right)  }\right\Vert }{2^{1+\left\Vert
\zeta^{\left(  1\right)  }\right\Vert }}\times\ldots\times\frac{1+\left\Vert
\zeta^{\left(  l\right)  }\right\Vert }{2^{\left\Vert \zeta^{\left(  l\right)
}\right\Vert }}\\
& =\frac{2+\left\Vert \zeta^{\left(  1\right)  }\right\Vert }{2^{1+\left\Vert
\zeta^{\left(  1\right)  }\right\Vert }}\times\ldots\times\frac{1+\left\Vert
\zeta^{\left(  l\right)  }\right\Vert }{2^{\left\Vert \zeta^{\left(  l\right)
}\right\Vert }}\\
& =\frac{1}{2}\frac{2+\left\Vert \zeta^{\left(  1\right)  }\right\Vert
}{1+\left\Vert \zeta^{\left(  1\right)  }\right\Vert }\frac{1+\left\Vert
\zeta^{\left(  1\right)  }\right\Vert }{2^{\left\Vert \zeta^{\left(  1\right)
}\right\Vert }}\times\ldots\times\frac{1+\left\Vert \zeta^{\left(  l\right)
}\right\Vert }{2^{\left\Vert \zeta^{\left(  l\right)  }\right\Vert }}\\
& =\frac{1}{2}\frac{2+\left\Vert \zeta^{\left(  1\right)  }\right\Vert
}{1+\left\Vert \zeta^{\left(  1\right)  }\right\Vert }c_{\gamma}^{\left(
n\right)  }\\
& =\frac{1+\left\Vert \zeta^{\left(  1\right)  }\right\Vert /2}{1+\left\Vert
\zeta^{\left(  1\right)  }\right\Vert }c_{\gamma}^{\left(  n\right)  }\text{,}%
\end{align*}

and hence%
\[
\frac{1}{2}c_{\beta}^{\left(  m\right)  }c_{\gamma}^{\left(  n\right)
}<c_{\left(  \beta,\gamma\right)  }^{\left(  m+n\right)  }<c_{\beta}^{\left(
m\right)  }c_{\gamma}^{\left(  n\right)  }.
\]

\fbox{\textbf{Case 3} $\beta$ has a non-trivial zero sequence but $\gamma$ has
only isolated zeros.} In a similar manner to Case 2 we obtain
\[
\frac{1}{2}c_{\beta}^{\left(  m\right)  }c_{\gamma}^{\left(  n\right)
}<c_{\left(  \beta,\gamma\right)  }^{\left(  m+n\right)  }<c_{\beta}^{\left(
m\right)  }c_{\gamma}^{\left(  n\right)  }.
\]
\medskip

\fbox{\textbf{Case 4} $\beta$,$\gamma$ both have non-trivial zero sequences.}
Noting symmetry, there are five cases:\medskip

\qquad\fbox{\textbf{Subcase 4.1} $\beta=\left(  \ldots,1\right)  $} The zero
sequences of $\left(  \beta,\gamma\right)  $ are $z^{\left(  1\right)
},\ldots,z^{\left(  k\right)  },\zeta^{\left(  1\right)  },\ldots
,\zeta^{\left(  k\right)  }$ and hence
\begin{align*}
c_{\left(  \beta,\gamma\right)  }^{\left(  m+n\right)  }  & =\frac
{1+\left\Vert z^{\left(  1\right)  }\right\Vert }{2^{\left\Vert z^{\left(
1\right)  }\right\Vert }}\times\ldots\times\frac{1+\left\Vert z^{\left(
k\right)  }\right\Vert }{2^{\left\Vert z^{\left(  k\right)  }\right\Vert }%
}\frac{1+\left\Vert \zeta^{\left(  1\right)  }\right\Vert }{2^{\left\Vert
\zeta^{\left(  1\right)  }\right\Vert }}\times\ldots\times\frac{1+\left\Vert
\zeta^{\left(  l\right)  }\right\Vert }{2^{\left\Vert \zeta^{\left(  l\right)
}\right\Vert }}\\
& =c_{\beta}^{\left(  m\right)  }c_{\gamma}^{\left(  n\right)  }.
\end{align*}
\medskip

\qquad\fbox{\textbf{Subcase 4.2} $\beta=\left(  \ldots,1,0\right)  $ and
$\gamma=\left(  1,\ldots\right)  $}

As in Subcase 4.1 the zero sequences of $\left(  \beta,\gamma\right)  $ are
$z^{\left(  1\right)  },\ldots,z^{\left(  k\right)  },\zeta^{\left(  1\right)
},\ldots,\zeta^{\left(  k\right)  }$ and hence
\[
c_{\left(  \beta,\gamma\right)  }^{\left(  m+n\right)  }=c_{\beta}^{\left(
m\right)  }c_{\gamma}^{\left(  n\right)  }.
\]
\medskip

\qquad\fbox{\textbf{Subcase 4.3} $\beta=\left(  \ldots,1,0\right)  $ and
$\gamma=\left(  0,1,\ldots\right)  $}

The zero sequences of $\left(  \beta,\gamma\right)  $ are $z^{\left(
1\right)  },\ldots,z^{\left(  k\right)  },\left(  0,0\right)  ,$
$\zeta^{\left(  1\right)  },\ldots,\zeta^{\left(  l\right)  }$ and hence%
\begin{align*}
c_{\left(  \beta,\gamma\right)  }^{\left(  m+n\right)  }  & =\frac
{1+\left\Vert z^{\left(  1\right)  }\right\Vert }{2^{\left\Vert z^{\left(
1\right)  }\right\Vert }}\times\ldots\times\frac{1+\left\Vert z^{\left(
k\right)  }\right\Vert }{2^{\left\Vert z^{\left(  k\right)  }\right\Vert }%
}\frac{3}{4}\frac{1+\left\Vert \zeta^{\left(  1\right)  }\right\Vert
}{2^{\left\Vert \zeta^{\left(  1\right)  }\right\Vert }}\times\ldots
\times\frac{1+\left\Vert \zeta^{\left(  l\right)  }\right\Vert }{2^{\left\Vert
\zeta^{\left(  l\right)  }\right\Vert }}\\
& =\frac{3}{4}c_{\beta}^{\left(  m\right)  }c_{\gamma}^{\left(  n\right)  },
\end{align*}
\medskip

\qquad\fbox{\textbf{Subcase 4.4} $\beta=\left(  \ldots,1,0\right)  $ and
$\gamma=\left(  0,0,\ldots\right)  $}

The zero sequences of $\left(  \beta,\gamma\right)  $ are $z^{\left(
1\right)  },\ldots,z^{\left(  k\right)  },\left(  0,\zeta^{\left(  1\right)
}\right)  ,\ldots,\zeta^{\left(  l\right)  }$ and hence%
\begin{align*}
c_{\left(  \beta,\gamma\right)  }^{\left(  m+n\right)  }  & =\frac
{1+\left\Vert z^{\left(  1\right)  }\right\Vert }{2^{\left\Vert z^{\left(
1\right)  }\right\Vert }}\times\ldots\times\frac{1+\left\Vert z^{\left(
k\right)  }\right\Vert }{2^{\left\Vert z^{\left(  k\right)  }\right\Vert }%
}\frac{2+\left\Vert \zeta^{\left(  1\right)  }\right\Vert }{2^{1+\left\Vert
\zeta^{\left(  1\right)  }\right\Vert }}\times\ldots\times\frac{1+\left\Vert
\zeta^{\left(  l\right)  }\right\Vert }{2^{\left\Vert \zeta^{\left(  l\right)
}\right\Vert }}\\
& =\frac{1}{2}\frac{2+\left\Vert \zeta^{\left(  1\right)  }\right\Vert
}{1+\left\Vert \zeta^{\left(  1\right)  }\right\Vert }\frac{1+\left\Vert
z^{\left(  1\right)  }\right\Vert }{2^{\left\Vert z^{\left(  1\right)
}\right\Vert }}\times\ldots\times\frac{1+\left\Vert z^{\left(  k\right)
}\right\Vert }{2^{\left\Vert z^{\left(  k\right)  }\right\Vert }}%
\frac{1+\left\Vert \zeta^{\left(  1\right)  }\right\Vert }{2^{1+\left\Vert
\zeta^{\left(  1\right)  }\right\Vert }}\times\ldots\times\frac{1+\left\Vert
\zeta^{\left(  l\right)  }\right\Vert }{2^{\left\Vert \zeta^{\left(  l\right)
}\right\Vert }}\\
& =\frac{1}{2}\frac{2+\left\Vert \zeta^{\left(  1\right)  }\right\Vert
}{1+\left\Vert \zeta^{\left(  1\right)  }\right\Vert }c_{\beta}^{\left(
m\right)  }c_{\gamma}^{\left(  n\right)  },
\end{align*}

so that%
\[
\frac{1}{2}c_{\beta}^{\left(  m\right)  }c_{\gamma}^{\left(  n\right)
}<c_{\left(  \beta,\gamma\right)  }^{\left(  m+n\right)  }<c_{\beta}^{\left(
m\right)  }c_{\gamma}^{\left(  n\right)  }.
\]
\medskip

\qquad\fbox{\textbf{Subcase 4.5} $\beta=\left(  \ldots,0,0\right)  $ and
$\gamma=\left(  0,0,\ldots\right)  $}

The zero sequences of $\left(  \beta,\gamma\right)  $ are $z^{\left(
1\right)  },\ldots,z^{\left(  k-1\right)  },\left(  z^{\left(  k\right)
},\zeta^{\left(  1\right)  }\right)  ,\zeta^{\left(  2\right)  },\ldots
,\zeta^{\left(  l\right)  }$ and hence
\begin{align*}
c_{\left(  \beta,\gamma\right)  }^{\left(  m+n\right)  }  & =\frac
{1+\left\Vert z^{\left(  1\right)  }\right\Vert }{2^{\left\Vert z^{\left(
1\right)  }\right\Vert }}\times\ldots\times\frac{1+\left\Vert z^{\left(
k-1\right)  }\right\Vert }{2^{\left\Vert z^{\left(  k-1\right)  }\right\Vert
}}\frac{1+\left\Vert z^{\left(  k\right)  },\zeta^{\left(  1\right)
}\right\Vert }{2^{\left\Vert z^{\left(  k\right)  },\zeta^{\left(  1\right)
}\right\Vert }}\frac{1+\left\Vert \zeta^{\left(  2\right)  }\right\Vert
}{2^{\left\Vert \zeta^{\left(  2\right)  }\right\Vert }}\times\ldots
\times\frac{1+\left\Vert \zeta^{\left(  l\right)  }\right\Vert }{2^{\left\Vert
\zeta^{\left(  l\right)  }\right\Vert }}\\
& =\frac{1+\left\Vert z^{\left(  1\right)  }\right\Vert }{2^{\left\Vert
z^{\left(  1\right)  }\right\Vert }}\times\ldots\times\frac{1+\left\Vert
z^{\left(  k-1\right)  }\right\Vert }{2^{\left\Vert z^{\left(  k-1\right)
}\right\Vert }}\frac{1+\left\Vert z^{\left(  k\right)  }\right\Vert
+\left\Vert \zeta^{\left(  1\right)  }\right\Vert }{2^{\left\Vert z^{\left(
k\right)  }\right\Vert }2^{\left\Vert \zeta^{\left(  1\right)  }\right\Vert }%
}\frac{1+\left\Vert \zeta^{\left(  2\right)  }\right\Vert }{2^{\left\Vert
\zeta^{\left(  2\right)  }\right\Vert }}\times\ldots\times\frac{1+\left\Vert
\zeta^{\left(  l\right)  }\right\Vert }{2^{\left\Vert \zeta^{\left(  l\right)
}\right\Vert }}\\
& =\frac{1+\left\Vert z^{\left(  k\right)  }\right\Vert +\left\Vert
\zeta^{\left(  1\right)  }\right\Vert }{\left(  1+\left\Vert z^{\left(
k\right)  }\right\Vert \right)  \left(  1+\left\Vert \zeta^{\left(  1\right)
}\right\Vert \right)  }\frac{1+\left\Vert z^{\left(  1\right)  }\right\Vert
}{2^{\left\Vert z^{\left(  1\right)  }\right\Vert }}\times\ldots\\
& \ldots\times\frac{1+\left\Vert z^{\left(  k-1\right)  }\right\Vert
}{2^{\left\Vert z^{\left(  k-1\right)  }\right\Vert }}\frac{1+\left\Vert
z^{\left(  k\right)  }\right\Vert }{2^{\left\Vert z^{\left(  k\right)
}\right\Vert }}\frac{1+\left\Vert \zeta^{\left(  1\right)  }\right\Vert
}{2^{\left\Vert \zeta^{\left(  1\right)  }\right\Vert }}\frac{1+\left\Vert
\zeta^{\left(  2\right)  }\right\Vert }{2^{\left\Vert \zeta^{\left(  2\right)
}\right\Vert }}\times\ldots\\
& \ldots\times\frac{1+\left\Vert \zeta^{\left(  l\right)  }\right\Vert
}{2^{\left\Vert \zeta^{\left(  l\right)  }\right\Vert }}\\
& =\frac{1+\left\Vert z^{\left(  k\right)  }\right\Vert +\left\Vert
\zeta^{\left(  1\right)  }\right\Vert }{\left(  1+\left\Vert z^{\left(
k\right)  }\right\Vert \right)  \left(  1+\left\Vert \zeta^{\left(  1\right)
}\right\Vert \right)  }c_{\beta}^{\left(  m\right)  }c_{\gamma}^{\left(
n\right)  }.
\end{align*}

Now $2\leq\left\Vert z^{\left(  k\right)  }\right\Vert \leq m-\left\vert
\beta\right\vert =\left\vert \mathbf{1}-\beta\right\vert $ and $2\leq
\left\Vert \zeta^{\left(  1\right)  }\right\Vert \leq n-\left\vert
\gamma\right\vert =\left\vert \mathbf{1}-\gamma\right\vert $ so
\begin{align*}
\frac{1}{1+m+n}<\frac{1+m+n}{\left(  1+m\right)  \left(  1+n\right)  }  &
\leq\frac{1+m-\left\vert \beta\right\vert +n-\left\vert \gamma\right\vert
}{\left(  1+m-\left\vert \beta\right\vert \right)  \left(  1+n-\left\vert
\gamma\right\vert \right)  }\\
& \leq\frac{1+\left\Vert z^{\left(  k\right)  }\right\Vert +\left\Vert
\zeta^{\left(  1\right)  }\right\Vert }{\left(  1+\left\Vert z^{\left(
k\right)  }\right\Vert \right)  \left(  1+\left\Vert \zeta^{\left(  1\right)
}\right\Vert \right)  }\\
& \leq\frac{5}{9},
\end{align*}

and%
\begin{align*}
\frac{1}{1+m+n}\leq\frac{1}{1+m-\left\vert \beta\right\vert +n-\left\vert
\gamma\right\vert }  & <\frac{1+m-\left\vert \beta\right\vert +n-\left\vert
\gamma\right\vert }{\left(  1+m-\left\vert \beta\right\vert \right)  \left(
1+n-\left\vert \gamma\right\vert \right)  }\\
& \leq\frac{1+\left\Vert z^{\left(  k\right)  }\right\Vert +\left\Vert
\zeta^{\left(  1\right)  }\right\Vert }{\left(  1+\left\Vert z^{\left(
k\right)  }\right\Vert \right)  \left(  1+\left\Vert \zeta^{\left(  1\right)
}\right\Vert \right)  }\\
& \leq\frac{5}{9},
\end{align*}

or%
\begin{align*}
\frac{1}{1+m+n}\leq\frac{1}{1+\left\vert \mathbf{1}-\left(  \beta
,\gamma\right)  \right\vert }  & <\frac{1+\left\vert \mathbf{1}-\left(
\beta,\gamma\right)  \right\vert }{\left(  1+\left\vert \mathbf{1}%
-\beta\right\vert \right)  \left(  1+\left\vert \mathbf{1}-\gamma\right\vert
\right)  }\\
& \leq\frac{1+\left\Vert z^{\left(  k\right)  }\right\Vert +\left\Vert
\zeta^{\left(  1\right)  }\right\Vert }{\left(  1+\left\Vert z^{\left(
k\right)  }\right\Vert \right)  \left(  1+\left\Vert \zeta^{\left(  1\right)
}\right\Vert \right)  }\\
& \leq\frac{5}{9}.
\end{align*}

Thus%
\begin{equation}
\frac{1+\left\vert \mathbf{1}-\left(  \beta,\gamma\right)  \right\vert
}{\left(  1+\left\vert \mathbf{1}-\beta\right\vert \right)  \left(
1+\left\vert \mathbf{1}-\gamma\right\vert \right)  }c_{\beta}^{\left(
m\right)  }c_{\gamma}^{\left(  n\right)  }\leq c_{\left(  \beta,\gamma\right)
}^{\left(  m+n\right)  }\leq\frac{5}{9}c_{\beta}^{\left(  m\right)  }%
c_{\gamma}^{\left(  n\right)  }.\label{av805}%
\end{equation}

These inequalities are sharp: equation \ref{av598} i.e. $c_{\mathbf{0}%
}^{\left(  m\right)  }=\frac{m+1}{2^{m}}$, implies that the left inequality of
\ref{av805} is an equality whenever $\beta=\left(  \mathbf{1},0_{p}\right)  $,
$2\leq p\leq m$ and $\gamma=\left(  \mathbf{1},0_{q}\right)  $, $2\leq q\leq n
$.
\end{proof}

\begin{remark}
Inequality \ref{av848} can be written%
\[
\frac{1}{1+m-\left\vert \beta\right\vert }c_{\beta}^{\left(  m\right)  }%
\frac{1}{1+n-\left\vert \gamma\right\vert }c_{\gamma}^{\left(  n\right)  }%
\leq\frac{1}{1+\left(  m+n\right)  -\left\vert \left(  \beta,\gamma\right)
\right\vert }c_{\left(  \beta,\gamma\right)  }^{\left(  m+n\right)  }%
\]

\end{remark}

Now we can prove:

\begin{theorem}
\label{vThm_sum_Cmb_UVb}\textbf{Separation of variables} Suppose
$u\in\mathbb{R}_{\oplus}^{k^{\prime}}$, $v\in\mathbb{R}_{\oplus}%
^{k^{\prime\prime}}$ and $p=k^{\prime}+k^{\prime\prime}$. Then%
\begin{equation}
\sum\limits_{\substack{\beta\leq\mathbf{1} \\\left\vert \beta\right\vert
=q}}c_{\beta}^{\left(  p\right)  }\left(  u,v\right)  ^{\beta}\geq
\sum\limits_{\substack{l= \\\max\left\{  q-k^{\prime\prime},0\right\}
}}^{\min\left\{  q,k^{\prime}\right\}  }\frac{1+p-q}{\left(  1+k^{\prime
}-l\right)  \left(  1+k^{\prime\prime}-q+l\right)  }\sum
\limits_{\substack{\beta^{\prime}\leq\mathbf{1} \\\left\vert \beta^{\prime
}\right\vert =l }}c_{\beta^{\prime}}^{\left(  k^{\prime}\right)  }%
u^{\beta^{\prime}}\sum\limits_{\substack{\beta^{\prime\prime}\leq\mathbf{1}
\\\left\vert \beta^{\prime\prime}\right\vert =q-l}}c_{\beta^{\prime\prime}%
}^{\left(  k^{\prime\prime}\right)  }v^{\beta^{\prime\prime}},\label{av030}%
\end{equation}

when $0\leq q\leq p$.
\end{theorem}

\begin{proof}
Applying inequality \ref{av842} of Theorem
\ref{vThm_estim_prod_C^m_and_b_C^n_c} to \ref{av863} below:
\begin{align}
\sum\limits_{\substack{\beta\leq\mathbf{1} \\\left\vert \beta\right\vert =q}}
&  c_{\beta}^{\left(  p\right)  }\left(  u,v\right)  ^{\beta}\nonumber\\
&  =\sum\limits_{\substack{\beta^{\prime}\leq1_{k^{\prime}} \\\beta
^{\prime\prime}\leq1_{k^{\prime\prime}} \\\left\vert \beta\right\vert
=q}}c_{\beta^{\prime},\beta^{\prime\prime}}^{\left(  p\right)  }%
u^{\beta^{\prime}}v^{\beta^{\prime\prime}}\label{av053}\\
&  =\sum\limits_{\substack{\beta^{\prime}\leq1_{k^{\prime}} \\\beta
^{\prime\prime}\leq1_{k^{\prime\prime}} \\\left\vert \beta^{\prime}\right\vert
+\left\vert \beta^{\prime\prime}\right\vert =q}}c_{\beta^{\prime}%
,\beta^{\prime\prime}}^{\left(  p\right)  }u^{\beta^{\prime}}v^{\beta
^{\prime\prime}}\nonumber\\
&  =\sum\limits_{\substack{\beta^{\prime}\leq1_{k^{\prime}} \\\beta
^{\prime\prime}\leq1_{k^{\prime\prime}} \\\left\vert \beta^{\prime\prime
}\right\vert =q-\left\vert \beta^{\prime}\right\vert }}c_{\beta^{\prime}%
,\beta^{\prime\prime}}^{\left(  p\right)  }u^{\beta^{\prime}}v^{\beta
^{\prime\prime}}\nonumber\\
&  =\sum\limits_{\substack{\beta^{\prime}\leq1_{k^{\prime}} \\0\leq\left\vert
\beta^{\prime}\right\vert \leq q}}u^{\beta^{\prime}}\sum
\limits_{\substack{\beta^{\prime\prime}\leq1_{k^{\prime\prime}} \\\left\vert
\beta^{\prime\prime}\right\vert =q-\left\vert \beta^{\prime}\right\vert
}}c_{\beta^{\prime},\beta^{\prime\prime}}^{\left(  p\right)  }v^{\beta
^{\prime\prime}}\nonumber\\
&  =\sum\limits_{l=0}^{q}\sum\limits_{\substack{\beta^{\prime}\leq
1_{k^{\prime}} \\\left\vert \beta^{\prime}\right\vert =l}}u^{\beta^{\prime}%
}\sum\limits_{\substack{\beta^{\prime\prime}\leq1_{k^{\prime\prime}}
\\\left\vert \beta^{\prime\prime}\right\vert =q-\left\vert \beta^{\prime
}\right\vert }}c_{\beta^{\prime},\beta^{\prime\prime}}^{\left(  p\right)
}v^{\beta^{\prime\prime}}\nonumber\\
&  =\sum\limits_{l=0}^{q}\sum\limits_{\substack{\beta^{\prime}\leq
1_{k^{\prime}} \\\left\vert \beta^{\prime}\right\vert =l}}u^{\beta^{\prime}%
}\sum\limits_{\substack{\beta^{\prime\prime}\leq1_{k^{\prime\prime}}
\\\left\vert \beta^{\prime\prime}\right\vert =q-l}}c_{\beta^{\prime}%
,\beta^{\prime\prime}}^{\left(  p\right)  }v^{\beta^{\prime\prime}}\nonumber\\
&  =\sum\limits_{l=0}^{\min\left\{  q,k^{\prime}\right\}  }\sum
\limits_{\substack{\beta^{\prime}\leq1_{k^{\prime}} \\\left\vert \beta
^{\prime}\right\vert =l}}u^{\beta^{\prime}}\sum\limits_{\substack{\beta
^{\prime\prime}\leq1_{k^{\prime\prime}} \\\left\vert \beta^{\prime\prime
}\right\vert =q-l}}c_{\beta^{\prime},\beta^{\prime\prime}}^{\left(  p\right)
}v^{\beta^{\prime\prime}}\label{av863}\\
&  \geq\sum\limits_{l=0}^{\min\left\{  q,k^{\prime}\right\}  }\sum
\limits_{\substack{\beta^{\prime}\leq1_{k^{\prime}} \\\left\vert \beta
^{\prime}\right\vert =l}}u^{\beta^{\prime}}\sum\limits_{\substack{\beta
^{\prime\prime}\leq1_{k^{\prime\prime}} \\\left\vert \beta^{\prime\prime
}\right\vert =q-l}}\frac{1+k^{\prime}-\left\vert \beta^{\prime}\right\vert
+k^{\prime\prime}-\left\vert \beta^{\prime\prime}\right\vert }{\left(
1+k^{\prime}-\left\vert \beta^{\prime}\right\vert \right)  \left(
1+k^{\prime\prime}-\left\vert \beta^{\prime\prime}\right\vert \right)
}c_{\beta^{\prime}}^{\left(  k^{\prime}\right)  }c_{\beta^{\prime\prime}%
}^{\left(  k^{\prime\prime}\right)  }v^{\beta^{\prime\prime}}\label{av866}\\
&  =\sum\limits_{l=0}^{\min\left\{  q,k^{\prime}\right\}  }\sum
\limits_{\substack{\beta^{\prime}\leq\mathbf{1} \\\left\vert \beta^{\prime
}\right\vert =l}}c_{\beta^{\prime}}^{\left(  k^{\prime}\right)  }%
u^{\beta^{\prime}}\sum\limits_{\substack{\beta^{\prime\prime}\leq\mathbf{1}
\\\left\vert \beta^{\prime\prime}\right\vert =q-l}}\frac{1+p-q}{\left(
1+k^{\prime}-l\right)  \left(  1+k^{\prime\prime}-q+l\right)  }c_{\beta
^{\prime\prime}}^{\left(  k^{\prime\prime}\right)  }v^{\beta^{\prime\prime}%
}\nonumber\\
&  =\sum\limits_{l=0}^{\min\left\{  q,k^{\prime}\right\}  }\frac
{1+p-q}{\left(  1+k^{\prime}-l\right)  \left(  1+k^{\prime\prime}-q+l\right)
}\sum\limits_{\substack{\beta^{\prime}\leq\mathbf{1} \\\left\vert
\beta^{\prime}\right\vert =l}}c_{\beta^{\prime}}^{\left(  k^{\prime}\right)
}u^{\beta^{\prime}}\sum\limits_{\substack{\beta^{\prime\prime}\leq\mathbf{1}
\\\left\vert \beta^{\prime\prime}\right\vert =q-l}}c_{\beta^{\prime\prime}%
}^{\left(  k^{\prime\prime}\right)  }v^{\beta^{\prime\prime}}.\label{av869}%
\end{align}

But we must have $q-l\leq k^{\prime\prime}$ i.e. $l\geq q-k^{\prime\prime}$.
Hence%
\[
\sum\limits_{\substack{\beta\leq\mathbf{1} \\\left\vert \beta\right\vert
=q}}c_{\beta}^{\left(  p\right)  }\left(  u,v\right)  ^{\beta}\geq
\sum\limits_{l=\max\left\{  q-k^{\prime\prime},0\right\}  }^{\min\left\{
q,k^{\prime}\right\}  }\frac{1+p-q}{\left(  1+k^{\prime}-l\right)  \left(
1+k^{\prime\prime}-q+l\right)  }\sum\limits_{\substack{\beta^{\prime}%
\leq\mathbf{1} \\\left\vert \beta^{\prime}\right\vert =l}}c_{\beta^{\prime}%
}^{\left(  k^{\prime}\right)  }u^{\beta^{\prime}}\sum\limits_{\substack{\beta
^{\prime\prime}\leq\mathbf{1} \\\left\vert \beta^{\prime\prime}\right\vert
=q-l}}c_{\beta^{\prime\prime}}^{\left(  k^{\prime\prime}\right)  }%
v^{\beta^{\prime\prime}},
\]

which is \ref{av030}.
\end{proof}

\begin{remark}
(not referenced) Inequality \ref{av030} can be written:%
\begin{align*}
\frac{1}{1+p-q} &  \sum\limits_{\substack{\beta\leq\mathbf{1} \\\left\vert
\beta\right\vert =q}}c_{\beta}^{\left(  p\right)  }\left(  u,v\right)
^{\beta}\\
&  \geq\sum\limits_{l=\max\left\{  q-k^{\prime\prime},0\right\}  }%
^{\min\left\{  q,k^{\prime}\right\}  }\frac{1}{\left(  1+k^{\prime}-l\right)
\left(  1+k^{\prime\prime}-(q-l\right)  }\sum\limits_{\substack{\beta^{\prime
}\leq\mathbf{1} \\\left\vert \beta^{\prime}\right\vert =l}}c_{\beta^{\prime}%
}^{\left(  k^{\prime}\right)  }u^{\beta^{\prime}}\sum\limits_{\substack{\beta
^{\prime\prime}\leq\mathbf{1} \\\left\vert \beta^{\prime\prime}\right\vert
=q-l}}c_{\beta^{\prime\prime}}^{\left(  k^{\prime\prime}\right)  }%
v^{\beta^{\prime\prime}}\\
&  =\sum\limits_{l=\max\left\{  q-k^{\prime\prime},0\right\}  }^{\min\left\{
q,k^{\prime}\right\}  }\frac{1}{1+k^{\prime}-l}\sum\limits_{\substack{\beta
^{\prime}\leq\geq\mathbf{1} \\\left\vert \beta^{\prime}\right\vert
=l}}c_{\beta^{\prime}}^{\left(  k^{\prime}\right)  }u^{\beta^{\prime}}\frac
{1}{1+k^{\prime\prime}-\left(  q-l\right)  }\sum\limits_{\substack{\beta
^{\prime\prime}\leq\mathbf{1} \\\left\vert \beta^{\prime\prime}\right\vert
=q-l}}c_{\beta^{\prime\prime}}^{\left(  k^{\prime\prime}\right)  }%
v^{\beta^{\prime\prime}}\\
&  =\sum\limits_{l=\max\left\{  q-k^{\prime\prime},0\right\}  }^{\min\left\{
q,k^{\prime}\right\}  }\frac{1}{1+k^{\prime}-l}\sum\limits_{\substack{\beta
^{\prime}\leq\mathbf{1} \\\left\vert \beta^{\prime}\right\vert =l}%
}c_{\beta^{\prime}}^{\left(  k^{\prime}\right)  }u^{\beta^{\prime}}\frac
{1}{1+k^{\prime\prime}-\left(  q-l\right)  }\sum\limits_{\substack{\beta
^{\prime\prime}\leq\mathbf{1} \\\left\vert \beta^{\prime\prime}\right\vert
=q-l}}c_{\beta^{\prime\prime}}^{\left(  k^{\prime\prime}\right)  }%
v^{\beta^{\prime\prime}},
\end{align*}

so that if%
\[
\kappa_{l}^{k}\left(  z\right)  :=\frac{1}{1+k-l}\sum\limits_{\substack{\beta
\leq\mathbf{1} \\\left\vert \beta\right\vert =l}}c_{\beta}^{\left(  k\right)
}z^{\beta},\quad z\in\mathbb{R}^{k},\text{ }0\leq l\leq k,
\]

the estimate \ref{av030} can be also written%
\[
\kappa_{q}^{p}\left(  u,v\right)  \geq\sum\limits_{l=\max\left\{
q-k^{\prime\prime},0\right\}  }^{\min\left\{  q,k^{\prime}\right\}  }%
\kappa_{l}^{k^{\prime}}\left(  u\right)  \kappa_{q-l}^{k^{\prime\prime}%
}\left(  v\right)  ,\quad0\leq q\leq p,\text{ }u\in\mathbb{R}_{\oplus
}^{k^{\prime}},\text{ }v\in\mathbb{R}_{\oplus}^{k^{\prime\prime}}.
\]

\end{remark}

When now derive our upper and lower bounds which remove leading and trailing
zeros from the arguments of the determinants:

\begin{theorem}
\label{vThm_low_bnd_remove_zeros}(\ref{av935}) For $n\geq1$,%
\begin{equation}
\left\vert T\left(  0_{n}\right)  \right\vert =c_{\mathbf{0}}^{\left(
n\right)  }=\frac{1+n}{2^{n}},\label{av877}%
\end{equation}

and for $p\geq1$:%
\begin{equation}
\left\vert T\left(  0_{n},x\right)  \right\vert \geq\frac{1+n+p}{1+p}\frac
{1}{2^{n}}\left\vert T\left(  x\right)  \right\vert ,\quad x\in\mathbb{R}%
_{\oplus}^{p},\label{av876}%
\end{equation}

which is exact when $x=0_{p}$. Also%
\begin{equation}
\left\vert T\left(  x,0_{n}\right)  \right\vert \geq\frac{1+n+p}{1+p}\frac
{1}{2^{n}}\left\vert T\left(  x\right)  \right\vert ,\quad x\in\mathbb{R}%
_{\oplus}^{p},\label{av878}%
\end{equation}

which is exact when $x=0_{p},$ Finally%
\begin{equation}
\left\vert T\left(  0_{n},x\right)  \right\vert <\left\vert T\left(  x\right)
\right\vert ,\quad x\in\mathbb{R}_{\oplus}^{p},\label{av915}%
\end{equation}

and%
\begin{equation}
\left\vert T\left(  x,0_{n}\right)  \right\vert <\left\vert T\left(  x\right)
\right\vert ,\quad x\in\mathbb{R}_{\oplus}^{p}.\label{av885}%
\end{equation}

\end{theorem}

\begin{proof}
From \ref{av599},%
\[
T\left(  x\right)  =\sum_{k=0}^{m}\sum\limits_{\substack{\beta\leq\mathbf{1}
\\\left\vert \beta\right\vert =k}}c_{\beta}^{\left(  m\right)  }x^{\beta
}=c_{\mathbf{0}}^{\left(  m\right)  }+\sum_{k=1}^{m}\sum
\limits_{\substack{\beta\leq\mathbf{1} \\\left\vert \beta\right\vert
=k}}c_{\beta}^{\left(  m\right)  }x^{\beta},\quad x\in\mathbb{R}^{m},
\]

so that $T\left(  0_{n}\right)  =c_{\mathbf{0}}^{\left(  n\right)  }$.
Further, from Theorem \ref{vThm_estim_prod_C^m_and_b_C^n_c},%
\[
c_{\left(  \beta^{\prime},\beta^{\prime\prime}\right)  }^{\left(  n+p\right)
}\geq\frac{1+n-\left\vert \beta^{\prime}\right\vert +p-\left\vert
\beta^{\prime\prime}\right\vert }{\left(  1+n-\left\vert \beta^{\prime
}\right\vert \right)  \left(  1+p-\left\vert \beta^{\prime\prime}\right\vert
\right)  }c_{\beta^{\prime}}^{\left(  n\right)  }c_{\beta^{\prime\prime}%
}^{\left(  p\right)  },
\]

so%
\begin{align*}
\left\vert T\left(  0_{n},x\right)  \right\vert  & =\sum_{k=0}^{n+p}%
\sum\limits_{\substack{\beta^{\prime}\leq1_{n}\mathbf{,}\beta^{\prime\prime
}\leq1_{p} \\\left\vert \beta^{\prime}\right\vert +\left\vert \beta
^{\prime\prime}\right\vert =k}}c_{\beta^{\prime},\beta^{\prime\prime}%
}^{\left(  n+p\right)  }\left(  0_{n},x\right)  ^{\left(  \beta^{\prime}%
,\beta^{\prime\prime}\right)  }\\
& =\sum_{k=0}^{n+p}\sum\limits_{\substack{\beta^{\prime\prime}\leq1_{p}
\\\left\vert \beta^{\prime\prime}\right\vert =k}}c_{0_{n},\beta^{\prime\prime
}}^{\left(  n+p\right)  }x^{\beta^{\prime\prime}}\\
& \geq\sum_{k=0}^{n+p}\sum\limits_{\substack{\beta^{\prime\prime}%
\leq\mathbf{1}_{p} \\\left\vert \beta^{\prime\prime}\right\vert =k}%
}\frac{1+n+p-\left\vert \beta^{\prime\prime}\right\vert }{\left(  1+n\right)
\left(  1+p-\left\vert \beta^{\prime\prime}\right\vert \right)  }%
c_{\mathbf{0}}^{\left(  n\right)  }c_{\beta^{\prime\prime}}^{\left(  p\right)
}x^{\beta^{\prime\prime}}\\
& =c_{\mathbf{0}}^{\left(  n\right)  }\sum_{k=0}^{n+p}\sum
\limits_{\substack{\beta^{\prime\prime}\leq\mathbf{1} \\\left\vert
\beta^{\prime\prime}\right\vert =k}}\frac{1+n+p-k}{\left(  1+n\right)  \left(
1+p-k\right)  }c_{\beta^{\prime\prime}}^{\left(  p\right)  }x^{\beta
^{\prime\prime}}\\
& =c_{\mathbf{0}}^{\left(  n\right)  }\sum_{k=0}^{n+p}\frac{1+n+p-k}{\left(
1+n\right)  \left(  1+p-k\right)  }\sum\limits_{\substack{\beta^{\prime\prime
}\leq\mathbf{1} \\\left\vert \beta^{\prime\prime}\right\vert =k }%
}c_{\beta^{\prime\prime}}^{\left(  p\right)  }x^{\beta^{\prime\prime}}\\
& =\frac{c_{\mathbf{0}}^{\left(  n\right)  }}{1+n}\sum_{k=0}^{p}\frac
{1+n+p-k}{1+p-k}\sum\limits_{\substack{\beta^{\prime\prime}\leq\mathbf{1}
\\\left\vert \beta^{\prime\prime}\right\vert =k}}c_{\beta^{\prime\prime}%
}^{\left(  p\right)  }x^{\beta^{\prime\prime}}\\
& =\frac{c_{\mathbf{0}}^{\left(  n\right)  }}{1+n}\sum_{k=0}^{p}\left(
1+\frac{n}{1+p-k}\right)  \sum\limits_{\substack{\beta^{\prime\prime}%
\leq\mathbf{1} \\\left\vert \beta^{\prime\prime}\right\vert =k}}c_{\beta
^{\prime\prime}}^{\left(  p\right)  }x^{\beta^{\prime\prime}}\\
& \geq\frac{c_{\mathbf{0}}^{\left(  n\right)  }}{1+n}\sum_{k=0}^{p}\left(
1+\frac{n}{1+p}\right)  \sum\limits_{\substack{\beta^{\prime\prime}%
\leq\mathbf{1} \\\left\vert \beta^{\prime\prime}\right\vert =k}}c_{\beta
^{\prime\prime}}^{\left(  p\right)  }x^{\beta^{\prime\prime}}=\\
& =\frac{1+n+p}{\left(  1+n\right)  \left(  1+p\right)  }c_{\mathbf{0}%
}^{\left(  n\right)  }\sum_{k=0}^{p}\sum\limits_{\substack{\beta^{\prime
\prime}\leq\mathbf{1} \\\left\vert \beta^{\prime\prime}\right\vert
=k}}c_{\beta^{\prime\prime}}^{\left(  p\right)  }x^{\beta^{\prime\prime}}\\
& =\frac{1+n+p}{\left(  1+n\right)  \left(  1+p\right)  }c_{\mathbf{0}%
}^{\left(  n\right)  }\left\vert T\left(  x\right)  \right\vert \\
& =\frac{1+n+p}{\left(  1+n\right)  \left(  1+p\right)  }\frac{1+n}{2^{n}%
}\left\vert T\left(  x\right)  \right\vert \\
& =\frac{1+n+p}{1+p}\frac{1}{2^{n}}\left\vert T\left(  x\right)  \right\vert .
\end{align*}

Regarding \ref{av878}, from Corollary \ref{vCor_1_formula_for_C^(m)_gam} the
reverse order permutation $\mathcal{R}$ satisfies $\left\vert T\left(
x\right)  \right\vert =\left\vert T\left(  \mathcal{R}x\right)  \right\vert $
for all $x$ so%
\begin{align*}
\left\vert T\left(  x,0_{n}\right)  \right\vert =\left\vert T\left(
\mathcal{R}\left(  x,0_{n}\right)  \right)  \right\vert =\left\vert T\left(
0_{n},\mathcal{R}x\right)  \right\vert  & \geq\frac{1+n+p}{1+p}\frac{1}{2^{n}%
}\left\vert T\left(  \mathcal{R}x\right)  \right\vert \\
& =\frac{1+n+p}{1+p}\frac{1}{2^{n}}\left\vert T\left(  x\right)  \right\vert .
\end{align*}

Regarding \ref{av915} and \ref{av885}: from \ref{av540}, $\left\vert T\left(
a\right)  \right\vert =\left(  1+a_{1}\right)  \left\vert T\left(
a_{2:m}\right)  \right\vert -\frac{1}{4}\left\vert T\left(  a_{3,m}\right)
\right\vert $ so
\[
\left\vert T\left(  0_{n},x\right)  \right\vert <\left\vert T\left(
0_{n-1},x\right)  \right\vert <\ldots<\left\vert T\left(  0,x\right)
\right\vert .
\]
If $\dim x>1$ then $\left\vert T\left(  0_{n},x\right)  \right\vert
<\left\vert T\left(  x\right)  \right\vert $ and if $\dim x=1$ then
$\left\vert T\left(  0,x\right)  \right\vert =%
\begin{vmatrix}
1 & -1/2\\
-1/2 & 1+x
\end{vmatrix}
=\frac{3}{4}+x<1+x=\left\vert T\left(  x\right)  \right\vert $. Further,
\[
\left\vert T\left(  x,0_{n}\right)  \right\vert =\left\vert T\left(
\mathcal{R}\left(  x,0_{n}\right)  \right)  \right\vert =\left\vert T\left(
0_{n},\mathcal{R}x\right)  \right\vert <\left\vert T\left(  \mathcal{R}%
x\right)  \right\vert =\left\vert T\left(  x\right)  \right\vert .
\]

\end{proof}

\begin{remark}
The upper bound \ref{av915} will be improved on below in \ref{av934}.
\end{remark}

\subsection{Upper bounds for $\frac{1}{v}\left\vert T^{-1}\Delta
_{X}f\right\vert _{\max}$ and $\left\vert \mu\left(  f\right)  \right\vert
_{\max}$\label{SbSecT_max_invT_deltaXFn}}

The next theorem will be derived by using the lower bounds \ref{av876} and
\ref{av878} to eliminate the zeros from the equations \ref{av039} and
\ref{av037} which have terms with factors such as $\left\vert T\left(
0_{k-j},x_{k+1:m}\right)  \right\vert $ and $\left\vert T\left(
x_{1:k-2},0_{l-k+1}\right)  \right\vert $.

\begin{theorem}
\label{vThm_low_bnd_detT}If $m\geq4$ and $x\in\mathbb{R}_{\oplus}^{m}$ then%
\begin{equation}%
\begin{array}
[c]{ll}%
\left\vert T\left(  x\right)  \right\vert \geq\left(  x_{1}+\frac{m+1}{m}%
\frac{1}{2}\right)  \left\vert T\left(  x_{2:m}\right)  \right\vert , & k=1,\\
\left\vert T\left(  x\right)  \right\vert \geq\left(
\begin{array}
[c]{r}%
\frac{m}{m-k+1}\frac{x_{1}}{2^{k-1}}+\sum\limits_{n=2}^{k}\frac{m-n+1}%
{m-k+1}\frac{x_{n}}{2^{k-n}}\left\vert T\left(  x_{1:n-1}\right)  \right\vert
+\\
+\frac{m+1}{m-k+1}\frac{1}{2^{k}}%
\end{array}
\right)  \left\vert T\left(  x_{k+1:m}\right)  \right\vert ,\quad & 2\leq
k\leq m-1,\\
\left\vert T\left(  x\right)  \right\vert =m\frac{x_{1}}{2^{m-1}}%
+\sum\limits_{n=2}^{m}\left(  m-n+1\right)  \left\vert T\left(  x_{1:n-1}%
\right)  \right\vert \frac{x_{n}}{2^{m-n}}+\frac{m+1}{2^{m}}, & k=m,
\end{array}
\label{av879}%
\end{equation}

and,%
\begin{equation}%
\begin{array}
[c]{ll}%
\left\vert T\left(  x\right)  \right\vert =\sum\limits_{n=1}^{m-1}n\frac
{x_{n}}{2^{n-1}}\left\vert T\left(  x_{n+1:m}\right)  \right\vert
+m\frac{x_{m}}{2^{m-1}}+\frac{m+1}{2^{m}}, & k=1,\\
\left\vert T\left(  x\right)  \right\vert \geq\left\vert T\left(
x_{1:k-1}\right)  \right\vert \left(  \sum\limits_{n=k}^{m-1}\frac{n}{k}%
\frac{x_{n}}{2^{n-k}}\left\vert T\left(  x_{n+1:m}\right)  \right\vert
+\frac{m}{k}\frac{x_{m}}{2^{m-k}}+\frac{m+1}{k}\frac{1}{2^{m-k+1}}\right)
,\quad & 2\leq k\leq m-1,\\
\left\vert T\left(  x\right)  \right\vert \geq\left\vert T\left(
x_{1:m-1}\right)  \right\vert \left(  x_{m}+\frac{m+1}{m}\frac{1}{2}\right)
, & k=m,
\end{array}
\label{av880}%
\end{equation}

\end{theorem}

\begin{proof}
We first prove the bounds \ref{av879} by applying the lower bound
\ref{av876}:
\[
\left\vert T\left(  0_{p},x_{l:m}\right)  \right\vert \geq\frac{2+p+m-l}%
{2+m-l}\frac{1}{2^{p}}\left\vert T\left(  x\right)  \right\vert ,
\]

to equation \ref{av039}.\medskip

\fbox{\textbf{Case} A1} When $k=1$ in \ref{av039}:%
\begin{align}
\left\vert T\left(  x\right)  \right\vert =x_{1}\left\vert T\left(
x_{2:m}\right)  \right\vert +\left\vert T\left(  0,x_{2:m}\right)
\right\vert  & \geq x_{1}\left\vert T\left(  x_{2:m}\right)  \right\vert
+\frac{1+m}{1+m-1}\frac{1}{2}\left\vert T\left(  x_{2:m}\right)  \right\vert
\nonumber\\
& =x_{1}\left\vert T\left(  x_{2:m}\right)  \right\vert +\frac{1+m}{m}\frac
{1}{2}\left\vert T\left(  x_{2:m}\right)  \right\vert \nonumber\\
& =\left(  x_{1}+\frac{m+1}{m}\frac{1}{2}\right)  \left\vert T\left(
x_{2:m}\right)  \right\vert .\label{av038}%
\end{align}

\fbox{\textbf{Case} B1} When $k=2$ in \ref{av039}:%
\begin{align*}
\left\vert T\left(  x\right)  \right\vert  & =x_{1}\left\vert T\left(
0,x_{3:m}\right)  \right\vert +\left\vert T\left(  x_{1}\right)  \right\vert
x_{2}\left\vert T\left(  x_{3:m}\right)  \right\vert +\left\vert T\left(
0_{2},x_{3:m}\right)  \right\vert \\
& \geq x_{1}\frac{2+m-2}{1+m-2}\frac{1}{2}\left\vert T\left(  x_{3:m}\right)
\right\vert +\left\vert T\left(  x_{1}\right)  \right\vert x_{2}\left\vert
T\left(  x_{3:m}\right)  \right\vert +\frac{1+m-2+1}{1+m-2}\frac{1}{2^{p}%
}\left\vert T\left(  x_{3:m}\right)  \right\vert \\
& =\left(  \frac{m}{m-1}\frac{1}{2}x_{1}+\left\vert T\left(  x_{1}\right)
\right\vert x_{2}+\frac{m}{m-1}\frac{1}{2^{2}}\right)  \left\vert T\left(
x_{3:m}\right)  \right\vert .
\end{align*}
\medskip

\fbox{\textbf{Case} C1} When $3\leq k\leq m-1$ in \ref{av039}:%
\begin{align*}
\left\vert T\left(  x\right)  \right\vert  & =x_{1}\left\vert T\left(
0_{k-1},x_{k+1:m}\right)  \right\vert +\sum\limits_{n=2}^{k-1}\left\vert
T\left(  x_{1:n-1}\right)  \right\vert x_{n}\left\vert T\left(  0_{k-n}%
,x_{k+1:m}\right)  \right\vert +\\
& \qquad\qquad+\left\vert T\left(  x_{1:k-1}\right)  \right\vert
x_{k}\left\vert T\left(  x_{k+1:m}\right)  \right\vert +\left\vert T\left(
0_{k},x_{k+1:m}\right)  \right\vert \\
& \geq\frac{m}{m-k+1}\frac{1}{2^{k-1}}x_{1}\left\vert T\left(  x_{k+1:m}%
\right)  \right\vert +\sum\limits_{n=2}^{k-1}\frac{m-n+1}{m-k+1}\frac
{1}{2^{k-n}}\left\vert T\left(  x_{1:n-1}\right)  \right\vert x_{n}\left\vert
T\left(  x_{k+1:m}\right)  \right\vert +\\
& \qquad\qquad+\left\vert T\left(  x_{1:k-1}\right)  \right\vert
x_{k}\left\vert T\left(  x_{k+1:m}\right)  \right\vert +\frac{m+1}{m-k+1}%
\frac{1}{2^{k}}\left\vert T\left(  x_{k+1:m}\right)  \right\vert \\
& =\left(  \sum\limits_{n=1}^{k-1}\frac{m-n+1}{m-k+1}\frac{1}{2^{k-n}%
}\left\vert T\left(  x_{1:n-1}\right)  \right\vert x_{n}+\left\vert T\left(
x_{1:k-1}\right)  \right\vert x_{k}+\frac{m+1}{m-k+1}\frac{1}{2^{k}}\right)
\left\vert T\left(  x_{k+1:m}\right)  \right\vert \\
& =\left(  \sum\limits_{n=1}^{k}\frac{m-n+1}{m-k+1}\frac{1}{2^{k-n}}\left\vert
T\left(  x_{1:n-1}\right)  \right\vert x_{n}+\frac{m+1}{m-k+1}\frac{1}{2^{k}%
}\right)  \left\vert T\left(  x_{k+1:m}\right)  \right\vert .
\end{align*}

\fbox{\textbf{Case} D1} We will use the equation of case $k=m$ of
\ref{av039}.\medskip

This proves \ref{av879}.

Next to prove \ref{av880} by applying the following consequence of
\ref{av878}:%
\begin{equation}
\left\vert T\left(  x_{1:p},0_{q}\right)  \right\vert \geq\frac{1+p+q}%
{1+p}\frac{1}{2^{q}}\left\vert T\left(  x_{1:p}\right)  \right\vert ,\quad
x_{1:p}>0,\label{av040}%
\end{equation}

to the equations \ref{av037}.\medskip

\fbox{\textbf{Case} A2} Starting with the case $k=1$ of \ref{av037}:%
\begin{align}
\left\vert T\left(  x\right)  \right\vert  & =\left\vert T\left(
x_{1:m-1}\right)  \right\vert x_{m}+\left\vert T\left(  x_{1:m-1},0\right)
\right\vert \nonumber\\
& \geq\left\vert T\left(  x_{1:m-1}\right)  \right\vert x_{m}+\frac{m+1}%
{m}\frac{1}{2}\left\vert T\left(  x_{1:m-1}\right)  \right\vert \nonumber\\
& =\left\vert T\left(  x_{1:m-1}\right)  \right\vert \left(  x_{m}+\frac
{m+1}{m}\frac{1}{2}\right)  .\label{av041}%
\end{align}

\fbox{\textbf{Case} B2} Starting with the case $k=2$ of \ref{av037}:%
\begin{align}
\left\vert T\left(  x\right)  \right\vert  & =\left\vert T\left(
x_{1:m-2}\right)  \right\vert x_{m-1}\left\vert T\left(  x_{m}\right)
\right\vert +\left\vert T\left(  x_{1:m-2},0\right)  \right\vert
x_{m}+\left\vert T\left(  x_{1:m-2},0_{2}\right)  \right\vert \nonumber\\
& \geq\left\vert T\left(  x_{1:m-2}\right)  \right\vert x_{m-1}\left\vert
T\left(  x_{m}\right)  \right\vert +\frac{m}{m-1}\frac{1}{2}\left\vert
T\left(  x_{1:m-2}\right)  \right\vert x_{m}+\frac{m+1}{m-1}\frac{1}{2^{2}%
}\left\vert T\left(  x_{1:m-2}\right)  \right\vert \nonumber\\
& =\left\vert T\left(  x_{1:m-2}\right)  \right\vert \left(  x_{m-1}\left\vert
T\left(  x_{m}\right)  \right\vert +\frac{m}{m-1}\frac{x_{m}}{2}+\frac
{m+1}{m-1}\frac{1}{2^{2}}\right)  .\label{av042}%
\end{align}

\fbox{\textbf{Case} C2} Starting with the case $3\leq k\leq m-1$ of
\ref{av037}:%
\begin{align*}
\left\vert T\left(  x\right)  \right\vert  & =\left\vert T\left(
x_{1:m-k}\right)  \right\vert x_{m-k+1}\left\vert T\left(  x_{m-k+2:m}\right)
\right\vert +\\
& \qquad+\sum\limits_{n=m-k+2}^{m-1}\left\vert T\left(  x_{1:m-k}%
,0_{n+k-m-1}\right)  \right\vert x_{n}\left\vert T\left(  x_{n+1:m}\right)
\right\vert +\\
& \qquad+\left\vert T\left(  x_{1:m-k},0_{k-1}\right)  \right\vert
x_{m}+\left\vert T\left(  x_{1:m-k},0_{k}\right)  \right\vert \\
& \geq\left\vert T\left(  x_{1:m-k}\right)  \right\vert x_{m-k+1}\left\vert
T\left(  x_{m-k+2:m}\right)  \right\vert +\\
& \qquad+\sum\limits_{n=m-k+2}^{m-1}\frac{n}{m-k+1}\frac{1}{2^{n+k-m-1}%
}\left\vert T\left(  x_{1:m-k}\right)  \right\vert x_{n}\left\vert T\left(
x_{n+1:m}\right)  \right\vert +\\
& \qquad+\frac{m}{m-k+1}\frac{1}{2^{k-1}}\left\vert T\left(  x_{1:m-k}\right)
\right\vert x_{m}+\frac{m+1}{m-k+1}\frac{1}{2^{k}}\left\vert T\left(
x_{1:m-k}\right)  \right\vert \\
& =\left\vert T\left(  x_{1:m-k}\right)  \right\vert \left(
\begin{array}
[c]{r}%
x_{m-k+1}\left\vert T\left(  x_{m-k+2:m}\right)  \right\vert +\sum
\limits_{n=m-k+2}^{m-1}\frac{n}{m-k+1}\frac{x_{n}}{2^{n+k-m-1}}\left\vert
T\left(  x_{n+1:m}\right)  \right\vert +\\
+\frac{m}{m-k+1}\frac{x_{m}}{2^{k-1}}+\frac{m+1}{m-k+1}\frac{1}{2^{k}}%
\end{array}
\right) \\
& =\left\vert T\left(  x_{1:m-k}\right)  \right\vert \left(
\begin{array}
[c]{r}%
\sum\limits_{n=m-k+1}^{m-1}\frac{n}{m-k+1}\frac{x_{n}}{2^{n+k-m-1}}\left\vert
T\left(  x_{n+1:m}\right)  \right\vert +\frac{m}{m-k+1}\frac{x_{m}}{2^{k-1}%
}+\\
+\frac{m+1}{m-k+1}\frac{1}{2^{k}}%
\end{array}
\right)  .
\end{align*}

Next apply the transformation $k\rightarrow m-k+1$ to give%
\begin{equation}
\left\vert T\left(  x\right)  \right\vert =\left\vert T\left(  x_{1:k-1}%
\right)  \right\vert \left(  \sum\limits_{n=k}^{m-1}\frac{n}{k}\frac{x_{n}%
}{2^{n-k}}\left\vert T\left(  x_{n+1:m}\right)  \right\vert +\frac{m}{k}%
\frac{x_{m}}{2^{m-k}}+\frac{m+1}{k}\frac{1}{2^{m-k+1}}\right)  ,\quad2\leq
k\leq m-2.\label{av043}%
\end{equation}

\fbox{\textbf{Case} D2} We will use the equation of case $k=m$ of \ref{av037}:%
\begin{equation}
\left\vert T\left(  x\right)  \right\vert =\sum\limits_{n=1}^{m-1}\frac
{n}{2^{n-1}}x_{n}\left\vert T\left(  x_{n+1:m}\right)  \right\vert +\frac
{m}{2^{m-1}}x_{m}+\frac{m+1}{2^{m}}.\label{av044}%
\end{equation}

Now if we substitute $k=1$ into \ref{av043} and apply the notation of
Definition \ref{vDef_T_minors} we obtain \ref{av044}. If we substitute $k=m-1
$ into \ref{av043} and apply the notation of Definition \ref{vDef_T_minors} we
obtain \ref{av042}. If we substitute $k=m$ into \ref{av043}, assume
$\sum\limits_{n=m}^{m-1}\ldots:=0$ and then apply the notation of Definition
\ref{vDef_T_minors} we obtain \ref{av041}.
\end{proof}

We now apply Theorem \ref{vThm_low_bnd_detT} to obtain:

\begin{corollary}
\label{vCor_1_Thm_low_bnd_detT}Suppose $x\in\mathbb{R}_{\oplus}^{m}$ and
$m\geq4$. Then:%
\begin{equation}%
\begin{array}
[c]{ll}%
\left\vert T\left(  x\right)  \right\vert =\sum\limits_{n=1}^{m-1}\frac{x_{n}%
}{2^{n-1}}\left\vert T\left(  x_{n+1:m}\right)  \right\vert +\frac{x_{m}%
}{2^{m-1}}+\frac{1}{2^{m}}, & k=1,\\
\left\vert T\left(  x\right)  \right\vert \geq\left(
\begin{array}
[c]{r}%
\frac{1}{2}\left\vert T\left(  x_{1:k-1}\right)  \right\vert \left(
\sum\limits_{n=k}^{m-1}\frac{x_{n}}{2^{n-k}}\left\vert T\left(  x_{n+1:m}%
\right)  \right\vert +\frac{x_{m}}{2^{m-k}}+\frac{1}{2^{m-k+1}}\right)  +\\
+\frac{1}{2}\left(  \frac{x_{1}}{2^{k-1}}+\sum\limits_{n=2}^{k}\left\vert
T\left(  x_{1:n-1}\right)  \right\vert \frac{x_{n}}{2^{k-n}}+\frac{1}{2^{k}%
}\right)  \left\vert T\left(  x_{k+1:m}\right)  \right\vert
\end{array}
\right)  ,\quad & 2\leq k\leq m-1,\\
\left\vert T\left(  x\right)  \right\vert =\frac{x_{1}}{2^{m-1}}%
+\sum\limits_{n=2}^{m}\left\vert T\left(  x_{1:n-1}\right)  \right\vert
\frac{x_{n}}{2^{m-n}}+\frac{1}{2^{m}}, & k=m.
\end{array}
\label{av047}%
\end{equation}

\end{corollary}

\begin{proof}
Suppose $m\geq4$. Then from \ref{av879},%
\begin{equation}%
\begin{array}
[c]{lll}%
\left\vert T\left(  x\right)  \right\vert \geq & \left(  x_{1}+\frac{1}%
{2}\right)  \left\vert T\left(  x_{2:m}\right)  \right\vert , & k=1,\\
\left\vert T\left(  x\right)  \right\vert \geq & \left(  \frac{x_{1}}{2^{k-1}%
}+\sum\limits_{n=2}^{k}\left\vert T\left(  x_{1:n-1}\right)  \right\vert
\frac{x_{n}}{2^{k-n}}+\frac{1}{2^{k}}\right)  \left\vert T\left(
x_{k+1:m}\right)  \right\vert ,\quad & 2\leq k\leq m-1,\\
\left\vert T\left(  x\right)  \right\vert = & \frac{x_{1}}{2^{m-1}}%
+\sum\limits_{n=2}^{m}\left\vert T\left(  x_{1:n-1}\right)  \right\vert
\frac{x_{n}}{2^{m-n}}+\frac{1}{2^{m}}, & k=m,
\end{array}
\label{av045}%
\end{equation}

and from \ref{av880},%
\begin{equation}%
\begin{array}
[c]{lll}%
\left\vert T\left(  x\right)  \right\vert = & \sum\limits_{n=1}^{m-1}%
\frac{x_{n}}{2^{n-1}}\left\vert T\left(  x_{n+1:m}\right)  \right\vert
+\frac{x_{m}}{2^{m-1}}+\frac{1}{2^{m}}, & k=1,\\
\left\vert T\left(  x\right)  \right\vert \geq & \left\vert T\left(
x_{1:k-1}\right)  \right\vert \left(  \sum\limits_{n=k}^{m-1}\frac{x_{n}%
}{2^{n-k}}\left\vert T\left(  x_{n+1:m}\right)  \right\vert +\frac{x_{m}%
}{2^{m-k}}+\frac{1}{2^{m-k+1}}\right)  ,\quad & 2\leq k\leq m-1,\\
\left\vert T\left(  x\right)  \right\vert \geq & \left\vert T\left(
x_{1:m-1}\right)  \right\vert \left(  x_{m}+\frac{1}{2}\right)  , & k=m,
\end{array}
\label{av046}%
\end{equation}

Adding terms with the same value of $k:2\leq k\leq m-1$ yields this corollary.
\end{proof}

\begin{corollary}
\label{vCor_2_Thm_low_bnd_detT}Recall that $T=T\left(  x;1\right)  $ where
$x=\frac{1}{v}\Delta_{1}X$ and $v=\rho N$. Then if $f\in C_{B}^{\left(
0\right)  }\left(  \Omega\right)  $ and $Df\in L^{\infty}\left(
\Omega\right)  $:%
\begin{equation}
\frac{1}{v}\left\vert T^{-1}\Delta_{X}f\right\vert _{\max}\leq2\left\Vert
Df\right\Vert _{\infty;\Omega},\label{av048}%
\end{equation}

and%
\begin{equation}
\left\vert \mu\left(  f\right)  \right\vert _{\max}\leq2\left\Vert
Df\right\Vert _{\infty;\Omega}+\left\vert \mu_{1}\left(  f\right)  \right\vert
,\label{av081}%
\end{equation}

when $N\geq4$ and $\rho>0$.
\end{corollary}

\begin{proof}
From \ref{av029},%
\begin{align*}
&  \frac{1}{v}\left\vert \left(  T^{-1}\Delta_{X}f\right)  _{k}\right\vert \\
&  \leq\left\Vert Df\right\Vert _{\infty;\Omega}\times\left\{
\begin{array}
[c]{ll}%
\frac{1}{\left\vert T\left(  x\right)  \right\vert }\left(  \sum
\limits_{n=1}^{N-2}\frac{x_{n}}{2^{n-1}}\left\vert T\left(  x_{n+1:N-1}%
\right)  \right\vert +\frac{x_{N-1}}{2^{N-2}}\right)  , & k=1,\\
& \\
\left.
\begin{array}
[c]{l}%
\frac{\left\vert T\left(  x_{k+1:N-1}\right)  \right\vert }{\left\vert
T\left(  x\right)  \right\vert }\left(  \frac{x_{1}}{2^{k-1}}+\sum
\limits_{n=2}^{k}\left\vert T\left(  x_{1:n-1}\right)  \right\vert \frac
{x_{n}}{2^{k-n}}\right)  +\\
+\frac{\left\vert T\left(  x_{1:k-1}\right)  \right\vert }{\left\vert T\left(
x\right)  \right\vert }\left(  \sum\limits_{n=k+1}^{N-2}\frac{x_{n}}{2^{n-k}%
}\left\vert T\left(  x_{n+1:N-1}\right)  \right\vert +\frac{x_{N-1}}%
{2^{N-1-k}}\right)
\end{array}
\right\}  , & 2\leq k\leq N-2,\\
& \\
\frac{1}{\left\vert T\left(  x\right)  \right\vert }\left(  \frac{x_{1}%
}{2^{N-2}}+\sum\limits_{n=2}^{N-1}\left\vert T\left(  x_{1:n-1}\right)
\right\vert \frac{x_{n}}{2^{N-1-n}}\right)  , & k=N-1,
\end{array}
\right.
\end{align*}

Set $m=N-1$. The case $k=1$ of \ref{av047} now implies%
\[
\sum\limits_{n=1}^{N-2}\frac{x_{n}}{2^{n-1}}\left\vert T\left(  x_{n+1:N-1}%
\right)  \right\vert +\frac{x_{N-1}}{2^{N-2}}=\sum\limits_{n=1}^{m-1}%
\frac{x_{n}}{2^{n-1}}\left\vert T\left(  x_{n+1:m}\right)  \right\vert
+\frac{x_{m}}{2^{m-1}}<\left\vert T\left(  x\right)  \right\vert ,
\]

so that%
\[
\frac{1}{v}\left\vert \left(  T^{-1}\Delta_{X}f\right)  _{1}\right\vert
\leq\frac{\left\Vert Df\right\Vert _{\infty;\Omega}}{\left\vert T\left(
x\right)  \right\vert }\left(  \sum\limits_{n=1}^{N-2}\frac{x_{n}}{2^{n-1}%
}\left\vert T\left(  x_{n+1:N-1}\right)  \right\vert +\frac{x_{N-1}}{2^{N-2}%
}\right)  <\left\Vert Df\right\Vert _{\infty;\Omega}.
\]

Similarly when $k=N-1$ we obtain%
\[
\frac{1}{v}\left\vert \left(  T^{-1}\Delta_{X}f\right)  _{N-1}\right\vert
<\left\Vert Df\right\Vert _{\infty;\Omega}.
\]

The case $2\leq k\leq m-1$ of \ref{av047} implies%
\begin{align*}
\left\vert T\left(  x\right)  \right\vert  & \geq\frac{1}{2}\left\vert
T\left(  x_{k+1:m}\right)  \right\vert \left(  \frac{x_{1}}{2^{k-1}}%
+\sum\limits_{n=2}^{k}\left\vert T\left(  x_{1:n-1}\right)  \right\vert
\frac{x_{n}}{2^{k-n}}\right)  +\\
& \qquad\qquad+\frac{1}{2}\left\vert T\left(  x_{1:k-1}\right)  \right\vert
\left(  \sum\limits_{n=k+1}^{m-1}\frac{x_{n}}{2^{n-k}}\left\vert T\left(
x_{n+1:m}\right)  \right\vert +\frac{x_{m}}{2^{m-k}}\right) \\
& =\frac{1}{2}\left\vert T\left(  x_{k+1:N-1}\right)  \right\vert \left(
\frac{x_{1}}{2^{k-1}}+\sum\limits_{n=2}^{k}\left\vert T\left(  x_{1:n-1}%
\right)  \right\vert \frac{x_{n}}{2^{k-n}}\right)  +\\
& \qquad\qquad+\frac{1}{2}\left\vert T\left(  x_{1:k-1}\right)  \right\vert
\left(  \sum\limits_{n=k+1}^{N-2}\frac{x_{n}}{2^{n-k}}\left\vert T\left(
x_{n+1:N-1}\right)  \right\vert +\frac{x_{N-1}}{2^{N-1-k}}\right)  ,
\end{align*}

so that%
\begin{multline*}
\frac{\left\vert T\left(  x_{k+1:N-1}\right)  \right\vert }{\left\vert
T\left(  x\right)  \right\vert }\left(  \frac{x_{1}}{2^{k-1}}+\sum
\limits_{n=2}^{k}\left\vert T\left(  x_{1:n-1}\right)  \right\vert \frac
{x_{n}}{2^{k-n}}\right)  +\\
+\frac{\left\vert T\left(  x_{1:k-1}\right)  \right\vert }{\left\vert T\left(
x\right)  \right\vert }\left(  \sum\limits_{n=k+1}^{N-2}\frac{x_{n}}{2^{n-k}%
}\left\vert T\left(  x_{n+1:N-1}\right)  \right\vert +\frac{x_{N-1}}%
{2^{N-1-k}}\right)  \leq2,
\end{multline*}

and thus%
\[
\frac{1}{v}\left\vert \left(  T^{-1}\Delta_{X}f\right)  _{k}\right\vert
\leq2\left\Vert Df\right\Vert _{\infty;\Omega},\quad2\leq k\leq N-2.
\]

Substituting this upper bound into \ref{av087} proves the second inequality
\ref{av081}.
\end{proof}

\section{Upper bounds for the Exact smoother derivative}

The goal of this section is to derive the upper bound $\left\Vert
Ds\right\Vert _{\infty;\Omega}\leq2\left\Vert f\right\Vert _{\infty;\Omega
}+\max\left\{  5,2+\rho N\right\}  \left\Vert Df\right\Vert _{\infty;\Omega}$,
on the assumption that $\operatorname*{diam}\Omega\leq1$, $f\in C_{B}^{\left(
0\right)  }\left(  \Omega\right)  $ and $Df\in L^{\infty}\left(
\Omega\right)  $. This is Theorem \ref{vThm_bound_deriv_hat_smth}. We also
derive a corresponding result for the scaled hat basis function.

The estimate \ref{av081}, namely $\left\vert \mu\left(  f\right)  \right\vert
_{\max}\leq2\left\Vert Df\right\Vert _{\infty;\Omega}+\left\vert \mu
_{1}\left(  f\right)  \right\vert $, means that if we can derive an upper
bound for $\left\vert \mu_{1}\left(  f\right)  \right\vert $ then we are done.
The approach is to first prove the preliminary estimates and equations of the
next two subsections.

In Subsection \ref{SbSect_bound_mu1} we start with formula \ref{av008}:%
\[
\mu_{1}\left(  f\right)  =\frac{f\left(  x^{\left(  1\right)  }\right)
+f\left(  x^{\left(  N\right)  }\right)  -\frac{1}{2}\beta^{T}T^{-1}\Delta
_{X}f}{v\left(  \left(  \frac{2}{v}-\mathbf{1}\cdot\frac{1}{v}\Delta
_{1}X+1\right)  -\frac{1}{4}\beta^{T}T^{-1}\beta\right)  }.
\]

and use the formulas \ref{av023} and \ref{av725} etc. to write $\beta
^{T}T^{-1}\beta$ and $\beta^{T}T^{-1}\Delta_{X}f$ in terms of $1/\left\vert
T\right\vert $ and the sub-determinants $\left\vert T_{i:j}\right\vert $.

We will show that the denominator exceeds $1$ so that%
\begin{align*}
\left\vert \mu_{1}\left(  f\right)  \right\vert  & \leq\frac{\left\vert
f\left(  x^{\left(  1\right)  }\right)  \right\vert +\left\vert f\left(
x^{\left(  N\right)  }\right)  \right\vert }{v\left(  \left(  \frac{2}%
{v}-\mathbf{1}\cdot\frac{1}{v}\Delta_{1}X+1\right)  -\frac{1}{4}\beta
^{T}T^{-1}\beta\right)  }+\frac{1}{2}\frac{\left\vert \beta^{T}T^{-1}%
\Delta_{X}f\right\vert }{v\left(  \left(  \frac{2}{v}-\mathbf{1}\cdot\frac
{1}{v}\Delta_{1}X+1\right)  -\frac{1}{4}\beta^{T}T^{-1}\beta\right)  }\\
& \leq2\left\Vert f\right\Vert _{\infty,\Omega}+\frac{1}{2}\frac{\left\vert
\beta^{T}T^{-1}\Delta_{X}f\right\vert }{v\left(  \left(  \frac{2}%
{v}-\mathbf{1}\cdot\frac{1}{v}\Delta_{1}X+1\right)  -\frac{1}{4}\beta
^{T}T^{-1}\beta\right)  },
\end{align*}

and the estimate \ref{av080}: $\left\vert \frac{\Delta_{1}f\left(  x^{\left(
n\right)  }\right)  }{\Delta_{1}x^{\left(  n\right)  }}\right\vert
\leq\left\Vert Df\right\Vert _{\infty;\Omega}$ is used to show the second term
is bounded above by

$\max\left\{  5,2+\rho N\right\}  \left\Vert Df\right\Vert _{\infty;\Omega}$.
But from \ref{av089}, $\left\Vert Ds\right\Vert _{\infty;\Omega}=\left\vert
\mu\left(  f\right)  \right\vert _{\max}$, so%
\[
\left\Vert Ds\right\Vert _{\infty;\Omega}\leq2\left\Vert f\right\Vert
_{\infty,\Omega}+\left\{  5,2+\rho N\right\}  \left\Vert Df\right\Vert
_{\infty;\Omega}.
\]

\subsection{Bounds for $T\left(  0_{n},x\right)  -T\left(  0_{n+m}\right)  $}

In this subsection we will derive the bounds \ref{av935} for $T\left(
0_{n},x\right)  -T\left(  0_{n+m}\right)  $ which are in effect better bounds
for $T\left(  0_{n},x\right)  $ than those of \ref{av915}. These results will
be used to obtain an upper bound for the numerator of the formula \ref{av008}
for $\mu_{1}\left(  f\right)  $.

From \ref{av915},%
\begin{equation}
\left\vert T\left(  0_{n},x\right)  \right\vert \leq\left\vert T\left(
x\right)  \right\vert ,\quad n\geq0,\label{av929}%
\end{equation}

where from Definition \ref{vDef_T_minors}, $0_{0}:=\left\{  {}\right\}  $ is
the empty sequence. The rest of this subsubsection will be devoted to
improving this estimate and obtaining some related lower bounds. From
\ref{av585}, for $x\in\mathbb{R}^{m}$, $m\geq2$,%
\begin{equation}
\left\vert T\left(  0_{n},x\right)  \right\vert =\left\{
\begin{array}
[c]{ll}%
\left\vert T\left(  0_{0},x\right)  \right\vert , & n=0,\\
\left\vert T\left(  x\right)  \right\vert -\frac{1}{4}\left\vert T\left(
x_{2:m}\right)  \right\vert , & n=1,\\
\left\vert T\left(  0_{1},x\right)  \right\vert -\frac{1}{4}\left\vert
T\left(  x\right)  \right\vert , & n=2,\\
\left\vert T\left(  0_{n-1},x\right)  \right\vert -\frac{1}{4}\left\vert
T\left(  0_{n-2},x\right)  \right\vert , & n\geq3.
\end{array}
\right. \label{av921}%
\end{equation}

This system of equations can be interpreted as a linear recurrence relation in
$\left\vert T\left(  0_{n},x\right)  \right\vert $. The auxiliary equation is
$u^{2}-u+1/4$ which has repeated root $1/2$ and so the general solution has
the form $\frac{1}{2^{n}}\left(  A+Bn\right)  $. Thus%
\begin{align*}
\left\vert T\left(  0_{0},x\right)  \right\vert  & =\left\vert T\left(
x\right)  \right\vert =A,\\
\left\vert T\left(  0_{1},x\right)  \right\vert  & =\left(  A+B\right)
\frac{1}{2}=\left(  \left\vert T\left(  x\right)  \right\vert +B\right)
\frac{1}{2}=\frac{1}{2}\left\vert T\left(  x\right)  \right\vert +\frac{1}%
{2}B,
\end{align*}

so that $A=\left\vert T\left(  x\right)  \right\vert $ and $B=2\left\vert
T\left(  0_{1},x\right)  \right\vert -\left\vert T\left(  x\right)
\right\vert $ and hence for $n\geq0$ and $m\geq2$,%
\begin{align}
\left\vert T\left(  0_{n},x\right)  \right\vert  & =\frac{1}{2^{n}}\left(
A+Bn\right) \nonumber\\
& =\frac{1}{2^{n}}\left(  \left\vert T\left(  x\right)  \right\vert +\left(
2\left\vert T\left(  0_{1},x\right)  \right\vert -\left\vert T\left(
x\right)  \right\vert \right)  n\right) \nonumber\\
& =\frac{1}{2^{n}}\left(  \left\vert T\left(  x\right)  \right\vert
+2n\left\vert T\left(  0_{1},x\right)  \right\vert -n\left\vert T\left(
x\right)  \right\vert \right) \nonumber\\
& =\frac{1}{2^{n}}\left(  2n\left\vert T\left(  0_{1},x\right)  \right\vert
-\left(  n-1\right)  \left\vert T\left(  x\right)  \right\vert \right)
\label{av925}\\
& =\frac{1}{2^{n}}\left(  2n\left(  \left\vert T\left(  x\right)  \right\vert
-\frac{1}{4}\left\vert T\left(  x_{2:m}\right)  \right\vert \right)  -\left(
n-1\right)  \left\vert T\left(  x\right)  \right\vert \right) \nonumber\\
& =\frac{1}{2^{n}}\left(  \left(  n+1\right)  \left\vert T\left(  x\right)
\right\vert -\frac{n}{2}\left\vert T\left(  x_{2:m}\right)  \right\vert
\right)  .\label{av922}%
\end{align}

From \ref{av922} and \ref{av929}:%
\begin{equation}
\left\vert T\left(  0_{n},x\right)  \right\vert <\frac{n+1}{2^{n}}\left\vert
T\left(  x\right)  \right\vert =\left\vert T\left(  0_{n}\right)  \right\vert
\text{ }\left\vert T\left(  x\right)  \right\vert ,\quad x\in\mathbb{R}%
_{\oplus}^{m};\text{ }m,n\geq1,\label{av926}%
\end{equation}

which improves on \ref{av929}. However, there is inequality when $x=0$.

We obtain equality when $x=0_{m}$ as follows. From \ref{av922}:%
\begin{align*}
&  \left\vert T\left(  0_{n},x\right)  \right\vert \\
&  =\frac{1}{2^{n}}\left(  \left(  n+1\right)  \left(  \left\vert T\left(
x\right)  \right\vert -\left\vert T\left(  0_{m}\right)  \right\vert
+\left\vert T\left(  0_{m}\right)  \right\vert \right)  -\frac{n}{2}\left(
\left\vert T\left(  x_{2:m}\right)  \right\vert -\left\vert T\left(
0_{m-1}\right)  \right\vert +\left\vert T\left(  0_{m-1}\right)  \right\vert
\right)  \right) \\
&  =\frac{1}{2^{n}}\left(
\begin{array}
[c]{r}%
\left(  n+1\right)  \left(  \left\vert T\left(  x\right)  \right\vert
-\left\vert T\left(  0_{m}\right)  \right\vert \right)  +\left(  n+1\right)
\left\vert T\left(  0_{m}\right)  \right\vert -\frac{n}{2}\left(  \left\vert
T\left(  x_{2:m}\right)  \right\vert -\left\vert T\left(  0_{m-1}\right)
\right\vert \right)  -\\
-\frac{n}{2}\left\vert T\left(  0_{m-1}\right)  \right\vert
\end{array}
\right) \\
&  =\frac{1}{2^{n}}\left(  \left(  n+1\right)  \left(  \left\vert T\left(
x\right)  \right\vert -\left\vert T\left(  0_{m}\right)  \right\vert \right)
-\frac{n}{2}\left(  \left\vert T\left(  x_{2:m}\right)  \right\vert
-\left\vert T\left(  0_{m-1}\right)  \right\vert \right)  \right)  -\\
&  \qquad\qquad-\frac{1}{2^{n}}\frac{n}{2}\left\vert T\left(  0_{m-1}\right)
\right\vert +\frac{1}{2^{n}}\left(  n+1\right)  \left\vert T\left(
0_{m}\right)  \right\vert \\
&  =\frac{1}{2^{n}}\left\{  \left(  n+1\right)  \left(  \left\vert T\left(
x\right)  \right\vert -\left\vert T\left(  0_{m}\right)  \right\vert \right)
-\frac{n}{2}\left(  \left\vert T\left(  x_{2:m}\right)  \right\vert
-\left\vert T\left(  0_{m-1}\right)  \right\vert \right)  \right\}  -\\
&  \qquad\qquad-\frac{1}{2^{n}}\frac{n}{2}\frac{m}{2^{m-1}}+\frac{1}{2^{n}%
}\left(  n+1\right)  \frac{m+1}{2^{m}},
\end{align*}

which becomes on rearranging%
\[
\left\vert T\left(  0_{n},x\right)  \right\vert -\left\vert T\left(
0_{m+n}\right)  \right\vert =\frac{1}{2^{n}}\left(  \left(  n+1\right)
\left(  \left\vert T\left(  x\right)  \right\vert -\left\vert T\left(
0_{m}\right)  \right\vert \right)  -\frac{n}{2}\left(  \left\vert T\left(
x_{2:m}\right)  \right\vert -\left\vert T\left(  0_{m-1}\right)  \right\vert
\right)  \right)  .
\]

Thus%
\begin{equation}
\left\vert T\left(  0_{n},x\right)  \right\vert -\left\vert T\left(
0_{m+n}\right)  \right\vert \leq\frac{n+1}{2^{n}}\left(  \left\vert T\left(
x\right)  \right\vert -\left\vert T\left(  0_{m}\right)  \right\vert \right)
,\quad x\in\mathbb{R}_{\oplus}^{m},\text{ }m\geq1,\label{av927}%
\end{equation}

with equality when $x=0_{m}$. This can be rewritten as%
\begin{equation}
\left\vert T\left(  0_{n},x\right)  \right\vert \leq\frac{n+1}{2^{n}%
}\left\vert T\left(  x\right)  \right\vert -\frac{mn}{2^{m+n}},\quad
x\in\mathbb{R}_{\oplus}^{m},\text{ }m\geq1,\label{av928}%
\end{equation}

with equality when $x=0_{m}$. These (equivalent) inequalities improve on both
\ref{av926} and \ref{av929}. Now%
\begin{align*}
\left\vert T\left(  0_{1},x\right)  \right\vert -\left\vert T\left(
0_{m+1}\right)  \right\vert  & =\left\vert T\left(  0,x\right)  \right\vert
-\left\vert T\left(  0_{m+1}\right)  \right\vert \\
& =\sum\limits_{k=1}^{m+1}\sum\limits_{\substack{\beta\leq\mathbf{1}_{m+1}
\\\left\vert \beta\right\vert =k}}c_{\beta}^{\left(  m+1\right)  }\left(
0,x\right)  ^{\beta}\\
& =\sum\limits_{k=1}^{m+1}\sum\limits_{\substack{\beta\leq\mathbf{1}_{m+1}
\\\beta_{1}+\left\vert \beta^{\prime}\right\vert =k}}c_{\beta_{1}%
,\beta^{\prime}}^{\left(  m+1\right)  }0^{\beta_{1}}x^{\beta^{\prime}}\\
& =\sum\limits_{k=1}^{m}\sum\limits_{\substack{\beta^{\prime}\leq
\mathbf{1}_{m} \\\left\vert \beta^{\prime}\right\vert =k}}c_{0,\beta^{\prime}%
}^{\left(  m+1\right)  }x^{\beta^{\prime}}.
\end{align*}

On the other hand, by part 6 of Theorem \ref{vThm_C^m_0n_gt_.5C^(m-1)_n}:
$\frac{1}{2}\frac{m+1}{m}c_{\beta^{\prime}}^{\left(  m\right)  }\leq
c_{0,\beta^{\prime}}^{\left(  m+1\right)  }$ when $\beta^{\prime}>0_{m}$ and
$m\geq1$, so that%
\begin{align*}
\left\vert T\left(  0_{1},x\right)  \right\vert -\left\vert T\left(
0_{m+1}\right)  \right\vert  & \geq\frac{1}{2}\frac{m+1}{m}\sum\limits_{k=1}%
^{m}\sum\limits_{\beta^{\prime}\leq\mathbf{1},\left\vert \beta^{\prime
}\right\vert =k}c_{\beta^{\prime}}^{\left(  m\right)  }x^{\beta^{\prime}}\\
& =\frac{1}{2}\frac{m+1}{m}\left(  \sum\limits_{k=0}^{m}\sum\limits_{\beta
^{\prime}\leq\mathbf{1},\left\vert \beta^{\prime}\right\vert =k}%
c_{\beta^{\prime}}^{\left(  m\right)  }x^{\beta^{\prime}}-\left\vert T\left(
0_{m}\right)  \right\vert \right) \\
& =\frac{1}{2}\frac{m+1}{m}\left(  \left\vert T\left(  x\right)  \right\vert
-\left\vert T\left(  0_{m}\right)  \right\vert \right)  ,\quad m\geq1,
\end{align*}

and $n$ applications yields%
\begin{align}
\left\vert T\left(  0_{n},x\right)  \right\vert -\left\vert T\left(
0_{m+n}\right)  \right\vert  & \geq\frac{1}{2}\frac{m+n}{m+n-1}\left(
\left\vert T\left(  0_{n-1},x\right)  \right\vert -\left\vert T\left(
0_{m+n-1}\right)  \right\vert \right) \nonumber\\
& \geq\frac{1}{2^{2}}\frac{m+n}{m+n-1}\frac{m+n-1}{m+n-2}\left(  \left\vert
T\left(  0_{n-2},x\right)  \right\vert -\left\vert T\left(  0_{m+n-2}\right)
\right\vert \right) \nonumber\\
& \geq\frac{1}{2^{n}}\left(  \frac{m+n}{m+n-1}\frac{m+n-1}{m+n-2}\ldots
\frac{m+1}{m}\right)  \left(  \left\vert T\left(  x\right)  \right\vert
-\left\vert T\left(  0_{m}\right)  \right\vert \right) \nonumber\\
& =\frac{1}{2^{n}}\frac{m+n}{m}\left(  \left\vert T\left(  x\right)
\right\vert -\left\vert T\left(  0_{m}\right)  \right\vert \right) \nonumber\\
& =\frac{1}{2^{n}}\left(  1+\frac{n}{m}\right)  \left(  \left\vert T\left(
x\right)  \right\vert -\left\vert T\left(  0_{m}\right)  \right\vert \right)
,\label{av933}%
\end{align}

or on rearranging%
\begin{align}
\left\vert T\left(  0_{n},x\right)  \right\vert  & \geq\frac{1}{2^{n}}\left(
1+\frac{n}{m}\right)  \left\vert T\left(  x\right)  \right\vert -\frac
{1}{2^{n}}\left(  1+\frac{n}{m}\right)  \left\vert T\left(  0_{m}\right)
\right\vert +\left\vert T\left(  0_{m+n}\right)  \right\vert \nonumber\\
& =\frac{1}{2^{n}}\left(  1+\frac{n}{m}\right)  \left\vert T\left(  x\right)
\right\vert -\left(  1+\frac{n}{m}\right)  \frac{m+1}{2^{m+n}}+\frac
{m+n+1}{2^{m+n}}\nonumber\\
& =\frac{1}{2^{n}}\left(  1+\frac{n}{m}\right)  \left\vert T\left(  x\right)
\right\vert -\frac{1}{2^{m+n}}\left\{  \left(  1+\frac{n}{m}\right)  \left(
m+1\right)  -\left(  m+n+1\right)  \right\} \nonumber\\
& =\frac{1}{2^{n}}\left(  1+\frac{n}{m}\right)  \left\vert T\left(  x\right)
\right\vert -\frac{1}{2^{m+n}}\left\{  m+n+1+\frac{n}{m}-\left(  m+n+1\right)
\right\} \nonumber\\
& =\frac{1}{2^{n}}\left(  1+\frac{n}{m}\right)  \left\vert T\left(  x\right)
\right\vert -\frac{1}{2^{m+n}}\frac{n}{m}.\label{av930}%
\end{align}

Combining \ref{av927} with \ref{av933} gives for $x\in\mathbb{R}_{\oplus}^{m}
$ and $m,n\geq1$,%
\begin{align}
0<\frac{1}{2^{n}}\left(  1+\frac{n}{m}\right)  \left(  \left\vert T\left(
x\right)  \right\vert -\left\vert T\left(  0_{m}\right)  \right\vert \right)
& \leq\left\vert T\left(  0_{n},x\right)  \right\vert -\left\vert T\left(
0_{m+n}\right)  \right\vert \nonumber\\
& \leq\frac{1}{2^{n}}\left(  1+n\right)  \left(  \left\vert T\left(  x\right)
\right\vert -\left\vert T\left(  0_{m}\right)  \right\vert \right)
,\label{av935}%
\end{align}

so \ref{av928} and \ref{av930} can be combined to give: if $x\in
\mathbb{R}_{\oplus}^{m}$ and $m,n\geq1$ then
\begin{equation}
0<\frac{1}{2^{n}}\left(  1+\frac{n}{m}\right)  \left\vert T\left(  x\right)
\right\vert -\frac{1}{2^{m+n}}\frac{n}{m}\leq\left\vert T\left(
0_{n},x\right)  \right\vert \leq\frac{1}{2^{n}}\left(  1+n\right)  \left\vert
T\left(  x\right)  \right\vert -\frac{mn}{2^{m+n}}.\label{av934}%
\end{equation}

\subsection{Lower bounds for $\left\vert T\left(  a\right)  \right\vert $}

This subsection is devoted to the derivation of the lower bounds \ref{av835}
and \ref{av896} for $T\left(  a\right)  $ which will be used to estimate the
numerator of the formula \ref{av008} for $\mu_{1}\left(  f\right)  $. The
bound for \ref{av835} involves linear combinations of terms like
$a_{k}\left\vert T\left(  a_{k+1:m}\right)  \right\vert $ and the bound
\ref{av896} involves linear combinations of terms like $\left\vert T\left(
a_{1:j-1}\right)  \right\vert a_{j}$. The approximations in this subsection
all derive from part 1 of Theorem \ref{vThm_C^m_0n_gt_.5C^(m-1)_n}.

From Definition \ref{vDef_T_minors} above%
\[
a_{l:m}:=\left(  a_{k}\right)  _{k=l}^{m}\in\mathbb{R}_{+}^{m-l+1}.
\]

If $m\geq2$:%
\begin{align*}
\sum\limits_{\substack{\beta\leq\mathbf{1} \\\left\vert \beta\right\vert
=k}}c_{\beta}^{\left(  m\right)  }a^{\beta}  & =\sum\limits_{\substack{\left(
\beta_{1},\beta^{\prime\prime}\right)  \leq\mathbf{1} \\\left\vert \beta
_{1},\beta^{\prime\prime}\right\vert =k}}c_{\beta_{1},\beta^{\prime\prime}%
}^{\left(  m\right)  }\left(  a_{1}a_{2:m}\right)  ^{\beta}\\
& =\sum\limits_{\substack{\left(  \beta_{1},\beta^{\prime\prime}\right)
\leq\mathbf{1} \\\left\vert \beta_{1},\beta^{\prime\prime}\right\vert
=k}}c_{\beta_{1},\beta^{\prime\prime}}^{\left(  m\right)  }a_{1}^{\beta_{1}%
}a_{2:m}^{\beta^{\prime\prime}}\\
& =\sum\limits_{\substack{\left(  1,\beta^{\prime\prime}\right)
\leq\mathbf{1} \\\left\vert 1,\beta^{\prime\prime}\right\vert =k}%
}c_{1,\beta^{\prime\prime}}^{\left(  m\right)  }a_{1}a_{2:m}^{\beta
^{\prime\prime}}+\sum\limits_{\substack{\left(  0,\beta^{\prime\prime}\right)
\leq\mathbf{1} \\\left\vert 0,\beta^{\prime\prime}\right\vert =k}%
}c_{0,\beta^{\prime\prime}}^{\left(  m\right)  }a_{2:m}^{\beta^{\prime\prime}%
}\\
& =a_{1}\sum\limits_{\substack{\beta^{\prime\prime}\leq\mathbf{1} \\\left\vert
\beta^{\prime\prime}\right\vert =k-1}}c_{1,\beta^{\prime\prime}}^{\left(
m\right)  }a_{2:m}^{\beta^{\prime\prime}}+\sum\limits_{\substack{\beta
^{\prime\prime}\leq\mathbf{1} \\\left\vert \beta^{\prime\prime}\right\vert
=k}}c_{0,\beta^{\prime\prime}}^{\left(  m\right)  }a_{2:m}^{\beta
^{\prime\prime}}.
\end{align*}

It is now more convenient to write%
\[
\sum\limits_{\substack{\beta\leq\mathbf{1} \\\left\vert \beta\right\vert
=k}}c_{\beta}^{\left(  m\right)  }a^{\beta}=a_{1}\sum\limits_{\substack{\beta
\leq\mathbf{1} \\\left\vert \beta\right\vert =k-1}}c_{1,\beta}^{\left(
m\right)  }a_{2:m}^{\beta}+\sum\limits_{\substack{\beta\leq\mathbf{1}
\\\left\vert \beta\right\vert =k}}c_{0,\beta}^{\left(  m\right)  }%
a_{2:m}^{\beta}.
\]

Since the non-trivial sequences of zeros in $\left(  1,\beta\right)  $ and
$\beta$ are identical, Theorem \ref{vThm_formula_for_C^(m)_idx} tells us that
$c_{1,\beta}^{\left(  m\right)  }=c_{\beta}^{\left(  m-1\right)  }$. Also by
part 1 of Theorem \ref{vThm_C^m_0n_gt_.5C^(m-1)_n}, $\frac{1}{2}c_{\beta
}^{\left(  m-1\right)  }<c_{0,\beta}^{\left(  m\right)  }<c_{\beta}^{\left(
m-1\right)  }$. Hence we have the upper and lower bounds: for $1\leq k\leq m-1
$ and $m\geq2$,
\begin{align}
a_{1}\sum\limits_{\substack{\beta\leq\mathbf{1} \\\left\vert \beta\right\vert
=k-1}}c_{\beta}^{\left(  m-1\right)  }a_{2:m}^{\beta}+\frac{1}{2}%
\sum\limits_{\substack{\beta\leq\mathbf{1} \\\left\vert \beta\right\vert
=k}}c_{\beta}^{\left(  m-1\right)  }a_{2:m}^{\beta} &  <\sum
\limits_{\substack{\beta\leq\mathbf{1} \\\left\vert \beta\right\vert
=k}}c_{\beta}^{\left(  m\right)  }a^{\beta}\nonumber\\
&  <a_{1}\sum\limits_{\substack{\beta\leq\mathbf{1} \\\left\vert
\beta\right\vert =k-1}}c_{\beta}^{\left(  m-1\right)  }a_{2:m}^{\beta}%
+\sum\limits_{\substack{\beta\leq\mathbf{1} \\\left\vert \beta\right\vert
=k}}c_{\beta}^{\left(  m-1\right)  }a_{2:m}^{\beta},\label{av827}%
\end{align}

and thus%
\begin{align*}
a_{1}\sum\limits_{k=1}^{m-1}\sum\limits_{\substack{\beta\leq\mathbf{1}
\\\left\vert \beta\right\vert =k-1}}c_{\beta}^{\left(  m-1\right)  }%
a_{2:m}^{\beta} &  +\frac{1}{2}\sum\limits_{k=1}^{m-1}\sum
\limits_{\substack{\beta\leq\mathbf{1} \\\left\vert \beta\right\vert
=k}}c_{\beta}^{\left(  m-1\right)  }a_{2:m}^{\beta}\\
&  <\sum\limits_{k=1}^{m-1}\sum\limits_{\substack{\beta\leq\mathbf{1}
\\\left\vert \beta\right\vert =k}}c_{\beta}^{\left(  m\right)  }a^{\beta}\\
&  <a_{1}\sum\limits_{k=1}^{m-1}\sum\limits_{\substack{\beta\leq\mathbf{1}
\\\left\vert \beta\right\vert =k-1}}c_{\beta}^{\left(  m-1\right)  }%
a_{2:m}^{\beta}+\sum\limits_{k=1}^{m-1}\sum\limits_{\substack{\beta
\leq\mathbf{1} \\\left\vert \beta\right\vert =k}}c_{\beta}^{\left(
m-1\right)  }a_{2:m}^{\beta},
\end{align*}

i.e.%
\begin{align}
a_{1}\sum\limits_{k=0}^{m-2}\sum\limits_{\substack{\beta\leq\mathbf{1}
\\\left\vert \beta\right\vert =k}}c_{\beta}^{\left(  m-1\right)  }%
a_{2:m}^{\beta} &  +\frac{1}{2}\sum\limits_{k=1}^{m-1}\sum
\limits_{\substack{\beta\leq\mathbf{1} \\\left\vert \beta\right\vert
=k}}c_{\beta}^{\left(  m-1\right)  }a_{2:m}^{\beta}\nonumber\\
&  <\sum\limits_{k=1}^{m-1}\sum\limits_{\substack{\beta\leq\mathbf{1}
\\\left\vert \beta\right\vert =k}}c_{\beta}^{\left(  m\right)  }a^{\beta
}\nonumber\\
&  <a_{1}\sum\limits_{k=0}^{m-2}\sum\limits_{\substack{\beta\leq\mathbf{1}
\\\left\vert \beta\right\vert =k}}c_{\beta}^{\left(  m-1\right)  }%
a_{2:m}^{\beta}+\sum\limits_{k=1}^{m-1}\sum\limits_{\substack{\beta
\leq\mathbf{1} \\\left\vert \beta\right\vert =k}}c_{\beta}^{\left(
m-1\right)  }a_{2:m}^{\beta}.\label{av831}%
\end{align}

Since we are dealing with $T_{m}$s we want to complete the outer summations in
the LHS of \ref{av831}. Noting from \ref{av598} that $c_{\mathbf{0}}^{\left(
m\right)  }=\frac{m+1}{2^{m}}$ and from \ref{av595} that $c_{\mathbf{1}%
}^{\left(  m\right)  }=1$, we are lead to the following sequence of
inequalities:%
\begin{multline*}
a_{1}\sum\limits_{k=0}^{m-1}\sum\limits_{\substack{\beta\leq\mathbf{1}
\\\left\vert \beta\right\vert =k}}c_{\beta}^{\left(  m-1\right)  }%
a_{2:m}^{\beta}-a_{1}\sum\limits_{\substack{\beta\leq\mathbf{1} \\\left\vert
\beta\right\vert =m-1}}c_{\beta}^{\left(  m-1\right)  }a_{2:m}^{\beta}%
+\frac{1}{2}\sum\limits_{k=0}^{m-1}\sum\limits_{\substack{\beta\leq\mathbf{1}
\\\left\vert \beta\right\vert =k}}c_{\beta}^{\left(  m-1\right)  }%
a_{2:m}^{\beta}-\frac{1}{2}\sum\limits_{\substack{\beta\leq\mathbf{1}
\\\left\vert \beta\right\vert =0}}c_{\beta}^{\left(  m-1\right)  }%
a_{2:m}^{\beta}\\
<\sum\limits_{k=0}^{m}\sum\limits_{\substack{\beta\leq\mathbf{1} \\\left\vert
\beta\right\vert =k}}c_{\beta}^{\left(  m\right)  }a^{\beta}-\sum
\limits_{\substack{\beta\leq\mathbf{1} \\\left\vert \beta\right\vert
=0}}c_{\beta}^{\left(  m\right)  }a^{\beta}-\sum\limits_{\substack{\beta
\leq\mathbf{1} \\\left\vert \beta\right\vert =m}}c_{\beta}^{\left(  m\right)
}a^{\beta},
\end{multline*}

\begin{multline*}
a_{1}\sum\limits_{k=0}^{m-1}\sum\limits_{\substack{\beta\leq\mathbf{1}
\\\left\vert \beta\right\vert =k}}c_{\beta}^{\left(  m-1\right)  }%
a_{2:m}^{\beta}-a_{1}\sum\limits_{\substack{\beta\leq\mathbf{1} \\\left\vert
\beta\right\vert =m-1}}c_{\beta}^{\left(  m-1\right)  }a_{2:m}^{\beta}%
+\frac{1}{2}\sum\limits_{k=0}^{m-1}\sum\limits_{\substack{\beta\leq\mathbf{1}
\\\left\vert \beta\right\vert =k}}c_{\beta}^{\left(  m-1\right)  }%
a_{2:m}^{\beta}-\frac{1}{2}\sum\limits_{\substack{\beta\leq\mathbf{1}
\\\left\vert \beta\right\vert =0}}c_{\beta}^{\left(  m-1\right)  }%
a_{2:m}^{\beta}\\
<\sum\limits_{k=0}^{m}\sum\limits_{\substack{\beta\leq\mathbf{1} \\\left\vert
\beta\right\vert =k}}c_{\beta}^{\left(  m\right)  }a^{\beta}-\sum
\limits_{\substack{\beta\leq\mathbf{1} \\\left\vert \beta\right\vert
=0}}c_{\beta}^{\left(  m\right)  }a^{\beta}-\sum\limits_{\substack{\beta
\leq\mathbf{1} \\\left\vert \beta\right\vert =m}}c_{\beta}^{\left(  m\right)
}a^{\beta},
\end{multline*}

\begin{multline*}
a_{1}\sum\limits_{k=0}^{m-1}\sum\limits_{\substack{\beta\leq\mathbf{1}
\\\left\vert \beta\right\vert =k}}c_{\beta}^{\left(  m-1\right)  }%
a_{2:m}^{\beta}-c_{\mathbf{1}}^{\left(  m-1\right)  }a_{1}a_{2:m}^{\mathbf{1}%
}+\frac{1}{2}\sum\limits_{k=0}^{m-1}\sum\limits_{\substack{\beta\leq\mathbf{1}
\\\left\vert \beta\right\vert =k}}c_{\beta}^{\left(  m-1\right)  }%
a_{2:m}^{\beta}-\frac{1}{2}c_{\mathbf{0}}^{\left(  m-1\right)  }\\
<\sum\limits_{k=0}^{m}\sum\limits_{\substack{\beta\leq\mathbf{1} \\\left\vert
\beta\right\vert =k}}c_{\beta}^{\left(  m\right)  }a^{\beta}-c_{\mathbf{0}%
}^{\left(  m\right)  }-c_{\mathbf{1}}^{\left(  m\right)  }a^{\mathbf{1}},
\end{multline*}

\begin{multline*}
a_{1}\sum\limits_{k=0}^{m-1}\sum\limits_{\substack{\beta\leq\mathbf{1}
\\\left\vert \beta\right\vert =k}}c_{\beta}^{\left(  m-1\right)  }%
a_{2:m}^{\beta}-a^{\mathbf{1}}+\frac{1}{2}\sum\limits_{k=0}^{m-1}%
\sum\limits_{\substack{\beta\leq\mathbf{1} \\\left\vert \beta\right\vert =k
}}c_{\beta}^{\left(  m-1\right)  }a_{2:m}^{\beta}-\frac{1}{2}\frac{m}{2^{m-1}%
}\\
<\sum\limits_{k=0}^{m}\sum\limits_{\substack{\beta\leq\mathbf{1} \\\left\vert
\beta\right\vert =k}}c_{\beta}^{\left(  m\right)  }a^{\beta}-\frac{m+1}{2^{m}%
}-a^{\mathbf{1}},
\end{multline*}

\[
a_{1}\sum\limits_{k=0}^{m-1}\sum\limits_{\substack{\beta\leq\mathbf{1}
\\\left\vert \beta\right\vert =k}}c_{\beta}^{\left(  m-1\right)  }%
a_{2:m}^{\beta}+\frac{1}{2}\sum\limits_{k=0}^{m-1}\sum\limits_{\substack{\beta
\leq\mathbf{1} \\\left\vert \beta\right\vert =k}}c_{\beta}^{\left(
m-1\right)  }a_{2:m}^{\beta}-\frac{m}{2^{m}}<\sum\limits_{k=0}^{m}%
\sum\limits_{\substack{\beta\leq\mathbf{1} \\\left\vert \beta\right\vert
=k}}c_{\beta}^{\left(  m\right)  }a^{\beta}-\frac{m+1}{2^{m}},
\]

\[
a_{1}\sum\limits_{k=0}^{m-1}\sum\limits_{\substack{\beta\leq\mathbf{1}
\\\left\vert \beta\right\vert =k}}c_{\beta}^{\left(  m-1\right)  }%
a_{2:m}^{\beta}+\frac{1}{2}\sum\limits_{k=0}^{m-1}\sum\limits_{\substack{\beta
\leq\mathbf{1} \\\left\vert \beta\right\vert =k}}c_{\beta}^{\left(
m-1\right)  }a_{2:m}^{\beta}<\sum\limits_{k=0}^{m}\sum\limits_{\substack{\beta
\leq\mathbf{1} \\\left\vert \beta\right\vert =k}}c_{\beta}^{\left(  m\right)
}a^{\beta}-\frac{1}{2^{m}},
\]%
\begin{align}
a_{1}\sum\limits_{k=0}^{m-1}\sum\limits_{\substack{\beta\leq\mathbf{1}
\\\left\vert \beta\right\vert =k}}c_{\beta}^{\left(  m-1\right)  }%
a_{2:m}^{\beta}+\frac{1}{2}\sum\limits_{k=0}^{m-1}\sum\limits_{\substack{\beta
\leq\mathbf{1} \\\left\vert \beta\right\vert =k}}c_{\beta}^{\left(
m-1\right)  }a_{2:m}^{\beta}+\frac{1}{2^{m}}  & <\sum\limits_{k=1}^{m}%
\sum\limits_{\substack{\beta\leq\mathbf{1} \\\left\vert \beta\right\vert =k
}}c_{\beta}^{\left(  m\right)  }a^{\beta},\nonumber\\
when\text{ }m  & \geq2,\label{av830}%
\end{align}

so that%
\begin{equation}
a_{1}\left\vert T\left(  a_{2:m}\right)  \right\vert +\frac{1}{2}\left\vert
T\left(  a_{2:m}\right)  \right\vert +\frac{1}{2^{m}}<\left\vert T\left(
a\right)  \right\vert ,\quad m\geq2.\label{av828}%
\end{equation}

When we replace $a$ by $a_{2:m}$ we get%
\[
a_{2}\left\vert T\left(  a_{3:m}\right)  \right\vert +\frac{1}{2}\left\vert
T\left(  a_{3:m}\right)  \right\vert +\frac{1}{2^{m-1}}<\left\vert T\left(
a_{2:m}\right)  \right\vert ,\quad m\geq3,
\]

and substituting for the second occurrence of $T_{m-2}\left(  a_{3:m}\right)
$ in \ref{av828} gives%
\begin{equation}
a_{1}\left\vert T\left(  a_{2:m}\right)  \right\vert +\frac{a_{2}}%
{2}\left\vert T\left(  a_{3:m}\right)  \right\vert +\frac{1}{2^{2}}\left\vert
T\left(  a_{3:m}\right)  \right\vert +\frac{2}{2^{m}}<\left\vert T\left(
a\right)  \right\vert ,\quad m\geq3.\label{av832}%
\end{equation}

When we replace $a$ by $a_{3:m}$ in \ref{av828} we get
\[
a_{3}\left\vert T\left(  a_{4:m}\right)  \right\vert +\frac{1}{2}\left\vert
T\left(  a_{4:m}\right)  \right\vert +\frac{1}{2^{m-2}}<\left\vert T\left(
a_{3:m}\right)  \right\vert ,\quad m\geq4,
\]

and substituting for the second occurrence of $\left\vert T\left(
a_{3:m}\right)  \right\vert $ in \ref{av832} gives for $m\geq4:$
\begin{align}
a_{1}\left\vert T\left(  a_{2:m}\right)  \right\vert  &  +\frac{a_{2}}%
{2}\left\vert T\left(  a_{3:m}\right)  \right\vert +\frac{1}{2^{2}}\left(
a_{3}\left\vert T\left(  a_{4:m}\right)  \right\vert +\frac{1}{2}\left\vert
T\left(  a_{4:m}\right)  \right\vert +\frac{1}{2^{m-2}}\right)  +\frac
{1}{2^{m-1}}\nonumber\\
&  =a_{1}\left\vert T\left(  a_{2:m}\right)  \right\vert +\frac{a_{2}}%
{2}\left\vert T\left(  a_{3:m}\right)  \right\vert +\frac{a_{3}}{2^{2}%
}\left\vert T\left(  a_{4:m}\right)  \right\vert +\frac{1}{2^{3}}\left\vert
T\left(  a_{4:m}\right)  \right\vert +\frac{1}{2^{m}}+\frac{1}{2^{m-1}%
}\nonumber\\
&  =a_{1}\left\vert T\left(  a_{2:m}\right)  \right\vert +\frac{a_{2}}%
{2}\left\vert T\left(  a_{3:m}\right)  \right\vert +\frac{a_{3}}{2^{2}%
}\left\vert T\left(  a_{4:m}\right)  \right\vert +\frac{1}{2^{3}}\left\vert
T\left(  a_{4:m}\right)  \right\vert +\frac{3}{2^{m}}\nonumber\\
&  <\left\vert T\left(  a\right)  \right\vert .\label{av833}%
\end{align}

Noting \ref{av828}, \ref{av832} and \ref{av833} we suspect that in general for
$m\geq j$,
\[
a_{1}\left\vert T\left(  a_{2:m}\right)  \right\vert +\frac{a_{2}}%
{2}\left\vert T\left(  a_{3:m}\right)  \right\vert +\ldots+\frac{a_{j-1}%
}{2^{j-2}}\left\vert T\left(  a_{j:m}\right)  \right\vert +\frac{1}{2^{j-1}%
}\left\vert T\left(  a_{j:m}\right)  \right\vert +\frac{j-1}{2^{m}}<\left\vert
T\left(  a\right)  \right\vert ,
\]

so that when $j=m$ and $m\geq2$:%
\begin{align*}
a_{1}\left\vert T\left(  a_{2:m}\right)  \right\vert  &  +\frac{a_{2}}%
{2}\left\vert T\left(  a_{3:m}\right)  \right\vert +\ldots+\frac{a_{m-1}%
}{2^{m-2}}\left\vert T\left(  a_{m,m}\right)  \right\vert +\frac{1}{2^{m-1}%
}\left\vert T\left(  a_{m,m}\right)  \right\vert +\frac{m-1}{2^{m}}\\
&  =a_{1}\left\vert T\left(  a_{2:m}\right)  \right\vert +\frac{a_{2}}%
{2}\left\vert T\left(  a_{3:m}\right)  \right\vert +\ldots+\frac{a_{m-1}%
}{2^{m-2}}\left\vert T\left(  a_{m,m}\right)  \right\vert +\frac{a_{m}%
+1}{2^{m-1}}+\frac{m-1}{2^{m}}\\
&  =a_{1}\left\vert T\left(  a_{2:m}\right)  \right\vert +\frac{a_{2}}%
{2}\left\vert T\left(  a_{3:m}\right)  \right\vert +\ldots+\frac{a_{m-1}%
}{2^{m-2}}\left\vert T\left(  a_{m,m}\right)  \right\vert +\frac{a_{m}%
}{2^{m-1}}+\frac{1}{2^{m-1}}+\frac{m-1}{2^{m}}\\
&  =a_{1}\left\vert T\left(  a_{2:m}\right)  \right\vert +\frac{a_{2}}%
{2}\left\vert T\left(  a_{3:m}\right)  \right\vert +\ldots+\frac{a_{m-1}%
}{2^{m-2}}\left\vert T\left(  a_{m,m}\right)  \right\vert +\frac{a_{m}%
}{2^{m-1}}+\frac{m+1}{2^{m}}\\
&  <\left\vert T\left(  a\right)  \right\vert
\end{align*}

or more compactly%
\begin{equation}
\sum\limits_{k=1}^{m}\frac{a_{k}}{2^{k-1}}\left\vert T\left(  a_{k+1:m}%
\right)  \right\vert <\left\vert T\left(  a\right)  \right\vert -\frac
{m+1}{2^{m}},\quad m\geq2,\text{ }a\in\mathbb{R}_{+}^{m}.\label{av835}%
\end{equation}

From Corollary \ref{vCor_1_formula_for_C^(m)_gam} the reverse order
permutation satisfies $\left(  \mathcal{R}a\right)  _{k}=a_{m+1-k}$ and so%
\begin{align*}
\sum\limits_{k=1}^{m}\frac{\left(  \mathcal{R}a\right)  _{k}}{2^{k-1}%
}\left\vert T\left(  \left(  \mathcal{R}a\right)  _{k+1:m}\right)  \right\vert
=\sum\limits_{k=1}^{m}\frac{a_{m+1-k}}{2^{k-1}}\left\vert T\left(
a_{1:m-k}\right)  \right\vert  & <\left\vert T\left(  \mathcal{R}a\right)
\right\vert -\frac{m+1}{2^{m}}\\
& =\left\vert T\left(  a\right)  \right\vert -\frac{m+1}{2^{m}},
\end{align*}

i.e.%
\[
\sum\limits_{k=1}^{m}\frac{a_{m+1-k}}{2^{k-1}}\left\vert T\left(
a_{1:m-k}\right)  \right\vert =\sum\limits_{j=1}^{m}\frac{a_{j}}{2^{m-j}%
}\left\vert T\left(  a_{1:j-1}\right)  \right\vert <\left\vert T\left(
a\right)  \right\vert -\frac{m+1}{2^{m}},
\]

and we have%
\begin{equation}
\sum\limits_{j=1}^{m}\frac{a_{j}}{2^{m-j}}\left\vert T\left(  a_{1:j-1}%
\right)  \right\vert <\left\vert T\left(  a\right)  \right\vert -\frac
{m+1}{2^{m}},\quad m\geq2,\text{ }a\in\mathbb{R}_{+}^{m}.\label{av896}%
\end{equation}

\subsection{Upper bounds for $\left\vert \mu_{1}\left(  f\right)  \right\vert
$ and $\left\vert \mu\left(  f\right)  \right\vert _{\max}$%
\label{SbSect_bound_mu1}}

We will now prove the uniform upper bound \ref{av302} for $\mu_{1}\left(
f\right)  $ and the upper bound \ref{av036} for $\left\vert \mu\left(
f\right)  \right\vert _{\max}$. To do so we start with equation \ref{av008}:%
\begin{equation}
\mu_{1}\left(  f\right)  =\frac{f\left(  x^{\left(  1\right)  }\right)
+f\left(  x^{\left(  N\right)  }\right)  -\frac{1}{2}\beta^{T}T^{-1}\Delta
_{X}f}{v\left(  \left(  \frac{2}{v}-\mathbf{1}\cdot\frac{1}{v}\Delta
_{1}X+1\right)  -\frac{1}{4}\beta^{T}T^{-1}\beta\right)  },\label{av006a}%
\end{equation}

where $x=\frac{1}{v}\Delta_{1}X$. Write%
\begin{align}
\mu_{1}\left(  f\right)   & =\frac{\mathcal{N}\left(  f\right)  }{\mathcal{D}%
},\label{av980}\\
\mathcal{N}\left(  f\right)   & =f\left(  x^{\left(  1\right)  }\right)
+f\left(  x^{\left(  N\right)  }\right)  -\frac{1}{2}\beta^{T}T^{-1}\Delta
_{X}f,\label{av981}\\
\mathcal{D}  & =v\left(  \left(  \frac{2}{v}-\mathbf{1\cdot}x+1\right)
-\frac{1}{4}\beta^{T}T^{-1}\beta\right)  .\label{av982}%
\end{align}
\medskip

\fbox{The denominator $\mathcal{D}$} Here we show that $\mathcal{D}>0$. First
note that we have assumed in the introduction that%
\[
0<\mathbf{1\cdot}x=x^{\left(  N\right)  }-x^{\left(  1\right)  }\leq1.
\]

Using equation \ref{av756} to write $\beta^{T}T^{-1}\beta$ in terms of the
sub-determinants $\left\vert T_{i:j}\right\vert $, \ref{av982} becomes%
\begin{align}
\frac{1}{v}\mathcal{D}\left\vert T\right\vert  & =\left(  \frac{2}%
{v}-\mathbf{1\cdot}\frac{1}{v}\Delta_{1}X+1\right)  \left\vert T\right\vert
-\frac{1}{4}\left(  \beta^{T}T^{-1}\beta\right)  \left\vert T\right\vert
\nonumber\\
& =\left(  \frac{2}{v}-\mathbf{1\cdot}\frac{1}{v}\Delta_{1}X+1\right)
\left\vert T\right\vert -\frac{1}{4}\left(  \left\vert T_{2:N-1}\right\vert
-\left(  \frac{1}{2}\right)  ^{N-3}+\left\vert T_{1:N-2}\right\vert \right)
\nonumber\\
& =\left(  \frac{2}{v}-\mathbf{1}\cdot x+1\right)  \left\vert T\right\vert
-\frac{1}{4}\left\vert T_{2:N-1}\right\vert -\frac{1}{4}\left\vert
T_{1:N-2}\right\vert +\frac{1}{2^{N-1}}\nonumber\\
& =\left(  1+z\right)  \left\vert T\right\vert -\frac{1}{4}\left\vert
T_{2:N-1}\right\vert -\frac{1}{4}\left\vert T_{1:N-2}\right\vert +\frac
{1}{2^{N-1}},\label{av939}%
\end{align}

where
\begin{equation}
z=\frac{1}{v}+\left(  \frac{1}{v}-\mathbf{1}\cdot x\right)  =\frac{1}{v}%
+\frac{1}{v}\left(  1-\mathbf{1}\cdot vx\right)  =\frac{1}{v}+\frac{1}%
{v}\left(  1-\left(  x^{\left(  N\right)  }-x^{\left(  1\right)  }\right)
\right)  >\frac{1}{v}.\label{av938}%
\end{equation}

Thus%
\begin{align}
\mathcal{D}\left\vert T\right\vert  & >v\left(  \left(  1+\frac{1}{v}\right)
\left\vert T\right\vert -\frac{1}{4}\left\vert T_{2:N-1}\right\vert -\frac
{1}{4}\left\vert T_{1:N-2}\right\vert +\frac{1}{2^{N-1}}\right) \nonumber\\
& =\left(  1+v\right)  \left\vert T\right\vert +v\left(  \frac{1}{2^{N-1}%
}-\frac{1}{4}\left\vert T_{2:N-1}\right\vert -\frac{1}{4}\left\vert
T_{1:N-2}\right\vert \right)  .\label{av034}%
\end{align}

With reference to \ref{av939}, using equations \ref{av936} and \ref{av840},%
\begin{align*}
&  \left(  1+z\right)  \left\vert T\right\vert -\frac{1}{4}\left\vert
T_{2:N-1}\right\vert -\frac{1}{4}\left\vert T_{1:N-2}\right\vert +\frac
{1}{2^{N-1}}\\
&  =\left(  1+z\right)  \left\vert T\left(  x\right)  \right\vert -\frac{1}%
{4}\left\vert T\left(  x_{2:N-1}\right)  \right\vert -\frac{1}{4}\left\vert
T\left(  x_{1:N-2}\right)  \right\vert +\frac{1}{2^{N-1}}\\
&  =\left(  \frac{1}{2}\left(  1+z\right)  \left\vert T\left(  x\right)
\right\vert -\frac{1}{4}\left\vert T\left(  x_{2:N-1}\right)  \right\vert
\right)  +\left(  \frac{1}{2}\left(  1+z\right)  \left\vert T\left(  x\right)
\right\vert -\frac{1}{4}\left\vert T\left(  x_{1:N-2}\right)  \right\vert
\right)  +\frac{1}{2^{N-1}}\\
&  =\left(  \frac{1}{2}\left(  1+z\right)  \left(  x_{1}\left\vert T\left(
x_{2:N-1}\right)  \right\vert +\left\vert T\left(  0,x_{2:N-1}\right)
\right\vert \right)  -\frac{1}{4}\left\vert T\left(  x_{2:N-1}\right)
\right\vert \right)  +\\
&  \qquad\qquad+\left(  \frac{1}{2}\left(  1+z\right)  \left(  x_{N-1}%
\left\vert T\left(  x_{1:N-2}\right)  \right\vert +\left\vert T\left(
0,x_{1:N-2}\right)  \right\vert \right)  -\frac{1}{4}\left\vert T\left(
x_{1:N-2}\right)  \right\vert \right)  +\frac{1}{2^{N-1}}\\
&  =\frac{1}{2}\left(  1+z\right)  x_{1}\left\vert T\left(  x_{2:N-1}\right)
\right\vert +\frac{1}{2}\left(  1+z\right)  \left\vert T\left(  0,x_{2:N-1}%
\right)  \right\vert -\frac{1}{4}\left\vert T\left(  x_{2:N-1}\right)
\right\vert +\\
&  \qquad+\frac{1}{2}\left(  1+z\right)  x_{N-1}\left\vert T\left(
x_{1:N-2}\right)  \right\vert +\frac{1}{2}\left(  1+z\right)  \left\vert
T\left(  x_{1:N-2},0\right)  \right\vert -\frac{1}{4}\left\vert T\left(
x_{1:N-2}\right)  \right\vert +\frac{1}{2^{N-1}}.
\end{align*}

Applying the lower bounds \ref{av876} and \ref{av878} i.e. $\left\vert
T\left(  0,y\right)  \right\vert \geq\frac{2+p}{1+p}\frac{1}{2}\left\vert
T\left(  y\right)  \right\vert $ and $\left\vert T\left(  y,0\right)
\right\vert \geq\frac{2+p}{1+p}\frac{1}{2}\left\vert T\left(  y\right)
\right\vert $ when $y\in\mathbb{R}_{\oplus}^{p}$, we get%
\begin{align}
\frac{1}{v}\mathcal{D}\left\vert T\right\vert  &  =\left(  1+z\right)
\left\vert T\right\vert -\frac{1}{4}\left\vert T_{2:N-1}\right\vert -\frac
{1}{4}\left\vert T_{1:N-2}\right\vert +\frac{1}{2^{N-1}}\nonumber\\
&  \geq\frac{1}{2}\left(  1+z\right)  x_{1}\left\vert T_{2:N-1}\right\vert
+\frac{1}{2}\left(  1+z\right)  \frac{N}{N-1}\frac{1}{2}\left\vert
T_{2:N-1}\right\vert -\frac{1}{4}\left\vert T_{2:N-1}\right\vert +\nonumber\\
&  \qquad+\frac{1}{2}\left(  1+z\right)  x_{N-1}\left\vert T_{1:N-2}%
\right\vert +\frac{1}{2}\left(  1+z\right)  \frac{N}{N-1}\frac{1}{2}\left\vert
T_{1:N-2}\right\vert -\frac{1}{4}\left\vert T_{1:N-2}\right\vert +\frac
{1}{2^{N-1}}\nonumber\\
&  =\left(  \frac{1}{2}\left(  1+z\right)  x_{1}+\frac{1}{2}\left(
1+z\right)  \frac{N}{N-1}\frac{1}{2}-\frac{1}{4}\right)  \left\vert
T_{2:N-1}\right\vert +\nonumber\\
&  \qquad+\left(  \frac{1}{2}\left(  1+z\right)  x_{N-1}+\frac{1}{2}\left(
1+z\right)  \frac{N}{N-1}\frac{1}{2}-\frac{1}{4}\right)  \left\vert
T_{1:N-2}\right\vert +\frac{1}{2^{N-1}}\nonumber\\
&  =\frac{1}{2}\left(  \left(  1+z\right)  x_{1}+\left(  1+z\right)  \left(
1+\frac{1}{N-1}\right)  -\frac{1}{2}\right)  \left\vert T_{2:N-1}\right\vert
+\nonumber\\
&  \qquad+\frac{1}{2}\left(  \left(  1+z\right)  x_{N-1}+\left(  1+z\right)
\left(  1+\frac{1}{N-1}\right)  -\frac{1}{2}\right)  \left\vert T_{1:N-2}%
\right\vert +\frac{1}{2^{N-1}}\nonumber\\
&  >\frac{1}{2}\left(  x_{1}+1-\frac{1}{2}\right)  \left\vert T_{2:N-1}%
\right\vert +\frac{1}{2}\left(  x_{N-1}+1-\frac{1}{2}\right)  \left\vert
T_{1:N-2}\right\vert \nonumber\\
&  =\frac{1}{2}\left(  x_{1}+\frac{1}{2}\right)  \left\vert T_{2:N-1}%
\right\vert +\frac{1}{2}\left(  x_{N-1}+\frac{1}{2}\right)  \left\vert
T_{1:N-2}\right\vert \label{av943}\\
&  >0,\nonumber
\end{align}

for all $v>0$ and $N\geq3$. We have shown that
\begin{equation}
\mathcal{D}>0\text{ }when\text{ }v>0\text{ }and\text{ }N\geq3.\label{av033}%
\end{equation}
\medskip

\fbox{The numerator $\mathcal{N}\left(  f\right)  $} From the minor formulas
\ref{av023},%
\[%
\begin{array}
[c]{ll}%
M_{j,n}=\left(  -\frac{1}{2}\right)  ^{n-j}\left\vert T_{1:j-1}\right\vert
\left\vert T_{n+1:N-1}\right\vert , & 1\leq j\leq n\leq N-1,\\
M_{n,j}=M_{j,n}, &
\end{array}
\]

where
\[
T_{N:j}=T_{0:j}=T_{j:N}=T_{j:0}=1,\quad1\leq j\leq N-1.
\]

Hence%
\[
M_{n,1}=M_{1,n}=\left(  -\frac{1}{2}\right)  ^{n-1}\left\vert T_{1:0}%
\right\vert \left\vert T_{n+1:N-1}\right\vert =\left(  -\frac{1}{2}\right)
^{n-1}\left\vert T_{n+1:N-1}\right\vert .
\]

and
\[
M_{n,N-1}=\left(  -\frac{1}{2}\right)  ^{N-1-n}\left\vert T_{1:n-1}\right\vert
\left\vert T_{N:N-1}\right\vert =\left(  -\frac{1}{2}\right)  ^{N-1-n}%
\left\vert T_{1:n}\right\vert .
\]

We use the equation \ref{av981} for $\mathcal{N}\left(  f\right)  $. Since
$x_{n}=\frac{1}{v}\Delta_{1}x^{\left(  n\right)  }$ we have from \ref{av757},%
\begin{align*}
&  \beta^{T}T^{-1}\Delta_{X}f\\
&  =\frac{1}{\left\vert T\right\vert }\sum\limits_{n=1}^{N-1}\left(
-1\right)  ^{n-1}\left(  \left(  -1\right)  ^{N}M_{n,N-1}-M_{n,1}\right)
\Delta_{1}f\left(  x^{\left(  n\right)  }\right) \\
&  =\frac{1}{\left\vert T\right\vert }\sum\limits_{n=1}^{N-1}\left(
-1\right)  ^{n-1}\left(  \left(  -1\right)  ^{N}\left(  -\frac{1}{2}\right)
^{N-1-n}\left\vert T_{1:n-1}\right\vert -\left(  -\frac{1}{2}\right)
^{n-1}\left\vert T_{n+1:N-1}\right\vert \right)  \Delta_{1}f\left(  x^{\left(
n\right)  }\right) \\
&  =\frac{1}{\left\vert T\right\vert }\sum\limits_{n=1}^{N-1}\left(  \frac
{1}{2^{N-1-n}}\left\vert T_{1:n-1}\right\vert -\frac{1}{2^{n-1}}\left\vert
T_{n+1:N-1}\right\vert \right)  \Delta_{1}f\left(  x^{\left(  n\right)
}\right) \\
&  =\frac{1}{\left\vert T\right\vert }\sum\limits_{n=1}^{N-1}\left(  \frac
{1}{2^{N-1-n}}\left\vert T_{1:n-1}\right\vert -\frac{1}{2^{n-1}}\left\vert
T_{n+1:N-1}\right\vert \right)  \Delta_{1}f\left(  x^{\left(  n\right)
}\right) \\
&  =\frac{1}{\left\vert T\right\vert }\sum\limits_{n=1}^{N-1}\left(
\frac{\Delta_{1}x^{\left(  n\right)  }}{2^{N-1-n}}\left\vert T_{1:n-1}%
\right\vert -\frac{\Delta_{1}x^{\left(  n\right)  }}{2^{n-1}}\left\vert
T_{n+1:N-1}\right\vert \right)  \frac{\Delta_{1}f\left(  x^{\left(  n\right)
}\right)  }{\Delta_{1}x^{\left(  n\right)  }}\\
&  =\frac{v}{\left\vert T\right\vert }\sum\limits_{n=1}^{N-1}\left(
\frac{x_{n}}{2^{N-1-n}}\left\vert T_{1:n-1}\right\vert -\frac{x_{n}}{2^{n-1}%
}\left\vert T_{n+1:N-1}\right\vert \right)  \frac{\Delta_{1}f\left(
x^{\left(  n\right)  }\right)  }{\Delta_{1}x^{\left(  n\right)  }},
\end{align*}

and applying \ref{av080} of Lemma \ref{vLem_Taylor_extension} gives%
\begin{align}
\left\vert \mathcal{N}\left(  f\right)  \right\vert \left\vert T\right\vert
&  \leq\left(  \left\vert f\left(  x^{\left(  1\right)  }\right)  +f\left(
x^{\left(  N\right)  }\right)  \right\vert +\frac{1}{2}\left\vert \beta
^{T}T^{-1}\Delta_{X}f\right\vert \right)  \left\vert T\right\vert \nonumber\\
&  \leq2\left\Vert f\right\Vert _{\infty;\Omega}\left\vert T\right\vert
+\frac{1}{2}\left\vert \beta^{T}T^{-1}\Delta_{X}f\right\vert \left\vert
T\right\vert \nonumber\\
&  \leq2\left\Vert f\right\Vert _{\infty;\Omega}\left\vert T\right\vert
+\frac{v}{2}\sum\limits_{n=1}^{N-1}\left(  \frac{x_{n}}{2^{N-1-n}}\left\vert
T_{1:n-1}\right\vert +\frac{x_{n}}{2^{n-1}}\left\vert T_{n+1:N-1}\right\vert
\right)  \left\vert \frac{\Delta_{1}f\left(  x^{\left(  n\right)  }\right)
}{\Delta_{1}x^{\left(  n\right)  }}\right\vert \nonumber\\
&  \leq2\left\Vert f\right\Vert _{\infty;\Omega}\left\vert T\right\vert
+\frac{v}{2}\left\Vert Df\right\Vert _{\infty;\Omega}\sum\limits_{n=1}%
^{N-1}\left(  \frac{x_{n}}{2^{N-1-n}}\left\vert T_{1:n-1}\right\vert
+\frac{x_{n}}{2^{n-1}}\left\vert T_{n+1:N-1}\right\vert \right)
.\label{av898}%
\end{align}

From \ref{av835}, \ref{av896}, when $m=N-1$ and $a=x$,
\begin{align*}
\sum\limits_{n=1}^{N-1}\frac{x_{n}}{2^{n-1}}\left\vert T_{n+1:N-1}\right\vert
& <\left\vert T\right\vert -\frac{N}{2^{N-1}},\\
\sum\limits_{n=1}^{N-1}\frac{x_{n}}{2^{N-1-n}}\left\vert T_{1:n-1}\right\vert
& <\left\vert T\right\vert -\frac{N}{2^{N-1}},
\end{align*}

so that \ref{av898} yields%
\begin{align}
\left\vert \mathcal{N}\left(  f\right)  \right\vert \left\vert T\right\vert  &
\leq2\left\Vert f\right\Vert _{\infty;\Omega}\left\vert T\right\vert +\frac
{v}{2}\left\Vert Df\right\Vert _{\infty;\Omega}\left(  \left\vert T\right\vert
-\frac{N}{2^{N-1}}+\left\vert T\right\vert -\frac{N}{2^{N-1}}\right)
\nonumber\\
& =2\left\Vert f\right\Vert _{\infty;\Omega}\left\vert T\right\vert
+v\left\Vert Df\right\Vert _{\infty;\Omega}\left(  \left\vert T\right\vert
-\frac{N}{2^{N-1}}\right)  ,\label{av940}%
\end{align}

and thus by \ref{av980},
\begin{align}
\left\vert \mu_{1}\left(  f\right)  \right\vert =\frac{\left\vert
\mathcal{N}\left(  f\right)  \right\vert \left\vert T\right\vert }%
{\mathcal{D}\left\vert T\right\vert }  & \leq\frac{2\left\Vert f\right\Vert
_{\infty;\Omega}\left\vert T\right\vert +v\left\Vert Df\right\Vert
_{\infty;\Omega}\left(  \left\vert T\right\vert -\frac{N}{2^{N-1}}\right)
}{\mathcal{D}\left\vert T\right\vert }\nonumber\\
& =2\left\Vert f\right\Vert _{\infty;\Omega}\frac{1}{\mathcal{D}}+\left\Vert
Df\right\Vert _{\infty;\Omega}\frac{v\left(  \left\vert T\right\vert -\frac
{N}{2^{N-1}}\right)  }{\mathcal{D}\left\vert T\right\vert }.\label{av944}%
\end{align}

\fbox{Component of \ref{av944}\ involving $\left\Vert f\right\Vert
_{\infty;\Omega}$} Starting with \ref{av034} we write%
\begin{align}
\mathcal{D}\left\vert T\right\vert  & >\left(  1+v\right)  \left\vert
T\right\vert +v\left(  \frac{1}{2^{N-1}}-\frac{1}{4}\left\vert T_{2:N-1}%
\right\vert -\frac{1}{4}\left\vert T_{1:N-2}\right\vert \right) \nonumber\\
& =\left\vert T\right\vert +v\left(  \left\vert T\right\vert -\frac{1}%
{4}\left\vert T_{2:N-1}\right\vert -\frac{1}{4}\left\vert T_{1:N-2}\right\vert
+\frac{1}{2^{N-1}}\right)  .\label{av945}%
\end{align}

The next result improves on \ref{av033} i.e. $\mathcal{D}>0$, by showing that
$\mathcal{D}>1$.

\begin{theorem}
\label{vThm_T_minus_quartT} $\left\vert T\left(  a\right)  \right\vert
-\frac{1}{4}\left\vert T\left(  a_{2:m}\right)  \right\vert -\frac{1}%
{4}\left\vert T\left(  a_{1:m-1}\right)  \right\vert >0$ when $a\in
\mathbb{R}_{+}^{m}$ and $m\geq2$.
\end{theorem}

\begin{proof}
From \ref{av828}: $\left(  a_{1}+\frac{1}{2}\right)  \left\vert T\left(
a_{2:m}\right)  \right\vert +\frac{1}{2^{m}}<\left\vert T\left(  a\right)
\right\vert $ for $m\geq2$. Thus
\[
\left(  a_{1}+\frac{1}{2}\right)  \left\vert T\left(  a_{2:m}\right)
\right\vert +\frac{1}{2^{m}}<\left\vert T\left(  a\right)  \right\vert ,\quad
N\geq3,
\]

i.e.%
\begin{equation}
\left\vert T\left(  a_{2:m}\right)  \right\vert <\frac{\left\vert T\left(
a\right)  \right\vert -\frac{1}{2^{m}}}{a_{1}+\frac{1}{2}}.\label{av04}%
\end{equation}

From Corollary \ref{vCor_1_formula_for_C^(m)_gam} the order reversal
permutation $\mathcal{R}$ satisfies $\left\vert T\left(  \mathcal{R}a\right)
\right\vert =\left\vert T\left(  a\right)  \right\vert $ so%
\[
\left\vert T\left(  \left(  \mathcal{R}a\right)  _{2:N-1}\right)  \right\vert
<\frac{\left\vert T\left(  \mathcal{R}a\right)  \right\vert -\frac{1}{2^{m}}%
}{a_{m}+\frac{1}{2}}=\frac{\left\vert T\left(  a\right)  \right\vert -\frac
{1}{2^{m}}}{a_{m}+\frac{1}{2}},
\]

i.e.%
\begin{equation}
\left\vert T\left(  a_{1:m-1}\right)  \right\vert <\frac{\left\vert T\left(
a\right)  \right\vert -\frac{1}{2^{m}}}{a_{1}+\frac{1}{2}}.\label{av486}%
\end{equation}

Hence%
\begin{align*}
\left\vert T\left(  a\right)  \right\vert  &  -\frac{1}{4}\left\vert T\left(
a_{2:m}\right)  \right\vert -\frac{1}{4}\left\vert T\left(  a_{1:m-1}\right)
\right\vert \\
&  >\left\vert T\left(  a\right)  \right\vert -\frac{1}{4}\left(
\frac{\left\vert T\left(  a\right)  \right\vert -\frac{1}{2^{m}}}{a_{1}%
+\frac{1}{2}}+\frac{\left\vert T\left(  a\right)  \right\vert -\frac{1}{2^{m}%
}}{a_{m}+\frac{1}{2}}\right) \\
&  =\left\vert T\left(  a\right)  \right\vert -\frac{1}{4}\left(  \left\vert
T\left(  a\right)  \right\vert -\frac{1}{2^{m}}\right)  \left(  \frac{1}%
{a_{1}+\frac{1}{2}}+\frac{1}{a_{m}+\frac{1}{2}}\right) \\
&  =\left\vert T\left(  a\right)  \right\vert -\frac{1}{4}\left\vert T\left(
a\right)  \right\vert \left(  \frac{1}{a_{1}+\frac{1}{2}}+\frac{1}{a_{m}%
+\frac{1}{2}}\right)  +\frac{1}{2^{m+2}}\left(  \frac{1}{a_{1}+\frac{1}{2}%
}+\frac{1}{a_{m}+\frac{1}{2}}\right) \\
&  =\left\vert T\left(  a\right)  \right\vert -\frac{1}{4}\left\vert T\left(
a\right)  \right\vert \left(  \frac{1}{a_{1}+\frac{1}{2}}+\frac{1}{a_{m}%
+\frac{1}{2}}\right)  +\frac{1}{2^{m+2}}\left(  \frac{1}{a_{1}+\frac{1}{2}%
}+\frac{1}{a_{m}+\frac{1}{2}}\right) \\
&  =\left\vert T\left(  a\right)  \right\vert -\frac{1}{4}\left\vert T\left(
a\right)  \right\vert \left(  \frac{1}{a_{1}+\frac{1}{2}}+\frac{1}{a_{m}%
+\frac{1}{2}}\right)  +\frac{1}{2^{m+2}}\left(  \frac{1}{a_{1}+\frac{1}{2}%
}+\frac{1}{a_{m}+\frac{1}{2}}\right) \\
&  =\left(  1-\frac{1}{4}\left(  \frac{1}{a_{1}+\frac{1}{2}}+\frac{1}%
{a_{m}+\frac{1}{2}}\right)  \right)  \left\vert T\left(  a\right)  \right\vert
+\frac{1}{2^{m+2}}\left(  \frac{1}{a_{1}+\frac{1}{2}}+\frac{1}{a_{m}+\frac
{1}{2}}\right) \\
&  \geq\left(  1-\frac{1}{4}\left(  \frac{1}{\frac{1}{2}}+\frac{1}{\frac{1}%
{2}}\right)  \right)  \left\vert T\left(  a\right)  \right\vert +\frac
{1}{2^{m+2}}\left(  \frac{1}{a_{1}+\frac{1}{2}}+\frac{1}{a_{m}+\frac{1}{2}%
}\right) \\
&  =\frac{1}{2^{m+2}}\left(  \frac{1}{a_{1}+\frac{1}{2}}+\frac{1}{a_{m}%
+\frac{1}{2}}\right) \\
&  >0.
\end{align*}

\end{proof}

Since $\mathcal{D}>1$, \ref{av944} becomes%
\begin{align}
\left\vert \mu_{1}\left(  f\right)  \right\vert  & \leq2\left\Vert
f\right\Vert _{\infty;\Omega}+\frac{v\left(  \left\vert T\right\vert -\frac
{N}{2^{N-1}}\right)  }{\mathcal{D}\left\vert T\right\vert }\left\Vert
Df\right\Vert _{\infty;\Omega}\label{av946}\\
& <2\left\Vert f\right\Vert _{\infty;\Omega}+v\frac{\left\vert T\right\vert
-\frac{N}{2^{N-1}}}{\left\vert T\right\vert }\left\Vert Df\right\Vert
_{\infty;\Omega}\nonumber\\
& <2\left\Vert f\right\Vert _{\infty;\Omega}+v\left\Vert Df\right\Vert
_{\infty;\Omega}.\label{av086}%
\end{align}

\fbox{Component of \ref{av946}\ involving $\left\Vert Df\right\Vert
_{\infty,\Omega}$} From \ref{av943},%
\[
\frac{v\left(  \left\vert T\right\vert -\frac{N}{2^{N-1}}\right)
}{\mathcal{D}\left\vert T\right\vert }=\frac{\left\vert T\right\vert -\frac
{N}{2^{N-1}}}{\frac{1}{v}\mathcal{D}\left\vert T\right\vert }<\frac{\left\vert
T\right\vert -\frac{N}{2^{N-1}}}{\left(  x_{1}+\frac{1}{2}\right)  \left\vert
T_{2:N-1}\right\vert +\left(  x_{N-1}+\frac{1}{2}\right)  \left\vert
T_{1:N-2}\right\vert }.
\]

From \ref{av936} and \ref{av840},%
\begin{align*}
\left\vert T\right\vert -\frac{N}{2^{N-1}}  & =x_{1}\left\vert T\left(
x_{2:N-1}\right)  \right\vert +\left\vert T\left(  0,x_{2:N-1}\right)
\right\vert -\frac{N}{2^{N-1}},\\
\left\vert T\right\vert -\frac{N}{2^{N-1}}  & =x_{N}\left\vert T\left(
x_{1:N-2}\right)  \right\vert +\left\vert T\left(  x_{1:N-2},0\right)
\right\vert -\frac{N}{2^{N-1}},
\end{align*}

so that%
\begin{align*}
&  \frac{v\left(  \left\vert T\right\vert -\frac{N}{2^{N-1}}\right)
}{\mathcal{D}\left\vert T\right\vert }\\
&  <\frac{x_{1}\left\vert T\left(  x_{2:N-1}\right)  \right\vert +\left\vert
T\left(  0,x_{2:N-1}\right)  \right\vert -\frac{N}{2^{N-1}}+x_{N}\left\vert
T\left(  x_{1:N-2}\right)  \right\vert +\left\vert T\left(  x_{1:N-2}%
,0\right)  \right\vert -\frac{N}{2^{N-1}}}{\left(  x_{1}+\frac{1}{2}\right)
\left\vert T_{2:N-1}\right\vert +\left(  x_{N-1}+\frac{1}{2}\right)
\left\vert T_{1:N-2}\right\vert }\\
& \\
&  =\frac{%
\begin{array}
[c]{r}%
x_{1}\left\vert T\left(  x_{2:N-1}\right)  \right\vert +x_{N}\left\vert
T\left(  x_{1:N-2}\right)  \right\vert +\left\vert T\left(  0,x_{2:N-1}%
\right)  \right\vert -\left\vert T\left(  0_{N-1}\right)  \right\vert
+\left\vert T\left(  x_{1:N-2},0\right)  \right\vert -\\
-\left\vert T\left(  0_{N-1}\right)  \right\vert
\end{array}
}{\left(  x_{1}+\frac{1}{2}\right)  \left\vert T_{2:N-1}\right\vert +\left(
x_{N-1}+\frac{1}{2}\right)  \left\vert T_{1:N-2}\right\vert }\\
& \\
&  =\frac{x_{1}\left\vert T\left(  x_{2:N-1}\right)  \right\vert
+x_{N}\left\vert T\left(  x_{1:N-2}\right)  \right\vert }{\left(  x_{1}%
+\frac{1}{2}\right)  \left\vert T_{2:N-1}\right\vert +\left(  x_{N-1}+\frac
{1}{2}\right)  \left\vert T_{1:N-2}\right\vert }+\\
&  \qquad\qquad\qquad\qquad\qquad\qquad+\frac{\left\vert T\left(
0,x_{2:N-1}\right)  \right\vert -\left\vert T\left(  0_{N-1}\right)
\right\vert +\left\vert T\left(  x_{1:N-2},0\right)  \right\vert -\left\vert
T\left(  0_{N-1}\right)  \right\vert }{\left(  x_{1}+\frac{1}{2}\right)
\left\vert T_{2:N-1}\right\vert +\left(  x_{N-1}+\frac{1}{2}\right)
\left\vert T_{1:N-2}\right\vert }\\
& \\
&  <1+\frac{\left\vert T\left(  0,x_{2:N-1}\right)  \right\vert -\left\vert
T\left(  0_{N-1}\right)  \right\vert +\left\vert T\left(  x_{1:N-2},0\right)
\right\vert -\left\vert T\left(  0_{N-1}\right)  \right\vert }{\left(
x_{1}+\frac{1}{2}\right)  \left\vert T_{2:N-1}\right\vert +\left(
x_{N-1}+\frac{1}{2}\right)  \left\vert T_{1:N-2}\right\vert }.
\end{align*}

But from \ref{av935}, $0<\left\vert T\left(  0_{1},u\right)  \right\vert
-\left\vert T\left(  0_{m+1}\right)  \right\vert \leq\left\vert T\left(
u\right)  \right\vert -\left\vert T\left(  0_{m}\right)  \right\vert
<\left\vert T\left(  u\right)  \right\vert $ when $u\in\mathbb{R}_{\oplus}%
^{m}$ and $m\geq1$, so that%
\begin{align*}
\frac{v\left(  \left\vert T\right\vert -\frac{N}{2^{N-1}}\right)
}{\mathcal{D}\left\vert T\right\vert } &  <1+\frac{\left\vert T\left(
0,x_{2:N-1}\right)  \right\vert -\left\vert T\left(  0_{N-1}\right)
\right\vert +\left\vert T\left(  x_{1:N-2},0\right)  \right\vert -\left\vert
T\left(  0_{N-1}\right)  \right\vert }{\left(  x_{1}+\frac{1}{2}\right)
\left\vert T_{2:N-1}\right\vert +\left(  x_{N-1}+\frac{1}{2}\right)
\left\vert T_{1:N-2}\right\vert }\\
&  <1+\frac{\left\vert T_{2:N-1}\right\vert +\left\vert T_{1:N-2}\right\vert
}{\left(  x_{1}+\frac{1}{2}\right)  \left\vert T_{2:N-1}\right\vert +\left(
x_{N-1}+\frac{1}{2}\right)  \left\vert T_{1:N-2}\right\vert }\\
&  <1+\frac{\left\vert T_{2:N-1}\right\vert +\left\vert T_{1:N-2}\right\vert
}{\frac{1}{2}\left\vert T_{2:N-1}\right\vert +\frac{1}{2}\left\vert
T_{1:N-2}\right\vert }\\
&  =3,
\end{align*}

and the estimate \ref{av946} becomes%
\begin{equation}
\left\vert \mu_{1}\left(  f\right)  \right\vert \leq2\left\Vert f\right\Vert
_{\infty;\Omega}+3\left\Vert Df\right\Vert _{\infty;\Omega},\quad
N\geq3,\text{ }\rho>0.\label{av302}%
\end{equation}

Combining this bound with that of \ref{av086} we\ get%
\[
\left\vert \mu_{1}\left(  f\right)  \right\vert \leq2\left\Vert f\right\Vert
_{\infty;\Omega}+\min\left\{  3,\rho N\right\}  \left\Vert Df\right\Vert
_{\infty;\Omega},\quad N\geq3,\text{ }\rho>0,
\]

so that from \ref{av081},%
\begin{align}
\left\vert \mu\left(  f\right)  \right\vert _{\max}  & \leq2\left\Vert
Df\right\Vert _{\infty;\Omega}+\left\vert \mu_{1}\left(  f\right)  \right\vert
\nonumber\\
& \leq2\left\Vert f\right\Vert _{\infty;\Omega}+\min\left\{  5,2+\rho
N\right\}  \left\Vert Df\right\Vert _{\infty;\Omega},\quad N\geq4,\text{ }%
\rho>0.\label{av036}%
\end{align}

\subsection{Upper bounds for the Exact smoother derivative}

This result is the smoother equivalent of the interpolation derivative
estimate given in Theorem \ref{vThm_bound_deriv_hat_interpol}:

\begin{theorem}
\label{vThm_bound_deriv_hat_smth}\textbf{Hat basis function smoother} Suppose
the data region satisfies $\operatorname*{diam}\Omega\leq1$ and the data
function $f$ lies in the space $\left\{  f\in C_{B}^{\left(  0\right)
}\left(  \Omega\right)  :Df\in L^{\infty}\left(  \Omega\right)  \right\}  $.
Then the basis function Exact smoother $s$ corresponding to the hat basis
function $\Lambda$ satisfies%
\[
\left\Vert Ds\right\Vert _{\infty;\Omega}\leq2\left\Vert f\right\Vert
_{\infty;\Omega}+\min\left\{  5,2+\rho N\right\}  \left\Vert Df\right\Vert
_{\infty;\Omega},\quad N\geq4,\text{ }\rho>0.
\]

\end{theorem}

\begin{proof}
From \ref{av089} $\left\Vert Ds\right\Vert _{\infty;\Omega}=\left\vert
\mu\left(  f\right)  \right\vert _{\max}$ so \ref{av036} gives the result.
\end{proof}

In the interpolation case $\rho\rightarrow0^{+}$ yielding $\left\vert
\mu\left(  f\right)  \right\vert _{\max}\leq2\left\Vert f\right\Vert
_{\infty;\Omega}+2\left\Vert Df\right\Vert _{\infty;\Omega}$ which is a weaker
estimate than that of \ref{av088} of Theorem
\ref{vThm_bound_deriv_hat_interpol}.

The next result is the smoother equivalent of the interpolation estimate
derived in Corollary \ref{vCor_bound_deriv_scal_hat_interpol}:

\begin{corollary}
\label{vCor_Thm_bound_deriv_scal_hat_smth}\textbf{Scaled hat basis function
smoother }Suppose the data region satisfies $\operatorname*{diam}\Omega
\leq\lambda$ and the data function $f$ lies in the space $\left\{  f\in
C_{B}^{\left(  0\right)  }\left(  \Omega\right)  :Df\in L^{\infty}\left(
\Omega\right)  \right\}  $.

Then the basis function smoother $s$ corresponding to the scaled hat basis
function $\Lambda\left(  \cdot/\lambda\right)  $ satisfies%
\[
\left\Vert Ds\right\Vert _{\infty;\Omega}\leq\frac{2}{\lambda}\left\Vert
f\right\Vert _{\infty;\Omega}+\min\left\{  5,2+\rho N\right\}  \left\Vert
Df\right\Vert _{\infty;\Omega}\quad N\geq4,\text{ }\rho>0.
\]

\end{corollary}

\begin{proof}
From \ref{av084} and \ref{av083} of the introduction to this chapter it
follows that $s_{1/\lambda}\left(  x\right)  =s\left(  \lambda x\right)  $ is
the hat basis function smoother associated with data function $f_{1/\lambda
}\left(  x\right)  =f\left(  \lambda x\right)  $ and data $X_{\lambda}%
\subset\Omega_{\lambda}$ satisfying $\operatorname*{diam}X_{\lambda}\leq1$. By
Theorem \ref{vThm_bound_deriv_hat_smth},%
\begin{align*}
\left\Vert Ds_{1/\lambda}\right\Vert _{\infty;\Omega_{\lambda}}  &
\leq2\left\Vert f_{1/\lambda}\right\Vert _{\infty;\Omega_{\lambda}}%
+\min\left\{  5,2+\rho N\right\}  \left\Vert Df_{1/\lambda}\right\Vert
_{\infty;\Omega_{\lambda}}\\
& =2\left\Vert f\right\Vert _{\infty;\Omega}+\min\left\{  5,2+\rho N\right\}
\left\Vert Df_{1/\lambda}\right\Vert _{\infty;\Omega_{\lambda}}\\
& =2\left\Vert f\right\Vert _{\infty;\Omega}+\lambda\min\left\{  5,2+\rho
N\right\}  \left\Vert \left(  Df\right)  _{1/\lambda}\right\Vert
_{\infty;\Omega_{\lambda}}\\
& =2\left\Vert f\right\Vert _{\infty;\Omega}+\lambda\min\left\{  5,2+\rho
N\right\}  \left\Vert Df\right\Vert _{\infty;\Omega},
\end{align*}

and since%
\[
\left\Vert Ds_{1/\lambda}\right\Vert _{\infty;\Omega_{\lambda}}=\lambda
\left\Vert \left(  Ds\right)  _{1/\lambda}\right\Vert _{\infty;\Omega
_{\lambda}}=\lambda\left\Vert Ds\right\Vert _{\infty;\Omega},
\]

it follows that%
\[
\left\Vert Ds\right\Vert _{\infty;\Omega}\leq\frac{2}{\lambda}\left\Vert
f\right\Vert _{\infty;\Omega}+\min\left\{  5,2+\rho N\right\}  \left\Vert
Df\right\Vert _{\infty;\Omega}.
\]

\end{proof}

\section{Some unused results}

These results are interesting but have not been used in this document.

\subsection{A Taylor series formula for $\left\vert T\left(  \cdot\right)
\right\vert $}

We start by observing that \ref{av596} implies that
\begin{equation}
c_{\beta}^{\left(  m\right)  }=D^{\beta}\left\vert T\left(  a_{1:m}\right)
\right\vert \left(  0_{m}\right)  ,\quad0\leq\beta\leq1_{m}.\label{av911}%
\end{equation}

Our approach will be to derive an explicit formula for the derivatives of the
determinate $\left\vert T\left(  a\right)  \right\vert $.

\begin{example}
When $a\in\mathbb{R}^{6}$, differentiating gives%
\[
\left\vert T\left(  a\right)  \right\vert =%
\begin{vmatrix}
1+a_{1} & -\frac{1}{2} & 0 & 0 & 0 & 0\\
-\frac{1}{2} & 1+a_{2} & -\frac{1}{2} & 0 & 0 & 0\\
0 & -\frac{1}{2} & 1+a_{3} & -\frac{1}{2} & 0 & 0\\
0 & 0 & -\frac{1}{2} & 1+a_{4} & -\frac{1}{2} & 0\\
0 & 0 & 0 & -\frac{1}{2} & 1+a_{5} & -\frac{1}{2}\\
0 & 0 & 0 & 0 & -\frac{1}{2} & 1+a_{6}%
\end{vmatrix}
.
\]

Thus%
\begin{align*}
D_{1}\left\vert T\left(  a\right)  \right\vert  & =%
\begin{vmatrix}
1 & -\frac{1}{2} & 0 & 0 & 0 & 0\\
0 & 1+a_{2} & -\frac{1}{2} & 0 & 0 & 0\\
0 & -\frac{1}{2} & 1+a_{3} & -\frac{1}{2} & 0 & 0\\
0 & 0 & -\frac{1}{2} & 1+a_{4} & -\frac{1}{2} & 0\\
0 & 0 & 0 & -\frac{1}{2} & 1+a_{5} & -\frac{1}{2}\\
0 & 0 & 0 & 0 & -\frac{1}{2} & 1+a_{6}%
\end{vmatrix}
=\left\vert T\left(  a_{2:6}\right)  \right\vert ,\\
& \\
D_{6}\left\vert T\left(  a\right)  \right\vert  & =%
\begin{vmatrix}
1+a_{1} & -\frac{1}{2} & 0 & 0 & 0 & 0\\
-\frac{1}{2} & 1+a_{2} & -\frac{1}{2} & 0 & 0 & 0\\
0 & -\frac{1}{2} & 1+a_{3} & -\frac{1}{2} & 0 & 0\\
0 & 0 & -\frac{1}{2} & 1+a_{4} & -\frac{1}{2} & 0\\
0 & 0 & 0 & -\frac{1}{2} & 1+a_{5} & 0\\
0 & 0 & 0 & 0 & 0 & 1
\end{vmatrix}
=\left\vert T\left(  a_{1:5}\right)  \right\vert ,\\
& \\
D_{4}\left\vert T\left(  a\right)  \right\vert  & =%
\begin{vmatrix}
1+a_{1} & -\frac{1}{2} & 0 & 0 & 0 & 0\\
-\frac{1}{2} & 1+a_{2} & -\frac{1}{2} & 0 & 0 & 0\\
0 & -\frac{1}{2} & 1+a_{3} & 0 & 0 & 0\\
0 & 0 & 0 & 1 & -\frac{1}{2} & 0\\
0 & 0 & 0 & 0 & 1+a_{5} & -\frac{1}{2}\\
0 & 0 & 0 & 0 & -\frac{1}{2} & 1+a_{6}%
\end{vmatrix}
=\left\vert T\left(  a_{1:3}\right)  \right\vert \text{ }\left\vert T\left(
a_{5:6}\right)  \right\vert .
\end{align*}

\end{example}

In general, if $m\geq3$,%
\[
D_{k}\left\vert T\left(  a_{1:m}\right)  \right\vert =\left\{
\begin{array}
[c]{ll}%
\left\vert T\left(  a_{2:m}\right)  \right\vert , & k=1,\\
\left\vert T\left(  a_{1:k-1}\right)  \right\vert \text{ }\left\vert T\left(
a_{k+1:m}\right)  \right\vert , & k=2,3,\ldots,m-1,\\
\left\vert T\left(  a_{1:m-1}\right)  \right\vert , & k=m,
\end{array}
\right.
\]

so that%
\begin{align*}
D_{k}\left\vert T\left(  \cdot\right)  \right\vert \left(  0\right)   &
=\left\{
\begin{array}
[c]{ll}%
\left\vert T\left(  0_{m-1}\right)  \right\vert , & k=1,\\
\left\vert T\left(  0_{k-1}\right)  \right\vert \text{ }\left\vert T\left(
0_{m-k}\right)  \right\vert , & k=2,3,\ldots,m-1,\\
\left\vert T\left(  0_{m-1}\right)  \right\vert , & k=m,
\end{array}
\right. \\
& =\left\{
\begin{array}
[c]{ll}%
\left\vert T\left(  0_{m-1}\right)  \right\vert , & k=1,m,\\
\left\vert T\left(  0_{k-1}\right)  \right\vert \text{ }\left\vert T\left(
0_{m-k}\right)  \right\vert , & k=2,3,\ldots,m-1.
\end{array}
\right.
\end{align*}

In general, if $2\leq k_{1}<k_{2}<\ldots<k_{n}\leq m-1$, then%
\begin{equation}
D_{k_{1}}D_{k_{2}\ldots}D_{k_{n}}\left\vert T\left(  a\right)  \right\vert
=\left\{
\begin{array}
[c]{ll}%
\left\vert T\left(  a_{1:k_{1}-1}\right)  \right\vert \text{ }\left\vert
T\left(  a_{k_{n}+1:m}\right)  \right\vert , & n=1,\\
\left\vert T\left(  a_{1:k_{1}-1}\right)  \right\vert \left(  \prod
\limits_{i=1}^{n-1}\left\vert T\left(  a_{k_{i}+1:k_{i+1}-1}\right)
\right\vert \right)  \left\vert T\left(  a_{k_{n}+1:m}\right)  \right\vert , &
2\leq n\leq m,
\end{array}
\right. \label{av910}%
\end{equation}

so that%
\[
D_{k_{1}}D_{k_{2}\ldots}D_{k_{n}}\left\vert T\left(  \cdot\right)  \right\vert
\left(  0_{m}\right)  =\left\{
\begin{array}
[c]{ll}%
\left\vert T\left(  0_{k_{1}-1}\right)  \right\vert \text{ }\left\vert
T\left(  0_{m-k_{n}}\right)  \right\vert , & n=1,\\
\left\vert T\left(  0_{k_{1}-1}\right)  \right\vert \left(  \prod
\limits_{i=1}^{n-1}\left\vert T\left(  0_{k_{i+1}-k_{i}}\right)  \right\vert
\right)  \left\vert T\left(  0_{m-k_{n}}\right)  \right\vert , & 2\leq n\leq
m.
\end{array}
\right.  .
\]

But from \ref{av598}, $T\left(  0_{j}\right)  =\frac{j+1}{2^{j}}$, $j\geq1$,
so that if $n\geq2$,
\begin{align}
D_{k_{1}}D_{k_{2}\ldots}D_{k_{n}}\left\vert T\left(  \cdot\right)  \right\vert
\left(  0_{m}\right)   & =\left\vert T\left(  0_{k_{1}-1}\right)  \right\vert
\left(  \prod\limits_{i=1}^{n-1}\left\vert T\left(  0_{k_{i+1}-k_{i}}\right)
\right\vert \right)  \left\vert T\left(  0_{m-k_{n}}\right)  \right\vert
\nonumber\\
& =\frac{1}{2^{k_{1}-1}}k_{1}\left(  \prod\limits_{i=1}^{n-1}\frac
{k_{i+1}-k_{i}+1}{2^{k_{i+1}-k_{i}}}\right)  \frac{m-k_{n}+1}{2^{m-k_{n}}%
}\nonumber\\
& =\frac{1}{2^{m-1}}k_{1}\left(  \prod\limits_{i=1}^{n-1}\left(  k_{i+1}%
-k_{i}+1\right)  \right)  \left(  m-k_{n}+1\right)  ,\label{av909}%
\end{align}

and if $n=1$,%
\[
D_{k_{1}}D_{k_{2}\ldots}D_{k_{n}}\left\vert T\left(  \cdot\right)  \right\vert
\left(  0_{m}\right)  =\frac{1}{2^{m-1}}k_{1}\left(  m-k_{n}+1\right)  ,
\]

so we have shown that%
\begin{equation}
D_{k_{1}}D_{k_{2}\ldots}D_{k_{n}}\left\vert T\left(  \cdot\right)  \right\vert
\left(  0_{m}\right)  =\left\{
\begin{array}
[c]{ll}%
\frac{1}{2^{m-1}}k_{1}\left(  m-k_{1}+1\right)  , & n=1,\\
\frac{1}{2^{m-1}}k_{1}\left(  \prod\limits_{i=1}^{n-1}\left(  k_{i+1}%
-k_{i}+1\right)  \right)  \left(  m-k_{n}+1\right)  , & 2\leq n\leq m.
\end{array}
\right. \label{av916}%
\end{equation}

We can now use \ref{av911} and \ref{av916} to get an explicit formula for the
$c_{\beta}^{\left(  m\right)  }$.

\subsection{Another lower bound for $\left\vert T\left(  \cdot\right)
\right\vert $}

Lower bounds for $\left\vert T\left(  a_{1:m}\right)  \right\vert $ were
obtained in \ref{av835} and \ref{av896} but these were in terms of the form
$\left\vert T\left(  a_{k+1:m}\right)  \right\vert $ and $\left\vert T\left(
a_{1:j-1}\right)  \right\vert $. Here we obtain the polynomial estimates
\ref{av819}.

Suppose $a\in\mathbb{R}_{+}^{m}$. From \ref{av731},%
\[
c_{\gamma}^{\left(  m\right)  }\geq\frac{1+m-\left\vert \gamma\right\vert
}{2^{m-\left\vert \gamma\right\vert }}\geq\frac{1+m}{2^{m}}=c_{\mathbf{0}%
}^{\left(  m\right)  },
\]

and consequently%
\begin{align}
\sum\limits_{\substack{\gamma\leq\mathbf{1} \\\left\vert \gamma\right\vert
=k}}c_{\gamma}^{\left(  m\right)  }a^{\mathbf{\gamma}}\geq\sum
\limits_{\substack{\gamma\leq\mathbf{1} \\\left\vert \gamma\right\vert
=k}}\frac{1+m-\left\vert \gamma\right\vert }{2^{m-\left\vert \gamma\right\vert
}}a^{\mathbf{\gamma}}=\sum\limits_{\substack{\gamma\leq\mathbf{1} \\\left\vert
\gamma\right\vert =k}}\frac{1+m-k}{2^{m-k}}a^{\mathbf{\gamma}}  &
=\frac{1+m-k}{2^{m-k}}\sum\limits_{\substack{\gamma\leq\mathbf{1} \\\left\vert
\gamma\right\vert =k}}a^{\mathbf{\gamma}}\nonumber\\
& =\frac{1+m-k}{2^{m-k}}\sigma_{k}\left(  a\right)  ,\label{av817}%
\end{align}

so that%
\begin{align}
\left\vert T\left(  a\right)  \right\vert =\sum\limits_{k=0}^{m}%
\sum\limits_{\substack{\gamma\leq\mathbf{1} \\\left\vert \gamma\right\vert
=k}}c_{\gamma}^{\left(  m\right)  }a^{\mathbf{\gamma}}  & \geq\sum
\limits_{k=0}^{m}\frac{1+m-k}{2^{m-k}}\sigma_{k}\left(  a\right) \nonumber\\
& =\sum\limits_{k=0}^{m}\frac{1}{2^{m-k}}\sigma_{k}\left(  a\right)
+\sum\limits_{k=0}^{m-1}\frac{m-k}{2^{m-k}}\sigma_{k}\left(  a\right)
.\label{av810}%
\end{align}

But if $s\in\mathbb{R}^{1}$,%
\begin{align}
\sigma_{m}\left(  s\mathbf{1}+a\right)   & =\left(  s\mathbf{1}+a\right)
^{\mathbf{1}}=\sum\limits_{\delta\leq\mathbf{1}}\tbinom{\mathbf{1}}{\delta
}\left(  s\mathbf{1}\right)  ^{\mathbf{1}-\delta}a^{\delta}=\sum
\limits_{\delta\leq\mathbf{1}}s^{m-\left\vert \delta\right\vert }a^{\delta
}=\nonumber\\
& =\sum\limits_{k=0}^{m}\sum\limits_{\substack{\delta\leq\mathbf{1}
\\\left\vert \delta\right\vert =k}}s^{m-k}a^{\delta}=\sum\limits_{k=0}%
^{m}s^{m-k}\sum\limits_{\substack{\delta\leq\mathbf{1} \\\left\vert
\delta\right\vert =k}}a^{\delta}=\sum\limits_{k=0}^{m}s^{m-k}\sigma_{k}\left(
a\right)  ,\label{av811}%
\end{align}

and%
\begin{equation}
D_{s}\left(  s\mathbf{1}+a\right)  ^{\mathbf{1}}=\sum\limits_{k=0}%
^{m-1}\left(  m-k\right)  s^{m-k-1}\sigma_{k}\left(  a\right)  =\frac{1}%
{s}\sum\limits_{k=0}^{m-1}\left(  m-k\right)  s^{m-k}\sigma_{k}\left(
a\right)  ,\label{av818}%
\end{equation}

so that \ref{av810} becomes%
\[
\left\vert T\left(  a\right)  \right\vert \geq\left(  \frac{1}{2v}%
\mathbf{1}+a\right)  ^{\mathbf{1}}+\frac{1}{2}\left(  D_{s}\left(
s\mathbf{1}+a\right)  ^{\mathbf{1}}\right)  \left(  s=\frac{1}{2}\right)  .
\]

However%
\begin{align*}
D_{s}\left(  s\mathbf{1}+a\right)  ^{\mathbf{1}}  & =D_{s}\left\{  \left(
s+a_{1}\right)  \left(  s+a_{2}\right)  \ldots\left(  s+a_{m}\right)  \right\}
\\
& =\sigma_{m-1}\left(  s\mathbf{1}+a\right) \\
& =\left(  s\mathbf{1}+a\right)  ^{\mathbf{1}}\left(  \frac{1}{s+a_{1}}%
+\ldots+\frac{1}{s+a_{m}}\right)  ,
\end{align*}

so%
\begin{align*}
\left\vert T\left(  a\right)  \right\vert  & \geq\left(  \frac{1}{2}+a\right)
^{\mathbf{1}}+\frac{1}{2}\sigma_{m-1}\left(  \frac{1}{2v}+a\right) \\
& =\left(  \frac{1}{2}+a\right)  ^{\mathbf{1}}+\frac{1}{2}\left(  \frac{1}%
{2}1+a\right)  ^{\mathbf{1}_{m}}\left(  \frac{1}{\frac{1}{2}+a_{1}}%
+\ldots+\frac{1}{\frac{1}{2}+a_{m}}\right) \\
& =\left(  \frac{1}{2}+a\right)  ^{\mathbf{1}}\left(  1+\frac{1}{2}\left(
\frac{1}{\frac{1}{2}+a_{1}}+\ldots+\frac{1}{\frac{1}{2}+a_{m}}\right)
\right)  ,
\end{align*}

and we have thus proved the polynomial lower bound%
\begin{equation}
\left\vert T\left(  a\right)  \right\vert \geq\left(  \frac{1}{2}+a\right)
^{\mathbf{1}}\left(  1+\frac{1}{2}\left(  \frac{1}{\frac{1}{2}+a_{1}}%
+\ldots+\frac{1}{\frac{1}{2}+a_{m}}\right)  \right)  .\label{av819}%
\end{equation}

When $a=0$ this inequality is exact since the right side becomes $\frac
{m+1}{2^{m}}=c_{\mathbf{0}}^{\left(  m\right)  }=T\left(  0\right)  $.

\subsection{Two upper bounds for $\left\vert T\left(  a;v\right)  \right\vert
$}

In this section we prove the polynomial upper bounds \ref{av05} and
\ref{av051} for $\left\vert T\left(  a\right)  \right\vert =\left\vert
T\left(  a;1\right)  \right\vert $.

Compare the following lemma with\ \ref{av004}.

\begin{lemma}
\label{vLem_T(s1m)}If $s\neq0$ then%
\begin{equation}
T\left(  s1_{m}\right)  =\frac{\beta^{m+1}-\alpha^{m+1}}{\beta-\alpha
},\label{av049}%
\end{equation}

where $\alpha=\frac{1+s-\sqrt{2s+s^{2}}}{2}$ and $\beta=\frac{1+s+\sqrt
{2s+s^{2}}}{2}$.
\end{lemma}

\begin{proof}
From \ref{av585} and \ref{av586} we have the general recurrence relation for
$\left\vert T\left(  a_{1:m}\right)  \right\vert $:%
\[%
\begin{array}
[c]{ll}%
\left\vert T\left(  a_{1:m}\right)  \right\vert =\left(  1+a_{m}\right)
\left\vert T\left(  a_{1:m-1}\right)  \right\vert -\frac{1}{4}\left\vert
T\left(  a_{1:m-2}\right)  \right\vert , & m\geq3,\\
\left\vert T\left(  a_{1:2}\right)  \right\vert =\left(  1+a_{1}\right)
\left(  1+a_{2}\right)  -\frac{1}{4}, & m=2,\\
\left\vert T\left(  a_{1}\right)  \right\vert =1+a_{1}, & m=1,
\end{array}
\]

and in particular when $a=s\mathbf{1}$ we have the simple linear difference
equation for $\left\vert T\left(  s1_{m}\right)  \right\vert $:%
\[%
\begin{array}
[c]{ll}%
\left\vert T\left(  s1_{m}\right)  \right\vert =\left(  1+s\right)  \left\vert
T\left(  s1_{m-1}\right)  \right\vert -\frac{1}{4}\left\vert T\left(
s1_{m-2}\right)  \right\vert , & m\geq3,\\
\left\vert T\left(  s1_{2}\right)  \right\vert =\left(  1+s\right)  ^{2}%
-\frac{1}{4}, & m=2,\\
\left\vert T\left(  s1\right)  \right\vert =1+s, & m=1.
\end{array}
\]

The auxiliary equation is $x^{2}-\left(  1+s\right)  x+\frac{1}{4}=0$ which
has the two roots $\alpha=\frac{1+s-\sqrt{\left(  1+s\right)  ^{2}-1}}%
{2}=\frac{1+s-\sqrt{2s+s^{2}}}{2}$ and $\beta=\frac{1+s+\sqrt{2s+s^{2}}}{2}$
which satisfy $\alpha\beta=\frac{1}{4}$, $\alpha+\beta=1+s$ and $\beta
-\alpha=\sqrt{2s+s^{2}}$.

Since the solution now has the form $A\alpha^{m}+B\beta^{m}$ for $m\geq1$ the
initial conditions imply $A\alpha+B\beta=1+s$ and $A\alpha^{2}+B\beta
^{2}=\left(  1+s\right)  ^{2}-\frac{1}{4}$ i.e.
\[%
\begin{pmatrix}
\alpha & \beta\\
\alpha^{2} & \beta^{2}%
\end{pmatrix}%
\begin{pmatrix}
A\\
B
\end{pmatrix}
=%
\begin{pmatrix}
1+s\\
\left(  1+s\right)  ^{2}-\frac{1}{4}%
\end{pmatrix}
=%
\begin{pmatrix}
\alpha+\beta\\
\left(  \alpha+\beta\right)  ^{2}-\frac{1}{4}%
\end{pmatrix}
=%
\begin{pmatrix}
\alpha+\beta\\
\alpha^{2}+\beta^{2}+\frac{1}{4}%
\end{pmatrix}
,
\]

so that%
\begin{align*}%
\begin{pmatrix}
A\\
B
\end{pmatrix}
& =%
\begin{pmatrix}
\alpha & \beta\\
\alpha^{2} & \beta^{2}%
\end{pmatrix}
^{-1}%
\begin{pmatrix}
\alpha+\beta\\
\alpha^{2}+\beta^{2}+\frac{1}{4}%
\end{pmatrix}
\\
& =\frac{1}{\alpha\beta\left(  \beta-\alpha\right)  }%
\begin{pmatrix}
\beta^{2} & -\beta\\
-\alpha^{2} & \alpha
\end{pmatrix}%
\begin{pmatrix}
\alpha+\beta\\
\alpha^{2}+\beta^{2}+\frac{1}{4}%
\end{pmatrix}
\\
& =\frac{4}{\beta-\alpha}%
\begin{pmatrix}
\beta^{2} & -\beta\\
-\alpha^{2} & \alpha
\end{pmatrix}%
\begin{pmatrix}
\alpha+\beta\\
\alpha^{2}+\beta^{2}+\frac{1}{4}%
\end{pmatrix}
\\
& =\frac{4}{\beta-\alpha}%
\begin{pmatrix}
\alpha\beta^{2}+\beta^{3}-\alpha^{2}\beta-\beta^{3}-\frac{1}{4}\beta\\
-\alpha^{3}-\alpha^{2}\beta+\alpha^{3}+\alpha\beta^{2}+\frac{1}{4}\alpha
\end{pmatrix}
\\
& =\frac{4}{\beta-\alpha}%
\begin{pmatrix}
\alpha\beta^{2}-\alpha^{2}\beta-\frac{1}{4}\beta\\
-\alpha^{2}\beta+\alpha\beta^{2}+\frac{1}{4}\alpha
\end{pmatrix}
\\
& =\frac{4}{\beta-\alpha}%
\begin{pmatrix}
\frac{1}{4}\beta-\frac{1}{4}\alpha-\frac{1}{4}\beta\\
-\frac{1}{4}\alpha+\frac{1}{4}\beta+\frac{1}{4}\alpha
\end{pmatrix}
\\
& =\frac{4}{\beta-\alpha}%
\begin{pmatrix}
-\frac{1}{4}\alpha\\
\frac{1}{4}\beta
\end{pmatrix}
\\
& =\frac{1}{\beta-\alpha}%
\begin{pmatrix}
-\alpha\\
\beta
\end{pmatrix}
,
\end{align*}

and hence%
\[
\left\vert T\left(  s1_{m}\right)  \right\vert =A\alpha^{m}+B\beta^{m}%
=\frac{\beta^{m+1}-\alpha^{m+1}}{\beta-\alpha}.
\]

\end{proof}

\fbox{Estimate 1} The homogeneity result \ref{av825} i.e. $\left\vert T\left(
ta;tv\right)  \right\vert =t^{m}\left\vert T\left(  a;1\right)  \right\vert $,
implies%
\[
\left\vert T\left(  a;v\right)  \right\vert =\left\vert T\left(  \frac
{a}{\left\vert \left(  a,v\right)  \right\vert };\frac{v}{\left\vert \left(
a,v\right)  \right\vert }\right)  \right\vert \left\vert \left(  a,v\right)
\right\vert ^{m},\quad a\in\mathbb{R}^{m},\text{ }v\in\mathbb{R}^{1},
\]

and so if $\left(  x,s\right)  \geq0_{m+1}$,
\begin{align*}
\left\vert T\left(  a;v\right)  \right\vert  & \leq\left(  \max_{\left\vert
x\right\vert ^{2}+s^{2}=1}\left\vert T\left(  x;s\right)  \right\vert \right)
\left\vert \left(  a,v\right)  \right\vert ^{m}\\
& =\max_{\left\vert x\right\vert \leq1,x\geq\mathbf{0}}\left\vert T\left(
x;\left(  1-\left\vert x\right\vert ^{2}\right)  ^{1/2}\right)  \right\vert
\,\left\vert \left(  a,v\right)  \right\vert ^{m}\\
& =\max_{\left\vert x\right\vert \leq1,x\geq\mathbf{0}}\left\vert T\left(
x;\left(  1-\left\vert x\right\vert ^{2}\right)  ^{1/2}\right)  \right\vert
\left(  \left\vert a\right\vert ^{2}+v^{2}\right)  ^{\frac{m}{2}}.
\end{align*}

But%
\begin{align*}
\left\vert T\left(  x;\left(  1-\left\vert x\right\vert ^{2}\right)
^{1/2}\right)  \right\vert  & =\sum\limits_{k=0}^{m}\left(  1-\left\vert
x\right\vert ^{2}\right)  ^{k/2}\sum\limits_{\gamma\leq\mathbf{1,}\left\vert
\gamma\right\vert =k}c_{\gamma}^{\left(  m\right)  }x^{\mathbf{\gamma}}\\
& =\sum\limits_{k=0}^{m}\left(  1-\left\vert x\right\vert ^{2}\right)
^{k/2}\sum\limits_{\gamma\leq\mathbf{1,}\left\vert \gamma\right\vert
=k}c_{\gamma}^{\left(  m\right)  }\left\vert x^{\mathbf{\gamma}}\right\vert \\
& \leq\sum\limits_{k=0}^{m}\left(  1-\left\vert x\right\vert ^{2}\right)
^{k/2}\sum\limits_{\gamma\leq\mathbf{1,}\left\vert \gamma\right\vert
=k}c_{\gamma}^{\left(  m\right)  }\left\vert x\right\vert ^{\left\vert
\gamma\right\vert }\\
& =\sum\limits_{k=0}^{m}\left(  1-\left\vert x\right\vert ^{2}\right)
^{k/2}\sum\limits_{\gamma\leq\mathbf{1,}\left\vert \gamma\right\vert
=k}c_{\gamma}^{\left(  m\right)  }\left\vert x\right\vert ^{k}\\
& =\sum\limits_{k=0}^{m}\left(  1-\left\vert x\right\vert ^{2}\right)
^{k/2}\left\vert x\right\vert ^{k}\sum\limits_{\gamma\leq\mathbf{1,}\left\vert
\gamma\right\vert =k}c_{\gamma}^{\left(  m\right)  }.
\end{align*}

But when $\left\vert x\right\vert \leq1$ and $x\geq0$:%
\[
\left(  1-\left\vert x\right\vert ^{2}\right)  ^{k/2}\left\vert x\right\vert
^{k}=\left(  \left(  1-\left\vert x\right\vert ^{2}\right)  \left\vert
x\right\vert ^{2}\right)  ^{k/2}\leq\left(  \frac{1}{4}\right)  ^{k/2}%
=\frac{1}{2^{k}},
\]

so that%
\[
\left\vert T\left(  x;\left(  1-\left\vert x\right\vert ^{2}\right)
^{1/2}\right)  \right\vert \leq\sum\limits_{k=0}^{m}\frac{1}{2^{k}}%
\sum\limits_{\substack{\gamma\leq\mathbf{1} \\\left\vert \gamma\right\vert
=k}}c_{\gamma}^{\left(  m\right)  }=\sum\limits_{k=0}^{m}\sum
\limits_{\substack{\gamma\leq\mathbf{1} \\\left\vert \gamma\right\vert
=k}}c_{\gamma}^{\left(  m\right)  }\left(  \frac{1}{2}1_{m}\right)  ^{\gamma
}=T\left(  \frac{1}{2}1_{m}\right)  .
\]

Now we can conclude that%
\[
\left\vert T\left(  a;v\right)  \right\vert \leq\max_{\substack{\left\vert
x\right\vert \leq1 \\x\geq0}}\left\vert T\left(  x;\left(  1-\left\vert
x\right\vert ^{2}\right)  ^{\frac{1}{2}}\right)  \right\vert \,\left(
\left\vert a\right\vert ^{2}+v^{2}\right)  ^{\frac{m}{2}}\leq\left\vert
T\left(  \frac{1}{2}1_{m}\right)  \right\vert \left(  \left\vert a\right\vert
^{2}+v^{2}\right)  ^{\frac{m}{2}},\quad a\in\mathbb{R}_{\oplus}^{m},
\]

where $\left\vert T\left(  \frac{1}{2}1_{m}\right)  \right\vert $ is given by
Lemma \ref{vLem_T(s1m)}. Hence%
\begin{equation}
\left\vert T\left(  a\right)  \right\vert \leq\left\vert T\left(  \frac{1}%
{2}1_{m}\right)  \right\vert \left(  \left\vert a\right\vert ^{2}+1\right)
^{\frac{m}{2}},\quad a\in\mathbb{R}_{\oplus}^{m},\label{av051}%
\end{equation}
\medskip

\fbox{Estimate 2} If $a\in\mathbb{R}_{\oplus}^{m}$ then%
\begin{align}
\left\vert T\left(  a\right)  \right\vert =\sum\limits_{k=0}^{m}%
\sum\limits_{\gamma\leq\mathbf{1},\left\vert \gamma\right\vert =k}c_{\gamma
}^{\left(  m\right)  }a^{\gamma}=\sum\limits_{k=0}^{m}\sum\limits_{\gamma
\leq\mathbf{1},\left\vert \gamma\right\vert =k}c_{\gamma}^{\left(  m\right)
}\left\vert a^{\gamma}\right\vert  & =\sum\limits_{k=0}^{m}\sum\limits_{\gamma
\leq\mathbf{1},\left\vert \gamma\right\vert =k}c_{\gamma}^{\left(  m\right)
}\left\vert a\right\vert ^{\left\vert \gamma\right\vert }\nonumber\\
& =\sum\limits_{k=0}^{m}\sum\limits_{\gamma\leq\mathbf{1},\left\vert
\gamma\right\vert =k}c_{\gamma}^{\left(  m\right)  }\left\vert a\right\vert
^{k}\nonumber\\
& \leq\left(  \sum\limits_{k=0}^{m}\sum\limits_{\gamma\leq\mathbf{1}%
,\left\vert \gamma\right\vert =k}c_{\gamma}^{\left(  m\right)  }\right)
\max_{k=0}^{m}\left\vert a\right\vert ^{k}\nonumber\\
& =\left(  \sum\limits_{k=0}^{m}\sum\limits_{\gamma\leq\mathbf{1},\left\vert
\gamma\right\vert =k}c_{\gamma}^{\left(  m\right)  }1_{m}^{\gamma}\right)
\max_{k=0}^{m}\left\vert a\right\vert ^{k}\nonumber\\
& =\left\vert T\left(  1_{m}\right)  \right\vert \max_{k=0}^{m}\left\vert
a\right\vert ^{k}\nonumber\\
& =\left\vert T\left(  1_{m}\right)  \right\vert \left(  \max\left\{
1,\left\vert a\right\vert \right\}  \right)  ^{m},\label{av05}%
\end{align}

where $\left\vert T\left(  1_{m}\right)  \right\vert $ is given by Lemma
\ref{vLem_T(s1m)}.

\chapter{Multivariate bounds for the hat function interpolant/exact smoother
?? incomplete ??}

From ?? we only need to show that each $D_{j}s$ has an appropriate bound.

?? The analogue of the approach of Section \ref{Sect_Cramer_1dim_deriv_smther}
is the take the derivative $D^{\mathbf{1}}=D_{1}D_{2}\ldots D_{d}$ of the
smoother. This again yields $N$ linear equations in the $\left(  \alpha
_{i}\right)  $ with coefficients that are $\pm1$. We will start by using the
unscaled tensor product hat basis function smoother \ref{av073} i.e.%
\[
s\left(  x\right)  =\sum_{k=1}^{N}\alpha_{k}\Lambda\left(  x-x^{\left(
k\right)  }\right)  ,\text{\quad}x\in\mathbb{R}^{d},\text{ }X=\left\{
x^{\left(  1\right)  },x^{\left(  2\right)  },\ldots,x^{\left(  N\right)
}\right\}  ,
\]

and assume that $\Lambda$ has large support w.r.t. the data region i.e. when
$x\in R\left(  x^{\left(  0\right)  },x^{\left(  N\right)  }\right)
\supseteq\Omega$,%
\[
\Lambda\left(  x\right)  =\prod_{i=1}^{d}\Lambda\left(  x_{i}\right)
=\prod_{i=1}^{d}\left(  1-\left\vert x_{i}\right\vert \right)  .
\]

This means that the smoother is a continuous, piecewise linear function with a
finite number of derivative values. Indeed%
\[
D^{\mathbf{1}}\Lambda\left(  x\right)  =\operatorname*{sgn}x\in\pm1,\quad
x\centerdot\neq\mathbf{0},
\]

and thus%
\begin{equation}
D^{\mathbf{1}}s\left(  x\right)  =\sum_{k=1}^{N}\alpha_{k}\operatorname*{sgn}%
\left(  x-x^{\left(  k\right)  }\right)  \text{,\quad}x\notin X,\label{av024}%
\end{equation}
%

\begin{figure}[ptb]%
\centering
\includegraphics[
natheight=5.263800in,
natwidth=5.316100in,
height=5.2638in,
width=5.3161in
]%
{C:/Math_SwBasisFunc/InterpolSmthDev/PapersMonog/ZeroOrd/ZeroOrdDev/graphics/figLargeSuppHatBasisFunc__21.pdf}%
\caption{Two-dimensional example with five data points.}%
\label{Fig_large_supp_hat_multivar}%
\end{figure}

Assume that%
\begin{equation}
j\neq k\Rightarrow x_{i}^{\left(  j\right)  }\neq x_{i}^{\left(  k\right)
}\text{ }\forall i,\label{av013}%
\end{equation}

and also that%
\[
x_{1}^{\left(  j\right)  }\neq x_{1}^{\left(  j+1\right)  }.
\]

If $R$ is an open rectangle which contains $\Omega$ let $x^{\left(  0\right)
}$ be the left-most point of $R$ and $x^{\left(  N+1\right)  }$ be the
right-most point of $R$.

Set%
\[
X^{+}=X\cup\left\{  x^{\left(  0\right)  },x^{\left(  N+1\right)  }\right\}  .
\]

For each dimension $k$ there exists a permutation $\sigma_{k}$ of $1:d$ which
places the components $x_{k}^{\left(  j\right)  }$ in increasing order i.e.%
\begin{align*}
\sigma_{1}  & =1,\\
z_{k}^{\left(  j\right)  }  & :=x_{k}^{\left(  \sigma_{k}\left(  j\right)
\right)  },\quad z_{k}^{\left(  j\right)  }<z_{k}^{\left(  j+1\right)  },\\
z_{k}^{\left(  0\right)  }  & :=x_{k}^{\left(  0\right)  },\quad
z_{k}^{\left(  N+1\right)  }:=x_{k}^{\left(  N+1\right)  }.
\end{align*}

Clearly%
\[
x^{\left(  k\right)  }=\left(  z_{1}^{\left(  k\right)  },\ldots
,z_{d}^{\left(  k\right)  }\right)  .
\]

The projections of $X^{+}$ onto each axis are denoted by%
\[
\mathcal{P}_{k}X^{+}:=\left\{  z_{k}^{\left(  j\right)  }\right\}
_{j=0}^{N+1},
\]

which in turn allows us to define the grid points $\left\{  z^{\left(
\alpha\right)  }\right\}  =\bigotimes_{k=1}^{d}\mathcal{P}_{k}X^{+}$;
\[
z^{\left(  \alpha\right)  }:=\left(  z_{1}^{\alpha_{1}},\ldots,z_{d}%
^{\alpha_{d}}\right)  ,\quad\mathbf{0}\leq\alpha\leq N+\mathbf{1}.
\]

The corresponding set of open rectangles is denoted
\[
R_{\alpha}:=\left\{  x:z^{\left(  \alpha\right)  }\centerdot<x\centerdot
<z^{\left(  \alpha+\mathbf{1}\right)  }\right\}  .
\]

Thus, from \ref{av024} $D^{\mathbf{1}}s$ is constant on each $R_{\alpha}$.

We also have%
\[
\int_{R\left[  a,b\right]  }D^{\mathbf{1}}s=s\left(  b\right)  -s\left(
a\right)  .
\]

If $0<\varepsilon<\min_{\alpha}\left\vert R_{\alpha}\right\vert $ and $u\in
\pm\mathbf{1}$ then%
\begin{align*}
D^{\mathbf{1}}s\left(  x^{\left(  j\right)  }+\varepsilon u\right)   &
=\sum_{\substack{k=1 \\k\neq j}}^{N}\alpha_{k}\operatorname*{sgn}\left(
x^{\left(  j\right)  }-x^{\left(  k\right)  }\right)  +\alpha_{j}%
\operatorname*{sgn}\left(  \varepsilon u\right) \\
& =\sum_{\substack{k=1 \\k\neq j}}^{N}\alpha_{k}\operatorname*{sgn}\left(
x^{\left(  j\right)  }-x^{\left(  k\right)  }\right)  +\alpha_{j}%
\operatorname*{sgn}u,
\end{align*}

and so the jumps at $x^{\left(  j\right)  }$ are%
\[
D^{\mathbf{1}}s\left(  x^{\left(  j\right)  }+\varepsilon u\right)
-D^{\mathbf{1}}s\left(  x^{\left(  j\right)  }+\varepsilon v\right)
=\alpha_{j}\left(  \operatorname*{sgn}u-\operatorname*{sgn}v\right)  ,\text{
}u,v\in\pm\mathbf{1}.
\]

Thus%
\[
D^{\mathbf{1}}s\left(  x^{\left(  j\right)  }+\varepsilon\right)
=\sum_{\substack{k=1 \\k\neq j}}^{N}\alpha_{k}\operatorname*{sgn}\left(
x^{\left(  j\right)  }-x^{\left(  k\right)  }\right)  +\alpha_{j}%
=:\Upsilon_{X,X}\alpha.
\]

\chapter{Explicit extension operators based on Wloka but using the rectangle
condition\label{Ch_exten_rect_condit}}

\section{Introduction}

I needed these explicit extension operators to characterize of the restriction
spaces $X_{w}^{0}\left(  \Omega\right)  $\textbf{\ }for several classes of
weight function, especially the tensor products in Subsection
\ref{Sect_local_data_space}, but $\Omega$ has ??? very awkward constraints.
Subsequently I discovered the very general theoretical extension operator
$r_{\Omega}^{\ast}$ (Theorem \ref{Thm_canon_exten_op}) and so the explicit
extension became unnecessary.

In this chapter, motivated by the work of Wloka \cite{Wloka87} and the Russian
mathematicians described in the Background section, we will derive
\textbf{four extension operators}. Technically:

\begin{enumerate}
\item We will assume that the bounded domain $\Omega$ satisfies the very
restrictive rectangle condition instead of the cone condition - a sphere does
not satisfy the rectangle condition.

\item I will modify the Calderon-Zygmund extension result used in Theorem 5.4,
Section 5.2 of Wloka \cite{Wloka87} by replacing the cone condition by the
rectangle condition but I will still use the integral representation and the
Fourier transform.

\item We adapt the (convolution) integral representation technique. This is in
the form of the integral operators $J_{\delta}\left[  v\right]  \left(
z\right)  =\int_{\mathcal{O}_{\mathbf{1}}}\lambda^{\delta}e^{-\lambda
\mathbf{1}}v\left(  \frac{z}{\lambda}\right)  d\lambda$ where $v$ is weakened
to be an $L^{1}$ function.

\item Develop an extension operator which extends a function from a single
orthant $\mathcal{O}_{\theta}\subset\mathbb{R}^{d}$ to the entire space.

\item Still uses a $C^{\infty}$ partition of unity to define an extension on a
domain $\Omega$.

\item The basic norms are the $L^{2}$ and $L^{\infty}$ norms.
\end{enumerate}

The four extension operators are:

\begin{enumerate}
\item In Section \ref{Sect_Exten_FromOrthant} I construct continuous extension
$\mathcal{E}_{\alpha}^{\mathbf{1}}:C_{0}^{\left(  \alpha\right)  }\left(
\overline{\mathcal{O}_{\mathbf{1}}}\right)  \rightarrow C_{B}^{\left(
\alpha\right)  }\left(  \mathbb{R}^{d}\right)  $ for each $\alpha
\geq\mathbf{1}$. This is then generalized to a continuous extension
$\mathcal{E}_{\alpha}^{\theta}:C_{0}^{\left(  \alpha\right)  }\left(
\overline{\mathcal{O}_{\theta}}\right)  \rightarrow C_{B}^{\left(
\alpha\right)  }\left(  \mathbb{R}^{d}\right)  $. The results of this section
are not used elsewhere and are not used to construct any of the other
extensions in this Chapter.

\item In Section \ref{Sect_ExtenLocWn1_to_Wn1_Fourier} a partition of unity
and the integral representation of Lemma \ref{Lem_SmthFuncIntegRepInOrthant}
to construct continuous convolution extension operators $E_{\Omega
}^{n\mathbf{1}}:W^{n\mathbf{1}}\left(  \Omega\right)  \rightarrow
W^{n\mathbf{1}}\left(  \mathbb{R}^{d}\right)  $ for $n\geq1$. This is done in
Theorem \ref{Thm_ExtenOrthantSobolFourier} and continuity is demonstrated
using a Fourier transform argument. Here $J_{\left(  n-2\right)  \mathbf{1}%
}^{\theta}\left[  v\right]  $ satisfies \ref{a030}.

\item In Section \ref{Sect_ExtenLocCatoBndCa} we generalize the $E_{\Omega
}^{n\mathbf{1}}$ to extension operators $E_{\Omega}^{\alpha}:C^{\left(
\alpha\right)  }\left(  \overline{\Omega}\right)  \rightarrow C_{B}^{\left(
\alpha\right)  }\left(  \mathbb{R}^{d}\right)  $ which are continuous under
the supremum norm. To do this we constrain the function $v$ to have bounded
support in $\mathcal{O}_{\mathbf{1}}$ and to be a tensor product function
which satisfies \ref{X551}. The main extension result is Theorem
\ref{Thm_ExtenContinFuncs_OrthantProp}.

\item In this section we show that the particular extension operators
$\left\{  E_{\Omega}^{n\mathbf{1}}\right\}  _{n\geq1}$ are such that each
$E_{\Omega}^{n\mathbf{1}}:W^{n\mathbf{1}}\left(  \Omega\right)  \rightarrow
W^{n\mathbf{1}}\left(  \mathbb{R}^{d}\right)  $ is continuous if we assume
that in the integral operator $J_{\left(  n-2\right)  \mathbf{1}}^{\theta
}\left[  v\right]  $ the function $v\in L_{0}^{1}\left(  \mathcal{O}%
_{1}\right)  $ is a tensor product with property \ref{X551}.
\end{enumerate}

\section{Background\label{Sect_ExtenBackground}}

\textbf{Background to the extension problem and the rectangle condition} From
Theorem 4.26 of Adams \cite{Adams75} there exists a continuous extension
operator $\mathcal{E}_{\Omega}^{\prime}:W^{m}\left(  \Omega\right)
\rightarrow W^{m}$ \textit{which is constructed using the \textbf{method of
reflections}} and assumes $\Omega$ has the \textbf{uniform }$C^{\left(
m\right)  }$\textbf{-regularity property} of Section 4.6. We have
$\mathcal{E}_{\Omega}^{\prime}:W^{k}\left(  \Omega\right)  \rightarrow W^{k}$
is continuous for $k\leq m$. The proof of Theorem 4.26 involves a
\textbf{two-step} construction. The \textbf{first (reflection) step} involves
an extension from the half-space $\Omega=\mathbb{R}_{+}^{d}$ to $\mathbb{R}%
^{d} $ by reflection and it is easy to show that this extension is also
continuous from $W^{m\mathbf{1}}\left(  \mathbb{R}_{+}^{d}\right)  $ to
$W^{m\mathbf{1}}\left(  \mathbb{R}^{d}\right)  $.

The \textbf{second (localization) step} uses the uniform $C^{\left(  m\right)
}$-regularity property to locally map neighborhoods of the boundary
$\partial\Omega$ into the unit ball $B_{1}$.

It turns out that if we try to modify the localization step by replacing the
$C^{\left(  m\right)  }$-regular map by one having smoothness $C^{\left(
m\mathbf{1}\right)  }$ then we encounter a basic problem and our attempt
fails. In fact, the equation
\begin{equation}
D_{i}\left(  u\left(  y_{1}\left(  x_{1},x_{2}\right)  ,y_{2}\left(
x_{1},x_{2}\right)  \right)  \right)  =\left(  D_{1}u\right)  \left(
y_{1},y_{2}\right)  D_{i}y_{1}\left(  x_{1},x_{2}\right)  +\left(
D_{2}u\right)  \left(  y_{1},y_{2}\right)  D_{i}y_{1}\left(  x_{1}%
,x_{2}\right)  ,\label{2.18}%
\end{equation}

implies that in general $D_{1}D_{2}\left(  u\left(  y_{1}\left(  x_{1}%
,x_{2}\right)  ,y_{2}\left(  x_{1},x_{2}\right)  \right)  \right)  $ always
involves the terms $D_{1}^{2}u,D_{1}D_{2}u,D_{2}^{2}u$ and so we have
continuity into $W^{2}$ but not into $W^{2\mathbf{1}}$.

Problems relating to \ref{2.18} were encountered by the Russian school of
mathematicians studying generalizations of the isotropic Sobolev spaces
$W^{m,??p}$ after World War 2. Prominent among this group and relevant to my
work include S. L. Sobolev, O. V. Besov, S. M. Nikol'ski\u{\i}, V. P. Il'in,
A. D. D\v{z}abrailov, P. L. Lizorkin. Their work employed function spaces
based on the spaces $L^{p}$, $1\leq p\leq\infty$, but latter I will only
interested in the $L^{2}$ base space and integral derivatives which simplifies
calculations considerably.

Note that the \textbf{embedding theorems} (problems) referred to below are
estimates of the form%
\begin{equation}
\left\Vert D^{\nu}u\right\Vert _{L^{p}\left(  \Omega\right)  }\leq
c_{\mathcal{E}}%
{\textstyle\sum\limits_{\alpha\in\mathcal{E}}}
\left\Vert D^{\alpha}u\right\Vert _{L^{p}\left(  \Omega\right)  }<\infty,\quad
u\in L_{loc}^{1}\left(  \Omega\right)  .\label{2.19}%
\end{equation}

where given a set of multi-indexes $\mathcal{E}$ we want to determine the set
of $\nu$ for which \ref{2.19} holds for some constant $c_{\mathcal{E}}$
independent of the $u$.\smallskip

(\textbf{a}) Nikol'ski\u{\i} \cite{Nikol62}, \cite{Nikol63} encountered
boundary geometry problems when he introduced the Banach spaces of
\textbf{mixed dominant derivatives} $S_{p}^{\alpha}\left(  \Omega\right)  $
which Nikol'ski\u{\i} \cite{Nikol6373} suggested be permanently
\textbf{renamed} $S_{p}^{\alpha}W\left(  \Omega\right)  $. Here $\alpha
.\mathbf{\geq0}$ is a multi-index. In particular%
\begin{equation}%
\begin{array}
[c]{l}%
\mathcal{A}_{\alpha}=\alpha.\left\{  0,1\right\}  ^{d}=%
{\textstyle\bigotimes\limits_{i=1}^{d}}
\left\{  0,\alpha_{i}\right\}  ,\\
S_{p}^{\alpha}\left(  \Omega\right)  \equiv S_{p}^{\alpha}W\left(
\Omega\right)  =\left\{  u\in L^{p}\left(  \Omega\right)  :D^{\beta}u\in
L^{p}\left(  \Omega\right)  \text{ }for\text{ }\beta\in\mathcal{A}_{\alpha
}\right\}  ,\\
\left\Vert u\right\Vert _{S_{p}^{\alpha}\left(  \Omega\right)  }^{2}=%
{\textstyle\sum_{\beta\in\mathcal{A}_{\alpha}}}
\left\Vert D^{\beta}u\right\Vert _{L^{p}\left(  \Omega\right)  }^{2}.
\end{array}
\label{X02}%
\end{equation}

have continuous boundary values (traces). This uses conditions on the boundary
that relate to\textbf{\ rectangles or cubes} $\Delta$ \textbf{with edges
parallel to the coordinate axes}.\smallskip

(\textbf{b}) Il'in \cite{Ilin63A} is a brief paper which discusses basic
embedding theorems, of the form \ref{2.19}, for domains that \textbf{are
rectangles} $\Delta$\textbf{\ with edges parallel to the coordinate axes}. He
used integral representations for smooth functions.\smallskip

(\textbf{c}) Lizorkin and Nikol'ski\u{\i} \cite{LizNikol6567} discussed the
spaces $S_{p}^{\alpha}\left(  \Delta\right)  \equiv S_{p}^{\alpha}W\left(
\Delta\right)  $ and
\[
S_{p}^{\alpha^{\left(  1\right)  },\ldots,\alpha^{\left(  N\right)  }}\left(
\mathbb{R}^{d}\right)  \equiv S_{p}^{\alpha^{\left(  1\right)  },\ldots
,\alpha^{\left(  N\right)  }}W\left(  \mathbb{R}^{d}\right)  =\cap_{i=1}%
^{N}S_{p}^{\alpha^{\left(  i\right)  }}W\left(  \mathbb{R}^{d}\right)
\equiv\cap_{i=1}^{N}S_{p}^{\alpha^{\left(  i\right)  }}\left(  \mathbb{R}%
^{d}\right)  .
\]

He derives some integral representations for smooth functions.\smallskip

(\textbf{d}) D\v{z}abrailov \cite{Dzhab67} used difference operators to
specify the multiple-integral H\"{o}lder condition. D\v{z}abrailov introduced
the spaces $S_{p}^{r}W\left(  \mathbb{R}^{d}\right)  $ of mixed dominant
derivatives(or smoothness) for arbitrary $r\in\mathbb{R}^{d}$, $r\geq
\mathbf{0}$. He uses integral representations for smooth functions to derive
embedding theorems as well as trace theorems into $\mathbb{R}^{m}$ where
$m\leq d$.\smallskip

(\textbf{e}) Il'in \cite{Ilin68} derives basic embedding results for domains
$\Omega$ which satisfy the \textbf{rectangle condition }$C\left(  H\right)  $
\textbf{which uses rectangles with edges parallel to the coordinate axes} but
assumes $\mathcal{E}$ (of inequality \ref{2.19}) contains $\leq d$
multi-indexes. In general this excludes $W^{m\mathbf{1}}$. Il'in used integral
representations for smooth functions.

This paper demonstrates the link between the geometry of the boundary and the
set $\mathcal{E}$ in the embedding \ref{2.19} inequality.\smallskip

(\textbf{f}) D\v{z}abrailov \cite{Dzhab69} announced in interpolation theorems
for $S_{p}^{r}W\left(  \Omega\right)  $ and

$S_{\left(  p_{i}\right)  }^{\left(  r^{\left(  1\right)  },\ldots,r^{\left(
d\right)  }\right)  }W\left(  \Omega\right)  =%
{\textstyle\bigcap_{i=1}^{d}}
S_{\left(  p_{i}\right)  }^{r^{\left(  i\right)  }}W\left(  \Omega\right)  $
where the domain $\Omega$ satisfies the \textbf{rectangle condition} $C\left(
H\right)  $ as described in Il'in \cite{Ilin68}.

(\textbf{g}) D\v{z}abrailov \cite{Dzhab72Fam} introduces the space $%
{\textstyle\bigcap_{i=1}^{2^{d}}}
L_{p_{i}}^{r^{\left(  i\right)  }}\left(  \Omega\right)  $ which generalizes
the space $S_{p}^{r}W\left(  \Omega\right)  $. The spaces are defined using
difference operators and $\forall i$: $1\leq p_{i}\leq\infty$, $r^{\left(
i\right)  }\in\mathbb{R}^{d}$ and $r^{\left(  i\right)  }\geq\mathbf{0}$.
Embedding results are derived in \S 2 using a new integral representation
formula where the domain $\Omega$ satisfies the \textbf{rectangle condition}
$C\left(  H\right)  $.

(\textbf{i}) D\v{z}abrailov \cite{Dzhab74} studies the continuous extension of
functions in $S_{p}^{r}W\left(  \Omega\right)  $ to $S_{p}^{r}W\left(
\mathbb{R}^{d}\right)  $ where $r\in\mathbb{R}^{d}$, $r\geq\mathbf{0}$ and
$\Omega$ satisfies the \textbf{rectangle condition} $C_{0,\varepsilon}\left(
H\right)  $ when $r$ is integral and the condition $C_{\alpha,\varepsilon
}\left(  H\right)  $ when $r$ is not integral. Section 2 derives estimates for
some singular integral operators. Theorem 1 in Section 3 is the extension
result.\smallskip

(\textbf{j}) The review of Besov, Il'in, Nikol'ski\u{\i} \cite{BesovEtAl90} by
Krantz \cite{Krantz80} mentions the simple general principle of the integral
representation technique.

(\textbf{k}) Besov et al. \cite{BesovEtAl90} mentions extensions and $SW$
spaces; D\v{z}abrailov \cite{Dzhab99} mentions extensions of $SW$ spaces;
Schmeisser \cite{Schmeis2006} present a \textbf{survey} of some results on
spaces of functions with dominating mixed smoothness - see Remark
\ref{Rem_Schmeis_DomMixDerivSobol} in appendix; Adams and Fournier
\cite{AdamFour2003} discuss imbedding theorems in Chapter 4 and extensions in
Chapter 5. Note paragraph 4.33.

See also, for example, equation 3.2 of the survey \cite{BungGrieb04} by
Bungartz and Griebel. Section 3 of Chernov \cite{Chern2012}, also Chapter 2 of
the book by Schmeisser and Triebel \cite{SchmeisTrieb87}. The Introduction in
Hansen and Vybiral \cite{HansVyb2009} gives uses of these spaces
$S_{p}^{\widetilde{r}}W\left(  \mathbb{R}^{n_{1}}\mathbb{\times}\ldots
\times\mathbb{R}^{n_{d}}\right)  $.

(\textbf{l}) When $n=1$, \ref{X02} implies $W^{\mathbf{1}}\left(
\Omega\right)  =S_{2}^{\mathbf{1}}W\left(  \Omega\right)  $ and so we can use
the extension results of D\v{z}abrailov \cite{Dzhab74}.\medskip

\textbf{With the above considerations in mind}, in Section
\ref{Ch_exten_rect_condit} I will derive a continuous extension operator
$\mathcal{E}_{\Omega}:W^{n\mathbf{1}}\left(  \Omega\right)  \rightarrow
W^{n\mathbf{1}}$ by \textbf{severely restricting the boundary} of $\Omega$ and
constrain it to satisfy the \textbf{rectangle condition }$C\left(  H\right)  $
of D\v{z}abrailov \cite{Dzhab69}, \cite{Dzhab72Fam} and \cite{Dzhab74}. This
restraint is not so surprising in the light of part (\textbf{e}) above.

I will modify the \textbf{Calderon-Zygmund extension} result used in Theorem
5.4, Section 5.2 of Wloka \cite{Wloka87} by replacing the cone condition by
the rectangle condition but I will still use the integral representation approach.

\section{Integral representation\label{Sect_ExtenIntegApproach}}

We first define various classes of continuous functions, the orthants
$\mathcal{O}_{\theta}$ and the rectangle condition. The function $J_{\delta} $
is introduced which here is defined here using a function $\sigma\in
L^{1}\left(  \mathcal{O}_{\mathbf{1}}\right)  $ instead of a $C_{0}^{\infty
}\left(  \mathcal{O}_{\mathbf{1}}\right)  $ function. The convolution function
$J_{\delta}$ is the basis for our integral representation in an orthant Lemma
\ref{Lem_SmthFuncIntegRepInOrthant}.

\begin{definition}
\textbf{Spaces of continuous functions }Suppose $\Omega$ is an open subset of
$\mathbb{R}^{d}$. Then for any multi-index $\alpha$:

\begin{enumerate}
\item $C^{\left(  \alpha\right)  }\left(  \Omega\right)  =\left\{  f\in
C^{\left(  0\right)  }\left(  \Omega\right)  :D^{\beta}f\in C^{\left(
0\right)  }\left(  \Omega\right)  \text{ }\forall\beta\leq\alpha\right\}  $;

\item $C^{\infty}\left(  \Omega\right)  =\left\{  f\in C^{\left(  0\right)
}\left(  \Omega\right)  :D^{\beta}f\in C^{\left(  0\right)  }\left(
\Omega\right)  \text{ }\forall\beta\right\}  $;

\item $C_{0}^{\left(  \alpha\right)  }\left(  \Omega\right)  =\left\{  f\in
C^{\left(  \alpha\right)  }\left(  \Omega\right)  :\operatorname*{supp}f\text{
}is\text{ }bounded\right\}  $;

\item $C_{0}^{\infty}\left(  \Omega\right)  =\left\{  f\in C^{\infty}\left(
\Omega\right)  :\operatorname*{supp}f\text{ }is\text{ }bounded\right\}  $;

\item For $\alpha\leq\infty$: $C^{\left(  \alpha\right)  }\left(
\overline{\Omega}\right)  $ denotes all functions in $C^{\left(
\alpha\right)  }\left(  \Omega\right)  $ which are uniformly continuous and
bounded on $\Omega$.
\end{enumerate}
\end{definition}

We now define orthants.

\begin{definition}
\label{Def_orthant}\textbf{Orthants} (or sometimes \textbf{hyperoctant) }An
orthant is the analogue in n-dimensional Euclidean space of a quadrant in the
plane or an octant in three dimensions.

The\textbf{\ positive open orthant} is $\mathcal{O}_{\mathbf{1}}:=\left\{
x\in\mathbb{R}^{d}:x_{j}>0\text{ }\forall j\right\}  =\left\{  x\in
\mathbb{R}^{d}:x.>\mathbf{0}\right\}  $.

$\overline{\mathcal{O}_{\mathbf{1}}}$ is known as the \textbf{positive closed
orthant}.

There are $2^{d}$ orthants $\mathcal{O}_{\theta}$, each one corresponding to a
member $\theta$ of the Cartesian product $\left\{  -1,1\right\}  ^{d}$.
Indeed,
\[
\mathcal{O}_{\theta}:=\theta.\mathcal{O}_{\mathbf{1}}:=\left\{  \theta
.x:x\in\mathcal{O}_{\mathbf{1}}\right\}  .
\]

Also, we will need to define the \textbf{unit orthant }rectangle%
\begin{align*}
\widehat{\mathcal{O}}_{\theta}  & :=\left\{  x\in\mathcal{O}_{\theta}:\left(
x_{+}\right)  .<\mathbf{1}\right\}  =\left\{  x\in\mathcal{O}_{\theta}:\left(
\theta.x\right)  .<\mathbf{1}\right\}  =\theta.\widehat{\mathcal{O}%
}_{\mathbf{1}}=\left\{  x\in\mathcal{O}_{\theta}:\left(  x.x\right)
.<\mathbf{1}\right\}  ,\\
\widehat{\mathcal{O}}_{\mathbf{1}}  & :=\left\{  x\in\mathcal{O}_{\mathbf{1}%
}:x.<\mathbf{1}\right\}  .
\end{align*}

where $x_{+}:=\left(  \left\vert x_{i}\right\vert \right)  $.

?? Will translations of orthants also be called orthants?

Write%
\begin{align*}
\int_{0}^{\infty}\int_{0}^{\infty}\ldots\int_{0}^{\infty}f  & =\int%
_{\mathbf{0}}^{\mathbf{\infty}}f=\int_{\mathcal{O}_{\mathbf{1}}}f.\\
\int_{\mathcal{O}_{\theta}}f  & =\int_{\mathcal{O}_{\mathbf{1}}}f\left(
\theta.x\right)  dx.
\end{align*}

\end{definition}

We will only use equation \ref{X30} derived in the next lemma. Using repeated
integration by parts we obtain:

\begin{lemma}
\label{Lem_integ_repres}??? \textbf{FIX}!

If $f\in C^{\left(  n\right)  }\left(  \overline{\mathbb{R}}_{+}^{1}\right)  $
and $D^{k}f\rightarrow0$ exponentially at $+\infty$ for $k\leq n$, then%
\begin{equation}
\frac{\left(  -1\right)  ^{n}}{\left(  n-1\right)  !}\int_{0}^{\infty}%
s^{n-1}D^{n}f\left(  s\right)  ds=f\left(  0\right)  .\label{X30}%
\end{equation}

Further, if $\alpha\geq\mathbf{1}$, $g\in C_{BP}^{\left(  \alpha\right)
}\left(  \overline{\mathcal{O}_{\mathbf{1}}}\right)  $ and $D^{\beta
}g\rightarrow0$ exponentially at $+\infty$ for $\beta\leq\alpha$, then%
\begin{equation}
\frac{\left(  -1\right)  ^{\left\vert \alpha\right\vert }}{\left(
\alpha-1\right)  !}\int_{\mathcal{O}_{\mathbf{1}}}x^{\alpha-1}D^{\alpha
}g\left(  x\right)  dx=g\left(  \mathbf{0}\right)  ,\label{X16}%
\end{equation}

and%
\begin{equation}
\frac{\left(  -1\right)  ^{\left\vert \alpha\right\vert }}{\left(
\alpha-1\right)  !}\int_{\mathcal{O}_{\mathbf{1}}}x^{\alpha-1}\left(
D^{\alpha}g\right)  \left(  y+x\right)  dx=g\left(  y\right)  ,\quad
y\in\overline{\mathcal{O}_{\mathbf{1}}}.\label{X20}%
\end{equation}

\end{lemma}

\begin{proof}
We prove \ref{X30} by repeated integration by parts: if $\varepsilon>0$ and
$n\geq1$,%
\begin{align*}
\int_{\varepsilon}^{\infty}s^{n-1}D^{n}f\left(  s\right)  ds  & =\int%
_{\varepsilon}^{\infty}s^{n-1}dD^{n-1}f\left(  s\right)  =\left[
s^{n-1}D^{n-1}f\left(  s\right)  \right]  _{\varepsilon}^{\infty}-\left(
n-1\right)  \int_{\varepsilon}^{\infty}s^{n-2}D^{n-1}f\left(  s\right)  =\\
& =\varepsilon^{n-1}D^{n-1}f\left(  \varepsilon\right)  -\left(  n-1\right)
\int_{\varepsilon}^{\infty}s^{n-2}D^{n-1}f\left(  s\right)  ds\\
& =\varepsilon^{n-1}D^{n-1}f\left(  \varepsilon\right)  -\left(  n-1\right)
\left(  \varepsilon^{n-2}D^{n-2}f\left(  \varepsilon\right)  -\left(
n-2\right)  \int_{\varepsilon}^{\infty}s^{n-3}D^{n-2}f\left(  s\right)
ds\right) \\
& =\varepsilon^{n-1}D^{n-1}f\left(  \varepsilon\right)  -\left(  n-1\right)
\varepsilon^{n-2}D^{n-2}f\left(  \varepsilon\right)  +\left(  n-1\right)
\left(  n-2\right)  \int_{\varepsilon}^{\infty}s^{n-3}D^{n-2}f\left(
s\right)  ds\\
& =\varepsilon^{n-1}D^{n-1}f\left(  \varepsilon\right)  -\left(  n-1\right)
\varepsilon^{n-2}D^{n-2}f\left(  \varepsilon\right)  +\left(  n-2\right)
\left(  n-1\right)  \varepsilon^{n-3}D^{n-3}f\left(  \varepsilon\right)
-\ldots\\
& -\left(  n-1\right)  \left(  n-2\right)  \ldots\left(  1\right)
\int_{\varepsilon}^{\infty}Df\left(  s\right)  ds\\
& =\varepsilon^{n-1}D^{n-1}f\left(  \varepsilon\right)  -\left(  n-1\right)
\varepsilon^{n-2}D^{n-2}f\left(  \varepsilon\right)  +\left(  n-2\right)
\left(  n-1\right)  \varepsilon^{n-3}D^{n-3}f\left(  \varepsilon\right)
-\ldots\\
& \qquad+\left(  -1\right)  ^{n}\left(  n-1\right)  !f\left(  \varepsilon
\right)  .
\end{align*}

and letting $\varepsilon\rightarrow0^{+}$ we get our first formula \ref{X30}.

??? \textbf{FIX}! \textbf{SEE INDUCTIVE proof IN MY BIT journal paper}!

Set $f_{1}\left(  x_{1}\right)  :=g\left(  x_{1},x^{\prime}\right)  $ so that
$f_{1}\in C^{\left(  \alpha_{1}\right)  }\left(  \overline{\mathbb{R}}_{+}%
^{1}\right)  $ and $D^{k}f_{1}\rightarrow0$ exponentially at $+\infty$ for
$k\leq\alpha_{1}$. By \ref{X30},
\[
\int_{0}^{\infty}x_{1}^{\alpha_{1}-1}D^{\alpha_{1}}f_{1}\left(  x_{1}\right)
dx_{1}=\left(  -1\right)  ^{\alpha_{1}}\left(  \alpha_{1}-1\right)
!f_{1}\left(  0\right)  .
\]

i.e.%
\begin{equation}
\int_{0}^{\infty}x_{1}^{\alpha_{1}-1}D_{1}^{\alpha_{1}}g\left(  x_{1}%
,x^{\prime}\right)  dx_{1}=\left(  -1\right)  ^{\alpha_{1}}\left(  \alpha
_{1}-1\right)  !g\left(  0,x^{\prime}\right)  .\label{X11}%
\end{equation}

Set $f_{2}\left(  x_{2}\right)  :=\left(  -1\right)  ^{\alpha_{1}}\left(
\alpha_{1}-1\right)  !g\left(  0,x_{2},x^{\prime\prime}\right)  $ so that
$f_{2}\in C^{\left(  \alpha_{2}\right)  }\left(  \overline{\mathbb{R}}_{+}%
^{1}\right)  $ and $D^{k}f_{2}\rightarrow0$ exponentially at $+\infty$ for
$k\leq\alpha_{2}$. Then by \ref{X30},%
\begin{equation}
\int_{0}^{\infty}x_{2}^{\alpha_{2}-1}D^{\alpha_{2}}f_{2}\left(  x_{2}\right)
dx_{2}=\left(  -1\right)  ^{\alpha_{2}}\left(  \alpha_{2}-1\right)
!f_{2}\left(  0\right)  .\label{X12}%
\end{equation}

Now%
\begin{equation}
\left(  -1\right)  ^{\alpha_{2}}\left(  \alpha_{2}-1\right)  !f_{2}\left(
0\right)  =\left(  -1\right)  ^{\alpha_{1}+\alpha_{2}}\left(  \alpha
_{1}-1\right)  !\left(  \alpha_{2}-1\right)  !g\left(  0,0,x^{\prime\prime
}\right)  .\label{X13}%
\end{equation}

Also, from \ref{X11},%
\begin{align*}
f_{2}\left(  x_{2}\right)   & =\left(  -1\right)  ^{\alpha_{1}}\left(
\alpha_{1}-1\right)  !g\left(  0,x_{2},x^{\prime\prime}\right) \\
& =\int_{0}^{\infty}x_{1}^{\alpha_{1}-1}D_{1}^{\alpha_{1}}g\left(  x_{1}%
,x_{2},x^{\prime\prime}\right)  dx_{1},
\end{align*}

and from \ref{X12},%
\[
\int_{0}^{\infty}x_{2}^{\alpha_{2}-1}D^{\alpha_{2}}f_{2}\left(  x_{2}\right)
dx_{2}=\int_{0}^{\infty}x_{2}^{\alpha_{2}-1}D_{x_{2}}^{\alpha_{2}}\left(
\int_{0}^{\infty}x_{1}^{\alpha_{1}-1}D_{1}^{\alpha_{1}}g\left(  x_{1}%
,x_{2},x^{\prime\prime}\right)  dx_{1}\right)  dx_{2}.
\]

We will use Lemma \ref{Lem_diff_under_integral_2} to justify differentiating
under the integral sign: clearly the integrand has sufficient smoothness and
also the condition $g\in C_{BP}^{\left(  \alpha\right)  }\left(
\overline{\mathcal{O}_{\mathbf{1}}}\right)  $ means that the integral is
absolutely convergent. Thus%
\[
\int_{0}^{\infty}x_{2}^{\alpha_{2}-1}D^{\alpha_{2}}f_{2}\left(  x_{2}\right)
dx_{2}=\int_{0}^{\infty}x_{2}^{\alpha_{2}-1}\int_{0}^{\infty}x_{1}^{\alpha
_{1}-1}D_{1}^{\alpha_{1}}D_{2}^{\alpha_{2}}g\left(  x_{1},x_{2},x^{\prime
\prime}\right)  dx_{1}dx_{2},
\]

and so \ref{X12} and \ref{X13} imply%
\begin{align*}
\left(  -1\right)  ^{\alpha_{1}+\alpha_{2}} &  \left(  \alpha_{1}-1\right)
!\left(  \alpha_{2}-1\right)  !g\left(  0,0,x^{\prime\prime}\right) \\
&  =\int_{0}^{\infty}\int_{0}^{\infty}x_{1}^{\alpha_{1}-1}x_{2}^{\alpha_{2}%
-1}D_{1}^{\alpha_{1}}D_{2}^{\alpha_{2}}g\left(  x_{1},x_{2},x^{\prime\prime
}\right)  dx_{1}dx_{2}.
\end{align*}

Repeating this process until exhaustion yields \ref{X16}.

Finally, for given $y\in\overline{\mathcal{O}_{\mathbf{1}}}$, $g\left(
y+\cdot\right)  \in C_{BP}^{\left(  \alpha\right)  }\left(  \overline
{\mathcal{O}_{\mathbf{1}}}\right)  $ and \ref{X20} holds.
\end{proof}

Following the Russian mathematicians (see Background subsection above) we
introduce the \textbf{rectangle condition} and then it's uniform version:

\begin{definition}
\label{Def_OthantCondit}\textbf{Rectangle condition} Suppose $\Omega
\subset\mathbb{R}^{d}$ is a bounded open set. Then $\Omega$ satisfies the
\textbf{rectangle condition} if for each $x\in\overline{\Omega}$ there exists
an open set $U_{x}$ containing $x$, an orthant $\mathcal{O}_{\theta_{x}}$ and
$h_{x}.>\mathbf{0}$ such that $\left(  U_{x}\cap\overline{\Omega}\right)
+\tau.h_{x}.\widehat{\mathcal{O}}_{\theta_{x}}\subset\Omega$ whenever
$\mathbf{0}\leq\tau\leq\mathbf{1}$.

We observe that the rectangle condition:

\begin{enumerate}
\item is invariant under translations and dilations - simply apply the
translation to $U_{x}$ and $\Omega$, and the dilation to $U_{x}$, $\Omega$ and
$h_{x}$;

\item implies the \textbf{segment property} - indeed choose $\tau=t\mathbf{1}$
and any $y_{x}\in h_{x}.\mathcal{O}_{\theta_{x}}$;

\item is not satisfied by an open ball.
\end{enumerate}
\end{definition}

\begin{definition}
\label{Def_UnifRectCondit}\textbf{Uniform rectangle condition} Suppose
$\Omega\subset\mathbb{R}^{d}$ is a \textbf{bounded} open set. Suppose there
exists a finite covering $\left\{  U_{i}\right\}  _{i=1}^{M}$ of $\Omega$ such
that for each $i$ there exists an orthant $\mathcal{O}_{\theta^{\left(
i\right)  }}$ and $h^{\left(  i\right)  }.>\mathbf{0}$ such that%
\[
\left(  U_{i}\cap\overline{\Omega}\right)  +\tau.h^{\left(  i\right)
}.\widehat{\mathcal{O}}_{\theta^{\left(  i\right)  }}\subset\Omega\text{
}whenever\text{ }\mathbf{0}\leq\tau\leq\mathbf{1}.
\]

We observe that the uniform rectangle condition:

\begin{enumerate}
\item is invariant under translations and dilations of $\Omega$ - simply apply
the translation to $U_{i}$ and $\Omega$, and the dilation to $U_{i}$, $\Omega$
and $h^{\left(  i\right)  }$;

\item implies the \textbf{segment property} - indeed choose $\tau=t\mathbf{1}$
and any $y^{\left(  i\right)  }\in h^{\left(  i\right)  }.\mathcal{O}%
_{\theta^{\left(  i\right)  }}$;

\item is not satisfied by an open ball;

\item on a bounded set the rectangle condition and the uniform rectangle
condition are equivalent.
\end{enumerate}
\end{definition}

\begin{definition}
\label{Def_funcs_and_orthants}\textbf{Spaces of continuous functions on
orthants} Assume $\theta\in\left\{  -1,1\right\}  ^{d}$. Then:

\begin{enumerate}
\item $f\in C_{0}^{\infty}\left(  \overline{\mathcal{O}_{\mathbf{1}}}\right)
$ iff $f\in C^{\infty}\left(  \overline{\mathcal{O}_{\mathbf{1}}}\right)  $
and $\operatorname*{supp}f\subset\overline{\mathcal{O}_{\mathbf{1}}}$ is bounded;

\item $f\in C_{0}^{\infty}\left(  \overline{\mathcal{O}_{\theta}}\right)  $
iff $f\left(  \theta.x\right)  \in C_{0}^{\infty}\left(  \overline
{\mathcal{O}_{\mathbf{1}}}\right)  $;

\item $f\in C_{0}^{\infty}\left(  \mathcal{O}_{\mathbf{1}}\right)  $ iff $f\in
C^{\infty}\left(  \mathcal{O}_{\mathbf{1}}\right)  $, $\operatorname*{supp}%
f\subset\mathcal{O}_{\mathbf{1}}$ and $\operatorname*{supp}f$ is bounded;

\item $f\in C_{0}^{\infty}\left(  \mathcal{O}_{\theta}\right)  $ iff $f\left(
\theta.x\right)  \in C_{0}^{\infty}\left(  \mathcal{O}_{\mathbf{1}}\right)  $.
\end{enumerate}
\end{definition}

\begin{definition}
\label{Def_L1o_PWinf}\ 

\textbf{1}) \textbf{The space }$L_{0}^{1}\left(  \Omega\right)  $: If
$\Omega\subset\mathbb{R}^{d}$ is an open set then%
\[
L_{0}^{1}\left(  \Omega\right)  :=\left\{  u\in L^{1}:\operatorname*{supp}%
u\Subset\Omega\right\}  .
\]

\textbf{2}) \textbf{Piecewise} $C^{\left(  \alpha\right)  }$ \textbf{function,
denoted }$PWC^{\left(  \alpha\right)  }$ or\textbf{\ }$PWC^{\left(
\alpha\right)  }\left(  \mathbb{R}^{d}\right)  $\textbf{:} Defined w.r.t. a
rectangular grid parallel to the coordinate axes. Some of the cells may be
unbounded and the union of the closed cells is $\mathbb{R}^{d}$. Each bounded
rectangle must intersect only a finite number of the grid cells.

The function must be $C^{\left(  \alpha\right)  }$ on each closed cell i.e. it
is $C^{\left(  \alpha\right)  }$, bounded and uniformly continuous on each
open cell so that it has a unique $C^{\left(  \alpha\right)  }$ extension to
the closed cell. Thus the function has a \textbf{jump} at the boundary of each
cell which includes at infinity.
\end{definition}

We now introduce our first definition of the function $J_{\delta}$. This will
weakened latter in Definition \ref{Def_Jd_2}. $J_{\delta}$ will be the basis
for our integral representations.

\begin{definition}
\label{Def_Jd}\textbf{The function }$J_{\delta}$\textbf{:} If $\sigma\in
L_{0}^{1}\left(  \mathcal{O}_{\mathbf{1}}\right)  $ and $\delta\in
\mathbb{R}^{d}$ then we define:%
\begin{equation}
J_{\delta}\left(  z\right)  =J_{\delta}\left[  \sigma\right]  \left(
z\right)  :=\int_{\mathcal{O}_{\mathbf{1}}}\lambda^{\delta}e^{-\lambda
\mathbf{1}}\sigma\left(  \frac{z}{\lambda}\right)  d\lambda,\text{\quad}%
z\in\mathbb{R}^{d}.\label{X39}%
\end{equation}

Observe that $J_{\delta}\left(  z\right)  =0$ when $z\in\mathbb{R}%
^{d}\setminus\overline{\mathcal{O}_{\mathbf{1}}}$ i.e.%
\[
J_{\delta}\left(  z\right)  =\left\{
\begin{array}
[c]{ll}%
\int_{\mathcal{O}_{\mathbf{1}}}\lambda^{\delta}e^{-\lambda\mathbf{1}}%
\sigma\left(  \frac{z}{\lambda}\right)  d\lambda, & z\in\mathcal{O}%
_{\mathbf{1}},\\
0 & z\notin\mathbb{R}^{d}\setminus\overline{\mathcal{O}_{\mathbf{1}}}.
\end{array}
\right.
\]

\end{definition}

This lemma is not used by the Fourier proof i.e. by Theorem
\ref{Thm_ExtenOrthantSobolFourier}.

\begin{lemma}
\label{Lem_property_J}\textbf{Properties of }$J_{\delta}=J_{\delta}\left[
\sigma\right]  $\textbf{\ when }$\sigma\in L_{0}^{1}\left(  \Omega\right)  $.
Assume $\operatorname*{supp}\sigma\subset\left[  a,b\right]  \Subset
\mathcal{O}_{\mathbf{1}}$.

\begin{enumerate}
\item
\begin{enumerate}
\item $z^{-\left(  \delta+1\right)  }J_{\delta}\left(  z\right)  \in
C^{\infty}\left(  \overline{\mathcal{O}_{\mathbf{1}}}\right)  $ when
$\delta\in\mathbb{R}^{d}$. In general $z^{-\left(  \delta+1\right)  }%
J_{\delta}\left(  z\right)  \neq0$ for some $z\in\partial\mathcal{O}%
_{\mathbf{1}}$.

\item For each multi-integer $\delta\geq-\mathbf{1}$, $J_{\delta}$ is a
"piecewise $C^{\infty}$" function ($PWC^{\infty}$) of exponential decrease
with $\operatorname*{supp}J_{\delta}\subseteq\overline{\mathcal{O}%
_{\mathbf{1}}}$.
\end{enumerate}

\item For multi-integer $\delta\geq\mathbf{0}$, $J_{\delta}\in C^{\left(
\delta\right)  }\left(  \mathbb{R}^{d}\right)  $ and has exponential decrease
with $\operatorname*{supp}J_{\delta}\subseteq\overline{\mathcal{O}%
_{\mathbf{1}}}$.

\item Suppose $\delta\geq-\mathbf{1}$ is a multi-integer. Then,
$\operatorname*{supp}\sigma\subset R\left(  a,b\right)  \subset\mathcal{O}%
_{\mathbf{1}}$ where $a.<b$, implies%
\begin{equation}
\left.
\begin{array}
[c]{ll}%
\left\vert J_{\delta}\left(  z\right)  \right\vert  & \leq\left(
\delta+1\right)  ^{\left(  \delta+1\right)  }e^{-\left(  \delta+1\right)
\mathbf{1}}\int\frac{\left\vert \sigma\left(  \mu\right)  \right\vert }%
{\mu^{\mathbf{1}}}d\mu\\
& \leq\left(  2\pi\right)  ^{-\frac{d}{2}}\frac{\left(  \delta+1\right)
!}{\left(  \delta+1\right)  ^{\mathbf{1}/2}}\int\frac{\left\vert \sigma\left(
\mu\right)  \right\vert }{\mu^{\mathbf{1}}}d\mu
\end{array}
\right\}  ,\text{\quad}z\in\mathcal{O}_{\mathbf{1}},\label{X96}%
\end{equation}

and%
\begin{equation}
\left\vert D^{\gamma}J_{\delta}\left(  z\right)  \right\vert \leq\left\Vert
\frac{\sigma\left(  x\right)  }{x^{\gamma+1}}\right\Vert _{1}q_{\gamma
}^{\delta}\left(  \frac{z}{a}\right)  e^{-\frac{z}{b}\mathbf{1}},\text{\quad
}\forall\gamma\text{ and }z\in\mathcal{O}_{\mathbf{1}},\label{X32}%
\end{equation}

where $q_{\gamma}^{\delta}\geq0$ is a tensor product polynomial of degree
$\delta+1$ which is independent of $\sigma$ and is given by \ref{X41} in the
proof. Also%
\[
q_{\gamma}^{\delta}\left(  x\right)  =q_{\gamma_{1}}^{\delta_{1}}\left(
x_{1}\right)  \ldots q_{\gamma_{d}}^{\delta_{d}}\left(  x_{d}\right)
,\text{\quad}q_{\mathbf{0}}^{\delta}\left(  x\right)  =x^{\delta+1}.
\]

\item Suppose $\delta\geq-\mathbf{1}$ is a multi-integer. Then $D^{\gamma
}J_{\delta}\in L^{1}\left(  \mathcal{O}_{\mathbf{1}}\right)  $ and%
\[
\int_{\mathcal{O}_{\mathbf{1}}}\left\vert D^{\gamma}J_{\delta}\left[
\sigma\right]  \right\vert \leq\left\Vert \frac{\sigma\left(  x\right)
}{x^{\gamma}}\right\Vert _{1,\mathcal{O}_{\mathbf{1}}}\int_{\mathcal{O}%
_{\mathbf{1}}}q_{\gamma}^{\delta}\left(  y\right)  e^{-y\mathbf{1}}%
dy,\quad\forall\gamma,
\]

where%
\begin{equation}
\int_{0}^{\infty}q_{\gamma_{i}}^{\delta_{i}}\left(  t\right)  e^{-t}dt=\left(
\delta_{i}+1\right)  !\left\{
\begin{array}
[c]{ll}%
2^{\gamma_{i}}, & \gamma_{i}\leq\delta_{i}+1,\\
\sum\limits_{k=0}^{\delta_{i}+1}\tbinom{\gamma_{i}}{k}, & \gamma_{i}%
>\delta_{i}+1.
\end{array}
\right\}  \leq2^{\gamma_{i}}\left(  \delta_{i}+1\right)  !\label{X77}%
\end{equation}

\item Suppose $\delta\geq-\mathbf{1}$ is a multi-integer. Suppose $\left(
p_{\gamma}^{\delta}\right)  _{\beta}$ denotes the coefficient of $x^{\beta}$
of the polynomial $p_{\gamma}^{\delta}\left(  x\right)  $ defined by \ref{X50}
in the proof. Then%
\[
D^{\gamma}J_{\delta}\left(  z\right)  =\frac{1}{z^{\gamma}}\sum_{\beta
\leq\delta+1}\left(  p_{\gamma}^{\delta}\right)  _{\beta}J_{\beta+\gamma
-1}\left(  z\right)  =\sum_{\beta\leq\delta+1}\left(  p_{\gamma}^{\delta
}\right)  _{\beta}\text{ }z^{\beta}g_{\beta+\gamma-1}\left(  z\right)
,\quad\gamma\leq\delta+1.
\]

where each $g_{\beta+\gamma-1}\left(  z\right)  :=\frac{J_{\beta+\gamma
-1}\left(  z\right)  }{z^{\beta+\gamma}}\in C^{\infty}\left(  \overline
{\mathcal{O}_{\mathbf{1}}}\right)  $ and in general $g_{\beta+\gamma-1}\neq0 $
on $\partial\mathcal{O}_{\mathbf{1}}$.

\item If $\delta\geq-\mathbf{1}$ then $D^{\gamma}J_{\delta}$ is a
$PWC^{\infty}$ function of exponential decrease at infinity when $\gamma
\leq\delta+1$.

Part 2 can be extended as follows: If $\delta\geq\mathbf{0}$ and
$\delta<\gamma\leq\delta+1$ then $D^{\gamma}J_{\delta}$ is a $PWC^{\infty}$
function of exponential decrease.

\item If $\mathbf{0}\leq\gamma\leq\delta$ then $\left(  p_{\gamma}^{\delta
}\right)  _{\beta}=0$ if $\beta<\mathbf{1}$.

If $\gamma>\delta$ then $\left(  p_{\gamma}^{\delta}\right)  _{\mathbf{0}%
}=\left(  -1\right)  ^{\left\vert \gamma-\delta-1\right\vert }\frac{\left(
\delta+1\right)  !}{\left(  \gamma-\delta-1\right)  !}\neq0$.
\end{enumerate}
\end{lemma}

\begin{proof}
We have $\operatorname*{supp}\sigma\subset\left[  a,b\right]  \subset
\mathcal{O}_{\mathbf{1}}$ for some $\mathbf{0}.<a.<b$. Now suppose
$z.>\mathbf{0}$. Then%
\[
J_{\delta}\left(  z\right)  =\int_{a\leq\frac{z}{\lambda}\leq b}%
\lambda^{\delta}e^{-\lambda\mathbf{1}}\sigma\left(  \frac{z}{\lambda}\right)
d\lambda=\int_{z./b}^{z./a}\lambda^{\delta}e^{-\lambda\mathbf{1}}\sigma\left(
\frac{z}{\lambda}\right)  d\lambda.
\]

The change of variables: $\mu=z./\lambda$, $d\lambda=\left\vert \left(
-z.\mu^{-2}\right)  ^{\mathbf{1}}\right\vert d\mu=\left(  z.\mu^{-2}\right)
^{\mathbf{1}}d\mu$, yields,%
\begin{equation}
J_{\delta}\left(  z\right)  =\int_{a}^{b}e^{-\frac{z}{\mu}\mathbf{1}}\left(
\frac{z}{\mu}\right)  ^{\delta}\sigma\left(  \mu\right)  \frac{\left(
z\right)  ^{\mathbf{1}}}{\mu^{2\mathbf{1}}}d\mu=z^{\delta+1}\int_{a}%
^{b}e^{-\frac{z}{\mu}\mathbf{1}}\frac{\sigma\left(  \mu\right)  }{\mu
^{\delta+2}}d\mu\in C^{\infty}\left(  \overline{\mathcal{O}_{\mathbf{1}}%
}\right)  ,\label{X31}%
\end{equation}

\textbf{which proves Part 1a}.

Thus if $\delta\geq\mathbf{0}$ then $J_{\delta}\left(  z\right)  =0$ when
$z\in\partial\mathcal{O}_{\mathbf{1}}$. Also, if $\delta\geq-\mathbf{1}$,%
\[
\left\vert J_{\delta}\left(  z\right)  \right\vert \leq\frac{z^{\delta+1}%
}{a^{\delta+1}}e^{-\frac{z}{b}}\int_{a}^{b}\frac{\left\vert \sigma\left(
\mu\right)  \right\vert }{\mu^{\mathbf{1}}}d\mu=\left\Vert \frac{\sigma\left(
x\right)  }{x^{\mathbf{1}}}\right\Vert _{1}\left(  \frac{z}{a}\right)
^{\delta+1}e^{-\frac{z}{b}\mathbf{1}},\quad z\geq\mathbf{0}.
\]

Also, if $\gamma\geq\mathbf{0}$ and $\delta\geq-\mathbf{1}$,%
\begin{align}
D^{\gamma}J_{\delta}\left(  z\right)   & =\int_{a}^{b}D_{z}^{\gamma}\left(
e^{-\frac{z}{\mu}\mathbf{1}}\left(  \frac{z}{\mu}\right)  ^{\delta+1}\right)
\frac{\sigma\left(  \mu\right)  }{\mu}d\mu\nonumber\\
& =\int_{a}^{b}\left(  D_{x}^{\gamma}\left(  x^{\delta+1}e^{-x\mathbf{1}%
}\right)  \right)  \left(  x=\frac{z}{\mu}\right)  \frac{\sigma\left(
\mu\right)  }{\mu^{\gamma+1}}d\mu\nonumber\\
& =\int_{a}^{b}\left(  e^{x\mathbf{1}}D_{x}^{\gamma}\left(  x^{\delta
+1}e^{-x\mathbf{1}}\right)  \right)  \left(  x=\frac{z}{\mu}\right)
\frac{\sigma\left(  \mu\right)  }{\mu^{\gamma+1}}e^{-\frac{z}{\mu}\mathbf{1}%
}d\mu\nonumber\\
& =\int_{a}^{b}p_{\gamma}^{\delta}\left(  \frac{z}{\mu}\right)  \frac
{\sigma\left(  \mu\right)  }{\mu^{\gamma+1}}e^{-\frac{z}{\mu}\mathbf{1}}%
d\mu,\label{X40}%
\end{align}

where%
\begin{equation}
\left.
\begin{array}
[c]{l}%
p_{\gamma}^{\delta}\left(  x\right)  :=e^{x\mathbf{1}}D^{\gamma}\left(
x^{\delta+1}e^{-x\mathbf{1}}\right)  =p_{\gamma_{1}}^{\delta_{1}}\left(
x_{1}\right)  \ldots p_{\gamma_{d}}^{\delta_{d}}\left(  x_{d}\right)  ,\\
\deg p_{\gamma}^{\delta}=\delta+1,\\
\operatorname*{coef}\left(  x^{\delta+1}\right)  =\left(  -1\right)
^{\left\vert \gamma\right\vert },
\end{array}
\right\} \label{X50}%
\end{equation}

?? Note that\ expressing $p_{\gamma}^{\delta}$ in terms of Laguerre
polynomials $L_{n}^{\left(  t\right)  }$ gives the rather cumbersome
expression:
\begin{align*}
p_{\gamma}^{\delta}\left(  x\right)   & =\gamma!x^{2\gamma-\delta-1}L_{\gamma
}^{\left(  \delta+1-\gamma\right)  }\left(  x\right)  .\\
L_{\gamma}^{\left(  \delta+1-\gamma\right)  }  & :=L_{\gamma_{1}}^{\left(
\delta_{1}+1-\gamma_{1}\right)  }\ldots L_{\gamma_{d}}^{\left(  \delta
_{d}+1-\gamma_{d}\right)  }.\\
L_{n}^{\left(  t\right)  }\left(  x\right)   & =\frac{1}{n!}x^{-t}e^{x}%
D^{n}\left(  x^{n+t}e^{-x}\right)  .
\end{align*}

For Laguerre polynomial theory see the online NIST Digital Library of
Mathematical Functions [??], Section 13.2 of Arfken \cite{Arfken70}, and many
other references.

Now for $m\geq0$ and $n\geq0$,%
\begin{align}
p_{n}^{m-1}\left(  t\right)   & =e^{t}D^{n}\left(  t^{m}e^{-t}\right)
\nonumber\\
& =e^{t}\sum_{k=0}^{n}\tbinom{n}{k}D^{k}t^{m}D^{n-k}e^{-t}=e^{t}\sum_{k=0}%
^{n}\tbinom{n}{k}D^{k}t^{m}\left(  -1\right)  ^{n-k}e^{-t}=\nonumber\\
& =\sum_{k=0}^{n}\left(  -1\right)  ^{n-k}\tbinom{n}{k}D^{k}t^{m}=m!\sum
_{k=0}^{\min\left\{  m,n\right\}  }\left(  -1\right)  ^{n-k}\tbinom{n}{k}%
\frac{t^{m-k}}{\left(  m-k\right)  !}=\nonumber\\
& =m!\sum_{j=m}^{m-\min\left\{  m,n\right\}  }\left(  -1\right)
^{n-m+j}\tbinom{n}{m-j}\frac{t^{j}}{j!}\nonumber\\
& =\left(  -1\right)  ^{n-m}m!\sum_{j=m-\min\left\{  m,n\right\}  }^{m}\left(
-1\right)  ^{j}\tbinom{n}{m-j}\frac{t^{j}}{j!}\nonumber\\
& =\left(  -1\right)  ^{n-m}m!\sum_{j=m+\max\left\{  -m,-n\right\}  }%
^{m}\left(  -1\right)  ^{j}\tbinom{n}{m-j}\frac{t^{j}}{j!}\nonumber\\
& =\left(  -1\right)  ^{n-m}m!\sum_{j=\max\left\{  0,m-n\right\}  }^{m}%
\tbinom{n}{m-j}\frac{\left(  -t\right)  ^{j}}{j!}.\label{X52}%
\end{align}

Thus, in the multivariate case we have%
\begin{equation}
\left.
\begin{array}
[c]{l}%
p_{\gamma}^{\delta-1}\left(  x\right)  =\sum\limits_{\beta=\max\left\{
0,\delta-\gamma\right\}  }^{m}\left(  p_{\gamma}^{\delta-1}\right)  _{\beta
}x^{\beta},\\
\left(  p_{\gamma}^{\delta-1}\right)  _{\beta}:=??.
\end{array}
\right\} \label{X72}%
\end{equation}

Equation \ref{X52} also implies%
\[
\left\vert p_{n}^{m}\left(  t\right)  \right\vert \leq\left(  m+1\right)
!\sum_{j=\max\left\{  0,m+1-n\right\}  }^{m+1}\tbinom{n}{m+1-j}\frac
{\left\vert t\right\vert ^{j}}{j!},
\]

and%
\begin{equation}
\left.
\begin{array}
[c]{ll}%
p_{n}^{m}\left(  0\right)  =0, & n<m+1,\\
p_{n}^{m}\left(  0\right)  =\left(  -1\right)  ^{n-m-1}\frac{\left(
m+1\right)  !}{\left(  n-m-1\right)  !}, & n\geq m+1.
\end{array}
\right\} \label{X43}%
\end{equation}

Hence, if $\gamma\geq\mathbf{0}$,%
\begin{align}
\left\vert D^{\gamma}J_{\delta}\left(  z\right)  \right\vert  & \leq\int%
_{a}^{b}\left\vert p_{\gamma}^{\delta+1}\left(  \frac{z}{\mu}\right)
\right\vert \frac{\left\vert \sigma\left(  \mu\right)  \right\vert }%
{\mu^{\gamma+1}}e^{-\frac{z}{\mu}\mathbf{1}}d\mu\nonumber\\
& \leq\left(  \int_{a}^{b}\left\vert p_{\gamma}^{\delta+1}\left(  \frac{z}%
{\mu}\right)  \right\vert \frac{\left\vert \sigma\left(  \mu\right)
\right\vert }{\mu^{\gamma+1}}d\mu\right)  e^{-\frac{z}{b}\mathbf{1}%
}\nonumber\\
& \leq\left(  \int_{a}^{b}\frac{\left\vert \sigma\left(  \mu\right)
\right\vert }{\mu^{\gamma+1}}d\mu\right)  \left(  \max_{\mu\in\left[
a,b\right]  }\left\vert p_{\gamma}^{\delta+1}\left(  \frac{z}{\mu}\right)
\right\vert \right)  e^{-\frac{z}{b}}\mathbf{1}\nonumber\\
& =\left\Vert \frac{\sigma\left(  x\right)  }{x^{\gamma+1}}\right\Vert
_{1}\left(  \max_{\mu\in\left[  a,b\right]  }\left\vert p_{\gamma}^{\delta
+1}\left(  \frac{z}{\mu}\right)  \right\vert \right)  e^{-\frac{z}%
{b}\mathbf{1}},\label{X38}%
\end{align}

which exists and is continuous for all $z\geq\mathbf{0}$. Further, for
$\delta\geq-\mathbf{1}$ and $z\geq\mathbf{0}$,%
\begin{align*}
\max_{\mu\in\left[  a,b\right]  }\left\vert p_{\gamma}^{\delta+1}\left(
\frac{z}{\mu}\right)  \right\vert  & \leq\prod\limits_{i=1}^{d}\max_{\mu
_{i}\in\left[  a_{i},b_{i}\right]  }\left(  \delta_{i}+1\right)  !\sum
_{j=\max\left\{  0,\delta_{i}-\gamma_{i}\right\}  }^{\delta_{i}}\tbinom
{\gamma_{i}}{\delta_{i}+1-j}\frac{1}{j!}\left(  \frac{z_{i}}{\mu_{i}}\right)
^{j}\\
& =\left(  \delta+1\right)  !\prod\limits_{i=1}^{d}\max_{\mu_{i}\in\left[
a_{i},b_{i}\right]  }\sum_{j=\max\left\{  0,\delta_{i}-\gamma_{i}\right\}
}^{\delta_{i}}\tbinom{\gamma_{i}}{\delta_{i}+1-j}\frac{1}{j!}\left(
\frac{z_{i}}{\mu_{i}}\right)  ^{j}\\
& \leq\left(  \delta+1\right)  !\prod\limits_{i=1}^{d}\sum_{j=\max\left\{
0,\delta_{i}+1-\gamma_{i}\right\}  }^{\delta_{i}+1}\tbinom{\gamma_{i}}%
{\delta_{i}+1-j}\frac{1}{j!}\left(  \frac{z_{i}}{a_{i}}\right)  ^{j}\\
& =q_{\gamma}^{\delta+1}\left(  \frac{z}{a}\right)  ,
\end{align*}

where we have defined%
\begin{equation}
q_{\gamma}^{\delta+1}\left(  \zeta\right)  :=\left(  \delta+1\right)
!\prod\limits_{i=1}^{d}\sum_{j=\max\left\{  0,\delta_{i}+1-\gamma_{i}\right\}
}^{\delta_{i}+1}\tbinom{\gamma_{i}}{\delta_{i}+1-j}\frac{1}{j!}\zeta_{i}%
^{j},\quad\zeta\geq\mathbf{0},\text{ }\delta\geq-\mathbf{1}.\label{X41}%
\end{equation}

\textbf{This proves Part 1b}.\medskip

\frame{\textbf{Part 2}} Further, if $z\in\partial\mathcal{O}_{\mathbf{1}}$
then $\frac{z}{\mu}\in\partial\mathcal{O}_{\mathbf{1}}$ and $p_{\gamma
}^{\delta+1}\left(  \frac{z}{\mu}\right)  =0$ when $\gamma\leq\delta$ and so
\ref{X38} implies that $D^{\gamma}J_{\delta}\left(  z\right)  =0$. Thus
$J_{\delta}\in C^{\left(  \delta\right)  }\left(  \mathbb{R}^{d}\right)
$.\medskip

\frame{\textbf{Part 3}} From \ref{X31}%
\[
J_{\delta}\left(  z\right)  =z^{\delta+1}\int e^{-\frac{z}{\mu}\mathbf{1}%
}\frac{\sigma\left(  \mu\right)  }{\mu^{\delta+2}}d\mu=\int\frac{z^{\delta+1}%
}{\mu^{\delta+1}}e^{-\frac{z}{\mu}\mathbf{1}}\frac{\sigma\left(  \mu\right)
}{\mu^{\mathbf{1}}}d\mu,
\]

so that%
\begin{align*}
\left\vert J_{\delta}\left(  z\right)  \right\vert \leq\int\frac{z^{\delta+1}%
}{\mu^{\delta+1}}e^{-\frac{z}{\mu}\mathbf{1}}\frac{\left\vert \sigma\left(
\mu\right)  \right\vert }{\mu^{\mathbf{1}}}d\mu &  \leq\left\Vert y^{\delta
+1}e^{-y\mathbf{1}}\right\Vert _{\infty,\mathcal{O}_{\mathbf{1}}}\int%
\frac{\left\vert \sigma\left(  \mu\right)  \right\vert }{\mu^{\mathbf{1}}}%
d\mu\\
&  =\left(  \delta+1\right)  ^{\left(  \delta+1\right)  }e^{-\left(
\delta+1\right)  \mathbf{1}}\int\frac{\left\vert \sigma\left(  \mu\right)
\right\vert }{\mu^{\mathbf{1}}}d\mu\\
&  \leq\left(  2\pi\right)  ^{-d/2}\frac{\left(  \delta+1\right)  !}{\left(
\delta+1\right)  ^{\mathbf{1}/2}}\int\frac{\left\vert \sigma\left(
\mu\right)  \right\vert }{\mu^{\mathbf{1}}}d\mu,
\end{align*}

where the last step used Lemma \ref{Lem_Jd[f]_bnd_f_in_L1}: $\left\Vert
x^{\lambda}e^{-x\mathbf{1}}\right\Vert _{\infty,\mathcal{O}_{\mathbf{1}}%
}=\lambda^{\lambda}e^{-\lambda\mathbf{1}}\leq\left(  2\pi\right)  ^{-d/2}%
\frac{\lambda!}{\lambda^{\mathbf{1}/2}}$.

Clearly%
\begin{equation}
\left.
\begin{array}
[c]{l}%
q_{\gamma}^{\delta+1}\left(  \zeta\right)  =q_{\gamma_{1}}^{\delta_{1}%
+1}\left(  \zeta_{1}\right)  \ldots q_{\gamma_{d}}^{\delta_{d}+1}\left(
\zeta_{d}\right)  .\\
q_{\gamma_{i}}^{\delta_{i}+1}\left(  0\right)  =p_{\gamma_{i}}^{\delta_{i}%
+1}\left(  0\right)  \text{ }\forall i.
\end{array}
\right\} \label{X42}%
\end{equation}
\medskip

\frame{\textbf{Part 4}} If $\sigma$ is a tensor product function then $J$ is
also.%
\begin{align*}
\int_{\mathcal{O}_{\mathbf{1}}}\left\vert D^{\gamma}J_{\delta}\left(
z\right)  \right\vert dz  & \leq\int_{\mathcal{O}_{\mathbf{1}}}\int_{a}%
^{b}\left\vert p_{\gamma}^{\delta+1}\left(  \frac{z}{\mu}\right)  \right\vert
\frac{\left\vert \sigma\left(  \mu\right)  \right\vert }{\mu^{\gamma+1}%
}e^{-\frac{z}{\mu}\mathbf{1}}d\mu\text{ }dz\\
& =\int_{a}^{b}\frac{\left\vert \sigma\left(  \mu\right)  \right\vert }%
{\mu^{\gamma+1}}\int_{\mathcal{O}_{\mathbf{1}}}\left\vert p_{\gamma}%
^{\delta+1}\left(  \frac{z}{\mu}\right)  \right\vert e^{-\frac{z}{\mu
}\mathbf{1}}dz\text{ }d\mu\\
chg.\text{ }var.  & :y=\frac{z}{\mu},\text{ }z=\mu.y,\text{ }dz=\mu
^{\mathbf{1}}dy\Rightarrow\\
& =\int_{a}^{b}\frac{\left\vert \sigma\left(  \mu\right)  \right\vert }%
{\mu^{\gamma}}d\mu\int_{\mathcal{O}_{\mathbf{1}}}\left\vert p_{\gamma}%
^{\delta+1}\left(  y\right)  \right\vert e^{-y\mathbf{1}}dy\\
& =\left\Vert \frac{\sigma\left(  x\right)  }{x^{\gamma}}\right\Vert _{1}%
\int_{\mathcal{O}_{\mathbf{1}}}\left\vert p_{\gamma}^{\delta+1}\left(
y\right)  \right\vert e^{-y\mathbf{1}}dy\\
& \leq\left\Vert \frac{\sigma\left(  x\right)  }{x^{\gamma}}\right\Vert
_{1}\int_{\mathcal{O}_{\mathbf{1}}}q_{\gamma}^{\delta+1}\left(  y\right)
e^{-y\mathbf{1}}dy,
\end{align*}

We now estimate $\int_{\mathcal{O}_{\mathbf{1}}}q_{\gamma}^{\delta+1}\left(
y\right)  e^{-y\mathbf{1}}dy$ using \ref{X41}.%
\begin{align*}
\int_{0}^{\infty}q_{\gamma_{i}}^{\delta_{i}+1}\left(  t\right)  e^{-t}dt  &
=\int_{0}^{\infty}\left(  \delta_{i}+1\right)  !\sum_{j=\max\left\{
0,\delta_{i}+1-\gamma_{i}\right\}  }^{\delta_{i}+1}\tbinom{\gamma_{i}}%
{\delta_{i}+1-j}\frac{1}{j!}t^{j}e^{-t}dt\\
& =\left(  \delta_{i}+1\right)  !\sum_{j=\max\left\{  0,\delta_{i}%
+1-\gamma_{i}\right\}  }^{\delta_{i}+1}\tbinom{\gamma_{i}}{\delta_{i}%
+1-j}\frac{1}{j!}\int_{0}^{\infty}t^{j}e^{-t}dt\\
& =\left(  \delta_{i}+1\right)  !\sum_{j=\max\left\{  0,\delta_{i}%
+1-\gamma_{i}\right\}  }^{\delta_{i}+1}\tbinom{\gamma_{i}}{\delta_{i}%
+1-j}\frac{1}{j!}j!\\
& =\left(  \delta_{i}+1\right)  !\sum_{j=\max\left\{  0,\delta_{i}%
+1-\gamma_{i}\right\}  }^{\delta_{i}+1}\tbinom{\gamma_{i}}{\delta_{i}+1-j}.
\end{align*}

If $\gamma_{i}\leq\delta_{i}+1$ then
\[
\int_{0}^{\infty}q_{\gamma_{i}}^{\delta_{i}+1}\left(  t\right)  e^{-t}%
dt=\left(  \delta_{i}+1\right)  !\sum_{j=\delta_{i}+1-\gamma_{i}}^{\delta
_{i}+1}\tbinom{\gamma_{i}}{\delta_{i}+1-j}=2^{\gamma_{i}}\left(  \delta
_{i}+1\right)  !,
\]

and if $\gamma_{i}>\delta_{i}+1$ then%
\[
\int_{0}^{\infty}q_{\gamma_{i}}^{\delta_{i}+1}\left(  t\right)  e^{-t}%
dt=\left(  \delta_{i}+1\right)  !\sum_{j=0}^{\delta_{i}+1}\tbinom{\gamma_{i}%
}{\delta_{i}+1-j}=\left(  \delta_{i}+1\right)  !\sum_{k=0}^{\delta_{i}%
+1}\tbinom{\gamma_{i}}{k},
\]

which means that%
\[
\int_{\mathcal{O}_{\mathbf{1}}}q_{\gamma}^{\delta+1}\left(  y\right)
e^{-y\mathbf{1}}dy\leq2^{\left\vert \gamma\right\vert }\left(  \delta
+1\right)  !
\]
\medskip

\textbf{Part 5} If $\gamma\geq\mathbf{0}$ and $\delta\geq-\mathbf{1}$ then
from \ref{X40} and then \ref{X50} and then \ref{X31},%
\begin{align*}
D^{\gamma}J_{\delta}\left(  z\right)   & =\int_{a}^{b}p_{\gamma}^{\delta
}\left(  \frac{z}{\mu}\right)  \frac{\sigma\left(  \mu\right)  }{\mu
^{\gamma+1}}e^{-\frac{z}{\mu}}d\mu=\int_{a}^{b}\left(  \sum_{\beta\leq
\delta+1}\left(  p_{\gamma}^{\delta}\right)  _{\beta}\left(  \frac{z}{\mu
}\right)  ^{\beta}\right)  \frac{\sigma\left(  \mu\right)  }{\mu^{\gamma+1}%
}e^{-\frac{z}{\mu}}d\mu=\\
& =\sum_{\beta\leq\delta+1}\left(  p_{\gamma}^{\delta}\right)  _{\beta
}z^{\beta}\int_{a}^{b}\frac{\sigma\left(  \mu\right)  }{\mu^{\beta+\gamma+1}%
}e^{-\frac{z}{\mu}}d\mu=\sum_{\beta\leq\delta+1}\left(  p_{\gamma}^{\delta
}\right)  _{\beta}z^{\beta}\frac{J_{\beta+\gamma-1}\left(  z\right)
}{z^{\beta+\gamma}}=\\
& =\sum_{\beta\leq\delta+1}\left(  p_{\gamma}^{\delta}\right)  _{\beta}%
\frac{1}{z^{\gamma}}J_{\beta+\gamma-1}\left(  z\right)  =\frac{1}{z^{\gamma}%
}\sum_{\beta\leq\delta+1}\left(  p_{\gamma}^{\delta}\right)  _{\beta}%
J_{\beta+\gamma-1}\left(  z\right)  .
\end{align*}

Also, these calculations imply
\[
D^{\gamma}J_{\delta}\left(  z\right)  =\sum_{\beta\leq\delta+1}\left(
p_{\gamma}^{\delta}\right)  _{\beta}z^{\beta}\frac{J_{\beta+\gamma-1}\left(
z\right)  }{z^{\beta+\gamma}}=\sum_{\beta\leq\delta+1}\left(  p_{\gamma
}^{\delta}\right)  _{\beta}z^{\beta}g_{\beta+\gamma-1}\left(  z\right)  ,
\]

and part 1a implies $g_{\beta+\gamma-1}\in C^{\infty}\left(  \overline
{\mathcal{O}_{\mathbf{1}}}\right)  $.\medskip

\textbf{Part 6} The first statement is a direct consequence of Part 5. If
$\delta\geq\mathbf{0}$ the first statement directly implies that $D^{\gamma
}J_{\delta}$ is a $PWC^{\infty}$ function of exponential decrease when
$\delta<\gamma\leq\delta+1$.\medskip

\textbf{Part 7} From \ref{X50}, $p_{\gamma}^{\delta}\left(  x\right)
:=p_{\gamma_{1}}^{\delta_{1}}\left(  x_{1}\right)  \ldots p_{\gamma_{d}%
}^{\delta_{d}}\left(  x_{d}\right)  $ and from \ref{X43},%
\[
p_{n}^{m}\left(  0\right)  =\left\{
\begin{array}
[c]{ll}%
0, & n\leq m,\\
\left(  -1\right)  ^{n-m-1}\frac{\left(  m+1\right)  !}{\left(  n-m-1\right)
!}, & n\geq m+1.
\end{array}
\right.
\]

Hence $\mathbf{0}\leq\gamma\leq\delta$ implies $p_{\gamma_{i}}^{\delta_{i}%
}\left(  0\right)  =0$ $\forall i$. Now when $\beta<\mathbf{0}$, $\beta_{i}=0$
for some $i$ and so
\[
\left(  p_{\gamma}^{\delta}\right)  _{\beta}=\left(  D^{\beta}p_{\gamma
}^{\delta}\right)  \left(  0\right)  =D_{1}^{\beta_{1}}p_{\gamma_{1}}%
^{\delta_{1}}\left(  0\right)  \ldots D_{d}^{\beta_{d}}p_{\gamma_{d}}%
^{\delta_{d}}\left(  0\right)  =0.
\]

On the other hand, when $\gamma>\delta$,
\[
\left(  p_{\gamma}^{\delta}\right)  _{\mathbf{0}}=p_{\gamma}^{\delta}\left(
0\right)  =p_{\gamma_{1}}^{\delta_{1}}\left(  0\right)  \ldots p_{\gamma_{d}%
}^{\delta_{d}}\left(  0\right)  =\left(  -1\right)  ^{\left\vert \gamma
-\delta-1\right\vert }\frac{\left(  \delta+1\right)  !}{\left(  \gamma
-\delta-1\right)  !}\neq0.
\]

\end{proof}

\begin{theorem}
If $\alpha\leq\delta+1$ and $\sigma\in L_{0}^{1}\left(  \mathcal{O}%
_{\mathbf{1}}\right)  $ then:

\begin{enumerate}
\item
\[
D^{\alpha}J_{\delta}\left[  \sigma\right]  =\left(  -1\right)  ^{\left\vert
\alpha\right\vert }\left(  \delta+1\right)  !\sum_{\beta=\mathbf{0}}^{\alpha
}\frac{\left(  -1\right)  ^{\left\vert \beta\right\vert }\tbinom{\alpha}%
{\beta}}{\left(  \delta-\beta+1\right)  !}J_{\delta-\beta}\left[  \frac
{\sigma}{y^{\alpha}}\right]  .
\]

\item
\[
\left\Vert D^{\alpha}J_{\delta}\left[  \sigma\right]  \right\Vert _{\infty
}\leq2^{\left\vert \alpha\right\vert }\left(  \delta+1\right)  !\left\Vert
\frac{\sigma\left(  x\right)  }{x^{\alpha+1}}\right\Vert _{1}.
\]

\item ??%
\[
\left\Vert D^{\alpha}J_{\delta}\left[  \sigma\right]  \right\Vert _{1}\leq??
\]

\end{enumerate}
\end{theorem}

\begin{proof}
\textbf{Part 1} See Theorem \ref{Thm_DJ_eq_sumJ} below.\medskip

\textbf{Part 2} From part 3 of Lemma \ref{Lem_property_J}, $\left\vert
J_{\delta}\left(  z\right)  \right\vert \leq\left(  \delta+1\right)
^{\delta+1}e^{-\left(  \delta+1\right)  \mathbf{1}}\left\Vert \frac
{\sigma\left(  x\right)  }{x^{\mathbf{1}}}\right\Vert _{1}$, so%
\begin{align*}
\left\Vert D^{\alpha}J_{\delta}\left[  \sigma\right]  \right\Vert _{\infty}  &
\leq\left(  \delta+1\right)  !\sum_{\beta=\mathbf{0}}^{\alpha}\frac
{\tbinom{\alpha}{\beta}}{\left(  \delta-\beta+1\right)  !}\left\Vert
J_{\delta-\beta}\left[  \frac{\sigma}{y^{\alpha}}\right]  \right\Vert
_{\infty}\\
& \leq\left(  \delta+1\right)  !\left(  \sum_{\beta=\mathbf{0}}^{\alpha
}\tbinom{\alpha}{\beta}\frac{\left(  \delta-\beta+1\right)  ^{\delta-\beta+1}%
}{\left(  \delta-\beta+1\right)  !}e^{-\left(  \delta-\beta+1\right)
\mathbf{1}}\right)  \left\Vert \frac{\sigma\left(  x\right)  }{x^{\alpha+1}%
}\right\Vert _{1}\\
& =\left(  \delta+1\right)  !\prod_{k=1}^{d}\left(  \sum_{m=0}^{\alpha_{k}%
}\tbinom{\alpha_{k}}{m}\frac{\left(  \delta_{k}-m+1\right)  ^{\delta_{k}-m+1}%
}{\left(  \delta_{k}-m+1\right)  !}e^{-\left(  \delta_{k}-m+1\right)
}\right)  \left\Vert \frac{\sigma\left(  x\right)  }{x^{\alpha+1}}\right\Vert
_{1}.
\end{align*}

Note that $\frac{\lambda^{\lambda}e^{-\lambda\mathbf{1}}}{\lambda!}\leq\left(
2\pi\right)  ^{-d/2}\frac{1}{\lambda^{\mathbf{1}/2}}$.

If $\alpha_{k}=0$ then%
\begin{align*}
\sum_{m=0}^{\alpha_{k}}\tbinom{\alpha_{k}}{m}\frac{\left(  \delta
_{k}-m+1\right)  ^{\delta_{k}-m+1}}{\left(  \delta_{k}-m+1\right)
!}e^{-\left(  \delta_{k}-m+1\right)  } &  =\frac{\left(  \delta_{k}+1\right)
^{\delta_{k}+1}}{\left(  \delta_{k}+1\right)  !}e^{-\left(  \delta
_{k}+1\right)  }\left(  2\pi\right)  ^{-1/2}\leq\\
&  \leq\frac{\left(  2\pi\right)  ^{-1/2}}{\left(  \delta_{k}+1\right)
^{1/2}}\leq\left(  2\pi\right)  ^{-1/2}<2^{\alpha_{k}}.
\end{align*}

If $1\leq\alpha_{k}<\delta_{k}+1$ then, ,
\begin{align*}
\sum_{m=0}^{\alpha_{k}}\tbinom{\alpha_{k}}{m}\frac{\left(  \delta
_{k}-m+1\right)  ^{\delta_{k}-m+1}}{\left(  \delta_{k}-m+1\right)
!}e^{-\left(  \delta_{k}-m+1\right)  }  & \leq\left(  2\pi\right)  ^{-1/2}%
\sum_{m=0}^{\alpha_{k}}\frac{\tbinom{\alpha_{k}}{m}}{\left(  \delta
_{k}-m+1\right)  ^{1/2}}\\
& \leq\left(  2\pi\right)  ^{-1/2}\sum_{m=0}^{\alpha_{k}}\tbinom{\alpha_{k}%
}{m}\\
& =\left(  2\pi\right)  ^{-1/2}2^{\alpha_{k}}.
\end{align*}

If $1\leq\alpha_{k}$ and $\alpha_{k}=\delta_{k}+1$ then%
\begin{align*}
\sum_{m=0}^{\alpha_{k}}\tbinom{\alpha_{k}}{m}\frac{\left(  \delta
_{k}-m+1\right)  ^{\delta_{k}-m+1}}{\left(  \delta_{k}-m+1\right)
!}e^{-\left(  \delta_{k}-m+1\right)  }  & =\sum_{m=0}^{\alpha_{k}}%
\tbinom{\alpha_{k}}{m}\frac{\left(  \alpha_{k}-m\right)  ^{\alpha_{k}-m}%
}{\left(  \alpha_{k}-m\right)  !}e^{-\left(  \alpha_{k}-m\right)  }\\
& \leq1+\sum_{m=0}^{\alpha_{k}-1}\tbinom{\alpha_{k}}{m}\frac{\left(
\alpha_{k}-m\right)  ^{\alpha_{k}-m}}{\left(  \alpha_{k}-m\right)
!}e^{-\left(  \alpha_{k}-m\right)  }\\
& \leq1+\left(  2\pi\right)  ^{-1/2}\sum_{m=0}^{\alpha_{k}-1}\tbinom
{\alpha_{k}}{m}\frac{1}{\left(  \alpha_{k}-m\right)  ^{\mathbf{1}/2}}\\
& \leq1+\left(  2\pi\right)  ^{-1/2}\sum_{m=0}^{\alpha_{k}-1}\tbinom
{\alpha_{k}}{m}\\
& =1+\frac{2^{\alpha_{k}-1}}{\sqrt{2\pi}}\\
& <2^{\alpha_{k}}.
\end{align*}

Thus%
\[
\left\Vert D^{\alpha}J_{\delta}\left[  \sigma\right]  \right\Vert _{\infty
}\leq\left(  \delta+1\right)  !\left(  \prod_{k=1}^{d}2^{\alpha_{k}}\right)
\left\Vert \frac{\sigma\left(  x\right)  }{x^{\alpha+1}}\right\Vert
_{1}=\left(  \delta+1\right)  !2^{\left\vert \alpha\right\vert }\left\Vert
\frac{\sigma\left(  x\right)  }{x^{\alpha+1}}\right\Vert _{1}.
\]
\medskip

\textbf{Part 3} ??
\end{proof}

In Lemma \ref{Lem_property_J} it was assumed that $\sigma\in L^{1}$ has
bounded support in $\mathcal{O}_{\mathbf{1}}$ i.e. $\sigma\in L_{0}^{1}\left(
\Omega\right)  $. We now weaken the constraint on $\sigma$:

\begin{definition}
\label{Def_Jd_2}\textbf{The function }$J_{\delta}$\textbf{:} If $\sigma\in
L_{\overline{\mathcal{O}_{\mathbf{1}}}}^{1}$ i.e. $\sigma\in L^{1}$,
$\operatorname*{supp}\sigma\subseteq\overline{\mathcal{O}_{\mathbf{1}}}$, and
$\delta\in\mathbb{R}^{d}$ then we define:%
\[
J_{\delta}\left(  z\right)  =J_{\delta}\left[  \sigma\right]  \left(
z\right)  :=\lim_{\substack{r\rightarrow0^{+} \\R\rightarrow\infty}%
}\int_{r\mathbf{1}}^{R\mathbf{1}}\lambda^{\delta}e^{-\lambda\mathbf{1}}%
\sigma\left(  \frac{z}{\lambda}\right)  d\lambda,\text{\quad}z\in
\mathbb{R}^{d}.
\]

Observe that $J_{\delta}\left(  z\right)  =0$ when $z\in\mathbb{R}%
^{d}\setminus\overline{\mathcal{O}_{\mathbf{1}}}$ i.e.%
\[
J_{\delta}\left(  z\right)  =\left\{
\begin{array}
[c]{ll}%
\int_{\mathcal{O}_{\mathbf{1}}}\lambda^{\delta}e^{-\lambda\mathbf{1}}%
\sigma\left(  \frac{z}{\lambda}\right)  d\lambda, & z\in\mathcal{O}%
_{\mathbf{1}},\\
0 & z\notin\mathbb{R}^{d}\setminus\overline{\mathcal{O}_{\mathbf{1}}}.
\end{array}
\right.
\]

\end{definition}

This is a weaker version of Lemma \ref{Lem_property_J} and only assumes
$\sigma\in L_{\overline{\mathcal{O}_{\mathbf{1}}}}^{1}$. This is used by the
Fourier proof i.e. by Theorem \ref{Thm_ExtenOrthantSobolFourier}.

\begin{lemma}
\label{Lem_Jd[f]_bnd_f_in_L1}Suppose $\delta\geq-2\mathbf{1}$ and
$z\in\mathcal{O}_{\mathbf{1}}$. Then:

\begin{enumerate}
\item $\left\Vert x^{\lambda}e^{-x\mathbf{1}}\right\Vert _{\infty
,\mathcal{O}_{\mathbf{1}}}=\lambda^{\lambda}e^{-\lambda\mathbf{1}}\leq\left(
2\pi\right)  ^{-d/2}\frac{\lambda!}{\lambda^{\mathbf{1}/2}}$ when $\lambda
\geq0$.

\item ($L^{\infty}$ result) If $\sigma\in L_{\overline{\mathcal{O}%
_{\mathbf{1}}}}^{1}$ then%
\[
J_{\delta}\left[  \sigma\right]  \left(  z\right)  =z^{\delta+1}%
\int_{\mathcal{O}_{\mathbf{1}}}e^{-\frac{z}{\mu}\mathbf{1}}\frac{\sigma\left(
\mu\right)  }{\mu^{\delta+2}}d\mu,\quad\delta\geq-2,
\]

and%
\[
\left\vert J_{\delta}\left[  \sigma\right]  \left(  z\right)  \right\vert
\leq\left(  \delta+2\right)  ^{\delta+2}e^{-\left(  \delta+2\right)
\mathbf{1}}\left\Vert \sigma\right\Vert _{1}\frac{1}{z^{\mathbf{1}}}%
\leq\left(  2\pi\right)  ^{-d/2}\frac{\left(  \delta+2\right)  !}{\left(
\delta+2\right)  ^{\mathbf{1}/2}}\left\Vert \sigma\right\Vert _{1}\frac
{1}{z^{\mathbf{1}}}.
\]

\item ($L^{1}$ result) When $\delta\geq-\mathbf{1}$, $J_{\delta}\in L^{1}$ and
$\left\Vert J_{\delta}\left[  \sigma\right]  \right\Vert _{1}\leq\left(
\delta+1\right)  !\left\Vert \sigma\right\Vert _{1}$.\medskip

\item ?? FIX! If $x^{\beta}\sigma\in L_{\overline{\mathcal{O}_{\mathbf{1}}}%
}^{1}$ for $\max\left\{  -\alpha,-\left(  \delta+2\right)  \right\}  \leq
\beta\leq\alpha$ then%
\[
\left\vert J_{\delta}\left[  \sigma\right]  \left(  z\right)  \right\vert
\leq\frac{1}{z^{\beta+1}}\left(  \int_{\mathcal{O}_{\mathbf{1}}}\mu^{\beta
}\sigma\left(  \mu\right)  d\mu\right)  \left(  \delta+2+\beta\right)
^{\delta+2+\beta}e^{-\left(  \delta+2+\beta\right)  \mathbf{1}}.
\]

\item ?? If $n\geq\mathbf{1}$, $\delta=\left(  n-2\right)  \mathbf{1}$ and
$\alpha=n\mathbf{1}$ then for $-n\mathbf{1}\leq\beta\leq n\mathbf{1}$:%
\[
\left\vert J_{\left(  n-2\right)  \mathbf{1}}\left[  \sigma\right]  \left(
z\right)  \right\vert \leq\frac{1}{z^{\beta+1}}\left(  \int_{\mathcal{O}%
_{\mathbf{1}}}\mu^{\beta}\sigma\left(  \mu\right)  d\mu\right)  \left(
n\mathbf{1}+\beta\right)  ^{n\mathbf{1}+\beta}e^{-\left(  nd+\beta
\mathbf{1}\right)  }.
\]

?? In fact, if $\beta_{i}\left(  z\right)  =1$ when $z_{i}>1$ and $\beta
_{i}\left(  z\right)  =-1$ when $z_{i}<1$ then%
\[
\left\vert J_{\left(  n-2\right)  \mathbf{1}}\left[  \sigma\right]  \left(
z\right)  \right\vert \leq\frac{1}{z^{n\beta\left(  z\right)  +1}}\left(
\int_{\mathcal{O}_{\mathbf{1}}}\mu^{n\beta\left(  z\right)  }\sigma\left(
\mu\right)  d\mu\right)  \left(  n\mathbf{1}+n\beta\left(  z\right)  \right)
^{n\mathbf{1}+n\beta\left(  z\right)  }e^{-\left(  nd+n\beta\left(  z\right)
\mathbf{1}\right)  }.
\]

\end{enumerate}
\end{lemma}

\begin{proof}
\textbf{Part 1} The maximum of $e^{-x\mathbf{1}}x^{\lambda}$ can be easily
derived by setting the derivative to zero etc. The second estimate follows
from Problem 27 of Section 9.E of Jones \cite{Jones2011}.\medskip

\textbf{Part 2} Assume $\delta\geq-2$. Suppose $z.>\mathbf{0}$. Then%
\[
J_{\delta}\left(  z\right)  =\lim_{\substack{r\rightarrow0^{+} \\R\rightarrow
\infty}}\int_{r\mathbf{1}}^{R\mathbf{1}}\lambda^{\delta}e^{-\lambda\mathbf{1}%
}\sigma\left(  \frac{z}{\lambda}\right)  d\lambda.
\]

The change of variables: $\mu=z./\lambda$, $d\lambda=\left\vert \left(
-z.\mu^{-2}\right)  ^{\mathbf{1}}\right\vert d\mu=\left(  z.\mu^{-2}\right)
^{\mathbf{1}}d\mu$, yields,%
\[
\int_{r\mathbf{1}}^{R\mathbf{1}}\lambda^{\delta}e^{-\lambda\mathbf{1}}%
\sigma\left(  \frac{z}{\lambda}\right)  d\lambda=\int_{z/R}^{z/r}e^{-\frac
{z}{\mu}\mathbf{1}}\left(  \frac{z}{\mu}\right)  ^{\delta}\sigma\left(
\mu\right)  \frac{\left(  z\right)  ^{\mathbf{1}}}{\mu^{2\mathbf{1}}}%
d\mu=z^{\delta+1}\int_{z/R}^{z/r}e^{-\frac{z}{\mu}\mathbf{1}}\frac
{\sigma\left(  \mu\right)  }{\mu^{\delta+2}}d\mu.
\]
\smallskip

Then%
\begin{align*}
z^{\delta+1}\int_{z/R}^{z/r}e^{-\frac{z}{\mu}\mathbf{1}}\frac{\left\vert
\sigma\left(  \mu\right)  \right\vert }{\mu^{\delta+2}}d\mu=\frac
{1}{z^{\mathbf{1}}}\int_{z/R}^{z/r}e^{-\frac{z}{\mu}\mathbf{1}}\frac
{z^{\delta+2}}{\mu^{\delta+2}}\left\vert \sigma\left(  \mu\right)  \right\vert
d\mu & \leq\frac{1}{z^{\mathbf{1}}}\left\Vert e^{-x\mathbf{1}}x^{\delta
+2}\right\Vert _{\infty}\int_{z/R}^{z/r}\left\vert \sigma\left(  \mu\right)
\right\vert d\mu\\
& \leq\frac{1}{z^{\mathbf{1}}}\left\Vert e^{-x\mathbf{1}}x^{\delta
+2}\right\Vert _{\infty}\left\Vert \sigma\right\Vert _{1}\\
& :part\text{ }1\Rightarrow\\
& =\frac{1}{z^{\mathbf{1}}}\left(  \delta+2\right)  ^{\delta+2}e^{-\left(
\delta+2\right)  \mathbf{1}}\left\Vert \sigma\right\Vert _{1}\\
& <\infty.
\end{align*}

Thus%
\[
J_{\delta}\left(  z\right)  =z^{\delta+1}\int_{\mathcal{O}_{\mathbf{1}}%
}e^{-\frac{z}{\mu}\mathbf{1}}\frac{\sigma\left(  \mu\right)  }{\mu^{\delta+2}%
}d\mu,\quad\delta\geq-2,
\]

and%
\[
\left\vert J_{\delta}\left(  z\right)  \right\vert \leq\left(  \delta
+2\right)  ^{\delta+2}e^{-\left(  \delta+2\right)  \mathbf{1}}\frac
{1}{z^{\mathbf{1}}}.
\]

\textbf{Part 3} Since $J_{\delta}\left(  z\right)  =z^{\delta+1}%
\int_{\mathcal{O}_{\mathbf{1}}}e^{-\frac{z}{\mu}\mathbf{1}}\frac{\sigma\left(
\mu\right)  }{\mu^{\delta+2}}d\mu$,%
\begin{align*}
& \int_{\mathcal{O}_{\mathbf{1}}}\left\vert J_{\delta}\right\vert \\
& \leq\int_{\mathcal{O}_{\mathbf{1}}}z^{\delta+1}\int_{\mathcal{O}%
_{\mathbf{1}}}e^{-\frac{z}{\mu}\mathbf{1}}\frac{\left\vert \sigma\left(
\mu\right)  \right\vert }{\mu^{\delta+2}}d\mu dz=\int_{\mathcal{O}%
_{\mathbf{1}}}\int_{\mathcal{O}_{\mathbf{1}}}z^{\delta+1}e^{-\frac{z}{\mu
}\mathbf{1}}dz\frac{\left\vert \sigma\left(  \mu\right)  \right\vert }%
{\mu^{\delta+2}}d\mu=\\
& :\zeta=z/\mu,\text{ }d\mu=-z.\zeta^{-2}d\zeta\Rightarrow\\
& =\int_{\mathcal{O}_{\mathbf{1}}}\int_{\mathcal{O}_{\mathbf{1}}}\left(
\mu.\zeta\right)  ^{\delta+1}e^{-\zeta\mathbf{1}}\mu^{\mathbf{1}}d\zeta
\frac{\left\vert \sigma\left(  \mu\right)  \right\vert }{\mu^{\delta+2}}%
d\mu=\int_{\mathcal{O}_{\mathbf{1}}}\int_{\mathcal{O}_{\mathbf{1}}}\mu
^{\delta+1}.\zeta^{\delta+1}e^{-\zeta\mathbf{1}}\mu^{\mathbf{1}}d\zeta
\frac{\left\vert \sigma\left(  \mu\right)  \right\vert }{\mu^{\delta+2}}%
d\mu=\\
& =\int_{\mathcal{O}_{\mathbf{1}}}\int_{\mathcal{O}_{\mathbf{1}}}\zeta
^{\delta+1}e^{-\zeta\mathbf{1}}d\zeta\left\vert \sigma\left(  \mu\right)
\right\vert d\mu=\left(  \delta+1\right)  !\left\Vert \sigma\right\Vert _{1}.
\end{align*}

Thus $J_{\delta}\in L^{1}$ (when $\delta>-2\mathbf{1}$).\medskip

\textbf{Part 4} Clearly $\sigma\in L^{1}$ so $J_{\delta}\left(  z\right)
=z^{\delta+1}\int_{\mathcal{O}_{\mathbf{1}}}e^{-\frac{z}{\mu}\mathbf{1}}%
\frac{\sigma\left(  \mu\right)  }{\mu^{\delta+2}}d\mu$ when $\delta\geq-2$.
Now if $\delta+2+\beta\geq\mathbf{0}$ i.e. $\alpha\geq\beta\geq\max\left\{
-\alpha,-\left(  \delta+2\right)  \right\}  $ we can write%
\[
J_{\delta}\left(  z\right)  =z^{\delta+1}\int_{\mathcal{O}_{\mathbf{1}}%
}e^{-\frac{z}{\mu}\mathbf{1}}\frac{\sigma\left(  \mu\right)  }{\mu^{\delta+2}%
}d\mu=\frac{1}{z^{\beta+1}}\int_{\mathcal{O}_{\mathbf{1}}}e^{-\frac{z}{\mu
}\mathbf{1}}\frac{z^{\delta+2+\beta}}{\mu^{\delta+2+\beta}}\mu^{\beta}%
\sigma\left(  \mu\right)  d\mu,
\]

which allows the estimate%
\begin{align*}
\left\vert J_{\delta}\left(  z\right)  \right\vert  & \leq\frac{1}{z^{\beta
+1}}\left(  \int_{\mathcal{O}_{\mathbf{1}}}\mu^{\beta}\sigma\left(
\mu\right)  d\mu\right)  \left\Vert e^{-x\mathbf{1}}x^{\delta+2+\beta
}\right\Vert _{\infty,\mathcal{O}_{\mathbf{1}}}\\
& =\frac{1}{z^{\beta+1}}\left(  \int_{\mathcal{O}_{\mathbf{1}}}\mu^{\beta
}\sigma\left(  \mu\right)  d\mu\right)  \left(  \delta+2+\beta\right)
^{\delta+2+\beta}e^{-\left(  \delta+2+\beta\right)  \mathbf{1}}.
\end{align*}

\textbf{Part 5} Specifically%
\[
\left\vert J_{\left(  n-2\right)  \mathbf{1}}\left(  z\right)  \right\vert
\leq\frac{1}{z^{\beta+1}}\left(  \int_{\mathcal{O}_{\mathbf{1}}}\mu^{\beta
}\sigma\left(  \mu\right)  d\mu\right)  \left(  n\mathbf{1}+\beta\right)
^{n\mathbf{1}+\beta}e^{-\left(  n\mathbf{1}+\beta\right)  \mathbf{1}},\text{
}\beta_{+}\leq n\mathbf{1}.
\]

\end{proof}

\begin{remark}
?? MOVE? \textbf{Idea} Proving $J_{\delta}^{\mathbf{1}}\left[  v\right]  $ is
continuous on each orthant by using the Fourier transform and the fact that
each has support on a separate orthant.
\begin{align*}
J_{\delta}^{\mathbf{1}}\left[  \sigma\right]  \left(  z\right)   &
=z^{\delta+1}\int_{\mathcal{O}_{\mathbf{1}}}e^{-\frac{z}{\mu}\mathbf{1}}%
\frac{\sigma\left(  \mu\right)  }{\mu^{\delta+2}}d\mu,\quad\delta
\geq-2\mathbf{1}.\\
& J_{\delta}^{\theta}\left[  \sigma\left(  \theta.\right)  \right]  \left(
\theta.z\right)  .
\end{align*}

\[
I_{\delta}\left[  \sigma\right]  \left(  x\right)  :=J_{\delta}^{\mathbf{1}%
}\left[  \sigma\right]  \left(  \theta.x\right)  ,\text{ }x\in\mathcal{O}%
_{\theta}.
\]

\begin{align*}
\widehat{I_{\delta}\left[  \sigma\right]  }\left(  \xi\right)   & =\left(
2\pi\right)  ^{-d/2}\int e^{-i\xi x}I_{\delta}\left[  \sigma\right]  \left(
x\right)  dx=\left(  2\pi\right)  ^{-d/2}\sum\limits_{\theta=-\mathbf{1}%
}^{\mathbf{1}}\int_{\mathcal{O}_{\theta}}e^{-i\xi x}I_{\delta}\left[
\sigma\right]  \left(  x\right)  dx=\\
& =\left(  2\pi\right)  ^{-d/2}\sum\limits_{\theta=-\mathbf{1}}^{\mathbf{1}%
}\int_{\mathcal{O}_{\theta}}e^{-i\xi x}J_{\delta}^{\mathbf{1}}\left[
\sigma\right]  \left(  \theta.x\right)  dx=\left(  2\pi\right)  ^{-d/2}%
\sum\limits_{\theta=-\mathbf{1}}^{\mathbf{1}}\int_{\mathcal{O}_{\mathbf{1}}%
}e^{-i\xi\theta.x}J_{\delta}^{\mathbf{1}}\left[  \sigma\right]  \left(
x\right)  dx=\\
& =\left(  2\pi\right)  ^{-d/2}\sum\limits_{\theta=-\mathbf{1}}^{\mathbf{1}%
}\int_{\mathcal{O}_{\mathbf{1}}}e^{-i\left(  \theta.\xi\right)  x}J_{\delta
}^{\mathbf{1}}\left[  \sigma\right]  \left(  x\right)  dx=\sum\limits_{\theta
=-\mathbf{1}}^{\mathbf{1}}\widehat{J_{\delta}^{\mathbf{1}}\left[
\sigma\right]  }\left(  \theta.\xi\right)  =\\
& =\sum\limits_{\theta=-\mathbf{1}}^{\mathbf{1}}\int_{\mathcal{O}_{\mathbf{1}%
}}\left(  H\left(  \zeta\right)  \zeta^{\delta+1}e^{-\zeta\mathbf{1}}\right)
^{\wedge}\left(  \theta.\xi.\mu\right)  \sigma\left(  \mu\right)  d\mu\\
& =\sum\limits_{\theta=-\mathbf{1}}^{\mathbf{1}}\int_{\mathcal{O}_{\mathbf{1}%
}}\left(  h\left(  \theta.x\right)  \right)  ^{\wedge}\left(  \xi.\mu\right)
\sigma\left(  \mu\right)  d\mu=\int_{\mathcal{O}_{\mathbf{1}}}\left(
\sum\limits_{\theta=-\mathbf{1}}^{\mathbf{1}}h\left(  \theta.x\right)
\right)  ^{\wedge}\left(  \xi.\mu\right)  \sigma\left(  \mu\right)  d\mu.
\end{align*}

Observe that $\sum\limits_{\theta=-\mathbf{1}}^{\mathbf{1}}h\left(
\theta.x\right)  \in PWC^{\infty}$.
\end{remark}

In the next lemma equation \ref{X06} expresses the smooth function value at
$x$ in terms of its derivatives in the orthant $x+\mathcal{O}_{\theta}$. This
is our basic integral representation in an orthant.

\begin{lemma}
\label{Lem_SmthFuncIntegRepInOrthant}\textbf{Orthant integral representation}
Suppose $\mathcal{O}_{\theta}=\theta.\mathcal{O}_{\mathbf{1}}$ is any orthant.
Suppose $\alpha\geq\mathbf{1}$, $f\in C_{B??0}^{\left(  \alpha\right)
}\left(  \overline{\mathcal{O}_{\theta}}\right)  $ and $v\in L_{\overline
{\mathcal{O}_{\theta}}}^{1}$ satisfies: $\int_{\mathcal{O}_{\theta}}v=1$ and
$x^{\beta}v\in L_{\overline{\mathcal{O}_{\theta}}}^{1}$ when $\mathbf{0}%
\leq\beta\leq\alpha$. Then:

\begin{enumerate}
\item When $\gamma\leq\alpha-1$, $D^{\gamma}f\in C_{0}^{\left(  \alpha
-\gamma\right)  }\left(  \overline{\mathcal{O}_{\theta}}\right)  $ and%
\begin{equation}
D^{\gamma}f\left(  x\right)  =\tfrac{1}{\left(  \alpha-\gamma-1\right)  !}%
\sum\limits_{\mathbf{0}\leq\beta\leq\alpha-\gamma}\tbinom{\alpha-\gamma}%
{\beta}\int\limits_{x+\mathcal{O}_{\theta}}J_{\alpha-\gamma-2}^{\theta}\left[
v_{\beta}\right]  \left(  y-x\right)  D^{\gamma+\beta}f\left(  y\right)
dy,\quad x\in\overline{\mathcal{O}_{\theta}},\label{X06}%
\end{equation}

where we define%
\begin{equation}
v_{\beta}\left(  x\right)  :=\left(  -x\right)  ^{\beta}v\left(  x\right)
,\quad\forall\beta\in\mathbb{Z}^{d},\label{X49}%
\end{equation}

and%
\begin{equation}
J_{\lambda}^{\theta}\left[  v\right]  \left(  z\right)  :=J_{\lambda}\left[
v\left(  \theta.\right)  \right]  \left(  \theta.z\right)  ,\quad\lambda
\geq-\mathbf{1},\label{X05}%
\end{equation}

and $J_{\lambda}^{\theta}\left[  v_{\beta}\right]  $ can be written in the
integral form \ref{X48}.\medskip

\textbf{Now further assume that} $v\in L_{0}^{1}\left(  \mathcal{O}_{\theta
}\right)  $. Then:\medskip

\item When $\lambda\geq-\mathbf{1}$ we have $J_{\lambda}^{\theta}\left[
v_{\beta}\right]  \in PWC^{\infty}$, $J_{\lambda}^{\theta}\left[  v_{\beta
}\right]  \in C^{\infty}\left(  -\overline{\mathcal{O}_{\theta}}\right)  $ and
each derivative decreases exponentially to zero at infinity.

\item When $\lambda\geq\mathbf{0}$ we have $J_{\lambda}^{\theta}\left[
v_{\beta}\right]  \in C_{B}^{\left(  \lambda\right)  }\left(  \mathbb{R}%
^{d}\right)  $.

??? \textbf{ADD more analogous results from Lemma} \ref{Lem_property_J}?

\item Suppose $\delta\geq-\mathbf{1}$ is a multi-integer and
$\operatorname*{supp}v\subset R\left[  c,d\right]  \subset\mathcal{O}_{\theta
}$ where $c.<d$. Then%
\begin{equation}
\left\vert D^{\gamma}J_{\delta}^{\theta}\left[  v\right]  \left(  z\right)
\right\vert \leq\left\Vert \frac{v\left(  y\right)  }{y^{\gamma+1}}\right\Vert
_{1}q_{\gamma}^{\delta}\left(  \frac{\theta.z}{\min\left\{  \theta
.c,\theta.d\right\}  }\right)  e^{-\frac{\theta z}{\max\left\{  \theta
.c,\theta.d\right\}  }},\quad z\in\mathcal{O}_{\theta},\label{X88}%
\end{equation}

where $q_{\gamma}^{\delta}\geq0$ is a tensor product polynomial of degree
$\delta+1$ which is independent of $v$ and is given by \ref{X41}.

\item Suppose $\delta\geq-\mathbf{1}$ is a multi-integer. Then $\forall\gamma
$, $D^{\gamma}J_{\delta}^{\theta}\in L^{1}\left(  \mathcal{O}_{\theta}\right)
$ and%
\[
\int_{\mathcal{O}_{\theta}}\left\vert D^{\gamma}J_{\delta}^{\theta}\left[
v\right]  \right\vert \leq\left\Vert \frac{v}{x^{\gamma}}\right\Vert _{1}%
\int_{\mathcal{O}_{\mathbf{1}}}q_{\gamma}^{\delta}\left(  y\right)
e^{-y\mathbf{1}}dy<\infty.
\]

\end{enumerate}
\end{lemma}

\begin{proof}
The approach is to first prove the lemma for the orthant $\theta=1$ and then
use Definition \ref{Def_funcs_and_orthants} to prove for the result for the
other orthants.\medskip

\fbox{We first prove the lemma for $\theta=\mathbf{1}$} and then use
Definition \ref{Def_funcs_and_orthants} for the other orthants.

Set $g\left(  \tau\right)  :=e^{-\tau\mathbf{1}}f\left(  x+\tau.\eta\right)  $
for fixed $x\in\overline{\mathcal{O}_{\mathbf{1}}}$ and $\eta\in
\mathcal{O}_{\mathbf{1}}$. Then, since ?? $f\in C_{??0}^{\left(
\alpha\right)  }\left(  \overline{\mathcal{O}_{\mathbf{1}}}\right)  $ implies
$f\in C_{BP}^{\left(  \alpha\right)  }\left(  \overline{\mathcal{O}%
_{\mathbf{1}}}\right)  $, we have $g\in C_{BP}^{\left(  \alpha\right)
}\left(  \overline{\mathcal{O}_{\mathbf{1}}}\right)  $ and $g\rightarrow0$ ??
exponentially at infinity equation \ref{X16} of Lemma \ref{Lem_integ_repres}
implies%
\begin{equation}
f\left(  x\right)  =\tfrac{\left(  -1\right)  ^{\left\vert \alpha\right\vert
}}{\left(  \alpha-1\right)  !}\int_{\mathcal{O}_{\mathbf{1}}}\tau^{\alpha
-1}D_{\tau}^{\alpha}\left(  e^{-\tau\mathbf{1}}f\left(  x+\tau.\eta\right)
\right)  d\tau,\text{\quad}\left\{
\begin{array}
[c]{l}%
\alpha\geq\mathbf{1},\\
x\in\overline{\mathcal{O}_{\mathbf{1}}},\text{ }\eta\in\mathcal{O}%
_{\mathbf{1}}.
\end{array}
\right. \label{X23}%
\end{equation}

Now%
\begin{align*}
D_{\tau}^{\alpha}\left(  e^{-\tau\mathbf{1}}f\left(  x+\tau.\eta\right)
\right)   &  =\sum_{\beta\leq\alpha}\tbinom{\alpha}{\beta}\left(  D_{\tau
}^{\alpha-\beta}e^{-\tau\mathbf{1}}\right)  D_{\tau}^{\beta}\left(  f\left(
x+\tau.\eta\right)  \right) \\
&  =\sum_{\beta\leq\alpha}\left(  -1\right)  ^{\left\vert \alpha
-\beta\right\vert }\tbinom{\alpha}{\beta}e^{-\tau\mathbf{1}}D_{\tau}^{\beta
}\left(  f\left(  x+\tau.\eta\right)  \right) \\
&  =\sum_{\beta\leq\alpha}\left(  -1\right)  ^{\left\vert \alpha\right\vert
-\left\vert \beta\right\vert }\tbinom{\alpha}{\beta}e^{-\tau\mathbf{1}}%
D_{\tau}^{\beta}\left(  f\left(  x+\tau.\eta\right)  \right) \\
&  =\left(  -1\right)  ^{\left\vert \alpha\right\vert }e^{-\tau\mathbf{1}}%
\sum_{\beta\leq\alpha}\left(  -1\right)  ^{\left\vert \beta\right\vert
}\tbinom{\alpha}{\beta}D_{\tau}^{\beta}\left(  f\left(  x+\tau.\eta\right)
\right) \\
&  =\left(  -1\right)  ^{\left\vert \alpha\right\vert }e^{-\tau\mathbf{1}}%
\sum_{\beta\leq\alpha}\left(  -1\right)  ^{\left\vert \beta\right\vert
}\tbinom{\alpha}{\beta}\eta^{\beta}\left(  D^{\beta}f\right)  \left(
x+\tau.\eta\right) \\
&  =\left(  -1\right)  ^{\left\vert \alpha\right\vert }e^{-\tau\mathbf{1}}%
\sum_{\beta\leq\alpha}\tbinom{\alpha}{\beta}\left(  -\eta\right)  ^{\beta
}\left(  D^{\beta}f\right)  \left(  x+\tau.\eta\right)  ,
\end{align*}

so that \ref{X23} becomes%
\begin{align}
f\left(  x\right)   & =\tfrac{\left(  -1\right)  ^{\left\vert \alpha
\right\vert }}{\left(  \alpha-1\right)  !}\int_{\mathcal{O}_{\mathbf{1}}}%
\tau^{\alpha-1}\left(  -1\right)  ^{\left\vert \alpha\right\vert }%
e^{-\tau\mathbf{1}}\sum_{\beta\leq\alpha}\tbinom{\alpha}{\beta}\left(
-\eta\right)  ^{\beta}\left(  D^{\beta}f\right)  \left(  x+\tau.\eta\right)
d\tau\nonumber\\
& =\tfrac{1}{\left(  \alpha-1\right)  !}\int_{\mathcal{O}_{\mathbf{1}}}%
\tau^{\alpha-1}e^{-\tau\mathbf{1}}\sum_{\beta\leq\alpha}\tbinom{\alpha}{\beta
}\left(  -\eta\right)  ^{\beta}\left(  D^{\beta}f\right)  \left(  x+\tau
.\eta\right)  d\tau\nonumber\\
& =\tfrac{1}{\left(  \alpha-1\right)  !}\sum_{\beta\leq\alpha}\tbinom{\alpha
}{\beta}\int_{\mathcal{O}_{\mathbf{1}}}\tau^{\alpha-1}e^{-\tau\mathbf{1}%
}\left(  -\eta\right)  ^{\beta}\left(  D^{\beta}f\right)  \left(  x+\tau
.\eta\right)  d\tau.\label{X011}%
\end{align}

Now choose a function $v\in L^{1}\left(  \mathcal{O}_{\mathbf{1}}\right)  $
such that $\int_{\mathcal{O}_{\mathbf{1}}}v=1$. Then using \ref{X23} we get%
\begin{align*}
f\left(  x\right)   & =f\left(  x\right)  \int_{\mathcal{O}_{\mathbf{1}}%
}v\left(  \eta\right)  d\eta=\int_{\mathcal{O}_{\mathbf{1}}}f\left(  x\right)
v\left(  \eta\right)  d\eta=\\
& =\int_{\mathcal{O}_{\mathbf{1}}}\tfrac{1}{\left(  \alpha-1\right)  !}%
\sum_{\beta\leq\alpha}\tbinom{\alpha}{\beta}\int_{\mathcal{O}_{\mathbf{1}}%
}\tau^{\alpha-1}e^{-\tau\mathbf{1}}\left(  -\eta\right)  ^{\beta}\left(
D^{\beta}f\right)  \left(  x+\tau.\eta\right)  d\tau\text{ }v\left(
\eta\right)  d\eta\\
& =\tfrac{1}{\left(  \alpha-1\right)  !}\sum_{\beta\leq\alpha}\tbinom{\alpha
}{\beta}\int_{\mathcal{O}_{\mathbf{1}}}\int_{\mathcal{O}_{\mathbf{1}}}%
\tau^{\alpha-1}e^{-\tau\mathbf{1}}v_{\beta}\left(  \eta\right)  \left(
D^{\beta}f\right)  \left(  x+\tau.\eta\right)  d\tau d\eta.
\end{align*}

The next step is to show that the integrand is absolutely convergent i.e. it
is $L^{1}$. In fact, since $x,\tau\in\mathcal{O}_{\mathbf{1}}$,%
\[
\int_{\mathcal{O}_{\mathbf{1}}}\left\vert v_{\beta}\left(  \eta\right)
\left(  D^{\beta}f\right)  \left(  x+\tau.\eta\right)  \right\vert d\eta
\leq\left\Vert v_{\beta}\right\Vert _{1,\mathcal{O}_{\mathbf{1}}}\left\Vert
D^{\beta}f\right\Vert _{\infty,\mathcal{O}_{\mathbf{1}}},
\]

and so%
\begin{align}
\int_{\mathcal{O}_{\mathbf{1}}}\int_{\mathcal{O}_{\mathbf{1}}}\left\vert
\tau^{\alpha-1}e^{-\tau\mathbf{1}}v_{\beta}\left(  \eta\right)  \left(
D^{\beta}f\right)  \left(  x+\tau.\eta\right)  \right\vert d\tau d\eta &
=\int_{\mathcal{O}_{\mathbf{1}}}\tau^{\alpha-1}e^{-\tau\mathbf{1}}%
\int_{\mathcal{O}_{\mathbf{1}}}\left\vert v_{\beta}\left(  \eta\right)
\left(  D^{\beta}f\right)  \left(  x+\tau.\eta\right)  \right\vert d\tau
d\eta\nonumber\\
& \leq\left(  \int_{\mathcal{O}_{\mathbf{1}}}\tau^{\alpha-1}e^{-\tau
\mathbf{1}}d\tau\right)  \left\Vert v_{\beta}\right\Vert _{1,\mathcal{O}%
_{\mathbf{1}}}\left\Vert D^{\beta}f\right\Vert _{\infty,\mathcal{O}%
_{\mathbf{1}}}\nonumber\\
& =\left(  \alpha-1\right)  !\left\Vert v_{\beta}\right\Vert _{1,\mathcal{O}%
_{\mathbf{1}}}\left\Vert D^{\beta}f\right\Vert _{\infty,\mathcal{O}%
_{\mathbf{1}}}\label{X71}\\
& <\infty,\nonumber
\end{align}

The change of variables $y=x+\tau.\eta$, $\lambda=\tau$ is a $C^{1}$ bijection
from $\mathcal{O}_{\mathbf{1}}\times\mathcal{O}_{\mathbf{1}}$ to $\left(
x+\mathcal{O}_{\mathbf{1}}\right)  \times\mathcal{O}_{\mathbf{1}}$ with
Jacobian $\lambda^{-\mathbf{1}}$ e.g. p502 in Chapter 15, Section J of Jones
\cite{Jones2011}, and since the integrand is $L^{1}$:%
\begin{align}
f\left(  x\right)   &  =\tfrac{1}{\left(  \alpha-1\right)  !}\sum_{\beta
\leq\alpha}\tbinom{\alpha}{\beta}\int_{\mathcal{O}_{\mathbf{1}}}%
\int_{\mathcal{O}_{\mathbf{1}}}\tau^{\alpha-1}e^{-\tau\mathbf{1}}v_{\beta
}\left(  \eta\right)  \left(  D^{\beta}f\right)  \left(  x+\tau.\eta\right)
d\tau d\eta\nonumber\\
&  =\tfrac{1}{\left(  \alpha-1\right)  !}\sum_{\beta\leq\alpha}\tbinom{\alpha
}{\beta}\int_{x+\mathcal{O}_{\mathbf{1}}}\int_{\mathcal{O}_{\mathbf{1}}%
}\lambda^{\alpha-1}e^{-\lambda\mathbf{1}}\left(  -\frac{y-x}{\lambda}\right)
^{\beta}D^{\beta}f\left(  y\right)  v\left(  \frac{y-x}{\lambda}\right)
\frac{d\lambda}{\lambda^{1}}dy\nonumber\\
&  =\tfrac{1}{\left(  \alpha-1\right)  !}\sum_{\beta\leq\alpha}\tbinom{\alpha
}{\beta}\int_{x+\mathcal{O}_{\mathbf{1}}}\int_{\mathcal{O}_{\mathbf{1}}%
}\lambda^{\alpha-2}e^{-\lambda\mathbf{1}}\left(  -\frac{y-x}{\lambda}\right)
^{\beta}v\left(  \frac{y-x}{\lambda}\right)  d\lambda\text{ }D^{\beta}f\left(
y\right)  dy\nonumber\\
&  =\tfrac{1}{\left(  \alpha-1\right)  !}\sum_{\beta\leq\alpha}\tbinom{\alpha
}{\beta}\int_{x+\mathcal{O}_{\mathbf{1}}}J_{\alpha-2}\left[  v_{\beta}\right]
\left(  y-x\right)  D^{\beta}f\left(  y\right)  dy,\label{X24}%
\end{align}

where $v_{\beta}\in L^{1}\left(  \mathcal{O}_{\mathbf{1}}\right)  $ is assumed
in the statement of this lemma and $J_{\alpha-2}^{\mathbf{1}}\left[  v_{\beta
}\right]  $ is given by \ref{X39} of Lemma \ref{Lem_property_J}.

If $f\in C_{0}^{\left(  \alpha\right)  }\left(  \overline{\mathcal{O}%
_{\mathbf{1}}}\right)  $ and $\gamma\leq\alpha$ then $D^{\gamma}f\in
C_{0}^{\left(  \alpha-\gamma\right)  }\left(  \overline{\mathcal{O}%
_{\mathbf{1}}}\right)  $ and so by \ref{X24},%
\begin{equation}
D^{\gamma}f\left(  x\right)  =\tfrac{1}{\left(  \alpha-\gamma-1\right)  !}%
\sum_{\beta\leq\alpha-\gamma}\tbinom{\alpha-\gamma}{\beta}\int_{x+\mathcal{O}%
_{\mathbf{1}}}J_{\alpha-\gamma-2}\left[  v_{\beta}\right]  \left(  y-x\right)
D^{\gamma+\beta}f\left(  y\right)  dy,\label{X46}%
\end{equation}

which is \ref{X06} when $\theta=\mathbf{1}$.\medskip

\fbox{Now consider arbitrary $\theta$.} From Definition
\ref{Def_funcs_and_orthants}, in general $f\in C_{0}^{\left(  \alpha\right)
}\left(  \overline{\mathcal{O}_{\theta}}\right)  $ and $v\in L^{1}\left(
\mathcal{O}_{\theta}\right)  $ iff $f\left(  \theta.\right)  \in
C_{0}^{\left(  \alpha\right)  }\left(  \overline{\mathcal{O}_{\mathbf{1}}%
}\right)  $ and $v\left(  \theta.\right)  \in L^{1}\left(  \mathcal{O}%
_{\mathbf{1}}\right)  $. Also $\int v\left(  \theta.\right)  =\int v=1$ and%
\begin{equation}
\left(  v\left(  \theta.\right)  \right)  _{\beta}\left(  x\right)  =\left(
-x\right)  ^{\beta}v\left(  \theta.x\right)  =\theta^{\beta}\left(
-\theta.x\right)  ^{\beta}v\left(  \theta.x\right)  =\theta^{\beta}v_{\beta
}\left(  \theta.x\right)  .\label{X81}%
\end{equation}

Thus, if $x\in\mathcal{O}_{\mathbf{1}}$, \ref{X46} implies%
\begin{align*}
f\left(  \theta.x\right)   & =\tfrac{1}{\left(  \alpha-1\right)  !}\sum
_{\beta\leq\alpha}\tbinom{\alpha}{\beta}\int_{x+\mathcal{O}_{\mathbf{1}}%
}J_{\alpha-2}\left[  \left(  v\left(  \theta.\right)  \right)  _{\beta
}\right]  \left(  y-x\right)  D^{\beta}\left(  f\left(  \theta.y\right)
\right)  dy.\\
& =\tfrac{1}{\left(  \alpha-1\right)  !}\sum_{\beta\leq\alpha}\tbinom{\alpha
}{\beta}\int_{x+\mathcal{O}_{\mathbf{1}}}\theta^{\beta}J_{\alpha-2}\left[
v_{\beta}\left(  \theta.\right)  \right]  \left(  y-x\right)  \left(
D^{\beta}f\right)  \left(  \theta.y\right)  dy.
\end{align*}

The change of variables: $z=\theta.y,$ $dz=\left\vert \theta^{\mathbf{1}%
}\right\vert dy=dy$ gives%
\[
f\left(  \theta.x\right)  =\tfrac{1}{\left(  \alpha-1\right)  !}\sum
_{\beta\leq\alpha}\tbinom{\alpha}{\beta}\int_{\theta.x+\mathcal{O}_{\theta}%
}\theta^{\beta}J_{\alpha-2}\left[  v_{\beta}\left(  \theta.\right)  \right]
\left(  \theta.z-x\right)  D^{\beta}f\left(  z\right)  dz,
\]

which is equivalent to: if $x\in\mathcal{O}_{\theta}$ then using \ref{X05},%
\begin{align}
f\left(  x\right)   & =\tfrac{1}{\left(  \alpha-1\right)  !}\sum_{\beta
\leq\alpha}\tbinom{\alpha}{\beta}\int_{x+\mathcal{O}_{\theta}}\theta^{\beta
}J_{\alpha-2}\left[  v_{\beta}\left(  \theta.\right)  \right]  \left(
\theta.z-\theta.x\right)  D^{\beta}f\left(  z\right)  dz\nonumber\\
& =\tfrac{1}{\left(  \alpha-1\right)  !}\sum_{\beta\leq\alpha}\tbinom{\alpha
}{\beta}\int_{x+\mathcal{O}_{\theta}}\theta^{\beta}J_{\alpha-2}\left[
v_{\beta}\left(  \theta.\right)  \right]  \left(  \theta.\left(  z-x\right)
\right)  D^{\beta}f\left(  z\right)  dz\nonumber\\
& =\tfrac{1}{\left(  \alpha-1\right)  !}\sum_{\beta\leq\alpha}\tbinom{\alpha
}{\beta}\int_{x+\mathcal{O}_{\theta}}J_{\alpha-2}^{\theta}\left[  v_{\beta
}\right]  \left(  z-x\right)  D^{\beta}f\left(  z\right)  dz,\label{X37}%
\end{align}

where, when $z\in-\mathcal{O}_{\theta}$, by part 2 of Lemma
\ref{Lem_Jd[f]_bnd_f_in_L1},%
\begin{align}
J_{\alpha-2}^{\theta}\left[  v_{\beta}\right]  \left(  z\right)
=\theta^{\beta}J_{\alpha-2}\left[  v_{\beta}\left(  \theta.\right)  \right]
\left(  \theta.z\right)   &  =\theta^{\beta}\int_{\mathcal{O}_{\mathbf{1}}%
}\lambda^{\alpha-2}e^{-\lambda\mathbf{1}}\left(  v_{\beta}\left(
\theta.\right)  \right)  \left(  \frac{\theta.z}{\lambda}\right)
d\lambda\nonumber\\
&  =\int_{\mathcal{O}_{\theta}}\lambda^{\alpha-2}e^{-\lambda\mathbf{1}%
}v_{\beta}\left(  \frac{z}{\lambda}\right)  d\lambda,\label{X47}%
\end{align}

for any $v\in L_{\overline{\mathcal{O}_{\theta}}}^{1}$. In particular when
$\beta=\mathbf{0}$ and $\lambda=\alpha-2\geq-\mathbf{1}$,%
\[
J_{\lambda}^{\theta}\left[  v\right]  \left(  z\right)  =J_{\lambda}\left[
v\left(  \theta.\right)  \right]  \left(  \theta.z\right)  .
\]

If $\mathbf{0}<\gamma\leq\alpha$ then $D^{\gamma}f\in C_{0}^{\left(
\alpha-\gamma\right)  }\left(  \overline{\mathcal{O}_{\theta}}\right)  $ and
equation \ref{X37} implies the derivative formula%
\[
D^{\gamma}f\left(  x\right)  =\tfrac{1}{\left(  \alpha-1\right)  !}\sum
_{\beta\leq\alpha-\gamma}\tbinom{\alpha-\gamma}{\beta}\int_{x+\mathcal{O}%
_{\theta}}J_{\alpha-\gamma-2}^{\theta}\left[  v_{\beta}\right]  \left(
z-x\right)  D^{\gamma+\beta}f\left(  z\right)  dz,
\]

which is \ref{X06} as claimed, and \ref{X47} implies%
\begin{equation}
J_{\alpha-\gamma-2}^{\theta}\left[  v_{\beta}\right]  \left(  z\right)
=\int_{\mathcal{O}_{\theta}}\lambda^{\alpha-\gamma-2}e^{-\lambda\mathbf{1}%
}v_{\beta}\left(  \frac{z}{\lambda}\right)  d\lambda,\quad\beta\leq
\alpha-\gamma,\label{X48}%
\end{equation}

which proves part 1 for all orthants.\medskip

\fbox{\textbf{Part 2}} Apply \ref{X05}, \ref{X49} to parts 2 of Lemma
\ref{Lem_property_J}.\medskip

\fbox{\textbf{Part 3}} Apply \ref{X05}, \ref{X49} to part 6 of Lemma
\ref{Lem_property_J}.\medskip

\fbox{\textbf{Part 4}} From part 1,%
\[
D^{\gamma}J_{\lambda}^{\theta}\left[  v\right]  \left(  z\right)
=\theta^{\gamma}\left(  D^{\gamma}J_{\lambda}^{\mathbf{1}}\left[  v\left(
\theta.\right)  \right]  \right)  \left(  \theta.z\right)  ,\quad
z\in\mathcal{O}_{\theta}.
\]

Now inequality \ref{X32} is: if $\operatorname*{supp}\sigma\subset R\left[
a,b\right]  \subset\mathcal{O}_{\mathbf{1}}$ then%
\[
\left\vert D^{\gamma}J_{\delta}^{\mathbf{1}}\left[  \sigma\right]  \left(
z\right)  \right\vert \leq\left\Vert \frac{\sigma}{x^{\gamma+1}}\right\Vert
_{1}q_{\gamma}^{\delta}\left(  \frac{z}{a}\right)  e^{-\frac{z}{b}\mathbf{1}%
},\quad z\in\mathcal{O}_{\mathbf{1}}.
\]

Now, noting \ref{X85} ?? BELOW ??,
\begin{align*}
\operatorname*{supp}v\left(  \theta.\right)   & =\theta.\operatorname*{supp}%
v\subseteq\theta.R\left[  c,d\right]  =R\left[  \theta.c,\theta.d\right]  =\\
& =R\left[  \min\left\{  \theta.c,\theta.d\right\}  ,\max\left\{
\theta.c,\theta.d\right\}  \right]  ,
\end{align*}

so for $z\in\mathcal{O}_{\theta}$,
\begin{align*}
\left\vert D^{\gamma}J_{\lambda}^{\theta}\left[  v\right]  \left(  z\right)
\right\vert  & =\left\vert \theta^{\gamma}\left(  D^{\gamma}J_{\lambda
}^{\mathbf{1}}\left[  v\left(  \theta.\right)  \right]  \right)  \left(
\theta.z\right)  \right\vert \\
& =\left\vert \left(  D^{\gamma}J_{\lambda}^{\mathbf{1}}\left[  v\left(
\theta.\right)  \right]  \right)  \left(  \theta.z\right)  \right\vert \\
& \leq\left\Vert \frac{v\left(  \theta.y\right)  }{y^{\gamma+1}}\right\Vert
_{1}q_{\gamma}^{\delta}\left(  \frac{\theta.z}{\min\left\{  \theta
.c,\theta.d\right\}  }\right)  e^{-\frac{\theta.z}{\max\left\{  \theta
.c,\theta.d\right\}  }\mathbf{1}}\\
& =\left\Vert \frac{v}{y^{\gamma+1}}\right\Vert _{1}q_{\gamma}^{\delta}\left(
\frac{\theta.z}{\min\left\{  \theta.c,\theta.d\right\}  }\right)
e^{-\frac{\theta z}{\max\left\{  \theta.c,\theta.d\right\}  }}.
\end{align*}
\medskip

\fbox{\textbf{Part 5}} From part 4 of Lemma \ref{Lem_property_J},%
\[
\int_{\mathcal{O}_{\mathbf{1}}}\left\vert D^{\gamma}J_{\lambda}^{\mathbf{1}%
}\left[  \sigma\right]  \right\vert \leq\left\Vert \frac{\sigma\left(
x\right)  }{x^{\gamma}}\right\Vert _{1,\mathcal{O}_{\mathbf{1}}}%
\int_{\mathcal{O}_{\mathbf{1}}}q_{\gamma}^{\lambda}\left(  y\right)
e^{-y\mathbf{1}}dy,\quad\forall\gamma.
\]

From \ref{X05}, $J_{\lambda}^{\theta}\left[  v\right]  \left(  z\right)
:=J_{\lambda}^{\mathbf{1}}\left[  v\left(  \theta.\right)  \right]  \left(
\theta.z\right)  $, so that%
\[
D^{\gamma}J_{\lambda}^{\theta}\left[  v\right]  \left(  z\right)
=\theta^{\gamma}\left(  D^{\gamma}J_{\lambda}^{\mathbf{1}}\left[  v\left(
\theta.\right)  \right]  \right)  \left(  \theta.z\right)  ,
\]

and hence%
\begin{align*}
\int_{\mathcal{O}_{\mathbf{1}}}\left\vert D^{\gamma}J_{\lambda}^{\theta
}\left[  v\right]  \right\vert  & =\int_{\mathcal{O}_{\mathbf{1}}}\left\vert
\theta^{\gamma}\left(  D^{\gamma}J_{\lambda}^{\mathbf{1}}\left[  v\left(
\theta.\right)  \right]  \right)  \left(  \theta.z\right)  \right\vert dz\\
& =\int_{\mathcal{O}_{\theta}}\left\vert D^{\gamma}J_{\lambda}^{\mathbf{1}%
}\left[  v\left(  \theta.\right)  \right]  \left(  z\right)  \right\vert dz\\
& \leq\left\Vert \frac{\sigma\left(  \theta.x\right)  }{x^{\gamma}}\right\Vert
_{1,\mathcal{O}_{\mathbf{1}}}\int_{\mathcal{O}_{\mathbf{1}}}q_{\gamma
}^{\lambda}\left(  y\right)  e^{-y\mathbf{1}}dy\\
& =\left\Vert \frac{\sigma\left(  x\right)  }{x^{\gamma}}\right\Vert
_{1,\mathcal{O}_{\theta}}\int_{\mathcal{O}_{\mathbf{1}}}q_{\gamma}^{\lambda
}\left(  y\right)  e^{-y\mathbf{1}}dy.
\end{align*}

\end{proof}

?? FINISH!

\begin{lemma}
\label{Lem_convol_PWsmth}Suppose $g\in C_{B}^{\left(  0\right)  }\left(
\overline{\mathcal{O}_{\mathbf{1}}}\right)  $ and $g=0$ on $\mathbb{R}%
^{d}\setminus\overline{\mathcal{O}_{\mathbf{1}}}$. Suppose $f\in C^{\left(
1\right)  }\left(  \overline{\mathcal{O}_{\mathbf{1}}}\right)  $, $f=0$ on
$\mathbb{R}^{d}\setminus\overline{\mathcal{O}_{\mathbf{1}}}$ and $f$ and all
$D_{i}f$ tend to zero exponentially at infinity on $\overline{\mathcal{O}%
_{\mathbf{1}}}$.

Then $f\ast g\in C_{B}^{\left(  0\right)  }\cap L^{1}$.
\end{lemma}

\begin{proof}
??? CONSIDER $f\in L_{loc}^{1}\cap L^{\infty}$ and $g\in C^{\left(  0\right)
}$ implies ?? $f\ast g\in C^{\left(  0\right)  }$ ??. Approx. discontin. $g$
using mollifier.

GIVE REF! $f$ is piecewise smooth. Write $g$ as the sum of piecewise smooth
functions with each term contributing to only a bounded region of the
convolution. ETC. When does sum of continuous functions converge to a
continuous function.

Clearly $f,g\in L^{1}$ so $f\ast g\in L^{1}$.%
\begin{align*}
\left\vert \left(  f\ast g\right)  \left(  x+h\right)  -\left(  f\ast
g\right)  \left(  x\right)  \right\vert  & =\left\vert \left(  f\left(
\cdot+h\right)  \ast g\right)  \left(  x\right)  -\left(  f\ast g\right)
\left(  x\right)  \right\vert \\
& =\left\vert \left(  \left(  f\left(  \cdot+h\right)  -f\right)  \ast
g\right)  \left(  x\right)  \right\vert \\
& \leq\left\Vert f\left(  \cdot+h\right)  -f\right\Vert _{1}\left\Vert
g\right\Vert _{\infty}\\
& \leq\left(  \left\Vert f\left(  \cdot+h\right)  \right\Vert _{1}+\left\Vert
f\right\Vert _{1}\right)  \left\Vert g\right\Vert _{\infty}\\
& =2\left\Vert f\right\Vert _{1}\left\Vert g\right\Vert _{\infty}.
\end{align*}

\begin{align*}
g\left(  h\right)   & =\int\left\vert f\left(  z+h\right)  -f\left(  z\right)
\right\vert dz=\int\left\vert f\left(  z+\frac{h}{2}\right)  -f\left(
z-\frac{h}{2}\right)  \right\vert dz.\\
g\left(  \pi h\right)   & =\int\left\vert f\left(  z+\pi\frac{h}{2}\right)
-f\left(  z-\pi\frac{h}{2}\right)  \right\vert dz=\int\left\vert f\left(
\pi\zeta+\pi\frac{h}{2}\right)  -f\left(  \pi\zeta-\pi\frac{h}{2}\right)
\right\vert d\zeta=\\
& =\int\left\vert \left(  f\pi\right)  \left(  z+\frac{h}{2}\right)  -\left(
f\pi\right)  \left(  z-\frac{h}{2}\right)  \right\vert dz.
\end{align*}

Consider first $h.>\mathbf{0}$:
\begin{align*}
\int &  \left\vert f\left(  z+h\right)  -f\left(  z\right)  \right\vert dz\\
&  =\int_{\mathcal{O}_{\mathbf{1}}}\left\vert f\left(  z+h\right)  -f\left(
z\right)  \right\vert dz+\int_{\left(  \mathcal{O}_{\mathbf{1}}-h\right)
\backslash\mathcal{O}_{\mathbf{1}}}\left\vert f\left(  z+h\right)  -f\left(
z\right)  \right\vert dz\\
&  =\int_{\mathcal{O}_{\mathbf{1}}}\left\vert f\left(  z+h\right)  -f\left(
z\right)  \right\vert dz+\int_{\left(  \mathcal{O}_{\mathbf{1}}-h\right)
\backslash\mathcal{O}_{\mathbf{1}}}\left\vert f\left(  z+h\right)  \right\vert
dz\\
&  =\int_{\mathcal{O}_{\mathbf{1}}}\left\vert f\left(  z+h\right)  -f\left(
z\right)  \right\vert dz+\int_{\mathcal{O}_{\mathbf{1}}\backslash\left(
\mathcal{O}_{\mathbf{1}}+h\right)  }\left\vert f\left(  z\right)  \right\vert
dz.
\end{align*}

But $\left\vert f\left(  z\right)  \right\vert \leq b^{\left(  0\right)
}e^{-a^{\left(  0\right)  }z},$ $\left\vert D_{i}f\left(  z\right)
\right\vert \leq b^{\left(  i\right)  }e^{-a^{\left(  i\right)  }z},$
$z\in\mathcal{O}_{\mathbf{1}}$ so%
\begin{align*}
\int &  \left\vert f\left(  z+h\right)  -f\left(  z\right)  \right\vert dz\\
&  \leq\int_{\mathcal{O}_{\mathbf{1}}}\left\vert f\left(  z+h\right)
-f\left(  z\right)  \right\vert dz+\int_{0}^{h_{1}}\int_{0^{\prime}}%
^{\infty^{\prime}}\left\vert f\left(  z\right)  \right\vert dz_{1}dz^{\prime
}+\int_{0}^{h_{2}}\int_{0^{\prime}}^{\infty^{\prime}}\left\vert f\left(
z\right)  \right\vert dz_{2}dz^{\prime}+\ldots\\
&  =\ldots+\int_{0}^{h_{1}}\int_{0^{\prime}}^{\infty^{\prime}}b^{\left(
0\right)  }e^{-a^{\left(  0\right)  }z}dz_{1}dz^{\prime}+\int_{0}^{h_{2}}%
\int_{0^{\prime}}^{\infty^{\prime}}b^{\left(  0\right)  }e^{-a^{\left(
0\right)  }z}dz_{2}dz^{\prime}+\ldots\\
&  =\ldots+b^{\left(  0\right)  }\int_{0}^{h_{1}}\int_{0^{\prime}}%
^{\infty^{\prime}}e^{-a^{\left(  0\right)  }z}dz_{1}dz^{\prime}+b^{\left(
0\right)  }\int_{0}^{h_{2}}\int_{0^{\prime}}^{\infty^{\prime}}e^{-a^{\left(
0\right)  }z}dz_{2}dz^{\prime}+\ldots\\
&  =\ldots+b^{\left(  0\right)  }\left(  \int_{0}^{h_{1}}e^{-a_{1}^{\left(
0\right)  }z_{1}}dz_{1}\right)  \int_{0^{\prime}}^{\infty^{\prime}}e^{-\left(
a^{\left(  0\right)  }\right)  ^{\prime}z^{\prime}}dz^{\prime}+b^{\left(
0\right)  }\int_{0}^{h_{2}}\int_{0^{\prime}}^{\infty^{\prime}}e^{-a^{\left(
0\right)  }z}dz_{2}dz^{\prime}+\ldots\\
&  <\ldots+\frac{b^{\left(  0\right)  }}{a_{2}^{\left(  0\right)  }\ldots
a_{d}^{\left(  0\right)  }}h_{1}+\frac{b^{\left(  0\right)  }}{a_{1}^{\left(
0\right)  }a_{3}^{\left(  0\right)  }\ldots a_{d}^{\left(  0\right)  }}%
h_{2}+\ldots\\
&  =\ldots+b^{\left(  0\right)  }a_{1}^{\left(  0\right)  }h_{1}+b^{\left(
0\right)  }a_{2}^{\left(  0\right)  }h_{2}+\ldots\\
&  =\ldots+b^{\left(  0\right)  }\left(  a^{\left(  0\right)  }h\right)  .
\end{align*}

Choose the permutation $\pi$ so that $\pi h=\left(  h^{\prime},h^{\prime
\prime}\right)  $ where $h^{\prime}.>0$ and $h^{\prime\prime}.<0$.

Two cases: $h^{\prime\prime}$ empty and not empty.

\textbf{First integral} Split $\mathcal{O}_{\mathbf{1}}$ into unit cubes
$\left\{  \mathcal{C}_{\alpha}\right\}  _{\alpha\geq\mathbf{0}}$ and use a
Taylor series expansion and remainder estimate \ref{1.34}:%
\begin{align*}
\left\vert f\left(  z+h\right)  -f\left(  z\right)  \right\vert  & \leq
\sum_{i}h_{i}\int_{0}^{1}\left\vert \left(  D_{i}f\right)  \left(
z+th\right)  \right\vert dt\leq\sum_{i}h_{i}\int_{0}^{1}\left\vert \left(
b^{\left(  i\right)  }e^{-a^{\left(  i\right)  }z}\right)  \left(
z+th\right)  \right\vert dt=\\
& =\sum_{i}h_{i}b^{\left(  i\right)  }\int_{0}^{1}e^{-a^{\left(  i\right)
}\left(  z+th\right)  }dt\leq\sum_{i}h_{i}b^{\left(  i\right)  }e^{-a^{\left(
i\right)  }z}\int_{0}^{1}e^{-ta^{\left(  i\right)  }h}dt\leq\\
& \leq\sum_{i}h_{i}b^{\left(  i\right)  }e^{-a^{\left(  i\right)  }z}.
\end{align*}

\[
\int_{\mathcal{O}_{\mathbf{1}}}\left\vert f\left(  z+h\right)  -f\left(
z\right)  \right\vert dz\leq\sum_{i}h_{i}b^{\left(  i\right)  }\int%
_{\mathcal{O}_{\mathbf{1}}}e^{-a^{\left(  i\right)  }z}dz=\sum_{i}h_{i}%
\frac{b^{\left(  i\right)  }}{\left(  a^{\left(  i\right)  }\right)
^{\mathbf{1}}}=h\left(  \frac{b^{\left(  i\right)  }}{\left(  a^{\left(
i\right)  }\right)  ^{\mathbf{1}}}\right)  .
\]

$\operatorname*{supp}f\subseteq\mathcal{O}_{\mathbf{1}}$. Suppose
$h\in\mathcal{O}_{\theta}$. Perhaps $h=\theta.h^{\prime}$ where $h^{\prime}%
\in\mathcal{O}_{\mathbf{1}}$.%
\[
\int\left\vert f\left(  z+h\right)  -f\left(  z\right)  \right\vert
dz=\int\int\left\vert f\left(  z^{\prime}+h^{\prime}.\mathbf{1}^{\prime
},z^{\prime\prime}-h^{\prime\prime}.\mathbf{1}^{\prime\prime}\right)
-f\left(  z^{\prime},z^{\prime\prime}\right)  \right\vert dz^{\prime
}dz^{\prime\prime}.
\]

Suppose: $h^{\prime}.>\mathbf{0}$ and $h^{\prime\prime}.<\mathbf{0}$.

??

Case 1 $z^{\prime}+h^{\prime}.\mathbf{1}^{\prime}\in\mathcal{O}_{\mathbf{1}%
^{\prime}}$ and $z^{\prime\prime}-h^{\prime\prime}.\mathbf{1}^{\prime\prime
}\in\mathcal{O}_{\mathbf{1}^{\prime\prime}}$ and $z^{\prime}\in\mathcal{O}%
_{\mathbf{1}^{\prime}}$ and $z^{\prime\prime}\in\mathcal{O}_{\mathbf{1}%
^{\prime\prime}}$.

====================

Case 2 $z^{\prime}+h^{\prime}.\mathbf{1}^{\prime}\in\mathcal{O}_{\mathbf{1}%
^{\prime}}$ and $z^{\prime\prime}-h^{\prime\prime}.\mathbf{1}^{\prime\prime
}\in\mathcal{O}_{\mathbf{1}^{\prime\prime}}$ and $z^{\prime}\notin%
\mathcal{O}_{\mathbf{1}^{\prime}}$ and $z^{\prime\prime}\in\mathcal{O}%
_{\mathbf{1}^{\prime\prime}}$.

Case 3 $z^{\prime}+h^{\prime}.\mathbf{1}^{\prime}\in\mathcal{O}_{\mathbf{1}%
^{\prime}}$ and $z^{\prime\prime}-h^{\prime\prime}.\mathbf{1}^{\prime\prime
}\in\mathcal{O}_{\mathbf{1}^{\prime\prime}}$ and $z^{\prime}\in\mathcal{O}%
_{\mathbf{1}^{\prime}}$ and $z^{\prime\prime}\notin\mathcal{O}_{\mathbf{1}%
^{\prime\prime}}$.

====================

Case 4 $z^{\prime}+h^{\prime}.\mathbf{1}^{\prime}\in\mathcal{O}_{\mathbf{1}%
^{\prime}}$ and $z^{\prime\prime}-h^{\prime\prime}.\mathbf{1}^{\prime\prime
}\notin\mathcal{O}_{\mathbf{1}^{\prime\prime}}$ and $z^{\prime}\in
\mathcal{O}_{\mathbf{1}^{\prime}}$ and $z^{\prime\prime}\in\mathcal{O}%
_{\mathbf{1}^{\prime\prime}}$.

Case 5 $z^{\prime}+h^{\prime}.\mathbf{1}^{\prime}\notin\mathcal{O}%
_{\mathbf{1}^{\prime}}$ and $z^{\prime\prime}-h^{\prime\prime}.\mathbf{1}%
^{\prime\prime}\in\mathcal{O}_{\mathbf{1}^{\prime\prime}}$ and $z^{\prime}%
\in\mathcal{O}_{\mathbf{1}^{\prime}}$ and $z^{\prime\prime}\in\mathcal{O}%
_{\mathbf{1}^{\prime\prime}}$.

??? FINISH! ====================

Want $z\in\left(  \mathcal{O}_{\mathbf{1}}-\theta.h\right)  \setminus
\mathcal{O}_{\mathbf{1}}$, $z\in\mathcal{O}_{\mathbf{1}}\setminus\left(
\mathcal{O}_{\mathbf{1}}-\theta.h\right)  $, $z\in\left(  \mathcal{O}%
_{\mathbf{1}}-\theta.h\right)  \setminus\mathcal{O}_{\mathbf{1}}%
\cap\mathcal{O}_{\mathbf{1}}\setminus\left(  \mathcal{O}_{\mathbf{1}}%
-\theta.h\right)  =\left\{  {}\right\}  $.

The resultant sum converges since $f$ decreases exponentially at infinity.

\textbf{Second integral}
\[
\int_{\left(  \partial\mathcal{O}_{\mathbf{1}}\right)  _{h_{+}}}\left\vert
f\left(  z+h\right)  -f\left(  z\right)  \right\vert dz=\int_{\left\vert
z_{1}\right\vert \leq\left\vert h_{1}\right\vert }\int_{\mathcal{O}%
_{\mathbf{1}^{\prime}}-\left\vert h_{1}\right\vert }\left\vert f\left(
z+h\right)  -f\left(  z\right)  \right\vert dz_{1}dz^{\prime}+\ldots
\]

\end{proof}

\section{An extension from an orthant for the $C_{0}^{\left(  \alpha\right)
}\left(  \overline{\mathcal{O}_{\theta}}\right)  $ spaces.
\label{Sect_Exten_FromOrthant}}

In this section I first derive Proposition \ref{Prop_exten_orthant_O1} which
uses the function $J_{\gamma}[v_{\beta}]$ to construct continuous extension
$\mathcal{E}_{\alpha}^{\mathbf{1}}:C_{0}^{\left(  \alpha\right)  }\left(
\overline{\mathcal{O}_{\mathbf{1}}}\right)  \rightarrow C_{B}^{\left(
\alpha\right)  }\left(  \mathbb{R}^{d}\right)  $ for each $\alpha
\geq\mathbf{1}$. Here $v$ has the property described in part 1 of Lemma
\ref{Lem_PWCinf_J1_b}. In Corollary \ref{Cor_Prop_exten_any_orthant} we
construct a continuous extension $\mathcal{E}_{\alpha}^{\theta}:C_{0}^{\left(
\alpha\right)  }\left(  \overline{\mathcal{O}_{\theta}}\right)  \rightarrow
C_{B}^{\left(  \alpha\right)  }\left(  \mathbb{R}^{d}\right)  $.

The results of this section are not used elsewhere and are not used to
construct any of the other extensions in this Chapter.

\begin{proposition}
\label{Prop_exten_orthant_O1}\textbf{Continuous} $C_{0}^{\left(
\alpha\right)  }\left(  \overline{\mathcal{O}_{\mathbf{1}}}\right)  $
\textbf{extension operator }$\mathcal{E}_{\alpha}^{\mathbf{1}}$. Noting part 1
of Lemma \ref{Lem_SmthFuncIntegRepInOrthant}, for each $\alpha\geq\mathbf{1}$
we associate the operator $\mathcal{E}_{\alpha}^{\mathbf{1}}$ defined on
$C_{0}^{\left(  \alpha\right)  }\left(  \overline{\mathcal{O}_{\mathbf{1}}%
}\right)  $ by%
\begin{equation}
\mathcal{E}_{\alpha}^{\mathbf{1}}f\left(  x\right)  :=\tfrac{1}{\left(
\alpha-1\right)  !}\sum\limits_{\beta\leq\alpha}\tbinom{\alpha}{\beta}%
\int\limits_{\mathcal{O}_{\mathbf{1}}}J_{\alpha-2}\left[  v_{\beta}\right]
\left(  y-x\right)  D^{\beta}f\left(  y\right)  dy,\quad x\in\mathbb{R}%
^{d}.\label{X80}%
\end{equation}

Sometimes we will write $\mathcal{E}_{\alpha}^{\mathbf{1}}\left[  v\right]  $. Then:

\begin{enumerate}
\item For each $\alpha\geq\mathbf{1}$, $\mathcal{E}_{\alpha}^{\mathbf{1}}$ is
a continuous extension from $C_{0}^{\left(  \alpha\right)  }\left(
\overline{\mathcal{O}_{\mathbf{1}}}\right)  $ to $C_{B}^{\left(
\alpha-1\right)  }\left(  \mathbb{R}^{d}\right)  $ when these spaces are
endowed with their respective supremum norms.\smallskip

Now suppose $v$ is a function with the property described in part 1 of Lemma
\ref{Lem_PWCinf_J1_b}. Then:\medskip

\item
\[
D^{\gamma}\mathcal{E}_{\alpha}^{\mathbf{1}}f\left(  x\right)  =\tfrac{\left(
-1\right)  ^{\left\vert \gamma\right\vert }}{\left(  \alpha-1\right)  !}%
\sum\limits_{\beta\leq\alpha}\tbinom{\alpha}{\beta}\int\left(  D^{\gamma
}J_{\alpha-2}\left[  v_{\beta}\right]  \right)  \left(  y-x\right)  D^{\beta
}f\left(  y\right)  dy,\quad\gamma\leq\alpha.
\]

\item The operator $\mathcal{E}_{\alpha}^{\mathbf{1}}$ is a continuous
extension from $C_{0}^{\left(  \alpha\right)  }\left(  \overline
{\mathcal{O}_{\mathbf{1}}}\right)  $ to $C_{B}^{\left(  \alpha\right)
}\left(  \mathbb{R}^{d}\right)  $ when endowed with their respective supremum
norms. In fact
\[
\max_{\gamma\leq\alpha}\left\Vert D^{\gamma}\mathcal{E}_{\alpha}^{\mathbf{1}%
}f\right\Vert _{\infty}\leq k_{v,\alpha}\max_{\beta\leq\alpha}\left\Vert
D^{\beta}f\right\Vert _{\infty,\overline{\mathcal{O}_{\mathbf{1}}}},
\]

where $k_{v,\alpha}$ is given by \ref{X35} or equivalently \ref{X36}.
\end{enumerate}
\end{proposition}

\begin{proof}
\textbf{Part 1} Part 4 of Lemma \ref{Lem_property_J} implies $J_{\alpha
-2}\left[  v_{\beta}\right]  \in L^{1}$ and clearly $\left\{  D^{\beta
}f\right\}  \in L_{0}^{1}\cap L^{\infty}$. Hence Young's convolution theorem
\ref{a1.022} implies
\begin{equation}
\mathcal{E}_{\alpha}^{\mathbf{1}}f=\tfrac{\left(  2\pi\right)  ^{d/2}}{\left(
\alpha-1\right)  !}\sum\limits_{\beta\leq\alpha}\tbinom{\alpha}{\beta
}J_{\alpha-2}\left[  v_{\beta}\right]  _{-}\ast\left\{  D^{\beta}f\right\}
\in L^{1}\cap L^{\infty}.\label{X73}%
\end{equation}

From part 1 of Lemma \ref{Lem_SmthFuncIntegRepInOrthant}, if $x\in
\overline{\mathcal{O}_{\mathbf{1}}}$,%
\begin{align*}
\mathcal{E}_{\alpha}^{\mathbf{1}}f\left(  x\right)   & =\tfrac{1}{\left(
\alpha-1\right)  !}\sum\limits_{\beta\leq\alpha}\tbinom{\alpha}{\beta}%
\int_{\mathcal{O}_{\mathbf{1}}}J_{\alpha-2}\left[  v_{\beta}\right]  \left(
y-x\right)  D^{\beta}f\left(  y\right)  dy\\
& =\tfrac{1}{\left(  \alpha-1\right)  !}\sum\limits_{\beta\leq\alpha}%
\tbinom{\alpha}{\beta}\int_{x+\mathcal{O}_{\mathbf{1}}}J_{\alpha-2}\left[
v_{\beta}\right]  \left(  y-x\right)  D^{\beta}f\left(  y\right)  dy\\
& =f\left(  x\right)  .
\end{align*}

so $\mathcal{E}_{\alpha}^{\mathbf{1}}$ is an extension operator. Define%
\[
\left\{  D^{\beta}f\right\}  \left(  x\right)  :=\left\{
\begin{array}
[c]{ll}%
D^{\beta}f\left(  x\right)  , & x\in\overline{\mathcal{O}_{1}},\\
0, & x\notin\overline{\mathcal{O}_{1}}.
\end{array}
\right.
\]

and since $J_{\alpha-2}\left[  v_{\beta}\right]  \in S^{\prime}$ and $\left\{
D^{\beta}f\right\}  \in\mathcal{E}^{\prime}$ we can use the $S^{\prime}%
\ast\mathcal{E}^{\prime}$ convolution (see, for example, Subsection 2.9.6 of
Vladimirov \cite{Vladimirov}) to write%
\begin{equation}
D^{\gamma}\mathcal{E}_{\alpha}^{\mathbf{1}}f=\left(  -1\right)  ^{\left\vert
\gamma\right\vert }\tfrac{\left(  2\pi\right)  ^{d/2}}{\left(  \alpha
-1\right)  !}\sum\limits_{\beta\leq\alpha}\tbinom{\alpha}{\beta}\left(
D^{\gamma}J_{\alpha-2}\left[  v_{\beta}\right]  \right)  _{-}\ast\left\{
D^{\beta}f\right\}  ,\quad\forall\gamma.\label{X51}%
\end{equation}

Also $J_{\alpha-2}^{\mathbf{1}}\left[  v_{\beta}\right]  _{-}\in C^{\left(
0\right)  }\left(  \overline{\mathcal{O}_{-\mathbf{1}}}\right)  $,
$J_{\alpha-2}^{\mathbf{1}}\left[  v_{\beta}\right]  _{-}=0$ on $\mathbb{R}%
^{d}\setminus\overline{\mathcal{O}_{-\mathbf{1}}}$ and $J_{\alpha
-2}^{\mathbf{1}}\left[  v_{\beta}\right]  _{-}$ tends to zero exponentially at
infinity on $\overline{\mathcal{O}_{-\mathbf{1}}}$. Since $\left\{  D^{\beta
}f\right\}  $ can be translated so that it's bounded support lies in
$\overline{\mathcal{O}_{-\mathbf{1}}}$ we can apply Lemma
\ref{Lem_convol_PWsmth} and, since part 6 of Lemma \ref{Lem_property_J}
implies $D^{\gamma}J_{\alpha-2}\left[  v_{\beta}\right]  $ is $PWC^{\infty}$
when $\mathbf{0}\leq\gamma\leq\alpha-1$, we can conclude that $D^{\gamma
}\mathcal{E}_{\alpha}^{\mathbf{1}}f\in C_{B}^{\left(  0\right)  }$ when
$\mathbf{0}\leq\gamma\leq\alpha-1$ i.e. $\mathcal{E}_{\alpha}^{\mathbf{1}}f\in
C_{B}^{\left(  \alpha-1\right)  } $.

We now use part 6 of Definition \ref{Def_convol} concerning Young's inequality
for convolutions. Here $1/p+1/q=1+1/r$ with $r=\infty$, $p=1$, $q=\infty$, so
for $\gamma\leq\alpha-1$,%
\begin{equation}
D^{\gamma}\mathcal{E}_{\alpha}^{\mathbf{1}}f\left(  x\right)  =\tfrac{\left(
-1\right)  ^{\left\vert \gamma\right\vert }}{\left(  \alpha-1\right)  !}%
\sum\limits_{\beta\leq\alpha}\tbinom{\alpha}{\beta}\int\left(  D^{\gamma
}J_{\alpha-2}\left[  v_{\beta}\right]  \right)  \left(  y-x\right)  D^{\beta
}f\left(  y\right)  dy,\quad\gamma\leq\alpha-1,\label{X82}%
\end{equation}

and%
\begin{align}
\left\Vert D^{\gamma}\mathcal{E}_{\alpha}^{\mathbf{1}}f\right\Vert
_{\infty,\mathbb{R}^{d}}  & \leq\tfrac{1}{\left(  \alpha-1\right)  !}%
\sum\limits_{\beta\leq\alpha}\tbinom{\alpha}{\beta}\left\Vert \left(
D^{\gamma}J_{\alpha-2}\left[  v_{\beta}\right]  \right)  _{-}\right\Vert
_{1}\left\Vert \left\{  D^{\beta}f\right\}  \right\Vert _{\infty}\nonumber\\
& =\tfrac{1}{\left(  \alpha-1\right)  !}\sum\limits_{\beta\leq\alpha}%
\tbinom{\alpha}{\beta}\left\Vert D^{\gamma}J_{\alpha-2}\left[  v_{\beta
}\right]  \right\Vert _{1}\left\Vert D^{\beta}f\right\Vert _{\infty
,\mathcal{O}_{\mathbf{1}}}\nonumber\\
& \leq\tfrac{1}{\left(  \alpha-1\right)  !}\left(  \sum\limits_{\beta
\leq\alpha}\tbinom{\alpha}{\beta}\left\Vert D^{\gamma}J_{\alpha-2}\left[
v_{\beta}\right]  \right\Vert _{1}\right)  \max_{\beta\leq\alpha}\left\Vert
D^{\beta}f\right\Vert _{\infty,\mathcal{O}_{\mathbf{1}}},\label{X75}%
\end{align}

and the $L^{1}$ estimate of part 4 of Lemma \ref{Lem_property_J} implies that
$\left\Vert D^{\gamma}\mathcal{E}_{\alpha}^{\mathbf{1}}f\right\Vert
_{\infty,\mathbb{R}^{d}}<\infty$. Thus
\begin{equation}
\mathcal{E}_{\alpha}^{\mathbf{1}}f\in C_{B}^{\left(  \alpha-1\right)  },\quad
f\in C_{0}^{\left(  \alpha\right)  }\left(  \mathcal{O}_{\mathbf{1}}\right)
,\label{X74}%
\end{equation}

and%
\begin{equation}
\max_{\gamma\leq\alpha-1}\left\Vert D^{\gamma}\mathcal{E}_{\alpha}%
^{\mathbf{1}}f\right\Vert _{\infty,\mathbb{R}^{d}}\leq\tfrac{1}{\left(
\alpha-1\right)  !}\left(  \max_{\gamma\leq\alpha-1}\sum\limits_{\beta
\leq\alpha}\tbinom{\alpha}{\beta}\left\Vert D^{\gamma}J_{\alpha-2}\left[
v_{\beta}\right]  \right\Vert _{1}\right)  \max_{\beta\leq\alpha}\left\Vert
D^{\beta}f\right\Vert _{\infty,\mathcal{O}_{\mathbf{1}}},\label{X79}%
\end{equation}

which confirms that $\mathcal{E}_{\alpha}^{\mathbf{1}}:C_{0}^{\left(
\alpha\right)  }\left(  \mathcal{O}_{\mathbf{1}}\right)  \rightarrow
C_{B}^{\left(  \alpha-1\right)  }\left(  \mathbb{R}^{d}\right)  $ is
continuous.\medskip

\textbf{Parts 2 and 3} To have $\mathcal{E}_{\alpha}^{\mathbf{1}}f\in
C_{B}^{\left(  \alpha\right)  }$ we will assume that the function $\sigma$
used to define $J_{\alpha-2}$ is a tensor product. From \ref{X51},%
\[
D^{\gamma}\mathcal{E}_{\alpha}^{\mathbf{1}}f=\tfrac{\left(  -1\right)
^{\left\vert \alpha\right\vert }\left(  2\pi\right)  ^{d/2}}{\left(
\alpha-1\right)  !}\sum\limits_{\beta\leq\alpha}\tbinom{\alpha}{\beta
}D^{\gamma}\left(  J_{\alpha-2}\left[  v_{\beta}\right]  \right)  _{-}%
\ast\left\{  D^{\beta}f\right\}  ,\quad\forall\gamma,
\]

and from part 1 we know that $\mathcal{E}_{\alpha}^{\mathbf{1}}f\in
C_{B}^{\left(  \alpha-1\right)  }$. So we must prove that $D^{\gamma
}\mathcal{E}_{\alpha}^{\mathbf{1}}f\in C_{B}^{\left(  0\right)  }$ when
$\alpha-1<\gamma\leq\alpha$. Hence Lemma \ref{Lem_PWCinf_Jtheta_b} implies
$D^{\gamma}J_{\alpha-2}\left[  v_{\beta}\right]  \in PWC^{\infty}\left(
\mathbb{R}^{d}\right)  $ when $\beta,\gamma\leq\alpha$ and since $D^{\gamma
}J_{\alpha-2}\left[  v_{\beta}\right]  $ is exponentially decreasing at
infinity we have $D^{\gamma}J_{\alpha-2}\left[  v_{\beta}\right]  \in L^{1}$.
Young's convolution theorem (part 6 of Definition \ref{Def_convol}) and
\ref{X82} now imply that%
\[
D^{\gamma}\mathcal{E}_{\alpha}^{\mathbf{1}}f=\tfrac{\left(  -1\right)
^{\left\vert \gamma\right\vert }}{\left(  \alpha-1\right)  !}\sum
\limits_{\beta\leq\alpha}\tbinom{\alpha}{\beta}\int\left(  D^{\gamma}%
J_{\alpha-2}\left[  v_{\beta}\right]  \right)  \left(  y-x\right)  D^{\beta
}f\left(  y\right)  dy,\quad\gamma\leq\alpha,
\]

and%
\begin{align}
\left\Vert D^{\gamma}\mathcal{E}_{\alpha}^{\mathbf{1}}f\right\Vert _{\infty}
& \leq\tfrac{1}{\left(  \alpha-1\right)  !}\sum\limits_{\beta\leq\alpha
}\tbinom{\alpha}{\beta}\left\Vert D^{\gamma}\left(  J_{\alpha-2}\left[
v_{\beta}\right]  \right)  _{-}\right\Vert _{1}\left\Vert D^{\beta
}f\right\Vert _{\infty,\mathcal{O}_{\mathbf{1}}}\nonumber\\
& =\tfrac{1}{\left(  \alpha-1\right)  !}\sum\limits_{\beta\leq\alpha}%
\tbinom{\alpha}{\beta}\left\Vert D^{\gamma}J_{\alpha-2}\left[  v_{\beta
}\right]  \right\Vert _{1}\left\Vert D^{\beta}f\right\Vert _{\infty
,\mathcal{O}_{\mathbf{1}}}\nonumber\\
& \leq\tfrac{1}{\left(  \alpha-1\right)  !}\left(  \sum\limits_{\beta
\leq\alpha}\tbinom{\alpha}{\beta}\left\Vert D^{\gamma}J_{\alpha-2}\left[
v_{\beta}\right]  \right\Vert _{1}\right)  \max_{\beta\leq\alpha}\left\Vert
D^{\beta}f\right\Vert _{\infty,\mathcal{O}_{\mathbf{1}}}\label{X76}\\
& <\infty.\nonumber
\end{align}

Thus $D^{\gamma}\mathcal{E}_{\alpha}^{\mathbf{1}}f\in C_{B}^{\left(  0\right)
}$ and we can conclude that $\mathcal{E}_{\alpha}^{\mathbf{1}}f\in
C_{B}^{\left(  \alpha\right)  }$. Finally, \ref{X76} yields%
\begin{align}
\max_{\gamma\leq\alpha}\left\Vert D^{\gamma}\mathcal{E}_{\alpha}^{\mathbf{1}%
}f\right\Vert _{\infty}  & \leq\tfrac{1}{\left(  \alpha-1\right)  !}%
\max_{\gamma\leq\alpha}\left(  \sum\limits_{\beta\leq\alpha}\tbinom{\alpha
}{\beta}\left\Vert D^{\gamma}J_{\alpha-2}\left[  v_{\beta}\right]  \right\Vert
_{1}\right)  \max_{\beta\leq\alpha}\left\Vert D^{\beta}f\right\Vert
_{\infty,\mathcal{O}_{\mathbf{1}}}\nonumber\\
& =\tfrac{1}{\left(  \alpha-1\right)  !}\max_{\gamma\leq\alpha}\left(
\sum\limits_{\beta\leq\alpha}\tbinom{\alpha}{\beta}\left\Vert D^{\gamma
}J_{\alpha-2}\left[  v_{\beta}\right]  \right\Vert _{1,\mathcal{O}_{1}%
}\right)  \max_{\beta\leq\alpha}\left\Vert D^{\beta}f\right\Vert
_{\infty,\mathcal{O}_{\mathbf{1}}}.\label{X78}%
\end{align}

Thus $\mathcal{E}_{\alpha}^{\mathbf{1}}$ is a continuous extension from
$C_{0}^{\left(  \alpha\right)  }\left(  \overline{\mathcal{O}_{\mathbf{1}}%
}\right)  $ to $C_{B}^{\left(  \alpha\right)  }\left(  \mathbb{R}^{d}\right)
$ under the supremum norm.

Now, from ??,%
\[
\left\Vert D^{\gamma}J_{\alpha-2}\left[  v_{\beta}\right]  \right\Vert
_{1,\mathcal{O}_{1}}\leq\left\Vert \frac{v_{\beta}\left(  x\right)
}{x^{\gamma}}\right\Vert _{1}\int_{\mathcal{O}_{\mathbf{1}}}q_{\gamma}%
^{\alpha-2}\left(  y\right)  e^{-y\mathbf{1}}dy=\left\Vert \frac{x^{\beta
}v\left(  x\right)  }{x^{\gamma}}\right\Vert _{1}\int_{\mathcal{O}%
_{\mathbf{1}}}q_{\gamma}^{\alpha-2}\left(  y\right)  e^{-y\mathbf{1}}dy,
\]

so that%
\begin{align*}
\sum\limits_{\beta\leq\alpha}\tbinom{\alpha}{\beta}\left\Vert D^{\gamma
}J_{\alpha-2}\left[  v_{\beta}\right]  \right\Vert _{1,\mathcal{O}_{1}}  &
\leq\sum\limits_{\beta\leq\alpha}\tbinom{\alpha}{\beta}\left\Vert
\frac{x^{\beta}v\left(  x\right)  }{x^{\gamma}}\right\Vert _{1}\int%
_{\mathcal{O}_{\mathbf{1}}}q_{\gamma}^{\alpha-2}\left(  y\right)
e^{-y\mathbf{1}}dy\\
& =\sum\limits_{\beta\leq\alpha}\tbinom{\alpha}{\beta}\int\frac{x^{\beta
}\left\vert v\left(  x\right)  \right\vert }{\left\vert x^{\gamma}\right\vert
}dx\int_{\mathcal{O}_{\mathbf{1}}}q_{\gamma}^{\alpha-2}\left(  y\right)
e^{-y\mathbf{1}}dy\\
& =\int_{a}^{b}\frac{\left(  1+x\right)  ^{\alpha}\left\vert v\left(
x\right)  \right\vert }{x^{\gamma}}dx\int_{\mathcal{O}_{\mathbf{1}}}q_{\gamma
}^{\alpha-2}\left(  y\right)  e^{-y\mathbf{1}}dy.
\end{align*}

From part 4 of Lemma \ref{Lem_property_J},%
\[
\int_{0}^{\infty}q_{\gamma_{i}}^{\delta_{i}}\left(  t\right)  e^{-t}dt=\left(
\delta_{i}+1\right)  !\times\left\{
\begin{array}
[c]{ll}%
2^{\gamma_{i}}, & \gamma_{i}\leq\delta_{i}+1,\\
\sum\limits_{k=0}^{\delta_{i}+1}\tbinom{\gamma_{i}}{k}, & \gamma_{i}\geq
\delta_{i}+2,
\end{array}
\right.
\]

so that%
\[
\int_{0}^{\infty}q_{\gamma_{i}}^{\alpha_{i}-2}\left(  t\right)  e^{-t}%
dt=\left(  \alpha_{i}-1\right)  !\times\left\{
\begin{array}
[c]{ll}%
2^{\gamma_{i}}, & \gamma_{i}\leq\alpha_{i}-1,\\
\sum\limits_{k=0}^{\alpha_{i}-1}\tbinom{\gamma_{i}}{k}, & \gamma_{i}\geq
\alpha_{i},
\end{array}
\right.
\]

and hence%
\[
\int_{0}^{\infty}q_{\gamma_{i}}^{\alpha_{i}-2}\left(  t\right)  e^{-t}%
dt=\left(  \alpha_{i}-1\right)  !\times\left\{
\begin{array}
[c]{ll}%
2^{\gamma_{i}}, & \gamma_{i}\leq\alpha_{i}-1,\\
2^{\gamma_{i}}-1, & \gamma_{i}=\alpha_{i},
\end{array}
\right.
\]

which implies%
\begin{equation}
\int_{0}^{\infty}q_{\gamma_{i}}^{\alpha_{i}-2}\left(  t\right)  e^{-t}%
dt\leq\left(  \alpha_{i}-1\right)  !2^{\gamma_{i}},\quad\gamma_{i}\leq
\alpha_{i}.\label{X33}%
\end{equation}

The estimate \ref{X78} can now be written%
\[
\max_{\gamma\leq\alpha}\left\Vert D^{\gamma}\mathcal{E}_{\alpha}^{\mathbf{1}%
}f\right\Vert _{\infty}\leq k_{v,\alpha}\max_{\beta\leq\alpha}\left\Vert
D^{\beta}f\right\Vert _{\infty,\overline{\mathcal{O}_{\mathbf{1}}}},
\]

where%
\begin{align}
k_{v,\alpha}  & =\max_{\gamma\leq\alpha}\left(  \int_{a}^{b}\frac{\left(
1+x\right)  ^{\alpha}}{\left(  x/2\right)  ^{\gamma}}\left\vert v\left(
x\right)  \right\vert dx\right) \label{X35}\\
& =\max_{\gamma\leq\alpha}\left(  \prod\limits_{i=1}^{d}\int_{a_{i}}^{b_{i}%
}\frac{\left(  1+x_{i}\right)  ^{\alpha_{i}}}{\left(  x_{i}/2\right)
^{\gamma_{i}}}\left\vert v_{i}\left(  x_{i}\right)  \right\vert dx_{i}\right)
\nonumber\\
& =\prod\limits_{i=1}^{d}\left(  \max_{\gamma_{i}\leq\alpha_{i}}\int_{a_{i}%
}^{b_{i}}\frac{\left(  1+x_{i}\right)  ^{\alpha_{i}}}{\left(  x_{i}/2\right)
^{\gamma_{i}}}\left\vert v_{i}\left(  x_{i}\right)  \right\vert dx_{i}\right)
\nonumber\\
& =\prod\limits_{i=1}^{d}\left(  \max_{k=0}^{\alpha_{i}}\int_{a_{i}}^{b_{i}%
}\frac{\left(  1+t\right)  ^{\alpha_{i}}}{\left(  t/2\right)  ^{k}}\left\vert
v_{i}\left(  t\right)  \right\vert dt\right)  .\label{X36}%
\end{align}

\end{proof}

\begin{remark}
Suppose $\alpha\geq\mathbf{1}$. Then from part 4 of Lemma
\ref{Cor_Prop_exten_any_orthant},%
\[
D^{\gamma}\left(  \mathcal{E}_{\alpha}^{\mathbf{1}}\left[  v\right]  f\right)
\left(  x\right)  =\tfrac{\left(  -1\right)  ^{\left\vert \gamma\right\vert }%
}{\left(  \alpha-1\right)  !}\sum\limits_{\beta\leq\alpha}\tbinom{\alpha
}{\beta}\int_{\mathcal{O}_{\mathbf{1}}}\left(  D^{\gamma}J_{\alpha
-2}^{\mathbf{1}}\left[  v_{\beta}\right]  \right)  \left(  y-x\right)
D^{\beta}f\left(  y\right)  dy,
\]

and from part 5 of Lemma \ref{Lem_property_J},%
\begin{align*}
D^{\gamma}J_{\alpha-2}^{\mathbf{1}}\left[  v_{\beta}\right]  \left(  z\right)
& =\frac{1}{z^{\gamma}}\sum_{\beta\leq\alpha-1}\left(  p_{\gamma}^{\alpha
-2}\right)  _{\beta}J_{\beta+\gamma-1}^{\mathbf{1}}\left[  v_{\beta}\right]
\left(  z\right) \\
& =\sum_{\beta\leq\alpha-1}\left(  p_{\gamma}^{\alpha-2}\right)  _{\beta
}z^{\beta}\frac{J_{\beta+\gamma-1}^{\mathbf{1}}\left[  v_{\beta}\right]
\left(  z\right)  }{z^{\beta+\gamma}},\quad\gamma\leq\alpha-1,
\end{align*}

where $\frac{J_{\beta+\gamma-1}^{\mathbf{1}}\left[  v_{\beta}\right]  \left(
z\right)  }{z^{\beta+\gamma}}\in C^{\infty}\left(  \overline{\mathcal{O}%
_{\mathbf{1}}}\right)  $.
\end{remark}

Using the equations and definitions of part 1 of Lemma
\ref{Lem_SmthFuncIntegRepInOrthant} we now generalize Proposition
\ref{Prop_exten_orthant_O1} by defining an extension from an arbitrary orthant.

\begin{definition}
\label{Def_exten_arbitrary_orthant}By definition, $f\in C_{0}^{\left(
\alpha\right)  }\left(  \overline{\mathcal{O}_{\theta}}\right)  $ iff
$f\left(  \theta.\right)  \in C_{0}^{\left(  \alpha\right)  }\left(
\overline{\mathcal{O}_{\mathbf{1}}}\right)  $ and $v\in L_{0}^{1}\left(
\mathcal{O}_{\theta}\right)  $ iff $v\left(  \theta.\right)  \in L_{0}%
^{1}\left(  \mathcal{O}_{\mathbf{1}}\right)  $ so noting the formula \ref{X80}
for $\mathcal{E}_{\alpha}^{\mathbf{1}}=\mathcal{E}_{\alpha}$ we define
$\mathcal{E}_{\alpha}^{\theta}=\mathcal{E}_{\alpha}^{\theta}\left[  v\right]
$ by%
\begin{equation}
\left(  \mathcal{E}_{\alpha}^{\theta}\left[  v\right]  f\right)  \left(
x\right)  :=\left(  \mathcal{E}_{\alpha}^{\mathbf{1}}\left[  v\left(
\theta.\right)  \right]  f\left(  \theta.\right)  \right)  \left(
\theta.x\right)  ,\label{X87}%
\end{equation}

and noting \ref{X73} it can be described by the following commutative diagram%
\begin{equation}
\left.
\begin{array}
[c]{ccc}%
C_{0}^{\left(  \alpha\right)  }\left(  \overline{\mathcal{O}_{\mathbf{1}}%
}\right)  & \overset{\mathcal{E}_{\alpha}^{\mathbf{1}}\left[  v\left(
\theta.\right)  \right]  }{\longrightarrow\longrightarrow\longrightarrow} &
L^{1}\cap L^{\infty}\\
\uparrow &  & \downarrow\\
\quad\uparrow\theta. &  & \quad\downarrow\theta.\\
\uparrow &  & \downarrow\\
C_{0}^{\left(  \alpha\right)  }\left(  \overline{\mathcal{O}_{\theta}}\right)
& \longrightarrow\underset{\mathcal{E}_{\alpha}^{\theta}\left[  v\right]
}{\longrightarrow}\longrightarrow & L^{1}\cap L^{\infty}%
\end{array}
\right\} \label{X86}%
\end{equation}

because the operators $\theta.:C_{0}^{\left(  \alpha\right)  }\left(
\overline{\mathcal{O}_{\theta}}\right)  \rightarrow C_{0}^{\left(
\alpha\right)  }\left(  \overline{\mathcal{O}_{\mathbf{1}}}\right)  $ and
$\theta.:L^{1}\rightarrow L^{1}\cap L^{\infty}$ are both isomorphisms.

Further, noting part 4 of topology Definition \ref{Def_topol_on_Rd}, if
$\operatorname*{supp}v\subset R\left[  c,d\right]  \subset\mathcal{O}_{\theta
}$ and $c.<d$ then
\begin{equation}
\left.
\begin{array}
[c]{ll}%
\operatorname*{supp}v\left(  \theta.\right)  & =\theta.\operatorname*{supp}%
v\subset\theta.R\left[  c,d\right]  =R\left[  \theta.c,\theta.d\right]  =\\
& =R\left[  \min\left(  \theta.c,\theta.d\right)  ,\max\left(  \theta
.c,\theta.d\right)  \right]  \subset\mathcal{O}_{\mathbf{1}}.
\end{array}
\right\} \label{X85}%
\end{equation}

\end{definition}

We now show that $\mathcal{E}_{\alpha}^{\theta}$ is a continuous extension
operator from $C_{0}^{\left(  \alpha\right)  }\left(  \overline{\mathcal{O}%
_{\theta}}\right)  $ to $C_{B}^{\left(  \alpha\right)  }\left(  \mathbb{R}%
^{d}\right)  $.

\begin{corollary}
\label{Cor_Prop_exten_any_orthant}\textbf{Properties of }$\mathcal{E}_{\alpha
}^{\theta}$. Noting part 1 of Lemma \ref{Lem_SmthFuncIntegRepInOrthant} we have:

\begin{enumerate}
\item For $f\in C_{0}^{\left(  \alpha\right)  }\left(  \overline
{\mathcal{O}_{\theta}}\right)  $ and $\alpha\geq-\mathbf{1}$:%
\[
\mathcal{E}_{\alpha}^{\theta}f\left(  x\right)  :=\tfrac{1}{\left(
\alpha-1\right)  !}\sum\limits_{\beta\leq\alpha}\tbinom{\alpha}{\beta}%
\int_{\mathcal{O}_{\theta}}J_{\alpha-2}^{\theta}\left[  v_{\beta}\right]
\left(  y-x\right)  D^{\beta}f\left(  y\right)  dy,\quad x\in\mathbb{R}^{d}.
\]

\item For each $\alpha\geq\mathbf{1}$, $\mathcal{E}_{\alpha}^{\theta}$ is a
continuous extension from $C_{0}^{\left(  \alpha\right)  }\left(
\overline{\mathcal{O}_{\theta}}\right)  $ to $C_{B}^{\left(  \alpha-1\right)
}\left(  \mathbb{R}^{d}\right)  $ when these spaces are endowed with their
respective supremum norms. We have the continuous commutative diagram:%
\[%
\begin{array}
[c]{ccc}%
C_{0}^{\left(  \alpha\right)  }\left(  \overline{\mathcal{O}_{\mathbf{1}}%
}\right)  & \overset{\mathcal{E}_{\alpha}^{\mathbf{1}}\left[  v\left(
\theta.\right)  \right]  }{\longrightarrow\longrightarrow\longrightarrow} &
C_{B}^{\left(  \alpha-1\right)  }\left(  \mathbb{R}^{d}\right) \\
\uparrow &  & \downarrow\\
\quad\uparrow\theta. &  & \quad\downarrow\theta.\\
\uparrow &  & \downarrow\\
C_{0}^{\left(  \alpha\right)  }\left(  \overline{\mathcal{O}_{\theta}}\right)
& \longrightarrow\underset{\mathcal{E}_{\alpha}^{\theta}\left[  v\right]
}{\longrightarrow}\longrightarrow & C_{B}^{\left(  \alpha-1\right)  }\left(
\mathbb{R}^{d}\right)
\end{array}
\]

\item ?? FIX!

\item For $\gamma\leq\alpha$,%
\[
D^{\gamma}\left(  \mathcal{E}_{\alpha}^{\theta}\left[  v\right]  f\right)
\left(  x\right)  =\tfrac{\left(  -1\right)  ^{\left\vert \gamma\right\vert }%
}{\left(  \alpha-1\right)  !}\sum\limits_{\beta\leq\alpha}\tbinom{\alpha
}{\beta}\int_{\mathcal{O}_{\theta}}\left(  D^{\gamma}J_{\alpha-2}^{\theta
}\left[  v_{\beta}\right]  \right)  \left(  y-x\right)  D^{\beta}f\left(
y\right)  dy.
\]

\item If $v$ is a function with the property described in part 3 then
$\mathcal{E}_{\alpha}^{\theta}$ is a continuous extension from $C_{0}^{\left(
\alpha\right)  }\left(  \overline{\mathcal{O}_{\theta}}\right)  $ to
$C_{B}^{\left(  \alpha\right)  }\left(  \mathbb{R}^{d}\right)  $ under the
supremum norm. In fact, if $\operatorname*{supp}v\subset R\left[  c,d\right]
\subset\mathcal{O}_{\theta}$ with $c.<d$ then
\begin{align*}
\max_{\gamma\leq\alpha} &  \left\Vert D^{\gamma}\mathcal{E}_{\alpha}^{\theta
}f\right\Vert _{\infty}\\
&  \leq\max_{\gamma\leq\alpha}\left(  \frac{2}{\min\left(  \theta
.c,\theta.d\right)  }\right)  ^{\gamma}\left(  \int_{c}^{d}\left(
1+\theta.z\right)  ^{\alpha}\left\vert v\left(  z\right)  \right\vert
dz\right)  \max_{\beta\leq\alpha}\left\Vert D^{\beta}f\right\Vert
_{\infty,\overline{\mathcal{O}_{\theta}}},
\end{align*}

and we have the continuous commutative diagram:%
\[%
\begin{array}
[c]{ccc}%
C_{0}^{\left(  \alpha\right)  }\left(  \overline{\mathcal{O}_{\mathbf{1}}%
}\right)  & \overset{\mathcal{E}_{\alpha}^{\mathbf{1}}\left[  v\left(
\theta.\right)  \right]  }{\longrightarrow\longrightarrow\longrightarrow} &
C_{B}^{\left(  \alpha\right)  }\left(  \mathbb{R}^{d}\right) \\
\uparrow &  & \downarrow\\
\quad\uparrow\theta. &  & \quad\downarrow\theta.\\
\uparrow &  & \downarrow\\
C_{0}^{\left(  \alpha\right)  }\left(  \overline{\mathcal{O}_{\theta}}\right)
& \longrightarrow\underset{\mathcal{E}_{\alpha}^{\theta}\left[  v\right]
}{\longrightarrow}\longrightarrow & C_{B}^{\left(  \alpha\right)  }\left(
\mathbb{R}^{d}\right)
\end{array}
\]

\item If $v$ is a function with the property described in part 3 then
$D^{\gamma}\left(  \mathcal{E}_{\alpha}^{\theta}\left[  v\right]  f\right)
\in L^{1}$ $\forall\gamma\leq\alpha$. Further%
\[
\left\Vert D^{\gamma}\left(  \mathcal{E}_{\alpha}^{\theta}\left[  v\right]
f\right)  \right\Vert _{1}=\left\Vert D^{\gamma}\left(  \mathcal{E}_{\alpha
}^{\theta}\left[  v\right]  f\right)  \right\Vert _{1,\mathcal{O}_{\theta}},
\]

and the upper bound given in part 5 of Lemma
\ref{Lem_SmthFuncIntegRepInOrthant} applies.\medskip

Now we consider the behavior of $\mathcal{E}_{\alpha}^{\theta}\left[
v\right]  f$ at infinity:\medskip

\item At infinity $\mathcal{E}_{\alpha}^{\theta}f$ is a rapidly decreasing
function and for all $\mu>0$,%
\[
\left(  1+\left\vert x\right\vert \right)  ^{\mu}\left\vert \mathcal{E}%
_{\alpha}^{\theta}f\left(  x\right)  \right\vert \leq\tfrac{1}{\left(
\alpha-1\right)  !}\sum\limits_{\beta\leq\alpha}\tbinom{\alpha}{\beta
}\left\Vert \left(  1+\left\vert \cdot\right\vert \right)  ^{\mu}J_{\alpha
-2}^{\theta}\left[  v_{\beta}\right]  \right\Vert _{1,\overline{\mathcal{O}%
_{\theta}}}\left\Vert \left(  1+\left\vert \cdot\right\vert \right)  ^{\mu
}D^{\beta}f\right\Vert _{\infty,\overline{\mathcal{O}_{\theta}}}.
\]
\medskip

In fact, $\mathcal{E}_{\alpha}^{\theta}f$ decreases exponentially at
infinity:\medskip

\item Suppose $\operatorname*{supp}v\subset R\left[  c,d\right]
\subset\mathcal{O}_{\theta}$ with $c.<d$. Then for any real $0<s<1$ there
exists real $C_{s,f,\alpha}>0$ given by \ref{X92} such that%
\[
\left\vert \mathcal{E}_{\alpha}^{\theta}f\left(  x\right)  \right\vert \leq
C_{s,f,\alpha}e^{-\frac{s}{\max\left\{  \theta.c,\theta.d\right\}  }x_{+}%
},\quad x\in\mathbb{R}^{d},\text{ }f\in C_{0}^{\left(  \alpha\right)  }\left(
\overline{\mathcal{O}_{\theta}}\right)  .
\]

Here $x_{+}=\left(  \left\vert x_{i}\right\vert \right)  $ and $C_{s,f,\alpha}
$ is independent of $x\in\mathbb{R}^{d}$.
\end{enumerate}
\end{corollary}

\begin{proof}
\textbf{Part 1} Using formula \ref{X80} for the action of $\mathcal{E}%
_{\alpha}^{\mathbf{1}}$, i.e.%
\[
\mathcal{E}_{\alpha}^{\mathbf{1}}f\left(  x\right)  :=\tfrac{1}{\left(
\alpha-1\right)  !}\sum\limits_{\beta\leq\alpha}\tbinom{\alpha}{\beta}%
\int_{\mathcal{O}_{\mathbf{1}}}J_{\alpha-2}\left[  v_{\beta}\right]  \left(
y-x\right)  D^{\beta}f\left(  y\right)  dy,\quad x\in\mathbb{R}^{d},
\]

we calculate%
\begin{align*}
\left(  \mathcal{E}_{\alpha}^{\theta}\left[  v\right]  f\right)  \left(
x\right)   & =\left(  \mathcal{E}_{\alpha}^{\mathbf{1}}\left[  v\left(
\theta.\right)  \right]  f\left(  \theta.\right)  \right)  \left(
\theta.x\right) \\
& =\tfrac{1}{\left(  \alpha-1\right)  !}\sum\limits_{\beta\leq\alpha}%
\tbinom{\alpha}{\beta}\int\limits_{\mathcal{O}_{\mathbf{1}}}J_{\alpha
-2}\left[  \left(  v\left(  \theta.\right)  \right)  _{\beta}\right]  \left(
y-\theta.x\right)  D^{\beta}\left(  f\left(  \theta.\right)  \right)  \left(
y\right)  dy\\
& =\tfrac{1}{\left(  \alpha-1\right)  !}\sum\limits_{\beta\leq\alpha}%
\tbinom{\alpha}{\beta}\int\limits_{\mathcal{O}_{\mathbf{1}}}J_{\alpha
-2}\left[  \left(  v\left(  \theta.\right)  \right)  _{\beta}\right]  \left(
y-\theta.x\right)  \theta^{\beta}\left(  D^{\beta}f\right)  \left(
\theta.y\right)  dy\\
& =\tfrac{1}{\left(  \alpha-1\right)  !}\sum\limits_{\beta\leq\alpha}%
\tbinom{\alpha}{\beta}\int\limits_{\mathcal{O}_{\theta}}\theta^{\beta
}J_{\alpha-2}\left[  \left(  v\left(  \theta.\right)  \right)  _{\beta
}\right]  \left(  \theta.y-\theta.x\right)  D^{\beta}f\left(  y\right)  dy\\
& =\tfrac{1}{\left(  \alpha-1\right)  !}\sum\limits_{\beta\leq\alpha}%
\tbinom{\alpha}{\beta}\int\limits_{\mathcal{O}_{\theta}}\theta^{\beta
}J_{\alpha-2}\left[  \left(  v\left(  \theta.\right)  \right)  _{\beta
}\right]  \left(  \theta.y-\theta.x\right)  D^{\beta}f\left(  y\right)  dy,
\end{align*}

but from \ref{X05},%
\[
J_{\lambda}^{\theta}\left[  v_{\beta}\right]  \left(  z\right)  :=\theta
^{\beta}J_{\lambda}\left[  \left(  v\left(  \theta.\right)  \right)  _{\beta
}\right]  \left(  \theta.z\right)  ,\quad\lambda\geq-\mathbf{1},
\]

so that%
\[
\left(  \mathcal{E}_{\alpha}^{\theta}\left[  v\right]  f\right)  \left(
x\right)  =\tfrac{1}{\left(  \alpha-1\right)  !}\sum\limits_{\beta\leq\alpha
}\tbinom{\alpha}{\beta}\int\limits_{\mathcal{O}_{\theta}}J_{\alpha-2}^{\theta
}\left[  v_{\beta}\right]  \left(  y-x\right)  D^{\beta}f\left(  y\right)  dy,
\]

as claimed.\medskip

\textbf{Part 2} The operators $\theta.:C_{0}^{\left(  \alpha-1\right)
}\left(  \overline{\mathcal{O}_{\theta}}\right)  \rightarrow C_{0}^{\left(
\alpha-1\right)  }\left(  \overline{\mathcal{O}_{\mathbf{1}}}\right)  $ and
$\theta.:C_{B}^{\left(  \alpha-1\right)  }\left(  \mathbb{R}^{d}\right)
\rightarrow C_{B}^{\left(  \alpha-1\right)  }\left(  \mathbb{R}^{d}\right)  $
are both isometric isomorphisms so it is clear from Proposition
\ref{Prop_exten_orthant_O1} that $\mathcal{E}_{\alpha}^{\theta}\left[
v\right]  :C_{0}^{\left(  \alpha\right)  }\left(  \overline{\mathcal{O}%
_{\theta}}\right)  \rightarrow C_{B}^{\left(  \alpha-1\right)  }\left(
\mathbb{R}^{d}\right)  $ is continuous.\medskip

\textbf{Part ??3} \medskip

\textbf{Part ??4 } From the definition of $\mathcal{E}_{\alpha}^{\theta}$,%
\begin{equation}
D^{\gamma}\left(  \mathcal{E}_{\alpha}^{\theta}\left[  v\right]  f\right)
\left(  x\right)  =\theta^{\gamma}D^{\gamma}\left\{  \mathcal{E}_{\alpha
}^{\mathbf{1}}\left[  v\left(  \theta.\right)  \right]  f\left(
\theta.\right)  \right\}  \left(  \theta.x\right)  ,\label{X84}%
\end{equation}

and from part 3 of Proposition \ref{Prop_exten_orthant_O1},%
\[
D^{\gamma}\mathcal{E}_{\alpha}^{\mathbf{1}}f\left(  x\right)  =\tfrac{\left(
-1\right)  ^{\left\vert \gamma\right\vert }}{\left(  \alpha-1\right)  !}%
\sum\limits_{\beta\leq\alpha}\tbinom{\alpha}{\beta}\int\left(  D^{\gamma
}J_{\alpha-2}\left[  v_{\beta}\right]  \right)  \left(  y-x\right)  D^{\beta
}f\left(  y\right)  dy,\quad\gamma\leq\alpha.
\]

Thus%
\begin{align*}
D^{\gamma} &  \left(  \mathcal{E}_{\alpha}^{\theta}f\right)  \left(  x\right)
\\
&  =\theta^{\gamma}\left(  \tfrac{\left(  -1\right)  ^{\left\vert
\gamma\right\vert }}{\left(  \alpha-1\right)  !}\sum\limits_{\beta\leq\alpha
}\tbinom{\alpha}{\beta}\int\left(  D^{\gamma}J_{\alpha-2}\left[  \left(
v\left(  \theta.\right)  \right)  _{\beta}\right]  \right)  \left(
y-x\right)  D^{\beta}\left(  f\left(  \theta.\right)  \right)  \left(
y\right)  dy\right)  \left(  \theta.x\right) \\
&  =\left(  \tfrac{\left(  -1\right)  ^{\left\vert \gamma\right\vert }%
}{\left(  \alpha-1\right)  !}\sum\limits_{\beta\leq\alpha}\tbinom{\alpha
}{\beta}\int\theta^{\gamma+\beta}\left(  D^{\gamma}J_{\alpha-2}\left[  \left(
v\left(  \theta.\right)  \right)  _{\beta}\right]  \right)  \left(
y-x\right)  \left(  D^{\beta}f\right)  \left(  \theta.y\right)  dy\right)
\left(  \theta.x\right) \\
&  =\tfrac{\left(  -1\right)  ^{\left\vert \gamma\right\vert }}{\left(
\alpha-1\right)  !}\sum\limits_{\beta\leq\alpha}\tbinom{\alpha}{\beta}%
\int\theta^{\gamma+\beta}\left(  D^{\gamma}J_{\alpha-2}\left[  \left(
v\left(  \theta.\right)  \right)  _{\beta}\right]  \right)  \left(
y-\theta.x\right)  \left(  D^{\beta}f\right)  \left(  \theta.y\right)  dy\\
&  =\tfrac{\left(  -1\right)  ^{\left\vert \gamma\right\vert }}{\left(
\alpha-1\right)  !}\sum\limits_{\beta\leq\alpha}\tbinom{\alpha}{\beta}%
\int\theta^{\gamma+\beta}\left(  D^{\gamma}J_{\alpha-2}\left[  \left(
v\left(  \theta.\right)  \right)  _{\beta}\right]  \right)  \left(
\theta.z-\theta.x\right)  D^{\beta}f\left(  z\right)  dz\\
&  =\tfrac{\left(  -1\right)  ^{\left\vert \gamma\right\vert }}{\left(
\alpha-1\right)  !}\sum\limits_{\beta\leq\alpha}\tbinom{\alpha}{\beta}%
\int\theta^{\gamma+\beta}\left(  D^{\gamma}J_{\alpha-2}\left[  \left(
v\left(  \theta.\right)  \right)  _{\beta}\right]  \right)  \left(
\theta.\left(  z-x\right)  \right)  D^{\beta}f\left(  z\right)  dz,
\end{align*}

but from \ref{X05},%
\[
J_{\alpha-2}^{\theta}\left[  v_{\beta}\right]  \left(  z\right)
:=\theta^{\beta}J_{\alpha-2}\left[  \left(  v\left(  \theta.\right)  \right)
_{\beta}\right]  \left(  \theta.z\right)  ,\quad\lambda\geq-\mathbf{1},
\]

so%
\[
D^{\gamma}J_{\alpha-2}^{\theta}\left[  v_{\beta}\right]  \left(  z\right)
:=\theta^{\beta+\gamma}D^{\gamma}\left(  J_{\alpha-2}\left[  \left(  v\left(
\theta.\right)  \right)  _{\beta}\right]  \right)  \left(  \theta.z\right)  ,
\]

and hence%
\[
D^{\gamma}\left(  \mathcal{E}_{\alpha}^{\theta}\left[  v\right]  f\right)
\left(  x\right)  =\tfrac{\left(  -1\right)  ^{\left\vert \gamma\right\vert }%
}{\left(  \alpha-1\right)  !}\sum\limits_{\beta\leq\alpha}\tbinom{\alpha
}{\beta}\int\left(  D^{\gamma}J_{\alpha-2}^{\theta}\left[  v_{\beta}\right]
\right)  \left(  y-x\right)  D^{\beta}f\left(  y\right)  dy.
\]
\medskip

\textbf{Part ??5} We begin with the inequality of part 4 of Proposition
\ref{Prop_exten_orthant_O1} i.e. if $\operatorname*{supp}v\subseteq R\left[
a,b\right]  \subset\mathcal{O}_{\mathbf{1}}$ where $a.<b$ then%
\[
\max_{\gamma\leq\alpha}\left\Vert D^{\gamma}\mathcal{E}_{\alpha}^{\mathbf{1}%
}f\right\Vert _{\infty}\leq\max_{\gamma\leq\alpha}\left(  \frac{2}{a}\right)
^{\gamma}\left(  \int_{a}^{b}\left(  1+x\right)  ^{\alpha}\left\vert v\left(
x\right)  \right\vert dx\right)  \max_{\beta\leq\alpha}\left\Vert D^{\beta
}f\right\Vert _{\infty,\mathcal{O}_{\mathbf{1}}}.
\]

From \ref{X84},%
\[
D^{\gamma}\left(  \mathcal{E}_{\alpha}^{\theta}\left[  v\right]  f\right)
\left(  x\right)  =\theta^{\gamma}D^{\gamma}\left\{  \mathcal{E}_{\alpha
}^{\mathbf{1}}\left[  v\left(  \theta.\right)  \right]  f\left(
\theta.\right)  \right\}  \left(  \theta.x\right)  ,
\]

where $\operatorname*{supp}v\subseteq R\left[  c,d\right]  \subset
\mathcal{O}_{\theta}$ and $c.<d$. But from \ref{X85}, $\operatorname*{supp}%
v\left(  \theta.\right)  \subseteq R\left[  \theta.c,\theta.d\right]
=R\left[  a,b\right]  \subset\mathcal{O}_{\mathbf{1}}$ where $a:=\min\left(
\theta.c,\theta.d\right)  $ and $b:=\max\left(  \theta.c,\theta.d\right)  $
uniquely satisfy $a.<b$, so that%
\begin{align*}
\max_{\gamma\leq\alpha} &  \left\Vert D^{\gamma}\mathcal{E}_{\alpha}^{\theta
}\left[  v\right]  f\right\Vert _{\infty}\\
&  =\max_{\gamma\leq\alpha}\left\Vert \theta^{\gamma}D^{\gamma}\left\{
\mathcal{E}_{\alpha}^{\mathbf{1}}\left[  v\left(  \theta.\right)  \right]
f\left(  \theta.\right)  \right\}  \left(  \theta.x\right)  \right\Vert
_{\infty}\\
&  =\max_{\gamma\leq\alpha}\left\Vert D^{\gamma}\left\{  \mathcal{E}_{\alpha
}^{\mathbf{1}}\left[  v\left(  \theta.\right)  \right]  f\left(
\theta.\right)  \right\}  \right\Vert _{\infty}\\
&  \leq\max_{\gamma\leq\alpha}\left(  \frac{2}{a}\right)  ^{\gamma}\left(
\int_{a}^{b}\left(  1+x\right)  ^{\alpha}\left\vert v\left(  \theta.x\right)
\right\vert dx\right)  \max_{\beta\leq\alpha}\left\Vert D^{\beta}\left(
f\left(  \theta.\right)  \right)  \right\Vert _{\infty,\mathcal{O}%
_{\mathbf{1}}}\\
&  =\max_{\gamma\leq\alpha}\left(  \frac{2}{a}\right)  ^{\gamma}\left(
\int_{a}^{b}\left(  1+x\right)  ^{\alpha}\left\vert v\left(  \theta.x\right)
\right\vert dx\right)  \max_{\beta\leq\alpha}\left\Vert \theta^{\beta}\left(
D^{\beta}f\right)  \left(  \theta.\right)  \right\Vert _{\infty,\mathcal{O}%
_{\mathbf{1}}}\\
&  =\max_{\gamma\leq\alpha}\left(  \frac{2}{a}\right)  ^{\gamma}\left(
\int_{a}^{b}\left(  1+x\right)  ^{\alpha}\left\vert v\left(  \theta.x\right)
\right\vert dx\right)  \max_{\beta\leq\alpha}\left\Vert \left(  D^{\beta
}f\right)  \left(  \theta.\right)  \right\Vert _{\infty,\mathcal{O}%
_{\mathbf{1}}}\\
&  =\max_{\gamma\leq\alpha}\left(  \frac{2}{a}\right)  ^{\gamma}\left(
\int_{a}^{b}\left(  1+x\right)  ^{\alpha}\left\vert v\left(  \theta.x\right)
\right\vert dx\right)  \max_{\beta\leq\alpha}\left\Vert D^{\beta}f\right\Vert
_{\infty,\mathcal{O}_{\theta}}\\
&  =\max_{\gamma\leq\alpha}\left(  \frac{2}{\min\left(  \theta.c,\theta
.d\right)  }\right)  ^{\gamma}\left(  \int_{c}^{d}\left(  1+x\right)
^{\alpha}\left\vert v\left(  \theta.x\right)  \right\vert dx\right)
\max_{\beta\leq\alpha}\left\Vert D^{\beta}f\right\Vert _{\infty,\mathcal{O}%
_{\theta}}\\
&  =\max_{\gamma\leq\alpha}\left(  \frac{2}{\min\left(  \theta.c,\theta
.d\right)  }\right)  ^{\gamma}\left(  \int_{c}^{d}\left(  1+\theta.z\right)
^{\alpha}\left\vert v\left(  z\right)  \right\vert dz\right)  \max_{\beta
\leq\alpha}\left\Vert D^{\beta}f\right\Vert _{\infty,\mathcal{O}_{\theta}},
\end{align*}

where the last step used the fact that Definition \ref{Def_topol_on_Rd}
implies $R\left[  \theta.a,\theta.b\right]  =R\left[  c,d\right]  $, and as a
consequence the change of variables%
\[
\int_{a}^{b}\left(  1+x\right)  ^{\alpha}\left\vert v\left(  \theta.x\right)
\right\vert dx=\int_{c}^{d}\left(  1+\theta.z\right)  ^{\alpha}\left\vert
v\left(  z\right)  \right\vert dz,
\]

is valid.\medskip

\textbf{Part ??6} ?? FINISH!\medskip

\textbf{Part ??7} From part 1%
\[
\left\vert \mathcal{E}_{\alpha}^{\theta}f\left(  x\right)  \right\vert
\leq\tfrac{1}{\left(  \alpha-1\right)  !}\sum\limits_{\beta\leq\alpha}%
\tbinom{\alpha}{\beta}\int_{\mathcal{O}_{\theta}}\left\vert J_{\alpha
-2}^{\theta}\left[  v_{\beta}\right]  \left(  y-x\right)  \right\vert
\left\vert D^{\beta}f\left(  y\right)  \right\vert dy.
\]

Using Peetre's inequality \ref{p33}: for any $\mu>0$,%
\begin{align*}
\left(  1+\left\vert x\right\vert \right)  ^{\mu}\left\vert \mathcal{E}%
_{\alpha}^{\theta}f\left(  x\right)  \right\vert  & \leq\tfrac{1}{\left(
\alpha-1\right)  !}\sum\limits_{\beta\leq\alpha}\tbinom{\alpha}{\beta}%
\int_{\mathcal{O}_{\theta}}\left(  1+\left\vert x\right\vert \right)  ^{\mu
}\left\vert J_{\alpha-2}^{\theta}\left[  v_{\beta}\right]  \left(  y-x\right)
\right\vert \left\vert D^{\beta}f\left(  y\right)  \right\vert dy\\
& =\tfrac{1}{\left(  \alpha-1\right)  !}\sum\limits_{\beta\leq\alpha}%
\tbinom{\alpha}{\beta}\int_{\mathcal{O}_{\theta}}\frac{\left(  1+\left\vert
x\right\vert \right)  ^{\mu}}{\left(  1+\left\vert y\right\vert \right)
^{\mu}}\left\vert J_{\alpha-2}^{\theta}\left[  v_{\beta}\right]  \left(
y-x\right)  \right\vert \left(  1+\left\vert y\right\vert \right)  ^{\mu
}\left\vert D^{\beta}f\left(  y\right)  \right\vert dy\\
& \leq\tfrac{1}{\left(  \alpha-1\right)  !}\sum\limits_{\beta\leq\alpha
}\tbinom{\alpha}{\beta}\int_{\mathcal{O}_{\theta}}\left(  1+\left\vert
y-x\right\vert \right)  ^{\mu}\left\vert J_{\alpha-2}^{\theta}\left[
v_{\beta}\right]  \left(  y-x\right)  \right\vert \left(  1+\left\vert
y\right\vert \right)  ^{\mu}\left\vert D^{\beta}f\left(  y\right)  \right\vert
dy\\
& \leq\tfrac{1}{\left(  \alpha-1\right)  !}\sum\limits_{\beta\leq\alpha
}\tbinom{\alpha}{\beta}\left\Vert \left(  1+\left\vert \cdot\right\vert
\right)  ^{\mu}J_{\alpha-2}^{\theta}\left[  v_{\beta}\right]  \right\Vert
_{1}\left\Vert \left(  1+\left\vert \cdot\right\vert \right)  ^{\mu}D^{\beta
}f\right\Vert _{\infty}\\
& =\tfrac{1}{\left(  \alpha-1\right)  !}\sum\limits_{\beta\leq\alpha}%
\tbinom{\alpha}{\beta}\left\Vert \left(  1+\left\vert \cdot\right\vert
\right)  ^{\mu}J_{\alpha-2}^{\theta}\left[  v_{\beta}\right]  \right\Vert
_{1,\overline{\mathcal{O}_{\theta}}}\left\Vert \left(  1+\left\vert
\cdot\right\vert \right)  ^{\mu}D^{\beta}f\right\Vert _{\infty,\overline
{\mathcal{O}_{\theta}}}\\
& <\infty,
\end{align*}

where the last inequality follows from the fact that $f\in C_{0}^{\left(
\alpha\right)  }\left(  \overline{\mathcal{O}_{\theta}}\right)  $ and
\ref{X32} implies $J_{\alpha-2}^{\theta}\left[  v_{\beta}\right]  $ decreases
exponentially at infinity.\medskip

\textbf{Part ??7} Suppose $x\in\mathcal{O}_{\omega}$ and $\omega\neq\theta$.%
\begin{align*}
e^{px}\left\vert \mathcal{E}_{\alpha}^{\theta}f\left(  x\right)  \right\vert
& \leq\tfrac{1}{\left(  \alpha-1\right)  !}\sum\limits_{\beta\leq\alpha
}\tbinom{\alpha}{\beta}\int_{\mathcal{O}_{\theta}}e^{px}J_{\alpha-2}^{\theta
}\left[  v_{\beta}\right]  \left(  y-x\right)  D^{\beta}f\left(  y\right)
dy\\
& =\tfrac{1}{\left(  \alpha-1\right)  !}\sum\limits_{\beta\leq\alpha}%
\tbinom{\alpha}{\beta}\int e^{p\left(  x-y\right)  }\left\vert J_{\alpha
-2}^{\theta}\left[  v_{\beta}\right]  \left(  y-x\right)  \right\vert
e^{py}\left\vert D^{\beta}f\left(  y\right)  \right\vert dy\\
& =\tfrac{1}{\left(  \alpha-1\right)  !}\sum\limits_{\beta\leq\alpha}%
\tbinom{\alpha}{\beta}\int e^{p\left(  x-y\right)  }\left\vert J_{\alpha
-2}^{\theta}\left[  v_{\beta}\right]  _{-}\left(  x-y\right)  \right\vert
e^{py}\left\vert D^{\beta}f\left(  y\right)  \right\vert dy\\
& \Rightarrow Young^{\prime}s\Rightarrow\\
& \leq\tfrac{1}{\left(  \alpha-1\right)  !}\sum\limits_{\beta\leq\alpha
}\tbinom{\alpha}{\beta}\left\Vert e^{pz}J_{\alpha-2}^{\theta}\left[  v_{\beta
}\right]  _{-}\right\Vert _{\infty}\left\Vert e^{py}D^{\beta}f\left(
y\right)  \right\Vert _{1}\\
& =\tfrac{1}{\left(  \alpha-1\right)  !}\sum\limits_{\beta\leq\alpha}%
\tbinom{\alpha}{\beta}\left\Vert e^{-pz}J_{\alpha-2}^{\theta}\left[  v_{\beta
}\right]  \right\Vert _{\infty}\left\Vert e^{py}D^{\beta}f\left(  y\right)
\right\Vert _{1,\mathcal{O}_{\theta}}\\
& =\tfrac{1}{\left(  \alpha-1\right)  !}\sum\limits_{\beta\leq\alpha}%
\tbinom{\alpha}{\beta}\left\Vert e^{-pz}J_{\alpha-2}^{\theta}\left[  v_{\beta
}\right]  \right\Vert _{\infty,\mathcal{O}_{\theta}}\left\Vert e^{py}D^{\beta
}f\left(  y\right)  \right\Vert _{1,\mathcal{O}_{\theta}},
\end{align*}

but from \ref{X88},%
\[
\left\vert J_{\delta}^{\theta}\left[  v_{\beta}\right]  \left(  z\right)
\right\vert \leq\left\Vert \frac{v_{\beta}}{x^{\mathbf{1}}}\right\Vert
_{1}q_{\mathbf{0}}^{\delta}\left(  \frac{\theta.z}{\min\left\{  \theta
.c,\theta.d\right\}  }\right)  e^{-\frac{\theta z}{\max\left\{  \theta
.c,\theta.d\right\}  }},\quad z\in\mathcal{O}_{\theta},
\]

so for any $x\in\mathbb{R}^{d}$,%
\begin{align*}
&  e^{px}\left\vert \mathcal{E}_{\alpha}^{\theta}f\left(  x\right)
\right\vert \\
&  \leq\tfrac{1}{\left(  \alpha-1\right)  !}\sum\limits_{\beta\leq\alpha
}\tbinom{\alpha}{\beta}\left\Vert e^{-pz}\left\Vert \frac{v_{\beta}%
}{x^{\mathbf{1}}}\right\Vert _{1}q_{\mathbf{0}}^{\alpha-2}\left(  \frac
{\theta.z}{\min\left\{  \theta.c,\theta.d\right\}  }\right)  e^{-\frac{\theta
z}{\max\left\{  \theta.c,\theta.d\right\}  }}\right\Vert _{\infty
,\mathcal{O}_{\theta}}\left\Vert e^{py}D^{\beta}f\left(  y\right)  \right\Vert
_{1,\mathcal{O}_{\theta}}\\
&  =\tfrac{1}{\left(  \alpha-1\right)  !}\sum\limits_{\beta\leq\alpha}%
\tbinom{\alpha}{\beta}\left\Vert \frac{v_{\beta}}{x^{\mathbf{1}}}\right\Vert
_{1}\left\Vert q_{\mathbf{0}}^{\alpha-2}\left(  \frac{\theta.z}{\min\left\{
\theta.c,\theta.d\right\}  }\right)  e^{-pz}e^{-\frac{\theta z}{\max\left\{
\theta.c,\theta.d\right\}  }}\right\Vert _{\infty,\mathcal{O}_{\theta}%
}\left\Vert e^{py}D^{\beta}f\left(  y\right)  \right\Vert _{1,\mathcal{O}%
_{\theta}}\\
&  =\tfrac{1}{\left(  \alpha-1\right)  !}\sum\limits_{\beta\leq\alpha}%
\tbinom{\alpha}{\beta}\left\Vert \frac{v_{\beta}}{x^{\mathbf{1}}}\right\Vert
_{1}\left\Vert q_{\mathbf{0}}^{\alpha-2}\left(  \frac{y}{\min\left\{
\theta.c,\theta.d\right\}  }\right)  e^{-p\text{ }\theta.y}e^{-\frac
{\mathbf{1}y}{\max\left\{  \theta.c,\theta.d\right\}  }}\right\Vert
_{\infty,\mathcal{O}_{\mathbf{1}}}\left\Vert e^{py}D^{\beta}f\left(  y\right)
\right\Vert _{1,\mathcal{O}_{\theta}}\\
&  =\tfrac{1}{\left(  \alpha-1\right)  !}\sum\limits_{\beta\leq\alpha}%
\tbinom{\alpha}{\beta}\left\Vert \frac{v_{\beta}}{x^{\mathbf{1}}}\right\Vert
_{1}\left\Vert q_{\mathbf{0}}^{\alpha-2}\left(  \frac{y}{\min\left\{
\theta.c,\theta.d\right\}  }\right)  e^{-p.\theta\text{ }y}e^{-\frac
{\mathbf{1}y}{\max\left\{  \theta.c,\theta.d\right\}  }}\right\Vert
_{\infty,\mathcal{O}_{\mathbf{1}}}\left\Vert e^{py}D^{\beta}f\left(  y\right)
\right\Vert _{1,\mathcal{O}_{\theta}}\\
&  =\tfrac{1}{\left(  \alpha-1\right)  !}\sum\limits_{\beta\leq\alpha}%
\tbinom{\alpha}{\beta}\left\Vert \frac{v_{\beta}}{x^{\mathbf{1}}}\right\Vert
_{1}\left\Vert q_{\mathbf{0}}^{\alpha-2}\left(  \frac{y}{\min\left\{
\theta.c,\theta.d\right\}  }\right)  e^{-\left(  p.\theta+\frac{\mathbf{1}%
}{\max\left\{  \theta.c,\theta.d\right\}  }\right)  y}\right\Vert
_{\infty,\mathcal{O}_{\mathbf{1}}}\left\Vert e^{py}D^{\beta}f\left(  y\right)
\right\Vert _{1,\mathcal{O}_{\theta}}\\
&  \leq\tfrac{1}{\left(  \alpha-1\right)  !}\left(  \sum\limits_{\beta
\leq\alpha}\tbinom{\alpha}{\beta}\left\Vert \frac{v_{\beta}}{x^{\mathbf{1}}%
}\right\Vert _{1}\left\Vert e^{py}D^{\beta}f\left(  y\right)  \right\Vert
_{1,\mathcal{O}_{\theta}}\right)  \left\Vert q_{\mathbf{0}}^{\alpha-2}\left(
\frac{y}{\min\left\{  \theta.c,\theta.d\right\}  }\right)  e^{-\left(
p.\theta+\frac{\mathbf{1}}{\max\left\{  \theta.c,\theta.d\right\}  }\right)
y}\right\Vert _{\infty,\mathcal{O}_{\mathbf{1}}}.
\end{align*}

Assume $x\in\overline{\mathcal{O}_{\omega}}$. Choose $s\in\mathbb{R}^{1}$ such
that $0<s<1$, and set
\begin{equation}
p^{\omega}:=\frac{s\omega}{\max\left\{  \theta.c,\theta.d\right\}  }%
\in\mathcal{O}_{\omega}.\label{X89}%
\end{equation}

Then%
\begin{align*}
p^{\omega}.\theta+\frac{\mathbf{1}}{\max\left\{  \theta.c,\theta.d\right\}
}=\frac{s\omega}{\max\left\{  \theta.c,\theta.d\right\}  }.\theta
+\frac{\mathbf{1}}{\max\left\{  \theta.c,\theta.d\right\}  } &  =\frac
{\mathbf{1}+s\omega.\theta}{\max\left\{  \theta.c,\theta.d\right\}  }\\
. &  >\frac{1-s}{\max\left\{  \theta.c,\theta.d\right\}  }\\
. &  >\mathbf{0},
\end{align*}

and so
\begin{align*}
\left\Vert q_{\mathbf{0}}^{\alpha-2}\left(  \frac{y}{\min\left\{
\theta.c,\theta.d\right\}  }\right)  e^{-\left(  p^{\omega}.\theta
+\frac{\mathbf{1}}{\max\left\{  \theta.c,\theta.d\right\}  }\right)
y}\right\Vert _{\infty,\mathcal{O}_{\mathbf{1}}} &  \leq\left\Vert
q_{\mathbf{0}}^{\alpha-2}\left(  \frac{y}{\min\left\{  \theta.c,\theta
.d\right\}  }\right)  e^{-\frac{\left(  1-s\right)  y}{\max\left\{
\theta.c,\theta.d\right\}  }}\right\Vert _{\infty,\mathcal{O}_{\mathbf{1}}}\\
&  <\infty.
\end{align*}

Consequently%
\begin{equation}
\left\vert \mathcal{E}_{\alpha}^{\theta}f\left(  x\right)  \right\vert \leq
C_{\omega}e^{-p^{\omega}x},\quad x\in\overline{\mathcal{O}_{\omega}%
},\label{X91}%
\end{equation}

where%
\begin{align}
C_{\omega}:=\tfrac{1}{\left(  \alpha-1\right)  !} &  \left(  \sum
\limits_{\beta\leq\alpha}\tbinom{\alpha}{\beta}\left\Vert \frac{v_{\beta}%
}{x^{\mathbf{1}}}\right\Vert _{1,\mathcal{O}_{\theta}}\left\Vert e^{p^{\omega
}y}D^{\beta}f\left(  y\right)  \right\Vert _{1,\mathcal{O}_{\theta}}\right)
\times\nonumber\\
&  \times\left\Vert q_{\mathbf{0}}^{\alpha-2}\left(  \frac{y}{\min\left\{
\theta.c,\theta.d\right\}  }\right)  e^{-\frac{\left(  1-s\right)  y}%
{\max\left\{  \theta.c,\theta.d\right\}  }}\right\Vert _{\infty,\mathcal{O}%
_{\mathbf{1}}}.\label{X90}%
\end{align}

Regarding \ref{X89},%
\[
p^{\omega}\leq\frac{s}{\max\left\{  \theta.c,\theta.d\right\}  }.
\]

Regarding \ref{X91}, since $x\in\overline{\mathcal{O}_{\omega}}$,%
\[
\left\vert \mathcal{E}_{\alpha}^{\theta}f\left(  x\right)  \right\vert \leq
C_{\omega}e^{-p^{\omega}x}=C_{\omega}e^{-\left(  \omega.p^{\omega}\right)
\left(  \omega.x\right)  }=C_{\omega}e^{-\frac{sx_{+}}{\max\left\{
\theta.c,\theta.d\right\}  }}.
\]

Hence%
\begin{align*}
C_{\omega} &  \leq\tfrac{1}{\left(  \alpha-1\right)  !}\left(  \sum
\limits_{\beta\leq\alpha}\tbinom{\alpha}{\beta}\left\Vert \frac{v_{\beta}%
}{x^{\mathbf{1}}}\right\Vert _{1,\mathcal{O}_{\theta}}\left\Vert e^{\frac
{sy}{\max\left\{  \theta.c,\theta.d\right\}  }}D^{\beta}f\left(  y\right)
\right\Vert _{1,\mathcal{O}_{\theta}}\right)  \times\\
&  \qquad\qquad\times\left\Vert q_{\mathbf{0}}^{\alpha-2}\left(  \frac{y}%
{\min\left\{  \theta.c,\theta.d\right\}  }\right)  e^{-\frac{\left(
1-s\right)  y}{\max\left\{  \theta.c,\theta.d\right\}  }}\right\Vert
_{\infty,\mathcal{O}_{\mathbf{1}}}\\
&  =\tfrac{\left(  \min\left\{  \theta.c,\theta.d\right\}  \right)
^{\mathbf{1}}}{\left(  \alpha-1\right)  !}\left(  \sum\limits_{\beta\leq
\alpha}\tbinom{\alpha}{\beta}\left\Vert \frac{v_{\beta}}{x^{\mathbf{1}}%
}\right\Vert _{1,\mathcal{O}_{\theta}}\left\Vert e^{\frac{sy}{\max\left\{
\theta.c,\theta.d\right\}  }}D^{\beta}f\left(  y\right)  \right\Vert
_{1,\mathcal{O}_{\theta}}\right)  \times\\
&  \qquad\qquad\times\left\Vert q_{\mathbf{0}}^{\alpha-2}\left(  z\right)
e^{-\left(  1-s\right)  \frac{\min\left\{  \theta.c,\theta.d\right\}  }%
{\max\left\{  \theta.c,\theta.d\right\}  }z}\right\Vert _{\infty
,\mathcal{O}_{\mathbf{1}}},
\end{align*}

which is independent of the orthant $\mathcal{O}_{\omega}$. Thus%
\[
\left\vert \mathcal{E}_{\alpha}^{\theta}f\left(  x\right)  \right\vert \leq
C_{s,f,\alpha}e^{-\frac{sx_{+}}{\max\left\{  \theta.c,\theta.d\right\}  }},
\]

where%
\begin{align}
C_{s,f,\alpha}:=\tfrac{\left(  \min\left\{  \theta.c,\theta.d\right\}
\right)  ^{\mathbf{1}}}{\left(  \alpha-1\right)  !} &  \left(  \sum
\limits_{\beta\leq\alpha}\tbinom{\alpha}{\beta}\left\Vert \frac{v_{\beta}%
}{x^{\mathbf{1}}}\right\Vert _{1,\mathcal{O}_{\theta}}\left\Vert e^{\frac
{sy}{\max\left\{  \theta.c,\theta.d\right\}  }}D^{\beta}f\left(  y\right)
\right\Vert _{1,\mathcal{O}_{\theta}}\right)  \times\nonumber\\
&  \times\left\Vert q_{\mathbf{0}}^{\alpha-2}\left(  z\right)  e^{-\left(
1-s\right)  \frac{\min\left\{  \theta.c,\theta.d\right\}  }{\max\left\{
\theta.c,\theta.d\right\}  }z}\right\Vert _{\infty,\mathcal{O}_{\mathbf{1}}%
}.\label{X92}%
\end{align}

Here $q_{\mathbf{0}}^{\alpha-2}$ is the tensor product polynomial in $z$ given
by \ref{X41}.
\end{proof}

\section{A continuous extension from $W^{n\mathbf{1}}\left(  \Omega\right)  $
to $W^{n\mathbf{1}}\mathbb{\ }\left(  \mathbb{R}^{d}\right)  $ - a Fourier
transform proof\label{Sect_ExtenLocWn1_to_Wn1_Fourier}}

In the next theorem we use a partition of unity and the integral
representation of Lemma \ref{Lem_SmthFuncIntegRepInOrthant} to construct a
continuous convolution extension operator $E_{\Omega}$ from $W^{n\mathbf{1}%
}\left(  \Omega\right)  $ to $W^{n\mathbf{1}}\left(  \mathbb{R}^{d}\right)  $.
Here $\Omega$ is bounded and satisfies the rectangle condition. This is done
in Theorem \ref{Thm_ExtenOrthantSobolFourier} and continuity is demonstrated
using a Fourier transform argument. With regard to the integral representation
lemma, the argument $v$ of the $J_{\delta}$ function is assumed to satisfy the
condition \ref{a030}.

First we need:

\begin{lemma}
\label{Lem_property_g_d}\textbf{Properties of the function: }%
\begin{equation}
g_{\delta}\mathbf{:=}\left(  H\left(  \zeta\right)  \zeta^{\delta}%
e^{-\zeta\mathbf{1}}\right)  ^{\wedge}\mathbf{,}\text{ }\delta\in
\mathbb{Z}^{d},\text{ }\delta\geq\mathbf{0},\label{X93}%
\end{equation}

where $H$ is the tensor product Heavyside step function.

\begin{enumerate}
\item $D^{j}\left(  H\left(  s\right)  s^{k}e^{-s}\right)  =H\left(  s\right)
p_{j}^{k-1}\left(  s\right)  e^{-s}\in L^{1}$ for $j\leq k$, and

$D^{k+1}\left(  H\left(  s\right)  s^{k}e^{-s}\right)  =k!\delta+H\left(
s\right)  p_{k+1}^{k-1}\left(  s\right)  e^{-s}$, where $p_{k}^{m}$ is defined
in \ref{X50}.

\begin{enumerate}
\item $H\left(  s\right)  s^{k}e^{-s}\in C_{B}^{\infty}\left(  \overline
{\mathcal{O}_{\mathbf{1}}}\right)  $;

\item $H\left(  s\right)  s^{k}e^{-s}\in C_{B}^{\left(  k-1\right)  }\cap
PWC_{B}^{\infty}$ $\forall k\geq\mathbf{1}$;

\item $H\left(  s\right)  s^{k}e^{-s}\in L^{p}$ for $1\leq p\leq\infty$.
\end{enumerate}

\item For all $\delta$: $g_{\delta}\in C_{BP}^{\infty}$ and $D^{\alpha
}g_{\delta}=\left(  -i\right)  ^{\left\vert \alpha\right\vert }g_{\alpha
+\delta}$ for all $\alpha$.

\item In one dimension: for $k=0,1,2,\ldots$%
\[
g_{k}\left(  s\right)  =\frac{1}{\sqrt{2\pi}}\frac{k!}{\left(  1+is\right)
^{k+1}},\quad\left\vert g_{k}\left(  s\right)  \right\vert =\frac{1}%
{\sqrt{2\pi}}\frac{k!}{\left(  1+s^{2}\right)  ^{\frac{k+1}{2}}}.
\]

Clearly $s^{t}g_{k}\in C_{B}^{\left(  0\right)  }$ when $0\leq t\leq k+1$.

\item If $k\geq0$ and $0\leq t\leq k+1$ then%
\begin{align}
\left\Vert s^{t}g_{k}\left(  s\right)  \right\Vert _{\infty}  & =\frac
{k!}{\sqrt{2\pi}}\left(  \frac{t}{k+1}\right)  ^{\frac{t}{2}}\left(
1-\frac{t}{k+1}\right)  ^{\frac{k+1-t}{2}}\label{X95}\\
& \leq\frac{k!}{\sqrt{2\pi}},\label{X97}%
\end{align}

where $0^{0}:=1$, and%
\[
\operatorname*{argmax}\limits_{s\geq0}\left\vert s^{t}g_{k}\left(  s\right)
\right\vert =\sqrt{\frac{t}{k+1-t}}.
\]

\item If $k\geq0$ and $0\leq t<k$ then%
\[
\int\left\vert s^{t}g_{k}\left(  s\right)  \right\vert =\frac{k!}{\sqrt{2\pi}%
}\frac{\Gamma\left(  \frac{t+1}{2}\right)  \Gamma\left(  \frac{k-t}{2}\right)
}{\Gamma\left(  \frac{k+1}{2}\right)  }=\frac{1}{\sqrt{2\pi}}\frac
{k!}{B\left(  \frac{t+1}{2},\frac{k-t-1}{2}\right)  }.
\]

\item $g_{\delta}\in C_{B}^{\infty}\cap H^{\infty}$ for all $\delta$.

\item If $\delta\geq\mathbf{0}$ then $g_{\delta}\in L^{1}$.

\item Alternative proof that $x^{\alpha}g_{\delta}\in C_{B}^{\left(  0\right)
}\left(  \mathbb{R}^{d}\right)  $ when $\alpha\leq\delta+1$.

\item Express $x^{\alpha}g_{\delta}$ in terms of a linear combination of the
$g_{\beta}$.%
\[
\xi^{\alpha}g_{\lambda}=\left(  -i\right)  ^{\left\vert \alpha\right\vert
}\sum_{\beta=0}^{\lambda-1}\left(  p_{\alpha}^{\lambda-2}\right)  _{\beta
}g_{\beta},\quad\alpha\leq\lambda-1,
\]

and in one dimension%
\[
t^{k}g_{k}\left(  t\right)  =\left(  -i\right)  ^{k+1}\left(  \frac{\left(
k-1\right)  !}{\left(  2\pi\right)  ^{d/2}}+\sum_{j=0}^{k-1}\left(
p_{k+1}^{k-1}\right)  _{j}g_{j}\right)  .
\]

\item For all multi-indexes $\alpha:\mathbf{0}\leq\alpha\leq\lambda$,%
\[
\tau^{\alpha}g_{\lambda}\left(  \tau\right)  =i^{\left\vert \alpha\right\vert
}\lambda!\sum_{\beta=0}^{\alpha}\left(  -1\right)  ^{\left\vert \beta
\right\vert }\tbinom{\alpha}{\beta}\frac{g_{\lambda-\beta}\left(  \tau\right)
}{\left(  \lambda-\beta\right)  !}.
\]

\end{enumerate}
\end{lemma}

\begin{proof}
\textbf{Part 1} ?? FINISH!.\medskip

\textbf{Part 2 }In one dimension, for all integer $m\geq0$, $s^{m}\left(
H\left(  s\right)  s^{k}e^{-s}\right)  \in L^{1}$ and so $D^{m}g_{k}%
=D^{m}\left(  H\left(  s\right)  s^{k}e^{-s}\right)  ^{\wedge}=??i^{m}\left(
s^{m}H\left(  s\right)  s^{k}e^{-s}\right)  ^{\wedge}\in C_{B}^{\left(
0\right)  }$ i.e. $g_{k}\in C_{B}^{\infty}$.

Also, for any derivative
\begin{align*}
D^{\alpha}g_{\delta}  & =D^{\alpha}\left(  H\left(  \zeta\right)
\zeta^{\delta}e^{-\zeta\mathbf{1}}\right)  ^{\wedge}=\left(  \left(
-i\right)  ^{\left\vert \alpha\right\vert }\zeta^{\alpha}H\left(
\zeta\right)  \zeta^{\delta}e^{-\zeta\mathbf{1}}\right)  ^{\wedge}=\\
& =\left(  -i\right)  ^{\left\vert \alpha\right\vert }\left(  H\left(
\zeta\right)  \zeta^{\delta+\alpha}e^{-\zeta\mathbf{1}}\right)  ^{\wedge
}=\left(  -i\right)  ^{\left\vert \alpha\right\vert }g_{\delta+\alpha}.
\end{align*}
\medskip

\textbf{Part 3} In one dimension%
\begin{align*}
\sqrt{2\pi}g_{0}\left(  s\right)   & =\int e^{-ist}H\left(  t\right)
e^{-t}dt=\int_{0}^{\infty}e^{-ist}e^{-t}dt=\int_{0}^{\infty}e^{-\left(
1+is\right)  t}dt=\\
& =\left[  \frac{-e^{-\left(  1+is\right)  t}}{\left(  1+is\right)  }\right]
_{0}^{\infty}=\frac{1}{1+is},
\end{align*}

and so from part 2,%
\begin{align*}
g_{k}\left(  s\right)  =i^{k}D^{k}g_{0}\left(  s\right)  =\frac{i^{k}}%
{\sqrt{2\pi}}D^{k}\frac{1}{1+is} &  =\frac{i^{k}}{\sqrt{2\pi}}\frac{\left(
-1\right)  ^{k}k!}{\left(  1+is\right)  ^{k+1}}i^{k}\\
&  =\frac{1}{\sqrt{2\pi}}\frac{k!}{\left(  1+is\right)  ^{k+1}}.
\end{align*}

Thus, in one dimension%
\[
\left\vert g_{k}\left(  s\right)  \right\vert =\frac{k!}{\sqrt{2\pi}}\frac
{1}{\left(  1+s^{2}\right)  ^{\frac{k+1}{2}}},
\]

and in higher dimensions%
\[
\left\vert g_{\delta}\left(  x\right)  \right\vert =\frac{\delta!}{\left(
2\pi\right)  ^{d/2}}\frac{1}{\left(  1+x.x\right)  ^{\frac{\delta+1}{2}}}.
\]
\medskip

\textbf{Part 4} If $0\leq t\leq k+1$,%
\begin{align*}
\max_{s}\left\vert s^{t}g_{k}\left(  s\right)  \right\vert  & =\max_{s}%
\frac{k!}{\sqrt{2\pi}}\frac{s^{t}}{\left(  1+s^{2}\right)  ^{\frac{k+1}{2}}%
}=\frac{k!}{\sqrt{2\pi}}\max_{s}\frac{s^{t}}{\left(  1+s^{2}\right)
^{\frac{k+1}{2}}}\leq\\
& \leq\frac{k!}{\sqrt{2\pi}}\max_{s}\frac{\left(  1+s^{2}\right)  ^{\frac
{t}{2}}}{\left(  1+s^{2}\right)  ^{\frac{k+1}{2}}}=\frac{k!}{\sqrt{2\pi}}%
\max_{s}\frac{1}{\left(  1+s^{2}\right)  ^{\frac{k+1-t}{2}}}=\\
& =\frac{k!}{\sqrt{2\pi}},
\end{align*}

which proves \ref{X97} directly.

If $0\leq t\leq k$,%
\begin{align*}
D_{s}\frac{s^{t}}{\left(  1+s^{2}\right)  ^{\frac{k+1}{2}}}  & =\frac{\left(
1+s^{2}\right)  ^{\frac{k+1}{2}}Ds^{t}-s^{t}D\left(  1+s^{2}\right)
^{\frac{k+1}{2}}}{\left(  1+s^{2}\right)  ^{k+1}}\\
& =\frac{\left(  1+s^{2}\right)  ^{\frac{k+1}{2}}ts^{t-1}-s^{t}\frac{k+1}%
{2}\left(  1+s^{2}\right)  ^{\frac{k-1}{2}}2s}{\left(  1+s^{2}\right)  ^{k+1}%
}\\
& =\frac{\left(  1+s^{2}\right)  ^{\frac{k+1}{2}}ts^{t-1}-s^{t+1}\left(
k+1\right)  \left(  1+s^{2}\right)  ^{\frac{k-1}{2}}}{\left(  1+s^{2}\right)
^{k+1}}\\
& =0,
\end{align*}

at a maximum.%
\begin{align*}
\left(  1+s^{2}\right)  ^{\frac{k+1}{2}}ts_{\max}^{t-1}  & =s_{\max}%
^{t+1}\left(  k+1\right)  \left(  1+s_{\max}^{2}\right)  ^{\frac{k-1}{2}},\\
\left(  1+s_{\max}^{2}\right)  t  & =\left(  k+1\right)  s_{\max}^{2}.\\
s_{\max}^{2}  & =\frac{t}{k+1-t}.
\end{align*}

More precisely,%
\begin{align*}
\max_{s}\left\vert s^{t}g_{k}\left(  s\right)  \right\vert  & =\max_{s}%
\frac{k!}{\sqrt{2\pi}}\frac{s^{t}}{\left(  1+s^{2}\right)  ^{\frac{k+1}{2}}%
}=\frac{k!}{\sqrt{2\pi}}\frac{\left(  \frac{t}{k+1-t}\right)  ^{\frac{t}{2}}%
}{\left(  1+\frac{t}{k+1-t}\right)  ^{\frac{k+1}{2}}}=\\
& =\frac{k!}{\sqrt{2\pi}}\frac{\left(  \frac{t}{k+1-t}\right)  ^{\frac{t}{2}}%
}{\left(  \frac{k+1}{k+1-t}\right)  ^{\frac{k+1}{2}}}=\frac{k!}{\sqrt{2\pi}%
}\left(  \frac{t}{k+1-t}\right)  ^{\frac{t}{2}}\left(  \frac{k+1-t}%
{k+1}\right)  ^{\frac{k+1}{2}}=\\
& =\frac{k!}{\sqrt{2\pi}}\left(  \frac{t}{k+1-t}\right)  ^{\frac{t}{2}}\left(
\frac{k+1-t}{k+1}\right)  ^{\frac{k+1}{2}}\\
& =\frac{k!}{\sqrt{2\pi}}t^{\frac{t}{2}}\left(  k+1-t\right)  ^{\frac
{k+1-t}{2}}\left(  k+2\right)  ^{-\frac{k+1}{2}}\\
& =\frac{k!}{\sqrt{2\pi}}t^{\frac{t}{2}}\left(  k+1\right)  ^{\frac{k+11-t}%
{2}}\left(  1-\frac{t}{k+1}\right)  ^{\frac{k+1-t}{2}}\left(  k+1\right)
^{-\frac{k+1}{2}}\\
& =\frac{k!}{\sqrt{2\pi}}\left(  \frac{t}{k+1}\right)  ^{\frac{t}{2}}\left(
1-\frac{t}{k+1}\right)  ^{\frac{k+1-t}{2}}.
\end{align*}

The formula also makes sense when $t=k+1$ since $0^{0}:=1$.\medskip

\textbf{Part 5} Suppose $0\leq t<k$. Then from part 3,%
\[
\int\left\vert s^{t}g_{k}\left(  s\right)  \right\vert =\int\frac{k!}%
{\sqrt{2\pi}}\frac{\left\vert s^{t}\right\vert ds}{\left(  1+s^{2}\right)
^{\frac{k+1}{2}}}=\frac{2k!}{\sqrt{2\pi}}\int_{0}^{\infty}\frac{s^{t}%
ds}{\left(  1+s^{2}\right)  ^{\frac{k+1}{2}}}.
\]

The change of variables: $s=\tan\theta$, $ds=\sec^{2}\theta d\theta$, and then
858.515 of Dwight \cite{Dwight61} gives%
\begin{align*}
\int\left\vert s^{t}g_{k}\left(  s\right)  \right\vert  & =\frac{2k!}%
{\sqrt{2\pi}}\int_{0}^{\pi/2}\frac{\tan^{t}\theta}{\sec^{k+1}\theta}\sec
^{2}\theta d\theta=\frac{2k!}{\sqrt{2\pi}}\int_{0}^{\pi/2}\tan^{t}\theta
\cos^{k-1}\theta d\theta=\\
& =\frac{2k!}{\sqrt{2\pi}}\int_{0}^{\pi/2}\sin^{t}\theta\cos^{k-1-t}\theta
d\theta=\frac{2k!}{\sqrt{2\pi}}\frac{\Gamma\left(  \frac{t+1}{2}\right)
\Gamma\left(  \frac{k-t}{2}\right)  }{2\Gamma\left(  \frac{k+1}{2}\right)
}=\\
& =\frac{k!}{\sqrt{2\pi}}\frac{\Gamma\left(  \frac{t+1}{2}\right)
\Gamma\left(  \frac{k-t}{2}\right)  }{\Gamma\left(  \frac{k+1}{2}\right)
}=\frac{1}{\sqrt{2\pi}}\frac{k!}{B\left(  \frac{t+1}{2},\frac{k-t}{2}\right)
}.
\end{align*}
\medskip

\textbf{Part 6} For all multi-indexes $\alpha$, $H\left(  \zeta\right)
\zeta^{\delta+\alpha}e^{-\zeta\mathbf{1}}\in L^{1}\cap L^{2}$. Hence $\left(
H\left(  \zeta\right)  \zeta^{\delta+\alpha}e^{-\zeta\mathbf{1}}\right)
^{\wedge}\in C_{B}^{\left(  0\right)  }\cap L^{2}$ and

$D^{\alpha}\left(  H\left(  \zeta\right)  \zeta^{\delta}e^{-\zeta\mathbf{1}%
}\right)  ^{\wedge}\in C_{B}^{\left(  0\right)  }\cap L^{2}$.\medskip

\textbf{Part 7} Follows directly from part 5.\medskip

\textbf{Part 8} ?? Noting that $x^{\alpha}g_{\delta}\left(  x\right)  $ is a
tensor product function we first consider one dimension.

Observe that $H\left(  s\right)  s^{k}\in C_{BP}^{\left(  k-1\right)  }\left(
\mathbb{R}^{1}\right)  $ when $k\geq1$, $D^{k}\left(  H\left(  s\right)
s^{k}\right)  =k!H\left(  s\right)  $ when $k\geq0$, and $D^{k+1}\left(
H\left(  s\right)  s^{k}\right)  =k!\delta$ when $k\geq0$.

Thus $D^{j}\left(  H\left(  s\right)  s^{k}e^{-s}\right)  \in C_{B}^{\left(
0\right)  }\cap L^{1}$ when $j\leq k-1$, $D^{k}\left(  H\left(  s\right)
s^{k}e^{-s}\right)  \in\left(  k!H\left(  s\right)  +C_{B}^{\left(
k-1\right)  }\right)  e^{-s}\subset L^{1}$ and

$D^{k+1}\left(  H\left(  s\right)  s^{k}e^{-s}\right)  \in k!\delta+\left(
cH+C_{B}^{\left(  k-1\right)  }\right)  e^{-s}\subset k!\delta+L^{1}$.

Thus $\left(  D^{j}\left(  H\left(  s\right)  s^{k}e^{-s}\right)  \right)
^{\wedge}\in C_{B}^{\left(  0\right)  }$ for $j\leq k+1$ i.e. $t^{j}\left(
H\left(  s\right)  s^{k}e^{-s}\right)  ^{\wedge}\left(  t\right)  \in
C_{B}^{\left(  0\right)  }$ when $j\leq k+1$.

But $\left(  H\left(  s\right)  s^{k}e^{-s}\right)  ^{\wedge}\left(  t\right)
\in C_{BP}^{\infty}$ so we must have $\left(  H\left(  s\right)  s^{k}%
e^{-s}\right)  ^{\wedge}\left(  t\right)  \in C_{B}^{\left(  0\right)  }$ when
$j\leq k+1$.\medskip

\textbf{Part 9} (Part 5 of Lemma \ref{Lem_property_J}) From part 1, in one
dimension $D^{j}\left(  H\left(  s\right)  s^{k}e^{-s}\right)  =H\left(
s\right)  p_{j}^{k-1}\left(  s\right)  e^{-s}$ when $j\leq k$, so%
\[
t^{j}g_{k}\left(  t\right)  =t^{j}\left(  H\left(  s\right)  s^{k}%
e^{-s}\right)  ^{\wedge}=\left(  -i\right)  ^{j}\left(  D^{j}\left(  H\left(
s\right)  s^{k}e^{-s}\right)  \right)  ^{\wedge}=\left(  -i\right)
^{j}\left(  H\left(  s\right)  p_{j}^{k-1}\left(  s\right)  e^{-s}\right)
^{\wedge},
\]

and now noting \ref{X50} and \ref{X72} we get for $\alpha\leq\lambda$,%
\begin{align*}
\tau^{\alpha}g_{\lambda}  & =\left(  -i\right)  ^{\left\vert \alpha\right\vert
}\left(  H\left(  x\right)  p_{\alpha}^{\lambda-1}\left(  x\right)
e^{-x\mathbf{1}}\right)  ^{\wedge}=\left(  -i\right)  ^{\left\vert
\alpha\right\vert }\left(  H\left(  x\right)  p_{\alpha}^{\lambda-1}\left(
x\right)  e^{-x\mathbf{1}}\right)  ^{\wedge}=\\
& =\left(  -i\right)  ^{\left\vert \alpha\right\vert }\sum_{\beta=0}^{\lambda
}\left(  p_{\alpha}^{\lambda-1}\right)  _{\beta}\left(  H\left(  x\right)
x^{\beta}e^{-x\mathbf{1}}\right)  ^{\wedge}=\left(  -i\right)  ^{\left\vert
\alpha\right\vert }\sum_{\beta=0}^{\lambda}\left(  p_{\alpha}^{\lambda
-1}\right)  _{\beta}g_{\beta}.
\end{align*}

When $j=k+1$, $D^{k+1}\left(  H\left(  s\right)  s^{k}e^{-s}\right)
=k!\delta+H\left(  s\right)  p_{k+1}^{k-1}\left(  s\right)  e^{-s}$ and%
\begin{align*}
t^{k+1}g_{k}\left(  t\right)  =t^{k+1}\left(  H\left(  s\right)  s^{k}%
e^{-s}\right)  ^{\wedge} &  =\left(  -i\right)  ^{k+1}\left(  D^{k+1}\left(
H\left(  s\right)  s^{k}e^{-s}\right)  \right)  ^{\wedge}\\
&  =\left(  -i\right)  ^{k+1}\left(  k!\delta+H\left(  s\right)  p_{k+1}%
^{k-1}\left(  s\right)  e^{-s}\right)  ^{\wedge}\\
&  =\left(  -i\right)  ^{k+1}\left(  \frac{k!}{\left(  2\pi\right)  ^{d/2}%
}+\sum_{j=0}^{k}\left(  p_{k+1}^{k-1}\right)  _{j}g_{j}\right)  .
\end{align*}
\medskip

\textbf{Part 10}
\begin{align*}
\frac{\tau^{\alpha}}{\left(  1+i\tau\right)  ^{\lambda}}=\frac{\left(
-i\right)  ^{\left\vert \alpha\right\vert }\left(  i\tau\right)  ^{\alpha}%
}{\left(  1+i\tau\right)  ^{\lambda}} &  =\left(  -i\right)  ^{\left\vert
\alpha\right\vert }\frac{\left(  1+i\tau-1\right)  ^{\alpha}}{\left(
1+i\tau\right)  ^{\lambda}}\\
&  =\left(  -i\right)  ^{\left\vert \alpha\right\vert }\sum_{\beta=0}^{\alpha
}\tbinom{\alpha}{\beta}\frac{\left(  -1\right)  ^{\left\vert \alpha
-\beta\right\vert }\left(  1+i\tau\right)  ^{\beta}}{\left(  1+i\tau\right)
^{\lambda}}\\
&  =i^{\left\vert \alpha\right\vert }\sum_{\beta=0}^{\alpha}\frac{\left(
-1\right)  ^{\left\vert \beta\right\vert }\tbinom{\alpha}{\beta}}{\left(
1+i\tau\right)  ^{\lambda-\beta}},
\end{align*}

and since from part 3, $g_{k}\left(  s\right)  =\frac{1}{\sqrt{2\pi}}\frac
{k!}{\left(  1+is\right)  ^{k+1}}$, it follows that
\begin{align*}
\tau^{\alpha}g_{\lambda}\left(  \tau\right)  =\frac{\lambda!}{\left(
2\pi\right)  ^{d/2}}\frac{\tau^{\alpha}}{\left(  1+i\tau\right)  ^{\lambda+1}}
&  =\frac{\lambda!}{\left(  2\pi\right)  ^{d/2}}i^{\left\vert \alpha
\right\vert }\sum_{\beta=0}^{\alpha}\frac{\left(  -1\right)  ^{\left\vert
\beta\right\vert }\tbinom{\alpha}{\beta}}{\left(  1+i\tau\right)
^{\lambda-\beta+1}}\\
&  =i^{\left\vert \alpha\right\vert }\lambda!\sum_{\beta=0}^{\alpha}\left(
-1\right)  ^{\left\vert \beta\right\vert }\tbinom{\alpha}{\beta}%
\frac{g_{\lambda-\beta}\left(  \tau\right)  }{\left(  \lambda-\beta\right)
!}\\
&  =i^{\left\vert \alpha\right\vert }\lambda!\sum_{\beta=0}^{\alpha}\left(
-1\right)  ^{\left\vert \beta\right\vert }\tbinom{\alpha}{\beta}%
\frac{g_{\lambda-\beta}\left(  \tau\right)  }{\left(  \lambda-\beta\right)
!}.
\end{align*}

\end{proof}

Regarding $J_{\delta}$, I will avoid differentiation under the integral sign
and instead use the Fourier transform.

\begin{lemma}
\label{Lem_FourierJd[u]_formulas} If $\sigma\in L_{\overline{\mathcal{O}%
_{\mathbf{1}}}}^{1}$ and $\delta\geq-\mathbf{1}$ then%
\[
\widehat{J_{\delta}\left[  \sigma\right]  }\left(  \xi\right)  =\int%
_{\mathcal{O}_{\mathbf{1}}}\widehat{\sigma}\left(  \xi.\zeta\right)
\zeta^{\delta+1}e^{-\zeta\mathbf{1}}d\zeta=\int_{\mathcal{O}_{\mathbf{1}}%
}g_{\delta+1}\left(  \xi.\mu\right)  \sigma\left(  \mu\right)  d\mu,
\]

where the function $g_{\delta}$ was introduced in \ref{X93}.
\end{lemma}

\begin{proof}
From part 3 of Lemma \ref{Lem_Jd[f]_bnd_f_in_L1}, $J_{\delta}\left[
\sigma\right]  \in L^{1}$ and so $\widehat{J_{\delta}\left[  \sigma\right]
}\in C_{B}^{\left(  0\right)  }$. Consequently%
\begin{align*}
\widehat{J_{\delta}\left[  \sigma\right]  }\left(  \xi\right)   & =\left(
2\pi\right)  ^{-d/2}\int_{\mathcal{O}_{\mathbf{1}}}e^{-i\xi z}J_{\delta
}\left[  \sigma\right]  \left(  z\right)  dz\\
& =\left(  2\pi\right)  ^{-d/2}\int_{\mathcal{O}_{\mathbf{1}}}e^{-i\xi
z}z^{\delta+1}\int_{\mathcal{O}_{\mathbf{1}}}\frac{1}{\mu^{\delta+2}}%
e^{-\frac{z}{\mu}\mathbf{1}}\sigma\left(  \mu\right)  d\mu\text{ }dz\\
& =\left(  2\pi\right)  ^{-d/2}\int_{\mathcal{O}_{\mathbf{1}}}\int%
_{\mathcal{O}_{\mathbf{1}}}e^{-i\xi z}z^{\delta+1}\frac{e^{-\frac{z}{\mu
}\mathbf{1}}}{\mu^{\delta+2}}\sigma\left(  \mu\right)  d\mu\text{ }dz,
\end{align*}

and this integral is absolutely convergent since,%
\begin{align*}
\int_{\mathcal{O}_{\mathbf{1}}}\int_{\mathcal{O}_{\mathbf{1}}}\left\vert
e^{-i\xi z}z^{\delta+1}\frac{e^{-\frac{z}{\mu}\mathbf{1}}}{\mu^{\delta+2}%
}\sigma\left(  \mu\right)  \right\vert dz\text{ }d\mu & \leq\int%
_{\mathcal{O}_{\mathbf{1}}}\int_{\mathcal{O}_{\mathbf{1}}}z^{\delta
+1}e^{-\frac{z}{\mu}\mathbf{1}}dz\frac{\left\vert \sigma\left(  \mu\right)
\right\vert }{\mu^{\delta+2}}d\mu\\
& :chg.\text{ }var.\text{ }\zeta=z./\mu,\text{ }dz=\mu^{\mathbf{1}}%
d\zeta\Rightarrow\\
& =\int_{\mathcal{O}_{\mathbf{1}}}\int_{\mathcal{O}_{\mathbf{1}}}\mu
^{\delta+2}.\zeta^{\delta+1}e^{-\zeta\mathbf{1}}dz\frac{\left\vert
\sigma\left(  \mu\right)  \right\vert }{\mu^{\delta+2}}d\mu\\
& =\int_{\mathcal{O}_{\mathbf{1}}}\int_{\mathcal{O}_{\mathbf{1}}}\zeta
^{\delta+1}e^{-\zeta\mathbf{1}}d\zeta\left\vert \sigma\left(  \mu\right)
\right\vert d\mu\\
& =\left(  \delta+1\right)  !\left\Vert \sigma\right\Vert _{1}.
\end{align*}

This allows us to use Fubini's theorem to justify changing the order of
integration. Indeed, for all multi-integers $\delta\geq-\mathbf{1}$,%
\begin{align}
\widehat{J_{\delta}\left[  \sigma\right]  }\left(  \xi\right)   & =\left(
2\pi\right)  ^{-d/2}\int_{\mathcal{O}_{\mathbf{1}}}\int_{\mathcal{O}%
_{\mathbf{1}}}e^{-i\xi z}z^{\delta+1}\frac{e^{-\frac{z}{\mu}\mathbf{1}}}%
{\mu^{\delta+2}}\sigma\left(  \mu\right)  d\mu\text{ }dz\nonumber\\
& =\left(  2\pi\right)  ^{-d/2}\int_{\mathcal{O}_{\mathbf{1}}}\int%
_{\mathcal{O}_{\mathbf{1}}}e^{-i\xi z}z^{\delta+1}e^{-\frac{z}{\mu}\mathbf{1}%
}dz\text{ }\frac{\sigma\left(  \mu\right)  }{\mu^{\delta+2}}d\mu\nonumber\\
& \Rightarrow\zeta=z./\mu,\text{ }z=\mu.\zeta,\text{ }dz=\mu^{\mathbf{1}%
}d\zeta\Rightarrow\nonumber\\
\widehat{J_{\delta}\left[  \sigma\right]  }\left(  \xi\right)   & =\left(
2\pi\right)  ^{-d/2}\int_{\mathcal{O}_{\mathbf{1}}}\int_{\mathcal{O}%
_{\mathbf{1}}}e^{-i\xi\left(  \mu.\zeta\right)  }\left(  \mu.\zeta\right)
^{\delta+1}e^{-\zeta\mathbf{1}}\mu^{\mathbf{1}}d\zeta\text{ }\frac
{\sigma\left(  \mu\right)  }{\mu^{\delta+2}}d\mu\nonumber\\
& =\left(  2\pi\right)  ^{-d/2}\int_{\mathcal{O}_{\mathbf{1}}}\int%
_{\mathcal{O}_{\mathbf{1}}}\mu^{\delta+2}e^{-i\left(  \xi.\zeta\right)  \mu
}\zeta^{\delta+1}e^{-\zeta\mathbf{1}}d\zeta\text{ }\frac{\sigma\left(
\mu\right)  }{\mu^{\delta+2}}d\mu\nonumber\\
& =\left(  2\pi\right)  ^{-d/2}\int_{\mathcal{O}_{\mathbf{1}}}\left(
\int_{\mathcal{O}_{\mathbf{1}}}e^{-i\left(  \xi.\zeta\right)  \mu}%
\zeta^{\delta+1}e^{-\zeta\mathbf{1}}d\zeta\right)  \sigma\left(  \mu\right)
d\mu\label{X18}\\
& \Rightarrow clearly\text{ }absolutely\text{ }convergent\Rightarrow
\nonumber\\
& =\int_{\mathcal{O}_{\mathbf{1}}}\left(  \left(  2\pi\right)  ^{-d/2}%
\int_{\mathcal{O}_{\mathbf{1}}}e^{-i\left(  \xi.\zeta\right)  \mu}%
\sigma\left(  \mu\right)  d\mu\right)  \zeta^{\delta+1}e^{-\zeta\mathbf{1}%
}d\zeta\nonumber\\
& =\int_{\mathcal{O}_{\mathbf{1}}}\widehat{\sigma}\left(  \xi.\zeta\right)
\zeta^{\delta+1}e^{-\zeta\mathbf{1}}d\zeta,\nonumber
\end{align}

as claimed. Also, from \ref{X18},%
\begin{align*}
\widehat{J_{\delta}\left[  \sigma\right]  }\left(  \xi\right)   & =\left(
2\pi\right)  ^{-d/2}\int_{\mathcal{O}_{\mathbf{1}}}\left(  \int_{\mathcal{O}%
_{\mathbf{1}}}e^{-i\left(  \xi.\zeta\right)  \mu}\zeta^{\delta+1}%
e^{-\zeta\mathbf{1}}d\zeta\right)  \sigma\left(  \mu\right)  d\mu\\
& =\left(  2\pi\right)  ^{-d/2}\int_{\mathcal{O}_{\mathbf{1}}}\left(
\int_{\mathcal{O}_{\mathbf{1}}}e^{-i\left(  \xi.\mu\right)  \zeta}%
\zeta^{\delta+1}e^{-\zeta\mathbf{1}}d\zeta\right)  \sigma\left(  \mu\right)
d\mu\\
& =\int_{\mathcal{O}_{\mathbf{1}}}\left(  \left(  2\pi\right)  ^{-d/2}%
\int_{\mathcal{O}_{\mathbf{1}}}e^{-i\left(  \xi.\mu\right)  \zeta}%
\zeta^{\delta+1}e^{-\zeta\mathbf{1}}d\zeta\right)  \sigma\left(  \mu\right)
d\mu\\
& =\int_{\mathcal{O}_{\mathbf{1}}}\left(  \left(  2\pi\right)  ^{-d/2}%
\int_{\mathcal{O}_{\mathbf{1}}}e^{-i\left(  \xi.\mu\right)  \zeta}H\left(
\zeta\right)  \zeta^{\delta+1}e^{-\zeta\mathbf{1}}d\zeta\right)  \sigma\left(
\mu\right)  d\mu\\
& =\int_{\mathcal{O}_{\mathbf{1}}}\left(  H\left(  \zeta\right)  \zeta
^{\delta+1}e^{-\zeta\mathbf{1}}\right)  ^{\wedge}\left(  \xi.\mu\right)
\sigma\left(  \mu\right)  d\mu\\
& =\int_{\mathcal{O}_{\mathbf{1}}}g_{\delta+1}\left(  \xi.\mu\right)
\sigma\left(  \mu\right)  d\mu.
\end{align*}

\end{proof}

?? \textbf{NEED TO ADD }a result for $v\in L_{0}^{1}\left(  \mathcal{O}%
_{1}\right)  $. This is needed in Theorem \ref{Thm_ExtenOrthantSobolFourier}.

\begin{lemma}
\label{Lem_property_xF[Jb]} \textbf{Properties of} $\widehat{J_{\delta}}$:
Suppose $\delta\geq-\mathbf{1}$. Then:

\begin{enumerate}
\item $L^{\infty}$ \textbf{bound}: Suppose $x^{-\left(  \delta+1\right)  }v\in
L_{\overline{\mathcal{O}_{\mathbf{1}}}}^{1}$ and $x^{\delta+1}v\in
L_{\overline{\mathcal{O}_{\mathbf{1}}}}^{1}$. Then $x^{\gamma}v\in
L_{\overline{\mathcal{O}_{\mathbf{1}}}}^{1}$ for all multi-integers $-\left(
\delta+1\right)  \leq\gamma\leq\delta+1$ i.e. $\gamma_{+}\leq\delta+1$.

Further%
\begin{align*}
\left\vert \xi^{\beta}\widehat{J_{\delta}\left[  x^{\alpha}v\right]  }\left(
\xi\right)  \right\vert \leq\left\Vert \eta^{\beta}g_{\delta+1}\right\Vert
_{\infty}\left\Vert y^{\alpha-\beta}v\right\Vert _{1,\mathcal{O}_{\mathbf{1}}%
}\leq\frac{\left(  \delta+1\right)  !}{\left(  2\pi\right)  ^{d/2}} &
\left\Vert y^{\alpha-\beta}v\right\Vert _{1,\mathcal{O}_{\mathbf{1}}},\\
\mathbf{0} &  \leq\alpha,\beta\leq\delta+1.
\end{align*}

\item $L^{1}$ \textbf{bound}: Suppose $\mathbf{0}\leq\beta\leq\delta$. Then
$\sigma\in L_{\overline{\mathcal{O}_{\mathbf{1}}}}^{1}$ and $y^{-\left(
\beta+1\right)  }\sigma\in L_{\overline{\mathcal{O}_{\mathbf{1}}}}^{1}$ imply%
\[
\int\left\vert \xi^{\beta}\widehat{J_{\delta}\left[  \sigma\right]  }\left(
\xi\right)  \right\vert d\xi\leq\left\Vert \eta^{\beta}g_{\delta+1}\right\Vert
_{1}\left\Vert \frac{\sigma}{y^{\beta+1}}\right\Vert _{1,\mathcal{O}%
_{\mathbf{1}}}\leq\frac{\left(  \delta+1\right)  !}{\left(  2\pi\right)
^{d/2}}\left\Vert \frac{\sigma}{y^{\beta+1}}\right\Vert _{1,\mathcal{O}%
_{\mathbf{1}}}.
\]

\item ?? If $\sigma\in??L^{1}\left(  \mathcal{O}_{\mathbf{1}}\right)  $ and
$y^{-\left(  \delta+2\right)  }\sigma\in??L^{1}\left(  \mathcal{O}%
_{\mathbf{1}}\right)  $ then $J_{\delta}\in C_{B}^{\left(  \delta+1\right)
}\left(  \mathbb{R}^{d}\right)  $.

\item ?? Suppose $\delta\geq\mathbf{0}$. Then $\sigma\in L^{1}??\left(
\mathcal{O}_{\mathbf{1}}\right)  $ and $y^{-\left(  \delta+1\right)  }%
\sigma\in L^{1}\left(  \mathcal{O}_{\mathbf{1}}\right)  $ imply $J_{\delta}\in
C_{B}^{\left(  \delta\right)  }\left(  \mathbb{R}^{d}\right)  $.

\item $L^{\infty}$ \textbf{bound}: Suppose $x^{-\left(  \delta+2\right)  }v\in
L_{\overline{\mathcal{O}_{\mathbf{1}}}}^{1}$ and $x^{\delta+2}v\in
L_{\overline{\mathcal{O}_{\mathbf{1}}}}^{1}$. Then $x^{\gamma}v\in
L_{\overline{\mathcal{O}_{\mathbf{1}}}}^{1}$ for all multi-integers $-\left(
\delta+2\right)  \leq\gamma\leq\delta+2$ i.e. $\gamma_{+}\leq\delta+2$.

Further%
\begin{align*}
\left\vert \xi^{\beta}\widehat{J_{\delta}\left[  x^{\alpha}v\right]  }\left(
\xi\right)  \right\vert \leq\left\Vert \eta^{\beta}g_{\delta+1}\right\Vert
_{\infty}\left\Vert y^{\alpha-\beta}v\right\Vert _{1,\mathcal{O}_{\mathbf{1}}%
}\leq\frac{\left(  \delta+1\right)  !}{\left(  2\pi\right)  ^{d/2}} &
\left\Vert y^{\alpha-\beta}v\right\Vert _{1,\mathcal{O}_{\mathbf{1}}},\\
\mathbf{0} &  \leq\alpha,\beta\leq\delta+2.
\end{align*}

\end{enumerate}
\end{lemma}

\begin{proof}
\textbf{Part 1} From Lemma \ref{Lem_FourierJd[u]_formulas}, since $x^{\alpha
}v\in L_{\overline{\mathcal{O}_{\mathbf{1}}}}^{1}$,%
\begin{align*}
\xi^{\beta}\widehat{J_{\delta}\left[  x^{\alpha}v\right]  }\left(  \xi\right)
=\xi^{\beta}\int_{\mathcal{O}_{\mathbf{1}}}y^{\alpha}v\left(  y\right)
g_{\delta+1}\left(  y.\xi\right)  dy &  =\int_{\mathcal{O}_{\mathbf{1}}%
}y^{\alpha}v\left(  y\right)  \xi^{\beta}g_{\delta+1}\left(  y.\xi\right)
dy\\
&  =\int_{\mathcal{O}_{\mathbf{1}}}y^{\alpha-\beta}v\left(  y\right)  \left(
y.\xi\right)  ^{\beta}g_{\delta+1}\left(  y.\xi\right)  dy,
\end{align*}

so that%
\begin{align*}
\left\vert \xi^{\beta}\widehat{J_{\delta}\left[  x^{\alpha}v\right]  }\left(
\xi\right)  \right\vert  & \leq\int_{\mathcal{O}_{\mathbf{1}}}\left\vert
y^{\alpha-\beta}v\left(  y\right)  \right\vert \left\vert \left(
y.\xi\right)  ^{\beta}g_{\delta+1}\left(  y.\xi\right)  \right\vert dy\\
& \leq\left\Vert \eta^{\beta}g_{\delta+1}\right\Vert _{\infty}\left\Vert
y^{\alpha-\beta}v\right\Vert _{1,\mathcal{O}_{\mathbf{1}}}\\
& \leq\frac{\left(  \delta+1\right)  !}{\left(  2\pi\right)  ^{d/2}}\left\Vert
y^{\alpha-\beta}v\right\Vert _{1,\mathcal{O}_{\mathbf{1}}},
\end{align*}

by part 4 of Lemma \ref{Lem_property_g_d} (Note: $\left(  \alpha-\beta\right)
_{+}\leq\delta+1$ and $\beta\leq\delta+1<\delta+2$).\medskip

\textbf{Part 2}%
\begin{align*}
\int\left\vert \xi^{\beta}\widehat{J_{\delta}\left[  \sigma\right]  }\left(
\xi\right)  \right\vert d\xi & \leq\int\int_{\mathcal{O}_{\mathbf{1}}}%
\frac{\left\vert \sigma\left(  y\right)  \right\vert }{y^{\beta}}\left\vert
\left(  y.\xi\right)  ^{\beta}g_{\delta+1}\left(  y.\xi\right)  \right\vert
dyd\xi\\
& =\int_{\mathcal{O}_{\mathbf{1}}}\frac{\left\vert \sigma\left(  y\right)
\right\vert }{y^{\beta}}\int\left\vert \left(  y.\xi\right)  ^{\beta}%
g_{\delta+1}\left(  y.\xi\right)  \right\vert d\xi\text{ }dy\\
& \Rightarrow\eta=y.\xi,\text{ }d\eta=\left\vert y^{\mathbf{1}}\right\vert
d\xi\Rightarrow\\
& =\int_{\mathcal{O}_{\mathbf{1}}}\frac{\left\vert \sigma\left(  y\right)
\right\vert }{y^{\beta+1}}\int\left\vert \eta^{\beta}g_{\delta+1}\left(
\eta\right)  \right\vert d\eta\text{ }dy\\
& =\int\left\vert \eta^{\beta}g_{\delta+1}\left(  \eta\right)  \right\vert
d\eta\left(  \int_{\mathcal{O}_{\mathbf{1}}}\frac{\left\vert \sigma\left(
y\right)  \right\vert }{y^{\beta+1}}dy\right) \\
& \leq\frac{\left(  \delta+1\right)  !}{\left(  2\pi\right)  ^{d/2}}%
\int_{\mathcal{O}_{\mathbf{1}}}\frac{\left\vert \sigma\left(  y\right)
\right\vert }{y^{\beta+1}}dy,
\end{align*}

by part 4 of Lemma \ref{Lem_property_g_d}.\medskip

\textbf{Part 3} ?? Suppose $\mathbf{0}\leq\beta\leq\delta+3$. Then $\sigma\in
L^{1}\left(  \mathcal{O}_{\mathbf{1}}\right)  $ and $y^{-\left(
\beta+1\right)  }\sigma\in L^{1}\left(  \mathcal{O}_{\mathbf{1}}\right)  $.
Part 2 now implies%
\[
\int\left\vert \xi^{\beta}\widehat{J_{\delta}\left[  \sigma\right]  }\left(
\xi\right)  \right\vert d\xi\leq\left\Vert \frac{\sigma}{y^{\beta+1}%
}\right\Vert _{1,\mathcal{O}_{\mathbf{1}}}\left\Vert \eta^{\beta}g_{\delta
+1}\right\Vert _{1}.
\]

which means that $\xi^{\beta}\widehat{J_{\delta}\left[  \sigma\right]  }\in
L^{1}$ i.e. $D^{\beta}J_{\delta}\left[  \sigma\right]  \in C_{B}^{\left(
0\right)  }$ i.e. $J_{\delta}\left[  \sigma\right]  \in C_{B}^{\left(
\delta+1\right)  }$.\medskip

\textbf{Part 4} Suppose $\mathbf{0}\leq\beta\leq\delta$. Then $\sigma\in
L^{1}\left(  \mathcal{O}_{\mathbf{1}}\right)  $ and $y^{-\left(
\beta+1\right)  }\sigma\in L^{1}\left(  \mathcal{O}_{\mathbf{1}}\right)  $.
Hence part 2 implies%
\[
\int\left\vert \xi^{\beta}\widehat{J_{\delta}\left[  \sigma\right]  }\left(
\xi\right)  \right\vert d\xi\leq\left\Vert \frac{\sigma}{y^{\beta+1}%
}\right\Vert _{1,\mathcal{O}_{\mathbf{1}}}\left\Vert \eta^{\beta}g_{\delta
+1}\right\Vert _{1},
\]

which implies $\widehat{D^{\beta}J_{\delta}\left[  \sigma\right]  }\in L^{1} $
for $\mathbf{0}\leq\beta\leq\delta$, which implies $D^{\beta}J_{\delta}\left[
\sigma\right]  \in C_{B}^{\left(  0\right)  }$ for $\mathbf{0}\leq\beta
\leq\delta$, and so $J_{\delta}\left[  \sigma\right]  \in C_{B}^{\left(
\delta\right)  }$.\medskip

\textbf{Part 5} Slight modification of Part 1.
\end{proof}

We now derive some formulas for $D^{\alpha}J_{\delta}\left[  \sigma\right]  $
using the Fourier transform.

\begin{theorem}
\label{Thm_DJ_eq_sumJ}If $\alpha\leq\delta$ and $\sigma,\sigma/y^{\alpha}\in
L_{\overline{\mathcal{O}_{\mathbf{1}}}}^{1}$ then%
\[
D^{\alpha}J_{\delta}\left[  \sigma\right]  =\left(  -1\right)  ^{\left\vert
\alpha\right\vert }\left(  \delta+1\right)  !\sum_{\beta=\mathbf{0}}^{\alpha
}\frac{\left(  -1\right)  ^{\left\vert \beta\right\vert }\tbinom{\alpha}%
{\beta}}{\left(  \delta-\beta+1\right)  !}J_{\delta-\beta}\left[  \frac
{\sigma}{y^{\alpha}}\right]  .
\]

\end{theorem}

\begin{proof}
From Lemma \ref{Lem_FourierJd[u]_formulas},
\begin{align*}
\widehat{J_{\delta}\left[  \sigma\right]  }\left(  \xi\right)   &
=\int_{\mathcal{O}_{\mathbf{1}}}\sigma\left(  y\right)  g_{\delta+1}\left(
y.\xi\right)  dy.\\
\xi^{\alpha}\widehat{J_{\delta}\left[  \sigma\right]  }\left(  \xi\right)   &
=\int_{\mathcal{O}_{\mathbf{1}}}\frac{\sigma\left(  y\right)  }{y^{\alpha}%
}\left(  y.\xi\right)  ^{\alpha}g_{\delta+1}\left(  y.\xi\right)  dy.
\end{align*}

Since $\alpha\leq\delta$, from part 10 of Lemma \ref{Lem_property_g_d},%
\[
\tau^{\alpha}g_{\delta}\left(  \tau\right)  =i^{\left\vert \alpha\right\vert
}\delta!\sum_{\beta=\mathbf{0}}^{\alpha}\left(  -1\right)  ^{\left\vert
\beta\right\vert }\tbinom{\alpha}{\beta}\frac{g_{\delta-\beta}\left(
\tau\right)  }{\left(  \delta-\beta\right)  !},\quad\mathbf{0}\leq\alpha
\leq\delta.
\]

Hence%
\begin{align*}
\xi^{\alpha}\widehat{J_{\delta}\left[  \sigma\right]  }\left(  \xi\right)   &
=\int_{\mathcal{O}_{\mathbf{1}}}\frac{\sigma\left(  y\right)  }{y^{\alpha}%
}\left(  i^{\left\vert \alpha\right\vert }\left(  \delta+1\right)
!\sum_{\beta=\mathbf{0}}^{\alpha}\left(  -1\right)  ^{\left\vert
\beta\right\vert }\tbinom{\alpha}{\beta}\frac{g_{\delta+1-\beta}\left(
y.\xi\right)  }{\left(  \delta+1-\beta\right)  !}\right)  dy\\
& =i^{\left\vert \alpha\right\vert }\left(  \delta+1\right)  !\sum
_{\beta=\mathbf{0}}^{\alpha}\frac{\left(  -1\right)  ^{\left\vert
\beta\right\vert }\tbinom{\alpha}{\beta}}{\left(  \delta-\beta+1\right)
!}\int_{\mathcal{O}_{\mathbf{1}}}\frac{\sigma\left(  y\right)  }{y^{\alpha}%
}g_{\delta+1-\beta}\left(  y.\xi\right)  dy\\
& =i^{\left\vert \alpha\right\vert }\left(  \delta+1\right)  !\sum
_{\beta=\mathbf{0}}^{\alpha}\frac{\left(  -1\right)  ^{\left\vert
\beta\right\vert }\tbinom{\alpha}{\beta}}{\left(  \delta-\beta+1\right)
!}\widehat{J_{\delta-\beta}\left[  \frac{\sigma}{y^{\alpha}}\right]  }\left(
\xi\right)  .
\end{align*}

Finally, since $\xi^{\alpha}\widehat{J_{\delta}\left[  \sigma\right]
}=\left(  -i\right)  ^{\left\vert \alpha\right\vert }\widehat{D^{\alpha
}J_{\delta}\left[  \sigma\right]  }$,%
\[
D^{\alpha}J_{\delta}\left[  \sigma\right]  =\left(  -1\right)  ^{\left\vert
\alpha\right\vert }\left(  \delta+1\right)  !\sum_{\beta=\mathbf{0}}^{\alpha
}\frac{\left(  -1\right)  ^{\left\vert \beta\right\vert }\tbinom{\alpha}%
{\beta}}{\left(  \delta-\beta+1\right)  !}J_{\delta-\beta}\left[  \frac
{\sigma}{y^{\alpha}}\right]  ,
\]

as claimed.
\end{proof}

\begin{remark}
\textbf{CHECK}! From ??%
\[
J_{\delta}\left[  \sigma\right]  \left(  z\right)  =z^{\delta+1}%
\int_{\mathcal{O}_{\mathbf{1}}}e^{-\frac{z}{\mu}\mathbf{1}}\frac{\sigma\left(
\mu\right)  }{\mu^{\delta+2}}d\mu,\quad\delta\geq-2,
\]

so that in one-dimension%
\begin{align*}
DJ_{\delta}\left[  \sigma\right]  \left(  z\right)   & =\left(  \delta
+1\right)  z^{\delta}\int_{\mathcal{O}_{\mathbf{1}}}e^{-\frac{z}{\mu
}\mathbf{1}}\frac{\sigma\left(  \mu\right)  }{\mu^{\delta+2}}d\mu-z^{\delta
+1}\int_{\mathcal{O}_{\mathbf{1}}}e^{-\frac{z}{\mu}\mathbf{1}}\frac
{\sigma\left(  \mu\right)  }{\mu^{\delta+3}}d\mu\\
& =\left(  \delta+1\right)  z^{\delta}\int_{\mathcal{O}_{\mathbf{1}}}%
e^{-\frac{z}{\mu}\mathbf{1}}\frac{\mu^{-1}\sigma\left(  \mu\right)  }%
{\mu^{\delta+1}}d\mu-z^{\delta+1}\int_{\mathcal{O}_{\mathbf{1}}}e^{-\frac
{z}{\mu}\mathbf{1}}\frac{\mu^{-1}\sigma\left(  \mu\right)  }{\mu^{\delta+2}%
}d\mu\\
& =\left(  \delta+1\right)  J_{\delta-1}\left[  \mu^{-1}\sigma\right]  \left(
z\right)  -J_{\delta-1}\left[  \mu^{-1}\sigma\right]  \left(  z\right)  .
\end{align*}

The result of the theorem yields:%
\begin{align*}
DJ_{\delta}\left[  \sigma\right]   & =-\left(  \delta+1\right)  !\sum
_{\beta=0}^{1}\frac{\left(  -1\right)  ^{\left\vert \beta\right\vert }%
}{\left(  \delta-\beta+1\right)  !}J_{\delta-\beta}\left[  y^{-1}\sigma\right]
\\
& =-\left(  \delta+1\right)  !\frac{1}{\left(  \delta+1\right)  !}J_{\delta
}\left[  y^{-1}\sigma\right]  -\left(  \delta+1\right)  !\sum_{\beta=0}%
^{1}\frac{-1}{\left(  \delta-1+1\right)  !}J_{\delta-1}\left[  y^{-1}%
\sigma\right] \\
& =-J_{\delta}\left[  y^{-1}\sigma\right]  +\left(  \delta+1\right)
J_{\delta-1}\left[  y^{-1}\sigma\right]  .
\end{align*}

\end{remark}

??

\begin{remark}%
\begin{align*}
f\left(  \xi\right)   & =\tfrac{\left(  2\pi\right)  ^{d/2}}{\left(
\alpha-1\right)  !}\sum\limits_{\beta\leq\alpha}\tbinom{\alpha}{\beta
}\widehat{J_{\alpha-2}\left[  v_{\beta}\right]  _{-}}\text{ }\widehat{\left\{
D^{\beta}f\right\}  }\\
& =\tfrac{\left(  2\pi\right)  ^{d/2}}{\left(  \alpha-1\right)  !}%
\sum\limits_{\beta\leq\alpha}\tbinom{\alpha}{\beta}\int_{\mathcal{O}%
_{\mathbf{1}}}\widehat{v_{\beta}}\left(  \xi.\zeta\right)  \zeta^{\alpha
-1}e^{-\zeta\mathbf{1}}d\zeta\text{ }\widehat{\left\{  D^{\beta}f\right\}  }\\
& =\tfrac{\left(  2\pi\right)  ^{d/2}}{\left(  \alpha-1\right)  !}%
\sum\limits_{\beta\leq\alpha}\tbinom{\alpha}{\beta}\int_{\mathcal{O}%
_{\mathbf{1}}}\widehat{v_{\beta}}\left(  \xi.\zeta\right)  \zeta^{\alpha
-1}e^{-\zeta\mathbf{1}}d\zeta\text{ }\widehat{\left\{  D^{\beta}f\right\}
}\left(  \xi\right) \\
& =\tfrac{\left(  2\pi\right)  ^{d/2}}{\left(  \alpha-1\right)  !}%
\sum\limits_{\beta\leq\alpha}\tbinom{\alpha}{\beta}\int_{\mathcal{O}%
_{\mathbf{1}}}\widehat{v_{\beta}}\left(  \xi.\zeta\right)  \widehat{\left\{
D^{\beta}f\right\}  }\left(  \xi\right)  \zeta^{\alpha-1}e^{-\zeta\mathbf{1}%
}d\zeta
\end{align*}

?? \textbf{Question} Given $g\in W^{??\alpha}$ can $f\in W^{\alpha}\left(
\mathcal{O}_{\mathbf{1}}\right)  $ be found such that%
\[
\widehat{u}=\tfrac{\left(  2\pi\right)  ^{d/2}}{\left(  \alpha-1\right)
!}\sum\limits_{\beta\leq\alpha}\tbinom{\alpha}{\beta}\int_{\mathcal{O}%
_{\mathbf{1}}}\widehat{v_{\beta}}\left(  \xi.\zeta\right)  \widehat{\left\{
D^{\beta}f\right\}  }\left(  \xi\right)  \zeta^{\alpha-1}e^{-\zeta\mathbf{1}%
}d\zeta?
\]

\textbf{Question} If $u$ is a tensor product function can a tensor product $f
$ be found?
\end{remark}

In the next theorem we use a partition of unity to construct a continuous
convolution extension operator from $W^{n\mathbf{1}}\left(  \Omega\right)  $
to $W^{n\mathbf{1}}\left(  \mathbb{R}^{d}\right)  $ where $\Omega$ satisfies
the rectangle condition.

?? \textbf{ERROR? Must assume that} $v\in L_{0}^{1}\left(  \mathcal{O}%
_{1}\right)  $?

\begin{theorem}
\label{Thm_ExtenOrthantSobolFourier}(\textbf{Extension proof using Fourier
transform}) Let $\Omega$ be a bounded open set that fulfils the
\textbf{uniform rectangle condition} of Definition \ref{Def_OthantCondit}.
Then for each $n=1,2,3,\ldots$ there exists a continuous linear
\textbf{extension operator}%
\[
E_{\Omega}:W^{n\mathbf{1}}\left(  \Omega\right)  \rightarrow W^{n\mathbf{1}%
}\left(  \mathbb{R}^{d}\right)  ,
\]

defined by \ref{X08}, \ref{X09}, \ref{X141}, \ref{X00}, \ref{X151} where $v$
satisfies
\begin{equation}
x^{\beta}v\in L_{\overline{\mathcal{O}_{\mathbf{1}}}}^{1}\text{ for all
multi-integers }-n\mathbf{1}\leq\beta\leq n\mathbf{1}.\label{a030}%
\end{equation}

and for simplicity we have assumed that%
\[
v^{k}\left(  x\right)  =v\left(  \theta^{k}.x\right)  ,\quad x\in
\mathcal{O}^{k}.
\]

We also obtain the continuity bound%
\begin{align*}
&  \left\Vert E_{\Omega}u\right\Vert _{W^{n\mathbf{1}}}\\
&  \leq2^{nd/2}\left(  \int_{\mathcal{O}_{\mathbf{1}}}\frac{\left(
1+x\right)  ^{2n\mathbf{1}}}{x^{n\mathbf{1}}}\left\vert v\left(  x\right)
\right\vert dx\right)  \left(  \sum_{k=1}^{M}\left\Vert \phi_{k}\right\Vert
_{W^{n\mathbf{1},\infty}\left(  \Omega\right)  }\right)  \left\Vert
u\right\Vert _{W^{n\mathbf{1}}\left(  \Omega\right)  },\quad u\in
W^{n\mathbf{1}}\left(  \Omega\right)  .
\end{align*}

\end{theorem}

\begin{proof}
\textbf{Step 1} $C^{\infty}$\textbf{\ density} Since the uniform rectangle
condition is equivalent to the rectangle condition when $\Omega$\ is bounded,
it follows that the segment property holds and so part 4 of Remark
\ref{Rem_SobolevSpace2} implies that $r_{\Omega}C^{\infty}$ is dense in
$W^{n\mathbf{1}}\left(  \Omega\right)  $. Thus it is sufficient to prove this
theorem for functions in $C^{\infty}\left(  \overline{\Omega}\right)  $
because $r_{\Omega}C^{\infty}\subset C^{\infty}\left(  \overline{\Omega
}\right)  \subset W^{n\mathbf{1}}\left(  \Omega\right)  $.\medskip

\textbf{Step 2} \textbf{Localization} Fix $u\in C^{\infty}\left(
\overline{\Omega}\right)  $. Since $\overline{\Omega}$ is compact we can cover
it by a finite number of orthant neighborhoods $\left\{  U_{x^{\left(
k\right)  }}\right\}  _{k=1}^{M}$, choose a $C^{\infty}$ partition of unity
$\left\{  \phi_{k}\right\}  _{k=1}^{M}$ subordinate to this cover and, noting
the definition of the uniform rectangle condition, set%
\[
U_{k}:=U_{x^{\left(  k\right)  }},\text{ }\theta^{k}=\theta_{x^{\left(
k\right)  }},\text{ }h^{k}:=h_{x^{\left(  k\right)  }},\text{ }\mathcal{O}%
^{k}=\theta^{k}.\mathcal{O}_{\mathbf{1}},
\]

where each $\theta^{k}\in\left\{  -1,1\right\}  ^{d}$, and also set%
\begin{equation}
u:=\sum\limits_{k=1}^{M}u_{k},\quad u_{k}:=\left(  r_{\Omega}\phi_{k}\right)
u.\label{X09}%
\end{equation}

Now $u_{k}\in C^{\infty}\left(  \overline{U_{k}\cap\Omega}\right)  $ and in
fact, since $\operatorname*{supp}\phi_{k}\subset U_{k}$ implies
$\operatorname*{dist}\left(  \operatorname*{supp}\phi_{k},\partial
U_{k}\right)  >0$, we have $\operatorname*{supp}u_{k}\subset U_{k}%
\cap\overline{\Omega}$.\medskip

\textbf{Step 3} \textbf{Local extension} We now define extension operators for
the local functions $u_{k}$ introduced in step 2.

Suppose $x\in U_{k}\cap\overline{\Omega}$. Then the uniform rectangle
condition implies $x+\tau.h^{k}.\overline{\widehat{\mathcal{O}^{k}}}%
\subset\Omega$ whenever $\mathbf{0}\leq\tau\leq\mathbf{1}$ and hence $u_{k}\in
C^{\infty}\left(  \overline{\Omega}\right)  $ implies $u_{k}\in C^{\infty
}\left(  x+\tau.h^{k}.\overline{\widehat{\mathcal{O}^{k}}}\right)  $ when
$\mathbf{0}\leq\tau\leq\mathbf{1}$. Thus $u_{k}\in C_{0}^{\infty}\left(
x+\overline{\mathcal{O}^{k}}\right)  $ and so by part 1 of Lemma
\ref{Lem_SmthFuncIntegRepInOrthant} with $\gamma=0$,%
\[
u_{k}\left(  x\right)  =\tfrac{1}{\left(  \left(  n-1\right)  !\right)  ^{d}%
}\sum_{\beta\leq n\mathbf{1}}\tbinom{n\mathbf{1}}{\beta}\int_{x+\mathcal{O}%
^{k}}I_{k,\beta}\left(  x-y\right)  D^{\beta}u_{k}\left(  y\right)  dy,
\]

where
\begin{equation}
I_{k,\beta}:=J_{\left(  n-2\right)  \mathbf{1}}^{\theta^{k}}\left[  \left(
v^{k}\right)  _{\beta}\right]  _{-},\quad v^{k}\in L_{\overline{\mathcal{O}%
^{k}}}^{1}.\label{X00}%
\end{equation}

Now define%
\begin{equation}
u_{k,\beta}\left(  y\right)  :=\left\{
\begin{array}
[c]{ll}%
D^{\beta}u_{k}\left(  y\right)  , & y\in U_{k}\cap\overline{\Omega},\\
0, & otherwise,
\end{array}
\right. \label{X151}%
\end{equation}

and the operator $E_{k}$ by%
\begin{equation}
E_{k}u_{k}\left(  x\right)  :=\tfrac{1}{\left(  \left(  n-1\right)  !\right)
^{d}}\sum_{\beta\leq n\mathbf{1}}\tbinom{n\mathbf{1}}{\beta}\int I_{k,\beta
}\left(  x-y\right)  u_{k,\beta}\left(  y\right)  dy,\quad x\in\mathbb{R}%
^{d}.\label{X141}%
\end{equation}

Thus, if $x\in U_{k}\cap\overline{\Omega}$,%
\[
E_{k}u_{k}\left(  x\right)  =\tfrac{1}{\left(  \left(  n-1\right)  !\right)
^{d}}\sum_{\beta\leq n\mathbf{1}}\tbinom{n\mathbf{1}}{\beta}\int%
_{x+\mathcal{O}^{k}}I_{k,\beta}\left(  x-y\right)  u_{k,\beta}\left(
y\right)  dy=u_{k}\left(  x\right)  ,
\]

confirming that $E_{k}$ is an extension operator from $U_{k}\cap
\overline{\Omega}$ to $\mathbb{R}^{d}$.

Formally, we can write the expression \ref{X141} for this operator in terms of
convolutions as:%
\[
E_{k}u_{k}=\left(  \tfrac{\sqrt{2\pi}}{\left(  n-1\right)  !}\right)  ^{d}%
\sum_{\beta\leq n\mathbf{1}}\tbinom{n\mathbf{1}}{\beta}I_{k,\beta}\ast
u_{k,\beta}.
\]

Clearly $u_{k,\beta}\in L^{2}$ since $\operatorname*{supp}u_{k,\beta
}=\operatorname*{supp}D^{\beta}u_{k}$ is bounded and from part 5 of Lemma
\ref{Lem_SmthFuncIntegRepInOrthant}, $I_{k,\beta}\in L^{1}$ when $\beta\leq
n\mathbf{1}$. Thus, by Young's convolution inequality, $I_{k,\beta}\ast
u_{k,\beta}\in L^{2}$ and from, for example, Exercise 4.7 of Petersen
\cite{Petersen83}, $\left(  I_{k,\beta}\ast u_{k,\beta}\right)  ^{\wedge
}=\widehat{I_{k,\beta}}$ $\widehat{u_{k,\beta}}\in L^{2}$ allowing us to
conclude that%
\[
\widehat{E_{k}u_{k}}=\left(  \tfrac{\sqrt{2\pi}}{\left(  n-1\right)
!}\right)  ^{d}\sum_{\beta\leq n\mathbf{1}}\tbinom{n\mathbf{1}}{\beta
}\widehat{I_{k,\beta}}\text{ }\widehat{u_{k,\beta}}\in L^{2},
\]

where each $\widehat{I_{k,\beta}}\in C_{B}^{\left(  0\right)  }$. This means
$E_{k}u_{k}\in L^{2}$. The Cauchy-Schwartz inequality now yields%
\[
\left\vert \widehat{E_{k}u_{k}}\right\vert \leq\left(  \tfrac{\sqrt{2\pi}%
}{\left(  n-1\right)  !}\right)  ^{d}\left(  \sum_{\beta\leq n\mathbf{1}%
}\tbinom{n\mathbf{1}}{\beta}\left\vert \widehat{I_{k,\beta}}\right\vert
^{2}\right)  ^{1/2}\left(  \sum_{\beta\leq n\mathbf{1}}\tbinom{n\mathbf{1}%
}{\beta}\left\vert \widehat{u_{k,\beta}}\right\vert ^{2}\right)  ^{1/2},
\]

so that from Definition \ref{Def_SobolevSpace},%
\begin{align*}
\left\Vert E_{k}u_{k}\right\Vert _{W^{n\mathbf{1}}}^{2}  & =\sum_{\gamma\leq
n\mathbf{1}}\tbinom{n\mathbf{1}}{\gamma}\int\xi^{2\gamma}\left\vert
\widehat{E_{k}u_{k}}\right\vert ^{2}\\
& \leq\sum_{\gamma\leq n\mathbf{1}}\tbinom{n\mathbf{1}}{\gamma}\int%
\xi^{2\gamma}\left(  \tfrac{\sqrt{2\pi}}{\left(  n-1\right)  !}\right)
^{2d}\left\vert \left(  \sum_{\beta\leq n\mathbf{1}}\tbinom{n\mathbf{1}}%
{\beta}\left\vert \widehat{I_{k,\beta}}\right\vert ^{2}\right)  ^{1/2}\left(
\sum_{\beta\leq n\mathbf{1}}\tbinom{n\mathbf{1}}{\beta}\left\vert
\widehat{u_{k,\beta}}\right\vert ^{2}\right)  ^{1/2}\right\vert ^{2}\\
& =\left(  \tfrac{\sqrt{2\pi}}{\left(  n-1\right)  !}\right)  ^{2d}\int%
\sum_{\gamma\leq n\mathbf{1}}\tbinom{n\mathbf{1}}{\gamma}\xi^{2\gamma}\left(
\sum_{\beta\leq n\mathbf{1}}\tbinom{n\mathbf{1}}{\beta}\left\vert
\widehat{I_{k,\beta}}\right\vert ^{2}\right)  \sum_{\beta\leq n\mathbf{1}%
}\tbinom{n\mathbf{1}}{\beta}\left\vert \widehat{u_{k,\beta}}\right\vert ^{2}\\
& =\left(  \tfrac{\sqrt{2\pi}}{\left(  n-1\right)  !}\right)  ^{2d}\int\left(
\sum_{\beta,\gamma\leq n\mathbf{1}}\tbinom{n\mathbf{1}}{\beta}\tbinom
{n\mathbf{1}}{\gamma}\xi^{2\gamma}\left\vert \widehat{I_{k,\beta}}\right\vert
^{2}\right)  \sum_{\beta\leq n\mathbf{1}}\tbinom{n\mathbf{1}}{\beta}\left\vert
\widehat{u_{k,\beta}}\right\vert ^{2}\\
& =\left(  \tfrac{\sqrt{2\pi}}{\left(  n-1\right)  !}\right)  ^{2d}\int\left(
\sum_{\beta,\gamma\leq n\mathbf{1}}\tbinom{n\mathbf{1}}{\beta}\tbinom
{n\mathbf{1}}{\gamma}\left\vert \xi^{\gamma}\widehat{I_{k,\beta}}\right\vert
^{2}\right)  \sum_{\beta\leq n\mathbf{1}}\tbinom{n\mathbf{1}}{\beta}\left\vert
\widehat{u_{k,\beta}}\right\vert ^{2}\\
& \leq\left(  \tfrac{\sqrt{2\pi}}{\left(  n-1\right)  !}\right)  ^{2d}\left(
\sum_{\beta,\gamma\leq n\mathbf{1}}\tbinom{n\mathbf{1}}{\beta}\tbinom
{n\mathbf{1}}{\gamma}\left\Vert \xi^{\gamma}\widehat{I_{k,\beta}}\right\Vert
_{\infty}^{2}\right)  \int\sum_{\beta\leq n\mathbf{1}}\tbinom{n\mathbf{1}%
}{\beta}\left\vert \widehat{u_{k,\beta}}\right\vert ^{2}\\
& =\left(  \tfrac{\sqrt{2\pi}}{\left(  n-1\right)  !}\right)  ^{2d}\left(
\sum_{\beta,\gamma\leq n\mathbf{1}}\tbinom{n\mathbf{1}}{\beta}\tbinom
{n\mathbf{1}}{\gamma}\left\Vert \xi^{\gamma}\widehat{I_{k,\beta}}\right\Vert
_{\infty}^{2}\right)  \int\sum_{\beta\leq n\mathbf{1}}\tbinom{n\mathbf{1}%
}{\beta}\left\vert \widehat{u_{k,\beta}}\right\vert ^{2}.
\end{align*}

But from \ref{X00},%
\begin{align*}
\left\Vert \xi^{\gamma}\widehat{I_{k,\beta}}\right\Vert _{\infty}%
^{2}=\left\Vert \xi^{\gamma}\widehat{J_{\left(  n-2\right)  \mathbf{1}%
}^{\theta^{k}}\left[  \left(  v^{k}\right)  _{\beta}\right]  _{-}}\left(
\xi\right)  \right\Vert _{\infty}^{2} &  =\left\Vert \xi^{\gamma
}\widehat{J_{\left(  n-2\right)  \mathbf{1}}^{\theta^{k}}\left[  \left(
v^{k}\right)  _{\beta}\right]  }\left(  -\xi\right)  \right\Vert _{\infty}%
^{2}\\
&  =\left\Vert \xi^{\gamma}\widehat{J_{\left(  n-2\right)  \mathbf{1}}%
^{\theta^{k}}\left[  \left(  v^{k}\right)  _{\beta}\right]  }\right\Vert
_{\infty}^{2},
\end{align*}

and since $\widehat{u_{k,\beta}}\in L^{2}$,
\begin{align}
\int\sum_{\beta\leq n\mathbf{1}}\tbinom{n\mathbf{1}}{\beta}\left\vert
\widehat{u_{k,\beta}}\right\vert ^{2}=\sum_{\beta\leq n\mathbf{1}}%
\tbinom{n\mathbf{1}}{\beta}\left\Vert \widehat{u_{k,\beta}}\right\Vert
_{2}^{2}=\sum_{\beta\leq n\mathbf{1}}\tbinom{n\mathbf{1}}{\beta}\left\Vert
u_{k,\beta}\right\Vert _{2}^{2} &  =\sum_{\beta\leq n\mathbf{1}}%
\tbinom{n\mathbf{1}}{\beta}\int_{U_{k}\cap\Omega}\left\vert D^{\beta}%
u_{k}\right\vert ^{2}\nonumber\\
&  =\left\Vert u_{k}\right\Vert _{W^{n\mathbf{1}}\left(  U_{k}\cap
\Omega\right)  }^{2},\label{X26}%
\end{align}

which means that%
\[
\left\Vert E_{k}u_{k}\right\Vert _{W^{n\mathbf{1}}}^{2}\leq\left(
\tfrac{\sqrt{2\pi}}{\left(  n-1\right)  !}\right)  ^{2d}\left(  \sum
_{\beta,\gamma\leq n\mathbf{1}}\tbinom{n\mathbf{1}}{\beta}\tbinom{n\mathbf{1}%
}{\gamma}\left\Vert \xi^{\gamma}\widehat{J_{\left(  n-2\right)  \mathbf{1}%
}^{\theta^{k}}\left[  \left(  v^{k}\right)  _{\beta}\right]  }\right\Vert
_{\infty}^{2}\right)  \left\Vert u_{k}\right\Vert _{W^{n\mathbf{1}}\left(
U_{k}\cap\Omega\right)  }^{2},
\]

and hence%
\begin{equation}
\left\Vert E_{k}u_{k}\right\Vert _{W^{n\mathbf{1}}}\leq\left(  \tfrac
{\sqrt{2\pi}}{\left(  n-1\right)  !}\right)  ^{d}\left(  \sum_{\beta
,\gamma\leq n\mathbf{1}}\tbinom{n\mathbf{1}}{\beta}^{\frac{1}{2}}%
\tbinom{n\mathbf{1}}{\gamma}^{\frac{1}{2}}\left\Vert \xi^{\gamma
}\widehat{J_{\left(  n-2\right)  \mathbf{1}}^{\theta^{k}}\left[  \left(
v^{k}\right)  _{\beta}\right]  }\right\Vert _{\infty}\right)  \left\Vert
u_{k}\right\Vert _{W^{n\mathbf{1}}\left(  U_{k}\cap\Omega\right)
}.\label{X70}%
\end{equation}

Since $v^{k}\left(  x\right)  =v\left(  \theta^{k}.x\right)  $, \ref{X05}
implies%
\[
\left\Vert \xi^{\gamma}\widehat{J_{\left(  n-2\right)  \mathbf{1}}^{\theta
^{k}}\left[  \left(  v^{k}\right)  _{\beta}\right]  }\right\Vert _{\infty
}=\left\Vert \xi^{\gamma}\widehat{J_{\left(  n-2\right)  \mathbf{1}}\left[
v_{\beta}\right]  }\right\Vert _{\infty},
\]

and from part 5 of Lemma \ref{Lem_property_xF[Jb]} and then part 4 of Lemma
\ref{Lem_property_g_d},%
\[
\left\Vert \xi^{\gamma}\widehat{J_{\left(  n-2\right)  \mathbf{1}}^{\theta
^{k}}\left[  \left(  v^{k}\right)  _{\beta}\right]  }\right\Vert _{\infty}%
\leq\left(  \frac{\left(  n-1\right)  !}{\sqrt{2\pi}}\right)  ^{d}\left\Vert
x^{\beta-\gamma}v\right\Vert _{1,\mathcal{O}_{\mathbf{1}}},\quad\mathbf{0}%
\leq\beta,\gamma\leq n\mathbf{1}.
\]

Consequently \ref{X70} becomes%
\[
\left\Vert E_{k}u_{k}\right\Vert _{W^{n\mathbf{1}}}\leq\left(  \sum
_{\beta,\gamma\leq n\mathbf{1}}\tbinom{n\mathbf{1}}{\beta}^{\frac{1}{2}%
}\tbinom{n\mathbf{1}}{\gamma}^{\frac{1}{2}}\left\Vert x^{\beta-\gamma
}v\right\Vert _{1,\mathcal{O}_{\mathbf{1}}}\right)  \left\Vert u_{k}%
\right\Vert _{W^{n\mathbf{1}}\left(  U_{k}\cap\Omega\right)  }.
\]

Regarding the summation term:%
\begin{align}
&  \sum_{\beta,\gamma\leq n\mathbf{1}}\tbinom{n\mathbf{1}}{\beta}^{\frac{1}%
{2}}\tbinom{n\mathbf{1}}{\gamma}^{\frac{1}{2}}\left\Vert x^{\beta-\gamma
}v\right\Vert _{1,\mathcal{O}_{\mathbf{1}}}\nonumber\\
&  =\sum_{\beta,\gamma\leq n\mathbf{1}}\tbinom{n\mathbf{1}}{\beta}^{\frac
{1}{2}}\tbinom{n\mathbf{1}}{\gamma}^{\frac{1}{2}}\int\frac{x^{\beta}%
}{x^{\gamma}}\left\vert v\left(  x\right)  \right\vert dx\nonumber\\
&  =\int\left(  \sum_{\gamma\leq n\mathbf{1}}\tbinom{n\mathbf{1}}{\gamma
}^{\frac{1}{2}}\mathbf{1}^{n\mathbf{1}-\gamma}\frac{1}{x^{\gamma}}\right)
\left(  \sum\limits_{\beta\leq n\mathbf{1}}\tbinom{n\mathbf{1}}{\beta}%
^{\frac{1}{2}}\mathbf{1}^{n\mathbf{1}-\beta}x^{\beta}\right)  \left\vert
v\left(  x\right)  \right\vert dx\nonumber\\
&  =\int\left(  \sum_{\gamma\leq n\mathbf{1}}\tbinom{n\mathbf{1}}{\gamma
}^{\frac{1}{2}}\mathbf{1}^{n\mathbf{1}-\gamma}\left(  \frac{\mathbf{1}}%
{x}\right)  ^{\gamma}\right)  \left(  \sum\limits_{\beta\leq n\mathbf{1}%
}\tbinom{n\mathbf{1}}{\beta}^{\frac{1}{2}}\mathbf{1}^{n\mathbf{1}-\beta
}x^{\beta}\right)  \left\vert v\left(  x\right)  \right\vert dx\nonumber\\
&  \leq\int\left(  \sum_{\gamma\leq n\mathbf{1}}\tbinom{n\mathbf{1}}{\gamma
}\mathbf{1}^{n\mathbf{1}-\gamma}\left(  \frac{\mathbf{1}}{x}\right)  ^{\gamma
}\right)  \left(  \sum\limits_{\beta\leq n\mathbf{1}}\tbinom{n\mathbf{1}%
}{\beta}\mathbf{1}^{n\mathbf{1}-\beta}x^{\beta}\right)  \left\vert v\left(
x\right)  \right\vert dx\nonumber\\
&  =\int_{\mathcal{O}_{\mathbf{1}}}\left(  1+\frac{\mathbf{1}}{x}\right)
^{n\mathbf{1}}\left(  1+x\right)  ^{n\mathbf{1}}\left\vert v\left(  x\right)
\right\vert dx\nonumber\\
&  =\int_{\mathcal{O}_{\mathbf{1}}}\frac{\left(  1+x\right)  ^{2n\mathbf{1}}%
}{x^{n\mathbf{1}}}\left\vert v\left(  x\right)  \right\vert dx,\label{X94}%
\end{align}

and so we have the local continuity estimate%
\begin{equation}
\left\Vert E_{k}u_{k}\right\Vert _{W^{n\mathbf{1}}}\leq\left(  \int%
_{\mathcal{O}_{\mathbf{1}}}\frac{\left(  1+x\right)  ^{2n\mathbf{1}}%
}{x^{n\mathbf{1}}}\left\vert v\left(  x\right)  \right\vert dx\right)
\left\Vert u_{k}\right\Vert _{W^{n\mathbf{1}}\left(  U_{k}\cap\Omega\right)
}.\label{X19}%
\end{equation}

Further%
\begin{align}
\left\Vert u_{k}\right\Vert _{W^{n\mathbf{1}}\left(  U_{k}\cap\Omega\right)
}=\left\Vert \left(  r_{\Omega}\phi_{k}\right)  u\right\Vert _{W^{n\mathbf{1}%
}\left(  U_{k}\cap\Omega\right)  } &  =\left\Vert \phi_{k}u\right\Vert
_{W^{n\mathbf{1}}\left(  U_{k}\cap\Omega\right)  }\nonumber\\
&  \leq\left\Vert \phi_{k}u\right\Vert _{W^{n\mathbf{1}}\left(  \Omega\right)
}\nonumber\\
&  \leq2^{nd/2}\left\Vert \phi_{k}\right\Vert _{W^{n\mathbf{1},\infty}\left(
\Omega\right)  }\left\Vert u\right\Vert _{W^{n\mathbf{1}}\left(
\Omega\right)  },\label{X27}%
\end{align}

where the last estimate is from Theorem \ref{Thm_norm_Cm1xHm1_loc} in the
Appendix. We now introduce our candidate for the extension operator
$E_{\Omega}$, namely%
\begin{equation}
E_{\Omega}u:=\sum_{k=1}^{M}E_{k}u_{k}.\label{X08}%
\end{equation}

Observe that $E_{\Omega}$ is an extension because each $E_{k}$ is an extension
and $u=\sum\limits_{k=1}^{M}u_{k}$. Finally, applying \ref{X27} to the
continuity estimate \ref{X70} for $E_{k}u_{k}$ we get%
\begin{align*}
&  \left\Vert E_{\Omega}u\right\Vert _{W^{n\mathbf{1}}}=\left\Vert \sum
_{k=1}^{M}E_{k}u_{k}\right\Vert _{W^{n\mathbf{1}}}\leq\sum_{k=1}^{M}\left\Vert
E_{k}u_{k}\right\Vert _{W^{n\mathbf{1}}}\leq\\
&  \leq2^{nd/2}\left(  \int_{\mathcal{O}_{\mathbf{1}}}\frac{\left(
1+x\right)  ^{2n\mathbf{1}}}{x^{n\mathbf{1}}}\left\vert v\left(  x\right)
\right\vert dx\right)  \left(  \sum_{k=1}^{M}\left\Vert \phi_{k}\right\Vert
_{W^{n\mathbf{1},\infty}\left(  \Omega\right)  }\right)  \left\Vert
u\right\Vert _{W^{n\mathbf{1}}\left(  \Omega\right)  }.
\end{align*}

\end{proof}

\begin{remark}
\ 

\begin{enumerate}
\item $x^{\beta}v\in L_{\overline{\mathcal{O}_{\mathbf{1}}}}^{1}$ for all
multi-integers $\beta_{+}\leq\left(  n-1\right)  \mathbf{1}$ iff $x^{-\left(
n-1\right)  \mathbf{1}}v,x^{\left(  n-1\right)  \mathbf{1}}v\in L_{\overline
{\mathcal{O}_{\mathbf{1}}}}^{1}$.

\item
\[
\int_{\mathcal{O}_{\mathbf{1}}}\frac{\left(  1+x\right)  ^{2n\mathbf{1}}%
}{x^{n\mathbf{1}}}\left\vert v\left(  x\right)  \right\vert dx=\int%
_{\widehat{\mathcal{O}}_{\mathbf{1}}}\frac{\left(  1+x\right)  ^{2n\mathbf{1}%
}}{x^{n\mathbf{1}}}\left(  \left\vert v\left(  x\right)  \right\vert +\frac
{1}{x^{2\mathbf{1}}}\left\vert v\left(  \frac{\mathbf{1}}{x}\right)
\right\vert \right)  dx.
\]

\item If $x^{\left(  n-1\right)  \mathbf{1}}v\in L_{\overline
{\widehat{\mathcal{O}}_{\mathbf{1}}}}^{1}$ and $v\in L_{\overline
{\widehat{\mathcal{O}}_{\mathbf{1}}}}^{1}$ then $x^{\beta}v\in L_{\overline
{\widehat{\mathcal{O}}_{\mathbf{1}}}}^{1}$ for all $\mathbf{0}\leq\beta
\leq\left(  n-1\right)  \mathbf{1}$. Further, if we extend $v$ to
$\mathcal{O}_{\mathbf{1}}$ as%
\[
v\left(  x\right)  :=\frac{1}{x^{2\mathbf{1}}}\left\vert v\left(
\frac{\mathbf{1}}{x}\right)  \right\vert ,\quad x\in\mathcal{O}_{\mathbf{1}%
}\setminus\widehat{\mathcal{O}}_{\mathbf{1}},
\]

then \ref{a030} holds i.e. $x^{\beta}v\in L_{\overline{\mathcal{O}%
_{\mathbf{1}}}}^{1}$ for all multi-integers $\beta_{+}\leq\left(  n-1\right)
\mathbf{1}$. Here $\widehat{\mathcal{O}}_{\mathbf{1}}:=\left[  0,1\right]
^{d}$.
\end{enumerate}
\end{remark}

\section{A continuous extension $E_{\Omega}^{\alpha}:C^{\left(  \alpha\right)
}\left(  \overline{\Omega}\right)  \rightarrow C_{B}^{\left(  \alpha\right)
}$.\label{Sect_ExtenLocCatoBndCa}}

In this section will use the functional $J_{\gamma}\left[  v\right]  $ to
construct a set of extension operators $E_{\Omega}^{\alpha}:C^{\left(
\alpha\right)  }\left(  \overline{\Omega}\right)  \rightarrow C_{B}^{\left(
\alpha\right)  }$, $\alpha\geq\mathbf{1}$, which are continuous under the
supremum norm. To do this we will constrain the function $v$ to have bounded
support in $\mathcal{O}_{\mathbf{1}}$ and to be a tensor product function
which satisfies \ref{X551}. The main extension result is Theorem
\ref{Thm_ExtenContinFuncs_OrthantProp}.

But first we need:

\begin{lemma}
\label{Lem_PWCinf_J1_b}??? Suppose $\alpha\geq\mathbf{1}$. Then:

\begin{enumerate}
\item There exists a tensor product function $v\left(  x\right)  =v_{1}\left(
x_{1}\right)  \ldots v_{d}\left(  x_{d}\right)  \in L_{0}^{1}\left(
\mathcal{O}_{\mathbf{1}}\right)  $ with support in some rectangle $R\left[
a,b\right]  \subset\mathcal{O}_{\mathbf{1}}$, $a.<b$, such that%
\begin{equation}
\forall i:\int_{a_{i}}^{b_{i}}\frac{v_{i}\left(  t\right)  }{t^{\lambda_{i}}%
}dt=\delta_{\lambda_{i},0},\quad\lambda_{i}=0,1,2,\ldots,\alpha_{i}%
+1,\label{X55}%
\end{equation}

where $\alpha\geq\mathbf{1}$ and $\delta_{\lambda_{i},0}$ denotes the
Kronecker delta.\medskip

If $v$ satisfies condition \ref{X55} then:\medskip

\item
\begin{equation}
D_{i}^{\gamma_{i}}J_{\alpha_{i}-2}\left[  \xi_{i}^{\beta_{i}}v_{i}\right]  \in
PWC^{\infty}\left(  \mathbb{R}^{1}\right)  ,\quad\gamma_{i}=\alpha
_{i}-1,\alpha_{i};\text{ }\beta_{i}=0,1,\ldots,\alpha_{i},\label{X53}%
\end{equation}

?? with jumps at zero.

\item $D^{\gamma}J_{\alpha-2}\left[  \xi^{\beta}v\right]  \in PWC^{\infty
}\left(  \mathbb{R}^{d}\right)  $ when $\beta,\gamma\leq\alpha$.?? When
$\alpha-1<\gamma\leq\alpha$ there are jumps at $\partial\mathcal{O}%
_{\mathbf{1}}$? ??

\item $D^{\gamma}J_{\alpha-2}\left[  \xi^{\beta}v\right]  \in L^{1}$ when
$\beta,\gamma\leq\alpha$, and $\left\Vert D^{\gamma}J_{\alpha-2}\left[
\xi^{\beta}v\right]  \right\Vert _{1}=\left\Vert D^{\gamma}J_{\alpha-2}\left[
\xi^{\beta}v\right]  \right\Vert _{1,\mathcal{O}_{\mathbf{1}}}$. The estimate
of part 4 in Lemma \ref{Lem_property_J} applies.
\end{enumerate}
\end{lemma}

\begin{proof}
\textbf{Part 1} Set $v_{i}\left(  t\right)  =t^{\alpha_{i}+1}u_{i}\left(
t\right)  $. Condition \ref{X55} now becomes:%
\[
\forall i:\int_{a_{i}}^{b_{i}}t^{\mu_{i}}u_{i}\left(  t\right)  dt=\delta
_{\mu_{i},\alpha_{i}+1},\quad\mu_{i}=0,1,2,\ldots,\alpha_{i}+1.
\]

We assume that each $u_{i}$\ is a step function which is constant on each unit
open interval $\left(  k,k+1\right)  \subset\mathbb{R}_{+}^{1}$.

?? Formulate the problem??: assume $g$\ is a step function on the interval
$\left[  1,m+1\right]  $. Assume $m$ steps with $g=a_{j}$ on $\left(
j,j+1\right)  $. Thus%
\begin{align*}
\sum_{j=1}^{m}\int_{j}^{j+1}t^{i-1}a_{j}dt  & =0,\quad i=1,\ldots,m-1.\\
\sum_{j=1}^{m}\int_{j}^{j+1}t^{m-1}a_{j}dt  & =1.
\end{align*}

i.e.%
\begin{align*}
\frac{1}{i}\sum_{j=1}^{m}\left(  \left(  j+1\right)  ^{i}-j^{i}\right)  a_{j}
& =0,\quad i=1,\ldots,m-1.\\
\frac{1}{m}\sum_{j=1}^{m}\left(  \left(  j+1\right)  ^{m}-j^{m}\right)  a_{j}
& =1.
\end{align*}

i.e. we have $m$ equations in the $m$ unknowns $a_{j}$:%
\begin{align*}
\sum_{j=1}^{m}\left(  \left(  j+1\right)  ^{i}-j^{i}\right)  a_{j}  & =0,\quad
i=1,\ldots,m-1.\\
\sum_{j=1}^{m}\left(  \left(  j+1\right)  ^{m}-j^{m}\right)  a_{j}  & =m.
\end{align*}

In matrix form:%
\[%
\begin{pmatrix}
2^{1}-1^{1} & 3^{1}-2^{1} & \cdots & \left(  m+1\right)  ^{1}-m^{1}\\
2^{2}-1^{2} & 3^{2}-2^{2} & \cdots & \left(  m+1\right)  ^{2}-m^{2}\\
\vdots & \vdots & \ddots & \vdots\\
2^{m-1}-1^{m-1} & 3^{m-1}-2^{m-1} & \cdots & \left(  m+1\right)
^{m-1}-m^{m-1}\\
2^{m}-1^{m} & 3^{m}-2^{m} & \cdots & \left(  m+1\right)  ^{m}-m^{m}%
\end{pmatrix}%
\begin{pmatrix}
a_{1}\\
a_{2}\\
a_{3}\\
\vdots\\
a_{m-1}\\
a_{m}%
\end{pmatrix}
=%
\begin{pmatrix}
0\\
0\\
0\\
\vdots\\
0\\
m
\end{pmatrix}
.
\]

Show the determinant is non-zero:%
\begin{align*}
&  \left\vert
\begin{array}
[c]{cccc}%
2^{1}-1^{1} & 3^{1}-2^{1} & \cdots & \left(  m+1\right)  ^{1}-m^{1}\\
2^{2}-1^{2} & 3^{2}-2^{2} & \cdots & \left(  m+1\right)  ^{2}-m^{2}\\
\vdots & \vdots & \ddots & \vdots\\
2^{m-1}-1^{m-1} & 3^{m-1}-2^{m-1} & \cdots & \left(  m+1\right)
^{m-1}-m^{m-1}\\
2^{m}-1^{m} & 3^{m}-2^{m} & \cdots & \left(  m+1\right)  ^{m}-m^{m}%
\end{array}
\right\vert \\
&  =%
\begin{vmatrix}
1 & 1 & 1 & \cdots & 1\\
1^{1} & 2^{1} & 3^{1} & \cdots & \left(  m+1\right)  ^{1}\\
1^{2} & 2^{2} & 3^{2} & \cdots & \left(  m+1\right)  ^{2}\\
\vdots & \vdots & \vdots & \ddots & \vdots\\
1^{m-1} & 2^{m-1} & 3^{m-1} & \cdots & \left(  m+1\right)  ^{m-1}\\
1^{m} & 2^{m} & 3^{m} & \cdots & \left(  m+1\right)  ^{m}%
\end{vmatrix}
\\
&  =Vandemonde\left(  1,2,3,\ldots,m,m+1\right) \\
&  =1!2!3!\ldots m!\\
&  \neq0.
\end{align*}

We conclude that the $a_{i}$ always exist and are unique and so $g$ exists and
hence a $v_{i}$ exists $\forall i$.\medskip

\textbf{Part 2} If $v$ is a tensor product:%
\[
D^{\gamma}J_{\alpha-2}\left[  v_{\beta}\right]  =\left(  -1\right)
^{\left\vert \beta\right\vert }D_{1}^{\gamma_{1}}J_{\alpha_{1}-2}\left[
\xi^{\beta_{1}}v_{1}\right]  D_{2}^{\gamma_{2}}J_{\alpha_{2}-2}\left[
\xi^{\beta_{2}}v_{2}\right]  \ldots D_{d}^{\gamma_{d}}J_{\alpha_{d}-2}\left[
\xi^{\beta_{d}}v_{d}\right]  .
\]

Then we want $v_{i}\in L_{0}^{1}\left(  0,\infty\right)  $ such that $\forall
i$ \ref{X53} holds, and since by Lemma \ref{Lem_property_J}, $D_{i}%
^{\alpha_{i}-1}J_{\alpha_{i}-2}\left[  \xi_{i}^{\beta_{i}}v_{i}\right]  \in
PWC^{\infty}\left(  \mathbb{R}^{1}\right)  $, this can be ensured by requiring
that $\forall i$,%
\begin{equation}
D_{i}^{\alpha_{i}}J_{\alpha_{i}-2}\left[  \xi_{i}^{\beta_{i}}v_{i}\right]
\left(  0^{+}\right)  =0,\quad\beta_{i}=0,1,\ldots,\alpha_{i}.\label{X57}%
\end{equation}

But from \ref{X40}, $\forall\gamma\geq\mathbf{0}$,
\begin{equation}
D^{\gamma}J_{\delta}\left[  \sigma\right]  \left(  z\right)  =\int_{a}%
^{b}p_{\gamma}^{\delta}\left(  \frac{z}{\mu}\right)  \frac{\sigma\left(
\mu\right)  }{\mu^{\gamma+1}}e^{-\frac{z}{\mu}\mathbf{1}}d\mu,\quad
z\in\overline{\mathcal{O}_{\mathbf{1}}},\label{X58}%
\end{equation}

so that $\forall i$,%
\begin{equation}
D_{i}^{\alpha_{i}}J_{\alpha_{i}-2}\left[  \xi_{i}^{\beta_{i}}v_{i}\left(
\xi_{i}\right)  \right]  \left(  s\right)  =\int_{a_{i}}^{b_{i}}p_{\alpha_{i}%
}^{\alpha_{i}-2}\left(  \frac{s}{t}\right)  \frac{t^{\beta_{i}}v_{i}\left(
t\right)  }{t^{\alpha_{i}+1}}e^{-\frac{s}{t}}dt,\text{\quad}s\geq0,\text{
}\beta_{i}\leq\alpha_{i},\label{X54}%
\end{equation}

and thus condition \ref{X57} becomes: for all $i$,
\[
p_{\alpha_{i}}^{\alpha_{i}-2}\left(  0\right)  \int_{a_{i}}^{b_{i}}%
\frac{t^{\beta_{i}}v_{i}\left(  t\right)  }{t^{\alpha_{i}+1}}dt=0,\quad
\beta_{i}\leq\alpha_{i},
\]

or%
\[
p_{\alpha_{i}}^{\alpha_{i}-2}\left(  0\right)  \int_{a_{i}}^{b_{i}}\frac
{v_{i}\left(  t\right)  }{t^{\alpha_{i}-\beta_{i}+1}}dt=0,\quad\beta_{i}%
\leq\alpha_{i},
\]

or%
\[
p_{\alpha_{i}}^{\alpha_{i}-2}\left(  0\right)  \int_{a_{i}}^{b_{i}}\frac
{v_{i}\left(  t\right)  }{t^{k}}dt=0,\quad k=1,\ldots,\alpha_{i}+1.
\]

But from \ref{X43}: $p_{n}^{m}\left(  0\right)  =\left(  -1\right)
^{n-1-m}\frac{\left(  m+1\right)  !}{\left(  n-1-m\right)  !}$ when $n\geq
m+1$, so that in particular $p_{n}^{n-2}\left(  0\right)  =-\left(
n-1\right)  !$ and our condition now becomes
\[
\forall i:\int_{a_{i}}^{b_{i}}\frac{v_{i}\left(  t\right)  }{t^{k}}dt=0,\quad
k=1,2,\ldots,\alpha_{i}+1.
\]

We also require that%
\[
\forall i:\int_{a_{i}}^{b_{i}}v_{i}\left(  t\right)  dt=1,
\]

so that the conditions can be written:
\[
\forall i:\int_{a_{i}}^{b_{i}}\frac{v_{i}\left(  t\right)  }{t^{\lambda_{i}}%
}dt=\delta_{\lambda_{i},0},\quad\lambda_{i}\leq\alpha_{i}+1,
\]

which is \ref{X55}.\medskip

\textbf{Part 3} From part 6 of Lemma \ref{Lem_property_J}, $D^{\gamma
}J_{\alpha-2}\left[  \xi^{\beta}v\right]  \in PWC^{\infty}\left(
\mathbb{R}^{d}\right)  $ when $\beta,\gamma\leq\alpha-1$, and so the result of
part 2 directly implies part 3.\medskip

\textbf{Part 4} Since $\alpha\geq\mathbf{1}$, part 5 of Lemma
\ref{Lem_SmthFuncIntegRepInOrthant} implies $D^{\gamma}J_{\alpha-2}\left[
\xi^{\beta}v\right]  \in L^{1}\left(  \mathcal{O}_{\mathbf{1}}\right)  $ for
all $\beta,\gamma$. From part 3, $D^{\gamma}J_{\alpha-2}\left[  \xi^{\beta
}v\right]  \in PWC^{\infty}\left(  \mathbb{R}^{d}\right)  $ when $\beta
,\gamma\leq\alpha$. Thus, since $\operatorname*{supp}D^{\gamma}J_{\alpha
-2}\left[  \xi^{\beta}v\right]  \subseteq\overline{\mathcal{O}_{\mathbf{1}}}$,
we can conclude that $D^{\gamma}J_{\alpha-2}\left[  \xi^{\beta}v\right]  \in
L^{1} $ when $\beta,\gamma\leq\alpha$.
\end{proof}

This is the arbitrary rectangle version of Lemma \ref{Lem_PWCinf_J1_b}.

\begin{lemma}
\label{Lem_PWCinf_Jtheta_b}Suppose $\alpha\geq\mathbf{1}$. Then:

\begin{enumerate}
\item There exists a tensor product function $v\left(  x\right)  =v_{1}\left(
x_{1}\right)  \ldots v_{d}\left(  x_{d}\right)  \in L_{0}^{1}\left(
\mathcal{O}_{\theta}\right)  $ with bounded support in the rectangle $R\left[
c,d\right]  \subset\mathcal{O}_{\theta},$ $c\centerdot<d$, which satisfies%
\begin{equation}
\forall i:\int_{c_{i}}^{d_{i}}\frac{v_{i}\left(  t\right)  }{t^{\lambda_{i}}%
}dt=\delta_{\lambda_{i},0},\quad\lambda_{i}=0,1,2,\ldots,\alpha_{i}%
+1.\label{X551}%
\end{equation}

Here $\delta_{\lambda_{i},0}$ denotes the Kronecker delta.

\item If $v$ satisfies condition \ref{X551} then:%
\[
D_{i}^{\gamma_{i}}J_{\alpha_{i}-2}^{\theta_{i}}\left[  \xi_{i}^{\beta_{i}%
}v_{i}\right]  \in PWC^{\infty}\left(  \mathbb{R}^{1}\right)  ,\quad\gamma
_{i}=\alpha_{i}-1,\alpha_{i};\text{ }\beta_{i}=0,1,\ldots,\alpha_{i},
\]
\medskip

If $v$ satisfies condition \ref{X551} then:\medskip

\item $D_{i}^{\gamma_{i}}J_{\alpha_{i}-2}^{\theta_{i}}\left[  \xi_{i}%
^{\beta_{i}}v_{i}\right]  \in PWC^{\infty}\left(  \mathbb{R}^{1}\right)
,\quad\gamma_{i}=\alpha_{i}-1,\alpha_{i};$ $\beta_{i}=0,1,\ldots,\alpha_{i}$,
?? with a jump at zero.

\item $L^{1}$ \textbf{estimate}: $D^{\gamma}J_{\alpha-2}^{\theta}\left[
\xi^{\beta}v\right]  \in L^{1}$ when $\beta,\gamma\leq\alpha$, and
\[
\left\Vert D^{\gamma}J_{\alpha-2}^{\theta}\left[  \xi^{\beta}v\right]
\right\Vert _{1}\leq??\left\Vert \frac{x^{\beta}v}{x^{\gamma}}\right\Vert
_{1,\mathcal{O}_{\theta}}\int_{\mathcal{O}_{\mathbf{1}}}q_{\gamma}^{\alpha
-2}\left(  y\right)  e^{-y\mathbf{1}}dy.
\]

\end{enumerate}
\end{lemma}

\begin{proof}
?? \textbf{Part 1} Concerning rectangles see Definition \ref{Def_topol_on_Rd}.
By definition, $v\in L_{0}^{1}\left(  \mathcal{O}_{\theta}\right)  $ iff
$v\left(  \theta.\right)  \in L_{0}^{1}\left(  \mathcal{O}_{\mathbf{1}%
}\right)  $ and by part 1 of Lemma \ref{Lem_PWCinf_J1_b} there exists a closed
rectangle $R\left[  a,b\right]  \subset\mathcal{O}_{\mathbf{1}}$ such that%
\[
\forall i:\int_{a_{i}}^{b_{i}}\frac{v_{i}\left(  \theta_{i}t\right)
}{t^{\lambda_{i}}}dt=\delta_{\lambda_{i},0},\quad\lambda_{i}=0,1,2,\ldots
,\alpha_{i}+1,
\]

which can be written: for all $i$,%
\[%
\begin{array}
[c]{l}%
\int_{a_{i}}^{b_{i}}\frac{v_{i}\left(  \theta_{i}t\right)  }{t^{\lambda_{i}}%
}dt=\delta_{\lambda_{i},0},\quad\lambda_{i}=1,2,\ldots,\alpha_{i}+1,\\
\int_{a_{i}}^{b_{i}}v_{i}\left(  \theta_{i}t\right)  dt=1.
\end{array}
\]

Applying the change of variables $s=\theta_{i}t$, $t=\theta_{i}s$,
$dt=\theta_{i}ds$ we get%
\[%
\begin{array}
[c]{l}%
\int_{\theta_{i}a_{i}}^{\theta_{i}b_{i}}\frac{v_{i}\left(  s\right)
}{s^{\lambda_{i}}}ds=0,\quad\lambda_{i}=1,2,\ldots,\alpha_{i}+1,\\
\theta_{i}\int_{\theta_{i}a_{i}}^{\theta_{i}b_{i}}v_{i}\left(  s\right)  ds=1,
\end{array}
\]

which can be written%
\[%
\begin{array}
[c]{l}%
\int_{\min\left\{  \theta_{i}a_{i},\theta_{i}b_{i}\right\}  }^{\max\left\{
\theta_{i}a_{i},\theta_{i}b_{i}\right\}  }\frac{v_{i}\left(  s\right)
}{s^{\lambda_{i}}}ds=0,\quad\lambda_{i}=1,2,\ldots,\alpha_{i}+1,\\
\int_{\min\left\{  \theta_{i}a_{i},\theta_{i}b_{i}\right\}  }^{\max\left\{
\theta_{i}a_{i},\theta_{i}b_{i}\right\}  }v_{i}\left(  s\right)  ds=1.
\end{array}
\]

Thus there exists a rectangle $R\left[  c,d\right]  \subset\mathcal{O}%
_{\theta}$, $c=\min\left\{  \theta.a,\theta.b\right\}  $, $d=\max\left\{
\theta.a,\theta.b\right\}  $ such that%
\[
\forall i:\int_{c_{i}}^{d_{i}}\frac{v_{i}\left(  t\right)  }{t^{\lambda_{i}}%
}dt=\delta_{\lambda_{i},0},\quad\lambda_{i}=0,1,2,\ldots,\alpha_{i}+1.
\]
\medskip

\textbf{Part 2} ??\medskip

\textbf{Part 3} ??\medskip

\textbf{Part 4} ??%
\[
\left\Vert D^{\gamma}J_{\alpha-2}^{\theta}\left[  \xi^{\beta}v\right]
\right\Vert _{1}=\left\Vert D^{\gamma}J_{\alpha-2}^{\theta}\left[  \xi^{\beta
}v\right]  \right\Vert _{1,\mathcal{O}_{\theta}}\leq\left\Vert \frac{x^{\beta
}v}{x^{\gamma}}\right\Vert _{1,\mathcal{O}_{\theta}}??\int_{\mathcal{O}%
_{\mathbf{1}}}q_{\gamma}^{\alpha-2}\left(  y\right)  e^{-y\mathbf{1}}dy\text{
}??,
\]

where the estimate of part 5 in Lemma \ref{Lem_SmthFuncIntegRepInOrthant} is used.
\end{proof}

This is the main extension theorem of this section.

\begin{theorem}
\label{Thm_ExtenContinFuncs_OrthantProp}Assume $\alpha\geq\mathbf{1}$ and
$\Omega$ has the uniform rectangle property. We construct a continuous
extension operator $E_{\Omega}^{\alpha}:C_{B}^{\left(  \alpha\right)  }\left(
\overline{\Omega}\right)  \rightarrow C_{B}^{\left(  \alpha\right)  }\left(
\mathbb{R}^{d}\right)  $. Assume $v\in L_{0}^{1}\left(  \mathcal{O}%
_{\mathbf{1}}\right)  $ is a tensor product function satisfying condition
\ref{X551}. Assume (for simplicity) that for some fixed $v\in L_{0}^{1}\left(
\mathcal{O}_{\mathbf{1}}\right)  $ and all $k$,
\begin{equation}
v^{k}\left(  x\right)  :=v\left(  \theta^{k}.x\right)  ,\quad x\in
\mathcal{O}_{\theta}\text{.}\label{X68}%
\end{equation}

Then we have the continuity estimate%
\[
\left\Vert E_{\Omega}^{\alpha}u\right\Vert _{\alpha,\infty;\Omega}\leq
k_{\alpha,v}\left(  \sum_{k=1}^{M}\left\Vert \phi_{k}\right\Vert
_{\alpha,\infty;\Omega}\right)  \left\Vert u\right\Vert _{\alpha,\infty
;\Omega},
\]

where $k_{v,\alpha}$ is given by \ref{X34} and $\left\{  \phi_{k}\right\}  $
is a special partition of unity.
\end{theorem}

\begin{proof}
\textbf{Step 1}: \textbf{Localization} Fix $u\in C^{\left(  \alpha\right)
}\left(  \overline{\Omega}\right)  $. Since $\overline{\Omega}$ is compact we
can cover it by a finite number of orthant neighborhoods $\left\{
U_{x^{\left(  k\right)  }}\right\}  _{k=1}^{M}$, choose a $C^{\infty}$
\textbf{partition of unity} $\left\{  \phi_{k}\right\}  _{k=1}^{M}$
subordinate to this cover and, noting the definition of the uniform rectangle
condition, set%
\[
U_{k}:=U_{x^{\left(  k\right)  }},\text{ }\theta^{k}=\theta_{x^{\left(
k\right)  }},\text{ }h^{k}:=h_{x^{\left(  k\right)  }},\text{ }\mathcal{O}%
^{k}=\theta^{k}.\mathcal{O}_{\mathbf{1}},
\]

where each $\theta^{k}\in\left\{  -1,1\right\}  ^{d}$, and also set%
\[
u:=\sum\limits_{k=1}^{M}u_{k},\quad u_{k}:=\left(  r_{\Omega}\phi_{k}\right)
u.
\]

Now $u_{k}\in C^{\left(  \alpha\right)  }\left(  \overline{U_{k}\cap\Omega
}\right)  $ and in fact, since $\operatorname*{supp}\phi_{k}\subset U_{k}$
implies $\operatorname*{dist}\left(  \operatorname*{supp}\phi_{k},\partial
U_{k}\right)  >0$, we have $\operatorname*{supp}u_{k}\subset U_{k}%
\cap\overline{\Omega}$.\medskip

\textbf{Step 2}: \textbf{Local extension} We now define extension operators
for the local functions $u_{k}$ introduced in step 2.

Suppose $x\in U_{k}\cap\overline{\Omega}$. Then the uniform rectangle
condition implies $x+\tau.h^{k}.\overline{\widehat{\mathcal{O}^{k}}}%
\subset\Omega$ whenever $\mathbf{0}\leq\tau\leq\mathbf{1}$ and hence $u_{k}\in
C^{\left(  \alpha\right)  }\left(  \overline{\Omega}\right)  $ implies
$u_{k}\in C^{\left(  \alpha\right)  }\left(  x+\tau.h^{k}.\overline
{\widehat{\mathcal{O}^{k}}}\right)  $ when $\mathbf{0}\leq\tau\leq\mathbf{1}$.
Thus $u_{k}\in C_{0}^{\left(  \alpha\right)  }\left(  x+\overline
{\mathcal{O}^{k}}\right)  $ and so by part 1 of Lemma
\ref{Lem_SmthFuncIntegRepInOrthant} with $\gamma=0$,%
\[
u_{k}\left(  x\right)  =\tfrac{1}{\left(  \alpha-1\right)  !}\sum_{\beta
\leq\alpha}\tbinom{\alpha}{\beta}\int_{x+\mathcal{O}^{k}}I_{k,\beta}\left(
x-y\right)  D^{\beta}u_{k}\left(  y\right)  dy,
\]

where
\begin{equation}
I_{k,\beta}:=J_{\alpha-2}^{\theta^{k}}\left[  \left(  v^{k}\right)  _{\beta
}\right]  _{-},\quad v^{k}\in L_{0}^{1}\left(  \mathcal{O}^{k}\right)
,\label{X01}%
\end{equation}

and $v^{k}$ has property \ref{X551}. Now define%
\begin{equation}
u_{k,\beta}\left(  y\right)  :=\left\{
\begin{array}
[c]{ll}%
D^{\beta}u_{k}\left(  y\right)  , & y\in U_{k}\cap\overline{\Omega},\\
0, & otherwise,
\end{array}
\right. \label{X15}%
\end{equation}

and the operator $E_{k}$ by%
\begin{equation}
E_{k}u_{k}\left(  x\right)  :=\tfrac{1}{\left(  \alpha-1\right)  !}\sum
_{\beta\leq\alpha}\tbinom{\alpha}{\beta}\int I_{k,\beta}\left(  x-y\right)
u_{k,\beta}\left(  y\right)  dy,\quad x\in\mathbb{R}^{d}.\label{X14}%
\end{equation}

Thus, if $x\in U_{k}\cap\overline{\Omega}$,%
\[
E_{k}u_{k}\left(  x\right)  =\tfrac{1}{\left(  \alpha-1\right)  !}\sum
_{\beta\leq\alpha}\tbinom{\alpha}{\beta}\int_{x+\mathcal{O}^{k}}I_{k,\beta
}\left(  x-y\right)  u_{k,\beta}\left(  y\right)  dy=u_{k}\left(  x\right)  ,
\]

confirming that $E_{k}$ is an extension operator from $U_{k}\cap
\overline{\Omega}$ to $\mathbb{R}^{d}$.

We now introduce our candidate for the extension operator $E_{\Omega}$ from
$\Omega$:%
\begin{equation}
E_{\Omega}^{\alpha}u:=\sum_{k=1}^{M}E_{k}u_{k}.\label{X22}%
\end{equation}

Observe that $E_{\Omega}^{\alpha}$ is an extension because each $E_{k}$ is an
extension and $u=\sum\limits_{k=1}^{M}u_{k}$.

Since $u_{k,\beta}\in C_{0}^{\left(  \alpha-\beta\right)  }\left(
x+\overline{\mathcal{O}^{k}}\right)  \in\mathcal{E}^{\prime}$ and $I_{k,\beta
}\in S^{\prime}$ we can use the $S^{\prime}\ast\mathcal{E}^{\prime}$
convolution (see, for example, Subsection 2:9.6 of Vladimirov
\cite{Vladimirov} OR Corollary 2:6.9 and Theorem 2:6.10 of Petersen
\cite{Petersen83}) to write%
\[
E_{k}u_{k}\left(  x\right)  =\tfrac{\left(  2\pi\right)  ^{d/2}}{\left(
\alpha-1\right)  !}\sum_{\beta\leq\alpha}\tbinom{\alpha}{\beta}I_{k,\beta}\ast
u_{k,\beta},
\]

and%
\begin{equation}
D^{\gamma}E_{k}u_{k}=\tfrac{\left(  2\pi\right)  ^{d/2}}{\left(
\alpha-1\right)  !}\sum_{\beta\leq\alpha}\tbinom{\alpha}{\beta}\left(
D^{\gamma}I_{k,\beta}\right)  \ast u_{k,\beta},\quad\forall\gamma.\label{X25}%
\end{equation}

Noting \ref{X01}, part 4 of Lemma \ref{Lem_PWCinf_Jtheta_b} implies
$D^{\gamma}I_{k,\beta}\in L^{1}$ when $\beta,\gamma\leq n\mathbf{1}$. Also,
$u_{k}\in C_{0}^{\left(  \alpha\right)  }\left(  x+\overline{\mathcal{O}^{k}%
}\right)  $ implies $u_{k}\in L^{\infty}$. Thus, by Young's convolution
inequality, $I_{k,\beta}\ast u_{k,\beta}\in L^{\infty}$ and so $D^{\gamma
}E_{k}u_{k}\in L^{\infty}$. Further%
\begin{align*}
\left\Vert D^{\gamma}E_{k}u_{k}\right\Vert _{\infty}  & \leq\tfrac{\left(
2\pi\right)  ^{d/2}}{\left(  \alpha-1\right)  !}\sum_{\beta\leq\alpha}%
\tbinom{\alpha}{\beta}\left\Vert \left(  D^{\gamma}I_{k,\beta}\right)  \ast
u_{k,\beta}\right\Vert _{\infty}\\
& \leq\tfrac{1}{\left(  \alpha-1\right)  !}\sum_{\beta\leq\alpha}%
\tbinom{\alpha}{\beta}\left\Vert D^{\gamma}I_{k,\beta}\right\Vert
_{1}\left\Vert u_{k,\beta}\right\Vert _{\infty,\Omega}\\
& =\tfrac{1}{\left(  \alpha-1\right)  !}\sum_{\beta\leq\alpha}\tbinom{\alpha
}{\beta}\left\Vert D^{\gamma}I_{k,\beta}\right\Vert _{1}\left\Vert u_{k,\beta
}\right\Vert _{\infty,\Omega}\\
& =\tfrac{1}{\left(  \alpha-1\right)  !}\sum_{\beta\leq\alpha}\tbinom{\alpha
}{\beta}\left\Vert D^{\gamma}I_{k,\beta}\right\Vert _{1}\left\Vert D^{\beta
}u_{k}\right\Vert _{\infty,\Omega},
\end{align*}

so that%
\begin{align*}
\left\Vert E_{k}u_{k}\right\Vert _{\alpha;\infty;\Omega}  & =\sum_{\gamma
\leq\alpha}\tbinom{\alpha}{\gamma}\left\Vert D^{\gamma}E_{k}u_{k}\right\Vert
_{\infty}\\
& \leq\tfrac{1}{\left(  \alpha-1\right)  !}\sum_{\gamma\leq\alpha}%
\tbinom{\alpha}{\gamma}\sum_{\beta\leq\alpha}\tbinom{\alpha}{\beta}\left\Vert
D^{\gamma}I_{k,\beta}\right\Vert _{1}\left\Vert D^{\beta}u_{k}\right\Vert
_{\infty,\Omega}\\
& \leq\tfrac{1}{\left(  \alpha-1\right)  !}\sum_{\gamma\leq\alpha}%
\tbinom{\alpha}{\gamma}\max_{\beta\leq\alpha}\left\Vert D^{\beta}%
u_{k}\right\Vert _{\infty,\Omega}\sum_{\beta\leq\alpha}\tbinom{\alpha}{\beta
}\left\Vert D^{\gamma}I_{k,\beta}\right\Vert _{1}.
\end{align*}

Noting \ref{X01}, part 4 of Lemma \ref{Lem_property_J} and\ \ref{X05} imply ??%
\begin{align*}
\left\Vert D^{\gamma}I_{k,\beta}\right\Vert _{1}=\left\Vert D^{\gamma
}J_{\alpha-2}^{\theta^{k}}\left[  \left(  v^{k}\right)  _{\beta}\right]
_{-}\right\Vert _{1} &  =\left\Vert D^{\gamma}J_{\alpha-2}^{\theta^{k}}\left[
\left(  v^{k}\right)  _{\beta}\right]  \right\Vert _{1}\\
&  =\left\Vert D^{\gamma}J_{\alpha-2}\left[  v_{\beta}\right]  \right\Vert
_{1}\\
&  \leq2^{\left\vert \gamma\right\vert }\left(  \alpha-1\right)  !\left\Vert
\frac{x^{\beta}v}{x^{\gamma}}\right\Vert _{1,\mathcal{O}^{k}},
\end{align*}

so that%
\begin{align*}
\sum_{\beta\leq\alpha}\tbinom{\alpha}{\beta}\left\Vert D^{\gamma}I_{k,\beta
}\right\Vert _{1}  & \leq2^{\left\vert \gamma\right\vert }\left(
\alpha-1\right)  !\sum_{\beta\leq\alpha}\tbinom{\alpha}{\beta}\left\Vert
\frac{x^{\beta}}{x^{\gamma}}v\right\Vert _{1,\mathcal{O}_{\mathbf{1}}}\\
& =2^{\left\vert \gamma\right\vert }\left(  \alpha-1\right)  !\sum_{\beta
\leq\alpha}\int_{\mathcal{O}_{\mathbf{1}}}\tbinom{\alpha}{\beta}\frac
{x^{\beta}}{x^{\gamma}}\left\vert v\left(  x\right)  \right\vert dx\\
& =2^{\left\vert \gamma\right\vert }\left(  \alpha-1\right)  !\int%
_{\mathcal{O}_{\mathbf{1}}}\frac{\sum\limits_{\beta\leq\alpha}\tbinom{\alpha
}{\beta}x^{\beta}\mathbf{1}^{\alpha-\beta}}{\left\vert x^{\gamma}\right\vert
}\left\vert v\left(  x\right)  \right\vert dx\\
& =2^{\left\vert \gamma\right\vert }\left(  \alpha-1\right)  !\left(
\int_{\mathcal{O}_{\mathbf{1}}}\frac{\left(  1+x\right)  ^{\alpha}}{x^{\gamma
}}\left\vert v\left(  x\right)  \right\vert dx\right)  ,
\end{align*}

and hence%
\begin{align*}
\left\Vert E_{k}u_{k}\right\Vert _{\alpha;\infty;\Omega}  & \leq\tfrac
{1}{\left(  \alpha-1\right)  !}\sum_{\gamma\leq\alpha}\tbinom{\alpha}{\gamma
}\max_{\beta\leq\alpha}\left\Vert D^{\beta}u_{k}\right\Vert _{\infty,\Omega
}2^{\left\vert \gamma\right\vert }\left(  \alpha-1\right)  !\int%
_{\mathcal{O}_{\mathbf{1}}}\frac{\left(  1+x\right)  ^{\alpha}}{x^{\gamma}%
}\left\vert v\left(  x\right)  \right\vert dx\\
& =\sum_{\gamma\leq\alpha}\tbinom{\alpha}{\gamma}\max_{\beta\leq\alpha
}\left\Vert D^{\beta}u_{k}\right\Vert _{\infty,\Omega}2^{\left\vert
\gamma\right\vert }\int_{\mathcal{O}_{\mathbf{1}}}\frac{\left(  1+x\right)
^{\alpha}}{x^{\gamma}}\left\vert v\left(  x\right)  \right\vert dx\\
& =\max_{\beta\leq\alpha}\left\Vert D^{\beta}u_{k}\right\Vert _{\infty,\Omega
}\sum_{\gamma\leq\alpha}\tbinom{\alpha}{\gamma}2^{\left\vert \gamma\right\vert
}\int_{\mathcal{O}_{\mathbf{1}}}\frac{\left(  1+x\right)  ^{\alpha}}%
{x^{\gamma}}\left\vert v\left(  x\right)  \right\vert dx\\
& =\max_{\beta\leq\alpha}\left\Vert D^{\beta}u_{k}\right\Vert _{\infty,\Omega
}\int_{\mathcal{O}_{\mathbf{1}}}\sum_{\gamma\leq\alpha}\tbinom{\alpha}{\gamma
}\left(  \frac{2}{x}\right)  ^{\gamma}\left(  1+x\right)  ^{\alpha}\left\vert
v\left(  x\right)  \right\vert dx\\
& =\max_{\beta\leq\alpha}\left\Vert D^{\beta}u_{k}\right\Vert _{\infty,\Omega
}\int_{\mathcal{O}_{\mathbf{1}}}\left(  1+\frac{2}{x}\right)  ^{\alpha}\left(
1+x\right)  ^{\alpha}\left\vert v\left(  x\right)  \right\vert dx\\
& =\max_{\beta\leq\alpha}\left\Vert D^{\beta}u_{k}\right\Vert _{\infty,\Omega
}\int_{\mathcal{O}_{\mathbf{1}}}\frac{\left(  1+x\right)  ^{\alpha}\left(
2+x\right)  ^{\alpha}}{x^{\alpha}}\left\vert v\left(  x\right)  \right\vert
dx.
\end{align*}

Thus, from \ref{X22},
\[
\left\Vert E_{\Omega}^{\alpha}u\right\Vert _{\alpha;\infty;\Omega}\leq
\sum_{k=1}^{M}\left\Vert E_{k}u_{k}\right\Vert _{\alpha;\infty;\Omega}%
\leq\left(  \sum_{k=1}^{M}\max_{\beta\leq\alpha}\left\Vert D^{\beta}%
u_{k}\right\Vert _{\infty,\Omega}\right)  \int_{\mathcal{O}_{\mathbf{1}}}%
\frac{\left(  1+x\right)  ^{\alpha}\left(  2+x\right)  ^{\alpha}}{x^{\alpha}%
}\left\vert v\left(  x\right)  \right\vert dx.
\]

But%
\begin{align*}
\max_{\beta\leq\alpha}\left\Vert D^{\beta}u_{k}\right\Vert _{\infty,\Omega}  &
=\max_{\beta\leq\alpha}\left\Vert D^{\beta}\left(  \phi_{k}u\right)
\right\Vert _{\infty,\Omega}=\max_{\beta\leq\alpha}\left\Vert \sum_{\gamma
\leq\beta}\tbinom{\beta}{\gamma}D^{\beta-\gamma}\phi_{k}D^{\gamma}u\right\Vert
_{\infty,\Omega}\leq\\
& \leq\max_{\beta\leq\alpha}\sum_{\gamma\leq\beta}\tbinom{\beta}{\gamma
}\left\Vert D^{\beta-\gamma}\phi_{k}\right\Vert _{\infty,\Omega}\left\Vert
D^{\gamma}u\right\Vert _{\infty,\Omega}\\
& \leq\max_{\beta\leq\alpha}\left(  \max_{\gamma\leq\beta}\left\Vert
D^{\gamma}\phi_{k}\right\Vert _{\infty,\Omega}\right)  \left\Vert D^{\alpha
}u\right\Vert _{\beta,\infty,\Omega}\\
& =\left(  \max_{\gamma\leq\alpha}\left\Vert D^{\gamma}\phi_{k}\right\Vert
_{\infty,\Omega}\right)  \left\Vert D^{\alpha}u\right\Vert _{\beta
,\infty,\Omega}\\
& \leq\left\Vert D^{\alpha}\phi_{k}\right\Vert _{\beta,\infty,\Omega
}\left\Vert D^{\alpha}u\right\Vert _{\beta,\infty,\Omega},
\end{align*}

so we can write%
\begin{equation}
\left\Vert E_{\Omega}^{\alpha}u\right\Vert _{\alpha;\infty;\Omega}\leq
k_{\alpha,v}\left(  \sum_{k=1}^{M}\left\Vert \phi_{k}\right\Vert
_{\alpha,\infty,\Omega}\right)  \left\Vert u\right\Vert _{\alpha,\infty
,\Omega},\label{X59}%
\end{equation}

where%
\begin{equation}
\left.
\begin{array}
[c]{ll}%
k_{\alpha,v} & =\int_{\mathcal{O}_{\mathbf{1}}}\frac{\left(  1+x\right)
^{\alpha}\left(  2+x\right)  ^{\alpha}}{x^{\alpha}}\left\vert v\left(
x\right)  \right\vert dx,\\
& =\prod\limits_{i=1}^{d}\int_{a_{i}}^{b_{i}}\left(  \frac{\left(  1+t\right)
\left(  2+t\right)  }{t}\right)  ^{\alpha_{i}}\left\vert v_{i}\left(
t\right)  \right\vert dt.
\end{array}
\right\} \label{X34}%
\end{equation}

This establishes the continuity of the extension operator i.e. $E_{\Omega
}:C_{B}^{\alpha}\left(  \overline{\Omega}\right)  \rightarrow C_{B}^{\alpha
}\left(  \mathbb{R}^{d}\right)  $ is continuous under supremum norms.
\end{proof}

\section{A contin. exten. $E_{\Omega}^{n\mathbf{1}}:W^{n\mathbf{1}}\left(
\Omega\right)  \rightarrow W^{n\mathbf{1}}$ derived$\mathbb{\ }$using spaces
of contin. funcs.\label{Sect_ExtenLocWn1toWn1_ContinFunc}}

In Section \ref{Sect_ExtenLocCatoBndCa} we constructed the continuous
extension operators $\left\{  E_{\Omega}^{\alpha}\right\}  _{\alpha
\geq\mathbf{1}}$ such that each $E_{\Omega}^{\alpha}:C^{\left(  \alpha\right)
}\left(  \overline{\Omega}\right)  \rightarrow C_{B}^{\left(  \alpha\right)
}$ is continuous under the supremum norm.

In this section we show that the particular extension operators $\left\{
E_{\Omega}^{n\mathbf{1}}\right\}  _{n\geq1}$ are such that each $E_{\Omega
}^{n\mathbf{1}}:W^{n\mathbf{1}}\left(  \Omega\right)  \rightarrow
W^{n\mathbf{1}}\left(  \mathbb{R}^{d}\right)  $ is continuous if we assume
that the integral operator $J_{\left(  n-2\right)  \mathbf{1}}^{\theta}\left[
v\right]  $ is a tensor product function with property \ref{X551}.

\begin{theorem}
\label{Thm_exten_Hn1(open)_to_Hn1}Suppose $n\geq1$ and $\Omega$ has the
uniform rectangle property. Suppose the function $v\in L_{0}^{1}\left(
\mathcal{O}_{\mathbf{1}}\right)  $ used to define the operator $J_{\left(
n-2\right)  \mathbf{1}}^{\theta}\left[  v\right]  $ is a tensor product
function with property \ref{X551}. Then the extension operator $E_{\Omega
}^{n\mathbf{1}}$ of Theorem \ref{Thm_ExtenContinFuncs_OrthantProp} satisfies
$E_{\Omega}^{n\mathbf{1}}:W^{n\mathbf{1}}\left(  \Omega\right)
\overset{c}{\rightarrow}W^{n\mathbf{1}}$. In fact%
\[
\left\Vert E_{\Omega}^{n\mathbf{1}}u\right\Vert _{W^{n\mathbf{1}}}\leq
k_{v,n}\left(  \sum_{k=1}^{M}\left\Vert \phi_{k}\right\Vert _{W^{n\mathbf{1}%
,\infty}\left(  \Omega\right)  }\right)  \left\Vert u\right\Vert
_{W^{n\mathbf{1}}\left(  \Omega\right)  },
\]

where $k_{v,n}$ is given by \ref{X28} in the proof .
\end{theorem}

\begin{proof}
The uniform rectangle property implies the segment property holds and so part
4 of Remark \ref{Rem_SobolevSpace2} implies that $C^{\left(  n\mathbf{1}%
\right)  }\left(  \overline{\Omega}\right)  $ is dense in $W^{n\mathbf{1}%
}\left(  \Omega\right)  $. We will follow the proof of Theorem
\ref{Thm_ExtenOrthantSobolFourier} and use a partition of unity $\left\{
\phi_{k}\right\}  _{k=1}^{M}$ to define each $u_{k}=\phi_{k}u$ and set
$E_{\Omega}^{n\mathbf{1}}u:=\sum_{k=1}^{M}E_{k}u_{k}$ so that%
\begin{equation}
\left\Vert E_{\Omega}^{n\mathbf{1}}u\right\Vert _{W^{n\mathbf{1}}}\leq
\sum_{k=1}^{M}\left\Vert E_{k}u_{k}\right\Vert _{W^{n\mathbf{1}}},\quad u\in
C^{\left(  n\mathbf{1}\right)  }\left(  \overline{\Omega}\right)
.\label{X102}%
\end{equation}

From \ref{X25} of the previous section, if $u\in C^{\left(  n\mathbf{1}%
\right)  }\left(  \overline{\Omega}\right)  $ then $u_{k,\beta}\in
C_{0}^{\left(  n\mathbf{1}-\beta\right)  }\left(  \Omega\right)  \subset
L^{2}$ and%
\[
D^{\gamma}E_{k}u_{k}=\tfrac{\left(  2\pi\right)  ^{d/2}}{\left(  \left(
n-1\right)  !\right)  ^{d}}\sum_{\beta\leq n\mathbf{1}}\tbinom{n\mathbf{1}%
}{\beta}\left(  D^{\gamma}I_{k,\beta}\right)  \ast u_{k,\beta},\quad
\forall\gamma,
\]

where $I_{k,\beta}=J_{\alpha-2}^{\theta^{k}}\left[  \left(  v^{k}\right)
_{\beta}\right]  _{-}$ for $v^{k}\in L_{0}^{1}\left(  \mathcal{O}^{k}\right)
$ is given by \ref{X01}.

But $D^{\gamma}I_{k,\beta}\in L^{1}$ and so by Young's theorem%
\begin{align}
\left\Vert D^{\gamma}E_{k}u_{k}\right\Vert _{2}  & \leq\tfrac{\left(
2\pi\right)  ^{d/2}}{\left(  \left(  n-1\right)  !\right)  ^{d}}\sum
_{\beta\leq n\mathbf{1}}\tbinom{n\mathbf{1}}{\beta}\left\Vert \left(
D^{\gamma}I_{k,\beta}\right)  \ast u_{k,\beta}\right\Vert _{2}\nonumber\\
& \leq\tfrac{1}{\left(  \left(  n-1\right)  !\right)  ^{d}}\sum_{\beta\leq
n\mathbf{1}}\tbinom{n\mathbf{1}}{\beta}\left\Vert D^{\gamma}I_{k,\beta
}\right\Vert _{1}\left\Vert u_{k,\beta}\right\Vert _{2}\nonumber\\
& \leq\tfrac{1}{\left(  \left(  n-1\right)  !\right)  ^{d}}\left(  \sum
_{\beta\leq n\mathbf{1}}\tbinom{n\mathbf{1}}{\beta}\left\Vert D^{\gamma
}I_{k,\beta}\right\Vert _{1}^{2}\right)  ^{1/2}\left(  \sum_{\beta\leq
n\mathbf{1}}\tbinom{n\mathbf{1}}{\beta}\left\Vert u_{k,\beta}\right\Vert
_{2}^{2}\right)  ^{1/2},\label{1.095}%
\end{align}

and from \ref{X26} and \ref{X27} of the proof of Theorem
\ref{Thm_ExtenOrthantSobolFourier},%
\[
\left(  \sum_{\beta\leq n\mathbf{1}}\tbinom{n\mathbf{1}}{\beta}\left\Vert
u_{k,\beta}\right\Vert _{2}^{2}\right)  ^{1/2}=\left\Vert u_{k}\right\Vert
_{W^{n\mathbf{1}}\left(  U_{k}\cap\Omega\right)  }\leq2^{nd/2}\left\Vert
\phi_{k}\right\Vert _{W^{n\mathbf{1},\infty}\left(  \Omega\right)  }\left\Vert
u\right\Vert _{W^{n\mathbf{1}}\left(  \Omega\right)  },
\]

where the constant $2^{nd/2}$ comes from Theorem \ref{Thm_norm_Cm1xHm1_loc} in
the Appendix. Thus \ref{1.095} becomes%
\[
\left\Vert D^{\gamma}E_{k}u_{k}\right\Vert _{2}\leq\tfrac{2^{nd/2}}{\left(
\left(  n-1\right)  !\right)  ^{d}}\left(  \sum_{\beta\leq n\mathbf{1}}%
\tbinom{n\mathbf{1}}{\beta}\left\Vert D^{\gamma}I_{k,\beta}\right\Vert
_{1}^{2}\right)  ^{1/2}\left\Vert \phi_{k}\right\Vert _{W^{n\mathbf{1},\infty
}\left(  \Omega\right)  }\left\Vert u\right\Vert _{W^{n\mathbf{1}}\left(
\Omega\right)  }.
\]

\begin{align}
&  \left\Vert E_{k}u_{k}\right\Vert _{W^{n\mathbf{1}}\left(  \Omega\right)
}=\left(  \sum_{\gamma\leq n\mathbf{1}}\tbinom{n\mathbf{1}}{\gamma}\left\Vert
D^{\gamma}E_{k}u_{k}\right\Vert _{2}^{2}\right)  ^{1/2}\leq\nonumber\\
&  \leq\tfrac{2^{nd/2}}{\left(  \left(  n-1\right)  !\right)  ^{d}}\left(
\sum_{\beta,\gamma\leq n\mathbf{1}}\tbinom{n\mathbf{1}}{\beta}\tbinom
{n\mathbf{1}}{\gamma}\left\Vert D^{\gamma}I_{k,\beta}\right\Vert _{1}%
^{2}\right)  ^{1/2}\left\Vert \phi_{k}\right\Vert _{W^{n\mathbf{1},\infty
}\left(  \Omega\right)  }\left\Vert u\right\Vert _{W^{n\mathbf{1}}\left(
\Omega\right)  }\label{X101}%
\end{align}

From \ref{X00} and then part 4 of Lemma \ref{Lem_PWCinf_Jtheta_b} (?? CHECK!)%
\begin{align*}
\sum_{\beta,\gamma\leq n\mathbf{1}}\tbinom{n\mathbf{1}}{\beta} &
\tbinom{n\mathbf{1}}{\gamma}\left\Vert D^{\gamma}I_{k,\beta}\right\Vert
_{1}^{2}\\
&  =\sum_{\beta,\gamma\leq n\mathbf{1}}\tbinom{n\mathbf{1}}{\beta}%
\tbinom{n\mathbf{1}}{\gamma}\left\Vert D^{\gamma}J_{\left(  n-2\right)
\mathbf{1}}^{\theta^{k}}\left[  \left(  v^{k}\right)  _{\beta}\right]
_{-}\right\Vert _{1}^{2}\\
&  =\sum_{\beta,\gamma\leq n\mathbf{1}}\tbinom{n\mathbf{1}}{\beta}%
\tbinom{n\mathbf{1}}{\gamma}\left\Vert D^{\gamma}J_{\left(  n-2\right)
\mathbf{1}}^{\theta^{k}}\left[  \left(  v^{k}\right)  _{\beta}\right]
\right\Vert _{1}^{2}\\
&  \leq\sum_{\beta,\gamma\leq n\mathbf{1}}\tbinom{n\mathbf{1}}{\beta}%
\tbinom{n\mathbf{1}}{\gamma}\left(  \left\Vert \frac{\left(  v^{k}\right)
_{\beta}}{x^{\gamma}}\right\Vert _{1,\mathcal{O}_{\theta^{k}}}\int%
\limits_{\mathcal{O}_{\mathbf{1}}}q_{\gamma}^{\left(  n-2\right)  \mathbf{1}%
}\left(  y\right)  e^{-y\mathbf{1}}dy\right)  ^{2},
\end{align*}

but we have assumed that $v^{k}\left(  x\right)  =v\left(  \theta
^{k}.x\right)  $ so $\left\Vert \frac{x^{\beta}v^{k}}{x^{\gamma}}\right\Vert
_{1,\mathcal{O}_{\theta^{k}}}=\left\Vert \frac{x^{\beta}v}{x^{\gamma}%
}\right\Vert _{1,\mathcal{O}_{\mathbf{1}}}$ and%
\[
\sum_{\beta\leq n\mathbf{1}}\tbinom{n\mathbf{1}}{\beta}\tbinom{n\mathbf{1}%
}{\gamma}\left\Vert D^{\gamma}I_{k,\beta}\right\Vert _{1}^{2}\leq\sum
_{\beta,\gamma\leq n\mathbf{1}}\tbinom{n\mathbf{1}}{\beta}\tbinom{n\mathbf{1}%
}{\gamma}\left(  \left\Vert \frac{x^{\beta}}{x^{\gamma}}v\right\Vert
_{1,\mathcal{O}_{\mathbf{1}}}\int\limits_{\mathcal{O}_{\mathbf{1}}}q_{\gamma
}^{\left(  n-2\right)  \mathbf{1}}\left(  y\right)  e^{-y\mathbf{1}}dy\right)
^{2},
\]

From \ref{X77},%
\[
\int_{\mathcal{O}_{\mathbf{1}}}q_{\gamma}^{\left(  n-2\right)  \mathbf{1}%
}\left(  y\right)  e^{-y\mathbf{1}}dy\leq\left(  \left(  n-1\right)  !\right)
^{d}2^{\left\vert \gamma\right\vert },
\]

so%
\begin{align*}
\sum_{\beta\leq n\mathbf{1}} &  \tbinom{n\mathbf{1}}{\beta}\tbinom
{n\mathbf{1}}{\gamma}\left\Vert D^{\gamma}I_{k,\beta}\right\Vert _{1}^{2}\\
&  \leq\sum_{\beta,\gamma\leq n\mathbf{1}}\tbinom{n\mathbf{1}}{\beta}%
\tbinom{n\mathbf{1}}{\gamma}\left(  \left\Vert \frac{x^{\beta}}{x^{\gamma}%
}v\right\Vert _{1,\mathcal{O}_{\mathbf{1}}}\left(  \left(  n-1\right)
!\right)  ^{d}2^{\left\vert \gamma\right\vert }\right)  ^{2}\\
&  =\left(  \left(  n-1\right)  !\right)  ^{2d}\sum_{\beta,\gamma\leq
n\mathbf{1}}\tbinom{n\mathbf{1}}{\beta}\tbinom{n\mathbf{1}}{\gamma}\left\Vert
2^{\left\vert \gamma\right\vert }\frac{x^{\beta}}{x^{\gamma}}v\right\Vert
_{1,\mathcal{O}_{\mathbf{1}}}^{2}\\
&  \leq\left(  \left(  n-1\right)  !\right)  ^{2d}\max_{\beta,\gamma\leq
n\mathbf{1}}\left(  2^{\left\vert \gamma\right\vert }\left\Vert \frac
{x^{\beta}}{x^{\gamma}}v\right\Vert _{1,\mathcal{O}_{\mathbf{1}}}\right)
\sum_{\beta,\gamma\leq n\mathbf{1}}\tbinom{n\mathbf{1}}{\beta}\tbinom
{n\mathbf{1}}{\gamma}\left\Vert 2^{\left\vert \gamma\right\vert }%
\frac{x^{\beta}}{x^{\gamma}}v\right\Vert _{1,\mathcal{O}_{\mathbf{1}}}\\
&  =\left(  \left(  n-1\right)  !\right)  ^{2d}\left(  \max_{\beta,\gamma\leq
n\mathbf{1}}2^{\left\vert \gamma\right\vert }\left\Vert x^{\beta-\gamma
}v\right\Vert _{1,\mathcal{O}_{\mathbf{1}}}\right)  \sum_{\beta,\gamma\leq
n\mathbf{1}}\tbinom{n\mathbf{1}}{\beta}\tbinom{n\mathbf{1}}{\gamma}\left\Vert
\frac{x^{\beta}}{\left(  x/2\right)  ^{\gamma}}v\right\Vert _{1,\mathcal{O}%
_{\mathbf{1}}}\\
&  =\left(  \left(  n-1\right)  !\right)  ^{2d}\left(  \max_{\beta,\gamma\leq
n\mathbf{1}}2^{\left\vert \gamma\right\vert }\left\Vert x^{\beta-\gamma
}v\right\Vert _{1,\mathcal{O}_{\mathbf{1}}}\right)  \sum_{\beta,\gamma\leq
n\mathbf{1}}\tbinom{n\mathbf{1}}{\beta}\tbinom{n\mathbf{1}}{\gamma}%
\int\limits_{\mathcal{O}_{\mathbf{1}}}\frac{x^{\beta}}{\left(  x/2\right)
^{\gamma}}\left\vert v\left(  x\right)  \right\vert dx\\
&  =\left(  \left(  n-1\right)  !\right)  ^{2d}\left(  \max_{\beta,\gamma\leq
n\mathbf{1}}2^{\left\vert \gamma\right\vert }\left\Vert x^{\beta-\gamma
}v\right\Vert _{1,\mathcal{O}_{\mathbf{1}}}\right)  \times\\
&  \qquad\qquad\times\int\limits_{\mathcal{O}_{\mathbf{1}}}\left(
\sum_{\gamma\leq n\mathbf{1}}\tbinom{n\mathbf{1}}{\gamma}\left(  \frac{1}%
{x/2}\right)  ^{\gamma}\right)  \left(  \sum_{\beta\leq n\mathbf{1}}%
\tbinom{n\mathbf{1}}{\beta}x^{\beta}\right)  \left\vert v\left(  x\right)
\right\vert dx\\
&  =\left(  \left(  n-1\right)  !\right)  ^{2d}\left(  \max_{\beta,\gamma\leq
n\mathbf{1}}2^{\left\vert \gamma\right\vert }\left\Vert x^{\beta-\gamma
}v\right\Vert _{1,\mathcal{O}_{\mathbf{1}}}\right)  \int\limits_{\mathcal{O}%
_{\mathbf{1}}}\left(  1+\frac{1}{x/2}\right)  ^{n\mathbf{1}}\left(
1+x\right)  ^{n\mathbf{1}}\left\vert v\left(  x\right)  \right\vert dx\\
&  =\left(  \left(  n-1\right)  !\right)  ^{2d}\left(  \max_{\beta,\gamma\leq
n\mathbf{1}}2^{\left\vert \gamma\right\vert }\int\limits_{\mathcal{O}%
_{\mathbf{1}}}\frac{x^{\beta}}{x^{\gamma}}\left\vert v\left(  x\right)
\right\vert dx\right)  \int\limits_{\mathcal{O}_{\mathbf{1}}}\frac{\left(
2+x\right)  ^{n\mathbf{1}}\left(  1+x\right)  ^{n\mathbf{1}}}{x^{n\mathbf{1}}%
}\left\vert v\left(  x\right)  \right\vert dx.
\end{align*}

and consequently \ref{X101} becomes%
\begin{align*}
&  \left\Vert E_{k}u_{k}\right\Vert _{W^{n\mathbf{1}}\left(  \Omega\right)
}\\
&  \leq\tfrac{2^{nd/2}}{\left(  \left(  n-1\right)  !\right)  ^{d}}\left(
\left(  \left(  n-1\right)  !\right)  ^{2d}\left(  \max_{\beta,\gamma\leq
n\mathbf{1}}2^{\left\vert \gamma\right\vert }\int\limits_{\mathcal{O}%
_{\mathbf{1}}}\frac{x^{\beta}}{x^{\gamma}}\left\vert v\left(  x\right)
\right\vert dx\right)  \int\limits_{\mathcal{O}_{\mathbf{1}}}\frac{\left(
2+x\right)  ^{n\mathbf{1}}\left(  1+x\right)  ^{n\mathbf{1}}}{x^{n\mathbf{1}}%
}\left\vert v\left(  x\right)  \right\vert dx\right)  ^{1/2}\times\\
&  \qquad\qquad\times\left\Vert \phi_{k}\right\Vert _{W^{n\mathbf{1},\infty
}\left(  \Omega\right)  }\left\Vert u\right\Vert _{W^{n\mathbf{1}}\left(
\Omega\right)  }\\
&  =k_{v,n}\left\Vert \phi_{k}\right\Vert _{W^{n\mathbf{1},\infty}\left(
\Omega\right)  }\left\Vert u\right\Vert _{W^{n\mathbf{1}}\left(
\Omega\right)  },
\end{align*}

where%
\begin{equation}
k_{v,n}=2^{nd/2}\left(  \left(  \max_{\beta,\gamma\leq n\mathbf{1}%
}2^{\left\vert \gamma\right\vert }\int\limits_{\mathcal{O}_{\mathbf{1}}}%
\frac{x^{\beta}}{x^{\gamma}}\left\vert v\left(  x\right)  \right\vert
dx\right)  \int\limits_{\mathcal{O}_{\mathbf{1}}}\frac{\left(  2+x\right)
^{n\mathbf{1}}\left(  1+x\right)  ^{n\mathbf{1}}}{x^{n\mathbf{1}}}\left\vert
v\left(  x\right)  \right\vert dx\right)  ^{1/2}.\label{X28}%
\end{equation}

Inequality \ref{X102} now yields our theorem.
\end{proof}

\begin{remark}
\label{Rem_Thm_exten_Hn1(open)_to_Hn1}We should be able to show that if $v$
has property of \ref{X551} for $\alpha$ then the extension operator
$E_{\Omega}^{\alpha}$ of Section \ref{Sect_ExtenLocCatoBndCa} is in fact a
continuous extension $E_{\Omega}^{\alpha}:W^{\alpha}\left(  \Omega\right)
\rightarrow W^{\alpha}\left(  \mathbb{R}^{d}\right)  $.

Following the Russian mathematicians in the Background section we would need
to define
\[
W^{\alpha}\left(  \Omega\right)  :=\left\{  u\in L^{2}\left(  \Omega\right)
:D^{\beta}u\in L^{2}\left(  \Omega\right)  \text{ when }\beta\leq
\alpha\right\}  .
\]

\end{remark}

\section{Some unused estimates\label{Sect_UnusedEstim}}

?? In this section ??

The next result shows that to prove $f\in C^{\left(  \alpha\right)  }$ we may
need only calculate one derivative, namely $D^{\alpha}$.

\begin{lemma}
\label{Lem_intermed_deriv_contin}??\textbf{FINISH}!

\begin{enumerate}
\item If $f\in\mathcal{D}^{\prime}$ and $D^{\alpha}f\in C^{\left(  0\right)
}$ then $D^{\beta}f\in C^{\left(  0\right)  }$ when $\beta\leq\alpha$ i.e.
$f\in C^{\left(  \alpha\right)  }$.

\item May need to assume boundedness e.g. If $f\in L^{\infty}$ and in the
distribution sense $D^{\alpha}f\in C_{B}^{\left(  0\right)  }$ then $D^{\beta
}f\in C_{B}^{\left(  0\right)  }$ when $\beta\leq\alpha$ i.e. $f\in
C_{B}^{\left(  \alpha\right)  }$.
\end{enumerate}
\end{lemma}

\begin{proof}
\textbf{Part 1} Use induction on $\max\alpha$? ?? First prove true for
$\max\alpha=1$: suppose $f\in??\mathcal{D}^{\prime}$ and $D_{1}f\in C^{\left(
0\right)  }$. Define%
\[
g\left(  x\right)  :=\int_{0}^{x_{1}}D_{1}f\left(  s,x^{\prime}\right)  ds,
\]

which exists since $D_{1}f\in L_{loc}^{1}$. Then%
\[
\left\vert g\left(  x\right)  -g\left(  y\right)  \right\vert \leq\int%
_{0}^{x_{1}}\left\vert D_{1}f\left(  s,x^{\prime}\right)  ds-D_{1}f\left(
s,y^{\prime}\right)  \right\vert ds+\left\vert \int_{x_{1}}^{y_{1}}\left\vert
D_{1}f\left(  s,y^{\prime}\right)  \right\vert ds\right\vert .
\]

Since $D_{1}f$ is continuous on $\mathbb{R}^{d}$ it is uniformly continuous on
any bounded closed rectangle. Suppose $\mathcal{R}$ is a rectangle which
contains $x$ and $\mathbf{0}$. etc. which proves $g$ is continuous.

Show that for each $x^{\prime}$, $g\in C^{\left(  1\right)  }\left(
\mathbb{R}^{1}\right)  $ and $D_{1}g\left(  x\right)  =D_{1}f\left(  x\right)
$ a.e. and $D_{1}g=D_{1}f$ in $\mathcal{D}^{\prime}$ so that $f-g$ is a
constant and hence $f$ is also continuous.

?? Generalize the following argument!: Suppose $f\in\mathcal{D}^{\prime}$ and
$D_{1}D_{2}f\in C^{\left(  0\right)  }$. Then $D_{2}f\in\mathcal{D}^{\prime}$
and $D_{1}D_{2}f\in C^{\left(  0\right)  }$ so $D_{2}f\in C^{\left(  0\right)
}$. Similarly $D_{1}f\in C^{\left(  0\right)  }$.

\textbf{Part 2} ?? Suppose $f\in L^{\infty}$ and $D_{1}D_{2}f\in
C_{B}^{\left(  0\right)  }$. From part 1, $f\in C^{\left(  \mathbf{1}\right)
}\cap L^{\infty}$.%
\[
D_{1}f\left(  x\right)  =\int_{0}^{x_{2}}D_{1}D_{2}f\left(  x_{1}%
,y_{2}\right)  dy_{2}.
\]

=====================

Suppose $D^{\mathbf{1}}f\in C_{B}^{\left(  0\right)  }\cap L^{1}$ and $f\in
C_{B}^{\left(  \mathbf{1}\right)  }\cap L^{1}$. Then%
\[
\left(  2\pi\right)  ^{d/2}H\ast D^{\mathbf{1}}f=\left(  2\pi\right)
^{d/2}D^{\mathbf{1}}H\ast f=\left(  2\pi\right)  ^{d/2}\delta\ast f=f.
\]

\end{proof}

\begin{lemma}
\label{Lem_f_inL1_D1f_inL1}?? \textbf{NEEDED}? Prove that%
\[
f\in L^{1},\text{ }D^{\mathbf{1}}f\in L^{1}\Rightarrow D^{\gamma}f\in
L^{1}\text{ }for\text{ }all\text{ }\gamma\leq\mathbf{1},
\]

and%
\[
\left\Vert D^{\gamma}f\right\Vert _{1}\leq\left(  2\pi\right)  ^{\left(
d-\left\vert \gamma\right\vert \right)  /2}\left\Vert D^{\mathbf{1}%
}f\right\Vert _{1},\quad\gamma\leq\mathbf{1}.
\]

\end{lemma}

\begin{proof}
First assume that $f\in L_{0}^{1}\subset\mathcal{E}^{\prime}$. Then we can use
the $\mathcal{D}^{\prime}\ast\mathcal{E}^{\prime}$ convolution. Using the
tensor product Heavyside step function $H\in\mathcal{D}^{\prime}$, for $f\in
L^{1}$ we can write%
\[
\left(  2\pi\right)  ^{d/2}H\ast D^{\mathbf{1}}f=\left(  2\pi\right)
^{d/2}D^{\mathbf{1}}H\ast f=\left(  2\pi\right)  ^{d/2}\delta\ast f=f.
\]

Hence, using the concepts of partial delta functions and partial convolutions
?? ADD REF!,%
\begin{align*}
D^{\gamma}f  & =\left(  2\pi\right)  ^{d/2}\left(  D^{\gamma}H\right)  \ast
D^{\mathbf{1}}f\\
& =\left(  2\pi\right)  ^{d/2}\left(  \delta^{\prime}H^{\prime\prime}\right)
\ast D^{\mathbf{1}}f\\
& =\left(  2\pi\right)  ^{d/2}H^{\prime\prime}\ast^{\prime\prime}\left(
\delta^{\prime}\ast^{\prime}D^{\mathbf{1}}f\right) \\
& =\left(  2\pi\right)  ^{\left(  d-\left\vert \gamma\right\vert \right)
/2}H^{\prime\prime}\ast^{\prime\prime}D^{\mathbf{1}}f.
\end{align*}

Thus%
\[
\left\Vert D^{\gamma}f\right\Vert _{1}^{\prime\prime}\leq\left(  2\pi\right)
^{\left(  d-\left\vert \gamma\right\vert \right)  /2}\left\Vert H^{\prime
\prime}\right\Vert _{\infty}\left\Vert D^{\mathbf{1}}f\right\Vert _{1}%
^{\prime\prime}=\left(  2\pi\right)  ^{\left(  d-\left\vert \gamma\right\vert
\right)  /2}\left\Vert D^{\mathbf{1}}f\right\Vert _{1}^{\prime\prime},
\]

and so when $f\in L_{0}^{1}$,%
\[
\left\Vert D^{\gamma}f\right\Vert _{1}\leq\left(  2\pi\right)  ^{\left(
d-\left\vert \gamma\right\vert \right)  /2}\left\Vert D^{\mathbf{1}%
}f\right\Vert _{1}.
\]

Now suppose $f\in L^{1}$. Choose the smooth cut-off function $\phi_{0}\in
C_{0}^{\infty}$ such that $\phi_{0}=1$ on the closed rectangle $R\left[
-\mathbf{1},\mathbf{1}\right]  $ and $\phi_{0}=0$ outside the open rectangle
$R\left(  -2\mathbf{1},2\mathbf{1}\right)  $. Then $\phi_{0}\left(
x/k\right)  f\in L_{0}^{1}$ and $\phi_{0}\left(  x/k\right)  f\rightarrow f$
in $L^{1}$ as integer $k\rightarrow\infty$. Also ?? $D^{\mathbf{1}}\left(
\phi_{0}\left(  x/k\right)  f\right)  =\ldots+\phi_{0}\left(  x/k\right)
D^{\mathbf{1}}f\rightarrow0+D^{\mathbf{1}}f$ in $L^{1}$.
\end{proof}

?? A distribution approach to the last lemma is:

\begin{lemma}
\label{Lem_integ_repres_distrib}?? REMOVE? Suppose $g\in\mathcal{E}^{\prime}$,
$D^{\alpha}g\in C_{0}^{\left(  0\right)  }$ for some $\alpha\geq\mathbf{1}$. Then:

\begin{enumerate}
\item $g\in C_{0}^{\left(  \alpha\right)  }$;

\item When $\gamma\leq\alpha-1$,%
\begin{align*}
\tfrac{1}{\left(  \alpha-\gamma-1\right)  !}\int\left(  x-y\right)
^{\alpha-\gamma-1}H\left(  x-y\right)  D^{\alpha}g\left(  y\right)  dy  &
=D^{\gamma}g\left(  x\right)  .\\
\int\chi_{+}^{\alpha-\gamma-1}\left(  x-y\right)  D^{\alpha}g\left(  y\right)
dy  & =D^{\gamma}g\left(  x\right)  .
\end{align*}

\end{enumerate}
\end{lemma}

\begin{proof}
\textbf{Part 1} Show: if $g\in\mathcal{D}^{\prime}$ and $D^{\alpha}g\in
C^{\left(  0\right)  }$ then $g\in C^{\left(  \alpha\right)  }$.\medskip

\textbf{Part 2} We begin with the identity%
\[
D^{\beta}\left(  \frac{x^{\beta-1}}{\left(  \beta-1\right)  !}H\right)
=\delta,\quad\beta\geq\mathbf{1}.
\]

From Chapter 13 of Duistermaat and Kolk \cite{DuistKolk2010},%
\[
\chi_{+}^{\beta}:=\frac{x^{\beta-1}}{\left(  \beta-1\right)  !}H=\frac
{x^{\beta-1}}{\Gamma\left(  \beta\right)  }H,\quad\chi_{+}^{\mathbf{0}%
}:=\delta,
\]

where $H\in S^{\prime}$ is the Heavyside step function. Then using the
$S^{\prime},\mathcal{E}^{\prime}$ convolution we have%
\begin{equation}
\left(  2\pi\right)  ^{d/2}\left(  \frac{x^{\beta-1}}{\left(  \beta-1\right)
!}H\right)  \ast D^{\beta}g=g,\quad\beta\geq\mathbf{1}.\label{X83}%
\end{equation}

From DuKo \cite{DuistKolk2010}, $D^{\gamma}\chi_{+}^{\alpha}=\chi_{+}%
^{\alpha-\gamma}$ when $\gamma\leq\alpha-1$ and so that from \ref{X83},%
\[
\left(  2\pi\right)  ^{d/2}\left(  \frac{x^{\alpha-\gamma-1}}{\left(
\alpha-\gamma-1\right)  !}H\right)  \ast D^{\alpha}g=D^{\gamma}g,\quad
\alpha\geq\mathbf{1},
\]

and using the $L_{loc}^{1},C_{0}^{\left(  0\right)  }$ convolution gives%
\[
\tfrac{1}{\left(  \alpha-\gamma-1\right)  !}\int\left(  x-y\right)
^{\alpha-\gamma-1}H\left(  x-y\right)  D^{\alpha}g\left(  y\right)
dy=D^{\gamma}g\left(  x\right)  .
\]

\end{proof}

\begin{remark}
\label{Rem_Lem_integ_repres}?? REMOVE? What about generalizing \ref{X20} using
the tempered distribution Taylor series expansions of Chapter
\ref{Ch_interp_err_Taylor_temper_distrib}?

Assume $D^{n}f\in L_{loc}^{1}\left(  \left[  0,\infty\right)  \right)  $ and
$f\in S^{\prime}$. This implies ?? $D^{k}f\in C^{\left(  n-k-1\right)
}\left(  \left[  0,\infty\right)  \right)  $ for $k\leq n-1$. Further assume
that $D^{k}f\rightarrow0$ exponentially for $k\leq n$.
\end{remark}

\begin{lemma}
\label{Lem_bnd_FourTran_J[v]_1}\textbf{Bounds for} $\left\Vert \xi^{\gamma
}\widehat{J_{\delta}\left[  \sigma\right]  }\right\Vert _{\infty}$
\textbf{when} $\delta\geq-\mathbf{1}$, $\gamma\leq\delta+2$.\medskip

\textbf{Part 1} If $\delta\geq-\mathbf{1}$ and $\sigma\in L_{0}^{1}\left(
\mathcal{O}_{\mathbf{1}}\right)  $ then $\widehat{J_{\delta}\left[
\sigma\right]  }\in C_{B}^{\left(  0\right)  }$ and%
\begin{equation}
\widehat{J_{\delta}\left[  \sigma\right]  }\left(  \xi\right)  =\int%
_{\mathcal{O}_{\mathbf{1}}}\widehat{\sigma}\left(  \xi.\zeta\right)
\zeta^{\delta+1}e^{-\zeta\mathbf{1}}d\zeta.\label{X171}%
\end{equation}

\end{lemma}

\begin{proof}
\textbf{Part 1} From Lemma \ref{Lem_SmthFuncIntegRepInOrthant}, $J_{\delta
}\left[  \sigma\right]  \in C^{\infty}\left(  \overline{\mathcal{O}_{1}%
}\right)  $, $\operatorname*{supp}J_{\delta}\left[  \sigma\right]
\subseteq\overline{\mathcal{O}_{1}}$ and the estimate \ref{X32} implies it
decreases exponentially at infinity. Thus $J_{\delta}\left[  \sigma\right]
\in L^{1}$ and so $\widehat{J_{\delta}\left[  \sigma\right]  }\in
C_{B}^{\left(  0\right)  }$, and from \ref{X31},%
\begin{align*}
\widehat{J_{\delta}\left[  \sigma\right]  }\left(  \xi\right)   & =\left(
2\pi\right)  ^{-d/2}\int e^{-i\xi z}J_{\delta}\left(  z\right)  dz\\
& =\left(  2\pi\right)  ^{-d/2}\int_{\mathcal{O}_{\mathbf{1}}}e^{-i\xi
z}z^{\delta+1}\int_{a}^{b}\frac{1}{\mu^{\delta+2}}e^{-\frac{z}{\mu}\mathbf{1}%
}\sigma\left(  \mu\right)  d\mu\text{ }dz\\
& =\left(  2\pi\right)  ^{-d/2}\int_{\mathcal{O}_{\mathbf{1}}}\int_{a}%
^{b}e^{-i\xi z}z^{\delta+1}\frac{e^{-\frac{z}{\mu}\mathbf{1}}}{\mu^{\delta+2}%
}\sigma\left(  \mu\right)  d\mu\text{ }dz,
\end{align*}

and this integral is absolutely convergent since,%
\begin{align*}
\int_{\mathcal{O}_{\mathbf{1}}}\int_{a}^{b}\left\vert e^{-i\xi z}z^{\delta
+1}\frac{e^{-\frac{z}{\mu}\mathbf{1}}}{\mu^{\delta+2}}\sigma\left(
\mu\right)  \right\vert d\mu\text{ }dz  & \leq\int_{\mathcal{O}_{\mathbf{1}}%
}\int_{a}^{b}z^{\delta+1}\frac{e^{-\frac{z}{\mu}\mathbf{1}}}{\mu^{\delta+2}%
}\left\vert \sigma\left(  \mu\right)  \right\vert d\mu\text{ }dz\\
& <\int_{\mathcal{O}_{\mathbf{1}}}\int_{a}^{b}z^{\delta+1}\frac{e^{-\frac
{z}{b}\mathbf{1}}}{a^{\delta+2}}\left\vert \sigma\left(  \mu\right)
\right\vert d\mu\text{ }dz\\
& <\infty,
\end{align*}

which allows us to use Fubini's theorem to justify changing the order of
integration. Indeed, for all multi-integers $\delta\geq-\mathbf{1}$,%
\begin{align*}
\widehat{J_{\delta}\left[  \sigma\right]  }\left(  \xi\right)   & =\left(
2\pi\right)  ^{-d/2}\int_{\mathcal{O}_{\mathbf{1}}}\int_{a}^{b}e^{-i\xi
z}z^{\delta+1}\frac{e^{-\frac{z}{\mu}\mathbf{1}}}{\mu^{\delta+2}}\sigma\left(
\mu\right)  d\mu\text{ }dz\\
& =\left(  2\pi\right)  ^{-d/2}\int_{a}^{b}\int_{\mathcal{O}_{\mathbf{1}}%
}e^{-i\xi z}z^{\delta+1}e^{-\frac{z}{\mu}\mathbf{1}}dz\text{ }\frac
{\sigma\left(  \mu\right)  }{\mu^{\delta+2}}d\mu\\
& \Rightarrow\zeta=z/\mu,\text{ }z=\mu.\zeta,\text{ }dz=\mu^{\mathbf{1}}%
d\zeta\Rightarrow\\
\widehat{J_{\delta}\left[  \sigma\right]  }\left(  \xi\right)   & =\left(
2\pi\right)  ^{-d/2}\int_{a}^{b}\int_{\mathcal{O}_{\mathbf{1}}}e^{-i\xi\left(
\mu.\zeta\right)  }\left(  \mu.\zeta\right)  ^{\delta+1}e^{-\zeta\mathbf{1}%
}\mu^{\mathbf{1}}d\zeta\text{ }\frac{\sigma\left(  \mu\right)  }{\mu
^{\delta+2}}d\mu\\
& =\left(  2\pi\right)  ^{-d/2}\int_{a}^{b}\mu^{\delta+1}\mu^{\mathbf{1}}%
\int_{\mathcal{O}_{\mathbf{1}}}e^{-i\left(  \xi.\zeta\right)  \mu}%
\zeta^{\delta+1}e^{-\zeta\mathbf{1}}d\zeta\text{ }\frac{\sigma\left(
\mu\right)  }{\mu^{\delta+2}}d\mu\\
& =\left(  2\pi\right)  ^{-d/2}\int_{a}^{b}\left(  \int_{\mathcal{O}%
_{\mathbf{1}}}e^{-i\left(  \xi.\zeta\right)  \mu}\zeta^{\delta+1}%
e^{-\zeta\mathbf{1}}d\zeta\right)  \sigma\left(  \mu\right)  d\mu\\
& \Rightarrow clearly\text{ }absolutely\text{ }convergent\Rightarrow\\
& =\left(  2\pi\right)  ^{-d/2}\int_{\mathcal{O}_{\mathbf{1}}}\left(  \int%
_{a}^{b}e^{-i\left(  \xi.\zeta\right)  \mu}\sigma\left(  \mu\right)
d\mu\right)  \zeta^{\delta+1}e^{-\zeta\mathbf{1}}d\zeta\\
& =\int_{\mathcal{O}_{\mathbf{1}}}\widehat{\sigma}\left(  \xi.\zeta\right)
\zeta^{\delta+1}e^{-\zeta\mathbf{1}}d\zeta,
\end{align*}

which proves \ref{X17}.
\end{proof}

\begin{lemma}
?? If $\sigma\in L^{1}\left(  \mathcal{O}_{\mathbf{1}}\right)  $ and
$\delta\geq-\mathbf{1}$ then%
\[
\widehat{J_{\delta}\left[  \sigma\right]  }\left(  \xi\right)  =\int%
_{\mathcal{O}_{\mathbf{1}}}\widehat{\sigma}\left(  \xi.\zeta\right)
\zeta^{\delta+1}e^{-\zeta\mathbf{1}}d\zeta.
\]

\end{lemma}

\begin{proof}
?? From ?? $J_{\delta}\left[  \sigma\right]  \in L^{1}$ and so
$\widehat{J_{\delta}\left[  \sigma\right]  }\in C_{B}^{\left(  0\right)  }$.
Thus%
\begin{align*}
\widehat{J_{\delta}\left[  \sigma\right]  }\left(  \xi\right)   & =\left(
2\pi\right)  ^{-d/2}\int e^{-i\xi z}J_{\delta}\left(  z\right)  dz\\
& =\left(  2\pi\right)  ^{-d/2}\int_{\mathcal{O}_{\mathbf{1}}}e^{-i\xi
z}z^{\delta+1}\int_{\mathcal{O}_{\mathbf{1}}}\frac{1}{\mu^{\delta+2}}%
e^{-\frac{z}{\mu}\mathbf{1}}\sigma\left(  \mu\right)  d\mu\text{ }dz\\
& =\left(  2\pi\right)  ^{-d/2}\int_{\mathcal{O}_{\mathbf{1}}}\int%
_{\mathcal{O}_{\mathbf{1}}}e^{-i\xi z}z^{\delta+1}\frac{e^{-\frac{z}{\mu
}\mathbf{1}}}{\mu^{\delta+2}}\sigma\left(  \mu\right)  d\mu\text{ }dz,
\end{align*}

and this integral is absolutely convergent since,%
\begin{align*}
\int_{\mathcal{O}_{\mathbf{1}}}\int_{\mathcal{O}_{\mathbf{1}}}\left\vert
e^{-i\xi z}z^{\delta+1}\frac{e^{-\frac{z}{\mu}\mathbf{1}}}{\mu^{\delta+2}%
}\sigma\left(  \mu\right)  \right\vert dz\text{ }d\mu & \leq\int%
_{\mathcal{O}_{\mathbf{1}}}\int_{\mathcal{O}_{\mathbf{1}}}z^{\delta+1}%
\frac{e^{-\frac{z}{\mu}\mathbf{1}}}{\mu^{\delta+2}}\left\vert \sigma\left(
\mu\right)  \right\vert dz\text{ }d\mu\\
& :\zeta=z./\mu\Rightarrow\\
& \text{?? FIX! ??}\\
& =\int_{\mathcal{O}_{\mathbf{1}}}\int_{\mathcal{O}_{\mathbf{1}}}\left(
\mu.\zeta\right)  ^{\delta+1}\frac{e^{-\zeta\mathbf{1}}}{\mu^{\delta+2}%
}\left\vert \sigma\left(  \mu\right)  \right\vert d\mu\text{ }\mu^{\mathbf{1}%
}d\zeta\\
& =\int_{\mathcal{O}_{\mathbf{1}}}\int_{\mathcal{O}_{\mathbf{1}}}\zeta
^{\delta+1}e^{-\zeta\mathbf{1}}\left\vert \sigma\left(  \mu\right)
\right\vert d\mu d\zeta\\
& =\left(  \delta+1\right)  !\left\Vert \sigma\right\Vert _{1}.
\end{align*}

This allows us to use Fubini's theorem to justify changing the order of
integration. Indeed, for all multi-integers $\delta\geq-\mathbf{1}$,%
\begin{align*}
\widehat{J_{\delta}\left[  \sigma\right]  }\left(  \xi\right)   & =\left(
2\pi\right)  ^{-d/2}\int_{\mathcal{O}_{\mathbf{1}}}\int_{\mathcal{O}%
_{\mathbf{1}}}e^{-i\xi z}z^{\delta+1}\frac{e^{-\frac{z}{\mu}\mathbf{1}}}%
{\mu^{\delta+2}}\sigma\left(  \mu\right)  d\mu\text{ }dz\\
& =\left(  2\pi\right)  ^{-d/2}\int_{\mathcal{O}_{\mathbf{1}}}\int%
_{\mathcal{O}_{\mathbf{1}}}e^{-i\xi z}z^{\delta+1}e^{-\frac{z}{\mu}\mathbf{1}%
}dz\text{ }\frac{\sigma\left(  \mu\right)  }{\mu^{\delta+2}}d\mu\\
& \Rightarrow\zeta=z./\mu,\text{ }z=\mu.\zeta,\text{ }dz=\mu^{\mathbf{1}%
}d\zeta\Rightarrow\\
\widehat{J_{\delta}\left[  \sigma\right]  }\left(  \xi\right)   & =\left(
2\pi\right)  ^{-d/2}\int_{\mathcal{O}_{\mathbf{1}}}\int_{\mathcal{O}%
_{\mathbf{1}}}e^{-i\xi\left(  \mu.\zeta\right)  }\left(  \mu.\zeta\right)
^{\delta+1}e^{-\zeta\mathbf{1}}\mu^{\mathbf{1}}d\zeta\text{ }\frac
{\sigma\left(  \mu\right)  }{\mu^{\delta+2}}d\mu\\
& =\left(  2\pi\right)  ^{-d/2}\int_{\mathcal{O}_{\mathbf{1}}}\mu^{\delta
+1}\mu^{\mathbf{1}}\int_{\mathcal{O}_{\mathbf{1}}}e^{-i\left(  \xi
.\zeta\right)  \mu}\zeta^{\delta+1}e^{-\zeta\mathbf{1}}d\zeta\text{ }%
\frac{\sigma\left(  \mu\right)  }{\mu^{\delta+2}}d\mu\\
& =\left(  2\pi\right)  ^{-d/2}\int_{\mathcal{O}_{\mathbf{1}}}\left(
\int_{\mathcal{O}_{\mathbf{1}}}e^{-i\left(  \xi.\zeta\right)  \mu}%
\zeta^{\delta+1}e^{-\zeta\mathbf{1}}d\zeta\right)  \sigma\left(  \mu\right)
d\mu\\
& \Rightarrow clearly\text{ }absolutely\text{ }convergent\Rightarrow\\
& =\left(  2\pi\right)  ^{-d/2}\int_{\mathcal{O}_{\mathbf{1}}}\left(
\int_{\mathcal{O}_{\mathbf{1}}}e^{-i\left(  \xi.\zeta\right)  \mu}%
\sigma\left(  \mu\right)  d\mu\right)  \zeta^{\delta+1}e^{-\zeta\mathbf{1}%
}d\zeta\\
& =\int_{\mathcal{O}_{\mathbf{1}}}\widehat{\sigma}\left(  \xi.\zeta\right)
\zeta^{\delta+1}e^{-\zeta\mathbf{1}}d\zeta,
\end{align*}

as claimed.
\end{proof}

The estimates below were used in Section \ref{Sect_ExtenLocWn1_to_Wn1_Fourier}
but were replaced by those of Lemma ?? \ref{Lem_bnd_FourTran_J[v]_1}? ??.
However, I decided to keep these estimates here.

We will need estimates for $\left\Vert \xi^{\gamma}\widehat{J_{\left(
n-2\right)  \mathbf{1}}^{\theta^{k}}\left[  \left(  v^{k}\right)  _{\beta
}\right]  }\left(  \xi\right)  \right\Vert _{\infty}$ for $\gamma\leq
n\mathbf{1}$ and $n\geq1$. The next lemma will do the job but first some notation:

\begin{notation}
\label{Not_partial_integ}Suppose $\gamma\in\left\{  0,1\right\}  ^{d}$ and
suppose $\gamma$ has $1$s at positions $i_{1},i_{2},i_{3},\ldots,i_{\left\vert
\gamma\right\vert }$.\medskip

\textbf{Partial integration} We define $d^{\gamma}\eta=d\eta_{i_{1}}%
d\eta_{i_{2}}\ldots d\eta_{i_{\left\vert \phi\right\vert }}$ and write $\int
f\left(  \eta\right)  d^{\gamma}\eta$.

\textbf{Projection operator} $\pi_{\gamma}\left(  x\right)  :=\left(
x_{i_{1}},x_{i_{2}},\ldots,x_{i_{\left\vert \gamma\right\vert }}\right)
\in\mathbb{R}^{\left\vert \gamma\right\vert }$ when $x\in\mathbb{R}^{d}$.

Set $\left\Vert \int\left\vert f\left(  \eta\right)  \right\vert
d^{\mathbf{0}}\eta\right\Vert _{\infty}:=\left\Vert f\right\Vert _{\infty}$.
\end{notation}

\begin{lemma}
\label{Lem_bnd_FourTran_J[v]}\textbf{Bounds for} $\left\Vert \xi^{\gamma
}\widehat{J_{\delta}\left[  \sigma\right]  }\right\Vert _{\infty}$
\textbf{when} $\delta\geq-\mathbf{1}$, $\gamma\leq\delta+2$.\medskip

\textbf{Part 1} If $\delta\geq-\mathbf{1}$ and $\sigma\in L_{0}^{1}\left(
\mathcal{O}_{\mathbf{1}}\right)  $ then $\widehat{J_{\delta}\left[
\sigma\right]  }\in C_{B}^{\left(  0\right)  }$ and%
\begin{equation}
\widehat{J_{\delta}\left[  \sigma\right]  }\left(  \xi\right)  =\int%
_{\mathcal{O}_{\mathbf{1}}}\widehat{\sigma}\left(  \xi.\zeta\right)
\zeta^{\delta+1}e^{-\zeta\mathbf{1}}d\zeta.\label{X17}%
\end{equation}

?? USE FOURIER\ TRANSFORM TO\ CALC SMOOTHNESS\ PROPERTIES?\medskip

\fbox{Suppose $\delta\geq-\mathbf{1}$.} If $\mathbf{1}\leq\gamma\leq\delta+2$
then for $\left\Vert \xi^{\gamma}\widehat{J_{\delta}\left[  \sigma\right]
}\right\Vert _{\infty}$ we have the upper bound \ref{X63}. If $\mathbf{0}%
\leq\gamma\leq\delta+1$ then for $\left\Vert \xi^{\gamma}\widehat{J_{\delta
}\left[  \sigma\right]  }\right\Vert _{\infty}$ we have the upper bounds
\ref{X64} and \ref{X62}.\medskip

\fbox{When $\delta=\mathbf{0}$ and $\gamma\leq\mathbf{1}$,} we have the upper
bounds \ref{X65} and \ref{X66} for $\left\Vert \xi^{\gamma}%
\widehat{J_{\mathbf{0}}\left[  \sigma\right]  }\right\Vert _{\infty}$.\medskip

\fbox{When $\delta=-\mathbf{1}$ and $\gamma\leq\mathbf{1}$,} if%
\begin{equation}
\left\Vert \int\left\vert \widehat{\sigma}\left(  \eta\right)  \right\vert
d^{\gamma}\eta\right\Vert _{\infty}<\infty,\label{X56}%
\end{equation}

we have the upper bound \ref{X67} for $\left\Vert \xi^{\gamma}%
\widehat{J_{\mathbf{-1}}\left[  \sigma\right]  }\right\Vert _{\infty}%
$.\medskip

\textbf{Part 2} This is the analogue of part 1 for $J_{\delta}^{\theta}\left[
\sigma\right]  $ which is in \ref{X05} by%
\[
J_{\delta}^{\theta}\left[  \sigma\right]  \left(  z\right)  :=J_{\delta
}\left[  \sigma\left(  \theta.\right)  \right]  \left(  \theta.z\right)
,\text{\quad}\sigma\in L_{0}^{1}\left(  \mathcal{O}_{\theta}\right)  .
\]

If $\delta\geq-\mathbf{1}$ and $\sigma\in L_{0}^{1}\left(  \mathcal{O}%
_{\theta}\right)  $ then $\widehat{J_{\delta}^{\theta}\left[  \sigma\right]
}\in C_{B}^{\left(  0\right)  }$ and
\[
\widehat{J_{\delta}^{\theta}\left[  \sigma\right]  }=\int_{\mathcal{O}%
_{\mathbf{1}}}\widehat{\sigma}\left(  \xi.\zeta\right)  \zeta^{\delta
+1}e^{-\zeta\mathbf{1}}d\zeta.
\]

Further, all the bounds described in part 1 also hold here.
\end{lemma}

\begin{proof}
\textbf{Part 1} From Lemma \ref{Lem_SmthFuncIntegRepInOrthant}, $J_{\delta
}\left[  \sigma\right]  \in C^{\infty}\left(  \overline{\mathcal{O}_{1}%
}\right)  $, $\operatorname*{supp}J_{\delta}\left[  \sigma\right]
\subset\overline{\mathcal{O}_{1}}$ and the estimate \ref{X32} implies it
decreases exponentially at infinity. Thus $J_{\delta}\left[  \sigma\right]
\in L^{1}$ and so $\widehat{J_{\delta}\left[  \sigma\right]  }\in
C_{B}^{\left(  0\right)  }$, and from \ref{X31},%
\begin{align*}
\widehat{J_{\delta}\left[  \sigma\right]  }\left(  \xi\right)   & =\left(
2\pi\right)  ^{-d/2}\int e^{-i\xi z}J_{\delta}\left(  z\right)  dz\\
& =\left(  2\pi\right)  ^{-d/2}\int_{\mathcal{O}_{\mathbf{1}}}e^{-i\xi
z}z^{\delta+1}\int_{a}^{b}\frac{1}{\mu^{\delta+2}}e^{-\frac{z}{\mu}\mathbf{1}%
}\sigma\left(  \mu\right)  d\mu\text{ }dz\\
& =\left(  2\pi\right)  ^{-d/2}\int_{\mathcal{O}_{\mathbf{1}}}\int_{a}%
^{b}e^{-i\xi z}z^{\delta+1}\frac{e^{-\frac{z}{\mu}\mathbf{1}}}{\mu^{\delta+2}%
}\sigma\left(  \mu\right)  d\mu\text{ }dz,
\end{align*}

and this integral is absolutely convergent since,%
\begin{align*}
\int_{\mathcal{O}_{\mathbf{1}}}\int_{a}^{b}\left\vert e^{-i\xi z}z^{\delta
+1}\frac{e^{-\frac{z}{\mu}\mathbf{1}}}{\mu^{\delta+2}}\sigma\left(
\mu\right)  \right\vert d\mu\text{ }dz  & \leq\int_{\mathcal{O}_{\mathbf{1}}%
}\int_{a}^{b}z^{\delta+1}\frac{e^{-\frac{z}{\mu}\mathbf{1}}}{\mu^{\delta+2}%
}\left\vert \sigma\left(  \mu\right)  \right\vert d\mu\text{ }dz\\
& <\int_{\mathcal{O}_{\mathbf{1}}}\int_{a}^{b}z^{\delta+1}\frac{e^{-\frac
{z}{b}\mathbf{1}}}{a^{\delta+2}}\left\vert \sigma\left(  \mu\right)
\right\vert d\mu\text{ }dz\\
& <\infty,
\end{align*}

which allows us to use Fubini's theorem to justify changing the order of
integration. Indeed, for all multi-integers $\delta\geq-\mathbf{1}$,%
\begin{align*}
\widehat{J_{\delta}\left[  \sigma\right]  }\left(  \xi\right)   & =\left(
2\pi\right)  ^{-d/2}\int_{\mathcal{O}_{\mathbf{1}}}\int_{a}^{b}e^{-i\xi
z}z^{\delta+1}\frac{e^{-\frac{z}{\mu}\mathbf{1}}}{\mu^{\delta+2}}\sigma\left(
\mu\right)  d\mu\text{ }dz\\
& =\left(  2\pi\right)  ^{-d/2}\int_{a}^{b}\int_{\mathcal{O}_{\mathbf{1}}%
}e^{-i\xi z}z^{\delta+1}e^{-\frac{z}{\mu}\mathbf{1}}dz\text{ }\frac
{\sigma\left(  \mu\right)  }{\mu^{\delta+2}}d\mu\\
& \Rightarrow\zeta=z/\mu,\text{ }z=\mu.\zeta,\text{ }dz=\mu^{\mathbf{1}}%
d\zeta\Rightarrow\\
\widehat{J_{\delta}\left[  \sigma\right]  }\left(  \xi\right)   & =\left(
2\pi\right)  ^{-d/2}\int_{a}^{b}\int_{\mathcal{O}_{\mathbf{1}}}e^{-i\xi\left(
\mu.\zeta\right)  }\left(  \mu.\zeta\right)  ^{\delta+1}e^{-\zeta\mathbf{1}%
}\mu^{\mathbf{1}}d\zeta\text{ }\frac{\sigma\left(  \mu\right)  }{\mu
^{\delta+2}}d\mu\\
& =\left(  2\pi\right)  ^{-d/2}\int_{a}^{b}\mu^{\delta+1}\mu^{\mathbf{1}}%
\int_{\mathcal{O}_{\mathbf{1}}}e^{-i\left(  \xi.\zeta\right)  \mu}%
\zeta^{\delta+1}e^{-\zeta\mathbf{1}}d\zeta\text{ }\frac{\sigma\left(
\mu\right)  }{\mu^{\delta+2}}d\mu\\
& =\left(  2\pi\right)  ^{-d/2}\int_{a}^{b}\left(  \int_{\mathcal{O}%
_{\mathbf{1}}}e^{-i\left(  \xi.\zeta\right)  \mu}\zeta^{\delta+1}%
e^{-\zeta\mathbf{1}}d\zeta\right)  \sigma\left(  \mu\right)  d\mu\\
& \Rightarrow clearly\text{ }absolutely\text{ }convergent\Rightarrow\\
& =\left(  2\pi\right)  ^{-d/2}\int_{\mathcal{O}_{\mathbf{1}}}\left(  \int%
_{a}^{b}e^{-i\left(  \xi.\zeta\right)  \mu}\sigma\left(  \mu\right)
d\mu\right)  \zeta^{\delta+1}e^{-\zeta\mathbf{1}}d\zeta\\
& =\int_{\mathcal{O}_{\mathbf{1}}}\widehat{\sigma}\left(  \xi.\zeta\right)
\zeta^{\delta+1}e^{-\zeta\mathbf{1}}d\zeta,
\end{align*}

which proves \ref{X17}.

If $\beta\leq\delta+1$, from \ref{X17},
\begin{align*}
\xi^{\beta}\widehat{J_{\delta}\left[  \sigma\right]  }\left(  \xi\right)   &
=\int_{\mathcal{O}_{\mathbf{1}}}\widehat{\sigma}\left(  \xi.\zeta\right)
\xi^{\beta}\zeta^{\delta+1}e^{-\zeta\mathbf{1}}d\zeta\\
& =\int_{\mathcal{O}_{\mathbf{1}}}\widehat{\sigma}\left(  \xi.\zeta\right)
\xi^{\beta}\zeta^{\beta}\zeta^{\delta+1-\beta}e^{-\zeta\mathbf{1}}d\zeta\\
& =\int_{\mathcal{O}_{\mathbf{1}}}\widehat{\sigma}\left(  \xi.\zeta\right)
\left(  \xi.\zeta\right)  ^{\beta}\zeta^{\delta+1-\beta}e^{-\zeta\mathbf{1}%
}d\zeta\\
& =\left(  -i\right)  ^{\left\vert \beta\right\vert }\int_{\mathcal{O}%
_{\mathbf{1}}}\widehat{D^{\beta}\sigma}\left(  \xi.\zeta\right)  \zeta
^{\delta+1-\beta}e^{-\zeta\mathbf{1}}d\zeta.
\end{align*}

If $\xi\in\mathcal{O}_{\theta}$, the change of variables $\chi=\xi.\zeta$
yields%
\begin{align}
\left\vert \xi^{\beta}\widehat{J_{\delta}\left[  \sigma\right]  }\left(
\xi\right)  \right\vert  & \leq\int_{\mathcal{O}_{\mathbf{1}}}\left\vert
\widehat{D^{\beta}\sigma}\left(  \xi.\zeta\right)  \right\vert \zeta
^{\delta+1-\beta}e^{-\zeta\mathbf{1}}d\zeta\label{X61}\\
& \leq\left\Vert \zeta^{\delta+1-\beta}e^{-\zeta\mathbf{1}}\right\Vert
_{\infty,\overline{\mathcal{O}_{\mathbf{1}}}}\int_{\mathcal{O}_{\mathbf{1}}%
}\left\vert \widehat{D^{\beta}\sigma}\left(  \xi.\zeta\right)  \right\vert
d\zeta\nonumber\\
& =\left\Vert \zeta^{\delta+1-\beta}e^{-\zeta\mathbf{1}}\right\Vert
_{\infty,\overline{\mathcal{O}_{\mathbf{1}}}}\frac{1}{\left\vert
\xi^{\mathbf{1}}\right\vert }\int_{\mathcal{O}_{\mathbf{1}}}\left\vert
\widehat{D^{\beta}\sigma}\left(  \chi\right)  \right\vert d\chi,\nonumber
\end{align}

so that%
\begin{align*}
\left\Vert \xi^{\beta+1}\widehat{J_{\delta}\left[  \sigma\right]  }\right\Vert
_{\infty}  & \leq\left\Vert \zeta^{\delta+1-\beta}e^{-\zeta\mathbf{1}%
}\right\Vert _{\infty,\overline{\mathcal{O}_{\mathbf{1}}}}\left\Vert
\widehat{D^{\beta}\sigma}\right\Vert _{1}\\
& \leq\left(  \delta+1-\beta\right)  ^{\delta+1-\beta}e^{-\left(
\delta+1-\beta\right)  \mathbf{1}}\left\Vert \widehat{D^{\beta}\sigma
}\right\Vert _{1},\quad\mathbf{0}\leq\beta\leq\delta+1,
\end{align*}

or equivalently \textbf{the upper bound}: when $\widehat{D^{\gamma}\sigma}\in
L^{1}$,%
\begin{equation}
\left\Vert \xi^{\gamma}\widehat{J_{\delta}\left[  \sigma\right]  }\right\Vert
_{\infty}\leq\left(  \delta+2-\gamma\right)  ^{\delta+2-\gamma}e^{-\left(
\delta+2-\gamma\right)  \mathbf{1}}\left\Vert \widehat{D^{\gamma-1}\sigma
}\right\Vert _{1},\quad\mathbf{1}\leq\gamma\leq\delta+2.\label{X63}%
\end{equation}

Alternatively, starting with \ref{X61}, if $\xi\in\mathcal{O}_{\theta}$,%
\begin{align*}
\left\vert \xi^{\beta}\widehat{J_{\delta}\left[  \sigma\right]  }\left(
\xi\right)  \right\vert  & \leq\int_{\mathcal{O}_{\mathbf{1}}}\left\vert
\widehat{D^{\beta}\sigma}\left(  \xi.\zeta\right)  \right\vert \zeta
^{\delta+1-\beta}e^{-\zeta\mathbf{1}}d\zeta\\
& \leq\left\Vert \widehat{D^{\beta}\sigma}\right\Vert _{\infty,\overline
{\mathcal{O}_{\theta}}}\int_{\mathcal{O}_{\mathbf{1}}}\zeta^{\delta+1-\beta
}e^{-\zeta\mathbf{1}}d\zeta\\
& =\left(  \delta+1-\beta\right)  !\left\Vert \widehat{D^{\beta}\sigma
}\right\Vert _{\infty,\overline{\mathcal{O}_{\theta}}}\\
& \leq\tfrac{\left(  \delta+1-\beta\right)  !}{\left(  2\pi\right)  ^{d/2}%
}\left\Vert D^{\beta}\sigma\right\Vert _{1},\quad\beta\leq\delta+1,
\end{align*}

and we have \textbf{the two upper bounds}: when $\widehat{D^{\beta}\sigma}\in
L^{\infty}$,%
\begin{equation}
\left\Vert \xi^{\beta}\widehat{J_{\delta}\left[  \sigma\right]  }\right\Vert
_{\infty}\leq\left(  \delta+1-\beta\right)  !\left\Vert \widehat{D^{\beta
}\sigma}\right\Vert _{\infty},\quad\beta\leq\delta+1,\text{ }\delta
\geq-\mathbf{1};\label{X64}%
\end{equation}

and in addition, if $D^{\beta}\sigma\in L^{1}$,%
\begin{equation}
\left\Vert \xi^{\beta}\widehat{J_{\delta}\left[  \sigma\right]  }\right\Vert
_{\infty}\leq\tfrac{\left(  \delta+1-\beta\right)  !}{\left(  2\pi\right)
^{d/2}}\left\Vert D^{\beta}\sigma\right\Vert _{1},\quad\beta\leq
\delta+1,\text{ }\delta\geq-\mathbf{1}.\label{X62}%
\end{equation}
\medskip

\fbox{Case: $\gamma\leq\mathbf{1}$ and $\delta=\mathbf{0}$.} From \ref{X17},
if $\widehat{D^{\gamma}\sigma}\in L^{\infty}$,\smallskip%
\begin{align*}
\xi^{\gamma}\widehat{J_{\mathbf{0}}\left[  \sigma\right]  }\left(  \xi\right)
=\int_{\mathcal{O}_{\mathbf{1}}}\widehat{\sigma}\left(  \xi.\zeta\right)
\xi^{\gamma}\zeta^{1}e^{-\zeta\mathbf{1}}d\zeta &  =\int_{\mathcal{O}%
_{\mathbf{1}}}\widehat{\sigma}\left(  \xi.\zeta\right)  \left(  \xi
.\zeta\right)  ^{\gamma}\zeta^{1-\gamma}e^{-\zeta\mathbf{1}}d\zeta=\\
&  =\left(  -1\right)  ^{\left\vert \gamma\right\vert }\int_{\mathcal{O}%
_{\mathbf{1}}}\widehat{D^{\gamma}\sigma}\left(  \xi.\zeta\right)
\zeta^{1-\gamma}e^{-\zeta\mathbf{1}}d\zeta,
\end{align*}

so that if $\xi\in\mathcal{O}_{\theta}$,%
\begin{align*}
\left\vert \xi^{\gamma}\widehat{J_{\mathbf{0}}\left[  \sigma\right]  }\left(
\xi\right)  \right\vert \leq\int_{\mathcal{O}_{\mathbf{1}}}\left\vert
\widehat{D^{\gamma}\sigma}\left(  \xi.\zeta\right)  \right\vert \zeta
^{1-\gamma}e^{-\zeta\mathbf{1}}d\zeta &  \leq\left\Vert \widehat{D^{\gamma
}\sigma}\right\Vert _{\infty,\overline{\mathcal{O}_{\theta}}}\int%
_{\mathcal{O}_{\mathbf{1}}}\zeta^{1-\gamma}e^{-\zeta\mathbf{1}}d\zeta\\
&  =\left(  1-\gamma\right)  !\left\Vert \widehat{D^{\gamma}\sigma}\right\Vert
_{\infty,\overline{\mathcal{O}_{\theta}}}\\
&  \leq\left(  1-\gamma\right)  !\left\Vert \widehat{D^{\gamma}\sigma
}\right\Vert _{\infty}\\
&  \leq\left\Vert D^{\gamma}\sigma\right\Vert _{1}.
\end{align*}

Hence: when $\widehat{D^{\gamma}\sigma}\in L^{\infty}$,%
\begin{equation}
\left\Vert \xi^{\gamma}\widehat{J_{\mathbf{0}}\left[  \sigma\right]
}\right\Vert _{\infty}\leq\left\Vert \widehat{D^{\gamma}\sigma}\right\Vert
_{\infty},\quad\gamma\leq\mathbf{1},\label{X65}%
\end{equation}

and if, in addition, $D^{\gamma}\sigma\in L^{1}$ then%
\begin{equation}
\left\Vert \xi^{\gamma}\widehat{J_{\mathbf{0}}\left[  \sigma\right]
}\right\Vert _{\infty}\leq\left(  2\pi\right)  ^{-d/2}\left\Vert D^{\gamma
}\sigma\right\Vert _{1},\quad\gamma\leq\mathbf{1}.\label{X66}%
\end{equation}
\medskip

\fbox{Case: $\gamma<\mathbf{1}$ and $\delta=-\mathbf{1}$.} If $\gamma
=\mathbf{0}$ then%
\begin{equation}
\left\Vert \widehat{J_{-1}\left[  \sigma\right]  }\right\Vert _{\infty}%
\leq\left\Vert \widehat{\sigma}\right\Vert _{\infty}.\label{X69}%
\end{equation}

Now suppose $\mathbf{0}<\gamma<\mathbf{1}$. Choose a permutation $\pi_{\gamma
}$ such that $\pi_{\gamma}\gamma=$ $\left(  \mathbf{1}^{\prime},\mathbf{0}%
^{\prime\prime}\right)  $ i.e. a sequence of ones then some zeros. Denote the
corresponding partition of variable by $\left(  \cdot^{\prime},\cdot
^{\prime\prime}\right)  $.Then%
\[
\left(  \pi_{\gamma}\xi\right)  ^{\gamma}\widehat{J_{\delta}\left[
\sigma\right]  }\left(  \pi_{\gamma}\xi\right)  =\left(  \xi^{\prime}\right)
^{\mathbf{1}^{\prime}}\widehat{J_{\delta}\left[  \sigma\right]  }\left(
\pi_{\gamma}\xi\right)  ,
\]

and%
\begin{align*}
\widehat{J_{\delta}\left[  \sigma\right]  }\left(  \pi_{\gamma}\xi\right)   &
=\int_{\mathcal{O}_{\mathbf{1}}}\widehat{\sigma}\left(  \left(  \pi_{\gamma
}\xi\right)  .\zeta\right)  \zeta^{\delta+1}e^{-\zeta\mathbf{1}}d\zeta\\
& \Rightarrow\chi=\pi_{\gamma}\zeta\Rightarrow\\
& =\int_{\mathcal{O}_{\mathbf{1}}}\widehat{\sigma}\left(  \left(  \pi_{\gamma
}\xi\right)  .\left(  \pi_{\gamma}\chi\right)  \right)  \left(  \pi_{\gamma
}\chi\right)  ^{\delta+1}e^{-\chi\mathbf{1}}d\chi\\
& =\int_{\mathcal{O}_{\mathbf{1}}}\widehat{\sigma\pi_{\gamma}}\left(  \xi
.\chi\right)  \chi^{\pi_{\gamma}\delta+1}e^{-\chi\mathbf{1}}d\chi,
\end{align*}

so that%
\[
\left(  \pi_{\gamma}\xi\right)  ^{\gamma}\widehat{J_{\pi_{\gamma}\delta
}\left[  \sigma\pi_{\gamma}\right]  }\left(  \pi_{\gamma}\xi\right)
=\int_{\mathcal{O}_{\mathbf{1}}}\widehat{\sigma}\left(  \xi.\chi\right)
\left(  \xi^{\prime}\right)  ^{\mathbf{1}^{\prime}}\chi^{\delta+1}%
e^{-\chi\mathbf{1}}d\chi.
\]

\[
\left(  \pi_{\gamma}\xi\right)  ^{\gamma}\widehat{J_{-\mathbf{1}}\left[
\sigma\pi_{\gamma}\right]  }\left(  \pi_{\gamma}\xi\right)  =\int%
_{\mathcal{O}_{\mathbf{1}^{\prime\prime}}}\int_{\mathcal{O}_{\mathbf{1}%
^{\prime}}}\widehat{\sigma}\left(  \xi.\chi\right)  \left(  \xi^{\prime
}\right)  ^{\mathbf{1}^{\prime}}e^{-\chi\mathbf{1}}d\chi.
\]

Suppose $\chi^{\prime}\in\mathcal{O}_{\theta}$. Now apply the change of
variables: $\omega^{\prime}=\xi^{\prime}.\chi^{\prime}$, $d\omega^{\prime
}=\left\vert \left(  \xi^{\prime}\right)  ^{\mathbf{1}^{\prime}}\right\vert
d\chi^{\prime}$, to yield%
\[
\left(  \pi_{\gamma}\xi\right)  ^{\gamma}\widehat{J_{-\mathbf{1}}\left[
\sigma\pi_{\gamma}\right]  }\left(  \pi_{\gamma}\xi\right)  =\int%
_{\mathcal{O}_{\mathbf{1}^{\prime\prime}}}\int_{\mathcal{O}_{\theta^{\prime}}%
}\widehat{\sigma}\left(  \omega^{\prime},\xi^{\prime\prime}.\chi^{\prime
\prime}\right)  \left(  \xi^{\prime}\right)  ^{\mathbf{1}^{\prime}}%
e^{-\frac{\omega^{\prime}}{\xi^{\prime}}\mathbf{1}^{\prime}}e^{-\chi
^{\prime\prime}\mathbf{1}^{\prime\prime}}\frac{d\omega^{\prime}}{\left\vert
\left(  \xi^{\prime}\right)  ^{\mathbf{1}^{\prime}}\right\vert }d\chi
^{\prime\prime},
\]

which implies%
\begin{align*}
\left\vert \left(  \pi_{\gamma}\xi\right)  ^{\gamma}\widehat{J_{-\mathbf{1}%
}\left[  \sigma\pi_{\gamma}\right]  }\left(  \pi_{\gamma}\xi\right)
\right\vert  & \leq\int_{\mathcal{O}_{\mathbf{1}^{\prime\prime}}}%
\int_{\mathcal{O}_{\theta^{\prime}}}\left\vert \widehat{\sigma}\left(
\omega^{\prime},\xi^{\prime\prime}.\chi^{\prime\prime}\right)  \right\vert
d\omega^{\prime}e^{-\chi^{\prime\prime}\mathbf{1}^{\prime\prime}}d\chi
^{\prime\prime}\\
& \leq\int_{\mathcal{O}_{\mathbf{1}^{\prime\prime}}}\left(  \int%
_{\mathbb{R}^{d^{\prime}}}\left\vert \widehat{\sigma}\left(  \omega^{\prime
},\xi^{\prime\prime}.\chi^{\prime\prime}\right)  \right\vert d\omega^{\prime
}\right)  e^{-\chi^{\prime\prime}\mathbf{1}^{\prime\prime}}d\chi^{\prime
\prime},
\end{align*}

and assumption \ref{X56} now implies that for all $\xi$,%
\begin{align*}
\left\vert \left(  \pi_{\gamma}\xi\right)  ^{\gamma}\widehat{J_{-\mathbf{1}%
}\left[  \sigma\pi_{\gamma}\right]  }\left(  \pi_{\gamma}\xi\right)
\right\vert  & \leq\left\Vert \int\left\vert \widehat{\sigma}\left(
\omega\right)  \right\vert d^{\pi_{\gamma}\gamma}\omega\right\Vert _{\infty
}\int_{\mathcal{O}_{\mathbf{1}^{\prime\prime}}}e^{-\chi^{\prime\prime
}\mathbf{1}^{\prime\prime}}d\chi^{\prime\prime}\\
& =\left\Vert \int\left\vert \widehat{\sigma}\left(  \omega\right)
\right\vert d^{\pi_{\gamma}\gamma}\omega\right\Vert _{\infty}.
\end{align*}

i.e.%
\[
\left\vert \xi^{\gamma}\widehat{J_{-\mathbf{1}}\left[  \sigma\right]  }\left(
\xi\right)  \right\vert \leq\left\Vert \int\left\vert \widehat{\sigma
\pi_{\gamma}}\left(  \omega\right)  \right\vert d^{\pi_{\gamma}\gamma}%
\omega\right\Vert _{\infty},\quad\gamma<\mathbf{1}.
\]

But%
\[
\left\Vert \int\left\vert \widehat{\sigma\pi_{\gamma}}\left(  \omega\right)
\right\vert d^{\pi_{\gamma}\gamma}\omega\right\Vert _{\infty}=\left\Vert
\int\left\vert \widehat{\sigma}\left(  \pi_{\gamma}\omega\right)  \right\vert
d^{\pi_{\gamma}\gamma}\omega\right\Vert _{\infty}=\left\Vert \int\left\vert
\widehat{\sigma}\left(  \omega\right)  \right\vert d^{\gamma}\omega\right\Vert
_{\infty},
\]

so, on noting Notation \ref{Not_partial_integ} and \ref{X69},%
\begin{equation}
\left\Vert \xi^{\gamma}\widehat{J_{-\mathbf{1}}\left[  \sigma\right]
}\right\Vert _{\infty}\leq\left\Vert \int\left\vert \widehat{\sigma}\left(
\omega\right)  \right\vert d^{\gamma}\omega\right\Vert _{\infty},\quad
\gamma\leq\mathbf{1}.\label{X67}%
\end{equation}
\medskip

\textbf{Part 2} Since $J_{\delta}^{\theta}\left[  \sigma\right]  \left(
z\right)  =J_{\delta}\left[  \sigma\left(  \theta.\right)  \right]  \left(
\theta.z\right)  $, from \ref{X17},%
\[
\widehat{J_{\delta}\left[  \sigma\right]  }\left(  \xi\right)  =\int%
_{\mathcal{O}_{\mathbf{1}}}\widehat{\sigma}\left(  \xi.\zeta\right)
\zeta^{\delta+1}e^{-\zeta\mathbf{1}}d\zeta.
\]

Part 13 of Definition \ref{Def_Fourier} implies $\widehat{g\left(
\theta.\right)  }=\widehat{g}\left(  \theta.\right)  $ so that%
\begin{align*}
\widehat{J_{\delta}^{\theta}\left[  \sigma\right]  }\left(  \xi\right)   &
=F_{z}\left[  J_{\delta}\left[  \sigma\left(  \theta.\right)  \right]  \left(
\theta.z\right)  \right]  \left(  \xi\right) \\
& =F_{z}\left[  J_{\delta}\left[  \sigma\left(  \theta.\right)  \right]
\right]  \left(  \theta.\xi\right) \\
& =\int_{\mathcal{O}_{\mathbf{1}}}\widehat{\sigma\left(  \theta.\right)
}\left(  \theta.\xi.\zeta\right)  \zeta^{\delta+1}e^{-\zeta\mathbf{1}}d\zeta\\
& =\int_{\mathcal{O}_{\mathbf{1}}}\widehat{\sigma}\left(  \xi.\zeta\right)
\zeta^{\delta+1}e^{-\zeta\mathbf{1}}d\zeta.
\end{align*}

The estimates for $\left\Vert \xi^{\beta}\widehat{J_{\delta}^{\theta}\left[
\sigma\right]  }\right\Vert _{\infty}$ now follow easily.
\end{proof}

\begin{lemma}
If $\delta\geq-\mathbf{1}$ and $\sigma\in L_{0}^{1}\left(  \mathcal{O}%
_{\mathbf{1}}\right)  $ then for all multi-indexes $\mathbf{0}\leq\alpha
\leq\delta+2$,%
\[
\left\vert \xi^{\alpha}\widehat{J_{\delta}\left[  \sigma\right]  }\left(
\xi\right)  \right\vert \leq\left\Vert \frac{\sigma}{y^{\alpha}}\right\Vert
_{1,\mathcal{O}_{\mathbf{1}}}\left\Vert x^{\alpha}g_{\delta+1}\right\Vert
_{\infty},
\]

where%
\[
g_{\delta}:=\left(  H\left(  \zeta\right)  \zeta^{\delta}e^{-\zeta\mathbf{1}%
}\right)  ^{\vee},\quad\delta\geq\mathbf{0},
\]

and $H$ is the tensor product Heavyside step function. The function
$g_{\delta}$ has the following properties:

\begin{enumerate}
\item $\widehat{g_{\delta}}\in C_{B}^{\infty}\left(  \overline{\mathcal{O}%
_{\mathbf{1}}}\right)  $, $\widehat{g_{\delta}}\in C_{B}^{\left(
\delta-1\right)  }$ $\delta\geq\mathbf{1}$, $\widehat{g_{\delta}}\in L^{p}$
for $1\leq p\leq\infty$.

\item $g_{\delta}\in C_{B}^{\infty}\cap H^{\infty}$.

\item $x^{\alpha}g_{\delta}\left(  x\right)  \in C_{B}^{\left(  0\right)  }$
for all $\alpha\leq\delta+1$.

\item $\left\vert g_{\delta}\left(  x\right)  \right\vert \leq\frac{C}{\left(
1+x_{+}\right)  ^{\delta+1}}=\frac{C}{\left(  1+\left\vert x_{1}\right\vert
\right)  ^{\delta_{1}+1}\ldots\left(  1+\left\vert x_{d}\right\vert \right)
^{\delta_{d}+1}}$.
\end{enumerate}
\end{lemma}

\begin{proof}
?? Since $\sigma\in L_{0}^{1}\subset L_{0}^{2}$ and $H\left(  \zeta\right)
\zeta^{\delta+1}e^{-\zeta\mathbf{1}}\in L^{2}\cap L^{1}$, \ref{X17} holds and
\begin{align*}
\widehat{J_{\delta}\left[  \sigma\right]  }\left(  \xi\right)   &
=\int_{\mathcal{O}_{\mathbf{1}}}\widehat{\sigma}\left(  \xi.\zeta\right)
\zeta^{\delta+1}e^{-\zeta\mathbf{1}}d\zeta\\
& =\int\widehat{\sigma}\left(  \xi.\zeta\right)  H\left(  \zeta\right)
\zeta^{\delta+1}e^{-\zeta\mathbf{1}}d\zeta\\
& =\frac{1}{\left\vert \xi^{\mathbf{1}}\right\vert }\int\widehat{\sigma\left(
\frac{\cdot}{\xi}\right)  }\left(  \zeta\right)  H\left(  \zeta\right)
\zeta^{\delta+1}e^{-\zeta\mathbf{1}}d\zeta\\
& =\frac{1}{\left\vert \xi^{\mathbf{1}}\right\vert }\int\sigma\left(  \frac
{x}{\xi}\right)  \left(  H\left(  \zeta\right)  \zeta^{\delta+1}%
e^{-\zeta\mathbf{1}}\right)  ^{\vee}\left(  x\right)  dx\\
& =\int_{\mathcal{O}_{\mathbf{1}}}\sigma\left(  y\right)  \left(  H\left(
\zeta\right)  \zeta^{\delta+1}e^{-\zeta\mathbf{1}}\right)  ^{\vee}\left(
y.\xi\right)  dy.
\end{align*}

Further%
\begin{align}
\xi^{\alpha}\widehat{J_{\delta}\left[  \sigma\right]  }\left(  \xi\right)   &
=\int_{\mathcal{O}_{\mathbf{1}}}\sigma\left(  y\right)  \xi^{\alpha}\left(
H\left(  \zeta\right)  \zeta^{\delta+1}e^{-\zeta\mathbf{1}}\right)  ^{\vee
}\left(  y.\xi\right)  dy\nonumber\\
& =\int_{\mathcal{O}_{\mathbf{1}}}\frac{\sigma\left(  y\right)  }{y^{\alpha}%
}\left(  y.\xi\right)  ^{\alpha}\left(  H\left(  \zeta\right)  \zeta
^{\delta+1}e^{-\zeta\mathbf{1}}\right)  ^{\vee}\left(  y.\xi\right)
dy\nonumber\\
& =\int_{\mathcal{O}_{\mathbf{1}}}\frac{\sigma\left(  y\right)  }{y^{\alpha}%
}\left(  y.\xi\right)  ^{\alpha}g_{\delta+1}\left(  y.\xi\right)
dy.\label{X44}%
\end{align}

\[
x^{\alpha}g_{\delta}\left(  x\right)  =x^{\alpha}\left(  H\left(
\zeta\right)  \zeta^{\delta}e^{-\zeta\mathbf{1}}\right)  ^{\vee}=\left(
D^{\alpha}\left(  H\left(  \zeta\right)  \zeta^{\delta}e^{-\zeta\mathbf{1}%
}\right)  \right)  ^{\vee}.
\]

??=================

When $s>0$, $D^{??\alpha}\left(  H\left(  s\right)  s^{\delta??}e^{-s}\right)
=D^{\alpha}\left(  s^{\delta}e^{-s}\right)  =D^{\alpha}\left(  s^{\delta
}e^{-s}\right)  =e^{-s}e^{s}D^{\alpha}\left(  s^{\delta}e^{-s}\right)
=p_{n}^{m-1}\left(  s\right)  e^{-s}.$

??==================

?? USE $p_{n}^{m-1}$! Noting that $x^{\alpha}g_{\delta}\left(  x\right)  $ is
a tensor product function we first consider one dimension. Since $H\left(
s\right)  s^{k}\in C_{B}^{\left(  k-1\right)  }$ when $k\geq1$, $D^{k}\left(
H\left(  s\right)  s^{k}\right)  =k!H\left(  s\right)  $ when $k\geq0$, and
$D^{k+1}\left(  H\left(  s\right)  s^{k}\right)  =k!\delta$ when $k\geq0$, we have:

$D^{j}\left(  H\left(  s\right)  s^{k}e^{-s}\right)  \in C_{B}^{\left(
0\right)  }\cap L^{1}$ when $j\leq k-1$, $D^{k}\left(  H\left(  s\right)
s^{k}e^{-s}\right)  \in\left(  k!H\left(  s\right)  +C_{B}^{\left(
k-1\right)  }\right)  e^{-s}\subset L^{1}$ and

$D^{k+1}\left(  H\left(  s\right)  s^{k}e^{-s}\right)  \in k!\delta+\left(
cH+C_{B}^{\left(  k-1\right)  }\right)  e^{-s}\subset k!\delta+L^{1}$.

Thus $\left(  D^{j}\left(  H\left(  s\right)  s^{k}e^{-s}\right)  \right)
^{\vee}\in C_{B}^{\left(  0\right)  }$ for $j\leq k+1$ i.e. $t^{j}\left(
H\left(  s\right)  s^{k}e^{-s}\right)  ^{\vee}\left(  t\right)  \in
C_{B}^{\left(  0\right)  }$ when $j\leq k+1$.

In $d$ dimensions we can now conclude that $x^{\alpha}\left(  H\left(
\xi\right)  \xi^{\delta}e^{-\xi\mathbf{1}}\right)  ^{\vee}\in C_{B}^{\left(
0\right)  }$ when $\alpha\leq\delta+1$ \textbf{which proves part 3 and hence
part 4}.

Returning to \ref{X44} : when $\alpha\leq\delta+2$,%
\begin{align*}
\left\vert \xi^{\alpha}\widehat{J_{\delta}\left[  \sigma\right]  }\left(
\xi\right)  \right\vert  & \leq\int_{\mathcal{O}_{\mathbf{1}}}\frac{\left\vert
\sigma\left(  y\right)  \right\vert }{\left\vert y^{\alpha}\right\vert
}\left\vert \left(  y.\xi\right)  ^{\alpha}g_{\delta+1}\left(  y.\xi\right)
\right\vert dy\\
& \leq\left\Vert \frac{\sigma}{y^{\alpha}}\right\Vert _{1,\mathcal{O}%
_{\mathbf{1}}}\left\Vert x^{\alpha}g_{\delta+1}\right\Vert _{\infty}.
\end{align*}

\textbf{Part 2} $\forall\alpha:D^{\alpha}g_{\delta}\in\left(  L^{1}\right)
^{\vee}\in C_{B}^{\left(  0\right)  }\Rightarrow g_{\delta}\in C_{B}^{\infty}%
$; $\forall\alpha:D^{\alpha}g\in\left(  L^{2}\right)  ^{\vee}=L^{2}\Rightarrow
g\in H^{\infty}$.
\end{proof}

\begin{lemma}
\label{Lem_bnd_FourTran_Jn1[v]}Bounds for $\left\Vert \xi^{\gamma
}\widehat{J_{n\mathbf{1}}^{\theta}\left[  \sigma\right]  }\right\Vert
_{\infty}$ when $n\geq-1$:

\begin{enumerate}
\item \fbox{$n=-1$}%
\[
\left\Vert \xi^{\gamma}\widehat{J_{-\mathbf{1}}^{\theta}\left[  \sigma\right]
}\right\Vert _{\infty}\leq\left\Vert \int\left\vert \widehat{\sigma}\left(
\omega\right)  \right\vert d^{\gamma}\omega\right\Vert _{\infty},\quad
\gamma\leq\mathbf{1},
\]

\item \fbox{$n=0$}%
\[
\left\Vert \xi^{\gamma}\widehat{J_{\mathbf{0}}^{\theta}\left[  \sigma\right]
}\right\Vert _{\infty}\leq\left\{
\begin{array}
[c]{ll}%
\left\Vert \widehat{D^{\gamma}\sigma}\right\Vert _{\infty}, & \gamma
\leq\mathbf{1},\\
e^{-\left(  2-\gamma\right)  \mathbf{1}}\left\Vert \widehat{D^{\gamma}\sigma
}\right\Vert _{1}, & \mathbf{1}\leq\gamma\leq2\mathbf{1}.
\end{array}
\right.
\]

Note that not all the derivatives bounds for $\gamma\leq2\mathbf{1}$ are given here.

\item \fbox{$n\geq1$}%
\begin{equation}
\left\Vert \xi^{\gamma}\widehat{J_{n\mathbf{1}}^{\theta}\left[  \sigma\right]
}\right\Vert _{\infty}\leq\left\{
\begin{array}
[c]{ll}%
\left(  \left(  n+1\right)  \mathbf{1}-\gamma\right)  !\left\Vert
\widehat{D^{\gamma}\sigma}\right\Vert _{\infty}, & \mathbf{0}\leq\gamma
\leq\left(  n+1\right)  \mathbf{1},\\
\left(  \left(  n+2\right)  \mathbf{1}-\gamma\right)  ^{\left(  n+2\right)
\mathbf{1}-\gamma}e^{-\left(  \left(  n+2\right)  \mathbf{1}-\gamma\right)
}\left\Vert \widehat{D^{\gamma}\sigma}\right\Vert _{1}, & \mathbf{1}\leq
\gamma\leq\left(  n+2\right)  \mathbf{1}.
\end{array}
\right. \label{X60}%
\end{equation}

Note that not all the derivative bounds for $\gamma\leq\left(  n+2\right)
\mathbf{1}$ are given here.
\end{enumerate}

NOTE: Parts 2 and 3 can be amalgamated i.e. the estimate of part 3 holds for
$n=0$. Also, in \ref{X60} of part 3 the ranges of $\gamma$ in the two
estimates can overlap substantially.
\end{lemma}

\begin{proof}
\textbf{Part 1} From Lemma \ref{Lem_bnd_FourTran_J[v]}: when $\gamma
=\mathbf{1}$ we have the upper bound \ref{X63} i.e. when $\widehat{D^{\gamma
}\sigma}\in L^{1}$,%
\[
\left\Vert \xi^{\gamma}\widehat{J_{\delta}^{\theta}\left[  \sigma\right]
}\right\Vert _{\infty}\leq\left(  \delta+2-\gamma\right)  ^{\delta+2-\gamma
}e^{-\left(  \delta+2-\gamma\right)  \mathbf{1}}\left\Vert \widehat{D^{\gamma
-1}\sigma}\right\Vert _{1},\quad\left\{
\begin{array}
[c]{l}%
\mathbf{1}\leq\gamma\leq\delta+2,\\
\delta\geq-\mathbf{1},
\end{array}
\right.
\]

so that when $\delta=-\mathbf{1}$ and $\gamma=\mathbf{1}$,`%
\[
\left\Vert \xi^{\mathbf{1}}\widehat{J_{-\mathbf{1}}^{\theta}\left[
\sigma\right]  }\right\Vert _{\infty}\leq\left\Vert \widehat{\sigma
}\right\Vert _{1}.
\]

Also, from \ref{X67},%
\[
\left\Vert \xi^{\gamma}\widehat{J_{-\mathbf{1}}\left[  \sigma\right]
}\right\Vert _{\infty}\leq\left\Vert \int\left\vert \widehat{\sigma}\left(
\omega\right)  \right\vert d^{\gamma}\omega\right\Vert _{\infty},\quad
\gamma\leq\mathbf{1}.
\]
\medskip

\textbf{Part 2} From part 2 of Lemma \ref{Lem_bnd_FourTran_J[v]}: when
$\delta=\mathbf{0}$ and $\gamma<\mathbf{1}$, we have the upper bound \ref{X65}
for $\left\Vert \xi^{\gamma}\widehat{J_{\mathbf{0}}^{\theta}\left[
\sigma\right]  }\right\Vert _{\infty}$ i.e. if $\widehat{D^{\gamma}\sigma}\in
L^{\infty}$ for $\gamma<\mathbf{1}$ then%
\[
\left\Vert \xi^{\gamma}\widehat{J_{\mathbf{0}}^{\theta}\left[  \sigma\right]
}\right\Vert _{\infty}\leq\left\Vert \widehat{D^{\gamma}\sigma}\right\Vert
_{\infty},\quad\gamma<\mathbf{1}.
\]

Also from \ref{X63} of part 2 of Lemma \ref{Lem_bnd_FourTran_J[v]}: if
$\widehat{D^{\gamma}\sigma}\in L^{1}$ for $\mathbf{1}\leq\gamma\leq\delta+2$
then%
\[
\left\Vert \xi^{\gamma}\widehat{J_{\delta}^{\theta}\left[  \sigma\right]
}\right\Vert _{\infty}\leq\left(  \delta+2-\gamma\right)  ^{\delta+2-\gamma
}e^{-\left(  \delta+2-\gamma\right)  \mathbf{1}}\left\Vert \widehat{D^{\gamma
}\sigma}\right\Vert _{1},\quad\mathbf{1}\leq\gamma\leq\delta+2,\text{ }%
\delta\geq-\mathbf{1},
\]

i.e. if $\widehat{D^{\gamma}\sigma}\in L^{1}$ for $\mathbf{1}\leq\gamma
\leq2\mathbf{1}$ then%
\[
\left\Vert \xi^{\gamma}\widehat{J_{\mathbf{0}}^{\theta}\left[  \sigma\right]
}\right\Vert _{\infty}\leq\left(  2-\gamma\right)  ^{2-\gamma}e^{-\left(
2-\gamma\right)  \mathbf{1}}\left\Vert \widehat{D^{\gamma}\sigma}\right\Vert
_{1}=e^{-\left(  2-\gamma\right)  \mathbf{1}}\left\Vert \widehat{D^{\gamma
}\sigma}\right\Vert _{1}.
\]
\medskip

\textbf{Part 3} The estimates of this part are \ref{X63} and \ref{X64}.
\end{proof}

\begin{remark}
The conditions of Lemma \ref{Lem_bnd_FourTran_Jn1[v]} are satisfied for any
$\sigma\in C_{0}^{\infty}\left(  \mathcal{O}_{\theta}\right)  $.
\end{remark}

A special case of Lemma \ref{Lem_bnd_FourTran_Jn1[v]} is the 1-dimensional case:

\begin{lemma}
\fbox{\textbf{The case} $d=1$}. Bounds for $\left\Vert \xi^{k}%
\widehat{J_{n\mathbf{1}}^{\theta}\left[  \sigma\right]  }\right\Vert _{\infty
}$ when $n\geq-1$: if $\sigma\in L_{0}^{1}\left(  \mathcal{O}_{\theta}\right)
$ then:

\begin{enumerate}
\item \fbox{$n=-1$}%
\[
\left\Vert \xi^{k}\widehat{J_{-\mathbf{1}}^{\theta}\left[  \sigma\right]
}\right\Vert _{\infty}\leq\left\{
\begin{array}
[c]{ll}%
\left\Vert \widehat{\sigma}\right\Vert _{\infty}, & k=0,\\
\left\Vert \widehat{\sigma}\right\Vert _{1}, & k=1.
\end{array}
\right.
\]

\item \fbox{$n=0$}%
\[
\left\Vert \xi^{k}\widehat{J_{0}^{\theta}\left[  \sigma\right]  }\right\Vert
_{\infty}\leq\left\{
\begin{array}
[c]{ll}%
\left\Vert \widehat{D^{k}\sigma}\right\Vert _{\infty}, & 0\leq k\leq1,\\
e^{-\left(  2-k\right)  }\left\Vert \widehat{D^{k}\sigma}\right\Vert _{1}, &
1\leq k\leq2.
\end{array}
\right.
\]
\medskip

Parts 1 and 2 are included in the following estimate:\medskip

\item \fbox{$n\geq-1$}%
\[
\left\Vert \xi^{k}\widehat{J_{n}^{\theta}\left[  \sigma\right]  }\right\Vert
_{\infty}\leq\left\{
\begin{array}
[c]{ll}%
\left(  n+1-k\right)  !\left\Vert \widehat{D^{k}\sigma}\right\Vert _{\infty
}, & 0\leq k\leq n+1,\\
\left(  n+2-k\right)  ^{\left(  n+2-k\right)  }e^{-\left(  n+2-k\right)
}\left\Vert \widehat{D^{k}\sigma}\right\Vert _{1}, & 1\leq k\leq n+2.
\end{array}
\right.
\]

\end{enumerate}
\end{lemma}

\begin{remark}
Regarding part 3 of the previous lemma: It is well-known that $s^{s}e^{-s}%
\leq\frac{s!}{\sqrt{2\pi s}}$ when $s\geq1$, so for $n\geq-1$ part 3 can be
estimated by%
\begin{equation}
\left\Vert \xi^{k}\widehat{J_{n}^{\theta}\left[  \sigma\right]  }\right\Vert
_{\infty}\leq\left\{
\begin{array}
[c]{ll}%
\left(  n+1-k\right)  !\left\Vert \widehat{D^{k}\sigma}\right\Vert _{\infty
}, & 0\leq k\leq n+1,\\
\left(  n+2-k\right)  !\left\Vert \widehat{D^{k}\sigma}\right\Vert _{1}, &
1\leq k\leq n+2,
\end{array}
\right. \label{X45}%
\end{equation}

\end{remark}

??

\begin{remark}
???
\end{remark}

\appendix

\chapter{Basic notation, definitions and
symbols\label{Ch_Appendx_basic_notation}}

\section{Basic function spaces}

\begin{definition}
\label{Def_Some_basic_spaces}\textbf{Basic function spaces}

All spaces below consist of\textbf{\ complex-valued} functions.

\begin{itemize}
\item $P_{0}=\left\{  0\right\}  $. For $n\geq0$, $P_{n}$ denotes the
polynomials of (total) order at most $n$ i.e. degree at most $n-1$ when
$n\geq1$. These polynomials have the form $\sum_{\left\vert \alpha\right\vert
<n}a_{\alpha}\xi^{\alpha}$, where $a_{\alpha}\in\mathbb{C}$ and $\xi
\in\mathbb{R}^{d}$.

\item $C^{\left(  0\right)  }$ is the space of continuous functions.

\item $C^{\left(  \alpha\right)  }=\left\{  f\in C^{\left(  0\right)
}:D^{\beta}f\in C^{\left(  0\right)  },\text{ }when\text{ }\beta\leq
\alpha\right\}  $.

\item $C_{B}^{\left(  0\right)  }$ is the space of bounded continuous functions.

\item $C_{BP}^{\left(  0\right)  }$ is the space of continuous functions
bounded by a polynomial.

\item $C^{\left(  m\right)  }=\left\{  f\in C^{\left(  0\right)  }:D^{\alpha
}f\in C^{\left(  0\right)  },\text{ }when\text{ }\left\vert \alpha\right\vert
=m\right\}  $.

\item $C_{B}^{\left(  m\right)  }=\left\{  f\in C_{B}^{\left(  0\right)
}:D^{\alpha}f\in C_{B}^{\left(  0\right)  },\text{ }when\text{ }\left\vert
\alpha\right\vert \leq m\right\}  $.

\item $C_{BP}^{\left(  m\right)  }=\left\{  f\in C_{BP}^{\left(  0\right)
}:D^{\alpha}f\in C_{BP}^{\left(  0\right)  },\text{ }when\text{ }\left\vert
\alpha\right\vert \leq m\right\}  $.

\item $C^{\infty}=\bigcap\limits_{m\geq0}C^{\left(  m\right)  }$%
\textbf{;}$\quad C_{B}^{\infty}=\bigcap\limits_{m\geq0}C_{B}^{\left(
m\right)  }$\textbf{;}\quad$C_{BP}^{\infty}=\bigcap\limits_{m\geq0}%
C_{BP}^{\left(  m\right)  }$.

\item $C_{0}^{\infty}$ is the space of $C^{\infty}$ functions that have
compact support.

These are the test functions for the space of distributions.

\item $S$ is the $C^{\infty}$ space of rapidly decreasing functions. These are
the test functions for the tempered distributions of Definition
\ref{Def_Distributions} below.

\item For a bounded continuous function on $\mathbb{R}^{d}$, $\left\Vert
f\right\Vert _{\infty}=\sup_{x\in\mathbb{R}^{d}}\left\vert f\left(  x\right)
\right\vert $.

\item $L_{loc}^{1}$ is the space of measurable functions which are absolutely
integrable on any compact set i.e. any closed bounded set.

I will try to avoid confusion between this notation and Sobolev norms:

\item $L^{1}$ is the complete normed vector space of measurable functions $f$
such that $\int\left\vert f\right\vert <\infty$. Norm is $\left\Vert
f\right\Vert _{1}=\int\left\vert f\right\vert $.

\item $L^{2}$ is the Hilbert space of measurable functions $f$ such that
$\int\left\vert f\right\vert ^{2}<\infty$. Norm is $\left\Vert f\right\Vert
_{2}=\left(  \int\left\vert f\right\vert ^{2}\right)  ^{1/2}$ and inner
product is $\left(  f,g\right)  _{2}=\int f\overline{g}$.

\item $L^{\infty}$ is the complete normed vector space of essentially bounded functions.

The norm is $\left\Vert f\right\Vert _{\infty}=\operatorname*{essup}%
_{x\in\mathbb{R}^{d}}\left\vert f\left(  x\right)  \right\vert $.
\end{itemize}
\end{definition}

\section{Multi-index and vector notation}

\begin{definition}
\label{Def_multi_id}\textbf{Notation}

\begin{enumerate}
\item Multi-indexes are vectors with non-negative integer components.

Let $x=\left(  x_{1},x_{2},\ldots,x_{d}\right)  $ and $y=\left(  y_{1}%
,y_{2},\ldots,y_{d}\right)  $.

Suppose $\symbol{126}$ is one of the binary operations $\leq,=,\geq$.

?? Write $\beta\symbol{126}a$ if $\beta_{i}\symbol{126}\alpha_{i}$ for all $i$.

For $s\in\mathbb{R}$ write $\beta\symbol{126}s$ if $\beta_{i}\symbol{126}s $
for all $i$.

$x\neq y$ if $x_{i}\neq y_{i}$ for some $i$.

$x<y$ if $x\leq y$ and $x\neq y$.

\textbf{Component-wise} inequality is denoted by $x.\neq y$.

$x\centerdot<y$ means $x_{i}<y_{i}$ for all $i$.

$x\centerdot>y$ means $x_{i}>y_{i}$ for all $i$.

The component-wise product is denoted $x.y=\left(  x_{i}y_{i}\right)  $.
Component-wise division is denoted $x\centerdot/y=\left(  x_{i}/y_{i}\right)
$ or $\frac{x}{y}=\left(  \frac{x_{i}}{y_{i}}\right)  $.

Denote $\mathbf{1}=1_{d}=\left(  1,1,1,\ldots,1\right)  \in\mathbb{R}^{d}$ and
let $\left\{  \mathbf{e}_{k}\right\}  _{k=1}^{d}$ be the canonical (or
standard) basis of $\mathbb{R}^{d}$.

\item Denote $\left\vert \alpha\right\vert =\sum\limits_{i=1}^{d}\alpha_{i} $.
The $D^{\alpha}f\left(  x\right)  $ is the derivative of the function $f$ of
degree $\alpha$
\[
D^{0}f\left(  x\right)  =f\left(  x\right)  ,\text{\quad}D^{\alpha}f\left(
x\right)  =\frac{D^{\left\vert \alpha\right\vert }f\left(  x_{1},x_{2}%
,\ldots,x_{d}\right)  }{D_{1}^{\alpha_{1}}x_{1}D_{2}^{\alpha_{2}}x_{2}\ldots
D_{d}^{\alpha_{d}}x_{d}},
\]

and%
\begin{equation}
\frac{1}{\beta!}D^{\beta}x^{\alpha}=\left\{
\begin{array}
[c]{ll}%
\binom{\alpha}{\beta}x^{\alpha-\beta}, & \beta\leq\alpha,\\
0, & otherwise.
\end{array}
\right. \label{a1.021}%
\end{equation}

\item
\[%
\begin{array}
[c]{ll}%
\text{monomial} & x^{\alpha}=x_{1}^{\alpha_{1}}x_{2}^{\alpha_{2}}\ldots
x_{d}^{\alpha_{d}},\\
\text{factorial} & \alpha!=\alpha_{1}!\alpha_{2}!\ldots\alpha_{d}!\text{
}and\text{ }0!=1,\\
\text{binomial} & \dbinom{\alpha}{\beta}=\dfrac{\alpha!}{\left(  \alpha
-\beta\right)  !\beta!},\quad if\text{ }\beta\leq\alpha.
\end{array}
\]

\item The inequality $\left\vert x^{\alpha}\right\vert \leq\left\vert
x\right\vert ^{\left\vert \alpha\right\vert }$ is used often.

\item Important identities are%
\begin{equation}
\frac{\left(  x,y\right)  ^{k}}{k!}=\sum_{\left\vert \alpha\right\vert
=k}\frac{x^{\alpha}y^{\alpha}}{\alpha!},\text{\quad}\frac{\left\vert
x\right\vert ^{2k}}{k!}=\sum_{\left\vert \alpha\right\vert =k}\frac
{x^{2\alpha}}{\alpha!}.\label{1.57}%
\end{equation}

\item If $\binom{k}{\alpha}=\frac{k!}{\alpha!\left(  k-\left\vert
\alpha\right\vert \right)  !}$ then%
\[
\left(  1+xy\right)  ^{k}=\sum_{\left\vert \alpha\right\vert \leq k}\binom
{k}{\alpha}x^{\alpha}y^{\alpha}.
\]

\item The binomial expansion is%
\begin{equation}
\left(  x+y\right)  ^{\alpha}=\sum_{\beta\leq\alpha}\binom{\alpha}{\beta
}x^{\beta}y^{\alpha-\beta},\label{1.070}%
\end{equation}

and a direct consequence is%
\[
\sum_{\beta\leq\alpha}\binom{\alpha}{\beta}=2^{\left\vert \alpha\right\vert
}\text{.}%
\]

\item Useful identities are%
\[
\sum\limits_{\left\vert \alpha\right\vert =k}1=\binom{d+k-1}{k},\quad k\geq1,
\]

and%
\begin{equation}
\sum\limits_{\left\vert \alpha\right\vert \leq k}1=\binom{d+k}{d},\quad
k\geq0,\label{a1.11}%
\end{equation}

and%
\begin{equation}
\sum\limits_{\mathbf{0}\leq\alpha\leq\mathbf{1}}f\left(  \left\vert
\alpha\right\vert \right)  =\sum\limits_{k=0}^{d}\binom{d}{k}f\left(
k\right)  .\label{a1.009}%
\end{equation}

\item Leibniz's rule is%
\[
D^{\alpha}\left(  uv\right)  =\sum_{\beta\leq\alpha}\dbinom{\alpha}{\beta
}D^{\beta}uD^{\alpha-\beta}v.
\]

\end{enumerate}
\end{definition}

\section{Topology}

\begin{definition}
\label{Def_topol_on_Rd}\textbf{Topology on }$\mathbb{R}^{d}$

\begin{enumerate}
\item The Euclidean norm is denoted by $\left\vert x\right\vert $ and the
inner product by $xy$ or $\left(  x,y\right)  $. Other norms are $\left\vert
x\right\vert _{1}=\sum\limits_{i=1}^{d}\left\vert x_{i}\right\vert $ and the
equivalent notation: $\left\vert x\right\vert _{\infty}=\left\vert
x\right\vert _{\max}=\max\limits_{i}\left\vert x_{i}\right\vert $.

\item The (open) ball $B\left(  x;r\right)  =\left\{  y:\left\vert
x-y\right\vert <r\right\}  $.

\item The $\varepsilon\mathbf{-}$neighborhood of a set: for $\varepsilon>0$
the $\varepsilon-$neighborhood of a set $S$ is $S_{\varepsilon}=\bigcup
\limits_{x\in S}B\left(  x;\varepsilon\right)  $.

\item If $a.\neq b$ then $R\left(  a,b\right)  $ will denote an \textbf{open
rectangle with extreme points} $a$ and $b$ given by:%
\[
R\left(  a,b\right)  :=\left\{  \tau.a+\left(  1-\tau\right)  .b:\mathbf{0}%
<\tau<\mathbf{1}\right\}  ,
\]

and%
\[
R\left[  a,b\right]  :=\left\{  \tau.a+\left(  1-\tau\right)  .b:\mathbf{0}%
\leq\tau\leq\mathbf{1}\right\}  ,
\]

will denote the \textbf{closed rectangle}.

\begin{enumerate}
\item We have $R\left(  a,b\right)  =R\left(  \min\left\{  a,b\right\}
,\max\left\{  a,b\right\}  \right)  $.

\item If $\theta\in\left\{  0,1\right\}  ^{d}$ then%
\begin{align*}
R\left(  a,b\right)   & =R\left(  \omega.a+\left(  1-\omega\right)  .b,\left(
1-\omega\right)  .a+\omega.b\right) \\
& =R\left(  \omega.a+\left(  1-\omega\right)  .b,\omega^{\prime}.a+\left(
1-\omega^{\prime}\right)  .b\right)  ,
\end{align*}

where $\omega+\omega^{\prime}=1$, $\omega^{\prime}:=\symbol{126}\omega$.

\item Also, $R\left(  a,b\right)  =R\left(  c,d\right)  $ iff $c=\omega
.a+\left(  1-\omega\right)  .b$ and $d=\left(  1-\omega\right)  .a+\omega.b$
for some $\omega\in\left\{  0,1\right\}  ^{d}$.

\item If $\omega=\left(  a<b\right)  $ and $\omega^{\prime}=\left(
a>b\right)  $ then $\min\left\{  a,b\right\}  =\omega.a+\left(  1-\omega
\right)  .b$ and $\max\left\{  a,b\right\}  =\omega^{\prime}.a+\left(
1-\omega^{\prime}\right)  .b$.

\item If $R\left(  a,b\right)  =R\left(  c,d\right)  $ then there exist unique
$c,d$ such that $c.<d$. In fact $c=\min\left\{  a,b\right\}  $ and
$d=\max\left\{  a,b\right\}  $.
\end{enumerate}

\item $R\left[  a,b\right)  =\left\{  \tau.a+\left(  1-\tau\right)
.b:\mathbf{0}\leq\tau.<\mathbf{1}\right\}  $ will denote
the\textbf{\ partially open rectangle} with extreme points $a.\neq b$.
\end{enumerate}
\end{definition}

\section{Tempered distributions}

\begin{definition}
\label{Def_Distributions}\textbf{Tempered distributions (or generalized
functions of slow growth)}

\begin{enumerate}
\item $S$ is the space of rapidly decreasing $C^{\infty}$ functions We endow
$S$ with the topology defined using the following equivalent sets of seminorms%
\begin{equation}
\left\Vert \left(  1+\left\vert \cdot\right\vert \right)  ^{n}D^{\alpha}%
\psi\right\Vert _{\infty},\text{\quad}n=0,1,2,\ldots;\text{ }\alpha
\geq0.\label{a1.2}%
\end{equation}

\[
\max_{\left\vert \alpha\right\vert \leq n}\left\Vert \left(  1+\left\vert
\cdot\right\vert \right)  ^{n}D^{\alpha}\psi\right\Vert _{\infty},\text{\quad
}n=0,1,2,\ldots.
\]

\begin{equation}
\left\Vert x^{\alpha}D^{\beta}\psi\right\Vert _{\infty},\text{\quad}%
\alpha,\beta\geq0.\label{a1.5}%
\end{equation}

\begin{equation}
\left\Vert x^{\alpha}D^{\beta}\psi\right\Vert _{2},\text{\quad}\alpha
,\beta\geq0.\label{a1.3}%
\end{equation}

\item $S^{\prime}$ denotes the space of tempered distributions or generalized
functions of slow growth. It is the set of all continuous linear functionals
on $S$ under the seminorm topology of part 1.

\item If $f\in S^{\prime}$ and $\phi\in S$ then $\left[  f,\phi\right]
\in\mathbb{C}$ will represent the action of $f$ on the test function $\phi$.
The functions in $S$ are called the test functions of $S^{\prime}$.
\end{enumerate}
\end{definition}

\section{Fourier Transforms\label{Sect_Four_transf}}

\begin{definition}
\label{Def_Fourier}\textbf{Fourier and Inverse Fourier transforms on }$S$
\textbf{and} $S^{\prime}$.

\begin{enumerate}
\item This document uses the two Fourier transform notations%
\[
F\left[  f\right]  =\widehat{f}(\xi)=\left(  2\pi\right)  ^{-d/2}%
\int\limits_{\mathbb{R}^{d}}e^{-ix\xi}f(x)dx,\quad f\in S,
\]

so that $F:S\rightarrow S$. The inverse Fourier transform $F^{-1}:S\rightarrow
S$ is
\[
F^{-1}\left[  f\right]  =\overset{\vee}{f}(\xi)=\left(  2\pi\right)
^{-d/2}\int\limits_{\mathbb{R}^{d}}e^{ix\xi}f(x)dx,\quad f\in S.
\]

\item The Fourier/Inverse Fourier transforms are extended to mappings from
$S^{\prime}$ to $S^{\prime}$ by
\[
\left[  \widehat{f},\phi\right]  =\left[  f,\overset{\vee}{\phi}\right]
;\qquad\left[  \overset{\vee}{f},\phi\right]  =\left[  f,\overset{\vee}{\phi
}\right]  ,\quad f\in S^{\prime},\phi\in S.
\]

\end{enumerate}
\end{definition}

Some important properties are: if $f\in S^{\prime}$ then:

\begin{summary}
\label{Sum_FourTransf}
\end{summary}

\begin{enumerate}
\item $\overset{\vee}{f}(\xi)=\widehat{f}(-\xi)$, $\widehat{\widehat{f}}%
(\xi)=f\left(  -\xi\right)  $, $\overline{\widehat{f}}=\overset{\vee}{\bar{f}%
}$, $\widehat{\overline{f}}=\overline{\overset{\vee}{f}}$.

\item For translations: $\left(  \tau_{c}f\right)  ^{\wedge}=\left(  f\left(
x-c\right)  \right)  ^{\wedge}=e^{-ic\xi}\widehat{f}$,$\quad c\in
\mathbb{R}^{d}$.

\item $\left(  e^{iax}\right)  ^{\wedge}=\left(  2\pi\right)  ^{d/2}%
\delta\left(  \xi-a\right)  $,$\quad a\in\mathbb{R}^{d}$.

\item $\left(  e^{iax}f\right)  ^{\wedge}=\widehat{f}\left(  \cdot-a\right)  $.

\item $D^{\alpha}\widehat{f}=\widehat{\left(  -ix\right)  ^{\alpha}f}=\left(
-i\right)  ^{\left\vert \alpha\right\vert }\widehat{x^{\alpha}f}$.

\item $\left(  x^{\alpha}f\right)  ^{\wedge}=\left(  iD\right)  ^{\alpha
}\widehat{f}=i^{\left\vert \alpha\right\vert }D^{\alpha}\widehat{f}$.

\item $\widehat{D^{\alpha}f}=i^{\left\vert \alpha\right\vert }\xi^{\alpha
}\widehat{f}=\left(  i\xi\right)  ^{\alpha}\widehat{f}$.

\item $\xi^{\alpha}\widehat{f}=\left(  -i\right)  ^{\left\vert \alpha
\right\vert }\widehat{D^{\alpha}f}$.

\item $\left[  \delta,\phi\right]  =\phi\left(  0\right)  $ and
$\widehat{\delta}=\left(  2\pi\right)  ^{-d/2}$ and $\widehat{1}=\left(
2\pi\right)  ^{d/2}\delta$.

\item $\widehat{x^{\alpha}}=\left(  2\pi\right)  ^{-d/2}\left(  iD\right)
^{\alpha}\delta$ and $\widehat{D^{\alpha}\delta}=\left(  2\pi\right)
^{-d/2}\left(  -i\right)  ^{\left\vert \alpha\right\vert }\xi^{\alpha}$.

\item If $p$ is a polynomial then $\widehat{p}=\left(  2\pi\right)
^{d/2}p(iD)\delta$ and $\widehat{p\left(  D\right)  f}=p\left(  i\xi\right)
\widehat{f}$.

\item For dilations: if $\lambda.\neq\mathbf{0}$ then $\left(  \sigma
_{\lambda}f\right)  ^{\wedge}=\left\vert \lambda^{\mathbf{1}}\right\vert
\sigma_{\mathbf{1}./\lambda}\widehat{f}$. Here $\sigma_{\lambda}f\left(
x\right)  :=f\left(  x./\lambda\right)  $.

\item More generally, if $A$ is a regular matrix then $\left(  f\circ
A\right)  ^{\wedge}=\left\vert \det A^{-1}\right\vert \left(  \widehat{f}\circ
A^{-1}\right)  $.
\end{enumerate}

\section{Convolutions}

\begin{definition}
\label{Def_convol}\textbf{The convolution}

\begin{enumerate}
\item If $f\in C_{BP}^{\left(  0\right)  }$ and $\phi\in S$ then the
convolution of $f$ and $\phi$ is denoted by $f\ast\phi$ or $\phi\ast f$ and is
defined by
\begin{equation}
\left(  \phi\ast f\right)  \left(  x\right)  =\left(  2\pi\right)  ^{-d/2}\int
f\left(  x-y\right)  \phi\left(  y\right)  dy=\left(  2\pi\right)  ^{-d/2}\int
f\left(  y\right)  \phi\left(  x-y\right)  dy.\label{a1.72}%
\end{equation}

\item This definition is easily extended to $f\in S^{\prime}$ by the formulas
\begin{equation}
\phi\ast f=\left(  \widehat{\phi}\widehat{f}\right)  ^{\vee}=\left(
2\pi\right)  ^{-d/2}\left[  f_{y},\phi\left(  \cdot-y\right)  \right]
,\quad\phi\in S,\label{1.20}%
\end{equation}

with $\phi\ast f\in C_{BP}^{\infty}$ and for all $\alpha=\beta+\gamma$%
\begin{equation}
D^{\alpha}\left(  \phi\ast f\right)  =\left(  D^{\beta}\phi\right)  \ast
D^{\gamma}f.\label{2.40}%
\end{equation}

NOTE that%
\[
\phi\ast\delta=\left(  \widehat{\phi}\widehat{\delta}\right)  ^{\vee}=\left(
2\pi\right)  ^{-d/2}\phi.
\]

\item In fact, if $f\in S^{\prime}$ is defined by seminorms of type \ref{a1.3}
(or \ref{a1.5}) then (Theorem 6.2 Petersen \cite{Petersen83})%
\[
\left\vert \left(  \phi\ast f\right)  \left(  x\right)  \right\vert
\leq\left(  1+\left\vert x\right\vert \right)  ^{m}\left\vert \phi\right\vert
_{\sigma},\quad\phi\in S,\text{ }x\in\mathbb{R}^{d},
\]

here $\left\vert \cdot\right\vert _{\sigma}$ is defined by seminorms of type
\ref{a1.3} (or \ref{a1.5}), and $m=\max\left\vert \alpha_{i}\right\vert $.

\item If $\psi\in S$ we can write $\left[  \phi\ast f,\psi\right]  =\left[
\left(  \widehat{\phi}\widehat{f}\right)  ^{\vee},\psi\right]  =\left[
\widehat{f},\overset{\vee}{\psi}\widehat{\phi}\right]  $ and the observation
that the last term makes sense when $\overset{\vee}{\phi}\in C_{BP}^{\infty}$
leads to the following extension of the definition of a convolution e.g.
Section 2.6 Petersen \cite{Petersen83}:%
\begin{equation}
\left[  g\ast f,\phi\right]  =\left[  \widehat{f},\overset{\vee}{\phi
}\widehat{g}\right]  =\left[  f,\phi\ast g^{\checkmark}\right]  ,\quad\phi\in
S,\text{ }g\in\left(  C_{BP}^{\infty}\right)  ^{\vee}=O_{c}^{\prime},\text{
}f\in S^{\prime}.\label{a1.1}%
\end{equation}

See Section \ref{Sect_temp_convol_multip} below.

\item If $f\in C_{BP}^{\left(  0\right)  }$ and $g\in C_{0}^{\left(  0\right)
}$ then $f\ast g$ is a regular tempered distribution and%
\[
\left(  f\ast g\right)  \left(  x\right)  =\left(  g\ast f\right)  \left(
x\right)  =\left(  2\pi\right)  ^{-\frac{d}{2}}\int f\left(  x-y\right)
g\left(  y\right)  dy=\left(  2\pi\right)  ^{-\frac{d}{2}}\int f\left(
y\right)  g\left(  x-y\right)  dy.
\]

This can be proved using Theorem 2.7.5 of Vladimirov \cite{Vladimirov}.

\item If $f\in L^{p}$, $g\in L^{q}$ satisfy $\frac{1}{p}+\frac{1}{q}%
=1+\frac{1}{r}$ and $1\leq p,q,r\leq\infty$, then $f\ast g\in L^{r}$ where the
convolution is given by the formulas of part 5.

Further, we have \textbf{Young's inequality}%
\begin{equation}
\left\Vert f\ast g\right\Vert _{r}\leq\left\Vert f\right\Vert _{p}\left\Vert
g\right\Vert _{q}.\label{a1.022}%
\end{equation}

\item If $A$ is a regular $d\times d$ matrix then%
\[
f\ast\left(  g\circ A\right)  =\frac{1}{\left\vert \det A\right\vert }\left(
\left(  f\circ A^{-1}\right)  \ast g\right)  \circ A.
\]

\item If $f\in S^{\prime}$ and $g\in S$ then using \ref{1.20} it can be shown
that%
\[
x^{\alpha}\left(  f\ast g\right)  =\sum\limits_{\beta\leq\alpha}\tbinom
{\alpha}{\beta}\left(  z^{\beta}f\right)  \ast\left(  z^{\alpha-\beta
}g\right)  .
\]

One consequence is that $S\ast S\subset S$.

\item $f,g\in L^{2}$ implies $fg\in L^{1}$ and $\widehat{fg}=\widehat{f}%
\ast\widehat{g}\in C_{B}^{\left(  0\right)  }$.
\end{enumerate}
\end{definition}

\section{Convolutors and multipliers\ in $S$ and $S^{\prime}$%
\label{Sect_temp_convol_multip}}

\begin{definition}
\label{Def_convol_multip_S}\textbf{Multipliers }$O_{M}$\textbf{\ in
}$S^{\prime}$: See for example Sections 2.5, 2.6 and 2.7 of Petersen
\cite{Petersen83}. See also Chapter 6 of Kinani and Oudadess \cite{KinOud2010}
and Section 3.9 of the Appendix of Dautray and Lions \cite{DautLion88}.

\begin{enumerate}
\item \textbf{Definition of} $O_{M}$:%
\begin{align*}
O_{M}  & :=\left\{  \psi\in C^{\infty}:\psi\phi\in S\text{ }\forall\phi\in
S\right\} \\
& =\left\{  \psi\in C^{\infty}:\psi f\in S^{\prime}\text{ }\forall f\in
S^{\prime}\right\}  .
\end{align*}

\textbf{Properties}

\item (Lemma 5.13 Petersen) If $\psi\in O_{M}$ then multiplication by $\psi$
maps $S\rightarrow S$ continuously and $S^{\prime}$ into $S^{\prime}$ continuously.

\item (Exercise 5.14 Petersen) $g\in\mathcal{D}^{\prime}$ and $\psi g\in S$
$\forall$ $\psi\in S$ implies $g\in O_{M}$.

\item (Theorem 5.15 Petersen) $\psi\in O_{M}$ iff $\psi\in C^{\infty}$ and all
its derivatives are slowly increasing i.e.%
\[
\left\vert D^{\alpha}\psi\left(  x\right)  \right\vert \leq C_{\alpha}\left(
1+\left\vert x\right\vert ^{2}\right)  ^{m_{\alpha}},\quad\forall\alpha,\text{
}and\text{ }x\in\mathbb{R}^{d}.
\]

Thus $O_{M}\equiv C_{BP}^{\infty}$.

\item (Corollary 5.16 Petersen) $S\hookrightarrow O_{M}\hookrightarrow
S^{\prime}$.
\end{enumerate}
\end{definition}

\begin{definition}
\label{Def_convolutor_convol}\textbf{Convolutors }$O_{c}^{\prime}$ \textbf{and
the convolution on }$S^{\prime}\times O_{c}^{\prime}$

\begin{enumerate}
\item $g\in S^{\prime}$ is a \textbf{convolutor} if $\phi\ast g\in S$ for all
$\phi\in S$ and the mapping $\phi\rightarrow\phi\ast g:S\rightarrow S $ is continuous.

\item $g\in O_{c}^{\prime}$ implies $g^{\checkmark}\in O_{c}^{\prime}$.

\item Noting \ref{a1.1} we define the following \textbf{convolution}: if $f\in
S^{\prime}$ and $g\in O_{c}^{\prime}$ then%
\begin{align*}
\left[  f\ast g,\phi\right]   & :=\left[  f,g^{\checkmark}\ast\phi\right]
,\quad\phi\in S.\\
g\ast f  & :=f\ast g.
\end{align*}

\item \textbf{Some properties:}

\begin{enumerate}
\item When $g\in O_{c}^{\prime}$ and $x\in\mathbb{R}^{d}$: $\tau
_{x}g,y^{\alpha}g,D^{\alpha}g,\widehat{g}\in O_{c}^{\prime}$ $\forall\alpha$.

\item $D^{\alpha}\left(  f\ast g\right)  =D^{\alpha}f\ast g=f\ast D^{\alpha}g$
$\forall\alpha$.

\item $\tau_{x}\left(  f\ast g\right)  =\tau_{x}f\ast g=f\ast\tau_{x}g$.

\item $f\in S^{\prime}$ and $\psi\in S$ implies $\psi f\in O_{c}^{\prime}$.

\item $\mathcal{E}^{\prime}\left(  \mathbb{R}^{d}\right)  \subseteq
O_{c}^{\prime}$.

\item $O_{c}^{\prime}$ and $O_{M}$ are not duals (see appendix of Dautray and
Lions \cite{DautLion88}).
\end{enumerate}

\item If $f\in S^{\prime}$ and $g,h\in O_{c}^{\prime}$ then $g\ast h=h\ast
g\in O_{c}^{\prime}$ and $f\ast\left(  g\ast h\right)  =\left(  f\ast
g\right)  \ast h$.

\item All the convolutions defined so far agree on common domains.

\item $\delta\in O_{c}^{\prime}$ because $\delta\in\mathcal{E}^{\prime}$. Also
$f\ast\delta=\left(  2\pi\right)  ^{-d/2}f$ when $f\in S^{\prime}$.
\end{enumerate}
\end{definition}

The extension of the Fourier transform to tempered distributions in Section
2.7 of Petersen allows the proof of the following results:

\begin{theorem}
\textbf{Exchange formulas} (Theorem 7.2 Petersen \cite{Petersen83})

\begin{enumerate}
\item If $g\in O_{c}^{\prime}$ then $\widehat{g}\in O_{M}$ and%
\[
\left(  f\ast g\right)  ^{\wedge}=\widehat{g}\widehat{f},\quad f\in S^{\prime
}.
\]

\item If $\psi\in O_{M}$ then $\widehat{\psi}\in O_{c}^{\prime}$ and%
\[
\widehat{\psi f}=\widehat{\psi}\ast\widehat{f},\quad f\in S^{\prime}.
\]

\item In particular, the Fourier transform is an \textbf{isomorphism} from
$O_{c}^{\prime}$ to $O_{M}$ and from $O_{M}$ to $O_{c}^{\prime}$.
\end{enumerate}
\end{theorem}

A concrete description of $O_{c}^{\prime}$ is:

\begin{theorem}
(Exercise 7.3 of Petersen \cite{Petersen83} and Proposition 8.7 of Kinani and
Oudadess \cite{KinOud2010})

$f\in O_{c}^{\prime}$ iff for every integer $m>0$ there is an integer
$k_{m}\geq0$ such that there exists a finite family of continuous functions
$\left\{  f_{\alpha}\right\}  _{\left\vert \alpha\right\vert \leq k_{m}}$
which satisfies $f=\sum_{\left\vert \alpha\right\vert \leq k_{m}}D^{\alpha
}f_{\alpha}$ and for which $\left(  1+\left\vert \cdot\right\vert ^{2}\right)
^{m}f_{\alpha}\in L^{\infty}$ for all $\left\vert \alpha\right\vert \leq
k_{m}$ or, alternatively, for which $\left(  1+\left\vert \cdot\right\vert
^{2}\right)  ^{m}f_{\alpha}\rightarrow0$ at infinity.
\end{theorem}

\section{Taylor series expansions with integral
remainder\label{Sect_apx_TaylorSeries}}

Suppose $u\in C^{\left(  n\right)  }\left(  \mathbb{R}^{d}\right)  $ for some
$n\geq1$. Then the Taylor series expansion about $z$ is given by
\begin{equation}
u\left(  z+b\right)  =\sum_{\left\vert \beta\right\vert <n}\frac{b^{\beta}%
}{\beta!}\left(  D^{\beta}u\right)  (z)+\left(  \mathcal{R}_{n}u\right)
\left(  z,b\right)  ,\label{p87}%
\end{equation}

where the integral remainder term
\begin{align}
\left(  \mathcal{R}_{n}u\right)  \left(  z,b\right)   & =n\sum_{\left\vert
\beta\right\vert =n}\frac{b^{\beta}}{\beta!}\int_{0}^{1}t^{n-1}\left(
D^{\beta}u\right)  \left(  z+\left(  1-t\right)  b\right)  dt\nonumber\\
& =n\sum_{\left\vert \beta\right\vert =n}\frac{b^{\beta}}{\beta!}\int_{0}%
^{1}\left(  1-t\right)  ^{n-1}\left(  D^{\beta}u\right)  \left(  z+tb\right)
dt\label{a1.12}\\
& =\frac{1}{\left(  n-1\right)  !}\int_{0}^{1}\left(  1-t\right)
^{n-1}\left(  \left(  bD\right)  ^{n}u\right)  \left(  z+tb\right)
dt\nonumber\\
& =\frac{\left\vert b\right\vert ^{n}}{\left(  n-1\right)  !}\int_{0}%
^{1}\left(  1-t\right)  ^{n-1}\left(  \left(  \widehat{b}D\right)
^{n}u\right)  \left(  z+tb\right)  ds,\label{Ap130}%
\end{align}

satisfies the estimate
\begin{equation}
\left\vert \left(  \mathcal{R}_{n}u\right)  \left(  z,b\right)  \right\vert
\leq\frac{\left\vert b\right\vert ^{n}}{n!}\max_{t\in\left[  z,z+b\right]
}\left\vert \left(  \widehat{b}D\right)  ^{n}u\left(  t\right)  \right\vert
,\label{2.68}%
\end{equation}

and the estimate (prove using \ref{1.57}),%
\begin{equation}
\left\vert \left(  \mathcal{R}_{n}u\right)  \left(  z,b\right)  \right\vert
\leq\frac{d^{\frac{n}{2}}\left\vert b\right\vert ^{n}}{n!}\max
_{\substack{\left\vert \beta\right\vert =n \\t\in\left[  z,z+b\right]
}}\left\vert D^{\beta}u\left(  t\right)  \right\vert ,\label{1.34}%
\end{equation}

and the estimate%
\begin{equation}
\left\vert \left(  \mathcal{R}_{n}u\right)  \left(  z,b\right)  \right\vert
\leq\left\vert b\right\vert ^{n}\sum_{\left\vert \beta\right\vert =n}%
\max_{t\in\left[  z,z+b\right]  }\left\vert \left(  D^{\beta}u\right)
(t)\right\vert ,\label{a1.37}%
\end{equation}

and%
\begin{equation}
\left\vert \left(  \mathcal{R}_{n}u\right)  \left(  z,b\right)  \right\vert
\leq\sum_{\left\vert \beta\right\vert =n}\frac{\left\vert b^{\beta}\right\vert
}{\beta!}\max_{y\in\left[  z,z+b\right]  }\left\vert \left(  D^{\beta
}u\right)  \left(  y\right)  \right\vert \leq\left\vert b\right\vert ^{n}%
\sum_{\left\vert \beta\right\vert =n}\frac{1}{\beta!}\max_{y\in\left[
z,z+b\right]  }\left\vert \left(  D^{\beta}u\right)  \left(  y\right)
\right\vert .\label{p20}%
\end{equation}

The remainder term can also be expressed in the form%
\begin{equation}
\left(  \mathcal{R}_{n}u\right)  \left(  z,b\right)  =\sum_{\left\vert
\beta\right\vert =n}\frac{b^{\beta}}{\beta!}D^{\beta}u\left(  z+t\left(
z,b\right)  b\right)  ,\label{a1.35}%
\end{equation}

where the function $t$ satisfies $0<t\left(  z,b\right)  <1$.

Also%
\begin{align}
\left\vert \left(  \mathcal{R}_{n}u\right)  \left(  0,x\right)  \right\vert  &
\leq n\sum_{\left\vert \beta\right\vert =n}\frac{x_{+}^{\beta}}{\sqrt{\beta!}%
}\frac{1}{\sqrt{\beta!}}\int_{0}^{1}\left(  1-t\right)  ^{n-1}\left\vert
\left(  D^{\beta}u\right)  \left(  tx\right)  \right\vert dt\nonumber\\
& \leq n\left(  \frac{\left\vert x\right\vert ^{2n}}{n!}\right)  ^{1/2}\left(
\sum_{\left\vert \beta\right\vert =n}\frac{1}{\beta!}\left(  \int_{0}%
^{1}\left(  1-t\right)  ^{n-1}\left\vert \left(  D^{\beta}u\right)  \left(
tx\right)  \right\vert dt\right)  ^{2}\right)  ^{1/2}\nonumber\\
& \leq\left(  \frac{\left\vert x\right\vert ^{2n}}{n!}\right)  ^{1/2}\left(
\sum_{\left\vert \beta\right\vert =n}\frac{1}{\beta!}\left\Vert D^{\beta
}u\right\Vert _{\infty,\overline{B}_{\left\vert x\right\vert }}^{2}\right)
^{1/2}\nonumber\\
& =\frac{1}{\sqrt{n!}}\left(  \sum_{\left\vert \beta\right\vert =n}\frac
{1}{\beta!}\left\Vert D^{\beta}u\right\Vert _{\infty,\overline{B}_{\left\vert
x\right\vert }}^{2}\right)  ^{1/2}\left\vert x\right\vert ^{n}.\label{a1.73}%
\end{align}

\begin{theorem}
\label{Thm_Taylor_rad_remain}Regarding \ref{Ap130}:

\begin{enumerate}
\item When $z=0$ we can write%
\begin{equation}
\left(  \mathcal{R}_{n}u\right)  \left(  0,b\right)  =\frac{\left\vert
b\right\vert ^{n}}{\left(  n-1\right)  !}\int_{0}^{1}\left(  1-s\right)
^{n-1}\left(  \left(  \widehat{\cdot}D\right)  ^{n}u\right)  \left(
sb\right)  ds,\label{Ap134}%
\end{equation}

with the bound%
\begin{equation}
\left\vert \left(  \mathcal{R}_{n}u\right)  \left(  0,b\right)  \right\vert
\leq\frac{\left\vert b\right\vert ^{n}}{n!}\max_{s\in\left[  0,1\right]
}\left\vert \left(  \left(  \widehat{\cdot}D\right)  ^{n}u\right)  \left(
sb\right)  \right\vert .\label{Ap144}%
\end{equation}

When $u$ is radial $\left(  \mathcal{R}_{n}u\right)  \left(  0,\cdot\right)  $
is also radial. In fact, if $u_{\odot}\left(  \left\vert \cdot\right\vert
\right)  =u\left(  \cdot\right)  $ then%
\begin{equation}
\left(  \mathcal{R}_{n}u\right)  \left(  0,b\right)  =\frac{\left\vert
b\right\vert ^{n}}{\left(  n-1\right)  !}\int_{0}^{1}\left(  1-s\right)
^{n-1}\left(  D^{n}u_{\odot}\right)  \left(  s\left\vert b\right\vert \right)
ds,\label{Ap135}%
\end{equation}

with upper bound%
\begin{equation}
\left\vert \left(  \mathcal{R}_{n}u\right)  \left(  0,b\right)  \right\vert
\leq\frac{\left\vert b\right\vert ^{n}}{n!}\max_{s\in\left[  0,1\right]
}\left\vert \left(  D^{n}u_{\odot}\right)  \left(  sb\right)  \right\vert
.\label{Ap146}%
\end{equation}

\item If $u$ is radial then $\frac{1}{n!}\left(  \left(  bD\right)
^{n}u\right)  \left(  0\right)  =\sum\limits_{\left\vert \alpha\right\vert
=n}\frac{b^{\alpha}}{\alpha!}D^{\alpha}u\left(  0\right)  $ is a radial
function of $b$ for $n\geq0$.

\item In fact, if $u$ is radial then
\begin{equation}
\sum\limits_{\left\vert \alpha\right\vert =n}\frac{b^{\alpha}}{\alpha
!}D^{\alpha}u\left(  0\right)  =\frac{\left\vert b\right\vert ^{n}}{n!}%
D^{n}u_{\odot}\left(  0\right)  =\frac{\left\vert b\right\vert ^{n}}%
{n!}\left(  \left(  \widehat{\cdot}D\right)  ^{n}u\right)  \left(  0\right)
.\label{Ap133}%
\end{equation}

\item If $\left\vert b\right\vert =r$ then%
\begin{align*}
u_{\odot}\left(  r\right)   & =\sum_{k=0}^{n-1}\frac{r^{k}}{k!}D^{k}u_{\odot
}\left(  0\right)  +\mathcal{R}_{n}u_{\odot}\left(  0,r\right)  ,\\
\mathcal{R}_{n}u_{\odot}\left(  0,r\right)   & =\mathcal{R}_{n}u\left(
0,b\right)  .
\end{align*}

\end{enumerate}
\end{theorem}

\begin{proof}
\textbf{Part 1} When $z=0$ \ref{Ap130} becomes%
\begin{align*}
\left(  \mathcal{R}_{n}u\right)  \left(  0,b\right)   & =\frac{\left\vert
b\right\vert ^{n}}{\left(  n-1\right)  !}\int_{0}^{1}\left(  1-t\right)
^{n-1}\left(  \left(  \widehat{b}D\right)  ^{n}u\right)  \left(  tb\right)
ds\\
& =\frac{\left\vert b\right\vert ^{n}}{\left(  n-1\right)  !}\int_{0}%
^{1}\left(  1-s\right)  ^{n-1}\left(  \left(  \widehat{sb}D\right)
^{n}u\right)  \left(  sb\right)  ds\\
& =\frac{\left\vert b\right\vert ^{n}}{\left(  n-1\right)  !}\int_{0}%
^{1}\left(  1-s\right)  ^{n-1}\left(  \left(  \widehat{\cdot}D\right)
^{n}u\right)  \left(  sb\right)  ds.
\end{align*}

Now suppose $u\left(  \cdot\right)  =u_{\odot}\left(  \left\vert
\cdot\right\vert \right)  $. Then from \ref{Ap027} below, $\left(  \left(
\widehat{\cdot}D\right)  ^{n}u\right)  \left(  x\right)  =\left(
D^{n}u_{\odot}\right)  \left(  \left\vert x\right\vert \right)  $ and so%
\[
\left(  \mathcal{R}_{n}u\right)  \left(  0,b\right)  =\frac{\left\vert
b\right\vert ^{n}}{\left(  n-1\right)  !}\int_{0}^{1}\left(  1-s\right)
^{n-1}\left(  D^{n}u_{\odot}\right)  \left(  s\left\vert b\right\vert \right)
ds
\]
$\medskip$

\textbf{Part 2} From \ref{p87} we have\textbf{\ }%
\begin{align*}
u\left(  b\right)   & =\sum_{\left\vert \alpha\right\vert \leq n-1}%
\frac{b^{\alpha}}{\alpha!}\left(  D^{\alpha}u\right)  \left(  0\right)
+\left(  \mathcal{R}_{n}u\right)  \left(  0,b\right)  ,\\
u\left(  b\right)   & =\sum_{\left\vert \alpha\right\vert \leq n-1}%
\frac{b^{\alpha}}{\alpha!}\left(  D^{\alpha}u\right)  \left(  0\right)
+\sum_{\left\vert \alpha\right\vert =n}\frac{b^{\alpha}}{\alpha!}\left(
D^{\alpha}u\right)  \left(  0\right)  +\left(  \mathcal{R}_{n+1}u\right)
\left(  0,b\right)  ,
\end{align*}

so that by inspection%
\begin{equation}
\sum_{\left\vert \alpha\right\vert =n}\frac{b^{\alpha}}{\alpha!}\left(
D^{\alpha}u\right)  \left(  0\right)  =\left(  \mathcal{R}_{n}u\right)
\left(  0,b\right)  -\left(  \mathcal{R}_{n+1}u\right)  \left(  0,b\right)
.\label{Ap132}%
\end{equation}

This result follows since from part 1 we know that $\left(  \mathcal{R}%
_{n+1}u\right)  \left(  0,b\right)  $ and $\left(  \mathcal{R}_{n}u\right)
\left(  0,b\right)  $ are both radial.\medskip

\textbf{Part 3} Equation \ref{Ap132} can be written as%
\[
\frac{1}{n!}\left(  \left(  bD\right)  ^{n}u\right)  \left(  0\right)
=\left(  \mathcal{R}_{n}u\right)  \left(  0,b\right)  -\left(  \mathcal{R}%
_{n+1}u\right)  \left(  0,b\right)  ,
\]

and using \ref{Ap135} and applying integration by parts gives
\begin{align*}
& \left(  \mathcal{R}_{n}u\right)  \left(  0,b\right)  -\left(  \mathcal{R}%
_{n+1}u\right)  \left(  0,b\right) \\
& =\frac{1}{\left(  n-1\right)  !}\int_{0}^{1}\left(  1-s\right)
^{n-1}\left(  D^{n}u_{\odot}\right)  \left(  s\right)  ds-\frac{1}{n!}\int%
_{0}^{1}\left(  1-s\right)  ^{n}\left(  D^{n+1}u_{\odot}\right)  \left(
s\right)  ds\\
& =\frac{1}{\left(  n-1\right)  !}\int_{0}^{1}\left(  1-s\right)
^{n-1}\left(  D^{n}u_{\odot}\right)  \left(  s\right)  ds-\\
& \qquad\qquad-\frac{1}{n!}\left(  \left[  \left(  1-s\right)  ^{n}%
D^{n}u_{\odot}\left(  s\right)  \right]  _{0}^{1}-\int_{0}^{1}\left(
D_{s}\left(  1-s\right)  ^{n}\right)  D^{n}u_{\odot}\left(  s\right)
ds\right) \\
& =\frac{1}{\left(  n-1\right)  !}\int_{0}^{1}\left(  1-s\right)
^{n-1}\left(  D^{n}u_{\odot}\right)  \left(  s\right)  ds-\\
& \qquad\qquad-\frac{1}{n!}\left(  -D^{n}u_{\odot}\left(  0\right)  +\int%
_{0}^{1}\left(  1-s\right)  ^{n-1}D^{n}u_{\odot}\left(  s\right)  ds\right) \\
& =\frac{1}{n!}D^{n}u_{\odot}\left(  0\right)  .
\end{align*}

The second equality of \ref{Ap133} follows from \ref{Ap027} below.
\end{proof}

\begin{theorem}
\label{Thm_Tay_rem_zeros}Suppose $u\in C_{B}^{\left(  n\right)  }\left(
\mathbb{R}^{d}\right)  $ and $D^{\alpha}u\left(  0\right)  =0$ when
$\left\vert \alpha\right\vert <n$. Then for all $r>0$,%
\begin{equation}
\left\vert u\left(  x\right)  \right\vert \leq\frac{\left\vert x\right\vert
^{n}}{n!}\left\Vert \left(  \widehat{\cdot}D\right)  ^{n}u\right\Vert
_{\infty;\overline{B}_{r}},\quad\left\vert x\right\vert \leq r,\label{Ap019}%
\end{equation}

where $\left(  \left(  \widehat{\cdot}D\right)  u\right)  \left(  x\right)
:=\sum\limits_{k=1}^{d}\widehat{x}_{k}D_{k}u\left(  x\right)  =\frac
{1}{\left\vert x\right\vert }\sum\limits_{k=1}^{d}x_{k}D_{k}u\left(  x\right)
$. Further%
\begin{equation}
\left\vert u\left(  x\right)  \right\vert \leq\frac{\left\vert x\right\vert
^{n}}{n!}\left\Vert \left(  \widehat{\cdot}D\right)  ^{n}u\right\Vert
_{\infty},\quad x\in\mathbb{R}^{d}.\label{Ap01}%
\end{equation}

Things are especially nice when $u$ is radial. In fact, if $u\left(  x\right)
=u_{\odot}\left(  \left\vert x\right\vert \right)  $ then $\left(
\widehat{\cdot}D\right)  ^{n}u$ is radial and%
\begin{equation}
\left(  \left(  \widehat{\cdot}D\right)  ^{n}u\right)  \left(  x\right)
=\left(  D^{n}u_{\odot}\right)  \left(  \left\vert x\right\vert \right)
.\label{Ap027}%
\end{equation}

\end{theorem}

\begin{proof}
Inequality \ref{Ap144} with $b=x$ and $\left\vert x\right\vert \leq r$ gives%
\begin{align*}
\left\vert u\left(  x\right)  \right\vert  & \leq\frac{\left\vert x\right\vert
^{n}}{n!}\max_{t\in\left[  0,1\right]  }\left\vert \left(  \left(
\widehat{x}D\right)  ^{n}u\right)  \left(  tx\right)  \right\vert
=\frac{\left\vert x\right\vert ^{n}}{n!}\max_{t\in\left[  0,1\right]
}\left\vert \left(  \left(  \widehat{tx}D\right)  ^{n}u\right)  \left(
tx\right)  \right\vert =\\
& =\frac{\left\vert x\right\vert ^{n}}{n!}\max_{t\in\left[  0,1\right]
}\left\vert \left(  \left(  \widehat{\cdot}D\right)  ^{n}u\right)  \left(
tx\right)  \right\vert \leq\frac{\left\vert x\right\vert ^{n}}{n!}\max
_{z\in\overline{B}_{r}}\left\vert \left(  \left(  \widehat{\cdot}D\right)
^{n}u\right)  \left(  z\right)  \right\vert \leq\\
& =\frac{\left\vert x\right\vert ^{n}}{n!}\left\Vert \left(  \widehat{\cdot
}D\right)  ^{n}u\right\Vert _{\infty;\leq r}\leq\frac{\left\vert x\right\vert
^{n}}{n!}\left\Vert \left(  \widehat{\cdot}D\right)  ^{n}u\right\Vert
_{\infty}.
\end{align*}

Equation \ref{Ap027} follows from ??.
\end{proof}

The next result may or may not be useful.

\begin{theorem}
\label{Thm_Taylor_rem_1divBk}Set%
\[
f_{n}\left(  \beta\right)  =\max_{y\in\left[  z,z+b\right]  }\left\vert
\left(  D^{\beta}u\right)  \left(  y\right)  \right\vert .
\]

Then%
\[
\left\vert \left(  \mathcal{R}_{n}u\right)  \left(  z,b\right)  \right\vert
\leq\left\vert b\right\vert ^{n}\frac{\left(  \chi_{n}d\right)  ^{n}}{n!},
\]

where%
\[
\chi_{n}=\max_{\left\vert \beta\right\vert =n}\max_{k}\left\vert f_{n}\left(
\beta\right)  \right\vert ^{\frac{1}{d\beta_{k}}}=\left(  \max_{\left\vert
\beta\right\vert =n}\max_{k}\left\vert f_{n}\left(  \beta\right)  \right\vert
^{\frac{1}{\beta_{k}}}\right)  ^{\frac{1}{d}}.
\]

\end{theorem}

\begin{proof}
From \ref{p20},%
\begin{align*}
\left\vert \left(  \mathcal{R}_{n}u\right)  \left(  z,b\right)  \right\vert  &
\leq\left\vert b\right\vert ^{n}\sum\limits_{\left\vert \beta\right\vert
=n}\frac{1}{\beta!}\max_{y\in\left[  z,z+b\right]  }\left\vert \left(
D^{\beta}u\right)  \left(  y\right)  \right\vert =\left\vert b\right\vert
^{n}\sum\limits_{\left\vert \beta\right\vert =n}\frac{1}{\beta!}f_{n}\left(
\beta\right)  =\\
& =\left\vert b\right\vert ^{n}\sum\limits_{\left\vert \beta\right\vert
=n}\frac{1}{\beta!}\mathbf{1}^{\beta}\left(  \left\vert f_{n}\left(
\beta\right)  \right\vert ^{\frac{1}{d\beta_{k}}}\right)  ^{\beta}%
\leq\left\vert b\right\vert ^{n}\sum\limits_{\left\vert \beta\right\vert
=n}\frac{1}{\beta!}\mathbf{1}^{\beta}\left(  \chi_{n}\mathbf{1}\right)
^{\beta}=\\
& =\left\vert b\right\vert ^{n}\left(  \chi_{n}\right)  ^{n}\sum
\limits_{\left\vert \beta\right\vert =n}\frac{\mathbf{1}^{2\beta}}{\beta
!}=\left\vert b\right\vert ^{n}\left(  \chi_{n}\right)  ^{n}\frac{1}%
{n!}\left\vert \mathbf{1}\right\vert ^{2n}=\left\vert b\right\vert ^{n}%
d^{n}\frac{\left(  \chi_{n}\right)  ^{n}}{n!}=\\
& =\left\vert b\right\vert ^{n}\frac{\left(  \chi_{n}d\right)  ^{n}}{n!}.
\end{align*}

\end{proof}

\section{Limits and differentiation under the integral
sign\label{Sect_lim_under_integ}}

Not used below but Bauer \cite{Bauer2001} studies integrals with parameters in
$\mathbb{R}^{n}$ in Section 16 of Chapter II - specifically in Lemma 16.1 and
Corollary 16.3.

\begin{lemma}
\label{Lem_lim_under_integ}Suppose $f:\mathbb{R}^{m+n}\rightarrow\mathbb{C}$
and we write $f\left(  \xi,x\right)  $ where $\xi\in\mathbb{R}^{m}$ and
$x\in\mathbb{R}^{n}$. Suppose that for all $\xi$, $f\left(  \xi,x\right)  $ is
continuous at $x_{0}$, and that there exists a function $g$ satisfying:

\begin{enumerate}
\item $\left\vert f\left(  \xi,x\right)  \right\vert \leq g\left(  \xi\right)
$ $\forall\xi$ and $x$;

\item $g\in L^{1}$.
\end{enumerate}

Then $\int f\left(  \xi,x\right)  d\xi$ is continuous at $x_{0}$.
\end{lemma}

\begin{proof}
Use the Lebesgue dominated convergence theorem with dominating function $g$.
\end{proof}

\begin{lemma}
\label{Lem_diff_under_integ}Suppose $f:\mathbb{R}^{m+n}\rightarrow\mathbb{C} $
and we write $f\left(  \xi,x\right)  $ where $\xi\in\mathbb{R}^{m}$ and
$x\in\mathbb{R}^{n}$. Further, suppose that \textbf{for all} $\xi$: $D_{x_{j}%
}f\left(  \xi,x\right)  $ exists at $x_{0}$ and that there exists a function
$g\left(  \xi,x_{0}\right)  $ and real $h_{0}>0$ satisfying:

\begin{enumerate}
\item $\left\vert f\left(  \xi,x_{0}+h\right)  -f\left(  \xi,x_{0}\right)
\right\vert \leq g\left(  \xi,x_{0}\right)  h$, for all $x$ and $\left\vert
h\right\vert \leq h_{0}$;

\item $g\left(  \xi,x_{0}\right)  \in L^{1}$.
\end{enumerate}

Then at $x=x_{0}$,%
\[
D_{x_{j}}\int f\left(  \xi,x\right)  d\xi=\int D_{x_{j}}f\left(  \xi,x\right)
d\xi.
\]

\end{lemma}

\begin{proof}
Use Lemma \ref{Lem_lim_under_integ}.
\end{proof}

\begin{lemma}
\label{Lem_diff_under_integral_2}Suppose $f:\mathbb{R}^{m+n}\rightarrow
\mathbb{C}$ and we write $f\left(  \xi,x\right)  $ where $\xi\in\mathbb{R}%
^{m}$ and $x\in\mathbb{R}^{n}$. Further suppose that:

\begin{enumerate}
\item For each $\xi$, $f\left(  \xi,\cdot\right)  \in C^{\left(  k\right)
}\left(  \mathbb{R}^{n}\right)  $.

\item For each $x$, $\int\left\vert D_{x}^{\alpha}f\left(  \xi,x\right)
\right\vert d\xi<\infty$ for $\left\vert \alpha\right\vert \leq k$.
\end{enumerate}

Then we have
\[
D_{x}^{\alpha}\int f\left(  \xi,x\right)  d\xi=\int D_{x}^{\alpha}f\left(
\xi,x\right)  d\xi,\;when\text{ }\left\vert \alpha\right\vert \leq k,
\]

and $\int f\left(  \xi,\cdot\right)  d\xi\in C^{\left(  k\right)  }\left(
\mathbb{R}^{n}\right)  $.
\end{lemma}

\begin{proof}
Induction on $k$ using Lemma \ref{Lem_diff_under_integ}.
\end{proof}

\section{Spherical coordinate integrals\label{Sect_spher_coord_integrals}}

Since I derived the results of this section I have discovered the references
Gradshteyn and Ryzhik \cite{GradRyz07} and Prudnikov, Brychkov and Marichev
\cite{PrudBryMar86} for multivariate integration.

These results were proved in the section \textit{Application of general
multivariate spherical coordinates} in the appendix of the positive order
weight function document. This lemma will employ the technique of Section 4.1,
Stein and Weiss \cite{SteinWeiss71} which defines radial functions in terms of
orthogonal transformations.

Stein and Weiss observe that a function $f$ is radial if and only if $f\left(
\mathcal{O}x\right)  =f\left(  x\right)  $ for any linear, orthogonal
transformation $\mathcal{O}:\mathbb{R}^{d}\rightarrow\mathbb{R}^{d}$ and any
$x\in\mathbb{R}^{d}$.

Note that an orthogonal transformation $\mathcal{O}$ satisfies $\mathcal{O}%
^{T}=\mathcal{O}^{-1}$ where $\left(  \mathcal{O}x,y\right)  =\left(
x,\mathcal{O}^{T}y\right)  $ for the Euclidean inner product. An orthogonal
transformation has a Jacobian of one.

\begin{theorem}
\label{Thm_Integ_u(xy)f(|x|)dx}Suppose $\xi,x\in\mathbb{R}^{d}$ and $\xi
x=\left(  \xi,x\right)  $ denotes the Euclidean inner product. Then:

\begin{enumerate}
\item $\int_{r\leq\left\vert x\right\vert \leq R}u\left(  \xi\widehat{x}%
\right)  f\left(  \left\vert x\right\vert \right)  dx$ is a radial function of
$\xi$ and%
\begin{equation}
\int\limits_{r\leq\left\vert x\right\vert \leq R}u\left(  \xi x\right)
f\left(  \left\vert x\right\vert \right)  dx=\int\limits_{r\leq\left\vert
x\right\vert \leq R}u\left(  \left\vert \xi\right\vert x_{k}\right)  f\left(
\left\vert x\right\vert \right)  dx,\quad1\leq k\leq d.\label{Ap016}%
\end{equation}

\item When $d\geq2$, $s\in\mathbb{R}^{1}$ and $1\leq k\leq d$,%
\begin{equation}
\int\limits_{r\leq\left\vert x\right\vert \leq R}u\left(  sx_{k}\right)
f\left(  \left\vert x\right\vert \right)  dx=\omega_{d-1}\int_{r}^{R}\left(
\int_{0}^{\pi}u\left(  s\rho\cos t\right)  \sin^{d-2}tdt\right)  f\left(
\rho\right)  \rho^{d-1}d\rho,\label{Ap001}%
\end{equation}

where $\omega_{k}=\frac{2\pi^{k/2}}{\Gamma\left(  k/2\right)  }$ for $k\geq1$.
\end{enumerate}
\end{theorem}

\begin{theorem}
\label{Thm_Integ_u(xyhat)f(|x|)dx}Suppose $\xi,x\in\mathbb{R}^{d}$ and $\xi
x=\left(  \xi,x\right)  $ denotes the Euclidean inner product. Then:

\begin{enumerate}
\item $\int_{\left\vert x\right\vert \leq r}u\left(  \xi\widehat{x}\right)
f\left(  \left\vert x\right\vert \right)  dx$ is a radial function of $\xi$
and%
\[
\int_{\left\vert x\right\vert \leq r}u\left(  \xi\widehat{x}\right)  f\left(
\left\vert x\right\vert \right)  dx=\int_{\left\vert x\right\vert \leq
r}u\left(  \left\vert \xi\right\vert \widehat{x}_{k}\right)  f\left(
\left\vert x\right\vert \right)  dx,\quad1\leq k\leq d.
\]

\item When $d\geq2$, $s\in\mathbb{R}^{1}$ and $1\leq k\leq d$,%
\[
\int\limits_{\left\vert x\right\vert \leq r}u\left(  s\widehat{x}_{k}\right)
f\left(  \left\vert x\right\vert \right)  dx=\omega_{d-1}\left(  \int_{0}%
^{\pi}u\left(  s\cos t\right)  \sin^{d-2}tdt\right)  \int_{0}^{r}f\left(
\rho\right)  \rho^{d-1}d\rho,
\]

where $\omega_{k}=\frac{2\pi^{k/2}}{\Gamma\left(  k/2\right)  }$ for $k\geq1$.
\end{enumerate}
\end{theorem}

\begin{theorem}
\label{Thm_Integ_u(xy,|x|)dx}When $d\geq2$, from Prudnikov \cite{PrudBryMar86}%
,%
\begin{align*}
\int_{\left\vert x\right\vert \leq r}\Phi\left(  \left\vert x\right\vert
^{2},\xi x\right)  dx  & =\omega_{d-1}\int_{0}^{r}\rho^{d-1}\int_{0}^{\pi}%
\Phi\left(  \rho^{2},\left\vert \xi\right\vert \rho\cos\theta\right)
\sin^{d-2}\theta d\theta d\rho,\\
\omega_{t}  & :=\frac{2\pi^{t/2}}{\Gamma\left(  t/2\right)  },\text{ }t>0.
\end{align*}

Formula 858.515 of Dwight \cite{Dwight61} is%
\[%
\begin{array}
[c]{c}%
\int\limits_{0}^{\pi/2}\sin^{p}\theta\cos^{q}\theta d\theta=\frac{1}{2}%
\frac{\Gamma\left(  \frac{p+1}{2}\right)  \Gamma\left(  \frac{q+1}{2}\right)
}{\Gamma\left(  \frac{p+q}{2}+1\right)  }=\frac{1}{2}B\left(  \frac{p+1}%
{2},\frac{q+1}{2}\right)  ,\\
where\text{ }p,q>-\frac{1}{2}.
\end{array}
\]

so that%
\[
\int_{0}^{\pi}\sin^{d-2}\theta d\theta=\frac{\omega_{d}}{\omega_{d-1}},\quad
d\geq2.
\]

Also from Prudnikov \cite{PrudBryMar86},%
\[
\int_{\left\vert x\right\vert \leq r}\Phi\left(  \left\vert x\right\vert
^{2},\xi x,\eta x\right)  dx=\omega_{d-2}\int\limits_{-r}^{r}\int%
\limits_{-\sqrt{r^{2}-t^{2}}}^{\sqrt{r^{2}-t^{2}}}\int\limits_{0}^{\sqrt
{r^{2}-t^{2}-u^{2}}}v^{d-3}\Phi\left(  t^{2}+u^{2}+v^{2},\left\vert
\xi\right\vert t,\beta t+\gamma u\right)  \text{ }dvdudt,
\]

where%
\begin{align*}
\beta & =\widehat{\xi}\eta=\left\vert \eta\right\vert \widehat{\xi
}\widehat{\eta},\\
\gamma & =\frac{\left\vert \xi\right\vert ^{2}\left\vert \eta\right\vert
^{2}-\left\vert \xi\eta\right\vert ^{2}}{\left\vert \xi\right\vert ^{2}%
}=\left\vert \eta\right\vert ^{2}\left(  1-\left\vert \widehat{\xi
}\widehat{\eta}\right\vert ^{2}\right)  =\left\vert \eta\right\vert ^{2}%
-\beta^{2}.
\end{align*}

\end{theorem}

\begin{corollary}
\label{Cor_Thm_IntegBr_exp(ixy)f(|x|)dx}Suppose $d\geq2$. Then for $\xi
\in\mathbb{R}^{d}$ and $1\leq k\leq d$:
\begin{align*}
\left(  2\pi\right)  ^{-d/2}\int\limits_{r\leq\left\vert x\right\vert \leq
R}e^{i\xi x}f\left(  \left\vert x\right\vert \right)  dx  & =\left(
2\pi\right)  ^{-d/2}\int\limits_{r\leq\left\vert x\right\vert \leq
R}e^{i\left\vert \xi\right\vert x_{k}}f\left(  \left\vert x\right\vert
\right)  dx\\
& =\left(  2\pi\right)  ^{-d/2}\Gamma\left(  \frac{d}{2}\right)
\int\limits_{r\leq\left\vert x\right\vert \leq R}\frac{J_{\frac{d-2}{2}%
}\left(  \left\vert \xi\right\vert \left\vert x\right\vert \right)  }{\left(
\frac{\left\vert \xi\right\vert \left\vert x\right\vert }{2}\right)
^{\frac{d-2}{2}}}f\left(  \left\vert x\right\vert \right)  dx=\\
& =\left(  2\pi\right)  ^{-d/2}\Gamma\left(  \frac{d}{2}\right)  \omega
_{d}\int\limits_{r\leq\left\vert x\right\vert \leq R}\frac{J_{\frac{d-2}{2}%
}\left(  \left\vert \xi\right\vert r\right)  }{\left(  \frac{\left\vert
\xi\right\vert r}{2}\right)  ^{\frac{d-2}{2}}}r^{d-1}f\left(  r\right)  dr\\
& =\left(  2\pi\right)  ^{-d/2}2\pi^{d/2}\int\limits_{r\leq\left\vert
x\right\vert \leq R}\frac{J_{\frac{d-2}{2}}\left(  \left\vert \xi\right\vert
r\right)  }{\left(  \frac{\left\vert \xi\right\vert r}{2}\right)  ^{\frac
{d-2}{2}}}r^{d-1}f\left(  r\right)  dr\\
& =\frac{1}{2^{\left(  d-2\right)  /2}}\int\limits_{r\leq\left\vert
x\right\vert \leq R}\frac{J_{\frac{d-2}{2}}\left(  \left\vert \xi\right\vert
r\right)  }{\left(  \frac{\left\vert \xi\right\vert r}{2}\right)  ^{\frac
{d-2}{2}}}r^{d-1}f\left(  r\right)  dr\\
& =\int\limits_{r\leq\left\vert x\right\vert \leq R}\frac{J_{\frac{d-2}{2}%
}\left(  \left\vert \xi\right\vert r\right)  }{\left(  \left\vert
\xi\right\vert r\right)  ^{\frac{d-2}{2}}}r^{d-1}f\left(  r\right)  dr,
\end{align*}

where $\omega_{d}=\frac{2\pi^{d/2}}{\Gamma\left(  d/2\right)  }$. ?? If $f\in
L^{1}$ then this is the Fourier transform.
\end{corollary}

\begin{corollary}
\label{Cor_Thm_Integ_u(xy)f(|x|)dx}If the integral exists then for $R\geq
r\geq0$:

\begin{enumerate}
\item
\begin{equation}
\int_{r\leq\left\vert x\right\vert \leq R}f\left(  \left\vert x\right\vert
\right)  dx=\omega_{d}\int_{r}^{R}t^{d-1}f\left(  t\right)  dt,\label{Ap148}%
\end{equation}

where $\omega_{d}=\frac{2\pi^{d/2}}{\Gamma\left(  d/2\right)  }$.

\item Assuming $p>-1/2$,
\begin{equation}
\int\limits_{\left\vert x\right\vert \leq r}\left\vert x_{k}\right\vert
^{p}f\left(  \left\vert x\right\vert \right)  dx=\frac{B\left(  \frac{d-1}%
{2},\frac{p+1}{2}\right)  }{B\left(  \frac{d-1}{2},\frac{1}{2}\right)  }%
\omega_{d}\int_{0}^{r}t^{p+d-1}f\left(  t\right)  dt.\label{Ap017}%
\end{equation}

\end{enumerate}
\end{corollary}

?? ADD ref!

\begin{theorem}
\label{Thm_convol_rad_fn}\textbf{Convolution of radial functions} The
expression $\int\limits_{\left\vert x\right\vert \leq R}f\left(  \left\vert
x\right\vert \right)  g\left(  \left\vert y-x\right\vert \right)  dx$, is an
even function of $y$ when $d=1$ and a radial function of $y$ when $d>1$. In
fact, if $f,g\in L^{1}$,
\begin{align*}
&  \int\limits_{\left\vert x\right\vert \leq R}f\left(  \left\vert
x\right\vert \right)  g\left(  \left\vert y-x\right\vert \right)  dx\\
&  =\left\{
\begin{array}
[c]{ll}%
\int\limits_{0}^{R}f\left(  \left\vert x\right\vert \right)  \left(  g\left(
\left\vert \left\vert y\right\vert +x\right\vert \right)  +g\left(  \left\vert
\left\vert y\right\vert -x\right\vert \right)  \right)  dx, & d=1,\\
\left\{
\begin{array}
[c]{l}%
\omega_{d-1}\int\limits_{0}^{\pi}\int\limits_{0}^{R}\rho^{d-1}f\left(
\rho\right)  g\left(  \sqrt{\left\vert y\right\vert ^{2}-2\left\vert
y\right\vert \rho\cos\theta+\rho^{2}}\right)  \sin^{d-2}\theta\text{ }d\rho
d\theta,\\
\omega_{d-1}\int\limits_{-1}^{1}\int\limits_{0}^{R}\rho^{d-1}f\left(
\rho\right)  g\left(  \sqrt{\left\vert y\right\vert ^{2}-2\left\vert
y\right\vert \rho t+\rho^{2}}\right)  d\rho\text{ }\left(  1-t^{2}\right)
^{\frac{d-3}{2}}dt,\\
\omega_{d-1}\left\vert y\right\vert ^{d}\int\limits_{-1}^{1}\int%
\limits_{0}^{R/\left\vert y\right\vert }\rho^{d-1}f\left(  \left\vert
y\right\vert \rho\right)  g\left(  \left\vert y\right\vert \sqrt{1-2\rho
t+\rho^{2}}\right)  d\rho\text{ }\left(  1-t^{2}\right)  ^{\frac{d-3}{2}}dt,
\end{array}
\right.  & d\geq2,
\end{array}
\right.
\end{align*}

where%
\[
\omega_{k}=\frac{2\pi^{k/2}}{\Gamma\left(  k/2\right)  },\quad k\geq1.
\]

\end{theorem}

\ 

\begin{corollary}
\label{Cor_Thm_convol_rad_fn}Suppose $f,g\in L^{1}$ are both radial, say
$f_{\odot}\left(  \left\vert x\right\vert \right)  :=f\left(  x\right)  $ and
$g_{\odot}\left(  \left\vert x\right\vert \right)  :=g\left(  x\right)  $.
Then%
\begin{align*}
&  \int\limits_{\left\vert x\right\vert \leq R}f\left(  x\right)  g\left(
y-x\right)  dx\\
&  =\frac{\omega_{d-1}}{2^{d-3}}\frac{1}{\left\vert y\right\vert ^{d-2}}%
\times\left\{
\begin{array}
[c]{l}%
\int\limits_{0}^{R}\rho f_{\odot}\left(  \rho\right)  \int\limits_{\left(
\left\vert y\right\vert -\rho\right)  ^{2}}^{\left(  \left\vert y\right\vert
+\rho\right)  ^{2}}g_{\odot}\left(  \sqrt{s}\right)  \left(  s-\left(
\left\vert y\right\vert -\rho\right)  ^{2}\right)  ^{\frac{d-3}{2}}\left(
\left(  \left\vert y\right\vert +\rho\right)  ^{2}-s\right)  ^{\frac{d-3}{2}%
}ds\text{ }d\rho,\\
or\\
\int\limits_{0}^{R}\rho f_{\odot}\left(  \rho\right)  \int\limits_{\left\vert
\left\vert y\right\vert -\rho\right\vert }^{\left\vert y\right\vert +\rho
}g_{\odot}\left(  t\right)  \left(  t^{2}-\left(  \left\vert y\right\vert
-\rho\right)  ^{2}\right)  ^{\frac{d-3}{2}}\left(  \left(  \left\vert
y\right\vert +\rho\right)  ^{2}-t^{2}\right)  ^{\frac{d-3}{2}}tdt\text{ }%
d\rho,
\end{array}
\right.
\end{align*}

if one of the integrals exists. Further%
\begin{align*}
&  \left(  \int f\left(  x\right)  g\left(  \cdot-x\right)  dx\right)
_{\odot}\left(  r\right) \\
&  =\frac{\omega_{d-1}}{2^{d-3}}\frac{1}{r^{d-2}}\left(
\begin{array}
[c]{l}%
\int\limits_{0}^{r}\rho f_{\odot}\left(  \rho\right)  \int\limits_{r-\rho
}^{r+\rho}tg_{\odot}\left(  t\right)  \left(  t^{2}-\left(  r-\rho\right)
^{2}\right)  ^{\frac{d-3}{2}}\left(  \left(  r+\rho\right)  ^{2}-t^{2}\right)
^{\frac{d-3}{2}}dt\text{ }d\rho+\\
\quad+\int\limits_{r}^{\infty}\rho f_{\odot}\left(  \rho\right)
\int\limits_{\rho-r}^{\rho+r}tg_{\odot}\left(  t\right)  \left(  t^{2}-\left(
\rho-r\right)  ^{2}\right)  ^{\frac{d-3}{2}}\left(  \left(  \rho+r\right)
^{2}-t^{2}\right)  ^{\frac{d-3}{2}}dt\text{ }d\rho
\end{array}
\right)
\end{align*}

\end{corollary}

Define%
\[
\theta\left(  r,\rho\right)  :=\int\limits_{r-\rho}^{r+\rho}tg_{\odot}\left(
t\right)  \left(  t^{2}-\left(  r-\rho\right)  ^{2}\right)  ^{c}\left(
\left(  r+\rho\right)  ^{2}-t^{2}\right)  ^{c}dt,
\]

so that%
\[
\left(  \int f\left(  x\right)  g\left(  \cdot-x\right)  dx\right)  _{\odot
}\left(  r\right)  =\frac{\omega_{d-1}}{2^{d-3}}\frac{1}{r^{d-2}}\left(
\int\limits_{0}^{r}\rho f_{\odot}\left(  \rho\right)  \theta\left(
r,\rho\right)  d\rho+\int\limits_{r}^{\infty}\rho f_{\odot}\left(
\rho\right)  \theta\left(  \rho,r\right)  d\rho\right)  .
\]

We first remove the variables $\rho$ and $r$ from the range of the inner
integral by using the rectangular function $\Pi$ with support $\left[
-1,1\right]  $ and non-zero value $1$:%
\begin{align*}
\theta\left(  r,\rho\right)   &  =\int_{r-\rho}^{r+\rho}tg_{\odot}\left(
t\right)  \left(  t^{2}-\left(  r-\rho\right)  ^{2}\right)  ^{c}\left(
\left(  r+\rho\right)  ^{2}-t^{2}\right)  ^{c}dt\\
&  =\int\Pi\left(  \frac{t-r}{\rho}\right)  tg_{\odot}\left(  t\right)
\left(  t^{2}-\left(  r-\rho\right)  ^{2}\right)  ^{c}\left(  \left(
r+\rho\right)  ^{2}-t^{2}\right)  ^{c}dt\\
&  =\int\Pi\left(  \frac{t-r}{\rho}\right)  tg_{\odot}\left(  t\right)
\left(  t-\left(  r-\rho\right)  \right)  ^{c}\left(  t+\left(  r-\rho\right)
\right)  ^{c}\left(  \left(  r+\rho\right)  -t\right)  ^{c}\left(  \left(
r+\rho\right)  +t\right)  ^{c}dt\\
&  =\int\Pi\left(  \frac{t-r}{\rho}\right)  tg_{\odot}\left(  t\right)
\left(  t-\left(  r-\rho\right)  \right)  ^{c}\left(  \left(  r+\rho\right)
-t\right)  ^{c}\left\{  \left(  t+\left(  r-\rho\right)  \right)  ^{c}\left(
\left(  r+\rho\right)  +t\right)  ^{c}\right\}  dt\\
&  =\rho^{2c}\int\Pi\left(  \frac{t-r}{\rho}\right)  tg_{\odot}\left(
t\right)  \left(  \frac{t-r}{\rho}+1\right)  ^{c}\left(  1-\frac{t-r}{\rho
}\right)  ^{c}\left\{  \left(  t-\rho+r\right)  ^{c}\left(  t+\rho+r\right)
^{c}\right\}  dt.
\end{align*}

Now apply the change of variables: $s=\frac{t-r}{\rho}$, $\rho ds=dt$ to get%
\begin{align*}
\theta\left(  r,\rho\right)   &  =\rho^{d-2}\int_{-1}^{1}\left(  \rho
s+r\right)  g_{\odot}\left(  \rho s+r\right)  \left(  s+1\right)  ^{c}\left(
1-s\right)  ^{c}\left\{  \left(  \rho s+r-\rho+r\right)  ^{c}\left(  \rho
s+r+\rho+r\right)  ^{c}\right\}  ds\\
&  =\rho^{d-2}\int_{-1}^{1}\left(  1-s^{2}\right)  ^{c}\left(  \rho
s+r\right)  g_{\odot}\left(  \rho s+r\right)  \left\{  \left(  \rho\left(
s-1\right)  +2r\right)  ^{c}\left(  \rho\left(  s+1\right)  +2r\right)
^{c}\right\}  ds.
\end{align*}

The next step is to apply the binomial theorem to express the right side as a
polynomial in $r$.%
\begin{align*}
&  \theta\left(  r,\rho\right) \\
&  =\rho^{d-2}\int_{-1}^{1}\left(  1-s^{2}\right)  ^{c}\left(  \rho
s+r\right)  g_{\odot}\left(  \rho s+r\right)  \sum_{i_{1}=0}^{c}\tbinom
{c}{i_{1}}\left(  \rho\left(  s-1\right)  \right)  ^{c-i_{1}}\left(
2r\right)  ^{i_{1}}\sum_{i_{2}=0}^{c}\tbinom{c}{i_{2}}\left(  \rho\left(
s+1\right)  \right)  ^{c-i_{2}}\left(  2r\right)  ^{i_{2}}ds\\
&  =\rho^{d-2}\sum_{i_{1},i_{2}=0}^{c}\tbinom{c}{i_{1}}\tbinom{c}{i_{2}%
}\left(  \int_{-1}^{1}\left(  1-s^{2}\right)  ^{c}\left(  \rho s+r\right)
g_{\odot}\left(  \rho s+r\right)  \left(  \rho\left(  s-1\right)  \right)
^{c-i_{1}}\left(  \rho\left(  s+1\right)  \right)  ^{c-i_{2}}ds\right)
\left(  2r\right)  ^{i_{1}}\left(  2r\right)  ^{i_{2}}\\
&  =\rho^{d-2}\sum_{i_{1},i_{2}=0}^{c}\tbinom{c}{i_{1}}\tbinom{c}{i_{2}}%
\rho^{2c-i_{1}-i_{2}}\left(  \int_{-1}^{1}\left(  1-s^{2}\right)  ^{c}\left(
\rho s+r\right)  g_{\odot}\left(  \rho s+r\right)  \left(  s-1\right)
^{c-i_{1}}\left(  s+1\right)  ^{c-i_{2}}ds\right)  \left(  2r\right)
^{i_{1}+i_{2}}.
\end{align*}

When $\operatorname*{supp}f\subset\overline{B}_{R}$ the estimate%
\begin{align*}
\left\vert \theta\left(  r,\rho\right)  \right\vert  & \leq\rho^{d-2}%
\sum_{i_{1},i_{2}=0}^{c}\tbinom{c}{i_{1}}\tbinom{c}{i_{2}}\rho^{2c-i_{1}%
-i_{2}}4\left(  \rho+r\right)  ^{2n-d+1}2^{c-i_{1}}2^{c-i_{2}}\left(
2r\right)  ^{i_{1}+i_{2}}\\
& =42^{2c}\left(  \rho+r\right)  ^{2n-d+1}\rho^{d-2}\sum_{i_{1},i_{2}=0}%
^{c}\tbinom{c}{i_{1}}\tbinom{c}{i_{2}}\rho^{2c-i_{1}-i_{2}}r^{i_{1}+i_{2}}\\
& =2^{2c+2}\left(  \rho+r\right)  ^{2n+2c-d+1}\rho^{d-2},
\end{align*}

implies%
\begin{align*}
\left\vert \left(  \int f\left(  x\right)  g\left(  \cdot-x\right)  dx\right)
_{\odot}\left(  r\right)  \right\vert  & \leq\frac{\omega_{d-1}}{2^{d-3}}%
\frac{1}{r^{d-2}}\left\vert \int\limits_{0}^{R}\rho f_{\odot}\left(
\rho\right)  \theta\left(  r,\rho\right)  d\rho\right\vert \\
& \leq\frac{\omega_{d-1}}{2^{d-3}}\frac{1}{r^{d-2}}\int\limits_{0}^{R}%
\rho\left\vert f_{\odot}\left(  \rho\right)  \right\vert \left\vert
\theta\left(  r,\rho\right)  \right\vert d\rho\\
& \leq\frac{\omega_{d-1}}{2^{d-3}}\frac{1}{r^{d-2}}\int\limits_{0}^{R}%
\rho\left\vert f_{\odot}\left(  \rho\right)  \right\vert 2^{2c+2}\left(
\rho+r\right)  ^{2n+2c-d+1}\rho^{d-2}d\rho\\
& =\frac{\omega_{d-1}}{2^{d-3}}\frac{1}{r^{d-2}}\int\limits_{0}^{R}%
\rho\left\vert f_{\odot}\left(  \rho\right)  \right\vert 2^{2c+2}\left(
\rho+r\right)  ^{2n-2}\rho^{d-2}d\rho\\
& =4\omega_{d-1}\frac{1}{r^{d-2}}\int\limits_{0}^{R}\rho^{d-1}\left\vert
f_{\odot}\left(  \rho\right)  \right\vert \left(  \rho+r\right)  ^{2n-2}%
d\rho\\
& <4\omega_{d-1}\frac{\left(  R+r\right)  ^{2n-2}}{r^{d-2}}\int\limits_{0}%
^{R}\rho^{d-1}\left\vert f_{\odot}\left(  \rho\right)  \right\vert d\rho,
\end{align*}

and note that $2n-d>0$.

\section{Some (combinatorial) binomial sums}

From ??%
\[
\sum\limits_{j=-l}^{l}\left(  -1\right)  ^{j}\tbinom{2l}{j+l}j^{k}=0.\quad
k=0,1,2,\ldots,??
\]

?? Derive using complex variable techniques e.g. Section I.8 of Zwillinger
\cite{Zwill92}?%

\begin{align}
\sum\limits_{j=1}^{l}\left(  -1\right)  ^{j}\tbinom{2l}{j+l}j^{2m}  &
=\sum\limits_{j=1}^{l}\left(  -1\right)  ^{j}\tbinom{2l}{l-j}j^{2m}\nonumber\\
& =\sum\limits_{k=0}^{l-1}\left(  -1\right)  ^{l-k}\tbinom{2l}{k}\left(
l-k\right)  ^{2m}\nonumber\\
& =\left(  -1\right)  ^{l}\sum\limits_{k=0}^{l-1}\left(  -1\right)
^{k}\tbinom{2l}{k}\left(  l-k\right)  ^{2m}\nonumber\\
& =\left(  -1\right)  ^{l}\sum\limits_{k=0}^{l}\left(  -1\right)  ^{k}%
\tbinom{2l}{k}\left(  l-k\right)  ^{2m}\nonumber\\
& =\left\{
\begin{array}
[c]{ll}%
0, & 1\leq m\leq l-1,\\
\left(  -1\right)  ^{l}\frac{\left(  2l\right)  !}{2}, & m=l,
\end{array}
\right.  ,\label{a020}%
\end{align}

from formulas 33 and 34 of Subsection 4.2.2 of Prudnikov \cite{PrudBryMar86}.
Also, the identity%
\begin{equation}
\sum\limits_{j=0}^{l-1}\left(  -1\right)  ^{j}\tbinom{2l}{j}=\frac{\left(
-1\right)  ^{l-1}}{2}\tbinom{2l}{l},\label{a019}%
\end{equation}

yields%
\begin{align}
\sum\limits_{j=1}^{l}\left(  -1\right)  ^{j}\tbinom{2l}{j+l}  & =\sum
\limits_{j=1}^{l}\left(  -1\right)  ^{j}\tbinom{2l}{l-j}=\sum\limits_{k=0}%
^{l-1}\left(  -1\right)  ^{l-k}\tbinom{2l}{k}=\nonumber\\
& =\left(  -1\right)  ^{l}\sum\limits_{k=0}^{l-1}\left(  -1\right)
^{k}\tbinom{2l}{k}=\left(  -1\right)  ^{l}\frac{\left(  -1\right)  ^{l-1}}%
{2}\tbinom{2l}{l}=\nonumber\\
& =-\frac{1}{2}\tbinom{2l}{l},\label{a021}%
\end{align}

and so \ref{a020} extends to%
\begin{equation}
\sum\limits_{j=1}^{l}\left(  -1\right)  ^{j}\tbinom{2l}{j+l}j^{2m}=\left\{
\begin{array}
[c]{ll}%
-\frac{1}{2}\tbinom{2l}{l}, & m=0,\\
0, & 1\leq m\leq l-1,\\
\left(  -1\right)  ^{l}\frac{\left(  2l\right)  !}{2}, & m=l.
\end{array}
\right. \label{a022}%
\end{equation}

From formula 32 of Subsection 4.2.2 of Prudnikov \cite{PrudBryMar86},
\begin{equation}
\sum\limits_{k=0}^{l-1}\left(  -1\right)  ^{k}\tbinom{2l}{k}\left(
l-k\right)  ^{s}=\left(  -1\right)  ^{l+1}\frac{1}{\pi}s!2^{2l-s}\sin\frac{\pi
s}{2}\int_{0}^{\infty}\frac{\sin^{2l}t}{t^{s+1}}dt,\quad0<s<2l,\label{a026}%
\end{equation}

\[
\int_{0}^{\infty}\frac{\sin^{2l}t}{t^{s+1}}dt=\left(  -1\right)  ^{l+1}%
\pi\frac{\sum\limits_{k=0}^{l-1}\left(  -1\right)  ^{k}\tbinom{2l}{k}\left(
l-k\right)  ^{s}}{s!2^{2l-s}\sin\frac{\pi s}{2}},\quad0<s<2l.
\]

But%
\begin{align*}
&  \sum\limits_{k=0}^{l-1}\left(  -1\right)  ^{k}\tbinom{2l}{k}\left(
l-k\right)  ^{s}\\
&  :j=l-k;\quad0:l-1\rightarrow l:1\\
&  =\sum\limits_{k=l}^{1}\left(  -1\right)  ^{l-j}\tbinom{2l}{l-j}%
j^{s}=\left(  -1\right)  ^{l}\sum\limits_{j=1}^{l}\left(  -1\right)
^{j}\tbinom{2l}{l-j}j^{s}=\left(  -1\right)  ^{l}\sum\limits_{j=1}^{l}\left(
-1\right)  ^{j}\tbinom{2l}{j+l}j^{s},
\end{align*}

and so%
\begin{align}
\left(  -1\right)  ^{l}\sum\limits_{j=1}^{l}\left(  -1\right)  ^{j}\tbinom
{2l}{j+l}j^{s}  & =\left(  -1\right)  ^{l+1}\frac{1}{\pi}s!2^{2l-s}\sin
\frac{\pi s}{2}\int_{0}^{\infty}\frac{\sin^{2l}t}{t^{s+1}}dt,\label{a027}\\
when\text{ }0  & <s<2l.\nonumber
\end{align}

Thus%
\begin{equation}
\int_{0}^{\infty}\frac{\sin^{2l}t}{t^{s+1}}dt=\frac{-\pi}{s!2^{2l-s}\sin
\frac{\pi s}{2}}\sum\limits_{j=1}^{l}\left(  -1\right)  ^{j}\tbinom{2l}%
{j+l}j^{s},\quad0<s<2l.\label{a053}%
\end{equation}

Clearly this formula implies
\begin{equation}
\operatorname{sgn}\sum\limits_{j=1}^{l}\left(  -1\right)  ^{j}\tbinom{2l}%
{j+l}j^{s}=-\operatorname{sgn}\sin\frac{\pi s}{2},\quad1\leq n\leq l,\text{
}l\geq2,\label{a056}%
\end{equation}

and when $s=2n-1$,%
\begin{equation}
\int_{0}^{\infty}\frac{\sin^{2l}t}{t^{2n}}dt=\frac{\left(  -1\right)  ^{n}\pi
}{\left(  2n-1\right)  !2^{2l-2n+1}}\sum\limits_{j=1}^{l}\left(  -1\right)
^{j}\tbinom{2l}{j+l}j^{2n-1},\quad1\leq n\leq l,\label{a058}%
\end{equation}

or%
\begin{equation}
\int_{0}^{\infty}\frac{\sin^{2l}t}{t^{m}}dt=\frac{\left(  -1\right)
^{\frac{m}{2}}\pi}{\left(  m-1\right)  !2^{2l-m+1}}\sum\limits_{j=1}%
^{l}\left(  -1\right)  ^{j}\tbinom{2l}{j+l}j^{m-1},\quad m=2,4,\ldots
,2l.\label{a052}%
\end{equation}

Suppose $s=2n+u$ in \ref{a053}. Then for $1\leq n<l$,%
\[
\int_{0}^{\infty}\frac{\sin^{2l}t}{t^{2n+1+u}}dt=\frac{-\pi}{\left(
2n+u\right)  !2^{2l-2n-u}}\frac{\sum\limits_{j=1}^{l}\left(  -1\right)
^{j}\tbinom{2l}{j+l}j^{2n+u}}{\sin\frac{\pi\left(  2n+u\right)  }{2}}%
,\quad0<s<2l.
\]

and using L'H\^{o}pital's rule,%
\begin{align*}
\lim_{u\rightarrow0}\frac{\sum\limits_{j=1}^{l}\left(  -1\right)  ^{j}%
\tbinom{2l}{j+l}j^{2n+u}}{\sin\frac{\pi\left(  2n+u\right)  }{2}}  &
=\lim_{u\rightarrow0}\frac{\sum\limits_{j=1}^{l}\left(  -1\right)  ^{j}%
\tbinom{2l}{j+l}j^{2n}j^{u}}{\sin\left(  \pi n+\pi\frac{u}{2}\right)  }%
=\lim_{u\rightarrow0}\frac{\sum\limits_{j=1}^{l}\left(  -1\right)  ^{j}%
\tbinom{2l}{j+l}j^{2n}j^{u}}{\left(  -1\right)  ^{n}\sin\frac{\pi u}{2}}=\\
& =\lim_{u\rightarrow0}\frac{\sum\limits_{j=1}^{l}\left(  -1\right)
^{j}\tbinom{2l}{j+l}j^{2n}j^{u}\ln j}{\left(  -1\right)  ^{n}\frac{\pi}{2}%
\cos\frac{u}{2}}=\lim_{u\rightarrow0}\frac{\sum\limits_{j=1}^{l}\left(
-1\right)  ^{j}\tbinom{2l}{j+l}j^{2n}\ln j}{\left(  -1\right)  ^{n}\frac{\pi
}{2}}=\\
& =\left(  -1\right)  ^{n}\frac{2}{\pi}\sum\limits_{j=1}^{l}\left(  -1\right)
^{j}\tbinom{2l}{j+l}j^{2n}\ln j,
\end{align*}

and so%
\begin{align*}
\int_{0}^{\infty}\frac{\sin^{2l}t}{t^{2n+1}}dt  & =\lim_{u\rightarrow0}%
\frac{-\pi}{\left(  2n+u\right)  !2^{2l-2n-u}}\frac{\sum\limits_{j=1}%
^{l}\left(  -1\right)  ^{j}\tbinom{2l}{j+l}j^{2n+u}}{\sin\frac{\pi\left(
2n+u\right)  }{2}}\\
& =\frac{-\pi}{\left(  2n\right)  !2^{2l-2n}}\lim_{u\rightarrow0}\frac
{\sum\limits_{j=1}^{l}\left(  -1\right)  ^{j}\tbinom{2l}{j+l}j^{2n+u}}%
{\sin\frac{\pi\left(  2n+u\right)  }{2}}\\
& =\frac{-\pi}{\left(  2n\right)  !2^{2l-2n}}\left(  -1\right)  ^{n}\frac
{2}{\pi}\sum\limits_{j=1}^{l}\left(  -1\right)  ^{j}\tbinom{2l}{j+l}j^{2n}\ln
j\\
& =\frac{\left(  -1\right)  ^{n+1}}{\left(  2n\right)  !2^{2l-2n-1}}%
\sum\limits_{j=1}^{l}\left(  -1\right)  ^{j}\tbinom{2l}{j+l}j^{2n}\ln j,
\end{align*}

and hence (checked using Matlab)%
\begin{equation}
\int_{0}^{\infty}\frac{\sin^{2l}t}{t^{2n-1}}dt=\frac{\left(  -1\right)  ^{n}%
}{\left(  2n-2\right)  !2^{2l-2n+1}}\sum\limits_{j=1}^{l}\left(  -1\right)
^{j}\tbinom{2l}{j+l}j^{2n-2}\ln j,\quad2\leq n\leq l,\label{a055}%
\end{equation}

i.e.%
\begin{equation}
\int_{0}^{\infty}\frac{\sin^{2l}t}{t^{m}}dt=\frac{\left(  -1\right)
^{\frac{m+1}{2}}}{\left(  m-1\right)  !2^{2l-m}}\sum\limits_{j=1}^{l}\left(
-1\right)  ^{j}\tbinom{2l}{j+l}j^{m-1}\ln j,\quad m=3,5,\ldots
,2l-1.\label{a057}%
\end{equation}

We summarize \ref{a057} and \ref{a052} as,%
\begin{equation}
\int_{0}^{\infty}\frac{\sin^{2l}t}{t^{m}}dt=\tfrac{1}{\left(  m-1\right)
!2^{2l-m+1}}\times\left\{
\begin{array}
[c]{ll}%
\left(  -1\right)  ^{\frac{m+1}{2}}2\sum\limits_{j=1}^{l}\left(  -1\right)
^{j}\tbinom{2l}{j+l}j^{m-1}\ln j, & m=3,5,\ldots,2l-1,\\
\left(  -1\right)  ^{\frac{m}{2}}\pi\sum\limits_{j=1}^{l}\left(  -1\right)
^{j}\tbinom{2l}{j+l}j^{m-1}, & m=2,4,\ldots,2l.
\end{array}
\right. \label{a059}%
\end{equation}

This can be written: \textbf{if }$p$\textbf{\ is even} then%
\[
\int_{0}^{\infty}\frac{\sin^{p}t}{t^{m}}dt=\tfrac{1}{\left(  m-1\right)
!2^{p-m+1}}\times\left\{
\begin{array}
[c]{ll}%
\left(  -1\right)  ^{\frac{m+1}{2}}2\sum\limits_{j=1}^{p/2}\left(  -1\right)
^{j}\tbinom{p}{j+p/2}j^{m-1}\ln j, & m=3,5,\ldots,p-1,\\
\left(  -1\right)  ^{\frac{m}{2}}\pi\sum\limits_{j=1}^{p/2}\left(  -1\right)
^{j}\tbinom{p}{j+p/2}j^{m-1}, & m=2,4,\ldots,p.
\end{array}
\right.
\]

\section{Approximating the gamma function and binomials}

A well known gamma function approximation is
\begin{equation}
\sqrt{2\pi}e^{-x+\frac{1}{1+12x}}x^{x+\frac{1}{2}}<\Gamma\left(  x+1\right)
=x!<\sqrt{2\pi}e^{-x+\frac{1}{12x}}x^{x+\frac{1}{2}},\quad x\geq
0,\label{Ap029}%
\end{equation}

\begin{align}
1  & \leq\frac{\sin x}{x}\leq\frac{2}{\pi},\quad\quad\quad0\leq x\leq
\pi/2.\label{Ap000}\\
1-\frac{1}{6}x^{2}  & \leq\frac{\sin x}{x}\leq1-\frac{1}{7}x^{2},\quad0\leq
x\leq\pi/2.\label{Ap002}%
\end{align}

and also that%
\begin{equation}
\Gamma\left(  x\right)  \Gamma\left(  1-x\right)  =\frac{\pi}{\sin\pi x}%
,\quad\Gamma\left(  1/2\right)  =\sqrt{\pi}.\label{Ap039}%
\end{equation}

Now $\Gamma\left(  x+1\right)  =x\Gamma\left(  x\right)  $ so \ref{Ap029}
becomes%
\begin{equation}
\sqrt{2\pi}e^{-x+\frac{1}{1+12x}}x^{x-\frac{1}{2}}<\Gamma\left(  x\right)
<\sqrt{2\pi}e^{-x+\frac{1}{12x}}x^{x-\frac{1}{2}},\quad x>0.\label{Ap049}%
\end{equation}

The estimates \ref{Ap049} are OK for $x\geq1/2$ ?? SUPPLY\ NUM\ EVIDENCE?.

From \ref{Ap029},%
\begin{align*}
\frac{\left(  2n\right)  !}{n!n!}  & <\frac{\sqrt{2\pi}e^{-2n+\frac{1}{24n}%
}\left(  2n\right)  ^{2n+\frac{1}{2}}}{\left(  \sqrt{2\pi}e^{-n+\frac
{1}{1+12n}}n^{n+\frac{1}{2}}\right)  ^{2}}=\frac{\sqrt{2\pi}e^{-2n+\frac
{1}{24n}}\left(  2n\right)  ^{2n+\frac{1}{2}}}{\left(  \sqrt{2\pi}\right)
^{2}e^{-2n+\frac{2}{1+12n}}n^{2n+1}}=\\
& =\frac{e^{\frac{1}{24n}}\left(  2n\right)  ^{2n+\frac{1}{2}}}{\sqrt{2\pi
}e^{\frac{2}{1+12n}}n^{2n+1}}=\frac{e^{\frac{1}{24n}}2^{2n+\frac{1}{2}%
}n^{2n+\frac{1}{2}}}{\sqrt{2\pi}e^{\frac{2}{1+12n}}n^{2n+1}}=\frac{e^{\frac
{1}{24n}}2^{2n}}{\sqrt{\pi}e^{\frac{2}{1+12n}}\sqrt{n}}=\frac{e^{\frac{1}%
{24n}-\frac{2}{1+12n}}2^{2n}}{\sqrt{\pi}\sqrt{n}}<\\
& <\frac{2^{2n}}{\sqrt{\pi}\sqrt{n}},
\end{align*}

and so%
\[
2^{-2n}\binom{2n}{n}<\frac{1}{\sqrt{\pi n}}.
\]

Also%
\begin{align*}
\frac{\left(  2n\right)  !}{n!n!}  & >\frac{\sqrt{2\pi}e^{-2n+\frac{1}{1+24n}%
}\left(  2n\right)  ^{2n+\frac{1}{2}}}{\left(  \sqrt{2\pi}e^{-n+\frac{1}{12n}%
}n^{n+\frac{1}{2}}\right)  ^{2}}=\frac{\sqrt{2\pi}e^{-2n+\frac{1}{1+24n}%
}2^{2n+\frac{1}{2}}n^{2n+\frac{1}{2}}}{\left(  \sqrt{2\pi}\right)
^{2}e^{-2n+\frac{1}{6n}}n^{2n+1}}=\\
& =\frac{e^{\frac{1}{1+24n}}2^{2n+\frac{1}{2}}n^{\frac{1}{2}}}{\sqrt{2\pi
}e^{\frac{1}{6n}}n}=\frac{e^{\frac{1}{1+24n}-\frac{1}{6n}}2^{2n+\frac{1}{2}}%
}{\sqrt{2\pi n}}=\frac{2^{2n}}{\sqrt{\pi n}}e^{-\left(  \frac{1}{6n}-\frac
{1}{1+24n}\right)  }>\\
& >\frac{2^{2n}}{\sqrt{\pi n}}e^{-\frac{1}{6n}},
\end{align*}

so that%
\begin{equation}
\frac{^{e^{-\frac{1}{6n}}}}{\sqrt{\pi n}}<2^{-2n}\binom{2n}{n}<\frac{1}%
{\sqrt{\pi n}}.\label{Ap004}%
\end{equation}

Also%
\[
\binom{2n-1}{n-1}=\frac{\left(  2n-1\right)  !}{\left(  n-1\right)  !n!}%
=\frac{n}{2n}\frac{\left(  2n\right)  !}{n!n!}=\frac{1}{2}\frac{\left(
2n\right)  !}{n!n!},
\]

so%
\begin{equation}
\frac{^{e^{-\frac{1}{6n}}}}{\sqrt{\pi n}}<2^{-2n+1}\binom{2n-1}{n-1}<\frac
{1}{\sqrt{\pi n}}.\label{Ap005}%
\end{equation}
\medskip

\fbox{For $0<x<1/2$} we can use \ref{Ap039}. In fact%
\begin{align*}
x\Gamma\left(  x\right)   & =\frac{\pi x}{\sin\pi x}\frac{1}{\Gamma\left(
1-x\right)  }\leq\frac{\pi}{2}\frac{1}{\Gamma\left(  1-x\right)  }\leq
\frac{\pi}{2}\frac{1}{\Gamma\left(  1\right)  }=\frac{\pi}{2}.\\
x\Gamma\left(  x\right)   & =\frac{\pi x}{\sin\pi x}\frac{1}{\Gamma\left(
1-x\right)  }\geq\frac{1}{\Gamma\left(  1-x\right)  }\geq\frac{1}{\sqrt{\pi}}.
\end{align*}

From \ref{Ap049},%
\[
\frac{1}{\Gamma\left(  1-x\right)  }\leq\frac{1}{\sqrt{2\pi}e^{-\left(
1-x\right)  +\frac{1}{1+12\left(  1-x\right)  }}\left(  1-x\right)  ^{\left(
1-x\right)  -\frac{1}{2}}}\leq\frac{e^{1-x}}{\sqrt{2\pi}\left(  1-x\right)
^{\frac{1}{2}-x}}%
\]

Numerically%
\begin{align*}
x\Gamma\left(  x\right)  \leq\left(  x\Gamma\left(  x\right)  \right)  \left(
0\right)  -\frac{\left(  x\Gamma\left(  x\right)  \right)  \left(  0\right)
-\left(  x\Gamma\left(  x\right)  \right)  \left(  1/2\right)  }{1/2}x  &
=1-\frac{1-\sqrt{\pi}/2}{1/2}x\\
& =1-\left(  2-\sqrt{\pi}\right)  x.
\end{align*}

Also, by adjusting the formula%
\[
1-\gamma x+0.5\left(  \frac{1}{2}\gamma^{2}+\frac{\pi^{2}}{12}\right)
x^{2}\leq x\Gamma\left(  x\right)  \leq1-\gamma x+\left(  \frac{1}{2}%
\gamma^{2}+\frac{\pi^{2}}{12}\right)  x^{2},
\]

by numerical observation we obtain the tighter%
\[
1-\gamma x+0.7\left(  \frac{1}{2}\gamma^{2}+\frac{\pi^{2}}{12}\right)
x^{2}\leq x\Gamma\left(  x\right)  \leq1-\gamma x+\left(  \frac{1}{2}%
\gamma^{2}+\frac{\pi^{2}}{12}\right)  x^{2}-0.525x^{3}.
\]

\section{Generalized Legendre polynomials}

?? ADD blah!

\begin{lemma}
\label{Lem_Legendre}Suppose $\lambda,t\in\mathbb{R}^{1}$, $\lambda
\notin??\mathbb{Z}_{+}^{1}??$ and $0\leq t<1$. Then%
\[
\left\vert 1+te^{i\theta}\right\vert ^{2\lambda}=\sum_{n=0}^{m-1}%
P_{n}^{\left(  \lambda\right)  }\left(  \cos\theta\right)  t^{n}+\sum
_{n=m}^{2m-2}Q_{n}^{\left(  \lambda\right)  }\left(  m,\cos\theta\right)
t^{n}+\mathcal{S}_{m}^{\left(  \lambda\right)  }\left(  te^{i\theta
},te^{-i\theta}\right)  ,
\]

where $P_{n}^{\left(  \lambda\right)  }\left(  \cos\theta\right)  $ is a
polynomial of degree $n$ in $\cos\theta$ given by \ref{7.78}, $Q_{n}^{\left(
\lambda\right)  }\left(  m,\cos\theta\right)  $ is a polynomial of degree
$n-m+1$ in $\cos\theta$ given by \ref{a2.27}, and $\mathcal{S}_{m}^{\left(
\lambda\right)  }$ is given by \ref{7.76} and \ref{7.79}.

The remainder $\mathcal{S}_{m}^{\left(  \lambda\right)  }$ satisfies%
\begin{equation}
\left.
\begin{array}
[c]{ll}%
\left\vert \mathcal{S}_{m}^{\left(  \lambda\right)  }\left(  te^{i\theta
},te^{-i\theta}\right)  \right\vert  & \leq\left\vert \mathcal{R}_{m}^{\left(
\lambda\right)  }\left(  te^{i\theta}\right)  \right\vert \left(  2\left\vert
\sum\limits_{k=0}^{m-1}\tbinom{\lambda}{k}t^{k}\right\vert +\left\vert
\mathcal{R}_{m}^{\left(  \lambda\right)  }\left(  te^{i\theta}\right)
\right\vert \right)  ,\\
\left\vert \mathcal{R}_{m}^{\left(  \lambda\right)  }\left(  te^{i\theta
}\right)  \right\vert  & \leq\left\{
\begin{array}
[c]{ll}%
\tbinom{\lambda}{m}t^{m}\left(  1+t\right)  ^{\lambda-m}, & \lambda\geq m,\\
\tbinom{\lambda}{m}t^{m}\left(  1-t\right)  ^{\lambda-m}, & \lambda\leq m,
\end{array}
\right.
\end{array}
\right\} \label{a50.9}%
\end{equation}

which in turn imply the upper bound%
\begin{align}
\left\vert \mathcal{S}_{m}^{\left(  \lambda\right)  }\left(  te^{i\theta
},te^{-i\theta}\right)  \right\vert  & \leq\left\{
\begin{array}
[c]{ll}%
\tbinom{\lambda}{m}t^{m}\left(  1-t\right)  ^{2\left(  \lambda-m\right)
}\left(  2\left(  1+t\right)  ^{m}+3\tbinom{\lambda}{m}t^{m}\right)  , &
\lambda\leq m,\\
\tbinom{\lambda}{m}t^{m}\left(  1+t\right)  ^{2\lambda-m}\left(
2+3\tbinom{\lambda}{m}t^{m}\right)  , & \lambda\geq m,
\end{array}
\right. \label{a4.9}\\
when\text{ }0  & \leq t<1.\nonumber
\end{align}

\end{lemma}

\begin{proof}
Since%
\[
D^{n}\left(  1+s\right)  ^{\lambda}=\lambda\left(  \lambda-1\right)
\ldots\left(  \lambda-n+1\right)  \left(  1+s\right)  ^{\lambda-n},
\]

we have the Taylor series binomial expansion about $1$:%
\[
\left(  1+z\right)  ^{\lambda}=\sum_{k=0}^{m-1}\tbinom{\lambda}{k}%
z^{k}+\mathcal{R}_{m}^{\left(  \lambda\right)  }\left(  z\right)  ,\quad
z\in\mathbb{C},
\]

where the%
\[
\tbinom{\lambda}{k}=\frac{\lambda\left(  \lambda-1\right)  \ldots\left(
\lambda-k+1\right)  }{k!},
\]

are the generalized binomial coefficients and%
\begin{align}
\mathcal{R}_{m}^{\left(  \lambda\right)  }\left(  z\right)   & =m\tbinom
{\lambda}{m}z^{m}\int_{0}^{1}\left(  1-\tau\right)  ^{m-1}\left(  1+\tau
z\right)  ^{\lambda-m}d\tau,\label{7.79}\\
\overline{\mathcal{R}_{m}^{\left(  \lambda\right)  }\left(  z\right)  }  &
=\mathcal{R}_{m}^{\left(  \lambda\right)  }\left(  \overline{z}\right)
,\nonumber\\
\left\vert \mathcal{R}_{m}^{\left(  \lambda\right)  }\left(  z\right)
\right\vert  & \leq m\tbinom{\lambda}{m}\left\vert z\right\vert ^{m}\int%
_{0}^{1}\left(  1-\tau\right)  ^{m-1}\left\vert 1+\tau z\right\vert
^{\lambda-m}d\tau\nonumber\\
& \leq\tbinom{\lambda}{m}\left\vert z\right\vert ^{m}\max_{0\leq\tau\leq
1}\left\vert 1+\tau z\right\vert ^{\lambda-m}\nonumber\\
& \leq\left\{
\begin{array}
[c]{ll}%
\tbinom{\lambda}{m}\left\vert z\right\vert ^{m}\left(  1+\left\vert
z\right\vert \right)  ^{\lambda-m}, & \lambda\geq m,\\
\tbinom{\lambda}{m}\left\vert z\right\vert ^{m}\left(  1-\left\vert
z\right\vert \right)  ^{\lambda-m}, & \lambda<m,
\end{array}
\right. \nonumber
\end{align}

which is the second estimate of \ref{a50.9}. Thus%
\begin{align*}
\left\vert 1+te^{i\theta}\right\vert ^{2\lambda} &  =\left(  1+z\right)
^{\lambda}\left(  1+\overline{z}\right)  ^{\lambda}\\
&  =\left(  \sum_{k=0}^{m-1}\tbinom{\lambda}{k}z^{k}+\mathcal{R}_{m}\left(
z\right)  \right)  \left(  \sum_{j=0}^{m-1}\tbinom{\lambda}{j}\overline{z}%
^{j}+\mathcal{R}_{m}^{\left(  \lambda\right)  }\left(  \overline{z}\right)
\right) \\
&  =\sum_{j,k=0}^{m-1}\tbinom{\lambda}{k}\tbinom{\lambda}{j}z^{k}\overline
{z}^{j}+\mathcal{R}_{m}^{\left(  \lambda\right)  }\left(  \overline{z}\right)
\sum_{k=0}^{m-1}\tbinom{\lambda}{k}z^{k}+\mathcal{R}_{m}^{\left(
\lambda\right)  }\left(  z\right)  \sum_{j=0}^{m-1}\tbinom{\lambda}%
{j}\overline{z}^{j}+\\
&  \qquad\qquad+\mathcal{R}_{m}^{\left(  \lambda\right)  }\left(  z\right)
\mathcal{R}_{m}^{\left(  \lambda\right)  }\left(  \overline{z}\right) \\
&  =\sum_{j,k=0}^{m-1}\tbinom{\lambda}{k}\tbinom{\lambda}{j}z^{k}\overline
{z}^{j}+\mathcal{S}_{m}^{\left(  \lambda\right)  }\left(  z,\overline
{z}\right)  ,
\end{align*}

where%
\begin{equation}
\mathcal{S}_{m}^{\left(  \lambda\right)  }\left(  z,\overline{z}\right)
=\mathcal{R}_{m}^{\left(  \lambda\right)  }\left(  \overline{z}\right)
\sum_{k=0}^{m-1}\tbinom{\lambda}{k}z^{k}+\mathcal{R}_{m}^{\left(
\lambda\right)  }\left(  z\right)  \sum_{j=0}^{m-1}\tbinom{\lambda}%
{j}\overline{z}^{j}+\mathcal{R}_{m}^{\left(  \lambda\right)  }\left(
z\right)  \mathcal{R}_{m}^{\left(  \lambda\right)  }\left(  \overline
{z}\right)  .\label{7.76}%
\end{equation}

In particular%
\begin{equation}
\left\vert 1+te^{i\theta}\right\vert ^{2\lambda}=\sum_{j,k=0}^{m-1}%
\tbinom{\lambda}{k}\tbinom{\lambda}{j}t^{k+j}e^{i\left(  k-j\right)  \theta
}+\mathcal{S}_{m}^{\left(  \lambda\right)  }\left(  te^{i\theta},te^{-i\theta
}\right)  .\label{7.77}%
\end{equation}

\begin{align*}
\left\vert \mathcal{S}_{m}^{\left(  \lambda\right)  }\left(  te^{i\theta
},te^{-i\theta}\right)  \right\vert  & \leq\left\vert \mathcal{R}_{m}^{\left(
\lambda\right)  }\left(  \overline{z}\right)  \right\vert \left\vert
\sum_{k=0}^{m-1}\tbinom{\lambda}{k}z^{k}\right\vert +\left\vert \mathcal{R}%
_{m}^{\left(  \lambda\right)  }\left(  z\right)  \right\vert \left\vert
\sum_{j=0}^{m-1}\tbinom{\lambda}{j}\overline{z}^{j}\right\vert +\left\vert
\mathcal{R}_{m}^{\left(  \lambda\right)  }\left(  z\right)  \right\vert
\left\vert \mathcal{R}_{m}^{\left(  \lambda\right)  }\left(  \overline
{z}\right)  \right\vert \\
& =\left\vert \mathcal{R}_{m}^{\left(  \lambda\right)  }\left(  z\right)
\right\vert \left\vert \sum_{k=0}^{m-1}\tbinom{\lambda}{k}z^{k}\right\vert
+\left\vert \mathcal{R}_{m}^{\left(  \lambda\right)  }\left(  z\right)
\right\vert \left\vert \sum_{j=0}^{m-1}\tbinom{\lambda}{j}z^{j}\right\vert
+\left\vert \mathcal{R}_{m}^{\left(  \lambda\right)  }\left(  z\right)
\right\vert ^{2}\\
& =\left\vert \mathcal{R}_{m}^{\left(  \lambda\right)  }\left(  z\right)
\right\vert \left(  2\left\vert \sum_{k=0}^{m-1}\tbinom{\lambda}{k}%
z^{k}\right\vert +\left\vert \mathcal{R}_{m}^{\left(  \lambda\right)  }\left(
z\right)  \right\vert \right)
\end{align*}

Partitioning the square $0\leq j,k\leq m-1$ using the main off-diagonal
$k+j=m-1$ gives%
\begin{align}
&  \left\vert 1+te^{i\theta}\right\vert ^{2\lambda}\nonumber\\
&  =\underset{below/on\text{ }diag.}{\underbrace{\sum_{n=0}^{m-1}\sum
_{k+j=n}\tbinom{\lambda}{k}\tbinom{\lambda}{j}t^{k+j}e^{i\left(  k-j\right)
\theta}}}+\underset{above\text{ }diag.}{\underbrace{\sum_{n=m}^{2m-2}%
\sum_{\substack{k+j=n \\k,j\leq m-1}}\tbinom{\lambda}{k}\tbinom{\lambda}%
{j}t^{k+j}e^{i\left(  k-j\right)  \theta}}}+\mathcal{S}_{m}^{\left(
\lambda\right)  }\nonumber\\
&  =\sum_{n=0}^{m-1}\sum_{k=0}^{n}\tbinom{\lambda}{k}\tbinom{\lambda}%
{n-k}t^{n}e^{i\left(  2k-n\right)  \theta}+\sum_{n=m}^{2m-2}\left(
\ldots\right)  t^{n}+\mathcal{S}_{m}^{\left(  \lambda\right)  }\nonumber\\
&  =\sum_{n=0}^{m-1}\left(  \sum_{k=0}^{n}\tbinom{\lambda}{k}\tbinom{\lambda
}{n-k}e^{i\left(  -n+2k\right)  \theta}\right)  t^{n}+\sum_{n=m}^{2m-2}\left(
\ldots\right)  t^{n}+\mathcal{S}_{m}^{\left(  \lambda\right)  }\nonumber\\
&  =\sum_{\substack{n=0 \\n\text{ }odd}}^{m-1}\left(  \sum_{k=0}^{n}%
\tbinom{\lambda}{k}\tbinom{\lambda}{n-k}e^{i\left(  -n+2k\right)  \theta
}\right)  t^{n}+\sum_{\substack{n=0 \\n\text{ }even}}^{m-1}\left(
\sum_{\substack{k+j=n \\k,j\leq m-1}}\tbinom{\lambda}{k}\tbinom{\lambda}%
{n-k}e^{i\left(  2k-n\right)  \theta}\right)  t^{n}+\sum_{n=m}^{2m-2}\left(
\ldots\right)  t^{n}+\mathcal{\ldots}\nonumber\\
&  =\sum_{\substack{n=0 \\n\text{ }odd}}^{m-1}\left(  \sum_{\substack{l=-n
\\l\text{ }odd}}^{n}\tbinom{\lambda}{k}\tbinom{\lambda}{n-k}e^{il\theta
}\right)  t^{n}+\sum_{\substack{n=0 \\n\text{ }even}}^{m-1}\left(
\sum_{\substack{l=-n \\l\text{ }even}}^{n}\tbinom{\lambda}{k}\tbinom{\lambda
}{n-k}e^{il\theta}\right)  t^{n}+\sum_{n=m}^{2m-2}\left(  \ldots\right)
t^{n}+\mathcal{\ldots}\nonumber\\
&  =\sum_{\substack{n=0 \\n\text{ }odd}}^{m-1}\left(  \sum_{\substack{l=-n
\\l\text{ }odd}}^{n}\tbinom{\lambda}{\frac{n+l}{2}}\tbinom{\lambda}{\frac
{n-l}{2}}e^{il\theta}\right)  t^{n}+\sum_{\substack{n=0 \\n\text{ }%
even}}^{m-1}\left(  \sum_{\substack{l=-n \\l\text{ }even}}^{n}\tbinom{\lambda
}{\frac{n+l}{2}}\tbinom{\lambda}{\frac{n-l}{2}}e^{il\theta}\right)  t^{n}%
+\sum_{n=m}^{2m-2}\left(  \ldots\right)  t^{n}+\mathcal{\ldots}\nonumber\\
&  =\sum_{\substack{n=0 \\n\text{ }odd}}^{m-1}\left(  \sum_{\substack{l=1
\\l\text{ }odd}}^{n}\tbinom{\lambda}{\frac{n+l}{2}}\tbinom{\lambda}{\frac
{n-l}{2}}\left(  e^{il\theta}+e^{-il\theta}\right)  \right)  t^{n}%
+\sum_{\substack{n=0 \\n\text{ }even}}^{m-1}\left(  \tbinom{\lambda}{\frac
{n}{2}}^{2}+\sum_{\substack{l=2 \\l\text{ }even}}^{n}\tbinom{\lambda}%
{\frac{n+l}{2}}\tbinom{\lambda}{\frac{n-l}{2}}\left(  e^{il\theta
}+e^{-il\theta}\right)  \right)  t^{n}+\ldots\nonumber\\
&  =2\sum_{\substack{n=0 \\n\text{ }odd}}^{m-1}\left(  \sum_{\substack{l=1
\\l\text{ }odd}}^{n}\tbinom{\lambda}{\frac{n+l}{2}}\tbinom{\lambda}{\frac
{n-l}{2}}\cos l\theta\right)  t^{n}+\sum_{\substack{n=0 \\n\text{ }%
even}}^{m-1}\left(  \tbinom{\lambda}{\frac{n}{2}}^{2}+2\sum_{\substack{l=2
\\l\text{ }even}}^{n}\tbinom{\lambda}{\frac{n+l}{2}}\tbinom{\lambda}%
{\frac{n-l}{2}}\cos l\theta\right)  t^{n}+\ldots\nonumber\\
&  =\sum_{n=0}^{m-1}P_{n}^{\left(  \lambda\right)  }\left(  \cos\theta\right)
t^{n}+\sum_{n=m}^{2m-2}\left(  \sum_{\substack{k+j=n \\k,j\leq m-1}%
}\tbinom{\lambda}{k}\tbinom{\lambda}{n-k}e^{i\left(  2k-n\right)  \theta
}\right)  t^{n}+\mathcal{S}_{m}^{\left(  \lambda\right)  }\left(  te^{i\theta
},te^{-i\theta}\right)  ,\label{a1.78}%
\end{align}

where the polynomials $P_{n}^{\left(  \lambda\right)  }$ of degree $n$ are
defined by
\begin{equation}
P_{n}^{\left(  \lambda\right)  }\left(  \cos\theta\right)  :=\left\{
\begin{array}
[c]{ll}%
2\sum\limits_{\substack{l=1 \\l\text{ }odd}}^{n}\tbinom{\lambda}{\frac{n+l}%
{2}}\tbinom{\lambda}{\frac{n-l}{2}}\cos l\theta, & n\text{ }odd,\\
\tbinom{\lambda}{\frac{n}{2}}^{2}+2\sum\limits_{\substack{l=2 \\l\text{ }%
even}}^{n}\tbinom{\lambda}{\frac{n+l}{2}}\tbinom{\lambda}{\frac{n-l}{2}}\cos
l\theta, & n\text{ }even.
\end{array}
\right\} \label{7.78}%
\end{equation}

Next we study the second summation in \ref{a1.78}. First write%
\[
\sum_{n=m}^{2m-2}\sum_{\substack{k+j=n \\k,j\leq m-1}}\tbinom{\lambda}%
{k}\tbinom{\lambda}{j}t^{k+j}e^{i\left(  k-j\right)  \theta}=\sum_{n=m}%
^{2m-2}\left(  \sum_{\substack{k-j=n-m-1 \\k\leq m-1}}\tbinom{\lambda}%
{k}\tbinom{\lambda}{j}e^{i\left(  k-j\right)  \theta}\right)  t^{n}.
\]

Now partition the region above the main off-diagonal using the main diagonal:%
\begin{align*}
& \sum_{n=m}^{2m-2}\left(  \sum_{\substack{k+j=n \\k,j\leq m-1}}\tbinom
{\lambda}{k}\tbinom{\lambda}{j}e^{i\left(  k-j\right)  \theta}\right)  t^{n}\\
& =\sum_{\substack{n=m \\n\text{ }even}}^{2m-2}\left(  \sum_{\substack{k+j=n
\\k,j\leq m-1}}\tbinom{\lambda}{k}\tbinom{\lambda}{j}e^{i\left(  k-j\right)
\theta}\right)  t^{n}+\sum_{\substack{n=m \\n\text{ }odd}}^{2m-2}\left(
\sum_{\substack{k+j=n \\k,j\leq m-1}}\tbinom{\lambda}{k}\tbinom{\lambda}%
{j}e^{i\left(  k-j\right)  \theta}\right)  t^{n}\\
& =\sum_{\substack{n=m \\n\text{ }even}}^{2m-2}\left(
\begin{array}
[c]{c}%
\sum\limits_{k=\frac{n}{2}}^{m-1}\sum\limits_{\substack{k-j \\=n-m+1}%
}\tbinom{\lambda}{k}\tbinom{\lambda}{j}e^{i\left(  k-j\right)  \theta}+\\
+\sum\limits_{k=\frac{n}{2}}^{m-1}\sum\limits_{k=j}\tbinom{\lambda}{k}%
\tbinom{\lambda}{j}e^{i\left(  k-j\right)  \theta}+\\
+\sum\limits_{j=\frac{n}{2}}^{m-1}\sum\limits_{\substack{k-j \\=m-n-1}%
}\tbinom{\lambda}{k}\tbinom{\lambda}{j}e^{i\left(  k-j\right)  \theta}%
\end{array}
\right)  t^{n}+\sum_{\substack{n=m \\n\text{ }odd}}^{2m-2}\left(
\begin{array}
[c]{c}%
\sum\limits_{k=\frac{n-1}{2}}^{m-1}\sum\limits_{\substack{k-j \\=n-m+1}%
}\tbinom{\lambda}{k}\tbinom{\lambda}{j}e^{i\left(  k-j\right)  \theta}+\\
+\sum\limits_{k=\frac{n-1}{2}}^{m-1}\sum\limits_{k=j}\tbinom{\lambda}%
{k}\tbinom{\lambda}{j}e^{i\left(  k-j\right)  \theta}+\\
+\sum\limits_{j=\frac{n-1}{2}}^{m-1}\sum\limits_{\substack{k-j \\=m-n-1}%
}\tbinom{\lambda}{k}\tbinom{\lambda}{j}e^{i\left(  k-j\right)  \theta}%
\end{array}
\right)  t^{n}\\
& =\sum_{\substack{n=m \\n\text{ }even}}^{2m-2}\left(
\begin{array}
[c]{c}%
\left(  \sum\limits_{k=\frac{n}{2}}^{m-1}\sum\limits_{k-j=n-m+1}%
\tbinom{\lambda}{k}\tbinom{\lambda}{k-\left(  n-m-1\right)  }\right)
e^{i\left(  n-m-1\right)  \theta}+\sum\limits_{k=\frac{n}{2}}^{m-1}%
\tbinom{\lambda}{k}^{2}+\\
+\left(  \sum\limits_{j=\frac{n}{2}}^{m-1}\sum\limits_{k-j=m-n-1}%
\tbinom{\lambda}{j-\left(  n-m-1\right)  }\tbinom{\lambda}{j}\right)
e^{-i\left(  n-m-1\right)  \theta}%
\end{array}
\right)  t^{n}+\\
& \qquad+\sum_{\substack{n=m \\n\text{ }odd}}^{2m-2}\left(
\begin{array}
[c]{c}%
\left(  \sum\limits_{k=\frac{n-1}{2}}^{m-1}\sum\limits_{k-j=n-m+1}%
\tbinom{\lambda}{k}\tbinom{\lambda}{j}\right)  e^{i\left(  n-m-1\right)
\theta}+\sum\limits_{k=\frac{n-1}{2}}^{m-1}\tbinom{\lambda}{k}^{2}+\\
+\left(  \sum\limits_{j=\frac{n-1}{2}}^{m-1}\sum\limits_{k-j=m-n-1}%
\tbinom{\lambda}{k}\tbinom{\lambda}{j}\right)  e^{-i\left(  n-m-1\right)
\theta}%
\end{array}
\right)  t^{n}\\
& =\sum_{\substack{n=m \\n\text{ }even}}^{2m-2}\sum\limits_{k=\frac{n}{2}%
}^{m-1}\left(  \tbinom{\lambda}{k}^{2}+2\tbinom{\lambda}{k}\tbinom{\lambda
}{k-\left(  n-m+1\right)  }\right)  \cos\left(  n-m+1\right)  \theta\text{
}t^{n}+\\
& \qquad+\sum_{\substack{n=m \\n\text{ }odd}}^{2m-2}\sum\limits_{k=\frac
{n-1}{2}}^{m-1}\left(  \tbinom{\lambda}{k}^{2}+2\tbinom{\lambda}{k}%
\tbinom{\lambda}{k-\left(  n-m+1\right)  }\right)  \cos\left(  n-m-1\right)
\theta\text{ }t^{n}\\
& =\sum_{n=m}^{2m-2}Q_{n}^{\left(  \lambda\right)  }\left(  m,\cos
\theta\right)  t^{n},
\end{align*}

where%
\begin{align}
& Q_{n}^{\left(  \lambda\right)  }\left(  m,\cos\theta\right) \nonumber\\
& :=\left\{
\begin{array}
[c]{ll}%
\sum\limits_{k=\frac{n}{2}}^{m-1}\left(  \tbinom{\lambda}{k}^{2}%
+2\tbinom{\lambda}{k}\tbinom{\lambda}{k-\left(  n-m+1\right)  }\right)
\cos\left(  n-m+1\right)  \theta, & \left\{
\begin{array}
[c]{l}%
n\text{ }even,\\
m\leq n\leq2m-2.
\end{array}
\right. \\
\sum\limits_{k=\frac{n-1}{2}}^{m-1}\left(  \tbinom{\lambda}{k}^{2}%
+2\tbinom{\lambda}{k}\tbinom{\lambda}{k-\left(  n-m+1\right)  }\right)
\cos\left(  n-m+1\right)  \theta, & \left\{
\begin{array}
[c]{l}%
n\text{ }odd,\\
m\leq n\leq2m-3.
\end{array}
\right.
\end{array}
\right. \label{a2.27}%
\end{align}

It remains to prove the upper bound \ref{a4.9}. From \ref{a50.9},%
\[
\left\vert \mathcal{S}_{m}^{\left(  \lambda\right)  }\left(  te^{i\theta
},te^{-i\theta}\right)  \right\vert \leq\left\vert \mathcal{R}_{m}^{\left(
\lambda\right)  }\left(  te^{i\theta}\right)  \right\vert \left(  2\left\vert
\sum_{k=0}^{m-1}\tbinom{\lambda}{k}t^{k}\right\vert +\left\vert \mathcal{R}%
_{m}^{\left(  \lambda\right)  }\left(  te^{i\theta}\right)  \right\vert
\right)  .
\]

But%
\[
\left(  1+z\right)  ^{\lambda}=\sum_{k=0}^{m-1}\tbinom{\lambda}{k}%
z^{k}+\mathcal{R}_{m}^{\left(  \lambda\right)  }\left(  z\right)  ,
\]

so%
\begin{align*}
\left\vert \sum_{k=0}^{m-1}\tbinom{\lambda}{k}t^{k}\right\vert  &
\leq\left\vert \left(  1+te^{i\theta}\right)  ^{\lambda}\right\vert
+\left\vert \mathcal{R}_{m}^{\left(  \lambda\right)  }\left(  te^{i\theta
}\right)  \right\vert \\
& =\left\vert 1+te^{i\theta}\right\vert ^{\lambda}+\left\vert \mathcal{R}%
_{m}^{\left(  \lambda\right)  }\left(  te^{i\theta}\right)  \right\vert ,
\end{align*}

and since $0\leq t<1$, the estimates \ref{a50.9} imply%
\begin{align*}
&  \left\vert \mathcal{S}_{m}^{\left(  \lambda\right)  }\left(  te^{i\theta
},te^{-i\theta}\right)  \right\vert \\
&  \leq\left\vert \mathcal{R}_{m}^{\left(  \lambda\right)  }\left(
te^{i\theta}\right)  \right\vert \left(  2\left\vert 1+te^{i\theta}\right\vert
^{\lambda}+3\left\vert \mathcal{R}_{m}^{\left(  \lambda\right)  }\left(
te^{i\theta}\right)  \right\vert \right) \\
&  =\left\vert \mathcal{R}_{m}^{\left(  \lambda\right)  }\left(  te^{i\theta
}\right)  \right\vert \left(  2\left\vert 1+te^{i\theta}\right\vert
^{\lambda-m}\left\vert 1+te^{i\theta}\right\vert ^{m}+3\left\vert
\mathcal{R}_{m}^{\left(  \lambda\right)  }\left(  te^{i\theta}\right)
\right\vert \right) \\
&  \leq\left\{
\begin{array}
[c]{ll}%
\left\vert \mathcal{R}_{m}^{\left(  \lambda\right)  }\left(  te^{i\theta
}\right)  \right\vert \left(  2\left(  1-t\right)  ^{\lambda-m}\left(
1+t\right)  ^{m}+3\left\vert \mathcal{R}_{m}^{\left(  \lambda\right)  }\left(
te^{i\theta}\right)  \right\vert \right)  , & \lambda\leq m,\\
\left\vert \mathcal{R}_{m}^{\left(  \lambda\right)  }\left(  te^{i\theta
}\right)  \right\vert \left(  2\left(  1+t\right)  ^{\lambda-m}\left(
1+t\right)  ^{m}+3\left\vert \mathcal{R}_{m}^{\left(  \lambda\right)  }\left(
te^{i\theta}\right)  \right\vert \right)  , & \lambda\geq m,
\end{array}
\right. \\
&  \leq\left\{
\begin{array}
[c]{ll}%
\left\vert \mathcal{R}_{m}^{\left(  \lambda\right)  }\left(  te^{i\theta
}\right)  \right\vert \left(  2\left(  1-t\right)  ^{\lambda-m}\left(
1+t\right)  ^{m}+3\left\vert \mathcal{R}_{m}^{\left(  \lambda\right)  }\left(
te^{i\theta}\right)  \right\vert \right)  , & \lambda\leq m,\\
\left\vert \mathcal{R}_{m}^{\left(  \lambda\right)  }\left(  te^{i\theta
}\right)  \right\vert \left(  2\left(  1+t\right)  ^{\lambda}+3\left\vert
\mathcal{R}_{m}^{\left(  \lambda\right)  }\left(  te^{i\theta}\right)
\right\vert \right)  , & \lambda\geq m,
\end{array}
\right. \\
&  =\left\{
\begin{array}
[c]{ll}%
\tbinom{\lambda}{m}t^{m}\left(  1-t\right)  ^{\lambda-m}\left\{  2\left(
1-t\right)  ^{\lambda-m}\left(  1+t\right)  ^{m}+3\tbinom{\lambda}{m}%
t^{m}\left(  1-t\right)  ^{\lambda-m}\right\}  , & \lambda\leq m,\\
\tbinom{\lambda}{m}t^{m}\left(  1+t\right)  ^{\lambda-m}\left\{  2\left(
1+t\right)  ^{\lambda}+3\tbinom{\lambda}{m}t^{m}\left(  1+t\right)
^{\lambda-m}\right\}  , & \lambda\geq m,
\end{array}
\right. \\
&  =\left\{
\begin{array}
[c]{ll}%
\tbinom{\lambda}{m}t^{m}\left(  1-t\right)  ^{2\left(  \lambda-m\right)
}\left\{  2\left(  1+t\right)  ^{m}+3\tbinom{\lambda}{m}t^{m}\right\}  , &
\lambda\leq m,\\
\tbinom{\lambda}{m}t^{m}\left(  1+t\right)  ^{2\lambda-m}\left\{
2+3\tbinom{\lambda}{m}t^{m}\left(  1+t\right)  ^{-m}\right\}  , & \lambda\geq
m,
\end{array}
\right. \\
&  \leq\left\{
\begin{array}
[c]{ll}%
\tbinom{\lambda}{m}t^{m}\left(  1-t\right)  ^{2\left(  \lambda-m\right)
}\left(  2\left(  1+t\right)  ^{m}+3\tbinom{\lambda}{m}t^{m}\right)  , &
\lambda\leq m,\\
\tbinom{\lambda}{m}t^{m}\left(  1+t\right)  ^{2\lambda-m}\left(
2+3\tbinom{\lambda}{m}t^{m}\right)  , & \lambda\geq m,
\end{array}
\right.
\end{align*}

and we have derived the upper bound%
\begin{align*}
\left\vert \mathcal{S}_{m}^{\left(  \lambda\right)  }\left(  te^{i\theta
},te^{-i\theta}\right)  \right\vert  & \leq\left\{
\begin{array}
[c]{ll}%
\tbinom{\lambda}{m}t^{m}\left(  1-t\right)  ^{2\left(  \lambda-m\right)
}\left(  2\left(  1+t\right)  ^{m}+3\tbinom{\lambda}{m}t^{m}\right)  , &
\lambda\leq m,\\
\tbinom{\lambda}{m}t^{m}\left(  1+t\right)  ^{2\lambda-m}\left(
2+3\tbinom{\lambda}{m}t^{m}\right)  , & \lambda\geq m,
\end{array}
\right. \\
when\text{ }0  & \leq t<1.
\end{align*}

Estimate of $\mathcal{S}_{m}^{\left(  \lambda\right)  }\left(  te^{i\theta
},te^{-i\theta}\right)  $: from Lemma \ref{Lem_Legendre},%
\[
\left\vert \mathcal{S}_{m}^{\left(  \lambda\right)  }\left(  te^{i\theta
},te^{-i\theta}\right)  \right\vert \leq\left\vert \mathcal{R}_{m}^{\left(
\lambda\right)  }\left(  te^{i\theta}\right)  \right\vert \left(  2\left\vert
\sum_{k=0}^{m-1}\tbinom{\lambda}{k}t^{k}\right\vert +\left\vert \mathcal{R}%
_{m}^{\left(  \lambda\right)  }\left(  te^{i\theta}\right)  \right\vert
\right)  .
\]

But%
\[
\left(  1+z\right)  ^{\lambda}=\sum_{k=0}^{m-1}\tbinom{\lambda}{k}%
z^{k}+\mathcal{R}_{m}^{\left(  \lambda\right)  }\left(  z\right)  ,
\]

so%
\begin{align*}
\left\vert \sum_{k=0}^{m-1}\tbinom{\lambda}{k}t^{k}\right\vert  &
\leq\left\vert \left(  1+te^{i\theta}\right)  ^{\lambda}\right\vert
+\left\vert \mathcal{R}_{m}^{\left(  \lambda\right)  }\left(  te^{i\theta
}\right)  \right\vert \\
& =\left\vert 1+te^{i\theta}\right\vert ^{\lambda}+\left\vert \mathcal{R}%
_{m}^{\left(  \lambda\right)  }\left(  te^{i\theta}\right)  \right\vert ,
\end{align*}

and, since $0\leq t<1$, the estimates \ref{a50.9} imply%
\begin{align*}
&  \left\vert \mathcal{S}_{m}^{\left(  \lambda\right)  }\left(  te^{i\theta
},te^{-i\theta}\right)  \right\vert \\
&  \leq\left\vert \mathcal{R}_{m}^{\left(  \lambda\right)  }\left(
te^{i\theta}\right)  \right\vert \left(  2\left\vert 1+te^{i\theta}\right\vert
^{\lambda}+3\left\vert \mathcal{R}_{m}^{\left(  \lambda\right)  }\left(
te^{i\theta}\right)  \right\vert \right) \\
&  =\left\vert \mathcal{R}_{m}^{\left(  \lambda\right)  }\left(  te^{i\theta
}\right)  \right\vert \left(  2\left\vert 1+te^{i\theta}\right\vert
^{\lambda-m}\left\vert 1+te^{i\theta}\right\vert ^{m}+3\left\vert
\mathcal{R}_{m}^{\left(  \lambda\right)  }\left(  te^{i\theta}\right)
\right\vert \right) \\
&  \leq\left\{
\begin{array}
[c]{ll}%
\left\vert \mathcal{R}_{m}^{\left(  \lambda\right)  }\left(  te^{i\theta
}\right)  \right\vert \left(  2\left(  1-t\right)  ^{\lambda-m}\left(
1+t\right)  ^{m}+3\left\vert \mathcal{R}_{m}^{\left(  \lambda\right)  }\left(
te^{i\theta}\right)  \right\vert \right)  , & \lambda\leq m,\\
\left\vert \mathcal{R}_{m}^{\left(  \lambda\right)  }\left(  te^{i\theta
}\right)  \right\vert \left(  2\left(  1+t\right)  ^{\lambda-m}\left(
1+t\right)  ^{m}+3\left\vert \mathcal{R}_{m}^{\left(  \lambda\right)  }\left(
te^{i\theta}\right)  \right\vert \right)  , & \lambda\geq m,
\end{array}
\right. \\
&  \leq\left\{
\begin{array}
[c]{ll}%
\left\vert \mathcal{R}_{m}^{\left(  \lambda\right)  }\left(  te^{i\theta
}\right)  \right\vert \left(  2\left(  1-t\right)  ^{\lambda-m}\left(
1+t\right)  ^{m}+3\left\vert \mathcal{R}_{m}^{\left(  \lambda\right)  }\left(
te^{i\theta}\right)  \right\vert \right)  , & \lambda\leq m,\\
\left\vert \mathcal{R}_{m}^{\left(  \lambda\right)  }\left(  te^{i\theta
}\right)  \right\vert \left(  2\left(  1+t\right)  ^{\lambda}+3\left\vert
\mathcal{R}_{m}^{\left(  \lambda\right)  }\left(  te^{i\theta}\right)
\right\vert \right)  , & \lambda\geq m,
\end{array}
\right. \\
&  =\left\{
\begin{array}
[c]{ll}%
\tbinom{\lambda}{m}t^{m}\left(  1-t\right)  ^{\lambda-m}\left\{  2\left(
1-t\right)  ^{\lambda-m}\left(  1+t\right)  ^{m}+3\tbinom{\lambda}{m}%
t^{m}\left(  1-t\right)  ^{\lambda-m}\right\}  , & \lambda\leq m,\\
\tbinom{\lambda}{m}t^{m}\left(  1+t\right)  ^{\lambda-m}\left\{  2\left(
1+t\right)  ^{\lambda}+3\tbinom{\lambda}{m}t^{m}\left(  1+t\right)
^{\lambda-m}\right\}  , & \lambda\geq m,
\end{array}
\right. \\
&  =\left\{
\begin{array}
[c]{ll}%
\tbinom{\lambda}{m}t^{m}\left(  1-t\right)  ^{2\left(  \lambda-m\right)
}\left\{  2\left(  1+t\right)  ^{m}+3\tbinom{\lambda}{m}t^{m}\right\}  , &
\lambda\leq m,\\
\tbinom{\lambda}{m}t^{m}\left(  1+t\right)  ^{2\lambda-m}\left\{
2+3\tbinom{\lambda}{m}t^{m}\left(  1+t\right)  ^{-m}\right\}  , & \lambda\geq
m,
\end{array}
\right. \\
&  \leq\left\{
\begin{array}
[c]{ll}%
\tbinom{\lambda}{m}t^{m}\left(  1-t\right)  ^{2\left(  \lambda-m\right)
}\left(  2\left(  1+t\right)  ^{m}+3\tbinom{\lambda}{m}t^{m}\right)  , &
\lambda\leq m,\\
\tbinom{\lambda}{m}t^{m}\left(  1+t\right)  ^{2\lambda-m}\left(
2+3\tbinom{\lambda}{m}t^{m}\right)  , & \lambda\geq m,
\end{array}
\right.
\end{align*}

which is \ref{a4.9}.
\end{proof}

\chapter{Quotient spaces and reproducing kernels\label{Ch_quot_sp_reprod_kern}%
}

We present the basic theory of Banach space quotient spaces and then the
general theory of reproducing kernels. Finally we characterize restriction
spaces in terms of a quotient space.

\section{General quotient space theory}

This section is from Chapter 3 "Quotient Spaces" of Schwartz \cite{Schwartz09}
and also Chapter 5 of Packel \cite{Packel74}. We present without proofs the
basic theory of Banach space quotient spaces including annihilators, linear
functionals and adjoints.

\begin{definition}
\label{Def_QuotientSpace}\textbf{Quotient space} Suppose $\left(  X,\left\vert
\cdot\right\vert \right)  $ is a \textbf{seminorm space} and $M\subset X$ is a
linear subspace. Denote the \textbf{quotient space} by $X/M$, where%
\[
X/M=\left\{  \left[  x\right]  =x+M:x\in X\right\}  ,
\]

where $\left[  x\right]  $ is termed a \textbf{coset}. $\left[  x\right]  $
can be associated with an equivalence relation $\sim$ defined on $X$ by%
\[
x\sim y\Leftrightarrow x-y\in M,
\]

so that%
\[
\left[  x\right]  =\left\{  y:y\sim x\right\}  .
\]

Endow $\left[  x\right]  $ with the \textbf{quotient seminorm}%
\[
\left\Vert \left[  x\right]  \right\Vert ^{\prime}:=\inf_{y\sim x}\left\vert
y\right\vert =\inf\left\{  \left\vert x+m\right\vert :m\in M\right\}
=\operatorname*{dist}\left(  \left\{  x\right\}  ,M\right)  .
\]

\end{definition}

\begin{theorem}
\label{Thm_QuotientProperties}\ 

\begin{enumerate}
\item $\left\Vert \cdot\right\Vert ^{\prime}$ defines a seminorm on $X/M$.

\item $\left\Vert \cdot\right\Vert ^{\prime}$ is a norm iff $M$ is closed.

\item If $M$ is closed the \textbf{quotient map} $x\rightarrow\left[
x\right]  $ is a reducing norm w.r.t. $\left\Vert \cdot\right\Vert ^{\prime}%
$as well as being continuous, onto and open.

\item $X/M$ is complete iff $X$ is complete.

\item If $X$ is complete and $M$ is closed, then $X/M$ is a Banach space.
\end{enumerate}
\end{theorem}

Next we consider the case where $M$ is the kernel of the seminorm defined on
$X$.

\begin{definition}
\label{Def_KerSeminorm}\textbf{Kernel of a seminorm} is denoted by $k\left(
X\right)  =\ker\left\vert \cdot\right\vert =\left\{  x:\left\vert x\right\vert
=0\right\}  $.
\end{definition}

\begin{theorem}
\label{Thm_KerSeminorm}\ 

\begin{enumerate}
\item $k\left(  X\right)  $ is a closed subspace of $X$.

\item The quotient map $x\rightarrow\left[  x\right]  $ is an isometry.

\item $X/k\left(  X\right)  $ is a normed space. It is also a Banach space iff
$X$ is complete.
\end{enumerate}
\end{theorem}

\begin{definition}
\label{Def_Annihilator}\textbf{Annihilator} If $K\subset X$, the annihilator
of $K$ (in $X^{\prime}$) is defined to be%
\[
K^{0}=\left\{  x^{\prime}\in X^{\prime}:x^{\prime}\left(  k\right)  =0\text{
}\forall\text{ }k\in K\right\}  .
\]

\end{definition}

\begin{theorem}
\label{Thm_FnalQuotSp_AnnhilDenom}Suppose $M$ is a closed subspace of $X$.
Then $\left(  X/M\right)  ^{\prime}$ and $M^{0}$ are isometrically isomorphic
via the map $V:\left(  X/M\right)  ^{\prime}\rightarrow M^{0}$ defined by%
\[
\left(  Vz^{\prime}\right)  \left(  x\right)  =z^{\prime}\left(  \left[
x\right]  \right)  ,\text{\quad}z^{\prime}\in\left(  X/M\right)  ^{\prime
},\text{ }x\in X.
\]

\end{theorem}

\begin{proof}
??
\end{proof}

\begin{theorem}
\label{Thm_exercises}\textbf{Results from the exercises in Schwartz
\cite{Schwartz09}:}

\begin{enumerate}
\item Suppose $X$ and $Y$ are seminormed spaces and $T\in B\left(  X,Y\right)
$.

Then the induced map $\widehat{T}:X/\ker T\rightarrow Y$ is linear, continuous
and satisfies $\left\Vert T\right\Vert =\left\Vert \widehat{T}\right\Vert $.

\item If $K\subset X$ then $K^{0}$ is a closed subspace of $X^{\prime}$.

\item Suppose $M$ is a linear subspace of $X$.

Then $X^{\prime}/M^{0}$ and $M^{\prime}$ are isometrically isomorphic under
the map $U:X^{\prime}/M^{0}\rightarrow M^{\prime}$ defined by $U\left[
x^{\prime}\right]  =x_{M}^{\prime}$, where $x_{M}^{\prime}$ is the restriction
of $x^{\prime}$ to $M$.

\item $B\subset X/M$ is bounded iff there exists a bounded subset $A\subset X
$ such that $B\subset\left[  A\right]  =A+M$.
\end{enumerate}
\end{theorem}

Standard basic result:

\begin{theorem}
\label{Thm_adjoint_op}\textbf{Adjoint operator} Suppose $X,Y$ are normed
vector spaces and $S:X\rightarrow Y$ is a homeomorphism. Then the adjoint
operator $S^{T}:Y^{\prime}\rightarrow X^{\prime}$ defined by
\[
\left(  S^{T}y^{\prime}\right)  x=y^{\prime}\left(  Sx\right)  ,\text{\quad
}y^{\prime}\in Y^{\prime},\text{ }x\in X,
\]

is also a homeomorphism.
\end{theorem}

\begin{proof}
??
\end{proof}

\section{Reproducing kernels}

We present the general theory of reproducing kernel Hilbert spaces from, for
example, Part I of Aronszajn \cite{Aronszajn50} and Chapter 10 of Wendland
\cite{Wendland05}.

\begin{definition}
\label{Def_reprod_ker}\textbf{Reproducing kernel} (I.1 Aronszajn
\cite{Aronszajn50}) Suppose $\mathcal{H}$ is a class of functions defined on a
set $E\subseteq\mathbb{R}^{d}$ which forms a complex Hilbert space $\left(
\mathcal{H},\left\Vert \cdot\right\Vert ,\left(  \cdot,\cdot\right)  \right)
$. The function $\Phi:E\times E\rightarrow\mathbb{C}$ is called a reproducing
kernel (r.k.) of $\mathcal{H}$ if:

\begin{enumerate}
\item For all $y\in E$, $\Phi\left(  \cdot,y\right)  \in\mathcal{H}$;

\item \textbf{The reproducing property}:%
\[
f\left(  y\right)  =\left(  f,\Phi\left(  \cdot,y\right)  \right)
,\text{\quad}y\in E,\text{ }f\in\mathcal{H}.
\]

\end{enumerate}
\end{definition}

\begin{remark}
\label{Rem_proper_reprod_ker}\textbf{Properties of reproducing kernels} (I.2
of Aronszajn \cite{Aronszajn50})

\begin{enumerate}
\item \textbf{Uniqueness} If a reproducing kernel exists it is unique.

\item \textbf{Existence} A r.k. exists iff for each $y\in E$ there exists a
constant $C_{y}\geq0$ such that
\[
\left\vert f\left(  y\right)  \right\vert \leq C_{y}\left\Vert f\right\Vert
\text{ }\forall f\in\mathcal{H}.
\]

\item $\Phi$ is a \textbf{positive matrix} in the sense that for any $n\geq1$,%
\[
\left\Vert \sum\limits_{k=1}^{n}\Phi\left(  \cdot,y^{\left(  k\right)
}\right)  \xi_{k}\right\Vert ^{2}=\sum\limits_{i,j=1}^{n}\Phi\left(
y^{\left(  i\right)  },y^{\left(  j\right)  }\right)  \overline{\xi_{i}}%
\xi_{j}\geq0,\text{\quad}\xi\in\mathbb{C}^{n},\text{ }y^{\left(  i\right)
}\in E.
\]

Also%
\[
\Phi\left(  x,x\right)  \geq0,\text{\quad}\Phi\left(  x,y\right)
=\overline{\Phi\left(  y,x\right)  },\text{\quad}\left\vert \Phi\left(
x,y\right)  \right\vert ^{2}\leq\Phi\left(  x,x\right)  \Phi\left(
y,y\right)  .
\]

\item \textbf{Inverse of part 3 }To every positive matrix $\Phi$ there
corresponds one and only one class of functions with a uniquely determined
quadratic form in it, forming a Hilbert space and admitting $\Phi\left(
x,y\right)  $ as a reproducing kernel.

This class of functions is generated by all functions of the form
$\sum\limits_{k=1}^{n}\Phi\left(  \cdot,y^{\left(  k\right)  }\right)  \xi
_{k}$ with quadratic norm defined by
\[
\left\Vert \sum\limits_{k=1}^{n}\Phi\left(  \cdot,y^{\left(  k\right)
}\right)  \xi_{k}\right\Vert ^{2}=\sum\limits_{i,j=1}^{n}\Phi\left(
y^{\left(  i\right)  },y^{\left(  j\right)  }\right)  \overline{\xi_{i}}%
\xi_{j}.
\]

This norm has an obvious inner product associated with it and can be extended
to a (closed) Hilbert space using Cauchy sequences.

See the \textbf{native spaces} of Chapter 10 of Wendland \cite{Wendland05}.

\item \textbf{Convergence} If $f,f_{n}\in\mathcal{H}$,%
\[
\left\vert f\left(  x\right)  -f_{n}\left(  x\right)  \right\vert
\leq\left\Vert f-f_{n}\right\Vert \sqrt{\Phi\left(  x,x\right)  },\text{\quad
}x\in E.
\]

?? Add more?

\item Suppose $\mathcal{H}_{1}$ has r.k. $\Phi_{1}$ and is a subspace of the
Hilbert space $\mathcal{H}_{2}$ i.e. $\mathcal{H}_{1}\subset\mathcal{H}_{2}$
and $f\in\mathcal{H}_{1}$ implies $\left\Vert f\right\Vert _{1}=\left\Vert
f\right\Vert _{2}$. Then $\mathcal{H}_{1}$ is a closed subspace of
$\mathcal{H}_{2}$ and (Ch. 5 Packel \cite{Packel74}) we have the
\textbf{unique direct sum decomposition} $\mathcal{H}_{2}=\mathcal{H}%
_{1}\oplus\mathcal{H}_{1}^{\bot}$ where \textbf{orthogonal complement}
$\mathcal{H}_{1}^{\bot}$ is closed subspace of $\mathcal{H}_{2}$. In fact,
given $h\in\mathcal{H}_{2}$ there exist unique $f\in\mathcal{H}_{1}$ and
$g\in\mathcal{H}_{1}^{\bot}$ such that $h=f+g$ and $\left(  f,g\right)  _{2}=0
$. The mapping $h\rightarrow f$ defines the \textbf{orthogonal projection} of
$\mathcal{H}_{2}$ onto $\mathcal{H}_{1}$. Since $\Phi\left(  \cdot,y\right)
\in\mathcal{H}_{1}\subset\mathcal{H}_{2}$ we have for $y\in E$,%
\begin{align*}
\left(  h,\Phi\left(  \cdot,y\right)  \right)  _{2}  & =\left(  f+g,\Phi
\left(  \cdot,y\right)  \right)  _{2}=\left(  f,\Phi\left(  \cdot,y\right)
\right)  _{2}+\left(  g,\Phi\left(  \cdot,y\right)  \right)  _{2}=\\
& =\left(  f,\Phi\left(  \cdot,y\right)  \right)  _{2}=\left(  f,\Phi\left(
\cdot,y\right)  \right)  _{1}=f\left(  y\right)  .
\end{align*}

\item If $\mathcal{H}$ has a r.k. then part 2 implies that the same is true
for all closed subspaces of $\mathcal{H}$.

Suppose $\mathcal{H}=\mathcal{H}_{1}\oplus\mathcal{H}_{2}$ where
$\mathcal{H}_{i}$ has r.k. $\Phi_{i}$ and $\mathcal{H}$ has r.k. $\Phi$.
Suppose $f\in\mathcal{H}$.

Then $f=f_{1}+f_{2}$ for unique $f_{i}\in\mathcal{H}_{1}$. Since
$\mathcal{H}_{1}=\mathcal{H}_{2}^{\bot}$ and $\mathcal{H}_{2}=\mathcal{H}%
_{1}^{\bot}$ we have for all $y\in E$,%
\begin{align*}
f_{1}\left(  y\right)   & =\left(  f_{1},\Phi_{1}\left(  \cdot,y\right)
\right)  =\left(  f_{1}+f_{2},\Phi_{1}\left(  \cdot,y\right)  \right)
=\left(  f,\Phi_{1}\left(  \cdot,y\right)  \right)  .\\
f_{2}\left(  y\right)   & =\left(  f_{2},\Phi_{2}\left(  \cdot,y\right)
\right)  =\left(  f_{1}+f_{2},\Phi_{2}\left(  \cdot,y\right)  \right)
=\left(  f,\Phi_{2}\left(  \cdot,y\right)  \right)  .\\
f\left(  y\right)   & =f_{1}\left(  y\right)  +f_{2}\left(  y\right)  =\left(
f,\Phi_{1}\left(  \cdot,y\right)  \right)  +\left(  f,\Phi_{2}\left(
\cdot,y\right)  \right)  =\left(  f,\Phi_{1}\left(  \cdot,y\right)  +\Phi
_{2}\left(  \cdot,y\right)  \right)  .
\end{align*}

But $f\left(  y\right)  =\left(  f,\Phi\left(  \cdot,y\right)  \right)  $ so
$\left(  f,\Phi\left(  \cdot,y\right)  \right)  =\left(  f,\Phi_{1}\left(
\cdot,y\right)  +\Phi_{2}\left(  \cdot,y\right)  \right)  $ for all $y\in E$
and $f\in\mathcal{H}$, and hence $\Phi=\Phi_{1}+\Phi_{2}$ on $E\times E$.

\item \textbf{Orthonormal systems} ?? Add more?
\end{enumerate}
\end{remark}

??

\begin{theorem}
\textbf{Pointwise continuity} Suppose $E$ is open. Then:

\begin{enumerate}
\item For all $f\in\mathcal{H}$,
\[
\left\vert f\left(  x\right)  -f\left(  y\right)  \right\vert \leq\left\Vert
f\right\Vert \left\{  \Phi\left(  x,x\right)  -\Phi\left(  x,y\right)
-\Phi\left(  y,x\right)  +\Phi\left(  y,y\right)  \right\}  ^{1/2},\quad
x,y\in E.
\]

\item If $\lim\limits_{y\rightarrow x,y\in E}\Phi\left(  x,y\right)
=\Phi\left(  x,x\right)  $ for each $x\in E$ and $\Phi\left(  x,x\right)  \in
C^{\left(  0\right)  }\left(  E\right)  $ then $f\in C^{\left(  0\right)
}\left(  E\right)  $.
\end{enumerate}
\end{theorem}

\begin{proof}
\textbf{Part 1}%
\begin{align*}
&  \left\vert f\left(  x\right)  -f\left(  y\right)  \right\vert ^{2}\\
&  =\left\vert \left(  f,\Phi\left(  \cdot,x\right)  \right)  -\left(
f,\Phi\left(  \cdot,y\right)  \right)  \right\vert ^{2}\\
&  =\left\vert \left(  f,\Phi\left(  \cdot,x\right)  -\Phi\left(
\cdot,y\right)  \right)  \right\vert ^{2}\\
&  \leq\left\Vert f\right\Vert ^{2}\left\Vert \Phi\left(  \cdot,x\right)
-\Phi\left(  \cdot,y\right)  \right\Vert ^{2}\\
&  =\left\Vert f\right\Vert ^{2}\left(  \Phi\left(  \cdot,x\right)
-\Phi\left(  \cdot,y\right)  ,\Phi\left(  \cdot,x\right)  -\Phi\left(
\cdot,y\right)  \right) \\
&  =\left\Vert f\right\Vert ^{2}\left\{  \left(  \Phi\left(  \cdot,x\right)
,\Phi\left(  \cdot,x\right)  \right)  -\left(  \Phi\left(  \cdot,x\right)
,\Phi\left(  \cdot,y\right)  \right)  -\left(  \Phi\left(  \cdot,y\right)
,\Phi\left(  \cdot,x\right)  \right)  +\left(  \Phi\left(  \cdot,y\right)
,\Phi\left(  \cdot,y\right)  \right)  \right\} \\
&  =\left\Vert f\right\Vert ^{2}\left\{  \Phi\left(  x,x\right)  -\Phi\left(
y,x\right)  -\Phi\left(  x,y\right)  +\Phi\left(  y,y\right)  \right\}  .
\end{align*}
\medskip

\textbf{Part 2} Write%
\begin{align*}
\Phi\left(  x,x\right)   &  -\Phi\left(  x,y\right)  -\Phi\left(  y,x\right)
+\Phi\left(  y,y\right) \\
&  =\Phi\left(  x,x\right)  -\Phi\left(  x,y\right)  +\Phi\left(  y,y\right)
-\overline{\Phi\left(  x,y\right)  }\\
&  =\Phi\left(  x,x\right)  -\Phi\left(  x,y\right)  +\left(  \Phi\left(
y,y\right)  -\Phi\left(  x,x\right)  \right)  +\overline{\Phi\left(
x,x\right)  -\Phi\left(  x,y\right)  }.
\end{align*}

\end{proof}

From Section I.5 of Aronszajn \cite{Aronszajn50}:

\begin{theorem}
\label{Thm_restrict_rkhs}\textbf{Restrictions of r.k. Hilbert spaces} Suppose
$\left(  \mathcal{H},\left(  \cdot,\cdot\right)  ,\left\Vert \cdot\right\Vert
,\Phi\right)  $ is a r.k. Hilbert space of functions on a set $E\subseteq
\mathbb{R}^{d}$. Define the restriction subspace $\mathcal{H}\left(
E_{1}\right)  :=r_{E_{1}}\mathcal{H}$ where $E_{1}\subset E$.

Then $\left(  \mathcal{H}\left(  E_{1}\right)  ,\left(  \cdot,\cdot\right)
_{1},\left\Vert \cdot\right\Vert _{1},\Phi_{1}\right)  $ is also a r.k.
Hilbert space where:%
\begin{equation}
\left\Vert f_{1}\right\Vert _{1}:=\operatorname*{dist}\left(  \left\{
0\right\}  ,r_{E_{1}}^{-1}f_{1}\right)  ,\quad f_{1}\in\mathcal{H}\left(
E_{1}\right)  ,\label{Q013}%
\end{equation}

is a norm satisfying the parallelogram law and the reproducing kernel is
$\Phi_{1}=r_{E_{1}\times E_{1}}\Phi$. Finally $r_{E_{1}}:\mathcal{H}%
\rightarrow\mathcal{H}\left(  E_{1}\right)  $ is continuous.
\end{theorem}

\begin{proof}
Define the subspace%
\begin{equation}
\mathcal{H}_{E_{1}}:=\left\{  f\in\mathcal{H}:f=0\text{ }on\text{ }%
E_{1}\right\}  .\label{Q014}%
\end{equation}

Suppose $\left\{  f_{k}\right\}  $ is a Cauchy sequence in $\mathcal{H}%
_{E_{1}} $. Then $\exists!$ $f\in\mathcal{H}$ such that $f_{k}\rightarrow f $
in $\mathcal{H}$. But $\left(  f,\Phi\left(  \cdot,x\right)  \right)
=f\left(  x\right)  $ $\forall x\in E$ and $f\in\mathcal{H}$. So if $x\in
E_{1}$,
\[
\left\vert f\left(  x\right)  \right\vert =\left\vert f\left(  x\right)
-f_{k}\left(  x\right)  \right\vert =\left\vert \left(  f-f_{k},\Phi\left(
\cdot,x\right)  \right)  \right\vert \leq\sqrt{\Phi\left(  x,x\right)
}\left\Vert f-f_{k}\right\Vert \rightarrow0,
\]

as $k\rightarrow\infty$. Thus $f\left(  x\right)  =0$ and hence $f\in
\mathcal{H}_{E_{1}}$ which confirms that $\mathcal{H}_{E_{1}}$ is closed and
thus a Hilbert space. Hence we have the unique\textbf{\ direct sum
decomposition}%
\begin{equation}
\mathcal{H}=\mathcal{H}_{E_{1}}\oplus\mathcal{H}_{E_{1}}^{\perp},\label{Q010}%
\end{equation}

where the \textbf{orthogonal complement} $\mathcal{H}_{E_{1}}^{\perp}$ is
unique and is also a closed linear subspace of $\mathcal{H}$.

From part 7 of Remark \ref{Rem_proper_reprod_ker} $\mathcal{H}_{E_{1}}$ and
$\mathcal{H}_{E_{1}}^{\perp}$ have unique reproducing kernels, say $\Phi_{0} $
and $\Phi_{\perp}$, such that%
\[
\Phi=\Phi_{0}+\Phi_{\perp}\text{ }on\text{ }E\times E.
\]

For every fixed $y$, $\Phi_{0}\left(  \cdot,y\right)  \in\mathcal{H}_{E_{1}}
$. Thus $\Phi_{0}\left(  x,y\right)  =0$ when $x\in E_{1}$ and so%
\begin{equation}
\Phi\left(  x,y\right)  =\Phi_{\perp}\left(  x,y\right)  \text{ }on\text{
}E_{1}\times E.\label{Q012}%
\end{equation}

Next consider the definition \ref{Q013} of $\left\Vert f_{1}\right\Vert _{1}$.
Since $r_{E_{1}}^{-1}f_{1}$ is a linear subspace of $\mathcal{H}$ part 1 of
Theorem \ref{Thm_QuotientProperties} implies that $\left\Vert \cdot\right\Vert
_{1}$ is at least a seminorm. But $\left\Vert f_{1}\right\Vert _{1}=0$ implies
$0\in r_{E_{1}}^{-1}f_{1}$ which in turn implies $f_{1}=0$ and so $\left\Vert
\cdot\right\Vert _{1}$ is a norm.

Since $r_{E_{1}}:\mathcal{H}\rightarrow\mathcal{H}\left(  E_{1}\right)  $ is
now continuous, $\mathcal{H}\left(  E_{1}\right)  $ is a Banach space.
Further, $r_{E_{1}}^{-1}f_{1}$ is closed and there exists a unique
$f\in\mathcal{H}$ such that%
\[
\left\Vert f_{1}\right\Vert _{1}=\left\Vert f\right\Vert .
\]

Indeed, if $f\in\mathcal{H}$ then from \ref{Q010} there is a unique
decomposition $f=f_{0}+f_{\perp}$ where $f_{0}\in\mathcal{H}_{E_{1}}$ and
$f_{\perp}\in\mathcal{H}_{E_{1}}^{\perp}$ and further $f\rightarrow f_{\perp}$
defines the \textbf{orthogonal projection} onto $\mathcal{H}_{E_{1}}^{\perp}$.
Also $r_{E_{1}}f=r_{E_{1}}f_{\perp}=f_{1}$ so $f_{\perp}$ is an extension of
$f_{1}$ and hence $\left\Vert f_{1}\right\Vert _{1}\leq\left\Vert f_{\perp
}\right\Vert $. But $\left\Vert f\right\Vert ^{2}=\left\Vert f_{0}+f_{\perp
}\right\Vert ^{2}=\left\Vert f_{0}\right\Vert ^{2}+\left\Vert f_{\perp
}\right\Vert ^{2}$, so $f_{0}=0$ and $f=f_{\perp}$. Thus
\[
\left\Vert f_{1}\right\Vert _{1}=\left\Vert f_{\perp}\right\Vert ,
\]

which establishes an isometric isomorphism between $\mathcal{H}\left(
E_{1}\right)  $ and $\mathcal{H}_{E_{1}}^{\perp}$. This result allows us to
easily prove that the \textbf{parallelogram law} is satisfied:%
\begin{align*}
\left\Vert f_{1}-g_{1}\right\Vert _{1}^{2}+\left\Vert f_{1}+g_{1}\right\Vert
_{1}^{2}=\left\Vert f_{\perp}-g_{\perp}\right\Vert ^{2}+\left\Vert f_{\perp
}+g_{\perp}\right\Vert ^{2}  & =2\left\Vert f_{\perp}\right\Vert
^{2}+2\left\Vert g_{\perp}\right\Vert ^{2}\\
& =2\left\Vert f_{1}\right\Vert _{1}^{2}+2\left\Vert g_{1}\right\Vert _{1}%
^{2},
\end{align*}

and so the norm $\left\Vert \cdot\right\Vert _{1}$ generates the inner product%
\[
\left(  f_{1},g_{1}\right)  _{1}=\left(  f_{\perp},g_{\perp}\right)
\]

and $\left(  \mathcal{H}\left(  E_{1}\right)  ,\left(  \cdot,\cdot\right)
_{1},\left\Vert \cdot\right\Vert _{1}\right)  $ is a Hilbert space.

For $y\in E_{1}$, $f=f_{E_{1}}+f_{\perp}$ implies
\begin{equation}
f_{1}\left(  y\right)  =f\left(  y\right)  =f_{E_{1}}\left(  y\right)
+f_{\perp}\left(  y\right)  =f_{\perp}\left(  y\right)  =\left(  f_{\perp
},\Phi_{\perp}\left(  \cdot,y\right)  \right)  ,\label{Q011}%
\end{equation}

so that for each fixed $y\in E_{1}$ there exists a constant $c_{y}^{\perp}>0 $
such that%
\[
\left\vert f_{1}\left(  y\right)  \right\vert =\left\vert f_{\perp}\left(
y\right)  \right\vert \leq c_{y}^{\perp}\left\Vert f_{\perp}\right\Vert
=c_{y}^{\prime}\left\Vert f_{1}\right\Vert _{1},\quad\forall f_{1}%
\in\mathcal{H}\left(  E_{1}\right)  .
\]

part 2 of Remark \ref{Rem_proper_reprod_ker} now implies that $\mathcal{H}%
\left(  E_{1}\right)  $ has a r.k., say $\Phi_{1}\left(  x,y\right)  $, on
$E_{1}\times E_{1}$. Thus from \ref{Q011},%
\[
f_{1}\left(  y\right)  =\left(  f_{\perp},\Phi_{\perp}\left(  \cdot,y\right)
\right)  =\left(  f_{1},\Phi_{1}\left(  \cdot,y\right)  \right)  _{1},\quad
y\in\mathcal{H}\left(  E_{1}\right)  ,
\]

and then from \ref{Q011} and \ref{Q012},%
\[
\left(  f_{1},\Phi_{1}\left(  \cdot,y\right)  \right)  _{1}=\left(  f_{\perp
},\Phi_{\perp}\left(  \cdot,y\right)  \right)  =\left(  f_{1},\Phi\left(
\cdot,y\right)  \right)  ,\quad y\in E_{1},f_{1}\in\mathcal{H}\left(
E_{1}\right)  .
\]

so that
\[
\Phi_{1}=\Phi\text{ }on\text{ }E_{1}\times E_{1}.
\]

\end{proof}

\subsection{Restrictions of r.k. Hilbert spaces as quotient spaces}

Theorem \ref{Thm_restrict_rkhs} is the relevant result. We now characterize
the \textbf{restriction space} $\mathcal{H}\left(  E_{1}\right)  $ as a
quotient space . This will require an extension operator.

\begin{theorem}
\label{Thm_restrict_sp_homeo_quot_space}Assume we have the notation and
results of the last two sections and suppose $E=\mathbb{R}^{d}$ and that
$E_{1}$ is open.

Then the operator
\[
\Psi:\mathcal{H}\diagup\mathcal{H}_{E_{1}^{c}}\rightarrow\mathcal{H}\left(
E_{1}\right)  ,
\]

defined by%
\[
\Psi\left[  f\right]  =r_{E_{1}}f,
\]

is an isometry. Here $\mathcal{H}_{E_{1}^{c}}$ is defined by \ref{Q014}.

Further, if there exists a continuous linear extension operator $\mathcal{E}%
:\mathcal{H}\left(  E_{1}\right)  \rightarrow\mathcal{H}$ then $\Psi$ is also
an isomorphism with inverse%
\[
\Psi^{-1}F=\left[  \mathcal{E}F\right]  ,\quad F\in\mathcal{H}\left(
E_{1}\right)  .
\]

\end{theorem}

\begin{proof}
From \ref{Q014} and the sequel, $\mathcal{H}_{\mathbb{R}^{d}\setminus E_{1}}$
is a closed subspace of $\mathcal{H}$ and their quotient space is defined.
Regarding the quotient space coset%
\begin{align*}
f\sim g\Leftrightarrow f-g\in\mathcal{H}_{E\setminus E_{1}}\Leftrightarrow
\operatorname*{supp}\left(  f-g\right)  \subset E\setminus E_{1} &
\Leftrightarrow f=g\text{ }on\text{ }E_{1}\\
&  \Leftrightarrow r_{E_{1}}f=r_{E_{1}}g,
\end{align*}

and so $\Psi$ is well-defined.

From Definition \ref{Def_QuotientSpace}, $\left\Vert \left[  f\right]
\right\Vert =\min_{h\sim f}\left\Vert h\right\Vert _{\mathcal{H}}$ so that
from \ref{Q013},%
\begin{align*}
\left\Vert \Psi\left[  f\right]  \right\Vert _{\mathcal{H}\left(
E_{1}\right)  }=\operatorname*{dist}\left(  \left\{  0\right\}  ,r_{E_{1}%
}^{-1}\Psi\left[  f\right]  \right)  =\operatorname*{dist}\left(  \left\{
0\right\}  ,r_{E_{1}}^{-1}r_{E_{1}}f\right)   & =\min_{r_{E_{1}}g=r_{E_{1}}%
f}\left\Vert g\right\Vert _{\mathcal{H}}\\
& =\min_{g\sim f}\left\Vert g\right\Vert _{\mathcal{H}}\\
& =\left\Vert \left[  f\right]  \right\Vert ,
\end{align*}

and hence $\Psi$ is an isometry. Now assume the extension exists. Since $\Psi$
is an isometry it is already 1-1 and since $\mathcal{E}$ is an extension
\[
\Psi\Psi^{-1}F=\Psi\left[  \mathcal{E}F\right]  =r_{E_{1}}\mathcal{E}F=F,
\]

which means that $\Psi$ is onto and hence an isomorphism.
\end{proof}

\chapter{Notes on Schmeisser \cite{Schmeis2006}\label{Ch_NotesOnSchmeisser}}

?? This chapter contains... ADD MORE! ??

\begin{remark}
\label{Rem_Schmeis_DomMixDerivSobol}\textbf{From the 2006 survey concerning
Sobolev spaces with dominating mixed derivatives by Schmeisser
\cite{Schmeis2006}.\medskip}

\fbox{\textbf{Subsection 1.1}} Suppose $1<p<\infty$ and $\widetilde{r}%
,\widetilde{n}\geq\mathbf{0}$ are multi-indexes on $\mathbb{R}^{d}$. Here
$\widetilde{r}$ will partition a given multi-index $\alpha$ by $\alpha=\left(
\alpha^{j}\right)  $ where $\left\vert \alpha^{j}\right\vert =r_{j}$.

Then in \textbf{Definition 1.2} Schneisser introduces the Sobolev spaces with
dominating mixed derivatives given by
\begin{align*}
S_{p}^{\widetilde{r}}W\left(  \mathbb{R}^{n_{1}}\mathbb{\times}\ldots
\times\mathbb{R}^{n_{d}}\right)   & =\left\{  f\in L^{p}\left(  \mathbb{R}%
^{n_{1}}\mathbb{\times}\ldots\times\mathbb{R}^{n_{d}}\right)  :D^{\alpha}f\in
L^{p}\left(  \mathbb{R}^{n_{1}}\mathbb{\times}\ldots\times\mathbb{R}^{n_{d}%
}\right)  ,\text{ }\left\vert \alpha^{j}\right\vert \leq r_{j}\text{ }\forall
j\right\}  ,\\
\left\Vert f\right\Vert _{S_{p}^{\widetilde{r}}W}  & =\sum_{\substack{j=1
\\\left\vert \alpha^{j}\right\vert \leq r_{j}}}^{d}\left\Vert D^{\alpha
}f\right\Vert _{p}.
\end{align*}

In \textbf{Remark 1.3}: Sobolev spaces with dominating mixed derivatives are
well adapted to tensor products of functions.

\textbf{After equation 1.3}: the special case where $\widetilde{r}%
=r\mathbf{1}$ and $\widetilde{n}=\mathbf{1}$ and defines
\begin{align*}
S_{p}^{r}W\left(  \mathbb{R}^{d}\right)   & :=S_{p}^{r\mathbf{1}}W\left(
\underset{d\text{ }times}{\underbrace{\mathbb{R}^{1}\mathbb{\times}%
\ldots\times\mathbb{R}^{1}}}\right) \\
& =\left\{  f\in L^{p}\left(  \mathbb{R}^{d}\right)  :D^{\alpha}f\in
L^{p}\left(  \mathbb{R}^{d}\right)  \text{ }\forall\alpha\leq r\mathbf{1}%
\right\} \\
& =\left\{  f\in L^{p}:D^{\alpha}f\in L^{p}\text{ }\forall\alpha\leq
r\mathbf{1}\right\}  .
\end{align*}

for compactness.

Schneisser mentions $S_{p}^{r}W\left(  \mathbb{\Omega}\right)  $ spaces have
been studied on the torus $\mathbb{T}^{d}=\left[  0,1\right]  ^{d}$ and can be
easily defined on domains $\Omega\subset\mathbb{R}^{d}$.\textbf{\medskip}

\fbox{\textbf{Subsection 1.2}} By using the Fourier transform
\textbf{Definition 1.5} introduces \textbf{fractional Sobolev spaces with
dominating mixed smoothness }on $\mathbb{R}^{d}$.\textbf{\medskip}

\fbox{\textbf{Subsection 1.3}} \textbf{Multivariate approximation}. Schneisser
discusses how we can benefit from dominating mixed smoothness. To this end he
considers an example related to multivariate approximation.\textbf{\medskip}

\fbox{\textbf{Subsection 1.4}} \textbf{Sobolev embeddings}.\textbf{\medskip}

\fbox{\textbf{Subsection 1.5}} \textbf{Traces} e.g. the \textbf{diagonal
trace} $f\left(  x\mathbf{1}\right)  $.\textbf{\medskip}

\fbox{\textbf{Subsection 1.6}} \textbf{Nikol'skii-Besov spaces} Dominating
mixed smoothness properties can also be described in a classical way using
mixed differences or mixed moduli of continuity. \textbf{Definition 1.11}
introduces \textbf{Nikol'skii spaces} with dominating mixed smoothness.
\textbf{Definition 1.12} concerns Besov spaces with dominating mixed
smoothness.\medskip

There are 5 sections in all but no table of contents.\textbf{\medskip}

\textbf{Section 2} \textbf{Littlewood-Paley analysis} -- the Fourier
analytical approach. \textbf{2.1} Littlewood-Paley theorems; \textbf{2.2}
Dominating mixed derivatives and differences; \textbf{2.3}. Besov and
Lizorkin-Triebel spaces.

\textbf{Section 3} \textbf{Approximation with respect to hyperbolic crosses}.
\textbf{3.1} Approximation spaces; \textbf{3.2} Comparison theorems.

\textbf{Section 4} \textbf{New instruments}. Recent developments in the theory
of function spaces are characterized by the use of new instruments such as
\textbf{local means}, \textbf{atoms and wavelets}, which are successfully
applied to tackle \textbf{key problems} such as \textbf{limiting embeddings},
\textbf{traces, extensions, entropy numbers or sampling numbers}. Moreover,
\textbf{far-reaching applications} to \textbf{fractals} and to
\textbf{spectral theory} have been achieved. The aim of this section is to
develop these new tools for spaces with dominating mixed smoothness. It yields
the basis for applications which will be described in the next section.

\textbf{Section 5} \textbf{Applications} to \textbf{traces} and
\textbf{entropy numbers}.
\end{remark}

\chapter{Proofs of the claims made in Remark \ref{Rem_SobolevSpace2}
\label{Ch_Hm1_proofs_Adams}}

This chapter contains the proofs of the claims made for the spaces
$W^{m\mathbf{1}}\left(  \Omega\right)  $, $W_{\overline{\Omega}}^{m\mathbf{1}%
}$ in Remark \ref{Rem_SobolevSpace2}.

\begin{lemma}
\label{Lem_Thm_2.29_AdamF_mollifier}\textbf{Theorem 2.29 of AF
\cite{AdamFour2003}:} \textbf{properties of mollification} Suppose
$\operatorname*{supp}u\subseteq\overline{\Omega}$ where $\Omega$ is an open
subset of $\mathbb{R}^{d}$. Then:

\begin{description}
\item[(a)] If $u\in L_{loc}^{1}$ then $J_{\varepsilon}\ast u\in C^{\infty}$.

\item[(b)] If $u\in L_{loc}^{1}\left(  \Omega\right)  $ and
$\operatorname*{supp}u\Subset\Omega$ then $J_{\varepsilon}\ast u\in
C_{0}^{\infty}\left(  \Omega\right)  $ when $\varepsilon<\operatorname*{dist}%
\left(  \operatorname*{supp}u,\partial\Omega\right)  $.

\item[(c)] If $u\in L^{p}\left(  \Omega\right)  $ where $1\leq p<\infty$, then
$J_{\varepsilon}\ast u\in L^{p}\left(  \Omega\right)  $. Also%
\[
\left\Vert J_{\varepsilon}\ast u\right\Vert _{p}\leq\left\Vert u\right\Vert
_{p},\quad\lim_{\varepsilon\rightarrow0^{+}}\left\Vert J_{\varepsilon}\ast
u-u\right\Vert _{p}=0.
\]

\item[(d)] If $u\in C\left(  \Omega\right)  $ and $G\Subset\Omega$, then
$\lim\limits_{\varepsilon\rightarrow0^{+}}\left(  J_{\varepsilon}\ast
u\right)  \left(  x\right)  =u\left(  x\right)  $ uniformly on $G$.

\item[(e)] If $u\in C\left(  \overline{\Omega}\right)  $ then $\lim
\limits_{\varepsilon\rightarrow0^{+}}\left(  J_{\varepsilon}\ast u\right)
\left(  x\right)  =u\left(  x\right)  $ uniformly on $\Omega$.
\end{description}
\end{lemma}

??

\begin{lemma}
\label{Lem_3.16_AdamsF}\textbf{Mollification in }$W^{m\mathbf{1}}\left(
\Omega\right)  $ (\textbf{Modification of Lemma 3.16 of AF \cite{AdamFour2003}%
}) Let $J_{\varepsilon}$ be defined as in Paragraph 2.28 of AF
\cite{AdamFour2003} and let $u\in W^{m\mathbf{1}}\left(  \Omega\right)  $.

Then if $\Omega^{\prime}\subset\Omega$ is a subdomain then $\lim
\limits_{\varepsilon\rightarrow0^{+}}J_{\varepsilon}\ast u=u$ in
$W^{m\mathbf{1}}\left(  \Omega^{\prime}\right)  $.
\end{lemma}

\begin{proof}
Easy modification of the proof in Adams and Fournier.
\end{proof}

\subsection{Part 1 of Remark \ref{Rem_SobolevSpace2}}

??

\subsection{Part 2 of Remark \ref{Rem_SobolevSpace2}}

??

\subsection{Part 3 of Remark \ref{Rem_SobolevSpace2}}

??

\subsection{Part 4 of Remark \ref{Rem_SobolevSpace2}}

??

\begin{theorem}
\label{Thm_3.17_AdamsF}\textbf{Modified Theorem 3.17 of AF}
\textbf{\cite{AdamFour2003}} For any open set $\Omega\subset\mathbb{R}^{d}$ we
have $C^{\left(  m\mathbf{1}\right)  }\left(  \Omega\right)  \cap
W^{m\mathbf{1}}\left(  \Omega\right)  $ is dense in $W^{m\mathbf{1}}\left(
\Omega\right)  $. In fact we show that $C^{\infty}\left(  \Omega\right)  \cap
W^{m\mathbf{1}}\left(  \Omega\right)  $ is dense in $W^{m\mathbf{1}}\left(
\Omega\right)  $.
\end{theorem}

\begin{proof}
The proof will use Theorem 2.29 of \cite{AdamFour2003}\textbf{\ }which gives
some \textbf{mollifier properties} and Lemma 3.15 of \cite{AdamFour2003} which
concerns $C^{\infty}$ \textbf{subordinate partitions of unity}. These theorems
are independent of the Sobolev spaces.

For $k=1,2,3,\ldots$ let%
\[
\Omega_{k}=\left\{  x\in\Omega:\operatorname*{dist}\left(  x,\partial
\Omega\right)  >1/k\right\}  ,
\]

and let $\Omega_{0}=\Omega_{-1}=\emptyset$, the empty set. Then%
\[
\mathcal{O=}\left\{  U_{k}:U_{k}=\Omega_{k+1}\cap\left(  \overline
{\Omega_{k-1}}\right)  ^{c},\text{ }k=1,2,3,\ldots\right\}  ,
\]

is a collection of open subsets of $\Omega$ that covers $\Omega$. Let $\Psi$
be a $C^{\infty}$-partition of unity of $\Omega$ subordinate to $\mathcal{O}$.
Let $\psi_{k}$ denote the sum of the finitely many functions $\psi\in\Psi$
whose supports are contained in $U_{k}$. Then $\psi_{k}\in C_{0}^{\infty
}\left(  U_{k}\right)  $ and $\sum_{k=1}^{\infty}\psi_{k}\left(  x\right)  =1$
on $\Omega$.

If $0<\varepsilon<1/\left(  k+1\right)  \left(  k+2\right)  $ then from
Theorem 2.29 \cite{AdamFour2003}, $J_{\varepsilon}\ast\left(  \psi
_{k}u\right)  $ has support in $V_{k}=\Omega_{k+2}\cap\Omega_{k-2}^{c}%
\Subset\Omega$. Since $\psi_{k}u\in W^{m\mathbf{1}}\left(  \Omega\right)  $,
Lemma \ref{Lem_3.16_AdamsF} allows us to choose $\varepsilon_{k}$ satisfying
$0<\varepsilon_{k}<1/\left(  k+1\right)  \left(  k+2\right)  $ such that%
\[
\left\Vert J_{\varepsilon_{k}}\ast\left(  \psi_{k}u\right)  -\psi
_{k}u\right\Vert _{m\mathbf{1},2,\Omega}=\left\Vert J_{\varepsilon_{k}}%
\ast\left(  \psi_{k}u\right)  -\psi_{k}u\right\Vert _{m\mathbf{1},2,V_{k}%
}<2^{-k}\varepsilon.
\]

Let $\phi=\sum_{k=1}^{\infty}J_{\varepsilon_{k}}\ast\left(  \psi_{k}u\right)
$. On any $\Omega^{\prime}\Subset\Omega$ only finitely many terms in the sum
can be nonzero. Thus by Theorem 2.29 \cite{AdamFour2003} $\phi\in C^{\infty
}\left(  \Omega\right)  $. For $x\in\Omega_{k}$, we have%
\[
u\left(  x\right)  =\sum_{j=1}^{k+2}\psi_{j}\left(  x\right)  u\left(
x\right)  ,\quad\phi\left(  x\right)  =\sum_{j=1}^{k+2}\left(  J_{\varepsilon
_{j}}\ast\left(  \psi_{j}u\right)  \right)  \left(  x\right)  ,
\]

and thus%
\[
\left\Vert u-\phi\right\Vert _{m\mathbf{1},2,\Omega_{k}}\leq\sum_{j=1}%
^{k+2}\left\Vert J_{\varepsilon_{k}}\ast\left(  \psi_{j}u\right)  -\psi
_{j}u\right\Vert _{m\mathbf{1},2,\Omega}<\varepsilon.
\]

By the monotone convergence theorem, $\left\Vert u-\phi\right\Vert
_{m\mathbf{1},2,\Omega}<\varepsilon$.
\end{proof}

\begin{theorem}
\label{Thm_3.22_AdamsF}\textbf{Modified Theorem 3.22 of AF}
\textbf{\cite{AdamFour2003}} If $\Omega$ is bounded and satisfies the segment
condition of Definition \ref{Def_SegCondit} then $r_{\Omega}C_{0}^{\infty}$ is
dense in $W^{m\mathbf{1}}\left(  \Omega\right)  $ and $r_{\Omega}C_{0}%
^{\infty}=r_{\Omega}C^{\infty}=r_{\Omega}C_{B}^{\infty}$.
\end{theorem}

\begin{proof}
Choose $u\in W^{m\mathbf{1}}\left(  \Omega\right)  $. Let
$K=\operatorname*{supp}u\cap\Omega$. The set $F=\overline{K}\setminus
\bigcup_{x\in\partial\Omega}U_{x}$ is compact and contained in $\Omega$,
$\left\{  U_{x}\right\}  $ being the collection of open sets referred to in
the definition of the segment condition. There also exists an open set $U_{0}$
such that $F\Subset U_{0}\Subset\Omega$. Since $\overline{K}$ is compact,
there exist finitely many of the sets $U_{x}$, let us rename them $\left\{
U_{i}\right\}  _{i=1}^{k}$, such that $\overline{K}\subset U_{0}\cup U_{1}%
\cup\ldots\cup U_{k}$. Moreover, there are other open sets $V_{0},V_{1}%
,\ldots,V_{k}$ such that $V_{j}\Subset U_{j}$ for $0\leq j\leq k$ but still
$\overline{K}\subset V_{0}\cup V_{1}\cup\ldots\cup V_{k}$.

Let $\Psi$ be a $C^{\infty}$-partition of unity subordinate to $\left\{
V_{j}\right\}  _{j=1}^{k}$, and let $\psi_{j}$ be the sum of the finitely many
functions $\psi\in\Psi$ which have supports in $V_{j}$. Let $u_{j}=\psi_{j}u$.
Suppose that for each $j$ we can find $\phi_{j}\in C_{0}^{\infty}\left(
\mathbb{R}^{d}\right)  $ such that%
\begin{equation}
\left\Vert u_{j}-\phi_{j}\right\Vert _{m\mathbf{1},2,\Omega}<\frac
{\varepsilon}{k+1}.\label{X03}%
\end{equation}

Then, putting $\phi=\sum_{j=0}^{k}\phi_{j}$, we would obtain%
\[
\left\Vert u-\phi\right\Vert _{m\mathbf{1},2,\Omega}\leq\sum_{j=0}%
^{k}\left\Vert u_{j}-\phi_{j}\right\Vert _{m\mathbf{1},2,\Omega}<\varepsilon.
\]

Since $\operatorname*{supp}u_{0}\subset V_{0}\Subset\Omega$, a function
$\phi_{0}\in C_{0}^{\infty}$ satisfying \ref{X03} for $j=0$ can be found via
the mollification Lemma \ref{Lem_3.16_AdamsF}.

It remains, therefore, to find $\phi_{j}$ satisfying \ref{X03} for $1\leq
j\leq k$. For fixed such $j$ we extend $u_{j}$ to be identically zero outside
$\Omega$. Thus $u_{j}\in W^{m\mathbf{1}}\left(  \mathbb{R}^{d}\setminus
\Gamma\right)  $, where $\Gamma=\overline{V_{j}}\cap\operatorname*{bdry}%
\Omega$. Let $y^{j}$ be the nonzero vector associated with the set $U_{j}$ in
the definition of the segment condition. Let $\Gamma_{t}=\left\{
x-ty^{j}:x\in\Gamma\right\}  $, where $t$ is so chosen that%
\[
0<t<\min\left\{  1,\operatorname*{dist}\left(  V_{j},\mathbb{R}^{d}%
\setminus\Gamma\right)  /\left\vert y^{j}\right\vert \right\}  .
\]

Then $\Gamma_{t}\subset U_{j}$ and $\Gamma_{t}\cap\overline{\Omega}$ is empty
by the segment condition. Now we use the segment condition to define the
\textbf{outwardly displaced} functions
\begin{equation}
u_{j,t}\left(  x\right)  :=u_{j}\left(  x+ty^{j}\right)  ,\label{X04}%
\end{equation}

and observe that $u_{j,t}\in W^{m\mathbf{1}}\left(  \mathbb{R}^{d}%
\setminus\Gamma_{t}\right)  $. Translation is continuous in $L^{2}$ so
$D^{\alpha}u_{j,t}\rightarrow D^{\alpha}u_{j}$ in $L^{2}$ as $t\rightarrow
0^{+}$ for $\alpha\leq m\mathbf{1}$, and so it is sufficient to find $\phi
_{j}\in C_{0}^{\infty}$ such that $\left\Vert u_{j,t}-\phi_{j}\right\Vert
_{m\mathbf{1},2}$ is sufficiently small. However, $\Omega\cap U_{j}%
\Subset\mathbb{R}^{d}\setminus\Gamma_{t}$, and so by the modified
mollification Lemma \ref{Lem_3.16_AdamsF} we set $\phi_{j}=J_{\delta}\ast
u_{j,t}$ for suitably small $\delta>0$.

\textbf{It remains to be proven that} $r_{\Omega}C_{0}^{\infty}=r_{\Omega
}C^{\infty}=r_{\Omega}C_{B}^{\infty}$. Clearly $r_{\Omega}C_{0}^{\infty
}\subset r_{\Omega}C^{\infty}$. Conversely, suppose $\psi\in r_{\Omega}\phi$
with $\phi\in C^{\infty}$. Since $\Omega$ is bounded there exists a ball
$B_{R}$ such that $\Omega\subset B_{R}$ and we can choose $\theta\in
C_{0}^{\infty}$ such that $\theta=1$ on $B_{R}$ and $\theta=0$ and
$\mathbb{R}^{d}\setminus B_{2R}$. Hence $\theta\phi\in C_{0}^{\infty}$ and
$\psi\in r_{\Omega}\left(  \theta\phi\right)  $ and so $r_{\Omega}%
C_{0}^{\infty}=r_{\Omega}C^{\infty}$ which in turn directly yields $r_{\Omega
}C^{\infty}=r_{\Omega}C_{B}^{\infty}$.
\end{proof}

\subsection{Part 6 of Remark \ref{Rem_SobolevSpace2}}

Recall that from part 9 of Remark \ref{Rem_SobolevSpace2},%
\[
W_{\overline{\Omega}}^{m\mathbf{1}}=\left\{  u\in W^{m\mathbf{1}%
}:\operatorname*{supp}u\subset\overline{\Omega}\right\}  ,\quad\left\Vert
u\right\Vert _{W_{\overline{\Omega}}^{m\mathbf{1}}}=\left\Vert u\right\Vert
_{W^{m\mathbf{1}}}.
\]

\begin{lemma}
\label{Lem_3.27_AdamsF}\textbf{Modified Lemma 3.27 of Adams and Fournier}
\textbf{\cite{AdamFour2003}} Let $\mathcal{E}_{0}u$ denote the \textbf{zero
extension} of $u\in W_{0}^{m\mathbf{1}}\left(  \Omega\right)  $ outside
$\Omega$. Then for $\alpha\leq m\mathbf{1}$, $D^{\alpha}\mathcal{E}%
_{0}u=\mathcal{E}_{0}D^{\alpha}u$ in $\mathcal{D}^{\prime}$, and
$\mathcal{E}_{0}u\in W^{m\mathbf{1}}\left(  \mathbb{R}^{d}\right)  $ with
$\left\Vert \mathcal{E}_{0}u\right\Vert _{W^{m\mathbf{1}}}=\left\Vert
u\right\Vert _{W_{0}^{m\mathbf{1}}\left(  \Omega\right)  }$ and
$\operatorname*{supp}u\subset\overline{\Omega}$.

In other words, $\mathcal{E}_{0}:W_{0}^{m\mathbf{1}}\left(  \Omega\right)
\rightarrow W_{\overline{\Omega}}^{m\mathbf{1}}$\ is an isometry.
\end{lemma}

\begin{proof}
Easy modification of the proof of Lemma 3.27.
\end{proof}

\subsection{Part 8 of Remark \ref{Rem_SobolevSpace2}}

\begin{theorem}
\label{Thm_norm_Cm1xHm1_loc} Suppose $\Omega$ is a bounded open set. Then
$\phi\in W^{m\mathbf{1},\infty}\left(  \Omega\right)  $ and $u\in
W^{m\mathbf{1}}\left(  \Omega\right)  $ implies $\phi u\in W^{m\mathbf{1}%
}\left(  \Omega\right)  $ and
\[
\left\Vert \phi u\right\Vert _{m\mathbf{1},2,\Omega}\leq2^{md/2}\left\Vert
\phi\right\Vert _{m\mathbf{1},\infty,\Omega}\left\Vert u\right\Vert
_{m\mathbf{1},2,\Omega}.
\]

\end{theorem}

\begin{proof}
Part 2 of Remark \ref{Rem_Thm_data_func} discusses weakening the distribution
Leibniz theorem to: $u,Du\in L^{p}\left(  \Omega\right)  $ and $v,Dv\in
L^{q}\left(  \Omega\right)  $ where $p^{-1}+q^{-1}=1$ and $q<\infty$, implies
$D\left(  uv\right)  =\left(  Du\right)  v+uDv\in L^{2}\left(  \Omega\right)
$.

Suppose $p=\infty$ and $q=1$. Suppose $\Omega$ is bounded. Then, since
$\Omega$ is bounded, $v,Dv\in L^{2}\left(  \Omega\right)  $ implies $v,Dv\in
L^{1}\left(  \Omega\right)  $ and so Leibniz theorem holds for $u,Du\in
L^{\infty}\left(  \Omega\right)  $ and $v,Dv\in L^{2}\left(  \Omega\right)  $.
Consequently Leibniz theorem holds for $u\in W^{m\mathbf{1},\infty}\left(
\Omega\right)  $ and $v\in W^{m\mathbf{1}}\left(  \Omega\right)  $. Now, if
$\tbinom{\mathbf{0}}{\mathbf{0}}:=1$,%
\begin{align*}
\left\Vert \phi u\right\Vert _{m\mathbf{1},2,\Omega}^{2}  & =\sum_{\alpha\leq
m\mathbf{1}}\tbinom{m\mathbf{1}}{\alpha}\left\Vert D^{\alpha}\left(  \phi
u\right)  \right\Vert _{L^{2}\left(  \Omega\right)  }^{2}\\
& =\sum_{\alpha\leq m\mathbf{1}}\tbinom{m\mathbf{1}}{\alpha}\left\Vert
\sum_{\beta\leq\alpha}\tbinom{\alpha}{\beta}D^{\alpha-\beta}\phi D^{\beta
}u\right\Vert _{L^{2}\left(  \Omega\right)  }^{2}\\
& \leq\sum_{\alpha\leq m\mathbf{1}}\tbinom{m\mathbf{1}}{\alpha}\left\Vert
\sum_{\beta\leq\alpha}\tbinom{\alpha}{\beta}\left\vert D^{\alpha-\beta}%
\phi\right\vert \text{ }\left\vert D^{\beta}u\right\vert \right\Vert
_{L^{2}\left(  \Omega\right)  }^{2}\\
& \leq\sum_{\alpha\leq m\mathbf{1}}\tbinom{m\mathbf{1}}{\alpha}\left\Vert
\left(  \sum_{\beta\leq\alpha}\tbinom{\alpha}{\beta}\left\vert D^{\alpha
-\beta}\phi\right\vert ^{2}\right)  ^{1/2}\left(  \sum_{\beta\leq\alpha
}\tbinom{\alpha}{\beta}\left\vert D^{\beta}u\right\vert ^{2}\right)
^{1/2}\right\Vert _{L^{2}\left(  \Omega\right)  }^{2}\\
& =\sum_{\alpha\leq m\mathbf{1}}\tbinom{m\mathbf{1}}{\alpha}\left\Vert \left(
\sum_{\beta\leq\alpha}\tbinom{\alpha}{\beta}\left\vert D^{\beta}%
\phi\right\vert ^{2}\right)  ^{1/2}\left(  \sum_{\beta\leq\alpha}%
\tbinom{\alpha}{\beta}\left\vert D^{\beta}u\right\vert ^{2}\right)
^{1/2}\right\Vert _{L^{2}\left(  \Omega\right)  }^{2}\\
& =\sum_{\alpha\leq m\mathbf{1}}\tbinom{m\mathbf{1}}{\alpha}\int_{\Omega
}\left(  \sum_{\beta\leq\alpha}\tbinom{\alpha}{\beta}\left\vert D^{\beta}%
\phi\right\vert ^{2}\right)  \left(  \sum_{\beta\leq\alpha}\tbinom{\alpha
}{\beta}\left\vert D^{\beta}u\right\vert ^{2}\right) \\
& \leq\sum_{\alpha\leq m\mathbf{1}}\tbinom{m\mathbf{1}}{\alpha}\int_{\Omega
}\left(  \sum_{\beta\leq\alpha}\tbinom{\alpha}{\beta}\left\Vert D^{\beta}%
\phi\right\Vert _{\infty}^{2}\right)  \left(  \sum_{\beta\leq\alpha}%
\tbinom{\alpha}{\beta}\left\vert D^{\beta}u\right\vert ^{2}\right) \\
& =\sum_{\alpha\leq m\mathbf{1}}\tbinom{m\mathbf{1}}{\alpha}\left(
\sum_{\beta\leq\alpha}\tbinom{\alpha}{\beta}\left\Vert D^{\beta}%
\phi\right\Vert _{\infty}^{2}\right)  \left(  \sum_{\beta\leq\alpha}%
\tbinom{\alpha}{\beta}\int_{\Omega}\left\vert D^{\beta}u\right\vert
^{2}\right) \\
& \leq\left(  \sum_{\alpha\leq m\mathbf{1}}\tbinom{m\mathbf{1}}{\alpha
}\right)  \left(  \sum_{\beta\leq m\mathbf{1}}\tbinom{m\mathbf{1}}{\beta
}\left\Vert D^{\beta}\phi\right\Vert _{\infty}^{2}\right)  \left(  \sum
_{\beta\leq m\mathbf{1}}\tbinom{m\mathbf{1}}{\beta}\int_{\Omega}\left\vert
D^{\beta}u\right\vert ^{2}\right) \\
& =2^{md}\left(  \sum_{\beta\leq m\mathbf{1}}\tbinom{m\mathbf{1}}{\beta
}\left\Vert D^{\beta}\phi\right\Vert _{\infty}^{2}\right)  \left(  \sum
_{\beta\leq m\mathbf{1}}\tbinom{m\mathbf{1}}{\beta}\int_{\Omega}\left\vert
D^{\beta}u\right\vert ^{2}\right) \\
& =2^{md}\left\Vert \phi\right\Vert _{m\mathbf{1},\infty,\Omega}^{2}\left\Vert
u\right\Vert _{m\mathbf{1},2,\Omega}^{2}.
\end{align*}

\end{proof}

\subsection{Part 9 of Remark \ref{Rem_SobolevSpace2}}

??

\begin{theorem}
\label{Thm_Thm3.7_Wloka}\textbf{Modified Theorem 3.7 of Wloka \cite{Wloka87}}
If $\Omega$ has the segment property of Definition \ref{Def_SegCondit} then
$r_{\Omega}W_{\overline{\Omega}}^{m\mathbf{1}}=W_{0}^{m\mathbf{1}}\left(
\Omega\right)  $ as sets.
\end{theorem}

\begin{proof}
Recall that%
\[
W_{0}^{m\mathbf{1}}\left(  \Omega\right)  =closure\text{ }of\text{ }%
C_{0}^{\infty}\left(  \Omega\right)  \text{ }in\text{ }W^{m\mathbf{1}}\left(
\Omega\right)  .
\]

Suppose $v\in W_{0}^{m\mathbf{1}}\left(  \Omega\right)  $. Then by part 6 of
Remark \ref{Rem_SobolevSpace2}, $\mathcal{E}_{0}v\in W^{m\mathbf{1}}$ and
$\operatorname*{supp}\mathcal{E}_{0}v\subset\overline{\Omega}$ i.e.
$\mathcal{E}_{0}v\in W_{\overline{\Omega}}^{m\mathbf{1}}$. Thus $v=r_{\Omega
}\mathcal{E}_{0}v\in r_{\Omega}W_{\overline{\Omega}}^{m\mathbf{1}}$ and so
$W_{0}^{m\mathbf{1}}\left(  \Omega\right)  \subset r_{\Omega}W_{\overline
{\Omega}}^{m\mathbf{1}}$. Now to prove the converse.

Since $r_{\Omega}W_{\overline{\Omega}}^{m\mathbf{1}}\subset W^{m\mathbf{1}%
}\left(  \Omega\right)  $ we need to show that given $v\in r_{\Omega
}W_{\overline{\Omega}}^{m\mathbf{1}}$ there is a sequence in $C_{0}^{\infty
}\left(  \Omega\right)  $ which converges to $v$ in $W^{m\mathbf{1}}\left(
\Omega\right)  $. Suppose $v=r_{\Omega}u$ for some $u\in W_{\overline{\Omega}%
}^{m\mathbf{1}}$. We will modify the proof of Theorem \ref{Thm_3.22_AdamsF} by
not using the open sets $V_{j}$ and by using the segment condition to push the
function inwards.

Let $K=\operatorname*{supp}u\cap\Omega$. The set $F=\overline{K}%
\setminus\bigcup_{x\in\partial\Omega}U_{x}$ is compact and contained in
$\Omega$, $\left\{  U_{x}\right\}  $ being the collection of open sets
referred to in the definition of the segment condition. There also exists an
open set $U_{0} $ such that $F\Subset U_{0}\Subset\Omega$. Since $\overline
{K}$ is compact, there exist finitely many of the sets $U_{x}$, let us rename
them $\left\{  U_{i}\right\}  _{i=1}^{k}$, such that $\overline{K}\subset
U_{0}\cup U_{1}\cup\ldots\cup U_{k}$.

Let $\Psi$ be a $C^{\infty}$-partition of unity subordinate to $\left\{
U_{j}\right\}  _{j=1}^{k}$ and let $\psi_{j}$ be the sum of the finitely many
functions $\psi\in\Psi$ which have supports in $U_{j}$. Let $u_{j}=\psi
_{j}u\in W^{m\mathbf{1}}\left(  \Omega\right)  $ by part 8 of Remark
\ref{Rem_SobolevSpace2}. Suppose that for each $j$ we can find $\phi_{j}\in
C_{0}^{\infty}\left(  \mathbb{R}^{d}\right)  $ such that%
\begin{equation}
\left\Vert u_{j}-\phi_{j}\right\Vert _{m\mathbf{1},2,\Omega}<\frac
{\varepsilon}{k+1}.\label{X07}%
\end{equation}

Then, putting $\phi=\sum_{j=0}^{k}\phi_{j}$, we would obtain%
\[
\left\Vert u-\phi\right\Vert _{n\mathbf{1},2,\Omega}\leq\sum_{j=0}%
^{k}\left\Vert u_{j}-\phi_{j}\right\Vert _{n\mathbf{1},2,\Omega}<\varepsilon.
\]

Since $\operatorname*{supp}u_{0}\subset U_{0}\Subset\Omega$, a function
$\phi_{0}\in C_{0}^{\infty}$ satisfying \ref{X07} for $j=0$ can be found via
the mollification Lemma \ref{Lem_3.16_AdamsF}. It remains, therefore, to find
$\phi_{j}$ satisfying \ref{X07} for $1\leq j\leq k$. Because of the segment
property the\textbf{\ inwardly displaced} functions%
\begin{equation}
u_{j,t}:=u_{j}\left(  \cdot-ty^{j}\right)  ,\quad0<t<1,\text{ }1\leq j\leq
k,\label{X10}%
\end{equation}

have their supports in $\Omega$. Here each $y^{j}\neq\mathbf{0}$ is the
translation vector required by the segment condition. Translation is
continuous in $L^{2}$ so for $\alpha\leq m\mathbf{1}$, $D^{\alpha}%
u_{j,t}\rightarrow D^{\alpha}u_{j}$ in $L^{2}\left(  \Omega\right)  $ as
$t\rightarrow0^{+}$, and we can find some $t_{\varepsilon}$ with%
\[
\left\Vert u_{j}-u_{j,t_{\varepsilon}}\right\Vert _{m\mathbf{1},2,\Omega
}<\frac{\varepsilon/2}{k+1},\quad1\leq j\leq k.
\]

Further, $\operatorname*{supp}u_{j,t_{\varepsilon}}\subset\Omega$ so if
$\eta=\frac{1}{2}\operatorname*{dist}\left(  \operatorname*{supp}%
u_{j,t_{\varepsilon}},\partial\Omega\right)  $ then $\Omega^{\prime}=\left(
\operatorname*{supp}u_{j,t_{\varepsilon}}\right)  _{\eta}\Subset\Omega$ and
$\operatorname*{dist}\left(  \Omega^{\prime},\partial\Omega\right)  \leq\eta$.
Since $u_{j,t_{\varepsilon}}\in L_{loc}^{1}\left(  \Omega\right)  $ and
$\operatorname*{supp}u_{j,t_{\varepsilon}}\Subset\Omega^{\prime}$, part (b)
Theorem 2.29 of Adams and F. \cite{AdamFour2003} implies $J_{\delta}\ast
u_{j,t_{\varepsilon}}\in C_{0}^{\infty}\left(  \Omega^{\prime}\right)  $ when
$\delta<\eta$. Also, Lemma \ref{Lem_3.16_AdamsF} implies $J_{\delta}\ast
u_{j,t_{\varepsilon}}\rightarrow u_{j,t_{\varepsilon}}$ in $W^{m\mathbf{1}%
}\left(  \Omega^{\prime}\right)  $ as $\delta\rightarrow0^{+}$. Thus there
exists $\phi_{j}\in C_{0}^{\infty}\left(  \Omega\right)  $ such that
$\left\Vert u_{j,t_{\varepsilon}}-\phi_{j}\right\Vert <\frac{\varepsilon
/2}{k+1}$ for $1\leq j\leq k$ and the proof is complete.
\end{proof}

\bibliographystyle{amsplain}
\bibliography{acompat,BasisFunc}

\newif\ifabfull\abfulltrue
\providecommand{\bysame}{\leavevmode\hbox to3em{\hrulefill}\thinspace}
\providecommand{\MR}{\relax\ifhmode\unskip\space\fi MR }
\providecommand{\MRhref}[2]{%
  \href{http://www.ams.org/mathscinet-getitem?mr=#1}{#2}
}
\providecommand{\href}[2]{#2}
\begin{thebibliography}{10}

\bibitem{KinOud2010}
M.~Oudadess A.~Kinani, \emph{Distribution theory and applications}, Series on
  Concrete and Applicable Mathematics, vol.~9, World Scientific Publishing,
  2010.

\bibitem{PrudBryMar86}
Yu.~Brychkov A.~Prudnikov and O.~Marichevl, \emph{Integrals and series}, vol.
  1: Elementary Functions, Gordon And Breach Science Publishers, New York,
  1986.

\bibitem{AbramowStegun70}
M.~Abramowitz and I.~Stegun, \emph{Handbook of mathematical functions}, Dover
  Publications, 1970.

\bibitem{Adams75}
R.~A. Adams, \emph{Sobolev spaces}, 1st ed., Pure and applied mathematics
  series, vol.~65, Academic Press, London, 1975.

\bibitem{AdamFour2003}
R.~A. Adams and J.~J. Fournier, \emph{Sobolev spaces}, 2nd ed., Pure and
  applied mathematics series, vol. 140, Academic Press, London, 2003.

\bibitem{Arfken70}
G.~Arfken, \emph{Mathematical methods for physicists}, 2nd ed., Academic Press,
  New York, London, 1970.

\bibitem{Aronszajn50}
N.~Aronszajn, \emph{Theory of reproducing kernels}, Trans. Amer. Math. Soc.
  \textbf{68} (1950), no.~3, 337--404.

\bibitem{Bauer2001}
H.~Bauer, \emph{Measure and integration theory}, De Gruyter studies in
  mathematics, vol.~26, Walter de Gruyter, Berlin, New York, 2001.

\bibitem{BesovEtAl90}
O.V. Besov, L.~Kudrayavtsev, P.~Lizorkin, and S.~Nikol'ski\u{i},
  \emph{Investigations in the theory of spaces of differentiable functions of
  several variables}, Proc. Steklov Inst. Math. \textbf{182} (1990), no.~1,
  73--139.

\bibitem{Bojarski2011}
B.~Bojarski, \emph{Taylor expansions and sobolev spaces}, Bull. Georg. Natl.
  Acad. Sci. \textbf{5} (2011), no.~2, 5--10.

\bibitem{Bojarski2012}
\bysame, \emph{Sobolev spaces and {L}agrange interpolation}, arXiv:1201.4708v1
  [math.AP], 2012.

\bibitem{BojarHaj93}
B.~Bojarski and P.~Hajlasz, \emph{Pointwise inequalities for {S}obolev
  functions and some applications}, Studia Math. \textbf{106} (1993), no.~1,
  77--92.

\bibitem{BojHajStr2002}
B.~Bojarski, P.~Hajlasz, and P.~Strzelecki, \emph{Improved approximation of
  higher order {S}obolev functions in norm and capacity}, Indiana Univ. Math.
  J. \textbf{51} (2002), 507--540.

\bibitem{BungGrieb04}
H-J. Bungartz and M.~Griebel, \emph{Sparse grids}, Acta. Num. \textbf{13}
  (2004), 147--269.

\bibitem{Chern2012}
A.~Chernov, \emph{Sparse polynomial approximation in positive order {S}obolev
  spaces with bounded mixed derivatives and applications to elliptic problems
  with random loading}, Appl. Num. Math. \textbf{62} (2012), no.~4, 360--377.

\bibitem{DautLion88}
R.~Dautray and J-L Lions, \emph{Functional and variational methods},
  Mathematical Analysis and Numerical Methods for Science and Technology,
  vol.~2, Springer-Verlag (Berlin), 1988.

\bibitem{Duchon77}
J.~Duchon, \emph{Splines minimising rotation-invariant seminorms in {S}obolev
  spaces}, Constructive theory of functions of several variables (W.~Schempp
  and K.~Zeller, eds.), Lecture Notes in Mathematics, vol. 571, Springer-Verlag
  (Berlin), 1977, pp.~85--100.

\bibitem{DuistKolk2010}
J.~Duistermaat and J.~Kolk, \emph{Distributions: theory and applications},
  Cornerstones, vol.~16, Springer Science+Business Media, NY, 2010.

\bibitem{Dzhab67}
A.~D. D\u{z}abrailov, \emph{On the theory of "imbedding theorems"}, Proc.
  Steklov Inst. Math. \textbf{89} (1967), 92--135 (English. Russian original).

\bibitem{Dzhab69}
\bysame, \emph{Interpolation spaces of differentiable functions defined in an
  n-dimensional domain}, Sov. Math. Dokl. \textbf{10} (1969), no.~4, 833--837
  (English. Russian original).

\bibitem{Dzhab72Fam}
\bysame, \emph{Families of spaces of functions whose mixed derivatives satisfy
  a multiple-integral {H}older condition}, Proc. Steklov Inst. Math.
  \textbf{117} (1972), 165--188 ({E}nglish. {R}ussian original).

\bibitem{Dzhab74}
\bysame, \emph{Theorems on the extension of functions in the spaces
  ${S}_{p}^{\mathbf{r}}{W}$ and ${S}_{p,\theta}^{\mathbf{r}}{B}$ beyond the
  boundary of a domain}, Proc. Steklov Inst. Math. \textbf{131} (1974), 85--97
  (English. Russian original).

\bibitem{Dzhab99}
\bysame, \emph{The method of integral representation in the theory of spaces of
  functions of several groups of variables}, Recent Developments in Complex
  Analysis and Computer Algebra (R.~Gilbert, J.~Kajiwara, and Y.~Xu, eds.),
  ISAAC 2001, Kluwer Academic Publishers, Dordrecht, Netherlands, 1999,
  pp.~69--80.

\bibitem{Dwight61}
H.~B. Dwight, \emph{Tables of integrals and other mathematical data}, 4th ed.,
  Macmillan, 1961.

\bibitem{Dyn89}
N.~Dyn, \emph{Interpolation and approximation by radial and related functions},
  Approximation Theory VI (C.~Chui., L.~Schumaker., and J.~Ward, eds.), Acad.
  Press, 1989, pp.~211--234.

\bibitem{HarLitPol52}
G.~Polya G.~Hardy, J.~Littlewood, \emph{Inequalities}, 2nd ed., Univ. Press,
  Cambridge, 1952.

\bibitem{GarckGrieb05}
J.~Garke and M.~Griebel, \emph{Semi-supervised learning with sparse grids},
  Proceedings of the 22nd International Conference on Machine Learning, 2005.

\bibitem{GradRyz07}
I.~S. Gradshteyn. and I.~M. Ryzhik, \emph{Tables of integrals, series and
  products}, 7th ed., Acad. Press, New York, 2007.

\bibitem{Grisvard85}
P.~Grisvard, \emph{Elliptic problems in nonsmooth domains}, Monographs and
  studies in mathematics, vol.~24, Pitman Advanced Pub. Program, 1985.

\bibitem{Hajl96}
P.~Hajlasz, \emph{Sobolev spaces on an arbitrary metric space}, Potential Anal.
  (1996), no.~5, 403--415.

\bibitem{Hajl2003}
\bysame, \emph{A new characterization of the {S}obolev space}, Studia Math.
  \textbf{159} (2003), 263--275.

\bibitem{HansVyb2009}
M.~Hansen and J.~Vybiral, \emph{The {J}awerth-{F}ranke embedding of spaces with
  dominating mixed smoothness}, Georgian Math. J. \textbf{16} (2009), no.~4,
  667--682.

\bibitem{Hollig2003}
K.~Hollig, \emph{Finite element methods with b-splines}, Frontiers in applied
  mathematics, vol.~3, SIAM, Philadelphia, PA, 2003.

\bibitem{Hormand63}
L.~H\"{o}rmander, \emph{Linear partial differential operators}, Grundlehren der
  mathematischen Wissenschaften, vol. 116, Springer, 1969.

\bibitem{Ilin63A}
V.~P. Il'in, \emph{On inequalities between norms of partial derivatives of
  functions of many variables}, Soviet Math. Dokl. (1963), no.~3, 784--787
  (English. Russian original).

\bibitem{Ilin68}
\bysame, \emph{Conditions of validity of inequalities between ${L}_{p}$ norms
  of partial derivatives of functions of several variables}, Proc. Steklov
  Inst. Math. \textbf{96} (1968), 259--305 (English. Russian original).

\bibitem{Jones2011}
F.~Jones, \emph{Lebesgue integration on euclidean space}, Pure and Applied
  Mathematics, Jones and Bartlett Mathematics, 2001.

\bibitem{Krantz80}
S.~G. Krantz, \emph{Review: "{I}ntegral representations of functions and
  imbedding theorems" by {B}esov, {I}l'in, {N}ikol'ski\u{i}}, Bull. Amer. Math.
  Soc. (N.S.) \textbf{2} (1980), no.~1, 216--222.

\bibitem{Kuller1969}
R.~G. Kuller, \emph{Topics in modern analysis}, Prentice-Hall, Englewood
  Cliffs, NJ, 1969.

\bibitem{LightWayne95ErrEst}
W.~Light and H.~Wayne, \emph{Error estimates for approximation by radial basis
  functions.}, Approximation theory, wavelets and applications (Maratea, 1994),
  NATO Adv. Sci. Inst. Ser. C Math. Phys. Sci., vol. 454, Kluwer Acad. Publ.,
  Dordrecht, 1995, pp.~215--246.

\bibitem{LightWayne98PowFunc}
\bysame, \emph{On power functions and error estimates for radial basis function
  interpolation}, J. Approx. Th. \textbf{92} (1998), no.~2, 245--266.

\bibitem{LightWayneX98Weight}
\bysame, \emph{Spaces of distributions, interpolation by translates of a basis
  function and error estimates}, Numer. Math. \textbf{81} (1999), no.~3,
  415--450.

\bibitem{LizNikol6567}
P.I. Lizorkin and S.M. Nikol'ski\u{i}, \emph{Classification of differentiable
  functions on the basis of spaces with a dominant derivative}, Proc. Steklov.
  Inst. Math. \textbf{77} (1967), 160--187 (English. Russian original).

\bibitem{MagOberSoni66}
W.~Magnus, F.~Oberhettinger, and R.~Soni, \emph{Formulas and theorems for the
  special functions of mathematical physics}, 3rd ed., Springer-Verlag, Berlin,
  Heidelberg, New York, 1966.

\bibitem{Malliavin95}
P.~Malliavin, H.~Airault, L.~Kay, and G.~Letac, \emph{Integration and
  probability}, Springer-Verlag, Berlin, 1995.

\bibitem{NarcWardWend2004}
F.~Narcowich, J.~Ward, and H.~Wendland, \emph{Sobolev bounds on functions with
  scattered zeros, with applications to radial basis function surface fitting},
  Math. Comput. \textbf{74} (2004), no.~250, 743--763.

\bibitem{Nikol62}
S.M. Nikol'ski\u{i}, \emph{Properties of differentiable functions of several
  variables at the boundary of their domain of definition}, Sov. Math. Dokl.
  \textbf{62} (1962), 1357--1360 (English. Russian original).

\bibitem{Nikol63}
\bysame, \emph{On stable boundary values of differentiable functions of several
  variables}, Mat. Sb. \textbf{61} (1963), no.~2, 224--252 (Russian).

\bibitem{Nikol6373}
\bysame, \emph{Functions with dominant mixed derivatives satisfying a multiple
  {H}older condition}, Amer. Math. Soc. Transl. \textbf{102} (1973), 27--51
  (English. Russian original).

\bibitem{Nikol75}
\bysame, \emph{Approximation of functions of several variables and imbedding
  theorems}, Grundlehren der mathematischen Wissenschaften, vol. 205,
  Springer-Verlag, Berlin, New York, 1975.

\bibitem{Packel74}
E.~W. Packel, \emph{Functional analysis - a short course}, Intext series in
  mathematics, Intext Educational Publishers, New York, London, 1974.

\bibitem{Petersen83}
B.~Petersen, \emph{Introduction to the fourier transform and
  pseudo-differential operators}, Monographs and studies in mathematics,
  vol.~19, Pitman Advanced Pub. Program, Boston, 1983.

\bibitem{ReedSimon72}
M.~Reed and B.~Simon, \emph{Methods of modern mathematical physics i:
  Functional analysis}, Academic Press, New York, 1972.

\bibitem{Schmeis2006}
H.~J. Schmeisser, \emph{Recent developments in the theory of function spaces
  with dominating mixed smoothness (english)}, Nonlin. Anal., Func. Spaces and
  Applic. \textbf{8} (2006), 145--204.

\bibitem{SchmeisTrieb87}
H.~J. Schmeisser and H.~Triebel, \emph{Topics in {F}ourier analysis and
  function spaces}, John Wiley, Chichester, 1987.

\bibitem{Schwartz09}
C.~Schwartz, \emph{Elementary functional analysis}, World Scientific Pub., New
  Jersey, 2009.

\bibitem{SteinWeiss71}
E.~M. Stein and G.~Weiss, \emph{Fourier analysis and {E}uclidean spaces},
  Princeton University Press, Princeton, 1971.

\bibitem{Trev75}
F.~Treves, \emph{Basic linear partial differential equations}, Pure and Applied
  Mathematics, vol.~62, Academic Press, 1975.

\bibitem{Vladimirov}
V.~Vladimirov, \emph{Equations of mathematical physics}, Pure and Applied
  Mathematics series, vol.~3, Marcel Dekker, New York, 1971.

\bibitem{Watson95}
G.~N. Watson, \emph{A treatise on the theory of {B}essel functions}, Cambridge
  Univ. Press, Cambridge, New York, Melbourne, 1995.

\bibitem{Wendland05}
H.~Wendland, \emph{Scattered data approximation}, Cambridge Monographs on
  Computational and Applied Mathematics, vol.~17, Cambridge Univ. Press,
  Cambridge, 2005.

\bibitem{WilliamsPosOrdSmthV3}
P.~Y. Williams, \emph{A weight function theory of positive order basis function
  interpolants and smoothers}, unpublished arXiv:0708.0795v3 [math.NA], 2014.

\bibitem{Wloka87}
J.~Wloka, \emph{Partial differential equations}, Cambridge Univ. Press,
  Cambridge, 1987.

\bibitem{Yosida58}
K.~Yosida, \emph{Functional analysis}, Springer-Verlag, Berlin, 1958.

\bibitem{Zwill89}
D.~Zwillinger, \emph{Handbook of differential equations}, Academic Press, 1989.

\bibitem{Zwill92}
\bysame, \emph{Handbook of integration}, Jones and Bartlett Publishers, London,
  Boston, 1992.

\end{thebibliography}

\end{document}